\documentclass[letterpaper]{article}
\usepackage{exscale,relsize}
\usepackage[dvips]{graphicx,color}
\usepackage{paralist}
\usepackage{pst-3d}
\usepackage{pst-all}
\usepackage{pst-text}
\usepackage{pst-coil}
\usepackage{pst-eps}
\usepackage{pst-fill}
\usepackage{pst-grad}
\usepackage{pst-node}
\usepackage{pst-plot}
\usepackage{pst-text}
\usepackage{pst-tree}
\usepackage{fancybox}

\usepackage{amsmath}
\usepackage{amssymb}
\usepackage{bm}
\usepackage[mathscr]{eucal}
\usepackage{amsfonts}
\usepackage{latexsym}
\usepackage{amsxtra}
\usepackage{amsbsy}
\usepackage{amsthm}
\usepackage{amscd}
\usepackage{amsopn}
\usepackage{amstext}
\usepackage{upref}
\usepackage{oldgerm}
\usepackage[varumlaut]{yfonts}
\newcounter{z0}

\newtheorem{dddd}{Theorem}[subsection]

\newtheorem{bbbbb}{Proposition}[section]
\newtheorem{ccccc}{Lemma}[section]

\newtheorem{fffff}{Corollary}[section]

\newtheorem*{rhp2}{$\bm{\mathrm{Lemma}} \, \, 
\bm{\mathrm{RHP}_{\mathrm{MPC}}}$}

\theoremstyle{definition}

\theoremstyle{definition}
\newtheorem{aaaa}{Definition}[subsection]

\theoremstyle{definition}
\newtheorem{aaaaa}{Definition}[section]

\theoremstyle{definition}

\theoremstyle{definition}
\newtheorem{eeee}{Remark}[subsection]

\theoremstyle{definition}
\newtheorem{eeeee}{Remark}[section]

\theoremstyle{definition}

\newtheorem{notetolem}{Remark to Lemma}[section]

\newcommand{\me}{\mathrm{e}}
\newcommand{\mi}{\mathrm{i}}
\newcommand{\md}{\mathrm{d}}
\renewcommand{\Im}{\mathrm{Im}}
\renewcommand{\Re}{\mathrm{Re}}
\newcommand{\id}{\pmb{\mi \md}}
\newcommand{\norm}[1]{\lVert#1\rVert}
\allowdisplaybreaks[4]

\newcommand{\pvi}{\ensuremath{\int \hspace{-2.75mm} \rule[2.5pt]{2mm}{0.25mm}}}

\usepackage{stmaryrd}
\usepackage{txfonts}
\DeclareMathOperator*{\supp}{supp}

\DeclareMathOperator*{\ad}{ad}
\let \Im \relax
\let \Re \relax
\DeclareMathOperator{\Im}{Im}
\DeclareMathOperator{\Re}{Re}
\usepackage{epic}
\usepackage{eepic}
\definecolor{pink}{rgb}{1,0.75,0.8}
\definecolor{lightblue}{rgb}{0.68,0.85,0.9}
\begin{document}
\fontsize{10pt}{11pt}\selectfont
\fontencoding{T1}\selectfont
\baselineskip=11pt
\frenchspacing
\title{Riemann-Hilbert Characterisation of Rational Functions with a General 
Distribution of Poles on the Extended Real Line Orthogonal with Respect to 
Varying Exponential Weights: Multi-Point Pad\'{e} Approximants and Asymptotics}
\author{A.~Vartanian \and X.~Zhou}
\date{\today}
\maketitle
\begin{abstract}
\noindent
Given $K$ poles $\alpha_{1},\alpha_{2},\dotsc,\alpha_{K}$ on the extended real 
line, $\overline{\mathbb{R}}$, which are neither necessarily distinct nor bounded, 
denote by $\Lambda^{\mathbb{R}}$ the linear space over $\mathbb{R}$ spanned 
by a constant and the rational functions $\mathscr{S}^{n}_{k}(z)$, $(n,k) \! \in \! 
\mathbb{N} \times \lbrace 1,2,\dotsc,K \rbrace$, where $\mathscr{S}^{n}_{k}(z) 
\! = \! z^{\varkappa_{nk}}$ if $\alpha_{k} \! = \! \infty$ and $\mathscr{S}^{n}_{k}
(z) \! = \! (z \! - \! \alpha_{k})^{-\varkappa_{nk}}$ if $\alpha_{k} \! \neq \! \infty$, 
with $\varkappa_{nk} \! = \! (n \! - \! 1) \gamma_{k} \! + \! \varrho_{k}$, where 
$\gamma_{k}$ counts the number of repetitions of the pole $\alpha_{k}$ in 
the (entire) set $\lbrace \alpha_{1},\alpha_{2},\dotsc,\alpha_{K} \rbrace$, and 
$\varrho_{k}$ enumerates the number of times the pole $\alpha_{k}$ appears 
up to and including the $k$th position in the set $\lbrace \alpha_{1},\alpha_{2},
\dotsc,\alpha_{K} \rbrace$. Define the real inner product $\left\langle 
\boldsymbol{\cdot},\boldsymbol{\cdot} \right\rangle_{\mathscr{L}} \colon 
\Lambda^{\mathbb{R}} \times \Lambda^{\mathbb{R}} \! \to \! \mathbb{R}$, 
$(f,g) \! \mapsto \! \left\langle f,g \right\rangle_{\mathscr{L}} \! := \! 
\int_{\mathbb{R}}f(\xi)g(\xi) \exp (-\mathscr{N}V(\xi)) \, \md \xi$, 
$\mathscr{N} \! \in \! \mathbb{N}$, where the external field, $V$, is real 
analytic on $\overline{\mathbb{R}} \setminus \lbrace \alpha_{1},\alpha_{2},
\dotsc,\alpha_{K} \rbrace$ and satisfies the growth conditions $\lim_{x \to 
\alpha_{i}}(V(x)/\ln (x^{2} \! + \! 1)) \! = \! +\infty$, $i \! \in \! \lbrace 
\mathstrut k \! \in \! \lbrace 1,2,\dotsc,K \rbrace; \, \alpha_{k} \! = \! \infty 
\rbrace$ and $\lim_{x \to \alpha_{j}}(V(x)/\ln ((x \! - \! \alpha_{j})^{-2} \! + 
\! 1)) \! = \! +\infty$, $j \! \in \! \lbrace \mathstrut k \! \in \! \lbrace 1,2,
\dotsc,K \rbrace; \, \alpha_{k} \! \neq \! \infty \rbrace$. For $(n,k) \! \in 
\! \mathbb{N} \times \lbrace 1,2,\dotsc,K \rbrace$, corresponding to 
the repeated pole sequence $\lbrace \text{no pole} \rbrace \cup 
\lbrace \alpha_{1},\alpha_{2},\dotsc,\alpha_{K},\alpha_{1},\alpha_{2},
\dotsc,\alpha_{K},\dotsc,\alpha_{1},\alpha_{2},\linebreak[4]
\dotsc,\alpha_{k} \rbrace$, orthogonalisation of the associated ordered 
base of rational functions $\lbrace \text{const.} \rbrace \cup \lbrace 
\mathscr{S}^{1}_{1}(z),\mathscr{S}^{1}_{2}(z),\dotsc,\mathscr{S}^{1}_{K}(z),
\linebreak[4]
\mathscr{S}^{2}_{1}(z),\mathscr{S}^{2}_{2}(z),\dotsc,\mathscr{S}^{2}_{K}(z),
\dotsc,\mathscr{S}^{n}_{1}(z),\mathscr{S}^{n}_{2}(z),\dotsc,\mathscr{S}^{n}_{k}(z) 
\rbrace$ with respect to $\left\langle \boldsymbol{\cdot},\boldsymbol{\cdot} 
\right\rangle_{\mathscr{L}}$ gives rise to $K$ families of monic orthogonal 
rational functions (ORFs), denoted by $\lbrace 1 \rbrace \cup \lbrace 
\pmb{\pi}_{k}^{n}(z) \rbrace_{\underset{k \in \lbrace 1,2,\dotsc,K \rbrace}{n 
\in \mathbb{N}}}$, where one subfamily consists of $\hat{K} \! := \! \# \lbrace 
\mathstrut k \! \in \! \lbrace 1,2,\dotsc,K \rbrace; \, \alpha_{k} \! = \! \infty 
\rbrace$ subfamilies of ORFs, and another subfamily consists of $\tilde{K} 
\! := \! \# \lbrace \mathstrut k \! \in \! \lbrace 1,2,\dotsc,K \rbrace; \, 
\alpha_{k} \! \neq \! \infty \rbrace$ subfamilies of ORFs, with $K \! = \! \hat{K} 
\! + \! \tilde{K}$. The monic ORF problem is characterised as a family of $K$ 
matrix Riemann-Hilbert problems (RHPs) on $\overline{\mathbb{R}}$, and 
a corresponding family of $K$ energy minimisation (variational) problems 
containing external fields with singular points is presented, and the existence, 
uniqueness, and regularity properties of the associated family of equilibrium 
measures is established. The family of $K$ equilibrium measures, in conjunction 
with a corresponding family of complex potentials $(g$-functions), is used to 
derive a family of $K$ `model' matrix RHPs on $\overline{\mathbb{R}}$ which 
are amenable to asymptotic analysis, in the double-scaling limit $\mathscr{N},
n \! \to \! \infty$ such that $\mathscr{N}/n \! \to \! 1$, via the Deift-Zhou 
non-linear steepest-descent method: this is used to derive uniform asymptotics 
of the monic ORF, $\pmb{\pi}_{k}^{n}(z)$, its leading coefficient, and the ORF, 
denoted by $\phi_{k}^{n}(z)$, as well as related, important objects (zeros of ORFs, 
normalised zero counting measures, etc.), in the entire complex plane. Furthermore, 
for $(n,k) \! \in \! \mathbb{N} \times \lbrace 1,2,\dotsc,K \rbrace$, a family 
of $K$ multi-point Pad\'{e} approximants (MPAs) of type $((n \! - \! 1)K \! 
+ \! k \! - \! 1,(n \! - \! 1)K \! + \! k)$ for the Markov-Stieltjes transform, 
$\int_{\mathbb{R}}(z \! - \! \xi)^{-1} \exp (-\mathscr{N}V(\xi)) \, \md \xi$, is 
introduced, and uniform asymptotics, in the double-scaling limit $\mathscr{N},n 
\! \to \! \infty$ such that $\mathscr{N}/n \! \to \! 1$, are obtained for the 
corresponding MPAs and their associated errors in approximation (MPA error 
terms) in the entire complex plane.

\vspace{0.35cm}

30E15, 30E20, 30E25

\vspace{0.20cm}
\textbf{Abbreviated Title.} Orthogonal Rational Functions: Characterisation, 
Multi-Point Pad\'{e} Approximants, 

and Asymptotics

\vspace{0.20cm}

Equilibrium Measures, Asymptotics
\end{abstract}
\clearpage
\tableofcontents
\listoffigures
\clearpage
\section{Introduction} \label{intro} 
\subsection{Orthogonal Rational Functions and Multi-Point Pad\'{e} 
Approximants: Motivation and Background} \label{subsec1.1} 
In approximation problems, whether complex or real, there are typically 
three seminal questions which need to be addressed: (Q1) for which class 
of functions is/are the proposed approximation scheme(s)/method(s) 
(approximant(s)) applicable?; (Q2) with what degree of accuracy, 
theoretical or numerical, is/are the proposed approximant(s) 
constructed?; and (Q3) what is/are the rate(s) of convergence(s)? With 
respect to question (Q1), consider the following problem: given $K$ 
$(\in \! \mathbb{N}$ and finite) poles $\alpha_{1},\alpha_{2},\dotsc,
\alpha_{K}$, which are neither necessarily distinct nor bounded, on 
the extended real line, $\overline{\mathbb{R}} \! := \! \mathbb{R} 
\cup \lbrace \pm \infty \rbrace$, and the positive probability measure 
with varying exponential weight, $\md \widetilde{\mu}(z) \! = \! 
\exp (-n \widetilde{V}(z)) \, \md z$, $n \! \in \! \mathbb{N}$, with 
$\widetilde{V}(z) \! := \! z_{o}V(z)$, where $z_{o} \! = \! z_{o}
(\mathscr{N},n) \! > \! 0$,\footnote{The explicit formula for $z_{o}$ 
is not imperative at this stage; rather, the relevant fact is that it is a 
positive, real-valued function of the two natural numbers $\mathscr{N}$ 
and $n$.} $\mathscr{N} \! \in \! \mathbb{N}$, and $V$ satisfies 
(see Subsection~\ref{subsec1.2}) the analyticity and growth 
conditions~\eqref{eq3}--\eqref{eq5}, what are asymptotics, in the 
double-scaling limit $\mathscr{N},n \! \to \! \infty$ such that $z_{o} 
\! = \! 1 \! + \! o(1)$, for the class of functions (the so-called 
Markov-Stieltjes transform of the probability measure 
$\widetilde{\mu})$ with representation $\mathrm{F}_{\tilde{\mu}}(z) 
\! := \! \int_{\mathbb{R}}(z \! - \! \xi)^{-1} \, \md \widetilde{\mu}
(\xi)$ for $z \! \in \! \mathbb{C} \setminus \mathbb{R}$? It turns 
out that a set of objects colloquially referred to as Orthogonal Rational 
Functions (ORFs) is the key analytical tool needed in order to solve 
the above-mentioned `asymptotic problem'; in fact, this monograph addresses 
a number of questions related to ORFs, precise definitions of which are given 
in Subsection~\ref{subsec1.2}. ORFs may be thought of as a generalisation 
of Orthogonal Polynomials (OPs); and, just as OPs are obtained via the 
Gram-Schmidt orthogonalisation procedure (with respect to a given measure) 
applied to the monomial sequence $\lbrace 1,z,z^{2},\dotsc \rbrace$, ORFs, 
too, may be obtained via the Gram-Schmidt procedure, but applied to a more 
general, pre-assigned sequence of (basis) rational functions.

In order to elucidate, further, the notion of ORFs, it is pedagogically advantageous 
to present, for illustrative purposes, a basic, albeit contrived, example. Using, 
say, the following measure defined on $\overline{\mathbb{R}}$, $\md 
\mathfrak{m}(z) \! = \! \me^{-z^{-2}} \me^{-(z-1)^{-2}} \me^{-z^{2}} 
\me^{-(z-\sqrt{2} \,)^{-2}} \me^{-(z-\pi)^{-2}} \, \md z$, one may apply 
the Gram-Schmidt procedure (with respect to this measure) to an ordered, 
initial sequence of---simple---rational functions with poles at $0,1,\infty,
\sqrt{2}$, and $\pi$. An infinite sequence of rational functions will now 
be constructed by first specifying the (finite) pole sequence of `length' $K \! 
= \! 7$, namely, $\lbrace \text{no pole} \rbrace \cup \lbrace 0,1,\infty,1,
\sqrt{2},\pi,\infty \rbrace$. Two running indices, $n \! \in \! \mathbb{N}$ and 
$k \! \in \! \lbrace 1,2,\dotsc,K \rbrace$, are now described: $n$ indexes 
the number of times the pole sequence is repeated, and $k$ indexes where 
in the pole sequence the pole $\alpha_{k}$ resides; e.g., for the $K \! = \! 7$ 
pole sequence above, that is, $\alpha_{1} \! = \! 0$, $\alpha_{2} \! = \! 1$, 
$\alpha_{3} \! = \! \infty$, $\alpha_{4} \! = \! 1$, $\alpha_{5} \! = \! \sqrt{2}$, 
$\alpha_{6} \! = \! \pi$, and $\alpha_{7} \! = \! \infty$, $n \! = \! 1$ corresponds 
to the pole sequence $\lbrace 0,1,\infty,1,\sqrt{2},\pi,\infty \rbrace$, $n \! 
= \! 2$ corresponds to the (repeated) pole sequence $\lbrace 0,1,\infty,
1,\sqrt{2},\pi,\infty \rbrace \cup \lbrace 0,1,\infty,1,\sqrt{2},\pi,\infty 
\rbrace$, etc. The first appearance of the pole sequence (corresponding 
to $n \! = \! 1)$ $\lbrace \text{no pole} \rbrace \cup \lbrace 0,1,\infty,
1,\sqrt{2},\pi,\infty \rbrace$ gives rise to the eight rational functions 
$\lbrace 1 \rbrace \cup \lbrace z^{-1},(z \! - \! 1)^{-1},z,(z \! - \! 1)^{-2},
(z \! - \! \sqrt{2} \,)^{-1},(z \! - \! \pi)^{-1},z^{2} \rbrace$. The infinite 
sequence of rational functions is now specified by repeating the pole sequence, 
and incrementing, by one, the order of the pole at $0,1,\infty,\sqrt{2}$, and 
$\pi$ every time that pole is encountered; e.g., after the first eight rational 
functions specified above, one repeats the pole sequence again (corresponding 
to $n \! = \! 2)$ so that the next seven members of the infinite sequence are 
$\lbrace z^{-2},(z \! - \! 1)^{-3},z^{3},(z \! - \! 1)^{-4},(z \! - \! \sqrt{2} \,)^{-2},
(z \! - \! \pi)^{-2},z^{4} \rbrace$, then, after specifying the latter fifteen rational 
functions, one repeats, again, the pole sequence (corresponding to $n \! = \! 3)$ 
so that the next seven members of the infinite sequence are $\lbrace z^{-3},
(z \! - \! 1)^{-5},z^{5},(z \! - \! 1)^{-6},(z \! - \! \sqrt{2} \,)^{-3},(z \! - \! \pi)^{-3},
z^{6} \rbrace$, etc.,\footnote{Of course, the complete, cyclically repeated $n \! 
= \! 3$ pole sequence is $\lbrace \text{no pole} \rbrace \cup \lbrace 0,1,\infty,
1,\sqrt{2},\pi,\infty \rbrace \cup \lbrace 0,1,\infty,1,\sqrt{2},\pi,\infty \rbrace 
\cup \lbrace 0,1,\infty,1,\sqrt{2},\pi,\infty \rbrace$, and the corresponding 
sequence of rational functions reads $\lbrace 1 \rbrace \cup \lbrace z^{-1},
(z \! - \! 1)^{-1},z,(z \! - \! 1)^{-2},(z \! - \! \sqrt{2} \,)^{-1},(z \! - \! \pi)^{-1},
z^{2} \rbrace \cup \lbrace z^{-2},(z \! - \! 1)^{-3},z^{3},(z \! - \! 1)^{-4},
(z \! - \! \sqrt{2} \,)^{-2},(z \! - \! \pi)^{-2},z^{4} \rbrace \cup \lbrace z^{-3},
(z \! - \! 1)^{-5},z^{5},(z \! - \! 1)^{-6},(z \! - \! \sqrt{2} \,)^{-3},(z \! - \! \pi)^{-3},
z^{6} \rbrace$.} and this continues $(n \! = \! 4,5,6,\dotsc)$, with the rational 
functions appearing in seven-tuples. The result of applying the Gram-Schmidt 
procedure to this ordered sequence of rational functions, using the aforementioned 
measure $\md \mathfrak{m}$, is a sequence of ORFs; more precisely, the 
$((n \! - \! 1)K \! + \! k)$th term in the new (ORF) sequence is in the linear span 
of the first $(n \! - \! 1)K \! + \! k$ members of the original sequence of rational 
functions and is orthogonal to all `previous' members of the sequence.

The general definition of a sequence of ORFs is particularly cumbersome: 
the reason being that one must specify the ordered pole sequence, and 
if a pole should happen to repeat, then the associated `simple rational 
function' must be linearly independent from all previous members of 
the ordered sequence; therefore, typically, the order of the pole 
increases with each occurrence of that pole. Although one might merely 
describe the sequence of rational functions abstractly, it is useful, and 
indeed necessary for a Riemann-Hilbert characterisation of the ORF 
problem, to explicitly annotate this level of detail. This monograph deals 
exclusively with the case of a finite sequence of poles on $\overline{
\mathbb{R}}$ which are repeated cyclically: \textbf{the situation in which 
the pole sequence is infinite is not considered in this monograph} (for 
this latter case, see the monograph \cite{a15} for complete details). In 
Subsection~\ref{subsec1.2}, then, the reader will find a detailed exposition 
of the repeating pole sequence: it is imperative to note that the family of 
Riemann-Hilbert problems (RHPs) which characterise the ORFs described 
herein depend on this level of detail regarding the ordered pole sequence. 
More abstractly, then, and without dwelving into technical particulars, 
given $K$ $(\in \! \mathbb{N}$ and finite) poles $\alpha_{1},\alpha_{2},
\dotsc,\alpha_{K}$ on $\overline{\mathbb{R}}$, which are neither necessarily 
distinct nor bounded, denote by $\Lambda^{\mathbb{R}}$ the linear space 
over $\mathbb{R}$ spanned by a constant and the rational functions 
$\mathscr{S}^{n}_{k}(z)$, $(n,k) \! \in \! \mathbb{N} \times \lbrace 1,2,
\dotsc,K \rbrace$, where $\mathscr{S}^{n}_{k}(z) \! = \! z^{\varkappa_{nk}}$ 
if $\alpha_{k} \! = \! \infty$ and $\mathscr{S}^{n}_{k}(z) \! = \! (z \! - \! 
\alpha_{k})^{-\varkappa_{nk}}$ if $\alpha_{k} \! \neq \! \infty$, with 
$\varkappa_{nk} \colon \mathbb{N} \times \lbrace 1,2,\dotsc,K \rbrace \! 
\to \! \mathbb{N}$, $(n,k) \! \mapsto \! (n \! - \! 1) \gamma_{k} \! + \! 
\varrho_{k} \! =: \! \varkappa_{nk}$, where $\gamma_{k}$ counts the 
number of repetitions of the pole $\alpha_{k}$ in the (entire) set $\lbrace 
\alpha_{1},\alpha_{2},\dotsc,\alpha_{K} \rbrace$, and $\varrho_{k}$ 
enumerates the number of times the pole $\alpha_{k}$ appears up to and 
including the $k$th position in the set $\lbrace \alpha_{1},\alpha_{2},
\dotsc,\alpha_{K} \rbrace$; e.g., for the $K \! = \! 7$ pole set $\lbrace 
\alpha_{1},\alpha_{2},\alpha_{3},\alpha_{4},\alpha_{5},\alpha_{6},\alpha_{7} 
\rbrace \! = \! \lbrace 0,1,\infty,1,\sqrt{2},\pi,\infty \rbrace$, $\gamma_{1} 
\! = \! 1$, $\gamma_{2} \! = \! 2$, $\gamma_{3} \! = \! 2$, $\gamma_{4} 
\! = \! 2$, $\gamma_{5} \! = \! 1$, $\gamma_{6} \! = \! 1$, $\gamma_{7} 
\! = \! 2$, $\varrho_{1} \! = \! 1$, $\varrho_{2} \! = \! 1$, $\varrho_{3} \! = 
\! 1$, $\varrho_{4} \! = \! 2$, $\varrho_{5} \! = \! 1$, $\varrho_{6} \! = \! 1$, 
$\varrho_{7} \! = \! 2$, and $\varkappa_{n1} \! = \! n$, $\varkappa_{n2} 
\! = \! 2n \! - \! 1$, $\varkappa_{n3} \! = \! 2n \! - \! 1$, $\varkappa_{n4} 
\! = \! 2n$, $\varkappa_{n5} \! = \! n$, $\varkappa_{n6} \! = \! n$, and 
$\varkappa_{n7} \! = \! 2n$. Define the real inner product $\left\langle 
\boldsymbol{\cdot},\boldsymbol{\cdot} \right\rangle_{\mathscr{L}} \colon 
\Lambda^{\mathbb{R}} \times \Lambda^{\mathbb{R}} \! \to \! \mathbb{R}$, 
$(f,g) \! \mapsto \! \left\langle f,g \right\rangle_{\mathscr{L}} \! := \! 
\int_{\mathbb{R}}f(\xi)g(\xi) \me^{-n \widetilde{V}(\xi)} \, \md \xi$, $n \! \in 
\! \mathbb{N}$, where the `scaled' external field, $\widetilde{V}$, satisfies 
(see Subsection~\ref{subsub2}) the analyticity and growth 
conditions~\eqref{eq20}--\eqref{eq22}. Corresponding to the repeated pole 
sequence $\lbrace \text{no pole} \rbrace \cup \lbrace \alpha_{1},\alpha_{2},
\dotsc,\alpha_{K},\alpha_{1},\alpha_{2},\dotsc,\alpha_{K},\dotsc,\alpha_{1},
\alpha_{2},\dotsc,\alpha_{k} \rbrace$, orthogonalisation of the associated, 
ordered base of rational functions $\lbrace \text{const.} \rbrace \cup \lbrace 
\mathscr{S}^{1}_{1}(z),\mathscr{S}^{1}_{2}(z),\dotsc,\mathscr{S}^{1}_{K}(z),
\mathscr{S}^{2}_{1}(z),\mathscr{S}^{2}_{2}(z),\dotsc,\mathscr{S}^{2}_{K}(z),
\dotsc,\mathscr{S}^{n}_{1}(z),\mathscr{S}^{n}_{2}(z),\dotsc,\mathscr{S}^{n}_{k}(z) 
\rbrace$ with respect to $\left\langle \boldsymbol{\cdot},\boldsymbol{\cdot} 
\right\rangle_{\mathscr{L}}$ gives rise to $K$ families of monic ORFs, denoted 
by $\lbrace 1 \rbrace \cup \lbrace \pmb{\pi}_{k}^{n}(z) \rbrace_{\underset{k 
\in \lbrace 1,2,\dotsc,K \rbrace}{n \in \mathbb{N}}}$, where one subfamily 
consists of $\hat{K} \! := \! \# \lbrace \mathstrut k \! \in \! \lbrace 1,2,
\dotsc,K \rbrace; \, \alpha_{k} \! = \! \infty \rbrace$ subfamilies of ORFs, 
and another subfamily consists of $\tilde{K} \! := \! \# \lbrace \mathstrut 
k \! \in \! \lbrace 1,2,\dotsc,K \rbrace; \, \alpha_{k} \! \neq \! \infty \rbrace$ 
subfamilies of ORFs, with $K \! = \! \hat{K} \! + \! \tilde{K}$. An important 
feature worth noting is that the Riemann-Hilbert characterisation of the 
(monic) ORF problem requires the consideration of two disjoint classes 
of RHPs: one class corresponding to those indices associated with the poles 
at the point at infinity, $\hat{\mathfrak{k}} \! := \! \lbrace \mathstrut k \! 
\in \! \lbrace 1,2,\dotsc,K \rbrace; \, \alpha_{k} \! = \! \infty \rbrace$, and 
another class corresponding to those indices associated with the bounded 
poles, $\tilde{\mathfrak{k}} \! := \! \lbrace \mathstrut k \! \in \! \lbrace 
1,2,\dotsc,K \rbrace; \, \alpha_{k} \! \neq \! \infty \rbrace$ (of course, 
$\hat{\mathfrak{k}} \cap \tilde{\mathfrak{k}} \! = \! \varnothing)$; e.g., for the 
$K \! = \! 7$ pole set $\lbrace \alpha_{1},\alpha_{2},\alpha_{3},\alpha_{4},
\alpha_{5},\alpha_{6},\alpha_{7} \rbrace \! = \! \lbrace 0,1,\infty,1,\sqrt{2},
\pi,\infty \rbrace$, $\hat{\mathfrak{k}} \! = \! \lbrace 3,7 \rbrace$ and 
$\tilde{\mathfrak{k}} \! = \! \lbrace 1,2,4,5,6 \rbrace$. Furthermore, even 
though the expression for the monic ORF is different for each member of 
the above-mentioned two classes, the generic notation $\pmb{\pi}_{k}^{n}(z)$ 
is retained throughout this monograph, irrespective of class.

With respect to question (Q2), consider the following discussion 
drawn {}from rational approximation with fixed poles (a general reference 
for the discussion that follows is the monograph \cite{gabpgm}). To 
approximate analytic/meromorphic functions in the neighbourhood(s) 
of their pole(s), Taylor series are not suitable approximants due to 
the fact that their partial sums are polynomials; therefore, Pad\'{e} 
approximants (PAs) were introduced: PAs are rational functions which 
interpolate a given function at a point of interest (typically the point 
at infinity) with maximal possible degree. Furthermore, there are more 
general multi-point Pad\'{e} approximants (MPAs), which, too, are 
rational functions, but which interpolate a given function, or its 
asymptotic development, at several points of interest simultaneously 
with maximal possible degrees. While PAs and MPAs are suitable for 
the uniform approximation of certain, special classes of functions, 
they can not, in general, guarantee convergence anywhere except 
at the point(s) of interpolation: this disadvantage of PAs and MPAs 
is due mainly to the fact that it is often not possible to control 
the spurious distribution of their poles (see, for example, 
\cite{mppta28,mppta6}; see, also, the recent  study \cite{scvnfsky} for a 
discussion of the generic properties of poles (and zeros) of Pad\'{e} 
approximants). A further impetus for the use of PAs is the following 
well-known problem: a classic theorem of Markov \cite{markovpa} states 
that the so-called `diagonal' PAs (constructed at the point at infinity) to 
a Markov-Stieltjes function/transform $\mathrm{u}_{\mathrm{MS}}(z) 
\! := \! \int_{\mathrm{D}_{\mathrm{MS}}}(z \! - \! \xi)^{-1} \, \md 
\mathfrak{u}(\xi)$, $z \! \in \! \overline{\mathbb{C}} \setminus 
\mathrm{D}_{\mathrm{MS}}$, where $\mathfrak{u}$ is a positive measure 
with compact support on $\mathbb{R}$, $\mathrm{D}_{\mathrm{MS}}$ is 
the minimal subset of $\mathbb{R}$ containing $\supp (\mathfrak{u})$, 
and $\overline{\mathbb{C}} \! := \! \mathbb{C} \cup \lbrace \infty \rbrace$ 
is the extended complex plane, converge locally uniformly to $\mathrm{u}_{
\mathrm{MS}}(z)$ for $z \! \in \! \overline{\mathbb{C}} \setminus 
\mathrm{D}_{\mathrm{MS}}$. Rakhmanov (see, for example, 
\cite{newearakh}) later showed that the analogue of Markov's classic 
theorem does not necessarily hold true for a function of the form 
(a rational perturbation of $\mathrm{u}_{\mathrm{MS}}(z))$ 
$\mathrm{f}_{\mathrm{MS}}(z) \! = \! \mathrm{u}_{\mathrm{MS}}(z) \! 
+ \! \mathrm{r}_{\mathrm{MS}}(z)$, where $\mathrm{r}_{\mathrm{MS}}
(z)$ is a rational function, because one does not, in general, have 
local uniform convergence for $z \! \in \! \overline{\mathbb{C}} 
\setminus (\mathrm{D}_{\mathrm{MS}} \cup \lbrace \mathstrut 
z^{\prime} \! \in \! \mathbb{C}; \, \mathrm{r}_{\mathrm{MS}}
(z^{\prime}) \! = \! \infty \rbrace)$ (unless `severe' restrictions 
are imposed on the measure $\mathfrak{u}$ or on the rational function 
$\mathrm{r}_{\mathrm{MS}}(z))$, due to the fact that the poles of the 
corresponding PAs may cluster in the domain of analyticity of 
$\mathrm{f}_{\mathrm{MS}}(z)$. Whilst the following topic is not addressed 
in this monograph, it is worth mentioning that, in an attempt to deal with 
the non-local-uniform-convergence (spurious pole distribution) issue  of 
the corresponding PAs or MPAs, a type of rational approximant(s) with 
pre-assigned poles, known as Pad\'{e}-type approximants (PTAs), were 
introduced (see, for example, \cite{mppta34,mppta33,mppta37,mppta36,
mppta7,mppta14,mppta8,mppta9,mppta10,mppta11,mppta13,mppta12,
mppta15,mppta20,mppta16}; see, also, \cite{eaovyevgen} for further 
generalisations and conjectures).

In order to elucidate, further, the notion of MPAs, it is pedagogically 
advantageous to present, for illustrative purposes, an extension of the 
contrived example considered earlier. Suppose, for example, that one 
wants to find an MPA to the Markov-Stieltjes function/transform of the 
measure $\md \mathfrak{m}(z) \! = \! \me^{-z^{-2}} \me^{-(z-1)^{-2}} 
\me^{-z^{2}} \me^{-(z-\sqrt{2} \,)^{-2}} \me^{-(z-\pi)^{-2}} \, \md z$, 
corresponding to the $K \! = \! 7$ pole set $\lbrace \alpha_{1},\alpha_{2},
\alpha_{3},\alpha_{4},\alpha_{5},\alpha_{6},\alpha_{7} \rbrace \! = \! \lbrace 
0,1,\infty,1,\sqrt{2},\pi,\infty \rbrace$, that is, $\mathrm{u}_{\mathrm{MS}}
(z)  \! = \! \int_{\mathbb{R}}(z \! - \! \xi)^{-1} \, \md \mathfrak{m}(\xi)$, 
$z \! \in \! \mathbb{C} \setminus \mathbb{R}$. Now, in open neighbourhoods 
of each of the points $z \! = \! 0,1,\sqrt{2},\pi$, and the point at infinity, 
$\mathrm{u}_{\mathrm{MS}}(z)$ possesses a complete asymptotic expansion 
(with respect to the indeterminate $z)$. For $n \! \in \! \mathbb{N}$ and, say, 
$k \! = \! 3$ (that is, $\alpha_{3} \! = \! \infty)$, the ($5$-point) MPA of type 
$((n \! - 1)K \! + \! k \! - \! 1,(n \! - \! 1)K \! + \! k) \! = \! (7n \! - \! 5,
7n \! - \! 4)$ for $\mathrm{u}_{\mathrm{MS}}(z)$ is a---proper---rational 
function of the form $\mathcal{U}(z)(\mathcal{V}(z))^{-1}$, where 
$\mathcal{U}(z)$ is a polynomial of the form $\sum_{m=0}^{7n-5}u_{m}
z^{m}$ with $\deg (\mathcal{U}) \! = \! 7n \! - \! 5$, $\mathcal{V}(z)$ 
$(\not\equiv \! 0)$ is a monic polynomial of the form $\sum_{m=0}^{7n-4}
v_{m}z^{m}$, $v_{7n-4} \! = \! 1$, with $\deg (\mathcal{V}) \! = \! 7n \! - 
\! 4$, and $(\mathcal{U}(z),\mathcal{V}(z))$ coprime,\footnote{Two polynomials 
are coprime if they have no common roots.} interpolating $\mathrm{u}_{
\mathrm{MS}}(z)$ in the sense that it `matches' each of the asymptotic 
expansions of $\mathrm{u}_{\mathrm{MS}}(z)$ to specified degrees at each 
of the points $0,1,\sqrt{2},\pi$, and the point at infinity, with the specified 
degrees chosen so that their sum equals $2(7n \! - \! 4)$, which leads precisely 
to the number of conditions necessary in order to determine, uniquely, the 
coefficients $\lbrace u_{m},v_{m} \rbrace_{m=0}^{7n-5}$; e.g., for the 
$K \! = \! 7$ pole set just discussed, for $n \! \in \! \mathbb{N}$ and, say, 
$k \! = \! 3$ $(\alpha_{3} \! = \! \infty)$, the corresponding interpolation 
conditions read $\mathcal{U}(z)(\mathcal{V}(z))^{-1} \! - \! \sum_{j=0}^{
\mathfrak{d}(z_{k})-1}c^{(k)}_{j}(z_{k})(z \! - \! z_{k})^{j} \! =_{z \to z_{k}} 
\! \mathcal{O}((z \! - \! z_{k})^{\mathfrak{d}(z_{k})})$, $k \! = \! 1,2,3,4$, 
where $z_{1} \! = \! 0$, $z_{2} \! = \! 1$, $z_{3} \! = \! \sqrt{2}$, $z_{4} 
\! = \! \pi$, and $c^{(k)}_{j}(z_{k}) \! := \! -\int_{\mathbb{R}}(\xi \! - \! 
z_{k})^{-(1+j)} \, \md \mathfrak{m}(\xi)$, and $\mathcal{U}(z)(\mathcal{V}
(z))^{-1} \! - \! \sum_{j=1}^{\mathfrak{d}(\infty)}c^{(\infty)}_{j}z^{-j} \! 
=_{z \to \infty} \! \mathcal{O}(z^{-(\mathfrak{d}(\infty)+1)})$, where 
$c^{(\infty)}_{j} \! := \! \int_{\mathbb{R}} \xi^{j-1} \, \md \mathfrak{m}(\xi)$, 
with $\mathfrak{d}(z_{k}),\mathfrak{d}(\infty) \! \in \! \mathbb{N}_{0}$, $k 
\! = \! 1,2,3,4$, chosen so that $\sum_{k=1}^{4} \mathfrak{d}(z_{k}) \! + \! 
\mathfrak{d}(\infty) \! = \! 2(7n \! - \! 4)$. The related, technically challenging 
problem is, of course, to study the asymptotic behaviour, as $n \! \to \! 
\infty$, of the sequence of MPAs; the remarkable connection, though, is 
this: the MPAs themselves, as well as the `error in approximation' (e.g., 
for the explicit $\alpha_{3} \! = \! \infty$ example just considered, 
the `error' is defined as $\mathcal{U}(z)(\mathcal{V}(z))^{-1} \! - \! 
\mathrm{u}_{\mathrm{MS}}(z))$, are described explicitly in terms of 
the monic ORFs and their---appropriately normalised---Cauchy 
transforms! The goal of describing the asymptotic behaviour of the 
sequence of MPAs can be achieved if one has a complete asymptotic, 
as $n \! \to \! \infty$, description of the monic ORFs and their 
corresponding Cauchy transforms; for further descriptions of these, 
and more, connections, the interested reader is referred, for example, to 
\cite{mppta38,mppta39,mppta40,mppta42,mppta35,mppta41,mppta32,
mppta22,mppta24,mppta18,mppta21,mppta61,mppta612,mpptanbglles,
mpptabosuwsino,mpptagllsmp,inopolgapam,newufpgllsmp1,mppta17,
mppta17b}). More generally, though, in this monograph, given $K$ 
$(\in \! \mathbb{N}$ and finite) poles $\alpha_{1},\alpha_{2},\dotsc,
\alpha_{K}$ on $\overline{\mathbb{R}}$, which are neither necessarily distinct 
nor bounded, and the positive probability measure with varying exponential 
weight $\widetilde{\mu}$, for $n \! \in \! \mathbb{N}$ and $k \! \in \! 
\lbrace 1,2,\dotsc,K \rbrace$, asymptotics, in the double-scaling limit 
$\mathscr{N},n \! \to \! \infty$ such that $z_{o} \! = \! 1 \! + \! o(1)$, 
are obtained for the corresponding MPAs of type $((n \! - 1)K \! + \! 
k \! - \! 1,(n \! - \! 1)K \! + \! k)$ and their associated `error(s) in 
approximation' for the so-called Markov-Stieltjes class of functions with 
representation $\mathrm{F}_{\tilde{\mu}}(z) \! := \! \int_{\mathbb{R}}
(z \! - \! \xi)^{-1} \, \md \widetilde{\mu}(\xi)$ in the entire complex 
plane; note, however, that this latter class of asymptotic problems requires 
the careful consideration of two sub-families of problems: one set of 
asymptotic analyses for those $k \! \in \! \lbrace 1,2,\dotsc,K \rbrace$ for 
which $\alpha_{k} \! = \! \infty$, and another set of asymptotic analyses 
for those $k \! \in \! \lbrace 1,2,\dotsc,K \rbrace$ for which $\alpha_{k} 
\! \neq \! \infty$. The papers which most closely capture the sense of the 
MPA problems addressed in this monograph are the seminal works of 
Nj\r{a}stad \cite{n2,n1}. With respect to question (Q3), it, and related 
rate(s)-of-convergence(s) problems, will be addressed elsewhere.

As further representative applications of why a generalisation {}from OPs 
to ORFs warrants consideration, the reader's attention is drawn to the 
following studies (by no means an exhaustive list!) {}from moment problems, 
in which one is given a list of pole locations and corresponding infinite 
sequences of real numbers and asks for the existence and the uniqueness 
of a measure whose moments relative to each pole are prescribed by one 
of the infinite sequences of real numbers, and `reconstruction problems' 
(see the monograph \cite{a15} for a more complete list of applications): \textbf{(i)} 
The \emph{Extended Hamburger Moment Problem} (EHMP) \cite{a1,a2}: given 
$p$ distinct real numbers $\tilde{a}_{1},\tilde{a}_{2},\dotsc,\tilde{a}_{p}$ 
and $p$ sequences of bounded real numbers $\lbrace 
\overset{\text{\tiny EH}}{c}^{\raise-1.0ex\hbox{$\scriptstyle (i)$}}_{k} 
\rbrace_{k \in \mathbb{N}}$, $i \! = \! 1,2,\dotsc,p$, find necessary and 
sufficient conditions for the existence of a distribution function (a real-valued, 
bounded, non-decreasing function with infinitely many points of increase 
on its domain of definition) $\mu^{\text{\tiny EH}}_{\text{\tiny MP}}$ on 
$\mathbb{R}$ such that 
$\int_{\mathbb{R}} \md \mu^{\text{\tiny EH}}_{\text{\tiny MP}}(\xi) \! = \! 1$ 
and $\overset{\text{\tiny EH}}{c}^{\raise-1.0ex\hbox{$\scriptstyle (i)$}}_{k} 
\! = \! \int_{\mathbb{R}}(\xi \! - \! \tilde{a}_{i})^{-k} \, \md 
\mu^{\text{\tiny EH}}_{\text{\tiny MP}}(\xi)$, $i \! = \! 1,2,\dotsc,p$, 
$k \! \in \! \mathbb{N}$; \textbf{(ii)} The \emph{Extended Stieltjes Moment Problem} 
(ESMP) \cite{a3}: given $q$ distinct real numbers $\hat{a}_{1},\hat{a}_{2},\dotsc,
\hat{a}_{q}$ ordered by size (e.g., $\hat{a}_{1} \! < \! \hat{a}_{2} \! < \! \dotsb \! 
< \! \hat{a}_{q})$, agree to call the real interval $[c,d]$ a \emph{Stieltjes interval} 
for the finite point set $\lbrace \hat{a}_{1},\hat{a}_{2},\dotsc,\hat{a}_{q} \rbrace$ 
if $(c,d) \cap \lbrace \hat{a}_{1},\hat{a}_{2},\dotsc,\hat{a}_{q} \rbrace \! = \! 
\varnothing$. Given $q$ sequences of bounded real numbers $\lbrace 
\overset{\text{\tiny ES}}{c}^{\raise-1.0ex\hbox{$\scriptstyle (r)$}}_{j} 
\rbrace_{j \in \mathbb{N}}$, $r \! = \! 1,2,\dotsc,q$, and a bounded real 
number $\overset{\text{\tiny ES}}{c}_{0}$, find necessary and sufficient 
conditions for the existence of a distribution function 
$\mu^{\text{\tiny ES}}_{\text{\tiny MP}}$, with all its points of 
increase on a given Stieltjes interval $[c,d]$, such that 
$\overset{\text{\tiny ES}}{c}_{0} \! = \! \int_{c}^{d} \md 
\mu^{\text{\tiny ES}}_{\text{\tiny MP}}(\xi)$ and 
$\overset{\text{\tiny ES}}{c}^{\raise-1.0ex\hbox{$\scriptstyle (r)$}}_{j} 
\! = \! \int_{c}^{d}(\xi \! - \! \hat{a}_{r})^{-j} \, \md 
\mu^{\text{\tiny ES}}_{\text{\tiny MP}}(\xi)$, $r \! = \! 1,2,\dotsc,q$, $j \! 
\in \! \mathbb{N}$; \textbf{(iii)} The \emph{Pick-Nevanlinna Problem} (PNP) \cite{a4}: 
consider a sequence of complex points $\lbrace z_{n} \rbrace_{n \in \mathbb{N}}$ 
such that $z_{n}$ coalesce into a finite number, $p_{o}$, say, of distinct 
real points $\check{a}_{1},\check{a}_{2},\dotsc,\check{a}_{p_{o}}$ according 
to the prescription $z_{p_{o}q_{o}+1} \! = \! \check{a}_{1}$, $z_{p_{o}q_{o}
+2} \! = \! \check{a}_{2}$, $\dotsc$, $z_{p_{o}q_{o}+p_{o}} \! = \! 
\check{a}_{p_{o}}$, $q_{o} \! \in \! \mathbb{N}_{0} \! := \! \mathbb{N} \cup 
\lbrace 0 \rbrace$. The corresponding PNP can be formulated thus: given 
the $p_{o}$ sequences of numbers $\lbrace \tilde{\gamma}^{(m)}_{l} 
\rbrace_{l \in \mathbb{N}_{0}}$, $m \! = \! 1,2,\dotsc,p_{o}$, find a 
Nevanlinna function {}\footnote{A function $\mathfrak{X}_{\Delta}(z)$ which is 
analytic for $z \! \in \! \mathbb{C}_{+} \! := \! \lbrace \mathstrut z \! \in \! 
\mathbb{C}; \, \Im (z) \! > \! 0 \rbrace$ with $\Im (\mathfrak{X}_{\Delta}(z)) \! 
\geqslant \! 0$ is called a Nevanlinna function.} $\mathscr{X}_{\Delta}(z)$ which 
has asymptotic expansions $\mathscr{X}_{\Delta}(z) \! =_{R_{m,\delta} \ni z \to 
\check{a}_{m}} \! \sum_{l \in \mathbb{N}_{0}} \tilde{\gamma}_{l}^{(m)}(z \! - \! 
\check{a}_{m})^{l}$, $m \! = \! 1,2,\dotsc,p_{o}$, where $R_{m,\delta} \! := \! 
\lbrace \mathstrut z \! \in \! \mathbb{C}; \, \delta \! < \! \text{Arg}(z \! - \! 
\check{a}_{m}) \! < \! \pi \! - \! \delta \rbrace$, $\delta \! > \! 0$ (as shown 
in \cite{a4}, this PNP is related to the EHMP and certain (weak) MPA problems); 
and \textbf{(iv)} The \emph{Frequency Analysis Problem} (FAP) (see, for example, 
the review article \cite{a5}): the determination of the unknown frequencies 
$\omega_{j}$, $j \! = \! 1,2,\dotsc,J$, {}from a discrete time signal $\lbrace 
x_{N}(m) \! = \! \sum_{j=-J}^{J} \alpha_{j} \me^{\mi m \omega_{j}} 
\rbrace_{m=0}^{N-1}$ of observed values, where $N$ is the `sample size', 
$\alpha_{0} \! \geqslant \! 0$, $0 \! \neq \! \alpha_{-j} \! = \! \overline{\alpha_{j}} 
\! \in \! \mathbb{C}$, $0 \! =: \! \omega_{0} \! < \! \omega_{1} \! < \! \dotsb 
\! < \! \omega_{J} \! < \! \pi$, and $\omega_{-j} \! = \! -\omega_{j} \! \in \! 
\mathbb{R}$, $j \! = \! 1,2,\dotsc,J$. The FAP has been dealt with by exploiting 
the fact that, under various conditions, certain roots of the Szeg\H{o} OPs 
\cite{a6} on the unit circle, $\mathbb{T} \! := \! \lbrace \mathstrut z \! 
\in \! \mathbb{C}; \, \lvert z \rvert \! = \! 1 \rbrace$, converge, as $N \! 
\to \! \infty$, to the `frequency points' $\me^{\mi \omega_{j}}$ and $1$, $j \! 
= \! \pm 1,\pm 2,\dotsc,\pm J$, and the `remaining (uninteresting) roots' are 
bounded away {}from $\mathbb{T}$ as $N \! \to \! \infty$. Recently, however, 
the question of the generalisation, or extension, of this theory, where, in lieu of 
Szeg\H{o} OPs, certain ORFs replace polynomials, has been raised \cite{a7,a8,a9}.

As shown in \cite{n2,n1,a1,a2,a3,a4,a7,a8,a9}, the principal technical 
observation subsumed in the analyses of the above-mentioned PA, MPA, 
moment, and reconstruction problems is to consider, in lieu of (Szeg\H{o}) 
OPs \cite{a6}, suitably orthogonalised rational functions, with no poles in 
$\overline{\mathbb{C}}$ outside of a pre-assigned pole set whose elements 
are neither necessarily distinct nor bounded (the point at infinity is included), 
called, generically, ORFs. Historically, to the best of the authors' knowledge 
as at the time of the presents, it seems that M.~M.~Djrbashian instigated 
the study of ORFs (see, for example, \cite{a10,a11,a12,a13,a14}). Since 
then, a monumental study of ORFs has been undertaken by A.~Bultheel, 
P.~Gonz\'{a}lez-Vera, E.~Hendriksen, and O.~Nj\r{a}stad (see the 
monograph \cite{a15} and the plethora of references therein; see, also, 
\cite{p1,p2,a18}): there has also been concomitant progress for the matrix 
generalisation of the scalar ORF theory on $\mathbb{T}$ (see, for example, 
\cite{a19,a20}). ORFs have applications, and potential applications, to 
diverse problems in generalised moment theory (see, for example, 
\cite{a1,a2,a3,a4,a21,a22}), interpolation theory (see, for example, 
\cite{a23}), control theory (see, for example, \cite{a29}), MPAs (see, for 
example, \cite{n2,n1,mppta59,ap2,ap1,a31,mppta23,a30,pz20}), an 
extended Toda lattice (see, for example, \cite{a38}), uniform approximation 
of $\text{sgn}(x)$ (see, for example, \cite{a39}; see, also, \cite{mppta25,
mppta29}), spectral theory (see, for example, \cite{a40}), the search for a 
rational variant of the Khrushchev formula (see, for example, \cite{a41}), 
and Christoffel functions and universality limits (see, for example, 
\cite{kddsl1,kddsl2}).

Thus far, the bulk of the analyses, asymptotic or otherwise, of ORFs on 
$\mathbb{T}$ assume that the poles of the ORFs lie in the interior of 
the open unit disc, $\mathbb{O} \! := \! \lbrace \mathstrut z \! \in \! 
\mathbb{C}; \, \lvert z \rvert \! < \! 1 \rbrace$ (see, however, Chapter~11 
of \cite{a15}, and \cite{a17,a42,l1,bklo}), and, for the scarce number of 
analyses, asymptotic or otherwise, of ORFs on the extended real line, 
$\overline{\mathbb{R}}$, the poles of the ORFs are assumed to be 
real and disjoint {}from the support of the orthogonality measure 
\cite{a43,a16,abpgvehon} (see, however, Chapter~11 of \cite{a15}), to 
lie in $\overline{\mathbb{C}} \setminus \overline{\mathbb{R}}$, or, due 
to so-called `special technical considerations', some `forbidden value' 
of a pole must be excluded {}from the analyses. Since, in the vast 
majority of the analyses, asymptotic or otherwise, the generalisation of 
OPs on $\mathbb{T}$ and on $\overline{\mathbb{R}}$ requires that the 
poles be in the exterior of the closed unit disc, $\overline{\mathbb{C}} 
\setminus (\mathbb{O} \cup \mathbb{T})$, or in $\overline{\mathbb{R}}$, 
respectively, the ORF pole location problem eluded to must be addressed 
in order to have a comprehensive understanding of the ORF theory 
(see, also, \cite{eelosoasar}). Most poignantly for the case of ORFs on 
$\overline{\mathbb{R}}$, which is the focus of this monograph, there is 
a dearth of analysis: the technically challenging aspect of the analysis, 
asymptotic or otherwise, is when the poles of the ORFs lie in the support 
of the orthogonality measure; in particular, to quote the words of Van 
Deun \emph{et al.} \cite{a44}: ``The computation of orthogonal rational 
functions with poles very close to the boundary is a complicated matter''. 
(In \cite{a44}, Van Deun \emph{et al.} consider the case of a finite interval.)

The present monograph, which is the second installment of a planned 
series of works dedicated to ORFs and related quantities, constitutes a 
generic study of ORFs when the prescribed (but otherwise arbitrary!), 
finite-in-number and not necessarily distinct poles lie on $\overline{
\mathbb{R}}$ (the support of the orthogonality measure): \textbf{these 
ORFs will, hereafter, be referred to as `mixed-pole case ORFs' (MPC 
ORFs)}.\footnote{In the theoretical study \cite{a45}, the so-called 
`fixed-pole case ORFs' (FPC ORFs) were studied, where the pre-assigned, 
but otherwise arbitrary, finite-in-number and not necessarily distinct 
poles of the FPC ORFs are all bounded and lie in $\mathbb{R}$.} More 
precisely, given the $K$ $(\in \! \mathbb{N}$ and finite) poles $\alpha_{1},
\alpha_{2},\dotsc,\alpha_{K}$ on $\overline{\mathbb{R}}$, which are neither 
necessarily distinct nor bounded, and the positive probability measure 
with varying exponential weight, $\md \widetilde{\mu}(z) \! = \! \exp 
(-n \widetilde{V}(z)) \, \md z$, $n \! \in \! \mathbb{N}$, where, with 
$z_{o} \! = \! z_{o}(\mathscr{N},n) \! > \! 0$, $\mathscr{N} \! \in \! 
\mathbb{N}$, $\widetilde{V}(z) \! := \! z_{o}V(z)$ satisfies the analyticity 
and growth conditions~\eqref{eq20}--\eqref{eq22}, define, for $n \! 
\in \! \mathbb{N}$ and $k \! \in \! \lbrace 1,2,\dotsc,K \rbrace$, the 
monic MPC ORF as per the---technical---discussion above (via the 
class-independent notation), $\pmb{\pi}^{n}_{k}(z)$, $z \! \in \! 
\mathbb{C}$: for $k \! \in \! \lbrace 1,2,\dotsc,K \rbrace$ such that 
$\alpha_{k} \! = \! \infty$, denote the associated `norming constant' 
(leading coefficient) by $\mu_{n,\varkappa_{nk}}^{\infty}(n,k)$ $(> \! 0)$, 
and define the corresponding MPC ORF as $\phi_{k}^{n}(z) \! := \! 
\mu_{n,\varkappa_{nk}}^{\infty}(n,k) \pmb{\pi}^{n}_{k}(z)$, $z \! \in \! 
\mathbb{C}$; and, for $k \! \in \! \lbrace 1,2,\dotsc,K \rbrace$ such that 
$\alpha_{k} \! \neq \! \infty$, denote the associated `norming constant' 
(leading coefficient) by $\mu_{n,\varkappa_{nk}}^{f}(n,k)$ $(> \! 0)$, and define 
the corresponding MPC ORF as $\phi_{k}^{n}(z) \! := \! \mu_{n,\varkappa_{nk}}^{f}
(n,k) \pmb{\pi}^{n}_{k}(z)$, $z \! \in \! \mathbb{C}$.\footnote{Despite the fact 
that the expression for the MPC ORF is different for those $k \! \in \! \lbrace 
1,2,\dotsc,K \rbrace$ for which $\alpha_{k} \! = \! \infty$ or $\alpha_{k} \! 
\neq \! \infty$, the generic notation $\phi_{k}^{n}(z)$ is retained throughout 
this monograph, irrespective of class.} Furthermore, for $n \! \in \! \mathbb{N}$ 
and $k \! \in \! \lbrace 1,2,\dotsc,K \rbrace$ such that $\alpha_{k} \! = \! 
\infty$, denote the corresponding MPA of type $((n \! - \! 1)K \! + \! k \! - \! 
1,(n \! - \! 1)K \! + \! k)$ for the Markov-Stieltjes transform of the probability 
measure $\widetilde{\mu}$, that is, $\mathrm{F}_{\tilde{\mu}}(z)$, by $\widehat{
\mathrm{U}}_{\tilde{\mu}}(z)(\widehat{\mathrm{V}}_{\tilde{\mu}}(z))^{-1}$, 
where $\widehat{\mathrm{U}}_{\tilde{\mu}}(z)$ is a polynomial of degree 
$(n \! - \! 1)K \! + \! k \! - \! 1$, $\widehat{\mathrm{V}}_{\tilde{\mu}}(z)$ 
$(\not\equiv \! 0)$ is a monic polynomial of degree $(n \! - \! 1)K \! + \! 
k$, and $(\widehat{\mathrm{U}}_{\tilde{\mu}}(z),\widehat{\mathrm{V}}_{
\tilde{\mu}}(z))$ coprime, with associated MPA error term define as 
$\widehat{\pmb{\mathrm{E}}}_{\tilde{\mu}}(z) \! := \! \widehat{\mathrm{U}}_{
\tilde{\mu}}(z)(\widehat{\mathrm{V}}_{\tilde{\mu}}(z))^{-1} \! - \! \mathrm{F}_{
\tilde{\mu}}(z)$, and, for $n \! \in \! \mathbb{N}$ and $k \! \in \! \lbrace 1,2,
\dotsc,K \rbrace$ such that $\alpha_{k} \! \neq \! \infty$, denote the corresponding 
MPA of type $((n \! - \! 1)K \! + \! k \! - \! 1,(n \! - \! 1)K \! + \! k)$ for 
$\mathrm{F}_{\tilde{\mu}}(z)$ by $\widetilde{\mathrm{U}}_{\tilde{\mu}}
(z)(\widetilde{\mathrm{V}}_{\tilde{\mu}}(z))^{-1}$, where $\widetilde{
\mathrm{U}}_{\tilde{\mu}}(z)$ is a polynomial of degree $(n \! - \! 1)K \! 
+ \! k \! - \! 1$, $\widetilde{\mathrm{V}}_{\tilde{\mu}}(z)$ $(\not\equiv 
\! 0)$ is a non-monic polynomial of degree $(n \! - \! 1)K \! + \! k$, 
and $(\widetilde{\mathrm{U}}_{\tilde{\mu}}(z),\widetilde{\mathrm{V}}_{
\tilde{\mu}}(z))$ coprime, with associated MPA error term define as 
$\widetilde{\pmb{\mathrm{E}}}_{\tilde{\mu}}(z) \! := \! \widetilde{\mathrm{U}}_{
\tilde{\mu}}(z)(\widetilde{\mathrm{V}}_{\tilde{\mu}}(z))^{-1} \! - \! \mathrm{F}_{
\tilde{\mu}}(z)$. Succinctly stated: for $n \! \in \! \mathbb{N}$ and $k 
\! \in \! \lbrace 1,2,\dotsc,K \rbrace$, the monic MPC ORF problem is 
characterised as a family of $K$ matrix RHPs on $\overline{\mathbb{R}}$, 
where one subfamily, corresponding to those $k \! \in \! \lbrace 1,2,
\dotsc,K \rbrace$ for which $\alpha_{k} \! = \! \infty$, consists of 
$\hat{K} \! := \! \# \lbrace \mathstrut k \! \in \! \lbrace 1,2,\dotsc,
K \rbrace; \, \alpha_{k} \! = \! \infty \rbrace$ matrix RHPs on 
$\overline{\mathbb{R}}$, and another subfamily, corresponding to those 
$k \! \in \! \lbrace 1,2,\dotsc,K \rbrace$ for which $\alpha_{k} \! \neq \! 
\infty$, consists of $\tilde{K} \! := \! \# \lbrace \mathstrut k \! \in \! \lbrace 
1,2,\dotsc,K \rbrace; \, \alpha_{k} \! \neq \! \infty \rbrace$ matrix RHPs on 
$\overline{\mathbb{R}}$, with $K \! = \! \hat{K} \! + \! \tilde{K}$, which, 
after the formal introduction of a family of $K$ variational problems 
containing external fields with singular points and the corresponding 
establishment of the existence, uniqueness, and regularity of the 
associated family of $K$ \emph{equilibrium measures}, is transformed 
into an equivalent family of $K$ `model' matrix RHPs on $\overline{
\mathbb{R}}$, whence uniform asymptotics, in the double-scaling limit 
$\mathscr{N},n \! \to \! \infty$ such that $z_{o} \! = \! 1 \! + \! o(1)$, 
in the entire complex plane for $\pmb{\pi}_{k}^{n}(z)$, $\mu_{n,
\varkappa_{nk}}^{\infty}(n,k)$, $\mu_{n,\varkappa_{nk}}^{f}(n,k)$, 
$\phi_{k}^{n}(z)$, $\widehat{\mathrm{U}}_{\tilde{\mu}}(z)
(\widehat{\mathrm{V}}_{\tilde{\mu}}(z))^{-1}$, $\widetilde{\mathrm{U}}_{
\tilde{\mu}}(z)(\widetilde{\mathrm{V}}_{\tilde{\mu}}(z))^{-1}$, 
$\widehat{\pmb{\mathrm{E}}}_{\tilde{\mu}}(z)$, and 
$\widetilde{\pmb{\mathrm{E}}}_{\tilde{\mu}}(z)$, as well as related, 
important (in their own right!) objects (zeros of monic MPC ORFs, 
normalised zero counting measures, etc.), are obtained.

This monograph is organised as follows. In Subsection~\ref{subsec1.2}, 
for $K \! \in \! \mathbb{N}$ and finite, the pole sequence (on 
$\overline{\mathbb{R}})$ $\alpha_{1},\alpha_{2},\dotsc,\alpha_{K}$ 
is partitioned according to whether $\alpha_{k} \! = \! \infty$ (see 
Subsection~\ref{subsubsec1.2.1}) or $\alpha_{k} \! \neq \! \infty$ (see 
Subsection~\ref{subsubsec1.2.2}), $k \! \in \! \lbrace 1,2,\dotsc,K 
\rbrace$, and, for $(n,k) \! \in \! \mathbb{N} \times \lbrace 1,2,\dotsc,
K \rbrace$, the corresponding monic MPC ORFs, $\pmb{\pi}_{k}^{n}(z)$, 
$z \! \in \! \mathbb{C}$, `norming constants' (leading coefficients), 
$\mu_{n,\varkappa_{nk}}^{r}(n,k)$, $r \! \in \! \lbrace \infty,f \rbrace$, 
and MPC ORFs, $\phi_{k}^{n}(z) \! = \! \mu_{n,\varkappa_{nk}}^{r}(n,k) 
\pmb{\pi}_{k}^{n}(z)$, $z \! \in \! \mathbb{C}$, are defined; furthermore, 
the corresponding MPAs of type $((n \! - \! 1)K \! + \! k \! - \! 1,(n \! - 
\! 1)K \! + \! k)$ and the associated MPA error terms are introduced. 
In Subsection~\ref{subsec1.3}, a self-contained synopsis of relevant 
facts {}from the theory of compact Riemann Surfaces which are 
necessary for the analysis machinery of this monograph is given (see 
Subsection~\ref{subsub1}), and the asymptotic, in the double-scaling 
limit $\mathscr{N},n \! \to \! \infty$ such that $z_{o} \! = \! 1 \! + \! o(1)$, 
results of this monograph are discussed, and, for the reader's convenience, 
summarised as Theorems~\ref{maintheoforinf1}--\ref{maintheoforfin2} 
(see Subsection~\ref{subsub2}). In Section~\ref{sec2}, the monic 
MPC ORF problem is formulated as a family of $K$ matrix RHPs on 
$\overline{\mathbb{R}}$, and the corresponding existence and 
uniqueness of this family of RHPs is established; en route, various 
integral representations for $\pmb{\pi}_{k}^{n}(z)$, $z \! \in \! 
\mathbb{C}$, $\mu_{n,\varkappa_{nk}}^{r}(n,k)$, $r \! \in \! \lbrace 
\infty,f \rbrace$, and, subsequently, $\phi_{k}^{n}(z) \! = \! \mu_{n,
\varkappa_{nk}}^{r}(n,k) \pmb{\pi}_{k}^{n}(z)$, $z \! \in \! \mathbb{C}$, 
are derived. In Section~\ref{sec3}, a family of $K$ variational (energy 
minimisation) problems containing external fields with singular points 
is introduced, and the existence, uniqueness, and regularity for the 
corresponding family of $K$ equilibrium measures is established; moreover, 
the asymptotic, in the double-scaling limit $\mathscr{N},n \! \to \! \infty$ 
such that $z_{o} \! = \! 1 \! + \! o(1)$, distribution of the zeros of the monic 
MPC ORFs is addressed. In Section~\ref{sec4}, the family of $K$ equilibrium 
measures, in conjunction with a corresponding family of complex potentials 
$(g$-functions), is used to transform the family of $K$ matrix RHPs on 
$\overline{\mathbb{R}}$ characterising the monic MPC ORF problem into an 
equivalent family of $K$ `model' matrix RHPs on $\overline{\mathbb{R}}$, 
which, in the double-scaling limit $\mathscr{N},n \! \to \! \infty$ such that 
$z_{o} \! = \! 1 \! + \! o(1)$, are solved explicitly in terms of Riemann theta 
functions and Airy functions. In Section~\ref{sek5}, which is lengthy, several 
key technical results in the form of lemmata are established, whence all the 
asymptotic, in the double-scaling limit $\mathscr{N},n \! \to \! \infty$ such 
that $z_{o} \! = \! 1 \! + \! o(1)$, results of this monograph are derived; 
moreover, a detailed analysis of the corresponding MPAs and the associated 
MPA error terms is presented. This monograph concludes with an appendix 
which discusses, concisely, several representative `first moments' of the 
family of $K$ equilibrium measures.
\begin{eeee} \label{aboutrecrels} 
\textsl{{\rm MPC ORFs} satisfy a cumbersome system of rational recurrence 
relations: their study and corresponding asymptotics, in the double-scaling 
limit $\mathscr{N},n \! \to \! \infty$ such that $z_{o} \! = \! 1 \! + \! o(1)$, 
will be presented elsewhere.}
\end{eeee}
\subsection{Partitions, MPC ORFs, and MPAs} \label{subsec1.2} 
There is a plethora of definitions which must be presented in this 
Subsection~\ref{subsec1.2} in order to completely describe the family of 
MPC ORFs. In order to parse the definitions in a digestible manner, an 
arbitrary, basic example of an MPC ORF (real) pole sequence is interspersed 
at various stages.

One begins by prescribing a sequence of $K$ $(\in \! \mathbb{N}$ and 
finite) poles on $\overline{\mathbb{R}}$, denoted by $\alpha_{1},\alpha_{2},
\dotsc,\alpha_{K}$, the elements of which are neither necessarily distinct 
(that is, $\lbrace \mathstrut \alpha_{i} \! = \! \alpha_{j}, \, i \! \neq \! j \! 
\in \! \lbrace 1,2,\dotsc,K \rbrace \rbrace$ is not necessarily empty) nor 
bounded; e.g., the MPC ORF pole set that will be referred to frequently 
throughout this Subsection~\ref{subsec1.2}, unless stated otherwise, for 
notational and illustrative purposes only is the following one of `length' 
$K \! = \! 7$: $\lbrace \alpha_{1},\alpha_{2},\alpha_{3},\alpha_{4},
\alpha_{5},\alpha_{6},\alpha_{7} \rbrace \! = \! \lbrace 0,1,\infty,1,
\sqrt{2},\pi,\infty \rbrace$.
\begin{eeee} \label{rem1.2.1} 
\textsl{The reader should not mistakenly interpret the notation $\lbrace 
\alpha_{1},\alpha_{2},\dotsc,\alpha_{K} \rbrace$ to be standard notation 
{}from set theory; rather, and as it is used throughout this monograph, 
unless stated otherwise, it should be thought of as a lexicographic listing 
of elements (individual pole locations) which may, therefore, repeat, and 
repetitions, if any, must be enumerated accordingly (see the examples 
below).}
\end{eeee}
\begin{eeee} \label{rem1.2.2} 
\textsl{Although $\overline{\mathbb{C}}$ (resp., $\overline{\mathbb{R}})$ 
is the standard notation for the (closed) Riemann sphere (resp., closed real 
line), the simplified, and somewhat abusive, notation $\mathbb{C}$ (resp., 
$\mathbb{R})$ is used to denote both the (closed) Riemann sphere (resp., 
closed real line) and the (open) complex field (resp., open real line), and 
the context(s) should make clear which object(s) the notation $\mathbb{C}$ 
(resp., $\mathbb{R})$ represents.}
\end{eeee}
Prescribe a measure of orthogonality, which will be taken to be a probability 
measure $\mu \! \in \! \mathscr{M}_{1}(\mathbb{R})$ {}\footnote{$\mathscr{M}_{1}
(\mathbb{R})$ denotes the set of all positive unit Borel measures on $\mathbb{R}$ 
for which all moments at the point at infinity and at $\alpha_{k} \! \neq \! \infty$, 
$k \! \in \! \lbrace 1,2,\dotsc,K \rbrace$, exist, that is, $\mathscr{M}_{1}
(\mathbb{R}) \! := \! \left\lbrace \mathstrut \mu; \int_{\mathbb{R}} \md \mu 
(\xi) \! = \! 1, \, \int_{\mathbb{R}} \xi^{m} \, \md \mu (\xi) \! < \! \infty, \, m 
\! \in \! \mathbb{N}, \, \int_{\mathbb{R}}(\xi \! - \! \alpha_{k})^{-m} \, \md 
\mu (\xi) \! < \! \infty, \, \alpha_{k} \! \neq \! \infty, \, k \! \in \! \lbrace 1,2,
\dotsc,K \rbrace \right\rbrace$. Throughout this monograph, $\mathbb{R}$ 
is oriented {}from $-\infty$ to $+\infty$, unless stated otherwise.} of the form
\begin{equation}
\md \mu (z) \! = \! \widetilde{w}(z) \, \md z, \label{eq1}
\end{equation}
with varying exponential weight function
\begin{equation}
\widetilde{w}(z) \! = \! \exp (-\mathscr{N}V(z)), \quad \mathscr{N} 
\! \in \! \mathbb{N}, \label{eq2}
\end{equation}
where the \emph{external field} $V \colon \overline{\mathbb{R}} \setminus 
\lbrace \alpha_{1},\alpha_{2},\dotsc,\alpha_{K} \rbrace \! \to \! \mathbb{R}$ 
satisfies the following conditions:
\begin{gather}
\text{$V(z)$ is real analytic on $\overline{\mathbb{R}} \setminus 
\lbrace \alpha_{1},\alpha_{2},\dotsc,\alpha_{K} \rbrace$}; \label{eq3} \\
\lim_{x \to \alpha_{i}} \left(\dfrac{V(x)}{\ln (x^{2} \! + \! 1)} \right) \! 
= \! +\infty, \quad i \! \in \! \lbrace \mathstrut k \! \in \! \lbrace 
1,2,\dotsc,K \rbrace; \, \alpha_{k} \! = \! \infty \rbrace; \label{eq4} \\
\lim_{x \to \alpha_{j}} \left(\dfrac{V(x)}{\ln ((x \! - \! \alpha_{j})^{-2} \! 
+ \! 1)} \right) \! = \! +\infty, \quad j \! \in \! \lbrace \mathstrut k 
\! \in \! \lbrace 1,2,\dotsc,K \rbrace; \, \alpha_{k} \! \neq \! \infty 
\rbrace. \label{eq5}
\end{gather}
For any pair $(n,k) \! \in \! \mathbb{N} \times \lbrace 1,2,\dotsc,K 
\rbrace$, denote by $\Lambda^{\mathbb{R}}_{n,k}$ the set of all 
rational functions with poles restricted to the pole set $\lbrace 
\alpha_{1},\alpha_{2},\dotsc,\alpha_{K} \rbrace$, that 
is,\footnote{Note the convention $\sum_{m=1}^{0} \boldsymbol{\ast} 
\! := \! 0$.} $\Lambda^{\mathbb{R}}_{n,k} \! := \! \lbrace \mathstrut 
f \colon \mathbb{N} \times \lbrace 1,2,\dotsc,K \rbrace \times 
\overline{\mathbb{C}} \setminus \lbrace \alpha_{1},\alpha_{2},
\dotsc,\alpha_{K} \rbrace \! \to \! \mathbb{C}; \, f(z) \! = \! d^{\flat}_{0}
(n,k) \! + \! \sum_{p=1}^{n-1} \sum_{q=1}^{K}d^{\sharp}_{p,q}(n,k) 
\mathscr{S}^{p}_{q}(z) \! + \! \sum_{r=1}^{k}d^{\natural}_{n,r}(n,k) 
\mathscr{S}^{n}_{r}(z), \, d^{\flat}_{0}(n,k), \, d^{\sharp}_{p,q}(n,k), \, 
d^{\natural}_{n,r}(n,k) \! \in \! \mathbb{R} \rbrace$, with
\begin{equation*}
\mathscr{S}^{n}_{k}(z) \! = \! 
\begin{cases}
z^{\varkappa_{nk}}, &\text{$\alpha_{k} \! = \! \infty$,} \\
(z \! - \! \alpha_{k})^{-\varkappa_{nk}}, &\text{$\alpha_{k} 
\! \neq \! \infty$,}
\end{cases}
\end{equation*}
where
\begin{equation*}
\varkappa_{nk} \colon \mathbb{N} \times \lbrace 1,2,\dotsc,K \rbrace 
\! \to \! \mathbb{N}, \, \, (n,k) \! \mapsto \! \varkappa_{nk} \! = \! 
(n \! - \! 1) \gamma_{k} \! + \! \varrho_{k}
\end{equation*}
denotes the multiplicity of the pole $\alpha_{k}$ in the \emph{repeated 
pole sequence}
\begin{equation*}
\overset{1}{\lbrace \alpha_{1},\alpha_{2},\dotsc,\alpha_{K} \rbrace} \cup 
\dotsb \cup \overset{n-1}{\lbrace \alpha_{1},\alpha_{2},\dotsc,\alpha_{K} 
\rbrace} \cup \overset{n}{\lbrace \alpha_{1},\alpha_{2},\dotsc,\alpha_{k} 
\rbrace},
\end{equation*}
with $\gamma_{k}$ the `repeating number' of the pole $\alpha_{k}$ in the set 
$\lbrace \alpha_{1},\alpha_{2},\dotsc,\alpha_{K} \rbrace$, and $\varrho_{k}$ 
the `repeating index' of the pole $\alpha_{k}$ up to and including the 
$k$th position in the set $\lbrace \alpha_{1},\alpha_{2},\dotsc,\alpha_{K} 
\rbrace$; e.g., for the $K \! = \! 7$ pole set $\lbrace \alpha_{1},\alpha_{2},
\alpha_{3},\alpha_{4},\alpha_{5},\alpha_{6},\alpha_{7} \rbrace \linebreak[4] 
\! = \! \lbrace 0,1,\infty,1,\sqrt{2},\pi,\infty \rbrace$,
\begin{equation*}
\gamma_{1} \! = \! 1, \, \gamma_{2} \! = \! \gamma_{4} \! = \! 2, \, 
\gamma_{3} \! = \! \gamma_{7} \! = \! 2, \, \gamma_{5} \! = \! 1, \, 
\gamma_{6} \! = \! 1, \, \varrho_{1} \! = \! 1, \, \varrho_{2} \! = \! 1, 
\, \varrho_{3} \! = \! 1, \, \varrho_{4} \! = \! 2, \, \varrho_{5} \! 
= \! 1, \, \varrho_{6} \! = \! 1, \, \varrho_{7} \! = \! 2,
\end{equation*}
and
\begin{align*}
\varkappa_{n1}=& \, (n \! - \! 1) \gamma_{1} \! + \! \varrho_{1} \! = \! 
(n \! - \! 1) \! + \! 1 \! = \! n, \\
\varkappa_{n2}=& \, (n \! - \! 1) \gamma_{2} \! + \! \varrho_{2} \! = \! 
2(n \! - \! 1) \! + \! 1 \! = \! 2n \! - \! 1, \\
\varkappa_{n3}=& \, (n \! - \! 1) \gamma_{3} \! + \! \varrho_{3} \! = \! 
2(n \! - \! 1) \! + \! 1 \! = \! 2n \! - \! 1, \\
\varkappa_{n4}=& \, (n \! - \! 1) \gamma_{4} \! + \! \varrho_{4} \! = \! 
2(n \! - \! 1) \! + \! 2 \! = \! 2n, \\
\varkappa_{n5}=& \, (n \! - \! 1) \gamma_{5} \! + \! \varrho_{5} \! = \! 
(n \! - \! 1) \! + \! 1 \! = \! n, \\
\varkappa_{n6}=& \, (n \! - \! 1) \gamma_{6} \! + \! \varrho_{6} \! = \! 
(n \! - \! 1) \! + \! 1 \! = \! n, \\
\varkappa_{n7}=& \, (n \! - \! 1) \gamma_{7} \! + \! \varrho_{7} \! = \! 
2(n \! - \! 1) \! + \! 2 \! = \! 2n.
\end{align*}
Denote by $\Lambda^{\mathbb{R}} \! := \! \cup_{n \in \mathbb{N}} 
\cup_{k=1}^{K} \Lambda_{n,k}^{\mathbb{R}}$ the linear space 
over the field $\mathbb{R}$ spanned by a constant and the 
set of rational functions $\lbrace \mathscr{S}^{n}_{k}(z) 
\rbrace_{\underset{k=1,2,\dotsc,K}{n \in \mathbb{N}}}$.\footnote{
Note that the function theoretic sets $\Lambda_{n,k}^{\mathbb{R}}$, 
$(n,k) \! \in \! \mathbb{N} \times \lbrace 1,2,\dotsc,K \rbrace$, and 
$\Lambda^{\mathbb{R}}$ form linear spaces over $\mathbb{R}$ with 
respect to the binary operations of addition and multiplication by a scalar.} 
A function element $0 \! \neq \! f \! \in \! \Lambda^{\mathbb{R}}$ is 
called a \emph{rational function corresponding to the pole set} $\lbrace 
\alpha_{1},\alpha_{2},\dotsc,\alpha_{K} \rbrace$. The ordered base of 
rational functions for $\Lambda^{\mathbb{R}}$ is
\begin{equation*}
\mathcal{B} \sim \left\{\text{const}., \, \overbrace{\underbrace{
\mathscr{S}^{1}_{1}(z),\mathscr{S}^{1}_{2}(z),\dotsc,\mathscr{S}^{1}_{K}
(z)}_{k=1,2,\dotsc,K}}^{n=1}, \, \overbrace{\underbrace{\mathscr{S}^{2}_{1}
(z),\mathscr{S}^{2}_{2}(z),\dotsc,\mathscr{S}^{2}_{K}(z)}_{k=1,2,\dotsc,
K}}^{n=2}, \, \dotsc, \, \overbrace{\underbrace{\mathscr{S}^{m}_{1}(z),
\mathscr{S}^{m}_{2}(z),\dotsc,\mathscr{S}^{m}_{K}(z)}_{k=1,2,\dotsc,K}}^{n=m}, 
\, \dotsc \right\},
\end{equation*}
corresponding, respectively, to the \emph{cyclically repeated pole sequence}
\begin{equation*}
\mathcal{P} \sim \left\{\text{no pole}, \, \overbrace{\underbrace{\alpha_{1},
\alpha_{2},\dotsc,\alpha_{K}}_{k=1,2,\dotsc,K}}^{n=1}, \, \overbrace{
\underbrace{\alpha_{1},\alpha_{2},\dotsc,\alpha_{K}}_{k=1,2,\dotsc,K}}^{n=2}, 
\, \dotsc, \, \overbrace{\underbrace{\alpha_{1},\alpha_{2},\dotsc,
\alpha_{K}}_{k=1,2,\dotsc,K}}^{n=m}, \, \dotsc \right\}.
\end{equation*}
For $0 \! \neq \! f \! \in \! \Lambda^{\mathbb{R}}$, there exists a unique 
pair $(n,k) \! \in \! \mathbb{N} \times \lbrace 1,2,\dotsc,K \rbrace$ 
such that $f \! \in \! \Lambda^{\mathbb{R}}_{n,k}$. For $(n,k) \! \in \! 
\mathbb{N} \times \lbrace 1,2,\dotsc,K \rbrace$ and $0 \! \neq \! f \! 
\in \! \Lambda^{\mathbb{R}}_{n,k}$, define the \emph{leading coefficient} 
of $f$, symbolically $\operatorname{LC}(f)$, as $\operatorname{LC}(f) 
\! := \! d^{\natural}_{n,k}(n,k)$. For $(n,k) \! \in \! \mathbb{N} \times 
\lbrace 1,2,\dotsc,K \rbrace$, $0 \! \neq \! f \! \in \! \Lambda^{
\mathbb{R}}_{n,k}$ is called \emph{monic} if $\operatorname{LC}(f) \! = 
\! 1$. For $(n,k) \! \in \! \mathbb{N} \times \lbrace 1,2,\dotsc,K \rbrace$, 
define the linear functional $\mathscr{L}$ by its action on the (rational) 
basis elements of $\Lambda^{\mathbb{R}}$ as follows: $\mathscr{L} 
\colon \Lambda^{\mathbb{R}} \! \to \! \mathbb{R}$, $f \! = \! d^{\flat}_{0}
(n,k) \! + \! \sum_{p=1}^{n-1} \sum_{q=1}^{K}d^{\sharp}_{p,q}(n,k) 
\mathscr{S}^{p}_{q}(z) \! + \! \sum_{r=1}^{k}d^{\natural}_{n,r}(n,k) 
\mathscr{S}^{n}_{r}(z) \! \mapsto \! \mathscr{L}(f) \! = \! d^{\flat}_{0}
(n,k) \! + \! \sum_{p=1}^{n-1} \sum_{q=1}^{K}d^{\sharp}_{p,q}(n,k)
c_{p,q} \! + \! \sum_{r=1}^{k}d^{\natural}_{n,r}(n,k)c_{n,r}$, where
\begin{equation*}
c_{n,k} \! := \! \mathscr{L}(\mathscr{S}^{n}_{k}(z)) \! = \! 
\begin{cases}
\int_{\mathbb{R}} \xi^{\varkappa_{nk}} \, \md \mu (\xi), 
&\text{$\alpha_{k} \! = \! \infty$,} \\
\int_{\mathbb{R}}(\xi \! - \! \alpha_{k})^{-\varkappa_{nk}} \, \md \mu (\xi), 
&\text{$\alpha_{k} \! \neq \! \infty$.}
\end{cases}
\end{equation*}
(Of course, since $\mu \! \in \! \mathscr{M}_{1}(\mathbb{R})$, $c_{0,0} \! 
:= \! \mathscr{L}(1) \! = \! \int_{\mathbb{R}} \md \mu (\xi) \! = \! 1$.) 
Define the real bilinear form $\langle \pmb{\cdot},\pmb{\cdot} 
\rangle_{\mathscr{L}}$ as follows: $\langle \pmb{\cdot},\pmb{\cdot} 
\rangle_{\mathscr{L}} \colon \Lambda^{\mathbb{R}} \times \Lambda^{
\mathbb{R}} \! \to \! \mathbb{R}$, $(f,g) \! \mapsto \! \langle f,g 
\rangle_{\mathscr{L}} \! := \! \mathscr{L}(f(z)g(z)) \! = \! \int_{\mathbb{R}}
f(\xi)g(\xi) \, \md \mu (\xi)$. It is a fact that the bilinear form $\langle 
\pmb{\cdot},\pmb{\cdot} \rangle_{\mathscr{L}}$ thus defined is an inner 
product (see Section~\ref{sec2}, the proof of Lemma~\ref{lem2.1}), and 
this fact is used, with little or no further reference, throughout this lecture 
note.

If $f \! \in \! \Lambda^{\mathbb{R}}$, then $\lvert \lvert f(\pmb{\cdot}) 
\rvert \rvert_{\mathscr{L}} \! := \! (\langle f,f \rangle)^{1/2}$ is called 
the \emph{norm of $f$ with respect to $\mathscr{L}$} (note that $\lvert 
\lvert f(\pmb{\cdot}) \rvert \rvert_{\mathscr{L}} \! \geqslant \! 0$ for all 
$f \! \in \! \Lambda^{\mathbb{R}}$, and $\lvert \lvert f(\pmb{\cdot}) 
\rvert \rvert_{\mathscr{L}} \! > \! 0$ if $0 \! \neq \! f \! \in \! 
\Lambda^{\mathbb{R}})$. A sequence of rational functions will now be 
defined: $\lbrace \hat{\phi}^{n}_{k}(z) \rbrace_{\underset{k=1,2,\dotsc,
K}{n \in \mathbb{N}}}$ is called a---real---orthonormal rational function 
sequence with respect to $\mathscr{L}$ if, for $(n,k) \! \in \! \mathbb{N} 
\times \lbrace 1,2,\dotsc,K \rbrace$: (i) $\hat{\phi}^{n}_{k}(z) \! \in 
\! \Lambda^{\mathbb{R}}_{n,k}$; (ii) $\langle \hat{\phi}^{n}_{k},
\hat{\phi}^{n^{\prime}}_{k^{\prime}} \rangle_{\mathscr{L}} \! = \! 
\delta_{nn^{\prime}} \delta_{kk^{\prime}}$, where $\delta_{ij} \colon 
\mathbb{N} \times \mathbb{N} \! \to \! \lbrace 0,1 \rbrace$ is the 
Kronecker delta; and (iii) $\langle \hat{\phi}^{n}_{k},\hat{\phi}^{n}_{k} 
\rangle_{\mathscr{L}} \! =: \! \lvert \lvert \hat{\phi}^{n}_{k}(\pmb{\cdot}) 
\rvert \rvert_{\mathscr{L}}^{2} \! = \! 1$ (for consistency of notation, set 
$\hat{\phi}^{0}_{0}(z) \! \equiv \! 1$). In order to elucidate the precise 
structure of the orthogonality conditions for the MPC ORFs and to state 
the results of this monograph (see Subsection~\ref{subsub2}), the 
following notational preamble is necessary.

What follows next is an ordered (disjoint) partitioning of the index set 
$\lbrace 1,2,\dotsc,K \rbrace$ and the pole set $\lbrace \alpha_{1},
\alpha_{2},\dotsc,\linebreak[4]
\alpha_{K} \rbrace$. In order to proceed, though, a parameter, denoted 
by $\mathfrak{s}$, must be defined.
\begin{aaaa} \label{def1.2.1} 
\emph{For the general pole set $\lbrace \alpha_{1},\alpha_{2},\dotsc,\alpha_{K} 
\rbrace$ described above, let $\mathfrak{s}$ denote the number of distinct 
poles (see Remark~\ref{commentonss} below$)$$;$ e.g., for the $K \! = \! 7$ 
pole set $\lbrace \alpha_{1},\alpha_{2},\alpha_{3},\alpha_{4},\alpha_{5},
\alpha_{6},\alpha_{7} \rbrace \! = \! \lbrace 0,1,\infty,1,\sqrt{2},\pi,\infty 
\rbrace$, $\mathfrak{s} \! = \! 5$.}
\end{aaaa}
\begin{eeee} \label{commentonss} 
\textsl{More precisely, let $\hat{\mathfrak{s}}(j) \! := \! \lbrace \mathstrut 
k \! \in \! \lbrace 1,2,\dotsc,K \rbrace \setminus \lbrace 1,2,\dotsc,j \rbrace; 
\, \alpha_{k} \! = \! \alpha_{j} \rbrace$, $j \! = \! 1,2,\dotsc,K$$;$ then, 
$\mathfrak{s} \! = \! K \! - \! \sum_{j=1}^{K} \# \hat{\mathfrak{s}}(j)$, with 
the caveat that the---standard---cardinal number operator, $\#$, obeys 
the following, additional rule: if $\hat{\mathfrak{s}}(k) \! \subseteq \! 
\hat{\mathfrak{s}}(j)$ for $j \! < \! k$, then $\# \hat{\mathfrak{s}}(k) \! := \! 0$ 
(this excludes repeated elements; of course, if $\hat{\mathfrak{s}}(j) \! = \! 
\varnothing$, $j \! \in \! \lbrace 1,2,\dotsc,K \rbrace$, then $\# 
\hat{\mathfrak{s}}(j) \! = \! 0)$. One can also let $\mathfrak{s} \! = \! K \! - \! 
\# (\cup_{j=1}^{K} \hat{\mathfrak{s}}(j))$, provided that in the set-theoretic union, 
$\cup_{j=1}^{K} \hat{\mathfrak{s}}(j)$, repeated elements are not counted: of 
course, with either of the above formulae, one arrives at the unique positive 
integer $\mathfrak{s}$$;$ e.g., for the $K \! = \! 13$ pole set $\lbrace \alpha_{1},
\alpha_{2},\alpha_{3},\alpha_{4},\alpha_{5},\alpha_{6},\alpha_{7},\alpha_{8},
\alpha_{9},\alpha_{10},\alpha_{11},\alpha_{12},\alpha_{13} \rbrace \! = \! \lbrace 
0,1,\infty,1,\sqrt{2},\pi,\infty,-2,-1,1,\infty,0,0 \rbrace$, $\hat{\mathfrak{s}}(1) 
\! = \! \lbrace 12,13 \rbrace$ $\Rightarrow$ $\# \hat{\mathfrak{s}}(1) \! = 
\! 2$, $\hat{\mathfrak{s}}(2) \! = \! \lbrace 4,10 \rbrace$ $\Rightarrow$ 
$\# \hat{\mathfrak{s}}(2) \! = \! 2$, $\hat{\mathfrak{s}}(3) \! = \! \lbrace 
7,11 \rbrace$ $\Rightarrow$ $\# \hat{\mathfrak{s}}(3) \! = \! 2$, 
$\hat{\mathfrak{s}}(4) \! = \! \lbrace 10 \rbrace$ $\Rightarrow$ $\# 
\hat{\mathfrak{s}}(4) \! = \! 0$ (since $\hat{\mathfrak{s}}(4) \! \subseteq \! 
\hat{\mathfrak{s}}(2))$, $\hat{\mathfrak{s}}(5) \! = \! \varnothing$ $\Rightarrow$ 
$\# \hat{\mathfrak{s}}(5) \! = \! 0$, $\hat{\mathfrak{s}}(6) \! = \! \varnothing$ 
$\Rightarrow$ $\# \hat{\mathfrak{s}}(6) \! = \! 0$, $\hat{\mathfrak{s}}(7) \! = 
\lbrace 11 \rbrace$ $\Rightarrow$ $\# \hat{\mathfrak{s}}(7) \! = \! 0$ (since 
$\hat{\mathfrak{s}}(7) \! \subseteq \! \hat{\mathfrak{s}}(3))$, $\hat{\mathfrak{s}}
(8) \! = \! \varnothing$ $\Rightarrow$ $\# \hat{\mathfrak{s}}(8) \! = \! 0$, 
$\hat{\mathfrak{s}}(9) \! = \! \varnothing$ $\Rightarrow$ $\# \hat{\mathfrak{s}}
(9) \! = \! 0$, $\hat{\mathfrak{s}}(10) \! = \! \varnothing$ $\Rightarrow$ $\# 
\hat{\mathfrak{s}}(10) \! = \! 0$, $\hat{\mathfrak{s}}(11) \! = \! \varnothing$ 
$\Rightarrow$ $\# \hat{\mathfrak{s}}(11) \! = \! 0$, $\hat{\mathfrak{s}}(12) 
\! = \! \lbrace 13 \rbrace$ $\Rightarrow$ $\# \hat{\mathfrak{s}}(12) \! = \! 0$ 
(since $\hat{\mathfrak{s}}(12) \! \subseteq \! \hat{\mathfrak{s}}(1))$, and 
$\hat{\mathfrak{s}}(13) \! = \! \varnothing$ $\Rightarrow$ $\# \hat{\mathfrak{s}}
(13) \! = \! 0$, hence $\mathfrak{s} \! = \! 13 \! - \! \sum_{j=1}^{13} \# 
\hat{\mathfrak{s}}(j) \! = \! 13 \! - \! (2 \! + \! 2 \! + \! 2 \! + \! 0 \! + \! 0 \! 
+ \! 0 \! + \! 0 \! + \! 0 \! + \! 0 \! + \! 0 \! + \! 0 \! + \! 0 \! + \! 0) \! = \! 
13 \! - \! 6 \! = \! 7$ (also, $\mathfrak{s} \! = \! 13 \! - \! \# (\cup_{j=1}^{13} 
\hat{\mathfrak{s}}(j)) \! = \! 13 \! - \! \# \lbrace 4,7,10,11,12,13 \rbrace \! = 
\! 13 \! - \! 6 \! = \! 7)$.}
\end{eeee}
\begin{eeee} \label{specialKScases} 
\textsl{The cases $K \! = \! \mathfrak{s} \! = \! 1$ and $K \! = \! \mathfrak{s} 
\! = \! 2$ correspond, respectively, after a possible M\"{o}bius transformation 
argument and a modification of the associated weight function, to {\rm OPs} 
(see, for example, {\rm \cite{a6,bs1,bs2,hstahlvtotikop,mehi,a51}}, and the 
references therein) and Orthogonal Laurent Polynomials (see, for example, 
{\rm \cite{a15}}, and the references therein$)$$;$ of course, the cases $K \! = 
\! \mathfrak{s} \! = \! 1$ and $K \! = \! \mathfrak{s} \! = \! 2$ are important, 
special reductions of the general case considered in this monograph.}
\end{eeee}
For $k \! = \! 1,2,\dotsc,K$, a decomposition of the index set 
corresponding to the poles distinct from $\alpha_{k}$ will be needed, 
that is, $\lbrace \mathstrut k^{\prime} \! \in \! \lbrace 1,2,\dotsc,K \rbrace; 
\, \alpha_{k^{\prime}} \! \neq \! \alpha_{k} \rbrace$.\footnote{Note that, 
strictly speaking, this set decomposition is a function of $k$; however, 
for simplicity of notation, this particular $k$ dependence, as well as the 
explicit $k$ dependence of all decompositions, partitions, indices, etc., 
henceforth, will, unless absolutely necessary, be suppressed, and the 
reader should be cognizant of this fact.} In order to decompose this set, 
one needs to consider the (possibly smaller) collection of distinct poles, 
ordered consistently with the original pole sequence, and with the pole 
$\alpha_{k}$ excised: the `size' (see below) of this set is $\mathfrak{s} \! 
- \! 1$. For the $j$th member of the reduced collection of poles, referred 
to as the `residual' pole set, the number of times that that pole appears 
in the original pole sequence shall be denoted $k_{j}$, $j \! = \! 1,2,
\dotsc,\mathfrak{s} \! - \! 1$. One then decomposes the set of integers 
into a disjoint union, so that the first subset is the collection of all 
integers corresponding to the first pole in this reduced collection of 
poles, the second subset is the collection of all integers corresponding 
to the second pole in this reduced collection of poles, etc., until the 
$(\mathfrak{s} \! - \! 1)$th subset. The precise definition of this 
decompsition, which exhausts the remainder of this 
Subsection~\ref{subsec1.2}, follows; however, for ease of presentation 
and readability, it is decomposed into two cases: (i) $n \! \in \! \mathbb{N}$ 
and $k \! \in \! \lbrace 1,2,\dotsc,K \rbrace$ such that $\alpha_{k} \! 
= \! \infty$ (see Subsection~\ref{subsubsec1.2.1}); and (ii) $n \! \in \! 
\mathbb{N}$ and $k \! \in \! \lbrace 1,2,\dotsc,K \rbrace$ such that 
$\alpha_{k} \! \neq \! \infty$ (see Subsection~\ref{subsubsec1.2.2}).
\subsubsection{$n \! \in \! \mathbb{N}$ and $k \! \in \! \lbrace 1,2,\dotsc,
K \rbrace$ such that $\alpha_{k} \! = \! \infty$} \label{subsubsec1.2.1} 
For $k \! \in \! \lbrace 1,2,\dotsc,K \rbrace$ such that $\alpha_{k} \! = 
\! \infty$, write the ordered disjoint integer partition
\begin{align*}
\lbrace \mathstrut k^{\prime} \! \in \! \lbrace 1,2,\dotsc,K \rbrace; \, 
\alpha_{k^{\prime}} \! \neq \! \alpha_{k} \! = \! \infty \rbrace \! :=& \, 
\lbrace \underbrace{i(1)_{1},i(1)_{2},\dotsc,i(1)_{k_{1}}}_{k_{1}} \rbrace 
\cup \lbrace \underbrace{i(2)_{1},i(2)_{2},\dotsc,i(2)_{k_{2}}}_{k_{2}} 
\rbrace \cup \dotsb \\
\dotsb& \, \cup \lbrace \underbrace{i(\mathfrak{s} \! - \! 1)_{1},
i(\mathfrak{s} \! - \! 1)_{2},\dotsc,
i(\mathfrak{s} \! - \! 1)_{k_{\mathfrak{s}-1}}}_{k_{\mathfrak{s}-1}} \rbrace 
\! = \bigcup_{q=1}^{\mathfrak{s}-1} \lbrace \underbrace{i(q)_{1},i(q)_{2},
\dotsc,i(q)_{k_{q}}}_{k_{q}} \rbrace,
\end{align*}
where $1 \! \leqslant \! i(q)_{1} \! < \! i(q)_{2} \! < \! \dotsb \! < \! 
i(q)_{k_{q}} \! \leqslant \! K$, $q \! \in \! \lbrace 1,2,\dotsc,\mathfrak{s} 
\! - \! 1 \rbrace$, $\lbrace i(j)_{1},i(j)_{2},\dotsc,i(j)_{k_{j}} \rbrace 
\cap \lbrace i(l)_{1},i(l)_{2},\dotsc,i(l)_{k_{l}} \rbrace \! = \! 
\varnothing$ $\forall$ $l \! \neq \! j \! \in \! \lbrace 1,2,\dotsc,
\mathfrak{s} \! - \! 1 \rbrace$, and $\# \lbrace \mathstrut k^{\prime} \! \in 
\! \lbrace 1,2,\dotsc,K \rbrace; \, \alpha_{k^{\prime}} \! \neq \! \alpha_{k} 
\! = \! \infty \rbrace \! = \! \sum_{q=1}^{\mathfrak{s}-1}k_{q} \! = \! 
K \! - \! \gamma_{k}$, which induces, on the pole set $\lbrace \alpha_{1},
\alpha_{2},\dotsc,\alpha_{K} \rbrace$, the following disjoint ordering,
\begin{align*}
\lbrace \mathstrut \alpha_{k^{\prime}}, \, k^{\prime} \! \in \! \lbrace 1,2,
\dotsc,K \rbrace; \, \alpha_{k^{\prime}} \! \neq \! \alpha_{k} \! = \! \infty 
\rbrace \! :=& \, \lbrace \underbrace{\alpha_{i(1)_{1}},\alpha_{i(1)_{2}},
\dotsc,\alpha_{i(1)_{k_{1}}}}_{k_{1}} \rbrace \cup \lbrace 
\underbrace{\alpha_{i(2)_{1}},\alpha_{i(2)_{2}},\dotsc,
\alpha_{i(2)_{k_{2}}}}_{k_{2}} \rbrace \cup \dotsb \\
\dotsb& \, \cup \lbrace \underbrace{\alpha_{i(\mathfrak{s}-1)_{1}},\alpha_{i
(\mathfrak{s}-1)_{2}},\dotsc,\alpha_{i(\mathfrak{s}-1)_{k_{\mathfrak{s}-1}}}
}_{k_{\mathfrak{s}-1}} \rbrace \! = \bigcup_{q=1}^{\mathfrak{s}-1} \lbrace 
\underbrace{\alpha_{i(q)_{1}},\alpha_{i(q)_{2}},\dotsc,
\alpha_{i(q)_{k_{q}}}}_{k_{q}} \rbrace,
\end{align*}
with $\alpha_{i(q)_{1}} \! \prec \! \alpha_{i(q)_{2}} \! \prec \! \dotsb \! 
\prec \! \alpha_{i(q)_{k_{q}}}$, $q \! \in \! \lbrace 1,2,\dotsc,\mathfrak{s} 
\! - \! 1 \rbrace$, where the notation $a \! \prec \! b$ means `$a$ 
precedes $b$' or `$a$ is to the left of $b$', and $\lbrace \alpha_{i(j)_{1}},
\alpha_{i(j)_{2}},\dotsc,\alpha_{i(j)_{k_{j}}} \rbrace \cap \lbrace 
\alpha_{i(l)_{1}},\alpha_{i(l)_{2}},\dotsc,\alpha_{i(l)_{k_{l}}} \rbrace \! 
= \! \varnothing$ $\forall$ $l \! \neq \! j \! \in \! \lbrace 1,2,\dotsc,
\mathfrak{s} \! - \! 1 \rbrace$, such that {}\footnote{If all the poles are 
distinct, that is, $\alpha_{i} \! \neq \! \alpha_{j}$ $\forall$ $i \! \neq 
\! j \! \in \! \lbrace 1,2,\dotsc,K \rbrace$, then, for $k \! \in \! \lbrace 
1,2,\dotsc,K \rbrace$ such that $\alpha_{k} \! = \! \infty$, $\lbrace 
\mathstrut k^{\prime} \! \in \! \lbrace 1,2,\dotsc,K \rbrace; 
\alpha_{k^{\prime}} \! \neq \! \alpha_{k} \! = \! \infty \rbrace$ is the 
ordered disjoint union of singletons, that is, $\cup_{q=1}^{K-1} \lbrace 
i(q)_{k_{q}} \rbrace$, with $k_{q} \! = \! 1$, $q \! = \! 1,2,\dotsc,K \! - \! 
1$, $1 \! \leqslant \! i(1)_{1} \! < \! i(2)_{1} \! < \! \dotsb \! < \! i(K \! 
- \! 1)_{1} \! \leqslant \! K$, and $\lbrace i(q)_{k_{q}} \rbrace \cap \lbrace 
i(r)_{k_{r}} \rbrace \! = \! \varnothing$ $\forall$ $q \! \neq \! r \! \in 
\! \lbrace 1,2,\dotsc,K \! - \! 1 \rbrace$, which induces, on the pole set 
$\lbrace \alpha_{1},\alpha_{2},\dotsc,\alpha_{K} \rbrace$, the following 
disjoint ordering, $\lbrace \mathstrut \alpha_{k^{\prime}}, \, k^{\prime} 
\! \in \! \lbrace 1,2,\dotsc,K \rbrace; \, \alpha_{k^{\prime}} \! \neq \! 
\alpha_{k} \! = \! \infty \rbrace \! := \! \cup_{q=1}^{K-1} \lbrace 
\alpha_{i(q)_{k_{q}}} \rbrace$, with $\alpha_{i(1)_{1}} \! \prec \! 
\alpha_{i(2)_{1}} \! \prec \! \dotsb \! \prec \! \alpha_{i(K-1)_{1}}$, $\# 
\lbrace \alpha_{i(q)_{k_{q}}} \rbrace \! = \! k_{q} \! = \! 1$, $q \! = \! 
1,2,\dotsc,K \! - \! 1$, $\lbrace \alpha_{i(q)_{k_{q}}} \rbrace \cap \lbrace 
\alpha_{i(r)_{k_{r}}} \rbrace \! = \! \varnothing$ $\forall$ $q \! \neq \! r 
\! \in \! \lbrace 1,2,\dotsc,K \! - \! 1 \rbrace$, and $\# \lbrace \mathstrut 
\alpha_{k^{\prime}}, \, k^{\prime} \! \in \! \lbrace 1,2,\dotsc,K \rbrace; 
\, \alpha_{k^{\prime}} \! \neq \! \alpha_{k} \! = \! \infty \rbrace \! = \! 
\sum_{q=1}^{K-1}k_{q} \! = \! K \! - \! 1$.}
\begin{gather*}
\alpha_{i(q)_{1}} \! = \! \alpha_{i(q)_{2}} \! = \! \dotsb \! = \! \alpha_{i
(q)_{k_{q}}} \, (\neq \! \infty), \quad q \! = \! 1,2,\dotsc,\mathfrak{s} \! 
- \! 1, \\
\# \lbrace \alpha_{i(q)_{1}},\alpha_{i(q)_{2}},\dotsc,\alpha_{i(q)_{k_{q}}} 
\rbrace \! = \! k_{q}, \quad q \! = \! 1,2,\dotsc,\mathfrak{s} \! - \! 1.
\end{gather*}
In order to illustrate this notation, consider the $K \! = \! 7$ pole set 
$\lbrace \alpha_{1},\alpha_{2},\alpha_{3},\alpha_{4},\alpha_{5},\alpha_{6},
\alpha_{7} \rbrace \! = \! \lbrace 0,1,\infty,1,\sqrt{2},\pi,\infty \rbrace$ 
for which $\mathfrak{s} \! = \! 5$ and $\alpha_{k} \! = \! \infty$, 
$k \! = \! 3,7$:
\begin{enumerate}
\item[(i)] \fbox{$k \! = \! 3$}
\begin{align*}
\lbrace \mathstrut k^{\prime} \! \in \! \lbrace 1,2,\dotsc,7 \rbrace; \, 
\alpha_{k^{\prime}} \! \neq \! \alpha_{3} \! = \! \infty \rbrace =& \, 
\lbrace 1,2,4,5,6 \rbrace \! = \! \lbrace 1 \rbrace \cup \lbrace 2,4 
\rbrace \cup \lbrace 5 \rbrace \cup \lbrace 6 \rbrace \\
:=& \, \lbrace i(1)_{k_{1}} \rbrace \cup \lbrace i(2)_{1},i(2)_{k_{2}} 
\rbrace \cup \lbrace i(3)_{k_{3}} \rbrace \cup \lbrace i(4)_{k_{4}} 
\rbrace \, \Rightarrow \\
k_{1} \! = \! 1, \, \, i(1)_{1} \! = \! 1, \, \, k_{2} \! = \! 2, \, \, i(2)_{1} 
\! = \! 2, \, \, &i(2)_{2} \! = \! 4, \, \, k_{3} \! = \! 1, \, \, i(3)_{1} \! 
= \! 5, \, \, k_{4} \! = \! 1, \, \, i(4)_{1} \! = \! 6
\end{align*}
(note that $1 \! \leqslant \! i(1)_{1} \! \leqslant \! 7$, $1 \! \leqslant \! 
i(2)_{1} \! < \! i(2)_{2} \! \leqslant \! 7$, $1 \! \leqslant \! i(3)_{1} \! 
\leqslant \! 7$, $1 \! \leqslant \! i(4)_{1} \! \leqslant \! 7$, $\lbrace 
i(1)_{1} \rbrace \cap \lbrace i(2)_{1},i(2)_{2} \rbrace \! = \! \lbrace 
i(1)_{1} \rbrace \cap \lbrace i(3)_{1} \rbrace \! = \! \lbrace i(1)_{1} 
\rbrace \cap \lbrace i(4)_{1} \rbrace \! = \! \lbrace i(2)_{1},i(2)_{2} 
\rbrace \cap \lbrace i(3)_{1} \rbrace \! = \! \lbrace i(2)_{1},i(2)_{2} 
\rbrace \cap \lbrace i(4)_{1} \rbrace \! = \! \lbrace i(3)_{1} \rbrace \cap 
\lbrace i(4)_{1} \rbrace \! = \! \varnothing$, and $\# \lbrace \mathstrut 
k^{\prime} \! \in \! \lbrace 1,2,\dotsc,7 \rbrace; \, \alpha_{k^{\prime}} \! 
\neq \! \alpha_{3} \! = \! \infty \rbrace \! = \! \sum_{q=1}^{4}k_{q} \! = \! 
1 \! + \! 2 \! + \! 1 \! + \! 1 = \! 5 \! = \! K \! - \! \gamma_{3} \! = \! 7 
\! - \! 2)$, which induces the ordering (on the `residual' pole set)
\begin{align*}
\lbrace \mathstrut \alpha_{k^{\prime}}, \, k^{\prime} \! \in \! \lbrace 
1,2,\dotsc,7 \rbrace; \, \alpha_{k^{\prime}} \! \neq \! \alpha_{3} \! = \! 
\infty \rbrace :=& \, \lbrace \alpha_{i(1)_{k_{1}}} \rbrace \cup \lbrace 
\alpha_{i(2)_{1}},\alpha_{i(2)_{k_{2}}} \rbrace \cup \lbrace 
\alpha_{i(3)_{k_{3}}} \rbrace \cup \lbrace \alpha_{i(4)_{k_{4}}} \rbrace \\
=& \,\lbrace \alpha_{1} \rbrace \cup \lbrace \alpha_{2},\alpha_{4} \rbrace 
\cup \lbrace \alpha_{5} \rbrace \cup \lbrace \alpha_{6} \rbrace \! = \! 
\lbrace 0 \rbrace \cup \lbrace 1,1 \rbrace \cup \lbrace \! \sqrt{2} \rbrace 
\cup \lbrace \pi \rbrace
\end{align*}
(note that $\alpha_{i(2)_{1}} \! \prec \! \alpha_{i(2)_{2}}$, $\lbrace 
\alpha_{i(1)_{1}} \rbrace \cap \lbrace \alpha_{i(2)_{1}},\alpha_{i(2)_{2}} 
\rbrace \! = \! \lbrace \alpha_{i(1)_{1}} \rbrace \cap \lbrace 
\alpha_{i(3)_{1}} \rbrace \! = \! \lbrace \alpha_{i(1)_{1}} \rbrace \cap 
\lbrace \alpha_{i(4)_{1}} \rbrace \! = \! \lbrace \alpha_{i(2)_{1}},
\alpha_{i(2)_{2}} \rbrace \cap \lbrace \alpha_{i(3)_{1}} \rbrace \! = 
\! \lbrace \alpha_{i(2)_{1}},\alpha_{i(2)_{2}} \rbrace \cap \lbrace 
\alpha_{i(4)_{1}} \rbrace \! = \! \lbrace \alpha_{i(3)_{1}} \rbrace \cap 
\lbrace \alpha_{i(4)_{1}} \rbrace \! = \! \varnothing$, $\# \lbrace 
\alpha_{i(1)_{1}} \rbrace \! = \! k_{1} \! = \! 1$, $\# \lbrace 
\alpha_{i(2)_{1}},\alpha_{i(2)_{2}} \rbrace \! = \! k_{2} \! = \! 2$, $\# 
\lbrace \alpha_{i(3)_{1}} \rbrace \! = \! k_{3} \! = \! 1$, and $\# \lbrace 
\alpha_{i(4)_{1}} \rbrace \! = \! k_{4} \! = \! 1)$;
\item[(ii)] \fbox{$k \! = \! 7$}
\begin{align*}
\lbrace \mathstrut k^{\prime} \! \in \! \lbrace 1,2,\dotsc,7 \rbrace; \, 
\alpha_{k^{\prime}} \! \neq \! \alpha_{7} \! = \! \infty \rbrace =& \, 
\lbrace 1,2,4,5,6 \rbrace \! = \! \lbrace 1 \rbrace \cup \lbrace 2,4 \rbrace 
\cup \lbrace 5 \rbrace \cup \lbrace 6 \rbrace \\
:=& \, \lbrace i(1)_{k_{1}} \rbrace \cup \lbrace i(2)_{1},i(2)_{k_{2}} 
\rbrace \cup \lbrace i(3)_{k_{3}} \rbrace \cup \lbrace i(4)_{k_{4}} \rbrace \, 
\Rightarrow \\
k_{1} \! = \! 1, \, \, i(1)_{1} \! = \! 1, \, \, k_{2} \! = \! 2, \, \, i(2)_{1} 
\! = \! 2, \, \, &i(2)_{2} \! = \! 4, \, \, k_{3} \! = \! 1, \, \, i(3)_{1} \! 
= \! 5, \, \, k_{4} \! = \! 1, \, \, i(4)_{1} \! = \! 6
\end{align*}
(note that $1 \! \leqslant \! i(1)_{1} \! \leqslant \! 7$, $1 \! \leqslant \! 
i(2)_{1} \! < \! i(2)_{2} \! \leqslant \! 7$, $1 \! \leqslant \! i(3)_{1} \! 
\leqslant \! 7$, $1 \! \leqslant \! i(4)_{1} \! \leqslant \! 7$, $\lbrace 
i(1)_{1} \rbrace \cap \lbrace i(2)_{1},i(2)_{2} \rbrace \! = \! \lbrace 
i(1)_{1} \rbrace \cap \lbrace i(3)_{1} \rbrace \! = \! \lbrace i(1)_{1} 
\rbrace \cap \lbrace i(4)_{1} \rbrace \! = \! \lbrace i(2)_{1},i(2)_{2} 
\rbrace \cap \lbrace i(3)_{1} \rbrace \! = \! \lbrace i(2)_{1},i(2)_{2} 
\rbrace \cap \lbrace i(4)_{1} \rbrace \! = \! \lbrace i(3)_{1} \rbrace \cap 
\lbrace i(4)_{1} \rbrace \! = \! \varnothing$, and $\# \lbrace \mathstrut 
k^{\prime} \! \in \! \lbrace 1,2,\dotsc,7 \rbrace; \, \alpha_{k^{\prime}} \! 
\neq \! \alpha_{7} \! = \! \infty \rbrace \! = \! \sum_{q=1}^{4}k_{q} \! = \! 
1 \! + \! 2 \! + \! 1 \! + \! 1 = \! 5 \! = \! K \! - \! \gamma_{7} \! = \! 
7 \! - \! 2)$, which induces the ordering (on the `residual' pole set)
\begin{align*}
\lbrace \mathstrut \alpha_{k^{\prime}}, \, k^{\prime} \! \in \! \lbrace 1,2,
\dotsc,7 \rbrace; \, \alpha_{k^{\prime}} \! \neq \! \alpha_{7} \! = \! 
\infty \rbrace :=& \, \lbrace \alpha_{i(1)_{k_{1}}} \rbrace \cup \lbrace 
\alpha_{i(2)_{1}},\alpha_{i(2)_{k_{2}}} \rbrace \cup \lbrace 
\alpha_{i(3)_{k_{3}}} \rbrace \cup \lbrace \alpha_{i(4)_{k_{4}}} \rbrace \\
=& \, \lbrace \alpha_{1} \rbrace \cup \lbrace \alpha_{2},\alpha_{4} \rbrace 
\cup \lbrace \alpha_{5} \rbrace \cup \lbrace \alpha_{6} \rbrace \! = \! 
\lbrace 0 \rbrace \cup \lbrace 1,1 \rbrace \cup \lbrace \! \sqrt{2} \rbrace 
\cup \lbrace \pi \rbrace
\end{align*}
(note that $\alpha_{i(2)_{1}} \! \prec \! \alpha_{i(2)_{2}}$, $\lbrace 
\alpha_{i(1)_{1}} \rbrace \cap \lbrace \alpha_{i(2)_{1}},\alpha_{i(2)_{2}} 
\rbrace \! = \! \lbrace \alpha_{i(1)_{1}} \rbrace \cap \lbrace 
\alpha_{i(3)_{1}} \rbrace \! = \! \lbrace \alpha_{i(1)_{1}} \rbrace \cap 
\lbrace \alpha_{i(4)_{1}} \rbrace \! = \! \lbrace \alpha_{i(2)_{1}},
\alpha_{i(2)_{2}} \rbrace \cap \lbrace \alpha_{i(3)_{1}} \rbrace \! = \! 
\lbrace \alpha_{i(2)_{1}},\alpha_{i(2)_{2}} \rbrace \cap \lbrace 
\alpha_{i(4)_{1}} \rbrace \! = \! \lbrace \alpha_{i(3)_{1}} \rbrace \cap 
\lbrace \alpha_{i(4)_{1}} \rbrace \! = \! \varnothing$, $\# \lbrace 
\alpha_{i(1)_{1}} \rbrace \! = \! k_{1} \! = \! 1$, $\# \lbrace 
\alpha_{i(2)_{1}},\alpha_{i(2)_{2}} \rbrace \! = \! k_{2} \! = \! 2$, $\# 
\lbrace \alpha_{i(3)_{1}} \rbrace \! = \! k_{3} \! = \! 1$, and $\# \lbrace 
\alpha_{i(4)_{1}} \rbrace \! = \! k_{4} \! = \! 1)$.
\end{enumerate}
This concludes the example.

In order to proceed further, it will be important to know if, within each of 
the---decomposed---index sets $\lbrace i(q)_{1},i(q)_{2},\dotsc,i(q)_{k_{q}} 
\rbrace$, $q \! = \! 1,2,\dotsc,\mathfrak{s} \! - \! 1$, there is a positive 
integer less than $k$. Towards this end, the following notational preamble is 
helpful.
\begin{aaaa} \label{def1.2.2} 
\emph{For $k \! = \! 1,2,\dotsc,K$ and a given set of positive integers $i_{1},
i_{2},\dotsc,i_{M}$, let $\operatorname{ind} \lbrace i_{1},i_{2},\dotsc,i_{M} 
\vert k \rbrace$ denote the largest positive integer $\hat{\jmath}$, if it 
exists, {}from the collection $i_{1},i_{2},\dotsc,i_{M}$ that is less than 
$k$$;$ e.g., $\operatorname{ind} \lbrace 1,3,7,11 \vert 8 \rbrace$ $\Rightarrow$ 
$\hat{\jmath} \! = \! 7$, and $\operatorname{ind} \lbrace 5,6,9 \vert 3 
\rbrace$ $\Rightarrow$ no such $\hat{\jmath}$ exists.}
\end{aaaa}
For $k \! \in \! \lbrace 1,2,\dotsc,K \rbrace$ such that $\alpha_{k} \! = \! 
\infty$, let
\begin{equation*}
\mathfrak{J}_{q}(k) \! := \! \lbrace \operatorname{ind} \lbrace i(q)_{1},
i(q)_{2},\dotsc,i(q)_{k_{q}} \vert k \rbrace \rbrace, \quad q \! = \! 1,2,
\dotsc,\mathfrak{s} \! - \! 1:
\end{equation*}
if, for a given index set $\lbrace i(q)_{1},i(q)_{2},\dotsc,i(q)_{k_{q}} \rbrace$, 
$q \! = \! 1,2,\dotsc,\mathfrak{s} \! - \! 1$, there does not exist a positive 
integer {}\footnote{This positive integer, if it exists, depends not only on $k$, 
but on $q$, too; however, for notational simplicity, this additional $q$ 
dependence is also suppressed.} $\hat{\jmath}$ such that $\hat{\jmath} 
\! = \! \operatorname{ind} \lbrace i(q)_{1},i(q)_{2},\dotsc,i(q)_{k_{q}} \vert 
k \rbrace$, then put $\mathfrak{J}_{q}(k) \! = \! \varnothing$; otherwise, 
denote by $\widetilde{m}_{q}(k)$, $q \! \in \! \lbrace 1,2,\dotsc,\mathfrak{s} 
\! - \! 1 \rbrace$, the unique element of the (non-empty) set $\mathfrak{J}_{q}
(k)$.\footnote{Note that $\mathfrak{J}_{q}(k)$ is either empty or a singleton.} 
For $n \! \in \! \mathbb{N}$ and $k \! \in \! \lbrace 1,2,\dotsc,K \rbrace$ 
such that $\alpha_{k} \! = \! \infty$, set, with the above orderings and 
definitions,
\begin{equation*}
\varkappa_{nk \tilde{k}_{q}} \! := \! 
\begin{cases}
(n \! - \! 1) \gamma_{i(q)_{k_{q}}}, &\text{$\mathfrak{J}_{q}(k) \! = \! 
\varnothing, \quad q \! \in \! \lbrace 1,2,\dotsc,\mathfrak{s} \! - \! 1 
\rbrace$,} \\
(n \! - \! 1) \gamma_{\widetilde{m}_{q}(k)} \! + \! 
\varrho_{\widetilde{m}_{q}(k)}, &\text{$\mathfrak{J}_{q}(k) \! \neq \! 
\varnothing, \quad q \! \in \! \lbrace 1,2,\dotsc,\mathfrak{s} \! - \! 1 
\rbrace$,}
\end{cases}
\end{equation*}
where $\varkappa_{nk \tilde{k}_{q}} \colon \mathbb{N} \times \lbrace 
1,2,\dotsc,K \rbrace \! \to \! \mathbb{N}_{0}$, $q \! = \! 1,2,\dotsc,
\mathfrak{s} \! - \! 1$, is the \emph{residual multiplicity} of the pole 
$\alpha_{i(q)_{k_{q}}}$ in the repeated pole sequence
\begin{equation*}
\mathcal{P}_{n,k} \! := \! \lbrace \overset{1}{\underbrace{\alpha_{1},
\alpha_{2},\dotsc,\alpha_{K}}_{K}} \rbrace \cup \dotsb \cup \lbrace 
\overset{n-1}{\underbrace{\alpha_{1},\alpha_{2},\dotsc,\alpha_{K}}_{K}} 
\rbrace \cup \lbrace \overset{n}{\underbrace{\alpha_{1},\alpha_{2},
\dotsc,\alpha_{k}}_{k}} \rbrace,
\end{equation*}
that is, for $n \! \in \! \mathbb{N}$ and $k \! \in \! \lbrace 1,2,\dotsc,
K \rbrace$ such that $\alpha_{k} \! = \! \infty$, as all occurrences of the 
pole (the point at infinity) $\alpha_{k} \! = \! \infty$ are excised {}from 
the repeated pole sequence $\mathcal{P}_{n,k}$, where the multiplicity, 
or number of occurrences, of the pole $\alpha_{k} \! = \! \infty$ is 
$\varkappa_{nk} \! = \! (n \! - \! 1) \gamma_{k} \! + \! \varrho_{k}$, one 
is left with the `residual' pole set (via the above induced ordering on the 
poles, as one walks across the pole sequence {}from left to right)
\begin{align*}
\mathcal{P}_{n,k} \setminus \lbrace \underbrace{\alpha_{k},\alpha_{k},\dotsc,
\alpha_{k}}_{\varkappa_{nk}} \rbrace =& \, \mathcal{P}_{n,k} \setminus 
\lbrace \underbrace{\infty,\infty,\dotsc,\infty}_{\varkappa_{nk}} \rbrace 
:= \lbrace \underbrace{\alpha_{i(1)_{k_{1}}},\alpha_{i(1)_{k_{1}}},\dotsc,
\alpha_{i(1)_{k_{1}}}}_{\varkappa_{nk \tilde{k}_{1}}} \rbrace \cup 
\lbrace \underbrace{\alpha_{i(2)_{k_{2}}},\alpha_{i(2)_{k_{2}}},\dotsc,
\alpha_{i(2)_{k_{2}}}}_{\varkappa_{nk \tilde{k}_{2}}} \rbrace \cup \dotsb \\
\dotsb& \, \cup \lbrace \underbrace{\alpha_{i(\mathfrak{s}-1)_{k_{\mathfrak{s}
-1}}},\alpha_{i(\mathfrak{s}-1)_{k_{\mathfrak{s}-1}}},\dotsc,\alpha_{i
(\mathfrak{s}-1)_{k_{\mathfrak{s}-1}}}}_{\varkappa_{nk \tilde{k}_{\mathfrak{s}
-1}}} \rbrace = \bigcup_{q=1}^{\mathfrak{s}-1} \lbrace \underbrace{
\alpha_{i(q)_{k_{q}}},\alpha_{i(q)_{k_{q}}},\dotsc,\alpha_{i(q)_{k_{q}}}}_{
\varkappa_{nk \tilde{k}_{q}}} \rbrace,
\end{align*}
where the number of times the pole $\alpha_{i(q)_{k_{q}}}$ $(\neq \! 
\alpha_{k} \! = \! \infty)$ occurs is (its multiplicity) 
$\varkappa_{nk \tilde{k}_{q}}$, $q \! = \! 1,2,\dotsc,\mathfrak{s} \! - \! 
1$.\footnote{Note: for $n \! = \! 1$ and $k \! \in \! \lbrace 1,2,\dotsc,K 
\rbrace$ such that $\alpha_{k} \! = \! \infty$, it can happen that, for some 
values of $q \! \in \! \lbrace 1,2,\dotsc,\mathfrak{s} \! - \! 1 \rbrace$, 
$\varkappa_{1k \tilde{k}_{q}} \! = \! 0$, in which case, one defines $\lbrace 
\alpha_{i(q)_{k_{q}}},\alpha_{i(q)_{k_{q}}},\dotsc,\alpha_{i(q)_{k_{q}}} 
\rbrace \! := \! \varnothing$, $q \! \in \! \lbrace 1,2,\dotsc,\mathfrak{s} 
\! - \! 1 \rbrace$; however, for $(\mathbb{N} \! \ni)$ $n \! \geqslant \! 
2$ and $k \! \in \! \lbrace 1,2,\dotsc,K \rbrace$ such that $\alpha_{k} 
\! = \! \infty$, $\varkappa_{nk \tilde{k}_{q}} \! \geqslant \! 1$, 
$q \! \in \! \lbrace 1,2,\dotsc,\mathfrak{s} \! - \! 1 \rbrace$.}

For $n \! \in \! \mathbb{N}$ and $k \! \in \! \lbrace 1,2,\dotsc,K \rbrace$ 
such that $\alpha_{k} \! = \! \infty$, a counting-of-residual-multiplicities 
argument gives rise to the following ordered sum formula:\footnote{It should 
be noted that, in this context, as well as for further reference, for any 
function $g \colon \mathbb{N} \times \lbrace 1,2,\dotsc,K \rbrace \times 
\overline{\mathbb{C}} \! \to \! \mathbb{C}$, the following ordered sums 
and products will be used, with little or no further reference, throughout 
this monograph: $\sum_{q=1}^{\mathfrak{s}-1}(g(n,k,\tilde{k}_{q}))^{
\varkappa_{nk \tilde{k}_{q}}} \! := \! (g(n,k,\tilde{k}_{1}))^{\varkappa_{n
k \tilde{k}_{1}}} \! + \! (g(n,k,\tilde{k}_{2}))^{\varkappa_{nk \tilde{k}_{2}}} 
\! + \! \dotsb \! + \! (g(n,k,\tilde{k}_{\mathfrak{s}-1}))^{\varkappa_{nk 
\tilde{k}_{\mathfrak{s}-1}}}$ and $\prod_{q=1}^{\mathfrak{s}-1}(g(n,k,
\tilde{k}_{q}))^{\varkappa_{nk \tilde{k}_{q}}} \! := \! (g(n,k,\tilde{k}_{1}))^{
\varkappa_{nk \tilde{k}_{1}}}(g(n,k,\tilde{k}_{2}))^{\varkappa_{nk 
\tilde{k}_{2}}} \dotsb (g(n,k,\tilde{k}_{\mathfrak{s}-1}))^{\varkappa_{nk 
\tilde{k}_{\mathfrak{s}-1}}}$.}
\begin{equation} \label{infcount} 
\sum_{q=1}^{\mathfrak{s}-1} \varkappa_{nk \tilde{k}_{q}} \! := \! 
\varkappa_{nk \tilde{k}_{1}} \! + \! \varkappa_{nk \tilde{k}_{2}} 
\! + \! \dotsb \! + \! \varkappa_{nk \tilde{k}_{\mathfrak{s}-1}} 
\! = \! (n \! - \! 1)K \! + \! k \! - \! \varkappa_{nk},
\end{equation}
whence $\sum_{q=1}^{\mathfrak{s}-1} \varkappa_{nk \tilde{k}_{q}} \! + \! 
\varkappa_{nk} \! = \! (n \! - \! 1)K \! + \! k$.

In order to illustrate the latter notation, consider, again, the $K \! = \! 7$ 
pole set $\lbrace \alpha_{1},\alpha_{2},\alpha_{3},\alpha_{4},\alpha_{5},
\alpha_{6},\alpha_{7} \rbrace \! = \! \lbrace 0,1,\infty,1,\sqrt{2},\pi,\infty 
\rbrace$ for which $\mathfrak{s} \! = \! 5$ and $\alpha_{k} \! = \! \infty$, 
$k \! = \! 3,7$:
\begin{enumerate}
\item[(i)] \fbox{$k \! = \! 3$} 
\begin{gather*}
\mathfrak{J}_{1}(3) \! := \! \lbrace \operatorname{ind} \lbrace i(1)_{k_{1}} 
\vert 3 \rbrace \rbrace \! = \! \lbrace \operatorname{ind} \lbrace 1 \vert 3 
\rbrace \rbrace \! = \! \lbrace 1 \rbrace \Rightarrow \widetilde{m}_{1}(3) \! 
= \! 1, \\
\mathfrak{J}_{2}(3) \! := \! \lbrace \operatorname{ind} \lbrace i(2)_{1},
i(2)_{k_{2}} \vert 3 \rbrace \rbrace \! = \! \lbrace \operatorname{ind} 
\lbrace 2,4 \vert 3 \rbrace \rbrace \! = \! \lbrace 2 \rbrace \Rightarrow 
\widetilde{m}_{2}(3) \! = \! 2, \\
\mathfrak{J}_{3}(3) \! := \! \lbrace \operatorname{ind} \lbrace i(3)_{k_{3}} 
\vert 3 \rbrace \rbrace \! = \! \lbrace \operatorname{ind} \lbrace 5 \vert 3 
\rbrace \rbrace \! = \! \varnothing, \\
\mathfrak{J}_{4}(3) \! := \! \lbrace \operatorname{ind} \lbrace i(4)_{k_{4}} 
\vert 3 \rbrace \rbrace \! = \! \lbrace \operatorname{ind} \lbrace 6 \vert 3 
\rbrace \rbrace \! = \! \varnothing,
\end{gather*}
hence
\begin{gather*}
\varkappa_{n3 \tilde{k}_{1}} \! = \! (n \! - \! 1) \gamma_{\widetilde{m}_{1}
(3)} \! + \! \varrho_{\widetilde{m}_{1}(3)} \! = \! (n \! - \! 1) \gamma_{1} 
\! + \! \varrho_{1} \! = \! (n \! - \! 1) \! + \! 1 \! = \! n, \\
\varkappa_{n3 \tilde{k}_{2}} \! = \! (n \! - \! 1) \gamma_{\widetilde{m}_{2}
(3)} \! + \! \varrho_{\widetilde{m}_{2}(3)} \! = \! (n \! - \! 1) \gamma_{2} 
\! + \! \varrho_{2} \! = \! 2(n \! - \! 1) \! + \! 1 \! = \! 2n \! - \! 1, \\
\varkappa_{n3 \tilde{k}_{3}} \! = \! (n \! - \! 1) \gamma_{i(3)_{k_{3}}} \! = 
\! (n \! - \! 1) \gamma_{5} \! = \! n \! - \! 1, \\
\varkappa_{n3 \tilde{k}_{4}} \! = \! (n \! - \! 1) \gamma_{i(4)_{k_{4}}} \! = 
\! (n \! - \! 1) \gamma_{6} \! = \! n \! - \! 1,
\end{gather*}
that is, as one moves {}from left to right across the repeated pole sequence
\begin{align*}
\mathcal{P}_{n,3} =& \, \lbrace \overset{1}{\underbrace{\alpha_{1},\alpha_{2},
\dotsc,\alpha_{7}}_{7}} \rbrace \cup \dotsb \cup \lbrace \overset{n-1}{
\underbrace{\alpha_{1},\alpha_{2},\dotsc,\alpha_{7}}_{7}} \rbrace \cup \lbrace 
\overset{n}{\underbrace{\alpha_{1},\alpha_{2},\alpha_{3}}_{3}} \rbrace \\
=& \, \lbrace \overset{1}{\underbrace{0,1,\infty,1,\sqrt{2},\pi,\infty}_{7}} 
\rbrace \cup \dotsb \cup \lbrace \overset{n-1}{\underbrace{0,1,\infty,1,
\sqrt{2},\pi,\infty}_{7}} \rbrace \cup \lbrace \overset{n}{\underbrace{0,1,
\infty}_{3}} \rbrace
\end{align*}
and removes all occurrences of the pole $\alpha_{3} \! = \! \infty$, which 
occurs $\varkappa_{n3} \! = \! (n \! - \! 1) \gamma_{3} \! + \! \varrho_{3} 
\! = \! 2(n \! - \! 1) \! + \! 1 \! = \! 2n \! - \! 1$ times, one is left with 
the residual pole set (via the above induced ordering)
\begin{align*}
\mathcal{P}_{n,3} \setminus \lbrace \underbrace{\alpha_{3},\alpha_{3},
\dotsc,\alpha_{3}}_{\varkappa_{n3}} \rbrace =& \, \mathcal{P}_{n,3} 
\setminus \lbrace \underbrace{\infty,\infty,\dotsc,\infty}_{2n-1} \rbrace 
\! = \bigcup_{q=1}^{4} \lbrace \underbrace{\alpha_{i(q)_{k_{q}}},
\alpha_{i(q)_{k_{q}}},\dotsc,\alpha_{i(q)_{k_{q}}}}_{\varkappa_{n3 
\tilde{k}_{q}}} \rbrace \\
:=& \, \lbrace \underbrace{\alpha_{i(1)_{k_{1}}},\alpha_{i(1)_{k_{1}}},\dotsc,
\alpha_{i(1)_{k_{1}}}}_{\varkappa_{n3 \tilde{k}_{1}}} \rbrace \cup \lbrace 
\underbrace{\alpha_{i(2)_{k_{2}}},\alpha_{i(2)_{k_{2}}},\dotsc,
\alpha_{i(2)_{k_{2}}}}_{\varkappa_{n3 \tilde{k}_{2}}} \rbrace \\
\cup& \, \lbrace \underbrace{\alpha_{i(3)_{k_{3}}},\alpha_{i(3)_{k_{3}}},
\dotsc,\alpha_{i(3)_{k_{3}}}}_{\varkappa_{n3 \tilde{k}_{3}}} \rbrace \cup 
\lbrace \underbrace{\alpha_{i(4)_{k_{4}}},\alpha_{i(4)_{k_{4}}},\dotsc,
\alpha_{i(4)_{k_{4}}}}_{\varkappa_{n3 \tilde{k}_{4}}} \rbrace \\
=& \, \lbrace \underbrace{0,0,\dotsc,0}_{n} \rbrace \cup \lbrace \underbrace{
1,1,\dotsc,1}_{2n-1} \rbrace \cup \lbrace \underbrace{\sqrt{2},\sqrt{2},
\dotsc,\sqrt{2}}_{n-1} \rbrace \cup \lbrace \underbrace{\pi,\pi,\dotsc,
\pi}_{n-1} \rbrace,
\end{align*}
where the number of times the pole $\alpha_{i(1)_{k_{1}}} \! = \! 0$ $(\neq 
\! \alpha_{3} \! = \! \infty)$ occurs is $\varkappa_{n3 \tilde{k}_{1}} \! = \! 
n$, the number of times the pole $\alpha_{i(2)_{k_{2}}} \! = \! 1$ $(\neq \! 
\alpha_{3} \! = \! \infty)$ occurs is $\varkappa_{n3 \tilde{k}_{2}} \! = \! 
2n \! - \! 1$, the number of times the pole $\alpha_{i(3)_{k_{3}}} \! = \! 
\sqrt{2}$ $(\neq \! \alpha_{3} \! = \! \infty)$ occurs is $\varkappa_{n3 
\tilde{k}_{3}} \! = \! n \! - \! 1$, and the number of times the pole 
$\alpha_{i(4)_{k_{4}}} \! = \! \pi$ $(\neq \! \alpha_{3} \! = \! \infty)$ 
occurs is $\varkappa_{n3 \tilde{k}_{4}} \! = \! n \! - \! 1$.\footnote{Note: 
for $n \! = \! 1$, since $\varkappa_{13 \tilde{k}_{3}} \! = \! \varkappa_{13 
\tilde{k}_{4}} \! = \! 0$, one sets, as per the convention above, $\lbrace 
\alpha_{i(q)_{k_{q}}},\alpha_{i(q)_{k_{q}}},\dotsc,\alpha_{i(q)_{k_{q}}} 
\rbrace \! := \! \varnothing$, $q \! = \! 3,4$, in which case, as 
$\varkappa_{13 \tilde{k}_{1}} \! = \! \varkappa_{13 \tilde{k}_{2}} \! = \! 
1$ and $\varkappa_{13} \! = \! \varrho_{3} \! = \! 1$, $\mathcal{P}_{1,3} 
\setminus \lbrace \alpha_{3} \rbrace \! = \! \mathcal{P}_{1,3} \setminus 
\lbrace \infty \rbrace \! = \! \lbrace 0 \rbrace \cup \lbrace 1 \rbrace 
\cup \varnothing \cup \varnothing \! = \! \lbrace 0 \rbrace \cup \lbrace 
1 \rbrace$.} In this case, the ordered sum formula reads
\begin{equation*}
\sum_{q=1}^{4} \varkappa_{n3 \tilde{k}_{q}} \! := \! \varkappa_{n3 
\tilde{k}_{1}} \! + \! \varkappa_{n3 \tilde{k}_{2}} \! + \! \varkappa_{n3 
\tilde{k}_{3}} \! + \! \varkappa_{n3 \tilde{k}_{4}} \! = \! n \! + \! (2n \! 
- \! 1) \! + \! (n \! - \! 1) \! + \! (n \! - \! 1) \! = \! 5n \! - \! 3,
\end{equation*}
whence $\sum_{q=1}^{4} \varkappa_{n3 \tilde{k}_{q}} \! + \! \varkappa_{n3} 
\! = \! (5n \! - \! 3) \! + \! (2n \! - \! 1) \! = \! 7(n \! - \! 1) \! + \! 3$;
\item[(ii)] \fbox{$k \! = \! 7$}
\begin{gather*}
\mathfrak{J}_{1}(7) \! := \! \lbrace \operatorname{ind} \lbrace i(1)_{k_{1}} 
\vert 7 \rbrace \rbrace \! = \! \lbrace \operatorname{ind} \lbrace 1 \vert 7 
\rbrace \rbrace \! = \! \lbrace 1 \rbrace \Rightarrow \widetilde{m}_{1}(7) \! 
= \! 1, \\
\mathfrak{J}_{2}(7) \! := \! \lbrace \operatorname{ind} \lbrace i(2)_{1},
i(2)_{k_{2}} \vert 7 \rbrace \rbrace \! = \! \lbrace \operatorname{ind} 
\lbrace 2,4 \vert 7 \rbrace \rbrace \! = \! \lbrace 4 \rbrace \Rightarrow 
\widetilde{m}_{2}(7) \! = \! 4, \\
\mathfrak{J}_{3}(7) \! := \! \lbrace \operatorname{ind} \lbrace i(3)_{k_{3}} 
\vert 7 \rbrace \rbrace \! = \! \lbrace \operatorname{ind} \lbrace 5 \vert 7 
\rbrace \rbrace \! = \! \lbrace 5 \rbrace \Rightarrow \widetilde{m}_{3}(7) \! 
= \! 5, \\
\mathfrak{J}_{4}(7) \! := \! \lbrace \operatorname{ind} \lbrace i(4)_{k_{4}} 
\vert 7 \rbrace \rbrace \! = \! \lbrace \operatorname{ind} \lbrace 6 \vert 7 
\rbrace \rbrace \! = \! \lbrace 6 \rbrace \Rightarrow \widetilde{m}_{4}(7) \! 
= \! 6,
\end{gather*}
hence
\begin{gather*}
\varkappa_{n7 \tilde{k}_{1}} \! = \! (n \! - \! 1) \gamma_{\widetilde{m}_{1}
(7)} \! + \! \varrho_{\widetilde{m}_{1}(7)} \! = \! (n \! - \! 1) \gamma_{1} 
\! + \! \varrho_{1} \! = \! (n \! - \! 1) \! + \! 1 \! = \! n, \\
\varkappa_{n7 \tilde{k}_{2}} \! = \! (n \! - \! 1) \gamma_{\widetilde{m}_{2}
(7)} \! + \! \varrho_{\widetilde{m}_{2}(7)} \! = \! (n \! - \! 1) \gamma_{4} 
\! + \! \varrho_{4} \! = \! 2(n \! - \! 1) \! + \! 2 \! = \! 2n, \\
\varkappa_{n7 \tilde{k}_{3}} \! = \! (n \! - \! 1) \gamma_{\widetilde{m}_{3}
(7)} \! + \! \varrho_{\widetilde{m}_{3}(7)} \! = \! (n \! - \! 1) \gamma_{5} 
\! + \! \varrho_{5} \! = \! (n \! - \! 1) \! + \! 1 \! = \! n, \\
\varkappa_{n7 \tilde{k}_{4}} \! = \! (n \! - \! 1) \gamma_{\widetilde{m}_{4}
(7)} \! + \! \varrho_{\widetilde{m}_{4}(7)} \! = \! (n \! - \! 1) \gamma_{6} 
\! + \! \varrho_{6} \! = \! (n \! - \! 1) \! + \! 1 \! = \! n,
\end{gather*}
that is, as one moves {}from left to right across the repeated pole sequence
\begin{align*}
\mathcal{P}_{n,7} =& \, \lbrace \overset{1}{\underbrace{\alpha_{1},\alpha_{2},
\dotsc,\alpha_{7}}_{7}} \rbrace \cup \dotsb \cup \lbrace \overset{n-1}{
\underbrace{\alpha_{1},\alpha_{2},\dotsc,\alpha_{7}}_{7}} \rbrace \cup \lbrace 
\overset{n}{\underbrace{\alpha_{1},\alpha_{2},\dotsc,\alpha_{7}}_{7}} \rbrace 
\\
=& \, \lbrace \overset{1}{\underbrace{0,1,\infty,1,\sqrt{2},\pi,\infty}_{7}} 
\rbrace \cup \dotsb \cup \lbrace \overset{n-1}{\underbrace{0,1,\infty,1,
\sqrt{2},\pi,\infty}_{7}} \rbrace \cup \lbrace \overset{n}{\underbrace{0,1,
\infty,1,\sqrt{2},\pi,\infty}_{7}} \rbrace
\end{align*}
and removes all occurrences of the pole $\alpha_{7} \! = \! \infty$, which 
occurs $\varkappa_{n7} \! = \! (n \! - \! 1) \gamma_{7} \! + \! \varrho_{7} 
\! = \! 2(n \! - \! 1) \! + \! 2 \! = \! 2n$ times, one is left with the 
residual pole set (via the above induced ordering)
\begin{align*}
\mathcal{P}_{n,7} \setminus \lbrace \underbrace{\alpha_{7},\alpha_{7},\dotsc,
\alpha_{7}}_{\varkappa_{n7}} \rbrace =& \, \mathcal{P}_{n,7} \setminus \lbrace 
\underbrace{\infty,\infty,\dotsc,\infty}_{2n} \rbrace \! = \bigcup_{q=1}^{4} 
\lbrace \underbrace{\alpha_{i(q)_{k_{q}}},\alpha_{i(q)_{k_{q}}},\dotsc,
\alpha_{i(q)_{k_{q}}}}_{\varkappa_{n7 \tilde{k}_{q}}} \rbrace \\
:=& \, \lbrace \underbrace{\alpha_{i(1)_{k_{1}}},\alpha_{i(1)_{k_{1}}},\dotsc,
\alpha_{i(1)_{k_{1}}}}_{\varkappa_{n7 \tilde{k}_{1}}} \rbrace \cup \lbrace 
\underbrace{\alpha_{i(2)_{k_{2}}},\alpha_{i(2)_{k_{2}}},\dotsc,
\alpha_{i(2)_{k_{2}}}}_{\varkappa_{n7 \tilde{k}_{2}}} \rbrace \\
\cup& \, \lbrace \underbrace{\alpha_{i(3)_{k_{3}}},\alpha_{i(3)_{k_{3}}},
\dotsc,\alpha_{i(3)_{k_{3}}}}_{\varkappa_{n7 \tilde{k}_{3}}} \rbrace \cup 
\lbrace \underbrace{\alpha_{i(4)_{k_{4}}},\alpha_{i(4)_{k_{4}}},\dotsc,
\alpha_{i(4)_{k_{4}}}}_{\varkappa_{n7 \tilde{k}_{4}}} \rbrace \\
=& \, \lbrace \underbrace{0,0,\dotsc,0}_{n} \rbrace \cup \lbrace \underbrace{
1,1,\dotsc,1}_{2n} \rbrace \cup \lbrace \underbrace{\sqrt{2},\sqrt{2},\dotsc,
\sqrt{2}}_{n} \rbrace \cup \lbrace \underbrace{\pi,\pi,\dotsc,\pi}_{n} \rbrace,
\end{align*}
where the number of times the pole $\alpha_{i(1)_{k_{1}}} \! = \! 0$ $(\neq 
\! \alpha_{7} \! = \! \infty)$ occurs is $\varkappa_{n7 \tilde{k}_{1}} \! = \! 
n$, the number of times the pole $\alpha_{i(2)_{k_{2}}} \! = \! 1$ $(\neq \! 
\alpha_{7} \! = \! \infty)$ occurs is $\varkappa_{n7 \tilde{k}_{2}} \! = \! 
2n$, the number of times the pole $\alpha_{i(3)_{k_{3}}} \! = \! \sqrt{2}$ 
$(\neq \! \alpha_{7} \! = \! \infty)$ occurs is $\varkappa_{n7 \tilde{k}_{3}} 
\! = \! n$, and the number of times the pole $\alpha_{i(4)_{k_{4}}} \! = \! 
\pi$ $(\neq \! \alpha_{7} \! = \! \infty)$ occurs is $\varkappa_{n7 
\tilde{k}_{4}} \! = \! n$.\footnote{Note: for $n \! = \! 1$, as $\varkappa_{17 
\tilde{k}_{1}} \! = \! 1$, $\varkappa_{17 \tilde{k}_{2}} \! = \! 2$, 
$\varkappa_{17 \tilde{k}_{3}} \! = \! 1$, $\varkappa_{17 \tilde{k}_{4}} \! = 
\! 1$, and $\varkappa_{17} \! = \! \varrho_{7} \! = \! 2$, it follows that 
$\mathcal{P}_{1,7} \setminus \lbrace \alpha_{7},\alpha_{7} \rbrace \! = \! 
\mathcal{P}_{1,7} \setminus \lbrace \infty,\infty \rbrace \! = \! \lbrace 
0 \rbrace \cup \lbrace 1,1 \rbrace \cup \lbrace \! \sqrt{2} \rbrace \cup 
\lbrace \pi \rbrace$.} In this case, the ordered sum formula reads
\begin{equation*}
\sum_{q=1}^{4} \varkappa_{n7 \tilde{k}_{q}} \! := \! \varkappa_{n7 
\tilde{k}_{1}} \! + \! \varkappa_{n7 \tilde{k}_{2}} \! + \! 
\varkappa_{n7 \tilde{k}_{3}} \! + \! \varkappa_{n7 \tilde{k}_{4}} \! 
= \! n \! + \! 2n \! + \! n \! + \! n \! = \! 5n,
\end{equation*}
whence $\sum_{q=1}^{4} \varkappa_{n7 \tilde{k}_{q}} \! + \! \varkappa_{n7} 
\! = \! 5n \! + \! 2n \! = \! 7(n \! - \! 1) \! + \! 7$.
\end{enumerate}
This concludes the example.

For simplicity of notation, set, hereafter,
\begin{equation*}
\alpha_{i(q)_{k_{q}}} \! := \! \alpha_{p_{q}}, \quad q \! = \! 1,2,\dotsc,
\mathfrak{s} \! - \! 1.
\end{equation*}

The next subset in this ordering contains all positive integers, {}from $1$ 
to $k$, corresponding to the pole $\alpha_{k} \! = \! \infty$. For $k \! \in 
\! \lbrace 1,2,\dotsc,K \rbrace$ such that $\alpha_{k} \! = \! \infty$, write 
the ordered integer partition
\begin{equation*}
\lbrace \mathstrut k^{\prime} \! \in \! \lbrace 1,2,\dotsc,K \rbrace; \, 
k^{\prime} \! \leqslant \! k, \, \alpha_{k^{\prime}} \! = \! \alpha_{k} 
\! = \! \infty \rbrace \! := \! \lbrace \underbrace{i(\mathfrak{s})_{1},
i(\mathfrak{s})_{2},\dotsc,
i(\mathfrak{s})_{k_{\mathfrak{s}}}}_{k_{\mathfrak{s}}} \rbrace,
\end{equation*}
with $i(\mathfrak{s})_{k_{\mathfrak{s}}} \! = \! k$, $1 \! \leqslant \! 
i(\mathfrak{s})_{1} \! < \! i(\mathfrak{s})_{2} \! < \! \dotsb \! < \! 
i(\mathfrak{s})_{k_{\mathfrak{s}}} \! \leqslant \! K$, $\# \lbrace 
i(\mathfrak{s})_{1},i(\mathfrak{s})_{2},\dotsc,
i(\mathfrak{s})_{k_{\mathfrak{s}}} \rbrace \! = \! k_{\mathfrak{s}}$ $(= \! 
\varkappa_{1i(\mathfrak{s})_{k_{\mathfrak{s}}}} \! = \! \varkappa_{1k} \! 
= \! \varrho_{i(\mathfrak{s})_{k_{\mathfrak{s}}}})$ $= \! \varrho_{k}$, 
and $\lbrace i(j)_{1},i(j)_{2},\dotsc,i(j)_{k_{j}} \rbrace \cap 
\lbrace i(\mathfrak{s})_{1},i(\mathfrak{s})_{2},\dotsc,
i(\mathfrak{s})_{k_{\mathfrak{s}}} \rbrace \! = \! \varnothing$, 
$j \! = \! 1,2,\dotsc,\mathfrak{s} \! - \! 1$, which induces, by the 
above definition, the following pole ordering,
\begin{equation*}
\lbrace \mathstrut \alpha_{k^{\prime}}, \, k^{\prime} \! \in \! \lbrace 
1,2,\dotsc,K \rbrace; \, k^{\prime} \! \leqslant \! k, \alpha_{k^{\prime}} 
\! = \! \alpha_{k} \! = \! \infty \rbrace \! := \! \lbrace \underbrace{
\alpha_{i(\mathfrak{s})_{1}},\alpha_{i(\mathfrak{s})_{2}},\dotsc,
\alpha_{i(\mathfrak{s})_{k_{\mathfrak{s}}}}}_{k_{\mathfrak{s}}} \rbrace,
\end{equation*}
with $\alpha_{i(\mathfrak{s})_{k_{\mathfrak{s}}}} \! = \! \alpha_{k} 
\! = \! \infty$, $\alpha_{i(\mathfrak{s})_{1}} \! \prec \! 
\alpha_{i(\mathfrak{s})_{2}} \! \prec \! \dotsb \! \prec \! 
\alpha_{i(\mathfrak{s})_{k_{\mathfrak{s}}}}$, and $\lbrace \alpha_{i(j)_{1}},
\alpha_{i(j)_{2}},\dotsc,\alpha_{i(j)_{k_{j}}} \rbrace \cap \lbrace 
\alpha_{i(\mathfrak{s})_{1}},\alpha_{i(\mathfrak{s})_{2}},\dotsc,
\alpha_{i(\mathfrak{s})_{k_{\mathfrak{s}}}} \rbrace \! = \! \varnothing$, 
$j \! = \! 1,2,\dotsc,\mathfrak{s} \! - \! 1$, such that
\begin{equation*}
\alpha_{i(\mathfrak{s})_{1}} \! = \! \alpha_{i(\mathfrak{s})_{2}} 
\! = \! \dotsb \! = \! \alpha_{i(\mathfrak{s})_{k_{\mathfrak{s}}}} 
\! = \! \alpha_{k} \! =: \! \alpha_{p_{\mathfrak{s}}} \! = \! \infty, 
\quad \quad \# \lbrace \mathstrut \alpha_{i(\mathfrak{s})_{1}},
\alpha_{i(\mathfrak{s})_{2}},\dotsc,\alpha_{i(\mathfrak{s})_{k_{\mathfrak{s}}}} 
\rbrace \! = \! k_{\mathfrak{s}} \! = \! \varrho_{k}.
\end{equation*}
For $(\mathbb{N} \! \ni)$ $n \! \geqslant \! 2$ and $k \! \in \! \lbrace 
1,2,\dotsc,K \rbrace$ such that $\alpha_{k} \! = \! \infty$,
\begin{equation*}
\# \lbrace \mathstrut \alpha_{i(\mathfrak{s})_{k_{\mathfrak{s}}}},
\alpha_{i(\mathfrak{s})_{k_{\mathfrak{s}}}},\dotsc,
\alpha_{i(\mathfrak{s})_{k_{\mathfrak{s}}}} \rbrace \! = \! 
\varkappa_{ni(\mathfrak{s})_{k_{\mathfrak{s}}}} \! = \! \varkappa_{nk} 
\! = \! (n \! - \! 1) \gamma_{k} \! + \! \varrho_{k},
\end{equation*}
that is, as one moves {}from left to right across the repeated pole sequence 
$\mathcal{P}_{n,k}$ and removes the residual pole set (cf. the discussion 
above) $\lbrace \mathstrut \alpha_{k^{\prime}}, \, k^{\prime} \! \in \! 
\lbrace 1,2,\dotsc,K \rbrace; \, \alpha_{k^{\prime}} \! \neq \! \alpha_{k} 
\! = \! \infty \rbrace$, one is left with the set (via the above induced 
ordering) that consists of all occurrences of the pole $\alpha_{k} \! = \! 
\infty$, which occurs $\varkappa_{nk} \! = \! (n \! - \! 1) \gamma_{k} \! + 
\! \varrho_{k}$ times, namely,
\begin{equation*}
\mathcal{P}_{n,k} \setminus \bigcup_{q=1}^{\mathfrak{s}-1} \lbrace 
\underbrace{\alpha_{p_{q}},\alpha_{p_{q}},\dotsc,
\alpha_{p_{q}}}_{\varkappa_{nk \tilde{k}_{q}}} \rbrace \! := \! 
\lbrace \underbrace{\alpha_{i(\mathfrak{s})_{k_{\mathfrak{s}}}},
\alpha_{i(\mathfrak{s})_{k_{\mathfrak{s}}}},\dotsc,
\alpha_{i(\mathfrak{s})_{k_{\mathfrak{s}}}}}_{\varkappa_{ni(\mathfrak{s})_{
k_{\mathfrak{s}}}}} \rbrace \! = \! \lbrace \underbrace{\alpha_{k},
\alpha_{k},\dotsc,\alpha_{k}}_{\varkappa_{nk}} \rbrace \! = \! \lbrace 
\underbrace{\infty,\infty,\dotsc,\infty}_{\varkappa_{nk}} \rbrace.
\end{equation*}
In order to illustrate this latter notation, consider, again, the $K \! = \! 
7$ pole set $\lbrace \alpha_{1},\alpha_{2},\alpha_{3},\alpha_{4},\alpha_{5},
\alpha_{6},\alpha_{7} \rbrace \! = \! \lbrace 0,1,\infty,1,\linebreak[4]
\sqrt{2},\pi,\infty \rbrace$ for which $\mathfrak{s} \! = \! 5$ and 
$\alpha_{k} \! = \! \infty$, $k \! = \! 3,7$:
\begin{enumerate}
\item[(i)] \fbox{$k \! = \! 3$}
\begin{equation*}
\lbrace \mathstrut k^{\prime} \! \in \! \lbrace 1,2,\dotsc,7 \rbrace; \, 
k^{\prime} \! \leqslant \! 3, \, \alpha_{k^{\prime}} \! = \! \alpha_{3} \! 
= \! \infty \rbrace \! = \! \lbrace 3 \rbrace \! := \! \lbrace i(5)_{k_{5}} 
\rbrace \, \Rightarrow \, k_{5} \! = \! 1, \, \, i(5)_{1} \! = \! 3,
\end{equation*}
which induces the pole ordering
\begin{equation*}
\lbrace \mathstrut \alpha_{k^{\prime}}, \, k^{\prime} \! \in \! \lbrace 1,2,
\dotsc,7 \rbrace; \, k^{\prime} \! \leqslant \! 3, \, \alpha_{k^{\prime}} \! 
= \! \alpha_{3} \! = \! \infty \rbrace \! := \! \lbrace \alpha_{i(5)_{k_{5}}} 
\rbrace \! = \! \lbrace \alpha_{3} \rbrace \! = \! \lbrace \infty \rbrace,
\end{equation*}
hence
\begin{equation*}
\varkappa_{ni(\mathfrak{s})_{k_{\mathfrak{s}}}} \! = \! \varkappa_{ni(5)_{1}} 
\! = \! \varkappa_{n3} \! = \! (n \! - \! 1) \gamma_{3} \! + \! \varrho_{3} 
\! = \! 2(n \! - \! 1) \! + \! 1 \! = \! 2n \! - \! 1,
\end{equation*}
that is, as one moves {}from left to right across the repeated pole sequence
\begin{align*}
\mathcal{P}_{n,3} =& \, \lbrace \overset{1}{\underbrace{\alpha_{1},\alpha_{2},
\dotsc,\alpha_{7}}_{7}} \rbrace \cup \dotsb \cup \lbrace \overset{n-1}{
\underbrace{\alpha_{1},\alpha_{2},\dotsc,\alpha_{7}}_{7}} \rbrace \cup \lbrace 
\overset{n}{\underbrace{\alpha_{1},\alpha_{2},\alpha_{3}}_{3}} \rbrace \\
=& \, \lbrace \overset{1}{\underbrace{0,1,\infty,1,\sqrt{2},\pi,\infty}_{7}} 
\rbrace \cup \dotsb \cup \lbrace \overset{n-1}{\underbrace{0,1,\infty,1,
\sqrt{2},\pi,\infty}_{7}} \rbrace \cup \lbrace \overset{n}{\underbrace{0,1,
\infty}_{3}} \rbrace
\end{align*}
and removes the residual pole set (cf. the discussion and examples above) 
$\lbrace \mathstrut \alpha_{k^{\prime}}, \, k^{\prime} \! \in \! \lbrace 
1,2,\dotsc,7 \rbrace; \, \alpha_{k^{\prime}} \! \neq \! \alpha_{3} \! = \! 
\infty \rbrace$, one is left with the set (via the above induced ordering) 
that consists of all occurrences of the pole $\alpha_{3} \! = \! \infty$, 
which occurs $\varkappa_{n3} \! = \! 2n \! - \! 1$ times, namely,
\begin{equation*}
\mathcal{P}_{n,3} \setminus \bigcup_{q=1}^{4} \lbrace 
\underbrace{\alpha_{p_{q}},\alpha_{p_{q}},\dotsc,
\alpha_{p_{q}}}_{\varkappa_{n3 \tilde{k}_{q}}} \rbrace \! := \! \lbrace 
\underbrace{\alpha_{i(5)_{k_{5}}},\alpha_{i(5)_{k_{5}}},\dotsc,
\alpha_{i(5)_{k_{5}}}}_{\varkappa_{ni(5)_{k_{5}}}} \rbrace \! = \! \lbrace 
\underbrace{\alpha_{3},\alpha_{3},\dotsc,\alpha_{3}}_{\varkappa_{n3}} \rbrace 
\! = \! \lbrace \underbrace{\infty,\infty,\dotsc,\infty}_{2n-1} \rbrace;
\end{equation*}
\item[(ii)] \fbox{$k \! = \! 7$}
\begin{equation*}
\lbrace \mathstrut k^{\prime} \! \in \! \lbrace 1,2,\dotsc,7 \rbrace; \, 
k^{\prime} \! \leqslant \! 7, \, \alpha_{k^{\prime}} \! = \! \alpha_{7} \! 
= \! \infty \rbrace \! = \! \lbrace 3,7 \rbrace \! := \! \lbrace i(5)_{1},
i(5)_{k_{5}} \rbrace \, \Rightarrow \, k_{5} \! = \! 2, \, \, i(5)_{1} 
\! = \! 3, \, \, i(5)_{2} \! = \! 7,
\end{equation*}
which induces the pole ordering
\begin{equation*}
\lbrace \mathstrut \alpha_{k^{\prime}}, \, k^{\prime} \! \in \! \lbrace 1,2,
\dotsc,7 \rbrace; \, k^{\prime} \! \leqslant \! 7, \, \alpha_{k^{\prime}} \! 
= \! \alpha_{7} \! = \! \infty \rbrace \! := \! \lbrace \alpha_{i(5)_{1}},
\alpha_{i(5)_{k_{5}}} \rbrace \! = \! \lbrace \alpha_{3},\alpha_{7} \rbrace 
\! = \! \lbrace \infty,\infty \rbrace,
\end{equation*}
hence
\begin{equation*}
\varkappa_{ni(\mathfrak{s})_{k_{\mathfrak{s}}}} \! = \! \varkappa_{ni(5)_{2}} 
\! = \! \varkappa_{n7} \! = \! (n \! - \! 1) \gamma_{7} \! + \! \varrho_{7} 
\! = \! 2(n \! - \! 1) \! + \! 2 \! = \! 2n,
\end{equation*}
that is, as one moves {}from left to right across the repeated pole sequence
\begin{align*}
\mathcal{P}_{n,7} =& \, \lbrace \overset{1}{\underbrace{\alpha_{1},\alpha_{2},
\dotsc,\alpha_{7}}_{7}} \rbrace \cup \dotsb \cup \lbrace \overset{n-1}{
\underbrace{\alpha_{1},\alpha_{2},\dotsc,\alpha_{7}}_{7}} \rbrace \cup \lbrace 
\overset{n}{\underbrace{\alpha_{1},\alpha_{2},\dotsc,\alpha_{7}}_{7}} \rbrace 
\\
=& \, \lbrace \overset{1}{\underbrace{0,1,\infty,1,\sqrt{2},\pi,\infty}_{7}} 
\rbrace \cup \dotsb \cup \lbrace \overset{n-1}{\underbrace{0,1,\infty,1,
\sqrt{2},\pi,\infty}_{7}} \rbrace \cup \lbrace \overset{n}{\underbrace{0,1,
\infty,1,\sqrt{2},\pi,\infty}_{7}} \rbrace
\end{align*}
and removes the residual pole set (cf. the discussion and examples above) 
$\lbrace \mathstrut \alpha_{k^{\prime}}, \, k^{\prime} \! \in \! \lbrace 
1,2,\dotsc,7 \rbrace; \, \alpha_{k^{\prime}} \! \neq \! \alpha_{7} \! = \! 
\infty \rbrace$, one is left with the set (via the above induced ordering) 
that consists of all occurrences of the pole $\alpha_{7} \! = \! \infty$, 
which occurs $\varkappa_{n7} \! = \! 2n$ times, namely,
\begin{equation*}
\mathcal{P}_{n,7} \setminus \bigcup_{q=1}^{4} \lbrace 
\underbrace{\alpha_{p_{q}},\alpha_{p_{q}},\dotsc,
\alpha_{p_{q}}}_{\varkappa_{n7 \tilde{k}_{q}}} \rbrace \! := \! \lbrace 
\underbrace{\alpha_{i(5)_{k_{5}}},\alpha_{i(5)_{k_{5}}},\dotsc,
\alpha_{i(5)_{k_{5}}}}_{\varkappa_{ni(5)_{k_{5}}}} \rbrace \! = \! \lbrace 
\underbrace{\alpha_{7},\alpha_{7},\dotsc,\alpha_{7}}_{\varkappa_{n7}} 
\rbrace \! = \! \lbrace \underbrace{\infty,\infty,\dotsc,\infty}_{2n} 
\rbrace.
\end{equation*}
\end{enumerate}
This concludes the example.

With the above conventions and ordered disjoint partitions, one writes, for 
$n \! \in \! \mathbb{N}$ and $k \! \in \! \lbrace 1,2,\dotsc,K \rbrace$ 
such that $\alpha_{k} \! = \! \infty$, the repeated pole sequence 
$\mathcal{P}_{n,k}$ as the following ordered union of disjoint partitions:
\begin{align*}
\mathcal{P}_{n,k} =& \, \bigcup_{q=1}^{\mathfrak{s}-1} \lbrace 
\underbrace{\alpha_{i(q)_{k_{q}}},\alpha_{i(q)_{k_{q}}},\dotsc,
\alpha_{i(q)_{k_{q}}}}_{\varkappa_{nk \tilde{k}_{q}}} \rbrace \cup 
\lbrace \underbrace{\alpha_{i(\mathfrak{s})_{k_{\mathfrak{s}}}},
\alpha_{i(\mathfrak{s})_{k_{\mathfrak{s}}}},\dotsc,
\alpha_{i(\mathfrak{s})_{k_{\mathfrak{s}}}}}_{\varkappa_{ni(\mathfrak{s})_{
k_{\mathfrak{s}}}}} \rbrace \! := \! \bigcup_{q=1}^{\mathfrak{s}-1} 
\lbrace \underbrace{\alpha_{p_{q}},\alpha_{p_{q}},\dotsc,
\alpha_{p_{q}}}_{\varkappa_{nk \tilde{k}_{q}}} \rbrace \cup \lbrace 
\underbrace{\alpha_{k},\alpha_{k},\dotsc,\alpha_{k}}_{\varkappa_{nk}} 
\rbrace \\
=& \, \bigcup_{q=1}^{\mathfrak{s}-1} \lbrace \underbrace{\alpha_{p_{q}},
\alpha_{p_{q}},\dotsc,\alpha_{p_{q}}}_{\varkappa_{nk \tilde{k}_{q}}} \rbrace 
\cup \lbrace \underbrace{\infty,\infty,\dotsc,\infty}_{\varkappa_{nk}} 
\rbrace,
\end{align*}
where, by convention, the set $\lbrace \mathstrut \alpha_{k},\alpha_{k},
\dotsc,\alpha_{k} \rbrace$ $(= \! \lbrace \infty,\infty,\dotsc,\infty \rbrace)$ 
is written as the right-most set.

With the above notational preamble concluded, one now returns to the 
precise formulation of the orthogonality conditions for the MPC ORFs. One 
has a nested sequence of rational base sets. For $n \! \in \! \mathbb{N}$ 
and $k \! \in \! \lbrace 1,2,\dotsc,K \rbrace$ such that $\alpha_{k} \! = 
\! \infty$, one member of this nested sequence of rational base sets is
\begin{align*}
&\left\{\mathstrut \text{const.},\, \overset{1}{\underbrace{\mathscr{S}^{1}_{1}
(z),\mathscr{S}^{1}_{2}(z),\dotsc,\mathscr{S}^{1}_{K}(z)}_{K}},\, \overset{2}{
\underbrace{\mathscr{S}^{2}_{1}(z),\mathscr{S}^{2}_{2}(z),\dotsc,
\mathscr{S}^{2}_{K}(z)}_{K}},\dotsc,\, \overset{n}{\underbrace{\mathscr{S}^{n}_{1}
(z),\mathscr{S}^{n}_{2}(z),\dotsc,\mathscr{S}^{n}_{k}(z)}_{k}} \, \right\} \\
&:= \lbrace \text{const.} \rbrace \bigcup \cup_{q=1}^{\mathfrak{s}-1} 
\cup_{r=1}^{\varkappa_{nk \tilde{k}_{q}}} \left\{(z \! - \! 
\alpha_{p_{q}})^{-r} \right\} \bigcup \cup_{m=1}^{\varkappa_{nk}} 
\left\{z^{m} \right\} \\
&= \lbrace \text{const.} \rbrace \bigcup \cup_{q=1}^{\mathfrak{s}-1} \left\{
(z \! - \! \alpha_{p_{q}})^{-1},(z \! - \! \alpha_{p_{q}})^{-2},\dotsc,
(z \! - \! \alpha_{p_{q}})^{-\varkappa_{nk \tilde{k}_{q}}} \right\} \bigcup 
\left\{z,z^{2},\dotsc,z^{\varkappa_{nk}} \right\},
\end{align*}
corresponding, respectively, to the ordered repeated pole sequence
\begin{align*}
&\left\{\mathstrut \text{no pole},\, \overset{1}{\underbrace{\alpha_{1},
\alpha_{2},\dotsc,\alpha_{K}}_{K}},\, \overset{2}{\underbrace{\alpha_{1},
\alpha_{2},\dotsc,\alpha_{K}}_{K}},\, \dotsc,\, \overset{n}{\underbrace{
\alpha_{1},\alpha_{2},\dotsc,\alpha_{k}}_{k}} \, \right\} \\
&:= \lbrace \text{no pole} \rbrace \bigcup \cup_{q=1}^{\mathfrak{s}-1} 
\lbrace \underbrace{\alpha_{p_{q}},\alpha_{p_{q}},\dotsc,
\alpha_{p_{q}}}_{\varkappa_{nk \tilde{k}_{q}}} \rbrace \bigcup \lbrace 
\underbrace{\alpha_{k},\alpha_{k},\dotsc,\alpha_{k}}_{\varkappa_{nk}} 
\rbrace \\
&= \lbrace \text{no pole} \rbrace \bigcup \cup_{q=1}^{\mathfrak{s}-1} 
\lbrace \underbrace{\alpha_{p_{q}},\alpha_{p_{q}},\dotsc,
\alpha_{p_{q}}}_{\varkappa_{nk \tilde{k}_{q}}} \rbrace \bigcup \lbrace 
\underbrace{\infty,\infty,\dotsc,\infty}_{\varkappa_{nk}} \rbrace.
\end{align*}
Orthonormalisation with respect to $\langle \pmb{\cdot},\pmb{\cdot} 
\rangle_{\mathscr{L}}$, via the Gram-Schmidt orthogonalisation 
method, leads to the MPC ORFs, $\lbrace \phi^{n}_{k}(z) 
\rbrace_{\underset{k=1,2,\dotsc,K}{n \in \mathbb{N}}}$ (for consistency 
of notation, set $\phi^{0}_{0}(z) \! \equiv \! 1)$, which, by suitable 
normalisation (see below), can be written as, with the above orderings,
\begin{align*}
\phi^{n}_{k} \colon& \, \mathbb{N} \times \lbrace 1,2,\dotsc,K \rbrace 
\times \overline{\mathbb{C}} \setminus \lbrace \alpha_{1},\alpha_{2},
\dotsc,\alpha_{K} \rbrace \! \to \! \mathbb{C}, \, \, (n,k,z) \! \mapsto \\
\phi^{n}_{k}(z) :=& \, \phi_{0}^{\infty}(n,k) \! + \! \sum_{q=1}^{
\mathfrak{s}-1} \sum_{r=1}^{\varkappa_{nk \tilde{k}_{q}}} 
\dfrac{\nu^{\infty}_{r,q}(n,k)}{(z \! - \! \alpha_{p_{q}})^{r}} \! + \! 
\sum_{m=1}^{\varkappa_{nk}} \mu^{\infty}_{n,m}(n,k)z^{m},
\end{align*}
where the $\phi^{n}_{k}$'s are normalised so that they all have real 
coefficients; in particular, for $n \! \in \! \mathbb{N}$ and $k \! \in \! 
\lbrace 1,2,\dotsc,K \rbrace$ such that $\alpha_{k} \! = \! \infty$,
\begin{equation*}
\operatorname{LC}(\phi^{n}_{k}) \! = \! \mu^{\infty}_{n,\varkappa_{nk}}
(n,k) \! > \! 0.
\end{equation*}
(For consistency of notation, set $\phi^{\infty}_{0}(0,0) \! \equiv \! 1$.) 
Furthermore, for $n \! \in \! \mathbb{N}$ and $k \! \in \! \lbrace 1,2,
\dotsc,K \rbrace$ such that $\alpha_{k} \! = \! \infty$, note that, by 
construction:
\begin{gather}
\left\langle \phi^{n}_{k},z^{j} \right\rangle_{\mathscr{L}} \! = \! 
\int_{\mathbb{R}} \phi^{n}_{k}(\xi) \xi^{j} \, \md \mu (\xi) \! = \! 0, 
\quad j \! = \! 0,1,\dotsc,\varkappa_{nk} \! - \! 1, \label{eq6} \\
\left\langle \phi^{n}_{k},\mu^{\infty}_{n,\varkappa_{nk}}(n,k)
z^{\varkappa_{nk}} \right\rangle_{\mathscr{L}} \! = \! 
\mu^{\infty}_{n,\varkappa_{nk}}(n,k) \int_{\mathbb{R}} \phi^{n}_{k}
(\xi) \xi^{\varkappa_{nk}} \, \md \mu (\xi) \! = \! 1, \label{eq7} \\
\left\langle \phi^{n}_{k},(z \! - \! \alpha_{p_{q}})^{-r} 
\right\rangle_{\mathscr{L}} \! = \! \int_{\mathbb{R}} \phi^{n}_{k}(\xi)
(\xi \! - \! \alpha_{p_{q}})^{-r} \, \md \mu (\xi) \! = \! 0, \quad q 
\! = \! 1,2,\dotsc,\mathfrak{s} \! - \! 1, \quad r \! = \! 1,2,\dotsc,
\varkappa_{nk \tilde{k}_{q}}. \label{eq8}
\end{gather}
(Note: if, for $k \! \in \! \lbrace 1,2,\dotsc,K \rbrace$ such that 
$\alpha_{k} \! = \! \infty$, the residual pole set $\lbrace \mathstrut 
\alpha_{k^{\prime}}, \, k^{\prime} \! \in \! \lbrace 1,2,\dotsc,K \rbrace; 
\, \alpha_{k^{\prime}} \! \neq \! \alpha_{k} \! = \! \infty \rbrace \! 
= \! \varnothing$, then the orthogonality condition~\eqref{eq8} is 
vacuous; actually, this can only occur for $n \! = \! 1$.) For $n \! 
\in \! \mathbb{N}$ and $k \! \in \! \lbrace 1,2,\dotsc,K \rbrace$ 
such that $\alpha_{k} \! = \! \infty$, the orthogonality 
conditions~\eqref{eq6}--\eqref{eq8} give rise to a total 
of (cf. Equation~\eqref{infcount})
\begin{equation*}
\sum_{q=1}^{\mathfrak{s}-1} \varkappa_{nk \tilde{k}_{q}} \! + \! 
\varkappa_{nk} \! + \! 1 \! = \! (n \! - \! 1)K \! + \! k \! + \! 1
\end{equation*}
linear equations determining the $(n \! - \! 1)K \! + \! k \! + \! 1$ real 
$(n$- and $k$-dependent) coefficients.

It is convenient to introduce, at this stage, the main object of study of this 
monograph, that is, the corresponding \emph{monic} MPC ORFs, $\lbrace 
\pmb{\pi}^{n}_{k}(z) \rbrace_{\underset{k=1,2,\dotsc,K}{n \in \mathbb{N}}}$ 
(for consistency of notation, set $\pmb{\pi}^{0}_{0}(z) \! \equiv \! 1)$. For 
$n \! \in \! \mathbb{N}$ and $k \! \in \! \lbrace 1,2,\dotsc,K \rbrace$ such 
that $\alpha_{k} \! = \! \infty$,
\begin{align*}
\pmb{\pi}^{n}_{k} \colon& \, \mathbb{N} \times \lbrace 1,2,\dotsc,K \rbrace 
\times \overline{\mathbb{C}} \setminus \lbrace \alpha_{1},\alpha_{2},
\dotsc,\alpha_{K} \rbrace \! \to \! \mathbb{C}, \, \, (n,k,z) \! \mapsto \\
\pmb{\pi}^{n}_{k}(z) :=& \, \dfrac{\phi^{n}_{k}(z)}{\operatorname{LC}
(\phi^{n}_{k})} \! = \! \dfrac{\phi^{\infty}_{0}(n,k)}{\mu^{\infty}_{n,
\varkappa_{nk}}(n,k)} \! + \! \dfrac{1}{\mu^{\infty}_{n,\varkappa_{nk}}(n,k)} 
\sum_{q=1}^{\mathfrak{s}-1} \sum_{r=1}^{\varkappa_{nk \tilde{k}_{q}}} 
\dfrac{\nu^{\infty}_{r,q}(n,k)}{(z \! - \! \alpha_{p_{q}})^{r}} \! + \! 
\dfrac{1}{\mu^{\infty}_{n,\varkappa_{nk}}(n,k)} \sum_{m=1}^{\varkappa_{nk}-1} 
\mu^{\infty}_{n,m}(n,k)z^{m} \! + \! z^{\varkappa_{nk}}.
\end{align*}
The monic MPC ORFs, $\lbrace \pmb{\pi}^{n}_{k}(z) \rbrace_{
\underset{k=1,2,\dotsc,K}{n \in \mathbb{N}}}$, possess the following 
orthogonality properties:
\begin{gather}
\left\langle \pmb{\pi}^{n}_{k},z^{j} \right\rangle_{\mathscr{L}} \! = \! 
\int_{\mathbb{R}} \pmb{\pi}^{n}_{k}(\xi) \xi^{j} \, \md \mu (\xi) \! = \! 
0, \quad j \! = \! 0,1,\dotsc,\varkappa_{nk} \! - \! 1, \label{eq9} \\
\left\langle \pmb{\pi}^{n}_{k},z^{\varkappa_{nk}} \right\rangle_{\mathscr{L}} 
\! = \! \int_{\mathbb{R}} \pmb{\pi}^{n}_{k}(\xi) \xi^{\varkappa_{nk}} \, 
\md \mu (\xi) \! = \! (\operatorname{LC}(\phi^{n}_{k}))^{-2} \! = \! 
(\mu^{\infty}_{n,\varkappa_{nk}}(n,k))^{-2}, \label{eq10} \\
\left\langle \pmb{\pi}^{n}_{k},(z \! - \! \alpha_{p_{q}})^{-r} 
\right\rangle_{\mathscr{L}} \! = \! \int_{\mathbb{R}} \pmb{\pi}^{n}_{k}(\xi)
(\xi \! - \! \alpha_{p_{q}})^{-r} \, \md \mu (\xi) \! = \! 0, \quad q 
\! = \! 1,2,\dotsc,\mathfrak{s} \! - \! 1, \quad r \! = \! 1,2,\dotsc,
\varkappa_{nk \tilde{k}_{q}}. \label{eq11}
\end{gather}
(Note: if, for $k \! \in \! \lbrace 1,2,\dotsc,K \rbrace$ such that 
$\alpha_{k} \! = \! \infty$, the residual pole set $\lbrace \mathstrut 
\alpha_{k^{\prime}}, \, k^{\prime} \! \in \! \lbrace 1,2,\dotsc,K \rbrace; 
\, \alpha_{k^{\prime}} \! \neq \! \alpha_{k} \! = \! \infty \rbrace \! = 
\! \varnothing$, then the orthogonality condition~\eqref{eq11} is 
vacuous; actually, this can occur only for $n \! = \! 1$.) For $n \! \in \! 
\mathbb{N}$ and $k \! \in \! \lbrace 1,2,\dotsc,K \rbrace$ such that 
$\alpha_{k} \! = \! \infty$, it follows {}from the monic MPC ORF 
orthogonality conditions~\eqref{eq9}--\eqref{eq11} that
\begin{equation*}
\left\langle \pmb{\pi}^{n}_{k},\pmb{\pi}^{n}_{k} \right\rangle_{
\mathscr{L}} \! =: \! \lvert \lvert \pmb{\pi}^{n}_{k}(\pmb{\cdot}) \rvert 
\rvert^{2}_{\mathscr{L}} \! = \! (\operatorname{LC}(\phi^{n}_{k}))^{-2} 
\! = \! (\mu_{n,\varkappa_{nk}}^{\infty}(n,k))^{-2},
\end{equation*}
whence $\lvert \lvert \pmb{\pi}^{n}_{k}(\pmb{\cdot}) \rvert \rvert_{
\mathscr{L}} \! = \! (\mu^{\infty}_{n,\varkappa_{nk}}(n,k))^{-1} \! > 
\! 0$.\footnote{See, in particular, Section~\ref{sec2}, Lemma~\ref{lem2.1}, 
Equations~\eqref{nhormmconstatinf}, \eqref{eq37}, and~\eqref{eq39}, and 
Corollary~\ref{cor2.1}, Equation~\eqref{eqnmctinf1}; incidentally, this also 
establishes the positive definiteness of the linear functional $\mathscr{L} 
\colon \Lambda^{\mathbb{R}} \! \to \! \mathbb{R}$, and hence 
the---real---bilinear form $\langle \pmb{\cdot},\pmb{\cdot} 
\rangle_{\mathscr{L}} \colon \Lambda^{\mathbb{R}} \times 
\Lambda^{\mathbb{R}} \! \to \! \mathbb{R}$, defined at the beginning 
of Subsection~\ref{subsec1.2}.}
\begin{eeee} \label{rem1.2.3} 
\textsl{For $n \! \in \! \mathbb{N}$ and $k \! \in \! \lbrace 1,2,\dotsc,
K \rbrace$ such that $\alpha_{p_{\mathfrak{s}}} \! := \! \alpha_{k} \! 
= \! \infty$, {}from the above partial fraction decomposition for the monic 
{\rm MPC ORFs}, $\lbrace \pmb{\pi}^{n}_{k}(z) \rbrace_{\underset{k=1,2,
\dotsc,K}{n \in \mathbb{N}}}$, one writes
\begin{equation*}
\pmb{\pi}^{n}_{k} \colon \mathbb{N} \times \lbrace 1,2,\dotsc,K \rbrace 
\times \overline{\mathbb{C}} \setminus \lbrace \alpha_{1},\alpha_{2},
\dotsc,\alpha_{K} \rbrace \! \ni \! (n,k,z) \! \mapsto \! \dfrac{1}{
\mu^{\infty}_{n,\varkappa_{nk}}(n,k)} \dfrac{\hat{\mathfrak{P}}^{n}_{k}
(z)}{\prod_{q=1}^{\mathfrak{s}-1}(z \! - \! \alpha_{p_{q}})^{\varkappa_{nk 
\tilde{k}_{q}}}},
\end{equation*}
where
\begin{align*}
\hat{\mathfrak{P}}^{n}_{k}(z) :=& \, \phi^{\infty}_{0}(n,k) 
\prod_{q=1}^{\mathfrak{s}-1}(z \! - \! \alpha_{p_{q}})^{\varkappa_{nk 
\tilde{k}_{q}}} \! + \! \prod_{q=1}^{\mathfrak{s}-1}(z \! - \! 
\alpha_{p_{q}})^{\varkappa_{nk \tilde{k}_{q}}} \sum_{r=1}^{\varkappa_{nk}} 
\mu^{\infty}_{n,r}(n,k)z^{r} \\
+& \, \sum_{j=1}^{\mathfrak{s}-1} 
\prod_{\substack{q=1\\q \neq j}}^{\mathfrak{s}-1}
(z \! - \! \alpha_{p_{q}})^{\varkappa_{nk \tilde{k}_{q}}} 
\sum_{r=1}^{\varkappa_{nk \tilde{k}_{j}}} \nu^{\infty}_{r,j}(n,k)
(z \! - \! \alpha_{p_{j}})^{\varkappa_{nk \tilde{k}_{j}}-r}.
\end{align*}
Since, for $n \! \in \! \mathbb{N}$ and $k \! \in \! \lbrace 1,2,\dotsc,K 
\rbrace$ such that $\alpha_{p_{\mathfrak{s}}} \! := \! \alpha_{k} \! = \! 
\infty$, $\mu^{\infty}_{n,\varkappa_{nk}}(n,k) \! > \! 0$, it follows that 
$\hat{\mathfrak{P}}^{n}_{k}(z)$ is a polynomial of degree at most 
$\sum_{q=1}^{\mathfrak{s}-1} \varkappa_{nk \tilde{k}_{q}} \! + \! 
\varkappa_{nk} \! = \! (n \! - \! 1)K \! + \! k$. The following discussion 
constitutes a brief synopsis of the classification theory for the zeros of the 
monic {\rm MPC ORF}, $\pmb{\pi}^{n}_{k}(z)$,\footnote{Similar statements 
apply, \emph{verbatim}, for the corresponding MPC ORF, $\phi^{n}_{k}(z) \! 
= \! \mu^{\infty}_{n,\varkappa_{nk}}(n,k) \pmb{\pi}^{n}_{k}(z)$.} for the case 
$n \! \in \! \mathbb{N}$ and $k \! \in \! \lbrace 1,2,\dotsc,K \rbrace$ such 
that $\alpha_{p_{\mathfrak{s}}} \! := \! \alpha_{k} \! = \! \infty$, that is based 
on the seminal works of Nj\r{a}stad {\rm \cite{a1,a2,ol1}}. Via the above 
definition of $\hat{\mathfrak{P}}^{n}_{k}(z)$, a straightforward calculation 
shows that $\hat{\mathfrak{P}}^{n}_{k}(\alpha_{p_{q}}) \! = \! \nu^{\infty}_{
\varkappa_{nk \tilde{k}_{q}},q}(n,k) \prod_{\underset{q^{\prime} \neq q}{
q^{\prime}=1}}^{\mathfrak{s}-1}(\alpha_{p_{q}} \! - \! \alpha_{p_{q^{\prime}}})^{
\varkappa_{nk \tilde{k}_{q^{\prime}}}}$, $q \! \in \! \lbrace 1,2,\dotsc,\mathfrak{s} 
\! - \! 1 \rbrace$, and, since $\mu^{\infty}_{n,\varkappa_{nk}}(n,k) \! > \! 0$, 
$\hat{\mathfrak{P}}^{n}_{k}(\alpha_{k})$ $(=: \! \hat{\mathfrak{P}}^{n}_{k}
(\alpha_{p_{\mathfrak{s}}}) \! = \! \hat{\mathfrak{P}}^{n}_{k}(\infty))$ 
$\neq \! 0$.\footnote{Note: $\hat{\mathfrak{P}}^{n}_{k}(z)$ has a pole of 
order $(n \! - \! 1)K \! + \! k$ at (the point at infinity) $\alpha_{p_{\mathfrak{s}}} 
\! := \! \alpha_{k} \! = \! \infty$ with positive coefficient $\mu^{\infty}_{n,
\varkappa_{nk}}(n,k)$, that is, $\hat{\mathfrak{P}}^{n}_{k}(z) \! 
=_{\overline{\mathbb{C}} \ni z \to \alpha_{p_{\mathfrak{s}}} := \alpha_{k} = \infty} 
\! \mu^{\infty}_{n,\varkappa_{nk}}(n,k)z^{(n-1)K+k}(1 \! + \! \mathcal{O}(z^{-1}))$.} 
One calls the monic {\rm MPC ORF}, $\pmb{\pi}^{n}_{k}(z)$, and the corresponding 
index two-tuple $(n,k)$ \emph{degenerate} if $\# \lbrace \mathstrut z \! \in \! 
\overline{\mathbb{C}}; \, \hat{\mathfrak{P}}^{n}_{k}(z) \! = \! 0 \rbrace \! \leqslant 
\! (n \! - \! 1)K \! + \! k \! - \! 1$, that is, $\deg (\hat{\mathfrak{P}}^{n}_{k}) \! 
< \! (n \! - \! 1)K \! + \! k$. One calls the monic {\rm MPC ORF}, $\pmb{\pi}^{n}_{k}
(z)$, and the corresponding index two-tuple $(n,k)$ $\alpha_{p_{q}}$-\emph{defective}, 
$q \! \in \! \lbrace 1,2,\dotsc,\mathfrak{s} \rbrace$, if $\alpha_{p_{q}}$ is a zero 
of $\hat{\mathfrak{P}}^{n}_{k}(z)$, namely: {\rm (i)} if $\nu^{\infty}_{\varkappa_{nk 
\tilde{k}_{q}},q}(n,k) \! = \! 0$, $q \! \in \! \lbrace 1,2,\dotsc,\mathfrak{s} \! - \! 
1 \rbrace$, then $\hat{\mathfrak{P}}^{n}_{k}(\alpha_{p_{q}}) \! = \! 0$, $q \! \in 
\! \lbrace 1,2,\dotsc,\mathfrak{s} \! - \! 1 \rbrace$, thus $\pmb{\pi}^{n}_{k}(z)$ 
and the corresponding index two-tuple $(n,k)$ are $\alpha_{p_{q}}$-defective$;$ 
and {\rm (ii)} since $\mu^{\infty}_{n,\varkappa_{nk}}(n,k) \! > \! 0$, that is, 
$\hat{\mathfrak{P}}^{n}_{k}(\alpha_{k})$ $(=: \! \hat{\mathfrak{P}}^{n}_{k}
(\alpha_{p_{\mathfrak{s}}}) \! = \! \hat{\mathfrak{P}}^{n}_{k}(\infty))$ $\neq 
\! 0$, it follows that $\pmb{\pi}^{n}_{k}(z)$ and the corresponding index 
two-tuple $(n,k)$ can not be $\alpha_{k}$ $(=: \! \alpha_{p_{\mathfrak{s}}} \! = 
\! \infty)$-defective. One calls the monic {\rm MPC ORF}, $\pmb{\pi}^{n}_{k}(z)$, 
(and the corresponding index two-tuple $(n,k))$ \emph{defective} if it is 
$\alpha_{p_{q}}$-defective for at least one $q \! \in \! \lbrace 1,2,\dotsc,
\mathfrak{s} \! - \! 1 \rbrace$ or \emph{maximally defective} if it is 
$\alpha_{p_{q}}$-defective $\forall$ $q \! \in \! \lbrace 1,2,\dotsc,\mathfrak{s} 
\! - \! 1 \rbrace$. One calls the monic {\rm MPC ORF}, $\pmb{\pi}^{n}_{k}(z)$, 
(and the corresponding index two-tuple $(n,k))$ \emph{singular} if it is 
degenerate or defective; otherwise, it is non-singular. In this monograph, it is 
assumed that $\pmb{\pi}^{n}_{k}(z)$ (and the corresponding index two-tuple 
$(n,k))$ is neither degenerate nor defective, that is, $\deg (\hat{\mathfrak{P}}^{n}_{k}) 
\! = \! (n \! - \! 1)K \! + \! k$ and $\nu^{\infty}_{\varkappa_{nk \tilde{k}_{q}},q}
(n,k) \! \neq \! 0$ $\forall$ $q \! \in \! \lbrace 1,2,\dotsc,\mathfrak{s} \! - \! 1 
\rbrace$ (recall that, since $\mu^{\infty}_{n,\varkappa_{nk}}(n,k) \! > \! 0$, 
$\hat{\mathfrak{P}}^{n}_{k}(\alpha_{k})$ $(=: \! \hat{\mathfrak{P}}^{n}_{k}
(\alpha_{p_{\mathfrak{s}}}) \! = \! \hat{\mathfrak{P}}^{n}_{k}(\infty))$ 
$\neq \! 0)$$;$ therefore, hereafter, it is to be understood that `monic 
{\rm MPC ORF}' and `non-singular monic {\rm MPC ORF}' are synonymous 
objects (the singular case(s) will be considered elsewhere$)$$;$ for simplicity, 
though, only the phrase {\rm MPC ORF} (monic or not) will be used. For 
$n \! \in \! \mathbb{N}$ and $k \! \in \! \lbrace 1,2,\dotsc,K \rbrace$ such that 
$\alpha_{p_{\mathfrak{s}}} \! := \! \alpha_{k} \! = \! \infty$, writing the factorisation 
$\hat{\mathfrak{P}}^{n}_{k}(z) \! := \! \mu^{\infty}_{n,\varkappa_{nk}}(n,k) 
\prod_{j=1}^{(n-1)K+k}(z \! - \! \hat{\mathfrak{z}}^{n}_{k}(j))$,\footnote{Of course, 
$\hat{\mathfrak{z}}^{n}_{k}(j)$ also depends on $\alpha_{p_{q}}$, $q \! = \! 1,2,
\dotsc,\mathfrak{s} \! - \! 1$; but, for simplicity of notation, this extraneous 
dependence is suppressed.} where, counting multiplicities, $\lbrace \mathstrut 
\hat{\mathfrak{z}}^{n}_{k}(j) \rbrace_{j=1}^{(n-1)K+k} \! := \! \lbrace \mathstrut 
z \! \in \! \overline{\mathbb{C}}; \, \hat{\mathfrak{P}}^{n}_{k}(z) \! = \! 0 \rbrace$ 
$(= \! \lbrace \mathstrut z \! \in \! \overline{\mathbb{C}}; \, \pmb{\pi}^{n}_{k}(z) 
\! = \! 0 \rbrace)$, it is shown in Section~\ref{sec3} (see, in particular, the 
corresponding items of Lemmata~\ref{lemrootz} and~\ref{lemetatomu} for precise 
statements) that, in the double-scaling limit $\mathscr{N},n \! \to \! \infty$ 
such that $\mathscr{N}/n \! \to \! 1$, the zeros (counting multiplicities) 
$\hat{\mathfrak{z}}^{n}_{k}(j) \! \in \! \overline{\mathbb{R}} \setminus \lbrace 
\alpha_{1},\alpha_{2},\dotsc,\alpha_{K} \rbrace$, $j \! = \! 1,2,\dotsc,
(n \! - \! 1)K \! + \! k$$;$ more precisely, in the double-scaling limit 
$\mathscr{N},n \! \to \! \infty$ such that $\mathscr{N}/n \! \to \! 1$, 
$\lbrace \hat{\mathfrak{z}}^{n}_{k}(j) \rbrace_{j=1}^{(n-1)K+k}$ accumulates 
on the real compact set $J_{\infty}$, where $J_{\infty}$ is defined in 
Subsection~\ref{subsub2}. (See, also,  {\rm \cite{a1,a2,ol1}}$;$ and, 
in a related context, see {\rm \cite{jvdab,davd,pz14,pz11,pz12,pz7,
pz16,pz13,pz15,pz18,pz17}.)}}
\end{eeee}
The brief discussion that follows is motivated, in part, by the seminal 
works of Nj\r{a}stad \cite{n2,n1} related to MPAs and ORFs: more 
precise statements can be located in Section~\ref{sek5}; see, in 
particular, Lemma~\ref{lemmpainffin}. Furthermore, in order to mitigate 
notational encumbrances, explicit $n$- and $k$-dependencies will 
be temporarily suppressed, except where absolutely necessary. The 
Markov-Stieltjes transform of the probability measure $\mu$ (cf. 
Equations~\eqref{eq1}--\eqref{eq5}) is defined as
\begin{equation} \label{mvssinf1} 
\mathrm{F}_{\mu}(z) \! := \! \int_{\mathbb{R}}(z \! - \! \xi)^{-1} 
\, \md \mu (\xi);
\end{equation}
for $n \! \in \! \mathbb{N}$ and $k \! \in \! \lbrace 1,2,\dotsc,K 
\rbrace$ such that $\alpha_{p_{\mathfrak{s}}} \! := \! \alpha_{k} 
\! = \! \infty$, a straightforward calculation shows that 
$\mathrm{F}_{\mu}(z)$ has the formal asymptotic expansions
\begin{equation} \label{mvssinf2} 
\mathrm{F}_{\mu}(z) \underset{z \to \alpha_{p_{q}}}{=} \sum_{j=
0}^{\infty}c^{(q)}_{j}(\alpha_{p_{q}})(z \! - \! \alpha_{p_{q}})^{j}, 
\quad q \! = \! 1,2,\dotsc,\mathfrak{s} \! - \! 1,
\end{equation}
where $c^{(q)}_{j}(\alpha_{p_{q}}) \! := \! -\int_{\mathbb{R}}
(\xi \! - \! \alpha_{p_{q}})^{-(1+j)} \, \md \mu (\xi)$, $(j,q) \! 
\in \! \mathbb{N}_{0} \times \lbrace 1,2,\dotsc,\mathfrak{s} 
\! - \! 1 \rbrace$, and
\begin{equation} \label{mvssinf3} 
\mathrm{F}_{\mu}(z) \underset{z \to \alpha_{k}}{=} 
\sum_{j=1}^{\infty}c^{(\infty)}_{j}z^{-j},
\end{equation}
where $c^{(\infty)}_{j} \! := \! \int_{\mathbb{R}} \xi^{j-1} \, \md 
\mu (\xi)$, $j \! \in \! \mathbb{N}$, with $c^{(\infty)}_{1} \! = \! 1$. 
For $n \! \in \! \mathbb{N}$ and $k \! \in \! \lbrace 1,2,\dotsc,K 
\rbrace$ such that $\alpha_{p_{\mathfrak{s}}} \! := \! \alpha_{k} 
\! = \! \infty$, define the associated $\mathrm{R}$-function 
\cite{n2,n1} as follows:
\begin{equation} \label{mvssinf4} 
\widehat{\pmb{\mathrm{R}}}_{\mu} \colon \mathbb{N} \times 
\lbrace 1,2,\dotsc,K \rbrace \times \overline{\mathbb{C}} 
\setminus \lbrace \alpha_{1},\alpha_{2},\dotsc,\alpha_{K} 
\rbrace \! \ni \! (n,k,z) \! \mapsto \! \int_{\mathbb{R}} 
\left(\dfrac{\pmb{\pi}^{n}_{k}(\xi) \! - \! \pmb{\pi}^{n}_{k}
(z)}{\xi \! - \! z} \right) \md \mu (\xi) \! =: \! 
\widehat{\pmb{\mathrm{R}}}_{\mu}(z):
\end{equation}
{}from Equation~\eqref{mvssinf4} and the fact that $\mu \! 
\in \! \mathscr{M}_{1}(\mathbb{R})$, one shows, via the identity 
$y_{1}^{m} \! - \! y_{2}^{m} \! = \! (y_{1} \! - \! y_{2})(y_{1}^{m-1} 
\! + \! y_{1}^{m-2}y_{2} \! + \! \dotsb \! + \! y_{1}y_{2}^{m-2} \! + 
\! y_{2}^{m-1})$, that $\widehat{\pmb{\mathrm{R}}}_{\mu}(z)$ can 
be presented as the improper fraction
\begin{equation} \label{mvssinf5} 
\widehat{\pmb{\mathrm{R}}}_{\mu}(z) \! = \! \dfrac{\widehat{
\mathrm{U}}_{\mu}(z)}{\prod_{q=1}^{\mathfrak{s}-1}(z \! - \! 
\alpha_{p_{q}})^{\varkappa_{nk \tilde{k}_{q}}}},
\end{equation}
where $\widehat{\mathrm{U}}_{\mu}(z) \! := \! \sum_{j=0}^{(n-1)K
+k-1} \hat{r}_{j}z^{j}$, with $\operatorname{deg}(\widehat{\mathrm{U}}_{
\mu}(z)) \! = \! (n \! - \! 1)K \! + \! k \! - \! 1$ (since $\hat{r}_{(n-1)K
+k-1} \! \neq \! 0)$, and (cf. Equation~\eqref{infcount}) 
$\operatorname{deg}(\prod_{q=1}^{\mathfrak{s}-1}(z \! - \! 
\alpha_{p_{q}})^{\varkappa_{nk \tilde{k}_{q}}}) \! = \! (n \! - \! 1)K 
\! + \! k \! - \! \varkappa_{nk}$; moreover, note that, for $n \! \in 
\! \mathbb{N}$ and $k \! \in \! \lbrace 1,2,\dotsc,K \rbrace$ such 
that $\alpha_{p_{\mathfrak{s}}} \! := \! \alpha_{k} \! = \! \infty$, 
the corresponding monic MPC ORF, $\pmb{\pi}^{n}_{k}(z)$, can 
also be presented as the improper fraction
\begin{equation} \label{mvssinf6} 
\pmb{\pi}^{n}_{k}(z) \! = \! \dfrac{\widehat{\mathrm{V}}_{\mu}(z)}{
\prod_{q=1}^{\mathfrak{s}-1}(z \! - \! \alpha_{p_{q}})^{\varkappa_{nk 
\tilde{k}_{q}}}},
\end{equation}
where $\widehat{\mathrm{V}}_{\mu}(z) \! := \! \sum_{j=0}^{(n-1)K+k} 
\hat{t}_{j}z^{j}$, with $\operatorname{deg}(\widehat{\mathrm{V}}_{\mu}
(z)) \! = \! (n \! - \! 1)K \! + \! k$ (since $\hat{t}_{(n-1)K+k} \! = \! 1)$, 
that is, $\widehat{\mathrm{V}}_{\mu}(z)$ is a monic polynomial of 
degree $(n \! - \! 1)K \! + \! k$. It turns out that, for $n \! \in \! 
\mathbb{N}$ and $k \! \in \! \lbrace 1,2,\dotsc,K \rbrace$ such 
that $\alpha_{p_{\mathfrak{s}}} \! := \! \alpha_{k} \! = \! \infty$, 
$\widehat{\mathrm{U}}_{\mu}(z)/\widehat{\mathrm{V}}_{\mu}(z)$ 
$(= \! \widehat{\pmb{\mathrm{R}}}_{\mu}(z)/\pmb{\pi}^{n}_{k}(z))$ 
is the MPA of type $((n \! - \! 1)K \! + \! k \! - \! 1,(n \! - \! 1)
K \! + \! k)$ for the Markov-Stieltjes transform, that is, it is 
the---unique---(proper) rational function with $\operatorname{deg}
(\widehat{\mathrm{U}}_{\mu}(z)) \! = \! (n \! - \! 1)K \! + \! k \! - \! 1$, 
$\operatorname{deg}(\widehat{\mathrm{V}}_{\mu}(z)) \! = \! (n \! - \! 1)K 
\! + \! k$, and $(\widehat{\mathrm{U}}_{\mu}(z),\widehat{\mathrm{V}}_{
\mu}(z))$ coprime interpolating $\mathrm{F}_{\mu}(z)$ and satisfying the 
interpolation conditions
\begin{gather}
\dfrac{\widehat{\mathrm{U}}_{\mu}(z)}{\widehat{\mathrm{V}}_{\mu}
(z)} \! - \! \sum_{j=0}^{2 \varkappa_{nk \tilde{k}_{q}}-1}c^{(q)}_{j}
(\alpha_{p_{q}})(z \! - \! \alpha_{p_{q}})^{j} \underset{z \to \alpha_{p_{q}}}{=} 
\mathcal{O} \left((z \! - \! \alpha_{p_{q}})^{2 \varkappa_{nk \tilde{k}_{q}}} 
\right), \quad q \! = \! 1,2,\dotsc,\mathfrak{s} \! - \! 1, \label{mvssinf7} \\
\dfrac{\widehat{\mathrm{U}}_{\mu}(z)}{\widehat{\mathrm{V}}_{\mu}(z)} 
\! - \! \sum_{j=1}^{2 \varkappa_{nk}}c^{(\infty)}_{j}z^{-j} \underset{z 
\to \alpha_{k}}{=} \mathcal{O} \left(z^{-(2 \varkappa_{nk}+1)} \right). 
\label{mvssinf8} 
\end{gather}
\begin{eeee} \label{mvssreminf} 
\textsl{Note that Equations~\eqref{mvssinf7} and~\eqref{mvssinf8} give 
rise to $2 \sum_{q=1}^{\mathfrak{s}-1} \varkappa_{nk \tilde{k}_{q}}$ 
and $2 \varkappa_{nk}$ interpolation conditions, respectively, for a 
combined total of $2(\sum_{q=1}^{\mathfrak{s}-1} \varkappa_{nk 
\tilde{k}_{q}} \! + \! \varkappa_{nk}) \! = \! 2((n \! - \! 1)K \! + \! k)$ 
conditions, which is precisely the number necessary in order 
to determine the coefficients $\lbrace \hat{r}_{j},\hat{t}_{j} 
\rbrace_{j=0}^{(n-1)K+k-1}$}.\footnote{Strictly speaking, $\hat{r}_{j} 
\! = \! \hat{r}_{j}(n,k),\hat{t}_{j} \! = \! \hat{t}_{j}(n,k)$, $j \! = \! 0,1,
\dotsc,(n \! - \! 1)K \! + \! k \! - \! 1$, and $\hat{t}_{(n-1)K+k} \! 
= \! \hat{t}_{(n-1)K+k}(n,k) \! = \! 1$.}
\end{eeee}
For $n \! \in \! \mathbb{N}$ and $k \! \in \! \lbrace 1,2,\dotsc,K \rbrace$ 
such that $\alpha_{p_{\mathfrak{s}}} \! := \! \alpha_{k} \! = \! \infty$, 
define the corresponding MPA \emph{error term} as follows:\footnote{In 
this context, \emph{error term} means that, for $n \! \in \! \mathbb{N}$ 
and $k \! \in \! \lbrace 1,2,\dotsc,K \rbrace$ such that $\alpha_{
p_{\mathfrak{s}}} \! := \! \alpha_{k} \! = \! \infty$, $\widehat{\pmb{
\mathrm{E}}}_{\mu}(z) \! =_{z \to \alpha_{p_{q}}} \! \mathcal{O}((z \! - \! 
\alpha_{p_{q}})^{2 \varkappa_{nk \tilde{k}_{q}}})$, $q \! = \! 1,2,\dotsc,
\mathfrak{s} \! - \! 1$, and $\widehat{\pmb{\mathrm{E}}}_{\mu}(z) \! 
=_{z \to \alpha_{k}} \! \mathcal{O}(z^{-(2 \varkappa_{nk}+1)})$.}
\begin{equation} \label{mvssinf9} 
\widehat{\pmb{\mathrm{E}}}_{\mu}(z) \! := \! \dfrac{\widehat{
\pmb{\mathrm{R}}}_{\mu}(z)}{\pmb{\pi}^{n}_{k}(z)} \! - \! 
\mathrm{F}_{\mu}(z);
\end{equation}
a calculation based on Equations~\eqref{mvssinf1} and~\eqref{mvssinf4} 
reveals that, in fact,
\begin{equation} \label{mvssinf10}
\widehat{\pmb{\mathrm{E}}}_{\mu}(z) \! = \! \dfrac{1}{\pmb{\pi}^{n}_{k}
(z)} \int_{\mathbb{R}} \dfrac{\pmb{\pi}^{n}_{k}(\xi)}{\xi \! - \! z} \, 
\md \mu (\xi).
\end{equation}
\subsubsection{$n \! \in \! \mathbb{N}$ and $k \! \in \! \lbrace 1,2,\dotsc,
K \rbrace$ such that $\alpha_{k} \! \neq \! \infty$} \label{subsubsec1.2.2} 
For $k \! \in \! \lbrace 1,2,\dotsc,K \rbrace$ such that $\alpha_{k} \! \neq 
\! \infty$, write the ordered disjoint integer partition
\begin{align*}
\lbrace \mathstrut k^{\prime} \! \in \! \lbrace 1,2,\dotsc,K \rbrace; \, 
\alpha_{k^{\prime}} \! \neq \! \alpha_{k}, \, \alpha_{k} \! \neq \! 
\infty \rbrace \! :=& \, \lbrace \underbrace{i(1)_{1},i(1)_{2},\dotsc,
i(1)_{k_{1}}}_{k_{1}} \rbrace \cup \dotsb \cup \lbrace \underbrace{
i(\mathfrak{s} \! - \! 2)_{1},i(\mathfrak{s} \! - \! 2)_{2},\dotsc,
i(\mathfrak{s} \! - \! 2)_{k_{\mathfrak{s}-2}}}_{k_{\mathfrak{s}-2}} 
\rbrace \\
\cup& \, \lbrace \underbrace{i(\mathfrak{s} \! - \! 1)_{1},
i(\mathfrak{s} \! - \! 1)_{2},\dotsc,
i(\mathfrak{s} \! - \! 1)_{k_{\mathfrak{s}-1}}}_{k_{\mathfrak{s}-1}} \rbrace 
\! = \bigcup_{q=1}^{\mathfrak{s}-2} \lbrace \underbrace{i(q)_{1},
i(q)_{2},\dotsc,i(q)_{k_{q}}}_{k_{q}} \rbrace \\
\cup& \, \lbrace \underbrace{i(\mathfrak{s} \! - \! 1)_{1},
i(\mathfrak{s} \! - \! 1)_{2},\dotsc,
i(\mathfrak{s} \! - \! 1)_{k_{\mathfrak{s}-1}}}_{k_{\mathfrak{s}-1}} \rbrace,
\end{align*}
where $1 \! \leqslant \! i(q)_{1} \! < \! i(q)_{2} \! < \! \dotsb \! < \! 
i(q)_{k_{q}} \! \leqslant \! K$, $q \! \in \! \lbrace 1,2,\dotsc,\mathfrak{s} 
\! - \! 1 \rbrace$, $\lbrace i(j)_{1},i(j)_{2},\dotsc,i(j)_{k_{j}} \rbrace 
\cap \lbrace i(l)_{1},i(l)_{2},\dotsc,i(l)_{k_{l}} \rbrace \! = \! 
\varnothing$ $\forall$ $l \! \neq \! j \! \in \! \lbrace 1,2,\dotsc,
\mathfrak{s} \! - \! 1 \rbrace$, and $\# \lbrace \mathstrut k^{\prime} \! \in 
\! \lbrace 1,2,\dotsc,K \rbrace; \, \alpha_{k^{\prime}} \! \neq \! \alpha_{k}, 
\, \alpha_{k} \! \neq \! \infty \rbrace \! = \! \sum_{q=1}^{\mathfrak{s}-1}
k_{q} \! = \! K \! - \! \gamma_{k}$, which induces, on the pole set $\lbrace 
\alpha_{1},\alpha_{2},\dotsc,\alpha_{K} \rbrace$, the following disjoint 
ordering,
\begin{align*}
\lbrace \mathstrut \alpha_{k^{\prime}}, \, k^{\prime} \! \in \! \lbrace 1,2,
\dotsc,K \rbrace; \, \alpha_{k^{\prime}} \! \neq \! \alpha_{k}, \, \alpha_{k} 
\! \neq \! \infty \rbrace \! :=& \, \lbrace \underbrace{\alpha_{i(1)_{1}},
\alpha_{i(1)_{2}},\dotsc,\alpha_{i(1)_{k_{1}}}}_{k_{1}} \rbrace 
\cup \dotsb \cup \lbrace \underbrace{\alpha_{i(\mathfrak{s}-2)_{1}},
\alpha_{i(\mathfrak{s}-2)_{2}},\dotsc,
\alpha_{i(\mathfrak{s}-2)_{k_{\mathfrak{s}-2}}}}_{k_{\mathfrak{s}-2}} 
\rbrace \\
\cup& \, \lbrace \underbrace{\alpha_{i(\mathfrak{s}-1)_{1}},
\alpha_{i(\mathfrak{s}-1)_{2}},\dotsc,
\alpha_{i(\mathfrak{s}-1)_{k_{\mathfrak{s}-1}}}}_{k_{\mathfrak{s}-1}} \rbrace 
\! = \bigcup_{q=1}^{\mathfrak{s}-2} \lbrace \underbrace{\alpha_{i(q)_{1}},
\alpha_{i(q)_{2}},\dotsc,\alpha_{i(q)_{k_{q}}}}_{k_{q}} \rbrace \\
\cup& \, \lbrace \underbrace{\alpha_{i(\mathfrak{s}-1)_{1}},
\alpha_{i(\mathfrak{s}-1)_{2}},\dotsc,
\alpha_{i(\mathfrak{s}-1)_{k_{\mathfrak{s}-1}}}}_{k_{\mathfrak{s}-1}} \rbrace,
\end{align*}
with $\alpha_{i(q)_{1}} \! \prec \! \alpha_{i(q)_{2}} \! \prec \! \dotsb \! 
\prec \! \alpha_{i(q)_{k_{q}}}$, where, as above, the notation $a \! \prec 
\! b$ means `$a$ precedes $b$' or `$a$ is to the left of $b$', and $\lbrace 
\alpha_{i(j)_{1}},\alpha_{i(j)_{2}},\dotsc,\alpha_{i(j)_{k_{j}}} \rbrace \cap 
\lbrace \alpha_{i(l)_{1}},\alpha_{i(l)_{2}},\dotsc,\alpha_{i(l)_{k_{l}}} 
\rbrace \! = \! \varnothing$ $\forall$ $l \! \neq \! j \! \in \! \lbrace 1,2,
\dotsc,\mathfrak{s} \! - \! 1 \rbrace$, such that {}\footnote{If all the poles 
are distinct, that is, $\alpha_{i} \! \neq \! \alpha_{j}$ $\forall$ $i \! 
\neq \! j \! \in \! \lbrace 1,2,\dotsc,K \rbrace$, then, for $k \! \in \! 
\lbrace 1,2,\dotsc,K \rbrace$ such that $\alpha_{k} \! \neq \! \infty$, 
$\lbrace \mathstrut k^{\prime} \! \in \! \lbrace 1,2,\dotsc,K \rbrace; \, 
\alpha_{k^{\prime}} \! \neq \! \alpha_{k}, \, \alpha_{k} \! \neq \! 
\infty \rbrace$ is the ordered disjoint union of singletons, that is, 
$\cup_{q=1}^{K-1} \lbrace i(q)_{k_{q}} \rbrace$, with $k_{q} \! = \! 1$, 
$q \! = \! 1,2,\dotsc,K \! - \! 1$, $1 \! \leqslant \! i(1)_{1} \! < \! i(2)_{1} 
\! < \! \dotsb \! < \! i(K \! - \! 1)_{1} \! \leqslant \! K$, and $\lbrace 
i(q)_{k_{q}} \rbrace \cap \lbrace i(r)_{k_{r}} \rbrace \! = \! \varnothing$ 
$\forall$ $q \! \neq \! r \! \in \! \lbrace 1,2,\dotsc,K \! - \! 1 \rbrace$, 
which induces, on the pole set $\lbrace \alpha_{1},\alpha_{2},\dotsc,
\alpha_{K} \rbrace$, the following disjoint ordering, $\lbrace \mathstrut 
\alpha_{k^{\prime}}, \, k^{\prime} \! \in \! \lbrace 1,2,\dotsc,K \rbrace; \, 
\alpha_{k^{\prime}} \! \neq \! \alpha_{k}, \, \alpha_{k} \! \neq \! \infty 
\rbrace \! := \! \cup_{q=1}^{K-1} \lbrace \alpha_{i(q)_{k_{q}}} \rbrace \! 
= \! \cup_{q=1}^{K-2} \lbrace \alpha_{i(q)_{k_{q}}} \rbrace \cup \lbrace 
\alpha_{i(K-1)_{k_{K-1}}} \rbrace$, where $\alpha_{i(q)_{k_{q}}} \! \neq \! 
\infty$, $q \! = \! 1,2,\dotsc,K\! - \! 2$, and $\alpha_{i(K-1)_{k_{K-1}}} \! 
= \! \infty$, with $\alpha_{i(1)_{1}} \! \prec \! \alpha_{i(2)_{1}} \! \prec 
\! \dotsb \! \prec \! \alpha_{i(K-1)_{1}}$, $\# \lbrace \alpha_{i(q)_{k_{q}}} 
\rbrace \! = \! k_{q} \! = \! 1$, $q \! = \! 1,2,\dotsc,K \! - \! 1$, $\lbrace 
\alpha_{i(q)_{k_{q}}} \rbrace \cap \lbrace \alpha_{i(r)_{k_{r}}} \rbrace \! = 
\! \varnothing$ $\forall$ $q \! \neq \! r \! \in \! \lbrace 1,2,\dotsc,K \! - 
\! 1 \rbrace$, and $\# \lbrace \mathstrut \alpha_{k^{\prime}}, \, k^{\prime} 
\! \in \! \lbrace 1,2,\dotsc,K \rbrace; \, \alpha_{k^{\prime}} \! \neq \! 
\alpha_{k}, \, \alpha_{k} \! \neq \! \infty \rbrace \! = \! \sum_{q=1}^{K-1}
k_{q} \! = \! K \! - \! 1$.}
\begin{gather*}
\alpha_{i(q)_{1}} \! = \! \alpha_{i(q)_{2}} \! = \! \dotsb \! = \! 
\alpha_{i(q)_{k_{q}}} \, (\neq \! \infty), \quad q \! = \! 1,2,\dotsc,
\mathfrak{s} \! - \! 2, \\
\alpha_{i(\mathfrak{s}-1)_{1}} \! = \! \alpha_{i(\mathfrak{s}-1)_{2}} \! = \! 
\dotsb \! = \! \alpha_{i(\mathfrak{s}-1)_{k_{\mathfrak{s}-1}}} \! = \! \infty, 
\\
\# \lbrace \alpha_{i(q)_{1}},\alpha_{i(q)_{2}},\dotsc,\alpha_{i(q)_{k_{q}}} 
\rbrace \! = \! k_{q}, \quad q \! = \! 1,2,\dotsc,\mathfrak{s} \! - \! 1.
\end{gather*}
In order to illustrate this notation, consider the $K \! = \! 7$ pole set 
$\lbrace \alpha_{1},\alpha_{2},\alpha_{3},\alpha_{4},\alpha_{5},\alpha_{6},
\alpha_{7} \rbrace \! = \! \lbrace 0,1,\infty,1,\sqrt{2},\pi,\infty \rbrace$ 
for which $\mathfrak{s} \! = \! 5$ and $\alpha_{k} \! \neq \! \infty$, 
$k \! = \! 1,2,4,5,6$:
\begin{enumerate}
\item[(i)] \fbox{$k \! = \! 1$}
\begin{align*}
\lbrace \mathstrut k^{\prime} \! \in \! \lbrace 1,2,\dotsc,7 \rbrace; \, 
\alpha_{k^{\prime}} \! \neq \! \alpha_{1} \! = \! 0 \rbrace =& \, \lbrace 
2,3,4,5,6,7 \rbrace \! = \! \lbrace 2,4 \rbrace \cup \lbrace 5 \rbrace \cup 
\lbrace 6 \rbrace \cup \lbrace 3,7 \rbrace \\
:=& \, \lbrace i(1)_{1},i(1)_{k_{1}} \rbrace \cup \lbrace i(2)_{k_{2}} \rbrace 
\cup \lbrace i(3)_{k_{3}} \rbrace \cup \lbrace i(4)_{1},i(4)_{k_{4}} \rbrace 
\, \Rightarrow \\
k_{1} \! = \! 2, \, \, i(1)_{1} \! = \! 2, \, \, i(1)_{2} \! = \! 4, \, \, k_{2} 
\! = \! 1, \, \, &i(2)_{1} \! = \! 5, \, \, k_{3} \! = \! 1, \, \, i(3)_{1} \! = 
\! 6, \, \, k_{4} \! = \! 2, \, \, i(4)_{1} \! = \! 3, \, \, i(4)_{2} \! = \! 7
\end{align*}
(note that $1 \! \leqslant \! i(1)_{1} \! < \! i(1)_{2} \! \leqslant \! 7$, 
$1 \! \leqslant \! i(2)_{1} \! \leqslant \! 7$, $1 \! \leqslant \! i(3)_{1} 
\! \leqslant \! 7$, $1 \! \leqslant \! i(4)_{1} \! < \! i(4)_{2} \! \leqslant 
\! 7$, $\lbrace i(1)_{1},i(1)_{2} \rbrace \cap \lbrace i(2)_{1} \rbrace \! = 
\! \lbrace i(1)_{1},i(1)_{2} \rbrace \cap \lbrace i(3)_{1} \rbrace \! = \! 
\lbrace i(1)_{1},i(1)_{2} \rbrace \cap \lbrace i(4)_{1},i(4)_{2} \rbrace \! = 
\! \lbrace i(2)_{1} \rbrace \cap \lbrace i(3)_{1} \rbrace \! = \! \lbrace 
i(2)_{1} \rbrace \cap \lbrace i(4)_{1},i(4)_{2} \rbrace \! = \! \lbrace 
i(3)_{1} \rbrace \cap \lbrace i(4)_{1},i(4)_{2} \rbrace \! = \! \varnothing$, 
and $\# \lbrace \mathstrut k^{\prime} \! \in \! \lbrace 1,2,\dotsc,7 \rbrace; 
\, \alpha_{k^{\prime}} \! \neq \! \alpha_{1} \! = \! 0 \rbrace \! = \! 
\sum_{q=1}^{4}k_{q} \! = \! 2 \! + \! 1 \! + \! 1 \! + \! 2 = \! 6 \! = \! K 
\! - \! \gamma_{1} \! = \! 7 \! - \! 1)$, which induces the ordering (on the 
`residual' pole set)
\begin{align*}
\lbrace \mathstrut \alpha_{k^{\prime}}, \, k^{\prime} \! \in \! \lbrace 1,
2,\dotsc,7 \rbrace; \, \alpha_{k^{\prime}} \! \neq \! \alpha_{1} \! = \! 0 
\rbrace :=& \, \lbrace \alpha_{i(1)_{1}},\alpha_{i(1)_{k_{1}}} \rbrace \cup 
\lbrace \alpha_{i(2)_{k_{2}}} \rbrace \cup \lbrace \alpha_{i(3)_{k_{3}}} 
\rbrace \cup \lbrace \alpha_{i(4)_{1}},\alpha_{i(4)_{k_{4}}} \rbrace \\
=& \, \lbrace \alpha_{2},\alpha_{4} \rbrace \cup \lbrace \alpha_{5} \rbrace 
\cup \lbrace \alpha_{6} \rbrace \cup \lbrace \alpha_{3},\alpha_{7} \rbrace 
\! = \! \lbrace 1,1 \rbrace \cup \lbrace \sqrt{2} \rbrace \cup \lbrace \pi 
\rbrace \cup \lbrace \infty,\infty \rbrace
\end{align*}
(note that $\alpha_{i(1)_{1}} \! \prec \! \alpha_{i(1)_{2}}$, 
$\alpha_{i(4)_{1}} \! \prec \! \alpha_{i(4)_{2}}$, $\lbrace \alpha_{i(1)_{1}},
\alpha_{i(1)_{2}} \rbrace \cap \lbrace \alpha_{i(2)_{1}} \rbrace \! = \! 
\lbrace \alpha_{i(1)_{1}},\alpha_{i(1)_{2}} \rbrace \cap \lbrace 
\alpha_{i(3)_{1}} \rbrace \! = \! \lbrace \alpha_{i(1)_{1}},\alpha_{i(1)_{2}} 
\rbrace \cap \lbrace \alpha_{i(4)_{1}},\alpha_{i(4)_{2}} \rbrace \! = \! 
\lbrace \alpha_{i(2)_{1}} \rbrace \cap \lbrace \alpha_{i(3)_{1}} \rbrace 
\! = \! \lbrace \alpha_{i(2)_{1}} \rbrace \cap \lbrace \alpha_{i(4)_{1}},
\alpha_{i(4)_{2}} \rbrace \! = \! \lbrace \alpha_{i(3)_{1}} \rbrace \cap 
\lbrace \alpha_{i(4)_{1}},\alpha_{i(4)_{2}} \rbrace \! = \! \varnothing$, 
$\# \lbrace \alpha_{i(1)_{1}},\alpha_{i(1)_{2}} \rbrace \! = \! k_{1} \! = 
\! 2$, $\# \lbrace \alpha_{i(2)_{1}} \rbrace \! = \! k_{2} \! = \! 1$, $\# 
\lbrace \alpha_{i(3)_{1}} \rbrace \! = \! k_{3} \! = \! 1$, and $\# \lbrace 
\alpha_{i(4)_{1}},\alpha_{i(4)_{2}} \rbrace \! = \! k_{4} \! = \! 2)$;
\item[(ii)] \fbox{$k \! = \! 2$}
\begin{align*}
\lbrace \mathstrut k^{\prime} \! \in \! \lbrace 1,2,\dotsc,7 \rbrace; \, 
\alpha_{k^{\prime}} \! \neq \! \alpha_{2} \! = \! 1 \rbrace =& \, \lbrace 
1,3,5,6,7 \rbrace \! = \! \lbrace 1 \rbrace \cup \lbrace 5 \rbrace \cup 
\lbrace 6 \rbrace \cup \lbrace 3,7 \rbrace \\
:=& \, \lbrace i(1)_{k_{1}} \rbrace \cup \lbrace i(2)_{k_{2}} \rbrace \cup 
\lbrace i(3)_{k_{3}} \rbrace \cup \lbrace i(4)_{1},i(4)_{k_{4}} \rbrace \, 
\Rightarrow \\
k_{1} \! = \! 1, \, \, i(1)_{1} \! = \! 1, \, \, k_{2} \! = \! 1, \, \, 
&i(2)_{1} \! = \! 5, \, \, k_{3} \! = \! 1, \, \, i(3)_{1} \! = \! 6, \, \, 
k_{4} \! = \! 2, \, \, i(4)_{1} \! = \! 3, \, \, i(4)_{2} \! = \! 7
\end{align*}
(note that $1 \! \leqslant \! i(1)_{1} \! \leqslant \! 7$, $1 \! \leqslant \! 
i(2)_{1} \! \leqslant \! 7$, $1 \! \leqslant \! i(3)_{1} \! \leqslant \! 7$, 
$1 \! \leqslant \! i(4)_{1} \! < \! i(4)_{2} \! \leqslant \! 7$, $\lbrace 
i(1)_{1} \rbrace \cap \lbrace i(2)_{1} \rbrace \! = \! \lbrace i(1)_{1} 
\rbrace \cap \lbrace i(3)_{1} \rbrace \! = \! \lbrace i(1)_{1} \rbrace \cap 
\lbrace i(4)_{1},i(4)_{2} \rbrace \! = \! \lbrace i(2)_{1} \rbrace \cap 
\lbrace i(3)_{1} \rbrace \! = \! \lbrace i(2)_{1} \rbrace \cap \lbrace 
i(4)_{1},i(4)_{2} \rbrace \! = \! \lbrace i(3)_{1} \rbrace \cap \lbrace 
i(4)_{1},i(4)_{2} \rbrace \! = \! \varnothing$, and $\# \lbrace \mathstrut 
k^{\prime} \! \in \! \lbrace 1,2,\dotsc,7 \rbrace; \, \alpha_{k^{\prime}} \! 
\neq \! \alpha_{2} \! = \! 1 \rbrace \! = \! \sum_{q=1}^{4}k_{q} \! = \! 1 
\! + \! 1 \! + \! 1 \! + \! 2 = \! 5 \! = \! K \! - \! \gamma_{2} \! = \! 7 
\! - \! 2)$, which induces the ordering (on the `residual' pole set)
\begin{align*}
\lbrace \mathstrut \alpha_{k^{\prime}}, \, k^{\prime} \! \in \! \lbrace 1,2,
\dotsc,7 \rbrace; \, \alpha_{k^{\prime}} \! \neq \! \alpha_{2} \! = \! 1 
\rbrace :=& \, \lbrace \alpha_{i(1)_{k_{1}}} \rbrace \cup \lbrace 
\alpha_{i(2)_{k_{2}}} \rbrace \cup \lbrace \alpha_{i(3)_{k_{3}}} \rbrace \cup 
\lbrace \alpha_{i(4)_{1}},\alpha_{i(4)_{k_{4}}} \rbrace \\
=& \, \lbrace \alpha_{1} \rbrace \cup \lbrace \alpha_{5} \rbrace \cup \lbrace 
\alpha_{6} \rbrace \cup \lbrace \alpha_{3},\alpha_{7} \rbrace \! = \! \lbrace 
0 \rbrace \cup \lbrace \sqrt{2} \rbrace \cup \lbrace \pi \rbrace \cup \lbrace 
\infty,\infty \rbrace
\end{align*}
(note that $\alpha_{i(4)_{1}} \! \prec \! \alpha_{i(4)_{2}}$, $\lbrace 
\alpha_{i(1)_{1}} \rbrace \cap \lbrace \alpha_{i(2)_{1}} \rbrace \! = \! 
\lbrace \alpha_{i(1)_{1}} \rbrace \cap \lbrace \alpha_{i(3)_{1}} \rbrace 
\! = \! \lbrace \alpha_{i(1)_{1}} \rbrace \cap \lbrace \alpha_{i(4)_{1}},
\alpha_{i(4)_{2}} \rbrace \! = \! \lbrace \alpha_{i(2)_{1}} \rbrace \cap 
\lbrace \alpha_{i(3)_{1}} \rbrace \! = \! \lbrace \alpha_{i(2)_{1}} \rbrace 
\cap \lbrace \alpha_{i(4)_{1}},\alpha_{i(4)_{2}} \rbrace \! = \! \lbrace 
\alpha_{i(3)_{1}} \rbrace \cap \lbrace \alpha_{i(4)_{1}},\alpha_{i(4)_{2}} 
\rbrace \! = \! \varnothing$, $\# \lbrace \alpha_{i(1)_{1}} \rbrace \! = \! 
k_{1} \! = \! 1$, $\# \lbrace \alpha_{i(2)_{1}} \rbrace \! = \! k_{2} \! = \! 
1$, $\# \lbrace \alpha_{i(3)_{1}} \rbrace \! = \! k_{3} \! = \! 1$, and $\# 
\lbrace \alpha_{i(4)_{1}},\alpha_{i(4)_{2}} \rbrace \! = \! k_{4} \! = \! 2)$;
\item[(iii)] \fbox{$k \! = \! 4$}
\begin{align*}
\lbrace \mathstrut k^{\prime} \! \in \! \lbrace 1,2,\dotsc,7 \rbrace; \, 
\alpha_{k^{\prime}} \! \neq \! \alpha_{4} \! = \! 1 \rbrace =& \, \lbrace 
1,3,5,6,7 \rbrace \! = \! \lbrace 1 \rbrace \cup \lbrace 5 \rbrace \cup 
\lbrace 6 \rbrace \cup \lbrace 3,7 \rbrace \\
:=& \, \lbrace i(1)_{k_{1}} \rbrace \cup \lbrace i(2)_{k_{2}} \rbrace \cup 
\lbrace i(3)_{k_{3}} \rbrace \cup \lbrace i(4)_{1},i(4)_{k_{4}} \rbrace \, 
\Rightarrow \\
k_{1} \! = \! 1, \, \, i(1)_{1} \! = \! 1, \, \, k_{2} \! = \! 1, \, \, 
&i(2)_{1} \! = \! 5, \, \, k_{3} \! = \! 1, \, \, i(3)_{1} \! = \! 6, \, \, 
k_{4} \! = \! 2, \, \, i(4)_{1} \! = \! 3, \, \, i(4)_{2} \! = \! 7
\end{align*}
(note that $1 \! \leqslant \! i(1)_{1} \! \leqslant \! 7$, $1 \! \leqslant \! 
i(2)_{1} \! \leqslant \! 7$, $1 \! \leqslant \! i(3)_{1} \! \leqslant \! 7$, 
$1 \! \leqslant \! i(4)_{1} \! < \! i(4)_{2} \! \leqslant \! 7$, $\lbrace 
i(1)_{1} \rbrace \cap \lbrace i(2)_{1} \rbrace \! = \! \lbrace i(1)_{1} 
\rbrace \cap \lbrace i(3)_{1} \rbrace \! = \! \lbrace i(1)_{1} \rbrace \cap 
\lbrace i(4)_{1},i(4)_{2} \rbrace \! = \! \lbrace i(2)_{1} \rbrace \cap 
\lbrace i(3)_{1} \rbrace \! = \! \lbrace i(2)_{1} \rbrace \cap \lbrace 
i(4)_{1},i(4)_{2} \rbrace \! = \! \lbrace i(3)_{1} \rbrace \cap \lbrace 
i(4)_{1},i(4)_{2} \rbrace \! = \! \varnothing$, and $\# \lbrace \mathstrut 
k^{\prime} \! \in \! \lbrace 1,2,\dotsc,7 \rbrace; \, \alpha_{k^{\prime}} \! 
\neq \! \alpha_{4} \! = \! 1 \rbrace \! = \! \sum_{q=1}^{4}k_{q} \! = \! 1 
\! + \! 1 \! + \! 1 \! + \! 2 = \! 5 \! = \! K \! - \! \gamma_{4} \! = \! 7 
\! - \! 2)$, which induces the ordering (on the `residual' pole set)
\begin{align*}
\lbrace \mathstrut \alpha_{k^{\prime}}, \, k^{\prime} \! \in \! \lbrace 1,2,
\dotsc,7 \rbrace; \, \alpha_{k^{\prime}} \! \neq \! \alpha_{4} \! = \! 1 
\rbrace :=& \, \lbrace \alpha_{i(1)_{k_{1}}} \rbrace \cup \lbrace 
\alpha_{i(2)_{k_{2}}} \rbrace \cup \lbrace \alpha_{i(3)_{k_{3}}} \rbrace 
\cup \lbrace \alpha_{i(4)_{1}},\alpha_{i(4)_{k_{4}}} \rbrace \\
=& \, \lbrace \alpha_{1} \rbrace \cup \lbrace \alpha_{5} \rbrace \cup \lbrace 
\alpha_{6} \rbrace \cup \lbrace \alpha_{3},\alpha_{7} \rbrace \! = \! \lbrace 
0 \rbrace \cup \lbrace \sqrt{2} \rbrace \cup \lbrace \pi \rbrace \cup \lbrace 
\infty,\infty \rbrace
\end{align*}
(note that $\alpha_{i(4)_{1}} \! \prec \! \alpha_{i(4)_{2}}$, $\lbrace 
\alpha_{i(1)_{1}} \rbrace \cap \lbrace \alpha_{i(2)_{1}} \rbrace \! = \! 
\lbrace \alpha_{i(1)_{1}} \rbrace \cap \lbrace \alpha_{i(3)_{1}} \rbrace 
\! = \! \lbrace \alpha_{i(1)_{1}} \rbrace \cap \lbrace \alpha_{i(4)_{1}},
\alpha_{i(4)_{2}} \rbrace \! = \! \lbrace \alpha_{i(2)_{1}} \rbrace \cap 
\lbrace \alpha_{i(3)_{1}} \rbrace \! = \! \lbrace \alpha_{i(2)_{1}} \rbrace 
\cap \lbrace \alpha_{i(4)_{1}},\alpha_{i(4)_{2}} \rbrace \! = \! \lbrace 
\alpha_{i(3)_{1}} \rbrace \cap \lbrace \alpha_{i(4)_{1}},\alpha_{i(4)_{2}} 
\rbrace \! = \! \varnothing$, $\# \lbrace \alpha_{i(1)_{1}} \rbrace \! = \! 
k_{1} \! = \! 1$, $\# \lbrace \alpha_{i(2)_{1}} \rbrace \! = \! k_{2} \! = 
\! 1$, $\# \lbrace \alpha_{i(3)_{1}} \rbrace \! = \! k_{3} \! = \! 1$, and 
$\# \lbrace \alpha_{i(4)_{1}},\alpha_{i(4)_{2}} \rbrace \! = \! k_{4} \! = 
\! 2)$;
\item[(iv)] \fbox{$k \! = \! 5$}
\begin{align*}
\lbrace \mathstrut k^{\prime} \! \in \! \lbrace 1,2,\dotsc,7 \rbrace; \, 
\alpha_{k^{\prime}} \! \neq \! \alpha_{5} \! = \! \! \sqrt{2} \rbrace =& 
\, \lbrace 1,2,3,4,6,7 \rbrace \! = \! \lbrace 1 \rbrace \cup \lbrace 2,4 
\rbrace \cup \lbrace 6 \rbrace \cup \lbrace 3,7 \rbrace \\
:=& \, \lbrace i(1)_{k_{1}} \rbrace \cup \lbrace i(2)_{1},i(2)_{k_{2}} 
\rbrace \cup \lbrace i(3)_{k_{3}} \rbrace \cup \lbrace i(4)_{1},i(4)_{k_{4}} 
\rbrace \, \Rightarrow \\
k_{1} \! = \! 1, \, \, i(1)_{1} \! = \! 1, \, \, k_{2} \! = \! 2, \, \, 
i(2)_{1} \! = \! 2, \, \, &i(2)_{2} \! = \! 4, \, \, k_{3} \! = \! 1, 
\, \, i(3)_{1} \! = \! 6, \, \, k_{4} \! = \! 2, \, \, i(4)_{1} \! = \! 
3, \, \, i(4)_{2} \! = \! 7
\end{align*}
(note that $1 \! \leqslant \! i(1)_{1} \! \leqslant \! 7$, $1 \! \leqslant \! 
i(2)_{1} \! < \! i(2)_{2} \! \leqslant \! 7$, $1 \! \leqslant \! i(3)_{1} \! 
\leqslant \! 7$, $1 \! \leqslant \! i(4)_{1} \! < \! i(4)_{2} \! \leqslant 
\! 7$, $\lbrace i(1)_{1} \rbrace \cap \lbrace i(2)_{1},i(2)_{2} \rbrace \! = 
\! \lbrace i(1)_{1} \rbrace \cap \lbrace i(3)_{1} \rbrace \! = \! \lbrace 
i(1)_{1} \rbrace \cap \lbrace i(4)_{1},i(4)_{2} \rbrace \! = \! \lbrace 
i(2)_{1},i(2)_{2} \rbrace \cap \lbrace i(3)_{1} \rbrace \! = \! \lbrace 
i(2)_{1},i(2)_{2} \rbrace \cap \lbrace i(4)_{1},i(4)_{2} \rbrace \! = \! 
\lbrace i(3)_{1} \rbrace \cap \lbrace i(4)_{1},i(4)_{2} \rbrace \! = \! 
\varnothing$, and $\# \lbrace \mathstrut k^{\prime} \! \in \! \lbrace 1,2,
\dotsc,7 \rbrace; \, \alpha_{k^{\prime}} \! \neq \! \alpha_{5} \! = \! 
\sqrt{2} \rbrace \! = \! \sum_{q=1}^{4}k_{q} \! = \! 1 \! + \! 2 \! + \! 1 
\! + \! 2 = \! 6 \! = \! K \! - \! \gamma_{5} \! = \! 7 \! - \! 1)$, which 
induces the ordering (on the `residual' pole set)
\begin{align*}
\lbrace \mathstrut \alpha_{k^{\prime}}, \, k^{\prime} \! \in \! \lbrace 
1,2,\dotsc,7 \rbrace; \, \alpha_{k^{\prime}} \! \neq \! \alpha_{5} \! = 
\! \sqrt{2} \rbrace :=& \, \lbrace \alpha_{i(1)_{k_{1}}} \rbrace \cup 
\lbrace \alpha_{i(2)_{1}},\alpha_{i(2)_{k_{2}}} \rbrace \cup \lbrace 
\alpha_{i(3)_{k_{3}}} \rbrace \cup \lbrace \alpha_{i(4)_{1}},
\alpha_{i(4)_{k_{4}}} \rbrace \\
=& \, \lbrace \alpha_{1} \rbrace \cup \lbrace \alpha_{2},\alpha_{4} \rbrace 
\cup \lbrace \alpha_{6} \rbrace \cup \lbrace \alpha_{3},\alpha_{7} \rbrace 
\! = \! \lbrace 0 \rbrace \cup \lbrace 1,1 \rbrace \cup \lbrace \pi \rbrace 
\cup \lbrace \infty,\infty \rbrace
\end{align*}
(note that $\alpha_{i(2)_{1}} \! \prec \! \alpha_{i(2)_{2}}$, 
$\alpha_{i(4)_{1}} \! \prec \! \alpha_{i(4)_{2}}$, $\lbrace \alpha_{i(1)_{1}} 
\rbrace \cap \lbrace \alpha_{i(2)_{1}},\alpha_{i(2)_{2}} \rbrace \! = \! 
\lbrace \alpha_{i(1)_{1}} \rbrace \cap \lbrace \alpha_{i(3)_{1}} \rbrace 
\! = \! \lbrace \alpha_{i(1)_{1}} \rbrace \cap \lbrace \alpha_{i(4)_{1}},
\alpha_{i(4)_{2}} \rbrace \! = \! \lbrace \alpha_{i(2)_{1}},\alpha_{i(2)_{2}} 
\rbrace \cap \lbrace \alpha_{i(3)_{1}} \rbrace \! = \! \lbrace 
\alpha_{i(2)_{1}},\alpha_{i(2)_{2}} \rbrace \cap \lbrace \alpha_{i(4)_{1}},
\alpha_{i(4)_{2}} \rbrace \! = \! \lbrace \alpha_{i(3)_{1}} \rbrace \cap 
\lbrace \alpha_{i(4)_{1}},\alpha_{i(4)_{2}} \rbrace \! = \! \varnothing$, 
$\# \lbrace \alpha_{i(1)_{1}} \rbrace \! = \! k_{1} \! = \! 1$, $\# \lbrace 
\alpha_{i(2)_{1}},\alpha_{i(2)_{2}} \rbrace \! = \! k_{2} \! = \! 2$, $\# 
\lbrace \alpha_{i(3)_{1}} \rbrace \! = \! k_{3} \! = \! 1$, and $\# \lbrace 
\alpha_{i(4)_{1}},\alpha_{i(4)_{2}} \rbrace \! = \! k_{4} \! = \! 2)$;
\item[(v)] \fbox{$k \! = \! 6$}
\begin{align*}
\lbrace \mathstrut k^{\prime} \! \in \! \lbrace 1,2,\dotsc,7 \rbrace; \, 
\alpha_{k^{\prime}} \! \neq \! \alpha_{6} \! = \! \pi \rbrace =& \, \lbrace 
1,2,3,4,5,7 \rbrace \! = \! \lbrace 1 \rbrace \cup \lbrace 2,4 \rbrace \cup 
\lbrace 5 \rbrace \cup \lbrace 3,7 \rbrace \\
:=& \, \lbrace i(1)_{k_{1}} \rbrace \cup \lbrace i(2)_{1},i(2)_{k_{2}} \rbrace 
\cup \lbrace i(3)_{k_{3}} \rbrace \cup \lbrace i(4)_{1},i(4)_{k_{4}} \rbrace 
\, \Rightarrow \\
k_{1} \! = \! 1, \, \, i(1)_{1} \! = \! 1, \, \, k_{2} \! = \! 2, \, \, 
i(2)_{1} \! = \! 2, \, \, &i(2)_{2} \! = \! 4, \, \, k_{3} \! = \! 1, 
\, \, i(3)_{1} \! = \! 5, \, \, k_{4} \! = \! 2, \, \, i(4)_{1} \! = \! 
3, \, \, i(4)_{2} \! = \! 7
\end{align*}
(note that $1 \! \leqslant \! i(1)_{1} \! \leqslant \! 7$, $1 \! \leqslant \! 
i(2)_{1} \! < \! i(2)_{2} \! \leqslant \! 7$, $1 \! \leqslant \! i(3)_{1} \! 
\leqslant \! 7$, $1 \! \leqslant \! i(4)_{1} \! < \! i(4)_{2} \! \leqslant \! 
7$, $\lbrace i(1)_{1} \rbrace \cap \lbrace i(2)_{1},i(2)_{2} \rbrace \! = \! 
\lbrace i(1)_{1} \rbrace \cap \lbrace i(3)_{1} \rbrace \! = \! \lbrace 
i(1)_{1} \rbrace \cap \lbrace i(4)_{1},i(4)_{2} \rbrace \! = \! \lbrace 
i(2)_{1},i(2)_{2} \rbrace \cap \lbrace i(3)_{1} \rbrace \! = \! \lbrace 
i(2)_{1},i(2)_{2} \rbrace \cap \lbrace i(4)_{1},i(4)_{2} \rbrace \! = \! 
\lbrace i(3)_{1} \rbrace \cap \lbrace i(4)_{1},i(4)_{2} \rbrace \! = \! 
\varnothing$, and $\# \lbrace \mathstrut k^{\prime} \! \in \! \lbrace 1,2,
\dotsc,7 \rbrace; \, \alpha_{k^{\prime}} \! \neq \! \alpha_{6} \! = \! \pi 
\rbrace \! = \! \sum_{q=1}^{4}k_{q} \! = \! 1 \! + \! 2 \! + \! 1 \! + \! 2 
= \! 6 \! = \! K \! - \! \gamma_{6} \! = \! 7 \! - \! 1)$, which induces the 
ordering (on the `residual' pole set)
\begin{align*}
\lbrace \mathstrut \alpha_{k^{\prime}}, \, k^{\prime} \! \in \! \lbrace 
1,2,\dotsc,7 \rbrace; \, \alpha_{k^{\prime}} \! \neq \! \alpha_{6} \! = 
\! \pi \rbrace :=& \, \lbrace \alpha_{i(1)_{k_{1}}} \rbrace \cup \lbrace 
\alpha_{i(2)_{1}},\alpha_{i(2)_{k_{2}}} \rbrace \cup \lbrace 
\alpha_{i(3)_{k_{3}}} \rbrace \cup \lbrace \alpha_{i(4)_{1}},
\alpha_{i(4)_{k_{4}}} \rbrace \\
=& \, \lbrace \alpha_{1} \rbrace \cup \lbrace \alpha_{2},\alpha_{4} \rbrace 
\cup \lbrace \alpha_{5} \rbrace \cup \lbrace \alpha_{3},\alpha_{7} \rbrace 
\! = \! \lbrace 0 \rbrace \cup \lbrace 1,1 \rbrace \cup \lbrace \sqrt{2} 
\rbrace \cup \lbrace \infty,\infty \rbrace
\end{align*}
(note that $\alpha_{i(2)_{1}} \! \prec \! \alpha_{i(2)_{2}}$, 
$\alpha_{i(4)_{1}} \! \prec \! \alpha_{i(4)_{2}}$, $\lbrace \alpha_{i(1)_{1}} 
\rbrace \cap \lbrace \alpha_{i(2)_{1}},\alpha_{i(2)_{2}} \rbrace \! = \! 
\lbrace \alpha_{i(1)_{1}} \rbrace \cap \lbrace \alpha_{i(3)_{1}} \rbrace 
\! = \! \lbrace \alpha_{i(1)_{1}} \rbrace \cap \lbrace \alpha_{i(4)_{1}},
\alpha_{i(4)_{2}} \rbrace \! = \! \lbrace \alpha_{i(2)_{1}},\alpha_{i(2)_{2}} 
\rbrace \cap \lbrace \alpha_{i(3)_{1}} \rbrace \! = \! \lbrace 
\alpha_{i(2)_{1}},\alpha_{i(2)_{2}} \rbrace \cap \lbrace \alpha_{i(4)_{1}},
\alpha_{i(4)_{2}} \rbrace \! = \! \lbrace \alpha_{i(3)_{1}} \rbrace \cap 
\lbrace \alpha_{i(4)_{1}},\alpha_{i(4)_{2}} \rbrace \! = \! \varnothing$, 
$\# \lbrace \alpha_{i(1)_{1}} \rbrace \! = \! k_{1} \! = \! 1$, $\# \lbrace 
\alpha_{i(2)_{1}},\alpha_{i(2)_{2}} \rbrace \! = \! k_{2} \! = \! 2$, $\# 
\lbrace \alpha_{i(3)_{1}} \rbrace \! = \! k_{3} \! = \! 1$, and $\# \lbrace 
\alpha_{i(4)_{1}},\alpha_{i(4)_{2}} \rbrace \! = \! k_{4} \! = \! 2)$.
\end{enumerate}
This concludes the example.

For $k \! \in \! \lbrace 1,2,\dotsc,K \rbrace$ such that $\alpha_{k} 
\! \neq \! \infty$, let (recall Definition~\ref{def1.2.2})
\begin{equation*}
\hat{\mathfrak{J}}_{q}(k) \! := \! \lbrace \operatorname{ind} \lbrace 
i(q)_{1},i(q)_{2},\dotsc,i(q)_{k_{q}} \vert k \rbrace \rbrace, \quad q 
\! = \! 1,2,\dotsc,\mathfrak{s} \! - \! 1:
\end{equation*}
if, for a given index set $\lbrace i(q)_{1},i(q)_{2},\dotsc,i(q)_{k_{q}} \rbrace$, 
$q \! = \! 1,2,\dotsc,\mathfrak{s} \! - \! 1$, there does not exist a positive 
integer $\tilde{\jmath}$ \footnote{This positive integer, if it exists, depends 
not only on $k$, but on $q$, too; however, for notational simplicity, this 
additional $q$ dependence is also suppressed.} such that $\tilde{\jmath} 
\! = \! \operatorname{ind} \lbrace i(q)_{1},i(q)_{2},\dotsc,i(q)_{k_{q}} \vert 
k \rbrace$, then put $\hat{\mathfrak{J}}_{q}(k) \! = \! \varnothing$; 
otherwise, denote by $\hat{m}_{q}(k)$, $q \! \in \! \lbrace 1,2,\dotsc,
\mathfrak{s} \! - \! 1 \rbrace$, the unique element of the (non-empty) set 
$\hat{\mathfrak{J}}_{q}(k)$.\footnote{Note that $\hat{\mathfrak{J}}_{q}(k)$ 
is either empty or a singleton.} For $n \! \in \! \mathbb{N}$ and $k \! \in 
\! \lbrace 1,2,\dotsc,K \rbrace$ such that $\alpha_{k} \! \neq \! \infty$, 
set, with the above orderings and definitions,
\begin{align*}
\varkappa_{nk \tilde{k}_{q}} \! :=& \, 
\begin{cases}
(n \! - \! 1) \gamma_{i(q)_{k_{q}}}, &\text{$\hat{\mathfrak{J}}_{q}(k) \! = 
\! \varnothing, \quad q \! \in \! \lbrace 1,2,\dotsc,\mathfrak{s} \! - \! 
2 \rbrace$,} \\
(n \! - \! 1) \gamma_{\hat{m}_{q}(k)} \! + \! \varrho_{\hat{m}_{q}(k)}, 
&\text{$\hat{\mathfrak{J}}_{q}(k) \! \neq \! \varnothing, \quad q \! \in \! 
\lbrace 1,2,\dotsc,\mathfrak{s} \! - \! 2 \rbrace$,}
\end{cases} \\
\varkappa^{\infty}_{nk \tilde{k}_{\mathfrak{s}-1}} \! :=& \, 
\begin{cases}
(n \! - \! 1) \gamma_{i(\mathfrak{s}-1)_{k_{\mathfrak{s}-1}}}, 
&\text{$\hat{\mathfrak{J}}_{\mathfrak{s}-1}(k) \! = \! \varnothing$,} \\
(n \! - \! 1) \gamma_{\hat{m}_{\mathfrak{s}-1}(k)} \! + \! 
\varrho_{\hat{m}_{\mathfrak{s}-1}(k)}, 
&\text{$\hat{\mathfrak{J}}_{\mathfrak{s}-1}(k) \! \neq \! \varnothing$,}
\end{cases}
\end{align*}
where $\varkappa_{nk \tilde{k}_{q}}$ (resp., $\varkappa^{\infty}_{nk 
\tilde{k}_{\mathfrak{s}-1}})$ $\colon \mathbb{N} \times \lbrace 1,2,
\dotsc,K \rbrace \! \to \! \mathbb{N}_{0}$, $q \! = \! 1,2,\dotsc,
\mathfrak{s} \! - \! 2$, is the residual multiplicity of the pole 
$\alpha_{i(q)_{k_{q}}}$ (resp., the point at infinity $\alpha_{i(\mathfrak{s}
-1)_{k_{\mathfrak{s}-1}}} \! = \! \infty)$, with $\alpha_{i(q)_{k_{q}}} \! 
\neq \! \alpha_{k}$ and $\alpha_{i(q)_{k_{q}}} \! \neq \! \infty$, in the 
repeated pole sequence
\begin{equation*}
\mathcal{P}_{n,k} \! := \! \lbrace \overset{1}{\underbrace{\alpha_{1},
\alpha_{2},\dotsc,\alpha_{K}}_{K}} \rbrace \cup \dotsb \cup \lbrace 
\overset{n-1}{\underbrace{\alpha_{1},\alpha_{2},\dotsc,\alpha_{K}}_{K}} 
\rbrace \cup \lbrace \overset{n}{\underbrace{\alpha_{1},\alpha_{2},\dotsc,
\alpha_{k}}_{k}} \rbrace,
\end{equation*}
that is, for $n \! \in \! \mathbb{N}$ and $k \! \in \! \lbrace 1,2,\dotsc,K 
\rbrace$ such that $\alpha_{k} \! \neq \! \infty$, as all occurrences of the 
pole $\alpha_{k}$ $(\neq \! \infty)$ are excised {}from the repeated pole 
sequence $\mathcal{P}_{n,k}$, where the multiplicity, or number of 
occurrences, of the pole $\alpha_{k}$ is $\varkappa_{nk} \! = \! (n \! - \! 
1) \gamma_{k} \! + \! \varrho_{k}$, one is left with the `residual' pole set 
(via the above induced ordering on the poles, as one walks across the pole 
sequence {}from left to right)
\begin{align*}
\mathcal{P}_{n,k} \setminus \lbrace \underbrace{\alpha_{k},\alpha_{k},
\dotsc,\alpha_{k}}_{\varkappa_{nk}} \rbrace \! :=& \, \bigcup_{q=1}^{
\mathfrak{s}-2} \lbrace \underbrace{\alpha_{i(q)_{k_{q}}},\alpha_{i(q)_{k_{q}}},
\dotsc,\alpha_{i(q)_{k_{q}}}}_{\varkappa_{nk \tilde{k}_{q}}} \rbrace \cup 
\lbrace \underbrace{\alpha_{i(\mathfrak{s}-1)_{k_{\mathfrak{s}-1}}},
\alpha_{i(\mathfrak{s}-1)_{k_{\mathfrak{s}-1}}},\dotsc,
\alpha_{i(\mathfrak{s}-1)_{k_{\mathfrak{s}-1}}}}_{\varkappa^{\infty}_{nk 
\tilde{k}_{\mathfrak{s}-1}}} \rbrace \\
:=& \, \lbrace \underbrace{\alpha_{i(1)_{k_{1}}},\alpha_{i(1)_{k_{1}}},\dotsc,
\alpha_{i(1)_{k_{1}}}}_{\varkappa_{nk \tilde{k}_{1}}} \rbrace \cup \dotsb 
\cup \lbrace \underbrace{\alpha_{i(\mathfrak{s}-2)_{k_{\mathfrak{s}-2}}},
\alpha_{i(\mathfrak{s}-2)_{k_{\mathfrak{s}-2}}},\dotsc,
\alpha_{i(\mathfrak{s}-2)_{k_{\mathfrak{s}-2}}}}_{\varkappa_{nk 
\tilde{k}_{\mathfrak{s}-2}}} \rbrace \\
\cup& \, \lbrace \underbrace{\infty,\infty,\dotsc,
\infty}_{\varkappa^{\infty}_{nk \tilde{k}_{\mathfrak{s}-1}}} \rbrace,
\end{align*}
where the number of times the pole $\alpha_{i(q)_{k_{q}}}$ (with 
$\alpha_{i(q)_{k_{q}}} \! \neq \! \alpha_{k}$ and $\alpha_{i(q)_{k_{q}}} \! 
\neq \! \infty)$ occurs is (its multiplicity) $\varkappa_{nk \tilde{k}_{q}}$, 
$q \! = \! 1,2,\dotsc,\mathfrak{s} \! - \! 2$, and the number of times the 
pole (the point at infinity) $\alpha_{i(\mathfrak{s}-1)_{k_{\mathfrak{s}
-1}}} \! = \! \infty$ $(\neq \! \alpha_{k})$ occurs is (its multiplicity) 
$\varkappa^{\infty}_{nk \tilde{k}_{\mathfrak{s}-1}}$.\footnote{Note: for 
$n \! = \! 1$ and $k \! \in \! \lbrace 1,2,\dotsc,K \rbrace$ such that 
$\alpha_{k} \! \neq \! \infty$, it can happen that $\varkappa_{1k \tilde{k}_{q}} 
\! = \! 0$, $q \! \in \! \lbrace 1,2,\dotsc,\mathfrak{s} \! - \! 2 \rbrace$, 
and $\varkappa^{\infty}_{1k \tilde{k}_{\mathfrak{s}-1}} \! = \! 0$, in 
which case, one defines, $\lbrace \alpha_{i(q^{\prime})_{k_{q^{\prime}}}},
\alpha_{i(q^{\prime})_{k_{q^{\prime}}}},\dotsc,\alpha_{i(q^{\prime})_{
k_{q^{\prime}}}} \rbrace \! := \! \varnothing$, $q^{\prime} \! \in \! 
\lbrace 1,2,\dotsc,\mathfrak{s} \! - \! 1 \rbrace$; however, for 
$(\mathbb{N} \! \ni)$ $n \! \geqslant \! 2$ and $k \! \in \! \lbrace 
1,2,\dotsc,K \rbrace$ such that $\alpha_{k} \! \neq \! \infty$, 
$\varkappa_{nk \tilde{k}_{q}} \! \geqslant \! 1$, $q \! \in \! \lbrace 
1,2,\dotsc,\mathfrak{s} \! - \! 2 \rbrace$, and $\varkappa^{\infty}_{nk 
\tilde{k}_{\mathfrak{s}-1}} \! \geqslant \! 1$.}

For $n \! \in \! \mathbb{N}$ and $k \! \in \! \lbrace 1,2,\dotsc,K 
\rbrace$ such that $\alpha_{k} \! \neq \! \infty$, a 
counting-of-residual-multiplicities argument gives rise to the following 
ordered sum formula:
\begin{equation} \label{fincount} 
\sum_{q=1}^{\mathfrak{s}-2} \varkappa_{nk \tilde{k}_{q}} \! 
+ \! \varkappa^{\infty}_{nk \tilde{k}_{\mathfrak{s}-1}} \! 
:= \! \varkappa_{nk \tilde{k}_{1}} \! + \! \dotsb \! + \! 
\varkappa_{nk \tilde{k}_{\mathfrak{s}-2}} \! + \! 
\varkappa^{\infty}_{nk \tilde{k}_{\mathfrak{s}-1}} \! = \! 
(n \! - \! 1)K \! + \! k \! - \! \varkappa_{nk},
\end{equation}
whence $\sum_{q=1}^{\mathfrak{s}-2} \varkappa_{nk \tilde{k}_{q}} \! + \! 
\varkappa^{\infty}_{nk \tilde{k}_{\mathfrak{s}-1}} \! + \! \varkappa_{nk} \! 
= \! (n \! - \! 1)K \! + \! k$.

In order to illustrate the latter notation, consider, again, the $K \! = \! 7$ 
pole set $\lbrace \alpha_{1},\alpha_{2},\alpha_{3},\alpha_{4},\alpha_{5},
\alpha_{6},\alpha_{7} \rbrace \! = \! \lbrace 0,1,\infty,1,\sqrt{2},\pi,\infty 
\rbrace$ for which $\mathfrak{s} \! = \! 5$ and $\alpha_{k} \! \neq \! 
\infty$, $k \! = \! 1,2,4,5,6$:
\begin{enumerate}
\item[(i)] \fbox{$k \! = \! 1$}
\begin{gather*}
\hat{\mathfrak{J}}_{1}(1) \! := \! \lbrace \operatorname{ind} \lbrace 
i(1)_{1},i(1)_{k_{1}} \vert 1 \rbrace \rbrace \! = \! \lbrace 
\operatorname{ind} \lbrace 2,4 \vert 1 \rbrace \rbrace \! = \! \varnothing, \\
\hat{\mathfrak{J}}_{2}(1) \! := \! \lbrace \operatorname{ind} \lbrace 
i(2)_{k_{2}} \vert 1 \rbrace \rbrace \! = \! \lbrace \operatorname{ind} 
\lbrace 5 \vert 1 \rbrace \rbrace \! = \! \varnothing, \\
\hat{\mathfrak{J}}_{3}(1) \! := \! \lbrace \operatorname{ind} \lbrace 
i(3)_{k_{3}} \vert 1 \rbrace \rbrace \! = \! \lbrace \operatorname{ind} 
\lbrace 6 \vert 1 \rbrace \rbrace \! = \! \varnothing, \\
\hat{\mathfrak{J}}_{4}(1) \! := \! \lbrace \operatorname{ind} \lbrace 
i(4)_{1},i(4)_{k_{4}} \vert 1 \rbrace \rbrace \! = \! \lbrace 
\operatorname{ind} \lbrace 3,7 \vert 1 \rbrace \rbrace \! = \! \varnothing,
\end{gather*}
hence
\begin{gather*}
\varkappa_{n1 \tilde{k}_{1}} \! = \! (n \! - \! 1) \gamma_{i(1)_{k_{1}}} \! 
= \! (n \! - \! 1) \gamma_{4} \! = \! 2(n \! - \! 1), \\
\varkappa_{n1 \tilde{k}_{2}} \! = \! (n \! - \! 1) \gamma_{i(2)_{k_{2}}} \! 
= \! (n \! - \! 1) \gamma_{5} \! = \! n \! - \! 1, \\
\varkappa_{n1 \tilde{k}_{3}} \! = \! (n \! - \! 1) \gamma_{i(3)_{k_{3}}} \! 
= \! (n \! - \! 1) \gamma_{6} \! = \! n \! - \! 1, \\
\varkappa^{\infty}_{n1 \tilde{k}_{4}} \! = \! (n \! - \! 1) 
\gamma_{i(4)_{k_{4}}} \! = \! (n \! - \! 1) \gamma_{7} \! = \! 2(n \! - \! 1),
\end{gather*}
that is, as one moves {}from left to right across the repeated pole sequence
\begin{align*}
\mathcal{P}_{n,1} =& \, \lbrace \overset{1}{\underbrace{\alpha_{1},
\alpha_{2},\dotsc,\alpha_{7}}_{7}} \rbrace \cup \dotsb \cup \lbrace 
\overset{n-1}{\underbrace{\alpha_{1},\alpha_{2},\dotsc,\alpha_{7}}_{7}} 
\rbrace \cup \lbrace \overset{n}{\underbrace{\alpha_{1}}_{1}} \rbrace \\
=& \, \lbrace \overset{1}{\underbrace{0,1,\infty,1,\sqrt{2},\pi,\infty}_{7}} 
\rbrace \cup \dotsb \cup \lbrace \overset{n-1}{\underbrace{0,1,\infty,1,
\sqrt{2},\pi,\infty}_{7}} \rbrace \cup \lbrace \overset{n}{\underbrace{0}_{1}} 
\rbrace
\end{align*}
and removes all occurrences of the pole $\alpha_{1} \! = \! 0$, which occurs 
$\varkappa_{n1} \! = \! (n \! - \! 1) \gamma_{1} \! + \! \varrho_{1} \! = \! 
(n \! - \! 1) \! + \! 1 \! = \! n$ times, one is left with the residual pole 
set (via the above induced ordering)
\begin{align*}
\mathcal{P}_{n,1} \setminus \lbrace \underbrace{\alpha_{1},\alpha_{1},\dotsc,
\alpha_{1}}_{\varkappa_{n1}} \rbrace \! =& \, \mathcal{P}_{n,1} \setminus 
\lbrace \underbrace{0,0,\dotsc,0}_{n} \rbrace \! := \bigcup_{q=1}^{3} 
\lbrace \underbrace{\alpha_{i(q)_{k_{q}}},\alpha_{i(q)_{k_{q}}},\dotsc,
\alpha_{i(q)_{k_{q}}}}_{\varkappa_{n1 \tilde{k}_{q}}} \rbrace \cup \lbrace 
\underbrace{\alpha_{i(4)_{k_{4}}},\alpha_{i(4)_{k_{4}}},\dotsc,
\alpha_{i(4)_{k_{4}}}}_{\varkappa^{\infty}_{n1 \tilde{k}_{4}}} \rbrace \\
:=& \, \lbrace \underbrace{\alpha_{i(1)_{k_{1}}},\alpha_{i(1)_{k_{1}}},\dotsc,
\alpha_{i(1)_{k_{1}}}}_{\varkappa_{n1 \tilde{k}_{1}}} \rbrace \cup \lbrace 
\underbrace{\alpha_{i(2)_{k_{2}}},\alpha_{i(2)_{k_{2}}},\dotsc,
\alpha_{i(2)_{k_{2}}}}_{\varkappa_{n1 \tilde{k}_{2}}} \rbrace \\
\cup& \, \lbrace \underbrace{\alpha_{i(3)_{k_{3}}},\alpha_{i(3)_{k_{3}}},
\dotsc,\alpha_{i(3)_{k_{3}}}}_{\varkappa_{n1 \tilde{k}_{3}}} \rbrace \cup 
\lbrace \underbrace{\alpha_{i(4)_{k_{4}}},\alpha_{i(4)_{k_{4}}},\dotsc,
\alpha_{i(4)_{k_{4}}}}_{\varkappa^{\infty}_{n1 \tilde{k}_{4}}} \rbrace \\
=& \, \lbrace \underbrace{1,1,\dotsc,1}_{2(n-1)} \rbrace \cup \lbrace 
\underbrace{\sqrt{2},\sqrt{2},\dotsc,\sqrt{2}}_{n-1} \rbrace \cup 
\lbrace \underbrace{\pi,\pi,\dotsc,\pi}_{n-1} \rbrace \cup \lbrace 
\underbrace{\infty,\infty,\dotsc,\infty}_{2(n-1)} \rbrace,
\end{align*}
where the number of times the pole $\alpha_{i(1)_{k_{1}}} \! = \! 1$ 
$(\neq \! \alpha_{1} \! = \! 0)$ occurs is $\varkappa_{n1 \tilde{k}_{1}} 
\! = \! 2(n \! - \! 1)$, the number of times the pole $\alpha_{i(2)_{k_{2}}} 
\! = \! \sqrt{2}$ $(\neq \! \alpha_{1} \! = \! 0)$ occurs is 
$\varkappa_{n1 \tilde{k}_{2}} \! = \! n \! - \! 1$, the number of times the 
pole $\alpha_{i(3)_{k_{3}}} \! = \! \pi$ $(\neq \! \alpha_{1} \! = \! 0)$ 
occurs is $\varkappa_{n1 \tilde{k}_{3}} \! = \! n \! - \! 1$, and the number 
of times the pole $\alpha_{i(4)_{k_{4}}} \! = \! \infty$ $(\neq \! \alpha_{1} 
\! = \! 0)$ occurs is $\varkappa^{\infty}_{n1 \tilde{k}_{4}} \! = \! 2(n \! - 
\! 1)$.\footnote{Note: for $n \! = \! 1$, since $\varkappa_{11 \tilde{k}_{1}} 
\! = \! \varkappa_{11 \tilde{k}_{2}} \! = \! \varkappa_{11 \tilde{k}_{3}} \! 
= \! \varkappa^{\infty}_{11 \tilde{k}_{4}} \! = \! 0$, one sets, as per the 
convention above, $\lbrace \alpha_{i(q)_{k_{q}}},\alpha_{i(q)_{k_{q}}},
\dotsc,\alpha_{i(q)_{k_{q}}} \rbrace \! := \! \varnothing$, $q \! = \! 
1,2,3,4$, in which case, as $\varkappa_{11} \! = \! \varrho_{1} \! = \! 
1$, $\mathcal{P}_{1,1} \setminus \lbrace \alpha_{1} \rbrace \! = \! 
\mathcal{P}_{1,1} \setminus \lbrace 0 \rbrace \! = \! \varnothing \cup 
\varnothing \cup \varnothing \cup \varnothing \! = \! \varnothing$.} 
In this case, the ordered sum formula reads
\begin{equation*}
\sum_{q=1}^{3} \varkappa_{n1 \tilde{k}_{q}} \! + \! \varkappa^{\infty}_{n1 
\tilde{k}_{4}} \! := \! \varkappa_{n1 \tilde{k}_{1}} \! + \! \varkappa_{n1 
\tilde{k}_{2}} \! + \! \varkappa_{n1 \tilde{k}_{3}} \! + \! 
\varkappa^{\infty}_{n1 \tilde{k}_{4}} \! = \! 2(n \! - \! 1) \! + \! (n \! - 
\! 1) \! + \! (n \! - \! 1) \! + \! 2(n \! - \! 1) \! = \! 6(n \! - \! 1),
\end{equation*}
whence $\sum_{q=1}^{3} \varkappa_{n1 \tilde{k}_{q}} \! + \! 
\varkappa^{\infty}_{n1 \tilde{k}_{4}} \! + \! \varkappa_{n1} \! = \! 
6(n \! - \! 1) \! + \! n \! = \! 7(n \! - \! 1) \! + \! 1$;
\item[(ii)] \fbox{$k \! = \! 2$}
\begin{gather*}
\hat{\mathfrak{J}}_{1}(2) \! := \! \lbrace \operatorname{ind} \lbrace 
i(1)_{k_{1}} \vert 2 \rbrace \rbrace \! = \! \lbrace \operatorname{ind} 
\lbrace 1 \vert 2 \rbrace \rbrace \! = \! \lbrace 1 \rbrace \Rightarrow 
\hat{m}_{1}(2) \! = \! 1, \\
\hat{\mathfrak{J}}_{2}(2) \! := \! \lbrace \operatorname{ind} \lbrace 
i(2)_{k_{2}} \vert 2 \rbrace \rbrace \! = \! \lbrace \operatorname{ind} 
\lbrace 5 \vert 2 \rbrace \rbrace \! = \! \varnothing, \\
\hat{\mathfrak{J}}_{3}(2) \! := \! \lbrace \operatorname{ind} \lbrace 
i(3)_{k_{3}} \vert 2 \rbrace \rbrace \! = \! \lbrace \operatorname{ind} 
\lbrace 6 \vert 2 \rbrace \rbrace \! = \! \varnothing, \\
\hat{\mathfrak{J}}_{4}(2) \! := \! \lbrace \operatorname{ind} \lbrace 
i(4)_{1},i(4)_{k_{4}} \vert 2 \rbrace \rbrace \! = \! \lbrace 
\operatorname{ind} \lbrace 3,7 \vert 2 \rbrace \rbrace \! = \! \varnothing,
\end{gather*}
hence
\begin{gather*}
\varkappa_{n2 \tilde{k}_{1}} \! = \! (n \! - \! 1) \gamma_{\hat{m}_{1}(2)} 
\! + \! \varrho_{\hat{m}_{1}(2)} \! = \! (n \! - \! 1) \gamma_{1} \! + \! 
\varrho_{1} \! = \! (n \! - \! 1) \! + \! 1 \! = \! n, \\
\varkappa_{n2 \tilde{k}_{2}} \! = \! (n \! - \! 1) \gamma_{i(2)_{k_{2}}} \! 
= \! (n \! - \! 1) \gamma_{5} \! = \! n \! - \! 1, \\
\varkappa_{n2 \tilde{k}_{3}} \! = \! (n \! - \! 1) \gamma_{i(3)_{k_{3}}} \! 
= \! (n \! - \! 1) \gamma_{6} \! = \! n \! - \! 1, \\
\varkappa^{\infty}_{n2 \tilde{k}_{4}} \! = \! (n \! - \! 1) 
\gamma_{i(4)_{k_{4}}} \! = \! (n \! - \! 1) \gamma_{7} \! = \! 2(n \! - \! 1),
\end{gather*}
that is, as one moves {}from left to right across the repeated pole sequence
\begin{align*}
\mathcal{P}_{n,2} =& \, \lbrace \overset{1}{\underbrace{\alpha_{1},
\alpha_{2},\dotsc,\alpha_{7}}_{7}} \rbrace \cup \dotsb \cup \lbrace 
\overset{n-1}{\underbrace{\alpha_{1},\alpha_{2},\dotsc,\alpha_{7}}_{7}} 
\rbrace \cup \lbrace \overset{n}{\underbrace{\alpha_{1},\alpha_{2}}_{2}} 
\rbrace \\
=& \, \lbrace \overset{1}{\underbrace{0,1,\infty,1,\sqrt{2},\pi,\infty}_{7}} 
\rbrace \cup \dotsb \cup \lbrace \overset{n-1}{\underbrace{0,1,\infty,1,
\sqrt{2},\pi,\infty}_{7}} \rbrace \cup \lbrace 
\overset{n}{\underbrace{0,1}_{2}} \rbrace
\end{align*}
and removes all occurrences of the pole $\alpha_{2} \! = \! 1$, which occurs 
$\varkappa_{n2} \! = \! (n \! - \! 1) \gamma_{2} \! + \! \varrho_{2} \! = \! 
2(n \! - \! 1) \! + \! 1 \! = \! 2n \! - \! 1$ times, one is left with the 
residual pole set (via the above induced ordering)
\begin{align*}
\mathcal{P}_{n,2} \setminus \lbrace \underbrace{\alpha_{2},\alpha_{2},\dotsc,
\alpha_{2}}_{\varkappa_{n2}} \rbrace \! =& \, \mathcal{P}_{n,2} \setminus 
\lbrace \underbrace{1,1,\dotsc,1}_{2n-1} \rbrace \! := \bigcup_{q=1}^{3} 
\lbrace \underbrace{\alpha_{i(q)_{k_{q}}},\alpha_{i(q)_{k_{q}}},\dotsc,
\alpha_{i(q)_{k_{q}}}}_{\varkappa_{n2 \tilde{k}_{q}}} \rbrace \cup \lbrace 
\underbrace{\alpha_{i(4)_{k_{4}}},\alpha_{i(4)_{k_{4}}},\dotsc,
\alpha_{i(4)_{k_{4}}}}_{\varkappa^{\infty}_{n2 \tilde{k}_{4}}} \rbrace \\
:=& \, \lbrace \underbrace{\alpha_{i(1)_{k_{1}}},\alpha_{i(1)_{k_{1}}},
\dotsc,\alpha_{i(1)_{k_{1}}}}_{\varkappa_{n2 \tilde{k}_{1}}} \rbrace \cup 
\lbrace \underbrace{\alpha_{i(2)_{k_{2}}},\alpha_{i(2)_{k_{2}}},\dotsc,
\alpha_{i(2)_{k_{2}}}}_{\varkappa_{n2 \tilde{k}_{2}}} \rbrace \\
\cup& \, \lbrace \underbrace{\alpha_{i(3)_{k_{3}}},\alpha_{i(3)_{k_{3}}},
\dotsc,\alpha_{i(3)_{k_{3}}}}_{\varkappa_{n2 \tilde{k}_{3}}} \rbrace \cup 
\lbrace \underbrace{\alpha_{i(4)_{k_{4}}},\alpha_{i(4)_{k_{4}}},\dotsc,
\alpha_{i(4)_{k_{4}}}}_{\varkappa^{\infty}_{n2 \tilde{k}_{4}}} \rbrace \\
=& \, \lbrace \underbrace{0,0,\dotsc,0}_{n} \rbrace \cup \lbrace 
\underbrace{\sqrt{2},\sqrt{2},\dotsc,\sqrt{2}}_{n-1} \rbrace \cup 
\lbrace \underbrace{\pi,\pi,\dotsc,\pi}_{n-1} \rbrace \cup \lbrace 
\underbrace{\infty,\infty,\dotsc,\infty}_{2(n-1)} \rbrace,
\end{align*}
where the number of times the pole $\alpha_{i(1)_{k_{1}}} \! = \! 0$ 
$(\neq \! \alpha_{2} \! = \! 1)$ occurs is $\varkappa_{n2 \tilde{k}_{1}} \! = 
\! n$, the number of times the pole $\alpha_{i(2)_{k_{2}}} \! = \! \sqrt{2}$ 
$(\neq \! \alpha_{2} \! = \! 1)$ occurs is $\varkappa_{n2 \tilde{k}_{2}} \! 
= \! n \! - \! 1$, the number of times the pole $\alpha_{i(3)_{k_{3}}} \! 
= \! \pi$ $(\neq \! \alpha_{2} \! = \! 1)$ occurs is $\varkappa_{n2 
\tilde{k}_{3}} \! = \! n \! - \! 1$, and the number of times the pole 
$\alpha_{i(4)_{k_{4}}} \! = \! \infty$ $(\neq \! \alpha_{2} \! = \! 1)$ 
occurs is $\varkappa^{\infty}_{n2 \widetilde{k}_{4}} \! = \! 2(n \! - 
\! 1)$.\footnote{Note: for $n \! = \! 1$, since $\varkappa_{12 
\tilde{k}_{2}} \! = \! \varkappa_{12 \tilde{k}_{3}} \! = \! 
\varkappa^{\infty}_{12 \tilde{k}_{4}} \! = \! 0$, one sets, as per the 
convention above, $\lbrace \alpha_{i(q)_{k_{q}}},\alpha_{i(q)_{k_{q}}},
\dotsc,\alpha_{i(q)_{k_{q}}} \rbrace \! := \! \varnothing$, $q \! = \! 
2,3,4$, in which case, as $\varkappa_{12 \tilde{k}_{1}} \! = \! 1$ and 
$\varkappa_{12} \! = \! \varrho_{2} \! = \! 1$, $\mathcal{P}_{1,2} 
\setminus \lbrace \alpha_{2} \rbrace \! = \! \mathcal{P}_{1,2} 
\setminus \lbrace 1 \rbrace \! = \! \lbrace 0 \rbrace \cup \varnothing 
\cup \varnothing \cup \varnothing \! = \! \lbrace 0 \rbrace$.} In this 
case, the ordered sum formula reads
\begin{equation*}
\sum_{q=1}^{3} \varkappa_{n2 \tilde{k}_{q}} \! + \! \varkappa^{\infty}_{n2 
\tilde{k}_{4}} \! := \! \varkappa_{n2 \tilde{k}_{1}} \! + \! \varkappa_{n2 
\tilde{k}_{2}} \! + \! \varkappa_{n2 \tilde{k}_{3}} \! + \! 
\varkappa^{\infty}_{n2 \tilde{k}_{4}} \! = \! n \! + \! (n \! - \! 1) \! + \! 
(n \! - \! 1) \! + \! 2(n \! - \! 1) \! = \! 5n \! - \! 4,
\end{equation*}
whence $\sum_{q=1}^{3} \varkappa_{n2 \tilde{k}_{q}} \! + \! 
\varkappa^{\infty}_{n2 \tilde{k}_{4}} \! + \! \varkappa_{n2} \! = \! 
(5n \! - \! 4) \! + \! (2n \! - \! 1) \! = \! 7(n \! - \! 1) \! + \! 2$;
\item[(iii)] \fbox{$k \! = \! 4$}
\begin{gather*}
\hat{\mathfrak{J}}_{1}(4) \! := \! \lbrace \operatorname{ind} \lbrace 
i(1)_{k_{1}} \vert 4 \rbrace \rbrace \! = \! \lbrace \operatorname{ind} 
\lbrace 1 \vert 4 \rbrace \rbrace \! = \! \lbrace 1 \rbrace \Rightarrow 
\hat{m}_{1}(4) \! = \! 1, \\
\hat{\mathfrak{J}}_{2}(4) \! := \! \lbrace \operatorname{ind} \lbrace 
i(2)_{k_{2}} \vert 4 \rbrace \rbrace \! = \! \lbrace \operatorname{ind} 
\lbrace 5 \vert 4 \rbrace \rbrace \! = \! \varnothing, \\
\hat{\mathfrak{J}}_{3}(4) \! := \! \lbrace \operatorname{ind} \lbrace 
i(3)_{k_{3}} \vert 4 \rbrace \rbrace \! = \! \lbrace \operatorname{ind} 
\lbrace 6 \vert 4 \rbrace \rbrace \! = \! \varnothing, \\
\hat{\mathfrak{J}}_{4}(4) \! := \! \lbrace \operatorname{ind} \lbrace 
i(4)_{1},i(4)_{k_{4}} \vert 4 \rbrace \rbrace \! = \! \lbrace 
\operatorname{ind} \lbrace 3,7 \vert 4 \rbrace \rbrace \! = \! \lbrace 3 
\rbrace \Rightarrow \hat{m}_{4}(4) \! = \! 3,
\end{gather*}
hence
\begin{gather*}
\varkappa_{n4 \tilde{k}_{1}} \! = \! (n \! - \! 1) \gamma_{\hat{m}_{1}(4)} 
\! + \! \varrho_{\hat{m}_{1}(4)} \! = \! (n \! - \! 1) \gamma_{1} \! + \! 
\varrho_{1} \! = \! (n \! - \! 1) \! + \! 1 \! = \! n, \\
\varkappa_{n4 \tilde{k}_{2}} \! = \! (n \! - \! 1) \gamma_{i(2)_{k_{2}}} \! 
= \! (n \! - \! 1) \gamma_{5} \! = \! n \! - \! 1, \\
\varkappa_{n4 \tilde{k}_{3}} \! = \! (n \! - \! 1) \gamma_{i(3)_{k_{3}}} \! 
= \! (n \! - \! 1) \gamma_{6} \! = \! n \! - \! 1, \\
\varkappa^{\infty}_{n4 \tilde{k}_{4}} \! = \! (n \! - \! 1) 
\gamma_{\hat{m}_{4}(4)} \! + \! \varrho_{\hat{m}_{4}(4)} \! = \! 
(n \! - \! 1) \gamma_{3} \! + \! \varrho_{3} \! = \! 2(n \! - \! 1) 
\! + \! 1 \! = \! 2n \! - \! 1,
\end{gather*}
that is, as one moves {}from left to right across the repeated pole sequence
\begin{align*}
\mathcal{P}_{n,4} =& \, \lbrace \overset{1}{\underbrace{\alpha_{1},
\alpha_{2},\dotsc,\alpha_{7}}_{7}} \rbrace \cup \dotsb \cup \lbrace 
\overset{n-1}{\underbrace{\alpha_{1},\alpha_{2},\dotsc,\alpha_{7}}_{7}} 
\rbrace \cup \lbrace \overset{n}{\underbrace{\alpha_{1},\alpha_{2},\alpha_{3},
\alpha_{4}}_{4}} \rbrace \\
=& \, \lbrace \overset{1}{\underbrace{0,1,\infty,1,\sqrt{2},\pi,\infty}_{7}} 
\rbrace \cup \dotsb \cup \lbrace \overset{n-1}{\underbrace{0,1,\infty,1,
\sqrt{2},\pi,\infty}_{7}} \rbrace \cup \lbrace \overset{n}{\underbrace{0,1,
\infty,1}_{4}} \rbrace
\end{align*}
and removes all occurrences of the pole $\alpha_{4} \! = \! 1$, which occurs 
$\varkappa_{n4} \! = \! (n \! - \! 1) \gamma_{4} \! + \! \varrho_{4} \! = \! 
2(n \! - \! 1) \! + \! 2 \! = \! 2n$ times, one is left with the residual pole 
set (via the above induced ordering)
\begin{align*}
\mathcal{P}_{n,4} \setminus \lbrace \underbrace{\alpha_{4},\alpha_{4},
\dotsc,\alpha_{4}}_{\varkappa_{n4}} \rbrace \! =& \, \mathcal{P}_{n,4} 
\setminus \lbrace \underbrace{1,1,\dotsc,1}_{2n} \rbrace \! := 
\bigcup_{q=1}^{3} \lbrace \underbrace{\alpha_{i(q)_{k_{q}}},
\alpha_{i(q)_{k_{q}}},\dotsc,\alpha_{i(q)_{k_{q}}}}_{\varkappa_{n4 
\tilde{k}_{q}}} \rbrace \cup \lbrace \underbrace{\alpha_{i(4)_{k_{4}}},
\alpha_{i(4)_{k_{4}}},\dotsc,\alpha_{i(4)_{k_{4}}}}_{\varkappa^{\infty}_{n4 
\tilde{k}_{4}}} \rbrace \\
:=& \, \lbrace \underbrace{\alpha_{i(1)_{k_{1}}},\alpha_{i(1)_{k_{1}}},\dotsc,
\alpha_{i(1)_{k_{1}}}}_{\varkappa_{n4 \tilde{k}_{1}}} \rbrace \cup \lbrace 
\underbrace{\alpha_{i(2)_{k_{2}}},\alpha_{i(2)_{k_{2}}},\dotsc,
\alpha_{i(2)_{k_{2}}}}_{\varkappa_{n4 \tilde{k}_{2}}} \rbrace \\
\cup& \, \lbrace \underbrace{\alpha_{i(3)_{k_{3}}},\alpha_{i(3)_{k_{3}}},
\dotsc,\alpha_{i(3)_{k_{3}}}}_{\varkappa_{n4 \tilde{k}_{3}}} \rbrace \cup 
\lbrace \underbrace{\alpha_{i(4)_{k_{4}}},\alpha_{i(4)_{k_{4}}},\dotsc,
\alpha_{i(4)_{k_{4}}}}_{\varkappa^{\infty}_{n4 \tilde{k}_{4}}} \rbrace \\
=& \, \lbrace \underbrace{0,0,\dotsc,0}_{n} \rbrace \cup \lbrace 
\underbrace{\sqrt{2},\sqrt{2},\dotsc,\sqrt{2}}_{n-1} \rbrace \cup 
\lbrace \underbrace{\pi,\pi,\dotsc,\pi}_{n-1} \rbrace \cup \lbrace 
\underbrace{\infty,\infty,\dotsc,\infty}_{2n-1} \rbrace,
\end{align*}
where the number of times the pole $\alpha_{i(1)_{k_{1}}} \! = \! 0$ 
$(\neq \! \alpha_{4} \! = \! 1)$ occurs is $\varkappa_{n4 \tilde{k}_{1}} \! = 
\! n$, the number of times the pole $\alpha_{i(2)_{k_{2}}} \! = \! \sqrt{2}$ 
$(\neq \! \alpha_{4} \! = \! 1)$ occurs is $\varkappa_{n4 \tilde{k}_{2}} \! = 
\! n \! - \! 1$, the number of times the pole $\alpha_{i(3)_{k_{3}}} \! = \! 
\pi$ $(\neq \! \alpha_{4} \! = \! 1)$ occurs is $\varkappa_{n4 \tilde{k}_{3}} 
\! = \! n \! - \! 1$, and the number of times the pole $\alpha_{i(4)_{k_{4}}} 
\! = \! \infty$ $(\neq \! \alpha_{4} \! = \! 1)$ occurs is 
$\varkappa^{\infty}_{n4 \tilde{k}_{4}} \! = \! 2n \! - \! 1$.\footnote{Note: 
for $n \! = \! 1$, since $\varkappa_{14 \tilde{k}_{2}} \! = \! \varkappa_{14 
\tilde{k}_{3}} \! = \! 0$, one sets, as per the convention above, $\lbrace 
\alpha_{i(q)_{k_{q}}},\alpha_{i(q)_{k_{q}}},\dotsc,\alpha_{i(q)_{k_{q}}} 
\rbrace \! := \! \varnothing$, $q \! = \! 2,3$, in which case, as 
$\varkappa_{14 \tilde{k}_{1}} \! = \! \varkappa^{\infty}_{14 \tilde{k}_{4}} 
\! = \! 1$ and $\varkappa_{14} \! = \! \varrho_{4} \! = \! 2$, 
$\mathcal{P}_{1,4} \setminus \lbrace \alpha_{4},\alpha_{4} \rbrace \! = \! 
\mathcal{P}_{1,4} \setminus \lbrace 1,1 \rbrace \! = \! \lbrace 0 \rbrace \cup 
\varnothing \cup \varnothing \cup \lbrace \infty \rbrace \! = \! \lbrace 0 
\rbrace \cup \lbrace \infty \rbrace$.} In this case, the ordered sum formula 
reads
\begin{equation*}
\sum_{q=1}^{3} \varkappa_{n4 \tilde{k}_{q}} \! + \! \varkappa^{\infty}_{n4 
\tilde{k}_{4}} \! := \! \varkappa_{n4 \tilde{k}_{1}} \! + \! \varkappa_{n4 
\tilde{k}_{2}} \! + \! \varkappa_{n4 \tilde{k}_{3}} \! + \! 
\varkappa^{\infty}_{n4 \tilde{k}_{4}} \! = \! n \! + \! (n \! - \! 1) \! + \! 
(n \! - \! 1) \! + \! (2n \! - \! 1) \! = \! 5n \! - \! 3,
\end{equation*}
whence $\sum_{q=1}^{3} \varkappa_{n4 \tilde{k}_{q}} \! + \! 
\varkappa^{\infty}_{n4 \tilde{k}_{4}} \! + \! \varkappa_{n4} \! = \! 
(5n \! - \! 3) \! + \! 2n \! = \! 7(n \! - \! 1) \! + \! 4$;
\item[(iv)] \fbox{$k \! = \! 5$}
\begin{gather*}
\hat{\mathfrak{J}}_{1}(5) \! := \! \lbrace \operatorname{ind} \lbrace 
i(1)_{k_{1}} \vert 5 \rbrace \rbrace \! = \! \lbrace \operatorname{ind} 
\lbrace 1 \vert 5 \rbrace \rbrace \! = \! \lbrace 1 \rbrace \Rightarrow 
\hat{m}_{1}(5) \! = \! 1, \\
\hat{\mathfrak{J}}_{2}(5) \! := \! \lbrace \operatorname{ind} \lbrace 
i(2)_{1},i(2)_{k_{2}} \vert 5 \rbrace \rbrace \! = \! \lbrace 
\operatorname{ind} \lbrace 2,4 \vert 5 \rbrace \rbrace \! = \! \lbrace 
4 \rbrace \Rightarrow \hat{m}_{2}(5) \! = \! 4, \\
\hat{\mathfrak{J}}_{3}(5) \! := \! \lbrace \operatorname{ind} \lbrace 
i(3)_{k_{3}} \vert 5 \rbrace \rbrace \! = \! \lbrace \operatorname{ind} 
\lbrace 6 \vert 5 \rbrace \rbrace \! = \! \varnothing, \\
\hat{\mathfrak{J}}_{4}(5) \! := \! \lbrace \operatorname{ind} \lbrace 
i(4)_{1},i(4)_{k_{4}} \vert 5 \rbrace \rbrace \! = \! \lbrace 
\operatorname{ind} \lbrace 3,7 \vert 5 \rbrace \rbrace \! = \! \lbrace 
3 \rbrace \Rightarrow \hat{m}_{4}(5) \! = \! 3,
\end{gather*}
hence
\begin{gather*}
\varkappa_{n5 \tilde{k}_{1}} \! = \! (n \! - \! 1) \gamma_{\hat{m}_{1}(5)} 
\! + \! \varrho_{\hat{m}_{1}(5)} \! = \! (n \! - \! 1) \gamma_{1} \! + \! 
\varrho_{1} \! = \! (n \! - \! 1) \! + \! 1 \! = \! n, \\
\varkappa_{n5 \tilde{k}_{2}} \! = \! (n \! - \! 1) \gamma_{\hat{m}_{2}(5)} 
\! + \! \varrho_{\hat{m}_{2}(5)} \! = \! (n \! - \! 1) \gamma_{4} \! + \! 
\varrho_{4} \! = \! 2(n \! - \! 1) \! + \! 2 \! = \! 2n, \\
\varkappa_{n5 \tilde{k}_{3}} \! = \! (n \! - \! 1) \gamma_{i(3)_{k_{3}}} 
\! = \! (n \! - \! 1) \gamma_{6} \! = \! n \! - \! 1, \\
\varkappa^{\infty}_{n5 \tilde{k}_{4}} \! = \! (n \! - \! 1) 
\gamma_{\hat{m}_{4}(5)} \! + \! \varrho_{\hat{m}_{4}(5)} \! = \! 
(n \! - \! 1) \gamma_{3} \! + \! \varrho_{3} \! = \! 2(n \! - \! 1) 
\! + \! 1 \! = \! 2n \! - \! 1,
\end{gather*}
that is, as one moves {}from left to right across the repeated pole sequence
\begin{align*}
\mathcal{P}_{n,5} =& \, \lbrace \overset{1}{\underbrace{\alpha_{1},
\alpha_{2},\dotsc,\alpha_{7}}_{7}} \rbrace \cup \dotsb \cup \lbrace 
\overset{n-1}{\underbrace{\alpha_{1},\alpha_{2},\dotsc,\alpha_{7}}_{7}} 
\rbrace \cup \lbrace \overset{n}{\underbrace{\alpha_{1},\alpha_{2},
\alpha_{3},\alpha_{4},\alpha_{5}}_{5}} \rbrace \\
=& \, \lbrace \overset{1}{\underbrace{0,1,\infty,1,\sqrt{2},\pi,\infty}_{7}} 
\rbrace \cup \dotsb \cup \lbrace \overset{n-1}{\underbrace{0,1,\infty,1,
\sqrt{2},\pi,\infty}_{7}} \rbrace \cup \lbrace \overset{n}{\underbrace{0,1,
\infty,1,\sqrt{2}}_{5}} \rbrace
\end{align*}
and removes all occurrences of the pole $\alpha_{5} \! = \! \sqrt{2}$, which 
occurs $\varkappa_{n5} \! = \! (n \! - \! 1) \gamma_{5} \! + \! \varrho_{5} \! 
= \! (n \! - \! 1) \! + \! 1 \! = \! n$ times, one is left with the residual 
pole set (via the above induced ordering)
\begin{align*}
\mathcal{P}_{n,5} \setminus \lbrace \underbrace{\alpha_{5},\alpha_{5},\dotsc,
\alpha_{5}}_{\varkappa_{n5}} \rbrace \! =& \, \mathcal{P}_{n,5} \setminus 
\lbrace \underbrace{\sqrt{2},\sqrt{2},\dotsc,\sqrt{2}}_{n} \rbrace \! := 
\bigcup_{q=1}^{3} \lbrace \underbrace{\alpha_{i(q)_{k_{q}}},
\alpha_{i(q)_{k_{q}}},\dotsc,\alpha_{i(q)_{k_{q}}}}_{\varkappa_{n5 
\tilde{k}_{q}}} \rbrace \cup \lbrace \underbrace{\alpha_{i(4)_{k_{4}}},
\alpha_{i(4)_{k_{4}}},\dotsc,\alpha_{i(4)_{k_{4}}}}_{\varkappa^{\infty}_{n5 
\tilde{k}_{4}}} \rbrace \\
:=& \, \lbrace \underbrace{\alpha_{i(1)_{k_{1}}},\alpha_{i(1)_{k_{1}}},
\dotsc,\alpha_{i(1)_{k_{1}}}}_{\varkappa_{n5 \tilde{k}_{1}}} \rbrace \cup 
\lbrace \underbrace{\alpha_{i(2)_{k_{2}}},\alpha_{i(2)_{k_{2}}},\dotsc,
\alpha_{i(2)_{k_{2}}}}_{\varkappa_{n5 \tilde{k}_{2}}} \rbrace \\
\cup& \, \lbrace \underbrace{\alpha_{i(3)_{k_{3}}},\alpha_{i(3)_{k_{3}}},
\dotsc,\alpha_{i(3)_{k_{3}}}}_{\varkappa_{n5 \tilde{k}_{3}}} \rbrace \cup 
\lbrace \underbrace{\alpha_{i(4)_{k_{4}}},\alpha_{i(4)_{k_{4}}},\dotsc,
\alpha_{i(4)_{k_{4}}}}_{\varkappa^{\infty}_{n5 \tilde{k}_{4}}} \rbrace \\
=& \, \lbrace \underbrace{0,0,\dotsc,0}_{n} \rbrace \cup \lbrace 
\underbrace{1,1,\dotsc,1}_{2n} \rbrace \cup \lbrace \underbrace{\pi,\pi,
\dotsc,\pi}_{n-1} \rbrace \cup \lbrace \underbrace{\infty,\infty,\dotsc,
\infty}_{2n-1} \rbrace,
\end{align*}
where the number of times the pole $\alpha_{i(1)_{k_{1}}} \! = \! 0$ 
$(\neq \! \alpha_{5} \! = \! \sqrt{2})$ occurs is $\varkappa_{n5 
\tilde{k}_{1}} \! = \! n$, the number of times the pole $\alpha_{i(2)_{k_{2}}} 
\! = \! 1$ $(\neq \! \alpha_{5} \! = \! \sqrt{2})$ occurs is $\varkappa_{n5 
\tilde{k}_{2}} \! = \! 2n$, the number of times the pole 
$\alpha_{i(3)_{k_{3}}} \! = \! \pi$ $(\neq \! \alpha_{5} \! = \! \sqrt{2})$ 
occurs is $\varkappa_{n5 \tilde{k}_{3}} \! = \! n \! - \! 1$, and the number 
of times the pole $\alpha_{i(4)_{k_{4}}} \! = \! \infty$ $(\neq \! \alpha_{5} 
\! = \! \sqrt{2})$ occurs is $\varkappa^{\infty}_{n5 \tilde{k}_{4}} \! = \! 2n 
\! - \! 1$.\footnote{Note: for $n \! = \! 1$, since $\varkappa_{15 \tilde{k}_{3}} 
\! = \! 0$, one sets, as per the convention above, $\lbrace \alpha_{i(q)_{k_{q}}},
\alpha_{i(q)_{k_{q}}},\dotsc,\alpha_{i(q)_{k_{q}}} \rbrace \! := \! 
\varnothing$, $q \! = \! 3$, in which case, as $\varkappa_{15 \tilde{k}_{1}} 
\! = \! 1$, $\varkappa_{15 \tilde{k}_{2}} \! = \! 2$, $\varkappa^{\infty}_{15 
\tilde{k}_{4}} \! = \! 1$, and $\varkappa_{15} \! = \! \varrho_{5} \! = 
\! 1$, $\mathcal{P}_{1,5} \setminus \lbrace \alpha_{5} \rbrace \! = \! 
\mathcal{P}_{1,5} \setminus \lbrace \sqrt{2} \rbrace \! = \! \lbrace 0 \rbrace 
\cup \lbrace 1,1 \rbrace \cup \varnothing \cup \lbrace \infty \rbrace \! = 
\! \lbrace 0 \rbrace \cup \lbrace 1,1 \rbrace \cup \lbrace \infty \rbrace$.} 
In this case, the ordered sum formula reads
\begin{equation*}
\sum_{q=1}^{3} \varkappa_{n5 \tilde{k}_{q}} \! + \! \varkappa^{\infty}_{n5 
\tilde{k}_{4}} \! := \! \varkappa_{n5 \tilde{k}_{1}} \! + \! \varkappa_{n5 
\tilde{k}_{2}} \! + \! \varkappa_{n5 \tilde{k}_{3}} \! + \! 
\varkappa^{\infty}_{n5 \tilde{k}_{4}} \! = \! n \! + \! 2n \! + \! 
(n \! - \! 1) \! + \! (2n \! - \! 1) \! = \! 6n \! - \! 2,
\end{equation*}
whence $\sum_{q=1}^{3} \varkappa_{n5 \tilde{k}_{q}} \! + \!
\varkappa^{\infty}_{n5 \tilde{k}_{4}} \! + \! \varkappa_{n5} \! = \! 
(6n \! - \! 2) \! + \! n \! = \! 7(n \! - \! 1) \! + \! 5$;
\item[(v)] \fbox{$k \! = \! 6$}
\begin{gather*}
\hat{\mathfrak{J}}_{1}(6) \! := \! \lbrace \operatorname{ind} \lbrace 
i(1)_{k_{1}} \vert 6 \rbrace \rbrace \! = \! \lbrace \operatorname{ind} 
\lbrace 1 \vert 6 \rbrace \rbrace \! = \! \lbrace 1 \rbrace \Rightarrow 
\hat{m}_{1}(6) \! = \! 1, \\
\hat{\mathfrak{J}}_{2}(6) \! := \! \lbrace \operatorname{ind} \lbrace 
i(2)_{1},i(2)_{k_{2}} \vert 6 \rbrace \rbrace \! = \! \lbrace 
\operatorname{ind} \lbrace 2,4 \vert 6 \rbrace \rbrace \! = \! \lbrace 4 
\rbrace \Rightarrow \hat{m}_{2}(6) \! = \! 4, \\
\hat{\mathfrak{J}}_{3}(6) \! := \! \lbrace \operatorname{ind} \lbrace 
i(3)_{k_{3}} \vert 6 \rbrace \rbrace \! = \! \lbrace \operatorname{ind} 
\lbrace 5 \vert 6 \rbrace \rbrace \! = \! \lbrace 5 \rbrace \Rightarrow 
\hat{m}_{3}(6) \! = \! 5, \\
\hat{\mathfrak{J}}_{4}(6) \! := \! \lbrace \operatorname{ind} \lbrace 
i(4)_{1},i(4)_{k_{4}} \vert 6 \rbrace \rbrace \! = \! \lbrace 
\operatorname{ind} \lbrace 3,7 \vert 6 \rbrace \rbrace \! = \! \lbrace 3 
\rbrace \Rightarrow \hat{m}_{4}(6) \! = \! 3,
\end{gather*}
hence
\begin{gather*}
\varkappa_{n6 \tilde{k}_{1}} \! = \! (n \! - \! 1) \gamma_{\hat{m}_{1}(6)} 
\! + \! \varrho_{\hat{m}_{1}(6)} \! = \! (n \! - \! 1) \gamma_{1} \! + \! 
\varrho_{1} \! = \! (n \! - \! 1) \! + \! 1 \! = \! n, \\
\varkappa_{n6 \tilde{k}_{2}} \! = \! (n \! - \! 1) \gamma_{\hat{m}_{2}(6)} 
\! + \! \varrho_{\hat{m}_{2}(6)} \! = \! (n \! - \! 1) \gamma_{4} \! + \! 
\varrho_{4} \! = \! 2(n \! - \! 1) \! + \! 2 \! = \! 2n, \\
\varkappa_{n6 \tilde{k}_{3}} \! = \! (n \! - \! 1) \gamma_{\hat{m}_{3}(6)} 
\! + \! \varrho_{\hat{m}_{3}(6)} \! = \! (n \! - \! 1) \gamma_{5} \! + \! 
\varrho_{5} \! = \! (n \! - \! 1) \! + \! 1 \! = \! n, \\
\varkappa^{\infty}_{n6 \tilde{k}_{4}} \! = \! (n \! - \! 1) \gamma_{
\hat{m}_{4}(6)} \! + \! \varrho_{\hat{m}_{4}(6)} \! = \! (n \! - \! 1) 
\gamma_{3} \! + \! \varrho_{3} \! = \! 2(n \! - \! 1) \! + \! 1 \! = \! 
2n \! - \! 1,
\end{gather*}
that is, as one moves {}from left to right across the repeated pole sequence
\begin{align*}
\mathcal{P}_{n,6} =& \, \lbrace \overset{1}{\underbrace{\alpha_{1},
\alpha_{2},\dotsc,\alpha_{7}}_{7}} \rbrace \cup \dotsb \cup \lbrace 
\overset{n-1}{\underbrace{\alpha_{1},\alpha_{2},\dotsc,\alpha_{7}}_{7}} 
\rbrace \cup \lbrace \overset{n}{\underbrace{\alpha_{1},\alpha_{2},
\alpha_{3},\alpha_{4},\alpha_{5},\alpha_{6}}_{6}} \rbrace \\
=& \, \lbrace \overset{1}{\underbrace{0,1,\infty,1,\sqrt{2},\pi,\infty}_{7}} 
\rbrace \cup \dotsb \cup \lbrace \overset{n-1}{\underbrace{0,1,\infty,1,
\sqrt{2},\pi,\infty}_{7}} \rbrace \cup \lbrace \overset{n}{\underbrace{0,1,
\infty,1,\! \sqrt{2},\pi}_{6}} \rbrace
\end{align*}
and removes all occurrences of the pole $\alpha_{6} \! = \! \pi$, which occurs 
$\varkappa_{n6} \! = \! (n \! - \! 1) \gamma_{6} \! + \! \varrho_{6} \! = \! 
(n \! - \! 1) \! + \! 1 \! = \! n$ times, one is left with the residual pole 
set (via the above induced ordering)
\begin{align*}
\mathcal{P}_{n,6} \setminus \lbrace \underbrace{\alpha_{6},\alpha_{6},\dotsc,
\alpha_{6}}_{\varkappa_{n6}} \rbrace \! =& \, \mathcal{P}_{n,6} \setminus 
\lbrace \underbrace{\pi,\pi,\dotsc,\pi}_{n} \rbrace \! := \bigcup_{q=1}^{3} 
\lbrace \underbrace{\alpha_{i(q)_{k_{q}}},\alpha_{i(q)_{k_{q}}},\dotsc,
\alpha_{i(q)_{k_{q}}}}_{\varkappa_{n6 \tilde{k}_{q}}} \rbrace \cup 
\lbrace \underbrace{\alpha_{i(4)_{k_{4}}},\alpha_{i(4)_{k_{4}}},\dotsc,
\alpha_{i(4)_{k_{4}}}}_{\varkappa^{\infty}_{n6 \tilde{k}_{4}}} \rbrace \\
:=& \, \lbrace \underbrace{\alpha_{i(1)_{k_{1}}},\alpha_{i(1)_{k_{1}}},
\dotsc,\alpha_{i(1)_{k_{1}}}}_{\varkappa_{n6 \tilde{k}_{1}}} \rbrace \cup 
\lbrace \underbrace{\alpha_{i(2)_{k_{2}}},\alpha_{i(2)_{k_{2}}},\dotsc,
\alpha_{i(2)_{k_{2}}}}_{\varkappa_{n6 \tilde{k}_{2}}} \rbrace \\
\cup& \, \lbrace \underbrace{\alpha_{i(3)_{k_{3}}},\alpha_{i(3)_{k_{3}}},
\dotsc,\alpha_{i(3)_{k_{3}}}}_{\varkappa_{n6 \tilde{k}_{3}}} \rbrace \cup 
\lbrace \underbrace{\alpha_{i(4)_{k_{4}}},\alpha_{i(4)_{k_{4}}},\dotsc,
\alpha_{i(4)_{k_{4}}}}_{\varkappa^{\infty}_{n6 \tilde{k}_{4}}} \rbrace \\
=& \, \lbrace \underbrace{0,0,\dotsc,0}_{n} \rbrace \cup \lbrace 
\underbrace{1,1,\dotsc,1}_{2n} \rbrace \cup \lbrace \underbrace{\sqrt{2},
\sqrt{2},\dotsc,\sqrt{2}}_{n} \rbrace \cup \lbrace \underbrace{\infty,\infty,
\dotsc,\infty}_{2n-1} \rbrace,
\end{align*}
where the number of times the pole $\alpha_{i(1)_{k_{1}}} \! = \! 0$ 
$(\neq \! \alpha_{6} \! = \! \pi)$ occurs is $\varkappa_{n6 \tilde{k}_{1}} 
\! = \! n$, the number of times the pole $\alpha_{i(2)_{k_{2}}} \! = \! 1$ 
$(\neq \! \alpha_{6} \! = \! \pi)$ occurs is $\varkappa_{n6 \tilde{k}_{2}} 
\! = \! 2n$, the number of times the pole $\alpha_{i(3)_{k_{3}}} \! = \! 
\sqrt{2}$ $(\neq \! \alpha_{6} \! = \! \pi)$ occurs is $\varkappa_{n6 
\tilde{k}_{3}} \! = \! n$, and the number of times the pole 
$\alpha_{i(4)_{k_{4}}} \! = \! \infty$ $(\neq \! \alpha_{6} \! = \! \pi)$ 
occurs is $\varkappa^{\infty}_{n6 \tilde{k}_{4}} \! = \! 2n \! - \! 
1$.\footnote{Note: for $n \! = \! 1$, as $\varkappa_{16 \tilde{k}_{1}} \! = 
\! 1$, $\varkappa_{16 \tilde{k}_{2}} \! = \! 2$, $\varkappa_{16 \tilde{k}_{3}} 
\! = \! 1$, $\varkappa^{\infty}_{16 \widetilde{k}_{4}} \! = \! 1$, and 
$\varkappa_{16} \! = \! \varrho_{6} \! = \! 1$, it follows that 
$\mathcal{P}_{1,6} \setminus \lbrace \alpha_{6} \rbrace \! = \! 
\mathcal{P}_{1,6} \setminus \lbrace \pi \rbrace \! = \! \lbrace 0 \rbrace 
\cup \lbrace 1,1 \rbrace \cup \lbrace \! \sqrt{2} \rbrace \cup \lbrace 
\infty \rbrace$.} In this case, the ordered sum formula reads
\begin{equation*}
\sum_{q=1}^{3} \varkappa_{n6 \tilde{k}_{q}} \! + \! \varkappa^{\infty}_{n6 
\tilde{k}_{4}} \! := \! \varkappa_{n6 \tilde{k}_{1}} \! + \! \varkappa_{n6 
\tilde{k}_{2}} \! + \! \varkappa_{n6 \tilde{k}_{3}} \! + \! 
\varkappa^{\infty}_{n6 \tilde{k}_{4}} \! = \! n \! + \! 2n \! + \! n \! + \! 
(2n \! - \! 1) \! = \! 6n \! - \! 1,
\end{equation*}
whence $\sum_{q=1}^{3} \varkappa_{n6 \tilde{k}_{q}} \! + \! 
\varkappa^{\infty}_{n6 \tilde{k}_{4}} \! + \! \varkappa_{n6} \! = \! 
(6n \! - \! 1) \! + \! n \! = \! 7(n \! - \! 1) \! + \! 6$.
\end{enumerate}
This concludes the example.

For simplicity of notation, set, hereafter,
\begin{equation*}
\alpha_{i(q)_{k_{q}}} \! := \! \alpha_{p_{q}} \, \, (\neq \! \infty), 
\quad q \! = \! 1,2,\dotsc,\mathfrak{s} \! - \! 2, \quad \quad 
\alpha_{i(\mathfrak{s}-1)_{k_{\mathfrak{s}-1}}} \! := \! 
\alpha_{p_{\mathfrak{s}-1}} \, \, (= \! \infty), \quad q \! = \! 
\mathfrak{s} \! - \! 1.
\end{equation*}

The next subset in this ordering contains all positive integers, {}from 
$1$ to $k$, corresponding to the pole $\alpha_{k}$ $(\neq \! \infty)$. 
For $k \! \in \! \lbrace 1,2,\dotsc,K \rbrace$ such that $\alpha_{k} \! 
\neq \! \infty$, write the ordered disjoint integer partition
\begin{equation*}
\lbrace \mathstrut k^{\prime} \! \in \! \lbrace 1,2,\dotsc,K \rbrace; \, 
k^{\prime} \! \leqslant \! k, \, \alpha_{k^{\prime}} \! = \! \alpha_{k} 
\! \neq \! \infty \rbrace \! := \! \lbrace \underbrace{i(\mathfrak{s})_{1},
i(\mathfrak{s})_{2},\dotsc,
i(\mathfrak{s})_{k_{\mathfrak{s}}}}_{k_{\mathfrak{s}}} \rbrace,
\end{equation*}
with $i(\mathfrak{s})_{k_{\mathfrak{s}}} \! = \! k$, $1 \! \leqslant \! 
i(\mathfrak{s})_{1} \! < \! i(\mathfrak{s})_{2} \! < \! \dotsb \! < \! 
i(\mathfrak{s})_{k_{\mathfrak{s}}} \! \leqslant \! K$, $\# \lbrace 
i(\mathfrak{s})_{1},i(\mathfrak{s})_{2},\dotsc,
i(\mathfrak{s})_{k_{\mathfrak{s}}} \rbrace \! = \! k_{\mathfrak{s}}$ $(= \! 
\varkappa_{1i(\mathfrak{s})_{k_{\mathfrak{s}}}} \! = \! \varkappa_{1k} 
\! = \! \varrho_{i(\mathfrak{s})_{k_{\mathfrak{s}}}})$ $= \! \varrho_{k}$, 
and $\lbrace i(j)_{1},i(j)_{2},\dotsc,i(j)_{k_{j}} \rbrace \cap \lbrace 
i(\mathfrak{s})_{1},i(\mathfrak{s})_{2},\dotsc,
i(\mathfrak{s})_{k_{\mathfrak{s}}} \rbrace \! = \! \varnothing$, 
$j \! = \! 1,2,\dotsc,\mathfrak{s} \! - \! 1$, which induces, by the 
above definition, the following pole ordering,
\begin{equation*}
\lbrace \mathstrut \alpha_{k^{\prime}}, \, k^{\prime} \! \in \! \lbrace 
1,2,\dotsc,K \rbrace; \, k^{\prime} \! \leqslant \! k, \, \alpha_{k^{\prime}} 
\! = \! \alpha_{k} \! \neq \! \infty \rbrace \! := \! \lbrace 
\underbrace{\alpha_{i(\mathfrak{s})_{1}},\alpha_{i(\mathfrak{s})_{2}},
\dotsc,\alpha_{i(\mathfrak{s})_{k_{\mathfrak{s}}}}}_{k_{\mathfrak{s}}} 
\rbrace,
\end{equation*}
with $\alpha_{i(\mathfrak{s})_{k_{\mathfrak{s}}}} \! = \! \alpha_{k}$ 
$(\neq \! \infty)$, $\alpha_{i(\mathfrak{s})_{1}} \! \prec \! 
\alpha_{i(\mathfrak{s})_{2}} \! \prec \! \dotsb \! \prec \! 
\alpha_{i(\mathfrak{s})_{k_{\mathfrak{s}}}}$, and $\lbrace \alpha_{i(j)_{1}},
\alpha_{i(j)_{2}},\dotsc,\alpha_{i(j)_{k_{j}}} \rbrace \cap \lbrace 
\alpha_{i(\mathfrak{s})_{1}},\alpha_{i(\mathfrak{s})_{2}},\dotsc,
\alpha_{i(\mathfrak{s})_{k_{\mathfrak{s}}}} \rbrace \! = \! \varnothing$, 
$j \! = \! 1,2,\dotsc,\mathfrak{s} \! - \! 1$, such that
\begin{equation*}
\alpha_{i(\mathfrak{s})_{1}} \! = \! \alpha_{i(\mathfrak{s})_{2}} \! = 
\! \dotsb \! = \! \alpha_{i(\mathfrak{s})_{k_{\mathfrak{s}}}} \! = \! 
\alpha_{k} \! =: \! \alpha_{p_{\mathfrak{s}}} \! \neq \! \infty, \quad \quad 
\# \lbrace \mathstrut \alpha_{i(\mathfrak{s})_{1}},\alpha_{i(\mathfrak{s})_{2}},
\dotsc,\alpha_{i(\mathfrak{s})_{k_{\mathfrak{s}}}} \rbrace \! = \! 
k_{\mathfrak{s}} \! = \! \varrho_{k}.
\end{equation*}
For $(\mathbb{N} \! \ni)$ $n \! \geqslant \! 2$ and $k \! \in \! \lbrace 
1,2,\dotsc,K \rbrace$ such that $\alpha_{k} \! \neq \! \infty$,
\begin{equation*}
\# \lbrace \mathstrut \alpha_{i(\mathfrak{s})_{k_{\mathfrak{s}}}},
\alpha_{i(\mathfrak{s})_{k_{\mathfrak{s}}}},\dotsc,
\alpha_{i(\mathfrak{s})_{k_{\mathfrak{s}}}} \rbrace \! = \! 
\varkappa_{ni(\mathfrak{s})_{k_{\mathfrak{s}}}} \! = \! \varkappa_{nk} 
\! = \! (n \! - \! 1) \gamma_{k} \! + \! \varrho_{k},
\end{equation*}
that is, as one moves {}from left to right across the repeated pole sequence 
$\mathcal{P}_{n,k}$ and removes the residual pole set (cf. the discussion 
above) $\lbrace \mathstrut \alpha_{k^{\prime}}, \, k^{\prime} \! \in \! 
\lbrace 1,2,\dotsc,K \rbrace; \, \alpha_{k^{\prime}} \! \neq \! \alpha_{k}, 
\, \alpha_{k} \! \neq \! \infty \rbrace$, one is left with the set (via 
the above induced ordering) that consists of all occurrences of the pole 
$\alpha_{k} \! \neq \! \infty$, which occurs $\varkappa_{nk} \! = \! 
(n \! - \! 1) \gamma_{k} \! + \! \varrho_{k}$ times, namely,
\begin{align*}
\mathcal{P}_{n,k} \setminus \left(\bigcup_{q=1}^{\mathfrak{s}-2} 
\lbrace \underbrace{\alpha_{p_{q}},\alpha_{p_{q}},\dotsc,
\alpha_{p_{q}}}_{\varkappa_{nk \tilde{k}_{q}}} \rbrace \cup \lbrace 
\underbrace{\alpha_{p_{\mathfrak{s}-1}},\alpha_{p_{\mathfrak{s}-1}},
\dotsc,\alpha_{p_{\mathfrak{s}-1}}}_{\varkappa^{\infty}_{nk 
\tilde{k}_{\mathfrak{s}-1}}} \rbrace \right) =& \, \mathcal{P}_{n,k} 
\setminus \left(\bigcup_{q=1}^{\mathfrak{s}-2} \lbrace 
\underbrace{\alpha_{p_{q}},\alpha_{p_{q}},\dotsc,
\alpha_{p_{q}}}_{\varkappa_{nk \tilde{k}_{q}}} \rbrace \cup \lbrace 
\underbrace{\infty,\infty,\dotsc,\infty}_{\varkappa^{\infty}_{nk 
\tilde{k}_{\mathfrak{s}-1}}} \rbrace \right) \\
:=& \, \lbrace \underbrace{\alpha_{i(\mathfrak{s})_{k_{\mathfrak{s}}}},
\alpha_{i(\mathfrak{s})_{k_{\mathfrak{s}}}},\dotsc,
\alpha_{i(\mathfrak{s})_{k_{\mathfrak{s}}}}}_{\varkappa_{ni(\mathfrak{s})_{
k_{\mathfrak{s}}}}} \rbrace \! = \! \lbrace \underbrace{\alpha_{k},
\alpha_{k},\dotsc,\alpha_{k}}_{\varkappa_{nk}} \rbrace.
\end{align*}
In order to illustrate this latter notation, consider, again, the $K \! = \! 
7$ pole set $\lbrace \alpha_{1},\alpha_{2},\alpha_{3},\alpha_{4},\alpha_{5},
\alpha_{6},\alpha_{7} \rbrace \! = \! \lbrace 0,1,\infty,1,\linebreak[4]
\sqrt{2},\pi,\infty \rbrace$ for which $\mathfrak{s} \! = \! 5$ and 
$\alpha_{k} \! \neq \! \infty$, $k \! = \! 1,2,4,5,6$:
\begin{enumerate}
\item[(i)] \fbox{$k \! = \! 1$}
\begin{equation*}
\lbrace \mathstrut k^{\prime} \! \in \! \lbrace 1,2,\dotsc,7 \rbrace; \, 
k^{\prime} \! \leqslant \! 1, \, \alpha_{k^{\prime}} \! = \! \alpha_{1} 
\! = \! 0 \rbrace \! = \! \lbrace 1 \rbrace \! := \! \lbrace i(5)_{k_{5}} 
\rbrace \, \Rightarrow \, k_{5} \! = \! 1, \, \, i(5)_{1} \! = \! 1,
\end{equation*}
which induces the pole ordering
\begin{equation*}
\lbrace \alpha_{k^{\prime}}, \, k^{\prime} \! \in \! \lbrace 1,2,\dotsc,7 
\rbrace; \, k^{\prime} \! \leqslant \! 1, \, \alpha_{k^{\prime}} \! = \! 
\alpha_{1} \! = \! 0 \rbrace \! := \! \lbrace \alpha_{i(5)_{k_{5}}} \rbrace 
\! = \! \lbrace \alpha_{1} \rbrace \! = \! \lbrace 0 \rbrace,
\end{equation*}
hence
\begin{equation*}
\varkappa_{ni(\mathfrak{s})_{k_{\mathfrak{s}}}} \! = \! \varkappa_{ni(5)_{1}} 
\! = \! \varkappa_{n1} \! = \! (n \! - \! 1) \gamma_{1} \! + \! \varrho_{1} 
\! = \! (n \! - \! 1) \! + \! 1 \! = \! n,
\end{equation*}
that is, as one moves {}from left to right across the repeated pole sequence
\begin{align*}
\mathcal{P}_{n,1} =& \, \lbrace \overset{1}{\underbrace{\alpha_{1},
\alpha_{2},\dotsc,\alpha_{7}}_{7}} \rbrace \cup \dotsb \cup \lbrace 
\overset{n-1}{\underbrace{\alpha_{1},\alpha_{2},\dotsc,\alpha_{7}}_{7}} 
\rbrace \cup \lbrace \overset{n}{\underbrace{\alpha_{1}}_{1}} \rbrace \\
=& \, \lbrace \overset{1}{\underbrace{0,1,\infty,1,\sqrt{2},\pi,\infty}_{7}} 
\rbrace \cup \dotsb \cup \lbrace \overset{n-1}{\underbrace{0,1,\infty,1,
\sqrt{2},\pi,\infty}_{7}} \rbrace \cup \lbrace \overset{n}{\underbrace{0}_{1}} 
\rbrace
\end{align*}
and removes the residual pole set (cf. the discussion and examples above) 
$\lbrace \mathstrut \alpha_{k^{\prime}}, \, k^{\prime} \! \in \! \lbrace 
1,2,\dotsc,7 \rbrace; \, \alpha_{k^{\prime}} \! \neq \! \alpha_{1} \! = \! 
0 \rbrace$, one is left with the set (via the above induced ordering) that 
consists of all occurrences of the pole $\alpha_{1} \! = \! 0$, which occurs 
$\varkappa_{n1} \! = \! n$ times, namely,
\begin{align*}
\mathcal{P}_{n,1} \setminus \left(\bigcup_{q=1}^{3} \lbrace 
\underbrace{\alpha_{p_{q}},\alpha_{p_{q}},\dotsc,
\alpha_{p_{q}}}_{\varkappa_{n1 \tilde{k}_{q}}} \rbrace \cup 
\lbrace \underbrace{\alpha_{p_{4}},\alpha_{p_{4}},\dotsc,
\alpha_{p_{4}}}_{\varkappa^{\infty}_{n1 \tilde{k}_{4}}} \rbrace \right) 
=& \, \mathcal{P}_{n,1} \setminus \left(\bigcup_{q=1}^{3} \lbrace 
\underbrace{\alpha_{p_{q}},\alpha_{p_{q}},\dotsc,
\alpha_{p_{q}}}_{\varkappa_{n1 \tilde{k}_{q}}} \rbrace \cup \lbrace 
\underbrace{\infty,\infty,\dotsc,\infty}_{2(n-1)} \rbrace \right) \\
:=& \, \lbrace \underbrace{\alpha_{i(5)_{k_{5}}},\alpha_{i(5)_{k_{5}}},
\dotsc,\alpha_{i(5)_{k_{5}}}}_{\varkappa_{ni(5)_{k_{5}}}} \rbrace 
\! = \! \lbrace \underbrace{\alpha_{1},\alpha_{1},\dotsc,
\alpha_{1}}_{\varkappa_{n1}} \rbrace \\
=& \, \lbrace \underbrace{0,0,\dotsc,0}_{n} \rbrace;
\end{align*}
\item[(ii)] \fbox{$k \! = \! 2$}
\begin{equation*}
\lbrace \mathstrut k^{\prime} \! \in \! \lbrace 1,2,\dotsc,7 \rbrace; \, 
k^{\prime} \! \leqslant \! 2, \, \alpha_{k^{\prime}} \! = \! \alpha_{2} 
\! = \! 1 \rbrace \! = \! \lbrace 2 \rbrace \! := \! \lbrace i(5)_{k_{5}} 
\rbrace \, \Rightarrow \, k_{5} \! = \! 1, \, \, i(5)_{1} \! = \! 2,
\end{equation*}
which induces the pole ordering
\begin{equation*}
\lbrace \alpha_{k^{\prime}}, \, k^{\prime} \! \in \! \lbrace 1,2,\dotsc,7 
\rbrace; \, k^{\prime} \! \leqslant \! 2, \, \alpha_{k^{\prime}} \! = \! 
\alpha_{2} \! = \! 1 \rbrace \! := \! \lbrace \alpha_{i(5)_{k_{5}}} \rbrace 
\! = \! \lbrace \alpha_{2} \rbrace \! = \! \lbrace 1 \rbrace,
\end{equation*}
hence
\begin{equation*}
\varkappa_{ni(\mathfrak{s})_{k_{\mathfrak{s}}}} \! = \! \varkappa_{ni(5)_{1}} 
\! = \! \varkappa_{n2} \! = \! (n \! - \! 1) \gamma_{2} \! + \! \varrho_{2} 
\! = \! 2(n \! - \! 1) \! + \! 1 \! = \! 2n \! - \! 1,
\end{equation*}
that is, as one moves {}from left to right across the repeated pole sequence
\begin{align*}
\mathcal{P}_{n,2} =& \, \lbrace \overset{1}{\underbrace{\alpha_{1},
\alpha_{2},\dotsc,\alpha_{7}}_{7}} \rbrace \cup \dotsb \cup \lbrace 
\overset{n-1}{\underbrace{\alpha_{1},\alpha_{2},\dotsc,\alpha_{7}}_{7}} 
\rbrace \cup \lbrace \overset{n}{\underbrace{\alpha_{1},\alpha_{2}}_{2}} 
\rbrace \\
=& \, \lbrace \overset{1}{\underbrace{0,1,\infty,1,\sqrt{2},\pi,\infty}_{7}} 
\rbrace \cup \dotsb \cup \lbrace \overset{n-1}{\underbrace{0,1,\infty,1,
\sqrt{2},\pi,\infty}_{7}} \rbrace \cup \lbrace \overset{n}{\underbrace{0,
1}_{2}} \rbrace
\end{align*}
and removes the residual pole set (cf. the discussion and examples above) 
$\lbrace \mathstrut \alpha_{k^{\prime}}, \, k^{\prime} \! \in \! \lbrace 
1,2,\dotsc,7 \rbrace; \, \alpha_{k^{\prime}} \! \neq \! \alpha_{2} \! = \! 
1 \rbrace$, one is left with the set (via the above induced ordering) that 
consists of all occurrences of the pole $\alpha_{2} \! = \! 1$, which occurs 
$\varkappa_{n2} \! = \! 2n \! - \! 1$ times, namely,
\begin{align*}
\mathcal{P}_{n,2} \setminus \left(\bigcup_{q=1}^{3} \lbrace 
\underbrace{\alpha_{p_{q}},\alpha_{p_{q}},\dotsc,
\alpha_{p_{q}}}_{\varkappa_{n2 \tilde{k}_{q}}} \rbrace \cup 
\lbrace \underbrace{\alpha_{p_{4}},\alpha_{p_{4}},\dotsc,
\alpha_{p_{4}}}_{\varkappa^{\infty}_{n2 \tilde{k}_{4}}} \rbrace \right) 
=& \, \mathcal{P}_{n,2} \setminus \left(\bigcup_{q=1}^{3} \lbrace 
\underbrace{\alpha_{p_{q}},\alpha_{p_{q}},\dotsc,
\alpha_{p_{q}}}_{\varkappa_{n2 \tilde{k}_{q}}} \rbrace \cup \lbrace 
\underbrace{\infty,\infty,\dotsc,\infty}_{2(n-1)} \rbrace \right) \\
:=& \, \lbrace \underbrace{\alpha_{i(5)_{k_{5}}},\alpha_{i(5)_{k_{5}}},
\dotsc,\alpha_{i(5)_{k_{5}}}}_{\varkappa_{ni(5)_{k_{5}}}} \rbrace 
\! = \! \lbrace \underbrace{\alpha_{2},\alpha_{2},\dotsc,
\alpha_{2}}_{\varkappa_{n2}} \rbrace \\
=& \, \lbrace \underbrace{1,1,\dotsc,1}_{2n-1} \rbrace;
\end{align*}
\item[(iii)] \fbox{$k \! = \! 4$}
\begin{equation*}
\lbrace \mathstrut k^{\prime} \! \in \! \lbrace 1,2,\dotsc,7 \rbrace; \, 
k^{\prime} \! \leqslant \! 4, \, \alpha_{k^{\prime}} \! = \! \alpha_{4} 
\! = \! 1 \rbrace \! = \! \lbrace 2,4 \rbrace \! := \! \lbrace i(5)_{1},
i(5)_{k_{5}} \rbrace \, \Rightarrow \, k_{5} \! = \! 2, \, \, i(5)_{1} \! 
= \! 2, \, \, i(5)_{2} \! = \! 4,
\end{equation*}
which induces the pole ordering
\begin{equation*}
\lbrace \alpha_{k^{\prime}}, \, k^{\prime} \! \in \! \lbrace 1,2,\dotsc,
7 \rbrace; \, k^{\prime} \! \leqslant \! 4, \, \alpha_{k^{\prime}} \! 
= \! \alpha_{4} \! = \! 1 \rbrace \! := \! \lbrace \alpha_{i(5)_{1}},
\alpha_{i(5)_{k_{5}}} \rbrace \! = \! \lbrace \alpha_{2},\alpha_{4} 
\rbrace \! = \! \lbrace 1,1 \rbrace,
\end{equation*}
hence
\begin{equation*}
\varkappa_{ni(\mathfrak{s})_{k_{\mathfrak{s}}}} \! = \! \varkappa_{ni(5)_{2}} 
\! = \! \varkappa_{n4} \! = \! (n \! - \! 1) \gamma_{4} \! + \! \varrho_{4} 
\! = \! 2(n \! - \! 1) \! + \! 2 \! = \! 2n,
\end{equation*}
that is, as one moves {}from left to right across the repeated pole sequence
\begin{align*}
\mathcal{P}_{n,4} =& \, \lbrace \overset{1}{\underbrace{\alpha_{1},
\alpha_{2},\dotsc,\alpha_{7}}_{7}} \rbrace \cup \dotsb \cup \lbrace 
\overset{n-1}{\underbrace{\alpha_{1},\alpha_{2},\dotsc,\alpha_{7}}_{7}} 
\rbrace \cup \lbrace \overset{n}{\underbrace{\alpha_{1},\alpha_{2},
\alpha_{3},\alpha_{4}}_{4}} \rbrace \\
=& \, \lbrace \overset{1}{\underbrace{0,1,\infty,1,\sqrt{2},\pi,\infty}_{7}} 
\rbrace \cup \dotsb \cup \lbrace \overset{n-1}{\underbrace{0,1,\infty,1,
\sqrt{2},\pi,\infty}_{7}} \rbrace \cup \lbrace \overset{n}{\underbrace{0,1,
\infty,1}_{4}} \rbrace
\end{align*}
and removes the residual pole set (cf. the discussion and examples above) 
$\lbrace \mathstrut \alpha_{k^{\prime}}, \, k^{\prime} \! \in \! \lbrace 
1,2,\dotsc,7 \rbrace; \, \alpha_{k^{\prime}} \! \neq \! \alpha_{4} \! = \! 
1 \rbrace$, one is left with the set (via the above induced ordering) that 
consists of all occurrences of the pole $\alpha_{4} \! = \! 1$, which occurs 
$\varkappa_{n4} \! = \! 2n$ times, namely,
\begin{align*}
\mathcal{P}_{n,4} \setminus \left(\bigcup_{q=1}^{3} \lbrace 
\underbrace{\alpha_{p_{q}},\alpha_{p_{q}},\dotsc,
\alpha_{p_{q}}}_{\varkappa_{n4 \tilde{k}_{q}}} \rbrace \cup 
\lbrace \underbrace{\alpha_{p_{4}},\alpha_{p_{4}},\dotsc,
\alpha_{p_{4}}}_{\varkappa^{\infty}_{n4 \tilde{k}_{4}}} \rbrace \right) 
=& \, \mathcal{P}_{n,4} \setminus \left(\bigcup_{q=1}^{3} \lbrace 
\underbrace{\alpha_{p_{q}},\alpha_{p_{q}},\dotsc,
\alpha_{p_{q}}}_{\varkappa_{n4 \tilde{k}_{q}}} \rbrace \cup \lbrace 
\underbrace{\infty,\infty,\dotsc,\infty}_{2n-1} \rbrace \right) \\
:=& \, \lbrace \underbrace{\alpha_{i(5)_{k_{5}}},\alpha_{i(5)_{k_{5}}},
\dotsc,\alpha_{i(5)_{k_{5}}}}_{\varkappa_{ni(5)_{k_{5}}}} \rbrace 
\! = \! \lbrace \underbrace{\alpha_{4},\alpha_{4},\dotsc,
\alpha_{4}}_{\varkappa_{n4}} \rbrace \\
=& \, \lbrace \underbrace{1,1,\dotsc,1}_{2n} \rbrace;
\end{align*}
\item[(iv)] \fbox{$k \! = \! 5$}
\begin{equation*}
\lbrace \mathstrut k^{\prime} \! \in \! \lbrace 1,2,\dotsc,7 \rbrace; \, 
k^{\prime} \! \leqslant \! 5, \, \alpha_{k^{\prime}} \! = \! \alpha_{5} 
\! = \! \sqrt{2} \rbrace \! = \! \lbrace 5 \rbrace \! := \! \lbrace 
i(5)_{k_{5}} \rbrace \, \Rightarrow \, k_{5} \! = \! 1, \, \, i(5)_{1} \! = \! 5,
\end{equation*}
which induces the pole ordering
\begin{equation*}
\lbrace \alpha_{k^{\prime}}, \, k^{\prime} \! \in \! \lbrace 1,2,\dotsc,
7 \rbrace; \, k^{\prime} \! \leqslant \! 5, \, \alpha_{k^{\prime}} \! = \! 
\alpha_{5} \! = \! \sqrt{2} \rbrace \! := \! \lbrace \alpha_{i(5)_{k_{5}}} 
\rbrace \! = \! \lbrace \alpha_{5} \rbrace \! = \! \lbrace \sqrt{2} \rbrace,
\end{equation*}
hence
\begin{equation*}
\varkappa_{ni(\mathfrak{s})_{k_{\mathfrak{s}}}} \! = \! \varkappa_{ni(5)_{1}} 
\! = \! \varkappa_{n5} \! = \! (n \! - \! 1) \gamma_{5} \! + \! \varrho_{5} 
\! = \! (n \! - \! 1) \! + \! 1 \! = \! n,
\end{equation*}
that is, as one moves {}from left to right across the repeated pole sequence
\begin{align*}
\mathcal{P}_{n,5} =& \, \lbrace \overset{1}{\underbrace{\alpha_{1},
\alpha_{2},\dotsc,\alpha_{7}}_{7}} \rbrace \cup \dotsb \cup \lbrace 
\overset{n-1}{\underbrace{\alpha_{1},\alpha_{2},\dotsc,\alpha_{7}}_{7}} 
\rbrace \cup \lbrace \overset{n}{\underbrace{\alpha_{1},\alpha_{2},
\alpha_{3},\alpha_{4},\alpha_{5}}_{5}} \rbrace \\
=& \, \lbrace \overset{1}{\underbrace{0,1,\infty,1,\sqrt{2},\pi,\infty}_{7}} 
\rbrace \cup \dotsb \cup \lbrace \overset{n-1}{\underbrace{0,1,\infty,1,
\sqrt{2},\pi,\infty}_{7}} \rbrace \cup \lbrace \overset{n}{\underbrace{0,1,
\infty,1,\sqrt{2}}_{5}} \rbrace
\end{align*}
and removes the residual pole set (cf. the discussion and examples above) 
$\lbrace \mathstrut \alpha_{k^{\prime}}, \, k^{\prime} \! \in \! \lbrace 
1,2,\dotsc,7 \rbrace; \, \alpha_{k^{\prime}} \! \neq \! \alpha_{5} \! = \! 
\sqrt{2} \rbrace$, one is left with the set (via the above induced ordering) 
that consists of all occurrences of the pole $\alpha_{5} \! = \! \sqrt{2}$, 
which occurs $\varkappa_{n5} \! = \! n$ times, namely,
\begin{align*}
\mathcal{P}_{n,5} \setminus \left(\bigcup_{q=1}^{3} \lbrace 
\underbrace{\alpha_{p_{q}},\alpha_{p_{q}},\dotsc,
\alpha_{p_{q}}}_{\varkappa_{n5 \tilde{k}_{q}}} \rbrace \cup 
\lbrace \underbrace{\alpha_{p_{4}},\alpha_{p_{4}},\dotsc,
\alpha_{p_{4}}}_{\varkappa^{\infty}_{n5 \tilde{k}_{4}}} \rbrace \right) 
=& \, \mathcal{P}_{n,5} \setminus \left(\bigcup_{q=1}^{3} \lbrace 
\underbrace{\alpha_{p_{q}},\alpha_{p_{q}},\dotsc,
\alpha_{p_{q}}}_{\varkappa_{n5 \tilde{k}_{q}}} \rbrace \cup \lbrace 
\underbrace{\infty,\infty,\dotsc,\infty}_{2n-1} \rbrace \right) \\
:=& \, \lbrace \underbrace{\alpha_{i(5)_{k_{5}}},\alpha_{i(5)_{k_{5}}},
\dotsc,\alpha_{i(5)_{k_{5}}}}_{\varkappa_{ni(5)_{k_{5}}}} \rbrace 
\! = \! \lbrace \underbrace{\alpha_{5},\alpha_{5},\dotsc,
\alpha_{5}}_{\varkappa_{n5}} \rbrace \\
=& \, \lbrace \underbrace{\sqrt{2},\sqrt{2},\dotsc,\sqrt{2}}_{n} \rbrace;
\end{align*}
\item[(v)] \fbox{$k \! = \! 6$}
\begin{equation*}
\lbrace \mathstrut k^{\prime} \! \in \! \lbrace 1,2,\dotsc,7 \rbrace; \, 
k^{\prime} \! \leqslant \! 6, \, \alpha_{k^{\prime}} \! = \! \alpha_{6} 
\! = \! \pi \rbrace \! = \! \lbrace 6 \rbrace \! := \! \lbrace i(5)_{k_{5}} 
\rbrace \, \Rightarrow \, k_{5} \! = \! 1, \, \, i(5)_{1} \! = \! 6,
\end{equation*}
which induces the pole ordering
\begin{equation*}
\lbrace \alpha_{k^{\prime}}, \, k^{\prime} \! \in \! \lbrace 1,2,\dotsc,
7 \rbrace; \, k^{\prime} \! \leqslant \! 6, \, \alpha_{k^{\prime}} \! = 
\! \alpha_{6} \! = \! \pi \rbrace \! := \! \lbrace \alpha_{i(5)_{k_{5}}} 
\rbrace \! = \! \lbrace \alpha_{6} \rbrace \! = \! \lbrace \pi \rbrace,
\end{equation*}
hence
\begin{equation*}
\varkappa_{ni(\mathfrak{s})_{k_{\mathfrak{s}}}} \! = \! \varkappa_{ni(5)_{1}} 
\! = \! \varkappa_{n6} \! = \! (n \! - \! 1) \gamma_{6} \! + \! \varrho_{6} 
\! = \! (n \! - \! 1) \! + \! 1 \! = \! n,
\end{equation*}
that is, as one moves {}from left to right across the repeated pole sequence
\begin{align*}
\mathcal{P}_{n,6} =& \, \lbrace \overset{1}{\underbrace{\alpha_{1},
\alpha_{2},\dotsc,\alpha_{7}}_{7}} \rbrace \cup \dotsb \cup \lbrace 
\overset{n-1}{\underbrace{\alpha_{1},\alpha_{2},\dotsc,\alpha_{7}}_{7}} 
\rbrace \cup \lbrace \overset{n}{\underbrace{\alpha_{1},\alpha_{2},
\alpha_{3},\alpha_{4},\alpha_{5},\alpha_{6}}_{6}} \rbrace \\
=& \, \lbrace \overset{1}{\underbrace{0,1,\infty,1,\sqrt{2},\pi,\infty}_{7}} 
\rbrace \cup \dotsb \cup \lbrace \overset{n-1}{\underbrace{0,1,\infty,1,
\sqrt{2},\pi,\infty}_{7}} \rbrace \cup \lbrace \overset{n}{\underbrace{0,1,
\infty,1,\sqrt{2},\pi}_{6}} \rbrace
\end{align*}
and removes the residual pole set (cf. the discussion and examples above) 
$\lbrace \mathstrut \alpha_{k^{\prime}}, \, k^{\prime} \! \in \! \lbrace 
1,2,\dotsc,7 \rbrace; \, \alpha_{k^{\prime}} \! \neq \! \alpha_{6} \! = \! 
\pi \rbrace$, one is left with the set (via the above induced ordering) that 
consists of all occurrences of the pole $\alpha_{6} \! = \! \pi$, which occurs 
$\varkappa_{n6} \! = \! n$ times, namely,
\begin{align*}
\mathcal{P}_{n,6} \setminus \left(\bigcup_{q=1}^{3} \lbrace 
\underbrace{\alpha_{p_{q}},\alpha_{p_{q}},\dotsc,
\alpha_{p_{q}}}_{\varkappa_{n6 \tilde{k}_{q}}} \rbrace \cup 
\lbrace \underbrace{\alpha_{p_{4}},\alpha_{p_{4}},\dotsc,
\alpha_{p_{4}}}_{\varkappa^{\infty}_{n6 \tilde{k}_{4}}} \rbrace \right) 
=& \, \mathcal{P}_{n,6} \setminus \left(\bigcup_{q=1}^{3} \lbrace 
\underbrace{\alpha_{p_{q}},\alpha_{p_{q}},\dotsc,
\alpha_{p_{q}}}_{\varkappa_{n6 \tilde{k}_{q}}} \rbrace \cup \lbrace 
\underbrace{\infty,\infty,\dotsc,\infty}_{2n-1} \rbrace \right) \\
:=& \, \lbrace \underbrace{\alpha_{i(5)_{k_{5}}},\alpha_{i(5)_{k_{5}}},
\dotsc,\alpha_{i(5)_{k_{5}}}}_{\varkappa_{ni(5)_{k_{5}}}} \rbrace 
\! = \! \lbrace \underbrace{\alpha_{6},\alpha_{6},\dotsc,
\alpha_{6}}_{\varkappa_{n6}} \rbrace \\
=& \, \lbrace \underbrace{\pi,\pi,\dotsc,\pi}_{n} \rbrace.
\end{align*}
\end{enumerate}
This concludes the example.

With the above conventions and ordered disjoint partitions, one writes, for 
$n \! \in \! \mathbb{N}$ and $k \! \in \! \lbrace 1,2,\dotsc,K \rbrace$ 
such that $\alpha_{k} \! \neq \! \infty$, the repeated pole sequence 
$\mathcal{P}_{n,k}$ as the following ordered union of disjoint partitions:
\begin{align*}
\mathcal{P}_{n,k} =& \, \bigcup_{q=1}^{\mathfrak{s}-2} \lbrace 
\underbrace{\alpha_{i(q)_{k_{q}}},\alpha_{i(q)_{k_{q}}},\dotsc,
\alpha_{i(q)_{k_{q}}}}_{\varkappa_{nk \tilde{k}_{q}}} \rbrace \cup \lbrace 
\underbrace{\alpha_{i(\mathfrak{s}-1)_{k_{\mathfrak{s}-1}}},
\alpha_{i(\mathfrak{s}-1)_{k_{\mathfrak{s}-1}}},\dotsc,
\alpha_{i(\mathfrak{s}-1)_{k_{\mathfrak{s}-1}}}}_{\varkappa^{\infty}_{nk 
\tilde{k}_{\mathfrak{s}-1}}} \rbrace \cup \lbrace \underbrace{\alpha_{i
(\mathfrak{s})_{k_{\mathfrak{s}}}},\alpha_{i(\mathfrak{s})_{
k_{\mathfrak{s}}}},\dotsc,\alpha_{i(\mathfrak{s})_{k_{\mathfrak{s}}}}}_{
\varkappa_{ni(\mathfrak{s})_{k_{\mathfrak{s}}}}} \rbrace \\
:=& \, \bigcup_{q=1}^{\mathfrak{s}-2} \lbrace \underbrace{\alpha_{p_{q}},
\alpha_{p_{q}},\dotsc,\alpha_{p_{q}}}_{\varkappa_{nk \tilde{k}_{q}}} \rbrace 
\cup \lbrace \underbrace{\alpha_{p_{\mathfrak{s}-1}},\alpha_{p_{\mathfrak{s}
-1}},\dotsc,\alpha_{p_{\mathfrak{s}-1}}}_{\varkappa^{\infty}_{nk 
\tilde{k}_{\mathfrak{s}-1}}} \rbrace \cup \lbrace \underbrace{\alpha_{k},
\alpha_{k},\dotsc,\alpha_{k}}_{\varkappa_{nk}} \rbrace \\
=& \, \bigcup_{q=1}^{\mathfrak{s}-2} \lbrace \underbrace{\alpha_{p_{q}},
\alpha_{p_{q}},\dotsc,\alpha_{p_{q}}}_{\varkappa_{nk \tilde{k}_{q}}} \rbrace 
\cup \lbrace \underbrace{\infty,\infty,\dotsc,\infty}_{\varkappa^{\infty}_{nk 
\tilde{k}_{\mathfrak{s}-1}}} \rbrace \cup \lbrace \underbrace{\alpha_{k},
\alpha_{k},\dotsc,\alpha_{k}}_{\varkappa_{nk}} \rbrace,
\end{align*}
where, by convention, the set $\lbrace \alpha_{k},\alpha_{k},\dotsc,\alpha_{k} 
\rbrace$ is written as the right-most set.

With the above notational preamble concluded, one now returns to the precise 
formulation of the orthogonality conditions for the MPC ORFs. One has a nested 
sequence of rational base sets. For $n \! \in \! \mathbb{N}$ and $k \! \in 
\! \lbrace 1,2,\dotsc,K \rbrace$ such that $\alpha_{k} \! \neq \! \infty$, 
one member of this nested sequence of rational base sets is
\begin{align*}
& \left\{\mathstrut \text{const.},\, \overset{1}{\underbrace{
\mathscr{S}^{1}_{1}(z),\mathscr{S}^{1}_{2}(z),\dotsc,
\mathscr{S}^{1}_{K}(z)}_{K}},\, \overset{2}{\underbrace{\mathscr{S}^{2}_{1}
(z),\mathscr{S}^{2}_{2}(z),\dotsc,\mathscr{S}^{2}_{K}(z)}_{K}},\dotsc,\, 
\overset{n}{\underbrace{\mathscr{S}^{n}_{1}(z),\mathscr{S}^{n}_{2}(z),
\dotsc,\mathscr{S}^{n}_{k}(z)}_{k}} \, \right\} \\
& \, := \lbrace \text{const.} \rbrace \bigcup \cup_{q=1}^{\mathfrak{s}-2} 
\cup_{r=1}^{\varkappa_{nk \tilde{k}_{q}}} \left\{(z \! - \! 
\alpha_{p_{q}})^{-r} \right\} \bigcup \cup_{l=1}^{\varkappa^{\infty}_{nk 
\tilde{k}_{\mathfrak{s}-1}}} \left\{z^{l} \right\} \bigcup 
\cup_{m=1}^{\varkappa_{nk}} \left\{(z \! - \! \alpha_{k})^{-m} \right\} \\
& \, = \lbrace \text{const.} \rbrace \bigcup \cup_{q=1}^{\mathfrak{s}-2} 
\left\{(z \! - \! \alpha_{p_{q}})^{-1},(z \! - \! \alpha_{p_{q}})^{-2},
\dotsc,(z \! - \! \alpha_{p_{q}})^{-\varkappa_{nk \tilde{k}_{q}}} \right\} 
\bigcup \left\{z,z^{2},\dotsc,z^{\varkappa^{\infty}_{nk \tilde{k}_{
\mathfrak{s}-1}}} \right\} \\
& \, \bigcup \left\{(z \! - \! \alpha_{k})^{-1},(z \! - \! \alpha_{k})^{-2},
\dotsc,(z \! - \! \alpha_{k})^{-\varkappa_{nk}} \right\},
\end{align*}
corresponding, respectively, to the ordered repeated pole sequence
\begin{align*}
& \left\{\mathstrut \text{no pole},\, \overset{1}{\underbrace{\alpha_{1},
\alpha_{2},\dotsc,\alpha_{K}}_{K}},\, \overset{2}{\underbrace{\alpha_{1},
\alpha_{2},\dotsc,\alpha_{K}}_{K}},\, \dotsc,\, \overset{n}{\underbrace{
\alpha_{1},\alpha_{2},\dotsc,\alpha_{k}}_{k}} \, \right\} \\
& \, := \lbrace \text{no pole} \rbrace \bigcup \cup_{q=1}^{\mathfrak{s}-2} 
\lbrace \underbrace{\alpha_{p_{q}},\alpha_{p_{q}},\dotsc,
\alpha_{p_{q}}}_{\varkappa_{nk \tilde{k}_{q}}} \rbrace \bigcup \lbrace 
\underbrace{\alpha_{p_{\mathfrak{s}-1}},\alpha_{p_{\mathfrak{s}-1}},\dotsc,
\alpha_{p_{\mathfrak{s}-1}}}_{\varkappa^{\infty}_{nk \tilde{k}_{\mathfrak{s}
-1}}} \rbrace \bigcup \lbrace \underbrace{\alpha_{k},\alpha_{k},\dotsc,
\alpha_{k}}_{\varkappa_{nk}} \rbrace \\
& \, = \lbrace \text{no pole} \rbrace \bigcup \cup_{q=1}^{\mathfrak{s}-2} 
\lbrace \underbrace{\alpha_{p_{q}},\alpha_{p_{q}},\dotsc,
\alpha_{p_{q}}}_{\varkappa_{nk \tilde{k}_{q}}} \rbrace \bigcup \lbrace 
\underbrace{\infty,\infty,\dotsc,\infty}_{\varkappa^{\infty}_{nk 
\tilde{k}_{\mathfrak{s}-1}}} \rbrace \bigcup \lbrace \underbrace{\alpha_{k},
\alpha_{k},\dotsc,\alpha_{k}}_{\varkappa_{nk}} \rbrace.
\end{align*}
Orthonormalisation with respect to $\langle \pmb{\cdot},\pmb{\cdot} 
\rangle_{\mathscr{L}}$, via the Gram-Schmidt orthogonalisation method, 
leads to the MPC ORFs, $\lbrace \phi^{n}_{k}(z) 
\rbrace_{\underset{k=1,2,\dotsc,K}{n \in \mathbb{N}}}$ (for consistency 
of notation, set $\phi^{0}_{0}(z) \! \equiv \! 1)$, which, by suitable 
normalisation (see below), can be written as, with the above orderings,
\begin{align*}
\phi^{n}_{k} \colon& \, \mathbb{N} \times \lbrace 1,2,\dotsc,K \rbrace 
\times \overline{\mathbb{C}} \setminus \lbrace \alpha_{1},\alpha_{2},
\dotsc,\alpha_{K} \rbrace \! \to \! \mathbb{C}, \, \, (n,k,z) \! \mapsto \\
\phi^{n}_{k}(z) :=& \, \phi_{0}^{f}(n,k) \! + \! \sum_{q=1}^{\mathfrak{s}
-2} \sum_{r=1}^{\varkappa_{nk \tilde{k}_{q}}} \dfrac{\tilde{\nu}^{f}_{r,q}
(n,k)}{(z \! - \! \alpha_{p_{q}})^{r}} \! + \! \sum_{l=1}^{\varkappa^{
\infty}_{nk \tilde{k}_{\mathfrak{s}-1}}} \hat{\nu}^{f}_{n,l}(n,k)z^{l} \! + 
\! \sum_{m=1}^{\varkappa_{nk}} \dfrac{\mu^{f}_{n,m}(n,k)}{(z \! - \! 
\alpha_{k})^{m}},
\end{align*}
where the $\phi^{n}_{k}$'s are normalised so that they all have real 
coefficients; in particular, for $n \! \in \! \mathbb{N}$ and $k \! \in \! 
\lbrace 1,2,\dotsc,K \rbrace$ such that $\alpha_{k} \! \neq \! \infty$,
\begin{equation*}
\operatorname{LC}(\phi^{n}_{k}) \! = \! \mu^{f}_{n,\varkappa_{nk}}
(n,k) \! > \! 0.
\end{equation*}
(For consistency of notation, set $\phi^{f}_{0}(0,0) \! \equiv \! 1$.) 
Furthermore, for $n \! \in \! \mathbb{N}$ and $k \! \in \! \lbrace 1,2,
\dotsc,K \rbrace$ such that $\alpha_{k} \! \neq \! \infty$, note that, 
by construction:
\begin{gather}
\left\langle \phi^{n}_{k},(z \! - \! \alpha_{k})^{-j} 
\right\rangle_{\mathscr{L}} \! = \! \int_{\mathbb{R}} \phi^{n}_{k}(\xi)
(\xi \! - \! \alpha_{k})^{-j} \, \md \mu (\xi) \! = \! 0, \quad j 
\! = \! 0,1,\dotsc,\varkappa_{nk} \! - \! 1, \label{eq12} \\
\left\langle \phi^{n}_{k},\mu^{f}_{n,\varkappa_{nk}}(n,k)
(z \! - \! \alpha_{k})^{-\varkappa_{nk}} \right\rangle_{\mathscr{L}} \! 
= \! \mu^{f}_{n,\varkappa_{nk}}(n,k) \int_{\mathbb{R}} \phi^{n}_{k}(\xi)
(\xi \! - \! \alpha_{k})^{-\varkappa_{nk}} \, \md \mu (\xi) \! = \! 1, 
\label{eq13} \\
\left\langle \phi^{n}_{k},(z \! - \! \alpha_{p_{q}})^{-r_{1}} 
\right\rangle_{\mathscr{L}} \! = \! \int_{\mathbb{R}} \phi^{n}_{k}(\xi)
(\xi \! - \! \alpha_{p_{q}})^{-r_{1}} \, \md \mu (\xi) \! = \! 0, \quad 
q \! = \! 1,2,\dotsc,\mathfrak{s} \! - \! 2, \quad r_{1} \! = \! 1,2,
\dotsc,\varkappa_{nk \tilde{k}_{q}}, \label{eq14} \\
\left\langle \phi^{n}_{k},z^{r_{2}} \right\rangle_{\mathscr{L}} \! = \! 
\int_{\mathbb{R}} \phi^{n}_{k}(\xi) \xi^{r_{2}} \, \md \mu (\xi) 
\! = \! 0, \quad r_{2} \! = \! 1,2,\dotsc,
\varkappa^{\infty}_{nk \tilde{k}_{\mathfrak{s}-1}}. \label{eq15} 
\end{gather}
(Note: if, for $k \! \in \! \lbrace 1,2,\dotsc,K \rbrace$ such that 
$\alpha_{k} \! \neq \! \infty$, the residual pole set $\lbrace \mathstrut 
\alpha_{k^{\prime}}, \, k^{\prime} \! \in \! \lbrace 1,2,\dotsc,K \rbrace; 
\, \alpha_{k^{\prime}} \! \neq \! \alpha_{k}, \, \alpha_{k} \! \neq \! 
\infty \rbrace \! = \! \varnothing$, then the orthogonality 
conditions~\eqref{eq14} and~\eqref{eq15} are vacuous; actually, this can 
only occur for $n \! = \! 1$.) For $n \! \in \! \mathbb{N}$ and $k \! \in \! 
\lbrace 1,2,\dotsc,K \rbrace$ such that $\alpha_{k} \! \neq \! \infty$, the 
orthogonality conditions~\eqref{eq12}--\eqref{eq15} give rise to a total 
of (cf. Equation~\eqref{fincount})
\begin{equation*}
\sum_{q=1}^{\mathfrak{s}-2} \varkappa_{nk \tilde{k}_{q}} \! 
+ \! \varkappa^{\infty}_{nk \tilde{k}_{\mathfrak{s}-1}} \! + \! 
\varkappa_{nk} \! + \! 1 \! = \! (n \! - \! 1)K \! + \! k \! + \! 1
\end{equation*}
linear equations determining the $(n \! - \! 1)K \! + \! k \! + \! 1$ real 
$(n$- and $k$-dependent) coefficients.

It is convenient to introduce, at this stage, the main object of study of 
this monograph, namely, the corresponding \emph{monic} MPC ORFs, 
$\lbrace \pmb{\pi}^{n}_{k}(z) \rbrace_{\underset{k=1,2,\dotsc,K}{n \in 
\mathbb{N}}}$ (for consistency of notation, set $\pmb{\pi}^{0}_{0}(z) 
\! \equiv \! 1)$. For $n \! \in \! \mathbb{N}$ and $k \! \in \! \lbrace 
1,2,\dotsc,K \rbrace$ such that $\alpha_{k} \! \neq \! \infty$,
\begin{align*}
\pmb{\pi}^{n}_{k} \colon& \, \mathbb{N} \times \lbrace 1,2,\dotsc,K \rbrace 
\times \overline{\mathbb{C}} \setminus \lbrace \alpha_{1},\alpha_{2},\dotsc,
\alpha_{K} \rbrace \! \to \! \mathbb{C}, \, \, (n,k,z) \! \mapsto \\
\pmb{\pi}^{n}_{k}(z) :=& \, \dfrac{\phi^{n}_{k}(z)}{\operatorname{LC}
(\phi^{n}_{k})} \! = \! \dfrac{\phi^{f}_{0}(n,k)}{\mu^{f}_{n,\varkappa_{nk}}(n,k)} 
\! + \! \dfrac{1}{\mu^{f}_{n,\varkappa_{nk}}(n,k)} \sum_{q=1}^{\mathfrak{s}-2} 
\sum_{r=1}^{\varkappa_{nk \tilde{k}_{q}}} \dfrac{\tilde{\nu}^{f}_{r,q}(n,k)}{(z \! 
- \! \alpha_{p_{q}})^{r}} \! + \! \dfrac{1}{\mu^{f}_{n,\varkappa_{nk}}(n,k)} 
\sum_{l=1}^{\varkappa^{\infty}_{nk \tilde{k}_{\mathfrak{s}-1}}} 
\hat{\nu}^{f}_{n,l}(n,k)z^{l} \\
+& \, \dfrac{1}{\mu^{f}_{n,\varkappa_{nk}}(n,k)} \sum_{m=1}^{
\varkappa_{nk}-1} \dfrac{\mu^{f}_{n,m}(n,k)}{(z \! - \! \alpha_{k})^{m}} 
\! + \! \dfrac{1}{(z \! - \! \alpha_{k})^{\varkappa_{nk}}}.
\end{align*}
The monic MPC ORFs, $\lbrace \pmb{\pi}^{n}_{k}(z) 
\rbrace_{\underset{k=1,2,\dotsc,K}{n \in \mathbb{N}}}$, 
possess the following orthogonality properties:
\begin{gather}
\left\langle \pmb{\pi}^{n}_{k},(z \! - \! \alpha_{k})^{-j} 
\right\rangle_{\mathscr{L}} \! = \! \int_{\mathbb{R}} \pmb{\pi}^{n}_{k}
(\xi)(\xi \! - \! \alpha_{k})^{-j} \, \md \mu (\xi) \! = \! 0, \quad 
j \! = \! 0,1,\dotsc,\varkappa_{nk} \! - \! 1, \label{eq16} \\
\left\langle \pmb{\pi}^{n}_{k},(z \! - \! \alpha_{k})^{-\varkappa_{nk}} 
\right\rangle_{\mathscr{L}} \! = \! \int_{\mathbb{R}} \pmb{\pi}^{n}_{k}
(\xi)(\xi \! - \! \alpha_{k})^{-\varkappa_{nk}} \, \md \mu (\xi) 
\! = \! (\operatorname{LC}(\phi^{n}_{k}))^{-2} \! = \! 
(\mu^{f}_{n,\varkappa_{nk}}(n,k))^{-2}, \label{eq17} \\
\left\langle \pmb{\pi}^{n}_{k},(z \! - \! \alpha_{p_{q}})^{-r_{1}} 
\right\rangle_{\mathscr{L}} \! = \! \int_{\mathbb{R}} \pmb{\pi}^{n}_{k}
(\xi)(\xi \! - \! \alpha_{p_{q}})^{-r_{1}} \, \md \mu (\xi) \! = \! 0, 
\quad q \! = \! 1,2,\dotsc,\mathfrak{s} \! - \! 2, \quad r_{1} \! = \! 
1,2,\dotsc,\varkappa_{nk \tilde{k}_{q}}, \label{eq18} \\
\left\langle \pmb{\pi}^{n}_{k},z^{r_{2}} \right\rangle_{\mathscr{L}} 
\! = \! \int_{\mathbb{R}} \pmb{\pi}^{n}_{k}(\xi) \xi^{r_{2}} \, 
\md \mu (\xi) \! = \! 0, \quad r_{2} \! = \! 1,2,\dotsc,
\varkappa^{\infty}_{nk \tilde{k}_{\mathfrak{s}-1}}. \label{eq19} 
\end{gather}
(Note: if, for $k \! \in \! \lbrace 1,2,\dotsc,K \rbrace$ such that 
$\alpha_{k} \! \neq \! \infty$, the residual pole set $\lbrace \mathstrut 
\alpha_{k^{\prime}}, \, k^{\prime} \! \in \! \lbrace 1,2,\dotsc,K \rbrace; 
\, \alpha_{k^{\prime}} \! \neq \! \alpha_{k}, \, \alpha_{k} \! \neq \! 
\infty \rbrace \! = \! \varnothing$, then the orthogonality 
conditions~\eqref{eq18} and~\eqref{eq19} are vacuous; actually, this can 
occur only for $n \! = \! 1$.) For $n \! \in \! \mathbb{N}$ and $k \! \in 
\! \lbrace 1,2,\dotsc,K \rbrace$ such that $\alpha_{k} \! \neq \! \infty$, 
it follows {}from the monic MPC ORF orthogonality 
conditions~\eqref{eq16}--\eqref{eq19} that
\begin{equation*}
\left\langle \pmb{\pi}^{n}_{k},\pmb{\pi}^{n}_{k} \right\rangle_{\mathscr{L}} 
\! =: \! \lvert \lvert \pmb{\pi}^{n}_{k}(\pmb{\cdot}) \rvert 
\rvert^{2}_{\mathscr{L}} \! = \! (\operatorname{LC}(\phi^{n}_{k}))^{-2} 
\! = \! (\mu_{n,\varkappa_{nk}}^{f}(n,k))^{-2},
\end{equation*}
whence $\lvert \lvert \pmb{\pi}^{n}_{k}(\pmb{\cdot}) \rvert \rvert_{\mathscr{L}} 
\! = \! (\mu^{f}_{n,\varkappa_{nk}}(n,k))^{-1} \! > \! 0$.\footnote{See, in particular, 
Section~\ref{sec2}, Lemma~\ref{lem2.1}, Equations~\eqref{nhormmconstatfin}, 
\eqref{eq63}, and~\eqref{eq65}, and Corollary~\ref{cor2.1}, 
Equation~\eqref{eqnmctfin1}; incidentally, this also establishes the positive 
definiteness of the linear functional $\mathscr{L} \colon \Lambda^{\mathbb{R}} 
\! \to \! \mathbb{R}$, and hence the---real---bilinear form $\langle \pmb{\cdot},
\pmb{\cdot} \rangle_{\mathscr{L}} \colon \Lambda^{\mathbb{R}} \times 
\Lambda^{\mathbb{R}} \! \to \! \mathbb{R}$, defined at the beginning of 
Subsection~\ref{subsec1.2}.} 
\begin{eeee} \label{rem1.2.4} 
\textsl{For $n \! \in \! \mathbb{N}$ and $k \! \in \! \lbrace 1,2,\dotsc,
K \rbrace$ such that $\alpha_{p_{\mathfrak{s}}} \! := \! \alpha_{k} \! 
\neq \! \infty$, {}from the above partial fraction decomposition 
for the monic {\rm MPC ORFs}, $\lbrace \pmb{\pi}^{n}_{k}(z) 
\rbrace_{\underset{k=1,2,\dotsc,K}{n \in \mathbb{N}}}$, one writes
\begin{equation*}
\pmb{\pi}^{n}_{k} \colon \mathbb{N} \times \lbrace 1,2,\dotsc,K \rbrace 
\times \overline{\mathbb{C}} \setminus \lbrace \alpha_{1},\alpha_{2},\dotsc,
\alpha_{K} \rbrace \! \ni \! (n,k,z) \! \mapsto \! \dfrac{1}{\mu^{f}_{n,
\varkappa_{nk}}(n,k)} \dfrac{\tilde{\mathfrak{P}}^{n}_{k}(z)}{\prod_{q=1}^{
\mathfrak{s}-2}(z \! - \! \alpha_{p_{q}})^{\varkappa_{nk \tilde{k}_{q}}}
(z \! - \! \alpha_{k})^{\varkappa_{nk}}},
\end{equation*}
where
\begin{align*}
\tilde{\mathfrak{P}}^{n}_{k}(z) :=& \, \phi^{f}_{0}(n,k) 
\prod_{q=1}^{\mathfrak{s}-2}(z \! - \! \alpha_{p_{q}})^{\varkappa_{nk 
\tilde{k}_{q}}}(z \! - \! \alpha_{k})^{\varkappa_{nk}} \! + \! 
\prod_{q=1}^{\mathfrak{s}-2}(z \! - \! \alpha_{p_{q}})^{\varkappa_{nk 
\tilde{k}_{q}}}(z \! - \! \alpha_{k})^{\varkappa_{nk}} 
\sum_{r=1}^{\varkappa^{\infty}_{nk \tilde{k}_{\mathfrak{s}-1}}} 
\hat{\nu}^{f}_{n,r}(n,k)z^{r} \\
+& \, \prod_{q=1}^{\mathfrak{s}-2}(z \! - \! \alpha_{p_{q}})^{\varkappa_{nk 
\tilde{k}_{q}}} \sum_{r=1}^{\varkappa_{nk}} \mu^{f}_{n,r}(n,k)
(z \! - \! \alpha_{k})^{\varkappa_{nk}-r} \! + \! \sum_{j=1}^{\mathfrak{s}-2} 
\prod_{\substack{q=1\\q \neq j}}^{\mathfrak{s}-2}
(z \! - \! \alpha_{p_{q}})^{\varkappa_{nk \tilde{k}_{q}}}
(z \! - \! \alpha_{k})^{\varkappa_{nk}} \sum_{r=1}^{\varkappa_{nk 
\tilde{k}_{j}}} \tilde{\nu}^{f}_{r,j}(n,k)
(z \! - \! \alpha_{p_{j}})^{\varkappa_{nk \tilde{k}_{j}}-r}.
\end{align*}
Via the above definition of $\tilde{\mathfrak{P}}^{n}_{k}(z)$, it follows that, 
for $n \! \in \! \mathbb{N}$ and $k \! \in \! \lbrace 1,2,\dotsc,K \rbrace$ 
such that $\alpha_{p_{\mathfrak{s}}} \! := \! \alpha_{k} \! \neq \! \infty$, 
$\tilde{\mathfrak{P}}^{n}_{k}(z)$ is a polynomial of degree at most 
$\sum_{q=1}^{\mathfrak{s}-2} \varkappa_{nk \tilde{k}_{q}} \! + \! 
\varkappa^{\infty}_{nk \tilde{k}_{\mathfrak{s}-1}} \! + \! \varkappa_{nk} \! 
= \! (n \! - \! 1)K \! + \! k$. The following discussion constitutes a succinct 
synopsis of the classification theory for the zeros of the monic {\rm MPC ORF}, 
$\pmb{\pi}^{n}_{k}(z)$,\footnote{Similar statements apply, \emph{verbatim}, for 
the corresponding MPC ORF, $\phi^{n}_{k}(z) \! = \! \mu^{f}_{n,\varkappa_{nk}}
(n,k) \pmb{\pi}^{n}_{k}(z)$.} for the case $n \! \in \! \mathbb{N}$ and $k \! \in 
\! \lbrace 1,2,\dotsc,K \rbrace$ such that $\alpha_{p_{\mathfrak{s}}} \! := \! 
\alpha_{k} \! \neq \! \infty$, that is based on the seminal works of Nj\r{a}stad 
{\rm \cite{a1,a2,ol1}}. Via the above definition of $\tilde{\mathfrak{P}}^{n}_{k}
(z)$, a straightforward calculation shows that $\tilde{\mathfrak{P}}^{n}_{k}
(\alpha_{p_{q}}) \! = \! \tilde{\nu}^{f}_{\varkappa_{nk \tilde{k}_{q}},q}(n,k) 
\prod_{\underset{q^{\prime} \neq q}{q^{\prime}=1}}^{\mathfrak{s}-2}
(\alpha_{p_{q}} \! - \! \alpha_{p_{q^{\prime}}})^{\varkappa_{nk 
\tilde{k}_{q^{\prime}}}}(\alpha_{p_{q}} \! - \! \alpha_{k})^{\varkappa_{nk}}$, 
$q \! \in \! \lbrace 1,2,\dotsc,\mathfrak{s} \! - \! 2 \rbrace$, 
$\tilde{\mathfrak{P}}^{n}_{k}(z)$ has a pole of order $(n \! - \! 1)K \! + \! k$ 
at (the point at infinity) $\alpha_{p_{\mathfrak{s}-1}} \! = \! \infty$ with 
coefficient $\hat{\nu}^{f}_{n,\varkappa^{\infty}_{nk \tilde{k}_{\mathfrak{s}-1}}}
(n,k)$, that is, $\tilde{\mathfrak{P}}^{n}_{k}(z) \! =_{\overline{\mathbb{C}} \ni z 
\to \alpha_{p_{\mathfrak{s}-1}} = \infty} \! \hat{\nu}^{f}_{n,\varkappa^{\infty}_{n
k \tilde{k}_{\mathfrak{s}-1}}}(n,k)z^{(n-1)K+k}(1 \! + \! \mathcal{O}(z^{-1}))$, 
and, since $\mu^{f}_{n,\varkappa_{nk}}(n,k) \! > \! 0$, $\tilde{\mathfrak{P}}^{n}_{k}
(\alpha_{k})$ $(=: \! \tilde{\mathfrak{P}}^{n}_{k}(\alpha_{p_{\mathfrak{s}}}))$ 
$= \! \mu^{f}_{n,\varkappa_{nk}}(n,k) \prod_{q=1}^{\mathfrak{s}-2}(\alpha_{k} 
\! - \! \alpha_{p_{q}})^{\varkappa_{nk \tilde{k}_{q}}} \! \neq \! 0$. One calls the 
monic {\rm MPC ORF}, $\pmb{\pi}^{n}_{k}(z)$, and the corresponding index 
two-tuple $(n,k)$ \emph{degenerate} if $\# \lbrace \mathstrut z \! \in \! 
\overline{\mathbb{C}}; \, \tilde{\mathfrak{P}}^{n}_{k}(z) \! = \! 0 \rbrace \! \leqslant 
\! (n \! - \! 1)K \! + \! k \! - \! 1$, that is, $\deg (\tilde{\mathfrak{P}}^{n}_{k}) \! < \! 
(n \! - \! 1)K \! + \! k$. One calls the monic {\rm MPC ORF}, $\pmb{\pi}^{n}_{k}(z)$, 
and the corresponding index two-tuple $(n,k)$ $\alpha_{p_{q}}$-\emph{defective}, 
$q \! \in \! \lbrace 1,2,\dotsc,\mathfrak{s} \rbrace$, if $\alpha_{p_{q}}$ is a zero of 
$\tilde{\mathfrak{P}}^{n}_{k}(z)$, namely: {\rm (i)} if $\tilde{\nu}^{f}_{\varkappa_{nk 
\tilde{k}_{q}},q}(n,k) \! = \! 0$, $q \! \in \! \lbrace 1,2,\dotsc,\mathfrak{s} \! - \! 2 
\rbrace$, then $\tilde{\mathfrak{P}}^{n}_{k}(\alpha_{p_{q}}) \! = \! 0$, $q \! \in \! 
\lbrace 1,2,\dotsc,\mathfrak{s} \! - \! 2 \rbrace$, thus $\pmb{\pi}^{n}_{k}(z)$ and 
the corresponding index two-tuple $(n,k)$ are $\alpha_{p_{q}}$-defective$;$ 
{\rm (ii)} if $\hat{\nu}^{f}_{n,\varkappa^{\infty}_{nk \tilde{k}_{\mathfrak{s}-1}}}
(n,k) \! = \! 0$, then $\tilde{\mathfrak{P}}^{n}_{k}(z)$ has a pole of order 
$(n \! - \! 1)K \! + \! k \! - \! 1$ at (the point at infinity) $\alpha_{p_{\mathfrak{s}-1}} 
\! = \! \infty$ with coefficient $\hat{\nu}^{f}_{n,\varkappa^{\infty}_{nk \tilde{k}_{
\mathfrak{s}-1}}-1}(n,k)$, that is, $\tilde{\mathfrak{P}}^{n}_{k}(z) \! 
=_{\overline{\mathbb{C}} \ni z \to \alpha_{p_{\mathfrak{s}-1}} = \infty} \! 
\hat{\nu}^{f}_{n,\varkappa^{\infty}_{nk \tilde{k}_{\mathfrak{s}-1}}-1}(n,k)
z^{(n-1)K+k-1}(1 \! + \! \mathcal{O}(z^{-1}))$, thus $\pmb{\pi}^{n}_{k}(z)$ and 
the corresponding index two-tuple $(n,k)$ are $\alpha_{p_{\mathfrak{s}-1}}$ 
$(= \! \infty)$-defective; and {\rm (iii)} since $\mu^{f}_{n,\varkappa_{nk}}(n,k) \! > \! 
0$, that is, $\tilde{\mathfrak{P}}^{n}_{k}(\alpha_{k})$ $(=: \! \hat{\mathfrak{P}}^{n}_{k}
(\alpha_{p_{\mathfrak{s}}}))$ $\neq \! 0$, it follows that $\pmb{\pi}^{n}_{k}(z)$ 
and the corresponding index two-tuple $(n,k)$ can not be $\alpha_{k}$ $(=: \! 
\alpha_{p_{\mathfrak{s}}})$-defective. One calls the monic {\rm MPC ORF}, 
$\pmb{\pi}^{n}_{k}(z)$, (and the corresponding index two-tuple $(n,k))$ 
\emph{defective} if it is $\alpha_{p_{q}}$-defective for at least one 
$q \! \in \! \lbrace 1,2,\dotsc,\mathfrak{s} \! - \! 1 \rbrace$ or 
\emph{maximally defective} if it is $\alpha_{p_{q}}$-defective $\forall$ 
$q \! \in \! \lbrace 1,2,\dotsc,\mathfrak{s} \! - \! 1 \rbrace$. One calls the 
monic {\rm MPC ORF}, $\pmb{\pi}^{n}_{k}(z)$, (and the corresponding index 
two-tuple $(n,k))$ \emph{singular} if it is degenerate or defective; otherwise, 
it is non-singular. In this monograph, it is assumed that $\pmb{\pi}^{n}_{k}(z)$ 
(and the corresponding index two-tuple $(n,k))$ is neither degenerate nor 
defective, that is, $\deg (\tilde{\mathfrak{P}}^{n}_{k}) \! = \! (n \! - \! 1)K 
\! + \! k$, $\tilde{\nu}^{f}_{\varkappa_{nk \tilde{k}_{q}},q}(n,k) \! \neq \! 0$ 
$\forall$ $q \! \in \! \lbrace 1,2,\dotsc,\mathfrak{s} \! - \! 2 \rbrace$, and 
$\hat{\nu}^{f}_{n,\varkappa^{\infty}_{nk \tilde{k}_{\mathfrak{s}-1}}}(n,k) 
\! \neq \! 0$ (recall that, since $\mu^{f}_{n,\varkappa_{nk}}(n,k) \! > \! 0$, 
$\tilde{\mathfrak{P}}^{n}_{k}(\alpha_{k})$ $(=: \! \tilde{\mathfrak{P}}^{n}_{k}
(\alpha_{p_{\mathfrak{s}}}))$ $\neq \! 0)$$;$ therefore, hereafter, it is to be 
understood that `monic {\rm MPC ORF}' and `non-singular monic {\rm MPC ORF}' 
are synonymous objects (the singular case(s) will be considered elsewhere$)$$;$ 
for simplicity, though, only the phrase {\rm MPC ORF} (monic or not) will be used. 
For $n \! \in \! \mathbb{N}$ and $k \! \in \! \lbrace 1,2,\dotsc,K \rbrace$ such 
that $\alpha_{p_{\mathfrak{s}}} \! := \! \alpha_{k} \! \neq \! \infty$, writing 
the factorisation $\tilde{\mathfrak{P}}^{n}_{k}(z) \! := \! \hat{\nu}^{f}_{n,
\varkappa^{\infty}_{nk \tilde{k}_{\mathfrak{s}-1}}}(n,k) \prod_{j=1}^{(n-1)K+k}
(z \! - \! \tilde{\mathfrak{z}}^{n}_{k}(j))$,\footnote{Of course, 
$\tilde{\mathfrak{z}}^{n}_{k}(j)$ also depends on $\alpha_{p_{q}}$, 
$q \! = \! 1,\dotsc,\mathfrak{s} \! - \! 2,\mathfrak{s}$; but, for simplicity 
of notation, this extraneous dependence is suppressed.} where, counting 
multiplicities, $\lbrace \mathstrut \tilde{\mathfrak{z}}^{n}_{k}(j) 
\rbrace_{j=1}^{(n-1)K+k} \! := \! \lbrace \mathstrut z \! \in \! 
\overline{\mathbb{C}}; \, \tilde{\mathfrak{P}}^{n}_{k}(z) \! = \! 0 \rbrace$ 
$(= \! \lbrace \mathstrut z \! \in \! \overline{\mathbb{C}}; \, \pmb{\pi}^{n}_{k}
(z) \! = \! 0 \rbrace)$, it is shown in Section~\ref{sec3} (see, in particular, the 
corresponding items of Lemmata~\ref{lemrootz} and~\ref{lemetatomu} for 
precise statements) that, in the double-scaling limit $\mathscr{N},n \! \to \! 
\infty$ such that $\mathscr{N}/n \! \to \! 1$, the zeros (counting multiplicities) 
$\tilde{\mathfrak{z}}^{n}_{k}(j) \! \in \! \overline{\mathbb{R}} \setminus \lbrace 
\alpha_{1},\alpha_{2},\dotsc,\alpha_{K} \rbrace$, $j \! = \! 1,2,\dotsc,
(n \! - \! 1)K \! + \! k$$;$ more precisely, in the double-scaling limit 
$\mathscr{N},n \! \to \! \infty$ such that $\mathscr{N}/n \! \to \! 1$, 
$\lbrace \tilde{\mathfrak{z}}^{n}_{k}(j) \rbrace_{j=1}^{(n-1)K+k}$ 
accumulates on the real compact set $J_{f}$, where $J_{f}$ is defined in 
Subsection~\ref{subsub2}. (See, also,  {\rm \cite{a1,a2,ol1}}.$)$}
\end{eeee}
The brief discussion that follows is motivated, in part, by the seminal 
works of Nj\r{a}stad \cite{n2,n1} related to MPAs and ORFs: more 
precise statements can be located in Section~\ref{sek5}; see, in 
particular, Lemma~\ref{lemmpainffin}. Furthermore, in order to mitigate 
notational encumbrances, explicit $n$- and $k$-dependencies will be 
temporarily suppressed, except where absolutely necessary. Recall (cf. 
Subsection~\ref{subsubsec1.2.1}) the definition of the Markov-Stieltjes 
transform given by Equation~\eqref{mvssinf1}; for $n \! \in \! 
\mathbb{N}$ and $k \! \in \! \lbrace 1,2,\dotsc,K \rbrace$ such 
that $\alpha_{p_{\mathfrak{s}}} \! := \! \alpha_{k} \! \neq \! \infty$, a 
straightforward calculation shows that $\mathrm{F}_{\mu}(z)$ has the 
formal asymptotic expansions
\begin{equation} \label{mvssfin1} 
\mathrm{F}_{\mu}(z) \underset{z \to \alpha_{p_{q}}}{=} \sum_{j=
0}^{\infty}c^{(q)}_{j}(\alpha_{p_{q}})(z \! - \! \alpha_{p_{q}})^{j}, 
\quad q \! = \! 1,\dotsc,\mathfrak{s} \! - \! 2,\mathfrak{s},
\end{equation}
where $c^{(q)}_{j}(\alpha_{p_{q}}) \! := \! -\int_{\mathbb{R}}(\xi 
\! - \! \alpha_{p_{q}})^{-(1+j)} \, \md \mu (\xi)$, $(j,q) \! \in \! 
\mathbb{N}_{0} \times \lbrace 1,\dotsc,\mathfrak{s} \! - \! 2,
\mathfrak{s} \rbrace$, and
\begin{equation} \label{mvssfin2} 
\mathrm{F}_{\mu}(z) \underset{z \to \alpha_{p_{\mathfrak{s}-1}} 
= \infty}{=} \sum_{j=1}^{\infty}c^{(\infty)}_{j}z^{-j},
\end{equation}
where $c^{(\infty)}_{j} \! := \! \int_{\mathbb{R}} \xi^{j-1} \, \md 
\mu (\xi)$, $j \! \in \! \mathbb{N}$, with $c^{(\infty)}_{1} \! = \! 1$. 
For $n \! \in \! \mathbb{N}$ and $k \! \in \! \lbrace 1,2,\dotsc,K 
\rbrace$ such that $\alpha_{p_{\mathfrak{s}}} \! := \! \alpha_{k} 
\! \neq \! \infty$, define the associated $\mathrm{R}$-function 
\cite{n2,n1} as follows:
\begin{equation} \label{mvssfin3} 
\widetilde{\pmb{\mathrm{R}}}_{\mu} \colon \mathbb{N} \times \lbrace 
1,2,\dotsc,K \rbrace \times \overline{\mathbb{C}} \setminus \lbrace 
\alpha_{1},\alpha_{2},\dotsc,\alpha_{K} \rbrace \! \ni \! (n,k,z) \! 
\mapsto \! \int_{\mathbb{R}} \left(\dfrac{\pmb{\pi}^{n}_{k}(\xi) 
\! - \! \pmb{\pi}^{n}_{k}(z)}{\xi \! - \! z} \right) \md \mu (\xi) \! 
=: \! \widetilde{\pmb{\mathrm{R}}}_{\mu}(z):
\end{equation}
{}from Equation~\eqref{mvssfin3} and the fact that $\mu \! \in 
\! \mathscr{M}_{1}(\mathbb{R})$, one shows, via the identity 
$y_{1}^{m} \! - \! y_{2}^{m} \! = \! (y_{1} \! - \! y_{2})(y_{1}^{m-1} 
\! + \! y_{1}^{m-2}y_{2} \! + \! \dotsb \! + \! y_{1}y_{2}^{m-2} \! + 
\! y_{2}^{m-1})$, that $\widetilde{\pmb{\mathrm{R}}}_{\mu}(z)$ 
can be presented as the improper fraction
\begin{equation} \label{mvssfin4} 
\widetilde{\pmb{\mathrm{R}}}_{\mu}(z) \! = \! \dfrac{\widetilde{
\mathrm{U}}_{\mu}(z)}{\prod_{q=1}^{\mathfrak{s}-2}(z \! - 
\! \alpha_{p_{q}})^{\varkappa_{nk \tilde{k}_{q}}}(z \! - \! 
\alpha_{k})^{\varkappa_{nk}}},
\end{equation}
where $\widetilde{\mathrm{U}}_{\mu}(z) \! := \! \sum_{j=0}^{(n-1)K+
k-1} \tilde{r}_{j}z^{j}$, with $\operatorname{deg}(\widetilde{\mathrm{
U}}_{\mu}(z)) \! = \! (n \! - \! 1)K \! + \! k \! - \! 1$ (since $\tilde{r}_{
(n-1)K+k-1} \! \neq \! 0)$, and $\operatorname{deg}(\prod_{q=1}^{
\mathfrak{s}-2}(z \! - \! \alpha_{p_{q}})^{\varkappa_{nk \tilde{k}_{q}}}
(z \! - \! \alpha_{k})^{\varkappa_{nk}}) \! = \! (n \! - \! 1)K \! + \! k \! 
- \! \varkappa^{\infty}_{nk \tilde{k}_{\mathfrak{s}-1}}$; moreover, note 
that, for $n \! \in \! \mathbb{N}$ and $k \! \in \! \lbrace 1,2,\dotsc,K 
\rbrace$ such that $\alpha_{p_{\mathfrak{s}}} \! := \! \alpha_{k} \! \neq 
\! \infty$, the corresponding monic MPC ORF, $\pmb{\pi}^{n}_{k}(z)$, 
can also be presented as the improper fraction
\begin{equation} \label{mvssfin5} 
\pmb{\pi}^{n}_{k}(z) \! = \! \dfrac{\widetilde{\mathrm{V}}_{\mu}(z)}{
\prod_{q=1}^{\mathfrak{s}-2}(z \! - \! \alpha_{p_{q}})^{\varkappa_{nk 
\tilde{k}_{q}}}(z \! - \! \alpha_{k})^{\varkappa_{nk}}},
\end{equation}
where $\widetilde{\mathrm{V}}_{\mu}(z) \! := \! \sum_{j=0}^{(n-1)K+k} 
\tilde{t}_{j}z^{j}$, with $\operatorname{deg}(\widetilde{\mathrm{V}}_{
\mu}(z)) \! = \! (n \! - \! 1)K \! + \! k$ (since $\tilde{t}_{(n-1)K+k} \! 
\neq \! 0)$.\footnote{For $n \! \in \! \mathbb{N}$ and $k \! \in \! 
\lbrace 1,2,\dotsc,K \rbrace$ such that $\alpha_{p_{\mathfrak{s}}} \! 
:= \! \alpha_{k} \! \neq \! \infty$, note that $\widetilde{\mathrm{V}}_{
\mu}(z)$ is not monic.} It turns out that, for $n \! \in \! \mathbb{N}$ 
and $k \! \in \! \lbrace 1,2,\dotsc,K \rbrace$ such that $\alpha_{
p_{\mathfrak{s}}} \! := \! \alpha_{k} \! \neq \! \infty$, $\widetilde{
\mathrm{U}}_{\mu}(z)/\widetilde{\mathrm{V}}_{\mu}(z)$ $(= \! 
\widetilde{\pmb{\mathrm{R}}}_{\mu}(z)/\pmb{\pi}^{n}_{k}(z))$ is the 
MPA of type $((n \! - \! 1)K \! + \! k \! - \! 1,(n \! - \! 1)K \! + \! k)$ 
for the Markov-Stieltjes transform, that is, it is the---unique---(proper) 
rational function with $\operatorname{deg}(\widetilde{\mathrm{U}}_{
\mu}(z)) \! = \! (n \! - \! 1)K \! + \! k \! - \! 1$, $\operatorname{deg}
(\widetilde{\mathrm{V}}_{\mu}(z)) \! = \! (n \! - \! 1)K \! + \! k$, and 
$(\widetilde{\mathrm{U}}_{\mu}(z),\widetilde{\mathrm{V}}_{\mu}(z))$ 
coprime interpolating $\mathrm{F}_{\mu}(z)$ and satisfying the 
interpolation conditions
\begin{gather}
\dfrac{\widetilde{\mathrm{U}}_{\mu}(z)}{\widetilde{\mathrm{V}}_{\mu}
(z)} \! - \! \sum_{j=0}^{2 \varkappa_{nk \tilde{k}_{q}}-1}c^{(q)}_{j}
(\alpha_{p_{q}})(z \! - \! \alpha_{p_{q}})^{j} \underset{z \to \alpha_{p_{q}}}{=} 
\mathcal{O} \left((z \! - \! \alpha_{p_{q}})^{2 \varkappa_{nk \tilde{k}_{q}}} 
\right), \quad q \! = \! 1,2,\dotsc,\mathfrak{s} \! - \! 2, \label{mvssfin6} \\
\dfrac{\widetilde{\mathrm{U}}_{\mu}(z)}{\widetilde{\mathrm{V}}_{\mu}
(z)} \! - \! \sum_{j=0}^{2 \varkappa_{nk}-1}c^{(\mathfrak{s})}_{j}
(\alpha_{k})(z \! - \! \alpha_{k})^{j} \underset{z \to \alpha_{k}}{=} 
\mathcal{O} \left((z \! - \! \alpha_{k})^{2 \varkappa_{nk}} \right), 
\label{mvssfin7} \\
\dfrac{\widetilde{\mathrm{U}}_{\mu}(z)}{\widetilde{\mathrm{V}}_{\mu}(z)} 
\! - \! \sum_{j=1}^{2 \varkappa^{\infty}_{nk \tilde{k}_{\mathfrak{s}-1}}}
c^{(\infty)}_{j}z^{-j} \underset{z \to \alpha_{p_{\mathfrak{s}-1}} =\infty}{=} 
\mathcal{O} \left(z^{-(2 \varkappa^{\infty}_{nk \tilde{k}_{\mathfrak{s}-1}}
+1)} \right). \label{mvssfin8} 
\end{gather}
\begin{eeee} \label{mvssremfin} 
\textsl{Note that Equations~\eqref{mvssfin6}, \eqref{mvssfin7}, 
and~\eqref{mvssfin8} give rise to $2 \sum_{q=1}^{\mathfrak{s}-2} 
\varkappa_{nk \tilde{k}_{q}}$, $2 \varkappa_{nk}$, and $2 \varkappa^{
\infty}_{nk \tilde{k}_{\mathfrak{s}-1}}$ interpolation conditions, respectively, 
for a combined total of $2(\sum_{q=1}^{\mathfrak{s}-2} \varkappa_{nk 
\tilde{k}_{q}} \! + \! \varkappa^{\infty}_{nk \tilde{k}_{\mathfrak{s}-1}} \! 
+ \! \varkappa_{nk}) \! = \! 2((n \! - \! 1)K \! + \! k)$ conditions, which 
is precisely the number necessary in order to determine the coefficients 
$\lbrace \tilde{r}_{j},\tilde{t}_{j} \rbrace_{j=0}^{(n-1)K+k-1}$ (the `leading 
coefficient', $\tilde{t}_{(n-1)K+k}$ $(\neq \! 0)$, is determined independently, 
via a detailed study of the $z \! \to \! \alpha_{p_{\mathfrak{s}-1}} 
\! = \! \infty$ asymptotics of the corresponding monic 
{\rm MPC ORF}$)$.}\footnote{Strictly speaking, $\tilde{r}_{j} \! = \! 
\tilde{r}_{j}(n,k),\tilde{t}_{j} \! = \! \tilde{t}_{j}(n,k)$, $j \! = \! 0,1,\dotsc,
(n \! - \! 1)K \! + \! k \! - \! 1$, and $\tilde{t}_{(n-1)K+k} \! = \! 
\tilde{t}_{(n-1)K+k}(n,k)$.}
\end{eeee}
For $n \! \in \! \mathbb{N}$ and $k \! \in \! \lbrace 1,2,\dotsc,K \rbrace$ 
such that $\alpha_{p_{\mathfrak{s}}} \! := \! \alpha_{k} \! \neq \! \infty$, 
define the corresponding MPA \emph{error term} as follows:\footnote{In 
this context, \emph{error term} means that, for $n \! \in \! \mathbb{N}$ 
and $k \! \in \! \lbrace 1,2,\dotsc,K \rbrace$ such that $\alpha_{
p_{\mathfrak{s}}} \! := \! \alpha_{k} \! \neq \! \infty$, $\widetilde{
\pmb{\mathrm{E}}}_{\mu}(z) \! =_{z \to \alpha_{p_{q}}} \! \mathcal{O}
((z \! - \! \alpha_{p_{q}})^{2 \varkappa_{nk \tilde{k}_{q}}})$, $q \! = \! 1,2,
\dotsc,\mathfrak{s} \! - \! 2$, $\widetilde{\pmb{\mathrm{E}}}_{\mu}(z) \! 
=_{z \to \alpha_{k}} \! \mathcal{O}((z \! - \! \alpha_{k})^{2 \varkappa_{nk}})$, 
and $\widetilde{\pmb{\mathrm{E}}}_{\mu}(z) \! =_{z \to \alpha_{
p_{\mathfrak{s}-1}}:=\infty} \! \mathcal{O} \left(z^{-(2 
\varkappa^{\infty}_{nk \tilde{k}_{\mathfrak{s}-1}}+1)} \right)$.}
\begin{equation} \label{mvssfin9} 
\widetilde{\pmb{\mathrm{E}}}_{\mu}(z) \! := \! \dfrac{\widetilde{
\pmb{\mathrm{R}}}_{\mu}(z)}{\pmb{\pi}^{n}_{k}(z)} \! - \! 
\mathrm{F}_{\mu}(z);
\end{equation}
a calculation based on Equations~\eqref{mvssinf1} and~\eqref{mvssfin3} 
reveals that, in fact,
\begin{equation} \label{mvssfin10}
\widetilde{\pmb{\mathrm{E}}}_{\mu}(z) \! = \! \dfrac{1}{(z \! - \! \alpha_{k}) 
\pmb{\pi}^{n}_{k}(z)} \left((z \! - \! \alpha_{k}) \int_{\mathbb{R}} 
\dfrac{((\xi \! - \! \alpha_{k}) \pmb{\pi}^{n}_{k}(\xi))}{(\xi \! - \! 
\alpha_{k})(\xi \! - \! z)} \, \md \mu (\xi) \right).
\end{equation}
\subsection{Hyperelliptic Riemann Surfaces and Summary of Results} 
\label{subsec1.3} 
Having defined, heretofore, in considerable detail and with examples, 
the principal objects of study of this monograph, that is, for $n \! \in 
\! \mathbb{N}$ and $k \! \in \! \lbrace 1,2,\dotsc,K \rbrace$ such 
that $\alpha_{p_{\mathfrak{s}}} \! := \! \alpha_{k} \! = \! \infty$ (resp., 
$\alpha_{p_{\mathfrak{s}}} \! := \! \alpha_{k} \! \neq \! \infty)$, the 
monic MPC ORF, $\pmb{\pi}^{n}_{k}(z)$, $z \! \in \! \mathbb{C}$, the 
`norming constant', $\mu^{\infty}_{n,\varkappa_{nk}}(n,k)$ (resp., 
$\mu^{f}_{n,\varkappa_{nk}}(n,k))$, the MPC ORF, $\phi^{n}_{k}(z) \! 
:= \! \mu^{\infty}_{n,\varkappa_{nk}}(n,k) \pmb{\pi}^{n}_{k}(z)$ (resp., 
$\phi^{n}_{k}(z) \! := \! \mu^{f}_{n,\varkappa_{nk}}(n,k) \pmb{\pi}^{n}_{k}
(z))$, $z \! \in \! \mathbb{C}$, and the corresponding MPA 
error term, $\widehat{\pmb{\mathrm{E}}}_{\mu}(z)$ (resp., 
$\widetilde{\pmb{\mathrm{E}}}_{\mu}(z))$, $z \! \in \! 
\mathbb{C}$, it must be mentioned that the ultimate goal of this lecture 
note is to obtain precise, and uniform, asymptotics, in the double-scaling 
limit {}\footnote{Similar asymptotics can also be considered in the 
double-scaling limit $\mathscr{N},n \! \to \! \infty$ such that $0 \! < 
\! \hat{c}_{\downarrow} \! \leqslant \! \mathscr{N}/n \! \leqslant \! 
\hat{c}_{\uparrow} \! < \! +\infty$, for arbitrary positive constants 
$\hat{c}_{\downarrow}$ and $\hat{c}_{\uparrow}$.} $\mathscr{N},n 
\! \to \! \infty$ such that $\mathscr{N}/n \! = \! 1 \! + \! o(1)$, of 
$\pmb{\pi}^{n}_{k}(z)$, $z \! \in \! \mathbb{C}$, $\mu^{p}_{n,\varkappa_{nk}}
(n,k)$, $p \! \in \! \lbrace \infty,f \rbrace$, $\phi^{n}_{k}(z) \! := \! 
\mu^{p}_{n,\varkappa_{nk}}(n,k) \pmb{\pi}^{n}_{k}(z)$, $z \! \in \! 
\mathbb{C}$, and $\widehat{\pmb{\mathrm{E}}}_{\mu}(z)$ and 
$\widetilde{\pmb{\mathrm{E}}}_{\mu}(z)$, $z \! \in \! \mathbb{C}$. 
In order to follow through with the above-described asymptotic programme, 
however, a formulation of the monic MPC ORF problem, for an \emph{a priori} 
prescribed sequence of neither necessarily distinct nor bounded poles 
$\alpha_{1},\alpha_{2},\dotsc,\alpha_{K}$ $(K \! \in \! \mathbb{N}$ and 
finite) lying on the support of the orthogonality measure, is a seminal 
necessity; in fact, the present monograph on MPC ORFs (see \cite{a45} 
for a detailed account of the FPC ORF theory) serves to address both 
the `formulation problem' and the `asymptotic analysis problem' (in 
the double-scaling limit $\mathscr{N},n \! \to \! \infty$ such that 
$\mathscr{N}/n \! = \! 1 \! + \! o(1))$.

In order to address the latter two problems, though, a primer on compact 
Riemman Surfaces, which is the subject matter of the following 
Subsection~\ref{subsub1}, is an absolute necessity: the monic MPC ORF 
problem will be re-addressed in Subsection~\ref{subsub2}, subsequent 
to this mathematical prelude.
\begin{eeee} \label{rem2.a} 
\textsl{\textbf{Henceforth, the notations and definitions of 
Subsections~\ref{subsubsec1.2.1} and~\ref{subsubsec1.2.2} will be used 
extensively, with little, or no, explanation(s).}}
\end{eeee}
\subsubsection{Riemann Surfaces: Preliminaries} \label{subsub1} 
In this succinct, self-contained subsection, the basic elements associated 
with the construction of hyperelliptic and finite genus (compact) Riemann 
Surfaces are presented; for further details and proofs, see, for example, 
\cite{gsp,hmfik}.
\begin{eeee} \label{rem2.1} 
\textsl{The superscripts ${}^{\pm}$, and sometimes the subscripts ${}_{\pm}$, 
should not be confused with the subscripts ${}_{\pm}$ appearing in the 
various---scalar and matrix---{\rm RHPs} (this is a general comment which 
applies throughout the entire monograph, unless stated otherwise).}
\end{eeee}
Let $N \! \in \! \mathbb{N}$,\footnote{It is assumed that $N \! < \! \infty$.} 
and let $\varsigma_{j} \! \in \! \overline{\mathbb{R}} \setminus \lbrace 
\alpha_{1},\alpha_{2},\dotsc,\alpha_{K} \rbrace$, $j \! = \! 1,2,\dotsc,
2N \! + \! 2$, be such that $\varsigma_{i} \! \not= \! \varsigma_{j}$  
$\forall$ $i \! \not= \! j \! \in \! \lbrace 1,2,\dotsc,2N \! + \! 2 \rbrace$, 
and enumerated/ordered according to $-\infty \! < \! \varsigma_{1} \! 
< \! \varsigma_{2} \! < \! \cdots \! < \! \varsigma_{2N+2} \! < \! +\infty$. 
Let $R(z) \! := \! \prod_{j=1}^{N}(z \! - \! \varsigma_{2j-1})(z \! - \! 
\varsigma_{2j})$ be the unital polynomial of even degree, $\mathrm{deg}(R) 
\! = \! 2N \! + \! 2$ $(\mathrm{deg}(R) \! = \! 0$ $(\mathrm{mod} \, 2))$, 
whose simple roots are $\lbrace \varsigma_{j} \rbrace_{j=1}^{2N+2}$. Denote 
by $\mathscr{R}$ the hyperelliptic Riemann surface of genus $N$ defined by 
the equation $y^{2} \! = \! R(z)$ and realised as a two-sheeted branched 
(ramified) covering of the Riemann sphere such that its two sheets are 
two identical copies of $\mathbb{C}$ with branch cuts along the intervals 
$(\varsigma_{1},\varsigma_{2})$, $(\varsigma_{3},\varsigma_{4})$, 
$\dotsc$, $(\varsigma_{2N+1},\varsigma_{2N+2})$, and glued to each 
other `crosswise' along the opposite banks of the corresponding cuts 
$(\varsigma_{2j-1},\varsigma_{2j})$, $j \! = \! 1,2,\dotsc,N \! + \! 1$. 
Denote the two sheets of $\mathscr{R}$ by $\mathscr{R}^{+}$ (the 
first/upper sheet) and  $\mathscr{R}^{-}$ (the second/lower sheet): 
to indicate that $z$ lies on the first (resp., second) sheet, one writes 
$z^{+}$ (resp., $z^{-})$; of course, as points on the complex plane 
$\mathbb{C}$, $z^{+} \! = \! z^{-} \! = \! z$. For points $z$ on the 
first (resp., second) sheet $\mathscr{R}^{+}$ (resp., $\mathscr{R}^{-})$, 
one has that $z^{+} \! = \! (z,+(R(z))^{1/2})$ (resp., $z^{-} \! = \! 
(z,-(R(z))^{1/2}))$, where the single-valued branch for the square root is 
chosen such that $z^{-(N+1)}(R(z))^{1/2} \! \sim_{\mathscr{R}^{\pm} 
\ni z \to \infty} \! \pm 1$.

Let $\mathscr{E}_{j} \! := \! (\varsigma_{2j-1},\varsigma_{2j})$, 
$j \! = \! 1,2,\dotsc,N \! + \! 1$, and set $\mathscr{E} \! = \! 
\cup_{j=1}^{N+1} \mathscr{E}_{j}$ (note that $\mathscr{E}_{i} \cap 
\mathscr{E}_{j} \! = \! \varnothing$, $i \! \not= \! j \! \in \! \lbrace 
1,2,\dotsc,N \! + \! 1 \rbrace)$. Denote by $\mathscr{E}_{j}^{+}$ 
$(\subset \mathscr{R}^{+})$ (resp., $\mathscr{E}_{j}^{-}$ $(\subset 
\mathscr{R}^{-}))$ the upper (resp., lower) bank of the interval 
$\mathscr{E}_{j}$, $j \! = \! 1,2,\dotsc,N \! + \! 1$, forming 
$\mathscr{E}$, and oriented in accordance with the orientation 
of $\mathscr{E}$ as the boundary of $\mathbb{C} \setminus 
\mathscr{E}$, that is, the domain $\mathbb{C} \setminus \mathscr{E}$ 
is on the left as one proceeds along the upper bank of the $j$th 
interval {}from $\varsigma_{2j-1}$ to the point $\varsigma_{2j}$ 
and back along the lower bank from $\varsigma_{2j}$ to 
$\varsigma_{2j-1}$; thus, $\mathscr{E}_{j}^{\pm} \! := \! 
(\varsigma_{2j-1},\varsigma_{2j})^{\pm}$, $j \! = \! 1,2,\dotsc,N \! 
+ \! 1$, are two (identical) copies of $(\varsigma_{2j-1},\varsigma_{2j}) 
\subset \mathbb{R}$ `lifted' to $\mathscr{R}^{\pm}$. Set $\Gamma \! 
:= \! \cup_{j=1}^{N+1} \Gamma_{j}$ $(\subset \mathscr{R})$, where 
$\Gamma_{j} \! := \! \mathscr{E}_{j}^{+} \cup \mathscr{E}_{j}^{-}$, 
$j \! = \! 1,2,\dotsc,N \! + \! 1$ $(\Gamma \! = \! \mathscr{E}^{+} 
\cup \mathscr{E}^{-})$. {}From the above construction, it is clear 
that $\mathscr{R} \! = \! \mathscr{R}^{+} \cup \mathscr{R}^{-} \cup 
\Gamma$; furthermore, the canonical projection of $\Gamma$ onto 
$\mathbb{C}$ $(\boldsymbol{\mathrm{pr}} \colon \mathscr{R} \! \to 
\! \mathbb{C})$ is $\mathscr{E}$, that is, $\boldsymbol{\mathrm{pr}}
(\Gamma) \! = \! \mathscr{E}$ (also, $\boldsymbol{\mathrm{pr}}
(\mathscr{R}^{+}) \! = \! \boldsymbol{\mathrm{pr}}(\mathscr{R}^{-}) 
\! = \! \mathbb{C} \setminus \mathscr{E}$, or, alternately, 
$\boldsymbol{\mathrm{pr}}(z^{+}) \! = \! \boldsymbol{\mathrm{pr}}
(z^{-}) \! = \! z)$. One moves in the `positive $(+)$' (resp., `negative 
$(-)$') direction along the (closed) contour $\Gamma \subset 
\mathscr{R}$ if the domain $\mathscr{R}^{+}$ is on the left (resp., 
right) and the domain $\mathscr{R}^{-}$ is on the right (resp., left): 
the corresponding notation is $\Gamma^{+}$ (resp., $\Gamma^{-})$.

One takes the first $N$ contours among the (closed) contours 
$\Gamma_{j}$ for basis $\boldsymbol{\alpha}$-cycles $\lbrace 
\boldsymbol{\alpha}_{j}, \, j \! = \! 1,2,\dotsc,N \rbrace$ and then 
supplements this in the standard way with $\boldsymbol{\beta}$-cycles 
$\lbrace \boldsymbol{\beta}_{j}, \, j \! = \!1,2,\dotsc,N \rbrace$ so 
that the \emph{intersection matrix} has the canonical form 
$\boldsymbol{\alpha}_{i} \circ \boldsymbol{\alpha}_{j} \! = \! 
\boldsymbol{\beta}_{i} \circ \boldsymbol{\beta}_{j} \! = \! 0$ 
$\forall$ $i \! \not= \! j \! \in \! \lbrace 1,2,\dotsc,N \rbrace$, and 
$\boldsymbol{\alpha}_{i} \circ \boldsymbol{\beta}_{j} \! = \! \delta_{ij}$: 
the cycles $\lbrace \boldsymbol{\alpha}_{j},\boldsymbol{\beta}_{j} 
\rbrace$, $j \! = \! 1,2,\dotsc,N$, form the \emph{canonical 
$\boldsymbol{1}$-homology basis} on $\mathscr{R}$, namely, any 
cycle $\boldsymbol{\gamma}^{\ast} \subset \mathscr{R}$ is 
homologous to an integral linear combination of $\lbrace 
\boldsymbol{\alpha}_{j},\boldsymbol{\beta}_{j} \rbrace$, that is, 
$\boldsymbol{\gamma}^{\ast} \! = \! \sum_{j=1}^{N}(n_{j} 
\boldsymbol{\alpha}_{j} \! + \! m_{j} \boldsymbol{\beta}_{j})$, where 
$(n_{j},m_{j}) \! \in \! \mathbb{Z} \times \mathbb{Z}$, $j \! = \! 
1,2,\dotsc,N$. The $\boldsymbol{\alpha}$-cycles $\lbrace 
\boldsymbol{\alpha}_{j}, \, j \! = \! 1,2,\dotsc,N \rbrace$, in the 
present case, are the intervals $(\varsigma_{2j-1},\varsigma_{2j})$, 
$j \! = \! 1,2,\dotsc,N$, `going twice', that is, along the upper 
({}from $\varsigma_{2j-1}$ to $\varsigma_{2j})$ and lower ({}from 
$\varsigma_{2j}$ to $\varsigma_{2j-1})$ banks $(\alpha_{j} \! = \! 
\mathscr{E}_{j}^{+} \cup \mathscr{E}_{j}^{-}$, $j \! = \! 1,2,\dotsc,N)$, 
and the $\boldsymbol{\beta}$-cycles $\lbrace \boldsymbol{\beta}_{j}, 
\, j \! = \! 1,2,\dotsc,N \rbrace$ are as follows: the $j$th 
$\boldsymbol{\beta}$-cycle consists of the $\boldsymbol{\alpha}$-cycles 
$\boldsymbol{\alpha}_{i}$, $i \! = \! j \! + \! 1,j \! + \! 2,\dotsc,N$, 
and the cycles `linked' with them and consisting of (the gaps) 
$(\varsigma_{2j},\varsigma_{2j+1})$, $j \! = \! 1,2,\dotsc,N$, 
`going twice', that is, {}from $\varsigma_{2j}$ to $\varsigma_{2j+1}$ 
on the first sheet and in the reverse direction on the second sheet. 
For an arbitrary holomorphic Abelian differential (one-form) 
$\boldsymbol{\omega}$ on $\mathscr{R}$, the function $\int^{z} 
\boldsymbol{\omega}$ is defined uniquely modulo its 
$\boldsymbol{\alpha}$- and $\boldsymbol{\beta}$-periods, 
$\oint_{\boldsymbol{\alpha}_{j}} \boldsymbol{\omega}$ and 
$\oint_{\boldsymbol{\beta}_{j}} \boldsymbol{\omega}$, $j \! = \! 
1,2,\dotsc,N$, respectively. It is well known that the canonical 
$\boldsymbol{1}$-homology basis $\lbrace \boldsymbol{\alpha}_{j},
\boldsymbol{\beta}_{j} \rbrace$, $j \! = \! 1,2,\dotsc,N$, constructed 
above `generates', on $\mathscr{R}$, the corresponding 
$\boldsymbol{\alpha}$-normalised basis of holomorphic Abelian 
differentials (one-forms) $\lbrace \omega_{1},\omega_{2},\dotsc,
\omega_{N} \rbrace$, where $\omega_{j} \! := \! \sum_{i=1}^{N}c_{ji}
(R(z))^{-1/2}z^{N-i} \, \md z$, $j \! = \! 1,2,\dotsc,N$, $c_{ji} \! \in 
\! \mathbb{C}$, and $\oint_{\boldsymbol{\alpha}_{i}} \omega_{j} \! 
= \! \delta_{ij}$, $i,j \! = \! 1,2,\dotsc,N$: the associated $N \times 
N$ matrix of $\boldsymbol{\beta}$-periods, $\tau \! = \! (\tau_{ij})_{i,j
=1,2,\dotsc,N} \! := \! \left(\oint_{\boldsymbol{\beta}_{j}} \omega_{i} 
\right)_{i,j=1,2,\dotsc,N}$, is a \emph{Riemann matrix}, that is, it is 
symmetric, pure imaginary, and $-\mi \tau$ is positive definite; 
moreover, $\tau$ is non-degenerate.

Fix the `standard basis' $\boldsymbol{e}_{1},\boldsymbol{e}_{2},\dotsc,
\boldsymbol{e}_{N}$ in $\mathbb{R}^{N}$, that is, $(\boldsymbol{e}_{i})_{j}
\! = \! \delta_{ij}$, $i,j \! = \! 1,2,\dotsc,N$:\footnote{These standard 
basis vectors should be viewed as column vectors.} the vectors 
$\boldsymbol{e}_{1},\boldsymbol{e}_{2},\dotsc,\boldsymbol{e}_{N},
\tau \boldsymbol{e}_{1},\linebreak[4]
\tau\boldsymbol{e}_{2},\dotsc,\tau \boldsymbol{e}_{N}$ are 
linearly independent over $\mathbb{R}$, and form a `basis' in 
$\mathbb{C}^{N}$. The quotient space $\mathbb{C}^{N}/\lbrace 
\mathrm{N} \! + \! \tau \mathrm{M} \rbrace$, $(\mathrm{N},\mathrm{M}) 
\! \in \! \mathbb{Z}^{N} \times \mathbb{Z}^{N}$, where $\mathbb{Z}^{N} 
\! := \! \lbrace \mathstrut (m_{1},m_{2},\dotsc,m_{N}); \, m_{j} \! \in \! 
\mathbb{Z}, \, j \! = \! 1,2,\dotsc,N \rbrace$, is a $2N$-dimensional real 
torus $\mathbb{T}^{2N}$, and is referred to as the \emph{Jacobi variety}, 
symbolically $\mathrm{Jac}(\mathscr{R})$, of the two-sheeted (hyperelliptic) 
Riemann surface $\mathscr{R}$ of genus $N$. Let $z_{0}^{\ast}$ be a fixed 
point in $\mathscr{R}$. A vector-valued function $\boldsymbol{\mathscr{A}}
(z) \! = \! (\mathscr{A}_{1}(z),\mathscr{A}_{2}(z),\dotsc,\mathscr{A}_{N}(z)) 
\! \in \! \mathrm{Jac}(\mathscr{R})$ with co-ordinates $\mathscr{A}_{j}(z) \! 
\equiv \! \int_{z_{0}^{\ast}}^{z} \omega_{j}$, $j \! = \! 1,2,\dotsc,N$, where, 
hereafter, unless stated otherwise and/or where confusion may arise, 
$\equiv$ denotes `congruence modulo the period lattice', defines the 
\emph{Abel map} $\boldsymbol{\mathscr{A}} \colon \mathscr{R} \! \to \! 
\mathrm{Jac}(\mathscr{R})$. The unordered set of points $z_{1},z_{2},
\dotsc,z_{N}$, with $z_{j} \! \in \! \mathscr{R}$, form the $N$th symmetric 
power of $\mathscr{R}$, symbolically $\mathscr{R}^{N}_{\mathrm{symm}}$. 
The vector function $\boldsymbol{\mathfrak{U}} \! = \! (\mathfrak{U}_{1},
\mathfrak{U}_{2},\dotsc,\mathfrak{U}_{N})$ with co-ordinates $\mathfrak{U}_{j} 
\! = \! \mathfrak{U}_{j}(z_{1},z_{2},\dotsc,z_{N}) \! \equiv \! \sum_{i=1}^{N} 
\mathscr{A}_{j}(z_{i}) \! \equiv \! \sum_{i=1}^{N} \int_{z_{0}^{\ast}}^{
z_{i}} \omega_{j}$, $j \! = \! 1,2,\dotsc,N$, that is, $(z_{1},z_{2},\dotsc,
z_{N}) \! \to \! (\sum_{i=1}^{N} \int_{z_{0}^{\ast}}^{z_{i}} \omega_{1},
\sum_{i=1}^{N} \int_{z_{0}^{\ast}}^{z_{i}} \omega_{2},\dotsc,
\sum_{i=1}^{N} \int_{z_{0}^{\ast}}^{z_{i}} \omega_{N})$, is also referred 
to as the \emph{Abel map}, $\boldsymbol{\mathfrak{U}} \colon 
\mathscr{R}^{N}_{\mathrm{symm}} \! \to \! \mathrm{Jac}(\mathscr{R})$. 
The \emph{dissected} Riemann surface, symbolically $\mathscr{R}^{\ast}$, 
is obtained {}from $\mathscr{R}$ by `cutting' (canonical dissection) 
along the cycles of the canonical $\boldsymbol{1}$-homology basis 
$\boldsymbol{\alpha}_{j},\boldsymbol{\beta}_{j}$, $j \! = \! 1,2,\dotsc,N$, 
of the original surface, namely, $\mathscr{R}^{\ast} \! = \! \mathscr{R} 
\setminus (\cup_{j=1}^{N}(\boldsymbol{\alpha}_{j} \cup \boldsymbol{
\beta}_{j}))$; the surface $\mathscr{R}^{\ast}$ is not only connected, as 
one can `pass' {}from one sheet to the other `across' $\Gamma_{N+1}$, 
but also simply connected (a $4N$-sided polygon ($4N$-gon) of a 
canonical dissection of $\mathscr{R}$ associated with the given canonical 
$\boldsymbol{1}$-homology basis for $\mathscr{R})$. For a given vector 
$\vec{\boldsymbol{v}} \! = \! (\upsilon_{1},\upsilon_{2},\dotsc,\upsilon_{N}) 
\! \in \! \mathrm{Jac}(\mathscr{R})$, the problem of finding an unordered 
collection of points $z_{1},z_{2},\dotsc,z_{N}$, $z_{j} \! \in \! \mathscr{R}$, 
$j \! = \! 1,2,\dotsc,N$, for which $\mathfrak{U}_{i}(z_{1},z_{2},\dotsc,
z_{N}) \! \equiv \! \upsilon_{i}$, $i \! = \! 1,2,\dotsc,N$, is called the 
\emph{Jacobi inversion problem} for Abelian integrals: as is well known, 
the Jacobi inversion problem is always solvable, but not, in general, 
uniquely.

By a \emph{divisor} on the Riemann surface $\mathscr{R}$ is meant a 
formal`symbol' $\boldsymbol{\mathrm{d}} \! = \! z_{1}^{n_{f}(z_{1})}
z_{2}^{n_{f}(z_{2})} \cdots z_{m}^{n_{f}(z_{m})}$, where $z_{j} \! \in \! 
\mathscr{R}$ and $n_{f}(z_{j}) \! \in \! \mathbb{Z}$, $j \! = \! 1,2,
\dotsc,m$: the number $\vert \boldsymbol{\mathrm{d}} \vert \! := \! 
\sum_{j=1}^{m}n_{f}(z_{j})$ is called the \emph{degree} of the divisor 
$\boldsymbol{\mathrm{d}}$: if $z_{i} \! \not= \! z_{j}$ $\forall$ $i \! 
\not= \! j \! \in \! \lbrace 1,2,\dotsc,m \rbrace$ and if $n_{f}(z_{j}) 
\!\geqslant \! 0$, $j \! = \! 1,2,\dotsc,m$, then the divisor 
$\boldsymbol{\mathrm{d}}$ is said to be \emph{integral}. Let $g$ be 
a meromorphic function defined on $\mathscr{R}$: for an arbitrary 
point $a \! \in \! \mathscr{R}$, one denotes by $n_{g}(a)$ (resp., 
$p_{g}(a))$ the multiplicity of the zero (resp., pole) of the function 
$g$ at this point if $a$ is a zero (resp., pole), and sets $n_{g}(a) \! 
= \! 0$ (resp., $p_{g}(a) \! = \! 0)$ otherwise; thus, $n_{g}(a),p_{g}(a) 
\! \geqslant \! 0$. To a meromorphic function $g$ on $\mathscr{R}$, 
one assigns the divisor of zeros and poles of this function as $(g) 
\! = \! z_{1}^{n_{g}(z_{1})}z_{2}^{n_{g}(z_{2})} \cdots z_{l_{1}}^{n_{g}
(z_{l_{1}})} \lambda_{1}^{-p_{g}(\lambda_{1})} \lambda_{2}^{-p_{g}
(\lambda_{2})} \cdots \lambda_{l_{2}}^{-p_{g}(\lambda_{l_{2}})}$, 
where $z_{i},\lambda_{j} \! \in \! \mathscr{R}$, $i \! = \! 1,2,\dotsc,
l_{1}$, $j \! = \! 1,2,\dotsc,l_{2}$, are the zeros and poles of $g$ on 
$\mathscr{R}$, and $n_{g}(z_{i}),p_{g}(\lambda_{j}) \! \geqslant \! 0$ 
are their multiplicities: these divisors are said to be \emph{principal}.

Associated with the Riemann matrix of $\boldsymbol{\beta}$-periods, 
$\tau$, is the \emph{Riemann theta function}, defined by
\begin{equation*}
\boldsymbol{\theta}(z;\tau) \! =: \! \boldsymbol{\theta}(z) \! = \! 
\sum_{m \in \mathbb{Z}^{N}} \me^{2 \pi \mi (m,z)+\mi \pi 
(m,\tau m)}, \quad z \! \in\! \mathbb{C}^{N},
\end{equation*}
where $(\boldsymbol{\cdot},\boldsymbol{\cdot})$ denotes 
the---real---Euclidean scalar product,\footnote{For $\mathbf{A} \! = \! 
(A_{1},A_{2},\dotsc,A_{N}) \! \in \! \mathbb{E}^{N}$ and $\mathbf{B} \! 
= \! (B_{1},B_{2},\dotsc,B_{N}) \! \in \! \mathbb{E}^{N}$, $(A,B) \! := \! 
\sum_{j=1}^{N}A_{j}B_{j}$.} with the following evenness and (quasi-) 
periodicity properties,
\begin{equation*}
\boldsymbol{\theta}(-z) \! = \! \boldsymbol{\theta}(z), \qquad 
\boldsymbol{\theta}(z \! + \! e_{j}) \! = \! \boldsymbol{\theta}(z), \qquad 
\mathrm{and} \qquad \boldsymbol{\theta}(z \! \pm \! \tau_{j}) \! = \! 
\me^{\mp 2 \pi \mi z_{j}-\mi \pi \tau_{jj}} \boldsymbol{\theta}(z),
\end{equation*}
where $e_{j}$ is the standard (basis) column vector in $\mathbb{C}^{N}$ 
with $1$ in the $j$th entry and $0$ elsewhere, and $\tau_{j} \! := \! 
\tau e_{j}$ $(\in \! \mathbb{C}^{N})$, $j \! = \! 1,2,\dotsc,N$.

It turns out that, for the analysis in this monograph, the following 
multi-valued functions are essential:
\begin{enumerate}
\item[$\boldsymbol{\bullet}$] $(\overset{\ast}{R}(z))^{1/2} \! 
:= \! (\prod_{j=0}^{N}(z \! - \! \overset{\ast}{b}_{j})(z \! - \! 
\overset{\ast}{a}_{j+1}))^{1/2}$, $\ast \! \in \! \lbrace \, \hat{} \, ,\, 
\tilde{} \, \rbrace$, where, with the identification $\overset{\ast}{a}_{N+1} 
\! \equiv \! \overset{\ast}{a}_{0}$ (as points on the complex sphere, 
$\overline{\mathbb{C}})$ and with the point at infinity lying on the (open) 
interval $(\overset{\ast}{a}_{0},\overset{\ast}{b}_{0})$, $-\infty \! < \! 
\overset{\ast}{a}_{0} \! < \! \overset{\ast}{b}_{0} \! < \! \overset{\ast}{a}_{1} 
\! < \! \overset{\ast}{b}_{1} \! < \! \cdots \! < \! \overset{\ast}{a}_{N} 
\! < \! \overset{\ast}{b}_{N} \! < \! +\infty$ (see 
Figure~\ref{drawone}).\footnote{It is important to point out that the 
multi-valued function $(\hat{R}(z))^{1/2}$ (resp., $(\tilde{R}(z))^{1/2})$ 
is associated with the asymptotic analysis, in the double-scaling limit 
$\mathscr{N},n \! \to \! \infty$ such that $\mathscr{N}/n \! = \! 1 \! 
+ \! o(1)$, corresponding to the case $n \! \in \! \mathbb{N}$ and $k \! 
\in \! \lbrace 1,2,\dotsc,K \rbrace$ such that $\alpha_{p_{\mathfrak{s}}} 
\! := \! \alpha_{k} \! = \! \infty$ (resp., $\alpha_{p_{\mathfrak{s}}} \! 
:= \! \alpha_{k} \! \neq \! \infty)$; see Subsection~\ref{subsub2}.}
\begin{figure}[tbh]
\begin{center}
\vspace{0.35cm}
\begin{pspicture}(0,0)(12,3)
\psset{xunit=1cm,yunit=1cm}
\psline[linewidth=0.9pt,linestyle=solid,linecolor=black]{o-o}(0.5,2)(2,2)
\psline[linewidth=0.9pt,linestyle=solid,linecolor=black]{o-o}(3,2)(4.5,2)
\psline[linewidth=0.9pt,linestyle=solid,linecolor=black]{o-o}(6.5,2)(8,2)
\psline[linewidth=0.9pt,linestyle=solid,linecolor=black]{o-o}(10,2)(11.5,2)
\psline[linewidth=0.7pt,linestyle=dotted,linecolor=darkgray](4.65,2)(6.35,2)
\psline[linewidth=0.7pt,linestyle=dotted,linecolor=darkgray](8.15,2)(9.85,2)
\rput(0.5,1.7){\makebox(0,0){$\overset{\ast}{a}_{0}$}}
\rput(0.5,0.9){\makebox(0,0){$\overset{\ast}{a}_{N+1}$}}
\rput{90}(0.5,1.3){\makebox(0,0){$\equiv$}}
\rput(1.25,2){\makebox(0,0){$\pmb{\times}$}}
\rput(1.25,2.3){\makebox(0,0){$\infty$}}
\rput(2,1.7){\makebox(0,0){$\overset{\ast}{b}_{0}$}}
\rput(3,1.7){\makebox(0,0){$\overset{\ast}{a}_{1}$}}
\rput(4.5,1.7){\makebox(0,0){$\overset{\ast}{b}_{1}$}}
\rput(6.5,1.7){\makebox(0,0){$\overset{\ast}{a}_{j}$}}
\rput(8,1.7){\makebox(0,0){$\overset{\ast}{b}_{j}$}}
\rput(10,1.7){\makebox(0,0){$\overset{\ast}{a}_{N}$}}
\rput(11.5,1.7){\makebox(0,0){$\overset{\ast}{b}_{N}$}}
\end{pspicture}
\end{center}
\vspace{-0.95cm}
\caption{Union of open intervals in the complex $z$-plane 
$(\ast \! \in \! \lbrace \, \hat{} \, , \, \tilde{} \, \rbrace)$.} 
\label{drawone}
\end{figure}
\end{enumerate}
The functions $\hat{R}(z)$ and $\tilde{R}(z)$, respectively, are unital 
polynomials of even degree $(\deg (\hat{R}(z)) \! = \! \deg (\tilde{R}(z)) 
\! = \! 2(N \! + \! 1))$ whose simple roots are $\lbrace \hat{b}_{j-1},
\hat{a}_{j} \rbrace_{j=1}^{N+1}$ and $\lbrace \tilde{b}_{j-1},\tilde{a}_{j} 
\rbrace_{j=1}^{N+1}$. The basic ingredients associated with the 
construction of the hyperelliptic Riemann surfaces of genus $N$ 
corresponding, respectively, to the multi-valued functions $y^{2} \! 
= \! \hat{R}(z)$ and $y^{2} \! = \! \tilde{R}(z)$ was given above. One 
now uses the above construction, but particularised to the cases 
of the polynomials $\hat{R}(z)$ and $\tilde{R}(z)$, to arrive at the 
following constructions.
\begin{enumerate}
\item[\shadowbox{$\sqrt{\hat{R}(z)}$}] This discussion applies to the 
case $n \! \in \! \mathbb{N}$ and $k \! \in \! \lbrace 1,2,\dotsc,K 
\rbrace$ such that $\alpha_{p_{\mathfrak{s}}} \! := \! \alpha_{k} \! 
= \! \infty$ (cf. Subsection~\ref{subsubsec1.2.1}). Let $\hat{\mathcal{Y}}$ 
denote the two-sheeted Riemann surface of genus $N$ associated with 
$y^{2} \! = \! \hat{R}(z)$, where $\hat{R}(z)$ is characterised above: the 
first/upper (resp., second/lower) sheet of $\hat{\mathcal{Y}}$ is denoted 
by $\hat{\mathcal{Y}}^{+}$ (resp., $\hat{\mathcal{Y}}^{-})$, points on 
the first/upper (resp., second/lower) sheet are represented as $z^{+} 
\! := \! (z,+(\hat{R}(z))^{1/2})$ (resp., $z^{-} \! := \! (z,-(\hat{R}(z))^{1/2}))$, 
where, as points on the plane $\mathbb{C}$, $z^{+} \! = \! z^{-} \! = \! 
z$, and the single-valued branch for the square root of the multi-valued 
function $(\hat{R}(z))^{1/2}$ is chosen such that $z^{-(N+1)}
(\hat{R}(z))^{1/2} \! \sim_{\hat{\mathcal{Y}}^{\pm} \ni z \to 
\alpha_{k}} \! \pm 1$. $\hat{\mathcal{Y}}$ is realised as a 
(two-sheeted) ramified covering of the Riemann sphere such that its two 
sheets are two identical copies of $\mathbb{C}$ with branch cuts (slits) 
along the intervals $(\hat{a}_{0},\hat{b}_{0}),(\hat{a}_{1},\hat{b}_{1}),\dotsc,
(\hat{a}_{N},\hat{b}_{N})$ and glued together along $\cup_{j=1}^{N+1}
(\hat{a}_{j-1},\hat{b}_{j-1})$ in such a way that the cycles 
$\hat{\boldsymbol{\alpha}}_{0}$ and $\lbrace \hat{\boldsymbol{\alpha}}_{j},
\hat{\boldsymbol{\beta}}_{j} \rbrace$, $j \! = \! 1,2,\dotsc,N$, where the 
latter forms the canonical $\mathbf{1}$-homology basis for $\hat{\mathcal{Y}}$, 
are characterised by the fact that $\hat{\boldsymbol{\alpha}}_{j}$, $j \! = \! 0,
1,\dotsc,N$, lie on $\hat{\mathcal{Y}}^{+}$, and $\hat{\boldsymbol{\beta}}_{j}$, 
$j \! = \! 1,2,\dotsc,N$, pass {}from $\hat{\mathcal{Y}}^{+}$ (starting {}from 
the slit $(\hat{a}_{j},\hat{b}_{j}))$, through the slit $(\hat{a}_{0},\hat{b}_{0})$ 
to $\hat{\mathcal{Y}}^{-}$, and back again to $\hat{\mathcal{Y}}^{+}$ through 
the slit $(\hat{a}_{j},\hat{b}_{j})$ (see Figure~\ref{drawtwo}).
\begin{figure}[tbh]
\begin{center}
\vspace{0.45cm}

\end{center}
\vspace{-0.55cm}
\caption{The Riemann surface $\hat{\mathcal{Y}}$ of $y^{2} \! = \! 
\prod_{j=0}^{N}(z \! - \! \hat{b}_{j})(z \! - \! \hat{a}_{j+1})$. The solid 
(resp., dashed) lines are on the first/upper (resp., second/lower) sheet 
of $\hat{\mathcal{Y}}$, denoted $\hat{\mathcal{Y}}^{+}$ (resp., 
$\hat{\mathcal{Y}}^{-})$.}
\label{drawtwo}
\end{figure}

\hspace*{0.50cm}
The canonical $\mathbf{1}$-homology basis $\lbrace \hat{
\boldsymbol{\alpha}}_{j},\hat{\boldsymbol{\beta}}_{j} \rbrace$, 
$j \! = \! 1,2,\dotsc,N$, generates, on $\hat{\mathcal{Y}}$, the 
corresponding $\hat{\boldsymbol{\alpha}}$-norma\-l\-i\-s\-e\-d 
basis of holomorphic Abelian differentials (one-forms) $\lbrace 
\hat{\omega}_{1},\hat{\omega}_{2},\dotsc,\hat{\omega}_{N} 
\rbrace$, where $\hat{\omega}_{j} \! := \! \sum_{i=1}^{N} 
\hat{c}_{ji}(\hat{R}(z))^{-1/2} \linebreak[4]
\pmb{\cdot} z^{N-i} \, \md z$, $j \! = \! 1,2,\dotsc,N$, 
$\hat{c}_{ji} \! \in \! \mathbb{C}$, and $\oint_{\hat{\boldsymbol{
\alpha}}_{i}} \hat{\omega}_{j} \! = \! \delta_{ij}$, $i,j \! = \! 1,2,
\dotsc,N$. Let $\hat{\boldsymbol{\omega}} \! := \! (\hat{\omega}_{1},
\hat{\omega}_{2},\dotsc,\hat{\omega}_{N})$ denote the basis of 
holomorphic one-forms on $\hat{\mathcal{Y}}$ normalised as 
above with the associated $N \! \times \! N$ Riemann matrix 
of $\hat{\boldsymbol{\beta}}$-periods, $\hat{\tau} \! = \! 
(\hat{\tau})_{i,j=1,2,\dotsc,N} \! := \! (\oint_{\hat{\boldsymbol{\beta}}_{j}} 
\hat{\omega}_{i})_{i,j=1,2,\dotsc,N}$: the Riemann matrix, $\hat{\tau}$, 
is symmetric and pure imaginary, $-\mi \hat{\tau}$ is positive definite, 
and $\det (\hat{\tau}) \! \not= \! 0$. For the holomorphic Abelian 
differential (one-form) $\hat{\boldsymbol{\omega}}$ defined above, 
choose $\hat{a}_{N+1}$, say, as the \emph{base point}, and set 
$\hat{\boldsymbol{u}} \colon \hat{\mathcal{Y}} \! \to \! 
\operatorname{Jac}(\hat{\mathcal{Y}})$ $(:= \! \mathbb{C}^{N}/\lbrace 
N_{0} \! + \! \hat{\tau}M_{0} \rbrace$, $(N_{0},M_{0}) \! \in \! 
\mathbb{Z}^{N} \! \times \! \mathbb{Z}^{N})$, $z \! \mapsto \! \hat{
\boldsymbol{u}}(z) \! := \! \int_{\hat{a}_{N+1}}^{z} \hat{\boldsymbol{
\omega}}$, where the integration {}from $\hat{a}_{N+1}$ to $z$ $(\in 
\hat{\mathcal{Y}})$ is taken along any path on $\hat{\mathcal{Y}}^{+}$.
\begin{eeee} \label{remchat} 
\textsl{{}From the representation $\hat{\omega}_{j} \! = \! \sum_{i=1}^{N} 
\hat{c}_{ji}(\hat{R}(z))^{-1/2}z^{N-i} \, \md z$, $j \! = \! 1,2,\dotsc,N$, 
and the normalisation condition $\oint_{\hat{\boldsymbol{\alpha}}_{i}} 
\hat{\omega}_{j} \! = \! \delta_{ij}$, $i,j \! = \! 1,2,\dotsc,N$, one shows 
that $\hat{c}_{i_{1}i_{2}}$, $i_{1},i_{2} \! = \! 1,2,\dotsc,N$, are derived {}from
\begin{equation} \label{E1}
\begin{pmatrix}
\hat{c}_{11} & \hat{c}_{12} & \dotsb & \hat{c}_{1N} \\
\hat{c}_{21} & \hat{c}_{22} & \dotsb & \hat{c}_{2N} \\
\vdots & \vdots & \ddots & \vdots \\
\hat{c}_{N1} & \hat{c}_{N2} & \dotsb & \hat{c}_{NN}
\end{pmatrix} \! =: \! \hat{\mathfrak{S}}^{-1},
\end{equation}
where
\begin{equation} \label{E2} 
\hat{\mathfrak{S}} \! = \! 
\begin{pmatrix}
\oint_{\hat{\boldsymbol{\alpha}}_{1}} \frac{\xi_{1}^{N-1} \md \xi_{1}}{
\sqrt{\smash[b]{\hat{R}(\xi_{1})}}} & \oint_{\hat{\boldsymbol{\alpha}}_{2}} 
\frac{\xi_{2}^{N-1} \md \xi_{2}}{\sqrt{\smash[b]{\hat{R}(\xi_{2})}}} & 
\dotsb & \oint_{\hat{\boldsymbol{\alpha}}_{N}} \frac{\xi_{N}^{N-1} 
\md \xi_{N}}{\sqrt{\smash[b]{\hat{R}(\xi_{N})}}} \\
\vdots & \vdots & \ddots & \vdots \\
\oint_{\hat{\boldsymbol{\alpha}}_{1}} \frac{\xi_{1} \md \xi_{1}}{
\sqrt{\smash[b]{\hat{R}(\xi_{1})}}} & \oint_{\hat{\boldsymbol{\alpha}}_{2}} 
\frac{\xi_{2} \md \xi_{2}}{\sqrt{\smash[b]{\hat{R}(\xi_{2})}}} & \dotsb & 
\oint_{\hat{\boldsymbol{\alpha}}_{N}} \frac{\xi_{N} \md \xi_{N}}{
\sqrt{\smash[b]{\hat{R}(\xi_{N})}}} \\
\oint_{\hat{\boldsymbol{\alpha}}_{1}} \frac{\md \xi_{1}}{\sqrt{\smash[b]{
\hat{R}(\xi_{1})}}} & \oint_{\hat{\boldsymbol{\alpha}}_{2}} \frac{\md 
\xi_{2}}{\sqrt{\smash[b]{\hat{R}(\xi_{2})}}} & \dotsb & 
\oint_{\hat{\boldsymbol{\alpha}}_{N}} \frac{\md \xi_{N}}{\sqrt{\smash[b]{
\hat{R}(\xi_{N})}}}
\end{pmatrix}.
\end{equation}
Using the multi-linearity property of the determinant, it follows, via a 
Vandermonde argument applied to Equation~\eqref{E2}, that
\begin{equation*}
\det (\hat{\mathfrak{S}}) \! = \! (-1)^{[N/2]+N(N+1)/2}(2 \mi)^{N} 
\int_{\hat{a}_{1}}^{\hat{b}_{1}} \int_{\hat{a}_{2}}^{\hat{b}_{2}} \dotsb 
\int_{\hat{a}_{N}}^{\hat{b}_{N}} \prod_{m=1}^{N}(\lvert \hat{R}_{m}
(\xi_{m}) \rvert)^{-1/2} \prod_{\substack{i,j=1\\j<i}}^{N}(\xi_{i} \! - 
\! \xi_{j}) \, \md \xi_{1} \, \md \xi_{2} \, \dotsb \, \md \xi_{N},
\end{equation*}
where $[N/2]$ denotes the integer part of $N/2$, and, for 
$\xi \! \in \! (\hat{a}_{j},\hat{b}_{j})$, $j \! = \! 1,2,\dotsc,N$,
\begin{align*}
(\lvert \hat{R}_{j}(\xi) \rvert)^{1/2} :=& \, (\xi \! - \! \hat{b}_{0})^{1/2}
(\xi \! - \! \hat{a}_{0})^{1/2}(\hat{b}_{j} \! - \! \xi)^{1/2}(\xi \! - \! 
\hat{a}_{j})^{1/2} \prod_{m=1}^{j-1}(\xi \! - \! \hat{b}_{m})^{1/2}
(\xi \! - \! \hat{a}_{m})^{1/2} \\
\times& \, \prod_{m^{\prime}=j+1}^{N}(\hat{b}_{m^{\prime}} 
\! - \! \xi)^{1/2}(\hat{a}_{m^{\prime}} \! - \! \xi)^{1/2},
\end{align*}
where all the square roots are positive and all the end-point integrations 
have removable $1/2$-root singularities: since $(\hat{a}_{i},\hat{b}_{i}) 
\cap (\hat{a}_{j},\hat{b}_{j}) \! = \! \varnothing$ $\forall$ $i \! \neq \! j 
\! \in \! \lbrace 1,2,\dotsc,N \rbrace$, it follows that, for $i \! < \! j$, 
$((\hat{a}_{i},\hat{b}_{i}) \! \ni)$ $\xi_{i} \! < \! \xi_{j}$ $(\in \! (\hat{a}_{j},
\hat{b}_{j}))$, which implies that $\prod_{\underset{j<i}{i,j=1}}^{N}
(\xi_{i} \! - \! \xi_{j}) \! > \! 0$$;$ moreover, since, for $\xi \! \in \! 
(\hat{a}_{j},\hat{b}_{j})$, $j \! = \! 1,2,\dotsc,N$, $(\lvert \hat{R}_{j}(\xi) 
\rvert)^{1/2} \! > \! 0$, that is, $(\lvert \hat{R}_{j}(\xi) \rvert)^{1/2}$ does 
not change sign, it follows that the integrand, $\prod_{m=1}^{N}(\lvert 
\hat{R}_{m}(\xi_{m}) \rvert)^{-1/2} \prod_{\underset{j<i}{i,j=1}}^{N}(\xi_{i} 
\! - \! \xi_{j})$, is positive (it, too, has constant sign) in the domain of 
definition of multiple integration; thus, $\det (\hat{\mathfrak{S}}) \! 
\neq \! 0$, namely, $\mathfrak{S}$ is invertible (non-singular). For a 
representation-independent proof of the fact that $\det (\hat{\mathfrak{S}}) 
\! \not= \! 0$, see, for example, Chapter~{\rm 10}, Section~{\rm 10}--{\rm 2}, 
of {\rm \cite{gsp}}}.
\end{eeee}

\hspace*{0.50cm}
Set (see Section~\ref{sec4}), for $z \! \in \! \mathbb{C}_{+}$, 
$\hat{\gamma}(z) \! := \! (\prod_{j=1}^{N+1}(z \! - \! \hat{b}_{j-1})
(z \! - \! \hat{a}_{j})^{-1})^{1/4}$, and, for $z \! \in \! \mathbb{C}_{-}$, 
$\hat{\gamma}(z) \! := \! -\mi (\prod_{j=1}^{N+1}(z \! - \! \hat{b}_{j-1})
(z \! - \! \hat{a}_{j})^{-1})^{1/4}$. It is shown in Lemma~\ref{lem4.4} 
that $\hat{\gamma}(z) \! =_{\hat{\mathcal{Y}}^{\pm} \ni z \to \alpha_{k}} 
\! (-\mi)^{(1 \mp 1)/2}(1 \! + \! \mathcal{O}(z^{-1}))$, $\hat{\gamma}
(z) \! =_{\hat{\mathcal{Y}}^{\pm} \ni z \to \alpha_{p_{q}}} 
\! (-\mi)^{(1 \mp 1)/2} \linebreak[4] 
\pmb{\cdot} \hat{\gamma}(\alpha_{p_{q}})(1 \! + \! \mathcal{O}(z \! 
- \! \alpha_{p_{q}}))$, $q \! = \! 1,2,\dotsc,\mathfrak{s} \! - \! 1$, where 
$\hat{\gamma}(\alpha_{p_{q}})$ is defined by Equation~\eqref{eqssabra1}, 
and
\begin{equation*}
\lbrace \hat{z}_{j}^{\pm} \rbrace_{j=1}^{N} \! := \! \lbrace \mathstrut 
z^{\pm} \! \in \! \hat{\mathcal{Y}}^{\pm}; \, (\hat{\gamma}(z) \! 
\mp \! (\hat{\gamma}(z))^{-1}) \vert_{z=z^{\pm}} \! = \! 0 \rbrace,
\end{equation*}
with $\hat{z}_{j}^{\pm} \! \in \! (\hat{a}_{j},\hat{b}_{j})^{\pm}$ 
$(\subset \! \hat{\mathcal{Y}}^{\pm})$, $j \! = \! 1,2,\dotsc,N$, where, 
as points on the plane, $\hat{z}_{j}^{+} \! = \! \hat{z}_{j}^{-} \! = \! 
\hat{z}_{j} \! \in \! (\hat{a}_{j},\hat{b}_{j})$, $j \! = \! 1,2,\dotsc,
N$.\footnote{More precisely, $\boldsymbol{\mathrm{pr}}(\hat{z}_{j}^{+}) 
\! = \! \boldsymbol{\mathrm{pr}}(\hat{z}_{j}^{-}) \! = \! \hat{z}_{j} 
\! \in \! (\hat{a}_{j},\hat{b}_{j})$, $j \! = \! 1,2,\dotsc,N$.}

\hspace*{0.50cm}
Corresponding to $\hat{\mathcal{Y}}$, define $\hat{\boldsymbol{d}} \! 
:= \!-\hat{\boldsymbol{K}} \! - \! \sum_{j=1}^{N} \int_{\hat{a}_{N+1}}^{
\hat{z}_{j}^{-}} \hat{\boldsymbol{\omega}}$ $(\in \! \mathbb{C}^{N})$, 
where $\hat{\boldsymbol{K}}$ is the associated vector of Riemann 
constants, and the integration {}from $\hat{a}_{N+1}$ to $\hat{z}_{j}^{-}$, 
$j \! = \! 1,2,\dotsc,N$, is taken along a fixed path in $\hat{\mathcal{Y}}^{-}$. 
It is shown in Chapter~VII of \cite{hmfik} that $\hat{\boldsymbol{K}} \! = \! 
\sum_{j=1}^{N} \int_{\hat{a}_{j}}^{\hat{a}_{N+1}} \hat{\boldsymbol{\omega}}$; 
furthermore, $\hat{\boldsymbol{K}}$ is a point of order $2$, that is, $2 
\hat{\boldsymbol{K}} \! = \! 0$ and $p \hat{\boldsymbol{K}} \! \not= \! 0$ 
for $0 \! < \! p \! < \! 2$. Recalling the definition of $\hat{\boldsymbol{
\omega}}$, and that $z^{-(N+1)}(\hat{R}(z))^{1/2} \! \sim_{\hat{
\mathcal{Y}}^{\pm }\ni z \to \alpha_{k}} \! \pm 1$, using the fact 
that $\hat{\boldsymbol{K}}$ is a point of order $2$, one arrives at
\begin{align*}
\hat{\boldsymbol{d}} =& \, -\hat{\boldsymbol{K}} \! - \! \sum_{j=1}^{N} 
\int_{\hat{a}_{N+1}}^{\hat{z}_{j}^{-}} \hat{\boldsymbol{\omega}} \! = \! 
\hat{\boldsymbol{K}} \! -\! \sum_{j=1}^{N} \int_{\hat{a}_{N+1}}^{\hat{z}_{j}^{-}} 
\hat{\boldsymbol{\omega}} \! = \! -\hat{\boldsymbol{K}} \! + \! \sum_{j=1}^{N} 
\int_{\hat{a}_{N+1}}^{\hat{z}_{j}^{+}} \hat{\boldsymbol{\omega}} \! = \! 
\hat{\boldsymbol{K}} \! + \! \sum_{j=1}^{N} \int_{\hat{a}_{N+1}}^{\hat{z}_{j}^{+}} 
\hat{\boldsymbol{\omega}} \\
=& \, -\sum_{j=1}^{N} \int_{\hat{a}_{j}}^{\hat{z}_{j}^{-}} \hat{\boldsymbol{
\omega}} \! = \! \sum_{j=1}^{N} \int_{\hat{a}_{j}}^{\hat{z}_{j}^{+}} 
\hat{\boldsymbol{\omega}}.
\end{align*}

\hspace*{0.50cm}
Associated with the Riemann matrix of $\hat{\boldsymbol{\beta}}$-periods, 
$\hat{\tau}$, is the Riemann theta function
\begin{equation} \label{eqrmthetainf} 
\hat{\boldsymbol{\theta}}(z;\hat{\tau}) \! =: \! \hat{\boldsymbol{\theta}}
(z) \! = \! \sum_{m \in \mathbb{Z}^{N}} \me^{2 \pi \mi (m,z)+\mi \pi 
(m,\hat{\tau}m)}, \quad z \! \in \! \mathbb{C}^{N};
\end{equation}
$\hat{\boldsymbol{\theta}}(z)$ has the following evenness and (quasi-) 
periodicity properties,
\begin{equation*}
\hat{\boldsymbol{\theta}}(-z) \! = \! \hat{\boldsymbol{\theta}}(z), \qquad 
\hat{\boldsymbol{\theta}}(z \! + \! e_{j}) \! = \! \hat{\boldsymbol{\theta}}(z), 
\qquad \text{and} \qquad \hat{\boldsymbol{\theta}}(z \! \pm \! \hat{\tau}_{j}) 
\! = \! \me^{\mp 2 \pi \mi z_{j}-\mi \pi \hat{\tau}_{jj}} \hat{\boldsymbol{\theta}}
(z),
\end{equation*}
where $\hat{\tau}_{j} \! := \! \hat{\tau}e_{j}$ $(\in \! \mathbb{C}^{N})$, 
$j \! = \! 1,2,\dotsc,N$.

\item[\shadowbox{$\sqrt{\tilde{R}(z)}$}] This discussion applies to the case 
$n \! \in \! \mathbb{N}$ and $k \! \in \! \lbrace 1,2,\dotsc,K \rbrace$ 
such that $\alpha_{p_{\mathfrak{s}}} \! := \! \alpha_{k} \! \neq \! \infty$ 
(cf. Subsection~\ref{subsubsec1.2.2}). Let $\tilde{\mathcal{Y}}$ denote 
the two-sheeted Riemann surface of genus $N$ associated with $y^{2} 
\! = \! \tilde{R}(z)$, where $\tilde{R}(z)$ is characterised above: the 
first/upper (resp., second/lower) sheet of $\tilde{\mathcal{Y}}$ is denoted 
by $\tilde{\mathcal{Y}}^{+}$ (resp., $\tilde{\mathcal{Y}}^{-})$, points on 
the first/upper (resp., second/lower) sheet are represented as $z^{+} 
\! := \! (z,+(\tilde{R}(z))^{1/2})$ (resp., $z^{-} \! := \! (z,-(\tilde{R}(z))^{1/2}))$, 
where, as points on the plane $\mathbb{C}$, $z^{+} \! = \! z^{-} \! = \! 
z$, and the single-valued branch for the square root of the (multi-valued) 
function $(\tilde{R}(z))^{1/2}$ is chosen such that $z^{-(N+1)}
(\tilde{R}(z))^{1/2} \! \sim_{\tilde{\mathcal{Y}}^{\pm} \ni z \to 
\alpha_{p_{\mathfrak{s}-1}}} \! \pm 1$ ($\alpha_{p_{\mathfrak{s}-1}} \! 
:= \! \infty$; cf. Subsection~\ref{subsubsec1.2.2}). $\tilde{\mathcal{Y}}$ 
is realised as a (two-sheeted) ramified covering of the Riemann sphere 
such that its two sheets are two identical copies of $\mathbb{C}$ with 
branch cuts along the intervals $(\tilde{a}_{0},\tilde{b}_{0}),(\tilde{a}_{1},
\tilde{b}_{1}),\dotsc,(\tilde{a}_{N},\tilde{b}_{N})$ and glued together 
along $\cup_{j=1}^{N+1}(\tilde{a}_{j-1},\tilde{b}_{j-1})$ in such a way 
that the cycles $\tilde{\boldsymbol{\alpha}}_{0}$ and $\lbrace 
\tilde{\boldsymbol{\alpha}}_{j},\tilde{\boldsymbol{\beta}}_{j} \rbrace$, 
$j \! = \! 1,2,\dotsc,N$, where the latter forms the canonical 
$\mathbf{1}$-homology basis for $\tilde{\mathcal{Y}}$, are characterised 
by the fact that $\tilde{\boldsymbol{\alpha}}_{j}$, $j \! = \! 0,1,\dotsc,N$, 
lie on $\tilde{\mathcal{Y}}^{+}$, and $\tilde{\boldsymbol{\beta}}_{j}$, $j \! 
= \! 1,2,\dotsc,N$, pass {}from $\tilde{\mathcal{Y}}^{+}$ (starting {}from the 
slit $(\tilde{a}_{j},\tilde{b}_{j}))$, through the slit $(\tilde{a}_{0},\tilde{b}_{0})$ 
to $\tilde{\mathcal{Y}}^{-}$, and back again to $\tilde{\mathcal{Y}}^{+}$ 
through the slit $(\tilde{a}_{j},\tilde{b}_{j})$ (see Figure~\ref{drawthree}).
\begin{figure}[tbh]
\begin{center}
\vspace{0.45cm}

\end{center}
\vspace{-0.55cm}
\caption{The Riemann surface $\tilde{\mathcal{Y}}$ of $y^{2} \! = 
\! \prod_{j=0}^{N}(z \! - \! \tilde{b}_{j})(z \! - \! \tilde{a}_{j+1})$. 
The solid (resp., dashed) lines are on the first/upper (resp., 
second/lower) sheet of $\tilde{\mathcal{Y}}$, denoted $\tilde{
\mathcal{Y}}^{+}$ (resp., $\tilde{\mathcal{Y}}^{-})$.}
\label{drawthree}
\end{figure}

\hspace*{0.50cm}
The canonical $\mathbf{1}$-homology basis $\lbrace 
\tilde{\boldsymbol{\alpha}}_{j},\tilde{\boldsymbol{\beta}}_{j} \rbrace$, 
$j \! = \! 1,2,\dotsc,N$, generates, on $\tilde{\mathcal{Y}}$, the 
corresponding $\tilde{\boldsymbol{\alpha}}$-norma\-l\-i\-s\-e\-d 
basis of holomorphic Abelian differentials (one-forms) $\lbrace 
\tilde{\omega}_{1},\tilde{\omega}_{2},\dotsc,\tilde{\omega}_{N} 
\rbrace$, where $\tilde{\omega}_{j} \! := \! \sum_{i=1}^{N} 
\tilde{c}_{ji}(\tilde{R}(z))^{-1/2} \linebreak[4]
\pmb{\cdot} z^{N-i} \, \md z$, $j \! = \! 1,2,
\dotsc,N$, $\tilde{c}_{ji} \! \in \! \mathbb{C}$, and $\oint_{\tilde{
\boldsymbol{\alpha}}_{i}} \tilde{\omega}_{j} \! = \! \delta_{ij}$, 
$i,j \! = \! 1,2,\dotsc,N$. Let $\tilde{\boldsymbol{\omega}} \! := \! 
(\tilde{\omega}_{1},\tilde{\omega}_{2},\dotsc,\tilde{\omega}_{N})$ 
denote the basis of holomorphic one-forms on $\tilde{\mathcal{Y}}$ 
normalised as above with the associated $N \! \times \! N$ Riemann 
matrix of $\tilde{\boldsymbol{\beta}}$-periods, $\tilde{\tau} \! = \! 
(\tilde{\tau})_{i,j=1,2,\dotsc,N} \! := \! (\oint_{\tilde{\boldsymbol{
\beta}}_{j}} \tilde{\omega}_{i})_{i,j=1,2,\dotsc,N}$: the Riemann matrix, 
$\tilde{\tau}$, is symmetric and pure imaginary, $-\mi \tilde{\tau}$ 
is positive definite, and $\det (\tilde{\tau}) \! \not= \! 0$. For the 
holomorphic Abelian differential (one-form) $\tilde{\boldsymbol{
\omega}}$ defined above, choose $\tilde{a}_{N+1}$, say, as the base 
point, and set $\tilde{\boldsymbol{u}} \colon \tilde{\mathcal{Y}} \! \to \! 
\operatorname{Jac}(\tilde{\mathcal{Y}})$ $(:= \! \mathbb{C}^{N}/\lbrace 
N_{0} \! + \! \tilde{\tau}M_{0} \rbrace$, $(N_{0},M_{0}) \! \in \! 
\mathbb{Z}^{N} \! \times \! \mathbb{Z}^{N})$, $z \! \mapsto \! 
\tilde{\boldsymbol{u}}(z) \! := \! \int_{\tilde{a}_{N+1}}^{z} 
\tilde{\boldsymbol{\omega}}$, where the integration {}from 
$\tilde{a}_{N+1}$ to $z$ $(\in \tilde{\mathcal{Y}})$ is taken along any 
path on $\tilde{\mathcal{Y}}^{+}$.
\begin{eeee} \label{remctilde} 
\textsl{{}From the representation $\tilde{\omega}_{j} \! = \! 
\sum_{i=1}^{N} \tilde{c}_{ji}(\tilde{R}(z))^{-1/2}z^{N-i} \, \md z$, 
$j \! = \! 1,2,\dotsc,N$, and the normalisation condition 
$\oint_{\tilde{\boldsymbol{\alpha}}_{i}} \tilde{\omega}_{j} \! = \! 
\delta_{ij}$, $i,j \! = \! 1,2,\dotsc,N$, one shows that $\tilde{c}_{i_{1}i_{2}}$, 
$i_{1},i_{2} \! = \! 1,2,\dotsc,N$, are derived {}from
\begin{equation} \label{O1}
\begin{pmatrix}
\tilde{c}_{11} & \tilde{c}_{12} & \dotsb & \tilde{c}_{1N} \\
\tilde{c}_{21} & \tilde{c}_{22} & \dotsb & \tilde{c}_{2N} \\
\vdots & \vdots & \ddots & \vdots \\
\tilde{c}_{N1} & \tilde{c}_{N2} & \dotsb & \tilde{c}_{NN}
\end{pmatrix} \! =: \! \tilde{\mathfrak{S}}^{-1},
\end{equation}
where
\begin{equation} \label{O2} 
\tilde{\mathfrak{S}} \! = \! 
\begin{pmatrix}
\oint_{\tilde{\boldsymbol{\alpha}}_{1}} \frac{\xi_{1}^{N-1} \md \xi_{1}}{
\sqrt{\smash[b]{\tilde{R}(\xi_{1})}}} & \oint_{\tilde{\boldsymbol{\alpha}}_{2}} 
\frac{\xi_{2}^{N-1} \md \xi_{2}}{\sqrt{\smash[b]{\tilde{R}(\xi_{2})}}} & 
\dotsb & \oint_{\tilde{\boldsymbol{\alpha}}_{N}} \frac{\xi_{N}^{N-1} 
\md \xi_{N}}{\sqrt{\smash[b]{\tilde{R}(\xi_{N})}}} \\
\vdots & \vdots & \ddots & \vdots \\
\oint_{\tilde{\boldsymbol{\alpha}}_{1}} \frac{\xi_{1} \md \xi_{1}}{
\sqrt{\smash[b]{\tilde{R}(\xi_{1})}}} & \oint_{\tilde{\boldsymbol{
\alpha}}_{2}} \frac{\xi_{2} \md \xi_{2}}{\sqrt{\smash[b]{\tilde{R}
(\xi_{2})}}} & \dotsb & \oint_{\tilde{\boldsymbol{\alpha}}_{N}} 
\frac{\xi_{N} \md \xi_{N}}{\sqrt{\smash[b]{\tilde{R}(\xi_{N})}}} \\
\oint_{\tilde{\boldsymbol{\alpha}}_{1}} \frac{\md \xi_{1}}{\sqrt{
\smash[b]{\tilde{R}(\xi_{1})}}} & \oint_{\tilde{\boldsymbol{\alpha}}_{2}} 
\frac{\md \xi_{2}}{\sqrt{\smash[b]{\tilde{R}(\xi_{2})}}} & \dotsb & 
\oint_{\tilde{\boldsymbol{\alpha}}_{N}} \frac{\md \xi_{N}}{
\sqrt{\smash[b]{\tilde{R}(\xi_{N})}}}
\end{pmatrix}.
\end{equation}
Arguing as in Remark~\ref{remchat}, one shows that
\begin{align*}
\det (\tilde{\mathfrak{S}}) =& \, (-1)^{[N/2]+N(N+1)/2}(2 \mi)^{N} 
\int_{\tilde{a}_{1}}^{\tilde{b}_{1}} \int_{\tilde{a}_{2}}^{\tilde{b}_{2}} 
\dotsb \int_{\tilde{a}_{N}}^{\tilde{b}_{N}} \prod_{m=1}^{N} \left(
(\xi_{m} \! - \! \tilde{b}_{0})^{1/2}(\xi_{m} \! - \! \tilde{a}_{0})^{1/2}
(\tilde{b}_{m} \! - \! \xi_{m})^{1/2}(\xi_{m} \! - \! \tilde{a}_{m})^{1/2} 
\right. \\
\times&\left. \, \prod_{l=1}^{m-1}(\xi_{m} \! - \! \tilde{b}_{l})^{1/2}
(\xi_{m} \! - \! \tilde{a}_{l})^{1/2} \prod_{l^{\prime}=m+1}^{N}
(\tilde{b}_{l^{\prime}} \! - \! \xi_{m})^{1/2}(\tilde{a}_{l^{\prime}} \! - \! 
\xi_{m})^{1/2} \right)^{-1} \prod_{\substack{i,j=1\\j<i}}^{N}(\xi_{i} \! - 
\! \xi_{j}) \, \md \xi_{1} \, \md \xi_{2} \, \dotsb \, \md \xi_{N} \! \neq \! 0,
\end{align*}
that is, $\tilde{\mathfrak{S}}$ is invertible (non-singular).}
\end{eeee}

\hspace*{0.50cm}
Set (see Section~\ref{sec4}), for $z \! \in \! \mathbb{C}_{+}$, 
$\tilde{\gamma}(z) \! := \! (\prod_{j=1}^{N+1}(z \! - \! \tilde{b}_{j-1})
(z \! - \! \tilde{a}_{j})^{-1})^{1/4}$, and, for $z \! \in \! \mathbb{C}_{-}$, 
$\tilde{\gamma}(z) \! := \! -\mi (\prod_{j=1}^{N+1}(z \! - \! \tilde{b}_{j-1})
(z \! - \! \tilde{a}_{j})^{-1})^{1/4}$. It is shown in Lemma~\ref{lem4.4} 
that $\tilde{\gamma}(z) \! =_{\tilde{\mathcal{Y}}^{\pm} \ni z \to 
\alpha_{p_{q}}} \! (-\mi)^{(1 \mp 1)/2} \tilde{\gamma}(\alpha_{p_{q}})
(1 \! + \! \mathcal{O}(z \! - \! \alpha_{p_{q}}))$, $q \! = \! 1,\dotsc,
\mathfrak{s} \! - \! 2,\mathfrak{s}$, where $\tilde{\gamma}
(\alpha_{p_{q}})$ is defined by Equation~\eqref{eqssabra2}, and 
$\tilde{\gamma}(z) \! =_{\tilde{\mathcal{Y}}^{\pm} \ni z \to 
\alpha_{p_{\mathfrak{s}-1}}} \! (-\mi)^{(1 \mp 1)/2}(1 \! + \! 
\mathcal{O}(z^{-1}))$, and
\begin{equation*}
\lbrace \tilde{z}_{j}^{\pm} \rbrace_{j=1}^{N} \! := \! \lbrace \mathstrut 
z^{\pm} \! \in \! \tilde{\mathcal{Y}}^{\pm}; \, (\tilde{\gamma}(z) \! 
\mp \! (\tilde{\gamma}(z))^{-1}) \vert_{z=z^{\pm}} \! = \! 0 \rbrace,
\end{equation*}
with $\tilde{z}_{j}^{\pm} \! \in \! (\tilde{a}_{j},\tilde{b}_{j})^{\pm}$ 
$(\subset \! \tilde{\mathcal{Y}}^{\pm})$, $j \! = \! 1,2,\dotsc,N$, 
where, as points on the plane, $\tilde{z}_{j}^{+} \! = \! \tilde{z}_{j}^{-} 
\! = \! \tilde{z}_{j} \! \in \! (\tilde{a}_{j},\tilde{b}_{j})$, $j \! = \! 1,2,
\dotsc,N$.\footnote{More precisely, $\boldsymbol{\mathrm{pr}}
(\tilde{z}_{j}^{+}) \! = \! \boldsymbol{\mathrm{pr}}(\tilde{z}_{j}^{-}) 
\! = \! \tilde{z}_{j} \! \in \! (\tilde{a}_{j},\tilde{b}_{j})$, $j \! = \! 1,
2,\dotsc,N$.}

\hspace*{0.50cm}
Corresponding to $\tilde{\mathcal{Y}}$, define $\tilde{\boldsymbol{d}} \! 
:= \!-\tilde{\boldsymbol{K}} \! - \! \sum_{j=1}^{N} \int_{\tilde{a}_{N+1}}^{
\tilde{z}_{j}^{-}} \tilde{\boldsymbol{\omega}}$ $(\in \! \mathbb{C}^{N})$, 
where $\tilde{\boldsymbol{K}}$, which is a point of order $2$, is the associated 
vector of Riemann constants, and the integration {}from $\tilde{a}_{N+1}$ to 
$\tilde{z}_{j}^{-}$, $j \! = \! 1,2,\dotsc,N$, is taken along a fixed path in 
$\tilde{\mathcal{Y}}^{-}$. It is shown in Chapter~VII of \cite{hmfik} that 
$\tilde{\boldsymbol{K}} \! = \! \sum_{j=1}^{N} \int_{\tilde{a}_{j}}^{
\tilde{a}_{N+1}} \tilde{\boldsymbol{\omega}}$. Recalling the definition of 
$\tilde{\boldsymbol{\omega}}$, and that $z^{-(N+1)}(\tilde{R}(z))^{1/2} \! 
\sim_{\tilde{\mathcal{Y}}^{\pm} \ni z \to \alpha_{p_{\mathfrak{s}-1}}} \! 
\pm 1$, using the fact that $\tilde{\boldsymbol{K}}$ is a point of order 
$2$, one arrives at
\begin{align*}
\tilde{\boldsymbol{d}} =& \, -\tilde{\boldsymbol{K}} \! - \! \sum_{j=1}^{N} 
\int_{\tilde{a}_{N+1}}^{\tilde{z}_{j}^{-}} \tilde{\boldsymbol{\omega}} \! = \! 
\tilde{\boldsymbol{K}} \! -\! \sum_{j=1}^{N} \int_{\tilde{a}_{N+1}}^{\tilde{z}_{j}^{-}} 
\tilde{\boldsymbol{\omega}} \! = \! -\tilde{\boldsymbol{K}} \! + \! \sum_{j=1}^{N} 
\int_{\tilde{a}_{N+1}}^{\tilde{z}_{j}^{+}} \tilde{\boldsymbol{\omega}} \! = \! 
\tilde{\boldsymbol{K}} \! + \! \sum_{j=1}^{N} \int_{\tilde{a}_{N+1}}^{\tilde{z}_{j}^{+}} 
\tilde{\boldsymbol{\omega}} \\
=& \, -\sum_{j=1}^{N} \int_{\tilde{a}_{j}}^{\tilde{z}_{j}^{-}} \tilde{\boldsymbol{
\omega}} \! = \! \sum_{j=1}^{N} \int_{\tilde{a}_{j}}^{\tilde{z}_{j}^{+}} 
\tilde{\boldsymbol{\omega}}.
\end{align*}

\hspace*{0.50cm}
Associated with the Riemann matrix of $\tilde{\boldsymbol{\beta}}$-periods, 
$\tilde{\tau}$, is the Riemann theta function
\begin{equation} \label{eqrmthetafin} 
\tilde{\boldsymbol{\theta}}(z;\tilde{\tau}) \! =: \! \tilde{\boldsymbol{\theta}}
(z) \! = \! \sum_{m \in \mathbb{Z}^{N}} \me^{2 \pi \mi (m,z)+\mi \pi 
(m,\tilde{\tau}m)}, \quad z \! \in \! \mathbb{C}^{N};
\end{equation}
$\tilde{\boldsymbol{\theta}}(z)$ has the following evenness and (quasi-) 
periodicity properties,
\begin{equation*}
\tilde{\boldsymbol{\theta}}(-z) \! = \! \tilde{\boldsymbol{\theta}}(z), 
\qquad \tilde{\boldsymbol{\theta}}(z \! + \! e_{j}) \! = \! 
\tilde{\boldsymbol{\theta}}(z), \qquad \text{and} \qquad 
\tilde{\boldsymbol{\theta}}(z \! \pm \! \tilde{\tau}_{j}) \! = \! 
\me^{\mp 2 \pi \mi z_{j}-\mi \pi \tilde{\tau}_{jj}} 
\tilde{\boldsymbol{\theta}}(z),
\end{equation*}
where $\tilde{\tau}_{j} \! := \! \tilde{\tau}e_{j}$ $(\in \! \mathbb{C}^{N})$, 
$j \! = \! 1,2,\dotsc,N$.
\end{enumerate}
\subsubsection{MPC ORF and MPA Asymptotics} \label{subsub2} 
The genesis of the MPC ORF problem consists in reformulating it, in the 
spirit of Fokas-Its-Kitaev \cite{a49,a50}, as $K$ families of matrix RHPs 
on $\overline{\mathbb{R}}$, where one subfamily consists of $\hat{K} \! 
:= \! \# \lbrace \mathstrut k \! \in \! \lbrace 1,2,\dotsc,K \rbrace; \, 
\alpha_{k} \! = \! \infty \rbrace$ matrix RHPs on $\overline{\mathbb{R}}$, 
and another subfamily consists of $\tilde{K} \! := \! \# \lbrace \mathstrut 
k \! \in \! \lbrace 1,2,\dotsc,K \rbrace; \, \alpha_{k} \! \neq \! \infty \rbrace$ 
matrix RHPs on $\overline{\mathbb{R}}$, with $K \! = \! \hat{K} \! + \! 
\tilde{K}$, and then to study both the finite-$n$ and large-$n$ asymptotic 
(in the double-scaling limit $\mathscr{N},n \! \to \! \infty$ such that 
$\mathscr{N}/n \! = \! 1 \! + \! o(1))$ behaviours of the corresponding 
$K$ solution families, wherein the latter family of $K$ asymptotic analyses 
consists of an amalgamation of the Deift-Zhou non-linear steepest-descent 
method for undulatory matrix RHPs \cite{a52} (see, also, \cite{pz1,pz2,pz3,
qydo,ydptgsman}) and the extension of Deift-Venakides-Zhou \cite{a54} 
(see, also, \cite{pz19,pz19a}).

For $n \! \in \! \mathbb{N}$ and $k \! \in \! \lbrace 1,2,\dotsc,K \rbrace$, 
given the $K$ neither necessarily distinct nor bounded, but otherwise 
arbitrary, poles $\alpha_{1},\alpha_{2},\dotsc,\alpha_{K}$ lying on the 
support of the orthogonality measure (cf. Equations~\eqref{eq1} 
and~\eqref{eq2}) $\md \mu (z) \! = \! \exp (-\mathscr{N}V(z)) \, \md z$, 
$\mathscr{N} \! \in \! \mathbb{N}$, where the external field $V \colon 
\overline{\mathbb{R}} \setminus \lbrace \alpha_{1},\alpha_{2},
\dotsc,\alpha_{K} \rbrace \! \to \! \mathbb{R}$ satisfies 
conditions~\eqref{eq3}--\eqref{eq5}, the monic MPC ORF problem is formulated 
(see Section~\ref{sec2}, Lemma~$\bm{\mathrm{RHP}_{\mathrm{MPC}}})$ as a 
family of $K$ matrix RHPs on $\overline{\mathbb{R}}$, where one subfamily, 
corresponding to those $k \! \in \! \lbrace 1,2,\dotsc,K \rbrace$ for which 
$\alpha_{k} \! = \! \infty$, consists of $\hat{K}$ matrix RHPs on 
$\overline{\mathbb{R}}$, and another subfamily, corresponding to those 
$k \! \in \! \lbrace 1,2,\dotsc,K \rbrace$ for which $\alpha_{k} \! \neq \! \infty$, 
consists of $\tilde{K}$ matrix RHPs on $\overline{\mathbb{R}}$, with $K \! = \! 
\hat{K} \! + \! \tilde{K}$.

For $n \! \in \! \mathbb{N}$ and $k \! \in \! \lbrace 1,2,\dotsc,K \rbrace$ 
such that $\alpha_{p_{\mathfrak{s}}} \! := \! \alpha_{k} \! = \! \infty$ or 
$\alpha_{p_{\mathfrak{s}}} \! := \! \alpha_{k} \! \neq \! \infty$, the finite-$n$ 
analysis for the monic MPC ORFs, $\pmb{\pi}^{n}_{k}(z)$, $z \! \in \! \mathbb{C}$, 
the norming constants, $\mu^{r}_{n,\varkappa_{nk}}(n,k)$, $r \! \in \! \lbrace 
\infty,f \rbrace$, and the MPC ORFs, $\phi^{n}_{k}(z) \! := \! \mu^{r}_{n,
\varkappa_{nk}}(n,k) \pmb{\pi}^{n}_{k}(z)$, $z \! \in \! \mathbb{C}$, is 
presented in Section~\ref{sec2}.
\begin{eeee} \label{rem1.3.1} 
\textsl{It must be stressed that the discussion thus far, as well as that 
which follows, while valid in its own right for finite $n$ $(\in \! 
\mathbb{N})$, is germane, principally, to transforming the family 
of $K$ matrix {\rm RHPs} on $\overline{\mathbb{R}}$ characterising 
the monic {\rm MPC ORFs} into a family of $K$ equivalent, or `model', 
matrix {\rm RHPs} on $\overline{\mathbb{R}}$ suitable for asymptotic 
analysis (in the double-scaling limit $\mathscr{N},n \! \to \! \infty$ 
such that $\mathscr{N}/n \! = \! 1 \! + \! o(1))$.}
\end{eeee}
For $n \! \in \! \mathbb{N}$ and $k \! \in \! \lbrace 1,2,\dotsc,K \rbrace$ 
such that $\alpha_{p_{\mathfrak{s}}} \! := \! \alpha_{k} \! = \! \infty$ or 
$\alpha_{p_{\mathfrak{s}}} \! := \! \alpha_{k} \! \neq \! \infty$, the 
large-$n$ asymptotic (in the double-scaling limit $\mathscr{N},n \! \to 
\! \infty$ such that $\mathscr{N}/n \! = \! 1 \! + \! o(1))$ analysis for 
the monic MPC ORFs, $\pmb{\pi}^{n}_{k}(z)$, $z \! \in \! \mathbb{C}$, 
the norming constants, $\mu^{r}_{n,\varkappa_{nk}}(n,k)$, $r \! \in \! 
\lbrace \infty,f \rbrace$, the MPC ORFs, $\phi^{n}_{k}(z) \! := \! \mu^{r}_{n,
\varkappa_{nk}}(n,k) \pmb{\pi}^{n}_{k}(z)$, $z \! \in \! \mathbb{C}$, and the 
corresponding MPA error terms, $\widehat{\pmb{\mathrm{E}}}_{\mu}(z)$ 
(for $\alpha_{p_{\mathfrak{s}}} \! := \! \alpha_{k} \! = \! \infty)$ and 
$\widetilde{\pmb{\mathrm{E}}}_{\mu}(z)$ (for $\alpha_{p_{\mathfrak{s}}} \! 
:= \! \alpha_{k} \! \neq \! \infty)$, $z \! \in \! \mathbb{C}$, requires the 
consideration of a family of $K$ variational conditions for $K$ suitably posed 
energy minimisation problems, and the establishment of the existence, 
uniqueness, and regularity properties of the corresponding family of 
\emph{equilibrium measures} which solve the family of energy minimisation 
problems. It is a well-known mathematical fact that variational conditions for 
energy minimisation problems in logarithmic potential theory, via the equilibrium 
measure (see, for example, \cite{a51,a55,a56,a57,a58,ear,aagear,mppta43,
glgyhasz}), play an absolutely crucial r\^{o}le in asymptotic analyses of RHPs associated 
with continuous and discrete OPs, their roots, and the corresponding recurrence 
relations (see, for example \cite{a61,pz22,toconwaypaft,a60,mvl6,dnhysakneln,
mvl7,mvl4,apkh,pelerkhy1,pelerkhy2,pbkl,avkomlovuein,hambuckty,kal1,kltydwang1,
lywi,ghoargtizremd,coramsibvot,fbmbsylkm,amfs,afmfvls,nikishphr,airevgllamf,
aptekarevlagomasfinkels,mondth3,alfodeno,afdoakaasporn,barhou,seungmenglee,
xbowylsxxyqzo}), the Painlev\'{e} Transcendents (see, for example, \cite{asfikyn,a53,
arsolvyakhv,astiaorov,aiakjo2,aiakjo1,sachobotsingp2,pz21,pz8,pz5,mvl3,tcaiik,tmcl2,tltg4,
tmcl3,mdabjk1,daik,peetsd11,mbat1,mbat2,marcobot,rentll,trredghsy,mxrtin,robbybuck,
rbukpdmllp1,rbukpdmllpp2,sxxuoahzqyy}), Random Matrix Theory (see, for example, 
\cite{a67,a59,pgkv,pddg2,pddg3,dikz1,dikz3,pddg1,dikz2,pdikv,ditskz,pdariik,
pdaiik,pftgmntson,briou,abtkytaets,itik1,mrots1,peetsd7,stiarkkkoz,pz26,pbai,prabj1,
prabj3,pvmbl,bdku,pbku,thbotas3,pvlado,mlevaflgu,rehlbeadvel,peetsd8,tmsyclfahsbi,
petadzhaf,bfhsyvso,tcysisky,prabj2,mppta4,ivkras3,ivkras1,krischschven,mvl9,mvl8,
kmfw,kuijsabeylz,arktac1,abjkuijlaarsatovbis,galomdl,amfroear,riejszara,ieaezgmlavl,
peetsd12,peetsd9,guihlermo,claeysgravakdtrmcl,nekeshystiny,christopheiertmeys,
nikreilhoz,maxratcesfri,tcabjk1,mvl2,mvl1,tmcl1,pz6,dk1,dgk1,dkmo1,mddg,baiktbothren,
jbdw1,jbdw2,nmekmcl2,nmekmcl1,sdkz,stvdldgdlzg,stvdl1,dkrlz,kdak2,peetsd2,
mbertolamcaf,mbeoamcafo,marcobotts,bblp,mbat3,mymo1,mymo4,mymo,mymo2,
mymo3,terstaivoroh,rhotagik,renott,tuberct,tascyterirky,thbotas,thbotas2,pgavkoolyy,
mvl5,mgirotti,chtopheomys,peetsd6,mondth2,peetsd5,peetsd4,peetsd3,rlxyppee1,
rlxypqee2,ccharr,alfie,ccharlieradeano,ccharlieradoeraene,nikulasaksbb,cwemdg,
nberecwebmdg,yzlhodi,sxxudiyqzo,sxxuyqzhao,sxxudandaiyzao,yzzhengxuyzhao1,
yzzhengxuyzhao,jianxueuifanyncn,wuxuao,sxxudeedai,ddaisxxulzhang,iaduxg,peetsd10,
mchenychenegfhan,buckyiech,yuxlapr,pz10}), Combinatorics (see, for example, 
\cite{jbpedets,jbpdkj,pz23,pz25,pz24,jhbrj,jobkzl}), and Pad\'{e} Approximants (see, for 
example, \cite{mppta1,pz4,a62,mppta2,aiavibamfsps,mptpaabo6,mppta5,mpptaarkkai,
aiekvevfietin,pmpat6677,mppta44,mppta45,mppta47,lbmly1,mlyattslv,mlyattslv2,
mlyattslv232,simplyrakhman,aagearsps2,aagearsps1,amartelinnovet,mppta48,ksvaw,
mppta26,mppta27,mppta31,mppta54,mppta55,mppta56,mppta53,mppta52,etinps,
ssuetin1,spmpeay1,seyina,seyinb,peedst1,veubsps,mlovveltiemc,mondth1,nomovchevaue1,
nrmovrkevaspin1,nrmovrkevaspin,mppta3,lyuy,tcfwmpta,mppta50,mppta49}). The situation 
with respect to the large-$n$ asymptotic analysis (in the double-scaling limit $\mathscr{N},n 
\! \to \! \infty$ such that $\mathscr{N}/n \! = \! 1 \! + \! o(1))$ for the monic MPC ORFs is, 
philosophically, analogous; however, unlike the asymptotic analyses for the OPs case, the 
asymptotic analysis for the monic MPC ORFs requires the consideration of $K$ different 
families of matrix RHPs on $\overline{\mathbb{R}}$, which, at the technical level, 
complicates the associated asymptotic analysis.
\begin{eeee} \label{rem1.3.2} 
\textsl{Before proceeding, an important notational preamble is requisite. 
Re-write (cf. Equations~\eqref{eq1} and~\eqref{eq2}$)$ $\md \mu (z) \! 
= \! \exp (-\mathscr{N}V(z)) \, \md z \! = \! \exp (-n \widetilde{V}(z)) \, 
\md z \! =: \! \md \widetilde{\mu}(z)$, $n \! \in \! \mathbb{N}$, where 
$\widetilde{V}(z) \! = \! z_{o}V(z)$, with $z_{o} \colon \mathbb{N} \times 
\mathbb{N} \! \to \! \mathbb{R}_{+}, \, (\mathscr{N},n) \! \mapsto \! z_{o} 
\! := \! \mathscr{N}/n$,\footnote{$\mathbb{R}_{\pm} \! := \! \lbrace x \! 
\in \! \mathbb{R}; \, \pm x \! > \! 0 \rbrace$.} and where the `scaled' external 
field $\widetilde{V} \colon \overline{\mathbb{R}} \setminus \lbrace \alpha_{1},
\alpha_{2},\dotsc,\alpha_{K} \rbrace \! \to \! \mathbb{R}$ satisfies the following 
conditions:\footnote{Since the double-scaling limit of interest is $\mathscr{N},n 
\! \to \! \infty$ such that $z_{o} \! = \! 1 \! + \! o(1)$, the monic MPC ORFs, 
$\lbrace \pmb{\pi}^{n}_{k}(z) \rbrace_{\underset{k=1,2,\dotsc,K}{n \in \mathbb{N}}}$, 
$z \! \in \! \mathbb{C}$, are now orthogonal with respect to the varying exponential 
measure $\md \widetilde{\mu}(z) \! = \! \exp (-n \widetilde{V}(z)) \, \md z$, with 
$\widetilde{V}(z)$ satisfying conditions~\eqref{eq20}--\eqref{eq22}, where the 
`large parameter', $n$, enters simultaneously into the order, $(n \! - \! 1)K \! 
+ \! k$, of the monic MPC ORFs and into the varying exponential weight; thus, 
asymptotics of the monic MPC ORFs are studied along a `diagonal strip' of a 
doubly-indexed sequence.}
\begin{gather}
\text{$\widetilde{V}(z)$ is real analytic on $\overline{\mathbb{R}} 
\setminus \lbrace \alpha_{1},\alpha_{2},\dotsc,\alpha_{K} \rbrace$}, 
\label{eq20} \\
\lim_{x \to \alpha_{i}} \left(\dfrac{\widetilde{V}(x)}{\ln (x^{2} \! + \! 1)} 
\right) \! = \! +\infty, \quad i \! \in \! \lbrace \mathstrut k \! \in \! 
\lbrace 1,2,\dotsc,K \rbrace; \, \alpha_{k} \! = \! \infty \rbrace, 
\label{eq21} \\
\lim_{x \to \alpha_{j}} \left(\dfrac{\widetilde{V}(x)}{\ln ((x \! - \! 
\alpha_{j})^{-2} \! + \! 1)} \right) \! = \! +\infty, \quad j \! \in \! 
\lbrace \mathstrut k \! \in \! \lbrace 1,2,\dotsc,K \rbrace; \, 
\alpha_{k} \! \neq \! \infty \rbrace; \label{eq22}
\end{gather}
for example, a rational function of the form (counting multiplicities of 
poles) $\widetilde{V} \colon \overline{\mathbb{R}} \setminus \lbrace 
\alpha_{1},\alpha_{2},\dotsc,\alpha_{K} \rbrace \! \ni \! z \! \mapsto 
\! \sum_{i \in \lbrace \mathstrut k \in \lbrace 1,2,\dotsc,K \rbrace; \, 
\alpha_{k} \neq \infty \rbrace} \sum_{j=-2m_{i}}^{-1} \tilde{\varsigma}_{j,i}
(z \! - \! \alpha_{i})^{j} \! + \! \sum_{l=0}^{2m_{\infty}} \tilde{\varsigma}_{l,
\infty}z^{l}$, where $m_{i} \! \in \! \mathbb{N}$, $m_{\infty} \! \in 
\! \mathbb{N}$, $\tilde{\varsigma}_{-2m_{i},i} \! > \! 0$, and 
$\tilde{\varsigma}_{2m_{\infty},\infty} \! > \! 0$ would satisfy 
conditions~\eqref{eq20}--\eqref{eq22}. As a consequence of the 
above-given transformation $V \! \to \! \widetilde{V}$, one must 
make the change $\md \mu \! \to \! \md \widetilde{\mu}$ (cf. 
Subsections~\ref{subsubsec1.2.1} and~\ref{subsubsec1.2.2}$)$ in 
Equations~\eqref{mvssinf1}, \eqref{mvssinf4}, \eqref{mvssinf10}, 
\eqref{mvssfin3}, and~\eqref{mvssfin10} for the Markov-Stieltjes 
transform, and the corresponding associated $\mathrm{R}$-functions 
and {\rm MPA} error terms (see, in particular, Section~\ref{sek5}, 
Lemma~\ref{lemmpainffin}, Equations~\eqref{eqmvsstildemu}, 
\eqref{eqlemmvssinfmpa1}, \eqref{eqlemmvssfinmpa1}, 
\eqref{eqlemmvssinfmpa3}, 
and~\eqref{eqlemmvssfinmpa7}$)$.\footnote{Henceforth, all 
quantities also depend on the doubling-scaling parameter, $z_{o}$; 
but, for notational simplicity, unless where absolutely necessary, 
this additional dependence is suppressed.}}
\end{eeee}

Prior to stating, for $n \! \in \! \mathbb{N}$ and $k \! \in \! \lbrace 
1,2,\dotsc,K \rbrace$, the asymptotics, in the double-scaling limit 
$\mathscr{N},n \! \to \! \infty$ such that $z_{o} \! = \! 1 \! + \! o(1)$, 
for the monic MPC ORFs, $\pmb{\pi}^{n}_{k}(z)$, $z \! \in \! \mathbb{C}$, 
the norming constants, $\mu^{r}_{n,\varkappa_{nk}}(n,k)$, $r \! \in 
\! \lbrace \infty,f \rbrace$, the MPC ORFs, $\phi^{n}_{k}(z) \! := \! 
\mu^{r}_{n,\varkappa_{nk}}(n,k) \pmb{\pi}^{n}_{k}(z)$, $z \! \in \! 
\mathbb{C}$, and the corresponding MPA error terms (cf. 
Equations~\eqref{eqlemmvssinfmpa3} and~\eqref{eqlemmvssfinmpa7}, 
respectively), $\widehat{\pmb{\mathrm{E}}}_{\tilde{\mu}}(z)$ (for 
$\alpha_{p_{\mathfrak{s}}} \! := \! \alpha_{k} \! = \! \infty)$ and 
$\widetilde{\pmb{\mathrm{E}}}_{\tilde{\mu}}(z)$ (for $\alpha_{
p_{\mathfrak{s}}} \! := \! \alpha_{k} \! \neq \! \infty)$, $z \! \in \! 
\mathbb{C}$, the following discussion highlights, succinctly, several 
key results of this monograph, all of which are seminal ingredients 
subsumed in the asymptotic analysis (see Sections~\ref{sec3}--\ref{sek5} 
for complete details).
\begin{enumerate}
\item[\textbullet] For $n \! \in \! \mathbb{N}$ and $k \! \in \! \lbrace 1,2,
\dotsc,K \rbrace$ such that $\alpha_{p_{\mathfrak{s}}} \! := \! \alpha_{k} 
\! = \! \infty$ (resp., $\alpha_{p_{\mathfrak{s}}} \! := \! \alpha_{k} \! \neq 
\! \infty)$, let $\mathrm{I}^{\infty}_{\widetilde{V}}[\mu_{1}^{\text{\tiny EQ}}]$ 
(resp., $\mathrm{I}_{\widetilde{V}}^{f}[\mu_{2}^{\text{\tiny EQ}}])$ be 
the energy functional defined by Equation~\eqref{eqlm3.1a} (resp., 
\eqref{eqlm3.1b}), and consider the associated minimisation problem 
$E^{\infty}_{\widetilde{V}} \! := \! \inf \lbrace \mathstrut 
\mathrm{I}^{\infty}_{\widetilde{V}}[\mu_{1}^{\text{\tiny EQ}}]; \, 
\mu_{1}^{\text{\tiny EQ}} \! \in \! \mathscr{M}_{1}(\mathbb{R}) \rbrace$ 
(resp., $E^{f}_{\widetilde{V}} \! := \! \inf \lbrace \mathstrut \mathrm{I}^{f}_{
\widetilde{V}}[\mu_{2}^{\text{\tiny EQ}}]; \, \mu_{2}^{\text{\tiny EQ}} 
\! \in \! \mathscr{M}_{1}(\mathbb{R}) \rbrace)$. For $n \! \in \! \mathbb{N}$ 
and $k \! \in \! \lbrace 1,2,\dotsc,K \rbrace$ such that $\alpha_{p_{\mathfrak{s}}} 
\! := \! \alpha_{k} \! = \! \infty$ (resp., $\alpha_{p_{\mathfrak{s}}} \! := \! 
\alpha_{k} \! \neq \! \infty)$, the infimum is finite, and there exists an associated 
unique measure $\mu^{\infty}_{\widetilde{V}}$ (resp., $\mu^{f}_{\widetilde{V}})$, 
called the equilibrium measure, achieving this minimum, that is, $\mathscr{M}_{1}
(\mathbb{R}) \! \ni \! \mu^{\infty}_{\widetilde{V}}$ (resp., $\mathscr{M}_{1}
(\mathbb{R}) \! \ni \! \mu^{f}_{\widetilde{V}})$ and $\mathrm{I}^{\infty}_{
\widetilde{V}}[\mu^{\infty}_{\widetilde{V}}] \! = \! E^{\infty}_{\widetilde{V}}$ 
(resp., $\mathrm{I}^{f}_{\widetilde{V}}[\mu^{f}_{\widetilde{V}}] \! = \! 
E^{f}_{\widetilde{V}})$. (See Section~\ref{sec3}, the corresponding items of 
Lemmata~\ref{lem3.1}--\ref{lem3.3}.)
\end{enumerate}
\begin{enumerate}
\item[\textbullet] For $n \! \in \! \mathbb{N}$ and $k \! \in \! \lbrace 
1,2,\dotsc,K \rbrace$ such that $\alpha_{p_{\mathfrak{s}}} \! := \! 
\alpha_{k} \! = \! \infty$ (resp., $\alpha_{p_{\mathfrak{s}}} \! := \! 
\alpha_{k} \! \neq \! \infty)$, the associated equilibrium measure, 
$\mu_{\widetilde{V}}^{\infty}$ (resp., $\mu_{\widetilde{V}}^{f})$, 
is absolutely continuous with respect to Lebesgue measure. (See 
Section~\ref{sec3}, the corresponding items of Lemmata~\ref{lem3.5} 
and~\ref{lem3.6}.)
\item[\textbullet] For $n \! \in \! \mathbb{N}$ and $k \! \in \! \lbrace 
1,2,\dotsc,K \rbrace$ such that $\alpha_{p_{\mathfrak{s}}} \! := \! 
\alpha_{k} \! = \! \infty$ (resp., $\alpha_{p_{\mathfrak{s}}} \! := \! 
\alpha_{k} \! \neq \! \infty)$, the associated equilibrium measure, 
$\mu_{\widetilde{V}}^{\infty}$ (resp., $\mu_{\widetilde{V}}^{f})$, has 
support which consists of the disjoint union of a finite number, $N \! + \! 1$ 
$(\in \! \mathbb{N}$: see below), of bounded real (compact) intervals; 
in fact, $\supp (\mu^{\infty}_{\widetilde{V}}) \! =: \! J_{\infty} \! = \! 
\cup_{j=1}^{N+1}[\hat{b}_{j-1},\hat{a}_{j}]$ (resp., $\supp (\mu^{f}_{
\widetilde{V}}) \! =: \! J_{f} \! = \! \cup_{j=1}^{N+1}[\tilde{b}_{j-1},
\tilde{a}_{j}])$, where, with $N(n,k) \! := \! N \! \in \! \mathbb{N}_{0}$ 
and finite, and $\hat{b}_{j-1}(n,k) \! =: \! \hat{b}_{j-1}$ and $\hat{a}_{j}
(n,k) \! =: \! \hat{a}_{j}$ (resp., $\tilde{b}_{j-1}(n,k) \! =: \! \tilde{b}_{j-1}$ 
and $\tilde{a}_{j}(n,k) \! =: \! \tilde{a}_{j})$, $j \! = \! 1,2,\dotsc,N \! + \! 1$, 
the end-points of the intervals, $\lbrace \hat{b}_{j-1},\hat{a}_{j} \rbrace_{j=
1}^{N+1}$ (resp., $\lbrace \tilde{b}_{j-1},\tilde{a}_{j} \rbrace_{j=1}^{N+1})$, 
of the support, $J_{\infty}$ (resp., $J_{f})$, of $\mu^{\infty}_{\widetilde{V}}$ 
(resp., $\mu_{\widetilde{V}}^{f})$ have the property that $[\hat{b}_{j-1},
\hat{a}_{j}] \cap \lbrace \alpha_{1},\alpha_{2},\dotsc,\alpha_{K} \rbrace \! = 
\! \varnothing$ (resp., $[\tilde{b}_{j-1},\tilde{a}_{j}] \cap \lbrace \alpha_{1},
\alpha_{2},\dotsc,\alpha_{K} \rbrace \! = \! \varnothing)$, $j \! = \! 1,2,
\dotsc,N \! + \! 1$, and have been enumerated so that $-\infty \! < \! 
\hat{b}_{0} \! < \! \hat{a}_{1} \! < \! \hat{b}_{1} \! < \! \hat{a}_{2} \! < \! 
\dotsb \! < \! \hat{b}_{N} \! < \! \hat{a}_{N+1} \! < \! +\infty$ (resp., 
$-\infty \! < \! \tilde{b}_{0} \! < \! \tilde{a}_{1} \! < \! \tilde{b}_{1} \! < 
\! \tilde{a}_{2} \! < \! \dotsb \! < \! \tilde{b}_{N} \! < \! \tilde{a}_{N+1} 
\! < \! +\infty)$. (See Section~\ref{sec3}, the corresponding items of 
Lemma~\ref{lem3.7}.)
\item[\textbullet] For $n \! \in \! \mathbb{N}$ and $k \! \in \! \lbrace 
1,2,\dotsc,K \rbrace$ such that $\alpha_{p_{\mathfrak{s}}} \! := \! 
\alpha_{k} \! = \! \infty$ (resp., $\alpha_{p_{\mathfrak{s}}} \! := \! 
\alpha_{k} \! \neq \! \infty)$, the associated end-points, $\lbrace 
\hat{b}_{j-1},\hat{a}_{j} \rbrace_{j=1}^{N+1}$ (resp., $\lbrace 
\tilde{b}_{j-1},\tilde{a}_{j} \rbrace_{j=1}^{N+1})$, are not arbitrary; 
rather, they satisfy the---locally solvable---system of $2(N \! 
+ \! 1)$ \emph{moment equations} (transcendental equations) 
\eqref{eql3.7a}--\eqref{eql3.7c} (resp., 
\eqref{eql3.7g}--\eqref{eql3.7i}).\footnote{In the double-scaling 
limit $\mathscr{N},n \! \to \! \infty$ such that $z_{o} \! = \! 1 \! + 
\! o(1)$, the associated end-points, $\lbrace \hat{b}_{j-1},\hat{a}_{j} 
\rbrace_{j=1}^{N+1}$ (resp., $\lbrace \tilde{b}_{j-1},\tilde{a}_{j} 
\rbrace_{j=1}^{N+1})$, are real analytic functions of $z_{o}$ (except 
for a denumerable set of singularities, which are the \emph{critical 
values} of $z_{o}$; see, for example, \cite{a57}).} (See 
Section~\ref{sec3}, the corresponding items of Lemma~\ref{lem3.7}.)
\item[\textbullet] 
For $n \! \in \! \mathbb{N}$ and $k \! \in \! \lbrace 1,2,\dotsc,K \rbrace$ 
such that $\alpha_{p_{\mathfrak{s}}} \! := \! \alpha_{k} \! = \! \infty$, the 
\emph{density} of the associated equilibrium measure is given by
\begin{equation*}
\md \mu^{\infty}_{\widetilde{V}}(x) \! = \! \psi^{\infty}_{\widetilde{V}}
(x) \, \md x \! = \! (2 \pi \mi)^{-1}(\hat{R}(x))^{1/2}_{+} 
\hat{h}_{\widetilde{V}}(x) \chi_{J_{\infty}}(x) \, \md x,
\end{equation*}
where $(\hat{R}(z))^{1/2}$ is defined by Equation~\eqref{eql3.7d}, 
with $(\hat{R}(z))^{1/2}_{\pm} \! := \! \lim_{\varepsilon \downarrow 0}
(\hat{R}(z \! \pm \! \mi \varepsilon))^{1/2}$, $\hat{h}_{\widetilde{V}}(z)$ 
is defined by Equation~\eqref{eql3.7f}, and $\chi_{J_{\infty}}(x)$ is the 
characteristic function of the compact set $J_{\infty}$. For $n \! \in \! 
\mathbb{N}$ and $k \! \in \! \lbrace 1,2,\dotsc,K \rbrace$ such that 
$\alpha_{p_{\mathfrak{s}}} \! := \! \alpha_{k} \! \neq \! \infty$, the 
density of the associated equilibrium measure is given by
\begin{equation*}
\md \mu^{f}_{\widetilde{V}}(x) \! = \! \psi^{f}_{\widetilde{V}}(x) \, \md x 
\! = \! (2 \pi \mi)^{-1}(\tilde{R}(x))^{1/2}_{+} \tilde{h}_{\widetilde{V}}(x) 
\chi_{J_{f}}(x) \, \md x,
\end{equation*}
where $(\tilde{R}(z))^{1/2}$ is defined by Equation~\eqref{eql3.7j}, with 
$(\tilde{R}(z))^{1/2}_{\pm} \! := \! \lim_{\varepsilon \downarrow 0}
(\tilde{R}(z \! \pm \! \mi \varepsilon))^{1/2}$, $\tilde{h}_{\widetilde{V}}(z)$ is 
defined by Equation~\eqref{eql3.7l}, and $\chi_{J_{f}}(x)$ is the characteristic 
function of the compact set $J_{f}$.\footnote{For $n \! \in \! \mathbb{N}$ and 
$k \! \in \! \lbrace 1,2,\dotsc,K \rbrace$ such that $\alpha_{p_{\mathfrak{s}}} 
\! := \! \alpha_{k} \! = \! \infty$, $\psi^{\infty}_{\widetilde{V}}(x) \! 
\geqslant \! 0$ (resp., $\psi^{\infty}_{\widetilde{V}}(x) \! > \! 0)$ for 
$x \! \in \! J_{\infty}$ (resp., $x \! \in \! \operatorname{int}(J_{\infty}))$; 
in fact, $\psi^{\infty}_{\widetilde{V}}(x)$ behaves like a square root at 
the end-points of the intervals, $\lbrace \hat{b}_{j-1},\hat{a}_{j} 
\rbrace_{j=1}^{N+1}$, of the support, $J_{\infty}$, of the associated 
equilibrium measure, $\mu_{\widetilde{V}}^{\infty}$, that is, 
$\psi^{\infty}_{\widetilde{V}}(x) \! =_{x \downarrow \hat{b}_{j-1}} 
\! \mathcal{O}((x \! - \! \hat{b}_{j-1})^{1/2})$ and $\psi^{\infty}_{
\widetilde{V}}(x) \! =_{x \uparrow \hat{a}_{j}} \! \mathcal{O}((\hat{a}_{j} 
\! - \! x)^{1/2})$, $j \! = \! 1,2,\dotsc,N \! + \! 1$. For $n \! \in \! 
\mathbb{N}$ and $k \! \in \! \lbrace 1,2,\dotsc,K \rbrace$ such that 
$\alpha_{p_{\mathfrak{s}}} \! := \! \alpha_{k} \! \neq \! \infty$, 
$\psi^{f}_{\widetilde{V}}(x) \! \geqslant \! 0$ (resp., $\psi^{f}_{\widetilde{V}}
(x) \! > \! 0)$ for $x \! \in \! J_{f}$ (resp., $x \! \in \! \operatorname{int}(J_{f}))$; 
in fact, $\psi^{f}_{\widetilde{V}}(x)$ behaves like a square root at the end-points 
of the intervals, $\lbrace \tilde{b}_{j-1},\tilde{a}_{j} \rbrace_{j=1}^{N+1}$, 
of the support, $J_{f}$, of the associated equilibrium measure, 
$\mu_{\widetilde{V}}^{f}$, that is, $\psi^{f}_{\widetilde{V}}(x) \! =_{x 
\downarrow \tilde{b}_{j-1}} \! \mathcal{O}((x \! - \! \tilde{b}_{j-1})^{1/2})$ 
and $\psi^{f}_{\widetilde{V}}(x) \! =_{x \uparrow \tilde{a}_{j}} \! 
\mathcal{O}((\tilde{a}_{j} \! - \! x)^{1/2})$, $j \! = \! 1,2,\dotsc,N \! + \! 1$.} 
(See Section~\ref{sec3}, the corresponding items of Lemma~\ref{lem3.7}.)
\item[\textbullet] For $n \! \in \! \mathbb{N}$ and $k \! \in \! \lbrace 
1,2,\dotsc,K \rbrace$ such that $\alpha_{p_{\mathfrak{s}}} \! := \! 
\alpha_{k} \! = \! \infty$ (resp., $\alpha_{p_{\mathfrak{s}}} \! := \! 
\alpha_{k} \! \neq \! \infty)$, in the double-scaling limit $\mathscr{N},n 
\! \to \! \infty$ such that $z_{o} \! = \! 1 \! + \! o(1)$, $x \! \in \! J_{\infty} 
\Leftrightarrow \operatorname{dist}(x,\lbrace \hat{\mathfrak{z}}_{k}^{n}(j) 
\rbrace_{j=1}^{(n-1)K+k}) \! \to \! 0$ (resp., $x \! \in \! J_{f} \Leftrightarrow 
\operatorname{dist}(x,\lbrace \tilde{\mathfrak{z}}_{k}^{n}(j) 
\rbrace_{j=1}^{(n-1)K+k}) \! \to \! 0)$,\footnote{For $x \! \in \! \mathbb{R}$ 
and $\mathbb{X} \subset \mathbb{R}$, $\operatorname{dist}(x,\mathbb{X}) 
\! := \! \inf \lbrace \mathstrut \lvert x \! - \! y \rvert; \, y \! \in \! 
\mathbb{X} \rbrace$.} or, equivalently, $\lbrace \hat{\mathfrak{z}}_{k}^{n}
(j) \rbrace_{j=1}^{(n-1)K+k} \subset J_{\infty}$ (resp., $\lbrace 
\tilde{\mathfrak{z}}_{k}^{n}(j) \rbrace_{j=1}^{(n-1)K+k} \subset J_{f})$, that 
is, the zeros (counting multiplicities) $\hat{\mathfrak{z}}_{k}^{n}(j)$ (resp., 
$\tilde{\mathfrak{z}}_{k}^{n}(j))$, $j \! = \! 1,2,\dotsc,(n \! - \! 1)K \! + \! 
k$, accumulate on the real compact set $J_{\infty}$ (resp., $J_{f})$. (See 
Section~\ref{sec3}, the corresponding items of Lemmata~\ref{lemrootz} 
and~\ref{lemetatomu}.)
\item[\textbullet] For $n \! \in \! \mathbb{N}$ and $k \! \in \! \lbrace 
1,2,\dotsc,K \rbrace$ such that $\alpha_{p_{\mathfrak{s}}} \! := \! 
\alpha_{k} \! = \! \infty$ (resp., $\alpha_{p_{\mathfrak{s}}} \! := \! 
\alpha_{k} \! \neq \! \infty)$, define the associated normalised zero 
counting measure as follows: $\hat{\eta}_{\hat{\mathfrak{z}}}(x) \! 
:= \! ((n \! - \! 1)K \! + \! k)^{-1} \sum_{j=1}^{(n-1)K+k} \delta_{
\hat{\mathfrak{z}}^{n}_{k}(j)}(x)$ (resp., $\tilde{\eta}_{\tilde{\mathfrak{z}}}
(x) \! := \! ((n \! - \! 1)K \! + \! k)^{-1} \sum_{j=1}^{(n-1)K+k} \delta_{
\tilde{\mathfrak{z}}^{n}_{k}(j)}(x))$, where $\delta_{\hat{\mathfrak{z}}^{
n}_{k}(j)}(x)$ (resp., $\delta_{\tilde{\mathfrak{z}}^{n}_{k}(j)}(x))$ is the 
Dirac delta (atomic) mass concentrated at $\hat{\mathfrak{z}}^{n}_{k}
(j)$ (resp., $\tilde{\mathfrak{z}}^{n}_{k}(j))$, $j \! = \! 1,2,\dotsc,(n 
\! - \! 1)K \! + \! k$; then, in the double-scaling limit $\mathscr{N},n 
\! \to \! \infty$ such that $z_{o} \! = \! 1 \! + \! o(1)$, $\hat{\eta}_{
\hat{\mathfrak{z}}}$ (resp., $\tilde{\eta}_{\tilde{\mathfrak{z}}})$ converges 
weakly {}\footnote{A sequence of probability measures $\lbrace \mu_{m} 
\rbrace_{m \in \mathbb{N}}$ in $\mathscr{M}_{1}(\mathscr{D})$ is said 
to \emph{converge weakly} as $m \! \to \! \infty$ to $\mu \! \in \! 
\mathscr{M}_{1}(\mathscr{D})$, symbolically $\mu_{m} \! 
\overset{\ast}{\to} \! \mu$, if $\mu_{m}(f) \! := \! \int_{\mathscr{D}}
f(\xi) \, \md \mu_{m}(\xi) \! \to \! \int_{\mathscr{D}}f(\xi) \, \md 
\mu (\xi) \! =: \! \mu (f)$ as $m \! \to \! \infty$ for all $f \! \in 
\! \pmb{\operatorname{C}}_{\text{b}}^{0}(\mathscr{D})$, where 
$\pmb{\operatorname{C}}_{\text{b}}^{0}(\mathscr{D})$ denotes 
the set of all bounded continuous functions on $\mathscr{D}$ 
with compact support (see, for example, \cite{a55}).} to 
$\mu_{\widetilde{V}}^{\infty}$ (resp., $\mu_{\widetilde{V}}^{f})$, 
that is, $\hat{\eta}_{\hat{\mathfrak{z}}} \! \overset{\ast}{\to} \! 
\mu_{\widetilde{V}}^{\infty}$ (resp., $\tilde{\eta}_{\tilde{\mathfrak{z}}} 
\! \overset{\ast}{\to} \! \mu_{\widetilde{V}}^{f})$. (See Section~\ref{sec3}, 
the corresponding items of Lemma~\ref{lemetatomu}.)
\item[\textbullet] For $n \! \in \! \mathbb{N}$ and $k \! \in \! \lbrace 
1,2,\dotsc,K \rbrace$ such that $\alpha_{p_{\mathfrak{s}}} \! := \! 
\alpha_{k} \! = \! \infty$ (resp., $\alpha_{p_{\mathfrak{s}}} \! := \! 
\alpha_{k} \! \neq \! \infty)$, the associated equilibrium measure, 
$\mu^{\infty}_{\widetilde{V}}$ (resp., $\mu^{f}_{\widetilde{V}})$, 
and its compact support, $J_{\infty}$ (resp., $J_{f})$, are uniquely 
characterised by the following Euler-Lagrange variational equations: 
there exists $\hat{\ell} \! \in \! \mathbb{R}$ (resp., $\tilde{\ell} 
\! \in \! \mathbb{R})$, the associated Lagrange multiplier, such that 
the variational conditions~\eqref{eql3.8a} (resp., \eqref{eql3.8b}) 
are valid. (See Section~\ref{sec3}, the corresponding items of 
Lemma~\ref{lem3.8}.)
\item[\textbullet] For $n \! \in \! \mathbb{N}$ and $k \! \in \! \lbrace 
1,2,\dotsc,K \rbrace$ such that $\alpha_{p_{\mathfrak{s}}} \! := \! 
\alpha_{k} \! = \! \infty$, the associated Euler-Lagrange variational 
conditions can be conveniently recast in terms of the complex potential 
(the associated `$g$-function'), $g^{\infty}(z)$, of $\mu^{\infty}_{
\widetilde{V}}$, defined by Equation~\eqref{eql3.4gee1}. The associated 
`$g$-function' satisfies (see Section~\ref{sec3}, the corresponding item 
of Lemma~\ref{lem3.4}):
\begin{enumerate}
\item[$(G1)_{\infty}$] $g^{\infty}(z)$ is analytic for 
$z \! \in \! \mathbb{C} \setminus (-\infty,\max \lbrace 
\max_{q=1,2,\dotsc,\mathfrak{s}-1} \lbrace \alpha_{p_{q}} \rbrace,
\max \lbrace \supp (\mu^{\infty}_{\widetilde{V}}) \rbrace \rbrace)$;
\item[$(G2)_{\infty}$] $g^{\infty}(z) \! =_{\mathbb{C}_{\pm} 
\ni z \to \alpha_{k}} \! \tfrac{\varkappa_{nk}}{n} \ln z \! + 
\! \tilde{\mathscr{P}}_{0}(n,k) \! + \! \mathcal{O}(z^{-1})$, 
where $\tilde{\mathscr{P}}_{0}(n,k)$ is defined by 
Equation~\eqref{eql3.4gee2};
\item[$(G3)_{\infty}$] $g^{\infty}_{+}(z) \! + \! g^{\infty}_{-}(z) \! - 
\! 2 \tilde{\mathscr{P}}_{0} \! - \! \widetilde{V}(z) \! - \! \hat{\ell} 
\! = \! 0$, $z \! \in \! J_{\infty}$, where $g^{\infty}_{\pm}(z) \! 
:= \! \lim_{\varepsilon \downarrow 0}g^{\infty}(z \! \pm \! \mi 
\varepsilon)$;
\item[$(G4)_{\infty}$] $g^{\infty}_{+}(z) \! + \! g^{\infty}_{-}(z) \! 
- \! 2 \tilde{\mathscr{P}}_{0} \! - \! \widetilde{V}(z) \! - \! \hat{\ell} 
\! \leqslant \! 0$, $z \! \in \! \mathbb{R} \setminus J_{\infty}$, 
where equality holds for at most a finite number of points;
\item[$(G5)_{\infty}$] $g^{\infty}_{+}(z) \! - \! g^{\infty}_{-}(z) \! = 
\! \mi f_{g^{\infty}}^{\mathbb{R}}(z)$, $z \! \in \! \mathbb{R}$, where 
$f_{g^{\infty}}^{\mathbb{R}}$ is a piecewise-continuous, real-valued, 
bounded function;
\item[$(G6)_{\infty}$] $\mi (g^{\infty}_{+}(z) \! - \! g^{\infty}_{-}(z) 
\! + \! 2 \pi \mi \sum_{q \in \lbrace \mathstrut j \in \lbrace 1,2,
\dotsc,\mathfrak{s}-1 \rbrace; \, \alpha_{p_{j}} >z \rbrace} 
\tfrac{\varkappa_{nk \tilde{k}_{q}}}{n})^{\prime} \! = \! 
\digamma_{\infty}(z)$, where the prime denotes differentiation 
with respect to $z$,
\begin{equation*}
\digamma_{\infty}(z) \! = \! 
\begin{cases}
2 \pi \left(\tfrac{(n-1)K+k}{n} \right) \psi^{\infty}_{\widetilde{V}}
(z) \! \geqslant \! 0, &\text{$z \! \in \! J_{\infty}$,} \\
0, &\text{$z \in \mathbb{R} \setminus J_{\infty}$,}
\end{cases}
\end{equation*}
and equality holds for at most a finite number of points.
\end{enumerate}
For $n \! \in \! \mathbb{N}$ and $k \! \in \! \lbrace 1,2,\dotsc,K 
\rbrace$ such that $\alpha_{p_{\mathfrak{s}}} \! := \! \alpha_{k} 
\! \neq \! \infty$, the associated Euler-Lagrange variational 
conditions can be conveniently recast in terms of the complex 
potential (the associated `$g$-function'), $g^{f}(z)$, of $\mu^{f}_{
\widetilde{V}}$, defined by Equation~\eqref{eql3.4gee3}. The 
associated `$g$-function' satisfies (see Section~\ref{sec3}, the 
corresponding item of Lemma~\ref{lem3.4}):
\begin{enumerate}
\item[$(G1)_{f}$] $g^{f}(z)$ is analytic for $z \! \in \! 
\mathbb{C} \setminus (-\infty,\max \lbrace \max_{q=1,\dotsc,
\mathfrak{s}-2,\mathfrak{s}} \lbrace \alpha_{p_{q}} \rbrace,
\max \lbrace \supp (\mu^{f}_{\widetilde{V}}) \rbrace \rbrace)$;
\item[$(G2)_{f}$] $g^{f}(z) \! =_{\mathbb{C}_{\pm} \ni z \to 
\alpha_{k}} \! -\left(\tfrac{\varkappa_{nk}-1}{n} \right) \ln (z \! - \! 
\alpha_{k}) \! + \! \hat{\mathscr{P}}_{0}^{\pm}(n,k) \! + \! \mathcal{O}
(z \! - \! \alpha_{k})$, where $\hat{\mathscr{P}}_{0}^{\pm}(n,k)$ is 
defined by Equations~\eqref{eql3.4gee4} and~\eqref{eql3.4gee5};
\item[$(G3)_{f}$] $g^{f}_{+}(z) \! + \! g^{f}_{-}(z) \! - \! \hat{
\mathscr{P}}_{0}^{+} \! - \! \hat{\mathscr{P}}_{0}^{-} \! - \! 
\widetilde{V}(z) \! - \! \tilde{\ell} \! = \! 0$, $z \! \in \! J_{f}$, 
where $g^{f}_{\pm}(z) \! := \! \lim_{\varepsilon \downarrow 0}
g^{f}(z \! \pm \! \mi \varepsilon)$;
\item[$(G4)_{f}$] $g^{f}_{+}(z) \! + \! g^{f}_{-}(z) \! - \! \hat{
\mathscr{P}}_{0}^{+} \! - \! \hat{\mathscr{P}}_{0}^{-} \! - \! 
\widetilde{V}(z) \! - \! \tilde{\ell} \! \leqslant \! 0$, $z \! \in \! 
\mathbb{R} \setminus J_{f}$, where equality holds for at most a 
finite number of points;
\item[$(G5)_{f}$] $g^{f}_{+}(z) \! - \! g^{f}_{-}(z) \! + \! \hat{
\mathscr{P}}_{0}^{-} \! - \! \hat{\mathscr{P}}_{0}^{+} \! = \! 
\mi f_{g^{f}}^{\mathbb{R}}(z)$, $z \! \in \! \mathbb{R}$, 
where $f_{g^{f}}^{\mathbb{R}}$ is a piecewise-continuous, 
real-valued, bounded function;
\item[$(G6)_{f}$] $\mi (g^{f}_{+}(z) \! - \! g^{f}_{-}(z) \! + \! \hat{
\mathscr{P}}_{0}^{-} \! - \! \hat{\mathscr{P}}_{0}^{+} \! + \! 2 \pi 
\mi \left(\tfrac{\varkappa_{nk}-1}{n} \right) \chi_{\lbrace \mathstrut 
z \in \mathbb{R}; \, z< \alpha_{k} \rbrace}(z) \! + \! 2 \pi \mi 
\sum_{q \in \lbrace \mathstrut j \in \lbrace 1,2,\dotsc,\mathfrak{s}
-2 \rbrace; \, \alpha_{p_{j}} >z \rbrace} \tfrac{\varkappa_{nk 
\tilde{k}_{q}}}{n})^{\prime} \! = \! \digamma_{f}(z)$, where the prime 
denotes differentiation with respect to $z$,
\begin{equation*}
\digamma_{f}(z) \! = \! 
\begin{cases}
2 \pi \left(\tfrac{(n-1)K+k}{n} \right) \psi^{f}_{\widetilde{V}}
(z) \! \geqslant \! 0, &\text{$z \! \in \! J_{f}$,} \\
0, &\text{$z \in \mathbb{R} \setminus J_{f}$,}
\end{cases}
\end{equation*}
and equality holds for at most a finite number of points.
\end{enumerate}
\end{enumerate}

In this monograph on the characterisation and asymptotics, in the 
double-scaling limit $\mathscr{N},n \! \to \! \infty$ such that 
$z_{o} \! = \! 1 \! + \! o(1)$, of MPC ORFs, the so-called `regular 
case' (see Section~\ref{sec3}, \pmb{Remark to Lemma~\ref{lem3.1}}) 
is studied, namely: (i) for $n \! \in \! \mathbb{N}$ and $k \! \in \! 
\lbrace 1,2,\dotsc,K \rbrace$ such that $\alpha_{p_{\mathfrak{s}}} 
\! := \! \alpha_{k} \! = \! \infty$, $\md \mu^{\infty}_{\widetilde{V}}$ 
(or $\widetilde{V} \colon \overline{\mathbb{R}} \setminus \lbrace 
\alpha_{1},\alpha_{2},\dotsc,\alpha_{K} \rbrace \! \to \! \mathbb{R})$ 
is \emph{regular} if
\begin{enumerate}
\item[$(\mathrm{i})_{\infty}$] $\hat{h}_{\widetilde{V}}(x) \! \neq \! 0$, 
$x \! \in \! J_{\infty}$,
\item[$(\mathrm{ii})_{\infty}$] $2(\tfrac{(n-1)K+k}{n}) \int_{J_{\infty}} 
\ln (\lvert z \! - \! \xi \rvert) \, \md \mu^{\infty}_{\widetilde{V}}(\xi) 
\! - \! 2 \sum_{q=1}^{\mathfrak{s}-1} \tfrac{\varkappa_{nk \tilde{k}_{
q}}}{n} \ln \lvert z \! - \! \alpha_{p_{q}} \rvert \! - \! \widetilde{V}(z) \! 
- \! \hat{\ell} \! < \! 0$, $z \! \in \! \mathbb{R} \setminus J_{\infty}$,
\item[$(\mathrm{iii})_{\infty}$] the inequalities~$(G4)_{\infty}$ 
and~$(G6)_{\infty}$ are strict, that is, $\leqslant$ (in $(G4)_{\infty})$ 
and $\geqslant$ (in $(G6)_{\infty})$ are replaced by $<$ and $>$, 
respectively;
\end{enumerate}
and (ii) for $n \! \in \! \mathbb{N}$ and $k \! \in \! \lbrace 1,2,\dotsc,
K \rbrace$ such that $\alpha_{p_{\mathfrak{s}}} \! := \! \alpha_{k} \! 
\neq \! \infty$, $\md \mu^{f}_{\widetilde{V}}$ (or $\widetilde{V} \colon 
\overline{\mathbb{R}} \setminus \lbrace \alpha_{1},\alpha_{2},\dotsc,
\alpha_{K} \rbrace \! \to \! \mathbb{R})$ is regular if
\begin{enumerate}
\item[$(\mathrm{i})_{f}$] $\tilde{h}_{\widetilde{V}}(x) \! \neq \! 0$, 
$x \! \in \! J_{f}$,
\item[$(\mathrm{ii})_{f}$] $2(\tfrac{(n-1)K+k}{n}) \int_{J_{f}} \ln (\lvert 
(z \! - \! \xi)/(\xi \! - \! \alpha_{k}) \rvert) \, \md\mu^{f}_{\widetilde{V}}
(\xi) \! - \! 2 \sum_{q=1}^{\mathfrak{s}-2} \tfrac{\varkappa_{nk 
\tilde{k}_{q}}}{n} \ln \lvert (z \! - \! \alpha_{p_{q}})/(\alpha_{p_{q}} \! 
- \! \alpha_{k}) \rvert \! - \! 2 \left(\tfrac{\varkappa_{nk}-1}{n} \right) 
\ln \lvert z \! - \! \alpha_{k} \rvert \! - \! \widetilde{V}(z) \! - \! 
\tilde{\ell} \! < \! 0$, $z \! \in \! \mathbb{R} \setminus J_{f}$,
\item[$(\mathrm{iii})_{f}$] the inequalities~$(G4)_{f}$ and~$(G6)_{f}$ are 
strict, that is, $\leqslant$ (in $(G4)_{f})$ and $\geqslant$ (in $(G6)_{f})$ 
are replaced by $<$ and $>$, respectively.
\end{enumerate}
\begin{eeee} \label{possbreakdown} 
\textsl{For $n \! \in \! \mathbb{N}$ and $k \! \in \! \lbrace 1,2,\dotsc,
K \rbrace$ such that $\alpha_{p_{\mathfrak{s}}} \! := \! \alpha_{k} 
\! = \! \infty$ (resp., $\alpha_{p_{\mathfrak{s}}} \! := \! \alpha_{k} 
\! \neq \! \infty)$, there are three distinct situations in which 
conditions~$(\mathrm{i})_{\infty}$--$(\mathrm{iii})_{\infty}$ (resp., 
$(\mathrm{i})_{f}$--$(\mathrm{iii})_{f})$ may fail: {\rm (1)} for at least one 
$\tilde{z}_{1} \! \in \! \mathbb{R} \setminus J_{\infty}$ (resp., $\tilde{z}_{2} 
\! \in \! \mathbb{R} \setminus J_{f})$, $2(\tfrac{(n-1)K+k}{n}) \int_{J_{\infty}} 
\ln (\lvert \tilde{z}_{1} \! - \! \xi \rvert) \, \md \mu^{\infty}_{\widetilde{V}}
(\xi) \! - \! 2 \sum_{q=1}^{\mathfrak{s}-1} \tfrac{\varkappa_{nk \tilde{k}_{q}}}{n} 
\ln \lvert \tilde{z}_{1} \! - \! \alpha_{p_{q}} \rvert \! - \! \widetilde{V}(\tilde{z}_{1}) 
\! - \! \hat{\ell} \! = \! 0$ (resp., $2(\tfrac{(n-1)K+k}{n}) \int_{J_{f}} \ln (\lvert 
(\tilde{z}_{2} \! - \! \xi)/(\xi \! - \! \alpha_{k}) \rvert) \, \md \mu^{f}_{\widetilde{V}}
(\xi) \! - \! 2 \sum_{q=1}^{\mathfrak{s}-2} \tfrac{\varkappa_{nk \tilde{k}_{q}}}{n} 
\ln \lvert (\tilde{z}_{2} \! - \! \alpha_{p_{q}})/(\alpha_{p_{q}} \! - \! \alpha_{k}) 
\rvert \! - \! 2 \left(\tfrac{\varkappa_{nk}-1}{n} \right) \ln \lvert \tilde{z}_{2} 
\! - \! \alpha_{k} \rvert \! - \! \widetilde{V}(\tilde{z}_{2}) \! - \! \tilde{\ell} \! 
= \! 0)$, that is, equality is attained for at least one point $\tilde{z}_{1}$ 
(resp., $\tilde{z}_{2})$ in the complement of the support of the associated 
equilibrium measure, which corresponds to the situation in which a `band' has 
just closed, or is just about to open, about $\tilde{z}_{1}$ (resp., $\tilde{z}_{2})$, 
a \emph{singular point of type $\mathrm{I}$} {\rm \cite{a57}}$;$ {\rm (2)} for at 
least one $\hat{z}_{1} \! \in \! \operatorname{int}(J_{\infty})$ (resp., $\hat{z}_{2} \! 
\in \! \operatorname{int}(J_{f}))$, $\hat{h}_{\widetilde{V}}(\hat{z}_{1}) \! = \! 0$ 
(resp., $\tilde{h}_{\widetilde{V}}(\hat{z}_{2}) \! = \! 0)$, that is, the function 
$\hat{h}_{\widetilde{V}}(z)$ (resp., $\tilde{h}_{\widetilde{V}}(z))$ vanishes for 
at least one point $\hat{z}_{1}$ (resp., $\hat{z}_{2})$ within the support of the 
associated equilibrium measure, which corresponds to the situation in which a 
`gap' is about to open, or close, about $\hat{z}_{1}$ (resp., $\hat{z}_{2})$, a 
\emph{singular point of type $\mathrm{II}$} {\rm \cite{a57}}$;$ and {\rm (3)} there 
exists at least one $j \! \in \! \lbrace 1,2,\dotsc,N \! + \! 1 \rbrace$, denoted 
$j_{\ast}$ (resp., $j^{\ast})$, such that $\hat{h}_{\widetilde{V}}(\hat{b}_{j_{\ast}-1}) 
\! = \! 0$ or $\hat{h}_{\widetilde{V}}(\hat{a}_{j_{\ast}}) \! = \! 0$ (resp., $\tilde{h}_{
\widetilde{V}}(\tilde{b}_{j^{\ast}-1}) \! = \! 0$ or $\tilde{h}_{\widetilde{V}}
(\tilde{a}_{j^{\ast}}) \! = \! 0)$, a \emph{singular point of type $\mathrm{III}$} 
{\rm \cite{a57}}. Each of these three cases can only occur a finite number of 
times due to the fact that $\widetilde{V} \colon \overline{\mathbb{R}} \setminus 
\lbrace \alpha_{1},\alpha_{2},\dotsc,\alpha_{K} \rbrace \! \to \! \mathbb{R}$ 
satisfies conditions~\eqref{eq20}--\eqref{eq22} {\rm \cite{a57,a58,a59,
amfroear,riejszara}}.}
\end{eeee}
\begin{eeee} \label{rem1.3.3} 
\textsl{The following correspondences should be noted: {\rm (i)} 
for $n \! \in \! \mathbb{N}$ and $k \! \in \! \lbrace 1,2,\dotsc,K 
\rbrace$ such that $\alpha_{p_{\mathfrak{s}}} \! := \! \alpha_{k} 
\! = \! \infty$, $g^{\infty} \colon \mathbb{N} \times \lbrace 1,
2,\dotsc,K \rbrace \times \mathbb{C} \setminus (-\infty,\max 
\lbrace \max_{q=1,2,\dotsc,\mathfrak{s}-1} \lbrace \alpha_{p_{q}} 
\rbrace,\max \lbrace \supp (\mu^{\infty}_{\widetilde{V}}) 
\rbrace \rbrace) \! \to \! \mathbb{C}$ solves the phase 
conditions~$(G1)_{\infty}$--$(G6)_{\infty}$ $\Leftrightarrow$ 
$\mathscr{M}_{1}(\mathbb{R}) \! \ni \! \mu^{\infty}_{\widetilde{V}}$ 
solves the variational conditions~\eqref{eql3.8a}$;$ and {\rm (ii)} 
for $n \! \in \! \mathbb{N}$ and $k \! \in \! \lbrace 1,2,\dotsc,K 
\rbrace$ such that $\alpha_{p_{\mathfrak{s}}} \! := \! \alpha_{k} \! 
\neq \! \infty$, $g^{f} \colon \mathbb{N} \times \lbrace 1,2,\dotsc,
K \rbrace \times \mathbb{C} \setminus (-\infty,\max \lbrace 
\max_{q=1,\dotsc,\mathfrak{s}-2,\mathfrak{s}} \lbrace 
\alpha_{p_{q}} \rbrace,\max \lbrace \supp (\mu^{f}_{\widetilde{V}}) 
\rbrace \rbrace) \! \to \! \mathbb{C}$ solves the phase 
conditions~$(G1)_{f}$--$(G6)_{f}$ $\Leftrightarrow$ 
$\mathscr{M}_{1}(\mathbb{R}) \! \ni \! \mu^{f}_{\widetilde{V}}$ 
solves the variational conditions~\eqref{eql3.8b}.}
\end{eeee}

For $n \! \in \! \mathbb{N}$ and $k \! \in \! \lbrace 1,2,\dotsc,K 
\rbrace$ such that $\alpha_{p_{\mathfrak{s}}} \! := \! \alpha_{k} \! 
= \! \infty$ (resp., $\alpha_{p_{\mathfrak{s}}} \! := \! \alpha_{k} \! 
\neq \! \infty)$, the final results of this monograph, that is, 
asymptotics, in the double-scaling limit $\mathscr{N},n \! \to 
\! \infty$ such that $z_{o} \! = \! 1 \! + \! o(1)$, of the monic 
MPC ORFs, $\pmb{\pi}_{k}^{n}(z)$, $z \! \in \! \mathbb{C}$, the 
corresponding MPA error term, $\widehat{\pmb{\mathrm{E}}}_{
\tilde{\mu}}(z)$ (resp., $\widetilde{\pmb{\mathrm{E}}}_{\tilde{\mu}}
(z))$, $z \! \in \! \mathbb{C}$, the associated norming constants, 
$\mu_{n,\varkappa_{nk}}^{r}(n,k)$, $r \! \in \! \lbrace \infty,f 
\rbrace$, and the MPC ORFs, $\phi_{k}^{n}(z) \! = \! \mu_{n,
\varkappa_{nk}}^{r}(n,k) \pmb{\pi}_{k}^{n}(z)$, $z \! \in \! 
\mathbb{C}$, are presented in 
Theorems~\ref{maintheoforinf1}--\ref{maintheoforinf2} (resp., 
Theorems~\ref{maintheoforfin1}--\ref{maintheoforfin2}) below 
(the complete details of the above-mentioned asymptotic analysis 
are given in Sections~\ref{sec3}--\ref{sek5}); but before doing so, 
however, some more notational preamble is necessary.

For $n \! \in \! \mathbb{N}$ and $k \! \in \! \lbrace 1,2,\dotsc,K \rbrace$ 
such that $\alpha_{p_{\mathfrak{s}}} \! := \! \alpha_{k} \! = \! \infty$, and 
$j \! = \! 1,2,\dotsc,N \! + \! 1$, let
\begin{gather}
\hat{\Phi}_{\hat{b}_{j-1}}(z) \! := \! \left(-\dfrac{3}{4}((n \! - \! 1)K \! + \! k) 
\int_{z}^{\hat{b}_{j-1}}(\hat{R}(\xi))^{1/2} \hat{h}_{\widetilde{V}}(\xi) \, \md 
\xi \right)^{2/3}, \label{eqmaininf77} \\
\hat{\Phi}_{\hat{a}_{j}}(z) \! := \! \left(\dfrac{3}{4}((n \! - \! 1)K \! + \! k) 
\int_{\hat{a}_{j}}^{z}(\hat{R}(\xi))^{1/2} \hat{h}_{\widetilde{V}}(\xi) \, \md 
\xi \right)^{2/3}, \label{eqmaininf84}
\end{gather}
and define the mutually disjoint open discs about $\hat{b}_{j-1}$ and 
$\hat{a}_{j}$, respectively, as $\hat{\mathbb{U}}_{\hat{\delta}_{\hat{b}_{j-1}}} 
\! := \! \lbrace \mathstrut z \! \in \! \mathbb{C}; \, \lvert z \! - \! \hat{b}_{j
-1} \rvert \! < \! \hat{\delta}_{\hat{b}_{j-1}} \rbrace$ and $\hat{\mathbb{
U}}_{\hat{\delta}_{\hat{a}_{j}}} \! := \! \lbrace \mathstrut z \! \in \! 
\mathbb{C}; \, \lvert z \! - \! \hat{a}_{j} \rvert \! < \! \hat{\delta}_{\hat{a}_{j}} 
\rbrace$, where $\hat{\delta}_{\hat{b}_{j-1}},\hat{\delta}_{\hat{a}_{j}} \! \in 
\! (0,1)$ are sufficiently small, positive real numbers chosen so that $\hat{
\mathbb{U}}_{\hat{\delta}_{\hat{b}_{i-1}}} \cap \hat{\mathbb{U}}_{\hat{\delta}_{
\hat{a}_{j}}} \! = \! \varnothing$ $\forall$ $i,j \! \in \! \lbrace 1,2,\dotsc,N \! 
+ \! 1 \rbrace$, $\hat{\mathbb{U}}_{\hat{\delta}_{\hat{b}_{j-1}}} \cap \lbrace 
\alpha_{p_{1}},\alpha_{p_{2}},\dotsc,\alpha_{p_{\mathfrak{s}}} \rbrace \! = 
\! \varnothing \! = \! \hat{\mathbb{U}}_{\hat{\delta}_{\hat{a}_{j}}} \cap 
\lbrace \alpha_{p_{1}},\alpha_{p_{2}},\dotsc,\alpha_{p_{\mathfrak{s}}} 
\rbrace$, and $\hat{\Phi}_{\hat{b}_{j-1}}(z)$ (resp., $\hat{\Phi}_{\hat{a}_{j}}
(z))$, which are bi-holomorphic, conformal, and non-orientation preserving 
(resp., orientation preserving), map (see Figure~\ref{forbhat}) 
$\hat{\mathbb{U}}_{\hat{\delta}_{\hat{b}_{j-1}}}$ (resp., (see 
Figure~\ref{forahat}) $\hat{\mathbb{U}}_{\hat{\delta}_{\hat{a}_{j}}})$, and, 
thus, the oriented contours $\hat{\Sigma}_{\hat{b}_{j-1}} \! := \! \cup_{m=
1}^{4} \hat{\Sigma}_{\hat{b}_{j-1}}^{m}$ (resp., $\hat{\Sigma}_{\hat{a}_{j}} \! 
:= \! \cup_{m=1}^{4} \hat{\Sigma}_{\hat{a}_{j}}^{m})$, injectively onto open, 
$n$- and $k$-dependent, neighbourhoods $\hat{\mathbb{U}}_{\hat{
\delta}_{\hat{b}_{j-1}}}^{\ast}$ (resp., $\hat{\mathbb{U}}_{\hat{\delta}_{
\hat{a}_{j}}}^{\ast})$ of the origin, such that $\hat{\Phi}_{\hat{b}_{j-1}}
(\hat{b}_{j-1}) \! = \! 0$, $\hat{\Phi}_{\hat{b}_{j-1}} \colon \hat{\mathbb{
U}}_{\hat{\delta}_{\hat{b}_{j-1}}} \! \to \! \hat{\mathbb{U}}_{\hat{\delta}_{
\hat{b}_{j-1}}}^{\ast} \! := \! \hat{\Phi}_{\hat{b}_{j-1}}(\hat{\mathbb{U}}_{
\hat{\delta}_{\hat{b}_{j-1}}})$, $\hat{\Phi}_{\hat{b}_{j-1}}(\hat{\mathbb{U}}_{
\hat{\delta}_{\hat{b}_{j-1}}} \cap \hat{\Sigma}_{\hat{b}_{j-1}}^{m}) \! = \! 
\hat{\Phi}_{\hat{b}_{j-1}}(\hat{\mathbb{U}}_{\hat{\delta}_{\hat{b}_{j-1}}}) 
\cap \gamma_{\hat{b}_{j-1}}^{\ast,m}$, and $\hat{\Phi}_{\hat{b}_{j-1}}
(\hat{\mathbb{U}}_{\hat{\delta}_{\hat{b}_{j-1}}} \cap \hat{\Omega}_{
\hat{b}_{j-1}}^{m}) \! = \! \hat{\Phi}_{\hat{b}_{j-1}}(\hat{\mathbb{U}}_{
\hat{\delta}_{\hat{b}_{j-1}}}) \cap \hat{\Omega}_{\hat{b}_{j-1}}^{\ast,m}$, 
$m \! = \! 1,2,3,4$, with $\hat{\Omega}_{\hat{b}_{j-1}}^{\ast,1} \! = \! 
\lbrace \mathstrut \zeta \! \in \! \mathbb{C}; \, \arg (\zeta) \! \in \! (0,
2 \pi/3) \rbrace$, $\hat{\Omega}_{\hat{b}_{j-1}}^{\ast,2} \! = \! \lbrace 
\mathstrut \zeta \! \in \! \mathbb{C}; \, \arg (\zeta) \! \in \! (2 \pi/3,\pi) 
\rbrace$, $\hat{\Omega}_{\hat{b}_{j-1}}^{\ast,3} \! = \! \lbrace \mathstrut 
\zeta \! \in \! \mathbb{C}; \, \arg (\zeta) \! \in \! (-\pi,-2 \pi/3) \rbrace$, 
and $\hat{\Omega}_{\hat{b}_{j-1}}^{\ast,4} \! = \! \lbrace \mathstrut 
\zeta \! \in \! \mathbb{C}; \, \arg (\zeta) \! \in \! (-2 \pi/3,0) \rbrace$ 
(resp., $\hat{\Phi}_{\hat{a}_{j}}(\hat{a}_{j}) \! = \! 0$, $\hat{\Phi}_{\hat{a}_{j}} 
\colon \hat{\mathbb{U}}_{\hat{\delta}_{\hat{a}_{j}}} \! \to \! \hat{\mathbb{
U}}_{\hat{\delta}_{\hat{a}_{j}}}^{\ast} \! := \! \hat{\Phi}_{\hat{a}_{j}}(\hat{
\mathbb{U}}_{\hat{\delta}_{\hat{a}_{j}}})$, $\hat{\Phi}_{\hat{a}_{j}}(\hat{
\mathbb{U}}_{\hat{\delta}_{\hat{a}_{j}}} \cap \hat{\Sigma}_{\hat{a}_{j}}^{m}) 
\! = \! \hat{\Phi}_{\hat{a}_{j}}(\hat{\mathbb{U}}_{\hat{\delta}_{\hat{a}_{j}}}) 
\cap \gamma_{\hat{a}_{j}}^{\ast,m}$, and $\hat{\Phi}_{\hat{a}_{j}}(\hat{
\mathbb{U}}_{\hat{\delta}_{\hat{a}_{j}}} \cap \hat{\Omega}_{\hat{a}_{j}}^{m}) 
\! = \! \hat{\Phi}_{\hat{a}_{j}}(\hat{\mathbb{U}}_{\hat{\delta}_{\hat{a}_{j}}}) 
\cap \hat{\Omega}_{\hat{a}_{j}}^{\ast,m}$, $m \! = \! 1,2,3,4$, with $\hat{
\Omega}_{\hat{a}_{j}}^{\ast,1} \! = \! \lbrace \mathstrut \zeta \! \in \! 
\mathbb{C}; \, \arg (\zeta) \! \in \! (0,2 \pi/3) \rbrace$, $\hat{\Omega}_{
\hat{a}_{j}}^{\ast,2} \! = \! \lbrace \mathstrut \zeta \! \in \! \mathbb{C}; \, 
\arg (\zeta) \! \in \! (2 \pi/3,\pi) \rbrace$, $\hat{\Omega}_{\hat{a}_{j}}^{
\ast,3} \! = \! \lbrace \mathstrut \zeta \! \in \! \mathbb{C}; \, \arg (\zeta) 
\! \in \! (-\pi,-2 \pi/3) \rbrace$, and $\hat{\Omega}_{\hat{a}_{j}}^{\ast,4} 
\! = \! \lbrace \mathstrut \zeta \! \in \! \mathbb{C}; \, \arg (\zeta) \! \in \! 
(-2 \pi/3,0) \rbrace)$. Consider, also, the decomposition of $\mathbb{C}$ 
(see Figure~\ref{figsectorhat}) into bounded and unbounded regions, and 
open neighbourhoods surrounding the end-points of the intervals, $\lbrace 
\hat{b}_{j-1},\hat{a}_{j} \rbrace_{j=1}^{N+1}$, of the support, $J_{\infty}$, 
of the associated equilibrium measure, $\mu_{\widetilde{V}}^{\infty}$.
\begin{figure}[bht]
\begin{center}
\vspace{0.55cm}

\end{center}
\caption{For $n \! \in \! \mathbb{N}$ and $k \! \in \! \lbrace 1,2,\dotsc,
K \rbrace$ such that $\alpha_{p_{\mathfrak{s}}} \! := \! \alpha_{k} \! = \! 
\infty$, the conformal mapping $\zeta \! = \! \hat{\Phi}_{\hat{b}_{j-1}}(z)$, 
$j \! = \! 1,2,\dotsc,N \! + \! 1$, where $(\hat{\Phi}_{\hat{b}_{j-1}})^{-1}$ 
denotes the inverse mapping.}\label{forbhat}
\end{figure}
\begin{figure}[tbh]
\begin{center}
\vspace{0.55cm}
%
\end{center}
\caption{For $n \! \in \! \mathbb{N}$ and $k \! \in \! \lbrace 1,2,\dotsc,
K \rbrace$ such that $\alpha_{p_{\mathfrak{s}}} \! := \! \alpha_{k} \! = \! 
\infty$, the conformal mapping $\zeta \! = \! \hat{\Phi}_{\hat{a}_{j}}(z)$, 
$j \! = \! 1,2,\dotsc,N \! + \! 1$, where $(\hat{\Phi}_{\hat{a}_{j}})^{-1}$ 
denotes the inverse mapping.}\label{forahat}
\end{figure}
\begin{figure}[tbh]
\begin{center}
\vspace{0.55cm}
%
\end{center}
\caption{The open neighbourhoods surrounding the end-points of the 
intervals, $\lbrace \hat{b}_{j-1},\hat{a}_{j} \rbrace_{j=1}^{N+1}$, of 
the support, $J_{\infty}$, of the associated equilibrium measure, 
$\mu_{\widetilde{V}}^{\infty}$.}\label{figsectorhat}
\end{figure}

Introduce, now, the Airy function, $\operatorname{Ai}(z)$, which appears 
in the final results of this monograph: $\operatorname{Ai}(z)$ is determined 
(uniquely) as the solution of the second-order non-constant coefficient 
homogeneous ODE (see, for example, Chapter~10 of \cite{abramsteg})
\begin{equation*}
\operatorname{Ai}^{\prime \prime}(z) \! - \! z \operatorname{Ai}(z) \! = \! 0,
\end{equation*}
where the primes denote differentiation with respect to $z$, with asymptotics
\begin{equation} \label{eqairy} 
\begin{split}
\operatorname{Ai}(z) \underset{\underset{\vert \arg z \vert < \pi}{z \to 
\infty}}{\sim}& \, \dfrac{1}{2 \sqrt{\pi}}z^{-1/4} \, \me^{-\frac{2}{3}z^{3/2}} 
\sum_{m=0}^{\infty}(-1)^{m}s_{m} \left(\dfrac{2}{3} z^{3/2} \right)^{-m}, \\
\operatorname{Ai}^{\prime}(z) \underset{\underset{\vert \arg z \vert < 
\pi}{z \to \infty}}{\sim}& \, -\dfrac{1}{2 \sqrt{\pi}}z^{1/4} \, \me^{-
\frac{2}{3}z^{3/2}} \sum_{m=0}^{\infty}(-1)^{m}t_{m} \left(\dfrac{2}{3} 
z^{3/2} \right)^{-m},
\end{split}
\end{equation}
where $s_{0} \! = \! t_{0} \! = \! 1$,
\begin{equation*}
s_{m} \! = \! \dfrac{\Gamma (3m \! + \! \frac{1}{2})}{(54)^{m}m! \Gamma 
(m \! + \! \frac{1}{2})} \! = \! \dfrac{(2m \! + \! 1)(2m \! + \! 3) \cdots 
(6m \! - \! 1)}{(216)^{m}m!}, \quad \quad t_{m} \! = \! -\left(\dfrac{6m 
\! + \! 1}{6m \! -\! 1} \right) s_{m}, \quad m \! \in \! \mathbb{N},
\end{equation*}
and $\Gamma (\pmb{\cdot})$ is the (Euler) gamma function.
\begin{eeee} \label{errocons} 
\textsl{In order to eschew a flood of superfluous notation, the `simplified' 
and `generic' notation $\mathcal{O}(\mathfrak{c}(n,k,\linebreak[4] 
z_{o})((n \! - \! 1)K 
\! + \! k)^{-2})$, where $\mathfrak{c}(n,k,z_{o}) \! =_{\underset{z_{o}=1+
o(1)}{\mathscr{N},n \to \infty}} \! \mathcal{O}(1)$, is maintained throughout 
Theorems~\ref{maintheoforinf1} and~\ref{maintheoforfin1} (see below), 
and it is to be understood in the following, `normal' sense: for $n \! 
\in \! \mathbb{N}$ and $k \! \in \! \lbrace 1,2,\dotsc,K \rbrace$ such 
that $\alpha_{p_{\mathfrak{s}}} \! := \! \alpha_{k} \! = \! \infty$ or 
$\alpha_{p_{\mathfrak{s}}} \! := \! \alpha_{k} \! \neq \! \infty$, in the 
double-scaling limit $\mathscr{N},n \! \to \! \infty$ such that $z_{o} \! 
= \! 1 \! + \! o(1)$, for a compact subset $D$ of $\mathbb{C}$, and 
uniformly with respect to $z \! \in \! D^{\ast} \! := \! D \setminus 
\lbrace \alpha_{p_{1}},\alpha_{p_{2}},\dotsc,\alpha_{p_{\mathfrak{s}}} 
\rbrace$, $\mathcal{O}(\mathfrak{c}(n,k,z_{o})((n \! - \! 1)K \! + \! 
k)^{-2}) \! := \! \mathcal{O}(\mathfrak{c}^{\ast}(n,k,z_{o};z)((n \! - 
\! 1)K \! + \! k)^{-2})$, where $\norm{\mathfrak{c}^{\ast}(n,k,z_{o};
\pmb{\cdot})}_{\mathcal{L}^{r}(D^{\ast})} \! =_{\underset{z_{o}=1+
o(1)}{\mathscr{N},n \to \infty}} \! \mathcal{O}(1)$, $r \! \in \! \lbrace 
2,\infty \rbrace$.}
\end{eeee}
\begin{eeee} \label{remmatelem} 
\textsl{For a $2 \times 2$ matrix-valued function $\mathfrak{I}
(\pmb{\cdot})$, $(\mathfrak{I}(\pmb{\cdot}))_{ij}$ or $\mathfrak{I}_{ij}
(\pmb{\cdot})$, $i,j \!= \! 1,2$, denotes the $(i \, j)$-element of 
$\mathfrak{I}(\pmb{\cdot})$.}
\end{eeee}
\begin{dddd} \label{maintheoforinf1} 
Let the external field $\widetilde{V} \colon \overline{\mathbb{R}} 
\setminus \lbrace \alpha_{1},\alpha_{2},\dotsc,\alpha_{K} \rbrace 
\! \to \! \mathbb{R}$ satisfy conditions~\eqref{eq20}--\eqref{eq22}. 
For $n \! \in \! \mathbb{N}$ and $k \! \in \! \lbrace 1,2,\dotsc,K 
\rbrace$ such that $\alpha_{p_{\mathfrak{s}}} \! := \! \alpha_{k} 
\! = \! \infty$, let $\md \mu_{\widetilde{V}}^{\infty}(x) \! = \! 
\psi_{\widetilde{V}}^{\infty}(x) \, \md x \! = \! (2 \pi \mi)^{-1}
(\hat{R}(x))^{1/2}_{+} \hat{h}_{\widetilde{V}}(x) \chi_{J_{\infty}}
(x) \, \md x$, where 
$(\hat{R}(z))^{1/2}$ is defined by Equation~\eqref{eql3.7d}, with 
$(\hat{R}(x))^{1/2}_{\pm} \! := \! \lim_{\varepsilon \downarrow 0}(\hat{R}
(x \! \pm \! \mi \varepsilon))^{1/2}$,\footnote{The branch of the square 
root is chosen so that $z^{-(N+1)}(\hat{R}(z))^{1/2} \! \sim_{\mathbb{C}_{
\pm} \ni z \to \alpha_{k}} \! \pm 1$.} $J_{\infty} \! := \! \supp 
(\mu_{\widetilde{V}}^{\infty}) \! = \! \cup_{j=1}^{N+1}[\hat{b}_{j-1},
\hat{a}_{j}]$, where $N \! \in \! \mathbb{N}_{0}$ and is finite, $[\hat{b}_{i-1},
\hat{a}_{i}] \cap [\hat{b}_{j-1},\hat{a}_{j}] \! = \! \varnothing$ $\forall$ 
$i \! \neq \! j \! \in \! \lbrace 1,2,\dotsc,N \! + \! 1 \rbrace$, $[\hat{b}_{j-1},
\hat{a}_{j}] \cap \lbrace \alpha_{p_{1}},\alpha_{p_{2}},\dotsc,\alpha_{
p_{\mathfrak{s}}} \rbrace \! = \! \varnothing$, and $-\infty \! < \! 
\hat{b}_{0} \! < \! \hat{a}_{1} \! < \! \hat{b}_{1} \! < \! \hat{a}_{2} \! 
< \! \dotsb \! < \! \hat{b}_{N} \! < \! \hat{a}_{N+1} \! < \! +\infty$, 
with $\lbrace \hat{b}_{j-1},\hat{a}_{j} \rbrace_{j=1}^{N+1}$ satisfying 
the associated $n$- and $k$-dependent system of $2(N \! + \! 1)$ 
real moment Equations~\eqref{eql3.7a}--\eqref{eql3.7c}, $\hat{h}_{
\widetilde{V}}(z)$, which is real analytic for $z \! \in \! \overline{
\mathbb{R}} \setminus \lbrace \alpha_{p_{1}},\alpha_{p_{2}},
\dotsc,\alpha_{p_{\mathfrak{s}}} \rbrace$, is defined by 
Equation~\eqref{eql3.7f},\footnote{Note: in the definition of 
$\hat{h}_{\widetilde{V}}(z)$ given by Equation~\eqref{eql3.7f}, there 
appears the contour integral $\oint_{\hat{C}_{\widetilde{V}}}$, the 
detailed description of which is given in the corresponding item 
of Lemma~\ref{lem3.7}.} and $\chi_{J_{\infty}}(x)$ is the characteristic 
function of the compact set $J_{\infty}$. Furthermore, suppose that, 
for $n \! \in \! \mathbb{N}$ and $k \! \in \! \lbrace 1,2,\dotsc,K 
\rbrace$ such that $\alpha_{p_{\mathfrak{s}}} \! := \! \alpha_{k} 
\! = \! \infty$, the external field, $\widetilde{V}$, is regular:
\begin{description}
\item[{\rm (i)}] $\hat{h}_{\widetilde{V}}(x) \! \neq \! 0$, 
$x \! \in \! J_{\infty}$$;$
\item[{\rm (ii)}]
\begin{equation*}
2 \left(\dfrac{(n \! - \! 1)K \! + \! k}{n} \right) \int_{J_{\infty}} \ln (\lvert 
x \! - \! \xi \rvert) \, \md \mu^{\infty}_{\widetilde{V}}(\xi) \! - \! 2 
\sum_{q=1}^{\mathfrak{s}-1} \dfrac{\varkappa_{nk \tilde{k}_{q}}}{n} 
\ln \lvert x \! - \! \alpha_{p_{q}} \rvert \! - \! \widetilde{V}(x) \! - \! 
\hat{\ell} \! = \! 0, \quad x \! \in \! J_{\infty},
\end{equation*}
which defines the associated variational constant $\hat{\ell}$ $(\in \! 
\mathbb{R})$,\footnote{Note that $\hat{\ell}$ is the same on each 
compact real interval $[\hat{b}_{j-1},\hat{a}_{j}]$, $j \! = \! 1,2,\dotsc,
N \! + \! 1$.} and
\begin{equation*}
2 \left(\dfrac{(n \! - \! 1)K \! + \! k}{n} \right) \int_{J_{\infty}} \ln (\lvert 
x \! - \! \xi \rvert) \, \md \mu^{\infty}_{\widetilde{V}}(\xi) \! - \! 2 
\sum_{q=1}^{\mathfrak{s}-1} \dfrac{\varkappa_{nk \tilde{k}_{q}}}{n} 
\ln \lvert x \! - \! \alpha_{p_{q}} \rvert \! - \! \widetilde{V}(x) \! - \! 
\hat{\ell} \! < \! 0, \quad x \! \in \! \mathbb{R} \setminus J_{\infty};
\end{equation*}
\item[{\rm (iii)}]
\begin{equation*}
g^{\infty}_{+}(x) \! + \! g^{\infty}_{-}(x) \! - \! 2 \tilde{\mathscr{P}}_{0} 
\! - \! \widetilde{V}(x) \! - \! \hat{\ell} \! < \! 0, \quad x \! \in \! 
\mathbb{R} \setminus J_{\infty},
\end{equation*}
where, for $z \! \in \! \mathbb{C} \setminus (-\infty,\max \lbrace 
\max_{q=1,2,\dotsc,\mathfrak{s}-1} \lbrace \alpha_{p_{q}} \rbrace,
\max \lbrace J_{\infty} \rbrace \rbrace)$, $g^{\infty}(z)$ is defined 
by Equation~\eqref{eql3.4gee1}, with $g^{\infty}_{\pm}(x) \! := \! 
\lim_{\varepsilon \downarrow 0}g^{\infty}(x \! \pm \! \mi \varepsilon)$, 
and $\tilde{\mathscr{P}}_{0}$ is defined by Equation~\eqref{eql3.4gee2}$;$
\item[{\rm (iv)}]
\begin{equation*}
\mi \left(g^{\infty}_{+}(z) \! - \! g^{\infty}_{-}(z) \! + \! 2 \pi \mi 
\sum_{q \in \lbrace \mathstrut j \in \lbrace 1,2,\dotsc,\mathfrak{s}
-1 \rbrace; \, \alpha_{p_{j}} >z \rbrace} \dfrac{\varkappa_{nk 
\tilde{k}_{q}}}{n} \right)^{\prime} \! > \! 0, \quad z \! \in \! J_{\infty},
\end{equation*}
where the prime denotes differentiation with respect to $z$.
\end{description}

For $n \! \in \! \mathbb{N}$ and $k \! \in \! \lbrace 1,2,\dotsc,K 
\rbrace$ such that $\alpha_{p_{\mathfrak{s}}} \! := \! \alpha_{k} 
\! = \! \infty$, let
\begin{gather} 
\tilde{\mathfrak{m}}^{\raise-0.5ex\hbox{$\scriptstyle \infty$}} 
\! = \! \operatorname{diag}
(\tilde{\mathfrak{m}}^{\raise-0.5ex\hbox{$\scriptstyle \infty$}}_{11},
\tilde{\mathfrak{m}}^{\raise-0.5ex\hbox{$\scriptstyle \infty$}}_{22}) 
\! = \! 
\begin{pmatrix}
\frac{\hat{\boldsymbol{\theta}}(\hat{\boldsymbol{u}}_{+}(\infty)+
\hat{\boldsymbol{d}})}{\hat{\boldsymbol{\theta}}(\hat{\boldsymbol{u}}_{+}
(\infty)-\frac{1}{2 \pi}((n-1)K+k) \hat{\boldsymbol{\Omega}}+\hat{
\boldsymbol{d}})} & 0 \\
0 & \frac{\hat{\boldsymbol{\theta}}(\hat{\boldsymbol{u}}_{+}(\infty)+
\hat{\boldsymbol{d}})}{\hat{\boldsymbol{\theta}}(-\hat{\boldsymbol{u}}_{
+}(\infty)-\frac{1}{2 \pi}((n-1)K+k) \hat{\boldsymbol{\Omega}}-
\hat{\boldsymbol{d}})}
\end{pmatrix}, \label{eqmaininf8} \\
\hat{\mathbb{M}}(z) \! = \! 
\begin{pmatrix}
\frac{1}{2}(\hat{\gamma}(z) \! + \! (\hat{\gamma}(z))^{-1}) \frac{
\hat{\boldsymbol{\theta}}(\hat{\boldsymbol{u}}(z)-\frac{1}{2 \pi}
((n-1)K+k) \hat{\boldsymbol{\Omega}}+\hat{\boldsymbol{d}})}{
\hat{\boldsymbol{\theta}}(\hat{\boldsymbol{u}}(z)+\hat{
\boldsymbol{d}})} & -\frac{1}{2 \mi}(\hat{\gamma}(z) \! - \! 
(\hat{\gamma}(z))^{-1}) \frac{\hat{\boldsymbol{\theta}}(-\hat{
\boldsymbol{u}}(z)-\frac{1}{2 \pi}((n-1)K+k) \hat{\boldsymbol{\Omega}}
+\hat{\boldsymbol{d}})}{\hat{\boldsymbol{\theta}}(-\hat{\boldsymbol{u}}
(z)+\hat{\boldsymbol{d}})} \\
\frac{1}{2 \mi}(\hat{\gamma}(z) \! - \! (\hat{\gamma}(z))^{-1}) \frac{
\hat{\boldsymbol{\theta}}(\hat{\boldsymbol{u}}(z)-\frac{1}{2 \pi}((n-1)K
+k) \hat{\boldsymbol{\Omega}}-\hat{\boldsymbol{d}})}{\hat{\boldsymbol{
\theta}}(\hat{\boldsymbol{u}}(z)-\hat{\boldsymbol{d}})} & 
\frac{1}{2}(\hat{\gamma}(z) \! + \! (\hat{\gamma}(z))^{-1}) \frac{\hat{
\boldsymbol{\theta}}(-\hat{\boldsymbol{u}}(z)-\frac{1}{2 \pi}((n-1)K+k) 
\hat{\boldsymbol{\Omega}}-\hat{\boldsymbol{d}})}{\hat{\boldsymbol{
\theta}}(-\hat{\boldsymbol{u}}(z)-\hat{\boldsymbol{d}})}
\end{pmatrix}, \label{eqmaininf9}
\end{gather}
where {}\footnote{Note: 
$\det (\tilde{\mathfrak{m}}^{\raise-0.5ex\hbox{$\scriptstyle \infty$}} 
\hat{\mathbb{M}}(z)) \! = \! 1$.} $\hat{\gamma}(z)$ is defined by 
Equation~\eqref{eqmaininf10}, $\hat{\boldsymbol{u}}(z) \! 
:= \! \int_{\hat{a}_{N+1}}^{z} \hat{\boldsymbol{\omega}}$, 
$\hat{\boldsymbol{u}}_{+}(\infty) \! := \! \int_{\hat{a}_{N+1}}^{
\infty^{+}} \hat{\boldsymbol{\omega}}$, with $\hat{\boldsymbol{\omega}}$ 
the associated normalised basis of holomorphic one-forms on $\hat{
\mathcal{Y}}$, $\hat{\boldsymbol{\Omega}} \! := \! (\hat{\Omega}_{1},
\hat{\Omega}_{2},\dotsc,\hat{\Omega}_{N})^{\mathrm{T}}$ $(\in \! 
\mathbb{R}^{N})$,\footnote{Note: $\operatorname{T}$ denotes 
transposition.} where $\hat{\Omega}_{j} \! = \! 2 \pi \int_{
\hat{b}_{j}}^{\hat{a}_{N+1}} \psi_{\widetilde{V}}^{\infty}(\xi) \, 
\md \xi$, $j \! = \! 1,2,\dotsc,N$, and $\hat{\boldsymbol{d}} 
\! \equiv \! -\sum_{j=1}^{N} \int_{\hat{a}_{j}}^{\hat{z}_{j}^{-}} 
\hat{\boldsymbol{\omega}}$ $(= \! \sum_{j=1}^{N} \int_{\hat{a}_{j}}^{
\hat{z}_{j}^{+}} \hat{\boldsymbol{\omega}})$, where $\lbrace 
\hat{z}_{1}^{\pm},\hat{z}_{2}^{\pm},\dotsc,\hat{z}_{N}^{\pm} \rbrace \! := 
\! \lbrace \mathstrut z \! \in \! \mathbb{C}_{\pm}; \, \hat{\gamma}(z) \! 
\mp \! (\hat{\gamma}(z))^{-1} \! = \! 0 \rbrace$, with $\hat{z}_{j}^{\pm} 
\! \in \! (\hat{a}_{j},\hat{b}_{j})^{\pm}$ $(\subset \mathbb{C}_{\pm})$, 
$j \! = \! 1,2,\dotsc,N$.\footnote{As points on the plane, $\boldsymbol{
\mathrm{pr}}(\hat{z}_{j}^{+}) \! = \! \boldsymbol{\mathrm{pr}}
(\hat{z}_{j}^{-}) \! = \! \hat{z}_{j} \! \in \! (\hat{a}_{j},\hat{b}_{j})$, 
$j \! = \! 1,2,\dotsc,N$.}

For $n \! \in \! \mathbb{N}$ and $k \! \in \! \lbrace 1,2,\dotsc,K 
\rbrace$ such that $\alpha_{p_{\mathfrak{s}}} \! := \! \alpha_{k} \! 
= \! \infty$, let $\mathcal{X} \colon \overline{\mathbb{C}} \setminus 
\overline{\mathbb{R}} \! \to \! \mathrm{SL}_{2}(\mathbb{C})$ be the 
unique solution of the corresponding monic {\rm MPC ORF RHP} 
$(\mathcal{X}(z),\upsilon (z),\overline{\mathbb{R}})$ stated in 
Lemma~$\bm{\mathrm{RHP}_{\mathrm{MPC}}}$, with integral 
representation given by Equation~\eqref{intrepinf} (see 
Lemma~\ref{lem2.1}$)$$;$ in particular,
\begin{equation*}
(\mathcal{X}(z))_{11} \! = \! \pmb{\pi}_{k}^{n}(z) \quad \quad \, 
\text{and} \, \quad \quad (\mathcal{X}(z))_{12} \! = \! \int_{\mathbb{R}} 
\dfrac{\pmb{\pi}_{k}^{n}(\xi) \me^{-n \widetilde{V}(\xi)}}{\xi \! - \! z} 
\, \dfrac{\md \xi}{2 \pi \mi}.
\end{equation*}
Then, for $n \! \in \! \mathbb{N}$ and $k \! \in \! \lbrace 1,2,\dotsc,K 
\rbrace$ such that $\alpha_{p_{\mathfrak{s}}} \! := \! \alpha_{k} \! = \! 
\infty$$:$\\
{\rm \pmb{(1)}} for $z \! \in \! \hat{\Upsilon}_{1}$,
\begin{align}
\pmb{\pi}_{k}^{n}(z) \underset{\underset{z_{o}=1+o(1)}{
\mathscr{N},n \to \infty}}{=}& \, \left(
\tilde{\mathfrak{m}}^{\raise-0.5ex\hbox{$\scriptstyle \infty$}}_{11} 
\hat{\mathbb{M}}_{11}(z) \left(1 \! + \! \dfrac{1}{(n \! - \! 1)K \! + \! k} 
\hat{\mathcal{R}}_{11}^{\sharp}(z) \! + \! \mathcal{O} \left(\dfrac{
\mathfrak{c}(n,k,z_{o})}{((n \! - \! 1)K \! + \! k)^{2}} \right) \right) 
\right. \nonumber \\
+&\left. \, 
\tilde{\mathfrak{m}}^{\raise-0.5ex\hbox{$\scriptstyle \infty$}}_{22} \hat{
\mathbb{M}}_{21}(z) \left(\dfrac{1}{(n \! - \! 1)K \! + \! k} \hat{\mathcal{
R}}_{12}^{\sharp}(z) \! + \! \mathcal{O} \left(\dfrac{\mathfrak{c}(n,k,
z_{o})}{((n \! - \! 1)K \! + \! k)^{2}} \right) \right) \right) \nonumber \\
\times& \, \me^{n(g^{\infty}(z)-\tilde{\mathscr{P}}_{0})}, \label{eqmaininf11}
\end{align}
and
\begin{align}
\int_{\mathbb{R}} \dfrac{\pmb{\pi}_{k}^{n}(\xi) 
\me^{-n \widetilde{V}(\xi)}}{\xi \! - \! z} \, \dfrac{\md \xi}{2 \pi \mi} 
\underset{\underset{z_{o}=1+o(1)}{\mathscr{N},n \to \infty}}{=}& \,
\left(\tilde{\mathfrak{m}}^{\raise-0.5ex\hbox{$\scriptstyle \infty$}}_{11} 
\hat{\mathbb{M}}_{12}(z) \left(1 \! + \! \dfrac{1}{(n \! - \! 1)K \! + \! k} 
\hat{\mathcal{R}}_{11}^{\sharp}(z) \! + \! \mathcal{O} \left(\dfrac{
\mathfrak{c}(n,k,z_{o})}{((n \! - \! 1)K \! + \! k)^{2}} \right) \right) 
\right. \nonumber \\
+&\left. \, 
\tilde{\mathfrak{m}}^{\raise-0.5ex\hbox{$\scriptstyle \infty$}}_{22} \hat{
\mathbb{M}}_{22}(z) \left(\dfrac{1}{(n \! - \! 1)K \! + \! k} \hat{\mathcal{
R}}_{12}^{\sharp}(z) \! + \! \mathcal{O} \left(\dfrac{\mathfrak{c}(n,k,
z_{o})}{((n \! - \! 1)K \! + \! k)^{2}} \right) \right) \right) \nonumber \\
\times& \, \me^{n \hat{\ell}} \me^{-n(g^{\infty}(z)-\tilde{\mathscr{P}}_{0})}, 
\label{eqmaininf12}
\end{align}
with
\begin{align}
\hat{\mathcal{R}}^{\sharp}(z) =& \, \sum_{j=1}^{N+1} \left(\dfrac{
(\hat{\alpha}_{0}(\hat{b}_{j-1}))^{-1}}{(z \! - \! \hat{b}_{j-1})^{2}} 
\hat{\boldsymbol{\mathrm{A}}}(\hat{b}_{j-1}) \! + \! \dfrac{(\hat{
\alpha}_{0}(\hat{b}_{j-1}))^{-2}}{z \! - \! \hat{b}_{j-1}} \left(\hat{
\alpha}_{0}(\hat{b}_{j-1}) \hat{\boldsymbol{\mathrm{B}}}(\hat{b}_{j-1}) 
\! - \! \hat{\alpha}_{1}(\hat{b}_{j-1}) \hat{\boldsymbol{\mathrm{A}}}
(\hat{b}_{j-1}) \right) \right. \nonumber \\
+&\left. \, \dfrac{(\hat{\alpha}_{0}(\hat{a}_{j}))^{-1}}{(z \! - \! 
\hat{a}_{j})^{2}} \hat{\boldsymbol{\mathrm{A}}}(\hat{a}_{j}) \! + \! 
\dfrac{(\hat{\alpha}_{0}(\hat{a}_{j}))^{-2}}{z \! - \! \hat{a}_{j}} \left(
\hat{\alpha}_{0}(\hat{a}_{j}) \hat{\boldsymbol{\mathrm{B}}}(\hat{a}_{j}) 
\! - \! \hat{\alpha}_{1}(\hat{a}_{j}) \hat{\boldsymbol{\mathrm{A}}}
(\hat{a}_{j}) \right) \right), \label{eqmaininf13}
\end{align}
where, for $j \! = \! 1,2,\dotsc,N \! + \! 1$,
\begin{gather}
\hat{\boldsymbol{\mathrm{A}}}(\hat{b}_{j-1}) \! := \! 
 \me^{\mi ((n-1)K+k) \hat{\mho}_{j}}, \label{eqmaininf17}
\end{gather}
with
\begin{gather}
\hat{\mathbb{A}}_{11}(\hat{b}_{j-1}) \! = \! -\hat{\mathbb{A}}_{22}
(\hat{b}_{j-1}) \! = \! -s_{1} \hat{\kappa}_{1}(\hat{b}_{j-1}) \hat{\kappa}_{2}
(\hat{b}_{j-1})(\hat{\mathfrak{Q}}_{0}(\hat{b}_{j-1}))^{-1}, \label{eqmaininf18} \\
\hat{\mathbb{A}}_{12}(\hat{b}_{j-1}) \! = \! -\mi s_{1}(\hat{\kappa}_{1}
(\hat{b}_{j-1}))^{2}(\hat{\mathfrak{Q}}_{0}(\hat{b}_{j-1}))^{-1}, \qquad 
\hat{\mathbb{A}}_{21}(\hat{b}_{j-1}) \! = \! -\mi s_{1}(\hat{\kappa}_{2}
(\hat{b}_{j-1}))^{2}(\hat{\mathfrak{Q}}_{0}(\hat{b}_{j-1}))^{-1}, 
\label{eqmaininf19}
\end{gather}
\begin{align}
\hat{\mathbb{B}}_{11}(\hat{b}_{j-1}) \! = \! -\hat{\mathbb{B}}_{22}
(\hat{b}_{j-1}) =& \, \hat{\kappa}_{1}(\hat{b}_{j-1}) \hat{\kappa}_{2}
(\hat{b}_{j-1}) \left(-s_{1}(\hat{\mathfrak{Q}}_{0}(\hat{b}_{j-1}))^{-1} 
\left(\hat{\daleth}^{1}_{1}(\hat{b}_{j-1}) \! + \! \hat{\daleth}^{1}_{-1}
(\hat{b}_{j-1}) \! - \! \hat{\mathfrak{Q}}_{1}(\hat{b}_{j-1})
(\hat{\mathfrak{Q}}_{0}(\hat{b}_{j-1}))^{-1} \right) \right. \nonumber \\
-&\left. \, t_{1} \left(\hat{\mathfrak{Q}}_{0}(\hat{b}_{j-1}) \! + \! 
(\hat{\mathfrak{Q}}_{0}(\hat{b}_{j-1}))^{-1} \hat{\aleph}^{1}_{1}
(\hat{b}_{j-1}) \hat{\aleph}^{1}_{-1}(\hat{b}_{j-1}) \right) \! + \! \mi 
(s_{1} \! + \! t_{1}) \left(\hat{\aleph}^{1}_{-1}(\hat{b}_{j-1}) \! - \! 
\hat{\aleph}^{1}_{1}(\hat{b}_{j-1}) \right) \right), \label{eqmaininf20} \\
\hat{\mathbb{B}}_{12}(\hat{b}_{j-1}) =& \, (\hat{\kappa}_{1}
(\hat{b}_{j-1}))^{2} \left(-\mi s_{1}(\hat{\mathfrak{Q}}_{0}(\hat{b}_{j-1}))^{-1} 
\left(2 \hat{\daleth}^{1}_{1}(\hat{b}_{j-1}) \! - \! \hat{\mathfrak{Q}}_{1}
(\hat{b}_{j-1})(\hat{\mathfrak{Q}}_{0}(\hat{b}_{j-1}))^{-1} \right) \right. 
\nonumber \\
+&\left. \, \mi t_{1} \left(\hat{\mathfrak{Q}}_{0}(\hat{b}_{j-1}) \! - \! 
(\hat{\mathfrak{Q}}_{0}(\hat{b}_{j-1}))^{-1}(\hat{\aleph}^{1}_{1}(\hat{b}_{j
-1}))^{2} \right) \! + \! 2(s_{1} \! - \! t_{1}) \hat{\aleph}^{1}_{1}(\hat{b}_{j-1}) 
\right), \label{eqmaininf21} \\
\hat{\mathbb{B}}_{21}(\hat{b}_{j-1}) =& \, (\hat{\kappa}_{2}
(\hat{b}_{j-1}))^{2} \left(-\mi s_{1}(\hat{\mathfrak{Q}}_{0}(\hat{b}_{j-1}))^{-1} 
\left(2 \hat{\daleth}^{1}_{-1}(\hat{b}_{j-1}) \! - \! \hat{\mathfrak{Q}}_{1}
(\hat{b}_{j-1})(\hat{\mathfrak{Q}}_{0}(\hat{b}_{j-1}))^{-1} \right) \right. 
\nonumber \\
+&\left. \, \mi t_{1} \left(\hat{\mathfrak{Q}}_{0}(\hat{b}_{j-1}) \! - \! 
(\hat{\mathfrak{Q}}_{0}(\hat{b}_{j-1}))^{-1}(\hat{\aleph}^{1}_{-1}(\hat{b}_{j
-1}))^{2} \right) \! - \! 2(s_{1} \! - \! t_{1}) \hat{\aleph}^{1}_{-1}(\hat{b}_{j-1}) 
\right), \label{eqmaininf22}
\end{align}
\begin{gather}
\hat{\mathbb{A}}_{11}(\hat{a}_{j}) \! = \! -\hat{\mathbb{A}}_{22}(\hat{a}_{j}) 
\! = \! -s_{1} \hat{\kappa}_{1}(\hat{a}_{j}) \hat{\kappa}_{2}(\hat{a}_{j}) 
\hat{\mathfrak{Q}}_{0}(\hat{a}_{j}), \label{eqmaininf23} \\
\hat{\mathbb{A}}_{12}(\hat{a}_{j}) \! = \! \mi s_{1}(\hat{\kappa}_{1}
(\hat{a}_{j}))^{2} \hat{\mathfrak{Q}}_{0}(\hat{a}_{j}), \qquad 
\hat{\mathbb{A}}_{21}(\hat{a}_{j}) \! = \! \mi s_{1}(\hat{\kappa}_{2}
(\hat{a}_{j}))^{2} \hat{\mathfrak{Q}}_{0}(\hat{a}_{j}), \label{eqmaininf24}
\end{gather}
\begin{align}
\hat{\mathbb{B}}_{11}(\hat{a}_{j}) \! = \! -\hat{\mathbb{B}}_{22}
(\hat{a}_{j}) =& \, \hat{\kappa}_{1}(\hat{a}_{j}) \hat{\kappa}_{2}
(\hat{a}_{j}) \left(-s_{1} \left(\hat{\mathfrak{Q}}_{1}(\hat{a}_{j}) \! + \! 
\hat{\mathfrak{Q}}_{0}(\hat{a}_{j}) \left(\hat{\daleth}^{1}_{1}(\hat{a}_{j}) \! 
+ \! \hat{\daleth}^{1}_{-1}(\hat{a}_{j}) \right) \right) \right. \nonumber \\
-&\left. \, t_{1} \left((\hat{\mathfrak{Q}}_{0}(\hat{a}_{j}))^{-1} \! + \! 
\hat{\mathfrak{Q}}_{0}(\hat{a}_{j}) \hat{\aleph}^{1}_{1}(\hat{a}_{j}) 
\hat{\aleph}^{1}_{-1}(\hat{a}_{j}) \right) \! + \! \mi (s_{1} \! + \! t_{1}) 
\left(\hat{\aleph}^{1}_{-1}(\hat{a}_{j}) \! - \! \hat{\aleph}^{1}_{1}(\hat{a}_{j}) 
\right) \right), \label{eqmaininf25} \\
\hat{\mathbb{B}}_{12}(\hat{a}_{j}) =& \, (\hat{\kappa}_{1}(\hat{a}_{j}))^{2} 
\left(\mi s_{1} \left(\hat{\mathfrak{Q}}_{1}(\hat{a}_{j}) \! + \! 2 \hat{
\mathfrak{Q}}_{0}(\hat{a}_{j}) \hat{\daleth}^{1}_{1}(\hat{a}_{j}) \right) 
\right. \nonumber \\
+&\left. \, \mi t_{1} \left(\hat{\mathfrak{Q}}_{0}(\hat{a}_{j})(\hat{\aleph}^{
1}_{1}(\hat{a}_{j}))^{2} \! - \! (\hat{\mathfrak{Q}}_{0}(\hat{a}_{j}))^{-1} 
\right) \! - \! 2(s_{1} \! - \! t_{1}) \hat{\aleph}^{1}_{1}(\hat{a}_{j}) \right), 
\label{eqmaininf26} \\
\hat{\mathbb{B}}_{21}(\hat{a}_{j}) =& \, (\hat{\kappa}_{2}(\hat{a}_{j}))^{2} 
\left(\mi s_{1} \left(\hat{\mathfrak{Q}}_{1}(\hat{a}_{j}) \! + \! 2 \hat{
\mathfrak{Q}}_{0}(\hat{a}_{j}) \hat{\daleth}^{1}_{-1}(\hat{a}_{j}) \right) 
\right. \nonumber \\
+&\left. \, \mi t_{1} \left(\hat{\mathfrak{Q}}_{0}(\hat{a}_{j})(\hat{\aleph}^{
1}_{-1}(\hat{a}_{j}))^{2} \! - \! (\hat{\mathfrak{Q}}_{0}(\hat{a}_{j}))^{-1} 
\right) \! + \! 2(s_{1} \! - \! t_{1}) \hat{\aleph}^{1}_{-1}(\hat{a}_{j}) 
\right), \label{eqmaininf27}
\end{align}
where, for $\varepsilon_{1},\varepsilon_{2} \! = \! \pm 1$,\footnote{Note: 
$\boldsymbol{0} \! := \! (0,0,\dotsc,0)^{\operatorname{T}}$ $(\in \! 
\mathbb{R}^{N})$.}
\begin{equation} \label{eqmaininf28} 
s_{1} \! = \! 5/72, \qquad \qquad \quad t_{1} \! = \! -7/72, \qquad 
\qquad \quad \hat{\mho}_{i} \! = \!
\begin{cases}
\hat{\Omega}_{i}, &\text{$i \! = \! 1,2,\dotsc,N$,} \\
0, &\text{$i \! = \! 0,N \! + \! 1$,}
\end{cases}
\end{equation}
\begin{align}
\hat{\kappa}_{1}(\varsigma) =& \, \dfrac{\hat{\boldsymbol{\theta}}
(\hat{\boldsymbol{u}}_{+}(\varsigma) \! - \! \frac{1}{2 \pi}((n \! - \! 1)K \! 
+ \! k) \hat{\boldsymbol{\Omega}} \! + \! \hat{\boldsymbol{d}})}{\hat{
\boldsymbol{\theta}}(\hat{\boldsymbol{u}}_{+}(\varsigma) \! + \! \hat{
\boldsymbol{d}})}, \qquad \qquad \hat{\kappa}_{2}(\varsigma) = \dfrac{
\hat{\boldsymbol{\theta}}(\hat{\boldsymbol{u}}_{+}(\varsigma) \! - \! 
\frac{1}{2 \pi}((n \! - \! 1)K \! + \! k) \hat{\boldsymbol{\Omega}} \! - \! 
\hat{\boldsymbol{d}})}{\hat{\boldsymbol{\theta}}(\hat{\boldsymbol{u}}_{+}
(\varsigma) \! - \! \hat{\boldsymbol{d}})}, \label{eqmaininf29} \\
\hat{\aleph}^{\varepsilon_{1}}_{\varepsilon_{2}}(\varsigma) =& \, 
-\dfrac{\hat{\mathfrak{u}}(\varepsilon_{1},\varepsilon_{2},\bm{0};
\varsigma)}{\hat{\boldsymbol{\theta}}(\varepsilon_{1} 
\hat{\boldsymbol{u}}_{+}(\varsigma) \! + \! \varepsilon_{2} \hat{
\boldsymbol{d}})} \! + \! \dfrac{\hat{\mathfrak{u}}(\varepsilon_{1},
\varepsilon_{2},\hat{\boldsymbol{\Omega}};\varsigma)}{\hat{\boldsymbol{
\theta}}(\varepsilon_{1} \hat{\boldsymbol{u}}_{+}(\varsigma) \! - \! 
\frac{1}{2 \pi}((n \! - \! 1)K \! + \! k) \hat{\boldsymbol{\Omega}} \! + 
\! \varepsilon_{2} \hat{\boldsymbol{d}})}, \label{eqmaininf30} \\
\hat{\daleth}^{\varepsilon_{1}}_{\varepsilon_{2}}(\varsigma) =& \, 
-\dfrac{\hat{\mathfrak{v}}(\varepsilon_{1},\varepsilon_{2},\bm{0};
\varsigma)}{\hat{\boldsymbol{\theta}}(\varepsilon_{1} 
\hat{\boldsymbol{u}}_{+}(\varsigma) \! + \! \varepsilon_{2} \hat{
\boldsymbol{d}})} \! + \! \dfrac{\hat{\mathfrak{v}}(\varepsilon_{1},
\varepsilon_{2},\hat{\boldsymbol{\Omega}};\varsigma)}{\hat{\boldsymbol{
\theta}}(\varepsilon_{1} \hat{\boldsymbol{u}}_{+}(\varsigma) \! - \! 
\frac{1}{2 \pi}((n \! - \! 1)K \! + \! k) \hat{\boldsymbol{\Omega}} \! + \! 
\varepsilon_{2} \hat{\boldsymbol{d}})} \! - \! \left(\dfrac{\hat{\mathfrak{u}}
(\varepsilon_{1},\varepsilon_{2},\bm{0};\varsigma)}{\hat{\boldsymbol{\theta}}
(\varepsilon_{1} \hat{\boldsymbol{u}}_{+}(\varsigma) \! + \! \varepsilon_{2} 
\hat{\boldsymbol{d}})} \right)^{2} \nonumber \\
+& \, \dfrac{\hat{\mathfrak{u}}(\varepsilon_{1},\varepsilon_{2},\bm{0};\varsigma) 
\hat{\mathfrak{u}}(\varepsilon_{1},\varepsilon_{2},\hat{\boldsymbol{\Omega}};
\varsigma)}{\hat{\boldsymbol{\theta}}(\varepsilon_{1} \hat{\boldsymbol{u}}_{+}
(\varsigma) \! + \! \varepsilon_{2} \hat{\boldsymbol{d}}) \hat{\boldsymbol{
\theta}}(\varepsilon_{1} \hat{\boldsymbol{u}}_{+}(\varsigma) \! - \! \frac{1}{2 
\pi}((n \! - \! 1)K \! + \! k) \hat{\boldsymbol{\Omega}} \! + \! \varepsilon_{2} 
\hat{\boldsymbol{d}})}, \label{eqmaininf31}
\end{align}
with
\begin{gather}
\hat{\mathfrak{u}}(\varepsilon_{1},\varepsilon_{2},\hat{\boldsymbol{\Omega}};
\varsigma) \! := \! 2 \pi \hat{\Lambda}^{\raise-1.0ex\hbox{$\scriptstyle 1$}}_{0}
(\varepsilon_{1},\varepsilon_{2},\hat{\boldsymbol{\Omega}};\varsigma), \qquad 
\qquad \hat{\mathfrak{v}}(\varepsilon_{1},\varepsilon_{2},\hat{\boldsymbol{\Omega}};
\varsigma) \! := \! -2 \pi^{2} \hat{\Lambda}^{\raise-1.0ex\hbox{$\scriptstyle 2$}}_{0}
(\varepsilon_{1},\varepsilon_{2},\hat{\boldsymbol{\Omega}};\varsigma), 
\label{eqmaininf32} \\
\hat{\Lambda}^{\raise-1.0ex\hbox{$\scriptstyle j_{1}$}}_{0}(\varepsilon_{1},
\varepsilon_{2},\hat{\boldsymbol{\Omega}};\varsigma) \! = \! \sum_{m \in
\mathbb{Z}^{N}}(\hat{\mathfrak{r}}_{1}(\varsigma))^{j_{1}} \me^{2 \pi \mi 
(m,\varepsilon_{1} \hat{\boldsymbol{u}}_{+}(\varsigma)-\frac{1}{2 \pi}((n-1)K
+k) \hat{\boldsymbol{\Omega}}+ \varepsilon_{2} \hat{\boldsymbol{d}})+ 
\mi \pi (m,\hat{\boldsymbol{\tau}}m)}, \quad j_{1} \! = \! 1,2, \label{eqmaininf33} \\
\hat{\mathfrak{r}}_{1}(\varsigma) \! := \! \dfrac{2}{\hat{\leftthreetimes}
(\varsigma)} \sum_{i=1}^{N} \sum_{j=1}^{N}m_{i} \hat{c}_{ij} \varsigma^{N-j}, 
\label{eqmaininf34}
\end{gather}
where $\hat{c}_{i_{1}i_{2}}$, $i_{1},i_{2} \! = \! 1,2,\dotsc,N$, are obtained 
{}from Equations~\eqref{E1} and~\eqref{E2},
\begin{equation}
\hat{\leftthreetimes}(\hat{b}_{0}) \! = \! \mi (-1)^{N} \hat{\eta}_{\hat{b}_{0}}, 
\quad \hat{\leftthreetimes}(\hat{b}_{j}) \! = \! \mi (-1)^{N-j} \hat{\eta}_{
\hat{b}_{j}}, \quad \hat{\leftthreetimes}(\hat{a}_{N+1}) \! = \! \hat{\eta}_{
\hat{a}_{N+1}}, \quad \hat{\leftthreetimes}(\hat{a}_{j}) \! = \! (-1)^{N+1-j} 
\hat{\eta}_{\hat{a}_{j}}, \quad j \! = \! 1,2,\dotsc,N, \label{eqmaininf35}
\end{equation}
with {}\footnote{All square roots are positive.}
\begin{gather}
\hat{\eta}_{\hat{b}_{0}} \! := \! (\hat{a}_{N+1} \! - \! \hat{b}_{0})^{1/2} 
\prod_{m=1}^{N}(\hat{b}_{m} \! - \! \hat{b}_{0})^{1/2}(\hat{a}_{m} \! - \! 
\hat{b}_{0})^{1/2}, \label{eqmaininf55} \\
\hat{\eta}_{\hat{b}_{j}} \! := \! (\hat{b}_{j} \! - \! \hat{a}_{j})^{1/2}
(\hat{a}_{N+1} \! - \! \hat{b}_{j})^{1/2}(\hat{b}_{j} \! - \! \hat{b}_{0})^{1/2} 
\prod_{m=1}^{j-1}(\hat{b}_{j} \! - \! \hat{b}_{m})^{1/2}(\hat{b}_{j} \! - \! 
\hat{a}_{m})^{1/2} \prod_{m^{\prime}=j+1}^{N}(\hat{b}_{m^{\prime}} \! 
- \! \hat{b}_{j})^{1/2}(\hat{a}_{m^{\prime}} \! - \! \hat{b}_{j})^{1/2}, 
\label{eqmaininf56} \\
\hat{\eta}_{\hat{a}_{N+1}} \! := \! (\hat{a}_{N+1} \! - \! \hat{b}_{0})^{1/2} 
\prod_{m=1}^{N}(\hat{a}_{N+1} \! - \! \hat{b}_{m})^{1/2}(\hat{a}_{N+1} 
\! - \! \hat{a}_{m})^{1/2}, \label{eqmaininf57} \\
\hat{\eta}_{\hat{a}_{j}} \! := \! (\hat{b}_{j} \! - \! \hat{a}_{j})^{1/2}
(\hat{a}_{N+1} \! - \! \hat{a}_{j})^{1/2}(\hat{a}_{j} \! - \! \hat{b}_{0})^{1/2} 
\prod_{m=1}^{j-1}(\hat{a}_{j} \! - \! \hat{b}_{m})^{1/2}(\hat{a}_{j} 
\! - \! \hat{a}_{m})^{1/2} \prod_{m^{\prime}=j+1}^{N}(\hat{b}_{m^{\prime}} 
\! - \! \hat{a}_{j})^{1/2}(\hat{a}_{m^{\prime}} \! - \! \hat{a}_{j})^{1/2}, 
\label{eqmaininf58}
\end{gather}
and, for $j \! = \! 1,2,\dotsc,N$,
\begin{gather}
\hat{\mathfrak{Q}}_{0}(\hat{b}_{0}) = -\mi (\hat{a}_{N+1} \! - \! 
\hat{b}_{0})^{-1/2} \prod_{m=1}^{N} \dfrac{(\hat{b}_{m} \! - \! 
\hat{b}_{0})^{1/2}}{(\hat{a}_{m} \! - \! \hat{b}_{0})^{1/2}}, \label{eqmaininf39} \\
\hat{\mathfrak{Q}}_{1}(\hat{b}_{0}) = \dfrac{1}{2} \hat{\mathfrak{Q}}_{0}
(\hat{b}_{0}) \left(\sum_{m=1}^{N} \left(\dfrac{1}{\hat{b}_{0} \! - \! 
\hat{b}_{m}} \! - \! \dfrac{1}{\hat{b}_{0} \! - \! \hat{a}_{m}} \right) \! - \! 
\dfrac{1}{\hat{b}_{0} \! - \! \hat{a}_{N+1}} \right), \label{eqmaininf40} \\
\hat{\mathfrak{Q}}_{0}(\hat{b}_{j}) = -\dfrac{\mi (\hat{b}_{j} \! - \! 
\hat{b}_{0})^{1/2}}{(\hat{a}_{N+1} \! - \! \hat{b}_{j})^{1/2}(\hat{b}_{j} \! 
- \! \hat{a}_{j})^{1/2}} \prod_{m=1}^{j-1} \dfrac{(\hat{b}_{j} \! - \! 
\hat{b}_{m})^{1/2}}{(\hat{b}_{j} \! - \! \hat{a}_{m})^{1/2}} \prod_{
m^{\prime}=j+1}^{N} \dfrac{(\hat{b}_{m^{\prime}} \! - \! 
\hat{b}_{j})^{1/2}}{(\hat{a}_{m^{\prime}} \! - \! \hat{b}_{j})^{1/2}}, 
\label{eqmaininf41} \\
\hat{\mathfrak{Q}}_{1}(\hat{b}_{j}) = \dfrac{1}{2} \hat{\mathfrak{Q}}_{0}
(\hat{b}_{j}) \left(\sum_{\substack{m=1\\m \neq j}}^{N} \left(\dfrac{1}{
\hat{b}_{j} \! - \! \hat{b}_{m}} \! - \! \dfrac{1}{\hat{b}_{j} \! - \! \hat{a}_{m}} 
\right) \! + \! \dfrac{1}{\hat{b}_{j} \! - \! \hat{b}_{0}} \! - \! \dfrac{1}{
\hat{b}_{j} \! - \! \hat{a}_{N+1}} \! - \! \dfrac{1}{\hat{b}_{j} \! - \! 
\hat{a}_{j}} \right), \label{eqmaininf42} \\
\hat{\mathfrak{Q}}_{0}(\hat{a}_{N+1}) = (\hat{a}_{N+1} \! - \! 
\hat{b}_{0})^{1/2} \prod_{m=1}^{N} \dfrac{(\hat{a}_{N+1} \! - \! \hat{b}_{
m})^{1/2}}{(\hat{a}_{N+1} \! - \! \hat{a}_{m})^{1/2}}, \label{eqmaininf43} \\
\hat{\mathfrak{Q}}_{1}(\hat{a}_{N+1}) = \dfrac{1}{2} \hat{\mathfrak{Q}}_{0}
(\hat{a}_{N+1}) \left(\sum_{m=1}^{N} \left(\dfrac{1}{\hat{a}_{N+1} \! - \! 
\hat{b}_{m}} \! - \! \dfrac{1}{\hat{a}_{N+1} \! - \! \hat{a}_{m}} \right) \! + \! 
\dfrac{1}{\hat{a}_{N+1} \! - \! \hat{b}_{0}} \right), \label{eqmaininf44} \\
\hat{\mathfrak{Q}}_{0}(\hat{a}_{j}) = \dfrac{(\hat{a}_{j} \! - \! \hat{b}_{0})^{
1/2}(\hat{b}_{j} \! - \! \hat{a}_{j})^{1/2}}{(\hat{a}_{N+1} \! - \! \hat{a}_{j})^{1/2}} 
\prod_{m=1}^{j-1} \dfrac{(\hat{a}_{j} \! - \! \hat{b}_{m})^{1/2}}{(\hat{a}_{j} 
\! - \! \hat{a}_{m})^{1/2}} \prod_{m^{\prime}=j+1}^{N} 
\dfrac{(\hat{b}_{m^{\prime}} \! - \! \hat{a}_{j})^{1/2}}{(\hat{a}_{m^{\prime}} 
\! - \! \hat{a}_{j})^{1/2}}, \label{eqmaininf45} \\
\hat{\mathfrak{Q}}_{1}(\hat{a}_{j}) = \dfrac{1}{2} \hat{\mathfrak{Q}}_{0}
(\hat{a}_{j}) \left(\sum_{\substack{m=1\\m \neq j}}^{N} \left(\dfrac{1}{
\hat{a}_{j} \! - \! \hat{b}_{m}} \! - \! \dfrac{1}{\hat{a}_{j} \! - \! \hat{a}_{m}} 
\right) \! + \! \dfrac{1}{\hat{a}_{j} \! - \! \hat{b}_{0}} \! - \! \dfrac{1}{
\hat{a}_{j} \! - \! \hat{a}_{N+1}} \! + \! \dfrac{1}{\hat{a}_{j} \! - \! \hat{b}_{j}} 
\right), \label{eqmaiinf46} \\
\hat{\alpha}_{0}(\hat{b}_{0}) = \mi (-1)^{N} \dfrac{2}{3} \hat{h}_{\widetilde{V}}
(\hat{b}_{0}) \hat{\eta}_{\hat{b}_{0}}, \label{eqmaininf47} \\
\hat{\alpha}_{0}(\hat{b}_{j}) = \mi (-1)^{N-j} \dfrac{2}{3} \hat{h}_{\widetilde{V}}
(\hat{b}_{j}) \hat{\eta}_{\hat{b}_{j}}, \label{eqmaininf48} \\
\hat{\alpha}_{1}(\hat{b}_{0}) = \mi (-1)^{N} \hat{\eta}_{\hat{b}_{0}} 
\dfrac{2}{5} \left(\dfrac{1}{2} \hat{h}_{\widetilde{V}}(\hat{b}_{0}) \left(\sum_{m=
1}^{N} \left(\dfrac{1}{\hat{b}_{0} \! - \! \hat{b}_{m}} \! + \! \dfrac{1}{\hat{b}_{0} 
\! - \! \hat{a}_{m}} \right) \! + \! \dfrac{1}{\hat{b}_{0} \! - \! \hat{a}_{N+1}} 
\right) \! + \! (\hat{h}_{\widetilde{V}}(\hat{b}_{0}))^{\prime} \right), 
\label{eqmaininf49} \\
\hat{\alpha}_{1}(\hat{b}_{j}) = \mi (-1)^{N-j} \hat{\eta}_{\hat{b}_{j}} \dfrac{2}{5} 
\left(\dfrac{1}{2} \hat{h}_{\widetilde{V}}(\hat{b}_{j}) \left(\sum_{\substack{m=
1\\m \neq j}}^{N} \left(\dfrac{1}{\hat{b}_{j} \! - \! \hat{b}_{m}} \! + \! \dfrac{1}{
\hat{b}_{j} \! - \! \hat{a}_{m}} \right) \! + \! \dfrac{1}{\hat{b}_{j} \! - \! \hat{a}_{j}} 
\! + \! \dfrac{1}{\hat{b}_{j} \! - \! \hat{a}_{N+1}} \! + \! \dfrac{1}{\hat{b}_{j} 
\! - \! \hat{b}_{0}} \right) \! + \! (\hat{h}_{\widetilde{V}}(\hat{b}_{j}))^{\prime} 
\right), \label{eqmaininf50} \\
\hat{\alpha}_{0}(\hat{a}_{N+1}) = \dfrac{2}{3} \hat{h}_{\widetilde{V}}(\hat{a}_{N
+1}) \hat{\eta}_{\hat{a}_{N+1}}, \label{eqmaininf51} \\
\hat{\alpha}_{0}(\hat{a}_{j}) = (-1)^{N+1-j} \dfrac{2}{3} \hat{h}_{\widetilde{V}}
(\hat{a}_{j}) \hat{\eta}_{\hat{a}_{j}}, \label{eqmaininf52} \\
\hat{\alpha}_{1}(\hat{a}_{N+1}) = \hat{\eta}_{\hat{a}_{N+1}} \dfrac{2}{5} \left(
\dfrac{1}{2} \hat{h}_{\widetilde{V}}(\hat{a}_{N+1}) \left(\sum_{m=1}^{N} \left(
\dfrac{1}{\hat{a}_{N+1} \! - \! \hat{b}_{m}} \! + \! \dfrac{1}{\hat{a}_{N+1} 
\! - \! \hat{a}_{m}} \right) \! + \! \dfrac{1}{\hat{a}_{N+1} \! - \! \hat{b}_{0}} 
\right) \! + \! (\hat{h}_{\widetilde{V}}(\hat{a}_{N+1}))^{\prime} \right), 
\label{eqmaininf53} \\
\hat{\alpha}_{1}(\hat{a}_{j}) = (-1)^{N+1-j} \hat{\eta}_{\hat{a}_{j}} \dfrac{2}{5} 
\left(\dfrac{1}{2} \hat{h}_{\widetilde{V}}(\hat{a}_{j}) \left(\sum_{\substack{m
=1\\m \neq j}}^{N} \left(\dfrac{1}{\hat{a}_{j} \! - \! \hat{b}_{m}} \! + \! 
\dfrac{1}{\hat{a}_{j} \! - \! \hat{a}_{m}} \right) \! + \! \dfrac{1}{\hat{a}_{j} 
\! - \! \hat{b}_{j}} \! + \! \dfrac{1}{\hat{a}_{j} \! - \! \hat{a}_{N+1}} \! + \! 
\dfrac{1}{\hat{a}_{j} \! - \! \hat{b}_{0}} \right) \! + \! (\hat{h}_{\widetilde{V}}
(\hat{a}_{j}))^{\prime} \right); \label{eqmaininf54}
\end{gather}
{\rm \pmb{(2)}} for $z \! \in \! \hat{\Upsilon}_{2}$,
\begin{align}
\pmb{\pi}_{k}^{n}(z) \underset{\underset{z_{o}=1+o(1)}{\mathscr{N},
n \to \infty}}{=}& \, \left(
\tilde{\mathfrak{m}}^{\raise-0.5ex\hbox{$\scriptstyle \infty$}}_{11} 
\hat{\mathbb{M}}_{12}(z) \left(1 \! + \! \dfrac{1}{(n \! - \! 1)K \! + \! k} 
\hat{\mathcal{R}}_{11}^{\sharp}(z) \! + \! \mathcal{O} \left(\dfrac{
\mathfrak{c}(n,k,z_{o})}{((n \! - \! 1)K \! + \! k)^{2}} \right) \right) 
\right. \nonumber \\
+&\left. \, 
\tilde{\mathfrak{m}}^{\raise-0.5ex\hbox{$\scriptstyle \infty$}}_{22} \hat{
\mathbb{M}}_{22}(z) \left(\dfrac{1}{(n \! - \! 1)K \! + \! k} \hat{\mathcal{
R}}_{12}^{\sharp}(z) \! + \! \mathcal{O} \left(\dfrac{\mathfrak{c}(n,k,
z_{o})}{((n \! - \! 1)K \! + \! k)^{2}} \right) \right) \right) \nonumber \\
\times& \, \me^{n(g^{\infty}(z)-\tilde{\mathscr{P}}_{0})}, 
\label{eqmaininf59}
\end{align}
and
\begin{align}
\int_{\mathbb{R}} \dfrac{\pmb{\pi}_{k}^{n}(\xi) 
\me^{-n \widetilde{V}(\xi)}}{\xi \! - \! z} \, \dfrac{\md \xi}{2 \pi \mi} 
\underset{\underset{z_{o}=1+o(1)}{\mathscr{N},n \to \infty}}{=}& \, 
-\left(\tilde{\mathfrak{m}}^{\raise-0.5ex\hbox{$\scriptstyle \infty$}}_{11} 
\hat{\mathbb{M}}_{11}(z) \left(1 \! + \! \dfrac{1}{(n \! - \! 1)K \! + \! k} 
\hat{\mathcal{R}}_{11}^{\sharp}(z) \! + \! \mathcal{O} \left(\dfrac{
\mathfrak{c}(n,k,z_{o})}{((n \! - \! 1)K \! + \! k)^{2}} \right) \right) 
\right. \nonumber \\
+&\left. \, 
\tilde{\mathfrak{m}}^{\raise-0.5ex\hbox{$\scriptstyle \infty$}}_{22} \hat{
\mathbb{M}}_{21}(z) \left(\dfrac{1}{(n \! - \! 1)K \! + \! k} \hat{\mathcal{
R}}_{12}^{\sharp}(z) \! + \! \mathcal{O} \left(\dfrac{\mathfrak{c}(n,k,
z_{o})}{((n \! - \! 1)K \! + \! k)^{2}} \right) \right) \right) \nonumber \\
\times& \, \me^{n \hat{\ell}} \me^{-n(g^{\infty}(z)-
\tilde{\mathscr{P}}_{0})}; \label{eqmaininf60}
\end{align}
{\rm \pmb{(3)}} for $z \! \in \! \hat{\Upsilon}_{3}$,
\begin{align}
\pmb{\pi}_{k}^{n}(z) \underset{\underset{z_{o}=1+o(1)}{\mathscr{N},
n \to \infty}}{=}& \, \left(
\tilde{\mathfrak{m}}^{\raise-0.5ex\hbox{$\scriptstyle \infty$}}_{11} \left(
\hat{\mathbb{M}}_{11}(z) \! + \! \hat{\mathbb{M}}_{12}(z) \me^{-2 \pi \mi 
((n-1)K+k) \int_{z}^{\hat{a}_{N+1}} \psi_{\widetilde{V}}^{\infty}(\xi) \, \md 
\xi} \right) \left(1 \! + \! \dfrac{1}{(n \! - \! 1)K \! + \! k} \hat{\mathcal{
R}}_{11}^{\sharp}(z) \! + \! \mathcal{O} \left(\dfrac{\mathfrak{c}(n,k,
z_{o})}{((n \! - \! 1)K \! + \! k)^{2}} \right) \right) \right. \nonumber \\
+&\left. \, 
\tilde{\mathfrak{m}}^{\raise-0.5ex\hbox{$\scriptstyle \infty$}}_{22} \left(
\hat{\mathbb{M}}_{21}(z) \! + \! \hat{\mathbb{M}}_{22}(z) \me^{-2 \pi \mi 
((n-1)K+k) \int_{z}^{\hat{a}_{N+1}} \psi_{\widetilde{V}}^{\infty}(\xi) \, \md 
\xi} \right) \left(\dfrac{1}{(n \! - \! 1)K \! + \! k} \hat{\mathcal{R}}_{12}^{
\sharp}(z) \! + \! \mathcal{O} \left(\dfrac{\mathfrak{c}(n,k,z_{o})}{((n 
\! - \! 1)K \! + \! k)^{2}} \right) \right) \right) \nonumber \\
\times& \, \me^{n(g^{\infty}(z)-\tilde{\mathscr{P}}_{0})}, 
\label{eqmaininf61}
\end{align}
and
\begin{align}
\int_{\mathbb{R}} \dfrac{\pmb{\pi}_{k}^{n}(\xi) 
\me^{-n \widetilde{V}(\xi)}}{\xi \! - \! z} \, \dfrac{\md \xi}{2 \pi \mi} 
\underset{\underset{z_{o}=1+o(1)}{\mathscr{N},n \to \infty}}{=}& \,
\left(\tilde{\mathfrak{m}}^{\raise-0.5ex\hbox{$\scriptstyle \infty$}}_{11} 
\hat{\mathbb{M}}_{12}(z) \left(1 \! + \! \dfrac{1}{(n \! - \! 1)K \! + \! k} 
\hat{\mathcal{R}}_{11}^{\sharp}(z) \! + \! \mathcal{O} \left(\dfrac{
\mathfrak{c}(n,k,z_{o})}{((n \! - \! 1)K \! + \! k)^{2}} \right) \right) 
\right. \nonumber \\
+&\left. \, 
\tilde{\mathfrak{m}}^{\raise-0.5ex\hbox{$\scriptstyle \infty$}}_{22} \hat{
\mathbb{M}}_{22}(z) \left(\dfrac{1}{(n \! - \! 1)K \! + \! k} \hat{\mathcal{
R}}_{12}^{\sharp}(z) \! + \! \mathcal{O} \left(\dfrac{\mathfrak{c}(n,k,
z_{o})}{((n \! - \! 1)K \! + \! k)^{2}} \right) \right) \right) \nonumber \\
\times& \, \me^{n \hat{\ell}} \me^{-n(g^{\infty}(z)-
\tilde{\mathscr{P}}_{0})}; \label{eqmaininf62}
\end{align}
{\rm \pmb{(4)}} for $z \! \in \! \hat{\Upsilon}_{4}$,
\begin{align}
\pmb{\pi}_{k}^{n}(z) \underset{\underset{z_{o}=1+o(1)}{\mathscr{N},
n \to \infty}}{=}& \, \left(
\tilde{\mathfrak{m}}^{\raise-0.5ex\hbox{$\scriptstyle \infty$}}_{11} \left(
\hat{\mathbb{M}}_{12}(z) \! + \! \hat{\mathbb{M}}_{11}(z) \me^{2 \pi \mi 
((n-1)K+k) \int_{z}^{\hat{a}_{N+1}} \psi_{\widetilde{V}}^{\infty}(\xi) \, \md 
\xi} \right) \left(1 \! + \! \dfrac{1}{(n \! - \! 1)K \! + \! k} \hat{\mathcal{
R}}_{11}^{\sharp}(z) \! + \! \mathcal{O} \left(\dfrac{\mathfrak{c}(n,k,
z_{o})}{((n \! - \! 1)K \! + \! k)^{2}} \right) \right) \right. \nonumber \\
+&\left. \, 
\tilde{\mathfrak{m}}^{\raise-0.5ex\hbox{$\scriptstyle \infty$}}_{22} \left(
\hat{\mathbb{M}}_{22}(z) \! + \! \hat{\mathbb{M}}_{21}(z) \me^{2 \pi \mi 
((n-1)K+k) \int_{z}^{\hat{a}_{N+1}} \psi_{\widetilde{V}}^{\infty}(\xi) \, \md 
\xi} \right) \left(\dfrac{1}{(n \! - \! 1)K \! + \! k} \hat{\mathcal{R}}_{12}^{
\sharp}(z) \! + \! \mathcal{O} \left(\dfrac{\mathfrak{c}(n,k,z_{o})}{((n 
\! - \! 1)K \! + \! k)^{2}} \right) \right) \right) \nonumber \\
\times& \, \me^{n(g^{\infty}(z)-\tilde{\mathscr{P}}_{0})}, 
\label{eqmaininf63}
\end{align}
and
\begin{align}
\int_{\mathbb{R}} \dfrac{\pmb{\pi}_{k}^{n}(\xi) 
\me^{-n \widetilde{V}(\xi)}}{\xi \! - \! z} \, \dfrac{\md \xi}{2 \pi \mi} 
\underset{\underset{z_{o}=1+o(1)}{\mathscr{N},n \to \infty}}{=}& \, 
-\left(\tilde{\mathfrak{m}}^{\raise-0.5ex\hbox{$\scriptstyle \infty$}}_{11} 
\hat{\mathbb{M}}_{11}(z) \left(1 \! + \! \dfrac{1}{(n \! - \! 1)K \! + \! k} 
\hat{\mathcal{R}}_{11}^{\sharp}(z) \! + \! \mathcal{O} \left(\dfrac{
\mathfrak{c}(n,k,z_{o})}{((n \! - \! 1)K \! + \! k)^{2}} \right) \right) 
\right. \nonumber \\
+&\left. \, 
\tilde{\mathfrak{m}}^{\raise-0.5ex\hbox{$\scriptstyle \infty$}}_{22} \hat{
\mathbb{M}}_{21}(z) \left(\dfrac{1}{(n \! - \! 1)K \! + \! k} \hat{\mathcal{
R}}_{12}^{\sharp}(z) \! + \! \mathcal{O} \left(\dfrac{\mathfrak{c}(n,k,
z_{o})}{((n \! - \! 1)K \! + \! k)^{2}} \right) \right) \right) \nonumber \\
\times& \, \me^{n \hat{\ell}} \me^{-n(g^{\infty}(z)-
\tilde{\mathscr{P}}_{0})}; \label{eqmaininf64}
\end{align}
{\rm \pmb{(5)}} for $z \! \in \! \hat{\Omega}^{1}_{\hat{b}_{j-1}}$, 
$j \! = \! 1,2,\dotsc,N \! + \! 1$,
\begin{align}
\pmb{\pi}_{k}^{n}(z) \underset{\underset{z_{o}=1+o(1)}{
\mathscr{N},n \to \infty}}{=}& \, \left(\hat{\mathcal{X}}_{11}^{\hat{b},1}(z) 
\left(1 \! + \! \dfrac{1}{(n \! - \! 1)K \! + \! k} \left(\hat{\mathcal{R}}_{11}^{
\sharp}(z) \! - \! \hat{\mathcal{R}}_{11}^{\natural}(z) \right) \! + \! \mathcal{O} 
\left(\dfrac{\mathfrak{c}(n,k,z_{o})}{((n \! - \! 1)K \! + \! k)^{2}} \right) \right) 
\right. \nonumber \\
+&\left. \, \hat{\mathcal{X}}_{21}^{\hat{b},1}(z) \left(\dfrac{1}{(n \! - \! 1)K \! 
+ \! k} \left(\hat{\mathcal{R}}_{12}^{\sharp}(z) \! - \! \hat{\mathcal{R}}_{12}^{
\natural}(z) \right) \! + \! \mathcal{O} \left(\dfrac{\mathfrak{c}(n,k,z_{o})}{((n 
\! - \! 1)K \! + \! k)^{2}} \right) \right) \right) \nonumber \\
\times& \, \me^{n(g^{\infty}(z)-\tilde{\mathscr{P}}_{0})}, 
\label{eqmaininf65}
\end{align}
and
\begin{align}
\int_{\mathbb{R}} \dfrac{\pmb{\pi}_{k}^{n}(\xi) 
\me^{-n \widetilde{V}(\xi)}}{\xi \! - \! z} \, \dfrac{\md \xi}{2 \pi \mi} 
\underset{\underset{z_{o}=1+o(1)}{\mathscr{N},n \to \infty}}{=}& \, 
\left(\hat{\mathcal{X}}_{12}^{\hat{b},1}(z) \left(1 \! + \! \dfrac{1}{(n \! 
- \! 1)K \! + \! k} \left(\hat{\mathcal{R}}_{11}^{\sharp}(z) \! - \! \hat{
\mathcal{R}}_{11}^{\natural}(z) \right) \! + \! \mathcal{O} \left(\dfrac{
\mathfrak{c}(n,k,z_{o})}{((n \! - \! 1)K \! + \! k)^{2}} \right) \right) 
\right. \nonumber \\
+&\left. \, \hat{\mathcal{X}}_{22}^{\hat{b},1}(z) \left(\dfrac{1}{(n \! - 
\! 1)K \! + \! k} \left(\hat{\mathcal{R}}_{12}^{\sharp}(z) \! - \! \hat{
\mathcal{R}}_{12}^{\natural}(z) \right) \! + \! \mathcal{O} \left(\dfrac{
\mathfrak{c}(n,k,z_{o})}{((n \! - \! 1)K \! + \! k)^{2}} \right) \right) 
\right) \nonumber \\
\times& \, \me^{n \hat{\ell}} \me^{-n(g^{\infty}(z)-
\tilde{\mathscr{P}}_{0})}, \label{eqmaininf66}
\end{align}
where
\begin{equation} \label{eqmaininf67} 
\hat{\mathcal{R}}^{\natural}(z) \! = \! \sum_{j=1}^{N+1} \left(\hat{\mathbb{
Y}}_{\hat{b}_{j-1}}(z) \chi_{\hat{\mathbb{U}}_{\hat{\delta}_{\hat{b}_{j-1}}}}(z) 
\! + \! \hat{\mathbb{Y}}_{\hat{a}_{j}}(z) \chi_{\hat{\mathbb{U}}_{\hat{\delta}_{
\hat{a}_{j}}}}(z) \right),
\end{equation}
with $\chi_{\hat{\mathbb{U}}_{\hat{\delta}_{\hat{b}_{j-1}}}}(z)$ (resp., 
$\chi_{\hat{\mathbb{U}}_{\hat{\delta}_{\hat{a}_{j}}}}(z))$ the characteristic 
function of the set $\hat{\mathbb{U}}_{\hat{\delta}_{\hat{b}_{j-1}}}$ (resp., 
$\hat{\mathbb{U}}_{\hat{\delta}_{\hat{a}_{j}}})$,
\begin{gather}
\hat{\mathbb{Y}}_{\hat{b}_{j-1}}(z) := \dfrac{1}{\hat{\xi}_{\hat{b}_{j-1}}
(z)} \tilde{\mathfrak{m}}^{\raise-0.5ex\hbox{$\scriptstyle \infty$}} 
\hat{\mathbb{M}}(z) 
\begin{pmatrix}
-(s_{1} \! + \! t_{1}) & -\mi (s_{1} \! - \! t_{1}) \me^{\mi ((n-1)K+k) 
\hat{\mho}_{j-1}} \\
-\mi (s_{1} \! - \! t_{1}) \me^{-\mi ((n-1)K+k) \hat{\mho}_{j-1}} & 
(s_{1} \! + \! t_{1})
\end{pmatrix}
(\tilde{\mathfrak{m}}^{\raise-0.5ex\hbox{$\scriptstyle \infty$}} 
\hat{\mathbb{M}}(z))^{-1}, \label{eqmaininf68} \\
\hat{\mathbb{Y}}_{\hat{a}_{j}}(z) := \dfrac{1}{\hat{\xi}_{\hat{a}_{j}}
(z)} \tilde{\mathfrak{m}}^{\raise-0.5ex\hbox{$\scriptstyle \infty$}} 
\hat{\mathbb{M}}(z) 
\begin{pmatrix}
-(s_{1} \! + \! t_{1}) & \mi (s_{1} \! - \! t_{1}) \me^{\mi ((n-1)K+k) 
\hat{\mho}_{j}} \\
\mi (s_{1} \! - \! t_{1}) \me^{-\mi ((n-1)K+k) \hat{\mho}_{j}} & 
(s_{1} \! + \! t_{1})
\end{pmatrix}
(\tilde{\mathfrak{m}}^{\raise-0.5ex\hbox{$\scriptstyle \infty$}} 
\hat{\mathbb{M}}(z))^{-1}, \label{eqmaininf69}
\end{gather}
where
\begin{gather}
\hat{\xi}_{\hat{b}_{j-1}}(z) \! = \! -\int_{z}^{\hat{b}_{j-1}}(\hat{R}(\xi))^{1/2} 
\hat{h}_{\widetilde{V}}(\xi) \, \md \xi, \label{eqmaininf70} \\
\hat{\xi}_{\hat{a}_{j}}(z) \! = \! \int_{\hat{a}_{j}}^{z}(\hat{R}(\xi))^{1/2} 
\hat{h}_{\widetilde{V}}(\xi) \, \md \xi, \label{eqmaininf71}
\end{gather}
and
\begin{align}
\hat{\mathcal{X}}_{11}^{\hat{b},1}(z) =& \, -\mi \sqrt{\pi} \me^{\frac{1}{2} 
((n-1)K+k) \hat{\xi}_{\hat{b}_{j-1}}(z)} 
\tilde{\mathfrak{m}}^{\raise-0.5ex\hbox{$\scriptstyle \infty$}}_{11} \left(
\mi \hat{\mathbb{M}}_{11}(z) \Theta_{0,0}^{0,0,-}(\hat{\Phi}_{\hat{b}_{j-1}}(z)) 
\! - \! \hat{\mathbb{M}}_{12}(z) \Theta_{0,0}^{0,0,+}(\hat{\Phi}_{\hat{b}_{j-1}}
(z)) \me^{-\mi ((n-1)K+k) \hat{\mho}_{j-1}} \right), \label{eqmaininf72} \\
\hat{\mathcal{X}}_{12}^{\hat{b},1}(z) =& \, \sqrt{\pi} \me^{-\frac{\mi \pi}{6}} 
\me^{-\frac{1}{2} ((n-1)K+k) \hat{\xi}_{\hat{b}_{j-1}}(z)} 
\tilde{\mathfrak{m}}^{\raise-0.5ex\hbox{$\scriptstyle \infty$}}_{11} \left(
\hat{\mathbb{M}}_{12}(z) \Theta_{2,2}^{0,2,+}(\hat{\Phi}_{\hat{b}_{j-1}}(z)) 
\! - \! \mi \hat{\mathbb{M}}_{11}(z) \Theta_{2,2}^{0,2,-}(\hat{\Phi}_{\hat{b}_{j-1}}
(z)) \me^{\mi ((n-1)K+k) \hat{\mho}_{j-1}} \right), \label{eqmaininf73} \\
\hat{\mathcal{X}}_{21}^{\hat{b},1}(z) =& \, -\mi \sqrt{\pi} \me^{\frac{1}{2} 
((n-1)K+k) \hat{\xi}_{\hat{b}_{j-1}}(z)} 
\tilde{\mathfrak{m}}^{\raise-0.5ex\hbox{$\scriptstyle \infty$}}_{22} \left(
\mi \hat{\mathbb{M}}_{21}(z) \Theta_{0,0}^{0,0,-}(\hat{\Phi}_{\hat{b}_{j-1}}(z)) 
\! - \! \hat{\mathbb{M}}_{22}(z) \Theta_{0,0}^{0,0,+}(\hat{\Phi}_{\hat{b}_{j-1}}
(z)) \me^{-\mi ((n-1)K+k) \hat{\mho}_{j-1}} \right), \label{eqmaininf74} \\
\hat{\mathcal{X}}_{22}^{\hat{b},1}(z) =& \, \sqrt{\pi} \me^{-\frac{\mi \pi}{6}} 
\me^{-\frac{1}{2} ((n-1)K+k) \hat{\xi}_{\hat{b}_{j-1}}(z)} 
\tilde{\mathfrak{m}}^{\raise-0.5ex\hbox{$\scriptstyle \infty$}}_{22} \left(
\hat{\mathbb{M}}_{22}(z) \Theta_{2,2}^{0,2,+}(\hat{\Phi}_{\hat{b}_{j-1}}(z)) \! - 
\! \mi \hat{\mathbb{M}}_{21}(z) \Theta_{2,2}^{0,2,-}(\hat{\Phi}_{\hat{b}_{j-1}}
(z)) \me^{\mi ((n-1)K+k) \hat{\mho}_{j-1}} \right), \label{eqmaininf75}
\end{align}
where $\hat{\Phi}_{\hat{b}_{j-1}}(z)$ is defined by Equation~\eqref{eqmaininf77}, 
and
\begin{equation} \label{eqmaininf76} 
\Theta_{r_{3},r_{4}}^{r_{1},r_{2},\pm}(\zeta) \! := \! \zeta^{1/4} \omega^{r_{1}} 
\operatorname{Ai}(\omega^{r_{2}} \zeta) \! \pm \! \zeta^{-1/4} \omega^{r_{3}} 
\operatorname{Ai}^{\prime}(\omega^{r_{4}} \zeta), \quad r_{1},r_{2},r_{3},r_{4} 
\! = \! 0,1,2,
\end{equation}
with $\omega \! = \! \exp (2 \pi \mi/3)$$;$\\
{\rm \pmb{(6)}} for $z \! \in \! \hat{\Omega}_{\hat{a}_{j}}^{1}$, $j \! = \! 1,2,
\dotsc,N \! + \!1$,
\begin{align}
\pmb{\pi}_{k}^{n}(z) \underset{\underset{z_{o}=1+o(1)}{\mathscr{N},
n \to \infty}}{=}& \, \left(\hat{\mathcal{X}}_{11}^{\hat{a},1}(z) \left(1 \! + \! 
\dfrac{1}{(n \! - \! 1)K \! + \! k} \left(\hat{\mathcal{R}}_{11}^{\sharp}(z) \! - \! 
\hat{\mathcal{R}}_{11}^{\natural}(z) \right) \! + \! \mathcal{O} \left(\dfrac{
\mathfrak{c}(n,k,z_{o})}{((n \! - \! 1)K \! + \! k)^{2}} \right) \right) \right. 
\nonumber \\
+&\left. \, \hat{\mathcal{X}}_{21}^{\hat{a},1}(z) \left(\dfrac{1}{(n \! - \! 1)K \! 
+ \! k} \left(\hat{\mathcal{R}}_{12}^{\sharp}(z) \! - \! \hat{\mathcal{R}}_{12}^{
\natural}(z) \right) \! + \! \mathcal{O} \left(\dfrac{\mathfrak{c}(n,k,z_{o})}{((n 
\! - \! 1)K \! + \! k)^{2}} \right) \right) \right) \nonumber \\
\times& \, \me^{n(g^{\infty}(z)-\tilde{\mathscr{P}}_{0})}, \label{eqmaininf78}
\end{align}
and
\begin{align}
\int_{\mathbb{R}} \dfrac{\pmb{\pi}_{k}^{n}(\xi) 
\me^{-n \widetilde{V}(\xi)}}{\xi \! - \! z} \, \dfrac{\md \xi}{2 \pi \mi} 
\underset{\underset{z_{o}=1+o(1)}{\mathscr{N},n \to \infty}}{=}& \, \left(
\hat{\mathcal{X}}_{12}^{\hat{a},1}(z) \left(1 \! + \! \dfrac{1}{(n \! - \! 1)K \! + 
\! k} \left(\hat{\mathcal{R}}_{11}^{\sharp}(z) \! - \! \hat{\mathcal{R}}_{11}^{
\natural}(z) \right) \! + \! \mathcal{O} \left(\dfrac{\mathfrak{c}(n,k,z_{o})}{((n 
\! - \! 1)K \! + \! k)^{2}} \right) \right) \right. \nonumber \\
+&\left. \, \hat{\mathcal{X}}_{22}^{\hat{a},1}(z) \left(\dfrac{1}{(n \! - \! 1)K \! 
+ \! k} \left(\hat{\mathcal{R}}_{12}^{\sharp}(z) \! - \! \hat{\mathcal{R}}_{12}^{
\natural}(z) \right) \! + \! \mathcal{O} \left(\dfrac{\mathfrak{c}(n,k,z_{o})}{((
n \! - \! 1)K \! + \! k)^{2}} \right) \right) \right) \nonumber \\
\times& \, \me^{n \hat{\ell}} \me^{-n(g^{\infty}(z)-\tilde{\mathscr{P}}_{0})}, 
\label{eqmaininf79}
\end{align}
where
\begin{align}
\hat{\mathcal{X}}_{11}^{\hat{a},1}(z) =& \, -\mi \sqrt{\pi} \me^{\frac{1}{2} 
((n-1)K+k) \hat{\xi}_{\hat{a}_{j}}(z)} 
\tilde{\mathfrak{m}}^{\raise-0.5ex\hbox{$\scriptstyle \infty$}}_{11} \left(\mi 
\hat{\mathbb{M}}_{11}(z) \Theta_{0,0}^{0,0,-}(\hat{\Phi}_{\hat{a}_{j}}(z)) \! + 
\! \hat{\mathbb{M}}_{12}(z) \Theta_{0,0}^{0,0,+}(\hat{\Phi}_{\hat{a}_{j}}(z)) 
\me^{-\mi ((n-1)K+k) \hat{\mho}_{j}} \right), \label{eqmaininf80} \\
\hat{\mathcal{X}}_{12}^{\hat{a},1}(z) =& \, \sqrt{\pi} \me^{-\frac{\mi \pi}{6}} 
\me^{-\frac{1}{2} ((n-1)K+k) \hat{\xi}_{\hat{a}_{j}}(z)} 
\tilde{\mathfrak{m}}^{\raise-0.5ex\hbox{$\scriptstyle \infty$}}_{11} \left(
\hat{\mathbb{M}}_{12}(z) \Theta_{2,2}^{0,2,+}(\hat{\Phi}_{\hat{a}_{j}}(z)) \! + \! 
\mi \hat{\mathbb{M}}_{11}(z) \Theta_{2,2}^{0,2,-}(\hat{\Phi}_{\hat{a}_{j}}(z)) 
\me^{\mi ((n-1)K+k) \hat{\mho}_{j}} \right), \label{eqmaininf81} \\
\hat{\mathcal{X}}_{21}^{\hat{a},1}(z) =& \, -\mi \sqrt{\pi} \me^{\frac{1}{2} 
((n-1)K+k) \hat{\xi}_{\hat{a}_{j}}(z)} 
\tilde{\mathfrak{m}}^{\raise-0.5ex\hbox{$\scriptstyle \infty$}}_{22} \left(
\mi \hat{\mathbb{M}}_{21}(z) \Theta_{0,0}^{0,0,-}(\hat{\Phi}_{\hat{a}_{j}}(z)) \! 
+ \! \hat{\mathbb{M}}_{22}(z) \Theta_{0,0}^{0,0,+}(\hat{\Phi}_{\hat{a}_{j}}(z)) 
\me^{-\mi ((n-1)K+k) \hat{\mho}_{j}} \right), \label{eqmaininf82} \\
\hat{\mathcal{X}}_{22}^{\hat{a},1}(z) =& \, \sqrt{\pi} \me^{-\frac{\mi \pi}{6}} 
\me^{-\frac{1}{2} ((n-1)K+k) \hat{\xi}_{\hat{a}_{j}}(z)} 
\tilde{\mathfrak{m}}^{\raise-0.5ex\hbox{$\scriptstyle \infty$}}_{22} \left(
\hat{\mathbb{M}}_{22}(z) \Theta_{2,2}^{0,2,+}(\hat{\Phi}_{\hat{a}_{j}}(z)) \! + 
\! \mi \hat{\mathbb{M}}_{21}(z) \Theta_{2,2}^{0,2,-}(\hat{\Phi}_{\hat{a}_{j}}(z)) 
\me^{\mi ((n-1)K+k) \hat{\mho}_{j}} \right), \label{eqmaininf83}
\end{align}
with $\hat{\Phi}_{\hat{a}_{j}}(z)$ defined by Equation~\eqref{eqmaininf84}$;$\\
{\rm \pmb{(7)}} for $z \! \in \! \hat{\Omega}^{2}_{\hat{b}_{j-1}}$, $j \! = \! 
1,2,\dotsc,N \! + \!1$,
\begin{align}
\pmb{\pi}_{k}^{n}(z) \underset{\underset{z_{o}=1+o(1)}{\mathscr{N},
n \to \infty}}{=}& \, \left(\left(\hat{\mathcal{X}}_{11}^{\hat{b},2}(z) \! + \! 
\hat{\mathcal{X}}_{12}^{\hat{b},2}(z) \me^{-2 \pi \mi ((n-1)K+k) \int_{z}^{
\hat{a}_{N+1}} \psi_{\widetilde{V}}^{\infty}(\xi) \, \md \xi} \right) \left(1 \! + 
\! \dfrac{1}{(n \! - \! 1)K \! + \! k} \left(\hat{\mathcal{R}}_{11}^{\sharp}(z) 
\! - \! \hat{\mathcal{R}}_{11}^{\natural}(z) \right) \! + \! \mathcal{O} \left(
\dfrac{\mathfrak{c}(n,k,z_{o})}{((n \! - \! 1)K \! + \! k)^{2}} \right) \right) 
\right. \nonumber \\
+&\left. \, \left(\hat{\mathcal{X}}_{21}^{\hat{b},2}(z) \! + \! \hat{\mathcal{
X}}_{22}^{\hat{b},2}(z) \me^{-2 \pi \mi ((n-1)K+k) \int_{z}^{\hat{a}_{N+1}} 
\psi_{\widetilde{V}}^{\infty}(\xi) \, \md \xi} \right) \left(\dfrac{1}{(n \! - \! 
1)K \! + \! k} \left(\hat{\mathcal{R}}_{12}^{\sharp}(z) \! - \! \hat{\mathcal{
R}}_{12}^{\natural}(z) \right) \! + \! \mathcal{O} \left(\dfrac{\mathfrak{c}
(n,k,z_{o})}{((n \! - \! 1)K \! + \! k)^{2}} \right) \right) \right) \nonumber \\
\times& \, \me^{n(g^{\infty}(z)-\tilde{\mathscr{P}}_{0})}, \label{eqmaininf85}
\end{align}
and
\begin{align}
\int_{\mathbb{R}} \dfrac{\pmb{\pi}_{k}^{n}(\xi) 
\me^{-n \widetilde{V}(\xi)}}{\xi \! - \! z} \, \dfrac{\md \xi}{2 \pi \mi} 
\underset{\underset{z_{o}=1+o(1)}{\mathscr{N},n \to \infty}}{=}& \,
\left(\hat{\mathcal{X}}_{12}^{\hat{b},2}(z) \left(1 \! + \! \dfrac{1}{(n \! - \! 
1)K \! + \! k} \left(\hat{\mathcal{R}}_{11}^{\sharp}(z) \! - \! \hat{\mathcal{
R}}_{11}^{\natural}(z) \right) \! + \! \mathcal{O} \left(\dfrac{\mathfrak{c}
(n,k,z_{o})}{((n \! - \! 1)K \! + \! k)^{2}} \right) \right) \right. \nonumber \\
+&\left. \, \hat{\mathcal{X}}_{22}^{\hat{b},2}(z) \left(\dfrac{1}{(n \! - \! 1)K \! 
+ \! k} \left(\hat{\mathcal{R}}_{12}^{\sharp}(z) \! - \! \hat{\mathcal{R}}_{12}^{
\natural}(z) \right) \! + \! \mathcal{O} \left(\dfrac{\mathfrak{c}(n,k,z_{o})}{(
(n \! - \! 1)K \! + \! k)^{2}} \right) \right) \right) \nonumber \\
\times& \, \me^{n \hat{\ell}} \me^{-n(g^{\infty}(z)-\tilde{\mathscr{P}}_{0})}, 
\label{eqmaininf86}
\end{align}
where
\begin{align}
\hat{\mathcal{X}}_{11}^{\hat{b},2}(z) =& \, -\mi \sqrt{\pi} \me^{\frac{1}{2} 
((n-1)K+k) \hat{\xi}_{\hat{b}_{j-1}}(z)} 
\tilde{\mathfrak{m}}^{\raise-0.5ex\hbox{$\scriptstyle \infty$}}_{11} \left(\mi 
\hat{\mathbb{M}}_{11}(z) \varTheta_{2,0}^{0,2,-}(\hat{\Phi}_{\hat{b}_{j-1}}(z)) \! 
- \! \hat{\mathbb{M}}_{12}(z) \varTheta_{2,0}^{0,2,+}(\hat{\Phi}_{\hat{b}_{j-1}}
(z)) \me^{-\mi ((n-1)K+k) \hat{\mho}_{j-1}} \right), \label{eqmaininf87} \\
\hat{\mathcal{X}}_{12}^{\hat{b},2}(z) =& \, \sqrt{\pi} \me^{-\frac{\mi \pi}{6}} 
\me^{-\frac{1}{2} ((n-1)K+k) \hat{\xi}_{\hat{b}_{j-1}}(z)} 
\tilde{\mathfrak{m}}^{\raise-0.5ex\hbox{$\scriptstyle \infty$}}_{11} \left(
\hat{\mathbb{M}}_{12}(z) \Theta_{2,2}^{0,2,+}(\hat{\Phi}_{\hat{b}_{j-1}}(z)) \! 
- \! \mi \hat{\mathbb{M}}_{11}(z) \Theta_{2,2}^{0,2,-}(\hat{\Phi}_{\hat{b}_{j-
1}}(z)) \me^{\mi ((n-1)K+k) \hat{\mho}_{j-1}} \right), \label{eqmaininf88} \\
\hat{\mathcal{X}}_{21}^{\hat{b},2}(z) =& \, -\mi \sqrt{\pi} \me^{\frac{1}{2} 
((n-1)K+k) \hat{\xi}_{\hat{b}_{j-1}}(z)} 
\tilde{\mathfrak{m}}^{\raise-0.5ex\hbox{$\scriptstyle \infty$}}_{22} \left(
\mi \hat{\mathbb{M}}_{21}(z) \varTheta_{2,0}^{0,2,-}(\hat{\Phi}_{\hat{b}_{j-1}}
(z)) \! - \! \hat{\mathbb{M}}_{22}(z) \varTheta_{2,0}^{0,2,+}(\hat{\Phi}_{
\hat{b}_{j-1}}(z)) \me^{-\mi ((n-1)K+k) \hat{\mho}_{j-1}} \right), 
\label{eqmaininf89} \\
\hat{\mathcal{X}}_{22}^{\hat{b},2}(z) =& \, \sqrt{\pi} \me^{-\frac{\mi \pi}{6}} 
\me^{-\frac{1}{2} ((n-1)K+k) \hat{\xi}_{\hat{b}_{j-1}}(z)} 
\tilde{\mathfrak{m}}^{\raise-0.5ex\hbox{$\scriptstyle \infty$}}_{22} \left(
\hat{\mathbb{M}}_{22}(z) \Theta_{2,2}^{0,2,+}(\hat{\Phi}_{\hat{b}_{j-1}}(z)) \! - 
\! \mi \hat{\mathbb{M}}_{21}(z) \Theta_{2,2}^{0,2,-}(\hat{\Phi}_{\hat{b}_{j-1}}
(z)) \me^{\mi ((n-1)K+k) \hat{\mho}_{j-1}} \right), \label{eqmaininf90}
\end{align}
with
\begin{equation} \label{eqmaininf91} 
\varTheta_{r_{3},r_{4}}^{r_{1},r_{2},\pm}(\zeta) \! := \! \zeta^{1/4} \left(
\operatorname{Ai}(\zeta) \! - \! \me^{\frac{\mi \pi}{3}} \omega^{r_{1}} 
\operatorname{Ai}(\omega^{r_{2}} \zeta) \right) \! \pm \! \zeta^{-1/4} 
\left(\operatorname{Ai}^{\prime}(\zeta) \! - \! \me^{\frac{\mi \pi}{3}}
\omega^{r_{3}} \operatorname{Ai}^{\prime}(\omega^{r_{4}} \zeta) \right), 
\quad r_{1},r_{2},r_{3},r_{4} \! = \! 0,1,2;
\end{equation}
{\rm \pmb{(8)}} for $z \! \in \! \hat{\Omega}^{2}_{\hat{a}_{j}}$, $j \! = \! 
1,2,\dotsc,N \! + \!1$,
\begin{align}
\pmb{\pi}_{k}^{n}(z) \underset{\underset{z_{o}=1+o(1)}{\mathscr{N},
n \to \infty}}{=}& \, \left(\left(\hat{\mathcal{X}}_{11}^{\hat{a},2}(z) \! + \! 
\hat{\mathcal{X}}_{12}^{\hat{a},2}(z) \me^{-2 \pi \mi ((n-1)K+k) \int_{z}^{
\hat{a}_{N+1}} \psi_{\widetilde{V}}^{\infty}(\xi) \, \md \xi} \right) \left(1 \! 
+ \! \dfrac{1}{(n \! - \! 1)K \! + \! k} \left(\hat{\mathcal{R}}_{11}^{\sharp}(z) 
\! - \! \hat{\mathcal{R}}_{11}^{\natural}(z) \right) \! + \! \mathcal{O} \left(
\dfrac{\mathfrak{c}(n,k,z_{o})}{((n \! - \! 1)K \! + \! k)^{2}} \right) \right) 
\right. \nonumber \\
+&\left. \, \left(\hat{\mathcal{X}}_{21}^{\hat{a},2}(z) \! + \! \hat{\mathcal{
X}}_{22}^{\hat{a},2}(z) \me^{-2 \pi \mi ((n-1)K+k) \int_{z}^{\hat{a}_{N+1}} 
\psi_{\widetilde{V}}^{\infty}(\xi) \, \md \xi} \right) \left(\dfrac{1}{(n \! - \! 
1)K \! + \! k} \left(\hat{\mathcal{R}}_{12}^{\sharp}(z) \! - \! \hat{\mathcal{
R}}_{12}^{\natural}(z) \right) \! + \! \mathcal{O} \left(\dfrac{\mathfrak{c}
(n,k,z_{o})}{((n \! - \! 1)K \! + \! k)^{2}} \right) \right) \right) \nonumber \\
\times& \, \me^{n(g^{\infty}(z)-\tilde{\mathscr{P}}_{0})}, \label{eqmaininf92}
\end{align}
and
\begin{align}
\int_{\mathbb{R}} \dfrac{\pmb{\pi}_{k}^{n}(\xi) 
\me^{-n \widetilde{V}(\xi)}}{\xi \! - \! z} \, \dfrac{\md \xi}{2 \pi \mi} 
\underset{\underset{z_{o}=1+o(1)}{\mathscr{N},n \to \infty}}{=}& \,
\left(\hat{\mathcal{X}}_{12}^{\hat{a},2}(z) \left(1 \! + \! \dfrac{1}{(n \! - \! 
1)K \! + \! k} \left(\hat{\mathcal{R}}_{11}^{\sharp}(z) \! - \! \hat{\mathcal{
R}}_{11}^{\natural}(z) \right) \! + \! \mathcal{O} \left(\dfrac{\mathfrak{c}
(n,k,z_{o})}{((n \! - \! 1)K \! + \! k)^{2}} \right) \right) \right. \nonumber \\
+&\left. \, \hat{\mathcal{X}}_{22}^{\hat{a},2}(z) \left(\dfrac{1}{(n \! - \! 1)K \! 
+ \! k} \left(\hat{\mathcal{R}}_{12}^{\sharp}(z) \! - \! \hat{\mathcal{R}}_{12}^{
\natural}(z) \right) \! + \! \mathcal{O} \left(\dfrac{\mathfrak{c}(n,k,z_{o})}{(
(n \! - \! 1)K \! + \! k)^{2}} \right) \right) \right) \nonumber \\
\times& \, \me^{n \hat{\ell}} \me^{-n(g^{\infty}(z)-\tilde{\mathscr{P}}_{0})}, 
\label{eqmaininf93}
\end{align}
where
\begin{align}
\hat{\mathcal{X}}_{11}^{\hat{a},2}(z) =& \, -\mi \sqrt{\pi} \me^{\frac{1}{2} 
((n-1)K+k) \hat{\xi}_{\hat{a}_{j}}(z)} 
\tilde{\mathfrak{m}}^{\raise-0.5ex\hbox{$\scriptstyle \infty$}}_{11} \left(\mi 
\hat{\mathbb{M}}_{11}(z) \varTheta_{2,0}^{0,2,-}(\hat{\Phi}_{\hat{a}_{j}}(z)) \! + 
\! \hat{\mathbb{M}}_{12}(z) \varTheta_{2,0}^{0,2,+}(\hat{\Phi}_{\hat{a}_{j}}(z)) 
\me^{-\mi ((n-1)K+k) \hat{\mho}_{j}} \right), \label{eqmaininf94} \\
\hat{\mathcal{X}}_{12}^{\hat{a},2}(z) =& \, \sqrt{\pi} \me^{-\frac{\mi \pi}{6}} 
\me^{-\frac{1}{2} ((n-1)K+k) \hat{\xi}_{\hat{a}_{j}}(z)} 
\tilde{\mathfrak{m}}^{\raise-0.5ex\hbox{$\scriptstyle \infty$}}_{11} \left(
\hat{\mathbb{M}}_{12}(z) \Theta_{2,2}^{0,2,+}(\hat{\Phi}_{\hat{a}_{j}}(z)) \! + 
\! \mi \hat{\mathbb{M}}_{11}(z) \Theta_{2,2}^{0,2,-}(\hat{\Phi}_{\hat{a}_{j}}(z)) 
\me^{\mi ((n-1)K+k) \hat{\mho}_{j}} \right), \label{eqmaininf95} \\
\hat{\mathcal{X}}_{21}^{\hat{a},2}(z) =& \, -\mi \sqrt{\pi} \me^{\frac{1}{2} 
((n-1)K+k) \hat{\xi}_{\hat{a}_{j}}(z)} 
\tilde{\mathfrak{m}}^{\raise-0.5ex\hbox{$\scriptstyle \infty$}}_{22} \left(
\mi \hat{\mathbb{M}}_{21}(z) \varTheta_{2,0}^{0,2,-}(\hat{\Phi}_{\hat{a}_{j}}(z)) 
\! + \! \hat{\mathbb{M}}_{22}(z) \varTheta_{2,0}^{0,2,+}(\hat{\Phi}_{\hat{a}_{j}}
(z)) \me^{-\mi ((n-1)K+k) \hat{\mho}_{j}} \right), \label{eqmaininf96} \\
\hat{\mathcal{X}}_{22}^{\hat{a},2}(z) =& \, \sqrt{\pi} \me^{-\frac{\mi \pi}{6}} 
\me^{-\frac{1}{2} ((n-1)K+k) \hat{\xi}_{\hat{a}_{j}}(z)} 
\tilde{\mathfrak{m}}^{\raise-0.5ex\hbox{$\scriptstyle \infty$}}_{22} \left(
\hat{\mathbb{M}}_{22}(z) \Theta_{2,2}^{0,2,+}(\hat{\Phi}_{\hat{a}_{j}}(z)) \! + 
\! \mi \hat{\mathbb{M}}_{21}(z) \Theta_{2,2}^{0,2,-}(\hat{\Phi}_{\hat{a}_{j}}(z)) 
\me^{\mi ((n-1)K+k) \hat{\mho}_{j}} \right); \label{eqmaininf97}
\end{align}
{\rm \pmb{(9)}} for $z \! \in \! \hat{\Omega}^{3}_{\hat{b}_{j-1}}$, $j \! = \! 
1,2,\dotsc,N \! + \! 1$,
\begin{align}
\pmb{\pi}_{k}^{n}(z) \underset{\underset{z_{o}=1+o(1)}{\mathscr{N},
n \to \infty}}{=}& \, \left(\left(\hat{\mathcal{X}}_{11}^{\hat{b},3}(z) \! - \! 
\hat{\mathcal{X}}_{12}^{\hat{b},3}(z) \me^{2 \pi \mi ((n-1)K+k) \int_{z}^{
\hat{a}_{N+1}} \psi_{\widetilde{V}}^{\infty}(\xi) \, \md \xi} \right) \left(1 \! + 
\! \dfrac{1}{(n \! - \! 1)K \! + \! k} \left(\hat{\mathcal{R}}_{11}^{\sharp}(z) 
\! - \! \hat{\mathcal{R}}_{11}^{\natural}(z) \right) \! + \! \mathcal{O} \left(
\dfrac{\mathfrak{c}(n,k,z_{o})}{((n \! - \! 1)K \! + \! k)^{2}} \right) \right) 
\right. \nonumber \\
+&\left. \, \left(\hat{\mathcal{X}}_{21}^{\hat{b},3}(z) \! - \! \hat{\mathcal{
X}}_{22}^{\hat{b},3}(z) \me^{2 \pi \mi ((n-1)K+k) \int_{z}^{\hat{a}_{N+1}} 
\psi_{\widetilde{V}}^{\infty}(\xi) \, \md \xi} \right) \left(\dfrac{1}{(n \! - \! 
1)K \! + \! k} \left(\hat{\mathcal{R}}_{12}^{\sharp}(z) \! - \! \hat{\mathcal{
R}}_{12}^{\natural}(z) \right) \! + \! \mathcal{O} \left(\dfrac{\mathfrak{c}
(n,k,z_{o})}{((n \! - \! 1)K \! + \! k)^{2}} \right) \right) \right) \nonumber \\
\times& \, \me^{n(g^{\infty}(z)-\tilde{\mathscr{P}}_{0})}, \label{eqmaininf98}
\end{align}
and
\begin{align}
\int_{\mathbb{R}} \dfrac{\pmb{\pi}_{k}^{n}(\xi) 
\me^{-n \widetilde{V}(\xi)}}{\xi \! - \! z} \, \dfrac{\md \xi}{2 \pi \mi} 
\underset{\underset{z_{o}=1+o(1)}{\mathscr{N},n \to \infty}}{=}& \,
\left(\hat{\mathcal{X}}_{12}^{\hat{b},3}(z) \left(1 \! + \! \dfrac{1}{(n \! - \! 
1)K \! + \! k} \left(\hat{\mathcal{R}}_{11}^{\sharp}(z) \! - \! \hat{\mathcal{
R}}_{11}^{\natural}(z) \right) \! + \! \mathcal{O} \left(\dfrac{\mathfrak{c}
(n,k,z_{o})}{((n \! - \! 1)K \! + \! k)^{2}} \right) \right) \right. \nonumber \\
+&\left. \, \hat{\mathcal{X}}_{22}^{\hat{b},3}(z) \left(\dfrac{1}{(n \! - \! 1)K \! 
+ \! k} \left(\hat{\mathcal{R}}_{12}^{\sharp}(z) \! - \! \hat{\mathcal{R}}_{12}^{
\natural}(z) \right) \! + \! \mathcal{O} \left(\dfrac{\mathfrak{c}(n,k,z_{o})}{(
(n \! - \! 1)K \! + \! k)^{2}} \right) \right) \right) \nonumber \\
\times& \, \me^{n \hat{\ell}} \me^{-n(g^{\infty}(z)-\tilde{\mathscr{P}}_{0})}, 
\label{eqmaininf99}
\end{align}
where
\begin{align}
\hat{\mathcal{X}}_{11}^{\hat{b},3}(z) =& \, -\mi \sqrt{\pi} \me^{\frac{1}{2} 
((n-1)K+k) \hat{\xi}_{\hat{b}_{j-1}}(z)} 
\tilde{\mathfrak{m}}^{\raise-0.5ex\hbox{$\scriptstyle \infty$}}_{11} \left(\mi 
\hat{\mathbb{M}}_{12}(z) \varTheta_{0,1}^{2,1,-}(\hat{\Phi}_{\hat{b}_{j-1}}(z)) \! 
+ \! \hat{\mathbb{M}}_{11}(z) \varTheta_{0,1}^{2,1,+}(\hat{\Phi}_{\hat{b}_{j-1}}
(z)) \me^{\mi ((n-1)K+k) \hat{\mho}_{j-1}} \right), \label{eqmaininf100} \\
\hat{\mathcal{X}}_{12}^{\hat{b},3}(z) =& \, \sqrt{\pi} \me^{-\frac{\mi \pi}{6}} 
\me^{-\frac{1}{2} ((n-1)K+k) \hat{\xi}_{\hat{b}_{j-1}}(z)} 
\tilde{\mathfrak{m}}^{\raise-0.5ex\hbox{$\scriptstyle \infty$}}_{11} \left(
\hat{\mathbb{M}}_{11}(z) \Theta_{0,1}^{2,1,+}(\hat{\Phi}_{\hat{b}_{j-1}}(z)) \! 
+ \! \mi \hat{\mathbb{M}}_{12}(z) \Theta_{0,1}^{2,1,-}(\hat{\Phi}_{\hat{b}_{j-
1}}(z)) \me^{-\mi ((n-1)K+k) \hat{\mho}_{j-1}} \right), \label{eqmaininf101} \\
\hat{\mathcal{X}}_{21}^{\hat{b},3}(z) =& \, -\mi \sqrt{\pi} \me^{\frac{1}{2} 
((n-1)K+k) \hat{\xi}_{\hat{b}_{j-1}}(z)} 
\tilde{\mathfrak{m}}^{\raise-0.5ex\hbox{$\scriptstyle \infty$}}_{22} \left(
\mi \hat{\mathbb{M}}_{22}(z) \varTheta_{0,1}^{2,1,-}(\hat{\Phi}_{\hat{b}_{j-1}}
(z)) \! + \! \hat{\mathbb{M}}_{21}(z) \varTheta_{0,1}^{2,1,+}(\hat{\Phi}_{
\hat{b}_{j-1}}(z)) \me^{\mi ((n-1)K+k) \hat{\mho}_{j-1}} \right), 
\label{eqmaininf102} \\
\hat{\mathcal{X}}_{22}^{\hat{b},3}(z) =& \, \sqrt{\pi} \me^{-\frac{\mi \pi}{6}} 
\me^{-\frac{1}{2}((n-1)K+k) \hat{\xi}_{\hat{b}_{j-1}}(z)} 
\tilde{\mathfrak{m}}^{\raise-0.5ex\hbox{$\scriptstyle \infty$}}_{22} \left(
\hat{\mathbb{M}}_{21}(z) \Theta_{0,1}^{2,1,+}(\hat{\Phi}_{\hat{b}_{j-1}}(z)) \! + 
\! \mi \hat{\mathbb{M}}_{22}(z) \Theta_{0,1}^{2,1,-}(\hat{\Phi}_{\hat{b}_{j-1}}
(z)) \me^{-\mi ((n-1)K+k) \hat{\mho}_{j-1}} \right); \label{eqmaininf103}
\end{align}
{\rm \pmb{(10)}} for $z \! \in \! \hat{\Omega}^{3}_{\hat{a}_{j}}$, $j \! = \! 
1,2,\dotsc,N \! + \!1$,
\begin{align}
\pmb{\pi}_{k}^{n}(z) \underset{\underset{z_{o}=1+o(1)}{\mathscr{N},
n \to \infty}}{=}& \, \left(\left(\hat{\mathcal{X}}_{11}^{\hat{a},3}(z) \! - \! 
\hat{\mathcal{X}}_{12}^{\hat{a},3}(z) \me^{2 \pi \mi ((n-1)K+k) \int_{z}^{
\hat{a}_{N+1}} \psi_{\widetilde{V}}^{\infty}(\xi) \, \md \xi} \right) \left(1 \! + 
\! \dfrac{1}{(n \! - \! 1)K \! + \! k} \left(\hat{\mathcal{R}}_{11}^{\sharp}(z) 
\! - \! \hat{\mathcal{R}}_{11}^{\natural}(z) \right) \! + \! \mathcal{O} \left(
\dfrac{\mathfrak{c}(n,k,z_{o})}{((n \! - \! 1)K \! + \! k)^{2}} \right) \right) 
\right. \nonumber \\
+&\left. \, \left(\hat{\mathcal{X}}_{21}^{\hat{a},3}(z) \! - \! \hat{\mathcal{
X}}_{22}^{\hat{a},3}(z) \me^{2 \pi \mi ((n-1)K+k) \int_{z}^{\hat{a}_{N+1}} 
\psi_{\widetilde{V}}^{\infty}(\xi) \, \md \xi} \right) \left(\dfrac{1}{(n \! - \! 
1)K \! + \! k} \left(\hat{\mathcal{R}}_{12}^{\sharp}(z) \! - \! \hat{\mathcal{
R}}_{12}^{\natural}(z) \right) \! + \! \mathcal{O} \left(\dfrac{\mathfrak{c}
(n,k,z_{o})}{((n \! - \! 1)K \! + \! k)^{2}} \right) \right) \right) \nonumber \\
\times& \, \me^{n(g^{\infty}(z)-\tilde{\mathscr{P}}_{0})}, \label{eqmaininf104}
\end{align}
and
\begin{align}
\int_{\mathbb{R}} \dfrac{\pmb{\pi}_{k}^{n}(\xi) 
\me^{-n \widetilde{V}(\xi)}}{\xi \! - \! z} \, \dfrac{\md \xi}{2 \pi \mi} 
\underset{\underset{z_{o}=1+o(1)}{\mathscr{N},n \to \infty}}{=}& \,
\left(\hat{\mathcal{X}}_{12}^{\hat{a},3}(z) \left(1 \! + \! \dfrac{1}{(n \! - \! 
1)K \! + \! k} \left(\hat{\mathcal{R}}_{11}^{\sharp}(z) \! - \! \hat{\mathcal{
R}}_{11}^{\natural}(z) \right) \! + \! \mathcal{O} \left(\dfrac{\mathfrak{c}
(n,k,z_{o})}{((n \! - \! 1)K \! + \! k)^{2}} \right) \right) \right. \nonumber \\
+&\left. \, \hat{\mathcal{X}}_{22}^{\hat{a},3}(z) \left(\dfrac{1}{(n \! - \! 1)K \! 
+ \! k} \left(\hat{\mathcal{R}}_{12}^{\sharp}(z) \! - \! \hat{\mathcal{R}}_{12}^{
\natural}(z) \right) \! + \! \mathcal{O} \left(\dfrac{\mathfrak{c}(n,k,z_{o})}{(
(n \! - \! 1)K \! + \! k)^{2}} \right) \right) \right) \nonumber \\
\times& \, \me^{n \hat{\ell}} \me^{-n(g^{\infty}(z)-\tilde{\mathscr{P}}_{0})}, 
\label{eqmaininf105}
\end{align}
where
\begin{align}
\hat{\mathcal{X}}_{11}^{\hat{a},3}(z) =& \, -\mi \sqrt{\pi} \me^{\frac{1}{2} 
((n-1)K+k) \hat{\xi}_{\hat{a}_{j}}(z)} 
\tilde{\mathfrak{m}}^{\raise-0.5ex\hbox{$\scriptstyle \infty$}}_{11} \left(\mi 
\hat{\mathbb{M}}_{12}(z) \varTheta_{0,1}^{2,1,-}(\hat{\Phi}_{\hat{a}_{j}}(z)) \! - 
\! \hat{\mathbb{M}}_{11}(z) \varTheta_{0,1}^{2,1,+}(\hat{\Phi}_{\hat{a}_{j}}(z)) 
\me^{\mi ((n-1)K+k) \hat{\mho}_{j}} \right), \label{eqmaininf106} \\
\hat{\mathcal{X}}_{12}^{\hat{a},3}(z) =& \, \sqrt{\pi} \me^{-\frac{\mi \pi}{6}} 
\me^{-\frac{1}{2} ((n-1)K+k) \hat{\xi}_{\hat{a}_{j}}(z)} 
\tilde{\mathfrak{m}}^{\raise-0.5ex\hbox{$\scriptstyle \infty$}}_{11} \left(
\hat{\mathbb{M}}_{11}(z) \Theta_{0,1}^{2,1,+}(\hat{\Phi}_{\hat{a}_{j}}(z)) \! - 
\! \mi \hat{\mathbb{M}}_{12}(z) \Theta_{0,1}^{2,1,-}(\hat{\Phi}_{\hat{a}_{j}}
(z)) \me^{-\mi ((n-1)K+k) \hat{\mho}_{j}} \right), \label{eqmaininf107} \\
\hat{\mathcal{X}}_{21}^{\hat{a},3}(z) =& \, -\mi \sqrt{\pi} \me^{\frac{1}{2} 
((n-1)K+k) \hat{\xi}_{\hat{a}_{j}}(z)} 
\tilde{\mathfrak{m}}^{\raise-0.5ex\hbox{$\scriptstyle \infty$}}_{22} \left(
\mi \hat{\mathbb{M}}_{22}(z) \varTheta_{0,1}^{2,1,-}(\hat{\Phi}_{\hat{a}_{j}}(z)) 
\! - \! \hat{\mathbb{M}}_{21}(z) \varTheta_{0,1}^{2,1,+}(\hat{\Phi}_{\hat{a}_{j}}
(z)) \me^{\mi ((n-1)K+k) \hat{\mho}_{j}} \right), \label{eqmaininf108} \\
\hat{\mathcal{X}}_{22}^{\hat{a},3}(z) =& \, \sqrt{\pi} \me^{-\frac{\mi \pi}{6}} 
\me^{-\frac{1}{2}((n-1)K+k) \hat{\xi}_{\hat{a}_{j}}(z)} 
\tilde{\mathfrak{m}}^{\raise-0.5ex\hbox{$\scriptstyle \infty$}}_{22} \left(
\hat{\mathbb{M}}_{21}(z) \Theta_{0,1}^{2,1,+}(\hat{\Phi}_{\hat{a}_{j}}(z)) \! - 
\! \mi \hat{\mathbb{M}}_{22}(z) \Theta_{0,1}^{2,1,-}(\hat{\Phi}_{\hat{a}_{j}}
(z)) \me^{-\mi ((n-1)K+k) \hat{\mho}_{j}} \right); \label{eqmaininf109}
\end{align}
{\rm \pmb{(11)}} for $z \! \in \! \hat{\Omega}^{4}_{\hat{b}_{j-1}}$, 
$j \! = \! 1,2,\dotsc,N \! + \!1$,
\begin{align}
\pmb{\pi}_{k}^{n}(z) \underset{\underset{z_{o}=1+o(1)}{\mathscr{N},
n \to \infty}}{=}& \, \left(\hat{\mathcal{X}}_{11}^{\hat{b},4}(z) \left(1 \! + \! 
\dfrac{1}{(n \! - \! 1)K \! + \! k} \left(\hat{\mathcal{R}}_{11}^{\sharp}(z) \! - \! 
\hat{\mathcal{R}}_{11}^{\natural}(z) \right) \! + \! \mathcal{O} \left(\dfrac{
\mathfrak{c}(n,k,z_{o})}{((n \! - \! 1)K \! + \! k)^{2}} \right) \right) \right. 
\nonumber \\
+&\left. \, \hat{\mathcal{X}}_{21}^{\hat{b},4}(z) \left(\dfrac{1}{(n \! - \! 1)K \! 
+ \! k} \left(\hat{\mathcal{R}}_{12}^{\sharp}(z) \! - \! \hat{\mathcal{R}}_{12}^{
\natural}(z) \right) \! + \! \mathcal{O} \left(\dfrac{\mathfrak{c}(n,k,z_{o})}{((n 
\! - \! 1)K \! + \! k)^{2}} \right) \right) \right) \nonumber \\
\times& \, \me^{n(g^{\infty}(z)-\tilde{\mathscr{P}}_{0})}, 
\label{eqmaininf110}
\end{align}
and
\begin{align}
\int_{\mathbb{R}} \dfrac{\pmb{\pi}_{k}^{n}(\xi) 
\me^{-n \widetilde{V}(\xi)}}{\xi \! - \! z} \, \dfrac{\md \xi}{2 \pi \mi} 
\underset{\underset{z_{o}=1+o(1)}{\mathscr{N},n \to \infty}}{=}& \, \left(
\hat{\mathcal{X}}_{12}^{\hat{b},4}(z) \left(1 \! + \! \dfrac{1}{(n \! - \! 1)K \! + 
\! k} \left(\hat{\mathcal{R}}_{11}^{\sharp}(z) \! - \! \hat{\mathcal{R}}_{11}^{
\natural}(z) \right) \! + \! \mathcal{O} \left(\dfrac{\mathfrak{c}(n,k,z_{o})}{((n 
\! - \! 1)K \! + \! k)^{2}} \right) \right) \right. \nonumber \\
+&\left. \, \hat{\mathcal{X}}_{22}^{\hat{b},4}(z) \left(\dfrac{1}{(n \! - \! 1)K \! 
+ \! k} \left(\hat{\mathcal{R}}_{12}^{\sharp}(z) \! - \! \hat{\mathcal{R}}_{12}^{
\natural}(z) \right) \! + \! \mathcal{O} \left(\dfrac{\mathfrak{c}(n,k,z_{o})}{((n 
\! - \! 1)K \! + \! k)^{2}} \right) \right) \right) \nonumber \\
\times& \, \me^{n \hat{\ell}} \me^{-n(g^{\infty}(z)-\tilde{\mathscr{P}}_{0})}, 
\label{eqmaininf111}
\end{align}
where
\begin{align}
\hat{\mathcal{X}}_{11}^{\hat{b},4}(z) =& \, -\mi \sqrt{\pi} \me^{\frac{1}{2} 
((n-1)K+k) \hat{\xi}_{\hat{b}_{j-1}}(z)} 
\tilde{\mathfrak{m}}^{\raise-0.5ex\hbox{$\scriptstyle \infty$}}_{11} \left(
\mi \hat{\mathbb{M}}_{12}(z) \Theta_{0,0}^{0,0,-}(\hat{\Phi}_{\hat{b}_{j-1}}(z)) 
\! + \! \hat{\mathbb{M}}_{11}(z) \Theta_{0,0}^{0,0,+}(\hat{\Phi}_{\hat{b}_{j-1}}
(z)) \me^{\mi ((n-1)K+k) \hat{\mho}_{j-1}} \right), \label{eqmaininf112} \\
\hat{\mathcal{X}}_{12}^{\hat{b},4}(z) =& \, \sqrt{\pi} \me^{-\frac{\mi \pi}{6}} 
\me^{-\frac{1}{2} ((n-1)K+k) \hat{\xi}_{\hat{b}_{j-1}}(z)} 
\tilde{\mathfrak{m}}^{\raise-0.5ex\hbox{$\scriptstyle \infty$}}_{11} \left(
\hat{\mathbb{M}}_{11}(z) \Theta_{0,1}^{2,1,+}(\hat{\Phi}_{\hat{b}_{j-1}}(z)) \! 
+ \! \mi \hat{\mathbb{M}}_{12}(z) \Theta_{0,1}^{2,1,-}(\hat{\Phi}_{\hat{b}_{j-1}}
(z)) \me^{-\mi ((n-1)K+k) \hat{\mho}_{j-1}} \right), \label{eqmaininf113} \\
\hat{\mathcal{X}}_{21}^{\hat{b},4}(z) =& \, -\mi \sqrt{\pi} \me^{\frac{1}{2} 
((n-1)K+k) \hat{\xi}_{\hat{b}_{j-1}}(z)} 
\tilde{\mathfrak{m}}^{\raise-0.5ex\hbox{$\scriptstyle \infty$}}_{22} \left(
\mi \hat{\mathbb{M}}_{22}(z) \Theta_{0,0}^{0,0,-}(\hat{\Phi}_{\hat{b}_{j-1}}(z)) 
\! + \! \hat{\mathbb{M}}_{21}(z) \Theta_{0,0}^{0,0,+}(\hat{\Phi}_{\hat{b}_{j-1}}
(z)) \me^{\mi ((n-1)K+k) \hat{\mho}_{j-1}} \right), \label{eqmaininf114} \\
\hat{\mathcal{X}}_{22}^{\hat{b},4}(z) =& \, \sqrt{\pi} \me^{-\frac{\mi \pi}{6}} 
\me^{-\frac{1}{2} ((n-1)K+k) \hat{\xi}_{\hat{b}_{j-1}}(z)} 
\tilde{\mathfrak{m}}^{\raise-0.5ex\hbox{$\scriptstyle \infty$}}_{22} \left(
\hat{\mathbb{M}}_{21}(z) \Theta_{0,1}^{2,1,+}(\hat{\Phi}_{\hat{b}_{j-1}}(z)) \! + 
\! \mi \hat{\mathbb{M}}_{22}(z) \Theta_{0,1}^{2,1,-}(\hat{\Phi}_{\hat{b}_{j-1}}
(z)) \me^{-\mi ((n-1)K+k) \hat{\mho}_{j-1}} \right); \label{eqmaininf115}
\end{align}
and {\rm \pmb{(12)}} for $z \! \in \! \hat{\Omega}^{4}_{\hat{a}_{j}}$, 
$j \! = \! 1,2,\dotsc,N \! + \! 1$,
\begin{align}
\pmb{\pi}_{k}^{n}(z) \underset{\underset{z_{o}=1+o(1)}{\mathscr{N},
n \to \infty}}{=}& \, \left(\hat{\mathcal{X}}_{11}^{\hat{a},4}(z) \left(1 \! + \! 
\dfrac{1}{(n \! - \! 1)K \! + \! k} \left(\hat{\mathcal{R}}_{11}^{\sharp}(z) \! - \! 
\hat{\mathcal{R}}_{11}^{\natural}(z) \right) \! + \! \mathcal{O} \left(\dfrac{
\mathfrak{c}(n,k,z_{o})}{((n \! - \! 1)K \! + \! k)^{2}} \right) \right) \right. 
\nonumber \\
+&\left. \, \hat{\mathcal{X}}_{21}^{\hat{a},4}(z) \left(\dfrac{1}{(n \! - \! 1)K \! 
+ \! k} \left(\hat{\mathcal{R}}_{12}^{\sharp}(z) \! - \! \hat{\mathcal{R}}_{12}^{
\natural}(z) \right) \! + \! \mathcal{O} \left(\dfrac{\mathfrak{c}(n,k,z_{o})}{((n 
\! - \! 1)K \! + \! k)^{2}} \right) \right) \right) \nonumber \\
\times& \, \me^{n(g^{\infty}(z)-\tilde{\mathscr{P}}_{0})}, 
\label{eqmaininf116}
\end{align}
and
\begin{align}
\int_{\mathbb{R}} \dfrac{\pmb{\pi}_{k}^{n}(\xi) 
\me^{-n \widetilde{V}(\xi)}}{\xi \! - \! z} \, \dfrac{\md \xi}{2 \pi \mi} 
\underset{\underset{z_{o}=1+o(1)}{\mathscr{N},n \to \infty}}{=}& \, \left(
\hat{\mathcal{X}}_{12}^{\hat{a},4}(z) \left(1 \! + \! \dfrac{1}{(n \! - \! 1)K \! + 
\! k} \left(\hat{\mathcal{R}}_{11}^{\sharp}(z) \! - \! \hat{\mathcal{R}}_{11}^{
\natural}(z) \right) \! + \! \mathcal{O} \left(\dfrac{\mathfrak{c}(n,k,z_{o})}{((n 
\! - \! 1)K \! + \! k)^{2}} \right) \right) \right. \nonumber \\
+&\left. \, \hat{\mathcal{X}}_{22}^{\hat{a},4}(z) \left(\dfrac{1}{(n \! - \! 1)K \! 
+ \! k} \left(\hat{\mathcal{R}}_{12}^{\sharp}(z) \! - \! \hat{\mathcal{R}}_{12}^{
\natural}(z) \right) \! + \! \mathcal{O} \left(\dfrac{\mathfrak{c}(n,k,z_{o})}{((n 
\! - \! 1)K \! + \! k)^{2}} \right) \right) \right) \nonumber \\
\times& \, \me^{n \hat{\ell}} \me^{-n(g^{\infty}(z)-\tilde{\mathscr{P}}_{0})}, 
\label{eqmaininf117}
\end{align}
where
\begin{align}
\hat{\mathcal{X}}_{11}^{\hat{a},4}(z) =& \, -\mi \sqrt{\pi} \me^{\frac{1}{2} 
((n-1)K+k) \hat{\xi}_{\hat{a}_{j}}(z)} 
\tilde{\mathfrak{m}}^{\raise-0.5ex\hbox{$\scriptstyle \infty$}}_{11} \left(
\mi \hat{\mathbb{M}}_{12}(z) \Theta_{0,0}^{0,0,-}(\hat{\Phi}_{\hat{a}_{j}}(z)) 
\! - \! \hat{\mathbb{M}}_{11}(z) \Theta_{0,0}^{0,0,+}(\hat{\Phi}_{\hat{a}_{j}}
(z)) \me^{\mi ((n-1)K+k) \hat{\mho}_{j}} \right), \label{eqmaininf118} \\
\hat{\mathcal{X}}_{12}^{\hat{a},4}(z) =& \, \sqrt{\pi} \me^{-\frac{\mi \pi}{6}} 
\me^{-\frac{1}{2} ((n-1)K+k) \hat{\xi}_{\hat{a}_{j}}(z)} 
\tilde{\mathfrak{m}}^{\raise-0.5ex\hbox{$\scriptstyle \infty$}}_{11} \left(
\hat{\mathbb{M}}_{11}(z) \Theta_{0,1}^{2,1,+}(\hat{\Phi}_{\hat{a}_{j}}(z)) \! - 
\! \mi \hat{\mathbb{M}}_{12}(z) \Theta_{0,1}^{2,1,-}(\hat{\Phi}_{\hat{a}_{j}}
(z)) \me^{-\mi ((n-1)K+k) \hat{\mho}_{j}} \right), \label{eqmaininf119} \\
\hat{\mathcal{X}}_{21}^{\hat{a},4}(z) =& \, -\mi \sqrt{\pi} \me^{\frac{1}{2} 
((n-1)K+k) \hat{\xi}_{\hat{a}_{j}}(z)} 
\tilde{\mathfrak{m}}^{\raise-0.5ex\hbox{$\scriptstyle \infty$}}_{22} \left(
\mi \hat{\mathbb{M}}_{22}(z) \Theta_{0,0}^{0,0,-}(\hat{\Phi}_{\hat{a}_{j}}(z)) \! 
- \! \hat{\mathbb{M}}_{21}(z) \Theta_{0,0}^{0,0,+}(\hat{\Phi}_{\hat{a}_{j}}(z)) 
\me^{\mi ((n-1)K+k) \hat{\mho}_{j}} \right), \label{eqmaininf120} \\
\hat{\mathcal{X}}_{22}^{\hat{a},4}(z) =& \, \sqrt{\pi} \me^{-\frac{\mi \pi}{6}} 
\me^{-\frac{1}{2} ((n-1)K+k) \hat{\xi}_{\hat{a}_{j}}(z)} 
\tilde{\mathfrak{m}}^{\raise-0.5ex\hbox{$\scriptstyle \infty$}}_{22} \left(
\hat{\mathbb{M}}_{21}(z) \Theta_{0,1}^{2,1,+}(\hat{\Phi}_{\hat{a}_{j}}(z)) \! - 
\! \mi \hat{\mathbb{M}}_{22}(z) \Theta_{0,1}^{2,1,-}(\hat{\Phi}_{\hat{a}_{j}}(z)) 
\me^{-\mi ((n-1)K+k) \hat{\mho}_{j}} \right). \label{eqmaininf121}
\end{align}
\end{dddd}
\begin{eeee} \label{remvalinf} 
\textsl{For $n \! \in \! \mathbb{N}$ and $k \! \in \! \lbrace 1,2,\dotsc,
K \rbrace$ such that $\alpha_{p_{\mathfrak{s}}} \! := \! \alpha_{k} \! = 
\! \infty$, using limiting values, if necessary, the asymptotics, in the 
double-scaling $\mathscr{N},n \! \to \! \infty$ such that $z_{o} \! = \! 1 
\! + \! o(1)$, in Theorem~\ref{maintheoforinf1} for $\pmb{\pi}_{k}^{n}
(z)$ and $\int_{\mathbb{R}} \pmb{\pi}_{k}^{n}(\xi)(\xi \! - \! z)^{-1} 
\me^{-n \widetilde{V}(\xi)} \, \md \xi/2 \pi \mi$ have a natural 
interpretation on $\mathbb{R} \setminus \lbrace \alpha_{1},\alpha_{2},
\dotsc,\alpha_{K} \rbrace$.}
\end{eeee}
\begin{eeee} \label{remfortheobelow1} 
\textsl{The bulk of the parameters appearing in 
Theorems~\ref{maintheompainf} and~\ref{maintheoforinf2} below 
have been defined heretofore in Theorem~\ref{maintheoforinf1}.}
\end{eeee}
\begin{dddd} \label{maintheompainf} 
Let the external field $\widetilde{V} \colon \overline{\mathbb{R}} \setminus 
\lbrace \alpha_{1},\alpha_{2},\dotsc,\alpha_{K} \rbrace \! \to \! \mathbb{R}$ 
satisfy conditions~\eqref{eq20}--\eqref{eq22}, and suppose that 
$\widetilde{V}$ is regular. For $n \! \in \! \mathbb{N}$ and $k \! \in \! 
\lbrace 1,2,\dotsc,K \rbrace$ such that $\alpha_{p_{\mathfrak{s}}} \! := \! 
\alpha_{k} \! = \! \infty$, let $\mathcal{X} \colon \overline{\mathbb{C}} 
\setminus \overline{\mathbb{R}} \! \to \! \mathrm{SL}_{2}(\mathbb{C})$ 
be the unique solution of the corresponding monic {\rm MPC ORF RHP} 
$(\mathcal{X}(z),\upsilon (z),\overline{\mathbb{R}})$ stated in 
Lemma~$\bm{\mathrm{RHP}_{\mathrm{MPC}}}$ with integral 
representation given by Equation~\eqref{intrepinf}, where, in particular,
\begin{equation*}
(\mathcal{X}(z))_{11} \! = \! \pmb{\pi}_{k}^{n}(z) \quad \quad \, 
\text{and} \, \quad \quad (\mathcal{X}(z))_{12} \! = \! \int_{\mathbb{R}} 
\dfrac{\pmb{\pi}_{k}^{n}(\xi) \me^{-n \widetilde{V}(\xi)}}{\xi \! - \! z} 
\, \dfrac{\md \xi}{2 \pi \mi}.
\end{equation*}
For $n \! \in \! \mathbb{N}$ and $k \! \in \! \lbrace 1,2,\dotsc,K 
\rbrace$ such that $\alpha_{p_{\mathfrak{s}}} \! := \! \alpha_{k} 
\! = \! \infty$, let the Markov-Stieltjes transform, $\mathrm{F}_{
\tilde{\mu}}(z)$, the associated $\mathrm{R}$-function, 
$\widehat{\pmb{\mathrm{R}}}_{\tilde{\mu}}(z)$, and the corresponding 
{\rm MPA} error term, $\widehat{\pmb{\mathrm{E}}}_{\tilde{\mu}}
(z)$, be defined by Equations~\eqref{eqmvsstildemu}, 
\eqref{eqlemmvssinfmpa1}, and~\eqref{eqlemmvssinfmpa3}, 
respectively,\footnote{Note that (cf. Subsection~\ref{subsubsec1.2.1}, 
Equations~\eqref{mvssinf1}, \eqref{mvssinf4}, and~\eqref{mvssinf9}) 
$\mathrm{F}_{\tilde{\mu}}(z)$, $\widehat{\pmb{\mathrm{R}}}_{\tilde{\mu}}(z)$, 
and $\widehat{\pmb{\mathrm{E}}}_{\tilde{\mu}}(z)$ are $\mathrm{F}_{\mu}(z)$, 
$\widehat{\pmb{\mathrm{R}}}_{\mu}(z)$, and $\widehat{\pmb{\mathrm{E}}}_{\mu}
(z)$, respectively, under the transformation (cf. Remark~\ref{rem1.3.2}) 
$\md \mu \! \to \! \md \widetilde{\mu}$.} where, in particular,
\begin{equation} \label{eqlemmpaerrorinf1} 
\widehat{\pmb{\mathrm{E}}}_{\tilde{\mu}}(z) \! = \! 2 \pi 
\mi \dfrac{(\mathcal{X}(z))_{12}}{(\mathcal{X}(z))_{11}}.
\end{equation} 
Then, for $n \! \in \! \mathbb{N}$ and $k \! \in \! \lbrace 1,2,\dotsc,
K \rbrace$ such that $\alpha_{p_{\mathfrak{s}}} \! := \! \alpha_{k} 
\! = \! \infty$, in the double-scaling limit $\mathscr{N},n \! \to 
\! \infty$ such that $z_{o} \! = \! 1 \! + \! o(1)$, $\widehat{
\pmb{\mathrm{E}}}_{\tilde{\mu}}(z)$ has asymptotics derivable 
{}from Equation~\eqref{eqlemmpaerrorinf1}, where (cf. 
Theorem~\ref{maintheoforinf1}$)$$:$ {\rm \pmb{(1)}} for $z \! \in 
\! \hat{\Upsilon}_{1}$, $(\mathcal{X}(z))_{11}$ and $(\mathcal{X}
(z))_{12}$ have, respectively, asymptotics~\eqref{eqmaininf11} 
and~\eqref{eqmaininf12}$;$ {\rm \pmb{(2)}} for $z \! \in \! 
\hat{\Upsilon}_{2}$, $(\mathcal{X}(z))_{11}$ and $(\mathcal{X}
(z))_{12}$ have, respectively, asymptotics~\eqref{eqmaininf59} 
and~\eqref{eqmaininf60}$;$ {\rm \pmb{(3)}} for $z \! \in \! 
\hat{\Upsilon}_{3}$, $(\mathcal{X}(z))_{11}$ and $(\mathcal{X}
(z))_{12}$ have, respectively, asymptotics~\eqref{eqmaininf61} 
and~\eqref{eqmaininf62}$;$ {\rm \pmb{(4)}} for $z \! \in \! 
\hat{\Upsilon}_{4}$, $(\mathcal{X}(z))_{11}$ and $(\mathcal{X}
(z))_{12}$ have, respectively, asymptotics~\eqref{eqmaininf63} 
and~\eqref{eqmaininf64}$;$ {\rm \pmb{(5)}} for $z \! \in \! 
\hat{\Omega}^{1}_{\hat{b}_{j-1}}$, $j \! = \! 1,2,\dotsc,N \! + \! 1$, 
$(\mathcal{X}(z))_{11}$ and $(\mathcal{X}(z))_{12}$ have, respectively, 
asymptotics~\eqref{eqmaininf65} and~\eqref{eqmaininf66}$;$ 
{\rm \pmb{(6)}} for $z \! \in \! \hat{\Omega}^{1}_{\hat{a}_{j}}$, 
$j \! = \! 1,2,\dotsc,N \! + \! 1$, $(\mathcal{X}(z))_{11}$ and 
$(\mathcal{X}(z))_{12}$ have, respectively, 
asymptotics~\eqref{eqmaininf78} and~\eqref{eqmaininf79}$;$ 
{\rm \pmb{(7)}} for $z \! \in \! \hat{\Omega}^{2}_{\hat{b}_{j-1}}$, 
$j \! = \! 1,2,\dotsc,N \! + \! 1$, $(\mathcal{X}(z))_{11}$ and 
$(\mathcal{X}(z))_{12}$ have, respectively, 
asymptotics~\eqref{eqmaininf85} and~\eqref{eqmaininf86}$;$ 
{\rm \pmb{(8)}} for $z \! \in \! \hat{\Omega}^{2}_{\hat{a}_{j}}$, 
$j \! = \! 1,2,\dotsc,N \! + \! 1$, $(\mathcal{X}(z))_{11}$ and 
$(\mathcal{X}(z))_{12}$ have, respectively, 
asymptotics~\eqref{eqmaininf92} and~\eqref{eqmaininf93}$;$ 
{\rm \pmb{(9)}} for $z \! \in \! \hat{\Omega}^{3}_{\hat{b}_{j-1}}$, 
$j \! = \! 1,2,\dotsc,N \! + \! 1$, $(\mathcal{X}(z))_{11}$ and 
$(\mathcal{X}(z))_{12}$ have, respectively, 
asymptotics~\eqref{eqmaininf98} and~\eqref{eqmaininf99}$;$ 
{\rm \pmb{(10)}} for $z \! \in \! \hat{\Omega}^{3}_{\hat{a}_{j}}$, 
$j \! = \! 1,2,\dotsc,N \! + \! 1$, $(\mathcal{X}(z))_{11}$ and 
$(\mathcal{X}(z))_{12}$ have, respectively, 
asymptotics~\eqref{eqmaininf104} and~\eqref{eqmaininf105}$;$ 
{\rm \pmb{(11)}} for $z \! \in \! \hat{\Omega}^{4}_{\hat{b}_{j-1}}$, 
$j \! = \! 1,2,\dotsc,N \! + \! 1$, $(\mathcal{X}(z))_{11}$ and 
$(\mathcal{X}(z))_{12}$ have, respectively, 
asymptotics~\eqref{eqmaininf110} and~\eqref{eqmaininf111}$;$ 
and {\rm \pmb{(12)}} for $z \! \in \! \hat{\Omega}^{4}_{\hat{a}_{j}}$, 
$j \! = \! 1,2,\dotsc,N \! + \! 1$, $(\mathcal{X}(z))_{11}$ and 
$(\mathcal{X}(z))_{12}$ have, respectively, 
asymptotics~\eqref{eqmaininf116} and~\eqref{eqmaininf117}.
\end{dddd}
\begin{eeee} \label{remmpaasyinf}
\textsl{It is shown in Section~\ref{sec3} (see, in particular, the corresponding 
items of Lemmata~\ref{lemrootz} and~\ref{lemetatomu}$)$ that, for 
$n \! \in \! \mathbb{N}$ and $k \! \in \! \lbrace 1,2,\dotsc,K \rbrace$ 
such that $\alpha_{p_{\mathfrak{s}}} \! := \! \alpha_{k} \! = \! \infty$, 
$\lbrace \mathstrut z \! \in \! \overline{\mathbb{C}}; \, \pmb{\pi}^{n}_{k}
(z) \! = \! 0 \rbrace \! =: \! \lbrace \hat{\mathfrak{z}}^{n}_{k}(j) 
\rbrace_{j=1}^{(n-1)K+k} \! \subset \! \cup_{j=1}^{N+1}[\hat{b}_{j-1},
\hat{a}_{j}] \! =: \! J_{\infty}$$;$ therefore, using limiting values, if necessary, 
for $\mathbb{C} \cap (\hat{\Upsilon}_{3} \cup \hat{\Upsilon}_{4} \cup 
\cup_{j=1}^{N+1}(\hat{\Omega}^{2}_{\hat{b}_{j-1}} \cup \hat{\Omega}^{3}_{
\hat{b}_{j-1}} \cup \hat{\Omega}^{2}_{\hat{a}_{j}} \cup \hat{\Omega}^{3}_{
\hat{a}_{j}})) \! \ni \! z \! \to \! x \! \in \! J_{\infty}$, the asymptotics, 
in the double-scaling limit $\mathscr{N},n \! \to \! \infty$ such that 
$z_{o} \! = \! 1 \! + \! o(1)$, in Theorem~\ref{maintheompainf} for 
$\widehat{\pmb{\mathrm{E}}}_{\tilde{\mu}}(z)$ have a natural interpretation 
on $J_{\infty} \setminus \lbrace \hat{\mathfrak{z}}^{n}_{k}(j) 
\rbrace_{j=1}^{(n-1)K+k}$.}
\end{eeee}
\begin{dddd} \label{maintheoforinf2} 
Let the external field $\widetilde{V} \colon \overline{\mathbb{R}} \setminus 
\lbrace \alpha_{1},\alpha_{2},\dotsc,\alpha_{K} \rbrace \! \to \! \mathbb{R}$ 
satisfy conditions~\eqref{eq20}--\eqref{eq22}, and suppose that 
$\widetilde{V}$ is regular. For $n \! \in \! \mathbb{N}$ and $k \! \in \! 
\lbrace 1,2,\dotsc,K \rbrace$ such that $\alpha_{p_{\mathfrak{s}}} \! := \! 
\alpha_{k} \! = \! \infty$, let $\mathcal{X} \colon \overline{\mathbb{C}} 
\setminus \overline{\mathbb{R}} \! \to \! \mathrm{SL}_{2}(\mathbb{C})$ 
be the unique solution of the corresponding monic {\rm MPC ORF RHP} 
$(\mathcal{X}(z),\upsilon (z),\overline{\mathbb{R}})$ stated in 
Lemma~$\bm{\mathrm{RHP}_{\mathrm{MPC}}}$, and let the monic 
{\rm MPC ORF} have the representation $\pmb{\pi}_{k}^{n}(z) \! = \! 
(\mathcal{X}(z))_{11}$, $z \! \in \! \mathbb{C}$, with asymptotics, in the 
double-scaling limit $\mathscr{N},n \! \to \! \infty$ such that $z_{o} \! 
= \! 1 \! + \! o(1)$, stated in Theorem~\ref{maintheoforinf1}$;$ moreover, let 
the associated norming constant be given by $\mu_{n,\varkappa_{nk}}^{\infty}
(n,k) \! = \! \lvert \lvert \pmb{\pi}_{k}^{n}(\pmb{\cdot}) \rvert \rvert_{
\mathscr{L}}^{-1}$. Then, for $n \! \in \! \mathbb{N}$ and $k \! \in \! 
\lbrace 1,2,\dotsc,K \rbrace$ such that $\alpha_{p_{\mathfrak{s}}} \! := 
\! \alpha_{k} \! = \! \infty$,{}\footnote{The asymptotics for $\mu_{n,
\varkappa_{nk}}^{\infty}(n,k)$ is obtained by taking the positive square 
root of the expression on the right-hand side of Equation~\eqref{eqcoefinff1}.}
\begin{equation} \label{eqcoefinff1} 
(\mu_{n,\varkappa_{nk}}^{\infty}(n,k))^{2} \underset{\underset{z_{o}=
1+o(1)}{\mathscr{N},n \to \infty}}{=} \dfrac{2 \me^{-n \hat{\ell}}}{\pi} 
\hat{\mathbb{X}}^{\sharp} \left(1 \! + \! \dfrac{1}{(n \! - \! 1)K \! + \!  k} 
\hat{\mathbb{X}}^{\sharp} \left(\hat{\mathbb{X}}^{\natural} \right)_{12} 
\! + \! \mathcal{O} \left(\dfrac{\hat{\mathfrak{c}}(n,k,z_{o})}{((n \! - \! 1)
K \! + \! k)^{2}} \right) \right),
\end{equation}
where
\begin{equation} \label{eqcoefinff2} 
\hat{\mathbb{X}}^{\sharp} \! = \! \left(\sum_{j=1}^{N+1}(\hat{a}_{j} 
\! - \! \hat{b}_{j-1}) \right)^{-1} \dfrac{\hat{\boldsymbol{\theta}}
(-\hat{\boldsymbol{u}}_{+}(\infty) \! + \! \hat{\boldsymbol{d}}) 
\hat{\boldsymbol{\theta}}(\hat{\boldsymbol{u}}_{+}(\infty) \! - \! 
\frac{1}{2 \pi}((n \! - \! 1)K \! + \! k) \hat{\boldsymbol{\Omega}} 
\! + \! \hat{\boldsymbol{d}})}{\hat{\boldsymbol{\theta}}
(\hat{\boldsymbol{u}}_{+}(\infty) \! + \! \hat{\boldsymbol{d}}) 
\hat{\boldsymbol{\theta}}(-\hat{\boldsymbol{u}}_{+}(\infty) \! - \! 
\frac{1}{2 \pi}((n \! - \! 1)K \! + \! k) \hat{\boldsymbol{\Omega}} 
\! + \! \hat{\boldsymbol{d}})},
\end{equation}
and
\begin{equation} \label{eqcoefinff3} 
\hat{\mathbb{X}}^{\natural} \! = \! 4 \mi \sum_{j=1}^{N+1} \left(
(\hat{\alpha}_{0}(\hat{b}_{j-1}))^{-2} \left(\hat{\alpha}_{0}(\hat{b}_{j-1}) 
\hat{\boldsymbol{\mathrm{B}}}(\hat{b}_{j-1}) \! - \! \hat{\alpha}_{1}
(\hat{b}_{j-1}) \hat{\boldsymbol{\mathrm{A}}}(\hat{b}_{j-1}) \right) \! 
+ \! (\hat{\alpha}_{0}(\hat{a}_{j}))^{-2} \left(\hat{\alpha}_{0}(\hat{a}_{j}) 
\hat{\boldsymbol{\mathrm{B}}}(\hat{a}_{j}) \! - \! \hat{\alpha}_{1}
(\hat{a}_{j}) \hat{\boldsymbol{\mathrm{A}}}(\hat{a}_{j}) \right) \right),
\end{equation}
with $\hat{\mathfrak{c}}(n,k,z_{o}) \! =_{\underset{z_{o}=1+o(1)}{\mathscr{N},
n \to \infty}} \! \mathcal{O}(1)$.

For $n \! \in \! \mathbb{N}$ and $k \! \in \! \lbrace 1,2,\dotsc,K \rbrace$ 
such that $\alpha_{p_{\mathfrak{s}}} \! := \! \alpha_{k} \! = \! \infty$, let 
the associated {\rm MPC ORF} be given by $\phi_{k}^{n}(z) \! = \! \mu_{n,
\varkappa_{nk}}^{\infty}(n,k) \pmb{\pi}_{k}^{n}(z)$, $z \! \in \! \mathbb{C}$$:$ 
the asymptotics, in the double-scaling limit $\mathscr{N},n \! \to \! \infty$ 
such that $z_{o} \! = \! 1 \! + \! o(1)$, for $\phi_{k}^{n}(z)$ are derived via 
the scalar multiplication of the asymptotics for $\pmb{\pi}_{k}^{n}(z)$ 
and $\mu_{n,\varkappa_{nk}}^{\infty}(n,k)$ stated, respectively, in 
Theorem~\ref{maintheoforinf1} and Equation~\eqref{eqcoefinff1}.
\end{dddd}

For $n \! \in \! \mathbb{N}$ and $k \! \in \! \lbrace 1,2,\dotsc,K \rbrace$ 
such that $\alpha_{p_{\mathfrak{s}}} \! := \! \alpha_{k} \! \neq \! \infty$, 
and $j \! = \! 1,2,\dotsc,N \! + \! 1$, let
\begin{gather}
\tilde{\Phi}_{\tilde{b}_{j-1}}(z) \! := \! \left(-\dfrac{3}{4}((n \! - \! 1)K \! + \! 
k) \int_{z}^{\tilde{b}_{j-1}}(\tilde{R}(\xi))^{1/2} \tilde{h}_{\widetilde{V}}(\xi) 
\, \md \xi \right)^{2/3}, \label{eqmainfin78} \\
\tilde{\Phi}_{\tilde{a}_{j}}(z) \! := \! \left(\dfrac{3}{4}((n \! - \! 1)K \! + \! k) 
\int_{\tilde{a}_{j}}^{z}(\tilde{R}(\xi))^{1/2} \tilde{h}_{\widetilde{V}}(\xi) \, 
\md \xi \right)^{2/3}, \label{eqmainfin85}
\end{gather}
and define the mutually disjoint open discs about $\tilde{b}_{j-1}$ and 
$\tilde{a}_{j}$, respectively, as $\tilde{\mathbb{U}}_{\tilde{\delta}_{\tilde{
b}_{j-1}}} \! := \! \lbrace \mathstrut z \! \in \! \mathbb{C}; \, \lvert z \! - 
\! \tilde{b}_{j-1} \rvert \! < \! \tilde{\delta}_{\tilde{b}_{j-1}} \rbrace$ and 
$\tilde{\mathbb{U}}_{\tilde{\delta}_{\tilde{a}_{j}}} \! := \! \lbrace \mathstrut 
z \! \in \! \mathbb{C}; \, \lvert z \! - \! \tilde{a}_{j} \rvert \! < \! \tilde{
\delta}_{\tilde{a}_{j}} \rbrace$, where $\tilde{\delta}_{\tilde{b}_{j-1}},\tilde{
\delta}_{\tilde{a}_{j}} \! \in \! (0,1)$ are sufficiently small, positive real 
numbers chosen so that $\tilde{\mathbb{U}}_{\tilde{\delta}_{\tilde{b}_{i-1}}} 
\cap \tilde{\mathbb{U}}_{\tilde{\delta}_{\tilde{a}_{j}}} \! = \! \varnothing$ 
$\forall$ $i,j \! \in \! \lbrace 1,2,\dotsc,N \! + \! 1 \rbrace$, $\tilde{
\mathbb{U}}_{\tilde{\delta}_{\tilde{b}_{j-1}}} \cap \lbrace \alpha_{p_{1}},
\alpha_{p_{2}},\dotsc,\alpha_{p_{\mathfrak{s}}} \rbrace \! = \! \varnothing 
\! = \! \tilde{\mathbb{U}}_{\tilde{\delta}_{\tilde{a}_{j}}} \cap \lbrace 
\alpha_{p_{1}},\alpha_{p_{2}},\dotsc,\alpha_{p_{\mathfrak{s}}} \rbrace$, and 
$\tilde{\Phi}_{\tilde{b}_{j-1}}(z)$ (resp., $\tilde{\Phi}_{\tilde{a}_{j}}(z))$, 
which are bi-holomorphic, conformal, and non-orientation preserving 
(resp., orientation preserving), map (see Figure~\ref{forbtil}) $\tilde{
\mathbb{U}}_{\tilde{\delta}_{\tilde{b}_{j-1}}}$ (resp., (see Figure~\ref{foratil}) 
$\tilde{\mathbb{U}}_{\tilde{\delta}_{\tilde{a}_{j}}})$, and, thus, the oriented 
contours $\tilde{\Sigma}_{\tilde{b}_{j-1}} \! := \! \cup_{m=1}^{4} 
\tilde{\Sigma}_{\tilde{b}_{j-1}}^{m}$ (resp., $\tilde{\Sigma}_{\tilde{a}_{j}} 
\! := \! \cup_{m=1}^{4} \tilde{\Sigma}_{\tilde{a}_{j}}^{m})$, injectively onto 
open, $n$- and $k$-dependent, neighbourhoods $\tilde{\mathbb{U}}_{
\tilde{\delta}_{\tilde{b}_{j-1}}}^{\ast}$ (resp., $\tilde{\mathbb{U}}_{\tilde{
\delta}_{\tilde{a}_{j}}}^{\ast})$ of the origin, such that $\tilde{\Phi}_{
\tilde{b}_{j-1}}(\tilde{b}_{j-1}) \! = \! 0$, $\tilde{\Phi}_{\tilde{b}_{j-1}} 
\colon \tilde{\mathbb{U}}_{\tilde{\delta}_{\tilde{b}_{j-1}}} \! \to \! 
\tilde{\mathbb{U}}_{\tilde{\delta}_{\tilde{b}_{j-1}}}^{\ast} \! := \! 
\tilde{\Phi}_{\tilde{b}_{j-1}}(\tilde{\mathbb{U}}_{\tilde{\delta}_{\tilde{b}_{j-1}}})$, 
$\tilde{\Phi}_{\tilde{b}_{j-1}}(\tilde{\mathbb{U}}_{\tilde{\delta}_{\tilde{b}_{j-1}}} 
\cap \tilde{\Sigma}_{\tilde{b}_{j-1}}^{m}) \! = \! \tilde{\Phi}_{\tilde{b}_{j-1}}
(\tilde{\mathbb{U}}_{\tilde{\delta}_{\tilde{b}_{j-1}}}) \cap \gamma_{
\tilde{b}_{j-1}}^{\ast,m}$, and $\tilde{\Phi}_{\tilde{b}_{j-1}}
(\tilde{\mathbb{U}}_{\tilde{\delta}_{\tilde{b}_{j-1}}} \cap \tilde{\Omega}_{
\tilde{b}_{j-1}}^{m}) \! = \! \tilde{\Phi}_{\tilde{b}_{j-1}}(\tilde{\mathbb{U}}_{
\tilde{\delta}_{\tilde{b}_{j-1}}}) \cap \tilde{\Omega}_{\tilde{b}_{j-1}}^{\ast,m}$, 
$m \! = \! 1,2,3,4$, with $\tilde{\Omega}_{\tilde{b}_{j-1}}^{\ast,1} \! = \! 
\lbrace \mathstrut \zeta \! \in \! \mathbb{C}; \, \arg (\zeta) \! \in \! (0,
2 \pi/3) \rbrace$, $\tilde{\Omega}_{\tilde{b}_{j-1}}^{\ast,2} \! = \! \lbrace 
\mathstrut \zeta \! \in \! \mathbb{C}; \, \arg (\zeta) \! \in \! (2 \pi/3,\pi) 
\rbrace$, $\tilde{\Omega}_{\tilde{b}_{j-1}}^{\ast,3} \! = \! \lbrace \mathstrut 
\zeta \! \in \! \mathbb{C}; \, \arg (\zeta) \! \in \! (-\pi,-2 \pi/3) \rbrace$, 
and $\tilde{\Omega}_{\tilde{b}_{j-1}}^{\ast,4} \! = \! \lbrace \mathstrut 
\zeta \! \in \! \mathbb{C}; \, \arg (\zeta) \! \in \! (-2 \pi/3,0) \rbrace$ 
(resp., $\tilde{\Phi}_{\tilde{a}_{j}}(\tilde{a}_{j}) \! = \! 0$, $\tilde{\Phi}_{
\tilde{a}_{j}} \colon \tilde{\mathbb{U}}_{\tilde{\delta}_{\tilde{a}_{j}}} \! \to 
\! \tilde{\mathbb{U}}_{\tilde{\delta}_{\tilde{a}_{j}}}^{\ast} \! := \! \tilde{
\Phi}_{\tilde{a}_{j}}(\tilde{\mathbb{U}}_{\tilde{\delta}_{\tilde{a}_{j}}})$, 
$\tilde{\Phi}_{\tilde{a}_{j}}(\tilde{\mathbb{U}}_{\tilde{\delta}_{\tilde{a}_{j}}} 
\cap \tilde{\Sigma}_{\tilde{a}_{j}}^{m}) \! = \! \tilde{\Phi}_{\tilde{a}_{j}}
(\tilde{\mathbb{U}}_{\tilde{\delta}_{\tilde{a}_{j}}}) \cap \gamma_{\tilde{a}_{j}}^{
\ast,m}$, and $\tilde{\Phi}_{\tilde{a}_{j}}(\tilde{\mathbb{U}}_{\tilde{\delta}_{
\tilde{a}_{j}}} \cap \tilde{\Omega}_{\tilde{a}_{j}}^{m}) \! = \! \tilde{\Phi}_{
\tilde{a}_{j}}(\tilde{\mathbb{U}}_{\tilde{\delta}_{\tilde{a}_{j}}}) \cap \tilde{
\Omega}_{\tilde{a}_{j}}^{\ast,m}$, $m \! = \! 1,2,3,4$, with $\tilde{\Omega}_{
\tilde{a}_{j}}^{\ast,1} \! = \! \lbrace \mathstrut \zeta \! \in \! \mathbb{C}; \, 
\arg (\zeta) \! \in \! (0,2 \pi/3) \rbrace$, $\tilde{\Omega}_{\tilde{a}_{j}}^{\ast,2} 
\! = \! \lbrace \mathstrut \zeta \! \in \! \mathbb{C}; \, \arg (\zeta) \! \in \! 
(2 \pi/3,\pi) \rbrace$, $\tilde{\Omega}_{\tilde{a}_{j}}^{\ast,3} \! = \! \lbrace 
\mathstrut \zeta \! \in \! \mathbb{C}; \, \arg (\zeta) \! \in \! (-\pi,-2 \pi/3) 
\rbrace$, and $\tilde{\Omega}_{\tilde{a}_{j}}^{\ast,4} \! = \! \lbrace \mathstrut 
\zeta \! \in \! \mathbb{C}; \, \arg (\zeta) \! \in \! (-2 \pi/3,0) \rbrace)$. Consider, 
also, the decomposition of $\mathbb{C}$ (see Figure~\ref{figsectortil}) into 
bounded and unbounded regions, and open neighbourhoods surrounding the 
end-points of the intervals, $\lbrace \tilde{b}_{j-1},\tilde{a}_{j} \rbrace_{j=
1}^{N+1}$, of the support, $J_{f}$, of the associated equilibrium measure, 
$\mu_{\widetilde{V}}^{f}$.
\begin{figure}[bht]
\begin{center}
\vspace{0.55cm}

\end{center}
\caption{For $n \! \in \! \mathbb{N}$ and $k \! \in \! \lbrace 1,2,\dotsc,
K \rbrace$ such that $\alpha_{p_{\mathfrak{s}}} \! := \! \alpha_{k} \! \neq \! 
\infty$, the conformal mapping $\zeta \! = \! \tilde{\Phi}_{\tilde{b}_{j-1}}(z)$, 
$j \! = \! 1,2,\dotsc,N \! + \! 1$, where $(\tilde{\Phi}_{\tilde{b}_{j-1}})^{-1}$ 
denotes the inverse mapping.}\label{forbtil}
\end{figure}
\begin{figure}[tbh]
\begin{center}
\vspace{0.55cm}
%
\end{center}
\caption{For $n \! \in \! \mathbb{N}$ and $k \! \in \! \lbrace 1,2,\dotsc,
K \rbrace$ such that $\alpha_{p_{\mathfrak{s}}} \! := \! \alpha_{k} \! \neq \! 
\infty$, the conformal mapping $\zeta \! = \! \tilde{\Phi}_{\tilde{a}_{j}}(z)$, 
$j \! = \! 1,2,\dotsc,N \! + \! 1$, where $(\tilde{\Phi}_{\tilde{a}_{j}})^{-1}$ 
denotes the inverse mapping.}\label{foratil}
\end{figure}
\begin{figure}[tbh]
\begin{center}
\vspace{0.55cm}
%
\end{center}
\caption{The open neighbourhoods surrounding the end-points of the 
intervals, $\lbrace \tilde{b}_{j-1},\tilde{a}_{j} \rbrace_{j=1}^{N+1}$, 
of  the support, $J_{f}$, of the associated equilibrium measure, 
$\mu_{\widetilde{V}}^{f}$.}\label{figsectortil}
\end{figure}
\begin{dddd} \label{maintheoforfin1} 
Let the external field $\widetilde{V} \colon \overline{\mathbb{R}} \setminus 
\lbrace \alpha_{1},\alpha_{2},\dotsc,\alpha_{K} \rbrace \! \to \! \mathbb{R}$ 
satisfy conditions~\eqref{eq20}--\eqref{eq22}. For $n \! \in \! \mathbb{N}$ 
and $k \! \in \! \lbrace 1,2,\dotsc,K \rbrace$ such that $\alpha_{
p_{\mathfrak{s}}} \! := \! \alpha_{k} \! \neq \! \infty$, let $\md 
\mu_{\widetilde{V}}^{f}(x) \! = \! \psi_{\widetilde{V}}^{f}(x) \, \md x \! 
= \! (2 \pi \mi)^{-1}(\tilde{R}(x))^{1/2}_{+} \tilde{h}_{\widetilde{V}}
(x) \chi_{J_{f}}(x) \, \md x$, where $(\tilde{R}(z))^{1/2}$ is defined 
by Equation~\eqref{eql3.7j}, with $(\tilde{R}(x))^{1/2}_{\pm} \! 
:= \! \lim_{\varepsilon \downarrow 0}(\tilde{R}(x \! \pm \! \mi 
\varepsilon))^{1/2}$,\footnote{The branch of the square root is chosen 
so that $z^{-(N+1)}(\tilde{R}(z))^{1/2} \! \sim_{\mathbb{C}_{\pm} \ni 
z \to \alpha_{p_{\mathfrak{s}-1}}} \! \pm 1$ $(\alpha_{p_{\mathfrak{s}-1}} 
\! = \! \infty)$.} $J_{f} \! := \! \supp (\mu_{\widetilde{V}}^{f}) \! = \! 
\cup_{j=1}^{N+1}[\tilde{b}_{j-1},\tilde{a}_{j}]$, where $N \! \in \! \mathbb{N}_{0}$ 
and is finite, $[\tilde{b}_{i-1},\tilde{a}_{i}] \cap [\tilde{b}_{j-1},\tilde{a}_{j}] \! = 
\! \varnothing$ $\forall$ $i \! \neq \! j \! \in \! \lbrace 1,2,\dotsc,N \! + \! 1 
\rbrace$, $[\tilde{b}_{j-1},\tilde{a}_{j}] \cap \lbrace \alpha_{p_{1}},\alpha_{p_{2}},
\dotsc,\alpha_{p_{\mathfrak{s}}} \rbrace \! = \! \varnothing$, and $-\infty \! 
< \! \tilde{b}_{0} \! < \! \tilde{a}_{1} \! < \! \tilde{b}_{1} \! < \! \tilde{a}_{2} \! 
< \! \dotsb \! < \! \tilde{b}_{N} \! < \! \tilde{a}_{N+1} \! < \! +\infty$, with 
$\lbrace \tilde{b}_{j-1},\tilde{a}_{j} \rbrace_{j=1}^{N+1}$ satisfying the 
associated $n$- and $k$-dependent system of $2(N \! + \! 1)$ real moment 
Equations~\eqref{eql3.7g}--\eqref{eql3.7i}, $\tilde{h}_{\widetilde{V}}(z)$, 
which is real analytic for $z \! \in \! \overline{\mathbb{R}} \setminus \lbrace 
\alpha_{p_{1}},\alpha_{p_{2}},\dotsc,\alpha_{p_{\mathfrak{s}}} \rbrace$, is 
defined by Equation~\eqref{eql3.7l},\footnote{Note: in the definition of 
$\tilde{h}_{\widetilde{V}}(z)$ given by Equation~\eqref{eql3.7l}, there appears 
the contour integral $\oint_{\tilde{C}_{\widetilde{V}}}$, the detailed description 
of which is given in the corresponding item of Lemma~\ref{lem3.7}.} and 
$\chi_{J_{f}}(x)$ is the characteristic function of the compact set $J_{f}$. 
Furthermore, suppose that, for $n \! \in \! \mathbb{N}$ and $k \! \in 
\! \lbrace 1,2,\dotsc,K \rbrace$ such that $\alpha_{p_{\mathfrak{s}}} 
\! := \! \alpha_{k} \! \neq \! \infty$, the external field, $\widetilde{V}$, 
is regular:
\begin{description}
\item[{\rm (i)}] $\tilde{h}_{\widetilde{V}}(x) \! \neq \! 0$, 
$x \! \in \! J_{f}$$;$
\item[{\rm (ii)}]
\begin{equation*}
2 \left(\dfrac{(n \! - \! 1)K \! + \! k}{n} \right) \int_{J_{f}} \ln \left(\left\vert 
\dfrac{x \! - \! \xi}{\xi \! - \! \alpha_{k}} \right\vert \right) \md 
\mu^{f}_{\widetilde{V}}(\xi) \! - \! 2 \sum_{q=1}^{\mathfrak{s}-2} 
\dfrac{\varkappa_{nk \tilde{k}_{q}}}{n} \ln \left\lvert \dfrac{x \! - \! 
\alpha_{p_{q}}}{\alpha_{p_{q}} \! - \! \alpha_{k}} \right\rvert \! - \! 2 
\left(\dfrac{\varkappa_{nk} \! - \! 1}{n} \right) \ln \lvert x \! - \! 
\alpha_{k} \rvert \! - \! \widetilde{V}(x) \! - \! \tilde{\ell} \! = \! 0, 
\quad x \! \in \! J_{f},
\end{equation*}
which defines the associated variational constant $\tilde{\ell}$ 
$(\in \! \mathbb{R})$,\footnote{Note that $\tilde{\ell}$ is the same 
on each compact real interval $[\tilde{b}_{j-1},\tilde{a}_{j}]$, 
$j \! = \! 1,2,\dotsc,N \! + \! 1$.} and
\begin{equation*}
2 \left(\dfrac{(n \! - \! 1)K \! + \! k}{n} \right) \int_{J_{f}} \ln \left(\left\vert 
\dfrac{x \! - \! \xi}{\xi \! - \! \alpha_{k}} \right\vert \right) \md 
\mu^{f}_{\widetilde{V}}(\xi) \! - \! 2 \sum_{q=1}^{\mathfrak{s}-2} 
\dfrac{\varkappa_{nk \tilde{k}_{q}}}{n} \ln \left\lvert \dfrac{x \! - \! 
\alpha_{p_{q}}}{\alpha_{p_{q}} \! - \! \alpha_{k}} \right\rvert \! - \! 2 
\left(\dfrac{\varkappa_{nk} \! - \! 1}{n} \right) \ln \lvert x \! - \! 
\alpha_{k} \rvert \! - \! \widetilde{V}(x) \! - \! \tilde{\ell} \! < \! 0, 
\quad x \! \in \! \mathbb{R} \setminus J_{f};
\end{equation*}
\item[{\rm (iii)}]
\begin{equation*}
g^{f}_{+}(x) \! + \! g^{f}_{-}(x) \! - \! \hat{\mathscr{P}}_{0}^{+} \! - \! 
\hat{\mathscr{P}}_{0}^{-} \! - \! \widetilde{V}(x) \! - \! \tilde{\ell} \! 
< \! 0, \quad x \! \in \! \mathbb{R} \setminus J_{f},
\end{equation*}
where, for $z \! \in \! \mathbb{C} \setminus (-\infty,\max \lbrace 
\max_{q=1,\dotsc,\mathfrak{s}-2,\mathfrak{s}} \lbrace \alpha_{p_{q}} 
\rbrace,\max \lbrace J_{f} \rbrace \rbrace)$, $g^{f}(z)$ is 
defined by Equation~\eqref{eql3.4gee3}, $g^{f}_{\pm}(x) \! 
:= \! \lim_{\varepsilon \downarrow 0}g^{f}(x \! \pm \! \mi 
\varepsilon)$, and $\hat{\mathscr{P}}_{0}^{\pm}$ is defined 
by Equation~\eqref{eql3.4gee5}$;$
\item[{\rm (iv)}]
\begin{equation*}
\mi \left(g^{f}_{+}(z) \! - \! g^{f}_{-}(z) \! + \! \hat{\mathscr{P}}_{0}^{-} 
\! - \! \hat{\mathscr{P}}_{0}^{+} \! + \! 2 \pi \mi \left(\dfrac{\varkappa_{nk}
-1}{n} \right) \chi_{\lbrace \mathstrut z \in \mathbb{R}; \, z< \alpha_{k} 
\rbrace}(z) \! + \! 2 \pi \mi \sum_{q \in \lbrace \mathstrut j \in \lbrace 
1,2,\dotsc,\mathfrak{s}-2 \rbrace; \, \alpha_{p_{j}} >z \rbrace} \dfrac{
\varkappa_{nk \tilde{k}_{q}}}{n} \right)^{\prime} \! > \! 0, \quad z \! 
\in \! J_{f},
\end{equation*}
where the prime denotes differentiation with respect to $z$.
\end{description}

For $n \! \in \! \mathbb{N}$ and $k \! \in \! \lbrace 1,2,\dotsc,K 
\rbrace$ such that $\alpha_{p_{\mathfrak{s}}} \! := \! \alpha_{k} 
\! \neq \! \infty$, let
\begin{gather}
\tilde{\mathbb{K}} \! = \! \dfrac{1}{\mathfrak{k}_{0}} 
, \label{eqmainfin9}
\end{gather}
where {}\footnote{Note: $\det (\tilde{\mathbb{K}} \tilde{\mathbb{M}}(z)) 
\! = \! 1$.}
\begin{align} \label{eqconsfinn1} 
\mathfrak{k}_{0} :=& \, \dfrac{1}{4} \left(\tilde{\gamma}(\alpha_{k}) 
\! + \! (\tilde{\gamma}(\alpha_{k}))^{-1} \right)^{2} \dfrac{\tilde{
\boldsymbol{\theta}}(\tilde{\boldsymbol{u}}_{+}(\alpha_{k}) \! - \! 
\frac{1}{2 \pi}((n \! - \! 1)K \! + \! k) \tilde{\boldsymbol{\Omega}} 
\! + \! \tilde{\boldsymbol{d}})}{\tilde{\boldsymbol{\theta}}(\tilde{
\boldsymbol{u}}_{+}(\alpha_{k}) \! + \! \tilde{\boldsymbol{d}})} 
\dfrac{\tilde{\boldsymbol{\theta}}(-\tilde{\boldsymbol{u}}_{+}
(\alpha_{k}) \! - \! \frac{1}{2 \pi}((n \! - \! 1)K \! + \! k) \tilde{
\boldsymbol{\Omega}} \! - \! \tilde{\boldsymbol{d}})}{\tilde{
\boldsymbol{\theta}}(-\tilde{\boldsymbol{u}}_{+}(\alpha_{k}) 
\! - \! \tilde{\boldsymbol{d}})} \nonumber \\
-& \, \dfrac{1}{4} \left(\tilde{\gamma}(\alpha_{k}) \! - \! (\tilde{\gamma}
(\alpha_{k}))^{-1} \right)^{2} \dfrac{\tilde{\boldsymbol{\theta}}(\tilde{
\boldsymbol{u}}_{+}(\alpha_{k}) \! - \! \frac{1}{2 \pi}((n \! - \! 1)K \! 
+ \! k) \tilde{\boldsymbol{\Omega}} \! - \! \tilde{\boldsymbol{d}})}{
\tilde{\boldsymbol{\theta}}(\tilde{\boldsymbol{u}}_{+}(\alpha_{k}) \! 
- \! \tilde{\boldsymbol{d}})} \dfrac{\tilde{\boldsymbol{\theta}}(-
\tilde{\boldsymbol{u}}_{+}(\alpha_{k}) \! - \! \frac{1}{2 \pi}((n \! - \! 
1)K \! + \! k) \tilde{\boldsymbol{\Omega}} \! + \! \tilde{\boldsymbol{
d}})}{\tilde{\boldsymbol{\theta}}(-\tilde{\boldsymbol{u}}_{+}(\alpha_{k}) 
\! + \! \tilde{\boldsymbol{d}})},
\end{align}
$\tilde{\gamma}(z)$ and $\tilde{\gamma}(\alpha_{k})$ are defined by 
Equations~\eqref{eqmainfin10} and~\eqref{eqssabra2}, respectively, 
$\tilde{\boldsymbol{u}}(z) \! := \! \int_{\tilde{a}_{N+1}}^{z} \tilde{
\boldsymbol{\omega}}$, $\tilde{\boldsymbol{u}}_{+}(\alpha_{k}) \! := \! 
\int_{\tilde{a}_{N+1}}^{\alpha_{k}^{+}} \tilde{\boldsymbol{\omega}}$, 
with $\tilde{\boldsymbol{\omega}}$ the associated normalised basis 
of holomorphic one-forms on $\tilde{\mathcal{Y}}$, 
$\tilde{\boldsymbol{\Omega}} \! := \! (\tilde{\Omega}_{1},
\tilde{\Omega}_{2},\dotsc,\tilde{\Omega}_{N})^{\operatorname{T}}$ 
$(\in \! \mathbb{R}^{N})$, where $\tilde{\Omega}_{j} \! = \! 2 \pi 
\int_{\tilde{b}_{j}}^{\tilde{a}_{N+1}} \psi_{\widetilde{V}}^{f}(\xi) \, 
\md \xi$, $j \! = \! 1,2,\dotsc,N$, and $\tilde{\boldsymbol{d}} \! 
\equiv \! -\sum_{j=1}^{N} \int_{\tilde{a}_{j}}^{\tilde{z}_{j}^{-}} 
\tilde{\boldsymbol{\omega}}$ $(= \! \sum_{j=1}^{N} \int_{\tilde{a}_{j}}^{
\tilde{z}_{j}^{+}} \tilde{\boldsymbol{\omega}})$, where $\lbrace 
\tilde{z}_{1}^{\pm},\tilde{z}_{2}^{\pm},\dotsc,\tilde{z}_{N}^{\pm} \rbrace 
\! := \! \lbrace \mathstrut z \! \in \! \mathbb{C}_{\pm}; \, \tilde{\gamma}
(z) \! \mp \! (\tilde{\gamma}(z))^{-1} \! = \! 0 \rbrace$, with 
$\tilde{z}_{j}^{\pm} \! \in \! (\tilde{a}_{j},\tilde{b}_{j})^{\pm}$ $(\subset 
\mathbb{C}_{\pm})$, $j \! = \! 1,2,\dotsc,N$.\footnote{As points on the 
plane, $\boldsymbol{\mathrm{pr}}(\tilde{z}_{j}^{+}) \! = \! \boldsymbol{
\mathrm{pr}}(\tilde{z}_{j}^{-}) \! = \! \tilde{z}_{j} \! \in \! (\tilde{a}_{j},
\tilde{b}_{j})$, $j \! = \! 1,2,\dotsc,N$.}

For $n \! \in \! \mathbb{N}$ and $k \! \in \! \lbrace 1,2,\dotsc,K \rbrace$ 
such that $\alpha_{p_{\mathfrak{s}}} \! := \! \alpha_{k} \! \neq \! \infty$, 
let $\mathcal{X} \colon \overline{\mathbb{C}} \setminus \overline{\mathbb{
R}} \! \to \! \mathrm{SL}_{2}(\mathbb{C})$ be the unique solution of the 
corresponding monic {\rm MPC ORF RHP} $(\mathcal{X}(z),\upsilon (z),
\overline{\mathbb{R}})$ stated in Lemma~$\bm{\mathrm{RHP}_{\mathrm{MPC}}}$, 
with integral representation given by Equation~\eqref{intrepfin} 
(see Lemma~\ref{lem2.1}$)$$;$ in particular,
\begin{equation*}
(\mathcal{X}(z))_{11} \! = \! (z \! - \! \alpha_{k}) \pmb{\pi}_{k}^{n}(z) 
\quad \quad \, \text{and} \, \quad \quad (\mathcal{X}(z))_{12} \! = \! 
(z \! - \! \alpha_{k}) \int_{\mathbb{R}} \dfrac{((\xi \! - \! \alpha_{k}) 
\pmb{\pi}_{k}^{n}(\xi)) \me^{-n \widetilde{V}(\xi)}}{(\xi \! - \! 
\alpha_{k})(\xi \! - \! z)} \, \dfrac{\md \xi}{2 \pi \mi}.
\end{equation*}
Then, for $n \! \in \! \mathbb{N}$ and $k \! \in \! \lbrace 1,2,\dotsc,K 
\rbrace$ such that $\alpha_{p_{\mathfrak{s}}} \! := \! \alpha_{k} \! \neq \! 
\infty$$:$\\
{\rm \pmb{(1)}} for $z \! \in \! \tilde{\Upsilon}_{1}$,
\begin{align}
(z \! - \! \alpha_{k}) \pmb{\pi}_{k}^{n}(z) \underset{\underset{z_{o}=
1+o(1)}{\mathscr{N},n \to \infty}}{=}& \, \left(\left(\tilde{\mathbb{K}}_{11} 
\tilde{\mathbb{M}}_{11}(z) \! + \! \tilde{\mathbb{K}}_{12} \tilde{\mathbb{
M}}_{21}(z) \right) \left(1 \! + \! \dfrac{1}{(n \! - \! 1)K \! + \! k} \tilde{
\mathcal{R}}_{11}^{\sharp}(z) \! + \! \mathcal{O} \left(\dfrac{\mathfrak{c}
(n,k,z_{o})}{((n \! - \! 1)K \! + \! k)^{2}} \right) \right) \right. \nonumber \\
+&\left. \, \mathscr{E}^{2} \left(\tilde{\mathbb{K}}_{21} \tilde{\mathbb{
M}}_{11}(z) \! + \! \tilde{\mathbb{K}}_{22} \tilde{\mathbb{M}}_{21}(z) 
\right) \left(\dfrac{1}{(n \! - \! 1)K \! + \! k} \tilde{\mathcal{R}}_{12}^{
\sharp}(z) \! + \! \mathcal{O} \left(\dfrac{\mathfrak{c}(n,k,z_{o})}{
((n \! - \! 1)K \! + \! k)^{2}} \right) \right) \right) \nonumber \\
\times& \, \me^{n(g^{f}(z)-\hat{\mathscr{P}}_{0}^{+})}, 
\label{eqmainfin11}
\end{align}
and
\begin{align}
(z \! - \! \alpha_{k}) \int_{\mathbb{R}} \dfrac{((\xi \! - \! \alpha_{k}) 
\pmb{\pi}_{k}^{n}(\xi)) \me^{-n \widetilde{V}(\xi)}}{(\xi \! - \! \alpha_{k})
(\xi \! - \! z)} \, \dfrac{\md \xi}{2 \pi \mi} \underset{\underset{z_{o}=1+o(1)}{
\mathscr{N},n \to \infty}}{=}& \, \left(\mathscr{E}^{-2} \left(\tilde{\mathbb{
K}}_{11} \tilde{\mathbb{M}}_{12}(z) \! + \! \tilde{\mathbb{K}}_{12} \tilde{
\mathbb{M}}_{22}(z) \right) \left(1 \! + \! \dfrac{1}{(n \! - \! 1)K \! + \! k} 
\tilde{\mathcal{R}}_{11}^{\sharp}(z) \right. \right. \nonumber \\
+&\left. \left. \, \mathcal{O} \left(\dfrac{\mathfrak{c}(n,k,z_{o})}{((n \! - 
\! 1)K \! + \! k)^{2}} \right) \right) \! + \! \left(\tilde{\mathbb{K}}_{21} 
\tilde{\mathbb{M}}_{12}(z) \! + \! \tilde{\mathbb{K}}_{22} 
\tilde{\mathbb{M}}_{22}(z) \right) \right. \nonumber \\
\times&\left. \, \left(\dfrac{1}{(n \! - \! 1)K \! + \! k} 
\tilde{\mathcal{R}}_{12}^{\sharp}(z) \! + \! \mathcal{O} \left(
\dfrac{\mathfrak{c}(n,k,z_{o})}{((n \! - \! 1)K \! + \! k)^{2}} \right) 
\right) \right) \nonumber \\
\times& \, \me^{n \tilde{\ell}} \me^{-n(g^{f}(z)-\hat{\mathscr{P}}_{0}^{+})}, 
\label{eqmainfin12}
\end{align}
with $\mathscr{E}$ defined by Equation~\eqref{eqmainfin13}, and
\begin{align}
\tilde{\mathcal{R}}^{\sharp}(z) =& \, \sum_{j=1}^{N+1} \left(\dfrac{(\tilde{
\alpha}_{0}(\tilde{b}_{j-1}))^{-1}}{(z \! - \! \tilde{b}_{j-1})^{2}} \tilde{
\boldsymbol{\mathrm{A}}}(\tilde{b}_{j-1}) \! + \! \dfrac{(\tilde{\alpha}_{0}
(\tilde{b}_{j-1}))^{-2}}{z \! - \! \tilde{b}_{j-1}} \left(\tilde{\alpha}_{0}
(\tilde{b}_{j-1}) \tilde{\boldsymbol{\mathrm{B}}}(\tilde{b}_{j-1}) \! - \! 
\tilde{\alpha}_{1}(\tilde{b}_{j-1}) \tilde{\boldsymbol{\mathrm{A}}}
(\tilde{b}_{j-1}) \right) \! + \! \dfrac{(\tilde{\alpha}_{0}(\tilde{a}_{j}))^{-
1}}{(z \! - \! \tilde{a}_{j})^{2}} \tilde{\boldsymbol{\mathrm{A}}}(\tilde{a}_{j}) 
\right. \nonumber \\
+&\left. \, \dfrac{(\tilde{\alpha}_{0}(\tilde{a}_{j}))^{-2}}{z \! - \! \tilde{a}_{j}} 
\left(\tilde{\alpha}_{0}(\tilde{a}_{j}) \tilde{\boldsymbol{\mathrm{B}}}
(\tilde{a}_{j}) \! - \! \tilde{\alpha}_{1}(\tilde{a}_{j}) \tilde{\boldsymbol{
\mathrm{A}}}(\tilde{a}_{j}) \right) \! + \! \dfrac{((\alpha_{k} \! - \! 
\tilde{b}_{j-1}) \tilde{\alpha}_{1}(\tilde{b}_{j-1}) \! - \! \tilde{\alpha}_{0}
(\tilde{b}_{j-1}))}{(\alpha_{k} \! - \! \tilde{b}_{j-1})^{2}(\tilde{\alpha}_{0}
(\tilde{b}_{j-1}))^{2}} \tilde{\boldsymbol{\mathrm{A}}}(\tilde{b}_{j-1}) 
\right. \nonumber \\
+&\left. \, \dfrac{((\alpha_{k} \! - \! \tilde{a}_{j}) \tilde{\alpha}_{1}
(\tilde{a}_{j}) \! - \! \tilde{\alpha}_{0}(\tilde{a}_{j}))}{(\alpha_{k} 
\! - \! \tilde{a}_{j})^{2}(\tilde{\alpha}_{0}(\tilde{a}_{j}))^{2}} 
\tilde{\boldsymbol{\mathrm{A}}}(\tilde{a}_{j}) \! - \! 
\dfrac{(\tilde{\alpha}_{0}(\tilde{b}_{j-1}))^{-1}}{\alpha_{k} \! - \! 
\tilde{b}_{j-1}} \tilde{\boldsymbol{\mathrm{B}}}(\tilde{b}_{j-1}) \! 
- \! \dfrac{(\tilde{\alpha}_{0}(\tilde{a}_{j}))^{-1}}{\alpha_{k} \! - \! 
\tilde{a}_{j}} \tilde{\boldsymbol{\mathrm{B}}}(\tilde{a}_{j}) \right), 
\label{eqmainfin14}
\end{align}
where, for $j \! = \! 1,2,\dotsc,N \! + \! 1$,
\begin{gather}
\tilde{\boldsymbol{\mathrm{A}}}(\tilde{b}_{j-1}) \! := \! 
\mathscr{E}^{-\sigma_{3}} 
 \mathscr{E}^{\sigma_{3}} 
\me^{\mi ((n-1)K+k) \tilde{\mho}_{j}}, 
\label{eqmainfin18}
\end{gather}
with
\begin{gather}
\tilde{\mathbb{A}}_{11}(\tilde{b}_{j-1}) \! = \! -\tilde{\mathbb{A}}_{22}
(\tilde{b}_{j-1}) \! = \! -s_{1} \tilde{\kappa}_{1}(\tilde{b}_{j-1}) \tilde{
\kappa}_{2}(\tilde{b}_{j-1})(\tilde{\mathfrak{Q}}_{0}(\tilde{b}_{j-1}))^{-1}, 
\label{eqmainfin20} \\
\tilde{\mathbb{A}}_{12}(\tilde{b}_{j-1}) \! = \! -\mi s_{1}(\tilde{\kappa}_{1}
(\tilde{b}_{j-1}))^{2}(\tilde{\mathfrak{Q}}_{0}(\tilde{b}_{j-1}))^{-1}, \qquad 
\tilde{\mathbb{A}}_{21}(\tilde{b}_{j-1}) \! = \! -\mi s_{1}(\tilde{\kappa}_{2}
(\tilde{b}_{j-1}))^{2}(\tilde{\mathfrak{Q}}_{0}(\tilde{b}_{j-1}))^{-1}, 
\label{eqmainfin21}
\end{gather}
\begin{align}
\tilde{\mathbb{B}}_{11}(\tilde{b}_{j-1}) \! = \! -\tilde{\mathbb{B}}_{22}
(\tilde{b}_{j-1}) =& \, \tilde{\kappa}_{1}(\tilde{b}_{j-1}) \tilde{\kappa}_{2}
(\tilde{b}_{j-1}) \left(-s_{1}(\tilde{\mathfrak{Q}}_{0}(\tilde{b}_{j-1}))^{-1} 
\left(\tilde{\daleth}^{1}_{1}(\tilde{b}_{j-1}) \! + \! \tilde{\daleth}^{1}_{-1}
(\tilde{b}_{j-1}) \! - \! \tilde{\mathfrak{Q}}_{1}(\tilde{b}_{j-1})
(\tilde{\mathfrak{Q}}_{0}(\tilde{b}_{j-1}))^{-1} \right) \right. \nonumber \\
-&\left. \, t_{1} \left(\tilde{\mathfrak{Q}}_{0}(\tilde{b}_{j-1}) \! + \! 
(\tilde{\mathfrak{Q}}_{0}(\tilde{b}_{j-1}))^{-1} \tilde{\aleph}^{1}_{1}
(\tilde{b}_{j-1}) \tilde{\aleph}^{1}_{-1}(\tilde{b}_{j-1}) \right) \! + \! 
\mi (s_{1} \! + \! t_{1}) \left(\tilde{\aleph}^{1}_{-1}(\tilde{b}_{j-1}) \! - \! 
\tilde{\aleph}^{1}_{1}(\tilde{b}_{j-1}) \right) \right), \label{eqmainfin22} \\
\tilde{\mathbb{B}}_{12}(\tilde{b}_{j-1}) =& \, (\tilde{\kappa}_{1}
(\tilde{b}_{j-1}))^{2} \left(-\mi s_{1}(\tilde{\mathfrak{Q}}_{0}(\tilde{b}_{j
-1}))^{-1} \left(2 \tilde{\daleth}^{1}_{1}(\tilde{b}_{j-1}) \! - \! \tilde{
\mathfrak{Q}}_{1}(\tilde{b}_{j-1})(\tilde{\mathfrak{Q}}_{0}(\tilde{b}_{j-
1}))^{-1} \right) \right. \nonumber \\
+&\left. \, \mi t_{1} \left(\tilde{\mathfrak{Q}}_{0}(\tilde{b}_{j-1}) \! - \! 
(\tilde{\mathfrak{Q}}_{0}(\tilde{b}_{j-1}))^{-1}(\tilde{\aleph}^{1}_{1}
(\tilde{b}_{j-1}))^{2} \right) \! + \! 2(s_{1} \! - \! t_{1}) \tilde{\aleph}^{1}_{1}
(\tilde{b}_{j-1}) \right), \label{eqmainfin23} \\
\tilde{\mathbb{B}}_{21}(\tilde{b}_{j-1}) =& \, (\tilde{\kappa}_{2}
(\tilde{b}_{j-1}))^{2} \left(-\mi s_{1}(\tilde{\mathfrak{Q}}_{0}(\tilde{b}_{j
-1}))^{-1} \left(2 \tilde{\daleth}^{1}_{-1}(\tilde{b}_{j-1}) \! - \! \tilde{
\mathfrak{Q}}_{1}(\tilde{b}_{j-1})(\tilde{\mathfrak{Q}}_{0}(\tilde{b}_{j-
1}))^{-1} \right) \right. \nonumber \\
+&\left. \, \mi t_{1} \left(\tilde{\mathfrak{Q}}_{0}(\tilde{b}_{j-1}) \! - \! 
(\tilde{\mathfrak{Q}}_{0}(\tilde{b}_{j-1}))^{-1}(\tilde{\aleph}^{1}_{-1}
(\tilde{b}_{j-1}))^{2} \right) \! - \! 2(s_{1} \! - \! t_{1}) \tilde{\aleph}^{
1}_{-1}(\tilde{b}_{j-1}) \right), \label{eqmainfin24}
\end{align}
\begin{gather}
\tilde{\mathbb{A}}_{11}(\tilde{a}_{j}) \! = \! -\tilde{\mathbb{A}}_{22}
(\tilde{a}_{j}) \! = \! -s_{1} \tilde{\kappa}_{1}(\tilde{a}_{j}) \tilde{\kappa}_{2}
(\tilde{a}_{j}) \tilde{\mathfrak{Q}}_{0}(\tilde{a}_{j}), \label{eqmainfin25} \\
\tilde{\mathbb{A}}_{12}(\tilde{a}_{j}) \! = \! \mi s_{1}(\tilde{\kappa}_{1}
(\tilde{a}_{j}))^{2} \tilde{\mathfrak{Q}}_{0}(\tilde{a}_{j}), \qquad 
\tilde{\mathbb{A}}_{21}(\tilde{a}_{j}) \! = \! \mi s_{1}(\tilde{\kappa}_{2}
(\tilde{a}_{j}))^{2} \tilde{\mathfrak{Q}}_{0}(\tilde{a}_{j}), \label{eqmainfin26}
\end{gather}
\begin{align}
\tilde{\mathbb{B}}_{11}(\tilde{a}_{j}) \! = \! -\tilde{\mathbb{B}}_{22}
(\tilde{a}_{j}) =& \, \tilde{\kappa}_{1}(\tilde{a}_{j}) \tilde{\kappa}_{2}
(\tilde{a}_{j}) \left(-s_{1} \left(\tilde{\mathfrak{Q}}_{1}(\tilde{a}_{j}) \! + \! 
\tilde{\mathfrak{Q}}_{0}(\tilde{a}_{j}) \left(\tilde{\daleth}^{1}_{1}(\tilde{a}_{j}) 
\! + \! \tilde{\daleth}^{1}_{-1}(\tilde{a}_{j}) \right) \right) \right. \nonumber \\
-&\left. \, t_{1} \left((\tilde{\mathfrak{Q}}_{0}(\tilde{a}_{j}))^{-1} \! + \! 
\tilde{\mathfrak{Q}}_{0}(\tilde{a}_{j}) \tilde{\aleph}^{1}_{1}(\tilde{a}_{j}) 
\tilde{\aleph}^{1}_{-1}(\tilde{a}_{j}) \right) \! + \! \mi (s_{1} \! + \! t_{1}) 
\left(\tilde{\aleph}^{1}_{-1}(\tilde{a}_{j}) \! - \! \tilde{\aleph}^{1}_{1}
(\tilde{a}_{j}) \right) \right), \label{eqmainfin27} \\
\tilde{\mathbb{B}}_{12}(\tilde{a}_{j}) =& \, (\tilde{\kappa}_{1}(\tilde{a}_{j}))^{2} 
\left(\mi s_{1} \left(\tilde{\mathfrak{Q}}_{1}(\tilde{a}_{j}) \! + \! 2 \tilde{
\mathfrak{Q}}_{0}(\tilde{a}_{j}) \tilde{\daleth}^{1}_{1}(\tilde{a}_{j}) \right) 
\right. \nonumber \\
+&\left. \, \mi t_{1} \left(\tilde{\mathfrak{Q}}_{0}(\tilde{a}_{j})(\tilde{
\aleph}^{1}_{1}(\tilde{a}_{j}))^{2} \! - \! (\tilde{\mathfrak{Q}}_{0}
(\tilde{a}_{j}))^{-1} \right) \! - \! 2(s_{1} \! - \! t_{1}) \tilde{\aleph}^{1}_{1}
(\tilde{a}_{j}) \right), \label{eqmainfin28} \\
\tilde{\mathbb{B}}_{21}(\tilde{a}_{j}) =& \, (\tilde{\kappa}_{2}(\tilde{a}_{j}))^{2} 
\left(\mi s_{1} \left(\tilde{\mathfrak{Q}}_{1}(\tilde{a}_{j}) \! + \! 2 \tilde{
\mathfrak{Q}}_{0}(\tilde{a}_{j}) \tilde{\daleth}^{1}_{-1}(\tilde{a}_{j}) \right) 
\right. \nonumber \\
+&\left. \, \mi t_{1} \left(\tilde{\mathfrak{Q}}_{0}(\tilde{a}_{j})(\tilde{
\aleph}^{1}_{-1}(\tilde{a}_{j}))^{2} \! - \! (\tilde{\mathfrak{Q}}_{0}
(\tilde{a}_{j}))^{-1} \right) \! + \! 2(s_{1} \! - \! t_{1}) \tilde{\aleph}^{1}_{-1}
(\tilde{a}_{j}) \right), \label{eqmainfin29}
\end{align}
where, for $\varepsilon_{1},\varepsilon_{2} \! = \! \pm 1$,
\begin{equation} \label{eqmainfin30} 
s_{1} \! = \! 5/72, \qquad \qquad \quad t_{1} \! = \! -7/72, \qquad 
\qquad \quad \tilde{\mho}_{i} \! = \!
\begin{cases}
\tilde{\Omega}_{i}, &\text{$i \! = \! 1,2,\dotsc,N$,} \\
0, &\text{$i \! = \! 0,N \! + \! 1$,}
\end{cases}
\end{equation}
\begin{align}
\tilde{\kappa}_{1}(\varsigma) =& \, \dfrac{\tilde{\boldsymbol{\theta}}
(\tilde{\boldsymbol{u}}_{+}(\varsigma) \! - \! \frac{1}{2 \pi}((n \! - \! 1)K \! 
+ \! k) \tilde{\boldsymbol{\Omega}} \! + \! \tilde{\boldsymbol{d}})}{\tilde{
\boldsymbol{\theta}}(\tilde{\boldsymbol{u}}_{+}(\varsigma) \! + \! \tilde{
\boldsymbol{d}})}, \qquad \qquad \tilde{\kappa}_{2}(\varsigma) = \dfrac{
\tilde{\boldsymbol{\theta}}(\tilde{\boldsymbol{u}}_{+}(\varsigma) \! - \! 
\frac{1}{2 \pi}((n \! - \! 1)K \! + \! k) \tilde{\boldsymbol{\Omega}} \! - \! 
\tilde{\boldsymbol{d}})}{\tilde{\boldsymbol{\theta}}(\tilde{\boldsymbol{
u}}_{+}(\varsigma) \! - \! \tilde{\boldsymbol{d}})}, \label{eqmainfin31} \\
\tilde{\aleph}^{\varepsilon_{1}}_{\varepsilon_{2}}(\varsigma) =& \, 
-\dfrac{\tilde{\mathfrak{u}}(\varepsilon_{1},\varepsilon_{2},\bm{0};
\varsigma)}{\tilde{\boldsymbol{\theta}}(\varepsilon_{1} \tilde{
\boldsymbol{u}}_{+}(\varsigma) \! + \! \varepsilon_{2} \tilde{
\boldsymbol{d}})} \! + \! \dfrac{\tilde{\mathfrak{u}}(\varepsilon_{1},
\varepsilon_{2},\tilde{\boldsymbol{\Omega}};\varsigma)}{\tilde{
\boldsymbol{\theta}}(\varepsilon_{1} \tilde{\boldsymbol{u}}_{+}(\varsigma) 
\! - \! \frac{1}{2 \pi}((n \! - \! 1)K \! + \! k) \tilde{\boldsymbol{\Omega}} 
\! + \! \varepsilon_{2} \tilde{\boldsymbol{d}})}, \label{eqmainfin32} \\
\tilde{\daleth}^{\varepsilon_{1}}_{\varepsilon_{2}}(\varsigma) =& \, 
-\dfrac{\tilde{\mathfrak{v}}(\varepsilon_{1},\varepsilon_{2},\bm{0};
\varsigma)}{\tilde{\boldsymbol{\theta}}(\varepsilon_{1} 
\tilde{\boldsymbol{u}}_{+}(\varsigma) \! + \! \varepsilon_{2} \tilde{
\boldsymbol{d}})} \! + \! \dfrac{\tilde{\mathfrak{v}}(\varepsilon_{1},
\varepsilon_{2},\tilde{\boldsymbol{\Omega}};\varsigma)}{\tilde{
\boldsymbol{\theta}}(\varepsilon_{1} \tilde{\boldsymbol{u}}_{+}(\varsigma) 
\! - \! \frac{1}{2 \pi}((n \! - \! 1)K \! + \! k) \tilde{\boldsymbol{\Omega}} 
\! + \! \varepsilon_{2} \tilde{\boldsymbol{d}})} \! - \! \left(\dfrac{\tilde{
\mathfrak{u}}(\varepsilon_{1},\varepsilon_{2},\bm{0};\varsigma)}{\tilde{
\boldsymbol{\theta}}(\varepsilon_{1} \tilde{\boldsymbol{u}}_{+}(\varsigma) 
\! + \! \varepsilon_{2} \tilde{\boldsymbol{d}})} \right)^{2} \nonumber \\
+& \, \dfrac{\tilde{\mathfrak{u}}(\varepsilon_{1},\varepsilon_{2},\bm{0};
\varsigma) \tilde{\mathfrak{u}}(\varepsilon_{1},\varepsilon_{2},\tilde{
\boldsymbol{\Omega}};\varsigma)}{\tilde{\boldsymbol{\theta}}
(\varepsilon_{1} \tilde{\boldsymbol{u}}_{+}(\varsigma) \! + \! \varepsilon_{2} 
\tilde{\boldsymbol{d}}) \tilde{\boldsymbol{\theta}}(\varepsilon_{1} \tilde{
\boldsymbol{u}}_{+}(\varsigma) \! - \! \frac{1}{2 \pi}((n \! - \! 1)K \! + \! k) 
\tilde{\boldsymbol{\Omega}} \! + \! \varepsilon_{2} \tilde{\boldsymbol{d}})}, 
\label{eqmainfin33}
\end{align}
with
\begin{gather}
\tilde{\mathfrak{u}}(\varepsilon_{1},\varepsilon_{2},\tilde{\boldsymbol{\Omega}};
\varsigma) \! := \! 2 \pi \tilde{\Lambda}^{\raise-1.0ex\hbox{$\scriptstyle 1$}}_{0}
(\varepsilon_{1},\varepsilon_{2},\tilde{\boldsymbol{\Omega}};\varsigma), \qquad 
\qquad \tilde{\mathfrak{v}}(\varepsilon_{1},\varepsilon_{2},\tilde{\boldsymbol{
\Omega}};\varsigma) \! := \! -2 \pi^{2} 
\tilde{\Lambda}^{\raise-1.0ex\hbox{$\scriptstyle 2$}}_{0}(\varepsilon_{1},
\varepsilon_{2},\tilde{\boldsymbol{\Omega}};\varsigma), \label{eqmainfin34} \\
\tilde{\Lambda}^{\raise-1.0ex\hbox{$\scriptstyle j_{1}$}}_{0}(\varepsilon_{1},
\varepsilon_{2},\tilde{\boldsymbol{\Omega}};\varsigma) \! = \! \sum_{m \in
\mathbb{Z}^{N}}(\tilde{\mathfrak{r}}_{1}(\varsigma))^{j_{1}} \me^{2 \pi \mi 
(m,\varepsilon_{1} \tilde{\boldsymbol{u}}_{+}(\varsigma)-\frac{1}{2 \pi}((n-1)K
+k) \tilde{\boldsymbol{\Omega}}+ \varepsilon_{2} \tilde{\boldsymbol{d}})+ \mi 
\pi (m,\tilde{\boldsymbol{\tau}}m)}, \quad j_{1} \! = \! 1,2, \label{eqmainfin35} \\
\tilde{\mathfrak{r}}_{1}(\varsigma) \! := \! \dfrac{2}{\tilde{\leftthreetimes}
(\varsigma)} \sum_{i=1}^{N} \sum_{j=1}^{N}m_{i} \tilde{c}_{ij} \varsigma^{N-j}, 
\label{eqmainfin36}
\end{gather}
where $\tilde{c}_{i_{1}i_{2}}$, $i_{1},i_{2} \! = \! 1,2,\dotsc,N$, are obtained 
{}from Equations~\eqref{O1} and~\eqref{O2},
\begin{equation}
\tilde{\leftthreetimes}(\tilde{b}_{0}) \! = \! \mi (-1)^{N} \tilde{\eta}_{\tilde{b}_{0}}, 
\quad \tilde{\leftthreetimes}(\tilde{b}_{j}) \! = \! \mi (-1)^{N-j} \tilde{\eta}_{
\tilde{b}_{j}}, \quad \tilde{\leftthreetimes}(\tilde{a}_{N+1}) \! = \! \tilde{\eta}_{
\tilde{a}_{N+1}}, \quad \tilde{\leftthreetimes}(\tilde{a}_{j}) \! = \! (-1)^{N+1-j} 
\tilde{\eta}_{\tilde{a}_{j}}, \quad j \! = \! 1,2,\dotsc,N, \label{eqmainfin37}
\end{equation}
with {}\footnote{All square roots are positive.}
\begin{gather}
\tilde{\eta}_{\tilde{b}_{0}} \! := \! (\tilde{a}_{N+1} \! - \! \tilde{b}_{0})^{1/2} 
\prod_{m=1}^{N}(\tilde{b}_{m} \! - \! \tilde{b}_{0})^{1/2}(\tilde{a}_{m} 
\! - \! \tilde{b}_{0})^{1/2}, \label{eqmainfin57} \\
\tilde{\eta}_{\tilde{b}_{j}} \! := \! (\tilde{b}_{j} \! - \! \tilde{a}_{j})^{1/2}
(\tilde{a}_{N+1} \! - \! \tilde{b}_{j})^{1/2}(\tilde{b}_{j} \! - \! \tilde{b}_{0})^{1/2} 
\prod_{m=1}^{j-1}(\tilde{b}_{j} \! - \! \tilde{b}_{m})^{1/2}(\tilde{b}_{j} \! - \! 
\tilde{a}_{m})^{1/2} \prod_{m^{\prime}=j+1}^{N}(\tilde{b}_{m^{\prime}} \! 
- \! \tilde{b}_{j})^{1/2}(\tilde{a}_{m^{\prime}} \! - \! \tilde{b}_{j})^{1/2}, 
\label{eqmainfin58} \\
\tilde{\eta}_{\tilde{a}_{N+1}} \! := \! (\tilde{a}_{N+1} \! - \! \tilde{b}_{0})^{1/2} 
\prod_{m=1}^{N}(\tilde{a}_{N+1} \! - \! \tilde{b}_{m})^{1/2}(\tilde{a}_{N+1} \! 
- \! \tilde{a}_{m})^{1/2}, \label{eqmainfin59} \\
\tilde{\eta}_{\tilde{a}_{j}} \! := \! (\tilde{b}_{j} \! - \! \tilde{a}_{j})^{1/2}
(\tilde{a}_{N+1} \! - \! \tilde{a}_{j})^{1/2}(\tilde{a}_{j} \! - \! \tilde{b}_{0})^{1/2} 
\prod_{m=1}^{j-1}(\tilde{a}_{j} \! - \! \tilde{b}_{m})^{1/2}(\tilde{a}_{j} \! - \! 
\tilde{a}_{m})^{1/2} \prod_{m^{\prime}=j+1}^{N}(\tilde{b}_{m^{\prime}} \! 
- \! \tilde{a}_{j})^{1/2}(\tilde{a}_{m^{\prime}} \! - \! \tilde{a}_{j})^{1/2}, 
\label{eqmainfin60}
\end{gather}
and, for $j \! = \! 1,2,\dotsc,N$,
\begin{gather}
\tilde{\mathfrak{Q}}_{0}(\tilde{b}_{0}) = -\mi (\tilde{a}_{N+1} \! - \! 
\tilde{b}_{0})^{-1/2} \prod_{m=1}^{N} \dfrac{(\tilde{b}_{m} \! - \! \tilde{b}_{
0})^{1/2}}{(\tilde{a}_{m} \! - \! \tilde{b}_{0})^{1/2}}, \label{eqmainfin41} \\
\tilde{\mathfrak{Q}}_{1}(\tilde{b}_{0}) = \dfrac{1}{2} \tilde{\mathfrak{
Q}}_{0}(\tilde{b}_{0}) \left(\sum_{m=1}^{N} \left(\dfrac{1}{\tilde{b}_{0} \! - 
\! \tilde{b}_{m}} \! - \! \dfrac{1}{\tilde{b}_{0} \! - \! \tilde{a}_{m}} \right) \! 
- \! \dfrac{1}{\tilde{b}_{0} \! - \! \tilde{a}_{N+1}} \right), \label{eqmainfin42} \\
\tilde{\mathfrak{Q}}_{0}(\tilde{b}_{j})= -\dfrac{\mi (\tilde{b}_{j} \! - \! 
\tilde{b}_{0})^{1/2}}{(\tilde{a}_{N+1} \! - \! \tilde{b}_{j})^{1/2}(\tilde{b}_{j} 
\! - \! \tilde{a}_{j})^{1/2}} \prod_{m=1}^{j-1} \dfrac{(\tilde{b}_{j} \! 
- \! \tilde{b}_{m})^{1/2}}{(\tilde{b}_{j} \! - \! \tilde{a}_{m})^{1/2}} 
\prod_{m^{\prime}=j+1}^{N} \dfrac{(\tilde{b}_{m^{\prime}} \! - \! 
\tilde{b}_{j})^{1/2}}{(\tilde{a}_{m^{\prime}} \! - \! \tilde{b}_{j})^{1/2}}, 
\label{eqmainfin43} \\
\tilde{\mathfrak{Q}}_{1}(\tilde{b}_{j}) = \dfrac{1}{2} \tilde{\mathfrak{
Q}}_{0}(\tilde{b}_{j}) \left(\sum_{\substack{m=1\\m \neq j}}^{N} \left(
\dfrac{1}{\tilde{b}_{j} \! - \! \tilde{b}_{m}} \! - \! \dfrac{1}{\tilde{b}_{j} \! - 
\! \tilde{a}_{m}} \right) \! + \! \dfrac{1}{\tilde{b}_{j} \! - \! \tilde{b}_{0}} \! 
- \! \dfrac{1}{\tilde{b}_{j} \! - \! \tilde{a}_{N+1}} \! - \! \dfrac{1}{\tilde{b}_{j} 
\! - \! \tilde{a}_{j}} \right), \label{eqmainfin44} \\
\tilde{\mathfrak{Q}}_{0}(\tilde{a}_{N+1}) = (\tilde{a}_{N+1} \! - \! 
\tilde{b}_{0})^{1/2} \prod_{m=1}^{N} \dfrac{(\tilde{a}_{N+1} \! - \! \tilde{b}_{
m})^{1/2}}{(\tilde{a}_{N+1} \! - \! \tilde{a}_{m})^{1/2}}, \label{eqmainfin45} \\
\tilde{\mathfrak{Q}}_{1}(\tilde{a}_{N+1}) = \dfrac{1}{2} \tilde{\mathfrak{Q}}_{0}
(\tilde{a}_{N+1}) \left(\sum_{m=1}^{N} \left(\dfrac{1}{\tilde{a}_{N+1} \! - \! 
\tilde{b}_{m}} \! - \! \dfrac{1}{\tilde{a}_{N+1} \! - \! \tilde{a}_{m}} \right) \! + 
\! \dfrac{1}{\tilde{a}_{N+1} \! - \! \tilde{b}_{0}} \right), \label{eqmainfin46} \\
\tilde{\mathfrak{Q}}_{0}(\tilde{a}_{j}) = \dfrac{(\tilde{a}_{j} \! - \! 
\tilde{b}_{0})^{1/2}(\tilde{b}_{j} \! - \! \tilde{a}_{j})^{1/2}}{(\tilde{a}_{N+1} 
\! - \! \tilde{a}_{j})^{1/2}} \prod_{m=1}^{j-1} \dfrac{(\tilde{a}_{j} \! 
- \! \tilde{b}_{m})^{1/2}}{(\tilde{a}_{j} \! - \! \tilde{a}_{m})^{1/2}} 
\prod_{m^{\prime}=j+1}^{N} \dfrac{(\tilde{b}_{m^{\prime}} \! - \! 
\tilde{a}_{j})^{1/2}}{(\tilde{a}_{m^{\prime}} \! - \! \tilde{a}_{j})^{1/2}}, 
\label{eqmainfin47} \\
\tilde{\mathfrak{Q}}_{1}(\tilde{a}_{j}) = \dfrac{1}{2} \tilde{\mathfrak{Q}}_{0}
(\tilde{a}_{j}) \left(\sum_{\substack{m=1\\m \neq j}}^{N} \left(\dfrac{1}{
\tilde{a}_{j} \! - \! \tilde{b}_{m}} \! - \! \dfrac{1}{\tilde{a}_{j} \! - \! \tilde{a}_{m}} 
\right) \! + \! \dfrac{1}{\tilde{a}_{j} \! - \! \tilde{b}_{0}} \! - \! \dfrac{1}{\tilde{a}_{j} 
\! - \! \tilde{a}_{N+1}} \! + \! \dfrac{1}{\tilde{a}_{j} \! - \! \tilde{b}_{j}} \right), 
\label{eqmainfin48} \\
\tilde{\alpha}_{0}(\tilde{b}_{0}) = \mi (-1)^{N} \dfrac{2}{3} \tilde{h}_{\widetilde{V}}
(\tilde{b}_{0}) \tilde{\eta}_{\tilde{b}_{0}}, \label{eqmainfin49} \\
\tilde{\alpha}_{0}(\tilde{b}_{j}) = \mi (-1)^{N-j} \dfrac{2}{3} \tilde{h}_{\widetilde{V}}
(\tilde{b}_{j}) \tilde{\eta}_{\tilde{b}_{j}}, \label{eqmainfin50} \\
\tilde{\alpha}_{1}(\tilde{b}_{0}) = \mi (-1)^{N} \tilde{\eta}_{\tilde{b}_{0}} \dfrac{2}{5} 
\left(\dfrac{1}{2} \tilde{h}_{\widetilde{V}}(\tilde{b}_{0}) \left(\sum_{m=1}^{N} 
\left(\dfrac{1}{\tilde{b}_{0} \! - \! \tilde{b}_{m}} \! + \! \dfrac{1}{\tilde{b}_{0} \! 
- \! \tilde{a}_{m}} \right) \! + \! \dfrac{1}{\tilde{b}_{0} \! - \! \tilde{a}_{N+1}} 
\right) \! + \! (\tilde{h}_{\widetilde{V}}(\tilde{b}_{0}))^{\prime} \right), 
\label{eqmainfin51} \\
\tilde{\alpha}_{1}(\tilde{b}_{j}) = \mi (-1)^{N-j} \tilde{\eta}_{\tilde{b}_{j}} \dfrac{2}{5} 
\left(\dfrac{1}{2} \tilde{h}_{\widetilde{V}}(\tilde{b}_{j}) \left(\sum_{\substack{m=
1\\m \neq j}}^{N} \left(\dfrac{1}{\tilde{b}_{j} \! - \! \tilde{b}_{m}} \! + \! \dfrac{
1}{\tilde{b}_{j} \! - \! \tilde{a}_{m}} \right) \! + \! \dfrac{1}{\tilde{b}_{j} \! - \! 
\tilde{a}_{j}} \! + \! \dfrac{1}{\tilde{b}_{j} \! - \! \tilde{a}_{N+1}} \! + \! \dfrac{
1}{\tilde{b}_{j} \! - \! \tilde{b}_{0}} \right) \! + \! (\tilde{h}_{\widetilde{V}}
(\tilde{b}_{j}))^{\prime} \right), \label{eqmainfin52} \\
\tilde{\alpha}_{0}(\tilde{a}_{N+1}) = \dfrac{2}{3} \tilde{h}_{\widetilde{V}}
(\tilde{a}_{N+1}) \tilde{\eta}_{\tilde{a}_{N+1}}, \label{eqmainfin53} \\
\tilde{\alpha}_{0}(\tilde{a}_{j}) = (-1)^{N+1-j} \dfrac{2}{3} \tilde{h}_{\widetilde{V}}
(\tilde{a}_{j}) \tilde{\eta}_{\tilde{a}_{j}}, \label{eqmainfin54} \\
\tilde{\alpha}_{1}(\tilde{a}_{N+1}) = \tilde{\eta}_{\tilde{a}_{N+1}} \dfrac{2}{5} 
\left(\dfrac{1}{2} \tilde{h}_{\widetilde{V}}(\tilde{a}_{N+1}) \left(\sum_{m=1}^{N} 
\left(\dfrac{1}{\tilde{a}_{N+1} \! - \! \tilde{b}_{m}} \! + \! \dfrac{1}{\tilde{a}_{N+1} 
\! - \! \tilde{a}_{m}} \right) \! + \! \dfrac{1}{\tilde{a}_{N+1} \! - \! \tilde{b}_{0}} 
\right) \! + \! (\tilde{h}_{\widetilde{V}}(\tilde{a}_{N+1}))^{\prime} \right), 
\label{eqmainfin55} \\
\tilde{\alpha}_{1}(\tilde{a}_{j}) = (-1)^{N+1-j} \tilde{\eta}_{\tilde{a}_{j}} \dfrac{2}{5} 
\left(\dfrac{1}{2} \tilde{h}_{\widetilde{V}}(\tilde{a}_{j}) \left(\sum_{\substack{m
=1\\m \neq j}}^{N} \left(\dfrac{1}{\tilde{a}_{j} \! - \! \tilde{b}_{m}} \! + \! 
\dfrac{1}{\tilde{a}_{j} \! - \! \tilde{a}_{m}} \right) \! + \! \dfrac{1}{\tilde{a}_{j} 
\! - \! \tilde{b}_{j}} \! + \! \dfrac{1}{\tilde{a}_{j} \! - \! \tilde{a}_{N+1}} \! + \! 
\dfrac{1}{\tilde{a}_{j} \! - \! \tilde{b}_{0}} \right) \! + \! (\tilde{h}_{\widetilde{V}}
(\tilde{a}_{j}))^{\prime} \right); \label{eqmainfin56}
\end{gather}
{\rm \pmb{(2)}} for $z \! \in \! \tilde{\Upsilon}_{2}$,
\begin{align}
(z \! - \! \alpha_{k}) \pmb{\pi}_{k}^{n}(z) \underset{\underset{z_{o}
=1+o(1)}{\mathscr{N},n \to \infty}}{=}& \, \left(\left(\tilde{\mathbb{K}}_{11} 
\tilde{\mathbb{M}}_{12}(z) \! + \! \tilde{\mathbb{K}}_{12} \tilde{\mathbb{
M}}_{22}(z) \right) \left(1 \! + \! \dfrac{1}{(n \! - \! 1)K \! + \! k} \tilde{
\mathcal{R}}_{11}^{\sharp}(z) \! + \! \mathcal{O} \left(\dfrac{\mathfrak{c}
(n,k,z_{o})}{((n \! - \! 1)K \! + \! k)^{2}} \right) \right) \right. \nonumber \\
+&\left. \, \mathscr{E}^{-2} \left(\tilde{\mathbb{K}}_{21} \tilde{\mathbb{
M}}_{12}(z) \! + \! \tilde{\mathbb{K}}_{22} \tilde{\mathbb{M}}_{22}(z) 
\right) \left(\dfrac{1}{(n \! - \! 1)K \! + \! k} \tilde{\mathcal{R}}_{12}^{
\sharp}(z) \! + \! \mathcal{O} \left(\dfrac{\mathfrak{c}(n,k,z_{o})}{
((n \! - \! 1)K \! + \! k)^{2}} \right) \right) \right) \nonumber \\
\times& \, \me^{n(g^{f}(z)-\hat{\mathscr{P}}_{0}^{-})}, 
\label{eqmainfin61}
\end{align}
and
\begin{align}
(z \! - \! \alpha_{k}) \int_{\mathbb{R}} \dfrac{((\xi \! - \! \alpha_{k}) 
\pmb{\pi}_{k}^{n}(\xi)) \me^{-n \widetilde{V}(\xi)}}{(\xi \! - \! \alpha_{k})
(\xi \! - \! z)} \, \dfrac{\md \xi}{2 \pi \mi} \underset{\underset{z_{o}=1+o(1)}{
\mathscr{N},n \to \infty}}{=}& \, -\left(\mathscr{E}^{2} \left(\tilde{\mathbb{K}}_{
11} \tilde{\mathbb{M}}_{11}(z) \! + \! \tilde{\mathbb{K}}_{12} \tilde{\mathbb{
M}}_{21}(z) \right) \left(1 \! + \! \dfrac{1}{(n \! - \! 1)K \! + \! k} \tilde{
\mathcal{R}}_{11}^{\sharp}(z) \! + \! \mathcal{O} \left(\dfrac{\mathfrak{c}
(n,k,z_{o})}{((n \! - \! 1)K \! + \! k)^{2}} \right) \right) \right. \nonumber \\
+&\left. \, \left(\tilde{\mathbb{K}}_{21} \tilde{\mathbb{M}}_{11}(z) \! + \! 
\tilde{\mathbb{K}}_{22} \tilde{\mathbb{M}}_{21}(z) \right) \left(\dfrac{1}{(n 
\! - \! 1)K \! + \! k} \tilde{\mathcal{R}}_{12}^{\sharp}(z) \! + \! \mathcal{O} 
\left(\dfrac{\mathfrak{c}(n,k,z_{o})}{((n \! - \! 1)K \! + \! k)^{2}} \right) \right) 
\right) \nonumber \\
\times& \, \me^{n \tilde{\ell}} \me^{-n(g^{f}(z)-\hat{\mathscr{P}}_{0}^{-})}; 
\label{eqmainfin62}
\end{align}
{\rm \pmb{(3)}} for $z \! \in \! \tilde{\Upsilon}_{3}$,
\begin{align}
(z \! - \! \alpha_{k}) \pmb{\pi}_{k}^{n}(z) \underset{\underset{z_{o}
=1+o(1)}{\mathscr{N},n \to \infty}}{=}& \, \left(\left(\tilde{\mathbb{K}}_{11} 
\left(\tilde{\mathbb{M}}_{11}(z) \! + \! \tilde{\mathbb{M}}_{12}(z) \me^{-2 
\pi \mi ((n-1)K+k) \int_{z}^{\tilde{a}_{N+1}} \psi_{\widetilde{V}}^{f}(\xi) \, 
\md \xi} \right) \! + \! \tilde{\mathbb{K}}_{12} \left(\tilde{\mathbb{M}}_{21}
(z) \! + \! \tilde{\mathbb{M}}_{22}(z) \me^{-2 \pi \mi ((n-1)K+k) \int_{z}^{
\tilde{a}_{N+1}} \psi_{\widetilde{V}}^{f}(\xi) \, \md \xi} \right) \right) \right.
\nonumber \\
\times&\left. \, \left(1 \! + \! \dfrac{1}{(n \! - \! 1)K \! + \! k} \tilde{
\mathcal{R}}_{11}^{\sharp}(z) \! + \! \mathcal{O} \left(\dfrac{\mathfrak{c}
(n,k,z_{o})}{((n \! - \! 1)K \! + \! k)^{2}} \right) \right) \! + \! \mathscr{E}^{2} 
\left(\tilde{\mathbb{K}}_{21} \left(\tilde{\mathbb{M}}_{11}(z) \! + \! \tilde{
\mathbb{M}}_{12}(z) \me^{-2 \pi \mi ((n-1)K+k) \int_{z}^{\tilde{a}_{N+1}} 
\psi_{\widetilde{V}}^{f}(\xi) \, \md \xi} \right) \right. \right. \nonumber \\
+&\left. \left. \, \tilde{\mathbb{K}}_{22} \left(\tilde{\mathbb{M}}_{21}(z) 
\! + \! \tilde{\mathbb{M}}_{22}(z) \me^{-2 \pi \mi ((n-1)K+k) \int_{z}^{
\tilde{a}_{N+1}} \psi_{\widetilde{V}}^{f}(\xi) \, \md \xi} \right) \right) 
\left(\dfrac{1}{(n \! - \! 1)K \! + \! k} \tilde{\mathcal{R}}_{12}^{\sharp}(z) 
\! + \! \mathcal{O} \left(\dfrac{\mathfrak{c}(n,k,z_{o})}{((n \! - \! 1)K 
\! + \! k)^{2}} \right) \right) \right) \nonumber \\
\times& \, \me^{n(g^{f}(z)-\hat{\mathscr{P}}_{0}^{+})}, \label{eqmainfin63}
\end{align}
and
\begin{align}
(z \! - \! \alpha_{k}) \int_{\mathbb{R}} \dfrac{((\xi \! - \! \alpha_{k}) 
\pmb{\pi}_{k}^{n}(\xi)) \me^{-n \widetilde{V}(\xi)}}{(\xi \! - \! \alpha_{k})
(\xi \! - \! z)} \, \dfrac{\md \xi}{2 \pi \mi} \underset{\underset{z_{o}=1+o(1)}{
\mathscr{N},n \to \infty}}{=}& \, \left(\mathscr{E}^{-2} \left(\tilde{\mathbb{K}}_{
11} \tilde{\mathbb{M}}_{12}(z) \! + \! \tilde{\mathbb{K}}_{12} \tilde{\mathbb{
M}}_{22}(z) \right) \left(1 \! + \! \dfrac{1}{(n \! - \! 1)K \! + \! k} \tilde{
\mathcal{R}}_{11}^{\sharp}(z) \! + \! \mathcal{O} \left(\dfrac{\mathfrak{c}
(n,k,z_{o})}{((n \! - \! 1)K \! + \! k)^{2}} \right) \right) \right. \nonumber \\
+&\left. \, \left(\tilde{\mathbb{K}}_{21} \tilde{\mathbb{M}}_{12}(z) \! + \! 
\tilde{\mathbb{K}}_{22} \tilde{\mathbb{M}}_{22}(z) \right) \left(\dfrac{1}{(n \! - 
\! 1)K \! + \! k} \tilde{\mathcal{R}}_{12}^{\sharp}(z) \! + \! \mathcal{O} \left(
\dfrac{\mathfrak{c}(n,k,z_{o})}{((n \! - \! 1)K \! + \! k)^{2}} \right) \right) \right) 
\nonumber \\
\times& \, \me^{n \tilde{\ell}} \me^{-n(g^{f}(z)-\hat{\mathscr{P}}_{0}^{+})}; 
\label{eqmainfin64}
\end{align}
{\rm \pmb{(4)}} for $z \! \in \! \tilde{\Upsilon}_{4}$,
\begin{align}
(z \! - \! \alpha_{k}) \pmb{\pi}_{k}^{n}(z) \underset{\underset{z_{o}
=1+o(1)}{\mathscr{N},n \to \infty}}{=}& \, \left(\left(\tilde{\mathbb{K}}_{11} 
\left(\tilde{\mathbb{M}}_{12}(z) \! + \! \tilde{\mathbb{M}}_{11}(z) \me^{2 \pi 
\mi ((n-1)K+k) \int_{z}^{\tilde{a}_{N+1}} \psi_{\widetilde{V}}^{f}(\xi) \, \md 
\xi} \right) \! + \! \tilde{\mathbb{K}}_{12} \left(\tilde{\mathbb{M}}_{22}(z) 
\! + \! \tilde{\mathbb{M}}_{21}(z) \me^{2 \pi \mi ((n-1)K+k) \int_{z}^{
\tilde{a}_{N+1}} \psi_{\widetilde{V}}^{f}(\xi) \, \md \xi} \right) \right) \right.
\nonumber \\
\times&\left. \, \left(1 \! + \! \dfrac{1}{(n \! - \! 1)K \! + \! k} \tilde{
\mathcal{R}}_{11}^{\sharp}(z) \! + \! \mathcal{O} \left(\dfrac{\mathfrak{c}
(n,k,z_{o})}{((n \! - \! 1)K \! + \! k)^{2}} \right) \right) \! + \! \mathscr{E}^{-2} 
\left(\tilde{\mathbb{K}}_{21} \left(\tilde{\mathbb{M}}_{12}(z) \! + \! \tilde{
\mathbb{M}}_{11}(z) \me^{2 \pi \mi ((n-1)K+k) \int_{z}^{\tilde{a}_{N+1}} 
\psi_{\widetilde{V}}^{f}(\xi) \, \md \xi} \right) \right. \right. \nonumber \\
+&\left. \left. \, \tilde{\mathbb{K}}_{22} \left(\tilde{\mathbb{M}}_{22}(z) 
\! + \! \tilde{\mathbb{M}}_{21}(z) \me^{2 \pi \mi ((n-1)K+k) \int_{z}^{
\tilde{a}_{N+1}} \psi_{\widetilde{V}}^{f}(\xi) \, \md \xi} \right) \right) 
\left(\dfrac{1}{(n \! - \! 1)K \! + \! k} \tilde{\mathcal{R}}_{12}^{\sharp}(z) 
\! + \! \mathcal{O} \left(\dfrac{\mathfrak{c}(n,k,z_{o})}{((n \! - \! 1)K 
\! + \! k)^{2}} \right) \right) \right) \nonumber \\
\times& \, \me^{n(g^{f}(z)-\hat{\mathscr{P}}_{0}^{-})}, \label{eqmainfin65}
\end{align}
and
\begin{align}
(z \! - \! \alpha_{k}) \int_{\mathbb{R}} \dfrac{((\xi \! - \! \alpha_{k}) 
\pmb{\pi}_{k}^{n}(\xi)) \me^{-n \widetilde{V}(\xi)}}{(\xi \! - \! \alpha_{k})
(\xi \! - \! z)} \, \dfrac{\md \xi}{2 \pi \mi} \underset{\underset{z_{o}=1+o(1)}{
\mathscr{N},n \to \infty}}{=}& \, -\left(\mathscr{E}^{2} \left(\tilde{\mathbb{K}}_{
11} \tilde{\mathbb{M}}_{11}(z) \! + \! \tilde{\mathbb{K}}_{12} \tilde{\mathbb{
M}}_{21}(z) \right) \left(1 \! + \! \dfrac{1}{(n \! - \! 1)K \! + \! k} \tilde{
\mathcal{R}}_{11}^{\sharp}(z) \! + \! \mathcal{O} \left(\dfrac{\mathfrak{c}
(n,k,z_{o})}{((n \! - \! 1)K \! + \! k)^{2}} \right) \right) \right. \nonumber \\
+&\left. \, \left(\tilde{\mathbb{K}}_{21} \tilde{\mathbb{M}}_{11}(z) \! + \! 
\tilde{\mathbb{K}}_{22} \tilde{\mathbb{M}}_{21}(z) \right) \left(\dfrac{1}{(n 
\! - \! 1)K \! + \! k} \tilde{\mathcal{R}}_{12}^{\sharp}(z) \! + \! \mathcal{O} 
\left(\dfrac{\mathfrak{c}(n,k,z_{o})}{((n \! - \! 1)K \! + \! k)^{2}} \right) 
\right) \right) \nonumber \\
\times& \, \me^{n \tilde{\ell}} \me^{-n(g^{f}(z)-\hat{\mathscr{P}}_{0}^{-})}; 
\label{eqmainfin66}
\end{align}
{\rm \pmb{(5)}} for $z \! \in \! \tilde{\Omega}^{1}_{\tilde{b}_{j-1}}$, $j \! = 
\! 1,2,\dotsc,N \! + \! 1$,
\begin{align}
(z \! - \! \alpha_{k}) \pmb{\pi}_{k}^{n}(z) \underset{\underset{z_{o}=1+
o(1)}{\mathscr{N},n \to \infty}}{=}& \, \left(\tilde{\mathcal{X}}_{11}^{\tilde{b},
1}(z) \left(1 \! + \! \dfrac{1}{(n \! - \! 1)K \! + \! k} \left(\tilde{\mathcal{R}}_{
11}^{\sharp}(z) \! - \! \tilde{\mathcal{R}}_{11}^{\natural}(z) \right) \! + \! 
\mathcal{O} \left(\dfrac{\mathfrak{c}(n,k,z_{o})}{((n \! - \! 1)K \! + \! k)^{2}} 
\right) \right) \right. \nonumber \\
+&\left. \, \tilde{\mathcal{X}}_{21}^{\tilde{b},1}(z) \left(\dfrac{1}{(n \! - \! 1)K \! 
+ \! k} \left(\tilde{\mathcal{R}}_{12}^{\sharp}(z) \! - \! \tilde{\mathcal{R}}_{12}^{
\natural}(z) \right) \! + \! \mathcal{O} \left(\dfrac{\mathfrak{c}(n,k,z_{o})}{
((n \! - \! 1)K \! + \! k)^{2}} \right) \right) \right) \nonumber \\
\times& \, \mathscr{E} \me^{n(g^{f}(z)-\hat{\mathscr{P}}_{0}^{+})}, 
\label{eqmainfin67}
\end{align}
and
\begin{align}
(z \! - \! \alpha_{k}) \int_{\mathbb{R}} \dfrac{((\xi \! - \! \alpha_{k}) 
\pmb{\pi}_{k}^{n}(\xi)) \me^{-n \widetilde{V}(\xi)}}{(\xi \! - \! \alpha_{k})
(\xi \! - \! z)} \, \dfrac{\md \xi}{2 \pi \mi} \underset{\underset{z_{o}=1+o(1)}{
\mathscr{N},n \to \infty}}{=}& \, \left(\tilde{\mathcal{X}}_{12}^{\tilde{b},1}(z) 
\left(1 \! + \! \dfrac{1}{(n \! - \! 1)K \! + \! k} \left(\tilde{\mathcal{R}}_{11}^{
\sharp}(z) \! - \! \tilde{\mathcal{R}}_{11}^{\natural}(z) \right) \! + \! \mathcal{O} 
\left(\dfrac{\mathfrak{c}(n,k,z_{o})}{((n \! - \! 1)K \! + \! k)^{2}} \right) \right) 
\right. \nonumber \\
+&\left. \, \tilde{\mathcal{X}}_{22}^{\tilde{b},1}(z) \left(\dfrac{1}{(n \! - \! 1)K \! 
+ \! k} \left(\tilde{\mathcal{R}}_{12}^{\sharp}(z) \! - \! \tilde{\mathcal{R}}_{12}^{
\natural}(z) \right) \! + \! \mathcal{O} \left(\dfrac{\mathfrak{c}(n,k,z_{o})}{
((n \! - \! 1)K \! + \! k)^{2}} \right) \right) \right) \nonumber \\
\times& \, \mathscr{E}^{-1} \me^{n \tilde{\ell}} \me^{-n(g^{f}(z)-
\hat{\mathscr{P}}_{0}^{+})}, \label{eqmainfin68}
\end{align}
where
\begin{equation} \label{eqmainfin69} 
\tilde{\mathcal{R}}^{\natural}(z) \! = \! \sum_{j=1}^{N+1} \left(\tilde{\mathbb{
Y}}_{\tilde{b}_{j-1}}(z) \chi_{\tilde{\mathbb{U}}_{\tilde{\delta}_{\tilde{b}_{j-1}}}}
(z) \! + \! \tilde{\mathbb{Y}}_{\tilde{a}_{j}}(z) \chi_{\tilde{\mathbb{U}}_{\tilde{
\delta}_{\hat{a}_{j}}}}(z) \right),
\end{equation}
with $\chi_{\tilde{\mathbb{U}}_{\tilde{\delta}_{\tilde{b}_{j-1}}}}(z)$ (resp., 
$\chi_{\tilde{\mathbb{U}}_{\tilde{\delta}_{\tilde{a}_{j}}}}(z))$ the characteristic 
function of the set $\tilde{\mathbb{U}}_{\tilde{\delta}_{\tilde{b}_{j-1}}}$ (resp., 
$\tilde{\mathbb{U}}_{\tilde{\delta}_{\tilde{a}_{j}}})$,
\begin{gather}
\tilde{\mathbb{Y}}_{\tilde{b}_{j-1}}(z) := \dfrac{1}{\tilde{\xi}_{\tilde{b}_{j-1}}(z)} 
\mathscr{E}^{-\sigma_{3}} \tilde{\mathbb{K}} \tilde{\mathbb{M}}(z) 
\begin{pmatrix}
-(s_{1} \! + \! t_{1}) & -\mi (s_{1} \! - \! t_{1}) \me^{\mi ((n-1)K+k) 
\tilde{\mho}_{j-1}} \\
-\mi (s_{1} \! - \! t_{1}) \me^{-\mi ((n-1)K+k) \tilde{\mho}_{j-1}} & 
(s_{1} \! + \! t_{1})
\end{pmatrix}
(\mathscr{E}^{-\sigma_{3}} \tilde{\mathbb{K}} \tilde{\mathbb{M}}(z))^{-1}, 
\label{eqmainfin70} \\
\tilde{\mathbb{Y}}_{\tilde{a}_{j}}(z) := \dfrac{1}{\tilde{\xi}_{\tilde{a}_{j}}(z)} 
\mathscr{E}^{-\sigma_{3}} \tilde{\mathbb{K}} \tilde{\mathbb{M}}(z) 
\begin{pmatrix}
-(s_{1} \! + \! t_{1}) & \mi (s_{1} \! - \! t_{1}) \me^{\mi ((n-1)K+k) 
\tilde{\mho}_{j}} \\
\mi (s_{1} \! - \! t_{1}) \me^{-\mi ((n-1)K+k) \tilde{\mho}_{j}} & 
(s_{1} \! + \! t_{1})
\end{pmatrix}
(\mathscr{E}^{-\sigma_{3}} \tilde{\mathbb{K}} \tilde{\mathbb{M}}(z))^{-1}, 
\label{eqmainfin71}
\end{gather}
where
\begin{gather}
\tilde{\xi}_{\tilde{b}_{j-1}}(z) \! = \! -\int_{z}^{\tilde{b}_{j-1}}(\tilde{R}
(\xi))^{1/2} \tilde{h}_{\widetilde{V}}(\xi) \, \md \xi, \label{eqmainfin72} \\
\tilde{\xi}_{\tilde{a}_{j}}(z) \! = \! \int_{\tilde{a}_{j}}^{z}(\tilde{R}(\xi))^{1/2} 
\tilde{h}_{\widetilde{V}}(\xi) \, \md \xi, \label{eqmainfin73}
\end{gather}
and
\begin{align}
\tilde{\mathcal{X}}_{11}^{\tilde{b},1}(z) =& \, -\mi \sqrt{\pi} \me^{\frac{1}{2} 
((n-1)K+k) \tilde{\xi}_{\tilde{b}_{j-1}}(z)} \mathscr{E}^{-1} \left(\mi \left(
\tilde{\mathbb{K}}_{11} \tilde{\mathbb{M}}_{11}(z) \! + \! \tilde{\mathbb{
K}}_{12} \tilde{\mathbb{M}}_{21}(z) \right) \Theta_{0,0}^{0,0,-}(\tilde{\Phi}_{
\tilde{b}_{j-1}}(z)) \! - \! \left(\tilde{\mathbb{K}}_{11} \tilde{\mathbb{M}}_{12}
(z) \! + \! \tilde{\mathbb{K}}_{12} \tilde{\mathbb{M}}_{22}(z) \right) \right. 
\nonumber \\
\times&\left. \, \Theta_{0,0}^{0,0,+}(\tilde{\Phi}_{\tilde{b}_{j-1}}(z)) 
\me^{-\mi ((n-1)K+k) \tilde{\mho}_{j-1}} \right), \label{eqmainfin74} \\
\tilde{\mathcal{X}}_{12}^{\tilde{b},1}(z) =& \, \sqrt{\pi} \me^{-\frac{\mi \pi}{6}} 
\me^{-\frac{1}{2}((n-1)K+k) \tilde{\xi}_{\tilde{b}_{j-1}}(z)} \mathscr{E}^{-1} 
\left(\left(\tilde{\mathbb{K}}_{11} \tilde{\mathbb{M}}_{12}(z) \! + \! \tilde{
\mathbb{K}}_{12} \tilde{\mathbb{M}}_{22}(z) \right) \Theta_{2,2}^{0,2,+}
(\tilde{\Phi}_{\tilde{b}_{j-1}}(z)) \! - \! \mi \left(\tilde{\mathbb{K}}_{11} 
\tilde{\mathbb{M}}_{11}(z) \! + \! \tilde{\mathbb{K}}_{12} \tilde{\mathbb{
M}}_{21}(z) \right) \right. \nonumber \\
\times&\left. \, \Theta_{2,2}^{0,2,-}(\tilde{\Phi}_{\tilde{b}_{j-1}}(z)) 
\me^{\mi ((n-1)K+k) \tilde{\mho}_{j-1}} \right), \label{eqmainfin75} \\
\tilde{\mathcal{X}}_{21}^{\tilde{b},1}(z) =& \, -\mi \sqrt{\pi} \me^{\frac{1}{2} 
((n-1)K+k) \tilde{\xi}_{\hat{b}_{j-1}}(z)} \mathscr{E} \left(\mi \left(
\tilde{\mathbb{K}}_{21} \tilde{\mathbb{M}}_{11}(z) \! + \! \tilde{\mathbb{
K}}_{22} \tilde{\mathbb{M}}_{21}(z) \right) \Theta_{0,0}^{0,0,-}(\tilde{\Phi}_{
\tilde{b}_{j-1}}(z)) \! - \! \left(\tilde{\mathbb{K}}_{21} \tilde{\mathbb{M}}_{12}
(z) \! + \! \tilde{\mathbb{K}}_{22} \tilde{\mathbb{M}}_{22}(z) \right) \right. 
\nonumber \\
\times&\left. \, \Theta_{0,0}^{0,0,+}(\tilde{\Phi}_{\tilde{b}_{j-1}}(z)) 
\me^{-\mi ((n-1)K+k) \tilde{\mho}_{j-1}} \right), \label{eqmainfin76} \\
\tilde{\mathcal{X}}_{22}^{\tilde{b},1}(z) =& \, \sqrt{\pi} \me^{-\frac{\mi \pi}{6}} 
\me^{-\frac{1}{2}((n-1)K+k) \tilde{\xi}_{\hat{b}_{j-1}}(z)} \mathscr{E} \left(
\left(\tilde{\mathbb{K}}_{21} \tilde{\mathbb{M}}_{12}(z) \! + \! \tilde{
\mathbb{K}}_{22} \tilde{\mathbb{M}}_{22}(z) \right) \Theta_{2,2}^{0,2,+}
(\tilde{\Phi}_{\tilde{b}_{j-1}}(z)) \! - \! \mi \left(\tilde{\mathbb{K}}_{21} 
\tilde{\mathbb{M}}_{11}(z) \! + \! \tilde{\mathbb{K}}_{22} \tilde{\mathbb{
M}}_{21}(z) \right) \right. \nonumber \\
\times&\left. \, \Theta_{2,2}^{0,2,-}(\tilde{\Phi}_{\tilde{b}_{j-1}}(z)) 
\me^{\mi ((n-1)K+k) \tilde{\mho}_{j-1}} \right), \label{eqmainfin77}
\end{align}
where $\Theta_{r_{3},r_{4}}^{r_{1},r_{2},\pm}(\zeta)$, $r_{1},r_{2},r_{3},r_{4} 
\! = \! 0,1,2$, is defined by Equation~\eqref{eqmaininf76}, $\omega \! = 
\! \exp (2 \pi \mi/3)$, and $\tilde{\Phi}_{\tilde{b}_{j-1}}(z)$ is defined by  
Equation~\eqref{eqmainfin78}$;$\\
{\rm \pmb{(6)}} for $z \! \in \! \tilde{\Omega}_{\tilde{a}_{j}}^{1}$, $j \! = \! 
1,2,\dotsc,N \! + \!1$,
\begin{align}
(z \! - \! \alpha_{k}) \pmb{\pi}_{k}^{n}(z) \underset{\underset{z_{o}=1+
o(1)}{\mathscr{N},n \to \infty}}{=}& \, \left(\tilde{\mathcal{X}}_{11}^{\tilde{a},1}
(z) \left(1 \! + \! \dfrac{1}{(n \! - \! 1)K \! + \! k} \left(\tilde{\mathcal{R}}_{11}^{
\sharp}(z) \! - \! \tilde{\mathcal{R}}_{11}^{\natural}(z) \right) \! + \! \mathcal{O} 
\left(\dfrac{\mathfrak{c}(n,k,z_{o})}{((n \! - \! 1)K \! + \! k)^{2}} \right) \right) 
\right. \nonumber \\
+&\left. \, \tilde{\mathcal{X}}_{21}^{\tilde{a},1}(z) \left(\dfrac{1}{(n \! - \! 1)K \! 
+ \! k} \left(\tilde{\mathcal{R}}_{12}^{\sharp}(z) \! - \! \tilde{\mathcal{R}}_{12}^{
\natural}(z) \right) \! + \! \mathcal{O} \left(\dfrac{\mathfrak{c}(n,k,z_{o})}{((n 
\! - \! 1)K \! + \! k)^{2}} \right) \right) \right) \nonumber \\
\times& \, \mathscr{E} \me^{n(g^{f}(z)-\hat{\mathscr{P}}_{0}^{+})}, 
\label{eqmainfin79}
\end{align}
and
\begin{align}
(z \! - \! \alpha_{k}) \int_{\mathbb{R}} \dfrac{((\xi \! - \! \alpha_{k}) 
\pmb{\pi}_{k}^{n}(\xi)) \me^{-n \widetilde{V}(\xi)}}{(\xi \! - \! \alpha_{k})
(\xi \! - \! z)} \, \dfrac{\md \xi}{2 \pi \mi} \underset{\underset{z_{o}=1+o(1)}{
\mathscr{N},n \to \infty}}{=}& \, \left(\tilde{\mathcal{X}}_{12}^{\tilde{a},1}(z) 
\left(1 \! + \! \dfrac{1}{(n \! - \! 1)K \! + \! k} \left(\tilde{\mathcal{R}}_{11}^{
\sharp}(z) \! - \! \tilde{\mathcal{R}}_{11}^{\natural}(z) \right) \! + \! \mathcal{O} 
\left(\dfrac{\mathfrak{c}(n,k,z_{o})}{((n \! - \! 1)K \! + \! k)^{2}} \right) \right) 
\right. \nonumber \\
+&\left. \, \tilde{\mathcal{X}}_{22}^{\tilde{a},1}(z) \left(\dfrac{1}{(n \! - \! 1)K \! 
+ \! k} \left(\tilde{\mathcal{R}}_{12}^{\sharp}(z) \! - \! \tilde{\mathcal{R}}_{12}^{
\natural}(z) \right) \! + \! \mathcal{O} \left(\dfrac{\mathfrak{c}(n,k,z_{o})}{((
n \! - \! 1)K \! + \! k)^{2}} \right) \right) \right) \nonumber \\
\times& \, \mathscr{E}^{-1} \me^{n \tilde{\ell}} \me^{-n(g^{f}(z)-
\hat{\mathscr{P}}_{0}^{+})}, \label{eqmainfin80}
\end{align}
where
\begin{align}
\tilde{\mathcal{X}}_{11}^{\tilde{a},1}(z) =& \, -\mi \sqrt{\pi} \me^{\frac{1}{2} 
((n-1)K+k) \tilde{\xi}_{\tilde{a}_{j}}(z)} \mathscr{E}^{-1} \left(\mi \left(
\tilde{\mathbb{K}}_{11} \tilde{\mathbb{M}}_{11}(z) \! + \! \tilde{\mathbb{
K}}_{12} \tilde{\mathbb{M}}_{21}(z) \right) \Theta_{0,0}^{0,0,-}(\tilde{\Phi}_{
\tilde{a}_{j}}(z)) \! + \! \left(\tilde{\mathbb{K}}_{11} \tilde{\mathbb{M}}_{12}
(z) \! + \! \tilde{\mathbb{K}}_{12} \tilde{\mathbb{M}}_{22}(z) \right) \right. 
\nonumber \\
\times&\left. \, \Theta_{0,0}^{0,0,+}(\tilde{\Phi}_{\tilde{a}_{j}}(z)) 
\me^{-\mi ((n-1)K+k) \tilde{\mho}_{j}} \right), \label{eqmainfin81} \\
\tilde{\mathcal{X}}_{12}^{\tilde{a},1}(z) =& \, \sqrt{\pi} \me^{-\frac{\mi \pi}{6}} 
\me^{-\frac{1}{2}((n-1)K+k) \tilde{\xi}_{\tilde{a}_{j}}(z)} \mathscr{E}^{-1} \left(
\left(\tilde{\mathbb{K}}_{11} \tilde{\mathbb{M}}_{12}(z) \! + \! \tilde{
\mathbb{K}}_{12} \tilde{\mathbb{M}}_{22}(z) \right) \Theta_{2,2}^{0,2,+}
(\tilde{\Phi}_{\tilde{a}_{j}}(z)) \! + \! \mi \left(\tilde{\mathbb{K}}_{11} 
\tilde{\mathbb{M}}_{11}(z) \! + \! \tilde{\mathbb{K}}_{12} \tilde{\mathbb{
M}}_{21}(z) \right) \right. \nonumber \\
\times&\left. \, \Theta_{2,2}^{0,2,-}(\tilde{\Phi}_{\tilde{a}_{j}}(z)) 
\me^{\mi ((n-1)K+k) \tilde{\mho}_{j}} \right), \label{eqmainfin82} \\
\tilde{\mathcal{X}}_{21}^{\tilde{a},1}(z) =& \, -\mi \sqrt{\pi} \me^{\frac{1}{2} 
((n-1)K+k) \tilde{\xi}_{\tilde{a}_{j}}(z)} \mathscr{E} \left(\mi \left(\tilde{
\mathbb{K}}_{21} \tilde{\mathbb{M}}_{11}(z) \! + \! \tilde{\mathbb{K}}_{22} 
\tilde{\mathbb{M}}_{21}(z) \right) \Theta_{0,0}^{0,0,-}(\tilde{\Phi}_{
\tilde{a}_{j}}(z)) \! + \! \left(\tilde{\mathbb{K}}_{21} \tilde{\mathbb{M}}_{12}
(z) \! + \! \tilde{\mathbb{K}}_{22} \tilde{\mathbb{M}}_{22}(z) \right) \right. 
\nonumber \\
\times&\left. \, \Theta_{0,0}^{0,0,+}(\tilde{\Phi}_{\tilde{a}_{j}}(z)) 
\me^{-\mi ((n-1)K+k) \tilde{\mho}_{j}} \right), \label{eqmainfin83} \\
\tilde{\mathcal{X}}_{22}^{\tilde{a},1}(z) =& \, \sqrt{\pi} \me^{-\frac{\mi \pi}{6}} 
\me^{-\frac{1}{2}((n-1)K+k) \tilde{\xi}_{\tilde{a}_{j}}(z)} \mathscr{E} \left(
\left(\tilde{\mathbb{K}}_{21} \tilde{\mathbb{M}}_{12}(z) \! + \! \tilde{
\mathbb{K}}_{22} \tilde{\mathbb{M}}_{22}(z) \right) \Theta_{2,2}^{0,2,+}
(\tilde{\Phi}_{\tilde{a}_{j}}(z)) \! + \! \mi \left(\tilde{\mathbb{K}}_{21} 
\tilde{\mathbb{M}}_{11}(z) \! + \! \tilde{\mathbb{K}}_{22} \tilde{\mathbb{
M}}_{21}(z) \right) \right. \nonumber \\
\times&\left. \, \Theta_{2,2}^{0,2,-}(\tilde{\Phi}_{\tilde{a}_{j}}(z)) 
\me^{\mi ((n-1)K+k) \tilde{\mho}_{j}} \right), \label{eqmainfin84}
\end{align}
with $\tilde{\Phi}_{\tilde{a}_{j}}(z)$ defined by Equation~\eqref{eqmainfin85}$;$\\
{\rm \pmb{(7)}} for $z \! \in \! \tilde{\Omega}^{2}_{\tilde{b}_{j-1}}$, 
$j \! = \! 1,2,\dotsc,N \! + \!1$,
\begin{align}
(z \! - \! \alpha_{k}) \pmb{\pi}_{k}^{n}(z) \underset{\underset{z_{o}=
1+o(1)}{\mathscr{N},n \to \infty}}{=}& \, \left(\left(\tilde{\mathcal{X}}_{11}^{
\tilde{b},2}(z) \! + \! \tilde{\mathcal{X}}_{12}^{\tilde{b},2}(z) \me^{-2 \pi \mi 
((n-1)K+k) \int_{z}^{\tilde{a}_{N+1}} \psi_{\widetilde{V}}^{f}(\xi) \, \md \xi} 
\right) \left(1 \! + \! \dfrac{1}{(n \! - \! 1)K \! + \! k} \left(\tilde{\mathcal{
R}}_{11}^{\sharp}(z) \! - \! \tilde{\mathcal{R}}_{11}^{\natural}(z) \right) \! + \! 
\mathcal{O} \left(\dfrac{\mathfrak{c}(n,k,z_{o})}{((n \! - \! 1)K \! + \! k)^{2}} 
\right) \right) \right. \nonumber \\
+&\left. \, \left(\tilde{\mathcal{X}}_{21}^{\tilde{b},2}(z) \! + \! \tilde{\mathcal{
X}}_{22}^{\tilde{b},2}(z) \me^{-2 \pi \mi ((n-1)K+k) \int_{z}^{\tilde{a}_{N+1}} 
\psi_{\widetilde{V}}^{f}(\xi) \, \md \xi} \right) \left(\dfrac{1}{(n \! - \! 1)K \! + 
\! k} \left(\tilde{\mathcal{R}}_{12}^{\sharp}(z) \! - \! \tilde{\mathcal{R}}_{12}^{
\natural}(z) \right) \! + \! \mathcal{O} \left(\dfrac{\mathfrak{c}(n,k,z_{o})}{((n 
\! - \! 1)K \! + \! k)^{2}} \right) \right) \right) \nonumber \\
\times& \, \mathscr{E} \me^{n(g^{f}(z)-\hat{\mathscr{P}}_{0}^{+})}, 
\label{eqmainfin86}
\end{align}
and
\begin{align}
(z \! - \! \alpha_{k}) \int_{\mathbb{R}} \dfrac{((\xi \! - \! \alpha_{k}) 
\pmb{\pi}_{k}^{n}(\xi)) \me^{-n \widetilde{V}(\xi)}}{(\xi \! - \! \alpha_{k})
(\xi \! - \! z)} \, \dfrac{\md \xi}{2 \pi \mi} \underset{\underset{z_{o}=1+o(1)}{
\mathscr{N},n \to \infty}}{=}& \, \left(\tilde{\mathcal{X}}_{12}^{\tilde{b},2}(z) 
\left(1 \! + \! \dfrac{1}{(n \! - \! 1)K \! + \! k} \left(\tilde{\mathcal{R}}_{11}^{
\sharp}(z) \! - \! \tilde{\mathcal{R}}_{11}^{\natural}(z) \right) \! + \! \mathcal{O} 
\left(\dfrac{\mathfrak{c}(n,k,z_{o})}{((n \! - \! 1)K \! + \! k)^{2}} \right) \right) 
\right. \nonumber \\
+&\left. \, \tilde{\mathcal{X}}_{22}^{\tilde{b},2}(z) \left(\dfrac{1}{(n \! - \! 1)K \! 
+ \! k} \left(\tilde{\mathcal{R}}_{12}^{\sharp}(z) \! - \! \tilde{\mathcal{R}}_{12}^{
\natural}(z) \right) \! + \! \mathcal{O} \left(\dfrac{\mathfrak{c}(n,k,z_{o})}{((n \! 
- \! 1)K \! + \! k)^{2}} \right) \right) \right) \nonumber \\
\times& \, \mathscr{E}^{-1} \me^{n \tilde{\ell}} \me^{-n(g^{f}(z)-
\hat{\mathscr{P}}_{0}^{+})}, \label{eqmainfin87}
\end{align}
where
\begin{align}
\tilde{\mathcal{X}}_{11}^{\tilde{b},2}(z) =& \, -\mi \sqrt{\pi} \me^{\frac{1}{2} 
((n-1)K+k) \tilde{\xi}_{\tilde{b}_{j-1}}(z)} \mathscr{E}^{-1}  \left(\mi \left(
\tilde{\mathbb{K}}_{11} \tilde{\mathbb{M}}_{11}(z) \! + \! \tilde{\mathbb{
K}}_{12} \tilde{\mathbb{M}}_{21}(z) \right) \varTheta_{2,0}^{0,2,-}(\tilde{\Phi}_{
\tilde{b}_{j-1}}(z)) \! - \! \left(\tilde{\mathbb{K}}_{11} \tilde{\mathbb{M}}_{12}
(z) \! + \! \tilde{\mathbb{K}}_{12} \tilde{\mathbb{M}}_{22}(z) \right) \right. 
\nonumber \\
\times&\left. \, \varTheta_{2,0}^{0,2,+}(\tilde{\Phi}_{\tilde{b}_{j-1}}(z)) 
\me^{-\mi ((n-1)K+k) \tilde{\mho}_{j-1}} \right), \label{eqmainfin88} \\
\tilde{\mathcal{X}}_{12}^{\tilde{b},2}(z) =& \, \sqrt{\pi} \me^{-\frac{\mi \pi}{6}} 
\me^{-\frac{1}{2}((n-1)K+k) \tilde{\xi}_{\tilde{b}_{j-1}}(z)} \mathscr{E}^{-1} 
\left(\left(\tilde{\mathbb{K}}_{11} \tilde{\mathbb{M}}_{12}(z) \! + \! \tilde{
\mathbb{K}}_{12} \tilde{\mathbb{M}}_{22}(z) \right) \Theta_{2,2}^{0,2,+}
(\tilde{\Phi}_{\tilde{b}_{j-1}}(z)) \! - \! \mi \left(\tilde{\mathbb{K}}_{11} 
\tilde{\mathbb{M}}_{11}(z) \! + \! \tilde{\mathbb{K}}_{12} \tilde{\mathbb{
M}}_{21}(z) \right) \right. \nonumber \\
\times&\left. \, \Theta_{2,2}^{0,2,-}(\tilde{\Phi}_{\tilde{b}_{j-1}}(z)) 
\me^{\mi ((n-1)K+k) \tilde{\mho}_{j-1}} \right), \label{eqmainfin89} \\
\tilde{\mathcal{X}}_{21}^{\tilde{b},2}(z) =& \, -\mi \sqrt{\pi} \me^{\frac{1}{2} 
((n-1)K+k) \tilde{\xi}_{\tilde{b}_{j-1}}(z)} \mathscr{E} \left(\mi \left(\tilde{
\mathbb{K}}_{21} \tilde{\mathbb{M}}_{11}(z) \! + \! \tilde{\mathbb{K}}_{22} 
\tilde{\mathbb{M}}_{21}(z) \right) \varTheta_{2,0}^{0,2,-}(\tilde{\Phi}_{
\tilde{b}_{j-1}}(z)) \! - \! \left(\tilde{\mathbb{K}}_{21} \tilde{\mathbb{M}}_{12}
(z) \! + \! \tilde{\mathbb{K}}_{22} \tilde{\mathbb{M}}_{22}(z) \right) \right. 
\nonumber \\
\times&\left. \, \varTheta_{2,0}^{0,2,+}(\tilde{\Phi}_{\tilde{b}_{j-1}}(z)) 
\me^{-\mi ((n-1)K+k) \tilde{\mho}_{j-1}} \right), \label{eqmainfin90} \\
\tilde{\mathcal{X}}_{22}^{\tilde{b},2}(z) =& \, \sqrt{\pi} \me^{-\frac{\mi \pi}{6}} 
\me^{-\frac{1}{2} ((n-1)K+k) \tilde{\xi}_{\tilde{b}_{j-1}}(z)} \mathscr{E} \left(
\left(\tilde{\mathbb{K}}_{21} \tilde{\mathbb{M}}_{12}(z) \! + \! \tilde{\mathbb{
K}}_{22} \tilde{\mathbb{M}}_{22}(z) \right) \Theta_{2,2}^{0,2,+}(\tilde{\Phi}_{
\tilde{b}_{j-1}}(z)) \! - \! \mi \left(\tilde{\mathbb{K}}_{21} \tilde{\mathbb{M}}_{
11}(z) \! + \! \tilde{\mathbb{K}}_{22} \tilde{\mathbb{M}}_{21}(z) \right) \right. 
\nonumber \\
\times&\left. \, \Theta_{2,2}^{0,2,-}(\tilde{\Phi}_{\tilde{b}_{j-1}}(z)) 
\me^{\mi ((n-1)K+k) \tilde{\mho}_{j-1}} \right), \label{eqmainfin91}
\end{align}
with $\varTheta_{r_{3},r_{4}}^{r_{1},r_{2},\pm}(\zeta)$, $r_{1},r_{2},r_{3},r_{4} 
\! = \! 0,1,2$, defined by Equation~\eqref{eqmaininf91}$;$\\
{\rm \pmb{(8)}} for $z \! \in \! \tilde{\Omega}^{2}_{\tilde{a}_{j}}$, $j \! = \! 
1,2,\dotsc,N \! + \!1$,
\begin{align}
(z \! - \! \alpha_{k}) \pmb{\pi}_{k}^{n}(z) \underset{\underset{z_{o}=
1+o(1)}{\mathscr{N},n \to \infty}}{=}& \, \left(\left(\tilde{\mathcal{X}}_{11}^{
\tilde{a},2}(z) \! + \! \tilde{\mathcal{X}}_{12}^{\tilde{a},2}(z) \me^{-2 \pi \mi 
((n-1)K+k) \int_{z}^{\tilde{a}_{N+1}} \psi_{\widetilde{V}}^{f}(\xi) \, \md \xi} \right) 
\left(1 \! + \! \dfrac{1}{(n \! - \! 1)K \! + \! k} \left(\tilde{\mathcal{R}}_{11}^{
\sharp}(z) \! - \! \tilde{\mathcal{R}}_{11}^{\natural}(z) \right) \! + \! \mathcal{O} 
\left(\dfrac{\mathfrak{c}(n,k,z_{o})}{((n \! - \! 1)K \! + \! k)^{2}} \right) \right) 
\right. \nonumber \\
+&\left. \, \left(\tilde{\mathcal{X}}_{21}^{\tilde{a},2}(z) \! + \! \tilde{\mathcal{
X}}_{22}^{\tilde{a},2}(z) \me^{-2 \pi \mi ((n-1)K+k) \int_{z}^{\tilde{a}_{N+1}} 
\psi_{\widetilde{V}}^{f}(\xi) \, \md \xi} \right) \left(\dfrac{1}{(n \! - \! 1)K \! + 
\! k} \left(\tilde{\mathcal{R}}_{12}^{\sharp}(z) \! - \! \tilde{\mathcal{R}}_{12}^{
\natural}(z) \right) \! + \! \mathcal{O} \left(\dfrac{\mathfrak{c}(n,k,z_{o})}{((n 
\! - \! 1)K \! + \! k)^{2}} \right) \right) \right) \nonumber \\
\times& \, \mathscr{E} \me^{n(g^{f}(z)-\hat{\mathscr{P}}_{0}^{+})}, 
\label{eqmainfin92}
\end{align}
and
\begin{align}
(z \! - \! \alpha_{k}) \int_{\mathbb{R}} \dfrac{((\xi \! - \! \alpha_{k}) 
\pmb{\pi}_{k}^{n}(\xi)) \me^{-n \widetilde{V}(\xi)}}{(\xi \! - \! \alpha_{k})
(\xi \! - \! z)} \, \dfrac{\md \xi}{2 \pi \mi} \underset{\underset{z_{o}=1+o(1)}{
\mathscr{N},n \to \infty}}{=}& \, \left(\tilde{\mathcal{X}}_{12}^{\tilde{a},2}(z) 
\left(1 \! + \! \dfrac{1}{(n \! - \! 1)K \! + \! k} \left(\tilde{\mathcal{R}}_{11}^{
\sharp}(z) \! - \! \tilde{\mathcal{R}}_{11}^{\natural}(z) \right) \! + \! \mathcal{O} 
\left(\dfrac{\mathfrak{c}(n,k,z_{o})}{((n \! - \! 1)K \! + \! k)^{2}} \right) \right) 
\right. \nonumber \\
+&\left. \, \tilde{\mathcal{X}}_{22}^{\tilde{a},2}(z) \left(\dfrac{1}{(n \! - \! 1)K \! 
+ \! k} \left(\tilde{\mathcal{R}}_{12}^{\sharp}(z) \! - \! \tilde{\mathcal{R}}_{12}^{
\natural}(z) \right) \! + \! \mathcal{O} \left(\dfrac{\mathfrak{c}(n,k,z_{o})}{((n \! 
- \! 1)K \! + \! k)^{2}} \right) \right) \right) \nonumber \\
\times& \, \mathscr{E}^{-1} \me^{n \tilde{\ell}} 
\me^{-n(g^{f}(z)-\hat{\mathscr{P}}_{0}^{+})}, \label{eqmainfin93}
\end{align}
where
\begin{align}
\tilde{\mathcal{X}}_{11}^{\tilde{a},2}(z) =& \, -\mi \sqrt{\pi} \me^{\frac{1}{2} 
((n-1)K+k) \tilde{\xi}_{\tilde{a}_{j}}(z)} \mathscr{E}^{-1} \left(\mi \left(\tilde{
\mathbb{K}}_{11} \tilde{\mathbb{M}}_{11}(z) \! + \! \tilde{\mathbb{K}}_{12} 
\tilde{\mathbb{M}}_{21}(z) \right) \varTheta_{2,0}^{0,2,-}(\tilde{\Phi}_{\tilde{a}_{j}}
(z)) \! + \! \left(\tilde{\mathbb{K}}_{11} \tilde{\mathbb{M}}_{12}(z) \! + \! \tilde{
\mathbb{K}}_{12} \tilde{\mathbb{M}}_{22}(z) \right) \right. \nonumber \\
\times&\left. \, \varTheta_{2,0}^{0,2,+}(\tilde{\Phi}_{\tilde{a}_{j}}(z)) 
\me^{-\mi ((n-1)K+k) \tilde{\mho}_{j}} \right), \label{eqmainfin94} \\
\tilde{\mathcal{X}}_{12}^{\tilde{a},2}(z) =& \, \sqrt{\pi} \me^{-\frac{\mi \pi}{6}} 
\me^{-\frac{1}{2}((n-1)K+k) \tilde{\xi}_{\tilde{a}_{j}}(z)} \mathscr{E}^{-1} \left(
\left(\tilde{\mathbb{K}}_{11} \tilde{\mathbb{M}}_{12}(z) \! + \! \tilde{\mathbb{
K}}_{12} \tilde{\mathbb{M}}_{22}(z) \right) \Theta_{2,2}^{0,2,+}(\tilde{\Phi}_{
\tilde{a}_{j}}(z)) \! + \! \mi \left(\tilde{\mathbb{K}}_{11} \tilde{\mathbb{M}}_{11}
(z) \! + \! \tilde{\mathbb{K}}_{12} \tilde{\mathbb{M}}_{21}(z) \right) \right. 
\nonumber \\
\times&\left. \, \Theta_{2,2}^{0,2,-}(\tilde{\Phi}_{\tilde{a}_{j}}(z)) 
\me^{\mi ((n-1)K+k) \tilde{\mho}_{j}} \right), \label{eqmainfin95} \\
\tilde{\mathcal{X}}_{21}^{\tilde{a},2}(z) =& \, -\mi \sqrt{\pi} \me^{\frac{1}{2} 
((n-1)K+k) \tilde{\xi}_{\tilde{a}_{j}}(z)} \mathscr{E} \left(\mi \left(\tilde{
\mathbb{K}}_{21} \tilde{\mathbb{M}}_{11}(z) \! + \! \tilde{\mathbb{K}}_{22} 
\tilde{\mathbb{M}}_{21}(z) \right) \varTheta_{2,0}^{0,2,-}(\tilde{\Phi}_{\tilde{a}_{j}}
(z)) \! + \! \left(\tilde{\mathbb{K}}_{21} \tilde{\mathbb{M}}_{12}(z) \! + \! \tilde{
\mathbb{K}}_{22} \tilde{\mathbb{M}}_{22}(z) \right) \right. \nonumber \\
\times&\left. \, \varTheta_{2,0}^{0,2,+}(\tilde{\Phi}_{\tilde{a}_{j}}(z)) 
\me^{-\mi ((n-1)K+k) \tilde{\mho}_{j}} \right), \label{eqmainfin96} \\
\tilde{\mathcal{X}}_{22}^{\tilde{a},2}(z) =& \, \sqrt{\pi} \me^{-\frac{\mi \pi}{6}} 
\me^{-\frac{1}{2}((n-1)K+k) \tilde{\xi}_{\tilde{a}_{j}}(z)} \mathscr{E} \left(\left(
\tilde{\mathbb{K}}_{21} \tilde{\mathbb{M}}_{12}(z) \! + \! \tilde{\mathbb{K}}_{22} 
\tilde{\mathbb{M}}_{22}(z) \right) \Theta_{2,2}^{0,2,+}(\tilde{\Phi}_{\tilde{a}_{j}}
(z)) \! + \! \mi \left(\tilde{\mathbb{K}}_{21} \tilde{\mathbb{M}}_{11}(z) \! + \! 
\tilde{\mathbb{K}}_{22} \tilde{\mathbb{M}}_{21}(z) \right) \right. \nonumber \\
\times&\left. \, \Theta_{2,2}^{0,2,-}(\tilde{\Phi}_{\tilde{a}_{j}}(z)) 
\me^{\mi ((n-1)K+k) \tilde{\mho}_{j}} \right); \label{eqmainfin97}
\end{align}
{\rm \pmb{(9)}} for $z \! \in \! \tilde{\Omega}^{3}_{\tilde{b}_{j-1}}$, $j \! = \! 
1,2,\dotsc,N \! + \! 1$,
\begin{align}
(z \! - \! \alpha_{k}) \pmb{\pi}_{k}^{n}(z) \underset{\underset{z_{o}=
1+o(1)}{\mathscr{N},n \to \infty}}{=}& \, \left(\left(\tilde{\mathcal{X}}_{11}^{
\tilde{b},3}(z) \! - \! \tilde{\mathcal{X}}_{12}^{\tilde{b},3}(z) \me^{2 \pi \mi 
((n-1)K+k) \int_{z}^{\tilde{a}_{N+1}} \psi_{\widetilde{V}}^{f}(\xi) \, \md \xi} \right) 
\left(1 \! + \! \dfrac{1}{(n \! - \! 1)K \! + \! k} \left(\tilde{\mathcal{R}}_{11}^{
\sharp}(z) \! - \! \tilde{\mathcal{R}}_{11}^{\natural}(z) \right) \! + \! \mathcal{O} 
\left(\dfrac{\mathfrak{c}(n,k,z_{o})}{((n \! - \! 1)K \! + \! k)^{2}} \right) \right) 
\right. \nonumber \\
+&\left. \, \left(\tilde{\mathcal{X}}_{21}^{\tilde{b},3}(z) \! - \! \tilde{\mathcal{
X}}_{22}^{\tilde{b},3}(z) \me^{2 \pi \mi ((n-1)K+k) \int_{z}^{\tilde{a}_{N+1}} 
\psi_{\widetilde{V}}^{f}(\xi) \, \md \xi} \right) \left(\dfrac{1}{(n \! - \! 1)K \! + 
\! k} \left(\tilde{\mathcal{R}}_{12}^{\sharp}(z) \! - \! \tilde{\mathcal{R}}_{12}^{
\natural}(z) \right) \! + \! \mathcal{O} \left(\dfrac{\mathfrak{c}(n,k,z_{o})}{((n 
\! - \! 1)K \! + \! k)^{2}} \right) \right) \right) \nonumber \\
\times& \, \mathscr{E}^{-1} \me^{n(g^{f}(z)-\hat{\mathscr{P}}_{0}^{-})}, 
\label{eqmainfin98}
\end{align}
and
\begin{align}
(z \! - \! \alpha_{k}) \int_{\mathbb{R}} \dfrac{((\xi \! - \! \alpha_{k}) 
\pmb{\pi}_{k}^{n}(\xi)) \me^{-n \widetilde{V}(\xi)}}{(\xi \! - \! \alpha_{k})
(\xi \! - \! z)} \, \dfrac{\md \xi}{2 \pi \mi} \underset{\underset{z_{o}=1+o(1)}{
\mathscr{N},n \to \infty}}{=}& \,\left(\tilde{\mathcal{X}}_{12}^{\tilde{b},3}(z) 
\left(1 \! + \! \dfrac{1}{(n \! - \! 1)K \! + \! k} \left(\tilde{\mathcal{R}}_{11}^{
\sharp}(z) \! - \! \tilde{\mathcal{R}}_{11}^{\natural}(z) \right) \! + \! \mathcal{O} 
\left(\dfrac{\mathfrak{c}(n,k,z_{o})}{((n \! - \! 1)K \! + \! k)^{2}} \right) \right) 
\right. \nonumber \\
+&\left. \, \tilde{\mathcal{X}}_{22}^{\tilde{b},3}(z) \left(\dfrac{1}{(n \! - \! 1)K \! 
+ \! k} \left(\tilde{\mathcal{R}}_{12}^{\sharp}(z) \! - \! \tilde{\mathcal{R}}_{12}^{
\natural}(z) \right) \! + \! \mathcal{O} \left(\dfrac{\mathfrak{c}(n,k,z_{o})}{((n \! 
- \! 1)K \! + \! k)^{2}} \right) \right) \right) \nonumber \\
\times& \, \mathscr{E} \me^{n \tilde{\ell}} \me^{-n(g^{f}(z)-
\hat{\mathscr{P}}_{0}^{-})}, \label{eqmainfin99}
\end{align}
where
\begin{align}
\tilde{\mathcal{X}}_{11}^{\tilde{b},3}(z) =& \, -\mi \sqrt{\pi} \me^{\frac{1}{2} 
((n-1)K+k) \tilde{\xi}_{\tilde{b}_{j-1}}(z)} \mathscr{E}^{-1} \left(\mi \left(
\tilde{\mathbb{K}}_{11} \tilde{\mathbb{M}}_{12}(z) \! + \! \tilde{\mathbb{
K}}_{12} \tilde{\mathbb{M}}_{22}(z) \right) \varTheta_{0,1}^{2,1,-}(\tilde{
\Phi}_{\tilde{b}_{j-1}}(z)) \! + \! \left(\tilde{\mathbb{K}}_{11} \tilde{\mathbb{
M}}_{11}(z) \! + \! \tilde{\mathbb{K}}_{12} \tilde{\mathbb{M}}_{21}(z) \right) 
\right. \nonumber \\
\times&\left. \, \varTheta_{0,1}^{2,1,+}(\tilde{\Phi}_{\tilde{b}_{j-1}}(z)) 
\me^{\mi ((n-1)K+k) \tilde{\mho}_{j-1}} \right), \label{eqmainfin100} \\
\tilde{\mathcal{X}}_{12}^{\tilde{b},3}(z) =& \, \sqrt{\pi} \me^{-\frac{\mi \pi}{6}} 
\me^{-\frac{1}{2} ((n-1)K+k) \tilde{\xi}_{\tilde{b}_{j-1}}(z)} \mathscr{E}^{-1}  
\left(\left(\tilde{\mathbb{K}}_{11} \tilde{\mathbb{M}}_{11}(z) \! + \! \tilde{
\mathbb{K}}_{21} \tilde{\mathbb{M}}_{21}(z) \right) \Theta_{0,1}^{2,1,+}
(\tilde{\Phi}_{\tilde{b}_{j-1}}(z)) \! + \! \mi \left(\tilde{\mathbb{K}}_{11} 
\tilde{\mathbb{M}}_{12}(z) \! + \! \tilde{\mathbb{K}}_{12} \tilde{\mathbb{
M}}_{22}(z) \right) \right. \nonumber \\
\times&\left. \, \Theta_{0,1}^{2,1,-}(\tilde{\Phi}_{\tilde{b}_{j-1}}(z)) 
\me^{-\mi ((n-1)K+k) \tilde{\mho}_{j-1}} \right), \label{eqmainfin101} \\
\tilde{\mathcal{X}}_{21}^{\tilde{b},3}(z) =& \, -\mi \sqrt{\pi} \me^{\frac{1}{2} 
((n-1)K+k) \tilde{\xi}_{\tilde{b}_{j-1}}(z)} \mathscr{E} \left(\mi \left(\tilde{
\mathbb{K}}_{21} \tilde{\mathbb{M}}_{12}(z) \! + \! \tilde{\mathbb{K}}_{22} 
\tilde{\mathbb{M}}_{22}(z) \right) \varTheta_{0,1}^{2,1,-}(\tilde{\Phi}_{
\tilde{b}_{j-1}}(z)) \! + \! \left(\tilde{\mathbb{K}}_{21} \tilde{\mathbb{M}}_{11}
(z) \! + \! \tilde{\mathbb{K}}_{22} \tilde{\mathbb{M}}_{21}(z) \right) \right. 
\nonumber \\
\times&\left. \, \varTheta_{0,1}^{2,1,+}(\tilde{\Phi}_{\tilde{b}_{j-1}}(z)) 
\me^{\mi ((n-1)K+k) \tilde{\mho}_{j-1}} \right), \label{eqmainfin102} \\
\tilde{\mathcal{X}}_{22}^{\tilde{b},3}(z) =& \, \sqrt{\pi} \me^{-\frac{\mi \pi}{6}} 
\me^{-\frac{1}{2}((n-1)K+k) \tilde{\xi}_{\tilde{b}_{j-1}}(z)} \mathscr{E} \left(
\left(\tilde{\mathbb{K}}_{21} \tilde{\mathbb{M}}_{11}(z) \! + \! \tilde{\mathbb{
K}}_{22} \tilde{\mathbb{M}}_{21}(z) \right) \Theta_{0,1}^{2,1,+}(\tilde{\Phi}_{
\tilde{b}_{j-1}}(z)) \! + \! \mi \left(\tilde{\mathbb{K}}_{21} \tilde{\mathbb{
M}}_{12}(z) \! + \! \tilde{\mathbb{K}}_{22} \tilde{\mathbb{M}}_{22}(z) \right) 
\right. \nonumber \\
\times&\left. \, \Theta_{0,1}^{2,1,-}(\tilde{\Phi}_{\tilde{b}_{j-1}}(z)) 
\me^{-\mi ((n-1)K+k) \tilde{\mho}_{j-1}} \right); \label{eqmainfin103}
\end{align}
{\rm \pmb{(10)}} for $z \! \in \! \tilde{\Omega}^{3}_{\tilde{a}_{j}}$, $j \! = 
\! 1,2,\dotsc,N \! + \!1$,
\begin{align}
(z \! - \! \alpha_{k}) \pmb{\pi}_{k}^{n}(z) \underset{\underset{z_{o}=
1+o(1)}{\mathscr{N},n \to \infty}}{=}& \, \left(\left(\tilde{\mathcal{X}}_{11}^{
\tilde{a},3}(z) \! - \! \tilde{\mathcal{X}}_{12}^{\tilde{a},3}(z) \me^{2 \pi \mi 
((n-1)K+k) \int_{z}^{\tilde{a}_{N+1}} \psi_{\widetilde{V}}^{f}(\xi) \, \md \xi} 
\right) \left(1 \! + \! \dfrac{1}{(n \! - \! 1)K \! + \! k} \left(\tilde{\mathcal{
R}}_{11}^{\sharp}(z) \! - \! \tilde{\mathcal{R}}_{11}^{\natural}(z) \right) \! + \! 
\mathcal{O} \left(\dfrac{\mathfrak{c}(n,k,z_{o})}{((n \! - \! 1)K \! + \! k)^{2}} 
\right) \right) \right. \nonumber \\
+&\left. \, \left(\tilde{\mathcal{X}}_{21}^{\tilde{a},3}(z) \! - \! \tilde{\mathcal{
X}}_{22}^{\tilde{a},3}(z) \me^{2 \pi \mi ((n-1)K+k) \int_{z}^{\tilde{a}_{N+1}} 
\psi_{\widetilde{V}}^{f}(\xi) \, \md \xi} \right) \left(\dfrac{1}{(n \! - \! 1)K \! + 
\! k} \left(\tilde{\mathcal{R}}_{12}^{\sharp}(z) \! - \! \tilde{\mathcal{R}}_{12}^{
\natural}(z) \right) \! + \! \mathcal{O} \left(\dfrac{\mathfrak{c}(n,k,z_{o})}{
((n \! - \! 1)K \! + \! k)^{2}} \right) \right) \right) \nonumber \\
\times& \, \mathscr{E}^{-1} \me^{n(g^{f}(z)-\hat{\mathscr{P}}_{0}^{-})}, 
\label{eqmainfin104}
\end{align}
and
\begin{align}
(z \! - \! \alpha_{k}) \int_{\mathbb{R}} \dfrac{((\xi \! - \! \alpha_{k}) 
\pmb{\pi}_{k}^{n}(\xi)) \me^{-n \widetilde{V}(\xi)}}{(\xi \! - \! \alpha_{k})
(\xi \! - \! z)} \, \dfrac{\md \xi}{2 \pi \mi} \underset{\underset{z_{o}=1+o(1)}{
\mathscr{N},n \to \infty}}{=}& \, \left(\tilde{\mathcal{X}}_{12}^{\tilde{a},3}(z) 
\left(1 \! + \! \dfrac{1}{(n \! - \! 1)K \! + \! k} \left(\tilde{\mathcal{R}}_{11}^{
\sharp}(z) \! - \! \tilde{\mathcal{R}}_{11}^{\natural}(z) \right) \! + \! \mathcal{O} 
\left(\dfrac{\mathfrak{c}(n,k,z_{o})}{((n \! - \! 1)K \! + \! k)^{2}} \right) \right) 
\right. \nonumber \\
+&\left. \, \tilde{\mathcal{X}}_{22}^{\tilde{a},3}(z) \left(\dfrac{1}{(n \! - \! 1)K \! 
+ \! k} \left(\tilde{\mathcal{R}}_{12}^{\sharp}(z) \! - \! \tilde{\mathcal{R}}_{12}^{
\natural}(z) \right) \! + \! \mathcal{O} \left(\dfrac{\mathfrak{c}(n,k,z_{o})}{((n \! 
- \! 1)K \! + \! k)^{2}} \right) \right) \right) \nonumber \\
\times& \, \mathscr{E} \me^{n \tilde{\ell}} \me^{-n(g^{f}(z)-
\hat{\mathscr{P}}_{0}^{-})}, \label{eqmainfin105}
\end{align}
where
\begin{align}
\tilde{\mathcal{X}}_{11}^{\tilde{a},3}(z) =& \, -\mi \sqrt{\pi} \me^{\frac{1}{2} 
((n-1)K+k) \tilde{\xi}_{\tilde{a}_{j}}(z)} \mathscr{E}^{-1} \left(\mi \left(\tilde{
\mathbb{K}}_{11} \tilde{\mathbb{M}}_{12}(z) \! + \! \tilde{\mathbb{K}}_{12} 
\tilde{\mathbb{M}}_{22}(z) \right) \varTheta_{0,1}^{2,1,-}(\tilde{\Phi}_{\tilde{a}_{j}}
(z)) \! - \! \left(\tilde{\mathbb{K}}_{11} \tilde{\mathbb{M}}_{11}(z) \! + \! \tilde{
\mathbb{K}}_{12} \tilde{\mathbb{M}}_{21}(z) \right) \right. \nonumber \\
\times&\left. \, \varTheta_{0,1}^{2,1,+}(\tilde{\Phi}_{\tilde{a}_{j}}(z)) 
\me^{\mi ((n-1)K+k) \tilde{\mho}_{j}} \right), \label{eqmainfin106} \\
\tilde{\mathcal{X}}_{12}^{\tilde{a},3}(z) =& \, \sqrt{\pi} \me^{-\frac{\mi \pi}{6}} 
\me^{-\frac{1}{2}((n-1)K+k) \tilde{\xi}_{\tilde{a}_{j}}(z)} \mathscr{E}^{-1} \left(
\left(\tilde{\mathbb{K}}_{11} \tilde{\mathbb{M}}_{11}(z) \! + \! \tilde{\mathbb{
K}}_{21} \tilde{\mathbb{M}}_{21}(z) \right) \Theta_{0,1}^{2,1,+}(\tilde{\Phi}_{
\tilde{a}_{j}}(z)) \! - \! \mi \left(\tilde{\mathbb{K}}_{11} \tilde{\mathbb{M}}_{12}
(z) \! + \! \tilde{\mathbb{K}}_{12} \tilde{\mathbb{M}}_{22}(z) \right) \right. 
\nonumber \\
\times&\left. \, \Theta_{0,1}^{2,1,-}(\tilde{\Phi}_{\tilde{a}_{j}}(z)) 
\me^{-\mi ((n-1)K+k) \tilde{\mho}_{j}} \right), \label{eqmainfin107} \\
\tilde{\mathcal{X}}_{21}^{\tilde{a},3}(z) =& \, -\mi \sqrt{\pi} \me^{\frac{1}{2} 
((n-1)K+k) \tilde{\xi}_{\tilde{a}_{j}}(z)} \mathscr{E} \left(\mi \left(\tilde{
\mathbb{K}}_{21} \tilde{\mathbb{M}}_{12}(z) \! + \! \tilde{\mathbb{K}}_{22} 
\tilde{\mathbb{M}}_{22}(z) \right) \varTheta_{0,1}^{2,1,-}(\tilde{\Phi}_{\tilde{a}_{j}}
(z)) \! - \! \left(\tilde{\mathbb{K}}_{21} \tilde{\mathbb{M}}_{11}(z) \! + \! \tilde{
\mathbb{K}}_{22} \tilde{\mathbb{M}}_{21}(z) \right) \right. \nonumber \\
\times&\left. \, \varTheta_{0,1}^{2,1,+}(\tilde{\Phi}_{\tilde{a}_{j}}(z)) 
\me^{\mi ((n-1)K+k) \tilde{\mho}_{j}} \right), \label{eqmainfin108} \\
\tilde{\mathcal{X}}_{22}^{\tilde{a},3}(z) =& \, \sqrt{\pi} \me^{-\frac{\mi \pi}{6}} 
\me^{-\frac{1}{2}((n-1)K+k) \tilde{\xi}_{\tilde{a}_{j}}(z)} \mathscr{E} \left(\left(
\tilde{\mathbb{K}}_{21} \tilde{\mathbb{M}}_{11}(z) \! + \! \tilde{\mathbb{K}}_{22} 
\tilde{\mathbb{M}}_{21}(z) \right) \Theta_{0,1}^{2,1,+}(\tilde{\Phi}_{\tilde{a}_{j}}(z)) 
\! - \! \mi \left(\tilde{\mathbb{K}}_{21} \tilde{\mathbb{M}}_{12}(z) \! + \! \tilde{
\mathbb{K}}_{22} \tilde{\mathbb{M}}_{22}(z) \right) \right. \nonumber \\
\times&\left. \, \Theta_{0,1}^{2,1,-}(\tilde{\Phi}_{\tilde{a}_{j}}(z)) 
\me^{-\mi ((n-1)K+k) \tilde{\mho}_{j}} \right); \label{eqmainfin109}
\end{align}
{\rm \pmb{(11)}} for $z \! \in \! \tilde{\Omega}^{4}_{\tilde{b}_{j-1}}$, $j \! 
= \! 1,2,\dotsc,N \! + \!1$,
\begin{align}
(z \! - \! \alpha_{k}) \pmb{\pi}_{k}^{n}(z) \underset{\underset{z_{o}=
1+o(1)}{\mathscr{N},n \to \infty}}{=}& \, \left(\tilde{\mathcal{X}}_{11}^{\tilde{b},4}
(z) \left(1 \! + \! \dfrac{1}{(n \! - \! 1)K \! + \! k} \left(\tilde{\mathcal{R}}_{11}^{
\sharp}(z) \! - \! \tilde{\mathcal{R}}_{11}^{\natural}(z) \right) \! + \! \mathcal{O} 
\left(\dfrac{\mathfrak{c}(n,k,z_{o})}{((n \! - \! 1)K \! + \! k)^{2}} \right) \right) 
\right. \nonumber \\
+&\left. \, \tilde{\mathcal{X}}_{21}^{\tilde{b},4}(z) \left(\dfrac{1}{(n \! - \! 1)K \! 
+ \! k} \left(\tilde{\mathcal{R}}_{12}^{\sharp}(z) \! - \! \tilde{\mathcal{R}}_{12}^{
\natural}(z) \right) \! + \! \mathcal{O} \left(\dfrac{\mathfrak{c}(n,k,z_{o})}{((n 
\! - \! 1)K \! + \! k)^{2}} \right) \right) \right) \nonumber \\
\times& \, \mathscr{E}^{-1} \me^{n(g^{f}(z)-\hat{\mathscr{P}}_{0}^{-})}, 
\label{eqmainfin110}
\end{align}
and
\begin{align}
(z \! - \! \alpha_{k}) \int_{\mathbb{R}} \dfrac{((\xi \! - \! \alpha_{k}) 
\pmb{\pi}_{k}^{n}(\xi)) \me^{-n \widetilde{V}(\xi)}}{(\xi \! - \! \alpha_{k})
(\xi \! - \! z)} \, \dfrac{\md \xi}{2 \pi \mi} \underset{\underset{z_{o}=1+o(1)}{
\mathscr{N},n \to \infty}}{=}& \, \left(\tilde{\mathcal{X}}_{12}^{\tilde{b},4}(z) 
\left(1 \! + \! \dfrac{1}{(n \! - \! 1)K \! + \! k} \left(\tilde{\mathcal{R}}_{11}^{
\sharp}(z) \! - \! \tilde{\mathcal{R}}_{11}^{\natural}(z) \right) \! + \! \mathcal{O} 
\left(\dfrac{\mathfrak{c}(n,k,z_{o})}{((n \! - \! 1)K \! + \! k)^{2}} \right) \right) 
\right. \nonumber \\
+&\left. \, \tilde{\mathcal{X}}_{22}^{\tilde{b},4}(z) \left(\dfrac{1}{(n \! - \! 1)K \! 
+ \! k} \left(\tilde{\mathcal{R}}_{12}^{\sharp}(z) \! - \! \tilde{\mathcal{R}}_{12}^{
\natural}(z) \right) \! + \! \mathcal{O} \left(\dfrac{\mathfrak{c}(n,k,z_{o})}{((n \! 
- \! 1)K \! + \! k)^{2}} \right) \right) \right) \nonumber \\
\times& \, \mathscr{E} \me^{n \tilde{\ell}} \me^{-n(g^{f}(z)-
\hat{\mathscr{P}}_{0}^{-})}, \label{eqmainfin111}
\end{align}
where
\begin{align}
\tilde{\mathcal{X}}_{11}^{\tilde{b},4}(z) =& \, -\mi \sqrt{\pi} \me^{\frac{1}{2} 
((n-1)K+k) \tilde{\xi}_{\tilde{b}_{j-1}}(z)} \mathscr{E}^{-1} \left(\mi \left(
\tilde{\mathbb{K}}_{11} \tilde{\mathbb{M}}_{12}(z) \! + \! \tilde{\mathbb{
K}}_{12} \tilde{\mathbb{M}}_{22}(z) \right) \Theta_{0,0}^{0,0,-}(\tilde{\Phi}_{
\tilde{b}_{j-1}}(z)) \! + \! \left(\tilde{\mathbb{K}}_{11} \tilde{\mathbb{M}}_{11}
(z) \! + \! \tilde{\mathbb{K}}_{12} \tilde{\mathbb{M}}_{21}(z) \right) \right. 
\nonumber \\
\times&\left. \, \Theta_{0,0}^{0,0,+}(\tilde{\Phi}_{\tilde{b}_{j-1}}(z)) 
\me^{\mi ((n-1)K+k) \tilde{\mho}_{j-1}} \right), \label{eqmainfin112} \\
\tilde{\mathcal{X}}_{12}^{\tilde{b},4}(z) =& \, \sqrt{\pi} \me^{-\frac{\mi \pi}{6}} 
\me^{-\frac{1}{2}((n-1)K+k) \tilde{\xi}_{\tilde{b}_{j-1}}(z)} \mathscr{E}^{-1} 
\left(\left(\tilde{\mathbb{K}}_{11} \tilde{\mathbb{M}}_{11}(z) \! + \! \tilde{
\mathbb{K}}_{12} \tilde{\mathbb{M}}_{21}(z) \right) \Theta_{0,1}^{2,1,+}
(\tilde{\Phi}_{\tilde{b}_{j-1}}(z)) \! + \! \mi \left(\tilde{\mathbb{K}}_{11} 
\tilde{\mathbb{M}}_{12}(z) \! + \! \tilde{\mathbb{K}}_{12} \tilde{\mathbb{
M}}_{22}(z) \right) \right. \nonumber \\
\times&\left. \, \Theta_{0,1}^{2,1,-}(\tilde{\Phi}_{\tilde{b}_{j-1}}(z)) 
\me^{-\mi ((n-1)K+k) \tilde{\mho}_{j-1}} \right), \label{eqmainfin113} \\
\tilde{\mathcal{X}}_{21}^{\tilde{b},4}(z) =& \, -\mi \sqrt{\pi} \me^{\frac{1}{2} 
((n-1)K+k) \tilde{\xi}_{\tilde{b}_{j-1}}(z)} \mathscr{E} \left(\mi \left(\tilde{
\mathbb{K}}_{21} \tilde{\mathbb{M}}_{12}(z) \! + \! \tilde{\mathbb{K}}_{22} 
\tilde{\mathbb{M}}_{22}(z) \right) \Theta_{0,0}^{0,0,-}(\tilde{\Phi}_{\tilde{b}_{j-1}}
(z)) \! + \! \left(\tilde{\mathbb{K}}_{21} \tilde{\mathbb{M}}_{11}(z) \! + \! \tilde{
\mathbb{K}}_{22} \tilde{\mathbb{M}}_{21}(z) \right) \right. \nonumber \\
\times&\left. \, \Theta_{0,0}^{0,0,+}(\tilde{\Phi}_{\tilde{b}_{j-1}}(z)) 
\me^{\mi ((n-1)K+k) \tilde{\mho}_{j-1}} \right), \label{eqmainfin114} \\
\tilde{\mathcal{X}}_{22}^{\tilde{b},4}(z) =& \, \sqrt{\pi} \me^{-\frac{\mi \pi}{6}} 
\me^{-\frac{1}{2} ((n-1)K+k) \tilde{\xi}_{\tilde{b}_{j-1}}(z)} \mathscr{E} \left(
\left(\tilde{\mathbb{K}}_{21} \tilde{\mathbb{M}}_{11}(z) \! + \! \tilde{\mathbb{
K}}_{22} \tilde{\mathbb{M}}_{21}(z) \right) \Theta_{0,1}^{2,1,+}(\tilde{\Phi}_{
\tilde{b}_{j-1}}(z)) \! + \! \mi \left(\tilde{\mathbb{K}}_{21} \tilde{\mathbb{M}}_{
12}(z) \! + \! \tilde{\mathbb{K}}_{22} \tilde{\mathbb{M}}_{22}(z) \right) \right. 
\nonumber \\
\times&\left. \, \Theta_{0,1}^{2,1,-}(\tilde{\Phi}_{\tilde{b}_{j-1}}(z)) 
\me^{-\mi ((n-1)K+k) \tilde{\mho}_{j-1}} \right); \label{eqmainfin115}
\end{align}
and {\rm \pmb{(12)}} for $z \! \in \! \tilde{\Omega}^{4}_{\tilde{a}_{j}}$, $j \! 
= \! 1,2,\dotsc,N \! + \! 1$,
\begin{align}
(z \! - \! \alpha_{k}) \pmb{\pi}_{k}^{n}(z) \underset{\underset{z_{o}=
1+o(1)}{\mathscr{N},n \to \infty}}{=}& \, \left(\tilde{\mathcal{X}}_{11}^{\tilde{a},4}
(z) \left(1 \! + \! \dfrac{1}{(n \! - \! 1)K \! + \! k} \left(\tilde{\mathcal{R}}_{11}^{
\sharp}(z) \! - \! \tilde{\mathcal{R}}_{11}^{\natural}(z) \right) \! + \! \mathcal{O} 
\left(\dfrac{\mathfrak{c}(n,k,z_{o})}{((n \! - \! 1)K \! + \! k)^{2}} \right) \right) 
\right. \nonumber \\
+&\left. \, \tilde{\mathcal{X}}_{21}^{\tilde{a},4}(z) \left(\dfrac{1}{(n \! - \! 1)K \! 
+ \! k} \left(\tilde{\mathcal{R}}_{12}^{\sharp}(z) \! - \! \tilde{\mathcal{R}}_{12}^{
\natural}(z) \right) \! + \! \mathcal{O} \left(\dfrac{\mathfrak{c}(n,k,z_{o})}{((n 
\! - \! 1)K \! + \! k)^{2}} \right) \right) \right) \nonumber \\
\times& \, \mathscr{E}^{-1} \me^{n(g^{f}(z)-\hat{\mathscr{P}}_{0}^{-})}, 
\label{eqmainfin116}
\end{align}
and
\begin{align}
(z \! - \! \alpha_{k}) \int_{\mathbb{R}} \dfrac{((\xi \! - \! \alpha_{k}) 
\pmb{\pi}_{k}^{n}(\xi)) \me^{-n \widetilde{V}(\xi)}}{(\xi \! - \! \alpha_{k})
(\xi \! - \! z)} \, \dfrac{\md \xi}{2 \pi \mi} \underset{\underset{z_{o}=1+o(1)}{
\mathscr{N},n \to \infty}}{=}& \, \left(\tilde{\mathcal{X}}_{12}^{\tilde{a},4}(z) 
\left(1 \! + \! \dfrac{1}{(n \! - \! 1)K \! + \! k} \left(\tilde{\mathcal{R}}_{11}^{
\sharp}(z) \! - \! \tilde{\mathcal{R}}_{11}^{\natural}(z) \right) \! + \! \mathcal{O} 
\left(\dfrac{\mathfrak{c}(n,k,z_{o})}{((n \! - \! 1)K \! + \! k)^{2}} \right) \right) 
\right. \nonumber \\
+&\left. \, \tilde{\mathcal{X}}_{22}^{\tilde{a},4}(z) \left(\dfrac{1}{(n \! - \! 1)K \! 
+ \! k} \left(\tilde{\mathcal{R}}_{12}^{\sharp}(z) \! - \! \tilde{\mathcal{R}}_{12}^{
\natural}(z) \right) \! + \! \mathcal{O} \left(\dfrac{\mathfrak{c}(n,k,z_{o})}{((n 
\! - \! 1)K \! + \! k)^{2}} \right) \right) \right) \nonumber \\
\times& \, \mathscr{E} \me^{n \tilde{\ell}} \me^{-n(g^{f}(z)-
\hat{\mathscr{P}}_{0}^{-})}, \label{eqmainfin117}
\end{align}
where
\begin{align}
\tilde{\mathcal{X}}_{11}^{\tilde{a},4}(z) =& \, -\mi \sqrt{\pi} \me^{\frac{1}{2} 
((n-1)K+k) \tilde{\xi}_{\tilde{a}_{j}}(z)} \mathscr{E}^{-1} \left(\mi \left(\tilde{
\mathbb{K}}_{11} \tilde{\mathbb{M}}_{12}(z) \! + \! \tilde{\mathbb{K}}_{12} 
\tilde{\mathbb{M}}_{22}(z) \right) \Theta_{0,0}^{0,0,-}(\tilde{\Phi}_{\tilde{a}_{j}}
(z)) \! - \! \left(\tilde{\mathbb{K}}_{11} \tilde{\mathbb{M}}_{11}(z) \! + \! 
\tilde{\mathbb{K}}_{12} \tilde{\mathbb{M}}_{21}(z) \right) \right. \nonumber \\
\times&\left. \, \Theta_{0,0}^{0,0,+}(\tilde{\Phi}_{\tilde{a}_{j}}(z)) 
\me^{\mi ((n-1)K+k) \tilde{\mho}_{j}} \right), \label{eqmainfin118} \\
\tilde{\mathcal{X}}_{12}^{\tilde{a},4}(z) =& \, \sqrt{\pi} \me^{-\frac{\mi \pi}{6}} 
\me^{-\frac{1}{2}((n-1)K+k) \tilde{\xi}_{\tilde{a}_{j}}(z)} \mathscr{E}^{-1} \left(
\left(\tilde{\mathbb{K}}_{11} \tilde{\mathbb{M}}_{11}(z) \! + \! \tilde{\mathbb{
K}}_{12} \tilde{\mathbb{M}}_{21}(z) \right) \Theta_{0,1}^{2,1,+}(\tilde{\Phi}_{
\tilde{a}_{j}}(z)) \! - \! \mi \left(\tilde{\mathbb{K}}_{11} \tilde{\mathbb{M}}_{12}
(z) \! + \! \tilde{\mathbb{K}}_{12} \tilde{\mathbb{M}}_{22}(z) \right) \right. 
\nonumber \\
\times&\left. \, \Theta_{0,1}^{2,1,-}(\tilde{\Phi}_{\tilde{a}_{j}}(z)) 
\me^{-\mi ((n-1)K+k) \tilde{\mho}_{j}} \right), \label{eqmainfin119} \\
\tilde{\mathcal{X}}_{21}^{\tilde{a},4}(z) =& \, -\mi \sqrt{\pi} \me^{\frac{1}{2} 
((n-1)K+k) \tilde{\xi}_{\tilde{a}_{j}}(z)} \mathscr{E} \left(\mi \left(\tilde{
\mathbb{K}}_{21} \tilde{\mathbb{M}}_{12}(z) \! + \! \tilde{\mathbb{K}}_{22} 
\tilde{\mathbb{M}}_{22}(z) \right) \Theta_{0,0}^{0,0,-}(\tilde{\Phi}_{\tilde{a}_{j}}
(z)) \! - \! \left(\tilde{\mathbb{K}}_{21} \tilde{\mathbb{M}}_{11}(z) \! + \! 
\tilde{\mathbb{K}}_{22} \tilde{\mathbb{M}}_{21}(z) \right) \right. \nonumber \\
\times&\left. \, \Theta_{0,0}^{0,0,+}(\tilde{\Phi}_{\tilde{a}_{j}}(z)) 
\me^{\mi ((n-1)K+k) \tilde{\mho}_{j}} \right), \label{eqmainfin120} \\
\tilde{\mathcal{X}}_{22}^{\tilde{a},4}(z) =& \, \sqrt{\pi} \me^{-\frac{\mi \pi}{6}} 
\me^{-\frac{1}{2}((n-1)K+k) \tilde{\xi}_{\tilde{a}_{j}}(z)} \mathscr{E} \left(\left(
\tilde{\mathbb{K}}_{21} \tilde{\mathbb{M}}_{11}(z) \! + \! \tilde{\mathbb{K}}_{22} 
\tilde{\mathbb{M}}_{21}(z) \right) \Theta_{0,1}^{2,1,+}(\tilde{\Phi}_{\tilde{a}_{j}}
(z)) \! - \! \mi \left(\tilde{\mathbb{K}}_{21} \tilde{\mathbb{M}}_{12}(z) \! + \! 
\tilde{\mathbb{K}}_{22} \tilde{\mathbb{M}}_{22}(z) \right) \right. \nonumber \\
\times&\left. \, \Theta_{0,1}^{2,1,-}(\tilde{\Phi}_{\tilde{a}_{j}}(z)) 
\me^{-\mi ((n-1)K+k) \tilde{\mho}_{j}} \right). \label{eqmainfin121}
\end{align}
\end{dddd}
\begin{eeee} \label{remvalfin} 
\textsl{For $n \! \in \! \mathbb{N}$ and $k \! \in \! \lbrace 1,2,\dotsc,K 
\rbrace$ such that $\alpha_{p_{\mathfrak{s}}} \! := \! \alpha_{k} \! \neq 
\! \infty$, using limiting values, if necessary, the asymptotics, in the 
double-scaling $\mathscr{N},n \! \to \! \infty$ such that $z_{o} \! = \! 1 
\! + \! o(1)$, in Theorem~\ref{maintheoforfin1} for $(z \! - \! \alpha_{k}) 
\pmb{\pi}_{k}^{n}(z)$ and $(z \! - \! \alpha_{k}) \int_{\mathbb{R}}((\xi \! 
- \! \alpha_{k}) \pmb{\pi}_{k}^{n}(\xi))((\xi \! - \! \alpha_{k})(\xi \! - \! 
z))^{-1} \me^{-n \widetilde{V}(\xi)} \, \md \xi/2 \pi \mi$ have a natural 
interpretation on $\mathbb{R} \setminus \lbrace \alpha_{1},\alpha_{2},
\dotsc,\alpha_{K} \rbrace$.}
\end{eeee}
\begin{eeee} \label{remfortheobelow2} 
\textsl{The bulk of the parameters appearing in 
Theorems~\ref{maintheompafin} and~\ref{maintheoforfin2} below 
have been defined heretofore in Theorem~\ref{maintheoforfin1}.}
\end{eeee}
\begin{dddd} \label{maintheompafin} 
Let the external field $\widetilde{V} \colon \overline{\mathbb{R}} 
\setminus \lbrace \alpha_{1},\alpha_{2},\dotsc,\alpha_{K} \rbrace \! 
\to \! \mathbb{R}$ satisfy conditions~\eqref{eq20}--\eqref{eq22}, and 
suppose that $\widetilde{V}$ is regular. For $n \! \in \! \mathbb{N}$ 
and $k \! \in \! \lbrace 1,2,\dotsc,K \rbrace$ such that $\alpha_{
p_{\mathfrak{s}}} \! := \! \alpha_{k} \! \neq \! \infty$, let 
$\mathcal{X} \colon \overline{\mathbb{C}} \setminus \overline{
\mathbb{R}} \! \to \! \mathrm{SL}_{2}(\mathbb{C})$ be the unique 
solution of the corresponding monic {\rm MPC ORF RHP} 
$(\mathcal{X}(z),\upsilon (z),\overline{\mathbb{R}})$ stated in 
Lemma~$\bm{\mathrm{RHP}_{\mathrm{MPC}}}$ with integral 
representation given by Equation~\eqref{intrepfin}, where, in 
particular,
\begin{equation*}
(\mathcal{X}(z))_{11} \! = \! (z \! - \! \alpha_{k}) \pmb{\pi}_{k}^{n}(z) 
\quad \quad \, \text{and} \, \quad \quad (\mathcal{X}(z))_{12} \! 
= \! (z \! - \! \alpha_{k}) \int_{\mathbb{R}} \dfrac{((\xi \! - \! 
\alpha_{k}) \pmb{\pi}_{k}^{n}(\xi)) \me^{-n \widetilde{V}(\xi)}}{(\xi 
\! - \! \alpha_{k})(\xi \! - \! z)} \, \dfrac{\md \xi}{2 \pi \mi}.
\end{equation*}
For $n \! \in \! \mathbb{N}$ and $k \! \in \! \lbrace 1,2,\dotsc,K 
\rbrace$ such that $\alpha_{p_{\mathfrak{s}}} \! := \! \alpha_{k} 
\! \neq \! \infty$, let the Markov-Stieltjes transform, $\mathrm{F}_{
\tilde{\mu}}(z)$, the associated $\mathrm{R}$-function, 
$\widetilde{\pmb{\mathrm{R}}}_{\tilde{\mu}}(z)$, and the corresponding 
{\rm MPA} error term, $\widetilde{\pmb{\mathrm{E}}}_{\tilde{\mu}}
(z)$, be defined by Equations~\eqref{eqmvsstildemu}, 
\eqref{eqlemmvssfinmpa1}, and~\eqref{eqlemmvssfinmpa7}, 
respectively,\footnote{Note that (cf. Subsections~\ref{subsubsec1.2.1} 
and~\ref{subsubsec1.2.2}, Equations~\eqref{mvssinf1}, \eqref{mvssfin3}, 
and~\eqref{mvssfin9}) $\mathrm{F}_{\tilde{\mu}}(z)$, 
$\widetilde{\pmb{\mathrm{R}}}_{\tilde{\mu}}(z)$, and 
$\widetilde{\pmb{\mathrm{E}}}_{\tilde{\mu}}(z)$ are $\mathrm{F}_{\mu}
(z)$, $\widetilde{\pmb{\mathrm{R}}}_{\mu}(z)$, and 
$\widetilde{\pmb{\mathrm{E}}}_{\mu}(z)$, respectively, under the 
transformation (cf. Remark~\ref{rem1.3.2}) $\md \mu \! \to \! \md 
\widetilde{\mu}$.} where, in particular,
\begin{equation} \label{eqlemmpaerrorfin1} 
\widetilde{\pmb{\mathrm{E}}}_{\tilde{\mu}}(z) \! = \! 2 \pi 
\mi \dfrac{(\mathcal{X}(z))_{12}}{(\mathcal{X}(z))_{11}}.
\end{equation} 
Then, for $n \! \in \! \mathbb{N}$ and $k \! \in \! \lbrace 1,2,\dotsc,
K \rbrace$ such that $\alpha_{p_{\mathfrak{s}}} \! := \! \alpha_{k} 
\! \neq \! \infty$, in the double-scaling limit $\mathscr{N},n \! \to 
\! \infty$ such that $z_{o} \! = \! 1 \! + \! o(1)$, $\widetilde{
\pmb{\mathrm{E}}}_{\tilde{\mu}}(z)$ has asymptotics derivable 
{}from Equation~\eqref{eqlemmpaerrorfin1}, where (cf. 
Theorem~\ref{maintheoforfin1}$)$$:$ {\rm \pmb{(1)}} for $z \! \in 
\! \tilde{\Upsilon}_{1}$, $(\mathcal{X}(z))_{11}$ and $(\mathcal{X}
(z))_{12}$ have, respectively, asymptotics~\eqref{eqmainfin11} 
and~\eqref{eqmainfin12}$;$ {\rm \pmb{(2)}} for $z \! \in \! 
\tilde{\Upsilon}_{2}$, $(\mathcal{X}(z))_{11}$ and $(\mathcal{X}
(z))_{12}$ have, respectively, asymptotics~\eqref{eqmainfin61} 
and~\eqref{eqmainfin62}$;$ {\rm \pmb{(3)}} for $z \! \in \! 
\tilde{\Upsilon}_{3}$, $(\mathcal{X}(z))_{11}$ and $(\mathcal{X}
(z))_{12}$ have, respectively, asymptotics~\eqref{eqmainfin63} 
and~\eqref{eqmainfin64}$;$ {\rm \pmb{(4)}} for $z \! \in \! 
\tilde{\Upsilon}_{4}$, $(\mathcal{X}(z))_{11}$ and $(\mathcal{X}
(z))_{12}$ have, respectively, asymptotics~\eqref{eqmainfin65} 
and~\eqref{eqmainfin66}$;$ {\rm \pmb{(5)}} for $z \! \in \! 
\tilde{\Omega}^{1}_{\tilde{b}_{j-1}}$, $j \! = \! 1,2,\dotsc,N \! + \! 1$, 
$(\mathcal{X}(z))_{11}$ and $(\mathcal{X}(z))_{12}$ have, respectively, 
asymptotics~\eqref{eqmainfin67} and~\eqref{eqmainfin68}$;$ 
{\rm \pmb{(6)}} for $z \! \in \! \tilde{\Omega}^{1}_{\tilde{a}_{j}}$, 
$j \! = \! 1,2,\dotsc,N \! + \! 1$, $(\mathcal{X}(z))_{11}$ 
and $(\mathcal{X}(z))_{12}$ have, respectively, 
asymptotics~\eqref{eqmainfin79} and~\eqref{eqmainfin80}$;$ 
{\rm \pmb{(7)}} for $z \! \in \! \tilde{\Omega}^{2}_{\tilde{b}_{j-1}}$, 
$j \! = \! 1,2,\dotsc,N \! + \! 1$, $(\mathcal{X}(z))_{11}$ and 
$(\mathcal{X}(z))_{12}$ have, respectively, 
asymptotics~\eqref{eqmainfin86} and~\eqref{eqmainfin87}$;$ 
{\rm \pmb{(8)}} for $z \! \in \! \tilde{\Omega}^{2}_{\tilde{a}_{j}}$, 
$j \! = \! 1,2,\dotsc,N \! + \! 1$, $(\mathcal{X}(z))_{11}$ 
and $(\mathcal{X}(z))_{12}$ have, respectively, 
asymptotics~\eqref{eqmainfin92} and~\eqref{eqmainfin93}$;$ 
{\rm \pmb{(9)}} for $z \! \in \! \tilde{\Omega}^{3}_{\tilde{b}_{j-1}}$, 
$j \! = \! 1,2,\dotsc,N \! + \! 1$, $(\mathcal{X}(z))_{11}$ 
and $(\mathcal{X}(z))_{12}$ have, respectively, 
asymptotics~\eqref{eqmainfin98} and~\eqref{eqmainfin99}$;$ 
{\rm \pmb{(10)}} for $z \! \in \! \tilde{\Omega}^{3}_{\tilde{a}_{j}}$, 
$j \! = \! 1,2,\dotsc,N \! + \! 1$, $(\mathcal{X}(z))_{11}$ 
and $(\mathcal{X}(z))_{12}$ have, respectively, 
asymptotics~\eqref{eqmainfin104} and~\eqref{eqmainfin105}$;$ 
{\rm \pmb{(11)}} for $z \! \in \! \tilde{\Omega}^{4}_{\tilde{b}_{j-1}}$, 
$j \! = \! 1,2,\dotsc,N \! + \! 1$, $(\mathcal{X}(z))_{11}$ 
and $(\mathcal{X}(z))_{12}$ have, respectively, 
asymptotics~\eqref{eqmainfin110} and~\eqref{eqmainfin111}$;$ and 
{\rm \pmb{(12)}} for $z \! \in \! \tilde{\Omega}^{4}_{\tilde{a}_{j}}$, 
$j \! = \! 1,2,\dotsc,N \! + \! 1$, $(\mathcal{X}(z))_{11}$ 
and $(\mathcal{X}(z))_{12}$ have, respectively, 
asymptotics~\eqref{eqmainfin116} and~\eqref{eqmainfin117}.
\end{dddd}
\begin{eeee} \label{remmpaasyfin}
\textsl{It is shown in Section~\ref{sec3} (see, in particular, the corresponding 
items of Lemmata~\ref{lemrootz} and~\ref{lemetatomu}$)$ that, for $n \! 
\in \! \mathbb{N}$ and $k \! \in \! \lbrace 1,2,\dotsc,K \rbrace$ such that 
$\alpha_{p_{\mathfrak{s}}} \! := \! \alpha_{k} \! \neq \! \infty$, $\lbrace 
\mathstrut z \! \in \! \overline{\mathbb{C}}; \, \pmb{\pi}^{n}_{k}(z) \! = \! 0 
\rbrace \! =: \! \lbrace \tilde{\mathfrak{z}}^{n}_{k}(j) \rbrace_{j=1}^{(n-1)K+k} 
\! \subset \! \cup_{j=1}^{N+1}[\tilde{b}_{j-1},\tilde{a}_{j}] \! =: \! J_{f}$$;$ 
therefore, using limiting values, if necessary, for $\mathbb{C} \cap 
(\tilde{\Upsilon}_{3} \cup \tilde{\Upsilon}_{4} \cup \cup_{j=1}^{N+1}
(\tilde{\Omega}^{2}_{\tilde{b}_{j-1}} \cup \tilde{\Omega}^{3}_{\tilde{b}_{j-1}} 
\cup \tilde{\Omega}^{2}_{\tilde{a}_{j}} \cup \tilde{\Omega}^{3}_{\tilde{a}_{j}})) 
\! \ni \! z \! \to \! x \! \in \! J_{f}$, the asymptotics, in the double-scaling 
limit $\mathscr{N},n \! \to \! \infty$ such that $z_{o} \! = \! 1 \! + \! o(1)$, in 
Theorem~\ref{maintheompafin} for $\widetilde{\pmb{\mathrm{E}}}_{\tilde{\mu}}(z)$ 
have a natural interpretation on $J_{f} \setminus \lbrace \tilde{\mathfrak{z}}^{n}_{k}
(j) \rbrace_{j=1}^{(n-1)K+k}$.}
\end{eeee}
\begin{dddd} \label{maintheoforfin2} 
Let the external field $\widetilde{V} \colon \overline{\mathbb{R}} 
\setminus \lbrace \alpha_{1},\alpha_{2},\dotsc,\alpha_{K} \rbrace \! 
\to \! \mathbb{R}$ satisfy conditions~\eqref{eq20}--\eqref{eq22}, and 
suppose that $\widetilde{V}$ is regular. For $n \! \in \! \mathbb{N}$ 
and $k \! \in \! \lbrace 1,2,\dotsc,K \rbrace$ such that $\alpha_{
p_{\mathfrak{s}}} \! := \! \alpha_{k} \! \neq \! \infty$, let $\mathcal{X} 
\colon \overline{\mathbb{C}} \setminus \overline{\mathbb{R}} 
\! \to \! \mathrm{SL}_{2}(\mathbb{C})$ be the unique solution 
of the corresponding monic {\rm MPC ORF RHP} 
$(\mathcal{X}(z),\upsilon (z),\overline{\mathbb{R}})$ stated in 
Lemma~$\bm{\mathrm{RHP}_{\mathrm{MPC}}}$, and let the monic 
{\rm MPC ORF} have the representation $\pmb{\pi}_{k}^{n}(z) \! 
= \! (z \! - \! \alpha_{k})^{-1}(\mathcal{X}(z))_{11}$, $z \! \in \! 
\mathbb{C}$, with asymptotics, in the double-scaling limit 
$\mathscr{N},n \! \to \! \infty$ such that $z_{o} \! = \! 1 \! + \! 
o(1)$, stated in Theorem~\ref{maintheoforfin1}$;$ moreover, let the 
associated norming constant be given by $\mu_{n,\varkappa_{nk}}^{f}
(n,k) \! = \! \lvert \lvert \pmb{\pi}_{k}^{n}(\pmb{\cdot}) \rvert 
\rvert_{\mathscr{L}}^{-1}$. Then, for $n \! \in \! \mathbb{N}$ 
and $k \! \in \! \lbrace 1,2,\dotsc,K \rbrace$ such that 
$\alpha_{p_{\mathfrak{s}}} \! := \! \alpha_{k} \! \neq \! 
\infty$,\footnote{The asymptotics for $\mu_{n,\varkappa_{nk}}^{f}(n,k)$ 
is obtained by taking the positive square root of the expression on the 
right-hand side of Equation~\eqref{eqcoeffinf1}.}
\begin{equation} \label{eqcoeffinf1} 
(\mu_{n,\varkappa_{nk}}^{f}(n,k))^{2} \underset{\underset{z_{o}=1+
o(1)}{\mathscr{N},n \to \infty}}{=} \dfrac{\me^{-n \tilde{\ell}}}{2 \pi} 
\tilde{\mathbb{X}}^{\sharp} \left(1 \! + \! \dfrac{1}{(n \! - \! 1)K \! + \!  k} 
\tilde{\mathbb{X}}^{\sharp} \left(\tilde{\mathbb{X}}^{\natural} \right)_{12} 
\! + \! \mathcal{O} \left(\dfrac{\tilde{\mathfrak{c}}(n,k,z_{o})}{((n \! - \! 1)
K \! + \! k)^{2}} \right) \right),
\end{equation}
where
\begin{equation} \label{eqcoeffinf2} 
\tilde{\mathbb{X}}^{\sharp} \! = \! \dfrac{\mathscr{E}^{2}}{\mi \tilde{
\mathbb{K}}_{11} \mathbb{H}(1) \! + \! \tilde{\mathbb{K}}_{12} 
\mathbb{H}(-1)},
\end{equation}
with
\begin{align} \label{eqcoeffinf3} 
\mathbb{H}(\varepsilon_{2}) :=& \, \dfrac{\mi}{2} \dfrac{\tilde{\boldsymbol{
\theta}}(-\tilde{\boldsymbol{u}}_{+}(\alpha_{k}) \! - \! \frac{1}{2 \pi}((n \! 
- \! 1)K \! + \! k) \tilde{\boldsymbol{\Omega}} \! + \! \varepsilon_{2} 
\tilde{\boldsymbol{d}})}{\tilde{\boldsymbol{\theta}}(-\tilde{\boldsymbol{
u}}_{+}(\alpha_{k}) \! + \! \varepsilon_{2} \tilde{\boldsymbol{d}})} \left(
\dfrac{1}{4} \left(\tilde{\gamma}(\alpha_{k}) \! + \! \varepsilon_{2}
(\tilde{\gamma}(\alpha_{k}))^{-1} \right) \sum_{j=1}^{N+1} \left(
\dfrac{1}{\tilde{a}_{j} \! - \! \alpha_{k}} \! - \! \dfrac{1}{\tilde{b}_{j-1} \! 
- \! \alpha_{k}} \right) \right. \nonumber \\
-&\left. \, \left(\tilde{\gamma}(\alpha_{k}) \! - \! \varepsilon_{2}(\tilde{
\gamma}(\alpha_{k}))^{-1} \right) \left(\mathbb{L}(\varepsilon_{2};
\tilde{\boldsymbol{\Omega}}) \! - \! \mathbb{L}(\varepsilon_{2};
\boldsymbol{0}) \right) \right), \quad \varepsilon_{2} \! = \! \pm 1,
\end{align}
where
\begin{gather}
\mathbb{L}(\varepsilon_{2};\tilde{\boldsymbol{\Omega}}) \! = \! \dfrac{1}{
\tilde{\boldsymbol{\theta}}(-\tilde{\boldsymbol{u}}_{+}(\alpha_{k}) \! - 
\! \frac{1}{2 \pi}((n \! - \! 1)K \! + \! k) \tilde{\boldsymbol{\Omega}} \! + 
\! \varepsilon_{2} \tilde{\boldsymbol{d}})} \sum_{m \in \mathbb{Z}^{N}} 
\tilde{\boldsymbol{\lambda}}(\boldsymbol{m}) \me^{2 \pi \mi (m,
-\tilde{\boldsymbol{u}}_{+}(\alpha_{k})-\frac{1}{2 \pi}((n-1)K+k) 
\tilde{\boldsymbol{\Omega}}+ \varepsilon_{2} \tilde{\boldsymbol{d}})
+\mi \pi (m,\tilde{\boldsymbol{\tau}}m)}, \label{eqcoeffinf4} \\
\tilde{\boldsymbol{\lambda}}(\boldsymbol{m}) \! := \! 2 \pi \mi 
\sum_{j_{1}=1}^{N} \sum_{j_{2}=1}^{N} \dfrac{m_{j_{1}} \tilde{c}_{j_{1}j_{2}}
(\alpha_{k})^{N-j_{2}}}{(-1)^{\tilde{\mathfrak{n}}(\alpha_{k})} \prod_{j_{3}=
1}^{N+1}(\lvert \alpha_{k} \! - \! \tilde{b}_{j_{3}-1} \rvert \lvert \alpha_{k} 
\! - \! \tilde{a}_{j_{3}} \rvert)^{1/2}}, \label{eqcoeffinf5} \\
\tilde{\mathfrak{n}}(\alpha_{k}) \! = \! 
\begin{cases}
0, &\text{$\alpha_{k} \! \in \! (\tilde{a}_{N+1},+\infty)$,} \\
N \! + \! 1, &\text{$\alpha_{k} \! \in \! (-\infty,\tilde{b}_{0})$,} \\
N \! - \! j \! + \! 1, &\text{$\alpha_{k} \! \in \! (\tilde{a}_{j},\tilde{b}_{j}), 
\quad j \! = \! 1,2,\dotsc,N$,}
\end{cases} \label{eqcoeffinf6}
\end{gather}
and $\tilde{c}_{j_{1}j_{2}}$, $j_{1},j_{2} \! = \! 1,2,\dotsc,N$, are obtained 
{}from Equations~\eqref{O1} and~\eqref{O2}, and
\begin{align}
\tilde{\mathbb{X}}^{\natural} =& \, \mi \sum_{j=1}^{N+1} \left(\dfrac{
(\tilde{\alpha}_{0}(\tilde{b}_{j-1}))^{-2}}{(\tilde{b}_{j-1} \! - \! \alpha_{k})^{2}} 
\left(\tilde{\alpha}_{0}(\tilde{b}_{j-1}) \tilde{\boldsymbol{\mathrm{B}}}
(\tilde{b}_{j-1}) \! - \! \tilde{\alpha}_{1}(\tilde{b}_{j-1}) \tilde{\boldsymbol{
\mathrm{A}}}(\tilde{b}_{j-1}) \right) \! + \! \dfrac{(\tilde{\alpha}_{0}
(\tilde{a}_{j}))^{-2}}{(\tilde{a}_{j} \! - \! \alpha_{k})^{2}} \left(\tilde{
\alpha}_{0}(\tilde{a}_{j}) \tilde{\boldsymbol{\mathrm{B}}}(\tilde{a}_{j}) 
\! - \! \tilde{\alpha}_{1}(\tilde{a}_{j}) \tilde{\boldsymbol{\mathrm{A}}}
(\tilde{a}_{j}) \right) \right. \nonumber \\
-&\left. \, \dfrac{2(\tilde{\alpha}_{0}(\tilde{b}_{j-1}))^{-1}}{(\tilde{b}_{j-1} 
\! - \! \alpha_{k})^{3}} \tilde{\boldsymbol{\mathrm{A}}}(\tilde{b}_{j-1}) \! 
- \! \dfrac{2(\tilde{\alpha}_{0}(\tilde{a}_{j}))^{-1}}{(\tilde{a}_{j} \! - \! 
\alpha_{k})^{3}} \tilde{\boldsymbol{\mathrm{A}}}(\tilde{a}_{j}) \right), 
\label{eqcoeffinf7}
\end{align}
with $\tilde{\mathfrak{c}}(n,k,z_{o}) \! =_{\underset{z_{o}=1+o(1)}{
\mathscr{N},n \to \infty}} \! \mathcal{O}(1)$.

For $n \! \in \! \mathbb{N}$ and $k \! \in \! \lbrace 1,2,\dotsc,K 
\rbrace$ such that $\alpha_{p_{\mathfrak{s}}} \! := \! \alpha_{k} 
\! \neq \! \infty$, let the associated {\rm MPC ORF} be given 
by $\phi_{k}^{n}(z) \! = \! \mu_{n,\varkappa_{nk}}^{f}(n,k) 
\pmb{\pi}_{k}^{n}(z)$, $z \! \in \! \mathbb{C}$$:$ the asymptotics, 
in the double-scaling limit $\mathscr{N},n \! \to \! \infty$ such that 
$z_{o} \! = \! 1 \! + \! o(1)$, for $\phi_{k}^{n}(z)$ are derived via the 
scalar multiplication of the asymptotics for $\pmb{\pi}_{k}^{n}(z)$ 
and $\mu_{n,\varkappa_{nk}}^{f}(n,k)$ stated, respectively, in 
Theorem~\ref{maintheoforfin1} and Equation~\eqref{eqcoeffinf1}.
\end{dddd}
\section{The Monic MPC ORF Families of RHPs: Existence, Uniqueness, and 
Multi-Integral Representations} \label{sec2} 
In this section, the family (which consists of two subfamilies) of $K$ 
${}$\footnote{Recall that $K \! = \! \hat{K} \! + \! \tilde{K}$, where 
$\hat{K} \! := \! \# \lbrace \mathstrut k \! \in \! \lbrace 1,2,\dotsc,K 
\rbrace; \, \alpha_{k} \! = \! \infty \rbrace$, and $\tilde{K} \! := \! \# 
\lbrace \mathstrut k \! \in \! \lbrace 1,2,\dotsc,K \rbrace; \, \alpha_{k} 
\! \neq \! \infty \rbrace$.} matrix RHPs on $\overline{\mathbb{R}}$ 
characterising the monic MPC ORFs, $\lbrace \pmb{\pi}^{n}_{k}(z) 
\rbrace_{\underset{k=1,2,\dotsc,K}{n \in \mathbb{N}}}$, 
$z \! \in \! \mathbb{C}$, is stated (see 
Lemma~$\bm{\mathrm{RHP}_{\mathrm{MPC}}}$), whence the existence 
and the uniqueness of the associated monic MPC ORFs and the 
corresponding norming constants, $\mu^{r}_{n,\varkappa_{nk}}
(n,k)$, $r \! \in \! \lbrace \infty,f \rbrace$, is established via 
a---generalised---Hankel determinant analysis (see Lemma~\ref{lem2.1}), 
and explicit multi-integral representations for the associated monic 
MPC ORFs (see Lemma~\ref{lem2.2}) and the corresponding norming 
constants (see Corollary~\ref{cor2.1}) are obtained (see, also, 
Corollary~\ref{cor2.2}).

Before launching into the theory proper, however, it is convenient to 
summarise the notation used throughout this monograph:
\begin{enumerate}
\item[(1)] $\mathrm{I} \! = \! 
\left(
\begin{smallmatrix}
1 & 0 \\
0 & 1
\end{smallmatrix}
\right)$ is the $2 \times 2$ identity matrix, $\sigma_{1} \! = \! 
\left(
\begin{smallmatrix}
0 & 1 \\
1 & 0
\end{smallmatrix}
\right)$, $\sigma_{2} \! = \! 
\left(
\begin{smallmatrix}
0 & -\mi \\
\mi & 0
\end{smallmatrix}
\right)$, and $\sigma_{3} \! = \! 
\left(
\begin{smallmatrix}
1 & 0 \\
0 & -1
\end{smallmatrix}
\right)$ are the Pauli matrices, $\sigma_{\pm} \! := \! \tfrac{1}{2}(\sigma_{1} 
\! \pm \! \mi \sigma_{2})$, $\mathbb{R}_{\pm} \! := \! \lbrace \mathstrut x 
\! \in \! \mathbb{R}; \, \pm x \! > \! 0 \rbrace$, $\mathbb{R}_{x_{0}}^{\gtrless} 
\! := \! \lbrace \mathstrut x \! \in \! \mathbb{R}; \, x \! \gtrless \! x_{0} 
\rbrace$, $\mathbb{C}_{\pm} \! := \! \lbrace \mathstrut z \! \in \! 
\mathbb{C}; \, \pm \Im (z) \! > \! 0 \rbrace$, and $\mathrm{sgn}(x) \! = \! 
\left\{
\begin{smallmatrix}
x \lvert x \rvert^{-1}, & x \neq 0, \\
0, & x=0;
\end{smallmatrix}
\right.$
\item[(2)] for $\omega \! \in \! \mathbb{C}$ and $\hat{\Upsilon} \! \in \! 
\operatorname{M}_{2}(\mathbb{C})$, $\omega^{\mathrm{ad}(\sigma_{3})} 
\hat{\Upsilon} \! := \! \omega^{\sigma_{3}} \hat{\Upsilon} \omega^{-\sigma_{3}}$;
\item[(3)] a contour $\mathcal{D}$ which is the finite union of piecewise-smooth, 
simple curves (as closed sets) is said to be \emph{orientable} if its complement 
$\overline{\mathbb{C}} \setminus \mathcal{D}$ can always be divided into two, 
possibly disconnected, disjoint open sets $\mho^{+}$ and $\mho^{-}$, either 
of which has finitely many components, such that $\mathcal{D}$ admits an 
orientation so that it can either be viewed as a positively oriented boundary 
$\mathcal{D}^{+}$ for $\mho^{+}$ or as a negatively oriented boundary 
$\mathcal{D}^{-}$ for $\mho^{-}$ \cite{a63}, that is, the (possibly disconnected) 
components of $\overline{\mathbb{C}} \setminus \mathcal{D}$ can be coloured 
by $+$ or by $-$ in such a way that the $+$ regions do not share boundary 
with the $-$ regions, except, possibly, at finitely many points \cite{a64};
\item[(4)] for each segment of an oriented contour $\mathcal{D}$, according 
to the given orientation, the ``+'' side is to the left and the ``-'' side 
is to the right as one traverses the contour in the direction of orientation, 
that is, for a $2 \times 2$ matrix $\mathcal{A}_{ij}(\pmb{\cdot})$, $i,j \! = 
\! 1,2$, $(\mathcal{A}_{ij}(\pmb{\cdot}))_{\pm}$ denote the non-tangential 
limits $(\mathcal{A}_{ij}(z))_{\pm} \! := \! \lim_{\genfrac{}{}{0pt}{2}{z^{\prime} 
\, \to \, z}{z^{\prime} \, \in \, \pm \, \mathrm{side} \, \mathrm{of} \, 
\mathcal{D}}} \mathcal{A}_{ij}(z^{\prime})$;
\item[(5)] for $1 \! \leqslant \! p \! < \! \infty$ and $\mathcal{D}$ 
some point set,
\begin{equation*}
\mathcal{L}^{p}_{\mathrm{M}_{2}(\mathbb{C})}(\mathcal{D}) \! := \! \left\{
\mathstrut f \colon \mathcal{D} \! \to \! \mathrm{M}_{2}(\mathbb{C}); \, 
\vert \vert f(\pmb{\cdot}) \vert \vert_{\mathcal{L}^{p}_{\mathrm{M}_{2}
(\mathbb{C})}(\mathcal{D})} \! := \! \left(\int\nolimits_{\mathcal{D}} \vert 
f(z) \vert^{p} \, \vert \md z \vert \right)^{1/p} \! < \! \infty \right\},
\end{equation*}
where, for $\mathcal{A}(\pmb{\cdot}) \! \in \! \operatorname{M}_{2}
(\mathbb{C})$, $\lvert \mathcal{A}(\pmb{\cdot}) \rvert \! := \! (\sum_{i,j=1}^{2} 
\overline{\mathcal{A}_{ij}(\pmb{\cdot})} \, \mathcal{A}_{ij}(\pmb{\cdot}))^{1/2}$ 
is the Hilbert-Schmidt norm, $\overline{\pmb{\ast}}$ denotes ${}$\footnote{Not 
to be confused with $\overline{\mathbb{C}}$ and $\overline{\mathbb{R}}$.} the 
complex conjugate of $\pmb{\ast}$, and $\lvert \md \mathbf{\bullet} \rvert$ 
denotes integration with respect to arc length, for $p \! = \! \infty$,
\begin{equation*}
\mathcal{L}^{\infty}_{\mathrm{M}_{2}(\mathbb{C})}(\mathcal{D}) \! := \! 
\left\{\mathstrut g \colon \mathcal{D} \! \to \! \mathrm{M}_{2}(\mathbb{C}); 
\, \vert \vert g(\pmb{\cdot}) \vert \vert_{\mathcal{L}^{\infty}_{\mathrm{M}_{2}
(\mathbb{C})}(\mathcal{D})} \! := \! \max_{i,j=1,2} \, \sup_{z \in \mathcal{
D}} \vert g_{ij}(z) \vert \! < \! \infty \right\},
\end{equation*}
and, for $f \! \in \! \mathrm{I} \! + \! \mathcal{L}^{2}_{\mathrm{M}_{2}
(\mathbb{C})}(\mathcal{D}) \! := \! \lbrace \mathstrut \mathrm{I} \! + \! 
h; \, h \! \in \! \mathcal{L}^{2}_{\mathrm{M}_{2}(\mathbb{C})}(\mathcal{D}) 
\rbrace$ and $k \! \in \! \lbrace 1,2,\dotsc,K \rbrace$,
\begin{equation*}
\vert \vert f(\pmb{\cdot}) \vert \vert_{\mathrm{I}+\mathcal{L}^{2}_{
\mathrm{M}_{2}(\mathbb{C})}(\mathcal{D})} \! := \! \left(\vert \vert 
f(\alpha_{k}) \vert \vert_{\mathcal{L}^{\infty}_{\mathrm{M}_{2}(\mathbb{C})}
(\mathcal{D})}^{2} \! + \! \vert \vert f(\pmb{\cdot}) \! - \! f(\alpha_{k}) \vert 
\vert_{\mathcal{L}^{2}_{\mathrm{M}_{2}(\mathbb{C})}(\mathcal{D})}^{2} 
\right)^{1/2};
\end{equation*}
\item[(6)] for a $2 \times 2$ matrix $\mathcal{A}_{ij}(\pmb{\cdot})$, $i,j \! 
= \! 1,2$, to have boundary values in the $\mathcal{L}^{2}_{\mathrm{M}_{2}
(\mathbb{C})}(\mathcal{D})$ sense on an oriented contour $\mathcal{D}$, 
it is meant that $\lim_{\genfrac{}{}{0pt}{2}{z^{\prime} \, \to \, z}{z^{\prime} 
\, \in \, \pm \, \mathrm{side} \, \mathrm{of} \, \mathcal{D}}} \int_{\mathcal{D}} 
\vert \mathcal{A}(z^{\prime}) \! - \! (\mathcal{A}(z))_{\pm} \vert^{2} \, \vert 
\md z \vert \! = \! 0$;\footnote{For example, if $\mathcal{D} \! = \! \mathbb{R}$ 
is oriented {}from $+\infty$ to $-\infty$, then $\mathcal{A}(\pmb{\cdot})$ has 
$\mathcal{L}^{2}_{\mathrm{M}_{2}(\mathbb{C})}(\mathcal{D})$ boundary values on 
$\mathcal{D}$ means that $\lim_{\varepsilon \downarrow 0} \int_{\mathbb{R}} 
\vert \mathcal{A}(x \! \mp \! \mi \varepsilon) \! - \! (\mathcal{A}(x))_{\pm} 
\vert^{2} \, \md x \! = \! 0$.}
\item[(7)] for $\mathrm{card}(J) \! < \! \infty$, $\vert \vert \mathscr{F}
(\pmb{\cdot}) \vert \vert_{\cap_{p \in J} \mathcal{L}^{p}_{\mathrm{M}_{2}
(\mathbb{C})}(\ast)} \! := \! \sum_{p \in J} \vert \vert \mathscr{F}
(\pmb{\cdot}) \vert \vert_{\mathcal{L}^{p}_{\mathrm{M}_{2}(\mathbb{C})}
(\ast)}$;
\item[(8)] for $(\nu_{1},\nu_{2}) \! \in \! \mathbb{R} \times \mathbb{R}$, 
the function $(\bullet \! - \! \nu_{1})^{\mi \nu_{2}} \colon \mathbb{C} 
\setminus (-\infty,\nu_{1}] \! \to \! \mathbb{C}$, $\bullet \! \mapsto \! 
\exp (\mi \nu_{2} \ln (\bullet -\nu_{1}))$, with the branch cut taken along 
$(-\infty,\nu_{1}]$, and with the principal branch of the logarithm chosen 
(that is, for $\nu_{0} \! \in \! \mathbb{R}$, $\ln (\nu \! - \! \nu_{0}) \! 
:= \! \ln \lvert \nu \! - \! \nu_{0} \rvert \! + \! \mi \arg (\nu \! - \! 
\nu_{0})$, where $\arg (\nu \! - \! \nu_{0}) \! \in \! (-(1 \! \mp \! 1) 
\pi/2,(1 \! \pm \! 1) \pi/2)$, $\pm \Im (\nu) \! > \! 0)$;
\item[(9)] for a $2 \times 2$ matrix-valued function $\mathfrak{T}(z)$, 
$\mathfrak{T}(z) \! =_{z \to z_{0}} \! \mathcal{O}(\ast)$ (resp., $o(\ast))$ 
means $\mathfrak{T}_{ij}(z) \! =_{z \to z_{0}} \! \mathcal{O}(\ast_{ij})$ 
(resp., $o(\ast_{ij}))$, $i,j \! = \! 1,2$.
\end{enumerate}
The family of matrix RHPs characterising the monic MPC ORFs, $\lbrace 
\pmb{\pi}^{n}_{k}(z) \rbrace_{\underset{k=1,2,\dotsc,K}{n \in \mathbb{N}}}$, 
$z \! \in \! \mathbb{C}$, is now stated.
\begin{rhp2} \label{rhpmpc} 
Let $\widetilde{V} \colon \overline{\mathbb{R}} \setminus \lbrace 
\alpha_{1},\alpha_{2},\dotsc,\alpha_{K} \rbrace \! \to \! \mathbb{R}$ satisfy 
conditions~\eqref{eq20}--\eqref{eq22}. For $n \! \in \! \mathbb{N}$ and $k 
\! \in \! \lbrace 1,2,\dotsc,K \rbrace$, find $\mathcal{X} \colon \mathbb{N} 
\times \lbrace 1,2,\dotsc,K \rbrace \times \overline{\mathbb{C}} \setminus 
\overline{\mathbb{R}} \! \to \! \operatorname{SL}_{2}(\mathbb{C})$ solving: 
{\rm (i)} $\mathcal{X}(n,k;z) \! := \! \mathcal{X}(z)$ is analytic for 
$z \! \in \! \overline{\mathbb{C}} \setminus \overline{\mathbb{R}};$ 
{\rm (ii)} the boundary values $\mathcal{X}_{\pm}(z) \! := \! 
\lim_{\varepsilon \downarrow 0} \mathcal{X}(z \! \pm \! \mi 
\varepsilon)$ satisfy the jump condition
\begin{equation*}
\mathcal{X}_{+}(z) \! = \! \mathcal{X}_{-}(z) \upsilon (z) \quad 
\mathrm{a.e.} \quad z \! \in \! \overline{\mathbb{R}},
\end{equation*}
where $\upsilon \colon \overline{\mathbb{R}} \! \ni \! z \! \mapsto \! 
\mathrm{I} \! + \! \exp (-n \widetilde{V}(z)) \sigma_{+}$$;$ {\rm (iii)} 
for $k \! \in \! \lbrace 1,2,\dotsc,K \rbrace$ such that $\alpha_{
p_{\mathfrak{s}}} \! := \! \alpha_{k} \! = \! \infty$,
\begin{gather*}
\mathcal{X}(z)z^{-\varkappa_{nk} \sigma_{3}} \underset{\overline{\mathbb{C}} 
\setminus \mathbb{R} \, \ni \, z \to \alpha_{k}}{=} \mathrm{I} \! + 
\mathcal{O}(z^{-1}), \\
\mathcal{X}(z)(z \! - \! \alpha_{p_{q}})^{\varkappa_{nk \tilde{k}_{q}} 
\sigma_{3}} \underset{\mathbb{C} \setminus \mathbb{R} \, \ni \, z 
\to \alpha_{p_{q}}}{=} \mathcal{O}(1), \quad q \! = \! 1,2,\dotsc,
\mathfrak{s} \! - \! 1;
\end{gather*}
and {\rm (iv)} for $k \! \in \! \lbrace 1,2,\dotsc,K \rbrace$ such that 
$\alpha_{p_{\mathfrak{s}}} \! := \! \alpha_{k} \! \neq \! \infty$,
\begin{gather*}
\mathcal{X}(z)(z \! - \! \alpha_{k})^{(\varkappa_{nk}-1) \sigma_{3}} 
\underset{\mathbb{C} \setminus \mathbb{R} \, \ni \, 
z \to \alpha_{k}}{=} \mathrm{I} \! + \mathcal{O}(z \! - \! \alpha_{k}), \\
\mathcal{X}(z)z^{-(\varkappa_{nk \tilde{k}_{\mathfrak{s}-1}}^{\infty}+1) 
\sigma_{3}} \underset{\overline{\mathbb{C}} \setminus \mathbb{R} \, 
\ni \, z \to \alpha_{p_{\mathfrak{s}-1}} = \infty}{=} \mathcal{O}(1), \\
\mathcal{X}(z)(z \! - \! \alpha_{p_{q}})^{\varkappa_{nk \tilde{k}_{q}} 
\sigma_{3}} \underset{\mathbb{C} \setminus \mathbb{R} \, \ni \, z 
\to \alpha_{p_{q}}}{=} \mathcal{O}(1), \quad q \! = \! 1,2,\dotsc,
\mathfrak{s} \! - \! 2.
\end{gather*}
\end{rhp2}
\begin{ccccc} \label{lem2.1} 
Let $\mathcal{X} \colon \mathbb{N} \times \lbrace 1,2,\dotsc,K \rbrace 
\times \overline{\mathbb{C}} \setminus \overline{\mathbb{R}} \! \to \! 
\operatorname{SL}_{2}(\mathbb{C})$ solve the {\rm RHP} stated in 
Lemma~{\rm $\bm{\mathrm{RHP}_{\mathrm{MPC}}}$}. For $n \! \in \! 
\mathbb{N}$ and $k \! \in \! \lbrace 1,2,\dotsc,K \rbrace$ such that 
$\alpha_{p_{\mathfrak{s}}} \! := \! \alpha_{k} \! = \! \infty$, the {\rm RHP} 
stated in Lemma~{\rm $\bm{\mathrm{RHP}_{\mathrm{MPC}}}$} 
possesses a unique solution given by
\begin{equation} \label{intrepinf} 
\mathcal{X}(z) \! = \! 
\begin{pmatrix}
\pmb{\pi}^{n}_{k}(z) & \int_{\mathbb{R}} \frac{\pmb{\pi}^{n}_{k}(\xi) 
\exp (-n \widetilde{V}(\xi))}{\xi -z} \, \frac{\md \xi}{2 \pi \mi} \\
\mathcal{X}_{21}(z) & \int_{\mathbb{R}} \frac{\mathcal{X}_{21}(\xi) 
\exp (-n \widetilde{V}(\xi))}{\xi -z} \, \frac{\md \xi}{2 \pi \mi}
\end{pmatrix}, \quad z \! \in \! \mathbb{C} \setminus \mathbb{R},
\end{equation}
where $\pmb{\pi}^{n}_{k} \colon \mathbb{N} \times \lbrace 1,2,\dotsc,K 
\rbrace \times \overline{\mathbb{C}} \setminus \lbrace \alpha_{p_{1}},
\alpha_{p_{2}},\dotsc,\alpha_{p_{\mathfrak{s}}} \rbrace \! \to \! \mathbb{C}$ 
is the corresponding monic {\rm MPC ORF} defined in 
Subsection~\ref{subsubsec1.2.1}, and $\mathcal{X}_{21} \colon \mathbb{N} 
\times \lbrace 1,2,\dotsc,K \rbrace \times \overline{\mathbb{C}} \setminus 
\lbrace \alpha_{p_{1}},\alpha_{p_{2}},\dotsc,\alpha_{p_{\mathfrak{s}}} 
\rbrace \! \to \! \mathbb{C}$$;$ and, for $n \! \in \! \mathbb{N}$ and $k 
\! \in \! \lbrace 1,2,\dotsc,K \rbrace$ such that $\alpha_{p_{\mathfrak{s}}} 
\! := \! \alpha_{k} \! \neq \! \infty$, the {\rm RHP} stated in 
Lemma~{\rm $\bm{\mathrm{RHP}_{\mathrm{MPC}}}$} possesses a unique 
solution given by
\begin{equation} \label{intrepfin} 
\mathcal{X}(z) \! = \! 
\begin{pmatrix}
(z \! - \! \alpha_{k}) \pmb{\pi}^{n}_{k}(z) & (z \! - \! \alpha_{k}) 
\int_{\mathbb{R}} \frac{((\xi -\alpha_{k}) \pmb{\pi}^{n}_{k}(\xi)) 
\exp (-n \widetilde{V}(\xi))}{(\xi - \alpha_{k})(\xi -z)} \, 
\frac{\md \xi}{2 \pi \mi} \\
\mathcal{X}_{21}(z) & (z \! - \! \alpha_{k}) \int_{\mathbb{R}} 
\frac{\mathcal{X}_{21}(\xi) \exp (-n \widetilde{V}(\xi))}{(\xi 
-\alpha_{k})(\xi -z)} \, \frac{\md \xi}{2 \pi \mi}
\end{pmatrix}, \quad z \! \in \! \mathbb{C} \setminus \mathbb{R},
\end{equation}
where 
$\pmb{\pi}^{n}_{k} \colon \mathbb{N} \times \lbrace 1,2,\dotsc,K \rbrace 
\times \overline{\mathbb{C}} \setminus \lbrace \alpha_{p_{1}},\alpha_{p_{2}},
\dotsc,\alpha_{p_{\mathfrak{s}}} \rbrace \! \to \! \mathbb{C}$ is the 
corresponding monic {\rm MPC ORF} defined in 
Subsection~\ref{subsubsec1.2.2}, and $\mathcal{X}_{21} \colon \mathbb{N} 
\times \lbrace 1,2,\dotsc,K \rbrace \times \overline{\mathbb{C}} \setminus 
\lbrace \alpha_{p_{1}},\alpha_{p_{2}},\dotsc,\alpha_{p_{\mathfrak{s}}} \rbrace 
\! \to \! \mathbb{C}$.
\end{ccccc}

\emph{Proof}. The proof consists of two cases: (i) $n \! \in \! \mathbb{N}$ 
and $k \! \in \! \lbrace 1,2,\dotsc,K \rbrace$ such that $\alpha_{
p_{\mathfrak{s}}} \! := \! \alpha_{k} \! = \! \infty$ (see case~\pmb{(1)} 
below); and (ii) $n \! \in \! \mathbb{N}$ and $k \! \in \! \lbrace 1,2,
\dotsc,K \rbrace$ such that $\alpha_{p_{\mathfrak{s}}} \! := \! \alpha_{k} 
\! \neq \! \infty$ (see case~\pmb{(2)} below). Let $\widetilde{w}(z) 
\! = \! \exp (-n \widetilde{V}(z))$, where $\widetilde{V} \colon 
\overline{\mathbb{R}} \setminus \lbrace \alpha_{1},\alpha_{2},
\dotsc,\alpha_{K} \rbrace \! \to \! \mathbb{R}$ satisfies 
conditions~\eqref{eq20}--\eqref{eq22}.

\pmb{(1)} If, for $n \! \in \! \mathbb{N}$ and $k \! \in \! \lbrace 1,2,
\dotsc,K \rbrace$ such that $\alpha_{p_{\mathfrak{s}}} \! := \! \alpha_{k} 
\! = \! \infty$,\footnote{Recall that $\mathcal{X}(n,k;z) \! := \! 
\mathcal{X}(z)$.} $\mathcal{X} \colon \overline{\mathbb{C}} \setminus 
\overline{\mathbb{R}} \! \to \! \operatorname{SL}_{2}(\mathbb{C})$ solves 
the RHP stated in Lemma~$\bm{\mathrm{RHP}_{\mathrm{MPC}}}$, then, 
{}from the jump condition~(ii) of Lemma~$\bm{\mathrm{RHP}_{\mathrm{MPC}}}$, 
it follows that, for the elements of the first column of $\mathcal{X}(z)$,
\begin{equation*}
(\mathcal{X}_{j1}(z))_{+} \! = \! (\mathcal{X}_{j1}(z))_{-} \! := \! 
\mathcal{X}_{j1}(z), \quad j \! = \! 1,2,
\end{equation*}
and, for the elements of the second column of $\mathcal{X}(z)$,
\begin{equation}
(\mathcal{X}_{j2}(z))_{+} \! - \! (\mathcal{X}_{j2}(z))_{-} \! = \! 
\mathcal{X}_{j1}(z) \widetilde{w}(z), \quad j \! = \! 1,2. \label{eq25}
\end{equation}
Via the normalisation and boundedness conditions~(iii) of 
Lemma~$\bm{\mathrm{RHP}_{\mathrm{MPC}}}$, it follows that, for 
$n \! \in \! \mathbb{N}$ and $k \! \in \! \lbrace 1,2,\dotsc,K \rbrace$ 
such that $\alpha_{p_{\mathfrak{s}}} \! := \! \alpha_{k} \! = \! \infty$,
\begin{equation}
\begin{split}
\mathcal{X}_{11}(z)z^{-\varkappa_{nk}} \underset{\overline{\mathbb{C}} 
\setminus \mathbb{R} \, \ni \, z \to \alpha_{k}}{=} 1 \! + \! \mathcal{O}
(z^{-1}), \qquad \quad \mathcal{X}_{12}(z)z^{\varkappa_{nk}} 
\underset{\overline{\mathbb{C}} \setminus \mathbb{R} \, \ni \, z \to 
\alpha_{k}}{=} \mathcal{O}(z^{-1}), \\
\mathcal{X}_{21}(z)z^{-\varkappa_{nk}} \underset{\overline{\mathbb{C}} 
\setminus \mathbb{R} \, \ni \, z \to \alpha_{k}}{=} \mathcal{O}(z^{-1}), 
\quad \qquad \mathcal{X}_{22}(z)z^{\varkappa_{nk}} 
\underset{\overline{\mathbb{C}} \setminus \mathbb{R} \, \ni \, z \to 
\alpha_{k}}{=} 1 \! + \! \mathcal{O}(z^{-1}),
\end{split} \label{eq26}
\end{equation}
and, for $q \! = \! 1,2,\dotsc,\mathfrak{s} \! - \! 1$,
\begin{equation}
\begin{split}
\mathcal{X}_{11}(z)(z \! - \! \alpha_{p_{q}})^{\varkappa_{nk \tilde{k}_{q}}} 
\underset{\mathbb{C} \setminus \mathbb{R} \, \ni \, 
z \to \alpha_{p_{q}}}{=} \mathcal{O}(1), \qquad \quad \mathcal{X}_{12}(z)
(z \! - \! \alpha_{p_{q}})^{-\varkappa_{nk \tilde{k}_{q}}} 
\underset{\mathbb{C} \setminus \mathbb{R} \, \ni \, 
z \to \alpha_{p_{q}}}{=} \mathcal{O}(1), \\
\mathcal{X}_{21}(z)(z \! - \! \alpha_{p_{q}})^{\varkappa_{nk \tilde{k}_{q}}} 
\underset{\mathbb{C} \setminus \mathbb{R} \, \ni \, 
z \to \alpha_{p_{q}}}{=} \mathcal{O}(1), \qquad \quad \mathcal{X}_{22}(z)
(z \! - \! \alpha_{p_{q}})^{-\varkappa_{nk \tilde{k}_{q}}} 
\underset{\mathbb{C} \setminus \mathbb{R} \, \ni \, 
z \to \alpha_{p_{q}}}{=} \mathcal{O}(1),
\end{split} \label{eq27}
\end{equation}
whence, temporarily re-inserting explicit $n$- and $k$-dependencies, 
that is, $\mathcal{X}_{ij}(n,k;z) \! := \! \mathcal{X}_{ij}(z)$, $i,j \! = \! 1,2$, 
one notes that, for $m \! = \! 1,2$, $\mathcal{X}_{m1} \colon \mathbb{N} 
\times \lbrace 1,2,\dotsc,K \rbrace \times \overline{\mathbb{C}} \setminus 
\lbrace \alpha_{p_{1}},\alpha_{p_{2}},\dotsc,\alpha_{p_{\mathfrak{s}}} \rbrace 
\! \to \! \mathbb{C}$ and $\mathcal{X}_{m2} \colon \mathbb{N} \times 
\lbrace 1,2,\dotsc,K \rbrace \times \overline{\mathbb{C}} \setminus 
\overline{\mathbb{R}} \! \to \! \mathbb{C}$; in particular, $\mathcal{X}_{11}(z)$ 
and $\mathcal{X}_{21}(z)$ have no jumps throughout the complex $z$-plane, 
$\mathcal{X}_{11}(z)$ is a monic (that is, $\operatorname{coeff} \lbrace 
z^{\varkappa_{nk}} \rbrace \! = \! 1)$ meromorphic function with poles 
at $\alpha_{p_{1}},\alpha_{p_{2}},\dotsc,\alpha_{p_{\mathfrak{s}}}$, and 
$\mathcal{X}_{21}(z)$ is a meromorphic function with poles at $\alpha_{p_{1}},
\alpha_{p_{2}},\dotsc,\alpha_{p_{\mathfrak{s}}}$. Application of the 
Sokhotski-Plemelj formula to the jump condition~\eqref{eq25}, with the Cauchy 
kernel normalised at $\alpha_{p_{\mathfrak{s}}} \! := \! \alpha_{k} \! = \! \infty$, 
$k \! \in \! \lbrace 1,2,\dotsc,K \rbrace$, gives rise to the following Cauchy 
integral representation:
\begin{equation}
\mathcal{X}_{j2}(z) \! = \! \int_{\mathbb{R}} \dfrac{\mathcal{X}_{j1}(\xi) 
\widetilde{w}(\xi)}{\xi \! - \! z} \, \dfrac{\md \xi}{2 \pi \mi}, \quad z 
\! \in \! \mathbb{C} \setminus \mathbb{R}, \quad j \! = \! 1,2; \label{eq28}
\end{equation}
hence, for $n \! \in \! \mathbb{N}$ and $k \! \in \! \lbrace 1,2,\dotsc,
K \rbrace$ such that $\alpha_{p_{\mathfrak{s}}} \! := \! \alpha_{k} \! 
= \! \infty$, $\mathcal{X} \colon \overline{\mathbb{C}} \setminus 
\overline{\mathbb{R}} \! \to \! \operatorname{SL}_{2}(\mathbb{C})$ 
has the integral representation
\begin{equation*}
\mathcal{X}(z) \! = \! 
\begin{pmatrix}
\mathcal{X}_{11}(z) & \int_{\mathbb{R}} \frac{\mathcal{X}_{11}(\xi) 
\widetilde{w}(\xi)}{\xi -z} \, \frac{\md \xi}{2 \pi \mi} \\
\mathcal{X}_{21}(z) & \int_{\mathbb{R}} \frac{\mathcal{X}_{21}(\xi) 
\widetilde{w}(\xi)}{\xi -z} \, \frac{\md \xi}{2 \pi \mi}
\end{pmatrix}, \quad z \! \in \! \mathbb{C} \setminus \mathbb{R}.
\end{equation*}
A detailed analysis for the elements of the first row of $\mathcal{X}(z)$ 
is presented first; then, the corresponding analysis for the elements 
of the second row is presented. Recalling that $\int_{\mathbb{R}} \xi^{m} 
\widetilde{w}(\xi) \, \md \xi \! < \! \infty$, $m \! \in \! \mathbb{N}_{0}$, 
and $\int_{\mathbb{R}}(\xi \! - \! \alpha_{k})^{-(m+1)} \widetilde{w}(\xi) 
\, \md \xi \! < \! \infty$, $\alpha_{k} \! \neq \! \infty$, $k \! \in \! \lbrace 
1,2,\dotsc,K \rbrace$, it follows, via the expansion $\tfrac{1}{z_{1}-z_{2}} 
\! = \! \sum_{i=0}^{l} \tfrac{z_{2}^{i}}{z_{1}^{i+1}} \! + \! 
\tfrac{z_{2}^{l+1}}{z_{1}^{l+1}(z_{1}-z_{2})}$, $l \! \in \! \mathbb{N}_{0}$, 
the integral representation (cf. Equation~\eqref{eq28}) $\mathcal{X}_{12}(z) \! 
= \! \int_{\mathbb{R}} \tfrac{\mathcal{X}_{11}(\xi) \widetilde{w}(\xi)}{\xi -z} 
\, \tfrac{\md \xi}{2 \pi \mi}$, $z \! \in \! \mathbb{C} \setminus \mathbb{R}$, 
and the first line of each of the asymptotic conditions~\eqref{eq26} 
and~\eqref{eq27}, that, for 
$n \! \in \! \mathbb{N}$ and $k \! \in \! \lbrace 1,2,\dotsc,K \rbrace$ 
such that $\alpha_{p_{\mathfrak{s}}} \! := \! \alpha_{k} \! = \! \infty$,
\begin{gather}
\int_{\mathbb{R}} \mathcal{X}_{11}(\xi) \xi^{p} \widetilde{w}(\xi) \, 
\md \xi \! = \! 0, \quad p \! = \! 0,1,\dotsc,\varkappa_{nk} \! - \! 1, 
\label{eq29} \\
\int_{\mathbb{R}} \mathcal{X}_{11}(\xi) \xi^{\varkappa_{nk}} \widetilde{w}
(\xi) \, \md \xi \! \neq \! 0, \label{eq30} \\
\int_{\mathbb{R}} \mathcal{X}_{11}(\xi)(\xi \! - \! \alpha_{p_{q}})^{-r} 
\widetilde{w}(\xi) \, \md \xi \! = \! 0, \quad q \! = \! 1,2,\dotsc,
\mathfrak{s} \! - \! 1, \quad r \! = \! 1,2,\dotsc,
\varkappa_{nk \tilde{k}_{q}}. \label{eq31}
\end{gather}
(Note: if, for $n \! \in \! \mathbb{N}$ and $k \! \in \! \lbrace 1,2,\dotsc,
K \rbrace$ such that $\alpha_{p_{\mathfrak{s}}} \! := \! \alpha_{k} \! = \! 
\infty$, the set $\lbrace \mathstrut \alpha_{k^{\prime}}, \, k^{\prime} 
\! \in \! \lbrace 1,2,\dotsc,K \rbrace; \, \alpha_{k^{\prime}} \! \neq 
\! \alpha_{k} \! = \! \infty \rbrace \! = \! \varnothing$, then 
Equation~\eqref{eq31} is vacuous; this can only occur if $n \! = \! 1$.) 
Recalling {}from the analysis preceding the integral 
representation~\eqref{eq28} that, for $n \! \in \! \mathbb{N}$ and $k \! 
\in \! \lbrace 1,2,\dotsc,K \rbrace$ such that $\alpha_{p_{\mathfrak{s}}} 
\! := \! \alpha_{k} \! = \! \infty$, $\mathcal{X}_{11} \colon \mathbb{N} 
\times \lbrace 1,2,\dotsc,K \rbrace \times \overline{\mathbb{C}} \setminus 
\lbrace \alpha_{p_{1}},\alpha_{p_{2}},\dotsc,\alpha_{p_{\mathfrak{s}}} 
\rbrace \! \to \! \mathbb{C}$ is a monic $(\operatorname{coeff} \lbrace 
z^{\varkappa_{nk}} \rbrace \! = \! 1)$ meromorphic function with pole set 
$\lbrace \alpha_{p_{1}},\alpha_{p_{2}},\dotsc,\alpha_{p_{\mathfrak{s}}} 
\rbrace$ and with no jumps throughout the $z$-plane, and that, 
for $n \! \in \! \mathbb{N}$ and $k \! \in \! \lbrace 1,2,\dotsc,K 
\rbrace$ such that $\alpha_{p_{\mathfrak{s}}} \! := \! \alpha_{k} \! = \! 
\infty$, the corresponding monic MPC ORF satisfies the orthogonality 
conditions~\eqref{eq9}--\eqref{eq11}, it follows {}from the latter two 
observations and Equations~\eqref{eq29}--\eqref{eq31} that, for 
$n \! \in \! \mathbb{N}$ and $k \! \in \! \lbrace 1,2,\dotsc,K \rbrace$ 
such that $\alpha_{p_{\mathfrak{s}}} \! := \! \alpha_{k} \! = \! \infty$,
\begin{equation*}
\pmb{\pi}^{n}_{k}(z) \! = \! \mathcal{X}_{11}(z).
\end{equation*}
Via Equations~\eqref{eq29} and~\eqref{eq31}, and this latter formula, one 
writes, for $n \! \in \! \mathbb{N}$ and $k \! \in \! \lbrace 1,2,\dotsc,K 
\rbrace$ such that $\alpha_{p_{\mathfrak{s}}} \! := \! \alpha_{k} \! = \! 
\infty$, Equation~\eqref{eq30} in a more transparent form:
\begin{align*}
&\int_{\mathbb{R}} \mathcal{X}_{11}(\xi) \xi^{\varkappa_{nk}} \widetilde{w}
(\xi) \, \md \xi \! = \! \int_{\mathbb{R}} \mathcal{X}_{11}(\xi) 
\mu^{\infty}_{n,\varkappa_{nk}}(n,k) \xi^{\varkappa_{nk}} 
\dfrac{\widetilde{w}(\xi)}{\mu^{\infty}_{n,\varkappa_{nk}}(n,k)} \, \md \xi \\
&= \, \int_{\mathbb{R}} \mathcal{X}_{11}(\xi) \underbrace{
\left(\phi_{0}^{\infty}(n,k) \! + \! \sum_{q=1}^{\mathfrak{s}-1} 
\sum_{r=1}^{\varkappa_{nk \tilde{k}_{q}}} \dfrac{\nu^{\infty}_{r,q}
(n,k)}{(z \! - \! \alpha_{p_{q}})^{r}} \! + \! \sum_{m=1}^{\varkappa_{nk}-1} 
\mu^{\infty}_{n,m}(n,k)z^{m} \! + \! \mu^{\infty}_{n,\varkappa_{nk}}(n,k)
z^{\varkappa_{nk}} \right)}_{= \, \phi^{n}_{k}(z)} \dfrac{\widetilde{w}
(\xi)}{\mu^{\infty}_{n,\varkappa_{nk}}(n,k)} \, \md \xi \\
&= \, \int_{\mathbb{R}} \underbrace{\mathcal{X}_{11}(\xi)}_{= \, 
(\mu^{\infty}_{n,\varkappa_{nk}}(n,k))^{-1} \phi^{n}_{k}(\xi)} \phi^{n}_{k}
(\xi) \dfrac{\widetilde{w}(\xi)}{\mu^{\infty}_{n,\varkappa_{nk}}(n,k)} \, 
\md \xi \! = \! (\mu^{\infty}_{n,\varkappa_{nk}}(n,k))^{-2} 
\underbrace{\int_{\mathbb{R}}(\phi^{n}_{k}(\xi))^{2} \, \widetilde{w}(\xi) 
\, \md \xi}_{= \, 1} \! = \! (\mu^{\infty}_{n,\varkappa_{nk}}(n,k))^{-2};
\end{align*}
hence, for $n \! \in \! \mathbb{N}$ and $k \! \in \! \lbrace 1,2,\dotsc,K 
\rbrace$ such that $\alpha_{p_{\mathfrak{s}}} \! := \! \alpha_{k} \! = \! 
\infty$, via this latter relation, Equations~\eqref{eq29}--\eqref{eq31} 
are written in the following, more convenient, form:
\begin{gather}
\int_{\mathbb{R}} \pmb{\pi}^{n}_{k}(\xi) \xi^{p} \widetilde{w}(\xi) \, 
\md \xi \! = \! 0, \quad p \! = \! 0,1,\dotsc,\varkappa_{nk} \! - \! 1, 
\label{eq32} \\
\int_{\mathbb{R}} \pmb{\pi}^{n}_{k}(\xi) \xi^{\varkappa_{nk}} \widetilde{w}
(\xi) \, \md \xi \! = \! (\mu^{\infty}_{n,\varkappa_{nk}}(n,k))^{-2}, 
\label{eq33} \\
\int_{\mathbb{R}} \pmb{\pi}^{n}_{k}(\xi)(\xi \! - \! \alpha_{p_{q}})^{-r} 
\widetilde{w}(\xi) \, \md \xi \! = \! 0, \quad q \! = \! 1,2,\dotsc,
\mathfrak{s} \! - \! 1, \quad r \! = \! 1,2,\dotsc,
\varkappa_{nk \tilde{k}_{q}}. \label{eq34}
\end{gather}
(Note: if, for $n \! \in \! \mathbb{N}$ and $k \! \in \! \lbrace 1,2,\dotsc,
K \rbrace$ such that $\alpha_{p_{\mathfrak{s}}} \! := \! \alpha_{k} \! = \! 
\infty$, the set $\lbrace \mathstrut \alpha_{k^{\prime}}, \, k^{\prime} 
\! \in \! \lbrace 1,2,\dotsc,K \rbrace; \, \alpha_{k^{\prime}} \! \neq \! 
\alpha_{k} \! = \! \infty \rbrace \! = \! \varnothing$, then 
Equation~\eqref{eq34} is vacuous; this can only occur if $n \! = \! 1$.) 
For $n \! \in \! \mathbb{N}$ and $k \! \in \! \lbrace 1,2,\dotsc,K \rbrace$ 
such that $\alpha_{p_{\mathfrak{s}}} \! := \! \alpha_{k} \! = \! \infty$, 
Equation~\eqref{eq32} gives rise to $\varkappa_{nk}$ conditions, 
Equation~\eqref{eq33} gives rise to $1$ condition, and Equation~\eqref{eq34} 
gives rise to (cf. Equation~\eqref{infcount}) $\sum_{q=1}^{\mathfrak{s}-1} 
\varkappa_{nk \tilde{k}_{q}} \! = \! (n \! - \! 1)K \! + \! k \! - \! 
\varkappa_{nk}$ conditions, for a total of $(n \! - \! 1)K \! + \! k \! + \! 1$ 
conditions, which is precisely the number necessary in order to determine, 
uniquely (see below), the associated $(n$- and $k$-dependent) norming 
constant, $\mu^{\infty}_{n,\varkappa_{nk}}(n,k)$.

One now examines, for $n \! \in \! \mathbb{N}$ and $k \! \in \! 
\lbrace 1,2,\dotsc,K \rbrace$ such that $\alpha_{p_{\mathfrak{s}}} 
\! := \! \alpha_{k} \! = \! \infty$, Equations~\eqref{eq32}--\eqref{eq34} 
in detail. Proceeding as per the detailed discussion of 
Subsection~\ref{subsubsec1.2.1}, write, for $n \! \in \! \mathbb{N}$ 
and $k \! \in \! \lbrace 1,2,\dotsc,K \rbrace$ such that 
$\alpha_{p_{\mathfrak{s}}} \! := \! \alpha_{k} \! = \! \infty$, the ordered 
disjoint partition for the repeated pole sequence:\footnote{Note the 
convention $\cup_{m=1}^{0} \lbrace \pmb{\ast},\pmb{\ast},\dotsc,
\pmb{\ast} \rbrace \! := \! \varnothing$.}
\begin{align*}
&\lbrace \overset{1}{\underbrace{\alpha_{1},\alpha_{2},\dotsc,\alpha_{K}}_{K}} 
\rbrace \cup \dotsb \cup \lbrace \overset{n-1}{\underbrace{\alpha_{1},
\alpha_{2},\dotsc,\alpha_{K}}_{K}} \rbrace \cup \lbrace 
\overset{n}{\underbrace{\alpha_{1},\alpha_{2},\dotsc,\alpha_{k}}_{k}} 
\rbrace \\
& \, := \bigcup_{q=1}^{\mathfrak{s}-1} \lbrace 
\underbrace{\alpha_{i(q)_{k_{q}}},\alpha_{i(q)_{k_{q}}},\dotsc,
\alpha_{i(q)_{k_{q}}}}_{l_{q}=\varkappa_{nk \tilde{k}_{q}}} \rbrace \cup 
\lbrace \underbrace{\alpha_{i(\mathfrak{s})_{k_{\mathfrak{s}}}},
\alpha_{i(\mathfrak{s})_{k_{\mathfrak{s}}}},\dotsc,
\alpha_{i(\mathfrak{s})_{k_{\mathfrak{s}}}}}_{l_{\mathfrak{s}}=
\varkappa_{ni(\mathfrak{s})_{k_{\mathfrak{s}}}}} \rbrace \\
& \, := \bigcup_{q=1}^{\mathfrak{s}-1} \lbrace \underbrace{\alpha_{p_{q}},
\alpha_{p_{q}},\dotsc,\alpha_{p_{q}}}_{l_{q}=\varkappa_{nk \tilde{k}_{q}}} 
\rbrace \cup \lbrace \underbrace{\alpha_{k},\alpha_{k},\dotsc,
\alpha_{k}}_{l_{\mathfrak{s}}=\varkappa_{nk}} \rbrace \! = 
\bigcup_{q=1}^{\mathfrak{s}-1} \lbrace \underbrace{\alpha_{p_{q}},
\alpha_{p_{q}},\dotsc,\alpha_{p_{q}}}_{l_{q}=\varkappa_{nk \tilde{k}_{q}}} 
\rbrace \cup \lbrace \underbrace{\infty,\infty,\dotsc,\infty}_{
l_{\mathfrak{s}}=\varkappa_{nk}} \rbrace,
\end{align*}
where $\sum_{q=1}^{\mathfrak{s}}l_{q} \! = \! \sum_{q=1}^{\mathfrak{s}-1}
l_{q} \! + \! l_{\mathfrak{s}} \! = \! \sum_{q=1}^{\mathfrak{s}-1} \varkappa_{nk 
\tilde{k}_{q}} \! + \! \varkappa_{nk} \! = \! (n \! - \! 1)K \! + \! k$. 
Hence, via this notational preamble, and the analysis leading to the 
orthogonality Equations~\eqref{eq32}--\eqref{eq34}, one writes, in the 
indicated order, for $n \! \in \! \mathbb{N}$ and $k \! \in \! \lbrace 
1,2,\dotsc,K \rbrace$ such that $\alpha_{p_{\mathfrak{s}}} \! := \! 
\alpha_{k} \! = \! \infty$,
\begin{align*}
\pmb{\pi}^{n}_{k}(z)=& \, \dfrac{\phi_{0}^{\infty}(n,k)}{\mu^{\infty}_{n,
\varkappa_{nk}}(n,k)} \! + \! \dfrac{1}{\mu^{\infty}_{n,\varkappa_{nk}}(n,k)} 
\sum_{q=1}^{\mathfrak{s}-1} \sum_{r=1}^{\varkappa_{nk \tilde{k}_{q}}} 
\dfrac{\nu^{\infty}_{r,q}(n,k)}{(z \! - \! \alpha_{i(q)_{k_{q}}})^{r}} 
\! + \! \dfrac{1}{\mu^{\infty}_{n,\varkappa_{nk}}(n,k)} 
\sum_{m=1}^{\varkappa_{ni(\mathfrak{s})_{k_{\mathfrak{s}}}}-1} 
\mu^{\infty}_{n,m}(n,k)z^{m} \! + \! z^{\varkappa_{nk}} \\
=& \, \dfrac{\phi^{\infty}_{0}(n,k)}{\mu^{\infty}_{n,\varkappa_{nk}}(n,k)} 
\! + \! \dfrac{1}{\mu^{\infty}_{n,\varkappa_{nk}}(n,k)} 
\sum_{q=1}^{\mathfrak{s}-1} \sum_{r=1}^{\varkappa_{nk \tilde{k}_{q}}} 
\dfrac{\nu^{\infty}_{r,q}(n,k)}{(z \! - \! \alpha_{p_{q}})^{r}} \! + \! 
\dfrac{1}{\mu^{\infty}_{n,\varkappa_{nk}}(n,k)} \sum_{m=1}^{\varkappa_{nk}-1} 
\mu^{\infty}_{n,m}(n,k)z^{m} \! + \! z^{\varkappa_{nk}} \\
:=& \, \widetilde{\phi}_{0}^{\raise-1.0ex\hbox{$\scriptstyle \infty$}}(n,k) 
\! + \! \sum_{m=1}^{\mathfrak{s}-1} 
\sum_{q=1}^{l_{m}=\varkappa_{nk \tilde{k}_{m}}} 
\dfrac{\widetilde{\nu}_{m,q}^{\raise-1.0ex\hbox{$\scriptstyle \infty$}}
(n,k)}{(z \! - \! \alpha_{p_{m}})^{q}} \! + \! 
\sum_{r=1}^{l_{\mathfrak{s}}=\varkappa_{nk}} 
\widetilde{\nu}_{\mathfrak{s},r}^{\raise-1.0ex\hbox{$\scriptstyle \infty$}}
(n,k)z^{r}, \qquad \widetilde{\nu}_{\mathfrak{s},
l_{\mathfrak{s}}}^{\raise-1.0ex\hbox{$\scriptstyle \infty$}}(n,k) 
\! \equiv \! 1.
\end{align*}
Substituting the latter partial fraction expansion for $\pmb{\pi}^{n}_{k}(z)$ 
into the orthogonality conditions~\eqref{eq32}--\eqref{eq34}, one arrives at, 
for $n \! \in \! \mathbb{N}$ and $k \! \in \! \lbrace 1,2,\dotsc,K \rbrace$ 
such that $\alpha_{p_{\mathfrak{s}}} \! := \! \alpha_{k} \! = \! \infty$, the 
orthogonality conditions
\begin{gather*}
\int_{\mathbb{R}} \left(
\widetilde{\phi}_{0}^{\raise-1.0ex\hbox{$\scriptstyle \infty$}}(n,k) \! + \! 
\sum_{m=1}^{\mathfrak{s}-1} \sum_{q=1}^{l_{m}=\varkappa_{nk \tilde{k}_{m}}} 
\dfrac{\widetilde{\nu}_{m,q}^{\raise-1.0ex\hbox{$\scriptstyle \infty$}}
(n,k)}{(\xi \! - \! \alpha_{p_{m}})^{q}} \! + \! 
\sum_{q=1}^{l_{\mathfrak{s}}=\varkappa_{nk}} 
\widetilde{\nu}_{\mathfrak{s},q}^{\raise-1.0ex\hbox{$\scriptstyle \infty$}}
(n,k) \xi^{q} \right) \xi^{r} \widetilde{w}(\xi) \, \md \xi \! = \! 0, \quad 
r \! = \! 0,1,\dotsc,\varkappa_{nk} \! - \! 1, \\
\int_{\mathbb{R}} \left(
\widetilde{\phi}_{0}^{\raise-1.0ex\hbox{$\scriptstyle \infty$}}(n,k) \! + \! 
\sum_{m=1}^{\mathfrak{s}-1} \sum_{q=1}^{l_{m}=\varkappa_{nk \tilde{k}_{m}}} 
\dfrac{\widetilde{\nu}_{m,q}^{\raise-1.0ex\hbox{$\scriptstyle \infty$}}
(n,k)}{(\xi \! - \! \alpha_{p_{m}})^{q}} \! + \! 
\sum_{q=1}^{l_{\mathfrak{s}}=\varkappa_{nk}} 
\widetilde{\nu}_{\mathfrak{s},q}^{\raise-1.0ex\hbox{$\scriptstyle \infty$}}
(n,k) \xi^{q} \right) \xi^{\varkappa_{nk}} \widetilde{w}(\xi) \, \md \xi \! 
= \! (\mu^{\infty}_{n,\varkappa_{nk}}(n,k))^{-2}, \\
\int_{\mathbb{R}} \left(
\widetilde{\phi}_{0}^{\raise-1.0ex\hbox{$\scriptstyle \infty$}}(n,k) \! + \! 
\sum_{m=1}^{\mathfrak{s}-1} \sum_{q=1}^{l_{m}=\varkappa_{nk \tilde{k}_{m}}} 
\dfrac{\widetilde{\nu}_{m,q}^{\raise-1.0ex\hbox{$\scriptstyle \infty$}}
(n,k)}{(\xi \! - \! \alpha_{p_{m}})^{q}} \! + \! 
\sum_{q=1}^{l_{\mathfrak{s}}=\varkappa_{nk}} 
\widetilde{\nu}_{\mathfrak{s},q}^{\raise-1.0ex\hbox{$\scriptstyle \infty$}}
(n,k) \xi^{q} \right) 
\dfrac{\widetilde{w}(\xi)}{(\xi \! - \! \alpha_{p_{i}})^{j}} \, \md \xi 
\! = \! 0, \quad i \! = \! 1,2,\dotsc,\mathfrak{s} \! - \! 1, \quad j \! 
= \! 1,2,\dotsc,l_{i}.
\end{gather*}
For $n \! \in \! \mathbb{N}$ and $k \! \in \! \lbrace 1,2,\dotsc,K \rbrace$ 
such that $\alpha_{p_{\mathfrak{s}}} \! := \! \alpha_{k} \! = \! \infty$, 
the above orthogonality conditions give rise to a total of $\sum_{r=1}^{
\mathfrak{s}}l_{r} \! + \! 1 \! = \! \sum_{r=1}^{\mathfrak{s}-1}l_{r} \! 
+ \! l_{\mathfrak{s}} \! + \! 1 \! = \! \sum_{r=1}^{\mathfrak{s}-1} 
\varkappa_{nk \tilde{k}_{r}} \! + \! \varkappa_{nk} \! + \! 1 \! = \! 
(n \! - \! 1)K \! + \! k \! + \! 1$ linear inhomogeneous algebraic 
equations for the $(n \! - \! 1)K \! + \! k \! + \! 1$ (real) unknowns 
$\widetilde{\phi}_{0}^{\raise-1.0ex\hbox{$\scriptstyle \infty$}}(n,k),
\widetilde{\nu}_{1,1}^{\raise-1.0ex\hbox{$\scriptstyle \infty$}}(n,k),
\dotsc,\widetilde{\nu}_{1,l_{1}}^{\raise-1.0ex\hbox{$\scriptstyle \infty$}}
(n,k),\dotsc,
\widetilde{\nu}_{\mathfrak{s}-1,1}^{\raise-1.0ex\hbox{$\scriptstyle \infty$}}
(n,k),\dotsc,\widetilde{\nu}_{\mathfrak{s}-1,
l_{\mathfrak{s}-1}}^{\raise-1.0ex\hbox{$\scriptstyle \infty$}}(n,k),
\widetilde{\nu}_{\mathfrak{s},1}^{\raise-1.0ex\hbox{$\scriptstyle \infty$}}
(n,k),\dotsc,\widetilde{\nu}_{\mathfrak{s},
l_{\mathfrak{s}}-1}^{\raise-1.0ex\hbox{$\scriptstyle \infty$}}(n,k),
\linebreak[4]
(\mu_{n,\varkappa_{nk}}^{\infty}(n,k))^{-2}$, that is, with $\md \widetilde{\mu}
(z) \! := \! \widetilde{w}(z) \, \md z$,
\begin{align}
\setcounter{MaxMatrixCols}{12}
&\left(

\right), \label{eq35}
\end{align}
where
\begin{equation*}
n_{1} \! = \! l_{1} \! + \! \dotsb \! + \! l_{\mathfrak{s}-2} \! + \! 1, \qquad 
n_{2} \! = \! (n \! - \! 1)K \! + \! k \! - \! \varkappa_{nk}, \qquad m_{1} \! 
= \! (n \! - \! 1)K \! + \! k \! - \! 1, \qquad \text{and} \qquad m_{2} \! = 
\! (n \! - \! 1)K \! + \! k.
\end{equation*}
For $n \! \in \! \mathbb{N}$ and $k \! \in \! \lbrace 1,2,\dotsc,K \rbrace$ 
such that $\alpha_{p_{\mathfrak{s}}} \! := \! \alpha_{k} \! = \! \infty$, 
the linear system~\eqref{eq35} of $(n \! - \! 1)K \! + \! k \! + \! 1$ 
inhomogeneous algebraic equations for the $(n \! - \! 1)K \! + \! k \! 
+ \! 1$ (real) unknowns 
$\widetilde{\phi}_{0}^{\raise-1.0ex\hbox{$\scriptstyle \infty$}}(n,k),
\widetilde{\nu}_{1,1}^{\raise-1.0ex\hbox{$\scriptstyle \infty$}}(n,k),
\dotsc,\widetilde{\nu}_{1,l_{1}}^{\raise-1.0ex\hbox{$\scriptstyle \infty$}}
(n,k),\dotsc,
\widetilde{\nu}_{\mathfrak{s}-1,1}^{\raise-1.0ex\hbox{$\scriptstyle \infty$}}
(n,k),\dotsc,\linebreak[4]
\widetilde{\nu}_{\mathfrak{s}-1,
l_{\mathfrak{s}-1}}^{\raise-1.0ex\hbox{$\scriptstyle \infty$}}(n,k),
\widetilde{\nu}_{\mathfrak{s},1}^{\raise-1.0ex\hbox{$\scriptstyle \infty$}}
(n,k),\dotsc,\widetilde{\nu}_{\mathfrak{s},
l_{\mathfrak{s}}-1}^{\raise-1.0ex\hbox{$\scriptstyle \infty$}}(n,k),
(\mu_{n,\varkappa_{nk}}^{\infty}(n,k))^{-2}$ admits a unique solution if, 
and only if, the determinant of the coefficient matrix is non-zero; this 
fact will now be established, and, en route, an explicit, multi-integral 
representation for the associated $(n$- and $k$-dependent) norming 
constant, $(\mu_{n,\varkappa_{nk}}^{\infty}(n,k))^{-2}$, will be derived. 
For $n \! \in \! \mathbb{N}$ and $k \! \in \! \lbrace 1,2,\dotsc,K \rbrace$ 
such that $\alpha_{p_{\mathfrak{s}}} \! := \! \alpha_{k} \! = \! \infty$, 
one uses Cramer's Rule and the multi-linearity property of the determinant 
to show that
\begin{equation} \label{nhormmconstatinf}
\left(\mu_{n,\varkappa_{nk}}^{\infty}(n,k) \right)^{-2} \! = \! 
\dfrac{\hat{\mathfrak{c}}_{N_{\infty}}}{\hat{\mathfrak{c}}_{D_{\infty}}},
\end{equation}
where
\begin{align*}
\hat{\mathfrak{c}}_{N_{\infty}} :=& \, \underbrace{\int_{\mathbb{R}} 
\int_{\mathbb{R}} \dotsb \int_{\mathbb{R}}}_{(n-1)K+k+1} \md \widetilde{\mu}
(\xi_{0}) \, \md \widetilde{\mu}(\xi_{1}) \, \dotsb \, \md \widetilde{\mu}(\xi_{l_{1}}) 
\, \dotsb \, \md \widetilde{\mu}(\xi_{n_{1}}) \, \dotsb \, \md \widetilde{\mu}
(\xi_{n_{2}}) \, \md \widetilde{\mu}(\xi_{m_{2}-\varkappa_{nk}+1}) \, \dotsb \, 
\md \widetilde{\mu} (\xi_{m_{2}}) \\
\times& \, \dfrac{1}{\dfrac{1}{\xi_{0}^{0}}(\xi_{1} \! - \! \alpha_{p_{1}})^{1} 
\dotsb (\xi_{l_{1}} \! - \! \alpha_{p_{1}})^{l_{1}} \dotsb (\xi_{n_{1}} \! - \! 
\alpha_{p_{\mathfrak{s}-1}})^{1} \dotsb (\xi_{n_{2}} \! - \! 
\alpha_{p_{\mathfrak{s}-1}})^{l_{\mathfrak{s}-1}} \dfrac{1}{(\xi_{m_{2}-
\varkappa_{nk}+1})^{1}} \dotsb \dfrac{1}{(\xi_{m_{2}})^{\varkappa_{nk}}}} \\
\times& \, 
\left\vert

\right\vert,
\end{align*}
and
\begin{align*}
\hat{\mathfrak{c}}_{D_{\infty}} :=& \, \underbrace{\int_{\mathbb{R}} \int_{\mathbb{R}} 
\dotsb \int_{\mathbb{R}}}_{(n-1)K+k} \md \widetilde{\mu}(\xi_{0}) \, \md 
\widetilde{\mu}(\xi_{1}) \, \dotsb \, \md \widetilde{\mu}(\xi_{l_{1}}) \, \dotsb \, 
\md \widetilde{\mu}(\xi_{n_{1}}) \, \dotsb \, \md \widetilde{\mu}(\xi_{n_{2}}) 
\, \md \widetilde{\mu}(\xi_{m_{2}-\varkappa_{nk}+1}) \, \dotsb \, \md 
\widetilde{\mu}(\xi_{m_{1}}) \\
\times& \, \dfrac{1}{\dfrac{1}{\xi_{0}^{0}}(\xi_{1} \! - \! 
\alpha_{p_{1}})^{1} \dotsb (\xi_{l_{1}} \! - \! \alpha_{p_{1}})^{l_{1}} 
\dotsb (\xi_{n_{1}} \! - \! \alpha_{p_{\mathfrak{s}-1}})^{1} \dotsb 
(\xi_{n_{2}} \! - \! \alpha_{p_{\mathfrak{s}-1}})^{l_{\mathfrak{s}-1}} 
\dfrac{1}{(\xi_{m_{2}-\varkappa_{nk}+1})^{1}} \dotsb 
\dfrac{1}{(\xi_{m_{1}})^{\varkappa_{nk}-1}}} \\
\times& \, 
\left\vert

\right\vert.
\end{align*}
For $n \! \in \! \mathbb{N}$ and $k \! \in \! \lbrace 1,2,\dotsc,K 
\rbrace$ such that $\alpha_{p_{\mathfrak{s}}} \! := \! \alpha_{k} \! 
= \! \infty$, $\hat{\mathfrak{c}}_{N_{\infty}}$ is studied first, and 
then $\hat{\mathfrak{c}}_{D_{\infty}}$. For $n \! \in \! \mathbb{N}$ 
and $k \! \in \! \lbrace 1,2,\dotsc,K \rbrace$ such that 
$\alpha_{p_{\mathfrak{s}}} \! := \! \alpha_{k} \! = \! \infty$, 
introduce the following notation (recall that $l_{q} \! = \! 
\varkappa_{nk \tilde{k}_{q}}$, $q \! = \! 1,2,\dotsc,\mathfrak{s} 
\! - \! 1)$:
\begin{equation*}
\widehat{\phi}_{0}(z) \! := \! \prod_{m=1}^{\mathfrak{s}-1}(z \! - \! 
\alpha_{p_{m}})^{l_{m}} \! =: \sum_{j=0}^{(n-1)K+k} 
\mathfrak{a}^{\spcheck}_{j,0}z^{j},
\end{equation*}
and (with some abuse of notation), for $r \! = \! 1,\dotsc,\mathfrak{s} 
\! - \! 1,\mathfrak{s}$, $q(r) \! = \! \sum_{i=1}^{r-1}l_{i} \! + \! 1,
\sum_{i=1}^{r-1}l_{i} \! + \! 2,\dotsc,\sum_{i=1}^{r-1}l_{i} \! + \! l_{r}$, 
and $m(r) \! = \! 1,2,\dotsc,l_{r}$,\footnote{With the convention $(z 
\! - \! \infty)^{0} \! := \! 1$.}
\begin{equation*}
\widehat{\phi}_{q(r)}(z) \! := \! \left\lbrace \mathstrut \widehat{\phi}_{0}
(z)(z \! - \! \alpha_{p_{r}})^{-m(r)(1-\delta_{r \mathfrak{s}})}
z^{m(r) \delta_{r \mathfrak{s}}} \! =: \sum_{j=0}^{(n-1)K+k} 
\mathfrak{a}^{\spcheck}_{j,q(r)}z^{j} \right\rbrace;
\end{equation*}
e.g., for $r \! = \! 1$, the notation $\widehat{\phi}_{q(1)}(z) \! := \! 
\lbrace \widehat{\phi}_{0}(z)(z \! - \! \alpha_{p_{1}})^{-m(1)} \! =: \! 
\sum_{j=0}^{(n-1)K+k} \mathfrak{a}^{\spcheck}_{j,q(1)}z^{j} \rbrace$, 
$q(1) \! = \! 1,2,\dotsc,l_{1}$, $m(1) \! = \! 1,2,\dotsc,l_{1}$, denotes 
the (set of) $l_{1} \! = \! \varkappa_{nk \tilde{k}_{1}}$ functions
\begin{gather*}
\widehat{\phi}_{1}(z) \! = \! \dfrac{\widehat{\phi}_{0}(z)}{z \! - \! 
\alpha_{p_{1}}} \! =: \sum_{j=0}^{(n-1)K+k} \mathfrak{a}^{\spcheck}_{j,1}
z^{j}, \, \widehat{\phi}_{2}(z) \! = \! \dfrac{\widehat{\phi}_{0}(z)}{(z 
\! - \! \alpha_{p_{1}})^{2}} \! =: \sum_{j=0}^{(n-1)K+k} 
\mathfrak{a}^{\spcheck}_{j,2}z^{j}, \, \dotsc, \, \widehat{\phi}_{l_{1}}(z) 
\! = \! \dfrac{\widehat{\phi}_{0}(z)}{(z \! - \! \alpha_{p_{1}})^{l_{1}}} 
\! =: \sum_{j=0}^{(n-1)K+k} \mathfrak{a}^{\spcheck}_{j,l_{1}}z^{j},
\end{gather*}
for $r \! = \! \mathfrak{s} \! - \! 1$, the notation 
$\widehat{\phi}_{q(\mathfrak{s}-1)}(z) \! := \! \lbrace \widehat{\phi}_{0}
(z)(z \! - \! \alpha_{p_{\mathfrak{s}-1}})^{-m(\mathfrak{s}-1)} \! =: \! 
\sum_{j=0}^{(n-1)K+k} \mathfrak{a}^{\spcheck}_{j,q(\mathfrak{s}-1)}z^{j} 
\rbrace$, $q(\mathfrak{s} \! - \! 1) \! = \! l_{1} \! + \! \dotsb \! + \! 
l_{\mathfrak{s}-2} \! + \! 1,l_{1} \! + \! \dotsb \! + \! l_{\mathfrak{s}-2} 
\! + \! 2,\dotsc,l_{1} \! + \! \dotsb \! + \! l_{\mathfrak{s}-2} \! + \! 
l_{\mathfrak{s}-1} \! = \! (n \! - \! 1)K \! + \! k \! - \! \varkappa_{nk}$, 
$m(\mathfrak{s} \! - \! 1) \! = \! 1,2,\dotsc,l_{\mathfrak{s}-1}$, 
denotes the (set of) $l_{\mathfrak{s}-1} \! = \! \varkappa_{nk 
\tilde{k}_{\mathfrak{s}-1}}$ functions
\begin{gather*}
\widehat{\phi}_{l_{1}+ \dotsb +l_{\mathfrak{s}-2}+1}(z) \! = \! 
\dfrac{\widehat{\phi}_{0}(z)}{z \! - \! \alpha_{p_{\mathfrak{s}-1}}} \! 
=: \sum_{j=0}^{(n-1)K+k} \mathfrak{a}^{\spcheck}_{j,l_{1}+\dotsb 
+l_{\mathfrak{s}-2}+1}z^{j}, \, \widehat{\phi}_{l_{1}+\dotsb 
+l_{\mathfrak{s}-2}+2}(z) \! = \! \dfrac{\widehat{\phi}_{0}(z)}{(z \! 
- \! \alpha_{p_{\mathfrak{s}-1}})^{2}} \! =: \sum_{j=0}^{(n-1)K+k} 
\mathfrak{a}^{\spcheck}_{j,l_{1}+\dotsb +l_{\mathfrak{s}-2}+2}
z^{j}, \, \dotsc \\\dotsc, \, \widehat{\phi}_{(n-1)K+k-
\varkappa_{nk}}(z) \! = \! \dfrac{\widehat{\phi}_{0}(z)}{(z \! 
- \! \alpha_{p_{\mathfrak{s}-1}})^{l_{\mathfrak{s}-1}}} \! =: 
\sum_{j=0}^{(n-1)K+k} \mathfrak{a}^{\spcheck}_{j,(n-1)K+k
-\varkappa_{nk}}z^{j},
\end{gather*}
etc., and, for $r \! = \! \mathfrak{s}$, the notation 
$\widehat{\phi}_{q(\mathfrak{s})}(z) \! := \! \lbrace \widehat{\phi}_{0}
(z)z^{m(\mathfrak{s})} \! =: \! \sum_{j=0}^{(n-1)K+k} 
\mathfrak{a}^{\spcheck}_{j,q(\mathfrak{s})}z^{j} \rbrace$, 
$q(\mathfrak{s}) \! = \! (n \! - \! 1)K \! + \! k \! - \! \varkappa_{nk} 
\! + \! 1,(n \! - \! 1)K \! + \! k \! - \! \varkappa_{nk} \! + \! 2,
\dotsc,(n \! - \! 1)K \! + \! k$, $m(\mathfrak{s}) \! = \! 1,2,
\dotsc,l_{\mathfrak{s}}$, denotes the (set of) $l_{\mathfrak{s}} \! 
= \! \varkappa_{nk}$ functions
\begin{gather*}
\widehat{\phi}_{(n-1)K+k-\varkappa_{nk}+1}(z) \! = \! 
\widehat{\phi}_{0}(z)z \! =: \sum_{j=0}^{(n-1)K+k} 
\mathfrak{a}^{\spcheck}_{j,(n-1)K+k-\varkappa_{nk}+1}z^{j}, \, 
\widehat{\phi}_{(n-1)K+k-\varkappa_{nk}+2}(z) \! = \! \widehat{\phi}_{0}
(z)z^{2} \\
=: \sum_{j=0}^{(n-1)K+k} \mathfrak{a}^{\spcheck}_{j,(n-1)K+k-\varkappa_{nk}+2}
z^{j}, \, \dotsc, \, \widehat{\phi}_{(n-1)K+k}(z) \! = \! \widehat{\phi}_{0}
(z)z^{\varkappa_{nk}} \! =: \sum_{j=0}^{(n-1)K+k} 
\mathfrak{a}^{\spcheck}_{j,(n-1)K+k}z^{j}.
\end{gather*}
(Note: $\# \lbrace \widehat{\phi}_{0}(z)(z \! - \! \alpha_{p_{r}})^{-m(r)
(1-\delta_{r \mathfrak{s}})}z^{m(r) \delta_{r \mathfrak{s}}} \rbrace \! 
= \! l_{r}$, $r \! = \! 1,2,\dotsc,\mathfrak{s}$, and $\# \cup_{r=1}^{
\mathfrak{s}} \lbrace \widehat{\phi}_{0}(z)(z \! - \! \alpha_{p_{r}})^{-
m(r)(1-\delta_{r \mathfrak{s}})}z^{m(r) \delta_{r \mathfrak{s}}} \rbrace 
\! = \! \sum_{r=1}^{\mathfrak{s}}l_{r} \! = \! (n \! - \! 1)K \! + \! k$.) 
One notes that, for $n \! \in \! \mathbb{N}$ and $k \! \in \! \lbrace 
1,2,\dotsc,K \rbrace$ such that $\alpha_{p_{\mathfrak{s}}} \! := \! 
\alpha_{k} \! = \! \infty$, the $l_{1} \! + \! \dotsb \! + \! 
l_{\mathfrak{s}-1} \! + \! l_{\mathfrak{s}} \! + \! 1 \! = \! 
(n \! - \! 1)K \! + \! k \! + \! 1$ functions $\widehat{\phi}_{0}(z),
\widehat{\phi}_{1}(z),\dotsc,\widehat{\phi}_{l_{1}}(z),\dotsc,
\widehat{\phi}_{l_{1}+\dotsb +l_{\mathfrak{s}-2}+1}(z),\dotsc,
\widehat{\phi}_{(n-1)K+k-\varkappa_{nk}}(z),
\widehat{\phi}_{(n-1)K+k-\varkappa_{nk}+1}(z),\dotsc,
\widehat{\phi}_{(n-1)K+k}(z)$ are linearly independent on $\mathbb{R}$, 
that is, for $z \! \in \! \mathbb{R}$, $\sum_{j=0}^{(n-1)K+k} 
\mathfrak{c}^{\spcheck}_{j} \widehat{\phi}_{j}(z) \! = \! 0$ $\Rightarrow$ 
(via a Vandermonde-type argument; see the $((n \! - \! 1)K \! + \! k \! 
+ \! 1) \times ((n \! - \! 1)K \! + \! k \! + \! 1)$ non-zero determinant 
$\mathbb{D}^{\spcheck}$ in Equation~\eqref{eq36} below) 
$\mathfrak{c}^{\spcheck}_{j} \! = \! 0$, $j \! = \! 0,1,\dotsc,
(n \! - \! 1)K \! + \! k$. For $n \! \in \! \mathbb{N}$ and $k \! \in \! 
\lbrace 1,2,\dotsc,K \rbrace$ such that $\alpha_{p_{\mathfrak{s}}} \! 
:= \! \alpha_{k} \! = \! \infty$, let $\mathfrak{S}_{(n-1)K+k+1}$ denote 
the $((n \! - \! 1)K \! + \! k \! + \! 1)!$ permutations of $\lbrace 0,1,
\dotsc,(n \! - \! 1)K \! + \! k \rbrace$. Using the above notation and the 
multi-linearity property of the determinant, one studies, thus, for $n \! 
\in \! \mathbb{N}$ and $k \! \in \! \lbrace 1,2,\dotsc,K \rbrace$ such 
that $\alpha_{p_{\mathfrak{s}}} \! := \! \alpha_{k} \! = \! \infty$, 
$\hat{\mathfrak{c}}_{N_{\infty}}$:
\begin{align*}
\hat{\mathfrak{c}}_{N_{\infty}} =& \, \underbrace{\int_{\mathbb{R}} 
\int_{\mathbb{R}} \dotsb \int_{\mathbb{R}}}_{(n-1)K+k+1} \md \widetilde{\mu}
(\xi_{0}) \, \md \widetilde{\mu}(\xi_{1}) \, \dotsb \, \md \widetilde{\mu}
(\xi_{l_{1}}) \, \dotsb \, \md \widetilde{\mu}(\xi_{n_{1}}) \, \dotsb \, \md 
\widetilde{\mu}(\xi_{n_{2}}) \, \md \widetilde{\mu}(\xi_{m_{2}-\varkappa_{nk}+1}) 
\, \dotsb \, \md \widetilde{\mu}(\xi_{m_{2}}) \\
\times& \, \dfrac{1}{\dfrac{1}{\xi_{0}^{0}}(\xi_{1} \! - \! 
\alpha_{p_{1}})^{1} \dotsb (\xi_{l_{1}} \! - \! \alpha_{p_{1}})^{l_{1}} 
\dotsb (\xi_{n_{1}} \! - \! \alpha_{p_{\mathfrak{s}-1}})^{1} \dotsb 
(\xi_{n_{2}} \! - \! \alpha_{p_{\mathfrak{s}-1}})^{l_{\mathfrak{s}-1}} 
\dfrac{1}{(\xi_{m_{2}-\varkappa_{nk}+1})^{1}} \dotsb 
\dfrac{1}{\xi_{m_{2}}^{\varkappa_{nk}}}} \\
\times& \, \dfrac{1}{\widehat{\phi}_{0}(\xi_{0}) \widehat{\phi}_{0}(\xi_{1}) 
\dotsb \widehat{\phi}_{0}(\xi_{l_{1}}) \dotsb \widehat{\phi}_{0}(\xi_{n_{1}}) 
\dotsb \widehat{\phi}_{0}(\xi_{n_{2}}) \widehat{\phi}_{0}(\xi_{m_{2}-
\varkappa_{nk}+1}) \dotsb \widehat{\phi}_{0}(\xi_{m_{2}})} \\
\times& \, 
\underbrace{\left\vert

\right\vert}_{=: \, \mathbb{G}^{\spcheck}(\xi_{0},\xi_{1},\dotsc,\xi_{l_{1}},
\dotsc,\xi_{n_{1}},\dotsc,\xi_{n_{2}},\xi_{m_{2}-\varkappa_{nk}+1},\dotsc,
\xi_{(n-1)K+k})} \\
=& \, \dfrac{1}{(m_{2} \! + \! 1)!} \sum_{\pmb{\sigma} \in \mathfrak{S}_{m_{2}
+1}} \underbrace{\int_{\mathbb{R}} \int_{\mathbb{R}} \dotsb \int_{\mathbb{R}}
}_{(n-1)K+k+1} \md \widetilde{\mu}(\xi_{\sigma (0)}) \, \md \widetilde{\mu}
(\xi_{\sigma (1)}) \, \dotsb \, \md \widetilde{\mu}(\xi_{\sigma (l_{1})}) \, \dotsb 
\, \md \widetilde{\mu}(\xi_{\sigma (n_{1})}) \, \dotsb \, \md \widetilde{\mu}
(\xi_{\sigma (n_{2})}) \, \dotsb \\
\times& \, \tfrac{\dotsb \, \md \widetilde{\mu}(\xi_{\sigma (m_{2}-\varkappa_{nk}+1)}) 
\, \dotsb \, \md \widetilde{\mu}(\xi_{\sigma (m_{2})})}{\dfrac{1}{\xi_{\sigma (0)}^{0}}
(\xi_{\sigma (1)}-\alpha_{p_{1}})^{1} \dotsb (\xi_{\sigma (l_{1})}-
\alpha_{p_{1}})^{l_{1}} \dotsb (\xi_{\sigma (n_{1})}-\alpha_{p_{\mathfrak{s}
-1}})^{1} \dotsb (\xi_{\sigma (n_{2})}-\alpha_{p_{\mathfrak{s}-1}})^{
l_{\mathfrak{s}-1}} \dfrac{1}{(\xi_{\sigma (m_{2}-\varkappa_{nk}+1)})^{1}} 
\dotsb \dfrac{1}{\xi_{\sigma (m_{2})}^{\varkappa_{nk}}}} \\
\times& \, \dfrac{1}{\widehat{\phi}_{0}(\xi_{\sigma (0)}) \widehat{\phi}_{0}
(\xi_{\sigma (1)}) \dotsb \widehat{\phi}_{0}(\xi_{\sigma (l_{1})}) \dotsb 
\widehat{\phi}_{0}(\xi_{\sigma (n_{1})}) \dotsb \widehat{\phi}_{0}(\xi_{\sigma 
(n_{2})}) \widehat{\phi}_{0}(\xi_{\sigma (m_{2}-\varkappa_{nk}+1)}) \dotsb 
\widehat{\phi}_{0}(\xi_{\sigma (m_{2})})} \\
\times& \, \mathbb{G}^{\spcheck}(\xi_{\sigma (0)},\xi_{\sigma (1)},\dotsc,
\xi_{\sigma (l_{1})},\dotsc,\xi_{\sigma (n_{1})},\dotsc,\xi_{\sigma (n_{2})},
\xi_{\sigma (m_{2}-\varkappa_{nk}+1)},\dotsc,\xi_{\sigma ((n-1)K+k)}) \\
=& \, \dfrac{1}{(m_{2} \! + \! 1)!} \underbrace{\int_{\mathbb{R}} 
\int_{\mathbb{R}} \dotsb \int_{\mathbb{R}}}_{(n-1)K+k+1} \md \widetilde{\mu}
(\xi_{0}) \, \md \widetilde{\mu}(\xi_{1}) \, \dotsb \, \md \widetilde{\mu}(\xi_{l_{1}}) 
\, \dotsb \, \md \widetilde{\mu}(\xi_{n_{1}}) \, \dotsb \, \md \widetilde{\mu}
(\xi_{n_{2}}) \, \md \widetilde{\mu}(\xi_{m_{2}-\varkappa_{nk}+1}) \, \dotsb \\
\times& \, \dotsb \, \md \widetilde{\mu}(\xi_{m_{2}}) \, \mathbb{G}^{\spcheck}
(\xi_{0},\xi_{1},\dotsc,\xi_{l_{1}},\dotsc,\xi_{n_{1}},\dotsc,\xi_{n_{2}},
\xi_{m_{2}-\varkappa_{nk}+1},\dotsc,\xi_{(n-1)K+k}) \\
\times& \, \sum_{\pmb{\sigma} \in \mathfrak{S}_{m_{2}+1}} \operatorname{sgn}
(\pmb{\pmb{\sigma}}) \tfrac{1}{\dfrac{1}{\xi_{\sigma (0)}^{0}}
(\xi_{\sigma (1)}-\alpha_{p_{1}})^{1} \dotsb (\xi_{\sigma (l_{1})}-
\alpha_{p_{1}})^{l_{1}} \dotsb (\xi_{\sigma (n_{1})}-
\alpha_{p_{\mathfrak{s}-1}})^{1} \dotsb (\xi_{\sigma (n_{2})}-
\alpha_{p_{\mathfrak{s}-1}})^{l_{\mathfrak{s}-1}} \dfrac{1}{(
\xi_{\sigma (m_{2}-\varkappa_{nk}+1)})^{1}} \dotsb \dfrac{1}{\xi_{\sigma 
(m_{2})}^{\varkappa_{nk}}}} \\
\times& \, \dfrac{1}{\widehat{\phi}_{0}(\xi_{\sigma (0)}) \widehat{\phi}_{0}
(\xi_{\sigma (1)}) \dotsb \widehat{\phi}_{0}(\xi_{\sigma (l_{1})}) \dotsb 
\widehat{\phi}_{0}(\xi_{\sigma (n_{1})}) \dotsb \widehat{\phi}_{0}(\xi_{\sigma 
(n_{2})}) \widehat{\phi}_{0}(\xi_{\sigma (m_{2}-\varkappa_{nk}+1)}) \dotsb 
\widehat{\phi}_{0}(\xi_{\sigma (m_{2})})} \\
=& \, \dfrac{1}{(m_{2} \! + \! 1)!} \underbrace{\int_{\mathbb{R}} \int_{
\mathbb{R}} \dotsb \int_{\mathbb{R}}}_{(n-1)K+k+1} \md \widetilde{\mu}(\xi_{0}) 
\, \md \widetilde{\mu}(\xi_{1}) \, \dotsb \, \md \widetilde{\mu}(\xi_{l_{1}}) \, 
\dotsb \, \md \widetilde{\mu}(\xi_{n_{1}}) \, \dotsb \, \md \widetilde{\mu}
(\xi_{n_{2}}) \, \md \widetilde{\mu}(\xi_{m_{2}-\varkappa_{nk}+1}) \, \dotsb \\
\times& \, \dotsb \, \md \widetilde{\mu}(\xi_{m_{2}}) \, \mathbb{G}^{\spcheck}
(\xi_{0},\xi_{1},\dotsc,\xi_{l_{1}},\dotsc,\xi_{n_{1}},\dotsc,\xi_{n_{2}},
\xi_{m_{2}-\varkappa_{nk}+1},\dotsc,\xi_{(n-1)K+k}) \\
\times& \, 
\left\lvert

\right\rvert \\
=& \, \dfrac{1}{(m_{2} \! + \! 1)!} \underbrace{\int_{\mathbb{R}} 
\int_{\mathbb{R}} \dotsb \int_{\mathbb{R}}}_{(n-1)K+k+1} \md \widetilde{\mu}
(\xi_{0}) \, \md \widetilde{\mu}(\xi_{1}) \, \dotsb \, \md \widetilde{\mu}
(\xi_{l_{1}}) \, \dotsb \, \md \widetilde{\mu}(\xi_{n_{1}}) \, \dotsb \, \md 
\widetilde{\mu}(\xi_{n_{2}}) \, \md \widetilde{\mu}(\xi_{m_{2}-\varkappa_{nk}+1}) 
\, \dotsb \\
\times& \, \dotsb \, \md \widetilde{\mu}(\xi_{m_{2}}) \, \dfrac{\mathbb{G}^{\spcheck}
(\xi_{0},\xi_{1},\dotsc,\xi_{l_{1}},\dotsc,\xi_{n_{1}},\dotsc,\xi_{n_{2}},
\xi_{m_{2}-\varkappa_{nk}+1},\dotsc,\xi_{(n-1)K+k})}{(\widehat{\phi}_{0}
(\xi_{0}) \widehat{\phi}_{0}(\xi_{1}) \dotsb \widehat{\phi}_{0}(\xi_{l_{1}}) 
\dotsb \widehat{\phi}_{0}(\xi_{n_{1}}) \dotsb \widehat{\phi}_{0}(\xi_{n_{2}}) 
\widehat{\phi}_{0}(\xi_{m_{2}-\varkappa_{nk}+1}) \dotsb \widehat{\phi}_{0}
(\xi_{m_{2}}))^{2}} \\
\times& \, 
\underbrace{\left\vert

\right\vert}_{= \, \mathbb{G}^{\spcheck}(\xi_{0},\xi_{1},\dotsc,\xi_{l_{1}},
\dotsc,\xi_{n_{1}},\dotsc,\xi_{n_{2}},\xi_{m_{2}-\varkappa_{nk}+1},\dotsc,
\xi_{(n-1)K+k})} \\
=& \, \dfrac{1}{(m_{2} \! + \! 1)!} \underbrace{\int_{\mathbb{R}} \int_{
\mathbb{R}} \dotsb \int_{\mathbb{R}}}_{(n-1)K+k+1} \md \widetilde{\mu}
(\xi_{0}) \, \md \widetilde{\mu}(\xi_{1}) \, \dotsb \, \md \widetilde{\mu}
(\xi_{l_{1}}) \, \dotsb \, \md \widetilde{\mu}(\xi_{n_{1}}) \, \dotsb \, 
\md \widetilde{\mu}(\xi_{n_{2}}) \, \md \widetilde{\mu}(\xi_{m_{2}-
\varkappa_{nk}+1}) \, \dotsb \\
\times& \, \dotsb \, \md \widetilde{\mu}(\xi_{m_{2}}) \left(\dfrac{\mathbb{G}^{\spcheck}
(\xi_{0},\xi_{1},\dotsc,\xi_{l_{1}},\dotsc,\xi_{n_{1}},\dotsc,\xi_{n_{2}},
\xi_{m_{2}-\varkappa_{nk}+1},\dotsc,\xi_{(n-1)K+k})}{\widehat{\phi}_{0}
(\xi_{0}) \widehat{\phi}_{0}(\xi_{1}) \dotsb \widehat{\phi}_{0}(\xi_{l_{1}}) 
\dotsb \widehat{\phi}_{0}(\xi_{n_{1}}) \dotsb \widehat{\phi}_{0}(\xi_{n_{2}}) 
\widehat{\phi}_{0}(\xi_{m_{2}-\varkappa_{nk}+1}) \dotsb \widehat{\phi}_{0}
(\xi_{m_{2}})} \right)^{2};
\end{align*}
but, noting the determinantal factorisation
\begin{equation*}
\mathbb{G}^{\spcheck}(\xi_{0},\xi_{1},\dotsc,\xi_{l_{1}},\dotsc,\xi_{n_{1}},
\dotsc,\xi_{n_{2}},\xi_{m_{2}-\varkappa_{nk}+1},\dotsc,\xi_{(n-1)K+k}) \! 
:= \! \overset{\infty}{\mathcal{V}}_{1}(\xi_{0},\xi_{1},\dotsc,\xi_{(n-1)K+k}) 
\mathbb{D}^{\spcheck},
\end{equation*}
where
\begin{equation*}
\overset{\infty}{\mathcal{V}}_{1}(\xi_{0},\xi_{1},\dotsc,\xi_{(n-1)K+k}) 
\! = \left\lvert
 
\right), \label{eq36}
\end{equation}
with {}\footnote{Note the convention $\sum_{m=1}^{0} \boldsymbol{\ast} \! := \! 0$.}
\begin{equation*}
(\mathscr{A}^{\spcheck}(0))_{1j} \! = \! 
\begin{cases}
\mathfrak{a}_{j-1}^{\spcheck}(\vec{\bm{\alpha}}), &\text{$j \! = \! 1,2,
\dotsc,(n \! - \! 1)K \! + \! k \! - \! \varkappa_{nk} \! + \! 1$,} \\
0, &\text{$j \! = \! (n \! - \! 1)K \! + \! k \! - \! \varkappa_{nk} 
\! + \! 2,\dotsc,(n \! - \! 1)K \! + \! k \! + \! 1$,}
\end{cases}
\end{equation*}
for $r \! = \! 1,2,\dotsc,\mathfrak{s} \! - \! 1$, with $l_{r} \! = \! 
\varkappa_{nk \tilde{k}_{r}}$ and $\sum_{m=1}^{\mathfrak{s}-1}l_{m} 
\! = \! (n \! - \! 1)K \! + \! k \! - \! \varkappa_{nk}$,
\begin{align*}
i(r)=& \, 2 \! + \! \sum_{m=1}^{r-1}l_{m},3 \! + \! \sum_{m=1}^{r-1}l_{m},
\dotsc,1 \! + \! l_{r} \! + \! \sum_{m=1}^{r-1}l_{m}, \qquad \qquad 
q(i(r),r) \! = \! 1,2,\dotsc,l_{r}, \\
(\mathscr{A}^{\spcheck}(r))_{i(r)j(r)}=& \, 
\begin{cases}
\dfrac{(-1)^{q(i(r),r)}}{\prod_{m=0}^{q(i(r),r)-1}(l_{r} \! - \! m)} 
\left(\dfrac{\partial}{\partial \alpha_{p_{r}}} \right)^{q(i(r),r)} 
\mathfrak{a}_{j(r)-1}^{\spcheck}(\vec{\bm{\alpha}}), &\text{$j(r) \! = 
\! 1,2,\dotsc,(n \! - \! 1)K \! + \! k \! - \! \varkappa_{nk} \! - \! 
q(i(r),r) \! + \! 1$,} \\
0, &\text{$j(r) \! = \! (n \! - \! 1)K \! + \! k \! - \! \varkappa_{nk} 
\! - \! q(i(r),r) \! + \! 2,\dotsc,(n \! - \! 1)K \! + \! k \! + \! 1$,}
\end{cases}
\end{align*}
and,\footnote{Notation such as, for example, ``for $r \! = \! 1,2,\dotsc,
\mathfrak{s} \! - \! 1$, with $l_{r} \! = \! \varkappa_{nk \tilde{k}_{r}}$ 
and $\sum_{m=1}^{\mathfrak{s}-1}l_{m} \! = \! (n \! - \! 1)K \! + \! k \! - 
\! \varkappa_{nk}$, $i(r) \! = \! 2 \! + \! \sum_{m=1}^{r-1}l_{m},3 \! + \! 
\sum_{m=1}^{r-1}l_{m},\dotsc,1 \! + \! l_{r} \! + \! \sum_{m=1}^{r-1}l_{m}$, 
$q(i(r),r) \! = \! 1,2,\dotsc,l_{r}$,'' may lead to confusion; it should be 
understood as follows: for $r \! = \! 1$, the running index $i(1)$ takes the 
$l_{1}$ $(= \! \varkappa_{nk \tilde{k}_{1}})$ values $i(1) \! = \! 2,3,\dotsc,
1 \! + \! l_{1}$, whilst the corresponding running index $q(i(1),1)$ takes 
the $l_{1}$ values $q(i(1),1) \! = \! 1,2,\dotsc,l_{1}$, where the latter 
pair of running indices, namely, $i(1)$ and $q(i(1),1)$, are sequenced in 
lexicographic unison, that is, when $i(1) \! = \! 2$ $\Rightarrow$ $q(i(1),
1) \! = \! 1$, when $i(1) \! = \! 3$ $\Rightarrow$ $q(i(1),1) \! = \! 2$, 
$\dotsc$, and when $i(1) \! = \! 1 \! + \! l_{1}$ $\Rightarrow$ $q(i(1),1) 
\! = \! l_{1}$, then, for $r \! = \! 2$, the running index $i(2)$ takes the 
$l_{2}$ $(= \! \varkappa_{nk \tilde{k}_{2}})$ values $i(2) \! = \! 2 \! + 
\! l_{1},3 \! + \! l_{1},\dotsc,1 \! + \! l_{1} \! + \! l_{2}$, whilst the 
corresponding running index $q(i(2),2)$ takes the $l_{2}$ values $q(i(2),2) 
\! = \! 1,2,\dotsc,l_{2}$, where the latter pair of running indices, namely, 
$i(2)$ and $q(i(2),2)$, are sequenced in lexicographic unison, that is, when 
$i(2) \! = \! 2 \! + \! l_{1}$ $\Rightarrow$ $q(i(2),2) \! = \! 1$, when $i(2) 
\! = \! 3 \! + \! l_{1}$ $\Rightarrow$ $q(i(2),2) \! = \! 2$, $\dotsc$, and 
when $i(2) \! = \! 1 \! + \! l_{1} \! + \! l_{2}$ $\Rightarrow$ $q(i(2),2) \! 
= \! l_{2}$, etc. All the notational analogues/variants appearing throughout 
the remainder of this monograph should, \emph{mutatis mutandis}, be 
understood in the above-defined sense.} for $r \! = \! \mathfrak{s}$, 
with $l_{\mathfrak{s}} \! = \! \varkappa_{nk}$,
\begin{align*}
i(\mathfrak{s})=& \, (n \! - \! 1)K \! + \! k \! - \! \varkappa_{nk} 
\! + \! 2,(n \! - \! 1)K \! + \! k \! - \! \varkappa_{nk} \! + \! 3,
\dotsc,(n \! - \! 1)K \! + \! k \! + \! 1, \quad \quad q(i(\mathfrak{s}),
\mathfrak{s}) \! = \! 1,2,\dotsc,\varkappa_{nk}, \\
(\mathscr{A}^{\spcheck}(\mathfrak{s}))_{i(\mathfrak{s})j(\mathfrak{s})}=& \, 
\begin{cases}
0, &\text{$j(\mathfrak{s}) \! = \! 1,\dotsc,q(i(\mathfrak{s}),
\mathfrak{s})$,} \\
\mathfrak{a}_{j(\mathfrak{s})-q(i(\mathfrak{s}),\mathfrak{s})-1}^{\spcheck}
(\vec{\bm{\alpha}}), &\text{$j(\mathfrak{s}) \! = \! q(i(\mathfrak{s}),
\mathfrak{s}) \! + \! 1,\dotsc,i(\mathfrak{s})$,} \\
0, &\text{$j(\mathfrak{s}) \! = \! i(\mathfrak{s}) \! + \! 1,\dotsc,
(n \! - \! 1)K \! + \! k \! + \! 1$,}
\end{cases}
\end{align*}
where, for $\tilde{m}_{1} \! = \! 0,1,\dotsc,(n \! - \! 1)K \! + \! k \! - \! 
\varkappa_{nk}$,
\begin{equation*}
\mathfrak{a}_{\tilde{m}_{1}}^{\spcheck}(\vec{\bm{\alpha}}) := \, 
\mathlarger{\sum_{\underset{\underset{\sum_{m=1}^{\mathfrak{s}-1}
i_{m}=(n-1)K+k-\varkappa_{nk}-\tilde{m}_{1}}{p \in \lbrace 1,2,\dotsc,
\mathfrak{s}-1 \rbrace}}{i_{p}=0,1,\dotsc,l_{p}}}}(-1)^{(n-1)K+k-
\varkappa_{nk}-\tilde{m}_{1}} \prod_{j=1}^{\mathfrak{s}-1} 
\binom{l_{j}}{i_{j}} \prod_{m=1}^{\mathfrak{s}-1}(\alpha_{p_{m}})^{i_{m}},
\end{equation*}
it follows that, for $n \! \in \! \mathbb{N}$ and $k \! \in \! \lbrace 
1,2,\dotsc,K \rbrace$ such that $\alpha_{p_{\mathfrak{s}}} \! := \! 
\alpha_{k} \! = \! \infty$,
\begin{align}
\hat{\mathfrak{c}}_{N_{\infty}} =& \, \dfrac{1}{(m_{2} \! + \! 1)!} 
\underbrace{\int_{\mathbb{R}} \int_{\mathbb{R}} \dotsb \int_{\mathbb{R}}}_{(n-
1)K+k+1} \md \widetilde{\mu}(\xi_{0}) \, \md \widetilde{\mu}(\xi_{1}) \, \dotsb 
\, \md \widetilde{\mu}(\xi_{l_{1}}) \, \dotsb \, \md \widetilde{\mu}(\xi_{n_{1}}) 
\, \dotsb \, \md \widetilde{\mu}(\xi_{n_{2}}) \, \md \widetilde{\mu}(\xi_{m_{2}-
\varkappa_{nk}+1}) \, \dotsb \nonumber \\
\times& \, \dotsb \, \md \widetilde{\mu}(\xi_{m_{2}}) \dfrac{(\mathbb{D}^{\spcheck})^{2} 
\prod_{\underset{j<i}{i,j=0}}^{(n-1)K+k}(\xi_{i} \! - \! \xi_{j})^{2}}{(
\widehat{\phi}_{0}(\xi_{0}) \widehat{\phi}_{0}(\xi_{1}) \dotsb \widehat{
\phi}_{0}(\xi_{l_{1}}) \dotsb \widehat{\phi}_{0}(\xi_{n_{1}}) \dotsb 
\widehat{\phi}_{0}(\xi_{n_{2}}) \widehat{\phi}_{0}(\xi_{m_{2}-
\varkappa_{nk}+1}) \dotsb \widehat{\phi}_{0}(\xi_{(n-1)K+k}))^{2}} \nonumber \\
=& \, \dfrac{(\mathbb{D}^{\spcheck})^{2}}{((n \! - \! 1)K \! + \! k \! + \! 1)!} 
\underbrace{\int_{\mathbb{R}} \int_{\mathbb{R}} \dotsb \int_{\mathbb{R}}}_{
(n-1)K+k+1} \md \widetilde{\mu}(\xi_{0}) \, \md \widetilde{\mu}(\xi_{1}) \, 
\dotsb \, \md \widetilde{\mu}(\xi_{l_{1}}) \, \dotsb \, \md \widetilde{\mu}
(\xi_{n_{1}}) \, \dotsb \, \md \widetilde{\mu}(\xi_{n_{2}}) \nonumber \\
\times& \, \md \widetilde{\mu}(\xi_{(n-1)K+k-\varkappa_{nk}+1}) \, \dotsb 
\, \md \widetilde{\mu}(\xi_{(n-1)K+k}) \left(\prod_{m=0}^{(n-1)K+k} 
\widehat{\phi}_{0}(\xi_{m}) \right)^{-2} \, \, \prod_{\substack{i,j=0\\j<i}}^{
(n-1)K+k}(\xi_{i} \! - \! \xi_{j})^{2} \nonumber \\
=& \, \dfrac{(\mathbb{D}^{\spcheck})^{2}}{((n \! - \! 1)K \! + \! k \! + \! 1)!} 
\underbrace{\int_{\mathbb{R}} \int_{\mathbb{R}} \dotsb \int_{\mathbb{R}}}_{
(n-1)K+k+1} \md \widetilde{\mu}(\xi_{0}) \, \md \widetilde{\mu}(\xi_{1}) \, 
\dotsb \, \md \widetilde{\mu}(\xi_{l_{1}}) \, \dotsb \, \md \widetilde{\mu}
(\xi_{n_{1}}) \, \dotsb \, \md \widetilde{\mu}(\xi_{n_{2}}) \nonumber \\
\times& \, \md \widetilde{\mu}(\xi_{(n-1)K+k-\varkappa_{nk}+1}) \, \dotsb 
\, \md \widetilde{\mu}(\xi_{(n-1)K+k}) \prod_{\substack{i,j=0\\j<i}}^{(n-1)K
+k}(\xi_{i} \! - \! \xi_{j})^{2} \left(\prod_{m=0}^{(n-1)K+k} \prod_{q=1}^{
\mathfrak{s}-1}(\xi_{m} \! - \! \alpha_{p_{q}})^{\varkappa_{nk \tilde{k}_{q}}} 
\right)^{-2} \, \Rightarrow \nonumber \\
\hat{\mathfrak{c}}_{N_{\infty}} =& \, \dfrac{(\mathbb{D}^{\spcheck})^{2}}{
((n \! - \! 1)K \! + \! k \! + \! 1)!} \underbrace{\int_{\mathbb{R}} 
\int_{\mathbb{R}} \dotsb \int_{\mathbb{R}}}_{(n-1)K+k+1} \md \widetilde{\mu}
(\xi_{0}) \, \md \widetilde{\mu}(\xi_{1}) \, \dotsb \, \md \widetilde{\mu}
(\xi_{(n-1)K+k}) \, \prod_{\substack{i,j=0\\j<i}}^{(n-1)K+k}(\xi_{i} \! - \! 
\xi_{j})^{2} \nonumber \\
\times& \, \left(\prod_{m=0}^{(n-1)K+k} \prod_{q=1}^{\mathfrak{s}-1}
(\xi_{m} \! - \! \alpha_{p_{q}})^{\varkappa_{nk \tilde{k}_{q}}} \right)^{-2} 
\quad (> \! 0). \label{eq37}
\end{align}

Now, $\hat{\mathfrak{c}}_{D_{\infty}}$ is studied. For $n \! \in \! 
\mathbb{N}$ and $k \! \in \! \lbrace 1,2,\dotsc,K \rbrace$ such 
that $\alpha_{p_{\mathfrak{s}}} \! := \! \alpha_{k} \! = \! \infty$, 
introduce the following notation:
\begin{equation*}
\psi_{0}(z) \! := \! \prod_{m=1}^{\mathfrak{s}-1}
(z \! - \! \alpha_{p_{m}})^{l_{m}} \! =: \sum_{j=0}^{(n-1)K+k-1} 
\mathfrak{a}^{\sphat}_{j,0}z^{j},
\end{equation*}
and (with some abuse of notation), for $r \! = \! 1,\dotsc,\mathfrak{s} 
\! - \! 1,\mathfrak{s}$, $q(r) \! = \! \sum_{i=1}^{r-1}l_{i} \! + \! 1,
\sum_{i=1}^{r-1}l_{i} \! + \! 2,\dotsc,\sum_{i=1}^{r-1}l_{i} \! + \! l_{r} 
\! - \! \delta_{r \mathfrak{s}}$, and $m(r) \! = \! 1,2,\dotsc,l_{r} \! 
- \! \delta_{r \mathfrak{s}}$,
\begin{equation*}
\psi_{q(r)}(z) \! := \! \left\lbrace \psi_{0}(z)
(z \! - \! \alpha_{p_{r}})^{-m(r)(1-\delta_{r \mathfrak{s}})}
z^{m(r) \delta_{r \mathfrak{s}}} \! =: \sum_{j=0}^{(n-1)K+k-1} 
\mathfrak{a}^{\sphat}_{j,q(r)}z^{j} \right\rbrace;
\end{equation*}
e.g., for $r \! = \! 1$, the notation $\psi_{q(1)}(z) \! := \! \lbrace 
\psi_{0}(z)(z \! - \! \alpha_{p_{1}})^{-m(1)} \! =: \! \sum_{j=0}^{(n-1)K+k-1} 
\mathfrak{a}^{\sphat}_{j,q(1)}z^{j} \rbrace$, $q(1) \! = \! 1,2,\dotsc,
l_{1}$, $m(1) \! = \! 1,2,\dotsc,l_{1}$, denotes the (set of) $l_{1} \! = \! 
\varkappa_{nk \tilde{k}_{1}}$ functions
\begin{gather*}
\psi_{1}(z) \! = \! \dfrac{\psi_{0}(z)}{z \! - \! \alpha_{p_{1}}} \! =: 
\sum_{j=0}^{(n-1)K+k-1} \mathfrak{a}^{\sphat}_{j,1}z^{j}, \, \psi_{2}
(z) \! = \! \dfrac{\psi_{0}(z)}{(z \! - \! \alpha_{p_{1}})^{2}} 
\! =: \sum_{j=0}^{(n-1)K+k-1} \mathfrak{a}^{\sphat}_{j,2}z^{j}, \, 
\dotsc, \, \psi_{l_{1}}(z) \! = \! 
\dfrac{\psi_{0}(z)}{(z \! - \! \alpha_{p_{1}})^{l_{1}}} \! =: 
\sum_{j=0}^{(n-1)K+k-1} \mathfrak{a}^{\sphat}_{j,l_{1}}z^{j},
\end{gather*}
for $r \! = \! \mathfrak{s} \! - \! 1$, the notation $\psi_{q(\mathfrak{s}
-1)}(z) \! := \! \lbrace \psi_{0}(z)(z \! - \! \alpha_{p_{\mathfrak{s}
-1}})^{-m(\mathfrak{s}-1)} \! =: \! \sum_{j=0}^{(n-1)K+k-1} 
\mathfrak{a}^{\sphat}_{j,q(\mathfrak{s}-1)}z^{j} \rbrace$, $q(\mathfrak{s} 
\! - \! 1) \! = \! l_{1} \! + \! \dotsb \! + \! l_{\mathfrak{s}-2} \! + \! 1,
l_{1} \! + \! \dotsb \! + \! l_{\mathfrak{s}-2} \! + \! 2,\dotsc,l_{1} \! + 
\! \dotsb \! + \! l_{\mathfrak{s}-2} \! + \! l_{\mathfrak{s}-1} \! = \! 
(n \! - \! 1)K \! + \! k \! - \! \varkappa_{nk}$, $m(\mathfrak{s} \! 
- \! 1) \! = \! 1,2,\dotsc,l_{\mathfrak{s}-1}$, denotes the (set of) 
$l_{\mathfrak{s}-1} \! = \! \varkappa_{nk \tilde{k}_{\mathfrak{s}-1}}$ 
functions
\begin{gather*}
\psi_{l_{1}+ \dotsb +l_{\mathfrak{s}-2}+1}(z) \! = \! 
\dfrac{\psi_{0}(z)}{z \! - \! \alpha_{p_{\mathfrak{s}-1}}} \! =: 
\sum_{j=0}^{(n-1)K+k-1} \mathfrak{a}^{\sphat}_{j,l_{1}+\dotsb 
+l_{\mathfrak{s}-2}+1}z^{j}, \, \psi_{l_{1}+\dotsb +l_{\mathfrak{s}-2}+2}(z) 
\! = \! \dfrac{\psi_{0}(z)}{(z \! - \! \alpha_{p_{\mathfrak{s}-1}})^{2}} \! 
=: \sum_{j=0}^{(n-1)K+k-1} \mathfrak{a}^{\sphat}_{j,l_{1}+\dotsb 
+l_{\mathfrak{s}-2}+2}z^{j}, \, \dotsc \\
\dotsc, \, \psi_{(n-1)K+k-\varkappa_{nk}}(z) \! = \! 
\dfrac{\psi_{0}(z)}{(z \! - \! \alpha_{p_{\mathfrak{s}-1}})^{
l_{\mathfrak{s}-1}}} \! =: \sum_{j=0}^{(n-1)K+k-1} 
\mathfrak{a}^{\sphat}_{j,(n-1)K+k-\varkappa_{nk}}z^{j},
\end{gather*}
etc., and, for $r \! = \! \mathfrak{s}$, the notation $\psi_{q(\mathfrak{s})}
(z) \! := \! \lbrace \psi_{0}(z)z^{m(\mathfrak{s})} \! =: \! 
\sum_{j=0}^{(n-1)K+k-1} \mathfrak{a}^{\sphat}_{j,q(\mathfrak{s})}z^{j} 
\rbrace$, $q(\mathfrak{s}) \! = \! (n \! - \! 1)K \! + \! k \! - \! 
\varkappa_{nk} \! + \! 1,(n \! - \! 1)K \! + \! k \! - \! \varkappa_{nk} \! + 
\! 2,\dotsc,(n \! - \! 1)K \! + \! k \! - \! 1$, $m(\mathfrak{s}) \! = \! 1,2,
\dotsc,l_{\mathfrak{s}} \! - \! 1$, denotes the (set of) $l_{\mathfrak{s}} 
\! - \! 1 \! = \! \varkappa_{nk} \! - \! 1$ functions
\begin{gather*}
\psi_{(n-1)K+k-\varkappa_{nk}+1}(z) \! = \! \psi_{0}(z)z \! =: 
\sum_{j=0}^{(n-1)K+k-1} \mathfrak{a}^{\sphat}_{j,(n-1)K+k-\varkappa_{nk}+1}
z^{j}, \, \psi_{(n-1)K+k-\varkappa_{nk}+2}(z) \! = \! \psi_{0}(z)z^{2} \\
=: \sum_{j=0}^{(n-1)K+k-1} \mathfrak{a}^{\sphat}_{j,(n-1)K+k-\varkappa_{nk}+2}
z^{j}, \, \dotsc, \, \psi_{(n-1)K+k-1}(z) \! = \! \psi_{0}(z)
z^{\varkappa_{nk}-1} \! =: \sum_{j=0}^{(n-1)K+k-1} 
\mathfrak{a}^{\sphat}_{j,(n-1)K+k-1}z^{j}.
\end{gather*}
(Note: $\# \lbrace \psi_{0}(z)(z \! - \! \alpha_{p_{r}})^{-m(r)
(1-\delta_{r \mathfrak{s}})}z^{m(r) \delta_{r \mathfrak{s}}} \rbrace \! = \! 
l_{r} \! - \! \delta_{r \mathfrak{s}}$, $r \! = \! 1,\dotsc,\mathfrak{s} \! - 
\! 1,\mathfrak{s}$, and $\# \cup_{r=1}^{\mathfrak{s}} \lbrace \psi_{0}(z)
(z \! - \! \alpha_{p_{r}})^{-m(r)(1-\delta_{r \mathfrak{s}})}z^{m(r) 
\delta_{r \mathfrak{s}}} \rbrace \! = \! \sum_{r=1}^{\mathfrak{s}}
l_{r} \! - \! 1 \! = \! (n \! - \! 1)K \! + \! k \! - \! 1$.) One notes that, for 
$n \! \in \! \mathbb{N}$ and $k \! \in \! \lbrace 1,2,\dotsc,K \rbrace$ 
such that $\alpha_{p_{\mathfrak{s}}} \! := \! \alpha_{k} \! = \! \infty$, the 
$l_{1} \! + \! \dotsb \! + \! l_{\mathfrak{s}-1} \! + \! l_{\mathfrak{s}} \! = 
\! (n \! - \! 1)K \! + \! k$ functions $\psi_{0}(z),\psi_{1}(z),\dotsc,\psi_{l_{1}}
(z),\dotsc,\psi_{l_{1}+\dotsb +l_{\mathfrak{s}-2}+1}(z),\dotsc,\psi_{(n-1)K+k
-\varkappa_{nk}}(z),\psi_{(n-1)K+k-\varkappa_{nk}+1}(z),\dotsc,\linebreak[4] 
\psi_{(n-1)K+k-1}(z)$ are linearly independent on $\mathbb{R}$, 
that is, for $z \! \in \! \mathbb{R}$, $\sum_{j=0}^{(n-1)K+k-1} 
\mathfrak{c}^{\sphat}_{j} \psi_{j}(z) \! = \! 0$ $\Rightarrow$ (via a 
Vandermonde-type argument; see the $((n \! - \! 1)K \! + \! k) \times 
((n \! - \! 1)K \! + \! k)$ non-zero determinant $\mathbb{D}^{\sphat}$ 
in Equation~\eqref{eq38} below) $\mathfrak{c}^{\sphat}_{j} \! = \! 0$, 
$j \! = \! 0,1,\dotsc,(n \! - \! 1)K \! + \! k \! - \! 1$. For $n \! \in \! 
\mathbb{N}$ and $k \! \in \! \lbrace 1,2,\dotsc,K \rbrace$ such that 
$\alpha_{p_{\mathfrak{s}}} \! := \! \alpha_{k} \! = \! \infty$, let 
$\mathfrak{S}_{(n-1)K+k}$ denote the $((n \! - \! 1)K \! + \! k)!$ 
permutations of $\lbrace 0,1,\dotsc,(n \! - \! 1)K \! + \! k \! - \! 1 
\rbrace$. Using the above notation and the multi-linearity property 
of the determinant, one proceeds to study, for $n \! \in \! \mathbb{N}$ 
and $k \! \in \! \lbrace 1,2,\dotsc,K \rbrace$ such that 
$\alpha_{p_{\mathfrak{s}}} \! := \! \alpha_{k} \! = \! \infty$, 
$\hat{\mathfrak{c}}_{D_{\infty}}$:
\begin{align*}
\hat{\mathfrak{c}}_{D_{\infty}} =& \, \underbrace{\int_{\mathbb{R}} 
\int_{\mathbb{R}} \dotsb \int_{\mathbb{R}}}_{(n-1)K+k} \md \widetilde{\mu}
(\xi_{0}) \, \md \widetilde{\mu}(\xi_{1}) \, \dotsb \, \md \widetilde{\mu}
(\xi_{l_{1}}) \, \dotsb \, \md \widetilde{\mu}(\xi_{n_{1}}) \, \dotsb \, 
\md \widetilde{\mu}(\xi_{n_{2}}) \, \md \widetilde{\mu}(\xi_{m_{2}-
\varkappa_{nk}+1}) \, \dotsb \, \md \widetilde{\mu}(\xi_{m_{1}}) \\
\times& \, \dfrac{1}{\dfrac{1}{\xi_{0}^{0}}(\xi_{1} \! - \! \alpha_{p_{1}})^{
1} \dotsb (\xi_{l_{1}} \! - \! \alpha_{p_{1}})^{l_{1}} \dotsb (\xi_{n_{1}} \! 
- \! \alpha_{p_{\mathfrak{s}-1}})^{1} \dotsb (\xi_{n_{2}} \! - \! \alpha_{p_{
\mathfrak{s}-1}})^{l_{\mathfrak{s}-1}} \dfrac{1}{(\xi_{m_{2}-\varkappa_{nk}
+1})^{1}} \dotsb \dfrac{1}{\xi_{m_{1}}^{\varkappa_{nk}-1}}} \\
\times& \, \dfrac{1}{\psi_{0}(\xi_{0}) \psi_{0}(\xi_{1}) \dotsb \psi_{0}
(\xi_{l_{1}}) \dotsb \psi_{0}(\xi_{n_{1}}) \dotsb \psi_{0}(\xi_{n_{2}}) 
\psi_{0}(\xi_{m_{2}-\varkappa_{nk}+1}) \dotsb \psi_{0}(\xi_{m_{1}})} \\
\times& \, 
\underbrace{\left\vert
 
\right\vert}_{=: \, \mathbb{G}^{\sphat}(\xi_{0},\xi_{1},\dotsc,\xi_{l_{1}},
\dotsc,\xi_{n_{1}},\dotsc,\xi_{n_{2}},\xi_{m_{2}-\varkappa_{nk}+1},\dotsc,
\xi_{(n-1)K+k-1})} \\
=& \, \dfrac{1}{m_{2}!} \sum_{\pmb{\sigma} \in \mathfrak{S}_{m_{2}}} 
\underbrace{\int_{\mathbb{R}} \int_{\mathbb{R}} \dotsb \int_{\mathbb{R}}}_{
(n-1)K+k} \md \widetilde{\mu}(\xi_{\sigma (0)}) \, \md \widetilde{\mu}
(\xi_{\sigma (1)}) \, \dotsb \, \md \widetilde{\mu}(\xi_{\sigma (l_{1})}) \, 
\dotsb \, \md \widetilde{\mu}(\xi_{\sigma (n_{1})}) \, \dotsb \, \md 
\widetilde{\mu}(\xi_{\sigma (n_{2})}) \, \dotsb \\
\times& \, \tfrac{\dotsb \, \md \widetilde{\mu}(\xi_{\sigma (m_{2}-\varkappa_{nk}+1)}) 
\, \dotsb \, \md \widetilde{\mu}(\xi_{\sigma (m_{1})})}{\dfrac{1}{\xi_{\sigma (0)}^{0}}
(\xi_{\sigma (1)}-\alpha_{p_{1}})^{1} \dotsb (\xi_{\sigma (l_{1})}-
\alpha_{p_{1}})^{l_{1}} \dotsb (\xi_{\sigma (n_{1})}-\alpha_{p_{\mathfrak{s}
-1}})^{1} \dotsb (\xi_{\sigma (n_{2})}-\alpha_{p_{\mathfrak{s}-1}})^{
l_{\mathfrak{s}-1}} \dfrac{1}{(\xi_{\sigma (m_{2}-\varkappa_{nk}+1)})^{1}} 
\dotsb \dfrac{1}{\xi_{\sigma (m_{1})}^{\varkappa_{nk}-1}}} \\
\times& \, \dfrac{1}{\psi_{0}(\xi_{\sigma (0)}) \psi_{0}(\xi_{\sigma (1)}) 
\dotsb \psi_{0}(\xi_{\sigma (l_{1})}) \dotsb \psi_{0}(\xi_{\sigma (n_{1})}) 
\dotsb \psi_{0}(\xi_{\sigma (n_{2})}) \psi_{0}(\xi_{\sigma (m_{2}-
\varkappa_{nk}+1)}) \dotsb \psi_{0}(\xi_{\sigma (m_{1})})} \\
\times& \, \mathbb{G}^{\sphat}(\xi_{\sigma (0)},\xi_{\sigma (1)},\dotsc,
\xi_{\sigma (l_{1})},\dotsc,\xi_{\sigma (n_{1})},\dotsc,\xi_{\sigma (n_{2})},
\xi_{\sigma (m_{2}-\varkappa_{nk}+1)},\dotsc,\xi_{\sigma ((n-1)K+k-1)}) \\
=& \, \dfrac{1}{m_{2}!} \underbrace{\int_{\mathbb{R}} \int_{\mathbb{R}} \dotsb 
\int_{\mathbb{R}}}_{(n-1)K+k} \md \widetilde{\mu}(\xi_{0}) \, \md \widetilde{\mu}
(\xi_{1}) \, \dotsb \, \md \widetilde{\mu}(\xi_{l_{1}}) \, \dotsb \, \md \widetilde{\mu}
(\xi_{n_{1}}) \, \dotsb \, \md \widetilde{\mu}(\xi_{n_{2}}) \, \md \widetilde{\mu}
(\xi_{m_{2}-\varkappa_{nk}+1}) \, \dotsb \\
\times& \, \dotsb \, \md \widetilde{\mu}(\xi_{m_{1}}) \, \mathbb{G}^{\sphat}
(\xi_{0},\xi_{1},\dotsc,\xi_{l_{1}},\dotsc,\xi_{n_{1}},\dotsc,\xi_{n_{2}},\xi_{m_{2}
-\varkappa_{nk}+1},\dotsc,\xi_{(n-1)K+k-1}) \\
\times& \, \sum_{\pmb{\sigma} \in \mathfrak{S}_{m_{2}}} \operatorname{sgn}
(\pmb{\pmb{\sigma}}) \tfrac{1}{\dfrac{1}{\xi_{\sigma (0)}^{0}}
(\xi_{\sigma (1)}-\alpha_{p_{1}})^{1} \dotsb (\xi_{\sigma (l_{1})}-
\alpha_{p_{1}})^{l_{1}} \dotsb (\xi_{\sigma (n_{1})}-\alpha_{p_{\mathfrak{s}
-1}})^{1} \dotsb (\xi_{\sigma (n_{2})}-\alpha_{p_{\mathfrak{s}-1}})^{
l_{\mathfrak{s}-1}} \dfrac{1}{(\xi_{\sigma (m_{2}-\varkappa_{nk}+1)})^{1}} 
\dotsb \dfrac{1}{\xi_{\sigma (m_{1})}^{\varkappa_{nk}-1}}} \\
\times& \, \dfrac{1}{\psi_{0}(\xi_{\sigma (0)}) \psi_{0}(\xi_{\sigma (1)}) 
\dotsb \psi_{0}(\xi_{\sigma (l_{1})}) \dotsb \psi_{0}(\xi_{\sigma (n_{1})}) 
\dotsb \psi_{0}(\xi_{\sigma (n_{2})}) \psi_{0}(\xi_{\sigma (m_{2}-\varkappa_{
nk}+1)}) \dotsb \psi_{0}(\xi_{\sigma (m_{1})})} \\
=& \, \dfrac{1}{m_{2}!} \underbrace{\int_{\mathbb{R}} \int_{\mathbb{R}} \dotsb 
\int_{\mathbb{R}}}_{(n-1)K+k} \md \widetilde{\mu}(\xi_{0}) \, \md \widetilde{\mu}
(\xi_{1}) \, \dotsb \, \md \widetilde{\mu}(\xi_{l_{1}}) \, \dotsb \, \md \widetilde{\mu}
(\xi_{n_{1}}) \, \dotsb \, \md \widetilde{\mu}(\xi_{n_{2}}) \, \md \widetilde{\mu}
(\xi_{m_{2}-\varkappa_{nk}+1}) \, \dotsb \\
\times& \, \dotsb \, \md \widetilde{\mu}(\xi_{m_{1}}) \, \mathbb{G}^{\sphat}
(\xi_{0},\xi_{1},\dotsc,\xi_{l_{1}},\dotsc,\xi_{n_{1}},\dotsc,\xi_{n_{2}},\xi_{m_{2}-
\varkappa_{nk}+1},\dotsc,\xi_{(n-1)K+k-1}) \\
\times& \, 
\left\lvert

\right\rvert \\
=& \, \dfrac{1}{m_{2}!} \underbrace{\int_{\mathbb{R}} \int_{\mathbb{R}} \dotsb 
\int_{\mathbb{R}}}_{(n-1)K+k} \md \widetilde{\mu}(\xi_{0}) \, \md \widetilde{\mu}
(\xi_{1}) \, \dotsb \, \md \widetilde{\mu}(\xi_{l_{1}}) \, \dotsb \, \md \widetilde{\mu}
(\xi_{n_{1}}) \, \dotsb \, \md \widetilde{\mu}(\xi_{n_{2}}) \, \md \widetilde{\mu}
(\xi_{m_{2}-\varkappa_{nk}+1}) \, \dotsb \\
\times& \, \dotsb \, \md \widetilde{\mu}(\xi_{m_{1}}) \, \dfrac{\mathbb{G}^{\sphat}
(\xi_{0},\xi_{1},\dotsc,\xi_{l_{1}},\dotsc,\xi_{n_{1}},\dotsc,\xi_{n_{2}},
\xi_{m_{2}-\varkappa_{nk}+1},\dotsc,\xi_{(n-1)K+k-1})}{(\psi_{0}(\xi_{0}) 
\psi_{0}(\xi_{1}) \dotsb \psi_{0}(\xi_{l_{1}}) \dotsb \psi_{0}(\xi_{n_{1}}) 
\dotsb \psi_{0}(\xi_{n_{2}}) \psi_{0}(\xi_{m_{2}-\varkappa_{nk}+1}) \dotsb 
\psi_{0}(\xi_{m_{1}}))^{2}} \\
\times& \, 
\underbrace{\left\vert

\right\vert}_{= \, \mathbb{G}^{\sphat}(\xi_{0},\xi_{1},\dotsc,\xi_{l_{1}},
\dotsc,\xi_{n_{1}},\dotsc,\xi_{n_{2}},\xi_{m_{2}-\varkappa_{nk}+1},\dotsc,
\xi_{(n-1)K+k-1})} \\
=& \, \dfrac{1}{m_{2}!} \underbrace{\int_{\mathbb{R}} \int_{\mathbb{R}} \dotsb 
\int_{\mathbb{R}}}_{(n-1)K+k} \md \widetilde{\mu}(\xi_{0}) \, \md \widetilde{\mu}
(\xi_{1}) \, \dotsb \, \md \widetilde{\mu}(\xi_{l_{1}}) \, \dotsb \, \md \widetilde{\mu}
(\xi_{n_{1}}) \, \dotsb \, \md \widetilde{\mu}(\xi_{n_{2}}) \, \md \widetilde{\mu}
(\xi_{m_{2}-\varkappa_{nk}+1}) \, \dotsb \\
\times& \, \dotsb \, \md \widetilde{\mu}(\xi_{m_{1}}) \left(\dfrac{\mathbb{G}^{\sphat}
(\xi_{0},\xi_{1},\dotsc,\xi_{l_{1}},\dotsc,\xi_{n_{1}},\dotsc,\xi_{n_{2}},
\xi_{m_{2}-\varkappa_{nk}+1},\dotsc,\xi_{(n-1)K+k-1})}{\psi_{0}(\xi_{0}) 
\psi_{0}(\xi_{1}) \dotsb \psi_{0}(\xi_{l_{1}}) \dotsb \psi_{0}(\xi_{n_{1}}) 
\dotsb \psi_{0}(\xi_{n_{2}}) \psi_{0}(\xi_{m_{2}-\varkappa_{nk}+1}) \dotsb 
\psi_{0}(\xi_{m_{1}})} \right)^{2};
\end{align*}
but, noting the determinantal factorisation
\begin{equation*}
\mathbb{G}^{\sphat}(\xi_{0},\xi_{1},\dotsc,\xi_{l_{1}},\dotsc,\xi_{n_{1}},
\dotsc,\xi_{n_{2}},\xi_{m_{2}-\varkappa_{nk}+1},\dotsc,\xi_{(n-1)K+k-1}) 
\! := \! \overset{\infty}{\mathcal{V}}_{2}(\xi_{0},\xi_{1},\dotsc,
\xi_{(n-1)K+k-1}) \mathbb{D}^{\sphat},
\end{equation*}
where
\begin{equation*}
\overset{\infty}{\mathcal{V}}_{2}(\xi_{0},\xi_{1},\dotsc,\xi_{(n-1)K+k-1}) \! = 
\left\lvert

\right), \label{eq38}
\end{equation}
with
\begin{equation*}
(\mathscr{A}^{\sphat}(0))_{1j} \! = \! 
\begin{cases}
\mathfrak{a}_{j-1}^{\sphat}(\vec{\bm{\alpha}}), &\text{$j \! = \! 1,2,
\dotsc,(n \! - \! 1)K \! + \! k \! - \! \varkappa_{nk} \! + \! 1$,} \\
0, &\text{$j \! = \! (n \! - \! 1)K \! + \! k \! - \! \varkappa_{nk} \! 
+ \! 2,\dotsc,(n \! - \! 1)K \! + \! k$,}
\end{cases}
\end{equation*}
for $r \! = \! 1,2,\dotsc,\mathfrak{s} \! - \! 1$,
\begin{align*}
i(r)=& \, 2 \! + \! \sum_{m=1}^{r-1}l_{m},3 \! + \! \sum_{m=1}^{r-1}l_{m},
\dotsc,1 \! + \! l_{r} \! + \! \sum_{m=1}^{r-1}l_{m}, \qquad \qquad q(i(r),r) 
\! = \! 1,2,\dotsc,l_{r}, \\
(\mathscr{A}^{\sphat}(r))_{i(r)j(r)}=& \, 
\begin{cases}
\dfrac{(-1)^{q(i(r),r)}}{\prod_{m=0}^{q(i(r),r)-1}(l_{r} \! - \! m)} 
\left(\dfrac{\partial}{\partial \alpha_{p_{r}}} \right)^{q(i(r),r)} 
\mathfrak{a}_{j(r)-1}^{\sphat}(\vec{\bm{\alpha}}), &\text{$j(r) \! = \! 
1,2,\dotsc,(n \! - \! 1)K \! + \! k \! - \! \varkappa_{nk} \! - \! q(i(r),r) 
\! + \! 1$,} \\
0, &\text{$j(r) \! = \! (n \! - \! 1)K \! + \! k \! - \! \varkappa_{nk} 
\! - \! q(i(r),r) \! + \! 2,\dotsc,(n \! - \! 1)K \! + \! k$,}
\end{cases}
\end{align*}
and, for $r \! = \! \mathfrak{s}$,
\begin{align*}
i(\mathfrak{s})=& \, (n \! - \! 1)K \! + \! k \! - \! \varkappa_{nk} \! 
+ \! 2,(n \! - \! 1)K \! + \! k \! - \! \varkappa_{nk} \! + \! 3,\dotsc,
(n \! - \! 1)K \! + \! k, \quad \quad q(i(\mathfrak{s}),\mathfrak{s}) 
\! = \! 1,2,\dotsc,\varkappa_{nk} \! - \! 1, \\
(\mathscr{A}^{\sphat}(\mathfrak{s}))_{i(\mathfrak{s})j(\mathfrak{s})}=& \, 
\begin{cases}
0, &\text{$j(\mathfrak{s}) \! = \! 1,\dotsc,q(i(\mathfrak{s}),
\mathfrak{s})$,} \\
\mathfrak{a}_{j(\mathfrak{s})-q(i(\mathfrak{s}),\mathfrak{s})-1}^{\sphat}
(\vec{\bm{\alpha}}), &\text{$j(\mathfrak{s}) \! = \! q(i(\mathfrak{s}),
\mathfrak{s}) \! + \! 1,\dotsc,i(\mathfrak{s})$,} \\
0, &\text{$j(\mathfrak{s}) \! = \! i(\mathfrak{s}) \! + \! 1,\dotsc,
(n \! - \! 1)K \! + \! k$,}
\end{cases}
\end{align*}
where, for $\tilde{m}_{1} \! = \! 0,1,\dotsc,(n \! - \! 1)K \! + \! k \! - \! 
\varkappa_{nk}$, $\mathfrak{a}_{\tilde{m}_{1}}^{\sphat}(\vec{\bm{\alpha}}) 
\! = \! \mathfrak{a}_{\tilde{m}_{1}}^{\spcheck}(\vec{\bm{\alpha}})$, it 
follows that, for $n \! \in \! \mathbb{N}$ and $k \! \in \! \lbrace 1,2,
\dotsc,K \rbrace$ such that $\alpha_{p_{\mathfrak{s}}} \! := \! \alpha_{k} 
\! = \! \infty$,
\begin{align}
\hat{\mathfrak{c}}_{D_{\infty}} =& \, \dfrac{1}{m_{2}!} \underbrace{
\int_{\mathbb{R}} \int_{\mathbb{R}} \dotsb \int_{\mathbb{R}}}_{(n-1)K+k} 
\md \widetilde{\mu}(\xi_{0}) \, \md \widetilde{\mu}(\xi_{1}) \, \dotsb \, \md 
\widetilde{\mu}(\xi_{l_{1}}) \, \dotsb \, \md \widetilde{\mu}(\xi_{n_{1}}) \, 
\dotsb \, \md \widetilde{\mu}(\xi_{n_{2}}) \, \md \widetilde{\mu}(\xi_{m_{2}
-\varkappa_{nk}+1}) \, \dotsb \nonumber \\
\times& \, \dotsb \, \md \widetilde{\mu}(\xi_{m_{1}}) \dfrac{(\mathbb{D}^{\sphat} 
\,)^{2} \prod_{\underset{j<i}{i,j=0}}^{(n-1)K+k-1}(\xi_{i} \! - \! \xi_{j})^{2}}{
(\psi_{0}(\xi_{0}) \psi_{0}(\xi_{1}) \dotsb \psi_{0}(\xi_{l_{1}}) \dotsb 
\psi_{0}(\xi_{n_{1}}) \dotsb \psi_{0}(\xi_{n_{2}}) \psi_{0}(\xi_{m_{2}-
\varkappa_{nk}+1}) \dotsb \psi_{0}(\xi_{(n-1)K+k-1}))^{2}} \nonumber \\
=& \, \dfrac{(\mathbb{D}^{\sphat} \,)^{2}}{((n \! - \! 1)K \! + \! k)!} 
\underbrace{\int_{\mathbb{R}} \int_{\mathbb{R}} \dotsb \int_{\mathbb{R}}}_{
(n-1)K+k} \md \widetilde{\mu}(\xi_{0}) \, \md \widetilde{\mu}(\xi_{1}) \, \dotsb 
\, \md \widetilde{\mu}(\xi_{l_{1}}) \, \dotsb \, \md \widetilde{\mu}(\xi_{n_{1}}) 
\, \dotsb \, \md \widetilde{\mu}(\xi_{n_{2}}) \nonumber \\
\times& \, \md \widetilde{\mu}(\xi_{(n-1)K+k-\varkappa_{nk}+1}) \, \dotsb \, 
\md \widetilde{\mu}(\xi_{(n-1)K+k-1}) \left(\prod_{m=0}^{(n-1)K+k-1} 
\psi_{0}(\xi_{m}) \right)^{-2} \prod_{\substack{i,j=0\\j<i}}^{(n-1)K+k-1}
(\xi_{i} \! - \! \xi_{j})^{2} \nonumber \\
=& \, \dfrac{(\mathbb{D}^{\sphat} \,)^{2}}{((n \! - \! 1)K \! + \! k)!} 
\underbrace{\int_{\mathbb{R}} \int_{\mathbb{R}} \dotsb \int_{\mathbb{R}}}_{
(n-1)K+k} \md \widetilde{\mu}(\xi_{0}) \, \md \widetilde{\mu}(\xi_{1}) \, 
\dotsb \, \md \widetilde{\mu}(\xi_{l_{1}}) \, \dotsb \, \md \widetilde{\mu}
(\xi_{n_{1}}) \, \dotsb \, \md \widetilde{\mu}(\xi_{n_{2}}) \nonumber \\
\times& \, \md \widetilde{\mu}(\xi_{(n-1)K+k-\varkappa_{nk}+1}) \, \dotsb \, 
\md \widetilde{\mu}(\xi_{(n-1)K+k-1}) \prod_{\substack{i,j=0\\j<i}}^{(n-1)K+k-1}
(\xi_{i} \! - \! \xi_{j})^{2} \left(\prod_{m=0}^{(n-1)K+k-1} \prod_{q=1}^{\mathfrak{s}
-1}(\xi_{m} \! - \! \alpha_{p_{q}})^{\varkappa_{nk \tilde{k}_{q}}} \right)^{-2} \, 
\Rightarrow  \nonumber \\
\hat{\mathfrak{c}}_{D_{\infty}} =& \, \dfrac{(\mathbb{D}^{\sphat} \,)^{2}}{
((n \! - \! 1)K \! + \! k)!} \underbrace{\int_{\mathbb{R}} \int_{\mathbb{R}} 
\dotsb \int_{\mathbb{R}}}_{(n-1)K+k} \md \widetilde{\mu}(\xi_{0}) \, \md 
\widetilde{\mu}(\xi_{1}) \, \dotsb \, \md \widetilde{\mu}(\xi_{(n-1)K+k-1}) 
\prod_{\substack{i,j=0\\j<i}}^{(n-1)K+k-1}(\xi_{i} \! - \! \xi_{j})^{2} \nonumber \\
\times& \, \left(\prod_{m=0}^{(n-1)K+k-1} \prod_{q=1}^{\mathfrak{s}-1}
(\xi_{m} \! - \! \alpha_{p_{q}})^{\varkappa_{nk \tilde{k}_{q}}} \right)^{-2} 
\quad (> \! 0). \label{eq39}
\end{align}
Hence, for $n \! \in \! \mathbb{N}$ and $k \! \in \! \lbrace 1,2,\dotsc,K 
\rbrace$ such that $\alpha_{p_{\mathfrak{s}}} \! := \! \alpha_{k} \! = \! 
\infty$, Equations~\eqref{eq37} and~\eqref{eq39} establish the existence 
and the uniqueness of the corresponding monic MPC ORF, $\pmb{\pi}^{n}_{k} 
\colon \mathbb{N} \times \lbrace 1,2,\dotsc,K \rbrace \times \overline{
\mathbb{C}} \setminus \lbrace \alpha_{p_{1}},\alpha_{p_{2}},\dotsc,
\alpha_{p_{\mathfrak{s}}} \rbrace \! \to \! \mathbb{C}$, $(n,k,z) 
\! \mapsto \! \pmb{\pi}^{n}_{k}(z) \! = \! \mathcal{X}_{11}(z)$.

There remains, still, the question of the existence and the uniqueness 
of $\mathcal{X}_{21} \colon \mathbb{N} \times \lbrace 1,2,\dotsc,K 
\rbrace \times \overline{\mathbb{C}} \setminus \lbrace \alpha_{p_{1}},
\alpha_{p_{2}},\dotsc,\linebreak[4] 
\alpha_{p_{\mathfrak{s}}} \rbrace \! \to \! \mathbb{C}$, which necessitates 
an analysis for the elements of the second row of $\mathcal{X}(z)$, that is, 
$\mathcal{X}_{2m}(z)$, $m \! = \! 1,2$ (thus the gist of the subsequent 
calculations). Recalling that $\int_{\mathbb{R}} \xi^{m} \widetilde{w}(\xi) 
\, \md \xi \! < \! \infty$, $m \! \in \! \mathbb{N}_{0}$, and 
$\int_{\mathbb{R}}(\xi \! - \! \alpha_{k})^{-(m+1)} \widetilde{w}(\xi) \, 
\md \xi \! < \! \infty$, $\alpha_{k} \! \neq \! \infty$, $k \! \in \! \lbrace 
1,2,\dotsc,K \rbrace$, it follows, via the expansion $\tfrac{1}{z_{1}-z_{2}} 
\! = \! \sum_{i=0}^{l} \tfrac{z_{2}^{i}}{z_{1}^{i+1}} \! + \! 
\tfrac{z_{2}^{l+1}}{z_{1}^{l+1}(z_{1}-z_{2})}$, $l \! \in \! \mathbb{N}_{0}$, 
the integral representation (cf. Equation~\eqref{eq28}) $\mathcal{X}_{22}
(z) \! = \! \int_{\mathbb{R}} \tfrac{\mathcal{X}_{21}(\xi) \widetilde{w}
(\xi)}{\xi -z} \, \tfrac{\md \xi}{2 \pi \mi}$, $z \! \in \! \mathbb{C} 
\setminus \mathbb{R}$, and the second line of each of the asymptotic 
conditions~\eqref{eq26} and~\eqref{eq27}, that, for $n \! \in \! 
\mathbb{N}$ and $k \! \in \! \lbrace 1,2,\dotsc,K \rbrace$ such 
that $\alpha_{p_{\mathfrak{s}}} \! := \! \alpha_{k} \! = \! \infty$,
\begin{gather}
\int_{\mathbb{R}} \mathcal{X}_{21}(\xi) \xi^{p} \widetilde{w}(\xi) \, \md \xi 
\! = \! 0, \quad p \! = \! 0,1,\dotsc,\varkappa_{nk} \! - \! 2, \label{eq40} \\
\int_{\mathbb{R}} \mathcal{X}_{21}(\xi) \xi^{\varkappa_{nk}-1} \widetilde{w}
(\xi) \, \md \xi \! = \! -2 \pi \mi, \label{eq41} \\
\int_{\mathbb{R}} \mathcal{X}_{21}(\xi)(\xi \! - \! \alpha_{p_{q}})^{-r} 
\widetilde{w}(\xi) \, \md \xi \! = \! 0, \quad q \! = \! 1,2,\dotsc,
\mathfrak{s} \! - \! 1, \quad r \! = \! 1,2,\dotsc,
\varkappa_{nk \tilde{k}_{q}}. \label{eq42}
\end{gather}
(Note: for $n \! = \! 1$ and $k \! \in \! \lbrace 1,2,\dotsc,K \rbrace$ 
such that $\alpha_{p_{\mathfrak{s}}} \! := \! \alpha_{k} \! = \! \infty$, 
it can happen that the corresponding $\varkappa_{1k} \! < \! 2$, in 
which case, Equation~\eqref{eq40} is vacuous (of course, for $n \! 
\geqslant \! 2$, $k \! \in \! \lbrace 1,2,\dotsc,K \rbrace$, $\varkappa_{nk} 
\! \geqslant \! 2)$; moreover, if, for $n \! \in \! \mathbb{N}$ 
and $k \! \in \! \lbrace 1,2,\dotsc,K \rbrace$ such that 
$\alpha_{p_{\mathfrak{s}}} \! := \! \alpha_{k} \! = \! \infty$, $\lbrace 
\mathstrut \alpha_{k^{\prime}}, \, k^{\prime} \! \in \! \lbrace 1,2,\dotsc,
K \rbrace; \, \alpha_{k^{\prime}} \! \neq \! \alpha_{k} \! = \! \infty \rbrace 
\! = \! \varnothing$, then the corresponding Equation~\eqref{eq42} is 
vacuous; this can only occur if $n \! = \! 1$.) Via the above ordered disjoint 
partition for the repeated pole sequence $\lbrace \alpha_{1},\alpha_{2},
\dotsc,\alpha_{K} \rbrace \cup \lbrace \alpha_{1},\alpha_{2},\dotsc,
\alpha_{K} \rbrace \cup \dotsb \cup \lbrace \alpha_{1},\alpha_{2},\dotsc,
\alpha_{k} \rbrace$, and Equations~\eqref{eq40}--\eqref{eq42}, one writes, 
in the indicated order,
\begin{equation*}
\mathcal{X}_{21}(z) \! = \! \sum_{q=1}^{\mathfrak{s}-1} \sum_{r=1}^{l_{q}} 
\dfrac{\widehat{\nu}^{\raise-1.0ex\hbox{$\scriptstyle \infty$}}_{q,r}
(n,k)}{(z \! - \! \alpha_{i(q)_{k_{q}}})^{r}} \! + \! 
\sum_{m=1}^{l_{\mathfrak{s}}} 
\widehat{\nu}^{\raise-1.0ex\hbox{$\scriptstyle \infty$}}_{\mathfrak{s},m}
(n,k)z^{m-1} \! = \! \sum_{q=1}^{\mathfrak{s}-1} \sum_{r=1}^{l_{q}=
\varkappa_{nk \tilde{k}_{q}}} 
\dfrac{\widehat{\nu}^{\raise-1.0ex\hbox{$\scriptstyle \infty$}}_{q,r}
(n,k)}{(z \! - \! \alpha_{p_{q}})^{r}} \! + \! \sum_{m=1}^{l_{\mathfrak{s}}
=\varkappa_{nk}} 
\widehat{\nu}^{\raise-1.0ex\hbox{$\scriptstyle \infty$}}_{\mathfrak{s},m}
(n,k)z^{m-1}, \quad 
\widehat{\nu}^{\raise-1.0ex\hbox{$\scriptstyle \infty$}}_{\mathfrak{s},
l_{\mathfrak{s}}}(n,k) \! \neq \! 0.
\end{equation*}
Substituting the latter partial fraction expansion for $\mathcal{X}_{21}(z)$ 
into Equations~\eqref{eq40}--\eqref{eq42}, one arrives at, for $n \! \in 
\! \mathbb{N}$ and $k \! \in \! \lbrace 1,2,\dotsc,K \rbrace$ such 
that $\alpha_{p_{\mathfrak{s}}} \! := \! \alpha_{k} \! = \! \infty$, the 
orthogonality conditions
\begin{gather}
\int_{\mathbb{R}} \left(\sum_{m=1}^{\mathfrak{s}-1} \sum_{q=1}^{l_{m}=
\varkappa_{nk \tilde{k}_{m}}} 
\dfrac{\widehat{\nu}^{\raise-1.0ex\hbox{$\scriptstyle \infty$}}_{m,q}(n,k)}{
(\xi \! - \! \alpha_{p_{m}})^{q}} \! + \! \sum_{q=1}^{l_{\mathfrak{s}}=
\varkappa_{nk}} 
\widehat{\nu}^{\raise-1.0ex\hbox{$\scriptstyle \infty$}}_{\mathfrak{s},q}
(n,k) \xi^{q-1} \right) \xi^{r} \widetilde{w}(\xi) \, \md \xi \! = \! 0, 
\quad r \! = \! 0,1,\dotsc,\varkappa_{nk} \! - \! 2, \label{eq43} \\
\int_{\mathbb{R}} \left(\sum_{m=1}^{\mathfrak{s}-1} \sum_{q=1}^{l_{m}=
\varkappa_{nk \tilde{k}_{m}}} 
\dfrac{\widehat{\nu}^{\raise-1.0ex\hbox{$\scriptstyle \infty$}}_{m,q}
(n,k)}{(\xi \! - \! \alpha_{p_{m}})^{q}} \! + \! \sum_{q=1}^{l_{\mathfrak{s}}
=\varkappa_{nk}} 
\widehat{\nu}^{\raise-1.0ex\hbox{$\scriptstyle \infty$}}_{\mathfrak{s},q}
(n,k) \xi^{q-1} \right) \xi^{\varkappa_{nk}-1} \widetilde{w}(\xi) \, \md 
\xi \! = \! -2 \pi \mi, \label{eq44} \\
\int_{\mathbb{R}} \left(\sum_{m=1}^{\mathfrak{s}-1} \sum_{q=1}^{l_{m}=
\varkappa_{nk \tilde{k}_{m}}} 
\dfrac{\widehat{\nu}^{\raise-1.0ex\hbox{$\scriptstyle \infty$}}_{m,q}
(n,k)}{(\xi \! - \! \alpha_{p_{m}})^{q}} \! + \! \sum_{q=1}^{l_{\mathfrak{s}}
=\varkappa_{nk}} 
\widehat{\nu}^{\raise-1.0ex\hbox{$\scriptstyle \infty$}}_{\mathfrak{s},q}
(n,k) \xi^{q-1} \right) \dfrac{\widetilde{w}(\xi)}{(\xi \! - \! 
\alpha_{p_{i}})^{j}} \, \md \xi \! = \! 0, \quad i \! = \! 1,2,\dotsc,
\mathfrak{s} \! - \! 1, \quad j \! = \! 1,2,\dotsc,l_{i}. \label{eq45}
\end{gather}
Incidentally, for $n \! \in \! \mathbb{N}$ and $k \! \in \! \lbrace 1,2,\dotsc,
K \rbrace$ such that $\alpha_{p_{\mathfrak{s}}} \! := \! \alpha_{k} \! = \! 
\infty$, via the orthogonality conditions~\eqref{eq43} and~\eqref{eq45}, 
the orthogonality condition~\eqref{eq44} can be manipulated thus:
\begin{align*}
-2 \pi \mi =& \, \int_{\mathbb{R}} \left(\sum_{m=1}^{\mathfrak{s}-1} 
\sum_{q=1}^{l_{m}=\varkappa_{nk \tilde{k}_{m}}} 
\dfrac{\widehat{\nu}^{\raise-1.0ex\hbox{$\scriptstyle \infty$}}_{m,q}
(n,k)}{(\xi \! - \! \alpha_{p_{m}})^{q}} \! + \! \sum_{q=1}^{l_{\mathfrak{s}}
=\varkappa_{nk}} 
\widehat{\nu}^{\raise-1.0ex\hbox{$\scriptstyle \infty$}}_{\mathfrak{s},q}
(n,k) \xi^{q-1} \right) 
\widehat{\nu}^{\raise-1.0ex\hbox{$\scriptstyle \infty$}}_{\mathfrak{s},
l_{\mathfrak{s}}}(n,k) \xi^{\varkappa_{nk}-1} \dfrac{\widetilde{w}
(\xi)}{\widehat{\nu}^{\raise-1.0ex\hbox{$\scriptstyle \infty$}}_{\mathfrak{s},
l_{\mathfrak{s}}}(n,k)} \, \md \xi \\
=& \, \int_{\mathbb{R}} \left(\sum_{m=1}^{\mathfrak{s}-1} \sum_{q=1}^{l_{m}
=\varkappa_{nk \tilde{k}_{m}}} 
\dfrac{\widehat{\nu}^{\raise-1.0ex\hbox{$\scriptstyle \infty$}}_{m,q}
(n,k)}{(\xi \! - \! \alpha_{p_{m}})^{q}} \! + \! \sum_{q=1}^{l_{\mathfrak{s}}
=\varkappa_{nk}} 
\widehat{\nu}^{\raise-1.0ex\hbox{$\scriptstyle \infty$}}_{\mathfrak{s},q}
(n,k) \xi^{q-1} \right) \left(\sum_{m=1}^{\mathfrak{s}-1} \sum_{q=1}^{l_{m}
=\varkappa_{nk \tilde{k}_{m}}} \dfrac{
\widehat{\nu}^{\raise-1.0ex\hbox{$\scriptstyle \infty$}}_{m,q}
(n,k)}{(\xi \! - \! \alpha_{p_{m}})^{q}} \right. \\
+&\left. \, \sum_{q=1}^{l_{\mathfrak{s}}-1=\varkappa_{nk}-1} 
\widehat{\nu}^{\raise-1.0ex\hbox{$\scriptstyle \infty$}}_{\mathfrak{s},q}
(n,k) \xi^{q-1} \! + \! 
\widehat{\nu}^{\raise-1.0ex\hbox{$\scriptstyle \infty$}}_{\mathfrak{s},
l_{\mathfrak{s}}}(n,k) \xi^{\varkappa_{nk}-1} \right) \dfrac{\widetilde{w}
(\xi)}{\widehat{\nu}^{\raise-1.0ex\hbox{$\scriptstyle \infty$}}_{\mathfrak{s},
l_{\mathfrak{s}}}(n,k)} \, \md \xi \\
=& \, \int_{\mathbb{R}} \underbrace{\left(\sum_{m=1}^{\mathfrak{s}-1} 
\sum_{q=1}^{l_{m}=\varkappa_{nk \tilde{k}_{m}}} 
\dfrac{\widehat{\nu}^{\raise-1.0ex\hbox{$\scriptstyle \infty$}}_{m,q}
(n,k)}{(\xi \! - \! \alpha_{p_{m}})^{q}} \! + \! \sum_{q=1}^{l_{\mathfrak{s}}
=\varkappa_{nk}} 
\widehat{\nu}^{\raise-1.0ex\hbox{$\scriptstyle \infty$}}_{\mathfrak{s},q}
(n,k) \xi^{q-1} \right)}_{= \, \mathcal{X}_{21}(\xi)} \, \underbrace{
\left(\sum_{m=1}^{\mathfrak{s}-1} \sum_{q=1}^{l_{m}=\varkappa_{nk 
\tilde{k}_{m}}} \dfrac{
\widehat{\nu}^{\raise-1.0ex\hbox{$\scriptstyle \infty$}}_{m,q}(n,k)}{
(\xi \! - \! \alpha_{p_{m}})^{q}} \! + \! \sum_{q=1}^{l_{\mathfrak{s}}=
\varkappa_{nk}} 
\widehat{\nu}^{\raise-1.0ex\hbox{$\scriptstyle \infty$}}_{\mathfrak{s},q}
(n,k) \xi^{q-1} \right)}_{= \, \mathcal{X}_{21}(\xi)} \, \dfrac{\widetilde{w}
(\xi)}{\widehat{\nu}^{\raise-1.0ex\hbox{$\scriptstyle \infty$}}_{\mathfrak{s},
l_{\mathfrak{s}}}(n,k)} \, \md \xi \\
=& \, \int_{\mathbb{R}}(\mathcal{X}_{21}(\xi))^{2} \dfrac{\widetilde{w}
(\xi)}{\widehat{\nu}^{\raise-1.0ex\hbox{$\scriptstyle \infty$}}_{\mathfrak{s},
l_{\mathfrak{s}}}(n,k)} \, \md \xi;
\end{align*}
hence, for $n \! \in \! \mathbb{N}$ and $k \! \in \! \lbrace 1,2,\dotsc,K 
\rbrace$ such that $\alpha_{p_{\mathfrak{s}}} \! := \! \alpha_{k} \! = \! 
\infty$, one arrives at the `normalisation formula'
\begin{equation*}
\int_{\mathbb{R}}(\mathcal{X}_{21}(\xi))^{2} \widetilde{w}(\xi) \, \md \xi \! 
= \! -2 \pi \mi \, \widehat{\nu}^{\raise-1.0ex\hbox{$\scriptstyle \infty$}}_{
\mathfrak{s},\varkappa_{nk}}(n,k).
\end{equation*}
For $n \! \in \! \mathbb{N}$ and $k \! \in \! \lbrace 1,2,\dotsc,K \rbrace$ 
such that $\alpha_{p_{\mathfrak{s}}} \! := \! \alpha_{k} \! = \! \infty$, the 
orthogonality conditions~\eqref{eq43}--\eqref{eq45} give rise to a total 
of (cf. Equation~\eqref{infcount}) $\sum_{r=1}^{\mathfrak{s}}l_{r} \! 
= \! \sum_{r=1}^{\mathfrak{s}-1}l_{r} \! + \! l_{\mathfrak{s}} \! = \! 
\sum_{r=1}^{\mathfrak{s}-1} \varkappa_{nk \tilde{k}_{r}} \! + \! 
\varkappa_{nk} \! = \! (n \! - \! 1)K \! + \! k$ linear inhomogeneous 
algebraic equations for the $(n \! - \! 1)K \! + \! k$ unknowns 
$\widehat{\nu}^{\raise-1.0ex\hbox{$\scriptstyle \infty$}}_{1,1}(n,k),\dotsc,
\widehat{\nu}^{\raise-1.0ex\hbox{$\scriptstyle \infty$}}_{1,l_{1}}(n,k),
\dotsc,\widehat{\nu}^{\raise-1.0ex\hbox{$\scriptstyle \infty$}}_{\mathfrak{s}
-1,1}(n,k),\dotsc,
\widehat{\nu}^{\raise-1.0ex\hbox{$\scriptstyle \infty$}}_{\mathfrak{s}-1,
l_{\mathfrak{s}-1}}(n,k),
\widehat{\nu}^{\raise-1.0ex\hbox{$\scriptstyle \infty$}}_{\mathfrak{s},1}
(n,k),\dotsc,\widehat{\nu}^{\raise-1.0ex\hbox{$\scriptstyle \infty$}}_{
\mathfrak{s},l_{\mathfrak{s}}}(n,k)$, that is,
\begin{align}
\setcounter{MaxMatrixCols}{12}
&\left(
, \label{eq46}
\end{align}
where (with abuse of notation)
\begin{equation*}
n_{1} \! = \! l_{1} \! + \! \dotsb \! + \! l_{\mathfrak{s}-2} \! + \! 1, 
\qquad n_{2} \! = \! (n \! - \! 1)K \! + \! k \! - \! \varkappa_{nk}, \qquad 
m_{1} \! = \! (n \! - \! 1)K \! + \! k \! - \! \varkappa_{nk} \! + \! 1, 
\qquad \text{and} \qquad m_{2} \! = \! (n \! - \! 1)K \! + \! k \! - \! 1.
\end{equation*}
For $n \! \in \! \mathbb{N}$ and $k \! \in \! \lbrace 1,2,\dotsc,K \rbrace$ 
such that $\alpha_{p_{\mathfrak{s}}} \! := \! \alpha_{k} \! = \! \infty$, the 
linear system~\eqref{eq46} of $(n \! - \! 1)K \! + \! k$ inhomogeneous 
algebraic equations for the $(n \! - \! 1)K \! + \! k$ unknowns 
$\widehat{\nu}^{\raise-1.0ex\hbox{$\scriptstyle \infty$}}_{1,1}(n,k),\dotsc,
\widehat{\nu}^{\raise-1.0ex\hbox{$\scriptstyle \infty$}}_{1,l_{1}}(n,k),
\dotsc,\widehat{\nu}^{\raise-1.0ex\hbox{$\scriptstyle \infty$}}_{\mathfrak{s}
-1,1}(n,k),\dotsc,
\widehat{\nu}^{\raise-1.0ex\hbox{$\scriptstyle \infty$}}_{\mathfrak{s}-1,
l_{\mathfrak{s}-1}}(n,k),
\widehat{\nu}^{\raise-1.0ex\hbox{$\scriptstyle \infty$}}_{\mathfrak{s},1}
(n,k),\linebreak[4]
\dotsc,\widehat{\nu}^{\raise-1.0ex\hbox{$\scriptstyle \infty$}}_{\mathfrak{s},
l_{\mathfrak{s}}}(n,k)$ admits a unique solution if, and only if, the 
determinant of the coefficient matrix, denoted $\mathcal{N}_{\infty}^{\sharp}
(n,k)$, is non-zero: this fact will now be established. Via the multi-linearity 
property of the determinant, one shows that, for $n \! \in \! 
\mathbb{N}$ and $k \! \in \! \lbrace 1,2,\dotsc,K \rbrace$ such that 
$\alpha_{p_{\mathfrak{s}}} \! := \! \alpha_{k} \! = \! \infty$,
\begin{align*}
\mathcal{N}_{\infty}^{\sharp}(n,k) =& \, \underbrace{\int_{\mathbb{R}} \int_{
\mathbb{R}} \dotsb \int_{\mathbb{R}}}_{(n-1)K+k} \md \widetilde{\mu}(\xi_{0}) 
\, \md \widetilde{\mu}(\xi_{1}) \, \dotsb \, \md \widetilde{\mu}(\xi_{l_{1}}) \, 
\dotsb \, \md \widetilde{\mu}(\xi_{n_{1}}) \, \dotsb \, \md \widetilde{\mu}
(\xi_{n_{2}}) \, \md \widetilde{\mu}(\xi_{m_{1}}) \, \dotsb \, \md \widetilde{\mu}
(\xi_{m_{2}}) \\
\times& \, \dfrac{1}{\dfrac{1}{\xi_{0}^{0}}(\xi_{1} \! - \! \alpha_{p_{1}})^{
1} \dotsb (\xi_{l_{1}} \! - \! \alpha_{p_{1}})^{l_{1}} \dotsb (\xi_{n_{1}} \! 
- \! \alpha_{p_{\mathfrak{s}-1}})^{1} \dotsb (\xi_{n_{2}} \! - \! \alpha_{p_{
\mathfrak{s}-1}})^{l_{\mathfrak{s}-1}} \dfrac{1}{(\xi_{m_{1}})^{1}} \dotsb 
\dfrac{1}{\xi_{m_{2}}^{\varkappa_{nk}-1}}} \\
\times& \, 
\left\vert

\right\vert.
\end{align*}
Recalling, for $n \! \in \! \mathbb{N}$ and $k \! \in \! \lbrace 1,2,\dotsc,
K \rbrace$ such that $\alpha_{p_{\mathfrak{s}}} \! := \! \alpha_{k} \! = \! 
\infty$, the $(n \! - \! 1)K \! + \! k$ linearly independent functions $\psi_{0}
(z),\psi_{1}(z),\dotsc,\psi_{l_{1}}(z),\dotsc,\psi_{n_{1}}(z),\dotsc,\psi_{n_{2}}
(z),\psi_{m_{1}}(z),\dotsc,\psi_{m_{2}}(z)$ introduced above for the analysis 
of $\hat{\mathfrak{c}}_{D_{\infty}}$, one shows that
\begin{align*}
\mathcal{N}_{\infty}^{\sharp}(n,k) =& \, \underbrace{\int_{\mathbb{R}} 
\int_{\mathbb{R}} \dotsb \int_{\mathbb{R}}}_{(n-1)K+k} \md \widetilde{\mu}
(\xi_{0}) \, \md \widetilde{\mu}(\xi_{1}) \, \dotsb \, \md \widetilde{\mu}
(\xi_{l_{1}}) \, \dotsb \, \md \widetilde{\mu}(\xi_{n_{1}}) \, \dotsb \, \md 
\widetilde{\mu}(\xi_{n_{2}}) \, \md \widetilde{\mu}(\xi_{m_{1}}) \, \dotsb 
\, \md \widetilde{\mu}(\xi_{m_{2}}) \\
\times& \, \dfrac{(-1)^{(n-1)K+k-\varkappa_{nk}}(\psi_{0}(\xi_{0}) \psi_{0}
(\xi_{1}) \dotsb \psi_{0}(\xi_{l_{1}}) \dotsb \psi_{0}(\xi_{n_{1}}) \dotsb 
\psi_{0}(\xi_{n_{2}}) \psi_{0}(\xi_{m_{1}}) \dotsb \psi_{0}(\xi_{m_{2}}))^{-
1}}{\dfrac{1}{\xi_{0}^{0}}(\xi_{1} \! - \! \alpha_{p_{1}})^{1} \dotsb (\xi_{
l_{1}} \! - \! \alpha_{p_{1}})^{l_{1}} \dotsb (\xi_{n_{1}} \! - \! \alpha_{
p_{\mathfrak{s}-1}})^{1} \dotsb (\xi_{n_{2}} \! - \! \alpha_{p_{\mathfrak{s}-
1}})^{l_{\mathfrak{s}-1}} \dfrac{1}{(\xi_{m_{1}})^{1}} \dotsb 
\dfrac{1}{\xi_{m_{2}}^{\varkappa_{nk}-1}}} \\
\times& \, 
\left\vert

\right\vert.
\end{align*}
Modulo the oscillatory factor $(-1)^{(n-1)K+k-\varkappa_{nk}}$, this 
is the same determinantal expression encountered while studying 
$\hat{\mathfrak{c}}_{D_{\infty}}$; therefore, for $n \! \in \! 
\mathbb{N}$ and $k \! \in \! \lbrace 1,2,\dotsc,K \rbrace$ such that 
$\alpha_{p_{\mathfrak{s}}} \! := \! \alpha_{k} \! = \! \infty$, it follows that
\begin{align}
\mathcal{N}_{\infty}^{\sharp}(n,k) =& \, \dfrac{(-1)^{(n-1)K+k-\varkappa_{nk}}
(\mathbb{D}^{\sphat} \,)^{2}}{((n \! - \! 1)K \! + \! k)!} \underbrace{
\int_{\mathbb{R}} \int_{\mathbb{R}} \dotsb \int_{\mathbb{R}}}_{(n-1)K+k} 
\md \widetilde{\mu}(\xi_{0}) \, \md \widetilde{\mu}(\xi_{1}) \, \dotsb \, 
\md \widetilde{\mu}(\xi_{(n-1)K+k-1}) \nonumber \\
\times& \, \prod_{\substack{i,j=0\\j<i}}^{(n-1)K+k-1}(\xi_{i} \! - \! 
\xi_{j})^{2} \left(\prod_{m=0}^{(n-1)K+k-1} \prod_{q=1}^{\mathfrak{s}-1}
(\xi_{m} \! - \! \alpha_{p_{q}})^{\varkappa_{nk \tilde{k}_{q}}} \right)^{-2} 
\quad (\neq 0), \label{eq47}
\end{align}
whence follows, for $n \! \in \! \mathbb{N}$ and $k \! \in \! \lbrace 1,2,
\dotsc,K \rbrace$ such that $\alpha_{p_{\mathfrak{s}}} \! := \! \alpha_{k} 
\! = \! \infty$, the existence and the uniqueness of $\mathcal{X}_{21}(z)$.

\textbf{(2)} If, for $n \! \in \! \mathbb{N}$ and $k \! \in \! \lbrace 1,2,
\dotsc,K \rbrace$ such that $\alpha_{p_{\mathfrak{s}}} \! := \! \alpha_{k} 
\! \neq \! \infty$, $\mathcal{X} \colon \overline{\mathbb{C}} \setminus 
\overline{\mathbb{R}} \! \to \! \operatorname{SL}_{2}(\mathbb{C})$ solves 
the RHP stated in Lemma~$\bm{\mathrm{RHP}_{\mathrm{MPC}}}$, then, 
{}from the jump condition~(ii) of Lemma~$\bm{\mathrm{RHP}_{\mathrm{MPC}}}$, 
it follows that, for the elements of the first column of $\mathcal{X}(z)$,
\begin{equation*}
(\mathcal{X}_{j1}(z))_{+} \! = \! (\mathcal{X}_{j1}(z))_{-} \! := \! 
\mathcal{X}_{j1}(z), \quad j \! = \! 1,2,
\end{equation*}
and, for the elements of the second column of $\mathcal{X}(z)$,
\begin{equation}
(\mathcal{X}_{j2}(z))_{+} \! - \! (\mathcal{X}_{j2}(z))_{-} \! = \! 
\mathcal{X}_{j1}(z) \widetilde{w}(z), \quad j \! = \! 1,2. \label{eq48}
\end{equation}
Via the normalisation and boundedness conditions~(iv) of 
Lemma~$\bm{\mathrm{RHP}_{\mathrm{MPC}}}$, it follows that, for 
$n \! \in \! \mathbb{N}$ and $k \! \in \! \lbrace 1,2,\dotsc,K \rbrace$ 
such that $\alpha_{p_{\mathfrak{s}}} \! := \! \alpha_{k} \! \neq \! \infty$,
\begin{equation}
\begin{split}
\mathcal{X}_{11}(z)(z \! - \! \alpha_{k})^{\varkappa_{nk}-1} 
\underset{\mathbb{C} \setminus \mathbb{R} \, \ni \, z \to 
\alpha_{k}}{=} 1 \! + \! \mathcal{O}(z \! - \! \alpha_{k}), \qquad \, 
\mathcal{X}_{12}(z)(z \! - \! \alpha_{k})^{-(\varkappa_{nk}-1)} 
\underset{\mathbb{C} \setminus \mathbb{R} \, \ni \, z \to 
\alpha_{k}}{=} \mathcal{O}(z \! - \! \alpha_{k}), \\
\mathcal{X}_{21}(z)(z \! - \! \alpha_{k})^{\varkappa_{nk}-1} 
\underset{\mathbb{C} \setminus \mathbb{R} \, \ni \, z \to 
\alpha_{k}}{=} \mathcal{O}(z \! - \! \alpha_{k}), \qquad \, \mathcal{X}_{22}
(z)(z \! - \! \alpha_{k})^{-(\varkappa_{nk}-1)} \underset{\mathbb{C} 
\setminus \mathbb{R} \, \ni \, z \to \alpha_{k}}{=} 1 \! + \! 
\mathcal{O}(z \! - \! \alpha_{k}),
\end{split} \label{eq49}
\end{equation}
\begin{equation}
\begin{split}
\mathcal{X}_{11}(z)z^{-(\varkappa_{nk \tilde{k}_{\mathfrak{s}-1}}^{\infty}+1)} 
\underset{\overline{\mathbb{C}} \setminus \mathbb{R} \, \ni \, z \to 
\alpha_{p_{\mathfrak{s}-1}} = \infty}{=} \mathcal{O}(1), \qquad \, 
\mathcal{X}_{12}(z)z^{\varkappa_{nk \tilde{k}_{\mathfrak{s}-1}}^{\infty}+1} 
\underset{\overline{\mathbb{C}} \setminus \mathbb{R} \, \ni \, z \to 
\alpha_{p_{\mathfrak{s}-1}} = \infty}{=} \mathcal{O}(1), \\
\mathcal{X}_{21}(z)z^{-(\varkappa_{nk \tilde{k}_{\mathfrak{s}-1}}^{\infty}+1)} 
\underset{\overline{\mathbb{C}} \setminus \mathbb{R} \, \ni \, z \to 
\alpha_{p_{\mathfrak{s}-1}} = \infty}{=} \mathcal{O}(1), \qquad \, 
\mathcal{X}_{22}(z)z^{\varkappa_{nk \tilde{k}_{\mathfrak{s}-1}}^{\infty}+1} 
\underset{\overline{\mathbb{C}} \setminus \mathbb{R} \, \ni \, z \to 
\alpha_{p_{\mathfrak{s}-1}} = \infty}{=} \mathcal{O}(1),
\end{split} \label{eq50}
\end{equation}
and, for $q \! = \! 1,2,\dotsc,\mathfrak{s} \! - \! 2$,
\begin{equation}
\begin{split}
\mathcal{X}_{11}(z)(z \! - \! \alpha_{p_{q}})^{\varkappa_{nk \tilde{k}_{q}}} 
\underset{\mathbb{C} \setminus \mathbb{R} \, \ni 
\, z \to \alpha_{p_{q}}}{=} \mathcal{O}(1), \qquad \, \mathcal{X}_{12}(z)
(z \! - \! \alpha_{p_{q}})^{-\varkappa_{nk \tilde{k}_{q}}} \underset{
\mathbb{C} \setminus \mathbb{R} \, \ni \, z \to 
\alpha_{p_{q}}}{=} \mathcal{O}(1), \\
\mathcal{X}_{21}(z)(z \! - \! \alpha_{p_{q}})^{\varkappa_{nk \tilde{k}_{q}}} 
\underset{\mathbb{C} \setminus \mathbb{R} \, \ni 
\, z \to \alpha_{p_{q}}}{=} \mathcal{O}(1), \qquad \, \mathcal{X}_{22}(z)
(z \! - \! \alpha_{p_{q}})^{-\varkappa_{nk \tilde{k}_{q}}} \underset{
\mathbb{C} \setminus \mathbb{R} \, \ni \, z \to 
\alpha_{p_{q}}}{=} \mathcal{O}(1),
\end{split} \label{eq51}
\end{equation}
whence, temporarily re-inserting explicit $n$- and $k$-dependencies, 
that is, $\mathcal{X}_{ij}(n,k;z) \! := \! \mathcal{X}_{ij}(z)$, $i,j \! = \! 1,2$, 
one notes that, for $m \! = \! 1,2$, $\mathcal{X}_{m1} \colon \mathbb{N} 
\times \lbrace 1,2,\dotsc,K \rbrace \times \overline{\mathbb{C}} \setminus 
\lbrace \alpha_{p_{1}},\alpha_{p_{2}},\dotsc,\alpha_{p_{\mathfrak{s}}} \rbrace 
\! \to \! \mathbb{C}$ and $\mathcal{X}_{m2} \colon \mathbb{N} \times 
\lbrace 1,2,\dotsc,K \rbrace \times \overline{\mathbb{C}} \setminus 
\overline{\mathbb{R}} \! \to \! \mathbb{C}$; in particular, $\mathcal{X}_{11}
(z)$ and $\mathcal{X}_{21}(z)$ have no jumps throughout the $z$-plane, 
$\mathcal{X}_{11}(z)$ is a monic (that is, $\operatorname{coeff} \lbrace 
(z \! - \! \alpha_{k})^{-\varkappa_{nk}} \rbrace \! = \! 1)$ meromorphic 
function with poles at $\alpha_{p_{1}},\alpha_{p_{2}},\dotsc,
\alpha_{p_{\mathfrak{s}}}$, and $\mathcal{X}_{21}(z)$ is a meromorphic 
function with poles at $\alpha_{p_{1}},\alpha_{p_{2}},\dotsc,
\alpha_{p_{\mathfrak{s}}}$. Application of the Sokhotski-Plemelj formula 
to the jump condition~\eqref{eq48}, with the Cauchy kernel normalised at 
$\alpha_{p_{\mathfrak{s}}} \! := \! \alpha_{k} \! \neq \! \infty$, $k \! \in 
\! \lbrace 1,2,\dotsc,K \rbrace$, gives rise to the following Cauchy-type 
integral representation:
\begin{equation}
\mathcal{X}_{j2}(z) \! = \! \int_{\mathbb{R}} \dfrac{(z \! - \! \alpha_{k}) 
\mathcal{X}_{j1}(\xi) \widetilde{w}(\xi)}{(\xi \! - \! \alpha_{k})
(\xi \! - \! z)} \, \dfrac{\md \xi}{2 \pi \mi}, \quad z \! \in \! \mathbb{C} 
\setminus \mathbb{R}, \quad j \! = \! 1,2; \label{eq52}
\end{equation}
hence, for $n \! \in \! \mathbb{N}$ and $k \! \in \! \lbrace 1,2,\dotsc,
K \rbrace$ such that $\alpha_{p_{\mathfrak{s}}} \! := \! \alpha_{k} \! 
\neq \! \infty$, $\mathcal{X} \colon \overline{\mathbb{C}} \setminus 
\overline{\mathbb{R}} \! \to \! \operatorname{SL}_{2}(\mathbb{C})$ has 
the integral representation
\begin{equation*}
\mathcal{X}(z) \! = \! 
\begin{pmatrix}
\mathcal{X}_{11}(z) & \int_{\mathbb{R}} \frac{(z-\alpha_{k}) \mathcal{X}_{11}
(\xi) \widetilde{w}(\xi)}{(\xi -\alpha_{k})(\xi -z)} \, \frac{\md \xi}{2 \pi 
\mi} \\
\mathcal{X}_{21}(z) & \int_{\mathbb{R}} \frac{(z-\alpha_{k}) \mathcal{X}_{21}
(\xi) \widetilde{w}(\xi)}{(\xi -\alpha_{k})(\xi -z)} \, \frac{\md \xi}{2 \pi 
\mi}
\end{pmatrix}, \quad z \! \in \! \mathbb{C} \setminus \mathbb{R}.
\end{equation*}
A detailed analysis for the elements of the first row of $\mathcal{X}(z)$ 
is presented first; then, the corresponding analysis for the elements 
of the second row is presented. Recalling that $\int_{\mathbb{R}} \xi^{m} 
\widetilde{w}(\xi) \, \md \xi \! < \! \infty$, $m \! \in \! \mathbb{N}_{0}$, 
and $\int_{\mathbb{R}}(\xi \! - \! \alpha_{k})^{-(m+1)} \widetilde{w}(\xi) 
\, \md \xi \! < \! \infty$, $\alpha_{k} \! \neq \! \infty$, $k \! \in \! \lbrace 
1,2,\dotsc,K \rbrace$, it follows, via the expansion $\tfrac{1}{z_{1}-z_{2}} 
\! = \! \sum_{i=0}^{l} \tfrac{z_{2}^{i}}{z_{1}^{i+1}} \! + \! 
\tfrac{z_{2}^{l+1}}{z_{1}^{l+1}(z_{1}-z_{2})}$, $l \! \in \! 
\mathbb{N}_{0}$, the integral representation (cf. Equation~\eqref{eq52}) 
$\mathcal{X}_{12}(z) \! = \! \int_{\mathbb{R}} \tfrac{(z-\alpha_{k}) 
\mathcal{X}_{11}(\xi) \widetilde{w}(\xi)}{(\xi -\alpha_{k})(\xi -z)} \, 
\tfrac{\md \xi}{2 \pi \mi}$, $z \! \in \! \mathbb{C} \setminus \mathbb{R}$, 
and the first line of each of the asymptotic conditions~\eqref{eq49}, 
\eqref{eq50}, and~\eqref{eq51}, that, for $n \! \in \! \mathbb{N}$ and $k 
\! \in \! \lbrace 1,2,\dotsc,K \rbrace$ such that $\alpha_{p_{\mathfrak{s}}} 
\! := \! \alpha_{k} \! \neq \! \infty$,
\begin{gather}
\int_{\mathbb{R}} \left(\dfrac{\mathcal{X}_{11}(\xi)}{\xi \! - \! \alpha_{k}} 
\right) \dfrac{\widetilde{w}(\xi)}{(\xi \! - \! \alpha_{k})^{p}} \, \md \xi 
\! = \! 0, \quad p \! = \! 0,1,\dotsc,\varkappa_{nk} \! - \! 1, \label{eq53} 
\\
\int_{\mathbb{R}} \left(\dfrac{\mathcal{X}_{11}(\xi)}{\xi \! - \! \alpha_{k}} 
\right) \dfrac{\widetilde{w}(\xi)}{(\xi \! - \! \alpha_{k})^{\varkappa_{nk}}} 
\, \md \xi \! \neq \! 0, \label{eq54} \\
\int_{\mathbb{R}} \left(\dfrac{\mathcal{X}_{11}(\xi)}{\xi \! - \! \alpha_{k}} 
\right) \dfrac{\widetilde{w}(\xi)}{(\xi \! - \! \alpha_{p_{q}})^{r}} \, \md 
\xi \! = \! 0, \quad q \! = \! 1,2,\dotsc,\mathfrak{s} \! - \! 2, \quad r \! 
= \! 1,2,\dotsc,\varkappa_{nk \tilde{k}_{q}}, \label{eq55} \\
\int_{\mathbb{R}} \left(\dfrac{\mathcal{X}_{11}(\xi)}{\xi \! - \! \alpha_{k}} 
\right) \xi^{r} \widetilde{w}(\xi) \, \md \xi \! = \! 0, \quad r \! = \! 1,2,
\dotsc,\varkappa_{nk \tilde{k}_{\mathfrak{s}-1}}^{\infty}. \label{eq56}
\end{gather}
(Note: if, for $n \! \in \! \mathbb{N}$ and $k \! \in \! \lbrace 1,2,\dotsc,K 
\rbrace$ such that $\alpha_{p_{\mathfrak{s}}} \! := \! \alpha_{k} \! \neq 
\! \infty$, the set $\lbrace \mathstrut \alpha_{k^{\prime}}, \, k^{\prime} 
\! \in \! \lbrace 1,2,\dotsc,K \rbrace; \, \alpha_{k^{\prime}} \! \neq \! 
\alpha_{k}, \, \alpha_{k} \! \neq \! \infty \rbrace \! = \! \varnothing$, 
then Equations~\eqref{eq55} and~\eqref{eq56} are vacuous; this can 
only occur if $n \! = \! 1$.) Recalling {}from the analysis preceding the 
integral representation~\eqref{eq52} that, for $n \! \in \! \mathbb{N}$ 
and $k \! \in \! \lbrace 1,2,\dotsc,K \rbrace$ such that 
$\alpha_{p_{\mathfrak{s}}} \! := \! \alpha_{k} \! \neq \! \infty$, 
$\mathcal{X}_{11} \colon \mathbb{N} \times \lbrace 1,2,\dotsc,K \rbrace 
\times \overline{\mathbb{C}} \setminus \lbrace \alpha_{p_{1}},\alpha_{p_{2}},
\dotsc,\alpha_{p_{\mathfrak{s}}} \rbrace \! \to \! \mathbb{C}$ is a monic 
$(\operatorname{coeff} \lbrace (z \! - \! \alpha_{k})^{-\varkappa_{nk}} 
\rbrace \! = \! 1)$ meromorphic function with pole set $\lbrace \alpha_{p_{1}},
\alpha_{p_{2}},\dotsc,\alpha_{p_{\mathfrak{s}}} \rbrace$ and with no jumps 
throughout the $z$-plane, and that, for $n \! \in \! \mathbb{N}$ and $k \! 
\in \! \lbrace 1,2,\dotsc,K \rbrace$ such that $\alpha_{p_{\mathfrak{s}}} 
\! := \! \alpha_{k} \! \neq \! \infty$, the monic MPC ORF satisfies the 
orthogonality conditions~\eqref{eq16}--\eqref{eq19}, it follows {}from 
the latter two observations and Equations~\eqref{eq53}--\eqref{eq56} 
that, for $n \! \in \! \mathbb{N}$ and $k \! \in \! \lbrace 1,2,\dotsc,
K \rbrace$ such that $\alpha_{p_{\mathfrak{s}}} \! := \! \alpha_{k} 
\! \neq \! \infty$,
\begin{equation*}
\pmb{\pi}^{n}_{k}(z) \! = \! \dfrac{\mathcal{X}_{11}(z)}{z \! - \! \alpha_{k}}.
\end{equation*}
Via Equations~\eqref{eq53}, \eqref{eq55}, and~\eqref{eq56}, and this 
latter formula, one writes, for $n \! \in \! \mathbb{N}$ and $k \! \in \! 
\lbrace 1,2,\dotsc,K \rbrace$ such that $\alpha_{p_{\mathfrak{s}}} \! 
:= \! \alpha_{k} \! \neq \! \infty$, Equation~\eqref{eq54} in a more 
transparent form:
\begin{align*}
&\int_{\mathbb{R}} \left(\dfrac{\mathcal{X}_{11}(\xi)}{\xi \! - \! \alpha_{k}} 
\right) \dfrac{\widetilde{w}(\xi)}{(\xi \! - \! \alpha_{k})^{\varkappa_{nk}}} 
\, \md \xi \! = \! \int_{\mathbb{R}} \left(\dfrac{\mathcal{X}_{11}(\xi)}{\xi 
\! - \! \alpha_{k}} \right) \dfrac{\mu^{f}_{n,\varkappa_{nk}}(n,k)}{(\xi 
\! - \! \alpha_{k})^{\varkappa_{nk}}} \dfrac{\widetilde{w}(\xi)}{\mu^{f}_{n,
\varkappa_{nk}}(n,k)} \, \md \xi \! = \! \int_{\mathbb{R}} \underbrace{
\left(\dfrac{\mathcal{X}_{11}(\xi)}{\xi \! - \! \alpha_{k}} \right)}_{= \, 
\pmb{\pi}^{n}_{k}(\xi)} \\
&\times \, \underbrace{\left(\phi_{0}^{f}(n,k) \! + \! 
\sum_{q=1}^{\mathfrak{s}-2} \sum_{r=1}^{\varkappa_{nk \tilde{k}_{q}}} 
\dfrac{\tilde{\nu}^{f}_{r,q}(n,k)}{(\xi \! - \! \alpha_{p_{q}})^{r}} 
\! + \! \sum_{l=1}^{\varkappa^{\infty}_{nk \tilde{k}_{\mathfrak{s}-1}}} 
\hat{\nu}^{f}_{n,l}(n,k) \xi^{l} \! + \! \sum_{m=1}^{\varkappa_{nk}-1} 
\dfrac{\mu^{f}_{n,m}(n,k)}{(\xi \! - \! \alpha_{k})^{m}} \! + \! 
\dfrac{\mu^{f}_{n,\varkappa_{nk}}(n,k)}{(\xi \! - \! \alpha_{k})^{
\varkappa_{nk}}} \right)}_{= \, \phi^{n}_{k}(\xi)} \, 
\dfrac{\widetilde{w}(\xi)}{\mu^{f}_{n,\varkappa_{nk}}(n,k)} \, \md \xi \\
&= \, \int_{\mathbb{R}} \underbrace{\pmb{\pi}^{n}_{k}(\xi)}_{= \, 
(\mu^{f}_{n,\varkappa_{nk}}(n,k))^{-1} \phi^{n}_{k}(\xi)} \phi^{n}_{k}(\xi) 
\dfrac{\widetilde{w}(\xi)}{\mu^{f}_{n,\varkappa_{nk}}(n,k)} \, \md \xi 
\! = \! (\mu^{f}_{n,\varkappa_{nk}}(n,k))^{-2} \int_{\mathbb{R}}
(\phi^{n}_{k}(\xi))^{2} \, \widetilde{w}(\xi) \, \md \xi \! = \! 
(\mu^{f}_{n,\varkappa_{nk}}(n,k))^{-2};
\end{align*}
hence, for $n \! \in \! \mathbb{N}$ and $k \! \in \! \lbrace 1,2,\dotsc,K 
\rbrace$ such that $\alpha_{p_{\mathfrak{s}}} \! := \! \alpha_{k} \! \neq 
\! \infty$, via this latter relation, Equations~\eqref{eq53}--\eqref{eq56} 
are written in the following, more convenient, form:
\begin{gather}
\int_{\mathbb{R}} \pmb{\pi}^{n}_{k}(\xi)(\xi \! - \! \alpha_{k})^{-p} \, 
\widetilde{w}(\xi) \, \md \xi \! = \! 0, \quad p \! = \! 0,1,\dotsc,
\varkappa_{nk} \! - \! 1, \label{eq57} \\
\int_{\mathbb{R}} \pmb{\pi}^{n}_{k}(\xi)(\xi \! - \! \alpha_{k})^{-
\varkappa_{nk}} \, \widetilde{w}(\xi) \, \md \xi \! = \! 
(\mu^{f}_{n,\varkappa_{nk}}(n,k))^{-2}, \label{eq58} \\
\int_{\mathbb{R}} \pmb{\pi}^{n}_{k}(\xi)(\xi \! - \! \alpha_{p_{q}})^{-r} 
\, \widetilde{w}(\xi) \, \md \xi \! = \! 0, \quad q \! = \! 1,2,\dotsc,
\mathfrak{s} \! - \! 2, \quad r \! = \! 1,2,\dotsc,\varkappa_{nk 
\tilde{k}_{q}}, \label{eq59} \\
\int_{\mathbb{R}} \pmb{\pi}^{n}_{k}(\xi) \xi^{r} \, \widetilde{w}
(\xi) \, \md \xi \! = \! 0, \quad r \! = \! 1,2,\dotsc,
\varkappa_{nk \tilde{k}_{\mathfrak{s}-1}}^{\infty}. \label{eq60}
\end{gather}
(Note: if, for $n \! \in \! \mathbb{N}$ and $k \! \in \! \lbrace 1,2,\dotsc,
K \rbrace$ such that $\alpha_{p_{\mathfrak{s}}} \! := \! \alpha_{k} \! 
\neq \! \infty$, the set $\lbrace \mathstrut \alpha_{k^{\prime}}, \, 
k^{\prime} \! \in \! \lbrace 1,2,\dotsc,K \rbrace; \, \alpha_{k^{\prime}} 
\! \neq \! \alpha_{k}, \, \alpha_{k} \! \neq \! \infty \rbrace \! = \! 
\varnothing$, then Equations~\eqref{eq59} and~\eqref{eq60} are 
vacuous; this can only occur if $n \! = \! 1$.) For $n \! \in \! 
\mathbb{N}$ and $k \! \in \! \lbrace 1,2,\dotsc,K \rbrace$ such 
that $\alpha_{p_{\mathfrak{s}}} \! := \! \alpha_{k} \! \neq \! 
\infty$, Equation~\eqref{eq57} gives rise to $\varkappa_{nk}$ 
conditions, Equation~\eqref{eq58} gives rise to $1$ condition, 
Equations~\eqref{eq59} and~\eqref{eq60} give rise to 
(cf. Equation~\eqref{fincount}) $\sum_{q=1}^{\mathfrak{s}-2} 
\varkappa_{nk \tilde{k}_{q}} \! + \! \varkappa_{nk \tilde{k}_{
\mathfrak{s}-1}}^{\infty} \! = \! (n \! - \! 1)K \! + \! k \! - \! \varkappa_{nk}$ 
conditions, for a total of $(n \! - \! 1)K \! + \! k \! + \! 1$ conditions, 
which is precisely the number necessary in order to determine, uniquely 
(see below), the associated $(n$- and $k$-dependent) norming constant, 
$\mu^{f}_{n,\varkappa_{nk}}(n,k)$.

One now examines, for $n \! \in \! \mathbb{N}$ and $k \! \in \! \lbrace 
1,2,\dotsc,K \rbrace$ such that $\alpha_{p_{\mathfrak{s}}} \! := \! 
\alpha_{k} \! \neq \! \infty$, Equations~\eqref{eq57}--\eqref{eq60} in 
detail. Proceeding as per the discussion of Subsection~\ref{subsubsec1.2.2}, 
write, for $n \! \in \! \mathbb{N}$ and $k \! \in \! \lbrace 1,2,\dotsc,K 
\rbrace$ such that $\alpha_{p_{\mathfrak{s}}} \! := \! \alpha_{k} \! \neq 
\! \infty$, the ordered disjoint partition for the repeated pole sequence:
\begin{align*}
&\overset{1}{\lbrace \underbrace{\alpha_{1},\alpha_{2},\dotsc,\alpha_{K}}_{K} 
\rbrace} \cup \dotsb \cup \overset{n-1}{\lbrace \underbrace{\alpha_{1},
\alpha_{2},\dotsc,\alpha_{K}}_{K} \rbrace} \cup \overset{n}{\lbrace 
\underbrace{\alpha_{1},\alpha_{2},\dotsc,\alpha_{k}}_{k} \rbrace} \\
&:= \, \bigcup_{q=1}^{\mathfrak{s}-2} \lbrace \underbrace{\alpha_{
i(q)_{k_{q}}},\alpha_{i(q)_{k_{q}}},\dotsc,\alpha_{i(q)_{k_{q}}}}_{l_{q}=
\varkappa_{nk \tilde{k}_{q}}} \rbrace \cup \lbrace \underbrace{\alpha_{i
(\mathfrak{s}-1)_{k_{\mathfrak{s}-1}}},\alpha_{i(\mathfrak{s}-1)_{
k_{\mathfrak{s}-1}}},\dotsc,\alpha_{i(\mathfrak{s}-1)_{k_{\mathfrak{s}
-1}}}}_{l_{\mathfrak{s}-1}=\varkappa^{\infty}_{nk \tilde{k}_{\mathfrak{s}-1}}} 
\rbrace \cup \lbrace \underbrace{\alpha_{i(\mathfrak{s})_{k_{\mathfrak{s}}}},
\alpha_{i(\mathfrak{s})_{k_{\mathfrak{s}}}},\dotsc,\alpha_{i(\mathfrak{s})_{
k_{\mathfrak{s}}}}}_{l_{\mathfrak{s}}=\varkappa_{ni(\mathfrak{s})_{
k_{\mathfrak{s}}}}} \rbrace \\
&:= \, \bigcup_{q=1}^{\mathfrak{s}-2} \lbrace \underbrace{\alpha_{p_{q}},
\alpha_{p_{q}},\dotsc,\alpha_{p_{q}}}_{l_{q}=\varkappa_{nk \tilde{k}_{q}}} 
\rbrace \cup \lbrace \underbrace{\alpha_{p_{\mathfrak{s}-1}},\alpha_{p_{
\mathfrak{s}-1}},\dotsc,\alpha_{p_{\mathfrak{s}-1}}}_{l_{\mathfrak{s}-1}=
\varkappa^{\infty}_{nk \tilde{k}_{\mathfrak{s}-1}}} \rbrace \cup \lbrace 
\underbrace{\alpha_{k},\alpha_{k},\dotsc,\alpha_{k}}_{l_{\mathfrak{s}}=
\varkappa_{nk}} \rbrace \\
&= \, \bigcup_{q=1}^{\mathfrak{s}-2} \lbrace \underbrace{\alpha_{p_{q}},
\alpha_{p_{q}},\dotsc,\alpha_{p_{q}}}_{l_{q}=\varkappa_{nk \tilde{k}_{q}}} 
\rbrace \cup \lbrace \underbrace{\infty,\infty,\dotsc,\infty}_{l_{\mathfrak{s}
-1}=\varkappa^{\infty}_{nk \tilde{k}_{\mathfrak{s}-1}}} \rbrace \cup \lbrace 
\underbrace{\alpha_{k},\alpha_{k},\dotsc,\alpha_{k}}_{l_{\mathfrak{s}}=
\varkappa_{nk}} \rbrace,
\end{align*}
where $\sum_{q=1}^{\mathfrak{s}}l_{q} \! = \! \sum_{q=1}^{
\mathfrak{s}-2}l_{q} \! + \! l_{\mathfrak{s}-1} \! + \! l_{\mathfrak{s}} 
\! = \! \sum_{q=1}^{\mathfrak{s}-2} \varkappa_{nk \tilde{k}_{q}} \! + \! 
\varkappa^{\infty}_{nk \tilde{k}_{\mathfrak{s}-1}} \! + \! \varkappa_{nk} 
\! = \! (n \! - \! 1)K \! + \! k$. Hence, via this notational 
preamble, and the analysis leading to the orthogonality 
Equations~\eqref{eq57}--\eqref{eq60}, one writes, in the indicated 
order, for $n \! \in \! \mathbb{N}$ and $k \! \in \! \lbrace 1,2,\dotsc,
K \rbrace$ such that $\alpha_{p_{\mathfrak{s}}} \! := \! \alpha_{k} \! 
\neq \! \infty$,
\begin{align*}
\pmb{\pi}^{n}_{k}(z)=& \, \dfrac{\phi_{0}^{f}(n,k)}{\mu^{f}_{n,\varkappa_{nk}}
(n,k)} \! + \! \dfrac{1}{\mu^{f}_{n,\varkappa_{nk}}(n,k)} \sum_{q=1}^{
\mathfrak{s}-2} \sum_{r=1}^{\varkappa_{nk \tilde{k}_{q}}} \dfrac{\tilde{
\nu}^{f}_{r,q}(n,k)}{(z \! - \! \alpha_{i(q)_{k_{q}}})^{r}} \! + \! \dfrac{
1}{\mu^{f}_{n,\varkappa_{nk}}(n,k)} \sum_{l=1}^{\varkappa^{\infty}_{nk 
\tilde{k}_{\mathfrak{s}-1}}} \hat{\nu}^{f}_{n,l}(n,k)z^{l} \\
+& \, \dfrac{1}{\mu^{f}_{n,\varkappa_{nk}}(n,k)} \sum_{m=1}^{\varkappa_{ni
(\mathfrak{s})_{k_{\mathfrak{s}}}}-1} \dfrac{\mu^{f}_{n,m}(n,k)}{(z \! - \! 
\alpha_{i(\mathfrak{s})_{k_{\mathfrak{s}}}})^{m}} \! + \! \dfrac{1}{(z \! - 
\! \alpha_{k})^{\varkappa_{nk}}} \\
=& \, \dfrac{\phi^{f}_{0}(n,k)}{\mu^{f}_{n,\varkappa_{nk}}(n,k)} \! + \! 
\dfrac{1}{\mu^{f}_{n,\varkappa_{nk}}(n,k)} \sum_{q=1}^{\mathfrak{s}-2} \sum_{r
=1}^{\varkappa_{nk \tilde{k}_{q}}} \dfrac{\tilde{\nu}^{f}_{r,q}(n,k)}{(z \! 
- \! \alpha_{p_{q}})^{r}} \! + \! \dfrac{1}{\mu^{f}_{n,\varkappa_{nk}}(n,k)} 
\sum_{l=1}^{\varkappa^{\infty}_{nk \tilde{k}_{\mathfrak{s}-1}}} \hat{\nu}^{
f}_{n,l}(n,k)z^{l} \\
+& \, \dfrac{1}{\mu^{f}_{n,\varkappa_{nk}}(n,k)} \sum_{m=1}^{\varkappa_{nk}-1} 
\dfrac{\mu^{f}_{n,m}(n,k)}{(z \! - \! \alpha_{k})^{m}} \! + \! \dfrac{1}{(z \! 
- \! \alpha_{k})^{\varkappa_{nk}}} \\
:=& \, \widetilde{\phi}_{0}^{\raise-1.0ex\hbox{$\scriptstyle f$}}(n,k) \! + \! 
\sum_{m=1}^{\mathfrak{s}-2} \, \sum_{q=1}^{l_{m}=\varkappa_{nk \tilde{k}_{m}}} 
\dfrac{\widetilde{\nu}_{m,q}^{\raise-1.0ex\hbox{$\scriptstyle f$}}(n,k)}{(z 
\! - \! \alpha_{p_{m}})^{q}} \! + \! \sum_{q=1}^{l_{\mathfrak{s}-1}=
\varkappa^{\infty}_{nk \tilde{k}_{\mathfrak{s}-1}}} \widetilde{\nu}_{
\mathfrak{s}-1,q}^{\raise-1.0ex\hbox{$\scriptstyle f$}}(n,k)z^{q} \! + \! 
\sum_{r=1}^{l_{\mathfrak{s}}=\varkappa_{nk}} \dfrac{\widetilde{\nu}_{
\mathfrak{s},r}^{\raise-1.0ex\hbox{$\scriptstyle f$}}(n,k)}{(z 
\! - \! \alpha_{k})^{r}}, \quad \widetilde{\nu}_{\mathfrak{s},
l_{\mathfrak{s}}}^{\raise-1.0ex\hbox{$\scriptstyle f$}}(n,k) \! \equiv \! 1.
\end{align*}
Substituting the latter partial fraction expansion of $\pmb{\pi}^{n}_{k}(z)$ 
into the orthogonality conditions~\eqref{eq57}--\eqref{eq60}, one arrives 
at, for $n \! \in \! \mathbb{N}$ and $k \! \in \! \lbrace 1,2,\dotsc,K 
\rbrace$ such that $\alpha_{p_{\mathfrak{s}}} \! := \! \alpha_{k} \! \neq 
\! \infty$, the orthogonality conditions
\begin{align*}
\int_{\mathbb{R}} &\left(
\widetilde{\phi}_{0}^{\raise-1.0ex\hbox{$\scriptstyle f$}}(n,k) \! + \! 
\sum_{m=1}^{\mathfrak{s}-2} \sum_{q=1}^{l_{m}=\varkappa_{nk \tilde{k}_{m}}} 
\dfrac{\widetilde{\nu}_{m,q}^{\raise-1.0ex\hbox{$\scriptstyle f$}}(n,k)}{
(\xi \! - \! \alpha_{p_{m}})^{q}} \! + \! \sum_{q=1}^{l_{\mathfrak{s}-1}=
\varkappa^{\infty}_{nk \tilde{k}_{\mathfrak{s}-1}}} 
\widetilde{\nu}_{\mathfrak{s}-1,q}^{\raise-1.0ex\hbox{$\scriptstyle f$}}
(n,k) \xi^{q} \! + \! \sum_{q=1}^{l_{\mathfrak{s}}=\varkappa_{nk}} \dfrac{
\widetilde{\nu}_{\mathfrak{s},q}^{\raise-1.0ex\hbox{$\scriptstyle f$}}(n,k)}{
(\xi \! - \! \alpha_{k})^{q}} \right) \dfrac{\widetilde{w}(\xi)}{(\xi \! - \! 
\alpha_{k})^{r_{1}}} \, \md \xi \! = \! 0, \\
&r_{1} \! = \! 0,1,\dotsc,\varkappa_{nk} \! - \! 1, \\
\int_{\mathbb{R}} &\left(
\widetilde{\phi}_{0}^{\raise-1.0ex\hbox{$\scriptstyle f$}}(n,k) \! + \! 
\sum_{m=1}^{\mathfrak{s}-2} \sum_{q=1}^{l_{m}=\varkappa_{nk \tilde{k}_{m}}} 
\dfrac{\widetilde{\nu}_{m,q}^{\raise-1.0ex\hbox{$\scriptstyle f$}}(n,k)}{
(\xi \! - \! \alpha_{p_{m}})^{q}} \! + \! \sum_{q=1}^{l_{\mathfrak{s}-1}=
\varkappa^{\infty}_{nk \tilde{k}_{\mathfrak{s}-1}}} \widetilde{\nu}_{
\mathfrak{s}-1,q}^{\raise-1.0ex\hbox{$\scriptstyle f$}}(n,k) \xi^{q} \! + \! 
\sum_{q=1}^{l_{\mathfrak{s}}=\varkappa_{nk}} \dfrac{\widetilde{\nu}_{
\mathfrak{s},q}^{\raise-1.0ex\hbox{$\scriptstyle f$}}(n,k)}{(\xi \! - \! 
\alpha_{k})^{q}} \right) \dfrac{\widetilde{w}(\xi)}{(\xi \! - \! 
\alpha_{k})^{\varkappa_{nk}}} \, \md \xi \! = \! (\mu^{f}_{n,\varkappa_{nk}}
(n,k))^{-2}, \\
\int_{\mathbb{R}} &\left(
\widetilde{\phi}_{0}^{\raise-1.0ex\hbox{$\scriptstyle f$}}(n,k) \! + \! 
\sum_{m=1}^{\mathfrak{s}-2} \sum_{q=1}^{l_{m}=\varkappa_{nk \tilde{k}_{m}}} 
\dfrac{\widetilde{\nu}_{m,q}^{\raise-1.0ex\hbox{$\scriptstyle f$}}(n,k)}{(
\xi \! - \! \alpha_{p_{m}})^{q}} \! + \! \sum_{q=1}^{l_{\mathfrak{s}-1}=
\varkappa^{\infty}_{nk \tilde{k}_{\mathfrak{s}-1}}} \widetilde{\nu}_{
\mathfrak{s}-1,q}^{\raise-1.0ex\hbox{$\scriptstyle f$}}(n,k) \xi^{q} \! + \! 
\sum_{q=1}^{l_{\mathfrak{s}}=\varkappa_{nk}} \dfrac{\widetilde{\nu}_{
\mathfrak{s},q}^{\raise-1.0ex\hbox{$\scriptstyle f$}}(n,k)}{(\xi \! - \! 
\alpha_{k})^{q}} \right) \dfrac{\widetilde{w}(\xi)}{(\xi \! - \! 
\alpha_{p_{i}})^{j}} \, \md \xi \! = \! 0, \\
&i \! = \! 1,2,\dotsc,\mathfrak{s} \! - \! 2, \quad j \! = \! 1,2,\dotsc,
l_{i}, \\
\int_{\mathbb{R}} &\left(
\widetilde{\phi}_{0}^{\raise-1.0ex\hbox{$\scriptstyle f$}}(n,k) \! + \! 
\sum_{m=1}^{\mathfrak{s}-2} \sum_{q=1}^{l_{m}=\varkappa_{nk \tilde{k}_{m}}} 
\dfrac{\widetilde{\nu}_{m,q}^{\raise-1.0ex\hbox{$\scriptstyle f$}}(n,k)}{(
\xi \! - \! \alpha_{p_{m}})^{q}} \! + \! \sum_{q=1}^{l_{\mathfrak{s}-1}=
\varkappa^{\infty}_{nk \tilde{k}_{\mathfrak{s}-1}}} \widetilde{\nu}_{
\mathfrak{s}-1,q}^{\raise-1.0ex\hbox{$\scriptstyle f$}}(n,k) \xi^{q} \! + \! 
\sum_{q=1}^{l_{\mathfrak{s}}=\varkappa_{nk}} \dfrac{\widetilde{\nu}_{
\mathfrak{s},q}^{\raise-1.0ex\hbox{$\scriptstyle f$}}(n,k)}{(\xi \! - \! 
\alpha_{k})^{q}} \right) \xi^{r_{2}} \, \widetilde{w}(\xi) \, \md \xi \! = 
\! 0, \\
&r_{2} \! = \! 1,2,\dotsc,l_{\mathfrak{s}-1}.
\end{align*}
For $n \! \in \! \mathbb{N}$ and $k \! \in \! \lbrace 1,2,\dotsc,K \rbrace$ 
such that $\alpha_{p_{\mathfrak{s}}} \! := \! \alpha_{k} \! \neq \! \infty$, 
the above orthogonality conditions give rise to a total of $\sum_{r=1}^{
\mathfrak{s}}l_{r} \! + \! 1 \! = \! \sum_{r=1}^{\mathfrak{s}-2}l_{r} 
\! + \! l_{\mathfrak{s}-1} \! + \! l_{\mathfrak{s}} \! + \! 1 \! = \! 
\sum_{r=1}^{\mathfrak{s}-2} \varkappa_{nk \tilde{k}_{r}} \! + \! 
\varkappa^{\infty}_{nk \tilde{k}_{\mathfrak{s}-1}} \! + \! \varkappa_{nk} 
\! + \! 1 \! = \! (n \! - \! 1)K \! + \! k \! + \! 1$ linear inhomogeneous 
algebraic equations for the $(n \! - \! 1)K \! + \! k \! + \! 1$ (real) 
unknowns $\widetilde{\phi}_{0}^{\raise-1.0ex\hbox{$\scriptstyle f$}}(n,k),
\widetilde{\nu}_{1,1}^{\raise-1.0ex\hbox{$\scriptstyle f$}}(n,k),\dotsc,
\widetilde{\nu}_{1,l_{1}}^{\raise-1.0ex\hbox{$\scriptstyle f$}}(n,k),\dotsc,
\widetilde{\nu}_{\mathfrak{s}-1,1}^{\raise-1.0ex\hbox{$\scriptstyle f$}}
(n,k),\dotsc,\widetilde{\nu}_{\mathfrak{s}-1,
l_{\mathfrak{s}-1}}^{\raise-1.0ex\hbox{$\scriptstyle f$}}(n,k),
\widetilde{\nu}_{\mathfrak{s},1}^{\raise-1.0ex\hbox{$\scriptstyle f$}}
(n,k),\dotsc,\linebreak[4]
\widetilde{\nu}_{\mathfrak{s},
l_{\mathfrak{s}}-1}^{\raise-1.0ex\hbox{$\scriptstyle f$}}(n,k),
(\mu^{f}_{n,\varkappa_{nk}}(n,k))^{-2}$, that is,
\begin{align}
\setcounter{MaxMatrixCols}{12}
&\left(

\right), \label{eq61}
\end{align}
where (with abuse of notation)
\begin{equation*}
n_{1} \! = \! l_{1} \! + \! \dotsb \! + \! l_{\mathfrak{s}-2} \! + \! 1, 
\qquad n_{2} \! = \! (n \! - \! 1)K \! + \! k \! - \! \varkappa_{nk}, \qquad 
m_{1} \! = \! (n \! - \! 1)K \! + \! k \! - \! 1, \qquad \text{and} \qquad 
m_{2} \! = \! (n \! - \! 1)K \! + \! k.
\end{equation*}
For $n \! \in \! \mathbb{N}$ and $k \! \in \! \lbrace 1,2,\dotsc,K \rbrace$ 
such that $\alpha_{p_{\mathfrak{s}}} \! := \! \alpha_{k} \! \neq \! \infty$, 
the linear system~\eqref{eq61} of $(n \! - \! 1)K \! + \! k \! + \! 1$ 
inhomogeneous algebraic equations for the $(n \! - \! 1)K \! + \! k \! 
+ \! 1$ (real) unknowns 
$\widetilde{\phi}_{0}^{\raise-1.0ex\hbox{$\scriptstyle f$}}(n,k),
\widetilde{\nu}_{1,1}^{\raise-1.0ex\hbox{$\scriptstyle f$}}(n,k),\dotsc,
\widetilde{\nu}_{1,l_{1}}^{\raise-1.0ex\hbox{$\scriptstyle f$}}(n,k),\dotsc,
\widetilde{\nu}_{\mathfrak{s}-1,1}^{\raise-1.0ex\hbox{$\scriptstyle f$}}
(n,k),\dotsc,\linebreak[4]
\widetilde{\nu}_{\mathfrak{s}-1,
l_{\mathfrak{s}-1}}^{\raise-1.0ex\hbox{$\scriptstyle f$}}(n,k),
\widetilde{\nu}_{\mathfrak{s},1}^{\raise-1.0ex\hbox{$\scriptstyle f$}}
(n,k),\dotsc,\widetilde{\nu}_{\mathfrak{s},
l_{\mathfrak{s}}-1}^{\raise-1.0ex\hbox{$\scriptstyle f$}}(n,k),
(\mu^{f}_{n,\varkappa_{nk}}(n,k))^{-2}$ admits a unique solution if, 
and only if, the determinant of the coefficient matrix is non-zero; this 
fact will now be established, and, en route, an explicit multi-integral 
representation for the associated $(n$- and $k$-dependent) norming 
constant, $(\mu^{f}_{n,\varkappa_{nk}}(n,k))^{-2}$, will be derived. For 
$n \! \in \! \mathbb{N}$ and $k \! \in \! \lbrace 1,2,\dotsc,K \rbrace$ 
such that $\alpha_{p_{\mathfrak{s}}} \! := \! \alpha_{k} \! \neq \! 
\infty$, one uses Cramer's Rule and the multi-linearity property of 
the determinant to show that
\begin{equation} \label{nhormmconstatfin} 
\left(\mu^{f}_{n,\varkappa_{nk}}(n,k) \right)^{-2} \! = \! 
\dfrac{\hat{\mathfrak{c}}_{N_{f}}}{\hat{\mathfrak{c}}_{D_{f}}},
\end{equation}
where
\begin{align*}
\hat{\mathfrak{c}}_{N_{f}} :=& \, \underbrace{\int_{\mathbb{R}} 
\int_{\mathbb{R}} \dotsb \int_{\mathbb{R}}}_{(n-1)K+k+1} \md \widetilde{\mu}
(\xi_{0}) \, \md \widetilde{\mu}(\xi_{1}) \, \dotsb \, \md \widetilde{\mu}(\xi_{l_{1}}) 
\, \dotsb \, \md \widetilde{\mu}(\xi_{n_{1}}) \, \dotsb \, \md \widetilde{\mu}
(\xi_{n_{2}}) \, \md \widetilde{\mu}(\xi_{m_{2}-\varkappa_{nk}+1}) \, \dotsb \, 
\md \widetilde{\mu}(\xi_{m_{2}}) \\
\times& \, \dfrac{1}{(\xi_{0} \! - \! \alpha_{k})^{0}(\xi_{1} \! - \! 
\alpha_{p_{1}})^{1} \dotsb (\xi_{l_{1}} \! - \! \alpha_{p_{1}})^{l_{1}} 
\dotsb \dfrac{1}{(\xi_{n_{1}})^{1}} \dotsb \dfrac{1}{(\xi_{n_{2}})^{
l_{\mathfrak{s}-1}}}(\xi_{m_{2}-\varkappa_{nk}+1}-\alpha_{k})^{1} \dotsb 
(\xi_{m_{2}} \! - \! \alpha_{k})^{\varkappa_{nk}}} \\
\times& \, 
\left\lvert 

\right\rvert,
\end{align*}
and
\begin{align*}
\hat{\mathfrak{c}}_{D_{f}} :=& \, \underbrace{\int_{\mathbb{R}} 
\int_{\mathbb{R}} \dotsb \int_{\mathbb{R}}}_{(n-1)K+k} \md \widetilde{\mu}
(\xi_{0}) \, \md \widetilde{\mu}(\xi_{1}) \, \dotsb \, \md \widetilde{\mu}
(\xi_{l_{1}}) \, \dotsb \, \md \widetilde{\mu}(\xi_{n_{1}}) \, \dotsb \, \md 
\widetilde{\mu}(\xi_{n_{2}}) \, \md \widetilde{\mu}(\xi_{m_{2}-\varkappa_{nk}+1}) 
\, \dotsb \, \md \widetilde{\mu}(\xi_{m_{1}}) \\
\times& \, \dfrac{1}{(\xi_{0} \! - \! \alpha_{k})^{0}(\xi_{1} \! - \! 
\alpha_{p_{1}})^{1} \dotsb (\xi_{l_{1}} \! - \! \alpha_{p_{1}})^{l_{1}} 
\dotsb \dfrac{1}{(\xi_{n_{1}})^{1}} \dotsb \dfrac{1}{(\xi_{n_{2}})^{
l_{\mathfrak{s}-1}}}(\xi_{m_{2}-\varkappa_{nk}+1} \! - \! \alpha_{k})^{1} 
\dotsb (\xi_{m_{1}} \! - \! \alpha_{k})^{\varkappa_{nk}-1}} \\
\times& \, 
\left\lvert 

\right\rvert.
\end{align*}
For $n \! \in \! \mathbb{N}$ and $k \! \in \! \lbrace 1,2,\dotsc,K \rbrace$ 
such that $\alpha_{p_{\mathfrak{s}}} \! := \! \alpha_{k} \! \neq \! \infty$, 
$\hat{\mathfrak{c}}_{N_{f}}$ is studied first, and then $\hat{\mathfrak{c}}_{D_{f}}$. 
For $n \! \in \! \mathbb{N}$ and $k \! \in \! \lbrace 1,2,\dotsc,K \rbrace$ 
such that $\alpha_{p_{\mathfrak{s}}} \! := \! \alpha_{k} \! \neq \! \infty$, 
introduce the following notation (recall that $l_{q} \! = \! \varkappa_{nk 
\tilde{k}_{q}}$, $q \! = \! 1,2,\dotsc,\mathfrak{s} \! - \! 2$, 
$l_{\mathfrak{s}-1} \! = \! \varkappa^{\infty}_{nk \tilde{k}_{\mathfrak{s}-1}}$, 
and $l_{\mathfrak{s}} \! = \! \varkappa_{nk})$:
\begin{equation*}
\phi_{0}(z) \! := \! \prod_{m=1}^{\mathfrak{s}-2}(z \! - \! 
\alpha_{p_{m}})^{l_{m}}(z \! - \! \alpha_{k})^{\varkappa_{nk}} 
\! =: \sum_{j=0}^{(n-1)K+k} \mathfrak{a}_{j,0}z^{j},
\end{equation*}
and (with some abuse of notation), for $r \! = \! 1,\dotsc,\mathfrak{s} 
\! - \! 1,\mathfrak{s}$, $q(r) \! = \! \sum_{i=1}^{r-1}l_{i} \! + \! 1,
\sum_{i=1}^{r-1}l_{i} \! + \! 2,\dotsc,\sum_{i=1}^{r-1}l_{i} \! + \! l_{r}$, 
and $m(r) \! = \! 1,2,\dotsc,l_{r}$,
\begin{equation*}
\phi_{q(r)}(z) \! := \! \left\lbrace \phi_{0}(z)(z \! - \! \alpha_{p_{r}})^{-
m(r)(1-\delta_{r \mathfrak{s}-1})}z^{m(r) \delta_{r \mathfrak{s}-1}} 
\! =: \sum_{j=0}^{(n-1)K+k} \mathfrak{a}_{j,q(r)}z^{j} \right\rbrace;
\end{equation*}
e.g., for $r \! = \! 1$, the notation $\phi_{q(1)}(z) \! := \! \lbrace \phi_{0}
(z)(z \! - \! \alpha_{p_{1}})^{-m(1)} \! =: \! \sum_{j=0}^{(n-1)K+k} 
\mathfrak{a}_{j,q(1)}z^{j} \rbrace$, $q(1) \! = \! 1,2,\dotsc,l_{1}$, 
$m(1) \! = \! 1,2,\dotsc,l_{1}$, denotes the (set of) $l_{1} \! = \! 
\varkappa_{nk \tilde{k}_{1}}$ functions
\begin{gather*}
\phi_{1}(z) \! = \! \dfrac{\phi_{0}(z)}{z \! - \! \alpha_{p_{1}}} \! =: 
\sum_{j=0}^{(n-1)K+k} \mathfrak{a}_{j,1}z^{j}, \, \phi_{2}(z) \! 
= \! \dfrac{\phi_{0}(z)}{(z \! - \! \alpha_{p_{1}})^{2}} \! =: 
\sum_{j=0}^{(n-1)K+k} \mathfrak{a}_{j,2}z^{j}, \, \dotsc, \, 
\phi_{l_{1}}(z) \! = \! \dfrac{\phi_{0}(z)}{(z \! - \! \alpha_{p_{1}})^{l_{1}}} 
\! =: \sum_{j=0}^{(n-1)K+k} \mathfrak{a}_{j,l_{1}}z^{j},
\end{gather*}
for $r \! = \! \mathfrak{s} \! - \! 1$ (recall that 
$\alpha_{p_{\mathfrak{s}-1}} \! = \! \infty)$, the notation 
$\phi_{q(\mathfrak{s}-1)}(z) \! := \! \lbrace \phi_{0}(z)z^{m(\mathfrak{s}-1)} 
\! =: \! \sum_{j=0}^{(n-1)K+k} \mathfrak{a}_{j,q(\mathfrak{s}-1)}z^{j} 
\rbrace$, $q(\mathfrak{s} \! - \! 1) \! = \! l_{1} \! + \! \dotsb \! + \! 
l_{\mathfrak{s}-2} \! + \! 1,l_{1} \! + \! \dotsb \! + \! l_{\mathfrak{s}-2} 
\! + \! 2,\dotsc,l_{1} \! + \! \dotsb \! + \! l_{\mathfrak{s}-2} \! + \! 
l_{\mathfrak{s}-1} \! = \! (n \! - \! 1)K \! + \! k \! - \! \varkappa_{nk}$, 
$m(\mathfrak{s} \! - \! 1) \! = \! 1,2,\dotsc,l_{\mathfrak{s}-1}$, denotes 
the (set of) $l_{\mathfrak{s}-1} \! = \! \varkappa^{\infty}_{nk 
\tilde{k}_{\mathfrak{s}-1}}$ functions
\begin{gather*}
\phi_{l_{1}+ \dotsb +l_{\mathfrak{s}-2}+1}(z) \! = \! \phi_{0}(z)z \! =: 
\sum_{j=0}^{(n-1)K+k} \mathfrak{a}_{j,l_{1}+\dotsb +l_{\mathfrak{s}-2}+1}
z^{j}, \, \phi_{l_{1}+\dotsb +l_{\mathfrak{s}-2}+2}(z) \! = \! \phi_{0}(z)
z^{2} \! =: \sum_{j=0}^{(n-1)K+k} \mathfrak{a}_{j,l_{1}+\dotsb +
l_{\mathfrak{s}-2}+2}z^{j}, \, \dotsc \\
\dotsc, \, \phi_{(n-1)K+k-\varkappa_{nk}}(z) \! = \! \phi_{0}(z)
z^{l_{\mathfrak{s}-1}} \! =: \sum_{j=0}^{(n-1)K+k} 
\mathfrak{a}_{j,(n-1)K+k-\varkappa_{nk}}z^{j},
\end{gather*}
etc., and, for $r \! = \! \mathfrak{s}$, the notation $\phi_{q(\mathfrak{s})}
(z) \! := \! \lbrace \phi_{0}(z)(z \! - \! \alpha_{k})^{-m(\mathfrak{s})} \! =: 
\! \sum_{j=0}^{(n-1)K+k} \mathfrak{a}_{j,q(\mathfrak{s})}z^{j} \rbrace$, 
$q(\mathfrak{s}) \! = \! (n \! - \! 1)K \! + \! k \! - \! \varkappa_{nk} \! 
+ \! 1,(n \! - \! 1)K \! + \! k \! - \! \varkappa_{nk} \! + \! 2,\dotsc,(n \! 
- \! 1)K \! + \! k$, $m(\mathfrak{s}) \! = \! 1,2,\dotsc,l_{\mathfrak{s}}$, 
denotes the (set of) $l_{\mathfrak{s}} \! = \! \varkappa_{nk}$ functions
\begin{gather*}
\phi_{(n-1)K+k-\varkappa_{nk}+1}(z) \! = \! \dfrac{\phi_{0}(z)}{z \! 
- \! \alpha_{k}} \! =: \sum_{j=0}^{(n-1)K+k} \mathfrak{a}_{j,(n-1)K+k-
\varkappa_{nk}+1}z^{j}, \, \phi_{(n-1)K+k-\varkappa_{nk}+2}(z) \! = \! 
\dfrac{\phi_{0}(z)}{(z \! - \! \alpha_{k})^{2}} \! =: \sum_{j=0}^{(n-1)K+k} 
\mathfrak{a}_{j,(n-1)K+k-\varkappa_{nk}+2}z^{j}, \, \dotsc \\
\dotsc, \, \phi_{(n-1)K+k}(z) \! = \! \dfrac{\phi_{0}(z)}{(z \! - \! 
\alpha_{k})^{\varkappa_{nk}}} \! =: \sum_{j=0}^{(n-1)K+k} 
\mathfrak{a}_{j,(n-1)K+k}z^{j}.
\end{gather*}
(Note: $\# \lbrace \phi_{0}(z)(z \! - \! \alpha_{p_{r}})^{-m(r)(1-\delta_{r 
\mathfrak{s}-1})}z^{m(r) \delta_{r \mathfrak{s}-1}} \rbrace \! = \! l_{r}$, 
$r \! = \! 1,2,\dotsc,\mathfrak{s}$, and $\# \cup_{r=1}^{\mathfrak{s}} 
\lbrace \phi_{0}(z)(z \! - \! \alpha_{p_{r}})^{-m(r)(1-\delta_{r \mathfrak{s}-1})}
z^{m(r) \delta_{r \mathfrak{s}-1}} \rbrace \! = \! \sum_{r=1}^{\mathfrak{s}}
l_{r} \! = \! (n \! - \! 1)K \! + \! k$.) One notes that, for $n \! \in \! 
\mathbb{N}$ and $k \! \in \! \lbrace 1,2,\dotsc,K \rbrace$ such that 
$\alpha_{p_{\mathfrak{s}}} \! := \! \alpha_{k} \! \neq \! \infty$, the 
$l_{1} \! + \! \dotsb \! + \! l_{\mathfrak{s}-1} \! + \! l_{\mathfrak{s}} \! 
+ \! 1 \! = \! (n \! - \! 1)K \! + \! k \! + \! 1$ functions $\phi_{0}(z),\phi_{1}
(z),\dotsc,\phi_{l_{1}}(z),\dotsc,\phi_{l_{1}+\dotsb +l_{\mathfrak{s}-2}+1}(z),
\dotsc,\phi_{(n-1)K+k-\varkappa_{nk}}(z),\phi_{(n-1)K+k-\varkappa_{nk}+1}
(z),\dotsc,\phi_{(n-1)K+k}(z)$ are linearly independent on $\mathbb{R}$, that 
is, for $z \! \in \! \mathbb{R}$, $\sum_{j=0}^{(n-1)K+k} \mathfrak{c}_{j} 
\phi_{j}(z) \! = \! 0$ $\Rightarrow$ (via a Vandermonde-type argument; see 
the $((n \! - \! 1)K \! + \! k \! + \! 1) \times ((n \! - \! 1)K \! + \! k \! + \! 1)$ 
non-zero determinant $\mathbb{D}^{\lozenge}$ in Equation~\eqref{eq62} 
below) $\mathfrak{c}_{j} \! = \! 0$, $j \! = \! 0,1,\dotsc,(n \! - \! 1)K \! + \! k$. 
For $n \! \in \! \mathbb{N}$ and $k \! \in \! \lbrace 1,2,\dotsc,K \rbrace$ 
such that $\alpha_{p_{\mathfrak{s}}} \! := \! \alpha_{k} \! \neq \! \infty$, let 
$\mathfrak{S}_{(n-1)K+k+1}$ denote the $((n \! - \! 1)K \! + \! k \! + \! 1)!$ 
permutations of $\lbrace 0,1,\dotsc,(n \! - \! 1)K \! + \! k \rbrace$. Using the 
above notation and the multi-linearity property of the determinant, one studies, 
thus, for $n \! \in \! \mathbb{N}$ and $k \! \in \! \lbrace 1,2,\dotsc,K \rbrace$ 
such that $\alpha_{p_{\mathfrak{s}}} \! := \! \alpha_{k} \! \neq \! \infty$, 
$\hat{\mathfrak{c}}_{N_{f}}$:
\begin{align*}
\hat{\mathfrak{c}}_{N_{f}} =& \, \underbrace{\int_{\mathbb{R}} 
\int_{\mathbb{R}} \dotsb \int_{\mathbb{R}}}_{(n-1)K+k+1} \md \widetilde{\mu}
(\xi_{0}) \, \md \widetilde{\mu}(\xi_{1}) \, \dotsb \, \md \widetilde{\mu}(\xi_{l_{1}}) 
\, \dotsb \, \md \widetilde{\mu}(\xi_{n_{1}}) \, \dotsb \, \md \widetilde{\mu}
(\xi_{n_{2}}) \, \md \widetilde{\mu}(\xi_{m_{2}-\varkappa_{nk}+1}) \, \dotsb 
\, \md \widetilde{\mu}(\xi_{m_{2}}) \\
\times& \, \dfrac{1}{(\xi_{0} \! - \! \alpha_{k})^{0}(\xi_{1} \! - \! 
\alpha_{p_{1}})^{1} \dotsb (\xi_{l_{1}} \! - \! \alpha_{p_{1}})^{l_{1}} 
\dotsb \dfrac{1}{(\xi_{n_{1}})^{1}} \dotsb \dfrac{1}{(\xi_{n_{2}})^{
l_{\mathfrak{s}-1}}}(\xi_{m_{2}-\varkappa_{nk}+1} \! - \! \alpha_{k})^{1} 
\dotsb (\xi_{m_{2}} \! - \! \alpha_{k})^{\varkappa_{nk}}} \\
\times& \, \dfrac{1}{\phi_{0}(\xi_{0}) \phi_{0}(\xi_{1}) \dotsb \phi_{0}
(\xi_{l_{1}}) \dotsb \phi_{0}(\xi_{n_{1}}) \dotsb \phi_{0}(\xi_{n_{2}}) 
\phi_{0}(\xi_{m_{2}-\varkappa_{nk}+1}) \dotsb \phi_{0}(\xi_{m_{2}})} \\
\times& \, 
\underbrace{
\left\vert 

\right\vert}_{=: \, \mathbb{G}^{\lozenge}(\xi_{0},\xi_{1},\dotsc,\xi_{l_{1}},
\dotsc,\xi_{n_{1}},\dotsc,\xi_{n_{2}},\xi_{m_{2}-\varkappa_{nk}+1},\dotsc,
\xi_{(n-1)K+k})} \\
=& \, \dfrac{1}{(m_{2} \! + \! 1)!} \sum_{\pmb{\sigma} \in \mathfrak{S}_{m_{2}
+1}} \underbrace{\int_{\mathbb{R}} \int_{\mathbb{R}} \dotsb \int_{\mathbb{R}}
}_{(n-1)K+k+1} \md \widetilde{\mu}(\xi_{\sigma (0)}) \, \md \widetilde{\mu}
(\xi_{\sigma (1)}) \, \dotsb \, \md \widetilde{\mu}(\xi_{\sigma (l_{1})}) \, \dotsb 
\, \md \widetilde{\mu}(\xi_{\sigma (n_{1})}) \, \dotsb \, \md \widetilde{\mu}
(\xi_{\sigma (n_{2})}) \\
\times& \, \tfrac{\dotsb \, \md \widetilde{\mu}(\xi_{\sigma (m_{2}-\varkappa_{nk}+1)}) 
\, \dotsb \, \md \widetilde{\mu}(\xi_{\sigma (m_{2})})}{(\xi_{\sigma (0)}-\alpha_{k})^{0}
(\xi_{\sigma (1)}-\alpha_{p_{1}})^{1} \dotsb (\xi_{\sigma (l_{1})}-\alpha_{p_{1}})^{l_{1}} 
\dotsb \frac{1}{(\xi_{\sigma (n_{1})})^{1}} \dotsb \frac{1}{(\xi_{\sigma (n_{2})})^{
l_{\mathfrak{s}-1}}}(\xi_{\sigma (m_{2}-\varkappa_{nk}+1)}-\alpha_{k})^{1} \dotsb 
(\xi_{\sigma (m_{2})}-\alpha_{k})^{\varkappa_{nk}}} \\
\times& \, \dfrac{1}{\phi_{0}(\xi_{\sigma (0)}) \phi_{0}(\xi_{\sigma (1)}) 
\dotsb \phi_{0}(\xi_{\sigma (l_{1})}) \dotsb \phi_{0}(\xi_{\sigma (n_{1})}) 
\dotsb \phi_{0}(\xi_{\sigma (n_{2})}) \phi_{0}(\xi_{\sigma (m_{2}-\varkappa_{
nk}+1)}) \dotsb \phi_{0}(\xi_{\sigma (m_{2})})} \\
\times& \, \mathbb{G}^{\lozenge}(\xi_{\sigma (0)},\xi_{\sigma (1)},\dotsc,
\xi_{\sigma (l_{1})},\dotsc,\xi_{\sigma (n_{1})},\dotsc,\xi_{\sigma (n_{2})},
\xi_{\sigma (m_{2}-\varkappa_{nk}+1)},\dotsc,\xi_{\sigma ((n-1)K+k)}) \\
=& \, \dfrac{1}{(m_{2} \! + \! 1)!} \underbrace{\int_{\mathbb{R}} \int_{
\mathbb{R}} \dotsb \int_{\mathbb{R}}}_{(n-1)K+k+1} \md \widetilde{\mu}(\xi_{0}) 
\, \md \widetilde{\mu}(\xi_{1}) \, \dotsb \, \md \widetilde{\mu}(\xi_{l_{1}}) \, 
\dotsb \, \md \widetilde{\mu}(\xi_{n_{1}}) \, \dotsb \, \md \widetilde{\mu}
(\xi_{n_{2}}) \, \md \widetilde{\mu}(\xi_{m_{2}-\varkappa_{nk}+1}) \, \dotsb \\
\times& \, \dotsb \, \md \widetilde{\mu}(\xi_{m_{2}}) \, \mathbb{G}^{\lozenge}
(\xi_{0},\xi_{1},\dotsc,\xi_{l_{1}},\dotsc,\xi_{n_{1}},\dotsc,\xi_{n_{2}},
\xi_{m_{2}-\varkappa_{nk}+1},\dotsc,\xi_{(n-1)K+k}) \\
\times& \, \sum_{\pmb{\sigma} \in \mathfrak{S}_{m_{2}+1}} \operatorname{sgn}
(\pmb{\pmb{\sigma}}) \tfrac{1}{(\xi_{\sigma (0)}-\alpha_{k})^{0}
(\xi_{\sigma (1)}-\alpha_{p_{1}})^{1} \dotsb (\xi_{\sigma (l_{1})}-
\alpha_{p_{1}})^{l_{1}} \dotsb \frac{1}{(\xi_{\sigma (n_{1})})^{1}} \dotsb 
\frac{1}{(\xi_{\sigma (n_{2})})^{l_{\mathfrak{s}-1}}}(\xi_{\sigma (m_{2}-
\varkappa_{nk}+1)}-\alpha_{k})^{1} \dotsb (\xi_{\sigma (m_{2})}-
\alpha_{k})^{\varkappa_{nk}}} \\
\times& \, \dfrac{1}{\phi_{0}(\xi_{\sigma (0)}) \phi_{0}(\xi_{\sigma (1)}) 
\dotsb \phi_{0}(\xi_{\sigma (l_{1})}) \dotsb \phi_{0}(\xi_{\sigma (n_{1})}) 
\dotsb \phi_{0}(\xi_{\sigma (n_{2})}) \phi_{0}(\xi_{\sigma (m_{2}-\varkappa_{
nk}+1)}) \dotsb \phi_{0}(\xi_{\sigma (m_{2})})} \\
=& \, \dfrac{1}{(m_{2} \! + \! 1)!} \underbrace{\int_{\mathbb{R}} \int_{
\mathbb{R}} \dotsb \int_{\mathbb{R}}}_{(n-1)K+k+1} \md \widetilde{\mu}(\xi_{0}) \, 
\md \widetilde{\mu}(\xi_{1}) \, \dotsb \, \md \widetilde{\mu}(\xi_{l_{1}}) \, \dotsb \, 
\md \widetilde{\mu}(\xi_{n_{1}}) \, \dotsb \, \md \widetilde{\mu}(\xi_{n_{2}}) \, 
\md \widetilde{\mu}(\xi_{m_{2}-\varkappa_{nk}+1}) \, \dotsb \\
\times& \, \dotsb \, \md \widetilde{\mu}(\xi_{m_{2}}) \, \mathbb{G}^{\lozenge}
(\xi_{0},\xi_{1},\dotsc,\xi_{l_{1}},\dotsc,\xi_{n_{1}},\dotsc,\xi_{n_{2}},
\xi_{m_{2}-\varkappa_{nk}+1},\dotsc,\xi_{(n-1)K+k}) \\
\times& \, 
\left\lvert 
 
\right\rvert \\
=& \, \dfrac{1}{(m_{2} \! + \! 1)!} \underbrace{\int_{\mathbb{R}} \int_{
\mathbb{R}} \dotsb \int_{\mathbb{R}}}_{(n-1)K+k+1} \md \widetilde{\mu}(\xi_{0}) 
\, \md \widetilde{\mu}(\xi_{1}) \, \dotsb \, \md \widetilde{\mu}(\xi_{l_{1}}) \, 
\dotsb \, \md \widetilde{\mu}(\xi_{n_{1}}) \, \dotsb \, \md \widetilde{\mu}
(\xi_{n_{2}}) \, \md \widetilde{\mu}(\xi_{m_{2}-\varkappa_{nk}+1}) \, \dotsb \\
\times& \, \dotsb \, \md \widetilde{\mu}(\xi_{m_{2}}) \, \dfrac{\mathbb{G}^{\lozenge}
(\xi_{0},\xi_{1},\dotsc,\xi_{l_{1}},\dotsc,\xi_{n_{1}},\dotsc,\xi_{n_{2}},
\xi_{m_{2}-\varkappa_{nk}+1},\dotsc,\xi_{(n-1)K+k})}{(\phi_{0}(\xi_{0}) 
\phi_{0}(\xi_{1}) \dotsb \phi_{0}(\xi_{l_{1}}) \dotsb \phi_{0}(\xi_{n_{1}}) 
\dotsb \phi_{0}(\xi_{n_{2}}) \phi_{0}(\xi_{m_{2}-\varkappa_{nk}+1}) \dotsb 
\phi_{0}(\xi_{m_{2}}))^{2}} \\
\times& \, 
\underbrace{
\left\vert

\right\vert}_{= \, \mathbb{G}^{\lozenge}(\xi_{0},\xi_{1},\dotsc,\xi_{l_{1}},
\dotsc,\xi_{n_{1}},\dotsc,\xi_{n_{2}},\xi_{m_{2}-\varkappa_{nk}+1},\dotsc,
\xi_{(n-1)K+k})} \\
=& \, \dfrac{1}{(m_{2} \! + \! 1)!} \underbrace{\int_{\mathbb{R}} \int_{
\mathbb{R}} \dotsb \int_{\mathbb{R}}}_{(n-1)K+k+1} \md \widetilde{\mu}(\xi_{0}) 
\, \md \widetilde{\mu}(\xi_{1}) \, \dotsb \, \md \widetilde{\mu}(\xi_{l_{1}}) \, 
\dotsb \, \md \widetilde{\mu}(\xi_{n_{1}}) \, \dotsb \, \md \widetilde{\mu}
(\xi_{n_{2}}) \, \md \widetilde{\mu}(\xi_{m_{2}-\varkappa_{nk}+1}) \, \dotsb \\
\times& \, \dotsb \, \md \widetilde{\mu}(\xi_{m_{2}}) \left(\dfrac{\mathbb{G}^{\lozenge}
(\xi_{0},\xi_{1},\dotsc,\xi_{l_{1}},\dotsc,\xi_{n_{1}},\dotsc,\xi_{n_{2}},
\xi_{m_{2}-\varkappa_{nk}+1},\dotsc,\xi_{(n-1)K+k})}{\phi_{0}(\xi_{0}) \phi_{
0}(\xi_{1}) \dotsb \phi_{0}(\xi_{l_{1}}) \dotsb \phi_{0}(\xi_{n_{1}}) \dotsb 
\phi_{0}(\xi_{n_{2}}) \phi_{0}(\xi_{m_{2}-\varkappa_{nk}+1}) \dotsb \phi_{0}
(\xi_{m_{2}})} \right)^{2};
\end{align*}
but, noting the determinantal factorisation
\begin{equation*}
\mathbb{G}^{\lozenge}(\xi_{0},\xi_{1},\dotsc,\xi_{l_{1}},\dotsc,\xi_{n_{1}},
\dotsc,\xi_{n_{2}},\xi_{m_{2}-\varkappa_{nk}+1},\dotsc,\xi_{(n-1)K+k}) 
\! := \! \overset{f}{\mathcal{V}}_{1}(\xi_{0},\xi_{1},\dotsc,\xi_{(n-1)K+k}) 
\mathbb{D}^{\lozenge},
\end{equation*}
where
\begin{equation*}
\overset{f}{\mathcal{V}}_{1}(\xi_{0},\xi_{1},\dotsc,\xi_{(n-1)K+k}) 
\! = \! \left\lvert 

\right), \label{eq62}
\end{equation}
with
\begin{equation*}
(\mathfrak{A}^{\lozenge}(0))_{1j} \! = \! 
\begin{cases}
\mathfrak{a}_{j-1}^{\lozenge}(\vec{\bm{\alpha}}), &\text{$j \! = \! 1,2,
\dotsc,(n \! - \! 1)K \! + \! k \! - \! l_{\mathfrak{s}-1} \! + \! 1$,} \\
0, &\text{$j \! = \! (n \! - \! 1)K \! + \! k \! - \! l_{\mathfrak{s}-1} 
\! + \! 2,\dotsc,(n \! - \! 1)K \! + \! k \! + \! 1$,}
\end{cases}
\end{equation*}
for $r \! = \! 1,\dotsc,\mathfrak{s} \! - \! 2,\mathfrak{s}$, with $l_{r} \! 
= \! \varkappa_{nk \tilde{k}_{r}}$, $r \! = \! 1,2,\dotsc,\mathfrak{s}-2$, 
$l_{\mathfrak{s}-1} \! = \! \varkappa^{\infty}_{nk \tilde{k}_{\mathfrak{s}-1}}$, 
$l_{\mathfrak{s}} \! = \! \varkappa_{nk}$, and $\sum_{m=1}^{\mathfrak{s}-2}
l_{m} \! + \! l_{\mathfrak{s}-1} \! = \! (n \! - \! 1)K \! + \! k \! - \! \varkappa_{nk}$,
\begin{align*}
i(r)=& \, 2 \! + \! \sum_{m=1}^{r-1}l_{m},3 \! + \! \sum_{m=1}^{r-1}l_{m},
\dotsc,1 \! + \! l_{r} \! + \! \sum_{m=1}^{r-1}l_{m}, \qquad \qquad q(i(r),r) 
\! = \! 1,2,\dotsc,l_{r}, \\
(\mathfrak{A}^{\lozenge}(r))_{i(r)j(r)}=& \, 
\begin{cases}
\dfrac{(-1)^{q(i(r),r)}}{\prod_{m=0}^{q(i(r),r)-1}(l_{r} \! - \! m)} 
\left(\dfrac{\partial}{\partial \alpha_{p_{r}}} \right)^{q(i(r),r)} 
\mathfrak{a}_{j(r)-1}^{\lozenge}(\vec{\bm{\alpha}}), &\text{$j(r) \! = \! 
1,2,\dotsc,(n \! - \! 1)K \! + \! k \! - \! l_{\mathfrak{s}-1} \! - \! 
q(i(r),r) \! + \! 1$,} \\
0, &\text{$j(r) \! = \! (n \! - \! 1)K \! + \! k \! - \! l_{\mathfrak{s}-1} 
\! - \! q(i(r),r) \! + \! 2,\dotsc,(n \! - \! 1)K \! + \! k \! + \! 1$,}
\end{cases}
\end{align*}
and, for $r \! = \! \mathfrak{s} \! - \! 1$,
\begin{align*}
i(\mathfrak{s} \! - \! 1)=& \, 2 \! + \! \sum_{m=1}^{\mathfrak{s}-2}l_{m},3 
\! + \! \sum_{m=1}^{\mathfrak{s}-2}l_{m},\dotsc,1 \! + \! l_{\mathfrak{s}-1} 
\! + \! \sum_{m=1}^{\mathfrak{s}-2}l_{m}, \quad \quad q(i(\mathfrak{s} \! 
- \! 1),\mathfrak{s} \! - \! 1) \! = \! 1,2,\dotsc,l_{\mathfrak{s}-1}, \\
(\mathfrak{A}^{\lozenge}(\mathfrak{s} \! - \! 1))_{i(\mathfrak{s}-1)
j(\mathfrak{s}-1)}=& \, 
\begin{cases}
0, &\text{$j(\mathfrak{s} \! - \! 1) \! = \! 1,\dotsc,q(i(\mathfrak{s} 
\! - \! 1),\mathfrak{s} \! - \! 1)$,} \\
\mathfrak{a}_{j(\mathfrak{s}-1)-q(i(\mathfrak{s}-1),\mathfrak{s}-1)-1}^{
\lozenge}(\vec{\bm{\alpha}}), &\text{$j(\mathfrak{s} \! - \! 1) \! = \! 
q(i(\mathfrak{s} \! - \! 1),\mathfrak{s} \! - \! 1) \! + \! 1,\dotsc,
i(\mathfrak{s} \! - \! 1) \! + \! \varkappa_{nk}$,} \\
0, &\text{$j(\mathfrak{s} \! - \! 1) \! = \! i(\mathfrak{s} \! - \! 1) \! 
+ \! \varkappa_{nk} \! + \! 1,\dotsc,(n \! - \! 1)K \! + \! k \! + \! 1$,}
\end{cases}
\end{align*}
where, for $\tilde{m}_{1} \! = \! 0,1,\dotsc,(n \! - \! 1)K \! + \! k \! - \! 
l_{\mathfrak{s}-1}$,
\begin{equation*}
\mathfrak{a}_{\tilde{m}_{1}}^{\lozenge}(\vec{\bm{\alpha}}):= \, \, 
\mathlarger{\sum_{\underset{\underset{\underset{\sum_{\underset{m \neq 
\mathfrak{s}-1}{m=1}}^{\mathfrak{s}}i_{m}=(n-1)K+k-l_{\mathfrak{s}-1}-
\tilde{m}_{1}}{i_{\mathfrak{s}}=0,1,\dotsc,l_{\mathfrak{s}}}}{p \in 
\lbrace 1,2,\dotsc,\mathfrak{s}-2 \rbrace}}{i_{p}=0,1,\dotsc,
l_{p}}}}(-1)^{(n-1)K+k-l_{\mathfrak{s}-1}-\tilde{m}_{1}} 
\prod_{j=1}^{\mathfrak{s}-2} \binom{l_{j}}{i_{j}} 
\binom{l_{\mathfrak{s}}}{i_{\mathfrak{s}}} \prod_{m=1}^{\mathfrak{s}-2}
(\alpha_{p_{m}})^{i_{m}} \alpha_{k}^{i_{\mathfrak{s}}},
\end{equation*}
it follows that, for $n \! \in \! \mathbb{N}$ and $k \! \in \! \lbrace 
1,2,\dotsc,K \rbrace$ such that $\alpha_{p_{\mathfrak{s}}} \! := \! 
\alpha_{k} \! \neq \! \infty$,
\begin{align}
\hat{\mathfrak{c}}_{N_{f}} =& \, \dfrac{1}{(m_{2} \! + \! 1)!} \underbrace{
\int_{\mathbb{R}} \int_{\mathbb{R}} \dotsb \int_{\mathbb{R}}}_{(n-1)K+k+1} 
\md \widetilde{\mu}(\xi_{0}) \, \md \widetilde{\mu}(\xi_{1}) \, \dotsb \, \md 
\widetilde{\mu}(\xi_{l_{1}}) \, \dotsb \, \md \widetilde{\mu}(\xi_{n_{1}}) \, 
\dotsb \, \md \widetilde{\mu}(\xi_{n_{2}}) \, \md \widetilde{\mu}(\xi_{m_{2}-
\varkappa_{nk}+1}) \, \dotsb \nonumber \\
\times& \, \dotsb \, \md \widetilde{\mu}(\xi_{m_{2}}) \dfrac{(\mathbb{D}^{
\lozenge})^{2} \prod_{\underset{j<i}{i,j=0}}^{(n-1)K+k}(\xi_{i} \! - \! 
\xi_{j})^{2}}{(\phi_{0}(\xi_{0}) \phi_{0}(\xi_{1}) \dotsb \phi_{0}(\xi_{l_{1}}) 
\dotsb \phi_{0}(\xi_{n_{1}}) \dotsb \phi_{0}(\xi_{n_{2}}) \phi_{0}(\xi_{m_{2}-
\varkappa_{nk}+1}) \dotsb \phi_{0}(\xi_{(n-1)K+k}))^{2}} \nonumber \\
=& \, \dfrac{(\mathbb{D}^{\lozenge})^{2}}{((n \! - \! 1)K \! + \! k \! + \! 1)!} 
\underbrace{\int_{\mathbb{R}} \int_{\mathbb{R}} \dotsb \int_{\mathbb{R}}
}_{(n-1)K+k+1} \md \widetilde{\mu}(\xi_{0}) \, \md \widetilde{\mu}(\xi_{1}) 
\, \dotsb \, \md \widetilde{\mu}(\xi_{l_{1}}) \, \dotsb \, \md \widetilde{\mu}
(\xi_{n_{1}}) \, \dotsb \, \md \widetilde{\mu}(\xi_{n_{2}}) \nonumber \\
\times& \, \md \widetilde{\mu}(\xi_{(n-1)K+k-\varkappa_{nk}+1}) \, \dotsb \, 
\md \widetilde{\mu}(\xi_{(n-1)K+k}) \left(\prod_{m=0}^{(n-1)K+k} \phi_{0}
(\xi_{m}) \right)^{-2} \prod_{\substack{i,j=0\\j<i}}^{(n-1)K+k}(\xi_{i} \! - \! 
\xi_{j})^{2} \nonumber \\
=& \, \dfrac{(\mathbb{D}^{\lozenge})^{2}}{((n \! - \! 1)K \! + \! k \! + \! 1)!} 
\underbrace{\int_{\mathbb{R}} \int_{\mathbb{R}} \dotsb \int_{\mathbb{R}}
}_{(n-1)K+k+1} \md \widetilde{\mu}(\xi_{0}) \, \md \widetilde{\mu}(\xi_{1}) \, 
\dotsb \, \md \widetilde{\mu}(\xi_{l_{1}}) \, \dotsb \, \md \widetilde{\mu}
(\xi_{n_{1}}) \, \dotsb \, \md \widetilde{\mu}(\xi_{n_{2}}) \nonumber \\
\times& \, \md \widetilde{\mu}(\xi_{(n-1)K+k-\varkappa_{nk}+1}) \, \dotsb \, 
\md \widetilde{\mu}(\xi_{(n-1)K+k}) \prod_{\substack{i,j=0\\j<i}}^{(n-1)K+k}
(\xi_{i} \! - \! \xi_{j})^{2} \left(\prod_{m=0}^{(n-1)K+k} \prod_{q=1}^{
\mathfrak{s}-2}(\xi_{m} \! - \! \alpha_{p_{q}})^{\varkappa_{nk \tilde{k}_{q}}}
(\xi_{m} \! - \! \alpha_{k})^{\varkappa_{nk}} 
\vphantom{\prod_{A}^{B}} \right)^{-2} \, \Rightarrow \nonumber \\
\hat{\mathfrak{c}}_{N_{f}} =& \, \dfrac{(\mathbb{D}^{\lozenge})^{2}}{((n \! - 
\! 1)K \! + \! k \! + \! 1)!} \underbrace{\int_{\mathbb{R}} \int_{\mathbb{R}} 
\dotsb \int_{\mathbb{R}}}_{(n-1)K+k+1} \md \widetilde{\mu}(\xi_{0}) \, \md 
\widetilde{\mu}(\xi_{1}) \, \dotsb \, \md \widetilde{\mu}(\xi_{(n-1)K+k}) \, 
\prod_{\substack{i,j=0\\j<i}}^{(n-1)K+k}(\xi_{i} \! - \! \xi_{j})^{2} \nonumber \\
\times& \, \left(\prod_{m=0}^{(n-1)K+k} \prod_{q=1}^{\mathfrak{s}-2}
(\xi_{m} \! - \! \alpha_{p_{q}})^{\varkappa_{nk \tilde{k}_{q}}}(\xi_{m} 
\! - \! \alpha_{k})^{\varkappa_{nk}} \right)^{-2} \quad (> \! 0). \label{eq63}
\end{align}

Now, $\hat{\mathfrak{c}}_{D_{f}}$ is studied. For $n \! \in \! \mathbb{N}$ and 
$k \! \in \! \lbrace 1,2,\dotsc,K \rbrace$ such that $\alpha_{p_{\mathfrak{s}}} 
\! := \! \alpha_{k} \! \neq \! \infty$, introduce the following notation:
\begin{equation*}
\tilde{\phi}_{0}(z) \! := \! \prod_{m=1}^{\mathfrak{s}-2}(z \! - \! 
\alpha_{p_{m}})^{l_{m}}(z \! - \! \alpha_{k})^{\varkappa_{nk}-1} 
\! =: \sum_{j=0}^{(n-1)K+k-1} \mathfrak{a}^{\sharp}_{j,0}z^{j},
\end{equation*}
and (with some abuse of notation), for $r \! = \! 1,\dotsc,\mathfrak{s} 
\! - \! 1,\mathfrak{s}$, $q(r) \! = \! \sum_{i=1}^{r-1}l_{i} \! + \! 1,
\sum_{i=1}^{r-1}l_{i} \! + \! 2,\dotsc,\sum_{i=1}^{r-1}l_{i} \! + \! l_{r} 
\! - \! \delta_{r \mathfrak{s}}$, and $m(r) \! = \! 1,2,\dotsc,l_{r} \! - \! 
\delta_{r \mathfrak{s}}$,
\begin{equation*}
\tilde{\phi}_{q(r)}(z) \! := \! \left\lbrace \tilde{\phi}_{0}(z)(z \! - \! 
\alpha_{p_{r}})^{-m(r)(1-\delta_{r \mathfrak{s}-1})}z^{m(r) \delta_{r 
\mathfrak{s}-1}} \! =: \sum_{j=0}^{(n-1)K+k-1} \mathfrak{a}^{
\sharp}_{j,q(r)}z^{j} \right\rbrace;
\end{equation*}
e.g., for $r \! = \! 1$, the notation $\tilde{\phi}_{q(1)}(z) \! := \! \lbrace 
\tilde{\phi}_{0}(z)(z \! - \! \alpha_{p_{1}})^{-m(1)} \! =: \! \sum_{j=0}^{(
n-1)K+k-1} \mathfrak{a}^{\sharp}_{j,q(1)}z^{j} \rbrace$, $q(1) \! = \! 
1,2,\dotsc,l_{1}$, $m(1) \! = \! 1,2,\dotsc,l_{1}$, denotes the (set of) 
$l_{1} \! = \! 
\varkappa_{nk \tilde{k}_{1}}$ functions
\begin{gather*}
\tilde{\phi}_{1}(z) \! = \! \dfrac{\tilde{\phi}_{0}(z)}{z \! - \! \alpha_{
p_{1}}} \! =: \sum_{j=0}^{(n-1)K+k-1} \mathfrak{a}^{\sharp}_{j,1}z^{j}, 
\, \tilde{\phi}_{2}(z) \! = \! \dfrac{\tilde{\phi}_{0}(z)}{(z \! - \! \alpha_{
p_{1}})^{2}} \! =: \sum_{j=0}^{(n-1)K+k-1} \mathfrak{a}^{\sharp}_{j,2}
z^{j}, \, \dotsc, \, \tilde{\phi}_{l_{1}}(z) \! = \! \dfrac{\tilde{\phi}_{0}(z)}{(
z \! - \! \alpha_{p_{1}})^{l_{1}}} \! =: \sum_{j=0}^{(n-1)K+k-1} 
\mathfrak{a}^{\sharp}_{j,l_{1}}z^{j},
\end{gather*}
for $r \! = \! \mathfrak{s} \! - \! 1$, the notation $\tilde{\phi}_{q
(\mathfrak{s}-1)}(z) \! := \! \lbrace \tilde{\phi}_{0}(z)z^{m(\mathfrak{s}-1)} 
\! =: \! \sum_{j=0}^{(n-1)K+k-1} \mathfrak{a}^{\sharp}_{j,q(\mathfrak{s}-1)}
z^{j} \rbrace$, $q(\mathfrak{s} \! - \! 1) \! = \! l_{1} \! + \! \dotsb \! + 
\! l_{\mathfrak{s}-2} \! + \! 1,l_{1} \! + \! \dotsb \! + \! l_{\mathfrak{s}-
2} \! + \! 2,\dotsc,l_{1} \! + \! \dotsb \! + \! l_{\mathfrak{s}-2} \! + \! 
l_{\mathfrak{s}-1} \! = \! (n \! - \! 1)K \! + \! k \! - \! \varkappa_{nk}$, 
$m(\mathfrak{s} \! - \! 1) \! = \! 1,2,\dotsc,l_{\mathfrak{s}-1}$, denotes 
the (set of) $l_{\mathfrak{s}-1} \! = \! \varkappa^{\infty}_{nk 
\tilde{k}_{\mathfrak{s}-1}}$ functions
\begin{gather*}
\tilde{\phi}_{l_{1}+\dotsb +l_{\mathfrak{s}-2}+1}(z) \! = \! \tilde{\phi}_{0}
(z)z \! =: \sum_{j=0}^{(n-1)K+k-1} \mathfrak{a}^{\sharp}_{j,l_{1}+\dotsb +
l_{\mathfrak{s}-2}+1}z^{j}, \, \tilde{\phi}_{l_{1}+\dotsb +l_{\mathfrak{s}-2}
+2}(z) \! = \! \tilde{\phi}_{0}(z)z^{2} \! =: \sum_{j=0}^{(n-1)K+k-1} 
\mathfrak{a}^{\sharp}_{j,l_{1}+\dotsb +l_{\mathfrak{s}-2}+2}z^{j}, \, 
\dotsc \\\dotsc, \, \tilde{\phi}_{(n-1)K+k-\varkappa_{nk}}(z) \! = \! 
\tilde{\phi}_{0}(z)z^{l_{\mathfrak{s}-1}} \! =: \sum_{j=0}^{(n-1)K+k-1} 
\mathfrak{a}^{\sharp}_{j,(n-1)K+k-\varkappa_{nk}}z^{j},
\end{gather*}
etc., and, for $r \! = \! \mathfrak{s}$, the notation $\tilde{\phi}_{q
(\mathfrak{s})}(z) \! := \! \lbrace \tilde{\phi}_{0}(z)(z \! - \! \alpha_{k}
)^{-m(\mathfrak{s})} \! =: \! \sum_{j=0}^{(n-1)K+k-1} \mathfrak{a}^{\sharp}_{
j,q(\mathfrak{s})}z^{j} \rbrace$, $q(\mathfrak{s}) \! = \! (n \! - \! 1)K 
\! + \! k \! - \! \varkappa_{nk} \! + \! 1,(n \! - \! 1)K \! + \! k \! - \! 
\varkappa_{nk} \! + \! 2,\dotsc,(n \! - \! 1)K \! + \! k \! - \! 1$, $m
(\mathfrak{s}) \! = \! 1,2,\dotsc,l_{\mathfrak{s}} \! - \! 1$, denotes 
the (set of) $l_{\mathfrak{s}} \! - \! 1 \! = \! \varkappa_{nk} \! - \! 1$ 
functions
\begin{gather*}
\tilde{\phi}_{(n-1)K+k-\varkappa_{nk}+1}(z) \! = \! \dfrac{\tilde{\phi}_{0}
(z)}{z \! - \! \alpha_{k}} \! =: \sum_{j=0}^{(n-1)K+k-1} \mathfrak{a}^{
\sharp}_{j,(n-1)K+k-\varkappa_{nk}+1}z^{j}, \, \tilde{\phi}_{(n-1)K+k-
\varkappa_{nk}+2}(z) \! = \! \dfrac{\tilde{\phi}_{0}(z)}{(z \! - \! 
\alpha_{k})^{2}} \! =: \sum_{j=0}^{(n-1)K+k-1} \mathfrak{a}^{\sharp}_{j,
(n-1)K+k-\varkappa_{nk}+2}z^{j}, \, \dotsc \\
\dotsc, \, \tilde{\phi}_{(n-1)K+k-1}(z) \! = \! \dfrac{\tilde{\phi}_{0}(z)}{
(z \! - \! \alpha_{k})^{\varkappa_{nk}-1}} \! =: \sum_{j=0}^{(n-1)K+k-1} 
\mathfrak{a}^{\sharp}_{j,(n-1)K+k-1}z^{j}.
\end{gather*}
(Note: $\# \lbrace \tilde{\phi}_{0}(z)(z \! - \! \alpha_{p_{r}})^{-m(r)
(1-\delta_{r \mathfrak{s}-1})}z^{m(r) \delta_{r \mathfrak{s}-1}} \rbrace 
\! = \! l_{r} \! - \! \delta_{r \mathfrak{s}}$, $r \! = \! 1,2,\dotsc,
\mathfrak{s}$, and $\# \cup_{r=1}^{\mathfrak{s}} \lbrace \tilde{\phi}_{0}(z)
(z \! - \! \alpha_{p_{r}})^{-m(r)(1-\delta_{r \mathfrak{s}-1})}z^{m(r) 
\delta_{r \mathfrak{s}-1}} \rbrace \! = \! \sum_{r=1}^{\mathfrak{s}}l_{r} 
\! - \! 1 \! = \! (n \! - \! 1)K \! + \! k \! - \! 1$.) One notes that, for 
$n \! \in \! \mathbb{N}$ and $k \! \in \! \lbrace 1,2,\dotsc,K \rbrace$ 
such that $\alpha_{p_{\mathfrak{s}}} \! := \! \alpha_{k} \! \neq \! \infty$, 
the $l_{1} \! + \! \dotsb \! + \! l_{\mathfrak{s}-1} \! + \! l_{\mathfrak{s}} 
\! = \! (n \! - \! 1)K \! + \! k$ functions $\tilde{\phi}_{0}(z),\tilde{\phi}_{1}
(z),\dotsc,\tilde{\phi}_{l_{1}}(z),\dotsc,\tilde{\phi}_{l_{1}+\dotsb 
+l_{\mathfrak{s}-2}+1}(z),\dotsc,\tilde{\phi}_{(n-1)K+k-\varkappa_{nk}}(z),
\tilde{\phi}_{(n-1)K+k-\varkappa_{nk}+1}(z),\dotsc,\tilde{\phi}_{(n-1)K+k-1}
(z)$ are linearly independent on $\mathbb{R}$, that is, for $z \! \in \! 
\mathbb{R}$, $\sum_{j=0}^{(n-1)K+k-1} \mathfrak{c}^{\sharp}_{j} \, 
\tilde{\phi}_{j}(z) \! = \! 0$ $\Rightarrow$ (via a Vandermonde-type argument; 
see the $((n \! - \! 1)K \! + \! k) \times ((n \! - \! 1)K \! + \! k)$ non-zero 
determinant $\mathbb{D}^{\blacklozenge}$ in Equation~\eqref{eq64} below) 
$\mathfrak{c}^{\sharp}_{j} \! = \! 0$, $j \! = \! 0,1,\dotsc,(n \! - \! 1)K \! + \! 
k \! - \! 1$. For $n \! \in \! \mathbb{N}$ and $k \! \in \! \lbrace 1,2,\dotsc,K 
\rbrace$ such that $\alpha_{p_{\mathfrak{s}}} \! := \! \alpha_{k} \! \neq \! \infty$, 
let $\mathfrak{S}_{(n-1)K+k}$ denote the $((n \! - \! 1)K \! + \! k)!$ permutations 
of $\lbrace 0,1,\dotsc,(n \! - \! 1)K \! + \! k \! - \! 1 \rbrace$. Using the above 
notation and the multi-linearity property of the determinant, one proceeds to 
study, for $n \! \in \! \mathbb{N}$ and $k \! \in \! \lbrace 1,2,\dotsc,K \rbrace$ 
such that $\alpha_{p_{\mathfrak{s}}} \! := \! \alpha_{k} \! \neq \! \infty$, 
$\hat{\mathfrak{c}}_{D_{f}}$:
\begin{align*}
\hat{\mathfrak{c}}_{D_{f}} =& \, \underbrace{\int_{\mathbb{R}} \int_{\mathbb{
R}} \dotsb \int_{\mathbb{R}}}_{(n-1)K+k} \md \widetilde{\mu}(\xi_{0}) \, \md 
\widetilde{\mu}(\xi_{1}) \, \dotsb \, \md \widetilde{\mu}(\xi_{l_{1}}) \, \dotsb 
\, \md \widetilde{\mu}(\xi_{n_{1}}) \, \dotsb \, \md \widetilde{\mu}(\xi_{n_{2}}) 
\, \md \widetilde{\mu}(\xi_{m_{2}-\varkappa_{nk}+1}) \, \dotsb \, \md 
\widetilde{\mu}(\xi_{m_{1}}) \\
\times& \, \dfrac{1}{(\xi_{0} \! - \! \alpha_{k})^{0}(\xi_{1} \! - \! 
\alpha_{p_{1}})^{1} \dotsb (\xi_{l_{1}} \! - \! \alpha_{p_{1}})^{l_{1}} 
\dotsb \dfrac{1}{(\xi_{n_{1}})^{1}} \dotsb \dfrac{1}{(\xi_{n_{2}})^{
l_{\mathfrak{s}-1}}}(\xi_{m_{2}-\varkappa_{nk}+1} \! - \! \alpha_{k})^{1} 
\dotsb (\xi_{m_{1}} \! - \! \alpha_{k})^{\varkappa_{nk}-1}} \\
\times& \, \dfrac{1}{\tilde{\phi}_{0}(\xi_{0}) \tilde{\phi}_{0}(\xi_{1}) 
\dotsb \tilde{\phi}_{0}(\xi_{l_{1}}) \dotsb \tilde{\phi}_{0}(\xi_{n_{1}}) 
\dotsb \tilde{\phi}_{0}(\xi_{n_{2}}) \tilde{\phi}_{0}(\xi_{m_{2}-\varkappa_{n
k}+1}) \dotsb \tilde{\phi}_{0}(\xi_{m_{1}})} \\
\times& \, 
\underbrace{
\left\vert

\right\vert}_{=: \, \mathbb{G}^{\blacklozenge}(\xi_{0},\xi_{1},\dotsc,\xi_{l_{
1}},\dotsc,\xi_{n_{1}},\dotsc,\xi_{n_{2}},\xi_{m_{2}-\varkappa_{nk}+1},\dotsc,
\xi_{(n-1)K+k-1})} \\
=& \, \dfrac{1}{m_{2}!} \sum_{\pmb{\sigma} \in \mathfrak{S}_{m_{2}}} 
\underbrace{\int_{\mathbb{R}} \int_{\mathbb{R}} \dotsb \int_{\mathbb{R}}}_{(n-1)K
+k} \md \widetilde{\mu}(\xi_{\sigma (0)}) \, \md \widetilde{\mu}(\xi_{\sigma (1)}) 
\, \dotsb \, \md \widetilde{\mu}(\xi_{\sigma (l_{1})}) \, \dotsb \, \md \widetilde{\mu}
(\xi_{\sigma (n_{1})}) \, \dotsb \, \md \widetilde{\mu}(\xi_{\sigma (n_{2})}) \\
\times& \, \tfrac{\dotsb \, \md \widetilde{\mu}(\xi_{\sigma (m_{2}-\varkappa_{nk}+1)}) 
\, \dotsb \, \md \widetilde{\mu}(\xi_{\sigma (m_{1})})}{(\xi_{\sigma (0)}-\alpha_{k})^{0}
(\xi_{\sigma (1)}-\alpha_{p_{1}})^{1} \dotsb (\xi_{\sigma (l_{1})}-\alpha_{
p_{1}})^{l_{1}} \dotsb \frac{1}{(\xi_{\sigma (n_{1})})^{1}} \dotsb \frac{1}{
(\xi_{\sigma (n_{2})})^{l_{\mathfrak{s}-1}}}(\xi_{\sigma (m_{2}-\varkappa_{nk}
+1)}-\alpha_{k})^{1} \dotsb (\xi_{\sigma (m_{1})}-\alpha_{k})^{\varkappa_{nk}
-1}} \\
\times& \, \dfrac{1}{\tilde{\phi}_{0}(\xi_{\sigma (0)}) \tilde{\phi}_{0}(\xi_{
\sigma (1)}) \dotsb \tilde{\phi}_{0}(\xi_{\sigma (l_{1})}) \dotsb \tilde{
\phi}_{0}(\xi_{\sigma (n_{1})}) \dotsb \tilde{\phi}_{0}(\xi_{\sigma (n_{2})}) 
\tilde{\phi}_{0}(\xi_{\sigma (m_{2}-\varkappa_{nk}+1)}) \dotsb \tilde{\phi}_{
0}(\xi_{\sigma (m_{1})})} \\
\times& \, \mathbb{G}^{\blacklozenge}(\xi_{\sigma (0)},\xi_{\sigma (1)},\dotsc,
\xi_{\sigma (l_{1})},\dotsc,\xi_{\sigma (n_{1})},\dotsc,\xi_{\sigma (n_{2})},
\xi_{\sigma (m_{2}-\varkappa_{nk}+1)},\dotsc,\xi_{\sigma ((n-1)K+k-1)}) \\
=& \, \dfrac{1}{m_{2}!} \underbrace{\int_{\mathbb{R}} \int_{\mathbb{R}} \dotsb 
\int_{\mathbb{R}}}_{(n-1)K+k} \md \widetilde{\mu}(\xi_{0}) \, \md \widetilde{\mu}
(\xi_{1}) \, \dotsb \, \md \widetilde{\mu}(\xi_{l_{1}}) \, \dotsb \, \md \widetilde{\mu}
(\xi_{n_{1}}) \, \dotsb \, \md \widetilde{\mu}(\xi_{n_{2}}) \, \md \widetilde{\mu}
(\xi_{m_{2}-\varkappa_{nk}+1}) \, \dotsb \\
\times& \, \dotsb \, \md \widetilde{\mu}(\xi_{m_{1}}) \, \mathbb{G}^{\blacklozenge}
(\xi_{0},\xi_{1},\dotsc,\xi_{l_{1}},\dotsc,\xi_{n_{1}},\dotsc,\xi_{n_{2}},\xi_{m_{2}-
\varkappa_{nk}+1},\dotsc,\xi_{(n-1)K+k-1}) \\
\times& \, \sum_{\pmb{\sigma} \in \mathfrak{S}_{m_{2}}} \operatorname{sgn}
(\pmb{\pmb{\sigma}}) \tfrac{1}{(\xi_{\sigma (0)}-\alpha_{k})^{0}
(\xi_{\sigma (1)}-\alpha_{p_{1}})^{1} \dotsb (\xi_{\sigma (l_{1})}-
\alpha_{p_{1}})^{l_{1}} \dotsb \frac{1}{(\xi_{\sigma (n_{1})})^{1}} \dotsb 
\frac{1}{(\xi_{\sigma (n_{2})})^{l_{\mathfrak{s}-1}}}(\xi_{\sigma (m_{2}-
\varkappa_{nk}+1)}-\alpha_{k})^{1} \dotsb (\xi_{\sigma (m_{1})}-
\alpha_{k})^{\varkappa_{nk}-1}} \\
\times& \, \dfrac{1}{\tilde{\phi}_{0}(\xi_{\sigma (0)}) \tilde{\phi}_{0}(\xi_{
\sigma (1)}) \dotsb \tilde{\phi}_{0}(\xi_{\sigma (l_{1})}) \dotsb \tilde{
\phi}_{0}(\xi_{\sigma (n_{1})}) \dotsb \tilde{\phi}_{0}(\xi_{\sigma (n_{2})}) 
\tilde{\phi}_{0}(\xi_{\sigma (m_{2}-\varkappa_{nk}+1)}) \dotsb \tilde{\phi}_{
0}(\xi_{\sigma (m_{1})})} \\
=& \, \dfrac{1}{m_{2}!} \underbrace{\int_{\mathbb{R}} \int_{\mathbb{R}} \dotsb 
\int_{\mathbb{R}}}_{(n-1)K+k} \md \widetilde{\mu}(\xi_{0}) \, \md \widetilde{\mu}
(\xi_{1}) \, \dotsb \, \md \widetilde{\mu}(\xi_{l_{1}}) \, \dotsb \, \md \widetilde{\mu}
(\xi_{n_{1}}) \, \dotsb \, \md \widetilde{\mu}(\xi_{n_{2}}) \, \md \widetilde{\mu}
(\xi_{m_{2}-\varkappa_{nk}+1}) \, \dotsb \\
\times& \, \dotsb \, \md \widetilde{\mu}(\xi_{m_{1}}) \, \mathbb{G}^{\blacklozenge}
(\xi_{0},\xi_{1},\dotsc,\xi_{l_{1}},\dotsc,\xi_{n_{1}},\dotsc,\xi_{n_{2}},\xi_{m_{2}
-\varkappa_{nk}+1},\dotsc,\xi_{(n-1)K+k-1}) \\
\times& \, 
\left\lvert

\right\rvert \\
=& \, \dfrac{1}{m_{2}!} \underbrace{\int_{\mathbb{R}} \int_{\mathbb{R}} \dotsb 
\int_{\mathbb{R}}}_{(n-1)K+k} \md \widetilde{\mu}(\xi_{0}) \, \md \widetilde{\mu}
(\xi_{1}) \, \dotsb \, \md \widetilde{\mu}(\xi_{l_{1}}) \, \dotsb \, \md \widetilde{\mu}
(\xi_{n_{1}}) \, \dotsb \, \md \widetilde{\mu}(\xi_{n_{2}}) \, \md \widetilde{\mu}
(\xi_{m_{2}-\varkappa_{nk}+1}) \, \dotsb \\
\times& \, \dotsb \, \md \widetilde{\mu}(\xi_{m_{1}}) \, \dfrac{\mathbb{G}^{
\blacklozenge}(\xi_{0},\xi_{1},\dotsc,\xi_{l_{1}},\dotsc,\xi_{n_{1}},\dotsc,\xi_{n_{2}},
\xi_{m_{2}-\varkappa_{nk}+1},\dotsc,\xi_{(n-1)K+k-1})}{(\tilde{\phi}_{0}
(\xi_{0}) \tilde{\phi}_{0}(\xi_{1}) \dotsb \tilde{\phi}_{0}(\xi_{l_{1}}) 
\dotsb \tilde{\phi}_{0}(\xi_{n_{1}}) \dotsb \tilde{\phi}_{0}(\xi_{n_{2}}) 
\tilde{\phi}_{0}(\xi_{m_{2}-\varkappa_{nk}+1}) \dotsb \tilde{\phi}_{0}
(\xi_{m_{1}}))^{2}} \\
\times& \, 
\underbrace{
\left\vert

\right\vert}_{= \, \mathbb{G}^{\blacklozenge}(\xi_{0},\xi_{1},\dotsc,\xi_{l_{
1}},\dotsc,\xi_{n_{1}},\dotsc,\xi_{n_{2}},\xi_{m_{2}-\varkappa_{nk}+1},\dotsc,
\xi_{(n-1)K+k-1})} \\
=& \, \dfrac{1}{m_{2}!} \underbrace{\int_{\mathbb{R}} \int_{\mathbb{R}} \dotsb 
\int_{\mathbb{R}}}_{(n-1)K+k} \md \widetilde{\mu}(\xi_{0}) \, \md \widetilde{\mu}
(\xi_{1}) \, \dotsb \, \md \widetilde{\mu}(\xi_{l_{1}}) \, \dotsb \, \md \widetilde{\mu}
(\xi_{n_{1}}) \, \dotsb \, \md \widetilde{\mu}(\xi_{n_{2}}) \, \md \widetilde{\mu}
(\xi_{m_{2}-\varkappa_{nk}+1}) \, \dotsb \\
\times& \, \dotsb \, \md \widetilde{\mu}(\xi_{m_{1}}) \left(\dfrac{\mathbb{G}^{
\blacklozenge}(\xi_{0},\xi_{1},\dotsc,\xi_{l_{1}},\dotsc,\xi_{n_{1}},\dotsc,
\xi_{n_{2}},\xi_{m_{2}-\varkappa_{nk}+1},\dotsc,\xi_{(n-1)K+k-1})}{\tilde{
\phi}_{0}(\xi_{0}) \tilde{\phi}_{0}(\xi_{1}) \dotsb \tilde{\phi}_{0}(\xi_{l_{
1}}) \dotsb \tilde{\phi}_{0}(\xi_{n_{1}}) \dotsb \tilde{\phi}_{0}(\xi_{n_{2}}) 
\tilde{\phi}_{0}(\xi_{m_{2}-\varkappa_{nk}+1}) \dotsb \tilde{\phi}_{0}(\xi_{
m_{1}})} \right)^{2};
\end{align*}
but, noting the determinantal factorisation
\begin{equation*}
\mathbb{G}^{\blacklozenge}(\xi_{0},\xi_{1},\dotsc,\xi_{l_{1}},\dotsc,
\xi_{n_{1}},\dotsc,\xi_{n_{2}},\xi_{m_{2}-\varkappa_{nk}+1},\dotsc,
\xi_{(n-1)K+k-1}) \! := \! \overset{f}{\mathcal{V}}_{2}(\xi_{0},\xi_{1},
\dotsc,\xi_{(n-1)K+k-1}) \mathbb{D}^{\blacklozenge},
\end{equation*}
where
\begin{equation*}
\overset{f}{\mathcal{V}}_{2}(\xi_{0},\xi_{1},\dotsc,\xi_{(n-1)K+k-1}) 
= \left\lvert

\right), \label{eq64}
\end{equation}
with
\begin{equation*}
(\mathfrak{A}^{\blacklozenge}(0))_{1j} \! = \! 
\begin{cases}
\mathfrak{a}_{j-1}^{\blacklozenge}(\vec{\bm{\alpha}}), &\text{$j \! 
= \! 1,2,\dotsc,(n \! - \! 1)K \! + \! k \! - \! l_{\mathfrak{s}-1},$} \\
0, &\text{$j \! = \! (n \! - \! 1)K \! + \! k \! - \! l_{\mathfrak{s}-1} 
\! + \! 1,\dotsc,(n \! - \! 1)K \! + \! k,$}
\end{cases}
\end{equation*}
for $r \! = \! 1,\dotsc,\mathfrak{s} \! - \! 2,\mathfrak{s}$,
\begin{align*}
i(r)=& \, 2 \! + \! \sum_{m=1}^{r-1}l_{m},3 \! + \! \sum_{m=1}^{r-1}l_{m},
\dotsc,1 \! - \! \delta_{r \mathfrak{s}} \! + \! l_{r} \! + \! 
\sum_{m=1}^{r-1}l_{m}, \qquad \qquad q(i(r),r) \! = \! 1,2,\dotsc,
l_{r} \! - \! \delta_{r \mathfrak{s}}, \\
(\mathfrak{A}^{\blacklozenge}(r))_{i(r)j(r)}=& \, 
\begin{cases}
\dfrac{(-1)^{q(i(r),r)}}{\prod_{m=0}^{q(i(r),r)-1}(l_{r} \! - \! m \! - \! 
\delta_{r \mathfrak{s}})} \left(\dfrac{\partial}{\partial \alpha_{p_{r}}} 
\right)^{q(i(r),r)} \mathfrak{a}_{j(r)-1}^{\blacklozenge}(\vec{\bm{\alpha}}), 
&\text{$j(r) \! = \! 1,2,\dotsc,(n \! - \! 1)K \! + \! k \! - \! 
l_{\mathfrak{s}-1} \! - \! q(i(r),r)$,} \\
0, &\text{$j(r) \! = \! (n \! - \! 1)K \! + \! k \! - \! l_{\mathfrak{s}-1} 
\! - \! q(i(r),r) \! + \! 1,\dotsc,(n \! - \! 1)K \! + \! k$,}
\end{cases}
\end{align*}
and, for $r \! = \! \mathfrak{s} \! - \! 1$,
\begin{align*}
i(\mathfrak{s} \! - \! 1)=& \, 2 \! + \! \sum_{m=1}^{\mathfrak{s}-2}l_{m},3 
\! + \! \sum_{m=1}^{\mathfrak{s}-2}l_{m},\dotsc,1 \! + \! l_{\mathfrak{s}-1} 
\! + \! \sum_{m=1}^{\mathfrak{s}-2}l_{m}, \quad \quad q(i(\mathfrak{s} \! 
- \! 1),\mathfrak{s} \! - \! 1) \! = \! 1,2,\dotsc,l_{\mathfrak{s}-1}, \\
(\mathfrak{A}^{\blacklozenge}(\mathfrak{s} \! - \! 1))_{i(\mathfrak{s}-1)
j(\mathfrak{s}-1)}=& \, 
\begin{cases}
0, &\text{$j(\mathfrak{s} \! - \! 1) \! = \! 1,\dotsc,q(i(\mathfrak{s} 
\! - \! 1),\mathfrak{s} \! - \! 1)$,} \\
\mathfrak{a}_{j(\mathfrak{s}-1)-q(i(\mathfrak{s}-1),\mathfrak{s}-1)
-1}^{\blacklozenge}(\vec{\bm{\alpha}}), &\text{$j(\mathfrak{s} \! - \! 1) 
\! = \! q(i(\mathfrak{s} \! - \! 1),\mathfrak{s} \! - \! 1) \! + \! 1,
\dotsc,i(\mathfrak{s} \! - \! 1) \! + \! \varkappa_{nk} \! - \! 1$,} \\
0, &\text{$j(\mathfrak{s} \! - \! 1) \! = \! i(\mathfrak{s} \! - \! 1) 
\! + \! \varkappa_{nk},\dotsc,(n \! - \! 1)K \! + \! k$,}
\end{cases}
\end{align*}
where, for $\tilde{m}_{1} \! = \! 0,1,\dotsc,(n \! - \! 1)K \! + \! k \! - \! 
l_{\mathfrak{s}-1} \! - \! 1$,
\begin{equation*}
\mathfrak{a}_{\tilde{m}_{1}}^{\blacklozenge}(\vec{\bm{\alpha}}) := \, 
\mathlarger{\sum_{\underset{\underset{\underset{\sum_{\underset{m \neq 
\mathfrak{s}-1}{m=1}}^{\mathfrak{s}}i_{m}=(n-1)K+k-l_{\mathfrak{s}-1}-
\tilde{m}_{1}-1}{i_{\mathfrak{s}}=0,1,\dotsc,l_{\mathfrak{s}}-1}}{p \in 
\lbrace 1,2,\dotsc,\mathfrak{s}-2 \rbrace}}{i_{p}=0,1,\dotsc,l_{p}}}}
(-1)^{(n-1)K+k-l_{\mathfrak{s}-1}-\tilde{m}_{1}-1} \prod_{j=1}^{
\mathfrak{s}-2} \binom{l_{j}}{i_{j}} \binom{l_{\mathfrak{s}} \! - \! 1}{
i_{\mathfrak{s}}} \prod_{m=1}^{\mathfrak{s}-2}(\alpha_{p_{m}})^{i_{m}} 
\alpha_{k}^{i_{\mathfrak{s}}},
\end{equation*}
it follows that, for $n \! \in \! \mathbb{N}$ and $k \! \in \! \lbrace 
1,2,\dotsc,K \rbrace$ such that $\alpha_{p_{\mathfrak{s}}} \! := \! 
\alpha_{k} \! \neq \! \infty$,
\begin{align}
\hat{\mathfrak{c}}_{D_{f}} =& \, \dfrac{1}{m_{2}!} \underbrace{\int_{\mathbb{
R}} \int_{\mathbb{R}} \dotsb \int_{\mathbb{R}}}_{(n-1)K+k} \md \widetilde{\mu}
(\xi_{0}) \, \md \widetilde{\mu}(\xi_{1}) \, \dotsb \, \md \widetilde{\mu}
(\xi_{l_{1}}) \, \dotsb \, \md \widetilde{\mu}(\xi_{n_{1}}) \, \dotsb \, \md 
\widetilde{\mu}(\xi_{n_{2}}) \, \md \widetilde{\mu}(\xi_{m_{2}-\varkappa_{nk}+1}) 
\, \dotsb \nonumber \\
\times& \, \dotsb \, \md \widetilde{\mu}(\xi_{m_{1}}) \dfrac{(\mathbb{D}^{
\blacklozenge})^{2} \prod_{\underset{j<i}{i,j=0}}^{(n-1)K+k-1}(\xi_{i} \! - \! 
\xi_{j})^{2}}{(\tilde{\phi}_{0}(\xi_{0}) \tilde{\phi}_{0}(\xi_{1}) \dotsb \tilde{\phi}_{0}
(\xi_{l_{1}}) \dotsb \tilde{\phi}_{0}(\xi_{n_{1}}) \dotsb \tilde{\phi}_{0}
(\xi_{n_{2}}) \tilde{\phi}_{0}(\xi_{m_{2}-\varkappa_{nk}+1}) \dotsb \tilde{
\phi}_{0}(\xi_{(n-1)K+k-1}))^{2}} \nonumber \\
=& \, \dfrac{(\mathbb{D}^{\blacklozenge})^{2}}{((n \! - \! 1)K \! + \! k)!} 
\underbrace{\int_{\mathbb{R}} \int_{\mathbb{R}} \dotsb \int_{\mathbb{R}}}_{(
n-1)K+k} \md \widetilde{\mu}(\xi_{0}) \, \md \widetilde{\mu}(\xi_{1}) \, \dotsb 
\, \md \widetilde{\mu}(\xi_{l_{1}}) \, \dotsb \, \md \widetilde{\mu}(\xi_{n_{1}}) 
\, \dotsb \, \md \widetilde{\mu}(\xi_{n_{2}}) \nonumber \\
\times& \, \md \widetilde{\mu}(\xi_{(n-1)K+k-\varkappa_{nk}+1}) \, \dotsb 
\, \md \widetilde{\mu}(\xi_{(n-1)K+k-1}) \left(\prod_{m=0}^{(n-1)K+k-1} 
\tilde{\phi}_{0}(\xi_{m}) \right)^{-2} \prod_{\substack{i,j=0\\j<i}}^{(n-1)K+k-1}
(\xi_{i} \! - \! \xi_{j})^{2} \nonumber \\
=& \, \dfrac{(\mathbb{D}^{\blacklozenge})^{2}}{((n \! - \! 1)K \! + \! k)!} 
\underbrace{\int_{\mathbb{R}} \int_{\mathbb{R}} \dotsb \int_{\mathbb{R}}}_{(
n-1)K+k} \md \widetilde{\mu}(\xi_{0}) \, \md \widetilde{\mu}(\xi_{1}) \, \dotsb 
\, \md \widetilde{\mu}(\xi_{l_{1}}) \, \dotsb \, \md \widetilde{\mu}(\xi_{n_{1}}) 
\, \dotsb \, \md \widetilde{\mu}(\xi_{n_{2}}) \nonumber \\
\times& \, \md \widetilde{\mu}(\xi_{(n-1)K+k-\varkappa_{nk}+1}) \, \dotsb \, 
\md \widetilde{\mu}(\xi_{(n-1)K+k-1}) \prod_{\substack{i,j=0\\j<i}}^{(n-1)K+k-1}
(\xi_{i} \! - \! \xi_{j})^{2} \left(\prod_{m=0}^{(n-1)K+k-1} \prod_{q=1}^{
\mathfrak{s}-2}(\xi_{m} \! - \! \alpha_{p_{q}})^{\varkappa_{nk \tilde{k}_{q}}}
(\xi_{m} \! - \! \alpha_{k})^{\varkappa_{nk}-1} \vphantom{\prod_{A}^{B}} 
\right)^{-2} \, \Rightarrow \nonumber \\
\hat{\mathfrak{c}}_{D_{f}} =& \, \dfrac{(\mathbb{D}^{\blacklozenge})^{2}}{(
(n \! - \! 1)K \! + \! k)!} \underbrace{\int_{\mathbb{R}} \int_{\mathbb{R}} 
\dotsb \int_{\mathbb{R}}}_{(n-1)K+k} \md \widetilde{\mu}(\xi_{0}) \, \md 
\widetilde{\mu}(\xi_{1}) \, \dotsb \, \md \widetilde{\mu}(\xi_{(n-1)K+k-1}) 
\prod_{\substack{i,j=0\\j<i}}^{(n-1)K+k-1}(\xi_{i} \! - \! \xi_{j})^{2} \nonumber \\
\times& \, \left(\prod_{m=0}^{(n-1)K+k-1} \prod_{q=1}^{\mathfrak{s}
-2}(\xi_{m} \! - \! \alpha_{p_{q}})^{\varkappa_{nk \tilde{k}_{q}}}
(\xi_{m} \! - \! \alpha_{k})^{\varkappa_{nk}-1} \right)^{-2} \quad (> \! 0). 
\label{eq65}
\end{align}
Hence, for $n \! \in \! \mathbb{N}$ and $k \! \in \! \lbrace 1,2,\dotsc,K 
\rbrace$ such that $\alpha_{p_{\mathfrak{s}}} \! := \! \alpha_{k} \! \neq \! 
\infty$, Equations~\eqref{eq63} and~\eqref{eq65} establish the existence 
and the uniqueness of the corresponding monic MPC ORF, $\pmb{\pi}^{n}_{k} 
\colon \mathbb{N} \times \lbrace 1,2,\dotsc,K \rbrace \times \overline{\mathbb{C}} 
\setminus \lbrace \alpha_{p_{1}},\alpha_{p_{2}},\dotsc,\alpha_{p_{\mathfrak{s}}} 
\rbrace \! \to \! \mathbb{C}$, $(n,k,z) \! \mapsto \! \pmb{\pi}^{n}_{k}(z) 
\! = \! \mathcal{X}_{11}(z)/(z \! - \! \alpha_{k})$.

There remains, still, the question of the existence and the uniqueness of 
$\mathcal{X}_{21} \colon \mathbb{N} \times \lbrace 1,2,\dotsc,K \rbrace 
\times \overline{\mathbb{C}} \setminus \lbrace \alpha_{p_{1}},\alpha_{p_{2}},
\dotsc,
\linebreak[4] 
\alpha_{p_{\mathfrak{s}}} \rbrace \! \to \! \mathbb{C}$, which necessitates 
an analysis for the elements of the second row of $\mathcal{X}(z)$, 
that is, $\mathcal{X}_{2m}(z)$, $m \! = \! 1,2$ (thus the gist of the 
subsequent calculations). Recalling that $\int_{\mathbb{R}} \xi^{m} 
\widetilde{w}(\xi) \, \md \xi \! < \! \infty$, $m \! \in \! \mathbb{N}_{0}$, 
and $\int_{\mathbb{R}}(\xi \! - \! \alpha_{k})^{-(m+1)} \widetilde{w}(\xi) 
\, \md \xi \! < \! \infty$, $\alpha_{k} \! \neq \! \infty$, $k \! \in \! \lbrace 
1,2,\dotsc,K \rbrace$, it follows, via the expansion $\tfrac{1}{z_{1}-z_{2}} \! 
= \! \sum_{i=0}^{l} \tfrac{z_{2}^{i}}{z_{1}^{i+1}} \! + \! \tfrac{z_{2}^{l+1}}{
z_{1}^{l+1}(z_{1}-z_{2})}$, $l \! \in \! \mathbb{N}_{0}$, the integral representation 
(cf. Equation~\eqref{eq52}) $\mathcal{X}_{22}(z) \! = \! \int_{\mathbb{R}} 
\tfrac{(z-\alpha_{k}) \mathcal{X}_{21}(\xi) \widetilde{w}(\xi)}{(\xi -
\alpha_{k})(\xi -z)} \, \tfrac{\md \xi}{2 \pi \mi}$, $z \! \in \! \mathbb{C} 
\setminus \mathbb{R}$, and the second line of each of the asymptotic 
conditions~\eqref{eq49}, \eqref{eq50}, and~\eqref{eq51}, that, for $n \! 
\in \! \mathbb{N}$ and $k \! \in \! \lbrace 1,2,\dotsc,K \rbrace$ such 
that $\alpha_{p_{\mathfrak{s}}} \! := \! \alpha_{k} \! \neq \! \infty$,
\begin{gather}
\int_{\mathbb{R}} \left(\dfrac{\mathcal{X}_{21}(\xi)}{\xi \! - \! \alpha_{k}} 
\right) \dfrac{\widetilde{w}(\xi)}{(\xi \! - \! \alpha_{k})^{p}} \, \md \xi \! 
= \! 0, \quad p \! = \! 0,1,\dotsc,\varkappa_{nk} \! - \! 2, \label{eq66} \\
\int_{\mathbb{R}} \left(\dfrac{\mathcal{X}_{21}(\xi)}{\xi \! - \! \alpha_{k}} 
\right) \dfrac{\widetilde{w}(\xi)}{(\xi \! - \! \alpha_{k})^{\varkappa_{nk}-
1}} \, \md \xi \! = \! 2 \pi \mi, \label{eq67} \\
\int_{\mathbb{R}} \left(\dfrac{\mathcal{X}_{21}(\xi)}{\xi \! - \! \alpha_{k}} 
\right) \dfrac{\widetilde{w}(\xi)}{(\xi \! - \! \alpha_{p_{q}})^{r}} \, \md 
\xi \! = \! 0, \quad q \! = \! 1,2,\dotsc,\mathfrak{s} \! - \! 2, \quad r 
\! = \! 1,2,\dotsc,\varkappa_{nk \tilde{k}_{q}}, \label{eq68} \\
\int_{\mathbb{R}} \left(\dfrac{\mathcal{X}_{21}(\xi)}{\xi \! - \! \alpha_{k}} 
\right) \xi^{r} \widetilde{w}(\xi) \, \md \xi \! = \! 0, \quad r \! = \! 
1,2,\dotsc,\varkappa_{nk \tilde{k}_{\mathfrak{s}-1}}^{\infty}. \label{eq69}
\end{gather}
(Note: for $n \! = \! 1$ and $k \! \in \! \lbrace 1,2,\dotsc,K \rbrace$ such 
that $\alpha_{p_{\mathfrak{s}}} \! := \! \alpha_{k} \! \neq \! \infty$, it can 
happen that the corresponding $\varkappa_{1k} \! < \! 2$, in which case, 
Equation~\eqref{eq66} is vacuous (of course, for $n \! \geqslant \! 2$, $k 
\! \in \! \lbrace 1,2,\dotsc,K \rbrace$, $\varkappa_{nk} \! \geqslant \! 2$); 
moreover, if, for $n \! \in \! \mathbb{N}$ and $k \! \in \! \lbrace 1,2,\dotsc,K 
\rbrace$ such that $\alpha_{p_{\mathfrak{s}}} \! := \! \alpha_{k} \! \neq \! \infty$, 
$\lbrace \mathstrut \alpha_{k^{\prime}}, \, k^{\prime} \! \in \! \lbrace 1,2,\dotsc,
K \rbrace; \, \alpha_{k^{\prime}} \! \neq \! \alpha_{k}, \, \alpha_{k} \! \neq \! \infty 
\rbrace \! = \! \varnothing$, then the corresponding Equations~\eqref{eq68} 
and~\eqref{eq69} are vacuous; this can only happen if $n \! = \! 1$.) Via the 
above ordered disjoint partition for the repeated pole sequence $\lbrace \alpha_{1},
\alpha_{2},\dotsc,\alpha_{K} \rbrace \cup \lbrace \alpha_{1},\alpha_{2},\dotsc,
\alpha_{K} \rbrace \cup \dotsb \cup \lbrace \alpha_{1},\alpha_{2},\dotsc,\alpha_{k} 
\rbrace$, and Equations~\eqref{eq66}--\eqref{eq69}, one writes, in the indicated order,
\begin{align*}
\dfrac{\mathcal{X}_{21}(z)}{z \! - \! \alpha_{k}} =& \, 
\sum_{q=1}^{\mathfrak{s}-2} \sum_{r=1}^{l_{q}} \dfrac{
\widehat{\nu}_{q,r}^{\raise-1.0ex\hbox{$\scriptstyle f$}}(n,k)}{(z \! 
- \! \alpha_{i(q)_{k_{q}}})^{r}} \! + \! \sum_{q=1}^{l_{\mathfrak{s}-1}} 
\widehat{\nu}_{\mathfrak{s}-1,q}^{\raise-1.0ex\hbox{$\scriptstyle f$}}
(n,k)z^{q} \! + \! \sum_{q=1}^{l_{\mathfrak{s}}} \dfrac{
\widehat{\nu}_{\mathfrak{s},q}^{\raise-1.0ex\hbox{$\scriptstyle f$}}
(n,k)}{(z \! - \! \alpha_{i(\mathfrak{s})_{k_{\mathfrak{s}}}})^{q-1}} \\
=& \, \sum_{m=1}^{\mathfrak{s}-2} \sum_{q=1}^{l_{m}=\varkappa_{nk 
\tilde{k}_{m}}} \dfrac{
\widehat{\nu}_{m,q}^{\raise-1.0ex\hbox{$\scriptstyle f$}}(n,k)}{(z \! - \! 
\alpha_{p_{m}})^{q}} \! + \! \sum_{q=1}^{l_{\mathfrak{s}-1}=\varkappa_{nk 
\tilde{k}_{\mathfrak{s}-1}}^{\infty}} 
\widehat{\nu}_{\mathfrak{s}-1,q}^{\raise-1.0ex\hbox{$\scriptstyle f$}}
(n,k)z^{q} \! + \! \sum_{q=1}^{l_{\mathfrak{s}}=\varkappa_{nk}} \dfrac{
\widehat{\nu}_{\mathfrak{s},q}^{\raise-1.0ex\hbox{$\scriptstyle f$}}
(n,k)}{(z \! - \! \alpha_{k})^{q-1}}, \quad \widehat{\nu}_{\mathfrak{s},
l_{\mathfrak{s}}}^{\raise-1.0ex\hbox{$\scriptstyle f$}}(n,k) \! \neq \! 0.
\end{align*}
Substituting the latter partial fraction expansion for $\mathcal{X}_{21}
(z)/(z \! - \! \alpha_{k})$ into Equations~\eqref{eq66}--\eqref{eq69}, 
one arrives at, for $n \! \in \! \mathbb{N}$ and $k \! \in \! \lbrace 1,2,
\dotsc,K \rbrace$ such that $\alpha_{p_{\mathfrak{s}}} \! := \! \alpha_{k} 
\! \neq \! \infty$, the orthogonality conditions
\begin{align}
\int_{\mathbb{R}} &\left(\sum_{m=1}^{\mathfrak{s}-2} 
\sum_{q=1}^{l_{m}=\varkappa_{nk \tilde{k}_{m}}} \dfrac{\widehat{\nu}_{m,q}^{
\raise-1.0ex\hbox{$\scriptstyle f$}}(n,k)}{(\xi \! - \! \alpha_{p_{m}})^{q}} 
\! + \! \sum_{q=1}^{l_{\mathfrak{s}-1}=\varkappa^{\infty}_{nk 
\tilde{k}_{\mathfrak{s}-1}}} \widehat{\nu}_{\mathfrak{s}-1,q}^{
\raise-1.0ex\hbox{$\scriptstyle f$}}(n,k) \xi^{q} \! + \! 
\sum_{q=1}^{l_{\mathfrak{s}}=\varkappa_{nk}} \dfrac{\widehat{\nu}_{
\mathfrak{s},q}^{\raise-1.0ex\hbox{$\scriptstyle f$}}(n,k)}{(\xi \! - 
\! \alpha_{k})^{q-1}} \right) \dfrac{\widetilde{w}(\xi)}{(\xi \! - \! 
\alpha_{k})^{r}} \, \md \xi \! = \! 0, \quad r \! = \! 0,1,\dotsc,
\varkappa_{nk} \! - \! 2, \label{eq70} \\
\int_{\mathbb{R}} &\left(\sum_{m=1}^{\mathfrak{s}-2} 
\sum_{q=1}^{l_{m}=\varkappa_{nk \tilde{k}_{m}}} \dfrac{\widehat{\nu}_{m,q}^{
\raise-1.0ex\hbox{$\scriptstyle f$}}(n,k)}{(\xi \! - \! \alpha_{p_{m}})^{q}} 
\! + \! \sum_{q=1}^{l_{\mathfrak{s}-1}=\varkappa^{\infty}_{nk 
\tilde{k}_{\mathfrak{s}-1}}} \widehat{\nu}_{\mathfrak{s}-1,q}^{
\raise-1.0ex\hbox{$\scriptstyle f$}}(n,k) \xi^{q} \! + \! 
\sum_{q=1}^{l_{\mathfrak{s}}=\varkappa_{nk}} \dfrac{\widehat{\nu}_{
\mathfrak{s},q}^{\raise-1.0ex\hbox{$\scriptstyle f$}}(n,k)}{(\xi \! - 
\! \alpha_{k})^{q-1}} \right) \dfrac{\widetilde{w}(\xi)}{(\xi \! - \! 
\alpha_{k})^{\varkappa_{nk}-1}} \, \md \xi \! = \! 2 \pi \mi, \label{eq71} \\
\int_{\mathbb{R}} &\left(\sum_{m=1}^{\mathfrak{s}-2} 
\sum_{q=1}^{l_{m}=\varkappa_{nk \tilde{k}_{m}}} \dfrac{\widehat{\nu}_{m,q}^{
\raise-1.0ex\hbox{$\scriptstyle f$}}(n,k)}{(\xi \! - \! \alpha_{p_{m}})^{q}} 
\! + \! \sum_{q=1}^{l_{\mathfrak{s}-1}=\varkappa^{\infty}_{nk 
\tilde{k}_{\mathfrak{s}-1}}} \widehat{\nu}_{\mathfrak{s}-1,q}^{
\raise-1.0ex\hbox{$\scriptstyle f$}}(n,k) \xi^{q} \! + \! 
\sum_{q=1}^{l_{\mathfrak{s}}=\varkappa_{nk}} \dfrac{\widehat{\nu}_{
\mathfrak{s},q}^{\raise-1.0ex\hbox{$\scriptstyle f$}}(n,k)}{(\xi \! - 
\! \alpha_{k})^{q-1}} \right) \dfrac{\widetilde{w}(\xi)}{(\xi \! - \! 
\alpha_{p_{i}})^{j}} \, \md \xi \! = \! 0, \nonumber \\
& \, \quad i \! = \! 1,2,\dotsc,\mathfrak{s} \! - \! 2, \quad j \! = \! 
1,2,\dotsc,l_{i}, \label{eq72} \\
\int_{\mathbb{R}} &\left(\sum_{m=1}^{\mathfrak{s}-2} 
\sum_{q=1}^{l_{m}=\varkappa_{nk \tilde{k}_{m}}} \dfrac{\widehat{\nu}_{m,q}^{
\raise-1.0ex\hbox{$\scriptstyle f$}}(n,k)}{(\xi \! - \! \alpha_{p_{m}})^{q}} 
\! + \! \sum_{q=1}^{l_{\mathfrak{s}-1}=\varkappa^{\infty}_{nk 
\tilde{k}_{\mathfrak{s}-1}}} \widehat{\nu}_{\mathfrak{s}-1,q}^{
\raise-1.0ex\hbox{$\scriptstyle f$}}(n,k) \xi^{q} \! + \! 
\sum_{q=1}^{l_{\mathfrak{s}}=\varkappa_{nk}} \dfrac{\widehat{\nu}_{
\mathfrak{s},q}^{\raise-1.0ex\hbox{$\scriptstyle f$}}(n,k)}{(\xi \! - 
\! \alpha_{k})^{q-1}} \right) \xi^{r} \widetilde{w}(\xi) \, \md \xi \! 
= \! 0, \quad r \! = \! 1,2,\dotsc,l_{\mathfrak{s}-1}. \label{eq73}
\end{align}
Incidentally, for $n \! \in \! \mathbb{N}$ and $k \! \in \! \lbrace 
1,2,\dotsc,K \rbrace$ such that $\alpha_{p_{\mathfrak{s}}} 
\! := \! \alpha_{k} \! \neq \! \infty$, via the orthogonality 
conditions~\eqref{eq70}, \eqref{eq72}, and~\eqref{eq73}, the 
orthogonality condition~\eqref{eq71} can be manipulated thus:
\begin{align*}
2 \pi \mi =& \, \int_{\mathbb{R}} \left(\sum_{m=1}^{\mathfrak{s}-2} 
\sum_{q=1}^{l_{m}=\varkappa_{nk \tilde{k}_{m}}} \dfrac{\widehat{\nu}_{m,q}^{
\raise-1.0ex\hbox{$\scriptstyle f$}}(n,k)}{(\xi \! - \! \alpha_{p_{m}})^{q}} 
\! + \! \sum_{q=1}^{l_{\mathfrak{s}-1}=\varkappa^{\infty}_{nk \tilde{k}_{
\mathfrak{s}-1}}} \widehat{\nu}_{\mathfrak{s}-1,q}^{
\raise-1.0ex\hbox{$\scriptstyle f$}}(n,k) \xi^{q} \! + \! \sum_{q=1}^{
l_{\mathfrak{s}}=\varkappa_{nk}} \dfrac{\widehat{\nu}_{\mathfrak{s},q}^{
\raise-1.0ex\hbox{$\scriptstyle f$}}(n,k)}{(\xi \! - \! \alpha_{k})^{q-1}} 
\right) \dfrac{\widehat{\nu}_{\mathfrak{s},l_{\mathfrak{s}}}^{
\raise-1.0ex\hbox{$\scriptstyle f$}}(n,k)}{(\xi \! - \! \alpha_{k})^{
\varkappa_{nk}-1}} \\
\times& \, \dfrac{\widetilde{w}(\xi)}{\widehat{\nu}_{\mathfrak{s},
l_{\mathfrak{s}}}^{\raise-1.0ex\hbox{$\scriptstyle f$}}(n,k)} \, \md 
\xi \! = \! \int_{\mathbb{R}} \left(\sum_{m=1}^{\mathfrak{s}-2} 
\sum_{q=1}^{l_{m}=\varkappa_{nk \tilde{k}_{m}}} \dfrac{\widehat{\nu}_{m,q}^{
\raise-1.0ex\hbox{$\scriptstyle f$}}(n,k)}{(\xi \! - \! \alpha_{p_{m}})^{q}} 
\! + \! \sum_{q=1}^{l_{\mathfrak{s}-1}=\varkappa^{\infty}_{nk \tilde{k}_{
\mathfrak{s}-1}}} \widehat{\nu}_{\mathfrak{s}-1,q}^{
\raise-1.0ex\hbox{$\scriptstyle f$}}(n,k) \xi^{q} \! + \! \sum_{q=1}^{
l_{\mathfrak{s}}=\varkappa_{nk}} \dfrac{\widehat{\nu}_{\mathfrak{s},q}^{
\raise-1.0ex\hbox{$\scriptstyle f$}}(n,k)}{(\xi \! - \! \alpha_{k})^{q-1}} 
\right) \\
\times& \, \left(\sum_{m=1}^{\mathfrak{s}-2} \sum_{q=1}^{l_{m}=
\varkappa_{nk \tilde{k}_{m}}} \dfrac{\widehat{\nu}_{m,q}^{
\raise-1.0ex\hbox{$\scriptstyle f$}}(n,k)}{(\xi \! - \! \alpha_{p_{m}})^{q}} 
\! + \! \sum_{q=1}^{l_{\mathfrak{s}-1}=\varkappa^{\infty}_{nk \tilde{k}_{
\mathfrak{s}-1}}} \widehat{\nu}_{\mathfrak{s}-1,q}^{
\raise-1.0ex\hbox{$\scriptstyle f$}}(n,k) \xi^{q} \! + \! \sum_{q=1}^{
l_{\mathfrak{s}}-1=\varkappa_{nk}-1} \dfrac{\widehat{\nu}_{\mathfrak{s},q}^{
\raise-1.0ex\hbox{$\scriptstyle f$}}(n,k)}{(\xi \! - \! \alpha_{k})^{q-1}} 
\! + \! \dfrac{\widehat{\nu}_{\mathfrak{s},l_{\mathfrak{s}}}^{
\raise-1.0ex\hbox{$\scriptstyle f$}}(n,k)}{(\xi \! - \! \alpha_{k})^{
\varkappa_{nk}-1}} \right) \\
\times& \, \dfrac{\widetilde{w}(\xi)}{\widehat{\nu}_{\mathfrak{s},
l_{\mathfrak{s}}}^{\raise-1.0ex\hbox{$\scriptstyle f$}}(n,k)} \, \md \xi 
\! = \! \int_{\mathbb{R}} \underbrace{\left(\sum_{m=1}^{\mathfrak{s}-2} 
\sum_{q=1}^{l_{m}=\varkappa_{nk \tilde{k}_{m}}} \dfrac{\widehat{\nu}_{m,q}^{
\raise-1.0ex\hbox{$\scriptstyle f$}}(n,k)}{(\xi \! - \! \alpha_{p_{m}})^{q}} 
\! + \! \sum_{q=1}^{l_{\mathfrak{s}-1}=\varkappa^{\infty}_{nk \tilde{k}_{
\mathfrak{s}-1}}} \widehat{\nu}_{\mathfrak{s}-1,q}^{
\raise-1.0ex\hbox{$\scriptstyle f$}}(n,k) \xi^{q} \! + \! \sum_{q=1}^{
l_{\mathfrak{s}}=\varkappa_{nk}} \dfrac{\widehat{\nu}_{\mathfrak{s},q}^{
\raise-1.0ex\hbox{$\scriptstyle f$}}(n,k)}{(\xi \! - \! \alpha_{k})^{q-1}} 
\right)}_{= \, \mathcal{X}_{21}(\xi)/(\xi -\alpha_{k})} \\
\times& \, \underbrace{\left(\sum_{m=1}^{\mathfrak{s}-2} \sum_{q=1}^{
l_{m}=\varkappa_{nk \tilde{k}_{m}}} \dfrac{\widehat{\nu}_{m,q}^{
\raise-1.0ex\hbox{$\scriptstyle f$}}(n,k)}{(\xi \! - \! \alpha_{p_{m}})^{q}} 
\! + \! \sum_{q=1}^{l_{\mathfrak{s}-1}=\varkappa^{\infty}_{nk \tilde{k}_{
\mathfrak{s}-1}}} \widehat{\nu}_{\mathfrak{s}-1,q}^{
\raise-1.0ex\hbox{$\scriptstyle f$}}(n,k) \xi^{q} \! + \! \sum_{q=1}^{
l_{\mathfrak{s}}=\varkappa_{nk}} \dfrac{\widehat{\nu}_{\mathfrak{s},q}^{
\raise-1.0ex\hbox{$\scriptstyle f$}}(n,k)}{(\xi \! - \! \alpha_{k})^{q-1}} 
\right)}_{= \, \mathcal{X}_{21}(\xi)/(\xi -\alpha_{k})} \, 
\dfrac{\widetilde{w}(\xi)}{\widehat{\nu}_{\mathfrak{s},l_{\mathfrak{s}}}^{
\raise-1.0ex\hbox{$\scriptstyle f$}}(n,k)} \, \md \xi \\
=& \, \int_{\mathbb{R}} \left(\dfrac{\mathcal{X}_{21}(\xi)}{\xi \! - \! 
\alpha_{k}} \right)^{2} \dfrac{\widetilde{w}(\xi)}{\widehat{\nu}_{
\mathfrak{s},l_{\mathfrak{s}}}^{\raise-1.0ex\hbox{$\scriptstyle f$}}(n,k)} 
\, \md \xi;
\end{align*}
hence, for $n \! \in \! \mathbb{N}$ and $k \! \in \! \lbrace 1,2,\dotsc,
K \rbrace$ such that $\alpha_{p_{\mathfrak{s}}} \! := \! \alpha_{k} 
\! \neq \! \infty$, one arrives at the `normalisation formula'
\begin{equation*}
\int_{\mathbb{R}} \left(\dfrac{\mathcal{X}_{21}(\xi)}{\xi \! - \! \alpha_{k}} 
\right)^{2} \widetilde{w}(\xi) \, \md \xi \! = \! 2 \pi \mi \, \widehat{\nu}_{
\mathfrak{s},\varkappa_{nk}}^{\raise-1.0ex\hbox{$\scriptstyle f$}}(n,k).
\end{equation*}
For $n \! \in \! \mathbb{N}$ and $k \! \in \! \lbrace 1,2,\dotsc,K 
\rbrace$ such that $\alpha_{p_{\mathfrak{s}}} \! := \! \alpha_{k} \! \neq 
\! \infty$, the orthogonality conditions~\eqref{eq70}--\eqref{eq73} give 
rise to a total of (cf. Equation~\eqref{fincount}) $\sum_{r=1}^{\mathfrak{s}}
l_{r} \! = \! \sum_{r=1}^{\mathfrak{s}-2}l_{r} \! + \! l_{\mathfrak{s}-1} \! 
+ \! l_{\mathfrak{s}} \! = \! \sum_{r=1}^{\mathfrak{s}-2} \varkappa_{nk 
\tilde{k}_{r}} \! + \! \varkappa_{nk \tilde{k}_{\mathfrak{s}-1}}^{\infty} \! 
+ \! \varkappa_{nk} \! = \! (n \! - \! 1)K \! + \! k$ linear inhomogeneous 
algebraic equations for the $(n \! - \! 1)K \! + \! k$ unknowns 
$\widehat{\nu}_{1,1}^{\raise-1.0ex\hbox{$\scriptstyle f$}}(n,k),\dotsc,
\widehat{\nu}_{1,l_{1}}^{\raise-1.0ex\hbox{$\scriptstyle f$}}(n,k),\dotsc,
\widehat{\nu}_{\mathfrak{s}-1,1}^{\raise-1.0ex\hbox{$\scriptstyle f$}}(n,k),
\dotsc,\widehat{\nu}_{\mathfrak{s}-1,l_{\mathfrak{s}-1}}^{
\raise-1.0ex\hbox{$\scriptstyle f$}}(n,k),\widehat{\nu}_{\mathfrak{s},1}^{
\raise-1.0ex\hbox{$\scriptstyle f$}}(n,k),\dotsc,\widehat{\nu}_{\mathfrak{s},
l_{\mathfrak{s}}}^{\raise-1.0ex\hbox{$\scriptstyle f$}}(n,k)$, that is,
\begin{align}
\setcounter{MaxMatrixCols}{12}
&\left(
, \label{eq74}
\end{align}
where (with abuse of notation)
\begin{equation*}
n_{1} \! = \! l_{1} \! + \! \dotsb \! + \! l_{\mathfrak{s}-2} \! + \! 1, 
\qquad n_{2} \! = \! (n \! - \! 1)K \! + \! k \! - \! \varkappa_{nk}, 
\qquad \text{and} \qquad m_{1} \! = \! (n \! - \! 1)K \! + \! k \! - \! 1.
\end{equation*}
For $n \! \in \! \mathbb{N}$ and $k \! \in \! \lbrace 1,2,\dotsc,K 
\rbrace$ such that $\alpha_{p_{\mathfrak{s}}} \! := \! \alpha_{k} 
\! \neq \! \infty$, the linear system~\eqref{eq74} of $(n \! - 
\! 1)K \! + \! k$ inhomogeneous algebraic equations for the 
$(n \! - \! 1)K \! + \! k$ unknowns 
$\widehat{\nu}_{1,1}^{\raise-1.0ex\hbox{$\scriptstyle f$}}(n,k),\dotsc,
\widehat{\nu}_{1,l_{1}}^{\raise-1.0ex\hbox{$\scriptstyle f$}}(n,k),\dotsc,
\widehat{\nu}_{\mathfrak{s}-1,1}^{\raise-1.0ex\hbox{$\scriptstyle f$}}
(n,k),\dotsc,\widehat{\nu}_{\mathfrak{s}-1,l_{\mathfrak{s}-1}}^{
\raise-1.0ex\hbox{$\scriptstyle f$}}(n,k),
\widehat{\nu}_{\mathfrak{s},1}^{\raise-1.0ex\hbox{$\scriptstyle f$}}
(n,k),\linebreak[4] 
\dotsc,
\widehat{\nu}_{\mathfrak{s},l_{\mathfrak{s}}}^{\raise-1.0ex\hbox{$\scriptstyle f$}}
(n,k)$ admits a unique solution if, and only if, the determinant of the 
coefficient matrix, denoted $\mathcal{N}_{f}^{\sharp}(n,k)$, is non-zero: 
this fact will now be established. Via the multi-linearity property of the 
determinant, one shows that, for $n \! \in \! \mathbb{N}$ and $k \! \in 
\! \lbrace 1,2,\dotsc,K \rbrace$ such that $\alpha_{p_{\mathfrak{s}}} 
\! := \! \alpha_{k} \! \neq \! \infty$,
\begin{align*}
\mathcal{N}^{\sharp}_{f}(n,k) =& \, \underbrace{\int_{\mathbb{R}} 
\int_{\mathbb{R}} \dotsb \int_{\mathbb{R}}}_{(n-1)K+k} \md \widetilde{\mu}
(\xi_{0}) \, \md \widetilde{\mu}(\xi_{1}) \, \dotsb \, \md \widetilde{\mu}(\xi_{l_{1}}) 
\, \dotsb \, \md \widetilde{\mu}(\xi_{n_{1}}) \, \dotsb \, \md \widetilde{\mu}
(\xi_{n_{2}}) \, \md \widetilde{\mu}(\xi_{n_{2}+1}) \, \dotsb \, \md \widetilde{\mu}
(\xi_{m_{1}}) \\
\times& \, \dfrac{1}{(\xi_{0} \! - \! \alpha_{k})^{0}(\xi_{1} \! - \! \alpha_{p_{1}})^{1} 
\dotsb (\xi_{l_{1}} \! - \! \alpha_{p_{1}})^{l_{1}} \dotsb \dfrac{1}{(\xi_{n_{1}})^{1}} 
\dotsb \dfrac{1}{(\xi_{n_{2}})^{l_{\mathfrak{s}-1}}}(\xi_{n_{2}+1} \! - \! \alpha_{k})^{1} 
\dotsb (\xi_{m_{1}} \! - \! \alpha_{k})^{\varkappa_{nk}-1}} \\
\times& \, 
\left\vert

\right\vert.
\end{align*}
Recalling, for $n \! \in \! \mathbb{N}$ and $k \! \in \! \lbrace 1,2,\dotsc,K 
\rbrace$ such that $\alpha_{p_{\mathfrak{s}}} \! := \! \alpha_{k} \! \neq 
\! \infty$, the $(n \! - \! 1)K \! + \! k$ linearly independent functions 
$\tilde{\phi}_{0}(z),\tilde{\phi}_{1}(z),\dotsc,\tilde{\phi}_{l_{1}}(z),
\dotsc,\tilde{\phi}_{n_{1}}(z),\dotsc,\tilde{\phi}_{n_{2}}(z),\tilde{\phi}_{n_{2}+1}
(z),\dotsc,\tilde{\phi}_{m_{1}}(z)$ introduced above for the analysis of 
$\hat{\mathfrak{c}}_{D_{f}}$, one shows that
\begin{align*}
\mathcal{N}_{f}^{\sharp}(n,k) =& \, \underbrace{\int_{\mathbb{R}} 
\int_{\mathbb{R}} \dotsb \int_{\mathbb{R}}}_{(n-1)K+k} \md \widetilde{\mu}
(\xi_{0}) \, \md \widetilde{\mu}(\xi_{1}) \, \dotsb \, \md \widetilde{\mu}(\xi_{l_{1}}) 
\, \dotsb \, \md \widetilde{\mu}(\xi_{n_{1}}) \, \dotsb \, \md \widetilde{\mu}
(\xi_{n_{2}}) \, \md \widetilde{\mu}(\xi_{n_{2}+1}) \, \dotsb \, \md \widetilde{\mu}
(\xi_{m_{1}}) \\
\times& \, \dfrac{(-1)^{(n-1)K+k-\varkappa_{nk}}(\tilde{\phi}_{0}(\xi_{0}) 
\tilde{\phi}_{0}(\xi_{1}) \dotsb \tilde{\phi}_{0}(\xi_{l_{1}}) \dotsb \tilde{
\phi}_{0}(\xi_{n_{1}}) \dotsb \tilde{\phi}_{0}(\xi_{n_{2}}) \tilde{\phi}_{0}
(\xi_{n_{2}+1}) \dotsb \tilde{\phi}_{0}(\xi_{m_{1}}))^{-1}}{(\xi_{0} \! - 
\! \alpha_{k})^{0}(\xi_{1} \! - \! \alpha_{p_{1}})^{1} \dotsb (\xi_{l_{1}} 
\! - \! \alpha_{p_{1}})^{l_{1}} \dotsb \dfrac{1}{(\xi_{n_{1}})^{1}} 
\dotsb \dfrac{1}{(\xi_{n_{2}})^{l_{\mathfrak{s}-1}}}(\xi_{n_{2}+1} \! - \! 
\alpha_{k})^{1} \dotsb (\xi_{m_{1}} \! - \! \alpha_{k})^{\varkappa_{nk}-1}} \\
\times& \, 
\left\vert

\right\vert.
\end{align*}
Modulo the undulatory factor $(-1)^{(n-1)K+k-\varkappa_{nk}}$, this 
is the same determinantal expression encountered while studying 
$\hat{\mathfrak{c}}_{D_{f}}$; therefore, for $n \! \in \! \mathbb{N}$ 
and $k \! \in \! \lbrace 1,2,\dotsc,K \rbrace$ such that 
$\alpha_{p_{\mathfrak{s}}} \! := \! \alpha_{k} \! \neq \! \infty$, 
it follows that
\begin{align}
\mathcal{N}_{f}^{\sharp}(n,k) =& \, \dfrac{(-1)^{(n-1)K+k-\varkappa_{nk}}
(\mathbb{D}^{\blacklozenge})^{2}}{((n \! - \! 1)K \! + \! k)!} \underbrace{
\int_{\mathbb{R}} \int_{\mathbb{R}} \dotsb \int_{\mathbb{R}}}_{(n-1)K+k} 
\md \widetilde{\mu}(\xi_{0}) \, \md \widetilde{\mu}(\xi_{1}) \, \dotsb \, 
\md \widetilde{\mu}(\xi_{(n-1)K+k-1}) 
\nonumber \\
\times& \, \prod_{\substack{i,j=0\\j<i}}^{(n-1)K+k-1}(\xi_{i} \! - \! 
\xi_{j})^{2} \left(\prod_{m=0}^{(n-1)K+k-1} \prod_{q=1}^{\mathfrak{s}-2}
(\xi_{m} \! - \! \alpha_{p_{q}})^{\varkappa_{nk \tilde{k}_{q}}}(\xi_{m} \! - \! 
\alpha_{k})^{\varkappa_{nk}-1} \right)^{-2} \quad (\neq \! 0), 
\label{eq75}
\end{align}
whence follows, for $n \! \in \! \mathbb{N}$ and $k \! \in \! \lbrace 1,2,
\dotsc,K \rbrace$ such that $\alpha_{p_{\mathfrak{s}}} \! := \! \alpha_{k} 
\! \neq \! \infty$, the existence and the uniqueness of $\mathcal{X}_{21}
(z)/(z \! - \! \alpha_{k})$. \hfill $\qed$
\begin{eeeee} \label{remdets}
\textsl{A calculation shows that (cf. Equations~\eqref{eq36} and~\eqref{eq38}$)$, 
for $(n,k) \! \in \! \mathbb{N} \times \lbrace 1,2,\dotsc,K \rbrace$ such that 
$\alpha_{p_{\mathfrak{s}}} \! := \! \alpha_{k} \! = \! \infty$, $\mathbb{D}^{
\spcheck} \! = \! \mathbb{D}^{\sphat}$, whilst (cf. Equations~\eqref{eq62} 
and~\eqref{eq64}$)$ for $(n,k) \! \in \! \mathbb{N} \times \lbrace 1,2,\dotsc,K 
\rbrace$ such that $\alpha_{p_{\mathfrak{s}}} \! := \! \alpha_{k} \! \neq \! \infty$, 
$\mathbb{D}^{\lozenge}/\mathbb{D}^{\blacklozenge} \! = \! (-1)^{(n-1)K+k} 
\prod_{q=1}^{\mathfrak{s}-2}(\alpha_{k} \! - \! \alpha_{p_{q}})^{\varkappa_{nk 
\tilde{k}_{q}}}$.}
\end{eeeee}
\begin{ccccc} \label{lem2.2} 
Let $\widetilde{V} \colon \overline{\mathbb{R}} \setminus \lbrace \alpha_{1},
\alpha_{2},\dotsc,\alpha_{K} \rbrace \! \to \! \mathbb{R}$ satisfy 
conditions~\eqref{eq20}--\eqref{eq22}. Let $\mathcal{X} \colon \mathbb{N} 
\times \lbrace 1,2,\dotsc,K \rbrace \times \overline{\mathbb{C}} 
\setminus \overline{\mathbb{R}} \! \to \! \operatorname{SL}_{2}
(\mathbb{C})$ be the unique solution of the {\rm RHP} stated in 
Lemma~$\bm{\mathrm{RHP}_{\mathrm{MPC}}}$ with integral representation 
given by Equation~\eqref{intrepinf} for $\alpha_{p_{\mathfrak{s}}} \! 
:= \! \alpha_{k} \! = \! \infty$ (resp., Equation~\eqref{intrepfin} for 
$\alpha_{p_{\mathfrak{s}}} \! := \! \alpha_{k} \! \neq \! \infty)$, 
$k \! \in \! \lbrace 1,2,\dotsc,K \rbrace$. For $n \! \in \! \mathbb{N}$ and 
$k \! \in \! \lbrace 1,2,\dotsc,K \rbrace$ such that $\alpha_{p_{\mathfrak{s}}} 
\! := \! \alpha_{k} \! = \! \infty$ (resp., $\alpha_{p_{\mathfrak{s}}} \! := \! 
\alpha_{k} \! \neq \! \infty)$, let $\pmb{\pi}^{n}_{k}(z)$ be the corresponding 
monic {\rm MPC ORF} defined in Subsection~\ref{subsubsec1.2.1} (resp., 
Subsection~\ref{subsubsec1.2.2}$)$. Then$:$ {\rm (i)} for $n \! \in 
\! \mathbb{N}$ and $k \! \in \! \lbrace 1,2,\dotsc,K \rbrace$ such 
that $\alpha_{p_{\mathfrak{s}}} \! := \! \alpha_{k} \! = \! \infty$, 
$\pmb{\pi}^{n}_{k}(z)$ and $\mathcal{X}_{21}(z)$ have, respectively, 
the following integral representations,
\begin{equation} \label{eintreppinf1} 
\pmb{\pi}^{n}_{k}(z) \! = \! \mathcal{X}_{11}(z) \! = \! 
\dfrac{\mathbb{D}^{\spcheck}}{\mathbb{D}^{\sphat}} \left(
\prod_{q=1}^{\mathfrak{s}-1}(z \! - \! \alpha_{p_{q}})^{\varkappa_{nk 
\tilde{k}_{q}}} \right)^{-1} \dfrac{\Upsilon_{N_{\infty}}(n,k;z)}{
\Upsilon_{D_{\infty}}(n,k)},
\end{equation}
where $\mathbb{D}^{\spcheck}$ and $\mathbb{D}^{\sphat}$, with 
$\mathbb{D}^{\spcheck} \! = \! \mathbb{D}^{\sphat}$, are defined by 
Equations~\eqref{eq36} and~\eqref{eq38}, respectively,
\begin{align*}
\Upsilon_{N_{\infty}}(n,k;z) =& \, \underbrace{\int_{\mathbb{R}} 
\int_{\mathbb{R}} \dotsb \int_{\mathbb{R}}}_{(n-1)K+k} \md \widetilde{\mu}
(\xi_{0}) \, \md \widetilde{\mu}(\xi_{1}) \, \dotsb \, \md \widetilde{\mu}
(\xi_{(n-1)K+k-1}) \prod_{\substack{i,j=0\\j<i}}^{(n-1)K+k-1}(\xi_{i} \! - \! 
\xi_{j})^{2} \\
\times& \, \left(\prod_{m=0}^{(n-1)K+k-1} \prod_{q=1}^{\mathfrak{s}-1}
(\xi_{m} \! - \! \alpha_{p_{q}})^{\varkappa_{nk \tilde{k}_{q}}} \right)^{-2} 
\prod_{l=0}^{(n-1)K+k-1}(z \! - \! \xi_{l}), \\
\Upsilon_{D_{\infty}}(n,k) =& \, \underbrace{\int_{\mathbb{R}} 
\int_{\mathbb{R}} \dotsb \int_{\mathbb{R}}}_{(n-1)K+k} \md \widetilde{\mu}
(\tau_{0}) \, \md \widetilde{\mu}(\tau_{1}) \, \dotsb \, \md \widetilde{\mu}
(\tau_{(n-1)K+k-1}) \prod_{\substack{i,j=0\\j<i}}^{(n-1)K+k-1}(\tau_{i} \! - \! 
\tau_{j})^{2} \\
\times& \, \left(\prod_{m=0}^{(n-1)K+k-1} \prod_{q=1}^{\mathfrak{s}-1}
(\tau_{m} \! - \! \alpha_{p_{q}})^{\varkappa_{nk \tilde{k}_{q}}} \right)^{-2},
\end{align*}
with $\md \widetilde{\mu}(z) \! = \! \exp (-n \widetilde{V}(z)) \, \md z$, and
\begin{equation} \label{eintreppinf2} 
\mathcal{X}_{21}(z) \! = \! -2 \pi \mi ((n \! - \! 1)K \! + \! k) 
\dfrac{\mathbb{D}^{\spadesuit}}{\mathbb{D}^{\sphat}} \left(
\prod_{q=1}^{\mathfrak{s}-1}(z \! - \! \alpha_{p_{q}})^{\varkappa_{nk 
\tilde{k}_{q}}} \right)^{-1} \dfrac{\Lambda_{N_{\infty}}(n,k;z)}{
\Lambda_{D_{\infty}}(n,k)},
\end{equation}
where $\mathbb{D}^{\spadesuit}$ is defined by Equation~\eqref{eq76} 
(see below),
\begin{align*}
\Lambda_{N_{\infty}}(n,k;z) =& \, \underbrace{\int_{\mathbb{R}} 
\int_{\mathbb{R}} \dotsb \int_{\mathbb{R}}}_{(n-1)K+k-1} \md \widetilde{\mu}
(\xi_{0}) \, \md \widetilde{\mu}(\xi_{1}) \, \dotsb \, \md \widetilde{\mu}
(\xi_{(n-1)K+k-2}) \prod_{\substack{i,j=0\\j<i}}^{(n-1)K+k-2}(\xi_{i} \! - \! 
\xi_{j})^{2} \\
\times& \, \left(\prod_{m=0}^{(n-1)K+k-2} \prod_{q=1}^{\mathfrak{s}-1}
(\xi_{m} \! - \! \alpha_{p_{q}})^{\varkappa_{nk \tilde{k}_{q}}} \right)^{-2} 
\prod_{l=0}^{(n-1)K+k-2}(z \! - \! \xi_{l}), \\
\Lambda_{D_{\infty}}(n,k) =& \, \underbrace{\int_{\mathbb{R}} 
\int_{\mathbb{R}} \dotsb \int_{\mathbb{R}}}_{(n-1)K+k} \md \widetilde{\mu}
(\tau_{0}) \, \md \widetilde{\mu}(\tau_{1}) \, \dotsb \, \md \widetilde{\mu}
(\tau_{(n-1)K+k-1}) \prod_{\substack{i,j=0\\j<i}}^{(n-1)K+k-1}(\tau_{i} \! 
- \! \tau_{j})^{2} \\
\times& \, \left(\prod_{m=0}^{(n-1)K+k-1} \prod_{q=1}^{\mathfrak{s}-1}
(\tau_{m} \! - \! \alpha_{p_{q}})^{\varkappa_{nk \tilde{k}_{q}}} 
\right)^{-2};
\end{align*}
and {\rm (ii)} for $n \! \in \! \mathbb{N}$ and $k \! \in \! \lbrace 1,2,
\dotsc,K \rbrace$ such that $\alpha_{p_{\mathfrak{s}}} \! := \! \alpha_{k} 
\! \neq \! \infty$, $\pmb{\pi}^{n}_{k}(z)$ and $\mathcal{X}_{21}(z)$ 
have, respectively, the following integral representations,
\begin{equation} \label{eintreppfin1} 
\pmb{\pi}^{n}_{k}(z) \! = \! \dfrac{\mathcal{X}_{11}(z)}{z \! - \! \alpha_{k}} 
\! = \! \dfrac{\mathbb{D}^{\lozenge}}{\mathbb{D}^{\blacklozenge}} \left(
\prod_{q=1}^{\mathfrak{s}-2}(z \! - \! \alpha_{p_{q}})^{\varkappa_{nk 
\tilde{k}_{q}}}(z \! - \! \alpha_{k})^{\varkappa_{nk}} \right)^{-1} 
\dfrac{\Upsilon_{N_{f}}(n,k;z)}{\Upsilon_{D_{f}}(n,k)},
\end{equation}
where $\mathbb{D}^{\lozenge}$ and $\mathbb{D}^{\blacklozenge}$, with 
$\mathbb{D}^{\lozenge}/\mathbb{D}^{\blacklozenge} \! = \! (-1)^{(n-1)K+k} 
\prod_{q=1}^{\mathfrak{s}-2}(\alpha_{k} \! - \! \alpha_{p_{q}})^{\varkappa_{nk 
\tilde{k}_{q}}}$, are defined by Equations~\eqref{eq62} and~\eqref{eq64}, 
respectively,
\begin{align*}
\Upsilon_{N_{f}}(n,k;z) =& \, \underbrace{\int_{\mathbb{R}} \int_{\mathbb{R}} 
\dotsb \int_{\mathbb{R}}}_{(n-1)K+k} \md \widetilde{\mu}(\xi_{0}) \, \md 
\widetilde{\mu}(\xi_{1}) \, \dotsb \, \md \widetilde{\mu}(\xi_{(n-1)K+k-1}) 
\prod_{\substack{i,j=0\\j<i}}^{(n-1)K+k-1}(\xi_{i} \! - \! \xi_{j})^{2} \\
\times& \, \left(\prod_{m=0}^{(n-1)K+k-1} \prod_{q=1}^{\mathfrak{s}-2}
(\xi_{m} \! - \! \alpha_{p_{q}})^{\varkappa_{nk \tilde{k}_{q}}}(\xi_{m} \! - 
\! \alpha_{k})^{\varkappa_{nk}} \right)^{-2} \prod_{l=0}^{(n-1)K+k-1}
(\xi_{l} \! - \! \alpha_{k})(z \! - \! \xi_{l}), \\
\Upsilon_{D_{f}}(n,k) =& \, \underbrace{\int_{\mathbb{R}} \int_{\mathbb{R}} 
\dotsb \int_{\mathbb{R}}}_{(n-1)K+k} \md \widetilde{\mu}(\tau_{0}) \, \md 
\widetilde{\mu}(\tau_{1}) \, \dotsb \, \md \widetilde{\mu}(\tau_{(n-1)K+k-1}) 
\prod_{\substack{i,j=0\\j<i}}^{(n-1)K+k-1}(\tau_{i} \! - \! \tau_{j})^{2} \\
\times& \, \left(\prod_{m=0}^{(n-1)K+k-1} \prod_{q=1}^{\mathfrak{s}-2}
(\tau_{m} \! - \! \alpha_{p_{q}})^{\varkappa_{nk \tilde{k}_{q}}}(\tau_{m} \! 
- \! \alpha_{k})^{\varkappa_{nk}-1} \right)^{-2},
\end{align*}
and
\begin{equation} \label{eintreppfin2} 
\dfrac{\mathcal{X}_{21}(z)}{z \! - \! \alpha_{k}} \! = \! 2 \pi \mi ((n \! - 
\! 1)K \! + \! k) \dfrac{\mathbb{D}^{\clubsuit}}{\mathbb{D}^{\blacklozenge}} 
\left(\prod_{q=1}^{\mathfrak{s}-2}(z \! - \! \alpha_{p_{q}})^{\varkappa_{nk 
\tilde{k}_{q}}}(z \! - \! \alpha_{k})^{\varkappa_{nk}-1} \right)^{-1} 
\dfrac{\Lambda_{N_{f}}(n,k;z)}{\Lambda_{D_{f}}(n,k)},
\end{equation}
where $\mathbb{D}^{\clubsuit}$ is defined by Equation~\eqref{eq77} 
(see below),
\begin{align*}
\Lambda_{N_{f}}(n,k;z) =& \, \underbrace{\int_{\mathbb{R}} 
\int_{\mathbb{R}} \dotsb \int_{\mathbb{R}}}_{(n-1)K+k-1} \md \widetilde{\mu}
(\xi_{0}) \, \md \widetilde{\mu}(\xi_{1}) \, \dotsb \, \md \widetilde{\mu}
(\xi_{(n-1)K+k-2}) \prod_{\substack{i,j=0\\j<i}}^{(n-1)K+k-2}
(\xi_{i} \! - \! \xi_{j})^{2} \\
\times& \, \left(\prod_{m=0}^{(n-1)K+k-2} \prod_{q=1}^{\mathfrak{s}-2}
(\xi_{m} \! - \! \alpha_{p_{q}})^{\varkappa_{nk \tilde{k}_{q}}}
(\xi_{m} \! - \! \alpha_{k})^{\varkappa_{nk}-1} \right)^{-2} 
\prod_{l=0}^{(n-1)K+k-2}(\xi_{l} \! - \! \alpha_{k})(z \! - \! \xi_{l}), \\
\Lambda_{D_{f}}(n,k) =& \, \underbrace{\int_{\mathbb{R}} \int_{\mathbb{R}} 
\dotsb \int_{\mathbb{R}}}_{(n-1)K+k} \md \widetilde{\mu}(\tau_{0}) \, \md 
\widetilde{\mu}(\tau_{1}) \, \dotsb \, \md \widetilde{\mu}(\tau_{(n-1)K+k-1}) 
\prod_{\substack{i,j=0\\j<i}}^{(n-1)K+k-1}(\tau_{i} \! - \! \tau_{j})^{2} \\
\times& \, \left(\prod_{m=0}^{(n-1)K+k-1} \prod_{q=1}^{\mathfrak{s}-2}
(\tau_{m} \! - \! \alpha_{p_{q}})^{\varkappa_{nk \tilde{k}_{q}}}
(\tau_{m} \! - \! \alpha_{k})^{\varkappa_{nk}-1} \right)^{-2}.
\end{align*}
\end{ccccc}

\emph{Proof}. The proof consists of two cases: (i) $n \! \in \! 
\mathbb{N}$ and $k \! \in \! \lbrace 1,2,\dotsc,K \rbrace$ such that 
$\alpha_{p_{\mathfrak{s}}} \! := \! \alpha_{k} \! = \! \infty$ (see case 
\pmb{(1)} below); and (ii) $n \! \in \! \mathbb{N}$ and $k \! \in \! 
\lbrace 1,2,\dotsc,K \rbrace$ such that $\alpha_{p_{\mathfrak{s}}} 
\! := \! \alpha_{k} \! \neq \! \infty$ (see case \pmb{(2)} below).

\pmb{(1)} For $n \! \in \! \mathbb{N}$ and $k \! \in \! \lbrace 1,2,
\dotsc,K \rbrace$ such that $\alpha_{p_{\mathfrak{s}}} \! := \! 
\alpha_{k} \! = \! \infty$, recall {}from the proof of Lemma~\ref{lem2.1} 
(cf. case~\textbf{(1)}) the ordered disjoint partition 
for $\lbrace \alpha_{1},\alpha_{2},\dotsc,\alpha_{K} \rbrace \cup \lbrace 
\alpha_{1},\alpha_{2},\dotsc,\alpha_{K} \rbrace \cup \dotsb \cup \lbrace 
\alpha_{1},\alpha_{2},\dotsc,\alpha_{k} \rbrace$ and the associated 
formula for the corresponding monic MPC ORF,
\begin{equation*}
\pmb{\pi}^{n}_{k}(z) \! = \! \mathcal{X}_{11}(z) \! = \! 
\widetilde{\phi}_{0}^{\raise-1.0ex\hbox{$\scriptstyle \infty$}}(n,k) \! + \! 
\sum_{m=1}^{\mathfrak{s}-1} \sum_{q=1}^{l_{m}=\varkappa_{nk \tilde{k}_{m}}} 
\dfrac{\widetilde{\nu}^{\raise-1.0ex\hbox{$\scriptstyle \infty$}}_{m,q}
(n,k)}{(z \! - \! \alpha_{p_{m}})^{q}} \! + \! \sum_{q=1}^{l_{\mathfrak{s}}
-1=\varkappa_{nk}-1} 
\widetilde{\nu}^{\raise-1.0ex\hbox{$\scriptstyle \infty$}}_{\mathfrak{s},q}
(n,k)z^{q} \! + \! z^{\varkappa_{nk}},
\end{equation*}
where the $(n \! - \! 1)K \! + \! k$ coefficients 
$\widetilde{\phi}^{\raise-1.0ex\hbox{$\scriptstyle \infty$}}_{0}(n,k),
\widetilde{\nu}^{\raise-1.0ex\hbox{$\scriptstyle \infty$}}_{1,1}(n,k),
\dotsc,\widetilde{\nu}^{\raise-1.0ex\hbox{$\scriptstyle \infty$}}_{1,l_{1}}
(n,k),\dotsc,
\widetilde{\nu}^{\raise-1.0ex\hbox{$\scriptstyle \infty$}}_{\mathfrak{s}-1,1}
(n,k),\dotsc,
\widetilde{\nu}^{\raise-1.0ex\hbox{$\scriptstyle \infty$}}_{\mathfrak{s}-1,
l_{\mathfrak{s}-1}}(n,k),
\widetilde{\nu}^{\raise-1.0ex\hbox{$\scriptstyle \infty$}}_{\mathfrak{s},1}
(n,k),\linebreak[4] 
\dotsc,\widetilde{\nu}^{\raise-1.0ex\hbox{$\scriptstyle \infty$}}_{
\mathfrak{s},l_{\mathfrak{s}}-1}(n,k)$, $(\mu^{\infty}_{n,\varkappa_{nk}}
(n,k))^{-2}$, with $\mu^{\infty}_{n,\varkappa_{nk}}(n,k)$ being the associated 
norming constant (see Corollary~\ref{cor2.1}, Equation~\eqref{eqnmctinf1}, 
below), satisfy the linear inhomogeneous algebraic system of 
equations~\eqref{eq35}. Via the multi-linearity property of the determinant 
and an application of Cramer's Rule to system~\eqref{eq35}, one arrives at, 
for $n \! \in \! \mathbb{N}$ and $k \! \in \! \lbrace 1,2,\dotsc,K \rbrace$ 
such that $\alpha_{p_{\mathfrak{s}}} \! := \! \alpha_{k} \! = \! \infty$, the 
following---ordered---determinantal representation for the corresponding 
monic MPC ORF:
\begin{equation*}
\pmb{\pi}^{n}_{k}(z) \! = \! \mathcal{X}_{11}(z) \! = \! 
\hat{\mathfrak{c}}_{D_{\infty}}^{-1} 
\overset{\infty}{\Xi}^{\raise-1.0ex\hbox{$\scriptstyle n$}}_{k}(z),
\end{equation*}
where $\hat{\mathfrak{c}}_{D_{\infty}}$ is given by Equation~\eqref{eq39}, 
and, with
\begin{equation*}
n_{1} \! = \! l_{1} \! + \! \dotsb \! + \! l_{\mathfrak{s}-2} \! + \! 1, \qquad 
n_{2} \! = \! (n \! - \! 1)K \! + \! k \! - \! \varkappa_{nk}, \qquad m_{1} \! = 
\! (n \! - \! 1)K \! + \! k \! - \! 1, \qquad m_{2} \! = \! (n \! - \! 1)K \! + \! k,
\end{equation*}
and $\md \widetilde{\mu}(z) \! = \! \exp (-n \widetilde{V}(z)) \, \md z$,
\begin{align*}
\overset{\infty}{\Xi}^{\raise-1.0ex\hbox{$\scriptstyle n$}}_{k}(z) := 
\setcounter{MaxMatrixCols}{12}
&\left\vert

\right\vert \\
=& \, \underbrace{\int_{\mathbb{R}} \int_{\mathbb{R}} \dotsb \int_{\mathbb{R}
}}_{(n-1)K+k} \md \widetilde{\mu}(\xi_{0}) \, \md \widetilde{\mu}(\xi_{1}) \, \dotsb 
\, \md \widetilde{\mu}(\xi_{l_{1}}) \, \dotsb \, \md \widetilde{\mu}(\xi_{n_{1}}) \, 
\dotsb \, \md \widetilde{\mu}(\xi_{n_{2}}) \, \md \widetilde{\mu}(\xi_{m_{2}-
\varkappa_{nk}+1}) \, \dotsb \, \md \widetilde{\mu}(\xi_{m_{1}}) \\
\times& \, \dfrac{1}{\dfrac{1}{\xi_{0}^{0}}(\xi_{1} \! - \! \alpha_{p_{1}})^{
1} \dotsb (\xi_{l_{1}} \! - \! \alpha_{p_{1}})^{l_{1}} \dotsb (\xi_{n_{1}} \! 
- \! \alpha_{p_{\mathfrak{s}-1}})^{1} \dotsb (\xi_{n_{2}} \! - \! \alpha_{
p_{\mathfrak{s}-1}})^{l_{\mathfrak{s}-1}} \dfrac{1}{(\xi_{m_{2}-\varkappa_{nk}
+1})^{1}} \dotsb \dfrac{1}{(\xi_{m_{1}})^{\varkappa_{nk}-1}}} \\
\times& \, 
\left\vert

\right\vert.
\end{align*}
For $n \! \in \! \mathbb{N}$ and $k \! \in \! \lbrace 1,2,\dotsc,K \rbrace$ such 
that $\alpha_{p_{\mathfrak{s}}} \! := \! \alpha_{k} \! = \! \infty$, recalling {}from 
the proof of Lemma~\ref{lem2.1} (cf. case~\textbf{(1)}) the $(n \! - \! 1)K \! + \! 
k \! + \! 1$ linearly independent functions $\widehat{\phi}_{0}(z) \! := \! 
\prod_{m=1}^{\mathfrak{s}-1}(z \! - \! \alpha_{p_{m}})^{l_{m}} \! =: \! 
\sum_{j=0}^{(n-1)K+k} \mathfrak{a}^{\spcheck}_{j,0}z^{j}$, $\widehat{\phi}_{q(r_{1})}
(z) \! := \! \lbrace \widehat{\phi}_{0}(z)(z \! - \! \alpha_{p_{r_{1}}})^{-m(r_{1})
(1-\delta_{r_{1} \mathfrak{s}})}z^{m(r_{1}) \delta_{r_{1} \mathfrak{s}}} \! =: \! 
\sum_{j=0}^{(n-1)K+k} \mathfrak{a}^{\spcheck}_{j,q(r_{1})}z^{j} \rbrace$, 
$r_{1} \! = \! 1,\dotsc,\mathfrak{s} \! - \! 1,\mathfrak{s}$, $q(r_{1}) \! = \! 
\sum_{i=1}^{r_{1}-1}l_{i} \! + \! 1,\sum_{i=1}^{r_{1}-1}l_{i} \! + \! 2,\dotsc,
\sum_{i=1}^{r_{1}-1}l_{i} \! + \! l_{r_{1}}$, $m(r_{1}) \! = \! 1,2,\dotsc,
l_{r_{1}}$, and the $(n \! - \! 1)K \! + \! k$ linearly independent functions 
$\psi_{0}(z) \! := \! \prod_{m=1}^{\mathfrak{s}-1}(z \! - \! \alpha_{p_{m}})^{l_{m}} 
\! =: \! \sum_{j=0}^{(n-1)K+k-1} \mathfrak{a}^{\sphat}_{j,0}z^{j}$, $\psi_{q(r_{2})}
(z) \! := \! \lbrace \psi_{0}(z)(z \! - \! \alpha_{p_{r_{2}}})^{-m(r_{2})(1-\delta_{r_{2} 
\mathfrak{s}})}z^{m(r_{2}) \delta_{r_{2} \mathfrak{s}}} \! =: \! \sum_{j=0}^{(n-1)K+k-1} 
\mathfrak{a}^{\sphat}_{j,q(r_{2})}z^{j} \rbrace$, $r_{2} \! = \! 1,\dotsc,\mathfrak{s} 
\! - \! 1,\mathfrak{s}$, $q(r_{2}) \! = \! \sum_{i=1}^{r_{2}-1}l_{i} \! + \! 1,
\sum_{i=1}^{r_{2}-1}l_{i} \! + \! 2,\dotsc,\sum_{i=1}^{r_{2}-1}l_{i} \! + \! l_{r_{2}} 
\! - \! \delta_{r_{2} \mathfrak{s}}$, $m(r_{2}) \! = \! 1,2,\dotsc,l_{r_{2}} \! - \! 
\delta_{r_{2} \mathfrak{s}}$, one proceeds, via the latter determinantal expression, 
with the analysis of 
$\overset{\infty}{\Xi}^{\raise-1.0ex\hbox{$\scriptstyle n$}}_{k}(z)$:
\begin{align*}
\overset{\infty}{\Xi}^{\raise-1.0ex\hbox{$\scriptstyle n$}}_{k}(z) =& \, 
\underbrace{\int_{\mathbb{R}} \int_{\mathbb{R}} \dotsb \int_{\mathbb{R}}}_{(
n-1)K+k} \md \widetilde{\mu}(\xi_{0}) \, \md \widetilde{\mu}(\xi_{1}) \, \dotsb 
\, \md \widetilde{\mu}(\xi_{l_{1}}) \, \dotsb \, \md \widetilde{\mu}(\xi_{n_{1}}) 
\, \dotsb \, \md \widetilde{\mu}(\xi_{n_{2}}) \, \md \widetilde{\mu}(\xi_{m_{2}-
\varkappa_{nk}+1}) \, \dotsb \, \md \widetilde{\mu}(\xi_{m_{1}}) \\
\times& \, \dfrac{1}{\dfrac{1}{\xi_{0}^{0}}(\xi_{1} \! - \! \alpha_{p_{1}})^{
1} \dotsb (\xi_{l_{1}} \! - \! \alpha_{p_{1}})^{l_{1}} \dotsb (\xi_{n_{1}} \! 
- \! \alpha_{p_{\mathfrak{s}-1}})^{1} \dotsb (\xi_{n_{2}} \! - \! \alpha_{p_{
\mathfrak{s}-1}})^{l_{\mathfrak{s}-1}} \dfrac{1}{(\xi_{m_{2}-\varkappa_{nk}
+1})^{1}} \dotsb \dfrac{1}{\xi_{m_{1}}^{\varkappa_{nk}-1}}} \\
\times& \, \dfrac{(\widehat{\phi}_{0}(z))^{-1}}{\widehat{\phi}_{0}(\xi_{0}) 
\widehat{\phi}_{0}(\xi_{1}) \dotsb \widehat{\phi}_{0}(\xi_{l_{1}}) \dotsb 
\widehat{\phi}_{0}(\xi_{n_{1}}) \dotsb \widehat{\phi}_{0}(\xi_{n_{2}}) 
\widehat{\phi}_{0}(\xi_{m_{2}-\varkappa_{nk}+1}) \dotsb \widehat{\phi}_{0}
(\xi_{m_{1}})} \\
\times& \, 
\underbrace{\left\vert

\right\vert}_{=: \, \mathbb{G}^{\blacktriangledown}(\xi_{0},\xi_{1},\dotsc,
\xi_{l_{1}},\dotsc,\xi_{n_{1}},\dotsc,\xi_{n_{2}},\dotsc,\xi_{(n-1)K+k-1};z)} \\
=& \, \dfrac{(\, \widehat{\phi}_{0}(z))^{-1}}{m_{2}!} \sum_{\pmb{\sigma} \in 
\mathfrak{S}_{m_{2}}} \underbrace{\int_{\mathbb{R}} \int_{\mathbb{R}} \dotsb 
\int_{\mathbb{R}}}_{(n-1)K+k} \md \widetilde{\mu}(\xi_{\sigma (0)}) \, \md 
\widetilde{\mu}(\xi_{\sigma (1)}) \, \dotsb \, \md \widetilde{\mu}(\xi_{\sigma (l_{1})}) 
\, \dotsb \, \md \widetilde{\mu}(\xi_{\sigma (n_{1})}) \, \dotsb \, \md \widetilde{\mu}
(\xi_{\sigma (n_{2})}) \\
\times& \, \tfrac{\dotsb \, \md \widetilde{\mu}(\xi_{\sigma (m_{2}-\varkappa_{nk}+1)}) 
\, \dotsb \, \md \widetilde{\mu}(\xi_{\sigma (m_{1})})}{\dfrac{1}{\xi_{\sigma (0)}^{0}}
(\xi_{\sigma (1)}-\alpha_{p_{1}})^{1} \dotsb (\xi_{\sigma (l_{1})}-\alpha_{p_{1}})^{l_{1}} 
\dotsb (\xi_{\sigma (n_{1})}-\alpha_{p_{\mathfrak{s}-1}})^{1} \dotsb 
(\xi_{\sigma (n_{2})}-\alpha_{p_{\mathfrak{s}-1}})^{l_{\mathfrak{s}-1}} 
\dfrac{1}{(\xi_{\sigma (m_{2}-\varkappa_{nk}+1)})^{1}} \dotsb 
\dfrac{1}{\xi_{\sigma (m_{1})}^{\varkappa_{nk}-1}}} \\
\times& \, \dfrac{1}{\widehat{\phi}_{0}(\xi_{\sigma (0)}) \widehat{\phi}_{0}
(\xi_{\sigma (1)}) \dotsb \widehat{\phi}_{0}(\xi_{\sigma (l_{1})}) \dotsb 
\widehat{\phi}_{0}(\xi_{\sigma (n_{1})}) \dotsb \widehat{\phi}_{0}(\xi_{\sigma 
(n_{2})}) \widehat{\phi}_{0}(\xi_{\sigma (m_{2}-\varkappa_{nk}+1)}) \dotsb 
\widehat{\phi}_{0}(\xi_{\sigma (m_{1})})} \\
\times& \, \mathbb{G}^{\blacktriangledown}(\xi_{\sigma (0)},\xi_{\sigma (1)},
\dotsc,\xi_{\sigma (l_{1})},\dotsc,\xi_{\sigma (n_{1})},\dotsc,
\xi_{\sigma (n_{2})},\dotsc,\xi_{\sigma ((n-1)K+k-1)};z) \\
=& \, \dfrac{(\, \widehat{\phi}_{0}(z))^{-1}}{m_{2}!} \underbrace{\int_{\mathbb{R}} 
\int_{\mathbb{R}} \dotsb \int_{\mathbb{R}}}_{(n-1)K+k} \md \widetilde{\mu}
(\xi_{0}) \, \md \widetilde{\mu}(\xi_{1}) \, \dotsb \, \md \widetilde{\mu}(\xi_{l_{1}}) 
\, \dotsb \, \md \widetilde{\mu}(\xi_{n_{1}}) \, \dotsb \, \md \widetilde{\mu}
(\xi_{n_{2}}) \, \md \widetilde{\mu}(\xi_{m_{2}-\varkappa_{nk}+1}) \, \dotsb \\
\times& \, \dotsb \, \md \widetilde{\mu}(\xi_{m_{1}}) \, \mathbb{G}^{
\blacktriangledown}(\xi_{0},\xi_{1},\dotsc,\xi_{l_{1}},\dotsc,\xi_{n_{1}},
\dotsc,\xi_{n_{2}},\dotsc,\xi_{(n-1)K+k-1};z) \\
\times& \, \sum_{\pmb{\sigma} \in \mathfrak{S}_{m_{2}}} 
\operatorname{sgn}(\pmb{\pmb{\sigma}}) 
\tfrac{1}{\dfrac{1}{\xi_{\sigma (0)}^{0}}(\xi_{\sigma (1)}-\alpha_{p_{1}})^{1} 
\dotsb (\xi_{\sigma (l_{1})}-\alpha_{p_{1}})^{l_{1}} \dotsb (\xi_{\sigma 
(n_{1})}-\alpha_{p_{\mathfrak{s}-1}})^{1} \dotsb (\xi_{\sigma (n_{2})}-
\alpha_{p_{\mathfrak{s}-1}})^{l_{\mathfrak{s}-1}} \dfrac{1}{(\xi_{\sigma 
(m_{2}-\varkappa_{nk}+1)})^{1}} \dotsb \dfrac{1}{\xi_{\sigma (m_{1})}^{
\varkappa_{nk}-1}}} \\
\times& \, \dfrac{1}{\widehat{\phi}_{0}(\xi_{\sigma (0)}) \widehat{\phi}_{0}
(\xi_{\sigma (1)}) \dotsb \widehat{\phi}_{0}(\xi_{\sigma (l_{1})}) \dotsb 
\widehat{\phi}_{0}(\xi_{\sigma (n_{1})}) \dotsb \widehat{\phi}_{0}(\xi_{\sigma 
(n_{2})}) \widehat{\phi}_{0}(\xi_{\sigma (m_{2}-\varkappa_{nk}+1)}) \dotsb 
\widehat{\phi}_{0}(\xi_{\sigma (m_{1})})} \\
=& \, \dfrac{(\, \widehat{\phi}_{0}(z))^{-1}}{m_{2}!} \underbrace{\int_{\mathbb{
R}} \int_{\mathbb{R}} \dotsb \int_{\mathbb{R}}}_{(n-1)K+k} \md \widetilde{\mu}
(\xi_{0}) \, \md \widetilde{\mu}(\xi_{1}) \, \dotsb \, \md \widetilde{\mu}
(\xi_{l_{1}}) \, \dotsb \, \md \widetilde{\mu}(\xi_{n_{1}}) \, \dotsb \, \md 
\widetilde{\mu}(\xi_{n_{2}}) \, \md \widetilde{\mu}(\xi_{m_{2}-\varkappa_{nk}+1}) 
\, \dotsb \\
\times& \, \dotsb \, \md \widetilde{\mu}(\xi_{m_{1}}) \, \dfrac{\mathbb{G}^{
\blacktriangledown}(\xi_{0},\xi_{1},\dotsc,\xi_{l_{1}},\dotsc,\xi_{n_{1}},
\dotsc,\xi_{n_{2}},\dotsc,\xi_{(n-1)K+k-1};z)}{\widehat{\phi}_{0}(\xi_{0}) 
\widehat{\phi}_{0}(\xi_{1}) \dotsb \widehat{\phi}_{0}(\xi_{l_{1}}) \dotsb 
\widehat{\phi}_{0}(\xi_{n_{1}}) \dotsb \widehat{\phi}_{0}(\xi_{n_{2}}) 
\widehat{\phi}_{0}(\xi_{m_{2}-\varkappa_{nk}+1}) \dotsb \widehat{\phi}_{0}
(\xi_{m_{1}})} \\
\times& \, 
\left\vert 

\right\vert \\
=& \, \dfrac{(\, \widehat{\phi}_{0}(z))^{-1}}{m_{2}!} \underbrace{\int_{\mathbb{
R}} \int_{\mathbb{R}} \dotsb \int_{\mathbb{R}}}_{(n-1)K+k} \md \widetilde{\mu}
(\xi_{0}) \, \md \widetilde{\mu}(\xi_{1}) \, \dotsb \, \md \widetilde{\mu}(\xi_{l_{1}}) 
\, \dotsb \, \md \widetilde{\mu}(\xi_{n_{1}}) \, \dotsb \, \md \widetilde{\mu}
(\xi_{n_{2}}) \, \md \widetilde{\mu}(\xi_{m_{2}-\varkappa_{nk}+1}) \, \dotsb \\
\times& \, \dotsb \, \md \widetilde{\mu}(\xi_{m_{1}}) \, \dfrac{\mathbb{G}^{
\blacktriangledown}(\xi_{0},\xi_{1},\dotsc,\xi_{l_{1}},\dotsc,\xi_{n_{1}},
\dotsc,\xi_{n_{2}},\dotsc,\xi_{(n-1)K+k-1};z)}{\widehat{\phi}_{0}(\xi_{0}) 
\widehat{\phi}_{0}(\xi_{1}) \dotsb \widehat{\phi}_{0}(\xi_{l_{1}}) \dotsb 
\widehat{\phi}_{0}(\xi_{n_{1}}) \dotsb \widehat{\phi}_{0}(\xi_{n_{2}}) 
\widehat{\phi}_{0}(\xi_{m_{2}-\varkappa_{nk}+1}) \dotsb \widehat{\phi}_{0}
(\xi_{m_{1}})} \\
\times& \, \dfrac{1}{\psi_{0}(\xi_{0}) \psi_{0}(\xi_{1}) \dotsb \psi_{0}(\xi_{
l_{1}}) \dotsb \psi_{0}(\xi_{n_{1}}) \dotsb \psi_{0}(\xi_{n_{2}}) \psi_{0}
(\xi_{m_{2}-\varkappa_{nk}+1}) \dotsb \psi_{0}(\xi_{m_{1}})} \\
\times& \, 
\underbrace{\left\vert
\begin{smallmatrix}
\psi_{0}(\xi_{0}) & \psi_{1}(\xi_{0}) & \dotsb & \psi_{l_{1}}(\xi_{0}) & 
\dotsb & \psi_{n_{1}}(\xi_{0}) & \dotsb & \psi_{n_{2}}(\xi_{0}) & \psi_{m_{2}
-\varkappa_{nk}+1}(\xi_{0}) & \dotsb & \psi_{(n-1)K+k-1}(\xi_{0}) \\
\psi_{0}(\xi_{1}) & \psi_{1}(\xi_{1}) & \dotsb & \psi_{l_{1}}(\xi_{1}) & 
\dotsb & \psi_{n_{1}}(\xi_{1}) & \dotsb & \psi_{n_{2}}(\xi_{1}) & \psi_{m_{2}
-\varkappa_{nk}+1}(\xi_{1}) & \dotsb & \psi_{(n-1)K+k-1}(\xi_{1}) \\
\vdots & \vdots & \ddots & \vdots & \ddots & \vdots & \ddots & \vdots & \vdots 
& \ddots & \vdots \\
\psi_{0}(\xi_{l_{1}}) & \psi_{1}(\xi_{l_{1}}) & \dotsb & \psi_{l_{1}}(\xi_{l_{
1}}) & \dotsb & \psi_{n_{1}}(\xi_{l_{1}}) & \dotsb & \psi_{n_{2}}(\xi_{l_{1}}) 
& \psi_{m_{2}-\varkappa_{nk}+1}(\xi_{l_{1}}) & \dotsb & \psi_{(n-1)K+k-1}(\xi_{
l_{1}}) \\
\vdots & \vdots & \ddots & \vdots & \ddots & \vdots & \ddots & \vdots & \vdots 
& \ddots & \vdots \\
\psi_{0}(\xi_{n_{1}}) & \psi_{1}(\xi_{n_{1}}) & \dotsb & \psi_{l_{1}}(\xi_{n_{
1}}) & \dotsb & \psi_{n_{1}}(\xi_{n_{1}}) & \dotsb & \psi_{n_{2}}(\xi_{n_{1}}) 
& \psi_{m_{2}-\varkappa_{nk}+1}(\xi_{n_{1}}) & \dotsb & \psi_{(n-1)K+k-1}(\xi_{
n_{1}}) \\
\vdots & \vdots & \ddots & \vdots & \ddots & \vdots & \ddots & \vdots & \vdots 
& \ddots & \vdots \\
\psi_{0}(\xi_{n_{2}}) & \psi_{1}(\xi_{n_{2}}) & \dotsb & \psi_{l_{1}}(\xi_{n_{
2}}) & \dotsb & \psi_{n_{1}}(\xi_{n_{2}}) & \dotsb & \psi_{n_{2}}(\xi_{n_{2}}) 
& \psi_{m_{2}-\varkappa_{nk}+1}(\xi_{n_{2}}) & \dotsb & \psi_{(n-1)K+k-1}
(\xi_{n_{2}}) \\
\vdots & \vdots & \ddots & \vdots & \ddots & \vdots & \ddots & \vdots & \vdots 
& \ddots & \vdots \\
\psi_{0}(\xi_{m_{1}}) & \psi_{1}(\xi_{m_{1}}) & \dotsb & \psi_{l_{1}}(\xi_{m_{
1}}) & \dotsb & \psi_{n_{1}}(\xi_{m_{1}}) & \dotsb & \psi_{n_{2}}(\xi_{m_{1}}) 
& \psi_{m_{2}-\varkappa_{nk}+1}(\xi_{m_{1}}) & \dotsb & \psi_{(n-1)K+k-1}
(\xi_{m_{1}})
\end{smallmatrix}
\right\vert}_{= \, \mathbb{D}^{\sphat} \prod_{\underset{j<i}{i,j=0}}^{(n-1)
K+k-1}(\xi_{i}-\xi_{j}) \quad \text{(cf. proof of Lemma~\ref{lem2.1}, 
case~\pmb{(1)})}} \\
=& \, \dfrac{\mathbb{D}^{\sphat} \, (\, \widehat{\phi}_{0}(z))^{-1}}{m_{2}!} 
\underbrace{\int_{\mathbb{R}} \int_{\mathbb{R}} \dotsb \int_{\mathbb{R}}}_{(
n-1)K+k} \md \widetilde{\mu}(\xi_{0}) \, \md \widetilde{\mu}(\xi_{1}) \, \dotsb 
\, \md \widetilde{\mu}(\xi_{l_{1}}) \, \dotsb \, \md \widetilde{\mu}(\xi_{n_{1}}) 
\, \dotsb \, \md \widetilde{\mu}(\xi_{n_{2}}) \, \md \widetilde{\mu}
(\xi_{m_{2}-\varkappa_{nk}+1}) \, \dotsb \\
\times& \, \dotsb \, \md \widetilde{\mu}(\xi_{m_{1}}) \, 
\dfrac{\prod_{\underset{j<i}{i,j=0}}^{(n-1)K+k-1}(\xi_{i} \! - \! 
\xi_{j})}{\widehat{\phi}_{0}(\xi_{0}) \widehat{\phi}_{0}(\xi_{1}) \dotsb 
\widehat{\phi}_{0}(\xi_{l_{1}}) \dotsb \widehat{\phi}_{0}(\xi_{n_{1}}) 
\dotsb \widehat{\phi}_{0}(\xi_{n_{2}}) \widehat{\phi}_{0}(\xi_{m_{2}-
\varkappa_{nk}+1}) \dotsb \widehat{\phi}_{0}(\xi_{m_{1}})} \\
\times& \, \dfrac{1}{\psi_{0}(\xi_{0}) \psi_{0}(\xi_{1}) \dotsb \psi_{0}(\xi_{
l_{1}}) \dotsb \psi_{0}(\xi_{n_{1}}) \dotsb \psi_{0}(\xi_{n_{2}}) \psi_{0}
(\xi_{m_{2}-\varkappa_{nk}+1}) \dotsb \psi_{0}(\xi_{m_{1}})} \\
\times& \, 
\underbrace{\left\vert 

\right\vert}_{= \, \mathbb{D}^{\spcheck} \prod_{\underset{j<i}{i,j=0}}^{
(n-1)K+k-1}(\xi_{i}-\xi_{j}) \prod_{m=0}^{(n-1)K+k-1}(z-\xi_{m}) 
\quad \text{(cf. proof of Lemma~\ref{lem2.1}, case~\pmb{(1)})}} \\
=& \, \dfrac{\mathbb{D}^{\sphat} \, \mathbb{D}^{\spcheck}
(\, \widehat{\phi}_{0}(z))^{-1}}{((n \! - \! 1)K \! + \! k)!} \underbrace{\int_{
\mathbb{R}} \int_{\mathbb{R}} \dotsb \int_{\mathbb{R}}}_{(n-1)K+k} \md 
\widetilde{\mu}(\xi_{0}) \, \md \widetilde{\mu}(\xi_{1}) \, \dotsb \, \md 
\widetilde{\mu}(\xi_{l_{1}}) \, \dotsb \, \md \widetilde{\mu}(\xi_{n_{1}}) 
\, \dotsb \, \md \widetilde{\mu}(\xi_{n_{2}}) \, \md \widetilde{\mu}
(\xi_{m_{2}-\varkappa_{nk}+1}) \, \dotsb \\
\times& \, \dotsb \, \md \widetilde{\mu}(\xi_{m_{1}}) \, \dfrac{(\psi_{0}
(\xi_{0}) \psi_{0}(\xi_{1}) \dotsb \psi_{0}(\xi_{l_{1}}) \dotsb \psi_{0}(\xi_{n_{1}}) 
\dotsb \psi_{0}(\xi_{n_{2}}) \psi_{0}(\xi_{m_{2}-\varkappa_{nk}+1}) \dotsb 
\psi_{0}(\xi_{m_{1}}))^{-1}}{\widehat{\phi}_{0}(\xi_{0}) \widehat{\phi}_{0}
(\xi_{1}) \dotsb \widehat{\phi}_{0}(\xi_{l_{1}}) \dotsb \widehat{\phi}_{0}
(\xi_{n_{1}}) \dotsb \widehat{\phi}_{0}(\xi_{n_{2}}) \widehat{\phi}_{0}
(\xi_{m_{2}-\varkappa_{nk}+1}) \dotsb \widehat{\phi}_{0}(\xi_{m_{1}})} \\
\times& \, \prod_{\substack{i,j=0\\j<i}}^{(n-1)K+k-1}(\xi_{i} \! - \! 
\xi_{j})^{2} \, \prod_{m=0}^{(n-1)K+k-1}(z \! - \! \xi_{m});
\end{align*}
but, recalling that, for $n \! \in \! \mathbb{N}$ and $k \! \in \! \lbrace 
1,2,\dotsc,K \rbrace$ such that $\alpha_{p_{\mathfrak{s}}} \! := \! 
\alpha_{k} \! = \! \infty$, $\widehat{\phi}_{0}(z) \! = \! \psi_{0}(z)$, 
one arrives at
\begin{align*}
\overset{\infty}{\Xi}^{\raise-1.0ex\hbox{$\scriptstyle n$}}_{k}(z) =& \, 
\dfrac{\mathbb{D}^{\sphat} \, \mathbb{D}^{\spcheck}(\, \widehat{\phi}_{0}
(z))^{-1}}{((n \! - \! 1)K \! + \! k)!} \underbrace{\int_{\mathbb{R}} 
\int_{\mathbb{R}} \dotsb \int_{\mathbb{R}}}_{(n-1)K+k} \md \widetilde{\mu}
(\xi_{0}) \, \md \widetilde{\mu}(\xi_{1}) \, \dotsb \, \md \widetilde{\mu}
(\xi_{l_{1}}) \, \dotsb \, \md \widetilde{\mu}(\xi_{n_{1}}) \, \dotsb \, \md 
\widetilde{\mu}(\xi_{n_{2}}) \\
\times& \, \md \widetilde{\mu}(\xi_{(n-1)K+k-\varkappa_{nk}+1}) \, \dotsb 
\, \md \widetilde{\mu}(\xi_{(n-1)K+k-1}) \left(\prod_{l=0}^{(n-1)K+k-1} 
\widehat{\phi}_{0}(\xi_{l}) \right)^{-2} \prod_{\substack{i,j=0\\j<i}}^{(n-1)K
+k-1}(\xi_{i} \! - \! \xi_{j})^{2} \prod_{m=0}^{(n-1)K+k-1}(z \! - \! \xi_{m}) 
\, \, \Rightarrow \\
\overset{\infty}{\Xi}^{\raise-1.0ex\hbox{$\scriptstyle n$}}_{k}(z) =& \, 
\dfrac{\mathbb{D}^{\sphat} \, \mathbb{D}^{\spcheck}}{((n \! - \! 1)K \! + \! k)!} 
\left(\prod_{q=1}^{\mathfrak{s}-1}(z \! - \! \alpha_{p_{q}})^{\varkappa_{nk 
\tilde{k}_{q}}} \right)^{-1} \underbrace{\int_{\mathbb{R}} \int_{\mathbb{R}} 
\dotsb \int_{\mathbb{R}}}_{(n-1)K+k} \md \widetilde{\mu}(\xi_{0}) \, \md 
\widetilde{\mu}(\xi_{1}) \, \dotsb \, \md \widetilde{\mu}(\xi_{(n-1)K+k-1}) \\
\times& \, \left(\prod_{m=0}^{(n-1)K+k-1} \prod_{q=1}^{\mathfrak{s}-1}
(\xi_{m} \! - \! \alpha_{p_{q}})^{\varkappa_{nk \tilde{k}_{q}}} \right)^{-2} 
\prod_{\substack{i,j=0\\j<i}}^{(n-1)K+k-1}(\xi_{i} \! - \! \xi_{j})^{2} \, 
\prod_{l=0}^{(n-1)K+k-1}(z \! - \! \xi_{l}).
\end{align*}
Recalling that, for $n \! \in \! \mathbb{N}$ and $k \! \in \! \lbrace 1,2,
\dotsc,K \rbrace$ such that $\alpha_{p_{\mathfrak{s}}} \! := \! \alpha_{k} 
\! = \! \infty$, $\pmb{\pi}^{n}_{k}(z) \! = \! \mathcal{X}_{11}(z) \! = \! 
\overset{\infty}{\Xi}^{\raise-1.0ex\hbox{$\scriptstyle n$}}_{k}(z)/
\hat{\mathfrak{c}}_{D_{\infty}}$, where 
$\overset{\infty}{\Xi}^{\raise-1.0ex\hbox{$\scriptstyle n$}}_{k}(z)$ is 
given directly above, and $\hat{\mathfrak{c}}_{D_{\infty}}$ is given by 
Equation~\eqref{eq39}, one arrives at the integral representation for 
$\pmb{\pi}^{n}_{k}(z)$ given in Equation~\eqref{eintreppinf1}.

The determinantal representation for $\mathcal{X}_{21}(z)$ is now considered. 
For $n \! \in \! \mathbb{N}$ and $k \! \in \! \lbrace 1,2,\dotsc,K \rbrace$ 
such that $\alpha_{p_{\mathfrak{s}}} \! := \! \alpha_{k} \! = \! \infty$, recall 
the ordered disjoint partition for $\lbrace \alpha_{1},\alpha_{2},\dotsc,
\alpha_{K} \rbrace \cup \lbrace \alpha_{1},\alpha_{2},\dotsc,\alpha_{K} 
\rbrace \cup \dotsb \cup \lbrace \alpha_{1},\alpha_{2},\dotsc,\alpha_{k} 
\rbrace$ introduced in the proof of Lemma~\ref{lem2.1} (cf. case~\textbf{(1)}). 
For this ordered disjoint partition, introduce, for $n \! \in \! \mathbb{N}$ and 
$k \! \in \! \lbrace 1,2,\dotsc,K \rbrace$ such that $\alpha_{p_{\mathfrak{s}}} 
\! := \! \alpha_{k} \! = \! \infty$, the following notation (recall that $l_{q} \! 
= \! \varkappa_{nk \tilde{k}_{q}}$, $q \! = \! 1,2,\dotsc,\mathfrak{s} \! - \! 
1$, and $l_{\mathfrak{s}} \! = \! \varkappa_{nk})$:
\begin{equation*}
\widehat{\psi}_{0}(z) \! := \! \prod_{m=1}^{\mathfrak{s}-1}
(z \! - \! \alpha_{p_{m}})^{l_{m}} \! =: \sum_{j=0}^{(n-1)K+k-2} 
\mathfrak{a}^{\spadesuit}_{j,0}z^{j},
\end{equation*}
and, for $r \! = \! 1,\dotsc,\mathfrak{s} \! - \! 1,\mathfrak{s}$, $q(r) \! 
= \! \sum_{i=1}^{r-1}l_{i} \! + \! 1,\sum_{i=1}^{r-1}l_{i} \! + \! 2,\dotsc,
\sum_{i=1}^{r-1}l_{i} \! + \! l_{r} \! - \! 2 \delta_{r \mathfrak{s}}$, and 
$m(r) \! = \! 1,2,\dotsc,l_{r} \! - \! 2 \delta_{r \mathfrak{s}}$,\footnote{With 
the convention $(z \! - \! \infty)^{0} \! := \! 1$.}
\begin{equation*}
\widehat{\psi}_{q(r)}(z) \! := \! \left\lbrace \widehat{\psi}_{0}(z)
(z \! - \! \alpha_{p_{r}})^{-m(r)(1-\delta_{r \mathfrak{s}})}
z^{m(r) \delta_{r \mathfrak{s}}} \! =: \sum_{j=0}^{(n-1)K+k-2} 
\mathfrak{a}^{\spadesuit}_{j,q(r)}z^{j} \right\rbrace;
\end{equation*}
e.g., for $r \! = \! 1$, the notation $\widehat{\psi}_{q(1)}(z) \! := \! 
\lbrace \widehat{\psi}_{0}(z)(z \! - \! \alpha_{p_{1}})^{-m(1)} \! =: \! 
\sum_{j=0}^{(n-1)K+k-2} \mathfrak{a}^{\spadesuit}_{j,q(1)}z^{j} \rbrace$, 
$q(1) \! = \! 1,2,\dotsc,l_{1}$, $m(1) \! = \! 1,2,\dotsc,l_{1}$, denotes 
the (set of) $l_{1} \! = \! \varkappa_{nk \tilde{k}_{1}}$ functions
\begin{gather*}
\widehat{\psi}_{1}(z) \! = \! \dfrac{\widehat{\psi}_{0}(z)}{z \! - \! 
\alpha_{p_{1}}} \! =: \sum_{j=0}^{(n-1)K+k-2} \mathfrak{a}^{\spadesuit}_{j,1}
z^{j}, \, \widehat{\psi}_{2}(z) \! = \! \dfrac{\widehat{\psi}_{0}(z)}{
(z \! - \! \alpha_{p_{1}})^{2}} \! =: \sum_{j=0}^{(n-1)K+k-2} 
\mathfrak{a}^{\spadesuit}_{j,2}z^{j}, \, \dotsc, \, \widehat{\psi}_{l_{1}}
(z) \! = \! \dfrac{\widehat{\psi}_{0}(z)}{(z \! - \! \alpha_{p_{1}})^{l_{1}}} 
\! =: \sum_{j=0}^{(n-1)K+k-2} \mathfrak{a}^{\spadesuit}_{j,l_{1}}z^{j},
\end{gather*}
for $r \! = \! \mathfrak{s} \! - \! 1$, the notation 
$\widehat{\psi}_{q(\mathfrak{s}-1)}(z) \! := \! \lbrace \widehat{\psi}_{0}
(z)(z \! - \! \alpha_{p_{\mathfrak{s}-1}})^{-m(\mathfrak{s}-1)} \! =: \! 
\sum_{j=0}^{(n-1)K+k-2} \mathfrak{a}^{\spadesuit}_{j,q(\mathfrak{s}-1)}z^{j} 
\rbrace$, $q(\mathfrak{s} \! - \! 1) \! = \! l_{1} \! + \! \dotsb \! + \! 
l_{\mathfrak{s}-2} \! + \! 1,l_{1} \! + \! \dotsb \! + \! l_{\mathfrak{s}-2} 
\! + \! 2,\dotsc,l_{1} \! + \! \dotsb \! + \! l_{\mathfrak{s}-2} \! + \! 
l_{\mathfrak{s}-1} \! = \! (n \! - \! 1)K \! + \! k \! - \! \varkappa_{nk}$, 
$m(\mathfrak{s} \! - \! 1) \! = \! 1,2,\dotsc,l_{\mathfrak{s}-1}$, denotes 
the (set of) $l_{\mathfrak{s}-1} \! = \! \varkappa_{nk \tilde{k}_{\mathfrak{s}
-1}}$ functions
\begin{gather*}
\widehat{\psi}_{l_{1}+ \dotsb +l_{\mathfrak{s}-2}+1}(z) \! = \! 
\dfrac{\widehat{\psi}_{0}(z)}{z \! - \! \alpha_{p_{\mathfrak{s}-1}}} \! 
=: \sum_{j=0}^{(n-1)K+k-2} \mathfrak{a}^{\spadesuit}_{j,l_{1}+\dotsb +
l_{\mathfrak{s}-2}+1}z^{j}, \, \widehat{\psi}_{l_{1}+\dotsb +
l_{\mathfrak{s}-2}+2}(z) \! = \! \dfrac{\widehat{\psi}_{0}(z)}{(z \! - \! 
\alpha_{p_{\mathfrak{s}-1}})^{2}} \\
=: \sum_{j=0}^{(n-1)K+k-2} \mathfrak{a}^{\spadesuit}_{j,l_{1}+\dotsb 
+l_{\mathfrak{s}-2}+2}z^{j}, \, \dotsc, \, \widehat{\psi}_{(n-1)K+k-
\varkappa_{nk}}(z) \! = \! \dfrac{\widehat{\psi}_{0}(z)}{(z \! - \! 
\alpha_{p_{\mathfrak{s}-1}})^{l_{\mathfrak{s}-1}}} \! =: 
\sum_{j=0}^{(n-1)K+k-2} \mathfrak{a}^{\spadesuit}_{j,(n-1)K+k-
\varkappa_{nk}}z^{j},
\end{gather*}
etc., and, for $r \! = \! \mathfrak{s}$, the notation $\widehat{\psi}_{
q(\mathfrak{s})}(z) \! := \! \lbrace \widehat{\psi}_{0}(z)z^{m(\mathfrak{s})} 
\! =: \! \sum_{j=0}^{(n-1)K+k-2} \mathfrak{a}^{\spadesuit}_{j,q(\mathfrak{s})}
z^{j} \rbrace$, $q(\mathfrak{s}) \! = \! (n \! - \! 1)K \! + \! k \! - \! \varkappa_{nk} 
\! + \! 1,(n \! - \! 1)K \! + \! k \! - \! \varkappa_{nk} \! + \! 2,\dotsc,(n \! - \! 
1)K \! + \! k \! - \! 2$, $m(\mathfrak{s}) \! = \! 1,2,\dotsc,l_{\mathfrak{s}} \! 
- \! 2$, denotes the (set of) $l_{\mathfrak{s}} \! - \! 2 \! = \! \varkappa_{nk} 
\! - \! 2$ functions
\begin{gather*}
\widehat{\psi}_{(n-1)K+k-\varkappa_{nk}+1}(z) \! = \! \widehat{\psi}_{0}
(z)z \! =: \sum_{j=0}^{(n-1)K+k-2} \mathfrak{a}^{\spadesuit}_{j,(n-1)K+k-
\varkappa_{nk}+1}z^{j}, \, \widehat{\psi}_{(n-1)K+k-\varkappa_{nk}+2}(z) 
\! = \! \widehat{\psi}_{0}(z)z^{2} \\
=: \sum_{j=0}^{(n-1)K+k-2} \mathfrak{a}^{\spadesuit}_{j,(n-1)K+k-
\varkappa_{nk}+2}z^{j}, \, \dotsc, \, \widehat{\psi}_{(n-1)K+k-2}(z) \! 
= \! \widehat{\psi}_{0}(z)z^{\varkappa_{nk}-2} \! =: \sum_{j=0}^{(n-1)K+k-2} 
\mathfrak{a}^{\spadesuit}_{j,(n-1)K+k-2}z^{j}.
\end{gather*}
(Note: $\# \lbrace \widehat{\psi}_{0}(z)(z \! - \! \alpha_{p_{r}})^{-m(r)(1-
\delta_{r \mathfrak{s}})}z^{m(r) \delta_{r \mathfrak{s}}} \rbrace \! = \! l_{r} 
\! - \! 2 \delta_{r \mathfrak{s}}$, $r \! = \! 1,\dotsc,\mathfrak{s} \! - \! 1,
\mathfrak{s}$, and $\# \cup_{r=1}^{\mathfrak{s}} \lbrace \widehat{\psi}_{0}
(z)(z \! - \! \alpha_{p_{r}})^{-m(r)(1-\delta_{r \mathfrak{s}})}z^{m(r) 
\delta_{r \mathfrak{s}}} \rbrace \! = \! \sum_{r=1}^{\mathfrak{s}}l_{r} \! - \! 
2 \! = \! (n \! - \! 1)K \! + \! k \! - \! 2$.) One notes that, for $n \! \in \! 
\mathbb{N}$ and $k \! \in \! \lbrace 1,2,\dotsc,K \rbrace$ such that 
$\alpha_{p_{\mathfrak{s}}} \! := \! \alpha_{k} \! = \! \infty$, the $l_{1} \! + \! 
\dotsb \! + \! l_{\mathfrak{s}-1} \! + \! l_{\mathfrak{s}} \! - \! 1 \! = \! (n \! 
- \! 1)K \! + \! k \! - \! 1$ functions $\widehat{\psi}_{0}(z),\widehat{\psi}_{1}
(z),\dotsc,\widehat{\psi}_{l_{1}}(z),\dotsc,\widehat{\psi}_{l_{1}+\dotsb +
l_{\mathfrak{s}-2}+1}(z),\dotsc,\widehat{\psi}_{(n-1)K+k-\varkappa_{nk}}(z),
\widehat{\psi}_{(n-1)K+k-\varkappa_{nk}+1}(z),\dotsc,\widehat{\psi}_{(n-1)K
+k-2}(z)$ are linearly independent on $\mathbb{R}$, that is, for $z \! \in 
\! \mathbb{R}$, $\sum_{j=0}^{(n-1)K+k-2} \mathfrak{c}^{\spadesuit}_{j} 
\widehat{\psi}_{j}(z) \! = \! 0$ $\Rightarrow$ (via a Vandermonde-type argument; 
see the $((n \! - \! 1)K \! + \! k \! - \! 1) \times ((n \! - \! 1)K \! + \! k \! - \! 1)$ 
non-zero determinant $\mathbb{D}^{\spadesuit}$ in Equation~\eqref{eq76} below) 
$\mathfrak{c}^{\spadesuit}_{j} \! = \! 0$, $j \! = \! 0,1,\dotsc,(n \! - \! 1)K \! + 
\! k \! - \! 2$. For $n \! \in \! \mathbb{N}$ and $k \! \in \! \lbrace 1,2,\dotsc,K 
\rbrace$ such that $\alpha_{p_{\mathfrak{s}}} \! := \! \alpha_{k} \! = \! \infty$, 
let $\mathfrak{S}_{(n-1)K+k-1}$ denote the $((n \! - \! 1)K \! + \! k \! - \! 1)!$ 
permutations of $\lbrace 0,1,\dotsc,(n \! - \! 1)K \! + \! k \! - \! 2 \rbrace$. 
For $n \! \in \! \mathbb{N}$ and $k \! \in \! \lbrace 1,2,\dotsc,K \rbrace$ 
such that $\alpha_{p_{\mathfrak{s}}} \! := \! \alpha_{k} \! = \! \infty$, 
corresponding to the ordered disjoint partition above, recall the formula 
for $\mathcal{X}_{21}(z)$ given in the proof of Lemma~\ref{lem2.1} (cf. 
case~\textbf{(1)}):
\begin{equation*}
\mathcal{X}_{21}(z) \! = \! \sum_{m=1}^{\mathfrak{s}-1} \sum_{q=1}^{l_{m}
=\varkappa_{nk \tilde{k}_{m}}} 
\dfrac{\widehat{\nu}_{m,q}^{\raise-1.0ex\hbox{$\scriptstyle \infty$}}
(n,k)}{(z \! - \! \alpha_{p_{m}})^{q}} \! + \! \sum_{q=1}^{l_{\mathfrak{s}}
=\varkappa_{nk}} 
\widehat{\nu}_{\mathfrak{s},q}^{\raise-1.0ex\hbox{$\scriptstyle \infty$}}
(n,k)z^{q-1}, \quad \widehat{\nu}_{\mathfrak{s},\varkappa_{nk}}^{
\raise-1.0ex\hbox{$\scriptstyle \infty$}}(n,k) \! \neq \! 0,
\end{equation*}
where the $(n \! - \! 1)K \! + \! k$ coefficients 
$\widehat{\nu}_{1,1}^{\raise-1.0ex\hbox{$\scriptstyle \infty$}}(n,k),
\dotsc,\widehat{\nu}_{1,l_{1}}^{\raise-1.0ex\hbox{$\scriptstyle \infty$}}
(n,k),\dotsc,\widehat{\nu}_{\mathfrak{s}-1,1}^{
\raise-1.0ex\hbox{$\scriptstyle \infty$}}(n,k),\dotsc,
\widehat{\nu}_{\mathfrak{s}-1,l_{\mathfrak{s}-1}}^{
\raise-1.0ex\hbox{$\scriptstyle \infty$}}(n,k),\widehat{\nu}_{
\mathfrak{s},1}^{\raise-1.0ex\hbox{$\scriptstyle \infty$}}(n,k),
\dotsc,
\widehat{\nu}_{\mathfrak{s},\varkappa_{nk}}^{\raise-1.0ex\hbox{$\scriptstyle \infty$}}
(n,k)$, with $\widehat{\nu}_{\mathfrak{s},\varkappa_{nk}}^{
\raise-1.0ex\hbox{$\scriptstyle \infty$}}(n,k) \! \neq \! 0$, satisfy the 
linear inhomogeneous algebraic system of equations~\eqref{eq46}. Via the 
multi-linearity property of the determinant and an application of Cramer's 
Rule to system~\eqref{eq46}, one arrives at, for $n \! \in \! \mathbb{N}$ 
and $k \! \in \! \lbrace 1,2,\dotsc,K \rbrace$ such that $\alpha_{p_{\mathfrak{s}}} 
\! := \! \alpha_{k} \! = \! \infty$, the following---ordered---determinantal 
representation for $\mathcal{X}_{21}(z)$:
\begin{equation*}
\mathcal{X}_{21}(z) \! = \! -\dfrac{2 \pi \mi}{\mathcal{N}_{\infty}^{\sharp}(n,k)} 
\overset{\blacktriangledown}{\Xi}^{\raise-1.0ex\hbox{$\scriptstyle n$}}_{k}(z),
\end{equation*}
where $\mathcal{N}_{\infty}^{\sharp}(n,k)$ is given by Equation~\eqref{eq47}, 
(with abuse of notation) $m_{1} \! = \! (n \! - \! 1)K \! + \! k \! - \! 2$, and
\begin{align*}
\overset{\blacktriangledown}{\Xi}^{\raise-1.0ex\hbox{$\scriptstyle n$}}_{k}
(z) \! := 
\setcounter{MaxMatrixCols}{12}
&\left\vert 

\right\vert \\
=& \, \underbrace{\int_{\mathbb{R}} \int_{\mathbb{R}} \dotsb \int_{\mathbb{R}}}_{
(n-1)K+k-1} \md \widetilde{\mu}(\xi_{0}) \, \md \widetilde{\mu}(\xi_{1}) \, \dotsb 
\, \md \widetilde{\mu}(\xi_{l_{1}}) \, \dotsb \, \md \widetilde{\mu}(\xi_{n_{1}}) \, 
\dotsb \, \md \widetilde{\mu}(\xi_{n_{2}}) \, \md \widetilde{\mu}(\xi_{m_{2}-
\varkappa_{nk}+1}) \, \dotsb \, \md \widetilde{\mu} (\xi_{m_{1}}) \\
\times& \, \dfrac{(-1)^{(n-1)K+k-\varkappa_{nk}}}{\dfrac{1}{\xi_{0}^{0}}
(\xi_{1} \! - \! \alpha_{p_{1}})^{1} \dotsb (\xi_{l_{1}} \! - \! 
\alpha_{p_{1}})^{l_{1}} \dotsb (\xi_{n_{1}} \! - \! \alpha_{p_{\mathfrak{s}
-1}})^{1} \dotsb (\xi_{n_{2}} \! - \! \alpha_{p_{\mathfrak{s}-1}})^{
l_{\mathfrak{s}-1}} \dfrac{1}{(\xi_{m_{2}-\varkappa_{nk}+1})^{1}} \dotsb 
\dfrac{1}{(\xi_{m_{1}})^{\varkappa_{nk}-2}}} \\
\times& \, 
\left\vert

\right\vert.
\end{align*}
For $n \! \in \! \mathbb{N}$ and $k \! \in \! \lbrace 1,2,\dotsc,K \rbrace$ such 
that $\alpha_{p_{\mathfrak{s}}} \! := \! \alpha_{k} \! = \! \infty$, recalling the 
$(n \! - \! 1)K \! + \! k$ linearly independent functions $\psi_{0}(z) \! := \! 
\prod_{m=1}^{\mathfrak{s}-1}(z \! - \! \alpha_{p_{m}})^{l_{m}} \! =: \! 
\sum_{j=0}^{(n-1)K+k-1} \mathfrak{a}^{\sphat}_{j,0}z^{j}$, $\psi_{q(r_{1})}(z) 
\! := \! \lbrace \psi_{0}(z)(z \! - \! \alpha_{p_{r_{1}}})^{-m(r_{1})(1-\delta_{r_{1} 
\mathfrak{s}})}z^{m(r_{1}) \delta_{r_{1} \mathfrak{s}}} \! =: \! \sum_{j=0}^{(n-1)K+k-1} 
\mathfrak{a}^{\sphat}_{j,q(r_{1})}z^{j} \rbrace$, $r_{1} \! = \! 1,\dotsc,\mathfrak{s} 
\! - \! 1,\mathfrak{s}$, $q(r_{1}) \! = \! \sum_{i=1}^{r_{1}-1}l_{i} \! + \! 1,
\sum_{i=1}^{r_{1}-1}l_{i} \! + \! 2,\dotsc,\sum_{i=1}^{r_{1}-1}l_{i} \! + \! 
l_{r_{1}} \! - \! \delta_{r_{1} \mathfrak{s}}$, $m(r_{1}) \! = \! 1,2,\dotsc,
l_{r_{1}} \! - \! \delta_{r_{1} \mathfrak{s}}$, and the $(n \! - \! 1)K \! + \! 
k \! - \! 1$ linearly independent functions $\widehat{\psi}_{0}(z) \! := \! 
\prod_{m=1}^{\mathfrak{s}-1}(z \! - \! \alpha_{p_{m}})^{l_{m}} \! =: \! 
\sum_{j=0}^{(n-1)K+k-2} \mathfrak{a}^{\spadesuit}_{j,0}z^{j}$, 
$\widehat{\psi}_{q(r_{2})}(z) \! := \! \lbrace \widehat{\psi}_{0}(z)
(z \! - \! \alpha_{p_{r_{2}}})^{-m(r_{2})(1-\delta_{r_{2} \mathfrak{s}})}
z^{m(r_{2}) \delta_{r_{2} \mathfrak{s}}} \! =: \! \sum_{j=0}^{(n-1)K+k-2} 
\mathfrak{a}^{\spadesuit}_{j,q(r_{2})}z^{j} \rbrace$, $r_{2} \! = \! 1,
\dotsc,\mathfrak{s} \! - \! 1,\mathfrak{s}$, $q(r_{2}) \! = \! 
\sum_{i=1}^{r_{2}-1}l_{i} \! + \! 1,\sum_{i=1}^{r_{2}-1}l_{i} \! + \! 
2,\dotsc,\sum_{i=1}^{r_{2}-1}l_{i} \! + \! l_{r_{2}} \! - \! 2 \delta_{r_{2} 
\mathfrak{s}}$, $m(r_{2}) \! = \! 1,2,\dotsc,l_{r_{2}} \! - \! 2 \delta_{r_{2} 
\mathfrak{s}}$, one proceeds, via the latter determinantal expression, with 
the analysis of 
$\overset{\blacktriangledown}{\Xi}^{\raise-1.0ex\hbox{$\scriptstyle n$}}_{k}
(z) \! := \! (-1)^{(n-1)K+k-\varkappa_{nk}} \overset{\blacktriangle}{\Xi}^{
\raise-1.0ex\hbox{$\scriptstyle n$}}_{k}(z)$:
\begin{align*}
\overset{\blacktriangle}{\Xi}^{\raise-1.0ex\hbox{$\scriptstyle n$}}_{k}(z) =& 
\, \underbrace{\int_{\mathbb{R}} \int_{\mathbb{R}} \dotsb \int_{\mathbb{R}}}_{
(n-1)K+k-1} \md \widetilde{\mu}(\xi_{0}) \, \md \widetilde{\mu}(\xi_{1}) \, 
\dotsb \, \md \widetilde{\mu}(\xi_{l_{1}}) \, \dotsb \, \md \widetilde{\mu}
(\xi_{n_{1}}) \, \dotsb \, \md \widetilde{\mu}(\xi_{n_{2}}) \, \md \widetilde{\mu}
(\xi_{m_{2}-\varkappa_{nk}+1}) \, \dotsb \, \md \widetilde{\mu}(\xi_{m_{1}}) \\
\times& \, \dfrac{1}{\dfrac{1}{\xi_{0}^{0}}(\xi_{1} \! - \! \alpha_{p_{1}})^{
1} \dotsb (\xi_{l_{1}} \! - \! \alpha_{p_{1}})^{l_{1}} \dotsb (\xi_{n_{1}} \! 
- \! \alpha_{p_{\mathfrak{s}-1}})^{1} \dotsb (\xi_{n_{2}} \! - \! \alpha_{p_{
\mathfrak{s}-1}})^{l_{\mathfrak{s}-1}} \dfrac{1}{(\xi_{m_{2}-\varkappa_{nk}
+1})^{1}} \dotsb \dfrac{1}{\xi_{m_{1}}^{\varkappa_{nk}-2}}} \\
\times& \, \dfrac{(\psi_{0}(z))^{-1}}{\psi_{0}(\xi_{0}) \psi_{0}(\xi_{1}) 
\dotsb \psi_{0}(\xi_{l_{1}}) \dotsb \psi_{0}(\xi_{n_{1}}) \dotsb \psi_{0}
(\xi_{n_{2}}) \psi_{0}(\xi_{m_{2}-\varkappa_{nk}+1}) \dotsb \psi_{0}(\xi_{
m_{1}})} \\
\times& \, 
\underbrace{\left\vert 

\right\vert}_{=: \, \mathbb{G}^{\blacktriangle}(\xi_{0},\xi_{1},\dotsc,\xi_{
l_{1}},\dotsc,\xi_{n_{1}},\dotsc,\xi_{n_{2}},\dotsc,\xi_{(n-1)K+k-2};z)} \\
=& \, \dfrac{(\psi_{0}(z))^{-1}}{(m_{1} \! + \! 1)!} \sum_{\pmb{\sigma} \in 
\mathfrak{S}_{m_{1}+1}} \underbrace{\int_{\mathbb{R}} \int_{\mathbb{R}} 
\dotsb \int_{\mathbb{R}}}_{(n-1)K+k-1} \md \widetilde{\mu}(\xi_{\sigma (0)}) 
\, \md \widetilde{\mu}(\xi_{\sigma (1)}) \, \dotsb \, \md \widetilde{\mu}
(\xi_{\sigma (l_{1})}) \, \dotsb \, \md \widetilde{\mu}(\xi_{\sigma (n_{1})}) 
\, \dotsb \, \md \widetilde{\mu}(\xi_{\sigma (n_{2})}) \\
\times& \, \dfrac{\dotsb \, \md \widetilde{\mu}(\xi_{\sigma (m_{2}-
\varkappa_{nk}+1)}) \, \dotsb \, \md \widetilde{\mu}(\xi_{\sigma (m_{1})})}{
\dfrac{1}{\xi_{\sigma (0)}^{0}}(\xi_{\sigma (1)} \! - \! \alpha_{p_{1}})^{1} \dotsb 
(\xi_{\sigma (l_{1})} \! - \! \alpha_{p_{1}})^{l_{1}} \dotsb (\xi_{\sigma (n_{1})} 
\! - \! \alpha_{p_{\mathfrak{s}-1}})^{1} \dotsb (\xi_{\sigma (n_{2})} \! - \! 
\alpha_{p_{\mathfrak{s}-1}})^{l_{\mathfrak{s}-1}} \dfrac{1}{(\xi_{\sigma 
(m_{2}-\varkappa_{nk}+1)})^{1}} \dotsb \dfrac{1}{\xi_{\sigma (m_{1})}^{
\varkappa_{nk}-2}}} \\
\times& \, \dfrac{1}{\psi_{0}(\xi_{\sigma (0)}) \psi_{0}(\xi_{\sigma (1)}) 
\dotsb \psi_{0}(\xi_{\sigma (l_{1})}) \dotsb \psi_{0}(\xi_{\sigma (n_{1})}) 
\dotsb \psi_{0}(\xi_{\sigma (n_{2})}) \psi_{0}(\xi_{\sigma (m_{2}-
\varkappa_{nk}+1)}) \dotsb \psi_{0}(\xi_{\sigma (m_{1})})} \\
\times& \, \mathbb{G}^{\blacktriangle}(\xi_{\sigma (0)},\xi_{\sigma (1)},
\dotsc,\xi_{\sigma (l_{1})},\dotsc,\xi_{\sigma (n_{1})},\dotsc,
\xi_{\sigma (n_{2})},\dotsc,\xi_{\sigma ((n-1)K+k-2)};z) \\
=& \, \dfrac{(\psi_{0}(z))^{-1}}{(m_{1} \! + \! 1)!} \underbrace{\int_{\mathbb{R}} 
\int_{\mathbb{R}} \dotsb \int_{\mathbb{R}}}_{(n-1)K+k-1} \md \widetilde{\mu}
(\xi_{0}) \, \md \widetilde{\mu}(\xi_{1}) \, \dotsb \, \md \widetilde{\mu}(\xi_{l_{1}}) 
\, \dotsb \, \md \widetilde{\mu}(\xi_{n_{1}}) \, \dotsb \, \md \widetilde{\mu}
(\xi_{n_{2}}) \, \md \widetilde{\mu}(\xi_{m_{2}-\varkappa_{nk}+1}) \, \dotsb \\
\times& \, \dotsb \, \md \widetilde{\mu}(\xi_{m_{1}}) \, \mathbb{G}^{\blacktriangle}
(\xi_{0},\xi_{1},\dotsc,\xi_{l_{1}},\dotsc,\xi_{n_{1}},\dotsc,\xi_{n_{2}},\dotsc,
\xi_{(n-1)K+k-2};z) \\
\times& \, \sum_{\pmb{\sigma} \in \mathfrak{S}_{m_{1}+1}} \operatorname{sgn}
(\pmb{\pmb{\sigma}}) \tfrac{1}{\dfrac{1}{\xi_{\sigma (0)}^{0}}
(\xi_{\sigma (1)}-\alpha_{p_{1}})^{1} \dotsb (\xi_{\sigma (l_{1})}-
\alpha_{p_{1}})^{l_{1}} \dotsb (\xi_{\sigma (n_{1})}-\alpha_{p_{\mathfrak{s}
-1}})^{1} \dotsb (\xi_{\sigma (n_{2})}-\alpha_{p_{\mathfrak{s}-1}})^{
l_{\mathfrak{s}-1}} \dfrac{1}{(\xi_{\sigma (m_{2}-\varkappa_{nk}+1)})^{1}} 
\dotsb \dfrac{1}{\xi_{\sigma (m_{1})}^{\varkappa_{nk}-2}}} \\
\times& \, \dfrac{1}{\psi_{0}(\xi_{\sigma (0)}) \psi_{0}(\xi_{\sigma (1)}) 
\dotsb \psi_{0}(\xi_{\sigma (l_{1})}) \dotsb \psi_{0}(\xi_{\sigma (n_{1})}) 
\dotsb \psi_{0}(\xi_{\sigma (n_{2})}) \psi_{0}(\xi_{\sigma (m_{2}-\varkappa_{
nk}+1)}) \dotsb \psi_{0}(\xi_{\sigma (m_{1})})} \\
=& \, \dfrac{(\psi_{0}(z))^{-1}}{(m_{1} \! + \! 1)!} \underbrace{\int_{\mathbb{R}} 
\int_{\mathbb{R}} \dotsb \int_{\mathbb{R}}}_{(n-1)K+k-1} \md \widetilde{\mu}
(\xi_{0}) \, \md \widetilde{\mu}(\xi_{1}) \, \dotsb \, \md \widetilde{\mu}(\xi_{l_{1}}) 
\, \dotsb \, \md \widetilde{\mu}(\xi_{n_{1}}) \, \dotsb \, \md \widetilde{\mu}
(\xi_{n_{2}}) \, \md \widetilde{\mu}(\xi_{m_{2}-\varkappa_{nk}+1}) \, \dotsb \\
\times& \, \dotsb \, \md \widetilde{\mu}(\xi_{m_{1}}) \, \dfrac{\mathbb{G}^{
\blacktriangle}(\xi_{0},\xi_{1},\dotsc,\xi_{l_{1}},\dotsc,\xi_{n_{1}},\dotsc,
\xi_{n_{2}},\dotsc,\xi_{(n-1)K+k-2};z)}{\psi_{0}(\xi_{0}) \psi_{0}(\xi_{1}) 
\dotsb \psi_{0}(\xi_{l_{1}}) \dotsb \psi_{0}(\xi_{n_{1}}) \dotsb \psi_{0}
(\xi_{n_{2}}) \psi_{0}(\xi_{m_{2}-\varkappa_{nk}+1}) \dotsb \psi_{0}(\xi_{
m_{1}})} \\
\times& \, 
\left\vert

\right\vert \\
=& \, \dfrac{(\psi_{0}(z))^{-1}}{(m_{1} \! + \! 1)!} \underbrace{\int_{\mathbb{R}} 
\int_{\mathbb{R}} \dotsb \int_{\mathbb{R}}}_{(n-1)K+k-1} \md \widetilde{\mu}
(\xi_{0}) \, \md \widetilde{\mu}(\xi_{1}) \, \dotsb \, \md \widetilde{\mu}(\xi_{l_{1}}) 
\, \dotsb \, \md \widetilde{\mu}(\xi_{n_{1}}) \, \dotsb \, \md \widetilde{\mu}
(\xi_{n_{2}}) \, \md \widetilde{\mu}(\xi_{m_{2}-\varkappa_{nk}+1}) \, \dotsb \\
\times& \, \dotsb \, \md \widetilde{\mu}(\xi_{m_{1}}) \, \dfrac{\mathbb{G}^{
\blacktriangle}(\xi_{0},\xi_{1},\dotsc,\xi_{l_{1}},\dotsc,\xi_{n_{1}},\dotsc,
\xi_{n_{2}},\dotsc,\xi_{(n-1)K+k-2};z)}{\psi_{0}(\xi_{0}) \psi_{0}(\xi_{1}) 
\dotsb \psi_{0}(\xi_{l_{1}}) \dotsb \psi_{0}(\xi_{n_{1}}) \dotsb \psi_{0}
(\xi_{n_{2}}) \psi_{0}(\xi_{m_{2}-\varkappa_{nk}+1}) \dotsb \psi_{0}(\xi_{
m_{1}})} \\
\times& \, \dfrac{1}{\widehat{\psi}_{0}(\xi_{0}) \widehat{\psi}_{0}(\xi_{1}) 
\dotsb \widehat{\psi}_{0}(\xi_{l_{1}}) \dotsb \widehat{\psi}_{0}(\xi_{n_{1}}) 
\dotsb \widehat{\psi}_{0}(\xi_{n_{2}}) \widehat{\psi}_{0}(\xi_{m_{2}-
\varkappa_{nk}+1}) \dotsb \widehat{\psi}_{0}(\xi_{m_{1}})} \\
\times& \, 
\underbrace{
\left\vert 

\right\vert}_{=: \, \mathbb{G}^{\triangleright}(\xi_{0},\xi_{1},\dotsc,\xi_{
l_{1}},\dotsc,\xi_{n_{1}},\dotsc,\xi_{n_{2}},\dotsc,\xi_{(n-1)K+k-2})};
\end{align*}
but, noting the determinantal factorisation
\begin{equation*}
\mathbb{G}^{\triangleright}(\xi_{0},\xi_{1},\dotsc,\xi_{l_{1}},\dotsc,
\xi_{n_{1}},\dotsc,\xi_{n_{2}},\dotsc,\xi_{(n-1)K+k-2}) \! := \! 
\overset{\infty}{\mathcal{V}}_{3}(\xi_{0},\xi_{1},\dotsc,\xi_{(n-1)K+k-2}) 
\mathbb{D}^{\spadesuit},
\end{equation*}
where
\begin{equation*}
\overset{\infty}{\mathcal{V}}_{3}(\xi_{0},\xi_{1},\dotsc,\xi_{(n-1)K+k-2}) 
= \left\lvert

\right), \label{eq76}
\end{equation}
with
\begin{equation*}
(\mathscr{A}^{\spadesuit}(0))_{1j} \! = \! 
\begin{cases}
\mathfrak{a}_{j-1}^{\spadesuit}(\vec{\bm{\alpha}}), &\text{$j \! = \! 1,2,
\dotsc,(n \! - \! 1)K \! + \! k \! - \! \varkappa_{nk} \! + \! 1$,} \\
0, &\text{$j \! = \! (n \! - \! 1)K \! + \! k \! - \! \varkappa_{nk} \! 
+ \! 2,\dotsc,(n \! - \! 1)K \! + \! k \! - \! 1$,}
\end{cases}
\end{equation*}
for $r \! = \! 1,2,\dotsc,\mathfrak{s} \! - \! 1$, with $l_{r} \! = \! 
\varkappa_{nk \tilde{k}_{r}}$ and $\sum_{m=1}^{\mathfrak{s}-1}l_{m} \! 
= \! (n \! - \! 1)K \! + \! k \! - \! \varkappa_{nk}$,\footnote{Note the 
convention $\sum_{m=1}^{0} \pmb{\ast} \! := \! 0$.}
\begin{align*}
i(r)=& \, 2 \! + \! \sum_{m=1}^{r-1}l_{m},3 \! + \! \sum_{m=1}^{r-1}l_{m},
\dotsc,1 \! + \! l_{r} \! + \! \sum_{m=1}^{r-1}l_{m}, \qquad \qquad q(i(r),r) 
\! = \! 1,2,\dotsc,l_{r}, \\
(\mathscr{A}^{\spadesuit}(r))_{i(r)j(r)}=& \, 
\begin{cases}
\dfrac{(-1)^{q(i(r),r)}}{\prod_{m=0}^{q(i(r),r)-1}(l_{r} \! - \! m)} 
\left(\dfrac{\partial}{\partial \alpha_{p_{r}}} \right)^{q(i(r),r)} 
\mathfrak{a}_{j(r)-1}^{\spadesuit}(\vec{\bm{\alpha}}), &\text{$j(r) \! = \! 
1,2,\dotsc,(n \! - \! 1)K \! + \! k \! - \! \varkappa_{nk} \! - \! q(i(r),r) 
\! + \! 1$,} \\
0, &\text{$j(r) \! = \! (n \! - \! 1)K \! + \! k \! - \! \varkappa_{nk} 
\! - \! q(i(r),r) \! + \! 2,\dotsc,(n \! - \! 1)K \! + \! k \! - \! 1$,}
\end{cases}
\end{align*}
and, for $r \! = \! \mathfrak{s}$, with $l_{\mathfrak{s}} \! = \! \varkappa_{nk}$,
\begin{align*}
i(\mathfrak{s})=& \, (n \! - \! 1)K \! + \! k \! - \! \varkappa_{nk} \! + \! 
2,(n \! - \! 1)K \! + \! k \! - \! \varkappa_{nk} \! + \! 3,\dotsc,(n \! - \! 
1)K \! + \! k \! - \! 1, \quad \quad q(i(\mathfrak{s}),\mathfrak{s}) \! = \! 
1,2,\dotsc,\varkappa_{nk} \! - \! 2, \\
(\mathscr{A}^{\spadesuit}(\mathfrak{s}))_{i(\mathfrak{s})j(\mathfrak{s})}=& \, 
\begin{cases}
0, &\text{$j(\mathfrak{s}) \! = \! 1,\dotsc,q(i(\mathfrak{s}),\mathfrak{s})$,} 
\\
\mathfrak{a}_{j(\mathfrak{s})-q(i(\mathfrak{s}),\mathfrak{s})-1}^{\spadesuit}
(\vec{\bm{\alpha}}), &\text{$j(\mathfrak{s}) \! = \! q(i(\mathfrak{s}),
\mathfrak{s}) \! + \! 1,\dotsc,i(\mathfrak{s})$,} \\
0, &\text{$j(\mathfrak{s}) \! = \! i(\mathfrak{s}) \! + \! 1,\dotsc,
(n \! - \! 1)K \! + \! k \! - \! 1$,}
\end{cases}
\end{align*}
where, for $\tilde{m}_{1} \! = \! 0,1,\dotsc,(n \! - \! 1)K \! + \! k \! - \! 
\varkappa_{nk}$,
\begin{equation*}
\mathfrak{a}_{\tilde{m}_{1}}^{\spadesuit}(\vec{\bm{\alpha}}) := 
\mathlarger{\sum_{\underset{\underset{\sum_{m=1}^{\mathfrak{s}-1}
i_{m}=(n-1)K+k-\varkappa_{nk}-\tilde{m}_{1}}{p \in \lbrace 1,2,\dotsc,
\mathfrak{s}-1 \rbrace}}{i_{p}=0,1,\dotsc,l_{p}}}}(-1)^{(n-1)K+k-
\varkappa_{nk}-\tilde{m}_{1}} \prod_{j=1}^{\mathfrak{s}-1} \binom{l_{j}}{i_{j}} 
\prod_{m=1}^{\mathfrak{s}-1}(\alpha_{p_{m}})^{i_{m}},
\end{equation*}
it follows that, for $n \! \in \! \mathbb{N}$ and $k \! \in \! \lbrace 1,2,
\dotsc,K \rbrace$ such that $\alpha_{p_{\mathfrak{s}}} \! := \! \alpha_{k} 
\! = \! \infty$,
\begin{align*}
\overset{\blacktriangle}{\Xi}^{\raise-1.0ex\hbox{$\scriptstyle n$}}_{k}(z)
=& \, \dfrac{\mathbb{D}^{\spadesuit}(\psi_{0}(z))^{-1}}{(m_{1} \! + \! 1)!} 
\underbrace{\int_{\mathbb{R}} \int_{\mathbb{R}} \dotsb \int_{\mathbb{R}}}_{
(n-1)K+k-1} \md \widetilde{\mu}(\xi_{0}) \, \md \widetilde{\mu}(\xi_{1}) \, 
\dotsb \, \md \widetilde{\mu}(\xi_{l_{1}}) \, \dotsb \, \md \widetilde{\mu}
(\xi_{n_{1}}) \, \dotsb \, \md \widetilde{\mu}(\xi_{n_{2}}) \, \md \widetilde{\mu}
(\xi_{m_{2}-\varkappa_{nk}+1}) \, \dotsb \\
\times& \, \dotsb \, \md \widetilde{\mu}(\xi_{m_{1}}) \, 
\dfrac{\prod_{\underset{j<i}{i,j=0}}^{(n-1)K+k-2}
(\xi_{i} \! - \! \xi_{j})}{\psi_{0}(\xi_{0}) \psi_{0}(\xi_{1}) \dotsb \psi_{0}
(\xi_{l_{1}}) \dotsb \psi_{0}(\xi_{n_{1}}) \dotsb \psi_{0}(\xi_{n_{2}}) 
\psi_{0}(\xi_{m_{2}-\varkappa_{nk}+1}) \dotsb \psi_{0}(\xi_{m_{1}})} \\
\times& \, \dfrac{1}{\widehat{\psi}_{0}(\xi_{0}) \widehat{\psi}_{0}(\xi_{1}) 
\dotsb \widehat{\psi}_{0}(\xi_{l_{1}}) \dotsb \widehat{\psi}_{0}(\xi_{n_{1}}) 
\dotsb \widehat{\psi}_{0}(\xi_{n_{2}}) \widehat{\psi}_{0}(\xi_{m_{2}-
\varkappa_{nk}+1}) \dotsb \widehat{\psi}_{0}(\xi_{m_{1}})} \\
\times& \, 
\underbrace{\left\vert

\right\vert}_{= \, \mathbb{D}^{\sphat} \prod_{\underset{j<i}{i,j=0}}^{(n-1)
K+k-2}(\xi_{i}-\xi_{j}) \prod_{m=0}^{(n-1)K+k-2}(z-\xi_{m}) \quad 
\text{(cf. proof of Lemma~\ref{lem2.1}, case~\pmb{(1)})}} \\
=& \, \dfrac{\mathbb{D}^{\spadesuit} \mathbb{D}^{\sphat} \, 
(\psi_{0}(z))^{-1}}{((n \! - \! 1)K \! + \! k \! - \! 1)!} \underbrace{
\int_{\mathbb{R}} \int_{\mathbb{R}} \dotsb \int_{\mathbb{R}}}_{(n-1)K+k-1} 
\md \widetilde{\mu}(\xi_{0}) \, \md \widetilde{\mu}(\xi_{1}) \, \dotsb \, \md 
\widetilde{\mu}(\xi_{l_{1}}) \, \dotsb \, \md \widetilde{\mu}(\xi_{n_{1}}) \, 
\dotsb \, \md \widetilde{\mu}(\xi_{n_{2}}) \, \md \widetilde{\mu}
(\xi_{m_{2}-\varkappa_{nk}+1}) \, \dotsb \\
\times& \, \dotsb \, \md \widetilde{\mu}(\xi_{m_{1}}) \, \dfrac{(\widehat{\psi}_{0}
(\xi_{0}) \widehat{\psi}_{0}(\xi_{1}) \dotsb \widehat{\psi}_{0}(\xi_{l_{1}}) 
\dotsb \widehat{\psi}_{0}(\xi_{n_{1}}) \dotsb \widehat{\psi}_{0}(\xi_{n_{2}}) 
\widehat{\psi}_{0}(\xi_{m_{2}-\varkappa_{nk}+1}) \dotsb \widehat{\psi}_{0}
(\xi_{m_{1}}))^{-1}}{\psi_{0}(\xi_{0}) \psi_{0}(\xi_{1}) \dotsb \psi_{0}
(\xi_{l_{1}}) \dotsb \psi_{0}(\xi_{n_{1}}) \dotsb \psi_{0}(\xi_{n_{2}}) 
\psi_{0}(\xi_{m_{2}-\varkappa_{nk}+1}) \dotsb \psi_{0}(\xi_{m_{1}})} \\
\times& \, \prod_{\substack{i,j=0\\j<i}}^{(n-1)K+k-2}(\xi_{i} \! - \! 
\xi_{j})^{2} \, \prod_{m=0}^{(n-1)K+k-2}(z \! - \! \xi_{m});
\end{align*}
but, recalling that, for $n \! \in \! \mathbb{N}$ and $k \! \in \! \lbrace 
1,2,\dotsc,K \rbrace$ such that $\alpha_{p_{\mathfrak{s}}} \! := \! 
\alpha_{k} \! = \! \infty$, $\psi_{0}(z) \! = \! \widehat{\psi}_{0}(z)$, and 
$\overset{\blacktriangledown}{\Xi}^{\raise-1.0ex\hbox{$\scriptstyle n$}}_{k}
(z) \! := \! (-1)^{(n-1)K+k-\varkappa_{nk}} \linebreak[4] 
\pmb{\cdot} 
\overset{\blacktriangle}{\Xi}^{\raise-1.0ex\hbox{$\scriptstyle n$}}_{k}(z)$, 
one arrives at
\begin{align*}
\overset{\blacktriangledown}{\Xi}^{\raise-1.0ex\hbox{$\scriptstyle n$}}_{k}
(z) =& \, \dfrac{(-1)^{(n-1)K+k-\varkappa_{nk}} \mathbb{D}^{\spadesuit} 
\mathbb{D}^{\sphat} \, (\psi_{0}(z))^{-1}}{((n \! - \! 1)K \! + \! k \! - \! 1)!} 
\underbrace{\int_{\mathbb{R}} \int_{\mathbb{R}} \dotsb \int_{\mathbb{R}}}_{
(n-1)K+k-1} \md \widetilde{\mu}(\xi_{0}) \, \md \widetilde{\mu}(\xi_{1}) \, 
\dotsb \, \md \widetilde{\mu}(\xi_{l_{1}}) \, \dotsb \, \md \widetilde{\mu}
(\xi_{n_{1}}) \, \dotsb \, \md \widetilde{\mu}(\xi_{n_{2}}) \\
\times& \, \md \widetilde{\mu}(\xi_{(n-1)K+k-\varkappa_{nk}+1}) \, \dotsb 
\, \md \widetilde{\mu}(\xi_{(n-1)K+k-2}) \left(\prod_{l=0}^{(n-1)K+k-2} 
\psi_{0}(\xi_{l}) \right)^{-2} \prod_{\substack{i,j=0\\j<i}}^{(n-1)K+k-2}
(\xi_{i} \! - \! \xi_{j})^{2} \, \prod_{m=0}^{(n-1)K+k-2}(z \! - \! \xi_{m}) 
\, \, \Rightarrow \\
\overset{\blacktriangledown}{\Xi}^{\raise-1.0ex\hbox{$\scriptstyle n$}}_{k}
(z) =& \, \dfrac{(-1)^{(n-1)K+k-\varkappa_{nk}} \mathbb{D}^{\spadesuit} 
\mathbb{D}^{\sphat}}{((n \! - \! 1)K \! + \! k \! - \! 1)!} \left(
\prod_{q=1}^{\mathfrak{s}-1}(z \! - \! \alpha_{p_{q}})^{\varkappa_{nk 
\tilde{k}_{q}}} \right)^{-1} \underbrace{\int_{\mathbb{R}} \int_{\mathbb{R}} 
\dotsb \int_{\mathbb{R}}}_{(n-1)K+k-1} \md \widetilde{\mu}(\xi_{0}) \, \md 
\widetilde{\mu}(\xi_{1}) \, \dotsb \, \md \widetilde{\mu}(\xi_{(n-1)K+k-2}) \\
\times& \, \left(\prod_{m=0}^{(n-1)K+k-2} \prod_{q=1}^{\mathfrak{s}-1}
(\xi_{m} \! - \! \alpha_{p_{q}})^{\varkappa_{nk \tilde{k}_{q}}} \right)^{-2} 
\prod_{\substack{i,j=0\\j<i}}^{(n-1)K+k-2}(\xi_{i} \! - \! \xi_{j})^{2} \, 
\prod_{l=0}^{(n-1)K+k-2}(z \! - \! \xi_{l}).
\end{align*}
Recalling that, for $n \! \in \! \mathbb{N}$ and $k \! \in \! \lbrace 1,2,
\dotsc,K \rbrace$ such that $\alpha_{p_{\mathfrak{s}}} \! := \! \alpha_{k} 
\! = \! \infty$, $\mathcal{X}_{21}(z) \! = \! -2 \pi \mi 
\overset{\blacktriangledown}{\Xi}^{\raise-1.0ex\hbox{$\scriptstyle n$}}_{k}
(z)/\mathcal{N}_{\infty}^{\sharp}(n,k)$, where 
$\overset{\blacktriangledown}{\Xi}^{\raise-1.0ex\hbox{$\scriptstyle n$}}_{k}
(z)$ is given directly above, and $\mathcal{N}_{\infty}^{\sharp}(n,k)$ is given 
by Equation~\eqref{eq47}, one arrives at the integral representation for 
$\mathcal{X}_{21}(z)$ given in Equation~\eqref{eintreppinf2}.

\pmb{(2)} For $n \! \in \! \mathbb{N}$ and $k \! \in \! \lbrace 1,2,\dotsc,
K \rbrace$ such that $\alpha_{p_{\mathfrak{s}}} \! := \! \alpha_{k} \! \neq 
\! \infty$, recall {}from the proof of Lemma~\ref{lem2.1} (cf. case~\textbf{(2)}) 
the ordered disjoint partition for $\lbrace \alpha_{1},\alpha_{2},\dotsc,
\alpha_{K} \rbrace \cup \lbrace \alpha_{1},\alpha_{2},\dotsc,\alpha_{K} 
\rbrace \cup \dotsb \cup \lbrace \alpha_{1},\alpha_{2},\dotsc,\alpha_{k} 
\rbrace$ and the associated formula for the corresponding monic MPC ORF,
\begin{align*}
\pmb{\pi}^{n}_{k}(z) \! = \! \dfrac{\mathcal{X}_{11}(z)}{z \! - \! \alpha_{k}}
=& \, \widetilde{\phi}^{\raise-1.0ex\hbox{$\scriptstyle f$}}_{0}(n,k) \! + \! 
\sum_{m=1}^{\mathfrak{s}-2} \sum_{q=1}^{l_{m}=\varkappa_{nk \tilde{k}_{m}}} 
\dfrac{\widetilde{\nu}_{m,q}^{\raise-1.0ex\hbox{$\scriptstyle f$}}(n,k)}{
(z \! - \! \alpha_{p_{m}})^{q}} \! + \! \sum_{q=1}^{l_{\mathfrak{s}-1}=
\varkappa^{\infty}_{nk \tilde{k}_{\mathfrak{s}-1}}} \widetilde{\nu}_{
\mathfrak{s}-1,q}^{\raise-1.0ex\hbox{$\scriptstyle f$}}(n,k)z^{q} \! + \! 
\sum_{q=1}^{l_{\mathfrak{s}}-1=\varkappa_{nk}-1} \dfrac{\widetilde{\nu}_{
\mathfrak{s},q}^{\raise-1.0ex\hbox{$\scriptstyle f$}}(n,k)}{(z \! - \! 
\alpha_{k})^{q}} \! + \! \dfrac{1}{(z \! - \! \alpha_{k})^{\varkappa_{nk}}},
\end{align*}
where the $(n \! - \! 1)K \! + \! k$ coefficients 
$\widetilde{\phi}_{0}^{\raise-1.0ex\hbox{$\scriptstyle f$}}(n,k),
\widetilde{\nu}_{1,1}^{\raise-1.0ex\hbox{$\scriptstyle f$}}(n,k),\dotsc,
\widetilde{\nu}_{1,l_{1}}^{\raise-1.0ex\hbox{$\scriptstyle f$}}(n,k),\dotsc,
\widetilde{\nu}_{\mathfrak{s}-1,1}^{\raise-1.0ex\hbox{$\scriptstyle f$}}
(n,k),\dotsc,\widetilde{\nu}_{\mathfrak{s}-1,l_{\mathfrak{s}-1}}^{
\raise-1.0ex\hbox{$\scriptstyle f$}}(n,k),
\widetilde{\nu}_{\mathfrak{s},1}^{\raise-1.0ex\hbox{$\scriptstyle f$}}
(n,k),\linebreak[4] 
\dotsc,\widetilde{\nu}_{\mathfrak{s},l_{\mathfrak{s}}-1}^{
\raise-1.0ex\hbox{$\scriptstyle f$}}(n,k)$, $(\mu_{n,\varkappa_{nk}}^{f}
(n,k))^{-2}$, with $\mu_{n,\varkappa_{nk}}^{f}(n,k)$ being the associated 
norming constant (see Corollary~\ref{cor2.1}, Equation~\eqref{eqnmctfin1}, 
below), satisfy the linear inhomogeneous algebraic system~\eqref{eq61}. 
Via the multi-linearity property of the determinant and an application of 
Cramer's Rule to system~\eqref{eq61}, one arrives at, for $n \! \in \! 
\mathbb{N}$ and $k \! \in \! \lbrace 1,2,\dotsc,K \rbrace$ such that 
$\alpha_{p_{\mathfrak{s}}} \! := \! \alpha_{k} \! \neq \! \infty$, the 
following---ordered---determinantal representation for the corresponding 
monic MPC ORF:
\begin{equation*}
\pmb{\pi}^{n}_{k}(z) \! = \! \dfrac{\mathcal{X}_{11}(z)}{z \! - \! \alpha_{k}} 
\! = \! \hat{\mathfrak{c}}_{D_{f}}^{-1} 
\overset{\scriptscriptstyle f}{\Xi}_{k}^{\raise-1.0ex\hbox{$\scriptstyle n$}}(z),
\end{equation*}
where $\hat{\mathfrak{c}}_{D_{f}}$ is given by Equation~\eqref{eq65}, 
and, with (abuse of notation)
\begin{equation*}
n_{1} \! = \! l_{1} \! + \! \dotsb \! + \! l_{\mathfrak{s}-2} \! + \! 1, 
\qquad n_{2} \! = \! (n \! - \! 1)K \! + \! k \! - \! \varkappa_{nk}, \qquad 
m_{1} \! = \! (n \! - \! 1)K \! + \! k \! - \! 1, \qquad \text{and} \qquad 
m_{2} \! = \! (n \! - \! 1)K \! + \! k,
\end{equation*}
\begin{align*}
\overset{\scriptscriptstyle f}{\Xi}_{k}^{\raise-1.0ex\hbox{$\scriptstyle n$}}
(z) := 
\setcounter{MaxMatrixCols}{12}
&\left\vert

\right\vert \\
=& \, \underbrace{\int_{\mathbb{R}} \int_{\mathbb{R}} \dotsb \int_{\mathbb{R}}
}_{(n-1)K+k} \md \widetilde{\mu}(\xi_{0}) \, \md \widetilde{\mu}(\xi_{1}) \, 
\dotsb \, \md \widetilde{\mu}(\xi_{l_{1}}) \, \dotsb \, \md \widetilde{\mu}
(\xi_{n_{1}}) \, \dotsb \, \md \widetilde{\mu}(\xi_{n_{2}}) \, \md \widetilde{\mu}
(\xi_{m_{2}-\varkappa_{nk}+1}) \, \dotsb \, \md \widetilde{\mu}(\xi_{m_{1}}) \\
\times& \, \dfrac{1}{(\xi_{0} \! - \! \alpha_{k})^{0}(\xi_{1} \! - \! 
\alpha_{p_{1}})^{1} \dotsb (\xi_{l_{1}} \! - \! \alpha_{p_{1}})^{l_{1}} 
\dotsb \dfrac{1}{(\xi_{n_{1}})^{1}} \dotsb \dfrac{1}{(\xi_{n_{2}})^{
l_{\mathfrak{s}-1}}}(\xi_{m_{2}-\varkappa_{nk}+1} \! - \! \alpha_{k})^{1} 
\dotsb (\xi_{m_{1}} \! - \! \alpha_{k})^{\varkappa_{nk}-1}} \\
\times& \, 
\left\vert

\right\vert.
\end{align*}
For $n \! \in \! \mathbb{N}$ and $k \! \in \! \lbrace 1,2,\dotsc,K \rbrace$ 
such that $\alpha_{p_{\mathfrak{s}}} \! := \! \alpha_{k} \! \neq \! \infty$, 
recalling {}from the proof of Lemma~\ref{lem2.1} (cf. case~\textbf{(2)}) the 
$(n \! - \! 1)K \! + \! k \! + \! 1$ linearly independent functions $\phi_{0}
(z) \! := \! \prod_{m=1}^{\mathfrak{s}-2}(z \! - \! \alpha_{p_{m}})^{l_{m}}
(z \! - \! \alpha_{k})^{\varkappa_{nk}} \! =: \! \sum_{j=0}^{(n-1)K+k} 
\mathfrak{a}_{j,0}z^{j}$, $\phi_{q(r_{1})}(z) \! := \! \lbrace \phi_{0}(z)
(z \! - \! \alpha_{p_{r_{1}}})^{-m(r_{1})(1-\delta_{r_{1} \mathfrak{s}-1})}
z^{m(r_{1}) \delta_{r_{1} \mathfrak{s}-1}} \! =: \! \sum_{j=0}^{(n-1)K+k} 
\mathfrak{a}_{j,q(r_{1})}z^{j} \rbrace$, $r_{1} \! = \! 1,\dotsc,\mathfrak{s} 
\! - \! 1,\mathfrak{s}$, $q(r_{1}) \! = \! \sum_{i=1}^{r_{1}-1}l_{i} \! + \! 1,
\sum_{i=1}^{r_{1}-1}l_{i} \! + \! 2,\dotsc,\sum_{i=1}^{r_{1}-1}l_{i} \! + \! 
l_{r_{1}}$, $m(r_{1}) \! = \! 1,2,\dotsc,l_{r_{1}}$, and the $(n \! - \! 1)K \! 
+ \! k$ linearly independent functions $\tilde{\phi}_{0}(z) \! := \! \prod_{m
=1}^{\mathfrak{s}-2}(z \! - \! \alpha_{p_{m}})^{l_{m}}(z \! - \! \alpha_{k})^{
\varkappa_{nk}-1} \! =: \! \sum_{j=0}^{(n-1)K+k-1} \mathfrak{a}^{
\sharp}_{j,0}z^{j}$, $\tilde{\phi}_{q(r_{2})}(z) \! := \! \lbrace \tilde{\phi}_{0}
(z)(z \! - \! \alpha_{p_{r_{2}}})^{-m(r_{2})(1-\delta_{r_{2} \mathfrak{s}-1})}
z^{m(r_{2}) \delta_{r_{2} \mathfrak{s}-1}} \! =: \! \sum_{j=0}^{(n-1)K+k-1} 
\mathfrak{a}^{\sharp}_{j,q(r_{2})}z^{j} \rbrace$, $r_{2} \! = \! 1,\dotsc,
\mathfrak{s} \! - \! 1,\mathfrak{s}$, $q(r_{2}) \! = \! \sum_{i=1}^{r_{2}-1}
l_{i} \! + \! 1,\sum_{i=1}^{r_{2}-1}l_{i} \! + \! 2,\dotsc,\sum_{i=1}^{r_{2}-1}
l_{i} \! + \! l_{r_{2}} \! - \! \delta_{r_{2} \mathfrak{s}}$, $m(r_{2}) \! = \! 
1,2,\dotsc,l_{r_{2}} \! - \! \delta_{r_{2} \mathfrak{s}}$, one proceeds, 
via the latter determinantal expression, with the analysis of 
$\overset{\scriptscriptstyle f}{\Xi}^{\raise-1.0ex\hbox{$\scriptstyle n$}}_{k}(z)$:
\begin{align*}
\overset{\scriptscriptstyle f}{\Xi}^{\raise-1.0ex\hbox{$\scriptstyle n$}}_{k}
(z) =& \, \underbrace{\int_{\mathbb{R}} \int_{\mathbb{R}} \dotsb \int_{
\mathbb{R}}}_{(n-1)K+k} \md \widetilde{\mu}(\xi_{0}) \, \md \widetilde{\mu}
(\xi_{1}) \, \dotsb \, \md \widetilde{\mu}(\xi_{l_{1}}) \, \dotsb \, \md \widetilde{\mu}
(\xi_{n_{1}}) \, \dotsb \, \md \widetilde{\mu}(\xi_{n_{2}}) \, \md \widetilde{\mu}
(\xi_{m_{2}-\varkappa_{nk}+1}) \, \dotsb \, \md \widetilde{\mu}(\xi_{m_{1}}) \\
\times& \, \dfrac{1}{(\xi_{0} \! - \! \alpha_{k})^{0}(\xi_{1} \! - \! 
\alpha_{p_{1}})^{1} \dotsb (\xi_{l_{1}} \! - \! \alpha_{p_{1}})^{l_{1}} 
\dotsb \dfrac{1}{(\xi_{n_{1}})^{1}} \dotsb \dfrac{1}{(\xi_{n_{2}})^{
l_{\mathfrak{s}-1}}}(\xi_{m_{2}-\varkappa_{nk}+1} \! - \! \alpha_{k})^{1} 
\dotsb (\xi_{m_{1}} \! - \! \alpha_{k})^{\varkappa_{nk}-1}} \\
\times& \, \dfrac{(\phi_{0}(z))^{-1}}{\phi_{0}(\xi_{0}) \phi_{0}(\xi_{1}) 
\dotsb \phi_{0}(\xi_{l_{1}}) \dotsb \phi_{0}(\xi_{n_{1}}) \dotsb \phi_{0}
(\xi_{n_{2}}) \phi_{0}(\xi_{m_{2}-\varkappa_{nk}+1}) \dotsb \phi_{0}(\xi_{
m_{1}})} \\
\times& \, 
\underbrace{\left\vert

\right\vert}_{=: \, \mathbb{G}^{\spadesuit}(\xi_{0},\xi_{1},\dotsc,\xi_{l_{1}},
\dotsc,\xi_{n_{1}},\dotsc,\xi_{n_{2}},\xi_{m_{2}-\varkappa_{nk}+1},\dotsc,
\xi_{(n-1)K+k-1};z)} \\
=& \, \dfrac{(\phi_{0}(z))^{-1}}{m_{2}!} \sum_{\pmb{\sigma} \in \mathfrak{S}_{
m_{2}}} \underbrace{\int_{\mathbb{R}} \int_{\mathbb{R}} \dotsb \int_{\mathbb{
R}}}_{(n-1)K+k} \md \widetilde{\mu}(\xi_{\sigma (0)}) \, \md \widetilde{\mu}
(\xi_{\sigma (1)}) \, \dotsb \, \md \widetilde{\mu}(\xi_{\sigma (l_{1})}) \, \dotsb 
\, \md \widetilde{\mu}(\xi_{\sigma (n_{1})}) \, \dotsb \, \md \widetilde{\mu}
(\xi_{\sigma (n_{2})}) \\
\times& \, \tfrac{\dotsb \, \md \widetilde{\mu}(\xi_{\sigma (m_{2}-\varkappa_{nk}+1)}) 
\, \dotsb \, \md \widetilde{\mu}(\xi_{\sigma (m_{1})})}{(\xi_{\sigma (0)}-\alpha_{k})^{0}
(\xi_{\sigma (1)}-\alpha_{p_{1}})^{1} \dotsb (\xi_{\sigma (l_{1})}-\alpha_{p_{1}})^{l_{1}} 
\dotsb \frac{1}{(\xi_{\sigma (n_{1})})^{1}} \dotsb \frac{1}{(\xi_{\sigma (n_{2})})^{
l_{\mathfrak{s}-1}}}(\xi_{\sigma (m_{2}-\varkappa_{nk}+1)}-\alpha_{k})^{1} \dotsb 
(\xi_{\sigma (m_{1})}-\alpha_{k})^{\varkappa_{nk}-1}} \\
\times& \, \dfrac{1}{\phi_{0}(\xi_{\sigma (0)}) \phi_{0}(\xi_{\sigma (1)}) 
\dotsb \phi_{0}(\xi_{\sigma (l_{1})}) \dotsb \phi_{0}(\xi_{\sigma (n_{1})}) 
\dotsb \phi_{0}(\xi_{\sigma (n_{2})}) \phi_{0}(\xi_{\sigma (m_{2}-\varkappa_{
nk}+1)}) \dotsb \phi_{0}(\xi_{\sigma (m_{1})})} \\
\times& \, \mathbb{G}^{\spadesuit}(\xi_{\sigma (0)},\xi_{\sigma (1)},\dotsc,
\xi_{\sigma (l_{1})},\dotsc,\xi_{\sigma (n_{1})},\dotsc,\xi_{\sigma (n_{2})},
\xi_{\sigma (m_{2}-\varkappa_{nk}+1)},\dotsc,\xi_{\sigma ((n-1)K+k-1)};z) \\
=& \, \dfrac{(\phi_{0}(z))^{-1}}{m_{2}!} \underbrace{\int_{\mathbb{R}} 
\int_{\mathbb{R}} \dotsb \int_{\mathbb{R}}}_{(n-1)K+k} \md \widetilde{\mu}
(\xi_{0}) \, \md \widetilde{\mu}(\xi_{1}) \, \dotsb \, \md \widetilde{\mu}
(\xi_{l_{1}}) \, \dotsb \, \md \widetilde{\mu}(\xi_{n_{1}}) \, \dotsb \, \md 
\widetilde{\mu}(\xi_{n_{2}}) \, \md \widetilde{\mu}(\xi_{m_{2}-\varkappa_{nk}+1}) 
\, \dotsb \\
\times& \, \dotsb \, \md \widetilde{\mu}(\xi_{m_{1}}) \, \mathbb{G}^{\spadesuit}
(\xi_{0},\xi_{1},\dotsc,\xi_{l_{1}},\dotsc,\xi_{n_{1}},\dotsc,\xi_{n_{2}},
\xi_{m_{2}-\varkappa_{nk}+1},\dotsc,\xi_{(n-1)K+k-1};z) \\
\times& \, \sum_{\pmb{\sigma} \in \mathfrak{S}_{m_{2}}} \operatorname{sgn}
(\pmb{\pmb{\sigma}}) \tfrac{1}{(\xi_{\sigma (0)}-\alpha_{k})^{0}(\xi_{\sigma 
(1)}-\alpha_{p_{1}})^{1} \dotsb (\xi_{\sigma (l_{1})}-\alpha_{p_{1}})^{l_{1}} 
\dotsb \frac{1}{(\xi_{\sigma (n_{1})})^{1}} \dotsb \frac{1}{(\xi_{\sigma 
(n_{2})})^{l_{\mathfrak{s}-1}}}(\xi_{\sigma (m_{2}-\varkappa_{nk}+1)}-
\alpha_{k})^{1} \dotsb (\xi_{\sigma (m_{1})}-\alpha_{k})^{\varkappa_{nk}-1}} \\
\times& \, \dfrac{1}{\phi_{0}(\xi_{\sigma (0)}) \phi_{0}(\xi_{\sigma (1)}) 
\dotsb \phi_{0}(\xi_{\sigma (l_{1})}) \dotsb \phi_{0}(\xi_{\sigma (n_{1})}) 
\dotsb \phi_{0}(\xi_{\sigma (n_{2})}) \phi_{0}(\xi_{\sigma (m_{2}-\varkappa_{
nk}+1)}) \dotsb \phi_{0}(\xi_{\sigma (m_{1})})} \\
=& \, \dfrac{(\phi_{0}(z))^{-1}}{m_{2}!} \underbrace{\int_{\mathbb{R}} \int_{
\mathbb{R}} \dotsb \int_{\mathbb{R}}}_{(n-1)K+k} \md \widetilde{\mu}(\xi_{0}) 
\, \md \widetilde{\mu}(\xi_{1}) \, \dotsb \, \md \widetilde{\mu}(\xi_{l_{1}}) \, 
\dotsb \, \md \widetilde{\mu}(\xi_{n_{1}}) \, \dotsb \, \md \widetilde{\mu}
(\xi_{n_{2}}) \, \md \widetilde{\mu}(\xi_{m_{2}-\varkappa_{nk}+1}) \, \dotsb \\
\times& \, \dotsb \, \md \widetilde{\mu}(\xi_{m_{1}}) \, \dfrac{\mathbb{G}^{
\spadesuit}(\xi_{0},\xi_{1},\dotsc,\xi_{l_{1}},\dotsc,\xi_{n_{1}},\dotsc,\xi_{n_{2}},
\xi_{m_{2}-\varkappa_{nk}+1},\dotsc,\xi_{(n-1)K+k-1};z)}{\phi_{0}(\xi_{0}) 
\phi_{0}(\xi_{1}) \dotsb \phi_{0}(\xi_{l_{1}}) \dotsb \phi_{0}(\xi_{n_{1}}) 
\dotsb \phi_{0}(\xi_{n_{2}}) \phi_{0}(\xi_{m_{2}-\varkappa_{nk}+1}) \dotsb 
\phi_{0}(\xi_{m_{1}})} \\
\times& \, 
\left\vert

\right\vert \\
=& \, \dfrac{(\phi_{0}(z))^{-1}}{m_{2}!} \underbrace{\int_{\mathbb{R}} \int_{
\mathbb{R}} \dotsb \int_{\mathbb{R}}}_{(n-1)K+k} \md \widetilde{\mu}(\xi_{0}) 
\, \md \widetilde{\mu}(\xi_{1}) \, \dotsb \, \md \widetilde{\mu}(\xi_{l_{1}}) \, 
\dotsb \, \md \widetilde{\mu}(\xi_{n_{1}}) \, \dotsb \, \md \widetilde{\mu}
(\xi_{n_{2}}) \, \md \widetilde{\mu}(\xi_{m_{2}-\varkappa_{nk}+1}) \, \dotsb \\
\times& \, \dotsb \, \md \widetilde{\mu}(\xi_{m_{1}}) \, \dfrac{\mathbb{G}^{
\spadesuit}(\xi_{0},\xi_{1},\dotsc,\xi_{l_{1}},\dotsc,\xi_{n_{1}},\dotsc,\xi_{n_{2}},
\xi_{m_{2}-\varkappa_{nk}+1},\dotsc,\xi_{(n-1)K+k-1};z)}{\phi_{0}(\xi_{0}) 
\phi_{0}(\xi_{1}) \dotsb \phi_{0}(\xi_{l_{1}}) \dotsb \phi_{0}(\xi_{n_{1}}) 
\dotsb \phi_{0}(\xi_{n_{2}}) \phi_{0}(\xi_{m_{2}-\varkappa_{nk}+1}) \dotsb 
\phi_{0}(\xi_{m_{1}})} \\
\times& \, \dfrac{1}{\tilde{\phi}_{0}(\xi_{0}) \tilde{\phi}_{0}(\xi_{1}) 
\dotsb \tilde{\phi}_{0}(\xi_{l_{1}}) \dotsb \tilde{\phi}_{0}(\xi_{n_{1}}) 
\dotsb \tilde{\phi}_{0}(\xi_{n_{2}}) \tilde{\phi}_{0}(\xi_{m_{2}-\varkappa_{n
k}+1}) \dotsb \tilde{\phi}_{0}(\xi_{m_{1}})} \\
\times& \, 
\underbrace{\left\vert

\right\vert}_{= \, \mathbb{D}^{\blacklozenge} \prod_{\underset{j<i}{i,j=0}}^{
(n-1)K+k-1}(\xi_{i}-\xi_{j}) \quad \text{(cf. proof of Lemma~\ref{lem2.1}, 
case~\textbf{(2)})}} \\
=& \, \dfrac{\mathbb{D}^{\blacklozenge}(\phi_{0}(z))^{-1}}{((n \! - \! 1)K \! 
+ \! k)!} \underbrace{\int_{\mathbb{R}} \int_{\mathbb{R}} \dotsb \int_{
\mathbb{R}}}_{(n-1)K+k} \md \widetilde{\mu}(\xi_{0}) \, \md \widetilde{\mu}
(\xi_{1}) \, \dotsb \, \md \widetilde{\mu}(\xi_{l_{1}}) \, \dotsb \, \md \widetilde{\mu}
(\xi_{n_{1}}) \, \dotsb \, \md \widetilde{\mu}(\xi_{n_{2}}) \, \md \widetilde{\mu}
(\xi_{m_{2}-\varkappa_{nk}+1}) \, \dotsb \\
\times& \, \dotsb \, \md \widetilde{\mu}(\xi_{m_{1}}) \, 
\dfrac{\prod_{\underset{j<i}{i,j=0}}^{(n-1)K+k-1}(\xi_{i} \! - \! \xi_{j})}{
\phi_{0}(\xi_{0}) \phi_{0}(\xi_{1}) \dotsb \phi_{0}(\xi_{l_{1}}) \dotsb \phi_{0}
(\xi_{n_{1}}) \dotsb \phi_{0}(\xi_{n_{2}}) \phi_{0}(\xi_{m_{2}-\varkappa_{nk}+1}) 
\dotsb \phi_{0}(\xi_{m_{1}})} \\
\times& \, \dfrac{1}{\tilde{\phi}_{0}(\xi_{0}) \tilde{\phi}_{0}(\xi_{1}) 
\dotsb \tilde{\phi}_{0}(\xi_{l_{1}}) \dotsb \tilde{\phi}_{0}(\xi_{n_{1}}) 
\dotsb \tilde{\phi}_{0}(\xi_{n_{2}}) \tilde{\phi}_{0}(\xi_{m_{2}-\varkappa_{n
k}+1}) \dotsb \tilde{\phi}_{0}(\xi_{m_{1}})} \\
\times& \, 
\underbrace{\left\vert

\right\vert}_{= \, \mathbb{D}^{\lozenge} \prod_{\underset{j<i}{i,j=0}}^{(n-1)
K+k-1}(\xi_{i}-\xi_{j}) \prod_{m=0}^{(n-1)K+k-1}(z-\xi_{m}) \quad \text{(cf. 
proof of Lemma~\ref{lem2.1}, case~\textbf{(2)})}} \\
=& \, \dfrac{\mathbb{D}^{\blacklozenge} \mathbb{D}^{\lozenge}(\phi_{0}(z))^{-
1}}{((n \! - \! 1)K \! + \! k)!} \underbrace{\int_{\mathbb{R}} \int_{\mathbb{R}} 
\dotsb \int_{\mathbb{R}}}_{(n-1)K+k} \md \widetilde{\mu}(\xi_{0}) \, \md 
\widetilde{\mu}(\xi_{1}) \, \dotsb \, \md \widetilde{\mu}(\xi_{l_{1}}) \, \dotsb 
\, \md \widetilde{\mu}(\xi_{n_{1}}) \, \dotsb \, \md \widetilde{\mu}(\xi_{n_{2}}) 
\, \md \widetilde{\mu}(\xi_{m_{2}-\varkappa_{nk}+1}) \, \dotsb \\
\times& \, \dotsb \, \md \widetilde{\mu}(\xi_{m_{1}}) \, \dfrac{(\tilde{\phi}_{0}
(\xi_{0}) \tilde{\phi}_{0}(\xi_{1}) \dotsb \tilde{\phi}_{0}(\xi_{l_{1}}) \dotsb 
\tilde{\phi}_{0}(\xi_{n_{1}}) \dotsb \tilde{\phi}_{0}(\xi_{n_{2}}) \tilde{
\phi}_{0}(\xi_{m_{2}-\varkappa_{nk}+1}) \dotsb \tilde{\phi}_{0}(\xi_{m_{1}})
)^{-1}}{\phi_{0}(\xi_{0}) \phi_{0}(\xi_{1}) \dotsb \phi_{0}(\xi_{l_{1}}) 
\dotsb \phi_{0}(\xi_{n_{1}}) \dotsb \phi_{0}(\xi_{n_{2}}) \phi_{0}
(\xi_{m_{2}-\varkappa_{nk}+1}) \dotsb \phi_{0}(\xi_{m_{1}})} \\
\times& \, \prod_{\substack{i,j=0\\j<i}}^{(n-1)K+k-1}(\xi_{i} \! 
- \! \xi_{j})^{2} \, \prod_{m=0}^{(n-1)K+k-1}(z \! - \! \xi_{m});
\end{align*}
but, recalling that, for $n \! \in \! \mathbb{N}$ and $k \! \in \! \lbrace 
1,2,\dotsc,K \rbrace$ such that $\alpha_{p_{\mathfrak{s}}} \! := \! 
\alpha_{k} \! \neq \! \infty$, $\tilde{\phi}_{0}(z) \! = \! \phi_{0}
(z)/(z \! - \! \alpha_{k})$, one 
arrives at
\begin{align*}
\overset{\scriptscriptstyle f}{\Xi}_{k}^{\raise-1.0ex\hbox{$\scriptstyle n$}}
(z) =& \, \dfrac{\mathbb{D}^{\blacklozenge} \mathbb{D}^{\lozenge}(\phi_{0}
(z))^{-1}}{((n \! - \! 1)K \! + \! k)!} \underbrace{\int_{\mathbb{R}} \int_{
\mathbb{R}} \dotsb \int_{\mathbb{R}}}_{(n-1)K+k} \md \widetilde{\mu}(\xi_{0}) 
\, \md \widetilde{\mu}(\xi_{1}) \, \dotsb \, \md \widetilde{\mu}(\xi_{l_{1}}) \, 
\dotsb \, \md \widetilde{\mu}(\xi_{n_{1}}) \, \dotsb \, \md \widetilde{\mu}
(\xi_{n_{2}}) \, \dotsb \, \md \widetilde{\mu}(\xi_{(n-1)K+k-1}) \\
\times& \, \left(\prod_{m=0}^{(n-1)K+k-1} \phi_{0}(\xi_{m}) \right)^{-2} 
\prod_{\substack{i,j=0\\j<i}}^{(n-1)K+k-1}(\xi_{i} \! - \! \xi_{j})^{2} \, 
\prod_{l=0}^{(n-1)K+k-1}(\xi_{l} \! - \! \alpha_{k})(z \! - \! \xi_{l}) \, \, 
\Rightarrow \\
\overset{\scriptscriptstyle f}{\Xi}_{k}^{\raise-1.0ex\hbox{$\scriptstyle n$}}
(z) =& \, \dfrac{\mathbb{D}^{\blacklozenge} \mathbb{D}^{\lozenge}}{
((n \! - \! 1)K \! + \! k)!} \left(\prod_{q=1}^{\mathfrak{s}-2}
(z \! - \! \alpha_{p_{q}})^{\varkappa_{nk \tilde{k}_{q}}}(z \! - \! 
\alpha_{k})^{\varkappa_{nk}} \right)^{-1} \underbrace{\int_{\mathbb{R}} 
\int_{\mathbb{R}} \dotsb \int_{\mathbb{R}}}_{(n-1)K+k} \md \widetilde{\mu}
(\xi_{0}) \, \md \widetilde{\mu}(\xi_{1}) \, \dotsb \, \md \widetilde{\mu}
(\xi_{(n-1)K+k-1}) \\
\times& \, \left(\prod_{m=0}^{(n-1)K+k-1} \prod_{q=1}^{
\mathfrak{s}-2}(\xi_{m} \! - \! \alpha_{p_{q}})^{\varkappa_{nk \tilde{k}_{q}}} 
(\xi_{m} \! - \! \alpha_{k})^{\varkappa_{nk}} \right)^{-2} 
\prod_{\substack{i,j=0\\j<i}}^{(n-1)K+k-1}(\xi_{i} \! - \! \xi_{j})^{2} \, 
\prod_{l=0}^{(n-1)K+k-1}(\xi_{l} \! - \! \alpha_{k})(z \! - \! \xi_{l}).
\end{align*}
Recalling that, for $n \! \in \! \mathbb{N}$ and $k \! \in \! \lbrace 1,2,
\dotsc,K \rbrace$ such that $\alpha_{p_{\mathfrak{s}}} \! := \! \alpha_{k} 
\! \neq \! \infty$, $\pmb{\pi}^{n}_{k}(z) \! = \! \mathcal{X}_{11}
(z)/(z \! - \! \alpha_{k}) \! = \! 
\overset{\scriptscriptstyle f}{\Xi}_{k}^{\raise-1.0ex\hbox{$\scriptstyle n$}}
(z)/\hat{\mathfrak{c}}_{D_{f}}$, where 
$\overset{\scriptscriptstyle f}{\Xi}_{k}^{\raise-1.0ex\hbox{$\scriptstyle n$}}
(z)$ is given directly above, and $\hat{\mathfrak{c}}_{D_{f}}$ is given by 
Equation~\eqref{eq65}, one arrives at the integral representation for 
$\pmb{\pi}^{n}_{k}(z)$ given in Equation~\eqref{eintreppfin1}.

The determinantal representation for $\mathcal{X}_{21}(z)$ is now considered. 
For $n \! \in \! \mathbb{N}$ and $k \! \in \! \lbrace 1,2,\dotsc,K \rbrace$ 
such that $\alpha_{p_{\mathfrak{s}}} \! := \! \alpha_{k} \! \neq \! \infty$, 
recall the ordered disjoint partition for $\lbrace \alpha_{1},\alpha_{2},\dotsc,
\alpha_{K} \rbrace \cup \lbrace \alpha_{1},\alpha_{2},\dotsc,\alpha_{K} 
\rbrace \cup \dotsb \cup \lbrace \alpha_{1},\alpha_{2},\dotsc,\alpha_{k} 
\rbrace$ introduced in the proof of Lemma~\ref{lem2.1} (cf. case~\textbf{(2)}). 
For this ordered disjoint partition, introduce, for $n \! \in \! \mathbb{N}$ and 
$k \! \in \! \lbrace 1,2,\dotsc,K \rbrace$ such that $\alpha_{p_{\mathfrak{s}}} 
\! := \! \alpha_{k} \! \neq \! \infty$, the following notation (recall that $l_{q} 
\! = \! \varkappa_{nk \tilde{k}_{q}}$, $q \! = \! 1,2,\dotsc,\mathfrak{s} \! - \! 
2$, $l_{\mathfrak{s}-1} \! = \! \varkappa^{\infty}_{nk \tilde{k}_{\mathfrak{s}-1}}$, 
and $l_{\mathfrak{s}} \! = \! \varkappa_{nk})$:
\begin{equation*}
\phi_{0}^{\spcheck}(z) \! := \! \prod_{m=1}^{\mathfrak{s}-2}(z \! - \! 
\alpha_{p_{m}})^{l_{m}}(z \! - \! \alpha_{k})^{\varkappa_{nk}-2} \! =: 
\sum_{j=0}^{(n-1)K+k-2} \mathfrak{a}^{\clubsuit}_{j,0}z^{j},
\end{equation*}
and, for $r \! = \! 1,\dotsc,\mathfrak{s} \! - \! 1,\mathfrak{s}$, $q(r) \! = 
\! \sum_{i=1}^{r-1}l_{i} \! + \! 1,\sum_{i=1}^{r-1}l_{i} \! + \! 2,\dotsc,
\sum_{i=1}^{r-1}l_{i} \! + \! l_{r} \! - \! 2 \delta_{r \mathfrak{s}}$, and 
$m(r) \! = \! 1,2,\dotsc,l_{r} \! - \! 2 \delta_{r \mathfrak{s}}$,
\begin{equation*}
\phi_{q(r)}^{\spcheck}(z) \! := \! \left\lbrace \phi_{0}^{\spcheck}(z)
(z \! - \! \alpha_{p_{r}})^{-m(r)(1-\delta_{r \mathfrak{s}-1})}z^{m(r) 
\delta_{r \mathfrak{s}-1}} \! =: \sum_{j=0}^{(n-1)K+k-2} 
\mathfrak{a}^{\clubsuit}_{j,q(r)}z^{j} \right\rbrace;
\end{equation*}
e.g., for $r \! = \! 1$, the notation $\phi_{q(1)}^{\spcheck}(z) \! := \! 
\lbrace \phi_{0}^{\spcheck}(z)(z \! - \! \alpha_{p_{1}})^{-m(1)} \! =: \! 
\sum_{j=0}^{(n-1)K+k-2} \mathfrak{a}^{\clubsuit}_{j,q(1)}z^{j} \rbrace$, 
$q(1) \! = \! 1,2,\dotsc,l_{1}$, $m(1) \! = \! 1,2,\dotsc,l_{1}$, denotes 
the (set of) $l_{1} \! = \! \varkappa_{nk \tilde{k}_{1}}$ functions
\begin{gather*}
\phi_{1}^{\spcheck}(z) \! = \! \dfrac{\phi_{0}^{\spcheck}(z)}{z \! - \! 
\alpha_{p_{1}}} \! =: \sum_{j=0}^{(n-1)K+k-2} \mathfrak{a}^{\clubsuit}_{j,1}
z^{j}, \, \phi_{2}^{\spcheck}(z) \! = \! \dfrac{\phi_{0}^{\spcheck}(z)}{(z 
\! - \! \alpha_{p_{1}})^{2}} \! =: \sum_{j=0}^{(n-1)K+k-2} \mathfrak{a}^{
\clubsuit}_{j,2}z^{j}, \, \dotsc, \, \phi_{l_{1}}^{\spcheck}(z) \! = \! 
\dfrac{\phi_{0}^{\spcheck}(z)}{(z \! - \! \alpha_{p_{1}})^{l_{1}}} \! =: 
\sum_{j=0}^{(n-1)K+k-2} \mathfrak{a}^{\clubsuit}_{j,l_{1}}z^{j},
\end{gather*}
for $r \! = \! \mathfrak{s} \! - \! 1$ (recall that $\alpha_{p_{\mathfrak{s}
-1}} \! = \! \infty)$, the notation $\phi_{q(\mathfrak{s}-1)}^{\spcheck}(z) 
\! := \! \lbrace \phi_{0}^{\spcheck}(z)z^{m(\mathfrak{s}-1)} \! =: \! 
\sum_{j=0}^{(n-1)K+k-2} \mathfrak{a}^{\clubsuit}_{j,q(\mathfrak{s}-1)}z^{j} 
\rbrace$, $q(\mathfrak{s} \! - \! 1) \! = \! l_{1} \! + \! \dotsb \! + \! 
l_{\mathfrak{s}-2} \! + \! 1,l_{1} \! + \! \dotsb \! + \! l_{\mathfrak{s}-2} 
\! + \! 2,\dotsc,l_{1} \! + \! \dotsb \! + \! l_{\mathfrak{s}-2} \! + \! 
l_{\mathfrak{s}-1} \! = \! (n \! - \! 1)K \! + \! k \! - \! \varkappa_{nk}$, 
$m(\mathfrak{s} \! - \! 1) \! = \! 1,2,\dotsc,l_{\mathfrak{s}-1}$, denotes 
the (set of) $l_{\mathfrak{s}-1} \! = \! \varkappa^{\infty}_{nk \tilde{k}_{
\mathfrak{s}-1}}$ functions
\begin{gather*}
\phi_{l_{1}+\dotsb +l_{\mathfrak{s}-2}+1}^{\spcheck}(z) \! = \! \phi_{0}^{
\spcheck}(z)z \! =: \sum_{j=0}^{(n-1)K+k-2} \mathfrak{a}^{\clubsuit}_{j,l_{1}+
\dotsb +l_{\mathfrak{s}-2}+1}z^{j}, \, \phi_{l_{1}+\dotsb +l_{\mathfrak{s}-2}+
2}^{\spcheck}(z) \! = \! \phi_{0}^{\spcheck}(z)z^{2} \! =: \sum_{j=0}^{(n-1)K+
k-2} \mathfrak{a}^{\clubsuit}_{j,l_{1}+ \dotsb +l_{\mathfrak{s}-2}+2}z^{j}, \, 
\dotsc \\
\dotsc, \, \phi_{(n-1)K+k-\varkappa_{nk}}^{\spcheck}(z) \! = \! \phi_{0}^{
\spcheck}(z)z^{l_{\mathfrak{s}-1}} \! =: \sum_{j=0}^{(n-1)K+k-2} 
\mathfrak{a}^{\clubsuit}_{j,(n-1)K+k-\varkappa_{nk}}z^{j},
\end{gather*}
etc., and, for $r \! = \! \mathfrak{s}$, the notation $\phi_{q(\mathfrak{s})}^{
\spcheck}(z) \! := \! \lbrace \phi_{0}^{\spcheck}(z)(z \! - \! \alpha_{k})^{-
m(\mathfrak{s})} \! =: \! \sum_{j=0}^{(n-1)K+k-2} \mathfrak{a}^{\clubsuit}_{j,
q(\mathfrak{s})}z^{j} \rbrace$, $q(\mathfrak{s}) \! = \! (n \! - \! 1)K \! + \! 
k \! - \! \varkappa_{nk} \! + \! 1,(n \! - \! 1)K \! + \! k \! - \! \varkappa_{nk} 
\! + \! 2,\dotsc,(n \! - \! 1)K \! + \! k \! - \! 2$, $m(\mathfrak{s}) \! = \! 
1,2,\dotsc,l_{\mathfrak{s}} \! - \! 2$, denotes the (set of) $l_{\mathfrak{s}} 
\! - \! 2 \! = \! \varkappa_{nk} \! - \! 2$ functions
\begin{gather*}
\phi_{(n-1)K+k-\varkappa_{nk}+1}^{\spcheck}(z) \! = \! \dfrac{\phi_{0}^{
\spcheck}(z)}{z \! - \! \alpha_{k}} \! =: \sum_{j=0}^{(n-1)K+k-2} \mathfrak{
a}^{\clubsuit}_{j,(n-1)K+k-\varkappa_{nk}+1}z^{j},\, \phi_{(n-1)K+k-
\varkappa_{nk}+2}^{\spcheck}(z) \! = \! \dfrac{\phi_{0}^{\spcheck}(z)}{(z \! - 
\! \alpha_{k})^{2}} \! =: \sum_{j=0}^{(n-1)K+k-2} \mathfrak{a}^{\clubsuit}_{j,
(n-1)K+k-\varkappa_{nk}+2}z^{j}, \, \dotsc \\
\dotsc, \, \phi_{(n-1)K+k-2}^{\spcheck}(z) \! = \! \dfrac{\phi_{0}^{\spcheck}
(z)}{(z \! - \! \alpha_{k})^{\varkappa_{nk}-2}} \! =: \sum_{j=0}^{(n-1)K+k-2} 
\mathfrak{a}^{\clubsuit}_{j,(n-1)K+k-2}z^{j}.
\end{gather*}
(Note: $\# \lbrace \phi_{0}^{\spcheck}(z)(z \! - \! \alpha_{p_{r}})^{-m(r)
(1-\delta_{r \mathfrak{s}-1})}z^{m(r) \delta_{r \mathfrak{s}-1}} \rbrace 
\! = \! l_{r} \! - \! 2 \delta_{r \mathfrak{s}}$, $r \! = \! 1,2,\dotsc,
\mathfrak{s}$, and $\# \cup_{r=1}^{\mathfrak{s}} \lbrace \phi_{0}^{\spcheck}
(z)(z \! - \! \alpha_{p_{r}})^{-m(r)(1-\delta_{r \mathfrak{s}-1})}
z^{m(r) \delta_{r \mathfrak{s}-1}} \rbrace \! = \! \sum_{r=1}^{\mathfrak{s}}
l_{r} \! - \! 2 \! = \! (n \! - \! 1)K \! + \! k \! - \! 2$.) One notes that, 
for $n \! \in \! \mathbb{N}$ and $k \! \in \! \lbrace 1,2,\dotsc,K \rbrace$ 
such that $\alpha_{p_{\mathfrak{s}}} \! := \! \alpha_{k} \! \neq \! \infty$, 
the $l_{1} \! + \! \dotsb \! + \! l_{\mathfrak{s}-1} \! + \! l_{\mathfrak{s}} \! - 
\! 1 \! = \! (n \! - \! 1)K \! + \! k \! - \! 1$ functions $\phi_{0}^{\spcheck}(z),
\phi_{1}^{\spcheck}(z),\dotsc,\phi_{l_{1}}^{\spcheck}(z),\dotsc,\phi_{l_{1}+
\dotsb +l_{\mathfrak{s}-2}+1}^{\spcheck}(z),\dotsc,\phi_{(n-1)K+k-
\varkappa_{nk}}^{\spcheck}(z),\phi_{(n-1)K+k-\varkappa_{nk}+1}^{\spcheck}
(z),\dotsc,\linebreak[4] 
\phi_{(n-1)K+k-2}^{\spcheck}(z)$ are linearly independent on $\mathbb{R}$, 
that is, for $z \! \in \! \mathbb{R}$, $\sum_{j=0}^{(n-1)K+k-2} 
\mathfrak{c}^{\clubsuit}_{j} \phi_{j}^{\spcheck}(z) \! = \! 0$ $\Rightarrow$ 
(via a Vandermonde-type argument; see the $((n \! - \! 1)K \! + \! k \! - 
\! 1) \times ((n \! - \! 1)K \! + \! k \! - \! 1)$ non-zero determinant 
$\mathbb{D}^{\clubsuit}$ in Equation~\eqref{eq77} below) $\mathfrak{c}^{
\clubsuit}_{j} \! = \! 0$, $j \! = \! 0,1,\dotsc,(n \! - \! 1)K \! + \! k \! - \! 2$. 
For $n \! \in \! \mathbb{N}$ and $k \! \in \! \lbrace 1,2,\dotsc,K \rbrace$ 
such that $\alpha_{p_{\mathfrak{s}}} \! := \! \alpha_{k} \! \neq \! \infty$, let 
$\mathfrak{S}_{(n-1)K+k-1}$ denote the $((n \! - \! 1)K \! + \! k \! - \! 1)!$ 
permutations of $\lbrace 0,1,\dotsc,(n \! - \! 1)K \! + \! k \! - \! 2 \rbrace$. 
For $n \! \in \! \mathbb{N}$ and $k \! \in \! \lbrace 1,2,\dotsc,K \rbrace$ 
such that $\alpha_{p_{\mathfrak{s}}} \! := \! \alpha_{k} \! \neq \! \infty$, 
corresponding to the ordered disjoint partition above, recall the formula for 
$\mathcal{X}_{21}(z)/(z \! - \! \alpha_{k})$ given in the proof of 
Lemma~\ref{lem2.1} (cf. case~\textbf{(2)}):
\begin{equation*}
\dfrac{\mathcal{X}_{21}(z)}{z \! - \! \alpha_{k}} \! = \! \sum_{m=1}^{
\mathfrak{s}-2} \sum_{q=1}^{l_{m}=\varkappa_{nk \tilde{k}_{m}}} \dfrac{
\widehat{\nu}^{\raise-1.0ex\hbox{$\scriptstyle f$}}_{m,q}(n,k)}{(z \! - \! 
\alpha_{p_{m}})^{q}} \! + \! \sum_{q=1}^{l_{\mathfrak{s}-1}=\varkappa^{
\infty}_{nk \tilde{k}_{\mathfrak{s}-1}}} \widehat{\nu}^{
\raise-1.0ex\hbox{$\scriptstyle f$}}_{\mathfrak{s}-1,q}(n,k)z^{q} \! + \! 
\sum_{q=1}^{l_{\mathfrak{s}}=\varkappa_{nk}} \dfrac{\widehat{\nu}^{
\raise-1.0ex\hbox{$\scriptstyle f$}}_{\mathfrak{s},q}(n,k)}{(z \! - \! 
\alpha_{k})^{q-1}}, \quad 
\widehat{\nu}^{\raise-1.0ex\hbox{$\scriptstyle f$}}_{\mathfrak{s},
\varkappa_{nk}}(n,k) \! \neq \! 0,
\end{equation*}
where the $(n \! - \! 1)K \! + \! k$ coefficients $\widehat{\nu}^{
\raise-1.0ex\hbox{$\scriptstyle f$}}_{1,1}(n,k),\dotsc,\widehat{\nu}^{
\raise-1.0ex\hbox{$\scriptstyle f$}}_{1,l_{1}}(n,k),\dotsc,\widehat{\nu}^{
\raise-1.0ex\hbox{$\scriptstyle f$}}_{\mathfrak{s}-1,1}(n,k),\dotsc,\widehat{
\nu}^{\raise-1.0ex\hbox{$\scriptstyle f$}}_{\mathfrak{s}-1,l_{\mathfrak{s}-1}}
(n,k),\widehat{\nu}^{\raise-1.0ex\hbox{$\scriptstyle f$}}_{\mathfrak{s},1}(n,
k),\dotsc,\widehat{\nu}^{\raise-1.0ex\hbox{$\scriptstyle f$}}_{\mathfrak{s},
\varkappa_{nk}}(n,k)$, with $\widehat{\nu}^{\raise-1.0ex\hbox{$\scriptstyle 
f$}}_{\mathfrak{s},\varkappa_{nk}}(n,k) \! \neq \! 0$, satisfy the linear 
inhomogeneous algebraic system of equations~\eqref{eq74}. Via the 
multi-linearity property of the determinant and an application of Cramer's 
Rule to system~\eqref{eq74}, one arrives at, for $n \! \in \! \mathbb{N}$ 
and $k \! \in \! \lbrace 1,2,\dotsc,K \rbrace$ such that 
$\alpha_{p_{\mathfrak{s}}} \! := \! \alpha_{k} \! \neq \! \infty$, the 
following---ordered---determinantal representation for 
$\mathcal{X}_{21}(z)/(z \! - \! \alpha_{k})$:
\begin{equation*}
\dfrac{\mathcal{X}_{21}(z)}{z \! - \! \alpha_{k}} \! = \! \dfrac{2 \pi \mi}{
\mathcal{N}_{f}^{\sharp}(n,k)} 
\overset{\triangledown}{\Xi}^{\raise-1.0ex\hbox{$\scriptstyle n$}}_{k}(z),
\end{equation*}
where $\mathcal{N}_{f}^{\sharp}(n,k)$ is given by Equation~\eqref{eq75}, 
and, with (abuse of notation)
\begin{equation*}
n_{1} \! = \! l_{1} \! + \! \dotsb \! + \! l_{\mathfrak{s}-2} \! + \! 1, 
\qquad n_{2} \! = \! (n \! - \! 1)K \! + \! k \! - \! \varkappa_{nk}, \qquad 
m_{1} \! = \! (n \! - \! 1)K \! + \! k \! - \! \varkappa_{nk} \! + \! 1, 
\qquad \text{and} \qquad m_{2} \! = \! (n \! - \! 1)K \! + \! k \! - \! 2,
\end{equation*}
\begin{align*}
\overset{\triangledown}{\Xi}^{\raise-1.0ex\hbox{$\scriptstyle n$}}_{k}(z):=
\setcounter{MaxMatrixCols}{12}
&\left\vert

\right\vert \\
=& \, \underbrace{\int_{\mathbb{R}} \int_{\mathbb{R}} \dotsb \int_{\mathbb{R}
}}_{(n-1)K+k-1} \md \widetilde{\mu}(\xi_{0}) \, \md \widetilde{\mu}(\xi_{1}) \, 
\dotsb \, \md \widetilde{\mu}(\xi_{l_{1}}) \, \dotsb \, \md \widetilde{\mu}
(\xi_{n_{1}}) \, \dotsb \, \md \widetilde{\mu}(\xi_{n_{2}}) \, \md \widetilde{\mu}
(\xi_{m_{1}}) \, \dotsb \, \md \widetilde{\mu}(\xi_{m_{2}}) \\
\times& \, \dfrac{(-1)^{(n-1)K+k-\varkappa_{nk}}}{(\xi_{0} \! - \! \alpha_{k}
)^{0}(\xi_{1} \! - \! \alpha_{p_{1}})^{1} \dotsb (\xi_{l_{1}} \! - \! \alpha_{
p_{1}})^{l_{1}} \dotsb \dfrac{1}{(\xi_{n_{1}})^{1}} \dotsb \dfrac{1}{(\xi_{
n_{2}})^{l_{\mathfrak{s}-1}}}(\xi_{m_{1}} \! - \! \alpha_{k})^{1} \dotsb 
(\xi_{m_{2}} \! - \! \alpha_{k})^{\varkappa_{nk}-2}} \\
\times& \, 
\left\vert

\right\vert.
\end{align*}
For $n \! \in \! \mathbb{N}$ and $k \! \in \! \lbrace 1,2,\dotsc,K \rbrace$ 
such that $\alpha_{p_{\mathfrak{s}}} \! := \! \alpha_{k} \! \neq \! \infty$, 
recalling the $(n \! - \! 1)K \! + \! k$ linearly independent functions 
$\tilde{\phi}_{0}(z) \! := \! \prod_{m=1}^{\mathfrak{s}-2}(z \! - \! 
\alpha_{p_{m}})^{l_{m}}(z \! - \! \alpha_{k})^{\varkappa_{nk}-1} \! 
=: \! \sum_{j=0}^{(n-1)K+k-1} \mathfrak{a}_{j,0}^{\sharp}z^{j}$, 
$\tilde{\phi}_{q(r_{1})}(z) \! := \! \lbrace \tilde{\phi}_{0}(z)(z \! - \! 
\alpha_{p_{r_{1}}})^{-m(r_{1})(1-\delta_{r_{1} \mathfrak{s}-1})}z^{m
(r_{1}) \delta_{r_{1} \mathfrak{s}-1}} \! =: \! \sum_{j=0}^{(n-1)K+k-1} 
\mathfrak{a}_{j,q(r_{1})}^{\sharp}z^{j} \rbrace$, $r_{1} \! = \! 1,\dotsc,
\mathfrak{s} \! - \! 1,\mathfrak{s}$, $q(r_{1}) \! = \! \sum_{i=1}^{r_{1}-1}
l_{i} \! + \! 1,\sum_{i=1}^{r_{1}-1}l_{i} \! + \! 2,\dotsc,\sum_{i=1}^{r_{1}-1}
l_{i} \! + \! l_{r_{1}} \! - \! \delta_{r_{1} \mathfrak{s}}$, $m(r_{1}) \! = \! 1,2,
\dotsc,l_{r_{1}} \! - \! \delta_{r_{1} \mathfrak{s}}$, and the $(n \! - \! 1)
K \! + \! k \! - \! 1$ linearly independent functions $\phi_{0}^{\spcheck}
(z) \! := \! \prod_{m=1}^{\mathfrak{s}-2}(z \! - \! \alpha_{p_{m}})^{l_{m}}
(z \! - \! \alpha_{k})^{\varkappa_{nk}-2} \! =: \! \sum_{j=0}^{(n-1)K+k-2} 
\mathfrak{a}^{\clubsuit}_{j,0}z^{j}$, $\phi_{q(r_{2})}^{\spcheck}(z) \! := \! 
\lbrace \phi_{0}^{\spcheck}(z)(z \! - \! \alpha_{p_{r_{2}}})^{-m(r_{2})(1-
\delta_{r_{2} \mathfrak{s}-1})}z^{m(r_{2}) \delta_{r_{2} \mathfrak{s}-1}} \! 
=: \! \sum_{j=0}^{(n-1)K+k-2} \mathfrak{a}^{\clubsuit}_{j,q(r_{2})}z^{j} 
\rbrace$, $r_{2} \! = \! 1,\dotsc,\mathfrak{s} \! - \! 1,\mathfrak{s}$, 
$q(r_{2}) \! = \! \sum_{i=1}^{r_{2}-1}l_{i} \! + \! 1,\sum_{i=1}^{r_{2}-1}l_{i} 
\! + \! 2,\dotsc,\sum_{i=1}^{r_{2}-1}l_{i} \! + \! l_{r_{2}} \! - \! 2 \delta_{r_{2} 
\mathfrak{s}}$, $m(r_{2}) \! = \! 1,2,\dotsc,l_{r_{2}} \! - \! 2 \delta_{r_{2} 
\mathfrak{s}}$, one proceeds, via the latter determinantal expression, 
with the analysis of 
$\overset{\triangledown}{\Xi}^{\raise-1.0ex\hbox{$\scriptstyle n$}}_{k}
(z) \! := \! (-1)^{(n-1)K+k-\varkappa_{nk}} 
\overset{\triangle}{\Xi}^{\raise-1.0ex\hbox{$\scriptstyle n$}}_{k}(z)$:
\begin{align*}
\overset{\triangle}{\Xi}^{\raise-1.0ex\hbox{$\scriptstyle n$}}_{k}(z) =& \, 
\underbrace{\int_{\mathbb{R}} \int_{\mathbb{R}} \dotsb \int_{\mathbb{R}}}_{(
n-1)K+k-1} \md \widetilde{\mu}(\xi_{0}) \, \md \widetilde{\mu}(\xi_{1}) \, 
\dotsb \, \md \widetilde{\mu}(\xi_{l_{1}}) \, \dotsb \, \md \widetilde{\mu}
(\xi_{n_{1}}) \, \dotsb \, \md \widetilde{\mu}(\xi_{n_{2}}) \, \md \widetilde{\mu}
(\xi_{m_{1}}) \, \dotsb \, \md \widetilde{\mu}(\xi_{m_{2}}) \\
\times& \, \dfrac{1}{(\xi_{0} \! - \! \alpha_{k})^{0}(\xi_{1} \! - \! \alpha_{
p_{1}})^{1} \dotsb (\xi_{l_{1}} \! - \! \alpha_{p_{1}})^{l_{1}} \dotsb \dfrac{
1}{(\xi_{n_{1}})^{1}} \dotsb \dfrac{1}{(\xi_{n_{2}})^{l_{\mathfrak{s}-1}}}
(\xi_{m_{1}} \! - \! \alpha_{k})^{1} \dotsb (\xi_{m_{2}} \! - \! 
\alpha_{k})^{\varkappa_{nk}-2}} \\
\times& \, \dfrac{(\tilde{\phi}_{0}(z))^{-1}}{\tilde{\phi}_{0}(\xi_{0}) 
\tilde{\phi}_{0}(\xi_{1}) \dotsb \tilde{\phi}_{0}(\xi_{l_{1}}) \dotsb \tilde{
\phi}_{0}(\xi_{n_{1}}) \dotsb \tilde{\phi}_{0}(\xi_{n_{2}}) \tilde{\phi}_{0}
(\xi_{m_{1}}) \dotsb \tilde{\phi}_{0}(\xi_{m_{2}})} \\
\times& \, 
\underbrace{\left\vert

\right\vert}_{=: \, \mathbb{G}^{\triangle}(\xi_{0},\xi_{1},\dotsc,\xi_{l_{1}},
\dotsc,\xi_{n_{1}},\dotsc,\xi_{n_{2}},\xi_{m_{1}},\dotsc,\xi_{(n-1)K+k-2};z)} \\
=& \, \dfrac{(\tilde{\phi}_{0}(z))^{-1}}{(m_{2} \! + \! 1)!} \sum_{\pmb{\sigma} 
\in \mathfrak{S}_{m_{2}+1}} \underbrace{\int_{\mathbb{R}} \int_{\mathbb{R}} 
\dotsb \int_{\mathbb{R}}}_{(n-1)K+k-1} \md \widetilde{\mu}(\xi_{\sigma (0)}) 
\, \md \widetilde{\mu}(\xi_{\sigma (1)}) \, \dotsb \, \md \widetilde{\mu}
(\xi_{\sigma (l_{1})}) \, \dotsb \, \md \widetilde{\mu}(\xi_{\sigma (n_{1})}) \, 
\dotsb \, \md \widetilde{\mu}(\xi_{\sigma (n_{2})}) \\
\times& \, \tfrac{\dotsb \, \md \widetilde{\mu}(\xi_{\sigma (m_{1})}) \, \dotsb 
\, \md \widetilde{\mu}(\xi_{\sigma (m_{2})})}{(\xi_{\sigma (0)}-\alpha_{k})^{0}
(\xi_{\sigma (1)}-\alpha_{p_{1}})^{1} \dotsb (\xi_{\sigma (l_{1})}-\alpha_{p_{1}})^{l_{1}} 
\dotsb \frac{1}{(\xi_{\sigma (n_{1})})^{1}} \dotsb \frac{1}{(\xi_{\sigma (n_{2})})^{
l_{\mathfrak{s}-1}}}(\xi_{\sigma (m_{1})}-\alpha_{k})^{1} \dotsb (\xi_{\sigma 
(m_{2})}-\alpha_{k})^{\varkappa_{nk}-2}} \\
\times& \, \dfrac{1}{\tilde{\phi}_{0}(\xi_{\sigma (0)}) \tilde{\phi}_{0}(\xi_{
\sigma (1)}) \dotsb \tilde{\phi}_{0}(\xi_{\sigma (l_{1})}) \dotsb \tilde{
\phi}_{0}(\xi_{\sigma (n_{1})}) \dotsb \tilde{\phi}_{0}(\xi_{\sigma (n_{2})}) 
\tilde{\phi}_{0}(\xi_{\sigma (m_{1})}) \dotsb \tilde{\phi}_{0}(\xi_{\sigma 
(m_{2})})} \\
\times& \, \mathbb{G}^{\triangle}(\xi_{\sigma (0)},\xi_{\sigma (1)},\dotsc,
\xi_{\sigma (l_{1})},\dotsc,\xi_{\sigma (n_{1})},\dotsc,\xi_{\sigma (n_{2})},
\xi_{\sigma (m_{1})},\dotsc,\xi_{\sigma ((n-1)K+k-2)};z) \\
=& \, \dfrac{(\tilde{\phi}_{0}(z))^{-1}}{(m_{2} \! + \! 1)!} \underbrace{
\int_{\mathbb{R}} \int_{\mathbb{R}} \dotsb \int_{\mathbb{R}}}_{(n-1)K+k-1} 
\md \widetilde{\mu}(\xi_{0}) \, \md \widetilde{\mu}(\xi_{1}) \, \dotsb \, \md 
\widetilde{\mu}(\xi_{l_{1}}) \, \dotsb \, \md \widetilde{\mu}(\xi_{n_{1}}) \, 
\dotsb \, \md \widetilde{\mu}(\xi_{n_{2}}) \, \md \widetilde{\mu}(\xi_{m_{1}}) 
\, \dotsb \\
\times& \, \dotsb \, \md \widetilde{\mu}(\xi_{m_{2}}) \, \mathbb{G}^{\triangle}
(\xi_{0},\xi_{1},\dotsc,\xi_{l_{1}},\dotsc,\xi_{n_{1}},\dotsc,\xi_{n_{2}},\xi_{m_{1}},
\dotsc,\xi_{(n-1)K+k-2};z) \\
\times& \, \sum_{\pmb{\sigma} \in \mathfrak{S}_{m_{2}+1}} \operatorname{sgn}
(\pmb{\pmb{\sigma}}) \tfrac{1}{(\xi_{\sigma (0)}-\alpha_{k})^{0}(\xi_{\sigma 
(1)}-\alpha_{p_{1}})^{1} \dotsb (\xi_{\sigma (l_{1})}-\alpha_{p_{1}})^{l_{1}} 
\dotsb \frac{1}{(\xi_{\sigma (n_{1})})^{1}} \dotsb \frac{1}{(\xi_{\sigma 
(n_{2})})^{l_{\mathfrak{s}-1}}}(\xi_{\sigma (m_{1})}-\alpha_{k})^{1} \dotsb 
(\xi_{\sigma (m_{2})}-\alpha_{k})^{\varkappa_{nk}-2}} \\
\times& \, \dfrac{1}{\tilde{\phi}_{0}(\xi_{\sigma (0)}) \tilde{\phi}_{0}(\xi_{
\sigma (1)}) \dotsb \tilde{\phi}_{0}(\xi_{\sigma (l_{1})}) \dotsb \tilde{
\phi}_{0}(\xi_{\sigma (n_{1})}) \dotsb \tilde{\phi}_{0}(\xi_{\sigma (n_{2})}) 
\tilde{\phi}_{0}(\xi_{\sigma (m_{1})}) \dotsb \tilde{\phi}_{0}(\xi_{\sigma 
(m_{2})})} \\
=& \, \dfrac{(\tilde{\phi}_{0}(z))^{-1}}{(m_{2} \! + \! 1)!} \underbrace{
\int_{\mathbb{R}} \int_{\mathbb{R}} \dotsb \int_{\mathbb{R}}}_{(n-1)K+k-1} 
\md \widetilde{\mu}(\xi_{0}) \, \md \widetilde{\mu}(\xi_{1}) \, \dotsb \, 
\md \widetilde{\mu}(\xi_{l_{1}}) \, \dotsb \, \md \widetilde{\mu}(\xi_{n_{1}}) 
\, \dotsb \, \md \widetilde{\mu}(\xi_{n_{2}}) \, \md \widetilde{\mu}
(\xi_{m_{1}}) \, \dotsb \\
\times& \, \dotsb \, \md \widetilde{\mu}(\xi_{m_{2}}) \, \dfrac{\mathbb{G}^{
\triangle}(\xi_{0},\xi_{1},\dotsc,\xi_{l_{1}},\dotsc,\xi_{n_{1}},\dotsc,\xi_{n_{2}},
\xi_{m_{1}},\dotsc,\xi_{(n-1)K+k-2};z)}{\tilde{\phi}_{0}(\xi_{0}) \tilde{\phi}_{0}
(\xi_{1}) \dotsb \tilde{\phi}_{0}(\xi_{l_{1}}) \dotsb \tilde{\phi}_{0}(\xi_{n_{1}}) 
\dotsb \tilde{\phi}_{0}(\xi_{n_{2}}) \tilde{\phi}_{0}(\xi_{m_{1}}) \dotsb 
\tilde{\phi}_{0}(\xi_{m_{2}})} \\
\times& \, 
\left\vert 

\right\vert \\
=& \, \dfrac{(\tilde{\phi}_{0}(z))^{-1}}{(m_{2} \! + \! 1)!} \underbrace{
\int_{\mathbb{R}} \int_{\mathbb{R}} \dotsb \int_{\mathbb{R}}}_{(n-1)K+k-1} 
\md \widetilde{\mu}(\xi_{0}) \, \md \widetilde{\mu}(\xi_{1}) \, \dotsb \, \md 
\widetilde{\mu}(\xi_{l_{1}}) \, \dotsb \, \md \widetilde{\mu}(\xi_{n_{1}}) \, 
\dotsb \, \md \widetilde{\mu}(\xi_{n_{2}}) \, \md \widetilde{\mu}(\xi_{m_{1}}) 
\, \dotsb \\
\times& \, \dotsb \, \md \widetilde{\mu}(\xi_{m_{2}}) \, \dfrac{\mathbb{G}^{
\triangle}(\xi_{0},\xi_{1},\dotsc,\xi_{l_{1}},\dotsc,\xi_{n_{1}},\dotsc,\xi_{n_{2}},
\xi_{m_{1}},\dotsc,\xi_{(n-1)K+k-2};z)}{\tilde{\phi}_{0}(\xi_{0}) \tilde{
\phi}_{0}(\xi_{1}) \dotsb \tilde{\phi}_{0}(\xi_{l_{1}}) \dotsb \tilde{\phi}_{
0}(\xi_{n_{1}}) \dotsb \tilde{\phi}_{0}(\xi_{n_{2}}) \tilde{\phi}_{0}(\xi_{
m_{1}}) \dotsb \tilde{\phi}_{0}(\xi_{m_{2}})} \\
\times& \, \dfrac{1}{\phi_{0}^{\spcheck}(\xi_{0}) \phi_{0}^{\spcheck}(\xi_{1}) 
\dotsb \phi_{0}^{\spcheck}(\xi_{l_{1}}) \dotsb \phi_{0}^{\spcheck}(\xi_{n_{1}}
) \dotsb \phi_{0}^{\spcheck}(\xi_{n_{2}}) \phi_{0}^{\spcheck}(\xi_{m_{1}}) 
\dotsb \phi_{0}^{\spcheck}(\xi_{m_{2}})} \\
\times& \, 
\underbrace{\left\vert

\right\vert}_{=: \, \mathbb{G}^{\clubsuit}(\xi_{0},\xi_{1},\dotsc,\xi_{l_{1}},
\dotsc,\xi_{n_{1}},\dotsc,\xi_{n_{2}},\dotsc,\xi_{(n-1)K+k-2})};
\end{align*}
but, noting the determinantal factorisation
\begin{equation*}
\mathbb{G}^{\clubsuit}(\xi_{0},\xi_{1},\dotsc,\xi_{l_{1}},\dotsc,
\xi_{n_{1}},\dotsc,\xi_{n_{2}},\dotsc,\xi_{(n-1)K+k-2}) \! := \! 
\overset{f}{\mathcal{V}}_{3}(\xi_{0},\xi_{1},\dotsc,\xi_{(n-1)K+k-2}) 
\mathbb{D}^{\clubsuit},
\end{equation*}
where
\begin{equation*}
\overset{f}{\mathcal{V}}_{3}(\xi_{0},\xi_{1},\dotsc,\xi_{(n-1)K+k-2}) =
\left\lvert
 
\right), \label{eq77}
\end{equation}
with
\begin{equation*}
(\mathfrak{A}^{\clubsuit}(0))_{1j} \! = \! 
\begin{cases}
\mathfrak{a}_{j-1}^{\clubsuit}(\vec{\bm{\alpha}}), &\text{$j \! = \! 1,2,
\dotsc,(n \! - \! 1)K \! + \! k \! - \! l_{\mathfrak{s}-1} \! - \! 1$,} \\
0, &\text{$j \! = \! (n \! - \! 1)K \! + \! k \! - \! l_{\mathfrak{s}-1},
\dotsc,(n \! - \! 1)K \! + \! k \! - \! 1$,}
\end{cases}
\end{equation*}
for $r \! = \! 1,\dotsc,\mathfrak{s} \! - \! 2,\mathfrak{s}$, with $l_{r} \! 
= \! \varkappa_{nk \tilde{k}_{r}}$, $r \! = \! 1,2,\dotsc,\mathfrak{s}-2$, 
$l_{\mathfrak{s}-1} \! = \! \varkappa^{\infty}_{nk \tilde{k}_{\mathfrak{s}-1}}$, 
$l_{\mathfrak{s}} \! = \! \varkappa_{nk}$, and $\sum_{m=1}^{\mathfrak{s}-2}
l_{m} \! + \! l_{\mathfrak{s}-1} \! = \! (n \! - \! 1)K \! + \! k \! - \! \varkappa_{nk}$,
\begin{align*}
i(r)=& \, 2 \! + \! \sum_{m=1}^{r-1}l_{m},3 \! + \! \sum_{m=1}^{r-1}l_{m},
\dotsc,1 \! - \! 2 \delta_{r \mathfrak{s}} \! + \! l_{r} \! + \! 
\sum_{m=1}^{r-1}l_{m}, \qquad \qquad q(i(r),r) \! = \! 1,2,\dotsc,
l_{r} \! - \! 2 \delta_{r \mathfrak{s}}, \\
(\mathfrak{A}^{\clubsuit}(r))_{i(r)j(r)}=& \, 
\begin{cases}
\dfrac{(-1)^{q(i(r),r)}}{\prod_{m=0}^{q(i(r),r)-1}(l_{r} \! - \! m \! - \! 2 
\delta_{r \mathfrak{s}})} 
\left(\dfrac{\partial}{\partial \alpha_{p_{r}}} \right)^{q(i(r),r)} 
\mathfrak{a}_{j(r)-1}^{\clubsuit}(\vec{\bm{\alpha}}), &\text{$j(r) \! = \! 
1,2,\dotsc,(n \! - \! 1)K \! + \! k \! - \! l_{\mathfrak{s}-1} \! - \! 
q(i(r),r) \! - \! 1$,} \\
0, &\text{$j(r) \! = \! (n \! - \! 1)K \! + \! k \! - \! l_{\mathfrak{s}-1} 
\! - \! q(i(r),r),\dotsc,(n \! - \! 1)K \! + \! k \! - \! 1$,}
\end{cases}
\end{align*}
and, for $r \! = \! \mathfrak{s} \! - \! 1$,
\begin{align*}
i(\mathfrak{s} \! - \! 1)=& \, 2 \! + \! \sum_{m=1}^{\mathfrak{s}-2}l_{m},
3 \! + \! \sum_{m=1}^{\mathfrak{s}-2}l_{m},\dotsc,1 \! + \! l_{\mathfrak{s}-1} 
\! + \! \sum_{m=1}^{\mathfrak{s}-2}l_{m}, \qquad \quad q(i(\mathfrak{s} \! 
- \! 1),\mathfrak{s} \! - \! 1) \! = \! 1,2,\dotsc,l_{\mathfrak{s}-1}, \\
(\mathfrak{A}^{\clubsuit}(\mathfrak{s} \! - \! 1))_{i(\mathfrak{s}-1)
j(\mathfrak{s}-1)}=& \, 
\begin{cases}
0, &\text{$j(\mathfrak{s} \! - \! 1) \! = \! 1,\dotsc,
q(i(\mathfrak{s} \! - \! 1),\mathfrak{s} \! - \! 1)$,} \\
\mathfrak{a}_{j(\mathfrak{s}-1)-q(i(\mathfrak{s}-1),\mathfrak{s}-1)
-1}^{\clubsuit}(\vec{\bm{\alpha}}), &\text{$j(\mathfrak{s} \! - \! 1) 
\! = \! q(i(\mathfrak{s} \! - \! 1),\mathfrak{s} \! - \! 1) \! + \! 1,
\dotsc,i(\mathfrak{s} \! - \! 1) \! + \! \varkappa_{nk} \! - \! 2$,} \\
0, &\text{$j(\mathfrak{s} \! - \! 1) \! = \! i(\mathfrak{s} \! - \! 1) \! 
+ \! \varkappa_{nk} \! - \! 1,\dotsc,(n \! - \! 1)K \! + \! k \! - \! 1$,}
\end{cases}
\end{align*}
where, for $\tilde{m}_{1} \! = \! 0,1,\dotsc,(n \! - \! 1)K \! + \! k \! - \! 
l_{\mathfrak{s}-1} \! - \! 2$,
\begin{equation*}
\mathfrak{a}_{\tilde{m}_{1}}^{\clubsuit}(\vec{\bm{\alpha}}) := 
\mathlarger{\sum_{\underset{\underset{\underset{\sum_{\underset{m \neq 
\mathfrak{s}-1}{m=1}}^{\mathfrak{s}}i_{m}=(n-1)K+k-l_{\mathfrak{s}-1}-
\tilde{m}_{1}-2}{i_{\mathfrak{s}}=0,1,\dotsc,l_{\mathfrak{s}}-2}}{p \in 
\lbrace 1,2,\dotsc,\mathfrak{s}-2 \rbrace}}{i_{p}=0,1,\dotsc,l_{p}}}}
(-1)^{(n-1)K+k-l_{\mathfrak{s}-1}-\tilde{m}_{1}-2} \prod_{j=1}^{\mathfrak{s}-2} 
\binom{l_{j}}{i_{j}} \binom{l_{\mathfrak{s}} \! - \! 2}{i_{\mathfrak{s}}} \prod_{m
=1}^{\mathfrak{s}-2}(\alpha_{p_{m}})^{i_{m}} \alpha_{k}^{i_{\mathfrak{s}}},
\end{equation*}
it follows that, for $n \! \in \! \mathbb{N}$ and $k \! \in \! \lbrace 
1,2,\dotsc,K \rbrace$ such that $\alpha_{p_{\mathfrak{s}}} \! := \! 
\alpha_{k} \! \neq \! \infty$,
\begin{align*}
\overset{\triangle}{\Xi}^{\raise-1.0ex\hbox{$\scriptstyle n$}}_{k}(z) =& \, 
\dfrac{\mathbb{D}^{\clubsuit}(\tilde{\phi}_{0}(z))^{-1}}{((n \! - \! 1)K \! + 
\! k \! - \! 1)!} \underbrace{\int_{\mathbb{R}} \int_{\mathbb{R}} \dotsb 
\int_{\mathbb{R}}}_{(n-1)K+k-1} \md \widetilde{\mu}(\xi_{0}) \, \md \widetilde{\mu}
(\xi_{1}) \, \dotsb \, \md \widetilde{\mu}(\xi_{l_{1}}) \, \dotsb \, \md \widetilde{\mu}
(\xi_{n_{1}}) \, \dotsb \, \md \widetilde{\mu}(\xi_{n_{2}}) \, \md \widetilde{\mu}
(\xi_{m_{1}}) \, \dotsb \, \md \widetilde{\mu}(\xi_{m_{2}}) \\
\times& \, \dfrac{\prod_{\underset{j<i}{i,j=0}}^{(n-1)K+k-2}(\xi_{i} 
\! - \! \xi_{j})}{\tilde{\phi}_{0}(\xi_{0}) \tilde{\phi}_{0}(\xi_{1}) \dotsb 
\tilde{\phi}_{0}(\xi_{l_{1}}) \dotsb \tilde{\phi}_{0}(\xi_{n_{1}}) \dotsb 
\tilde{\phi}_{0}(\xi_{n_{2}}) \tilde{\phi}_{0}(\xi_{m_{1}}) \dotsb 
\tilde{\phi}_{0}(\xi_{m_{2}})} \\
\times& \, \dfrac{1}{\phi_{0}^{\spcheck}(\xi_{0}) \phi_{0}^{\spcheck}(\xi_{1}) 
\dotsb \phi_{0}^{\spcheck}(\xi_{l_{1}}) \dotsb \phi_{0}^{\spcheck}(\xi_{n_{1}}
) \dotsb \phi_{0}^{\spcheck}(\xi_{n_{2}}) \phi_{0}^{\spcheck}(\xi_{m_{1}}) 
\dotsb \phi_{0}^{\spcheck}(\xi_{m_{2}})} \\
\times& \, 
\underbrace{\left\vert

\right\vert}_{= \, \mathbb{D}^{\blacklozenge} \prod_{\underset{j<i}{i,j
=0}}^{(n-1)K+k-2}(\xi_{i}-\xi_{j}) \prod_{m=0}^{(n-1)K+k-2}(z-\xi_{m}) 
\quad \text{(cf. proof of Lemma~\ref{lem2.1}, case~\textbf{(2)})}} \\
=& \, \dfrac{\mathbb{D}^{\clubsuit} \mathbb{D}^{\blacklozenge}(\tilde{\phi}_{0}
(z))^{-1}}{((n \! - \! 1)K \! + \! k \! - \! 1)!} \underbrace{\int_{\mathbb{R}} 
\int_{\mathbb{R}} \dotsb \int_{\mathbb{R}}}_{(n-1)K+k-1} \md \widetilde{\mu}
(\xi_{0}) \, \md \widetilde{\mu}(\xi_{1}) \, \dotsb \, \md \widetilde{\mu}
(\xi_{l_{1}}) \, \dotsb \, \md \widetilde{\mu}(\xi_{n_{1}}) \, \dotsb \, \md 
\widetilde{\mu}(\xi_{n_{2}}) \, \md \widetilde{\mu}(\xi_{m_{1}}) \, \dotsb \\
\times& \, \dotsb \, \md \widetilde{\mu}(\xi_{m_{2}}) \, \dfrac{(\phi_{0}^{\spcheck}
(\xi_{0}) \phi_{0}^{\spcheck}(\xi_{1}) \dotsb \phi_{0}^{\spcheck}(\xi_{l_{1}}) 
\dotsb \phi_{0}^{\spcheck}(\xi_{n_{1}}) \dotsb \phi_{0}^{\spcheck}(\xi_{n_{2}}
) \phi_{0}^{\spcheck}(\xi_{m_{1}}) \dotsb \phi_{0}^{\spcheck}(\xi_{m_{2}}))^{-
1}}{\tilde{\phi}_{0}(\xi_{0}) \tilde{\phi}_{0}(\xi_{1}) \dotsb \tilde{\phi}_{
0}(\xi_{l_{1}}) \dotsb \tilde{\phi}_{0}(\xi_{n_{1}}) \dotsb \tilde{\phi}_{0}
(\xi_{n_{2}}) \tilde{\phi}_{0}(\xi_{m_{1}}) \dotsb \tilde{\phi}_{0}(\xi_{m_{
2}})} \\
\times& \, \prod_{\substack{i,j=0\\j<i}}^{(n-1)K+k-2}(\xi_{i} \! - \! 
\xi_{j})^{2} \, \prod_{m=0}^{(n-1)K+k-2}(z \! - \! \xi_{m});
\end{align*}
but, noting that, for $n \! \in \! \mathbb{N}$ and $k \! \in \! \lbrace 1,2,
\dotsc,K \rbrace$ such that $\alpha_{p_{\mathfrak{s}}} \! := \! \alpha_{k} 
\! \neq \! \infty$, $\phi_{0}^{\spcheck}(z) \! = \! \tilde{\phi}_{0}
(z)/(z \! - \! \alpha_{k})$, and 
$\overset{\triangledown}{\Xi}^{\raise-1.0ex\hbox{$\scriptstyle n$}}_{k}
(z) \! = \! (-1)^{(n-1)K+k-\varkappa_{nk}} 
\overset{\triangle}{\Xi}^{\raise-1.0ex\hbox{$\scriptstyle n$}}_{k}(z)$, 
one arrives at
\begin{align*}
\overset{\triangledown}{\Xi}^{\raise-1.0ex\hbox{$\scriptstyle n$}}_{k}(z) =& 
\, \dfrac{(-1)^{(n-1)K+k-\varkappa_{nk}} \mathbb{D}^{\clubsuit} \mathbb{D}^{
\blacklozenge}(\tilde{\phi}_{0}(z))^{-1}}{((n \! - \! 1)K \! + \! k \! - \! 1)!} 
\underbrace{\int_{\mathbb{R}} \int_{\mathbb{R}} \dotsb \int_{\mathbb{R}
}}_{(n-1)K+k-1} \md \widetilde{\mu}(\xi_{0}) \, \md \widetilde{\mu}(\xi_{1}) 
\, \dotsb \, \md \widetilde{\mu}(\xi_{l_{1}}) \, \dotsb \, \md \widetilde{\mu}
(\xi_{n_{1}}) \, \dotsb \, \md \widetilde{\mu}(\xi_{n_{2}}) \, \dotsb \\
\times& \, \dotsb \, \md \widetilde{\mu}(\xi_{(n-1)K+k-2}) \left(\prod_{m=
0}^{(n-1)K+k-2} \tilde{\phi}_{0}(\xi_{m}) \right)^{-2} \prod_{\substack{i,j=
0\\j<i}}^{(n-1)K+k-2}(\xi_{i} \! - \! \xi_{j})^{2} \, \prod_{l=0}^{(n-1)K+k-2}
(\xi_{l} \! - \! \alpha_{k})(z \! - \! \xi_{l}) \, \, \Rightarrow \\
\overset{\triangledown}{\Xi}^{\raise-1.0ex\hbox{$\scriptstyle n$}}_{k}(z) 
=& \, \dfrac{(-1)^{(n-1)K+k-\varkappa_{nk}} \mathbb{D}^{\clubsuit} 
\mathbb{D}^{\blacklozenge}}{((n \! - \! 1)K \! + \! k \! - \! 1)!} \left(
\prod_{q=1}^{\mathfrak{s}-2}(z \! - \! \alpha_{p_{q}})^{\varkappa_{nk 
\tilde{k}_{q}}}(z \! - \! \alpha_{k})^{\varkappa_{nk}-1} \right)^{-1} 
\underbrace{\int_{\mathbb{R}} \int_{\mathbb{R}} \dotsb \int_{\mathbb{R}}}_{
(n-1)K+k-1} \md \widetilde{\mu}(\xi_{0}) \, \md \widetilde{\mu}(\xi_{1}) 
\, \dotsb \, \md \widetilde{\mu}(\xi_{(n-1)K+k-2}) \\
\times& \, \left(\prod_{m=0}^{(n-1)K+k-2} \prod_{q=1}^{\mathfrak{s}-2}
(\xi_{m} \! - \! \alpha_{p_{q}})^{\varkappa_{nk \tilde{k}_{q}}}
(\xi_{m} \! - \! \alpha_{k})^{\varkappa_{nk}-1} \right)^{-2} 
\prod_{\substack{i,j=0\\j<i}}^{(n-1)K+k-2}(\xi_{i} \! - \! \xi_{j})^{2} \, 
\prod_{l=0}^{(n-1)K+k-2}(\xi_{l} \! - \! \alpha_{k})(z \! - \! \xi_{l}).
\end{align*}
Recalling that, for $n \! \in \! \mathbb{N}$ and $k \! \in \! \lbrace 1,2,\dotsc,
K \rbrace$ such that $\alpha_{p_{\mathfrak{s}}} \! := \! \alpha_{k} \! \neq 
\! \infty$, $\mathcal{X}_{21}(z)/(z \! - \! \alpha_{k}) \! = \! 2 \pi \mi 
\overset{\triangledown}{\Xi}^{\raise-1.0ex\hbox{$\scriptstyle n$}}_{k}(z)/
\mathcal{N}_{f}^{\sharp}(n,k)$, where 
$\overset{\triangledown}{\Xi}^{\raise-1.0ex\hbox{$\scriptstyle n$}}_{k}(z)$ 
is given directly above, and $\mathcal{N}_{f}^{\sharp}(n,k)$ is given by 
Equation~\eqref{eq75}, one arrives at the integral representation for 
$\mathcal{X}_{21}(z)/(z \! - \! \alpha_{k})$ given in Equation~\eqref{eintreppfin2}. 
\hfill $\qed$
\begin{fffff} \label{cor2.1} 
Let $\widetilde{V} \colon \overline{\mathbb{R}} \setminus \lbrace 
\alpha_{1},\alpha_{2},\dotsc,\alpha_{K} \rbrace \! \to \! \mathbb{R}$ 
satisfy conditions~\eqref{eq20}--\eqref{eq22}. Let $\mathcal{X} \colon 
\mathbb{N} \times \lbrace 1,2,\dotsc,K \rbrace \times \overline{\mathbb{C}} 
\setminus \overline{\mathbb{R}} \! \to \! \operatorname{SL}_{2}
(\mathbb{C})$ be the unique solution of the {\rm RHP} stated in 
Lemma~$\bm{\mathrm{RHP}_{\mathrm{MPC}}}$ with integral representation 
given by Equation~\eqref{intrepinf} for $\alpha_{p_{\mathfrak{s}}} \! 
:= \! \alpha_{k} \! = \! \infty$ (resp., Equation~\eqref{intrepfin} for 
$\alpha_{p_{\mathfrak{s}}} \! := \! \alpha_{k} \! \neq \! \infty)$, 
$k \! \in \! \lbrace 1,2,\dotsc,K \rbrace$. For $n \! \in \! \mathbb{N}$ 
and $k \! \in \! \lbrace 1,2,\dotsc,K \rbrace$ such that 
$\alpha_{p_{\mathfrak{s}}} \! := \! \alpha_{k} \! = \! \infty$ (resp., 
$\alpha_{p_{\mathfrak{s}}} \! := \! \alpha_{k} \! \neq \! \infty)$, let 
$\mu^{\infty}_{n,\varkappa_{nk}}(n,k)$ (resp., $\mu^{f}_{n,\varkappa_{nk}}
(n,k))$ be the associated {\rm MPC ORF} norming constant defined in 
Subsection~\ref{subsubsec1.2.1} (resp., Subsection~\ref{subsubsec1.2.2}$)$. 
Then{\rm :} {\rm (i)} for $n \! \in \! \mathbb{N}$ and $k \! \in \! \lbrace 
1,2,\dotsc,K \rbrace$ such that $\alpha_{p_{\mathfrak{s}}} \! := \! \alpha_{k} 
\! = \! \infty$,\footnote{One obtains $\mu^{\infty}_{n,\varkappa_{nk}}(n,k)$ 
by taking the positive square root of the expression on the right-hand side 
of Equation~\eqref{eqnmctinf1}.}
\begin{equation} \label{eqnmctinf1} 
\left(\mu^{\infty}_{n,\varkappa_{nk}}(n,k) \right)^{2} \! = \! 
\dfrac{\varsigma_{\infty}(n,k)}{\lambda_{\infty}(n,k)},
\end{equation}
where
\begin{align*}
\varsigma_{\infty}(n,k) =& \, \dfrac{(\mathbb{D}^{\sphat} \,)^{2}}{((n \! - 
\! 1)K \! + \! k)!} \underbrace{\int_{\mathbb{R}} \int_{\mathbb{R}} \dotsb 
\int_{\mathbb{R}}}_{(n-1)K+k} \md \widetilde{\mu}(\xi_{0}) \, \md \widetilde{\mu}
(\xi_{1}) \, \dotsb \, \md \widetilde{\mu}(\xi_{(n-1)K+k-1}) \, \prod_{\substack{i,
j=0\\j<i}}^{(n-1)K+k-1}(\xi_{i} \! - \! \xi_{j})^{2} \\
\times& \, \left(\prod_{m=0}^{(n-1)K+k-1} \prod_{q=1}^{\mathfrak{s}-1}
(\xi_{m} \! - \! \alpha_{p_{q}})^{\varkappa_{nk \tilde{k}_{q}}} \right)^{-2}, \\
\lambda_{\infty}(n,k) =& \, \dfrac{(\mathbb{D}^{\spcheck})^{2}}{((n \! - \! 
1)K \! + \! k \! + \! 1)!} \underbrace{\int_{\mathbb{R}} \int_{\mathbb{R}} 
\dotsb \int_{\mathbb{R}}}_{(n-1)K+k+1} \md \widetilde{\mu}(\tau_{0}) \, \md 
\widetilde{\mu}(\tau_{1}) \, \dotsb \, \md \widetilde{\mu}(\tau_{(n-1)K+k}) 
\, \prod_{\substack{i,j=0\\j<i}}^{(n-1)K+k}(\tau_{i} \! - \! \tau_{j})^{2} \\
\times& \, \left(\prod_{m=0}^{(n-1)K+k} \prod_{q=1}^{\mathfrak{s}-1}
(\tau_{m} \! - \! \alpha_{p_{q}})^{\varkappa_{nk \tilde{k}_{q}}} \right)^{-2},
\end{align*}
with $\md \widetilde{\mu}(z) \! = \! \exp (-n \widetilde{V}(z)) \, \md z$, 
and $\mathbb{D}^{\spcheck}$ and $\mathbb{D}^{\sphat}$, with 
$\mathbb{D}^{\spcheck} \! = \! \mathbb{D}^{\sphat}$, are defined by 
Equations~\eqref{eq36} and~\eqref{eq38}, respectively$;$ and {\rm (ii)} for 
$n \! \in \! \mathbb{N}$ and $k \! \in \! \lbrace 1,2,\dotsc,K \rbrace$ such that 
$\alpha_{p_{\mathfrak{s}}} \! := \! \alpha_{k} \! \neq \! \infty$,\footnote{One 
obtains $\mu^{f}_{n,\varkappa_{nk}}(n,k)$ by taking the positive square root 
of the expression on the right-hand side of Equation~\eqref{eqnmctfin1}.}
\begin{equation} \label{eqnmctfin1} 
\left(\mu^{f}_{n,\varkappa_{nk}}(n,k) \right)^{2} \! = \! 
\dfrac{\varsigma_{f}(n,k)}{\lambda_{f}(n,k)},
\end{equation}
where
\begin{align*}
\varsigma_{f}(n,k) =& \, \dfrac{(\mathbb{D}^{\blacklozenge})^{2}}{((n \! - 
\! 1)K \! + \! k)!} \underbrace{\int_{\mathbb{R}} \int_{\mathbb{R}} \dotsb 
\int_{\mathbb{R}}}_{(n-1)K+k} \md \widetilde{\mu}(\xi_{0}) \, \md \widetilde{\mu}
(\xi_{1}) \, \dotsb \, \md \widetilde{\mu}(\xi_{(n-1)K+k-1}) \prod_{\substack{i,j
=0\\j<i}}^{(n-1)K+k-1}(\xi_{i} \! - \! \xi_{j})^{2} \\
\times& \, \left(\prod_{m=0}^{(n-1)K+k-1} \prod_{q=1}^{\mathfrak{s}-2}
(\xi_{m} \! - \! \alpha_{p_{q}})^{\varkappa_{nk \tilde{k}_{q}}}
(\xi_{m} \! - \! \alpha_{k})^{\varkappa_{nk}-1} \right)^{-2}, \\
\lambda_{f}(n,k) =& \, \dfrac{(\mathbb{D}^{\lozenge})^{2}}{((n \! - \! 1)K 
\! + \! k \! + \! 1)!} \underbrace{\int_{\mathbb{R}} \int_{\mathbb{R}} \dotsb 
\int_{\mathbb{R}}}_{(n-1)K+k+1} \md \widetilde{\mu}(\tau_{0}) \, \md \widetilde{\mu}
(\tau_{1}) \, \dotsb \, \md \widetilde{\mu}(\tau_{(n-1)K+k}) \, \prod_{\substack{i,j
=0\\j<i}}^{(n-1)K+k}(\tau_{i} \! - \! \tau_{j})^{2} \\
\times& \, \left(\prod_{m=0}^{(n-1)K+k} \prod_{q=1}^{\mathfrak{s}-2}
(\tau_{m} \! - \! \alpha_{p_{q}})^{\varkappa_{nk \tilde{k}_{q}}}
(\tau_{m} \! - \! \alpha_{k})^{\varkappa_{nk}} \right)^{-2},
\end{align*}
and $\mathbb{D}^{\lozenge}$ and $\mathbb{D}^{\blacklozenge}$, with 
$\mathbb{D}^{\lozenge}/\mathbb{D}^{\blacklozenge} \! = \! (-1)^{(n-1)K+k} 
\prod_{q=1}^{\mathfrak{s}-2}(\alpha_{k} \! - \! \alpha_{p_{q}})^{\varkappa_{nk 
\tilde{k}_{q}}}$, are defined by Equations~\eqref{eq62} and~\eqref{eq64}, 
respectively.
\end{fffff}

\emph{Proof}. Recall {}from the proof of Lemma~\ref{lem2.1}, 
Equation~\eqref{nhormmconstatinf} (resp., \eqref{nhormmconstatfin}) that, 
for $n \! \in \! \mathbb{N}$ and $k \! \in \! \lbrace 1,2,\dotsc,K \rbrace$ 
such that $\alpha_{p_{\mathfrak{s}}} \! := \! \alpha_{k} \! = \! \infty$ 
(resp., $\alpha_{p_{\mathfrak{s}}} \! := \! \alpha_{k} \! \neq \! \infty)$, 
$(\mu^{\infty}_{n,\varkappa_{nk}}(n,k))^{-2} \! = \! \hat{\mathfrak{c}}_{
N_{\infty}}/\hat{\mathfrak{c}}_{D_{\infty}}$ (resp., $(\mu^{f}_{n,\varkappa_{nk}}
(n,k))^{-2} \! = \! \hat{\mathfrak{c}}_{N_{f}}/\hat{\mathfrak{c}}_{D_{f}})$, 
where $\hat{\mathfrak{c}}_{N_{\infty}}$ (resp., $\hat{\mathfrak{c}}_{N_{f}})$ 
is given in Equation~\eqref{eq37} (resp., \eqref{eq63}), with the associated, 
non-singular determinant $\mathbb{D}^{\spcheck}$ (resp., 
$\mathbb{D}^{\lozenge})$ defined in Equation~\eqref{eq36} 
(resp., \eqref{eq62}), and $\hat{\mathfrak{c}}_{D_{\infty}}$ (resp., 
$\hat{\mathfrak{c}}_{D_{f}})$ is given in Equation~\eqref{eq39} (resp., 
\eqref{eq65}), with the associated, non-singular determinant 
$\mathbb{D}^{\sphat}$ (resp., $\mathbb{D}^{\blacklozenge})$ defined 
by Equation~\eqref{eq38} (resp., \eqref{eq64}); hence, substituting 
Equations~\eqref{eq37} and~\eqref{eq39} (resp., \eqref{eq63} and~\eqref{eq65}) 
into the formula for $(\mu^{\infty}_{n,\varkappa_{nk}}(n,k))^{-2}$ (resp., 
$(\mu^{f}_{n,\varkappa_{nk}}(n,k))^{-2})$, one arrives at, for $n \! \in \! 
\mathbb{N}$ and $k \! \in \! \lbrace 1,2,\dotsc,K \rbrace$ such that 
$\alpha_{p_{\mathfrak{s}}} \! := \! \alpha_{k} \! = \! \infty$ (resp., 
$\alpha_{p_{\mathfrak{s}}} \! := \! \alpha_{k} \! \neq \! \infty)$, 
Equation~\eqref{eqnmctinf1} (resp., \eqref{eqnmctfin1}). \hfill $\qed$
\begin{fffff} \label{cor2.2} 
For $n \! \in \! \mathbb{N}$ and $k \! \in \! \lbrace 1,2,\dotsc,K \rbrace$ 
such that $\alpha_{p_{\mathfrak{s}}} \! := \! \alpha_{k} \! = \! \infty$ 
(resp., $\alpha_{p_{\mathfrak{s}}} \! := \! \alpha_{k} \! \neq \! \infty)$, 
the {\rm MPC ORF} is obtained via $\phi^{n}_{k}(z) \! = \! \mu^{\infty}_{n,
\varkappa_{nk}}(n,k) \pmb{\pi}^{n}_{k}(z)$ (resp., $\phi^{n}_{k}(z) \! = \! 
\mu^{f}_{n,\varkappa_{nk}}(n,k) \pmb{\pi}^{n}_{k}(z))$, $z \! \in \! 
\mathbb{C}$, where the corresponding monic {\rm MPC ORF}, 
$\pmb{\pi}^{n}_{k}(z)$, $z \! \in \! \mathbb{C}$, is given in 
Equation~\eqref{eintreppinf1} (resp., \eqref{eintreppfin1}$)$, and 
the associated norming constant $\mu^{\infty}_{n,\varkappa_{nk}}
(n,k)$ (resp., $\mu^{f}_{n,\varkappa_{nk}}(n,k))$ is given in 
Equation~\eqref{eqnmctinf1} (resp., \eqref{eqnmctfin1}$)$.
\end{fffff}

\emph{Proof}. Consequence of the definition of the corresponding MPC ORF in 
terms of the associated norming constant and the monic MPC ORF given in 
Subsection~\ref{subsubsec1.2.1} (resp., Subsection~\ref{subsubsec1.2.2}) 
for the case $\alpha_{p_{\mathfrak{s}}} \! := \! \alpha_{k} \! = \! \infty$ 
(resp., $\alpha_{p_{\mathfrak{s}}} \! := \! \alpha_{k} \! \neq \! \infty)$, 
$k \! \in \! \lbrace 1,2,\dotsc,K \rbrace$. \hfill $\qed$
\begin{eeeee} \label{remcoeffs} 
\textsl{Recall the following formulae for the corresponding monic {\rm MPC ORFs} 
(cf. the proof of Lemma~\ref{lem2.1}$)$$:$ {\rm (i)} for $n \! \in \! \mathbb{N}$ 
and $k \! \in \! \lbrace 1,2,\dotsc,K \rbrace$ such that $\alpha_{p_{\mathfrak{s}}} 
\! := \! \alpha_{k} \! = \! \infty$,
\begin{align*}
\pmb{\pi}^{n}_{k}(z) =& \, \dfrac{\phi^{\infty}_{0}(n,k)}{\mu^{\infty}_{n,
\varkappa_{nk}}(n,k)} \! + \! \dfrac{1}{\mu^{\infty}_{n,\varkappa_{nk}}(n,k)} 
\sum_{m=1}^{\mathfrak{s}-1} \sum_{q=1}^{l_{m}} \dfrac{\nu^{\infty}_{q,m}
(n,k)}{(z \! - \! \alpha_{p_{m}})^{q}} \! + \! \dfrac{1}{\mu^{\infty}_{n,
\varkappa_{nk}}(n,k)} \sum_{q=1}^{l_{\mathfrak{s}}-1} \mu^{\infty}_{n,q}
(n,k)z^{q} \! + \! z^{\varkappa_{nk}} \\
:=& \, \widetilde{\phi}_{0}^{\raise-1.0ex\hbox{$\scriptstyle \infty$}}
(n,k) \! + \! \sum_{m=1}^{\mathfrak{s}-1} \sum_{q=1}^{l_{m}} 
\dfrac{\widetilde{\nu}^{\raise-1.0ex\hbox{$\scriptstyle \infty$}}_{m,q}
(n,k)}{(z \! - \! \alpha_{p_{m}})^{q}} \! + \! \sum_{q=1}^{l_{\mathfrak{s}}-1} 
\widetilde{\nu}^{\raise-1.0ex\hbox{$\scriptstyle \infty$}}_{\mathfrak{s},q}
(n,k)z^{q} \! + \! z^{\varkappa_{nk}},
\end{align*}
where $\widetilde{\phi}_{0}^{\raise-1.0ex\hbox{$\scriptstyle \infty$}}(n,k) 
\! = \! \phi^{\infty}_{0}(n,k)/\mu^{\infty}_{n,\varkappa_{nk}}(n,k)$, 
$\widetilde{\nu}^{\raise-1.0ex\hbox{$\scriptstyle \infty$}}_{m,q}(n,k) \! 
= \! \nu^{\infty}_{q,m}(n,k)/\mu^{\infty}_{n,\varkappa_{nk}}(n,k)$, $q \! = 
\! 1,2,\dotsc,l_{m} \! = \! \varkappa_{nk \tilde{k}_{m}}$, $m \! = \! 1,2,
\dotsc,\mathfrak{s} \! - \! 1$, 
$\widetilde{\nu}^{\raise-1.0ex\hbox{$\scriptstyle \infty$}}_{\mathfrak{s},
q^{\prime}}(n,k) \! = \! \mu^{\infty}_{n,q^{\prime}}(n,k)/\mu^{\infty}_{n,
\varkappa_{nk}}(n,k)$, $q^{\prime} \! = \! 1,2,\dotsc,l_{\mathfrak{s}} \! - 
\! 1$, with $l_{\mathfrak{s}} \! = \! \varkappa_{nk}$, and $\mu^{\infty}_{n,
\varkappa_{nk}}(n,k)$ is given in Equation~\eqref{eqnmctinf1}$;$ and 
{\rm (ii)} for $n \! \in \! \mathbb{N}$ and $k \! \in \! \lbrace 1,2,\dotsc,
K \rbrace$ such that $\alpha_{p_{\mathfrak{s}}} \! := \! \alpha_{k} \! 
\neq \! \infty$,
\begin{align*}
\pmb{\pi}^{n}_{k}(z) =& \, \dfrac{\phi^{f}_{0}(n,k)}{\mu^{f}_{n,
\varkappa_{nk}}(n,k)} \! + \! \dfrac{1}{\mu^{f}_{n,\varkappa_{nk}}(n,k)} 
\sum_{m=1}^{\mathfrak{s}-2} \sum_{q=1}^{l_{m}} \dfrac{\tilde{\nu}^{f}_{q,m}
(n,k)}{(z \! - \! \alpha_{p_{m}})^{q}} \! + \! \dfrac{1}{\mu^{f}_{n,
\varkappa_{nk}}(n,k)} \sum_{q=1}^{l_{\mathfrak{s}-1}} \hat{\nu}^{f}_{n,q}
(n,k)z^{q} \! + \! \dfrac{1}{\mu^{f}_{n,\varkappa_{nk}}(n,k)} \sum_{q=1}^{
l_{\mathfrak{s}}-1} \dfrac{\mu^{f}_{n,q}(n,k)}{(z \! - \! \alpha_{k})^{q}} \! 
+ \! \dfrac{1}{(z \! - \! \alpha_{k})^{\varkappa_{nk}}} \\
:=& \, \widetilde{\phi}^{\raise-1.0ex\hbox{$\scriptstyle f$}}_{0}
(n,k) \! + \! \sum_{m=1}^{\mathfrak{s}-2} \sum_{q=1}^{l_{m}} 
\dfrac{\widetilde{\nu}_{m,q}^{\raise-1.0ex\hbox{$\scriptstyle f$}}(n,k)}{
(z \! - \! \alpha_{p_{m}})^{q}} \! + \! \sum_{q=1}^{l_{\mathfrak{s}-1}} 
\widetilde{\nu}_{\mathfrak{s}-1,q}^{\raise-1.0ex\hbox{$\scriptstyle f$}}
(n,k)z^{q} \! + \! \sum_{q=1}^{l_{\mathfrak{s}}-1} \dfrac{\widetilde{\nu}_{
\mathfrak{s},q}^{\raise-1.0ex\hbox{$\scriptstyle f$}}(n,k)}{(z \! - \! 
\alpha_{k})^{q}} \! + \! \dfrac{1}{(z \! - \! \alpha_{k})^{\varkappa_{nk}}},
\end{align*}
where $\widetilde{\phi}^{\raise-1.0ex\hbox{$\scriptstyle f$}}_{0}(n,k) 
\! = \! \phi^{f}_{0}(n,k)/\mu^{f}_{n,\varkappa_{nk}}(n,k)$, 
$\widetilde{\nu}_{m,q}^{\raise-1.0ex\hbox{$\scriptstyle f$}}(n,k) \! = \! 
\tilde{\nu}^{f}_{q,m}(n,k)/\mu^{f}_{n,\varkappa_{nk}}(n,k)$, $q \! = \! 1,2,
\dotsc,l_{m} \! = \! \varkappa_{nk \tilde{k}_{m}}$, $m \! = \! 1,2,\dotsc,
\mathfrak{s} \! - \! 2$, 
$\widetilde{\nu}_{\mathfrak{s}-1,q^{\prime}}^{\raise-1.0ex\hbox{$\scriptstyle f$}}
(n,k) \! = \! \hat{\nu}^{f}_{n,q^{\prime}}(n,k)/\mu^{f}_{n,\varkappa_{nk}}(n,k)$, 
$q^{\prime} \! = \! 1,2,\dotsc,l_{\mathfrak{s}-1} \! = \! \varkappa^{\infty}_{nk 
\tilde{k}_{\mathfrak{s}-1}}$, $\widetilde{\nu}_{\mathfrak{s},
q^{\prime \prime}}^{\raise-1.0ex\hbox{$\scriptstyle f$}}(n,k) \! = \! 
\mu^{f}_{n,q^{\prime \prime}}(n,k)/\mu^{f}_{n,\varkappa_{nk}}(n,k)$, 
$q^{\prime \prime} \! = \! 1,2,\dotsc,l_{\mathfrak{s}} \! - \! 1$, with 
$l_{\mathfrak{s}} \! = \! \varkappa_{nk}$, and $\mu^{f}_{n,\varkappa_{nk}}
(n,k)$ is given in Equation~\eqref{eqnmctfin1}. In conjunction with the above 
formulae, an algebraic analysis of the multi-integral representations for the 
corresponding monic {\rm MPC ORFs} given in Equations~\eqref{eintreppinf1} 
and~\eqref{eintreppfin1}, respectively, reveals that: {\rm (i)} for $n \! \in 
\! \mathbb{N}$ and $k \! \in \! \lbrace 1,2,\dotsc,K \rbrace$ such that 
$\alpha_{p_{\mathfrak{s}}} \! := \! \alpha_{k} \! = \! \infty$,
\begin{gather*}
\widetilde{\phi}_{0}^{\raise-1.0ex\hbox{$\scriptstyle \infty$}}(n,k) \! = 
\! \dfrac{\phi^{\infty}_{0}(n,k)}{\mu^{\infty}_{n,\varkappa_{nk}}(n,k)} \! 
= \! \dfrac{1}{\Upsilon_{D_{\infty}}(n,k)} \left(\mathfrak{b}_{(n-1)K+k-
l_{\mathfrak{s}}} \! + \! \sum_{m=1}^{l_{\mathfrak{s}}} 
\mathfrak{b}^{\natural}_{m} \mathfrak{b}_{m+(n-1)K+k-l_{\mathfrak{s}}} 
\right), \\
\widetilde{\nu}^{\raise-1.0ex\hbox{$\scriptstyle \infty$}}_{\mathfrak{s},
q^{\prime}}(n,k) \! = \! \dfrac{\mu^{\infty}_{n,q^{\prime}}(n,k)}{
\mu^{\infty}_{n,\varkappa_{nk}}(n,k)} = \! \dfrac{1}{\Upsilon_{D_{\infty}}
(n,k)} \left(\mathfrak{b}_{q^{\prime}+(n-1)K+k-l_{\mathfrak{s}}} \! + \! 
\mathlarger{\sum_{\underset{m_{1}-m_{2}=q^{\prime}}{m_{1},m_{2}=1,2,\dotsc,
l_{\mathfrak{s}}}}^{}} \mathfrak{b}^{\natural}_{m_{2}} \mathfrak{b}_{m_{1}
+(n-1)K+k-l_{\mathfrak{s}}} \right), \quad q^{\prime} \! = \! 1,2,\dotsc,
l_{\mathfrak{s}} \! - \! 1,
\end{gather*}
where $\Upsilon_{D_{\infty}}(n,k)$ is given in item~{\rm (i)} of 
Lemma~\ref{lem2.2},
\begin{align*}
\mathfrak{b}_{r} =& \, \underbrace{\int_{\mathbb{R}} \int_{\mathbb{R}} 
\dotsb \int_{\mathbb{R}}}_{(n-1)K+k} \md \widetilde{\mu}(\xi_{0}) \, \md 
\widetilde{\mu}(\xi_{1}) \, \dotsb \, \md \widetilde{\mu}(\xi_{(n-1)K+k-1}) 
\prod_{\substack{i,j=0\\j<i}}^{(n-1)K+k-1}(\xi_{i} \! - \! \xi_{j})^{2} \left(
\prod_{m=0}^{(n-1)K+k-1} \prod_{q=1}^{\mathfrak{s}-1}(\xi_{m} \! - \! 
\alpha_{p_{q}})^{l_{q}} \right)^{-2} \\
\times& \, \mathlarger{\sum_{\underset{\underset{\sum_{m^{\prime}=0}^{(n-1)K
+k-1}i_{m^{\prime}}=(n-1)K+k-r}{p=0,1,\dotsc,(n-1)K+k-1}}{i_{p} \in \lbrace 
0,1 \rbrace}}^{}} \prod_{q^{\prime}=0}^{(n-1)K+k-1} \binom{1}{i_{q^{\prime}}}
(-1)^{i_{q^{\prime}}} \xi_{q^{\prime}}^{i_{q^{\prime}}}, \quad r \! = \! 0,1,
\dotsc,(n \! - \! 1)K \! + \! k,
\end{align*}
\begin{gather*}
\mathfrak{b}^{\natural}_{p} \! = \! \mathlarger{\sum_{\underset{\underset{
\sum_{k^{\prime}=1}^{l_{\mathfrak{s}}}k^{\prime}m_{k^{\prime}}=p}{i=
1,2,\dotsc,l_{\mathfrak{s}}}}{m_{i} \in \mathbb{N}_{0}}}^{}} 
\dfrac{(\mathfrak{a}^{\natural}_{1})^{m_{1}}(\mathfrak{a}^{\natural}_{2})^{
m_{2}} \, \dotsb \, (\mathfrak{a}^{\natural}_{l_{\mathfrak{s}}})^{
m_{l_{\mathfrak{s}}}}}{m_{1}!m_{2}! \, \dotsb \, m_{l_{\mathfrak{s}}}!}, 
\quad p \! = \! 1,2,\dotsc,l_{\mathfrak{s}}, \\
\mathfrak{a}_{m}^{\natural} \! := \! \dfrac{1}{m} \sum_{j=1}^{\mathfrak{s}
-1}l_{j}(\alpha_{p_{j}})^{m}, \quad m \! = \! 1,2,\dotsc,l_{\mathfrak{s}},
\end{gather*}
with $\md \widetilde{\mu}(z) \! = \! \exp (-n \widetilde{V}(z)) \, \md z$, and
\begin{equation*}
\widetilde{\nu}^{\raise-1.0ex\hbox{$\scriptstyle \infty$}}_{q,m(q)}(n,k) \! = 
\! \dfrac{\nu^{\infty}_{m(q),q}(n,k)}{\mu^{\infty}_{n,\varkappa_{nk}}(n,k)} 
\! = \! \dfrac{1}{\Upsilon_{D_{\infty}}(n,k)} 
\mathlarger{\sum_{\underset{k^{\prime}-r-1=-m(q)}{k^{\prime},r=1,2,
\dotsc,l_{q}}}^{}} \dfrac{\mathfrak{f}^{(k^{\prime}-1)}(\alpha_{p_{q}}) 
\mathfrak{p}^{(q)}_{r}(n,k)}{(k^{\prime} \! - \! 1)!}, \quad q \! = \! 
1,2,\dotsc,\mathfrak{s} \! - \! 1, \quad m(q) \! = \! 1,2,\dotsc,l_{q},
\end{equation*}
where $\mathfrak{f}^{(k^{\prime}-1)}(\alpha_{p_{q}}) \! := \! \left. 
(\tfrac{\md}{\md z})^{k^{\prime}-1} \mathfrak{f}(z) \right\vert_{z=
\alpha_{p_{q}}}$, with $\mathfrak{f}(z) \! = \! \sum_{r=0}^{(n-1)K+k} 
\mathfrak{b}_{r}z^{r}$, and
\begin{gather*}
\mathfrak{p}^{(q)}_{l_{q}-r(q)}(n,k) \! := \! \mathfrak{b}^{\sharp}_{r(q)}
(q) \prod_{\substack{q^{\prime}=1\\q^{\prime} \neq q}}^{\mathfrak{s}
-1}(\alpha_{p_{q}} \! - \! \alpha_{p_{q^{\prime}}})^{-l_{q^{\prime}}}, 
\quad q \! = \! 1,2,\dotsc,\mathfrak{s} \! - \! 1, \quad r(q) \! = \! 
0,1,\dotsc,l_{q} \! - \! 1, \\
\mathfrak{b}^{\sharp}_{r(q)}(q) \! = \! \mathlarger{\sum_{\underset{
\underset{\sum_{k^{\prime}=1}^{l_{q}-1}k^{\prime}m_{k^{\prime}}=r(q)}{i=1,
2,\dotsc,l_{q}-1}}{m_{i} \in \mathbb{N}_{0}}}} \dfrac{(\mathfrak{a}_{
1}^{\sharp}(q))^{m_{1}}(\mathfrak{a}_{2}^{\sharp}(q))^{m_{2}} \, \dotsb \, 
(\mathfrak{a}_{l_{q}-1}^{\sharp}(q))^{m_{l_{q}-1}}}{m_{1}!m_{2}! \, \dotsb 
\, m_{l_{q}-1}!}, \\
\mathfrak{a}^{\sharp}_{m^{\prime}(q)}(q) \! := \! \dfrac{1}{m^{\prime}(q)} 
\sum_{\substack{j=1\\j \neq q}}^{\mathfrak{s}-1} \dfrac{l_{j}}{(\alpha_{p_{j}} 
\! - \! \alpha_{p_{q}})^{m^{\prime}(q)}}, \quad m^{\prime}(q) \! = \! 1,2,
\dotsc,l_{q} \! - \! 1;
\end{gather*}
and {\rm (ii)} for $n \! \in \! \mathbb{N}$ and $k \! \in \! \lbrace 1,2,
\dotsc,K \rbrace$ such that $\alpha_{p_{\mathfrak{s}}} \! := \! \alpha_{k} 
\! \neq \! \infty$,
\begin{gather*}
\widetilde{\phi}^{\raise-1.0ex\hbox{$\scriptstyle f$}}_{0}(n,k) \! = \! 
\dfrac{\phi^{f}_{0}(n,k)}{\mu^{f}_{n,\varkappa_{nk}}(n,k)} \! = \! 
\dfrac{(-1)^{(n-1)K+k} \prod_{q^{\prime}=1}^{\mathfrak{s}-2}(\alpha_{k} 
\! - \! \alpha_{p_{q^{\prime}}})^{l_{q^{\prime}}}}{\Upsilon_{D_{f}}(n,k)} 
\left(\mathfrak{b}^{o}_{(n-1)K+k-l_{\mathfrak{s}-1}} \! + \! \sum_{m=1}^{
l_{\mathfrak{s}-1}} \mathfrak{b}^{\flat}_{m} \mathfrak{b}^{o}_{m+(n-1)K+k
-l_{\mathfrak{s}-1}} \right), \\
\widetilde{\nu}_{\mathfrak{s}-1,
l_{\mathfrak{s}-1}}^{\raise-1.0ex\hbox{$\scriptstyle f$}}(n,k) \! = \! 
\dfrac{\hat{\nu}^{f}_{n,l_{\mathfrak{s}-1}}(n,k)}{\mu^{f}_{n,\varkappa_{nk}}
(n,k)} \! = \! \dfrac{(-1)^{(n-1)K+k} \mathfrak{b}^{o}_{(n-1)K+k} 
\prod_{q^{\prime}=1}^{\mathfrak{s}-2}(\alpha_{k} \! - \! 
\alpha_{p_{q^{\prime}}})^{l_{q^{\prime}}}}{\Upsilon_{D_{f}}(n,k)},
\end{gather*}
\begin{align*}
\widetilde{\nu}_{\mathfrak{s}-1,
q^{\prime \prime}}^{\raise-1.0ex\hbox{$\scriptstyle f$}}(n,k) \! = \! 
\dfrac{\hat{\nu}^{f}_{n,q^{\prime \prime}}(n,k)}{\mu^{f}_{n,\varkappa_{nk}}
(n,k)} =& \, \dfrac{(-1)^{(n-1)K+k} \prod_{q^{\prime}=1}^{\mathfrak{s}-2}
(\alpha_{k} \! - \! \alpha_{p_{q^{\prime}}})^{l_{q^{\prime}}}}{\Upsilon_{D_{f}}
(n,k)} \left(\vphantom{M^{M^{M^{M^{M^{M^{M}}}}}}} \mathfrak{b}^{o}_{
q^{\prime \prime}+(n-1)K+k-l_{\mathfrak{s}-1}} \right. \\
+&\left. \, \mathlarger{\sum_{\underset{m_{1}-m_{2}=q^{\prime \prime}}{m_{1},
m_{2}=1,2,\dotsc,l_{\mathfrak{s}-1}}}^{}} \mathfrak{b}^{\flat}_{m_{2}} 
\mathfrak{b}^{o}_{m_{1}+(n-1)K+k-l_{\mathfrak{s}-1}} \right), \quad 
q^{\prime \prime} \! = \! 1,2,\dotsc,l_{\mathfrak{s}-1} \! - \! 1,
\end{align*}
where $\Upsilon_{D_{f}}(n,k)$ is given in item~{\rm (ii)} of Lemma~\ref{lem2.2},
\begin{align*}
\mathfrak{b}^{o}_{r} =& \, \underbrace{\int_{\mathbb{R}} \int_{\mathbb{R}} 
\dotsb \int_{\mathbb{R}}}_{(n-1)K+k} \md \widetilde{\mu}(\xi_{0}) \, \md 
\widetilde{\mu}(\xi_{1}) \, \dotsb \, \md \widetilde{\mu}(\xi_{(n-1)K+k-1}) 
\prod_{\substack{i,j=0\\j<i}}^{(n-1)K+k-1}(\xi_{i} \! - \! \xi_{j})^{2} \left(
\prod_{m=0}^{(n-1)K+k-1} \prod_{q=1}^{\mathfrak{s}-2}(\xi_{m} \! - \! 
\alpha_{p_{q}})^{l_{q}}(\xi_{m} \! - \! \alpha_{k})^{l_{\mathfrak{s}}} \right)^{-2} \\
\times& \, \left(\prod_{l=0}^{(n-1)K+k-1}(\xi_{l} \! - \! \alpha_{k}) \right) 
\mathlarger{\sum_{\underset{\underset{\sum_{m^{\prime}=0}^{(n-1)K
+k-1}i_{m^{\prime}}=(n-1)K+k-r}{p=0,1,\dotsc,(n-1)K+k-1}}{i_{p} \in \lbrace 
0,1 \rbrace}}^{}} \prod_{q^{\prime}=0}^{(n-1)K+k-1} \binom{1}{i_{q^{\prime}}}
(-1)^{i_{q^{\prime}}} \xi_{q^{\prime}}^{i_{q^{\prime}}}, \quad r \! = \! 
0,1,\dotsc,(n \! - \! 1)K \! + \! k,
\end{align*}
\begin{gather*}
\mathfrak{b}^{\flat}_{p} \! = \! \mathlarger{\sum_{\underset{\underset{
\sum_{k^{\prime}=1}^{l_{\mathfrak{s}-1}}k^{\prime}m_{k^{\prime}}=p}{i
=1,2,\dotsc,l_{\mathfrak{s}-1}}}{m_{i} \in \mathbb{N}_{0}}}^{}} 
\dfrac{(\hat{\mathfrak{a}}^{\flat}_{1})^{m_{1}}(\hat{\mathfrak{a}}^{
\flat}_{2})^{m_{2}} \, \dotsb \, (\hat{\mathfrak{a}}^{\flat}_{
l_{\mathfrak{s}-1}})^{m_{l_{\mathfrak{s}-1}}}}{m_{1}!m_{2}! \, \dotsb \, 
m_{l_{\mathfrak{s}-1}}!}, \quad p \! = \! 1,2,\dotsc,l_{\mathfrak{s}-1}, \\
\hat{\mathfrak{a}}_{m}^{\flat} \! := \! \dfrac{1}{m} 
\sum_{\substack{j=1\\j \neq \mathfrak{s}-1}}^{\mathfrak{s}}l_{j}
(\alpha_{p_{j}})^{m}, \quad m \! = \! 1,2,\dotsc,l_{\mathfrak{s}-1},
\end{gather*}
and
\begin{align*}
\widetilde{\nu}_{q,m(q)}^{\raise-1.0ex\hbox{$\scriptstyle f$}}(n,k) \! = \! 
\dfrac{\tilde{\nu}^{f}_{m(q),q}(n,k)}{\mu^{f}_{n,\varkappa_{nk}}(n,k)} =& \, 
\dfrac{(-1)^{(n-1)K+k} \prod_{q^{\prime}=1}^{\mathfrak{s}-2}(\alpha_{k} 
\! - \! \alpha_{p_{q^{\prime}}})^{l_{q^{\prime}}}}{\Upsilon_{D_{f}}(n,k)} \\
\times& \, \mathlarger{\sum_{\underset{k^{\prime}-r-1=-m(q)}{k^{\prime}, 
r=1,2,\dotsc,l_{q}}}^{}} \dfrac{\mathfrak{f}_{o}^{(k^{\prime}-1)}(\alpha_{p_{q}}) 
\hat{\mathfrak{p}}^{(q)}_{r}(n,k)}{(k^{\prime} \! - \! 1)!}, \quad q \! = \! 1,2,
\dotsc,\mathfrak{s} \! - \! 2, \, m(q) \! = \! 1,2,\dotsc,l_{q}, \\
\widetilde{\nu}_{\mathfrak{s},
m(\mathfrak{s})}^{\raise-1.0ex\hbox{$\scriptstyle f$}}(n,k) \! = \! 
\dfrac{\mu^{f}_{n,m(\mathfrak{s})}(n,k)}{\mu^{f}_{n,\varkappa_{nk}}(n,k)} =& 
\, \dfrac{(-1)^{(n-1)K+k} \prod_{q^{\prime}=1}^{\mathfrak{s}-2}(\alpha_{k} 
\! - \! \alpha_{p_{q^{\prime}}})^{l_{q^{\prime}}}}{\Upsilon_{D_{f}}(n,k)} \\
\times& \, \mathlarger{\sum_{\underset{k^{\prime}-r-1=
-m(\mathfrak{s})}{k^{\prime},r=1,2,\dotsc,l_{\mathfrak{s}}}}^{}} 
\dfrac{\mathfrak{f}_{o}^{(k^{\prime}-1)}(\alpha_{k}) 
\hat{\mathfrak{p}}^{(\mathfrak{s})}_{r}(n,k)}{(k^{\prime} \! - \! 1)!}, 
\quad m(\mathfrak{s}) \! = \! 1,2,\dotsc,l_{\mathfrak{s}} \! - \! 1,
\end{align*}
where $\mathfrak{f}_{o}^{(k^{\prime}-1)}(\alpha_{p_{q}}) \! := \! \left. 
(\tfrac{\md}{\md z})^{k^{\prime}-1} \mathfrak{f}_{o}(z) \right\vert_{z=
\alpha_{p_{q}}}$, with $\mathfrak{f}_{o}(z) \! = \! \sum_{r=0}^{(n-1)K+k} 
\mathfrak{b}^{o}_{r}z^{r}$, and
\begin{gather*}
\hat{\mathfrak{p}}^{(q)}_{l_{q}-r(q)}(n,k) \! := \! \hat{\mathfrak{b}}^{o}_{r(q)}
(q) \prod_{\substack{q^{\prime}=1\\q^{\prime} \neq q,\mathfrak{s}-
1}}^{\mathfrak{s}}(\alpha_{p_{q}} \! - \! \alpha_{p_{q^{\prime}}})^{-
l_{q^{\prime}}}, \quad q \! = \! 1,\dotsc,\mathfrak{s} \! - \! 2,
\mathfrak{s}, \quad r(q) \! = \! 0,1,\dotsc,l_{q} \! - \! 1, \\
\hat{\mathfrak{b}}^{o}_{r(q)}(q) \! = \! \mathlarger{\sum_{\underset{
\underset{\sum_{k^{\prime}=1}^{l_{q}-1}k^{\prime}m_{k^{\prime}}
=r(q)}{i=1,2,\dotsc,l_{q}-1}}{m_{i} \in \mathbb{N}_{0}}}} 
\dfrac{(\hat{\mathfrak{a}}_{1}^{o}(q))^{m_{1}}(\hat{\mathfrak{a}}_{2}^{o}
(q))^{m_{2}} \, \dotsb \, (\hat{\mathfrak{a}}_{l_{q}-1}^{o}(q))^{m_{l_{q}
-1}}}{m_{1}!m_{2}! \, \dotsb \, m_{l_{q}-1}!}, \\
\hat{\mathfrak{a}}^{o}_{m^{\prime}(q)}(q) \! := \! \dfrac{1}{m^{\prime}(q)} 
\sum_{\substack{j=1\\j \neq q,\mathfrak{s}-1}}^{\mathfrak{s}} 
\dfrac{l_{j}}{(\alpha_{p_{j}} \! - \! \alpha_{p_{q}})^{m^{\prime}(q)}}, 
\quad m^{\prime}(q) \! = \! 1,2,\dotsc,l_{q} \! - \! 1.
\end{gather*}}
\end{eeeee}
\section{The Monic MPC ORF Families of Variational Problems: Existence, 
Uniqueness, Regularity, and the Transformed RHPs} \label{sec3} 
In this section, the analysis of the family of $K$ $(= \! \hat{K} \! + 
\! \tilde{K}$, where $\hat{K} \! = \! \# \lbrace \mathstrut k \! \in \! 
\lbrace 1,2,\dotsc,K \rbrace; \, \alpha_{k} \! = \! \infty \rbrace$, and 
$\tilde{K} \! = \! \# \lbrace \mathstrut k \! \in \! \lbrace 1,2,\dotsc,K 
\rbrace; \, \alpha_{k} \! \neq \! \infty \rbrace)$ variational problems 
corresponding to the associated monic MPC ORFs is undertaken. This family 
of energy minimisation problems requires an existence and regularity theory 
due to the presence of the poles $\alpha_{1},\alpha_{2},\dotsc,\alpha_{K}$ 
and the singular behaviour of $\widetilde{V}$ at each $\alpha_{k}$, $k \! \in 
\! \lbrace 1,2,\dotsc,K \rbrace$ (cf. conditions~\eqref{eq20}--\eqref{eq22}). 
More precisely, the two corresponding $n$-, $k$-, and $z_{o}$-dependent 
subfamilies of associated equilibrium measures, that is, their existence, 
uniqueness, and regularity theory, and generalised weighted Fekete sets are 
analysed; furthermore, by considering the two corresponding $n$-, $k$-, 
and $z_{o}$-dependent subfamilies of complex potentials (`$g$-functions'), 
Lemma~$\bm{\mathrm{RHP}_{\mathrm{MPC}}}$ is re-formulated as a family 
of $K$ equivalent {}\footnote{If there are two matrix RHPs, $(\mathcal{X}_{1}
(z),\upsilon_{1}(z),\varGamma_{1})$ and $(\mathcal{X}_{2}(z),\upsilon_{2}
(z),\varGamma_{2})$, say, where $\varGamma_{2} \subset \varGamma_{1}$ 
and $\upsilon_{1} \! \! \upharpoonright_{\varGamma_{1} \setminus 
\varGamma_{2}} \! =_{\underset{z_{o}=1+o(1)}{\mathscr{N},n \to \infty}} 
\! \mathrm{I} \! + \! o(1)$, then their solutions, $\mathcal{X}_{1}$ and 
$\mathcal{X}_{2}$, respectively, are asymptotically (modulo $o(1)$ 
estimates) equal.} matrix RHPs on $\overline{\mathbb{R}}$.
\begin{eeeee} \label{remm3.1} 
\textsl{The formulae, etc., of this Section~\ref{sec3}, while valid for finite 
$n$ $(\in \! \mathbb{N})$, are germane, principally, for the ensuing 
asymptotic analysis, in the double-scaling limit $\mathscr{N},n \! \to \! 
\infty$ such that $z_{o} \! = \! 1 \! + \! o(1)$.}
\end{eeeee} 
\begin{eeeee} \label{remm3.2} 
\textsl{The equilibrium measures, weighted logarithmic energy 
functionals, Fekete sets, etc., introduced in this monograph vary 
according to the parameters $n$ $(\in \! \mathbb{N})$, 
$k$ $(\in \! \lbrace 1,2,\dotsc,K \rbrace)$, and $z_{o}$ (cf. 
Remark~\ref{rem1.3.2}$)$, and thus, \emph{sensus strictu}, all 
associated variables, too, depend on $n$, $k$, and $z_{o}$, which 
would necessitate the introduction of an over-arching amount of 
additional, superfluous notation(s) to subsume these $n$-, $k$-, 
and $z_{o}$-dependencies; however, for simplicity of notation, 
such cumbersome $n$-, $k$-, and $z_{o}$-dependencies will be 
suppressed, unless where absolutely necessary, and the reader 
should be cognizant of this fact. This comment applies, 
\emph{mutatis mutandis}, throughout the remainder of this 
monograph.}
\end{eeeee}
\begin{eeeee} \label{rhemm3.3} 
\textsl{The `symbol(s)' (`notation(s)') $\mathfrak{c}_{1},
\mathfrak{c}_{2},\mathfrak{c}_{3},\dotsc$, with or without subscripts, 
superscripts, underscripts, overscripts, explicit reference(s) to 
their functional dependencies, etc., appearing in the various error 
estimates throughout this monograph are not equal, and should 
properly be denoted as (suppressing, again, all subscripts, 
superscripts, underscripts, and overscripts) $\mathfrak{c}_{1}(n,k,
z_{o}),\mathfrak{c}_{2}(n,k,z_{o}),\mathfrak{c}_{3}(n,k,z_{o}),\dotsc$. 
Since the principal concern of this monograph is not their explicit 
(functional) $n$-, $k$-, and $z_{o}$-dependencies, which, in every 
instance, is $\mathcal{O}(1)$ (in the double-scaling limit $\mathscr{N},
n \! \to \! \infty$ such that $z_{o} \! = \! 1 \! + \! o(1))$, but, rather, 
the explicit class(es) to which they belong, the simplified `notation(s)' 
$\mathfrak{c}_{1},\mathfrak{c}_{2},\mathfrak{c}_{3},\dotsc$ are 
retained throughout this monograph, unless stated otherwise.}
\end{eeeee}

One begins by establishing the existence and the uniqueness of the 
family of $K$ equilibrium measures.
\begin{ccccc} \label{lem3.1} 
Let the external field $\widetilde{V} \colon \overline{\mathbb{R}} 
\setminus \lbrace \alpha_{1},\alpha_{2},\dotsc,\alpha_{K} \rbrace \! 
\to \! \mathbb{R}$ satisfy conditions~\eqref{eq20}--\eqref{eq22}, 
and set $\varpi (z) \! := \! \exp (-\widetilde{V}(z)/2)$. For $n \! \in \! 
\mathbb{N}$ and $k \! \in \! \lbrace 1,2,\dotsc,K \rbrace$ such that 
$\alpha_{p_{\mathfrak{s}}} \! := \! \alpha_{k} \! = \! \infty$, with 
associated measure $\mu^{\text{\tiny $\mathrm{EQ}$}}_{1} \! \in 
\! \mathscr{M}_{1}(\mathbb{R})$, define the corresponding weighted 
logarithmic energy functional
\begin{align}
\mathrm{I}_{\widetilde{V}}^{\infty} &\colon \mathbb{N} \times \lbrace 
1,2,\dotsc,K \rbrace \times \mathscr{M}_{1}(\mathbb{R}) \! \ni \! 
(n,k,\mu^{\text{\tiny $\mathrm{EQ}$}}_{1}) \! \mapsto \! 
\mathrm{I}_{\widetilde{V}}^{\infty}[n,k;
\mu^{\text{\tiny $\mathrm{EQ}$}}_{1}] 
\! := \! \mathrm{I}_{\widetilde{V}}^{\infty}
[\mu^{\text{\tiny $\mathrm{EQ}$}}_{1}] \nonumber \\
&= \iint_{\mathbb{R}^{2}} \ln \left(\vert \xi \! - \! \tau 
\vert^{\frac{\varkappa_{nk}}{n}} \prod_{q=1}^{\mathfrak{s}-1} \left(
\dfrac{\vert \xi \! - \! \tau \vert}{\vert \xi \! - \! \alpha_{p_{q}} \vert 
\vert \tau \! - \! \alpha_{p_{q}} \vert} \right)^{\frac{\varkappa_{nk 
\tilde{k}_{q}}}{n}} \varpi (\xi) \varpi (\tau) \right)^{-1} \md 
\mu^{\text{\tiny $\mathrm{EQ}$}}_{1}(\xi) \, \md 
\mu^{\text{\tiny $\mathrm{EQ}$}}_{1}(\tau), \label{eqlm3.1a}
\end{align}
and consider its associated minimisation problem
\begin{equation*}
E_{\widetilde{V}}^{\infty} \! = \! \inf \left\lbrace \mathstrut \mathrm{I}_{
\widetilde{V}}^{\infty}[\mu^{\text{\tiny $\mathrm{EQ}$}}_{1}]; \, 
\mu^{\text{\tiny $\mathrm{EQ}$}}_{1} \! \in \! \mathscr{M}_{1}
(\mathbb{R}) \right\rbrace;
\end{equation*}
and, for $n \! \in \! \mathbb{N}$ and $k \! \in \! \lbrace 1,2,\dotsc,
K \rbrace$ such that $\alpha_{p_{\mathfrak{s}}} \! := \! \alpha_{k} 
\! \neq \! \infty$, with associated measure 
$\mu^{\text{\tiny $\mathrm{EQ}$}}_{2} \! \in \! \mathscr{M}_{1}
(\mathbb{R})$, define the corresponding weighted logarithmic energy 
functional
\begin{align}
\mathrm{I}_{\widetilde{V}}^{f} &\colon \mathbb{N} \times \lbrace 
1,2,\dotsc,K \rbrace \times \mathscr{M}_{1}(\mathbb{R}) \! \ni \! 
(n,k,\mu^{\text{\tiny $\mathrm{EQ}$}}_{2}) \! \mapsto \! 
\mathrm{I}_{\widetilde{V}}^{f}[n,k;\mu^{\text{\tiny $\mathrm{EQ}$}}_{2}] \! 
:= \! \mathrm{I}_{\widetilde{V}}^{f}[\mu^{\text{\tiny $\mathrm{EQ}$}}_{2}] 
\nonumber \\
&= \iint_{\mathbb{R}^{2}} \ln \left(\vert \xi \! - \! \tau 
\vert^{\frac{\varkappa_{nk \tilde{k}_{\mathfrak{s}-1}}^{\infty}+1}{n}} 
\left(\dfrac{\vert \xi \! - \! \tau \vert}{\vert \xi \! - \! \alpha_{k} \vert 
\vert \tau \! - \! \alpha_{k} \vert} \right)^{\frac{\varkappa_{nk}-1}{n}} 
\prod_{q=1}^{\mathfrak{s}-2} \left(\dfrac{\vert \xi \! - \! \tau \vert}{
\vert \xi \! - \! \alpha_{p_{q}} \vert \vert \tau \! - \! \alpha_{p_{q}} 
\vert} \right)^{\frac{\varkappa_{nk \tilde{k}_{q}}}{n}} \varpi (\xi) \varpi 
(\tau) \right)^{-1} \md \mu^{\text{\tiny $\mathrm{EQ}$}}_{2}(\xi) \, 
\md \mu^{\text{\tiny $\mathrm{EQ}$}}_{2}(\tau), \label{eqlm3.1b}
\end{align}
and consider its associated minimisation problem
\begin{equation*}
E_{\widetilde{V}}^{f} \! = \! \inf \left\lbrace \mathstrut \mathrm{I}_{
\widetilde{V}}^{f}[\mu^{\text{\tiny $\mathrm{EQ}$}}_{2}]; \, 
\mu^{\text{\tiny $\mathrm{EQ}$}}_{2} \! \in \! 
\mathscr{M}_{1}(\mathbb{R}) \right\rbrace.
\end{equation*}
Then, for $n \! \in \! \mathbb{N}$ and $k \! \in \! \lbrace 1,2,\dotsc,
K \rbrace$ such that $\alpha_{p_{\mathfrak{s}}} \! := \! \alpha_{k} 
\! = \! \infty$ (resp., $\alpha_{p_{\mathfrak{s}}} \! := \! \alpha_{k} 
\! \neq \! \infty)$$:$ {\rm (1)} $E_{\widetilde{V}}^{\infty}$ (resp., 
$E_{\widetilde{V}}^{f})$ is finite; {\rm (2)} the infimum is attained, that is, 
there exists $\mu_{\widetilde{V}}^{\infty} \! \in \! \mathscr{M}_{1}
(\mathbb{R})$ (resp., $\mu_{\widetilde{V}}^{f} \! \in \! \mathscr{M}_{1}
(\mathbb{R}))$, the associated equilibrium measure, such that 
$\mathrm{I}_{\widetilde{V}}^{\infty}[\mu_{\widetilde{V}}^{\infty}] \! 
= \! E_{\widetilde{V}}^{\infty}$ (resp., $\mathrm{I}_{\widetilde{V}}^{f}
[\mu_{\widetilde{V}}^{f}] \! = \! E_{\widetilde{V}}^{f})$, and 
$\mu_{\widetilde{V}}^{\infty}$ (resp., $\mu_{\widetilde{V}}^{f})$ 
has finite weighted logarithmic energy, that is, $-\infty \! < \! 
\mathrm{I}_{\widetilde{V}}^{\infty}[\mu_{\widetilde{V}}^{\infty}] \! 
< \! +\infty$ (resp., $-\infty \! < \! \mathrm{I}_{\widetilde{V}}^{f}
[\mu_{\widetilde{V}}^{f}] \! < \! +\infty)$$;$ and {\rm (3)} $J_{\infty} 
\! := \! \supp (\mu_{\widetilde{V}}^{\infty})$ (resp., $J_{f} \! := \! 
\supp (\mu_{\widetilde{V}}^{f}))$ is a proper compact subset of 
$\overline{\mathbb{R}} \setminus \lbrace \alpha_{1},\alpha_{2},
\dotsc,\alpha_{K} \rbrace$.
\end{ccccc}
\begin{notetolem} 
\textsl{The proof of Lemma~\ref{lem3.1} that is presented in this lecture 
note is modelled on the paradigm presented in Chapter~{\rm 6} of 
Deift's book {\rm \cite{a51}} (see, also, Section~{\rm 9} of {\rm \cite{ears}}, 
and {\rm \cite{abjkglfsa}):} there is, however, an alternate proof, which is 
based on the ideas used to prove Theorem~{\rm 1.2} in {\rm \cite{a57}} 
(see, in particular, p.~{\rm 740}, Section~{\rm 1}, Theorem~{\rm 1.2}, 
and pp.~{\rm 746}--{\rm 748}, Section~{\rm 3}, 
\textbf{Proof of Theorem~{\rm 1.2}} of {\rm \cite{a57})}. 
Adapting and suitably modifying the core ideas subsumed in the proof 
of Theorem~{\rm 1.2} of Kuijlaars-McLaughlin {\rm \cite{a57}} to the 
situation considered in this monograph, though, one must suppose 
that, for real analytic external fields $\widetilde{V}(z)$ and $\lbrace 
\widetilde{V}_{m}(z) \rbrace_{m \in \mathbb{N}}$ on $\overline{\mathbb{R}} 
\setminus \lbrace \alpha_{1},\alpha_{2},\dotsc,\alpha_{K} \rbrace$, the following 
hold: {\rm (i)} $\lim_{m^{\prime} \to \infty} \widetilde{V}_{m^{\prime}}^{(i)}
(x) \! = \! \widetilde{V}^{(i)}(x)$, $i \! = \! 0,1,2,3$, uniformly on compact 
subsets of $\overline{\mathbb{R}} \setminus \lbrace \alpha_{1},\alpha_{2},
\dotsc,\alpha_{K} \rbrace$$;$ and {\rm (ii)} the `growth conditions' 
$\lim_{x \to \alpha_{i}}(\widetilde{V}_{m}(x)/\ln \lvert x \rvert) \! = \! 
+\infty$, $i \! \in \! \lbrace \mathstrut k \! \in \! \lbrace 1,2,\dotsc,K 
\rbrace; \, \alpha_{k} \! = \! \infty \rbrace$, and $\lim_{x \to \alpha_{j}}
(\widetilde{V}_{m}(x)/\ln \lvert x \! - \! \alpha_{j} \rvert^{-1}) \! = \! 
+\infty$, $j \! \in \! \lbrace \mathstrut k \! \in \! \lbrace 1,2,\dotsc,K 
\rbrace; \, \alpha_{k} \! \neq \! \infty \rbrace$, hold uniformly in 
$m$. Furthermore, one must have \emph{a priori} knowledge 
of the analogue of their function $q_{V}(x)$ (see 
Equations~{\rm (1.12)} and~{\rm (3.8)} of {\rm \cite{a57};} see, also, 
{\rm \cite{a58})}, which, in this monograph, will not be available until 
(see Lemma~\ref{lem3.7} below) Equation~\eqref{eqveeinf}, for the case 
$\alpha_{p_{\mathfrak{s}}} \! := \! \alpha_{k} \! = \! \infty$, $k \! \in 
\! \lbrace 1,2,\dotsc,K \rbrace$, and Equation~\eqref{eqveefin}, for the 
case $\alpha_{p_{\mathfrak{s}}} \! := \! \alpha_{k} \! \neq \! \infty$, 
$k \! \in \! \lbrace 1,2,\dotsc,K \rbrace$, and the fact that the density 
of their equilibrium measure, denoted $\md \mu_{V}$, can be presented 
in the form $\md \mu_{V}(x) \! = \! \psi_{V}(x) \, \md x$ on $\supp 
(\mu_{V})$ (for $\psi_{V}$ continuous with compact support), which, in 
this monograph, is established only in Lemma~\ref{lem3.7} below (see 
Equation~\eqref{eql3.7e} for the case $\alpha_{p_{\mathfrak{s}}} \! := 
\! \alpha_{k} \! = \! \infty$, $k \! \in \! \lbrace 1,2,\dotsc,K \rbrace$, 
and Equation~\eqref{eql3.7k} for the case $\alpha_{p_{\mathfrak{s}}} \! 
:= \! \alpha_{k} \! \neq \! \infty$, $k \! \in \! \lbrace 1,2,\dotsc,K 
\rbrace)$. Since the above-mentioned suppositions are not made, and 
the above-mentioned facts are not established until well into the analysis 
of this Section~\ref{sec3}, it is more natural to present, at this stage, the 
proof of Lemma~\ref{lem3.1} based on the ideas of Chapter~{\rm 6} of 
Deift's book {\rm \cite{a51}} (see, also, {\rm \cite{a56})}. If, however, 
the above suppositions (facts) had been established, then, one could 
proceed, with suitable modifications and adaptations, as in the proof of 
Theorem~{\rm 1.2} of {\rm \cite{a57}}, to show that there exists some 
$\tilde{A} \! > \! 0$ (and bounded) such that, $\forall$ $m \! \in \! 
\mathbb{N}$, $\supp (\mu_{\widetilde{V}_{m}}) \subset [-\tilde{A},\tilde{A}] 
\setminus \lbrace \mathstrut \alpha_{k^{\prime}}, \, k^{\prime} \! \in \! 
\lbrace 1,2,\dotsc,K \rbrace; \, \alpha_{k^{\prime}} \! \neq \! \infty \rbrace$ 
and is compact; moreover, one could also show that, for regular $\widetilde{V}$, 
$\exists$ $m_{0} \! \in \! \mathbb{N}$ such that $\widetilde{V}_{m}$ is 
regular $\forall$ $m \! \geqslant \! m_{0}$. Incidentally, and as a 
by-product of the latter results, one can also establish that regular external 
fields are generic, that is, open and dense in the space $\widetilde{
\mathbb{V}}$ of real analytic functions $\widetilde{V}$ satisfying the 
conditions~\eqref{eq20}--\eqref{eq22}$;$ in particular, the issue of 
openness is addressed by providing the space $\widetilde{\mathbb{V}}$ with 
the---in the terminology of {\rm \cite{a57}}---``right topology'' (the question 
of denseness is not addressed in this monograph). Towards this end, one wants 
that a sequence $\lbrace \widetilde{V}_{m}(z) \rbrace_{m \in \mathbb{N}}$ 
converges to $\widetilde{V}(z)$ in $\widetilde{\mathbb{V}}$ if, and only if, 
the conditions~{\rm (i)} and~{\rm (ii)} above are satisfied. Define $\hat{M} \! 
:= \! \max \lbrace \mathstrut \lvert \alpha_{k^{\prime}} \rvert, \, k^{\prime} 
\! \in \! \lbrace 1,2,\dotsc,K \rbrace; \, \alpha_{k^{\prime}} \! \neq \! \infty 
\rbrace \! + \! \max \lbrace \mathstrut \lvert \alpha_{i} \! - \! \alpha_{j} 
\rvert, \, i \! \neq \! j \! \in \! \lbrace 1,2,\dotsc,K \rbrace; \, 
\alpha_{i},\alpha_{j} \! \neq \! \infty \rbrace$, and choose $\delta_{\hat{M}} 
\! > \! 0$ so that $3 \delta_{\hat{M}} \! < \! \hat{m}$, where $\hat{m} \! 
:= \! \min \lbrace \mathstrut \lvert \alpha_{i} \! - \! \alpha_{j} \rvert, \, 
i \! \neq \! j \! \in \! \lbrace 1,2,\dotsc,K \rbrace; \, \alpha_{i},\alpha_{j} 
\! \neq \! \infty \rbrace$. For $\alpha_{k^{\prime}} \! \neq \! \infty$, 
$k^{\prime} \! \in \! \lbrace 1,2,\dotsc,K \rbrace$, define 
$\hat{\mathscr{O}}_{\frac{\delta_{\hat{M}}}{K(j+1)}}(\alpha_{k^{\prime}}) 
\! := \! \lbrace \mathstrut x \! \in \! \mathbb{R}; \, \lvert x \! - \! 
\alpha_{k^{\prime}} \rvert \! \leqslant \! \delta_{\hat{M}}/K(j \! + \! 1) \rbrace$, 
$j \! \in \! \mathbb{N}$. For $j \! \in \! \mathbb{N}$ and $\alpha_{k^{\prime}} 
\! \neq \! \infty$, $k^{\prime} \! \in \! \lbrace 1,2,\dotsc,K \rbrace$, define 
the real-valued mappings (not the only ones possible!) $\hat{\mathcal{F}}_{j} 
\colon \widetilde{\mathbb{V}} \! \ni \! \widetilde{V} \! \mapsto \! 
\inf_{\lvert x \rvert \geqslant K \hat{M}+j}(\widetilde{V}(x)/\ln \lvert x 
\rvert)$ and $\hat{\mathcal{G}}_{k^{\prime},j} \colon \widetilde{\mathbb{V}} 
\! \ni \! \widetilde{V} \! \mapsto \! \inf_{\lvert x- \alpha_{k^{\prime}} 
\rvert \leqslant \delta_{\hat{M}}/K(j+1)}(\widetilde{V}(x)/\ln \lvert 
x \! - \! \alpha_{k^{\prime}} \rvert^{-1})$, and consider the space 
$\widetilde{\mathbb{V}}$ endowed with the metric (not the only choice 
possible!)
\begin{align*}
\widetilde{\sigma}(\widetilde{V},\widetilde{W}) :=& \, \sum_{j=1}^{\infty} 
\dfrac{1}{2^{j}} \dfrac{\lvert \hat{\mathcal{F}}_{j}(\widetilde{V}) \! - \! 
\hat{\mathcal{F}}_{j}(\widetilde{W}) \rvert}{1 \! + \! \lvert 
\hat{\mathcal{F}}_{j}(\widetilde{V}) \! - \! \hat{\mathcal{F}}_{j}
(\widetilde{W}) \rvert} \! + \! \sum_{\lbrace \mathstrut k^{\prime} \in 
\lbrace 1,2,\dotsc,K \rbrace, \, \alpha_{k^{\prime}} \neq \infty \rbrace} 
\sum_{j=1}^{\infty} \dfrac{1}{2^{j}} \dfrac{\lvert \hat{\mathcal{G}}_{
k^{\prime},j}(\widetilde{V}) \! - \! \hat{\mathcal{G}}_{k^{\prime},j}
(\widetilde{W}) \rvert}{1 \! + \! \lvert \hat{\mathcal{G}}_{k^{\prime},j}
(\widetilde{V}) \! - \! \hat{\mathcal{G}}_{k^{\prime},j}(\widetilde{W}) 
\rvert} \\
+& \, \sum_{i=0}^{3} \sum_{j=1}^{\infty} \dfrac{1}{2^{j}} \dfrac{\lvert 
\lvert \widetilde{V}^{(i)} \! - \! \widetilde{W}^{(i)} \rvert \rvert_{
\hat{A}(j)}}{1 \! + \! \lvert \lvert \widetilde{V}^{(i)} \! - \! 
\widetilde{W}^{(i)} \rvert \rvert_{\hat{A}(j)}},
\end{align*}
with $\hat{A}(j) \! := \! [-(K \hat{M} \! + \! j),K \hat{M} \! + \! j] \setminus 
\cup_{\lbrace \mathstrut \alpha_{k^{\prime}} \neq \infty, \, k^{\prime} 
\in \lbrace 1,2,\dotsc,K \rbrace \rbrace} \hat{\mathscr{O}}_{\frac{
\delta_{\hat{M}}}{K(j+1)}}(\alpha_{k^{\prime}})$, where the norm $\lvert 
\lvert \pmb{\mathbf{\cdot}} \rvert \rvert_{\hat{A}(j)}$ is the sup norm on 
$\hat{A}(j)$. One then notes that convergence of a sequence $\lbrace 
\widetilde{V}_{m}(z) \rbrace_{m \in \mathbb{N}}$ to $\widetilde{V}(z)$ 
in the metric space $(\widetilde{\mathbb{V}},\widetilde{\sigma})$ 
is equivalent to the convergence of $\widetilde{V}_{m}^{(i)}(z)$ to 
$\widetilde{V}^{(i)}(z)$, $i \! = \! 0,1,2,3$, uniformly on compact subsets 
of $\overline{\mathbb{R}} \setminus \lbrace \alpha_{1},\alpha_{2},\dotsc,
\alpha_{K} \rbrace$ (together with the uniform `growth conditions' at 
each pole $\alpha_{k}$, $k \! \in \! \lbrace 1,2,\dotsc,K \rbrace)$.}
\end{notetolem}

\emph{Proof}. The proof of this Lemma~\ref{lem3.1} consists of two cases: 
(i) $n \! \in \! \mathbb{N}$ and $k \! \in \! \lbrace 1,2,\dotsc,K \rbrace$ 
such that $\alpha_{p_{\mathfrak{s}}} \! := \! \alpha_{k} \! = \! \infty$; 
and (ii) $n \! \in \! \mathbb{N}$ and $k \! \in \! \lbrace 1,2,\dotsc,K 
\rbrace$ such that $\alpha_{p_{\mathfrak{s}}} \! := \! \alpha_{k} \! \neq 
\! \infty$. For the sake of eschewing redundancies, and without loss of 
generality, only case (ii) will be considered in detail (see $\pmb{(1)}$ 
below), whilst case (i) can be proved in an analogous manner (see 
$\pmb{(2)}$ below).

$\pmb{(1)}$ For $n \! \in \! \mathbb{N}$ and $k \! \in \! \lbrace 1,2,
\dotsc,K \rbrace$ such that $\alpha_{p_{\mathfrak{s}}} \! := \! \alpha_{k} 
\! \neq \! \infty$, let $\mu^{\text{\tiny $\mathrm{EQ}$}}_{2} \! \in \! 
\mathscr{M}_{1}(\mathbb{R})$, and set, as in the lemma, $\varpi (z) 
\! := \! \exp (-\widetilde{V}(z)/2)$, where $\widetilde{V} \colon 
\overline{\mathbb{R}} \setminus \lbrace \alpha_{1},\alpha_{2},
\dotsc,\alpha_{K} \rbrace \! \to \! \mathbb{R}$ satisfies 
conditions~\eqref{eq20}--\eqref{eq22}. For $n \! \in \! \mathbb{N}$ 
and $k \! \in \! \lbrace 1,2,\dotsc,K \rbrace$ such that 
$\alpha_{p_{\mathfrak{s}}} \! := \! \alpha_{k} \! \neq \! \infty$, {}from 
the definition of $\mathrm{I}_{\widetilde{V}}^{f} \colon \mathbb{N} 
\times \lbrace 1,2,\dotsc,K \rbrace \times \mathscr{M}_{1}(\mathbb{R}) 
\! \to \! \mathbb{R}$ stated in the lemma,\footnote{The definition of 
$\mathrm{I}_{\widetilde{V}}^{f}[\mu^{\text{\tiny $\mathrm{EQ}$}}_{2}]$ 
has sense provided the integrals exist and are finite.} one shows that
\begin{equation*}
\mathrm{I}_{\widetilde{V}}^{f}[\mu^{\text{\tiny $\mathrm{EQ}$}}_{2}] 
\! = \! \iint_{\mathbb{R}^{2}} \mathcal{K}_{\widetilde{V}}^{f}
(\xi,\tau) \, \md \mu^{\text{\tiny $\mathrm{EQ}$}}_{2}(\xi) \, 
\md \mu^{\text{\tiny $\mathrm{EQ}$}}_{2}(\tau),
\end{equation*}
where the symmetric kernel $\mathcal{K}_{\widetilde{V}}^{f}(\xi,\tau) 
\! = \! \mathcal{K}_{\widetilde{V}}^{f}(\tau,\xi)$ is given by
\begin{equation} \label{eqKvinf1} 
\begin{split}
\mathcal{K}_{\widetilde{V}}^{f}(\xi,\tau) =& \, \ln \left(\vert \xi \! - \! 
\tau \vert^{\frac{\varkappa_{nk \tilde{k}_{\mathfrak{s}-1}}^{\infty}+1}{n}} 
\left(\dfrac{\vert \xi \! - \! \tau \vert}{\vert \xi \! - \! \alpha_{k} \vert 
\vert \tau \! - \! \alpha_{k} \vert} \right)^{\frac{\varkappa_{nk}-1}{n}} 
\prod_{q=1}^{\mathfrak{s}-2} \left(\dfrac{\vert \xi \! - \! \tau \vert}{
\vert \xi \! - \! \alpha_{p_{q}} \vert \vert \tau \! - \! \alpha_{p_{q}} 
\vert} \right)^{\frac{\varkappa_{nk \tilde{k}_{q}}}{n}} \varpi (\xi) \varpi 
(\tau) \right)^{-1} \\
=& \, \left(\dfrac{\varkappa_{nk \tilde{k}_{\mathfrak{s}-1}}^{\infty} \! 
+ \! 1}{n} \right) \ln \left(\dfrac{1}{\lvert \xi \! - \! \tau \rvert} \right) 
\! + \! \left(\dfrac{\varkappa_{nk} \! - \! 1}{n} \right) \ln \left(\left\vert 
\dfrac{1}{\xi \! - \! \alpha_{k}} \! - \! \dfrac{1}{\tau \! - \! \alpha_{k}} 
\right\vert^{-1} \right) \\
+& \, \sum_{q=1}^{\mathfrak{s}-2} \dfrac{\varkappa_{nk \tilde{k}_{q}}}{n} 
\ln \left(\left\vert \dfrac{1}{\xi \! - \! \alpha_{p_{q}}} \! - \! \dfrac{1}{\tau 
\! - \! \alpha_{p_{q}}} \right\vert^{-1} \right) \! + \! \dfrac{1}{2} 
\widetilde{V}(\xi) \! + \! \dfrac{1}{2} \widetilde{V}(\tau) \\
=& \, \dfrac{1}{2} \widetilde{V}(\xi) \! + \! \dfrac{1}{2} \widetilde{V}(\tau) 
\! + \! \ln \left(\dfrac{\lvert \xi \! - \! \tau \rvert^{K}}{\prod_{q=1}^{
\mathfrak{s}-2}(\lvert \xi \! - \! \alpha_{p_{q}} \rvert \lvert \tau \! - \! 
\alpha_{p_{q}} \rvert)^{\gamma_{i(q)_{k_{q}}}}(\lvert \xi \! - \! \alpha_{k} 
\rvert \lvert \tau \! - \! \alpha_{k} \rvert)^{\gamma_{k}}} \right)^{-1} \\
+& \, 
\begin{cases} 
\dfrac{1}{n} \ln \left(\dfrac{\lvert \xi \! - \! \tau \rvert^{k} 
\prod_{j=1}^{\mathfrak{s}-2}(\lvert \xi \! - \! \alpha_{p_{j}} \rvert \lvert 
\tau \! - \! \alpha_{p_{j}} \rvert)^{\gamma_{i(j)_{k_{j}}}}(\lvert \xi \! - 
\! \alpha_{k} \rvert \lvert \tau \! - \! \alpha_{k} \rvert)^{\gamma_{k}}}{
\lvert \xi \! - \! \tau \rvert^{K}(\lvert \xi \! - \! \alpha_{k} \rvert 
\lvert \tau \! - \! \alpha_{k} \rvert)^{\varrho_{k}-1}} \right)^{-1}, 
&\text{$\hat{\mathfrak{J}}_{q}(k) \! = \! \varnothing, \quad q \! 
\in \! \lbrace 1,2,\dotsc,\mathfrak{s} \! - \! 2 \rbrace$,} \\
\dfrac{1}{n} \ln \left(\dfrac{\lvert \xi \! - \! \tau \rvert^{k} 
\prod_{j=1}^{\mathfrak{s}-2}(\lvert \xi \! - \! \alpha_{p_{j}} \rvert \lvert 
\tau \! - \! \alpha_{p_{j}} \rvert)^{\gamma_{i(j)_{k_{j}}}}(\lvert \xi \! - \! 
\alpha_{k} \rvert \lvert \tau \! - \! \alpha_{k} \rvert)^{\gamma_{k}}}{
\lvert \xi \! - \! \tau \rvert^{K} \prod_{j=1}^{\mathfrak{s}-2}(\lvert \xi 
\! - \! \alpha_{p_{j}} \rvert \lvert \tau \! - \! \alpha_{p_{j}} \rvert)^{
\varrho_{\hat{m}_{j}(k)}}(\lvert \xi \! - \! \alpha_{k} \rvert \lvert 
\tau \! - \! \alpha_{k} \rvert)^{\varrho_{k}-1}} \right)^{-1}, 
&\text{$\hat{\mathfrak{J}}_{q}(k) \! \neq \! \varnothing, \quad 
q \! \in \! \lbrace 1,2,\dotsc,\mathfrak{s} \! - \! 2 \rbrace$.}
\end{cases}
\end{split}
\end{equation}
Via the inequalities
\begin{gather*}
\ln \vert \xi \! - \! \tau \vert^{-1} \! \geqslant \! -\dfrac{1}{2} \ln 
(1 \! + \! \xi^{2}) \! - \! \dfrac{1}{2} \ln (1 \! + \! \tau^{2}), \\
\ln \vert (\xi \! - \! \alpha_{k})^{-1} \! - \! (\tau \! - \! 
\alpha_{k})^{-1} \vert^{-1} \! \geqslant \! -\dfrac{1}{2} \ln (1 \! + \! 
(\xi \! - \! \alpha_{k})^{-2}) \! - \! \dfrac{1}{2} \ln (1 \! + \! (\tau \! 
- \! \alpha_{k})^{-2}), \\
\ln \vert (\xi \! - \! \alpha_{p_{q}})^{-1} \! - \! (\tau \! - \! 
\alpha_{p_{q}})^{-1} \vert^{-1} \! \geqslant \! -\dfrac{1}{2} \ln (1 \! + \! 
(\xi \! - \! \alpha_{p_{q}})^{-2}) \! - \! \dfrac{1}{2} \ln (1 \! + \! (\tau 
\! - \! \alpha_{p_{q}})^{-2}),
\end{gather*}
one shows that, for $n \! \in \! \mathbb{N}$ and $k \! \in \! \lbrace 1,2,
\dotsc,K \rbrace$ such that $\alpha_{p_{\mathfrak{s}}} \! := \! \alpha_{k} 
\! \neq \! \infty$,
\begin{align*}
\mathcal{K}_{\widetilde{V}}^{f}(\xi,\tau) \geqslant& \, \dfrac{1}{2} \left(
\widetilde{V}(\xi) \! - \! \left(\dfrac{\varkappa_{nk \tilde{k}_{\mathfrak{s}-
1}}^{\infty} \! + \! 1}{n} \right) \ln (1 \! + \! \xi^{2}) \! - \! \left(
\dfrac{\varkappa_{nk} \! - \! 1}{n} \right) \ln (1 \! + \! (\xi \! - \! 
\alpha_{k})^{-2}) \! - \! \sum_{q=1}^{\mathfrak{s}-2} \dfrac{\varkappa_{nk 
\tilde{k}_{q}}}{n} \ln (1 \! + \! (\xi \! - \! \alpha_{p_{q}})^{-2}) \right) \\
+& \, \dfrac{1}{2} \left(\widetilde{V}(\tau) \! - \! \left(\dfrac{\varkappa_{nk 
\tilde{k}_{\mathfrak{s}-1}}^{\infty} \! + \! 1}{n} \right) \ln (1 \! + \! 
\tau^{2}) \! - \! \left(\dfrac{\varkappa_{nk} \! - \! 1}{n} \right) \ln (1 \! 
+ \! (\tau \! - \! \alpha_{k})^{-2}) \! - \! \sum_{q=1}^{\mathfrak{s}-2} 
\dfrac{\varkappa_{nk \tilde{k}_{q}}}{n} \ln (1 \! + \! (\tau \! - \! 
\alpha_{p_{q}})^{-2}) \right).
\end{align*}
Recalling conditions~\eqref{eq20}--\eqref{eq22} for the external field 
$\widetilde{V} \colon \overline{\mathbb{R}} \setminus \lbrace \alpha_{1},
\alpha_{2},\dotsc,\alpha_{K} \rbrace \! \to \! \mathbb{R}$, in particular, 
for $x \! \in \! \mathscr{O}_{\infty} \! := \! \lbrace \mathstrut x \! \in 
\! \mathbb{R}; \, \lvert x \rvert \! > \! \delta^{-1}_{\infty} \rbrace$, where 
$\delta_{\infty}$ $(= \! \delta_{\infty}(n,k,z_{o}))$ is some arbitrarily fixed, 
sufficiently small positive real number,{}\footnote{For example, $\delta_{\infty}^{-1} 
\! = \! K(1 \! + \! \max \lbrace \mathstrut \lvert \alpha_{p_{q}} \rvert, \, q \! 
= \! 1,\dotsc,\mathfrak{s} \! - \! 2,\mathfrak{s} \rbrace \! + \! \max \lbrace 
\mathstrut \lvert \alpha_{p_{i}} \! - \! \alpha_{p_{j}} \rvert, \, i \! \neq \! j \! 
\in \! \lbrace 1,\dotsc,\mathfrak{s} \! - \! 2,\mathfrak{s} \rbrace \rbrace)$.} 
$\widetilde{V}(x) \! \geqslant \! (1 \! + \! \mathfrak{c}_{\infty}) \ln (1 \! + \! 
x^{2})$, where $\mathfrak{c}_{\infty}$ $(= \! \mathfrak{c}_{\infty}(n,k,z_{o}))$ 
is some bounded, positive real number, and, for $x \! \in \! \mathscr{O}_{
\tilde{\delta}_{q}}(\alpha_{p_{q}}) \! := \! \lbrace \mathstrut x \! \in \! 
\mathbb{R}; \, \vert x \! - \! \alpha_{p_{q}} \vert \! < \! \tilde{\delta}_{q} 
\rbrace$, $q \! = \! 1,\dotsc,\mathfrak{s} \! - \! 2,\mathfrak{s}$, where 
$\tilde{\delta}_{q}$ $(= \! \tilde{\delta}_{q}(n,k,z_{o}))$, $q \! = \! 1,\dotsc,
\mathfrak{s} \! - \! 2,\mathfrak{s}$, are arbitrarily fixed, sufficiently small 
positive real numbers chosen so that {}\footnote{For example, $\tilde{\delta}_{q} 
\! < \! (3K)^{-1} \min \lbrace \mathstrut \lvert \alpha_{p_{i}} \! - \! \alpha_{p_{j}} 
\rvert, \, i \! \neq \! j \! \in \! \lbrace 1,\dotsc,\mathfrak{s} \! - \! 2,\mathfrak{s} 
\rbrace \rbrace$.} $\mathscr{O}_{\tilde{\delta}_{q^{\prime}}}(\alpha_{p_{q^{\prime}}}) 
\cap \mathscr{O}_{\tilde{\delta}_{q^{\prime \prime}}}(\alpha_{p_{q^{\prime \prime}}}) 
\! = \! \varnothing$ $\forall$ $q^{\prime} \! \neq \! q^{\prime \prime} \! \in \! 
\lbrace 1,\dotsc,\mathfrak{s} \! - \! 2,\mathfrak{s} \rbrace$, $\widetilde{V}(x) \! 
\geqslant \! (1 \! + \! \mathfrak{c}_{q}) \ln (1 \! + \! (x \! - \! \alpha_{p_{q}})^{-2})$, 
where $\mathfrak{c}_{q}$ $(= \! \mathfrak{c}_{q}(n,k,z_{o}))$, $q \! = \! 1,\dotsc,
\mathfrak{s} \! - \! 2,\mathfrak{s}$, are bounded, positive real numbers, it follows 
that, for $n \! \in \! \mathbb{N}$ and $k \! \in \! \lbrace 1,2,\dotsc,K \rbrace$ such 
that $\alpha_{p_{\mathfrak{s}}} \! := \! \alpha_{k} \! \neq \! \infty$,
\begin{equation*}
\widetilde{V}(x) \! - \! \left(\dfrac{\varkappa_{nk \tilde{k}_{\mathfrak{s}
-1}}^{\infty} \! + \! 1}{n} \right) \ln (1 \! + \! x^{2}) \! - \! \left(
\dfrac{\varkappa_{nk} \! - \! 1}{n} \right) \ln (1 \! + \! (x \! - \! 
\alpha_{k})^{-2}) \! - \! \sum_{q=1}^{\mathfrak{s}-2} \dfrac{
\varkappa_{nk \tilde{k}_{q}}}{n} \ln (1 \! + \! (x \! - \! \alpha_{p_{q}})^{-2}) 
\! \geqslant \! \hat{C}_{\widetilde{V}}^{f}(n,k,z_{o}) \! := \! \hat{C}_{
\widetilde{V}}^{f} \! > \! -\infty,
\end{equation*}
whence $\mathcal{K}_{\widetilde{V}}^{f}(\xi,\tau) \! \geqslant \! 
\hat{C}_{\widetilde{V}}^{f}$ $(> \! -\infty)$, and, consequently,
\begin{equation}
\mathrm{I}_{\widetilde{V}}^{f}[\mu^{\text{\tiny $\mathrm{EQ}$}}_{2}] 
\! \geqslant \! \iint_{\mathbb{R}^{2}} \hat{C}_{\widetilde{V}}^{f} 
\, \md \mu^{\text{\tiny $\mathrm{EQ}$}}_{2}(\xi) \, \md 
\mu^{\text{\tiny $\mathrm{EQ}$}}_{2}(\tau) \! \geqslant \! 
\hat{C}_{\widetilde{V}}^{f} \quad (> \! -\infty). \label{eqKvinf2}
\end{equation}
It follows {}from the above inequality and the definition of $E_{\widetilde{
V}}^{f}$ stated in the lemma that, for $n \! \in \! \mathbb{N}$ and $k \! 
\in \! \lbrace 1,2,\dotsc,K \rbrace$ such that $\alpha_{p_{\mathfrak{s}}} 
\! := \! \alpha_{k} \! \neq \! \infty$, and for all 
$\mu^{\text{\tiny $\mathrm{EQ}$}}_{2} \! \in \! \mathscr{M}_{1}
(\mathbb{R})$, $E_{\widetilde{V}}^{f} \! \geqslant \! \hat{C}_{
\widetilde{V}}^{f} \! > \! -\infty$, that is, $E_{\widetilde{V}}^{f}$ is 
bounded {}from below. Let $\varepsilon$ $(= \! \varepsilon (n,k,z_{o}))$ 
be an arbitrarily fixed, sufficiently small positive real number, and set 
$\hat{\sigma}_{\varepsilon} \! := \! \lbrace \mathstrut x \! \in \! 
\mathbb{R}; \, \varpi (x) \! \geqslant \! \varepsilon \rbrace$; then 
$\hat{\sigma}_{\varepsilon}$ is compact, and $\hat{\sigma}_{0} \! := 
\! \cup_{j=1}^{\infty} \hat{\sigma}_{1/j} \! = \! \lbrace \mathstrut 
x \! \in \! \mathbb{R}; \, \varpi (x) \! > \! 0 \rbrace$. Since, for 
$\widetilde{V} \colon \overline{\mathbb{R}} \setminus \lbrace 
\alpha_{1},\alpha_{2},\dotsc,\alpha_{K} \rbrace \! \to \! \mathbb{R}$ 
satisfying conditions~\eqref{eq20}--\eqref{eq22}, $\varpi$ is an 
\emph{admissible weight} {}\footnote{For $\widetilde{V} \colon 
\overline{\mathbb{R}} \setminus \lbrace \alpha_{1},\alpha_{2},
\dotsc,\alpha_{K} \rbrace \! \to \! \mathbb{R}$ satisfying 
conditions~\eqref{eq20}--\eqref{eq22}, $\varpi$ is an admissible 
weight means: (i) $\varpi (x)$ is upper semi-continuous (u.s.c.) 
on $\overline{\mathbb{R}} \setminus \lbrace \alpha_{1},\alpha_{2},
\dotsc,\alpha_{K} \rbrace$; (ii) $\lbrace \mathstrut x \! \in \! 
\mathbb{R}; \, \varpi (x) \! > \! 0 \rbrace$ has 
positive---weighted logarithmic---capacity, that is, $\operatorname{cap}
(\lbrace \mathstrut x \! \in \! \mathbb{R}; \, \varpi (x) \! > \! 0 \rbrace) 
\! := \! \exp (-\inf \lbrace \mathstrut \mathrm{I}_{\widetilde{V}}^{f}
[\mu^{\text{\tiny $\mathrm{EQ}$}}_{2}]; \, 
\mu^{\text{\tiny $\mathrm{EQ}$}}_{2} \! \in \! \mathscr{M}_{1}(\lbrace 
\mathstrut x \! \in \! \mathbb{R}; \, \varpi (x) \! > \! 0 \rbrace) 
\rbrace) \! > \! 0$; (iii) $\lvert x \rvert \varpi (x) \! \to \! 0$ 
as $\lvert x \rvert \! \to \! +\infty$; and (iv) $\lvert x \! - \! 
\alpha_{p_{q}} \rvert^{-1} \varpi (x) \! \to \! 0$ as $x \! \to \! 
\alpha_{p_{q}}$, $q \! = \! 1,\dotsc,\mathfrak{s} \! - \! 2,\mathfrak{s}$. 
(See, also, \cite{plsi,ebs,tbnlevvikfiel}.)} (see, for example, \cite{a55}), it 
follows that, for $n \! \in \! \mathbb{N}$ and $k \! \in \! \lbrace 1,2,
\dotsc,K \rbrace$ such that $\alpha_{p_{\mathfrak{s}}} \! := \! \alpha_{k} 
\! \neq \! \infty$, there exits $j^{\ast}$ $(= \! j^{\ast}(n,k,z_{o}))$ $\in \! 
\mathbb{N}$ such that $\operatorname{cap}(\hat{\sigma}_{1/j^{\ast}}) 
\! := \! \exp (-\inf \lbrace \mathstrut \mathrm{I}_{\widetilde{V}}^{f}
[\mu^{\text{\tiny $\mathrm{EQ}$}}_{2}]; \, 
\mu^{\text{\tiny $\mathrm{EQ}$}}_{2} \! \in \! \mathscr{M}_{1}
(\hat{\sigma}_{1/j^{\ast}}) \rbrace) \! > \! 0$, which, in turn, means 
that there exists an associated probability measure, denoted 
$\mu_{2,j^{\ast}}^{\text{\tiny $\mathrm{EQ}$}}$ $(= \! 
\mu_{2,j^{\ast}}^{\text{\tiny $\mathrm{EQ}$}}(n,k,z_{o}))$, with 
$\supp (\mu_{2,j^{\ast}}^{\text{\tiny $\mathrm{EQ}$}}) \subseteq 
\hat{\sigma}_{1/j^{\ast}}$, such that
\begin{equation*}
\iint_{\hat{\sigma}_{1/j^{\ast}} \, \times \, \hat{\sigma}_{1/j^{\ast}}} 
\ln \left(\vert \xi \! - \! \tau \vert^{\frac{\varkappa_{nk \tilde{k}_{
\mathfrak{s}-1}}^{\infty}+1}{n}} \left(\dfrac{\vert \xi \! - \! \tau 
\vert}{\vert \xi \! - \! \alpha_{k} \vert \vert \tau \! - \! \alpha_{k} \vert} 
\right)^{\frac{\varkappa_{nk}-1}{n}} \prod_{q=1}^{\mathfrak{s}-2} \left(
\dfrac{\vert \xi \! - \! \tau \vert}{\vert \xi \! - \! \alpha_{p_{q}} \vert 
\vert \tau \! - \! \alpha_{p_{q}} \vert} \right)^{\frac{\varkappa_{nk 
\tilde{k}_{q}}}{n}} \right)^{-1} \md 
\mu_{2,j^{\ast}}^{\text{\tiny $\mathrm{EQ}$}}(\xi) \, \md 
\mu_{2,j^{\ast}}^{\text{\tiny $\mathrm{EQ}$}}(\tau) \! < \! +\infty;
\end{equation*}
furthermore, for $x \! \in \! \supp 
(\mu_{2,j^{\ast}}^{\text{\tiny $\mathrm{EQ}$}}) \subseteq \hat{\sigma}_{
1/j^{\ast}}$, it follows that $\varpi (x) \! \geqslant \! 1/j^{\ast}$, whence
\begin{equation*}
\iint_{\hat{\sigma}_{1/j^{\ast}} \, \times \, \hat{\sigma}_{1/j^{\ast}}} 
\ln (\varpi (\xi) \varpi (\tau))^{-1} \, \md 
\mu_{2,j^{\ast}}^{\text{\tiny $\mathrm{EQ}$}}(\xi) \, \md 
\mu_{2,j^{\ast}}^{\text{\tiny $\mathrm{EQ}$}}(\tau) \! \leqslant \! 
2 \ln (j^{\ast}) \! < \! +\infty.
\end{equation*}
Hence, it follows that, for $n \! \in \! \mathbb{N}$ and $k \! \in \! \lbrace 
1,2,\dotsc,K \rbrace$ such that $\alpha_{p_{\mathfrak{s}}} \! := \! \alpha_{k} 
\! \neq \! \infty$,
\begin{align}
\mathrm{I}_{\widetilde{V}}^{f}[\mu_{2,j^{\ast}}^{\text{\tiny $\mathrm{EQ}$}}] 
=& \, \iint_{\hat{\sigma}_{1/j^{\ast}} \, \times \, \hat{\sigma}_{1/j^{\ast}}} 
\ln \left(\vert \xi \! - \! \tau \vert^{\frac{\varkappa_{nk \tilde{k}_{
\mathfrak{s}-1}}^{\infty}+1}{n}} \left(\dfrac{\vert \xi \! - \! \tau \vert}{
\vert \xi \! - \! \alpha_{k} \vert \vert \tau \! - \! \alpha_{k} \vert} 
\right)^{\frac{\varkappa_{nk}-1}{n}} \prod_{q=1}^{\mathfrak{s}-2} \left(
\dfrac{\vert \xi \! - \! \tau \vert}{\vert \xi \! - \! \alpha_{p_{q}} \vert 
\vert \tau \! - \! \alpha_{p_{q}} \vert} \right)^{\frac{\varkappa_{nk 
\tilde{k}_{q}}}{n}} \varpi (\xi) \varpi (\tau) \right)^{-1} \nonumber \\
\times& \, \md \mu_{2,j^{\ast}}^{\text{\tiny $\mathrm{EQ}$}}(\xi) \, \md 
\mu_{2,j^{\ast}}^{\text{\tiny $\mathrm{EQ}$}}(\tau) \! < \! +\infty, 
\label{eqKvinf3}
\end{align}
that is, via inequalities~\eqref{eqKvinf2} and~\eqref{eqKvinf3},
\begin{equation*}
-\infty \! < \! E_{\widetilde{V}}^{f} \! = \! \inf \lbrace \mathstrut 
\mathrm{I}_{\widetilde{V}}^{f}[\mu^{\text{\tiny $\mathrm{EQ}$}}_{2}]; \, 
\mu^{\text{\tiny $\mathrm{EQ}$}}_{2} \! \in \! \mathscr{M}_{1}
(\mathbb{R}) \rbrace \! < \! +\infty,
\end{equation*}
which establishes, for $n \! \in \! \mathbb{N}$ and $k \! \in \! \lbrace 
1,2,\dotsc,K \rbrace$ such that $\alpha_{p_{\mathfrak{s}}} \! := \! 
\alpha_{k} \! \neq \! \infty$, the corresponding claim~(1) of the lemma.

For $n \! \in \! \mathbb{N}$ and $k \! \in \! \lbrace 1,2,\dotsc,K 
\rbrace$ such that $\alpha_{p_{\mathfrak{s}}} \! := \! \alpha_{k} \! 
\neq \! \infty$, choose a sequence of probability measures $\lbrace 
\mu_{2,m}^{\text{\tiny $\mathrm{EQ}$}} \rbrace_{m=1}^{\infty}$ 
in $\mathscr{M}_{1}(\mathbb{R})$, with 
$\mu_{2,m}^{\text{\tiny $\mathrm{EQ}$}} \! = \! 
\mu_{2,m}^{\text{\tiny $\mathrm{EQ}$}}(n,k,z_{o})$, such that $\mathrm{I}_{
\widetilde{V}}^{f}[\mu_{2,m}^{\text{\tiny $\mathrm{EQ}$}}] \! \leqslant 
\! E_{\widetilde{V}}^{f} \! + \! \tfrac{1}{m}$. {}From the analysis above, 
it follows that, for $n \! \in \! \mathbb{N}$ and $k \! \in \! \lbrace 
1,2,\dotsc,K \rbrace$ such that $\alpha_{p_{\mathfrak{s}}} \! := \! 
\alpha_{k} \! \neq \! \infty$,
\begin{equation*}
\mathrm{I}_{\widetilde{V}}^{f}[\mu_{2,m}^{\text{\tiny $\mathrm{EQ}$}}] 
\! = \! \iint_{\mathbb{R}^{2}} \mathcal{K}_{\widetilde{V}}^{f}
(\xi,\tau) \, \md \mu_{2,m}^{\text{\tiny $\mathrm{EQ}$}}(\xi) \, \md 
\mu_{2,m}^{\text{\tiny $\mathrm{EQ}$}}(\tau) \! \geqslant \! \dfrac{1}{2} 
\iint_{\mathbb{R}^{2}} \left(\hat{\psi}_{\widetilde{V}}^{f}(\xi) 
\! + \! \hat{\psi}_{\widetilde{V}}^{f}(\tau) \right) \md 
\mu_{2,m}^{\text{\tiny $\mathrm{EQ}$}}(\xi) \, \md 
\mu_{2,m}^{\text{\tiny $\mathrm{EQ}$}}(\tau),
\end{equation*}
where
\begin{equation}
\hat{\psi}_{\widetilde{V}}^{f}(z) \! := \! \widetilde{V}(z) \! - \! \left(
\dfrac{\varkappa^{\infty}_{nk \tilde{k}_{\mathfrak{s}-1}} \! + \! 1}{n} 
\right) \ln (1 \! + \! z^{2}) \! - \! \left(\dfrac{\varkappa_{nk} \! - \! 1}{n} 
\right) \ln (1 \! + \! (z \! - \! \alpha_{k})^{-2}) \! - \! \sum_{q=1}^{
\mathfrak{s}-2} \dfrac{\varkappa_{nk \tilde{k}_{q}}}{n} \ln (1 \! + \! 
(z \! - \! \alpha_{p_{q}})^{-2}). \label{eqKvinf4}
\end{equation}
Then, for $n \! \in \! \mathbb{N}$ and $k \! \in \! \lbrace 
1,2,\dotsc,K \rbrace$ such that $\alpha_{p_{\mathfrak{s}}} 
\! := \! \alpha_{k} \! \neq \! \infty$, it follows that, for 
$\mu_{2,m}^{\text{\tiny $\mathrm{EQ}$}} \! \in \! \mathscr{M}_{1}
(\mathbb{R})$, $\mathrm{I}_{\widetilde{V}}^{f}
[\mu_{2,m}^{\text{\tiny $\mathrm{EQ}$}}] \! \geqslant \! 
\int_{\mathbb{R}} \hat{\psi}_{\widetilde{V}}^{f}(\tau) \, 
\md \mu_{2,m}^{\text{\tiny $\mathrm{EQ}$}}(\tau)$, whence
\begin{equation*}
\int_{\mathbb{R}} \hat{\psi}_{\widetilde{V}}^{f}(\lambda) \, \md 
\mu_{2,m}^{\text{\tiny $\mathrm{EQ}$}}(\lambda) \! \leqslant \! 
\mathrm{I}_{\widetilde{V}}^{f}[\mu_{2,m}^{\text{\tiny $\mathrm{EQ}$}}] 
\! \leqslant \! E_{\widetilde{V}}^{f} \! + \! \dfrac{1}{m}, 
\quad m \! \in \! \mathbb{N}.
\end{equation*}
For $n \! \in \! \mathbb{N}$ and $k \! \in \! \lbrace 1,2,\dotsc,K 
\rbrace$ such that $\alpha_{p_{\mathfrak{s}}} \! := \! \alpha_{k} \! 
\neq \! \infty$, recalling that there exists $\mathfrak{c}_{\infty} \! 
> \! 0$ and bounded such that, for $x \! \in \! \mathscr{O}_{\infty}$, 
$\widetilde{V}(x) \! \geqslant \! (1 \! + \! \mathfrak{c}_{\infty}) \ln 
(1 \! + \! x^{2})$, and, for $q \! = \! 1,\dotsc,\mathfrak{s} \! - \! 2,
\mathfrak{s}$, there exist $\mathfrak{c}_{q} \! > \! 0$ and bounded 
such that, for $x \! \in \! \mathscr{O}_{\tilde{\delta}_{q}}(\alpha_{p_{q}})$, 
$\widetilde{V}(x) \! \geqslant \! (1 \! + \! \mathfrak{c}_{q}) \ln (1 \! + \! 
(x \! - \! \alpha_{p_{q}})^{-2})$, it follows that, for any $\mathfrak{b}_{f}$ 
$(= \! \mathfrak{b}_{f}(n,k,z_{o}))$ $> \! 0$, there exists $M_{f}$ $(= \! 
M_{f}(n,k,z_{o}))$  $> \! 1$ {}\footnote{For example, $M_{f} \! = \! K(1 \! 
+ \! \max \lbrace \mathstrut \lvert \alpha_{p_{q}} \rvert, \, q \! = \! 1,
\dotsc,\mathfrak{s} \! - \! 2,\mathfrak{s} \rbrace \! + \! 3(\min \lbrace 
\mathstrut \lvert \alpha_{p_{i}} \! - \! \alpha_{p_{j}} \rvert, \, i \! \neq 
\! j \! \in \! \lbrace 1,\dotsc,\mathfrak{s} \! - \! 2,\mathfrak{s} \rbrace 
\rbrace)^{-1})$.} for which $\mathscr{O}_{\frac{1}{M_{f}}}(\alpha_{p_{i}}) 
\cap \mathscr{O}_{\frac{1}{M_{f}}}(\alpha_{p_{j}}) \! = \! \varnothing$, 
$i \! \neq \! j \! \in \! \lbrace 1,\dotsc,\mathfrak{s} \! - \! 2,\mathfrak{s} 
\rbrace$, such that, for $z \! \in \! \mathfrak{D}_{M_{f}} \! := \! \lbrace 
\lvert x \rvert \! \geqslant \! M_{f} \rbrace \cup \cup_{\underset{q \neq 
\mathfrak{s}-1}{q=1}}^{\mathfrak{s}} \operatorname{clos}(\mathscr{O}_{
\frac{1}{M_{f}}}(\alpha_{p_{q}}))$, $\hat{\psi}_{\widetilde{V}}^{f}(z) \! 
> \! b_{f}$, which implies that, upon writing $\mathbb{R} \! = \! 
(\mathbb{R} \setminus \mathfrak{D}_{M_{f}}) \cup \mathfrak{D}_{M_{f}}$ 
(with $(\mathbb{R} \setminus \mathfrak{D}_{M_{f}}) \cap \mathfrak{D}_{M_{f}} 
\! = \! \varnothing)$,
\begin{align*}
E_{\widetilde{V}}^{f} \! + \! \dfrac{1}{m} \geqslant& \int_{\mathbb{R}} 
\hat{\psi}_{\widetilde{V}}^{f}(\lambda) \, \md 
\mu_{2,m}^{\text{\tiny $\mathrm{EQ}$}}(\lambda) \! = \! \int_{
\mathbb{R} \setminus \mathfrak{D}_{M_{f}}} \underbrace{\hat{\psi}_{
\widetilde{V}}^{f}(\lambda)}_{\geqslant -\vert \hat{C}_{\widetilde{V}}^{f} 
\vert} \, \md \mu_{2,m}^{\text{\tiny $\mathrm{EQ}$}}(\lambda) \! + \! 
\int_{\mathfrak{D}_{M_{f}}} \underbrace{\hat{\psi}_{\widetilde{V}}^{f}
(\lambda)}_{>\mathfrak{b}_{f}} \, \md \mu_{2,m}^{\text{\tiny $\mathrm{EQ}$}}
(\lambda) \\
\geqslant& \, \mathfrak{b}_{f} \int_{\mathfrak{D}_{M_{f}}} \md 
\mu_{2,m}^{\text{\tiny $\mathrm{EQ}$}}(\lambda) \! - \! \vert \hat{C}_{
\widetilde{V}}^{f} \vert,
\end{align*}
whence, for $n \! \in \! \mathbb{N}$ and $k \! \in \! \lbrace 1,2,\dotsc,
K \rbrace$ such that $\alpha_{p_{\mathfrak{s}}} \! := \! \alpha_{k} \! 
\neq \! \infty$, choosing $\mathfrak{b}_{f}$ judiciously {}\footnote{For 
example, $\mathfrak{b}_{f} \! \gg \! E_{\widetilde{V}}^{f} \! + \! \vert 
\hat{C}_{\widetilde{V}}^{f} \vert \! + \! \tfrac{1}{m}$, $m \! \in \! 
\mathbb{N}$.} and setting $\hat{\varepsilon}$ $(= \! \hat{\varepsilon}
(n,k,z_{o}))$ $:= \! \limsup_{m \to \infty}((E_{\widetilde{V}}^{f} \! + \! 
\vert \hat{C}_{\widetilde{V}}^{f} \vert \! + \! m^{-1})/\mathfrak{b}_{f})$,
\begin{equation*}
\limsup_{m \to \infty} \int_{\mathfrak{D}_{M_{f}}} \md 
\mu_{2,m}^{\text{\tiny $\mathrm{EQ}$}}(\lambda) \! \leqslant \! 
\hat{\varepsilon},
\end{equation*}
that is, the sequence of probability measures $\lbrace 
\mu_{2,m}^{\text{\tiny $\mathrm{EQ}$}} \rbrace_{m=1}^{\infty}$ 
in $\mathscr{M}_{1}(\mathbb{R})$ is \emph{tight} \cite{a56}.\footnote{That 
is, for $n \! \in \! \mathbb{N}$ and $k \! \in \! \lbrace 1,2,\dotsc,K \rbrace$ 
such that $\alpha_{p_{\mathfrak{s}}} \! := \! \alpha_{k} \! \neq \! \infty$, 
given $\hat{\varepsilon}$ $(= \! \hat{\varepsilon}(n,k,z_{o}))$ $> \! 0$, 
there exists $M_{f}$ $(= \! M_{f}(n,k,z_{o}))$ $> \! 1$ for which 
$\mathscr{O}_{\frac{1}{M_{f}}}(\alpha_{p_{i}}) \cap \mathscr{O}_{\frac{1}{M_{f}}}
(\alpha_{p_{j}}) \! = \! \varnothing$, $i \! \neq \! j \! \in \! \lbrace 1,\dotsc,
\mathfrak{s} \! - \! 2,\mathfrak{s} \rbrace$, such that $\limsup_{m \to \infty} 
\mu_{2,m}^{\text{\tiny $\mathrm{EQ}$}}(\mathfrak{D}_{M_{f}}) \! 
\leqslant \! \hat{\varepsilon}$, where 
$\mu_{2,m}^{\text{\tiny $\mathrm{EQ}$}}(\mathfrak{D}_{M_{f}}) 
\! := \! \int_{\mathfrak{D}_{M_{f}}} \md \mu_{2,m}^{\text{\tiny 
$\mathrm{EQ}$}}(\lambda)$.} Since, for $n \! \in \! \mathbb{N}$ and $k \! 
\in \! \lbrace 1,2,\dotsc,K \rbrace$ such that $\alpha_{p_{\mathfrak{s}}} 
\! := \! \alpha_{k} \! \neq \! \infty$, the sequence of probabilty measures 
$\lbrace \mu_{2,m}^{\text{\tiny $\mathrm{EQ}$}} \rbrace_{m=1}^{\infty}$ 
in $\mathscr{M}_{1}(\mathbb{R})$ is tight, there exists, by a Helly selection 
theorem (see, for example, \cite{a55}), a weak-$\ast$ convergent subsequence 
of probability measures $\lbrace \mu_{2,m_{j}}^{\text{\tiny $\mathrm{EQ}$}} 
\rbrace_{j=1}^{\infty}$ in $\mathscr{M}_{1}(\mathbb{R})$ (with 
$\mu_{2,m_{j}}^{\text{\tiny $\mathrm{EQ}$}} \! = \! 
\mu_{2,m_{j}}^{\text{\tiny $\mathrm{EQ}$}}(n,k,z_{o})$ and $j \! = \! j(n,k,z_{o}))$ 
converging weakly to a probability measure $\mu_{f}^{\text{\tiny $\mathrm{EQ}$}} 
\! \in \! \mathscr{M}_{1}(\mathbb{R})$ (with $\mu_{f}^{\text{\tiny $\mathrm{EQ}$}} 
\! = \! \mu_{f}^{\text{\tiny $\mathrm{EQ}$}}(n,k,z_{o}))$, symbolically 
$\mu_{2,m_{j}}^{\text{\tiny $\mathrm{EQ}$}} \! \overset{\ast}{\to} \! 
\mu_{f}^{\text{\tiny $\mathrm{EQ}$}}$ as $j \! \to \! \infty$. One now shows 
that, if $\mu_{2,m}^{\text{\tiny $\mathrm{EQ}$}} \! \overset{\ast}{\to} 
\! \mu_{f}^{\text{\tiny $\mathrm{EQ}$}}$, then $\liminf_{m \to \infty} 
\mathrm{I}_{\widetilde{V}}^{f}[\mu_{2,m}^{\text{\tiny $\mathrm{EQ}$}}] 
\! \geqslant \! \mathrm{I}_{\widetilde{V}}^{f}
[\mu_{f}^{\text{\tiny $\mathrm{EQ}$}}]$. Since, for $n \! \in \! 
\mathbb{N}$ and $k \! \in \! \lbrace 1,2,\dotsc,K \rbrace$ such 
that $\alpha_{p_{\mathfrak{s}}} \! := \! \alpha_{k} \! \neq \! \infty$, 
$\varpi$ (resp., $\widetilde{V})$ is upper semi-continuous (resp., 
lower semi-continuous) on $\overline{\mathbb{R}} \setminus \lbrace 
\alpha_{1},\alpha_{2},\dotsc,\alpha_{K} \rbrace$, $\varpi$ (resp., 
$\widetilde{V})$ is the pointwise limit of a decreasing (resp., an 
increasing) sequence of positive {}\footnote{Adding a suitable 
constant, if necessary, which does not change 
$\mu_{2,m}^{\text{\tiny $\mathrm{EQ}$}}$ $(\in \! \mathscr{M}_{1}
(\mathbb{R}))$, or the regularity of $\widetilde{V} \colon 
\overline{\mathbb{R}} \setminus \lbrace \alpha_{1},\alpha_{2},\dotsc,
\alpha_{K} \rbrace \! \to \! \mathbb{R}$, one may assume that 
$\widetilde{V} \! \geqslant \! 0$ and $\widetilde{V}_{m} \! \geqslant 
\! 0$, $m \! \in \! \mathbb{N}$.} continuous functions $\lbrace 
\varpi_{m} \rbrace_{m=1}^{\infty}$ (resp., $\lbrace \widetilde{V}_{m} 
\rbrace_{m=1}^{\infty})$ on $\overline{\mathbb{R}} \setminus \lbrace 
\alpha_{1},\alpha_{2},\dotsc,\alpha_{K} \rbrace$, that is, $\varpi_{m}
(x) \! \searrow \! \varpi(x)$ (resp., $\widetilde{V}_{m}(x) \! \nearrow 
\! \widetilde{V}(x))$ as $m \! \to \! \infty$ for $x \! \in \! 
\overline{\mathbb{R}} \setminus \lbrace \alpha_{1},\alpha_{2},
\dotsc,\alpha_{K} \rbrace$. Noting that, for $n \! \in \! \mathbb{N}$ 
and $k \! \in \! \lbrace 1,2,\dotsc,K \rbrace$ such that 
$\alpha_{p_{\mathfrak{s}}} \! := \! \alpha_{k} \! \neq \! \infty$, 
$\mathrm{I}_{\widetilde{V}}^{f}[\mu_{2,m}^{\text{\tiny $\mathrm{EQ}$}}] 
\! = \! \iint_{\mathbb{R}^{2}} \mathcal{K}_{\widetilde{V}}^{f}
(\xi,\tau) \, \md \mu_{2,m}^{\text{\tiny $\mathrm{EQ}$}}(\xi) \, 
\md \mu_{2,m}^{\text{\tiny $\mathrm{EQ}$}}(\tau) \! \geqslant \! 
\iint_{\mathbb{R}^{2}} \mathcal{K}_{\widetilde{V}_{m}}^{f}(\xi,\tau) 
\, \md \mu_{2,m}^{\text{\tiny $\mathrm{EQ}$}}(\xi) \, \md 
\mu_{2,m}^{\text{\tiny $\mathrm{EQ}$}}(\tau)$, it follows that, for any 
$L_{f}$ $(= \! L_{f}(n,k,z_{o}))$ $\in \! \mathbb{R}$, $\mathrm{I}_{
\widetilde{V}}^{f}[\mu_{2,m}^{\text{\tiny $\mathrm{EQ}$}}] \! 
\geqslant \! \iint_{\mathbb{R}^{2}} \mathfrak{g}_{f}(\xi,\tau) \, \md 
\mu_{2,m}^{\text{\tiny $\mathrm{EQ}$}}(\xi) \, \md 
\mu_{2,m}^{\text{\tiny $\mathrm{EQ}$}}(\tau)$, where $\mathfrak{g}_{f} 
\colon \mathbb{N} \times \lbrace 1,2,\dotsc,K \rbrace \times 
\mathbb{R}^{2} \! \ni \! (n,k,\xi,\tau) \! \mapsto \! 
\mathfrak{g}_{f}(n,k,\xi,\tau) \! := \! \mathfrak{g}_{f}(\xi,\tau) 
\! = \! \mathfrak{g}_{f}(\tau,\xi) \! = \! \min \lbrace L_{f},
\mathcal{K}_{\widetilde{V}_{m}}^{f}(\xi,\tau) \rbrace$ is,\footnote{Recall 
that $\min \lbrace f_{1},f_{2} \rbrace \! = \! \tfrac{1}{2}(f_{1} \! + \! 
f_{2}) \! - \! \tfrac{1}{2} \lvert f_{1} \! - \! f_{2} \rvert$.} by virtue of 
conditions~\eqref{eq20}--\eqref{eq22} on the external field $\widetilde{V} 
\colon \overline{\mathbb{R}} \setminus \lbrace \alpha_{1},\alpha_{2},\dotsc,
\alpha_{K} \rbrace \! \to \! \mathbb{R}$, bounded and continuous on 
$\mathbb{R}^{2}$. Let $\hat{\varepsilon}$ $(= \! \hat{\varepsilon}(n,k,z_{o}))$ 
$> \! 0$ be given, and choose $M_{0}$ $(= \! M_{0}(n,k,z_{o}))$ $> \! 1$ 
for which $\mathscr{O}_{\frac{1}{M_{0}}}(\alpha_{p_{i}}) \cap 
\mathscr{O}_{\frac{1}{M_{0}}}(\alpha_{p_{j}}) \! = \! \varnothing$, $i \! 
\neq \! j \! \in \! \lbrace 1,\dotsc,\mathfrak{s} \! - \! 2,\mathfrak{s} 
\rbrace$ (e.g., $M_{0} \! = \! K(1 \! + \! \max \lbrace \mathstrut \lvert 
\alpha_{p_{q}} \rvert, \, q \! = \! 1,\dotsc,\mathfrak{s} \! - \! 2,
\mathfrak{s} \rbrace \! + \! 3(\min \lbrace \mathstrut \lvert \alpha_{p_{i}} 
\! - \! \alpha_{p_{j}} \rvert, \, i \! \neq \! j \! \in \! \lbrace 1,\dotsc,
\mathfrak{s} \! - \! 2,\mathfrak{s} \rbrace \rbrace)^{-1}))$, such that 
$\limsup_{m \to \infty} \int_{\mathfrak{D}_{M_{0}}} \md 
\mu_{2,m}^{\text{\tiny $\mathrm{EQ}$}}(\lambda) \! \leqslant \! 
\hat{\varepsilon}$, where $\mathfrak{D}_{M_{0}} \! := \! \lbrace \lvert x 
\rvert \! \geqslant \! M_{0} \rbrace \cup \cup_{\underset{q \neq \mathfrak{s}
-1}{q=1}}^{\mathfrak{s}} \operatorname{clos}(\mathscr{O}_{\frac{1}{M_{0}}}
(\alpha_{p_{q}}))$. For $n \! \in \! \mathbb{N}$ and $k \! \in \! \lbrace 
1,2,\dotsc,K \rbrace$ such that $\alpha_{p_{\mathfrak{s}}} 
\! := \! \alpha_{k} \! \neq \! \infty$, let the test function 
$\pmb{\operatorname{C}}_{\text{b}}^{0}(\mathbb{R}) \! \ni \! 
h_{M_{0}} \colon \mathbb{N} \times \lbrace 1,2,\dotsc,K \rbrace 
\times \mathbb{R} \! \to \! [0,1]$, $(n,k,x) \! \mapsto \! h_{M_{0}}
(n,k,x) \! := \! h_{M_{0}}(x)$ be chosen such that: (i) $0 \! \leqslant 
\! h_{M_{0}}(x) \! \leqslant \! 1$, $x \! \in \! \mathbb{R}$; (ii) 
$h_{M_{0}}(x) \! = \! 1$, $x \! \in \! \mathbb{R} \setminus 
\mathfrak{D}_{M_{0}}$; and (iii) $h_{M_{0}}(x) \! = \! 0$, $x \! \in \! 
\mathfrak{D}_{M_{0}+1}$. Write the following decomposition \cite{a51}
\begin{equation*}
\iint_{\mathbb{R}^{2}} \mathfrak{g}_{f}(\tau,\xi) \, 
\md \mu_{2,m}^{\text{\tiny $\mathrm{EQ}$}}(\tau) \, \md 
\mu_{2,m}^{\text{\tiny $\mathrm{EQ}$}}(\xi) \! = \! I_{a,f} \! + \! I_{b,f} 
\! + \! I_{c,f},
\end{equation*}
where (suppressing $n$-, $k$-, and $z_{o}$-dependencies)
\begin{align*}
I_{a,f} :=& \iint_{\mathbb{R}^{2}} \mathfrak{g}_{f}(\tau,\xi)
(1 \! - \! h_{M_{0}}(\xi)) \, \md \mu_{2,m}^{\text{\tiny $\mathrm{EQ}$}}
(\tau) \, \md \mu_{2,m}^{\text{\tiny $\mathrm{EQ}$}}(\xi), \\
I_{b,f} :=& \iint_{\mathbb{R}^{2}} \mathfrak{g}_{f}
(\tau,\xi)h_{M_{0}}(\xi)(1 \! - \! h_{M_{0}}(\tau)) 
\, \md \mu_{2,m}^{\text{\tiny $\mathrm{EQ}$}}(\tau) 
\, \md \mu_{2,m}^{\text{\tiny $\mathrm{EQ}$}}(\xi), \\
I_{c,f} :=& \iint_{\mathbb{R}^{2}} \mathfrak{g}_{f}(\tau,\xi)h_{M_{0}}
(\tau)h_{M_{0}}(\xi) \, \md \mu_{2,m}^{\text{\tiny $\mathrm{EQ}$}}
(\tau) \, \md \mu_{2,m}^{\text{\tiny $\mathrm{EQ}$}}(\xi).
\end{align*}
One shows that
\begin{equation*}
\limsup_{m \to \infty} \lvert I_{a,f} \rvert \! \leqslant \! \lvert \lvert 
\mathfrak{g}_{f} \rvert \rvert_{\infty} \limsup_{m \to \infty} 
\int_{\mathbb{R} \setminus \mathfrak{D}_{M_{0}+1}} \, \md 
\mu_{2,m}^{\text{\tiny $\mathrm{EQ}$}}(\xi) \! \leqslant \! 
\hat{\varepsilon} \lvert \lvert \mathfrak{g}_{f} \rvert \rvert_{\infty},
\end{equation*}
where $\lvert \lvert \mathfrak{g}_{f} \rvert \rvert_{\infty} \! := 
\! \sup_{(\tau,\xi) \in \mathbb{R}^{2}} \vert \mathfrak{g}_{f}
(\tau,\xi) \vert$, and, similarly,
\begin{equation*}
\limsup_{m \to \infty} \vert I_{b,f} \vert \! \leqslant \! \hat{\varepsilon} 
\lvert \lvert \mathfrak{g}_{f} \rvert \rvert_{\infty}.
\end{equation*}
By an appropriate generalisation of the (one-variable case) Stone-Weierstrass 
Theorem, there exists a (symmetric) polynomial in two variables with bounded 
$n$-, $k$-, and $z_{o}$-dependent coefficients, that is, $p_{f}(\tau,\xi) \! = 
\! \sum_{i \geqslant i_{0}} \sum_{j \geqslant j_{0}} \gamma^{f}_{ij} \tau^{i} 
\xi^{j}$, $i_{0},j_{0} \! \in \! \mathbb{N}$, where $\gamma^{f}_{ij} \! = \! 
\gamma^{f}_{ij}(n,k,z_{o})$ with $\lvert \gamma^{f}_{ij} \rvert \! < \! +\infty$, 
such that, $\forall$ $\xi,\tau \! \in \! [-(M_{0} \! + \! 1),M_{0} \! + \! 1] 
\setminus \cup_{\underset{q \neq \mathfrak{s}-1}{q=1}}^{\mathfrak{s}} 
\operatorname{clos}(\mathscr{O}_{\frac{1}{M_{0}+1}}(\alpha_{p_{q}}))$, with, 
say, $M_{0} \! = \! K(1 \! + \! \max \lbrace \mathstrut \lvert \alpha_{p_{q}} 
\rvert, \, q \! = \! 1,\dotsc,\mathfrak{s} \! - \! 2,\mathfrak{s} \rbrace \! 
+ \! 3(\min \lbrace \mathstrut \lvert \alpha_{p_{i}} \! - \! \alpha_{p_{j}} 
\rvert, \, i \! \neq \! j \! \in \! \lbrace 1,\dotsc,\mathfrak{s} \! - \! 2,
\mathfrak{s} \rbrace \rbrace)^{-1})$, $\lvert \mathfrak{g}_{f}(\tau,\xi) 
\! - \! p_{f}(\tau,\xi) \rvert \! \leqslant \! \tilde{\varepsilon}$ $(= \! 
\tilde{\varepsilon}(n,k,z_{o}))$, whence $\lvert h_{M_{0}}(\tau)h_{M_{0}}
(\xi)(\mathfrak{g}_{f}(\tau,\xi) \! - \! p_{f}(\tau,\xi)) \rvert \! \leqslant \! 
\tilde{\varepsilon}$, $\tau,\xi \! \in \! \mathbb{R}$. Split $I_{c,f}$ further 
\cite{a51}, that is, write $I_{c,f} \! = \! I_{c,f}^{\alpha} \! + \! I_{c,f}^{\beta}$, 
where
\begin{align*}
I_{c,f}^{\alpha} &:= \iint_{\mathbb{R}^{2}}h_{M_{0}}(\xi)
h_{M_{0}}(\tau)(\mathfrak{g}_{f}(\tau,\xi) \! - \! p_{f}(\tau,\xi)) 
\, \md \mu_{2,m}^{\text{\tiny $\mathrm{EQ}$}}(\tau) 
\, \md \mu_{2,m}^{\text{\tiny $\mathrm{EQ}$}}(\xi), \\
I_{c,f}^{\beta} &:= \iint_{\mathbb{R}^{2}}h_{M_{0}}(\xi)h_{M_{0}}
(\tau)p_{f}(\tau,\xi) \, \md \mu_{2,m}^{\text{\tiny $\mathrm{EQ}$}}(\tau) 
\, \md \mu_{2,m}^{\text{\tiny $\mathrm{EQ}$}}(\xi).
\end{align*}
Now, proceeding \emph{verbatim} as on pp.~270--271 of \cite{a45}, one shows 
that $\vert I_{c,f}^{\alpha} \vert \! \leqslant \! \tilde{\varepsilon}$, whence
\begin{equation*}
\limsup_{m \to \infty} \vert I_{c,f}^{\alpha} \vert \! \leqslant \! 
\tilde{\varepsilon},
\end{equation*}
and
\begin{equation*}
I_{c,f}^{\beta} \! \leqslant \! \iint_{\mathbb{R}^{2}} 
\mathfrak{g}_{f}(\tau,\xi) \, \md \mu_{f}^{\text{\tiny $\mathrm{EQ}$}}
(\tau) \, \md \mu_{f}^{\text{\tiny $\mathrm{EQ}$}}(\xi) \! + \! 
\tilde{\varepsilon} \! + \! 2 \hat{\varepsilon} \lvert \lvert 
\mathfrak{g}_{f} \rvert \rvert_{\infty} \! + \! \hat{\varepsilon}^{2} 
\lvert \lvert \mathfrak{g}_{f} \rvert \rvert_{\infty}.
\end{equation*}
Hence, assembling the above-derived bounds for $I_{a,f}$, $I_{b,f}$, 
$I_{c,f}^{\alpha}$, and $I_{c,f}^{\beta}$, one arrives at, for $n \! \in 
\! \mathbb{N}$ and $k \! \in \! \lbrace 1,2,\dotsc,K \rbrace$ such that 
$\alpha_{p_{\mathfrak{s}}} \! := \!\alpha_{k} \! \neq \! \infty$, upon 
setting $\varepsilon_{f}^{\natural} \! := \! 2(\tilde{\varepsilon} \! + \! 
2 \hat{\varepsilon} \lvert \lvert \mathfrak{g}_{f} \rvert \rvert_{\infty})$ 
and neglecting the $\mathcal{O}(\hat{\varepsilon}^{2})$ term,
\begin{equation*}
\iint_{\mathbb{R}^{2}} \mathfrak{g}_{f}(\tau,\xi) 
\, \md \mu_{2,m}^{\text{\tiny $\mathrm{EQ}$}}(\tau) \, \md 
\mu_{2,m}^{\text{\tiny $\mathrm{EQ}$}}(\xi) \! - \! \iint_{\mathbb{R}^{2}} 
\mathfrak{g}_{f}(\tau,\xi) \, \md \mu_{f}^{\text{\tiny $\mathrm{EQ}$}}
(\tau) \, \md \mu_{f}^{\text{\tiny $\mathrm{EQ}$}}(\xi) \! \leqslant \! 
\varepsilon_{f}^{\natural};
\end{equation*}
thus, as $\varepsilon_{f}^{\natural} \! > \! 0$ is arbitrarily small (since 
$\tilde{\varepsilon} \! > \! 0$ and $\hat{\varepsilon} \! > \! 0$ are), it 
follows that
\begin{equation*}
\iint_{\mathbb{R}^{2}} \mathfrak{g}_{f}(\tau,\xi) 
\, \md \mu_{2,m}^{\text{\tiny $\mathrm{EQ}$}}(\tau) \, \md 
\mu_{2,m}^{\text{\tiny $\mathrm{EQ}$}}(\xi) \! \to \! 
\iint_{\mathbb{R}^{2}} \mathfrak{g}_{f}(\tau,\xi) 
\, \md \mu_{f}^{\text{\tiny $\mathrm{EQ}$}}(\tau) \, \md 
\mu_{f}^{\text{\tiny $\mathrm{EQ}$}}(\xi) \quad \text{as} \, \, 
m \! \to \! \infty.
\end{equation*}
Recalling that, for $(L_{f},m) \! \in \! \mathbb{R} \times \mathbb{N}$, 
$\mathfrak{g}_{f}(\tau,\xi) \! = \! \min \lbrace L_{f},\mathcal{K}^{f}_{
\widetilde{V}_{m}}(\tau,\xi) \rbrace$, it follows that, for $n \! \in \! 
\mathbb{N}$ and $k \! \in \! \lbrace 1,2,\dotsc,K \rbrace$ such that 
$\alpha_{p_{\mathfrak{s}}} \! := \! \alpha_{k} \! \neq \! \infty$,
\begin{equation*}
\liminf_{m \to \infty} \mathrm{I}_{\widetilde{V}}^{f}
[\mu_{2,m}^{\text{\tiny $\mathrm{EQ}$}}] \! \geqslant \! \iint_{\mathbb{R}^{2}} 
\min \left\lbrace L_{f},\mathcal{K}^{f}_{\widetilde{V}_{m}}(\tau,\xi) 
\right\rbrace \md \mu_{f}^{\text{\tiny $\mathrm{EQ}$}}(\tau) 
\, \md \mu_{f}^{\text{\tiny $\mathrm{EQ}$}}(\xi):
\end{equation*}
letting $L_{f} \! \uparrow \! +\infty$ and noting that, as $m \! \to \! 
\infty$, $\min \lbrace L_{f},\mathcal{K}^{f}_{\widetilde{V}_{m}}(\tau,\xi) 
\rbrace \! \to \! \mathcal{K}^{f}_{\widetilde{V}}(\tau,\xi)$, one arrives at, 
for $n \! \in \! \mathbb{N}$ and $k \! \in \! \lbrace 1,2,\dotsc,K \rbrace$ 
such that $\alpha_{p_{\mathfrak{s}}} \! := \! \alpha_{k} \! \neq \! \infty$, 
via the Monotone Convergence Theorem,
\begin{equation*}
\liminf_{m \to \infty} \mathrm{I}_{\widetilde{V}}^{f}
[\mu_{2,m}^{\text{\tiny $\mathrm{EQ}$}}] \! \geqslant \! 
\iint_{\mathbb{R}^{2}} \mathcal{K}^{f}_{\widetilde{V}}
(\tau,\xi) \, \md \mu_{f}^{\text{\tiny $\mathrm{EQ}$}}(\tau) \, \md 
\mu_{f}^{\text{\tiny $\mathrm{EQ}$}}(\xi) \! = \! \mathrm{I}_{
\widetilde{V}}^{f}[\mu_{f}^{\text{\tiny $\mathrm{EQ}$}}].
\end{equation*}
Since, for $n \! \in \! \mathbb{N}$ and $k \! \in \! \lbrace 1,2,\dotsc,K 
\rbrace$ such that $\alpha_{p_{\mathfrak{s}}} \! := \! \alpha_{k} \! \neq 
\! \infty$, there exists a weakly convergent subsequence of probability 
measures $\lbrace 
\mu_{2,m_{j}}^{\text{\tiny $\mathrm{EQ}$}} \rbrace_{j=1}^{\infty}$ 
$(\subset \mathscr{M}_{1}(\mathbb{R}))$ of $\lbrace 
\mu_{2,m}^{\text{\tiny $\mathrm{EQ}$}} \rbrace_{m=1}^{\infty}$ 
$(\subset \mathscr{M}_{1}(\mathbb{R}))$ with a weak-$\ast$ 
limit $\mu_{f}^{\text{\tiny $\mathrm{EQ}$}} \! \in \! \mathscr{M}_{1}
(\mathbb{R})$, that is, $\mu_{2,m_{j}}^{\text{\tiny $\mathrm{EQ}$}} \! 
\overset{\ast}{\to} \! \mu_{f}^{\text{\tiny $\mathrm{EQ}$}}$ as $j \! 
\to \! \infty$, upon recalling that $\mathrm{I}_{\widetilde{V}}^{f}
[\mu_{2,m}^{\text{\tiny $\mathrm{EQ}$}}] \! \leqslant \! 
E_{\widetilde{V}}^{f} \! + \! \tfrac{1}{m}$, $m \! \in \! \mathbb{N}$, 
it follows that, for $n \! \in \! \mathbb{N}$ and $k \! \in \! \lbrace 
1,2,\dotsc,K \rbrace$ such that $\alpha_{p_{\mathfrak{s}}} \! := 
\! \alpha_{k} \! \neq \! \infty$, $\mathrm{I}^{f}_{\widetilde{V}}
[\mu_{f}^{\text{\tiny $\mathrm{EQ}$}}] \! \leqslant \! E^{f}_{\widetilde{V}} 
\! = \! \inf \lbrace \mathstrut \mathrm{I}^{f}_{\widetilde{V}}
[\mu_{f}^{\text{\tiny $\mathrm{EQ}$}}]; \, \mu_{f}^{\text{\tiny 
$\mathrm{EQ}$}} \! \in \! \mathscr{M}_{1}(\mathbb{R}) \rbrace$ as 
$m \! \to \! \infty$; hence, for $n \! \in \! \mathbb{N}$ and $k \! \in 
\! \lbrace 1,2,\dotsc,K \rbrace$ such that $\alpha_{p_{\mathfrak{s}}} 
\! := \! \alpha_{k} \! \neq \! \infty$, it follows that there 
exists $\mu_{f}^{\text{\tiny $\mathrm{EQ}$}}$ $(= \! 
\mu_{f}^{\text{\tiny $\mathrm{EQ}$}}(n,k,z_{o}))$ $:= \! 
\mu^{f}_{\widetilde{V}} \! \in \! \mathscr{M}_{1}(\mathbb{R})$, 
the associated equilibrium measure, such that $\mathrm{I}^{f}_{
\widetilde{V}}[\mu^{f}_{\widetilde{V}}] \! = \! \inf \lbrace \mathstrut 
\mathrm{I}^{f}_{\widetilde{V}}[\mu_{f}^{\text{\tiny $\mathrm{EQ}$}}]; \, 
\mu_{f}^{\text{\tiny $\mathrm{EQ}$}} \! \in \! \mathscr{M}_{1}(\mathbb{R}) 
\rbrace$, that is, the infimum is attained, which establishes, for 
$n \! \in \! \mathbb{N}$ and $k \! \in \! \lbrace 1,2,\dotsc,K \rbrace$ 
such that $\alpha_{p_{\mathfrak{s}}} \! := \! \alpha_{k} \! \neq \! \infty$, 
the corresponding claim~(2) of the lemma. (The unicity of the associated 
equilibrium measure is addressed in Lemma~\ref{lem3.3} below.)

Finally, it remains to show that, for $n \! \in \! \mathbb{N}$ and $k \! 
\in \! \lbrace 1,2,\dotsc,K \rbrace$ such that $\alpha_{p_{\mathfrak{s}}} 
\! := \! \alpha_{k} \! \neq \! \infty$, $(\lbrace \mathstrut x \! \in \! 
\mathbb{R}; \, \varpi (x) \! > \! 0 \rbrace \supset)$ $J_{f} \! := \! 
\supp (\mu^{f}_{\widetilde{V}}) \subset \overline{\mathbb{R}} \setminus 
\lbrace \alpha_{1},\alpha_{2},\dotsc,\alpha_{K} \rbrace$ is compact. (The 
following argument is valid for any $\mu_{2}^{\text{\tiny $\mathrm{EQ}$}}$ 
$(= \! \mu_{2}^{\text{\tiny $\mathrm{EQ}$}}(n,k,z_{o}))$ $\in \! 
\mathscr{M}_{1}(\mathbb{R})$ achieving the above infimum; in particular, 
for $\mu_{2}^{\text{\tiny $\mathrm{EQ}$}} \! = \! \mu^{f}_{\widetilde{V}}$.) 
For $n \! \in \! \mathbb{N}$ and $k \! \in \! \lbrace 1,2,\dotsc,K \rbrace$ 
such that $\alpha_{p_{\mathfrak{s}}} \! := \! \alpha_{k} \! \neq \! \infty$, 
let $\mathscr{M}_{1}(\mathbb{R}) \! \ni \! \mu_{\varpi,f}^{\text{\tiny 
$\mathrm{EQ}$}}$ $(= \! \mu_{\varpi,f}^{\text{\tiny $\mathrm{EQ}$}}
(n,k,z_{o}))$ be such that $\mathrm{I}^{f}_{\widetilde{V}}
[\mu_{\varpi,f}^{\text{\tiny $\mathrm{EQ}$}}] \! = \! E^{f}_{\widetilde{V}}$, 
and let $D_{f}$ $(= \! D_{f}(n,k,z_{o}))$ be any proper measurable subset of 
$\mathbb{R}$ for which $\mu_{\varpi,f}^{\text{\tiny $\mathrm{EQ}$}}(D_{f}) 
\! = \! \int_{D_{f}} \md \mu_{\varpi,f}^{\text{\tiny $\mathrm{EQ}$}}(\lambda) 
\! > \! 0$. For $n \! \in \! \mathbb{N}$ and $k \! \in \! \lbrace 1,2,\dotsc,K 
\rbrace$ such that $\alpha_{p_{\mathfrak{s}}} \! := \! \alpha_{k} \! \neq \! 
\infty$, set \cite{a56}
\begin{equation*}
\mu_{\varpi,f}^{\varepsilon}(z) \! := \! \left(1 \! + \! \varepsilon 
\mu_{\varpi,f}^{\text{\tiny $\mathrm{EQ}$}}(D_{f}) \right)^{-1} \left(
\mu_{\varpi,f}^{\text{\tiny $\mathrm{EQ}$}}(z) \! + \! \varepsilon 
(\mu_{\varpi,f}^{\text{\tiny $\mathrm{EQ}$}} \! \! \upharpoonright_{D_{f}})
(z) \right), \quad \varepsilon \! \in \! (-1,1),
\end{equation*}
where $\mu_{\varpi,f}^{\text{\tiny $\mathrm{EQ}$}} \! \! 
\upharpoonright_{D_{f}}$ denotes the restriction of 
$\mu_{\varpi,f}^{\text{\tiny $\mathrm{EQ}$}}$ to $D_{f}$ (note that 
$\int_{\mathbb{R}} \md \mu_{\varpi,f}^{\varepsilon}(\lambda) 
\! = \! 1$, $\varepsilon \! \in \! (-1,1)$; in particular, for $\varepsilon 
\! \in \! (-1,1)$, $\mu_{\varpi,f}^{\varepsilon} \! \in \! \mathscr{M}_{1}
(\mathbb{R}))$. Using the symmetry $\mathcal{K}^{f}_{\widetilde{V}}
(\tau,\xi) \! = \! \mathcal{K}^{f}_{\widetilde{V}}(\xi,\tau)$, one shows that, 
for $n \! \in \! \mathbb{N}$ and $k \! \in \! \lbrace 1,2,\dotsc,K \rbrace$ 
such that $\alpha_{p_{\mathfrak{s}}} \! := \! \alpha_{k} \! \neq \! \infty$,
\begin{align*}
\mathrm{I}^{f}_{\widetilde{V}}[\mu_{\varpi,f}^{\varepsilon}] =& \, 
\iint_{\mathbb{R}^{2}} \mathcal{K}^{f}_{\widetilde{V}}
(\xi,\tau) \, \md \mu^{\varepsilon}_{\varpi,f}(\xi) \, 
\md \mu^{\varepsilon}_{\varpi,f}(\tau) \! = \! \left(1 \! + \! \varepsilon 
\mu_{\varpi,f}^{\text{\tiny $\mathrm{EQ}$}}(D_{f}) \right)^{-2} \left(
\mathrm{I}^{f}_{\widetilde{V}}[\mu_{\varpi,f}^{\text{\tiny $\mathrm{EQ}$}}] 
\! + \! 2 \varepsilon \iint_{\mathbb{R}^{2}} \mathcal{K}^{f}_{\widetilde{V}}
(\xi,\tau) \, \md \mu_{\varpi,f}^{\text{\tiny $\mathrm{EQ}$}}(\xi) \right. \\
\times&\left. \, \md (\mu_{\varpi,f}^{\text{\tiny $\mathrm{EQ}$}} \! \! 
\upharpoonright_{D_{f}})(\tau) \! + \! \varepsilon^{2} \iint_{\mathbb{R}^{2}} 
\mathcal{K}^{f}_{\widetilde{V}}(\xi,\tau) \, \md 
(\mu_{\varpi,f}^{\text{\tiny $\mathrm{EQ}$}} \! \! \upharpoonright_{D_{f}})
(\xi) \, \md (\mu_{\varpi,f}^{\text{\tiny $\mathrm{EQ}$}} \! \! 
\upharpoonright_{D_{f}})(\tau) \right), \quad \varepsilon \! \in \! (-1,1).
\end{align*}
By the minimal property of $\mu_{\varpi,f}^{\text{\tiny $\mathrm{EQ}$}} 
\! \in \! \mathscr{M}_{1}(\mathbb{R})$, it follows that \cite{a56}, for 
$n \! \in \! \mathbb{N}$ and $k \! \in \! \lbrace 1,2,\dotsc,K \rbrace$ such 
that $\alpha_{p_{\mathfrak{s}}} \! := \! \alpha_{k} \! \neq \! \infty$, 
$\partial_{\varepsilon} \mathrm{I}^{f}_{\widetilde{V}}[\mu_{\varpi,f}^{
\varepsilon}] \! = \! 0$, whence
\begin{equation}
\iint_{\mathbb{R}^{2}} \left(\mathcal{K}^{f}_{\widetilde{V}}
(\xi,\tau) \! - \! \mathrm{I}^{f}_{\widetilde{V}}
[\mu_{\varpi,f}^{\text{\tiny $\mathrm{EQ}$}}] \right) \md 
\mu_{\varpi,f}^{\text{\tiny $\mathrm{EQ}$}}(\xi) \, \md 
(\mu_{\varpi,f}^{\text{\tiny $\mathrm{EQ}$}} \! \! \upharpoonright_{D_{f}})
(\tau) \! = \! 0; \label{eqKvinf5}
\end{equation}
but, recalling that, for $n \! \in \! \mathbb{N}$ and $k \! \in \! 
\lbrace 1,2,\dotsc,K \rbrace$ such that $\alpha_{p_{\mathfrak{s}}} 
\! := \!\alpha_{k} \! \neq \! \infty$, $\mathcal{K}^{f}_{\widetilde{V}}
(\xi,\tau) \! \geqslant \! (\hat{\psi}^{f}_{\widetilde{V}}(\xi) \! + \! 
\hat{\psi}^{f}_{\widetilde{V}}(\tau))/2$, where $\hat{\psi}^{f}_{
\widetilde{V}}(z)$ is defined by Equation~\eqref{eqKvinf4}, it follows 
{}from the minimisation condition~\eqref{eqKvinf5} that
\begin{gather}
\iint_{\mathbb{R}^{2}} \mathrm{I}^{f}_{\widetilde{V}}
[\mu_{\varpi,f}^{\text{\tiny $\mathrm{EQ}$}}] \, \md 
\mu_{\varpi,f}^{\text{\tiny $\mathrm{EQ}$}}(\xi) \, \md 
(\mu_{\varpi,f}^{\text{\tiny $\mathrm{EQ}$}} \! \! \upharpoonright_{D_{f}})
(\tau) \! \geqslant \! \dfrac{1}{2} \iint_{\mathbb{R}^{2}}
(\hat{\psi}^{f}_{\widetilde{V}}(\xi) \! + \! \hat{\psi}^{f}_{\widetilde{V}}
(\tau)) \, \md \mu_{\varpi,f}^{\text{\tiny $\mathrm{EQ}$}}(\xi) \, \md 
(\mu_{\varpi,f}^{\text{\tiny $\mathrm{EQ}$}} \! \! \upharpoonright_{D_{f}})
(\tau) \, \, \Rightarrow \nonumber \\
0 \! \geqslant \! \int_{\mathbb{R}} \left(\hat{\psi}^{f}_{\widetilde{V}}
(\tau) \! + \! \int_{\mathbb{R}} \hat{\psi}^{f}_{\widetilde{V}}(\xi) \, 
\md \mu_{\varpi,f}^{\text{\tiny $\mathrm{EQ}$}}(\xi) \! - \! 2 
\mathrm{I}^{f}_{\widetilde{V}}[\mu_{\varpi,f}^{\text{\tiny $\mathrm{EQ}$}}] 
\right) \md (\mu_{\varpi,f}^{\text{\tiny $\mathrm{EQ}$}} 
\! \! \upharpoonright_{D_{f}})(\tau). \label{eqKvinf6}
\end{gather}
But, for $n \! \in \! \mathbb{N}$ and $k \! \in \! \lbrace 1,2,\dotsc,K 
\rbrace$ such that $\alpha_{p_{\mathfrak{s}}} \! := \! \alpha_{k} \! \neq 
\! \infty$, via the growth conditions on $\widetilde{V} \colon \overline{
\mathbb{R}} \setminus \lbrace \alpha_{1},\alpha_{2},\dotsc,\alpha_{K} 
\rbrace \! \to \! \mathbb{R}$ satisfying 
conditions~\eqref{eq20}--\eqref{eq22}, that is, there exists $\mathfrak{c}_{
\infty} \! > \! 0$ and bounded such that, for $x \! \in \! \mathscr{O}_{
\infty}$, $\widetilde{V}(x) \! \geqslant \! (1 \! + \! \mathfrak{c}_{\infty}) 
\ln (1 \! + \! x^{2})$, and, for $q \! = \! 1,\dotsc,\mathfrak{s} \! - \! 2,
\mathfrak{s}$, there exist $\mathfrak{c}_{q} \! > \! 0$ and bounded 
such that, for $x \! \in \! \mathscr{O}_{\tilde{\delta}_{q}}(\alpha_{p_{q}})$, 
$\widetilde{V}(x) \! \geqslant \! (1 \! + \! \mathfrak{c}_{q}) \ln (1 \! 
+ \! (x \! - \! \alpha_{p_{q}})^{-2})$, it follows that, for $n \! \in \! 
\mathbb{N}$ and $k \! \in \! \lbrace 1,2,\dotsc,K \rbrace$ such that 
$\alpha_{p_{\mathfrak{s}}} \! := \! \alpha_{k} \! \neq \! \infty$, there 
exists (some) $T_{M_{f}}$ $(= \! T_{M_{f}}(n,k,z_{o}))$ $> \! 1$ (e.g., 
$T_{M_{f}} \! = \! K(1 \! + \! \max \lbrace \lvert \alpha_{p_{q}} \rvert, 
\, q \! = \! 1,\dotsc,\mathfrak{s} \! - \! 2,\mathfrak{s} \rbrace \! + \! 
3(\min \lbrace \mathstrut \lvert \alpha_{p_{i}} \! - \! \alpha_{p_{j}} 
\rvert, \, i \! \neq \! j \! \in \! \lbrace 1,\dotsc,\mathfrak{s} \! - \! 2,
\mathfrak{s} \rbrace \rbrace)^{-1}))$, with $\mathscr{O}_{\frac{1}{T_{M_{f}}}}
(\alpha_{p_{i}}) \cap \mathscr{O}_{\frac{1}{T_{M_{f}}}}(\alpha_{p_{j}}) \! = 
\! \varnothing$, $i \! \neq \! j \! \in \! \lbrace 1,\dotsc,\mathfrak{s} \! 
- \! 2,\mathfrak{s} \rbrace$, such that, for $\tau \! \in \! \lbrace \lvert 
x \rvert \! \geqslant \! T_{M_{f}} \rbrace \cup \cup_{\underset{q \neq 
\mathfrak{s}-1}{q=1}}^{\mathfrak{s}} \operatorname{clos}
(\mathscr{O}_{\frac{1}{T_{M_{f}}}}(\alpha_{p_{q}}))$,
\begin{equation}
\hat{\psi}^{f}_{\widetilde{V}}(\tau) \! + \! \int_{\mathbb{R}} 
\hat{\psi}^{f}_{\widetilde{V}}(\xi) \, \md 
\mu_{\varpi,f}^{\text{\tiny $\mathrm{EQ}$}}(\xi) \! - \! 2 \mathrm{I}^{f}_{
\widetilde{V}}[\mu_{\varpi,f}^{\text{\tiny $\mathrm{EQ}$}}] \! \geqslant \! 
1; \label{eqKvinf7}
\end{equation}
hence, for $D_{f} \subset \lbrace \lvert x \rvert \! \geqslant \! T_{M_{f}} 
\rbrace \cup \cup_{\underset{q \neq \mathfrak{s}-1}{q=1}}^{\mathfrak{s}} 
\operatorname{clos}(\mathscr{O}_{\frac{1}{T_{M_{f}}}}(\alpha_{p_{q}}))$, 
with $T_{M_{f}} \! > \! 1$ (chosen as above), it follows that, for $n \! 
\in \! \mathbb{N}$ and $k \! \in \! \lbrace 1,2,\dotsc,K \rbrace$ such 
that $\alpha_{p_{\mathfrak{s}}} \! := \! \alpha_{k} \! \neq \! \infty$, via 
inequalities~\eqref{eqKvinf6} and~\eqref{eqKvinf7},
\begin{equation*}
0 \! \geqslant \! \int_{\mathbb{R}} \left(\hat{\psi}^{f}_{
\widetilde{V}}(\tau) \! + \! \int_{\mathbb{R}} \hat{\psi}^{f}_{
\widetilde{V}}(\xi) \, \md \mu_{\varpi,f}^{\text{\tiny $\mathrm{EQ}$}}
(\xi) \! - \! 2 \mathrm{I}^{f}_{\widetilde{V}}
[\mu_{\varpi,f}^{\text{\tiny $\mathrm{EQ}$}}] \right) \md 
(\mu_{\varpi,f}^{\text{\tiny $\mathrm{EQ}$}} \! \! \upharpoonright_{D_{f}})
(\tau) \! \geqslant \! 1,
\end{equation*}
which is a contradiction; hence, for $n \! \in \! \mathbb{N}$ and $k \! 
\in \! \lbrace 1,2,\dotsc,K \rbrace$ such that $\alpha_{p_{\mathfrak{s}}} 
\! := \! \alpha_{k} \! \neq \! \infty$, $\supp 
(\mu_{\varpi,f}^{\text{\tiny $\mathrm{EQ}$}}) \subseteq \mathbb{R} 
\setminus (\lbrace \lvert x \rvert \! \geqslant \! T_{M_{f}} \rbrace 
\cup \cup_{\underset{q \neq \mathfrak{s}-1}{q=1}}^{\mathfrak{s}} 
\operatorname{clos}(\mathscr{O}_{\frac{1}{T_{M_{f}}}}(\alpha_{p_{q}})))$; 
in particular, $J_{f} \! := \! \supp (\mu^{f}_{\widetilde{V}}) \subseteq 
\mathbb{R} \setminus (\lbrace \lvert x \rvert \! \geqslant \! T_{M_{f}} 
\rbrace \cup \cup_{\underset{q \neq \mathfrak{s}-1}{q=1}}^{
\mathfrak{s}} \operatorname{clos}(\mathscr{O}_{\frac{1}{T_{M_{f}}}}
(\alpha_{p_{q}})))$, which establishes the compactness of the support of 
the associated equilibrium measure $\mu^{f}_{\widetilde{V}} \! \in \! 
\mathscr{M}_{1}(\mathbb{R})$ (it necessarily follows that $\supp 
(\mu^{f}_{\widetilde{V}}) \cap \lbrace \alpha_{1},\alpha_{2},\dotsc,
\alpha_{K} \rbrace \! = \! \varnothing)$. Since $\widetilde{V} \colon 
\overline{\mathbb{R}} \setminus \lbrace \alpha_{1},\alpha_{2},\dotsc,
\alpha_{K} \rbrace \! \to \! \mathbb{R}$ is real analytic and $\supp 
(\mu^{f}_{\widetilde{V}}) \cap \lbrace \alpha_{1},\alpha_{2},\dotsc,
\alpha_{K} \rbrace \! = \! \varnothing$, in which case, for $\supp 
(\mu^{f}_{\widetilde{V}}) \! \ni \! x$, $\inf_{x \in J_{f}} \widetilde{V}(x) 
\! \leqslant \! \widetilde{V}(x) \! \leqslant \! \sup_{x \in J_{f}} \widetilde{V}
(x)$, and $-\infty \! < \! \int_{J_{f}} \widetilde{V}(\lambda) \, \md 
\mu^{f}_{\widetilde{V}}(\lambda) \! < \! +\infty$, it follows that, for 
$n \! \in \! \mathbb{N}$ and $k \! \in \! \lbrace 1,2,\dotsc,K \rbrace$ 
such that $\alpha_{p_{\mathfrak{s}}} \! := \! \alpha_{k} \! \neq \! \infty$, 
both the weighted logarithmic energy and logarithmic energy of 
$\mu^{f}_{\widetilde{V}} \! \in \! \mathscr{M}_{1}(\mathbb{R})$, 
respectively, are bounded, that is,
\begin{align*}
-\infty \! < \! \mathrm{I}^{f}_{\widetilde{V}}[\mu^{f}_{\widetilde{V}}] 
=& \, \iint_{J_{f} \times J_{f}} \ln \left(\vert \xi \! - \! \tau \vert^{
\frac{\varkappa_{nk \tilde{k}_{\mathfrak{s}-1}}^{\infty}+1}{n}} \left(
\dfrac{\vert \xi \! - \! \tau \vert}{\vert \xi \! - \! \alpha_{k} \vert \vert 
\tau \! - \! \alpha_{k} \vert} \right)^{\frac{\varkappa_{nk}-1}{n}} 
\prod_{q=1}^{\mathfrak{s}-2} \left(\dfrac{\vert \xi \! - \! \tau \vert}{\vert 
\xi \! - \! \alpha_{p_{q}} \vert \vert \tau \! - \! \alpha_{p_{q}} \vert} 
\right)^{\frac{\varkappa_{nk \tilde{k}_{q}}}{n}} \varpi (\xi) \varpi (\tau) 
\right)^{-1} \\
\times& \, \md \mu^{f}_{\widetilde{V}}(\xi) \, \md \mu^{f}_{\widetilde{V}}
(\tau) \! < \! +\infty, \\
-\infty <& \, \iint_{J_{f} \times J_{f}} \ln \left(\vert \xi \! - \! \tau 
\vert^{\frac{\varkappa_{nk \tilde{k}_{\mathfrak{s}-1}}^{\infty}+1}{n}} \left(
\dfrac{\vert \xi \! - \! \tau \vert}{\vert \xi \! - \! \alpha_{k} \vert \vert 
\tau \! - \! \alpha_{k} \vert} \right)^{\frac{\varkappa_{nk}-1}{n}} 
\prod_{q=1}^{\mathfrak{s}-2} \left(\dfrac{\vert \xi \! - \! \tau \vert}{\vert 
\xi \! - \! \alpha_{p_{q}} \vert \vert \tau \! - \! \alpha_{p_{q}} \vert} 
\right)^{\frac{\varkappa_{nk \tilde{k}_{q}}}{n}} \right)^{-1} \md 
\mu^{f}_{\widetilde{V}}(\xi) \, \md \mu^{f}_{\widetilde{V}}(\tau) \! < \! 
+\infty,
\end{align*}
which establishes, for $n \! \in \! \mathbb{N}$ and $k \! \in \! \lbrace 
1,2,\dotsc,K \rbrace$ such that $\alpha_{p_{\mathfrak{s}}} \! := \! 
\alpha_{k} \! \neq \! \infty$, the corresponding claim~(3) of the lemma. 
Furthermore, it is a consequence of the facts established that, for $n \! 
\in \! \mathbb{N}$ and $k \! \in \! \lbrace 1,2,\dotsc,K \rbrace$ such 
that $\alpha_{p_{\mathfrak{s}}} \! := \! \alpha_{k} \! \neq \! \infty$, 
$J_{f} \! := \! \supp (\mu^{f}_{\widetilde{V}})$, which is a proper compact 
subset of $\overline{\mathbb{R}} \setminus \lbrace \alpha_{1},\alpha_{2},
\dotsc,\alpha_{K} \rbrace$, has positive---weighted logarithmic---capacity, 
that is, $\operatorname{cap}(J_{f}) \! = \! \exp (-E^{f}_{\widetilde{V}}) 
\! > \! 0$.

$\pmb{(2)}$ The proof of this case, that is, $n \! \in \! \mathbb{N}$ and $k 
\! \in \! \lbrace 1,2,\dotsc,K \rbrace$ such that $\alpha_{p_{\mathfrak{s}}} 
\! := \! \alpha_{k} \! = \! \infty$, is virtually identical to the proof presented 
in case $\pmb{(1)}$. One mimics, \emph{verbatim}, the scheme of the calculations 
presented in $\pmb{(1)}$ above in order to arrive at the corresponding 
claims~(1)--(3) of the lemma; in order to do so, however, one needs the 
analogues of definitions~\eqref{eqKvinf1} and~\eqref{eqKvinf4}, which, 
respectively, read: for $n \! \in \! \mathbb{N}$ and $k \! \in \! \lbrace 
1,2,\dotsc,K \rbrace$ such that $\alpha_{p_{\mathfrak{s}}} \! := \! 
\alpha_{k} \! = \! \infty$, with $\mathrm{I}^{\infty}_{\widetilde{V}}
[\mu^{\text{\tiny $\mathrm{EQ}$}}_{1}] \! = \! \iint_{\mathbb{R}^{2}} 
\mathcal{K}_{\widetilde{V}}^{\infty}(\xi,\tau) \, \md 
\mu^{\text{\tiny $\mathrm{EQ}$}}_{1}(\xi) \, \md 
\mu^{\text{\tiny $\mathrm{EQ}$}}_{1}(\tau)$,
\begin{equation} \label{eqKvinf8} 
\begin{split}
\mathcal{K}_{\widetilde{V}}^{\infty}(\xi,\tau) =& \, \mathcal{K}_{
\widetilde{V}}^{\infty}(\tau,\xi) \! = \! \ln \left(\vert \xi \! - \! 
\tau \vert^{\frac{\varkappa_{nk}}{n}} \prod_{q=1}^{\mathfrak{s}-1} 
\left(\dfrac{\vert \xi \! - \! \tau \vert}{\vert \xi \! - \! \alpha_{p_{q}} 
\vert \vert \tau \! - \! \alpha_{p_{q}} \vert} \right)^{\frac{\varkappa_{nk 
\tilde{k}_{q}}}{n}} \varpi (\xi) \varpi (\tau) \right)^{-1} \\
=& \, \dfrac{\varkappa_{nk}}{n} \ln \left(\dfrac{1}{\lvert \xi \! - \! \tau 
\rvert} \right) \! + \! \sum_{q=1}^{\mathfrak{s}-1} \dfrac{\varkappa_{nk 
\tilde{k}_{q}}}{n} \ln \left(\left\vert \dfrac{1}{\xi \! - \! \alpha_{p_{q}}} 
\! - \! \dfrac{1}{\tau \! - \! \alpha_{p_{q}}} \right\vert^{-1} \right) \! + 
\! \dfrac{1}{2} \widetilde{V}(\xi) \! + \! \dfrac{1}{2} \widetilde{V}(\tau) \\
=& \, \dfrac{1}{2} \widetilde{V}(\xi) \! + \! \dfrac{1}{2} \widetilde{V}(\tau) 
\! + \! \ln \left(\dfrac{\lvert \xi \! - \! \tau \rvert^{K}}{\prod_{q=1}^{
\mathfrak{s}-1}(\lvert \xi \! - \! \alpha_{p_{q}} \rvert \lvert \tau \! - \! 
\alpha_{p_{q}} \rvert)^{\gamma_{i(q)_{k_{q}}}}} \right)^{-1} \\
+& \, 
\begin{cases} 
\dfrac{1}{n} \ln \left(\dfrac{\lvert \xi \! - \! \tau \rvert^{k}}{\lvert 
\xi \! - \! \tau \rvert^{K}} \prod_{j=1}^{\mathfrak{s}-1}(\lvert \xi \! 
- \! \alpha_{p_{j}} \rvert \lvert \tau \! - \! \alpha_{p_{j}} \rvert)^{
\gamma_{i(j)_{k_{j}}}} \right)^{-1}, &\text{$\mathfrak{J}_{q}(k) \! = \! 
\varnothing, \quad q \! \in \! \lbrace 1,2,\dotsc,\mathfrak{s} \! - \! 1 
\rbrace$,} \\
\dfrac{1}{n} \ln \left(\dfrac{\lvert \xi \! - \! \tau \rvert^{k} 
\prod_{j=1}^{\mathfrak{s}-1}(\lvert \xi \! - \! \alpha_{p_{j}} \rvert 
\lvert \tau \! - \! \alpha_{p_{j}} \rvert)^{\gamma_{i(j)_{k_{j}}}}}{\lvert 
\xi \! - \! \tau \rvert^{K} \prod_{j=1}^{\mathfrak{s}-1}(\lvert \xi \! 
- \! \alpha_{p_{j}} \rvert \lvert \tau \! - \! \alpha_{p_{j}} \rvert)^{
\varrho_{\widetilde{m}_{j}(k)}}} \right)^{-1}, &\text{$\mathfrak{J}_{q}(k) 
\! \neq \! \varnothing, \quad q \! \in \! \lbrace 1,2,\dotsc,\mathfrak{s} 
\! - \! 1 \rbrace$,}
\end{cases}
\end{split}
\end{equation}
and
\begin{equation}
\hat{\psi}_{\widetilde{V}}^{\infty}(z) \! := \! \widetilde{V}(z) \! - \! 
\dfrac{\varkappa_{nk}}{n} \ln (1 \! + \! z^{2}) \! - \! \sum_{q=1}^{
\mathfrak{s}-1} \dfrac{\varkappa_{nk \tilde{k}_{q}}}{n} \ln 
(1 \! + \! (z \! - \! \alpha_{p_{q}})^{-2}). \label{eqKvinf9}
\end{equation}
This concludes the proof. \hfill $\qed$
\begin{eeeee} \label{rem1.3.4} 
\textsl{For $n \! \in \! \mathbb{N}$ and $k \! \in \! \lbrace 1,2,
\dotsc,K \rbrace$ such that $\alpha_{p_{\mathfrak{s}}} \! := \! 
\alpha_{k} \! = \! \infty$ or $\alpha_{p_{\mathfrak{s}}} \! := \! 
\alpha_{k} \! \neq \! \infty$, a combined family of $K$ energy 
minimisation problems has been presented in Lemma~\ref{lem3.1}. 
As the principal focus of this monograph is asymptotics, in the 
double-scaling limit $\mathscr{N},n \! \to \! \infty$ such that 
$z_{o} \! = \! 1 \! + \! o(1)$, of {\rm MPC ORFs} and related quantities, 
it is instructive to consider the `asymptotic structure' of this family 
of $K$ energy minimisation problems. In the double-scaling limit 
$\mathscr{N},n \! \to \! \infty$ such that $z_{o} \! = \! 1 \! + \! o(1)$, 
the family of $K$ energy minimisation problems `stabilizes', in the 
sense that one may write the associated energy functionals as `small' 
perturbations of certain `core' energy functionals: {\rm (i)} for $n \! 
\in \! \mathbb{N}$ and $k \! \in \! \lbrace 1,2,\dotsc,K \rbrace$ such 
that $\alpha_{p_{\mathfrak{s}}} \! := \! \alpha_{k} \! = \! \infty$,
\begin{equation*}
\mathrm{I}^{\infty}_{\widetilde{V}}[\mu^{\text{\tiny {\rm EQ}}}] 
\underset{\substack{\mathscr{N},n \to \infty\\z_{o}=1+o(1)}}{=} 
\mathfrak{X}^{\infty}_{\widetilde{V}}[\mu^{\text{\tiny {\rm EQ}}}] 
\! + \! \dfrac{1}{n} \hat{\mathfrak{Z}}^{\infty}_{\widetilde{V}}
[\mu^{\text{\tiny {\rm EQ}}}] \! + \! o(\mathfrak{e}),
\end{equation*}
with the associated `core' energy functional, $\mathfrak{X}^{\infty}_{
\widetilde{V}}[\mu^{\text{\tiny {\rm EQ}}}]$, given by
\begin{equation} \label{stochinf} 
\mathfrak{X}^{\infty}_{\widetilde{V}}[\mu^{\text{\tiny {\rm EQ}}}] 
\! = \! \iint_{\mathbb{R}^{2}} \ln \left(\dfrac{\lvert \xi \! - \! \tau 
\rvert^{K}}{\prod_{q=1}^{\mathfrak{s}-1}(\lvert \xi \! - \! \alpha_{p_{q}} 
\rvert \lvert \tau \! - \! \alpha_{p_{q}} \rvert)^{\gamma_{i(q)_{k_{q}}}}} 
\right)^{-1} \md \mu^{\text{\tiny {\rm EQ}}}(\xi) \, \md 
\mu^{\text{\tiny {\rm EQ}}}(\tau) \! + \! \int_{\mathbb{R}}V(\xi) 
\, \md \mu^{\text{\tiny {\rm EQ}}}(\xi),
\end{equation}
where, for $q \! \in \! \lbrace 1,2,\dotsc,\mathfrak{s} \! - \! 1 \rbrace$,
\begin{equation*}
\hat{\mathfrak{Z}}^{\infty}_{\widetilde{V}}[\mu^{\text{\tiny {\rm EQ}}}] 
\! = \! 
\begin{cases}
\iint_{\mathbb{R}^{2}} \ln \left(\dfrac{\lvert \xi \! - \! \tau 
\rvert^{k}}{\lvert \xi \! - \! \tau \rvert^{K}} \prod_{j=1}^{\mathfrak{s}-1}
(\lvert \xi \! - \! \alpha_{p_{j}} \rvert \lvert \tau \! - \! \alpha_{p_{j}} 
\rvert)^{\gamma_{i(j)_{k_{j}}}} \right)^{-1} \md 
\mu^{\text{\tiny {\rm EQ}}}(\xi) \, \md \mu^{\text{\tiny {\rm EQ}}}
(\tau), &\text{$\mathfrak{J}_{q}(k) \! = \! \varnothing$,} \\
\iint_{\mathbb{R}^{2}} \ln \left(\dfrac{\lvert \xi \! - \! \tau 
\rvert^{k} \prod_{j=1}^{\mathfrak{s}-1}(\lvert \xi \! - \! \alpha_{p_{j}} 
\rvert \lvert \tau \! - \! \alpha_{p_{j}} \rvert)^{\gamma_{i(j)_{k_{j}}}}}{
\lvert \xi \! - \! \tau \rvert^{K} \prod_{j=1}^{\mathfrak{s}-1}(\lvert 
\xi \! - \! \alpha_{p_{j}} \rvert \lvert \tau \! - \! \alpha_{p_{j}} 
\rvert)^{\varrho_{\widetilde{m}_{j}(k)}}} \right)^{-1} \md 
\mu^{\text{\tiny {\rm EQ}}}(\xi) \, \md \mu^{\text{\tiny {\rm EQ}}}
(\tau), &\text{$\mathfrak{J}_{q}(k) \! \neq \! \varnothing$,}
\end{cases}
\end{equation*}
and $\mathfrak{e} \! = \! \int_{\mathbb{R}}V(\xi) \, \md 
\mu^{\text{\tiny {\rm EQ}}}(\xi)$$;$\footnote{Via 
conditions~\eqref{eq20}--\eqref{eq22}, and the fact that, for 
$\mu^{\text{\tiny {\rm EQ}}} \! \in \! \mathscr{M}_{1}(\mathbb{R})$, 
$\supp (\mu^{\text{\tiny {\rm EQ}}})$ is a proper compact subset of 
$\overline{\mathbb{R}} \setminus \lbrace \alpha_{1},\alpha_{2},\dotsc,
\alpha_{K} \rbrace$, $\exists$ $\mathrm{M} \! > \! 0$ (and $\mathcal{O}
(1))$ such that $\lvert \mathfrak{e} \rvert \! \leqslant \! \mathrm{M}$.} 
and {\rm (ii)} for $n \! \in \! \mathbb{N}$ and $k \! \in \! \lbrace 1,2,
\dotsc,K \rbrace$ such that $\alpha_{p_{\mathfrak{s}}} \! := \! \alpha_{k} 
\! \neq \! \infty$,
\begin{equation*}
\mathrm{I}^{f}_{\widetilde{V}}[\mu^{\text{\tiny {\rm EQ}}}] 
\underset{\substack{\mathscr{N},n \to \infty\\z_{o}=1+o(1)}}{=} 
\mathfrak{X}^{f}_{\widetilde{V}}[\mu^{\text{\tiny {\rm EQ}}}] 
\! + \! \dfrac{1}{n} \tilde{\mathfrak{Z}}^{f}_{\widetilde{V}}
[\mu^{\text{\tiny {\rm EQ}}}] \! + \! o(\mathfrak{e}),
\end{equation*}
with the associated `core' energy functional, $\mathfrak{X}^{f}_{
\widetilde{V}}[\mu^{\text{\tiny {\rm EQ}}}]$, given by
\begin{equation} \label{stochfin} 
\mathfrak{X}^{f}_{\widetilde{V}}[\mu^{\text{\tiny {\rm EQ}}}] \! 
= \! \iint_{\mathbb{R}^{2}} \ln \left(\dfrac{\lvert \xi \! - \! \tau 
\rvert^{K}}{\prod_{q=1}^{\mathfrak{s}-2}(\lvert \xi \! - \! \alpha_{p_{q}} 
\rvert \lvert \tau \! - \! \alpha_{p_{q}} \rvert)^{\gamma_{i(q)_{k_{q}}}}
(\lvert \xi \! - \! \alpha_{k} \rvert \lvert \tau \! - \! \alpha_{k} 
\rvert)^{\gamma_{k}}} \right)^{-1} \md \mu^{\text{\tiny {\rm EQ}}}(\xi) 
\, \md \mu^{\text{\tiny {\rm EQ}}}(\tau) \! + \! \int_{\mathbb{R}}
V(\xi) \, \md \mu^{\text{\tiny {\rm EQ}}}(\xi),
\end{equation}
where, for $q \! \in \! \lbrace 1,2,\dotsc,\mathfrak{s} \! - \! 2 \rbrace$,
\begin{equation*}
\tilde{\mathfrak{Z}}^{f}_{\widetilde{V}}[\mu^{\text{\tiny {\rm EQ}}}] 
\! = \! 
\begin{cases}
\iint_{\mathbb{R}^{2}} \ln \left(\dfrac{\lvert \xi \! - \! \tau 
\rvert^{k} \prod_{j=1}^{\mathfrak{s}-2}(\lvert \xi \! - \! \alpha_{p_{j}} 
\rvert \lvert \tau \! - \! \alpha_{p_{j}} \rvert)^{\gamma_{i(j)_{k_{j}}}}
(\lvert \xi \! - \! \alpha_{k} \rvert \lvert \tau \! - \! \alpha_{k} 
\rvert)^{\gamma_{k}}}{\lvert \xi \! - \! \tau \rvert^{K}(\lvert \xi \! - \! 
\alpha_{k} \rvert \lvert \tau \! - \! \alpha_{k} \rvert)^{\varrho_{k}-1}} 
\right)^{-1} \md \mu^{\text{\tiny {\rm EQ}}}(\xi) \, 
\md \mu^{\text{\tiny {\rm EQ}}}(\tau), &\text{$\hat{\mathfrak{J}}_{q}(k) 
\! = \! \varnothing$,} \\
\iint_{\mathbb{R}^{2}} \ln \left(\dfrac{\lvert \xi \! - \! \tau 
\rvert^{k} \prod_{j=1}^{\mathfrak{s}-2}(\lvert \xi \! - \! \alpha_{p_{j}} 
\rvert \lvert \tau \! - \! \alpha_{p_{j}} \rvert)^{\gamma_{i(j)_{k_{j}}}}
(\lvert \xi \! - \! \alpha_{k} \rvert \lvert \tau \! - \! \alpha_{k} \rvert)^{
\gamma_{k}}}{\lvert \xi \! - \! \tau \rvert^{K} \prod_{j=1}^{\mathfrak{s}
-2}(\lvert \xi \! - \! \alpha_{p_{j}} \rvert \lvert \tau \! - \! \alpha_{p_{j}} 
\rvert)^{\varrho_{\hat{m}_{j}(k)}}(\lvert \xi \! - \! \alpha_{k} \rvert 
\lvert \tau \! - \! \alpha_{k} \rvert)^{\varrho_{k}-1}} \right)^{-1} \md 
\mu^{\text{\tiny {\rm EQ}}}(\xi) \, \md \mu^{\text{\tiny {\rm EQ}}}
(\tau), &\text{$\hat{\mathfrak{J}}_{q}(k) \! \neq \! \varnothing$.}
\end{cases}
\end{equation*}
Even though, at first glance, it might appear that the associated 
`core' energy functionals, $\mathfrak{X}^{\infty}_{\widetilde{V}}
[\mu^{\text{\tiny {\rm EQ}}}]$ and $\mathfrak{X}^{f}_{\widetilde{V}}
[\mu^{\text{\tiny {\rm EQ}}}]$, respectively, are different, this turns 
out not to be the case; in fact, a permutation argument based on 
the partitions introduced in Subsections~\ref{subsubsec1.2.1} 
and~\ref{subsubsec1.2.2} reveals that, in the double-scaling limit 
$\mathscr{N},n \! \to \! \infty$ such that $z_{o} \! = \! 1 \! + \! 
o(1)$, not only are the factors $\prod_{q=1}^{\mathfrak{s}-1}
(\lvert \xi \! - \! \alpha_{p_{q}} \rvert \lvert \tau \! - \! \alpha_{p_{q}} 
\rvert)^{\gamma_{i(q)_{k_{q}}}}$ and $\prod_{q=1}^{\mathfrak{s}-2}
(\lvert \xi \! - \! \alpha_{p_{q}} \rvert \lvert \tau \! - \! \alpha_{p_{q}} 
\rvert)^{\gamma_{i(q)_{k_{q}}}}(\lvert \xi \! - \! \alpha_{k} \rvert 
\lvert \tau \! - \! \alpha_{k} \rvert)^{\gamma_{k}}$ invariant with 
respect to $k$ $(\in \! \lbrace 1,2,\dotsc,K \rbrace)$, but they actually 
coincide, as a consequence of which $\mathfrak{X}^{\infty}_{\widetilde{V}}
[\mu^{\text{\tiny {\rm EQ}}}] \! = \! \mathfrak{X}^{f}_{\widetilde{V}}
[\mu^{\text{\tiny {\rm EQ}}}];$ e.g., for the $K \! = \! 7$ pole set 
$\lbrace \alpha_{1},\alpha_{2},\alpha_{3},\alpha_{4},\alpha_{5},
\alpha_{6},\alpha_{7} \rbrace \! = \! \lbrace 0,1,\infty,1,\sqrt{2},\pi,
\infty \rbrace$, for which $\mathfrak{s} \! = \! 5$, $\alpha_{k} \! = \! 
\infty$ for $k \! = \! 3,7$, and $\alpha_{k} \! \neq \! \infty$ for $k \! 
= \! 1,2,4,5,6$, the common expression for the associated `core' 
energy functionals reads
\begin{align*}
\mathfrak{X}^{\infty}_{\widetilde{V}}[\mu^{\text{\tiny {\rm EQ}}}] \! 
= \! \mathfrak{X}^{f}_{\widetilde{V}}[\mu^{\text{\tiny {\rm EQ}}}] 
=& \, \iint_{\mathbb{R}^{2}} \ln \left(\dfrac{\lvert \xi \! - \! \tau 
\rvert^{7}}{\lvert \xi \rvert \lvert \tau \rvert (\lvert \xi \! - \! 1 \rvert 
\lvert \tau \! - \! 1 \rvert)^{2} \lvert \xi \! - \! \sqrt{2} \rvert \lvert 
\tau \! - \! \sqrt{2} \rvert \lvert \xi \! - \! \pi \rvert \lvert \tau \! - \! 
\pi \rvert} \right)^{-1} \md \mu^{\text{\tiny {\rm EQ}}}(\xi) \, \md 
\mu^{\text{\tiny {\rm EQ}}}(\tau) \! + \! \int_{\mathbb{R}}V(\xi) \, 
\md \mu^{\text{\tiny {\rm EQ}}}(\xi).
\end{align*}
Furthermore, if all the poles in the sequence $\lbrace \alpha_{1},\alpha_{2},
\dotsc,\alpha_{K} \rbrace$ are distinct, that is, $\alpha_{i} \! \neq \! 
\alpha_{j}$ $\forall$ $i \! \neq \! j \! \in \! \lbrace 1,2,\dotsc,K 
\rbrace$, then
\begin{equation*}
\mathfrak{X}^{\infty}_{\widetilde{V}}[\mu^{\text{\tiny {\rm EQ}}}] \! 
= \! \mathfrak{X}^{f}_{\widetilde{V}}[\mu^{\text{\tiny {\rm EQ}}}] \! 
= \! \iint_{\mathbb{R}^{2}} \ln \left(\dfrac{\lvert \xi \! - \! \tau 
\rvert^{K}}{\prod_{\lbrace j \in \lbrace 1,2,\dotsc,K \rbrace, \, \alpha_{j} 
\neq \infty \rbrace} \lvert \xi \! - \! \alpha_{j} \rvert \lvert \tau \! - \! 
\alpha_{j} \rvert} \right)^{-1} \md \mu^{\text{\tiny {\rm EQ}}}(\xi) \, \md 
\mu^{\text{\tiny {\rm EQ}}}(\tau) \! + \! \int_{\mathbb{R}}V(\xi) \, \md 
\mu^{\text{\tiny {\rm EQ}}}(\xi).
\end{equation*}}
\end{eeeee}
The following Lemma~\ref{lem3.2}, which subsumes two families of 
variational inequalities, is necessary in order to prove, for $n \! \in \! 
\mathbb{N}$ and $k \! \in \! \lbrace 1,2,\dotsc,K \rbrace$ such that 
$\alpha_{p_{\mathfrak{s}}} \! := \! \alpha_{k} \! = \! \infty$ (resp., 
$\alpha_{p_{\mathfrak{s}}} \! := \! \alpha_{k} \! \neq \! \infty)$, the 
unicity of the associated equilibrium measure $\mathscr{M}_{1}
(\mathbb{R}) \! \ni \! \mu_{\widetilde{V}}^{\infty}$ (resp., 
$\mu_{\widetilde{V}}^{f})$, whose support, $J_{\infty} \! := \! 
\supp (\mu_{\widetilde{V}}^{\infty})$ (resp., $J_{f} \! := \! \supp 
(\mu_{\widetilde{V}}^{f}))$, is a proper compact subset of 
$\overline{\mathbb{R}} \setminus \lbrace\alpha_{1},\alpha_{2},
\dotsc,\alpha_{K} \rbrace$.
\begin{ccccc} \label{lem3.2} 
Let the external field $\widetilde{V} \colon \overline{\mathbb{R}} 
\setminus \lbrace \alpha_{1},\alpha_{2},\dotsc,\alpha_{K} \rbrace 
\! \to \! \mathbb{R}$ satisfy conditions~\eqref{eq20}--\eqref{eq22}, 
and set $\varpi (z) \! := \! \exp (-\widetilde{V}(z)/2)$. For $n \! \in \! 
\mathbb{N}$ and $k \! \in \! \lbrace 1,2,\dotsc,K \rbrace$ such that 
$\alpha_{p_{\mathfrak{s}}} \! := \! \alpha_{k} \! = \! \infty$ (resp., 
$\alpha_{p_{\mathfrak{s}}} \! := \! \alpha_{k} \! \neq \! \infty)$, let 
$\mu_{j}^{\infty}$ (resp., $\mu_{j}^{f})$ $\in \! \mathscr{M}_{1}(\mathbb{R})$, 
$j \! = \! 1,2$, be positive, finite-moment measures on $\mathbb{R}$ 
supported on distinct compact sets, that is, $\int_{\supp (\mu_{j}^{r})} 
\lambda^{m} \, \md \mu_{j}^{r}(\lambda) \! < \! \infty$ and 
$\int_{\supp (\mu_{j}^{\infty})}(\lambda \! - \! \alpha_{p_{q}})^{-m} 
\, \md \mu_{j}^{\infty}(\lambda) \! < \! \infty$, $q \! = \! 1,2,\dotsc,
\mathfrak{s} \! - \! 1$ (resp., $\int_{\supp (\mu_{j}^{f})}(\lambda \! - \! 
\alpha_{p_{q}})^{-m} \, \md \mu_{j}^{f}(\lambda) \! < \! \infty$, $q \! 
= \! 1,\dotsc,\mathfrak{s} \! - \! 2,\mathfrak{s})$, $j \! = \! 1,2$, $r \! 
\in \! \lbrace \infty,f \rbrace$, $m \! \in \! \mathbb{N}$, and $\supp 
(\mu_{1}^{r}) \cap \supp (\mu_{2}^{r}) \! = \! \varnothing$. For $n \! 
\in \! \mathbb{N}$ and $k \! \in \! \lbrace 1,2,\dotsc,K \rbrace$ such 
that $\alpha_{p_{\mathfrak{s}}} \! := \! \alpha_{k} \! = \! \infty$ (resp., 
$\alpha_{p_{\mathfrak{s}}} \! := \! \alpha_{k} \! \neq \! \infty)$, let 
$\mu^{\infty} \! := \! \mu_{1}^{\infty} \! - \! \mu_{2}^{\infty}$ (resp., 
$\mu^{f} \! := \! \mu_{1}^{f} \! - \! \mu_{2}^{f})$ be the unique Jordan 
decomposition of the finite-moment signed measure on $\mathbb{R}$ 
with compact support and mean zero, that is, $\int_{\supp (\mu^{r})} 
\md \mu^{r}(\lambda) \! = \! 0$, $r \! \in \! \lbrace \infty,f \rbrace$$;$ 
furthermore, suppose that, for $n \! \in \! \mathbb{N}$ and $k \! \in \! 
\lbrace 1,2,\dotsc,K \rbrace$ such that $\alpha_{p_{\mathfrak{s}}} \! 
:= \! \alpha_{k} \! = \! \infty$, the measures $\mu_{1}^{\infty}$ and 
$\mu_{2}^{\infty}$ have finite logarithmic energy,
\begin{equation*}
-\infty \! < \! \iint_{\mathbb{R}^{2}} \ln \left(\vert \xi \! - \! \tau 
\vert^{\frac{\varkappa_{nk}}{n}} \prod_{q=1}^{\mathfrak{s}-1} 
\left(\dfrac{\vert \xi \! - \! \tau \vert}{\vert \xi \! - \! \alpha_{p_{q}} 
\vert \vert \tau \! - \! \alpha_{p_{q}} \vert} \right)^{\frac{\varkappa_{nk 
\tilde{k}_{q}}}{n}} \right)^{-1} \md \mu^{\infty}_{j}(\xi) \, \md 
\mu^{\infty}_{j}(\tau) \! < \! +\infty, \quad j \! = \! 1,2,
\end{equation*}
and, for $n \! \in \! \mathbb{N}$ and $k \! \in \! \lbrace 1,2,\dotsc,
K \rbrace$ such that $\alpha_{p_{\mathfrak{s}}} \! := \! \alpha_{k} \! 
\neq \! \infty$, the measures $\mu_{1}^{f}$ and $\mu_{2}^{f}$ have 
finite logarithmic energy,
\begin{equation*}
-\infty \! < \! \iint_{\mathbb{R}^{2}} \ln \left(\vert \xi \! - \! \tau 
\vert^{\frac{\varkappa_{nk \tilde{k}_{\mathfrak{s}-1}}^{\infty}+1}{n}} 
\left(\dfrac{\vert \xi \! - \! \tau \vert}{\vert \xi \! - \! \alpha_{k} \vert 
\vert \tau \! - \! \alpha_{k} \vert} \right)^{\frac{\varkappa_{nk}-1}{n}} 
\prod_{q=1}^{\mathfrak{s}-2} \left(\dfrac{\vert \xi \! - \! \tau \vert}{\vert 
\xi \! - \! \alpha_{p_{q}} \vert \vert \tau \! - \! \alpha_{p_{q}} \vert} 
\right)^{\frac{\varkappa_{nk \tilde{k}_{q}}}{n}} \right)^{-1} \md \mu_{j}^{f}
(\xi) \, \md \mu_{j}^{f}(\tau) \! < \! +\infty, \quad j \! = \! 1,2.
\end{equation*}
Then$:$ {\rm (i)} for $n \! \in \! \mathbb{N}$ and $k \! \in \! \lbrace 1,2,
\dotsc,K \rbrace$ such that $\alpha_{p_{\mathfrak{s}}} \! := \! \alpha_{k} 
\! = \! \infty$,
\begin{align*}
& \, \iint_{\mathbb{R}^{2}} \ln \left(\vert \xi \! - \! \tau 
\vert^{\frac{\varkappa_{nk}}{n}} \prod_{q=1}^{\mathfrak{s}-1} \left(
\dfrac{\vert \xi \! - \! \tau \vert}{\vert \xi \! - \! \alpha_{p_{q}} \vert 
\vert \tau \! - \! \alpha_{p_{q}} \vert} \right)^{\frac{\varkappa_{nk 
\tilde{k}_{q}}}{n}} \right)^{-1} \md \mu^{\infty}(\xi) \, \md \mu^{\infty}
(\tau) \\
=& \, \iint_{\mathbb{R}^{2}} \ln \left(\vert \xi \! - \! \tau 
\vert^{\frac{\varkappa_{nk}}{n}} \prod_{q=1}^{\mathfrak{s}-1} \left(
\dfrac{\vert \xi \! - \! \tau \vert}{\vert \xi \! - \! \alpha_{p_{q}} \vert 
\vert \tau \! - \! \alpha_{p_{q}} \vert} \right)^{\frac{\varkappa_{nk 
\tilde{k}_{q}}}{n}} \varpi (\xi) \varpi (\tau) \right)^{-1} \md \mu^{\infty}
(\xi) \, \md \mu^{\infty}(\tau) \\
=& \, \left(\dfrac{(n \! - \! 1)K \! + \! k}{n} \right) \int_{0}^{+\infty} 
\lambda^{-1} \lvert \hat{\mu}_{\infty}(\lambda) \rvert^{2} \, \md 
\lambda \! \geqslant \! 0,
\end{align*}
where $\hat{\mu}_{\infty}(\lambda) \! = \! \widehat{\mu^{\infty}_{1}}
(\lambda) \! - \! \widehat{\mu^{\infty}_{2}}(\lambda)$, with 
$\widehat{\mu^{\infty}_{j}}(\lambda) \! := \! \int_{\mathbb{R}} 
\me^{\mi \lambda t} \, \md \mu^{\infty}_{j}(t)$, $j \! = \! 1,2$, and 
equality holds if, and only if, $\mu^{\infty} \! = \! 0$$;$ and {\rm (ii)} 
for $n \! \in \! \mathbb{N}$ and $k \! \in \! \lbrace 1,2,\dotsc,K \rbrace$ 
such that $\alpha_{p_{\mathfrak{s}}} \! := \! \alpha_{k} \! \neq \! \infty$,
\begin{align*}
& \, \iint_{\mathbb{R}^{2}} \ln \left(\vert \xi \! - \! \tau 
\vert^{\frac{\varkappa_{nk \tilde{k}_{\mathfrak{s}-1}}^{\infty}+1}{n}} 
\left(\dfrac{\vert \xi \! - \! \tau \vert}{\vert \xi \! - \! \alpha_{k} \vert 
\vert \tau \! - \! \alpha_{k} \vert} \right)^{\frac{\varkappa_{nk}-1}{n}} 
\prod_{q=1}^{\mathfrak{s}-2} \left(\dfrac{\vert \xi \! - \! \tau \vert}{\vert 
\xi \! - \! \alpha_{p_{q}} \vert \vert \tau \! - \! \alpha_{p_{q}} \vert} 
\right)^{\frac{\varkappa_{nk \tilde{k}_{q}}}{n}} \right)^{-1} \md \mu^{f}
(\xi) \, \md \mu^{f}(\tau) \\
=& \, \iint_{\mathbb{R}^{2}} \ln \left(\vert \xi \! - \! \tau 
\vert^{\frac{\varkappa_{nk \tilde{k}_{\mathfrak{s}-1}}^{\infty}+1}{n}} \left(
\dfrac{\vert \xi \! - \! \tau \vert}{\vert \xi \! - \! \alpha_{k} \vert 
\vert \tau \! - \! \alpha_{k} \vert} \right)^{\frac{\varkappa_{nk}-1}{n}} 
\prod_{q=1}^{\mathfrak{s}-2} \left(\dfrac{\vert \xi \! - \! \tau \vert}{\vert 
\xi \! - \! \alpha_{p_{q}} \vert \vert \tau \! - \! \alpha_{p_{q}} \vert} 
\right)^{\frac{\varkappa_{nk \tilde{k}_{q}}}{n}} \varpi (\xi) \varpi (\tau) 
\right)^{-1} \md \mu^{f}(\xi) \, \md \mu^{f}(\tau) \\
=& \, \left(\dfrac{(n \! - \! 1)K \! + \! k}{n} \right) \int_{0}^{+\infty} 
\lambda^{-1} \vert \hat{\mu}_{f}(\lambda) \vert^{2} \, \md \lambda \! 
\geqslant \! 0,
\end{align*}
where $\hat{\mu}_{f}(\lambda) \! = \! \widehat{\mu^{f}_{1}}(\lambda) \! - 
\! \widehat{\mu^{f}_{2}}(\lambda)$, with $\widehat{\mu^{f}_{j}}(\lambda) 
\! := \! \int_{\mathbb{R}} \me^{\mi \lambda t} \, \md \mu^{f}_{j}(t)$, 
$j \! = \! 1,2$, and equality holds if, and only if, $\mu^{f} \! = \! 0$.
\end{ccccc}
\begin{eeeee} \label{remm3.3} 
\textsl{Lemma~\ref{lem3.2} states that (upon noting the symmetry of 
the respective integrands under the interchange $\xi \! \leftrightarrow 
\! \tau)$$:$ {\rm (i)} for $n \! \in \! \mathbb{N}$ and $k \! \in \! \lbrace 
1,2,\dotsc,K \rbrace$ such that $\alpha_{p_{\mathfrak{s}}} \! := \! 
\alpha_{k} \! = \! \infty$,
\begin{align*}
& \, \iint_{\mathbb{R}^{2}} \ln \left(\vert \xi \! - \! \tau 
\vert^{\frac{\varkappa_{nk}}{n}} \prod_{q=1}^{\mathfrak{s}-1} \left(
\dfrac{\vert \xi \! - \! \tau \vert}{\vert \xi \! - \! \alpha_{p_{q}} \vert 
\vert \tau \! - \! \alpha_{p_{q}} \vert} \right)^{\frac{\varkappa_{nk 
\tilde{k}_{q}}}{n}} \varpi (\xi) \varpi (\tau) \right)^{-1} \left(\md 
\mu_{1}^{\infty}(\xi) \, \md \mu_{1}^{\infty}(\tau) \! + \! \md 
\mu_{2}^{\infty}(\xi) \, \md \mu_{2}^{\infty}(\tau) \right) \\
\geqslant& \, \iint_{\mathbb{R}^{2}} \ln \left(\vert \xi \! - \! 
\tau \vert^{\frac{\varkappa_{nk}}{n}} \prod_{q=1}^{\mathfrak{s}-1} \left(
\dfrac{\vert \xi \! - \! \tau \vert}{\vert \xi \! - \! \alpha_{p_{q}} \vert 
\vert \tau \! - \! \alpha_{p_{q}} \vert} \right)^{\frac{\varkappa_{nk 
\tilde{k}_{q}}}{n}} \varpi (\xi) \varpi (\tau) \right)^{-1} \left(\md 
\mu_{1}^{\infty}(\xi) \, \md \mu_{2}^{\infty}(\tau) \! + \! \md 
\mu_{2}^{\infty}(\xi) \, \md \mu_{1}^{\infty}(\tau) \right) \\
=& \, 2 \iint_{\mathbb{R}^{2}} \ln \left(\vert \xi \! - \! 
\tau \vert^{\frac{\varkappa_{nk}}{n}} \prod_{q=1}^{\mathfrak{s}-1} \left(
\dfrac{\vert \xi \! - \! \tau \vert}{\vert \xi \! - \! \alpha_{p_{q}} \vert 
\vert \tau \! - \! \alpha_{p_{q}} \vert} \right)^{\frac{\varkappa_{nk 
\tilde{k}_{q}}}{n}} \varpi (\xi) \varpi (\tau) \right)^{-1} \md 
\mu_{1}^{\infty}(\xi) \, \md \mu_{2}^{\infty}(\tau) \\
=& \, 2 \iint_{\mathbb{R}^{2}} \ln \left(\vert \xi \! - \! 
\tau \vert^{\frac{\varkappa_{nk}}{n}} \prod_{q=1}^{\mathfrak{s}-1} \left(
\dfrac{\vert \xi \! - \! \tau \vert}{\vert \xi \! - \! \alpha_{p_{q}} \vert 
\vert \tau \! - \! \alpha_{p_{q}} \vert} \right)^{\frac{\varkappa_{nk 
\tilde{k}_{q}}}{n}} \varpi (\xi) \varpi (\tau) \right)^{-1} \md 
\mu_{2}^{\infty}(\xi) \, \md \mu_{1}^{\infty}(\tau);
\end{align*}
and {\rm (ii)} for $n \! \in \! \mathbb{N}$ and $k \! \in \! \lbrace 
1,2,\dotsc,K \rbrace$ such that $\alpha_{p_{\mathfrak{s}}} \! := \! 
\alpha_{k} \! \neq \! \infty$,
\begin{align*}
& \, \iint_{\mathbb{R}^{2}} \ln \left(\vert \xi \! - \! \tau 
\vert^{\frac{\varkappa_{nk \tilde{k}_{\mathfrak{s}-1}}^{\infty}+1}{n}} 
\left(\dfrac{\vert \xi \! - \! \tau \vert}{\vert \xi \! - \! \alpha_{k} \vert 
\vert \tau \! - \! \alpha_{k} \vert} \right)^{\frac{\varkappa_{nk}-1}{n}} 
\prod_{q=1}^{\mathfrak{s}-2} \left(\dfrac{\vert \xi \! - \! \tau \vert}{\vert 
\xi \! - \! \alpha_{p_{q}} \vert \vert \tau \! - \! \alpha_{p_{q}} \vert} 
\right)^{\frac{\varkappa_{nk \tilde{k}_{q}}}{n}} \varpi (\xi) \varpi (\tau) 
\right)^{-1} \left(\md \mu_{1}^{f}(\xi) \, \md \mu_{1}^{f}(\tau) \right. \\
+&\left. \, \md \mu_{2}^{f}(\xi) \, \md \mu_{2}^{f}(\tau) \right) \! 
\geqslant \! \iint_{\mathbb{R}^{2}} \ln \left(\vert \xi \! - \! \tau 
\vert^{\frac{\varkappa_{nk \tilde{k}_{\mathfrak{s}-1}}^{\infty}+1}{n}} 
\left(\dfrac{\vert \xi \! - \! \tau \vert}{\vert \xi \! - \! \alpha_{k} \vert 
\vert \tau \! - \! \alpha_{k} \vert} \right)^{\frac{\varkappa_{nk}-1}{n}} 
\prod_{q=1}^{\mathfrak{s}-2} \left(\dfrac{\vert \xi \! - \! \tau \vert}{\vert 
\xi \! - \! \alpha_{p_{q}} \vert \vert \tau \! - \! \alpha_{p_{q}} \vert} 
\right)^{\frac{\varkappa_{nk \tilde{k}_{q}}}{n}} \varpi (\xi) \varpi (\tau) 
\right)^{-1} \\
\times& \, \left(\md \mu_{1}^{f}(\xi) \, \md \mu_{2}^{f}(\tau) \! + \! 
\md \mu_{2}^{f}(\xi) \, \md \mu_{1}^{f}(\tau) \right) \\
=& \, 2 \iint_{\mathbb{R}^{2}} \ln \left(\vert \xi \! - \! \tau 
\vert^{\frac{\varkappa_{nk \tilde{k}_{\mathfrak{s}-1}}^{\infty}+1}{n}} 
\left(\dfrac{\vert \xi \! - \! \tau \vert}{\vert \xi \! - \! \alpha_{k} \vert 
\vert \tau \! - \! \alpha_{k} \vert} \right)^{\frac{\varkappa_{nk}-1}{n}} 
\prod_{q=1}^{\mathfrak{s}-2} \left(\dfrac{\vert \xi \! - \! \tau \vert}{\vert 
\xi \! - \! \alpha_{p_{q}} \vert \vert \tau \! - \! \alpha_{p_{q}} \vert} 
\right)^{\frac{\varkappa_{nk \tilde{k}_{q}}}{n}} \varpi (\xi) \varpi (\tau) 
\right)^{-1} \md \mu_{1}^{f}(\xi) \, \md \mu_{2}^{f}(\tau) \\
=& \, 2 \iint_{\mathbb{R}^{2}} \ln \left(\vert \xi \! - \! \tau 
\vert^{\frac{\varkappa_{nk \tilde{k}_{\mathfrak{s}-1}}^{\infty}+1}{n}} 
\left(\dfrac{\vert \xi \! - \! \tau \vert}{\vert \xi \! - \! \alpha_{k} \vert 
\vert \tau \! - \! \alpha_{k} \vert} \right)^{\frac{\varkappa_{nk}-1}{n}} 
\prod_{q=1}^{\mathfrak{s}-2} \left(\dfrac{\vert \xi \! - \! \tau \vert}{\vert 
\xi \! - \! \alpha_{p_{q}} \vert \vert \tau \! - \! \alpha_{p_{q}} \vert} 
\right)^{\frac{\varkappa_{nk \tilde{k}_{q}}}{n}} \varpi (\xi) \varpi (\tau) 
\right)^{-1} \md \mu_{2}^{f}(\xi) \, \md \mu_{1}^{f}(\tau).
\end{align*}
In other words, if, for $n \! \in \! \mathbb{N}$ and $k \! \in \! \lbrace 1,2,
\dotsc,K \rbrace$ such that $\alpha_{p_{\mathfrak{s}}} \! := \! \alpha_{k} 
\! = \! \infty$ (resp., $\alpha_{p_{\mathfrak{s}}} \! := \! \alpha_{k} \! \neq 
\! \infty)$, $\ln (\vert \xi \! - \! \tau \vert^{\frac{\varkappa_{nk}}{n}} 
\prod_{q=1}^{\mathfrak{s}-1}(\vert \xi - \tau \vert/\vert \xi - \alpha_{p_{q}} 
\vert \vert \tau - \alpha_{p_{q}} \vert)^{\frac{\varkappa_{nk \tilde{k}_{q}}}{n}})$ 
(resp., $\ln (\vert \xi \! - \! \tau \vert^{\frac{\varkappa_{nk \tilde{k}_{
\mathfrak{s}-1}}^{\infty}+1}{n}}(\vert \xi - \tau \vert/\vert \xi - \alpha_{k} 
\vert \vert \tau - \alpha_{k} \vert)^{\frac{\varkappa_{nk}-1}{n}} \prod_{q=
1}^{\mathfrak{s}-2}(\vert \xi - \tau \vert/\vert \xi - \alpha_{p_{q}} \vert \vert 
\tau - \alpha_{p_{q}} \vert)^{\frac{\varkappa_{nk \tilde{k}_{q}}}{n}}))$ is integrable 
with respect to the product measures $\md \mu_{1}^{\infty}(\lambda) \, \md 
\mu_{1}^{\infty}(t)$ and $\md \mu_{2}^{\infty}(\lambda) \, \md \mu_{2}^{\infty}(t)$ 
(resp., $\md \mu_{1}^{f}(\lambda) \, \md \mu_{1}^{f}(t)$ and $\md \mu_{2}^{f}
(\lambda) \, \md \mu_{2}^{f}(t))$, then it is integrable with respect to the `mixed' 
product measures $\md \mu_{1}^{\infty}(\lambda) \, \md \mu_{2}^{\infty}(t)$ 
and $\md \mu_{2}^{\infty}(\lambda) \, \md \mu_{1}^{\infty}(t)$ (resp., $\md 
\mu_{1}^{f}(\lambda) \, \md \mu_{2}^{f}(t)$ and $\md \mu_{2}^{f}(\lambda) \, 
\md \mu_{1}^{f}(t))$.}
\end{eeeee}

\emph{Proof.} The proof of this Lemma~\ref{lem3.2} consists of two cases: (i) 
$n \! \in \! \mathbb{N}$ and $k \! \in \! \lbrace 1,2,\dotsc,K \rbrace$ such 
that $\alpha_{p_{\mathfrak{s}}} \! := \! \alpha_{k} \! = \! \infty$; and (ii) 
$n \! \in \! \mathbb{N}$ and $k \! \in \! \lbrace 1,2,\dotsc,K \rbrace$ such 
that $\alpha_{p_{\mathfrak{s}}} \! := \! \alpha_{k} \! \neq \! \infty$. The 
proof for case~(ii) will be considered in detail (see $\pmb{(1)}$ below), 
whilst case~(i) can be proved analogously (see $\pmb{(2)}$ below).

$\pmb{(1)}$ Recall the following identity (see \cite{a51}, Chapter~6, 
p.~147, Equation~(6.44)): for $\lambda \! \in \! \mathbb{R}$ and any 
$\varepsilon \! > \! 0$,
\begin{equation*}
\ln (\lambda^{2} \! + \! \varepsilon^{2}) \! = \! \ln (\varepsilon^{2}) \! + 
\! 2 \Im \left(\int_{0}^{+\infty}(\mi v)^{-1}(\me^{\mi \lambda v} \! - \! 1) 
\me^{-\varepsilon v} \, \md v \right).
\end{equation*}
For $n \! \in \! \mathbb{N}$ and $k \! \in \! \lbrace 1,2,\dotsc,K \rbrace$ 
such that $\alpha_{p_{\mathfrak{s}}} \! := \! \alpha_{k} \! \neq \! \infty$, 
let $\widetilde{V} \colon \overline{\mathbb{R}} \setminus \lbrace 
\alpha_{1},\alpha_{2},\dotsc,\alpha_{K} \rbrace \! \to \! \mathbb{R}$ 
satisfy conditions~\eqref{eq20}--\eqref{eq22}, set $\varpi (z) \! = \! \exp 
(-\widetilde{V}(z)/2)$, and let the measures $\mu^{f}_{1}$, $\mu^{f}_{2}$ 
and $\mu^{f} \! := \! \mu^{f}_{1} \! - \! \mu^{f}_{2}$ (and their respective 
supports) satisfy all of the conditions stated in the corresponding item of 
the lemma. Arguing as in the proof of Lemma~3.2 in \cite{a45} (see, in 
particular, pp.~276--279 of \cite{a45}), one shows, via the above identity, 
that, for $n \! \in \! \mathbb{N}$ and $k \! \in \! \lbrace 1,2,\dotsc,K 
\rbrace$ such that $\alpha_{p_{\mathfrak{s}}} \! := \! \alpha_{k} \! \neq 
\! \infty$,
\begin{gather*}
\dfrac{1}{2} \left(\dfrac{\varkappa^{\infty}_{nk \tilde{k}_{\mathfrak{s}-1}} 
\! + \! 1}{n} \right) \iint_{\mathbb{R}^{2}} \ln ((\xi \! - \! \tau)^{2} 
\! + \! \varepsilon^{2}) \, \md \mu^{f}(\xi) \, \md \mu^{f}(\tau) \! = \! 
\left(\dfrac{\varkappa^{\infty}_{nk \tilde{k}_{\mathfrak{s}-1}} 
\! + \! 1}{n} \right) \Im \left(\int_{0}^{+\infty}(\mi v)^{-1} \vert 
\hat{\mu}_{f}(v) \vert^{2} \me^{-\varepsilon v} \, \md v \right), \\
\dfrac{1}{2} \left(\dfrac{\varkappa_{nk} \! - \! 1}{n} \right) \iint_{
\mathbb{R}^{2}} \ln ((\xi \! - \! \tau)^{2} \! + \! \varepsilon^{2}) 
\, \md \mu^{f}(\xi) \, \md \mu^{f}(\tau) \! = \! \left(\dfrac{\varkappa_{nk} 
\! - \! 1}{n} \right) \Im \left(\int_{0}^{+\infty}(\mi v)^{-1} \vert 
\hat{\mu}_{f}(v) \vert^{2} \me^{-\varepsilon v} \, \md v \right), \\
\dfrac{1}{2} \dfrac{\varkappa_{nk \tilde{k}_{q}}}{n} \iint_{\mathbb{R}^{2}} 
\ln ((\xi \! - \! \tau)^{2} \! + \! \varepsilon^{2}) \, \md \mu^{f}(\xi) 
\, \md \mu^{f}(\tau) \! = \! \dfrac{\varkappa_{nk \tilde{k}_{q}}}{n} 
\Im \left(\int_{0}^{+\infty}(\mi v)^{-1} \vert \hat{\mu}_{f}(v) \vert^{2} 
\me^{-\varepsilon v} \, \md v \right), \quad q \! = \! 1,2,\dotsc,
\mathfrak{s} \! - \! 2, \\
\iint_{\mathbb{R}^{2}} \ln ((r \! - \! t)^{2} \! + \! 
\varepsilon^{2}) \, \md \mu^{f}(\xi) \, \md \mu^{f}(\tau) \! = \! 0, 
\quad r \! \in \! \lbrace \xi,\tau \rbrace, \quad t \! \in \! \lbrace 
\alpha_{p_{1}},\dotsc,\alpha_{p_{\mathfrak{s}-2}},\alpha_{k} \rbrace,
\end{gather*}
where $\hat{\mu}_{f}(\boldsymbol{\cdot})$ is defined in item~(ii) of 
the lemma. Noting that $\hat{\mu}_{f}(0) \! = \! \int_{\mathbb{R}} 
\md \mu^{f}(t) \! = \! 0$, a Taylor expansion about $\nu \! = \! 0$ shows 
that $\hat{\mu}_{f}(\nu) \! =_{\nu \to 0} \! \hat{\mu}_{f}^{\prime}(0) \nu 
\! + \! \mathcal{O}(\nu^{2})$, where $\hat{\mu}_{f}^{\prime}(0) \! := \! 
\partial_{\nu} \hat{\mu}_{f}(\nu) \vert_{\nu =0}$; thus, $\nu^{-1} \vert 
\hat{\mu}_{f}(\nu) \vert^{2} \! =_{\nu \to 0} \! \vert \hat{\mu}_{f}^{\prime}
(0) \vert^{2} \nu \! + \! \mathcal{O}(\nu^{2})$, which means that there is 
no singularity in the integrands at $\nu \! = \! 0$ (in fact, $\nu^{-1} \vert 
\hat{\mu}_{f}(\nu) \vert^{2}$ is real analytic in an open neighbourhood of 
the origin), whence, for $n \! \in \! \mathbb{N}$ and $k \! \in \! \lbrace 
1,2,\dotsc,K \rbrace$ such that $\alpha_{p_{\mathfrak{s}}} \! := \! 
\alpha_{k} \! \neq \! \infty$,
\begin{gather*}
\iint_{\mathbb{R}^{2}} \ln ((\xi \! - \! \tau)^{2} \! + \! 
\varepsilon^{2})^{-\frac{(\varkappa^{\infty}_{nk \tilde{k}_{\mathfrak{s}-1}} 
\! + \! 1)}{2n}} \, \md \mu^{f}(\xi) \, \md \mu^{f}(\tau) \! = \! \left(
\dfrac{\varkappa^{\infty}_{nk \tilde{k}_{\mathfrak{s}-1}} \! + \! 1}{n} 
\right) \int_{0}^{+\infty} \lambda^{-1} \vert \hat{\mu}_{f}(\lambda) 
\vert^{2} \me^{-\varepsilon \lambda} \, \md \lambda, \\
\iint_{\mathbb{R}^{2}} \ln ((\xi \! - \! \tau)^{2} \! + \! 
\varepsilon^{2})^{-\frac{(\varkappa_{nk}-1)}{2n}} \, \md \mu^{f}(\xi) \, 
\md \mu^{f}(\tau) \! = \! \left(\dfrac{\varkappa_{nk} \! - \! 1}{n} \right) 
\int_{0}^{+\infty} \lambda^{-1} \vert \hat{\mu}_{f}(\lambda) \vert^{2} 
\me^{-\varepsilon \lambda} \, \md \lambda, \\
\iint_{\mathbb{R}^{2}} \ln ((\xi \! - \! \tau)^{2} \! + \! 
\varepsilon^{2})^{-\frac{\varkappa_{nk \tilde{k}_{q}}}{2n}} \, \md \mu^{f}
(\xi) \, \md \mu^{f}(\tau) \! = \! \dfrac{\varkappa_{nk \tilde{k}_{q}}}{n} 
\int_{0}^{+\infty} \lambda^{-1} \vert \hat{\mu}_{f}(\lambda) \vert^{2} 
\me^{-\varepsilon \lambda} \, \md \lambda, \quad q \! = \! 1,2,\dotsc,
\mathfrak{s} \! - \! 2.
\end{gather*}
Now, adding the above, noting that $\iint_{\mathbb{R}^{2}} 
\ln ((r \! - \! t)^{2} \! + \! \varepsilon^{2}) \, \md \mu^{f}(\xi) \, \md 
\mu^{f}(\tau) \! = \! 0$, $r \! \in \! \lbrace \xi,\tau \rbrace$, $t \! \in 
\! \lbrace \alpha_{p_{1}},\dotsc,\alpha_{p_{\mathfrak{s}-2}},\alpha_{k} 
\rbrace$, and $\iint_{\mathbb{R}^{2}} \ln (\varpi (\xi) \varpi 
(\tau))^{-1} \, \md \mu^{f}(\xi) \, \md \mu^{f}(\tau) \! = \! 0$, 
recalling that $\sum_{q=1}^{\mathfrak{s}-2} \varkappa_{nk \tilde{k}_{q}} 
\! + \! \varkappa^{\infty}_{nk \tilde{k}_{\mathfrak{s}-1}} \! + \! 
\varkappa_{nk} \! = \! (n \! - \! 1)K \! + \! k$, noting that the respective 
integrands are uniformly bounded (some {}from above and some {}from 
below) with respect to $\varepsilon$ and that the associated measures 
have compact support, letting $\varepsilon \! \downarrow \! 0$ and 
using the Monotone Convergence and Dominated Convergence Theorems, 
one arrives at, for $n \! \in \! \mathbb{N}$ and $k \! \in \! \lbrace 1,2,
\dotsc,K \rbrace$ such that $\alpha_{p_{\mathfrak{s}}} \! := \! \alpha_{k} 
\! \neq \! \infty$,
\begin{align}
& \iint_{\mathbb{R}^{2}} \ln \left(\left((\xi \! - \! \tau)^{2} 
\! + \! \varepsilon^{2} \right)^{\frac{\varkappa_{nk \tilde{k}_{\mathfrak{s}
-1}}^{\infty}+1}{2n}} \left(\dfrac{(\xi \! - \! \tau)^{2} \! + \! 
\varepsilon^{2}}{((\xi \! - \! \alpha_{k})^{2} \! + \! \varepsilon^{2})
((\tau \! - \! \alpha_{k})^{2} \! + \! \varepsilon^{2})} \right)^{
\frac{\varkappa_{nk}-1}{2n}} \right. \nonumber \\
\times&\left. \, \prod_{q=1}^{\mathfrak{s}-2} \left(\dfrac{(\xi \! - \! 
\tau)^{2} \! + \! \varepsilon^{2}}{((\xi \! - \! \alpha_{p_{q}})^{2} \! + \! 
\varepsilon^{2})((\tau \! - \! \alpha_{p_{q}})^{2} \! + \! \varepsilon^{2})} 
\right)^{\frac{\varkappa_{nk \tilde{k}_{q}}}{2n}} \right)^{-1} \md \mu^{f}
(\xi) \, \md \mu^{f}(\tau) \nonumber \\
=& \, \iint_{\mathbb{R}^{2}} \ln \left(\left((\xi \! - \! \tau)^{2} 
\! + \! \varepsilon^{2} \right)^{\frac{\varkappa_{nk \tilde{k}_{\mathfrak{s}
-1}}^{\infty}+1}{2n}} \left(\dfrac{(\xi \! - \! \tau)^{2} \! + \! 
\varepsilon^{2}}{((\xi \! - \! \alpha_{k})^{2} \! + \! \varepsilon^{2})
((\tau \! - \! \alpha_{k})^{2} \! + \! \varepsilon^{2})} \right)^{
\frac{\varkappa_{nk}-1}{2n}} \right. \nonumber \\
\times&\left. \, \prod_{q=1}^{\mathfrak{s}-2} \left(\dfrac{(\xi \! - \! 
\tau)^{2} \! + \! \varepsilon^{2}}{((\xi \! - \! \alpha_{p_{q}})^{2} \! + \! 
\varepsilon^{2})((\tau \! - \! \alpha_{p_{q}})^{2} \! + \! \varepsilon^{2})} 
\right)^{\frac{\varkappa_{nk \tilde{k}_{q}}}{2n}} \varpi (\xi) \varpi (\tau) 
\right)^{-1} \md \mu^{f}(\xi) \, \md \mu^{f}(\tau) \nonumber \\
\underset{\varepsilon \downarrow 0}{=}& \, \iint_{\mathbb{R}^{2}} 
\ln \left(\vert \xi \! - \! \tau \vert^{\frac{\varkappa_{nk 
\tilde{k}_{\mathfrak{s}-1}}^{\infty}+1}{n}} \left(\dfrac{\vert \xi \! - \! 
\tau \vert}{\vert \xi \! - \! \alpha_{k} \vert \vert \tau \! - \! \alpha_{k} 
\vert} \right)^{\frac{\varkappa_{nk}-1}{n}} \prod_{q=1}^{\mathfrak{s}-2} 
\left(\dfrac{\vert \xi \! - \! \tau \vert}{\vert \xi \! - \! \alpha_{p_{q}} 
\vert \vert \tau \! - \! \alpha_{p_{q}} \vert} \right)^{\frac{\varkappa_{nk 
\tilde{k}_{q}}}{n}} \varpi (\xi) \varpi (\tau) \right)^{-1} \md \mu^{f}(\xi) 
\, \md \mu^{f}(\tau) \nonumber \\
=& \, \left(\dfrac{(n \! - \! 1)K \! + \! k}{n} \right) \int_{0}^{+\infty} 
\lambda^{-1} \vert \hat{\mu}_{f}(\lambda) \vert^{2} \, \md \lambda \! 
\geqslant \! 0, \label{eql3.2a}
\end{align}
where equality holds if, and only if, $\mu^{f} \! = \! 0$: this establishes 
the corresponding inequality of item~(ii) of the lemma.

$\pmb{(2)}$ The proof of this case, that is, $n \! \in \! \mathbb{N}$ 
and $k \! \in \! \lbrace 1,2,\dotsc,K \rbrace$ such that 
$\alpha_{p_{\mathfrak{s}}} \! := \! \alpha_{k} \! = \! \infty$, is virtually 
identical to the proof presented in $\pmb{(1)}$ above; one mimics, 
\emph{verbatim}, the scheme of the calculations presented in case $\pmb{(1)}$ 
in order to arrive at the associated analogue of inequality~\eqref{eql3.2a}, 
which reads
\begin{align}
& \iint_{\mathbb{R}^{2}} \ln \left(\left((\xi \! - \! \tau)^{2} 
\! + \! \varepsilon^{2} \right)^{\frac{\varkappa_{nk}}{2n}} 
\prod_{q=1}^{\mathfrak{s}-1} \left(\dfrac{(\xi \! - \! \tau)^{2} \! + \! 
\varepsilon^{2}}{((\xi \! - \! \alpha_{p_{q}})^{2} \! + \! \varepsilon^{2})
((\tau \! - \! \alpha_{p_{q}})^{2} \! + \! \varepsilon^{2})} \right)^{
\frac{\varkappa_{nk \tilde{k}_{q}}}{2n}} \right)^{-1} \md \mu^{\infty}
(\xi) \, \md \mu^{\infty}(\tau) \nonumber \\
=& \, \iint_{\mathbb{R}^{2}} \ln \left(\left((\xi \! - \! \tau)^{2} 
\! + \! \varepsilon^{2} \right)^{\frac{\varkappa_{nk}}{2n}} 
\prod_{q=1}^{\mathfrak{s}-1} \left(\dfrac{(\xi \! - \! \tau)^{2} \! + \! 
\varepsilon^{2}}{((\xi \! - \! \alpha_{p_{q}})^{2} \! + \! \varepsilon^{2})
((\tau \! - \! \alpha_{p_{q}})^{2} \! + \! \varepsilon^{2})} \right)^{
\frac{\varkappa_{nk \tilde{k}_{q}}}{2n}} \varpi (\xi) \varpi (\tau) 
\right)^{-1} \md \mu^{\infty}(\xi) \, \md \mu^{\infty}(\tau) \nonumber \\
\underset{\varepsilon \downarrow 0}{=}& \, \iint_{\mathbb{R}^{2}} 
\ln \left(\vert \xi \! - \! \tau \vert^{\frac{\varkappa_{nk}}{n}} 
\prod_{q=1}^{\mathfrak{s}-1} \left(\dfrac{\vert \xi \! - \! \tau \vert}{\vert 
\xi \! - \! \alpha_{p_{q}} \vert \vert \tau \! - \! \alpha_{p_{q}} \vert} 
\right)^{\frac{\varkappa_{nk \tilde{k}_{q}}}{n}} \varpi (\xi) \varpi (\tau) 
\right)^{-1} \md \mu^{\infty} (\xi) \, \md \mu^{\infty}(\tau) \nonumber \\
=& \, \left(\dfrac{(n \! - \! 1)K \! + \! k}{n} \right) \int_{0}^{+\infty} 
\lambda^{-1} \vert \hat{\mu}_{\infty}(\lambda) \vert^{2} \, \md \lambda 
\! \geqslant \! 0, \label{eql3.2b}
\end{align}
where $\hat{\mu}_{\infty}(\boldsymbol{\cdot})$ is defined in item~(i) of 
the lemma: this establishes the corresponding inequality of item~(i) 
of the lemma, and thus concludes the proof. \hfill $\qed$

Via the variational inequalities of Lemma~\ref{lem3.2}, one now 
establishes the unicity of the associated families of equilibrium 
measures $\mu_{\widetilde{V}}^{\infty}$ and $\mu_{\widetilde{V}}^{f}$.
\begin{ccccc} \label{lem3.3} 
Let the external field $\widetilde{V} \colon \overline{\mathbb{R}} 
\setminus \lbrace \alpha_{1},\alpha_{2},\dotsc,\alpha_{K} \rbrace 
\! \to \! \mathbb{R}$ satisfy conditions~\eqref{eq20}--\eqref{eq22}, 
and set $\varpi (z) \! := \! \exp (-\widetilde{V}(z)/2)$. For $n \! \in 
\! \mathbb{N}$ and $k \! \in \! \lbrace 1,2,\dotsc,K \rbrace$ such 
that $\alpha_{p_{\mathfrak{s}}} \! := \! \alpha_{k} \! = \! \infty$ (resp., 
$\alpha_{p_{\mathfrak{s}}} \! := \! \alpha_{k} \! \neq \! \infty)$, let 
$\mathrm{I}_{\widetilde{V}}^{\infty}$ (resp., $\mathrm{I}_{\widetilde{V}}^{f})$ 
be defined by Equation~\eqref{eqlm3.1a} (resp., Equation~\eqref{eqlm3.1b}$)$, 
and consider the associated minimisation problem $E_{\widetilde{V}}^{
\infty} \! = \! \inf \lbrace \mathstrut \mathrm{I}_{\widetilde{V}}^{\infty}
[\mu^{\text{\tiny $\mathrm{EQ}$}}_{1}]; \, 
\mu^{\text{\tiny $\mathrm{EQ}$}}_{1} \! \in \! \mathscr{M}_{1}
(\mathbb{R}) \rbrace$ (resp., $E_{\widetilde{V}}^{f} \! = \! 
\inf \lbrace \mathstrut \mathrm{I}_{\widetilde{V}}^{f}
[\mu^{\text{\tiny $\mathrm{EQ}$}}_{2}]; \, 
\mu^{\text{\tiny $\mathrm{EQ}$}}_{2} \! \in \! \mathscr{M}_{1}
(\mathbb{R}) \rbrace)$. Then, for $n \! \in \! \mathbb{N}$ and 
$k \! \in \! \lbrace 1,2,\dotsc,K \rbrace$ such that 
$\alpha_{p_{\mathfrak{s}}} \! := \! \alpha_{k} \! = \! \infty$ (resp., 
$\alpha_{p_{\mathfrak{s}}} \! := \! \alpha_{k} \! \neq \! \infty)$, 
there exists a unique $\mu_{\widetilde{V}}^{\infty}$ (resp., 
$\mu_{\widetilde{V}}^{f})$ $\in \! \mathscr{M}_{1}(\mathbb{R})$ such 
that $\mathrm{I}_{\widetilde{V}}^{\infty}[\mu_{\widetilde{V}}^{\infty}] 
\! = \! E_{\widetilde{V}}^{\infty}$ (resp., $\mathrm{I}_{\widetilde{V}}^{f}
[\mu_{\widetilde{V}}^{f}] \! = \! E_{\widetilde{V}}^{f})$.
\end{ccccc}

\emph{Proof}. The proof of this Lemma~\ref{lem3.3} consists of two 
cases: (i) $n \! \in \! \mathbb{N}$ and $k \! \in \! \lbrace 1,2,\dotsc,
K \rbrace$ such that $\alpha_{p_{\mathfrak{s}}} \! := \! \alpha_{k} \! 
= \! \infty$; and (ii) $n \! \in \! \mathbb{N}$ and $k \! \in \! \lbrace 
1,2,\dotsc,K \rbrace$ such that $\alpha_{p_{\mathfrak{s}}} \! := \! 
\alpha_{k} \! \neq \! \infty$. The proof for case~(ii) will be considered 
in detail (see $\pmb{(1)}$ below), whilst case~(i) can be proved 
analogously (see $\pmb{(2)}$ below).

$\pmb{(1)}$ Recall that, for $n \! \in \! \mathbb{N}$ and $k \! \in \! 
\lbrace 1,2,\dotsc,K \rbrace$ such that $\alpha_{p_{\mathfrak{s}}} \! 
:= \! \alpha_{k} \! \neq \! \infty$, it was shown in Lemma~\ref{lem3.1} 
that there exists $\mu^{f}_{\widetilde{V}} \! \in \! \mathscr{M}_{1}
(\mathbb{R})$, the associated equilibrium measure, such 
that $\mathrm{I}^{f}_{\widetilde{V}}[\mu^{f}_{\widetilde{V}}] \! = \! 
E^{f}_{\widetilde{V}}$ and $J_{f} \! := \! \supp (\mu^{f}_{\widetilde{V}})$ 
$(\subset \lbrace \mathstrut x \! \in \! \mathbb{R}; \, \varpi (x) \! > \! 
0 \rbrace)$ is a proper compact subset of $\overline{\mathbb{R}} \setminus 
\lbrace \alpha_{1},\alpha_{2},\dotsc,\alpha_{K} \rbrace$ (that is, for (some) 
$T_{M_{f}}$ $(= \! T_{M_{f}}(n,k,z_{o}))$ $> \! 1$, e.g., $T_{M_{f}} \! = \! 
K(1 \! + \! \max \lbrace \mathstrut \lvert \alpha_{p_{q}} \rvert, \, q \! = 
\! 1,\dotsc,\mathfrak{s} \! - \! 2,\mathfrak{s} \rbrace \! + \! 3(\min 
\lbrace \mathstrut \lvert \alpha_{p_{i}} \! - \! \alpha_{p_{j}} \rvert, \, i 
\! \neq \! j \! \in \! \lbrace 1,\dotsc,\mathfrak{s} \! - \! 2,\mathfrak{s} 
\rbrace \rbrace)^{-1})$ with $\mathscr{O}_{\frac{1}{T_{M_{f}}}}(\alpha_{p_{i}}) 
\cap \mathscr{O}_{\frac{1}{T_{M_{f}}}}(\alpha_{p_{j}}) \! = \! \varnothing$, 
$i \! \neq \! j \! \in \! \lbrace 1,\dotsc,\mathfrak{s} \! - \! 2,\mathfrak{s} 
\rbrace$, $\supp (\mu^{f}_{\widetilde{V}}) \subseteq \mathbb{R} 
\setminus (\lbrace \lvert x \rvert \! \geqslant \! T_{M_{f}} \rbrace \cup 
\cup_{\underset{q \neq \mathfrak{s}-1}{q=1}}^{\mathfrak{s}} 
\operatorname{clos}(\mathscr{O}_{\frac{1}{T_{M_{f}}}}(\alpha_{p_{q}}))))$; 
hence, it remains only to establish the uniqueness of the associated 
equilibrium measure $\mu^{f}_{\widetilde{V}}$. For $n \! \in \! 
\mathbb{N}$ and $k \! \in \! \lbrace 1,2,\dotsc,K \rbrace$ such that 
$\alpha_{p_{\mathfrak{s}}} \! := \! \alpha_{k} \! \neq \! \infty$, 
let $\widetilde{\mu}^{f}_{\widetilde{V}} \! \in \! \mathscr{M}_{1}
(\mathbb{R})$ be a second probability measure for which 
$\mathrm{I}^{f}_{\widetilde{V}}[\widetilde{\mu}^{f}_{\widetilde{V}}] 
\! = \! E^{f}_{\widetilde{V}}$: the analogue of the argument of 
Lemma~\ref{lem3.1} for the associated probability measure 
$\widetilde{\mu}^{f}_{\widetilde{V}}$ shows that $-\infty \! < \! 
\mathrm{I}^{f}_{\widetilde{V}}[\widetilde{\mu}^{f}_{\widetilde{V}}] 
\! < \! +\infty$ and $\widetilde{J}_{f} \! := \! \supp 
(\widetilde{\mu}^{f}_{\widetilde{V}})$ $(\subset \lbrace \mathstrut 
x \! \in \! \mathbb{R}; \, \varpi (x) \! > \! 0 \rbrace)$ is a proper 
compact (measurable) subset of $\overline{\mathbb{R}} \setminus 
\lbrace \alpha_{1},\alpha_{2},\dotsc,\alpha_{K} \rbrace$. For $n \! 
\in \! \mathbb{N}$ and $k \! \in \! \lbrace 1,2,\dotsc,K \rbrace$ 
such that $\alpha_{p_{\mathfrak{s}}} \! := \! \alpha_{k} \! \neq \! 
\infty$, define, for $\mu^{f}_{\widetilde{V}},\widetilde{\mu}^{f}_{
\widetilde{V}} \! \in \! \mathscr{M}_{1}(\mathbb{R})$, the finite-moment 
signed measure (on $\mathbb{R})$ $\mu^{\sharp}_{f} \! := \! 
\widetilde{\mu}^{f}_{\widetilde{V}} \! - \! \mu^{f}_{\widetilde{V}}$. For 
$n \! \in \! \mathbb{N}$ and $k \! \in \! \lbrace 1,2,\dotsc,K \rbrace$ 
such that $\alpha_{p_{\mathfrak{s}}} \! := \! \alpha_{k} \! \neq \! \infty$, 
{}from item~(ii) of Lemma~\ref{lem3.2}, with the change $\mu^{f} \! 
\to \! \mu^{\sharp}_{f}$, that is,
\begin{equation*}
\iint_{\mathbb{R}^{2}} \ln \left(\vert \xi \! - \! \tau 
\vert^{\frac{\varkappa^{\infty}_{nk \tilde{k}_{\mathfrak{s}-1}}+1}{n}} 
\left(\dfrac{\vert \xi \! - \! \tau \vert}{\vert \xi \! - \! \alpha_{k} \vert 
\vert \tau \! - \! \alpha_{k} \vert} \right)^{\frac{\varkappa_{nk}-1}{n}} 
\prod_{q=1}^{\mathfrak{s}-2} \left(\dfrac{\vert \xi \! - \! \tau \vert}{\vert 
\xi \! - \! \alpha_{p_{q}} \vert \vert \tau \! - \! \alpha_{p_{q}} \vert} 
\right)^{\frac{\varkappa_{nk \tilde{k}_{q}}}{n}} \varpi (\xi) \varpi (\tau) 
\right)^{-1} \md \mu^{\sharp}_{f}(\xi) \, \md \mu^{\sharp}_{f}(\tau) \! 
\geqslant \! 0,
\end{equation*}
it follows, via the symmetry of the integrand under the interchange 
$\xi \! \leftrightarrow \! \tau$, the definition of $\mu^{\sharp}_{f}$, 
and Remark~\ref{remm3.3}, that
\begin{align}
& \iint_{\mathbb{R}^{2}} \ln \left(\vert \xi \! - \! \tau 
\vert^{\frac{\varkappa^{\infty}_{nk \tilde{k}_{\mathfrak{s}-1}}+1}{n}} 
\left(\dfrac{\vert \xi \! - \! \tau \vert}{\vert \xi \! - \! \alpha_{k} \vert 
\vert \tau \! - \! \alpha_{k} \vert} \right)^{\frac{\varkappa_{nk}-1}{n}} 
\prod_{q=1}^{\mathfrak{s}-2} \left(\dfrac{\vert \xi \! - \! \tau \vert}{\vert 
\xi \! - \! \alpha_{p_{q}} \vert \vert \tau \! - \! \alpha_{p_{q}} \vert} 
\right)^{\frac{\varkappa_{nk \tilde{k}_{q}}}{n}} \varpi (\xi) \varpi (\tau) 
\right)^{-1} \nonumber \\
\times& \, \left(\md \widetilde{\mu}^{f}_{\widetilde{V}}(\xi) \, \md 
\widetilde{\mu}^{f}_{\widetilde{V}}(\tau) \! + \! \md \mu^{f}_{\widetilde{V}}
(\xi) \, \md \mu^{f}_{\widetilde{V}}(\tau) \right) \nonumber \\
\geqslant& \, \iint_{\mathbb{R}^{2}} \ln \left(\vert \xi \! - \! 
\tau \vert^{\frac{\varkappa^{\infty}_{nk \tilde{k}_{\mathfrak{s}-1}}+1}{n}} 
\left(\dfrac{\vert \xi \! - \! \tau \vert}{\vert \xi \! - \! \alpha_{k} \vert 
\vert \tau \! - \! \alpha_{k} \vert} \right)^{\frac{\varkappa_{nk}-1}{n}} 
\prod_{q=1}^{\mathfrak{s}-2} \left(\dfrac{\vert \xi \! - \! \tau \vert}{\vert 
\xi \! - \! \alpha_{p_{q}} \vert \vert \tau \! - \! \alpha_{p_{q}} \vert} 
\right)^{\frac{\varkappa_{nk \tilde{k}_{q}}}{n}} \varpi (\xi) \varpi (\tau) 
\right)^{-1} \nonumber \\
\times& \, \left(\md \widetilde{\mu}^{f}_{\widetilde{V}}(\xi) \, \md 
\mu^{f}_{\widetilde{V}}(\tau) \! + \! \md \mu^{f}_{\widetilde{V}}(\xi) 
\, \md \widetilde{\mu}^{f}_{\widetilde{V}}(\tau) \right) \nonumber \\
=& \, 2 \iint_{\mathbb{R}^{2}} \ln \left(\vert \xi \! - \! \tau 
\vert^{\frac{\varkappa^{\infty}_{nk \tilde{k}_{\mathfrak{s}-1}}+1}{n}} 
\left(\dfrac{\vert \xi \! - \! \tau \vert}{\vert \xi \! - \! \alpha_{k} \vert 
\vert \tau \! - \! \alpha_{k} \vert} \right)^{\frac{\varkappa_{nk}-1}{n}} 
\prod_{q=1}^{\mathfrak{s}-2} \left(\dfrac{\vert \xi \! - \! \tau \vert}{\vert 
\xi \! - \! \alpha_{p_{q}} \vert \vert \tau \! - \! \alpha_{p_{q}} \vert} 
\right)^{\frac{\varkappa_{nk \tilde{k}_{q}}}{n}} \varpi (\xi) \varpi 
(\tau) \right)^{-1} \md \widetilde{\mu}^{f}_{\widetilde{V}}(\xi) 
\, \md \mu^{f}_{\widetilde{V}}(\tau) \nonumber \\
=& \, 2 \iint_{\mathbb{R}^{2}} \ln \left(\vert \xi \! - \! \tau 
\vert^{\frac{\varkappa^{\infty}_{nk \tilde{k}_{\mathfrak{s}-1}}+1}{n}} 
\left(\dfrac{\vert \xi \! - \! \tau \vert}{\vert \xi \! - \! \alpha_{k} \vert 
\vert \tau \! - \! \alpha_{k} \vert} \right)^{\frac{\varkappa_{nk}-1}{n}} 
\prod_{q=1}^{\mathfrak{s}-2} \left(\dfrac{\vert \xi \! - \! \tau \vert}{\vert 
\xi \! - \! \alpha_{p_{q}} \vert \vert \tau \! - \! \alpha_{p_{q}} \vert} 
\right)^{\frac{\varkappa_{nk \tilde{k}_{q}}}{n}} \varpi (\xi) \varpi 
(\tau) \right)^{-1} \md \mu^{f}_{\widetilde{V}}(\xi) \, \md 
\widetilde{\mu}^{f}_{\widetilde{V}}(\tau). \label{eql3.3a} 
\end{align}
Since, for $n \! \in \! \mathbb{N}$ and $k \! \in \! \lbrace 1,2,\dotsc,
K \rbrace$ such that $\alpha_{p_{\mathfrak{s}}} \! := \! \alpha_{k} 
\! \neq \! \infty$, both $\mathrm{I}^{f}_{\widetilde{V}}[\mu^{f}_{
\widetilde{V}}]$ and $\mathrm{I}^{f}_{\widetilde{V}}[\widetilde{
\mu}^{f}_{\widetilde{V}}]$ are bounded, it follows {}from 
inequality~\eqref{eql3.3a} that $\ln (\lvert \xi \! - \! \tau 
\rvert^{\frac{\varkappa^{\infty}_{nk \tilde{k}_{\mathfrak{s}-1}}+1}{n}}
(\lvert \xi - \tau \rvert/\lvert \xi -\alpha_{k} \rvert \lvert \tau 
-\alpha_{k} \rvert)^{\frac{\varkappa_{nk}-1}{n}} \prod_{q=1}^{
\mathfrak{s}-2}(\lvert \xi -\tau \rvert/\lvert \xi -\alpha_{p_{q}} 
\rvert \lvert \tau -\alpha_{p_{q}} \rvert)^{\frac{\varkappa_{nk 
\tilde{k}_{q}}}{n}})$ is integrable with respect to both the `product' 
measures $\md \widetilde{\mu}^{f}_{\widetilde{V}}(\xi) \, \md 
\mu^{f}_{\widetilde{V}}(\tau)$ and $\md \mu^{f}_{\widetilde{V}}(\xi) 
\, \md \widetilde{\mu}^{f}_{\widetilde{V}}(\tau)$; hence, {}from an 
argument on p.~149 of \cite{a51}, it follows that, for $n \! \in \! 
\mathbb{N}$ and $k \! \in \! \lbrace 1,2,\dotsc,K \rbrace$ such 
that $\alpha_{p_{\mathfrak{s}}} \! := \! \alpha_{k} \! \neq \! \infty$, 
$\ln (\lvert \xi \! - \! \tau \rvert^{\frac{\varkappa^{\infty}_{nk 
\tilde{k}_{\mathfrak{s}-1}}+1}{n}}(\lvert \xi - \tau \rvert/\lvert \xi -
\alpha_{k} \rvert \lvert \tau -\alpha_{k} \rvert)^{\frac{\varkappa_{nk}
-1}{n}} \prod_{q=1}^{\mathfrak{s}-2}(\lvert \xi -\tau \rvert/\lvert 
\xi -\alpha_{p_{q}} \rvert \lvert \tau -\alpha_{p_{q}} \rvert)^{\frac{
\varkappa_{nk \tilde{k}_{q}}}{n}})$ is integrable with respect to the 
associated one-parameter family of `product' measures $\md 
\mu^{f}_{\lambda}(\pmb{\cdot}) \, \md \mu^{f}_{\lambda}
(\pmb{\cdot})$, where
\begin{equation}
\mu^{f}_{\lambda}(z) \! := \! \mu^{f}_{\widetilde{V}}(z) \! + \! \lambda 
(\widetilde{\mu}^{f}_{\widetilde{V}} \! - \! \mu^{f}_{\widetilde{V}})(z), 
\quad (z,\lambda) \! \in \! \mathbb{R} \times [0,1]. 
\label{eql3.3d}
\end{equation}
For $n \! \in \! \mathbb{N}$ and $k \! \in \! \lbrace 1,2,\dotsc,K 
\rbrace$ such that $\alpha_{p_{\mathfrak{s}}} \! := \! \alpha_{k} 
\! \neq \! \infty$, let $\mathscr{F}_{\mu}^{f} \colon \mathbb{N} 
\times \lbrace 1,2,\dotsc,K \rbrace \times \mathscr{M}_{1}
(\mathbb{R}) \times [0,1] \! \ni \! (n,k,\mu^{f}_{\lambda}) 
\! \mapsto \! \mathrm{I}^{f}_{\widetilde{V}}[\mu^{f}_{\lambda}] \! 
=: \! \mathscr{F}_{\mu}^{f}(\lambda)$, where
\begin{equation}
\mathscr{F}_{\mu}^{f}(\lambda) \! = \! \iint_{\mathbb{R}^{2}} 
\ln \left(\vert \xi \! - \! \tau \vert^{\frac{\varkappa^{\infty}_{nk 
\tilde{k}_{\mathfrak{s}-1}}+1}{n}} \left(\dfrac{\vert \xi \! - \! \tau 
\vert}{\vert \xi \! - \! \alpha_{k} \vert \vert \tau \! - \! \alpha_{k} \vert} 
\right)^{\frac{\varkappa_{nk}-1}{n}} \prod_{q=1}^{\mathfrak{s}-2} \left(
\dfrac{\vert \xi \! - \! \tau \vert}{\vert \xi \! - \! \alpha_{p_{q}} \vert 
\vert \tau \! - \! \alpha_{p_{q}} \vert} \right)^{\frac{\varkappa_{nk 
\tilde{k}_{q}}}{n}} \varpi (\xi) \varpi (\tau) \right)^{-1} \md 
\mu^{f}_{\lambda}(\xi) \, \md \mu^{f}_{\lambda}(\tau). \label{eql3.3b}
\end{equation}
Noting {}from definition~\eqref{eql3.3d} that $\md \mu^{f}_{\lambda}(\xi) \, 
\md \mu^{f}_{\lambda}(\tau) \! = \! \md \mu^{f}_{\widetilde{V}}(\xi) \, \md 
\mu^{f}_{\widetilde{V}}(\tau) \! + \! \lambda \md \mu^{f}_{\widetilde{V}}(\xi) 
\, \md (\widetilde{\mu}^{f}_{\widetilde{V}} \! - \! \mu^{f}_{\widetilde{V}})
(\tau) \! + \! \lambda \md \mu^{f}_{\widetilde{V}}(\tau) \, \md 
(\widetilde{\mu}^{f}_{\widetilde{V}} \! - \! \mu^{f}_{\widetilde{V}})(\xi) 
\! + \! \lambda^{2} \md (\widetilde{\mu}^{f}_{\widetilde{V}} \! - \! 
\mu^{f}_{\widetilde{V}})(\xi) \, \md (\widetilde{\mu}^{f}_{\widetilde{V}} \! 
- \! \mu^{f}_{\widetilde{V}})(\tau)$, it follows {}from Lemma~\ref{lem3.1} 
and Equation~\eqref{eql3.3b} that, for $n \! \in \! \mathbb{N}$ and $k \! 
\in \! \lbrace 1,2,\dotsc,K \rbrace$ such that $\alpha_{p_{\mathfrak{s}}} 
\! := \! \alpha_{k} \! \neq \! \infty$,
\begin{align}
\mathscr{F}_{\mu}^{f}(\lambda) =& \, \mathrm{I}^{f}_{\widetilde{V}}
[\mu^{f}_{\widetilde{V}}] \! + \! 2 \lambda \iint_{\mathbb{R}^{2}} 
\ln \left(\vert \xi \! - \! \tau \vert^{\frac{\varkappa^{\infty}_{nk 
\tilde{k}_{\mathfrak{s}-1}}+1}{n}} \left(\dfrac{\vert \xi \! - \! \tau 
\vert}{\vert \xi \! - \! \alpha_{k} \vert \vert \tau \! - \! \alpha_{k} \vert} 
\right)^{\frac{\varkappa_{nk}-1}{n}} \prod_{q=1}^{\mathfrak{s}-2} \left(
\dfrac{\vert \xi \! - \! \tau \vert}{\vert \xi \! - \! \alpha_{p_{q}} \vert 
\vert \tau \! - \! \alpha_{p_{q}} \vert} \right)^{\frac{\varkappa_{nk 
\tilde{k}_{q}}}{n}} \varpi (\xi) \varpi (\tau) \right)^{-1} \nonumber \\
\times& \, \md \mu^{f}_{\widetilde{V}}(\xi) \, \md 
(\widetilde{\mu}^{f}_{\widetilde{V}} \! - \! \mu^{f}_{\widetilde{V}})
(\tau) \! + \! \lambda^{2} \iint_{\mathbb{R}^{2}} \ln \left(
\vert \xi \! - \! \tau \vert^{\frac{\varkappa^{\infty}_{nk 
\tilde{k}_{\mathfrak{s}-1}}+1}{n}} \left(\dfrac{\vert \xi \! - \! \tau 
\vert}{\vert \xi \! - \! \alpha_{k} \vert \vert \tau \! - \! \alpha_{k} \vert} 
\right)^{\frac{\varkappa_{nk}-1}{n}} \prod_{q=1}^{\mathfrak{s}-2} \left(
\dfrac{\vert \xi \! - \! \tau \vert}{\vert \xi \! - \! \alpha_{p_{q}} \vert 
\vert \tau \! - \! \alpha_{p_{q}} \vert} \right)^{\frac{\varkappa_{nk 
\tilde{k}_{q}}}{n}} \right. \nonumber \\
\times&\left. \, \varpi (\xi) \varpi (\tau) \vphantom{M^{M^{M^{M^{M^{M}}}}}}
\right)^{-1} \md (\widetilde{\mu}^{f}_{\widetilde{V}} \! - \! 
\mu^{f}_{\widetilde{V}})(\xi) \, \md (\widetilde{\mu}^{f}_{\widetilde{V}} 
\! - \! \mu^{f}_{\widetilde{V}})(\tau). \label{eql3.3c}
\end{align}
Since $\mu_{f}^{\sharp} \! := \! \widetilde{\mu}^{f}_{\widetilde{V}} 
\! - \! \mu^{f}_{\widetilde{V}}$ $(\in \! \mathscr{M}_{1}(\mathbb{R}))$ 
is a finite-moment signed measure on $\mathbb{R}$ with compact 
support and mean zero (that is, $\int_{\mathbb{R}} \md \mu_{f}^{\sharp}
(t) \! = \! \int_{\mathbb{R}} \md (\widetilde{\mu}^{f}_{\widetilde{V}} \! - \! 
\mu^{f}_{\widetilde{V}})(t) \! = \! 0)$, it follows {}from Equation~\eqref{eql3.3c} 
and the corresponding variational inequality in item~(ii) of Lemma~\ref{lem3.2} 
that, for $n \! \in \! \mathbb{N}$ and $k \! \in \! \lbrace 1,2,\dotsc,K \rbrace$ 
such that $\alpha_{p_{\mathfrak{s}}} \! := \! \alpha_{k} \! \neq \! \infty$, 
$\mathscr{F}_{\mu}^{f}(\lambda)$ is a twice-differentiable, real-valued 
function of $\lambda$ on $[0,1]$, with
\begin{equation*}
\dfrac{1}{2} \dfrac{\md^{2} \mathscr{F}_{\mu}^{f}(\lambda)}{\md 
\lambda^{2}} \! = \! \iint_{\mathbb{R}^{2}} \ln \left(\vert \xi \! - \! \tau 
\vert^{\frac{\varkappa^{\infty}_{nk \tilde{k}_{\mathfrak{s}-1}}+1}{n}} 
\left(\dfrac{\vert \xi \! - \! \tau \vert}{\vert \xi \! - \! \alpha_{k} \vert 
\vert \tau \! - \! \alpha_{k} \vert} \right)^{\frac{\varkappa_{nk}-1}{n}} 
\prod_{q=1}^{\mathfrak{s}-2} \left(\dfrac{\vert \xi \! - \! \tau \vert}{\vert 
\xi \! - \! \alpha_{p_{q}} \vert \vert \tau \! - \! \alpha_{p_{q}} \vert} 
\right)^{\frac{\varkappa_{nk \tilde{k}_{q}}}{n}} \varpi (\xi) \varpi (\tau) 
\right)^{-1} \md \mu_{f}^{\sharp}(\xi) \, \md \mu_{f}^{\sharp}(\tau) \! 
\geqslant \! 0,
\end{equation*}
that is, $\mathscr{F}_{\mu}^{f}(\lambda)$ is a real-valued convex function 
of $\lambda$ on $[0,1]$; hence, for $n \! \in \! \mathbb{N}$ and $k \! \in 
\! \lbrace 1,2,\dotsc,K \rbrace$ such that $\alpha_{p_{\mathfrak{s}}} \! 
:= \! \alpha_{k} \! \neq \! \infty$, via Equation~\eqref{eql3.3b} and 
Definition~\eqref{eql3.3d},
\begin{align*}
\mathrm{I}^{f}_{\widetilde{V}}[\mu^{f}_{\widetilde{V}}] \leqslant& \, 
\mathscr{F}_{\mu}^{f}(\lambda) \! = \! \mathrm{I}^{f}_{\widetilde{V}}
[\mu^{f}_{\lambda}] \! = \! \mathscr{F}_{\mu}^{f}(\lambda \! + \! 
(1 \! - \! \lambda)0) \! \leqslant \! \lambda \mathscr{F}_{\mu}^{f}(1) 
\! + \! (1 \! - \! \lambda) \mathscr{F}_{\mu}^{f}(0) \\
=& \, \lambda \mathrm{I}^{f}_{\widetilde{V}}[\widetilde{\mu}^{f}_{
\widetilde{V}}] \! + \! (1 \! - \! \lambda) \mathrm{I}^{f}_{\widetilde{V}}
[\mu^{f}_{\widetilde{V}}] \! = \! \lambda \mathrm{I}^{f}_{\widetilde{V}}
[\mu^{f}_{\widetilde{V}}] \! + \! (1 \! - \! \lambda) \mathrm{I}^{f}_{
\widetilde{V}}[\mu^{f}_{\widetilde{V}}] \quad \Rightarrow \quad 
\mathrm{I}^{f}_{\widetilde{V}}[\mu^{f}_{\widetilde{V}}] \! \leqslant 
\! \mathrm{I}^{f}_{\widetilde{V}}[\mu^{f}_{\lambda}] \! \leqslant 
\! \mathrm{I}^{f}_{\widetilde{V}}[\mu^{f}_{\widetilde{V}}],
\end{align*}
whence $\mathrm{I}^{f}_{\widetilde{V}}[\mu^{f}_{\lambda}] \! = \! 
\mathrm{I}^{f}_{\widetilde{V}}[\mu^{f}_{\widetilde{V}}] \! = \! E^{f}_{
\widetilde{V}}$. Since, for $n \! \in \! \mathbb{N}$ and $k \! \in \! 
\lbrace 1,2,\dotsc,K \rbrace$ such that $\alpha_{p_{\mathfrak{s}}} 
\! := \! \alpha_{k} \! \neq \! \infty$, $\mathrm{I}^{f}_{\widetilde{V}}
[\mu^{f}_{\lambda}] \! = \! \mathscr{F}_{\mu}^{f}(\lambda) \! = \! 
E^{f}_{\widetilde{V}}$, it follows, in particular, that $\left. \tfrac{
\md^{2} \mathscr{F}_{\mu}^{f}(\lambda)}{\md \lambda^{2}} 
\right\vert_{\lambda =0} \! = \! 0$, in which case, via 
Equation~\eqref{eql3.3c}, the corresponding variational inequality in 
item~(ii) of Lemma~\ref{lem3.2}, and the fact that $\mu_{f}^{\sharp}$ 
$(\in \! \mathscr{M}_{1}(\mathbb{R}))$ is a finite-moment signed 
measure on $\mathbb{R}$ with compact support and mean zero,
\begin{align*}
0 =& \, \iint_{\mathbb{R}^{2}} \ln \left(\vert \xi \! - \! \tau 
\vert^{\frac{\varkappa^{\infty}_{nk \tilde{k}_{\mathfrak{s}-1}}+1}{n}} 
\left(\dfrac{\vert \xi \! - \! \tau \vert}{\vert \xi \! - \! \alpha_{k} \vert 
\vert \tau \! - \! \alpha_{k} \vert} \right)^{\frac{\varkappa_{nk}-1}{n}} 
\prod_{q=1}^{\mathfrak{s}-2} \left(\dfrac{\vert \xi \! - \! \tau \vert}{\vert 
\xi \! - \! \alpha_{p_{q}} \vert \vert \tau \! - \! \alpha_{p_{q}} \vert} 
\right)^{\frac{\varkappa_{nk \tilde{k}_{q}}}{n}} \right)^{-1} \md 
(\widetilde{\mu}^{f}_{\widetilde{V}} \! - \! \mu^{f}_{\widetilde{V}})(\xi) 
\, \md (\widetilde{\mu}^{f}_{\widetilde{V}} \! - \! \mu^{f}_{\widetilde{V}})
(\tau) \\
+& \, 2 \left(\int_{\mathbb{R}} \widetilde{V}(\xi) \, \md 
(\widetilde{\mu}^{f}_{\widetilde{V}} \! - \! \mu^{f}_{\widetilde{V}})(\xi) 
\right) \underbrace{\int_{\mathbb{R}} \md (\widetilde{\mu}^{f}_{
\widetilde{V}} \! - \! \mu^{f}_{\widetilde{V}})(\tau)}_{= \, 0} \\
=& \, \iint_{\mathbb{R}^{2}} \ln \left(\vert \xi \! - \! \tau 
\vert^{\frac{\varkappa^{\infty}_{nk \tilde{k}_{\mathfrak{s}-1}}+1}{n}} 
\left(\dfrac{\vert \xi \! - \! \tau \vert}{\vert \xi \! - \! \alpha_{k} \vert 
\vert \tau \! - \! \alpha_{k} \vert} \right)^{\frac{\varkappa_{nk}-1}{n}} 
\prod_{q=1}^{\mathfrak{s}-2} \left(\dfrac{\vert \xi \! - \! \tau \vert}{
\vert \xi \! - \! \alpha_{p_{q}} \vert \vert \tau \! - \! \alpha_{p_{q}} 
\vert} \right)^{\frac{\varkappa_{nk \tilde{k}_{q}}}{n}} \right)^{-1} \md 
(\widetilde{\mu}^{f}_{\widetilde{V}} \! - \! \mu^{f}_{\widetilde{V}})(\xi) 
\, \md (\widetilde{\mu}^{f}_{\widetilde{V}} \! - \! \mu^{f}_{\widetilde{V}})
(\tau) \\
=& \, \left(\dfrac{(n \! - \! 1)K \! + \! k}{n} \right) \int_{0}^{+\infty}
t^{-1} \lvert (\widehat{\widetilde{\mu}^{f}_{\widetilde{V}}} \! - 
\! \widehat{\mu^{f}_{\widetilde{V}}})(t) \rvert^{2} \, \md t \! 
\geqslant \! 0,
\end{align*}
whence
\begin{equation*}
\int_{0}^{+\infty}t^{-1} \lvert (\widehat{\widetilde{\mu}^{f}_{
\widetilde{V}}} \! - \! \widehat{\mu^{f}_{\widetilde{V}}})(t) 
\rvert^{2} \, \md t \! = \! 0 \quad \Rightarrow
\end{equation*}
$\widehat{\widetilde{\mu}^{f}_{\widetilde{V}}}(t) \! = \! \widehat{\mu^{f}_{
\widetilde{V}}}(t)$, $t \! \geqslant \! 0$. Noting that $\widehat{\widetilde{
\mu}^{f}_{\widetilde{V}}}(-t) \! = \! \int_{\mathbb{R}} \me^{-\mi t \zeta} 
\, \md \widetilde{\mu}^{f}_{\widetilde{V}}(\zeta) \! = \! 
\overline{\widehat{\widetilde{\mu}^{f}_{\widetilde{V}}}(t)}$ and 
$\widehat{\mu^{f}_{\widetilde{V}}}(-t) \! = \! \int_{\mathbb{R}} 
\me^{-\mi t \zeta} \, \md \mu^{f}_{\widetilde{V}}(\zeta) \! = \! \overline{
\widehat{\mu^{f}_{\widetilde{V}}}(t)}$, it follows {}from the relation 
$\widehat{\widetilde{\mu}^{f}_{\widetilde{V}}}(t) \! = \! \widehat{\mu^{f}_{
\widetilde{V}}}(t)$, $t \! \geqslant \! 0$, and a complex-conjugation 
argument, that $\widehat{\widetilde{\mu}^{f}_{\widetilde{V}}}(-t) \! = \! 
\widehat{\mu^{f}_{\widetilde{V}}}(-t)$, $t \! \geqslant \! 0$; hence, for 
$n \! \in \! \mathbb{N}$ and $k \! \in \! \lbrace 1,2,\dotsc,K \rbrace$ 
such that $\alpha_{p_{\mathfrak{s}}} \! := \! \alpha_{k} \! \neq \! \infty$, 
$\widehat{\widetilde{\mu}^{f}_{\widetilde{V}}}(t) \! = \! \widehat{
\mu^{f}_{\widetilde{V}}}(t)$, $t \! \in \! \mathbb{R}$, which shows that 
$\int_{\mathbb{R}} \me^{\mi t \zeta} \, \md (\widetilde{\mu}^{f}_{
\widetilde{V}} \! - \! \mu^{f}_{\widetilde{V}})(\zeta) \! = \! 0$ $\Rightarrow$ 
$\widetilde{\mu}^{f}_{\widetilde{V}} \! = \! \mu^{f}_{\widetilde{V}}$: 
thus the uniqueness of the associated equilibrium measure.

$\pmb{(2)}$ The proof of this case, that is, $n \! \in \! \mathbb{N}$ and 
$k \! \in \! \lbrace 1,2,\dotsc,K \rbrace$ such that $\alpha_{p_{\mathfrak{s}}} 
\! := \! \alpha_{k} \! = \! \infty$, is virtually identical to the proof presented 
in $\pmb{(1)}$ above; one mimics, \emph{verbatim}, the scheme of the calculations 
presented in case $\pmb{(1)}$ in order to arrive at the corresponding uniqueness 
claim for the associated equilibrium measure $\mu_{\widetilde{V}}^{\infty}$ 
stated in the lemma; in order to do so, however, one needs the analogues of 
Definition~\eqref{eql3.3d} and Equation~\eqref{eql3.3b}, which, in the present 
case, read: for $\widetilde{\mu}^{\infty}_{\widetilde{V}} \! \in \! \mathscr{M}_{1}
(\mathbb{R})$, a second probability measure for which $\mathrm{I}^{\infty}_{
\widetilde{V}}[\widetilde{\mu}^{\infty}_{\widetilde{V}}] \! = \! E^{\infty}_{
\widetilde{V}} \! = \! \mathrm{I}^{\infty}_{\widetilde{V}}[\mu^{\infty}_{
\widetilde{V}}]$,
\begin{equation}
\mu^{\infty}_{\lambda}(z) \! := \! \mu^{\infty}_{\widetilde{V}}(z) \! + \! 
\lambda (\widetilde{\mu}^{\infty}_{\widetilde{V}} \! - \! \mu^{\infty}_{
\widetilde{V}})(z), \quad (z,\lambda) \! \in \! \mathbb{R} \times [0,1]. 
\label{eql3.3e}
\end{equation}
and
\begin{equation}
\mathscr{F}_{\mu}^{\infty}(\lambda) \! = \! \iint_{\mathbb{R}^{2}} 
\ln \left(\vert \xi \! - \! \tau \vert^{\frac{\varkappa_{nk}}{n}} 
\prod_{q=1}^{\mathfrak{s}-1} \left(\dfrac{\vert \xi \! - \! \tau \vert}{\vert 
\xi \! - \! \alpha_{p_{q}} \vert \vert \tau \! - \! \alpha_{p_{q}} \vert} 
\right)^{\frac{\varkappa_{nk \tilde{k}_{q}}}{n}} \varpi (\xi) \varpi (\tau) 
\right)^{-1} \md \mu^{\infty}_{\lambda}(\xi) \, \md \mu^{\infty}_{\lambda}
(\tau). \label{eql3.3f}
\end{equation}
This concludes the proof. \hfill $\qed$

For $n \! \in \! \mathbb{N}$ and $k \! \in \! \lbrace 1,2,\dotsc,K \rbrace$ 
such that $\alpha_{p_{\mathfrak{s}}} \! := \! \alpha_{k} \! = \! \infty$ 
or $\alpha_{p_{\mathfrak{s}}} \! := \! \alpha_{k} \! \neq \! \infty$, 
Lemma~$\bm{\mathrm{RHP}_{\mathrm{MPC}}}$ is now reformulated 
as two families of equivalent, auxiliary matrix RHPs; before doing so, 
however, the following motivational preamble warrants consideration.
\begin{eeeee} \label{rem1.3.5} 
\textsl{In an attempt to motivate the derivation of the formulae for the 
associated complex logarithmic potentials $(g$-functions$)$ defined 
by (see Lemma~\ref{lem3.4} below) Equations~\eqref{eql3.4gee1} 
and~\eqref{eql3.4gee3}, heuristic calculations will now be presented$:$ 
as an illustration, consider, without loss of generality, the case $(n,k) 
\! \in \! \mathbb{N} \times \lbrace 1,2,\dotsc,K \rbrace$ such that 
$\alpha_{p_{\mathfrak{s}}} \! := \! \alpha_{k} \! \neq \! \infty$ (the 
calculation for the case $\alpha_{p_{\mathfrak{s}}} \! := \! \alpha_{k} 
\! = \! \infty$ is analogous). For $(n,k) \! \in \! \mathbb{N} \times 
\lbrace 1,2,\dotsc,K \rbrace$ such that $\alpha_{p_{\mathfrak{s}}} \! 
:= \!\alpha_{k} \! \neq \! \infty$, with $\mathcal{N} \! := \! (n \! - \! 
1)K \! + \! k$, denote by (cf. Remark~\ref{rem1.2.4}$)$ $\lbrace 
\tilde{\mathfrak{z}}^{n}_{k}(1),\tilde{\mathfrak{z}}^{n}_{k}(2),\dotsc,
\tilde{\mathfrak{z}}^{n}_{k}(\mathcal{N}) \rbrace$ $(\subset 
\overline{\mathbb{R}} \setminus \lbrace \alpha_{1},\alpha_{2},
\dotsc,\alpha_{K} \rbrace)$, with, say, the enumeration 
$\tilde{\mathfrak{z}}^{n}_{k}(1) \! < \! \tilde{\mathfrak{z}}^{n}_{k}(2) \! 
< \! \dotsb \! < \! \tilde{\mathfrak{z}}^{n}_{k}(\mathcal{N})$, the set 
of zeros (counting multiplicities; see the corresponding items of 
Lemmata~\ref{lemrootz} and~\ref{lemetatomu} below) of the 
{\rm MPC ORF}, $\phi^{n}_{k}(z)$, (also the set of zeros (counting 
multiplicities) of the monic {\rm MPC ORF}, $\pmb{\pi}^{n}_{k}(z))$ and 
define by $\tilde{\eta}_{\tilde{\mathfrak{z}}}(z) \! := \! \mathcal{N}^{-1} 
\sum_{j=1}^{\mathcal{N}} \delta_{\tilde{\mathfrak{z}}^{n}_{k}(j)}(z)$ the 
associated normalised $(\int_{\mathbb{R}} \md \tilde{\eta}_{\tilde{\mathfrak{z}}}
(\xi) \! = \! 1)$ zero counting measure, where $\delta_{\tilde{\mathfrak{z}}^{n}_{k}
(j)}(z)$, $j \! = \! 1,2,\dotsc,\mathcal{N}$, denotes the Dirac (atomic) mass 
concentrated at $\tilde{\mathfrak{z}}^{n}_{k}(j)$. For $(n,k) \! \in \! \mathbb{N} 
\times \lbrace 1,2,\dotsc,K \rbrace$ such that $\alpha_{p_{\mathfrak{s}}} 
\! := \! \alpha_{k} \! \neq \! \infty$, write $\phi^{n}_{k}(z) \! = \! \tilde{c}
(n,k,z_{o};\vec{\pmb{\alpha}}) \prod_{j=1}^{\mathcal{N}}(z \! - \! 
\tilde{\mathfrak{z}}^{n}_{k}(j)) \prod_{q=1}^{\mathfrak{s}-2}(z \! - \! 
\alpha_{p_{q}})^{-\varkappa_{nk \tilde{k}_{q}}}(z \! - \! \alpha_{k})^{-
\varkappa_{nk}}$, where $\tilde{c}(n,k,z_{o};\vec{\pmb{\alpha}})$ denotes 
an, in general, bounded $n$-, $k$-, $z_{o}$-, and $\lbrace \alpha_{1},
\alpha_{2},\dotsc,\alpha_{K} \rbrace$-dependent constant. One 
manipulates this expression thus: $(z \! - \! \alpha_{k}) \phi^{n}_{k}
(z) \! = \! \tilde{c}(n,k,z_{o};\vec{\pmb{\alpha}}) \exp (\mathcal{N} 
\int_{\mathbb{R}} \ln (z \! - \! \xi) \, \md \tilde{\eta}_{\tilde{
\mathfrak{z}}}(\xi)) \exp (-n \int_{\mathbb{R}} \ln (\prod_{q=1}^{
\mathfrak{s}-2}(z \! - \! \alpha_{p_{q}})^{\frac{\varkappa_{nk 
\tilde{k}_{q}}}{n}}(z \! - \! \alpha_{k})^{\frac{\varkappa_{nk}-1}{n}}) \, 
\md \tilde{\eta}_{\tilde{\mathfrak{z}}}(\xi))$, which, upon noting 
the `splitting' (cf. Equation~\eqref{fincount}$)$ $\mathcal{N} \! 
:= \! (n \! - \! 1)K \! + \! k \! = \! \sum_{q=1}^{\mathfrak{s}-2} 
\varkappa_{nk \tilde{k}_{q}} \! + \! (\varkappa^{\infty}_{nk \tilde{k}_{
\mathfrak{s}-1}} \! + \! 1) \! + \! (\varkappa_{nk} \! - \! 1)$, can be 
presented as $(z \! - \! \alpha_{k}) \phi^{n}_{k}(z) \! = \! \tilde{c}
(n,k,z_{o};\vec{\pmb{\alpha}}) \exp (n(\sum_{q=1}^{\mathfrak{s}-2} 
\tfrac{\varkappa_{nk \tilde{k}_{q}}}{n} \! + \! \tfrac{\varkappa^{\infty}_{nk 
\tilde{k}_{\mathfrak{s}-1}}+1}{n} \! + \! \tfrac{\varkappa_{nk}-1}{n}) 
\int_{\mathbb{R}} \ln (z \! - \! \xi) \, \md \tilde{\eta}_{\tilde{\mathfrak{z}}}
(\xi)) \exp (-n \int_{\mathbb{R}} \ln (\prod_{q=1}^{\mathfrak{s}-2}
(z \! - \! \alpha_{p_{q}})^{\frac{\varkappa_{nk \tilde{k}_{q}}}{n}}(z \! 
- \! \alpha_{k})^{\frac{\varkappa_{nk}-1}{n}}) \, \md \tilde{\eta}_{
\tilde{\mathfrak{z}}}(\xi))$, whence $(z \! - \! \alpha_{k}) \phi^{n}_{k}
(z) \! = \! \tilde{c}(n,k,z_{o};\vec{\pmb{\alpha}}) \exp (n \int_{\mathbb{R}} 
\ln ((z \! - \! \xi)^{\frac{\varkappa^{\infty}_{nk \tilde{k}_{\mathfrak{s}-1}}
+1}{n}}(z \! - \! \xi)^{\frac{\varkappa_{nk}-1}{n}} \prod_{q=1}^{\mathfrak{s}
-2}(z \! - \! \xi)^{\frac{\varkappa_{nk \tilde{k}_{q}}}{n}}) \, \md 
\tilde{\eta}_{\tilde{\mathfrak{z}}}(\xi)) \exp (-n \int_{\mathbb{R}} \ln 
(\prod_{q=1}^{\mathfrak{s}-2}(z \! - \! \alpha_{p_{q}})^{\frac{\varkappa_{nk 
\tilde{k}_{q}}}{n}}(z \! - \! \alpha_{k})^{\frac{\varkappa_{nk}-1}{n}}) \, \md 
\tilde{\eta}_{\tilde{\mathfrak{z}}}(\xi))$. As shown in the corresponding 
items of Lemmata~\ref{lemrootz} and~\ref{lemetatomu} below (see, 
also, Lemmata~\ref{lem3.5} and~\ref{lem3.6} below), $\tilde{\eta}_{
\tilde{\mathfrak{z}}}$ converges, in the weak-$\ast$ topology of 
measures, to the associated equilibrium measure, $\mu_{\widetilde{V}}^{f}$, 
that is, $\tilde{\eta}_{\tilde{\mathfrak{z}}} \! \overset{\ast}{\to} \! 
\mu_{\widetilde{V}}^{f}$, in the double-scaling limit $\mathscr{N},n 
\! \to \! \infty$ such that $z_{o} \! = \! 1 \! + \! o(1)$$;$ therefore, 
{}from the above calculations, recalling that $(z \! - \! \alpha_{k}) 
\phi^{n}_{k}(z) \! =_{z \to \alpha_{k}} \! \mu^{f}_{n,\varkappa_{nk}}(n,k)
(z \! - \! \alpha_{k})^{-(\varkappa_{nk}-1)}(1 \! + \! \mathcal{O}(z \! - \! 
\alpha_{k}))$, assuming that all limits (see below) commute, and using the 
fact that the complex logarithmic potential $(g$-function) can be specified 
up to a $z$-independent factor (so chosen in order to guarantee that 
$\pmb{\pi}^{n}_{k}(z)$ is monic) which can depend on $n$, $k$, 
$z_{o}$, and $\lbrace \alpha_{1},\alpha_{2},\dotsc,\alpha_{K} \rbrace$, 
via this latter degree of freedom in `fixing' the complex logarithmic 
potential in conjunction with the comments below, 
one can expect that, in the double-scaling limit $\mathscr{N},n \! \to \! 
\infty$ such that $z_{o} \! = \! 1 \! + \! o(1)$, for non-real $z \! \in \! 
\overline{\mathbb{C}} \setminus \lbrace \alpha_{1},\alpha_{2},\dotsc,
\alpha_{K} \rbrace$, $(z \! - \! \alpha_{k}) \phi^{n}_{k}(z) \! 
\sim_{z \to \alpha_{k}} \! \exp (n \mathfrak{g}(z))$,\footnote{Here and 
in the sequel, the notation $A \! \sim \! B$ is used to denote the fact 
that the quotient $A/B$ is fixed up to some $n$-, $k$-, $z_{o}$-, and 
$\lbrace \alpha_{1},\alpha_{2},\dotsc,\alpha_{K} \rbrace$-dependent 
factor.} where $\mathfrak{g}(z) \! := \! \int_{J_{f}} \ln ((z \! - \! \xi)^{
\frac{\varkappa^{\infty}_{nk \tilde{k}_{\mathfrak{s}-1}}+1}{n}}((z \! - 
\! \xi)/(z \! - \! \alpha_{k}))^{\frac{\varkappa_{nk}-1}{n}} \prod_{q=
1}^{\mathfrak{s}-2}((z \! - \! \xi)/(z \! - \! \alpha_{p_{q}}))^{\frac{
\varkappa_{nk \tilde{k}_{q}}}{n}}) \, \md \mu_{\widetilde{V}}^{f}(\xi)$ 
gives rise to, but does not yet fix (see the comments below), 
the associated complex logarithmic potential $g^{f}(z)$ 
given by Equation~\eqref{eql3.4gee3}.\footnote{Note: this is actually 
a family of $g$-functions, triply-indexed by $n$, $k$, and $z_{o}$.} 
Once the matrix {\rm RHP} for the associated monic {\rm MPC ORF} (cf. the 
corresponding item of Lemma~$\bm{\mathrm{RHP}_{\mathrm{MPC}}})$ 
and the associated $g$-function, whose properties are determined at 
a later stage, are available, one makes, after a careful analysis of the 
asymptotic behaviour of the integral representation of the `solution 
matrix' of the associated matrix {\rm RHP} (with the Cauchy kernel 
normalised at the pole $\alpha_{p_{\mathfrak{s}}} \! := \! \alpha_{k} 
\! \neq \! \infty)$ in an open neighbourhood of the pole 
$\alpha_{p_{\mathfrak{s}}} \! := \! \alpha_{k} \! \neq \! \infty$ 
(cf. Equation~\eqref{intrepfin}$)$, a transformation via this 
$g$-function, the result of which is a `new' matrix {\rm RHP} for 
the monic {\rm MPC ORF}. One then asks if there exists a choice 
of this $g$-function (which will also fix the indeterminacy for 
$\mathfrak{g}(z)$ above, so that one may arrive at the definition of 
$g^{f}(z))$ so that the `new', or `transformed', matrix {\rm RHP} for 
the monic {\rm MPC ORF} is in a form suitable for asymptotic analysis, 
in the double-scaling limit $\mathscr{N},n \! \to \! \infty$ such that 
$z_{o} \! = \! 1 \! + \! o(1)$. One then arrives at a collection of 
equations and inequalities, which, if satisfied, achieve the sought-after 
result of obtaining a matrix {\rm RHP} for the monic {\rm MPC ORF} 
having `canonical normalisation' (at the pole $\alpha_{p_{\mathfrak{s}}} 
\! := \! \alpha_{k} \! \neq \! \infty)$. These equations and inequalities, 
in turn, are shown to be equivalent to Euler-Lagrange variational 
conditions associated with the corresponding energy minimisation 
problems (see the corresponding items in Lemma~\ref{lem3.8} 
below).}
\end{eeeee}
\begin{eeeee} \label{remm3.4} 
\textsl{For completeness, the integrands appearing in the definitions of 
$g^{\infty}(z)$ and $g^{f}(z)$ given by Equations~\eqref{eql3.4gee1} 
and~\eqref{eql3.4gee3}, respectively, are defined as follows: for $n \! 
\in \! \mathbb{N}$ and $k \! \in \! \lbrace 1,2,\dotsc,K \rbrace$ such 
that $\alpha_{p_{\mathfrak{s}}} \! := \! \alpha_{k} \! = \! \infty$,
\begin{equation*}
\ln \left((z \! - \! \xi)^{\frac{\varkappa_{nk}}{n}} \prod_{q=1}^{
\mathfrak{s}-1} \left(\dfrac{(z \! - \! \xi)}{(z \! - \! \alpha_{p_{q}})
(\xi \! - \! \alpha_{p_{q}})} \right)^{\frac{\varkappa_{nk \tilde{k}_{q}}}{n}} 
\right) \! := \! \dfrac{\varkappa_{nk}}{n} \ln (z \! - \! \xi) \! + \! 
\sum_{q=1}^{\mathfrak{s}-1} \dfrac{\varkappa_{nk \tilde{k}_{q}}}{n} 
\left(\ln (z \! - \! \xi) \! - \! \ln (z \! - \! \alpha_{p_{q}}) \! - \! \ln 
(\xi \! - \! \alpha_{p_{q}}) \right),
\end{equation*}
and, for $n \! \in \! \mathbb{N}$ and $k \! \in \! \lbrace 1,2,\dotsc,
K \rbrace$ such that $\alpha_{p_{\mathfrak{s}}} \! := \! \alpha_{k} 
\! \neq \! \infty$,
\begin{align*}
& \ln \left((z \! - \! \xi)^{\frac{\varkappa_{nk \tilde{k}_{\mathfrak{s}-
1}}^{\infty}+1}{n}} \left(\dfrac{(z \! - \! \xi)}{(z \! - \! \alpha_{k})
(\xi \! - \! \alpha_{k})} \right)^{\frac{\varkappa_{nk}-1}{n}} \prod_{q=
1}^{\mathfrak{s}-2} \left(\dfrac{(z \! - \! \xi)}{(z \! - \! \alpha_{p_{q}})
(\xi \! - \! \alpha_{p_{q}})} \right)^{\frac{\varkappa_{nk \tilde{k}_{q}}}{n}} 
\right) \! := \! \left(\dfrac{\varkappa_{nk \tilde{k}_{\mathfrak{s}-1}}^{
\infty} \! + \! 1}{n} \right) \ln (z \! - \! \xi) \\
&+ \, \left(\dfrac{\varkappa_{nk} \! - \! 1}{n} \right) \left(\ln (z \! - \! \xi) 
\! - \! \ln (z \! - \! \alpha_{k}) \! - \! \ln (\xi \! - \! \alpha_{k}) \right) \! 
+ \! \sum_{q=1}^{\mathfrak{s}-2} \dfrac{\varkappa_{nk \tilde{k}_{q}}}{n} 
\left(\ln (z \! - \! \xi) \! - \! \ln (z \! - \! \alpha_{p_{q}}) \! - \! \ln (\xi 
\! - \! \alpha_{p_{q}}) \right),
\end{align*}
where, for any (generic) constant $c \! < \! 0$, $\ln c \! := \! \ln 
\lvert c \rvert \! + \! \mi \pi$.}
\end{eeeee}
\begin{ccccc} \label{lem3.4} 
Let the external field $\widetilde{V} \colon \overline{\mathbb{R}} 
\setminus \lbrace \alpha_{1},\alpha_{2},\dotsc,\alpha_{K} \rbrace \! 
\to \! \mathbb{R}$ satisfy conditions~\eqref{eq20}--\eqref{eq22}. 
For $n \! \in \! \mathbb{N}$ and $k \! \in \! \lbrace 1,2,\dotsc,K 
\rbrace$ such that $\alpha_{p_{\mathfrak{s}}} \! := \! \alpha_{k} \! 
= \! \infty$ (resp., $\alpha_{p_{\mathfrak{s}}} \! := \!\alpha_{k} 
\! \neq \! \infty)$, set, for the associated equilibrium measure 
$\mu_{\widetilde{V}}^{\infty}$ (resp., $\mu_{\widetilde{V}}^{f})$, 
$J_{\infty} \! := \! \supp (\mu_{\widetilde{V}}^{\infty})$ (resp., $J_{f} \! 
:= \! \supp (\mu_{\widetilde{V}}^{f}))$, where $J_{\infty}$ (resp., $J_{f})$ 
is a proper compact measurable subset of $\overline{\mathbb{R}} 
\setminus \lbrace \alpha_{1},\alpha_{2},\dotsc,\alpha_{K} \rbrace$. 
For $n \! \in \! \mathbb{N}$ and $k \! \in \! \lbrace 1,2,\dotsc,
K \rbrace$, let $\mathcal{X} \colon \mathbb{N} \times \lbrace 
1,2,\dotsc,K \rbrace \times \overline{\mathbb{C}} \setminus 
\overline{\mathbb{R}} \! \to \! \operatorname{SL}_{2}(\mathbb{C})$ 
be the unique solution of the monic {\rm MPC} {\rm ORF} {\rm RHP} 
stated in Lemma~{\rm $\bm{\mathrm{RHP}_{\mathrm{MPC}}}$}. For 
$n \! \in \! \mathbb{N}$ and $k \! \in \! \lbrace 1,2,\dotsc,K \rbrace$ 
such that $\alpha_{p_{\mathfrak{s}}} \! := \! \alpha_{k} \! = \! \infty$, 
set
\begin{equation} \label{eql3.4e} 
\mathscr{X}(n,k,z) \! = \! \mathscr{X}(z) \! := \! 
\me^{-\frac{n \hat{\ell}}{2} \ad (\sigma_{3})} \mathcal{X}(z) 
\me^{-n(g^{\infty}(z)-\tilde{\mathscr{P}}_{0}(n,k)) \sigma_{3}},
\end{equation}
where $\hat{\ell} \colon \mathbb{N} \times \lbrace 1,2,\dotsc,K 
\rbrace \! \to \! \mathbb{R}$, the associated variational constant, 
is given in the corresponding item of Lemma \ref{lem3.8} below, 
$g^{\infty} \colon \mathbb{N} \times \lbrace 1,2,\dotsc,K 
\rbrace \times \mathbb{C} \setminus (-\infty,\max \lbrace 
\max_{q=1,2,\dotsc,\mathfrak{s}-1} \lbrace \alpha_{p_{q}} \rbrace,
\max \lbrace J_{\infty} \rbrace \rbrace) \! \to \! \mathbb{C}$, 
$(n,k,z) \! \mapsto \! g^{\infty}(n,k,z) \! = \! g^{\infty}(z)$, with
\begin{equation} \label{eql3.4gee1} 
g^{\infty}(z) \! := \! \int_{J_{\infty}} \ln \left((z \! - \! \xi)^{\frac{
\varkappa_{nk}}{n}} \prod_{q=1}^{\mathfrak{s}-1} \left(\dfrac{(z \! - 
\! \xi)}{(z \! - \! \alpha_{p_{q}})(\xi \! - \! \alpha_{p_{q}})} 
\right)^{\frac{\varkappa_{nk \tilde{k}_{q}}}{n}} \right) \md 
\mu_{\widetilde{V}}^{\infty}(\xi),
\end{equation}
and
\begin{equation} \label{eql3.4gee2} 
\tilde{\mathscr{P}}_{0}(n,k) \! = \! \tilde{\mathscr{P}}_{0} \! := \! 
-\sum_{q=1}^{\mathfrak{s}-1} \dfrac{\varkappa_{nk \tilde{k}_{q}}}{n} 
\int_{J_{\infty}} \ln (\lvert \xi \! - \! \alpha_{p_{q}} \rvert) \, \md \mu_{
\widetilde{V}}^{\infty}(\xi) \! - \! \mi \pi \sum_{q=1}^{\mathfrak{s}-1} 
\dfrac{\varkappa_{nk \tilde{k}_{q}}}{n} \int_{J_{\infty} \cap \mathbb{R}_{
\alpha_{p_{q}}}^{<}} \, \md \mu_{\widetilde{V}}^{\infty}(\xi),
\end{equation}
and,\footnote{If $J_{\infty} \cap \mathbb{R}^{<}_{\alpha_{p_{q}}} \! = \! 
\varnothing$, then $\int_{J_{\infty} \cap \mathbb{R}^{<}_{\alpha_{p_{q}}}} 
\md \mu_{\widetilde{V}}^{\infty}(\xi) \! := \! 0$, $q \! \in \! \lbrace 1,2,
\dotsc,\mathfrak{s} \! - \! 1 \rbrace$.} for $n \! \in \! \mathbb{N}$ and 
$k \! \in \! \lbrace 1,2,\dotsc,K \rbrace$ such that $\alpha_{p_{\mathfrak{s}}} 
\! := \! \alpha_{k} \! \neq \! \infty$, set
\begin{equation} \label{eql3.4f} 
\mathscr{X}(n,k,z) \! = \! \mathscr{X}(z) \! := \! \me^{-\frac{n 
\tilde{\ell}}{2} \ad (\sigma_{3})} \mathcal{X}(z) \me^{-n(g^{f}(z)- 
\hat{\mathscr{P}}_{0}(n,k)) \sigma_{3}},
\end{equation}
where $\tilde{\ell} \colon \mathbb{N} \times \lbrace 1,2,\dotsc,K \rbrace 
\! \to \! \mathbb{R}$, the associated variational constant, is given in the 
corresponding item of Lemma \ref{lem3.8} below, $g^{f} \colon \mathbb{N} 
\times \lbrace 1,2,\dotsc,K \rbrace \times \mathbb{C} \setminus (-\infty,
\max \lbrace \max_{q=1,\dotsc,\mathfrak{s}-2,\mathfrak{s}} \lbrace 
\alpha_{p_{q}} \rbrace,\max \lbrace J_{f} \rbrace \rbrace) \! \to \! 
\mathbb{C}$, $(n,k,z) \! \mapsto \! g^{f}(n,k,z) \! = \! g^{f}(z)$, with
\begin{align} \label{eql3.4gee3} 
g^{f}(z) \! := \! \int_{J_{f}} \ln \left((z \! - \! \xi)^{\frac{\varkappa_{nk 
\tilde{k}_{\mathfrak{s}-1}}^{\infty}+1}{n}} \left(\dfrac{(z \! - \! \xi)}{(z 
\! - \! \alpha_{k})(\xi \! - \! \alpha_{k})} \right)^{\frac{\varkappa_{nk}-
1}{n}} \prod_{q=1}^{\mathfrak{s}-2} \left(\dfrac{(z \! - \! \xi)}{(z \! - \! 
\alpha_{p_{q}})(\xi \! - \! \alpha_{p_{q}})} \right)^{\frac{\varkappa_{nk 
\tilde{k}_{q}}}{n}} \right) \md \mu_{\widetilde{V}}^{f}(\xi),
\end{align}
and
\begin{equation} \label{eql3.4gee4} 
\hat{\mathscr{P}}_{0}(n,k) \! = \! 
\begin{cases}
\hat{\mathscr{P}}_{0}^{+}, &\text{$z \! \in \! \mathbb{C}_{+}$,} \\
\hat{\mathscr{P}}_{0}^{-}, &\text{$z \! \in \! \mathbb{C}_{-}$,}
\end{cases}
\end{equation}
with
\begin{align} \label{eql3.4gee5} 
\hat{\mathscr{P}}_{0}^{\pm} :=& \, \dfrac{1}{n} \int_{J_{f}} \ln \left(\lvert \xi \! 
- \! \alpha_{k} \rvert^{\varkappa_{nk \tilde{k}_{\mathfrak{s}-1}}^{\infty}+1} 
\prod_{q=1}^{\mathfrak{s}-2} \left(\dfrac{\lvert \xi \! - \! \alpha_{k} 
\rvert}{\lvert \xi \! - \! \alpha_{p_{q}} \rvert \lvert \alpha_{p_{q}} \! - \! \alpha_{k} 
\rvert} \right)^{\varkappa_{nk \tilde{k}_{q}}} \right) \md \mu_{\widetilde{V}}^{f}
(\xi) \! - \! \mi \pi \left(\dfrac{\varkappa_{nk} \! - \! 1}{n} \right) \int_{J_{f} \cap 
\mathbb{R}_{\alpha_{k}}^{<}} \md \mu_{\widetilde{V}}^{f}(\xi) \nonumber \\
-& \, \mi \pi \sum_{q=1}^{\mathfrak{s}-2} \dfrac{\varkappa_{nk 
\tilde{k}_{q}}}{n} \int_{J_{f} \cap \mathbb{R}_{\alpha_{p_{q}}}^{<}} \md 
\mu_{\widetilde{V}}^{f}(\xi) \! \pm \! \mi \pi \left(\dfrac{(n \! - \! 1)K 
\! + \! k}{n} \right) \int_{J_{f} \cap \mathbb{R}_{\alpha_{k}}^{>}} \md 
\mu_{\widetilde{V}}^{f}(\xi) \! \mp \! \mi \pi \sum_{j \in \hat{\Delta}_{f}
(k)} \dfrac{\varkappa_{nk \tilde{k}_{j}}}{n},
\end{align}
where $\hat{\Delta}_{f}(k) \! := \! \lbrace \mathstrut j \! \in \! \lbrace 
1,2,\dotsc,\mathfrak{s} \! - \! 2 \rbrace; \, \alpha_{p_{j}} \! > \! 
\alpha_{k} \rbrace$.\footnote{If $J_{f} \cap \mathbb{R}^{<}_{\alpha_{p_{q}}} 
\! = \! \varnothing$, then $\int_{J_{f} \cap \mathbb{R}^{<}_{\alpha_{p_{q}}}} 
\md \mu_{\widetilde{V}}^{f}(\xi) \! := \! 0$, $q \! \in \! \lbrace 1,\dotsc,
\mathfrak{s} \! - \! 2,\mathfrak{s} \rbrace$. If $\hat{\Delta}_{f}(k) \! = 
\! \varnothing$, then $\sum_{j \in \hat{\Delta}_{f}(k)} \varkappa_{nk 
\tilde{k}_{j}}/n \! := \! 0$.} Then, for $n \! \in \! \mathbb{N}$ and $k \! 
\in \! \lbrace 1,2,\dotsc,K \rbrace$, $\mathscr{X} \colon \mathbb{N} 
\times \lbrace 1,2,\dotsc,K \rbrace \times \overline{\mathbb{C}} 
\setminus \overline{\mathbb{R}} \! \to \! \operatorname{SL}_{2}
(\mathbb{C})$ solves the following {\rm RHP}$:$ {\rm \pmb{(i)}} 
$\mathscr{X}(z)$ is analytic for $z \! \in \! \overline{\mathbb{C}} 
\setminus \overline{\mathbb{R}}$$;$ {\rm \pmb{(ii)}} the boundary 
values $\mathscr{X}_{\pm}(z) \! := \! \lim_{\varepsilon \downarrow 0} 
\mathscr{X}(z \! \pm \! \mi \varepsilon)$ satisfy the jump condition
\begin{equation*}
\mathscr{X}_{+}(z) \! = \! \mathscr{X}_{-}(z) \mathscr{V}(z) \quad 
\mathrm{a.e.} \quad z \! \in \! \overline{\mathbb{R}},
\end{equation*}
where, for $n \! \in \! \mathbb{N}$ and $k \! \in \! \lbrace 1,2,\dotsc,K 
\rbrace$ such that $\alpha_{p_{\mathfrak{s}}} \! := \! \alpha_{k} \! = \! \infty$,
\begin{equation} \label{eql3.4g} 
\mathscr{V} \colon \mathbb{N} \times \lbrace 1,2,\dotsc,K \rbrace 
\times \overline{\mathbb{R}} \! \ni \! (n,k,z) \! \mapsto \! 
\begin{pmatrix}
\me^{-n(g_{+}^{\infty}(z)-g_{-}^{\infty}(z))} & \me^{n(g_{+}^{\infty}(z)+
g_{-}^{\infty}(z)-2 \tilde{\mathscr{P}}_{0}-\widetilde{V}(z)-\hat{\ell})} \\
0 & \me^{n(g_{+}^{\infty}(z)-g_{-}^{\infty}(z))}
\end{pmatrix} \! =: \! \mathscr{V}(n,k,z) \! = \! \mathscr{V}(z),
\end{equation}
and, for $n \! \in \! \mathbb{N}$ and $k \! \in \! \lbrace 1,2,\dotsc,K \rbrace$ 
such that $\alpha_{p_{\mathfrak{s}}} \! := \! \alpha_{k} \! \neq \! \infty$,
\begin{equation} \label{eql3.4h} 
\mathscr{V} \colon \mathbb{N} \times \lbrace 1,2,\dotsc,K \rbrace 
\times \overline{\mathbb{R}} \! \ni \! (n,k,z) \! \mapsto \! 
\begin{pmatrix}
\me^{-n(g_{+}^{f}(z)-g_{-}^{f}(z)+\hat{\mathscr{P}}_{0}^{-}-
\hat{\mathscr{P}}_{0}^{+})} & \me^{n(g_{+}^{f}(z)+g_{-}^{f}(z)-
\hat{\mathscr{P}}_{0}^{-}-\hat{\mathscr{P}}_{0}^{+}-\widetilde{V}(z)
-\tilde{\ell})} \\
0 & \me^{n(g_{+}^{f}(z)-g_{-}^{f}(z)+\hat{\mathscr{P}}_{0}^{-}-
\hat{\mathscr{P}}_{0}^{+})}
\end{pmatrix} \! =: \! \mathscr{V}(n,k,z) \! = \! \mathscr{V}(z),
\end{equation}
where $g_{\pm}^{r}(z) \! := \! \lim_{\varepsilon \downarrow 0}g^{r}
(z \! \pm \! \mi \varepsilon)$, $r \! \in \! \lbrace \infty,f \rbrace$$;$ 
{\rm \pmb{(iii)}} for $n \! \in \! \mathbb{N}$ and $k \! \in \! \lbrace 
1,2,\dotsc,K \rbrace$ such that $\alpha_{p_{\mathfrak{s}}} \! := \! 
\alpha_{k} \! = \! \infty$,
\begin{equation*}
\mathscr{X}(z) \underset{\overline{\mathbb{C}}_{\pm} \ni z \to 
\alpha_{k}}{=} \operatorname{I} \! + \mathcal{O}(z^{-1}), \quad 
\quad \mathscr{X}(z) \underset{\mathbb{C}_{\pm} \ni z \to 
\alpha_{p_{q}}}{=} \mathcal{O}(\me^{n(\tilde{\mathscr{P}}_{0}
-\tilde{\mathscr{P}}^{\pm}_{1}) \sigma_{3}}), \quad q \! = \! 
1,2,\dotsc,\mathfrak{s} \! - \! 1,
\end{equation*}
where $\tilde{\mathscr{P}}^{\pm}_{1}$ is defined by 
Equation~\eqref{eql3.4gee10}$;$ and {\rm \pmb{(iv)}} for $n \! \in \! 
\mathbb{N}$ and $k \! \in \! \lbrace 1,2,\dotsc,K \rbrace$ such that 
$\alpha_{p_{\mathfrak{s}}} \! := \! \alpha_{k} \! \neq \! \infty$,
\begin{gather*}
\mathscr{X}(z) \underset{\mathbb{C}_{\pm} \ni z \to \alpha_{k}}{=} 
\operatorname{I} \! + \mathcal{O}(z \! - \! \alpha_{k}), \qquad \quad 
\mathscr{X}(z) \underset{\overline{\mathbb{C}}_{\pm} \ni z \to 
\alpha_{p_{\mathfrak{s}-1}} = \infty}{=} \mathcal{O}(\me^{n
(\hat{\mathscr{P}}^{\pm}_{0}-\hat{\mathscr{P}}_{1}) \sigma_{3}}), \\
\qquad \mathscr{X}(z) \underset{\mathbb{C}_{\pm} \ni z \to 
\alpha_{p_{q}}}{=} \mathcal{O}(\me^{n(\hat{\mathscr{P}}^{\pm}_{0}
-\hat{\mathscr{P}}^{\pm}_{2}) \sigma_{3}}), \quad q \! = \! 1,2,\dotsc,
\mathfrak{s} \! - \! 2,
\end{gather*}
where $\hat{\mathscr{P}}_{1}$ and $\hat{\mathscr{P}}^{\pm}_{2}$ are 
defined by Equations~\eqref{eql3.4gee7} and~\eqref{eql3.4gee9}, 
respectively.
\end{ccccc}

\emph{Proof}. The proof of this Lemma~\ref{lem3.4} consists of two cases: 
(i) $n \! \in \! \mathbb{N}$ and $k \! \in \! \lbrace 1,2,\dotsc,K \rbrace$ 
such that $\alpha_{p_{\mathfrak{s}}} \! := \! \alpha_{k} \! = \! \infty$; and 
(ii) $n \! \in \! \mathbb{N}$ and $k \! \in \! \lbrace 1,2,\dotsc,K \rbrace$ 
such that $\alpha_{p_{\mathfrak{s}}} \! := \! \alpha_{k} \! \neq \! \infty$. 
The proof for the case $\alpha_{p_{\mathfrak{s}}} \! := \! \alpha_{k} \! \neq 
\! \infty$, $k \! \in \! \lbrace 1,2,\dotsc,K \rbrace$, will be considered in 
detail (see $\pmb{(1)}$ below), whilst the case $\alpha_{p_{\mathfrak{s}}} 
\! := \! \alpha_{k} \! = \! \infty$, $k \! \in \! \lbrace 1,2,\dotsc,K \rbrace$, 
can be proved analogously (see $\pmb{(2)}$ below).

$\pmb{(1)}$ For $n \! \in \! \mathbb{N}$ and $k \! \in \! \lbrace 1,2,
\dotsc,K \rbrace$ such that $\alpha_{p_{\mathfrak{s}}} \! := \! \alpha_{k} 
\! \neq \! \infty$, it follows {}from the definition of $g^{f}(z)$ (cf. 
Equation~\eqref{eql3.4gee3}) that, for arbitrary $z_{1},z_{2} \! \in \! 
\mathbb{C}_{\pm}$, $g^{f}(z_{2}) \! - \! g^{f}(z_{1}) \! = \! \mi 
\pi \int_{z_{1}}^{z_{2}} \mathfrak{F}_{f}(\xi) \, \md \xi$, where, with 
$\mathscr{D}_{f} \! := \! \mathbb{C} \setminus (J_{f} \cup \lbrace 
\alpha_{p_{1}},\dotsc,\alpha_{p_{\mathfrak{s}-2}},\alpha_{k} \rbrace)$, 
$\mathfrak{F}_{f} \colon \mathbb{N} \times \lbrace 1,2,\dotsc,K \rbrace 
\times \mathscr{D}_{f} \! \ni \! (n,k,z) \! \mapsto \! \mathfrak{F}_{f}
(n,k,z) \! = \! \mathfrak{F}_{f}(z)$, where
\begin{equation} \label{eql3.4a}
\mathfrak{F}_{f}(z) \! := \! -\dfrac{1}{\mi \pi} \left(\dfrac{(\varkappa_{nk} 
\! - \! 1)}{n(z \! - \! \alpha_{k})} \! + \! \sum_{q=1}^{\mathfrak{s}-2} 
\dfrac{\varkappa_{nk \tilde{k}_{q}}}{n(z \! - \! \alpha_{p_{q}})} 
\! + \! \left(\dfrac{(n \! - \! 1)K \! + \! k}{n} \right) \int_{J_{f}}
(\xi \! - \! z)^{-1} \, \md \mu^{f}_{\widetilde{V}}(\xi) \right).
\end{equation}
Since $J_{f} \cap \lbrace \alpha_{1},\alpha_{2},\dotsc,\alpha_{K} \rbrace \! 
= \! \varnothing$ and $\mu^{f}_{\widetilde{V}} \! \in \! \mathscr{M}_{1}
(\mathbb{R})$, it follows that (e.g., $0 \! < \! \vert z \! - \! \alpha_{k} \vert 
\! \ll \! \min \lbrace \min_{q=1,2,\dotsc,\mathfrak{s}-2} \lbrace \vert 
\alpha_{p_{q}} \! - \! \alpha_{k} \vert \rbrace,\min_{i \neq j \in \lbrace 
1,2,\dotsc,\mathfrak{s}-2 \rbrace} \lbrace \vert \alpha_{p_{i}} \! - \! 
\alpha_{p_{j}} \vert \rbrace,\inf_{t \in J_{f}} \lbrace \vert t \! - \! \alpha_{k} 
\rvert \rbrace \rbrace)$, 
\begin{equation*}
\mathfrak{F}_{f}(z) \! + \! \dfrac{1}{\mi \pi} \dfrac{(\varkappa_{nk} \! 
- \! 1)}{n(z \! - \! \alpha_{k})} \underset{\mathscr{D}_{f} \ni z \to 
\alpha_{k}}{=} \mathcal{O}(1),
\end{equation*}
and (e.g., $0 \! < \! \vert z \! - \! \alpha_{p_{q}} \vert \! \ll \! \min 
\lbrace \min_{q^{\prime}=1,2,\dotsc,\mathfrak{s}-2} \lbrace \vert 
\alpha_{p_{q^{\prime}}} \! - \! \alpha_{k} \vert \rbrace,\min_{i \neq j \in 
\lbrace 1,2,\dotsc,\mathfrak{s}-2 \rbrace} \lbrace \vert \alpha_{p_{i}} \! 
- \! \alpha_{p_{j}} \vert \rbrace,\inf_{\underset{q^{\prime}=1,2,\dotsc,
\mathfrak{s}-2}{t \in J_{f}}} \lbrace \vert t \! - \! \alpha_{p_{q^{\prime}}} 
\rvert \rbrace \rbrace$, $q \! = \! 1,2,\dotsc,\mathfrak{s} \! - \! 2)$,
\begin{equation*}
\mathfrak{F}_{f}(z) \! + \! \dfrac{1}{\mi \pi} \dfrac{\varkappa_{nk \tilde{
k}_{q}}}{n(z \! - \! \alpha_{p_{q}})} \underset{\mathscr{D}_{f} \ni z \to 
\alpha_{p_{q}}}{=} \mathcal{O}(1), \quad q \! = \! 1,2,\dotsc,\mathfrak{s} 
\! - \! 2;
\end{equation*}
hence, $\lvert g^{f}(z_{2}) \! - \! g^{f}(z_{1}) \rvert \! \leqslant \! 
\pi \lvert z_{2} \! - \! z_{1} \rvert \sup_{z \in \mathbb{C}_{\pm}} 
\lvert \mathfrak{F}_{f}(z) \rvert$. Thus, for $n \! \in \! \mathbb{N}$ 
and $k \! \in \! \lbrace 1,2,\dotsc,K \rbrace$ such that 
$\alpha_{p_{\mathfrak{s}}} \! := \! \alpha_{k} \! \neq \! \infty$, 
{}from the definition of $g^{f} \colon \mathbb{C} \setminus 
\mathfrak{D}_{f}^{\sharp} \! \to \! \mathbb{C}$, with $\mathfrak{D}_{f}^{
\sharp} \! := \! (-\infty,\max \lbrace \max_{q=1,\dotsc,\mathfrak{s}-2,
\mathfrak{s}} \lbrace \alpha_{p_{q}} \rbrace,\max \lbrace J_{f} \rbrace 
\rbrace)$, given in the lemma: (1) for $\xi \! \in \! J_{f}$ and $z \! \in \! 
\mathbb{C} \setminus \mathfrak{D}^{\sharp}_{f}$, with $0 \! < \! \vert 
z \! - \! \alpha_{k} \vert \! \ll \! \min \lbrace \min_{q=1,2,\dotsc,
\mathfrak{s}-2} \lbrace \vert \alpha_{p_{q}} \! - \! \alpha_{k} \vert \rbrace,
\min_{i \neq j \in \lbrace 1,2,\dotsc,\mathfrak{s}-2 \rbrace} \lbrace \vert 
\alpha_{p_{i}} \! - \! \alpha_{p_{j}} \vert \rbrace,\inf_{\xi \in J_{f}} \lbrace \vert 
\xi \! - \! \alpha_{k} \rvert \rbrace \rbrace$, and $\mu^{f}_{\widetilde{V}} 
\! \in \! \mathscr{M}_{1}(\mathbb{R})$, it follows {}from the expansions 
$\tfrac{1}{(z-\alpha_{k})-(\xi -\alpha_{k})} \! = \! -\sum_{j=0}^{l} 
\tfrac{(z-\alpha_{k})^{j}}{(\xi -\alpha_{k})^{j+1}} \! + \! 
\tfrac{(z-\alpha_{k})^{l+1}}{(\xi -\alpha_{k})^{l+1}(z-\xi)}$, $l \! \in 
\! \mathbb{N}_{0}$, and $\ln (1 \! - \! \lozenge) \! =_{\vert \lozenge 
\vert \to 0} \! -\sum_{j=1}^{\infty} \tfrac{\lozenge^{j}}{j}$, and a careful 
analysis of the branch cuts, that
\begin{align} \label{eql3.4b} 
g^{f}(z) \underset{\mathbb{C}_{\pm} \ni z \to \alpha_{k}}{=}& \, 
-\left(\dfrac{\varkappa_{nk} \! - \! 1}{n} \right) \ln (z \! - \! \alpha_{k}) 
\! + \! \hat{\mathscr{P}}_{0}(n,k) \! + \! \sum_{m=1}^{\infty} \dfrac{1}{m} 
\left(\sum_{q=1}^{\mathfrak{s}-2} \dfrac{\varkappa_{nk \tilde{k}_{q}}}{n
(\alpha_{p_{q}} \! - \! \alpha_{k})^{m}} \right. \nonumber \\
-&\left. \, \left(\dfrac{(n \! - \! 1)K \! + \! k}{n} \right) \int_{J_{f}}
(\xi \! - \! \alpha_{k})^{-m} \, \md \mu_{\widetilde{V}}^{f}(\xi) \right)
(z \! - \! \alpha_{k})^{m},
\end{align}
where $\hat{\mathscr{P}}_{0}(n,k)$ is defined by Equations~\eqref{eql3.4gee4} 
and~\eqref{eql3.4gee5}, and $\hat{\Delta}_{f}(k) \! := \! \lbrace 
\mathstrut j \! \in \! \lbrace 1,2,\dotsc,\mathfrak{s} \! - \! 2 \rbrace; \, 
\alpha_{p_{j}} \! > \! \alpha_{k} \rbrace$; (2) for $\xi \! \in \! J_{f}$ and 
$z \! \in \! \overline{\mathbb{C}} \setminus \mathfrak{D}^{\sharp}_{f}$, with 
$\vert z \vert \! \gg \! \max \lbrace \max_{i \neq j \in \lbrace 1,2,\dotsc,
\mathfrak{s}-2 \rbrace} \lbrace \vert \alpha_{p_{i}} \! - \! \alpha_{p_{j}} 
\vert \rbrace,\max_{q=1,2,\dotsc,\mathfrak{s}-2} \lbrace \vert 
\alpha_{p_{q}} \! - \! \alpha_{k} \vert \rbrace,\max_{q=1,\dotsc,
\mathfrak{s}-2,\mathfrak{s}} \lbrace \vert \alpha_{p_{q}} \vert \rbrace,
\max \lbrace J_{f} \rbrace \rbrace$, and $\mu^{f}_{\widetilde{V}} \! 
\in \! \mathscr{M}_{1}(\mathbb{R})$, it follows {}from the expansions 
$\tfrac{1}{\xi -z} \! = \! -\sum_{j=0}^{l} \tfrac{\xi^{j}}{z^{j+1}} \! + \! 
\tfrac{\xi^{l+1}}{z^{l+1}(\xi -z)}$, $l \! \in \! \mathbb{N}_{0}$, and 
$\ln (z \! - \! \lozenge) \! =_{\vert z \vert \to \infty} \! \ln (z) \! - \! 
\sum_{j=1}^{\infty} \tfrac{1}{j}(\tfrac{\lozenge}{z})^{j}$, and a careful 
analysis of the branch cuts, that
\begin{align} \label{eql3.4c} 
g^{f}(z) \underset{\overline{\mathbb{C}}_{\pm} \ni z \to \alpha_{p_{
\mathfrak{s}-1}} = \infty}{=}& \, \left(\dfrac{\varkappa_{nk \tilde{k}_{
\mathfrak{s}-1}}^{\infty} \! + \! 1}{n} \right) \ln (z) \! + \! 
\hat{\mathscr{P}}_{1}(n,k) \! + \! \sum_{m=1}^{\infty} \dfrac{1}{m} \left(
\left(\dfrac{\varkappa_{nk} \! - \! 1}{n} \right)(\alpha_{k})^{m} \! 
+ \! \sum_{q=1}^{\mathfrak{s}-2} \dfrac{\varkappa_{nk \tilde{k}_{q}}}{n}
(\alpha_{p_{q}})^{m} \right. \nonumber \\
-&\left. \, \left(\dfrac{(n \! - \! 1)K \! + \! k}{n} \right) \int_{J_{f}} 
\xi^{m} \, \md \mu_{\widetilde{V}}^{f}(\xi) \right)z^{-m},
\end{align}
where {}\footnote{If $J_{f} \cap \mathbb{R}^{<}_{\alpha_{p_{q}}} \! = \! 
\varnothing$, then $\int_{J_{f} \cap \mathbb{R}^{<}_{\alpha_{p_{q}}}} 
\md \mu_{\widetilde{V}}^{f}(\xi) \! := \! 0$, $q \! \in \! \lbrace 1,
\dotsc,\mathfrak{s} \! - \! 2,\mathfrak{s} \rbrace$.}
\begin{align} \label{eql3.4gee7} 
\hat{\mathscr{P}}_{1}(n,k) \! = \! \hat{\mathscr{P}}_{1} \! :=& \, -\left(
\dfrac{\varkappa_{nk} \! - \! 1}{n} \right) \int_{J_{f}} \ln (\lvert \xi 
\! - \! \alpha_{k} \rvert) \, \md \mu_{\widetilde{V}}^{f}(\xi) \! - \! 
\sum_{q=1}^{\mathfrak{s}-2} \dfrac{\varkappa_{nk \tilde{k}_{q}}}{n} 
\int_{J_{f}} \ln (\lvert \xi \! - \! \alpha_{p_{q}} \rvert) \, \md 
\mu_{\widetilde{V}}^{f}(\xi) \nonumber \\
-& \, \mi \pi \left(\dfrac{\varkappa_{nk} \! - \! 1}{n} \right) \int_{J_{f} 
\cap \mathbb{R}_{\alpha_{k}}^{<}} \md \mu_{\widetilde{V}}^{f}(\xi) \! - \! 
\mi \pi \sum_{q=1}^{\mathfrak{s}-2} \dfrac{\varkappa_{nk \tilde{k}_{q}}}{n} 
\int_{J_{f} \cap \mathbb{R}_{\alpha_{p_{q}}}^{<}} \md \mu_{\widetilde{V}}^{f}
(\xi);
\end{align}
and (3) for $\xi \! \in \! J_{f}$ and $z \! \in \! \mathbb{C} \setminus 
\mathfrak{D}^{\sharp}_{f}$, with $0 \! < \! \vert z \! - \! \alpha_{p_{q}} 
\vert \! \ll \! \min \lbrace \min_{q^{\prime}=1,2,\dotsc,\mathfrak{s}-2} 
\lbrace \vert \alpha_{p_{q^{\prime}}} \! - \! \alpha_{k} \vert \rbrace,
\min_{i \neq j \in \lbrace 1,2,\dotsc,\mathfrak{s}-2 \rbrace} 
\lbrace \vert \alpha_{p_{i}} \! - \! \alpha_{p_{j}} \vert \rbrace,
\inf_{\underset{q^{\prime}=1,2,\dotsc,\mathfrak{s}-2}{\xi \in J_{f}}} \lbrace 
\vert \xi \! - \! \alpha_{p_{q^{\prime}}} \rvert \rbrace \rbrace$, $q \! 
= \! 1,2,\dotsc,\mathfrak{s} \! - \! 2$, and $\mu^{f}_{\widetilde{V}} \! 
\in \! \mathscr{M}_{1}(\mathbb{R})$, it follows {}from the expansions 
$\tfrac{1}{(z-\alpha_{p_{q}})-(\xi -\alpha_{p_{q}})} \! = \! -\sum_{j=0}^{l} 
\tfrac{(z-\alpha_{p_{q}})^{j}}{(\xi -\alpha_{p_{q}})^{j+1}} \! + \! 
\tfrac{(z-\alpha_{p_{q}})^{l+1}}{(\xi -\alpha_{p_{q}})^{l+1}(z-\xi)}$, $l \! 
\in \! \mathbb{N}_{0}$, and $\ln (1 \! - \! \lozenge) \! =_{\vert \lozenge 
\vert \to 0} \! -\sum_{j=1}^{\infty} \tfrac{\lozenge^{j}}{j}$, and a careful 
analysis of the branch cuts, that
\begin{align} \label{eql3.4d} 
g^{f}(z) \underset{\mathbb{C}_{\pm} \ni z \to \alpha_{p_{q}}}{=}& \, 
-\dfrac{\varkappa_{nk \tilde{k}_{q}}}{n} \ln (z \! - \! \alpha_{p_{q}}) \! 
+ \! \hat{\mathscr{P}}_{2}(n,k) \! + \! \sum_{m=1}^{\infty} \dfrac{1}{m} 
\left(\dfrac{(\varkappa_{nk} \! - \! 1)}{n(\alpha_{k} \! - \! 
\alpha_{p_{q}})^{m}} \! + \! \sum_{\underset{j \neq q}{j=1}}^{\mathfrak{s}-2} 
\dfrac{\varkappa_{nk \tilde{k}_{j}}}{n(\alpha_{p_{j}} \! - \! \alpha_{p_{q}})^{m}} 
\right. \nonumber \\
-&\left. \, \left(\dfrac{(n \! - \! 1)K \! + \! k}{n} \right) \int_{J_{f}}
(\xi \! - \! \alpha_{p_{q}})^{-m} \, \md \mu^{f}_{\widetilde{V}}(\xi) \right)
(z \! - \! \alpha_{p_{q}})^{m}, \quad q \! = \! 1,2,\dotsc,\mathfrak{s} 
\! - \! 2,
\end{align}
where
\begin{equation} \label{eql3.4gee8} 
\hat{\mathscr{P}}_{2}(n,k) \! = \! 
\begin{cases}
\hat{\mathscr{P}}_{2}^{+}, &\text{$z \! \in \! \mathbb{C}_{+}$,} \\
\hat{\mathscr{P}}_{2}^{-}, &\text{$z \! \in \! \mathbb{C}_{-}$,}
\end{cases}
\end{equation}
with
\begin{align} \label{eql3.4gee9} 
\hat{\mathscr{P}}_{2}^{\pm} :=& \, \dfrac{1}{n} \int_{J_{f}} \ln \left(\dfrac{
\lvert \xi \! - \! \alpha_{p_{q}} \rvert^{\varkappa_{nk}+\varkappa_{nk 
\tilde{k}_{\mathfrak{s}-1}}^{\infty}}}{(\lvert \xi \! - \! \alpha_{k} \rvert 
\lvert \alpha_{k} \! - \! \alpha_{p_{q}} \rvert)^{\varkappa_{nk}-1}} \prod_{
\substack{j=1\\j \neq q}}^{\mathfrak{s}-2} \left(\dfrac{\lvert \xi \! - \! 
\alpha_{p_{q}} \rvert}{\lvert \xi \! - \! \alpha_{p_{j}} \rvert \lvert \alpha_{p_{j}} 
\! - \! \alpha_{p_{q}} \rvert} \right)^{\varkappa_{nk \tilde{k}_{j}}} \right) \md 
\mu_{\widetilde{V}}^{f}(\xi) \nonumber \\
-& \, \mi \pi \left(\dfrac{\varkappa_{nk} \! - \! 1}{n} \right) \int_{
J_{f} \cap \mathbb{R}_{\alpha_{k}}^{<}} \md \mu_{\widetilde{V}}^{f}(\xi) \! - 
\! \mi \pi \dfrac{\varkappa_{nk \tilde{k}_{q}}}{n} \int_{J_{f} \cap 
\mathbb{R}_{\alpha_{p_{q}}}^{<}} \md \mu_{\widetilde{V}}^{f}(\xi) \! - \! \mi 
\pi \sum_{\substack{j=1\\j \neq q}}^{\mathfrak{s}-2} \dfrac{\varkappa_{nk 
\tilde{k}_{j}}}{n} \int_{J_{f} \cap \mathbb{R}_{\alpha_{p_{j}}}^{<}} \md 
\mu_{\widetilde{V}}^{f}(\xi) \nonumber \\
\mp& \, \mi \pi \sum_{j \in \tilde{\Delta}_{f}(q)} \dfrac{\varkappa_{nk 
\tilde{k}_{j}}}{n} \! \mp \! \mi \pi \left(\dfrac{\varkappa_{nk} \! - \! 
1}{n} \right) \epsilon (k,q) \! \pm \! \mi \pi \left(\dfrac{(n \! - \! 
1)K \! + \! k}{n} \right) \int_{J_{f} \cap \mathbb{R}_{\alpha_{p_{q}}}^{>}} \md 
\mu_{\widetilde{V}}^{f}(\xi),
\end{align}
where $\tilde{\Delta}_{f}(q) \! := \! \lbrace \mathstrut j \! \in \! \lbrace 
1,2,\dotsc,\mathfrak{s} \! - \! 2 \rbrace \setminus \lbrace q \rbrace; \, 
\alpha_{p_{j}} \! > \! \alpha_{p_{q}} \rbrace$,\footnote{If $J_{f} \cap 
\mathbb{R}^{<}_{\alpha_{p_{j}}} \! = \! \varnothing$, then $\int_{J_{f} 
\cap \mathbb{R}^{<}_{\alpha_{p_{j}}}} \md \mu_{\widetilde{V}}^{f}(\xi) 
\! := \! 0$, $j \! \in \! \lbrace 1,2,\dotsc,\mathfrak{s} \! - \! 2 \rbrace 
\setminus \lbrace q \rbrace$, $q \! = \! 1,2,\dotsc,\mathfrak{s} \! - 
\! 2$. If $\tilde{\Delta}_{f}(q) \! = \! \varnothing$, then $\sum_{j \in 
\tilde{\Delta}_{f}(q)} \varkappa_{nk \tilde{k}_{j}}/n \! := \! 0$, $q \! 
\in \! \lbrace 1,2,\dotsc,\mathfrak{s} \! - \! 2 \rbrace$. If $J_{f} 
\cap \mathbb{R}^{>}_{\alpha_{p_{q}}} \! = \! \varnothing$, 
then $\int_{J_{f} \cap \mathbb{R}^{>}_{\alpha_{p_{q}}}} \md 
\mu_{\widetilde{V}}^{f}(\xi) \! := \! 0$, $q \! \in \! \lbrace 1,2,
\dotsc,\mathfrak{s} \! - \! 2 \rbrace$.} and $\epsilon (k,q) \! = \! 
\left\{
\begin{smallmatrix}
1, \, \, \, \alpha_{p_{q}}< \alpha_{k}, \\
0, \, \, \, \alpha_{p_{q}}> \alpha_{k}.
\end{smallmatrix}
\right.$ For $n \! \in \! \mathbb{N}$ and $k \! \in \! \lbrace 1,2,\dotsc,
K \rbrace$ such that $\alpha_{p_{\mathfrak{s}}} \! := \! \alpha_{k} \! 
\neq \! \infty$, recalling that $\mathcal{X} \colon \mathbb{N} \times 
\lbrace 1,2,\dotsc,K \rbrace \times \overline{\mathbb{C}} \setminus 
\overline{\mathbb{R}} \! \to \! \operatorname{SL}_{2}(\mathbb{C})$ 
solves, uniquely, the RHP stated in 
Lemma~{\rm $\bm{\mathrm{RHP}_{\mathrm{MPC}}}$}, it follows, via the 
Definition~\eqref{eql3.4f} and the asymptotic expansions~\eqref{eql3.4b}, 
\eqref{eql3.4c}, and~\eqref{eql3.4d}, that $\mathscr{X} \colon \mathbb{N} 
\times \lbrace 1,2,\dotsc,K \rbrace \times \overline{\mathbb{C}} 
\setminus \overline{\mathbb{R}} \! \to \! \operatorname{SL}_{2}
(\mathbb{C})$ is analytic for $z \! \in \! \overline{\mathbb{C}} \setminus 
\overline{\mathbb{R}}$ (cf. item~\pmb{(i)} of the lemma), its boundary 
values, $\mathscr{X}_{\pm}(z) \! := \! \lim_{\varepsilon \downarrow 0} 
\mathscr{X}(z \! \pm \! \mi \varepsilon)$, satisfy the jump condition 
$\mathscr{X}_{+}(z) \! = \! \mathscr{X}_{-}(z) \mathscr{V}(z)$ 
$\mathrm{a.e.}$ $z \! \in \! \overline{\mathbb{R}}$, where the 
associated jump matrix, $\mathscr{V}(z)$, is given by 
definition~\eqref{eql3.4h} (cf. item~\pmb{(ii)} of the lemma), and 
satisfies the asymptotic conditions stated in item~\pmb{(iv)} of the 
lemma.

$\pmb{(2)}$ The proof corresponding to the case $\alpha_{p_{\mathfrak{s}}} 
\! := \! \alpha_{k} \! = \! \infty$, $k \! \in \! \lbrace 1,2,\dotsc,K \rbrace$, 
is virtually identical to the proof presented in $\pmb{(1)}$ above; one mimics, 
\emph{verbatim}, the scheme of the calculations presented in case $\pmb{(1)}$ 
in order to arrive at the corresponding items of the RHP for $\mathscr{X} 
\colon \mathbb{N} \times \lbrace 1,2,\dotsc,K \rbrace \times 
\overline{\mathbb{C}} \setminus \overline{\mathbb{R}} \! \to \! 
\operatorname{SL}_{2}(\mathbb{C})$ stated in the lemma. More precisely, 
via the Definition~\eqref{eql3.4e}, the analogue of Equation~\eqref{eql3.4a}, 
that is, $\mathfrak{F}_{\infty} \colon \mathbb{N} \times \lbrace 1,2,
\dotsc,K \rbrace \times \mathscr{D}_{\infty} \! \ni \! (n,k,z) \! \mapsto 
\! \mathfrak{F}_{\infty}(n,k,z) \! = \! \mathfrak{F}_{\infty}(z)$, where, 
with $\mathscr{D}_{\infty} \! := \! \mathbb{C} \setminus (J_{\infty} \cup 
\lbrace \alpha_{p_{1}},\alpha_{p_{2}},\dotsc,\alpha_{p_{\mathfrak{s}-1}} 
\rbrace)$,
\begin{equation} \label{eql3.4i} 
\mathfrak{F}_{\infty}(z) \! := \! -\dfrac{1}{\mi \pi} \left(\sum_{q=1}^{
\mathfrak{s}-1} \dfrac{\varkappa_{nk \tilde{k}_{q}}}{n(z \! - \! 
\alpha_{p_{q}})} \! + \! \left(\dfrac{(n \! - \! 1)K \! + \! k}{n} \right) 
\int_{J_{\infty}}(\xi \! - \! z)^{-1} \, \md \mu_{\widetilde{V}}^{\infty}
(\xi) \right),
\end{equation}
and the asymptotic expansions
\begin{align} \label{eql3.4j} 
g^{\infty}(z) \underset{\overline{\mathbb{C}}_{\pm} \ni z \to 
\alpha_{k}}{=}& \, \dfrac{\varkappa_{nk}}{n} \ln (z) \! + \! \tilde{
\mathscr{P}}_{0}(n,k) \! + \! \sum_{m=1}^{\infty} \dfrac{1}{m} \left(
\sum_{q=1}^{\mathfrak{s}-1} \dfrac{\varkappa_{nk \tilde{k}_{q}}}{n}
(\alpha_{p_{q}})^{m} \! - \! \left(\dfrac{(n \! - \! 1)K \! + \! k}{n} 
\right) \int_{J_{\infty}} \xi^{m} \md \mu_{\widetilde{V}}^{\infty}(\xi) 
\right)z^{-m},
\end{align}
where $\tilde{\mathscr{P}}_{0}(n,k)$ is defined by 
Equation~\eqref{eql3.4gee2}, and
\begin{align} \label{eql3.4k} 
g^{\infty}(z) \underset{\mathbb{C}_{\pm} \ni z \to \alpha_{p_{q}}}{=}& \, 
-\dfrac{\varkappa_{nk \tilde{k}_{q}}}{n} \ln (z \! - \! \alpha_{p_{q}}) \! + 
\! \tilde{\mathscr{P}}_{1}(n,k) \! + \! \sum_{m=1}^{\infty} \dfrac{1}{m} 
\left(\sum_{\substack{j=1\\j \neq q}}^{\mathfrak{s}-1} \dfrac{\varkappa_{nk 
\tilde{k}_{j}}}{n(\alpha_{p_{j}} \! - \! \alpha_{p_{q}})^{m}} \! - \! 
\left(\dfrac{(n \! - \! 1)K \! + \! k}{n} \right) \right. \nonumber \\
\times&\left. \, \int_{J_{\infty}}(\xi \! - \! \alpha_{p_{q}})^{-m} \, \md 
\mu_{\widetilde{V}}^{\infty}(\xi) \right)(z \! - \! \alpha_{p_{q}})^{m}, 
\quad q \! = \! 1,2,\dotsc,\mathfrak{s} \! - \! 1,
\end{align}
where
\begin{equation} \label{eql3.4gee6} 
\tilde{\mathscr{P}}_{1}(n,k) \! = \! 
\begin{cases}
\tilde{\mathscr{P}}_{1}^{+}, &\text{$z \! \in \! \mathbb{C}_{+}$,} \\
\tilde{\mathscr{P}}_{1}^{-}, &\text{$z \! \in \! \mathbb{C}_{-}$,}
\end{cases}
\end{equation}
with
\begin{align} \label{eql3.4gee10} 
\tilde{\mathscr{P}}_{1}^{\pm} :=& \, \dfrac{1}{n} \int_{J_{\infty}} \ln \left(\lvert 
\xi \! - \! \alpha_{p_{q}} \rvert^{\varkappa_{nk}} \prod_{\substack{j=1\\j 
\neq q}}^{\mathfrak{s}-1} \left(\dfrac{\lvert \xi \! - \! \alpha_{p_{q}} \rvert}{\lvert 
\xi \! - \! \alpha_{p_{j}} \rvert \lvert \alpha_{p_{j}} \! - \! \alpha_{p_{q}} \rvert} 
\right)^{\varkappa_{nk \tilde{k}_{j}}} \right) \md \mu_{\widetilde{V}}^{\infty}
(\xi) \! - \! \mi \pi \dfrac{\varkappa_{nk \tilde{k}_{q}}}{n} \int_{J_{\infty} \cap 
\mathbb{R}_{\alpha_{p_{q}}}^{<}} \md \mu_{\widetilde{V}}^{\infty}(\xi) 
\nonumber \\
-& \, \mi \pi \sum_{\substack{j=1\\j \neq q}}^{\mathfrak{s}-1} 
\dfrac{\varkappa_{nk \tilde{k}_{j}}}{n} \int_{J_{\infty} \cap \mathbb{R}_{
\alpha_{p_{j}}}^{<}} \md \mu_{\widetilde{V}}^{\infty}(\xi) \! \mp \! \mi \pi 
\sum_{j \in \tilde{\Delta}_{\infty}(q)} \dfrac{\varkappa_{nk \tilde{k}_{j}}}{n} 
\! \pm \! \mi \pi \left(\dfrac{(n \! - \! 1)K \! + \! k}{n} \right) 
\int_{J_{\infty} \cap \mathbb{R}_{\alpha_{p_{q}}}^{>}} \md 
\mu_{\widetilde{V}}^{\infty}(\xi),
\end{align}
where $\tilde{\Delta}_{\infty}(q) \! := \! \lbrace \mathstrut j \! \in \! 
\lbrace 1,2,\dotsc,\mathfrak{s} \! - \! 1 \rbrace \setminus \lbrace q 
\rbrace; \, \alpha_{p_{j}} \! > \! \alpha_{p_{q}} \rbrace$,\footnote{If 
$J_{\infty} \cap \mathbb{R}^{<}_{\alpha_{p_{q}}} \! = \! \varnothing$, 
then $\int_{J_{\infty} \cap \mathbb{R}^{<}_{\alpha_{p_{q}}}} \md 
\mu_{\widetilde{V}}^{\infty}(\xi) \! := \! 0$, $q \! \in \! \lbrace 
1,2,\dotsc,\mathfrak{s} \! - \! 1 \rbrace$. If $J_{\infty} \cap 
\mathbb{R}^{<}_{\alpha_{p_{j}}} \! = \! \varnothing$, then $\int_{J_{\infty} 
\cap \mathbb{R}^{<}_{\alpha_{p_{j}}}} \md \mu_{\widetilde{V}}^{\infty}
(\xi) \! := \! 0$, $j \! \in \! \lbrace 1,2,\dotsc,\mathfrak{s} \! - \! 1 
\rbrace \setminus \lbrace q \rbrace$, $q \! = \! 1,2,\dotsc,
\mathfrak{s} \! - \! 1$. If $\tilde{\Delta}_{\infty}(q) \! = \! 
\varnothing$, then $\sum_{j \in \tilde{\Delta}_{\infty}(q)} 
\varkappa_{nk \tilde{k}_{j}}/n \! := \! 0$, $q \! \in \! \lbrace 
1,2,\dotsc,\mathfrak{s} \! - \! 1 \rbrace$. If $J_{\infty} \cap 
\mathbb{R}^{>}_{\alpha_{p_{q}}} \! = \! \varnothing$, then $\int_{J_{f} 
\cap \mathbb{R}^{>}_{\alpha_{p_{q}}}} \md \mu_{\widetilde{V}}^{\infty}
(\xi) \! := \! 0$, $q \! \in \! \lbrace 1,2,\dotsc,\mathfrak{s} \! - \! 1 
\rbrace$.} one verifies that $\mathscr{X} \colon \mathbb{N} \times 
\lbrace 1,2,\dotsc,K \rbrace \times \overline{\mathbb{C}} \setminus 
\overline{\mathbb{R}} \! \to \! \operatorname{SL}_{2}(\mathbb{C})$ is 
analytic for $z \! \in \! \overline{\mathbb{C}} \setminus \overline{
\mathbb{R}}$ (cf. item~\pmb{(i)} of the lemma), its boundary values 
satisfy the jump condition $\mathscr{X}_{+}(z) \! = \! \mathscr{X}_{-}
(z) \mathscr{V}(z)$ $\mathrm{a.e.}$ $z \! \in \! \overline{\mathbb{R}}$, 
where the associated jump matrix, $\mathscr{V}(z)$, is given by 
definition~\eqref{eql3.4g} (cf. item~\pmb{(ii)} of the lemma), and 
satisfies the asymptotic conditions stated in item~\pmb{(iii)} of the 
lemma. This concludes the proof. \hfill $\qed$

The following---technical---Lemma~\ref{lem3.5}, which is modelled on 
Theorem~6.6.2 of \cite{a51} (see, in particular, Chapter~6 of \cite{a51}), 
is necessary in order to prove the seminal Lemma~\ref{lem3.6} below.
\begin{ccccc} \label{lem3.5} 
Let $\widetilde{V} \colon \overline{\mathbb{R}} \setminus \lbrace 
\alpha_{1},\alpha_{2},\dotsc,\alpha_{K} \rbrace \! \to \! \mathbb{R}$ 
satisfy conditions~\eqref{eq20}--\eqref{eq22}. For $n \! \in \! 
\mathbb{N}$ and $k \! \in \! \lbrace 1,2,\dotsc,K \rbrace$ such that 
$\alpha_{p_{\mathfrak{s}}} \! := \! \alpha_{k} \! = \! \infty$, let
\begin{align*}
\mathfrak{d}^{\widetilde{V},\infty}_{(n-1)K+k} :=& \, \dfrac{1}{
((n \! - \! 1)K \! + \! k)((n \! - \! 1)K \! + \! k \! - \! 1)} \, 
\inf_{\{x_{1},x_{2},\dotsc,x_{(n-1)K+k}\} \subset \mathbb{R}} 
\left(\sum_{\substack{i,j=1\\i \neq j}}^{(n-1)K+k} \ln \left(
\vphantom{M^{M^{M^{M^{M^{M^{M}}}}}}} 
\lvert x_{i} \! - \! x_{j} \rvert^{\frac{\varkappa_{nk}}{n}} \right. \right. \\
\times&\left. \left. \, \prod_{q=1}^{\mathfrak{s}-1} \left(\dfrac{\lvert 
x_{i} \! - \! x_{j} \rvert}{\lvert x_{i} \! - \! \alpha_{p_{q}} \rvert 
\lvert x_{j} \! - \! \alpha_{p_{q}} \rvert} \right)^{\frac{\varkappa_{n
k \tilde{k}_{q}}}{n}} \right)^{-1} \! + \! ((n \! - \! 1)K \! + \! k \! - \! 
1) \sum_{m=1}^{(n-1)K+k} \widetilde{V}(x_{m}) \right),
\end{align*}
and, for $n \! \in \! \mathbb{N}$ and $k \! \in \! \lbrace 1,2,\dotsc,
K \rbrace$ such that $\alpha_{p_{\mathfrak{s}}} \! := \! \alpha_{k} \! 
\neq \! \infty$, let
\begin{align*}
\mathfrak{d}^{\widetilde{V},f}_{(n-1)K+k} :=& \, \dfrac{1}{((n \! - 
\! 1)K \! + \! k)((n \! - \! 1)K \! + \! k \! - \! 1)} \, \inf_{\{x_{1},
x_{2},\dotsc,x_{(n-1)K+k}\} \subset \mathbb{R}} \left(
\sum_{\substack{i,j=1\\i \neq j}}^{(n-1)K+k} \ln \left(\lvert x_{i} \! - \! 
x_{j} \rvert^{\frac{\varkappa_{nk \tilde{k}_{\mathfrak{s}-1}}^{\infty}+1}{n}} 
\left(\dfrac{\lvert x_{i} \! - \! x_{j} \rvert}{\lvert x_{i} \! - \! \alpha_{k} 
\rvert \lvert x_{j} \! - \! \alpha_{k} \rvert} \right)^{\frac{\varkappa_{nk}
-1}{n}} \right. \right. \\
\times&\left. \left. \, \prod_{q=1}^{\mathfrak{s}-2} \left(\dfrac{\lvert 
x_{i} \! - \! x_{j} \rvert}{\lvert x_{i} \! - \! \alpha_{p_{q}} \rvert 
\lvert x_{j} \! - \! \alpha_{p_{q}} \rvert} \right)^{\frac{\varkappa_{nk 
\tilde{k}_{q}}}{n}} \right)^{-1} \! + \! ((n \! - \! 1)K \! + \! k \! - \! 1) 
\sum_{m=1}^{(n-1)K+k} \widetilde{V}(x_{m}) \right).
\end{align*}
Then, $\lim_{\underset{z_{o}=1+o(1)}{\mathscr{N},n \to \infty}} 
\mathfrak{d}^{\widetilde{V},r}_{(n-1)K+k}$ exists, $r \! \in \! \lbrace 
\infty,f \rbrace$, that is, $\lim_{\underset{z_{o}=1+o(1)}{\mathscr{N},
n \to \infty}} \mathfrak{d}^{\widetilde{V},r}_{(n-1)K+k} \! = \! 
E_{\widetilde{V}}^{r} \! = \! \inf \lbrace \mathrm{I}_{\widetilde{V}}^{r}
[\mu^{\text{\tiny $\mathrm{EQ}$}}]; \, \mu^{\text{\tiny $\mathrm{EQ}$}} 
\! \in \! \mathscr{M}_{1}(\mathbb{R}) \rbrace \! = \! \mathrm{I}_{
\widetilde{V}}^{r}[\mu_{\widetilde{V}}^{r}]$, and $\lim_{\underset{z_{o}
=1+o(1)}{\mathscr{N},n \to \infty}} \exp (-\mathfrak{d}^{\widetilde{V},
r}_{(n-1)K+k}) \! = \! \exp (-E_{\widetilde{V}}^{r}) \! > \! 0$ and 
$\mathcal{O}(1)$.

For $n \! \in \! \mathbb{N}$ and $k \! \in \! \lbrace 1,2,\dotsc,K 
\rbrace$ such that $\alpha_{p_{\mathfrak{s}}} \! := \! \alpha_{k} 
\! = \! \infty$, let $\hat{x}_{1}^{\ast},\hat{x}_{2}^{\ast},\dotsc,
\hat{x}_{(n-1)K+k}^{\ast}$ denote the associated, weighted Fekete 
points, that is,
\begin{align*}
\mathfrak{d}^{\widetilde{V},\infty}_{(n-1)K+k} =& \, \dfrac{1}{
((n \! - \! 1)K \! + \! k)((n \! - \! 1)K \! + \! k \! - \! 1)} \left(
\sum_{\substack{i,j=1\\i \neq j}}^{(n-1)K+k} \ln \left(\lvert 
\hat{x}_{i}^{\ast} \! - \! \hat{x}_{j}^{\ast} \rvert^{\frac{\varkappa_{nk}}{n}} 
\prod_{q=1}^{\mathfrak{s}-1} \left(\dfrac{\lvert \hat{x}_{i}^{\ast} \! - 
\! \hat{x}_{j}^{\ast} \rvert}{\lvert \hat{x}_{i}^{\ast} \! - \! \alpha_{p_{q}} 
\rvert \lvert \hat{x}_{j}^{\ast} \! - \! \alpha_{p_{q}} \rvert} 
\right)^{\frac{\varkappa_{nk \tilde{k}_{q}}}{n}} \right)^{-1} \right. \\
+&\left. \, ((n \! - \! 1)K \! + \! k \! - \! 1) \sum_{m=1}^{(n-1)K+k} 
\widetilde{V}(\hat{x}_{m}^{\ast}) \right),
\end{align*}
and, for $n \! \in \! \mathbb{N}$ and $k \! \in \! \lbrace 1,2,\dotsc,K 
\rbrace$ such that $\alpha_{p_{\mathfrak{s}}} \! := \! \alpha_{k} \! 
\neq \! \infty$, let $\tilde{x}_{1}^{\ast},\tilde{x}_{2}^{\ast},\dotsc,
\tilde{x}_{(n-1)K+k}^{\ast}$ denote the associated, weighted Fekete 
points, that is,
\begin{align*}
\mathfrak{d}^{\widetilde{V},f}_{(n-1)K+k} =& \, \dfrac{1}{((n \! - \! 1)
K \! + \! k)((n \! - \! 1)K \! + \! k \! - \! 1)} \left(\sum_{\substack{i,j=
1\\i \neq j}}^{(n-1)K+k} \ln \left(\lvert \tilde{x}_{i}^{\ast} \! - \! 
\tilde{x}_{j}^{\ast} \rvert^{\frac{\varkappa_{nk \tilde{k}_{\mathfrak{s}
-1}}^{\infty}+1}{n}} \left(\dfrac{\lvert \tilde{x}_{i}^{\ast} \! - \! 
\tilde{x}_{j}^{\ast} \rvert}{\lvert \tilde{x}_{i}^{\ast} \! - \! \alpha_{k} 
\rvert \lvert \tilde{x}_{j}^{\ast} \! - \! \alpha_{k} \rvert} \right)^{
\frac{\varkappa_{nk}-1}{n}} \right. \right. \\
\times&\left. \left. \, \prod_{q=1}^{\mathfrak{s}-2} \left(\dfrac{
\lvert \tilde{x}_{i}^{\ast} \! - \! \tilde{x}_{j}^{\ast} \rvert}{\lvert 
\tilde{x}_{i}^{\ast} \! - \! \alpha_{p_{q}} \rvert \lvert \tilde{x}_{j}^{\ast} 
\! - \! \alpha_{p_{q}} \rvert} \right)^{\frac{\varkappa_{nk 
\tilde{k}_{q}}}{n}} \right)^{-1} \! + \! ((n \! - \! 1)K \! + \! k \! - \! 1) 
\sum_{m=1}^{(n-1)K+k} \widetilde{V}(\tilde{x}_{m}^{\ast}) 
\right).
\end{align*}
For $n \! \in \! \mathbb{N}$ and $k \! \in \! \lbrace 1,2,\dotsc,K 
\rbrace$ such that $\alpha_{p_{\mathfrak{s}}} \! := \! \alpha_{k} \! 
= \! \infty$ (resp., $\alpha_{p_{\mathfrak{s}}} \! := \!\alpha_{k} \! 
\neq \! \infty)$, with $\lbrace \hat{x}_{1}^{\ast},\hat{x}_{2}^{\ast},\dotsc,
\hat{x}_{(n-1)K+k}^{\ast} \rbrace$ (resp., $\lbrace \tilde{x}_{1}^{\ast},
\tilde{x}_{2}^{\ast},\dotsc,\tilde{x}_{(n-1)K+k}^{\ast} \rbrace)$ an 
associated, weighted $((n \! - \! 1)K \! + \! k)$-Fekete set, denote by
\begin{equation*}
\lambda_{(n-1)K+k}^{r} \! := \! \dfrac{1}{((n \! - \! 1)K \! + \! k)} 
\sum_{j=1}^{(n-1)K+k} \delta_{x_{j}(r)}, \quad r \! \in \! \lbrace \infty,
f \rbrace,
\end{equation*}
where $x_{j}(r) \! = \! \hat{x}_{j}^{\ast}$ for $r \! = \! \infty$ and $x_{j}(r) 
\! = \! \tilde{x}_{j}^{\ast}$ for $r \! = \! f$, and $\delta_{x_{j}(r)}$ is the Dirac 
delta (atomic) mass concentrated at $x_{j}(r)$, the associated normalised counting 
measure $(\int_{\mathbb{R}} \md \lambda_{(n-1)K+k}^{r}(\xi) \! = \! 1$, 
$r \! \in \! \lbrace \infty,f \rbrace)$. Then, for $n \! \in \! \mathbb{N}$ and 
$k \! \in \! \lbrace 1,2,\dotsc,K \rbrace$ such that $\alpha_{p_{\mathfrak{s}}} 
\! := \! \alpha_{k} \! = \! \infty$ (resp., $\alpha_{p_{\mathfrak{s}}} \! 
:= \! \alpha_{k} \! \neq \! \infty)$, $\lambda_{(n-1)K+k}^{\infty} \! 
\overset{\ast}{\to} \! \mu_{\widetilde{V}}^{\infty}$ (resp., $\lambda_{(n-1)K+k}^{f} 
\! \overset{\ast}{\to} \! \mu_{\widetilde{V}}^{f})$ in the double-scaling limit 
$\mathscr{N},n \! \to \! \infty$ such that $z_{o} \! = \! 1 \! + \! o(1)$.
\end{ccccc}
\begin{eeeee} \label{remm3.5} 
\textsl{In fact, as shown in the proof of Lemma~\ref{lem3.5} below, 
for $n \! \in \! \mathbb{N}$ and $k \! \in \! \lbrace 1,2,\dotsc,K 
\rbrace$ such that $\alpha_{p_{\mathfrak{s}}} \! := \! \alpha_{k} \! 
= \! \infty$ (resp., $\alpha_{p_{\mathfrak{s}}} \! := \! \alpha_{k} \! 
\neq \! \infty)$, $\lbrace \hat{x}_{1}^{\ast},\hat{x}_{2}^{\ast},\dotsc,
\hat{x}^{\ast}_{(n-1)K+k} \rbrace$ (resp., $\lbrace \tilde{x}_{1}^{\ast},
\tilde{x}_{2}^{\ast},\dotsc,\tilde{x}^{\ast}_{(n-1)K+k} \rbrace)$ $\subset 
\overline{\mathbb{R}} \setminus \lbrace \alpha_{1},\alpha_{2},\dotsc,
\alpha_{K} \rbrace$.}
\end{eeeee}

\emph{Proof}. The proof of this Lemma~\ref{lem3.5} consists of two 
cases: (i) $n \! \in \! \mathbb{N}$ and $k \! \in \! \lbrace 1,2,\dotsc,K 
\rbrace$ such that $\alpha_{p_{\mathfrak{s}}} \! := \! \alpha_{k} \! = \! 
\infty$; and (ii) $n \! \in \! \mathbb{N}$ and $k \! \in \! \lbrace 1,2,
\dotsc,K \rbrace$ such that $\alpha_{p_{\mathfrak{s}}} \! := \! \alpha_{k} 
\! \neq \! \infty$. The proof for the case $\alpha_{p_{\mathfrak{s}}} \! 
:= \! \alpha_{k} \! \neq \! \infty$, $k \! \in \! \lbrace 1,2,\dotsc,K \rbrace$, 
will be considered in detail (see $\pmb{(1)}$ below), whilst the case 
$\alpha_{p_{\mathfrak{s}}} \! := \! \alpha_{k} \! = \! \infty$, $k \! \in \! 
\lbrace 1,2,\dotsc,K \rbrace$, can be proved analogously (see $\pmb{(2)}$ 
below).

$\pmb{(1)}$ Let $\widetilde{V} \colon \overline{\mathbb{R}} \setminus 
\lbrace \alpha_{1},\alpha_{2},\dotsc,\alpha_{K} \rbrace \! \to \! \mathbb{R}$ 
satisfy conditions~\eqref{eq20}--\eqref{eq22}. For $n \! \in \! \mathbb{N}$ 
and $k \! \in \! \lbrace 1,2,\dotsc,K \rbrace$ such that 
$\alpha_{p_{\mathfrak{s}}} \! := \! \alpha_{k} \! \neq \! \infty$, with 
$\mathcal{N} \! := \! (n \! - \! 1)K \! + \! k$, set
\begin{equation} \label{eql3.5a} 
\delta^{\widetilde{V},f}_{\mathcal{N}} \! = \sup_{\lbrace x_{1},x_{2},
\dotsc,x_{\mathcal{N}} \rbrace \subset \mathbb{R}} \left(\prod_{i<j}
(\mathfrak{h}(x_{i},x_{j}))^{2} \me^{-\widetilde{V}(x_{i})} 
\me^{-\widetilde{V}(x_{j})} \right)^{\frac{1}{\mathcal{N}(\mathcal{N}-1)}} = 
\sup_{\lbrace x_{1},x_{2},\dotsc,x_{\mathcal{N}} \rbrace \subset \mathbb{R}} 
\left(\prod_{i \neq j} \mathfrak{h}(x_{i},x_{j}) \me^{-(\mathcal{N}-1) 
\sum_{j=1}^{\mathcal{N}} \widetilde{V}(x_{j})} \right)^{\frac{1}{\mathcal{N}
(\mathcal{N}-1)}},
\end{equation}
where $x_{j} \! = \! x_{j}(n,k,z_{o})$, $\prod_{i<j}(\boldsymbol{\ast}) \! := \! 
\prod_{i=1}^{\mathcal{N}-1} \prod_{j=i+1}^{\mathcal{N}}(\boldsymbol{\ast})$, 
and $\mathfrak{h} \colon \mathbb{N} \times \lbrace 1,2,\dotsc,K 
\rbrace \times \mathbb{R}^{2} \! \ni \! (n,k,x,y) \! \mapsto \! \mathfrak{h}
(n,k,x,y) \! := \! \mathfrak{h}(x,y)$ $(= \! \mathfrak{h}(y,x))$, with
\begin{equation} \label{eql3.5b} 
\mathfrak{h}(x,y) \! := \! \lvert x \! - \! y \rvert^{\frac{\varkappa_{nk 
\tilde{k}_{\mathfrak{s}-1}}^{\infty}+1}{n}} \left(\dfrac{\lvert x \! - \! y 
\rvert}{\lvert x \! - \! \alpha_{k} \rvert \lvert y \! - \! \alpha_{k} \rvert} 
\right)^{\frac{\varkappa_{nk}-1}{n}} \prod_{q=1}^{\mathfrak{s}-2} \left(
\dfrac{\lvert x \! - \! y \rvert}{\lvert x \! - \! \alpha_{p_{q}} \rvert 
\lvert y \! - \! \alpha_{p_{q}} \rvert} \right)^{\frac{\varkappa_{nk 
\tilde{k}_{q}}}{n}}.
\end{equation}
One begins by showing that, for $n \! \in \! \mathbb{N}$ and $k \! \in \! 
\lbrace 1,2,\dotsc,K \rbrace$ such that $\alpha_{p_{\mathfrak{s}}} \! := \! 
\alpha_{k} \! \neq \! \infty$, $\delta^{\widetilde{V},f}_{\mathcal{N}}$ is finite, 
and the supremum is attained at some (or many) finite point set(s) $\lbrace 
\tilde{x}_{1}^{\ast},\tilde{x}_{2}^{\ast},\dotsc,\tilde{x}_{\mathcal{N}}^{\ast} 
\rbrace \subset \mathbb{R}$, with $\tilde{x}_{j}^{\ast} \! = \! 
\tilde{x}_{j}^{\ast}(n,k,z_{o})$, $j \! = \! 1,2,\dotsc,\mathcal{N}$; in fact, 
as shown below, $\lbrace \tilde{x}_{1}^{\ast},\tilde{x}_{2}^{\ast},\dotsc,
\tilde{x}_{\mathcal{N}}^{\ast} \rbrace \subset \overline{\mathbb{R}} 
\setminus \lbrace \alpha_{1},\alpha_{2},\dotsc,\alpha_{K} \rbrace$. Via the 
inequality $\lvert t_{1} \! - \! t_{2} \rvert \! \leqslant \! (1 \! + \! t_{1}^{2})^{1/2}
(1 \! + \! t_{2}^{2})^{1/2}$, $t_{1},t_{2} \! \in \! \mathbb{R}$, one shows that, 
for $n \! \in \! \mathbb{N}$ and $k \! \in \! \lbrace 1,2,\dotsc,K \rbrace$ 
such that $\alpha_{p_{\mathfrak{s}}} \! := \! \alpha_{k} \! \neq \! \infty$,
\begin{align*}
\prod_{i \neq j} \lvert x_{i} \! - \! x_{j} \rvert^{\varkappa^{\infty}_{nk 
\tilde{k}_{\mathfrak{s}-1}}+1} \leqslant& \, \left(\prod_{j=1}^{\mathcal{N}}
(1 \! + \! x_{j}^{2})^{\varkappa^{\infty}_{nk \tilde{k}_{\mathfrak{s}-1}}+1} 
\right)^{\mathcal{N}-1} \\
\prod_{i \neq j} \left(\dfrac{\lvert x_{i} \! - \! x_{j} \rvert}{\lvert x_{i} 
\! - \! \alpha_{k} \rvert \lvert x_{j} \! - \! \alpha_{k} \rvert} \right)^{
\varkappa_{nk}-1} \leqslant& \left(\prod_{j=1}^{\mathcal{N}}(1 \! + \! (x_{j} 
\! - \! \alpha_{k})^{-2})^{\varkappa_{nk}-1} \right)^{\mathcal{N}-1} \\
\prod_{i \neq j} \prod_{q=1}^{\mathfrak{s}-2} \left(\dfrac{\lvert x_{i} \! - 
\! x_{j} \rvert}{\lvert x_{i} \! - \! \alpha_{p_{q}} \rvert \lvert x_{j} \! 
- \! \alpha_{p_{q}} \rvert} \right)^{\varkappa_{nk \widetilde{k}_{q}}} 
\leqslant& \, \left(\prod_{j=1}^{\mathcal{N}} \prod_{q=1}^{\mathfrak{s}-2}
(1 \! + \! (x_{j} \! - \! \alpha_{p_{q}})^{-2})^{\varkappa_{nk \tilde{k}_{q}}} 
\right)^{\mathcal{N}-1};
\end{align*}
hence, via the above inequalities, one arrives at, for $n \! \in \! 
\mathbb{N}$ and $k \! \in \! \lbrace 1,2,\dotsc,K \rbrace$ such that 
$\alpha_{p_{\mathfrak{s}}} \! := \! \alpha_{k} \! \neq \! \infty$,
\begin{align*}
\prod_{i<j}(\mathfrak{h}(x_{i},x_{j}))^{2} \me^{-\widetilde{V}(x_{i})} 
\me^{-\widetilde{V}(x_{j})} =& \, \prod_{i<j} \mathfrak{h}(x_{i},x_{j}) 
\me^{-(\mathcal{N}-1) \sum_{j=1}^{\mathcal{N}} \widetilde{V}(x_{j})} \! 
\leqslant \! \prod_{j=1}^{\mathcal{N}} \left((1 \! + \! x_{j}^{2})^{
\frac{\varkappa^{\infty}_{nk \tilde{k}_{\mathfrak{s}-1}}+1}{n}} \right. \\
\times&\left. \, (1 \! + \! (x_{j} \! - \! \alpha_{k})^{-2})^{\frac{
\varkappa_{nk}-1}{n}} \prod_{q=1}^{\mathfrak{s}-2}(1 \! + \! (x_{j} \! - \! 
\alpha_{p_{q}})^{-2})^{\frac{\varkappa_{nk \tilde{k}_{q}}}{n}} \, 
\me^{-\widetilde{V}(x_{j})} \right)^{\mathcal{N}-1}.
\end{align*}
For $\widetilde{V} \colon \overline{\mathbb{R}} \setminus \lbrace \alpha_{1},
\alpha_{2},\dotsc,\alpha_{K} \rbrace \! \to \! \mathbb{R}$ satisfying 
conditions~\eqref{eq20}--\eqref{eq22}, one shows that, for $n \! \in \! 
\mathbb{N}$ and $k \! \in \! \lbrace 1,2,\dotsc,K \rbrace$ such that 
$\alpha_{p_{\mathfrak{s}}} \! := \! \alpha_{k} \! \neq \! \infty$: (i) for 
$x_{j} \! \in \! \overline{\mathbb{R}} \setminus \lbrace \alpha_{1},
\alpha_{2},\dotsc,\alpha_{K} \rbrace$, $j \! = \! 1,2,\dotsc,\mathcal{N}$,
\begin{equation*}
0 \! \leqslant \! (1 \! + \! x_{j}^{2})^{\frac{\varkappa^{\infty}_{nk 
\tilde{k}_{\mathfrak{s}-1}}+1}{n}}(1 \! + \! (x_{j} \! - \! 
\alpha_{k})^{-2})^{\frac{\varkappa_{nk}-1}{n}} \prod_{q=1}^{\mathfrak{s}-2}
(1 \! + \! (x_{j} \! - \! \alpha_{p_{q}})^{-2})^{\frac{\varkappa_{nk 
\tilde{k}_{q}}}{n}} \me^{-\widetilde{V}(x_{j})} \! \leqslant \! \mathfrak{c}_{f}
(n,k,z_{o}) \! < \! +\infty;
\end{equation*}
(ii) for $x_{j} \! \in \! \mathscr{O}_{\infty} \! := \! \lbrace \mathstrut 
x \! \in \! \mathbb{R}; \, \lvert x \rvert \! > \! \delta_{\infty}^{-1} 
\rbrace$, $j \! = \! 1,2,\dotsc,\mathcal{N}$, with $\delta_{\infty}$ 
$(= \! \delta_{\infty}(n,k,z_{o}))$ some arbitrarily fixed, sufficiently small 
positive real number (e.g., $\delta_{\infty}^{-1} \! = \! K(1 \! + \! 
\max \lbrace \mathstrut \lvert \alpha_{p_{q}} \rvert, \, q \! = \! 1,
\dotsc,\mathfrak{s} \! - \! 2,\mathfrak{s} \rbrace \! + \! \max \lbrace 
\mathstrut \lvert \alpha_{p_{l}} \! - \! \alpha_{p_{m}} \rvert, \, l \! \neq 
\! m \! \in \! \lbrace 1,\dotsc,\mathfrak{s} \! - \! 2,\mathfrak{s} \rbrace 
\rbrace))$, $\widetilde{V}(x_{j}) \! \geqslant \! (1 \! + \! \mathfrak{c}_{
\infty}) \ln (1 \! + \! x_{j}^{2})$, where $\mathfrak{c}_{\infty}$ $(= \! 
\mathfrak{c}_{\infty}(n,k,z_{o}))$ is some bounded, positive real number, it 
follows that $(1 \! + \! x_{j}^{2})^{\frac{\varkappa^{\infty}_{nk \tilde{k}_{
\mathfrak{s}-1}}+1}{n}} \exp (-\widetilde{V}(x_{j})) \! \to \! 0$ as $x_{j} 
\! \to \! \alpha_{p_{\mathfrak{s}-1}} \! = \! \infty$, $j \! = \! 1,2,
\dotsc,\mathcal{N}$, whence
\begin{equation*}
0 \! \leqslant \! (1 \! + \! x_{j}^{2})^{\frac{\varkappa^{\infty}_{nk 
\tilde{k}_{\mathfrak{s}-1}}+1}{n}}(1 \! + \! (x_{j} \! - \! 
\alpha_{k})^{-2})^{\frac{\varkappa_{nk}-1}{n}} \prod_{q=1}^{\mathfrak{s}
-2}(1 \! + \! (x_{j} \! - \! \alpha_{p_{q}})^{-2})^{\frac{\varkappa_{nk 
\tilde{k}_{q}}}{n}} \me^{-\widetilde{V}(x_{j})} \! \leqslant \! 
\mathfrak{c}_{f}(n,k,z_{o}), \quad x_{j} \! \in \! \mathscr{O}_{\infty};
\end{equation*}
and (iii) for $x_{j} \! \in \! \mathscr{O}_{\tilde{\delta}_{q}}(\alpha_{p_{q}}) 
\! := \! \lbrace \mathstrut x \! \in \! \mathbb{R}; \, \lvert x \! - \! 
\alpha_{p_{q}} \rvert \! < \! \tilde{\delta}_{q} \rbrace$, $j \! = \! 1,
2,\dotsc,\mathcal{N}$, $q \! = \! 1,\dotsc,\mathfrak{s} \! - \! 2,
\mathfrak{s}$, with $\tilde{\delta}_{q}$ $(= \! \tilde{\delta}_{q}(n,k,z_{o}))$ 
some arbitrarily fixed, sufficiently small positive real number chosen 
so that $\mathscr{O}_{\tilde{\delta}_{q^{\prime}}}(\alpha_{p_{q^{\prime}}}) 
\cap \mathscr{O}_{\tilde{\delta}_{q^{\prime \prime}}}(\alpha_{p_{q^{\prime 
\prime}}}) \! = \! \varnothing$ $\forall$ $q^{\prime} \! \neq \! q^{\prime 
\prime} \! \in \! \lbrace 1,\dotsc,\mathfrak{s} \! - \! 2,\mathfrak{s} 
\rbrace$ (e.g., $\tilde{\delta}_{q} \! < \! (3K)^{-1} \min \lbrace \mathstrut 
\lvert \alpha_{p_{l}} \! - \! \alpha_{p_{m}} \rvert, \, l \! \neq \! m \! \in \! 
\lbrace 1,\dotsc,\mathfrak{s} \! - \! 2,\mathfrak{s} \rbrace \rbrace)$, 
$\widetilde{V}(x_{j}) \! \geqslant \! (1 \! + \! \mathfrak{c}_{q}) \ln (1 \! 
+ \! (x_{j} \! - \! \alpha_{p_{q}})^{-2})$, where $\mathfrak{c}_{q}$ $(= 
\! \mathfrak{c}_{q}(n,k,z_{o}))$ is some bounded, positive real number, it 
follows that $(1 \! + \! (x_{j} \! - \! \alpha_{k})^{-2})^{\frac{\varkappa_{nk}
-1}{n}} \exp (-\widetilde{V}(x_{j})) \! \to \! 0$ as $x_{j} \! \to \! \alpha_{k}$, 
and $(1 \! + \! (x_{j} \! - \! \alpha_{p_{q}})^{-2})^{\frac{\varkappa_{nk 
\tilde{k}_{q}}}{n}} \exp (-\widetilde{V}(x_{j})) \! \to \! 0$ as $x_{j} \! \to 
\! \alpha_{p_{q}}$, $j \! = \! 1,2,\dotsc,\mathcal{N}$, $q \! = \! 1,2,
\dotsc,\mathfrak{s} \! - \! 2$, whence
\begin{equation*}
0 \! \leqslant \! (1 \! + \! x_{j}^{2})^{\frac{\varkappa^{\infty}_{nk 
\tilde{k}_{\mathfrak{s}-1}}+1}{n}}(1 \! + \! (x_{j} \! - \! 
\alpha_{k})^{-2})^{\frac{\varkappa_{nk}-1}{n}} \prod_{q=1}^{\mathfrak{s}
-2}(1 \! + \! (x_{j} \! - \! \alpha_{p_{q}})^{-2})^{\frac{\varkappa_{nk 
\tilde{k}_{q}}}{n}} \me^{-\widetilde{V}(x_{j})} \! \leqslant \! \mathfrak{c}_{f}
(n,k,z_{o}), \quad x_{j} \! \in \! \mathscr{O}_{\tilde{\delta}_{q}}(\alpha_{p_{q}}).
\end{equation*}
Hence, for $n \! \in \! \mathbb{N}$ and $k \! \in \! \lbrace 1,2,\dotsc,K 
\rbrace$ such that $\alpha_{p_{\mathfrak{s}}} \! := \! \alpha_{k} \! \neq 
\! \infty$, with $x_{m} \! \in \! \mathbb{R}$, $m \! = \! 1,2,\dotsc,
\mathcal{N}$,
\begin{align*}
\prod_{i<j}(\mathfrak{h}(x_{i},x_{j}))^{2} \me^{-\widetilde{V}(x_{i})} 
\me^{-\widetilde{V}(x_{j})} =& \, \prod_{i \neq j} \mathfrak{h}(x_{i},x_{j}) 
\me^{-(\mathcal{N}-1) \sum_{j=1}^{\mathcal{N}} \widetilde{V}(x_{j})} 
\! \leqslant \! (\mathfrak{c}_{f}(n,k,z_{o}))^{(\mathcal{N}-1)^{2}} \left(
(1 \! + \! x_{m}^{2})^{\frac{\varkappa_{nk \tilde{k}_{\mathfrak{s}-1}}^{
\infty}+1}{n}} \right. \\
\times&\left. \, (1 \! + \! (x_{m} \! - \! \alpha_{k})^{-2})^{\frac{
\varkappa_{nk}-1}{n}} \prod_{q=1}^{\mathfrak{s}-2}(1 \! + \! (x_{m} 
\! - \! \alpha_{p_{q}})^{-2})^{\frac{\varkappa_{nk \tilde{k}_{q}}}{n}} 
\me^{-\widetilde{V}(x_{m})} \right)^{\mathcal{N}-1} \\
\to& \, 0 \, \, \, \text{as} \, \, \, x_{m} \! \to \! \alpha_{p_{\mathfrak{s}-1}} 
\! = \! \infty \, \, \, \text{and as} \, \, \, x_{m} \! \to \! \alpha_{p_{q}}, 
\, \, \, q \! = \! 1,\dotsc,\mathfrak{s} \! - \! 2,\mathfrak{s};
\end{align*}
hence, for $n \! \in \! \mathbb{N}$ and $k \! \in \! \lbrace 1,2,\dotsc,K 
\rbrace$ such that $\alpha_{p_{\mathfrak{s}}} \! := \! \alpha_{k} \! \neq 
\! \infty$, $\delta^{\widetilde{V},f}_{\mathcal{N}}$ is finite, and the 
supremum is attained at some (or many) finite point set(s) $\lbrace 
\tilde{x}_{1}^{\ast},\tilde{x}_{2}^{\ast},\dotsc,\tilde{x}_{\mathcal{N}}^{\ast} 
\rbrace \subset \mathbb{R}$, with $\tilde{x}_{j}^{\ast} \! = \! 
\tilde{x}_{j}^{\ast}(n,k,z_{o})$, $j \! = \! 1,2,\dotsc,\mathcal{N}$; in fact, 
as shown above, $\lbrace \tilde{x}_{j}^{\ast} \rbrace_{j=1}^{\mathcal{N}} 
\subset \overline{\mathbb{R}} \setminus \lbrace \alpha_{1},\alpha_{2},
\dotsc,\alpha_{K} \rbrace$. For $n \! \in \! \mathbb{N}$ and $k \! \in 
\! \lbrace 1,2,\dotsc,K \rbrace$ such that $\alpha_{p_{\mathfrak{s}}} \! 
:= \! \alpha_{k} \! \neq \! \infty$, one calls a maximising set, that is, 
a set $\lbrace \tilde{x}_{1}^{\ast},\tilde{x}_{2}^{\ast},\dotsc,\tilde{x}_{
\mathcal{N}}^{\ast} \rbrace$ for which
\begin{equation*}
\delta^{\widetilde{V},f}_{\mathcal{N}} \! = \! \left(\prod_{i<j}(\mathfrak{h}
(\tilde{x}_{i}^{\ast},\tilde{x}_{j}^{\ast}))^{2} \me^{-\widetilde{V}
(\tilde{x}_{i}^{\ast})} \me^{-\widetilde{V}(\tilde{x}_{j}^{\ast})} 
\right)^{\frac{1}{\mathcal{N}(\mathcal{N}-1)}} \! = \! \left(\prod_{i \neq j} 
\mathfrak{h}(\tilde{x}_{i}^{\ast},\tilde{x}_{j}^{\ast}) \me^{-(\mathcal{N}-1) 
\sum_{j=1}^{\mathcal{N}} \widetilde{V}(\tilde{x}_{j}^{\ast})} 
\right)^{\frac{1}{\mathcal{N}(\mathcal{N}-1)}},
\end{equation*}
an associated, weighted $\mathcal{N}$-Fekete set, and the $\mathcal{N}$ 
points $\tilde{x}_{1}^{\ast},\tilde{x}_{2}^{\ast},\dotsc,\tilde{x}_{\mathcal{N}}^{
\ast}$ will be called the associated, weighted Fekete points. For $n \! \in 
\! \mathbb{N}$ and $k \! \in \! \lbrace 1,2,\dotsc,K \rbrace$ such that 
$\alpha_{p_{\mathfrak{s}}} \! := \! \alpha_{k} \! \neq \! \infty$, let
\begin{equation} \label{eql3.5c} 
\mathscr{K}^{\widetilde{V},f}_{\mathcal{N}}(x_{1},x_{2},\dotsc,
x_{\mathcal{N}}) \! := \! \sum_{\substack{i,j=1\\i \neq j}}^{\mathcal{N}} 
\mathcal{K}_{\widetilde{V}}^{f}(x_{i},x_{j}),
\end{equation}
where $\mathcal{K}^{f}_{\widetilde{V}}(\xi,\tau)$ is given in 
Equation~\eqref{eqKvinf1}. A calculation shows that
\begin{align} \label{eql3.5d} 
\mathscr{K}^{\widetilde{V},f}_{\mathcal{N}}(x_{1},x_{2},\dotsc,
x_{\mathcal{N}}) =& \, \sum_{\substack{i,j=1\\i \neq j}}^{\mathcal{N}} 
\left(\left(\dfrac{\varkappa_{nk \tilde{k}_{\mathfrak{s}-1}}^{\infty} \! 
+ \! 1}{n} \right) \ln \left(\dfrac{1}{\lvert x_{i} \! - \! x_{j} \rvert} \right) \! 
+ \! \left(\dfrac{\varkappa_{nk} \! - \! 1}{n} \right) \ln \left(\dfrac{\lvert 
x_{i} \! - \! \alpha_{k} \rvert \lvert x_{j} \! - \! \alpha_{k} \rvert}{\lvert 
x_{i} \! - \! x_{j} \rvert} \right) \right. \nonumber \\
+&\left. \, \sum_{q=1}^{\mathfrak{s}-2} \dfrac{\varkappa_{nk 
\tilde{k}_{q}}}{n} \ln \left(\dfrac{\lvert x_{i} \! - \! \alpha_{p_{q}} \rvert 
\lvert x_{j} \! - \! \alpha_{p_{q}} \rvert}{\lvert x_{i} \! - \! x_{j} \rvert} 
\right) \right) \! + \! (\mathcal{N} \! - \! 1) \sum_{j=1}^{\mathcal{N}} 
\widetilde{V}(x_{j}).
\end{align}
For $n \! \in \! \mathbb{N}$ and $k \! \in \! \lbrace 1,2,\dotsc,K 
\rbrace$ such that $\alpha_{p_{\mathfrak{s}}} \! := \! \alpha_{k} \! 
\neq \! \infty$, set, as in the lemma,
\begin{equation} \label{eql3.5e} 
\mathfrak{d}^{\widetilde{V},f}_{\mathcal{N}} \! := \! (\mathcal{N}(\mathcal{N} 
\! - \! 1))^{-1} \inf_{\lbrace x_{1},x_{2},\dotsc,x_{\mathcal{N}} \rbrace 
\subset \mathbb{R}} \mathscr{K}^{\widetilde{V},f}_{\mathcal{N}}(x_{1},
x_{2},\dotsc,x_{\mathcal{N}}).
\end{equation}
One shows that, for $n \! \in \! \mathbb{N}$ and $k \! \in \! \lbrace 
1,2,\dotsc,K \rbrace$ such that $\alpha_{p_{\mathfrak{s}}} \! := \! 
\alpha_{k} \! \neq \! \infty$, via Equation~\eqref{eql3.5d},
\begin{align*}
\me^{-\mathscr{K}^{\widetilde{V},f}_{\mathcal{N}}(x_{1},x_{2},\dotsc,
x_{\mathcal{N}})} =& \, \prod_{i<j} \left(\lvert x_{i} \! - \! x_{j} 
\rvert^{\frac{\varkappa_{nk \tilde{k}_{\mathfrak{s}-1}}^{\infty}+1}{n}} 
\left(\dfrac{\lvert x_{i} \! - \! x_{j} \rvert}{\lvert x_{i} \! - \! 
\alpha_{k} \rvert \lvert x_{j} \! - \! \alpha_{k} \rvert} \right)^{
\frac{\varkappa_{nk}-1}{n}} \prod_{q=1}^{\mathfrak{s}-2} \left(
\dfrac{\lvert x_{i} \! - \! x_{j} \rvert}{\lvert x_{i} \! - \! \alpha_{p_{q}} 
\rvert \lvert x_{j} \! - \! \alpha_{p_{q}} \rvert} \right)^{\frac{
\varkappa_{nk \tilde{k}_{q}}}{n}} \right)^{2} \me^{-(\mathcal{N}-1) 
\sum_{j=1}^{\mathcal{N}} \widetilde{V}(x_{j})} \\
=& \, \prod_{i<j}(\mathfrak{h}(x_{i},x_{j}))^{2} \me^{-(\mathcal{N}-1) 
\sum_{j=1}^{\mathcal{N}} \widetilde{V}(x_{j})} \qquad \Rightarrow
\end{align*}
\begin{equation*}
\me^{-(\mathcal{N}(\mathcal{N}-1))^{-1} \mathscr{K}^{\widetilde{V},f}_{
\mathcal{N}}(x_{1},x_{2},\dotsc,x_{\mathcal{N}})} \! = \! \left(\prod_{i<j}
(\mathfrak{h}(x_{i},x_{j}))^{2} \me^{-(\mathcal{N}-1) \sum_{j=1}^{\mathcal{N}} 
\widetilde{V}(x_{j})} \right)^{(\mathcal{N}(\mathcal{N}-1))^{-1}},
\end{equation*}
where $\mathfrak{h}(x,y)$ is defined by Equation~\eqref{eql3.5b}, whence, 
via Equation~\eqref{eql3.5a}, the Definition~\eqref{eql3.5e}, and the fact 
that $\inf (\boldsymbol{\ast}) \! = \! -\sup (-\boldsymbol{\ast})$,
\begin{equation*}
\me^{-\mathfrak{d}^{\widetilde{V},f}_{\mathcal{N}}} \! = \! \sup_{\lbrace 
x_{1},x_{2},\dotsc,x_{\mathcal{N}} \rbrace \subset \mathbb{R}} \left(
\prod_{i<j}(\mathfrak{h}(x_{i},x_{j}))^{2} \me^{-(\mathcal{N}-1) 
\sum_{j=1}^{\mathcal{N}} \widetilde{V}(x_{j})} \right)^{(\mathcal{N}
(\mathcal{N}-1))^{-1}} \! = \! \delta^{\widetilde{V},f}_{\mathcal{N}}.
\end{equation*}
The above calculations also show that, for $n \! \in \! \mathbb{N}$ 
and $k \! \in \! \lbrace 1,2,\dotsc,K \rbrace$ such that 
$\alpha_{p_{\mathfrak{s}}} \! := \! \alpha_{k} \! \neq \! \infty$,
\begin{equation*}
\me^{-\mathscr{K}^{\widetilde{V},f}_{\mathcal{N}}(x_{1},x_{2},\dotsc,
x_{\mathcal{N}})} \! \leqslant \! \left(\prod_{j=1}^{\mathcal{N}}(1 \! + \! 
x_{j}^{2})^{\frac{\varkappa_{nk \tilde{k}_{\mathfrak{s}-1}}^{\infty}+1}{n}}
(1 \! + \! (x_{j} \! - \! \alpha_{k})^{-2})^{\frac{\varkappa_{nk}-1}{n}} 
\prod_{q=1}^{\mathfrak{s}-2}(1 \! + \! (x_{j} \! - \! \alpha_{p_{q}})^{-2})^{
\frac{\varkappa_{nk \tilde{k}_{q}}}{n}} \me^{-\widetilde{V}(x_{j})} 
\right)^{\mathcal{N}-1},
\end{equation*}
whence
\begin{align*}
\me^{-\mathfrak{d}^{\widetilde{V},f}_{\mathcal{N}}} \! \leqslant \! 
\sup_{\lbrace x_{1},x_{2},\dotsc,x_{\mathcal{N}} \rbrace \subset 
\mathbb{R}} \left(\prod_{j=1}^{\mathcal{N}} \left((1 \! + \! x_{j}^{2})^{
\frac{\varkappa_{nk \tilde{k}_{\mathfrak{s}-1}}^{\infty}+1}{n}}(1 \! + \! 
(x_{j} \! - \! \alpha_{k})^{-2})^{\frac{\varkappa_{nk}-1}{n}} \prod_{q=1}^{
\mathfrak{s}-2}(1 \! + \! (x_{j} \! - \! \alpha_{p_{q}})^{-2})^{\frac{
\varkappa_{nk \tilde{k}_{q}}}{n}} \me^{-\widetilde{V}(x_{j})} \right)^{
\mathcal{N}-1} \right)^{(\mathcal{N}(\mathcal{N}-1))^{-1}} \! < \! +\infty;
\end{align*}
hence, for $n \! \in \! \mathbb{N}$ and $k \! \in \! \lbrace 1,2,\dotsc,
K \rbrace$ such that $\alpha_{p_{\mathfrak{s}}} \! := \! \alpha_{k} \! \neq 
\! \infty$, $\mathfrak{d}^{\widetilde{V},f}_{\mathcal{N}} \! > \! -\infty$, 
that is, $(\mathcal{N}(\mathcal{N} \! - \! 1))^{-1} \inf_{\lbrace x_{1},
x_{2},\dotsc,x_{\mathcal{N}} \rbrace \subset \mathbb{R}} \mathscr{K}^{
\widetilde{V},f}_{\mathcal{N}}(x_{1},\linebreak[4]
x_{2},\dotsc,x_{\mathcal{N}}) \! > \! -\infty$.

Since it's been established above that, for $n \! \in \! \mathbb{N}$ 
and $k \! \in \! \lbrace 1,2,\dotsc,K \rbrace$ such that 
$\alpha_{p_{\mathfrak{s}}} \! := \! \alpha_{k} \! \neq \! \infty$, an 
associated, weighted $\mathcal{N}$-Fekete set for $\mathscr{K}^{
\widetilde{V},f}_{\mathcal{N}}(x_{1},x_{2},\dotsc,x_{\mathcal{N}})$ 
exists, that is, a set $\lbrace \tilde{x}_{1}^{\ast},\tilde{x}_{2}^{\ast},\dotsc,
\tilde{x}_{\mathcal{N}}^{\ast} \rbrace \subset \overline{\mathbb{R}} 
\setminus \lbrace \alpha_{1},\alpha_{2},\dotsc,\alpha_{K} \rbrace$ such 
that $\mathfrak{d}^{\widetilde{V},f}_{\mathcal{N}} \! = \! (\mathcal{N}
(\mathcal{N} \! - \! 1))^{-1} \mathscr{K}^{\widetilde{V},f}_{\mathcal{N}}
(\tilde{x}_{1}^{\ast},\tilde{x}_{2}^{\ast},\dotsc,\tilde{x}_{\mathcal{N}}^{
\ast})$, define the corresponding normalised counting measure for the 
associated, weighted $\mathcal{N}$-Fekete set as follows:
\begin{equation} \label{eql3.5f} 
\lambda_{\mathcal{N}}^{f} \! := \! \dfrac{1}{\mathcal{N}} 
\sum_{j=1}^{\mathcal{N}} \delta_{\tilde{x}_{j}^{\ast}},
\end{equation}
where $\delta_{\tilde{x}_{j}^{\ast}}$ is the Dirac delta (atomic) 
mass concentrated at $\tilde{x}_{j}^{\ast}$, $j \! = \! 1,2,\dotsc,
\mathcal{N}$.\footnote{Note that $\int_{\mathbb{R}} \md 
\lambda_{\mathcal{N}}^{f}(\xi) \! = \! \tfrac{1}{\mathcal{N}} 
\int_{\mathbb{R}} \sum_{j=1}^{\mathcal{N}} \delta (\xi \! - \! 
\tilde{x}_{j}^{\ast}) \, \md \xi \! = \! 1$.} One shows that, for $n \! 
\in \! \mathbb{N}$ and $k \! \in \! \lbrace 1,2,\dotsc,K \rbrace$ 
such that $\alpha_{p_{\mathfrak{s}}} \! := \! \alpha_{k} \! \neq \! 
\infty$, in the weak-$\ast$ topology of measures, $\lambda^{f}_{
\mathcal{N}}$ converges weakly to $\mu^{f}_{\widetilde{V}}$ 
$(\in \! \mathscr{M}_{1}(\mathbb{R}))$ in the double-scaling limit 
$\mathscr{N},n \! \to \! \infty$ such that $z_{o} \! = \! 1 \! + \! 
o(1)$. As in the proof of Lemma~3.5 in \cite{a45}, one proceeds via 
a contradiction argument, namely, one assumes that, for $k \! \in 
\! \lbrace 1,2,\dotsc,K \rbrace$ such that $\alpha_{p_{\mathfrak{s}}} 
\! := \! \alpha_{k} \! \neq \! \infty$ and (some) $\mathfrak{g} 
\! \in \! \pmb{\operatorname{C}}^{0}_{\text{b}}(\mathbb{R})$, 
$\int_{\mathbb{R}} \mathfrak{g}(\xi) \, \md \lambda^{f}_{
\mathcal{N}}(\xi) \! \nrightarrow \! \int_{\mathbb{R}} 
\mathfrak{g}(\xi) \, \md \mu^{f}_{\widetilde{V}}(\xi)$ in the 
double-scaling limit $\mathscr{N},n \! \to \! \infty$ such that 
$z_{o} \! = \! 1 \! + \! o(1)$ (for simplicity of notation, only 
$n \! \to \! \infty$ will be written below), that is, there exists 
$\varepsilon^{\flat}$ $(= \! \varepsilon^{\flat}(n,k,z_{o}))$ $> \! 0$ 
and a subsequence {}\footnote{Strictly speaking, subsequences 
$\mathscr{N}_{\hat{k}},n_{\hat{k}} \! \to \! \infty$ such that 
$z_{o,\hat{k}} \! := \! \mathscr{N}_{\hat{k}}/n_{\hat{k}} \! = \! 
1 \! + \! o(1)$.} $n_{\hat{k}} \! \to \! \infty$ (with $\mathcal{N}_{
\hat{k}} \! := \! (n_{\hat{k}} \! - \! 1)K \! + \! k)$ such that, for 
all $\hat{k} \! \in \! \mathbb{N}$, $\lvert \int_{\mathbb{R}} 
\mathfrak{g}(\xi) \, \md \lambda^{f}_{\mathcal{N}_{\hat{k}}}
(\xi) \! - \! \int_{\mathbb{R}} \mathfrak{g}(\xi) \, \md \mu^{f}_{
\widetilde{V}}(\xi) \rvert \! \geqslant \! \varepsilon^{\flat}$; but first, 
one must establish that, for $n \! \in \! \mathbb{N}$ and $k \! \in \! 
\lbrace 1,2,\dotsc,K \rbrace$ such that $\alpha_{p_{\mathfrak{s}}} 
\! := \! \alpha_{k} \! \neq \! \infty$, the corresponding sequence 
of probability measures $\lbrace \lambda^{f}_{\mathcal{N}} 
\rbrace_{\mathcal{N} \in \mathbb{N}}$ is tight (cf. the proof 
of Lemma~\ref{lem3.1}). Since, for $n \! \in \! \mathbb{N}$ 
and $k \! \in \! \lbrace 1,2,\dotsc,K \rbrace$ such that 
$\alpha_{p_{\mathfrak{s}}} \! := \! \alpha_{k} \! \neq \! \infty$, 
$(\mathcal{N}(\mathcal{N} \! - \! 1))^{-1} \inf_{\lbrace x_{1},
x_{2},\dotsc,x_{\mathcal{N}} \rbrace \subset \mathbb{R}} 
\mathscr{K}^{\widetilde{V},f}_{\mathcal{N}}(x_{1},x_{2},\dotsc,
x_{\mathcal{N}}) \! \leqslant \! (\mathcal{N}(\mathcal{N} \! - \! 
1))^{-1} \mathscr{K}^{\widetilde{V},f}_{\mathcal{N}}(x_{1},x_{2},
\dotsc,x_{\mathcal{N}}) \! \leqslant \! (\mathcal{N}(\mathcal{N} 
\! - \! 1))^{-1} \sup_{\lbrace x_{1},x_{2},\dotsc,x_{\mathcal{N}} 
\rbrace \subset \mathbb{R}} \mathscr{K}^{\widetilde{V},f}_{
\mathcal{N}}(x_{1},x_{2},\dotsc,x_{\mathcal{N}})$, it follows {}from 
the Definition~\eqref{eql3.5e} that $\mathfrak{d}^{\widetilde{V},f}_{
\mathcal{N}} \! \leqslant \! (\mathcal{N}(\mathcal{N} \! - \! 1))^{-1} 
\mathscr{K}^{\widetilde{V},f}_{\mathcal{N}}(x_{1},x_{2},\dotsc,
x_{\mathcal{N}})$: integrating both sides of this latter inequality 
with respect to the `product measure' (cf. Remark~\ref{remm3.3}) 
$\md \mu^{f}_{\widetilde{V}}(x_{1}) \, \md \mu^{f}_{\widetilde{V}}(x_{2}) 
\, \dotsb \, \md \mu^{f}_{\widetilde{V}}(x_{\mathcal{N}})$, one shows 
that, for $n \! \in \! \mathbb{N}$ and $k \! \in \! \lbrace 1,2,\dotsc,
K \rbrace$ such that $\alpha_{p_{\mathfrak{s}}} \! := \! \alpha_{k} 
\! \neq \! \infty$,
\begin{align*}
& \, \iint_{\mathbb{R}^{2}} \ln \left(\lvert \xi \! - \! \tau 
\rvert^{-\frac{(\varkappa_{nk \tilde{k}_{\mathfrak{s}-1}}^{\infty}
+1)}{n}} \left(\dfrac{\lvert \xi \! - \! \alpha_{k} \rvert \lvert \tau \! 
- \! \alpha_{k} \rvert}{\lvert \xi \! - \! \tau \rvert} \right)^{\frac{
\varkappa_{nk}-1}{n}} \prod_{q=1}^{\mathfrak{s}-2} \left(\dfrac{\lvert \xi 
\! - \! \alpha_{p_{q}} \rvert \lvert \tau \! - \! \alpha_{p_{q}} \rvert}{\lvert 
\xi \! - \! \tau \rvert} \right)^{\frac{\varkappa_{nk \tilde{k}_{q}}}{n}} 
\me^{\widetilde{V}(\xi)/2} \me^{\widetilde{V}(\tau)/2} \right) \md 
\mu^{f}_{\widetilde{V}}(\xi) \, \md \mu_{\widetilde{V}}^{f}(\tau) \\
&= \, \iint_{\mathbb{R}^{2}} \mathcal{K}_{\widetilde{V}}^{f}(\xi,\tau) 
\, \md \mu^{f}_{\widetilde{V}}(\xi) \, \md \mu^{f}_{\widetilde{V}}(\tau) 
\! = \! \mathrm{I}_{\widetilde{V}}^{f}[\mu_{\widetilde{V}}^{f}] \! = \! 
E_{\widetilde{V}}^{f} \! \geqslant \! \mathfrak{d}^{\widetilde{V},f}_{
\mathcal{N}} \quad (> \! -\infty).
\end{align*}
Via this latter inequality and the Definitions~\eqref{eqKvinf1} 
and~\eqref{eqKvinf4}, one proceeds thus; for $n \! \in \! 
\mathbb{N}$ and $k \! \in \! \lbrace 1,2,\dotsc,K \rbrace$ such 
that $\alpha_{p_{\mathfrak{s}}} \! := \! \alpha_{k} \! \neq \! \infty$:
\begin{align*}
E^{f}_{\widetilde{V}} \! \geqslant \! \mathfrak{d}^{\widetilde{V},f}_{
\mathcal{N}} =& \, (\mathcal{N}(\mathcal{N} \! - \! 1))^{-1} \inf_{
\lbrace x_{1},x_{2},\dotsc,x_{\mathcal{N}} \rbrace \subset \mathbb{R}} 
\mathscr{K}^{\widetilde{V},f}_{\mathcal{N}}(x_{1},x_{2},\dotsc,
x_{\mathcal{N}}) \! = \! (\mathcal{N}(\mathcal{N} \! - \! 1))^{-1} 
\mathscr{K}^{\widetilde{V},f}_{\mathcal{N}}(\tilde{x}_{1}^{\ast},
\tilde{x}_{2}^{\ast},\dotsc,\tilde{x}_{\mathcal{N}}^{\ast}) \\
=& \, (\mathcal{N}(\mathcal{N} \! - \! 1))^{-1} \sum_{\substack{i,j=
1\\i \neq j}}^{\mathcal{N}} \mathcal{K}^{f}_{\widetilde{V}}(\tilde{x}_{i}^{
\ast},\tilde{x}_{j}^{\ast}) \! \geqslant \! \dfrac{1}{2}(\mathcal{N}(\mathcal{N} 
\! - \! 1))^{-1} \underbrace{\sum_{\substack{i,j=1\\i \neq j}}^{\mathcal{N}} 
\left(\hat{\psi}^{f}_{\widetilde{V}}(\tilde{x}_{i}^{\ast}) \! + \! \hat{\psi}^{f}_{
\widetilde{V}}(\tilde{x}_{j}^{\ast}) \right)}_{= \, 2(\mathcal{N}-1) \sum_{j=
1}^{\mathcal{N}} \hat{\psi}^{f}_{\widetilde{V}}(\tilde{x}_{j}^{\ast})} \\
=& \, \mathcal{N}^{-1} \sum_{j=1}^{\mathcal{N}} \hat{\psi}^{f}_{
\widetilde{V}}(\tilde{x}_{j}^{\ast}) \! = \! \int_{\mathbb{R}} \hat{\psi}^{f}_{
\widetilde{V}}(\xi) \, \md \lambda^{f}_{\mathcal{N}}(\xi) \quad \Rightarrow 
\quad \mathrm{I}^{f}_{\widetilde{V}}[\mu^{f}_{\widetilde{V}}] \! = \! 
E^{f}_{\widetilde{V}} \! \geqslant \! \mathfrak{d}^{\widetilde{V},f}_{
\mathcal{N}} \! \geqslant \! \int_{\mathbb{R}} \hat{\psi}^{f}_{\widetilde{V}}
(\xi) \, \md \lambda^{f}_{\mathcal{N}}(\xi).
\end{align*}
One now proceeds as in the proof of Lemma~\ref{lem3.1} to show that, for 
$n \! \in \! \mathbb{N}$ and $k \! \in \! \lbrace 1,2,\dotsc,K \rbrace$ such 
that $\alpha_{p_{\mathfrak{s}}} \! := \! \alpha_{k} \! \neq \! \infty$, the 
sequence of probability measures $\lbrace \lambda^{f}_{\mathcal{N}} 
\rbrace_{\mathcal{N} \in \mathbb{N}}$ (in $\mathscr{M}_{1}(\mathbb{R}))$ 
is tight, that is, for (some) $\hat{\varepsilon}$ $(= \! \hat{\varepsilon}(n,k,
z_{o}))$ $\! > \! 0$ and sufficiently small, $\limsup_{\underset{z_{o}=
1+o(1)}{\mathscr{N},n \to \infty}} \int_{\mathfrak{D}_{M_{f}}} \md 
\lambda^{f}_{\mathcal{N}}(\xi) \! \leqslant \! \hat{\varepsilon}$, where 
$\mathfrak{D}_{M_{f}} \! := \! \lbrace \lvert x \rvert \! \geqslant \! M_{f} 
\rbrace \cup \cup_{\underset{q \neq \mathfrak{s}-1}{q=1}}^{\mathfrak{s}} 
\operatorname{clos}(\mathscr{O}_{\frac{1}{M_{f}}}(\alpha_{p_{q}}))$, with 
$M_{f}$ $(= \! M_{f}(n,k,z_{o}))$ $> \! 1$ chosen so that $\mathscr{O}_{
\frac{1}{M_{f}}}(\alpha_{p_{i}}) \cap \mathscr{O}_{\frac{1}{M_{f}}}(\alpha_{p_{j}}) 
\! = \! \varnothing$, $i \! \neq \! j \! \in \! \lbrace 1,\dotsc,\mathfrak{s} \! 
- \! 2,\mathfrak{s} \rbrace$ (e.g., $M_{f} \! = \! K(1 \! + \! \max \lbrace 
\mathstrut \lvert \alpha_{p_{q}} \rvert, \, q \! = \! 1,\dotsc,\mathfrak{s} \! - \! 
2,\mathfrak{s} \rbrace \! + \! 3(\min \lbrace \mathstrut \lvert \alpha_{p_{i}} 
\! - \! \alpha_{p_{j}} \rvert, \, i \! \neq \! j \! \in \! \lbrace 1,\dotsc,
\mathfrak{s} \! - \! 2,\mathfrak{s} \rbrace \rbrace)^{-1}))$. For $n \! \in 
\! \mathbb{N}$ and $k \! \in \! \lbrace 1,2,\dotsc,K \rbrace$ such that 
$\alpha_{p_{\mathfrak{s}}} \! := \! \alpha_{k} \! \neq \! \infty$, with 
$L_{f} \! \in \! \mathbb{R}$:
\begin{align} \label{eql3.5g} 
\mathfrak{d}^{\widetilde{V},f}_{\mathcal{N}} =& \, (\mathcal{N}
(\mathcal{N} \! - \! 1))^{-1} \sum_{\substack{i,j=1\\i \neq j}}^{
\mathcal{N}} \mathcal{K}^{f}_{\widetilde{V}}(\tilde{x}_{i}^{\ast},
\tilde{x}_{j}^{\ast}) \! \geqslant \! (\mathcal{N}(\mathcal{N} \! - \! 1))^{-1} 
\sum_{\substack{i,j=1\\i \neq j}}^{\mathcal{N}} \min \left\lbrace L_{f},
\mathcal{K}^{f}_{\widetilde{V}}(\tilde{x}_{i}^{\ast},\tilde{x}_{j}^{\ast}) 
\right\rbrace \nonumber \\
=& \, (\mathcal{N}(\mathcal{N} \! - \! 1))^{-1} \left(\mathcal{N}^{2} 
\iint_{\mathbb{R}^{2}} \min \left\lbrace L_{f},\mathcal{K}^{f}_{
\widetilde{V}}(\xi,\tau) \right\rbrace \underbrace{\dfrac{1}{\mathcal{N}} 
\sum_{i=1}^{\mathcal{N}} \delta (\xi \! - \! \tilde{x}_{i}^{\ast}) \, \md 
\xi}_{= \, \md \lambda^{f}_{\mathcal{N}}(\xi)} \, \underbrace{\dfrac{1}{
\mathcal{N}} \sum_{j=1}^{\mathcal{N}} \delta (\tau \! - \! \tilde{x}_{j}^{
\ast}) \, \md \tau}_{= \, \md \lambda^{f}_{\mathcal{N}}(\tau)} - 
\mathcal{N}L_{f} \right) \nonumber \\
=& \, \dfrac{\mathcal{N}^{2}}{\mathcal{N}(\mathcal{N} \! - \! 1)} 
\iint_{\mathbb{R}^{2}} \min \left\lbrace L_{f},\mathcal{K}^{f}_{
\widetilde{V}}(\xi,\tau) \right\rbrace \, \md \lambda^{f}_{\mathcal{N}}(\xi) 
\, \md \lambda^{f}_{\mathcal{N}}(\tau) \! - \! (\mathcal{N} \! - \! 1)^{-1}
L_{f}.
\end{align}
For $n \! \in \! \mathbb{N}$ and $k \! \in \! \lbrace 1,2,\dotsc,K \rbrace$ 
such that $\alpha_{p_{\mathfrak{s}}} \! := \! \alpha_{k} \! \neq \! \infty$, 
since, as a consequence of tightness, the subsequence of probability 
measures $\lbrace \lambda^{f}_{\mathcal{N}_{\hat{k}}} \rbrace_{\hat{k}
=1}^{\infty}$ (in $\mathscr{M}_{1}(\mathbb{R}))$, with 
$\mathcal{N}_{\hat{k}} \! := \! (n_{\hat{k}} \! - \! 1)K \! + \! k$, is also 
tight, there exists, by a Helly selection theorem, a further weak-$\ast$ 
convergent subsequence of probability measures $\lbrace \lambda^{f}_{
\mathcal{N}_{\hat{k}_{j}}} \rbrace_{j=1}^{\infty}$ (in $\mathscr{M}_{1}
(\mathbb{R}))$, with $\mathcal{N}_{\hat{k}_{j}} \! := \! (n_{\hat{k}_{j}} 
\! - \! 1)K \! + \! k$, converging {}\footnote{In the double-scaling limit 
$\mathscr{N}_{\hat{k}_{j}},n_{\hat{k}_{j}} \! \to \! \infty$ such that 
$z_{o,\hat{k}_{j}} \! := \! \mathscr{N}_{\hat{k}_{j}}/n_{\hat{k}_{j}} \! = \! 
1 \! + \! o(1)$.} weakly to a probability measure $\lambda_{f} \! \in \! 
\mathscr{M}_{1}(\mathbb{R})$. Proceeding, now, as in the proof of 
Lemma~\ref{lem3.1}, one shows, via inequality~\eqref{eql3.5g}, that 
$\mathcal{N}_{\hat{k}_{j}}^{2}(\mathcal{N}_{\hat{k}_{j}}(\mathcal{N}_{
\hat{k}_{j}} \! - \! 1))^{-1} \iint_{\mathbb{R}^{2}} \min \lbrace L_{f},
\mathcal{K}^{f}_{\widetilde{V}}(\xi,\tau) \rbrace \, \md \lambda^{f}_{
\mathcal{N}_{\hat{k}_{j}}}(\xi) \md \lambda^{f}_{\mathcal{N}_{\hat{k}_{j}}}
(\tau) \! - \! (\mathcal{N}_{\hat{k}_{j}} \! - \! 1)^{-1}L_{f} \! \to \! 
\iint_{\mathbb{R}^{2}} \min \lbrace L_{f},\mathcal{K}^{f}_{\widetilde{V}}
(\xi,\tau) \rbrace \, \md \lambda_{f}(\xi) \, \md \lambda_{f}(\tau)$ as 
$j \! \to \! \infty$, whence
\begin{equation*}
\liminf_{j \to \infty} \mathfrak{d}^{\widetilde{V},f}_{\mathcal{N}_{
\hat{k}_{j}}} \! \geqslant \! \iint_{\mathbb{R}^{2}} \min \left\lbrace 
L_{f},\mathcal{K}_{\widetilde{V}}^{f}(\xi,\tau) \right\rbrace \md 
\lambda_{f}(\xi) \, \md \lambda_{f}(\tau):
\end{equation*}
letting $L_{f} \! \uparrow \! +\infty$ and using monotone convergence, 
one arrives at
\begin{equation*}
\liminf_{j \to \infty} \mathfrak{d}^{\widetilde{V},f}_{\mathcal{N}_{
\hat{k}_{j}}} \! \geqslant \! \iint_{\mathbb{R}^{2}} \mathcal{K}_{
\widetilde{V}}^{f}(\xi,\tau) \, \md \lambda_{f}(\xi) \, \md \lambda_{f}
(\tau) \! = \! \mathrm{I}_{\widetilde{V}}^{f}[\lambda_{f}] \quad (> \! -\infty);
\end{equation*}
but,
\begin{equation*}
E_{\widetilde{V}}^{f} \! \geqslant \! \limsup_{j \to \infty} \mathfrak{d}^{
\widetilde{V},f}_{\mathcal{N}_{\hat{k}_{j}}} \! \geqslant \! \liminf_{j 
\to \infty} \mathfrak{d}^{\widetilde{V},f}_{\mathcal{N}_{\hat{k}_{j}}} \! 
\geqslant \! \mathrm{I}_{\widetilde{V}}^{f}[\lambda_{f}] \! \geqslant 
\! E_{\widetilde{V}}^{f} \! = \! \mathrm{I}_{\widetilde{V}}^{f}
[\mu_{\widetilde{V}}^{f}] \quad \Rightarrow \quad \mathrm{I}_{
\widetilde{V}}^{f}[\lambda_{f}] \! = \! \mathrm{I}_{\widetilde{V}}^{f}
[\mu_{\widetilde{V}}^{f}],
\end{equation*}
which, via Lemma~\ref{lem3.3}, implies that, for $n \! \in \! \mathbb{N}$ and 
$k \! \in \! \lbrace 1,2,\dotsc,K \rbrace$ such that $\alpha_{p_{\mathfrak{s}}} 
\! := \! \alpha_{k} \! \neq \! \infty$, $\lambda_{f} \! = \! 
\mu_{\widetilde{V}}^{f}$; in particular, for $\mathfrak{g} \! 
\in \! \pmb{\operatorname{C}}^{0}_{\text{b}}(\mathbb{R})$, 
$\int_{\mathbb{R}} \mathfrak{g}(\xi) \, \md 
\lambda^{f}_{\mathcal{N}_{\hat{k}_{j}}}(\xi) \! \to \! \int_{\mathbb{R}} 
\mathfrak{g}(\xi) \, \md \mu_{\widetilde{V}}^{f}(\xi)$ as $j \! \to \! 
\infty$, which is a contradiction; hence, $\lambda_{\mathcal{N}}^{f} 
\! \overset{\ast}{\to} \! \mu_{\widetilde{V}}^{f}$ in the double-scaling 
limit $\mathscr{N},n \! \to \! \infty$ such that $z_{o} \! = \! 1 \! + 
\! o(1)$. One now argues as in the final paragraph of the proof of 
Lemma~3.5 in \cite{a45} to show that, for $k \! \in \! \lbrace 
1,2,\dotsc,K \rbrace$ such that $\alpha_{p_{\mathfrak{s}}} \! 
:= \! \alpha_{k} \! \neq \! \infty$, if $\lbrace \mathfrak{d}^{
\widetilde{V},f}_{\mathcal{N}} \rbrace_{n=1}^{\infty}$, with 
$\mathcal{N} \! := \! (n \! - \! 1)K \! + \! k$, converges for some 
subsequence $(\lbrace \mathfrak{d}^{\widetilde{V},f}_{\mathcal{N}_{l}} 
\rbrace_{l=1}^{\infty}$ in $\mathscr{M}_{1}(\mathbb{R})$, say, with 
$\mathcal{N}_{l} \! := \! (n_{l} \! - \! 1)K \! + \! k)$, then the limit is 
always equal to $E_{\widetilde{V}}^{f}$, that is, $\lim_{\underset{z_{o}
=1+o(1)}{\mathscr{N},n \to \infty}} \mathfrak{d}^{\widetilde{V},f}_{
\mathcal{N}} \! = \! E_{\widetilde{V}}^{f}$.\footnote{Since, for $n \! 
\in \! \mathbb{N}$ and $k \! \in \! \lbrace 1,2,\dotsc,K \rbrace$ such 
that $\alpha_{p_{\mathfrak{s}}} \! := \! \alpha_{k} \! \neq \! \infty$, 
$\limsup_{j \to \infty} \mathfrak{d}^{\widetilde{V},f}_{\mathcal{N}_{
\hat{k}_{j}}}$ and $\liminf_{j \to \infty} \mathfrak{d}^{\widetilde{V},
f}_{\mathcal{N}_{\hat{k}_{j}}}$ exist, it follows that $\lim_{j \to \infty} 
\mathfrak{d}^{\widetilde{V},f}_{\mathcal{N}_{\hat{k}_{j}}}$ exists and 
equals $E_{\widetilde{V}}^{f} \! = \! \mathrm{I}_{\widetilde{V}}^{f}
[\mu_{\widetilde{V}}^{f}]$; in fact, $\lim_{\underset{z_{o}=1+
o(1)}{\mathscr{N},n \to \infty}} \mathfrak{d}^{\widetilde{V},f}_{
\mathcal{N}} \! = \! E_{\widetilde{V}}^{f} \! = \! \mathrm{I}_{
\widetilde{V}}^{f}[\mu_{\widetilde{V}}^{f}]$.}

$\pmb{(2)}$ The proof of this case, that is, $n \! \in \! \mathbb{N}$ 
and $k \! \in \! \lbrace 1,2,\dotsc,K \rbrace$ such that 
$\alpha_{p_{\mathfrak{s}}} \! := \! \alpha_{k} \! = \! \infty$, is 
virtually identical to the proof presented in $\pmb{(1)}$ above; one 
mimics, \emph{verbatim}, the scheme of the calculations presented 
in case $\pmb{(1)}$ in order to arrive at the corresponding claims 
stated in the lemma; in order to do so, however, the analogues of 
Equations~\eqref{eql3.5a}--\eqref{eql3.5f} are necessary, which, 
in the present case, read: for $n \! \in \! \mathbb{N}$ and $k \! \in 
\! \lbrace 1,2,\dotsc,K \rbrace$ such that $\alpha_{p_{\mathfrak{s}}} 
\! := \! \alpha_{k} \! = \! \infty$, with $\mathcal{N} \! := \! (n \! - 
\! 1)K \! + \! k$,
\begin{equation} \label{eql3.5h}
\delta^{\widetilde{V},\infty}_{\mathcal{N}} \! = \sup_{\lbrace x_{1},x_{2},
\dotsc,x_{\mathcal{N}} \rbrace \subset \mathbb{R}} \left(\prod_{i<j}
(h(x_{i},x_{j}))^{2} \me^{-\widetilde{V}(x_{i})} \me^{-\widetilde{V}(x_{j})} 
\right)^{\frac{1}{\mathcal{N}(\mathcal{N}-1)}} = \sup_{\lbrace x_{1},
x_{2},\dotsc,x_{\mathcal{N}} \rbrace \subset \mathbb{R}} \left(
\prod_{i \neq j}h(x_{i},x_{j}) \me^{-(\mathcal{N}-1) \sum_{j=1}^{
\mathcal{N}} \widetilde{V}(x_{j})} \right)^{\frac{1}{\mathcal{N}
(\mathcal{N}-1)}},
\end{equation}
\begin{equation} \label{eql3.5i} 
h \colon \mathbb{N} \times \lbrace 1,2,\dotsc,K \rbrace \times 
\mathbb{R}^{2} \! \ni \! (n,k,x,y) \! \mapsto \! \lvert x \! - \! y 
\rvert^{\frac{\varkappa_{nk}}{n}} \prod_{q=1}^{\mathfrak{s}-1} \left(
\dfrac{\lvert x \! - \! y \rvert}{\lvert x \! - \! \alpha_{p_{q}} \rvert 
\lvert y \! - \! \alpha_{p_{q}} \rvert} \right)^{\frac{\varkappa_{nk 
\tilde{k}_{q}}}{n}} \! = \! h(n,k,x,y) \! := \! h(x,y),
\end{equation}
\begin{align}
\mathscr{K}^{\widetilde{V},\infty}_{\mathcal{N}}(x_{1},x_{2},\dotsc,
x_{\mathcal{N}}) \! := \! \sum_{\substack{i,j=1\\i \neq j}}^{
\mathcal{N}} \mathcal{K}_{\widetilde{V}}^{\infty}(x_{i},x_{j}) =& \, 
\sum_{\substack{i,j=1\\i \neq j}}^{\mathcal{N}} \left(\dfrac{
\varkappa_{nk}}{n} \ln \left(\dfrac{1}{\lvert x_{i} \! - \! x_{j} \rvert} 
\right) \! + \! \sum_{q=1}^{\mathfrak{s}-1} \dfrac{\varkappa_{nk 
\tilde{k}_{q}}}{n} \ln \left(\dfrac{\lvert x_{i} \! - \! \alpha_{p_{q}} 
\rvert \lvert x_{j} \! - \! \alpha_{p_{q}} \rvert}{\lvert x_{i} \! - \! 
x_{j} \rvert} \right) \right) \nonumber \\
+& \, (\mathcal{N} \! - \! 1) \sum_{j=1}^{\mathcal{N}} \widetilde{V}
(x_{j}), \label{eql3.5k}
\end{align}
\begin{gather} 
\mathfrak{d}^{\widetilde{V},\infty}_{\mathcal{N}} \! := \! (\mathcal{N}
(\mathcal{N} \! - \! 1))^{-1} \inf_{\lbrace x_{1},x_{2},\dotsc,
x_{\mathcal{N}} \rbrace \subset \mathbb{R}} \mathscr{K}^{
\widetilde{V},\infty}_{\mathcal{N}}(x_{1},x_{2},\dotsc,
x_{\mathcal{N}}), \label{eql3.5l} \\
\lambda_{\mathcal{N}}^{\infty} \! := \! \dfrac{1}{\mathcal{N}} 
\sum_{j=1}^{\mathcal{N}} \delta_{\hat{x}_{j}^{\ast}}, \label{eql3.5m}
\end{gather}
where $\lbrace \hat{x}_{1}^{\ast},\hat{x}_{2}^{\ast},\dotsc,
\hat{x}_{\mathcal{N}}^{\ast} \rbrace$ $(\subset \overline{\mathbb{R}} 
\setminus \lbrace \alpha_{1},\alpha_{2},\dotsc,\alpha_{K} \rbrace)$, with 
$\hat{x}_{j}^{\ast} \! = \! \hat{x}_{j}^{\ast}(n,k,z_{o})$, $j \! = \! 1,2,\dotsc,
\mathcal{N}$, is an associated, weighted $\mathcal{N}$-Fekete set, that is, 
$\mathfrak{d}^{\widetilde{V},\infty}_{\mathcal{N}} \! := \! (\mathcal{N}
(\mathcal{N} \! - \! 1))^{-1} \mathscr{K}^{\widetilde{V},\infty}_{\mathcal{N}}
(\hat{x}_{1}^{\ast},\hat{x}_{2}^{\ast},\dotsc,\hat{x}_{\mathcal{N}}^{\ast})$. 
This concludes the proof. \hfill $\qed$

Via Lemma~\ref{lem3.5}, one now establishes, for $n \! \in \! 
\mathbb{N}$ and $k \! \in \! \lbrace 1,2,\dotsc,K \rbrace$ such that 
$\alpha_{p_{\mathfrak{s}}} \! := \! \alpha_{k} \! = \! \infty$ (resp., 
$\alpha_{p_{\mathfrak{s}}} \! := \!\alpha_{k} \! \neq \! \infty)$, the 
regularity of the associated family of $\hat{K} \! := \! \# \lbrace 
\mathstrut j \! \in \! \lbrace 1,2,\dotsc,K \rbrace; \, \alpha_{j} \! 
= \! \infty \rbrace$ (resp., $\tilde{K} \! := \! \# \lbrace \mathstrut 
j \! \in \! \lbrace 1,2,\dotsc,K \rbrace; \, \alpha_{j} \! \neq \! \infty 
\rbrace)$ equilibrium measures $\mu_{\widetilde{V}}^{\infty}$ (resp., 
$\mu_{\widetilde{V}}^{f})$. 
\begin{ccccc} \label{lem3.6} 
Let $\widetilde{V} \colon \overline{\mathbb{R}} \setminus \lbrace 
\alpha_{1},\alpha_{2},\dotsc,\alpha_{K} \rbrace \! \to \! \mathbb{R}$ 
satisfy conditions~\eqref{eq20}--\eqref{eq22}. Then, for $n \! \in \! 
\mathbb{N}$ and $k \! \in \! \lbrace 1,2,\dotsc,K \rbrace$ such that 
$\alpha_{p_{\mathfrak{s}}} \! := \! \alpha_{k} \! = \! \infty$ (resp., 
$\alpha_{p_{\mathfrak{s}}} \! := \! \alpha_{k} \! \neq \! \infty)$, the 
associated equilibrium measure $\mu_{\widetilde{V}}^{\infty}$ (resp., 
$\mu_{\widetilde{V}}^{f})$ $\in \! \mathscr{M}_{1}(\mathbb{R})$ is 
absolutely continuous with respect to Lebesgue measure.
\end{ccccc}

\emph{Proof}. The proof of this Lemma~\ref{lem3.6} consists of two 
cases: (i) $n \! \in \! \mathbb{N}$ and $k \! \in \! \lbrace 1,2,\dotsc,
K \rbrace$ such that $\alpha_{p_{\mathfrak{s}}} \! := \! \alpha_{k} \! 
= \! \infty$; and (ii) $n \! \in \! \mathbb{N}$ and $k \! \in \! \lbrace 
1,2,\dotsc,K \rbrace$ such that $\alpha_{p_{\mathfrak{s}}} \! 
:= \! \alpha_{k} \! \neq \! \infty$. The proof for the case 
$\alpha_{p_{\mathfrak{s}}} \! := \! \alpha_{k} \! \neq \! \infty$, $k \! 
\in \! \lbrace 1,2,\dotsc,K \rbrace$, will be considered in detail (see 
$\pmb{(1)}$ below), whilst the case $\alpha_{p_{\mathfrak{s}}} \! := 
\! \alpha_{k} \! = \! \infty$, $k \! \in \! \lbrace 1,2,\dotsc,K \rbrace$, 
can be proved, modulo technical and computational particulars, 
analogously (see $\pmb{(2)}$ below).

$\pmb{(1)}$ Let $\widetilde{V} \colon \overline{\mathbb{R}} \setminus 
\lbrace \alpha_{1},\alpha_{2},\dotsc,\alpha_{K} \rbrace \! \to \! 
\mathbb{R}$ satisfy conditions~\eqref{eq20}--\eqref{eq22}. For $n \! 
\in \! \mathbb{N}$ and $k \! \in \! \lbrace 1,2,\dotsc,K \rbrace$ such 
that $\alpha_{p_{\mathfrak{s}}} \! := \! \alpha_{k} \! \neq \! \infty$, 
with $\mathcal{N} \! := \! (n \! - \! 1)K \! + \! k$, let
\begin{equation} \label{eql3.6a} 
\tilde{f} \colon \mathbb{N} \times \lbrace 1,2,\dotsc,K \rbrace \times 
\overline{\mathbb{R}} \setminus \lbrace \alpha_{1},\alpha_{2},\dotsc,
\alpha_{K} \rbrace \! \ni \! (n,k,x) \! \mapsto \! \dfrac{\prod_{j=1}^{
\mathcal{N}}(x \! - \! \tilde{x}_{j}^{\ast})}{\prod_{\underset{q \neq 
\mathfrak{s}-1}{q=1}}^{\mathfrak{s}}(x \! - \! \alpha_{p_{q}})^{l_{q}}} 
\! =: \! \tilde{f}(x),
\end{equation}
where $\lbrace \tilde{x}_{1}^{\ast},\tilde{x}_{2}^{\ast},\dotsc,\tilde{x}_{
\mathcal{N}}^{\ast} \rbrace \! \subset \! \overline{\mathbb{R}} \setminus 
\lbrace \alpha_{1},\alpha_{2},\dotsc,\alpha_{K} \rbrace$, with the 
enumeration $\tilde{x}_{i}^{\ast} \! < \! \tilde{x}_{j}^{\ast}$ for 
$i \! < \! j$, is the associated, weighted Fekete set described in 
the corresponding item of Lemma~\ref{lem3.5}. For $n \! \in \! 
\mathbb{N}$ and $k \! \in \! \lbrace 1,2,\dotsc,K \rbrace$ such that 
$\alpha_{p_{\mathfrak{s}}} \! := \! \alpha_{k} \! \neq \! \infty$, one 
shows, via the identity
\begin{align*}
\left(\sum_{j=1}^{\mathcal{N}} \dfrac{1}{x \! - \! \tilde{x}_{j}^{\ast}} 
\right)^{2} \! - \! \sum_{j=1}^{\mathcal{N}} \dfrac{1}{(x \! - \! 
\tilde{x}_{j}^{\ast})^{2}}=& \, \sum_{j=1}^{\mathcal{N}} 
\sum_{\substack{k^{\prime}=1\\k^{\prime} \neq j}}^{\mathcal{N}} 
\left(\dfrac{1}{x \! - \! \tilde{x}_{j}^{\ast}} \! - \! \dfrac{1}{x \! - \! 
\tilde{x}_{k^{\prime}}^{\ast}} \right) \dfrac{1}{\tilde{x}_{j}^{\ast} 
\! - \! \tilde{x}_{k^{\prime}}^{\ast}} \! = \! \sum_{j=1}^{\mathcal{N}} 
\sum_{\substack{k^{\prime}=1\\k^{\prime} \neq j}}^{\mathcal{N}} 
\dfrac{1}{(x \! - \! \tilde{x}_{j}^{\ast})(x \! - \! 
\tilde{x}_{k^{\prime}}^{\ast})} \\
=& \, 2 \sum_{j=1}^{\mathcal{N}} \sum_{\substack{k^{\prime}=
1\\k^{\prime} \neq j}}^{\mathcal{N}} \dfrac{1}{(x \! - \! 
\tilde{x}_{j}^{\ast})(\tilde{x}_{j}^{\ast} \! - \! 
\tilde{x}_{k^{\prime}}^{\ast})},
\end{align*}
that
\begin{equation} \label{eql3.6b} 
\dfrac{\tilde{f}^{\prime \prime}(x)}{\tilde{f}^{\prime}(x)} 
\! = \! \dfrac{\tilde{f}_{N}(x)}{\tilde{f}_{D}(x)},
\end{equation}
where the prime denotes differentiation with respect to $x$,
\begin{align}
\tilde{f}_{N}(x) :=& \, 2 \sum_{j=1}^{\mathcal{N}} 
\prod_{\substack{l=1\\l \neq j}}^{\mathcal{N}}(x \! - \! \tilde{x}_{l}^{\ast}) 
\sum_{\substack{k^{\prime}=1\\k^{\prime} \neq j}}^{\mathcal{N}} 
\dfrac{1}{\tilde{x}_{j}^{\ast} \! - \! \tilde{x}_{k^{\prime}}^{\ast}} \! - \! 2 
\sum_{j=1}^{\mathcal{N}} \prod_{\substack{l=1\\l \neq j}}^{\mathcal{N}}
(x \! - \! \tilde{x}_{l}^{\ast}) \sum_{\substack{q=1\\q \neq \mathfrak{s}
-1}}^{\mathfrak{s}} \dfrac{l_{q}}{x \! - \! \alpha_{p_{q}}} \nonumber \\
+& \, \prod_{j=1}^{\mathcal{N}}(x \! - \! \tilde{x}_{j}^{\ast}) \left(\sum_{
\substack{q=1\\q \neq \mathfrak{s}-1}}^{\mathfrak{s}} \dfrac{l_{q}}{(x 
\! - \! \alpha_{p_{q}})^{2}} \! + \! \left(\sum_{\substack{q=1\\q \neq 
\mathfrak{s}-1}}^{\mathfrak{s}} \dfrac{l_{q}}{x \! - \! \alpha_{p_{q}}} 
\right)^{2} \right), \label{eql3.6c} \\
\tilde{f}_{D}(x) :=& \, \sum_{j=1}^{\mathcal{N}} \sum_{\substack{k^{\prime}
=1\\k^{\prime} \neq j}}^{\mathcal{N}}(x \! - \! \tilde{x}_{k^{\prime}}^{\ast}) 
\! - \! \prod_{j=1}^{\mathcal{N}}(x \! - \! \tilde{x}_{j}^{\ast}) 
\sum_{\substack{q=1\\q \neq \mathfrak{s}-1}}^{\mathfrak{s}} 
\dfrac{l_{q}}{x \! - \! \alpha_{p_{q}}}; \label{eql3.6d} 
\end{align}
hence, noting that, for $m \! = \! 1,2,\dotsc,\mathcal{N}$,
\begin{gather*}
\left. \sum_{j=1}^{\mathcal{N}} \prod_{\substack{l=1\\l \neq j}}^{
\mathcal{N}}(x \! - \! \tilde{x}_{l}^{\ast}) \sum_{\substack{k^{\prime}
=1\\k^{\prime} \neq j}}^{\mathcal{N}} \dfrac{1}{\tilde{x}_{j}^{\ast} \! - 
\! \tilde{x}_{k^{\prime}}^{\ast}} \right\vert_{x=\tilde{x}_{m}^{\ast}} = 
\prod_{\substack{l=1\\l \neq m}}^{\mathcal{N}}(\tilde{x}_{m}^{\ast} 
\! - \! \tilde{x}_{l}^{\ast}) \sum_{\substack{k^{\prime}=1\\k^{\prime} 
\neq m}}^{\mathcal{N}} \dfrac{1}{\tilde{x}_{m}^{\ast} \! - \! 
\tilde{x}_{k^{\prime}}^{\ast}}, \\
\left. \sum_{j=1}^{\mathcal{N}} \prod_{\substack{l=1\\l \neq j}}^{
\mathcal{N}}(x \! - \! \tilde{x}_{l}^{\ast}) \sum_{\substack{q=1\\
q \neq \mathfrak{s}-1}}^{\mathfrak{s}} \dfrac{l_{q}}{x \! - \! \alpha_{p_{q}}} 
\right\vert_{x=\tilde{x}_{m}^{\ast}} = \prod_{\substack{l=1\\l \neq 
m}}^{\mathcal{N}}(\tilde{x}_{m}^{\ast} \! - \! \tilde{x}_{l}^{\ast}) 
\sum_{\substack{q=1\\q \neq \mathfrak{s}-1}}^{\mathfrak{s}} 
\dfrac{l_{q}}{\tilde{x}_{m}^{\ast} \! - \! \alpha_{p_{q}}}, \\
\left. \sum_{j=1}^{\mathcal{N}} \prod_{\substack{k^{\prime}=
1\\k^{\prime} \neq j}}^{\mathcal{N}}(x \! - \! \tilde{x}_{k^{\prime}}^{\ast}) 
\! - \! \prod_{j=1}^{\mathcal{N}}(x \! - \! \tilde{x}_{j}^{\ast}) 
\sum_{\substack{q=1\\q \neq \mathfrak{s}-1}}^{\mathfrak{s}} 
\dfrac{l_{q}}{x \! - \! \alpha_{p_{q}}} \right\vert_{x=\tilde{x}_{m}^{\ast}}
=\prod_{\substack{k^{\prime}=1\\k^{\prime} \neq m}}^{\mathcal{N}}
(\tilde{x}_{m}^{\ast} \! - \! \tilde{x}_{k^{\prime}}^{\ast}),
\end{gather*}
one arrives at, upon recalling that $\lbrace \tilde{x}_{1}^{\ast},
\tilde{x}_{2}^{\ast},\dotsc,\tilde{x}_{\mathcal{N}}^{\ast} \rbrace 
\cap \lbrace \alpha_{1},\alpha_{2},\dotsc,\alpha_{K} \rbrace \! 
= \! \varnothing$,
\begin{equation} \label{eql3.6e} 
\dfrac{1}{2} \dfrac{\tilde{f}^{\prime \prime}(\tilde{x}_{m}^{\ast})}{
\tilde{f}^{\prime}(\tilde{x}_{m}^{\ast})} \! = \! \sum_{\substack{
k^{\prime}=1\\k^{\prime} \neq m}}^{\mathcal{N}} \dfrac{1}{
\tilde{x}_{m}^{\ast} \! - \! \tilde{x}_{k^{\prime}}^{\ast}} \! - \! 
\sum_{\substack{q=1\\q \neq \mathfrak{s}-1}}^{\mathfrak{s}} 
\dfrac{l_{q}}{\tilde{x}_{m}^{\ast} \! - \! \alpha_{p_{q}}}, \quad m 
\! = \! 1,2,\dotsc,\mathcal{N}.
\end{equation}
{}From Equations~\eqref{eql3.5a} and~\eqref{eql3.5b}, one shows that, 
for $n \! \in \! \mathbb{N}$ and $k \! \in \! \lbrace 1,2,\dotsc,K 
\rbrace$ such that $\alpha_{p_{\mathfrak{s}}} \! := \! \alpha_{k} \! 
\neq \! \infty$,
\begin{equation} \label{eql3.6f} 
-\dfrac{1}{2}(\mathcal{N} \! - \! 1) \widetilde{V}^{\prime}
(\tilde{x}_{m}^{\ast}) \! + \! \dfrac{\mathcal{N}}{n} 
\sum_{\substack{k^{\prime}=1\\k^{\prime} \neq m}}^{\mathcal{N}} 
\dfrac{1}{\tilde{x}_{m}^{\ast} \! - \! \tilde{x}_{k^{\prime}}^{\ast}} 
\! - \! (\mathcal{N} \! - \! 1) \left(\sum_{q=1}^{\mathfrak{s}-2} 
\dfrac{l_{q}/n}{\tilde{x}_{m}^{\ast} \! - \! \alpha_{p_{q}}} \! + \! 
\dfrac{(l_{\mathfrak{s}} \! - \! 1)/n}{\tilde{x}_{m}^{\ast} \! - \! 
\alpha_{p_{\mathfrak{s}}}} \right) \! = \! 0, \quad m \! = \! 
1,2,\dotsc,
\mathcal{N}.
\end{equation}
Combining Equations~\eqref{eql3.6e} and~\eqref{eql3.6f}, and noting 
that $\lbrace \mathstrut x \! \in \! \mathbb{R}; \, \tilde{f}(x) \! = \! 0 
\rbrace \! = \! \lbrace \tilde{x}_{1}^{\ast},\tilde{x}_{2}^{\ast},\dotsc,
\tilde{x}_{\mathcal{N}}^{\ast} \rbrace$ and $\lbrace \tilde{x}_{1}^{\ast},
\tilde{x}_{2}^{\ast},\dotsc,\tilde{x}_{\mathcal{N}}^{\ast} \rbrace \cap 
\lbrace \alpha_{1},\alpha_{2},\dotsc,\alpha_{K} \rbrace \! = \! \varnothing$, 
one proceeds as in the proof of Lemma~3.6 in \cite{a45} (see, also, 
\cite{rnveziia,mfr}) to show that, for $n \! \in \! \mathbb{N}$ and $k \! \in 
\! \lbrace 1,2,\dotsc,K \rbrace$ such that $\alpha_{p_{\mathfrak{s}}} \! := 
\! \alpha_{k} \! \neq \! \infty$,
\begin{equation} \label{eql3.6g} 
\tilde{f}^{\prime \prime}(x) \! + \! 2 \dfrac{n}{\mathcal{N}} \left(
-\dfrac{1}{2}(\mathcal{N} \! - \! 1) \widetilde{V}^{\prime}(x) \! 
+ \! \sum_{\substack{q=1\\q \neq \mathfrak{s}-1}}^{\mathfrak{s}} 
\dfrac{l_{q}/n}{x \! - \! \alpha_{p_{q}}} \! + \! \dfrac{(\mathcal{N} 
\! - \! 1)/n}{x \! - \! \alpha_{p_{\mathfrak{s}}}} \right) 
\tilde{f}^{\prime}(x) \! = \! \tilde{Q}(x) \tilde{f}(x),
\end{equation}
where (with $\tilde{Q}(\tilde{x}_{m}^{\ast}) \! \neq \! 0$, 
$m \! = \! 1,2,\dotsc,\mathcal{N})$
\begin{align}
\tilde{Q}(x)=& \, 2 \dfrac{(\mathcal{N} \! - \! 1)}{\mathcal{N}} 
\sum_{j=1}^{\mathcal{N}} \sum_{q=1}^{\mathfrak{s}-2} 
\dfrac{l_{q}}{(\tilde{x}_{j}^{\ast} \! - \! \alpha_{p_{q}})(x \! - \! 
\alpha_{p_{q}})} \! + \! 2 \dfrac{(\mathcal{N} \! - \! 1)}{\mathcal{N}} 
\sum_{j=1}^{\mathcal{N}} \dfrac{(l_{\mathfrak{s}} \! - \! 1)}{
(\tilde{x}_{j}^{\ast} \! - \! \alpha_{p_{\mathfrak{s}}})(x \! - \! 
\alpha_{p_{\mathfrak{s}}})} \! + \! \sum_{\substack{q=1\\q 
\neq \mathfrak{s}-1}}^{\mathfrak{s}} \dfrac{l_{q}}{(x \! - \! 
\alpha_{p_{q}})^{2}} \nonumber \\
+& \, \dfrac{(\mathcal{N} \! - \! 2)}{\mathcal{N}} \left(\sum_{
\substack{q=1\\q \neq \mathfrak{s}-1}}^{\mathfrak{s}} 
\dfrac{l_{q}}{x \! - \! \alpha_{p_{q}}} \right)^{2} \! - \! 2 
\dfrac{(\mathcal{N} \! - \! 1)}{\mathcal{N}} \dfrac{1}{x \! - \! 
\alpha_{p_{\mathfrak{s}}}} \sum_{\substack{q=1\\q \neq 
\mathfrak{s}-1}}^{\mathfrak{s}} \dfrac{l_{q}}{x \! - \! \alpha_{p_{q}}} 
\! + \! \dfrac{(\mathcal{N} \! - \! 1)}{\mathcal{N}}n \widetilde{V}^{
\prime}(x) \sum_{\substack{q=1\\q \neq \mathfrak{s}-1}}^{
\mathfrak{s}} \dfrac{l_{q}}{x \! - \! \alpha_{p_{q}}} \nonumber \\
-& \, \dfrac{n}{\mathcal{N}}(\mathcal{N} \! - \! 1) \sum_{j=1}^{\mathcal{N}} 
\dfrac{\widetilde{V}^{\prime}(x) \! - \! \widetilde{V}^{\prime}
(\tilde{x}_{j}^{\ast})}{x \! - \! \tilde{x}_{j}^{\ast}}. \label{eql3.6h} 
\end{align}
Since, for $n \! \in \! \mathbb{N}$ and $k \! \in \! \lbrace 1,2,
\dotsc,K \rbrace$ such that $\alpha_{p_{\mathfrak{s}}} \! := 
\! \alpha_{k} \! \neq \! \infty$, $\lbrace \tilde{x}_{1}^{\ast},
\tilde{x}_{2}^{\ast},\dotsc,\tilde{x}_{\mathcal{N}}^{\ast} \rbrace 
\subset \overline{\mathbb{R}} \setminus \lbrace \alpha_{1},
\alpha_{2},\dotsc,\alpha_{K} \rbrace$, choose, for $\tilde{N}$ $(= \! 
\tilde{N}(n,k,z_{o}))$ $\in \! \mathbb{N}_{0}$, $2(\tilde{N} \! + \! 1)$ real 
points $\lbrace \tilde{B}_{j-1},\tilde{A}_{j} \rbrace_{j=1}^{\tilde{N}+1}$ 
so that $\overline{\mathbb{R}} \setminus \lbrace \alpha_{1},\alpha_{2},
\dotsc,\alpha_{K} \rbrace \supset \cup_{j=1}^{\tilde{N}+1}[\tilde{B}_{j-1},
\tilde{A}_{j}] \supseteq \lbrace \tilde{x}_{1}^{\ast},\tilde{x}_{2}^{\ast},
\dotsc,\tilde{x}_{\mathcal{N}}^{\ast} \rbrace$, with enumeration 
$-\infty \! < \! \tilde{B}_{0} \! < \! \tilde{A}_{1} \! < \! \dotsb \! < \! 
\tilde{B}_{j-1} \! < \! \tilde{A}_{j} \! < \! \dotsb \! < \! \tilde{B}_{\tilde{N}} 
\! < \! \tilde{A}_{\tilde{N}+1} \! < \! +\infty$, and $[\tilde{B}_{j-1},
\tilde{A}_{j}] \cap \lbrace \alpha_{1},\alpha_{2},\dotsc,\alpha_{K} 
\rbrace \! = \! \varnothing$, $j \! = \! 1,2,\dotsc,\tilde{N} \! + \! 
1$.\footnote{$\cup_{j=1}^{\tilde{N}+1}[\tilde{B}_{j-1},\tilde{A}_{j}]$, 
which is the disjoint union of $\tilde{N} \! + \! 1$ compact real intervals, 
is the `pre-confinement domain' for the associated, weighted Fekete 
points $\tilde{x}_{m}^{\ast}$, $m \! = \! 1,2,\dotsc,\mathcal{N}$.} Let 
$\tilde{Z}(x) \! := \! \prod_{j=1}^{\tilde{N}+1}(z \! - \! \tilde{B}_{j-1})
(z \! - \! \tilde{A}_{j})$; then, a straightforward calculation shows that 
$\tilde{Z}(x)$ can be presented as $\tilde{Z}(x) \! = \! \sum_{j=2(\tilde{N}
+1)}^{0} \tilde{c}_{2(\tilde{N}+1)-j}x^{j}$,\footnote{For the purposes of 
this proof, the adopted convention is that $\sum_{j=2(\tilde{N}+1)}^{0} 
\tilde{\blacklozenge}(j) \! := \! \tilde{\blacklozenge}(2(\tilde{N} \! + \! 1)) 
\! + \! \dotsb \! + \! \tilde{\blacklozenge}(1) \! + \! \tilde{\blacklozenge}
(0)$.} where {}\footnote{The explicit $n$-, $k$-, 
and $z_{o}$-dependencies of $\tilde{B}_{j-1}$, $\tilde{A}_{j}$, 
$j \! = \! 1,2,\dotsc,\tilde{N} \! + \! 1$, and $\tilde{c}_{i}$, 
$i \! = \! 0,1,\dotsc,2(\tilde{N} \! + \! 1)$, are not essential here.}
\begin{equation*}
\tilde{c}_{j} \! = \! \mathlarger{\sum_{\underset{\underset{\sum_{i=1}^{
\tilde{N}+1}k^{\prime}_{i}+\sum_{m=1}^{\tilde{N}+1}l^{\prime}_{m}=j}{i,m 
\in \lbrace 1,2,\dotsc,\tilde{N}+1 \rbrace}}{k^{\prime}_{i},l^{\prime}_{m}=
0,1}}}(-1)^{j} \left(\prod_{i=1}^{\tilde{N}+1} \binom{1}{k^{\prime}_{i}}
(\tilde{B}_{i-1})^{k^{\prime}_{i}} \right) \left(\prod_{m=1}^{\tilde{N}+1} 
\binom{1}{l^{\prime}_{m}}(\tilde{A}_{m})^{l^{\prime}_{m}} \right), \quad 
j \! = \! 0,1,\dotsc,2(\tilde{N} \! + \! 1).
\end{equation*}
(Note: $\tilde{c}_{0} \! = \! 1$ and $\tilde{c}_{2(\tilde{N}+1)} \! = 
\! \prod_{j=1}^{\tilde{N}+1} \tilde{B}_{j-1} \tilde{A}_{j}$.) Via the 
above definition of $\tilde{Z}(x)$, one shows that, for $n \! \in \! 
\mathbb{N}$ and $k \! \in \! \lbrace 1,2,\dotsc,K \rbrace$ such 
that $\alpha_{p_{\mathfrak{s}}} \! := \! \alpha_{k} \! \neq \! \infty$, 
the term $\tfrac{n}{\mathcal{N}}(\mathcal{N} \! - \! 1) \sum_{j=1}^{
\mathcal{N}} \tfrac{\widetilde{V}^{\prime}(x)-\widetilde{V}^{\prime}
(\tilde{x}_{j}^{\ast})}{x-\tilde{x}_{j}^{\ast}}$, which appears in the 
expression for $\tilde{Q}(x)$ given in Equation~\eqref{eql3.6h}, 
can be presented as
\begin{align}
\dfrac{n}{\mathcal{N}}(\mathcal{N} \! - \! 1) \sum_{j=1}^{\mathcal{N}} 
\dfrac{\widetilde{V}^{\prime}(x) \! - \! \widetilde{V}^{\prime}
(\tilde{x}_{j}^{\ast})}{x \! - \! \tilde{x}_{j}^{\ast}} &= \, 
\dfrac{n}{\mathcal{N}} \dfrac{(\mathcal{N} \! - \! 1)}{\tilde{Z}(x)} 
\sum_{j=1}^{\mathcal{N}} \left(\dfrac{\widetilde{V}^{\prime}(x) \tilde{Z}
(x) \! - \! \widetilde{V}^{\prime}(\tilde{x}_{j}^{\ast}) \tilde{Z}
(\tilde{x}_{j}^{\ast})}{x \! - \! \tilde{x}_{j}^{\ast}} \right) \! - \! 
\dfrac{2}{\tilde{Z}(x)} \sum_{j=1}^{\mathcal{N}} \sum_{
\substack{k^{\prime}=1\\k^{\prime} \neq j}}^{\mathcal{N}} 
\dfrac{1}{\tilde{x}_{j}^{\ast} \! - \! \tilde{x}_{k^{\prime}}^{\ast}} \nonumber \\
\times& \, \left(\sum_{i=2(\tilde{N}+1)}^{3} \tilde{c}_{2(\tilde{N}+1)-i} 
\left(\sum_{r=0}^{i-1}x^{i-r-1}(\tilde{x}_{j}^{\ast})^{r} \right) \! 
+ \! \tilde{c}_{2 \tilde{N}}(x \! + \! \tilde{x}_{j}^{\ast}) \! + \! 
\tilde{c}_{2 \tilde{N}+1} \right) \nonumber \\
+& \, \dfrac{2}{\tilde{Z}(x)} \dfrac{(\mathcal{N} \! - \! 1)}{\mathcal{N}} 
\sum_{j=1}^{\mathcal{N}} \sum_{q=1}^{\mathfrak{s}-2} \dfrac{l_{q}}{
\tilde{x}_{j}^{\ast} \! - \! \alpha_{p_{q}}} \left(\sum_{i=2(\tilde{N}+1)}^{3} 
\tilde{c}_{2(\tilde{N}+1)-i} \left(\sum_{r=0}^{i-1}x^{i-r-1}
(\tilde{x}_{j}^{\ast})^{r} \right) \right. \nonumber \\
+&\left. \, \tilde{c}_{2 \tilde{N}}(x \! + \! \tilde{x}_{j}^{\ast}) \! 
+ \! \tilde{c}_{2 \tilde{N}+1} \vphantom{M^{M^{M^{M^{M}}}}} \right) \! + 
\! \dfrac{2}{\tilde{Z}(x)} \dfrac{(\mathcal{N} \! - \! 1)}{\mathcal{N}} 
\sum_{j=1}^{\mathcal{N}} \dfrac{(l_{\mathfrak{s}} \! - \! 1)}{
\tilde{x}_{j}^{\ast} \! - \! \alpha_{p_{\mathfrak{s}}}} \nonumber \\
\times& \, \left(\sum_{i=2(\tilde{N}+1)}^{3} \tilde{c}_{2(\tilde{N}+1)-i} 
\left(\sum_{r=0}^{i-1}x^{i-r-1}(\tilde{x}_{j}^{\ast})^{r} \right) \! 
+ \! \tilde{c}_{2 \tilde{N}}(x \! + \! \tilde{x}_{j}^{\ast}) \! + \! 
\tilde{c}_{2 \tilde{N}+1} \right): \label{eql3.6i}
\end{align}
noting that (cf. the proof of Lemma~\ref{lem3.5}), for an arbitrary 
function, $\tilde{h}$, say,
\begin{align*}
\sum_{j<k^{\prime}} \tilde{h}(\tilde{x}_{j}^{\ast}) \! := \! 
\sum_{j=1}^{\mathcal{N}-1} \sum_{k^{\prime}=j+1}^{\mathcal{N}} 
\tilde{h}(\tilde{x}_{j}^{\ast})=& \, (\mathcal{N} \! - \! 1) \tilde{h}
(\tilde{x}_{1}^{\ast}) \! + \! \sum_{m=2}^{\mathcal{N}-1}(\mathcal{N} 
\! - \! m) \tilde{h}(\tilde{x}_{m}^{\ast}), \\
\sum_{j<k^{\prime}} \tilde{h}(\tilde{x}_{k^{\prime}}^{\ast}) \! := \! 
\sum_{j=1}^{\mathcal{N}-1} \sum_{k^{\prime}=j+1}^{\mathcal{N}} \tilde{h}
(\tilde{x}_{k^{\prime}}^{\ast})=& \, (\mathcal{N} \! - \! 1) \tilde{h}
(\tilde{x}_{\mathcal{N}}^{\ast}) \! + \! \sum_{m=2}^{\mathcal{N}-1}
(m \! - \! 1) \tilde{h}(\tilde{x}_{m}^{\ast}), \\
\sum_{j<k^{\prime}} \tilde{h}(\tilde{x}_{j}^{\ast}) \! + \! 
\sum_{j<k^{\prime}} \tilde{h}(\tilde{x}_{k^{\prime}}^{\ast})=& \, 
(\mathcal{N} \! - \! 1) \sum_{m=1}^{\mathcal{N}} \tilde{h}
(\tilde{x}_{m}^{\ast}),
\end{align*}
one shows, via the identities $y_{1}^{n} \! - \! y_{2}^{n} \! = \! (y_{1} \! 
- \! y_{2})(y_{1}^{n-1} \! + \! y_{1}^{n-2}y_{2} \! + \! \dotsb \! + \! y_{1}
y_{2}^{n-2} \! + \! y_{2}^{n-1})$,
\begin{align*}
\sum_{j=1}^{\mathcal{N}} \sum_{\substack{k^{\prime}=1\\k^{\prime} 
\neq j}}^{\mathcal{N}} \dfrac{1}{\tilde{x}_{j}^{\ast} \! - \! 
\tilde{x}_{k^{\prime}}^{\ast}}=& \, 0, \\
\sum_{j=1}^{\mathcal{N}} \sum_{\substack{k^{\prime}=1\\k^{\prime} 
\neq j}}^{\mathcal{N}} \dfrac{x \! + \! \tilde{x}_{j}^{\ast}}{\tilde{x}_{j}^{
\ast} \! - \! \tilde{x}_{k^{\prime}}^{\ast}}=& \, \sum_{j=1}^{\mathcal{N}} 
\sum_{\substack{k^{\prime}=1\\k^{\prime} \neq j}}^{\mathcal{N}} 
\dfrac{\tilde{x}_{j}^{\ast}}{\tilde{x}_{j}^{\ast} \! - \! \tilde{x}_{k^{
\prime}}^{\ast}} \! = \! \sum_{j=1}^{\mathcal{N}-1} \sum_{k^{\prime}=
j+1}^{\mathcal{N}}1 \! = \! \mathcal{N}(\mathcal{N} \! - \! 1)/2, \\
\sum_{j=1}^{\mathcal{N}} \sum_{\substack{k^{\prime}=1\\k^{\prime} 
\neq j}}^{\mathcal{N}} \dfrac{(\tilde{x}_{j}^{\ast})^{r}}{\tilde{x}_{j}^{\ast} 
\! - \! \tilde{x}_{k^{\prime}}^{\ast}}=& \, \sum_{j=1}^{\mathcal{N}-1} 
\sum_{k^{\prime}=j+1}^{\mathcal{N}} \dfrac{(\tilde{x}_{j}^{\ast})^{r} 
\! - \! (\tilde{x}_{k^{\prime}}^{\ast})^{r}}{\tilde{x}_{j}^{\ast} \! - \! 
\tilde{x}_{k^{\prime}}^{\ast}} \! = \! \sum_{j=1}^{\mathcal{N}-1} 
\sum_{k^{\prime}=j+1}^{\mathcal{N}} \sum_{m=0}^{r-1}
(\tilde{x}_{j}^{\ast})^{r-m-1}(\tilde{x}_{k^{\prime}}^{\ast})^{m} \\
=& \, \sum_{j=1}^{\mathcal{N}-1} \sum_{k^{\prime}=j+1}^{\mathcal{N}}
(\tilde{x}_{j}^{\ast})^{r-1} \! + \! \sum_{j=1}^{\mathcal{N}-1} 
\sum_{k^{\prime}=j+1}^{\mathcal{N}}(\tilde{x}_{k^{\prime}}^{\ast})^{r-1} 
\! + \! \sum_{j=1}^{\mathcal{N}-1} \sum_{k^{\prime}=j+1}^{\mathcal{N}} 
\sum_{m=2}^{r-1}(\tilde{x}_{j}^{\ast})^{r-m}
(\tilde{x}_{k^{\prime}}^{\ast})^{m-1} \\
=& \, (\mathcal{N} \! - \! 1) \sum_{j=1}^{\mathcal{N}}
(\tilde{x}_{j}^{\ast})^{r-1} \! + \! \sum_{j=1}^{\mathcal{N}-1} 
\sum_{k^{\prime}=j+1}^{\mathcal{N}} \sum_{m=2}^{r-1}
(\tilde{x}_{j}^{\ast})^{r-m}(\tilde{x}_{k^{\prime}}^{\ast})^{m-1}, \quad 
r \! \geqslant \! 2,
\end{align*}
that Equation~\eqref{eql3.6i} can be re-written as
\begin{align}
\dfrac{n}{\mathcal{N}}(\mathcal{N} \! - \! 1) \sum_{j=1}^{\mathcal{N}} 
\dfrac{\widetilde{V}^{\prime}(x) \! - \! \widetilde{V}^{\prime}
(\tilde{x}_{j}^{\ast})}{x \! - \! \tilde{x}_{j}^{\ast}}=& \, 
\dfrac{n}{\mathcal{N}} \dfrac{(\mathcal{N} \! - \! 1)}{\tilde{Z}(x)} 
\sum_{j=1}^{\mathcal{N}} \left(\dfrac{\widetilde{V}^{\prime}(x) \tilde{Z}(x) 
\! - \! \widetilde{V}^{\prime}(\tilde{x}_{j}^{\ast}) \tilde{Z}
(\tilde{x}_{j}^{\ast})}{x \! - \! \tilde{x}_{j}^{\ast}} \right) 
\! - \! \dfrac{2}{\tilde{Z}(x)} \left(\sum_{i=2(\tilde{N}+1)}^{3} 
\tilde{c}_{2(\tilde{N}+1)-i} \right. \nonumber \\
\times&\left. \, \left(\dfrac{\mathcal{N}(\mathcal{N} \! - \! 1)}{2}x^{i-2} 
\! + \! \sum_{r=2}^{i-1} \left((\mathcal{N} \! - \! 1) \sum_{j=1}^{
\mathcal{N}}(\tilde{x}_{j}^{\ast})^{r-1} \! + \! \sum_{j=1}^{\mathcal{N}-1} 
\sum_{k^{\prime}=j+1}^{\mathcal{N}} \sum_{m=2}^{r-1}(\tilde{x}_{j}^{\ast})^{r
-m}(\tilde{x}_{k^{\prime}}^{\ast})^{m-1} \right) \right. \right. \nonumber \\
\times&\left. \left. \, x^{i-r-1} \right) \! + \! \dfrac{\mathcal{N}
(\mathcal{N} \! - \! 1)}{2} \tilde{c}_{2 \tilde{N}} \right) \! + \! 
\dfrac{2}{\tilde{Z}(x)} \dfrac{(\mathcal{N} \! - \! 1)}{\mathcal{N}} 
\left(\sum_{i=2(\tilde{N}+1)}^{3} \tilde{c}_{2(\tilde{N}+1)-i} \right. 
\nonumber \\
\times&\left. \, \left(\sum_{r=0}^{i-1} \sum_{j=1}^{\mathcal{N}} 
\sum_{q=1}^{\mathfrak{s}-2} \dfrac{l_{q}(\tilde{x}_{j}^{\ast})^{r}}{
\tilde{x}_{j}^{\ast} \! - \! \alpha_{p_{q}}}x^{i-r-1} \right) \! + \! 
\tilde{c}_{2 \tilde{N}} \sum_{j=1}^{\mathcal{N}} \sum_{q=1}^{\mathfrak{s}-2} 
\dfrac{l_{q}(x \! + \! \tilde{x}_{j}^{\ast})}{\tilde{x}_{j}^{\ast} \! - \! 
\alpha_{p_{q}}} \! + \! \tilde{c}_{2 \tilde{N}+1} \right. \nonumber \\
\times&\left. \, \sum_{j=1}^{\mathcal{N}} \sum_{q=1}^{\mathfrak{s}-2} 
\dfrac{l_{q}}{\tilde{x}_{j}^{\ast} \! - \! \alpha_{p_{q}}} \right) \! + \! 
\dfrac{2(l_{\mathfrak{s}} \! - \! 1)}{\tilde{Z}(x)} \dfrac{(\mathcal{N} \! 
- \! 1)}{\mathcal{N}} \left(\sum_{i=2(\tilde{N}+1)}^{3} \tilde{c}_{2(
\tilde{N}+1)-i} \left(\sum_{r=0}^{i-1} \sum_{j=1}^{\mathcal{N}} \dfrac{(
\tilde{x}_{j}^{\ast})^{r}}{\tilde{x}_{j}^{\ast} \! - \! \alpha_{p_{
\mathfrak{s}}}}x^{i-r-1} \right) \right. \nonumber \\
+&\left. \, \tilde{c}_{2 \tilde{N}} \sum_{j=1}^{\mathcal{N}} \dfrac{x \! + \! 
\tilde{x}_{j}^{\ast}}{\tilde{x}_{j}^{\ast} \! - \! \alpha_{p_{\mathfrak{s}}}} 
\! + \! \tilde{c}_{2 \tilde{N}+1} \sum_{j=1}^{\mathcal{N}} 
\dfrac{1}{\tilde{x}_{j}^{\ast} \! - \! \alpha_{p_{\mathfrak{s}}}} \right), 
\label{eql3.6j}
\end{align}
whence, via Equation~\eqref{eql3.6j}, one arrives at, for $n \! \in \! 
\mathbb{N}$ and $k \! \in \! \lbrace 1,2,\dotsc,K \rbrace$ such 
that $\alpha_{p_{\mathfrak{s}}} \! := \! \alpha_{k} \! \neq \! \infty$ 
(cf. Equation~\eqref{eql3.6h}),
\begin{align}
\tilde{Q}(x)=& \, \dfrac{1}{\tilde{Z}(x)} \left(\left(\sum_{i=2(\tilde{N}+
1)}^{0} \tilde{c}_{2(\tilde{N}+1)-i}x^{i} \right) \left(2 \dfrac{(\mathcal{N} 
\! - \! 1)}{\mathcal{N}} \sum_{j=1}^{\mathcal{N}} \sum_{q=1}^{
\mathfrak{s}-2} \dfrac{l_{q}}{(\tilde{x}_{j}^{\ast} \! - \! \alpha_{p_{q}})
(x \! - \! \alpha_{p_{q}})} \! + \! 2 \dfrac{(\mathcal{N} \! - \! 1)}{
\mathcal{N}} \sum_{j=1}^{\mathcal{N}} \dfrac{(l_{\mathfrak{s}} 
\! - \! 1)}{(\tilde{x}_{j}^{\ast} \! - \! \alpha_{p_{\mathfrak{s}}})
(x \! - \! \alpha_{p_{\mathfrak{s}}})} \right. \right. \nonumber \\
+&\left. \left. \, \sum_{\substack{q=1\\q \neq \mathfrak{s}-1}}^{
\mathfrak{s}} \dfrac{l_{q}}{(x \! - \! \alpha_{p_{q}})^{2}} \! + \! 
\dfrac{(\mathcal{N} \! - \! 2)}{\mathcal{N}} \left(\sum_{\substack{q
=1\\q \neq \mathfrak{s}-1}}^{\mathfrak{s}} \dfrac{l_{q}}{x \! - \! 
\alpha_{p_{q}}} \right)^{2} \! - \! 2 \dfrac{(\mathcal{N} \! - \! 1)}{
\mathcal{N}} \dfrac{1}{x \! - \! \alpha_{p_{\mathfrak{s}}}} \sum_{
\substack{q=1\\q \neq \mathfrak{s}-1}}^{\mathfrak{s}} 
\dfrac{l_{q}}{x \! - \! \alpha_{p_{q}}} \! + \! \dfrac{(\mathcal{N} 
\! - \! 1)}{\mathcal{N}}n \widetilde{V}^{\prime}(x) \sum_{
\substack{q=1\\q \neq \mathfrak{s}-1}}^{\mathfrak{s}} 
\dfrac{l_{q}}{x \! - \! \alpha_{p_{q}}} \right) \right. \nonumber \\
-&\left. \, n \dfrac{(\mathcal{N} \! - \! 1)}{\mathcal{N}} \sum_{j=1}^{
\mathcal{N}} \left(\dfrac{\widetilde{V}^{\prime}(x) \tilde{Z}(x) \! - \! 
\widetilde{V}^{\prime}(\tilde{x}_{j}^{\ast}) \tilde{Z}(\tilde{x}_{j}^{\ast})}{
x \! - \! \tilde{x}_{j}^{\ast}} \right) \! + \! 2 \left(\sum_{i=2(\tilde{N}+
1)}^{3} \tilde{c}_{2(\tilde{N}+1)-i} \left(\dfrac{\mathcal{N}(\mathcal{N} \! 
- \! 1)}{2}x^{i-2} \! + \! \sum_{r=2}^{i-1} \left(\vphantom{M^{M^{M^{M^{M}}}}}
(\mathcal{N} \! - \! 1) \right. \right. \right. \right. \nonumber \\
\times&\left. \left. \left. \left. \, \sum_{j=1}^{\mathcal{N}}(\tilde{x}_{
j}^{\ast})^{r-1} \! + \! \sum_{j=1}^{\mathcal{N}-1} \sum_{k^{\prime}=j+1}^{
\mathcal{N}} \sum_{m=2}^{r-1}(\tilde{x}_{j}^{\ast})^{r-m}(\tilde{x}_{
k^{\prime}}^{\ast})^{m-1} \right)x^{i-r-1} \right) \! + \! \dfrac{\mathcal{N}
(\mathcal{N} \! - \! 1)}{2} \tilde{c}_{2 \tilde{N}} \right) \! - \! 2 
\dfrac{(\mathcal{N} \! - \! 1)}{\mathcal{N}} \left(\sum_{i=2(\tilde{N}+1)}^{3} 
\tilde{c}_{2(\tilde{N}+1)-i} \right. \right. \nonumber \\
\times&\left. \left. \, \left(\sum_{r=0}^{i-1} \sum_{j=1}^{\mathcal{N}} 
\sum_{q=1}^{\mathfrak{s}-2} \dfrac{l_{q}(\tilde{x}_{j}^{\ast})^{r}}{
\tilde{x}_{j}^{\ast} \! - \! \alpha_{p_{q}}}x^{i-r-1} \right) \! + \! 
\tilde{c}_{2 \tilde{N}} \sum_{j=1}^{\mathcal{N}} \sum_{q=1}^{\mathfrak{s}-2} 
\dfrac{l_{q}(x \! + \! \tilde{x}_{j}^{\ast})}{\tilde{x}_{j}^{\ast} \! - \! 
\alpha_{p_{q}}} \! + \! \tilde{c}_{2 \tilde{N}+1} \sum_{j=1}^{\mathcal{N}} 
\sum_{q=1}^{\mathfrak{s}-2} \dfrac{l_{q}}{\tilde{x}_{j}^{\ast} \! - \! 
\alpha_{p_{q}}} \right) \! - \! 2(l_{\mathfrak{s}} \! - \! 1) \right. 
\nonumber \\
\times&\left. \, \dfrac{(\mathcal{N} \! - \! 1)}{\mathcal{N}} \left(\sum_{i=
2(\tilde{N}+1)}^{3} \tilde{c}_{2(\tilde{N}+1)-i} \left(\sum_{r=0}^{i-1} 
\sum_{j=1}^{\mathcal{N}} \dfrac{(\tilde{x}_{j}^{\ast})^{r}}{\tilde{x}_{
j}^{\ast} \! - \! \alpha_{p_{\mathfrak{s}}}}x^{i-r-1} \right) \! + \! 
\tilde{c}_{2 \tilde{N}} \sum_{j=1}^{\mathcal{N}} \dfrac{x \! + \! 
\tilde{x}_{j}^{\ast}}{\tilde{x}_{j}^{\ast} \! - \! \alpha_{p_{\mathfrak{s}}}} 
\! + \! \tilde{c}_{2 \tilde{N}+1} \sum_{j=1}^{\mathcal{N}} \dfrac{1}{
\tilde{x}_{j}^{\ast} \! - \! \alpha_{p_{\mathfrak{s}}}} \right) \right). 
\label{eql3.6k}
\end{align}
As in the proof of Lemma~2.15 in \cite{a58}, make, for $n \! \in \! 
\mathbb{N}$ and $k \! \in \! \lbrace 1,2,\dotsc,K \rbrace$ such that 
$\alpha_{p_{\mathfrak{s}}} \! := \! \alpha_{k} \! \neq \! \infty$, the 
following change of dependent variable (cf. Equations~\eqref{eql3.6g} 
and~\eqref{eql3.6k}),
\begin{equation} \label{eql3.6l}
\tilde{f}(x) \! = \! \tilde{F}(x) \exp \left(-\dfrac{1}{2} \int_{}^{x} 
\tilde{\mathcal{P}}(\xi) \, \md \xi \right),
\end{equation}
where
\begin{equation} \label{eql3.6m} 
\tilde{\mathcal{P}}(x) \! := \! 2 \dfrac{n}{\mathcal{N}} \left(-
\dfrac{1}{2}(\mathcal{N} \! - \! 1) \widetilde{V}^{\prime}(x) \! + \! 
\sum_{\substack{q=1\\q \neq \mathfrak{s}-1}}^{\mathfrak{s}} 
\dfrac{l_{q}/n}{x \! - \! \alpha_{p_{q}}} \! + \! \dfrac{(\mathcal{N} 
\! - \! 1)/n}{x \! - \! \alpha_{p_{\mathfrak{s}}}} \right).
\end{equation}
Via Equations~\eqref{eql3.6g}, \eqref{eql3.6k}, and~\eqref{eql3.6l}, 
and the Definition~\eqref{eql3.6m}, one arrives at, for $n \! \in \! 
\mathbb{N}$ and $k \! \in \! \lbrace 1,2,\dotsc,K \rbrace$ such that 
$\alpha_{p_{\mathfrak{s}}} \! := \! \alpha_{k} \! \neq \! \infty$, the 
following linear, second-order, rational-coefficient ODE for the function 
$\tilde{F}(x)$:
\begin{equation} \label{eql3.6n} 
\tilde{F}^{\prime \prime}(x) \! = \! \left(\tilde{Q}(x) \! + \! \dfrac{1}{2} 
\tilde{\mathcal{P}}^{\prime}(x) \! + \! \dfrac{1}{4}(\tilde{\mathcal{P}}
(x))^{2} \right) \tilde{F}(x),
\end{equation}
where,
\begin{align}
& \, \tilde{Q}(x) \! + \! \dfrac{1}{2} \tilde{\mathcal{P}}^{\prime}(x) 
\! + \! \dfrac{1}{4}(\tilde{\mathcal{P}}(x))^{2} \! = \! -\dfrac{1}{2} 
\dfrac{(\mathcal{N} \! - \! 1)}{\mathcal{N}}n \widetilde{V}^{\prime 
\prime}(x) \! - \! \dfrac{1}{\mathcal{N}} \dfrac{(\mathcal{N} \! - \! 
1)}{\mathcal{N}} \dfrac{1}{(x \! - \! \alpha_{p_{\mathfrak{s}}})^{2}} 
\! + \! \dfrac{(\mathcal{N} \! - \! 1)}{\mathcal{N}} \sum_{
\substack{q=1\\q \neq \mathfrak{s}-1}}^{\mathfrak{s}} 
\dfrac{l_{q}}{(x \! - \! \alpha_{p_{q}})^{2}} \nonumber \\
+& \, \dfrac{1}{4} \left(\dfrac{\mathcal{N} \! - \! 1}{\mathcal{N}} 
\right)^{2}(n \widetilde{V}^{\prime}(x))^{2} \! + \! n \widetilde{V}^{
\prime}(x) \left(\dfrac{\mathcal{N} \! - \! 1}{\mathcal{N}} \right)^{2} 
\sum_{\substack{q=1\\q \neq \mathfrak{s}-1}}^{\mathfrak{s}} 
\dfrac{l_{q}}{x \! - \! \alpha_{p_{q}}} \! - \! 2 \left(\dfrac{\mathcal{N} 
\! - \! 1}{\mathcal{N}} \right)^{2} \dfrac{1}{x \! - \! \alpha_{
p_{\mathfrak{s}}}} \sum_{\substack{q=1\\q \neq \mathfrak{s}
-1}}^{\mathfrak{s}} \dfrac{l_{q}}{x \! - \! \alpha_{p_{q}}} \! - \! \left(
\dfrac{\mathcal{N} \! - \! 1}{\mathcal{N}} \right)^{2} \dfrac{n 
\widetilde{V}^{\prime}(x)}{x \! - \! \alpha_{p_{\mathfrak{s}}}} \nonumber \\
+& \, \left(\dfrac{\mathcal{N} \! - \! 1}{\mathcal{N}} \right)^{2} \left(
\sum_{\substack{q=1\\q \neq \mathfrak{s}-1}}^{\mathfrak{s}} 
\dfrac{l_{q}}{x \! - \! \alpha_{p_{q}}} \right)^{2} \! + \! 2 \dfrac{(\mathcal{N} 
\! - \! 1)}{\mathcal{N}} \sum_{j=1}^{\mathcal{N}} \sum_{q=1}^{\mathfrak{s}
-2} \dfrac{l_{q}}{(\tilde{x}_{j}^{\ast} \! - \! \alpha_{p_{q}})(x \! - \! 
\alpha_{p_{q}})} \! + \! 2 \dfrac{(\mathcal{N} \! - \! 1)}{\mathcal{N}} 
\sum_{j=1}^{\mathcal{N}} \dfrac{(l_{\mathfrak{s}} \! - \! 1)}{(\tilde{x}_{j}^{\ast} 
\! - \! \alpha_{p_{\mathfrak{s}}})(x \! - \! \alpha_{p_{\mathfrak{s}}})} 
\nonumber \\
+& \, \dfrac{1}{\tilde{Z}(x)} \left(-n \dfrac{(\mathcal{N} \! - \! 1)}{
\mathcal{N}} \sum_{j=1}^{\mathcal{N}} \left(\dfrac{\widetilde{V}^{\prime}(x) 
\tilde{Z}(x) \! - \! \widetilde{V}^{\prime}(\tilde{x}_{j}^{\ast}) \tilde{Z}
(\tilde{x}_{j}^{\ast})}{x \! - \! \tilde{x}_{j}^{\ast}} \right) \! + \! 
2 \left(\sum_{i=2(\tilde{N}+1)}^{3} \tilde{c}_{2(\tilde{N}+1)-i} \left(
\dfrac{\mathcal{N}(\mathcal{N} \! - \! 1)}{2}x^{i-2} \! + \! \sum_{r=2}^{i-1} 
\left(\vphantom{M^{M^{M^{M^{M}}}}}(\mathcal{N} \! - \! 1) \right. \right. 
\right. \right. \nonumber \\
\times&\left. \left. \left. \left. \, \sum_{j=1}^{\mathcal{N}}(\tilde{x}_{
j}^{\ast})^{r-1} \! + \! \sum_{j=1}^{\mathcal{N}-1} \sum_{k^{\prime}=
j+1}^{\mathcal{N}} \sum_{m=2}^{r-1}(\tilde{x}_{j}^{\ast})^{r-m}(\tilde{x}_{
k^{\prime}}^{\ast})^{m-1} \right)x^{i-r-1} \right) \! + \! \dfrac{\mathcal{N}
(\mathcal{N} \! - \! 1)}{2} \tilde{c}_{2 \tilde{N}} \right) \! - \! 2 
\dfrac{(\mathcal{N} \! - \! 1)}{\mathcal{N}} \left(\sum_{i=2(\tilde{N}+1)}^{3} 
\tilde{c}_{2(\tilde{N}+1)-i} \right. \right. \nonumber \\
\times&\left. \left. \, \left(\sum_{r=0}^{i-1} \sum_{j=1}^{\mathcal{N}} 
\sum_{q=1}^{\mathfrak{s}-2} \dfrac{l_{q}(\tilde{x}_{j}^{\ast})^{r}}{
\tilde{x}_{j}^{\ast} \! - \! \alpha_{p_{q}}}x^{i-r-1} \right) \! + \! 
\tilde{c}_{2 \tilde{N}} \sum_{j=1}^{\mathcal{N}} \sum_{q=1}^{\mathfrak{s}-2} 
\dfrac{l_{q}(x \! + \! \tilde{x}_{j}^{\ast})}{\tilde{x}_{j}^{\ast} \! - \! 
\alpha_{p_{q}}} \! + \! \tilde{c}_{2 \tilde{N}+1} \sum_{j=1}^{\mathcal{N}} 
\sum_{q=1}^{\mathfrak{s}-2} \dfrac{l_{q}}{\tilde{x}_{j}^{\ast} \! - \! 
\alpha_{p_{q}}} \right) \! - \! 2(l_{\mathfrak{s}} \! - \! 1) 
\dfrac{(\mathcal{N} \! - \! 1)}{\mathcal{N}} \right. \nonumber \\
\times&\left. \, \left(\sum_{i=2(\tilde{N}+1)}^{3} \tilde{c}_{2(\tilde{N}+1)
-i} \left(\sum_{r=0}^{i-1} \sum_{j=1}^{\mathcal{N}} \dfrac{(\tilde{x}_{j}^{
\ast})^{r}}{\tilde{x}_{j}^{\ast} \! - \! \alpha_{p_{\mathfrak{s}}}}x^{i-r-1} 
\right) \! + \! \tilde{c}_{2 \tilde{N}} \sum_{j=1}^{\mathcal{N}} \dfrac{x \! 
+ \! \tilde{x}_{j}^{\ast}}{\tilde{x}_{j}^{\ast} \! - \! \alpha_{p_{\mathfrak{s}}}} 
\! + \! \tilde{c}_{2 \tilde{N}+1} \sum_{j=1}^{\mathcal{N}} 
\dfrac{1}{\tilde{x}_{j}^{\ast} \! - \! \alpha_{p_{\mathfrak{s}}}} \right) 
\right). \label{eql3.6o}
\end{align}
Since, for $n \! \in \! \mathbb{N}$ and $k \! \in \! \lbrace 1,2,\dotsc,K 
\rbrace$ such that $\alpha_{p_{\mathfrak{s}}} \! := \! \alpha_{k} \! \neq 
\! \infty$, $\overline{\mathbb{R}} \setminus \lbrace \alpha_{1},\alpha_{2},
\dotsc,\alpha_{K} \rbrace \supset \cup_{j=1}^{\tilde{N}+1}[\tilde{B}_{j-1},
\tilde{A}_{j}] \supseteq \lbrace \tilde{x}_{1}^{\ast},\tilde{x}_{2}^{\ast},
\dotsc,\tilde{x}_{\mathcal{N}}^{\ast} \rbrace$, with $\tilde{x}_{i}^{\ast} 
\! < \! \tilde{x}_{j}^{\ast}$ for $i \! < \! j$ $(\in \! \lbrace 
1,2,\dotsc,\mathcal{N} \rbrace)$, it follows {}from the proof of 
Lemma~\ref{lem3.5} that, for any two consecutive associated, weighted 
Fekete points $\tilde{x}_{j}^{\ast}$ and $\tilde{x}_{j+1}^{\ast}$, there 
are, for $r \! = \! 1,2,\dotsc,\tilde{N} \! + \! 1$, the following cases to 
consider: (1) $\tilde{B}_{r-1} \! \leqslant \! \tilde{x}_{j}^{\ast} \! < \! 
\tilde{x}_{j+1}^{\ast} \! \leqslant \! (\tilde{B}_{r-1} \! + \! \tilde{A}_{r})/2$; 
(2) $(\tilde{B}_{r-1} \! + \! \tilde{A}_{r})/2 \! \leqslant \! \tilde{x}_{j}^{\ast} 
\! < \! \tilde{x}_{j+1}^{\ast} \! \leqslant \! \tilde{A}_{r}$; and (3) 
$\tilde{B}_{r-1} \! \leqslant \! \tilde{x}_{j}^{\ast} \! < \! (\tilde{B}_{r-1} 
\! + \! \tilde{A}_{r})/2 \! < \! \tilde{x}_{j+1}^{\ast} \! \leqslant \! 
\tilde{A}_{r}$. This is, up to a linear scaling, the `one-interval-case' 
result given in Lemma~2.15 of \cite{a58}, wherein a lower bound for the 
distance between two consecutive Fekete points is estimated. In order to 
use the result of Lemma~2.15 of \cite{a58} and adapt it to the present 
situation,\footnote{If, for $n \! \in \! \mathbb{N}$ and $k \! \in \! 
\lbrace 1,2,\dotsc,K \rbrace$ such that $\alpha_{p_{\mathfrak{s}}} \! 
:= \! \alpha_{k} \! \neq \! \infty$, two consecutive associated, weighted 
Fekete points $\tilde{x}_{j}^{\ast}$ and $\tilde{x}_{j+1}^{\ast}$ lie, 
respectively, in the disjoint compact real intervals $[\tilde{B}_{r_{1}-1},
\tilde{A}_{r_{1}}]$ and $[\tilde{B}_{r_{2}-1},\tilde{A}_{r_{2}}]$, $r_{1} \! 
< \! r_{2}$ $(\in \! \lbrace 1,2,\dotsc,\tilde{N} \! + \! 1 \rbrace)$, 
then, clearly, $\tilde{x}_{j+1}^{\ast} \! - \! \tilde{x}_{j}^{\ast} \! 
\geqslant \! \lvert \tilde{B}_{r_{2}-1} \! - \! \tilde{A}_{r_{1}} \rvert \! 
> \! 0$ $(>_{\underset{z_{o}=1+o(1)}{\mathscr{N},n \to \infty}} 
\mathcal{O}(((n \! - \! 1)K \! + \! k)^{-1})$, by the Archimedean 
property).} one needs to, for $n \! \in \! \mathbb{N}$ and $k \! \in 
\! \lbrace 1,2,\dotsc,K \rbrace$ such that $\alpha_{p_{\mathfrak{s}}} 
\! := \! \alpha_{k} \! \neq \! \infty$: (i) map the compact real intervals 
$[\tilde{B}_{r-1},\tilde{A}_{r}]$, $r \! = \! 1,2,\dotsc,\tilde{N} \! + \! 1$, 
onto the compact real---symmetric---interval $[-1,1]$; and (ii) estimate, 
for $x \! \in \! [\tilde{B}_{r-1},\tilde{A}_{r}]$, $r \! = \! 1,2,\dotsc,\tilde{N} 
\! + \! 1$, an upper bound for $\tilde{Q}(x) \! + \! \tilde{\mathcal{P}}^{\prime}
(x)/2 \! + \!(\tilde{\mathcal{P}}(x))^{2}/4$. For the former problem~(i), one 
makes, for $r \! = \! 1,2,\dotsc,\tilde{N} \! + \! 1$, the linear change of 
variable $\tilde{\lambda}_{r} \colon \mathbb{C} \! \to \! \mathbb{C}$, 
$x \! \mapsto \! \tilde{\lambda}_{r}(x) \! := \! (2x \! - \! (\tilde{A}_{r} \! 
+ \! \tilde{B}_{r-1}))/(\tilde{A}_{r} \! - \! \tilde{B}_{r-1})$, which maps, 
one-to-one, the compact real interval $[\tilde{B}_{r-1},\tilde{A}_{r}]$ onto 
the compact real---symmetric---interval $[-1,1]$; and for the latter 
problem~(ii), noting that, for $n \! \in \! \mathbb{N}$ and $k \! \in \! 
\lbrace 1,2,\dotsc,K \rbrace$ such that $\alpha_{p_{\mathfrak{s}}} \! 
:= \! \alpha_{k} \! \neq \! \infty$, $\inf \lbrace \mathstrut \lvert 
\alpha_{p_{q}} \! - \! \tilde{x}_{j}^{\ast} \rvert; \, q \! = \! 1,\dotsc,
\mathfrak{s} \! - \! 2,\mathfrak{s}, \, j \! = \! 1,2,\dotsc,\mathcal{N} 
\rbrace \! > \! 0$ and $\inf \lbrace \mathstrut \lvert \alpha_{p_{q}} 
\! - \! x \rvert; \, q \! = \! 1,\dotsc,\mathfrak{s} \! - \! 2,\mathfrak{s}, 
\, x \! \in \! [\tilde{B}_{r-1},\tilde{A}_{r}], \, r \! = \! 1,2,\dotsc,\tilde{N} 
\! + \! 1 \rbrace \! > \! 0$, writing $\tilde{Z}(x) \! = \! (x \! - \! 
\tilde{B}_{r-1})(x \! - \! \tilde{A}_{r}) \prod_{\underset{j \neq r}{j=
1}}^{\tilde{N}+1}(x \! - \! \tilde{B}_{j-1})(x \! - \! \tilde{A}_{j})$, $r \! 
= \! 1,2,\dotsc,\tilde{N} \! + \! 1$, and using the real analyticity and 
regularity of $\widetilde{V} \colon \overline{\mathbb{R}} \setminus 
\lbrace \alpha_{1},\alpha_{2},\dotsc,\alpha_{K} \rbrace \! \to \! 
\mathbb{R}$, one shows {}from Equation~\eqref{eql3.6o} that, for 
$[\tilde{B}_{r-1},\tilde{A}_{r}] \! \ni \! x$, $r \! = \! 1,2,\dotsc,
\tilde{N} \! + \! 1$,
\begin{equation} \label{eql3.6p} 
\left\vert \tilde{Q}(x) \! + \! \dfrac{1}{2} \tilde{\mathcal{P}}^{\prime}(x) 
\! + \! \dfrac{1}{4}(\tilde{\mathcal{P}}(x))^{2} \right\vert \! \leqslant \! 
\dfrac{\tilde{\mathfrak{C}}_{r} \mathcal{N}^{2}}{(x \! - \! \tilde{B}_{r-1})
(\tilde{A}_{r} \! - \! x)},
\end{equation}
where $\tilde{\mathfrak{C}}_{r} \! = \! \tilde{\mathfrak{C}}_{r}(n,k,z_{o}) \! 
> \! 0$ and $\mathcal{O}(1)$ (in the double-scaling limit $\mathscr{N},
n \! \to \! \infty$ such that $z_{o} \! = \! 1 \! + \! o(1))$. One now 
uses this estimate in conjunction with the result of Lemma~2.15 of 
\cite{a58} to show that, for $n \! \in \! \mathbb{N}$ and $k \! \in \! 
\lbrace 1,2,\dotsc,K \rbrace$ such that $\alpha_{p_{\mathfrak{s}}} 
\! := \! \alpha_{k} \! \neq \! \infty$, two consecutive associated, 
weighted Fekete points $\tilde{x}_{j}^{\ast}$ and $\tilde{x}_{j+1}^{\ast}$ 
satisfy the following `nearest-neighbour-distance' inequality:
\begin{equation} \label{eql3.6q} 
\tilde{x}_{j+1}^{\ast} \! - \! \tilde{x}_{j}^{\ast} \! \geqslant \! 
\dfrac{\min \lbrace (3 \tilde{\mathfrak{C}}_{r})^{-1/2},1/4 
\rbrace}{((n \! - \! 1)K \! + \! k)} \left((\tilde{A}_{r} \! - \! 
\tilde{x}_{j}^{\ast})(\tilde{A}_{r} \! - \! \tilde{x}_{j+1}^{\ast}) 
\right)^{1/2}, \quad j \! = \! 1,2,\dotsc,(n \! - \! 1)K \! + \! 
k \! - \! 1, \quad r \! = \! 1,2,\dotsc,\tilde{N} \! + \! 1,
\end{equation}
where $((\tilde{A}_{r} \! - \! \tilde{x}_{j}^{\ast})(\tilde{A}_{r} \! - \! 
\tilde{x}_{j+1}^{\ast}))^{1/2} \! > \! 0$ and $\mathcal{O}(1)$ (in 
the double-scaling limit $\mathscr{N},n \! \to \! \infty$ such that 
$z_{o} \! = \! 1 \! + \! o(1))$. One now uses the estimate~\eqref{eql3.6q} 
in conjunction with the fact that (cf. Lemma~\ref{lem3.5}), for $k \! 
\in \! \lbrace 1,2,\dotsc,K \rbrace$ such that $\alpha_{p_{\mathfrak{s}}} 
\! := \! \alpha_{k} \! \neq \! \infty$, the associated normalised counting 
measure converges, in the double-scaling limit $\mathscr{N},n \! \to \! 
\infty$ such that $z_{o} \! = \! 1 \! + \! o(1)$, weakly (in the weak-$\ast$ 
topology of measures) to the corresponding equilibrium measure 
$\mu_{\widetilde{V}}^{f}$ $(\in \! \mathscr{M}_{1}(\mathbb{R}))$ in order 
to proceed, \emph{mutatis mutandis}, as in the proof of Lemma~2.26 of 
\cite{a58} to conclude that, for $n \! \in \! \mathbb{N}$ and $k \! \in \! 
\lbrace 1,2,\dotsc,K \rbrace$ such that $\alpha_{p_{\mathfrak{s}}} \! 
:= \! \alpha_{k} \! \neq \! \infty$, the associated family of $\tilde{K}$ 
equilibrium measures $\mu_{\widetilde{V}}^{f}$ is absolutely continuous 
with respect to Lebesgue measure, in the double-scaling limit 
$\mathscr{N},n \! \to \! \infty$ such that $z_{o} \! = \! 1 \! + \! o(1)$. 

$\pmb{(2)}$ The proof of this case, that is, $n \! \in \! \mathbb{N}$ 
and $k \! \in \! \lbrace 1,2,\dotsc,K \rbrace$ such that 
$\alpha_{p_{\mathfrak{s}}} \! := \! \alpha_{k} \! = \! \infty$, is virtually 
identical to the proof presented in $\pmb{(1)}$ above; one mimics, 
\emph{verbatim}, the scheme of the calculations presented in case $\pmb{(1)}$ 
in order to arrive at the corresponding claims stated in the lemma; in 
order to do so, however, the analogues of Equations~\eqref{eql3.6a} 
and~\eqref{eql3.6b}, the Definitions~\eqref{eql3.6c} and~\eqref{eql3.6d}, 
Equations~\eqref{eql3.6e}--\eqref{eql3.6h} 
and~\eqref{eql3.6j}--\eqref{eql3.6l}, the Definition~\eqref{eql3.6m}, 
Equations~\eqref{eql3.6n} and~\eqref{eql3.6o}, estimate~\eqref{eql3.6p}, 
and inequality~\eqref{eql3.6q} are necessary, which, in the present case, 
read: for $n \! \in \! \mathbb{N}$ and $k \! \in \! \lbrace 1,2,\dotsc,
K \rbrace$ such that $\alpha_{p_{\mathfrak{s}}} \! := \! \alpha_{k} 
\! = \! \infty$, with $\mathcal{N} \! := \! (n \! - \! 1)K \! + \! k$,
\begin{align} \label{eql3.6r} 
\hat{f} \colon \mathbb{N} \times \lbrace 1,2,\dotsc,K \rbrace \times 
\overline{\mathbb{R}} \setminus \lbrace \alpha_{1},\alpha_{2},\dotsc,
\alpha_{K} \rbrace \! \ni \! (n,k,x) \! \mapsto \! \dfrac{\prod_{j=1}^{
\mathcal{N}}(x \! - \! \hat{x}_{j}^{\ast})}{\prod_{q=1}^{\mathfrak{s}-1}
(x \! - \! \alpha_{p_{q}})^{l_{q}}} \! =: \! \hat{f}(x),
\end{align}
where $\lbrace \hat{x}_{1}^{\ast},\hat{x}_{2}^{\ast},\dotsc,\hat{x}_{
\mathcal{N}}^{\ast} \rbrace \subset \overline{\mathbb{R}} \setminus 
\lbrace \alpha_{1},\alpha_{2},\dotsc,\alpha_{K} \rbrace$, with the 
enumeration $\hat{x}_{i}^{\ast} \! < \! \hat{x}_{j}^{\ast}$ for $i \! < \! j$, 
is the associated, weighted Fekete set described in the corresponding 
item of Lemma~\ref{lem3.5},
\begin{equation} \label{eql3.6s} 
\dfrac{\hat{f}^{\prime \prime}(x)}{\hat{f}^{\prime}(x)} 
\! = \! \dfrac{\hat{f}_{N}(x)}{\hat{f}_{D}(x)},
\end{equation}
\begin{align}
\hat{f}_{N}(x) :=& \, 2 \sum_{j=1}^{\mathcal{N}} \prod_{\substack{l=
1\\l \neq j}}^{\mathcal{N}}(x \! - \! \hat{x}_{l}^{\ast}) \sum_{\substack{
k^{\prime}=1\\k^{\prime} \neq j}}^{\mathcal{N}} \dfrac{1}{\hat{x}_{j}^{
\ast} \! - \! \hat{x}_{k^{\prime}}^{\ast}} \! - \! 2 \sum_{j=1}^{\mathcal{N}} 
\prod_{\substack{l=1\\l \neq j}}^{\mathcal{N}}(x \! - \! \hat{x}_{l}^{\ast}) 
\sum_{q=1}^{\mathfrak{s}-1} \dfrac{l_{q}}{x \! - \! \alpha_{p_{q}}} 
\nonumber \\
+& \, \prod_{j=1}^{\mathcal{N}}(x \! - \! \hat{x}_{j}^{\ast}) \left(\sum_{q
=1}^{\mathfrak{s}-1} \dfrac{l_{q}}{(x \! - \! \alpha_{p_{q}})^{2}} \! + \! 
\left(\sum_{q=1}^{\mathfrak{s}-1} \dfrac{l_{q}}{x \! - \! \alpha_{p_{q}}} 
\right)^{2} \right), \label{eql3.6t}
\end{align}
\begin{gather}
\hat{f}_{D}(x) \! := \! \sum_{j=1}^{\mathcal{N}} \sum_{\substack{
k^{\prime}=1\\k^{\prime} \neq j}}^{\mathcal{N}}(x \! - \! \hat{x}_{
k^{\prime}}^{\ast}) \! - \! \prod_{j=1}^{\mathcal{N}}(x \! - \! 
\hat{x}_{j}^{\ast}) \sum_{q=1}^{\mathfrak{s}-1} \dfrac{l_{q}}{x \! 
- \! \alpha_{p_{q}}}, \label{eql3.6u} \\
\dfrac{1}{2} \dfrac{\hat{f}^{\prime \prime}(\hat{x}_{m}^{\ast})}{\hat{f}^{
\prime}(\hat{x}_{m}^{\ast})} \! = \! \sum_{\substack{k^{\prime}=1\\
k^{\prime} \neq m}}^{\mathcal{N}} \dfrac{1}{\hat{x}_{m}^{\ast} \! - 
\! \hat{x}_{k^{\prime}}^{\ast}} \! - \! \sum_{q=1}^{\mathfrak{s}-1} 
\dfrac{l_{q}}{\hat{x}_{m}^{\ast} \! - \! \alpha_{p_{q}}}, \quad 
m \! = \! 1,2,\dotsc,\mathcal{N}, \label{eql3.6v} \\
-\dfrac{1}{2}(\mathcal{N} \! - \! 1) \widetilde{V}^{\prime}
(\hat{x}_{m}^{\ast}) \! + \! \dfrac{\mathcal{N}}{n} \sum_{
\substack{k^{\prime}=1\\k^{\prime} \neq m}}^{\mathcal{N}} 
\dfrac{1}{\hat{x}_{m}^{\ast} \! - \! \hat{x}_{k^{\prime}}^{\ast}} 
\! - \! (\mathcal{N} \! - \! 1) \sum_{q=1}^{\mathfrak{s}-1} 
\dfrac{l_{q}/n}{\hat{x}_{m}^{\ast} \! - \! \alpha_{p_{q}}} \! = \! 0, 
\quad m \! = \! 1,2,\dotsc,\mathcal{N}, \label{eql3.6w} \\
\hat{f}^{\prime \prime}(x) \! + \! 2 \dfrac{n}{\mathcal{N}} \left(
-\dfrac{1}{2}(\mathcal{N} \! - \! 1) \widetilde{V}^{\prime}(x) \! + \! 
\sum_{q=1}^{\mathfrak{s}-1} \dfrac{l_{q}/n}{x \! - \! \alpha_{p_{q}}} 
\right) \hat{f}^{\prime}(x) \! = \! \hat{Q}(x) \hat{f}(x), \label{eql3.6x}
\end{gather}
where (with $\hat{Q}(\hat{x}_{m}^{\ast}) \! \neq \! 0$, $m \! = \! 1,2,
\dotsc,\mathcal{N})$
\begin{align}
\hat{Q}(x)=& \, 2 \dfrac{(\mathcal{N} \! - \! 1)}{\mathcal{N}} \sum_{j=1}^{
\mathcal{N}} \sum_{q=1}^{\mathfrak{s}-1} \dfrac{l_{q}}{(\hat{x}_{j}^{\ast} \! 
- \! \alpha_{p_{q}})(x \! - \! \alpha_{p_{q}})} \! + \! \sum_{q=1}^{\mathfrak{
s}-1} \dfrac{l_{q}}{(x \! - \! \alpha_{p_{q}})^{2}} \! + \! \dfrac{(\mathcal{
N} \! - \! 2)}{\mathcal{N}} \left(\sum_{q=1}^{\mathfrak{s}-1} \dfrac{l_{q}}{x 
\! - \! \alpha_{p_{q}}} \right)^{2} \nonumber \\
+& \, \dfrac{(\mathcal{N} \! - \! 1)}{\mathcal{N}}n \widetilde{V}^{\prime}
(x) \sum_{q=1}^{\mathfrak{s}-1} \dfrac{l_{q}}{x \! - \! \alpha_{p_{q}}} \! 
- \! \dfrac{n}{\mathcal{N}}(\mathcal{N} \! - \! 1) \sum_{j=1}^{\mathcal{N}} 
\dfrac{\widetilde{V}^{\prime}(x) \! - \! \widetilde{V}^{\prime}
(\hat{x}_{j}^{\ast})}{x \! - \! \hat{x}_{j}^{\ast}}, \label{eql3.6y}
\end{align}
\begin{align}
\dfrac{n}{\mathcal{N}}(\mathcal{N} \! - \! 1) \sum_{j=1}^{\mathcal{N}} 
\dfrac{\widetilde{V}^{\prime}(x) \! - \! \widetilde{V}^{\prime}(\hat{x}_{
j}^{\ast})}{x \! - \! \hat{x}_{j}^{\ast}}=& \, \dfrac{n}{\mathcal{N}} 
\dfrac{(\mathcal{N} \! - \! 1)}{\hat{Z}(x)} \sum_{j=1}^{\mathcal{N}} \left(
\dfrac{\widetilde{V}^{\prime}(x) \hat{Z}(x) \! - \! \widetilde{V}^{\prime}
(\hat{x}_{j}^{\ast}) \hat{Z}(\hat{x}_{j}^{\ast})}{x \! - \! \hat{x}_{j}^{
\ast}} \right) \! - \! \dfrac{2}{\hat{Z}(x)} \left(\sum_{i=2(\hat{N}+1)}^{3}
\hat{c}_{2(\hat{N}+1)-i} \right. \nonumber \\
\times&\left. \, \left(\dfrac{\mathcal{N}(\mathcal{N} \! - \! 1)}{2}x^{i-2} \! 
+ \! \sum_{r=2}^{i-1} \left((\mathcal{N} \! - \! 1) \sum_{j=1}^{\mathcal{N}}
(\hat{x}_{j}^{\ast})^{r-1} \! + \! \sum_{j=1}^{\mathcal{N}-1} \sum_{k^{\prime}
=j+1}^{\mathcal{N}} \sum_{m=2}^{r-1}(\hat{x}_{j}^{\ast})^{r-m}
(\hat{x}_{k^{\prime}}^{\ast})^{m-1} \right) \right. \right. \nonumber \\
\times&\left. \left. x^{i-r-1} \right) \! + \! \dfrac{\mathcal{N}(\mathcal{N} 
\! - \! 1)}{2} \hat{c}_{2 \hat{N}} \right) \! + \! \dfrac{2}{\hat{Z}(x)} 
\dfrac{(\mathcal{N} \! - \! 1)}{\mathcal{N}} \left(\sum_{i=2(\hat{N}+1)}^{3} 
\hat{c}_{2(\hat{N}+1)-i} \right. \nonumber \\
\times&\left. \left(\sum_{r=0}^{i-1} \sum_{j=1}^{\mathcal{N}} \sum_{q=1}^{
\mathfrak{s}-1} \dfrac{l_{q}(\hat{x}_{j}^{\ast})^{r}}{\hat{x}_{j}^{\ast} 
\! - \! \alpha_{p_{q}}}x^{i-r-1} \right) \! + \! \hat{c}_{2 \hat{N}} 
\sum_{j=1}^{\mathcal{N}} \sum_{q=1}^{\mathfrak{s}-1} \dfrac{l_{q}(x \! + \! 
\hat{x}_{j}^{\ast})}{\hat{x}_{j}^{\ast} \! - \! \alpha_{p_{q}}} \right. 
\nonumber \\
+&\left. \hat{c}_{2 \hat{N}+1} \sum_{j=1}^{\mathcal{N}} \sum_{q=1}^{
\mathfrak{s}-1} \dfrac{l_{q}}{\hat{x}_{j}^{\ast} \! - \! \alpha_{p_{q}}} 
\right), \label{eql3.6z}
\end{align}
where $\hat{Z}(x) \! := \! \prod_{j=1}^{\hat{N}+1}(z \! - \! \hat{B}_{j-1})
(z \! - \! \hat{A}_{j}) \! := \! \sum_{j=2(\hat{N}+1)}^{0} \hat{c}_{2(\hat{N}+1)-j}
x^{j}$, with the $2(\hat{N} \! + \! 1)$ real points $\lbrace \hat{B}_{j-1},
\hat{A}_{j} \rbrace_{j=1}^{\hat{N}+1}$, for $\hat{N}$ $(= \! \hat{N}(n,k,z_{o}))$ 
$\in \! \mathbb{N}_{0}$, chosen so that $\overline{\mathbb{R}} \setminus 
\lbrace \alpha_{1},\alpha_{2},\dotsc,\alpha_{K} \rbrace \supset \cup_{j=1}^{
\hat{N}+1}[\hat{B}_{j-1},\hat{A}_{j}] \supseteq \lbrace \hat{x}_{1}^{\ast},
\hat{x}_{2}^{\ast},\dotsc,\hat{x}_{\mathcal{N}}^{\ast} \rbrace$, with 
enumeration $-\infty \! < \! \hat{B}_{0} \! < \! \hat{A}_{1} \! < \! \dotsb \! 
< \! \hat{B}_{j-1} \! < \! \hat{A}_{j} \! < \! \dotsb \! < \! \hat{B}_{\hat{N}} \! 
< \! \hat{A}_{\hat{N}+1} \! < \! +\infty$, and $[\hat{B}_{j-1},\hat{A}_{j}] \cap 
\lbrace \alpha_{1},\alpha_{2},\dotsc,\alpha_{K} \rbrace \! = \! \varnothing$, 
$j \! = \! 1,2,\dotsc,\hat{N} \! + \! 1$, and
\begin{equation*}
\hat{c}_{j} \! = \! \mathlarger{\sum_{\underset{\underset{\sum_{i=1}^{\hat{N}
+1}k^{\prime}_{i}+\sum_{m=1}^{\hat{N}+1}l^{\prime}_{m}=j}{i,m \in \lbrace 1,2,
\dotsc,\hat{N}+1 \rbrace}}{k^{\prime}_{i},l^{\prime}_{m}=0,1}}}(-1)^{j} 
\left(\prod_{i=1}^{\hat{N}+1} \binom{1}{k^{\prime}_{i}}(\hat{B}_{i-1})^{
k^{\prime}_{i}} \right) \left(\prod_{m=1}^{\hat{N}+1} \binom{1}{l^{\prime}_{m}}
(\hat{A}_{m})^{l^{\prime}_{m}} \right), \quad j \! = \! 0,1,\dotsc,2(\hat{N} 
\! + \! 1),
\end{equation*}
with $\hat{c}_{0} \! = \! 1$ and $\hat{c}_{2(\hat{N}+1)} \! = \! 
\prod_{r=1}^{\hat{N}+1}\hat{B}_{r-1} \hat{A}_{r}$,
\begin{align}
\hat{Q}(x) =& \, \dfrac{1}{\hat{Z}(x)} \left(\left(\sum_{i=2(\hat{N}+1)}^{0} 
\hat{c}_{2(\hat{N}+1)-i}x^{i} \right) \left(2 \dfrac{(\mathcal{N} \! - \! 
1)}{\mathcal{N}} \sum_{j=1}^{\mathcal{N}} \sum_{q=1}^{\mathfrak{s}-1} 
\dfrac{l_{q}}{(\hat{x}_{j}^{\ast} \! - \! \alpha_{p_{q}})(x \! - \! 
\alpha_{p_{q}})} \! + \! \sum_{q=1}^{\mathfrak{s}-1} \dfrac{l_{q}}{(x \! - \! 
\alpha_{p_{q}})^{2}} \right. \right. \nonumber \\
+&\left. \left. \, \dfrac{(\mathcal{N} \! - \! 2)}{\mathcal{N}} \left(\sum_{q
=1}^{\mathfrak{s}-1} \dfrac{l_{q}}{x \! - \! \alpha_{p_{q}}} \right)^{2} \! + 
\! \dfrac{(\mathcal{N} \! - \! 1)}{\mathcal{N}}n \widetilde{V}^{\prime}(x) 
\sum_{q=1}^{\mathfrak{s}-1} \dfrac{l_{q}}{x \! - \! \alpha_{p_{q}}} \right) \! 
- \! n \dfrac{(\mathcal{N} \! - \! 1)}{\mathcal{N}} \sum_{j=1}^{\mathcal{N}} 
\left(\dfrac{\widetilde{V}^{\prime}(x) \hat{Z}(x) \! - \! \widetilde{V}^{
\prime}(\hat{x}_{j}^{\ast}) \hat{Z}(\hat{x}_{j}^{\ast})}{x \! - \! 
\hat{x}_{j}^{\ast}} \right) \right. \nonumber \\
+&\left. \, 2 \left(\sum_{i=2(\hat{N}+1)}^{3} \hat{c}_{2(\hat{N}+1)-i} \left(
\dfrac{\mathcal{N}(\mathcal{N} \! - \! 1)}{2}x^{i-2} \! + \! \sum_{r=2}^{i-1} 
\left((\mathcal{N} \! - \! 1) \sum_{j=1}^{\mathcal{N}}(\hat{x}_{j}^{\ast})^{r
-1} \! + \! \sum_{j=1}^{\mathcal{N}-1} \sum_{k^{\prime}=j+1}^{\mathcal{N}} 
\sum_{m=2}^{r-1}(\hat{x}_{j}^{\ast})^{r-m}(\hat{x}_{k^{\prime}}^{\ast})^{m-1} 
\right) \right. \right. \right. \nonumber \\
\times&\left. \left. \left. \, x^{i-r-1} \right) \! + \! \dfrac{\mathcal{N}
(\mathcal{N} \! - \! 1)}{2} \hat{c}_{2 \hat{N}} \right) \! - \! 2 \dfrac{
(\mathcal{N} \! - \! 1)}{\mathcal{N}} \left(\sum_{i=2(\hat{N}+1)}^{3} 
\hat{c}_{2(\hat{N}+1)-i} \left(\sum_{r=0}^{i-1} \sum_{j=1}^{\mathcal{N}} 
\sum_{q=1}^{\mathfrak{s}-1} \dfrac{l_{q}(\hat{x}_{j}^{\ast})^{r}}{
\hat{x}_{j}^{\ast} \! - \! \alpha_{p_{q}}}x^{i-r-1} \right) \right. \right. 
\nonumber \\
+&\left. \left. \, \hat{c}_{2 \hat{N}} \sum_{j=1}^{\mathcal{N}} \sum_{q=1}^{
\mathfrak{s}-1} \dfrac{l_{q}(x \! + \! \hat{x}_{j}^{\ast})}{\hat{x}_{j}^{\ast} 
\! - \! \alpha_{p_{q}}} \! + \! \hat{c}_{2 \hat{N}+1} \sum_{j=1}^{\mathcal{N}} 
\sum_{q=1}^{\mathfrak{s}-1} \dfrac{l_{q}}{\hat{x}_{j}^{\ast} \! - \! 
\alpha_{p_{q}}} \right) \right), \label{eql3.6aa}
\end{align}
\begin{equation} \label{eql3.6bb} 
\hat{f}(x) \! = \! \hat{F}(x) \exp \left(-\dfrac{1}{2} \int_{}^{x} 
\hat{\mathcal{P}}(\xi) \, \md \xi \right),
\end{equation}
where
\begin{equation} \label{eql3.6cc}
\hat{\mathcal{P}}(x) \! := \! 2 \dfrac{n}{\mathcal{N}} \left(-\dfrac{1}{2}
(\mathcal{N} \! - \! 1) \widetilde{V}^{\prime}(x) \! + \! 
\sum_{q=1}^{\mathfrak{s}-1} \dfrac{l_{q}/n}{x \! - \! \alpha_{p_{q}}} 
\right),
\end{equation}
\begin{equation} \label{eql3.6dd} 
\hat{F}^{\prime \prime}(x) \! = \! \left(\hat{Q}(x) \! + \! \dfrac{1}{2} 
\hat{\mathcal{P}}^{\prime}(x) \! + \! \dfrac{1}{4}(\hat{\mathcal{P}}(x))^{2} 
\right) \hat{F}(x),
\end{equation}
with
\begin{align}
& \hat{Q}(x) \! + \! \dfrac{1}{2} \hat{\mathcal{P}}^{\prime}(x) \! + \! 
\dfrac{1}{4}(\hat{\mathcal{P}}(x))^{2} \! = \! -\dfrac{1}{2} 
\dfrac{(\mathcal{N} \! - \! 1)}{\mathcal{N}}n \widetilde{V}^{\prime \prime}(x) 
\! + \! \dfrac{(\mathcal{N} \! - \! 1)}{\mathcal{N}} \sum_{q=1}^{\mathfrak{s}
-1} \dfrac{l_{q}}{(x \! - \! \alpha_{p_{q}})^{2}} \! + \! \dfrac{1}{4} 
\left(\dfrac{\mathcal{N} \! - \! 1}{\mathcal{N}} \right)^{2}(n 
\widetilde{V}^{\prime}(x))^{2} \nonumber \\
+& \, n \widetilde{V}^{\prime}(x) \left(\dfrac{\mathcal{N} \! - \! 1}{
\mathcal{N}} \right)^{2} \sum_{q=1}^{\mathfrak{s}-1} \dfrac{l_{q}}{x \! - 
\! \alpha_{p_{q}}} \! + \! \left(\dfrac{\mathcal{N} \! - \! 1}{\mathcal{N}} 
\right)^{2} \left(\sum_{q=1}^{\mathfrak{s}-1} \dfrac{l_{q}}{x \! - \! 
\alpha_{p_{q}}} \right)^{2} \! + \! 2 \dfrac{(\mathcal{N} \! - \! 1)}{
\mathcal{N}} \sum_{j=1}^{\mathcal{N}} \sum_{q=1}^{\mathfrak{s}-1} 
\dfrac{l_{q}}{(\hat{x}_{j}^{\ast} \! - \! \alpha_{p_{q}})(x \! - \! 
\alpha_{p_{q}})} \nonumber \\
+& \, \dfrac{1}{\hat{Z}(x)} \left(-n \dfrac{(\mathcal{N} \! - \! 1)}{
\mathcal{N}} \sum_{j=1}^{\mathcal{N}} \left(\dfrac{\widetilde{V}^{\prime}(x) 
\hat{Z}(x) \! - \! \widetilde{V}^{\prime}(\hat{x}_{j}^{\ast}) \hat{Z}
(\hat{x}_{j}^{\ast})}{x \! - \! \hat{x}_{j}^{\ast}} \right) \! + \! 2 \left(
\sum_{i=2(\hat{N}+1)}^{3} \hat{c}_{2(\hat{N}+1)-i} \left(\dfrac{\mathcal{N}
(\mathcal{N} \! - \! 1)}{2}x^{i-2} \! + \! \sum_{r=2}^{i-1} \left(
\vphantom{M^{M^{M^{M}}}}(\mathcal{N} \! - \! 1) \right. \right. \right. 
\right. \nonumber \\
\times&\left. \left. \left. \left. \, \sum_{j=1}^{\mathcal{N}}(\hat{x}_{j}^{
\ast})^{r-1} \! + \! \sum_{j=1}^{\mathcal{N}-1} \sum_{k^{\prime}=j+1}^{
\mathcal{N}} \sum_{m=2}^{r-1}(\hat{x}_{j}^{\ast})^{r-m}(\hat{x}_{k^{\prime}}^{
\ast})^{m-1} \right)x^{i-r-1} \right) \! + \! \dfrac{\mathcal{N}(\mathcal{N} 
\! - \! 1)}{2} \hat{c}_{2 \hat{N}} \right) \! - \! 2 \dfrac{(\mathcal{N} \! - 
\! 1)}{\mathcal{N}} \left(\sum_{i=2(\hat{N}+1)}^{3} \hat{c}_{2(\hat{N}+1)-i} 
\right. \right. \nonumber \\
\times&\left. \left. \, \left(\sum_{r=0}^{i-1} \sum_{j=1}^{\mathcal{N}} 
\sum_{q=1}^{\mathfrak{s}-1} \dfrac{l_{q}(\hat{x}_{j}^{\ast})^{r}}{\hat{x}_{
j}^{\ast} \! - \! \alpha_{p_{q}}}x^{i-r-1} \right) \! + \! \hat{c}_{2 \hat{N}} 
\sum_{j=1}^{\mathcal{N}} \sum_{q=1}^{\mathfrak{s}-1} \dfrac{l_{q}(x \! + \! 
\hat{x}_{j}^{\ast})}{\hat{x}_{j}^{\ast} \! - \! \alpha_{p_{q}}} \! + \! 
\hat{c}_{2 \hat{N}+1} \sum_{j=1}^{\mathcal{N}} \sum_{q=1}^{\mathfrak{s}-1} 
\dfrac{l_{q}}{\hat{x}_{j}^{\ast} \! - \! \alpha_{p_{q}}} \right) \right), 
\label{eql3.6ee}
\end{align}
\begin{equation} \label{eql3.6ff} 
\left\vert \hat{Q}(x) \! + \! \dfrac{1}{2} \hat{\mathcal{P}}^{\prime}(x) 
\! + \! \dfrac{1}{4}(\hat{\mathcal{P}}(x))^{2} \right\vert \! \leqslant 
\! \dfrac{\hat{\mathfrak{C}}_{r} \mathcal{N}^{2}}{(x \! - \! \hat{B}_{r-1})
(\hat{A}_{r} \! - \! x)}, \quad x \! \in \! [\hat{B}_{r-1},\hat{A}_{r}], 
\quad r \! = \! 1,2,\dotsc,\hat{N} \! + \! 1, 
\end{equation}
where $\hat{\mathfrak{C}}_{r} \! = \! \hat{\mathfrak{C}}_{r}(n,k,
z_{o}) \! > \! 0$ and $\mathcal{O}(1)$ (in the double-scaling limit 
$\mathscr{N},n \! \to \! \infty$ such that $z_{o} \! = \! 1 \! + \! o(1))$,
\begin{equation} \label{eql3.6gg} 
\hat{x}_{j+1}^{\ast} \! - \! \hat{x}_{j}^{\ast} \! \geqslant \! \dfrac{\min 
\lbrace (3 \hat{\mathfrak{C}}_{r})^{-1/2},1/4 \rbrace}{((n \! - \! 1)K 
\! + \! k)} \left((\hat{A}_{r} \! - \! \hat{x}_{j}^{\ast})(\hat{A}_{r} \! 
- \! \hat{x}_{j+1}^{\ast}) \right)^{1/2}, \quad j \! = \! 1,2,\dotsc,
\mathcal{N} \! - \! 1, \quad r \! = \! 1,2,\dotsc,\hat{N} \! + \! 1,
\end{equation}
with $((\hat{A}_{r} \! - \! \hat{x}_{j}^{\ast})(\hat{A}_{r} \! - \! 
\hat{x}_{j+1}^{\ast}))^{1/2} \! > \! 0$ and $\mathcal{O}(1)$ (in the 
double-scaling limit $\mathscr{N},n \! \to \! \infty$ such that $z_{o} 
\! = \! 1 \! + \! o(1))$. This concludes the proof. \hfill $\qed$

The following Lemma~\ref{lem3.7} provides more detailed 
information about the representation and support of the associated 
family of $\hat{K} \! := \! \# \lbrace \mathstrut j \! \in \! \lbrace 
1,2,\dotsc,K \rbrace; \, \alpha_{j} \! = \! \infty \rbrace$ (resp., 
$\tilde{K} \! := \! \# \lbrace \mathstrut j \! \in \! \lbrace 1,2,
\dotsc,K \rbrace; \, \alpha_{j} \! \neq \! \infty \rbrace)$ 
equilibrium measures $\mu^{\infty}_{\widetilde{V}}$ (resp., 
$\mu^{f}_{\widetilde{V}})$.
\begin{ccccc} \label{lem3.7} 
Let the external field $\widetilde{V} \colon \overline{\mathbb{R}} 
\setminus \lbrace \alpha_{1},\alpha_{2},\dotsc,\alpha_{K} \rbrace 
\! \to \! \mathbb{R}$ satisfy conditions~\eqref{eq20}--\eqref{eq22}$;$ 
suppose, furthermore, that $\widetilde{V}$ is regular.

$\pmb{(1)}$ For $n \! \in \! \mathbb{N}$ and $k \! \in \! \lbrace 
1,2,\dotsc,K \rbrace$ such that $\alpha_{p_{\mathfrak{s}}} \! := 
\! \alpha_{k} \! = \! \infty$, with associated equilibrium measure 
$\mu_{\widetilde{V}}^{\infty}$ $(\in \! \mathscr{M}_{1}(\mathbb{R}))$, 
set $J_{\infty} \! := \! \supp (\mu_{\widetilde{V}}^{\infty})$, which, 
{}from Lemma~\ref{lem3.1}, is a proper compact subset of 
$\overline{\mathbb{R}} \setminus \lbrace \alpha_{1},\alpha_{2},\dotsc,
\alpha_{K} \rbrace$. Then, for $n \! \in \! \mathbb{N}$ and $k \! \in \! 
\lbrace 1,2,\dotsc,K \rbrace$ such that $\alpha_{p_{\mathfrak{s}}} \! 
:= \! \alpha_{k} \! = \! \infty$$:$ $\pmb{{\rm (i)}}$ $J_{\infty} \! = \! 
\cup_{j=1}^{N+1}[\hat{b}_{j-1},\hat{a}_{j}]$,\footnote{Strictly speaking, 
the end-points of the support of the associated equilibrium measure, 
$\mu_{\widetilde{V}}^{\infty}$, depend on $n$, $k$, and $z_{o}$, that is, 
$\hat{b}_{j-1} \! = \! \hat{b}_{j-1}(n,k,z_{o})$ and $\hat{a}_{j} \! = \! 
\hat{a}_{j}(n,k,z_{o})$, $j \! = \! 1,2,\dotsc,N \! + \! 1$; but, for notational 
simplicity, unless where absolutely necessary, these additional $n$-, $k$-, 
and $z_{o}$-dependencies are suppressed. (There will be no discussion 
concerning the parametric dependence of the associated end-points 
$\hat{b}_{j-1},\hat{a}_{j}$, $j \! = \! 1,2,\dotsc,N \! + \! 1$, on the elements 
of the corresponding pole set $\lbrace \alpha_{p_{1}},\alpha_{p_{2}},\dotsc,
\alpha_{p_{\mathfrak{s}}} \rbrace$.)} with $N \! \in \! \mathbb{N}_{0}$ and 
finite, $[\hat{b}_{j-1},\hat{a}_{j}] \cap \lbrace \alpha_{p_{1}},\alpha_{p_{2}},
\dotsc,\alpha_{p_{\mathfrak{s}}} \rbrace \! = \! \varnothing$, $j \! = \! 1,2,
\dotsc,N \! + \! 1$, $[\hat{b}_{i-1},\hat{a}_{i}] \cap [\hat{b}_{j-1},\hat{a}_{j}] 
\! = \! \varnothing$, $i \! \neq \! j \! \in \! \lbrace 1,2,\dotsc,N \! + \! 1 
\rbrace$, and $-\infty \! < \! \hat{b}_{0} \! < \! \hat{a}_{1} \! < \! \hat{b}_{1} 
\! < \! \hat{a}_{2} \! < \! \dotsb \! < \! \hat{b}_{N} \! < \! \hat{a}_{N+1} \! 
< \! +\infty$, where $\{\hat{b}_{j-1},\hat{a}_{j}\}_{j=1}^{N+1}$ satisfy the 
associated, locally solvable system of $2(N \! + \! 1)$ real moment equations
\begin{align} 
&\int_{J_{\infty}} \dfrac{\xi^{j}}{(\hat{R}(\xi))^{1/2}_{+}} \left(
\dfrac{2}{\mi \pi} \sum_{q=1}^{\mathfrak{s}-1} \dfrac{\varkappa_{nk 
\tilde{k}_{q}}}{n(\xi \! - \! \alpha_{p_{q}})} \! + \! \dfrac{\widetilde{V}^{
\prime}(\xi)}{\mi \pi} \right) \md \xi \! = \! 0, \quad j \! = \! 0,1,
\dotsc,N, \label{eql3.7a} \\
&\int_{J_{\infty}} \dfrac{\xi^{N+1}}{(\hat{R}(\xi))^{1/2}_{+}} \left(
\dfrac{2}{\mi \pi} \sum_{q=1}^{\mathfrak{s}-1} \dfrac{\varkappa_{nk 
\tilde{k}_{q}}}{n(\xi \! - \! \alpha_{p_{q}})} \! + \! \dfrac{\widetilde{V}^{
\prime}(\xi)}{\mi \pi} \right) \md \xi \! = \! -2 \left(\dfrac{(n \! - \! 
1)K \! + \! k}{n} \right), \label{eql3.7b} \\
&\int_{\hat{a}_{j}}^{\hat{b}_{j}}(\hat{R}(\varsigma))^{1/2} \left(
\int_{J_{\infty}}(\hat{R}(\xi))^{-1/2}_{+} \left(\dfrac{2}{\mi \pi} 
\sum_{q=1}^{\mathfrak{s}-1} \dfrac{\varkappa_{nk \tilde{k}_{q}}}{n(\xi 
\! - \! \alpha_{p_{q}})} \! + \! \dfrac{\widetilde{V}^{\prime}(\xi)}{\mi \pi} 
\right) \dfrac{\md \xi}{\xi \! - \! \varsigma} \right) \md \varsigma 
\nonumber \\
&= \, 2 \sum_{q=1}^{\mathfrak{s}-1} \dfrac{\varkappa_{nk \tilde{k}_{q}}}{n} 
\ln \left\vert \dfrac{\hat{b}_{j} \! - \! \alpha_{p_{q}}}{\hat{a}_{j} \! - \! 
\alpha_{p_{q}}} \right\vert \! + \! \widetilde{V}(\hat{b}_{j}) \! - \! 
\widetilde{V}(\hat{a}_{j}), \quad j \! = \! 1,2,\dotsc,N, \label{eql3.7c}
\end{align}
where
\begin{equation} \label{eql3.7d} 
(\hat{R}(z))^{1/2} \! := \! \left(\prod_{j=1}^{N+1}(z \! - \! \hat{b}_{j-1})
(z \! - \! \hat{a}_{j}) \right)^{1/2},
\end{equation}
with $(\hat{R}(z))^{1/2}_{\pm} \! := \! \lim_{\varepsilon \downarrow 0}
(\hat{R}(z \! \pm \! \mi \varepsilon))^{1/2}$, and the branch of the 
square root is chosen so that $z^{-(N+1)}(\hat{R}(z))^{1/2} \! 
\sim_{\overline{\mathbb{C}}_{\pm} \ni z \to \alpha_{k}} \! \pm 1$$;$ 
and $\pmb{{\rm (ii)}}$ the density of the associated equilibrium measure, 
$\mu_{\widetilde{V}}^{\infty}$, which, via Lemma~\ref{lem3.6}, is 
absolutely continuous with respect to Lebesgue measure, is given by
\begin{equation} \label{eql3.7e} 
\md \mu_{\widetilde{V}}^{\infty}(x) \! = \! \psi_{\widetilde{V}}^{\infty}
(x) \, \md x \! = \! (2 \pi \mi)^{-1}(\hat{R}(x))^{1/2}_{+} 
\hat{h}_{\widetilde{V}}(x) \chi_{J_{\infty}}(x) \, \md x,
\end{equation}
where
\begin{equation} \label{eql3.7f} 
\hat{h}_{\widetilde{V}}(z) \! := \! \dfrac{1}{2} \left(\dfrac{(n \! - \! 1)K 
\! + \! k}{n} \right)^{-1} \oint_{\hat{C}_{\widetilde{V}}} \left(\dfrac{2 
\mi}{\pi} \sum_{q=1}^{\mathfrak{s}-1} \dfrac{\varkappa_{nk 
\tilde{k}_{q}}}{n(\xi \! - \! \alpha_{p_{q}})} \! + \! \dfrac{\mi 
\widetilde{V}^{\prime}(\xi)}{\pi} \right) \dfrac{(\hat{R}(\xi))^{-
1/2}}{\xi \! - \! z} \, \md \xi
\end{equation}
(real analytic for $z \! \in \! \overline{\mathbb{R}} \setminus \lbrace 
\alpha_{p_{1}},\alpha_{p_{2}},\dotsc,\alpha_{p_{\mathfrak{s}}} 
\rbrace)$, with $\overline{\mathbb{C}} \setminus \lbrace \alpha_{p_{1}},
\alpha_{p_{2}},\dotsc,\alpha_{p_{\mathfrak{s}}} \rbrace \supset 
\hat{C}_{\widetilde{V}}$ the simple boundary of any open 
$\mathfrak{s}$-connected punctured region of the type 
$\hat{C}_{\widetilde{V}} \! = \! \partial \, \mathbb{U}_{R_{1}}^{
\circlearrowright} \cup \cup_{q=1}^{\mathfrak{s}-1} \partial \, 
\mathbb{U}^{\circlearrowleft}_{\hat{\delta}_{q}}(\alpha_{p_{q}})$, 
where $\mathbb{U}_{R_{1}}^{\circlearrowright} \! = \! \lbrace \mathstrut 
z \! \in \! \mathbb{C}; \, \lvert z \rvert \! < \! R_{1} \rbrace$, with 
$\partial \, \mathbb{U}_{R_{1}}^{\circlearrowright} \! = \! \lbrace 
\mathstrut z \! \in \! \mathbb{C}; \, \lvert z \rvert \! = \! R_{1} 
\rbrace$ oriented clockwise, and, for $q \! = \! 1,2,\dotsc,\mathfrak{s} 
\! - \! 1$, $\mathbb{U}_{\hat{\delta}_{q}}^{\circlearrowleft}(\alpha_{p_{q}}) 
\! = \! \lbrace \mathstrut z \! \in \! \mathbb{C}; \, \lvert z \! - \! 
\alpha_{p_{q}} \rvert \! < \! \hat{\delta}_{q} \rbrace$, with $\partial 
\, \mathbb{U}_{\hat{\delta}_{q}}^{\circlearrowleft}(\alpha_{p_{q}}) \! = \! 
\lbrace \mathstrut z \! \in \! \mathbb{C}; \, \lvert z \! - \! \alpha_{p_{q}} 
\rvert \! = \! \hat{\delta}_{q} \rbrace$ oriented counter-clockwise, where 
the numbers $0 \! < \! \hat{\delta}_{q} \! < \! R_{1} \! < \! +\infty$, 
$q \! = \! 1,2,\dotsc,\mathfrak{s} \! - \! 1$, are chosen so that, for any 
non-real $z$ in the domain of analyticity of $\widetilde{V}$, $\partial 
\, \mathbb{U}_{\hat{\delta}_{i}}^{\circlearrowleft}(\alpha_{p_{i}}) 
\cap \partial \, \mathbb{U}_{\hat{\delta}_{j}}^{\circlearrowleft}
(\alpha_{p_{j}}) \! = \! \varnothing$, $i \! \neq \! j \! \in \! \lbrace 
1,2,\dotsc,\mathfrak{s} \! - \! 1 \rbrace$, $\partial \, \mathbb{U}_{
\hat{\delta}_{j}}^{\circlearrowleft}(\alpha_{p_{j}}) \cap \partial \, 
\mathbb{U}_{R_{1}}^{\circlearrowright} \! = \! \varnothing$, $j \! = 
\! 1,2,\dotsc,\mathfrak{s} \! - \! 1$, $(\mathbb{U}_{\hat{\delta}_{j}}^{
\circlearrowleft}(\alpha_{p_{j}}) \cup \partial \, \mathbb{U}_{\hat{
\delta}_{j}}^{\circlearrowleft}(\alpha_{p_{j}})) \cap J_{\infty} \! = \! 
\varnothing$, $j \! = \! 1,2,\dotsc,\mathfrak{s} \! - \! 1$, $\partial 
\, \mathbb{U}_{R_{1}}^{\circlearrowright} \cap J_{\infty} \! = \! 
\varnothing$, and $(\mathbb{U}_{R_{1}}^{\circlearrowright} \setminus 
\cup_{q=1}^{\mathfrak{s}-1}(\mathbb{U}_{\hat{\delta}_{q}}^{\circlearrowleft}
(\alpha_{p_{q}}) \cup \partial \, \mathbb{U}_{\hat{\delta}_{q}}^{
\circlearrowleft}(\alpha_{p_{q}})) \! =:)$ $\operatorname{int}
(\hat{C}_{\widetilde{V}}) \supset \{z\} \cup J_{\infty}$, $\chi_{J_{\infty}}
(x)$ is the characteristic function of the set $J_{\infty}$, and 
$\psi_{\widetilde{V}}^{\infty}(x) \! \geqslant \! 0$ (resp., 
$\psi_{\widetilde{V}}^{\infty}(x) \! > \! 0)$ for $x \! \in \! J_{\infty}$ 
(resp., $x \! \in \! \operatorname{int}(J_{\infty}))$.

$\pmb{(2)}$ For $n \! \in \! \mathbb{N}$ and $k \! \in \! \lbrace 
1,2,\dotsc,K \rbrace$ such that $\alpha_{p_{\mathfrak{s}}} \! := \! 
\alpha_{k} \! \neq \! \infty$, with associated equilibrium measure 
$\mu_{\widetilde{V}}^{f}$ $(\in \! \mathscr{M}_{1}(\mathbb{R}))$, 
set $J_{f} \! := \! \supp (\mu_{\widetilde{V}}^{f})$, which, 
{}from Lemma~\ref{lem3.1}, is a proper compact subset of 
$\overline{\mathbb{R}} \setminus \lbrace \alpha_{1},\alpha_{2},\dotsc,
\alpha_{K} \rbrace$. Then, for $n \! \in \! \mathbb{N}$ and $k \! \in 
\! \lbrace 1,2,\dotsc,K \rbrace$ such that $\alpha_{p_{\mathfrak{s}}} 
\! := \! \alpha_{k} \! \neq \! \infty$$:$ $\pmb{{\rm (i)}}$ $J_{f} \! = \! 
\cup_{j=1}^{N+1}[\tilde{b}_{j-1},\tilde{a}_{j}]$,\footnote{Strictly speaking, 
the end-points of the support of the associated equilibrium measure, 
$\mu_{\widetilde{V}}^{f}$, depend on $n$, $k$, and $z_{o}$, that is, 
$\tilde{b}_{j-1} \! = \! \tilde{b}_{j-1}(n,k,z_{o})$ and $\tilde{a}_{j} \! = \! 
\tilde{a}_{j}(n,k,z_{o})$, $j \! = \! 1,2,\dotsc,N \! + \! 1$; but, for notational 
simplicity, unless where absolutely necessary, these additional $n$-, $k$-, 
and $z_{o}$-dependencies are suppressed. (There will be no discussion 
concerning the parametric dependence of the associated end-points 
$\tilde{b}_{j-1},\tilde{a}_{j}$, $j \! = \! 1,2,\dotsc,N \! + \! 1$, on the 
elements of the corresponding pole set $\lbrace \alpha_{p_{1}},\alpha_{p_{2}},
\dotsc,\alpha_{p_{\mathfrak{s}}} \rbrace$.)} with $[\tilde{b}_{j-1},\tilde{a}_{j}] 
\cap \lbrace \alpha_{p_{1}},\alpha_{p_{2}},\dotsc,\alpha_{p_{\mathfrak{s}}} 
\rbrace \! = \! \varnothing$, $j \! = \! 1,2,\dotsc,N \! + \! 1$, $[\tilde{b}_{i-1},
\tilde{a}_{i}] \cap [\tilde{b}_{j-1},\tilde{a}_{j}] \! = \! \varnothing$, $i \! \neq 
\! j \! \in \! \lbrace 1,2,\dotsc,N \! + \! 1 \rbrace$, and $-\infty \! < \! 
\tilde{b}_{0} \! < \! \tilde{a}_{1} \! < \! \tilde{b}_{1} \! < \! \tilde{a}_{2} \! 
< \! \dotsb \! < \! \tilde{b}_{N} \! < \! \tilde{a}_{N+1} \! < \! +\infty$, 
where $\{\tilde{b}_{j-1},\tilde{a}_{j}\}_{j=1}^{N+1}$ satisfy the associated, 
locally solvable system of $2(N \! + \! 1)$ real moment equations
\begin{align}
&\int_{J_{f}} \dfrac{\xi^{j}}{(\tilde{R}(\xi))^{1/2}_{+}} \left(
\dfrac{2}{\mi \pi} \left(\dfrac{(\varkappa_{nk} \! - \! 1)}{n(\xi \! - \! 
\alpha_{k})} \! + \! \sum_{q=1}^{\mathfrak{s}-2} \dfrac{\varkappa_{nk 
\tilde{k}_{q}}}{n(\xi \! - \! \alpha_{p_{q}})} \right) \! + \! \dfrac{
\widetilde{V}^{\prime}(\xi)}{\mi \pi} \right) \md \xi \! = \! 0, 
\quad j \! = \! 0,1,\dotsc,N, \label{eql3.7g} \\
&\int_{J_{f}} \dfrac{\xi^{N+1}}{(\tilde{R}(\xi))^{1/2}_{+}} \left(
\dfrac{2}{\mi \pi} \left(\dfrac{(\varkappa_{nk} \! - \! 1)}{n(\xi \! - \! 
\alpha_{k})} \! + \! \sum_{q=1}^{\mathfrak{s}-2} \dfrac{\varkappa_{nk 
\tilde{k}_{q}}}{n(\xi \! - \! \alpha_{p_{q}})} \right) \! + \! \dfrac{
\widetilde{V}^{\prime}(\xi)}{\mi \pi} \right) \md \xi \! = \! -2 
\left(\dfrac{(n \! - \! 1)K \! + \! k}{n} \right), \label{eql3.7h} \\
&\int_{\tilde{a}_{j}}^{\tilde{b}_{j}}(\tilde{R}(\varsigma))^{1/2} \left(
\int_{J_{f}}(\tilde{R}(\xi))^{-1/2}_{+} \left(\dfrac{2}{\mi \pi} \left(
\dfrac{(\varkappa_{nk} \! - \! 1)}{n(\xi \! - \! \alpha_{k})} \! + \! 
\sum_{q=1}^{\mathfrak{s}-2} \dfrac{\varkappa_{nk \tilde{k}_{q}}}{n(\xi \! 
- \! \alpha_{p_{q}})} \right) \! + \! \dfrac{\widetilde{V}^{\prime}(\xi)}{\mi 
\pi} \right) \dfrac{\md \xi}{\xi \! - \! \varsigma} \right) \md \varsigma 
\nonumber \\
&= \, 2 \left(\dfrac{\varkappa_{nk} \! - \! 1}{n} \right) \ln \left\vert 
\dfrac{\tilde{b}_{j} \! - \! \alpha_{k}}{\tilde{a}_{j} \! - \! \alpha_{k}} 
\right\vert \! + \! 2 \sum_{q=1}^{\mathfrak{s}-2} \dfrac{\varkappa_{nk 
\tilde{k}_{q}}}{n} \ln \left\vert \dfrac{\tilde{b}_{j} \! - \! 
\alpha_{p_{q}}}{\tilde{a}_{j} \! - \! \alpha_{p_{q}}} \right\vert \! + \! 
\widetilde{V}(\tilde{b}_{j}) \! - \! \widetilde{V}(\tilde{a}_{j}), \quad 
j \! = \! 1,2,\dotsc,N, \label{eql3.7i}
\end{align}
where
\begin{equation} \label{eql3.7j} 
(\tilde{R}(z))^{1/2} \! := \! \left(\prod_{j=1}^{N+1}(z \! - \! 
\tilde{b}_{j-1})(z \! - \! \tilde{a}_{j}) \right)^{1/2},
\end{equation}
with $(\tilde{R}(z))^{1/2}_{\pm} \! := \! \lim_{\varepsilon \downarrow 0}
(\tilde{R}(z \! \pm \! \mi \varepsilon))^{1/2}$, and the branch of the 
square root is chosen so that $z^{-(N+1)}(\tilde{R}(z))^{1/2} \! 
\sim_{\overline{\mathbb{C}}_{\pm} \ni z \to \alpha_{p_{\mathfrak{s}-1}}
= \infty} \! \pm 1$$;$ and $\pmb{{\rm (ii)}}$ the density of the 
associated equilibrium measure, $\mu^{f}_{\widetilde{V}}$, which, via 
Lemma~\ref{lem3.6}, is absolutely continuous with respect to Lebesgue 
measure, is given by
\begin{equation} \label{eql3.7k} 
\md \mu_{\widetilde{V}}^{f}(x) \! = \! \psi_{\widetilde{V}}^{f}(x) \, \md x 
\! = \! (2 \pi \mi)^{-1}(\tilde{R}(x))^{1/2}_{+} \tilde{h}_{\widetilde{V}}(x) 
\chi_{J_{f}}(x) \, \md x,
\end{equation}
where
\begin{align} \label{eql3.7l} 
\tilde{h}_{\widetilde{V}}(z) :=& \, \dfrac{1}{2} \left(\dfrac{(n \! - \! 1)K 
\! + \! k}{n} \right)^{-1} \oint_{\tilde{C}_{\widetilde{V}}} \left(
\dfrac{2 \mi}{\pi} \left(\dfrac{(\varkappa_{nk} \! - \! 1)}{n(\xi \! - \! 
\alpha_{k})} \! + \! \sum_{q=1}^{\mathfrak{s}-2} \dfrac{\varkappa_{nk 
\tilde{k}_{q}}}{n(\xi \! - \! \alpha_{p_{q}})} \right) \! + \! \dfrac{\mi 
\widetilde{V}^{\prime}(\xi)}{\pi} \right) \dfrac{(\tilde{R}(\xi))^{-1/2}}{\xi 
\! - \! z} \, \md \xi
\end{align}
(real analytic for $z \! \in \! \overline{\mathbb{R}} \setminus \lbrace 
\alpha_{p_{1}},\alpha_{p_{2}},\dotsc,\alpha_{p_{\mathfrak{s}}} \rbrace)$, 
with $\overline{\mathbb{C}} \setminus \lbrace \alpha_{p_{1}},
\alpha_{p_{2}},\dotsc,\alpha_{p_{\mathfrak{s}}} \rbrace \supset 
\tilde{C}_{\widetilde{V}}$ the simple boundary of any open 
$\mathfrak{s}$-connected punctured region of the type 
$\tilde{C}_{\widetilde{V}} \! = \! \partial \, \mathbb{U}_{R_{2}}^{
\circlearrowright} \cup \cup_{\underset{q \neq \mathfrak{s}-1}{q=1}}^{
\mathfrak{s}} \partial \, \mathbb{U}^{\circlearrowleft}_{\tilde{\delta}_{q}}
(\alpha_{p_{q}})$, where $\mathbb{U}_{R_{2}}^{\circlearrowright} \! = \! 
\lbrace \mathstrut z \! \in \! \mathbb{C}; \, \lvert z \rvert \! < \! R_{2} 
\rbrace$, with $\partial \, \mathbb{U}_{R_{2}}^{\circlearrowright} \! = \! 
\lbrace \mathstrut z \! \in \! \mathbb{C}; \, \lvert z \rvert \! = \! R_{2} 
\rbrace$ oriented clockwise, and, for $q \! = \! 1,\dotsc,\mathfrak{s} 
\! - \! 2,\mathfrak{s}$, $\mathbb{U}_{\tilde{\delta}_{q}}^{\circlearrowleft}
(\alpha_{p_{q}}) \! = \! \lbrace \mathstrut z \! \in \! \mathbb{C}; \, \lvert 
z \! - \! \alpha_{p_{q}} \rvert \! < \! \tilde{\delta}_{q} \rbrace$, with 
$\partial \, \mathbb{U}_{\tilde{\delta}_{q}}^{\circlearrowleft}
(\alpha_{p_{q}}) \! = \! \lbrace \mathstrut z \! \in \! \mathbb{C}; \, \lvert 
z \! - \! \alpha_{p_{q}} \rvert \! = \! \tilde{\delta}_{q} \rbrace$ oriented 
counter-clockwise, where the numbers $0 \! < \! \tilde{\delta}_{q} \! 
< \! R_{2} \! < \! +\infty$, $q \! = \! 1,\dotsc,\mathfrak{s} \! - \! 2,
\mathfrak{s}$, are chosen so that, for any non-real $z$ in the domain of 
analyticity of $\widetilde{V}$, $\partial \, \mathbb{U}_{\tilde{\delta}_{
i}}^{\circlearrowleft}(\alpha_{p_{i}}) \cap \partial \, \mathbb{U}_{\tilde{
\delta}_{j}}^{\circlearrowleft}(\alpha_{p_{j}}) \! = \! \varnothing$, $i \! 
\neq \! j \! \in \! \lbrace 1,\dotsc,\mathfrak{s} \! - \! 2,\mathfrak{s} 
\rbrace$, $\partial \, \mathbb{U}_{\tilde{\delta}_{j}}^{\circlearrowleft}
(\alpha_{p_{j}}) \cap \partial \, \mathbb{U}_{R_{2}}^{\circlearrowright} 
\! = \! \varnothing$, $j \! = \! 1,\dotsc,\mathfrak{s} \! - \! 2,
\mathfrak{s}$, $(\mathbb{U}_{\tilde{\delta}_{j}}^{\circlearrowleft}
(\alpha_{p_{j}}) \cup \partial \, \mathbb{U}_{\tilde{\delta}_{j}}^{
\circlearrowleft}(\alpha_{p_{j}})) \cap J_{f} \! = \! \varnothing$, 
$j \! = \! 1,\dotsc,\mathfrak{s} \! - \! 2,\mathfrak{s}$, $\partial \, 
\mathbb{U}_{R_{2}}^{\circlearrowright} \cap J_{f} \! = \! \varnothing$, and 
$(\mathbb{U}_{R_{2}}^{\circlearrowright} \setminus \cup_{\underset{q \neq 
\mathfrak{s}-1}{q=1}}^{\mathfrak{s}}(\mathbb{U}_{\tilde{\delta}_{q}}^{
\circlearrowleft}(\alpha_{p_{q}}) \cup \partial \, \mathbb{U}_{\tilde{\delta}_{
q}}^{\circlearrowleft}(\alpha_{p_{q}})) \! =:)$ $\operatorname{int}
(\tilde{C}_{\widetilde{V}}) \supset \{z\} \cup J_{f}$, $\chi_{J_{f}}(x)$ is the 
characteristic function of the set $J_{f}$, and $\psi_{\widetilde{V}}^{f}(x) 
\! \geqslant \! 0$ (resp., $\psi_{\widetilde{V}}^{f}(x) \! > \! 0)$ for 
$x \! \in \! J_{f}$ (resp., $x \! \in \! \operatorname{int}(J_{f}))$.
\end{ccccc}

\emph{Proof}. The proof of this Lemma~\ref{lem3.7} consists 
of two---rather technical---cases: (i) $n \! \in \! \mathbb{N}$ 
and $k \! \in \! \lbrace 1,2,\dotsc,K \rbrace$ such that 
$\alpha_{p_{\mathfrak{s}}} \! := \! \alpha_{k} \! = \! \infty$; and 
(ii) $n \! \in \! \mathbb{N}$ and $k \! \in \! \lbrace 1,2,\dotsc,K 
\rbrace$ such that $\alpha_{p_{\mathfrak{s}}} \! := \! \alpha_{k} \! 
\neq \! \infty$. The proof for the case $\alpha_{p_{\mathfrak{s}}} 
\! := \! \alpha_{k} \! \neq \! \infty$, $k \! \in \! \lbrace 1,2,\dotsc,K 
\rbrace$, will be considered in detail (see $\pmb{(\mathrm{A})}$ below), 
whilst the case $\alpha_{p_{\mathfrak{s}}} \! := \! \alpha_{k} \! = \! 
\infty$, $k \! \in \! \lbrace 1,2,\dotsc,K \rbrace$, can be proved, 
modulo technical and computational particulars, analogously (see 
$\pmb{(\mathrm{B})}$ below).

$\pmb{(\mathrm{A})}$ Let $\widetilde{V} \colon \overline{\mathbb{R}} 
\setminus \lbrace \alpha_{1},\alpha_{2},\dotsc,\alpha_{K} \rbrace \! \to \! 
\mathbb{R}$ satisfy conditions~\eqref{eq20}--\eqref{eq22} and be regular. 
For $n \! \in \! \mathbb{N}$ and $k \! \in \! \lbrace 1,2,\dotsc,K \rbrace$ 
such that $\alpha_{p_{\mathfrak{s}}} \! := \! \alpha_{k} \! \neq \! \infty$, 
one begins by showing that the support, that is, $\supp (\mu^{f}_{
\widetilde{V}}) \! =: \! J_{f}$, of each member of the associated family 
of equilibrium measures consists of the union of a finite number of 
disjoint and bounded real intervals. Recalling {}from Lemma~\ref{lem3.1} 
that, for $n \! \in \! \mathbb{N}$ and $k \! \in \! \lbrace 1,2,\dotsc,K 
\rbrace$ such that $\alpha_{p_{\mathfrak{s}}} \! := \! \alpha_{k} \! \neq 
\! \infty$, $(\overline{\mathbb{R}} \setminus \lbrace \alpha_{1},
\alpha_{2},\dotsc,\alpha_{K} \rbrace \supset)$ $J_{f} \! \subseteq \! 
\mathbb{R} \setminus (\lbrace \lvert x \rvert \! \geqslant \! T_{M_{f}} 
\rbrace \cup \cup_{\underset{q \neq \mathfrak{s}-1}{q=1}}^{
\mathfrak{s}} \operatorname{clos}(\mathscr{O}_{\frac{1}{T_{M_{f}}}}
(\alpha_{p_{q}})))$, where $T_{M_{f}} \! > \! 1$ is chosen (e.g., $T_{M_{f}} 
\! = \! K(1 \! + \! \max \lbrace \mathstrut \lvert \alpha_{p_{q}} \rvert, \, 
q \! = \! 1,\dotsc,\mathfrak{s} \! - \! 2,\mathfrak{s} \rbrace \! + \! 3(\min 
\lbrace \mathstrut \lvert \alpha_{p_{i}} \! - \! \alpha_{p_{j}} \rvert, \, i \! 
\neq \! j \! \in \! \lbrace 1,\dotsc,\mathfrak{s} \! - \! 2,\mathfrak{s} 
\rbrace \rbrace)^{-1}))$ so that $\mathscr{O}_{\frac{1}{T_{M_{f}}}}(\alpha_{p_{i}}) 
\cap \mathscr{O}_{\frac{1}{T_{M_{f}}}}(\alpha_{p_{j}}) \! = \! \varnothing$, 
$i \! \neq \! j \! \in \! \lbrace 1,\dotsc,\mathfrak{s} \! - \! 2,\mathfrak{s} 
\rbrace$, and that $\widetilde{V}$ is real analytic on $J_{f}$ with an analytic 
extension to the following open neighbourhood (not the only choice 
possible!) of $J_{f}$, $\mathbb{U}_{f} \! = \! \lbrace \mathstrut z \! \in 
\! \overline{\mathbb{C}} \setminus \lbrace \alpha_{p_{1}},\alpha_{p_{2}},
\dotsc,\alpha_{p_{\mathfrak{s}}} \rbrace; \, \inf_{\xi \in J_{f}} \lvert z 
\! - \! \xi \rvert \! < \! r \rbrace$, with $r \! \in \! (0,1)$ chosen small 
enough so that $\mathbb{U}_{f} \cap \lbrace \alpha_{p_{1}},\alpha_{p_{2}},
\dotsc,\alpha_{p_{\mathfrak{s}}} \rbrace \! = \! \varnothing$, it follows 
{}from the corresponding result of Lemma~\ref{lem3.6}, that is, the 
elements of the associated family of equilibrium measures are 
absolutely continuous with respect to Lebesgue measure in the 
double-scaling limit $\mathscr{N},n \! \to \! \infty$ such that $z_{o} 
\! = \! 1 \! + \! o(1)$, and a calculation analogous to that subsumed 
in the proof of Lemma~2.26 of \cite{a58} that, for $n \! \in \! 
\mathbb{N}$ and $k \! \in \! \lbrace 1,2,\dotsc,K \rbrace$ such that 
$\alpha_{p_{\mathfrak{s}}} \! := \! \alpha_{k} \! \neq \! \infty$, the 
densities of the elements of the associated family of equilibrium 
measures have the representation $\md \mu_{\widetilde{V}}^{f}(x) 
\! = \! \psi_{\widetilde{V}}^{f}(x) \, \md x$, $x \! \in \! J_{f}$, with 
$\psi_{\widetilde{V}}^{f}$ determined below.

For $n \! \in \! \mathbb{N}$ and $k \! \in \! \lbrace 1,2,\dotsc,K 
\rbrace$ such that $\alpha_{p_{\mathfrak{s}}} \! := \! \alpha_{k} 
\! \neq \! \infty$, let
\begin{align} \label{eql3.7m}
\tilde{\mathcal{S}} \colon &\mathbb{N} \times \lbrace 1,2,\dotsc,K 
\rbrace \times \mathbb{C} \setminus (J_{f} \cup \lbrace \alpha_{p_{1}},
\dotsc,\alpha_{p_{\mathfrak{s}-2}},\alpha_{k} \rbrace) \! \ni \! (n,k,z) 
\! \mapsto \! 4 \mi \left(\dfrac{(n \! - \! 1)K \! + \! k}{n} \right)^{2}
(\mathcal{H} \psi_{\widetilde{V}}^{f})(z) \psi_{\widetilde{V}}^{f}(z) 
\nonumber \\
&-\dfrac{4 \mi}{\pi} \left(\dfrac{(n \! - \! 1)K \! + \! k}{n} \right) 
\left(\dfrac{(\varkappa_{nk} \! - \! 1)}{n(z \! - \! \alpha_{k})} \! + \! 
\sum_{q=1}^{\mathfrak{s}-2} \dfrac{\varkappa_{nk \tilde{k}_{q}}}{n
(z \! - \! \alpha_{p_{q}})} \right) \psi_{\widetilde{V}}^{f}(z) \! = \! 
\tilde{\mathcal{S}}(n,k,z) \! =: \! \tilde{\mathcal{S}}(z),
\end{align}
where $\mathcal{H} \colon \mathcal{L}^{2}(\mathbb{R}) \! \ni \! f \! 
\mapsto \! \pvi_{\raise-0.95ex\hbox{$\scriptstyle{}\mathbb{R}$}} 
\tfrac{f(\xi)}{z-\xi} \, \tfrac{\md \xi}{\pi} \! =: \! (\mathcal{H}f)(z)$ 
denotes the Hilbert transform, with $\pvi_{}$ the principle value 
integral, and
\begin{equation} \label{eql3.7n} 
\tilde{\mathfrak{H}} \colon \mathbb{N} \times \lbrace 1,2,\dotsc,K 
\rbrace \times \mathbb{C} \setminus (J_{f} \cup \lbrace \alpha_{p_{1}},
\dotsc,\alpha_{p_{\mathfrak{s}-2}},\alpha_{k} \rbrace) \! \ni \! 
(n,k,z) \! \mapsto \! (\mathfrak{F}_{f}(z))^{2} \! - \! \int_{J_{f}} 
\dfrac{\tilde{\mathcal{S}}(\xi)}{\xi \! - \! z} \, \dfrac{\md \xi}{2 \pi 
\mi} \! = \! \tilde{\mathfrak{H}}(n,k,z) \! =: \! \tilde{\mathfrak{H}}(z),
\end{equation}
where $\mathfrak{F}_{f}(z)$ is defined by Equation~\eqref{eql3.4a}. Via 
the distributional identities $(x \! - \! (x_{o} \! \pm \! \mi 0))^{-1} \! = 
\! (x \! - \! x_{o})^{-1} \! \pm \! \mi \pi \delta (x \! - \! x_{o})$, where 
$\delta (\pmb{\cdot})$ is the Dirac delta function, and $\int_{x_{1}}^{x_{2}}
f(\xi) \delta (\xi \! - \! x) \md \xi \! = \! 
\left\{
\begin{smallmatrix}
f(x), \, \, \, x \in (x_{1},x_{2}), \\
0, \, \, \, x \in \mathbb{R} \setminus (x_{1},x_{2}),
\end{smallmatrix}
\right.$ and the representation $\md \mu^{f}_{\widetilde{V}}(x) \! = \! 
\psi^{f}_{\widetilde{V}}(x) \, \md x$, $x \! \in \! J_{f}$, one shows that
\begin{equation} \label{eql3.7o} 
\tilde{\mathfrak{H}}_{\pm}(z) \! = \! 
\begin{cases}
(\mathfrak{F}_{f \, \pm}(z))^{2} \! - \! 
\pvi_{\raise-0.95ex\hbox{$\scriptstyle{}J_{f}$}} \tfrac{\tilde{\mathcal{S}}
(\xi)}{\xi -z} \, \tfrac{\md \xi}{2 \pi \mi} \! \mp \! \frac{1}{2} 
\tilde{\mathcal{S}}(z), &\text{$z \! \in \! J_{f}$,} \\
(\mathfrak{F}_{f}(z))^{2} \! - \! \int_{J_{f}} \tfrac{\tilde{\mathcal{S}}
(\xi)}{\xi -z} \, \tfrac{\md \xi}{2 \pi \mi}, &\text{$z \! \notin \! J_{f}$,}
\end{cases}
\end{equation}
and
\begin{equation} \label{eql3.7p} 
\mathfrak{F}_{f \, \pm}(z) \! = \! 
\begin{cases}
-\frac{1}{\mi \pi} \left(\frac{(\varkappa_{nk}-1)}{n(z-\alpha_{k})} \! + \! 
\sum\limits_{q=1}^{\mathfrak{s}-2} \frac{\varkappa_{nk \tilde{k}_{q}}}{n
(z-\alpha_{p_{q}})} \! - \! \pi \left(\frac{(n-1)K+k}{n} \right)(\mathcal{H} 
\psi^{f}_{\widetilde{V}})(z) \right) \! \mp \! \left(\frac{(n-1)K+k}{n} 
\right) \psi^{f}_{\widetilde{V}}(z), &\text{$z \! \in \! J_{f}$,} \\
-\frac{1}{\mi \pi} \left(\frac{(\varkappa_{nk}-1)}{n(z-\alpha_{k})} \! + \! 
\sum\limits_{q=1}^{\mathfrak{s}-2} \frac{\varkappa_{nk \tilde{k}_{q}}}{n
(z-\alpha_{p_{q}})} \! + \! \left(\frac{(n-1)K+k}{n} \right) \int_{J_{f}}
(\xi \! - \! z)^{-1} \psi^{f}_{\widetilde{V}}(\xi) \, \md \xi \right), 
&\text{$z \! \notin \! J_{f}$,}
\end{cases}
\end{equation}
where $\star_{\pm}(z) \! := \! \lim_{\varepsilon \downarrow 0} \star 
(z \! \pm \! \mi \varepsilon)$, $\star \! \in \! \lbrace \tilde{
\mathfrak{H}},\mathfrak{F}_{f} \rbrace$; thus, {}from the second 
line of Equation~\eqref{eql3.7p}, $\mathfrak{F}_{f \, +}(z) \! = \! 
\mathfrak{F}_{f \, -}(z) \! = \! \mathfrak{F}_{f}(z)$, $z \! \in \! 
\mathbb{R} \setminus (J_{f} \cup \lbrace \alpha_{p_{1}},\dotsc,
\alpha_{p_{\mathfrak{s}-2}},\alpha_{k} \rbrace)$, in which case, by 
Equations~\eqref{eql3.7n} and~\eqref{eql3.7o}, $\tilde{\mathfrak{H}}_{+}
(z) \! = \! \tilde{\mathfrak{H}}_{-}(z) \! = \! \tilde{\mathfrak{H}}(z)$. For 
$z \! \in \! J_{f}$, one notes {}from Equation~\eqref{eql3.7o} and the first 
line of Equation~\eqref{eql3.7p} that $\tilde{\mathfrak{H}}_{+}(z) \! - \! 
\tilde{\mathfrak{H}}_{-}(z) \! = \! (\mathfrak{F}_{f \, +}(z))^{2} \! - \! 
(\mathfrak{F}_{f \, -}(z))^{2} \! - \! \tilde{\mathcal{S}}(z)$, where
\begin{align*}
(\mathfrak{F}_{f \, \pm}(z))^{2} =& \, -\dfrac{1}{\pi^{2}} \left(\dfrac{
(\varkappa_{nk} \! - \! 1)}{n(z \! - \! \alpha_{k})} \! + \! \sum_{q=1}^{
\mathfrak{s}-2} \dfrac{\varkappa_{nk \tilde{k}_{q}}}{n(z \! - \! 
\alpha_{p_{q}})} \right)^{2} \! + \! \dfrac{2}{\pi} \left(\dfrac{(n \! - \! 
1)K \! + \! k}{n} \right) \left(\dfrac{(\varkappa_{nk} \! - \! 1)}{n(z \! - 
\! \alpha_{k})} \! + \! \sum_{q=1}^{\mathfrak{s}-2} \dfrac{\varkappa_{nk 
\tilde{k}_{q}}}{n(z \! - \! \alpha_{p_{q}})} \right)(\mathcal{H} 
\psi^{f}_{\widetilde{V}})(z) \\
\mp& \, \dfrac{2 \mi}{\pi} \left(\dfrac{(n \! - \! 1)K \! + \! k}{n} \right) 
\left(\dfrac{(\varkappa_{nk} \! - \! 1)}{n(z \! - \! \alpha_{k})} \! + \! 
\sum_{q=1}^{\mathfrak{s}-2} \dfrac{\varkappa_{nk \tilde{k}_{q}}}{n(z \! - 
\! \alpha_{p_{q}})} \right) \psi^{f}_{\widetilde{V}}(z) \! \pm \! 2 \mi 
\left(\dfrac{(n \! - \! 1)K \! + \! k}{n} \right)^{2}(\mathcal{H} 
\psi^{f}_{\widetilde{V}})(z) \psi^{f}_{\widetilde{V}}(z) \\
-& \, \left(\dfrac{(n \! - \! 1)K \! + \! k}{n} \right)^{2}((\mathcal{H} 
\psi^{f}_{\widetilde{V}})(z))^{2} \! + \! \left(\dfrac{(n \! - \! 1)K \! 
+ \! k}{n} \right)^{2} (\psi^{f}_{\widetilde{V}}(z))^{2},
\end{align*}
whence, via Equation~\eqref{eql3.7m},
\begin{align*}
(\mathfrak{F}_{f \, +}(z))^{2} \! - \! (\mathfrak{F}_{f \, -}(z))^{2} =& \, 
-\dfrac{4 \mi}{\pi} \left(\dfrac{(n \! - \! 1)K \! + \! k}{n} \right) 
\left(\dfrac{(\varkappa_{nk} \! - \! 1)}{n(z \! - \! \alpha_{k})} \! + \! 
\sum_{q=1}^{\mathfrak{s}-2} \dfrac{\varkappa_{nk \tilde{k}_{q}}}{n
(z \! - \! \alpha_{p_{q}})} \right) \psi^{f}_{\widetilde{V}}(z) \\
+& \, 4 \mi \left(\dfrac{(n \! - \! 1)K \! + \! k}{n} \right)^{2}
(\mathcal{H} \psi^{f}_{\widetilde{V}})(z) \psi^{f}_{\widetilde{V}}
(z) \! = \! \tilde{\mathcal{S}}(z),
\end{align*}
that is, $(\mathfrak{F}_{f \, +}(z))^{2} \! - \! (\mathfrak{F}_{f \, -}
(z))^{2} \! - \! \tilde{\mathcal{S}}(z) \! = \! 0$ $\Rightarrow$ 
$\tilde{\mathfrak{H}}_{+}(z) \! - \! \tilde{\mathfrak{H}}_{-}(z) \! 
= \! 0$; thus, for $z \! \in \! J_{f}$, via Equation~\eqref{eql3.7n}, 
$\tilde{\mathfrak{H}}_{+}(z) \! = \! \tilde{\mathfrak{H}}_{-}(z) \! = \! 
\tilde{\mathfrak{H}}(z)$. The above argument shows, in particular, 
that, for $n \! \in \! \mathbb{N}$ and $k \! \in \! \lbrace 1,2,\dotsc,
K \rbrace$ such that $\alpha_{p_{\mathfrak{s}}} \! := \! \alpha_{k} 
\! \neq \! \infty$, $\tilde{\mathfrak{H}}(z)$ is analytic across 
$\mathbb{R} \setminus \lbrace \alpha_{p_{1}},\dotsc,
\alpha_{p_{\mathfrak{s}-2}},\alpha_{k} \rbrace$; in fact, 
$\tilde{\mathfrak{H}}(z)$ is entire (resp., meromorphic) for 
$z \! \in \! \mathbb{C} \setminus \lbrace\alpha_{p_{1}},\dotsc,
\alpha_{p_{\mathfrak{s}-2}},\alpha_{k} \rbrace$ (resp., $z \! \in 
\! \mathbb{C})$. Recalling that, for $n \! \in \! \mathbb{N}$ 
and $k \! \in \! \lbrace 1,2,\dotsc,K \rbrace$ such that 
$\alpha_{p_{\mathfrak{s}}} \! := \! \alpha_{k} \! \neq \! \infty$, 
$\mu^{f}_{\widetilde{V}} \! \in \! \mathscr{M}_{1}(\mathbb{R})$, it follows 
that, for $\xi \! \in \! J_{f}$ and $z \! \notin \! J_{f}$ such that $\lvert 
(z \! - \! \alpha_{k})/(\xi \! - \! \alpha_{k}) \rvert  \! \ll \! 1$ (e.g., $0 \! 
< \! \lvert z \! - \! \alpha_{k} \rvert \! \ll \! \min \lbrace \min_{i \neq j 
\in \lbrace 1,\dotsc,\mathfrak{s}-2,\mathfrak{s} \rbrace} \lbrace \lvert 
\alpha_{p_{i}} \! - \! \alpha_{p_{j}} \rvert \rbrace,\inf_{\xi \in J_{f}} \lbrace 
\lvert \xi \! - \! \alpha_{k} \rvert \rbrace \rbrace)$, via the expansion 
$\tfrac{1}{(z-\alpha_{k})-(\xi -\alpha_{k})} \! = \! -\sum_{j=0}^{l} 
\tfrac{(z-\alpha_{k})^{j}}{(\xi -\alpha_{k})^{j+1}} \! + \! \tfrac{(z-
\alpha_{k})^{l+1}}{(\xi -\alpha_{k})^{l+1}(z-\xi)}$, $l \! \in \! 
\mathbb{N}_{0}$, and Equation~\eqref{eql3.7n},
\begin{align*}
\tilde{\mathfrak{H}}(z) \underset{z \to \alpha_{k}}{=}& \, 
-\dfrac{1}{(z \! - \! \alpha_{k})^{2}} \left(\dfrac{1}{\pi} \left(
\dfrac{\varkappa_{nk} \! - \! 1}{n} \right) \right)^{2} \! + \! 
\dfrac{1}{(z \! - \! \alpha_{k})} \left(\dfrac{2}{\pi^{2}} \left(
\dfrac{\varkappa_{nk} \! - \! 1}{n} \right) \left(\sum_{q=1}^{
\mathfrak{s}-2} \dfrac{\varkappa_{nk \tilde{k}_{q}}}{n(\alpha_{p_{q}} 
\! - \! \alpha_{k})} \right. \right. \\
-&\left. \left. \, \left(\dfrac{(n \! - \! 1)K \! + \! k}{n} \right) 
\int_{J_{f}}(\xi \! - \! \alpha_{k})^{-1} \psi^{f}_{\widetilde{V}}
(\xi) \, \md \xi \right) \right) \! + \! \mathcal{O}(1),
\end{align*}
which shows that $\tilde{\mathfrak{H}}(z)$ has a pole of order 
two at $\alpha_{p_{\mathfrak{s}}} \! := \! \alpha_{k}$, with 
$\operatorname{Res}(\tilde{\mathfrak{H}}(z);\alpha_{k}) \! = \! 
\tfrac{2}{\pi^{2}}(\tfrac{\varkappa_{nk}-1}{n})(\sum_{q=1}^{\mathfrak{s}
-2} \tfrac{\varkappa_{nk \tilde{k}_{q}}}{n(\alpha_{p_{q}}-\alpha_{k})} 
\! - \! (\tfrac{(n-1)K+k}{n}) \int_{J_{f}}(\xi \! - \! \alpha_{k})^{-1} 
\psi^{f}_{\widetilde{V}}(\xi) \, \md \xi)$. One learns {}from the above 
analysis that, for $n \! \in \! \mathbb{N}$ and $k \! \in \! \lbrace 
1,2,\dotsc,K \rbrace$ such that $\alpha_{p_{\mathfrak{s}}} \! := \! 
\alpha_{k} \! \neq \! \infty$, $\prod_{q=1}^{\mathfrak{s}-2}(z \! - \! 
\alpha_{p_{q}})^{2}(z \! - \! \alpha_{k})^{2} \tilde{\mathfrak{H}}(z)$ 
is entire;\footnote{It should be remarked that, for $n \! \in \! 
\mathbb{N}$ and $k \! \in \! \lbrace 1,2,\dotsc,K \rbrace$ such that 
$\alpha_{p_{\mathfrak{s}}} \! := \! \alpha_{k} \! \neq \! \infty$, there 
is no need to subsume the contribution of a `polar term' of the type 
(cf. the proof of Lemma~\ref{lem3.6}) $\prod_{j=1}^{\tilde{N}+1}
[\tilde{B}_{j-1},\tilde{A}_{j}]$ as an additional, multiplicative factor in the 
expression $\prod_{q=1}^{\mathfrak{s}-2}(z \! - \! \alpha_{p_{q}})^{2}
(z \! - \! \alpha_{k})^{2} \tilde{\mathfrak{H}}(z)$, since, in the 
double-scaling limit $\mathscr{N},n \! \to \! \infty$ such that 
$z_{o} \! = \! 1 \! + \! o(1)$, the associated, weighted Fekete points 
(cf. Lemma~\ref{lem3.6}) $\tilde{x}_{j}^{\ast}$, $j \! = \! 1,2,\dotsc,
(n \! - \! 1)K \! + \! k$, all lie in a proper compact subset of 
$\mathbb{R} \setminus \lbrace \alpha_{p_{1}},\dotsc,
\alpha_{p_{\mathfrak{s}-2}},\alpha_{k} \rbrace$. This additional 
observation stems {}from mimicking part of the calculations (but 
adapted to the present situation) in Section~4 of \cite{a58}; 
alternatively, one may apply suitably modified versions of the results 
of the present monograph to measures that are supported on a 
`large' real interval containing $J_{f}$; e.g., take $\tilde{B}_{0} \! 
\ll \! \min \lbrace J_{f} \rbrace$, $\tilde{A}_{\tilde{N}+1} \! \gg \! \max 
\lbrace J_{f} \rbrace$, exclude the factor $\prod_{j=1}^{\tilde{N}}
(z \! - \! \tilde{B}_{j})(z \! - \! \tilde{A}_{j})$, and consider the limits 
$\tilde{B}_{0} \! \to \! -\infty$ and $\tilde{A}_{\tilde{N}+1} \! \to \! 
+\infty$: the details are left to the interested reader.} in particular, 
look at the behaviour of $\prod_{q=1}^{\mathfrak{s}-2}(z \! - \! 
\alpha_{p_{q}})^{2}(z \! - \! \alpha_{k})^{2} \tilde{\mathfrak{H}}(z)$ as 
$z \! \to \! \alpha_{p_{\mathfrak{s}-1}} \! = \! \infty$. For $n \! \in \! 
\mathbb{N}$ and $k \! \in \! \lbrace 1,2,\dotsc,K \rbrace$ such that 
$\alpha_{p_{\mathfrak{s}}} \! := \! \alpha_{k} \! \neq \! \infty$, one has, 
via Equations~\eqref{eql3.4a}, \eqref{eql3.7m}, and~\eqref{eql3.7n},
\begin{align} \label{eql3.7q} 
\tilde{\mathfrak{H}}(z) =& \, -\dfrac{(\varkappa_{nk} \! - \! 1)^{2}}{
\pi^{2}n^{2}(z \! - \! \alpha_{k})^{2}} \! - \! \dfrac{2(\varkappa_{nk} 
\! - \! 1)}{\pi^{2}n(z \! - \! \alpha_{k})} \sum_{q=1}^{\mathfrak{s}-2} 
\dfrac{\varkappa_{nk \tilde{k}_{q}}}{n(z \! - \! \alpha_{p_{q}})} \! - \! 
\dfrac{1}{\pi^{2}} \left(\sum_{q=1}^{\mathfrak{s}-2} \dfrac{\varkappa_{nk 
\tilde{k}_{q}}}{n(z \! - \! \alpha_{p_{q}})} \right)^{2} \! - \! 
\dfrac{2(\varkappa_{nk} \! - \! 1)}{\pi^{2}n(z \! - \! \alpha_{k})} 
\left(\dfrac{(n \! - \! 1)K \! + \! k}{n} \right) \nonumber \\
\times& \, \int_{J_{f}} \dfrac{\md \mu^{f}_{\widetilde{V}}(\xi)}{\xi \! - \! 
z} \! - \! \dfrac{2}{\pi^{2}} \left(\dfrac{(n \! - \! 1)K \! + \! k}{n} \right) 
\sum_{q=1}^{\mathfrak{s}-2} \dfrac{\varkappa_{nk \tilde{k}_{q}}}{n
(z \! - \! \alpha_{p_{q}})} \int_{J_{f}} \dfrac{\md \mu^{f}_{\widetilde{V}}
(\xi)}{\xi \! - \! z} \! - \! \dfrac{1}{\pi^{2}} \left(\dfrac{(n \! - \! 1)K \! 
+ \! k}{n} \right)^{2} \left(\int_{J_{f}} \dfrac{\md \mu^{f}_{\widetilde{V}}
(\xi)}{\xi \! - \! z} \right)^{2} \nonumber \\
-& \, \int_{J_{f}} \dfrac{\tilde{\mathcal{S}}(\xi)}{\xi \! - \! z} \, 
\dfrac{\md \xi}{2 \pi \mi}, \quad z \! \in \! \mathbb{C} \setminus 
(J_{f} \cup \lbrace \alpha_{p_{1}},\dotsc,\alpha_{p_{\mathfrak{s}-2}},
\alpha_{k} \rbrace).
\end{align}
In order to proceed, a simplified expression for $\tilde{\mathcal{S}}(z)$ is 
requisite. Using the representation $\md \mu^{f}_{\widetilde{V}}(x) \! = \! 
\psi^{f}_{\widetilde{V}}(x) \, \md x$, $x \! \in \! J_{f}$, in the definition 
of $g^{f}(z)$ given by Equation~\eqref{eql3.4gee3}, one shows that 
(assuming differentiation commutes with taking boundary values)
\begin{equation*}
(g^{f}_{\pm}(z))^{\prime} \! := \! \lim_{\varepsilon \downarrow 0}
(g^{f})^{\prime}(z \! \pm \! \mi \varepsilon) \! = \! 
\begin{cases}
-\frac{(\varkappa_{nk}-1)}{n(z-\alpha_{k})} \! - \! \sum\limits_{q=1}^{
\mathfrak{s}-2} \frac{\varkappa_{nk \tilde{k}_{q}}}{n(z-\alpha_{p_{q}})} 
\! - \! \left(\frac{(n-1)K+k}{n} \right) 
\pvi_{\raise-1.05ex\hbox{$\scriptstyle{}J_{f}$}} 
\tfrac{\psi^{f}_{\widetilde{V}}(\xi)}{\xi -z} \, \md \xi \! \mp \! 
\mi \pi \left(\frac{(n-1)K+k}{n} \right) \psi^{f}_{\widetilde{V}}(z), 
&\text{$z \! \in \! J_{f}$,} \\
-\frac{(\varkappa_{nk}-1)}{n(z-\alpha_{k})} \! - \! \sum\limits_{q=1}^{
\mathfrak{s}-1} \frac{\varkappa_{nk \tilde{k}_{q}}}{n(z-\alpha_{p_{q}})} 
\! - \! \left(\frac{(n-1)K+k}{n} \right) 
\int_{\raise-1.05ex\hbox{$\scriptstyle{}J_{f}$}} 
\tfrac{\psi^{f}_{\widetilde{V}}(\xi)}{\xi -z} \, \md \xi, 
&\text{$z \! \notin \! J_{f}$,}
\end{cases}
\end{equation*}
where the prime denotes differentiation with respect to $z$, whence, 
by linearity,
\begin{equation*}
(g^{f}_{+}(z))^{\prime} \! + \! (g^{f}_{-}(z))^{\prime} \! = \! 
-\dfrac{2(\varkappa_{nk} \! - \! 1)}{n(z \! - \! \alpha_{k})} \! - \! 2 
\sum_{q=1}^{\mathfrak{s}-2} \dfrac{\varkappa_{nk \tilde{k}_{q}}}{n(z \! 
- \! \alpha_{p_{q}})} \! + \! 2 \pi \left(\dfrac{(n \! - \! 1)K \! + \! k}{n} 
\right)(\mathcal{H} \psi^{f}_{\widetilde{V}})(z), \quad z \! \in \! J_{f},
\end{equation*}
and
\begin{equation*}
(g^{f}_{+}(z))^{\prime} \! - \! (g^{f}_{-}(z))^{\prime} \! = \! 
\begin{cases}
-2 \pi \mi \left(\frac{(n-1)K+k}{n} \right) \psi^{f}_{\widetilde{V}}(z), 
&\text{$z \! \in \! J_{f}$,} \\
0, &\text{$z \! \notin \! J_{f}$.}
\end{cases}
\end{equation*}
Noting {}from item~$\pmb{(2)}$ of Lemma~\ref{lem3.8} below that, in a 
piecewise-continuous sense, $(g^{f}_{+}(z))^{\prime} \! + \! (g^{f}_{-}
(z))^{\prime} \! - \! \widetilde{V}^{\prime}(z) \! = \! 0$, $z \! \in \! 
J_{f}$, one shows {}from the first of the above pair of formulae that
\begin{equation} \label{eql3.7r} 
(\mathcal{H} \psi^{f}_{\widetilde{V}})(z) \! = \! \dfrac{1}{2 \pi} \left(
\dfrac{(n \! - \! 1)K \! + \! k}{n} \right)^{-1} \left(\dfrac{2(\varkappa_{nk} 
\! - \! 1)}{n(z \! - \! \alpha_{k})} \! + \! 2 \sum_{q=1}^{\mathfrak{s}-2} 
\dfrac{\varkappa_{nk \tilde{k}_{q}}}{n(z \! - \! \alpha_{p_{q}})} \! + \! 
\widetilde{V}^{\prime}(z) \right), \quad z \! \in \! J_{f}.
\end{equation}
Hence, via Equations~\eqref{eql3.7m} and~\eqref{eql3.7r}, one arrives at, 
for $n \! \in \! \mathbb{N}$ and $k \! \in \! \lbrace 1,2,\dotsc,K \rbrace$ 
such that $\alpha_{p_{\mathfrak{s}}} \! := \! \alpha_{k} \! \neq \! \infty$,
\begin{equation} \label{eql3.7s} 
\tilde{\mathcal{S}}(z) \! = \! \dfrac{2 \mi}{\pi} \left(\dfrac{(n \! - \! 1)K 
\! + \! k}{n} \right) \widetilde{V}^{\prime}(z) \psi^{f}_{\widetilde{V}}(z), 
\quad z \! \in \! J_{f}.
\end{equation}
{}From the observation $\deg (\prod_{q=1}^{\mathfrak{s}-2}(z \! - \! 
\alpha_{p_{q}})^{2}(z \! - \! \alpha_{k})^{2}) \! = \! 2(\mathfrak{s} \! - 
\! 1)$ and the fact that, for $n \! \in \! \mathbb{N}$ and $k \! \in \! 
\lbrace 1,2,\dotsc,K \rbrace$ such that $\alpha_{p_{\mathfrak{s}}} \! 
:= \! \alpha_{k} \! \neq \! \infty$, $\mu^{f}_{\widetilde{V}} \! \in \! 
\mathscr{M}_{1}(\mathbb{R})$, it follows that, for $\xi \! \in \! J_{f}$ 
and $z \! \notin \! J_{f}$ such that $\lvert \xi/z \rvert \! \ll \! 1$ 
(e.g., $\lvert z \rvert \! \gg \! \max \lbrace \max_{i \neq j \in 
\lbrace 1,\dotsc,\mathfrak{s}-2,\mathfrak{s} \rbrace} \lbrace 
\lvert \alpha_{p_{i}} \! - \! \alpha_{p_{j}} \rvert \rbrace,\max_{q=1,
\dotsc,\mathfrak{s}-2,\mathfrak{s}} \lbrace \lvert \alpha_{p_{q}} 
\rvert \rbrace,\max \lbrace J_{f} \rbrace \rbrace)$, via the expansion 
$\tfrac{1}{\xi -z} \! = \! -\sum_{j=0}^{l} \tfrac{\xi^{j}}{z^{j+1}} \! + 
\! \tfrac{\xi^{l+1}}{z^{l+1}(\xi -z)}$, $l \! \in \! \mathbb{N}_{0}$, 
$\prod_{q=1}^{\mathfrak{s}-2}(z \! - \! \alpha_{p_{q}})^{2}(z \! - 
\! \alpha_{k})^{2} \tilde{\mathfrak{H}}(z) \! - \! \sum_{m=0}^{2 
\mathfrak{s}-3} \tilde{\rho}_{m}(n,k)z^{m} \! =_{z \to \alpha_{
p_{\mathfrak{s}-1}} = \infty} \! \mathcal{O}(z^{-1})$, with 
$\tilde{\rho}_{m} \colon \mathbb{N} \times \lbrace 1,2,\dotsc,
K \rbrace \! \to \! \mathbb{R}$, $m \! = \! 0,1,\dotsc,2 
\mathfrak{s} \! - \! 3$, determined below, where, in particular, 
via Equations~\eqref{eql3.7q} and~\eqref{eql3.7s}, 
$\tilde{\rho}_{2 \mathfrak{s}-3}(n,k) \! = \! \tfrac{1}{\pi^{2}}
(\tfrac{(n-1)K+k}{n}) \int_{J_{f}} \widetilde{V}^{\prime}(\xi) 
\psi^{f}_{\widetilde{V}}(\xi) \, \md \xi$; hence, as $\prod_{q=1}^{
\mathfrak{s}-2}(z \! - \! \alpha_{p_{q}})^{2}(z \! - \! \alpha_{k})^{2} 
\tilde{\mathfrak{H}}(z)$ is entire, it follows, by a generalisation of 
Liouville's Theorem, that
\begin{equation} \label{eql3.7t} 
\tilde{\mathfrak{H}}(z) \! = \! \dfrac{\sum_{m=0}^{2 \mathfrak{s}-3} 
\tilde{\rho}_{m}(n,k)z^{m}}{\prod_{q=1}^{\mathfrak{s}-2}(z \! - \! 
\alpha_{p_{q}})^{2}(z \! - \! \alpha_{k})^{2}}, \qquad z \! \in \! 
\mathbb{C} \setminus \lbrace \alpha_{p_{1}},\dotsc,
\alpha_{p_{\mathfrak{s}-2}},\alpha_{k} \rbrace.
\end{equation}
It remains to determine $\tilde{\rho}_{m}(n,k)$, $m \! = \! 0,1,\dotsc,2 
\mathfrak{s} \! - \! 3$. Via Equations~\eqref{eql3.7q} and~\eqref{eql3.7s}, 
the representation $\md \mu^{f}_{\widetilde{V}}(x) \! = \! \psi^{f}_{
\widetilde{V}}(x) \, \md x$, $x \! \in \! J_{f}$, and the fact that 
$\mu^{f}_{\widetilde{V}} \! \in \! \mathscr{M}_{1}(\mathbb{R})$, one 
shows that, for $\xi \! \in \! J_{f}$ and $z \! \notin \! J_{f}$ such that 
$\lvert \xi/z \rvert \! \ll \! 1$ (e.g., $\lvert z \rvert \! \gg \! \max \lbrace 
\max_{i \neq j \in \lbrace 1,\dotsc,\mathfrak{s}-2,\mathfrak{s} \rbrace} 
\lbrace \lvert \alpha_{p_{i}} \! - \! \alpha_{p_{j}} \rvert \rbrace,\max_{q=
1,\dotsc,\mathfrak{s}-2,\mathfrak{s}} \lbrace \lvert \alpha_{p_{q}} \rvert 
\rbrace,\max \lbrace J_{f} \rbrace \rbrace)$, via the expansion 
$\tfrac{1}{\xi -z} \! = \! -\sum_{j=0}^{l} \tfrac{\xi^{j}}{z^{j+1}} \! 
+ \! \tfrac{\xi^{l+1}}{z^{l+1}(\xi -z)}$, $l \! \in \! \mathbb{N}_{0}$,
\begin{equation} \label{eql3.7u} 
-\prod_{q=1}^{\mathfrak{s}-2}(z \! - \! \alpha_{p_{q}})^{2}(z \! - \! 
\alpha_{k})^{2} \int_{J_{f}} \dfrac{\tilde{\mathcal{S}}(\xi)}{\xi \! - \! z} 
\, \dfrac{\md \xi}{2 \pi \mi} \underset{z \to \alpha_{p_{\mathfrak{s}
-1}} = \infty}{=} \sum_{m=0}^{2 \mathfrak{s}-3} \sum_{j=m}^{2 
\mathfrak{s}-3} \tilde{w}_{j-m}(n,k) \tilde{c}_{j+1}(n,k)z^{m} \! + \! 
\mathcal{O}(z^{-1}),
\end{equation}
where
\begin{gather}
\tilde{w}_{r}(n,k) \! = \! \dfrac{1}{\pi^{2}} \left(\dfrac{(n \! - \! 
1)K \! + \! k}{n} \right) \int_{J_{f}} \xi^{r} \widetilde{V}^{\prime}(\xi) 
\psi_{\widetilde{V}}^{f}(\xi) \, \md \xi, \quad r \! \in \! \mathbb{N}_{0}, 
\label{eql3.7v} \\
\tilde{c}_{l}(n,k) \! = \! \mathlarger{\sum_{\underset{\underset{
\sum_{\underset{r \neq \mathfrak{s}-1}{r=1}}^{\mathfrak{s}}i_{r}=2
(\mathfrak{s}-1)-l}{r \in \lbrace 1,2,\dotsc,\mathfrak{s} \rbrace \setminus 
\lbrace \mathfrak{s}-1 \rbrace}}{i_{r}=0,1,2}}}(-1)^{2(\mathfrak{s}-1)-l}
(2!)^{\mathfrak{s}-1} \prod_{\substack{m=1\\m \neq \mathfrak{s}
-1}}^{\mathfrak{s}} \dfrac{1}{i_{m}!(2 \! - \! i_{m})!} \prod_{j=1}^{
\mathfrak{s}-2}(\alpha_{p_{j}})^{i_{j}}(\alpha_{k})^{i_{\mathfrak{s}}}, 
\quad l \! = \! 0,1,\dotsc,2(\mathfrak{s} \! - \! 1), \label{eql3.7w}
\end{gather}
\begin{equation} \label{eql3.7x} 
-\dfrac{1}{\pi^{2}} \left(\dfrac{(n \! - \! 1)K \! + \! k}{n} \right)^{2} 
\prod_{q=1}^{\mathfrak{s}-2}(z \! - \! \alpha_{p_{q}})^{2}(z \! - \! 
\alpha_{k})^{2} \left(\int_{J_{f}} \dfrac{\md \mu^{f}_{\widetilde{V}}
(\xi)}{\xi \! - \! z} \right)^{2} \underset{z \to \alpha_{p_{\mathfrak{s}
-1}} = \infty}{=} \sum_{m=0}^{2(\mathfrak{s}-2)} \sum_{j=m}^{2
(\mathfrak{s}-2)} \tilde{u}_{j-m}(n,k) \tilde{c}_{j+2}^{\sharp}(n,k)z^{m} 
\! + \! \mathcal{O}(z^{-1}),
\end{equation}
where
\begin{gather}
\tilde{u}_{m}(n,k) \! = \! \sum_{j=0}^{m} \tilde{\nu}_{j}(n,k) 
\tilde{\nu}_{m-j}(n,k), \quad m \! = \! 0,1,\dotsc,2(\mathfrak{s} 
\! - \! 2), \label{eql3.7y} \\
\tilde{\nu}_{r}(n,k) \! = \! \int_{J_{f}} \xi^{r} \psi_{\widetilde{V}}^{f}
(\xi) \, \md \xi, \quad r \! \in \! \mathbb{N}_{0}, \label{eql3.7z} \\
\tilde{c}_{l}^{\sharp}(n,k) \! = \! -\dfrac{1}{\pi^{2}} \left(
\dfrac{(n \! - \! 1)K \! + \! k}{n} \right)^{2} \tilde{c}_{l}(n,k), 
\quad l \! = \! 0,1,\dotsc,2(\mathfrak{s} \! - \! 1), \label{eql3.7a1}
\end{gather}
\begin{align} \label{eql3.7a2} 
-& \, \prod_{q=1}^{\mathfrak{s}-2}(z \! - \! \alpha_{p_{q}})^{2}(z \! 
- \! \alpha_{k})^{2} \left(\dfrac{2(\varkappa_{nk} \! - \! 1)}{\pi^{2}n
(z \! - \! \alpha_{k})} \left(\dfrac{(n \! - \! 1)K \! + \! k}{n} \right) 
\int_{J_{f}} \dfrac{\md \mu^{f}_{\widetilde{V}}(\xi)}{\xi \! - \! z} \! 
+ \! \dfrac{2}{\pi^{2}} \left(\dfrac{(n \! - \! 1)K \! + \! k}{n} \right) 
\sum_{q=1}^{\mathfrak{s}-2} \dfrac{\varkappa_{nk \tilde{k}_{q}}}{n
(z \! - \! \alpha_{p_{q}})} \right. \nonumber \\
&\left. \times \, \int_{J_{f}} \dfrac{\md \mu^{f}_{\widetilde{V}}(\xi)}{\xi 
\! - \! z} \right) \underset{z \to \alpha_{p_{\mathfrak{s}-1}} = \infty}{=} 
\, \sum_{m=0}^{2(\mathfrak{s}-2)} \sum_{j=m}^{2(\mathfrak{s}-2)} 
\tilde{\nu}_{j-m}(n,k) \tilde{c}_{j+1}^{\flat}(n,k)z^{m} \! + \! \mathcal{O}
(z^{-1}),
\end{align}
where
\begin{align} \label{eql3.7a3} 
\tilde{c}_{l}^{\flat}(n,k) =& \, \dfrac{2}{\pi^{2}} \left(
\dfrac{\varkappa_{nk} \! - \! 1}{n} \right) \left(\dfrac{(n \! - \! 1)K \! 
+ \! k}{n} \right) \mathlarger{\sum_{\underset{\underset{\underset{
\sum_{\underset{m^{\prime} \neq \mathfrak{s}-1}{m^{\prime}=1}}^{
\mathfrak{s}}i_{m^{\prime}}=2(\mathfrak{s}-1)-l-1}{r \in \lbrace 1,2,
\dotsc,\mathfrak{s}-2 \rbrace}}{i_{r}=0,1,2}}{i_{\mathfrak{s}}=0,1}}} 
\dfrac{(-1)^{2(\mathfrak{s}-1)-l-1}(2!)^{\mathfrak{s}-2}}{i_{\mathfrak{s}}!
(1 \! - \! i_{\mathfrak{s}})!} \prod_{m=1}^{\mathfrak{s}-2} 
\dfrac{1}{i_{m}!(2 \! - \! i_{m})!} \nonumber \\
\times& \, \prod_{j=1}^{\mathfrak{s}-2}(\alpha_{p_{j}})^{i_{j}}
(\alpha_{k})^{i_{\mathfrak{s}}} \! + \! \dfrac{2}{\pi^{2}} \left(
\dfrac{(n \! - \! 1)K \! + \! k}{n} \right) \mathlarger{\sum_{q=1}^{
\mathfrak{s}-2}} \dfrac{\varkappa_{nk \tilde{k}_{q}}}{n} \mathlarger{
\sum_{\underset{\underset{\underset{\sum_{\underset{m^{\prime} 
\neq \mathfrak{s}-1}{m^{\prime}=1}}^{\mathfrak{s}}i_{m^{\prime}}=
2(\mathfrak{s}-1)-l-1}{r \in \lbrace 1,2,\dotsc,\mathfrak{s} \rbrace 
\setminus \lbrace q,\mathfrak{s}-1 \rbrace}}{i_{r}=0,1,2}}{i_{q}=0,1}}} 
\dfrac{(-1)^{2(\mathfrak{s}-1)-l-1}(2!)^{\mathfrak{s}-2}}{i_{q}!
(1 \! - \! i_{q})!} \nonumber \\
\times& \, \prod_{\substack{m=1\\m \neq q,\mathfrak{s}-
1}}^{\mathfrak{s}} \dfrac{1}{i_{m}!(2 \! - \! i_{m})!} \prod_{j=1}^{
\mathfrak{s}-2}(\alpha_{p_{j}})^{i_{j}}(\alpha_{k})^{i_{\mathfrak{s}}}, 
\quad l \! = \! 0,1,\dotsc,2 \mathfrak{s} \! - \! 3,
\end{align}
and
\begin{align} \label{eql3.7a4} 
-& \, \prod_{q=1}^{\mathfrak{s}-2}(z \! - \! \alpha_{p_{q}})^{2}(z \! 
- \! \alpha_{k})^{2} \left(\dfrac{(\varkappa_{nk} \! - \! 1)^{2}}{\pi^{2}
n^{2}(z \! - \! \alpha_{k})^{2}} \! + \! \dfrac{2(\varkappa_{nk} \! - \! 
1)}{\pi^{2}n(z \! - \! \alpha_{k})} \sum_{q=1}^{\mathfrak{s}-2} 
\dfrac{\varkappa_{nk \tilde{k}_{q}}}{n(z \! - \! \alpha_{p_{q}})} \! + \! 
\dfrac{1}{\pi^{2}} \left(\sum_{q=1}^{\mathfrak{s}-2} \dfrac{\varkappa_{nk 
\tilde{k}_{q}}}{n(z \! - \! \alpha_{p_{q}})} \right)^{2} \right) \nonumber \\
&\underset{z \to \alpha_{p_{\mathfrak{s}-1}} = \infty}{=} \, 
\sum_{m=0}^{2(\mathfrak{s}-2)} \tilde{c}_{m}^{\natural}(n,k)z^{m} 
\! + \! \mathcal{O}(z^{-1}),
\end{align}
where
\begin{align} \label{eql3.7a5} 
\tilde{c}_{l}^{\natural}(n,k) =& \, -\dfrac{1}{\pi^{2}} \left(
\dfrac{\varkappa_{nk} \! - \! 1}{n} \right)^{2} \mathlarger{
\sum_{\underset{\underset{\sum_{m^{\prime}=1}^{\mathfrak{s}-2}
i_{m^{\prime}}=2(\mathfrak{s}-2)-l}{r \in \lbrace 1,2,\dotsc,
\mathfrak{s}-2 \rbrace}}{i_{r}=0,1,2}}}(-1)^{2(\mathfrak{s}-2)-l}
(2!)^{\mathfrak{s}-2} \prod_{m=1}^{\mathfrak{s}-2} \dfrac{1}{i_{m}!
(2 \! - \! i_{m})!} \prod_{j=1}^{\mathfrak{s}-2}(\alpha_{p_{j}})^{i_{j}} 
\nonumber \\
-& \, \dfrac{2}{\pi^{2}} \left(\dfrac{\varkappa_{nk} \! - \! 1}{n} \right) 
\mathlarger{\sum_{q=1}^{\mathfrak{s}-2}} \dfrac{\varkappa_{nk 
\tilde{k}_{q}}}{n} \mathlarger{\sum_{\underset{\underset{\underset{
\sum_{\underset{m^{\prime} \neq \mathfrak{s}-1}{m^{\prime}=1}}^{
\mathfrak{s}}i_{m^{\prime}}=2(\mathfrak{s}-2)-l}{r \in \lbrace 1,2,
\dotsc,\mathfrak{s}-2 \rbrace \setminus \lbrace q \rbrace}}{i_{r}=0,
1,2}}{i_{q},i_{\mathfrak{s}}=0,1}}} \dfrac{(-1)^{2(\mathfrak{s}-2)-l}
(2!)^{\mathfrak{s}-3}}{i_{q}!(1 \! - \! i_{q})!i_{\mathfrak{s}}!(1 \! - 
\! i_{\mathfrak{s}})!} \prod_{\substack{m=1\\m \neq q}}^{
\mathfrak{s}-2} \dfrac{1}{i_{m}!(2 \! - \! i_{m})!} \nonumber \\
\times& \, \prod_{j=1}^{\mathfrak{s}-2}(\alpha_{p_{j}})^{i_{j}}
(\alpha_{k})^{i_{\mathfrak{s}}} \! - \! \dfrac{1}{\pi^{2}} \mathlarger{
\sum_{q=1}^{\mathfrak{s}-2}} \left(\dfrac{\varkappa_{nk \tilde{k}_{q}}}{n} 
\right)^{2} \mathlarger{\sum_{\underset{\underset{\sum_{\underset{
m^{\prime} \neq q,\mathfrak{s}-1}{m^{\prime}=1}}^{\mathfrak{s}}
i_{m^{\prime}}=2(\mathfrak{s}-2)-l}{r \in \lbrace 1,2,\dotsc,\mathfrak{s}
-2 \rbrace \setminus \lbrace q \rbrace}}{i_{r},i_{\mathfrak{s}}=0,1,2}}}
(-1)^{2(\mathfrak{s}-2)-l}(2!)^{\mathfrak{s}-2} \prod_{\substack{m=
1\\m \neq q,\mathfrak{s}-1}}^{\mathfrak{s}} \dfrac{1}{i_{m}!(2 \! - \! 
i_{m})!} \nonumber \\
\times& \, \prod_{\substack{j=1\\j \neq q}}^{\mathfrak{s}-2}
(\alpha_{p_{j}})^{i_{j}}(\alpha_{k})^{i_{\mathfrak{s}}} \! - \! 
\dfrac{2}{\pi^{2}} \mathlarger{\sum_{p^{\prime}=1}^{\mathfrak{s}-3}} 
\mathlarger{\sum_{q=p^{\prime}+1}^{\mathfrak{s}-2}} \dfrac{
\varkappa_{nk \tilde{k}_{p^{\prime}}}}{n} \dfrac{\varkappa_{nk \tilde{k}_{q}}}{n} 
\mathlarger{\sum_{\underset{\underset{\underset{i_{p^{\prime}}+i_{q}+
i_{\mathfrak{s}}+ \sum_{\underset{m^{\prime} \neq p^{\prime},q}{m^{\prime}
=1}}^{\mathfrak{s}-2}i_{m^{\prime}}=2(\mathfrak{s}-2)-l}{r \in \lbrace 
1,2,\dotsc,\mathfrak{s}-2 \rbrace \setminus \lbrace p^{\prime},q 
\rbrace}}{i_{r},i_{\mathfrak{s}}=0,1,2}}{i_{p^{\prime}},i_{q}=0,1}}} 
\dfrac{(-1)^{2(\mathfrak{s}-2)-l}}{i_{p^{\prime}}!(1 \! - \! i_{p^{\prime}})!
i_{q}!(1 \! - \! i_{q})!} \nonumber \\
\times& \, \dfrac{(2!)^{\mathfrak{s}-3}}{i_{\mathfrak{s}}!(2 \! - \! 
i_{\mathfrak{s}})!} \prod_{\substack{m=1\\m \neq p^{\prime},
q}}^{\mathfrak{s}-2} \dfrac{1}{i_{m}!(2 \! - \! i_{m})!} 
(\alpha_{p_{p^{\prime}}})^{i_{p^{\prime}}}(\alpha_{p_{q}})^{i_{q}} 
\prod_{\substack{j=1\\j \neq p^{\prime},q}}^{\mathfrak{s}-2}
(\alpha_{p_{j}})^{i_{j}}(\alpha_{k})^{i_{\mathfrak{s}}}, \quad l \! = \! 
0,1,\dotsc,2(\mathfrak{s} \! - \! 2);
\end{align}
hence, via the asymptotic representation $\prod_{q=1}^{\mathfrak{s}-2}
(z \! - \! \alpha_{p_{q}})^{2}(z \! - \! \alpha_{k})^{2} \tilde{\mathfrak{H}}
(z) \! - \! \sum_{m=0}^{2 \mathfrak{s}-3} \tilde{\rho}_{m}(n,k)z^{m} \! 
=_{z \to \alpha_{p_{\mathfrak{s}-1}} = \infty} \! \mathcal{O}(z^{-1})$, 
the expansions~\eqref{eql3.7u}, \eqref{eql3.7x}, \eqref{eql3.7a2}, 
and~\eqref{eql3.7a4}, and Equations~\eqref{eql3.7n}, \eqref{eql3.7s}, 
and~\eqref{eql3.7t}, one shows that, for $n \! \in \! \mathbb{N}$ 
and $k \! \in \! \lbrace 1,2,\dotsc,K \rbrace$ such that 
$\alpha_{p_{\mathfrak{s}}} \! := \! \alpha_{k} \! \neq \! \infty$,
\begin{equation} \label{eql3.7a6} 
(\mathfrak{F}_{f}(z))^{2} \! = \! \dfrac{1}{\pi^{2}} \left(\dfrac{(n \! - 
\! 1)K \! + \! k}{n} \right) \int_{J_{f}} \dfrac{\widetilde{V}^{\prime}
(\xi) \psi_{\widetilde{V}}^{f}(\xi)}{\xi \! - \! z} \, \md \xi \! + \! 
\dfrac{\sum_{m=0}^{2 \mathfrak{s}-3} \tilde{\rho}_{m}(n,k)z^{m}}{
\prod_{q=1}^{\mathfrak{s}-2}(z \! - \! \alpha_{p_{q}})^{2}(z \! - \! 
\alpha_{k})^{2}}, \quad z \! \in \! \mathbb{C} \setminus (J_{f} \cup 
\lbrace \alpha_{p_{1}},\dotsc,\alpha_{p_{\mathfrak{s}-2}},\alpha_{k} 
\rbrace),
\end{equation}
where
\begin{equation}
\tilde{\rho}_{2 \mathfrak{s}-3}(n,k) \! = \! \tilde{w}_{0}(n,k) \! = 
\! \dfrac{1}{\pi^{2}} \left(\dfrac{(n \! - \! 1)K \! + \! k}{n} \right) 
\int_{J_{f}} \widetilde{V}^{\prime}(\xi) \psi_{\widetilde{V}}^{f}(\xi) 
\, \md \xi, \label{eql3.7a7}
\end{equation}
\begin{align}
\tilde{\rho}_{m}(n,k) =& \, \tilde{w}_{2 \mathfrak{s}-m-3}(n,k) \! + \! 
\tilde{c}_{m}^{\natural}(n,k) \! + \! \sum_{j=m}^{2(\mathfrak{s}-2)} 
\left(\tilde{w}_{j-m}(n,k) \tilde{c}_{j+1}(n,k) \! + \! \tilde{\nu}_{j-m}
(n,k) \tilde{c}_{j+1}^{\flat}(n,k) \right. \nonumber \\
+&\left. \, \tilde{u}_{j-m}(n,k) \tilde{c}_{j+2}^{\sharp}(n,k) \right), 
\quad m \! = \! 0,1,\dotsc,2(\mathfrak{s} \! - \! 2), \label{eql3.7a8}
\end{align}
with $\tilde{w}_{r}(n,k)$, $\tilde{c}_{l}(n,k)$, $\tilde{u}_{r}
(n,k)$, $\tilde{\nu}_{r}(n,k)$, $\tilde{c}_{l}^{\sharp}(n,k)$, 
$\tilde{c}_{l}^{\flat}(n,k)$, and $\tilde{c}_{l}^{\natural}(n,k)$ 
given in Equations~\eqref{eql3.7v}, \eqref{eql3.7w}, \eqref{eql3.7y}, 
\eqref{eql3.7z}, \eqref{eql3.7a1}, \eqref{eql3.7a3}, and~\eqref{eql3.7a5}, 
respectively. Via Equation~\eqref{eql3.4a}, one writes
\begin{align*}
\dfrac{1}{\pi^{2}} \left(\dfrac{(n \! - \! 1)K \! + \! k}{n} \right) 
\int_{J_{f}} \dfrac{\widetilde{V}^{\prime}(\xi) \psi_{\widetilde{V}}^{f}
(\xi)}{\xi \! - \! z} \, \md \xi =& \, \dfrac{1}{\pi^{2}} \left(\dfrac{(n \! 
- \! 1)K \! + \! k}{n} \right) \int_{J_{f}} \dfrac{(\widetilde{V}^{\prime}
(\xi) \! - \! \widetilde{V}^{\prime}(z)) \psi^{f}_{\widetilde{V}}(\xi)}{\xi 
\! - \! z} \, \md \xi \\
-& \, \dfrac{\mi}{\pi} \widetilde{V}^{\prime}(z) \mathfrak{F}_{f}(z) \! - \! 
\dfrac{\widetilde{V}^{\prime}(z)}{\pi^{2}} \left(\dfrac{(\varkappa_{nk} \! 
- \! 1)}{n(z \! - \! \alpha_{k})} \! + \! \sum_{q=1}^{\mathfrak{s}-2} 
\dfrac{\varkappa_{nk \tilde{k}_{q}}}{n(z \! - \! \alpha_{p_{q}})} \right):
\end{align*}
substituting the latter relation into Equation~\eqref{eql3.7a6} and 
completing the square, one arrives at, for $n \! \in \! \mathbb{N}$ 
and $k \! \in \! \lbrace 1,2,\dotsc,K \rbrace$ such that 
$\alpha_{p_{\mathfrak{s}}} \! := \! \alpha_{k} \! \neq \! \infty$, 
upon re-arranging terms,
\begin{equation} \label{eql3.7a9} 
\left(\mathfrak{F}_{f}(z) \! + \! \dfrac{\mi \widetilde{V}^{\prime}(z)}{2 \pi} 
\right)^{2} \! + \! \dfrac{\tilde{q}_{\widetilde{V}}(z)}{\pi^{2}} \! = \! 0, 
\quad z \! \in \! \mathbb{C} \setminus (J_{f} \cup \lbrace \alpha_{p_{1}},
\dotsc,\alpha_{p_{\mathfrak{s}-2}},\alpha_{k} \rbrace),
\end{equation}
where
\begin{equation} \label{eqveefin}
\begin{aligned}
\tilde{q}_{\widetilde{V}}(z) :=& \, \left(\dfrac{\widetilde{V}^{\prime}(z)}{2} 
\right)^{2} \! + \! \widetilde{V}^{\prime}(z) \left(\dfrac{(\varkappa_{nk} \! 
- \! 1)}{n(z \! - \! \alpha_{k})} \! + \! \sum_{q=1}^{\mathfrak{s}-2} 
\dfrac{\varkappa_{nk \tilde{k}_{q}}}{n(z \! - \! \alpha_{p_{q}})} \right) \! - 
\! \dfrac{\pi^{2} \sum_{m=0}^{2 \mathfrak{s}-3} \tilde{\rho}_{m}(n,k)
z^{m}}{\prod_{q=1}^{\mathfrak{s}-2}(z \! - \! \alpha_{p_{q}})^{2}
(z \! - \! \alpha_{k})^{2}} \\
-& \, \left(\dfrac{(n \! - \! 1)K \! + \! k}{n} \right) \int_{J_{f}} 
\dfrac{(\widetilde{V}^{\prime}(\xi) \! - \! \widetilde{V}^{\prime}(z)) 
\psi_{\widetilde{V}}^{f}(\xi)}{\xi \! - \! z} \, \md \xi \\
=& \, \left(\dfrac{\widetilde{V}^{\prime}(z)}{2} \right)^{2} \! + \! 
\widetilde{V}^{\prime}(z) \left(\dfrac{(\varkappa_{nk} \! - \! 1)}{n(z \! - 
\! \alpha_{k})} \! + \! \sum_{q=1}^{\mathfrak{s}-2} \dfrac{\varkappa_{nk 
\tilde{k}_{q}}}{n(z \! - \! \alpha_{p_{q}})} \right) \! - \! \dfrac{\pi^{2} 
\sum_{m=0}^{2 \mathfrak{s}-3} \tilde{\rho}_{m}(n,k)z^{m}}{\prod_{q=
1}^{\mathfrak{s}-2}(z \! - \! \alpha_{p_{q}})^{2}(z \! - \! \alpha_{k})^{2}} \\
-& \, \left(\dfrac{(n \! - \! 1)K \! + \! k}{n} \right) \int_{J_{f}} 
\int_{0}^{1} \widetilde{V}^{\prime \prime}(t \xi \! + \! (1 \! - \! t)z) \, 
\md t \, \md \mu_{\widetilde{V}}^{f}(\xi).
\end{aligned}
\end{equation}
(Since $\widetilde{V} \colon \overline{\mathbb{R}} \setminus \lbrace 
\alpha_{1},\alpha_{2},\dotsc,\alpha_{K} \rbrace \! \to \! \mathbb{R}$ 
satisfies conditions~\eqref{eq20}--\eqref{eq22}, it follows {}from the 
identity $x_{1}^{m} \! - \! x_{2}^{m} \! = \! (x_{1} \! - \! x_{2})
(x_{1}^{m-1} \! + \! x_{1}^{m-2}x_{2} \! + \! \cdots \! + \! x_{1}x_{2}^{m-2} 
\! + \! x_{2}^{m-1})$, $m \! \in \! \mathbb{N}$, that, for $n \! \in \! 
\mathbb{N}$ and $k \! \in \! \lbrace 1,2,\dotsc,K \rbrace$ such that 
$\alpha_{p_{\mathfrak{s}}} \! := \! \alpha_{k} \! \neq \! \infty$, the 
function $\tilde{q}_{\widetilde{V}}(z)$ is, in particular, real analytic 
on $\mathbb{R} \setminus \lbrace \alpha_{p_{1}},\dotsc,
\alpha_{p_{\mathfrak{s}-2}},\alpha_{k} \rbrace$.) For $x \! \in \! J_{f}$, 
let $z \! = \! x \! \pm \! \mi \varepsilon$, and consider the $\varepsilon 
\! \downarrow \! 0$ limit of---the resulting---Equation~\eqref{eql3.7a9}: 
$\lim_{\varepsilon \downarrow 0}(\mathfrak{F}_{f}(x \! \pm \! \mi 
\varepsilon) \! + \! \tfrac{\mi \widetilde{V}^{\prime}(x \pm \mi 
\varepsilon)}{2 \pi})^{2} \! = \! (\mathfrak{F}_{f \, \pm}(x) \! + \! 
\tfrac{\mi \widetilde{V}^{\prime}(x)}{2 \pi})^{2}$, as $\widetilde{V}$ 
is real analytic on $J_{f}$; however, {}from the real analyticity of 
$\widetilde{V} \colon \overline{\mathbb{R}} \setminus \lbrace 
\alpha_{1},\alpha_{2},\dotsc,\alpha_{K} \rbrace \! \to \! \mathbb{R}$, 
Equation~\eqref{eql3.7r}, and the first line of Equation~\eqref{eql3.7p}, 
it follows that, for $n \! \in \! \mathbb{N}$ and $k \! \in \! \lbrace 
1,2,\dotsc,K \rbrace$ such that $\alpha_{p_{\mathfrak{s}}} \! := \! 
\alpha_{k} \! \neq \! \infty$,
\begin{equation} \label{eql3.7a10} 
\mathfrak{F}_{f \, \pm}(z) \! = \! \dfrac{\widetilde{V}^{\prime}(z)}{2 
\pi \mi} \! \mp \! \left(\dfrac{(n \! - \! 1)K \! + \! k}{n} \right) 
\psi^{f}_{\widetilde{V}}(z), \quad z \! \in \! J_{f},
\end{equation}
whence $(\mathfrak{F}_{f \, \pm}(x) \! + \! \tfrac{\mi \widetilde{V}^{\prime}
(x)}{2 \pi})^{2} \! = \! (\tfrac{(n-1)K+k}{n})^{2}(\psi^{f}_{\widetilde{V}}
(x))^{2}$, $x \! \in \! J_{f}$, which implies, via Equation~\eqref{eql3.7a9}, 
that, for $n \! \in \! \mathbb{N}$ and $k \! \in \! \lbrace 1,2,\dotsc,K 
\rbrace$ such that $\alpha_{p_{\mathfrak{s}}} \! := \! \alpha_{k} \! \neq 
\! \infty$, $(\pi (\tfrac{(n-1)K+k}{n}))^{2}(\psi^{f}_{\widetilde{V}}(x))^{2} 
\! = \! -\tilde{q}_{\widetilde{V}}(x)$, $x \! \in \! J_{f}$, whereupon, taking 
note of the fact that $\psi^{f}_{\widetilde{V}}(x) \! \geqslant \! 0$, $x \! 
\in \! J_{f}$, it follows that, for $x \! \in \! J_{f}$, $\tilde{q}_{\widetilde{V}}
(x) \! \leqslant \! 0$.\footnote{As a by-product, decomposing the function 
$\tilde{q}_{\widetilde{V}}(x)$, $x \! \in \! J_{f}$, into its positive and 
negative parts, that is, $\tilde{q}_{\widetilde{V}}(x) \! = \! (\tilde{q}_{
\widetilde{V}}(x))^{+} \! - \! (\tilde{q}_{\widetilde{V}}(x))^{-}$, $x \! 
\in \! J_{f}$, where $(\tilde{q}_{\widetilde{V}}(x))^{\pm} \! := \! \max 
\lbrace 0,\pm \tilde{q}_{\widetilde{V}}(x) \rbrace$ $(\geqslant \! 0)$, 
one learns that, for $x \! \in \! J_{f}$, $(\tilde{q}_{\widetilde{V}}(x))^{+} 
\! \equiv \! 0$ and $\psi^{f}_{\widetilde{V}}(x) \! = \! (\pi 
(\tfrac{(n-1)K+k}{n}))^{-1}((\tilde{q}_{\widetilde{V}}(x))^{-})^{1/2}$; as 
$\int_{J_{f}} \psi^{f}_{\widetilde{V}}(\xi) \, \md \xi \! = \! 1$, it follows 
{}from the latter relation that, for $n \! \in \! \mathbb{N}$ and $k \! \in 
\! \lbrace 1,2,\dotsc,K \rbrace$ such that $\alpha_{p_{\mathfrak{s}}} 
\! := \! \alpha_{k} \! \neq \! \infty$, $(\pi (\tfrac{(n-1)K+k}{n}))^{-1} 
\int_{J_{f}}((\tilde{q}_{\widetilde{V}}(\xi))^{-})^{1/2} \, \md \xi \! = \! 1$.} 
For $x \! \notin \! J_{f}$, let $z \! = \! x \! \pm \! \mi \varepsilon$, 
and study, again, the $\varepsilon \! \downarrow \! 0$ limit 
of---the resulting---Equation~\eqref{eql3.7a9}: in this case, however, 
$\lim_{\varepsilon \downarrow 0}(\mathfrak{F}_{f}(x \! \pm \! \mi 
\varepsilon) \! + \! \tfrac{\mi \widetilde{V}^{\prime}(x \pm \mi 
\varepsilon)}{2 \pi})^{2} \! = \! (\mathfrak{F}_{f \, \pm}(x) \! + \! 
\tfrac{\mi \widetilde{V}^{\prime}(x)}{2 \pi})^{2} \! = \! (\mathfrak{F}_{f}
(x) \! + \! \tfrac{\mi \widetilde{V}^{\prime}(x)}{2 \pi})^{2}$; via the 
second line of Equation~\eqref{eql3.7p}, that is, $\mathfrak{F}_{f}(x) 
\! = \! \tfrac{\mi}{\pi}(\tfrac{(\varkappa_{nk}-1)}{n(x-\alpha_{k})} \! + 
\! \sum_{q=1}^{\mathfrak{s}-2} \tfrac{\varkappa_{nk \tilde{k}_{q}}}{n
(x-\alpha_{p_{q}})}) \! - \! \mi (\tfrac{(n-1)K+k}{n})(\mathcal{H} 
\psi^{f}_{\widetilde{V}})(x)$, $x \! \notin \! J_{f}$, it follows {}from 
Equation~\eqref{eql3.7a9} that, for $n \! \in \! \mathbb{N}$ 
and $k \! \in \! \lbrace 1,2,\dotsc,K \rbrace$ such that 
$\alpha_{p_{\mathfrak{s}}} \! := \! \alpha_{k} \! \neq \! \infty$, 
$(\mathfrak{F}_{f}(x) \! + \! \tfrac{\mi \widetilde{V}^{\prime}(x)}{2 
\pi})^{2} \! = \! (\tfrac{1}{\pi}(\tfrac{(\varkappa_{nk}-1)}{n(x-\alpha_{k})} 
\! + \! \sum_{q=1}^{\mathfrak{s}-2} \tfrac{\varkappa_{nk \tilde{k}_{q}}}{n
(x-\alpha_{p_{q}})}) \! - \! (\tfrac{(n-1)K+k}{n})(\mathcal{H} \psi^{f}_{
\widetilde{V}})(x) \! + \! \tfrac{\widetilde{V}^{\prime}(x)}{2 \pi})^{2} 
\! = \! \tilde{q}_{\widetilde{V}}(x)/\pi^{2}$, $x \! \notin \! J_{f}$, 
whence $\tilde{q}_{\widetilde{V}}(x) \! > \! 0$, $x \! \notin \! 
J_{f}$.\footnote{For $x \! = \! \alpha_{p_{\mathfrak{s}}} \! := \! 
\alpha_{k} \! \neq \! \infty$, $k \! \in \! \lbrace 1,2,\dotsc,K \rbrace$, 
this latter relation merely states that $+\infty \! = \! +\infty$.} {}From 
the fact that $\widetilde{V} \colon \overline{\mathbb{R}} \setminus \lbrace 
\alpha_{1},\alpha_{2},\dotsc,\alpha_{K} \rbrace \! \to \! \mathbb{R}$ 
satisfies conditions~\eqref{eq20}--\eqref{eq22} (in particular, $\widetilde{V}$ 
is real analytic on $J_{f}$ and has an analytic extension to the open 
neighbourhood $\mathbb{U}_{f} \! = \! \lbrace \mathstrut z \! \in \! 
\overline{\mathbb{C}} \setminus \lbrace \alpha_{p_{1}},\alpha_{p_{2}},
\dotsc,\alpha_{p_{\mathfrak{s}}} \rbrace; \, \inf_{t \in J_{f}} \lvert z \! - \! 
t \rvert \! < \! r \rbrace$, with $r \! \in \! (0,1)$ chosen small enough 
so that $\mathbb{U}_{f} \cap \lbrace \alpha_{p_{1}},\alpha_{p_{2}},\dotsc,
\alpha_{p_{\mathfrak{s}}} \rbrace \! = \! \varnothing$), the associated 
equilibrium measure, $\mu^{f}_{\widetilde{V}}$, has compact support, 
and that on a proper compact subset of $\mathbb{R}$ a meromorphic/rational 
function changes sign an at most countable number of times, it follows {}from 
the above analysis for the function $\tilde{q}_{\widetilde{V}}(z)$ and mimicking 
a part of the calculations subsumed in the proof of Theorem~1.38 in \cite{a58} 
that, for $n \! \in \! \mathbb{N}$ and $k \! \in \! \lbrace 1,2,\dotsc,K \rbrace$ 
such that $\alpha_{p_{\mathfrak{s}}} \! := \! \alpha_{k} \! \neq \! \infty$, 
$\supp (\mu^{f}_{\widetilde{V}}) \! =: \! J_{f} \! = \! \lbrace \mathstrut 
x \! \in \! \overline{\mathbb{R}} \setminus \lbrace \alpha_{p_{1}},
\alpha_{p_{2}},\dotsc,\alpha_{p_{\mathfrak{s}}} \rbrace; \, \tilde{q}_{
\widetilde{V}}(x) \! \leqslant \! 0 \rbrace$ consists of the disjoint 
union of a finite number $\mathbb{N}_{0} \! \ni \! N$ (see 
Remark~\ref{remen} below) of bounded real (compact) intervals, which 
can be presented as $J_{f} \! = \! \cup_{j=1}^{N+1}[\tilde{b}_{j-1},
\tilde{a}_{j}]$, where, up to permutation, $\lbrace \tilde{b}_{0},
\tilde{a}_{1},\tilde{b}_{1},\tilde{a}_{2},\dotsc,\tilde{b}_{N},\tilde{a}_{N+1} 
\rbrace \! = \! \lbrace \mathstrut x \! \in \! \overline{\mathbb{R}} 
\setminus \lbrace \alpha_{p_{1}},\alpha_{p_{2}},\dotsc,
\alpha_{p_{\mathfrak{s}}} \rbrace; \, \tilde{q}_{\widetilde{V}}(x) \! = \! 0 
\rbrace$,\footnote{Recalling {}from the above analysis that, for $n \! 
\in \! \mathbb{N}$ and $k \! \in \! \lbrace 1,2,\dotsc,K \rbrace$ such 
that $\alpha_{p_{\mathfrak{s}}} \! := \! \alpha_{k} \! \neq \! \infty$, 
$\psi^{f}_{\widetilde{V}}(x) \! = \! (\pi (\tfrac{(n-1)K+k}{n}))^{-1}
((\tilde{q}_{\widetilde{V}}(x))^{-})^{1/2}$, $x \! \in \! J_{f}$, where 
$(\tilde{q}_{\widetilde{V}}(x))^{-} \! := \! \max \lbrace 0,-\tilde{q}_{
\widetilde{V}}(x) \rbrace$ $(\geqslant \! 0)$, it follows that the 
end-points of the intervals, $\lbrace \tilde{b}_{j-1},\tilde{a}_{j} 
\rbrace_{j=1}^{N+1}$, of the support, $J_{f}$, of the associated equilibrium 
measure, $\mu_{\widetilde{V}}^{f}$, are given by the zeros of the function 
$\tilde{q}_{\widetilde{V}}$ with odd multiplicity; in particular, in the present 
situation, due to the regularity of $\widetilde{V}$, they are the simple zeros 
of $\tilde{q}_{\widetilde{V}}$. Furthermore, one notes that $J_{f} \! = \! 
\lbrace \mathstrut x \! \in \! \overline{\mathbb{R}} \setminus \lbrace 
\alpha_{p_{1}},\alpha_{p_{2}},\dotsc,\alpha_{p_{\mathfrak{s}}} \rbrace; 
\, \psi^{f}_{\widetilde{V}}(x) \! \geqslant \! 0 \rbrace$ (resp., 
$\operatorname{int}(J_{f}) \! = \! \lbrace \mathstrut x \! \in \! 
\overline{\mathbb{R}} \setminus \lbrace \alpha_{p_{1}},\alpha_{p_{2}},
\dotsc,\alpha_{p_{\mathfrak{s}}} \rbrace; \, \psi^{f}_{\widetilde{V}}(x) \! 
> \! 0 \rbrace)$.} $[\tilde{b}_{j-1},\tilde{a}_{j}] \cap \lbrace \alpha_{p_{1}},
\alpha_{p_{2}},\dotsc,\alpha_{p_{\mathfrak{s}}} \rbrace \! = \! \varnothing$, 
$j \! = \! 1,2,\dotsc,N \! + \! 1$, $[\tilde{b}_{i-1},\tilde{a}_{i}] \cap 
[\tilde{b}_{j-1},\tilde{a}_{j}] \! = \! \varnothing$, $i \! \neq \! j \! \in \! 
\lbrace 1,2,\dotsc,N \! + \! 1 \rbrace$, and (with enumeration) $-\infty 
\! < \! \tilde{b}_{0} \! < \! \tilde{a}_{1} \! < \! \tilde{b}_{1} \! < \! 
\tilde{a}_{2} \! < \! \dotsb \! < \! \tilde{b}_{N} \! < \! \tilde{a}_{N+1} 
\! < \! +\infty$.\footnote{As a by-product, it follows that 
$\operatorname{meas}(J_{f}) \! = \! \sum_{j=1}^{N+1} \lvert \tilde{b}_{j-1} 
\! - \! \tilde{a}_{j} \rvert$ $(< \! +\infty)$.}

It remains to determine the $2(N \! + \! 1)$ equations satisfied by the 
end-points of the intervals, $\lbrace \tilde{b}_{j-1},\tilde{a}_{j} \rbrace_{j=
1}^{N+1}$, of the support, $J_{f}$, of the associated equilibrium measure, 
$\mu^{f}_{\widetilde{V}}$. {}From Equation~\eqref{eql3.4a} and the fact 
that $\mu^{f}_{\widetilde{V}} \! \in \! \mathscr{M}_{1}(\mathbb{R})$, it 
follows that, for $n \! \in \! \mathbb{N}$ and $k \! \in \! \lbrace 1,2,
\dotsc,K \rbrace$ such that $\alpha_{p_{\mathfrak{s}}} \! := \! \alpha_{k} 
\! \neq \! \infty$: (i) for $\xi \! \in \! J_{f}$ and $z \! \notin \! J_{f}$ such 
that $\vert \xi/z \vert \! \ll \! 1$ (e.g., $\lvert z \rvert \! \gg \! \max 
\lbrace \max_{i \neq j \in \lbrace 1,\dotsc,\mathfrak{s}-2,\mathfrak{s} 
\rbrace} \lbrace \lvert \alpha_{p_{i}} \! - \! \alpha_{p_{j}} \rvert \rbrace,
\max_{q=1,\dotsc,\mathfrak{s}-2,\mathfrak{s}} \lbrace \lvert 
\alpha_{p_{q}} \rvert \rbrace,\max_{j=1,2,\dotsc,N+1} \lbrace \lvert 
\tilde{b}_{j-1} \! - \! \tilde{a}_{j} \rvert \rbrace \rbrace)$, via the 
expansion $\tfrac{1}{\xi -z} \! = \! -\sum_{j=0}^{l} \tfrac{\xi^{j}}{z^{j+1}} 
\! + \! \tfrac{\xi^{l+1}}{z^{l+1}(\xi -z)}$, $l \! \in \! \mathbb{N}_{0}$,
\begin{equation} \label{eql3.7a11} 
\mathfrak{F}_{f}(z) \underset{z \to \alpha_{p_{\mathfrak{s}-1}} = \infty}{=} 
\dfrac{(\varkappa^{\infty}_{nk \tilde{k}_{\mathfrak{s}-1}} \! + \! 1)}{\mi 
\pi nz} \! - \! \dfrac{1}{\mi \pi z} \sum_{m=1}^{\infty} \left(\left(\dfrac{
\varkappa_{nk} \! - \! 1}{n} \right)(\alpha_{k})^{m} \! + \! \sum_{q=1}^{
\mathfrak{s}-2} \dfrac{\varkappa_{nk \tilde{k}_{q}}}{n}(\alpha_{p_{q}})^{m} 
\! - \! \left(\dfrac{(n \! - \! 1)K \! + \! k}{n} \right) \int_{J_{f}} \xi^{m} 
\psi^{f}_{\widetilde{V}}(\xi) \, \md \xi \right)z^{-m};
\end{equation}
(ii) for $\xi \! \in \! J_{f}$ and $z \! \notin \! J_{f}$ such that $\lvert 
(z \! - \! \alpha_{k})/(\xi \! - \! \alpha_{k}) \rvert \! \ll \! 1$ (e.g., $0 
\! < \! \lvert z \! - \! \alpha_{k} \rvert \! \ll \! \min \lbrace \min_{i 
\neq j \in \lbrace 1,\dotsc,\mathfrak{s}-2,\mathfrak{s} \rbrace} 
\lbrace \lvert \alpha_{p_{i}} \! - \! \alpha_{p_{j}} \rvert \rbrace,
\inf_{\xi \in J_{f}} \lbrace \lvert \xi \! - \! \alpha_{k} \rvert \rbrace,
\min_{j=1,2,\dotsc,N+1} \lbrace \lvert \lvert \tilde{b}_{j-1} \! - \! 
\tilde{a}_{j} \rvert \! - \! \alpha_{k} \rvert \rbrace \rbrace)$, via 
the expansion $\tfrac{1}{(z-\alpha_{k})-(\xi -\alpha_{k})} \! = \! 
-\sum_{j=0}^{l} \tfrac{(z-\alpha_{k})^{j}}{(\xi -\alpha_{k})^{j+1}} \! 
+ \! \tfrac{(z-\alpha_{k})^{l+1}}{(\xi -\alpha_{k})^{l+1}(z-\xi)}$, 
$l \! \in \! \mathbb{N}_{0}$,
\begin{equation*}
\mathfrak{F}_{f}(z) \! + \! \dfrac{1}{\mi \pi} \dfrac{(\varkappa_{nk} \! 
- \! 1)}{n(z-\alpha_{k})} \underset{z \to \alpha_{k}}{=} \mathcal{O}(1);
\end{equation*}
and (iii) for $\xi \! \in \! J_{f}$ and $z \! \notin \! J_{f}$ such that 
$\lvert (z \! - \! \alpha_{p_{q}})/(\xi \! - \! \alpha_{p_{q}}) \rvert \! 
\ll \! 1$, $q \! = \! 1,2,\dotsc,\mathfrak{s} \! - \! 2$ (e.g., $0 \! < 
\! \lvert z \! - \! \alpha_{p_{q}} \rvert \! \ll \! \min \lbrace \min_{i 
\neq j \in \lbrace 1,\dotsc,\mathfrak{s}-2,\mathfrak{s} \rbrace} 
\lbrace \lvert \alpha_{p_{i}} \! - \! \alpha_{p_{j}} \rvert \rbrace,
\inf_{\underset{q^{\prime}=1,2,\dotsc,\mathfrak{s}-2}{\xi \in J_{f}}} 
\lbrace \lvert \xi \! - \! \alpha_{p_{q^{\prime}}} \rvert \rbrace,
\min_{\underset{q^{\prime}=1,2,\dotsc,\mathfrak{s}-2}{j=1,2,\dotsc,
N+1}} \lbrace \lvert \lvert \tilde{b}_{j-1} \! - \! \tilde{a}_{j} \rvert \! 
- \! \alpha_{p_{q^{\prime}}} \rvert \rbrace \rbrace$, $q \! = \! 1,2,
\dotsc,\mathfrak{s} \! - \! 2$), via the expansion $\tfrac{1}{(z-
\alpha_{p_{q}})-(\xi -\alpha_{p_{q}})} \! = \! -\sum_{j=0}^{l} 
\tfrac{(z-\alpha_{p_{q}})^{j}}{(\xi -\alpha_{p_{q}})^{j+1}} \! + \! 
\tfrac{(z-\alpha_{p_{q}})^{l+1}}{(\xi -\alpha_{p_{q}})^{l+1}(z-\xi)}$, 
$l \! \in \! \mathbb{N}_{0}$,
\begin{equation*}
\mathfrak{F}_{f}(z) \! + \! \dfrac{1}{\mi \pi} \dfrac{\varkappa_{nk 
\tilde{k}_{q}}}{n(z-\alpha_{p_{q}})} \underset{z \to \alpha_{p_{q}}}{=} 
\mathcal{O}(1), \quad q \! = \! 1,2,\dotsc,\mathfrak{s} \! - \! 2. 
\end{equation*}
{}From the above asymptotic expansions and Equations~\eqref{eql3.7p}, 
\eqref{eql3.7r}, and~\eqref{eql3.7a10}, one deduces that $\mathfrak{F}_{f} 
\colon \mathbb{N} \times \lbrace 1,2,\dotsc,K \rbrace \linebreak[4] 
\times \mathbb{C} \setminus (J_{f} \cup \lbrace \alpha_{p_{1}},
\dotsc,\alpha_{p_{\mathfrak{s}-2}},\alpha_{k} \rbrace) \! \to \! \mathbb{C}$ 
solves the following scalar RHP: (1) $\mathfrak{F}_{f}(z)$ is analytic (resp., 
meromorphic) for $z \! \in \! \mathbb{C} \setminus (J_{f} \cup \lbrace 
\alpha_{p_{1}},\dotsc,\alpha_{p_{\mathfrak{s}-2}},\alpha_{k} \rbrace)$ 
(resp., $z \! \in \! \mathbb{C} \setminus J_{f})$; (2) $\mathfrak{F}_{f \, \pm}
(z) \! := \! \lim_{\varepsilon \downarrow 0} \mathfrak{F}_{f}(z \! \pm \! 
\mi \varepsilon)$ satisfy the boundary conditions $\mathfrak{F}_{f \, +}
(z) \! + \! \mathfrak{F}_{f \, -}(z) \! = \! \widetilde{V}^{\prime}(z)/\mi \pi$, 
$z \! \in \! J_{f}$, and $\mathfrak{F}_{f \, +}(z) \! = \! \mathfrak{F}_{f \, -}
(z) \! = \! \mathfrak{F}_{f}(z)$, $z \! \notin \! J_{f}$; (3) $\mathfrak{F}_{f}(z) 
\! =_{z \to \alpha_{p_{\mathfrak{s}-1}} = \infty} \! (\varkappa^{\infty}_{nk 
\tilde{k}_{\mathfrak{s}-1}} \! + \! 1)/\mi \pi nz \! + \! \mathcal{O}(z^{-2})$; 
and (4) $\operatorname{Res}(\mathfrak{F}_{f}(z);\alpha_{k}) \! = \! 
-(\varkappa_{nk} \! - \! 1)/\mi \pi n$, and $\operatorname{Res}
(\mathfrak{F}_{f}(z);\alpha_{p_{q}}) \! = \! -\varkappa_{nk \tilde{k}_{q}}
/\mi \pi n$, $q \! = \! 1,2,\dotsc,\mathfrak{s} \! - \! 2$. For $n \! \in \! 
\mathbb{N}$ and $k \! \in \! \lbrace 1,2,\dotsc,K \rbrace$ such that 
$\alpha_{p_{\mathfrak{s}}} \! := \! \alpha_{k} \! \neq \! \infty$, a 
representation for the solution of this scalar RHP is given by (see, for 
example, Chapter~VI of \cite{Ga})
\begin{align} \label{eql3.7a12} 
\mathfrak{F}_{f}(z) =& \, -\dfrac{1}{\mi \pi} \left(\dfrac{(\varkappa_{nk} 
\! - \! 1)}{n(z \! - \! \alpha_{k})} \! + \! \sum_{q=1}^{\mathfrak{s}-2} 
\dfrac{\varkappa_{nk \tilde{k}_{q}}}{n(z \! - \! \alpha_{p_{q}})} \right) \! 
+ \! (\tilde{R}(z))^{1/2} \int_{J_{f}} \left(\dfrac{2}{\mi \pi} \left(
\dfrac{(\varkappa_{nk} \! - \! 1)}{n(\xi \! - \! \alpha_{k})} \! + \! 
\sum_{q=1}^{\mathfrak{s}-2} \dfrac{\varkappa_{nk \tilde{k}_{q}}}{n(\xi \! - 
\! \alpha_{p_{q}})} \right) \! + \! \dfrac{\widetilde{V}^{\prime}(\xi)}{\mi 
\pi} \right) \nonumber \\
\times& \, \dfrac{(\tilde{R}(\xi))^{-1/2}_{+}}{\xi \! - \! z} \, 
\dfrac{\md \xi}{2 \pi \mi}, \quad z \! \in \! \mathbb{C} \setminus (J_{f} 
\cup \lbrace \alpha_{p_{1}},\dotsc,\alpha_{p_{\mathfrak{s}-2}},\alpha_{k} 
\rbrace),
\end{align}
where $(\tilde{R}(z))^{1/2}$ is defined by Equation~\eqref{eql3.7j}, 
$(\tilde{R}(z))^{1/2}_{\pm} \! := \! \lim_{\varepsilon \downarrow 0}
(\tilde{R}(z \! \pm \! \mi \varepsilon))^{1/2}$, and the branch of the 
square root is chosen so that $z^{-(N+1)}(\tilde{R}(z))^{1/2} \! 
\sim_{\overline{\mathbb{C}}_{\pm} \ni z \to \alpha_{p_{\mathfrak{s}-1}} 
= \infty} \! \pm 1$. One shows {}from the integral 
representation~\eqref{eql3.7a12} that, for $n \! \in \! \mathbb{N}$ 
and $k \! \in \! \lbrace 1,2,\dotsc,K \rbrace$ such that 
$\alpha_{p_{\mathfrak{s}}} \! := \! \alpha_{k} \! \neq \! \infty$, for 
$\xi \! \in \! J_{f}$ and $z \! \notin \! J_{f}$ such that $\vert \xi/z 
\vert \! \ll \! 1$ (e.g., $\lvert z \rvert \! \gg \! \max \lbrace \max_{i 
\neq j \in \lbrace 1,\dotsc,\mathfrak{s}-2,\mathfrak{s} \rbrace} \lbrace 
\lvert \alpha_{p_{i}} \! - \! \alpha_{p_{j}} \rvert \rbrace,\max_{q=1,\dotsc,
\mathfrak{s}-2,\mathfrak{s}} \lbrace \lvert \alpha_{p_{q}} \rvert \rbrace,
\max_{j=1,2,\dotsc,N+1} \lbrace \lvert \tilde{b}_{j-1} \! - \! \tilde{a}_{j} 
\rvert \rbrace \rbrace)$, via the expansion $\tfrac{1}{\xi -z} \! = \! 
-\sum_{j=0}^{l} \tfrac{\xi^{j}}{z^{j+1}} \! + \! \tfrac{\xi^{l+1}}{z^{l+1}
(\xi -z)}$, $l \! \in \! \mathbb{N}_{0}$,
\begin{align} \label{eql3.7a13} 
\mathfrak{F}_{f}(z) \underset{z \to \alpha_{p_{\mathfrak{s}-1}} = 
\infty}{=}& \, -\dfrac{(z^{N+1} \! + \! \dotsb)}{2 \pi \mi z} \int_{J_{f}} 
\dfrac{1}{(\tilde{R}(\xi))^{1/2}_{+}} \left(\dfrac{2}{\mi \pi} \left(
\dfrac{(\varkappa_{nk} \! - \! 1)}{n(\xi \! - \! \alpha_{k})} \! + \! 
\sum_{q=1}^{\mathfrak{s}-2} \dfrac{\varkappa_{nk \tilde{k}_{q}}}{n
(\xi \! - \! \alpha_{p_{q}})} \right) \! + \! \dfrac{\widetilde{V}^{\prime}
(\xi)}{\mi \pi} \right) \left(1 \! + \! \dfrac{\xi}{z} \! + \! 
\dfrac{\xi^{2}}{z^{2}} \! + \! \dotsb \! + \! \dfrac{\xi^{N}}{z^{N}} \right) 
\md \xi \nonumber \\
-& \, \dfrac{1}{z} \dfrac{1}{\mi \pi} \left(\dfrac{\varkappa_{nk} \! - \! 
1}{n} \right) \! - \! \dfrac{1}{z} \sum_{q=1}^{\mathfrak{s}-2} \dfrac{
\varkappa_{nk \tilde{k}_{q}}}{\mi \pi n} \! - \! \dfrac{(1 \! + \! \dotsb)}{2 
\pi \mi z} \int_{J_{f}} \dfrac{\xi^{N+1}}{(\tilde{R}(\xi))^{1/2}_{+}} \left(
\dfrac{2}{\mi \pi} \left(\dfrac{(\varkappa_{nk} \! - \! 1)}{n(\xi \! - \! 
\alpha_{k})} \! + \! \sum_{q=1}^{\mathfrak{s}-2} \dfrac{\varkappa_{nk 
\tilde{k}_{q}}}{n(\xi \! - \! \alpha_{p_{q}})} \right) \! + \! \dfrac{
\widetilde{V}^{\prime}(\xi)}{\mi \pi} \right) \md \xi \nonumber \\
-& \, \dfrac{(1 \! + \! \dotsb)}{2 \pi \mi z^{2}} \int_{J_{f}} 
\dfrac{\xi^{N+2}}{(\tilde{R}(\xi))^{1/2}_{+}} \left(\dfrac{2}{\mi \pi} 
\left(\dfrac{(\varkappa_{nk} \! - \! 1)}{n(\xi \! - \! \alpha_{k})} \! + \! 
\sum_{q=1}^{\mathfrak{s}-2} \dfrac{\varkappa_{nk \tilde{k}_{q}}}{n(\xi \! - \! 
\alpha_{p_{q}})} \right) \! + \! \dfrac{\widetilde{V}^{\prime}(\xi)}{\mi \pi} 
\right) \left(1 \! + \! \dfrac{\xi}{z} \! + \! \dotsb \right) \md \xi 
\nonumber \\
-& \, \dfrac{1}{\pi \mi z} \sum_{m=1}^{\infty} \left(\left(\dfrac{
\varkappa_{nk} \! - \! 1}{n} \right)(\alpha_{k})^{m} \! + \! \sum_{q=1}^{
\mathfrak{s}-2} \dfrac{\varkappa_{nk \tilde{k}_{q}}}{n}(\alpha_{p_{q}})^{m} 
\right)z^{-m}:
\end{align}
via the asymptotic expansions~\eqref{eql3.7a11} and~\eqref{eql3.7a13}, 
one arrives at, for $n \! \in \! \mathbb{N}$ and $k \! \in \! \lbrace 1,2,
\dotsc,K \rbrace$ such that $\alpha_{p_{\mathfrak{s}}} \! := \! \alpha_{k} 
\! \neq \! \infty$, upon removing secular terms, the $N \! + \! 1$ real 
moment equations~\eqref{eql3.7g}, and, upon equating $z^{-1}$ terms, 
the real moment equation~\eqref{eql3.7h}. It remains to determine the 
remaining $2(N \! + \! 1) \! - \! (N \! + \! 1) \! - \! 1 \! = \! N$ real 
moment equations. {}From Equation~\eqref{eql3.4a} and the second line 
of Equation~\eqref{eql3.7p}, it follows that, for $z \! \notin \! J_{f}$,
\begin{equation*}
\mathfrak{F}_{f}(z) \! + \! \dfrac{1}{\mi \pi} \left(\dfrac{(\varkappa_{nk} 
\! - \! 1)}{n(z \! - \! \alpha_{k})} \! + \! \sum_{q=1}^{\mathfrak{s}-2} 
\dfrac{\varkappa_{nk \tilde{k}_{q}}}{n(z \! - \! \alpha_{p_{q}})} \right) 
\! = \! -\mi \left(\dfrac{(n \! - \! 1)K \! + \! k}{n} \right)(\mathcal{H} 
\psi^{f}_{\widetilde{V}})(z),
\end{equation*}
which, in conjunction with the integral representation~\eqref{eql3.7a12}, 
gives rise to, for $z \! \notin \! J_{f}$,
\begin{equation} \label{eql3.7a14} 
(\mathcal{H} \psi^{f}_{\widetilde{V}})(z) \! = \! \mi \left(\dfrac{(n \! - 
\! 1)K \! + \! k}{n} \right)^{-1}(\tilde{R}(z))^{1/2} \int_{J_{f}} \left(
\dfrac{2}{\mi \pi} \left(\dfrac{(\varkappa_{nk} \! - \! 1)}{n(\xi \! - \! 
\alpha_{k})} \! + \! \sum_{q=1}^{\mathfrak{s}-2} \dfrac{\varkappa_{nk 
\tilde{k}_{q}}}{n(\xi \! - \! \alpha_{p_{q}})} \right) \! + \! 
\dfrac{\widetilde{V}^{\prime}(\xi)}{\mi \pi} \right) 
\dfrac{(\tilde{R}(\xi))^{-1/2}_{+}}{\xi \! - \! z} \, 
\dfrac{\md \xi}{2 \pi \mi}.
\end{equation}
A contour integration argument shows that (cf. Equation~\eqref{eql3.7r}), 
for $j \! = \! 1,2,\dotsc,N$,
\begin{equation} \label{eql3.7a15} 
\int_{\tilde{a}_{j}}^{\tilde{b}_{j}} \left((\mathcal{H} 
\psi^{f}_{\widetilde{V}})(\xi) \! - \! \dfrac{1}{2 \pi} \left(\dfrac{(n \! - 
\! 1)K \! + \! k}{n} \right)^{-1} \left(\dfrac{2(\varkappa_{nk} \! - \! 1)}{n
(\xi \! - \! \alpha_{k})} \! + \! 2 \sum_{q=1}^{\mathfrak{s}-2} \dfrac{
\varkappa_{nk \tilde{k}_{q}}}{n(\xi \! - \! \alpha_{p_{q}})} \! + \! 
\widetilde{V}^{\prime}(\xi) \right) \right) \md \xi \! = \! 0;
\end{equation}
hence, via Equations~\eqref{eql3.7a14} and~\eqref{eql3.7a15}, it follows 
that, for $n \! \in \! \mathbb{N}$ and $k \! \in \! \lbrace 1,2,\dotsc,
K \rbrace$ such that $\alpha_{p_{\mathfrak{s}}} \! := \! \alpha_{k} 
\! \neq \! \infty$,
\begin{align*}
&\int_{\tilde{a}_{j}}^{\tilde{b}_{j}}(\tilde{R}(\varsigma))^{1/2} \left(
\int_{J_{f}} \left(\dfrac{2}{\mi \pi} \left(\dfrac{(\varkappa_{nk} \! - \! 
1)}{n(\xi \! - \! \alpha_{k})} \! + \! \sum_{q=1}^{\mathfrak{s}-2} \dfrac{
\varkappa_{nk \tilde{k}_{q}}}{n(\xi \! - \! \alpha_{p_{q}})} \right) \! + 
\! \dfrac{\widetilde{V}^{\prime}(\xi)}{\mi \pi} \right) \dfrac{(\tilde{R}
(\xi))^{-1/2}_{+}}{\xi \! - \! \varsigma} \, \md \xi \right) \md \varsigma \\
=& \, \int_{\tilde{a}_{j}}^{\tilde{b}_{j}} \left(\dfrac{2(\varkappa_{nk} \! 
- \! 1)}{n(\xi \! - \! \alpha_{k})} \! + \! 2 \sum_{q=1}^{\mathfrak{s}-2} 
\dfrac{\varkappa_{nk \tilde{k}_{q}}}{n(\xi \! - \! \alpha_{p_{q}})} \! + \! 
\widetilde{V}^{\prime}(\xi) \right) \md \xi, \quad j \! = \! 1,2,\dotsc,N,
\end{align*}
which, via an integration argument,\footnote{If no $\alpha_{
p_{\mathfrak{s}}} \! := \! \alpha_{k} \! \neq \! \infty$, $k \! \in \! 
\lbrace 1,2,\dotsc,K \rbrace$, belongs to a `gap' $(\tilde{a}_{j},
\tilde{b}_{j})$, $j \! = \! 1,2,\dotsc,N$, then, via the real analyticity of 
$\widetilde{V} \colon \overline{\mathbb{R}} \setminus \lbrace \alpha_{1},
\alpha_{2},\dotsc,\alpha_{K} \rbrace \! \to \! \mathbb{R}$ satisfying 
conditions~\eqref{eq20}--\eqref{eq22} and an application of the 
Fundamental Theorem of Calculus, the integrals on the right-hand side 
are easily evaluated; if, however, an $\alpha_{p_{\mathfrak{s}}} \! := \! 
\alpha_{k} \! \neq \! \infty$, $k \! \in \! \lbrace 1,2,\dotsc,K \rbrace$, 
belongs to a gap $(\tilde{a}_{j},\tilde{b}_{j})$, $j \! = \! 1,2,\dotsc,N$, 
then the corresponding integrals on the right-hand side need to be 
evaluated in the Cauchy principal value sense; in either case, one arrives 
at, for $n \! \in \! \mathbb{N}$ and $k \! \in \! \lbrace 1,2,\dotsc,K 
\rbrace$ such that $\alpha_{p_{\mathfrak{s}}} \! := \! \alpha_{k} \! 
\neq \! \infty$, the $N$ real moment equations~\eqref{eql3.7i}.} 
gives rise to the remaining $N$ real moment 
equations~\eqref{eql3.7i}.\footnote{As $J_{f} \cap \lbrace \alpha_{p_{1}},
\alpha_{p_{2}},\dotsc,\alpha_{p_{\mathfrak{s}}} \rbrace \! = \! \varnothing$, 
$\widetilde{V} \colon \overline{\mathbb{R}} \setminus \lbrace \alpha_{1},
\alpha_{2},\dotsc,\alpha_{K} \rbrace \! \to \! \mathbb{R}$ satisfying 
conditions~\eqref{eq20}--\eqref{eq22} is real analytic on $J_{f}$, 
$(\tilde{R}(x))^{1/2} \! =_{x \downarrow \tilde{b}_{j-1}} \! \mathcal{O}
((x \! - \! \tilde{b}_{j-1})^{1/2})$ and $(\tilde{R}(x))^{1/2} \! 
=_{x \uparrow a_{j}} \! \mathcal{O}((\tilde{a}_{j} \! - \! x)^{1/2})$, 
$j \! = \! 1,2,\dotsc,N \! + \! 1$, it follows that the integrals 
constituting the system of $2(N \! + \! 1)$ real moment 
equations~\eqref{eql3.7g}--\eqref{eql3.7i} for the end-points of 
the intervals, $\lbrace \tilde{b}_{j-1},\tilde{a}_{j} \rbrace_{j=1}^{N+1}$, 
of the support, $J_{f}$, of the associated equilibrium measure, 
$\mu^{f}_{\widetilde{V}}$, have integrable singularities at 
$\tilde{b}_{j-1},\tilde{a}_{j}$, $j \! = \! 1,2,\dotsc,N \! + \! 1$.}

{}From the integral representation~\eqref{eql3.7a12}, a residue calculation 
shows that, for $n \! \in \! \mathbb{N}$ and $k \! \in \! \lbrace 1,2,\dotsc,
K \rbrace$ such that $\alpha_{p_{\mathfrak{s}}} \! := \! \alpha_{k} \! \neq 
\! \infty$, for non-real $z$ in the domain of analyticity of $\widetilde{V}$,
\begin{equation} \label{eql3.7a30} 
\mathfrak{F}_{f}(z) \! = \! \dfrac{\widetilde{V}^{\prime}(z)}{2 \pi \mi} \! 
+ \! \dfrac{1}{2}(\tilde{R}(z))^{1/2} \oint_{\tilde{C}_{\widetilde{V}}} 
\left(\dfrac{2}{\mi \pi} \left(\dfrac{(\varkappa_{nk} \! - \! 1)}{n(\xi \! 
- \! \alpha_{k})} \! + \! \sum_{q=1}^{\mathfrak{s}-2} \dfrac{\varkappa_{nk 
\tilde{k}_{q}}}{n(\xi \! - \! \alpha_{p_{q}})} \right) \! + \! 
\dfrac{\widetilde{V}^{\prime}(\xi)}{\mi \pi} \right) 
\dfrac{(\tilde{R}(\xi))^{-1/2}}{\xi \! - \! z} \, \dfrac{\md \xi}{2 \pi \mi},
\end{equation}
where the contour $\tilde{C}_{\widetilde{V}}$ is described in item~$\pmb{(2)}$ 
of the lemma: it follows {}from Equation~\eqref{eql3.7a30} and the real 
analyticity of $\widetilde{V} \colon \overline{\mathbb{R}} \setminus \lbrace 
\alpha_{1},\alpha_{2},\dotsc,\alpha_{K} \rbrace \! \to \! \mathbb{R}$ 
satisfying conditions~\eqref{eq20}--\eqref{eq22} that, for $z \! \in \! J_{f}$,
\begin{equation} \label{eql3.7a16} 
\mathfrak{F}_{f \, \pm}(z) \! = \! \dfrac{\widetilde{V}^{\prime}(z)}{2 
\pi \mi} \! + \! \dfrac{1}{2}(\tilde{R}(z))^{1/2}_{\pm} \oint_{\tilde{C}_{
\widetilde{V}}} \left(\dfrac{2}{\mi \pi} \left(\dfrac{(\varkappa_{nk} \! - 
\! 1)}{n(\xi \! - \! \alpha_{k})} \! + \! \sum_{q=1}^{\mathfrak{s}-2} 
\dfrac{\varkappa_{nk \tilde{k}_{q}}}{n(\xi \! - \! \alpha_{p_{q}})} \right) 
\! + \! \dfrac{\widetilde{V}^{\prime}(\xi)}{\mi \pi} \right) \dfrac{(\tilde{R}
(\xi))^{-1/2}}{\xi \! - \! z} \, \dfrac{\md \xi}{2 \pi \mi};
\end{equation}
hence, via Equations~\eqref{eql3.7a10} and~\eqref{eql3.7a16}, one 
arrives at, for $n \! \in \! \mathbb{N}$ and $k \! \in \! \lbrace 1,2,\dotsc,
K \rbrace$ such that $\alpha_{p_{\mathfrak{s}}} \! := \! \alpha_{k} \! 
\neq \! \infty$, the representation $\psi^{f}_{\widetilde{V}}(x) \! = \! 
(2 \pi \mi)^{-1}(\tilde{R}(x))^{1/2}_{+} \tilde{h}_{\widetilde{V}}(x) 
\chi_{J_{f}}(x)$, where $\tilde{h}_{\widetilde{V}}(z)$ is defined by 
Equation~\eqref{eql3.7l} (the integral representation~\eqref{eql3.7l} 
for $\tilde{h}_{\widetilde{V}}(z)$ shows that it is analytic in some open 
subset of $\mathbb{C} \setminus \lbrace \alpha_{p_{1}},\dotsc,
\alpha_{p_{\mathfrak{s}-2}},\alpha_{k} \rbrace$ containing $J_{f}$), 
and $\chi_{J_{f}}(x)$ is the characteristic function of the compact set 
$J_{f}$, which gives rise to the formula~\eqref{eql3.7k} for the density 
of the associated equilibrium measure.

Lastly, it will be shown that, for $n \! \in \! \mathbb{N}$ and $k \! \in 
\! \lbrace 1,2,\dotsc,K \rbrace$ such that $\alpha_{p_{\mathfrak{s}}} 
\! := \! \alpha_{k} \! \neq \! \infty$, if $J_{f} \! := \! \cup_{j=1}^{N+1}
[\tilde{b}_{j-1},\tilde{a}_{j}]$, then the end-points of the intervals, 
$\lbrace \tilde{b}_{j-1},\tilde{a}_{j} \rbrace_{j=1}^{N+1}$, of the support, 
$J_{f}$, of the associated equilibrium measure, $\mu^{f}_{\widetilde{V}}$, 
which satisfy the system of $2(N \! + \! 1)$ real moment 
equations~\eqref{eql3.7g}--\eqref{eql3.7i}, are real-analytic functions 
of $z_{o}$ (see, also, Proposition~6 in \cite{vltot}), thus establishing, in the 
double-scaling limit $\mathscr{N},n \! \to \! \infty$ such that $z_{o} \! = \! 
1 \! + \! o(1)$, the local {}\footnote{In an open (asymptotic) neighbourhood 
$\mathbb{U}_{o(1)}(1) \! := \! \lbrace \mathstrut z_{o} \! \in \! \mathbb{R}_{+}; 
\, \lvert z_{o} \! - \! 1 \rvert \! \lesssim \! o(1) \rbrace$.} solvability 
of the corresponding system of $2(N \! + \! 1)$ real moment 
equations~\eqref{eql3.7g}--\eqref{eql3.7i}. Towards this end, one follows 
closely the idea of the proof (but adapted to the present situation) of 
Theorem~1.3~(iii) in \cite{a57} (see, also, Section~8 of \cite{a58}).

It was shown in the above analysis that, for $n \! \in \! \mathbb{N}$ 
and $k \! \in \! \lbrace 1,2,\dotsc,K \rbrace$ such that 
$\alpha_{p_{\mathfrak{s}}} \! := \! \alpha_{k} \! \neq \! \infty$, the 
end-points of the intervals, $\lbrace \tilde{b}_{j-1},\tilde{a}_{j} 
\rbrace_{j=1}^{N+1}$, of the support, $J_{f}$, of the associated 
equilibrium measure, $\mu^{f}_{\widetilde{V}}$, are the simple 
zeros {}\footnote{These are the only zeros for the case of regular 
$\widetilde{V} \colon \overline{\mathbb{R}} \setminus \lbrace 
\alpha_{1},\alpha_{2},\dotsc,\alpha_{K} \rbrace \! \to \! \mathbb{R}$ 
satisfying conditions~\eqref{eq20}--\eqref{eq22} considered/studied in this 
monograph.} of the meromorphic function (cf. Equation~\eqref{eqveefin}) 
$\tilde{q}_{\widetilde{V}}(z)$, that is, with the enumeration $-\infty \! 
< \! \tilde{b}_{0} \! < \! \tilde{a}_{1} \! < \! \tilde{b}_{1} \! < \! 
\tilde{a}_{2} \! < \! \dotsb \! < \! \tilde{b}_{N} \! < \! \tilde{a}_{N+1} \! < 
\! +\infty$, $\lbrace \tilde{b}_{0},\tilde{a}_{1},\tilde{b}_{1},\tilde{a}_{2},
\dotsc,\tilde{b}_{N},\tilde{a}_{N+1} \rbrace \! = \! \lbrace \mathstrut 
x \! \in \! \mathbb{R} \setminus \lbrace \alpha_{p_{1}},\dotsc,
\alpha_{p_{\mathfrak{s}-2}},\alpha_{k} \rbrace; \, \tilde{q}_{\widetilde{V}}
(x) \! = \! 0 \rbrace$. The function $\tilde{q}_{\widetilde{V}}(z)$ is real 
rational (resp., real analytic) on $\mathbb{R}$ (resp., $\mathbb{R} 
\setminus \lbrace \alpha_{p_{1}},\dotsc,\alpha_{p_{\mathfrak{s}-2}},
\alpha_{k} \rbrace)$, it has analytic extension (independent of $z_{o})$ 
to the open neighbourhood $\tilde{\mathbb{U}}_{f} \! = \! \cup_{j=1}^{N+1} 
\tilde{\mathbb{U}}_{j}$, where $\tilde{\mathbb{U}}_{j} \! := \! \lbrace 
\mathstrut z \! \in \! \mathbb{C} \setminus \lbrace \alpha_{p_{1}},\dotsc,
\alpha_{p_{\mathfrak{s}-2}},\alpha_{k} \rbrace; \, \inf_{x \in (\tilde{b}_{j-1},
\tilde{a}_{j})} \lvert z \! - \! x \rvert \! < \! r_{j} \rbrace$, $j \! = \! 1,2,
\dotsc,N \! + \! 1$, with $r_{j} \! \in \! (0,1)$ chosen small enough so that 
$\tilde{\mathbb{U}}_{i} \cap  \tilde{\mathbb{U}}_{j} \! = \! \varnothing$, 
$i \! \neq \! j \! \in \! \lbrace 1,2,\dotsc,N \! + \! 1 \rbrace$, and depends 
continuously on $z_{o}$; thus, its simple zeros, that is, $\tilde{b}_{j-1} \! 
= \! \tilde{b}_{j-1}(z_{o})$ and $\tilde{a}_{j} \! = \! \tilde{a}_{j}(z_{o})$, 
$j \! = \! 1,2,\dotsc,N \! + \! 1$, are continuous functions of 
$z_{o}$.\footnote{The parametric dependence of the associated 
end-points of the intervals on the elements of the corresponding pole 
set $\lbrace \alpha_{p_{1}},\dotsc,\alpha_{p_{\mathfrak{s}-2}},\alpha_{k} 
\rbrace$ is not considered in this monograph.} For $n \! \in \! 
\mathbb{N}$ and $k \! \in \! \lbrace 1,2,\dotsc,K \rbrace$ such that 
$\alpha_{p_{\mathfrak{s}}} \! := \! \alpha_{k} \! \neq \! \infty$, the 
large-$z$ (e.g., $\lvert z \rvert \! \gg \! \max \lbrace \max_{i \neq j 
\in \lbrace 1,\dotsc,\mathfrak{s}-2,\mathfrak{s} \rbrace} \lbrace \lvert 
\alpha_{p_{i}} \! - \! \alpha_{p_{j}} \rvert \rbrace,\max_{q=1,\dotsc,
\mathfrak{s}-2,\mathfrak{s}} \lbrace \lvert \alpha_{p_{q}} \rvert \rbrace,
\max_{j=1,2,\dotsc,N+1} \lbrace \lvert \tilde{b}_{j-1} \! - \! \tilde{a}_{j} 
\rvert \rbrace \rbrace)$ asymptotic expansion for $\mathfrak{F}_{f}(z)$ 
given in Equation~\eqref{eql3.7a12} reads
\begin{align*}
\mathfrak{F}_{f}(z) \underset{z \to \alpha_{p_{\mathfrak{s}-1}} = \infty}{=}& 
\, -\dfrac{1}{\mi \pi z} \left(\left(\dfrac{\varkappa_{nk} \! - \! 1}{n} \right) 
\! + \! \sum_{q=1}^{\mathfrak{s}-2} \dfrac{\varkappa_{nk \tilde{k}_{q}}}{n} 
\right) \! - \! \dfrac{1}{\mi \pi z} \sum_{m=1}^{\infty} \left(\left(
\dfrac{\varkappa_{nk} \! - \! 1}{n} \right)(\alpha_{k})^{m} \! + \! 
\sum_{q=1}^{\mathfrak{s}-2} \dfrac{\varkappa_{nk \tilde{k}_{q}}}{n}
(\alpha_{p_{q}})^{m} \right)z^{-m} \! - \! \dfrac{(\tilde{R}(z))^{1/2}}{2 \pi \mi z} 
\sum_{j=1}^{\infty} \tilde{\mathcal{T}}_{j}z^{-j},
\end{align*}
where
\begin{equation} \label{eql3.7a17} 
\tilde{\mathcal{T}}_{j} \! := \! \int_{J_{f}} \dfrac{\xi^{j}}{(\tilde{R}
(\xi))^{1/2}_{+}} \left(\dfrac{2}{\mi \pi} \left(\dfrac{(\varkappa_{nk} \! 
- \! 1)}{n(\xi  \! - \! \alpha_{k})} \! + \! \sum_{q=1}^{\mathfrak{s}-2} 
\dfrac{\varkappa_{nk \tilde{k}_{q}}}{n(\xi \! - \! \alpha_{p_{q}})} \right) 
\! + \! \dfrac{\widetilde{V}^{\prime}(\xi)}{\mi \pi} \right) \md \xi, \quad 
j \! \in \! \mathbb{N}_{0}.
\end{equation}
For $n \! \in \! \mathbb{N}$ and $k \! \in \! \lbrace 1,2,\dotsc,K \rbrace$ 
such that $\alpha_{p_{\mathfrak{s}}} \! := \! \alpha_{k} \! \neq \! \infty$, 
set (cf. Equation~\eqref{eql3.7a15})\begin{equation} \label{eql3.7a18} 
\tilde{\mathcal{N}}_{j} \! := \! \int_{\tilde{a}_{j}}^{\tilde{b}_{j}} \left(
(\mathcal{H} \psi^{f}_{\widetilde{V}})(\xi) \! - \! \dfrac{1}{2 \pi} \left(
\dfrac{(n \! - \! 1)K \! + \! k}{n} \right)^{-1} \left(\dfrac{2(\varkappa_{nk} 
\! - \! 1)}{n(\xi  \! - \! \alpha_{k})} \! + \! 2 \sum_{q=1}^{\mathfrak{s}-2} 
\dfrac{\varkappa_{nk \tilde{k}_{q}}}{n(\xi \! - \! \alpha_{p_{q}})} \! + \! 
\widetilde{V}^{\prime}(\xi) \right) \right) \md \xi, \quad j \! = \! 1,2,
\dotsc,N.
\end{equation}
Via the Definitions~\eqref{eql3.7a17} and~\eqref{eql3.7a18}, it follows 
that, for $n \! \in \! \mathbb{N}$ and $k \! \in \! \lbrace 1,2,\dotsc,K 
\rbrace$ such that $\alpha_{p_{\mathfrak{s}}} \! := \! \alpha_{k} \! \neq 
\! \infty$, the associated system of $2(N \! + \! 1)$ real moment 
equations~\eqref{eql3.7g}--\eqref{eql3.7i} are equivalent to 
$\tilde{\mathcal{T}}_{j} \! = \! 0$, $j \! = \! 0,1,\dotsc,N$, $\tilde{
\mathcal{T}}_{N+1} \! = \! -2((n \! - \! 1)K \! + \! k)/n$, and $\tilde{
\mathcal{N}}_{j} \! = \! 0$, $j \! = \! 1,2,\dotsc,N$. It will first be shown 
that, for regular $\widetilde{V} \colon \overline{\mathbb{R}} \setminus 
\lbrace \alpha_{1},\alpha_{2},\dotsc,\alpha_{K} \rbrace \! \to \! \mathbb{R}$ 
satisfying conditions~\eqref{eq20}--\eqref{eq22}, the Jacobian of the 
transformation $\lbrace \tilde{b}_{0}(z_{o}),\tilde{b}_{1}(z_{o}),\dotsc,
\tilde{b}_{N}(z_{o}),\tilde{a}_{1}(z_{o}),\tilde{a}_{2}(z_{o}),\dotsc,
\tilde{a}_{N+1}(z_{o}) \rbrace \! \mapsto \! \lbrace \tilde{\mathcal{T}}_{1},
\tilde{\mathcal{T}}_{2},\dotsc,\tilde{\mathcal{T}}_{N+1},
\tilde{\mathcal{N}}_{1},\tilde{\mathcal{N}}_{2},\dotsc,\tilde{\mathcal{N}}_{N} 
\rbrace$ is non-zero whenever $\tilde{b}_{j-1} \! = \! \tilde{b}_{j-1}(z_{o})$ 
and $\tilde{a}_{j} \! = \! \tilde{a}_{j}(z_{o})$, $j \! = \! 1,2,\dotsc,
N \! + \! 1$, are chosen so that $J_{f} \! = \! \cup_{j=1}^{N+1}
[\tilde{b}_{j-1},\tilde{a}_{j}]$. Via Equation~\eqref{eql3.4a}, that is,
\begin{equation} \label{eql3.7a31} 
(\mathcal{H} \psi^{f}_{\widetilde{V}})(z) \! = \! \dfrac{1}{\pi} 
\left(\dfrac{(n \! - \! 1)K \! + \! k}{n} \right)^{-1} \left(\mi \pi 
\mathfrak{F}_{f}(z) \! + \! \dfrac{(\varkappa_{nk} \! - \! 1)}{n(z \! - 
\! \alpha_{k})} \! + \! \sum_{q=1}^{\mathfrak{s}-2} \dfrac{\varkappa_{nk 
\tilde{k}_{q}}}{n(z \! - \! \alpha_{p_{q}})} \right),
\end{equation}
and the integral representation~\eqref{eql3.7a12} for $\mathfrak{F}_{f}(z)$, 
one follows the analysis on pp.~778--779 of \cite{a57} (but adapted to 
the present situation; see, also, \cite{demcl}) to show that, for $n \! \in 
\! \mathbb{N}$ and $k \! \in \! \lbrace 1,2,\dotsc,K \rbrace$ such that 
$\alpha_{p_{\mathfrak{s}}} \! := \! \alpha_{k} \! \neq \! \infty$, and 
$i \! = \! 1,2,\dotsc,N \! + \! 1$:
\begin{gather}
\dfrac{\partial \tilde{\mathcal{T}}_{j}}{\partial \tilde{b}_{i-1}} \! = \! 
\tilde{b}_{i-1} \dfrac{\partial \tilde{\mathcal{T}}_{j-1}}{\partial 
\tilde{b}_{i-1}} \! + \! \dfrac{1}{2} \tilde{\mathcal{T}}_{j-1}, \quad j \! 
\in \! \mathbb{N}, \label{eql3.7a19} \\
\dfrac{\partial \tilde{\mathcal{T}}_{j}}{\partial \tilde{a}_{i}} \! = \! 
\tilde{a}_{i} \dfrac{\partial \tilde{\mathcal{T}}_{j-1}}{\partial 
\tilde{a}_{i}} \! + \! \dfrac{1}{2} \tilde{\mathcal{T}}_{j-1}, \quad j \! 
\in \! \mathbb{N}, \label{eql3.7a20} \\
\dfrac{\partial \mathfrak{F}_{f}(z)}{\partial \tilde{b}_{i-1}} \! = \! 
-\dfrac{1}{2 \pi \mi} \left(\dfrac{\partial \tilde{\mathcal{T}}_{0}}{\partial 
\tilde{b}_{i-1}} \right) \dfrac{(\tilde{R}(z))^{1/2}}{z \! - \! 
\tilde{b}_{i-1}}, \quad z \! \in \! \mathbb{C} \setminus (J_{f} \cup \lbrace 
\alpha_{p_{1}},\dotsc,\alpha_{p_{\mathfrak{s}-2}},\alpha_{k} \rbrace), 
\label{eql3.7a21} \\
\dfrac{\partial \mathfrak{F}_{f}(z)}{\partial \tilde{a}_{i}} \! = \! 
-\dfrac{1}{2 \pi \mi} \left(\dfrac{\partial \tilde{\mathcal{T}}_{0}}{\partial 
\tilde{a}_{i}} \right) \dfrac{(\tilde{R}(z))^{1/2}}{z \! - \! \tilde{a}_{i}}, 
\quad z \! \in \! \mathbb{C} \setminus (J_{f} \cup \lbrace \alpha_{p_{1}},
\dotsc,\alpha_{p_{\mathfrak{s}-2}},\alpha_{k} \rbrace), \label{eql3.7a22} \\
\dfrac{\partial \tilde{\mathcal{N}}_{j}}{\partial \tilde{b}_{i-1}} \! = \! 
-\dfrac{1}{2 \pi} \left(\dfrac{(n \! - \! 1)K \! + \! k}{n} \right)^{-1} 
\left(\dfrac{\partial \tilde{\mathcal{T}}_{0}}{\partial \tilde{b}_{i-1}} 
\right) \int_{\tilde{a}_{j}}^{\tilde{b}_{j}} 
\dfrac{(\tilde{R}(\xi))^{1/2}}{\xi \! - \! \tilde{b}_{i-1}} \, \md \xi, 
\quad j \! = \! 1,2,\dotsc,N, \label{eql3.7a23} \\
\dfrac{\partial \tilde{\mathcal{N}}_{j}}{\partial \tilde{a}_{i}} \! = \! 
-\dfrac{1}{2 \pi} \left(\dfrac{(n \! - \! 1)K \! + \! k}{n} \right)^{-1} 
\left(\dfrac{\partial \tilde{\mathcal{T}}_{0}}{\partial \tilde{a}_{i}} \right) 
\int_{\tilde{a}_{j}}^{\tilde{b}_{j}} \dfrac{(\tilde{R}(\xi))^{1/2}}{\xi \! 
- \! \tilde{a}_{i}} \, \md \xi, \quad j \! = \! 1,2,\dotsc,N. \label{eql3.7a24}
\end{gather}
Furthermore, evaluating Equations~\eqref{eql3.7a19} and~\eqref{eql3.7a20} 
on the solution of the associated system of $2(N \! + \! 1)$ real moment 
equations~\eqref{eql3.7g}--\eqref{eql3.7i}, that is, $\tilde{\mathcal{T}}_{j} 
\! = \! 0$, $j \! = \! 0,1,\dotsc,N$, $\tilde{\mathcal{T}}_{N+1} \! = \! 
-2((n \! - \! 1)K \! + \! k)/n$, and $\tilde{\mathcal{N}}_{j} \! = \! 0$, 
$j \! = \! 1,2,\dotsc,N$, one arrives at, for $i \! = \! 1,2,\dotsc,N \! + \! 1$,
\begin{equation} \label{eql3.7a25} 
\dfrac{\partial \tilde{\mathcal{T}}_{j}}{\partial \tilde{b}_{i-1}} \! = \! 
(\tilde{b}_{i-1})^{j} \dfrac{\partial \tilde{\mathcal{T}}_{0}}{\partial 
\tilde{b}_{i-1}} \quad \quad \text{and} \quad \quad \dfrac{\partial 
\tilde{\mathcal{T}}_{j}}{\partial \tilde{a}_{i}} \! = \! (\tilde{a}_{i})^{j} 
\dfrac{\partial \tilde{\mathcal{T}}_{0}}{\partial \tilde{a}_{i}}, \quad 
j \! = \! 0,1,\dotsc,N \! + \! 1.
\end{equation}
For $n \! \in \! \mathbb{N}$ and $k \! \in \! \lbrace 1,2,\dotsc,K 
\rbrace$ such that $\alpha_{p_{\mathfrak{s}}} \! := \! \alpha_{k} \! 
\neq \! \infty$, one evaluates, via Equations~\eqref{eql3.7a23}, 
\eqref{eql3.7a24}, and~\eqref{eql3.7a25}, and the multi-linearity 
property of the determinant, the Jacobian of the transformation 
$\lbrace \tilde{b}_{0}(z_{o}),\tilde{b}_{1}(z_{o}),\dotsc,\tilde{b}_{N}
(z_{o}),\tilde{a}_{1}(z_{o}),\linebreak[4]
\tilde{a}_{2}(z_{o}),\dotsc,\tilde{a}_{N+1}(z_{o}) \rbrace \! \mapsto 
\! \lbrace \tilde{\mathcal{T}}_{1},\tilde{\mathcal{T}}_{2},\dotsc,
\tilde{\mathcal{T}}_{N+1},\tilde{\mathcal{N}}_{1},\tilde{\mathcal{N}}_{2},
\dotsc,\tilde{\mathcal{N}}_{N} \rbrace$ on the solution of the associated 
system of $2(N \! + \! 1)$ real moment 
equations~\eqref{eql3.7g}--\eqref{eql3.7i}:
\begin{align}
&\operatorname{Jac}(\tilde{\mathcal{T}}_{0},\tilde{\mathcal{T}}_{1},\dotsc,
\tilde{\mathcal{T}}_{N+1},\tilde{\mathcal{N}}_{1},\tilde{\mathcal{N}}_{2},
\dotsc,\tilde{\mathcal{N}}_{N}) \! := \! \dfrac{\partial 
(\tilde{\mathcal{T}}_{0},\tilde{\mathcal{T}}_{1},\dotsc,
\tilde{\mathcal{T}}_{N+1},\tilde{\mathcal{N}}_{1},\tilde{\mathcal{N}}_{2},
\dotsc,\tilde{\mathcal{N}}_{N})}{\partial (\tilde{b}_{0},\tilde{b}_{1},\dotsc,
\tilde{b}_{N},\tilde{a}_{1},\tilde{a}_{2},\dotsc,\tilde{a}_{N+1})} \nonumber \\
&= 
.
\end{equation*}
The above determinant $\tilde{\Delta}_{d}(\vec{\xi})$ has, modulo notation, 
been evaluated on p.~780 of \cite{a57}:
\begin{equation*}
\tilde{\Delta}_{d}(\vec{\xi}) \! = \! \dfrac{\left(\prod_{j=1}^{N+1} 
\prod_{i=1}^{N+1}(\tilde{b}_{i-1} \! - \! \tilde{a}_{j}) \right) \left(
\prod_{\underset{j<i}{i,j=1}}^{N+1}(\tilde{b}_{i-1} \! - \! 
\tilde{b}_{j-1})(\tilde{a}_{i} \! - \! \tilde{a}_{j}) \right) \left(
\prod_{\underset{j<i}{i,j=1}}^{N}(\xi_{i} \! - \! \xi_{j}) \right)}{(-1)^{N} 
\prod_{j=1}^{N} \prod_{i=1}^{N+1}(\xi_{j} \! - \! \tilde{b}_{i-1})(\xi_{j} 
\! - \! \tilde{a}_{i})};
\end{equation*}
however, for $(-\infty \! < \! \tilde{b}_{0} \! <)$ $\tilde{a}_{1} \! < \! 
\xi_{1} \! < \! \tilde{b}_{1} \! < \! \tilde{a}_{2} \! < \! \xi_{2} \! < \! 
\tilde{b}_{2} \! < \! \dotsb \! < \! \tilde{b}_{N-1} \! < \! \tilde{a}_{N} \! 
< \! \xi_{N} \! < \! \tilde{b}_{N}$ $(< \! \tilde{a}_{N+1} \! < \! +\infty)$, 
$\tilde{\Delta}_{d}(\vec{\xi}) \! \not= \! 0$ (which means that it is of 
fixed sign), and $\int_{\tilde{a}_{j}}^{\tilde{b}_{j}}(\tilde{R}(\xi_{j}))^{1/2} 
\, \md \xi_{j} \! \neq \! 0$, $j \! = \! 1,2,\dotsc,N$, whence
\begin{equation} \label{eql3.7a28} 
\left(\prod_{j=1}^{N} \int_{\tilde{a}_{j}}^{\tilde{b}_{j}}(\tilde{R}
(\xi_{j}))^{1/2} \, \md \xi_{j} \right) \tilde{\Delta}_{d}(\vec{\xi}) 
\! \neq \! 0.
\end{equation}
Exploiting the fact that (cf. definition~\eqref{eql3.7a17}) 
$\tilde{\mathcal{T}}_{0}$ is independent of $z$, it follows {}from the 
integral representations~\eqref{eql3.7l} and~\eqref{eql3.7a12}, and 
Equations~\eqref{eql3.7a21} and~\eqref{eql3.7a22}, that, for 
$j \! = \! 1,2,\dotsc,N \! + \! 1$,
\begin{gather*}
\dfrac{(z \! - \! \tilde{b}_{j-1})}{(\tilde{R}(z))^{1/2}} \dfrac{\partial 
\mathfrak{F}_{f}(z)}{\partial \tilde{b}_{j-1}} \! = \! -\dfrac{1}{2 \pi \mi} 
\left(\dfrac{(n \! - \! 1)K \! + \! k}{n} \right) \left((z \! - \! 
\tilde{b}_{j-1}) \dfrac{\partial \tilde{h}_{\widetilde{V}}(z)}{\partial 
\tilde{b}_{j-1}} \! - \! \dfrac{1}{2} \tilde{h}_{\widetilde{V}}(z) \right), \\
\dfrac{(z \! - \! \tilde{a}_{j})}{(\tilde{R}(z))^{1/2}} \dfrac{\partial 
\mathfrak{F}_{f}(z)}{\partial \tilde{a}_{j}} \! = \! -\dfrac{1}{2 \pi \mi} 
\left(\dfrac{(n \! - \! 1)K \! + \! k}{n} \right) \left((z \! - \! 
\tilde{a}_{j}) \dfrac{\partial \tilde{h}_{\widetilde{V}}(z)}{\partial 
\tilde{a}_{j}} \! - \! \dfrac{1}{2} \tilde{h}_{\widetilde{V}}(z) \right):
\end{gather*}
via the latter two formulae, the $z$-independence of 
$\tilde{\mathcal{T}}_{0}$, and the fact that, for regular $\widetilde{V} 
\colon \overline{\mathbb{R}} \setminus \lbrace \alpha_{1},\alpha_{2},
\dotsc,\alpha_{K} \rbrace \! \to \! \mathbb{R}$ satisfying 
conditions~\eqref{eq20}--\eqref{eq22}, $\tilde{h}_{\widetilde{V}}(z) \! 
\neq \! 0$ for $z \! \in \! J_{f}$ (in particular, $\tilde{h}_{\widetilde{V}}
(\tilde{b}_{j-1}),\tilde{h}_{\widetilde{V}}(\tilde{a}_{j}) \! \neq \! 0$, $j \! 
= \! 1,2,\dotsc,N \! + \! 1)$, one shows that, for $j \! = \! 1,2,\dotsc,
N \! + \! 1$,
\begin{gather*}
\left. \dfrac{(z \! - \! \tilde{b}_{j-1})}{(\tilde{R}(z))^{1/2}} 
\dfrac{\partial \mathfrak{F}_{f}(z)}{\partial \tilde{b}_{j-1}} 
\right\vert_{z=\tilde{b}_{j-1}}=\dfrac{1}{4 \pi \mi} \left(\dfrac{(n \! - \! 
1)K \! + \! k}{n} \right) \tilde{h}_{\widetilde{V}}(\tilde{b}_{j-1}) \neq 0, \\
\left. \dfrac{(z \! - \! \tilde{a}_{j})}{(\tilde{R}(z))^{1/2}} \dfrac{\partial 
\mathfrak{F}_{f}(z)}{\partial \tilde{a}_{j}} \right\vert_{z=\tilde{a}_{j}}=
\dfrac{1}{4 \pi \mi} \left(\dfrac{(n \! - \! 1)K \! + \! k}{n} \right) 
\tilde{h}_{\widetilde{V}}(\tilde{a}_{j}) \neq 0,
\end{gather*}
which implies, via Equations~\eqref{eql3.7a21} and~\eqref{eql3.7a22}, 
that, for $j \! = \! 1,2,\dotsc,N \! + \! 1$,
\begin{equation*}
\dfrac{\partial \tilde{\mathcal{T}}_{0}}{\partial \tilde{b}_{j-1}} \! 
= \! -\dfrac{1}{2} \left(\dfrac{(n \! - \! 1)K \! + \! k}{n} \right) 
\tilde{h}_{\widetilde{V}}(\tilde{b}_{j-1}) \! \neq \! 0 \quad \quad 
\text{and} \quad \quad \dfrac{\partial \tilde{\mathcal{T}}_{0}}{\partial 
\tilde{a}_{j}} \! = \! -\dfrac{1}{2} \left(\dfrac{(n \! - \! 1)K \! + \! k}{n} 
\right) \tilde{h}_{\widetilde{V}}(\tilde{a}_{j}) \! \neq \! 0,
\end{equation*}
whence
\begin{equation} \label{eql3.7a29} 
\prod_{m=1}^{N+1} \dfrac{\partial \tilde{\mathcal{T}}_{0}}{\partial 
\tilde{b}_{m-1}} \dfrac{\partial \tilde{\mathcal{T}}_{0}}{\partial 
\tilde{a}_{m}} \! = \! \left(\dfrac{1}{2} \left(\dfrac{(n \! - \! 1)K \! + 
\! k}{n} \right) \right)^{2(N+1)} \prod_{j=1}^{N+1} \tilde{h}_{\widetilde{V}}
(\tilde{b}_{j-1}) \tilde{h}_{\widetilde{V}}(\tilde{a}_{j}) \! \neq \! 0;
\end{equation}
hence, via Equations~\eqref{eql3.7a26}, \eqref{eql3.7a28}, 
and~\eqref{eql3.7a29}, one concludes that, for $n \! \in \! 
\mathbb{N}$ and $k \! \in \! \lbrace 1,2,\dotsc,K \rbrace$ such 
that $\alpha_{p_{\mathfrak{s}}} \! := \! \alpha_{k} \! \neq \! \infty$, 
$\operatorname{Jac}(\tilde{\mathcal{T}}_{0},\tilde{\mathcal{T}}_{1},
\dotsc,\tilde{\mathcal{T}}_{N+1},\tilde{\mathcal{N}}_{1},
\tilde{\mathcal{N}}_{2},\dotsc,\tilde{\mathcal{N}}_{N}) \! \neq \!0$. In 
order to apply the Implicit Function Theorem, it remains to show that, 
for $n \! \in \! \mathbb{N}$ and $k \! \in \! \lbrace 1,2,\dotsc,K 
\rbrace$ such that $\alpha_{p_{\mathfrak{s}}} \! := \! \alpha_{k} \! 
\neq \! \infty$, $\tilde{\mathcal{T}}_{j}$, $j \! = \! 0,1,\dotsc,N \! 
+ \! 1$, and $\tilde{\mathcal{N}}_{i}$, $i \! = \! 1,2,\dotsc,N$, are 
real analytic functions of the end-points of the intervals, $\lbrace 
\tilde{b}_{j-1}(z_{o}),\tilde{a}_{j}(z_{o}) \rbrace_{j=1}^{N+1}$, of 
the support, $J_{f}$, of the associated equilibrium measure, 
$\mu^{f}_{\widetilde{V}}$. Noting the $z$-independence of 
$\tilde{\mathcal{T}}_{j}$, $j \! \in \! \mathbb{N}_{0}$, it follows 
via a residue calculation that, equivalently,
\begin{equation*}
\tilde{\mathcal{T}}_{j} \! = \! \dfrac{1}{2} \oint_{\tilde{C}_{\widetilde{V}}} 
\dfrac{\xi^{j}}{(\tilde{R}(\xi))^{1/2}} \left(\dfrac{2}{\mi \pi} \left(
\dfrac{(\varkappa_{nk} \! - \! 1)}{n(\xi \! - \! \alpha_{k})} \! + \! 
\sum_{q=1}^{\mathfrak{s}-2} \dfrac{\varkappa_{nk \tilde{k}_{q}}}{n(\xi \! - \! 
\alpha_{p_{q}})} \right) \! + \! \dfrac{\widetilde{V}^{\prime}(\xi)}{\mi \pi} 
\right) \md \xi, \quad j \! \in \! \mathbb{N}_{0},
\end{equation*}
where the contour $\tilde{C}_{\widetilde{V}}$ is described in 
item~$\pmb{(2)}$ of the lemma: the only factor depending on 
$\tilde{b}_{j-1},\tilde{a}_{j}$, $j \! = \! 1,2,\dotsc,N \! + \! 1$, is 
$(\tilde{R}(z))^{1/2}$. As $(\tilde{R}(z))^{1/2}$ is analytic for $z \! 
\in \! \mathbb{C} \setminus \cup_{j=1}^{N+1}[\tilde{b}_{j-1},
\tilde{a}_{j}]$, and since $\mathbb{C} \setminus \cup_{j=1}^{N+1}
[\tilde{b}_{j-1},\tilde{a}_{j}] \supset \tilde{C}_{\widetilde{V}}$ (with 
$\operatorname{int}(\tilde{C}_{\widetilde{V}}) \supset \lbrace z \rbrace 
\cup \cup_{j=1}^{N+1}[\tilde{b}_{j-1},\tilde{a}_{j}])$, it follows, in 
particular, that $(\tilde{R}(z))^{1/2} \! \! \upharpoonright_{\tilde{C}_{
\widetilde{V}}}$ is an analytic function of $\tilde{b}_{j-1},\tilde{a}_{j}$, 
$j \! = \! 1,2,\dotsc,N \! + \! 1$, which implies, via the above contour 
integral representation for $\tilde{\mathcal{T}}_{j}$, $j \! \in \! \mathbb{N}_{0}$, 
that, for $n \! \in \! \mathbb{N}$ and $k \! \in \! \lbrace 1,2,\dotsc,K \rbrace$ 
such that $\alpha_{p_{\mathfrak{s}}} \! := \! \alpha_{k} \! \neq \! \infty$, 
$\tilde{\mathcal{T}}_{i}$, $i \! = \! 0,1,\dotsc,N \! + \! 1$, are real analytic 
functions of $\tilde{b}_{j-1},\tilde{a}_{j}$, $j \! = \! 1,2,\dotsc,N \! + \! 1$. 
One shows {}from the integral representation~\eqref{eql3.7l} for 
$\tilde{h}_{\widetilde{V}}(z)$, Equations~\eqref{eql3.7a30} 
and~\eqref{eql3.7a31}, and the Definition~\eqref{eql3.7a18}, that
\begin{equation*}
\tilde{\mathcal{N}}_{j} \! = \! -\dfrac{1}{2 \pi} \int_{\tilde{a}_{j}}^{
\tilde{b}_{j}}(\tilde{R}(\xi))^{1/2} \tilde{h}_{\widetilde{V}}(\xi) \, 
\md \xi, \quad j \! = \! 1,2,\dotsc,N:
\end{equation*}
making the linear change of variable $\tilde{u}_{j} \colon \mathbb{C} 
\! \to \! \mathbb{C}$, $\xi \! \mapsto \! \tilde{u}_{j}(\xi) \! := \! 
(\tilde{b}_{j} \! - \! \tilde{a}_{j})^{-1}(\xi \! - \! \tilde{a}_{j})$, 
$j \! = \! 1,2,\dotsc,N$, one shows {}from the above expression for 
$\tilde{\mathcal{N}}_{j}$, $j \! = \! 1,2,\dotsc,N$, that, for $n \! \in 
\! \mathbb{N}$ and $k \! \in \! \lbrace 1,2,\dotsc,K \rbrace$ such 
that $\alpha_{p_{\mathfrak{s}}} \! := \! \alpha_{k} \! \neq \! \infty$,
\begin{equation*}
\tilde{\mathcal{N}}_{j} \! = \! -\dfrac{1}{2 \pi}(\tilde{b}_{j} \! - \! 
\tilde{a}_{j})^{2} \int_{0}^{1}(\tilde{\mathcal{R}}^{\ast}(\tilde{u}_{j}
(\tilde{b}_{j} \! - \! \tilde{a}_{j}) \! + \! \tilde{a}_{j}))^{1/2} 
\tilde{h}_{\widetilde{V}}(\tilde{u}_{j}(\tilde{b}_{j} \! - \! \tilde{a}_{j}) 
\! + \! \tilde{a}_{j})(\tilde{u}_{j} \lvert \tilde{u}_{j} \! - \! 1 
\rvert)^{1/2} \, \md \tilde{u}_{j}, \quad j \! = \! 1,2,\dotsc,N,
\end{equation*}
where $(\tilde{\mathcal{R}}^{\ast}(\xi))^{1/2} \! := \! (-1)^{N-j+1} 
\prod_{i_{1}=1}^{j} \lvert \xi \! - \! \tilde{b}_{i_{1}-1} \rvert^{1/2} 
\prod_{i_{2}=1}^{j-1} \lvert \xi \! - \! \tilde{a}_{i_{2}} \rvert^{1/2} 
\prod_{i_{3}=j+2}^{N+1} \lvert \xi \! - \! \tilde{b}_{i_{3}-1} \rvert^{1/2} 
\prod_{i_{4}=j+1}^{N+1} \lvert \xi \! - \! \tilde{a}_{i_{4}} \rvert^{1/2}$. 
Recalling that $\tilde{h}_{\widetilde{V}}(z)$ is real analytic on $\mathbb{R} 
\setminus \lbrace \alpha_{p_{1}},\dotsc,\alpha_{p_{\mathfrak{s}-2}},
\alpha_{k} \rbrace$, $\tilde{h}_{\widetilde{V}}(\tilde{b}_{j-1}),\tilde{h}_{
\widetilde{V}}(\tilde{a}_{j}) \! \neq \! 0$, $j \! = \! 1,2,\dotsc,N \! + \! 1$, 
and that it is an analytic function of $\tilde{b}_{j-1}(z_{o}),\tilde{a}_{j}
(z_{o})$, $j \! = \! 1,2,\dotsc,N \! + \! 1$, it follows {}from the above 
formula for $\tilde{\mathcal{N}}_{i}$, $i \! = \! 1,2,\dotsc,N$, that, for 
$n \! \in \! \mathbb{N}$ and $k \! \in \! \lbrace 1,2,\dotsc,K \rbrace$ 
such that $\alpha_{p_{\mathfrak{s}}} \! := \! \alpha_{k} \! \neq \! \infty$, 
it, too, is an analytic function of $\tilde{b}_{j-1}(z_{o}),\tilde{a}_{j}(z_{o})$, 
$j \! = \! 1,2,\dotsc,N \! + \! 1$. Thus, for $n \! \in \! \mathbb{N}$ 
and $k \! \in \! \lbrace 1,2,\dotsc,K \rbrace$ such that 
$\alpha_{p_{\mathfrak{s}}} \! := \! \alpha_{k} \! \neq \! \infty$, since 
the Jacobian of the transformation $\lbrace \tilde{b}_{0}(z_{o}),
\tilde{b}_{1}(z_{o}),\dotsc,\tilde{b}_{N}(z_{o}),\tilde{a}_{1}(z_{o}),\tilde{a}_{2}
(z_{o}),\dotsc,\tilde{a}_{N+1}(z_{o}) \rbrace \! \mapsto \! \lbrace 
\tilde{\mathcal{T}}_{0},\tilde{\mathcal{T}}_{1},\dotsc,\tilde{\mathcal{T}}_{N+1},
\tilde{\mathcal{N}}_{1},\tilde{\mathcal{N}}_{2},\dotsc,\tilde{\mathcal{N}}_{N} 
\rbrace$ is non-zero whenever $\lbrace \tilde{b}_{j-1}(z_{o}),\tilde{a}_{j}
(z_{o}) \rbrace_{j=1}^{N+1}$ is chosen so that, for regular $\widetilde{V} 
\colon \overline{\mathbb{R}} \setminus \lbrace \alpha_{1},\alpha_{2},
\dotsc,\alpha_{K} \rbrace \! \to \! \mathbb{R}$ satisfying 
conditions~\eqref{eq20}--\eqref{eq22}, $J_{f} \! := \! \cup_{j=1}^{N+1}
[\tilde{b}_{j-1}(z_{o}),\tilde{a}_{j}(z_{o})]$, and $\tilde{\mathcal{T}}_{j}$, 
$j \! = \! 0,1,\dotsc,N \! + \! 1$, and $\tilde{\mathcal{N}}_{i}$, $i \! = 
\! 1,2,\dotsc,N$, are real analytic functions of $\tilde{b}_{m-1}(z_{o}),
\tilde{a}_{m}(z_{o})$, $m \! = \! 1,2,\dotsc,N \! + \! 1$, it follows 
{}from the Implicit Function Theorem that, in the double-scaling limit 
$\mathscr{N},n\! \to \! \infty$ such that $z_{o} \! = \! 1 \! + \! o(1)$, 
$\tilde{b}_{j-1}(z_{o}),\tilde{a}_{j}(z_{o})$, $j \! = \! 1,2,\dotsc,N \! + \! 
1$, are, with the exception of a denumerable set of singularities (the 
`critical values' of $z_{o})$, real analytic functions of $z_{o}$.

$\pmb{(\mathrm{B})}$ The proof of this case, that is, $n \! \in \! \mathbb{N}$ 
and $k \! \in \! \lbrace 1,2,\dotsc,K \rbrace$ such that $\alpha_{p_{\mathfrak{s}}} 
\! := \! \alpha_{k} \! = \! \infty$, is virtually identical to the proof presented in 
$\pmb{(\mathrm{A})}$; one mimics, \emph{verbatim}, the scheme of the calculations 
presented in case $\pmb{(\mathrm{A})}$ in order to arrive at the corresponding 
claims stated in item~$\pmb{(1)}$ of the lemma; in order to do so, 
however, the analogues of Equations~\eqref{eql3.7m}, \eqref{eql3.7n}, 
\eqref{eql3.7r}--\eqref{eql3.7t}, \eqref{eql3.7a6}, 
\eqref{eql3.7a9}--\eqref{eql3.7a12}, \eqref{eql3.7a30}, \eqref{eql3.7a17}, 
\eqref{eql3.7a18}, and~\eqref{eql3.7a19}--\eqref{eql3.7a29}, respectively, 
are necessary, which, in the present case, read:
\begin{align} \label{eql3.7b1} 
\hat{\mathcal{S}} &\colon \mathbb{N} \times \lbrace 1,2,\dotsc,K \rbrace 
\times \mathbb{C} \setminus (J_{\infty} \cup \lbrace \alpha_{p_{1}},
\alpha_{p_{2}},\dotsc,\alpha_{p_{\mathfrak{s}-1}} \rbrace) \! \ni \! (n,k,z) 
\! \mapsto \! 4 \mi \left(\dfrac{(n \! - \! 1)K \! + \! k}{n} \right)^{2}
(\mathcal{H} \psi_{\widetilde{V}}^{\infty})(z) \psi_{\widetilde{V}}^{\infty}
(z) \nonumber \\
&-\dfrac{4 \mi}{\pi} \left(\dfrac{(n \! - \! 1)K \! + \! k}{n} \right) 
\left(\sum_{q=1}^{\mathfrak{s}-1} \dfrac{\varkappa_{nk \tilde{k}_{q}}}{n
(z \! - \! \alpha_{p_{q}})} \right) \psi_{\widetilde{V}}^{\infty}(z) \! = \! 
\hat{\mathcal{S}}(n,k,z) \! =: \! \hat{\mathcal{S}}(z),
\end{align}
where $J_{\infty} \! := \! \supp (\mu_{\widetilde{V}}^{\infty})$ and 
$\psi_{\widetilde{V}}^{\infty}(z)$ are described in item~$\pmb{(1)}$ 
of the lemma,
\begin{align} \label{eql3.7b2} 
\hat{\mathfrak{H}} &\colon \mathbb{N} \times \lbrace 1,2,\dotsc,K 
\rbrace \times \mathbb{C} \setminus (J_{\infty} \cup \lbrace \alpha_{p_{1}},
\alpha_{p_{2}},\dotsc,\alpha_{p_{\mathfrak{s}-1}} \rbrace) \! \ni \! 
(n,k,z) \! \mapsto \! (\mathfrak{F}_{\infty}(z))^{2} \! - \! \int_{J_{\infty}} 
\dfrac{\hat{\mathcal{S}}(\xi)}{\xi \! - \! z} \, \dfrac{\md \xi}{2 \pi \mi} 
\! = \! \hat{\mathfrak{H}}(n,k,z) \! =: \! \hat{\mathfrak{H}}(z),
\end{align}
where $\mathfrak{F}_{\infty}(z)$ is defined by Equation~\eqref{eql3.4i},
\begin{equation} \label{eql3.7b3} 
(\mathcal{H} \psi_{\widetilde{V}}^{\infty})(z) \! = \! \dfrac{1}{2 \pi} 
\left(\dfrac{(n \! - \! 1)K \! + \! k}{n} \right)^{-1} \left(2 
\sum_{q=1}^{\mathfrak{s}-1} \dfrac{\varkappa_{nk \tilde{k}_{q}}}{n
(z \! - \! \alpha_{p_{q}})} \! + \! \widetilde{V}^{\prime}(z) \right), 
\quad z \! \in \! J_{\infty},
\end{equation}
\begin{equation} \label{eql3.7b4} 
\hat{\mathcal{S}}(z) \! = \! \dfrac{2 \mi}{\pi} \left(\dfrac{(n \! - \! 1)K 
\! + \! k}{n} \right) \widetilde{V}^{\prime}(z) \psi^{\infty}_{\widetilde{V}}
(z), \quad z \! \in \! J_{\infty},
\end{equation}
\begin{equation} \label{eql3.7b5} 
\hat{\mathfrak{H}}(z) \! = \! \dfrac{\sum_{m=0}^{2 \mathfrak{s}-3} 
\hat{\rho}_{m}(n,k)z^{m}}{\prod_{q=1}^{\mathfrak{s}-1}(z \! - \! 
\alpha_{p_{q}})^{2}}, \qquad z \! \in \! \mathbb{C} \setminus (J_{\infty} 
\cup \lbrace \alpha_{p_{1}},\alpha_{p_{2}},\dotsc,\alpha_{p_{\mathfrak{s}
-1}} \rbrace),
\end{equation}
where
\begin{equation} \label{eql3.7b6} 
\hat{\rho}_{2 \mathfrak{s}-3} (n,k) \! = \! \dfrac{1}{\pi^{2}} 
\left(\dfrac{(n \! - \! 1)K \! + \! k}{n} \right) \int_{J_{\infty}} 
\widetilde{V}^{\prime}(\xi) \psi_{\widetilde{V}}^{\infty}(\xi) \, \md \xi,
\end{equation}
\begin{align} \label{eql3.7b7} 
\hat{\rho}_{m}(n,k) =& \, \hat{w}_{2 \mathfrak{s}-m-3}(n,k) \! + \! 
\hat{c}_{m}^{\natural}(n,k) \! + \! \sum_{j=m}^{2(\mathfrak{s}-2)} 
\left(\hat{w}_{j-m}(n,k) \hat{c}_{j+1}(n,k) \! + \! \hat{\nu}_{j-m}
(n,k) \hat{c}_{j+1}^{\flat}(n,k) \right. \nonumber \\
+&\left. \, \hat{u}_{j-m}(n,k) \hat{c}_{j+2}^{\sharp}(n,k) \right), 
\quad m \! = \! 0,1,\dotsc,2(\mathfrak{s} \! - \! 2),
\end{align}
with
\begin{equation} \label{eql3.7b8} 
\hat{w}_{r}(n,k) \! = \! \dfrac{1}{\pi^{2}} \left(\dfrac{(n \! - \! 1)K 
\! + \! k}{n} \right) \int_{J_{\infty}} \xi^{r} \widetilde{V}^{\prime}(\xi) 
\psi_{\widetilde{V}}^{\infty}(\xi) \, \md \xi, \quad r \! \in \! 
\mathbb{N}_{0},
\end{equation}
\begin{equation} \label{eql3.7b9} 
\hat{c}_{l}(n,k) \! = \! \mathlarger{\sum_{\underset{\underset{
\sum_{r=1}^{\mathfrak{s}-1}i_{r}=2(\mathfrak{s}-1)-l}{r \in \lbrace 1,2,
\dotsc,\mathfrak{s}-1 \rbrace}}{i_{r}=0,1,2}}}(-1)^{2(\mathfrak{s}-1)-l}
(2!)^{\mathfrak{s}-1} \prod_{m=1}^{\mathfrak{s}-1} \dfrac{1}{i_{m}!(2 
\! - \! i_{m})!} \prod_{j=1}^{\mathfrak{s}-1}(\alpha_{p_{j}})^{i_{j}}, \quad 
l \! = \! 0,1,\dotsc,2(\mathfrak{s} \! - \! 1),
\end{equation}
\begin{gather}
\hat{u}_{m}(n,k) \! = \! \sum_{j=0}^{m} \hat{\nu}_{j}(n,k) \hat{\nu}_{m-j}
(n,k), \quad m \! = \! 0,1,\dotsc,2(\mathfrak{s} \! - \! 2), 
\label{eql3.7b10} \\
\hat{\nu}_{r}(n,k) \! = \! \int_{J_{\infty}} \xi^{r} \psi_{\widetilde{V}}^{
\infty}(\xi) \, \md \xi, \quad r \! \in \! \mathbb{N}_{0}, \label{eql3.7b11}
\end{gather}
\begin{equation} \label{eql3.7b12} 
\hat{c}_{l}^{\sharp}(n,k) \! = \! -\dfrac{1}{\pi^{2}} \left(\dfrac{(n \! - \! 
1)K \! + \! k}{n} \right)^{2} \hat{c}_{l}(n,k), \quad l \! = \! 0,1,\dotsc,
2(\mathfrak{s} \! - \! 1),
\end{equation}
\begin{align}
\hat{c}_{l}^{\flat}(n,k) =& \, \dfrac{2}{\pi^{2}} \left(\dfrac{(n \! - \! 1)K 
\! + \! k}{n} \right) \mathlarger{\sum_{q=1}^{\mathfrak{s}-1}} \dfrac{
\varkappa_{nk \tilde{k}_{q}}}{n} \mathlarger{\sum_{\underset{\underset{
\underset{\sum_{m^{\prime}=1}^{\mathfrak{s}-1}i_{m^{\prime}}=2(\mathfrak{s}-
1)-l-1}{r \in \lbrace 1,2,\dotsc,\mathfrak{s}-1 \rbrace \setminus \lbrace 
q \rbrace}}{i_{r}=0,1,2}}{i_{q}=0,1}}} \dfrac{(-1)^{2(\mathfrak{s}-1)-l-1}
(2!)^{\mathfrak{s}-2}}{i_{q}!(1 \! - \! i_{q})!} \prod_{\substack{m=1\\
m \neq q}}^{\mathfrak{s}-1} \dfrac{1}{i_{m}!(2 \! - \! i_{m})!} \nonumber \\
\times& \, \prod_{j=1}^{\mathfrak{s}-1}(\alpha_{p_{j}})^{i_{j}}, \quad 
l \! = \! 0,1,\dotsc,2 \mathfrak{s} \! - \! 3, \label{eql3.7b13} 
\end{align}
\begin{align}
\hat{c}_{l}^{\natural}(n,k) =& \, -\dfrac{1}{\pi^{2}} \mathlarger{\sum_{q=
1}^{\mathfrak{s}-1}} \left(\dfrac{\varkappa_{nk \tilde{k}_{q}}}{n} \right)^{2} 
\mathlarger{\sum_{\underset{\underset{\sum_{\underset{m^{\prime} 
\neq q}{m^{\prime}=1}}^{\mathfrak{s}-1}i_{m^{\prime}}=2(\mathfrak{s}-
2)-l}{r \in \lbrace 1,2,\dotsc,\mathfrak{s}-1 \rbrace \setminus \lbrace q 
\rbrace}}{i_{r}=0,1,2}}}(-1)^{2(\mathfrak{s}-2)-l}(2!)^{\mathfrak{s}-2} 
\prod_{\substack{m=1\\m \neq q}}^{\mathfrak{s}-1} \dfrac{1}{i_{m}!
(2 \! - \! i_{m})!} \prod_{\substack{j=1\\j \neq q}}^{\mathfrak{s}-1}
(\alpha_{p_{j}})^{i_{j}} \nonumber \\
-& \, \dfrac{2}{\pi^{2}} \mathlarger{\sum_{p^{\prime}=1}^{\mathfrak{s}-2}} 
\mathlarger{\sum_{q=p^{\prime}+1}^{\mathfrak{s}-1}} \dfrac{\varkappa_{nk 
\tilde{k}_{p^{\prime}}}}{n} \dfrac{\varkappa_{nk \tilde{k}_{q}}}{n} 
\mathlarger{\sum_{\underset{\underset{\underset{i_{p^{\prime}}+i_{q}+
\sum_{\underset{m^{\prime} \neq p^{\prime},q}{m^{\prime}=1}}^{\mathfrak{s}
-1}i_{m^{\prime}}=2(\mathfrak{s}-2)-l}{r \in \lbrace 1,2,\dotsc,\mathfrak{s}-1 
\rbrace \setminus \lbrace p^{\prime},q \rbrace}}{i_{r}=0,1,2}}{i_{p^{\prime}},
i_{q}=0,1}}} \dfrac{(-1)^{2(\mathfrak{s}-2)-l}(2!)^{\mathfrak{s}-3}}{
i_{p^{\prime}}!(1 \! - \! i_{p^{\prime}})!i_{q}!(1 \! - \! i_{q})!} \nonumber \\
\times& \, \prod_{\substack{m=1\\m \neq p^{\prime},q} }^{
\mathfrak{s}-1} \dfrac{1}{i_{m}!(2 \! - \! i_{m})!}(\alpha_{p_{p^{\prime}}})^{
i_{p^{\prime}}}(\alpha_{p_{q}})^{i_{q}} \prod_{\substack{j=1\\j \neq p^{\prime},
q}}^{\mathfrak{s}-1}(\alpha_{p_{j}})^{i_{j}}, \quad l \! = \! 0,1,\dotsc,
2(\mathfrak{s} \! - \! 2), 
\label{eql3.7b14}
\end{align}
\begin{equation} \label{eql3.7b15} 
(\mathfrak{F}_{\infty}(z))^{2} \! = \! \dfrac{1}{\pi^{2}} \left(\dfrac{(n \! - 
\! 1)K \! + \! k}{n} \right) \int_{J_{\infty}} \dfrac{\widetilde{V}^{\prime}
(\xi) \psi_{\widetilde{V}}^{\infty}(\xi)}{\xi \! - \! z} \, \md \xi \! + \! 
\dfrac{\sum_{m=0}^{2 \mathfrak{s}-3} \hat{\rho}_{m}(n,k)z^{m}}{\prod_{q=1}^{
\mathfrak{s}-1}(z \! - \! \alpha_{p_{q}})^{2}}, \quad z \! \in \! \mathbb{C} 
\setminus (J_{\infty} \cup \lbrace \alpha_{p_{1}},\alpha_{p_{2}},\dotsc,
\alpha_{p_{\mathfrak{s}-1}} \rbrace),
\end{equation}
\begin{equation} \label{eql3.7b16} 
\left(\mathfrak{F}_{\infty}(z) \! + \! \dfrac{\mi \widetilde{V}^{\prime}(z)}{2 
\pi} \right)^{2} \! + \! \dfrac{\hat{q}_{\widetilde{V}}(z)}{\pi^{2}} \! = 
\! 0, \quad z \! \in \! \mathbb{C} \setminus (J_{\infty} \cup \lbrace 
\alpha_{p_{1}},\alpha_{p_{2}},\dotsc,\alpha_{p_{\mathfrak{s}-1}} \rbrace),
\end{equation}
where
\begin{equation} \label{eqveeinf} 
\begin{aligned}
\hat{q}_{\widetilde{V}}(z) :=& \, \left(\dfrac{\widetilde{V}^{\prime}(z)}{2} 
\right)^{2} \! + \! \widetilde{V}^{\prime}(z) \sum_{q=1}^{\mathfrak{s}-1} 
\dfrac{\varkappa_{nk \tilde{k}_{q}}}{n(z \! - \! \alpha_{p_{q}})} \! - \! 
\left(\dfrac{(n \! - \! 1)K \! + \! k}{n} \right) \int_{J_{\infty}} 
\dfrac{(\widetilde{V}^{\prime}(\xi) \! - \! \widetilde{V}^{\prime}(z)) 
\psi_{\widetilde{V}}^{\infty}(\xi)}{\xi \! - \! z} \, \md \xi \\
-& \, \dfrac{\pi^{2} \sum_{m=0}^{2 \mathfrak{s}-3} \hat{\rho}_{m}(n,k)
z^{m}}{\prod_{q=1}^{\mathfrak{s}-1}(z \! - \! \alpha_{p_{q}})^{2}} \\
=& \, \left(\dfrac{\widetilde{V}^{\prime}(z)}{2} \right)^{2} \! + \! 
\widetilde{V}^{\prime}(z) \sum_{q=1}^{\mathfrak{s}-1} \dfrac{\varkappa_{nk 
\tilde{k}_{q}}}{n(z \! - \! \alpha_{p_{q}})} \! - \! \left(\dfrac{(n \! - \! 
1)K \! + \! k}{n} \right) \int_{J_{\infty}} \int_{0}^{1} \widetilde{V}^{\prime 
\prime}(t \xi \! + \! (1 \! - \! t)z) \, \md t \, \md \mu_{\widetilde{V}}^{
\infty}(\xi) \\
-& \, \dfrac{\pi^{2} \sum_{m=0}^{2 \mathfrak{s}-3} \hat{\rho}_{m}(n,k)
z^{m}}{\prod_{q=1}^{\mathfrak{s}-1}(z \! - \! \alpha_{p_{q}})^{2}},
\end{aligned}
\end{equation}
with $\md \mu_{\widetilde{V}}^{\infty}$ described in item~$\pmb{(1)}$ 
of the lemma,
\begin{equation} \label{eql3.7b17} 
\mathfrak{F}_{\infty \, \pm}(z) \! = \! \dfrac{\widetilde{V}^{\prime}(z)}{2 
\pi \mi} \! \mp \! \left(\dfrac{(n \! - \! 1)K \! + \! k}{n} \right) 
\psi_{\widetilde{V}}^{\infty}(z), \quad z \! \in \! J_{\infty},
\end{equation}
where $\mathfrak{F}_{\infty \, \pm}(z) \! := \! 
\lim_{\varepsilon \downarrow 0} \mathfrak{F}_{\infty}(z \! \pm \! \mi 
\varepsilon)$,
\begin{equation} \label{eql3.7b18} 
\mathfrak{F}_{\infty}(z) \underset{z \to \alpha_{k}}{=} 
\dfrac{\varkappa_{nk}}{\mi \pi nz} \! - \! \dfrac{1}{\mi \pi z} 
\sum_{m=1}^{\infty} \left(\sum_{q=1}^{\mathfrak{s}-1} \dfrac{
\varkappa_{nk \tilde{k}_{q}}}{n}(\alpha_{p_{q}})^{m} \! - \! \left(
\dfrac{(n \! - \! 1)K \! + \! k}{n} \right) \int_{J_{\infty}} \xi^{m} 
\psi_{\widetilde{V}}^{\infty}(\xi) \, \md \xi \right)z^{-m},
\end{equation}
\begin{align}
\mathfrak{F}_{\infty}(z) =& \, -\dfrac{1}{\mi \pi} \sum_{q=1}^{\mathfrak{s}
-1} \dfrac{\varkappa_{nk \tilde{k}_{q}}}{n(z \! - \! \alpha_{p_{q}})} + 
(\hat{R}(z))^{1/2} \int_{J_{\infty}} \left(\dfrac{2}{\mi \pi} \sum_{q=
1}^{\mathfrak{s}-1} \dfrac{\varkappa_{nk \tilde{k}_{q}}}{n(\xi \! - \! 
\alpha_{p_{q}})} \! + \! \dfrac{\widetilde{V}^{\prime}(\xi)}{\mi \pi} 
\right) \nonumber \\
\times& \, \dfrac{(\hat{R}(\xi))^{-1/2}_{+}}{\xi \! - \! z} \, \dfrac{\md 
\xi}{2 \pi \mi}, \quad z \! \in \! \mathbb{C} \setminus (J_{\infty} \cup 
\lbrace \alpha_{p_{1}},\alpha_{p_{2}},\dotsc,\alpha_{p_{\mathfrak{s}-1}} 
\rbrace), \label{eql3.7b19}
\end{align}
where $(\hat{R}(z))^{1/2}$ is defined by Equation~\eqref{eql3.7d}, with 
$(\hat{R}(z))^{1/2}_{\pm} \! := \! \lim_{\varepsilon \downarrow 0}
(\hat{R}(z \! \pm \! \mi \varepsilon))^{1/2}$,
\begin{equation} \label{eql3.7b20} 
\mathfrak{F}_{\infty}(z) \! = \! \dfrac{\widetilde{V}^{\prime}(z)}{2 \pi \mi} 
\! + \! \dfrac{1}{2}(\hat{R}(z))^{1/2} \oint_{\hat{C}_{\widetilde{V}}} 
\left(\dfrac{2}{\mi \pi} \sum_{q=1}^{\mathfrak{s}-1} \dfrac{\varkappa_{nk 
\tilde{k}_{q}}}{n(\xi \! - \! \alpha_{p_{q}})} \! + \! \dfrac{\widetilde{V}^{
\prime}(\xi)}{\mi \pi} \right) \dfrac{(\hat{R}(\xi))^{-1/2}}{\xi \! - \! z} \, 
\dfrac{\md \xi}{2 \pi \mi},
\end{equation}
where the contour $\hat{C}_{\widetilde{V}}$ is described in 
item~$\pmb{(1)}$ of the lemma,
\begin{gather}
\hat{\mathcal{T}}_{j} \! := \! \int_{J_{\infty}} \dfrac{\xi^{j}}{(\hat{R}
(\xi))^{1/2}_{+}} \left(\dfrac{2}{\mi \pi} \sum_{q=1}^{\mathfrak{s}-1} 
\dfrac{\varkappa_{nk \tilde{k}_{q}}}{n(\xi \! - \! \alpha_{p_{q}})} \! + \! 
\dfrac{\widetilde{V}^{\prime}(\xi)}{\mi \pi} \right) \md \xi, \quad j \! 
\in \! \mathbb{N}_{0}, \label{eql3.7b21} \\
\hat{\mathcal{N}}_{j} \! := \! \int_{\hat{a}_{j}}^{\hat{b}_{j}} \left(
(\mathcal{H} \psi_{\widetilde{V}}^{\infty})(\xi) \! - \! \dfrac{1}{2 \pi} 
\left(\dfrac{(n \! - \! 1)K \! + \! k}{n} \right)^{-1} \left(2 \sum_{q=1}^{
\mathfrak{s}-1} \dfrac{\varkappa_{nk \tilde{k}_{q}}}{n(\xi \! - \! 
\alpha_{p_{q}})} \! + \! \widetilde{V}^{\prime}(\xi) \right) \right) \md \xi, 
\quad j \! = \! 1,2,\dotsc,N, \label{eql3.7b22} \\
\dfrac{\partial \hat{\mathcal{T}}_{j}}{\partial \hat{b}_{i-1}} \! = \! 
\hat{b}_{i-1} \dfrac{\partial \hat{\mathcal{T}}_{j-1}}{\partial \hat{b}_{i-1}} 
\! + \! \dfrac{1}{2} \hat{\mathcal{T}}_{j-1}, \quad j \! \in \! \mathbb{N}, 
\label{eql3.7b23} \\
\dfrac{\partial \hat{\mathcal{T}}_{j}}{\partial \hat{a}_{i}} \! = \! 
\hat{a}_{i} \dfrac{\partial \hat{\mathcal{T}}_{j-1}}{\partial \hat{a}_{i}} 
\! + \! \dfrac{1}{2} \hat{\mathcal{T}}_{j-1}, \quad j \! \in \! \mathbb{N}, 
\label{eql3.7b24} \\
\dfrac{\partial \mathfrak{F}_{\infty}(z)}{\partial \hat{b}_{i-1}} \! = \! 
-\dfrac{1}{2 \pi \mi} \left(\dfrac{\partial \hat{\mathcal{T}}_{0}}{\partial 
\hat{b}_{i-1}} \right) \dfrac{(\hat{R}(z))^{1/2}}{z \! - \! \hat{b}_{i-1}}, 
\quad z \! \in \! \mathbb{C} \setminus (J_{\infty} \cup \lbrace \alpha_{p_{1}},
\alpha_{p_{2}},\dots,\alpha_{p_{\mathfrak{s}-1}} \rbrace), \label{eql3.7b25} \\
\dfrac{\partial \mathfrak{F}_{\infty}(z)}{\partial \hat{a}_{i}} \! = \! 
-\dfrac{1}{2 \pi \mi} \left(\dfrac{\partial \hat{\mathcal{T}}_{0}}{\partial 
\hat{a}_{i}} \right) \dfrac{(\hat{R}(z))^{1/2}}{z \! - \! \hat{a}_{i}}, \quad 
z \! \in \! \mathbb{C} \setminus (J_{\infty} \cup \lbrace \alpha_{p_{1}},
\alpha_{p_{2}},\dots,\alpha_{p_{\mathfrak{s}-1}} \rbrace), \label{eql3.7b26} \\
\dfrac{\partial \hat{\mathcal{N}}_{j}}{\partial \hat{b}_{i-1}} \! = \! 
-\dfrac{1}{2 \pi} \left(\dfrac{(n \! - \! 1)K \! + \! k}{n} \right)^{-1} 
\left(\dfrac{\partial \hat{\mathcal{T}}_{0}}{\partial \hat{b}_{i-1}} \right) 
\int_{\hat{a}_{j}}^{\hat{b}_{j}} \dfrac{(\hat{R}(\xi))^{1/2}}{\xi \! - \! 
\hat{b}_{i-1}} \, \md \xi, \quad j \! = \! 1,2,\dotsc,N, \label{eql3.7b27} \\
\dfrac{\partial \hat{\mathcal{N}}_{j}}{\partial \hat{a}_{i}} \! = \! 
-\dfrac{1}{2 \pi} \left(\dfrac{(n \! - \! 1)K \! + \! k}{n} \right)^{-1} 
\left(\dfrac{\partial \hat{\mathcal{T}}_{0}}{\partial \hat{a}_{i}} \right) 
\int_{\hat{a}_{j}}^{\hat{b}_{j}} \dfrac{(\hat{R}(\xi))^{1/2}}{\xi \! - \! 
\hat{a}_{i}} \, \md \xi, \quad j \! = \! 1,2,\dotsc,N, \label{eql3.7b28} \\
\dfrac{\partial \hat{\mathcal{T}}_{j}}{\partial \hat{b}_{i-1}} \! = \! 
(\hat{b}_{i-1})^{j} \dfrac{\partial \hat{\mathcal{T}}_{0}}{\partial 
\hat{b}_{i-1}} \quad \quad \text{and} \quad \quad \dfrac{\partial 
\hat{\mathcal{T}}_{j}}{\partial \hat{a}_{i}} \! = \! (\hat{a}_{i})^{j} 
\dfrac{\partial \hat{\mathcal{T}}_{0}}{\partial \hat{a}_{i}}, \quad 
j \! = \! 0,1,\dotsc,N \! + \! 1, \label{eql3.7b29}
\end{gather}
\begin{align}
&\operatorname{Jac}(\hat{\mathcal{T}}_{0},\hat{\mathcal{T}}_{1},\dotsc,
\hat{\mathcal{T}}_{N+1},\hat{\mathcal{N}}_{1},\hat{\mathcal{N}}_{2},\dotsc,
\hat{\mathcal{N}}_{N}) \! := \! \dfrac{\partial (\hat{\mathcal{T}}_{0},
\hat{\mathcal{T}}_{1},\dotsc,\hat{\mathcal{T}}_{N+1},\hat{\mathcal{N}}_{1},
\hat{\mathcal{N}}_{2},\dotsc,\hat{\mathcal{N}}_{N})}{\partial (\hat{b}_{0},
\hat{b}_{1},\dotsc,\hat{b}_{N},\hat{a}_{1},\hat{a}_{2},\dotsc,\hat{a}_{N+1})} 
\nonumber \\
&= 
 \nonumber \\
&= \dfrac{(-1)^{N}}{(2 \pi)^{N}} \left(\dfrac{(n \! - \! 1)K \! + \! k}{n} 
\right)^{-N} \left(\prod_{m=1}^{N+1} \dfrac{\partial \hat{\mathcal{T}}_{0}}{
\partial \hat{b}_{m-1}} \dfrac{\partial \hat{\mathcal{T}}_{0}}{\partial 
\hat{a}_{m}} \right) \left(\prod_{j=1}^{N} \int_{\hat{a}_{j}}^{\hat{b}_{j}}
(\hat{R}(\xi_{j}))^{1/2} \, \md \xi_{j} \right) \hat{\Delta}_{d}(\vec{\xi}), 
\label{eql3.7b30}
\end{align} 
where
\begin{equation*}
\hat{\Delta}_{d}(\vec{\xi}) \! = \! \dfrac{\left(\prod_{j=1}^{N+1} 
\prod_{i=1}^{N+1}(\hat{b}_{i-1} \! - \! \hat{a}_{j}) \right) \left(
\prod_{\underset{j<i}{i,j=1}}^{N+1}(\hat{b}_{i-1} \! - \! \hat{b}_{j-1})
(\hat{a}_{i} \! - \! \hat{a}_{j}) \right) \left(\prod_{\underset{j<i}{i,j
=1}}^{N}(\xi_{i} \! - \! \xi_{j}) \right)}{(-1)^{N} \prod_{j=1}^{N} 
\prod_{i=1}^{N+1}(\xi_{j} \! - \! \hat{b}_{i-1})(\xi_{j} \! - \! \hat{a}_{i})},
\end{equation*}
with
\begin{equation} \label{eql3.7b31} 
\left(\prod_{j=1}^{N} \int_{\hat{a}_{j}}^{\hat{b}_{j}}(\hat{R}(\xi_{j}))^{1/2} 
\, \md \xi_{j} \right) \hat{\Delta}_{d}(\vec{\xi}) \! \neq \! 0,
\end{equation}
and
\begin{equation} \label{eql3.7b32} 
\prod_{m=1}^{N+1} \dfrac{\partial \hat{\mathcal{T}}_{0}}{\partial 
\hat{b}_{m-1}} \dfrac{\partial \hat{\mathcal{T}}_{0}}{\partial \hat{a}_{m}} 
\! = \! \left(\dfrac{1}{2} \left(\dfrac{(n \! - \! 1)K \! + \! k}{n} \right) 
\right)^{2(N+1)} \prod_{j=1}^{N+1} \hat{h}_{\widetilde{V}}(\hat{b}_{j-1}) 
\hat{h}_{\widetilde{V}}(\hat{a}_{j}) \! \neq \! 0,
\end{equation}
where $\hat{h}_{\widetilde{V}}(z)$ is defined by Equation~\eqref{eql3.7f}. 
This concludes the proof. \hfill $\qed$
\begin{eeeee} \label{remreg} 
\textsl{It is important to make note of the following facts, the rigorous 
details of which can be supplied by using the results of Lemma~\ref{lem3.7} 
in conjunction with the ideas (but adapted to the present situation) of 
Kuijlaars-McLaughlin {\rm \cite{a57}}$:$ the details are left to the interested 
reader. Without loss of generality, consider the case $n \! \in \! \mathbb{N}$ 
and $k \! \in \! \lbrace 1,2,\dotsc,K \rbrace$ such that $\alpha_{p_{\mathfrak{s}}} 
\! := \! \alpha_{k} \! \neq \! \infty$ in the double-scaling limit $\mathscr{N},
n \! \to \! \infty$ such that $z_{o} \! \to \! z_{o}^{\ast}$ $(= \! 1)$. 
(Even though, in this monograph, $z_{o}^{\ast} \! = \! 1$, the more 
general `notation' $z_{o}^{\ast}$ is retained in order to indicate that 
the following discussion is also applicable for any bounded $z_{o}^{\ast} 
\! \in \! \mathbb{R}_{+}$.$)$ One agrees to call $z_{o}^{\ast} \! > \! 0$ a 
\emph{regular value} for $V(z)$ (satisfying conditions~\eqref{eq3}--\eqref{eq5}$)$ 
if the---scaled---external field $\widetilde{V}(z) \! := \! z_{o}^{\ast}V(z)$ 
(satisfying conditions~\eqref{eq20}--\eqref{eq22}$)$ is regular. (It turns 
out that the external field $\widetilde{V}(z) \! := \! z_{o}V(z)$ is regular 
for every $z_{o} \! > \! 0$ except for an at most countable set 
of values without a finite accumulation point; see, also, the 
\textbf{Remark to Lemma~\ref{lem3.1}}.$)$ Let $z_{o}^{\ast} \! 
> \! 0$ be a regular value for $V$, and set $J_{f}(z_{o}^{\ast}) \! 
:= \! \cup_{j=1}^{N+1}[\tilde{b}_{j-1}(z_{o}^{\ast}),\tilde{a}_{j}
(z_{o}^{\ast})]$$;$ then, $z_{o}$ is a regular value for $V$ for 
every $z_{o}$ in an open neighbourhood of $z_{o}^{\ast}$ (e.g., 
$\mathbb{U}_{\varepsilon}(z_{o}^{\ast}) \! := \! \lbrace \mathstrut 
z_{o} \! \in \! \mathbb{R}_{+}; \, \lvert z_{o} \! - \! z_{o}^{\ast} \rvert 
\! \lesssim \! \varepsilon \rbrace$, where $\varepsilon \! = \! o(1)$ 
in the double-scaling limit $\mathscr{N},n \! \to \! \infty$ such that 
$z_{o} \! \to \! z_{o}^{\ast}$ $(= \! 1))$ and $J_{f}(z_{o}) \! = \! 
\cup_{j=1}^{N+1}[\tilde{b}_{j-1}(z_{o}),\tilde{a}_{j}(z_{o})]$ with the 
\textbf{same} number of intervals (see Remark~\ref{remen} below). 
As shown in the proof of Lemma~\ref{lem3.7}, the end-points of the 
intervals, $\tilde{b}_{j-1}(z_{o}),\tilde{a}_{j}(z_{o})$, $j \! = \! 1,2,
\dotsc,N \! + \! 1$, of the support, $J_{f}$ $(= \! J_{f}(z_{o}))$, of 
the associated equilibrium measure, $\mu_{\widetilde{V}}^{f}$ 
$(= \! \mu_{\widetilde{V}}^{f}(z_{o}))$, are real-analytic (thus 
continuous) functions of $z_{o}$$;$ more is true, namely, for 
$j \! = \! 1,2,\dotsc,N \! + \! 1$, the function $z_{o} \! \mapsto 
\! \tilde{b}_{j-1}(z_{o})$ is strictly decreasing, and the function 
$z_{o} \! \mapsto \! \tilde{a}_{j}(z_{o})$ is strictly increasing; in
particular, for $j \! = \! 1,2,\dotsc,N \! + \! 1$, $[\tilde{b}_{j-1}
(z_{o} \! - \! \varepsilon),\tilde{a}_{j}(z_{o} \! - \! \varepsilon)] \subset 
[\tilde{b}_{j-1}(z_{o}),\tilde{a}_{j}(z_{o})] \subset [\tilde{b}_{j-1}
(z_{o} \! + \! \varepsilon),\tilde{a}_{j}(z_{o} \! + \! \varepsilon)]$. 
For $j \! = \! 1,2,\dotsc,N \! + \! 1$, set $(J_{f}(z_{o}^{\ast} \! + 
\! r \varepsilon))_{j} \! := \! [\tilde{b}_{j-1}(z_{o}^{\ast} \! + \! r 
\varepsilon),\tilde{a}_{j}(z_{o}^{\ast} \! + \! r \varepsilon)]$, $r \! = 
\! 0,\pm 1$. It can be shown via the proof of Lemma~\ref{lem3.7} 
that, for $z_{o} \! \in \! \mathbb{U}_{\varepsilon}(z_{o}^{\ast})$ 
(re-introducing explicit $z_{o}$- and $z_{o}^{\ast}$-dependencies$)$$:$ 
{\rm (i)} since, for $j \! = \! 1,2,\dotsc,N \! + \! 1$, $\psi_{\widetilde{V}}^{f}
(z_{o}^{\ast},\tau) \! \sim \! (\tau \! - \! \tilde{b}_{j-1}(z_{o}^{\ast}))^{1/2}$ 
as $\tau \! \downarrow \! \tilde{b}_{j-1}(z_{o}^{\ast})$, and 
$\psi_{\widetilde{V}}^{f}(z_{o}^{\ast},\tau) \! \sim \! (\tilde{a}_{j}(z_{o}^{\ast}) 
\! - \! \tau)^{1/2}$ as $\tau \! \uparrow \! \tilde{a}_{j}(z_{o}^{\ast})$, it 
follows that $\tilde{b}_{j-1}(z_{o}^{\ast} \! - \! \varepsilon) \! - \! \tilde{b}_{j-1}
(z_{o}^{\ast}) \! \sim \! \varepsilon$ as $\varepsilon \! \downarrow \! 0$ (or 
$\tilde{b}_{j-1}(z_{o}) \! - \! \tilde{b}_{j-1}(z_{o}^{\ast}) \! \sim \! (z_{o}^{\ast} 
\! - \! z_{o})^{1}$ as $z_{o} \! \uparrow \! z_{o}^{\ast}$$)$ and $\tilde{a}_{j}
(z_{o}^{\ast}) \! - \! \tilde{a}_{j}(z_{o}^{\ast} \! - \! \varepsilon) \! 
\sim \! \varepsilon$ as $\varepsilon \! \downarrow \! 0$ (or $\tilde{a}_{j}
(z_{o}^{\ast}) \! - \! \tilde{a}_{j}(z_{o}) \! \sim \! (z_{o}^{\ast} \! 
- \! z_{o})^{1}$ as $z_{o} \! \uparrow \! z_{o}^{\ast}$$)$, respectively; 
and {\rm (ii)} since, for $j \! = \! 1,2,\dotsc,N \! + \! 1$, 
$\psi_{\widetilde{V}}^{f}(z_{o}^{\ast},\tau) \! \sim \! (\tilde{b}_{j-1}
(z_{o}^{\ast}) \! - \! \tau)^{1/2}$ as $\tau \! \uparrow \! \tilde{b}_{j-1}
(z_{o}^{\ast})$, and $\psi_{\widetilde{V}}^{f}(z_{o}^{\ast},\tau) \! \sim \! 
(\tau \! - \! \tilde{a}_{j}(z_{o}^{\ast}))^{1/2}$ as $\tau \! \downarrow \! 
\tilde{a}_{j}(z_{o}^{\ast})$, it follows that $\tilde{b}_{j-1}(z_{o}^{\ast}) 
\! - \! \tilde{b}_{j-1}(z_{o}^{\ast} \! + \! \varepsilon) \! \sim \! 
\varepsilon$ as $\varepsilon \! \downarrow \! 0$ (or $\tilde{b}_{j-1}
(z_{o}^{\ast}) \! - \! \tilde{b}_{j-1}(z_{o}) \! \sim \! (z_{o} \! - \! 
z_{o}^{\ast})^{1}$ as $z_{o} \! \downarrow \! z_{o}^{\ast}$$)$ and 
$\tilde{a}_{j}(z_{o}^{\ast} \! + \! \varepsilon) \! - \! \tilde{a}_{j}
(z_{o}^{\ast}) \! \sim \! \varepsilon$ as $\varepsilon \! \downarrow \! 0$ 
(or $\tilde{a}_{j}(z_{o}) \! - \! \tilde{a}_{j}(z_{o}^{\ast}) \! \sim \! 
(z_{o} \! - \! z_{o}^{\ast})^{1}$ as $z_{o} \! \downarrow \! z_{o}^{\ast}$$)$, 
respectively. As a consequence: {\rm (i)} since, for $j \! = \! 1,2,\dotsc,N 
\! + \! 1$, $(J_{f}(z_{o}^{\ast}))_{j} \setminus (J_{f}(z_{o}^{\ast} \! - \! 
\varepsilon))_{j} \! = \! [\tilde{b}_{j-1}(z_{o}^{\ast}),\tilde{b}_{j-1}
(z_{o}^{\ast} \! - \! \varepsilon)) \cup (\tilde{a}_{j}(z_{o}^{\ast} \! - \! 
\varepsilon),\tilde{a}_{j}(z_{o}^{\ast})]$, $J_{f}(z_{o}^{\ast}) \! = \! 
\cup_{j=1}^{N+1}(J_{f}(z_{o}^{\ast}))_{j}$, and $J_{f}(z_{o}^{\ast} \! - \! 
\varepsilon) \! = \! \cup_{j=1}^{N+1}(J_{f}(z_{o}^{\ast} \! - \! 
\varepsilon))_{j}$, it follows that $\int_{J_{f}(z_{o}^{\ast}) \setminus 
J_{f}(z_{o}^{\ast}-\varepsilon)} \psi_{\widetilde{V}}^{f}(z_{o}^{\ast},\tau) 
\, \md \tau \! \leqslant \! \sum_{j=1}^{N+1}(\int_{\tilde{b}_{j-1}
(z_{o}^{\ast})}^{\tilde{b}_{j-1}(z_{o}^{\ast}-\varepsilon)} 
\psi_{\widetilde{V}}^{f}(z_{o}^{\ast},\tau) \, \md \tau \! + \! 
\int_{\tilde{a}_{j}(z_{o}^{\ast}-\varepsilon)}^{\tilde{a}_{j}(z_{o}^{\ast})} 
\psi_{\widetilde{V}}^{f}(z_{o}^{\ast},\tau) \, \md \tau) \! \lesssim \! 
\sum_{j=1}^{N+1}(\int_{\tilde{b}_{j-1}(z_{o}^{\ast})}^{\tilde{b}_{j-1}
(z_{o}^{\ast})+ \varepsilon}(\tau \! - \! \tilde{b}_{j-1}(z_{o}^{\ast}))^{1/2} 
\, \md \tau \! + \! \int_{\tilde{a}_{j}(z_{o}^{\ast})-\varepsilon}^{
\tilde{a}_{j}(z_{o}^{\ast})}(\tilde{a}_{j}(z_{o}^{\ast}) \! - \! \tau)^{1/2} 
\, \md \tau) \! \lesssim \! \mathcal{O}(\varepsilon^{3/2}) \! = \! 
o(\varepsilon)$ as $\varepsilon \! \downarrow \! 0$$;$ and {\rm (ii)} 
since, for $j \! = \! 1,2,\dotsc,N \! + \! 1$, $(J_{f}(z_{o}^{\ast} \! + \! 
\varepsilon))_{j} \setminus (J_{f}(z_{o}^{\ast}))_{j} \! = \! [\tilde{b}_{j-1}
(z_{o}^{\ast} \! + \! \varepsilon),\tilde{b}_{j-1}(z_{o}^{\ast})) \cup 
(\tilde{a}_{j}(z_{o}^{\ast}),\tilde{a}_{j}(z_{o}^{\ast} \! + \! 
\varepsilon)]$, and $J_{f}(z_{o}^{\ast} \! + \! \varepsilon) \! = \! 
\cup_{j=1}^{N+1}(J_{f}(z_{o}^{\ast} \! + \! \varepsilon))_{j}$, it follows 
that $\int_{J_{f}(z_{o}^{\ast}+\varepsilon) \setminus J_{f}(z_{o}^{\ast})} 
\psi_{\widetilde{V}}^{f}(z_{o}^{\ast},\tau) \, \md \tau \! \leqslant \! 
\sum_{j=1}^{N+1}(\int_{\tilde{b}_{j-1}(z_{o}^{\ast}+\varepsilon)}^{
\tilde{b}_{j-1}(z_{o}^{\ast})} \psi_{\widetilde{V}}^{f}(z_{o}^{\ast},\tau) 
\, \md \tau \! + \! \int_{\tilde{a}_{j}(z_{o}^{\ast})}^{\tilde{a}_{j}
(z_{o}^{\ast}+\varepsilon)} \psi_{\widetilde{V}}^{f}(z_{o}^{\ast},\tau) \, 
\md \tau) \! \lesssim \! \sum_{j=1}^{N+1}(\int_{\tilde{b}_{j-1}(z_{o}^{\ast})
-\varepsilon}^{\tilde{b}_{j-1}(z_{o}^{\ast})}(\tilde{b}_{j-1}(z_{o}^{\ast}) 
\! - \! \tau)^{1/2} \, \md \tau \! + \! \int_{\tilde{a}_{j}(z_{o}^{\ast})}^{
\tilde{a}_{j}(z_{o}^{\ast})+\varepsilon}(\tau \! - \! \tilde{a}_{j}
(z_{o}^{\ast}))^{1/2} \, \md \tau) \! \lesssim \! \mathcal{O}
(\varepsilon^{3/2}) \! = \! o(\varepsilon)$ as $\varepsilon \! \downarrow 
\! 0$.\footnote{Within the context of the study of limit distributions of 
eigenvalues (equivalently, phase transitions) in---Hermitian---random 
matrix models with a polynomial (resp., rational) external field (potential), 
the authors in \cite{amfroear} (resp., \cite{riejszara}) study, in detail, the 
evolution (variation) of one-parameter families of equilibrium measures 
on $\mathbb{R}$ when the `total mass', $[0,+\infty) \! \ni \! t \! := \! 
\lambda_{t}(\mathbb{R})$, which is the principal parameter (see below), 
varies {}from $0$ to $+\infty$ in the presence of a polynomial (resp., 
rational) external field; in particular, they \cite{amfroear} (resp., 
\cite{riejszara}) introduce and investigate an autonomous non-linear 
system of ODEs governing the dynamics of the end-points of the support 
of the associated equilibrium measure as the parameter $t$ increases 
{}from $0$ to $+\infty$: also studied are the dynamics of the equilibrium 
measure, and its support, when parameters other than the total mass, 
$t$, are allowed to vary (see, also, the detailed analysis in Section~4 of 
\cite{ieaezgmlavl}).}}
\end{eeeee}
\begin{eeeee} \label{remen}
\textsl{While studying, in the double-scaling limit $\mathscr{N},n \! \to \! 
\infty$ such that $z_{o} \! = \! 1 \! + \! o(1)$, the respective numbers of 
intervals constituting the supports of the associated equilibrium measures 
$\mu_{\widetilde{V}}^{\infty}$ (for the case $n \! \in \! \mathbb{N}$ and 
$k \! \in \! \lbrace 1,2,\dotsc,K \rbrace$ such that $\alpha_{p_{\mathfrak{s}}} 
\! := \! \alpha_{k} \! = \! \infty)$ and $\mu_{\widetilde{V}}^{f}$ (for the case 
$n \! \in \! \mathbb{N}$ and $k \! \in \! \lbrace 1,2,\dotsc,K \rbrace$ such 
that $\alpha_{p_{\mathfrak{s}}} \! := \! \alpha_{k} \! \neq \! \infty)$, one does 
not know, \emph{a priori}, that they are the same (the common value for the 
number of intervals is denoted by $N \! + \! 1$ in this monograph$)$$;$ rather, 
this fact, which requires a lengthy---asymptotic---analysis similar to that 
contained in {\rm \cite{a57}}, must be established. Strictly speaking, for the case 
$n \! \in \! \mathbb{N}$ and $k \! \in \! \lbrace 1,2,\dotsc,K \rbrace$ such 
that $\alpha_{p_{\mathfrak{s}}} \! := \! \alpha_{k} \! = \! \infty$ (resp., 
$\alpha_{p_{\mathfrak{s}}} \! := \! \alpha_{k} \! \neq \! \infty)$, one denotes 
the corresponding number of intervals by $\hat{N} \! + \! 1$ (resp., $\tilde{N} 
\! + \! 1)$, in which case $J_{\infty} \! := \! \supp (\mu_{\widetilde{V}}^{\infty}) 
\! = \! \cup_{j=1}^{\hat{N}+1}[\hat{b}_{j-1}(z_{o}),\hat{a}_{j}(z_{o})]$ (resp., 
$J_{f} \! := \! \supp (\mu_{\widetilde{V}}^{f}) \! = \! \cup_{j=1}^{\tilde{N}
+1}[\tilde{b}_{j-1}(z_{o}),\tilde{a}_{j}(z_{o})])$, and proceeds as in the 
proof of Lemma~\ref{lem3.7} to show that the end-points of the intervals 
of the supports, that is, $\hat{b}_{j-1}(z_{o}),\hat{a}_{j}(z_{o})$, $j \! = \! 
1,2,\dotsc,\hat{N} \! + \! 1$ (resp., $\tilde{b}_{j-1}(z_{o}),\tilde{a}_{j}(z_{o})$, 
$j \! = \! 1,2,\dotsc,\tilde{N} \! + \! 1)$, satisfy the system of (transcendental) 
moment equations~\eqref{eql3.7a}--\eqref{eql3.7c} (resp., 
\eqref{eql3.7g}--\eqref{eql3.7i}$)$ with the change $N \! \to \! \hat{N}$ 
(resp., $N \! \to \! \tilde{N})$ everywhere. Writing out, in the double-scaling 
limit $\mathscr{N},n \! \to \! \infty$ such that $z_{o} \! = \! 1 \! + \! o(1)$, 
the above-described moment equations in full, one arrives at: {\rm (i)} for 
$n \! \in \! \mathbb{N}$ and $k \! \in \! \lbrace 1,2,\dotsc,K \rbrace$ such 
that $\alpha_{p_{\mathfrak{s}}} \! := \! \alpha_{k} \! = \! \infty$,
\begin{align*}
0 \underset{\underset{z_{o}=1+o(1)}{\mathscr{N},n \to \infty}}{=}& \, 
\int_{J_{\infty}} \dfrac{\xi^{j}}{(\hat{R}(\xi))^{1/2}_{+}} \left(
\dfrac{2}{\mi \pi} \sum_{q=1}^{\mathfrak{s}-1} \dfrac{\gamma_{i(q)_{k_{q}}}}{
\xi \! - \! \alpha_{p_{q}}} \! + \! \dfrac{V^{\prime}(\xi)}{\mi \pi} \right) 
\md \xi \! + \! \dfrac{1}{n} \left(\int_{J_{\infty}} \dfrac{\xi^{j}}{(\hat{R}(\xi))^{
1/2}_{+}} \dfrac{2}{\mi \pi} \sum_{q=1}^{\mathfrak{s}-1} \dfrac{\widehat{
\digamma}_{k}(q)}{\xi \! - \! \alpha_{p_{q}}} \, \md \xi \right) \\
+& \, o \left(\int_{J_{\infty}} \dfrac{\xi^{j}}{(\hat{R}(\xi))^{1/2}_{+}} 
\dfrac{V^{\prime}(\xi)}{\mi \pi} \, \md \xi \right), \quad j \! = \! 0,1,
\dotsc,\hat{N}, \\
-2K \! + \! \dfrac{2(K \! - \! k)}{n} \underset{\underset{z_{o}=1+o(1)}{
\mathscr{N},n \to \infty}}{=}& \, \int_{J_{\infty}} \dfrac{\xi^{\hat{N}+1}}{
(\hat{R}(\xi))^{1/2}_{+}} \left(\dfrac{2}{\mi \pi} \sum_{q=1}^{\mathfrak{s}
-1} \dfrac{\gamma_{i(q)_{k_{q}}}}{\xi \! - \! \alpha_{p_{q}}} \! + \! 
\dfrac{V^{\prime}(\xi)}{\mi \pi} \right) \md \xi \! + \! \dfrac{1}{n} \left(
\int_{J_{\infty}} \dfrac{\xi^{\hat{N}+1}}{(\hat{R}(\xi))^{1/2}_{+}} \dfrac{2}{\mi \pi} 
\sum_{q=1}^{\mathfrak{s}-1} \dfrac{\widehat{\digamma}_{k}(q)}{\xi \! - \! 
\alpha_{p_{q}}} \, \md \xi \right) \\
+& \, o \left(\int_{J_{\infty}} \dfrac{\xi^{\hat{N}+1}}{(\hat{R}(\xi))^{1/2}_{+}} 
\dfrac{V^{\prime}(\xi)}{\mi \pi} \, \md \xi \right), \\
& \, \int_{\hat{a}_{j}(z_{o})}^{\hat{b}_{j}(z_{o})}(\hat{R}(\varsigma))^{1/2} 
\left(\int_{J_{\infty}} \dfrac{(\hat{R}(\xi))^{-1/2}_{+}}{\xi \! - \! 
\varsigma} \left(\dfrac{2}{\mi \pi} \sum_{q=1}^{\mathfrak{s}-1} \dfrac{
\gamma_{i(q)_{k_{q}}}}{\xi \! - \! \alpha_{p_{q}}} \! + \! \dfrac{V^{\prime}
(\xi)}{\mi \pi} \right) \md \xi \right) \md \varsigma \\
+& \, \dfrac{1}{n} \left(\int_{\hat{a}_{j}(z_{o})}^{\hat{b}_{j}(z_{o})}
(\hat{R}(\varsigma))^{1/2} \left(\int_{J_{\infty}} \dfrac{(\hat{R}(\xi))^{-
1/2}_{+}}{\xi \! - \! \varsigma} \dfrac{2}{\mi \pi} \sum_{q=1}^{\mathfrak{s}
-1} \dfrac{\widehat{\digamma}_{k}(q)}{\xi \! - \! \alpha_{p_{q}}} \, 
\md \xi \right) \md \varsigma \right) \\
+& \, o \left(\int_{\hat{a}_{j}(z_{o})}^{\hat{b}_{j}(z_{o})}(\hat{R}
(\varsigma))^{1/2} \left(\int_{J_{\infty}} \dfrac{(\hat{R}(\xi))^{-1/2}_{+}}{
\xi \! - \! \varsigma} \dfrac{V^{\prime}(\xi)}{\mi \pi} \, \md \xi \right) 
\md \varsigma \right) \\
\, \underset{\underset{z_{o}=1+o(1)}{\mathscr{N},n \to \infty}}{=}& \, 
2 \sum_{q=1}^{\mathfrak{s}-1} \gamma_{i(q)_{k_{q}}} \ln \left\vert 
\dfrac{\hat{b}_{j}(z_{o}) \! - \! \alpha_{p_{q}}}{\hat{a}_{j}(z_{o}) \! 
- \! \alpha_{p_{q}}} \right\vert \! + \! V(\hat{b}_{j}(z_{o})) \! - \! 
V(\hat{a}_{j}(z_{o})) \\
+& \, \dfrac{2}{n} \sum_{q=1}^{\mathfrak{s}-1} \widehat{\digamma}_{k}
(q) \ln \left\vert \dfrac{\hat{b}_{j}(z_{o}) \! - \! \alpha_{p_{q}}}{\hat{a}_{j}
(z_{o}) \! - \! \alpha_{p_{q}}} \right\vert \! + \! o \left(V(\hat{b}_{j}(z_{o})) 
\! - \! V(\hat{a}_{j}(z_{o})) \right), \quad j \! = \! 1,2,\dotsc,\hat{N},
\end{align*}
where $(\hat{R}(z))^{1/2} \! = \! (\prod_{j=1}^{\hat{N}+1}(z \! - \! 
\hat{b}_{j-1}(z_{o}))(z \! - \! \hat{a}_{j}(z_{o})))^{1/2}$, with $(\hat{R}
(z))^{1/2}_{\pm} \! := \! \lim_{\varepsilon \downarrow 0}(\hat{R}
(z \! \pm \! \mi \varepsilon))^{1/2}$, and
\begin{equation*}
\widehat{\digamma}_{k}(q) \! = \! 
\begin{cases}
-\gamma_{i(q)_{k_{q}}}, &\text{$\mathfrak{J}_{q}(k) \! = \! \varnothing, 
\quad q \! \in \! \lbrace 1,2,\dotsc,\mathfrak{s} \! - \! 1 \rbrace$,} \\
\varrho_{\widetilde{m}_{q}(k)} \! - \! \gamma_{i(q)_{k_{q}}}, 
&\text{$\mathfrak{J}_{q}(k) \! \neq \! \varnothing, \quad q \! \in \! 
\lbrace 1,2,\dotsc,\mathfrak{s} \! - \! 1 \rbrace$$;$}
\end{cases}
\end{equation*}
and {\rm (ii)} for $n \! \in \! \mathbb{N}$ and $k \! \in \! \lbrace 1,2,\dotsc,K 
\rbrace$ such that $\alpha_{p_{\mathfrak{s}}} \! := \! \alpha_{k} \! \neq \! \infty$,
\begin{align*}
0 \underset{\underset{z_{o}=1+o(1)}{\mathscr{N},n \to \infty}}{=}& \, 
\int_{J_{f}} \dfrac{\xi^{j}}{(\tilde{R}(\xi))^{1/2}_{+}} \left(
\dfrac{2}{\mi \pi} \left(\dfrac{\gamma_{k}}{\xi \! - \! \alpha_{k}} \! + \! 
\sum_{q=1}^{\mathfrak{s}-2} \dfrac{\gamma_{i(q)_{k_{q}}}}{\xi \! - \! 
\alpha_{p_{q}}} \right) \! + \! \dfrac{V^{\prime}(\xi)}{\mi \pi} \right) 
\md \xi \! + \! o \left(\int_{J_{f}} \dfrac{\xi^{j}}{(\tilde{R}(\xi))^{1/2}_{+}} 
\dfrac{V^{\prime}(\xi)}{\mi \pi} \, \md \xi \right) \\
+& \, \dfrac{1}{n} \left(\int_{J_{f}} \dfrac{\xi^{j}}{(\tilde{R}(\xi))^{
1/2}_{+}} \dfrac{2}{\mi \pi} \left(\dfrac{\varrho_{k} \! - \! \gamma_{k} 
\! - \! 1}{\xi \! - \! \alpha_{k}} \! + \! \sum_{q=1}^{\mathfrak{s}-2} 
\dfrac{\widetilde{\digamma}_{k}(q)}{\xi \! - \! \alpha_{p_{q}}} \right) 
\md \xi \right), \quad j \! = \! 0,1,\dotsc,\tilde{N}, \\
-2K \! + \! \dfrac{2(K \! - \! k)}{n} \underset{\underset{z_{o}=1+o(1)}{
\mathscr{N},n \to \infty}}{=}& \, \int_{J_{f}} \dfrac{\xi^{\tilde{N}
+1}}{(\tilde{R}(\xi))^{1/2}_{+}} \left(\dfrac{2}{\mi \pi} \left(
\dfrac{\gamma_{k}}{\xi \! - \! \alpha_{k}} \! + \! \sum_{q=1}^{\mathfrak{s}-2} 
\dfrac{\gamma_{i(q)_{k_{q}}}}{\xi \! - \! \alpha_{p_{q}}} \right) \! + \! 
\dfrac{V^{\prime}(\xi)}{\mi \pi} \right) \md \xi \! + \! o \left(
\int_{J_{f}} \dfrac{\xi^{\tilde{N}+1}}{(\tilde{R}(\xi))^{1/2}_{+}} 
\dfrac{V^{\prime}(\xi)}{\mi \pi} \, \md \xi \right) \\
+& \, \dfrac{1}{n} \left(\int_{J_{f}} \dfrac{\xi^{\tilde{N}+1}}{(\tilde{R}
(\xi))^{1/2}_{+}} \dfrac{2}{\mi \pi} \left(\dfrac{\varrho_{k} \! - \! 
\gamma_{k} \! - \! 1}{\xi \! - \! \alpha_{k}} \! + \! \sum_{q=1}^{
\mathfrak{s}-2} \dfrac{\widetilde{\digamma}_{k}(q)}{\xi \! - \! 
\alpha_{p_{q}}} \right) \md \xi \right), \\
& \, \int_{\tilde{a}_{j}(z_{o})}^{\tilde{b}_{j}(z_{o})}(\tilde{R}
(\varsigma))^{1/2} \left(\int_{J_{f}} \dfrac{(\tilde{R}(\xi))^{-1/2}_{+}}{\xi 
\! - \! \varsigma} \left(\dfrac{2}{\mi \pi} \left(\dfrac{\gamma_{k}}{\xi \! - 
\! \alpha_{k}} \! + \! \sum_{q=1}^{\mathfrak{s}-2} \dfrac{\gamma_{i(q)_{
k_{q}}}}{\xi \! - \! \alpha_{p_{q}}} \right) \! + \! \dfrac{V^{\prime}
(\xi)}{\mi \pi} \right) \md \xi \right) \md \varsigma \\
+& \, \dfrac{1}{n} \left(\int_{\tilde{a}_{j}(z_{o})}^{\tilde{b}_{j}(z_{o})}
(\tilde{R}(\varsigma))^{1/2} \left(\int_{J_{f}} \dfrac{(\tilde{R}(\xi))^{-
1/2}_{+}}{\xi \! - \! \varsigma} \dfrac{2}{\mi \pi} \left(\dfrac{\varrho_{k} 
\! - \! \gamma_{k} \! - \! 1}{\xi \! - \! \alpha_{k}} \! + \! \sum_{q=1}^{
\mathfrak{s}-2} \dfrac{\widetilde{\digamma}_{k}(q)}{\xi \! - \! 
\alpha_{p_{q}}} \right) \md \xi \right) \md \varsigma \right) \\
+& \, o \left(\int_{\tilde{a}_{j}(z_{o})}^{\tilde{b}_{j}(z_{o})}(\tilde{R}
(\varsigma))^{1/2} \left(\int_{J_{f}} \dfrac{(\tilde{R}(\xi))^{-1/2}_{+}}{\xi 
\! - \! \varsigma} \dfrac{V^{\prime}(\xi)}{\mi \pi} \, \md \xi \right) \md 
\varsigma \right) \\
\, \underset{\underset{z_{o}=1+o(1)}{\mathscr{N},n \to \infty}}{=}& \, 
2 \gamma_{k} \ln \left\vert \dfrac{\tilde{b}_{j}(z_{o}) \! - \! \alpha_{k}}{
\tilde{a}_{j}(z_{o}) \! - \! \alpha_{k}} \right\vert \! + \! 2 \sum_{q=1}^{
\mathfrak{s}-2} \gamma_{i(q)_{k_{q}}} \ln \left\vert \dfrac{\tilde{b}_{j}
(z_{o}) \! - \! \alpha_{p_{q}}}{\tilde{a}_{j}(z_{o}) \! - \! \alpha_{p_{q}}} 
\right\vert \! + \! V(\tilde{b}_{j}(z_{o})) \! - \! V(\tilde{a}_{j}(z_{o})) \\
+& \, \dfrac{2(\varrho_{k} \! - \! \gamma_{k} \! - \! 1)}{n} \ln \left\vert 
\dfrac{\tilde{b}_{j}(z_{o}) \! - \! \alpha_{k}}{\tilde{a}_{j}(z_{o}) \! - \! 
\alpha_{k}} \right\vert \! + \! \dfrac{2}{n} \sum_{q=1}^{\mathfrak{s}-2} 
\widetilde{\digamma}_{k}(q) \ln \left\vert \dfrac{\tilde{b}_{j}(z_{o}) \! - \! 
\alpha_{p_{q}}}{\tilde{a}_{j}(z_{o}) \! - \! \alpha_{p_{q}}} \right\vert \\
+& \, o \left(V(\tilde{b}_{j}(z_{o})) \! - \! V(\tilde{a}_{j}(z_{o})) \right),  
\quad j \! = \! 1,2,\dotsc,\tilde{N},
\end{align*}
where $(\tilde{R}(z))^{1/2} \! = \! (\prod_{j=1}^{\tilde{N}+1}(z \! - \! 
\tilde{b}_{j-1}(z_{o}))(z \! - \! \tilde{a}_{j}(z_{o})))^{1/2}$, with 
$(\tilde{R}(z))^{1/2}_{\pm} \! := \! \lim_{\varepsilon \downarrow 0}
(\tilde{R}(z \! \pm \! \mi \varepsilon))^{1/2}$, and
\begin{equation*}
\widetilde{\digamma}_{k}(q) \! = \! 
\begin{cases}
-\gamma_{i(q)_{k_{q}}}, &\text{$\hat{\mathfrak{J}}_{q}(k) \! = \! \varnothing, 
\quad q \! \in \! \lbrace 1,2,\dotsc,\mathfrak{s} \! - \! 2 \rbrace$,} \\
\varrho_{\hat{m}_{q}(k)} \! - \! \gamma_{i(q)_{k_{q}}}, 
&\text{$\hat{\mathfrak{J}}_{q}(k) \! \neq \! \varnothing, \quad q \! \in \! 
\lbrace 1,2,\dotsc,\mathfrak{s} \! - \! 2 \rbrace$.}
\end{cases}
\end{equation*}
Concomitant with the results for the associated `core' energy functionals given 
in Remark~\ref{rem1.3.4}, and the above systems of $2(\hat{N} \! + \! 1)$ (for 
the case $n \! \in \! \mathbb{N}$ and $k \! \in \! \lbrace 1,2,\dotsc,K \rbrace$ 
such that $\alpha_{p_{\mathfrak{s}}} \! := \! \alpha_{k} \! = \! \infty)$ and 
$2(\tilde{N} \! + \! 1)$ (for the case $n \! \in \! \mathbb{N}$ and $k \! \in \! 
\lbrace 1,2,\dotsc,K \rbrace$ such that $\alpha_{p_{\mathfrak{s}}} \! := \! 
\alpha_{k} \! \neq \! \infty)$ associated moment equations, respectively, one 
mimicks (but adapted to the present situation) the---asymptotic---procedure 
of Kuijlaars-McLaughlin {\rm \cite{a57}} to show that, in the double-scaling limit 
$\mathscr{N},n \! \to \! \infty$ such that $z_{o} \! = \! 1 \! + \! o(1)$, $\hat{N} 
\! + \! 1 \! = \! \tilde{N} \! + \! 1 \! := \! N \! + \! 1$ $(\in \! \mathbb{N}$ and 
finite$)$$;$ furthermore, one also shows that (see, also, Subsections~{\rm 2.4} 
and~{\rm 7.3} of {\rm \cite{borotguionkozl})}, for $j \! = \! 1,2,\dotsc,N \! + \! 1$, 
$\hat{b}_{j-1}(z_{o}) \! =_{\underset{z_{o}=1+o(1)}{\mathscr{N},n \to \infty}} \! 
\hat{\mathfrak{b}}_{j-1}(1 \! + \! o(1))$, $\hat{a}_{j}(z_{o}) \! =_{\underset{z_{o}
=1+o(1)}{\mathscr{N},n \to \infty}} \! \hat{\mathfrak{a}}_{j}(1 \! + \! o(1))$, 
$\tilde{b}_{j-1}(z_{o}) \! =_{\underset{z_{o}=1+o(1)}{\mathscr{N},n \to \infty}} 
\! \tilde{\mathfrak{b}}_{j-1}(1 \! + \! o(1))$, and $\tilde{a}_{j}(z_{o}) \! 
=_{\underset{z_{o}=1+o(1)}{\mathscr{N},n \to \infty}} \! \tilde{\mathfrak{a}}_{j}
(1 \! + \! o(1))$, where $\hat{\mathfrak{b}}_{j-1}$, $\hat{\mathfrak{a}}_{j}$, 
$\tilde{\mathfrak{b}}_{j-1}$, and $\tilde{\mathfrak{a}}_{j}$ are $\mathcal{O}
(1)$ (the $o(1)$ terms which appear in the latter---asymptotic---expansions 
for the respective end-points are not the same as the $o(1)$ term which appears 
in the double-scaling limit$)$$:$ the details, which are over-arching, are left 
to the interested reader.}
\end{eeeee}
\begin{eeeee} \label{bhandeghape}
\textsl{For real analytic $V$ (satisfying conditions~\eqref{eq3}--\eqref{eq5}$)$, 
the band-gap structure can be obtained, at least theoretically, via finitely 
many pointwise (not differential) equations for the end-points of the intervals 
constituting the supports of the associated equilibrium measures; in reality, 
however: {\rm (i)} these equations are highly transcendental (they are usually 
moment conditions of integrations$)$$;$ {\rm (ii)} one does not know 
\emph{a priori} the number of bands (genus), and, therefore, one does not 
know how many equations there are; and {\rm (iii)} in the case of {\rm ORFs}, 
one also needs to know how many bands there are between each pair of 
poles in order to set up equations correctly. Thus, perspicacious guesswork 
is requisite in order to commence: the difficulty, however, is that insightful 
guesswork can not be made by simply looking at the graph of $V$$;$ 
for example, a `dent' in $V$ does not guarantee an interval of support 
of the corresponding equilibrium measure. These issues can be resolved 
through numerical simulations; for example, for the $K \! = \! 3$ pole 
set $\lbrace \alpha_{1},\alpha_{2},\alpha_{3} \rbrace \! = \! \lbrace 0,3,
\infty \rbrace$ with an external field of the form $V(z) \! = \! (z \! - 
\! 3)^{-2} \! + \! z^{-2} \! + \! z^{2}$, taking the queue {}from the 
associated `core' energy functionals given in Remark~\ref{rem1.3.4} 
(cf. Equations~\eqref{stochinf} and~\eqref{stochfin}$)$, that is, 
$\mathfrak{X}^{\infty}_{\widetilde{V}}[\mu^{\text{\tiny {\rm EQ}}}] \! 
= \! \mathfrak{X}^{f}_{\widetilde{V}}[\mu^{\text{\tiny {\rm EQ}}}]$, 
the corresponding finite-particle energy function reads
\begin{equation*}
-\dfrac{1}{(\mathfrak{n}(\#))^{2}} \sum_{\substack{i,j=1 \\i \neq j}}^{
\mathfrak{n}(\#)} \ln \left(\dfrac{\lvert z_{i} \! - \! z_{j} \rvert^{3}}{\lvert 
z_{i}z_{j}(z_{i} \! - \! 3)(z_{j} \! - \! 3) \rvert} \right)\! + \! \dfrac{1}{
\mathfrak{n}(\#)} \sum_{m=1}^{\mathfrak{n}(\#)}V(z_{m}),
\end{equation*}
where $\mathfrak{n}(\#)$ denotes the `number of particles'.\footnote{In 
this case, the discrete measure is $\md \mu_{\mathfrak{n}(\#)}(z) \! = \! 
\tfrac{1}{\mathfrak{n}(\#)} \sum_{j=1}^{\mathfrak{n}(\#)} \delta (z \! - \! 
z_{j}) \, \md z$; as $\mathfrak{n}(\#) \! \to \! \infty$, $\mu_{\mathfrak{n}
(\#)}$ converges to a Borel measure in the distribution(al) sense (or 
weak-$\ast$ topology over continuous functions with compact support).} 
Using $\mathfrak{n}(\#) \! = \! 500$, the simulated equilibrium measure 
is supported on two intervals, $[-2.107,-0.506] \cup [0.500,2.273]$, with 
global minimum energy equal to $5.408$$;$ but, modifying slightly the 
external field to $V(z) \! = \! (z \! - \! 3)^{-2} \! + \! z^{-2} \! + \! z^{2}/2$, 
the numerically simulated equilibrium measure is supported on three intervals, 
$[-2.800,-0.544] \cup [0.536,2.412] \cup [3.696,3.800]$, with global minimum 
energy equal to $4.43$. These results can be used to correctly set up equations 
in the form of moment conditions, and more accurate numerical results can be 
obtained {}from them; once the end-points of the intervals of the support are 
obtained, the corresponding equilibrium measure can be calculated via the 
so-called `$g$-function method'.}
\end{eeeee}
\begin{eeeee} \label{remstielt}
\textsl{For the case $n \! \in \! \mathbb{N}$ and $k \! \in \! \lbrace 
1,2,\dotsc,K \rbrace$ such that $\alpha_{p_{\mathfrak{s}}} \! := \! 
\alpha_{k} \! \neq \! \infty$, perusing Equations~\eqref{eql3.7a7} 
and~\eqref{eql3.7a8} (see, also, Equations~\eqref{eql3.7v}, 
\eqref{eql3.7w}, \eqref{eql3.7y}--\eqref{eql3.7a1}, \eqref{eql3.7a3}, 
and~\eqref{eql3.7a5}$)$ for the real-valued coefficients $\tilde{\rho}_{2 
\mathfrak{s}-3}(n,k)$ and $\tilde{\rho}_{m}(n,k)$, $m \! = \! 0,1,\dotsc,
2(\mathfrak{s} - \! 2)$, respectively, one notes that `moment integrals' 
of the type $\int_{J_{f}} \xi^{r_{1}}(\widetilde{V}^{\prime}(\xi))^{r_{2}} \, 
\md \mu_{\widetilde{V}}^{f}(\xi)$, $r_{1} \! \in \! \mathbb{N}_{0}$, 
$r_{2} \! = \! 0,1$, are encountered; in fact, one may use such moment 
integrals in order to evaluate the Stieltjes transform of the associated 
equilibrium measure, $\mu_{\widetilde{V}}^{f}$ $(\in \! \mathscr{M}_{1}
(\mathbb{R}))$, on the corresponding real pole set $\lbrace \alpha_{p_{1}},
\dotsc,\alpha_{p_{\mathfrak{s}-2}},\alpha_{p_{\mathfrak{s}}} \rbrace$, that 
is, $\int_{J_{f}}(\xi \! - \! \alpha_{p_{q}})^{-1} \, \md \mu_{\widetilde{V}}^{f}
(\xi)$, $q \! = \! 1,\dotsc,\mathfrak{s} \! - \! 2,\mathfrak{s}$. Via 
Equations~\eqref{eql3.7q} and~\eqref{eql3.7t}, proceeding as per 
the proof of Lemma~{\rm 3.7} in {\rm \cite{a45}} (see, in particular, 
pp.~{\rm 320}--{\rm 323}, the linear system~{\rm (62)} of {\rm \cite{a45})}, 
one shows that, for $n \! \in \! \mathbb{N}$ and $k \! \in \! \lbrace 
1,2,\dotsc,K \rbrace$ such that $\alpha_{p_{\mathfrak{s}}} \! := \! 
\alpha_{k} \! \neq \! \infty$,
\begin{equation} \label{eqremm1fin} 
\underbrace{
,
\end{equation}
where
\begin{align*}
\tilde{\Gamma}(n,k) :=& \, -\dfrac{2}{\pi^{2}} \left(\dfrac{\varkappa_{nk} 
\! - \! 1}{n} \right) \prod_{q=1}^{\mathfrak{s}-2}(\alpha_{k} \! - \! 
\alpha_{p_{q}})^{2} \left(\sum_{r=1}^{\mathfrak{s}-2} \dfrac{(\varkappa_{nk} 
\! + \! \varkappa_{nk \tilde{k}_{r}} \! - \! 1)}{n(\alpha_{k} \! - \! 
\alpha_{p_{r}})} \! + \! \left(\dfrac{(n \! - \! 1)K \! + \! k}{n} \right) 
\int_{J_{f}}(\xi \! - \! \alpha_{k})^{-1} \, \md \mu^{f}_{\widetilde{V}}(\xi) 
\right), \\
\tilde{\Xi}(n,k,q) :=& \, -\dfrac{2}{\pi^{2}} \dfrac{\varkappa_{nk 
\tilde{k}_{q}}}{n}(\alpha_{p_{q}} \! - \! \alpha_{k})^{2} \prod_{
\substack{r=1\\r \neq q}}^{\mathfrak{s}-2}(\alpha_{p_{q}} \! - \! 
\alpha_{p_{r}})^{2} \left(\sum_{\substack{r=1\\r \neq q}}^{\mathfrak{s}-2} 
\dfrac{(\varkappa_{nk \tilde{k}_{q}} \! + \! \varkappa_{nk \tilde{k}_{r}})}{n
(\alpha_{p_{q}} \! - \! \alpha_{p_{r}})} \! + \! \dfrac{(\varkappa_{nk} \! 
+ \! \varkappa_{nk \tilde{k}_{q}} \! - \! 1)}{n(\alpha_{p_{q}} \! - \! 
\alpha_{k})} \right. \\
+&\left. \, \left(\dfrac{(n \! - \! 1)K \! + \! k}{n} \right) \int_{J_{f}}
(\xi \! - \! \alpha_{p_{q}})^{-1} \, \md \mu_{\widetilde{V}}^{f}(\xi) \right), 
\quad q \! = \! 1,2,\dotsc,\mathfrak{s} \! - \! 2;
\end{align*}
one shows via Theorem~{\rm 20} of {\rm \cite{a66}} that the determinant of 
the $2(\mathfrak{s} \! - \! 1) \times 2(\mathfrak{s} \! - \! 1)$ coefficient 
matrix of system~\eqref{eqremm1fin}, that is, $\tilde{D}(n,k)$, is given by
\begin{equation*}
\det (\tilde{D}(n,k)) \! = \! (-1)^{\frac{(\mathfrak{s}-1)(\mathfrak{s}-
2)}{2}} \prod_{q=1}^{\mathfrak{s}-2}(\alpha_{p_{q}} \! - \! \alpha_{k})^{4} 
\prod_{\substack{i,j=1\\i<j}}^{\mathfrak{s}-2}(\alpha_{p_{j}} \! - \! 
\alpha_{p_{i}})^{4} \neq \! 0.
\end{equation*}
For $n \! \in \! \mathbb{N}$ and $k \! \in \! \lbrace 1,2,\dotsc,K \rbrace$ 
such that $\alpha_{p_{\mathfrak{s}}} \! := \! \alpha_{k} \! \neq \! \infty$, 
the first $\mathfrak{s} \! - \! 1$ equations of system~\eqref{eqremm1fin} 
give rise to some remarkable identities (exercise!), whilst the latter 
$\mathfrak{s} \! - \! 1$ equations of system~\eqref{eqremm1fin} allow 
one to express the associated Stieltjes transforms, $\int_{J_{f}}(\xi \! - \! 
\alpha_{p_{q}})^{-1} \, \md \mu_{\widetilde{V}}^{f}(\xi)$, $q \! = \! 1,\dotsc,
\mathfrak{s} \! - \! 2,\mathfrak{s}$, in terms of the corresponding moment 
integrals $\int_{J_{f}} \xi^{r_{1}} \widetilde{V}^{\prime}(\xi) \, \md 
\mu_{\widetilde{V}}^{f}(\xi)$, $r_{1} \! = \! 0,1,\dotsc,2 \mathfrak{s} \! 
- \! 3$, and $\int_{J_{f}} \xi^{r_{2}} \, \md \mu_{\widetilde{V}}^{f}(\xi)$, 
$r_{2} \! = \! 0,1,\dotsc,2(\mathfrak{s} \! - \! 2)$$:$ as a by-product, 
this shows, incidentally, that the associated Stieltjes transforms, 
$\int_{J_{f}}(\xi \! - \! \alpha_{p_{q}})^{-1} \, \md \mu_{\widetilde{V}}^{f}
(\xi)$, $q \! = \! 1,\dotsc,\mathfrak{s} \! - \! 2,\mathfrak{s}$, are real 
valued and bounded. Alternatively, given the associated Stieltjes transforms, 
$\int_{J_{f}}(\xi \! - \! \alpha_{p_{q}})^{-1} \, \md \mu_{\widetilde{V}}^{f}(\xi)$, 
$q \! = \! 1,\dotsc,\mathfrak{s} \! - \! 2,\mathfrak{s}$, one may regard 
system~\eqref{eqremm1fin} as a means by which real linear combinations 
of the corresponding moment integrals $\int_{J_{f}} \xi^{r_{1}} 
\widetilde{V}^{\prime}(\xi) \, \md \mu_{\widetilde{V}}^{f}(\xi)$, $r_{1} \! 
= \! 0,1,\dotsc,2 \mathfrak{s} \! - \! 3$, and $\int_{J_{f}} \xi^{r_{2}} \, 
\md \mu_{\widetilde{V}}^{f}(\xi)$, $r_{2} \! = \! 0,1,\dotsc,2(\mathfrak{s} 
\! - \! 2)$, may be expressed in terms of real linear combinations of the 
associated Stieltjes transforms. Similarly, for the case $n \! \in \! 
\mathbb{N}$ and $k \! \in \! \lbrace 1,2,\dotsc,K \rbrace$ such that 
$\alpha_{p_{\mathfrak{s}}} \! := \! \alpha_{k} \! = \! \infty$, the analogue 
of system~\eqref{eqremm1fin} reads 
(cf. Equations~\eqref{eql3.7b5}--\eqref{eql3.7b14}$):$
\begin{equation} \label{eqremm1inf} 
\underbrace{
,
\end{equation}
where, for $q \! = \! 1,2,\dotsc,\mathfrak{s} \! - \! 1$,
\begin{equation*}
\hat{\Xi}(n,k,q) \! := \! -\dfrac{2}{\pi^{2}} \dfrac{\varkappa_{nk 
\tilde{k}_{q}}}{n} \prod_{\substack{r=1\\r \neq q}}^{\mathfrak{s}-1}
(\alpha_{p_{q}} \! - \! \alpha_{p_{r}})^{2} \left(\sum_{\substack{r=
1\\r \neq q}}^{\mathfrak{s}-1} \dfrac{(\varkappa_{nk \tilde{k}_{q}} \! + 
\! \varkappa_{nk \tilde{k}_{r}})}{n(\alpha_{p_{q}} \! - \! \alpha_{p_{r}})} 
\! + \! \left(\dfrac{(n \! - \! 1)K \! + \! k}{n} \right) \int_{J_{\infty}}
(\xi \! - \! \alpha_{p_{q}})^{-1} \, \md \mu_{\widetilde{V}}^{\infty}(\xi) 
\right);
\end{equation*}
one shows via Theorem~{\rm 20} of {\rm \cite{a66}} that the determinant of 
the $2(\mathfrak{s} \! - \! 1) \times 2(\mathfrak{s} \! - \! 1)$ coefficient 
matrix of system~\eqref{eqremm1inf}, that is, $\hat{D}(n,k)$, is given by
\begin{equation*}
\det (\hat{D}(n,k)) \! = \! (-1)^{\frac{(\mathfrak{s}-1)(\mathfrak{s}-2)}{2}} 
\prod_{\substack{i,j=1\\i<j}}^{\mathfrak{s}-1}(\alpha_{p_{j}} \! - \! 
\alpha_{p_{i}})^{4} \neq 0.
\end{equation*}
For $n \! \in \! \mathbb{N}$ and $k \! \in \! \lbrace 1,2,\dotsc,K \rbrace$ 
such that $\alpha_{p_{\mathfrak{s}}} \! := \! \alpha_{k} \! = \! \infty$, the 
first $\mathfrak{s} \! - \! 1$ equations of system~\eqref{eqremm1inf} give 
rise to some remarkable identities (exercise!), whilst the latter $\mathfrak{s} 
\! - \! 1$ equations of system~\eqref{eqremm1inf} allow one to express the 
associated Stieltjes transforms, $\int_{J_{\infty}}(\xi \! - \! \alpha_{p_{q}})^{-1} 
\, \md \mu_{\widetilde{V}}^{\infty}(\xi)$, $q \! = \! 1,2,\dotsc,\mathfrak{s} 
\! - \! 1$, in terms of the corresponding real-valued and bounded moment 
integrals $\int_{J_{\infty}} \xi^{r_{1}} \widetilde{V}^{\prime}(\xi) \, \md 
\mu_{\widetilde{V}}^{\infty}(\xi)$, $r_{1} \! = \! 0,1,\dotsc,2 \mathfrak{s} 
\! - \! 3$, and $\int_{J_{\infty}} \xi^{r_{2}} \, \md \mu_{\widetilde{V}}^{
\infty}(\xi)$, $r_{2} \! = \! 0,1,\dotsc,2(\mathfrak{s} \! - \! 2)$. 
Alternatively, given the associated Stieltjes transforms, $\int_{J_{\infty}}
(\xi \! - \! \alpha_{p_{q}})^{-1} \, \md \mu_{\widetilde{V}}^{\infty}(\xi)$, 
$q \! = \! 1,2,\dotsc,\mathfrak{s} \! - \! 1$, one may regard 
system~\eqref{eqremm1inf} as a means by which real linear combinations 
of the corresponding moment integrals $\int_{J_{\infty}} \xi^{r_{1}} 
\widetilde{V}^{\prime}(\xi) \, \md \mu_{\widetilde{V}}^{\infty}(\xi)$, 
$r_{1} \! = \! 0,1,\dotsc,2 \mathfrak{s} \! - \! 3$, and $\int_{J_{\infty}} 
\xi^{r_{2}} \, \md \mu_{\widetilde{V}}^{\infty}(\xi)$, $r_{2} \! = \! 0,1,
\dotsc,2(\mathfrak{s} \! - \! 2)$, may be expressed in terms of real linear 
combinations of the associated Stieltjes transforms.}
\end{eeeee}
\begin{eeeee} \label{remext} 
\textsl{For regular $\widetilde{V} \colon \overline{\mathbb{R}} \setminus 
\lbrace \alpha_{1},\alpha_{2},\dotsc,\alpha_{K} \rbrace \! \to \! \mathbb{R}$ 
satisfying conditions~\eqref{eq20}--\eqref{eq22} of the form
\begin{equation} \label{eqremext1} 
\widetilde{V}(z) \! = \! \sum_{\lbrace \mathstrut q \in \lbrace 1,2,
\dotsc,\mathfrak{s} \rbrace; \, \alpha_{p_{q}} \neq \infty \rbrace} 
\sum_{j=-2m_{q}}^{-1} \tilde{\upsilon}_{j,q}(z \! - \! \alpha_{p_{q}})^{j} 
\! + \! \sum_{j=0}^{2m_{\infty}} \hat{\upsilon}_{j,\infty}z^{j},
\end{equation}
where {}\footnote{Equation~\eqref{eqremext1} should be understood as 
follows: (i) for the case $(n,k) \! \in \! \mathbb{N} \times \lbrace 1,2,
\dotsc,K \rbrace$ such that $\alpha_{p_{\mathfrak{s}}} \! := \! \alpha_{k} 
\! = \! \infty$, $\widetilde{V}(z) \! = \! \sum_{q=1}^{\mathfrak{s}-1} 
\sum_{j=-2m_{q}}^{-1} \tilde{\upsilon}_{j,q}(z \! - \! \alpha_{p_{q}})^{j} 
\! + \! \sum_{j=0}^{2m_{\infty}} \hat{\upsilon}_{j,\infty}z^{j}$; and (ii) 
for the case $(n,k) \! \in \! \mathbb{N} \times \lbrace 1,2,\dotsc,K 
\rbrace$ such that $\alpha_{p_{\mathfrak{s}}} \! := \! \alpha_{k} \! \neq 
\! \infty$, $\widetilde{V}(z) \! = \! \sum_{\underset{q \neq \mathfrak{s}
-1}{q=1}}^{\mathfrak{s}} \sum_{j=-2m_{q}}^{-1} \tilde{\upsilon}_{j,q}
(z \! - \! \alpha_{p_{q}})^{j} \! + \! \sum_{j=0}^{2m_{\infty}} 
\hat{\upsilon}_{j,\infty}z^{j}$.} $m_{q},m_{\infty} \! \in \! \mathbb{N}$ 
and $\tilde{\upsilon}_{-2m_{q},q},\hat{\upsilon}_{2m_{\infty},\infty} 
\! > \! 0$, closed form expressions for the associated functions 
$\hat{h}_{\widetilde{V}}(z)$ and $\tilde{h}_{\widetilde{V}}(z)$ defined 
by Equations~\eqref{eql3.7f} and~\eqref{eql3.7l}, respectively, can 
be derived. Via a residue calculation, one shows that, for $n \! \in \! 
\mathbb{N}$ and $k \! \in \! \lbrace 1,2,\dotsc,K \rbrace$ such that 
$\alpha_{p_{\mathfrak{s}}} \! := \! \alpha_{k} \! = \! \infty$, and $z 
\! \in \! \mathbb{C} \setminus \lbrace \alpha_{p_{1}},\alpha_{p_{2}},
\dotsc,\alpha_{p_{\mathfrak{s}-1}} \rbrace$,
\begin{align*}
\left(\dfrac{(n \! - \! 1)K \! + \! k}{n} \right) \hat{h}_{\widetilde{V}}(z) 
\! =& \, z^{2m_{\infty}-N-2} \sum_{j=0}^{2m_{\infty}-N-2}(2m_{\infty} \! - 
\! j) \hat{\upsilon}_{2m_{\infty}-j,\infty} \underset{\substack{0 \leqslant 
\vert \hat{k} \vert + \vert \hat{l} \vert \leqslant 2m_{\infty}-j-N-2\\k_{i} 
\in \mathbb{N}_{0}, \, \, l_{i} \in \mathbb{N}_{0}, \, \, i=0,1,
\dotsc,N}}{\sideset{}{'}{\sum}_{k_{0},k_{1},\dotsc,k_{N}} \, \, 
\sideset{}{'}{\sum}_{l_{0},l_{1},\dotsc,l_{N}}} \prod_{p_{1}=0}^{N} 
\prod_{j_{p_{1}}=0}^{k_{p_{1}}-1} \left(\dfrac{1}{2} \! + \! j_{p_{1}} 
\right) \prod_{p_{2}=0}^{N} \prod_{j_{p_{2}}=0}^{l_{p_{2}}-1} 
\left(\dfrac{1}{2} \! + \! j_{p_{2}} \right) \\
\times& \, \dfrac{\prod_{p_{3}=0}^{N}(\hat{b}_{p_{3}})^{k_{p_{3}}} 
\prod_{p_{4}=0}^{N}(\hat{a}_{p_{4}+1})^{l_{p_{4}}}}{\prod_{i=0}^{N}k_{i}! 
\prod_{j=0}^{N}l_{j}!}z^{-j-\vert \hat{k} \vert -\vert \hat{l} \vert} \! - \! 
\sum_{q=1}^{\mathfrak{s}-1} 
\dfrac{(-1)^{\hat{\mathfrak{n}}(\alpha_{p_{q}})}(\prod_{i=1}^{N+1} 
\lvert \hat{b}_{i-1} \! - \! \alpha_{p_{q}} \rvert \lvert \hat{a}_{i} \! - \! 
\alpha_{p_{q}} \rvert)^{-1/2}}{(z \! - \! \alpha_{p_{q}})^{2m_{q}+1}} \\
\times& \, \sum_{j=-2m_{q}+1}^{0}(2m_{q} \! + \! j) \tilde{\upsilon}_{-
2m_{q}-j,q} \underset{\substack{0 \leqslant \vert \hat{k} \vert + \vert 
\hat{l} \vert \leqslant 2m_{q}+j\\k_{i} \in \mathbb{N}_{0}, \, \, l_{i} \in 
\mathbb{N}_{0}, \, \, i=0,1,\dotsc,N}}{\sideset{}{''}{\sum}_{k_{0},k_{1},
\dotsc,k_{N}} \, \, \sideset{}{''}{\sum}_{l_{0},l_{1},\dotsc,l_{N}}} 
\dfrac{\prod_{p_{1}=0}^{N} \prod_{j_{p_{1}}=0}^{k_{p_{1}}-1} 
\left(\frac{1}{2} \! + \! j_{p_{1}} \right)}{\prod_{i_{1}=0}^{N}
(\hat{b}_{i_{1}} \! - \! \alpha_{p_{q}})^{k_{i_{1}}} \prod_{i_{2}=0}^{N}
(\hat{a}_{i_{2}+1} \! - \! \alpha_{p_{q}})^{l_{i_{2}}}} \\
\times& \, \dfrac{\prod_{p_{2}=0}^{N} \prod_{j_{p_{2}}=0}^{l_{p_{2}}-1} 
\left(\frac{1}{2} \! + \! j_{p_{2}} \right)}{\prod_{i_{3}=0}^{N}k_{i_{3}}! 
\prod_{i_{4}=0}^{N}l_{i_{4}}!}(z \! - \! \alpha_{p_{q}})^{-j+ \vert \hat{k} 
\vert + \vert \hat{l} \vert} \! + \! 2 \sum_{q=1}^{\mathfrak{s}-1} 
\dfrac{\varkappa_{nk \tilde{k}_{q}}}{n} \dfrac{(-1)^{\hat{\mathfrak{n}}
(\alpha_{p_{q}})}(\prod_{i=1}^{N+1} \lvert \hat{b}_{i-1} \! - \! 
\alpha_{p_{q}} \rvert \lvert \hat{a}_{i} \! - \! \alpha_{p_{q}} 
\rvert)^{-1/2}}{(z \! - \! \alpha_{p_{q}})},
\end{align*}
where $\lvert \hat{k} \rvert \! := \! k_{0} \! + \! k_{1} \! + \! \dotsb \! 
+ \! k_{N}$, $\lvert \hat{l} \rvert \! := \! l_{0} \! + \! l_{1} \! + \! 
\dotsb \! + \! l_{N}$, and, for $q \! = \! 1,2,\dotsc,\mathfrak{s} \! - \! 1$,
\begin{equation*}
\hat{\mathfrak{n}}(\alpha_{p_{q}}) \! = \! 
\begin{cases}
0, &\text{$\alpha_{p_{q}} \! \in \! (\hat{a}_{N+1},+\infty)$,} \\
N \! + \! 1, &\text{$\alpha_{p_{q}} \! \in \! (-\infty,\hat{b}_{0})$,} \\
N \! - \! j \! + \! 1, &\text{$\alpha_{p_{q}} \! \in \! (\hat{a}_{j},
\hat{b}_{j}), \quad j \! = \! 1,2,\dotsc,N$,}
\end{cases}
\end{equation*}
and, for $n \! \in \! \mathbb{N}$ and $k \! \in \! \lbrace 1,2,\dotsc,
K \rbrace$ such that $\alpha_{p_{\mathfrak{s}}} \! := \! \alpha_{k} 
\! \neq \! \infty$, and $z \! \in \! \mathbb{C} \setminus \lbrace 
\alpha_{p_{1}},\dotsc,\alpha_{p_{\mathfrak{s}-2}},\alpha_{k} \rbrace$,
\begin{align*}
\left(\dfrac{(n \! - \! 1)K \! + \! k}{n} \right) \tilde{h}_{\widetilde{V}}
(z) \! =& \, z^{2m_{\infty}-N-2} \sum_{j=0}^{2m_{\infty}-N-2}(2m_{\infty} \! 
- \! j) \hat{\upsilon}_{2m_{\infty}-j,\infty} \underset{\substack{0 \leqslant 
\vert \hat{k} \vert + \vert \hat{l} \vert \leqslant 2m_{\infty}-j-N-2\\k_{i} 
\in \mathbb{N}_{0}, \, \, l_{i} \in \mathbb{N}_{0}, \, \, i=0,1,\dotsc,N}}{
\sideset{}{'}{\sum}_{k_{0},k_{1},\dotsc,k_{N}} \, \, \sideset{}{'}{\sum}_{l_{0},
l_{1},\dotsc,l_{N}}} \prod_{p_{1}=0}^{N} \prod_{j_{p_{1}}=0}^{k_{p_{1}}-1} 
\left(\dfrac{1}{2} \! + \! j_{p_{1}} \right) \prod_{p_{2}=0}^{N} 
\prod_{j_{p_{2}}=0}^{l_{p_{2}}-1} \left(\dfrac{1}{2} \! + \! j_{p_{2}} \right) \\
\times& \, \dfrac{\prod_{p_{3}=0}^{N}(\tilde{b}_{p_{3}})^{k_{p_{3}}} 
\prod_{p_{4}=0}^{N}(\tilde{a}_{p_{4}+1})^{l_{p_{4}}}}{\prod_{i=0}^{N}k_{i}! 
\prod_{j=0}^{N}l_{j}!}z^{-j-\vert \hat{k} \vert -\vert \hat{l} \vert} \! - \! 
\sum_{\substack{q=1\\q \neq \mathfrak{s}-1}}^{\mathfrak{s}} 
\dfrac{(-1)^{\tilde{\mathfrak{n}}(\alpha_{p_{q}})}(\prod_{i=1}^{N+1} 
\lvert \tilde{b}_{i-1} \! - \! \alpha_{p_{q}} \rvert \lvert \tilde{a}_{i} \! - 
\! \alpha_{p_{q}} \rvert)^{-1/2}}{(z \! - \! \alpha_{p_{q}})^{2m_{q}+1}} \\
\times& \, \sum_{j=-2m_{q}+1}^{0}(2m_{q} \! + \! j) \tilde{\upsilon}_{-2
m_{q}-j,q} \underset{\substack{0 \leqslant \vert \hat{k} \vert + \vert 
\hat{l} \vert \leqslant 2m_{q}+j\\k_{i} \in \mathbb{N}_{0}, \, \, l_{i} \in 
\mathbb{N}_{0}, \, \, i=0,1,\dotsc,N}}{\sideset{}{''}{\sum}_{k_{0},k_{1},
\dotsc,k_{N}} \, \, \sideset{}{''}{\sum}_{l_{0},l_{1},\dotsc,l_{N}}} 
\dfrac{\prod_{p_{1}=0}^{N} \prod_{j_{p_{1}}=0}^{k_{p_{1}}-1} \left(
\frac{1}{2} \! + \! j_{p_{1}} \right)}{\prod_{i_{1}=0}^{N}(\tilde{b}_{i_{1}} 
\! - \! \alpha_{p_{q}})^{k_{i_{1}}} \prod_{i_{2}=0}^{N}(\tilde{a}_{i_{2}+1} 
\! - \! \alpha_{p_{q}})^{l_{i_{2}}}} \\
\times& \, \dfrac{\prod_{p_{2}=0}^{N} \prod_{j_{p_{2}}=0}^{l_{p_{2}}-1} 
\left(\frac{1}{2} \! + \! j_{p_{2}} \right)}{\prod_{i_{3}=0}^{N}k_{i_{3}}! 
\prod_{i_{4}=0}^{N}l_{i_{4}}!}(z \! - \! \alpha_{p_{q}})^{-j+ \vert \hat{k} 
\vert + \vert \hat{l} \vert} \! + \! 2 \sum_{q=1}^{\mathfrak{s}-2} 
\dfrac{\varkappa_{nk \tilde{k}_{q}}}{n} \dfrac{(-1)^{\tilde{\mathfrak{n}}
(\alpha_{p_{q}})}(\prod_{i=1}^{N+1} \lvert \tilde{b}_{i-1} \! - \! 
\alpha_{p_{q}} \rvert \lvert \tilde{a}_{i} \! - \! \alpha_{p_{q}} 
\rvert)^{-1/2}}{(z \! - \! \alpha_{p_{q}})} \\
+& \, 2 \left(\dfrac{\varkappa_{nk} \! - \! 1}{n} \right) 
\dfrac{(-1)^{\tilde{\mathfrak{n}}(\alpha_{k})}(\prod_{i=1}^{N+1} 
\lvert \tilde{b}_{i-1} \! - \! \alpha_{k} \rvert \lvert \tilde{a}_{i} \! 
- \! \alpha_{k} \rvert)^{-1/2}}{(z \! - \! \alpha_{k})},
\end{align*}
where, for $q \! = \! 1,\dotsc,\mathfrak{s} \! - \! 2,\mathfrak{s}$,
\begin{equation*}
\tilde{\mathfrak{n}}(\alpha_{p_{q}}) \! = \! 
\begin{cases}
0, &\text{$\alpha_{p_{q}} \! \in \! (\tilde{a}_{N+1},+\infty)$,} \\
N \! + \! 1, &\text{$\alpha_{p_{q}} \! \in \! (-\infty,\tilde{b}_{0})$,} \\
N \! - \! j \! + \! 1, &\text{$\alpha_{p_{q}} \! \in \! (\tilde{a}_{j},
\tilde{b}_{j}), \quad j \! = \! 1,2,\dotsc,N$,}
\end{cases}
\end{equation*}
and the primes (resp., double primes) on the respective summations mean 
that all possible sums over $\lbrace k_{0},k_{1},\dotsc,k_{N} \rbrace$ and 
$\lbrace l_{0},l_{1},\dotsc,l_{N} \rbrace$ must be taken for which $0 \! 
\leqslant \! \sum_{m=0}^{N}(k_{m} \! + \! l_{m}) \! \leqslant \! 2m_{\infty} 
\! - \! j \! - \! N \! - \! 2$, $j \! = \! 0,1,\dotsc,2m_{\infty} \! - \! N \! - \! 2$, 
$k_{i} \! \in \! \mathbb{N}_{0}$, $l_{i} \! \in \! \mathbb{N}_{0}$, $i \! = \! 0,1,
\dotsc,N$ (resp., $0 \! \leqslant \! \sum_{m=0}^{N}(k_{m} \! + \! l_{m}) \! 
\leqslant \! 2m_{q} \! + \! j$, $j \! = \! -2m_{q} \! + \! 1,-2m_{q} \! + \! 2,
\dotsc,0$, $k_{i} \! \in \! \mathbb{N}_{0}$, $l_{i} \! \in \! \mathbb{N}_{0}$, 
$i \! = \! 0,1,\dotsc,N$, where, for $k \! \in \! \lbrace 1,2,\dotsc,K \rbrace$ 
such that $\alpha_{p_{\mathfrak{s}}} \! := \! \alpha_{k} \! = \! \infty$, $q \! 
= \! 1,2,\dotsc,\mathfrak{s} \! - \! 1$, and, for $k \! \in \! \lbrace 1,2,\dotsc,
K \rbrace$ such that $\alpha_{p_{\mathfrak{s}}} \! := \! \alpha_{k} \! \neq \! 
\infty$, $q \! = \! 1,\dotsc,\mathfrak{s} \! - \! 2,\mathfrak{s})$. It is important 
to note that all of the above sums are finite sums: any sums for which the 
upper limit is less than the lower limit are defined to be equal to zero (e.g., 
$\sum_{j=0}^{-1} \boldsymbol{\ast}(j) \! := \! 0)$, and any products for 
which the upper limit is less than the lower limit are defined to be equal 
to one (e.g., $\prod_{j=0}^{-1} \boldsymbol{\ast}(j) \! := \! 1)$.}
\end{eeeee}

The following---lengthy and technical---Lemma~\ref{lemrootz}, which is modelled 
on the calculations in Chapter~6 of \cite{a51} and Section~4 of \cite{a56} (see, also, 
\cite{borotguionkozl,shr1,shr7,shrgbag,whylee,cfagtysgaw,peiiechr,shr2,shr3,shr6,
shr4,shr5,newmscba1}), shows that, for $n \! \in \! \mathbb{N}$ and $k \! \in \! 
\lbrace 1,2,\dotsc,K \rbrace$ such that $\alpha_{p_{\mathfrak{s}}} \! := \! \alpha_{k} 
\! = \! \infty$ (resp., $\alpha_{p_{\mathfrak{s}}} \! := \! \alpha_{k} \! \neq \! 
\infty)$, in the double-scaling limit $\mathscr{N},n \! \to \! \infty$ such that 
$z_{o} \! = \! 1 \! + \! o(1)$, there exists $\hat{\lambda}_{1} \! > \! 0$ and 
$\mathcal{O}(1)$ (resp., $\tilde{\lambda}_{1} \! > \! 0$ and $\mathcal{O}(1))$ 
such that the associated monic MPC ORF zeros (counting multiplicities) 
$\lbrace \hat{\mathfrak{z}}^{n}_{k}(j) \rbrace_{j=1}^{(n-1)K+k} \subseteq 
[-\hat{\lambda}_{1},\hat{\lambda}_{1}] \setminus \cup_{q=1}^{\mathfrak{s}-1} 
\mathscr{O}_{\frac{1}{\hat{\lambda}_{1}}}(\alpha_{p_{q}})$ (resp., $\lbrace 
\tilde{\mathfrak{z}}^{n}_{k}(j) \rbrace_{j=1}^{(n-1)K+k} \subseteq 
[-\tilde{\lambda}_{1},\tilde{\lambda}_{1}] \setminus \cup_{\underset{q \neq 
\mathfrak{s}-1}{q=1}}^{\mathfrak{s}} \mathscr{O}_{\frac{1}{\tilde{\lambda}_{1}}}
(\alpha_{p_{q}}))$.
\begin{eeeee} \label{commzero}
\textsl{Note that, for $n \! \in \! \mathbb{N}$ and $k \! \in \! \lbrace 1,2,\dotsc,
K \rbrace$ such that $\alpha_{p_{\mathfrak{s}}} \! := \! \alpha_{k} \! = \! \infty$ 
(resp., $\alpha_{p_{\mathfrak{s}}} \! := \! \alpha_{k} \! \neq \! \infty)$, in the 
double-scaling limit $\mathscr{N},n \! \to \! \infty$ such that $z_{o} \! = \! 1 
\! + \! o(1)$, $[-\hat{\lambda}_{1},\hat{\lambda}_{1}] \setminus \cup_{q=1}^{
\mathfrak{s}-1} \mathscr{O}_{\frac{1}{\hat{\lambda}_{1}}}(\alpha_{p_{q}}) 
\supseteq J_{\infty}$ (resp., $[-\tilde{\lambda}_{1},\tilde{\lambda}_{1}] 
\setminus \cup_{\underset{q \neq \mathfrak{s}-1}{q=1}}^{\mathfrak{s}} 
\mathscr{O}_{\frac{1}{\hat{\lambda}_{1}}}(\alpha_{p_{q}}) \supseteq J_{f})$.}
\end{eeeee}
\begin{ccccc} \label{lemrootz}
Let the external field $\widetilde{V} \colon \overline{\mathbb{R}} 
\setminus \lbrace \alpha_{1},\alpha_{2},\dotsc,\alpha_{K} 
\rbrace \! \to \! \mathbb{R}$ satisfy 
conditions~\eqref{eq20}--\eqref{eq22} and be regular. For $n \! \in \! 
\mathbb{N}$ and $k \! \in \! \lbrace 1,2,\dotsc,K \rbrace$ such that 
$\alpha_{p_{\mathfrak{s}}} \! := \! \alpha_{k} \! = \! \infty$ (resp., 
$\alpha_{p_{\mathfrak{s}}} \! := \! \alpha_{k} \! \neq \! \infty)$, let the 
associated equilibrium measure, $\mu_{\widetilde{V}}^{\infty}$ (resp., 
$\mu_{\widetilde{V}}^{f})$, and its support, $J_{\infty}$ (resp., $J_{f})$, 
be as described in item~$\pmb{(1)}$ (resp., item~$\pmb{(2)})$ of 
Lemma~\ref{lem3.7}. For $n \! \in \! \mathbb{N}$ and $k \! \in \! 
\lbrace 1,2,\dotsc,K \rbrace$ such that $\alpha_{p_{\mathfrak{s}}} 
\! := \!\alpha_{k} \! = \! \infty$ (resp., $\alpha_{p_{\mathfrak{s}}} 
\! := \!\alpha_{k} \! \neq \! \infty)$, denote the associated monic 
{\rm MPC ORF} zeros (counting multiplicities) by $\lbrace 
\hat{\mathfrak{z}}^{n}_{k}(j) \rbrace_{j=1}^{(n-1)K+k}$ (resp., $\lbrace 
\tilde{\mathfrak{z}}^{n}_{k}(j) \rbrace_{j=1}^{(n-1)K+k})$. Then, for 
$n \! \in \! \mathbb{N}$ and $k \! \in \! \lbrace 1,2,\dotsc,K \rbrace$ 
such that $\alpha_{p_{\mathfrak{s}}} \! := \! \alpha_{k} \! = \! \infty$ 
(resp., $\alpha_{p_{\mathfrak{s}}} \! := \! \alpha_{k} \! \neq \! \infty)$, 
in the double-scaling limit $\mathscr{N},n \! \to \! \infty$ such that 
$z_{o} \! = \! 1 \! + \! o(1)$, there exists $\hat{\lambda}_{1} \! = \! 
\hat{\lambda}_{1}(n,k,z_{o}) \! > \! 0$ and $\mathcal{O}(1)$ (resp., 
$\tilde{\lambda}_{1} \! = \! \tilde{\lambda}_{1}(n,k,z_{o}) \! > \! 0$ and 
$\mathcal{O}(1))$ such that $\lbrace \hat{\mathfrak{z}}^{n}_{k}(j) 
\rbrace_{j=1}^{(n-1)K+k} \subseteq [-\hat{\lambda}_{1},\hat{\lambda}_{1}] 
\setminus \cup_{q=1}^{\mathfrak{s}-1} \mathscr{O}_{\frac{1}{\hat{\lambda}_{1}}}
(\alpha_{p_{q}})$ (resp., $\lbrace \tilde{\mathfrak{z}}^{n}_{k}(j) 
\rbrace_{j=1}^{(n-1)K+k} \subseteq [-\tilde{\lambda}_{1},\tilde{\lambda}_{1}] 
\setminus \cup_{\underset{q \neq \mathfrak{s}-1}{q=1}}^{\mathfrak{s}} 
\mathscr{O}_{\frac{1}{\tilde{\lambda}_{1}}}(\alpha_{p_{q}}))$.
\end{ccccc}

\emph{Proof}. The proof of this Lemma~\ref{lemrootz} consists of two 
cases: (i) $n \! \in \! \mathbb{N}$ and $k \! \in \! \lbrace 1,2,\dotsc,K 
\rbrace$ such that $\alpha_{p_{\mathfrak{s}}} \! := \! \alpha_{k} \! = \! 
\infty$; and (ii) $n \! \in \! \mathbb{N}$ and $k \! \in \! \lbrace 1,2,
\dotsc,K \rbrace$ such that $\alpha_{p_{\mathfrak{s}}} \! := \! \alpha_{k} 
\! \neq \! \infty$. Notwithstanding the fact that the scheme of the 
proof is, \emph{mutatis mutandis}, similar for both cases, case~(ii) is 
the more technically challenging of the two; therefore, without loss of 
generality, its particulars are presented in detail (see \pmb{(1)} below), 
whilst case~(i) is proved analogously (see \pmb{(2)} below).

\pmb{(1)} For $n \! \in \! \mathbb{N}$ and $k \! \in \! \lbrace 1,2,\dotsc,
K \rbrace$ such that $\alpha_{p_{\mathfrak{s}}} \! := \! \alpha_{k} \! \neq 
\! \infty$, set, for any arbitrarily small $\tilde{\delta} \! > \! 
0$ {}\footnote{Strictly speaking, $\tilde{\delta} \! = \! \tilde{\delta}
(n,k,z_{o})$; in order to eschew a flood of superfluous notation, unless 
where absolutely necessary, explicit $n$-, $k$-, and $z_{o}$-dependencies 
are suppressed. The parametric dependence of $\tilde{\delta}$ 
on the corresponding pole set $\lbrace \alpha_{p_{1}},\dotsc,
\alpha_{p_{\mathfrak{s}-2}},\alpha_{p_{\mathfrak{s}}} \rbrace$ 
is not considered. This is a general comment which applies, 
\emph{mutatis mutandis}, throughout the proof of this 
Lemma~\ref{lemrootz} (as well as in Corollary~\ref{corol3.1} and 
Lemma~\ref{lemetatomu} below).} satisfying $[\tilde{b}_{j-1} \! - 
\! \tilde{\delta},\tilde{b}_{j-1}] \cap \lbrace \alpha_{p_{1}},\dotsc,
\alpha_{p_{\mathfrak{s}-2}},\alpha_{p_{\mathfrak{s}}} \rbrace \! = 
\! \varnothing \! = \! [\tilde{a}_{j},\tilde{a}_{j} \! + \! \tilde{\delta}] 
\cap \lbrace \alpha_{p_{1}},\dotsc,\alpha_{p_{\mathfrak{s}-2}},
\alpha_{p_{\mathfrak{s}}} \rbrace$, $j \! = \! 1,2,\dotsc,N \! + \! 1$, 
and $2 \tilde{\delta} \! \ll \! \min_{i=1,2,\dotsc,N+1} \lbrace \lvert 
\tilde{a}_{i} \! - \! \tilde{b}_{i-1} \rvert \rbrace$, $\tilde{\phi}_{
\tilde{\delta}}(x) \! := \!  (2 \tilde{\delta})^{-1} \int_{x-
\tilde{\delta}}^{x+\tilde{\delta}} \md \mu_{\widetilde{V}}^{f}
(\xi)$ $(= \! (2 \tilde{\delta})^{-1} \int_{x-\tilde{\delta}}^{x+
\tilde{\delta}} \psi_{\widetilde{V}}^{f}(\xi) \, \md \xi)$. Via 
the properties of the corresponding equilibrium measure, $\mu_{
\widetilde{V}}^{f}$, and its support, $\supp (\mu_{\widetilde{V}}^{f}) \! 
:= \! J_{f} \! = \! \cup_{j=1}^{N+1}[\tilde{b}_{j-1},\tilde{a}_{j}]$, stated 
in item~\pmb{(2)} of Lemma~\ref{lem3.7}, and the fact that $J_{f} 
\cap \lbrace \alpha_{p_{1}},\dotsc,\alpha_{p_{\mathfrak{s}-2}},
\alpha_{p_{\mathfrak{s}}} \rbrace \! = \! \varnothing$, it follows 
{}from the above definition of $\tilde{\phi}_{\tilde{\delta}}$, the proof 
of Lemma~\ref{lem3.1}, and an application of Fubini's Theorem, that: 
(i) $\supp (\tilde{\phi}_{\tilde{\delta}})$ is a proper compact subset 
of $\mathbb{R}$, with $\supp (\tilde{\phi}_{\tilde{\delta}}) \cap 
\lbrace \alpha_{p_{1}},\dotsc,\alpha_{p_{\mathfrak{s}-2}},
\alpha_{p_{\mathfrak{s}}} \rbrace \! = \! \varnothing$ (of course, 
$\supp (\tilde{\phi}_{\tilde{\delta}}) \cap \lbrace \alpha_{p_{
\mathfrak{s}-1}} \! = \! \infty \rbrace \! = \! \varnothing)$; (ii) 
$\tilde{\phi}_{\tilde{\delta}}$ is non-negative, bounded, and continuous 
on $\supp (\tilde{\phi}_{\tilde{\delta}})$; and (iii) $\int_{\mathbb{R}} 
\tilde{\phi}_{\tilde{\delta}}(\xi) \, \md \xi$ $(= \! \int_{\supp 
(\tilde{\phi}_{\tilde{\delta}})} \tilde{\phi}_{\tilde{\delta}}(\xi) \, \md \xi)$ 
$= \! 1$. These properties imply that (see, for example, \cite{a55}) 
$\tilde{\phi}_{\tilde{\delta}}(x) \, \md x \! \overset{\ast}{\to} \! \md 
\mu_{\widetilde{V}}^{f}(x)$ as $\tilde{\delta} \! \downarrow \! 0$. {}From 
the definition of the corresponding weighted logarithmic energy functional 
given by Equation~\eqref{eqlm3.1b}, one shows, via an application of the 
extended triangle inequality, that, for $n \! \in \! \mathbb{N}$ and $k \! 
\in \! \lbrace 1,2,\dotsc,K \rbrace$ such that $\alpha_{p_{\mathfrak{s}}} 
\! := \! \alpha_{k} \! \neq \! \infty$, with $\mathcal{N} \! = \! (n \! - \! 
1)K \! + \! k$,
\begin{align} \label{eqlmrtz1} 
\left\lvert \mathrm{I}_{\widetilde{V}}^{f}[\mu_{\widetilde{V}}^{f}] 
\! - \! \mathrm{I}_{\widetilde{V}}^{f}[\tilde{\phi}_{\tilde{\delta}}] 
\right\rvert \leqslant& \, \dfrac{\mathcal{N}}{n} \left\lvert 
\iint_{\mathbb{R}^{2}} \ln \lvert \xi \! - \! \tau \rvert \, 
\tilde{\phi}_{\tilde{\delta}}(\xi) \, \md \xi \, \tilde{\phi}_{\tilde{\delta}}
(\tau) \, \md \tau \! - \!\iint_{\mathbb{R}^{2}} 
\ln \lvert \xi \! - \! \tau \rvert \, \md \mu_{\widetilde{V}}^{f}(\xi) 
\, \md \mu_{\widetilde{V}}^{f}(\tau) \right\rvert \nonumber \\
+& \, 2 \left(\dfrac{\varkappa_{nk} \! - \! 1}{n} \right) \left\lvert 
\int_{\mathbb{R}} \ln \lvert \xi \! - \! \alpha_{k} \rvert \, \md 
\mu_{\widetilde{V}}^{f}(\xi) \! - \! \int_{\mathbb{R}} \ln \lvert 
\xi \! - \! \alpha_{k} \rvert \, \tilde{\phi}_{\tilde{\delta}}(\xi) \, 
\md \xi \right\rvert \nonumber \\
+& \, 2 \sum_{q=1}^{\mathfrak{s}-2} \dfrac{\varkappa_{nk 
\tilde{k}_{q}}}{n} \left\lvert \int_{\mathbb{R}} \ln \lvert \xi \! - \! 
\alpha_{p_{q}} \rvert \, \md \mu_{\widetilde{V}}^{f}(\xi) \! - \! 
\int_{\mathbb{R}} \ln \lvert \xi \! - \! \alpha_{p_{q}} \rvert \, 
\tilde{\phi}_{\tilde{\delta}}(\xi) \, \md \xi \right\rvert \nonumber \\
+& \, \left\lvert \int_{\mathbb{R}} \widetilde{V}(\xi) \, \md 
\mu_{\widetilde{V}}^{f}(\xi) \! - \! \int_{\mathbb{R}} \widetilde{V}
(\xi) \, \tilde{\phi}_{\tilde{\delta}}(\xi) \, \md \xi \right\rvert.
\end{align}
Via the properties of $\mu_{\widetilde{V}}^{f}$ $(\in \! \mathscr{M}_{1}
(\mathbb{R}))$ and $J_{f}$ stated in Lemma~\ref{lem3.1} and 
item~\pmb{(2)} of Lemma~\ref{lem3.7}, the fact that, for $q \! = \! 1,
\dotsc,\mathfrak{s} \! - \! 2,\mathfrak{s}$, $\inf \lbrace \mathstrut \lvert 
x \! - \! \alpha_{p_{q}} \rvert; \, x \! \in \! J_{f} \rbrace \! > \! 0$, and, for 
$0 \! \leqslant \! u \! \leqslant \! \tilde{\delta} \! < \! 1$ and $x \! \in 
\! J_{f}$, $0 \! \leqslant \! \ln (1 \! + \! u \lvert x \! - \! \alpha_{p_{q}} 
\rvert^{-1}) \! \leqslant \! \ln (2) \! + \! \lvert \ln (\inf \lbrace \mathstrut 
\lvert \xi \! - \! \alpha_{p_{q}} \rvert; \, \xi \! \in \! J_{f} \rbrace)^{-1} 
\rvert$ $(< \! +\infty)$, the expansion $\ln (1 \! - \! \lozenge) \! 
=_{\lvert \lozenge \rvert \to 0} -\sum_{m=1}^{\infty} \tfrac{\lozenge^{m}}{m}$, 
and the fact that $\widetilde{V} \colon \overline{\mathbb{R}} \setminus 
\lbrace \alpha_{1},\alpha_{2},\dotsc,\alpha_{K} \rbrace \! \to \! \mathbb{R}$ 
satisfies conditions~\eqref{eq20}--\eqref{eq22} and is regular, one shows, 
via an integration-by-parts argument, an application of Fubini's Theorem, 
and $\tilde{\phi}_{\tilde{\delta}}(x) \, \md x \! \overset{\ast}{\to} \! 
\md \mu_{\widetilde{V}}^{f}(x)$ as $\delta \! \downarrow \! 0$, that, 
for $n \! \in \! \mathbb{N}$ and $k \! \in \! \lbrace 1,2,\dotsc,K 
\rbrace$ such that $\alpha_{p_{\mathfrak{s}}} \! := \! \alpha_{k} \! 
\neq \! \infty$,
\begin{equation*}
\left\lvert \int_{\mathbb{R}} \widetilde{V}(\xi) \, \md 
\mu_{\widetilde{V}}^{f}(\xi) \! - \! \int_{\mathbb{R}} \widetilde{V}
(\xi) \, \tilde{\phi}_{\tilde{\delta}}(\xi) \, \md \xi \right\rvert 
\leqslant \int_{J_{f}} \left\lvert \dfrac{1}{2 \tilde{\delta}} \int_{\xi 
-\tilde{\delta}}^{\xi +\tilde{\delta}}(\widetilde{V}(\xi) \! - \! 
\widetilde{V}(\tau)) \, \md \tau \right\rvert \md \mu_{\widetilde{V}}^{f}
(\xi) \underset{\tilde{\delta} \downarrow 0}{=} \mathcal{O} \left(
\mathfrak{c}_{1}(n,k,z_{o}) \tilde{\delta}^{3/2} \right),
\end{equation*}
\begin{align*}
\left(\dfrac{\varkappa_{nk} \! - \! 1}{n} \right) \left\lvert 
\int_{\mathbb{R}} \ln \lvert \xi \! - \! \alpha_{k} \rvert \, \md 
\mu_{\widetilde{V}}^{f}(\xi) \! - \! \int_{\mathbb{R}} \ln \lvert 
\xi \! - \! \alpha_{k} \rvert \, \tilde{\phi}_{\tilde{\delta}}(\xi) \, \md 
\xi \right\rvert \leqslant& \, \left(\dfrac{\varkappa_{nk} \! - \! 1}{n} 
\right) \int_{J_{f}} \left\lvert \dfrac{1}{2 \tilde{\delta}} \int_{-
\tilde{\delta}}^{\tilde{\delta}} \ln \left(1 \! - \! \dfrac{\tau}{\lvert 
\xi \! - \! \alpha_{k} \rvert} \right) \md \tau \right\rvert \md 
\mu_{\widetilde{V}}^{f}(\xi) \\
\underset{\tilde{\delta} \downarrow 0}{=}& \, \mathcal{O} \left(
\mathfrak{c}_{2}(n,k,z_{o}) \tilde{\delta}^{2} \right), \\
\sum_{q=1}^{\mathfrak{s}-2} \dfrac{\varkappa_{nk \tilde{k}_{q}}}{n} 
\left\lvert \int_{\mathbb{R}} \ln \lvert \xi \! - \! \alpha_{p_{q}} 
\rvert \, \md \mu_{\widetilde{V}}^{f}(\xi) \! - \! \int_{\mathbb{R}} \ln 
\lvert \xi \! - \! \alpha_{p_{q}} \rvert \, \tilde{\phi}_{\tilde{\delta}}(\xi) 
\, \md \xi \right\rvert \leqslant& \, \sum_{q=1}^{\mathfrak{s}-2} 
\dfrac{\varkappa_{nk \tilde{k}_{q}}}{n} \int_{J_{f}} \left\lvert \dfrac{1}{2 
\tilde{\delta}} \int_{-\tilde{\delta}}^{\tilde{\delta}} \ln \left(1 \! - \! 
\dfrac{\tau}{\lvert \xi \! - \! \alpha_{p_{q}} \rvert} \right) \md \tau 
\right\rvert \md \mu_{\widetilde{V}}^{f}(\xi) \\
\underset{\tilde{\delta} \downarrow 0}{=}& \, \mathcal{O} \left(
\mathfrak{c}_{3}(n,k,z_{o}) \tilde{\delta}^{2} \right),
\end{align*}
where $\mathfrak{c}_{m}(n,k,z_{o}) \! =_{\underset{z_{o}=1+o(1)}{
\mathscr{N},n \to \infty}} \! \mathcal{O}(1)$, $m \! = \! 1,2,3$, and, 
via the monotonicity of $\ln (\pmb{\cdot})$ and a change-of-variables 
argument,
\begin{align*}
&\dfrac{\mathcal{N}}{n} \left\lvert \iint_{\mathbb{R}^{2}} \ln 
\lvert \xi \! - \! \tau \rvert \, \tilde{\phi}_{\tilde{\delta}}(\xi) 
\, \md \xi \, \tilde{\phi}_{\tilde{\delta}}(\tau) 
\, \md \tau \! - \! \iint_{\mathbb{R}^{2}} 
\ln \lvert \xi \! - \! \tau \rvert \, \md \mu_{\widetilde{V}}^{f}
(\xi) \, \md \mu_{\widetilde{V}}^{f}(\tau) \right\rvert \\
\leqslant& \, \dfrac{\mathcal{N}}{n} \iint_{J_{f}^{2}} \underbrace{\left\lvert 
\dfrac{\lvert \xi \! - \! \tau \rvert}{2 \tilde{\delta}} \int_{-2 \tilde{\delta} 
\lvert \xi -\tau \rvert^{-1}}^{2 \tilde{\delta} \lvert \xi - \tau \rvert^{-1}} 
\ln (1 \! + \! u) \, \md u \right\rvert}_{\leqslant \mathfrak{c}_{4}(n,k,z_{o}) 
\ln (1+ \tilde{\delta} \lvert \xi - \tau \rvert^{-1})} \md \mu_{\widetilde{
V}}^{f}(\xi) \, \md \mu_{\widetilde{V}}^{f}(\tau) \\
\leqslant& \, \mathfrak{c}_{4}(n,k,z_{o}) \dfrac{\mathcal{N}}{n} 
\iint_{J_{f}^{2}} (\ln (\lvert \xi \! - \! \tau \rvert \! + \! \tilde{\delta}) 
\! - \! \ln \lvert \xi \! - \! \tau \rvert) \, \md \mu_{\widetilde{V}}^{f}
(\xi) \, \md \mu_{\widetilde{V}}^{f}(\tau),
\end{align*} 
where $\mathfrak{c}_{4}(n,k,z_{o}) \! =_{\underset{z_{o}=1+o(1)}{
\mathscr{N},n \to \infty}} \! \mathcal{O}(1)$: arguing, now, as in the 
proof of Lemma~\ref{lem3.1}, one shows, via an application of the 
Dominated Convergence Theorem, that $\iint_{J_{f}^{2}} (\ln (\lvert 
\xi \! - \! \tau \rvert \! + \! \tilde{\delta}) \! - \! \ln \lvert \xi 
\! - \! \tau \rvert) \, \md \mu_{\widetilde{V}}^{f}(\xi) \, \md 
\mu_{\widetilde{V}}^{f}(\tau) \! \to \! 0$ as $\tilde{\delta} \! \downarrow 
\! 0$; hence, for $n \! \in \! \mathbb{N}$ and $k \! \in \! \lbrace 1,2,
\dotsc,K \rbrace$ such that $\alpha_{p_{\mathfrak{s}}} \! := \! \alpha_{k} \! 
\neq \! \infty$, $\mathrm{I}_{\widetilde{V}}^{f}[\tilde{\phi}_{\tilde{\delta}}] \! 
\to \! \mathrm{I}_{\widetilde{V}}^{f}[\mu_{\widetilde{V}}^{f}]$ as $\tilde{\delta} 
\! \downarrow \! 0$, that is, for $n \! \in \! \mathbb{N}$ and $k \! \in \! 
\lbrace 1,2,\dotsc,K \rbrace$ such that $\alpha_{p_{\mathfrak{s}}} \! := \! 
\alpha_{k} \! \neq \! \infty$, given $\tilde{\varepsilon} \! = \! \tilde{\varepsilon}
(n,k,z_{o}) \! > \! 0$ and sufficiently small, there exists sufficiently 
small $\tilde{\delta}(\tilde{\varepsilon}) \! := \! \tilde{\delta} \! 
> \! 0$ (as described above) such that one can choose 
$\tilde{\phi}_{\tilde{\varepsilon}}(\pmb{\cdot}) \! := \! 
\tilde{\phi}_{\tilde{\delta}(\tilde{\varepsilon})}(\pmb{\cdot})$ 
with finite entropy (that is, $\tilde{\phi}_{\tilde{\varepsilon}} 
\ln \tilde{\phi}_{\tilde{\varepsilon}}$ is integrable) so that 
$\mathrm{I}_{\widetilde{V}}^{f}[\tilde{\phi}_{\tilde{\varepsilon}}] 
\! \leqslant \! E_{\widetilde{V}}^{f} \! + \! \tilde{\varepsilon}/2$ 
(cf. the proof of Lemma~\ref{lem3.1}). For $n \! \in \! 
\mathbb{N}$ and $k \! \in \! \lbrace 1,2,\dotsc,K \rbrace$ 
such that $\alpha_{p_{\mathfrak{s}}} \! := \! \alpha_{k} \! 
\neq \! \infty$, denote the corresponding multi-dimensional 
probability measure on $\mathbb{R}^{\mathcal{N}}$ by 
$\tilde{\mathcal{P}}^{n}_{k}(x_{1},x_{2},\dotsc,x_{\mathcal{N}}) 
\, \md x_{1} \, \md x_{2} \, \dotsb \, \md x_{\mathcal{N}}$, 
with density {}\footnote{Note: $\idotsint\nolimits_{
\mathbb{R}^{\raise-0.7ex\hbox{$\scriptscriptstyle \mathcal{N}$}}} 
\tilde{\mathcal{P}}^{n}_{k}(\xi_{1},\xi_{2},\dotsc,\xi_{\mathcal{N}}) \, 
\md \xi_{1} \, \md \xi_{2} \, \dotsb \, \md \xi_{\mathcal{N}} 
\! = \! 1$.}
\begin{equation} \label{eqlmrtzp} 
\tilde{\mathcal{P}}^{n}_{k}(x_{1},x_{2},\dotsc,x_{\mathcal{N}}) \! := \! 
\dfrac{1}{\tilde{\mathscr{Z}}^{n}_{k}} \me^{-n \sum_{m_{1}=1}^{
\mathcal{N}} \widetilde{V}(x_{m_{1}})} \prod_{\substack{i,j=1\\j<i}}^{
\mathcal{N}}(x_{j} \! - \! x_{i})^{2} \left(\prod_{m=1}^{\mathcal{N}} 
\prod_{q=1}^{\mathfrak{s}-2}(x_{m} \! - \! \alpha_{p_{q}})^{\varkappa_{nk 
\tilde{k}_{q}}}(x_{m} \! - \! \alpha_{k})^{\varkappa_{nk}-1} \right)^{-2},
\end{equation}
where
\begin{equation} \label{eqlmrtzz} 
\tilde{\mathscr{Z}}^{n}_{k} \! = \! \idotsint\limits_{\mathbb{R}^{\mathcal{N}}} 
\me^{-n \sum_{m_{1}=1}^{\mathcal{N}} \widetilde{V}(\tau_{m_{1}})} 
\prod_{\substack{i,j=1\\j<i}}^{\mathcal{N}}(\tau_{j} \! - \! \tau_{i})^{2} 
\left(\prod_{m=1}^{\mathcal{N}} \prod_{q=1}^{\mathfrak{s}-2}
(\tau_{m} \! - \! \alpha_{p_{q}})^{\varkappa_{nk \tilde{k}_{q}}}
(\tau_{m} \! - \! \alpha_{k})^{\varkappa_{nk}-1} \right)^{-2} \md 
\tau_{1} \, \md \tau_{2} \, \dotsb \, \md \tau_{\mathcal{N}}.
\end{equation}
Recalling {}from the proof of Lemma~\ref{lem3.6} that, for a symmetric 
function $\tilde{f}(x,y)$, say, $\sum_{\underset{i \neq j}{i,j=1}}^{\mathcal{N}} 
\tilde{f}(x_{i},x_{j}) \! = \! 2 \sum_{\underset{i<j}{i,j=1}}^{\mathcal{N}} \tilde{f}
(x_{i},x_{j})$, and, for a function $\tilde{h}(x)$, say, $(\mathcal{N} \! - \! 1) 
\sum_{m=1}^{\mathcal{N}} \tilde{h}(x_{m}) \! = \! \sum_{\underset{j<m}{j,
m=1}}^{\mathcal{N}} \tilde{h}(x_{j}) \! + \! \sum_{\underset{j<m}{j,m=
1}}^{\mathcal{N}} \tilde{h}(x_{m})$, it follows via the decomposition 
(cf. Equation~\eqref{fincount}) $\sum_{q=1}^{\mathfrak{s}-2} 
\varkappa_{nk \tilde{k}_{q}} \! + \! \varkappa^{\infty}_{nk 
\tilde{k}_{\mathfrak{s}-1}} \! + \! \varkappa_{nk} \! = \! (n \! - \! 1)K \! 
+ \! k \! =: \! \mathcal{N}$, and Equation~\eqref{eql3.5d}, that, for $n 
\! \in \! \mathbb{N}$ and $k \! \in \! \lbrace 1,2,\dotsc,K \rbrace$ such 
that $\alpha_{p_{\mathfrak{s}}} \! := \! \alpha_{k} \! \neq \! \infty$,
\begin{equation*}
\tilde{\mathscr{Z}}^{n}_{k} \underset{\underset{z_{o}=1+o(1)}{
\mathscr{N},n \to \infty}}{=} \idotsint\limits_{\mathbb{R}^{\mathcal{N}}} 
\me^{-\frac{n}{\mathcal{N}}\left(\sum_{j=1}^{\mathcal{N}} \widetilde{V}
(\xi_{j})+ \mathscr{K}^{\widetilde{V},f}_{\mathcal{N}}(\xi_{1},\xi_{2},
\dotsc,\xi_{\mathcal{N}})+\mathcal{O}(\frac{1}{\mathcal{N}}) 
\sum_{\underset{i \neq j}{i,j=1}}^{\mathcal{N}} \tilde{\mathrm{F}}
(\xi_{i},\xi_{j}) \right)} \, \md \xi_{1} \, \md \xi_{2} \, \dotsb \, \md 
\xi_{\mathcal{N}},
\end{equation*}
where
\begin{equation} \label{eqlmrtzf} 
\tilde{\mathrm{F}}(x,y) \! := \! \left(\dfrac{\varkappa_{nk} \! 
- \! 1}{n} \right) \ln (\lvert x \! - \! \alpha_{k} \rvert \lvert y 
\! - \! \alpha_{k} \rvert) \! + \! \sum_{q=1}^{\mathfrak{s}-2} 
\dfrac{\varkappa_{nk \tilde{k}_{q}}}{n} \ln (\lvert x \! - \! 
\alpha_{p_{q}} \rvert \lvert y \! - \! \alpha_{p_{q}} \rvert).
\end{equation}
For $n \! \in \! \mathbb{N}$ and $k \! \in \! \lbrace 1,2,\dotsc,
K \rbrace$ such that $\alpha_{p_{\mathfrak{s}}} \! := \! 
\alpha_{k} \! \neq \! \infty$, set 
$\tilde{\mathbb{E}}^{\raise-0.0ex\hbox{$\scriptscriptstyle \mathcal{N}$}} 
\! := \! \lbrace \mathstrut (x_{1},x_{2},\dotsc,x_{\mathcal{N}}) \! 
\in \! \mathbb{R}^{\mathcal{N}}; \, \prod_{i=1}^{\mathcal{N}} 
\tilde{\phi}_{\tilde{\varepsilon}}(x_{i}) \! > \! 0 \rbrace$ $(\subset 
\mathbb{R}^{\mathcal{N}})$; hence, via this definition, and an 
application of the multi-dimensional version of Jensen's Inequality 
(see, for example, \cite{jbgtt}), it follows, via the monotonicty of 
$\ln (\pmb{\cdot})$, that
\begin{align*}
\tilde{\mathscr{Z}}^{n}_{k} \underset{\underset{z_{o}=1
+o(1)}{\mathscr{N},n \to \infty}}{\geqslant}& \, \idotsint\limits_{
\tilde{\mathbb{E}}^{\raise-0.0ex\hbox{$\scriptscriptstyle \mathcal{N}$}}} 
\me^{-\frac{n}{\mathcal{N}}\left(\mathscr{K}^{\widetilde{V},f}_{\mathcal{
N}}(\xi_{1},\xi_{2},\dotsc,\xi_{\mathcal{N}})+ \sum_{j=1}^{\mathcal{N}} 
\widetilde{V}(\xi_{j})+ \mathcal{O}(\frac{1}{\mathcal{N}}) 
\sum_{\underset{i \neq j}{i,j=1}}^{\mathcal{N}} \tilde{\mathrm{F}}
(\xi_{i},\xi_{j})+\frac{\mathcal{N}}{n} \sum_{i=1}^{\mathcal{N}} 
\ln (\tilde{\phi}_{\tilde{\varepsilon}}(\xi_{i})) \right)} \prod_{m
=1}^{\mathcal{N}} \tilde{\phi}_{\tilde{\varepsilon}}(\xi_{m}) \, 
\md \xi_{m} \\
\underset{\underset{z_{o}=1+o(1)}{\mathscr{N},n \to \infty}}{\geqslant}& 
\, \exp \left(-\dfrac{n}{\mathcal{N}} \idotsint\limits_{
\tilde{\mathbb{E}}^{\raise-0.0ex\hbox{$\scriptscriptstyle \mathcal{N}$}}} 
\left(\mathscr{K}^{\widetilde{V},f}_{\mathcal{N}}(\xi_{1},\xi_{2},\dotsc,
\xi_{\mathcal{N}}) \! + \! \sum_{j=1}^{\mathcal{N}} \widetilde{V}(\xi_{j}) \! 
+ \! \mathcal{O} \left(\dfrac{1}{\mathcal{N}} \right) \sum_{\substack{i,j
=1\\i \neq j}}^{\mathcal{N}} \tilde{\mathrm{F}}(\xi_{i},\xi_{j}) 
\! + \! \dfrac{\mathcal{N}}{n} \sum_{i=1}^{\mathcal{N}} \ln 
(\tilde{\phi}_{\tilde{\varepsilon}}(\xi_{i})) \right) \right. \\
\times&\left. \, \prod_{m=1}^{\mathcal{N}} \tilde{\phi}_{\tilde{
\varepsilon}}(\xi_{m}) \, \md \xi_{m} \right) \quad \Rightarrow \\
\ln (\tilde{\mathscr{Z}}^{n}_{k}) \underset{\underset{z_{o}=1+o(1)}{\mathscr{N},
n \to \infty}}{\geqslant}& \, -\dfrac{n}{\mathcal{N}}\idotsint\limits_{
\tilde{\mathbb{E}}^{\raise-0.0ex\hbox{$\scriptscriptstyle \mathcal{N}$}}} 
\mathscr{K}^{\widetilde{V},f}_{\mathcal{N}}(\xi_{1},\xi_{2},\dotsc,
\xi_{\mathcal{N}}) \prod_{i=1}^{\mathcal{N}} \tilde{\phi}_{\tilde{\varepsilon}}
(\xi_{i}) \, \md \xi_{i} \!  - \! \dfrac{n}{\mathcal{N}} \idotsint\limits_{
\tilde{\mathbb{E}}^{\raise-0.0ex\hbox{$\scriptscriptstyle \mathcal{N}$}}} 
\sum_{j=1}^{\mathcal{N}} \widetilde{V}(\xi_{j}) \prod_{i=1}^{\mathcal{N}} 
\tilde{\phi}_{\tilde{\varepsilon}}(\xi_{i}) \, \md \xi_{i} \\
-& \, \idotsint\limits_{
\tilde{\mathbb{E}}^{\raise-0.0ex\hbox{$\scriptscriptstyle \mathcal{N}$}}} 
\sum_{j=1}^{\mathcal{N}} \ln (\tilde{\phi}_{\tilde{\varepsilon}}(\xi_{j})) 
\prod_{i=1}^{\mathcal{N}} \tilde{\phi}_{\tilde{\varepsilon}}(\xi_{i}) \, \md 
\xi_{i} \! - \! \dfrac{n}{\mathcal{N}} \mathcal{O} \left(\dfrac{1}{\mathcal{N}} 
\right) \idotsint\limits_{
\tilde{\mathbb{E}}^{\raise-0.0ex\hbox{$\scriptscriptstyle \mathcal{N}$}}} 
\sum_{\substack{i,j=1\\i \neq j}}^{\mathcal{N}} \tilde{\mathrm{F}}(\xi_{i},
\xi_{j}) \prod_{m=1}^{\mathcal{N}} \tilde{\phi}_{\tilde{\varepsilon}}
(\xi_{m}) \, \md \xi_{m}.
\end{align*}
Using the fact that, for $i \! = \! 1,2,\dotsc,\mathcal{N}$, 
$\idotsint\nolimits_{
\tilde{\mathbb{E}}^{\raise-0.0ex\hbox{$\scriptscriptstyle \mathcal{N}-1$}}} 
\prod_{\underset{j \neq i}{j=1}}^{\mathcal{N}} 
\tilde{\phi}_{\tilde{\varepsilon}}(\xi_{j}) \, \md \xi_{j} 
\! = \! 1$, a straightforward argument shows that
\begin{gather}
-\dfrac{n}{\mathcal{N}} \idotsint\limits_{
\tilde{\mathbb{E}}^{\raise-0.0ex\hbox{$\scriptscriptstyle \mathcal{N}$}}} 
\sum_{j=1}^{\mathcal{N}} \widetilde{V}(\xi_{j}) \prod_{i=1}^{\mathcal{N}} 
\tilde{\phi}_{\tilde{\varepsilon}}(\xi_{i}) \, \md \xi_{i} \! = \! -n 
\int_{\tilde{\mathbb{E}}} \widetilde{V}(\xi) \, \tilde{\phi}_{
\tilde{\varepsilon}}(\xi) \, \md \xi, \label{eqlmrtz2} \\
-\idotsint\limits_{
\tilde{\mathbb{E}}^{\raise-0.0ex\hbox{$\scriptscriptstyle \mathcal{N}$}}} 
\sum_{j=1}^{\mathcal{N}} \ln (\tilde{\phi}_{\tilde{\varepsilon}}(\xi_{j})) 
\prod_{i=1}^{\mathcal{N}} \tilde{\phi}_{\tilde{\varepsilon}}(\xi_{i}) 
\, \md \xi_{i} \! = \! -\mathcal{N} \int_{\tilde{\mathbb{E}}} \ln 
(\tilde{\phi}_{\tilde{\varepsilon}}(\xi)) \, \tilde{\phi}_{\tilde{\varepsilon}}
(\xi) \, \md \xi, 
\label{eqlmrtz3}
\end{gather}
where $\tilde{\mathbb{E}} \! := \! \operatorname{int}
(\supp (\tilde{\phi}_{\tilde{\varepsilon}}))$ (with $\supp 
(\tilde{\phi}_{\tilde{\varepsilon}}) \cap \lbrace \alpha_{p_{1}},
\dotsc,\alpha_{p_{\mathfrak{s}-2}},\alpha_{p_{\mathfrak{s}}} 
\rbrace \! = \! \varnothing)$. It follows via 
Equation~\eqref{eql3.5d} that
\begin{align*}
-\dfrac{n}{\mathcal{N}} \idotsint\limits_{\tilde{\mathbb{E}}^{\mathcal{N}}} 
\mathscr{K}^{\widetilde{V},f}_{\mathcal{N}}(\xi_{1},\xi_{2},\dotsc,
\xi_{\mathcal{N}}) \prod_{i=1}^{\mathcal{N}} \tilde{\phi}_{\tilde{\varepsilon}}
(\xi_{i}) \, \md \xi_{i} =& \, -\dfrac{n}{\mathcal{N}} \idotsint\limits_{
\tilde{\mathbb{E}}^{\raise-0.0ex\hbox{$\scriptscriptstyle \mathcal{N}$}}} 
\sum_{\substack{i,j=1\\i \neq j}}^{\mathcal{N}} \left(\dfrac{\varkappa_{nk 
\tilde{k}_{\mathfrak{s}-1}}^{\infty} \! + \! 1}{n} \right) \ln \lvert \xi_{i} \! - \! 
\xi_{j} \rvert^{-1} \prod_{m=1}^{\mathcal{N}} \tilde{\phi}_{\tilde{\varepsilon}}
(\xi_{m}) \, \md \xi_{m} \\
-& \, \dfrac{n}{\mathcal{N}} \idotsint\limits_{
\tilde{\mathbb{E}}^{\raise-0.0ex\hbox{$\scriptscriptstyle \mathcal{N}$}}} 
\sum_{\substack{i,j=1\\i \neq j}}^{\mathcal{N}} \left(\dfrac{\varkappa_{nk} 
\! - \! 1}{n} \right) \ln \left(\dfrac{\lvert \xi_{i} \! - \! \alpha_{k} \rvert \lvert 
\xi_{j} \! - \! \alpha_{k} \rvert}{\lvert \xi_{i} \! - \! \xi_{j} \rvert} \right) 
\prod_{m=1}^{\mathcal{N}} \tilde{\phi}_{\tilde{\varepsilon}}(\xi_{m}) 
\, \md \xi_{m} \\
-& \, \dfrac{n}{\mathcal{N}} \idotsint\limits_{
\tilde{\mathbb{E}}^{\raise-0.0ex\hbox{$\scriptscriptstyle \mathcal{N}$}}} 
\sum_{\substack{i,j=1\\i \neq j}}^{\mathcal{N}} \sum_{q=1}^{\mathfrak{s}-2} 
\dfrac{\varkappa_{nk \tilde{k}_{q}}}{n} \ln \left(\dfrac{\lvert \xi_{i} \! - \! 
\alpha_{p_{q}} \rvert \lvert \xi_{j} \! - \! \alpha_{p_{q}} \rvert}{\lvert 
\xi_{i} \! - \! \xi_{j} \rvert} \right) \prod_{m=1}^{\mathcal{N}} 
\tilde{\phi}_{\tilde{\varepsilon}}(\xi_{m}) \, \md \xi_{m} \\
-& \, \dfrac{n}{\mathcal{N}}(\mathcal{N} \! - \! 1) 
\underbrace{\idotsint\limits_{
\tilde{\mathbb{E}}^{\raise-0.0ex\hbox{$\scriptscriptstyle \mathcal{N}$}}} 
\sum_{j=1}^{\mathcal{N}} \widetilde{V}(\xi_{j}) \prod_{m=1}^{\mathcal{N}} 
\tilde{\phi}_{\tilde{\varepsilon}}(\xi_{m}) \, \md \xi_{m}}_{= \, \mathcal{N} 
\int_{\tilde{\mathbb{E}}} \widetilde{V}(\xi) \tilde{\phi}_{\tilde{\varepsilon}}
(\xi) \, \md \xi}.
\end{align*}
Using the fact that, for $i \! \neq \! j \! \in \! \lbrace 
1,2,\dotsc,\mathcal{N} \rbrace$, $\idotsint\nolimits_{
\tilde{\mathbb{E}}^{\raise-0.0ex\hbox{$\scriptscriptstyle \mathcal{N}-2$}}} 
\prod_{\underset{m \neq i,j}{m=1}}^{\mathcal{N}} \tilde{\phi}_{
\tilde{\varepsilon}}(\xi_{m}) \, \md \xi_{m} \! = \! 1$, an enumeration 
argument shows that
\begin{align*}
-\dfrac{n}{\mathcal{N}} \idotsint\limits_{
\tilde{\mathbb{E}}^{\raise-0.0ex\hbox{$\scriptscriptstyle \mathcal{N}$}}} 
\sum_{\substack{i,j=1\\i \neq j}}^{\mathcal{N}} \left(\dfrac{\varkappa_{nk 
\tilde{k}_{\mathfrak{s}-1}}^{\infty} \! + \! 1}{n} \right) \ln \lvert 
\xi_{i} \! - \! \xi_{j} \rvert^{-1} \prod_{m=1}^{\mathcal{N}} 
\tilde{\phi}_{\tilde{\varepsilon}}(\xi_{m}) \, \md \xi_{m} =& \, 
-\dfrac{n}{\mathcal{N}} \mathcal{N}(\mathcal{N} \! - \! 1) \left(
\dfrac{\varkappa_{nk \tilde{k}_{\mathfrak{s}-1}}^{\infty} \! + \! 1}{n} 
\right) \\
\times& \, \iint_{\tilde{\mathbb{E}}^{2}} \ln (\lvert \xi 
\! - \! \tau \rvert^{-1}) \, \tilde{\phi}_{\tilde{\varepsilon}}(\xi) \, 
\md \xi \, \tilde{\phi}_{\tilde{\varepsilon}}(\tau) \, \md \tau, \\
-\dfrac{n}{\mathcal{N}} \idotsint\limits_{
\tilde{\mathbb{E}}^{\raise-0.0ex\hbox{$\scriptscriptstyle \mathcal{N}$}}} 
\sum_{\substack{i,j=1\\i \neq j}}^{\mathcal{N}} \left(\dfrac{
\varkappa_{nk} \! - \! 1}{n} \right) \ln \left(\dfrac{\lvert \xi_{i} \! - \! 
\alpha_{k} \rvert \lvert \xi_{j} \! - \! \alpha_{k} \rvert}{\lvert \xi_{i} 
\! - \! \xi_{j} \rvert} \right) \prod_{m=1}^{\mathcal{N}} 
\tilde{\phi}_{\tilde{\varepsilon}}(\xi_{m}) \, \md \xi_{m} =& \, 
-\dfrac{n}{\mathcal{N}} \mathcal{N}(\mathcal{N} \! - \! 1) \left(
\dfrac{\varkappa_{nk} \! - \! 1}{n} \right) \\
\times& \, \iint_{\tilde{\mathbb{E}}^{2}} \ln \left(
\dfrac{\lvert \xi \! - \! \alpha_{k} \rvert \lvert \tau \! - \! 
\alpha_{k} \rvert}{\lvert \xi \! - \! \tau \rvert} \right) \tilde{
\phi}_{\tilde{\varepsilon}}(\xi) \, \md \xi \, \tilde{\phi}_{
\tilde{\varepsilon}}(\tau) \, \md \tau, \\
-\dfrac{n}{\mathcal{N}} \idotsint\limits_{
\tilde{\mathbb{E}}^{\raise-0.0ex\hbox{$\scriptscriptstyle \mathcal{N}$}}} 
\sum_{\substack{i,j=1\\i \neq j}}^{\mathcal{N}} 
\sum_{q=1}^{\mathfrak{s}-2} \dfrac{\varkappa_{nk \tilde{k}_{q}}}{n} 
\ln \left(\dfrac{\lvert \xi_{i} \! - \! \alpha_{p_{q}} \rvert \lvert 
\xi_{j} \! - \! \alpha_{p_{q}} \rvert}{\lvert \xi_{i} \! - \! \xi_{j} 
\rvert} \right) \prod_{m=1}^{\mathcal{N}} \tilde{\phi}_{
\tilde{\varepsilon}}(\xi_{m}) \, \md \xi_{m} =& \, -\dfrac{n}{
\mathcal{N}} \mathcal{N}(\mathcal{N} \! - \! 1) \sum_{q=1}^{
\mathfrak{s}-2} \dfrac{\varkappa_{nk \tilde{k}_{q}}}{n} \\
\times& \, \iint_{\tilde{\mathbb{E}}^{2}} \ln \left(
\dfrac{\lvert \xi \! - \! \alpha_{p_{q}} \rvert \lvert \tau \! - 
\! \alpha_{p_{q}} \rvert}{\lvert \xi \! - \! \tau \rvert} \right) 
\tilde{\phi}_{\tilde{\varepsilon}}(\xi) \, \md \xi \, 
\tilde{\phi}_{\tilde{\varepsilon}}(\tau) \, \md \tau,
\end{align*}
whence, via Equation~\eqref{eqKvinf1} and the proof of Lemma~\ref{lem3.5},
\begin{align} \label{eqlmrtz4}
-\dfrac{n}{\mathcal{N}} \idotsint\limits_{
\tilde{\mathbb{E}}^{\raise-0.0ex\hbox{$\scriptscriptstyle \mathcal{N}$}}} 
\mathscr{K}^{\widetilde{V},f}_{\mathcal{N}}(\xi_{1},\xi_{2},\dotsc,
\xi_{\mathcal{N}}) \prod_{i=1}^{\mathcal{N}} \tilde{\phi}_{\tilde{\varepsilon}}
(\xi_{i}) \, \md \xi_{i} =& \, -\dfrac{n}{\mathcal{N}} \mathcal{N}
(\mathcal{N} \! - \! 1) \iint_{\tilde{\mathbb{E}}^{2}} \ln \left(\lvert 
\xi \! - \! \tau \rvert^{-\frac{(\varkappa_{nk \tilde{k}_{\mathfrak{s}
-1}}^{\infty}+1)}{n}} \left(\dfrac{\lvert \xi \! - \! \alpha_{k} 
\rvert \lvert \tau \! - \! \alpha_{k} \rvert}{\lvert \xi \! - \! \tau 
\rvert} \right)^{\frac{\varkappa_{nk}-1}{n}} \right. \nonumber \\
\times&\left. \, \prod_{q=1}^{\mathfrak{s}-2} \left(\dfrac{\lvert 
\xi \! - \! \alpha_{p_{q}} \rvert \lvert \tau \! - \! \alpha_{p_{q}} 
\rvert}{\lvert \xi \! - \! \tau \rvert} \right)^{\frac{\varkappa_{nk 
\tilde{k}_{q}}}{n}} \me^{\widetilde{V}(\xi)/2} \me^{\widetilde{V}
(\tau)/2} \right) \tilde{\phi}_{\tilde{\varepsilon}}(\xi) \, \md 
\xi \, \tilde{\phi}_{\tilde{\varepsilon}}(\tau) \, \md \tau 
\nonumber \\
=& \, -\dfrac{n}{\mathcal{N}} \mathcal{N}(\mathcal{N} \! - \! 1) 
\mathrm{I}_{\widetilde{V}}^{f}[\tilde{\phi}_{\tilde{\varepsilon}}].
\end{align}
Using the fact that, for $i \! \neq \! j \! \in \! \lbrace 1,2,
\dotsc,\mathcal{N} \rbrace$, $\idotsint\nolimits_{\tilde{
\mathbb{E}}^{\raise-0.0ex\hbox{$\scriptscriptstyle \mathcal{N}-2$}}} 
\prod_{\underset{m \neq i,j}{m=1}}^{\mathcal{N}} \tilde{\phi}_{
\tilde{\varepsilon}}(\xi_{m}) \, \md \xi_{m} \! = \! 1$, one shows, 
via Equation~\eqref{eqlmrtzf} and an enumeration argument, that
\begin{align*}
&\dfrac{1}{\mathcal{N}(\mathcal{N} \! - \! 1)} \idotsint\limits_{
\tilde{\mathbb{E}}^{\raise-0.0ex\hbox{$\scriptscriptstyle \mathcal{N}$}}} 
\sum_{\substack{i,j=1\\i \neq j}}^{\mathcal{N}} \tilde{\mathrm{F}}
(\xi_{i},\xi_{j}) \prod_{m=1}^{\mathcal{N}} \tilde{\phi}_{
\tilde{\varepsilon}}(\xi_{m}) \, \md \xi_{m} \! = \! \iint_{
\tilde{\mathbb{E}}^{2}} \tilde{\mathrm{F}}(\xi,\tau) \, 
\tilde{\phi}_{\tilde{\varepsilon}}(\xi) \, \md \xi \, 
\tilde{\phi}_{\tilde{\varepsilon}}(\tau) \, \md \tau \\
=& \, 2 \left(\dfrac{\varkappa_{nk} \! - \! 1}{n} \right) 
\left(\int_{\tilde{\mathbb{E}}} \ln \lvert \xi \! - \! 
\alpha_{k} \rvert \, \tilde{\phi}_{\tilde{\varepsilon}}(\xi) 
\, \md \xi \! - \! \int_{\tilde{\mathbb{E}}} \ln \lvert 
\xi \! - \! \alpha_{k} \rvert \, \md \mu_{\widetilde{V}}^{f}
(\xi) \right) \! + \! 2 \left(\dfrac{\varkappa_{nk} \! - \! 1}{n} 
\right) \int_{\tilde{\mathbb{E}}} \ln \lvert \xi \! - 
\! \alpha_{k} \rvert \, \md \mu_{\widetilde{V}}^{f}(\xi) \\
+& \, 2 \sum_{q=1}^{\mathfrak{s}-2} \dfrac{\varkappa_{nk 
\tilde{k}_{q}}}{n} \left(\int_{\tilde{\mathbb{E}}} \ln 
\lvert \xi \! - \! \alpha_{p_{q}} \rvert \, \tilde{\phi}_{
\tilde{\varepsilon}}(\xi) \, \md \xi \! - \! \int_{
\tilde{\mathbb{E}}} \ln \lvert \xi \! - \! \alpha_{p_{q}} 
\rvert \, \md \mu_{\widetilde{V}}^{f}(\xi) 
\right) \! + \! 2 \sum_{q=1}^{\mathfrak{s}-2} \dfrac{
\varkappa_{nk \tilde{k}_{q}}}{n} \int_{\tilde{
\mathbb{E}}} \ln \lvert \xi \! - \! \alpha_{p_{q}} \rvert 
\, \md \mu_{\widetilde{V}}^{f}(\xi):
\end{align*}
arguing, \emph{mutatis mutandis}, as at the beginning of the proof 
(cf. the calculations leading to the convergence $\mathrm{I}_{
\widetilde{V}}^{f}[\tilde{\phi}_{\tilde{\delta}(\tilde{\varepsilon})}] 
\linebreak[4] \! = \! 
\mathrm{I}_{\widetilde{V}}^{f}[\tilde{\phi}_{\tilde{\varepsilon}}] \! 
\to \! \mathrm{I}_{\widetilde{V}}^{f}[\mu_{\widetilde{V}}^{f}] \! = \! 
E_{\widetilde{V}}^{f}$ as $\tilde{\varepsilon} \! \downarrow \! 0)$, 
one shows that
\begin{gather*}
\left(\dfrac{\varkappa_{nk} \! - \! 1}{n} \right) \left(
\int_{\tilde{\mathbb{E}}} \ln \lvert \xi \! - \! \alpha_{k} 
\rvert \, \tilde{\phi}_{\tilde{\varepsilon}}(\xi) \, \md \xi \! - 
\! \int_{\tilde{\mathbb{E}}} \ln \lvert \xi \! - \! 
\alpha_{k} \rvert \, \md \mu_{\widetilde{V}}^{f}(\xi) 
\right) \underset{\tilde{\varepsilon} \downarrow 0}{=} 
\mathcal{O} \left(\mathfrak{c}_{5}(n,k,z_{o}) 
\tilde{\varepsilon}^{2} \right), \\
\sum_{q=1}^{\mathfrak{s}-2} \dfrac{\varkappa_{nk 
\tilde{k}_{q}}}{n} \left(\int_{\tilde{\mathbb{E}}} \ln 
\lvert \xi \! - \! \alpha_{p_{q}} \rvert \, \tilde{\phi}_{
\tilde{\varepsilon}}(\xi) \, \md \xi \! - \! \int_{
\tilde{\mathbb{E}}} \ln \lvert \xi \! - \! \alpha_{p_{q}} 
\rvert \, \md \mu_{\widetilde{V}}^{f}(\xi) \right) 
\underset{\tilde{\varepsilon} \downarrow 0}{=} 
\mathcal{O} \left(\mathfrak{c}_{6}(n,k,z_{o}) 
\tilde{\varepsilon}^{2} \right), \\
\left(\dfrac{\varkappa_{nk} \! - \! 1}{n} \right) 
\int_{\tilde{\mathbb{E}}} \ln \lvert \xi \! - \! 
\alpha_{k} \rvert \, \md \mu_{\widetilde{V}}^{f}(\xi) 
\! = \! \left(\dfrac{\varkappa_{nk} \! - \! 1}{n} 
\right) \int_{J_{f}} \ln \lvert \xi \! - \! \alpha_{k} 
\rvert \, \md \mu_{\widetilde{V}}^{f}(\xi) 
\underset{\tilde{\varepsilon} \downarrow 0}{=} 
\mathcal{O}(\mathfrak{c}_{7}(n,k,z_{o})), \\
\sum_{q=1}^{\mathfrak{s}-2} \dfrac{\varkappa_{nk 
\tilde{k}_{q}}}{n} \int_{\tilde{\mathbb{E}}} \ln 
\lvert \xi \! - \! \alpha_{p_{q}} \rvert \, \md \mu_{
\widetilde{V}}^{f}(\xi) \! = \! \sum_{q=1}^{\mathfrak{s}
-2} \dfrac{\varkappa_{nk \tilde{k}_{q}}}{n} \int_{J_{f}} 
\ln \lvert \xi \! - \! \alpha_{p_{q}} \rvert \, \md \mu_{
\widetilde{V}}^{f}(\xi) \underset{\tilde{\varepsilon} 
\downarrow 0}{=} \mathcal{O}(\mathfrak{c}_{8}
(n,k,z_{o})),
\end{gather*}
where $\mathfrak{c}_{m}(n,k,z_{o}) \! =_{\underset{z_{o}
=1+o(1)}{\mathscr{N},n \to \infty}} \! \mathcal{O}(1)$, 
$m \! = \! 5,6,7,8$; hence (assuming commutativity of 
limits),
\begin{equation} \label{eqlmrtz5} 
\iint_{\tilde{\mathbb{E}}^{2}} \tilde{\mathrm{F}}
(\xi,\tau) \, \tilde{\phi}_{\tilde{\varepsilon}}(\xi) \, \md 
\xi \, \tilde{\phi}_{\tilde{\varepsilon}}(\tau) \, \md 
\tau \underset{\tilde{\varepsilon} \downarrow 0}{=} 
\mathcal{O} \left(\mathfrak{c}_{9}(n,k,z_{o}) \! + \! 
\mathfrak{c}_{10}(n,k,z_{o}) \tilde{\varepsilon}^{2} \right),
\end{equation}
where $\mathfrak{c}_{m}(n,k,z_{o}) \! =_{\underset{z_{o}
=1+o(1)}{\mathscr{N},n \to \infty}} \! \mathcal{O}(1)$, 
$m \! = \! 9,10$. Hence, for $n \! \in \! \mathbb{N}$ 
and $k \! \in \! \lbrace 1,2,\dotsc,K \rbrace$ such 
that $\alpha_{p_{\mathfrak{s}}} \! := \! \alpha_{k} 
\! \neq \! \infty$, it follows via 
Equations~\eqref{eqlmrtz2}--\eqref{eqlmrtz5}, the 
weak-$\ast$ convergence $\tilde{\phi}_{\tilde{\varepsilon}}
(x) \, \md x \! \overset{\ast}{\to} \! \md 
\mu_{\widetilde{V}}^{f}(x)$ as $\tilde{\varepsilon} \! 
\downarrow \! 0$, the integrability of $\tilde{\phi}_{
\tilde{\varepsilon}} \ln (\tilde{\phi}_{\tilde{\varepsilon}})$ 
(on $\tilde{\mathbb{E}})$, the fact that $\mathrm{I}_{
\widetilde{V}}^{f}[\tilde{\phi}_{\tilde{\varepsilon}}] \! 
\leqslant \! E_{\widetilde{V}}^{f} \! + \! \tilde{\varepsilon}/2$, 
and the convergence $\mathrm{I}_{\widetilde{V}}^{f}
[\tilde{\phi}_{\tilde{\varepsilon}}] \! \to \! \mathrm{I}_{
\widetilde{V}}^{f}[\mu_{\widetilde{V}}^{f}] \! = \! 
E_{\widetilde{V}}^{f}$ as $\tilde{\varepsilon} \! \downarrow 
\! 0$, that (assuming commutativity of limits)
\begin{equation} \label{eqlmrtz6} 
\dfrac{\mathcal{N}}{n} \dfrac{1}{\mathcal{N}(\mathcal{N} 
\! - \! 1)} \ln (\tilde{\mathscr{Z}}^{n}_{k}) \underset{
\underset{z_{o}=1+o(1)}{\mathscr{N},n \to \infty}}{\geqslant} 
-\left(E_{\widetilde{V}}^{f} \! + \! \tilde{\varepsilon} \right) 
\left(1 \! + \! \mathcal{O} \left(\dfrac{\mathfrak{c}_{11}
(n,k,z_{o},\tilde{\varepsilon})}{n} \right) \right) \quad 
\text{as} \quad \tilde{\varepsilon} \! \downarrow \! 0,
\end{equation}
where $\mathfrak{c}_{11}(n,k,z_{o},\tilde{\varepsilon}) \! 
=_{\underset{z_{o}=1+o(1)}{\mathscr{N},n \to \infty}} \! 
\mathcal{O}(1)$ as $\tilde{\varepsilon} \! \downarrow \! 0$.

For $n \! \in \! \mathbb{N}$ and $k \! \in \! \lbrace 1,2,\dotsc,K 
\rbrace$ such that $\alpha_{p_{\mathfrak{s}}} \! := \! \alpha_{k} 
\! \neq \! \infty$, given $\tilde{\eta} \! = \! \tilde{\eta}(n,k,z_{o}) 
\! > \! 0$ and $\mathcal{O}(1)$, set $\tilde{\mathbb{A}}_{\mathcal{N},
\tilde{\eta}} \! := \! \lbrace \mathstrut (x_{1},x_{2},\dotsc,
x_{\mathcal{N}}) \! \in \! \mathbb{R}^{\mathcal{N}}; \, 
(\mathcal{N}(\mathcal{N} \! - \! 1))^{-1} \mathscr{K}_{
\mathcal{N}}^{\widetilde{V},f}(x_{1},x_{2},\dotsc,x_{\mathcal{N}}) 
\! \leqslant \! E_{\widetilde{V}}^{f} \! + \! \tilde{\eta} \rbrace$: 
then, via the Estimate~\eqref{eqlmrtz6}, it follows that, 
for any $\mathfrak{a} \! = \! \mathfrak{a}(n,k,z_{o}) \! \geqslant \! 0$, 
upon choosing $2 \tilde{\varepsilon} \! < \! \tilde{\eta}$, say,
\begin{align*}
\operatorname{Prob} \, (\mathbb{R}^{\mathcal{N}} \setminus 
\tilde{\mathbb{A}}_{\mathcal{N},\tilde{\eta}+\mathfrak{a}}) 
\underset{\underset{z_{o}=1+o(1)}{\mathscr{N},n \to \infty}}{\leqslant}& 
\, \me^{\frac{n}{\mathcal{N}} \mathcal{N}(\mathcal{N}-1)
(E_{\widetilde{V}}^{f}+ \tilde{\varepsilon})(1+O(n^{-1}))} 
\idotsint\limits_{\lbrace \mathstrut (\xi_{1},\xi_{2},\dotsc,
\xi_{\mathcal{N}}) \in \mathbb{R}^{\mathcal{N}}; \, 
(\mathcal{N}(\mathcal{N}-1))^{-1} \mathscr{K}_{\mathcal{N}}^{
\widetilde{V},f}(\xi_{1},\xi_{2},\dotsc,\xi_{\mathcal{N}})>
E_{\widetilde{V}}^{f}+ \tilde{\eta}+ \mathfrak{a} \rbrace} 
\me^{-\frac{n}{\mathcal{N}} \mathcal{N}(\mathcal{N}-1)
(E_{\widetilde{V}}^{f}+ \tilde{\eta}+ \mathfrak{a})} \\
\times& \, \me^{-\frac{n}{\mathcal{N}} \left(\sum_{j=1}^{
\mathcal{N}} \widetilde{V}(\xi_{j})+\mathcal{O}(\frac{1}{
\mathcal{N}}) \sum_{\underset{i \neq j}{i,j=1}}^{\mathcal{N}} 
\tilde{\mathrm{F}}(\xi_{i},\xi_{j}) \right)} \, \md \xi_{1} 
\, \md \xi_{2} \, \dotsb \, \md \xi_{\mathcal{N}} \\
\underset{\underset{z_{o}=1+o(1)}{\mathscr{N},n \to 
\infty}}{\leqslant}& \, \me^{-\frac{n}{\mathcal{N}} 
\mathcal{N}(\mathcal{N}-1) \mathfrak{a}} 
\me^{-\frac{1}{2} \frac{n}{\mathcal{N}} \mathcal{N}
(\mathcal{N}-1) \tilde{\eta}(1+\mathcal{O}(n^{-1}))} 
\idotsint\limits_{\mathbb{R}^{\mathcal{N}}} 
\me^{-\frac{n}{\mathcal{N}} \left(\sum_{j=1}^{\mathcal{N}} 
\widetilde{V}(\xi_{j})+\mathcal{O}(\frac{1}{\mathcal{N}}) 
\sum_{\underset{i \neq j}{i,j=1}}^{\mathcal{N}} \tilde{
\mathrm{F}}(\xi_{i},\xi_{j}) \right)} \, \md \xi_{1} \, \md 
\xi_{2} \, \dotsb \, \md \xi_{\mathcal{N}};
\end{align*}
one now shows, via Equation~\eqref{eqlmrtzf} and an enumeration 
argument, that
\begin{equation*}
\operatorname{Prob} \, (\mathbb{R}^{\mathcal{N}} \setminus 
\tilde{\mathbb{A}}_{\mathcal{N},\tilde{\eta}+\mathfrak{a}}) 
\underset{\underset{z_{o}=1+o(1)}{\mathscr{N},n \to 
\infty}}{\leqslant} \me^{-\frac{n}{\mathcal{N}} 
\mathcal{N}(\mathcal{N}-1) \mathfrak{a}} \left(
\int_{\mathbb{R}} \dfrac{\me^{-\frac{n}{\mathcal{N}} 
\widetilde{V}(\xi)}}{\prod_{\underset{q \neq \mathfrak{s}-1}{q
=1}}^{\mathfrak{s}} \lvert \xi \! - \! \alpha_{p_{q}} \rvert^{
\tilde{\gamma}_{q}(1+\mathcal{O}(n^{-1}))}} \md \xi 
\right)^{\mathcal{N}} \me^{-\frac{1}{2} \frac{n}{\mathcal{N}} 
\mathcal{N}(\mathcal{N}-1) \tilde{\eta}(1+\mathcal{O}(n^{-1}))},
\end{equation*}
where $\tilde{\gamma}_{q} \! = \! 2K^{-1} \gamma_{i(q)_{k_{q}}}$, 
$q \! \in \! \lbrace 1,2,\dotsc,\mathfrak{s} \! - \! 2 \rbrace$, and 
$\tilde{\gamma}_{q} \! = \! 2K^{-1} \gamma_{k}$, $q \! = \! 
\mathfrak{s}$. Recall {}from the proof of Lemma~\ref{lem3.1} 
that, for $n \! \in \! \mathbb{N}$ and $k \! \in \! \lbrace 
1,2,\dotsc,K \rbrace$ such that $\alpha_{p_{\mathfrak{s}}} 
\! := \! \alpha_{k} \! \neq \! \infty$, there exists $T_{M_{f}} 
\! > \! 1$ (e.g., $T_{M_{f}} \! = \! K(1 \! + \! \max_{q=1,\dotsc,
\mathfrak{s}-2,\mathfrak{s}} \lbrace \lvert \alpha_{p_{q}} 
\rvert \rbrace\! + \! 3(\min_{i \neq j \in \lbrace 1,\dotsc,
\mathfrak{s}-2,\mathfrak{s} \rbrace} \lbrace \lvert \alpha_{p_{i}} 
\! - \! \alpha_{p_{j}} \rvert \rbrace)^{-1}))$, with $\mathscr{O}_{
\frac{1}{T_{M_{f}}}}(\alpha_{p_{i}}) \cap \mathscr{O}_{\frac{1}{T_{
M_{f}}}}(\alpha_{p_{j}}) \! = \! \varnothing$ $\forall$ $i \! \neq \! 
j \! \in \! \lbrace 1,\dotsc,\mathfrak{s} \! - \! 2,\mathfrak{s} 
\rbrace$, such that $\mathbb{R} \setminus \mathfrak{D}_{
M_{f}} \supseteq J_{f}$, where $\mathfrak{D}_{M_{f}} \! := \! 
\lbrace \lvert x \rvert \! \geqslant \! T_{M_{f}} \rbrace \cup 
\cup_{\underset{q \neq \mathfrak{s}-1}{q=1}}^{\mathfrak{s}} 
\operatorname{clos}(\mathscr{O}_{\frac{1}{T_{M_{f}}}}
(\alpha_{p_{q}}))$. Write $\mathbb{R} \! = \! (\mathbb{R} 
\setminus \mathfrak{D}_{M_{f}}) \cup \mathfrak{D}_{M_{f}} \! 
= \! ((\mathbb{R} \setminus \mathfrak{D}_{M_{f}}) \setminus 
J_{f}) \cup J_{f} \cup \mathfrak{D}_{M_{f}}$, with $\mathbb{R} 
\setminus \mathfrak{D}_{M_{f}} \cap \mathfrak{D}_{M_{f}} 
\! = \! \varnothing \! = \! ((\mathbb{R} \setminus 
\mathfrak{D}_{M_{f}}) \setminus J_{f}) \cap J_{f}$; hence,
\begin{equation*}
\int_{\mathbb{R}} \dfrac{\me^{-\frac{n}{\mathcal{N}} 
\widetilde{V}(\xi)}}{\prod_{\underset{q \neq \mathfrak{s}-
1}{q=1}}^{\mathfrak{s}} \lvert \xi \! - \! \alpha_{p_{q}} 
\rvert^{\tilde{\gamma}_{q}(1+\mathcal{O}(n^{-1}))}} \, \md \xi 
\! = \! \left(\int_{(\mathbb{R} \setminus \mathfrak{D}_{M_{f}}) 
\setminus J_{f}} \! + \! \int_{J_{f}} \! + \! \int_{\lbrace \lvert 
\xi \rvert \geqslant T_{M_{f}} \rbrace} \! + \! \sum_{
\substack{q=1\\q \neq \mathfrak{s}-1}}^{\mathfrak{s}} 
\int_{\mathscr{O}_{\frac{1}{T_{M_{f}}}}(\alpha_{p_{q}})} \right) 
\dfrac{\me^{-\frac{n}{\mathcal{N}} \widetilde{V}(\xi)}}{
\prod_{\underset{q \neq \mathfrak{s}-1}{q=1}}^{\mathfrak{s}} 
\lvert \xi \! - \! \alpha_{p_{q}} \rvert^{\tilde{\gamma}_{q}
(1+\mathcal{O}(n^{-1}))}} \, \md \xi.
\end{equation*}
Recalling {}from the proof of Lemma~\ref{lem3.1} that, 
for $n \! \in \! \mathbb{N}$ and $k \! \in \! \lbrace 
1,2,\dotsc,K \rbrace$ such that $\alpha_{p_{\mathfrak{s}}} 
\! := \! \alpha_{k} \! \neq \! \infty$, in the double-scaling 
limit $\mathscr{N},n \! \to \! \infty$ such that 
$z_{o} \! = \! 1 \! + \! o(1)$, via 
Conditions~\eqref{eq20}--\eqref{eq22} for regular 
$\widetilde{V} \colon \overline{\mathbb{R}} \setminus 
\lbrace \alpha_{1},\alpha_{2},\dotsc,\alpha_{K} \rbrace 
\! \to \! \mathbb{R}$, there exists $\tilde{c}_{\infty} \! = \! 
\tilde{c}_{\infty}(n,k,z_{o}) \! > \! 0$ and $\mathcal{O}(1)$ 
such that, for $x \! \in \! \lbrace \lvert x \rvert \! \geqslant \! 
T_{M_{f}} \rbrace$, $\widetilde{V}(x) \! \geqslant \! (1 \! + \! 
\tilde{c}_{\infty}) \ln (1 \! + \! x^{2})$, and, for $q \! = \! 1,
\dotsc,\mathfrak{s} \! - \! 2,\mathfrak{s}$, there exists 
$\tilde{c}_{q} \! = \! \tilde{c}_{q}(n,k,z_{o}) \! > \! 0$ 
and $\mathcal{O}(1)$ such that, for $x \! \in \! 
\mathscr{O}_{\frac{1}{T_{M_{f}}}}(\alpha_{p_{q}})$, 
$\widetilde{V}(x) \! \geqslant \! (1 \! + \! \tilde{c}_{q}) 
\ln (1 \! + \! (x \! - \! \alpha_{p_{q}})^{-2})$, one shows, 
via item~\pmb{(2)} of Lemma~\ref{lem3.7}, and choosing 
$\tilde{c}_{q}$ so that $1 \! + \! \tilde{c}_{q} \! \geqslant 
\! \gamma_{i(q)_{k_{q}}}$, $q \! = \! 1,\dotsc,\mathfrak{s} 
\! - \! 2,\mathfrak{s}$, that, via the monotonicity of 
$\exp (\pmb{\cdot})$,
\begin{gather*}
\int_{(\mathbb{R} \setminus \mathfrak{D}_{M_{f}}) 
\setminus J_{f}} \dfrac{\me^{-\frac{n}{\mathcal{N}} 
\widetilde{V}(\xi)}}{\prod_{\underset{q \neq \mathfrak{s}-
1}{q=1}}^{\mathfrak{s}} \lvert \xi \! - \! \alpha_{p_{q}} 
\rvert^{\tilde{\gamma}_{q}(1+\mathcal{O}(n^{-1}))}} \, \md 
\xi \underset{\underset{z_{o}=1+o(1)}{\mathscr{N},n \to 
\infty}}{\leqslant} \dfrac{\me^{-K^{-1} \inf \lbrace \mathstrut 
\widetilde{V}(x); \, x \in (\mathbb{R} \setminus \mathfrak{D}_{M_{f}}) 
\setminus J_{f} \rbrace} \operatorname{meas}((\mathbb{R} 
\setminus \mathfrak{D}_{M_{f}}) \setminus J_{f})(1 \! + \! o(1))}{
\prod_{\underset{q \neq \mathfrak{s}-1}{q=1}}^{\mathfrak{s}}
(\inf \lbrace \mathstrut \lvert x \! - \! \alpha_{p_{q}} \rvert; 
\, x \! \in \! (\mathbb{R} \setminus \mathfrak{D}_{M_{f}}) 
\setminus J_{f} \rbrace)^{\tilde{\gamma}_{q}}} \! 
=: \! \mathfrak{c}_{12}(n,k,z_{o}), \\
\int_{J_{f}} \dfrac{\me^{-\frac{n}{\mathcal{N}} 
\widetilde{V}(\xi)}}{\prod_{\underset{q \neq \mathfrak{s}-
1}{q=1}}^{\mathfrak{s}} \lvert \xi \! - \! \alpha_{p_{q}} 
\rvert^{\tilde{\gamma}_{q}(1+\mathcal{O}(n^{-1}))}} \, 
\md \xi \underset{\underset{z_{o}=1+o(1)}{\mathscr{N},n 
\to \infty}}{\leqslant} \dfrac{\me^{-K^{-1} \inf \lbrace 
\mathstrut \widetilde{V}(x); \, x \in J_{f} \rbrace} 
\sum_{j=1}^{\mathcal{N}} \lvert \tilde{b}_{j-1} \! - \! 
\tilde{a}_{j} \rvert(1 \! + \! o(1))}{\prod_{\underset{q \neq 
\mathfrak{s}-1}{q=1}}^{\mathfrak{s}}(\inf \lbrace \mathstrut 
\lvert x \! - \! \alpha_{p_{q}} \rvert; \, x \! \in \! J_{f} 
\rbrace)^{\tilde{\gamma}_{q}}} \! =: \! \mathfrak{c}_{13}
(n,k,z_{o}), \\
\int_{\lbrace \lvert \xi \rvert \geqslant T_{M_{f}} \rbrace} 
\dfrac{\me^{-\frac{n}{\mathcal{N}} \widetilde{V}(\xi)}}{
\prod_{\underset{q \neq \mathfrak{s}-1}{q=1}}^{\mathfrak{s}} 
\lvert \xi \! - \! \alpha_{p_{q}} \rvert^{\tilde{\gamma}_{q}
(1+\mathcal{O}(n^{-1}))}} \, \md \xi \underset{\underset{
z_{o}=1+o(1)}{\mathscr{N},n \to \infty}}{\leqslant} 
\dfrac{2KT_{M_{f}}(K \! + \! 2(1 \! - \! \gamma_{i(\mathfrak{s}
-1)_{k_{\mathfrak{s}-1}}}))^{-1}(1 \! + \! o(1))}{(T_{M_{f}}^{2} 
\! + \! 1)^{\frac{\tilde{c}_{\infty}+1}{K}} \prod_{\underset{q 
\neq \mathfrak{s}-1}{q=1}}^{\mathfrak{s}}(T_{M_{f}} \! - \! 
\lvert \alpha_{p_{q}} \rvert)^{\tilde{\gamma}_{q}}} \! =: \! 
\mathfrak{c}_{14}(n,k,z_{o}), \\
\sum_{\substack{q=1\\q \neq \mathfrak{s}-1}}^{\mathfrak{s}} 
\int_{\mathscr{O}_{\frac{1}{T_{M_{f}}}}(\alpha_{p_{q}})} 
\dfrac{\me^{-\frac{n}{\mathcal{N}} \widetilde{V}(\xi)}}{
\prod_{\underset{q \neq \mathfrak{s}-1}{q=1}}^{\mathfrak{s}} 
\lvert \xi \! - \! \alpha_{p_{q}} \rvert^{\tilde{\gamma}_{q}
(1+\mathcal{O}(n^{-1}))}} \, \md \xi \underset{\underset{
z_{o}=1+o(1)}{\mathscr{N},n \to \infty}}{\leqslant} 
\sum_{\substack{q=1\\q \neq \mathfrak{s}-1}}^{\mathfrak{s}} 
\dfrac{2T_{M_{f}}^{-\left(1+2K^{-1}(1+ \tilde{c}_{q})-\tilde{
\gamma}_{q} \right)}(1 \! + \! o(1))}{(1 \! + \! 2K^{-1}(1 \! + 
\! \tilde{c}_{q}) \! - \! \tilde{\gamma}_{q}) \prod_{\underset{
q^{\prime} \neq q,\mathfrak{s}-1}{q^{\prime}=1}}^{\mathfrak{s}} 
\lvert \alpha_{p_{q}} \! - \! \alpha_{p_{q^{\prime}}} \rvert^{
\tilde{\gamma}_{q}}} \! =: \! \mathfrak{c}_{15}(n,k,z_{o}),
\end{gather*}
where $\mathfrak{c}_{m}(n,k,z_{o}) \! =_{\underset{z_{o}
=1+o(1)}{\mathscr{N},n \to \infty}} \! \mathcal{O}(1)$, 
$m \! = \! 12,13,14,15$, whence
\begin{equation*}
\int_{\mathbb{R}} \dfrac{\me^{-\frac{n}{\mathcal{N}} 
\widetilde{V}(\xi)}}{\prod_{\underset{q \neq \mathfrak{s}-1}{q
=1}}^{\mathfrak{s}} \lvert \xi \! - \! \alpha_{p_{q}} \rvert^{
\tilde{\gamma}_{q}(1+\mathcal{O}(n^{-1}))}} \, \md \xi 
\underset{\underset{z_{o}=1+o(1)}{\mathscr{N},n \to \infty}}{
\leqslant} \mathfrak{c}_{16}(n,k,z_{o}),
\end{equation*}
with $\mathfrak{c}_{16}(n,k,z_{o}) \! =_{\underset{z_{o}=1
+o(1)}{\mathscr{N},n \to \infty}} \! \mathcal{O}(1)$, which 
implies that
\begin{equation*}
\left(\int_{\mathbb{R}} \dfrac{\me^{-\frac{n}{\mathcal{N}} 
\widetilde{V}(\xi)}}{\prod_{\underset{q \neq 
\mathfrak{s}-1}{q=1}}^{\mathfrak{s}} \lvert \xi \! - \! 
\alpha_{p_{q}} \rvert^{\tilde{\gamma}_{q}(1+\mathcal{O}
(n^{-1}))}} \md \xi \right)^{\mathcal{N}} \me^{-\frac{1}{2} 
\frac{n}{\mathcal{N}} \mathcal{N}(\mathcal{N}-1) \tilde{\eta}
(1+\mathcal{O}(n^{-1}))} \underset{\underset{z_{o}=1+o(1)}{
\mathscr{N},n \to \infty}}{\leqslant} \mathfrak{c}_{17}(n,k,z_{o}),
\end{equation*}
where $\mathfrak{c}_{17}(n,k,z_{o}) \! \leqslant_{\underset{
z_{o}=1+o(1)}{\mathscr{N},n \to \infty}} \! 1$; hence, for $n 
\! \in \! \mathbb{N}$ and $k \! \in \! \lbrace 1,2,\dotsc,K 
\rbrace$ such that $\alpha_{p_{\mathfrak{s}}} \! := \! 
\alpha_{k} \! \neq \! \infty$, one arrives at
\begin{equation} \label{eqlmrtz7} 
\operatorname{Prob} \, (\mathbb{R}^{\mathcal{N}} \setminus 
\tilde{\mathbb{A}}_{\mathcal{N},\tilde{\eta}+\mathfrak{a}}) 
\underset{\underset{z_{o}=1+o(1)}{\mathscr{N},n \to \infty}}{
\leqslant} \mathfrak{c}_{17}(n,k,z_{o}) \me^{-\frac{n}{\mathcal{N}} 
\mathcal{N}(\mathcal{N}-1) \mathfrak{a}}.
\end{equation}

For $n \! \in \! \mathbb{N}$ and $k \! \in \! \lbrace 1,2,\dotsc,
K \rbrace$ such that $\alpha_{p_{\mathfrak{s}}} \! := \! \alpha_{k} 
\! \neq \! \infty$, take $\mathfrak{a} \! = \! \tilde{\eta}$, and 
consider the set $\tilde{\mathbb{A}}_{\mathcal{N},2 \tilde{\eta}} 
\! := \! \lbrace \mathstrut (x_{1},x_{2},\dotsc,x_{\mathcal{N}}) 
\! \in \! \mathbb{R}^{\mathcal{N}}; \, (\mathcal{N}
(\mathcal{N} \! - \! 1))^{-1} \mathscr{K}_{\mathcal{N}}^{
\widetilde{V},f}(x_{1},x_{2},\dotsc,x_{\mathcal{N}}) \! \leqslant 
\! E_{\widetilde{V}}^{f} \! + \! 2 \tilde{\eta} \rbrace$, and, for $m 
\! \in \! \mathbb{N}$, let $\phi_{0}$ $(= \! \phi_{0}(n,k,z_{o}))$ 
$\colon \mathbb{R}^{m} \! \to \! \mathbb{R}$ be bounded 
and continuous, and consider, in the double-scaling limit 
$\mathscr{N},n \! \to \! \infty$ such that $z_{o} \! = \! 
1 \! + \! o(1)$, the quantity $\mathcal{N}^{-1} \ln \tilde{
\mathcal{E}}^{n}_{k}(\exp (\mathcal{N}^{-(m-1)} \linebreak[4] 
\pmb{\cdot} \sum_{i_{1},i_{2},\dotsc,i_{m}} \phi_{0}(x_{i_{1}},
x_{i_{2}},\dotsc,x_{i_{m}})))$, where $i_{j} \! \in \! \lbrace 1,2,
\dotsc,\mathcal{N} \rbrace$, $j \! = \! 1,2,\dotsc,m$, and 
$\tilde{\mathcal{E}}^{n}_{k}$ denotes the expectation with 
respect to the multi-dimensional probability measure 
$\tilde{\mathcal{P}}^{n}_{k}(x_{1},x_{2},\dotsc,
x_{\mathcal{N}}) \, \md x_{1} \, \md x_{2} \, \dotsb \, \md 
x_{\mathcal{N}}$.\footnote{\, $\tilde{\mathcal{E}}^{n}_{k}
(F(x_{1},x_{2},\dotsc,x_{\mathcal{N}})) \! := \! \idotsint_{
\mathbb{R}^{\raise-0.7ex\hbox{$\scriptscriptstyle \mathcal{N}$}}}
F(\xi_{1},\xi_{2},\dotsc,\xi_{\mathcal{N}}) \tilde{\mathcal{P}}^{n}_{k}
(\xi_{1},\xi_{2},\dotsc,\xi_{\mathcal{N}}) \, \md \xi_{1} \, \md 
\xi_{2} \, \dotsb \, \md \xi_{\mathcal{N}}$.} Denote the characteristic 
function of the set $\tilde{\mathbb{A}}_{\mathcal{N},2 \tilde{\eta}}$ 
as $\chi_{\tilde{\mathbb{A}}_{\mathcal{N},2 \tilde{\eta}}}(\vec{x}) 
\! := \! \chi_{\tilde{\mathbb{A}}_{\mathcal{N},2 \tilde{\eta}}}
(x_{1},x_{2},\dotsc,x_{\mathcal{N}})$. Via the linearity of 
$\tilde{\mathcal{E}}^{n}_{k}$, one shows, for $n \! \in \! 
\mathbb{N}$ and $k \! \in \! \lbrace 1,2,\dotsc,K \rbrace$ such 
that $\alpha_{p_{\mathfrak{s}}} \! := \! \alpha_{k} \! \neq \! \infty$, 
that
\begin{align*}
\limsup_{\underset{z_{o}=1+o(1)}{\mathscr{N},n \to \infty}} 
\dfrac{1}{\mathcal{N}} \ln \tilde{\mathcal{E}}^{n}_{k} \left(
\me^{\mathcal{N}^{-(m-1)} \sum_{i_{1},i_{2},\dotsc,i_{m}} 
\phi_{0}(x_{i_{1}},x_{i_{2}},\dotsc,x_{i_{m}})} \right) =& \, 
\limsup_{\underset{z_{o}=1+o(1)}{\mathscr{N},n \to \infty}} 
\dfrac{1}{\mathcal{N}} \left(\vphantom{M^{M^{M^{M^{M^{M}}}}}} 
\ln \tilde{\mathcal{E}}^{n}_{k} \left(\chi_{\tilde{\mathbb{A}}_{
\mathcal{N},2 \tilde{\eta}}}(\vec{x}) \me^{\mathcal{N}^{-(m-1)} 
\sum_{i_{1},i_{2},\dotsc,i_{m}} \phi_{0}(x_{i_{1}},x_{i_{2}},\dotsc,
x_{i_{m}})} \right) \right. \\
+&\left. \, \ln \left(1 \! + \! \dfrac{\tilde{\mathcal{E}}^{n}_{k} 
\left((1 \! - \! \chi_{\tilde{\mathbb{A}}_{\mathcal{N},2 \tilde{
\eta}}}(\vec{x})) \me^{\mathcal{N}^{-(m-1)} \sum_{i_{1},i_{2},
\dotsc,i_{m}} \phi_{0}(x_{i_{1}},x_{i_{2}},\dotsc,x_{i_{m}})} 
\right)}{\tilde{\mathcal{E}}^{n}_{k} \left(\chi_{\tilde{\mathbb{
A}}_{\mathcal{N},2 \tilde{\eta}}}(\vec{x}) \me^{\mathcal{N}^{-
(m-1)} \sum_{i_{1},i_{2},\dotsc,i_{m}} \phi_{0}(x_{i_{1}},x_{i_{2}},
\dotsc,x_{i_{m}})} \right)} \right) \right).
\end{align*}
Since $\phi_{0}$ is bounded, there exists $\tilde{M}_{f} \! = \! 
\tilde{M}_{f}(n,k,z_{o}) \! > \! 0$ and $\mathcal{O}(1)$ such that 
$\lvert \phi_{0}(x_{i_{1}},x_{i_{2}},\dotsc,x_{i_{m}}) \rvert \! 
\leqslant \! \tilde{M}_{f}$; hence, via an application of the 
Fundamental Counting Principle, $-\mathcal{N} \tilde{M}_{f} \! 
\leqslant \! \mathcal{N}^{-(m-1)} \sum_{i_{1},i_{2},\dotsc,i_{m}} 
\phi_{0}(x_{i_{1}},x_{i_{2}},\dotsc,x_{i_{m}}) \! \leqslant \! 
\mathcal{N} \tilde{M}_{f}$, which implies, by the linearity of 
$\tilde{\mathcal{E}}^{n}_{k}$, the fact that $\tilde{\mathcal{E}}^{n}_{k}
(1) \! = \! 1$, and the monotonicity of $\exp (\pmb{\cdot})$ 
(that is, $\me^{-\mathcal{N} \tilde{M}_{f}} \! \leqslant \! 
\tilde{\mathcal{E}}^{n}_{k}(\me^{\mathcal{N}^{-(m-1)} 
\sum_{i_{1},i_{2},\dotsc,i_{m}} \phi_{0}(x_{i_{1}},x_{i_{2}},
\dotsc,x_{i_{m}})}) \! \leqslant \! \me^{\mathcal{N} 
\tilde{M}_{f}})$, that
\begin{equation*}
\dfrac{\tilde{\mathcal{E}}^{n}_{k} \left((1 \! - \! \chi_{
\tilde{\mathbb{A}}_{\mathcal{N},2 \tilde{\eta}}}(\vec{x})) 
\me^{\mathcal{N}^{-(m-1)} \sum_{i_{1},i_{2},\dotsc,i_{m}} 
\phi_{0}(x_{i_{1}},x_{i_{2}},\dotsc,x_{i_{m}})} \right)}{\tilde{
\mathcal{E}}^{n}_{k} \left(\chi_{\tilde{\mathbb{A}}_{\mathcal{N},
2 \tilde{\eta}}}(\vec{x}) \me^{\mathcal{N}^{-(m-1)} \sum_{
i_{1},i_{2},\dotsc,i_{m}} \phi_{0}(x_{i_{1}},x_{i_{2}},\dotsc,
x_{i_{m}})} \right)} \underset{\underset{z_{o}=1+o(1)}{
\mathscr{N},n \to \infty}}{\leqslant} \dfrac{\me^{2 \mathcal{N} 
\tilde{M}_{f}} \tilde{\mathcal{E}}^{n}_{k}(1 \! - \! \chi_{\tilde{
\mathbb{A}}_{\mathcal{N},2 \tilde{\eta}}}(\vec{x}))}{1 \! - \! 
\tilde{\mathcal{E}}^{n}_{k}(1 \! - \! \chi_{\tilde{\mathbb{A}}_{
\mathcal{N},2 \tilde{\eta}}}(\vec{x}))}.
\end{equation*}
Via the deomposition $\mathbb{R}^{\mathcal{N}} 
\! = \!  (\mathbb{R}^{\mathcal{N}} \setminus 
\tilde{\mathbb{A}}_{\mathcal{N},2 \tilde{\eta}}) \cup \tilde{\mathbb{
A}}_{\mathcal{N},2 \tilde{\eta}}$ (with $(\mathbb{R}^{\mathcal{N}} 
\setminus \tilde{\mathbb{A}}_{\mathcal{N},2 \tilde{\eta}}) \cap 
\tilde{\mathbb{A}}_{\mathcal{N},2 \tilde{\eta}} \! = \! \varnothing)$ 
and the Estimate~\eqref{eqlmrtz7}, it follows, via the linearity of 
$\tilde{\mathcal{E}}^{n}_{k}$, that
\begin{align*}
\tilde{\mathcal{E}}^{n}_{k}(1 \! - \! \chi_{\tilde{\mathbb{
A}}_{\mathcal{N},2 \tilde{\eta}}}(\vec{x})) =& \, 
\left(\idotsint\limits_{\mathbb{R}^{\mathcal{N}} \setminus 
\tilde{\mathbb{A}}_{\mathcal{N},2 \tilde{\eta}}} \! + \! 
\idotsint\limits_{\tilde{\mathbb{A}}_{\mathcal{N},2 
\tilde{\eta}}} \right)(1 \! - \! \chi_{\tilde{\mathbb{A}}_{
\mathcal{N},2 \tilde{\eta}}}(\vec{\xi})) \tilde{\mathcal{P}}^{n}_{k}
(\xi_{1},\xi_{2},\dotsc,\xi_{\mathcal{N}}) \, \md \xi_{1} 
\, \md \xi_{2} \, \dotsb \, \md \xi_{\mathcal{N}} \\
=& \, \idotsint\limits_{\lbrace \mathstrut (\xi_{1},\xi_{2},
\dotsc,\xi_{\mathcal{N}}) \in \mathbb{R}^{\mathcal{N}}; \, 
(\mathcal{N}(\mathcal{N}-1))^{-1} \mathscr{K}_{\mathcal{N}}^{
\widetilde{V},f}(\xi_{1},\xi_{2},\dotsc,\xi_{\mathcal{N}})>
E_{\widetilde{V}}^{f}+ 2 \tilde{\eta} \rbrace} \tilde{\mathcal{P}}^{n}_{k}
(\xi_{1},\xi_{2},\dotsc,\xi_{\mathcal{N}}) \, \md \xi_{1} \, \md 
\xi_{2} \, \dotsb \, \md \xi_{\mathcal{N}} \\
=& \, \operatorname{Prob} \, (\mathbb{R}^{\mathcal{N}} \setminus 
\tilde{\mathbb{A}}_{\mathcal{N},2\tilde{\eta}}) \underset{
\underset{z_{o}=1+o(1)}{\mathscr{N},n \to \infty}}{\leqslant} 
\mathfrak{c}_{17}(n,k,z_{o}) \me^{-\frac{n}{\mathcal{N}} 
\mathcal{N}(\mathcal{N} \! - \! 1) \tilde{\eta}} \quad \Rightarrow
\end{align*}
\begin{equation*}
\dfrac{\tilde{\mathcal{E}}^{n}_{k} \left((1 \! - \! \chi_{
\tilde{\mathbb{A}}_{\mathcal{N},2 \tilde{\eta}}}(\vec{x})) 
\me^{\mathcal{N}^{-(m-1)} \sum_{i_{1},i_{2},\dotsc,i_{m}} 
\phi_{0}(x_{i_{1}},x_{i_{2}},\dotsc,x_{i_{m}})} \right)}{\tilde{
\mathcal{E}}^{n}_{k} \left(\chi_{\tilde{\mathbb{A}}_{\mathcal{N},
2 \tilde{\eta}}}(\vec{x}) \me^{\mathcal{N}^{-(m-1)} \sum_{
i_{1},i_{2},\dotsc,i_{m}} \phi_{0}(x_{i_{1}},x_{i_{2}},\dotsc,
x_{i_{m}})} \right)} \underset{\underset{z_{o}=1+o(1)}{
\mathscr{N},n \to \infty}}{\leqslant} \mathfrak{c}_{18}
(n,k,z_{o}) \me^{-\frac{n}{\mathcal{N}} \mathcal{N}
(\mathcal{N}-1) \tilde{\eta}(1+o(1))},
\end{equation*}
where $\mathfrak{c}_{18}(n,k,z_{o}) \! =_{\underset{z_{o}=1+o(1)}{
\mathscr{N},n \to \infty}} \! \mathcal{O}(1)$; hence, via the expansion 
$\ln (1 \! - \! \lozenge) \! =_{\lvert \lozenge \rvert \to 0} \! -\sum_{m
=1}^{\infty} \tfrac{\lozenge^{m}}{m}$, one arrives at, for $n \! \in \! 
\mathbb{N}$ and $k \! \in \! \lbrace 1,2,\dotsc,K \rbrace$ such 
that $\alpha_{p_{\mathfrak{s}}} \! := \! \alpha_{k} \! \neq \! 
\infty$,\footnote{For an expression of the type $\lim_{\scriptscriptstyle 
\underset{z_{o}=1+o(1)}{\mathscr{N},n \to \infty}}f_{1}(n,k,z_{o}) \! 
=_{\scriptscriptstyle \underset{z_{o}=1+o(1)}{\mathscr{N},n \to \infty}} 
\! f_{2}(n,k,z_{o})$, it is meant that $\lvert (\lim_{\scriptscriptstyle 
\underset{z_{o}=1+o(1)}{\mathscr{N},n \to \infty}}f_{1}(n))f_{2}^{-1}(n) 
\rvert \! \leqslant \! \varepsilon_{\diamond}(n,k,z_{o})$, where 
$\varepsilon_{\diamond}(n,k,z_{o}) \! =_{\scriptscriptstyle \underset{z_{o}
=1+o(1)}{\mathscr{N},n \to \infty}} \! 1 \! + \! o(1)$ (the $o(1)$ term 
appearing in the latter asymptotic expansion is not the same as the $o(1)$ 
term in the double-scaling limit $\mathscr{N},n \! \to \! \infty$ such that 
$z_{o} \! = \! 1 \! + \! o(1))$. This notation applies, \emph{mutatis mutandis}, 
for asymptotic expressions of the type $\limsup_{\scriptscriptstyle 
\underset{z_{o}=1+o(1)}{\mathscr{N},n \to \infty}}f_{1}(n,k,z_{o}) \! 
=_{\scriptscriptstyle \underset{z_{o}=1+o(1)}{\mathscr{N},n \to \infty}} \! 
\limsup_{\scriptscriptstyle \underset{z_{o}=1+o(1)}{\mathscr{N},n \to \infty}}
f_{2}(n,k,z_{o})$, $\liminf_{\scriptscriptstyle \underset{z_{o}=1+o(1)}{\mathscr{N},
n \to \infty}}f_{1}(n,k,z_{o}) \! =_{\scriptscriptstyle \underset{z_{o}=1+o(1)}{
\mathscr{N},n \to \infty}} \! \liminf_{\scriptscriptstyle \underset{z_{o}=1+
o(1)}{\mathscr{N},n \to \infty}}f_{2}(n,k,z_{o})$, etc., which appear in the 
proof of this Lemma~\ref{lemrootz} (see, also, Corollary~\ref{corol3.1} 
and Lemma~\ref{lemetatomu} below).}
\begin{align*}
\limsup_{\underset{z_{o}=1+o(1)}{\mathscr{N},n \to \infty}} 
\dfrac{1}{\mathcal{N}} \ln \tilde{\mathcal{E}}^{n}_{k} \left(
\me^{\mathcal{N}^{-(m-1)} \sum_{i_{1},i_{2},\dotsc,i_{m}} 
\phi_{0}(x_{i_{1}},x_{i_{2}},\dotsc,x_{i_{m}})} \right) 
\underset{\underset{z_{o}=1+o(1)}{\mathscr{N},n \to 
\infty}}{\leqslant}& \, 
\limsup_{\underset{z_{o}=1+o(1)}{\mathscr{N},n \to \infty}} 
\dfrac{1}{\mathcal{N}} \ln \tilde{\mathcal{E}}^{n}_{k} \left(
\chi_{\tilde{\mathbb{A}}_{\mathcal{N},2 \tilde{\eta}}}(\vec{x}) 
\me^{\mathcal{N}^{-(m-1)} \sum_{i_{1},i_{2},\dotsc,i_{m}} 
\phi_{0}(x_{i_{1}},x_{i_{2}},\dotsc,x_{i_{m}})} \right) \\
+& \, \mathcal{O} \left(\dfrac{\mathfrak{c}_{19}(n,k,z_{o})}{
\mathcal{N}} \me^{-\frac{n}{\mathcal{N}} \mathcal{N}
(\mathcal{N}-1) \tilde{\eta}(1+o(1))} \right),
\end{align*}
where $\mathfrak{c}_{19}(n,k,z_{o}) \! =_{\underset{z_{o}
=1+o(1)}{\mathscr{N},n \to \infty}} \! \mathcal{O}(1)$. 
Since, for $(x_{1},x_{2},\dotsc,x_{\mathcal{N}}) \! \in \! 
\tilde{\mathbb{A}}_{\mathcal{N},2 \tilde{\eta}}$, 
$E_{\widetilde{V}}^{f} \! + \! 2 \tilde{\eta} \! \geqslant \! 
(\mathcal{N}(\mathcal{N} \! - \! 1))^{-1} \mathscr{K}_{
\mathcal{N}}^{\widetilde{V},f}(x_{1},x_{2},\dotsc,
x_{\mathcal{N}})$, it follows {}from the proofs of 
Lemmata~\ref{lem3.1} and~\ref{lem3.5}, 
Equation~\eqref{eqKvinf4}, and $\sum_{\underset{j \neq 
m}{j=1}}^{\mathcal{N}}1 \! = \! \mathcal{N} \! - \! 1$, 
$m \! = \! 1,2,\dotsc,\mathcal{N}$, that, for $n \! \in \! 
\mathbb{N}$ and $k \! \in \! \lbrace 1,2,\dotsc,K \rbrace$ 
such that $\alpha_{p_{\mathfrak{s}}} \! := \! \alpha_{k} \! 
\neq \! \infty$,
\begin{align*}
E_{\widetilde{V}}^{f} \! + \! 2 \tilde{\eta} \geqslant& \, 
\dfrac{1}{\mathcal{N}(\mathcal{N} \! - \! 1)} \mathscr{K}_{
\mathcal{N}}^{\widetilde{V},f}(x_{1},x_{2},\dotsc,x_{\mathcal{N}}) 
\! \geqslant \! \dfrac{1}{\mathcal{N}} \sum_{m=1}^{\mathcal{N}} 
\hat{\psi}_{\widetilde{V}}^{f}(x_{m}) \! \geqslant \! \dfrac{1}{
\mathcal{N}} \left(- \lvert \hat{C}_{\widetilde{V}}^{f} \rvert 
(\mathcal{N} \! - \! 1) \! + \! \widetilde{V}(x_{i}) \! - \! \left(
\dfrac{\varkappa_{nk \tilde{k}_{\mathfrak{s}-1}}^{\infty} \! + \! 1}{n} 
\right) \ln (1 \! + \! x_{i}^{2}) \right. \\
-&\left. \, \left(\dfrac{\varkappa_{nk} \! - \! 1}{n} \right) \ln (1 \! + 
\! (x_{i} \! - \! \alpha_{k})^{-2}) \! - \! \sum_{q=1}^{\mathfrak{s}-2} 
\dfrac{\varkappa_{nk \tilde{k}_{q}}}{n} \ln (1 \! + \! (x_{i} \! - \! 
\alpha_{p_{q}})^{-2}) \right), \quad i \! = \! 1,2,\dotsc,\mathcal{N},
\end{align*}
which shows that, for $i \! = \! 1,2,\dotsc,\mathcal{N}$, 
$\xi_{i}$ lies in the bounded set $\mathbb{R} \setminus 
\mathfrak{D}_{M_{f}}$ (in particular, $\xi_{i} \! \notin \! 
\lbrace \alpha_{p_{1}},\dotsc,\alpha_{p_{\mathfrak{s}-2}},
\alpha_{p_{\mathfrak{s}}} \rbrace$ (of course, $\xi_{i} \! 
\neq \! \alpha_{p_{\mathfrak{s}-1}} \! = \! \infty))$, that 
is, $\tilde{\mathbb{A}}_{\mathcal{N},2 \tilde{\eta}}$ is a 
bounded $\mathcal{N}$-dimensional proper subset of 
$\mathbb{R}^{\mathcal{N}}$; furthermore, since, with respect 
to any Euclidean-equivalent metric, $\mathbb{R}^{\mathcal{N}} 
\setminus \tilde{\mathbb{A}}_{\mathcal{N},2 \tilde{\eta}}$ is an 
open $\mathcal{N}$-dimensional proper subset of $\mathbb{R}^{
\mathcal{N}}$, it follows that $\tilde{\mathbb{A}}_{\mathcal{N},
2 \tilde{\eta}}$ is a closed $\mathcal{N}$-dimensional proper 
subset of $\mathbb{R}^{\mathcal{N}}$; hence, $\tilde{\mathbb{A}}_{
\mathcal{N},2 \tilde{\eta}}$ is a compact $\mathcal{N}$-dimensional 
set. Since $\phi_{0} \colon \mathbb{R}^{m} \! \to \! \mathbb{R}$, 
$m \! \in \! \mathbb{N}$, is bounded and continuous, it follows, via 
the compactness of $\tilde{\mathbb{A}}_{\mathcal{N},2 \tilde{\eta}}$, 
that there exist points $\tilde{\mathbf{x}}^{\sharp} \! := \! 
(x_{1}^{\sharp},x_{2}^{\sharp},\dotsc,x_{\mathcal{N}}^{\sharp}) 
\! \in \! \tilde{\mathbb{A}}_{\mathcal{N},2 \tilde{\eta}}$ and 
$\tilde{\mathbf{x}}^{\flat} \! := \! (x_{1}^{\flat},x_{2}^{\flat},\dotsc,
x_{\mathcal{N}}^{\flat}) \! \in \! \tilde{\mathbb{A}}_{\mathcal{N},
2 \tilde{\eta}}$ such that
\begin{align*}
\sum_{i_{1},i_{2},\dotsc,i_{m}} \phi_{0}(x_{i_{1}}^{\flat},
x_{i_{2}}^{\flat},\dotsc,x_{i_{m}}^{\flat}) =& \, \inf_{(x_{1},
x_{2},\dotsc,x_{\mathcal{N}}) \in \tilde{\mathbb{A}}_{
\mathcal{N},2 \tilde{\eta}}} \sum_{i_{1},i_{2},\dotsc,i_{m}} 
\phi_{0}(x_{i_{1}},x_{i_{2}},\dotsc,x_{i_{m}}) \! \leqslant \! 
\sum_{i_{1},i_{2},\dotsc,i_{m}} \phi_{0}(x_{i_{1}},x_{i_{2}},
\dotsc,x_{i_{m}}) \\
\leqslant& \, \sup_{(x_{1},x_{2},\dotsc,x_{\mathcal{N}}) \in 
\tilde{\mathbb{A}}_{\mathcal{N},2 \tilde{\eta}}} \sum_{i_{1},
i_{2},\dotsc,i_{m}} \phi_{0}(x_{i_{1}},x_{i_{2}},\dotsc,x_{i_{m}}) 
\! = \! \sum_{i_{1},i_{2},\dotsc,i_{m}} \phi_{0}(x_{i_{1}}^{\sharp},
x_{i_{2}}^{\sharp},\dotsc,x_{i_{m}}^{\sharp});
\end{align*}
hence, via the monotonicity of $\exp (\pmb{\cdot})$, the expansion 
$\ln (1 \! - \! \lozenge) \! =_{\lvert \lozenge \rvert \to 0} \! 
-\sum_{m=1}^{\infty} \tfrac{\lozenge^{m}}{m}$, the linearity of 
$\tilde{\mathcal{E}}^{n}_{k}$, the fact that $\tilde{\mathcal{E}}^{n}_{k}
(1) \! = \! 1$, the relation $\tilde{\mathcal{E}}^{n}_{k}(\chi_{
\tilde{\mathbb{A}}_{\mathcal{N},2 \tilde{\eta}}}(\vec{x})) \! 
= \! 1 \! - \! \operatorname{Prob} \, (\mathbb{R}^{\mathcal{N}} 
\setminus \tilde{\mathbb{A}}_{\mathcal{N},2 \tilde{\eta}})$, and 
the Estimate~\eqref{eqlmrtz7}, it follows that, for $n \! \in \! 
\mathbb{N}$ and $k \! \in \! \lbrace 1,2,\dotsc,K \rbrace$ such 
that $\alpha_{p_{\mathfrak{s}}} \! := \! \alpha_{k} \! \neq \! \infty$,
\begin{align*}
\limsup_{\underset{z_{o}=1+o(1)}{\mathscr{N},n \to \infty}} 
\dfrac{1}{\mathcal{N}} \ln \tilde{\mathcal{E}}^{n}_{k} \left(
\me^{\mathcal{N}^{-(m-1)} \sum_{i_{1},i_{2},\dotsc,i_{m}} 
\phi_{0}(x_{i_{1}},x_{i_{2}},\dotsc,x_{i_{m}})} \right) 
\underset{\underset{z_{o}=1+o(1)}{\mathscr{N},n \to 
\infty}}{\leqslant}& \, 
\limsup_{\underset{z_{o}=1+o(1)}{\mathscr{N},n \to \infty}} 
\dfrac{1}{\mathcal{N}^{m}} \sum_{i_{1},i_{2},\dotsc,i_{m}} \phi_{0}
(x_{i_{1}}^{\sharp},x_{i_{2}}^{\sharp},\dotsc,x_{i_{m}}^{\sharp}) \\
+& \, \mathcal{O} \left(\dfrac{\mathfrak{c}_{20}(n,k,z_{o})}{
\mathcal{N}} \me^{-\frac{n}{\mathcal{N}} \mathcal{N}
(\mathcal{N}-1) \tilde{\eta}(1+o(1))} \right),
\end{align*}
where $\mathfrak{c}_{20}(n,k,z_{o}) \! =_{\underset{z_{o}
=1+o(1)}{\mathscr{N},n \to \infty}} \! \mathcal{O}(1)$; in 
fact, an analogous argument shows that, for $n \! \in \! 
\mathbb{N}$ and $k \! \in \! \lbrace 1,2,\dotsc,K \rbrace$ 
such that $\alpha_{p_{\mathfrak{s}}} \! := \! \alpha_{k} \! 
\neq \! \infty$,
\begin{align*}
\liminf_{\underset{z_{o}=1+o(1)}{\mathscr{N},n \to \infty}} 
\dfrac{1}{\mathcal{N}} \ln \tilde{\mathcal{E}}^{n}_{k} \left(
\me^{\mathcal{N}^{-(m-1)} \sum_{i_{1},i_{2},\dotsc,i_{m}} 
\phi_{0}(x_{i_{1}},x_{i_{2}},\dotsc,x_{i_{m}})} \right) 
\underset{\underset{z_{o}=1+o(1)}{\mathscr{N},n \to \infty}}{\geqslant}& 
\, \liminf_{\underset{z_{o}=1+o(1)}{\mathscr{N},n \to \infty}} 
\dfrac{1}{\mathcal{N}^{m}} \sum_{i_{1},i_{2},\dotsc,i_{m}} \phi_{0}
(x_{i_{1}}^{\flat},x_{i_{2}}^{\flat},\dotsc,x_{i_{m}}^{\flat}) \\
+& \, \mathcal{O} \left(\dfrac{\mathfrak{c}_{21}(n,k,z_{o})}{
\mathcal{N}} \me^{-\frac{n}{\mathcal{N}} \mathcal{N}
(\mathcal{N}-1) \tilde{\eta}(1+o(1))} \right),
\end{align*}
where $\mathfrak{c}_{21}(n,k,z_{o}) \! =_{\underset{z_{o}
=1+o(1)}{\mathscr{N},n \to \infty}} \! \mathcal{O}(1)$. For 
$n \! \in \! \mathbb{N}$ and $k \! \in \! \lbrace 1,2,\dotsc,
K \rbrace$ such that $\alpha_{p_{\mathfrak{s}}} \! := \! 
\alpha_{k} \! \neq \! \infty$, let $\md \nu_{\mathcal{N}}^{r}
(\xi_{j}) \! := \! \mathcal{N}^{-1} \sum_{i_{j}=1}^{\mathcal{
N}} \delta (\xi_{j} \linebreak[4] 
\! - \! x_{i_{j}}^{r}) \, \md \xi_{j}$, $r \! \in \! \lbrace \sharp,\flat 
\rbrace$, $j \! = \! 1,2,\dotsc,m$, where $\delta (\xi_{j} \! - \! 
x_{i_{j}}^{r})$ is the Dirac delta (atomic) mass concentrated at 
$x_{i_{j}}^{r}$: a calculation shows that
\begin{align*}
\limsup_{\underset{z_{o}=1+o(1)}{\mathscr{N},n \to \infty}} 
\dfrac{1}{\mathcal{N}} \ln \tilde{\mathcal{E}}^{n}_{k} 
\left(\me^{\mathcal{N}^{-(m-1)} \sum_{i_{1},i_{2},\dotsc,
i_{m}} \phi_{0}(x_{i_{1}},x_{i_{2}},\dotsc,x_{i_{m}})} \right) 
\underset{\underset{z_{o}=1+o(1)}{\mathscr{N},n \to \infty}}{
\leqslant}& \, \limsup_{\underset{z_{o}=1+o(1)}{\mathscr{N},n \to \infty}} 
\idotsint\limits_{\mathbb{R}^{m}} \phi_{0}(\xi_{1},\xi_{2},\dotsc,\xi_{m}) 
\, \md \nu_{\mathcal{N}}^{\sharp}(\xi_{1}) \, \md \nu_{\mathcal{N}}^{\sharp}
(\xi_{2}) \, \dotsb \, \md \nu_{\mathcal{N}}^{\sharp}(\xi_{m}) \\
+& \, \mathcal{O} \left(\dfrac{\mathfrak{c}_{20}(n,k,z_{o})}{
\mathcal{N}} \me^{-\frac{n}{\mathcal{N}} \mathcal{N}
(\mathcal{N}-1) \tilde{\eta}(1+o(1))} \right), \\
\liminf_{\underset{z_{o}=1+o(1)}{\mathscr{N},n \to \infty}} 
\dfrac{1}{\mathcal{N}} \ln \tilde{\mathcal{E}}^{n}_{k} 
\left(\me^{\mathcal{N}^{-(m-1)} \sum_{i_{1},i_{2},\dotsc,
i_{m}} \phi_{0}(x_{i_{1}},x_{i_{2}},\dotsc,x_{i_{m}})} \right) 
\underset{\underset{z_{o}=1+o(1)}{\mathscr{N},n \to \infty}}{
\geqslant}& \, \liminf_{\underset{z_{o}=1+o(1)}{\mathscr{N},n \to \infty}} 
\idotsint\limits_{\mathbb{R}^{m}} \phi_{0}(\xi_{1},\xi_{2},\dotsc,\xi_{m}) 
\, \md \nu_{\mathcal{N}}^{\flat}(\xi_{1}) \, \md \nu_{\mathcal{N}}^{\flat}
(\xi_{2}) \, \dotsb \, \md \nu_{\mathcal{N}}^{\flat}(\xi_{m}) \\
+& \, \mathcal{O} \left(\dfrac{\mathfrak{c}_{21}(n,k,z_{o})}{
\mathcal{N}} \me^{-\frac{n}{\mathcal{N}} \mathcal{N}
(\mathcal{N}-1) \tilde{\eta}(1+o(1))} \right).
\end{align*}
Proceeding, \emph{mutatis mutandis}, as in the proof of 
Lemma~\ref{lem3.5}, that is, using tightness arguments 
for $\nu_{\mathcal{N}}^{r}$, $r \! \in \! \lbrace \sharp,
\flat \rbrace$, to extract, via a Helly Selection Theorem, 
subsequences of weakly convergent (in the weak-$\ast$ 
topology of measures) probability measures (cf. the 
calculations of Lemma~\ref{lem3.5} leading to 
Equation~\eqref{eql3.5g}, and the analogue of the 
discussion thereafter), one shows that, for $n \! \in \! 
\mathbb{N}$ and $k \! \in \! \lbrace 1,2,\dotsc,K \rbrace$ 
such that $\alpha_{p_{\mathfrak{s}}} \! := \! \alpha_{k} 
\! \neq \! \infty$, via the unicity of the corresponding 
equilibrium measure (cf. Lemma~\ref{lem3.3}) 
$\mu_{\widetilde{V}}^{f} \! \in \! \mathscr{M}_{1}
(\mathbb{R})$, and the latter two asymptotic inequalities, 
in the double-scaling limit $\mathscr{N},n \! \to \! \infty$ 
such that $z_{o} \! = \! 1 \! + \! o(1)$,
\begin{align*}
\limsup_{\underset{z_{o}=1+o(1)}{\mathscr{N},n \to \infty}} 
\dfrac{1}{\mathcal{N}} \ln \tilde{\mathcal{E}}^{n}_{k} 
\left(\me^{\mathcal{N}^{-(m-1)} \sum_{i_{1},i_{2},\dotsc,
i_{m}} \phi_{0}(x_{i_{1}},x_{i_{2}},\dotsc,x_{i_{m}})} \right) 
\underset{\underset{z_{o}=1+o(1)}{\mathscr{N},n \to 
\infty}}{\leqslant}& \, \idotsint\limits_{\mathbb{R}^{m}} \phi_{0}
(\xi_{1},\xi_{2},\dotsc,\xi_{m}) \, \md \mu_{\widetilde{V}}^{f}
(\xi_{1}) \, \md \mu_{\widetilde{V}}^{f}(\xi_{2}) \, \dotsb 
\, \md \mu_{\widetilde{V}}^{f}(\xi_{m}) \! + \! o(1), \\
\liminf_{\underset{z_{o}=1+o(1)}{\mathscr{N},n \to \infty}} 
\dfrac{1}{\mathcal{N}} \ln \tilde{\mathcal{E}}^{n}_{k} 
\left(\me^{\mathcal{N}^{-(m-1)} \sum_{i_{1},i_{2},\dotsc,
i_{m}} \phi_{0}(x_{i_{1}},x_{i_{2}},\dotsc,x_{i_{m}})} \right) 
\underset{\underset{z_{o}=1+o(1)}{\mathscr{N},n \to 
\infty}}{\geqslant}& \, \idotsint\limits_{\mathbb{R}^{m}} \phi_{0}
(\xi_{1},\xi_{2},\dotsc,\xi_{m}) \, \md \mu_{\widetilde{V}}^{f}
(\xi_{1}) \, \md \mu_{\widetilde{V}}^{f}(\xi_{2}) \, \dotsb 
\, \md \mu_{\widetilde{V}}^{f}(\xi_{m}) \! + \! o(1),
\end{align*}
whence
\begin{equation} \label{eqlmrtz8} 
\lim_{\underset{z_{o}=1+o(1)}{\mathscr{N},n \to \infty}} 
\dfrac{1}{\mathcal{N}} \ln \tilde{\mathcal{E}}^{n}_{k} 
\left(\me^{\mathcal{N}^{-(m-1)} \sum_{i_{1},i_{2},\dotsc,
i_{m}} \phi_{0}(x_{i_{1}},x_{i_{2}},\dotsc,x_{i_{m}})} \right) 
\underset{\underset{z_{o}=1+o(1)}{\mathscr{N},
n \to \infty}}{=} \idotsint\limits_{\mathbb{R}^{m}} \phi_{0}
(\xi_{1},\xi_{2},\dotsc,\xi_{m}) \, \md \mu_{\widetilde{V}}^{f}
(\xi_{1}) \, \md \mu_{\widetilde{V}}^{f}(\xi_{2}) \, 
\dotsb \, \md \mu_{\widetilde{V}}^{f}(\xi_{m}).
\end{equation}
Recalling that, for $n \! \in \! \mathbb{N}$ and $k \! \in 
\! \lbrace 1,2,\dotsc,K \rbrace$ such that $\alpha_{p_{
\mathfrak{s}}} \! := \! \alpha_{k} \! \neq \! \infty$, 
$\phi_{0} \colon \mathbb{R}^{m} \! \to \! \mathbb{R}$, 
$m \! \in \! \mathbb{N}$, is bounded and continuous, it 
follows via the argument on p.~162 of \cite{a51} that 
$\tilde{f}(t) \! := \! \ln \tilde{\mathcal{E}}^{n}_{k} (\exp (t 
\mathcal{N}^{-(m-1)} \sum_{i_{1},i_{2},\dotsc,i_{m}} \phi_{0}
(x_{i_{1}},x_{i_{2}},\dotsc.x_{i_{m}})))$ is a convex function 
of the real variable $t$, that is, $\tilde{f}^{\prime \prime}
(t) \! > \! 0$, where the prime denotes differentiation with 
respect to $t$; in particular, since $\tilde{f}(0) \! = \! \ln 
\tilde{\mathcal{E}}^{n}_{k}(1) \! = \! \ln 1 \! = \! 0$, the 
convexity of $\tilde{f}$ implies that $-\tilde{f}(-1) \! 
\leqslant \! \tilde{f}^{\prime}(0) \! \leqslant \! \tilde{f}(1)$, 
whence, for $m \! \in \! \mathbb{N}$,
\begin{align*}
-\dfrac{1}{\mathcal{N}} \ln \tilde{\mathcal{E}}^{n}_{k} \left(
\me^{-\mathcal{N}^{-(m-1)} \sum_{i_{1},i_{2},\dotsc,i_{m}} 
\phi_{0}(x_{i_{1}},x_{i_{2}},\dotsc,x_{i_{m}})} \right) \leqslant& 
\, \dfrac{1}{\mathcal{N}^{m}} \tilde{\mathcal{E}}^{n}_{k} 
\left(\sum_{i_{1},i_{2},\dotsc,i_{m}} \phi_{0}(x_{i_{1}},x_{i_{2}},
\dotsc,x_{i_{m}}) \right) \\
\leqslant& \, \dfrac{1}{\mathcal{N}} \ln \tilde{\mathcal{E}}^{n}_{k} 
\left(\me^{\mathcal{N}^{-(m-1)} \sum_{i_{1},i_{2},\dotsc,i_{m}} 
\phi_{0}(x_{i_{1}},x_{i_{2}},\dotsc,x_{i_{m}})} \right);
\end{align*}
hence, via the two estimates leading to Equation~\eqref{eqlmrtz8} 
and the latter estimate,
\begin{align*}
&\idotsint\limits_{\mathbb{R}^{m}} \phi_{0}(\xi_{1},\xi_{2},
\dotsc,\xi_{m}) \, \md \mu_{\widetilde{V}}^{f}(\xi_{1}) \, 
\md \mu_{\widetilde{V}}^{f}(\xi_{2}) \, \dotsb \, \md 
\mu_{\widetilde{V}}^{f}(\xi_{m}) \! + \! o(1) \underset{
\underset{z_{o}=1+o(1)}{\mathscr{N},n \to \infty}}{\leqslant} 
\lim_{\underset{z_{o}=1+o(1)}{\mathscr{N},n \to \infty}} 
\dfrac{1}{\mathcal{N}^{m}} \tilde{\mathcal{E}}^{n}_{k} 
\left(\sum_{i_{1},i_{2},\dotsc,i_{m}} \phi_{0}(x_{i_{1}},x_{i_{2}},
\dotsc,x_{i_{m}}) \right) \\
& \, \underset{\underset{z_{o}=1+o(1)}{\mathscr{N},n 
\to \infty}}{\leqslant} \idotsint\limits_{\mathbb{R}^{m}} \phi_{0}
(\xi_{1},\xi_{2},\dotsc,\xi_{m}) \, \md \mu_{\widetilde{V}}^{f}
(\xi_{1}) \, \md \mu_{\widetilde{V}}^{f}(\xi_{2}) \, \dotsb 
\, \md \mu_{\widetilde{V}}^{f}(\xi_{m}) \! + \! o(1),
\end{align*}
that is, for $n \! \in \! \mathbb{N}$ and $k \! \in \! \lbrace 
1,2,\dotsc,K \rbrace$ such that $\alpha_{p_{\mathfrak{s}}} 
\! := \! \alpha_{k} \! \neq \! \infty$,
\begin{equation} \label{eqlmrtz9} 
\lim_{\underset{z_{o}=1+o(1)}{\mathscr{N},n \to \infty}} 
\dfrac{1}{\mathcal{N}^{m}} \tilde{\mathcal{E}}^{n}_{k} \left(
\sum_{i_{1},i_{2},\dotsc,i_{m}} \phi_{0}(x_{i_{1}},x_{i_{2}},\dotsc,
x_{i_{m}}) \right) \underset{\underset{z_{o}=1+o(1)}{\mathscr{N},n 
\to \infty}}{=} \idotsint\limits_{\mathbb{R}^{m}} \phi_{0}(\xi_{1},
\xi_{2},\dotsc,\xi_{m}) \, \md \mu_{\widetilde{V}}^{f}(\xi_{1}) 
\, \md \mu_{\widetilde{V}}^{f}(\xi_{2}) \, \dotsb \, \md 
\mu_{\widetilde{V}}^{f}(\xi_{m}), \quad m \! \in \! \mathbb{N}.
\end{equation}
For $n \! \in \! \mathbb{N}$ and $k \! \in \! \lbrace 1,2,
\dotsc,K \rbrace$ such that $\alpha_{p_{\mathfrak{s}}} \! 
:= \! \alpha_{k} \! \neq \! \infty$, write, via the linearity of 
$\tilde{\mathcal{E}}^{n}_{k}$, $\mathcal{N}^{-m} 
\tilde{\mathcal{E}}^{n}_{k}(\sum_{i_{1},i_{2},\dotsc,i_{m}} 
\phi_{0}(x_{i_{1}},x_{i_{2}},\dotsc,\linebreak[4] 
x_{i_{m}})) \! = \! \mathcal{N}^{-m} \tilde{\mathcal{E}}^{n}_{k}
(\sum_{i_{1} \neq i_{2} \neq \dotsb \neq i_{m}} \phi_{0}(x_{i_{1}},
x_{i_{2}},\dotsc,x_{i_{m}})) \! + \! \mathcal{N}^{-m} \tilde{\mathcal{
E}}^{n}_{k}(\sum_{\underset{j \neq m}{i_{j}=i_{m}}} \phi_{0}
(x_{i_{1}},x_{i_{2}},\dotsc,x_{i_{m}}))$, $m \! \in \! \mathbb{N}$: 
an application of the Fundamental Counting Principle shows 
that the number of terms in the multi-dimensional sums 
$\sum_{i_{1},i_{2},\dotsc,i_{m}} \phi_{0}(x_{i_{1}},x_{i_{2}},
\dotsc,\linebreak[4] 
x_{i_{m}})$ and $\sum_{i_{1} \neq i_{2} \neq \dotsb 
\neq i_{m}} \phi_{0}(x_{i_{1}},x_{i_{2}},\dotsc,x_{i_{m}})$ are 
$\mathcal{N}^{m}$ $(= \! \mathcal{O}(\mathcal{N}^{m})$ as 
$\mathcal{N} \! \to \! \infty)$ and $\mathcal{N}(\mathcal{N} 
\! - \! 1) \dotsb (\mathcal{N} \! - \! (m \! - \! 1)) \! = \! 
\mathcal{N}!/(\mathcal{N} \! - \! m)!$ $(= \! \mathcal{O}
(\mathcal{N}^{m})$ as $\mathcal{N} \! \to \! \infty)$, 
respectively; hence, the number of terms in the multi-dimensional 
sum $\sum_{\underset{j \neq m}{i_{j}=i_{m}}} \phi_{0}(x_{i_{1}},
x_{i_{2}},\dotsc,x_{i_{m}})$ is $\mathcal{O}(\mathcal{N}^{m-1})$ 
as $\mathcal{N} \! \to \! \infty$, whence, via the boundedness 
and continuity of $\phi_{0} \colon \mathbb{R}^{m} \! \to \! 
\mathbb{R}$, $m \! \in \! \mathbb{N}$,
\begin{equation*}
\lim_{\underset{z_{o}=1+o(1)}{\mathscr{N},n \to \infty}} 
\dfrac{1}{\mathcal{N}^{m}} \tilde{\mathcal{E}}^{n}_{k} 
\left(\sum_{i_{1},i_{2},\dotsc,i_{m}} \phi_{0}(x_{i_{1}},x_{i_{2}},
\dotsc,x_{i_{m}}) \right) \underset{\underset{z_{o}=1+o(1)}{
\mathscr{N},n \to \infty}}{=} \lim_{\underset{z_{o}=1+o(1)}{
\mathscr{N},n \to \infty}} \dfrac{1}{\mathcal{N}^{m}} 
\tilde{\mathcal{E}}^{n}_{k} \left(\sum_{i_{1} \neq i_{2} \neq 
\dotsb \neq i_{m}} \phi_{0}(x_{i_{1}},x_{i_{2}},\dotsc,x_{i_{m}}) 
\right) \! + \! \mathcal{O} \left(\dfrac{\mathfrak{c}_{22}
(n,k,z_{o})}{\mathcal{N}} \right),
\end{equation*}
where $\mathfrak{c}_{22}(n,k,z_{o}) \! =_{\underset{z_{o}=
1+o(1)}{\mathscr{N},n \to \infty}} \! \mathcal{O}(1)$. Thus, 
via a symmetry argument (see, for example, Chapter~5 of 
\cite{a51}) and Equation~\eqref{eqlmrtz9}, one arrives at, 
for $n \! \in \! \mathbb{N}$ and $k \! \in \! \lbrace 1,2,
\dotsc,K \rbrace$ such that $\alpha_{p_{\mathfrak{s}}} \! 
:= \! \alpha_{k} \! \neq \! \infty$,
\begin{align} \label{eqlmrtz10} 
\lim_{\underset{z_{o}=1+o(1)}{\mathscr{N},n \to \infty}} 
\dfrac{1}{\mathcal{N}^{m}} \tilde{\mathcal{E}}^{n}_{k} 
\left(\sum_{i_{1},i_{2},\dotsc,i_{m}} \phi_{0}(x_{i_{1}},x_{i_{2}},
\dotsc,x_{i_{m}}) \right) \underset{\underset{z_{o}=1+
o(1)}{\mathscr{N},n \to \infty}}{=}& \, \lim_{\underset{z_{o}=
1+o(1)}{\mathscr{N},n \to \infty}} \dfrac{1}{\mathcal{N}^{m}} 
\idotsint\limits_{\mathbb{R}^{m}} \phi_{0}(\xi_{1},\xi_{2},\dotsc,
\xi_{m}) \tilde{\mathfrak{R}}^{n,k}_{m}(\xi_{1},\xi_{2},\dotsc,
\xi_{m}) \, \md \xi_{1} \, \md \xi_{2} \, \dotsb \, \md \xi_{m} 
\nonumber \\
\underset{\underset{z_{o}=1+o(1)}{\mathscr{N},n \to 
\infty}}{=}& \, \idotsint\limits_{\mathbb{R}^{m}} \phi_{0}
(\xi_{1},\xi_{2},\dotsc,\xi_{m}) \, \md \mu_{\widetilde{V}}^{f}
(\xi_{1}) \, \md \mu_{\widetilde{V}}^{f}(\xi_{2}) \, \dotsb 
\, \md \mu_{\widetilde{V}}^{f}(\xi_{m}), \quad m \! \in \! 
\mathbb{N},
\end{align}
where
\begin{equation} \label{eqlmrtz11} 
\tilde{\mathfrak{R}}^{n,k}_{m}(x_{1},x_{2},\dotsc,x_{m}) \! := \! 
\dfrac{\mathcal{N}!}{(\mathcal{N} \! - \! m)!} \idotsint\limits_{
\mathbb{R}^{
\raise-0.7ex\hbox{$\scriptscriptstyle \mathcal{N}-m$}}} \tilde{
\mathcal{P}}^{n}_{k}(x_{1},x_{2},\dotsc,x_{m},\xi_{m+1},\dotsc,
\xi_{\mathcal{N}}) \, \md \xi_{m+1} \, \dotsb \, \md \xi_{\mathcal{N}}, 
\quad m \! \in \! \mathbb{N},
\end{equation}
is the $m$-point correlation function (marginal 
density).\footnote{As a corollary, note that, for $m \! = \! 1$, 
Equation~\eqref{eqlmrtz10} reads: $\lim_{\underset{z_{o}=
1+o(1)}{\mathscr{N},n \to \infty}} \int_{\mathbb{R}} \phi_{0}(\tau) 
\mathcal{N}^{-1} \tilde{\mathfrak{R}}^{n,k}_{1}(\tau) \, \md \tau 
\! =_{\underset{z_{o}=1+o(1)}{\mathscr{N},n \to \infty}} \! 
\int_{\mathbb{R}} \phi_{0}(\tau) \, \md \mu_{\widetilde{V}}^{f}
(\tau)$, which implies that, for $n \! \in \! \mathbb{N}$ and 
$k \! \in \! \lbrace 1,2,\dotsc,K \rbrace$ such that $\alpha_{p_{
\mathfrak{s}}} \! := \! \alpha_{k} \! \neq \! \infty$, $\mathcal{N}^{-1} 
\tilde{\mathfrak{R}}^{n,k}_{1}(x) \, \md x \! \overset{\ast}{\to} \! \md 
\mu_{\widetilde{V}}^{f}(x)$ in the double-scaling limit $\mathscr{N},n 
\! \to \! \infty$ such that $z_{o} \! = \! 1 \! + \! o(1)$, that is, 
$\mathcal{N}^{-1} \tilde{\mathfrak{R}}^{n,k}_{1}(x) \, \md x$ 
converges weakly (in the weak-$\ast$ topology of measures) 
to the corresponding equilibrium measure $\mu_{\widetilde{
V}}^{f}$.} Consider, now, for $n \! \in \! \mathbb{N}$ and $k \! \in \! 
\lbrace 1,2,\dotsc,K \rbrace$ such that $\alpha_{p_{\mathfrak{s}}} \! 
:= \! \alpha_{k} \! \neq \! \infty$, the one-point correlation function:
\begin{equation*}
\dfrac{1}{\mathcal{N}} \tilde{\mathfrak{R}}^{n,k}_{1}
(\tau) \! = \! \idotsint\limits_{\mathbb{R}^{
\raise-0.7ex\hbox{$\scriptscriptstyle \mathcal{N}-1$}}} 
\tilde{\mathcal{P}}^{n}_{k}(\tau,\xi_{2},\dotsc,\xi_{\mathcal{N}}) 
\, \md \xi_{2} \, \md \xi_{3} \, \dotsb \, \md \xi_{\mathcal{N}};
\end{equation*}
via this formula and Equation~\eqref{eqlmrtzp}, a calculation shows that
\begin{equation} \label{eqlmrtzr} 
\dfrac{1}{\mathcal{N}} \tilde{\mathfrak{R}}^{n,k}_{1}(\tau) \! = \! 
\dfrac{\tilde{\mathscr{Z}}^{n,\sharp}_{k}}{\tilde{\mathscr{Z}}^{n}_{k}} 
\me^{-n \widetilde{V}(\tau)} \idotsint\limits_{\mathbb{R}^{\mathcal{N}-1}} 
\prod_{j=2}^{\mathcal{N}}(\tau \! - \! \xi_{j})^{2} \left(
\prod_{q=1}^{\mathfrak{s}-2}(\tau \! - \! \alpha_{p_{q}})^{\varkappa_{nk 
\tilde{k}_{q}}}(\tau \! - \! \alpha_{k})^{\varkappa_{nk}-1} \right)^{-2} 
\tilde{\mathcal{P}}^{n,\sharp}_{k}(\xi_{2},\xi_{3},\dotsc,\xi_{\mathcal{N}}) 
\, \md \xi_{2} \, \md \xi_{2} \, \dotsb \, \md \xi_{\mathcal{N}},
\end{equation}
where $\tilde{\mathscr{Z}}^{n}_{k}$ is given by 
Equation~\eqref{eqlmrtzz}, and $\tilde{\mathcal{P}}^{n,
\sharp}_{k}(x_{2},x_{3},\dotsc,x_{\mathcal{N}}) \, \md 
x_{2} \, \md x_{3} \, \dotsb \, \md x_{\mathcal{N}}$ is the 
multi-dimensional probability measure on $\mathbb{R}^{
\mathcal{N}-1}$ with density
\begin{equation} \label{eqlmrtzps} 
\tilde{\mathcal{P}}^{n,\sharp}_{k}(x_{2},x_{3},\dotsc,x_{\mathcal{N}}) 
\! := \! \dfrac{1}{\tilde{\mathscr{Z}}^{n,\sharp}_{k}} \me^{-n 
\sum_{m_{1}=2}^{\mathcal{N}} \widetilde{V}(x_{m_{1}})} \prod_{
\substack{i,j=2\\j<i}}^{\mathcal{N}}(x_{j} \! - \! x_{i})^{2} \left(
\prod_{m=2}^{\mathcal{N}} \prod_{q=1}^{\mathfrak{s}-2}
(x_{m} \! - \! \alpha_{p_{q}})^{\varkappa_{nk \tilde{k}_{q}}}
(x_{m} \! - \! \alpha_{k})^{\varkappa_{nk}-1} \right)^{-2},
\end{equation}
where
\begin{equation} \label{eqlmrtzzs} 
\tilde{\mathscr{Z}}^{n,\sharp}_{k} \! = \! 
\idotsint\limits_{\mathbb{R}^{\mathcal{N}-1}} 
\me^{-n \sum_{m_{1}=2}^{\mathcal{N}} \widetilde{V}(\tau_{m_{1}})} 
\prod_{\substack{i,j=2\\j<i}}^{\mathcal{N}}(\tau_{j} \! - \! \tau_{i})^{2} 
\left(\prod_{m=2}^{\mathcal{N}} \prod_{q=1}^{\mathfrak{s}-2}
(\tau_{m} \! - \! \alpha_{p_{q}})^{\varkappa_{nk \tilde{k}_{q}}}
(\tau_{m} \! - \! \alpha_{k})^{\varkappa_{nk}-1} \right)^{-2} \md 
\tau_{2} \, \md \tau_{3} \, \dotsb \, \md \tau_{\mathcal{N}}.
\end{equation}
(Note that $\idotsint\nolimits_{\mathbb{R}^{\mathcal{N}-1}} 
\tilde{\mathcal{P}}^{n,\sharp}_{k}(\xi_{2},\xi_{3},\dotsc,
\xi_{\mathcal{N}}) \, \md \xi_{2} \, \md \xi_{3} \, \dotsb \, \md 
\xi_{\mathcal{N}} \! = \! 1$.) {}From Equations~\eqref{eqlmrtzz} 
and~\eqref{eqlmrtzzs}, one shows that
\begin{equation} \label{eqlmrtz12} 
\dfrac{\tilde{\mathscr{Z}}^{n}_{k}}{\tilde{\mathscr{Z}}^{n,\sharp}_{k}} 
\! = \! \tilde{\mathcal{E}}^{n,\sharp}_{k} \left(\int_{\mathbb{R}} 
\me^{-n \widetilde{V}(\tau)} \prod_{j=2}^{\mathcal{N}}
(\tau \! - \! \xi_{j})^{2} \left(\prod_{q=1}^{\mathfrak{s}-2}
(\tau \! - \! \alpha_{p_{q}})^{\varkappa_{nk \tilde{k}_{q}}}(\tau \! 
- \! \alpha_{k})^{\varkappa_{nk}-1} \right)^{-2} \md \tau \right),
\end{equation}
where $\tilde{\mathcal{E}}^{n,\sharp}_{k}$ denotes the expectation 
with respect to the multi-dimensional probability measure $\tilde{
\mathcal{P}}^{n,\sharp}_{k}(x_{2},x_{3},\dotsc,x_{\mathcal{N}}) 
\linebreak[4] 
\pmb{\cdot} \, \md x_{2} \, \md x_{3} \, \dotsb \, \md 
x_{\mathcal{N}}$.\footnote{\, $\tilde{\mathcal{E}}^{n,\sharp}_{k}
(F(x_{2},x_{3},\dotsc,x_{\mathcal{N}})) \! := \! \idotsint\nolimits_{
\mathbb{R}^{\mathcal{N}-1}}F(\xi_{2},\xi_{3},\dotsc,\xi_{\mathcal{N}}) 
\tilde{\mathcal{P}}^{n,\sharp}_{k}(\xi_{2},\xi_{3},\dotsc,\xi_{\mathcal{N}}) 
\, \md \xi_{2} \, \md \xi_{3} \, \dotsb \, \md \xi_{\mathcal{N}}$.} 
Via an application of Jensen's Inequality to the argument of 
$\tilde{\mathcal{E}}^{n,\sharp}_{k}$ in Equation~\eqref{eqlmrtz12}, 
one shows that, for $n \! \in \! \mathbb{N}$ and $k \! \in \! \lbrace 
1,2,\dotsc,K \rbrace$ such that $\alpha_{p_{\mathfrak{s}}} \! := \! 
\alpha_{k} \! \neq \! \infty$,
\begin{align} \label{eqlmrtz13} 
&\int_{\mathbb{R}} \me^{-n \widetilde{V}(\tau)} \prod_{j=
2}^{\mathcal{N}}(\tau \! - \! \xi_{j})^{2} \left(\prod_{q=1}^{
\mathfrak{s}-2}(\tau \! - \! \alpha_{p_{q}})^{\varkappa_{nk 
\tilde{k}_{q}}}(\tau \! - \! \alpha_{k})^{\varkappa_{nk}-1} 
\right)^{-2} \md \tau \! \geqslant \! \tilde{\mathfrak{z}} \exp 
\left(\dfrac{1}{\tilde{\mathfrak{z}}} \left(-\dfrac{n}{\mathcal{N}}
(\mathcal{N} \! - \! 1) \int_{\mathbb{R}} \widetilde{V}(\tau) 
\me^{-\frac{n}{\mathcal{N}} \widetilde{V}(\tau)} \, \md \tau 
\right. \right. \nonumber \\
+&\left. \left. \, \int_{\mathbb{R}} \sum_{j=2}^{\mathcal{N}}
2 \ln (\lvert \tau \! - \! \xi_{j} \rvert) \me^{-\frac{n}{\mathcal{N}} 
\widetilde{V}(\tau)} \, \md \tau \! + \! \int_{\mathbb{R}} \left(
\sum_{q=1}^{\mathfrak{s}-2}-2 \varkappa_{nk \tilde{k}_{q}} \ln 
\lvert \tau \! - \! \alpha_{p_{q}} \rvert \! - \! 2(\varkappa_{nk} 
\! - \! 1) \ln \lvert \tau \! - \! \alpha_{k} \rvert \right) 
\me^{-\frac{n}{\mathcal{N}} \widetilde{V}(\tau)} \, \md \tau 
\right) \right),
\end{align}
where
\begin{equation*}
\tilde{\mathfrak{z}} \! = \! \tilde{\mathfrak{z}}(n,k,z_{o}) 
\! := \! \int_{\mathbb{R}} \me^{-\frac{n}{\mathcal{N}} 
\widetilde{V}(\tau)} \, \md \tau.
\end{equation*}
Recalling from the proof of Lemma~\ref{lem3.1} that, for $n \! \in 
\! \mathbb{N}$ and $k \! \in \! \lbrace 1,2,\dotsc,K \rbrace$ such 
that $\alpha_{p_{\mathfrak{s}}} \! := \! \alpha_{k} \! \neq \! \infty$, 
in the double-scaling limit $\mathscr{N},n \! \to \! \infty$ such 
that $z_{o} \! = \! 1 \! + \! o(1)$, via 
conditions~\eqref{eq20}--\eqref{eq22} for regular $\widetilde{V} 
\colon \overline{\mathbb{R}} \setminus \lbrace \alpha_{1},
\alpha_{2},\dotsc,\alpha_{K} \rbrace \! \to \! \mathbb{R}$, there 
exists $\tilde{c}_{\infty} \! = \! \tilde{c}_{\infty}(n,k,z_{o}) \! > \! 
0$ and $\mathcal{O}(1)$ such that, for $x \! \in \! \lbrace \lvert 
x \rvert \! \geqslant \! T_{M_{f}} \rbrace$, $\widetilde{V}(x) \! 
\geqslant \! (1 \! + \! \tilde{c}_{\infty}) \ln (1 \! + \! x^{2})$, and, 
for $q \! = \! 1,\dotsc,\mathfrak{s} \! - \! 2,\mathfrak{s}$, there 
exists $\tilde{c}_{q} \! = \! \tilde{c}_{q}(n,k,z_{o}) \! > \! 0$ and 
$\mathcal{O}(1)$ such that, for $x \! \in \! \mathscr{O}_{\frac{1}{
T_{M_{f}}}}(\alpha_{p_{q}})$, $\widetilde{V}(x) \! \geqslant \! (1 \! 
+ \! \tilde{c}_{q}) \ln (1 \! + \! (x \! - \! \alpha_{p_{q}})^{-2})$, 
and the definition $\mathfrak{D}_{M_{f}} \! := \! \lbrace \lvert 
x \rvert \! \geqslant \! T_{M_{f}} \rbrace \cup \cup_{\underset{q 
\neq \mathfrak{s}-1}{q=1}}^{\mathfrak{s}} \operatorname{clos}
(\mathscr{O}_{\frac{1}{T_{M_{f}}}}(\alpha_{p_{q}}))$, one arrives at, 
via the inequalities $\lvert y_{1} \! - \! y_{2} \rvert^{2} \! \leqslant 
\! (1 \! + \! y_{1}^{2})(1 \! + \! y_{2}^{2})$ and $\ln \lvert y_{1} 
\! - \! y_{2} \rvert^{-1} \! \geqslant \! -\tfrac{1}{2} \ln (1 \! + \! 
y_{1}^{2})(1 \! + \! y_{2}^{2})$, $y_{1},y_{2} \! \in \! \mathbb{R}$, the 
formula $\sum_{j=2}^{\mathcal{N}}1 \! = \! \mathcal{N} \! - \! 1$, 
and standard algebraic and logarithmic inequalities, the following 
estimates: (i)
\begin{equation*}
\tilde{\mathfrak{z}} \! = \! \left(\int_{(\mathbb{R} \setminus 
\mathfrak{D}_{M_{f}}) \setminus J_{f}} \! + \! \int_{J_{f}} 
\! + \! \int_{\lbrace \lvert \tau \rvert \geqslant 
T_{M_{f}} \rbrace} \! + \! \sum_{\substack{q=1\\q \neq 
\mathfrak{s}-1}}^{\mathfrak{s}} \int_{\mathscr{O}_{\frac{1}{
T_{M_{f}}}}(\alpha_{p_{q}})} \right) \me^{-\frac{n}{\mathcal{N}} 
\widetilde{V}(\tau)} \, \md \tau,
\end{equation*}
with
\begin{gather*}
\int_{(\mathbb{R} \setminus \mathfrak{D}_{M_{f}}) \setminus 
J_{f}} \me^{-\frac{n}{\mathcal{N}} \widetilde{V}(\tau)} 
\, \md \tau \underset{\underset{z_{o}=1+o(1)}{\mathscr{N},n 
\to \infty}}{\leqslant} \me^{-\frac{n}{\mathcal{N}} \inf \lbrace 
\mathstrut \widetilde{V}(\tau); \, \tau \in (\mathbb{R} \setminus 
\mathfrak{D}_{M_{f}}) \setminus J_{f} \rbrace} \operatorname{meas}
((\mathbb{R} \setminus \mathfrak{D}_{M_{f}}) \setminus J_{f}) 
\! =: \! \mathfrak{c}_{23}(n,k,z_{o}), \\
\int_{J_{f}} \me^{-\frac{n}{\mathcal{N}} \widetilde{V}(\tau)} 
\, \md \tau \underset{\underset{z_{o}=1+o(1)}{\mathscr{N},
n \to \infty}}{\leqslant} \me^{-\frac{n}{\mathcal{N}} \inf 
\lbrace \mathstrut \widetilde{V}(\tau); \, \tau \in J_{f} \rbrace} 
\sum_{j=1}^{N+1} \lvert \tilde{b}_{j-1} \! - \! \tilde{a}_{j} 
\rvert \! =: \! \mathfrak{c}_{24}(n,k,z_{o}), \\
\int_{\lbrace \lvert \tau \rvert \geqslant T_{M_{f}} \rbrace} 
\me^{-\frac{n}{\mathcal{N}} \widetilde{V}(\tau)} \, \md \tau 
\underset{\underset{z_{o}=1+o(1)}{\mathscr{N},n \to 
\infty}}{\leqslant} \dfrac{2T_{M_{f}}^{1-2 \frac{n}{\mathcal{N}}
(1+\tilde{c}_{\infty})}}{2 \frac{n}{\mathcal{N}}(1 \! + \! 
\tilde{c}_{\infty}) \! - \! 1} \! =: \! \mathfrak{c}_{25}
(n,k,z_{o}), \\
\sum_{\substack{q=1\\q \neq \mathfrak{s}-1}}^{\mathfrak{s}} 
\int_{\mathscr{O}_{\frac{1}{T_{M_{f}}}}(\alpha_{p_{q}})} \me^{-
\frac{n}{\mathcal{N}} \widetilde{V}(\tau)} \, \md \tau 
\underset{\underset{z_{o}=1+o(1)}{\mathscr{N},n \to \infty}}{
\leqslant} \sum_{\substack{q=1\\q \neq \mathfrak{s}-1}}^{
\mathfrak{s}} \dfrac{2T_{M_{f}}^{-\left(1+2 \frac{n}{\mathcal{N}}
(1+\tilde{c}_{q}) \right)}}{2 \frac{n}{\mathcal{N}}(1 \! + \! \tilde{c}_{q}) 
\! + \! 1} \! =: \! \mathfrak{c}_{26}(n,k,z_{o}),
\end{gather*}
where $2 \tfrac{n}{\mathcal{N}}(1 \! + \! \tilde{c}_{\infty}) \! > \! 
1$, and $\mathfrak{c}_{m}(n,k,z_{o}) \! =_{\underset{z_{o}=1+
o(1)}{\mathscr{N},n \to \infty}} \! \mathcal{O}(1)$, $m \! = \! 
23,24,25,26$, whence, gathering the above-derived bounds, it 
follows that $\tilde{\mathfrak{z}} \! =_{\underset{z_{o}=1+o(1)}{
\mathscr{N},n \to \infty}} \! \mathcal{O}(1)$ (and $\neq \! 0)$, 
which establishes the boundedness of $\tilde{\mathfrak{z}}$; (ii)
\begin{equation*}
-\dfrac{1}{\tilde{\mathfrak{z}}} \dfrac{n}{\mathcal{N}}
(\mathcal{N} \! - \! 1) \int_{\mathbb{R}} \widetilde{V}(\tau) 
\me^{-\frac{n}{\mathcal{N}} \widetilde{V}(\tau)} \, \md \tau 
\! = \! -\dfrac{1}{\tilde{\mathfrak{z}}} \dfrac{n}{\mathcal{N}}
(\mathcal{N} \! - \! 1) \left(\int_{(\mathbb{R} \setminus 
\mathfrak{D}_{M_{f}}) \setminus J_{f}} \! + \! \int_{J_{f}} \! + \! 
\int_{\lbrace \lvert \tau \rvert \geqslant T_{M_{f}} \rbrace} \! + \! 
\sum_{\substack{q=1\\q \neq \mathfrak{s}-1}}^{\mathfrak{s}} 
\int_{\mathscr{O}_{\frac{1}{T_{M_{f}}}}(\alpha_{p_{q}})} \right) 
\widetilde{V}(\tau) \me^{-\frac{n}{\mathcal{N}} \widetilde{V}(\tau)} 
\, \md \tau,
\end{equation*}
with
\begin{align*}
-\dfrac{1}{\tilde{\mathfrak{z}}} \dfrac{n}{\mathcal{N}}(\mathcal{N} 
\! - \! 1) \int_{(\mathbb{R} \setminus \mathfrak{D}_{M_{f}}) 
\setminus J_{f}} \widetilde{V}(\tau) \me^{-\frac{n}{\mathcal{N}} 
\widetilde{V}(\tau)} \, \md \tau \underset{\underset{z_{o}=1+
o(1)}{\mathscr{N},n \to \infty}}{\geqslant}& \, -\dfrac{n}{\mathcal{N}}
(\mathcal{N} \! - \! 1) \dfrac{\sup \lbrace \mathstrut \widetilde{V}
(x); \, x \! \in \! (\mathbb{R} \setminus \mathfrak{D}_{M_{f}}) 
\setminus J_{f} \rbrace}{\tilde{\mathfrak{z}}} \\
\times& \, \me^{-\frac{n}{\mathcal{N}} \sup \lbrace \mathstrut 
\widetilde{V}(x); \, x \in (\mathbb{R} \setminus \mathfrak{D}_{M_{f}}) 
\setminus J_{f} \rbrace} \operatorname{meas}((\mathbb{R} \setminus 
\mathfrak{D}_{M_{f}}) \setminus J_{f}) \\=:& \, -\dfrac{n}{\mathcal{N}}
(\mathcal{N} \! - \! 1) \mathfrak{c}_{27}(n,k,z_{o}),
\end{align*}
\begin{align*}
-\dfrac{1}{\tilde{\mathfrak{z}}} \dfrac{n}{\mathcal{N}}(\mathcal{N} 
\! - \! 1) \int_{J_{f}} \widetilde{V}(\tau) \me^{-\frac{n}{\mathcal{N}} 
\widetilde{V}(\tau)} \, \md \tau \underset{\underset{z_{o}=1+
o(1)}{\mathscr{N},n \to \infty}}{\geqslant}& \, 
-\dfrac{n}{\mathcal{N}}(\mathcal{N} \! - \! 1) \dfrac{\sup 
\lbrace \mathstrut \widetilde{V}(x); \, x \! \in \! J_{f} \rbrace}{
\tilde{\mathfrak{z}}} \me^{-\frac{n}{\mathcal{N}} \sup \lbrace 
\mathstrut \widetilde{V}(x); \, x \in J_{f} \rbrace} \sum_{j=1}^{
N+1} \lvert \tilde{b}_{j-1} \! - \! \tilde{a}_{j} \rvert \\
=:& \, -\dfrac{n}{\mathcal{N}}(\mathcal{N} \! - \! 1) 
\mathfrak{c}_{28}(n,k,z_{o}),
\end{align*}
\begin{equation*}
-\dfrac{1}{\tilde{\mathfrak{z}}} \dfrac{n}{\mathcal{N}}(\mathcal{N} 
\! - \! 1) \int_{\lbrace \lvert \tau \rvert \geqslant T_{M_{f}} \rbrace} 
\widetilde{V}(\tau) \me^{-\frac{n}{\mathcal{N}} \widetilde{V}(\tau)} 
\, \md \tau \underset{\underset{z_{o}=1+o(1)}{\mathscr{N},n \to 
\infty}}{\geqslant} -\dfrac{n}{\mathcal{N}}(\mathcal{N} \! - \! 1) 
\dfrac{(1 \! + \! \tilde{c}_{\infty})}{\tilde{\mathfrak{z}}}(\pi \! - \! 2 
\tan^{-1}(T_{M_{f}})) \! =: \! -\dfrac{n}{\mathcal{N}}(\mathcal{N} \! 
- \! 1) \mathfrak{c}_{29}(n,k,z_{o}),
\end{equation*}
\begin{align*}
-\dfrac{1}{\tilde{\mathfrak{z}}} \dfrac{n}{\mathcal{N}}(\mathcal{N} 
\! - \! 1) \sum_{\substack{q=1\\q \neq \mathfrak{s}-1}}^{\mathfrak{s}} 
\int_{\mathscr{O}_{\frac{1}{T_{M_{f}}}}(\alpha_{p_{q}})} \widetilde{V}
(\tau) \me^{-\frac{n}{\mathcal{N}} \widetilde{V}(\tau)} \, \md \tau 
\underset{\underset{z_{o}=1+o(1)}{\mathscr{N},n \to \infty}}{
\geqslant}& \, -\dfrac{n}{\mathcal{N}}(\mathcal{N} \! - \! 1) 
\left(\dfrac{2}{\tilde{\mathfrak{z}}} \dfrac{\mathcal{N}}{n} 
\sum_{\substack{q=1\\q \neq \mathfrak{s}-1}}^{\mathfrak{s}} 
\dfrac{\ln \left(T_{M_{f}}^{2 \frac{n}{\mathcal{N}}(1+\tilde{c}_{q})} 
\right)}{T_{M_{f}}^{1+2\frac{n}{\mathcal{N}}(1+\tilde{c}_{q})}} \right) \\
=:& \, -\dfrac{n}{\mathcal{N}}(\mathcal{N} \! - \! 1) 
\mathfrak{c}_{30}(n,k,z_{o}),
\end{align*}
where $\mathfrak{c}_{m}(n,k,z_{o}) \! =_{\underset{z_{o}=1+o(1)}{
\mathscr{N},n \to \infty}} \! \mathcal{O}(1)$, $m \! = \! 27,28,29,
30$, whence, gathering the above-derived bounds, it follows that 
$-\tilde{\mathfrak{z}}^{-1} \tfrac{n}{\mathcal{N}}(\mathcal{N} 
\! - \! 1) \int_{\mathbb{R}} \widetilde{V}(\tau) \me^{-\frac{n}{
\mathcal{N}} \widetilde{V}(\tau)} \, \md \tau \! \geqslant_{\underset{
z_{o}=1+o(1)}{\mathscr{N},n \to \infty}} \! -\tfrac{n}{\mathcal{N}}
(\mathcal{N} \! - \! 1) \mathfrak{c}_{31}(n,k,z_{o})$, where $\mathfrak{c}_{31}
(n,k,z_{o}) \! =_{\underset{z_{o}=1+o(1)}{\mathscr{N},n \to \infty}} \! 
\mathcal{O}(1)$; (iii)
\begin{equation*}
\dfrac{1}{\tilde{\mathfrak{z}}} \int_{\mathbb{R}} 
\sum_{j=2}^{\mathcal{N}}2 \ln (\lvert \tau \! - \! \xi_{j} 
\rvert) \me^{-\frac{n}{\mathcal{N}} \widetilde{V}(\tau)} \, 
\md \tau \! = \! \dfrac{1}{\tilde{\mathfrak{z}}} \left(
\int_{(\mathbb{R} \setminus \mathfrak{D}_{M_{f}}) 
\setminus J_{f}} \! + \! \int_{J_{f}} \! + \! \int_{\lbrace \lvert 
\tau \rvert \geqslant T_{M_{f}} \rbrace} \! + \! \sum_{
\substack{q=1\\q \neq \mathfrak{s}-1}}^{\mathfrak{s}} 
\int_{\mathscr{O}_{\frac{1}{T_{M_{f}}}}(\alpha_{p_{q}})} \right) 
\sum_{j=2}^{\mathcal{N}}2 \ln (\lvert \tau \! - \! \xi_{j} 
\rvert) \me^{-\frac{n}{\mathcal{N}} \widetilde{V}(\tau)} \, 
\md \tau,
\end{equation*}
with (using that $\sum_{j=2}^{\mathcal{N}} \ln (\xi_{j}^{2} \! 
+ \! 1) \! \geqslant \! \mathfrak{c}_{\ast}^{\natural}(n,k,z_{o})
(\mathcal{N} \! - \! 1)$, where $\mathfrak{c}_{\ast}^{\natural}
(n,k,z_{o}) \! =_{\underset{z_{o}=1+o(1)}{\mathscr{N},n \to \infty}} 
\! \mathcal{O}(1))$
\begin{align*}
\dfrac{1}{\tilde{\mathfrak{z}}} \int_{(\mathbb{R} \setminus 
\mathfrak{D}_{M_{f}}) \setminus J_{f}} &\sum_{j=2}^{
\mathcal{N}}2 \ln (\lvert \tau \! - \! \xi_{j} \rvert) 
\me^{-\frac{n}{\mathcal{N}} \widetilde{V}(\tau)} \, \md \tau 
\underset{\underset{z_{o}=1+o(1)}{\mathscr{N},n \to \infty}}{
\geqslant} -\dfrac{n}{\mathcal{N}}(\mathcal{N} \! - \! 1) 
\left(\mathfrak{c}_{\ast}^{\natural}(n,k,z_{o}) \! + \! \inf 
\lbrace \ln (1 \! + \! x^{2}); \, x \! \in \! (\mathbb{R} 
\setminus \mathfrak{D}_{M_{f}}) \setminus J_{f} \rbrace \right) \\
&\times \, \dfrac{1}{\tilde{\mathfrak{z}}} \dfrac{\mathcal{
N}}{n} \me^{-\frac{n}{\mathcal{N}} \sup \lbrace \mathstrut 
\widetilde{V}(x); \, x \in (\mathbb{R} \setminus 
\mathfrak{D}_{M_{f}}) \setminus J_{f} \rbrace} 
\operatorname{meas}((\mathbb{R} \setminus 
\mathfrak{D}_{M_{f}}) \setminus J_{f}) \! =: \! -\dfrac{n}{
\mathcal{N}}(\mathcal{N} \! - \! 1) \mathfrak{c}_{32}(n,k,z_{o}),
\end{align*}
\begin{align*}
\dfrac{1}{\tilde{\mathfrak{z}}} \int_{J_{f}} &\sum_{j=2}^{\mathcal{N}}
2 \ln (\lvert \tau \! - \! \xi_{j} \rvert) \me^{-\frac{n}{\mathcal{N}} 
\widetilde{V}(\tau)} \, \md \tau \underset{\underset{z_{o}=1+
o(1)}{\mathscr{N},n \to \infty}}{\geqslant} -\dfrac{n}{\mathcal{N}}
(\mathcal{N} \! - \! 1) \left(\mathfrak{c}_{\ast}^{\natural}(n,k,z_{o}) 
\! + \! \inf \lbrace \ln (1 \! + \! x^{2}); \, x \! \in \! J_{f} \rbrace \right) \\
&\times \, \dfrac{1}{\tilde{\mathfrak{z}}} \dfrac{\mathcal{
N}}{n} \me^{-\frac{n}{\mathcal{N}} \sup \lbrace \mathstrut 
\widetilde{V}(x); \, x \in J_{f} \rbrace} \sum_{j=1}^{N+1} 
\lvert \tilde{b}_{j-1} \! - \! \tilde{a}_{j} \rvert \! =: \! 
-\dfrac{n}{\mathcal{N}}(\mathcal{N} \! - \! 1) \mathfrak{c}_{33}
(n,k,z_{o}),
\end{align*}
\begin{align*}
\dfrac{1}{\tilde{\mathfrak{z}}} \int_{\lbrace \lvert \tau \rvert 
\geqslant T_{M_{f}} \rbrace} \sum_{j=2}^{\mathcal{N}}2 \ln 
(\lvert \tau \! - \! \xi_{j} \rvert) \me^{-\frac{n}{\mathcal{N}} 
\widetilde{V}(\tau)} \, \md \tau \underset{\underset{z_{o}=
1+o(1)}{\mathscr{N},n \to \infty}}{\geqslant}& \, -\dfrac{n}{
\mathcal{N}}(\mathcal{N} \! - \! 1) \left(\mathfrak{c}_{\ast}^{\natural}
(n,k,z_{o}) \! + \! \dfrac{\mathcal{N}}{n} \dfrac{1}{\tilde{c}_{
\infty}} \right) \dfrac{1}{\tilde{\mathfrak{z}}} \dfrac{\mathcal{
N}}{n}(\pi \! - \! 2 \tan^{-1}(T_{M_{f}})) \\
=:& \, -\dfrac{n}{\mathcal{N}}(\mathcal{N} \! - \! 1) 
\mathfrak{c}_{34}(n,k,z_{o}),
\end{align*}
\begin{align*}
\dfrac{1}{\tilde{\mathfrak{z}}} \sum_{\substack{q=1\\q \neq 
\mathfrak{s}-1}}^{\mathfrak{s}} \int_{\mathscr{O}_{\frac{1}{
T_{M_{f}}}}(\alpha_{p_{q}})} &\sum_{j=2}^{\mathcal{N}}2 \ln 
(\lvert \tau \! - \! \xi_{j} \rvert) \me^{-\frac{n}{\mathcal{N}} 
\widetilde{V}(\tau)} \, \md \tau \underset{\underset{z_{o}=
1+o(1)}{\mathscr{N},n \to \infty}}{\geqslant} -\dfrac{n}{
\mathcal{N}}(\mathcal{N} \! - \! 1) \left(\mathfrak{c}_{\ast}^{\natural}
(n,k,z_{o}) \sum_{\substack{q=1\\q \neq \mathfrak{s}-1}}^{
\mathfrak{s}} \dfrac{2T_{M_{f}}^{-\left(1+2\frac{n}{\mathcal{N}}
(1+\tilde{c}_{q}) \right)}}{2\frac{n}{\mathcal{N}}(1 \! + \! \tilde{c}_{q}) 
\! + \! 1} \right. \\
&\left. \, +2 \sum_{\substack{q=1\\q \neq \mathfrak{s}-1}}^{
\mathfrak{s}} \dfrac{\ln \left((\alpha_{p_{q}} \! - \! T_{M_{f}}^{-1})^{2} 
\! + \! 1 \right)}{ \left(2\frac{n}{\mathcal{N}}(1 \! + \! \tilde{c}_{q}) \! 
+ \! 1 \right)T_{M_{f}}^{1+2\frac{n}{\mathcal{N}}(1+\tilde{c}_{q})}} 
\right) \dfrac{1}{\tilde{\mathfrak{z}}} \dfrac{\mathcal{N}}{n} \! =: \! 
-\dfrac{n}{\mathcal{N}}(\mathcal{N} \! - \! 1) \mathfrak{c}_{35}
(n,k,z_{o}),
\end{align*}
where $\mathfrak{c}_{m}(n,k,z_{o}) \! =_{\underset{z_{o}=1+
o(1)}{\mathscr{N},n \to \infty}} \! \mathcal{O}(1)$, $m \! = \! 
32,33,34,35$, whence, gathering the above-derived bounds, 
it follows that $\tilde{\mathfrak{z}}^{-1} \int_{\mathbb{R}} 
\sum_{j=2}^{\mathcal{N}} 2 \ln (\lvert \tau \! - \! \xi_{j} \rvert) 
\me^{-\frac{n}{\mathcal{N}} \widetilde{V}(\tau)} \, \md \tau \! 
\geqslant_{\underset{z_{o}=1+o(1)}{\mathscr{N},n \to \infty}} 
\! -\tfrac{n}{\mathcal{N}}(\mathcal{N} \! - \! 1) \mathfrak{c}_{36}
(n,k,z_{o})$, where $\mathfrak{c}_{36}(n,k,z_{o}) \! 
=_{\underset{z_{o}=1+o(1)}{\mathscr{N},n \to \infty}} 
\! \mathcal{O}(1)$; and (iv)
\begin{align*}
\dfrac{1}{\tilde{\mathfrak{z}}} \int_{\mathbb{R}} \left(
\sum_{q=1}^{\mathfrak{s}-2}-2 \varkappa_{nk \tilde{k}_{q}} \ln 
\lvert \tau \! - \! \alpha_{p_{q}} \rvert \! - \! 2(\varkappa_{nk} 
\! - \! 1) \ln \lvert \tau \! - \! \alpha_{k} \rvert \right) \me^{-
\frac{n}{\mathcal{N}} \widetilde{V}(\tau)} \, \md \tau &= 
\dfrac{1}{\tilde{\mathfrak{z}}} \left(\int_{(\mathbb{R} \setminus 
\mathfrak{D}_{M_{f}}) \setminus J_{f}} \! + \! \int_{J_{f}} \! + \! 
\int_{\lbrace \lvert \tau \rvert \geqslant T_{M_{f}} \rbrace} \! + 
\! \sum_{\substack{q=1\\q \neq \mathfrak{s}-1}}^{\mathfrak{s}} 
\int_{\mathscr{O}_{\frac{1}{T_{M_{f}}}}(\alpha_{p_{q}})} \right) \\
&\pmb{\cdot} \left(\sum_{q=1}^{\mathfrak{s}-2}-2 \varkappa_{nk 
\tilde{k}_{q}} \ln \lvert \tau \! - \! \alpha_{p_{q}} \rvert \! - \! 
2(\varkappa_{nk} \! - \! 1) \ln \lvert \tau \! - \! \alpha_{k} 
\rvert \right) \me^{-\frac{n}{\mathcal{N}} \widetilde{V}(\tau)} 
\, \md \tau,
\end{align*}
with
\begin{align*}
& \, \dfrac{1}{\tilde{\mathfrak{z}}} \int_{(\mathbb{R} \setminus 
\mathfrak{D}_{M_{f}}) \setminus J_{f}} \left(\sum_{q=1}^{
\mathfrak{s}-2}-2 \varkappa_{nk \tilde{k}_{q}} \ln \lvert \tau 
\! - \! \alpha_{p_{q}} \rvert \! - \! 2(\varkappa_{nk} \! - \! 1) \ln 
\lvert \tau \! - \! \alpha_{k} \rvert \right) \me^{-\frac{n}{\mathcal{N}} 
\widetilde{V}(\tau)} \, \md \tau \underset{\underset{z_{o}=1+o(1)}{
\mathscr{N},n \to \infty}}{\geqslant} -\dfrac{n}{\mathcal{N}}
(\mathcal{N} \! - \! 1) \\
&\times \, \left(\sum_{q=1}^{\mathfrak{s}-2} \dfrac{2 \varkappa_{nk 
\tilde{k}_{q}}}{n} \inf \lbrace \ln \lvert x \! - \! \alpha_{p_{q}} \rvert; \, 
x \! \in \! (\mathbb{R} \setminus \mathfrak{D}_{M_{f}}) \setminus J_{f} 
\rbrace \! + \! 2 \left(\dfrac{\varkappa_{nk} \! - \! 1}{n} \right) \inf 
\lbrace \ln \lvert x \! - \! \alpha_{k} \rvert; \, x \! \in \! (\mathbb{R} 
\setminus \mathfrak{D}_{M_{f}}) \setminus J_{f} \rbrace \right) \\
&\times \, \dfrac{1}{\tilde{\mathfrak{z}}} \me^{-\frac{n}{\mathcal{N}} 
\sup \lbrace \mathstrut \widetilde{V}(x); \, x \in (\mathbb{R} 
\setminus \mathfrak{D}_{M_{f}}) \setminus J_{f} \rbrace} 
\operatorname{meas}((\mathbb{R} \setminus \mathfrak{D}_{M_{f}}) 
\setminus J_{f}) \! =: \! -\dfrac{n}{\mathcal{N}}(\mathcal{N} \! - \! 1) 
\mathfrak{c}_{37}(n,k,z_{o}),
\end{align*}
\begin{align*}
& \, \dfrac{1}{\tilde{\mathfrak{z}}} \int_{J_{f}} \left(\sum_{q=1}^{
\mathfrak{s}-2}-2 \varkappa_{nk \tilde{k}_{q}} \ln \lvert \tau \! - 
\! \alpha_{p_{q}} \rvert \! - \! 2(\varkappa_{nk} \! - \! 1) \ln \lvert 
\tau \! - \! \alpha_{k} \rvert \right) \me^{-\frac{n}{\mathcal{N}} 
\widetilde{V}(\tau)} \, \md \tau \underset{\underset{z_{o}=1+
o(1)}{\mathscr{N},n \to \infty}}{\geqslant} -\dfrac{n}{\mathcal{N}}
(\mathcal{N} \! - \! 1) \\
&\times \, \left(\sum_{q=1}^{\mathfrak{s}-2} \dfrac{2 
\varkappa_{nk \tilde{k}_{q}}}{n} \inf \lbrace \ln \lvert x \! - \! \alpha_{k} 
\rvert; \, x \! \in \! J_{f} \rbrace \! + \! 2 \left(\dfrac{\varkappa_{nk} \! 
- \! 1}{n} \right) \inf \lbrace \ln \lvert x \! - \! \alpha_{k} \rvert; \, x \! 
\in \! J_{f} \rbrace \right) \\&\times \, \dfrac{1}{\tilde{\mathfrak{z}}} 
\me^{-\frac{n}{\mathcal{N}} \sup \lbrace \mathstrut \widetilde{V}(x); 
\, x \in J_{f} \rbrace} \sum_{j=1}^{N+1} \lvert \tilde{b}_{j-1} \! - \! 
\tilde{a}_{j} \rvert \! =: \! -\dfrac{n}{\mathcal{N}}(\mathcal{N} \! - \! 1) 
\mathfrak{c}_{38}(n,k,z_{o}),
\end{align*}
\begin{align*}
& \, \dfrac{1}{\tilde{\mathfrak{z}}} \int_{\lbrace \lvert \tau \rvert 
\geqslant T_{M_{f}} \rbrace} \left(\sum_{q=1}^{\mathfrak{s}-2}
-2 \varkappa_{nk \tilde{k}_{q}} \ln \lvert \tau \! - \! \alpha_{p_{q}} 
\rvert \! - \! 2(\varkappa_{nk} \! - \! 1) \ln \lvert \tau \! - \! 
\alpha_{k} \rvert \right) \me^{-\frac{n}{\mathcal{N}} \widetilde{V}
(\tau)} \, \md \tau \underset{\underset{z_{o}=1+o(1)}{\mathscr{N},
n \to \infty}}{\geqslant} -\dfrac{n}{\mathcal{N}}(\mathcal{N} \! - \! 
1) \\
&\times \, \left(\sum_{q=1}^{\mathfrak{s}-2} \dfrac{2 
\varkappa_{nk \tilde{k}_{q}}}{n} \ln \lvert T_{M_{f}}^{2} \! - \! 
\alpha_{p_{q}}^{2} \rvert \! + \! 2 \left(\dfrac{\varkappa_{nk} \! - 
\! 1}{n} \right) \ln \lvert T_{M_{f}}^{2} \! - \! \alpha_{k}^{2} \rvert 
\right) \dfrac{1}{\tilde{\mathfrak{z}}} \left(\dfrac{1}{2} \sqrt{
\dfrac{\pi}{2 \tilde{c}_{\infty}}}- \! T_{M_{f}} \right) \! =: \! 
-\dfrac{n}{\mathcal{N}}(\mathcal{N} \! - \! 1) \mathfrak{c}_{39}
(n,k,z_{o}),
\end{align*}
\begin{align*}
& \, \dfrac{1}{\tilde{\mathfrak{z}}} \sum_{\substack{q=1\\q 
\neq \mathfrak{s}-1}}^{\mathfrak{s}} \int_{\mathscr{O}_{\frac{1}{
T_{M_{f}}}}(\alpha_{p_{q}})} \left(\sum_{q=1}^{\mathfrak{s}-2}
-2 \varkappa_{nk \tilde{k}_{q}} \ln \lvert \tau \! - \! \alpha_{
p_{q}} \rvert \! - \! 2(\varkappa_{nk} \! - \! 1) \ln \lvert \tau 
\! - \! \alpha_{k} \rvert \right) \me^{-\frac{n}{\mathcal{N}} 
\widetilde{V}(\tau)} \, \md \tau \underset{\underset{z_{o}=1+
o(1)}{\mathscr{N},n \to \infty}}{\geqslant} -\dfrac{n}{\mathcal{
N}}(\mathcal{N} \! - \! 1) \\
&\times \, \left(\dfrac{2T_{M_{f}}^{-\left(1+2\frac{n}{\mathcal{N}}
(1+\tilde{c}_{\mathfrak{s}}) \right)}}{\tilde{\mathfrak{z}}(2 \frac{n}{
\mathcal{N}}(1 \! + \! \tilde{c}_{\mathfrak{s}}) \! + \! 1)} \ln 
\left(\left(((\alpha_{k} \! - \! T_{M_{f}}^{-1})^{2} \! + \! 1)
(\alpha_{k}^{2} \! + \! 1) \right)^{\frac{\varkappa_{nk}-1}{n}} 
\prod_{q^{\prime}=1}^{\mathfrak{s}-2} \left(((\alpha_{k} \! - \! 
T_{M_{f}}^{-1})^{2} \! + \! 1)(\alpha_{p_{q^{\prime}}}^{2} \! + \! 
1) \right)^{\frac{\varkappa_{nk \tilde{k}_{q^{\prime}}}}{n}} 
\right) \right. \\
&\left. \, +\sum_{q=1}^{\mathfrak{s}-2} \dfrac{2T_{M_{f}}^{-
\left(1+2\frac{n}{\mathcal{N}}(1+\tilde{c}_{q}) \right)}}{\tilde{
\mathfrak{z}}(2 \frac{n}{\mathcal{N}}(1 \! + \! \tilde{c}_{q}) \! + \! 
1)} \ln \left(\left(((\alpha_{p_{q}} \! - \! T_{M_{f}}^{-1})^{2} \! + \! 
1)(\alpha_{k}^{2} \! + \! 1) \right)^{\frac{\varkappa_{nk}-1}{n}} 
\left(((\alpha_{p_{q}} \! - \! T_{M_{f}}^{-1})^{2} \! + \! 1)
(\alpha_{p_{q}}^{2} \! + \! 1) \right)^{\frac{\varkappa_{nk 
\tilde{k}_{q}}}{n}} \right. \right. \\
&\left. \left. \, \times \prod_{\substack{q^{\prime}=1\\q^{\prime} 
\neq q}}^{\mathfrak{s}-2} \left(((\alpha_{p_{q^{\prime}}} \! - \! 
T_{M_{f}}^{-1})^{2} \! + \! 1)(\alpha_{p_{q^{\prime}}}^{2} \! + \! 
1) \right)^{\frac{\varkappa_{nk \tilde{k}_{q^{\prime}}}}{n}} \right) 
\right) \! =: \! -\dfrac{n}{\mathcal{N}}(\mathcal{N} \! - \! 1) 
\mathfrak{c}_{40}(n,k,z_{o}),
\end{align*}
where $\mathfrak{c}_{m}(n,k,z_{o}) \! =_{\underset{z_{o}=1+
o(1)}{\mathscr{N},n \to \infty}} \! \mathcal{O}(1)$, $m \! = \! 
37,38,39,40$, whence, gathering the above-derived bounds, 
it follows that $\tilde{\mathfrak{z}}^{-1} \int_{\mathbb{R}}
(\sum_{q=1}^{\mathfrak{s}-2}-2 \varkappa_{nk \tilde{k}_{q}} \ln 
\lvert \tau \! - \! \alpha_{p_{q}} \rvert \! - \! 2(\varkappa_{nk} \! - \! 1) 
\ln \lvert \tau \! - \! \alpha_{k} \rvert) \me^{-\frac{n}{\mathcal{N}} 
\widetilde{V}(\tau)} \, \md \tau \! \geqslant_{\underset{z_{o}=1+
o(1)}{\mathscr{N},n \to \infty}} \! -\tfrac{n}{\mathcal{N}}(\mathcal{N} 
\! - \! 1) \mathfrak{c}_{41}(n,k,z_{o})$, where $\mathfrak{c}_{41}
(n,k,z_{o}) \! =_{\underset{z_{o}=1+o(1)}{\mathscr{N},n \to \infty}} 
\! \mathcal{O}(1)$. Assembling all of the above-derived estimates, 
one concludes, via Equation~\eqref{eqlmrtz13}, that, for $n \! 
\in \! \mathbb{N}$ and $k \! \in \! \lbrace 1,2,\dotsc,K \rbrace$ 
such that $\alpha_{p_{\mathfrak{s}}} \! := \! \alpha_{k} \! \neq \! 
\infty$,
\begin{equation*}
\int_{\mathbb{R}} \me^{-n \widetilde{V}(\tau)} \prod_{j=2}^{
\mathcal{N}}(\tau \! - \! \xi_{j})^{2} \left(\prod_{q=1}^{
\mathfrak{s}-2}(\tau \! - \! \alpha_{p_{q}})^{\varkappa_{nk 
\tilde{k}_{q}}}(\tau \! - \! \alpha_{k})^{\varkappa_{nk}-1} \right)^{-2} 
\md \tau \underset{\underset{z_{o}=1+o(1)}{\mathscr{N},n 
\to \infty}}{\geqslant} \tilde{\mathfrak{z}} \me^{-\frac{n}{
\mathcal{N}}(\mathcal{N} \! - \! 1) \mathfrak{c}_{42}(n,k,z_{o})},
\end{equation*}
where $\mathfrak{c}_{42}(n,k,z_{o}) \! =_{\underset{z_{o}=1+
o(1)}{\mathscr{N},n \to \infty}} \! \mathcal{O}(1)$; hence, via 
Equation~\eqref{eqlmrtz12} and the fact that $\tilde{\mathcal{
E}}^{n,\sharp}_{k}(1) \! = \! 1$, one arrives at
\begin{equation*}
\dfrac{\tilde{\mathscr{Z}}^{n}_{k}}{\tilde{\mathscr{Z}}^{n,
\sharp}_{k}} \underset{\underset{z_{o}=1+o(1)}{\mathscr{N},
n \to \infty}}{\geqslant} \tilde{\mathfrak{z}} \me^{-\frac{n}{
\mathcal{N}}(\mathcal{N}-1) \mathfrak{c}_{42}(n,k,z_{o})},
\end{equation*}
or, equivalently,
\begin{equation} \label{eqlmrtz14} 
\dfrac{\tilde{\mathscr{Z}}^{n,\sharp}_{k}}{\tilde{\mathscr{
Z}}^{n}_{k}} \underset{\underset{z_{o}=1+o(1)}{\mathscr{N},
n \to \infty}}{\leqslant} \dfrac{1}{\tilde{\mathfrak{z}}} 
\me^{\frac{n}{\mathcal{N}}(\mathcal{N}-1) \mathfrak{c}_{42}
(n,k,z_{o})}.
\end{equation}
Via Equations~\eqref{eqlmrtzr} and~\eqref{eqlmrtz14}, it follows that
\begin{align*}
\dfrac{1}{\mathcal{N}} \tilde{\mathfrak{R}}^{n,k}_{1}(\tau) 
\underset{\underset{z_{o}=1+o(1)}{\mathscr{N},n \to 
\infty}}{\leqslant}& \, \dfrac{\me^{-\frac{n}{\mathcal{N}}
(\mathcal{N}-1) \widetilde{V}(\tau)} \me^{\frac{n}{
\mathcal{N}}(\mathcal{N}-1) \mathfrak{c}_{42}(n,k,
z_{o})}}{\tilde{\mathfrak{z}}} \idotsint\limits_{\mathbb{R}^{
\mathcal{N}-1}} \prod_{j=2}^{\mathcal{N}}(\tau \! - \! 
\xi_{j})^{2} \left(\prod_{q=1}^{\mathfrak{s}-2}(\tau \! - \! 
\alpha_{p_{q}})^{\varkappa_{nk \tilde{k}_{q}}}(\tau \! - \! 
\alpha_{k})^{\varkappa_{nk}-1} \right)^{-2} \\
\times& \, \tilde{\mathcal{P}}^{n,\sharp}_{k}(\xi_{2},\xi_{3},
\dotsc,\xi_{\mathcal{N}}) \, \md \xi_{2} \, \md \xi_{2} \, 
\dotsb \, \md \xi_{\mathcal{N}}.
\end{align*}
Via the formula $\sum_{j=2}^{\mathcal{N}}1 \! = \! (\mathcal{N} \! - \! 1)$, 
the decomposition $\sum_{q=1}^{\mathfrak{s}-2} \varkappa_{nk \tilde{k}_{q}} 
\! + \! \varkappa_{nk \tilde{k}_{\mathfrak{s}-1}}^{\infty} \! + \! \varkappa_{nk} 
\! = \! (n \! - \! 1)K \! + \! k \! = \! \mathcal{N}$, and the inequality $\lvert 
y_{1} \! - \! y_{2} \rvert^{2} \! \leqslant \! (1 \! + \! y_{1}^{2})(1 \! + \! y_{2}^{2})$, 
$y_{1},y_{2} \! \in \! \mathbb{R}$, a calculation shows that
\begin{align*}
\prod_{j=2}^{\mathcal{N}}(\tau \! - \! \xi_{j})^{2} \left(
\prod_{q=1}^{\mathfrak{s}-2}(\tau \! - \! \alpha_{p_{q}})^{
\varkappa_{nk \tilde{k}_{q}}}(\tau \! - \! \alpha_{k})^{\varkappa_{nk}-1} 
\right)^{-2} \leqslant& \, (1 \! + \! \alpha_{k}^{2})^{\varkappa_{nk}-1}  
\prod_{q=1}^{\mathfrak{s}-2}(1 \! + \! \alpha_{p_{q}}^{2})^{
\varkappa_{nk \tilde{k}_{q}}}(1 \! + \! \tau^{2})^{\varkappa_{nk 
\tilde{k}_{\mathfrak{s}-1}}^{\infty}+1}(1 \! + \! (\tau \! - \! 
\alpha_{k})^{-2})^{\varkappa_{nk}-1} \\
\times& \, \prod_{q=1}^{\mathfrak{s}-2}(1 \! + \! (\tau 
\! - \! \alpha_{p_{q}})^{-2})^{\varkappa_{nk \tilde{k}_{q}}} 
\prod_{j=2}^{\mathcal{N}}(1 \! + \! \xi_{j})^{2}(1 \! + \! 
(\xi_{j} \! - \! \alpha_{k})^{-2})^{\frac{\varkappa_{nk}-1}{
\mathcal{N}-1}} \\
\times& \, \prod_{q=1}^{\mathfrak{s}-2}(1 \! + \! 
(\xi_{j} \! - \! \alpha_{p_{q}})^{-2})^{\frac{\varkappa_{nk 
\tilde{k}_{q}}}{\mathcal{N}-1}},
\end{align*}
whence, via the identity $y_{1}^{y_{2}} \! = \! \exp (y_{2} 
\ln y_{1})$, one arrives at, for $n \! \in \! \mathbb{N}$ 
and $k \! \in \! \lbrace 1,2,\dotsc,K \rbrace$ such that 
$\alpha_{p_{\mathfrak{s}}} \! := \! \alpha_{k} \! \neq \! 
\infty$,
\begin{align} \label{eqlmrtznr} 
\dfrac{1}{\mathcal{N}} \tilde{\mathfrak{R}}^{n,k}_{1}(\tau) 
\underset{\underset{z_{o}=1+o(1)}{\mathscr{N},n \to 
\infty}}{\leqslant}& \, \me^{-\frac{n}{\mathcal{N}}
(\mathcal{N}-1) \widetilde{V}(\tau)} \me^{\frac{n}{
\mathcal{N}}(\mathcal{N}-1) \mathfrak{c}_{43}(n,k,
z_{o})}(1 \! + \! \tau^{2})^{\varkappa_{nk \tilde{k}_{
\mathfrak{s}-1}}^{\infty}+1}(1 \! + \! (\tau \! - \! 
\alpha_{k})^{-2})^{\varkappa_{nk}-1} \prod_{q=1}^{
\mathfrak{s}-2}(1 \! + \! (\tau \! - \! \alpha_{p_{q}})^{-
2})^{\varkappa_{nk \tilde{k}_{q}}} \nonumber \\
\times& \, \idotsint\limits_{\mathbb{R}^{\mathcal{N}-1}} 
\prod_{j=2}^{\mathcal{N}}(1 \! + \! \xi_{j})^{2}(1 \! + \! 
(\xi_{j} \! - \! \alpha_{k})^{-2})^{\frac{\varkappa_{nk}-1}{
\mathcal{N}-1}} \prod_{q=1}^{\mathfrak{s}-2}(1 \! + \! 
(\xi_{j} \! - \! \alpha_{p_{q}})^{-2})^{\frac{\varkappa_{nk 
\tilde{k}_{q}}}{\mathcal{N}-1}} \tilde{\mathcal{P}}^{n,\sharp}_{k}
(\xi_{2},\xi_{3},\dotsc,\xi_{\mathcal{N}}) \nonumber \\
\times& \, \, \md \xi_{2} \, \md \xi_{2} \, \dotsb 
\, \md \xi_{\mathcal{N}},
\end{align}
where $\mathfrak{c}_{43}(n,k,z_{o}) \! =_{\underset{z_{o}
=1+o(1)}{\mathscr{N},n \to \infty}} \! \mathcal{O}(1)$. For 
$n \! \in \! \mathbb{N}$ and $k \! \in \! \lbrace 1,2,\dotsc,
K \rbrace$ such that $\alpha_{p_{\mathfrak{s}}} \! := \! 
\alpha_{k} \! \neq \! \infty$, let $\tilde{\mathbb{S}}^{n}_{k}$ 
denote the linear operator $\tilde{\mathbb{S}}^{n}_{k} 
\colon \widetilde{V} \! \to \! \widetilde{V} \! + \! \tilde{
\upsilon}_{0}$, where $\tilde{\upsilon}_{0} \! \in \! 
\mathbb{R}$: a calculation shows that 
(cf. Equation~\eqref{eqlmrtzr}) $\tilde{\mathbb{S}}^{n}_{k} 
\colon \widetilde{V} \! \to \! \widetilde{V} \! + \! \tilde{
\upsilon}_{0}$, $\mathcal{N}^{-1} \tilde{\mathfrak{R}}^{n,
k}_{1} \! \mapsto \! \tilde{\mathbb{S}}^{n}_{k} \mathcal{N}^{-1} 
\tilde{\mathfrak{R}}^{n,k}_{1} \! = \! \mathcal{N}^{-1} \tilde{
\mathfrak{R}}^{n,k}_{1}$, that is, $\mathcal{N}^{-1} \tilde{
\mathfrak{R}}^{n,k}_{1}$ is invariant under the action of 
$\tilde{\mathbb{S}}^{n}_{k}$. In analogy with 
Equation~\eqref{eql3.5c}, define, for $n \! \in \! \mathbb{N}$ 
and $k \! \in \! \lbrace 1,2,\dotsc,K \rbrace$ such that 
$\alpha_{p_{\mathfrak{s}}} \! := \! \alpha_{k} \! \neq \! \infty$,
\begin{equation*}
\mathscr{K}^{\widetilde{V},f}_{\mathcal{N}-1}(x_{2},x_{3},
\dotsc,x_{\mathcal{N}}) \! := \! \sum_{\substack{i,j=2\\i 
\neq j}}^{\mathcal{N}} \mathcal{K}_{\widetilde{V}}^{f}(x_{i},
x_{j}),
\end{equation*}
where the symmetric function $\mathcal{K}^{f}_{\widetilde{V}}(\xi,\tau)$ 
is defined by Equation~\eqref{eqKvinf1}: via the latter definition, and 
Equations~\eqref{eqKvinf1} and~\eqref{eqKvinf4}, a calculation shows that
\begin{align} \label{eqlmrtz15} 
\mathscr{K}^{\widetilde{V},f}_{\mathcal{N}-1}(x_{2},x_{3},
\dotsc,x_{\mathcal{N}}) =& \, \sum_{\substack{i,j=2\\i \neq 
j}}^{\mathcal{N}} \left(\left(\dfrac{\varkappa_{nk \tilde{k}_{
\mathfrak{s}-1}}^{\infty} \! + \! 1}{n} \right) \ln \left(\dfrac{
1}{\lvert x_{i} \! - \! x_{j} \rvert} \right) \! + \! \left(\dfrac{
\varkappa_{nk} \! - \! 1}{n} \right) \ln \left(\dfrac{\lvert x_{i} 
\! - \! \alpha_{k} \rvert \lvert x_{j} \! - \! \alpha_{k} \rvert}{
\lvert x_{i} \! - \! x_{j} \rvert} \right) \right. \nonumber \\
+&\left. \, \sum_{q=1}^{\mathfrak{s}-2} \dfrac{\varkappa_{n
k \tilde{k}_{q}}}{n} \ln \left(\dfrac{\lvert x_{i} \! - \! \alpha_{
p_{q}} \rvert \lvert x_{j} \! - \! \alpha_{p_{q}} \rvert}{\lvert 
x_{i} \! - \! x_{j} \rvert} \right) \right) \! + \! (\mathcal{N} \! 
- \! 2) \sum_{j=2}^{\mathcal{N}} \widetilde{V}(x_{j}) \nonumber \\
\geqslant& \, (\mathcal{N} \! - \! 2) \sum_{j=2}^{\mathcal{N}} 
\hat{\psi}^{f}_{\widetilde{V}}(x_{j}) \! \geqslant \! (\mathcal{N} 
\! - \! 2) \sum_{j=2}^{\mathcal{N}} \left(\widetilde{V}(x_{j}) 
\! - \! \left(\dfrac{\varkappa_{nk \tilde{k}_{\mathfrak{s}-1}}^{
\infty} \! + \! 1}{n} \right) \ln (1 \! + \! x_{j}^{2}) \right. 
\nonumber \\
-&\left. \, \left(\dfrac{\varkappa_{nk} \! - \! 1}{n} \right) \ln 
(1 \! + \! (x_{j} \! - \! \alpha_{k})^{-2}) \! - \! \sum_{q=1}^{
\mathfrak{s}-2} \dfrac{\varkappa_{nk \tilde{k}_{q}}}{n} \ln 
(1 \! + \! (x_{j} \! - \! \alpha_{p_{q}})^{-2}) \right),
\end{align}
whence, via the invariance of $\mathcal{N}^{-1} \tilde{\mathfrak{
R}}^{n,k}_{1}$ under the action of $\tilde{\mathbb{S}}^{n}_{k}$, 
that is,
\begin{equation*}
\widetilde{V}(x_{j}) \! \geqslant \! (1 \! + \! \tilde{c}_{\infty}) \ln 
(1 \! + \! x_{j}^{2}) \! + \! (1 \! + \! \tilde{c}_{\mathfrak{s}}) \ln 
(1 \! + \! (x_{j} \! - \! \alpha_{k})^{-2}) \! + \! \sum_{q=1}^{
\mathfrak{s}-2}(1 \! + \! \tilde{c}_{q}) \ln (1 \! + \! (x_{j} \! - \! 
\alpha_{p_{q}})^{-2}), \quad j \! = \! 2,3,\dotsc,\mathcal{N},
\end{equation*}
one arrives at
\begin{align*}
\mathscr{K}^{\widetilde{V},f}_{\mathcal{N}-1}(x_{2},x_{3},
\dotsc,x_{\mathcal{N}}) \geqslant& \, (\mathcal{N} \! - \! 2) 
\sum_{j=2}^{\mathcal{N}} \left(\left(1 \! + \! \tilde{c}_{\infty} 
\! - \! \left(\dfrac{\varkappa_{nk \tilde{k}_{\mathfrak{s}-1}}^{
\infty} \! + \! 1}{n} \right) \right) \ln (1 \! + \! x_{j}^{2}) \! + 
\! \left(1 \! + \! \tilde{c}_{\mathfrak{s}} \! - \! \left(\dfrac{
\varkappa_{nk} \! - \! 1}{n} \right) \right) \ln (1 \! + \! (x_{j} 
\! - \! \alpha_{k})^{-2}) \right. \\
+&\left. \, \sum_{q=1}^{\mathfrak{s}-2} \left(1 \! + \! 
\tilde{c}_{q} \! - \! \dfrac{\varkappa_{nk \tilde{k}_{q}}}{n} 
\right) \ln (1 \! + \! (x_{j} \! - \! \alpha_{p_{q}})^{-2}) \right);
\end{align*}
thus, for any $\mathfrak{c}_{m}(n,k,z_{o}) \! \geqslant_{
\underset{z_{o}=1+o(1)}{\mathscr{N},n \to \infty}} \! 0$ and 
$\mathcal{O}(1)$, $m \! = \! 44,45,46$, such that $\sum_{j
=2}^{\mathcal{N}} \ln (1 \! + \! x_{j}^{2}) \! \geqslant \! 
(\mathcal{N} \! - \! 1) \mathfrak{c}_{44}(n,k,z_{o})$, $\sum_{j
=2}^{\mathcal{N}} \ln (1 \! + \! (x_{j} \! - \! \alpha_{k})^{-2}) 
\! \geqslant \! (\mathcal{N} \! - \! 1) \mathfrak{c}_{45}(n,k,
z_{o})$, and $\sum_{j=2}^{\mathcal{N}} \sum_{q=1}^{
\mathfrak{s}-2} \ln (1 \! + \! (x_{j} \! - \! \alpha_{p_{q}})^{-2}) 
\! \geqslant \! (\mathcal{N} \! - \! 1) \mathfrak{c}_{46}(n,k,
z_{o})$, it follows that
\begin{equation} \label{eqlmrtz16} 
\mathscr{K}^{\widetilde{V},f}_{\mathcal{N}-1}(x_{2},x_{3},
\dotsc,x_{\mathcal{N}}) \! \geqslant \! (\mathcal{N} \! - \! 1)
(\mathcal{N} \! - \! 2) \mathfrak{c}_{47}(n,k,z_{o}),
\end{equation}
where $\mathfrak{c}_{47}(n,k,z_{o})$ $(:= \! \min \lbrace 
(1 \! + \! \tilde{c}_{\infty} \! - \! (\varkappa_{nk \tilde{k}_{
\mathfrak{s}-1}}^{\infty} \! + \! 1)/n) \mathfrak{c}_{44}
(n,k,z_{o}),(1 \! + \! \tilde{c}_{\mathfrak{s}} \! - \! 
(\varkappa_{nk} \! - \! 1)/n) \mathfrak{c}_{45}(n,k,z_{o}),
\min_{q=1,2,\dotsc,\mathfrak{s}-2} \lbrace (1 \! + \! 
\tilde{c}_{q} \! - \! \varkappa_{nk \tilde{k}_{q}}/n) 
\mathfrak{c}_{46}(n,k,z_{o}) \rbrace)$ $=_{\underset{z_{o}=
1+o(1)}{\mathscr{N},n \to \infty}} \! \mathcal{O}(1)$ and $> \! 
0$. Via Equations~\eqref{eqlmrtzzs} and~\eqref{eqlmrtz15}, 
one mimics the above analysis (leading to the 
Estimate~\eqref{eqlmrtz6}) to show that, for $n \! \in \! 
\mathbb{N}$ and $k \! \in \! \lbrace 1,2,\dotsc,K \rbrace$ 
such that $\alpha_{p_{\mathfrak{s}}} \! := \! \alpha_{k} \! 
\neq \! \infty$,
\begin{equation} \label{eqlmrtz17} 
\dfrac{\mathcal{N}}{n} \dfrac{1}{(\mathcal{N} \! - \! 1)
(\mathcal{N} \! - \! 2)} \ln (\tilde{\mathscr{Z}}^{n,\sharp}_{k}) 
\underset{\underset{z_{o}=1+o(1)}{\mathscr{N},n \to \infty}}{
\geqslant} -\left(E_{\widetilde{V}}^{f} \! + \! \tilde{\varepsilon} 
\right) \left(1 \! + \! \mathcal{O} \left(\dfrac{\mathfrak{c}_{48}
(n,k,z_{o},\tilde{\varepsilon})}{n} \right) \right) \quad \text{as} 
\quad \tilde{\varepsilon} \! \downarrow \! 0,
\end{equation}
where $\mathfrak{c}_{48}(n,k,z_{o},\tilde{\varepsilon}) \! =_{
\underset{z_{o}=1+o(1)}{\mathscr{N},n \to \infty}} \! \mathcal{O}
(1)$ as $\tilde{\varepsilon} \! \downarrow \! 0$. For $n \! \in \! 
\mathbb{N}$ and $k \! \in \! \lbrace 1,2,\dotsc,K \rbrace$ such 
that $\alpha_{p_{\mathfrak{s}}} \! := \! \alpha_{k} \! \neq \! \infty$, 
in the double-scaling limit $\mathscr{N},n \! \to \! \infty$ such that 
$z_{o} \! = \! 1 \! + \! o(1)$, let $\tilde{\mathbb{A}}_{\mathcal{N}
-1}^{\sharp} \! := \! \lbrace \mathstrut (x_{2},x_{3},\dotsc,
x_{\mathcal{N}}) \! \in \! \mathbb{R}^{\mathcal{N}-1}; \, 
((\mathcal{N} \! - \! 1)(\mathcal{N} \! - \! 2))^{-1} \mathscr{K}_{
\mathcal{N}-1}^{\widetilde{V},f}(x_{2},x_{3},\dotsc,\linebreak[4]
x_{\mathcal{N}}) \! \leqslant \! \mathfrak{c}_{47}(n,k,z_{o}) 
\rbrace$ $(\subset \mathbb{R}^{\mathcal{N}-1})$: via this definition, 
adjusting all $\mathcal{O}(1)$ parameters, if necessary, so that 
$\mathfrak{c}_{49}(n,k,z_{o}) \! := \! \mathfrak{c}_{47}(n,k,
z_{o}) \! - \! (E_{\widetilde{V}}^{f} \! + \! \tilde{\varepsilon}) \! 
>_{\underset{z_{o}=1+o(1)}{\mathscr{N},n \to \infty}} \! 0$ (and 
$\mathcal{O}(1))$, and Equations~\eqref{eqlmrtz15}---\eqref{eqlmrtz17}, 
one shows, by mimicking the above calculations (leading to the 
Estimate~\eqref{eqlmrtz7}), that, analogously,
\begin{equation*}
\operatorname{Prob} \, (\mathbb{R}^{\mathcal{N}-1} \setminus 
\tilde{\mathbb{A}}_{\mathcal{N}-1}^{\sharp}) \underset{
\underset{z_{o}=1+o(1)}{\mathscr{N},n \to \infty}}{\leqslant} 
\me^{-\frac{n}{\mathcal{N}}(\mathcal{N}-1)(\mathcal{N}-2) 
\mathfrak{c}_{49}(n,k,z_{o})} \underbrace{\left(
\int_{\mathbb{R}} \left(\dfrac{\me^{-\frac{n}{\mathcal{N}} 
\widetilde{V}(\xi)}}{\prod_{\underset{q \neq \mathfrak{s}-1}{q
=1}}^{\mathfrak{s}} \lvert \xi \! - \! \alpha_{p_{q}} \rvert^{
\tilde{\gamma}_{q}(1+\mathcal{O}(n^{-1}))}} \right)^{2} 
\md \xi \right)^{\mathcal{N}-1}}_{=: \, \mathfrak{c}_{50}
(n,k,z_{o})} \quad \Rightarrow
\end{equation*}
\begin{equation} \label{eqlmrtz18} 
\operatorname{Prob} \, (\mathbb{R}^{\mathcal{N}-1} 
\setminus \tilde{\mathbb{A}}_{\mathcal{N}-1}^{\sharp}) 
\underset{\underset{z_{o}=1+o(1)}{\mathscr{N},n \to \infty}}{
\leqslant} \mathfrak{c}_{50}(n,k,z_{o}) \me^{-\frac{n}{\mathcal{
N}}(\mathcal{N}-1)(\mathcal{N}-2) \mathfrak{c}_{49}(n,k,z_{o})},
\end{equation}
where $\mathfrak{c}_{50}(n,k,z_{o}) \! =_{\underset{z_{o}=1
+o(1)}{\mathscr{N},n \to \infty}} \! \mathcal{O}(1)$. For $n 
\! \in \! \mathbb{N}$ and $k \! \in \! \lbrace 1,2,\dotsc,K 
\rbrace$ such that $\alpha_{p_{\mathfrak{s}}} \! := \! 
\alpha_{k} \! \neq \! \infty$, let $\mathfrak{c}_{51}(n,k,z_{o}) 
\! := \! \mathfrak{c}_{44}(n,k,z_{o}) \! + \! (\mathcal{N} \! - 
\! 1)^{-1}(\varkappa_{nk} \! - \! 1) \mathfrak{c}_{45}(n,k,
z_{o}) \! + \! (\mathcal{N} \! - \! 1)^{-1} \mathfrak{c}_{46}(n,k,z_{o}) 
\min\limits_{q=1,2,\dotsc,\mathfrak{s}-2} \lbrace \varkappa_{nk 
\tilde{k}_{q}} \rbrace$, with $\mathfrak{c}_{51}(n,k,z_{o}) \! 
=_{\underset{z_{o}=1+o(1)}{\mathscr{N},n \to \infty}} \! 
\mathcal{O}(1)$ and $> \! 0$, and define the following 
$(\mathcal{N} \! - \! 1)$-dimensional subsets of 
$\mathbb{R}^{\mathcal{N}-1}$:
\begin{gather*}
\mathbb{S}^{n,\sharp}_{k,1} \! := \! \left\lbrace \mathstrut 
(x_{2},x_{3},\dotsc,x_{\mathcal{N}}) \! \in \! \mathbb{R}^{\mathcal{N}
-1}; \, \prod_{j=2}^{\mathcal{N}}(1 \! + \! x_{j}^{2})(1 \! + \! (x_{j} 
\! - \! \alpha_{k})^{-2})^{\frac{\varkappa_{nk}-1}{\mathcal{N}-1}} 
\prod_{q=1}^{\mathfrak{s}-2}(1 \! + \! (x_{j} \! - \! \alpha_{
p_{q}})^{-2})^{\frac{\varkappa_{nk \tilde{k}_{q}}}{\mathcal{N}-1}} 
\underset{\underset{z_{o}=1+o(1)}{\mathscr{N},n \to \infty}}{
\geqslant} \me^{(\mathcal{N}-1) \mathfrak{c}_{51}(n,k,z_{o})} 
\right\rbrace, \\
\mathbb{S}^{n,\sharp}_{k,2} \! := \! \left\lbrace \mathstrut 
(x_{2},x_{3},\dotsc,x_{\mathcal{N}}) \! \in \! \mathbb{R}^{
\mathcal{N}-1}; \, \prod_{j=2}^{\mathcal{N}}(1 \! + \! x_{j}^{2})
(1 \! + \! (x_{j} \! - \! \alpha_{k})^{-2})^{\frac{\varkappa_{nk}-1}{
\mathcal{N}-1}} \prod_{q=1}^{\mathfrak{s}-2}(1 \! + \! (x_{j} \! - \! 
\alpha_{p_{q}})^{-2})^{\frac{\varkappa_{nk \tilde{k}_{q}}}{\mathcal{N}
-1}} \underset{\underset{z_{o}=1+o(1)}{\mathscr{N},n \to \infty}}{
\leqslant} \me^{(\mathcal{N}-1) \mathfrak{c}_{51}(n,k,z_{o})} 
\right\rbrace.
\end{gather*}
(Note that $\mathbb{S}^{n,\sharp}_{k,1}$ and $\mathbb{S}^{
n,\sharp}_{k,2}$ may have a non-empty intersection.) Using 
the fact  that $\chi_{\mathbb{S}^{n,\sharp}_{k,1} \cup 
\mathbb{S}^{n,\sharp}_{k,2}}(\vec{x}) \! = \! \chi_{\mathbb{
S}^{n,\sharp}_{k,1}}(\vec{x}) \! + \! \chi_{\mathbb{S}^{n,
\sharp}_{k,2}}(\vec{x}) \! - \! \chi_{\mathbb{S}^{n,\sharp}_{k,
1} \cap \mathbb{S}^{n,\sharp}_{k,2}}(\vec{x})$,\footnote{If 
$\mathbb{S}^{n,\sharp}_{k,1} \cap \mathbb{S}^{n,\sharp}_{k,2} 
\! = \! \varnothing$, then, since $\chi_{\varnothing}
(\pmb{\cdot}) \! = \! 0$, $\chi_{\mathbb{S}^{n,\sharp}_{k,1} 
\cup \mathbb{S}^{n,\sharp}_{k,2}}(\vec{x}) \! = \! 
\chi_{\mathbb{S}^{n,\sharp}_{k,1}}(\vec{x}) \! + \! 
\chi_{\mathbb{S}^{n,\sharp}_{k,2}}(\vec{x})$.} and 
$\operatorname{meas}(\mathbb{S}^{n,\sharp}_{k,1} \cup 
\mathbb{S}^{n,\sharp}_{k,2}) \! \leqslant \! 
\operatorname{meas}(\mathbb{S}^{n,\sharp}_{k,1}) \! + \! 
\operatorname{meas}(\mathbb{S}^{n,\sharp}_{k,2}) \! - \! 
\operatorname{meas}(\mathbb{S}^{n,\sharp}_{k,1} \cap 
\mathbb{S}^{n,\sharp}_{k,2})$,\footnote{See, for example, 
\cite{hkestm}, p.~283, Exercise~\#74.} one estimates, 
via the linearity of $\tilde{\mathcal{E}}^{n,\sharp}_{k}$ 
and the analogue of the calculation leading to 
Equation~\eqref{eqlmrtz18}, the $(\mathcal{N} \! - \! 
1)$-dimensional multiple integral appearing in 
Equation~\eqref{eqlmrtznr} as follows:
\begin{align*}
& \, \idotsint\limits_{\mathbb{R}^{\mathcal{N}-1}} 
\prod_{j=2}^{\mathcal{N}}(1 \! + \! \xi_{j})^{2}(1 \! + \! 
(\xi_{j} \! - \! \alpha_{k})^{-2})^{\frac{\varkappa_{nk}-1}{
\mathcal{N}-1}} \prod_{q=1}^{\mathfrak{s}-2}(1 \! + \! 
(\xi_{j} \! - \! \alpha_{p_{q}})^{-2})^{\frac{\varkappa_{nk 
\tilde{k}_{q}}}{\mathcal{N}-1}} \tilde{\mathcal{P}}^{n,\sharp}_{k}
(\xi_{2},\xi_{3},\dotsc,\xi_{\mathcal{N}}) \, \md \xi_{2} 
\, \md \xi_{2} \, \dotsb \, \md \xi_{\mathcal{N}} \\
& \, \underset{\underset{z_{o}=1+o(1)}{\mathscr{N},n \to 
\infty}}{\leqslant} \underbrace{\tilde{\mathcal{E}}^{n,\sharp}_{k} 
\left(\chi_{\mathbb{S}^{n,\sharp}_{k,1}}(\vec{x}) \prod_{j=
2}^{\mathcal{N}}(1 \! + \! x_{j})^{2}(1 \! + \! (x_{j} \! - \! 
\alpha_{k})^{-2})^{\frac{\varkappa_{nk}-1}{\mathcal{N}-1}} 
\prod_{q=1}^{\mathfrak{s}-2}(1 \! + \! (x_{j} \! - \! 
\alpha_{p_{q}})^{-2})^{\frac{\varkappa_{nk \tilde{k}_{q}}}{
\mathcal{N}-1}} \right)}_{\leqslant_{\underset{z_{o}=1+
o(1)}{\mathscr{N},n \to \infty}} \, \mathfrak{c}_{52}(n,k,z_{o}) 
\exp (-\frac{n}{\mathcal{N}}(\mathcal{N}-1)(\mathcal{N}-2) 
\mathfrak{c}_{53}(n,k,z_{o}))} \\
& \, +\tilde{\mathcal{E}}^{n,\sharp}_{k} \left(\chi_{\mathbb{
S}^{n,\sharp}_{k,2}}(\vec{x}) \prod_{j=2}^{\mathcal{N}}(1 \! 
+ \! x_{j})^{2}(1 \! + \! (x_{j} \! - \! \alpha_{k})^{-2})^{\frac{
\varkappa_{nk}-1}{\mathcal{N}-1}} \prod_{q=1}^{\mathfrak{s}
-2}(1 \! + \! (x_{j} \! - \! \alpha_{p_{q}})^{-2})^{\frac{
\varkappa_{nk \tilde{k}_{q}}}{\mathcal{N}-1}} \right) 
\underset{\underset{z_{o}=1+o(1)}{\mathscr{N},n \to \infty}}{
\leqslant} \mathfrak{c}_{52}(n,k,z_{o}) \me^{-\frac{n}{
\mathcal{N}}(\mathcal{N}-1)(\mathcal{N}-2) \mathfrak{c}_{53}
(n,k,z_{o})} \\
& \, +\idotsint\limits_{\mathbb{S}^{n,\sharp}_{k,2}} 
\underbrace{\prod_{j=2}^{\mathcal{N}}(1 \! + \! \xi_{j})^{2}
(1 \! + \! (\xi_{j} \! - \! \alpha_{k})^{-2})^{\frac{\varkappa_{nk}
-1}{\mathcal{N}-1}} \prod_{q=1}^{\mathfrak{s}-2}(1 \! + \! 
(\xi_{j} \! - \! \alpha_{p_{q}})^{-2})^{\frac{\varkappa_{nk 
\tilde{k}_{q}}}{\mathcal{N}-1}}}_{\leqslant_{\underset{z_{o}=1
+o(1)}{\mathscr{N},n \to \infty}} \, \exp ((\mathcal{N}-1) 
\mathfrak{c}_{51}(n,k,z_{o}))} \tilde{\mathcal{P}}^{n,\sharp}_{k}
(\xi_{2},\xi_{3},\dotsc,\xi_{\mathcal{N}}) \, \md \xi_{2} \, \md 
\xi_{2} \, \dotsb \, \md \xi_{\mathcal{N}} \\
& \, \underset{\underset{z_{o}=1+o(1)}{\mathscr{N},n \to 
\infty}}{\leqslant} \mathfrak{c}_{52}(n,k,z_{o}) \me^{-\frac{n}{
\mathcal{N}}(\mathcal{N}-1)(\mathcal{N}-2) \mathfrak{c}_{53}
(n,k,z_{o})} \! + \! \me^{(\mathcal{N}-1) \mathfrak{c}_{51}
(n,k,z_{o})} \underbrace{\idotsint\limits_{\mathbb{R}^{
\mathcal{N}-1}} \tilde{\mathcal{P}}^{n,\sharp}_{k}(\xi_{2},
\xi_{3},\dotsc,\xi_{\mathcal{N}}) \, \md \xi_{2} \, \md 
\xi_{2} \, \dotsb \, \md \xi_{\mathcal{N}}}_{=1} \\
& \, \underset{\underset{z_{o}=1+o(1)}{\mathscr{N},n \to 
\infty}}{\leqslant} \me^{\frac{n}{\mathcal{N}}(\mathcal{N}-1) 
\mathfrak{c}_{54}(n,k,z_{o})} \left(1 \! + \! \mathcal{O} \left(
\mathfrak{c}_{55}(n,k,z_{o}) \me^{-\frac{n}{\mathcal{N}}
(\mathcal{N}-1)(\mathcal{N}-2) \mathfrak{c}_{53}(n,k,z_{o})} 
\right) \right),
\end{align*}
where $\mathfrak{c}_{m}(n,k,z_{o}) \! =_{\underset{z_{o}=1
+o(1)}{\mathscr{N},n \to \infty}} \! \mathcal{O}(1)$, $m \! = 
\! 52,53,54,55$; in particular, $\mathfrak{c}_{53}(n,k,z_{o}),
\mathfrak{c}_{54}(n,k,z_{o}) \! > \! 0$: via this latter estimate 
and Equation~\eqref{eqlmrtznr}, one arrives at, for $n \! \in 
\! \mathbb{N}$ and $k \! \in \! \lbrace 1,2,\dotsc,K \rbrace$ 
such that $\alpha_{p_{\mathfrak{s}}} \! := \! \alpha_{k} \! 
\neq \! \infty$,
\begin{align} \label{eqlmrtz19} 
\dfrac{1}{\mathcal{N}} \tilde{\mathfrak{R}}^{n,k}_{1}(\tau) 
\underset{\underset{z_{o}=1+o(1)}{\mathscr{N},n \to 
\infty}}{\leqslant}& \, \me^{-\frac{n}{\mathcal{N}}
(\mathcal{N}-1) \widetilde{V}(\tau)} \me^{\frac{n}{
\mathcal{N}}(\mathcal{N}-1) \mathfrak{c}_{56}(n,k,
z_{o})}(1 \! + \! \tau^{2})^{\varkappa_{nk \tilde{k}_{
\mathfrak{s}-1}}^{\infty}+1}(1 \! + \! (\tau \! - \! 
\alpha_{k})^{-2})^{\varkappa_{nk}-1} \prod_{q=1}^{
\mathfrak{s}-2}(1 \! + \! (\tau \! - \! \alpha_{p_{q}}
)^{-2})^{\varkappa_{nk \tilde{k}_{q}}} \nonumber \\
\times& \, \left(1 \! + \! \mathcal{O} \left(\mathfrak{c}_{55}
(n,k,z_{o}) \me^{-\frac{n}{\mathcal{N}}(\mathcal{N}-1)
(\mathcal{N}-2) \mathfrak{c}_{53}(n,k,z_{o})} \right) \right),
\end{align}
where $\mathfrak{c}_{56}(n,k,z_{o}) \! =_{\underset{z_{o}
=1+o(1)}{\mathscr{N},n \to \infty}} \! \mathcal{O}(1)$.

For $n \! \in \! \mathbb{N}$ and $k \! \in \! \lbrace 1,2,
\dotsc,K \rbrace$ such that $\alpha_{p_{\mathfrak{s}}} \! := 
\! \alpha_{k} \! \neq \! \infty$, choose, for the time being (see 
the refinements below), $\tilde{M}_{0} \! = \! \tilde{M}_{0}
(n,k,z_{o}) \! \gg_{\underset{z_{o}=1+o(1)}{\mathscr{N},n 
\to \infty}} \! K(1 \! + \! \max_{q=1,\dotsc,\mathfrak{s}-2,
\mathfrak{s}} \lbrace \lvert \alpha_{p_{q}} \rvert \rbrace \! 
+ \! 3(\min_{i \neq j \in \lbrace 1,\dotsc,\mathfrak{s}-2,
\mathfrak{s} \rbrace} \lbrace \lvert \alpha_{p_{i}} \! - \! 
\alpha_{p_{j}} \rvert \rbrace)^{-1})$ $(\gg \! 1)$ such that 
$\mathscr{O}_{\frac{1}{\tilde{M}_{0}}}(\alpha_{p_{i}}) \cap 
\mathscr{O}_{\frac{1}{\tilde{M}_{0}}}(\alpha_{p_{j}}) \! = \! 
\varnothing$ $\forall$ $i \! \neq \! j \! \in \! \lbrace 1,
\dotsc,\mathfrak{s} \! - \! 2,\mathfrak{s} \rbrace$, and 
$\tilde{\mathfrak{D}}(\tilde{M}_{0}) \! := \! [-\tilde{M}_{0},\tilde{M}_{0}] 
\setminus \cup_{\underset{q \neq \mathfrak{s}-1}{q=1}}^{\mathfrak{s}} 
\mathscr{O}_{\frac{1}{\tilde{M}_{0}}}(\alpha_{p_{q}}) \supseteq 
J_{f}$.\footnote{Incidentally, with this choice of $\tilde{M}_{0}$, 
it also follows that $[-\tilde{M}_{0},\alpha_{p_{1}} \! - \! 
\tilde{M}_{0}^{-1}] \cap [\alpha_{p_{i}} \! + \! \tilde{M}_{0}^{-
1},\alpha_{p_{i+1}} \! - \! \tilde{M}_{0}^{-1}] \! = \! \varnothing$, 
$i \! = \! 1,2,\dotsc,\mathfrak{s} \! - \! 3$, $[-\tilde{M}_{0},
\alpha_{p_{1}} \! - \! \tilde{M}_{0}^{-1}] \cap [\alpha_{p_{
\mathfrak{s}-2}} \! + \! \tilde{M}_{0}^{-1},\alpha_{k} \! - \! 
\tilde{M}_{0}^{-1}] \! = \! \varnothing$, $[-\tilde{M}_{0},
\alpha_{p_{1}} \! - \! \tilde{M}_{0}^{-1}] \cap [\alpha_{k} \! 
+ \! \tilde{M}_{0}^{-1},\tilde{M}_{0}] \! = \! \varnothing$, 
$[\alpha_{k} \! + \! \tilde{M}_{0}^{-1},\tilde{M}_{0}] \cap 
[\alpha_{p_{\mathfrak{s}-2}} \! + \! \tilde{M}_{0}^{-1},
\alpha_{k} \! - \! \tilde{M}_{0}^{-1}] \! = \! \varnothing$, 
$[\alpha_{k} \! + \! \tilde{M}_{0}^{-1},\tilde{M}_{0}] \cap 
[\alpha_{p_{j}} \! + \! \tilde{M}_{0}^{-1},\alpha_{p_{j+1}} \! - \! 
\tilde{M}_{0}^{-1}] \! = \! \varnothing$, $j \! = \! 1,2,\dotsc,
\mathfrak{s} \! - \! 3$, and $[\alpha_{p_{i}} \! + \! \tilde{M}_{
0}^{-1},\alpha_{p_{i+1}} \! - \! \tilde{M}_{0}^{-1}] \cap 
[\alpha_{p_{j}} \! + \! \tilde{M}_{0}^{-1},\alpha_{p_{j+1}} \! - \! 
\tilde{M}_{0}^{-1}] \! = \! \varnothing$ $\forall$ $i \! \neq \! j 
\! \in \! \lbrace 1,2,\dotsc,\mathfrak{s} \! - \! 3 \rbrace$.} 
Set $\tilde{\mathfrak{D}}^{c}(\tilde{M}_{0}) \! := \! \mathbb{R} 
\setminus \tilde{\mathfrak{D}}(\tilde{M}_{0}) \! = \! 
\mathbb{A}(\tilde{M}_{0}) \cup \cup_{\underset{q \neq 
\mathfrak{s}-1}{q=1}}^{\mathfrak{s}} \mathbb{B}_{q}(\tilde{M}_{
0})$ (the complement of $\tilde{\mathfrak{D}}(\tilde{M}_{0})$ 
relative to $\mathbb{R})$, where $\mathbb{A}(\tilde{M}_{0}) 
\! := \! \lbrace \mathstrut x \! \in \! \mathbb{R}; \, \lvert x 
\rvert \! > \! \tilde{M}_{0} \rbrace$, and $\mathbb{B}_{q}
(\tilde{M}_{0}) \! := \! \lbrace x \! \in \! \mathbb{R}; \, 
\lvert x \! - \! \alpha_{p_{q}} \rvert \! < \! \tilde{M}_{0}^{-1} 
\rbrace$, $q \! = \! 1,\dotsc,\mathfrak{s} \! - \! 2,
\mathfrak{s}$.\footnote{Note that $\mathbb{A}(\tilde{M}_{0}) 
\cap \mathbb{B}_{q}(\tilde{M}_{0}) \! = \! \varnothing$, $q \! = \! 1,
\dotsc,\mathfrak{s} \! - \! 2,\mathfrak{s}$, and $\mathbb{B}_{i}
(\tilde{M}_{0}) \cap \mathbb{B}_{j}(\tilde{M}_{0}) \! = \! \varnothing$ 
$\forall$ $i \! \neq \! j \! \in \! \lbrace 1,\dotsc,\mathfrak{s} \! 
- \! 2,\mathfrak{s} \rbrace$; hence, $\mathbb{A}(\tilde{M}_{0}) 
\cap (\cup_{\underset{q \neq \mathfrak{s}-1}{q=1}}^{
\mathfrak{s}} \mathbb{B}_{q}(\tilde{M}_{0})) \! = \! \cup_{
\underset{q \neq \mathfrak{s}-1}{q=1}}^{\mathfrak{s}}
(\mathbb{A}(\tilde{M}_{0}) \cap \mathbb{B}_{q}(\tilde{M}_{0})) 
\! = \! \cup_{\underset{q \neq \mathfrak{s}-1}{q=1}}^{
\mathfrak{s}} \varnothing \! = \! \varnothing$.} Denote the 
$\mathcal{N}$-fold Cartesian product of $\tilde{\mathfrak{
D}}(\tilde{M}_{0})$ by $\tilde{\mathfrak{D}}_{\mathcal{N}}
(\tilde{M}_{0}) \! = \! \tilde{\mathfrak{D}}(\tilde{M}_{0}) \times 
\tilde{\mathfrak{D}}(\tilde{M}_{0}) \times \dotsb \times 
\tilde{\mathfrak{D}}(\tilde{M}_{0})$ $(\subset \mathbb{R}^{
\mathcal{N}})$, and let $\tilde{\mathfrak{D}}_{\mathcal{N}}^{c}
(\tilde{M}_{0}) \! := \! \mathbb{R}^{\mathcal{N}} \setminus 
\tilde{\mathfrak{D}}_{\mathcal{N}}(\tilde{M}_{0})$ (the complement 
of $\tilde{\mathfrak{D}}_{\mathcal{N}}(\tilde{M}_{0})$ relative to 
$\mathbb{R}^{\mathcal{N}})$. A calculation shows that $\chi_{\tilde{
\mathfrak{D}}^{c}(\tilde{M}_{0})}(\pmb{\cdot}) \! = \! \chi_{
\mathbb{A}(\tilde{M}_{0})}(\pmb{\cdot}) \! + \! \sum_{\underset{q 
\neq \mathfrak{s}-1}{q=1}}^{\mathfrak{s}} \chi_{\mathbb{B}_{q}
(\tilde{M}_{0})}(\pmb{\cdot})$. For $\vec{x} \! := \! (x_{1},x_{2},
\dotsc,x_{\mathcal{N}}) \! \in \! \tilde{\mathfrak{D}}_{\mathcal{
N}}^{c}(\tilde{M}_{0})$, it follows that $\exists$ $j \! \in 
\! \lbrace 1,2,\dotsc,\mathcal{N} \rbrace$ such that 
$\lvert x_{j} \rvert \! > \! \tilde{M}_{0}$ or $\lvert x_{j} 
\! - \! \alpha_{p_{q}} \rvert \! < \! \tilde{M}_{0}^{-1}$, 
$q \! = \! 1,\dotsc,\mathfrak{s} \! - \! 2,\mathfrak{s}$; 
hence, via the analogue of the argument on p.~170 of 
\cite{a51}, and Equation~\eqref{eqlmrtz11},
\begin{align*}
\operatorname{Prob}(\tilde{\mathfrak{D}}_{\mathcal{N}}^{c}
(\tilde{M}_{0})) \leqslant& \, \idotsint\limits_{\mathbb{R}^{
\mathcal{N}}} \sum_{j=1}^{\mathcal{N}} \chi_{\tilde{\mathfrak{D}}^{c}
(\tilde{M}_{0})}(\xi_{j}) \, \tilde{\mathcal{P}}^{n}_{k}(\xi_{1},\xi_{2},
\dotsc,\xi_{\mathcal{N}}) \, \md \xi_{1} \, \md \xi_{2} \, \dotsb 
\, \md \xi_{\mathcal{N}} \\
\leqslant& \, \idotsint\limits_{\mathbb{R}^{\mathcal{N}}} \sum_{j
=1}^{\mathcal{N}} \left(\chi_{\mathbb{A}(\tilde{M}_{0})}(\xi_{j}) 
\! + \! \sum_{\substack{q=1\\q \neq \mathfrak{s}-1}}^{
\mathfrak{s}} \chi_{\mathbb{B}_{q}(\tilde{M}_{0})}(\xi_{j}) \right) 
\tilde{\mathcal{P}}^{n}_{k}(\xi_{1},\xi_{2},\dotsc,\xi_{\mathcal{N}}) 
\, \md \xi_{1} \, \md \xi_{2} \, \dotsb \, \md \xi_{\mathcal{N}} \\
\leqslant& \, \int\nolimits_{\mathbb{R}} \sum_{j=1}^{\mathcal{N}} 
\chi_{\mathbb{A}(\tilde{M}_{0})}(\xi_{j}) \dfrac{1}{\mathcal{N}} 
\tilde{\mathfrak{R}}^{n,k}_{1}(\xi_{j}) \, \md \xi_{j} \! + \! 
\int\nolimits_{\mathbb{R}} \sum_{j=1}^{\mathcal{N}} \sum_{
\substack{q=1\\q \neq \mathfrak{s}-1}}^{\mathfrak{s}} \chi_{
\mathbb{B}_{q}(\tilde{M}_{0})}(\xi_{j}) \dfrac{1}{\mathcal{N}} 
\tilde{\mathfrak{R}}^{n,k}_{1}(\xi_{j}) \, \md \xi_{j} \quad \Rightarrow
\end{align*}
\begin{equation} \label{eqlmrtz20} 
\operatorname{Prob}(\tilde{\mathfrak{D}}_{\mathcal{N}}^{c}
(\tilde{M}_{0})) \leqslant \int_{\lbrace \lvert \tau \rvert > 
\tilde{M}_{0} \rbrace} \tilde{\mathfrak{R}}^{n,k}_{1}(\tau) 
\, \md \tau \! + \! \sum_{\substack{q=1\\q \neq \mathfrak{s}
-1}}^{\mathfrak{s}} \int_{\lbrace \lvert \tau - \alpha_{p_{q}} 
\rvert < \tilde{M}_{0}^{-1} \rbrace} \tilde{\mathfrak{R}}^{n,
k}_{1}(\tau) \, \md \tau.
\end{equation}
Via Equation~\eqref{eqlmrtz19}, one shows, by mimicking 
the analogue of the calculations leading to the 
Estimate~\eqref{eqlmrtz14}, that, for $n \! \in \! \mathbb{N}$ 
and $k \! \in \! \lbrace 1,2,\dotsc,K \rbrace$ such that 
$\alpha_{p_{\mathfrak{s}}} \! := \! \alpha_{k} \! \neq \! \infty$, 
after an integration argument,
\begin{gather*}
\int_{\lbrace \lvert \tau \rvert > \tilde{M}_{0} \rbrace} \tilde{
\mathfrak{R}}^{n,k}_{1}(\tau) \, \md \tau \underset{\underset{
z_{o}=1+o(1)}{\mathscr{N},n \to \infty}}{\leqslant} \dfrac{2 
\me^{\frac{n}{\mathcal{N}}(\mathcal{N}-1) \mathfrak{c}_{57}(n,
k,z_{o})} \tilde{M}_{0}^{-n(2(1+\tilde{c}_{\infty}-\gamma_{
i(\mathfrak{s}-1)_{k_{\mathfrak{s}-1}}})+ \mathcal{O}(n^{-1}))}
(1 \! + \! o(1))}{n(2(1 \! + \! \tilde{c}_{\infty} \! - \! \gamma_{
i(\mathfrak{s}-1)_{k_{\mathfrak{s}-1}}}) \! + \! \mathcal{O}(n^{-1}))}, \\
\sum_{\substack{q=1\\q \neq \mathfrak{s}-1}}^{\mathfrak{s}} 
\int_{\lbrace \lvert \tau - \alpha_{p_{q}} \rvert < \tilde{M}_{
0}^{-1} \rbrace} \tilde{\mathfrak{R}}^{n,k}_{1}(\tau) \, \md \tau 
\underset{\underset{z_{o}=1+o(1)}{\mathscr{N},n \to \infty}}{
\leqslant} \sum_{\substack{q=1\\q \neq \mathfrak{s}-1}}^{
\mathfrak{s}} \dfrac{\mathfrak{c}_{58}(n,k,z_{o},q) \me^{\frac{n}{
\mathcal{N}}(\mathcal{N}-1) \mathfrak{c}_{59}(n,k,z_{o},q)} 
\tilde{M}_{0}^{-n(2(1+\tilde{c}_{q}-\gamma_{i(q)_{k_{q}}})+ 
\mathcal{O}(n^{-1}))}(1 \! + \! o(1))}{n(2(1 \! + \! \tilde{c}_{q} 
\! - \! \gamma_{i(q)_{k_{q}}}) \! + \! \mathcal{O}(n^{-1}))},
\end{gather*}
where $\mathfrak{c}_{57}(n,k,z_{o}) \! =_{\underset{z_{o}=1
+o(1)}{\mathscr{N},n \to \infty}} \! \mathcal{O}(1)$, $1 \! + 
\! \tilde{c}_{\infty} \! > \! \gamma_{i(\mathfrak{s}-1)_{k_{
\mathfrak{s}-1}}}$, $\mathfrak{c}_{m}(n,k,z_{o},q) \! =_{
\underset{z_{o}=1+o(1)}{\mathscr{N},n \to \infty}} \! \mathcal{O}
(1)$, $m \! = \! 58,59$, and $1 \! + \! \tilde{c}_{q} \! > \! 
\gamma_{i(q)_{k_{q}}}$, $q \! = \! 1,\dotsc,\mathfrak{s} \! - 
\! 2,\mathfrak{s}$: taking, if necessary, $\tilde{M}_{0}$ even 
larger so that
\begin{equation*}
\tilde{M}_{0} \! \gg \! \tilde{M}_{0}^{\sharp} \! := \! \max 
\left\lbrace \me^{\frac{2 \mathfrak{c}_{57}(n,k,z_{o})+1}{2
(1+\tilde{c}_{\infty}-\gamma_{i(\mathfrak{s}-1)_{k_{\mathfrak{s}
-1}}})}},\left\lbrace \me^{\frac{2 \mathfrak{c}_{59}(n,k,z_{o},q)+
1}{2(1+\tilde{c}_{q}-\gamma_{i(q)_{k_{q}}})}} \right\rbrace_{
\substack{q=1\\q \neq \mathfrak{s}-1}}^{\mathfrak{s}},K 
\left(1 \! + \! \max_{q=1,\dotsc,\mathfrak{s}-2,\mathfrak{s}} 
\lbrace \lvert \alpha_{q} \rvert \rbrace \! + \! 3 \left(\min_{i 
\neq j \in \lbrace 1,\dotsc,\mathfrak{s}-2,\mathfrak{s} \rbrace} 
\lbrace \lvert \alpha_{p_{i}} \! - \! \alpha_{p_{j}} \rvert 
\rbrace \right)^{-1} \right) \right\rbrace,
\end{equation*}
one shows, via the latter two estimates and 
Equation~\eqref{eqlmrtz20}, that, for $n \! \in \! \mathbb{N}$ 
and $k \! \in \! \lbrace 1,2,\dotsc,K \rbrace$ such that 
$\alpha_{p_{\mathfrak{s}}} \! := \! \alpha_{k} \! \neq \! \infty$,
\begin{equation} \label{eqlmrtz21} 
\operatorname{Prob}(\tilde{\mathfrak{D}}_{\mathcal{N}}^{c}
(\tilde{M}_{0})) \underset{\underset{z_{o}=1+o(1)}{\mathscr{N},
n \to \infty}}{\leqslant} \dfrac{\mathfrak{c}_{60}(n,k,z_{o})}{n} 
\me^{-\frac{n}{\mathcal{N}}(\mathcal{N}-1)} \me^{-\frac{n}{
\mathcal{N}}(\mathcal{N}-1) \mathfrak{c}_{61}(n,k,z_{o})}
(1 \! + \! o(1)),
\end{equation}
where $\mathfrak{c}_{60}(n,k,z_{o}) \! =_{\underset{z_{o}=
1+o(1)}{\mathscr{N},n \to \infty}} \! \mathcal{O}(1)$, 
and $\mathfrak{c}_{61}(n,k,z_{o})$ $(:= \! \min \lbrace 
\mathfrak{c}_{57}(n,k,z_{o}),\min_{q=1,\dotsc,\mathfrak{s}
-2,\mathfrak{s}} \lbrace \mathfrak{c}_{59}(n,k,z_{o},q) \rbrace 
\rbrace)$ $=_{\underset{z_{o}=1+o(1)}{\mathscr{N},n \to 
\infty}} \! \mathcal{O}(1)$. For $n \! \in \! \mathbb{N}$ and 
$k \! \in \! \lbrace 1,2,\dotsc,K \rbrace$ such that $\alpha_{
p_{\mathfrak{s}}} \! := \! \alpha_{k} \! \neq \! \infty$, one 
shows, via the decomposition $\sum_{q=1}^{\mathfrak{s}-2} 
\varkappa_{nk \tilde{k}_{q}} \! + \! \varkappa_{nk \tilde{k}_{
\mathfrak{s}-1}}^{\infty} \! + \! \varkappa_{nk} \! = \! (n \! 
- \! 1)K \! + \! k \! = \! \mathcal{N}$, Lemma~\ref{lem2.2}, 
and Equations~\eqref{eqlmrtzp} and~\eqref{eqlmrtzz}, 
that (temporarily changing $z$ to $\lambda$)
\begin{equation} \label{eqlmrtz22} 
\pmb{\pi}^{n}_{k}(\lambda) \! = \! \tilde{\mathcal{
E}}^{n}_{k} \left(\lambda^{\varkappa_{nk \tilde{k}_{\mathfrak{s}
-1}}^{\infty}} \tilde{\mathbb{P}}(x_{1},x_{2},\dotsc,
x_{\mathcal{N}};\lambda) \right),
\end{equation}
where
\begin{equation} \label{eqlmrtz23} 
\tilde{\mathbb{P}}(x_{1},x_{2},\dotsc,x_{\mathcal{N}};
\lambda) \! = \! \tilde{\mathbb{P}}_{0}(x_{1},x_{2},\dotsc,
x_{\mathcal{N}};\lambda) \tilde{\mathbb{P}}_{1}(x_{1},x_{2},
\dotsc,x_{\mathcal{N}};\lambda),
\end{equation}
with
\begin{align*}
\tilde{\mathbb{P}}_{0}(x_{1},x_{2},\dotsc,x_{\mathcal{N}};\lambda) 
:=& \, \prod_{j=1}^{\mathcal{N}} \left(1 \! - \! \dfrac{x_{j}}{
\lambda} \right)^{\frac{\varkappa_{nk \tilde{k}_{\mathfrak{s}-
1}}^{\infty}+1}{\mathcal{N}}} \left(\dfrac{1}{\frac{x_{j}-\alpha_{k}}{
\lambda}} \! - \! \dfrac{1}{1 \! - \! \frac{\alpha_{k}}{\lambda}} 
\right)^{\frac{\varkappa_{nk}-1}{\mathcal{N}}} \left(\dfrac{x_{j} 
\! - \! \alpha_{k}}{\lambda} \right)^{\frac{\varkappa_{nk}-1}{
\mathcal{N}}} \left(1 \! - \! \dfrac{\alpha_{k}}{\lambda} \right)^{-
\frac{1}{\mathcal{N}}} \\
\times&\, \prod_{q=1}^{\mathfrak{s}-2} \left(\dfrac{1}{
\frac{x_{j}-\alpha_{p_{q}}}{\lambda}} \! - \! \dfrac{1}{1 \! - \! 
\frac{\alpha_{p_{q}}}{\lambda}} \right)^{\frac{\varkappa_{nk 
\tilde{k}_{q}}}{\mathcal{N}}} \left(\dfrac{x_{j} \! - \! \alpha_{
p_{q}}}{\lambda} \right)^{\frac{\varkappa_{nk \tilde{k}_{q}}}{
\mathcal{N}}}(\alpha_{k} \! - \! \alpha_{p_{q}})^{\frac{
\varkappa_{nk \tilde{k}_{q}}}{\mathcal{N}}},
\end{align*}
and
\begin{equation*}
\tilde{\mathbb{P}}_{1}(x_{1},x_{2},\dotsc,x_{\mathcal{N}};\lambda) 
\! := \! \left(\prod_{j=1}^{\mathcal{N}}(\alpha_{k} \! - \! x_{j})^{
\frac{\varkappa_{nk \tilde{k}_{\mathfrak{s}-1}}^{\infty}+1}{
\mathcal{N}}}(\alpha_{k} \! - \! x_{j})^{\frac{\varkappa_{nk}-1}{
\mathcal{N}}} \prod_{q=1}^{\mathfrak{s}-2}(\alpha_{k} \! - \! 
x_{j})^{\frac{\varkappa_{nk \tilde{k}_{q}}}{\mathcal{N}}} \right)^{-1}.
\end{equation*}
Via the inequality $\lvert y_{1} \! - \! y_{2} \rvert^{2} \! \leqslant 
\! (1 \! + \! y_{1}^{2})(1 \! + \! y_{2}^{2})$, $y_{1},y_{2} \! \in \! 
\mathbb{R}$, a calculation shows that
\begin{align*}
\lvert \tilde{\mathbb{P}}_{0}(x_{1},x_{2},\dotsc,x_{\mathcal{N}};
\lambda) \rvert^{2} \leqslant& \, 2^{\varkappa_{nk \tilde{k}_{
\mathfrak{s}-1}}^{\infty}+2}(1 \! + \! \alpha_{k}^{2})^{\sum_{q
=1}^{\mathfrak{s}-2} \varkappa_{nk \tilde{k}_{q}}} \left(1 \! 
+ \! \left\lvert \dfrac{\alpha_{k}}{\lambda} \right\rvert^{2} 
\right)^{\varkappa_{nk}-1} \left(1 \! + \! \left\lvert 1 \! - \! 
\dfrac{\alpha_{k}}{\lambda} \right\rvert^{-2} \right)^{
\varkappa_{nk}} \\
\times& \, \prod_{q=1}^{\mathfrak{s}-2}(1 \! + \! \alpha_{
p_{q}}^{2})^{\varkappa_{nk \tilde{k}_{q}}} \left(1 \! + \! \left\lvert 
\dfrac{\alpha_{p_{q}}}{\lambda} \right\rvert^{2} \right)^{
\varkappa_{nk \tilde{k}_{q}}} \left(1 \! + \! \left\lvert 1 \! - \! 
\dfrac{\alpha_{p_{q}}}{\lambda} \right\rvert^{-2} \right)^{
\varkappa_{nk \tilde{k}_{q}}} \\
\times& \, \prod_{j=1}^{\mathcal{N}} \left(1 \! + \! \left\lvert 
\dfrac{x_{j}}{\lambda} \right\rvert^{2} \right) \left(1 \! + \! 
\left\lvert \dfrac{x_{j} \! - \! \alpha_{k}}{\lambda} \right\rvert^{
-2} \right)^{\frac{\varkappa_{nk}-1}{\mathcal{N}}} \prod_{q=
1}^{\mathfrak{s}-2} \left(1 \! + \! \left\lvert \dfrac{x_{j} \! - 
\! \alpha_{p_{q}}}{\lambda} \right\rvert^{-2} \right)^{\frac{
\varkappa_{nk \tilde{k}_{q}}}{\mathcal{N}}}.
\end{align*}
Now, choose $\lambda$ real and large enough, that is, $\lvert 
\lambda \rvert$ $(> \! \tilde{M}_{0}$ $(\gg \! \tilde{M}_{0}^{
\sharp}))$ $\gg \! 1$, so that $\lvert \alpha_{p_{q}}/\lambda \rvert 
\! \leqslant \! \tilde{\Delta}(q)$, $q \! = \! 1,\dotsc,\mathfrak{s} 
\! - \! 2,\mathfrak{s}$, where $\tilde{\Delta}(q) \! := \! \lvert 
\alpha_{p_{q}} \rvert (K(1 \! + \! \max_{m=1,\dotsc,\mathfrak{s}
-2,\mathfrak{s}} \lbrace \lvert \alpha_{p_{m}} \rvert \rbrace 
\! + \! 3(\min_{i \neq j \in \lbrace 1,\dotsc,\mathfrak{s}-2,
\mathfrak{s} \rbrace} \lbrace \lvert \alpha_{p_{i}} \! - \! 
\alpha_{p_{j}} \rvert \rbrace)^{-1}))^{-1}$ $(\in \! 
(0,1/2))$,\footnote{A more precise analysis shows, in fact,  
that $(0,1/3) \! \ni \! \tilde{\Delta}(q)$, $q \! = \! 1,\dotsc,
\mathfrak{s} \! - \! 2,\mathfrak{s}$; however, the open interval 
$(0,1/2)$ is sufficient for the purposes of this proof.} and 
$\lvert (x_{j} \! - \! \alpha_{p_{q}})/\lambda \rvert \! \geqslant 
\! \tilde{\mathfrak{d}}(q)$, $j \! = \! 1,2,\dotsc,\mathcal{N}$, 
$q \! = \! 1,\dotsc,\mathfrak{s} \! - \! 2,\mathfrak{s}$, where 
$\tilde{\mathfrak{d}}(q) \! := \! 4(1 \! - \! \tilde{\Delta}(q))
((\tilde{M}_{0} \! + \! 1)(2(\tilde{M}_{0} \! + \! 1) \! + \! 
3(\mathfrak{s} \! - \! 1)(1 \! + \! \max_{m=1,\dotsc,\mathfrak{s}
-2,\mathfrak{s}} \lbrace \lvert \alpha_{p_{m}} \rvert \rbrace)
(\min_{j=1,\dotsc,\mathfrak{s}-2,\mathfrak{s}} \lbrace 
\tilde{\Delta}(j) \rbrace)^{-1}))^{-1}$.\footnote{Upon writing $\lvert 
x_{j} \! - \! \alpha_{p_{q}} \rvert \! = \! \lvert \lambda (\tfrac{x_{j}-
\alpha_{p_{q}}}{\lambda}) \rvert$, $j \! = \! 1,2,\dotsc,\mathcal{N}$, 
$q \! = \! 1,\dotsc,\mathfrak{s} \! - \! 2,\mathfrak{s}$, one shows 
that  $\lvert x_{j} \! - \! \alpha_{p_{q}} \rvert \! \geqslant \! 4(1 \! - 
\! \tilde{\Delta}(q))(\tilde{M}_{0} \! + \! 1)^{-1} \! > \! 2(\tilde{M}_{0} 
\! + \! 1)^{-1} \! > \! (\tilde{M}_{0} \! + \! 1)^{-1}$, which will 
be important for the analysis below.} With this choice of 
$\lambda$,\footnote{Take, say, $\lvert \lambda \rvert \! \geqslant 
\! 2(\tilde{M}_{0} \! + \! 1) \! + \! 3(\mathfrak{s} \! - \! 1)(1 \! + 
\! \max_{q=1,\dotsc,\mathfrak{s}-2,\mathfrak{s}} \lbrace \lvert 
\alpha_{p_{q}} \rvert \rbrace)(\min_{q=1,\dotsc,\mathfrak{s}-2,
\mathfrak{s}} \lbrace \tilde{\Delta}(q) \rbrace)^{-1})$.} a calculation 
shows that, for $n \! \in \! \mathbb{N}$ and $k \! \in \! \lbrace 1,
2,\dotsc,K \rbrace$ such that $\alpha_{p_{\mathfrak{s}}} \! := \! 
\alpha_{k} \! \neq \! \infty$,
\begin{align*}
\lvert \tilde{\mathbb{P}}_{0}(x_{1},x_{2},\dotsc,x_{\mathcal{N}};
\lambda) \rvert^{2} \underset{\underset{z_{o}=1+o(1)}{
\mathscr{N},n \to \infty}}{\leqslant}& \, \mathfrak{c}_{62}
(n,k,z_{o}) \me^{\frac{n}{\mathcal{N}}(\mathcal{N}-1) 
\mathfrak{c}_{63}(n,k,z_{o})} \underbrace{\left(\dfrac{1 \! + \! 
(\tilde{\mathfrak{d}}(\mathfrak{s}))^{-2}}{1 \! + \! (\tilde{M}_{0} 
\tilde{\mathfrak{d}}(\mathfrak{s}))^{-2}} \right)^{\varkappa_{nk}
-1}}_{\leqslant_{\underset{z_{o}=1+o(1)}{\mathscr{N},n 
\to \infty}} \, \mathfrak{c}_{64}(n,k,z_{o}) \me^{\frac{n}{
\mathcal{N}}(\mathcal{N}-1) \mathfrak{c}_{65}(n,k,z_{o})}} \, \, 
\underbrace{\prod_{q=1}^{\mathfrak{s}-2} \left(\dfrac{1 \! + 
\! (\tilde{\mathfrak{d}}(q))^{-2}}{1 \! + \! (\tilde{M}_{0} \tilde{
\mathfrak{d}}(q))^{-2}} \right)^{\varkappa_{nk \tilde{k}_{q}}}}_{
\leqslant_{\underset{z_{o}=1+o(1)}{\mathscr{N},n \to 
\infty}} \, \mathfrak{c}_{66}(n,k,z_{o}) \me^{\frac{n}{
\mathcal{N}}(\mathcal{N}-1) \mathfrak{c}_{67}(n,k,z_{o})}} \\
\times& \, \prod_{j=1}^{\mathcal{N}}(1 \! + \! x_{j}^{2})
(1 \! + \! (x_{j} \! - \! \alpha_{k})^{-2})^{\frac{\varkappa_{nk}
-1}{\mathcal{N}}} \prod_{q=1}^{\mathfrak{s}-2}(1 \! + \! (x_{j} 
\! - \! \alpha_{p_{q}})^{-2})^{\frac{\varkappa_{nk \tilde{k}_{q}}}{
\mathcal{N}}},
\end{align*}
where $\mathfrak{c}_{m}(n,k,z_{o}) \! =_{\underset{z_{o}=1+o(1)}{
\mathscr{N},n \to \infty}} \! \mathcal{O}(1)$, $m \! = \! 62,63,64,
65,66,67$, \footnote{It is instructive to note that $\mathfrak{c}_{63}
(n,k,z_{o}) \! = \! \ln (2) \gamma_{i(\mathfrak{s}-1)_{k_{\mathfrak{s}-1}}} 
\! + \! \gamma_{k} \ln ((1 \! + \! (\tilde{\Delta}(\mathfrak{s}))^{2})
(1 \! + \! (1 \! - \! \tilde{\Delta}(\mathfrak{s}))^{-2})) \! + 
\! \sum_{q=1}^{\mathfrak{s}-2} \gamma_{i(q)_{k_{q}}} \ln 
((1 \! + \! \alpha_{k}^{2})(1 \! + \! \alpha_{p_{q}}^{2})(1 \! + \! 
(\tilde{\Delta}(q))^{2})(1 \! + \! (1 \! - \! \tilde{\Delta}(q))^{-2}))$.} whence
\begin{equation} \label{eqlmrtz24} 
\lvert \tilde{\mathbb{P}}_{0}(x_{1},x_{2},\dotsc,x_{\mathcal{N}};
\lambda) \rvert^{2} \underset{\underset{z_{o}=1+o(1)}{
\mathscr{N},n \to \infty}}{\leqslant} \mathfrak{c}_{68}
(n,k,z_{o}) \me^{\frac{n}{\mathcal{N}}(\mathcal{N}-1) 
\mathfrak{c}_{69}(n,k,z_{o})} \prod_{j=1}^{\mathcal{N}}(1 \! 
+ \! x_{j}^{2})(1 \! + \! (x_{j} \! - \! \alpha_{k})^{-2})^{\frac{
\varkappa_{nk}-1}{\mathcal{N}}} \prod_{q=1}^{\mathfrak{s}-2}
(1 \! + \! (x_{j} \! - \! \alpha_{p_{q}})^{-2})^{\frac{\varkappa_{n
k \tilde{k}_{q}}}{\mathcal{N}}},
\end{equation}
where $\mathfrak{c}_{m}(n,k,z_{o}) \! =_{\underset{z_{o}=1+
o(1)}{\mathscr{N},n \to \infty}} \! \mathcal{O}(1)$, $m \! = \! 
68,69$, and
\begin{equation} \label{eqlmrtz25} 
\lvert \tilde{\mathbb{P}}_{1}(x_{1},x_{2},\dotsc,x_{\mathcal{N}};
\lambda) \rvert^{2} \underset{\underset{z_{o}=1+o(1)}{
\mathscr{N},n \to \infty}}{\leqslant} \mathfrak{c}_{70}(n,k,z_{o}) 
\me^{\frac{n}{\mathcal{N}}(\mathcal{N}-1) \mathfrak{c}_{71}
(n,k,z_{o})},
\end{equation}
where $\mathfrak{c}_{m}(n,k,z_{o}) \! =_{\underset{z_{o}
=1+o(1)}{\mathscr{N},n \to \infty}} \! \mathcal{O}(1)$, 
$m \! = \! 70,71$; hence, via 
Equations~\eqref{eqlmrtz23}--\eqref{eqlmrtz25}, one arrives 
at, for $n \! \in \! \mathbb{N}$ and $k \! \in \! \lbrace 1,2,
\dotsc,K \rbrace$ such that $\alpha_{p_{\mathfrak{s}}} \! 
:= \! \alpha_{k} \! \neq \! \infty$,
\begin{equation} \label{eqlmrtz26} 
\lvert \tilde{\mathbb{P}}(x_{1},x_{2},\dotsc,x_{\mathcal{N}};
\lambda) \rvert^{2} \underset{\underset{z_{o}=1+o(1)}{
\mathscr{N},n \to \infty}}{\leqslant} \mathfrak{c}_{72}
(n,k,z_{o}) \me^{\frac{n}{\mathcal{N}}(\mathcal{N}-1) 
\mathfrak{c}_{73}(n,k,z_{o})} \prod_{j=1}^{\mathcal{N}}(1 \! 
+ \! x_{j}^{2})(1 \! + \! (x_{j} \! - \! \alpha_{k})^{-2})^{\frac{
\varkappa_{nk}-1}{\mathcal{N}}} \prod_{q=1}^{\mathfrak{s}-2}
(1 \! + \! (x_{j} \! - \! \alpha_{p_{q}})^{-2})^{\frac{\varkappa_{n
k \tilde{k}_{q}}}{\mathcal{N}}},
\end{equation}
where $\mathfrak{c}_{m}(n,k,z_{o}) \! =_{\underset{z_{o}
=1+o(1)}{\mathscr{N},n \to \infty}} \! \mathcal{O}(1)$, 
$m \! = \! 72,73$. Recalling the definition of the 
$\mathcal{N}$-dimensional sets $\tilde{\mathfrak{D}}_{
\mathcal{N}}(\tilde{M}_{0})$ and $\tilde{\mathfrak{D}}^{c}_{
\mathcal{N}}(\tilde{M}_{0})$ (with $\tilde{M}_{0}$ as refined 
above), it follows {}from Equation~\eqref{eqlmrtz22} and the 
linearity of $\tilde{\mathcal{E}}^{n}_{k}$ that
\begin{equation} \label{eqlmrtz27} 
\pmb{\pi}^{n}_{k}(\lambda) \! = \! \lambda^{
\varkappa_{nk \tilde{k}_{\mathfrak{s}-1}}^{\infty}} \left(
\tilde{\mathcal{E}}^{n}_{k} \left(\chi_{\tilde{\mathfrak{D}}_{
\mathcal{N}}(\tilde{M}_{0})}(\vec{x}) \tilde{\mathbb{P}}
(x_{1},x_{2},\dotsc,x_{\mathcal{N}};\lambda) \right) \! + \! 
\tilde{\mathcal{E}}^{n}_{k} \left(\chi_{\tilde{\mathfrak{D}}^{
c}_{\mathcal{N}}(\tilde{M}_{0})}(\vec{x}) \tilde{\mathbb{P}}
(x_{1},x_{2},\dotsc,x_{\mathcal{N}};\lambda) \right) \right):
\end{equation}
one now proceeds to analyse, in the double-scaling limit 
$\mathscr{N},n \! \to \! \infty$ such that $z_{o} \! = \! 1 \! 
+ \! o(1)$, the term $\tilde{\mathcal{E}}^{n}_{k}(\chi_{\tilde{
\mathfrak{D}}^{c}_{\mathcal{N}}(\tilde{M}_{0})}(\vec{x}) 
\linebreak[4] 
\cdot \tilde{\mathbb{P}}(x_{1},x_{2},\dotsc,x_{\mathcal{N}};
\lambda))$. An application of the Cauchy-Schwarz Inequality 
shows that
\begin{align*}
\left\lvert \tilde{\mathcal{E}}^{n}_{k} \left(\chi_{\tilde{
\mathfrak{D}}^{c}_{\mathcal{N}}(\tilde{M}_{0})}(\vec{x}) 
\tilde{\mathbb{P}}(x_{1},x_{2},\dotsc,x_{\mathcal{N}};
\lambda) \right) \right\vert^{2} \leqslant& \, \underbrace{
\tilde{\mathcal{E}}^{n}_{k} \left(\chi_{\tilde{\mathfrak{
D}}^{c}_{\mathcal{N}}(\tilde{M}_{0})}(\vec{x}) \right)}_{= \, 
\operatorname{Prob}(\tilde{\mathfrak{D}}^{c}_{\mathcal{N}}
(\tilde{M}_{0}))} \idotsint\limits_{\mathbb{R}^{\mathcal{N}}} 
\chi_{\tilde{\mathfrak{D}}^{c}_{\mathcal{N}}(\tilde{M}_{0})}
(\vec{\xi}) \lvert \tilde{\mathbb{P}}(\xi_{1},\xi_{2},\dotsc,
\xi_{\mathcal{N}};\lambda) \rvert^{2} \\
\times& \, \tilde{\mathcal{P}}^{n}_{k}(\xi_{1},\xi_{2},\dotsc,
\xi_{\mathcal{N}}) \, \md \xi_{1} \, \md \xi_{2} \, \dotsb 
\, \md \xi_{\mathcal{N}};
\end{align*}
hence, via the Estimates~\eqref{eqlmrtz21} 
and~\eqref{eqlmrtz26}, it follows that
\begin{align} \label{eqlmrtz28} 
\left\lvert \tilde{\mathcal{E}}^{n}_{k} \left(\chi_{\tilde{
\mathfrak{D}}^{c}_{\mathcal{N}}(\tilde{M}_{0})}(\vec{x}) 
\tilde{\mathbb{P}}(x_{1},x_{2},\dotsc,x_{\mathcal{N}};
\lambda) \right) \right\vert^{2} \underset{\underset{
z_{o}=1+o(1)}{\mathscr{N},n \to \infty}}{\leqslant}& \, 
\dfrac{\mathfrak{c}_{74}(n,k,z_{o})}{n} \me^{-\frac{n}{
\mathcal{N}}(\mathcal{N}-1)} \me^{-\frac{n}{\mathcal{N}}
(\mathcal{N}-1) \mathfrak{c}_{61}(n,k,z_{o})} \me^{
\frac{n}{\mathcal{N}}(\mathcal{N}-1) \mathfrak{c}_{69}
(n,k,z_{o})}(1 \! + \! o(1)) \nonumber \\
\times& \, \tilde{\mathcal{E}}^{n}_{k} \left(\chi_{\tilde{
\mathfrak{D}}^{c}_{\mathcal{N}}(\tilde{M}_{0})}(\vec{x}) 
\prod_{j=1}^{\mathcal{N}}(1 \! + \! x_{j}^{2})(1 \! + \! 
(x_{j} \! - \! \alpha_{k})^{-2})^{\frac{\varkappa_{nk}-1}{
\mathcal{N}}} \prod_{q=1}^{\mathfrak{s}-2}(1 \! + \! 
(x_{j} \! - \! \alpha_{p_{q}})^{-2})^{\frac{\varkappa_{nk 
\tilde{k}_{q}}}{\mathcal{N}}} \right),
\end{align}
where $\mathfrak{c}_{74}(n,k,z_{o}) \! =_{\underset{z_{o}=1+
o(1)}{\mathscr{N},n \to \infty}} \! \mathcal{O}(1)$. Define the 
following $\mathcal{N}$-dimensional subsets of $\mathbb{R}^{
\mathcal{N}}$ (cf. the $(\mathcal{N} \! - \! 1)$-dimensional sets 
$\mathbb{S}^{n,\sharp}_{k,j}$, $j \! = \! 1,2$, defined heretofore):
\begin{gather*}
\tilde{\mathbb{S}}^{n,\sharp}_{k,1} \! := \! \left\lbrace 
\mathstrut (x_{1},x_{2},\dotsc,x_{\mathcal{N}}) \! \in \! 
\mathbb{R}^{\mathcal{N}}; \, \prod_{j=1}^{\mathcal{N}}
(1 \! + \! x_{j}^{2})(1 \! + \! (x_{j} \! - \! \alpha_{k})^{-2})^{
\frac{\varkappa_{nk}-1}{\mathcal{N}}} \prod_{q=1}^{
\mathfrak{s}-2}(1 \! + \! (x_{j} \! - \! \alpha_{p_{q}})^{-2})^{
\frac{\varkappa_{nk \tilde{k}_{q}}}{\mathcal{N}}} \underset{
\underset{z_{o}=1+o(1)}{\mathscr{N},n \to \infty}}{\geqslant} 
\me^{\mathcal{N} \mathfrak{c}_{75}(n,k,z_{o})} \right\rbrace, \\
\tilde{\mathbb{S}}^{n,\sharp}_{k,2} \! := \! \left\lbrace 
\mathstrut (x_{1},x_{2},\dotsc,x_{\mathcal{N}}) \! \in \! 
\mathbb{R}^{\mathcal{N}}; \, \prod_{j=1}^{\mathcal{N}}
(1 \! + \! x_{j}^{2})(1 \! + \! (x_{j} \! - \! \alpha_{k})^{-2})^{
\frac{\varkappa_{nk}-1}{\mathcal{N}}} \prod_{q=1}^{\mathfrak{s}-2}
(1 \! + \! (x_{j} \! - \! \alpha_{p_{q}})^{-2})^{\frac{\varkappa_{nk 
\tilde{k}_{q}}}{\mathcal{N}}} \underset{\underset{z_{o}=1+o(1)}{
\mathscr{N},n \to \infty}}{\leqslant} \me^{\mathcal{N} 
\mathfrak{c}_{75}(n,k,z_{o})} \right\rbrace,
\end{gather*}
where $\mathfrak{c}_{75}(n,k,z_{o}) \! =_{\underset{z_{o}
=1+o(1)}{\mathscr{N},n \to \infty}} \! \mathcal{O}(1)$ (note 
that $\tilde{\mathbb{S}}^{n,\sharp}_{k,1}$ and $\tilde{
\mathbb{S}}^{n,\sharp}_{k,2}$ may have a non-empty 
intersection): proceeding, now, via the 
$\mathcal{N}$-dimensional analogue of the calculations 
subsumed in the derivation of the Estimate~\eqref{eqlmrtz19}, 
one shows that
\begin{align} \label{eqlmrtz29} 
& \, \tilde{\mathcal{E}}^{n}_{k} \left(\chi_{\tilde{\mathfrak{D}}^{
c}_{\mathcal{N}}(\tilde{M}_{0})}(\vec{x}) \prod_{j=1}^{\mathcal{
N}}(1 \! + \! x_{j}^{2})(1 \! + \! (x_{j} \! - \! \alpha_{k})^{-2})^{
\frac{\varkappa_{nk}-1}{\mathcal{N}}} \prod_{q=1}^{\mathfrak{s}
-2}(1 \! + \! (x_{j} \! - \! \alpha_{p_{q}})^{-2})^{\frac{\varkappa_{n
k \tilde{k}_{q}}}{\mathcal{N}}} \right) \nonumber \\
& \, \underset{\underset{z_{o}=1+o(1)}{\mathscr{N},n \to 
\infty}}{\leqslant} \underbrace{\tilde{\mathcal{E}}^{n}_{k} \left(
\chi_{\tilde{\mathbb{S}}^{n,\sharp}_{k,1}}(\vec{x}) \prod_{j
=1}^{\mathcal{N}}(1 \! + \! x_{j}^{2})(1 \! + \! (x_{j} \! - \! 
\alpha_{k})^{-2})^{\frac{\varkappa_{nk}-1}{\mathcal{N}}} 
\prod_{q=1}^{\mathfrak{s}-2}(1 \! + \! (x_{j} \! - \! \alpha_{
p_{q}})^{-2})^{\frac{\varkappa_{nk \tilde{k}_{q}}}{\mathcal{N}}} 
\right)}_{\leqslant_{\underset{z_{o}=1+o(1)}{\mathscr{N},n \to 
\infty}} \mathfrak{c}_{76}(n,k,z_{o}) \exp (-\frac{n}{\mathcal{N}} 
\mathcal{N}(\mathcal{N}-1) \mathfrak{c}_{77}(n,k,z_{o}))} 
\nonumber \\
& \, +\underbrace{\tilde{\mathcal{E}}^{n}_{k} \left(\chi_{\tilde{
\mathbb{S}}^{n,\sharp}_{k,2}}(\vec{x}) \prod_{j=1}^{\mathcal{N}}
(1 \! + \! x_{j}^{2})(1 \! + \! (x_{j} \! - \! \alpha_{k})^{-2})^{
\frac{\varkappa_{nk}-1}{\mathcal{N}}} \prod_{q=1}^{\mathfrak{s}
-2}(1 \! + \! (x_{j} \! - \! \alpha_{p_{q}})^{-2})^{\frac{\varkappa_{n
k \tilde{k}_{q}}}{\mathcal{N}}} \right)}_{\leqslant_{\underset{z_{o}
=1+o(1)}{\mathscr{N},n \to \infty}} \exp (\mathcal{N} 
\mathfrak{c}_{75}(n,k,z_{o}))} \nonumber \\
& \, \underset{\underset{z_{o}=1+o(1)}{\mathscr{N},n \to 
\infty}}{\leqslant} \mathfrak{c}_{78}(n,k,z_{o}) \me^{\frac{n}{
\mathcal{N}}(\mathcal{N}-1) \mathfrak{c}_{79}(n,k,z_{o})} 
\left(1 \! + \! \mathcal{O} \left(\mathfrak{c}_{80}(n,k,z_{o}) 
\me^{-\frac{n}{\mathcal{N}} \mathcal{N}(\mathcal{N}-1) 
\mathfrak{c}_{77}(n,k,z_{o})} \right) \right),
\end{align}  
where $\mathfrak{c}_{m}(n,k,z_{o}) \! =_{\underset{z_{o}=1+
o(1)}{\mathscr{N},n \to \infty}} \! \mathcal{O}(1)$, $m \! = \! 
76,77,78,79,80$; therefore, adjusting all $\mathcal{O}(1)$ 
parameters so that $\exp (\tfrac{n}{\mathcal{N}}(\mathcal{N} 
\! - \! 1)(\mathfrak{c}_{69}(n,k,z_{o}) \! + \! \mathfrak{c}_{79}
(n,k,z_{o}) \! - \! \mathfrak{c}_{61}(n,k,z_{o}))) \! =_{\underset{
z_{o}=1+o(1)}{\mathscr{N},n \to \infty}} \! \mathcal{O}(1)$, 
one shows, via Equation~\eqref{eqlmrtz28} and the 
Estimate~\eqref{eqlmrtz29}, that, for $n \! \in \! \mathbb{N}$ 
and $k \! \in \! \lbrace 1,2,\dotsc,K \rbrace$ such that 
$\alpha_{p_{\mathfrak{s}}} \! := \! \alpha_{k} \! \neq \! \infty$,
\begin{equation*}
\left\lvert \tilde{\mathcal{E}}^{n}_{k} \left(\chi_{\tilde{
\mathfrak{D}}^{c}_{\mathcal{N}}(\tilde{M}_{0})}(\vec{x}) 
\tilde{\mathbb{P}}(x_{1},x_{2},\dotsc,x_{\mathcal{N}};
\lambda) \right) \right\vert^{2} \underset{\underset{
z_{o}=1+o(1)}{\mathscr{N},n \to \infty}}{\leqslant} 
\dfrac{\mathfrak{c}_{81}(n,k,z_{o})}{n} \me^{-\frac{n}{
\mathcal{N}}(\mathcal{N}-1)}(1 \! + \! o(1)),
\end{equation*}
where $\mathfrak{c}_{81}(n,k,z_{o}) \! =_{\underset{z_{o}
=1+o(1)}{\mathscr{N},n \to \infty}} \! \mathcal{O}(1)$, 
whence
\begin{equation} \label{eqlmrtz30} 
\left\lvert \tilde{\mathcal{E}}^{n}_{k} \left(\chi_{\tilde{
\mathfrak{D}}^{c}_{\mathcal{N}}(\tilde{M}_{0})}(\vec{x}) 
\tilde{\mathbb{P}}(x_{1},x_{2},\dotsc,x_{\mathcal{N}};
\lambda) \right) \right\vert \underset{\underset{z_{o}=
1+o(1)}{\mathscr{N},n \to \infty}}{\leqslant} \dfrac{
\mathfrak{c}_{82}(n,k,z_{o})}{\sqrt{n}} \me^{-\frac{1}{2} 
\frac{n}{\mathcal{N}}(\mathcal{N}-1)}(1 \! + \! o(1)),
\end{equation}
where $\mathfrak{c}_{82}(n,k,z_{o}) \! =_{\underset{z_{o}
=1+o(1)}{\mathscr{N},n \to \infty}} \! \mathcal{O}(1)$. For 
$n \! \in \! \mathbb{N}$ and $k \! \in \! \lbrace 1,2,\dotsc,
K \rbrace$ such that $\alpha_{p_{\mathfrak{s}}} \! := \! 
\alpha_{k} \! \neq \! \infty$, with the above-refined choice 
of $\tilde{M}_{0}$ $(\gg \! 1)$, define, analogously as above, 
the following sets: (i) $\tilde{\mathfrak{D}}(\tilde{M}_{0} \! 
+ \! 1) \! := \! [-(\tilde{M}_{0} \! + \! 1),\tilde{M}_{0} \! + \! 
1] \setminus \cup_{\underset{q \neq \mathfrak{s}-1}{q=1}}^{
\mathfrak{s}} \mathscr{O}_{\frac{1}{\tilde{M}_{0}+1}}(\alpha_{
p_{q}})$ $(\supset \tilde{\mathfrak{D}}(\tilde{M}_{0}) \supseteq 
J_{f})$, with $\mathscr{O}_{\frac{1}{\tilde{M}_{0}+1}}(\alpha_{
p_{i}}) \cap \mathscr{O}_{\frac{1}{\tilde{M}_{0}+1}}(\alpha_{
p_{j}}) \! = \! \varnothing$ $\forall$ $i \! \neq \! j \! \in \! 
\lbrace 1,\dotsc,\mathfrak{s} \! - \! 2,\mathfrak{s} \rbrace$; 
(ii) the complement of $\tilde{\mathfrak{D}}(\tilde{M}_{0} 
\! + \! 1)$ relative to $\mathbb{R}$, $\tilde{\mathfrak{D}}^{c}
(\tilde{M}_{0} \! + \! 1) \! := \! \mathbb{R} \setminus \tilde{
\mathfrak{D}}(\tilde{M}_{0} \! + \! 1) \! = \! \mathbb{A}
(\tilde{M}_{0} \! + \! 1) \cup \cup_{\underset{q \neq 
\mathfrak{s}-1}{q=1}}^{\mathfrak{s}} \mathbb{B}_{q}
(\tilde{M}_{0} \! + \! 1)$, where $\mathbb{A}(\tilde{M}_{0} 
\! + \! 1) \! := \! \lbrace \mathstrut x \! \in \! \mathbb{R}; 
\, \lvert x \rvert \! \geqslant \! \tilde{M}_{0} \! + \! 1 \rbrace$, 
and $\mathbb{B}_{q}(\tilde{M}_{0} \! + \! 1) \! := \! \lbrace 
\mathstrut x \! \in \! \mathbb{R}; \, \lvert x \! - \! \alpha_{
p_{q}} \rvert \! < \! (\tilde{M}_{0} \! + \! 1)^{-1} \rbrace$, 
$q \! = \! 1,\dotsc,\mathfrak{s} \! - \! 2,\mathfrak{s}$, with 
$\mathbb{A}(\tilde{M}_{0} \! + \! 1) \cap \mathbb{B}_{q}
(\tilde{M}_{0} \! + \! 1) \! = \! \varnothing$, $q \! = \! 1,\dotsc,
\mathfrak{s} \! - \! 2,\mathfrak{s}$, and $\mathbb{B}_{i}
(\tilde{M}_{0} \! + \! 1) \cap \mathbb{B}_{j}(\tilde{M}_{0} \! + 
\! 1) \! = \! \varnothing$ $\forall$ $i \! \neq \! j \! \in \! 
\lbrace 1,\dotsc,\mathfrak{s} \! - \! 2,\mathfrak{s} \rbrace$ 
$(\Rightarrow \! \mathbb{A}(\tilde{M}_{0} \! + \! 1) \cap 
(\cup_{\underset{q \neq \mathfrak{s}-1}{q=1}}^{\mathfrak{s}} 
\mathbb{B}_{q}(\tilde{M}_{0} \! + \! 1)) \! = \! \cup_{
\underset{q \neq \mathfrak{s}-1}{q=1}}^{\mathfrak{s}}
(\mathbb{A}(\tilde{M}_{0} \! + \! 1) \cap \mathbb{B}_{q}
(\tilde{M}_{0} \! + \! 1)) \! = \! \cup_{\underset{q \neq 
\mathfrak{s}-1}{q=1}}^{\mathfrak{s}} \varnothing \! = \! 
\varnothing)$; (iii) the $\mathcal{N}$-fold Cartesian product 
of $\tilde{\mathfrak{D}}(\tilde{M}_{0} \! + \! 1)$, $\tilde{
\mathfrak{D}}_{\mathcal{N}}(\tilde{M}_{0} \! + \! 1) \! := \! 
\tilde{\mathfrak{D}}(\tilde{M}_{0} \! + \! 1) \! \times \! 
\tilde{\mathfrak{D}}(\tilde{M}_{0} \! + \! 1) \! \times 
\! \dotsb \! \times \tilde{\mathfrak{D}}(\tilde{M}_{0} \! + \! 1)$ 
$(\subset \mathbb{R}^{\mathcal{N}})$; and (iv) the complement 
of $\tilde{\mathfrak{D}}_{\mathcal{N}}(\tilde{M}_{0} \! + \! 1)$ 
relative to $\mathbb{R}^{\mathcal{N}}$, $\tilde{\mathfrak{D}}_{
\mathcal{N}}^{c}(\tilde{M}_{0} \! + \! 1) \! := \! \mathbb{R}^{
\mathcal{N}} \setminus \tilde{\mathfrak{D}}_{\mathcal{N}}(\tilde{M}_{0} 
\! + \! 1)$.\footnote{Note the nested-intervals property of these 
sets, that is, $\mathbb{R} \supset \tilde{\mathfrak{D}}(\tilde{M}_{0} 
\! + \! 1) \supset \tilde{\mathfrak{D}}(\tilde{M}_{0}) 
\supseteq J_{f}$ and $\mathbb{R}^{\mathcal{N}} \supset 
\tilde{\mathfrak{D}}_{\mathcal{N}}(\tilde{M}_{0} \! + \! 1) \supset 
\tilde{\mathfrak{D}}_{\mathcal{N}}(\tilde{M}_{0})$.} Recalling 
{}from the above-refined $\tilde{M}_{0}$ $(\gg \! 1)$ and the 
choice $(\mathbb{R} \! \ni)$ $\lvert \lambda \rvert \! \geqslant 
\! 2(\tilde{M}_{0} \! + \! 1) \! + \! 3(\mathfrak{s} \! - \! 1)(1 \! 
+ \! \max_{q=1,\dotsc,\mathfrak{s}-2,\mathfrak{s}} \lbrace 
\lvert \alpha_{p_{q}} \rvert \rbrace)(\min_{q=1,\dotsc,
\mathfrak{s}-2,\mathfrak{s}} \lbrace \tilde{\Delta}(q) 
\rbrace)^{-1}$, say, that, for $x \! \in \! \tilde{\mathfrak{D}}_{
\mathcal{N}}(\tilde{M}_{0} \! + \! 1)$, $\lvert \alpha_{p_{q}}/
\lambda \rvert \! \leqslant \! \tilde{\Delta}(q)$ $(\in \! (0,1/2))$, 
$q=1,\dotsc,\mathfrak{s}-2,\mathfrak{s}$, and $\lvert (x_{j} \! 
- \! \alpha_{p_{q}})/\lambda \rvert \! \geqslant \! 4(1 \! - \! 
\tilde{\Delta}(q))((\tilde{M}_{0} \! + \! 1)(2(\tilde{M}_{0} \! + \! 
1) \! + \! 3(\mathfrak{s} \! - \! 1)(1 \! + \! \max_{q=1,\dotsc,
\mathfrak{s}-2,\mathfrak{s}} \lbrace \lvert \alpha_{p_{q}} \rvert 
\rbrace)(\min_{q=1,\dotsc,\mathfrak{s}-2,\mathfrak{s}} \lbrace 
\tilde{\Delta}(q) \rbrace)^{-1}))^{-1}$ $(\Rightarrow \! \lvert x_{j} 
\! - \! \alpha_{p_{q}} \rvert \! > \! (\tilde{M}_{0} \! + \! 1)^{-1})$, 
$j \! = \! 1,2,\dotsc,\mathcal{N}$, $q=1,\dotsc,\mathfrak{s}-2,
\mathfrak{s}$, it also follows that $\max_{j=1,2,\dotsc,\mathcal{
N}} \lbrace \lvert x_{j}/\lambda \rvert \rbrace \! \leqslant \! 1/2$. 
For $n \! \in \! \mathbb{N}$ and $k \! \in \! \lbrace 1,2,\dotsc,K 
\rbrace$ such that $\alpha_{p_{\mathfrak{s}}} \! := \! \alpha_{k} 
\! \neq \! \infty$, let $\boldsymbol{\mathrm{C}}^{\infty}_{0}
(\mathbb{R}) \! \ni \! \tilde{F} \colon \mathbb{R} \! \to \! [0,1]$ 
be a test function satisfying: (i) $0 \! \leqslant \! \tilde{F}(x) \! 
\leqslant \! 1$, $x \! \in \! \mathbb{R}$; (ii) $\tilde{F}(x) \! = \! 1$, 
$x \! \in \! \tilde{\mathfrak{D}}(\tilde{M}_{0})$ $(= \! [-\tilde{M}_{0},
\tilde{M}_{0}] \setminus \cup_{\underset{q \neq \mathfrak{s}-1}{q=
1}}^{\mathfrak{s}} \mathbb{B}_{q}(\tilde{M}_{0}))$; and (iii) $\tilde{F}
(x) \! = \! 0$, $x \! \in \! \tilde{\mathfrak{D}}^{c}(\tilde{M}_{0} 
\! + \! 1)$ $(= \! \mathbb{A}(\tilde{M}_{0} \! + \! 1) \cup 
\cup_{\underset{q \neq \mathfrak{s}-1}{q=1}}^{\mathfrak{s}} 
\mathbb{B}_{q}(\tilde{M}_{0} \! + \! 1) \subset \mathbb{R} \setminus 
\tilde{\mathfrak{D}}(\tilde{M}_{0} \! + \! 1))$. For $n \! \in \! 
\mathbb{N}$ and $k \! \in \! \lbrace 1,2,\dotsc,K \rbrace$ such that 
$\alpha_{p_{\mathfrak{s}}} \! := \! \alpha_{k} \! \neq \! \infty$, let
\begin{equation} \label{eqlmrtz31} 
\tilde{\mathfrak{g}}(x) \! = \! \tilde{F}(x) \ln \tilde{\mathbb{G}}(x),
\end{equation}
where
\begin{align} \label{eqlmrtz32} 
\tilde{\mathbb{G}}(x) :=& \, \left(1 \! - \! \dfrac{x}{\lambda} 
\right)^{\frac{\varkappa_{nk \tilde{k}_{\mathfrak{s}-1}}^{\infty}
+1}{\mathcal{N}}} \left(\dfrac{1 \! - \! \frac{x}{\lambda}}{(1 \! - 
\! \frac{\alpha_{k}}{\lambda})(\frac{x-\alpha_{k}}{\lambda})} 
\right)^{\frac{\varkappa_{nk}-1}{\mathcal{N}}} \prod_{q=1}^{
\mathfrak{s}-2} \left(\dfrac{1 \! - \! \frac{x}{\lambda}}{(1 \! - 
\! \frac{\alpha_{p_{q}}}{\lambda})(\frac{x-\alpha_{p_{q}}}{
\lambda})} \right)^{\frac{\varkappa_{nk \tilde{k}_{q}}}{\mathcal{N}}} 
\left(\dfrac{x \! - \! \alpha_{k}}{\lambda} \right)^{\frac{\varkappa_{
nk}-1}{\mathcal{N}}} \prod_{q=1}^{\mathfrak{s}-2} \left(\dfrac{x \! 
- \! \alpha_{p_{q}}}{\lambda} \right)^{\frac{\varkappa_{nk \tilde{k}_{
q}}}{\mathcal{N}}} \nonumber \\
\times& \, \dfrac{\left(1 \! - \! \frac{\alpha_{k}}{\lambda} \right)^{-
\frac{1}{\mathcal{N}}} \prod_{q=1}^{\mathfrak{s}-2}(\alpha_{k} \! - 
\! \alpha_{p_{q}})^{\frac{\varkappa_{nk \tilde{k}_{q}}}{\mathcal{N}}}}{
(\alpha_{k} \! - \! x)^{\frac{\varkappa_{nk \tilde{k}_{\mathfrak{s}-1}}^{
\infty}+1}{\mathcal{N}}}(\alpha_{k} \! - \! x)^{\frac{\varkappa_{nk}-
1}{\mathcal{N}}} \prod_{q=1}^{\mathfrak{s}-2}(\alpha_{k} \! - \! x)^{
\frac{\varkappa_{nk \tilde{k}_{q}}}{\mathcal{N}}}}.
\end{align}
A tedious analysis shows that
\begin{align} \label{eqlmrtz33} 
-\ln \left(\dfrac{4}{3}(\tilde{M}_{0} \! + \! 1)6^{\mathfrak{s}} \left(
\min_{q=1,2,\dotsc,\mathfrak{s}-2} \lbrace \lvert \alpha_{k} \! - 
\! \alpha_{p_{q}} \rvert \rbrace \right)^{\mathfrak{s}-2} \right)
(1 \! + \! \mathcal{O}(n^{-1})) \underset{\underset{z_{o}=1+
o(1)}{\mathscr{N},n \to \infty}}{\leqslant}& \ln \left(\dfrac{\frac{1}{8}
(\tilde{M}_{0} \! + \! 1) \left(\min\limits_{q=1,2,\dotsc,\mathfrak{s}
-2} \lbrace \lvert \alpha_{k} \! - \! \alpha_{p_{q}} \rvert 
\rbrace \right)^{\sum_{q=1}^{\mathfrak{s}-2} \frac{\varkappa_{nk 
\tilde{k}_{q}}}{\mathcal{N}}}}{(1 \! + \! \tilde{\Delta}(\mathfrak{s}))^{
\frac{\varkappa_{nk}}{\mathcal{N}}} \prod_{q=1}^{\mathfrak{s}-2}
(1 \! + \! \tilde{\Delta}(q))^{\frac{\varkappa_{nk \tilde{k}_{q}}}{
\mathcal{N}}}} \right) \nonumber \\
\underset{\underset{z_{o}=1+o(1)}{\mathscr{N},n \to \infty}}{
\leqslant} \tilde{\mathfrak{g}}(x) \underset{\underset{z_{o}=1+
o(1)}{\mathscr{N},n \to \infty}}{\leqslant} \ln \left(\dfrac{\frac{3}{4}
(\tilde{M}_{0} \! + \! 1) \left(\max\limits_{q=1,2,\dotsc,\mathfrak{s}
-2} \lbrace \lvert \alpha_{k} \! - \! \alpha_{p_{q}} \rvert 
\rbrace \right)^{\sum_{q=1}^{\mathfrak{s}-2} \frac{\varkappa_{nk 
\tilde{k}_{q}}}{\mathcal{N}}}}{(1 \! - \! \tilde{\Delta}(\mathfrak{s}))^{
\frac{\varkappa_{nk}}{\mathcal{N}}} \prod_{q=1}^{\mathfrak{s}-2}
(1 \! - \! \tilde{\Delta}(q))^{\frac{\varkappa_{nk \tilde{k}_{q}}}{
\mathcal{N}}}} \right) 
&\underset{\underset{z_{o}=1+o(1)}{\mathscr{N},n \to \infty}}{\leqslant} 
\ln \left(6(\tilde{M}_{0} \! + \! 1)2^{\mathfrak{s}} \left(\max_{q=1,2,
\dotsc,\mathfrak{s}-2} \lbrace \lvert \alpha_{k} \! - \! \alpha_{p_{q}} 
\rvert \rbrace \right)^{\mathfrak{s}-2} \right) \nonumber \\
& \, \times (1 \! + \! \mathcal{O}(n^{-1})).
\end{align}
Take, if necessary, the above-refined $\tilde{M}_{0}$ $(\gg \! 1)$ even 
larger so that, also, $(\tilde{M}_{0} \! + \! 1)(\min_{q=1,2,\dotsc,
\mathfrak{s}-2} \lbrace \lvert \alpha_{k} \! - \! \alpha_{p_{q}} \rvert 
\rbrace)^{\mathfrak{s}-2} \! := \! 1 \! + \! \mathfrak{z}_{m} \! \gg 
\! 1$, where $\mathfrak{z}_{m} \! = \! \mathfrak{z}_{m}(n,k,z_{o}) \! 
\gg_{\underset{z_{o}=1+o(1)}{\mathscr{N},n \to \infty}} \! 1$ and 
$\mathcal{O}(1)$.\footnote{Of course, this implies that $(\tilde{M}_{0} 
\! + \! 1)(\max_{q=1,2,\dotsc,\mathfrak{s}-2} \lbrace \lvert 
\alpha_{k} \! - \! \alpha_{p_{q}} \rvert \rbrace)^{\mathfrak{s}-2} 
\! \geqslant \! 1 \! + \! \mathfrak{z}_{m} \! \gg \! 1$; furthermore, 
since, for $\vec{x} \! := \! (x_{1},x_{2},\dotsc,x_{\mathcal{N}}) \! 
\in \! \tilde{\mathfrak{D}}_{\mathcal{N}}(\tilde{M}_{0} \! + \! 1)$, 
$\lvert x_{j} \! - \! \alpha_{p_{q}} \rvert \! > \! (\tilde{M}_{0} \! + 
\! 1)^{-1}$, $j \! = \! 1,2,\dotsc,\mathcal{N}$, $q \! = \! 1,\dotsc,
\mathfrak{s} \! - \! 2,\mathfrak{s}$, it also follows that $(1 \! 
+ \! \mathfrak{z}_{m}) \lvert x_{j} \! - \! \alpha_{p_{q}} \rvert 
\! > \! (\min_{q=1,2,\dotsc,\mathfrak{s}-2} \lbrace \lvert 
\alpha_{k} \! - \! \alpha_{p_{q}} \rvert \rbrace)^{\mathfrak{s}-2}$.} Via 
Equations~\eqref{eqlmrtz23} and~\eqref{eqlmrtz32}, one shows that
\begin{equation*}
\tilde{\mathbb{P}}(x_{1},x_{2},\dotsc,x_{\mathcal{N}};\lambda) 
\! = \! \prod_{m=1}^{\mathcal{N}} \tilde{\mathbb{G}}(x_{m}):
\end{equation*}
for $\vec{x} \! \in \! \tilde{\mathfrak{D}}_{\mathcal{N}}(\tilde{M}_{0})$, 
it follows, via the latter relation and Equation~\eqref{eqlmrtz31}, that
\begin{equation*}
\tilde{\mathbb{P}}(x_{1},x_{2},\dotsc,x_{\mathcal{N}};\lambda) \! 
= \! \me^{\sum_{m=1}^{\mathcal{N}} \tilde{\mathfrak{g}}(x_{m})},
\end{equation*}
which implies, in turn, via Equation~\eqref{eqlmrtz27}, the 
Estimate~\eqref{eqlmrtz30}, and $0 \! \leqslant \! \tilde{F}(x_{m}) 
\! \leqslant \! 1$, $m \! = \! 1,2,\dotsc,\mathcal{N}$, that, for 
$\vec{x} \! \in \! \tilde{\mathfrak{D}}_{\mathcal{N}}(\tilde{M}_{0})$,
\begin{align*}
\pmb{\pi}^{n}_{k}(\lambda) \underset{\underset{z_{o}=1
+o(1)}{\mathscr{N},n \to \infty}}{=}& \, \lambda^{\varkappa_{nk 
\tilde{k}_{\mathfrak{s}-1}}^{\infty}} \left(\tilde{\mathcal{E}}^{n}_{k} 
\left(\chi_{\tilde{\mathfrak{D}}_{\mathcal{N}}(\tilde{M}_{0})}(\vec{x}) 
\me^{\sum_{m=1}^{\mathcal{N}} \tilde{\mathfrak{g}}(x_{m})} \right) 
\! + \! \mathcal{O} \left(\dfrac{\mathfrak{c}_{82}(n,k,z_{o})}{\sqrt{n}} 
\me^{-\frac{1}{2} \frac{n}{\mathcal{N}}(\mathcal{N}-1)} \right) 
\right) \\
\underset{\underset{z_{o}=1+o(1)}{\mathscr{N},n \to \infty}}{=}& \, 
\lambda^{\varkappa_{nk \tilde{k}_{\mathfrak{s}-1}}^{\infty}} \left(
\tilde{\mathcal{E}}^{n}_{k} \left(\me^{\sum_{m=1}^{\mathcal{N}} 
\tilde{\mathfrak{g}}(x_{m})} \right) \! - \! \underbrace{\tilde{
\mathcal{E}}^{n}_{k} \left(\chi_{\tilde{\mathfrak{D}}^{c}_{\mathcal{N}}
(\tilde{M}_{0})}(\vec{x}) \me^{\sum_{m=1}^{\mathcal{N}} \tilde{
\mathfrak{g}}(x_{m})} \right)}_{=_{\underset{z_{o}=1+o(1)}{
\mathscr{N},n \to \infty}} \mathcal{O}(\mathfrak{c}_{83}(n,k,z_{o})
n^{-1/2} \exp (-\frac{1}{2} \frac{n}{\mathcal{N}}(\mathcal{N}-1)))} 
\! + \mathcal{O} \left(\dfrac{\mathfrak{c}_{82}(n,k,z_{o})}{\sqrt{n}} 
\me^{-\frac{1}{2} \frac{n}{\mathcal{N}}(\mathcal{N}-1)} \right) 
\right) \\
\underset{\underset{z_{o}=1+o(1)}{\mathscr{N},n \to \infty}}{=}& \, 
\lambda^{\varkappa_{nk \tilde{k}_{\mathfrak{s}-1}}^{\infty}} \left(
\tilde{\mathcal{E}}^{n}_{k} \left(\me^{\sum_{m=1}^{\mathcal{N}} 
\tilde{\mathfrak{g}}(x_{m})} \right) \! + \! \mathcal{O} \left(
\dfrac{\mathfrak{c}_{84}(n,k,z_{o})}{\sqrt{n}} \me^{-\frac{1}{2} 
\frac{n}{\mathcal{N}}(\mathcal{N}-1)} \right) \right),
\end{align*}
where $\mathfrak{c}_{m}(n,k,z_{o}) \! =_{\underset{z_{o}=1+o(1)}{
\mathscr{N},n \to \infty}} \! \mathcal{O}(1)$, $m \! = \! 83,84$; 
hence, for $n \! \in \! \mathbb{N}$ and $k \! \in \! \lbrace 1,2,
\dotsc,K \rbrace$ such that $\alpha_{p_{\mathfrak{s}}} \! := \! 
\alpha_{k} \! \neq \! \infty$, one concludes that, uniformly for, 
say, $\lvert \lambda \rvert \! \geqslant \! 
2(\tilde{M}_{0} \! + \! 1) \! + \! 3(\mathfrak{s} \! - \! 1)(1 \! + \! 
\max_{q=1,\dotsc,\mathfrak{s}-2,\mathfrak{s}} \lbrace \lvert 
\alpha_{p_{q}} \rvert \rbrace)(\min_{q=1,\dotsc,\mathfrak{s}-2,
\mathfrak{s}} \lbrace \tilde{\Delta}(q) \rbrace)^{-1}$ $(\gg \! 1)$,
\begin{equation} \label{eqlmrtz34} 
\pmb{\pi}^{n}_{k}(\lambda) \underset{\underset{z_{o}=1
+o(1)}{\mathscr{N},n \to \infty}}{=} \lambda^{\varkappa_{nk 
\tilde{k}_{\mathfrak{s}-1}}^{\infty}} \left(\tilde{\mathcal{E}}^{n}_{k} 
\left(\me^{\sum_{m=1}^{\mathcal{N}} \tilde{\mathfrak{g}}(x_{m})} 
\right) \! + \! \mathcal{O} \left(\dfrac{\mathfrak{c}_{84}(n,k,z_{o})}{
\sqrt{n}} \me^{-\frac{1}{2} \frac{n}{\mathcal{N}}(\mathcal{N}-1)} 
\right) \right).
\end{equation}
If it can be established that, for $n \! \in \! \mathbb{N}$ and $k \! \in 
\! \lbrace 1,2,\dotsc,K \rbrace$ such that $\alpha_{p_{\mathfrak{s}}} 
\! := \! \alpha_{k} \! \neq \! \infty$, $\tilde{\mathcal{E}}^{n}_{k}
(\exp (\sum_{m=1}^{\mathcal{N}} \tilde{\mathfrak{g}}(x_{m}))) \! 
>_{\underset{z_{o}=1+o(1)}{\mathscr{N},n \to \infty}} \! 0$ (and 
$\mathcal{O}(1))$, then the proof of this part~\pmb{(1)} of the lemma will 
be complete. For $n \! \in \! \mathbb{N}$ and $k \! \in \! \lbrace 1,2,
\dotsc,K \rbrace$ such that $\alpha_{p_{\mathfrak{s}}} \! := \! \alpha_{k} 
\! \neq \! \infty$, with $\tilde{M}_{0}$ as refined above, choose, for $x \! 
= \! \tilde{M}_{0} \! + \! 1$, $\tilde{\lambda}_{0} \! > \! 2(\tilde{M}_{0} \! + 
\! 1) \! + \! 3(\mathfrak{s} \! - \! 1)(1 \! + \! \max_{q=1,\dotsc,\mathfrak{s}
-2,\mathfrak{s}} \lbrace \lvert \alpha_{p_{q}} \rvert \rbrace)(\min_{q=1,
\dotsc,\mathfrak{s}-2,\mathfrak{s}} \lbrace \tilde{\Delta}(q) \rbrace)^{-1}$ 
$(\gg \! 1)$ so that, in the double-scaling limit $\mathscr{N},n \! \to \! 
\infty$ such that $z_{o} \! = \! 1 \! + \! o(1)$,
\begin{equation} \label{eqlmrtzgt1} 
\left. \ln (\tilde{\mathbb{G}}(x)) \vphantom{M^{M^{M^{M^{M^{M}}}}}} 
\right\rvert_{\underset{\lambda = \tilde{\lambda}_{0}}{x=\tilde{M}_{0}+1}} 
\geqslant -\frac{n}{8 \mathcal{N} \varpi_{0}^{\ast}} \dfrac{\ln \left(8(\tilde{M}_{0} 
\! + \! 1)(1 \! + \! \tilde{\Delta}(\mathfrak{s}))^{\frac{\gamma_{k}}{K}} 
\left(\prod_{q=1}^{\mathfrak{s}-2}(1 \! + \! \tilde{\Delta}(q))^{\frac{
\gamma_{i(q)_{k_{q}}}}{K}} \right) \left(\min\limits_{q=1,2,\dotsc,\mathfrak{s}-2} 
\lbrace \lvert \alpha_{k} \! - \! \alpha_{p_{q}} \rvert \rbrace \right)^{\sum_{q=1}^{
\mathfrak{s}-2} \frac{\gamma_{i(q)_{k_{q}}}}{K}} \right)}{\ln \left(\frac{4}{3}
(\tilde{M}_{0} \! + \! 1)6^{\mathfrak{s}} \left(\min\limits_{q=1,2,\dotsc,
\mathfrak{s}-2} \lbrace \lvert \alpha_{k} \! - \! \alpha_{p_{q}} \rvert \rbrace 
\right)^{\mathfrak{s}-2} \right)},
\end{equation}
where $\varpi_{0}^{\ast}$ $(= \! \varpi_{0}^{\ast}(n,k,z_{o}) \! =_{\underset{z_{o}
=1+o(1)}{\mathscr{N},n \to \infty}} \! \mathcal{O}(1))$, whose exact 
form is not important for the purposes of this proof, is described in 
Equation~\eqref{eqlmrtzgt2} below: a tedious calculation reveals that, 
for $\lvert \lambda \rvert \! \geqslant \! \tilde{\lambda}_{0}$,
\begin{equation} \label{eqlmrtzgt2} 
\left. \ln (\tilde{\mathbb{G}}(x)) \vphantom{M^{M^{M^{M^{M^{M}}}}}} 
\right\vert_{\underset{\lambda = \lvert \lambda \rvert}{x=\tilde{M}_{0}+1}} 
\underset{\underset{z_{o}=1+o(1)}{\mathscr{N},n \to \infty}}{\geqslant} 
\varpi_{0}^{\ast} \left. \ln (\tilde{\mathbb{G}}(x)) 
\vphantom{M^{M^{M^{M^{M^{M}}}}}} \right\vert_{\underset{\lambda 
= \tilde{\lambda}_{0}}{x=\tilde{M}_{0}+1}}.
\end{equation}
It turns out that, in order to proceed further, the quantity 
$\min_{q=1,2,\dotsc,\mathfrak{s}-2} \lbrace \lvert \alpha_{k} \! 
- \! \alpha_{p_{q}} \rvert \rbrace$ requires careful consideration; 
in fact, the following two cases are distinguished: (i) $\min_{q=
1,2,\dotsc,\mathfrak{s}-2} \lbrace \lvert \alpha_{k} \! - \! 
\alpha_{p_{q}} \rvert \rbrace \! \in \! (0,1)$; and (ii) $\min_{q=
1,2,\dotsc,\mathfrak{s}-2} \lbrace \lvert \alpha_{k} \! - \! 
\alpha_{p_{q}} \rvert \rbrace \! \in \! [1,+\infty)$. A calculation 
shows that, for $k \! \in \! \lbrace 1,2,\dotsc,K \rbrace$ such 
that $\alpha_{p_{\mathfrak{s}}} \! := \! \alpha_{k} \! \neq \! \infty$, 
via the monotonicity of $\ln (\pmb{\cdot})$: (i) for $\min_{q=1,2,
\dotsc,\mathfrak{s}-2} \lbrace \lvert \alpha_{k} \! - \! \alpha_{p_{q}} 
\rvert \rbrace \! \in \! (0,1)$,
\begin{align} \label{eqlmrtzgamm1} 
0 <& \, \dfrac{\ln \left(8(\tilde{M}_{0} \! + \! 1) \left(\min\limits_{q
=1,2,\dotsc,\mathfrak{s}-2} \lbrace \lvert \alpha_{k} \! - \! \alpha_{p_{q}} 
\rvert \rbrace \right)^{(\mathfrak{s}-2)(1-\frac{\gamma_{k}}{K})} 
\right)}{\ln \left(\frac{4}{3}(\tilde{M}_{0} \! + \! 1)
6^{\mathfrak{s}} \left(\min\limits_{q=1,2,\dotsc,\mathfrak{s}-2} 
\lbrace \lvert \alpha_{k} \! - \! \alpha_{p_{q}} \rvert \rbrace 
\right)^{\mathfrak{s}-2} \right)} \nonumber \\
\leqslant& \, \dfrac{\ln \left(8(\tilde{M}_{0} \! + \! 1)(1 \! + \! 
\tilde{\Delta}(\mathfrak{s}))^{\frac{\gamma_{k}}{K}} \left(\prod_{q=
1}^{\mathfrak{s}-2}(1 \! + \! \tilde{\Delta}(q))^{\frac{\gamma_{i(q)_{
k_{q}}}}{K}} \right) \left(\min\limits_{q=1,2,\dotsc,\mathfrak{s}-2} 
\lbrace \lvert \alpha_{k} \! - \! \alpha_{p_{q}} \rvert \rbrace \right)^{
\sum_{q=1}^{\mathfrak{s}-2} \frac{\gamma_{i(q)_{k_{q}}}}{K}} 
\right)}{\ln \left(\frac{4}{3}(\tilde{M}_{0} \! + \! 1)6^{\mathfrak{s}} 
\left(\min\limits_{q=1,2,\dotsc,\mathfrak{s}-2} \lbrace \lvert 
\alpha_{k} \! - \! \alpha_{p_{q}} \rvert \rbrace \right)^{\mathfrak{s}-2} 
\right)} \nonumber \\
\leqslant& \, \dfrac{\ln \left(8(\tilde{M}_{0} \! + \! 1)
(\frac{3}{2})^{\mathfrak{s}-2}(\frac{3}{2})^{-\frac{(\mathfrak{s}-2) 
\gamma_{k}}{K}}(\frac{3}{2})^{\frac{\gamma_{k}}{K}} \left(
\min\limits_{q=1,2,\dotsc,\mathfrak{s}-2} \lbrace \lvert 
\alpha_{k}-\alpha_{p_{q}} \rvert \rbrace \right)^{\frac{
\mathfrak{s}-2}{K}}\right)}{\ln \left(\frac{4}{3}(\tilde{M}_{0} \! + \! 1)
6^{\mathfrak{s}} \left(\min\limits_{q=1,2,\dotsc,\mathfrak{s}-2} 
\lbrace \lvert \alpha_{k} \! - \! \alpha_{p_{q}} \rvert \rbrace 
\right)^{\mathfrak{s}-2} \right)} \! \leqslant \! 1 \! - \! 
\tilde{\Gamma}_{\mathrm{A}},
\end{align}
where
\begin{equation*}
\tilde{\Gamma}_{\mathrm{A}} \! := \! \dfrac{\ln \left(4^{\mathfrak{s}
-1} \left(\min\limits_{q=1,2,\dotsc,\mathfrak{s}-2} 
\lbrace \lvert \alpha_{k} \! - \! \alpha_{p_{q}} \rvert \rbrace 
\right)^{(\mathfrak{s}-2)(1-\frac{1}{K})} \right)}{\ln \left(8(\tilde{M}_{0} 
\! + \! 1)6^{\mathfrak{s}-1} \left(\min\limits_{q=1,2,\dotsc,\mathfrak{s}
-2} \lbrace \lvert \alpha_{k} \! - \! \alpha_{p_{q}} \rvert \rbrace 
\right)^{\mathfrak{s}-2} \right)} \quad (\in \! (0,1));
\end{equation*}
and (ii) for $\min_{q=1,2,\dotsc,\mathfrak{s}-2} \lbrace \lvert 
\alpha_{k} \! - \! \alpha_{p_{q}} \rvert \rbrace \! \in \! [1,+\infty)$,
\begin{align} \label{eqlmrtzgamm2} 
0 <& \, \dfrac{\ln \left(8(\tilde{M}_{0} \! + \! 1) \left(\min\limits_{q=
1,2,\dotsc,\mathfrak{s}-2} \lbrace \lvert \alpha_{k}-\alpha_{p_{q}} 
\rvert \rbrace \right)^{\frac{\mathfrak{s}-2}{K}} \right)}{\ln \left(
\frac{4}{3}(\tilde{M}_{0} \! + \! 1)6^{\mathfrak{s}} \left(\min\limits_{q=
1,2,\dotsc,\mathfrak{s}-2} \lbrace \lvert \alpha_{k} \! - \! \alpha_{p_{q}} 
\rvert \rbrace \right)^{\mathfrak{s}-2} \right)} \nonumber \\
\leqslant& \, \dfrac{\ln \left(8(\tilde{M}_{0} \! + \! 1)(1 \! + \! 
\tilde{\Delta}(\mathfrak{s}))^{\frac{\gamma_{k}}{K}} \left(\prod_{q=
1}^{\mathfrak{s}-2}(1 \! + \! \tilde{\Delta}(q))^{\frac{\gamma_{i(q)_{
k_{q}}}}{K}} \right) \left(\min\limits_{q=1,2,\dotsc,\mathfrak{s}-2} 
\lbrace \lvert \alpha_{k} \! - \! \alpha_{p_{q}} \rvert \rbrace \right)^{
\sum_{q=1}^{\mathfrak{s}-2} \frac{\gamma_{i(q)_{k_{q}}}}{K}} 
\right)}{\ln \left(\frac{4}{3}(\tilde{M}_{0} \! + \! 1)6^{\mathfrak{s}} 
\left(\min\limits_{q=1,2,\dotsc,\mathfrak{s}-2} \lbrace \lvert 
\alpha_{k} \! - \! \alpha_{p_{q}} \rvert \rbrace \right)^{
\mathfrak{s}-2} \right)} \nonumber \\
\leqslant& \, \dfrac{\ln \left(8(\tilde{M}_{0} \! + \! 1)(\frac{3}{2})^{
\mathfrak{s}-2}(\frac{3}{2})^{-\frac{(\mathfrak{s}-2) \gamma_{k}}{K}}
(\frac{3}{2})^{\frac{\gamma_{k}}{K}} \left(\min\limits_{q=1,2,
\dotsc,\mathfrak{s}-2} \lbrace \lvert \alpha_{k}-\alpha_{p_{q}} 
\rvert \rbrace \right)^{(\mathfrak{s}-2)(1-\frac{\gamma_{k}}{K})} 
\right)}{\ln \left(\frac{4}{3}(\tilde{M}_{0} \! + \! 1)6^{\mathfrak{s}} 
\left(\min\limits_{q=1,2,\dotsc,\mathfrak{s}-2} \lbrace \lvert \alpha_{k} 
\! - \! \alpha_{p_{q}} \rvert \rbrace \right)^{\mathfrak{s}-2} \right)} 
\! \leqslant \! 1 \! - \! \tilde{\Gamma}_{\mathrm{B}},
\end{align}
where
\begin{equation*}
\tilde{\Gamma}_{\mathrm{B}} \! := \! \dfrac{\ln \left(4^{\mathfrak{s}
-1} \left(\min\limits_{q=1,2,\dotsc,\mathfrak{s}-2} \lbrace \lvert 
\alpha_{k} \! - \! \alpha_{p_{q}} \rvert \rbrace \right)^{\frac{
(\mathfrak{s}-2) \gamma_{k}}{K}} \right)}{\ln \left(8(\tilde{M}_{0} \! 
+ \! 1)6^{\mathfrak{s}-1} \left(\min\limits_{q=1,2,\dotsc,\mathfrak{s}
-2} \lbrace \lvert \alpha_{k} \! - \! \alpha_{p_{q}} \rvert \rbrace 
\right)^{\mathfrak{s}-2} \right)} \quad (\in \! (0,1)).
\end{equation*}
Since $0 \! \leqslant \! \tilde{F}(x) \! \leqslant \! 1$, $x \! \in \! \mathbb{R}$, 
it follows that (cf. Equation~\eqref{eqlmrtzgt1}), for $k \! \in \! \lbrace 1,2,
\dotsc,K \rbrace$ such that $\alpha_{p_{\mathfrak{s}}} \! := \! \alpha_{k} 
\! \neq \! \infty$,
\begin{align} \label{eqlmrtzgt3} 
& \, -\dfrac{n}{8 \mathcal{N} \varpi_{0}^{\ast}} \dfrac{\ln \left(8(\tilde{M}_{0} 
\! + \! 1)(1 \! + \! \tilde{\Delta}(\mathfrak{s}))^{\frac{\gamma_{k}}{K}} \left(
\prod_{q=1}^{\mathfrak{s}-2}(1 \! + \! \tilde{\Delta}(q))^{\frac{\gamma_{
i(q)_{k_{q}}}}{K}} \right) \left(\min\limits_{q=1,2,\dotsc,\mathfrak{s}-2} 
\lbrace \lvert \alpha_{k} \! - \! \alpha_{p_{q}} \rvert \rbrace 
\right)^{\sum_{q=1}^{\mathfrak{s}-2} \frac{\gamma_{i(q)_{k_{q}}}}{K}} 
\right)}{\ln \left(\frac{4}{3}(\tilde{M}_{0} \! + \! 1)6^{\mathfrak{s}} 
\left(\min\limits_{q=1,2,\dotsc,\mathfrak{s}-2} \lbrace \lvert 
\alpha_{k} \! - \! \alpha_{p_{q}} \rvert \rbrace \right)^{\mathfrak{s}-2} 
\right)} \nonumber \\
& \, \leqslant -\dfrac{n}{8 \mathcal{N} \varpi_{0}^{\ast}} \dfrac{\ln \left(
8(\tilde{M}_{0} \! + \! 1)(1 \! + \! \tilde{\Delta}(\mathfrak{s}))^{\frac{
\gamma_{k}}{K}} \left(\prod_{q=1}^{\mathfrak{s}-2}(1 \! + \! \tilde{\Delta}
(q))^{\frac{\gamma_{i(q)_{k_{q}}}}{K}} \right) \left(\min\limits_{q=1,
2,\dotsc,\mathfrak{s}-2} \lbrace \lvert \alpha_{k} \! - \! \alpha_{
p_{q}} \rvert \rbrace \right)^{\sum_{q=1}^{\mathfrak{s}-2} \frac{
\gamma_{i(q)_{k_{q}}}}{K}} \right)}{\ln \left(\frac{4}{3}(\tilde{M}_{0} \! 
+ \! 1)6^{\mathfrak{s}} \left(\min\limits_{q=1,2,\dotsc,\mathfrak{s}-2} 
\lbrace \lvert \alpha_{k} \! - \! \alpha_{p_{q}} \rvert \rbrace 
\right)^{\mathfrak{s}-2} \right)} \tilde{F}(x) \! \leqslant \! 0, 
\quad x \! \in \! \mathbb{R},
\end{align} 
whence, for $k \! \in \! \lbrace 1,2,\dotsc,K \rbrace$ such that 
$\alpha_{p_{\mathfrak{s}}} \! := \! \alpha_{k} \! \neq \! \infty$ and 
uniformly for $\lvert \lambda \rvert \! \geqslant \! \tilde{\lambda}_{0} 
\! > \! 2(\tilde{M}_{0} \! + \! 1) \! + \! 3(\mathfrak{s} \! - \! 1)(1 \! + \! 
\max_{q=1,\dotsc,\mathfrak{s}-2,\mathfrak{s}} \lbrace \lvert \alpha_{p_{q}} 
\rvert \rbrace)(\min_{q=1,\dotsc,\mathfrak{s}-2,\mathfrak{s}} \lbrace 
\tilde{\Delta}(q) \rbrace)^{-1}$ $(\gg \! 1)$, via Equations~\eqref{eqlmrtz31}, 
\eqref{eqlmrtzgt2}, and~\eqref{eqlmrtzgt3}, in the double-scaling limit 
$\mathscr{N},n \! \to \! \infty$ such that $z_{o} \! = \! 1 \! + \! o(1)$,
\begin{align} \label{eqlmrtz35} 
\tilde{\mathfrak{g}}(x) \geqslant& \, \underbrace{-\dfrac{n}{8 
\mathcal{N}} \dfrac{\ln \left(8(\tilde{M}_{0} \! + \! 1)(1 \! + \! 
\tilde{\Delta}(\mathfrak{s}))^{\frac{\gamma_{k}}{K}} \left(\prod_{q=1}^{
\mathfrak{s}-2}(1 \! + \! \tilde{\Delta}(q))^{\frac{\gamma_{i(q)_{k_{q}}}}{K}} 
\right) \left(\min\limits_{q=1,2,\dotsc,\mathfrak{s}-2} \lbrace \lvert \alpha_{k} 
\! - \! \alpha_{p_{q}} \rvert \rbrace \right)^{\sum_{q=1}^{\mathfrak{s}-2} 
\frac{\gamma_{i(q)_{k_{q}}}}{K}} \right)}{\ln \left(\frac{4}{3}(\tilde{M}_{0} 
\! + \! 1)6^{\mathfrak{s}} \left(\min\limits_{q=1,2,\dotsc,\mathfrak{s}-2} 
\lbrace \lvert \alpha_{k} \! - \! \alpha_{p_{q}} \rvert \rbrace 
\right)^{\mathfrak{s}-2} \right)} \tilde{F}(x)}_{:= \, \tilde{F}^{\ast}(x)} 
\nonumber \\
\geqslant& \, -\underbrace{\dfrac{n}{8 \mathcal{N}} \dfrac{\ln \left(
8(\tilde{M}_{0} \! + \! 1)(1 \! + \! \tilde{\Delta}(\mathfrak{s}))^{
\frac{\gamma_{k}}{K}} \left(\prod_{q=1}^{\mathfrak{s}-2}(1 \! + \! 
\tilde{\Delta}(q))^{\frac{\gamma_{i(q)_{k_{q}}}}{K}} \right) \left(
\min\limits_{q=1,2,\dotsc,\mathfrak{s}-2} \lbrace \lvert \alpha_{k} \! 
- \! \alpha_{p_{q}} \rvert \rbrace \right)^{\sum_{q=1}^{\mathfrak{s}-2} 
\frac{\gamma_{i(q)_{k_{q}}}}{K}} \right)}{\ln \left(\frac{4}{3}(\tilde{M}_{0} 
\! + \! 1)6^{\mathfrak{s}} \left(\min\limits_{q=1,2,\dotsc,\mathfrak{s}-2} 
\lbrace \lvert \alpha_{k} \! - \! \alpha_{p_{q}} \rvert \rbrace \right)^{
\mathfrak{s}-2} \right)}}_{:= \, \tilde{\mathfrak{c}}^{\ast}} \, , 
\quad x \! \in \! \mathbb{R}.
\end{align}
(Note that $\tilde{\mathfrak{c}}^{\ast} \! =_{\underset{z_{o}=
1+o(1)}{\mathscr{N},n \to \infty}} \! \mathcal{O}(1)$ and $> 
\! 0$.) Equations~\eqref{eqlmrtz33}, \eqref{eqlmrtzgamm1}, 
\eqref{eqlmrtzgamm2}, and~\eqref{eqlmrtz35} imply, in particular, 
that, for $n \! \in \! \mathbb{N}$ and $k \! \in \! \lbrace 1,2,\dotsc,
K \rbrace$ such that $\alpha_{p_{\mathfrak{s}}} \! := \! \alpha_{k} 
\! \neq \! \infty$, in the double-scaling limit $\mathscr{N},n \! 
\to \! \infty$ such that $z_{o} \! = \! 1 \! + \! o(1)$, uniformly 
for $\lvert \lambda \rvert \! \geqslant \! \tilde{\lambda}_{0} \! > 
\! 2(\tilde{M}_{0} \! + \! 1) \! + \! 3(\mathfrak{s} \! - \! 1)(1 \! + 
\! \max_{q=1,\dotsc,\mathfrak{s}-2,\mathfrak{s}} \lbrace \lvert 
\alpha_{p_{q}} \rvert \rbrace)(\min_{q=1,\dotsc,\mathfrak{s}-2,
\mathfrak{s}} \lbrace \tilde{\Delta}(q) \rbrace)^{-1}$ $(\gg \! 1)$, 
the function $\tilde{F}^{\ast}$ is real-valued, bounded, and continuous 
on $\mathbb{R}$; hence, via Equation~\eqref{eqlmrtz8} and the conditions 
satisfied by $\tilde{F}$,
\begin{align*}
\lim_{\underset{z_{o}=1+o(1)}{\mathscr{N},n \to \infty}} 
\dfrac{1}{\mathcal{N}} \ln \tilde{\mathcal{E}}^{n}_{k} \left(
\me^{\sum_{m=1}^{\mathcal{N}} \tilde{F}^{\ast}(x_{m})} \right) 
\underset{\underset{z_{o}=1+o(1)}{\mathscr{N},n \to \infty}}{
=}& \, \int_{\mathbb{R}} \tilde{F}^{\ast}(\xi) \, \md \mu_{
\widetilde{V}}^{f}(\xi) \underset{\underset{z_{o}=1+o(1)}{
\mathscr{N},n \to \infty}}{=} \int_{\tilde{\mathfrak{D}}^{c}
(\tilde{M}_{0}+1)} \tilde{F}^{\ast}(\xi) \, \underbrace{\md 
\mu_{\widetilde{V}}^{f}(\xi)}_{= \, 0} \! + \! \int_{\tilde{
\mathfrak{D}}(\tilde{M}_{0}+1)} \tilde{F}^{\ast}(\xi) \, \md 
\mu_{\widetilde{V}}^{f}(\xi) \\
\underset{\underset{z_{o}=1+o(1)}{\mathscr{N},n \to \infty}}{=}& 
\, \int_{\tilde{\mathfrak{D}}(\tilde{M}_{0}+1) \setminus \tilde{
\mathfrak{D}}(\tilde{M}_{0})} \tilde{F}^{\ast}(\xi) \, \underbrace{
\md \mu_{\widetilde{V}}^{f}(\xi)}_{= \, 0} \! + \! \int_{\tilde{
\mathfrak{D}}(\tilde{M}_{0})} \tilde{F}^{\ast}(\xi) \, \md \mu_{
\widetilde{V}}^{f}(\xi) \\
\underset{\underset{z_{o}=1+o(1)}{\mathscr{N},n \to \infty}}{=}& 
\, \int_{\tilde{\mathfrak{D}}(\tilde{M}_{0}) \setminus J_{f}} 
\tilde{F}^{\ast}(\xi) \, \underbrace{\md \mu_{\widetilde{V}}^{f}
(\xi)}_{= \, 0} \! + \! \int_{J_{f}} \tilde{F}^{\ast}(\xi) \, \md 
\mu_{\widetilde{V}}^{f}(\xi) \\
\underset{\underset{z_{o}=1+o(1)}{\mathscr{N},n \to \infty}}{=}& 
\, -\tilde{\mathfrak{c}}^{\ast} \int_{J_{f}} \underbrace{\tilde{F}(\xi)}_{= 
\, 1} \, \md \mu_{\widetilde{V}}^{f}(\xi) \underset{\underset{z_{o}=
1+o(1)}{\mathscr{N},n \to \infty}}{=} -\tilde{\mathfrak{c}}^{\ast} 
\underbrace{\int_{J_{f}} \, \md \mu_{\widetilde{V}}^{f}(\xi)}_{= \, 1}
\underset{\underset{z_{o}=1+o(1)}{\mathscr{N},n \to \infty}}{=} 
-\tilde{\mathfrak{c}}^{\ast} \quad \Rightarrow
\end{align*}
\begin{equation} \label{eqlmrtz36} 
\lim_{\underset{z_{o}=1+o(1)}{\mathscr{N},n \to \infty}} \dfrac{1}{\mathcal{N}} 
\ln \tilde{\mathcal{E}}^{n}_{k} \left(\me^{\sum_{m=1}^{\mathcal{N}} 
\tilde{F}^{\ast}(x_{m})} \right) \underset{\underset{z_{o}=1+o(1)}{\mathscr{N},
n \to \infty}}{=} -\tilde{\mathfrak{c}}^{\ast}.
\end{equation}
Via Equation~\eqref{eqlmrtz35}, the linearity of $\tilde{\mathcal{E}}^{n}_{k}$, 
the fact that $\tilde{\mathcal{E}}^{n}_{k}(1) \! = \! 1$, and the monotonicity 
of $\exp (\pmb{\cdot})$, it follows that, for $m \! = \! 1,2,\dotsc,\mathcal{N}$, 
$\tilde{\mathfrak{g}}(x_{m}) \! \geq \! \tilde{F}^{\ast}(x_{m}) \! \geqslant \! 
-\tilde{\mathfrak{c}}^{\ast}$ $\Rightarrow$ $\sum_{m=1}^{\mathcal{N}} 
\tilde{\mathfrak{g}}(x_{m}) \! \geq \! \sum_{m=1}^{\mathcal{N}} \tilde{F}^{
\ast}(x_{m}) \! \geqslant \! -\mathcal{N} \, \tilde{\mathfrak{c}}^{\ast}$ 
$\Rightarrow$ $\tilde{\mathcal{E}}^{n}_{k}(\me^{\sum_{m=1}^{\mathcal{N}} 
\tilde{\mathfrak{g}}(x_{m})}) \! \geq \! \tilde{\mathcal{E}}^{n}_{k}(\me^{
\sum_{m=1}^{\mathcal{N}} \tilde{F}^{\ast}(x_{m})}) \! \geqslant \! 
\me^{-\mathcal{N} \, \tilde{\mathfrak{c}}^{\ast}}$; hence, for $n \! 
\in \! \mathbb{N}$ and $k \! \in \! \lbrace 1,2,\dotsc,K \rbrace$ such 
that $\alpha_{p_{\mathfrak{s}}} \! := \! \alpha_{k} \! \neq \! \infty$, 
uniformly for, say, $\lvert \lambda \rvert \! \geqslant \! \tilde{\lambda}_{0} 
\! > \! 2(\tilde{M}_{0} \! + \! 1) \! + \! 3(\mathfrak{s} \! - \! 1)(1 \! + \! 
\max_{q=1,\dotsc,\mathfrak{s}-2,\mathfrak{s}} \lbrace \lvert \alpha_{p_{q}} 
\rvert \rbrace)(\min_{q=1,\dotsc,\mathfrak{s}-2,\mathfrak{s}} \lbrace 
\tilde{\Delta}(q) \rbrace)^{-1}$ $(\gg \! 1)$, one arrives at (for 
$\vec{x} \! \in \! \tilde{\mathfrak{D}}_{\mathcal{N}}(\tilde{M}_{0}))$, 
via Equations~\eqref{eqlmrtzgamm1} and~\eqref{eqlmrtzgamm2}: 
(i) for $\min_{q=1,2,\dotsc,\mathfrak{s}-2} \lbrace \lvert \alpha_{k} 
\! - \! \alpha_{p_{q}} \rvert \rbrace \! \in \! (0,1)$,
\begin{align} \label{eqlmrtz37} 
\tilde{\mathcal{E}}^{n}_{k} \left(\me^{\sum_{m=1}^{\mathcal{N}} 
\tilde{\mathfrak{g}}(x_{m})} \right) \underset{\underset{z_{o}=1+
o(1)}{\mathscr{N},n \to \infty}}{\geqslant}& \, \exp \left(-\dfrac{n}{8} 
\dfrac{\ln \left(8(\tilde{M}_{0} \! + \! 1)(1 \! + \! \tilde{\Delta}
(\mathfrak{s}))^{\frac{\gamma_{k}}{K}} \left(\prod_{q=1}^{\mathfrak{s}
-2}(1 \! + \! \tilde{\Delta}(q))^{\frac{\gamma_{i(q)_{k_{q}}}}{K}} \right) 
\left(\min\limits_{q=1,2,\dotsc,\mathfrak{s}-2} \lbrace \lvert 
\alpha_{k} \! - \! \alpha_{p_{q}} \rvert \rbrace \right)^{\sum_{q=1}^{
\mathfrak{s}-2} \frac{\gamma_{i(q)_{k_{q}}}}{K}} \right)}{\ln \left(
\frac{4}{3}(\tilde{M}_{0} \! + \! 1)6^{\mathfrak{s}} \left(\min\limits_{q
=1,2,\dotsc,\mathfrak{s}-2} \lbrace \lvert \alpha_{k} \! - \! 
\alpha_{p_{q}} \rvert \rbrace \right)^{\mathfrak{s}-2} \right)} \right) 
\nonumber \\
\underset{\underset{z_{o}=1+o(1)}{\mathscr{N},n \to \infty}}{\geqslant}& 
\, \exp \left(-\dfrac{n}{8}(1 \! - \! \tilde{\Gamma}_{\mathrm{A}}) \right) 
> 0;
\end{align}
and (ii) for $\min_{q=1,2,\dotsc,\mathfrak{s}-2} \lbrace \lvert \alpha_{k} 
\! - \! \alpha_{p_{q}} \rvert \rbrace \! \in \! [1,+\infty)$,
\begin{align} \label{eqlmrtz38} 
\tilde{\mathcal{E}}^{n}_{k} \left(\me^{\sum_{m=1}^{\mathcal{N}} 
\tilde{\mathfrak{g}}(x_{m})} \right) \underset{\underset{z_{o}=1+
o(1)}{\mathscr{N},n \to \infty}}{\geqslant}& \, \exp \left(-\dfrac{n}{8} 
\dfrac{\ln \left(8(\tilde{M}_{0} \! + \! 1)(1 \! + \! \tilde{\Delta}
(\mathfrak{s}))^{\frac{\gamma_{k}}{K}} \left(\prod_{q=1}^{\mathfrak{s}-2}
(1 \! + \! \tilde{\Delta}(q))^{\frac{\gamma_{i(q)_{k_{q}}}}{K}} \right) 
\left(\min\limits_{q=1,2,\dotsc,\mathfrak{s}-2} \lbrace \lvert 
\alpha_{k} \! - \! \alpha_{p_{q}} \rvert \rbrace \right)^{\sum_{q=1}^{
\mathfrak{s}-2} \frac{\gamma_{i(q)_{k_{q}}}}{K}} \right)}{\ln \left(
\frac{4}{3}(\tilde{M}_{0} \! + \! 1)6^{\mathfrak{s}} \left(\min\limits_{q
=1,2,\dotsc,\mathfrak{s}-2} \lbrace \lvert \alpha_{k} \! - \! 
\alpha_{p_{q}} \rvert \rbrace \right)^{\mathfrak{s}-2} \right)} \right) 
\nonumber \\
\underset{\underset{z_{o}=1+o(1)}{\mathscr{N},n \to \infty}}{\geqslant}& 
\, \exp \left(-\dfrac{n}{8}(1 \! - \! \tilde{\Gamma}_{\mathrm{B}}) 
\right) > 0.
\end{align}
For $n \! \in \! \mathbb{N}$ and $k \! \in \! \lbrace 1,2,\dotsc,K \rbrace$ 
such that $\alpha_{p_{\mathfrak{s}}} \! := \! \alpha_{k} \! \neq \! \infty$, 
let $\tilde{\mathfrak{Z}}_{0} \! := \! \lbrace \mathstrut \lambda \! \in \! 
\overline{\mathbb{C}}; \, \pmb{\pi}^{n}_{k}(\lambda) \! = \! 0 \rbrace$: 
Equations~\eqref{eqlmrtz34}, \eqref{eqlmrtz37}, and~\eqref{eqlmrtz38} 
show that, in the double-scaling limit $\mathscr{N},n \! \to \! \infty$ such 
that $z_{o} \! = \! 1 \! + \! o(1)$, $\tilde{\mathfrak{Z}}_{0} \cap (\lbrace 
\alpha_{p_{1}},\dotsc,\alpha_{p_{\mathfrak{s}-2}},\alpha_{p_{\mathfrak{s}}} 
\rbrace \cup \lbrace \alpha_{p_{\mathfrak{s}-1}} \rbrace) \! = \! \varnothing$, 
that is (cf. Remark~\ref{rem1.2.4}), $\tilde{\mathfrak{Z}}_{0} \! = \! 
\lbrace \tilde{\mathfrak{z}}^{n}_{k}(j) \rbrace_{j=1}^{\mathcal{N}}$, with 
$\tilde{\mathfrak{z}}^{n}_{k}(j) \! \neq \! \alpha_{p_{q}}$, $j \! = \! 1,2,
\dotsc,\mathcal{N}$, $q \! = \! 1,2,\dotsc,\mathfrak{s}$. The above discussion, 
in conjunction with Equations~\eqref{eqlmrtz34}, \eqref{eqlmrtz37}, 
and~\eqref{eqlmrtz38}, implies, in particular, that, for $n \! \in \! \mathbb{N}$ 
and $k \! \in \! \lbrace 1,2,\dotsc,K \rbrace$ such that $\alpha_{p_{\mathfrak{s}}} 
\! := \! \alpha_{k} \! \neq \! \infty$, in the double-scaling limit $\mathscr{N},n \! 
\to \! \infty$ such that $z_{o} \! = \! 1 \! + \! o(1)$, there exists $\tilde{\lambda}_{1} 
\! = \! \tilde{\lambda}_{1}(n,k,z_{o}) \! > \! 0$ and $\mathcal{O}(1)$ {}\footnote{For 
$k \! \in \! \lbrace 1,2,\dotsc,K \rbrace$ such that $\alpha_{p_{\mathfrak{s}}} \! 
:= \! \alpha_{k} \! \neq \! \infty$, take, say, $\tilde{\lambda}_{1} \! = \! \max 
\left\lbrace \tilde{\lambda}_{1}^{\ast},2(\tilde{M}_{0} \! + \! 1) \! + \! 3(\mathfrak{s} 
\! - \! 1)(1 \! + \! \max_{q=1,\dotsc,\mathfrak{s}-2,\mathfrak{s}} \lbrace \lvert 
\alpha_{p_{q}} \rvert \rbrace)(\min_{q=1,\dotsc,\mathfrak{s}-2,\mathfrak{s}} 
\lbrace \tilde{\Delta}(q) \rbrace)^{-1} \right\rbrace$, where
\begin{equation*}
\tilde{\lambda}_{1}^{\ast} \! := \! (\tilde{M}_{0} \! + \! 1) \left(1 
\! - \! \left(\tfrac{(\max_{q=1,2,\dotsc,\mathfrak{s}-2} \lbrace 
\lvert \alpha_{k}-\alpha_{p_{q}} \rvert \rbrace)^{\sum_{q=1}^{
\mathfrak{s}-2} \frac{\gamma_{i(q)_{k_{q}}}}{K}}}{(\tilde{M}_{0}+1)-
\alpha_{k}} \right)^{\frac{K}{\mathfrak{s} \varpi_{0}^{\sharp}}} 
\exp \left(-\frac{\ln (\tilde{\mathrm{A}}_{0}^{\sharp})}{8 \mathfrak{s} 
\varpi_{0}^{\ast} \varpi_{0}^{\sharp} \ln (\tilde{\mathrm{B}}_{0}^{
\sharp})} \right) \right)^{-1},
\end{equation*}
with $\varpi_{0}^{\ast}$ discussed in Equation~\eqref{eqlmrtzgt2}, 
$\varpi_{0}^{\sharp} \! := \! \max \lbrace \gamma_{k},\max_{q=1,2,
\dotsc,\mathfrak{s}-2} \lbrace \gamma_{i(q)_{k_{q}}} \rbrace 
\rbrace$, $\tilde{\mathrm{A}}_{0}^{\sharp} \! := \! 8(\tilde{M}_{0} \! 
+ \! 1)(1 \! + \! \tilde{\Delta}(\mathfrak{s}))^{\frac{\gamma_{k}}{K}} 
\prod_{q=1}^{\mathfrak{s}-2}(1 \! + \! \tilde{\Delta}(q))^{\frac{
\gamma_{i(q)_{k_{q}}}}{K}}(\min_{q^{\prime}=1,2,\dotsc,\mathfrak{s}-2} 
\lbrace \lvert \alpha_{k} \! - \! \alpha_{p_{q^{\prime}}} \rvert \rbrace)^{
\sum_{q^{\prime \prime}=1}^{\mathfrak{s}-2} \frac{\gamma_{i
(q^{\prime \prime})_{k_{q^{\prime \prime}}}}}{K}}$, and 
$\tilde{\mathrm{B}}_{0}^{\sharp} \! := \! \tfrac{4}{3}(\tilde{M}_{0} 
\! + \! 1)6^{\mathfrak{s}}(\min_{q=1,2,\dotsc,\mathfrak{s}-2} 
\lbrace \lvert \alpha_{k} \! - \! \alpha_{p_{q}} \rvert \rbrace)^{
\mathfrak{s}-2}$.} such that $\pmb{\pi}^{n}_{k}(\lambda)$ has no 
zeros (counting multiplicities) in $\overline{\mathbb{R}} \setminus 
\tilde{\mathfrak{D}}(\tilde{\lambda}_{1})$, where $\tilde{\mathfrak{D}}
(\tilde{\lambda}_{1}) \! := \! [-\tilde{\lambda}_{1},\tilde{\lambda}_{1}] 
\setminus \cup_{\underset{q \neq \mathfrak{s}-1}{q=1}}^{\mathfrak{s}} 
\mathscr{O}_{\frac{1}{\tilde{\lambda}_{1}}}(\alpha_{p_{q}})$ (of course, 
$\mathscr{O}_{\frac{1}{\tilde{\lambda}_{1}}}(\alpha_{p_{i}}) \cap 
\mathscr{O}_{\frac{1}{\tilde{\lambda}_{1}}}(\alpha_{p_{j}}) \! = \! \varnothing$ 
$\forall$ $i \! \neq \! j \! \in \! \lbrace 1,\dotsc,\mathfrak{s} \! - \! 2,
\mathfrak{s} \rbrace$), that is, $\tilde{\mathfrak{Z}}_{0} \subseteq 
\tilde{\mathfrak{D}}(\tilde{\lambda}_{1})$. This concludes the proof 
of part~\pmb{(1)}.

\pmb{(2)} The proof of this case, that is, $n \! \in \! \mathbb{N}$ 
and $k \! \in \! \lbrace 1,2,\dotsc,K \rbrace$ such that $\alpha_{
p_{\mathfrak{s}}} \! := \! \alpha_{k} \! = \! \infty$, is virtually identical 
to the proof of \pmb{(1)} above; one mimics, \emph{verbatim}, the 
scheme of the calculations presented in \pmb{(1)} in order to arrive 
at the corresponding claim stated in the lemma; in order to do so, 
however, the analogues of Equations~\eqref{eqlmrtzp}, \eqref{eqlmrtzz}, 
\eqref{eqlmrtz6}--\eqref{eqlmrtz8}, \eqref{eqlmrtz10}, \eqref{eqlmrtz11}, 
\eqref{eqlmrtzps}, \eqref{eqlmrtzzs}, \eqref{eqlmrtz14}, \eqref{eqlmrtz22}, 
\eqref{eqlmrtz23}, \eqref{eqlmrtz27}, \eqref{eqlmrtz30}, \eqref{eqlmrtz31}, 
\eqref{eqlmrtz34}, \eqref{eqlmrtz37}, and~\eqref{eqlmrtz38}, respectively, 
are necessary, which, in the present case, read: with $\mathcal{N} \! = \! 
(n \! - \! 1)K \! + \! k \! = \! \sum_{q=1}^{\mathfrak{s}-1} \varkappa_{nk 
\tilde{k}_{q}} \! + \! \varkappa_{nk}$,
\begin{equation} \label{eqlmrtz39} 
\hat{\mathcal{P}}^{n}_{k}(x_{1},x_{2},\dotsc,x_{\mathcal{N}}) \! := \! 
\dfrac{1}{\hat{\mathscr{Z}}^{n}_{k}} \me^{-n \sum_{m_{1}=1}^{\mathcal{
N}} \widetilde{V}(x_{m_{1}})} \prod_{\substack{i,j=1\\j<i}}^{\mathcal{N}}
(x_{j} \! - \! x_{i})^{2} \left(\prod_{m=1}^{\mathcal{N}} \prod_{q=
1}^{\mathfrak{s}-1}(x_{m} \! - \! \alpha_{p_{q}})^{\varkappa_{nk 
\tilde{k}_{q}}} \right)^{-2},
\end{equation}
\begin{equation} \label{eqlmrtz40} 
\hat{\mathscr{Z}}^{n}_{k} \! = \! \idotsint\limits_{\mathbb{R}^{\mathcal{N}}} 
\me^{-n \sum_{m_{1}=1}^{\mathcal{N}} \widetilde{V}(\tau_{m_{1}})} 
\prod_{\substack{i,j=1\\j<i}}^{\mathcal{N}}(\tau_{j} \! - \! \tau_{i})^{2} 
\left(\prod_{m=1}^{\mathcal{N}} \prod_{q=1}^{\mathfrak{s}-1}(\tau_{m} 
\! - \! \alpha_{p_{q}})^{\varkappa_{nk \tilde{k}_{q}}} \right)^{-2} \md 
\tau_{1} \, \md \tau_{2} \, \dotsb \, \md \tau_{\mathcal{N}},
\end{equation}
\begin{equation} \label{eqlmrtz41} 
\dfrac{\mathcal{N}}{n} \dfrac{1}{\mathcal{N}(\mathcal{N} \! - \! 1)} 
\ln (\hat{\mathscr{Z}}^{n}_{k}) \underset{\underset{z_{o}=1+o(1)}{
\mathscr{N},n \to \infty}}{\geqslant} -\left(E_{\widetilde{V}}^{\infty} 
\! + \! \hat{\varepsilon} \right) \left(1 \! + \! \mathcal{O} \left(
\dfrac{\mathfrak{c}_{85}(n,k,z_{o},\hat{\varepsilon})}{n} \right) 
\right) \quad \text{as} \quad \hat{\varepsilon} \! \downarrow \! 0,
\end{equation}
where $\mathfrak{c}_{85}(n,k,z_{o},\hat{\varepsilon}) \! =_{\underset{
z_{o}=1+o(1)}{\mathscr{N},n \to \infty}} \! \mathcal{O}(1)$ 
as $\hat{\varepsilon} \! \downarrow \! 0$,
\begin{equation} \label{eqlmrtz42} 
\operatorname{Prob} \, (\mathbb{R}^{\mathcal{N}} \setminus 
\hat{\mathbb{A}}_{\mathcal{N},2 \hat{\eta}}) \underset{
\underset{z_{o}=1+o(1)}{\mathscr{N},n \to \infty}}{\leqslant} 
\mathfrak{c}_{86}(n,k,z_{o}) \me^{-\frac{n}{\mathcal{N}} 
\mathcal{N}(\mathcal{N}-1) \hat{\eta}},
\end{equation}
where $\mathfrak{c}_{86}(n,k,z_{o}) \! =_{\underset{z_{o}=1+o(1)}{
\mathscr{N},n \to \infty}} \! \mathcal{O}(1)$, and $\hat{\mathbb{A}}_{
\mathcal{N},2 \hat{\eta}} \! := \! \lbrace \mathstrut (x_{1},x_{2},
\dotsc,x_{\mathcal{N}}) \! \in \! \mathbb{R}^{\mathcal{N}}; \, 
(\mathcal{N}(\mathcal{N} \! - \! 1))^{-1} \mathscr{K}^{\widetilde{V},
\infty}_{\mathcal{N}}(x_{1},x_{2},\dotsc,x_{\mathcal{N}}) \! \leqslant 
\! E_{\widetilde{V}}^{\infty} \! + \! 2 \hat{\eta} \rbrace$, with 
$\hat{\eta} \! = \! \hat{\eta}(n,k,z_{o}) \! > \! 0$,
\begin{equation} \label{eqlmrtz43} 
\lim_{\underset{z_{o}=1+o(1)}{\mathscr{N},n \to \infty}} 
\dfrac{1}{\mathcal{N}} \ln \hat{\mathcal{E}}^{n}_{k} \left(\me^{
\mathcal{N}^{-(m-1)} \sum_{i_{1},i_{2},\dotsc,i_{m}} \phi_{0}(x_{i_{1}},
x_{i_{2}},\dotsc,x_{i_{m}})} \right) \underset{\underset{z_{o}=1+o(1)}{
\mathscr{N},n \to \infty}}{=} \idotsint\limits_{\mathbb{R}^{m}} \phi_{0}
(\xi_{1},\xi_{2},\dotsc,\xi_{m}) \, \md \mu_{\widetilde{V}}^{\infty}
(\xi_{1}) \, \md \mu_{\widetilde{V}}^{\infty}(\xi_{2}) \, \dotsb \, \md 
\mu_{\widetilde{V}}^{\infty}(\xi_{m}), \quad m \! \in \! \mathbb{N},
\end{equation}
where $\hat{\mathcal{E}}^{n}_{k}$ denotes the expectation with respect 
to the multi-dimensional probability measure $\hat{\mathcal{P}}^{n}_{k}
(x_{1},x_{2},\dotsc,x_{\mathcal{N}}) \, \md x_{1} \linebreak[4] 
\pmb{\cdot} \, \md x_{2} \, \dotsb \, \md x_{\mathcal{N}}$,
\begin{align} \label{eqlmrtz44} 
\lim_{\underset{z_{o}=1+o(1)}{\mathscr{N},n \to \infty}} 
\dfrac{1}{\mathcal{N}^{m}} \hat{\mathcal{E}}^{n}_{k} \left(
\sum_{i_{1},i_{2},\dotsc,i_{m}} \phi_{0}(x_{i_{1}},x_{i_{2}},\dotsc,x_{i_{m}}) 
\right) \underset{\underset{z_{o}=1+o(1)}{\mathscr{N},n \to \infty}}{=}& 
\, \lim_{\underset{z_{o}=1+o(1)}{\mathscr{N},n \to \infty}} 
\dfrac{1}{\mathcal{N}^{m}} \idotsint\limits_{\mathbb{R}^{m}} \phi_{0}
(\xi_{1},\xi_{2},\dotsc,\xi_{m}) \hat{\mathfrak{R}}^{n,k}_{m}(\xi_{1},
\xi_{2},\dotsc,\xi_{m}) \, \md \xi_{1} \, \md \xi_{2} \, \dotsb \, \md 
\xi_{m} \nonumber \\
\underset{\underset{z_{o}=1+o(1)}{\mathscr{N},n \to \infty}}{=}& 
\, \idotsint\limits_{\mathbb{R}^{m}} \phi_{0}(\xi_{1},\xi_{2},\dotsc,
\xi_{m}) \, \md \mu_{\widetilde{V}}^{\infty}(\xi_{1}) \, \md 
\mu_{\widetilde{V}}^{\infty}(\xi_{2}) \, \dotsb \, \md 
\mu_{\widetilde{V}}^{\infty}(\xi_{m}), \quad m \! \in \! \mathbb{N},
\end{align}
where, for $m \! \in \! \mathbb{N}$,
\begin{equation} \label{eqlmrtz45} 
\hat{\mathfrak{R}}^{n,k}_{m}(x_{1},x_{2},\dotsc,x_{m}) \! := \! 
\dfrac{\mathcal{N}!}{(\mathcal{N} \! - \! m)!} \idotsint\limits_{
\mathbb{R}^{\mathcal{N}-m}} \hat{\mathcal{P}}^{n}_{k}(x_{1},x_{2},
\dotsc,x_{m},\xi_{m+1},\dotsc,\xi_{\mathcal{N}}) \, \md \xi_{m+1} 
\, \dotsb \, \md \xi_{\mathcal{N}}
\end{equation}
is the $m$-point correlation function,
\begin{equation} \label{eqlmrtz46} 
\hat{\mathcal{P}}^{n,\sharp}_{k}(x_{2},x_{3},\dotsc,x_{\mathcal{N}}) 
\! := \! \dfrac{1}{\hat{\mathscr{Z}}^{n,\sharp}_{k}} \me^{-n \sum_{
m_{1}=2}^{\mathcal{N}} \widetilde{V}(x_{m_{1}})} \prod_{\substack{i,
j=2\\j<i}}^{\mathcal{N}}(x_{j} \! - \! x_{i})^{2} \left(\prod_{m=2}^{
\mathcal{N}} \prod_{q=1}^{\mathfrak{s}-1}(x_{m} \! - \! 
\alpha_{p_{q}})^{\varkappa_{nk \tilde{k}_{q}}} \right)^{-2},
\end{equation}
\begin{equation} \label{eqlmrtz47} 
\hat{\mathscr{Z}}^{n,\sharp}_{k} \! = \! \idotsint\limits_{
\mathbb{R}^{\mathcal{N}-1}} \me^{-n \sum_{m_{1}=2}^{\mathcal{N}} 
\widetilde{V}(\tau_{m_{1}})} \prod_{\substack{i,j=2\\j<i}}^{\mathcal{N}}
(\tau_{j} \! - \! \tau_{i})^{2} \left(\prod_{m=2}^{\mathcal{N}} \prod_{q
=1}^{\mathfrak{s}-1}(\tau_{m} \! - \! \alpha_{p_{q}})^{\varkappa_{nk 
\tilde{k}_{q}}} \right)^{-2} \md \tau_{2} \, \md \tau_{3} 
\, \dotsb \, \md \tau_{\mathcal{N}},
\end{equation}
\begin{equation} \label{eqlmrtz48} 
\dfrac{\hat{\mathscr{Z}}^{n,\sharp}_{k}}{\hat{\mathscr{Z}}^{n}_{k}} 
\underset{\underset{z_{o}=1+o(1)}{\mathscr{N},n \to \infty}}{\leqslant} 
\dfrac{1}{\hat{\mathfrak{z}}} \me^{\frac{n}{\mathcal{N}}(\mathcal{N}-1) 
\mathfrak{c}_{87}(n,k,z_{o})},
\end{equation}
with $\mathfrak{c}_{87}(n,k,z_{o}) \! =_{\underset{z_{o}=1+o(1)}{
\mathscr{N},n \to \infty}} \! \mathcal{O}(1)$,
\begin{equation} \label{eqlmrtz49} 
\pmb{\pi}^{n}_{k}(\lambda) \! = \! \hat{\mathcal{E}}^{n}_{k} 
\left(\lambda^{\varkappa_{nk}} \hat{\mathbb{P}}(x_{1},x_{2},\dotsc,
x_{\mathcal{N}};\lambda) \right),
\end{equation}
where
\begin{align} \label{eqlmrtz50} 
\hat{\mathbb{P}}(x_{1},x_{2},\dotsc,x_{\mathcal{N}};\lambda) :=& \, 
\prod_{j=1}^{\mathcal{N}} \left(1 \! - \! \dfrac{x_{j}}{\lambda} 
\right)^{\frac{\varkappa_{nk}}{\mathcal{N}}} \prod_{q=1}^{\mathfrak{
s}-1} \left(\dfrac{x_{j} \! - \! \alpha_{p_{q}}}{\lambda} \right)^{\frac{
\varkappa_{nk \tilde{k}_{q}}}{\mathcal{N}}} \prod_{q=1}^{\mathfrak{s}
-1} \left(\dfrac{1}{\frac{x_{j}-\alpha_{p_{q}}}{\lambda}} \! - \! \dfrac{
1}{1 \! - \! \frac{\alpha_{p_{q}}}{\lambda}} \right)^{\frac{\varkappa_{n
k \tilde{k}_{q}}}{\mathcal{N}}},
\end{align}
\begin{equation} \label{eqlmrtz51} 
\pmb{\pi}^{n}_{k}(\lambda) \! = \! \lambda^{\varkappa_{nk}} 
\left(\hat{\mathcal{E}}^{n}_{k} \left(\chi_{\hat{\mathfrak{D}}_{
\mathcal{N}}(\hat{M}_{0})}(\vec{x}) \hat{\mathbb{P}}(x_{1},x_{2},
\dotsc,x_{\mathcal{N}};\lambda) \right) \! + \! \hat{\mathcal{E}}^{
n}_{k} \left(\chi_{\hat{\mathfrak{D}}^{c}_{\mathcal{N}}(\hat{M}_{0})}
(\vec{x}) \hat{\mathbb{P}}(x_{1},x_{2},\dotsc,x_{\mathcal{N}};
\lambda) \right) \right),
\end{equation}
\begin{equation} \label{eqlmrtz52} 
\left\lvert \hat{\mathcal{E}}^{n}_{k} \left(\chi_{\hat{\mathfrak{D}}^{
c}_{\mathcal{N}}(\hat{M}_{0})}(\vec{x}) \hat{\mathbb{P}}(x_{1},x_{2},
\dotsc,x_{\mathcal{N}};\lambda) \right) \right\vert \underset{
\underset{z_{o}=1+o(1)}{\mathscr{N},n \to \infty}}{\leqslant} 
\dfrac{\mathfrak{c}_{88}(n,k,z_{o})}{\sqrt{n}} \me^{-\frac{1}{2} 
\frac{n}{\mathcal{N}}(\mathcal{N}-1)}(1 \! + \! o(1)),
\end{equation}
where $\mathfrak{c}_{88}(n,k,z_{o}) \! =_{\underset{z_{o}=1+o(1)}{
\mathscr{N},n \to \infty}} \! \mathcal{O}(1)$, $\hat{\mathfrak{D}}_{
\mathcal{N}}(\hat{M}_{0})$ is the $\mathcal{N}$-fold Cartesian 
product of $\hat{\mathfrak{D}}(\hat{M}_{0}) \! := \! [-\hat{M}_{0},
\hat{M}_{0}] \setminus \cup_{q=1}^{\mathfrak{s}-1} \mathscr{O}_{
\frac{1}{\hat{M}_{0}}}(\alpha_{p_{q}})$, that is, $\hat{\mathfrak{D}}_{
\mathcal{N}}(\hat{M}_{0}) \! = \! \hat{\mathfrak{D}}(\hat{M}_{0}) \! 
\times \! \hat{\mathfrak{D}}(\hat{M}_{0}) \! \times \! \dotsb \! 
\times \! \hat{\mathfrak{D}}(\hat{M}_{0})$, with $\hat{M}_{0} \! = \! 
\hat{M}_{0}(n,k,z_{o}) \! \gg_{\underset{z_{o}=1+o(1)}{\mathscr{N},
n \to \infty}} \! 1$ and bounded, and $\hat{\mathfrak{D}}_{
\mathcal{N}}^{c}(\hat{M}_{0}) \! := \! \mathbb{R}^{\mathcal{N}} 
\setminus \hat{\mathfrak{D}}_{\mathcal{N}}(\hat{M}_{0})$ (the 
complement of $\hat{\mathfrak{D}}_{\mathcal{N}}(\hat{M}_{0})$ 
relative to $\mathbb{R}^{\mathcal{N}})$,
\begin{equation} \label{eqlmrtz53} 
\hat{\mathfrak{g}}(x) \! := \! \hat{F}(x) \ln \left(\left(1 \! - \! 
\dfrac{x}{\lambda} \right)^{\frac{\varkappa_{nk}}{\mathcal{N}}} 
\prod_{q=1}^{\mathfrak{s}-1} \left(\dfrac{x \! - \! \alpha_{p_{q}}}{
\lambda} \right)^{\frac{\varkappa_{nk \tilde{k}_{q}}}{\mathcal{N}}} 
\prod_{q^{\prime}=1}^{\mathfrak{s}-1} \left(\dfrac{1 \! - \! 
\frac{x}{\lambda}}{(1 \! - \! \frac{\alpha_{p_{q^{\prime}}}}{\lambda})
(\frac{x-\alpha_{p_{q^{\prime}}}}{\lambda})} 
\right)^{\frac{\varkappa_{nk \tilde{k}_{q^{\prime}}}}{\mathcal{N}}} 
\right),
\end{equation}
where $\boldsymbol{\mathrm{C}}^{\infty}_{0}(\mathbb{R}) \! \ni \! 
\hat{F} \colon \mathbb{R} \! \to \! [0,1]$ is the test function satisfying 
$0 \! \leqslant \! \hat{F}(x) \! \leqslant \! 1$, $x \! \in \! \mathbb{R}$, 
$\hat{F}(x) \! = \! 1$, $x \! \in \! \hat{\mathfrak{D}}(\hat{M}_{0})$, and 
$\hat{F}(x) \! = \! 0$, $x \! \in \! \hat{\mathfrak{D}}^{c}(\hat{M}_{0} \! 
+ \! 1) \! := \! \lbrace \mathstrut x \! \in \! \mathbb{R}; \, \lvert x 
\rvert \! \geqslant \! \hat{M}_{0} \! + \! 1 \rbrace \cup \cup_{q=1}^{
\mathfrak{s}-1} \lbrace \mathstrut x \! \in \! \mathbb{R}; \, \lvert 
x \! - \! \alpha_{p_{q}} \rvert \! < \! (\hat{M}_{0} \! + \! 1)^{-1} 
\rbrace$,
\begin{equation} \label{eqlmrtz54} 
\pmb{\pi}^{n}_{k}(\lambda) \underset{\underset{z_{o}=1+
o(1)}{\mathscr{N},n \to \infty}}{=} \lambda^{\varkappa_{nk}} \left(
\hat{\mathcal{E}}^{n}_{k} \left(\me^{\sum_{m=1}^{\mathcal{N}} 
\hat{\mathfrak{g}}(x_{m})} \right) \! + \! \mathcal{O} \left(\dfrac{
\mathfrak{c}_{89}(n,k,z_{o})}{\sqrt{n}} \me^{-\frac{1}{2} \frac{n}{
\mathcal{N}}(\mathcal{N}-1)} \right) \right)
\end{equation}
uniformly for, say, $(\mathbb{R} \! \ni)$ $\lvert \lambda \rvert \! 
\geqslant \! 2(\hat{M}_{0} \! + \! 1) \! + \! 3(\mathfrak{s} \! - \! 
1)(1 \! + \! \max_{q=1,2,\dotsc,\mathfrak{s}-1} \lbrace \lvert 
\alpha_{p_{q}} \rvert \rbrace)(\min_{q=1,2,\dotsc,\mathfrak{s}-1} 
\lbrace \hat{\Delta}(q) \rbrace)^{-1}$ $(\gg \! 1)$, where $(0,1/2) 
\! \ni \! \hat{\Delta}(q) \! := \! \lvert \alpha_{p_{q}} \rvert (K(1 \! 
+ \! \max_{q=1,2,\dotsc,\mathfrak{s}-1} \lbrace \lvert \alpha_{
p_{q}} \rvert \rbrace \! + \! 3(\min_{i \neq j \in \lbrace 1,2,\dotsc,
\mathfrak{s}-1 \rbrace} \lbrace \lvert \alpha_{p_{i}} \! - \! \alpha_{
p_{j}} \rvert \rbrace)^{-1}))^{-1}$, $q \! = \! 1,2,\dotsc,\mathfrak{
s} \! - \! 1$, and $\mathfrak{c}_{89}(n,k,z_{o}) \! =_{\underset{z_{o}
=1+o(1)}{\mathscr{N},n \to \infty}} \! \mathcal{O}(1)$,
\begin{equation} \label{eqlmrtz55} 
\hat{\mathcal{E}}^{n}_{k} \left(\me^{\sum_{m=1}^{\mathcal{N}} 
\hat{\mathfrak{g}}(x_{m})} \right) \underset{\underset{z_{o}=1+
o(1)}{\mathscr{N},n \to \infty}}{\geqslant} \exp \left(-\dfrac{n}{8} 
\dfrac{\ln \left(2 \prod_{q=1}^{\mathfrak{s}-1}(1 \! + \! \hat{
\Delta}(q))^{\frac{\gamma_{i(q)_{k_{q}}}}{K}} \right)}{\ln \left(2
(3/2)^{\mathfrak{s}-1} \right)} \right) \underset{\underset{z_{o}=1
+o(1)}{\mathscr{N},n \to \infty}}{\geqslant} \exp \left(-\dfrac{n}{8}
(1 \! - \! \hat{\Gamma}_{0}) \right) > 0,
\end{equation}
where $(0,1) \! \ni \! \hat{\Gamma}_{0} \! := \! (\mathfrak{s} \! - 
\! 1) \gamma_{k} \ln (3/2)(K \ln (2(3/2)^{\mathfrak{s}-1}))^{-1}$.
For $n \! \in \! \mathbb{N}$ and $k \! \in \! \lbrace 1,2,\dotsc,K 
\rbrace$ such that $\alpha_{p_{\mathfrak{s}}} \! := \! \alpha_{k} \! 
= \! \infty$, let $\hat{\mathfrak{Z}}_{0} \! := \! \lbrace \mathstrut 
\lambda \! \in \! \overline{\mathbb{C}}; \, \pmb{\pi}^{n}_{k}(\lambda) 
\! = \! 0 \rbrace$: Equations~\eqref{eqlmrtz54} and~\eqref{eqlmrtz55} 
show that, in the double-scaling limit $\mathscr{N},n \! \to \! \infty$ 
such that $z_{o} \! = \! 1 \! + \! o(1)$, $\hat{\mathfrak{Z}}_{0} \cap 
(\lbrace \alpha_{p_{1}},\alpha_{p_{2}},\dotsc,\alpha_{p_{\mathfrak{s}-1}} 
\rbrace \cup \lbrace \alpha_{p_{\mathfrak{s}}} \rbrace) \! = \! 
\varnothing$, that is (cf. Remark~\ref{rem1.2.3}), $\hat{\mathfrak{Z}}_{0} 
\! = \! \lbrace \hat{\mathfrak{z}}^{n}_{k}(j) \rbrace_{j=1}^{\mathcal{N}}$, 
with $\hat{\mathfrak{z}}^{n}_{k}(j) \! \neq \! \alpha_{p_{q}}$, $j \! = \! 1,2,
\dotsc,\mathcal{N}$, $q \! = \! 1,2,\dotsc,\mathfrak{s}$. In conjunction 
with Equations~\eqref{eqlmrtz54} and~\eqref{eqlmrtz55}, this implies, in 
particular, that, for $n \! \in \! \mathbb{N}$ and $k \! \in \! \lbrace 1,2,
\dotsc,K \rbrace$ such that $\alpha_{p_{\mathfrak{s}}} \! := \! \alpha_{k} 
\! = \! \infty$, in the double-scaling limit $\mathscr{N},n \! \to \! \infty$ 
such that $z_{o} \! = \! 1 \! + \! o(1)$, there exists $\hat{\lambda}_{1} 
\! = \! \hat{\lambda}_{1}(n,k,z_{o}) \! > \! 0$ and $\mathcal{O}(1)$ 
{}\footnote{For $k \! \in \! \lbrace 1,2,\dotsc,K \rbrace$ such that 
$\alpha_{p_{\mathfrak{s}}} \! := \! \alpha_{k} \! = \! \infty$, take, say, 
$\hat{\lambda}_{1} \! = \! \max \left\lbrace \hat{\lambda}_{1}^{\ast},
2(\hat{M}_{0} \! + \! 1) \! + \! 3(\mathfrak{s} \! - \! 1)(1 \! + \! 
\max_{q=1,2,\dotsc,\mathfrak{s}-1} \lbrace \lvert \alpha_{p_{q}} \rvert 
\rbrace)(\min_{q=1,2,\dotsc,\mathfrak{s}-1} \lbrace \hat{\Delta}(q) 
\rbrace)^{-1} \right\rbrace$ $(\gg \! 1)$, where
\begin{equation*}
\hat{\lambda}_{1}^{\ast} \! := \! (\hat{M}_{0} \! + \! 1) \left(1 \! 
- \! \exp \left(-\left\lvert \dfrac{\ln \left(2 \prod_{q=1}^{\mathfrak{s}
-1}(1 \! + \! \hat{\Delta}(q))^{\frac{\gamma_{i(q)_{k_{q}}}}{K}} \right)}{8 
\mathfrak{s} \varpi_{0}^{\flat} \ln (2(3/2)^{\mathfrak{s}-1}) \max_{q=
1,2,\dotsc,\mathfrak{s}-1} \lbrace \gamma_{i(q)_{k_{q}}}\rbrace} 
\right\rvert \right) \right)^{-1},
\end{equation*}
with $\varpi_{0}^{\flat}$ defined via the following inequality: for 
$x \! = \! \hat{M}_{0} \! + \! 1$, $\hat{\lambda}_{0} \! > \! 
2(\hat{M}_{0} \! + \! 1) \! + \! 3(\mathfrak{s} \! - \! 1)(1 \! + \! 
\max_{q=1,2,\dotsc,\mathfrak{s}-1} \lbrace \lvert \alpha_{p_{q}} 
\rvert \rbrace)(\min_{q=1,2,\dotsc,\mathfrak{s}-1} \lbrace 
\hat{\Delta}(q) \rbrace)^{-1}$ $(\gg \! 1)$, and $\lvert \lambda 
\rvert \! \geqslant \! \hat{\lambda}_{0}$, one shows that
\begin{equation*}
\left. \ln \hat{\mathbb{G}}(x) \vphantom{M^{M^{M^{M}}}} 
\right\vert_{\underset{\lambda = \lvert \lambda \rvert}{x=
\hat{M}_{0}+1}} \underset{\underset{z_{o}=1+o(1)}{\mathscr{N},n 
\to \infty}}{\geqslant} \varpi_{0}^{\flat} \left. \ln \hat{\mathbb{G}}
(x) \vphantom{M^{M^{M^{M}}}} \right\vert_{\underset{\lambda 
= \hat{\lambda}_{0}}{x=\hat{M}_{0}+1}},
\end{equation*} 
where
\begin{equation*}
\hat{\mathbb{G}}(x) \! := \! \left(1 \! - \! \dfrac{x}{\lambda} \right)^{
\frac{\varkappa_{nk}}{\mathcal{N}}} \prod_{q=1}^{\mathfrak{s}-1} \left(
\dfrac{x \! - \! \alpha_{p_{q}}}{\lambda} \right)^{\frac{\varkappa_{nk 
\tilde{k}_{q}}}{\mathcal{N}}} \prod_{q^{\prime}=1}^{\mathfrak{s}-1} 
\left(\dfrac{1 \! - \! \frac{x}{\lambda}}{(1 \! - \! \frac{\alpha_{p_{
q^{\prime}}}}{\lambda})(\frac{x-\alpha_{p_{q^{\prime}}}}{\lambda})} 
\right)^{\frac{\varkappa_{nk \tilde{k}_{q^{\prime}}}}{\mathcal{N}}}.
\end{equation*}} such that $\pmb{\pi}^{n}_{k}(\lambda)$ has no 
zeros (counting multiplicities) in $\overline{\mathbb{R}} \setminus 
\hat{\mathfrak{D}}(\hat{\lambda}_{1})$, where $\hat{\mathfrak{D}}
(\hat{\lambda}_{1}) \! := \! [-\hat{\lambda}_{1},\hat{\lambda}_{1}] 
\setminus \cup_{q=1}^{\mathfrak{s}-1} \mathscr{O}_{\frac{1}{
\hat{\lambda}_{1}}}(\alpha_{p_{q}})$, that is, $\hat{\mathfrak{Z}}_{0} 
\subseteq \hat{\mathfrak{D}}(\hat{\lambda}_{1})$. \hfill $\qed$

The following corollary establishes, for $n \! \in \! \mathbb{N}$ and 
$k \! \in \! \lbrace 1,2,\dotsc,K \rbrace$ such that $\alpha_{p_{
\mathfrak{s}}} \! := \! \alpha_{k} \! = \! \infty$ (resp., $\alpha_{p_{
\mathfrak{s}}} \! := \! \alpha_{k} \! \neq \! \infty)$, in the double-scaling 
limit $\mathscr{N},n \! \to \! \infty$ such that $z_{o} \! = \! 1 \! + \! 
o(1)$, the existence of the corresponding monic MPC ORF `Free Energies'.
\begin{fffff} \label{corol3.1} 
Let the external field $\widetilde{V} \colon \overline{\mathbb{R}} 
\setminus \lbrace \alpha_{1},\alpha_{2},\dotsc,\alpha_{K} \rbrace 
\! \to \! \mathbb{R}$ satisfy conditions~\eqref{eq20}--\eqref{eq22} 
and be regular. For $n \! \in \! \mathbb{N}$ and $k \! \in \! \lbrace 
1,2,\dotsc,K \rbrace$ such that $\alpha_{p_{\mathfrak{s}}} \! := \! 
\alpha_{k} \! = \! \infty$ (resp., $\alpha_{p_{\mathfrak{s}}} \! := \! 
\alpha_{k} \! \neq \! \infty)$, let the associated equilibrium measure, 
$\mu_{\widetilde{V}}^{\infty}$ (resp., $\mu_{\widetilde{V}}^{f})$, and its 
support, $J_{\infty}$ (resp., $J_{f})$, be as described in item~$\pmb{(1)}$ 
(resp., item~$\pmb{(2)})$ of Lemma~\ref{lem3.7}. Then, for $n \! \in \! 
\mathbb{N}$ and $k \! \in \! \lbrace 1,2,\dotsc,K \rbrace$ such that 
$\alpha_{p_{\mathfrak{s}}} \! := \! \alpha_{k} \! = \! \infty$,
\begin{equation} \label{leedordvarinf} 
\lim_{\underset{z_{o}=1+o(1)}{\mathscr{N},n \to \infty}} \left(
-\dfrac{\mathcal{N}}{n} \dfrac{1}{\mathcal{N}(\mathcal{N} \! - \! 1)} 
\ln (\hat{\mathscr{Z}}^{n}_{k}) \right) \underset{\underset{z_{o}=
1+o(1)}{\mathscr{N},n \to \infty}}{=} E_{\widetilde{V}}^{\infty},
\end{equation}
with $\mathcal{N} \! = \! (n \! - \! 1)K \! + \! k$, where 
$\hat{\mathscr{Z}}^{n}_{k}$ is defined by Equation~\eqref{eqlmrtz40}, and, 
for $n \! \in \! \mathbb{N}$ and $k \! \in \! \lbrace 1,2,\dotsc,K \rbrace$ 
such that $\alpha_{p_{\mathfrak{s}}} \! := \! \alpha_{k} \! \neq \! \infty$,
\begin{equation} \label{leedordvarfin} 
\lim_{\underset{z_{o}=1+o(1)}{\mathscr{N},n \to \infty}} \left(
-\dfrac{\mathcal{N}}{n} \dfrac{1}{\mathcal{N}(\mathcal{N} \! - \! 1)} 
\ln (\tilde{\mathscr{Z}}^{n}_{k}) \right) \underset{\underset{z_{o}=
1+o(1)}{\mathscr{N},n \to \infty}}{=} E_{\widetilde{V}}^{f},
\end{equation}
where $\tilde{\mathscr{Z}}^{n}_{k}$ is defined by Equation~\eqref{eqlmrtzz}.
\end{fffff}

\emph{Proof}. The proof of this Corollary~\ref{corol3.1} consists of two 
cases: (i) $n \! \in \! \mathbb{N}$ and $k \! \in \! \lbrace 1,2,\dotsc,K 
\rbrace$ such that $\alpha_{p_{\mathfrak{s}}} \! := \! \alpha_{k} \! = \! 
\infty$; and (ii) $n \! \in \! \mathbb{N}$ and $k \! \in \! \lbrace 1,2,
\dotsc,K \rbrace$ such that $\alpha_{p_{\mathfrak{s}}} \! := \! \alpha_{k} 
\! \neq \! \infty$. Notwithstanding the fact that the scheme of the proof 
is, \emph{mutatis mutandis}, similar for both cases, case~(ii) is the more 
technically challenging of the two; therefore, without loss of generality, 
its particulars are presented in detail (see \pmb{(1)} below), whilst 
case~(i) is proved analogously (see \pmb{(2)} below).

\pmb{(1)} For $n \! \in \! \mathbb{N}$ and $k \! \in \! \lbrace 1,2,
\dotsc,K \rbrace$ such that $\alpha_{p_{\mathfrak{s}}} \! := \! 
\alpha_{k} \! \neq \! \infty$, it was shown in  Equation~\eqref{eqlmrtz6} 
that, for (arbitrary) $\tilde{\varepsilon} \! > \! 0$,
\begin{equation} \label{eqlmcorl1} 
\limsup_{\underset{z_{o}=1+o(1)}{\mathscr{N},n \to \infty}} \left(
-\dfrac{\mathcal{N}}{n} \dfrac{1}{\mathcal{N}(\mathcal{N} \! - \! 1)} 
\ln (\tilde{\mathscr{Z}}^{n}_{k}) \right) \underset{\underset{z_{o}=
1+o(1)}{\mathscr{N},n \to \infty}}{\leqslant} E_{\widetilde{V}}^{f} \! 
+ \! \tilde{\varepsilon};
\end{equation}
furthermore, it was shown in the proof of Lemma~\ref{lemrootz} 
(cf. case~\pmb{(1)}) that
\begin{equation} \label{eqlmcorl2} 
\tilde{\mathscr{Z}}^{n}_{k} \underset{\underset{z_{o}=1+o(1)}{
\mathscr{N},n \to \infty}}{=} \idotsint\limits_{\mathbb{R}^{\mathcal{N}}} 
\me^{-\frac{n}{\mathcal{N}}\left(\sum_{j=1}^{\mathcal{N}} \widetilde{V}
(\xi_{j})+ \mathscr{K}^{\widetilde{V},f}_{\mathcal{N}}(\xi_{1},\xi_{2},
\dotsc,\xi_{\mathcal{N}})+\mathcal{O}(\frac{1}{\mathcal{N}}) 
\sum_{\underset{i \neq j}{i,j=1}}^{\mathcal{N}} \tilde{\mathrm{F}}
(\xi_{i},\xi_{j}) \right)} \, \md \xi_{1} \, \md \xi_{2} \, \dotsb 
\, \md \xi_{\mathcal{N}},
\end{equation}
where the symmetric function $\tilde{\mathrm{F}}(x,y)$ is defined by 
Equation~\eqref{eqlmrtzf}. It follows {}from Equation~\eqref{eql3.5e} 
that $\mathfrak{d}_{\mathcal{N}}^{\widetilde{V},f} \! \leqslant \! 
(\mathcal{N}(\mathcal{N} \! - \! 1))^{-1} \mathscr{K}_{\mathcal{N}}^{
\widetilde{V},f}(x_{1},x_{2},\dotsc,x_{\mathcal{N}})$, which, via 
the monotonicity of $\exp (\pmb{\cdot})$, implies that $\exp 
(-\tfrac{n}{\mathcal{N}} \mathscr{K}_{\mathcal{N}}^{\widetilde{V},f}
(x_{1},x_{2},\dotsc,x_{\mathcal{N}})) \linebreak[4] 
\! \leqslant \! \exp (-\tfrac{n}{\mathcal{N}} \mathcal{N}(\mathcal{N} 
\! - \! 1) \mathfrak{d}_{\mathcal{N}}^{\widetilde{V},f})$; hence, via 
this latter relation and Equation~\eqref{eqlmcorl2}, one shows that
\begin{equation} \label{eqlmcorl3} 
\tilde{\mathscr{Z}}^{n}_{k} \underset{\underset{z_{o}=1+
o(1)}{\mathscr{N},n \to \infty}}{\leqslant} \me^{-\frac{n}{
\mathcal{N}} \mathcal{N}(\mathcal{N}-1) \mathfrak{d}_{
\mathcal{N}}^{\widetilde{V},f}} \idotsint\limits_{\mathbb{R}^{
\mathcal{N}}} \me^{-\frac{n}{\mathcal{N}}\left(\sum_{j=1}^{
\mathcal{N}} \widetilde{V}(\xi_{j})+\mathcal{O}(\frac{1}{\mathcal{N}}) 
\sum_{\underset{i \neq j}{i,j=1}}^{\mathcal{N}} \tilde{\mathrm{F}}
(\xi_{i},\xi_{j}) \right)} \, \md \xi_{1} \, \md \xi_{2} \, \dotsb 
\, \md \xi_{\mathcal{N}}.
\end{equation}
It was also shown in the proof of Lemma~\ref{lemrootz} 
(cf. case~\pmb{(1)}, the calculations leading to the 
Estimate~\eqref{eqlmrtz7}) that
\begin{align*}
\idotsint\limits_{\mathbb{R}^{\mathcal{N}}} \me^{-\frac{n}{\mathcal{N}} 
\left(\sum_{j=1}^{\mathcal{N}} \widetilde{V}(\xi_{j})+\mathcal{O}
(\frac{1}{\mathcal{N}}) \sum_{\underset{i \neq j}{i,j=1}}^{\mathcal{N}} 
\tilde{\mathrm{F}}(\xi_{i},\xi_{j}) \right)} \, \md \xi_{1} \, \md \xi_{2} 
\, \dotsb \, \md \xi_{\mathcal{N}} \underset{\underset{z_{o}=1+o
(1)}{\mathscr{N},n \to \infty}}{=}& \, \left(\int_{\mathbb{R}} \dfrac{
\me^{-\frac{n}{\mathcal{N}} \widetilde{V}(\tau)}}{\prod_{\underset{
q \neq \mathfrak{s}-1}{q=1}}^{\mathfrak{s}} \lvert \tau \! - \! 
\alpha_{p_{q}} \rvert^{\tilde{\gamma}_{q}(1+\mathcal{O}(n^{-1}))}} 
\, \md \tau \right)^{\mathcal{N}} \\
& \, \underset{\underset{z_{o}=1+o(1)}{\mathscr{N},n \to 
\infty}}{\leqslant}(\mathfrak{c}_{16}(n,k,z_{o}))^{\mathcal{N}},
\end{align*}
where $\mathfrak{c}_{16}(n,k,z_{o}) \! =_{\underset{z_{o}=1+
o(1)}{\mathscr{N},n \to \infty}} \! \mathcal{O}(1)$, which, via 
Equation~\eqref{eqlmcorl3}, implies that
\begin{equation} \label{eqlmcorl4} 
-\dfrac{\mathcal{N}}{n} \dfrac{1}{\mathcal{N}(\mathcal{N} \! - \! 
1)} \ln (\tilde{\mathscr{Z}}^{n}_{k}) \underset{\underset{z_{o}=1+
o(1)}{\mathscr{N},n \to \infty}}{\geqslant} \mathfrak{d}_{\mathcal{
N}}^{\widetilde{V},f} \! - \! \dfrac{\mathcal{N}}{n} \dfrac{1}{
(\mathcal{N} \! - \! 1)} \ln (\mathfrak{c}_{16}(n,k,z_{o})).
\end{equation}
Using the fact that (cf. Lemma~\ref{lem3.5}), for $n \! \in \! 
\mathbb{N}$ and $k \! \in \! \lbrace 1,2,\dotsc,K \rbrace$ 
such that $\alpha_{p_{\mathfrak{s}}} \! := \! \alpha_{k} \! \neq \! 
\infty$, $\lim_{\underset{z_{o}=1+o(1)}{\mathscr{N},n \to \infty}} 
\mathfrak{d}_{\mathcal{N}}^{\widetilde{V},f} \! = \! E_{\widetilde{
V}}^{f}$, it follows {}from Equation~\eqref{eqlmcorl4} that
\begin{equation} \label{eqlmcorl5} 
\liminf_{\underset{z_{o}=1+o(1)}{\mathscr{N},n \to \infty}} \left(
-\dfrac{\mathcal{N}}{n} \dfrac{1}{\mathcal{N}(\mathcal{N} \! - \! 1)} 
\ln (\tilde{\mathscr{Z}}^{n}_{k}) \right) \underset{\underset{z_{o}=
1+o(1)}{\mathscr{N},n \to \infty}}{\geqslant} E_{\widetilde{V}}^{f};
\end{equation}
via the $\tilde{\varepsilon} \! \downarrow \! 0$ limit of 
Equation~\eqref{eqlmcorl1} (since $\tilde{\varepsilon} \! > \! 0$ is 
arbitrary), and Equation~\eqref{eqlmcorl5}, one arrives at, for $n \! \in 
\! \mathbb{N}$ and $k \! \in \! \lbrace 1,2,\dotsc,K \rbrace$ such 
that $\alpha_{p_{\mathfrak{s}}} \! := \! \alpha_{k} \! \neq \! \infty$, 
in the double-scaling limit $\mathscr{N},n \! \to \! \infty$ such that 
$z_{o} \! = \! 1 \! + \! o(1)$, the corresponding result stated in 
Equation~\eqref{leedordvarfin}.

\pmb{(2)} The proof of this case, that is, $n \! \in \! \mathbb{N}$ 
and $k \! \in \! \lbrace 1,2,\dotsc,K \rbrace$ such that $\alpha_{
p_{\mathfrak{s}}} \! := \! \alpha_{k} \! = \! \infty$, is virtually identical 
to the proof of \pmb{(1)} above; one mimics, \emph{verbatim}, the 
scheme of the calculations presented in \pmb{(1)} in order to 
arrive at the corresponding claim stated in the corollary; in order 
to do so, however, the analogues of Equations~\eqref{eqlmcorl1} 
and~\eqref{eqlmcorl4}, respectively, are necessary, which, in the 
present case, read (cf. the proof of Lemma~\ref{lemrootz}, 
case~\pmb{(2)}): for (arbitrary) $\hat{\varepsilon} \! > \! 0$,
\begin{equation} \label{eqlmcorl6} 
\limsup_{\underset{z_{o}=1+o(1)}{\mathscr{N},n \to \infty}} 
\left(-\dfrac{\mathcal{N}}{n} \dfrac{1}{\mathcal{N}(\mathcal{N} 
\! - \! 1)} \ln (\hat{\mathscr{Z}}^{n}_{k}) \right) \underset{\underset{
z_{o}=1+o(1)}{\mathscr{N},n \to \infty}}{\leqslant} E_{\widetilde{V}}^{
\infty} \! + \! \hat{\varepsilon},
\end{equation}
and
\begin{equation} \label{eqlmcorl7} 
-\dfrac{\mathcal{N}}{n} \dfrac{1}{\mathcal{N}(\mathcal{N} \! - \! 
1)} \ln (\hat{\mathscr{Z}}^{n}_{k}) \underset{\underset{z_{o}=1+
o(1)}{\mathscr{N},n \to \infty}}{\geqslant} \mathfrak{d}_{\mathcal{
N}}^{\widetilde{V},\infty} \! - \! \dfrac{\mathcal{N}}{n} \dfrac{1}{
(\mathcal{N} \! - \! 1)} \ln (\mathfrak{c}_{\blacklozenge}^{\lozenge}
(n,k,z_{o})),\end{equation}
where $\mathfrak{c}_{\blacklozenge}^{\lozenge}(n,k,z_{o}) \! 
=_{\underset{z_{o}=1+o(1)}{\mathscr{N},n \to \infty}} \! \mathcal{O}
(1)$. Using the fact that (cf. Lemma~\ref{lem3.5}), for $n \! \in \! 
\mathbb{N}$ and $k \! \in \! \lbrace 1,2,\dotsc,K \rbrace$ such that 
$\alpha_{p_{\mathfrak{s}}} \! := \! \alpha_{k} \! = \! \infty$, $\lim_{
\underset{z_{o}=1+o(1)}{\mathscr{N},n \to \infty}} \mathfrak{d}_{
\mathcal{N}}^{\widetilde{V},\infty} \! = \! E_{\widetilde{V}}^{\infty}$, 
it follows {}from Equation~\eqref{eqlmcorl7} that
\begin{equation} \label{eqlmcorl8} 
\liminf_{\underset{z_{o}=1+o(1)}{\mathscr{N},n \to \infty}} 
\left(-\dfrac{\mathcal{N}}{n} \dfrac{1}{\mathcal{N}(\mathcal{N} 
\! - \! 1)} \ln (\hat{\mathscr{Z}}^{n}_{k}) \right) \underset{\underset{z_{o}
=1+o(1)}{\mathscr{N},n \to \infty}}{\geqslant} E_{\widetilde{V}}^{\infty};
\end{equation}
hence, via the $\hat{\varepsilon} \! \downarrow \! 0$ limit of 
Equation~\eqref{eqlmcorl6} (since $\hat{\varepsilon} \! > \! 0$ is 
arbitrary), and Equation~\eqref{eqlmcorl8}, one arrives at, for $n \! \in 
\! \mathbb{N}$ and $k \! \in \! \lbrace 1,2,\dotsc,K \rbrace$ such that 
$\alpha_{p_{\mathfrak{s}}} \! := \! \alpha_{k} \! = \! \infty$, in the 
double-scaling limit $\mathscr{N},n \! \to \! \infty$ such that 
$z_{o} \! = \! 1 \! + \! o(1)$, the corresponding result stated in 
Equation~\eqref{leedordvarinf}. \hfill $\qed$
\begin{eeeee} \label{freeenergyzeroord} 
\emph{To leading order in the double-scaling limit $\mathscr{N},n \! 
\to \! \infty$ such that $z_{o} \! = \! 1 \! + \! o(1)$, the right-hand side 
of Equation~\eqref{leedordvarinf} (resp., \eqref{leedordvarfin}$)$ reads 
$E_{\widetilde{V}}^{\infty} \! = \! \inf \lbrace \mathstrut \mathfrak{X}^{
\infty}_{\widetilde{V}}[\mu^{\text{\tiny {\rm EQ}}}]; \, \mu^{\text{\tiny {\rm EQ}}} 
\! \in \! \mathscr{M}_{1}(\mathbb{R}) \rbrace$ (resp., $E_{\widetilde{V}}^{f} 
\! = \! \inf \lbrace \mathstrut \mathfrak{X}^{f}_{\widetilde{V}}
[\mu^{\text{\tiny {\rm EQ}}}]; \, \mu^{\text{\tiny {\rm EQ}}} \! \in \! 
\mathscr{M}_{1}(\mathbb{R}) \rbrace)$, where $\mathfrak{X}^{\infty}_{
\widetilde{V}}[\mu^{\text{\tiny {\rm EQ}}}]$ (resp., $\mathfrak{X}^{f}_{
\widetilde{V}}[\mu^{\text{\tiny {\rm EQ}}}])$ is given in Equation~\eqref{stochinf} 
(resp., \eqref{stochfin}$)$$;$ recall, also, that (cf. Remark~\ref{rem1.3.4}$)$ 
$\mathfrak{X}^{\infty}_{\widetilde{V}}[\mu^{\text{\tiny {\rm EQ}}}] \! = \! 
\mathfrak{X}^{f}_{\widetilde{V}}[\mu^{\text{\tiny {\rm EQ}}}]$.}
\end{eeeee}
The following---crucial---Lemma~\ref{lemetatomu} establishes, for 
$n \! \in \! \mathbb{N}$ and $k \! \in \! \lbrace 1,2,\dotsc,K \rbrace$ 
such that $\alpha_{p_{\mathfrak{s}}} \! := \! \alpha_{k} \! = \! \infty$ 
(resp., $\alpha_{p_{\mathfrak{s}}} \! := \! \alpha_{k} \! \neq \! \infty)$, 
the weak-$\ast$ convergence, in the double-scaling limit $\mathscr{N},
n \! \to \! \infty$ such that $z_{o} \! = \! 1 \! + \! o(1)$, of the 
associated monic MPC ORF normalised zero counting measures.
\begin{ccccc} \label{lemetatomu} 
Let the external field $\widetilde{V} \colon \overline{\mathbb{R}} 
\setminus \lbrace \alpha_{1},\alpha_{2},\dotsc,\alpha_{K} \rbrace 
\! \to \! \mathbb{R}$ satisfy conditions~\eqref{eq20}--\eqref{eq22} 
and be regular. For $n \! \in \! \mathbb{N}$ and $k \! \in \! \lbrace 
1,2,\dotsc,K \rbrace$ such that $\alpha_{p_{\mathfrak{s}}} \! := \! 
\alpha_{k} \! = \! \infty$ (resp., $\alpha_{p_{\mathfrak{s}}} \! := \! 
\alpha_{k} \! \neq \! \infty)$, let the associated equilibrium measure, 
$\mu_{\widetilde{V}}^{\infty}$ (resp., $\mu_{\widetilde{V}}^{f})$, and its 
support, $J_{\infty}$ (resp., $J_{f})$, be as described in item~$\pmb{(1)}$ 
(resp., item~$\pmb{(2)})$ of Lemma~\ref{lem3.7}. For $n \! \in \! 
\mathbb{N}$ and $k \! \in \! \lbrace 1,2,\dotsc,K \rbrace$ such 
that $\alpha_{p_{\mathfrak{s}}} \! := \! \alpha_{k} \! = \! \infty$, let 
$\hat{\mathfrak{Z}}_{0} \! := \! \lbrace \mathstrut \lambda \! \in \! 
\overline{\mathbb{C}}; \, \pmb{\pi}^{n}_{k}(\lambda) \! = \! 0 \rbrace \! 
= \! \lbrace \hat{\mathfrak{z}}^{n}_{k}(j) \rbrace_{j=1}^{\mathcal{N}}$, 
where, with $\mathcal{N} \! = \! (n \! - \! 1)K \! + \! k$, $\lbrace 
\hat{\mathfrak{z}}^{n}_{k}(j) \rbrace_{j=1}^{\mathcal{N}}$ is described 
in the corresponding item of Lemma~\ref{lemrootz}, and let the 
associated normalised zero counting measure be defined as
\begin{equation*}
\hat{\eta}_{\hat{\mathfrak{z}}}(x) \! := \! \dfrac{1}{\mathcal{N}} 
\sum_{j=1}^{\mathcal{N}} \delta_{\hat{\mathfrak{z}}^{n}_{k}(j)}(x),
\end{equation*}
where $\delta_{\hat{\mathfrak{z}}^{n}_{k}(j)}(x)$ is the Dirac delta (atomic) 
mass concentrated at $\hat{\mathfrak{z}}^{n}_{k}(j)$, $j \! = \! 1,2,\dotsc,
\mathcal{N}$, and, for $n \! \in \! \mathbb{N}$ and $k \! \in \! \lbrace 
1,2,\dotsc,K \rbrace$ such that $\alpha_{p_{\mathfrak{s}}} \! := \! 
\alpha_{k} \! \neq \! \infty$, let $\tilde{\mathfrak{Z}}_{0} \! := \! \lbrace 
\mathstrut \lambda \! \in \! \overline{\mathbb{C}}; \, \pmb{\pi}^{n}_{k}
(\lambda) \! = \! 0 \rbrace \! = \! \lbrace \tilde{\mathfrak{z}}^{n}_{k}
(j) \rbrace_{j=1}^{\mathcal{N}}$, where $\lbrace \tilde{\mathfrak{z}}^{n}_{k}
(j) \rbrace_{j=1}^{\mathcal{N}}$ is described in the corresponding item 
of Lemma~\ref{lemrootz}, and let the associated normalised zero 
counting measure be defined as
\begin{equation*}
\tilde{\eta}_{\tilde{\mathfrak{z}}}(x) \! := \! \dfrac{1}{\mathcal{N}} 
\sum_{j=1}^{\mathcal{N}} \delta_{\tilde{\mathfrak{z}}^{n}_{k}(j)}(x),
\end{equation*}
where $\delta_{\tilde{\mathfrak{z}}^{n}_{k}(j)}(x)$ is the Dirac delta (atomic) 
mass concentrated at $\tilde{\mathfrak{z}}^{n}_{k}(j)$, $j \! = \! 1,2,
\dotsc,\mathcal{N}$. Then, for $n \! \in \! \mathbb{N}$ and $k \! \in 
\! \lbrace 1,2,\dotsc,K \rbrace$ such that $\alpha_{p_{\mathfrak{s}}} 
\! := \! \alpha_{k} \! = \! \infty$ (resp., $\alpha_{p_{\mathfrak{s}}} 
\! := \! \alpha_{k} \! \neq \! \infty)$, $\hat{\eta}_{\hat{\mathfrak{z}}} 
\! \overset{\ast}{\to} \! \mu_{\widetilde{V}}^{\infty}$ (resp., $\tilde{
\eta}_{\tilde{\mathfrak{z}}} \! \overset{\ast}{\to} \! \mu_{\widetilde{
V}}^{f})$ in the double-scaling limit $\mathscr{N},n \! \to \! \infty$ 
such that $z_{o} \! = \! 1 \! + \! o(1)$.
\end{ccccc}

\emph{Proof}. The proof of this Lemma~\ref{lemetatomu} consists 
of two cases: (i) $n \! \in \! \mathbb{N}$ and $k \! \in \! \lbrace 1,
2,\dotsc,K \rbrace$ such that $\alpha_{p_{\mathfrak{s}}} \! := \! 
\alpha_{k} \! = \! \infty$; and (ii) $n \! \in \! \mathbb{N}$ and $k \! \in 
\! \lbrace 1,2,\dotsc,K \rbrace$ such that $\alpha_{p_{\mathfrak{s}}} 
\! := \! \alpha_{k} \! \neq \! \infty$. Notwithstanding the fact that the 
scheme of the proof is, \emph{mutatis mutandis}, similar for both cases, 
case~(ii), nonetheless, is the more technically challenging of the two; 
therefore, without loss of generality, its particulars are presented in 
detail (see \pmb{(1)} below), whilst case~(i) is proved analogously (see 
\pmb{(2)} below).

\pmb{(1)} It was shown in the proof of Lemma~\ref{lemrootz} (cf. 
Equations~\eqref{eqlmrtz37} and~\eqref{eqlmrtz38}) that, for $n \! 
\in \! \mathbb{N}$ and $k \! \in \! \lbrace 1,2,\dotsc,K \rbrace$ 
such that $\alpha_{p_{\mathfrak{s}}} \! := \! \alpha_{k} \! \neq \! 
\infty$, for $(\mathbb{R} \! \ni)$ $\lvert \lambda \rvert \! > \! 
\tilde{\lambda}_{1} \! := \! \max \lbrace \tilde{\lambda}_{1}^{\ast},
2(\tilde{M}_{0} \! + \! 1) \! + \! 3(\mathfrak{s} \! - \! 1)(1 \! + \! 
\max_{q=1,\dotsc,\mathfrak{s}-2,\mathfrak{s}} \lbrace \lvert 
\alpha_{p_{q}} \rvert \rbrace)(\min_{q=1,\dotsc,\mathfrak{s}-2,
\mathfrak{s}} \lbrace \tilde{\Delta}(q) \rbrace)^{-1} \rbrace$ $(\gg \! 1)$,
\begin{equation*}
\tilde{\mathcal{E}}^{n}_{k} \left(\me^{\sum_{m=1}^{\mathcal{N}} 
\tilde{\mathfrak{g}}(x_{m})} \right) \underset{\underset{z_{o}=1+
o(1)}{\mathscr{N},n \to \infty}}{\geqslant} \exp \left(-\dfrac{n}{8}
(1 \! - \! \tilde{\Gamma}_{r}) \right) >0, \quad r \! \in \! \lbrace 
\mathrm{A},\mathrm{B} \rbrace,
\end{equation*}
where $\mathcal{N} \! = \! (n \! - \! 1)K \! + \! k$, $\tilde{
\mathfrak{g}}(x)$ is given in Equations~\eqref{eqlmrtz31} 
and~\eqref{eqlmrtz32}, $r \! = \! \mathrm{A} \! \leftrightarrow \! 
\min_{q=1,2,\dotsc,\mathfrak{s}-2} \lbrace \lvert \alpha_{k} \! - 
\! \alpha_{p_{q}} \rvert \rbrace \! \in \! (0,1)$, and $r \! = \! 
\mathrm{B} \! \leftrightarrow \! \min_{q=1,2,\dotsc,\mathfrak{s}-2} 
\lbrace \lvert \alpha_{k} \! - \! \alpha_{p_{q}} \rvert \rbrace \! \in \! 
[1,+\infty)$, which, via Equation~\eqref{eqlmrtz34}, implies that, 
for $\lvert \lambda \rvert \! > \! \tilde{\lambda}_{1}$,
\begin{equation} \label{eqlmeta1} 
\dfrac{1}{\mathcal{N}} \ln \pmb{\pi}^{n}_{k}(\lambda) 
\underset{\underset{z_{o}=1+o(1)}{\mathscr{N},n \to \infty}}{=} 
\dfrac{1}{\mathcal{N}} \ln \tilde{\mathcal{E}}^{n}_{k} \left(\me^{
\sum_{m=1}^{\mathcal{N}} \tilde{\mathfrak{g}}_{\blacklozenge}
(x_{m})} \right) \! + \! \mathcal{O} \left(\dfrac{\mathfrak{c}_{
\blacklozenge}^{1}(n,k,z_{o})}{n^{3/2}} \me^{-\frac{3n}{8}} 
\me^{-\frac{n \tilde{\Gamma}_{r}}{8}} \right), \quad r \! \in \! 
\lbrace \mathrm{A},\mathrm{B} \rbrace,
\end{equation}
where
\begin{equation*}
\tilde{\mathfrak{g}}_{\blacklozenge}(x) \! := \! \tilde{\mathfrak{g}}
(x) \! + \! \dfrac{\varkappa_{nk \tilde{k}_{\mathfrak{s}-1}}^{
\infty}}{\mathcal{N}} \ln \lambda,
\end{equation*}
and $\mathfrak{c}_{\blacklozenge}^{1}(n,k,z_{o}) \! =_{\underset{
z_{o}=1+o(1)}{\mathscr{N},n \to \infty}} \! \mathcal{O}(1)$; now, 
noting that, as a function of the independent variable $x$, $\tilde{
\mathfrak{g}}_{\blacklozenge}(x)$ is a real-valued, bounded, 
and continuous function on $\mathbb{R}$ (cf. the proof of 
case~\pmb{(1)} of Lemma~\ref{lemrootz}), it follows via 
Equations~\eqref{eqlmrtz8} and~\eqref{eqlmeta1} that, for $n \! 
\in \! \mathbb{N}$ and $k \! \in \! \lbrace 1,2,\dotsc,K \rbrace$ 
such that $\alpha_{p_{\mathfrak{s}}} \! := \! \alpha_{k} \! \neq \! 
\infty$, for $\lvert \lambda \rvert \! > \! \tilde{\lambda}_{1}$,
\begin{equation} \label{eqlmeta2} 
\dfrac{1}{\mathcal{N}} \ln \pmb{\pi}^{n}_{k}(\lambda) 
\underset{\underset{z_{o}=1+o(1)}{\mathscr{N},n \to \infty}}{=} 
\int_{\mathbb{R}} \tilde{\mathfrak{g}}_{\blacklozenge}(\xi) 
\, \md \mu_{\widetilde{V}}^{f}(\xi) \! + \! \mathcal{O} \left(
\dfrac{\mathfrak{c}_{\blacklozenge}^{1}(n,k,z_{o})}{n^{3/2}} 
\me^{-\frac{3n}{8}} \me^{-\frac{n \tilde{\Gamma}_{r}}{8}} \right), 
\quad r \! \in \! \lbrace \mathrm{A},\mathrm{B} \rbrace.
\end{equation}
One re-writes Equation~\eqref{eqlmrtz32} as follows:
\begin{equation} \label{eqlmeta3} 
\tilde{\mathbb{G}}(x) \! = \! \lambda^{-\frac{\varkappa^{
\infty}_{nk \tilde{k}_{\mathfrak{s}-1}}}{\mathcal{N}}} 
\tilde{\mathbb{G}}_{\blacklozenge}(x),
\end{equation}
where
\begin{align} \label{eqlmeta4} 
\tilde{\mathbb{G}}_{\blacklozenge}(x) :=& \, (\lambda \! - \! 
x)^{\frac{\varkappa_{nk \tilde{k}_{\mathfrak{s}-1}}^{\infty}+1}{
\mathcal{N}}} \left(\dfrac{1}{\lambda \! - \! \alpha_{k}} \right)^{
\frac{1}{\mathcal{N}}} \left(\dfrac{\lambda \! - \! x}{(\lambda \! 
- \! \alpha_{k})(x \! - \! \alpha_{k})} \right)^{\frac{\varkappa_{nk}
-1}{\mathcal{N}}} \prod_{q=1}^{\mathfrak{s}-2} \left(\dfrac{
\lambda \! - \! x}{(\lambda \! - \! \alpha_{p_{q}})(x \! - \! 
\alpha_{p_{q}})} \right)^{\frac{\varkappa_{nk \tilde{k}_{q}}}{
\mathcal{N}}} \nonumber \\
\times& \, \dfrac{(x \! - \! \alpha_{k})^{\frac{\varkappa_{nk}-1}{
\mathcal{N}}} \prod_{q=1}^{\mathfrak{s}-2}(x \! - \! \alpha_{p_{
q}})^{\frac{\varkappa_{nk \tilde{k}_{q}}}{\mathcal{N}}} \prod_{q=
1}^{\mathfrak{s}-2}(\alpha_{k} \! - \! \alpha_{p_{q}})^{\frac{
\varkappa_{nk \tilde{k}_{q}}}{\mathcal{N}}}}{(\alpha_{k} \! - \! x)}.
\end{align}
For $n \! \in \! \mathbb{N}$ and $k \! \in \! \lbrace 1,2,\dotsc,K 
\rbrace$ such that $\alpha_{p_{\mathfrak{s}}} \! := \! \alpha_{k} 
\! \neq \! \infty$, and $\lvert \lambda \rvert \! > \! \tilde{\lambda}_{1}$, 
say, via Equations~\eqref{eqlmeta1}--\eqref{eqlmeta3}, one proceeds 
thus: for $r \! \in \! \lbrace \mathrm{A},\mathrm{B} \rbrace$,
\begin{align*}
\dfrac{1}{\mathcal{N}} \ln \pmb{\pi}^{n}_{k}(\lambda) 
\underset{\underset{z_{o}=1+o(1)}{\mathscr{N},n \to 
\infty}}{=}& \, \int_{\mathbb{R}} \left(\dfrac{\varkappa_{nk 
\tilde{k}_{\mathfrak{s}-1}}^{\infty}}{\mathcal{N}} \ln \lambda \! 
- \! \tilde{F}(\xi) \dfrac{\varkappa_{nk \tilde{k}_{\mathfrak{s}-
1}}^{\infty}}{\mathcal{N}} \ln \lambda \! + \! \tilde{F}(\xi) \ln 
\tilde{\mathbb{G}}_{\blacklozenge}(\xi) \right) \md \mu_{
\widetilde{V}}^{f}(\xi) \! + \! \mathcal{O} \left(\dfrac{
\mathfrak{c}_{\blacklozenge}^{1}(n,k,z_{o})}{n^{3/2}} 
\me^{-\frac{3n}{8}} \me^{-\frac{n \tilde{\Gamma}_{r}}{8}} \right) \\
\underset{\underset{z_{o}=1+o(1)}{\mathscr{N},n \to \infty}}{
=}& \, \int_{\tilde{\mathfrak{D}}^{c}(\tilde{M}_{0}+1)} \left(
\dfrac{\varkappa_{nk \tilde{k}_{\mathfrak{s}-1}}^{\infty}}{
\mathcal{N}} \ln \lambda \! - \! \underbrace{\tilde{F}(\xi)}_{= 
\, 0} \dfrac{\varkappa_{nk \tilde{k}_{\mathfrak{s}-1}}^{\infty}}{
\mathcal{N}} \ln \lambda \! + \! \underbrace{\tilde{F}(\xi)}_{= \, 0} 
\ln \tilde{\mathbb{G}}_{\blacklozenge}(\xi) \right) \underbrace{
\md \mu_{\widetilde{V}}^{f}(\xi)}_{= \, 0} \\
+& \, \int_{\tilde{\mathfrak{D}}(\tilde{M}_{0}+1)} \left(\dfrac{
\varkappa_{nk \tilde{k}_{\mathfrak{s}-1}}^{\infty}}{\mathcal{N}} 
\ln \lambda \! - \! \tilde{F}(\xi) \dfrac{\varkappa_{nk \tilde{k}_{
\mathfrak{s}-1}}^{\infty}}{\mathcal{N}} \ln \lambda \! + \! 
\tilde{F}(\xi) \ln \tilde{\mathbb{G}}_{\blacklozenge}(\xi) \right) 
\md \mu_{\widetilde{V}}^{f}(\xi) \! + \! \mathcal{O} \left(\dfrac{
\mathfrak{c}_{\blacklozenge}^{1}(n,k,z_{o})}{n^{3/2}} \me^{-
\frac{3n}{8}} \me^{-\frac{n \tilde{\Gamma}_{r}}{8}} \right) \\
\underset{\underset{z_{o}=1+o(1)}{\mathscr{N},n \to \infty}}{
=}& \, \int_{\tilde{\mathfrak{D}}(\tilde{M}_{0}+1) \setminus 
\tilde{\mathfrak{D}}(\tilde{M}_{0})} \left(\dfrac{\varkappa_{nk 
\tilde{k}_{\mathfrak{s}-1}}^{\infty}}{\mathcal{N}} \ln \lambda \! - 
\! \underbrace{\tilde{F}(\xi)}_{\in \, [0,1]} \dfrac{\varkappa_{nk 
\tilde{k}_{\mathfrak{s}-1}}^{\infty}}{\mathcal{N}} \ln \lambda \! + 
\! \underbrace{\tilde{F}(\xi)}_{\in \, [0,1]} \ln \tilde{\mathbb{G}}_{
\blacklozenge}(\xi) \right) \underbrace{\md \mu_{\widetilde{
V}}^{f}(\xi)}_{= \, 0} \\
+& \, \int_{\tilde{\mathfrak{D}}(\tilde{M}_{0})} \left(\dfrac{
\varkappa_{nk \tilde{k}_{\mathfrak{s}-1}}^{\infty}}{\mathcal{N}} 
\ln \lambda \! - \! \underbrace{\tilde{F}(\xi)}_{= \, 1} \dfrac{
\varkappa_{nk \tilde{k}_{\mathfrak{s}-1}}^{\infty}}{\mathcal{N}} 
\ln \lambda \! + \! \underbrace{\tilde{F}(\xi)}_{= \, 1} \ln \tilde{
\mathbb{G}}_{\blacklozenge}(\xi) \right) \md \mu_{\widetilde{V}}^{f}
(\xi) \! + \! \mathcal{O} \left(\dfrac{\mathfrak{c}_{\blacklozenge}^{1}
(n,k,z_{o})}{n^{3/2}} \me^{-\frac{3n}{8}} \me^{-\frac{n \tilde{
\Gamma}_{r}}{8}} \right) \\
\underset{\underset{z_{o}=1+o(1)}{\mathscr{N},n \to \infty}}{
=}& \, \int_{\tilde{\mathfrak{D}}(\tilde{M}_{0}) \setminus J_{f}} 
\ln (\tilde{\mathbb{G}}_{\blacklozenge}(\xi)) \, \underbrace{\md 
\mu_{\widetilde{V}}^{f}(\xi)}_{= \, 0} \! + \! \int_{J_{f}} \ln (\tilde{
\mathbb{G}}_{\blacklozenge}(\xi)) \, \md \mu_{\widetilde{V}}^{f}(\xi) 
\! + \! \mathcal{O} \left(\dfrac{\mathfrak{c}_{\blacklozenge}^{1}
(n,k,z_{o})}{n^{3/2}} \me^{-\frac{3n}{8}} \me^{-\frac{n \tilde{
\Gamma}_{r}}{8}} \right) \\
\underset{\underset{z_{o}=1+o(1)}{\mathscr{N},n \to \infty}}{=}& 
\, \int_{J_{f}} \ln (\tilde{\mathbb{G}}_{\blacklozenge}(\xi)) \, \md 
\mu_{\widetilde{V}}^{f}(\xi) \! + \! \mathcal{O} \left(\dfrac{
\mathfrak{c}_{\blacklozenge}^{1}(n,k,z_{o})}{n^{3/2}} 
\me^{-\frac{3n}{8}} \me^{-\frac{n \tilde{\Gamma}_{r}}{8}} \right),
\end{align*}
whence, via Equation~\eqref{eqlmeta4}, for $\lvert \lambda \rvert \! 
> \! \tilde{\lambda}_{1}$,\footnote{Note that, since $\mathscr{M}_{1}
(\mathbb{R}) \! \ni \! \mu_{\widetilde{V}}^{f}$, a straightforward 
calculation shows that
\begin{equation*}
\int_{J_{f}} \ln \left(\dfrac{(\xi \! - \! \alpha_{k})^{\frac{
\varkappa_{nk}-1}{\mathcal{N}}} \prod_{q=1}^{\mathfrak{s}-2}
(\xi \! - \! \alpha_{p_{q}})^{\frac{\varkappa_{nk \tilde{k}_{q}}}{
\mathcal{N}}} \prod_{q^{\prime}=1}^{\mathfrak{s}-2}(\alpha_{k} \! - \! 
\alpha_{p_{q^{\prime}}})^{\frac{\varkappa_{nk \tilde{k}_{q^{\prime}}}}{
\mathcal{N}}}}{\alpha_{k} \! - \! \xi} \right) \md \mu_{\widetilde{V}}^{f}
(\xi) \underset{\underset{z_{o}=1+o(1)}{\mathscr{N},n \to \infty}}{=} 
\mathcal{O}(1).
\end{equation*}}
\begin{align} \label{eqlmeta5} 
\dfrac{1}{\mathcal{N}} \ln \pmb{\pi}^{n}_{k}(\lambda) 
\underset{\underset{z_{o}=1+o(1)}{\mathscr{N},n \to \infty}}{
=}& \, \int_{J_{f}} \ln \left((\lambda \! - \! \xi)^{\frac{\varkappa_{nk 
\tilde{k}_{\mathfrak{s}-1}}^{\infty}+1}{\mathcal{N}}} \left(\dfrac{1}{
\lambda \! - \! \alpha_{k}} \right)^{\frac{1}{\mathcal{N}}} \left(
\dfrac{\lambda \! - \! \xi}{(\lambda \! - \! \alpha_{k})(\xi \! - \! 
\alpha_{k})} \right)^{\frac{\varkappa_{nk}-1}{\mathcal{N}}} \prod_{q
=1}^{\mathfrak{s}-2} \left(\dfrac{\lambda \! - \! \xi}{(\lambda \! 
- \! \alpha_{p_{q}})(\xi \! - \! \alpha_{p_{q}})} \right)^{\frac{
\varkappa_{nk \tilde{k}_{q}}}{\mathcal{N}}} \right) \md \mu_{
\widetilde{V}}^{f}(\xi) \nonumber \\
+& \, \int_{J_{f}} \ln \left(\dfrac{(\xi \! - \! \alpha_{k})^{\frac{
\varkappa_{nk}-1}{\mathcal{N}}} \prod_{q=1}^{\mathfrak{s}-2}
(\xi \! - \! \alpha_{p_{q}})^{\frac{\varkappa_{nk \tilde{k}_{q}}}{
\mathcal{N}}} \prod_{q^{\prime}=1}^{\mathfrak{s}-2}(\alpha_{k} \! - \! 
\alpha_{p_{q^{\prime}}})^{\frac{\varkappa_{nk \tilde{k}_{q^{\prime}}}}{
\mathcal{N}}}}{\alpha_{k} \! - \! \xi} \right) \md \mu_{\widetilde{V}}^{f}
(\xi) \nonumber \\
+& \, \mathcal{O} \left(\dfrac{\mathfrak{c}_{\blacklozenge}^{1}
(n,k,z_{o})}{n^{3/2}} \me^{-\frac{3n}{8}} \me^{-\frac{n \tilde{
\Gamma}_{r}}{8}} \right), \quad r \! \in \! \lbrace \mathrm{A},
\mathrm{B} \rbrace.
\end{align}
For $n \! \in \! \mathbb{N}$ and $k \! \in \! \lbrace 1,2,\dotsc,K 
\rbrace$ such that $\alpha_{p_{\mathfrak{s}}} \! := \! \alpha_{k} 
\! \neq \! \infty$, it follows {}from the corresponding item of 
Lemma~\ref{lem2.2} that
\begin{equation} \label{eqlmeta6} 
\pmb{\pi}^{n}_{k}(\lambda) \left((\lambda \! - \! \alpha_{k})^{\varkappa_{nk}} 
\prod_{q=1}^{\mathfrak{s}-2} \left(\dfrac{\lambda \! - \! \alpha_{p_{q}}}{
\alpha_{k} \! - \! \alpha_{p_{q}}} \right)^{\varkappa_{nk \tilde{k}_{q}}} \right) 
\! = \! \tilde{\mathcal{E}}^{n}_{k} \left(\prod_{m=1}^{\mathcal{N}} \dfrac{
\lambda \! - \! x_{m}}{\alpha_{k} \! - \! x_{m}} \right).
\end{equation}
Recall {}from the proof of case~\pmb{(1)} of Lemma~\ref{lemrootz} 
that, for $n \! \in \! \mathbb{N}$ and $k \! \in \! \lbrace 
1,2,\dotsc,K \rbrace$ such that $\alpha_{p_{\mathfrak{s}}} \! 
:= \! \alpha_{k} \! \neq \! \infty$, in the double-scaling limit 
$\mathscr{N},n \! \to \! \infty$ such that $z_{o} \! = \! 1 \! + \! 
o(1)$, $\tilde{\mathfrak{Z}}_{0} \subseteq \tilde{\mathfrak{D}}
(\tilde{\lambda}_{1})$, where $\tilde{\mathfrak{Z}}_{0}$ is defined 
in the lemma, and $\tilde{\mathfrak{D}}(\tilde{\lambda}_{1}) \! := \! 
[-\tilde{\lambda}_{1},\tilde{\lambda}_{1}] \setminus \cup_{\underset{q 
\neq \mathfrak{s}-1}{q=1}}^{\mathfrak{s}} \mathscr{O}_{\frac{1}{\tilde{
\lambda}_{1}}}(\alpha_{p_{q}})$: as $\tilde{\mathfrak{Z}}_{0} \cap \tilde{
\mathfrak{D}}^{c}(\tilde{\lambda}_{1}) \! = \! \varnothing$, where 
$\tilde{\mathfrak{D}}^{c}(\tilde{\lambda}_{1}) \! := \! \mathbb{R} 
\setminus \tilde{\mathfrak{D}}(\tilde{\lambda}_{1})$, it follows that 
(cf. Equation~\eqref{eqlmeta6})
\begin{equation*}
\left. (\lambda \! - \! \alpha_{k})^{\varkappa_{nk}} \prod_{q=
1}^{\mathfrak{s}-2} \left(\dfrac{\lambda \! - \! \alpha_{p_{q}}}{
\alpha_{k} \! - \! \alpha_{p_{q}}} \right)^{\varkappa_{nk \tilde{
k}_{q}}} \right\vert_{\lambda \in \tilde{\mathfrak{Z}}_{0}} \! = \! 
(\tilde{\mathfrak{z}}^{n}_{k}(j) \! - \! \alpha_{k})^{\varkappa_{nk}} 
\prod_{q=1}^{\mathfrak{s}-2} \left(\dfrac{\tilde{\mathfrak{z}}^{n}_{k}
(j) \! - \! \alpha_{p_{q}}}{\alpha_{k} \! - \! \alpha_{p_{q}}} \right)^{
\varkappa_{nk \tilde{k}_{q}}} \! \neq \! 0, \quad j \! = \! 1,2,\dotsc,
\mathcal{N};
\end{equation*}
hence, since $\pmb{\pi}^{n}_{k}(\tilde{\mathfrak{z}}^{n}_{k}(j)) \! 
= \! 0$, $j \! = \! 1,2,\dotsc,\mathcal{N}$, and (via the linearity 
of $\tilde{\mathcal{E}}^{n}_{k}$ and the fact that $\tilde{\mathcal{
E}}^{n}_{k}(1) \! = \! 1)$ $\tilde{\mathcal{E}}^{n}_{k}(\prod_{m=
1}^{\mathcal{N}} \tfrac{\lambda -x_{m}}{\alpha_{k}-x_{m}})$ is 
a non-monic polynomial in $\lambda$ with non-zero leading 
coefficient $\tilde{\mathcal{E}}^{n}_{k}(\prod_{m=1}^{\mathcal{N}}
(\alpha_{k} \! - \! x_{m})^{-1}) \! := \! \tilde{\mathfrak{c}}_{
\mathcal{N}}(n,k,z_{o}) \! = \! \tilde{\mathfrak{c}}_{\mathcal{N}}$ (see 
Equation~\eqref{eqlmeta11} below), that is, $\operatorname{coeff}
(\lambda^{\mathcal{N}}) \! = \! \tilde{\mathfrak{c}}_{\mathcal{N}}$, 
it follows that
\begin{equation*}
\left. \tilde{\mathcal{E}}^{n}_{k} \left(\prod_{m=1}^{\mathcal{N}} 
\dfrac{\lambda \! - \! x_{m}}{\alpha_{k} \! - \! x_{m}} \right) 
\right\vert_{\lambda \in \tilde{\mathfrak{Z}}_{0}} \! = \! \tilde{
\mathcal{E}}^{n}_{k} \left(\prod_{m=1}^{\mathcal{N}} \dfrac{
\tilde{\mathfrak{z}}^{n}_{k}(j) \! - \! x_{m}}{\alpha_{k} \! - \! x_{m}} 
\right) \! = \! 0, \quad j \! = \! 1,2,\dotsc,\mathcal{N},
\end{equation*}
which implies that $\tilde{\mathcal{E}}^{n}_{k}(\prod_{m=1}^{
\mathcal{N}} \tfrac{\lambda -x_{m}}{\alpha_{k}-x_{m}})$ has 
the representation
\begin{equation*}
\tilde{\mathcal{E}}^{n}_{k} \left(\prod_{m=1}^{\mathcal{N}} 
\dfrac{\lambda \! - \! \xi_{m}}{\alpha_{k} \! - \! \xi_{m}} 
\right) \! = \! \tilde{\mathfrak{c}}_{\mathcal{N}} \prod_{j=
1}^{\mathcal{N}}(\lambda \! - \! \tilde{\mathfrak{z}}^{n}_{k}(j)),
\end{equation*}
which, via Equation~\eqref{eqlmeta6}, leads to
\begin{equation} \label{eqlmeta7} 
\pmb{\pi}^{n}_{k}(\lambda) \left((\lambda \! - \! \alpha_{k})^{\varkappa_{nk}} 
\prod_{q=1}^{\mathfrak{s}-2} \left(\dfrac{\lambda \! - \! \alpha_{p_{q}}}{
\alpha_{k} \! - \! \alpha_{p_{q}}} \right)^{\varkappa_{nk \tilde{k}_{q}}} 
\right) \! = \! \tilde{\mathfrak{c}}_{\mathcal{N}} \prod_{j=1}^{\mathcal{N}}
(\lambda \! - \! \tilde{\mathfrak{z}}^{n}_{k}(j)).
\end{equation}
For $n \! \in \! \mathbb{N}$ and $k \! \in \! \lbrace 1,2,\dotsc,K 
\rbrace$ such that $\alpha_{p_{\mathfrak{s}}} \! := \! \alpha_{k} 
\! \neq \! \infty$, define the associated normalised zero counting 
measure as follows:
\begin{equation} \label{eqlmeta8} 
\tilde{\eta}_{\tilde{\mathfrak{z}}}(x) \! := \! \dfrac{1}{\mathcal{N}} 
\sum_{j=1}^{\mathcal{N}} \delta_{\tilde{\mathfrak{z}}^{n}_{k}(j)}(x),
\end{equation}
where $\delta_{\tilde{\mathfrak{z}}^{n}_{k}(j)}(x)$ $(= \! \delta (x \! - \! 
\tilde{\mathfrak{z}}^{n}_{k}(j)))$ is the Dirac delta (atomic) mass concentrated 
at $\tilde{\mathfrak{z}}^{n}_{k}(j)$, $j \! = \! 1,2,\dotsc,\mathcal{N}$ (note 
that $\int_{\mathbb{R}} \md \tilde{\eta}_{\tilde{\mathfrak{z}}}(\xi) \! = \! 
\int_{\mathbb{R}} \tfrac{1}{\mathcal{N}} \sum_{j=1}^{\mathcal{N}} \delta 
(\xi \! - \! \tilde{\mathfrak{z}}^{n}_{k}(j)) \, \md \xi \! = \! 1)$. Via 
Equations~\eqref{eqlmeta7} and~\eqref{eqlmeta8}, and the 
decomposition (cf. Equation~\eqref{fincount}) $\sum_{q=1}^{\mathfrak{s}-2} 
\varkappa_{nk \tilde{k}_{q}} \! + \! \varkappa_{nk \tilde{k}_{\mathfrak{s}-1}}^{\infty} 
\! + \varkappa_{nk} \! = \! (n \! - \! 1)K \! + \! k \! = \! \mathcal{N}$, a calculation 
reveals that
\begin{align} \label{eqlmeta9} 
\dfrac{1}{\mathcal{N}} \ln \pmb{\pi}^{n}_{k}(\lambda) =& \, 
\int_{\mathbb{R}} \ln \left((\lambda \! - \! \xi)^{\frac{
\varkappa_{nk \tilde{k}_{\mathfrak{s}-1}}^{\infty}+1}{\mathcal{N}}} 
\left(\dfrac{1}{\lambda \! - \! \alpha_{k}} \right)^{\frac{1}{\mathcal{
N}}} \left(\dfrac{\lambda \! - \! \xi}{(\lambda \! - \! \alpha_{k})(\xi 
\! - \! \alpha_{k})} \right)^{\frac{\varkappa_{nk}-1}{\mathcal{N}}} 
\prod_{q=1}^{\mathfrak{s}-2} \left(\dfrac{\lambda \! - \! \xi}{(
\lambda \! - \! \alpha_{p_{q}})(\xi \! - \! \alpha_{p_{q}})} \right)^{
\frac{\varkappa_{nk \tilde{k}_{q}}}{\mathcal{N}}} \right) \md 
\tilde{\eta}_{\tilde{\mathfrak{z}}}(\xi) \nonumber \\
+& \, \int_{\mathbb{R}} \ln \left((\xi \! - \! \alpha_{k})^{\frac{
\varkappa_{nk}-1}{\mathcal{N}}} \prod_{q=1}^{\mathfrak{s}-2}
(\xi \! - \! \alpha_{p_{q}})^{\frac{\varkappa_{nk \tilde{k}_{q}}}{
\mathcal{N}}} \prod_{q^{\prime}=1}^{\mathfrak{s}-2}(\alpha_{k} \! - \! 
\alpha_{p_{q^{\prime}}})^{\frac{\varkappa_{nk \tilde{k}_{q^{\prime}}}}{
\mathcal{N}}} \right) \md \tilde{\eta}_{\tilde{\mathfrak{z}}}(\xi) \! 
+ \! \dfrac{1}{\mathcal{N}} \ln \tilde{\mathfrak{c}}_{\mathcal{N}}.
\end{align}
Now, for $n \! \in \! \mathbb{N}$ and $k \! \in \! \lbrace 1,2,
\dotsc,K \rbrace$ such that $\alpha_{p_{\mathfrak{s}}} \! := 
\! \alpha_{k} \! \neq \! \infty$, in the double-scaling limit 
$\mathscr{N},n \! \to \! \infty$ such that $z_{o} \! = \! 1 \! + \! 
o(1)$, for $\lambda$ in the `cut plane' $\mathbb{C} \setminus 
(-\infty,\tilde{\lambda}_{1}]$, because $\int_{\mathbb{R}} \ln 
(\lambda \! - \! \xi) \, \md \tilde{\eta}_{\tilde{\mathfrak{z}}}
(\xi)$ and $\ln (\lambda \! - \! \alpha_{p_{q}})$, $q \! = \! 
1,\dotsc,\mathfrak{s} \! - \! 2,\mathfrak{s}$, are analytic 
functions,\footnote{Note that, since $\tilde{\mathfrak{D}}(\tilde{
\lambda}_{1}) \! \supseteq \! \tilde{\mathfrak{Z}}_{0}$ and $\tilde{
\mathfrak{D}}(\tilde{\lambda}_{1}) \cap \lbrace \alpha_{p_{1}},\dotsc,
\alpha_{p_{\mathfrak{s}-2}},\alpha_{p_{\mathfrak{s}}} \rbrace \! = \! 
\varnothing$, $\int_{\mathbb{R}} \sum_{q=1}^{\mathfrak{s}-2} 
\tfrac{\varkappa_{nk \tilde{k}_{q}}}{\mathcal{N}} \ln (\xi \! - \! 
\alpha_{p_{q}}) \, \md \tilde{\eta}_{\tilde{\mathfrak{z}}}(\xi) \! = \!  
\sum_{j=1}^{\mathcal{N}} \sum_{q=1}^{\mathfrak{s}-2} \tfrac{
\varkappa_{nk \tilde{k}_{q}}}{\mathcal{N}} \tfrac{\ln (\tilde{\mathfrak{
z}}^{n}_{k}(j)-\alpha_{p_{q}})}{\mathcal{N}} \! =_{\underset{z_{o}=1+
o(1)}{\mathscr{N},n \to \infty}} \! \mathcal{O}(1)$, and $\int_{\mathbb{R}}
(\tfrac{\varkappa_{nk}-1}{\mathcal{N}}) \ln (\xi \! - \! \alpha_{k}) 
\, \md \tilde{\eta}_{\tilde{\mathfrak{z}}}(\xi) \! = \! 
(\tfrac{\varkappa_{nk}-1}{\mathcal{N}}) \sum_{j=1}^{\mathcal{N}} 
\tfrac{\ln (\tilde{\mathfrak{z}}^{n}_{k}(j)-\alpha_{k})}{\mathcal{N}} 
\! =_{\underset{z_{o}=1+o(1)}{\mathscr{N},n \to \infty}} \! 
\mathcal{O}(1)$.} it follows (cf. Equation~\eqref{eqlmeta9}) that 
$\mathcal{N}^{-1} \ln \pmb{\pi}^{n}_{k}(\lambda)$, too, is analytic 
for $\mathbb{C} \setminus (-\infty,\tilde{\lambda}_{1}] \! \ni \! 
\lambda$; furthermore, since $\tilde{\mathfrak{D}}(\tilde{\lambda}_{1}) 
\! \supseteq \! \tilde{\mathfrak{Z}}_{0}$ and $\tilde{\mathfrak{D}}
(\tilde{\lambda}_{1}) \cap \lbrace \alpha_{p_{1}},\dotsc,
\alpha_{p_{\mathfrak{s}-2}},\alpha_{p_{\mathfrak{s}}} \rbrace 
\! = \! \varnothing$,
\begin{equation*}
\dfrac{1}{\mathcal{N}} \ln \pmb{\pi}^{n}_{k}(\lambda) \! = \! 
\dfrac{1}{\mathcal{N}} \sum_{j=1}^{\mathcal{N}} \ln (\lambda \! - \! 
\tilde{\mathfrak{z}}^{n}_{k}(j)) \! - \! \dfrac{\varkappa_{nk}}{\mathcal{
N}} \ln (\lambda \! - \! \alpha_{k}) \! - \! \sum_{q=1}^{\mathfrak{s}
-2} \dfrac{\varkappa_{nk \tilde{k}_{q}}}{\mathcal{N}} \ln \left(\dfrac{
\lambda \! - \! \alpha_{p_{q}}}{\alpha_{k} \! - \! \alpha_{p_{q}}} 
\right) \! + \! \dfrac{1}{\mathcal{N}} \ln \tilde{\mathfrak{c}}_{
\mathcal{N}}
\end{equation*}
is bounded (locally) on compact subsets of the cut plane $\mathbb{C} 
\setminus (-\infty,\tilde{\lambda}_{1}]$ $(\ni \! \lambda)$; as $\int_{
J_{f}} \ln (\lambda \! - \! \xi) \, \md \mu_{\widetilde{V}}^{f}(\xi)$ 
and $\ln (\lambda \! - \! \alpha_{p_{q}})$, $q \! = \! 1,\dotsc,
\mathfrak{s} \! - \! 2,\mathfrak{s}$, too, are analytic functions 
for $\mathbb{C} \setminus (-\infty,\tilde{\lambda}_{1}] \! \ni \! 
\lambda$, it follows by Vitali's Theorem {}\footnote{See \cite{rremt}, 
p.~157: \emph{Let $G_{V}$ be a domain, and let $f_{0},f_{1},f_{2},
\dotsc \! \in \! \mathrm{Hol}(G_{V})$ be a sequence of functions that is 
locally bounded in $G_{V}$. Then the following statements are equivalent: 
{\rm (i)} the sequence $f_{n}$ is compactly convergent in $G_{V}$$;$ 
{\rm (ii)} there exists a point $c_{V} \! \in \! G_{V}$ such that for every 
$m \! \in \! \mathbb{N}$ the sequence of numbers $f_{0}^{(m)}(c_{V}),
f_{1}^{(m)}(c_{V}),f_{2}^{(m)}(c_{V}),\dotsc$ converges$;$ and {\rm (iii)} 
the set $A_{V} \! := \! \lbrace \mathstrut w \! \in \! G_{V}; \, 
\lim_{n \to \infty}f_{n}(w) \, \text{exists in $G_{V}$} \rbrace$ has an 
accumulation point in $G_{V}$.}} and Equation~\eqref{eqlmeta5} that, 
for $n \! \in \! \mathbb{N}$ and $k \! \in \! \lbrace 1,2,\dotsc,K \rbrace$ 
such that $\alpha_{p_{\mathfrak{s}}} \! := \! \alpha_{k} \! \neq \! \infty$, 
and $z \! \in \! \mathbb{C} \setminus (-\infty,\tilde{\lambda}_{1}]$,
\begin{align} \label{eqlmeta10} 
\lim_{\underset{z_{o}=1+o(1)}{\mathscr{N},n \to \infty}} \dfrac{1}{
\mathcal{N}} \ln \pmb{\pi}^{n}_{k}(z) \underset{\underset{z_{o}=1+
o(1)}{\mathscr{N},n \to \infty}}{=}& \, \int_{J_{f}} \ln \left((z \! - \! 
\xi)^{\frac{\varkappa_{nk \tilde{k}_{\mathfrak{s}-1}}^{\infty}+1}{
\mathcal{N}}} \left(\dfrac{1}{z \! - \! \alpha_{k}} \right)^{\frac{1}{
\mathcal{N}}} \left(\dfrac{z \! - \! \xi}{(z \! - \! \alpha_{k})(\xi \! - \! 
\alpha_{k})} \right)^{\frac{\varkappa_{nk}-1}{\mathcal{N}}} \prod_{q=
1}^{\mathfrak{s}-2} \left(\dfrac{z \! - \! \xi}{(z \! - \! \alpha_{p_{q}})
(\xi \! - \! \alpha_{p_{q}})} \right)^{\frac{\varkappa_{nk \tilde{k}_{q}}}{
\mathcal{N}}} \right) \md \mu_{\widetilde{V}}^{f}(\xi) \nonumber \\
+& \, \int_{J_{f}} \ln \left(\dfrac{(\xi \! - \! \alpha_{k})^{\frac{
\varkappa_{nk}-1}{\mathcal{N}}} \prod_{q=1}^{\mathfrak{s}-2}
(\xi \! - \! \alpha_{p_{q}})^{\frac{\varkappa_{nk \tilde{k}_{q}}}{
\mathcal{N}}} \prod_{q^{\prime}=1}^{\mathfrak{s}-2}(\alpha_{k} \! - \! 
\alpha_{p_{q^{\prime}}})^{\frac{\varkappa_{nk \tilde{k}_{q^{\prime}}}}{
\mathcal{N}}}}{\alpha_{k} \! - \! \xi} \right) \md \mu_{\widetilde{V}}^{f}
(\xi),
\end{align}
where the convergence, in the double-scaling limit $\mathscr{N},n 
\! \to \! \infty$ such that $z_{o} \! = \! 1 \! + \! o(1)$, is uniform on 
compact subsets of $\mathbb{C} \setminus (-\infty,\tilde{\lambda}_{1}]$. 
It follows {}from Equation~\eqref{eqlmrtz8} that
\begin{equation} \label{eqlmeta11} 
\dfrac{1}{\mathcal{N}} \ln \tilde{\mathfrak{c}}_{\mathcal{N}} \! 
= \! \dfrac{1}{\mathcal{N}} \ln \tilde{\mathcal{E}}^{n}_{k} \left(
\prod_{m=1}^{\mathcal{N}}(\alpha_{k} \! - \! \xi_{m})^{-1} \right) 
\underset{\underset{z_{o}=1+o(1)}{\mathscr{N},n \to \infty}}{=} 
\int_{\mathbb{R}} \ln (\alpha_{k} \! - \! \xi)^{-1} \, \md 
\mu_{\widetilde{V}}^{f}(\xi):
\end{equation}
via Equations~\eqref{eqlmeta9} and~\eqref{eqlmeta11}, a corollary 
of Equation~\eqref{eqlmeta10} is that, for $n \! \in \! \mathbb{N}$ 
and $k \! \in \! \lbrace 1,2,\dotsc,K \rbrace$ such that $\alpha_{p_{
\mathfrak{s}}} \! := \! \alpha_{k} \! \neq \! \infty$, for $\mathbb{C} 
\setminus (-\infty,\tilde{\lambda}_{1}] \! \ni \! z$,
\begin{align} \label{eqlmeta12} 
& \, \lim_{\underset{z_{o}=1+o(1)}{\mathscr{N},n \to \infty}} \left(
\int_{\mathbb{R}} \ln \left((z \! - \! \xi)^{\frac{\varkappa_{nk 
\tilde{k}_{\mathfrak{s}-1}}^{\infty}+1}{\mathcal{N}}} \left(
\dfrac{1}{z \! - \! \alpha_{k}} \right)^{\frac{1}{\mathcal{N}}} 
\left(\dfrac{z \! - \! \xi}{(z \! - \! \alpha_{k})(\xi \! - \! \alpha_{k})} 
\right)^{\frac{\varkappa_{nk}-1}{\mathcal{N}}} \prod_{q=1}^{
\mathfrak{s}-2} \left(\dfrac{z \! - \! \xi}{(z \! - \! \alpha_{p_{q}})
(\xi \! - \! \alpha_{p_{q}})} \right)^{\frac{\varkappa_{nk \tilde{k}_{
q}}}{\mathcal{N}}} \right) \md \tilde{\eta}_{\tilde{\mathfrak{z}}}(\xi) 
\right. \nonumber \\
&\left. \, + \int_{\mathbb{R}} \ln \left((\xi \! - \! \alpha_{k})^{
\frac{\varkappa_{nk}-1}{\mathcal{N}}} \prod_{q=1}^{\mathfrak{s}-2}
(\xi \! - \! \alpha_{p_{q}})^{\frac{\varkappa_{nk \tilde{k}_{q}}}{\mathcal{N}}} 
\prod_{q^{\prime}=1}^{\mathfrak{s}-2}(\alpha_{k} \! - \! \alpha_{p_{
q^{\prime}}})^{\frac{\varkappa_{nk \tilde{k}_{q^{\prime}}}}{\mathcal{N}}} 
\right) \md \tilde{\eta}_{\tilde{\mathfrak{z}}}(\xi) \right) \nonumber \\
&\underset{\underset{z_{o}=1+o(1)}{\mathscr{N},n \to \infty}}{=} 
\int_{J_{f}} \ln \left((z \! - \! \xi)^{\frac{\varkappa_{nk \tilde{k}_{
\mathfrak{s}-1}}^{\infty}+1}{\mathcal{N}}} \left(\dfrac{1}{z \! - \! 
\alpha_{k}} \right)^{\frac{1}{\mathcal{N}}} \left(\dfrac{z \! - \! 
\xi}{(z \! - \! \alpha_{k})(\xi \! - \! \alpha_{k})} \right)^{\frac{
\varkappa_{nk}-1}{\mathcal{N}}} \prod_{q=1}^{\mathfrak{s}-2} 
\left(\dfrac{z \! - \! \xi}{(z \! - \! \alpha_{p_{q}})(\xi \! - \! 
\alpha_{p_{q}})} \right)^{\frac{\varkappa_{nk \tilde{k}_{q}}}{
\mathcal{N}}} \right) \md \mu_{\widetilde{V}}^{f}(\xi) \nonumber \\
& \, + \int_{J_{f}} \ln \left((\xi \! - \! \alpha_{k})^{\frac{
\varkappa_{nk}-1}{\mathcal{N}}} \prod_{q=1}^{\mathfrak{s}
-2}(\xi \! - \! \alpha_{p_{q}})^{\frac{\varkappa_{nk \tilde{k}_{q}}}{
\mathcal{N}}} \prod_{q^{\prime}=1}^{\mathfrak{s}-2}(\alpha_{k} 
\! - \! \alpha_{p_{q^{\prime}}})^{\frac{\varkappa_{nk \tilde{k}_{
q^{\prime}}}}{\mathcal{N}}} \right) \md \mu_{\widetilde{V}}^{f}(\xi).
\end{align}
Since, for $n \! \in \! \mathbb{N}$ and $k \! \in \! \lbrace 1,2,
\dotsc,K \rbrace$ such that $\alpha_{p_{\mathfrak{s}}} \! := 
\! \alpha_{k} \! \neq \! \infty$, the convergence stated in 
Equation~\eqref{eqlmeta12} is (normal) uniform on compact subsets 
of $\mathbb{C} \setminus (-\infty,\tilde{\lambda}_{1}]$ $(\ni \! z)$, 
it follows by the equivalent Part~(ii) of Vitali's Theorem 
{}\footnote{In---the equivalent---Part~(ii) of Vitali's Theorem, take any 
`fixed' $z \! \in \! \mathbb{C} \setminus (-\infty,\tilde{\lambda}_{1}]$ 
that is real and greater than $\tilde{\lambda}_{1}$ to play the r\^{o}le 
of $c_{V}$.} and a successive differentiation (with respect to $z$) 
argument, that, for arbitrary $m \! \in \! \mathbb{N}$ and $z \! 
\in \! \mathbb{C} \setminus (-\infty,\tilde{\lambda}_{1}]$,
\begin{equation*}
\lim_{\underset{z_{o}=1+o(1)}{\mathscr{N},n \to \infty}} 
\int_{\mathbb{R}}(z \! - \! \xi)^{-m} \, \md \tilde{\eta}_{\tilde{\mathfrak{z}}}
(\xi) \underset{\underset{z_{o}=1+o(1)}{\mathscr{N},n \to \infty}}{=} 
\int_{J_{f}}(z \! - \! \xi)^{-m} \, \md \mu_{\widetilde{V}}^{f}(\xi):
\end{equation*} 
as $\md \tilde{\eta}_{\tilde{\mathfrak{z}}}$ and $\md \mu_{\widetilde{V}}^{f}$ 
are, for $n \! \in \! \mathbb{N}$ and $k \! \in \! \lbrace 1,2,\dotsc,K \rbrace$ 
such that $\alpha_{p_{\mathfrak{s}}} \! := \! \alpha_{k} \! \neq \! \infty$, 
supported on bounded (and measurable) compact sets, it follows via the 
complex form of the Stone-Weierstrass theorem {}\footnote{See, in particular, 
\cite{ehtksg}, Chapter~2, Section~6, Theorem~7.3.4.} and the latter 
relation that $\tilde{\eta}_{\tilde{\mathfrak{z}}} \! \overset{\ast}{\to} \! 
\mu_{\widetilde{V}}^{f}$ in the double-scaling limit $\mathscr{N},n \! \to \! 
\infty$ such that $z_{o} \! = \! 1 \! + \! o(1)$.\footnote{As a corollary, note 
that, for $n \! \in \! \mathbb{N}$ and $k \! \in \! \lbrace 1,2,\dotsc,K \rbrace$ 
such that $\alpha_{p_{\mathfrak{s}}} \! := \! \alpha_{k} \! \neq \! \infty$, $x \! 
\in \! J_{f} \! \Leftrightarrow \! \operatorname{dist}(x,\tilde{\mathfrak{Z}}_{0}) 
\! \to \! 0$ in the double-scaling limit $\mathscr{N},n \to \infty$ such that 
$z_{o} \! = \! 1+o(1)$.}

\pmb{(2)} The proof of this case, that is, $n \! \in \! \mathbb{N}$ 
and $k \! \in \! \lbrace 1,2,\dotsc,K \rbrace$ such that $\alpha_{
p_{\mathfrak{s}}} \! := \! \alpha_{k} \! = \! \infty$, is virtually identical 
to the proof of \pmb{(1)} above; one mimics, \emph{verbatim}, the 
scheme of the calculations presented in \pmb{(1)} in order to 
arrive at the corresponding claim stated in the lemma; in order to 
do so, however, the analogues of Equations~\eqref{eqlmeta5}, 
and \eqref{eqlmeta7}--\eqref{eqlmeta9}, respectively, are necessary, 
which, in the present case, read: for $\lvert \lambda \rvert \! > \! 
\hat{\lambda}_{1} \! := \! \max \lbrace \hat{\lambda}_{1}^{\ast},
2(\hat{M}_{0} \! + \! 1) \! + \! 3(\mathfrak{s} \! - \! 1)(1 \! + \! 
\max_{q=1,2,\dotsc,\mathfrak{s}-1} \lbrace \lvert \alpha_{p_{q}} 
\rvert \rbrace)(\min_{q=1,2,\dotsc,\mathfrak{s}-1} \lbrace \hat{\Delta}
(q) \rbrace)^{-1} \rbrace$ $(\gg \! 1)$,
\begin{align} \label{eqlmeta13} 
\dfrac{1}{\mathcal{N}} \ln \pmb{\pi}^{n}_{k}(\lambda) 
\underset{\underset{z_{o}=1+o(1)}{\mathscr{N},n \to \infty}}{
=}& \, \int_{J_{\infty}} \ln \left((\lambda \! - \! \xi)^{\frac{
\varkappa_{nk}}{\mathcal{N}}} \prod_{q=1}^{\mathfrak{s}-1} 
\left(\dfrac{\lambda \! - \! \xi}{(\lambda \! - \! \alpha_{p_{q}})
(\xi \! - \! \alpha_{p_{q}})} \right)^{\frac{\varkappa_{nk 
\tilde{k}_{q}}}{\mathcal{N}}} \right) \md \mu_{\widetilde{V}}^{
\infty}(\xi) \! + \! \int_{J_{\infty}} \ln \left(\prod_{q=1}^{\mathfrak{s}
-1}(\xi \! - \! \alpha_{p_{q}})^{\frac{\varkappa_{nk \tilde{k}_{q}}}{
\mathcal{N}}} \right) \md \mu_{\widetilde{V}}^{\infty}(\xi) 
\nonumber \\
+& \, \mathcal{O} \left(\dfrac{\mathfrak{c}_{\lozenge}^{1}(n,k,z_{o})}{
n^{3/2}} \me^{-\frac{3n}{8}} \me^{-\frac{n \hat{\Gamma}_{0}}{8}} 
\right),
\end{align}
\begin{equation} \label{eqlmeta14} 
\pmb{\pi}^{n}_{k}(\lambda) \prod_{q=1}^{\mathfrak{s}-1}
(\lambda \! - \! \alpha_{p_{q}})^{\varkappa_{nk \tilde{k}_{q}}} 
\! = \! \prod_{j=1}^{\mathcal{N}}(\lambda \! - \! 
\hat{\mathfrak{z}}^{n}_{k}(j)),
\end{equation}
with the associated normalised zero counting measure defined 
as
\begin{equation} \label{eqlmeta15} 
\hat{\eta}_{\hat{\mathfrak{z}}}(x) \! := \! \dfrac{1}{\mathcal{N}} 
\sum_{j=1}^{\mathcal{N}} \delta_{\hat{\mathfrak{z}}^{n}_{k}(j)}(x),
\end{equation}
and
\begin{equation} \label{eqlmeta16} 
\dfrac{1}{\mathcal{N}} \ln \pmb{\pi}^{n}_{k}(\lambda) \! = 
\!  \int_{\mathbb{R}} \ln \left((\lambda \! - \! \xi)^{\frac{
\varkappa_{nk}}{\mathcal{N}}} \prod_{q=1}^{\mathfrak{s}-1} \left(
\dfrac{\lambda \! - \! \xi}{(\lambda \! - \! \alpha_{p_{q}})(\xi \! 
- \! \alpha_{p_{q}})} \right)^{\frac{\varkappa_{nk \tilde{k}_{q}}}{
\mathcal{N}}} \right) \md \hat{\eta}_{\hat{\mathfrak{z}}}(\xi) \! 
+ \! \int_{\mathbb{R}} \ln \left(\prod_{q=1}^{\mathfrak{s}-1}
(\xi \! - \! \alpha_{p_{q}})^{\frac{\varkappa_{nk \tilde{k}_{q}}}{
\mathcal{N}}} \right) \md \hat{\eta}_{\hat{\mathfrak{z}}}(\xi).
\end{equation}
Now, for $n \! \in \! \mathbb{N}$ and $k \! \in \! \lbrace 1,2,
\dotsc,K \rbrace$ such that $\alpha_{p_{\mathfrak{s}}} \! := \! 
\alpha_{k} \! = \! \infty$, in the double-scaling limit $\mathscr{N},
n \! \to \! \infty$ such that $z_{o} \! = \! 1 \! + \! o(1)$, for 
$\lambda$ in the `cut plane' $\mathbb{C} \setminus (-\infty,
\hat{\lambda}_{1}]$, because $\int_{\mathbb{R}} \ln (\lambda 
\! - \! \xi) \, \md \hat{\eta}_{\hat{\mathfrak{z}}}(\xi)$ and 
$\ln (\lambda \! - \! \alpha_{p_{q}})$, $q \! = \! 1,2,\dotsc,
\mathfrak{s} \! - \! 1$, are analytic functions,\footnote{Note that, 
since $\hat{\mathfrak{D}}(\hat{\lambda}_{1}) \! \supseteq \! \hat{
\mathfrak{Z}}_{0}$ and $\hat{\mathfrak{D}}(\hat{\lambda}_{1}) \cap 
\lbrace \alpha_{p_{1}},\alpha_{p_{2}},\dotsc,\alpha_{p_{\mathfrak{s}
-1}} \rbrace \! = \! \varnothing$, $\int_{\mathbb{R}} 
\sum_{q=1}^{\mathfrak{s}-1} \tfrac{\varkappa_{nk \tilde{k}_{q}}}{
\mathcal{N}} \ln (\xi \! - \! \alpha_{p_{q}}) \, \md \hat{\eta}_{
\hat{\mathfrak{z}}}(\xi) \! = \!  \sum_{j=1}^{\mathcal{N}} \sum_{q
=1}^{\mathfrak{s}-1} \tfrac{\varkappa_{nk \tilde{k}_{q}}}{\mathcal{
N}} \tfrac{\ln (\hat{\mathfrak{z}}^{n}_{k}(j)-\alpha_{p_{q}})}{
\mathcal{N}} \! =_{\underset{z_{o}=1+o(1)}{\mathscr{N},n 
\to \infty}} \! \mathcal{O}(1)$.} it follows (cf. 
Equation~\eqref{eqlmeta16}) that $\mathcal{N}^{-1} \ln 
\pmb{\pi}^{n}_{k}(\lambda)$, too, is analytic for $\mathbb{C} 
\setminus (-\infty,\hat{\lambda}_{1}] \! \ni \! \lambda$; furthermore, 
since $\hat{\mathfrak{D}}(\hat{\lambda}_{1}) \! \supseteq \! 
\hat{\mathfrak{Z}}_{0}$ and $\hat{\mathfrak{D}}(\hat{\lambda}_{1}) 
\cap \lbrace \alpha_{p_{1}},\alpha_{p_{2}},\dotsc,\alpha_{p_{\mathfrak{s}-1}} 
\rbrace \! = \! \varnothing$,
\begin{equation*}
\dfrac{1}{\mathcal{N}} \ln \pmb{\pi}^{n}_{k}(\lambda) \! = \! 
\dfrac{1}{\mathcal{N}} \sum_{j=1}^{\mathcal{N}} \ln (\lambda \! - \! 
\hat{\mathfrak{z}}^{n}_{k}(j)) \! - \! \sum_{q=1}^{\mathfrak{s}-1} 
\dfrac{\varkappa_{nk \tilde{k}_{q}}}{\mathcal{N}} \ln (\lambda \! - \! 
\alpha_{p_{q}})
\end{equation*}
is bounded (locally) on compact subsets of the cut plane $\mathbb{C} 
\setminus (-\infty,\hat{\lambda}_{1}]$ $(\ni \! \lambda)$; as $\int_{
J_{\infty}} \ln (\lambda \! - \! \xi) \, \md \mu_{\widetilde{V}}^{\infty}
(\xi)$ and $\ln (\lambda \! - \! \alpha_{p_{q}})$, $q \! = \! 1,2,\dotsc,
\mathfrak{s} \! - \! 1$, too, are analytic functions for $\mathbb{C} 
\setminus (-\infty,\hat{\lambda}_{1}] \! \ni \! \lambda$, it follows by 
Vitali's Theorem and Equation~\eqref{eqlmeta13} that, for $n \! \in \! 
\mathbb{N}$ and $k \! \in \! \lbrace 1,2,\dotsc,K \rbrace$ such that 
$\alpha_{p_{\mathfrak{s}}} \! := \! \alpha_{k} \! = \! \infty$, and 
$z \! \in \! \mathbb{C} \setminus (-\infty,\hat{\lambda}_{1}]$,
\begin{equation} \label{eqlmeta17} 
\lim_{\underset{z_{o}=1+o(1)}{\mathscr{N},n \to \infty}} 
\dfrac{1}{\mathcal{N}} \ln \pmb{\pi}^{n}_{k}(z) \underset{
\underset{z_{o}=1+o(1)}{\mathscr{N},n \to \infty}}{=} \int_{J_{\infty}} 
\ln \left((z \! - \! \xi)^{\frac{\varkappa_{nk}}{\mathcal{N}}} \prod_{q=
1}^{\mathfrak{s}-1} \left(\dfrac{z \! - \! \xi}{(z \! - \! \alpha_{p_{q}})
(\xi \! - \! \alpha_{p_{q}})} \right)^{\frac{\varkappa_{nk \tilde{k}_{q}}}{
\mathcal{N}}} \right) \md \mu_{\widetilde{V}}^{\infty}(\xi) \! + \! 
\int_{J_{\infty}} \ln \left(\prod_{q=1}^{\mathfrak{s}-1}(\xi \! - \! 
\alpha_{p_{q}})^{\frac{\varkappa_{nk \tilde{k}_{q}}}{\mathcal{N}}} 
\right) \md \mu_{\widetilde{V}}^{\infty}(\xi),
\end{equation}
where the convergence, in the double-scaling limit $\mathscr{N},n 
\! \to \! \infty$ such that $z_{o} \! = \! 1 \! + \! o(1)$, is uniform 
on compact subsets of $\mathbb{C} \setminus (-\infty,
\hat{\lambda}_{1}]$. It follows {}from Equation~\eqref{eqlmeta16} 
that a corollary of Equation~\eqref{eqlmeta17} is that, for $n \! 
\in \! \mathbb{N}$ and $k \! \in \! \lbrace 1,2,\dotsc,K \rbrace$ 
such that $\alpha_{p_{\mathfrak{s}}} \! := \! \alpha_{k} \! = \! 
\infty$, for $\mathbb{C} \setminus (-\infty,\hat{\lambda}_{1}] \! 
\ni \! z$,
\begin{align} \label{eqlmeta18} 
& \, \lim_{\underset{z_{o}=1+o(1)}{\mathscr{N},n \to \infty}} \left(
\int_{\mathbb{R}} \ln \left((z \! - \! \xi)^{\frac{\varkappa_{nk}}{\mathcal{N}}} 
\prod_{q=1}^{\mathfrak{s}-1} \left(\dfrac{z \! - \! \xi}{(z \! - \! 
\alpha_{p_{q}})(\xi \! - \! \alpha_{p_{q}})} \right)^{\frac{\varkappa_{nk 
\tilde{k}_{q}}}{\mathcal{N}}} \right) \md \hat{\eta}_{\hat{\mathfrak{z}}}
(\xi) \! + \! \int_{\mathbb{R}} \ln \left(\prod_{q=1}^{\mathfrak{s}-1}
(\xi \! - \! \alpha_{p_{q}})^{\frac{\varkappa_{nk \tilde{k}_{q}}}{
\mathcal{N}}} \right) \md \hat{\eta}_{\hat{\mathfrak{z}}}(\xi) 
\right) \nonumber \\
&\underset{\underset{z_{o}=1+o(1)}{\mathscr{N},n \to \infty}}{=} 
\int_{J_{\infty}} \ln \left((z \! - \! \xi)^{\frac{\varkappa_{nk}}{
\mathcal{N}}} \prod_{q=1}^{\mathfrak{s}-1} \left(\dfrac{z \! 
- \! \xi}{(z \! - \! \alpha_{p_{q}})(\xi \! - \! \alpha_{p_{q}})} 
\right)^{\frac{\varkappa_{nk \tilde{k}_{q}}}{\mathcal{N}}} \right) 
\md \mu_{\widetilde{V}}^{\infty}(\xi) \! + \! \int_{J_{\infty}} \ln 
\left(\prod_{q=1}^{\mathfrak{s}-1}(\xi \! - \! \alpha_{p_{q}})^{
\frac{\varkappa_{nk \tilde{k}_{q}}}{\mathcal{N}}} \right) \md 
\mu_{\widetilde{V}}^{\infty}(\xi).
\end{align}
Since, for $n \! \in \! \mathbb{N}$ and $k \! \in \! \lbrace 1,
2,\dotsc,K \rbrace$ such that $\alpha_{p_{\mathfrak{s}}} \! 
:= \! \alpha_{k} \! = \! \infty$, the convergence stated in 
Equation~\eqref{eqlmeta18} is (normal) uniform on compact subsets 
of $\mathbb{C} \setminus (-\infty,\hat{\lambda}_{1}]$ $(\ni \! z)$, 
it follows by the equivalent Part~(ii) of Vitali's Theorem 
{}\footnote{In---the equivalent---Part~(ii) of Vitali's Theorem, take any 
`fixed' $z \! \in \! \mathbb{C} \setminus (-\infty,\hat{\lambda}_{1}]$ 
that is real and greater than $\hat{\lambda}_{1}$ to play the r\^{o}le 
of $c_{V}$.} and a successive differentiation (with respect to $z$) 
argument, that, for arbitrary $m \! \in \! \mathbb{N}$ and 
$z \! \in \! \mathbb{C} \setminus (-\infty,\hat{\lambda}_{1}]$,
\begin{equation*}
\lim_{\underset{z_{o}=1+o(1)}{\mathscr{N},n \to \infty}} 
\int_{\mathbb{R}}(z \! - \! \xi)^{-m} \, \md \hat{\eta}_{
\hat{\mathfrak{z}}}(\xi) \underset{\underset{z_{o}=1+o(1)}{
\mathscr{N},n \to \infty}}{=} \int_{J_{\infty}}(z \! - \! \xi)^{-m} 
\, \md \mu_{\widetilde{V}}^{\infty}(\xi):
\end{equation*} 
as $\md \hat{\eta}_{\hat{\mathfrak{z}}}$ and $\md \mu_{
\widetilde{V}}^{\infty}$ are, for $n \! \in \! \mathbb{N}$ and 
$k \! \in \! \lbrace 1,2,\dotsc,K \rbrace$ such that $\alpha_{
p_{\mathfrak{s}}} \! := \! \alpha_{k} \! = \! \infty$, supported 
on bounded compact sets, it follows via the complex form of 
the Stone-Weierstrass theorem and the latter relation that 
$\hat{\eta}_{\hat{\mathfrak{z}}} \! \overset{\ast}{\to} \! 
\mu_{\widetilde{V}}^{\infty}$ in the double-scaling limit 
$\mathscr{N},n \! \to \! \infty$ such that $z_{o} \! = \! 1 \! 
+ \! o(1)$.\footnote{As a corollary, note that, for $n \! \in \! 
\mathbb{N}$ and $k \! \in \! \lbrace 1,2,\dotsc,K \rbrace$ 
such that $\alpha_{p_{\mathfrak{s}}} \! := \! \alpha_{k} 
\! = \! \infty$, $x \! \in \! J_{\infty} \! \Leftrightarrow \! 
\operatorname{dist}(x,\hat{\mathfrak{Z}}_{0}) \! \to \! 0$ 
in the double-scaling limit $\mathscr{N},n \to \infty$ such 
that $z_{o} \! = \! 1+o(1)$.} \hfill $\qed$
\begin{ccccc} \label{lem3.8} 
Let the external field $\widetilde{V} \colon \overline{\mathbb{R}} 
\setminus \lbrace \alpha_{1},\alpha_{2},\dotsc,\alpha_{K} 
\rbrace \! \to \! \mathbb{R}$ satisfy 
conditions~\eqref{eq20}--\eqref{eq22} and be regular.

$\pmb{(1)}$ For $n \! \in \! \mathbb{N}$ and $k \! \in \! \lbrace 
1,2,\dotsc,K \rbrace$ such that $\alpha_{p_{\mathfrak{s}}} \! := \! 
\alpha_{k} \! = \! \infty$, let the associated equilibrium measure, 
$\mu_{\widetilde{V}}^{\infty}$, and its support, $J_{\infty}$, be as 
described in item~$\pmb{(1)}$ of Lemma~\ref{lem3.7}, and let there 
exist $\hat{\ell} \colon \mathbb{N} \times \lbrace 1,2,\dotsc,K \rbrace 
\! \to \! \mathbb{R}$, the associated variational constant, such that
\begin{gather} \label{eql3.8a}
\begin{split}
& \, 2 \left(\dfrac{(n \! - \! 1)K \! + \! k}{n} \right) \int_{J_{\infty}} 
\ln (\lvert z \! - \! \xi \rvert) \psi_{\widetilde{V}}^{\infty}(\xi) \, 
\md \xi \! - \! 2 \sum_{q=1}^{\mathfrak{s}-1} \dfrac{\varkappa_{nk 
\tilde{k}_{q}}}{n} \ln \lvert z \! - \! \alpha_{p_{q}} \rvert \! - \! 
\widetilde{V}(z) \! - \! \hat{\ell} \! = \! 0, \quad z \! \in \! J_{\infty}, \\
& \, 2 \left(\dfrac{(n \! - \! 1)K \! + \! k}{n} \right) \int_{J_{\infty}} 
\ln (\lvert z \! - \! \xi \rvert) \psi_{\widetilde{V}}^{\infty}(\xi) \, 
\md \xi \! - \! 2 \sum_{q=1}^{\mathfrak{s}-1} \dfrac{\varkappa_{nk 
\tilde{k}_{q}}}{n} \ln \lvert z \! - \! \alpha_{p_{q}} \rvert \! - \! 
\widetilde{V}(z) \! - \! \hat{\ell} \! \leqslant \! 0, \quad z \! \in \! 
\mathbb{R} \setminus J_{\infty},
\end{split}
\end{gather}
where, for regular $\widetilde{V}$, the inequality in the second of the 
variational conditions~\eqref{eql3.8a} is strict. Then, for $n \! \in \! 
\mathbb{N}$ and $k \! \in \! \lbrace 1,2,\dotsc,K \rbrace$ such that 
$\alpha_{p_{\mathfrak{s}}} \! := \! \alpha_{k} \! = \! \infty$$:$ {\rm (i)} 
$g_{+}^{\infty}(z) \! + \! g_{-}^{\infty}(z) \! - \! 2 \tilde{\mathscr{P}}_{0} 
\! - \! \widetilde{V}(z) \! - \! \hat{\ell} \! = \! 0$, $z \! \in \! J_{\infty}$, 
where $g^{\infty}(z)$ and $\tilde{\mathscr{P}}_{0}$ are defined by 
Equations~\eqref{eql3.4gee1} and~\eqref{eql3.4gee2}, respectively, 
and $g_{\pm}^{\infty}(z) \! := \! \lim_{\varepsilon \downarrow 0}
g^{\infty}(z \! \pm \! \mi \varepsilon)$$;$ {\rm (ii)} $g_{+}^{\infty}(z) \! 
+ \! g_{-}^{\infty}(z) \! - \! 2 \tilde{\mathscr{P}}_{0} \! - \! \widetilde{V}
(z) \! - \! \hat{\ell} \! \leqslant \! 0$, $z \! \in \! \mathbb{R} \setminus 
J_{\infty}$, where equality holds for at most a finite number of points, 
and, for regular $\widetilde{V}$, the inequality is strict$;$ {\rm (iii)}
\begin{equation*}
g_{+}^{\infty}(z) \! - \! g_{-}^{\infty}(z) \! = \! 
\begin{cases}
2 \pi \mi \hat{\mathfrak{G}}_{\Ydown}(z), &\text{$z \! \in \! 
[\hat{b}_{j-1},\hat{a}_{j}], \quad j \! = \! 1,2,\dotsc,N \! + \! 1$,} \\
2 \pi \mi \hat{\mathfrak{G}}_{\Yup}(z), &\text{$z \! \in \! 
(\hat{a}_{i},\hat{b}_{i}), \quad i \! = \! 1,2,\dotsc,N$,} \\
2 \pi \mi \hat{\mathfrak{G}}_{\Yright}(z), &\text{$z \! \in \! 
(\hat{a}_{N+1},+\infty)$,} \\
2 \pi \mi \hat{\mathfrak{G}}_{\Yleft}(z), &\text{$z \! \in \! 
(-\infty,\hat{b}_{0})$,}
\end{cases}
\end{equation*}
where
\begin{gather*}
\hat{\mathfrak{G}}_{\Ydown}(z) \! := \! \left(\dfrac{(n \! - \! 1)K \! + \! 
k}{n} \right) \int_{z}^{\hat{a}_{N+1}} \psi_{\widetilde{V}}^{\infty}(\xi) 
\, \md \xi \! - \! \sum_{q \in \hat{\Delta}_{k}(z)} \dfrac{\varkappa_{nk 
\tilde{k}_{q}}}{n}, \qquad \hat{\mathfrak{G}}_{\Yup}(z) \! := \! \left(
\dfrac{(n \! - \! 1)K \! + \! k}{n} \right) \int_{\hat{b}_{i}}^{\hat{a}_{N+1}} 
\psi_{\widetilde{V}}^{\infty}(\xi) \, \md \xi \! - \! \sum_{q \in 
\hat{\Delta}_{k}(z)} \dfrac{\varkappa_{nk \tilde{k}_{q}}}{n}, \\
\hat{\mathfrak{G}}_{\Yright}(z) \! := \! -\sum_{q \in \hat{\Delta}_{k}(z)} 
\dfrac{\varkappa_{nk \tilde{k}_{q}}}{n}, \quad \qquad \hat{\mathfrak{G}}_{
\Yleft}(z) \! := \! \left(\dfrac{(n \! - \! 1)K \! + \! k}{n} \right) \! - \! 
\sum_{q \in \hat{\Delta}_{k}(z)} \dfrac{\varkappa_{nk \tilde{k}_{q}}}{n},
\end{gather*}
with $\hat{\Delta}_{k}(z) \! := \! \lbrace \mathstrut j \! \in \! \lbrace 
1,2,\dotsc,\mathfrak{s} \! - \! 1 \rbrace; \, \alpha_{p_{j}} \! > \! z 
\rbrace$$;$\footnote{If $\hat{\Delta}_{k}(z) \! = \! \varnothing$, then 
$\sum_{q \in \hat{\Delta}_{k}(z)} \varkappa_{nk \tilde{k}_{q}}/n \! := 
\! 0$.} and {\rm (iv)}
\begin{equation*}
\mi \left(g_{+}^{\infty}(z) \! - \! g_{-}^{\infty}(z) \! + \! 2 \pi \mi 
\sum_{q \in \hat{\Delta}_{k}(z)} \dfrac{\varkappa_{nk \tilde{k}_{q}}}{n} 
\right)^{\prime} \! = \! 
\begin{cases}
2 \pi \left(\dfrac{(n \! - \! 1)K \! + \! k}{n} \right) \psi_{\widetilde{V}}^{
\infty}(z) \! \geqslant \! 0, &\text{$z \! \in \! J_{\infty}$,} \\
0, &\text{$z \! \in \! \mathbb{R} \setminus J_{\infty}$,}
\end{cases}
\end{equation*}
where the prime denotes differentiation with respect to $z$, and, for 
regular $\widetilde{V}$, equality holds for at most a finite number 
of points.

$\pmb{(2)}$ For $n \! \in \! \mathbb{N}$ and $k \! \in \! \lbrace 
1,2,\dotsc,K \rbrace$ such that $\alpha_{p_{\mathfrak{s}}} \! := \! 
\alpha_{k} \! \neq \! \infty$, let the associated equilibrium measure, 
$\mu_{\widetilde{V}}^{f}$, and its support, $J_{f}$, be as described 
in item~$\pmb{(2)}$ of Lemma~\ref{lem3.7}, and let there exist 
$\tilde{\ell} \colon \mathbb{N} \times \lbrace 1,2,\dotsc,K \rbrace 
\! \to \! \mathbb{R}$, the associated variational constant, such that
\begin{gather} \label{eql3.8b}
\begin{split}
&2 \left(\dfrac{(n \! - \! 1)K \! + \! k}{n} \right) \int_{J_{f}} \ln \left(
\left\lvert \dfrac{z \! - \! \xi}{\xi \! - \! \alpha_{k}} \right\rvert \right) 
\psi_{\widetilde{V}}^{f}(\xi) \, \md \xi \! - \! 2 \sum_{q=1}^{\mathfrak{s}-2} 
\dfrac{\varkappa_{nk \tilde{k}_{q}}}{n} \ln \left\lvert \dfrac{z \! - \! 
\alpha_{p_{q}}}{\alpha_{p_{q}} \! - \! \alpha_{k}} \right\rvert \! - \! 2 
\left(\dfrac{\varkappa_{nk} \! - \! 1}{n} \right) \ln \lvert z \! - \! 
\alpha_{k} \rvert \! - \! \widetilde{V}(z) \! - \! \tilde{\ell} \! = \! 0, 
\quad z \! \in \! J_{f}, \\
&2 \left(\dfrac{(n \! - \! 1)K \! + \! k}{n} \right) \int_{J_{f}} \ln \left(
\left\lvert \dfrac{z \! - \! \xi}{\xi \! - \! \alpha_{k}} \right\rvert \right) 
\psi_{\widetilde{V}}^{f}(\xi) \, \md \xi \! - \! 2 \sum_{q=1}^{\mathfrak{s}-2} 
\dfrac{\varkappa_{nk \tilde{k}_{q}}}{n} \ln \left\lvert \dfrac{z \! - \! 
\alpha_{p_{q}}}{\alpha_{p_{q}} \! - \! \alpha_{k}} \right\rvert \! - \! 2 
\left(\dfrac{\varkappa_{nk} \! - \! 1}{n} \right) \ln \lvert z \! - \! 
\alpha_{k} \rvert \! - \! \widetilde{V}(z) \! - \! \tilde{\ell} \! \leqslant 
\! 0, \quad z \! \in \! \mathbb{R} \setminus J_{f},
\end{split}
\end{gather}
where, for regular $\widetilde{V}$, the inequality in the second of the 
variational conditions~\eqref{eql3.8b} is strict. Then, for $n \! \in \! 
\mathbb{N}$ and $k \! \in \! \lbrace 1,2,\dotsc,K \rbrace$ such that 
$\alpha_{p_{\mathfrak{s}}} \! := \! \alpha_{k} \! \neq \! \infty$$:$ 
{\rm (i)} $g_{+}^{f}(z) \! + \! g_{-}^{f}(z) \! - \! \hat{\mathscr{P}}_{0}^{+} 
\! - \! \hat{\mathscr{P}}_{0}^{-} \! - \! \widetilde{V}(z) \! - \! \tilde{\ell} 
\! = \! 0$, $z \! \in \! J_{f}$, where $g^{f}(z)$ and $\hat{\mathscr{
P}}_{0}^{\pm}$ are defined by Equations~\eqref{eql3.4gee3} 
and~\eqref{eql3.4gee5}, respectively, and $g_{\pm}^{f}(z) \! 
:= \! \lim_{\varepsilon \downarrow 0}g^{f}(z \! \pm \! \mi 
\varepsilon)$$;$ {\rm (ii)} $g_{+}^{f}(z) \! + \! g_{-}^{f}(z) \! - \! 
\hat{\mathscr{P}}_{0}^{+} \! - \! \hat{\mathscr{P}}_{0}^{-} \! - \! 
\widetilde{V}(z) \! - \! \tilde{\ell} \! \leqslant \! 0$, $z \! \in \! 
\mathbb{R} \setminus J_{f}$, where equality holds for at most a finite 
number of points, and, for regular $\widetilde{V}$, the inequality is 
strict$;$ {\rm (iii)}
\begin{equation*}
g_{+}^{f}(z) \! - \! g_{-}^{f}(z) \! + \! \hat{\mathscr{P}}_{0}^{-} \! - \! 
\hat{\mathscr{P}}_{0}^{+} \! = \! 
\begin{cases}
2 \pi \mi \tilde{\mathfrak{G}}_{\Ydown}(z), &\text{$z \! \in \! 
[\tilde{b}_{j-1},\tilde{a}_{j}], \quad j \! = \! 1,2,\dotsc,N \! + \! 1$,} \\
2 \pi \mi \tilde{\mathfrak{G}}_{\Yup}(z), &\text{$z \! \in \! 
(\tilde{a}_{i},\tilde{b}_{i}), \quad i \! = \! 1,2,\dotsc,N$,} \\
2 \pi \mi \tilde{\mathfrak{G}}_{\Yright}(z), &\text{$z \! \in \! 
(\tilde{a}_{N+1},+\infty)$,} \\
2 \pi \mi \tilde{\mathfrak{G}}_{\Yleft}(z), &\text{$z \! \in \! 
(-\infty,\tilde{b}_{0})$,}
\end{cases}
\end{equation*}
where
\begin{align*}
\tilde{\mathfrak{G}}_{\Ydown}(z) :=& \, -\left(\dfrac{\varkappa_{nk} \! - \! 
1}{n} \right) \chi_{\mathbb{R}_{\alpha_{k}}^{<}}(z) \! + \! \left(
\dfrac{(n \! - \! 1)K \! + \! k}{n} \right) \int_{z}^{\tilde{a}_{N+1}} 
\psi_{\widetilde{V}}^{f}(\xi) \, \md \xi \! - \! \left(\dfrac{(n \! - \! 1)K 
\! + \! k}{n} \right) \int_{J_{f} \cap \mathbb{R}_{\alpha_{k}}^{>}} \psi_{
\widetilde{V}}^{f}(\xi) \, \md \xi \\
-& \, \sum_{q \in \tilde{\Delta}_{k}(z)} \dfrac{\varkappa_{nk 
\tilde{k}_{q}}}{n} \! + \! \sum_{q \in \tilde{\Delta}(k)} \dfrac{\varkappa_{nk 
\tilde{k}_{q}}}{n}, \\
\tilde{\mathfrak{G}}_{\Yup}(z) :=& \, -\left(\dfrac{\varkappa_{nk} \! - \! 
1}{n} \right) \chi_{\mathbb{R}_{\alpha_{k}}^{<}}(z) \! + \! 
\left(\dfrac{(n \! - \! 1)K \! + \! k}{n} \right) \int_{\tilde{b}_{i}}^{
\tilde{a}_{N+1}} \psi_{\widetilde{V}}^{f}(\xi) \, \md \xi \! - \! \left(
\dfrac{(n \! - \! 1)K \! + \! k}{n} \right) \int_{J_{f} \cap \mathbb{R}_{
\alpha_{k}}^{>}} \psi_{\widetilde{V}}^{f}(\xi) \, \md \xi \\
-& \, \sum_{q \in \tilde{\Delta}_{k}(z)} \dfrac{\varkappa_{nk 
\tilde{k}_{q}}}{n} \! + \! \sum_{q \in \tilde{\Delta}(k)} \dfrac{\varkappa_{nk 
\tilde{k}_{q}}}{n}, \\
\tilde{\mathfrak{G}}_{\Yright}(z) :=& \, -\left(\dfrac{\varkappa_{nk} \! - \! 
1}{n} \right) \chi_{\mathbb{R}_{\alpha_{k}}^{<}}(z) \! - \! \left(
\dfrac{(n \! - \! 1)K \! + \! k}{n} \right) \int_{J_{f} \cap \mathbb{R}_{
\alpha_{k}}^{>}} \psi_{\widetilde{V}}^{f}(\xi) \, \md \xi \! - \! \sum_{q 
\in \tilde{\Delta}_{k}(z)} \dfrac{\varkappa_{nk \tilde{k}_{q}}}{n} \! + \! 
\sum_{q \in \tilde{\Delta}(k)} \dfrac{\varkappa_{nk \tilde{k}_{q}}}{n}, \\
\tilde{\mathfrak{G}}_{\Yleft}(z) :=& \, -\left(\dfrac{\varkappa_{nk} \! - \! 
1}{n} \right) \chi_{\mathbb{R}_{\alpha_{k}}^{<}}(z) \! - \! \left(
\dfrac{(n \! - \! 1)K \! + \! k}{n} \right) \int_{J_{f} \cap \mathbb{R}_{
\alpha_{k}}^{>}} \psi_{\widetilde{V}}^{f}(\xi) \, \md \xi \! + \! \left(
\dfrac{(n \! - \! 1)K \! + \! k}{n} \right) \! - \! \sum_{q \in \tilde{
\Delta}_{k}(z)} \dfrac{\varkappa_{nk \tilde{k}_{q}}}{n} \! + \! \sum_{q \in 
\tilde{\Delta}(k)} \dfrac{\varkappa_{nk \tilde{k}_{q}}}{n},
\end{align*}
with {}\footnote{If $J_{f} \cap \mathbb{R}_{\alpha_{k}}^{>} \! = \! \varnothing$, 
then $\int_{J_{f} \cap \mathbb{R}_{\alpha_{k}}^{>}} \psi_{\widetilde{V}}^{f}
(\xi) \, \md \xi \! := \! 0$.} $\tilde{\Delta}_{k}(z) \! := \! \lbrace 
\mathstrut j \! \in \! \lbrace 1,2,\dotsc,\mathfrak{s} \! - \! 2 \rbrace; \, 
\alpha_{p_{j}} \! > \! z \rbrace$,\footnote{If $\tilde{\Delta}_{k}(z) \! = 
\! \varnothing$, then $\sum_{q \in \tilde{\Delta}_{k}(z)} \varkappa_{nk 
\tilde{k}_{q}}/n \! := \! 0$.} and $\tilde{\Delta}(k) \! := \! \lbrace 
\mathstrut j \! \in \! \lbrace 1,2,\dotsc,\mathfrak{s} \! - \! 2 \rbrace; \, 
\alpha_{p_{j}} \! > \! \alpha_{k} \rbrace$$;$\footnote{If $\tilde{\Delta}(k) 
\! = \! \varnothing$, then $\sum_{q \in \tilde{\Delta}(k)} \varkappa_{nk 
\tilde{k}_{q}}/n \! := \! 0$.} and {\rm (iv)}
\begin{equation*}
\mi \left(g_{+}^{f}(z) \! - \! g_{-}^{f}(z) \! + \! \hat{\mathscr{P}}_{0}^{-} 
\! - \! \hat{\mathscr{P}}_{0}^{+} \! + \! 2 \pi \mi \left(
\dfrac{\varkappa_{nk} \! - \! 1}{n} \right) \chi_{\mathbb{R}_{
\alpha_{k}}^{<}}(z) \! + \! 2 \pi \mi \sum_{q \in \tilde{\Delta}_{k}(z)} 
\dfrac{\varkappa_{nk \tilde{k}_{q}}}{n} \right)^{\prime} \! = \! 
\begin{cases}
2 \pi \left(\dfrac{(n \! - \! 1)K \! + \! k}{n} \right) \psi_{\widetilde{
V}}^{f}(z) \! \geqslant \! 0, &\text{$z \! \in \! J_{f}$,} \\
0, &\text{$z \! \in \! \mathbb{R} \setminus J_{f}$,}
\end{cases}
\end{equation*}
where, for regular $\widetilde{V}$, equality holds for at most a finite number 
of points.
\end{ccccc}

\emph{Proof}. The proof of this Lemma~\ref{lem3.8} consists of two cases: (i) 
$n \! \in \! \mathbb{N}$ and $k \! \in \! \lbrace 1,2,\dotsc,K \rbrace$ such 
that $\alpha_{p_{\mathfrak{s}}} \! := \! \alpha_{k} \! = \! \infty$; and (ii) $n 
\! \in \! \mathbb{N}$ and $k \! \in \! \lbrace 1,2,\dotsc,K \rbrace$ such that 
$\alpha_{p_{\mathfrak{s}}} \! := \! \alpha_{k} \! \neq \! \infty$. The proof 
for the case $\alpha_{p_{\mathfrak{s}}} \! := \! \alpha_{k} \! \neq \! \infty$, 
$k \! \in \! \lbrace 1,2,\dotsc,K \rbrace$, will be considered in detail (see 
$\pmb{(\mathrm{A})}$ below), whilst the case $\alpha_{p_{\mathfrak{s}}} \! 
:= \! \alpha_{k} \! = \! \infty$, $k \! \in \! \lbrace 1,2,\dotsc,K \rbrace$, can, 
modulo technical and computational particulars, be proved analogously (see 
$\pmb{(\mathrm{B})}$ below).

$\pmb{(\mathrm{A})}$ Let $\widetilde{V} \colon \overline{\mathbb{R}} 
\setminus \lbrace \alpha_{1},\alpha_{2},\dotsc,\alpha_{K} \rbrace \! \to 
\! \mathbb{R}$ satisfy conditions~\eqref{eq20}--\eqref{eq22} and be 
regular. Recall {}from item~$\pmb{(2)}$ of Lemma~\ref{lem3.7} that, for 
$n \! \in \! \mathbb{N}$ and $k \! \in \! \lbrace 1,2,\dotsc,K \rbrace$ 
such that $\alpha_{p_{\mathfrak{s}}} \! := \! \alpha_{k} \! \neq \! \infty$, 
$\supp (\mu^{f}_{\widetilde{V}}) \! =: \! J_{f} \! = \! \cup_{j=1}^{N+1} 
\tilde{J}_{j}$ $(\subset \overline{\mathbb{R}} \setminus \lbrace 
\alpha_{p_{1}},\alpha_{p_{2}},\dotsc,\alpha_{p_{\mathfrak{s}}} \rbrace)$, 
where $\tilde{J}_{j} \! := \! [\tilde{b}_{j-1},\tilde{a}_{j}]$, $j \! = \! 1,2,
\dotsc,N \! + \! 1$, and $\tilde{b}_{j-1},\tilde{a}_{j}$ satisfy, in the 
double-scaling limit $\mathscr{N},n \! \to \! \infty$ such that $z_{o} \! 
= \! 1 \! + \! o(1)$, the locally solvable system of $2(N \! + \! 1)$ real 
moment equations~\eqref{eql3.7g}--\eqref{eql3.7i}, with $N \! \in 
\! \mathbb{N}_{0}$ and finite, $\tilde{J}_{i} \cap \tilde{J}_{j} \! = \! 
\varnothing$, $i \! \neq \! j \! \in \! \lbrace 1,2,\dotsc,N \! + \! 1 
\rbrace$, and $-\infty \! < \! \tilde{b}_{0} \! < \! \tilde{a}_{1} \! < \! 
\tilde{b}_{1} \! < \! \tilde{a}_{2} \! < \! \dotsb \! < \! \tilde{b}_{N} \! 
< \! \tilde{a}_{N+1} \! < \! +\infty$. Write $\mathbb{R} \! = \! J_{f} 
\cup (\cup_{i=1}^{N}(\tilde{a}_{i},\tilde{b}_{i})) \cup (\tilde{a}_{N+1},
+\infty) \cup (-\infty,\tilde{b}_{0})$. Consider the following cases: 
$\pmb{(1)}$ $z \! \in \! \tilde{J}_{j}$, $j \! = \! 1,2,\dotsc,N \! + \! 1$; 
$\pmb{(2)}$ $z \! \in \! (\tilde{a}_{i},\tilde{b}_{i})$, $i \! = \! 1,2,\dotsc,
N$; $\pmb{(3)}$ $z \! \in \! (\tilde{a}_{N+1},+\infty)$; and $\pmb{(4)}$ 
$z \! \in \! (-\infty,\tilde{b}_{0})$.

$\pmb{(1)}$ Substituting the representation~\eqref{eql3.7k} for the 
density of the associated equilibrium measure into the definition of 
$g^{f}(z)$ given by Equation~\eqref{eql3.4gee3}, one shows that, for 
$z \! \in \! \tilde{J}_{j}$, $j \! = \! 1,2,\dotsc,N \! + \! 1$,
\begin{align} \label{eql3.8c} 
g^{f}_{\pm}(z) =& \, \left(\dfrac{(n \! - \! 1)K \! + \! k}{n} \right) 
\int_{J_{f}} \ln (\lvert z \! - \! \xi \rvert) \psi^{f}_{\widetilde{V}}(\xi) 
\, \md \xi \! - \! \left(\dfrac{\varkappa_{nk} \! - \! 1}{n} \right) 
\int_{J_{f}} \ln (\lvert \xi \! - \! \alpha_{k} \rvert) 
\psi^{f}_{\widetilde{V}}(\xi) \, \md \xi \nonumber \\
-& \, \mi \pi \left(\dfrac{\varkappa_{nk} \! - \! 1}{n} \right) \int_{J_{f} 
\cap \mathbb{R}_{\alpha_{k}}^{<}} \psi^{f}_{\widetilde{V}}(\xi) \, \md \xi 
\! - \! \sum_{q=1}^{\mathfrak{s}-2} \dfrac{\varkappa_{nk \tilde{k}_{q}}}{n} 
\int_{J_{f}} \ln (\lvert \xi \! - \! \alpha_{p_{q}} \rvert) \psi^{f}_{
\widetilde{V}}(\xi) \, \md \xi \nonumber \\
-& \, \mi \pi \sum_{q=1}^{\mathfrak{s}-2} \dfrac{\varkappa_{nk \tilde{k}_{
q}}}{n} \int_{J_{f} \cap \mathbb{R}_{\alpha_{p_{q}}}^{<}} \psi^{f}_{
\widetilde{V}}(\xi) \, \md \xi \! - \! \sum_{q=1}^{\mathfrak{s}-2} \dfrac{
\varkappa_{nk \tilde{k}_{q}}}{n} \ln \lvert z \! - \! \alpha_{p_{q}} \rvert 
\! \mp \! \mi \pi \sum_{q \in \tilde{\Delta}(j;z)} \dfrac{\varkappa_{nk 
\tilde{k}_{q}}}{n} \nonumber \\
\pm& \, \mi \pi \left(\dfrac{(n \! - \! 1)K \! + \! k}{n} \right) \int_{z}^{
\tilde{a}_{N+1}} \psi^{f}_{\widetilde{V}}(\xi) \, \md \xi \! - \! \left(
\dfrac{\varkappa_{nk} \! - \! 1}{n} \right) \left(\ln \lvert z \! - \! 
\alpha_{k} \rvert \! \pm \! \mi \pi \chi_{\mathbb{R}^{<}_{\alpha_{k}}}
(z) \right),
\end{align}
where {}\footnote{If $J_{f} \cap \mathbb{R}_{\alpha_{p_{q}}}^{<} \! = 
\! \varnothing$, $q \! \in \! \lbrace 1,\dotsc,\mathfrak{s} \! - \! 2,
\mathfrak{s} \rbrace$, then $\int_{J_{f} \cap \mathbb{R}_{\alpha_{
p_{q}}}^{<}} \psi^{f}_{\widetilde{V}}(\xi) \, \md \xi \! := \! 0$.} 
$g^{f}_{\pm}(z) \! := \! \lim_{\varepsilon \downarrow 0}g^{f}(z \! \pm 
\! \mi \varepsilon)$, and $\tilde{\Delta}(j;z) \! := \! \lbrace \mathstrut 
i \! \in \! \lbrace 1,2,\dotsc,\mathfrak{s} \! - \! 2 \rbrace; \, \alpha_{p_{i}} 
\! > \! z \rbrace$.\footnote{If $\tilde{\Delta}(j;z) \! = \! \varnothing$, 
$j \! \in \! \lbrace 1,2,\dotsc,N \! + \! 1 \rbrace$, then $\sum_{q 
\in \tilde{\Delta}(j;z)} \varkappa_{nk \tilde{k}_{q}}/n \! := \! 0$.} Via 
Equation~\eqref{eql3.8c} and the definition of $\hat{\mathscr{P}}_{
0}^{\pm}$ given by Equation~\eqref{eql3.4gee5}, one arrives at, for 
$n \! \in \! \mathbb{N}$ and $k \! \in \! \lbrace 1,2,\dotsc,K \rbrace$ 
such that $\alpha_{p_{\mathfrak{s}}} \! := \! \alpha_{k} \! \neq \! \infty$,
\begin{align} \label{eql3.8B10} 
(2 \pi \mi)^{-1}(g^{f}_{+}(z) \! - \! g^{f}_{-}(z) \! + \! 
\hat{\mathscr{P}}_{0}^{-} \! - \! \hat{\mathscr{P}}_{0}^{+}) =& \, \left(
\dfrac{(n \! - \! 1)K \! + \! k}{n} \right) \int_{z}^{\tilde{a}_{N+1}} 
\psi^{f}_{\widetilde{V}}(\xi) \, \md \xi \! - \! \left(\dfrac{(n \! - \! 
1)K \! + \! k}{n} \right) \int_{J_{f} \cap \mathbb{R}_{\alpha_{k}}^{>}} 
\psi^{f}_{\widetilde{V}}(\xi) \, \md \xi \nonumber \\
-& \, \sum_{q \in \tilde{\Delta}(j;z)} \dfrac{\varkappa_{nk \tilde{k}_{q}}}{n} 
\! + \! \sum_{q \in \tilde{\Delta}(k)} \dfrac{\varkappa_{nk \tilde{k}_{q}}}{n} 
\! - \! \left(\dfrac{\varkappa_{nk} \! - \! 1}{n} \right) 
\chi_{\mathbb{R}_{\alpha_{k}}^{<}}(z),
\end{align}
where $\tilde{\Delta}(k) \! := \! \lbrace \mathstrut j \! \in \! \lbrace 1,2,
\dotsc,\mathfrak{s} \! - \! 2 \rbrace; \, \alpha_{p_{j}} \! > \! \alpha_{k} 
\rbrace$,\footnote{If $\tilde{\Delta}(k) \! = \! \varnothing$, then $\sum_{q 
\in \tilde{\Delta}(k)} \varkappa_{nk \tilde{k}_{q}}/n \! := \! 0$.} which 
shows that, for $z \! \in \! \tilde{J}_{j}$, $j \! = \! 1,2,\dotsc,N \! + \! 
1$, $g^{f}_{+}(z) \! - \! g^{f}_{-}(z) \! + \! \hat{\mathscr{P}}_{0}^{-} \! - 
\! \hat{\mathscr{P}}_{0}^{+} \! \in \! \mi \mathbb{R}$; moreover, via 
Leibnitz's Rule and the fact that $\psi^{f}_{\widetilde{V}}(z) \! \geqslant 
\! 0$, $z \! \in \! \tilde{J}_{j}$, $j \! = \! 1,2,\dotsc,N \! + \! 1$, it also 
follows that
\begin{equation} \label{eql3.8D1} 
\mi \left(g^{f}_{+}(z) \! - \! g^{f}_{-}(z) \! + \! \hat{\mathscr{P}}_{0}^{-} 
\! - \! \hat{\mathscr{P}}_{0}^{+} \! + \! 2 \pi \mi \left(
\dfrac{\varkappa_{nk} \! - \! 1}{n} \right) \chi_{\mathbb{R}_{
\alpha_{k}}^{<}}(z) \! + \! 2 \pi \mi \sum_{q \in \tilde{\Delta}(j;z)} 
\dfrac{\varkappa_{nk \tilde{k}_{q}}}{n} \right)^{\prime} \! = \! 2 \pi 
\left(\dfrac{(n \! - \! 1)K \! + \! k}{n} \right) \psi^{f}_{\widetilde{V}}(z) 
\! \geqslant \! 0,
\end{equation}
where the prime denotes differentiation with respect to $z$. Via 
Equation~\eqref{eql3.8c}, the definition of $\hat{\mathscr{P}}_{0}^{\pm}$ 
given by Equation~\eqref{eql3.4gee5}, and the first of the variational 
conditions~\eqref{eql3.8b}, one also arrives at, for $n \! \in \! \mathbb{N}$ 
and $k \! \in \! \lbrace 1,2,\dotsc,K \rbrace$ such that $\alpha_{p_{\mathfrak{s}}} 
\! := \! \alpha_{k} \! \neq \! \infty$,
\begin{align} \label{eql3.8d} 
g^{f}_{+}(z) \! + \! g^{f}_{-}(z) \! - \! \hat{\mathscr{P}}_{0}^{-} \! - \! 
\hat{\mathscr{P}}_{0}^{+} \! - \! \widetilde{V}(z) \! - \! \tilde{\ell} =& \, 
2 \left(\dfrac{(n \! - \! 1)K \! + \! k}{n} \right) \int_{J_{f}} \ln \left(
\left\lvert \dfrac{z \! - \! \xi}{\xi \! - \! \alpha_{k}} \right\rvert \right) 
\psi^{f}_{\widetilde{V}}(\xi) \, \md \xi \! - \! 2 \sum_{q=1}^{\mathfrak{s}-2} 
\dfrac{\varkappa_{nk \tilde{k}_{q}}}{n} \ln \left\lvert \dfrac{z \! - \! 
\alpha_{p_{q}}}{\alpha_{p_{q}} \! - \! \alpha_{k}} \right\rvert \nonumber \\
-& \, 2 \left(\dfrac{\varkappa_{nk} \! - \! 1}{n} \right) \ln \lvert z \! - 
\! \alpha_{k} \rvert \! - \! \widetilde{V}(z) \! - \! \tilde{\ell} \! = \! 0,
\end{align}
which, as a by-product, gives the formula for the associated variational 
constant $\tilde{\ell} \colon \mathbb{N} \times \lbrace 1,2,\dotsc,K \rbrace 
\! \to \! \mathbb{R}$, which is the same \cite{a57} (see, also, Section~7 of 
\cite{a58}) for each compact interval $\tilde{J}_{j}$, $j \! = \! 1,2,\dotsc,
N \! + \! 1$; in particular,
\begin{align} \label{eql3.8ellwave} 
\tilde{\ell} =& \, 2 \left(\dfrac{(n \! - \! 1)K \! + \! k}{n} \right) 
\int_{J_{f}} \ln \left(\left\lvert \dfrac{\frac{1}{2}(\tilde{b}_{N} \! + \! 
\tilde{a}_{N+1}) \! - \! \xi}{\xi \! - \! \alpha_{k}} \right\rvert \right) 
\psi^{f}_{\widetilde{V}}(\xi) \, \md \xi \! - \! 2 \sum_{q=1}^{\mathfrak{s}-2} 
\dfrac{\varkappa_{nk \tilde{k}_{q}}}{n} \ln \left\lvert \dfrac{\tfrac{1}{2}
(\tilde{b}_{N} \! + \! \tilde{a}_{N+1}) \! - \! \alpha_{p_{q}}}{\alpha_{p_{q}} 
\! - \! \alpha_{k}} \right\rvert \nonumber \\
-& \, 2 \left(\dfrac{\varkappa_{nk} \! - \! 1}{n} \right) \ln \left\lvert 
\tfrac{1}{2}(\tilde{b}_{N} \! + \! \tilde{a}_{N+1}) \! - \! \alpha_{k} 
\right\rvert \! - \! \widetilde{V}(\tfrac{1}{2}(\tilde{b}_{N} \! + \! 
\tilde{a}_{N+1})).
\end{align}

$\pmb{(2)}$ Substituting the representation~\eqref{eql3.7k} for the density 
of the associated equilibrium measure into the definition of $g^{f}(z)$ 
given by Equation~\eqref{eql3.4gee3}, one shows that, for $z \! \in \! 
(\tilde{a}_{j},\tilde{b}_{j})$, $j \! = \! 1,2,\dotsc,N$,
\begin{align} \label{eql3.8e} 
g^{f}_{\pm}(z) =& \, \left(\dfrac{(n \! - \! 1)K \! + \! k}{n} \right) 
\int_{J_{f}} \ln (\lvert z \! - \! \xi \rvert) \psi^{f}_{\widetilde{V}}
(\xi) \, \md \xi \! - \! \left(\dfrac{\varkappa_{nk} \! - \! 1}{n} 
\right) \int_{J_{f}} \ln (\lvert \xi \! - \! \alpha_{k} \rvert) 
\psi^{f}_{\widetilde{V}}(\xi) \, \md \xi \nonumber \\
-& \, \mi \pi \left(\dfrac{\varkappa_{nk} \! - \! 1}{n} \right) \int_{J_{f} 
\cap \mathbb{R}_{\alpha_{k}}^{<}} \psi^{f}_{\widetilde{V}}(\xi) \, \md \xi 
\! - \! \sum_{q=1}^{\mathfrak{s}-2} \dfrac{\varkappa_{nk \tilde{k}_{q}}}{n} 
\int_{J_{f}} \ln (\lvert \xi \! - \! \alpha_{p_{q}} \rvert) 
\psi^{f}_{\widetilde{V}}(\xi) \, \md \xi \nonumber \\
-& \, \mi \pi \sum_{q=1}^{\mathfrak{s}-2} \dfrac{\varkappa_{nk 
\tilde{k}_{q}}}{n} \int_{J_{f} \cap \mathbb{R}_{\alpha_{p_{q}}}^{<}} 
\psi^{f}_{\widetilde{V}}(\xi) \, \md \xi \! - \! \sum_{q=1}^{\mathfrak{s}-2} 
\dfrac{\varkappa_{nk \tilde{k}_{q}}}{n} \ln \lvert z \! - \! \alpha_{p_{q}} 
\rvert \! \mp \! \mi \pi \sum_{q \in \tilde{\Delta}(j;z)} 
\dfrac{\varkappa_{nk \tilde{k}_{q}}}{n} \nonumber \\
\pm& \, \mi \pi \left(\dfrac{(n \! - \! 1)K \! + \! k}{n} \right) 
\int_{\tilde{b}_{j}}^{\tilde{a}_{N+1}} \psi^{f}_{\widetilde{V}}(\xi) 
\, \md \xi \! - \! \left(\dfrac{\varkappa_{nk} \! - \! 1}{n} \right) 
\left(\ln \lvert z \! - \! \alpha_{k} \rvert \! \pm \! \mi \pi 
\chi_{\mathbb{R}^{<}_{\alpha_{k}}}(z) \right);
\end{align}
hence, via Equation~\eqref{eql3.8e} and the definition of 
$\hat{\mathscr{P}}_{0}^{\pm}$ given by Equation~\eqref{eql3.4gee5}, 
one shows that, for $n \! \in \! \mathbb{N}$ and $k \! \in \! \lbrace 
1,2,\dotsc,K \rbrace$ such that $\alpha_{p_{\mathfrak{s}}} \! := \! 
\alpha_{k} \! \neq \! \infty$,
\begin{align} \label{eql3.8B11} 
(2 \pi \mi)^{-1}(g^{f}_{+}(z) \! - \! g^{f}_{-}(z) \! + \! 
\hat{\mathscr{P}}_{0}^{-} \! - \! \hat{\mathscr{P}}_{0}^{+}) =& \, 
\left(\dfrac{(n \! - \! 1)K \! + \! k}{n} \right) \int_{\tilde{b}_{j}}^{
\tilde{a}_{N+1}} \psi^{f}_{\widetilde{V}}(\xi) \, \md \xi \! - \! \left(
\dfrac{(n \! - \! 1)K \! + \! k}{n} \right) \int_{J_{f} \cap \mathbb{R}_{
\alpha_{k}}^{>}} \psi^{f}_{\widetilde{V}}(\xi) \, \md \xi \nonumber \\
-& \, \sum_{q \in \tilde{\Delta}(j;z)} \dfrac{\varkappa_{nk \tilde{k}_{q}}}{n} 
\! + \! \sum_{q \in \tilde{\Delta}(k)} \dfrac{\varkappa_{nk \tilde{k}_{q}}}{n} 
\! - \! \left(\dfrac{\varkappa_{nk} \! - \! 1}{n} \right) 
\chi_{\mathbb{R}_{\alpha_{k}}^{<}}(z),
\end{align}
whence $g^{f}_{+}(z) \! - \! g^{f}_{-}(z) \! + \! \hat{\mathscr{P}}_{0}^{-} 
\! - \! \hat{\mathscr{P}}_{0}^{+} \! \in \! \mi \mathbb{R}$; moreover, 
via Leibnitz's Rule, it follows that, for $z \! \in \! (\tilde{a}_{j},\tilde{b}_{j})$, 
$j \! = \! 1,2,\dotsc,N$,
\begin{equation} \label{eql3.8D2} 
\mi \left(g^{f}_{+}(z) \! - \! g^{f}_{-}(z) \! + \! \hat{\mathscr{P}}_{0}^{-} 
\! - \! \hat{\mathscr{P}}_{0}^{+} \! + \! 2 \pi \mi \left(
\dfrac{\varkappa_{nk} \! - \! 1}{n} \right) \chi_{\mathbb{R}_{\alpha_{k}}^{<}}
(z) \! + \! 2 \pi \mi \sum_{q \in \tilde{\Delta}(j;z)} \dfrac{\varkappa_{nk 
\tilde{k}_{q}}}{n} \right)^{\prime} \! = \! 0.
\end{equation}
Via Equation~\eqref{eql3.8e} and the definition of 
$\hat{\mathscr{P}}_{0}^{\pm}$ given by Equation~\eqref{eql3.4gee5}, one 
arrives at, for $n \! \in \! \mathbb{N}$ and $k \! \in \! \lbrace 1,2,\dotsc,K 
\rbrace$ such that $\alpha_{p_{\mathfrak{s}}} \! := \!\alpha_{k} \! \neq \! 
\infty$, and $z \! \in \! (\tilde{a}_{j},\tilde{b}_{j})$, $j \! = \! 1,2,\dotsc,N$,
\begin{align} \label{eql3.8f} 
g^{f}_{+}(z) \! + \! g^{f}_{-}(z) \! - \! \hat{\mathscr{P}}_{0}^{-} \! - \! 
\hat{\mathscr{P}}_{0}^{+} \! - \! \widetilde{V}(z) \! - \! \tilde{\ell} =& \, 
2 \left(\dfrac{(n \! - \! 1)K \! + \! k}{n} \right) \int_{J_{f}} \ln \left(
\left\lvert \dfrac{z \! - \! \xi}{\xi \! - \! \alpha_{k}} \right\rvert \right) 
\psi^{f}_{\widetilde{V}}(\xi) \, \md \xi \! - \! 2 \sum_{q=1}^{\mathfrak{s}-2} 
\dfrac{\varkappa_{nk \tilde{k}_{q}}}{n} \ln \left\lvert \dfrac{z \! - \! 
\alpha_{p_{q}}}{\alpha_{p_{q}} \! - \! \alpha_{k}} \right\rvert \nonumber \\
-& \, 2 \left(\dfrac{\varkappa_{nk} \! - \! 1}{n} \right) \ln \lvert z \! 
- \! \alpha_{k} \rvert \! - \! \widetilde{V}(z) \! - \! \tilde{\ell}.
\end{align}
Recalling {}from the proof of Lemma~\ref{lem3.7} that $\mathcal{H} \colon 
\mathcal{L}^{2}(\mathbb{R}) \! \ni \! f \! \mapsto \! \lim_{\varepsilon 
\downarrow 0} \int_{\lvert z-\xi \rvert > \varepsilon} \tfrac{f(\xi)}{z-\xi} 
\, \tfrac{\md \xi}{\pi} \! =: \! (\mathcal{H}f)(z)$ denotes the Hilbert 
transform, one shows that, for $z \! \in \! (\tilde{a}_{j},\tilde{b}_{j})$, 
$j \! = \! 1,2,\dotsc,N$,
\begin{equation} \label{eql3.8g} 
\int_{J_{f}} \ln (\lvert z \! - \! \xi \rvert) \psi_{\widetilde{V}}^{f}
(\xi) \, \md \xi \! = \! \pi \int_{\tilde{a}_{j}}^{z}(\mathcal{H} 
\psi_{\widetilde{V}}^{f})(\xi) \, \md \xi \! + \! \int_{J_{f}} \ln (\lvert 
\tilde{a}_{j} \! - \! \xi \rvert) \psi_{\widetilde{V}}^{f}(\xi) \, \md \xi;
\end{equation}
however, some of the poles $\alpha_{p_{q}}$, $q \! \in \! \lbrace 
1,\dotsc,\mathfrak{s} \! - \! 2,\mathfrak{s} \rbrace$, may lie between 
$\tilde{a}_{j}$ and $z$, in which case, the integral $\int_{\tilde{a}_{j}}^{z}
(\mathcal{H} \psi_{\widetilde{V}}^{f})(\xi) \, \md \xi$ has to be handled 
appropriately: proceeding as per the analysis on p.~347 of \cite{a45}, 
it turns out that, if there are any poles $\alpha_{p_{q}}$, $q \! \in \! 
\lbrace 1,\dotsc,\mathfrak{s} \! - \! 2,\mathfrak{s} \rbrace$, lying 
between $\tilde{a}_{j}$ and $z$, then the integral $\int_{\tilde{a}_{j}}^{z}
(\mathcal{H} \psi_{\widetilde{V}}^{f})(\xi) \, \md \xi$ gives rise to 
integrable logarithmic singularities, and the relation~\eqref{eql3.8g} 
holds true for $z \! \in \! (\tilde{a}_{j},\tilde{b}_{j})$, $j \! = \! 
1,2,\dotsc,N$. Substituting relation~\eqref{eql3.8g} into 
Equation~\eqref{eql3.8f}, one shows that, for $n \! \in \! 
\mathbb{N}$ and $k \! \in \! \lbrace 1,2,\dotsc,K \rbrace$ such that 
$\alpha_{p_{\mathfrak{s}}} \! := \! \alpha_{k} \! \neq \! \infty$, and 
$z \! \in \! (\tilde{a}_{j},\tilde{b}_{j})$, $j \! = \! 1,2,\dotsc,N$,
\begin{align*}
g^{f}_{+}(z) \! + \! g^{f}_{-}(z) \! - \! \hat{\mathscr{P}}_{0}^{-} \! - \! 
\hat{\mathscr{P}}_{0}^{+} \! - \! \widetilde{V}(z) \! - \! \tilde{\ell} 
=& \, 2 \pi \left(\dfrac{(n \! - \! 1)K \! + \! k}{n} \right) 
\int_{\tilde{a}_{j}}^{z}(\mathcal{H} \psi_{\widetilde{V}}^{f})(\xi) \, \md 
\xi \! + \! 2 \left(\dfrac{(n \! - \! 1)K \! + \! k}{n} \right) \int_{J_{f}} 
\ln (\lvert \tilde{a}_{j} \! - \! \xi \rvert) \psi_{\widetilde{V}}^{f}(\xi) 
\, \md \xi \nonumber \\
-& \, 2 \left(\dfrac{(n \! - \! 1)K \! + \! k}{n} \right) \int_{J_{f}} \ln 
(\lvert \xi \! - \! \alpha_{k} \rvert) \psi^{f}_{\widetilde{V}}(\xi) \, 
\md \xi \! - \ 2 \sum_{q=1}^{\mathfrak{s}-2} \dfrac{\varkappa_{nk 
\tilde{k}_{q}}}{n} \ln \left\lvert \dfrac{z \! - \! \alpha_{p_{q}}}{
\alpha_{p_{q}} \! - \! \alpha_{k}} \right\rvert \nonumber \\
-& \, 2 \left(\dfrac{\varkappa_{nk} \! - \! 1}{n} \right) \ln \lvert z 
\! - \! \alpha_{k} \rvert \! - \! \widetilde{V}(z) \! - \! \tilde{\ell} 
\! = \! 2 \pi \left(\dfrac{(n \! - \! 1)K \! + \! k}{n} \right) 
\int_{\tilde{a}_{j}}^{z}(\mathcal{H} \psi_{\widetilde{V}}^{f})(\xi) \, 
\md \xi \nonumber \\
-& \, \int_{\tilde{a}_{j}}^{z} \widetilde{V}^{\prime}(\xi) \, \md \xi \! - \! 
2 \left(\dfrac{\varkappa_{nk} \! - \! 1}{n} \right) \int_{\tilde{a}_{j}}^{z}
(\xi \! - \! \alpha_{k})^{-1} \, \md \xi \! - \! 2 \sum_{q=1}^{\mathfrak{s}-2} 
\dfrac{\varkappa_{nk \tilde{k}_{q}}}{n} \int_{\tilde{a}_{j}}^{z}(\xi \! - \! 
\alpha_{p_{q}})^{-1} \, \md \xi \nonumber \\
+& \, \left(2 \left(\dfrac{(n \! - \! 1)K \! + \! k}{n} \right) \int_{J_{f}} 
\ln \left(\left\lvert \dfrac{\tilde{a}_{j} \! - \! \xi}{\xi \! - \! 
\alpha_{k}} \right\rvert \right) \psi^{f}_{\widetilde{V}}(\xi) \, \md \xi \! 
- \! 2 \sum_{q=1}^{\mathfrak{s}-2} \dfrac{\varkappa_{nk \tilde{k}_{q}}}{n} \ln 
\left\lvert \dfrac{\tilde{a}_{j} \! - \! \alpha_{p_{q}}}{\alpha_{p_{q}} \! - 
\! \alpha_{k}} \right\rvert \right. \nonumber \\
-&\left. \, 2 \left(\dfrac{\varkappa_{nk} \! - \! 1}{n} \right) \ln \lvert 
\tilde{a}_{j} \! - \! \alpha_{k} \rvert \! - \! \widetilde{V}(\tilde{a}_{j}) 
\! - \! \tilde{\ell} \right),
\end{align*}
which, via Equation~\eqref{eql3.8d}, simplifies to
\begin{align} \label{eql3.8h} 
g^{f}_{+}(z) \! + \! g^{f}_{-}(z) \! - \! \hat{\mathscr{P}}_{0}^{-} \! - \! 
\hat{\mathscr{P}}_{0}^{+} \! - \! \widetilde{V}(z) \! - \! \tilde{\ell} =& \, 
\int_{\tilde{a}_{j}}^{z} \left(2 \pi \left(\dfrac{(n \! - \! 1)K \! + \! k}{n} 
\right)(\mathcal{H} \psi_{\widetilde{V}}^{f})(\xi) \! - \! \dfrac{2(
\varkappa_{nk} \! - \! 1)}{n(\xi \! - \! \alpha_{k})} \! - \! 2 \sum_{q=
1}^{\mathfrak{s}-2} \dfrac{\varkappa_{nk \tilde{k}_{q}}}{n(\xi \! - \! 
\alpha_{p_{q}})} \! - \! \widetilde{V}^{\prime}(\xi) \right) \md \xi.
\end{align}
{}From definition~\eqref{eql3.7l}, the second line of 
Equation~\eqref{eql3.7p}, and Equation~\eqref{eql3.7a30}, one shows that, for 
$z \! \in \! (\tilde{a}_{j},\tilde{b}_{j})$, $j \! = \! 1,2,\dotsc,N$,
\begin{equation} \label{eql3.8i} 
2 \pi \left(\dfrac{(n \! - \! 1)K \! + \! k}{n} \right)(\mathcal{H} 
\psi_{\widetilde{V}}^{f})(z) \! = \! \widetilde{V}^{\prime}(z) \! + \! 
\dfrac{2(\varkappa_{nk} \! - \! 1)}{n(z \! - \! \alpha_{k})} \! + \! 2 
\sum_{q=1}^{\mathfrak{s}-2} \dfrac{\varkappa_{nk \tilde{k}_{q}}}{n(z \! 
- \! \alpha_{p_{q}})} \! - \! \left(\dfrac{(n \! - \! 1)K \! + \! k}{n} 
\right)(\tilde{R}(z))^{1/2} \tilde{h}_{\widetilde{V}}(z):
\end{equation}
Substituting Equation~\eqref{eql3.8i} into Equation~\eqref{eql3.8h}, 
one arrives at, for $n \! \in \! \mathbb{N}$ and $k \! \in \! \lbrace 
1,2,\dotsc,K \rbrace$ such that $\alpha_{p_{\mathfrak{s}}} \! := \! 
\alpha_{k} \! \neq \! \infty$, and $z \! \in \! (\tilde{a}_{j},\tilde{b}_{j})$, 
$j \! = \! 1,2,\dotsc,N$,
\begin{equation} \label{eql3.8AA} 
g^{f}_{+}(z) \! + \! g^{f}_{-}(z) \! - \! \hat{\mathscr{P}}_{0}^{-} \! - \! 
\hat{\mathscr{P}}_{0}^{+} \! - \! \widetilde{V}(z) \! - \! \tilde{\ell} \! = 
\! -\left(\dfrac{(n \! - \! 1)K \! + \! k}{n} \right) \int_{\tilde{a}_{j}}^{z}
(\tilde{R}(\xi))^{1/2} \tilde{h}_{\widetilde{V}}(\xi) \, \md \xi \quad 
(\leqslant \! 0).
\end{equation}
(Since, for regular $\widetilde{V} \colon \overline{\mathbb{R}} 
\setminus \lbrace \alpha_{1},\alpha_{2},\dotsc,\alpha_{K} \rbrace \! 
\to \! \mathbb{R}$ satisfying conditions~\eqref{eq20}--\eqref{eq22}, 
$\tilde{h}_{\widetilde{V}} \colon \mathbb{R} \setminus \lbrace 
\alpha_{p_{1}},\dotsc,\alpha_{p_{\mathfrak{s}-2}},\alpha_{p_{\mathfrak{s}}} 
\rbrace \! \to \! \mathbb{R}$ is real analytic (and clearly not identically 
zero), it follows that one has equality only at points $\mathfrak{z}_{\tilde{h}} 
\! \in \! \cup_{j=1}^{N}(\tilde{a}_{j},\tilde{b}_{j})$ for which $\tilde{h}_{
\widetilde{V}}(\mathfrak{z}_{\tilde{h}}) \! = \! 0$, with $\# \lbrace 
\mathstrut \mathfrak{z}_{\tilde{h}} \! \in \! \cup_{j=1}^{N}(\tilde{a}_{j},
\tilde{b}_{j}); \, \tilde{h}_{\widetilde{V}}(\mathfrak{z}_{\tilde{h}}) \! = \! 
0 \rbrace \! < \! \infty$.) For $n \! \in \! \mathbb{N}$ and $k \! \in \! 
\lbrace 1,2,\dotsc,K \rbrace$ such that $\alpha_{p_{\mathfrak{s}}} \! := 
\! \alpha_{k} \! \neq \! \infty$, one shows {}from Equation~\eqref{eql3.8f} 
that, for $\xi \! \in \! J_{f}$ and $z \! \in \! (\tilde{a}_{j},\tilde{b}_{j})$, 
$j \! = \! 1,2,\dotsc,N$, such that $\lvert (z \! - \! \alpha_{k})/(\xi \! - \! 
\alpha_{k}) \rvert \! \ll \! 1$ (e.g., $0 \! < \! \lvert z \! - \! \alpha_{k} 
\rvert \! \ll \! \min \lbrace \min_{i \neq j \in \lbrace 1,\dotsc,
\mathfrak{s}-2,\mathfrak{s} \rbrace} \lbrace \lvert \alpha_{p_{i}} \! - \! 
\alpha_{p_{j}} \rvert \rbrace,\inf_{\xi \in J_{f}} \lbrace \lvert \xi \! - \! 
\alpha_{k} \rvert \rbrace,\min_{i=1,2,\dotsc,N+1} \lbrace \lvert \lvert 
\tilde{b}_{i-1} \! - \! \tilde{a}_{i} \rvert \! - \! \alpha_{k} \rvert \rbrace 
\rbrace)$, via the expansions $\tfrac{1}{(z-\alpha_{k})-(\xi -\alpha_{k})} 
\! = \! -\sum_{j=0}^{l} \tfrac{(z-\alpha_{k})^{j}}{(\xi -\alpha_{k})^{j+1}} 
\! + \! \tfrac{(z-\alpha_{k})^{l+1}}{(\xi -\alpha_{k})^{l+1}(z-\xi)}$, 
$l \! \in \! \mathbb{N}_{0}$, and $\ln (1 \! - \! \blacklozenge) \! = \! 
-\sum_{m=1}^{\infty} \blacklozenge^{m}/m$, $\lvert \blacklozenge 
\rvert \! \ll \! 1$,
\begin{equation*}
g^{f}_{+}(z) \! + \! g^{f}_{-}(z) \! - \! \hat{\mathscr{P}}_{0}^{-} \! - \! 
\hat{\mathscr{P}}_{0}^{+} \! - \! \widetilde{V}(z) \! - \! \tilde{\ell} 
\underset{z \to \alpha_{k}}{=} -\left(\widetilde{V}(z) \! - \! \left(
\dfrac{\varkappa_{nk} \! - \! 1}{n} \right) \ln (\lvert z \! - \! \alpha_{k} 
\rvert^{-2} \! + \! 1) \right) \! + \! \mathcal{O}(1),
\end{equation*}
which, for regular $\widetilde{V} \colon \overline{\mathbb{R}} \setminus 
\lbrace \alpha_{1},\alpha_{2},\dotsc,\alpha_{K} \rbrace \! \to \! \mathbb{R}$ 
satisfying conditions~\eqref{eq20}--\eqref{eq22}, shows that $g^{f}_{+}
(z) \! + \! g^{f}_{-}(z) \! - \! \hat{\mathscr{P}}_{0}^{-} \! - \! 
\hat{\mathscr{P}}_{0}^{+} \! - \! \widetilde{V}(z) \! - \! \tilde{\ell} \! < 
\! 0$, $z \! \in \! (\tilde{a}_{j},\tilde{b}_{j})$, $j \! = \! 1,2,\dotsc,N$.

$\pmb{(3)}$ Substituting the representation~\eqref{eql3.7k} for the density 
of the associated equilibrium measure into the definition of $g^{f}(z)$ 
given by Equation~\eqref{eql3.4gee3}, one shows that, for $z \! \in \! 
(\tilde{a}_{N+1},+\infty)$,
\begin{align} \label{eql3.8j} 
g^{f}_{\pm}(z) =& \, \left(\dfrac{(n \! - \! 1)K \! + \! k}{n} \right) 
\int_{J_{f}} \ln (\lvert z \! - \! \xi \rvert) \psi^{f}_{\widetilde{V}}
(\xi) \, \md \xi \! - \! \left(\dfrac{\varkappa_{nk} \! - \! 1}{n} 
\right) \int_{J_{f}} \ln (\lvert \xi \! - \! \alpha_{k} \rvert) 
\psi^{f}_{\widetilde{V}}(\xi) \, \md \xi \nonumber \\
-& \, \mi \pi \left(\dfrac{\varkappa_{nk} \! - \! 1}{n} \right) \int_{J_{f} 
\cap \mathbb{R}_{\alpha_{k}}^{<}} \psi^{f}_{\widetilde{V}}(\xi) \, \md \xi 
\! - \! \sum_{q=1}^{\mathfrak{s}-2} \dfrac{\varkappa_{nk \tilde{k}_{q}}}{n} 
\int_{J_{f}} \ln (\lvert \xi \! - \! \alpha_{p_{q}} \rvert) 
\psi^{f}_{\widetilde{V}}(\xi) \, \md \xi \nonumber \\
-& \, \mi \pi \sum_{q=1}^{\mathfrak{s}-2} \dfrac{\varkappa_{nk 
\tilde{k}_{q}}}{n} \int_{J_{f} \cap \mathbb{R}_{\alpha_{p_{q}}}^{<}} 
\psi^{f}_{\widetilde{V}}(\xi) \, \md \xi \! - \! \sum_{q=1}^{\mathfrak{s}-2} 
\dfrac{\varkappa_{nk \tilde{k}_{q}}}{n} \ln \lvert z \! - \! \alpha_{p_{q}} 
\rvert \! \mp \! \mi \pi \sum_{q \in \tilde{\Delta}(z)} \dfrac{\varkappa_{nk 
\tilde{k}_{q}}}{n} \nonumber \\
-& \, \left(\dfrac{\varkappa_{nk} \! - \! 1}{n} \right) 
\left(\ln \lvert z \! - \! \alpha_{k} \rvert \! \pm \! \mi \pi 
\chi_{\mathbb{R}^{<}_{\alpha_{k}}}(z) \right),
\end{align}
where $\tilde{\Delta}(z) \! := \! \lbrace \mathstrut j \! \in \! \lbrace 
1,2,\dotsc,\mathfrak{s} \! - \! 2 \rbrace; \, \alpha_{p_{j}} \! > \! z 
\rbrace$;\footnote{If, for $z \! \in \! (\tilde{a}_{N+1},+\infty)$, 
$\tilde{\Delta}(z) \! = \! \varnothing$, then $\sum_{q \in 
\tilde{\Delta}(z)} \varkappa_{nk \tilde{k}_{q}}/n \! := \! 0$.} hence, via 
Equation~\eqref{eql3.8j} and the definition of $\hat{\mathscr{P}}_{
0}^{\pm}$ given by Equation~\eqref{eql3.4gee5}, one shows that, for 
$n \! \in \! \mathbb{N}$ and $k \! \in \! \lbrace 1,2,\dotsc,K \rbrace$ 
such that $\alpha_{p_{\mathfrak{s}}} \! := \!\alpha_{k} \! \neq \! \infty$,
\begin{equation} \label{eql3.8B12} 
(2 \pi \mi)^{-1}(g^{f}_{+}(z) \! - \! g^{f}_{-}(z) \! + \! 
\hat{\mathscr{P}}_{0}^{-} \! - \! \hat{\mathscr{P}}_{0}^{+}) \! = \! -\left(
\dfrac{(n \! - \! 1)K \! + \! k}{n} \right) \int_{J_{f} \cap \mathbb{R}_{
\alpha_{k}}^{>}} \psi^{f}_{\widetilde{V}}(\xi) \, \md \xi \! - \! \sum_{q \in 
\tilde{\Delta}(z)} \dfrac{\varkappa_{nk \tilde{k}_{q}}}{n} \! + \! \sum_{q \in 
\tilde{\Delta}(k)} \dfrac{\varkappa_{nk \tilde{k}_{q}}}{n} \! - \! \left(
\dfrac{\varkappa_{nk} \! - \! 1}{n} \right) \chi_{\mathbb{R}_{
\alpha_{k}}^{<}}(z),
\end{equation}
whence $g^{f}_{+}(z) \! - \! g^{f}_{-}(z) \! + \! \hat{\mathscr{P}}_{0}^{-} 
\! - \! \hat{\mathscr{P}}_{0}^{+} \! \in \! \mi \mathbb{R}$; moreover, 
via Leibnitz's Rule, it follows that, for $z \! \in \! (\tilde{a}_{N+1},
+\infty)$,
\begin{equation} \label{eql3.8D3} 
\mi \left(g^{f}_{+}(z) \! - \! g^{f}_{-}(z) \! + \! \hat{\mathscr{P}}_{0}^{-} 
\! - \! \hat{\mathscr{P}}_{0}^{+} \! + \! 2 \pi \mi \left(
\dfrac{\varkappa_{nk} \! - \! 1}{n} \right) \chi_{\mathbb{R}_{
\alpha_{k}}^{<}}(z) \! + \! 2 \pi \mi \sum_{q \in \tilde{\Delta}(z)} 
\dfrac{\varkappa_{nk \tilde{k}_{q}}}{n} \right)^{\prime} \! = \! 0.
\end{equation}
Via Equation~\eqref{eql3.8j} and the definition of 
$\hat{\mathscr{P}}_{0}^{\pm}$ given by Equation~\eqref{eql3.4gee5}, 
one arrives at (cf. Equation~\eqref{eql3.4h}), for $n \! \in \! 
\mathbb{N}$ and $k \! \in \! \lbrace 1,2,\dotsc,K \rbrace$ such that 
$\alpha_{p_{\mathfrak{s}}} \! := \! \alpha_{k} \! \neq \! \infty$, and 
$z \! \in \! (\tilde{a}_{N+1},+\infty)$,
\begin{align} \label{eql3.8k} 
g^{f}_{+}(z) \! + \! g^{f}_{-}(z) \! - \! \hat{\mathscr{P}}_{0}^{-} \! - \! 
\hat{\mathscr{P}}_{0}^{+} \! - \! \widetilde{V}(z) \! - \! \tilde{\ell} =& \, 
2 \left(\dfrac{(n \! - \! 1)K \! + \! k}{n} \right) \int_{J_{f}} \ln \left(
\left\lvert \dfrac{z \! - \! \xi}{\xi \! - \! \alpha_{k}} \right\rvert \right) 
\psi^{f}_{\widetilde{V}}(\xi) \, \md \xi \! - \! 2 \sum_{q=1}^{\mathfrak{s}-2} 
\dfrac{\varkappa_{nk \tilde{k}_{q}}}{n} \ln \left\lvert \dfrac{z \! - \! 
\alpha_{p_{q}}}{\alpha_{p_{q}} \! - \! \alpha_{k}} \right\rvert \nonumber \\
-& \, 2 \left(\dfrac{\varkappa_{nk} \! - \! 1}{n} \right) \ln \lvert z \! 
- \! \alpha_{k} \rvert \! - \! \widetilde{V}(z) \! - \! \tilde{\ell}.
\end{align}
Proceeding as in case~$\pmb{(2)}$ above, one shows that, for $z \! \in \! 
(\tilde{a}_{N+1},+\infty)$,
\begin{equation} \label{eql3.8l} 
\int_{J_{f}} \ln (\lvert z \! - \! \xi \rvert) \psi_{\widetilde{V}}^{f}
(\xi) \, \md \xi \! = \! \pi \int_{\tilde{a}_{N+1}}^{z}(\mathcal{H} 
\psi_{\widetilde{V}}^{f})(\xi) \, \md \xi \! + \! \int_{J_{f}} \ln (\lvert 
\tilde{a}_{N+1} \! - \! \xi \rvert) \psi_{\widetilde{V}}^{f}(\xi) \, \md \xi:
\end{equation}
substituting relation~\eqref{eql3.8l} into Equation~\eqref{eql3.8k}, one 
shows that, for $n \! \in \! \mathbb{N}$ and $k \! \in \! \lbrace 1,2,
\dotsc,K \rbrace$ such that $\alpha_{p_{\mathfrak{s}}} \! := \! 
\alpha_{k} \! \neq \! \infty$, and $z \! \in \! (\tilde{a}_{N+1},+\infty)$,
\begin{align*}
g^{f}_{+}(z) \! + \! g^{f}_{-}(z) \! - \! \hat{\mathscr{P}}_{0}^{-} \! - \! 
\hat{\mathscr{P}}_{0}^{+} \! - \! \widetilde{V}(z) \! - \! \tilde{\ell} 
=& \, 2 \pi \left(\dfrac{(n \! - \! 1)K \! + \! k}{n} \right) 
\int_{\tilde{a}_{N+1}}^{z}(\mathcal{H} \psi_{\widetilde{V}}^{f})(\xi) \, 
\md \xi \! - \! \int_{\tilde{a}_{N+1}}^{z} \widetilde{V}^{\prime}(\xi) \, 
\md \xi \! - \! 2 \left(\dfrac{\varkappa_{nk} \! - \! 1}{n} \right) 
\int_{\tilde{a}_{N+1}}^{z}(\xi \! - \! \alpha_{k})^{-1} \, \md \xi \\
-& \, 2 \sum_{q=1}^{\mathfrak{s}-2} \dfrac{\varkappa_{nk \tilde{k}_{q}}}{n} 
\int_{\tilde{a}_{N+1}}^{z}(\xi \! - \! \alpha_{p_{q}})^{-1} \, \md \xi \! + 
\! \left(2 \left(\dfrac{(n \! - \! 1)K \! + \! k}{n} \right) \int_{J_{f}} 
\ln \left(\left\lvert \dfrac{\tilde{a}_{N+1} \! - \! \xi}{\xi \! - \! 
\alpha_{k}} \right\rvert \right) \psi^{f}_{\widetilde{V}}(\xi) \, \md \xi 
\right. \\
-&\left. \, 2 \sum_{q=1}^{\mathfrak{s}-2} \dfrac{\varkappa_{nk 
\tilde{k}_{q}}}{n} \ln \left\lvert \dfrac{\tilde{a}_{N+1} \! - \! 
\alpha_{p_{q}}}{\alpha_{p_{q}} \! - \! \alpha_{k}} \right\rvert \! - \! 2 
\left(\dfrac{\varkappa_{nk} \! - \! 1}{n} \right) \ln \lvert \tilde{a}_{N+1} 
\! - \! \alpha_{k} \rvert \! - \! \widetilde{V}(\tilde{a}_{N+1}) \! - \! 
\tilde{\ell} \right),
\end{align*}
which, via Equation~\eqref{eql3.8d}, simplifies to
\begin{align} \label{eql3.8m} 
g^{f}_{+}(z) \! + \! g^{f}_{-}(z) \! - \! \hat{\mathscr{P}}_{0}^{-} \! - \! 
\hat{\mathscr{P}}_{0}^{+} \! - \! \widetilde{V}(z) \! - \! \tilde{\ell} =& \, 
\int_{\tilde{a}_{N+1}}^{z} \left(2 \pi \left(\dfrac{(n \! - \! 1)K \! + \! 
k}{n} \right)(\mathcal{H} \psi_{\widetilde{V}}^{f})(\xi) \! - \! \dfrac{2
(\varkappa_{nk} \! - \! 1)}{n(\xi \! - \! \alpha_{k})} \! - \! 2 \sum_{q=
1}^{\mathfrak{s}-2} \dfrac{\varkappa_{nk \tilde{k}_{q}}}{n(\xi \! - \! 
\alpha_{p_{q}})} \! - \! \widetilde{V}^{\prime}(\xi) \right) \md \xi.
\end{align}
As a simple analysis shows, it turns out that Equation~\eqref{eql3.8i} is 
also valid for $z \! \in \! (\tilde{a}_{N+1},+\infty)$; hence, substituting 
Equation~\eqref{eql3.8i} into Equation~\eqref{eql3.8m}, one arrives 
at, for $n \! \in \! \mathbb{N}$ and $k \! \in \! \lbrace 1,2,\dotsc,
K \rbrace$ such that $\alpha_{p_{\mathfrak{s}}} \! := \! \alpha_{k} 
\! \neq \! \infty$, and $z \! \in \! (\tilde{a}_{N+1},+\infty)$,
\begin{equation} \label{eql3.8BB} 
g^{f}_{+}(z) \! + \! g^{f}_{-}(z) \! - \! \hat{\mathscr{P}}_{0}^{-} \! - 
\! \hat{\mathscr{P}}_{0}^{+} \! - \! \widetilde{V}(z) \! - \! \tilde{\ell} 
\! = \! -\left(\dfrac{(n \! - \! 1)K \! + \! k}{n} \right) 
\int_{\tilde{a}_{N+1}}^{z}(\tilde{R}(\xi))^{1/2} \tilde{h}_{\widetilde{V}}
(\xi) \, \md \xi \quad (\leqslant \! 0).
\end{equation}
(Since, for regular $\widetilde{V} \colon \overline{\mathbb{R}} \setminus 
\lbrace \alpha_{1},\alpha_{2},\dotsc,\alpha_{K} \rbrace \! \to \! \mathbb{R}$ 
satisfying conditions~\eqref{eq20}--\eqref{eq22}, $\tilde{h}_{\widetilde{V}} 
\colon \mathbb{R} \setminus \lbrace \alpha_{p_{1}},\dotsc,\alpha_{p_{
\mathfrak{s}-2}},\alpha_{p_{\mathfrak{s}}} \rbrace \! \to \! \mathbb{R}$ 
is real analytic, it follows that one has equality only at points 
$\mathfrak{z}^{\sharp}_{\tilde{h}} \! \in \! (\tilde{a}_{N+1},+\infty)$ for 
which $\tilde{h}_{\widetilde{V}}(\mathfrak{z}^{\sharp}_{\tilde{h}}) \! = \! 
0$, with $\# \lbrace \mathstrut \mathfrak{z}^{\sharp}_{\tilde{h}} \! \in \! 
(\tilde{a}_{N+1},+\infty); \, \tilde{h}_{\widetilde{V}}(\mathfrak{z}^{
\sharp}_{\tilde{h}}) \! = \! 0 \rbrace \! < \! \infty$.) For $n \! \in \! 
\mathbb{N}$ and $k \! \in \! \lbrace 1,2,\dotsc,K \rbrace$ such that 
$\alpha_{p_{\mathfrak{s}}} \! := \! \alpha_{k} \! \neq \! \infty$, one 
shows {}from Equation~\eqref{eql3.8k} that, for $\xi \! \in \! J_{f}$ 
and $z \! \in \! (\tilde{a}_{N+1},+\infty)$ such that $\lvert (z \! - \! 
\alpha_{k})/(\xi \! - \! \alpha_{k}) \rvert \! \ll \! 1$ (e.g., $0 \! < \! 
\lvert z \! - \! \alpha_{k} \rvert \! \ll \! \min \lbrace \min_{i \neq j \in 
\lbrace 1,\dotsc,\mathfrak{s}-2,\mathfrak{s} \rbrace} \lbrace \lvert 
\alpha_{p_{i}} \! - \! \alpha_{p_{j}} \rvert \rbrace,\inf_{\xi \in J_{f}} 
\lbrace \lvert \xi \! - \! \alpha_{k} \rvert \rbrace,\min_{i=1,2,\dotsc,
N+1} \lbrace \lvert \lvert \tilde{b}_{i-1} \! - \! \tilde{a}_{i} \rvert 
\! - \! \alpha_{k} \rvert \rbrace \rbrace)$, via the expansions 
$\tfrac{1}{(z-\alpha_{k})-(\xi -\alpha_{k})} \! = \! -\sum_{j=0}^{l} 
\tfrac{(z-\alpha_{k})^{j}}{(\xi -\alpha_{k})^{j+1}} \! + \! \tfrac{(z-
\alpha_{k})^{l+1}}{(\xi -\alpha_{k})^{l+1}(z-\xi)}$, $l \! \in \! 
\mathbb{N}_{0}$, and $\ln (1 \! - \! \blacklozenge) \! = \! 
-\sum_{m=1}^{\infty} \blacklozenge^{m}/m$, $\lvert \blacklozenge 
\rvert \! \ll \! 1$,
\begin{equation*}
g^{f}_{+}(z) \! + \! g^{f}_{-}(z) \! - \! \hat{\mathscr{P}}_{0}^{-} \! - \! 
\hat{\mathscr{P}}_{0}^{+} \! - \! \widetilde{V}(z) \! - \! \tilde{\ell} 
\underset{z \to \alpha_{k}}{=} -\left(\widetilde{V}(z) \! - \! \left(
\dfrac{\varkappa_{nk} \! - \! 1}{n} \right) \ln (\lvert z \! - \! \alpha_{k} 
\rvert^{-2} \! + \! 1) \right) \! + \! \mathcal{O}(1),
\end{equation*}
and, for $\xi \! \in \! J_{f}$ and $z \! \in \! (\tilde{a}_{N+1},+\infty)$ 
such that $\lvert \xi/z \rvert \! \ll \! 1$ (e.g., $\lvert z \rvert \! \gg 
\! \max \lbrace \max_{i \neq j \in \lbrace 1,\dotsc,\mathfrak{s}-2,
\mathfrak{s} \rbrace} \lbrace \lvert \alpha_{p_{i}} \! - \! \alpha_{p_{j}} 
\rvert \rbrace,\max_{q=1,\dotsc,\mathfrak{s}-2,\mathfrak{s}} 
\lbrace \lvert \alpha_{p_{q}} \rvert \rbrace,\linebreak[4] 
\max_{i=1,2,\dotsc,N+1} \lbrace \lvert \tilde{b}_{i-1} \! - \! 
\tilde{a}_{i} \rvert \rbrace \rbrace)$, via the expansions 
$\tfrac{1}{\xi -z} \! = \! -\sum_{j=0}^{l} \tfrac{\xi^{j}}{z^{j+1}} \! + 
\! \tfrac{\xi^{l+1}}{z^{l+1}(\xi -z)}$, $l \! \in \! \mathbb{N}_{0}$, 
and $\ln (1 \! - \! \blacklozenge) \! = \! -\sum_{m=1}^{\infty} 
\blacklozenge^{m}/m$, $\lvert \blacklozenge \rvert \! \ll \! 1$,
\begin{equation*}
g^{f}_{+}(z) \! + \! g^{f}_{-}(z) \! - \! \hat{\mathscr{P}}_{0}^{-} \! - \! 
\hat{\mathscr{P}}_{0}^{+} \! - \! \widetilde{V}(z) \! - \! \tilde{\ell} 
\underset{z \to +\alpha_{p_{\mathfrak{s}-1}} = +\infty}{=} 
-\left(\widetilde{V}(z) \! - \! \left(\dfrac{\varkappa^{\infty}_{nk 
\tilde{k}_{\mathfrak{s}-1}} \! + \! 1}{n} \right) \ln (z^{2} \! + \! 1) 
\right) \! + \! \mathcal{O}(1),
\end{equation*}
which, for regular $\widetilde{V} \colon \overline{\mathbb{R}} \setminus 
\lbrace \alpha_{1},\alpha_{2},\dotsc,\alpha_{K} \rbrace \! \to \! \mathbb{R}$ 
satisfying conditions~\eqref{eq20}--\eqref{eq22}, show that $g^{f}_{+}
(z) \! + \! g^{f}_{-}(z) \! - \! \hat{\mathscr{P}}_{0}^{-} \! - \! 
\hat{\mathscr{P}}_{0}^{+} \! - \! \widetilde{V}(z) \! - \! \tilde{\ell} \! 
< \! 0$, $z \! \in \! (\tilde{a}_{N+1},+\infty)$.

$\pmb{(4)}$ Substituting the representation~\eqref{eql3.7k} for the density 
of the associated equilibrium measure into the definition of $g^{f}(z)$ 
given by Equation~\eqref{eql3.4gee3}, one shows that, for $z \! \in \! 
(-\infty,\tilde{b}_{0})$,
\begin{align} \label{eql3.8n} 
g^{f}_{\pm}(z) =& \, \left(\dfrac{(n \! - \! 1)K \! + \! k}{n} \right) 
\int_{J_{f}} \ln (\lvert z \! - \! \xi \rvert) \psi^{f}_{\widetilde{V}}
(\xi) \, \md \xi \! - \! \left(\dfrac{\varkappa_{nk} \! - \! 1}{n} 
\right) \int_{J_{f}} \ln (\lvert \xi \! - \! \alpha_{k} \rvert) 
\psi^{f}_{\widetilde{V}}(\xi) \, \md \xi \nonumber \\
-& \, \mi \pi \left(\dfrac{\varkappa_{nk} \! - \! 1}{n} \right) \int_{J_{f} 
\cap \mathbb{R}_{\alpha_{k}}^{<}} \psi^{f}_{\widetilde{V}}(\xi) \, \md \xi 
\! - \! \sum_{q=1}^{\mathfrak{s}-2} \dfrac{\varkappa_{nk \tilde{k}_{q}}}{n} 
\int_{J_{f}} \ln (\lvert \xi \! - \! \alpha_{p_{q}} \rvert) 
\psi^{f}_{\widetilde{V}}(\xi) \, \md \xi \nonumber \\
-& \, \mi \pi \sum_{q=1}^{\mathfrak{s}-2} \dfrac{\varkappa_{nk 
\tilde{k}_{q}}}{n} \int_{J_{f} \cap \mathbb{R}_{\alpha_{p_{q}}}^{<}} 
\psi^{f}_{\widetilde{V}}(\xi) \, \md \xi \! - \! \sum_{q=1}^{\mathfrak{s}-2} 
\dfrac{\varkappa_{nk \tilde{k}_{q}}}{n} \ln \lvert z \! - \! \alpha_{p_{q}} 
\rvert \! \pm \! \mi \pi \left(\dfrac{(n \! - \! 1)K \! + \! k}{n} \right) 
\nonumber \\
\mp& \, \mi \pi \sum_{q \in \tilde{\Delta}(z)} \dfrac{\varkappa_{nk 
\tilde{k}_{q}}}{n} \! - \! \left(\dfrac{\varkappa_{nk} \! - \! 1}{n} 
\right) \left(\ln \lvert z \! - \! \alpha_{k} \rvert \! \pm \! \mi \pi 
\chi_{\mathbb{R}^{<}_{\alpha_{k}}}(z) \right);
\end{align}
hence, via Equation~\eqref{eql3.8n} and the definition of 
$\hat{\mathscr{P}}_{0}^{\pm}$ given by Equation~\eqref{eql3.4gee5}, 
one shows that, for $n \! \in \! \mathbb{N}$ and $k \! \in \! \lbrace 
1,2,\dotsc,K \rbrace$ such that $\alpha_{p_{\mathfrak{s}}} \! := \! 
\alpha_{k} \! \neq \! \infty$,
\begin{align} \label{eql3.8B13} 
(2 \pi \mi)^{-1}(g^{f}_{+}(z) \! - \! g^{f}_{-}(z) \! + \! 
\hat{\mathscr{P}}_{0}^{-} \! - \! \hat{\mathscr{P}}_{0}^{+}) =& \, 
\left(\dfrac{(n \! - \! 1)K \! + \! k}{n} \right) \! - \! \left(
\dfrac{(n \! - \! 1)K \! + \! k}{n} \right) \int_{J_{f} \cap \mathbb{R}_{
\alpha_{k}}^{>}} \psi^{f}_{\widetilde{V}}(\xi) \, \md \xi \! - \! \sum_{q \in 
\tilde{\Delta}(z)} \dfrac{\varkappa_{nk \tilde{k}_{q}}}{n} \! + \! \sum_{q 
\in \tilde{\Delta}(k)} \dfrac{\varkappa_{nk \tilde{k}_{q}}}{n} \nonumber \\
-& \, \left(\dfrac{\varkappa_{nk} \! - \! 1}{n} \right) 
\chi_{\mathbb{R}_{\alpha_{k}}^{<}}(z),
\end{align}
whence $g^{f}_{+}(z) \! - \! g^{f}_{-}(z) \! + \! \hat{\mathscr{P}}_{0}^{-} 
\! - \! \hat{\mathscr{P}}_{0}^{+} \! \in \! \mi \mathbb{R}$; moreover, 
via Leibnitz's Rule, it follows that, for $z \! \in \! (-\infty,\tilde{b}_{0})$,
\begin{equation} \label{eql3.8D4} 
\mi \left(g^{f}_{+}(z) \! - \! g^{f}_{-}(z) \! + \! \hat{\mathscr{P}}_{0}^{-} 
\! - \! \hat{\mathscr{P}}_{0}^{+} \! + \! 2 \pi \mi \left(
\dfrac{\varkappa_{nk} \! - \! 1}{n} \right) \chi_{\mathbb{R}_{\alpha_{k}}^{<}}
(z) \! + \! 2 \pi \mi \sum_{q \in \tilde{\Delta}(z)} \dfrac{\varkappa_{nk 
\tilde{k}_{q}}}{n} \right)^{\prime} \! = \! 0.
\end{equation}
Via Equation~\eqref{eql3.8n} and the definition of 
$\hat{\mathscr{P}}_{0}^{\pm}$ given by Equation~\eqref{eql3.4gee5}, 
one arrives at, for $n \! \in \! \mathbb{N}$ and $k \! \in \! \lbrace 1,2,
\dotsc,K \rbrace$ such that $\alpha_{p_{\mathfrak{s}}} \! := \! \alpha_{k} 
\! \neq \! \infty$ and $z \! \in \! (-\infty,\tilde{b}_{0})$,
\begin{align} \label{eql3.8o} 
g^{f}_{+}(z) \! + \! g^{f}_{-}(z) \! - \! \hat{\mathscr{P}}_{0}^{-} \! - \! 
\hat{\mathscr{P}}_{0}^{+} \! - \! \widetilde{V}(z) \! - \! \tilde{\ell} =& \, 
2 \left(\dfrac{(n \! - \! 1)K \! + \! k}{n} \right) \int_{J_{f}} \ln \left(
\left\lvert \dfrac{z \! - \! \xi}{\xi \! - \! \alpha_{k}} \right\rvert \right) 
\psi^{f}_{\widetilde{V}}(\xi) \, \md \xi \! - \! 2 \sum_{q=1}^{\mathfrak{s}-2} 
\dfrac{\varkappa_{nk \tilde{k}_{q}}}{n} \ln \left\lvert \dfrac{z \! - \! 
\alpha_{p_{q}}}{\alpha_{p_{q}} \! - \! \alpha_{k}} \right\rvert \nonumber \\
-& \, 2 \left(\dfrac{\varkappa_{nk} \! - \! 1}{n} \right) \ln \lvert z \! 
- \! \alpha_{k} \rvert \! - \! \widetilde{V}(z) \! - \! \tilde{\ell}.
\end{align}
Proceeding as in case~$\pmb{(2)}$ above, one shows that, for $z \! \in \! 
(-\infty,\tilde{b}_{0})$,
\begin{equation} \label{eql3.8p} 
\int_{J_{f}} \ln (\lvert z \! - \! \xi \rvert) \psi_{\widetilde{V}}^{f}
(\xi) \, \md \xi \! = \! -\pi \int_{z}^{\tilde{b}_{0}}(\mathcal{H} 
\psi_{\widetilde{V}}^{f})(\xi) \, \md \xi \! + \! \int_{J_{f}} \ln (\lvert 
\tilde{b}_{0} \! - \! \xi \rvert) \psi_{\widetilde{V}}^{f}(\xi) \, \md \xi.
\end{equation}
Substituting relation~\eqref{eql3.8p} into Equation~\eqref{eql3.8o}, 
and using Equation~\eqref{eql3.8i} (which remains valid for $z \! \in \! 
(-\infty,\tilde{b}_{0}))$, one shows, via Equation~\eqref{eql3.8d}, that, 
for $n \! \in \! \mathbb{N}$ and $k \! \in \! \lbrace 1,2,\dotsc,K 
\rbrace$ such that $\alpha_{p_{\mathfrak{s}}} \! := \! \alpha_{k} 
\! \neq \! \infty$, and $z \! \in \! (-\infty,\tilde{b}_{0})$,
\begin{equation} \label{eql3.8CC} 
g^{f}_{+}(z) \! + \! g^{f}_{-}(z) \! - \! \hat{\mathscr{P}}_{0}^{-} \! - 
\! \hat{\mathscr{P}}_{0}^{+} \! - \! \widetilde{V}(z) \! - \! \tilde{\ell} 
\! = \! \left(\dfrac{(n \! - \! 1)K \! + \! k}{n} \right) 
\int_{z}^{\tilde{b}_{0}}(\tilde{R}(\xi))^{1/2} \tilde{h}_{\widetilde{V}}(\xi) 
\, \md \xi \quad (\leqslant \! 0).
\end{equation}
(Since, for regular $\widetilde{V} \colon \overline{\mathbb{R}} \setminus 
\lbrace \alpha_{1},\alpha_{2},\dotsc,\alpha_{K} \rbrace \! \to \! \mathbb{R}$ 
satisfying conditions~\eqref{eq20}--\eqref{eq22}, $\tilde{h}_{\widetilde{V}} 
\colon \mathbb{R} \setminus \lbrace \alpha_{p_{1}},\dotsc,\alpha_{p_{
\mathfrak{s}-2}},\alpha_{p_{\mathfrak{s}}} \rbrace \! \to \! \mathbb{R}$ 
is real analytic, it follows that one has equality only at points 
$\mathfrak{z}^{\flat}_{\tilde{h}} \! \in \! (-\infty,\tilde{b}_{0})$ for which 
$\tilde{h}_{\widetilde{V}}(\mathfrak{z}^{\flat}_{\tilde{h}}) \! = \! 0$, with 
$\# \lbrace \mathstrut \mathfrak{z}^{\flat}_{\tilde{h}} \! \in \! (-\infty,
\tilde{b}_{0}); \, \tilde{h}_{\widetilde{V}}(\mathfrak{z}^{\flat}_{\tilde{h}}) 
\! = \! 0 \rbrace \! < \! \infty$.) For $n \! \in \! \mathbb{N}$ and $k \! 
\in \! \lbrace 1,2,\dotsc,K \rbrace$ such that $\alpha_{p_{\mathfrak{s}}} 
\! := \! \alpha_{k} \! \neq \! \infty$, one shows {}from 
Equation~\eqref{eql3.8o} that, for $\xi \! \in \! J_{f}$ and $z \! \in \! 
(-\infty,\tilde{b}_{0})$ such that $\lvert (z \! - \! \alpha_{k})/(\xi \! - \! 
\alpha_{k}) \rvert \! \ll \! 1$ (e.g., $0 \! < \! \lvert z \! - \! \alpha_{k} 
\rvert \! \ll \! \min \lbrace \min_{i \neq j \in \lbrace 1,\dotsc,
\mathfrak{s}-2,\mathfrak{s} \rbrace} \lbrace \lvert \alpha_{p_{i}} 
\! - \! \alpha_{p_{j}} \rvert \rbrace,\inf_{\xi \in J_{f}} \lbrace \lvert \xi \! 
- \! \alpha_{k} \rvert \rbrace,\min_{i=1,2,\dotsc,N+1} \lbrace \lvert \lvert 
\tilde{b}_{i-1} \! - \! \tilde{a}_{i} \rvert \! - \! \alpha_{k} \rvert \rbrace 
\rbrace)$, via the expansions $\tfrac{1}{(z-\alpha_{k})-(\xi -\alpha_{k})} 
\! = \! -\sum_{j=0}^{l} \tfrac{(z-\alpha_{k})^{j}}{(\xi -\alpha_{k})^{j+1}} 
\! + \! \tfrac{(z-\alpha_{k})^{l+1}}{(\xi -\alpha_{k})^{l+1}(z-\xi)}$, 
$l \! \in \! \mathbb{N}_{0}$, and $\ln (1 \! - \! \blacklozenge) \! = \! 
-\sum_{m=1}^{\infty} \blacklozenge^{m}/m$, $\lvert \blacklozenge 
\rvert \! \ll \! 1$,
\begin{equation*}
g^{f}_{+}(z) \! + \! g^{f}_{-}(z) \! - \! \hat{\mathscr{P}}_{0}^{-} \! - \! 
\hat{\mathscr{P}}_{0}^{+} \! - \! \widetilde{V}(z) \! - \! \tilde{\ell} 
\underset{z \to \alpha_{k}}{=} -\left(\widetilde{V}(z) \! - \! \left(
\dfrac{\varkappa_{nk} \! - \! 1}{n} \right) \ln (\lvert z \! - \! \alpha_{k} 
\rvert^{-2} \! + \! 1) \right) \! + \! \mathcal{O}(1),
\end{equation*}
and, for $\xi \! \in \! J_{f}$ and $z \! \in \! (-\infty,\tilde{b}_{0})$ such 
that $\lvert \xi/z \rvert \! \ll \! 1$ (e.g., $\lvert z \rvert \! \gg \! \max 
\lbrace \max_{i \neq j \in \lbrace 1,\dotsc,\mathfrak{s}-2,\mathfrak{s} 
\rbrace} \lbrace \lvert \alpha_{p_{i}} \! - \! \alpha_{p_{j}} \rvert 
\rbrace,\max_{q=1,\dotsc,\mathfrak{s}-2,\mathfrak{s}} \lbrace 
\lvert \alpha_{p_{q}} \rvert \rbrace,\linebreak[4]
\max_{i=1,2,\dotsc,N+1} \lbrace 
\lvert \tilde{b}_{i-1} \! - \! \tilde{a}_{i} \rvert \rbrace \rbrace)$, via the 
expansions $\tfrac{1}{\xi -z} \! = \! -\sum_{j=0}^{l} \tfrac{\xi^{j}}{
z^{j+1}} \! + \! \tfrac{\xi^{l+1}}{z^{l+1}(\xi -z)}$, $l \! \in \! 
\mathbb{N}_{0}$, and $\ln (1 \! - \! \blacklozenge) \! = \! 
-\sum_{m=1}^{\infty} \blacklozenge^{m}/m$, $\lvert \blacklozenge 
\rvert \! \ll \! 1$,
\begin{equation*}
g^{f}_{+}(z) \! + \! g^{f}_{-}(z) \! - \! \hat{\mathscr{P}}_{0}^{-} \! - \! 
\hat{\mathscr{P}}_{0}^{+} \! - \! \widetilde{V}(z) \! - \! \tilde{\ell} 
\underset{z \to -\alpha_{p_{\mathfrak{s}-1}} = -\infty}{=} -\left(
\widetilde{V}(z) \! - \! \left(\dfrac{\varkappa^{\infty}_{nk \tilde{k}_{
\mathfrak{s}-1}} \! + \! 1}{n} \right) \ln (z^{2} \! + \! 1) \right) \! 
+ \! \mathcal{O}(1),
\end{equation*}
which, for regular $\widetilde{V} \colon \overline{\mathbb{R}} 
\setminus \lbrace \alpha_{1},\alpha_{2},\dotsc,\alpha_{K} \rbrace \! 
\to \! \mathbb{R}$ satisfying conditions~\eqref{eq20}--\eqref{eq22}, 
show that $g^{f}_{+}(z) \! + \! g^{f}_{-}(z) \! - \! \hat{\mathscr{P}}_{0}^{-} 
\! - \! \hat{\mathscr{P}}_{0}^{+} \! - \! \widetilde{V}(z) \! - \! \tilde{\ell} 
\! < \! 0$, $z \! \in \! (-\infty,\tilde{b}_{0})$.

$\pmb{(\mathrm{B})}$ The proof of this case, that is, $n \! \in \! \mathbb{N}$ 
and $k \! \in \! \lbrace 1,2,\dotsc,K \rbrace$ such that $\alpha_{p_{\mathfrak{s}}} 
\! := \! \alpha_{k} \! = \! \infty$, is virtually identical to the proof presented in 
$\pmb{(\mathrm{A})}$ above; one mimics, \emph{verbatim}, the scheme of the 
calculations presented in case $\pmb{(\mathrm{A})}$ in order to arrive at the 
corresponding claims stated in item~$\pmb{(1)}$ of the lemma; in order to do 
so, however, the analogues of Equations~\eqref{eql3.8c}, \eqref{eql3.8B10}, 
\eqref{eql3.8d}--\eqref{eql3.8B11}, \eqref{eql3.8f}, 
\eqref{eql3.8AA}--\eqref{eql3.8B12}, \eqref{eql3.8k}, 
\eqref{eql3.8BB}--\eqref{eql3.8B13}, \eqref{eql3.8o}, and~\eqref{eql3.8CC}, 
respectively, are necessary, which, in the present case, and in conjunction with 
the definition of $g^{\infty}(z)$ given by Equation~\eqref{eql3.4gee1}, read:

$\pmb{(1)}$ for $z \! \in \! [\hat{b}_{j-1},\hat{a}_{j}]$, $j \! = \! 1,2,
\dotsc,N \! + \! 1$,
\begin{align} \label{eql3.8B1} 
g^{\infty}_{\pm}(z) =& \, \left(\dfrac{(n \! - \! 1)K \! + \! k}{n} 
\right) \int_{J_{\infty}} \ln (\lvert z \! - \! \xi \rvert) 
\psi^{\infty}_{\widetilde{V}}(\xi) \, \md \xi \! - \! 
\sum_{q=1}^{\mathfrak{s}-1} \dfrac{\varkappa_{nk \tilde{k}_{q}}}{n} 
\int_{J_{\infty}} \ln (\lvert \xi \! - \! \alpha_{p_{q}} \rvert) 
\psi^{\infty}_{\widetilde{V}}(\xi) \, \md \xi \nonumber \\
-& \, \mi \pi \sum_{q=1}^{\mathfrak{s}-1} \dfrac{\varkappa_{nk \tilde{k}_{
q}}}{n} \int_{J_{\infty} \cap \mathbb{R}_{\alpha_{p_{q}}}^{<}} \psi^{\infty}_{
\widetilde{V}}(\xi) \, \md \xi \! - \! \sum_{q=1}^{\mathfrak{s}-1} \dfrac{
\varkappa_{nk \tilde{k}_{q}}}{n} \ln \lvert z \! - \! \alpha_{p_{q}} \rvert 
\! \mp \! \mi \pi \sum_{q \in \hat{\Delta}(j;z)} \dfrac{\varkappa_{nk 
\tilde{k}_{q}}}{n} \nonumber \\
\pm& \, \mi \pi \left(\dfrac{(n \! - \! 1)K \! + \! k}{n} \right) \int_{z}^{
\hat{a}_{N+1}} \psi^{\infty}_{\widetilde{V}}(\xi) \, \md \xi,
\end{align}
that is,
\begin{equation} \label{eql3.8B14} 
(2 \pi \mi)^{-1}(g^{\infty}_{+}(z) \! - \! g^{\infty}_{-}(z)) \! = \! 
\left(\dfrac{(n \! - \! 1)K \! + \! k}{n} \right) \int_{z}^{\hat{a}_{N+1}} 
\psi^{\infty}_{\widetilde{V}}(\xi) \, \md \xi \! - \! \sum_{q \in 
\hat{\Delta}(j;z)} \dfrac{\varkappa_{nk \tilde{k}_{q}}}{n},
\end{equation}
where {}\footnote{If $J_{\infty} \cap \mathbb{R}_{\alpha_{p_{q}}}^{<} \! 
= \! \varnothing$, $q \! \in \! \lbrace 1,2,\dotsc,\mathfrak{s} \! - \! 1 
\rbrace$, then $\int_{J_{\infty} \cap \mathbb{R}_{\alpha_{p_{q}}}^{<}} 
\psi^{\infty}_{\widetilde{V}}(\xi) \, \md \xi \! := \! 0$.} $g^{\infty}_{\pm}
(z) \! := \! \lim_{\varepsilon \downarrow 0}g^{\infty}(z \! \pm \! \mi 
\varepsilon)$, and $\hat{\Delta}(j;z) \! := \! \lbrace \mathstrut i \! \in 
\! \lbrace 1,2,\dotsc,\mathfrak{s} \! - \! 1 \rbrace; \, \alpha_{p_{i}} 
\! > \! z \rbrace$, $j \! = \! 1,2,\dotsc,N \! + \! 1$,\footnote{If 
$\hat{\Delta}(j;z) \! = \! \varnothing$, $j \! \in \! \lbrace 1,2,
\dotsc,N \! + \! 1 \rbrace$, then $\sum_{q \in \hat{\Delta}(j;z)} 
\varkappa_{nk \tilde{k}_{q}}/n \! := \! 0$.}
\begin{equation} \label{eql3.8B2} 
g^{\infty}_{+}(z) \! + \! g^{\infty}_{-}(z) \! - \! 2 \tilde{\mathscr{P}}_{0} 
\! - \! \widetilde{V}(z) \! - \! \hat{\ell} \! = \! 2 \left(\dfrac{(n \! - \! 
1)K \! + \! k}{n} \right) \int_{J_{\infty}} \ln (\lvert z \! - \! \xi \rvert) 
\psi^{\infty}_{\widetilde{V}}(\xi) \, \md \xi \! - \! 2 \sum_{q=1}^{
\mathfrak{s}-1} \dfrac{\varkappa_{nk \tilde{k}_{q}}}{n} \ln \lvert z \! - \! 
\alpha_{p_{q}} \rvert \! - \! \widetilde{V}(z) \! - \! \hat{\ell} \! = \! 0,
\end{equation}
where $\tilde{\mathscr{P}}_{0}$ is defined by Equation~\eqref{eql3.4gee2}, 
and
\begin{align} \label{eql3.8B3} 
\hat{\ell} \! = \! 2 \left(\dfrac{(n \! - \! 1)K \! + \! k}{n} \right) 
\int_{J_{\infty}} \ln \left(\left\lvert \tfrac{1}{2}(\hat{b}_{N} \! + \! 
\hat{a}_{N+1}) \! - \! \xi \right\rvert \right) \psi^{\infty}_{\widetilde{V}}
(\xi) \, \md \xi \! - \! 2 \sum_{q=1}^{\mathfrak{s}-1} \dfrac{\varkappa_{nk 
\tilde{k}_{q}}}{n} \ln \left\lvert \tfrac{1}{2}(\hat{b}_{N} \! + \! 
\hat{a}_{N+1}) \! - \! \alpha_{p_{q}} \right\rvert \! - \! \widetilde{V}
(\tfrac{1}{2}(\hat{b}_{N} \! + \! \hat{a}_{N+1}));
\end{align}

$\pmb{(2)}$ for $z \! \in \! (\hat{a}_{j},\hat{b}_{j})$, $j \! = \! 1,2,
\dotsc,N$,
\begin{align} \label{eql3.8B4} 
g^{\infty}_{\pm}(z) =& \, \left(\dfrac{(n \! - \! 1)K \! + \! k}{n} \right) 
\int_{J_{\infty}} \ln (\lvert z \! - \! \xi \rvert) \psi^{\infty}_{
\widetilde{V}}(\xi) \, \md \xi \! - \! \sum_{q=1}^{\mathfrak{s}-1} 
\dfrac{\varkappa_{nk \tilde{k}_{q}}}{n} \int_{J_{\infty}} \ln (\lvert \xi \! 
- \! \alpha_{p_{q}} \rvert) \psi^{\infty}_{\widetilde{V}}(\xi) \, \md \xi 
\nonumber \\
-& \, \mi \pi \sum_{q=1}^{\mathfrak{s}-1} \dfrac{\varkappa_{nk 
\tilde{k}_{q}}}{n} \int_{J_{\infty} \cap \mathbb{R}_{\alpha_{p_{q}}}^{<}} 
\psi^{\infty}_{\widetilde{V}}(\xi) \, \md \xi \! - \! \sum_{q=1}^{
\mathfrak{s}-1} \dfrac{\varkappa_{nk \tilde{k}_{q}}}{n} \ln \lvert z \! - \! 
\alpha_{p_{q}} \rvert \! \mp \! \mi \pi \sum_{q \in \hat{\Delta}(j;z)} 
\dfrac{\varkappa_{nk \tilde{k}_{q}}}{n} \nonumber \\
\pm& \, \mi \pi \left(\dfrac{(n \! - \! 1)K \! + \! k}{n} \right) 
\int_{\hat{b}_{j}}^{\hat{a}_{N+1}} \psi^{\infty}_{\widetilde{V}}(\xi) 
\, \md \xi,
\end{align}
that is,
\begin{equation} \label{eql3.8B15} 
(2 \pi \mi)^{-1}(g^{\infty}_{+}(z) \! - \! g^{\infty}_{-}(z)) \! = \! 
\left(\dfrac{(n \! - \! 1)K \! + \! k}{n} \right) \int_{\hat{b}_{j}}^{
\hat{a}_{N+1}} \psi^{\infty}_{\widetilde{V}}(\xi) \, \md \xi \! - \! 
\sum_{q \in \hat{\Delta}(j;z)} \dfrac{\varkappa_{nk \tilde{k}_{q}}}{n},
\end{equation}
and
\begin{align} \label{eql3.8B5} 
g^{\infty}_{+}(z) \! + \! g^{\infty}_{-}(z) \! - \! 2 \tilde{\mathscr{P}}_{0} 
\! - \! \widetilde{V}(z) \! - \! \hat{\ell} =& \, 2 \left(\dfrac{(n \! - \! 
1)K \! + \! k}{n} \right) \int_{J_{\infty}} \ln (\lvert z \! - \! \xi \rvert) 
\psi^{\infty}_{\widetilde{V}}(\xi) \, \md \xi \! - \! 2 \sum_{q=1}^{
\mathfrak{s}-1} \dfrac{\varkappa_{nk \tilde{k}_{q}}}{n} \ln \lvert z \! - \! 
\alpha_{p_{q}} \rvert \! - \! \widetilde{V}(z) \! - \! \hat{\ell} \nonumber \\
=& \, -\left(\dfrac{(n \! - \! 1)K \! + \! k}{n} \right) \int_{\hat{a}_{j}}^{z}
(\hat{R}(\xi))^{1/2} \hat{h}_{\widetilde{V}}(\xi) \, \md \xi \quad (\leqslant 
\! 0);
\end{align}

$\pmb{(3)}$ for $z \! \in \! (\hat{a}_{N+1},+\infty)$,
\begin{align} \label{eql3.8B6} 
g^{\infty}_{\pm}(z) =& \, \left(\dfrac{(n \! - \! 1)K \! + \! k}{n} \right) 
\int_{J_{\infty}} \ln (\lvert z \! - \! \xi \rvert) \psi^{\infty}_{
\widetilde{V}}(\xi) \, \md \xi \! - \! \sum_{q=1}^{\mathfrak{s}-1} \dfrac{
\varkappa_{nk \tilde{k}_{q}}}{n} \int_{J_{\infty}} \ln (\lvert \xi \! - \! 
\alpha_{p_{q}} \rvert) \psi^{\infty}_{\widetilde{V}}(\xi) \, \md \xi 
\nonumber \\
-& \, \mi \pi \sum_{q=1}^{\mathfrak{s}-1} \dfrac{\varkappa_{nk 
\tilde{k}_{q}}}{n} \int_{J_{\infty} \cap \mathbb{R}_{\alpha_{p_{q}}}^{<}} 
\psi^{\infty}_{\widetilde{V}}(\xi) \, \md \xi \! - \! \sum_{q=1}^{
\mathfrak{s}-1} \dfrac{\varkappa_{nk \tilde{k}_{q}}}{n} \ln \lvert z \! - 
\! \alpha_{p_{q}} \rvert \! \mp \! \mi \pi \sum_{q \in \hat{\Delta}(z)} 
\dfrac{\varkappa_{nk \tilde{k}_{q}}}{n},
\end{align}
that is,
\begin{equation} \label{eql3.8B16} 
(2 \pi \mi)^{-1}(g^{\infty}_{+}(z) \! - \! g^{\infty}_{-}(z)) \! = \! 
-\sum_{q \in \hat{\Delta}(z)} \dfrac{\varkappa_{nk \tilde{k}_{q}}}{n},
\end{equation}
where $\hat{\Delta}(z) \! := \! \lbrace \mathstrut j \! \in \! \lbrace 
1,2,\dotsc,\mathfrak{s} \! - \! 1 \rbrace; \, \alpha_{p_{j}} \! > \! 
z \rbrace$,\footnote{If, for $z \! \in \! (\hat{a}_{N+1},+\infty)$, 
$\hat{\Delta}(z) \! = \! \varnothing$, then $\sum_{q \in \hat{\Delta}(z)} 
\varkappa_{nk \tilde{k}_{q}}/n \! := \! 0$.} and
\begin{align} \label{eql3.8B7} 
g^{\infty}_{+}(z) \! + \! g^{\infty}_{-}(z) \! - \! 2 \tilde{\mathscr{P}}_{0} 
\! - \! \widetilde{V}(z) \! - \! \hat{\ell} =& \, 2 \left(\dfrac{(n \! - \! 
1)K \! + \! k}{n} \right) \int_{J_{\infty}} \ln (\lvert z \! - \! \xi \rvert) 
\psi^{\infty}_{\widetilde{V}}(\xi) \, \md \xi \! - \! 2 \sum_{q=1}^{
\mathfrak{s}-1} \dfrac{\varkappa_{nk \tilde{k}_{q}}}{n} \ln \lvert z \! - \! 
\alpha_{p_{q}} \rvert \! - \! \widetilde{V}(z) \! - \! \hat{\ell} \nonumber \\
=& \, -\left(\dfrac{(n \! - \! 1)K \! + \! k}{n} \right) \int_{\hat{a}_{N+
1}}^{z}(\hat{R}(\xi))^{1/2} \hat{h}_{\widetilde{V}}(\xi) \, \md \xi \quad 
(\leqslant \! 0);
\end{align}

$\pmb{(4)}$ for $z \! \in \! (-\infty,\hat{b}_{0})$,
\begin{align} \label{eql3.8B8} 
g^{\infty}_{\pm}(z) =& \, \left(\dfrac{(n \! - \! 1)K \! + \! k}{n} \right) 
\int_{J_{\infty}} \ln (\lvert z \! - \! \xi \rvert) \psi^{\infty}_{
\widetilde{V}}(\xi) \, \md \xi \! - \! \sum_{q=1}^{\mathfrak{s}-1} 
\dfrac{\varkappa_{nk \tilde{k}_{q}}}{n} \int_{J_{\infty}} \ln (\lvert 
\xi \! - \! \alpha_{p_{q}} \rvert) \psi^{\infty}_{\widetilde{V}}(\xi) 
\, \md \xi \nonumber \\
-& \, \mi \pi \sum_{q=1}^{\mathfrak{s}-1} \dfrac{\varkappa_{nk 
\tilde{k}_{q}}}{n} \int_{J_{\infty} \cap \mathbb{R}_{\alpha_{p_{q}}}^{<}} 
\psi^{\infty}_{\widetilde{V}}(\xi) \, \md \xi \! - \! \sum_{q=1}^{
\mathfrak{s}-1} \dfrac{\varkappa_{nk \tilde{k}_{q}}}{n} \ln \lvert z 
\! - \! \alpha_{p_{q}} \rvert \! \pm \! \mi \pi \left(\dfrac{(n \! - \! 1)K 
\! + \! k}{n} \right) \nonumber \\
\mp& \, \mi \pi \sum_{q \in \hat{\Delta}(z)} \dfrac{\varkappa_{nk 
\tilde{k}_{q}}}{n},
\end{align}
that is,
\begin{equation} \label{eql3.8B17} 
(2 \pi \mi)^{-1}(g^{\infty}_{+}(z) \! - \! g^{\infty}_{-}(z)) \! = \! \left(
\dfrac{(n \! - \! 1)K \! + \! k}{n} \right) \! - \! \sum_{q \in \hat{\Delta}
(z)} \dfrac{\varkappa_{nk \tilde{k}_{q}}}{n},
\end{equation}
and
\begin{align} \label{eql3.8B9} 
g^{\infty}_{+}(z) \! + \! g^{\infty}_{-}(z) \! - \! 2 \tilde{\mathscr{P}}_{0} 
\! - \! \widetilde{V}(z) \! - \! \hat{\ell} =& \, 2 \left(\dfrac{(n \! - \! 
1)K \! + \! k}{n} \right) \int_{J_{\infty}} \ln (\lvert z \! - \! \xi \rvert) 
\psi^{\infty}_{\widetilde{V}}(\xi) \, \md \xi \! - \! 2 \sum_{q=1}^{
\mathfrak{s}-1} \dfrac{\varkappa_{nk \tilde{k}_{q}}}{n} \ln \lvert z \! - \! 
\alpha_{p_{q}} \rvert \! - \! \widetilde{V}(z) \! - \! \hat{\ell} \nonumber \\
=& \, \left(\dfrac{(n \! - \! 1)K \! + \! k}{n} \right) \int_{z}^{\hat{b}_{0}}
(\hat{R}(\xi))^{1/2} \hat{h}_{\widetilde{V}}(\xi) \, \md \xi \quad (\leqslant 
\! 0).
\end{align}
This concludes the proof. \hfill $\qed$
\section{The Monic MPC ORF Families of Model RHPs and Parametrices} 
\label{sec4}
In this section, the $K$ families of auxiliary matrix RHPs for the 
associated monic MPC ORFs formulated in Lemma~\ref{lem3.4} 
are augmented via contour deformations and transformations into 
simpler, `model' matrix RHPs which, for regular $\widetilde{V} 
\colon \overline{\mathbb{R}} \setminus \lbrace \alpha_{1},
\alpha_{2},\dotsc,\alpha_{K} \rbrace \! \to \! \mathbb{R}$ satisfying 
conditions~\eqref{eq20}--\eqref{eq22}, can be solved explicitly, in 
the double-scaling limit $\mathscr{N},n \! \to \! \infty$ such that 
$z_{o} \! = \! 1 \! + \! o(1)$, in terms of Riemann theta functions 
(associated with the underlying genus-$N$ hyperelliptic (compact) 
Riemann surfaces) and Airy functions (see, for example, 
\cite{a54,a60,a53,a67,a59}).
\begin{ccccc} \label{lem4.1} 
Let the external field $\widetilde{V} \colon \overline{\mathbb{R}} 
\setminus \lbrace \alpha_{1},\alpha_{2},\dotsc,\alpha_{K} \rbrace 
\! \to \! \mathbb{R}$ satisfy conditions~\eqref{eq20}--\eqref{eq22} 
and be regular. For $n \! \in \! \mathbb{N}$ and $k \! \in \! \lbrace 
1,2,\dotsc,K \rbrace$ such that $\alpha_{p_{\mathfrak{s}}} \! := \! 
\alpha_{k} \! = \! \infty$ (resp., $\alpha_{p_{\mathfrak{s}}} \! := \! 
\alpha_{k} \! \neq \! \infty)$, let the associated equilibrium measure, 
$\mu_{\widetilde{V}}^{\infty}$ (resp., $\mu_{\widetilde{V}}^{f})$, 
and its support, $J_{\infty}$ (resp., $J_{f})$, be as described in 
item~$\pmb{(1)}$ (resp., item~$\pmb{(2)})$ of Lemma~\ref{lem3.7}, 
and, along with the corresponding variational constant, $\hat{\ell}$ 
(resp., $\tilde{\ell})$, satisfy the variational conditions~\eqref{eql3.8a} 
(resp., conditions~\eqref{eql3.8b}$)$$;$ moreover, let the associated 
conditions~{\rm (i)}--{\rm (iv)} of item~$\pmb{(1)}$ (resp., 
item~$\pmb{(2)})$ of Lemma~\ref{lem3.8} be valid. For $n \! \in 
\! \mathbb{N}$ and $k \! \in \! \lbrace 1,2,\dotsc,K \rbrace$, let 
$\mathscr{X} \colon \mathbb{N} \times \lbrace 1,2,\dotsc,K \rbrace 
\times \overline{\mathbb{C}} \setminus \overline{\mathbb{R}} \! \to 
\! \operatorname{SL}_{2}(\mathbb{C})$ solve the monic {\rm MPC} 
{\rm ORF} {\rm RHP} $(\mathscr{X}(z),\mathscr{V}(z),\overline{\mathbb{R}})$ 
formulated in Lemma~\ref{lem3.4}, and, for $n \! \in \! \mathbb{N}$ and 
$k \! \in \! \lbrace 1,2,\dotsc,K \rbrace$ such that $\alpha_{p_{\mathfrak{s}}} 
\! := \! \alpha_{k} \! = \! \infty$, set
\begin{equation*}
\mathcal{M}(z) \! = \! \mathscr{X}(z),
\end{equation*}
and, for $n \! \in \! \mathbb{N}$ and $k \! \in \! \lbrace 1,2,\dotsc,K \rbrace$ 
such that $\alpha_{p_{\mathfrak{s}}} \! := \!  \alpha_{k} \! \neq \! \infty$, set
\begin{equation*}
\mathcal{M}(z) \! = \! 
\begin{cases}
\mathscr{X}(z) \mathscr{E}^{-\sigma_{3}}, &\text{$z \! \in \! 
\mathbb{C}_{+}$,} \\
\mathscr{X}(z) \mathscr{E}^{\sigma_{3}}, &\text{$z \! \in \! 
\mathbb{C}_{-}$,}
\end{cases}
\end{equation*}
where
\begin{equation} \label{eqmainfin13} 
\mathscr{E} \colon \mathbb{N} \times \lbrace 1,2,\dotsc,K \rbrace 
\times \mathscr{M}_{1}(\mathbb{R}) \! \ni \! (n,k,\mu_{\widetilde{V}}^{f}) 
\! \mapsto \! \exp \left(\mi \pi ((n \! - \! 1)K \! + \! k) \int_{J_{f} \cap 
\mathbb{R}_{\alpha_{k}}^{>}} \psi_{\widetilde{V}}^{f}(\xi) \, \md \xi \right) 
\! = \! \mathscr{E}[n,k,\mu_{\widetilde{V}}^{f}] \! =: \! \mathscr{E}.
\end{equation}
Then, for $n \! \in \! \mathbb{N}$ and $k \! \in \! \lbrace 1,2,\dotsc,K 
\rbrace$, $\mathcal{M} \colon \mathbb{N} \times \lbrace 1,2,\dotsc,K 
\rbrace \times \overline{\mathbb{C}} \setminus \overline{\mathbb{R}} \! \to 
\! \operatorname{SL}_{2}(\mathbb{C})$ solves the following matrix {\rm RHP:} 
$\pmb{{\rm (i)}}$ $\mathcal{M}(n,k,z) \! = \! \mathcal{M}(z)$ is analytic 
for $z \! \in \! \overline{\mathbb{C}} \setminus \overline{\mathbb{R}};$ 
$\pmb{{\rm (ii)}}$ the boundary values $\mathcal{M}_{\pm}(z) \! := \! 
\lim_{\varepsilon \downarrow 0} \mathcal{M}(z \! \pm \! \mi \varepsilon)$ 
satisfy the jump condition
\begin{equation*}
\mathcal{M}_{+}(z) \! = \! \mathcal{M}_{-}(z) \hat{\mathscr{V}}(z) \quad 
\mathrm{a.e.} \quad z \! \in \! \overline{\mathbb{R}},
\end{equation*}
where, for $n \! \in \! \mathbb{N}$ and $k \! \in \! \lbrace 1,2,\dotsc,K \rbrace$ 
such that $\alpha_{p_{\mathfrak{s}}} \! := \! \alpha_{k} \! = \! \infty$,
\begin{align*}
\hat{\mathscr{V}} &\colon \mathbb{N} \times \lbrace 1,2,\dotsc,K \rbrace 
\times \overline{\mathbb{R}} \! \to \! \operatorname{GL}_{2}(\mathbb{C}), \, 
(n,k,z) \! \mapsto \! \hat{\mathscr{V}}(n,k,z) \! =: \! \hat{\mathscr{V}}(z) 
\\
&= 
\begin{cases}
\begin{pmatrix}
\me^{-2 \pi \mi ((n-1)K+k) \int_{z}^{\hat{a}_{N+1}} 
\psi_{\widetilde{V}}^{\infty}(\xi) \, \md \xi} & 1 \\
0 & \me^{2 \pi \mi ((n-1)K+k) \int_{z}^{\hat{a}_{N+1}} 
\psi_{\widetilde{V}}^{\infty}(\xi) \, \md \xi}
\end{pmatrix}, &\text{$z \! \in \! (\hat{b}_{j-1},\hat{a}_{j}), \quad 
j \! = \! 1,2,\dotsc,N \! + \! 1$,} \\
\begin{pmatrix}
\me^{-2 \pi \mi ((n-1)K+k) \int_{\hat{b}_{i}}^{\hat{a}_{N+1}} 
\psi_{\widetilde{V}}^{\infty}(\xi) \, \md \xi} & \me^{n(g_{+}^{\infty}(z)
+g_{-}^{\infty}(z)-2 \tilde{\mathscr{P}}_{0}-\widetilde{V}(z)-\hat{\ell})} \\
0 & \me^{2 \pi \mi ((n-1)K+k) \int_{\hat{b}_{i}}^{\hat{a}_{N+1}} 
\psi_{\widetilde{V}}^{\infty}(\xi) \, \md \xi}
\end{pmatrix}, &\text{$z \! \in \! (\hat{a}_{i},\hat{b}_{i}), \quad i 
\! = \! 1,2,\dotsc,N$,} \\
\mathrm{I} \! + \! \me^{n(g_{+}^{\infty}(z)+g_{-}^{\infty}(z)-2 
\tilde{\mathscr{P}}_{0}-\widetilde{V}(z)-\hat{\ell})} \sigma_{+}, 
&\text{$z \! \in \! (-\infty,\hat{b}_{0}) \cup (\hat{a}_{N+1},+\infty)$,}
\end{cases}
\end{align*}
with $g^{\infty}(z)$ and $\tilde{\mathscr{P}}_{0}$ defined by 
Equations~\eqref{eql3.4gee1} and~\eqref{eql3.4gee2}, respectively, 
$g^{\infty}_{\pm}(z) \! := \! \lim_{\varepsilon \downarrow 0}g^{\infty}
(z \! \pm \! \mi \varepsilon)$, $g_{+}^{\infty}(z) \! + \! g_{-}^{\infty}
(z) \! - \! 2 \tilde{\mathscr{P}}_{0} \! - \! \widetilde{V}(z) \! - \! \hat{\ell} 
\! < \! 0$, $z \! \in \! (-\infty,\hat{b}_{0}) \cup \cup_{i=1}^{N}(\hat{a}_{i},
\hat{b}_{i}) \cup (\hat{a}_{N+1},+\infty)$, and $\pm \Re (\mi \int_{z}^{
\hat{a}_{N+1}} \psi_{\widetilde{V}}^{\infty}(\xi) \, \md \xi) \! > \! 0$, 
$z \! \in \! \mathbb{C}_{\pm} \cap (\cup_{j=1}^{N+1} \hat{\mathbb{U}}_{j})$, 
where, for $j \! = \! 1,2,\dotsc,N \! + \! 1$, $\hat{\mathbb{U}}_{j} \! := \! 
\lbrace \mathstrut z \! \in \! \mathbb{C} \setminus \lbrace \alpha_{p_{1}},
\alpha_{p_{2}},\dotsc,\alpha_{p_{\mathfrak{s}-1}} \rbrace; \, \inf_{\tau \in 
(\hat{b}_{j-1},\hat{a}_{j})} \lvert z \! - \! \tau \rvert \! < \! \hat{r}_{j} \rbrace$, 
with $\hat{r}_{j} \! \in \! (0,1)$ chosen small enough so that $\hat{\mathbb{U}}_{i} 
\cap \hat{\mathbb{U}}_{j} \! = \! \varnothing$, $i \! \neq \! j \! \in \! \lbrace 
1,2,\dotsc,N \! + \! 1 \rbrace$, and, for $n \! \in \! \mathbb{N}$ and $k \! \in 
\! \lbrace 1,2,\dotsc,K \rbrace$ such that $\alpha_{p_{\mathfrak{s}}} \! := \! 
\alpha_{k} \! \neq \! \infty$,
\begin{align*}
\hat{\mathscr{V}} &\colon \mathbb{N} \times \lbrace 1,2,\dotsc,K \rbrace 
\times \overline{\mathbb{R}} \! \to \! \operatorname{GL}_{2}(\mathbb{C}), \, 
(n,k,z) \! \mapsto \! \hat{\mathscr{V}}(n,k,z) \! =: \! \hat{\mathscr{V}}(z) 
\\
&= 
\begin{cases}
\begin{pmatrix}
\me^{-2 \pi \mi ((n-1)K+k) \int_{z}^{\tilde{a}_{N+1}} 
\psi_{\widetilde{V}}^{f}(\xi) \, \md \xi} & 1 \\
0 & \me^{2 \pi \mi ((n-1)K+k) \int_{z}^{\tilde{a}_{N+1}} 
\psi_{\widetilde{V}}^{f}(\xi) \, \md \xi}
\end{pmatrix}, &\text{$z \! \in \! (\tilde{b}_{j-1},\tilde{a}_{j}), \quad 
j \! = \! 1,2,\dotsc,N \! + \! 1$,} \\
\begin{pmatrix}
\me^{-2 \pi \mi ((n-1)K+k) \int_{\tilde{b}_{i}}^{\tilde{a}_{N+1}} 
\psi_{\widetilde{V}}^{f}(\xi) \, \md \xi} & \me^{n(g_{+}^{f}(z)+g_{-}^{f}
(z)-\hat{\mathscr{P}}_{0}^{+}-\hat{\mathscr{P}}_{0}^{-}-\widetilde{V}(z)-
\tilde{\ell})} \\
0 & \me^{2 \pi \mi ((n-1)K+k) \int_{\tilde{b}_{i}}^{\tilde{a}_{N+1}} 
\psi_{\widetilde{V}}^{f}(\xi) \, \md \xi}
\end{pmatrix}, &\text{$z \! \in \! (\tilde{a}_{i},\tilde{b}_{i}), \quad 
i \! = \! 1,2,\dotsc,N$,} \\
\mathrm{I} \! + \! \me^{n(g_{+}^{f}(z)+g_{-}^{f}(z)-\hat{\mathscr{P}}_{0}^{+}
-\hat{\mathscr{P}}_{0}^{-}-\widetilde{V}(z)-\tilde{\ell})} \sigma_{+}, 
&\text{$z \! \in \! (-\infty,\tilde{b}_{0}) \cup (\tilde{a}_{N+1},+\infty)$,}
\end{cases}
\end{align*}
with $g^{f}(z)$ and $\hat{\mathscr{P}}_{0}^{\pm}$ defined by 
Equations~\eqref{eql3.4gee3} and~\eqref{eql3.4gee5}, respectively, 
$g^{f}_{\pm}(z) \! := \! \lim_{\varepsilon \downarrow 0}g^{f}(z \! 
\pm \! \mi \varepsilon)$, $g_{+}^{f}(z) \! + \! g_{-}^{f}(z) \! - \! 
\hat{\mathscr{P}}_{0}^{+} \! - \! \hat{\mathscr{P}}_{0}^{-} \! - \! 
\widetilde{V}(z) \! - \! \tilde{\ell} \! < \! 0$, $z \! \in \! (-\infty,
\tilde{b}_{0}) \cup \cup_{i=1}^{N}(\tilde{a}_{i},\tilde{b}_{i}) \cup 
(\tilde{a}_{N+1},+\infty)$, and $\pm \Re (\mi \int_{z}^{\tilde{a}_{N+1}} 
\psi_{\widetilde{V}}^{f}(\xi) \, \md \xi) \! > \! 0$, $z \! \in \! 
\mathbb{C}_{\pm} \cap (\cup_{j=1}^{N+1} \tilde{\mathbb{U}}_{j})$, 
where, for $j \! = \! 1,2,\dotsc,N \! + \! 1$, $\tilde{\mathbb{U}}_{j} 
\! := \! \lbrace \mathstrut z \! \in \! \mathbb{C} \setminus \lbrace 
\alpha_{p_{1}},\dotsc,\alpha_{p_{\mathfrak{s}-2}},\alpha_{p_{\mathfrak{s}}} 
\rbrace; \, \inf_{\tau \in (\tilde{b}_{j-1},\tilde{a}_{j})} \lvert z \! - \! 
\tau \rvert \! < \! \tilde{r}_{j} \rbrace$, with $\tilde{r}_{j} \! \in \! 
(0,1)$ chosen small enough so that $\tilde{\mathbb{U}}_{i} \cap 
\tilde{\mathbb{U}}_{j} \! = \! \varnothing$, $i \! \neq \! j \! \in \! 
\lbrace 1,2,\dotsc,N \! + \! 1 \rbrace$$;$ $\pmb{{\rm (iii)}}$ for 
$n \! \in \! \mathbb{N}$ and $k \! \in \! \lbrace 1,2,\dotsc,K \rbrace$ 
such that $\alpha_{p_{\mathfrak{s}}} \! := \! \alpha_{k} \! = \! \infty$,
\begin{equation*}
\mathcal{M}(z) \underset{\overline{\mathbb{C}}_{\pm} \ni z \to 
\alpha_{k}}{=} \mathrm{I} \! + \! \mathcal{O}(z^{-1}), \, \quad 
\qquad \, \mathcal{M}(z) \underset{\mathbb{C}_{\pm} \ni z \to 
\alpha_{p_{q}}}{=} \mathcal{O}(\me^{n(\tilde{\mathscr{P}}_{0}-
\tilde{\mathscr{P}}^{\pm}_{1}) \sigma_{3}}), \quad q \! = \! 1,2,
\dotsc,\mathfrak{s} \! - \! 1,
\end{equation*}
where $\tilde{\mathscr{P}}^{\pm}_{1}$ is defined by 
Equation~\eqref{eql3.4gee10}$;$ and $\pmb{{\rm (iv)}}$ for $n \! \in 
\! \mathbb{N}$ and $k \! \in \! \lbrace 1,2,\dotsc,K \rbrace$ such that 
$\alpha_{p_{\mathfrak{s}}} \! := \! \alpha_{k} \! \neq \! \infty$,
\begin{gather*}
\mathcal{M}(z) \underset{\mathbb{C}_{\pm} \ni z \to \alpha_{k}}{=} 
(\mathrm{I} \! + \! \mathcal{O}(z \! - \! \alpha_{k})) \mathscr{E}^{\mp 
\sigma_{3}}, \, \quad \qquad \, \mathcal{M}(z) \underset{\overline{
\mathbb{C}}_{\pm} \ni z \to \alpha_{p_{\mathfrak{s}-1}} = \infty}{=} 
\mathcal{O}(\me^{n(\hat{\mathscr{P}}^{\pm}_{0}-\hat{\mathscr{P}}_{1}) 
\sigma_{3}}), \\
\mathcal{M}(z) \underset{\mathbb{C}_{\pm} \ni z \to \alpha_{p_{q}}}{=} 
\mathcal{O}(\me^{n(\hat{\mathscr{P}}^{\pm}_{0}-\hat{\mathscr{P}}^{
\pm}_{2}) \sigma_{3}}), \quad q \! = \! 1,2,\dotsc,\mathfrak{s} \! - \! 2,
\end{gather*}
where $\hat{\mathscr{P}}_{1}$ and $\hat{\mathscr{P}}^{\pm}_{2}$ are 
defined by Equations~\eqref{eql3.4gee7} and~\eqref{eql3.4gee9}, 
respectively.
\end{ccccc}

\emph{Proof}. The proof of this Lemma~\ref{lem4.1} consists of two cases: 
(i) $n \! \in \! \mathbb{N}$ and $k \! \in \! \lbrace 1,2,\dotsc,K \rbrace$ 
such that $\alpha_{p_{\mathfrak{s}}} \! := \! \alpha_{k} \! = \! \infty$; and 
(ii) $n \! \in \! \mathbb{N}$ and $k \! \in \! \lbrace 1,2,\dotsc,K \rbrace$ 
such that $\alpha_{p_{\mathfrak{s}}} \! := \! \alpha_{k} \! \neq \! \infty$. 
The proof for the case $\alpha_{p_{\mathfrak{s}}} \! := \! \alpha_{k} \! \neq 
\! \infty$, $k \! \in \! \lbrace 1,2,\dotsc,K \rbrace$, will be considered 
in detail (see $\pmb{(\mathrm{A})}$ below), whilst the case 
$\alpha_{p_{\mathfrak{s}}} \! := \! \alpha_{k} \! = \! \infty$, $k \! 
\in \! \lbrace 1,2,\dotsc,K \rbrace$, can be proved analogously (see 
$\pmb{(\mathrm{B})}$ below).

$\pmb{(\mathrm{A})}$ Let $\widetilde{V} \colon \overline{\mathbb{R}} 
\setminus \lbrace \alpha_{1},\alpha_{2},\dotsc,\alpha_{K} \rbrace \! \to 
\! \mathbb{R}$ satisfy conditions~\eqref{eq20}--\eqref{eq22} and be 
regular. For $n \! \in \! \mathbb{N}$ and $k \! \in \! \lbrace 1,2,\dotsc,
K \rbrace$ such that $\alpha_{p_{\mathfrak{s}}} \! := \! \alpha_{k} \! 
\neq \! \infty$, the corresponding item~$\pmb{{\rm (i)}}$ of the lemma 
follows {}from the definition of $\mathcal{M}(z)$ in terms of $\mathscr{X}
(z)$ stated in the formulation of the lemma and the corresponding 
item~$\pmb{{\rm (i)}}$ of Lemma~\ref{lem3.4}. Recall {}from the proof of 
Lemma~\ref{lem3.8} (cf. Equations~\eqref{eql3.8B10}, \eqref{eql3.8B11}, 
\eqref{eql3.8B12}, and~\eqref{eql3.8B13}) that, for $j \! = \! 1,2,\dotsc,
N \! + \! 1$ and $i \! = \! 1,2,\dotsc,N$,
\begin{equation*}
(2 \pi \mi)^{-1}(g_{+}^{f}(z) \! - \! g_{-}^{f}(z) \! + \! 
\hat{\mathscr{P}}_{0}^{-} \! - \! \hat{\mathscr{P}}_{0}^{+}) \! = \! 
\begin{cases}
{\fontsize{8pt}{9pt}\selectfont 
\begin{aligned}[b]
&\left(\frac{(n \! - \! 1)K \! + \! k}{n} \right) \int_{z}^{\tilde{a}_{N+1}} 
\psi_{\widetilde{V}}^{f}(\xi) \, \md \xi \! - \! \left(\frac{(n \! - \! 1)K 
\! + \! k}{n} \right) \int_{J_{f} \cap \mathbb{R}_{\alpha_{k}}^{>}} 
\psi_{\widetilde{V}}^{f}(\xi) \, \md \xi \\
&-\sum_{q \in \tilde{\Delta}(j;z)} \frac{\varkappa_{nk \tilde{k}_{q}}}{n} \! + 
\! \sum_{q \in \tilde{\Delta}(k)} \frac{\varkappa_{nk \tilde{k}_{q}}}{n} \! - 
\! \left(\frac{\varkappa_{nk} \! - \! 1}{n} \right) 
\chi_{\mathbb{R}_{\alpha_{k}}^{<}}(z),
\end{aligned}} \text{$z \! \in \! [\tilde{b}_{j-1},\tilde{a}_{j}]$,} \\
{\fontsize{8pt}{9pt}\selectfont 
\begin{aligned}[b]
&\left(\frac{(n \! - \! 1)K \! + \! k}{n} \right) \int_{\tilde{b}_{i}}^{
\tilde{a}_{N+1}} \psi_{\widetilde{V}}^{f}(\xi) \, \md \xi \! - \! \left(
\frac{(n \! - \! 1)K \! + \! k}{n} \right) \int_{J_{f} \cap \mathbb{R}_{
\alpha_{k}}^{>}} \psi_{\widetilde{V}}^{f}(\xi) \, \md \xi \\
&-\sum_{q \in \tilde{\Delta}(i;z)} \frac{\varkappa_{nk \tilde{k}_{q}}}{n} \! 
+ \! \sum_{q \in \tilde{\Delta}(k)} \frac{\varkappa_{nk \tilde{k}_{q}}}{n} 
\! - \! \left(\frac{\varkappa_{nk} \! - \! 1}{n} \right) 
\chi_{\mathbb{R}_{\alpha_{k}}^{<}}(z),
\end{aligned}} \text{$z \! \in \! (\tilde{a}_{i},\tilde{b}_{i})$,} \\
{\fontsize{8pt}{9pt}\selectfont 
\begin{aligned}[b]
&-\left(\frac{(n \! - \! 1)K \! + \! k}{n} \right) \int_{J_{f} \cap 
\mathbb{R}_{\alpha_{k}}^{>}} \psi_{\widetilde{V}}^{f}(\xi) \, \md \xi \! - 
\! \sum_{q \in \tilde{\Delta}(z)} \frac{\varkappa_{nk \tilde{k}_{q}}}{n} \! 
+ \! \sum_{q \in \tilde{\Delta}(k)} \frac{\varkappa_{nk \tilde{k}_{q}}}{n} \\
&-\left(\frac{\varkappa_{nk} \! - \! 1}{n} \right) 
\chi_{\mathbb{R}_{\alpha_{k}}^{<}}(z),
\end{aligned}} \text{$z \! \in \! (\tilde{a}_{N+1},+\infty)$,} \\
{\fontsize{8pt}{9pt}\selectfont 
\begin{aligned}[b]
&\left(\frac{(n \! - \! 1)K \! + \! k}{n} \right) \! - \! \left(
\frac{(n \! - \! 1)K \! + \! k}{n} \right) \int_{J_{f} \cap \mathbb{R}_{
\alpha_{k}}^{>}} \psi_{\widetilde{V}}^{f}(\xi) \, \md \xi \! - \! \sum_{q 
\in \tilde{\Delta}(z)} \frac{\varkappa_{nk \tilde{k}_{q}}}{n} \\
&+\sum_{q \in \tilde{\Delta}(k)} \frac{\varkappa_{nk \tilde{k}_{q}}}{n} 
\! - \! \left(\frac{\varkappa_{nk} \! - \! 1}{n} \right) 
\chi_{\mathbb{R}_{\alpha_{k}}^{<}}(z),
\end{aligned}} \text{$z \! \in \! (-\infty,\tilde{b}_{0})$,}
\end{cases}
\end{equation*}
where $\tilde{\Delta}(j;z)$, $\tilde{\Delta}(k)$, and $\tilde{\Delta}
(z)$ are defined in the proof of Lemma~\ref{lem3.8}, and (cf. 
Equations~\eqref{eql3.8d}, \eqref{eql3.8AA}, \eqref{eql3.8BB}, 
and~\eqref{eql3.8CC})
\begin{equation*}
g_{+}^{f}(z) \! + \! g_{-}^{f}(z) \! - \! \hat{\mathscr{P}}_{0}^{+} \! - \! 
\hat{\mathscr{P}}_{0}^{-} \! - \! \widetilde{V}(z) \! - \! \tilde{\ell} \! = 
\! 
\begin{cases}
0, &\text{$z \! \in \! [\tilde{b}_{j-1},\tilde{a}_{j}]$,} \\
-\left(\frac{(n-1)K+k}{n} \right) \int_{\tilde{a}_{i}}^{z}(\tilde{R}
(\xi))^{1/2} \tilde{h}_{\widetilde{V}}(\xi) \, \md \xi \, \, \, (< \! 0), 
&\text{$z \! \in \! (\tilde{a}_{i},\tilde{b}_{i})$,} \\
-\left(\frac{(n-1)K+k}{n} \right) \int_{\tilde{a}_{N+1}}^{z}(\tilde{R}
(\xi))^{1/2} \tilde{h}_{\widetilde{V}}(\xi) \, \md \xi \, \, \, (< \! 0), 
&\text{$z \! \in \! (\tilde{a}_{N+1},+\infty)$,} \\
\left(\frac{(n-1)K+k}{n} \right) \int_{z}^{\tilde{b}_{0}}(\tilde{R}
(\xi))^{1/2} \tilde{h}_{\widetilde{V}}(\xi) \, \md \xi \, \, \, (< \! 0), 
&\text{$z \! \in \! (-\infty,\tilde{b}_{0})$.}
\end{cases}
\end{equation*}
{}From the expression for the associated jump matrix given in 
Equation~\eqref{eql3.4h}, the above formulae, the fact that, for 
$n \! \in \! \mathbb{N}$ and $k \! \in \! \lbrace 1,2,\dotsc,K 
\rbrace$ such that $\alpha_{p_{\mathfrak{s}}} \! := \! \alpha_{k} 
\! \neq \! \infty$, $\varkappa_{nk} \! \in \! \mathbb{N}$ and 
$\varkappa_{nk \tilde{k}_{q}} \! \in \! \mathbb{N}_{0}$, $q \! = \! 
1,\dotsc,\mathfrak{s} \! - \! 2,\mathfrak{s}$, and the definition of 
$\mathcal{M}(z)$ in terms of $\mathscr{X}(z)$ stated in the formulation 
of the lemma, one arrives at the expression for $\hat{\mathscr{V}}(z)$ 
stated in the corresponding item~$\pmb{{\rm (ii)}}$ of the lemma; 
furthermore, the associated asymptotics in item~$\pmb{{\rm (iv)}}$ 
of the lemma are a consequence of the definition of $\mathcal{M}(z)$ 
in terms of $\mathscr{X}(z)$ stated in the formulation of the lemma 
and the asymptotics in the corresponding item~$\pmb{{\rm (iv)}}$ of 
Lemma~\ref{lem3.4}. It remains, therefore, to show that $\pm \Re 
(\mi \int_{z}^{\tilde{a}_{N+1}} \psi^{f}_{\widetilde{V}}(\xi) \, \md \xi) \! 
> \! 0$, $z \! \in \! \mathbb{C}_{\pm} \cap (\cup_{j=1}^{N+1} \tilde{
\mathbb{U}}_{j})$, where $\tilde{\mathbb{U}}_{j} \! := \! \lbrace \mathstrut 
z \! \in \! \mathbb{C} \setminus \lbrace \alpha_{p_{1}},\dotsc,\alpha_{
p_{\mathfrak{s}-2}},\alpha_{p_{\mathfrak{s}}} \rbrace; \, \inf_{\tau \in 
(\tilde{b}_{j-1},\tilde{a}_{j})} \lvert z \! - \! \tau \rvert \! < \! \tilde{r}_{j} 
\rbrace$, $j \! = \! 1,2,\dotsc,N \! + \! 1$, with $\tilde{r}_{j}$, which is 
an arbitrarily fixed, sufficiently small positive real number, chosen so that 
$\tilde{\mathbb{U}}_{i} \cap \tilde{\mathbb{U}}_{j} \! = \! \varnothing$, 
$i \! \neq \! j \! \in \! \lbrace 1,2,\dotsc,N \! + \! 1 \rbrace$. Via the 
uniform Lipschitz continuity of $g^{f}(z)$ for $z \! \in \! \mathbb{C}_{\pm}$ 
(cf. the proof of Lemma~\ref{lem3.4}), the Cauchy-Riemann conditions, and 
Equations~\eqref{eql3.8D1}, \eqref{eql3.8D2}, \eqref{eql3.8D3}, 
and~\eqref{eql3.8D4}, it follows that $g_{+}^{f}(z) \! - \! g_{-}^{f}(z) \! + \! 
\hat{\mathscr{P}}_{0}^{-} \! - \! \hat{\mathscr{P}}_{0}^{+}$ has an analytic 
extension,\footnote{Since, by piecewise continuity, $\exp (2 \pi \mi 
(\varkappa_{nk} \! - \! 1) \chi_{\mathbb{R}^{<}_{\alpha_{k}}}(z)) \! = \! 
\exp (2 \pi \mi \sum_{q \in \tilde{\Delta}(j;z)} \varkappa_{nk \tilde{k}_{q}}) \! 
= \! \exp (2 \pi \mi \sum_{q \in \tilde{\Delta}(z)} \varkappa_{nk \tilde{k}_{q}}) 
\! = \! 1$.} denoted $\tilde{\mathscr{G}}(z)$, to an open neighbourhood 
of $J_{f}$, namely, $\cup_{j=1}^{N+1} \tilde{\mathbb{U}}_{j}$, where 
$\tilde{\mathbb{U}}_{j}$, $j \! = \! 1,2,\dotsc,N \! + \! 1$, is defined above, with 
the property that, for $z \! \in \! \mathbb{C}_{\pm} \cap (\cup_{j=1}^{N+1} 
\tilde{\mathbb{U}}_{j})$, $\pm \Re (\tilde{\mathscr{G}}(z)) \! > \! 0$.

$\pmb{(\mathrm{B})}$ The proof of this case, that is, $n \! \in \! 
\mathbb{N}$ and $k \! \in \! \lbrace 1,2,\dotsc,K \rbrace$ such that 
$\alpha_{p_{\mathfrak{s}}} \! := \! \alpha_{k} \! = \! \infty$, is virtually 
identical to the proof presented in $\pmb{(\mathrm{A})}$ above; one 
mimics, \emph{verbatim}, the scheme of the calculations presented in 
case $\pmb{(\mathrm{A})}$ in order to arrive at the claims stated in 
the corresponding items~$\pmb{{\rm (i)}}$, $\pmb{{\rm (ii)}}$, 
and~$\pmb{{\rm (iii)}}$ of the lemma; in order to do so, however, 
the analogues of the formulae given above are requisite, which, in the 
present case, read (cf. Equations~\eqref{eql3.8B14}, \eqref{eql3.8B2}, 
\eqref{eql3.8B15}, \eqref{eql3.8B5}, \eqref{eql3.8B16}, \eqref{eql3.8B7}, 
\eqref{eql3.8B17}, and~\eqref{eql3.8B9}): for $j \! = \! 1,2,\dotsc,
N \! + \! 1$ and $i \! = \! 1,2,\dotsc,N$,
\begin{equation*}
(2 \pi \mi)^{-1}(g_{+}^{\infty}(z) \! - \! g_{-}^{\infty}(z)) \! = \! 
\begin{cases}
{\fontsize{8pt}{9pt}\selectfont 
\begin{aligned}[b]
&\left(\frac{(n \! - \! 1)K \! + \! k}{n} \right) \int_{z}^{\hat{a}_{N+1}} 
\psi_{\widetilde{V}}^{\infty}(\xi) \, \md \xi \! - \! \sum_{q \in 
\hat{\Delta}(j;z)} \frac{\varkappa_{nk \tilde{k}_{q}}}{n},
\end{aligned}} \text{$z \! \in \! [\hat{b}_{j-1},\hat{a}_{j}]$,} \\
{\fontsize{8pt}{9pt}\selectfont 
\begin{aligned}[b]
&\left(\frac{(n \! - \! 1)K \! + \! k}{n} \right) \int_{\hat{b}_{i}}^{
\hat{a}_{N+1}} \psi_{\widetilde{V}}^{\infty}(\xi) \, \md \xi \! - \! 
\sum_{q \in \hat{\Delta}(i;z)} \frac{\varkappa_{nk \tilde{k}_{q}}}{n},
\end{aligned}} \text{$z \! \in \! (\hat{a}_{i},\hat{b}_{i})$,} \\
{\fontsize{8pt}{9pt}\selectfont 
\begin{aligned}[b]
&-\sum_{q \in \hat{\Delta}(z)} \frac{\varkappa_{nk \tilde{k}_{q}}}{n},
\end{aligned}} \text{$z \! \in \! (\hat{a}_{N+1},+\infty)$,} \\
{\fontsize{8pt}{9pt}\selectfont 
\begin{aligned}[b]
&\left(\frac{(n \! - \! 1)K \! + \! k}{n} \right) \! - \! 
\sum_{q \in \hat{\Delta}(z)} \frac{\varkappa_{nk \tilde{k}_{q}}}{n},
\end{aligned}} \text{$z \! \in \! (-\infty,\hat{b}_{0})$,}
\end{cases}
\end{equation*}
where $\hat{\Delta}(j;z)$ and $\hat{\Delta}(z)$ are defined in the proof 
of Lemma~\ref{lem3.8}, and
\begin{equation*}
g_{+}^{\infty}(z) \! + \! g_{-}^{\infty}(z) \! - \! 2 \tilde{\mathscr{P}}_{0} 
\! - \! \widetilde{V}(z) \! - \! \hat{\ell} \! = \! 
\begin{cases}
0, &\text{$z \! \in \! [\hat{b}_{j-1},\hat{a}_{j}]$,} \\
-\left(\frac{(n-1)K+k}{n} \right) \int_{\hat{a}_{i}}^{z}(\hat{R}(\xi))^{1/2} 
\hat{h}_{\widetilde{V}}(\xi) \, \md \xi \, \, \, (< \! 0), 
&\text{$z \! \in \! (\hat{a}_{i},\hat{b}_{i})$,} \\
-\left( \frac{(n-1)K+k}{n} \right) \int_{\hat{a}_{N+1}}^{z}(\hat{R}
(\xi))^{1/2} \hat{h}_{\widetilde{V}}(\xi) \, \md \xi \, \, \, (< \! 0), 
&\text{$z \! \in \! (\hat{a}_{N+1},+\infty)$,} \\
\left(\frac{(n-1)K+k}{n} \right) \int_{z}^{\hat{b}_{0}}(\hat{R}(\xi))^{1/2} 
\hat{h}_{\widetilde{V}}(\xi) \, \md \xi \, \, \, (< \! 0), 
&\text{$z \! \in \! (-\infty,\hat{b}_{0})$.}
\end{cases}
\end{equation*}
This concludes the proof. \hfill $\qed$
\begin{eeeee} \label{rem4.1} 
\textsl{Recalling that the external field $\widetilde{V} \colon 
\overline{\mathbb{R}} \setminus \lbrace \alpha_{1},\alpha_{2},
\dotsc,\alpha_{K} \rbrace \! \to \! \mathbb{R}$ satisfying 
conditions~\eqref{eq20}--\eqref{eq22} is regular, and, for $n \! \in \! 
\mathbb{N}$ and $k \! \in \! \lbrace 1,2,\dotsc,K \rbrace$ such that 
$\alpha_{p_{\mathfrak{s}}} \! := \! \alpha_{k} \! = \! \infty$ (resp., 
$\alpha_{p_{\mathfrak{s}}} \! := \! \alpha_{k} \! \neq \! \infty)$, 
$g_{+}^{\infty}(z) \! + \! g_{-}^{\infty}(z) \! - \! 2 \tilde{\mathscr{P}}_{0} 
\! - \! \widetilde{V}(z) \! - \! \hat{\ell} \! < \! 0$, $z \! \in \! 
(-\infty,\hat{b}_{0}) \cup \cup_{i=1}^{N}(\hat{a}_{i},\hat{b}_{i}) \cup 
(\hat{a}_{N+1},+\infty)$ (resp., $g_{+}^{f}(z) \! + \! g_{-}^{f}(z) \! - 
\! \hat{\mathscr{P}}_{0}^{+} \! - \! \hat{\mathscr{P}}_{0}^{-} \! - \! 
\widetilde{V}(z) \! - \! \tilde{\ell} \! < \! 0$, $z \! \in \! (-\infty,
\tilde{b}_{0}) \cup \cup_{i=1}^{N}(\tilde{a}_{i},\tilde{b}_{i}) \cup 
(\tilde{a}_{N+1},+\infty))$, it follows that, for $n \! \in \! 
\mathbb{N}$ and $k \! \in \! \lbrace 1,2,\dotsc,K \rbrace$ such 
that $\alpha_{p_{\mathfrak{s}}} \! := \! \alpha_{k} \! = \! \infty$, 
\begin{equation*}
\hat{\mathscr{V}}(z) \underset{\underset{z_{o}=1+o(1)}{\mathscr{N},
n \to \infty}}{=} 
\begin{cases}
\me^{-(2 \pi \mi ((n-1)K+k) \int_{\hat{b}_{i}}^{\hat{a}_{N+1}} 
\psi_{\widetilde{V}}^{\infty}(\xi) \, \md \xi) \sigma_{3}}(\mathrm{I} 
\! + \! o(1) \sigma_{+}), &\text{$z \! \in \! (\hat{a}_{i},\hat{b}_{i}), 
\quad i \! = \! 1,2,\dotsc,N$,} \\
\mathrm{I} \! + \! o(1) \sigma_{+}, &\text{$z \! \in \! 
(-\infty,\hat{b}_{0}) \cup (\hat{a}_{N+1},+\infty)$,}
\end{cases}
\end{equation*}
where, here and below (in this section), $o(1)$ denotes terms that 
are exponentially small (in the double-scaling limit $\mathscr{N},n 
\! \to \! \infty$ such that $z_{o} \! = \! 1 \! + \! o(1))$, and, for 
$n \! \in \! \mathbb{N}$ and $k \! \in \! \lbrace 1,2,\dotsc,K \rbrace$ 
such that $\alpha_{p_{\mathfrak{s}}} \! := \! \alpha_{k} \! \neq \! \infty$,
\begin{equation*}
\hat{\mathscr{V}}(z) \underset{\underset{z_{o}=1+o(1)}{\mathscr{N},
n \to \infty}}{=} 
\begin{cases}
\me^{-(2 \pi \mi ((n-1)K+k) \int_{\tilde{b}_{i}}^{\tilde{a}_{N+1}} 
\psi_{\widetilde{V}}^{f}(\xi) \, \md \xi) \sigma_{3}}(\mathrm{I} \! + \! 
o(1) \sigma_{+}), &\text{$z \! \in \! (\tilde{a}_{i},\tilde{b}_{i}), \quad 
i \! = \! 1,2,\dotsc,N$,} \\
\mathrm{I} \! + \! o(1) \sigma_{+}, &\text{$z \! \in \! (-\infty,
\tilde{b}_{0}) \cup (\tilde{a}_{N+1},+\infty)$.}
\end{cases}
\end{equation*}}
\end{eeeee}
\begin{ccccc} \label{lem4.2} 
Let the external field $\widetilde{V} \colon \overline{\mathbb{R}} 
\setminus \lbrace \alpha_{1},\alpha_{2},\dotsc,\alpha_{K} \rbrace 
\! \to \! \mathbb{R}$ satisfy conditions~\eqref{eq20}--\eqref{eq22} 
and be regular. For $n \! \in \! \mathbb{N}$ and $k \! \in \! \lbrace 
1,2,\dotsc,K \rbrace$ such that $\alpha_{p_{\mathfrak{s}}} \! := \! 
\alpha_{k} \! = \! \infty$ (resp., $\alpha_{p_{\mathfrak{s}}} \! := \! 
\alpha_{k} \! \neq \! \infty)$, let the associated equilibrium measure, 
$\mu_{\widetilde{V}}^{\infty}$ (resp., $\mu_{\widetilde{V}}^{f})$, and 
its support, $J_{\infty}$ (resp., $J_{f})$, be as described in item~$\pmb{(1)}$ 
(resp., item~$\pmb{(2)})$ of Lemma~\ref{lem3.7}, and, along with the 
corresponding variational constant, $\hat{\ell}$ (resp., $\tilde{\ell})$, 
satisfy the variational conditions~\eqref{eql3.8a} (resp., 
conditions~\eqref{eql3.8b}$)$$;$ moreover, let the associated 
conditions~{\rm (i)}--{\rm (iv)} of item~$\pmb{(1)}$ (resp., 
item~$\pmb{(2)})$ of Lemma~\ref{lem3.8} be valid. For $n \! \in \! 
\mathbb{N}$ and $k \! \in \! \lbrace 1,2,\dotsc,K \rbrace$, let $\mathcal{M} 
\colon \mathbb{N} \times \lbrace 1,2,\dotsc,K \rbrace \times \overline{
\mathbb{C}} \setminus \overline{\mathbb{R}} \! \to \! \operatorname{SL}_{2}
(\mathbb{C})$ solve the matrix {\rm RHP} $(\mathcal{M}(z),\hat{\mathscr{V}}
(z),\overline{\mathbb{R}})$ formulated in Lemma~\ref{lem4.1}, and, for $n \! 
\in \! \mathbb{N}$ and $k \! \in \! \lbrace 1,2,\dotsc,K \rbrace$ such that 
$\alpha_{p_{\mathfrak{s}}} \! := \! \alpha_{k} \! = \! \infty$, let the deformed 
and oriented contour $\hat{\Sigma} \! := \! \mathbb{R} \cup \cup_{j=1}^{N+1}
(\hat{J}_{j}^{\smallfrown} \cup \hat{J}_{j}^{\smallsmile})$ be as in 
Figure~\ref{figgone}, with $\mathbb{R} \! = \! (-\infty,\hat{b}_{0}) \cup 
\cup_{j=1}^{N+1} \hat{J}_{j} \cup \cup_{i=1}^{N}(\hat{a}_{i},\hat{b}_{i}) \cup 
(\hat{a}_{N+1},+\infty)$, where $\hat{J}_{j} \! := \! [\hat{b}_{j-1},\hat{a}_{j}]$, 
$j \! = \! 1,2,\dotsc,N \! + \! 1$, and $\cup_{j=1}^{N+1}(\hat{\Omega}_{j}^{
\smallfrown} \cup \hat{\Omega}_{j}^{\smallsmile} \cup \hat{J}_{j}^{\smallfrown} 
\cup \hat{J}_{j}^{\smallsmile}) \subset \cup_{j=1}^{N+1} \hat{\mathbb{U}}_{j}$, 
where $\hat{\mathbb{U}}_{j}$, $j \! = \! 1,2,\dotsc,N \! + \! 1$, is defined in 
the corresponding item~$\pmb{{\rm (ii)}}$ of Lemma~\ref{lem4.1}, and set
\begin{equation*}
\mathfrak{M}(z) \! = \! 
\begin{cases}
\mathcal{M}(z), &\text{$z \! \in \! \mathbb{C} \setminus (\hat{\Sigma} 
\cup \cup_{j=1}^{N+1}(\hat{\Omega}_{j}^{\smallfrown} \cup 
\hat{\Omega}_{j}^{\smallsmile}))$,} \\
\mathcal{M}(z) \left(\mathrm{I} \! - \! \me^{-2 \pi \mi ((n-1)K+k) 
\int_{z}^{\hat{a}_{N+1}} \psi_{\widetilde{V}}^{\infty}(\xi) \, \md \xi} 
\, \sigma_{-} \right), &\text{$z \! \in \! \mathbb{C}_{+} \cap 
(\cup_{j=1}^{N+1} \hat{\Omega}_{j}^{\smallfrown})$,} \\
\mathcal{M}(z) \left(\mathrm{I} \! + \! \me^{2 \pi \mi ((n-1)K+k) 
\int_{z}^{\hat{a}_{N+1}} \psi_{\widetilde{V}}^{\infty}(\xi) \, \md \xi} 
\, \sigma_{-} \right), &\text{$z \! \in \! \mathbb{C}_{-} \cap 
(\cup_{j=1}^{N+1} \hat{\Omega}_{j}^{\smallsmile})$,}
\end{cases}
\end{equation*}
and, for $n \! \in \! \mathbb{N}$ and $k \! \in \! \lbrace 1,2,\dotsc,K \rbrace$ 
such that $\alpha_{p_{\mathfrak{s}}} \! := \! \alpha_{k} \! \neq \! \infty$, let 
the deformed and oriented contour $\tilde{\Sigma} \! := \! \mathbb{R} \cup 
\cup_{j=1}^{N+1}(\tilde{J}_{j}^{\smallfrown} \cup \tilde{J}_{j}^{\smallsmile})$ 
be as in Figure~\ref{figgtwo}, with $\mathbb{R} \! = \! (-\infty,\tilde{b}_{0}) \cup 
\cup_{j=1}^{N+1} \tilde{J}_{j} \cup \cup_{i=1}^{N}(\tilde{a}_{i},\tilde{b}_{i}) 
\cup (\tilde{a}_{N+1},+\infty)$, where $\tilde{J}_{j} \! := \! [\tilde{b}_{j-1},
\tilde{a}_{j}]$, $j \! = \! 1,2,\dotsc,N \! + \! 1$, and $\cup_{j=1}^{N+1}
(\tilde{\Omega}_{j}^{\smallfrown} \cup \tilde{\Omega}_{j}^{\smallsmile} 
\cup \tilde{J}_{j}^{\smallfrown} \cup \tilde{J}_{j}^{\smallsmile}) \subset 
\cup_{j=1}^{N+1} \tilde{\mathbb{U}}_{j}$, where $\tilde{\mathbb{U}}_{j}$, 
$j \! = \! 1,2,\dotsc,N \! + \! 1$, is defined in the corresponding 
item~$\pmb{{\rm (ii)}}$ of Lemma~\ref{lem4.1}, and set
\begin{equation*}
\mathfrak{M}(z) \! = \! 
\begin{cases}
\mathcal{M}(z), &\text{$z \! \in \! \mathbb{C} \setminus (\tilde{\Sigma} 
\cup \cup_{j=1}^{N+1}(\tilde{\Omega}_{j}^{\smallfrown} \cup 
\tilde{\Omega}_{j}^{\smallsmile}))$,} \\
\mathcal{M}(z) \left(\mathrm{I} \! - \! \me^{-2 \pi \mi ((n-1)K+k) 
\int_{z}^{\tilde{a}_{N+1}} \psi_{\widetilde{V}}^{f}(\xi) \, \md \xi} 
\, \sigma_{-} \right), &\text{$z \! \in \! \mathbb{C}_{+} \cap 
(\cup_{j=1}^{N+1} \tilde{\Omega}_{j}^{\smallfrown})$,} \\
\mathcal{M}(z) \left(\mathrm{I} \! + \! \me^{2 \pi \mi ((n-1)K+k) 
\int_{z}^{\tilde{a}_{N+1}} \psi_{\widetilde{V}}^{f}(\xi) \, \md \xi} 
\, \sigma_{-} \right), &\text{$z \! \in \! \mathbb{C}_{-} \cap 
(\cup_{j=1}^{N+1} \tilde{\Omega}_{j}^{\smallsmile})$.}
\end{cases}
\end{equation*}
Then: $\pmb{(1)}$ for $n \! \in \! \mathbb{N}$ and $k \! \in \! \lbrace 
1,2,\dotsc,K \rbrace$ such that $\alpha_{p_{\mathfrak{s}}} \! := \! 
\alpha_{k} \! = \! \infty$, $\mathfrak{M} \colon \mathbb{N} \times 
\lbrace 1,2,\dotsc,K \rbrace \times \mathbb{C} \setminus 
\hat{\Sigma} \! \to \! \operatorname{SL}_{2}(\mathbb{C})$ solves the 
following equivalent matrix {\rm RHP}$:$ {\rm (i)} $\mathfrak{M}
(n,k,z) \! = \! \mathfrak{M}(z)$ is analytic for $z \! \in \! \mathbb{C} 
\setminus \hat{\Sigma}$$;$ {\rm (ii)} the boundary values $\mathfrak{M}_{\pm}
(z) \! := \! \lim_{\underset{z^{\prime} \in \pm \, \mathrm{side} \, 
\mathrm{of} \, \hat{\Sigma}}{z^{\prime} \to z \in \hat{\Sigma}}} \mathfrak{M}
(z^{\prime})$ satisfy the jump condition
\begin{equation*}
\mathfrak{M}_{+}(z) \! = \! \mathfrak{M}_{-}(z) \mathfrak{v}(z) \quad 
\mathrm{a.e.} \quad z \! \in \! \hat{\Sigma},
\end{equation*}
where
\begin{align*}
\mathfrak{v} &\colon \mathbb{N} \times \lbrace 1,2,\dotsc,K \rbrace 
\times \hat{\Sigma} \! \to \! \operatorname{GL}_{2}(\mathbb{C}), \, 
(n,k,z) \! \mapsto \! \mathfrak{v}(n,k,z) \! =: \! \mathfrak{v}(z) \\
&= 
\begin{cases}
\mi \sigma_{2}, &\text{$z \! \in \! (\hat{b}_{j-1},\hat{a}_{j}), \quad 
j \! = \! 1,2,\dotsc,N \! + \! 1$,} \\
\mathrm{I} \! + \! \me^{-2 \pi \mi ((n-1)K+k) \int_{z}^{\hat{a}_{N+1}} 
\psi_{\widetilde{V}}^{\infty}(\xi) \, \md \xi} \, \sigma_{-}, 
&\text{$z \! \in \! \hat{J}_{j}^{\smallfrown}, \quad j \! = \! 1,2,\dotsc,
N \! + \! 1$,} \\
\mathrm{I} \! + \! \me^{2 \pi \mi ((n-1)K+k) \int_{z}^{\hat{a}_{N+1}} 
\psi_{\widetilde{V}}^{\infty}(\xi) \, \md \xi} \, \sigma_{-}, 
&\text{$z \! \in \! \hat{J}_{j}^{\smallsmile}, \quad j \! = \! 1,2,\dotsc,
N \! + \! 1$,} \\
\begin{pmatrix}
\me^{-2 \pi \mi ((n-1)K+k) \int_{\hat{b}_{i}}^{\hat{a}_{N+1}} 
\psi_{\widetilde{V}}^{\infty}(\xi) \, \md \xi} & \me^{n(g_{+}^{\infty}(z)
+g_{-}^{\infty}(z)-2 \tilde{\mathscr{P}}_{0}-\widetilde{V}(z)-\hat{\ell})} \\
0 & \me^{2 \pi \mi ((n-1)K+k) \int_{\hat{b}_{i}}^{\hat{a}_{N+1}} 
\psi_{\widetilde{V}}^{\infty}(\xi) \, \md \xi}
\end{pmatrix}, &\text{$z \! \in \! (\hat{a}_{i},\hat{b}_{i}), \quad 
i \! = \! 1,2,\dotsc,N$,} \\
\mathrm{I} \! + \! \me^{n(g_{+}^{\infty}(z)+g_{-}^{\infty}(z)-2 
\tilde{\mathscr{P}}_{0}-\widetilde{V}(z)-\hat{\ell})} \sigma_{+}, 
&\text{$z \! \in \! (-\infty,\hat{b}_{0}) \cup (\hat{a}_{N+1},+\infty)$,}
\end{cases}
\end{align*}
with $\Re (\mi \int_{z}^{\hat{a}_{N+1}} \psi_{\widetilde{V}}^{\infty}(\xi) \, 
\md \xi) \! > \! 0$, $z \! \in \! \hat{\Omega}_{j}^{\smallfrown}$ (resp., 
$\Re (\mi \int_{z}^{\hat{a}_{N+1}} \psi_{\widetilde{V}}^{\infty}(\xi) \, \md \xi) 
\! < \! 0$, $z \! \in \! \hat{\Omega}_{j}^{\smallsmile})$, $j \! = \! 1,2,
\dotsc,N \! + \! 1$$:$ {\rm (iii)} for $\mathfrak{p} \! \in \! \lbrace +,- 
\rbrace$, with $\mathfrak{q}(+) \! := \smallfrown$, and $\mathfrak{q}(-) 
\! := \smallsmile$,
\begin{equation*}
\mathfrak{M}(z) \underset{\overline{\mathbb{C}}_{\mathfrak{p}} 
\setminus \cup_{j=1}^{N+1}(\hat{J}_{j}^{\mathfrak{q}(\mathfrak{p})} 
\cup \hat{\Omega}_{j}^{\mathfrak{q}(\mathfrak{p})}) \ni z \to \alpha_{k}}{=} 
\mathrm{I} \! + \! \mathcal{O}(z^{-1}), \, \quad \quad \mathfrak{M}(z) 
\underset{\mathbb{C}_{\mathfrak{p}} \setminus \cup_{j=1}^{N+1}
(\hat{J}_{j}^{\mathfrak{q}(\mathfrak{p})} \cup \hat{\Omega}_{j}^{\mathfrak{q}
(\mathfrak{p})}) \ni z \to \alpha_{p_{q}}}{=} \mathcal{O}
(\me^{n(\tilde{\mathscr{P}}_{0}-\tilde{\mathscr{P}}^{\mathfrak{p}}_{1}) 
\sigma_{3}}), \quad q \! = \! 1,2,\dotsc,\mathfrak{s} \! - \! 1;
\end{equation*}
and $\pmb{(2)}$ for $n \! \in \! \mathbb{N}$ and $k \! \in \! \lbrace 1,2,
\dotsc,K \rbrace$ such that $\alpha_{p_{\mathfrak{s}}} \! := \! \alpha_{k} 
\! \neq \! \infty$, $\mathfrak{M} \colon \mathbb{N} \times \lbrace 
1,2,\dotsc,K \rbrace \times \mathbb{C} \setminus \tilde{\Sigma} \! \to \! 
\operatorname{SL}_{2}(\mathbb{C})$ solves the following equivalent matrix 
{\rm RHP}$:$ {\rm (i)} $\mathfrak{M}(n,k,z) \! = \! \mathfrak{M}(z)$ is 
analytic for $z \! \in \! \mathbb{C} \setminus \tilde{\Sigma}$$;$ {\rm (ii)} 
the boundary values $\mathfrak{M}_{\pm}(z) \! := \! 
\lim_{\underset{z^{\prime} \in \pm \, \mathrm{side} \, \mathrm{of} \, 
\tilde{\Sigma}}{z^{\prime} \to z \in \tilde{\Sigma}}} \mathfrak{M}
(z^{\prime})$ satisfy the jump condition
\begin{equation*}
\mathfrak{M}_{+}(z) \! = \! \mathfrak{M}_{-}(z) \mathfrak{v}(z) \quad 
\mathrm{a.e.} \quad z \! \in \! \tilde{\Sigma},
\end{equation*}
where
\begin{align*}
\mathfrak{v} &\colon \mathbb{N} \times \lbrace 1,2,\dotsc,K \rbrace \times 
\tilde{\Sigma} \! \to \! \operatorname{GL}_{2}(\mathbb{C}), \, (n,k,z) \! 
\mapsto \! \mathfrak{v}(n,k,z) \! =: \! \mathfrak{v}(z) \\
&= 
\begin{cases}
\mi \sigma_{2}, &\text{$z \! \in \! (\tilde{b}_{j-1},\tilde{a}_{j}), \quad 
j \! = \! 1,2,\dotsc,N \! + \! 1$,} \\
\mathrm{I} \! + \! \me^{-2 \pi \mi ((n-1)K+k) \int_{z}^{\tilde{a}_{N+1}} 
\psi_{\widetilde{V}}^{f}(\xi) \, \md \xi} \, \sigma_{-}, 
&\text{$z \! \in \! \tilde{J}_{j}^{\smallfrown}, \quad j \! = \! 1,2,\dotsc,
N \! + \! 1$,} \\
\mathrm{I} \! + \! \me^{2 \pi \mi ((n-1)K+k) \int_{z}^{\tilde{a}_{N+1}} 
\psi_{\widetilde{V}}^{f}(\xi) \, \md \xi} \, \sigma_{-}, 
&\text{$z \! \in \! \tilde{J}_{j}^{\smallsmile}, \quad j \! = \! 1,2,\dotsc,
N \! + \! 1$,} \\
\begin{pmatrix}
\me^{-2 \pi \mi ((n-1)K+k) \int_{\tilde{b}_{i}}^{\tilde{a}_{N+1}} 
\psi_{\widetilde{V}}^{f}(\xi) \, \md \xi} & \me^{n(g_{+}^{f}(z)+g_{-}^{f}
(z)-\hat{\mathscr{P}}_{0}^{+}-\hat{\mathscr{P}}_{0}^{-}-\widetilde{V}(z)-
\tilde{\ell})} \\
0 & \me^{2 \pi \mi ((n-1)K+k) \int_{\tilde{b}_{i}}^{\tilde{a}_{N+1}} 
\psi_{\widetilde{V}}^{f}(\xi) \, \md \xi}
\end{pmatrix}, &\text{$z \! \in \! (\tilde{a}_{i},\tilde{b}_{i}), \quad 
i \! = \! 1,2,\dotsc,N$,} \\
\mathrm{I} \! + \! \me^{n(g_{+}^{f}(z)+g_{-}^{f}(z)-\hat{\mathscr{P}}_{0}^{+}
-\hat{\mathscr{P}}_{0}^{-}-\widetilde{V}(z)-\tilde{\ell})} \sigma_{+}, 
&\text{$z \! \in \! (-\infty,\tilde{b}_{0}) \cup (\tilde{a}_{N+1},+\infty)$,}
\end{cases}
\end{align*}
with $\Re (\mi \int_{z}^{\tilde{a}_{N+1}} \psi_{\widetilde{V}}^{f}(\xi) \, 
\md \xi) \! > \! 0$, $z \! \in \! \tilde{\Omega}_{j}^{\smallfrown}$ (resp., 
$\Re (\mi \int_{z}^{\tilde{a}_{N+1}} \psi_{\widetilde{V}}^{f}(\xi) \, \md \xi) 
\! < \! 0$, $z \! \in \! \tilde{\Omega}_{j}^{\smallsmile})$, $j \! = \! 1,2,
\dotsc,N \! + \! 1$$:$ {\rm (iii)} for $\mathfrak{p} \! \in \! \lbrace +,- 
\rbrace$, with $\mathfrak{q}(+) \! := \smallfrown$, $\mathfrak{q}(-) \! := 
\smallsmile$, $\mathfrak{r}(+) \! := \! -1$, and $\mathfrak{r}(-) \! := \! 
+1$,
\begin{gather*}
\mathfrak{M}(z) \underset{\mathbb{C}_{\mathfrak{p}} \setminus 
\cup_{j=1}^{N+1}(\tilde{J}_{j}^{\mathfrak{q}(\mathfrak{p})} \cup 
\tilde{\Omega}_{j}^{\mathfrak{q}(\mathfrak{p})}) \ni z \to \alpha_{k}}{=} 
(\mathrm{I} \! + \! \mathcal{O}(z \! - \! \alpha_{k})) \mathscr{E}^{\mathfrak{r}
(\mathfrak{p}) \sigma_{3}}, \, \qquad \quad \mathfrak{M}(z) \underset{
\overline{\mathbb{C}}_{\mathfrak{p}} \setminus \cup_{j=1}^{N+1}
(\tilde{J}_{j}^{\mathfrak{q}(\mathfrak{p})} \cup \tilde{\Omega}_{j}^{
\mathfrak{q}(\mathfrak{p})}) \ni z \to \alpha_{p_{\mathfrak{s}-1}} = \infty}{=} 
\mathcal{O}(\me^{n(\hat{\mathscr{P}}_{0}^{\mathfrak{p}}-\hat{\mathscr{P}}_{1}) 
\sigma_{3}} \mathscr{E}^{\mathfrak{r}(\mathfrak{p}) \sigma_{3}}), \\
\mathfrak{M}(z) \underset{\mathbb{C}_{\mathfrak{p}} \setminus 
\cup_{j=1}^{N+1}(\tilde{J}_{j}^{\mathfrak{q}(\mathfrak{p})} \cup 
\tilde{\Omega}_{j}^{\mathfrak{q}(\mathfrak{p})}) \ni z \to \alpha_{p_{q}}}{=} 
\mathcal{O}(\me^{n(\hat{\mathscr{P}}_{0}^{\mathfrak{p}}-\hat{\mathscr{P}}^{
\mathfrak{p}}_{2}) \sigma_{3}} \mathscr{E}^{\mathfrak{r}(\mathfrak{p}) 
\sigma_{3}}), \quad q \! = \! 1,2,\dotsc,\mathfrak{s} \! - \! 2.
\end{gather*}
\end{ccccc}
\begin{figure}[tbh]
\begin{center}
\vspace{-0.40cm}

\end{center}
\vspace{-1.00cm}
\caption{Oriented and deformed contour $\hat{\Sigma} \! := \! 
\mathbb{R} \cup \cup_{j=1}^{N+1}(\hat{J}_{j}^{\smallfrown} \cup 
\hat{J}_{j}^{\smallsmile})$}\label{figgone}
\end{figure}
\begin{figure}[tbh]
\begin{center}
\vspace{-0.40cm}
%
\end{center}
\vspace{-1.00cm}
\caption{Oriented and deformed contour $\tilde{\Sigma} \! := \! 
\mathbb{R} \cup \cup_{j=1}^{N+1}(\tilde{J}_{j}^{\smallfrown} \cup 
\tilde{J}_{j}^{\smallsmile})$}\label{figgtwo}
\end{figure}

\emph{Proof}. The proof of this Lemma~\ref{lem4.2} consists of two cases: 
(i) $n \! \in \! \mathbb{N}$ and $k \! \in \! \lbrace 1,2,\dotsc,K \rbrace$ 
such that $\alpha_{p_{\mathfrak{s}}} \! := \! \alpha_{k} \! = \! \infty$; and 
(ii) $n \! \in \! \mathbb{N}$ and $k \! \in \! \lbrace 1,2,\dotsc,K \rbrace$ 
such that $\alpha_{p_{\mathfrak{s}}} \! := \! \alpha_{k} \! \neq \! \infty$. 
The proof for the case $\alpha_{p_{\mathfrak{s}}} \! := \! \alpha_{k} \! \neq 
\! \infty$, $k \! \in \! \lbrace 1,2,\dotsc,K \rbrace$, will be considered 
in detail (see $\pmb{(\mathrm{A})}$ below), whilst the case 
$\alpha_{p_{\mathfrak{s}}} \! := \! \alpha_{k} \! = \! \infty$, $k \! \in 
\! \lbrace 1,2,\dotsc,K \rbrace$, can be proved analogously (see 
$\pmb{(\mathrm{B})}$ below).

$\pmb{(\mathrm{A})}$ Let $\widetilde{V} \colon \overline{\mathbb{R}} 
\setminus \lbrace \alpha_{1},\alpha_{2},\dotsc,\alpha_{K} \rbrace \! \to 
\! \mathbb{R}$ satisfy conditions~\eqref{eq20}--\eqref{eq22} and be 
regular. For $n \! \in \! \mathbb{N}$ and $k \! \in \! \lbrace 1,2,\dotsc,K 
\rbrace$ such that $\alpha_{p_{\mathfrak{s}}} \! := \! \alpha_{k} \! \neq 
\! \infty$, the corresponding items~(i) and~(iii) for case~$\pmb{(2)}$ of 
the matrix RHP for $\mathfrak{M}(z)$ formulated in the lemma are a 
consequence of the definition of $\mathfrak{M}(z)$ in terms of 
$\mathcal{M}(z)$ stated in the lemma and the corresponding 
items~$\pmb{{\rm (i)}}$ and~$\pmb{{\rm (iii)}}$ of Lemma~\ref{lem4.1}; 
it remains to verify the corresponding item~(ii) of case~$\pmb{(2)}$, that is, 
the expression for the associated jump matrix, $\mathfrak{v}(z)$. Recall 
{}from the corresponding item~$\pmb{{\rm (ii)}}$ of Lemma~\ref{lem4.1} 
that, for $n \! \in \! \mathbb{N}$ and $k \! \in \! \lbrace 1,2,\dotsc,K 
\rbrace$ such that $\alpha_{p_{\mathfrak{s}}} \! := \! \alpha_{k} \! \neq 
\! \infty$, and $z \! \in \! (\tilde{b}_{j-1},\tilde{a}_{j})$, $j \! = \! 1,2,
\dotsc,N \! + \! 1$, $\mathcal{M}_{+}(z) \! = \! \mathcal{M}_{-}(z) 
\hat{\mathscr{V}}(z)$, where $\hat{\mathscr{V}}(z) \! = \! \sigma_{+} \! 
+ \! \exp (-(2 \pi \mi ((n \! - \! 1)K \! + \! k) \int_{z}^{\tilde{a}_{N+1}} 
\psi_{\widetilde{V}}^{f}(\xi) \, \md \xi) \sigma_{3})$: noting the matrix 
factorisation
\begin{equation*}
\sigma_{+} \! + \! \me^{-(2 \pi \mi ((n-1)K+k) \int_{z}^{\tilde{a}_{N+1}} 
\psi_{\widetilde{V}}^{f}(\xi) \, \md \xi) \sigma_{3}} = 
\begin{pmatrix}
1 & 0 \\
\me^{2 \pi \mi ((n-1)K+k) \int_{z}^{\tilde{a}_{N+1}} \psi_{\widetilde{V}}^{f}
(\xi) \, \md \xi} & 1
\end{pmatrix} 
\begin{pmatrix}
0 & 1 \\
-1 & 0
\end{pmatrix} 
\begin{pmatrix}
1 & 0 \\
\me^{-2 \pi \mi ((n-1)K+k) \int_{z}^{\tilde{a}_{N+1}} \psi_{\widetilde{V}}^{f}
(\xi) \, \md \xi} & 1
\end{pmatrix},
\end{equation*}
it follows that, for $z \! \in \! (\tilde{b}_{j-1},\tilde{a}_{j})$, $j \! = \! 
1,2,\dotsc,N \! + \! 1$,
\begin{equation*}
\mathcal{M}_{+}(z) \left(\mathrm{I} \! - \! \me^{-2 \pi \mi ((n-1)K+k) 
\int_{z}^{\tilde{a}_{N+1}} \psi_{\widetilde{V}}^{f}(\xi) \, \md \xi} 
\sigma_{-} \right) \! = \! \mathcal{M}_{-}(z) \left(\mathrm{I} \! + \! 
\me^{2 \pi \mi ((n-1)K+k) \int_{z}^{\tilde{a}_{N+1}} \psi_{\widetilde{V}}^{f}
(\xi) \, \md \xi} \sigma_{-} \right) \mi \sigma_{2}.
\end{equation*}
It was shown in the corresponding proof of Lemma~\ref{lem4.1} that, for 
$n \! \in \! \mathbb{N}$ and $k \! \in \! \lbrace 1,2,\dotsc,K \rbrace$ such 
that $\alpha_{p_{\mathfrak{s}}} \! := \! \alpha_{k} \! \neq \! \infty$, $\pm 
\Re (\mi \int_{z}^{\tilde{a}_{N+1}} \psi_{\widetilde{V}}^{f}(\xi) \, \md \xi) 
\! > \! 0$, $z \! \in \! \mathbb{C}_{\pm} \cap (\cup_{j=1}^{N+1} 
\tilde{\mathbb{U}}_{j})$: the terms $\pm 2 \pi \mi ((n \! - \! 1)K \! + \! k) 
\int_{z}^{\tilde{a}_{N+1}} \psi_{\widetilde{V}}^{f}(\xi) \, \md \xi$, which 
are pure imaginary for $z \! \in \! \mathbb{R}$, and corresponding 
to which $\exp (\pm 2 \pi \mi ((n \! - \! 1)K \! + \! k) 
\int_{z}^{\tilde{a}_{N+1}} \psi_{\widetilde{V}}^{f}(\xi) \, \md \xi)$ 
are oscillatory, are continued analytically to $\mathbb{C}_{\pm} \cap 
(\cup_{j=1}^{N+1} \tilde{\mathbb{U}}_{j})$, respectively, corresponding to 
which $\exp (\mp 2 \pi \mi ((n \! - \! 1)K \! + \! k) \int_{z}^{\tilde{a}_{N+1}} 
\psi_{\widetilde{V}}^{f}(\xi) \, \md \xi)$ are exponentially decreasing (in 
the double-scaling limit $\mathscr{N},n \! \to \! \infty$ such that $z_{o} 
\! = \! 1 \! + \! o(1))$. As per the Deift-Zhou non-linear steepest-descent 
method \cite{a52,a53} (see, also, \cite{a54}), one now `deforms' the 
original, oriented contour $\mathbb{R}$ to the extended, oriented 
contour/skeleton $\tilde{\Sigma} \! := \! \mathbb{R} \cup \cup_{j=1}^{N+1}
(\tilde{J}_{j}^{\smallfrown} \cup \tilde{J}_{j}^{\smallsmile})$ shown in 
Figure~\ref{figgtwo} in such a way that the upper (resp., lower) `lips' of the 
`lenses' $\tilde{J}_{j}^{\smallfrown}$ (resp., $\tilde{J}_{j}^{\smallsmile})$, 
$j \! = \! 1,2,\dotsc,N \! + \! 1$, which are the boundaries of 
$\tilde{\Omega}_{j}^{\smallfrown}$ (resp., $\tilde{\Omega}_{j}^{\smallsmile})$, 
lie in the domain of analytic continuation of $g^{f}_{+}(z) \! - \! g^{f}_{-}(z) 
\! + \! \hat{\mathscr{P}}_{0}^{-} \! - \! \hat{\mathscr{P}}_{0}^{+}$, namely, 
$\cup_{j=1}^{N+1}(\tilde{\Omega}_{j}^{\smallfrown} \cup \tilde{\Omega}_{j}^{
\smallsmile} \cup \tilde{J}_{j}^{\smallfrown} \cup \tilde{J}_{j}^{\smallsmile})$ 
$(\subset \cup_{j=1}^{N+1} \tilde{\mathbb{U}}_{j})$; more precisely, each 
oriented and bounded open interval $(\tilde{b}_{j-1},\tilde{a}_{j})$, $j \! = 
\! 1,2,\dotsc,N \! + \! 1$, in the original, oriented contour $\mathbb{R}$ 
is `split', or `branched', into three, and the resulting oriented contour 
$\tilde{\Sigma}$ is the old contour, $\mathbb{R}$, together with the oriented 
boundary of $N \! + \! 1$ lens-shaped regions, one region surrounding each 
bounded and oriented open interval $(\tilde{b}_{j-1},\tilde{a}_{j})$, $j \! = 
\! 1,2,\dotsc,N \! + \! 1$. Recalling, now, the corresponding definition of 
$\mathfrak{M}(z)$ in terms of $\mathcal{M}(z)$ stated in the formulation of 
the lemma, and the associated expression for $\hat{\mathscr{V}}(z)$ given in 
the corresponding item~$\pmb{{\rm (ii)}}$ of Lemma~\ref{lem4.1}, one arrives 
at, for $n \! \in \! \mathbb{N}$ and $k \! \in \! \lbrace 1,2,\dotsc,K \rbrace$ 
such that $\alpha_{p_{\mathfrak{s}}} \! := \! \alpha_{k} \! \neq \! \infty$, 
the expression for $\mathfrak{v}(z)$ given in the corresponding item~(ii) 
of case~$\pmb{(2)}$ of the lemma.

$\pmb{(\mathrm{B})}$ The proof of this case, that is, $n \! \in \! 
\mathbb{N}$ and $k \! \in \! \lbrace 1,2,\dotsc,K \rbrace$ such that 
$\alpha_{p_{\mathfrak{s}}} \! := \! \alpha_{k} \! = \! \infty$, is virtually 
identical to the proof presented in $\pmb{(\mathrm{A})}$ above; one 
mimics, \emph{verbatim}, the scheme of the calculations presented in 
case $\pmb{(\mathrm{A})}$ in order to arrive at the claims stated in the 
corresponding items~(i)--(iii) of case~$\pmb{(2)}$ of the lemma; in 
particular, all that is necessary in order to complete the proof is the 
associated definition of $\mathfrak{M}(z)$ in terms of $\mathcal{M}(z)$ 
stated in the formulation of the lemma and the following matrix factorisation: 
for $n \! \in \! \mathbb{N}$ and $k \! \in \! \lbrace 1,2,\dotsc,K \rbrace$ 
such that $\alpha_{p_{\mathfrak{s}}} \! := \! \alpha_{k} \! = \! \infty$, and 
$z \! \in \! (\hat{b}_{j-1},\hat{a}_{j})$, $j \! = \! 1,2,\dotsc,N \! + \! 1$,
\begin{equation*}
\sigma_{+} \! + \! \me^{-(2 \pi \mi ((n-1)K+k) \int_{z}^{\hat{a}_{N+1}} 
\psi_{\widetilde{V}}^{\infty}(\xi) \, \md \xi) \sigma_{3}} = 
\begin{pmatrix}
1 & 0 \\
\me^{2 \pi \mi ((n-1)K+k) \int_{z}^{\hat{a}_{N+1}} 
\psi_{\widetilde{V}}^{\infty}(\xi) \, \md \xi} & 1
\end{pmatrix} 
\begin{pmatrix}
0 & 1 \\
-1 & 0
\end{pmatrix} 
\begin{pmatrix}
1 & 0 \\
\me^{-2 \pi \mi ((n-1)K+k) \int_{z}^{\hat{a}_{N+1}} 
\psi_{\widetilde{V}}^{\infty}(\xi) \, \md \xi} & 1
\end{pmatrix}.
\end{equation*}
This concludes the proof. \hfill $\qed$
\begin{eeeee} \label{ref4.2} 
\textsl{Recalling {}from Lemma~\ref{lem4.1} that, for $n \! \in \! 
\mathbb{N}$ and $k \! \in \! \lbrace 1,2,\dotsc,K \rbrace$ such 
that $\alpha_{p_{\mathfrak{s}}} \! := \! \alpha_{k} \! = \! \infty$ 
(resp., $\alpha_{p_{\mathfrak{s}}} \! := \! \alpha_{k} \! \neq \! 
\infty)$, $g_{+}^{\infty}(z) \! + \! g_{-}^{\infty}(z) \! - \! 2 
\tilde{\mathscr{P}}_{0} \! - \! \widetilde{V}(z) \! - \! \hat{\ell} \! 
< \! 0$, $z \! \in \! (-\infty,\hat{b}_{0}) \cup \cup_{i=1}^{N}
(\hat{a}_{i},\hat{b}_{i}) \cup (\hat{a}_{N+1},+\infty)$ (resp., 
$g_{+}^{f}(z) \! + \! g_{-}^{f}(z) \! - \! \hat{\mathscr{P}}_{0}^{+} 
\! - \! \hat{\mathscr{P}}_{0}^{-} \! - \! \widetilde{V}(z) \! - \! 
\tilde{\ell} \! < \! 0$, $z \! \in \! (-\infty,\tilde{b}_{0}) \cup 
\cup_{i=1}^{N}(\tilde{a}_{i},\tilde{b}_{i}) \cup (\tilde{a}_{N+1},
+\infty))$, and $\pm \Re (\mi \int_{z}^{\hat{a}_{N+1}} 
\psi_{\widetilde{V}}^{\infty}(\xi) \, \md \xi) \! > \! 0$, $z \! 
\in \! \hat{J}_{j}^{\mathfrak{q}(\pm)}$ (resp., $\pm \Re (\mi 
\int_{z}^{\tilde{a}_{N+1}} \psi_{\widetilde{V}}^{f}(\xi) \, \md \xi) 
\! > \! 0$, $z \! \in \! \tilde{J}_{j}^{\mathfrak{q}(\pm)})$, $j \! 
= \! 1,2,\dotsc,N \! + \! 1$, $\mathfrak{q}(+) \! := \smallfrown$, 
$\mathfrak{q}(-) \! := \smallsmile$, one arrives at (see 
Section~\ref{sek5}, Lemma~\ref{lem5.1}$)$ the 
following---coarse---asymptotic behaviours (in the double-scaling 
limit $\mathscr{N},n \! \to \! \infty$ such that $z_{o} \! = \! 1 \! 
+ \! o(1))$ for the associated jump matrices (cf. items~$\pmb{(1)}$ 
and~$\pmb{(2)}$ of Lemma~\ref{lem4.2}$;$ see, also, 
Remark~\ref{rem4.1}$)$$:$ for $n \! \in \! \mathbb{N}$ and $k \! \in 
\! \lbrace 1,2,\dotsc,K \rbrace$ such that $\alpha_{p_{\mathfrak{s}}} 
\! := \! \alpha_{k} \! = \! \infty$, $j \! = \! 1,2,\dotsc,N \! + \! 1$, $i 
\! = \! 1,2,\dotsc,N$, and $q \! \in \! \lbrace 1,2,\dotsc,\mathfrak{s} 
\! - \! 1 \rbrace$,
\begin{equation*}
\mathfrak{v}(z) \underset{\underset{z_{o}=1+o(1)}{\mathscr{N},n 
\to \infty}}{=} 
\begin{cases}
\mi \sigma_{2}, &\text{$z \! \in \! (\hat{b}_{j-1},\hat{a}_{j})$,} \\
\mathrm{I} \! + \! \mathcal{O}(\me^{-((n-1)K+k) \hat{c} \lvert 
z-\hat{P} \rvert}) \sigma_{-}, &\text{$z \! \in \! \hat{J}_{j}^{
\smallfrown} \cup \hat{J}_{j}^{\smallsmile}$,} \\
\me^{-(2 \pi \mi ((n-1)K+k) \int_{\hat{b}_{i}}^{\hat{a}_{N+1}} 
\psi_{\widetilde{V}}^{\infty}(\xi) \, \md \xi) \sigma_{3}} \left(
\mathrm{I} \! + \! \mathcal{O}(\me^{-((n-1)K+k) \hat{c}(z-
\hat{a}_{i})}) \sigma_{+} \right), &\text{$z \! \in \! (\hat{a}_{i},
\hat{b}_{i}) \setminus \mathscr{O}_{\hat{\delta}_{p_{q}}}
(\alpha_{p_{q}})$,} \\
\me^{-(2 \pi \mi ((n-1)K+k) \int_{\hat{b}_{i}}^{\hat{a}_{N+1}} 
\psi_{\widetilde{V}}^{\infty}(\xi) \, \md \xi) \sigma_{3}} \left(
\mathrm{I} \! + \! \mathcal{O}(\lvert z \! - \! \alpha_{p_{q}} 
\rvert^{((n-1)K+k) \hat{c}}) \sigma_{+} \right), &\text{$z \! 
\in \! (\hat{a}_{i},\hat{b}_{i}) \cap \mathscr{O}_{\hat{\delta}_{p_{q}}}
(\alpha_{p_{q}})$,} \\
\mathrm{I} \! + \! \mathcal{O}(\me^{-((n-1)K+k) \hat{c}(z-
\hat{a}_{N+1})}) \sigma_{+}, &\text{$z \! \in \! (\hat{a}_{N+1},
+\infty) \setminus \mathscr{O}_{\hat{\delta}_{p_{q}}}
(\alpha_{p_{q}})$,} \\
\mathrm{I} \! + \! \mathcal{O}(\me^{-((n-1)K+k) \hat{c} 
\lvert z-\hat{b}_{0} \rvert}) \sigma_{+}, &\text{$z \! \in \! 
(-\infty,\hat{b}_{0}) \setminus \mathscr{O}_{\hat{\delta}_{p_{q}}}
(\alpha_{p_{q}})$,} \\
\mathrm{I} \! + \! \mathcal{O}(\lvert z \! - \! \alpha_{p_{q}} 
\rvert^{((n-1)K+k) \hat{c}}) \sigma_{+}, &\text{$z \! \in \! 
((-\infty,\hat{b}_{0}) \cup (\hat{a}_{N+1},+\infty)) \cap 
\mathscr{O}_{\hat{\delta}_{p_{q}}}(\alpha_{p_{q}})$,} \\
\mathrm{I} \! + \! \mathcal{O}(\me^{-((n-1)K+k) \hat{c} 
\ln \lvert z \rvert}) \sigma_{+}, &\text{$z \! \in \! 
((-\infty,\hat{b}_{0}) \cup (\hat{a}_{N+1},+\infty)) \cap 
\mathscr{O}_{\infty}$,}
\end{cases}
\end{equation*}
where $\hat{c} \! > \! 0$ is some generic constant (all of whose 
explicit dependencies are suppressed for notational simplicity), 
$\hat{P} \! = \! (\hat{a}_{j} \! + \! \hat{\varepsilon} \cos 
\hat{\theta},\hat{\varepsilon} \sin \hat{\theta})$, $\hat{\theta} 
\! \in \! (-\pi,-2 \pi/3) \cup (2 \pi/3,\pi)$, $\mathscr{O}_{
\hat{\delta}_{p_{q}}}(\alpha_{p_{q}}) \! := \! \lbrace \mathstrut z 
\! \in \! \mathbb{C}; \, \lvert z \! - \! \alpha_{p_{q}} \rvert  \! 
< \! \hat{\delta}_{p_{q}} \rbrace$, with $\hat{\varepsilon}$ and 
$\hat{\delta}_{p_{q}}$ some arbitrarily fixed, sufficiently small 
positive real numbers chosen so that, for $i_{1} \! \in \! \lbrace 
1,2,\dotsc,N \! + \! 1 \rbrace$ and $i_{2},i_{3} \! \in \! \lbrace 
1,2,\dotsc,\mathfrak{s} \! - \! 1 \rbrace$, $\mathscr{O}_{
\hat{\varepsilon}}(\hat{b}_{i_{1}-1}) \cap \mathscr{O}_{
\hat{\varepsilon}}(\hat{a}_{i_{1}}) \! = \! \mathscr{O}_{
\hat{\varepsilon}}(\hat{b}_{i_{1}-1}) \cap \mathscr{O}_{
\hat{\delta}_{p_{i_{2}}}}(\alpha_{p_{i_{2}}}) \! = \! 
\mathscr{O}_{\hat{\varepsilon}}(\hat{a}_{i_{1}}) \cap 
\mathscr{O}_{\hat{\delta}_{p_{i_{2}}}}(\alpha_{p_{i_{2}}}) \! = \! 
\varnothing$, and $\mathscr{O}_{\hat{\delta}_{p_{i_{2}}}}
(\alpha_{p_{i_{2}}}) \cap \mathscr{O}_{\hat{\delta}_{p_{i_{3}}}}
(\alpha_{p_{i_{3}}}) \! = \! \varnothing$ $\forall$ $i_{2} 
\neq \! i_{3} $, and $\mathscr{O}_{\infty}$ denotes 
the open neighbourhood of the point at infinity$;$ and, 
for $n \! \in \! \mathbb{N}$ and $k \! \in \! \lbrace 1,2,\dotsc,
K \rbrace$ such that $\alpha_{p_{\mathfrak{s}}} \! := \! \alpha_{k} 
\! \neq \! \infty$, $j \! = \! 1,2,\dotsc,N \! + \! 1$, $i \! = \! 1,2,
\dotsc,N$, and $q \! \in \! \lbrace 1,\dotsc,\mathfrak{s} \! - \! 2,
\mathfrak{s} \rbrace$,
\begin{equation*}
\mathfrak{v}(z) \underset{\underset{z_{o}=1+o(1)}{\mathscr{N},n 
\to \infty}}{=} 
\begin{cases}
\mi \sigma_{2}, &\text{$z \! \in \! (\tilde{b}_{j-1},\tilde{a}_{j})$,} \\
\mathrm{I} \! + \! \mathcal{O}(\me^{-((n-1)K+k) \tilde{c} \lvert 
z-\tilde{P} \rvert}) \sigma_{-}, &\text{$z \! \in \! \tilde{J}_{j}^{
\smallfrown} \cup \tilde{J}_{j}^{\smallsmile}$,} \\
\me^{-(2 \pi \mi ((n-1)K+k) \int_{\tilde{b}_{i}}^{\tilde{a}_{N+1}} 
\psi_{\widetilde{V}}^{f}(\xi) \, \md \xi) \sigma_{3}} \left(
\mathrm{I} \! + \! \mathcal{O}(\me^{-((n-1)K+k) \tilde{c}(z-
\tilde{a}_{i})}) \sigma_{+} \right), &\text{$z \! \in \! (\tilde{a}_{i},
\tilde{b}_{i}) \setminus \mathscr{O}_{\tilde{\delta}_{p_{q}}}
(\alpha_{p_{q}})$,} \\
\me^{-(2 \pi \mi ((n-1)K+k) \int_{\tilde{b}_{i}}^{\tilde{a}_{N+1}} 
\psi_{\widetilde{V}}^{f}(\xi) \, \md \xi) \sigma_{3}} \left(
\mathrm{I} \! + \! \mathcal{O}(\lvert z \! - \! \alpha_{p_{q}} 
\rvert^{((n-1)K+k) \tilde{c}}) \sigma_{+} \right), &\text{$z \! \in 
\! (\tilde{a}_{i},\tilde{b}_{i}) \cap \mathscr{O}_{\tilde{\delta}_{p_{q}}}
(\alpha_{p_{q}})$,} \\
\mathrm{I} \! + \! \mathcal{O}(\me^{-((n-1)K+k) \tilde{c}(z-
\tilde{a}_{N+1})}) \sigma_{+}, &\text{$z \! \in \! (\tilde{a}_{N+1},
+\infty) \setminus \mathscr{O}_{\tilde{\delta}_{p_{q}}}
(\alpha_{p_{q}})$,} \\
\mathrm{I} \! + \! \mathcal{O}(\me^{-((n-1)K+k) \tilde{c} 
\lvert z-\tilde{b}_{0} \rvert}) \sigma_{+}, &\text{$z \! \in \! 
(-\infty,\tilde{b}_{0}) \setminus \mathscr{O}_{\tilde{\delta}_{p_{q}}}
(\alpha_{p_{q}})$,} \\
\mathrm{I} \! + \! \mathcal{O}(\lvert z \! - \! \alpha_{p_{q}} 
\rvert^{((n-1)K+k) \tilde{c}}) \sigma_{+}, &\text{$z \! \in \! 
((-\infty,\tilde{b}_{0}) \cup (\tilde{a}_{N+1},+\infty)) \cap 
\mathscr{O}_{\tilde{\delta}_{p_{q}}}(\alpha_{p_{q}})$,} \\
\mathrm{I} \! + \! \mathcal{O}(\me^{-((n-1)K+k) \tilde{c} 
\ln \lvert z \rvert}) \sigma_{+}, &\text{$z \! \in \! 
((-\infty,\tilde{b}_{0}) \cup (\tilde{a}_{N+1},+\infty)) \cap 
\mathscr{O}_{\infty}$,}
\end{cases}
\end{equation*}
where $\tilde{c} \! > \! 0$ is some generic constant (all of whose 
explicit dependencies are suppressed for notational simplicity), 
$\tilde{P} \! = \! (\tilde{a}_{j} \! + \! \tilde{\varepsilon} \cos 
\tilde{\theta},\tilde{\varepsilon} \sin \tilde{\theta})$, $\tilde{\theta} 
\! \in \! (-\pi,-2 \pi/3) \cup (2 \pi/3,\pi)$, with $\tilde{\varepsilon}$ 
and $\tilde{\delta}_{p_{q}}$ some arbitrarily fixed, sufficiently small 
positive real numbers chosen so that, for $i_{1} \! \in \! \lbrace 
1,2,\dotsc,N \! + \! 1 \rbrace$ and $i_{2},i_{3} \! \in \! \lbrace 
1,\dotsc,\mathfrak{s} \! - \! 2,\mathfrak{s} \rbrace$, 
$\mathscr{O}_{\tilde{\varepsilon}}(\tilde{b}_{i_{1}-1}) \cap 
\mathscr{O}_{\tilde{\varepsilon}}(\tilde{a}_{i_{1}}) \! = \! 
\mathscr{O}_{\tilde{\varepsilon}}(\tilde{b}_{i_{1}-1}) \cap 
\mathscr{O}_{\tilde{\delta}_{p_{i_{2}}}}(\alpha_{p_{i_{2}}}) \! 
= \! \mathscr{O}_{\tilde{\varepsilon}}(\tilde{a}_{i_{1}}) \cap 
\mathscr{O}_{\tilde{\delta}_{p_{i_{2}}}}(\alpha_{p_{i_{2}}}) \! = 
\! \varnothing$, and $\mathscr{O}_{\tilde{\delta}_{p_{i_{2}}}}
(\alpha_{p_{i_{2}}}) \cap \mathscr{O}_{\tilde{\delta}_{p_{i_{3}}}}
(\alpha_{p_{i_{3}}}) \! = \! \varnothing$ $\forall$ $i_{2} \neq \! 
i_{3} $. All of the above-mentioned associated convergences are 
uniform in the respective compact subsets.}
\end{eeeee}
Recall {}from Lemma~2.56 of \cite{a52} that, for an oriented 
skeleton in $\mathbb{C}$ on which the jump matrix of a matrix 
RHP is defined, one may always choose to add or delete a portion 
of the skeleton on which the jump matrix equals $\mathrm{I}$ 
without altering the matrix RHP in the operator-theoretic sense. For 
$n \! \in \! \mathbb{N}$ and $k \! \in \! \lbrace 1,2,\dotsc,K \rbrace$ 
such that $\alpha_{p_{\mathfrak{s}}} \! := \! \alpha_{k} \! = \! \infty$ 
(resp., $\alpha_{p_{\mathfrak{s}}} \! := \! \alpha_{k} \! \neq \! \infty)$, 
neglecting those jumps on the associated oriented skeleton $\hat{
\Sigma}$ (resp., $\tilde{\Sigma}$) tending, in the double-scaling 
limit $\mathscr{N},n \! \to \! \infty$ such that $z_{o} \! = \! 1 \! + 
\! o(1)$, to $\mathrm{I}$, and deleting the corresponding oriented 
sub-skeletons from $\hat{\Sigma}$ (resp., $\tilde{\Sigma}$), it becomes 
more or less transparent how to construct associated parametrices 
(approximate solutions) for the matrix RHP $(\mathfrak{M}(z),
\mathfrak{v}(z),\hat{\Sigma})$ (resp., $(\mathfrak{M}(z),\mathfrak{v}(z),
\tilde{\Sigma}))$ formulated in Lemma~\ref{lem4.2}; in particular, for 
$n \! \in \! \mathbb{N}$ and $k \! \in \! \lbrace 1,2,\dotsc,K \rbrace$ 
such that $\alpha_{p_{\mathfrak{s}}} \! := \! \alpha_{k} \! = \! \infty$ 
(resp., $\alpha_{p_{\mathfrak{s}}} \! := \! \alpha_{k} \! \neq \! \infty$), 
the solution, in the double-scaling limit $\mathscr{N},n \! \to \! \infty$ 
such that $z_{o} \! = \! 1 \! + \! o(1)$, of the matrix RHP $\mathfrak{M} 
\colon \mathbb{N} \times \lbrace 1,2,\dotsc,K \rbrace \times \mathbb{C} 
\setminus \hat{\Sigma} \! \to \! \operatorname{SL}_{2}(\mathbb{C})$ (resp., 
$\mathfrak{M} \colon \mathbb{N} \times \lbrace 1,2,\dotsc,K \rbrace \times 
\mathbb{C} \setminus \tilde{\Sigma} \! \to \! \operatorname{SL}_{2}(\mathbb{C}))$ 
formulated in Lemma~\ref{lem4.2} should be `close to', in some appropriately 
defined operator-theoretic sense, the solution of the following family of $\hat{K} 
\! := \! \# \lbrace \mathstrut j \! \in \! \lbrace 1,2,\dotsc,K \rbrace; \, \alpha_{j} 
\! = \! \infty \rbrace$ (resp., $\tilde{K} \! := \! \# \lbrace \mathstrut j \! \in \! 
\lbrace 1,2,\dotsc,K \rbrace; \, \alpha_{j} \! \neq \! \infty \rbrace$) `model' 
matrix RHPs.
\begin{ccccc} \label{lem4.3} 
Let the external field $\widetilde{V} \colon \overline{\mathbb{R}} \setminus 
\lbrace \alpha_{1},\alpha_{2},\dotsc,\alpha_{K} \rbrace \! \to \! \mathbb{R}$ 
satisfy conditions~\eqref{eq20}--\eqref{eq22} and be regular. For $n \! \in 
\! \mathbb{N}$ and $k \! \in \! \lbrace 1,2,\dotsc,K \rbrace$ such that 
$\alpha_{p_{\mathfrak{s}}} \! := \! \alpha_{k} \! = \! \infty$ (resp., 
$\alpha_{p_{\mathfrak{s}}} \! := \! \alpha_{k} \! \neq \! \infty)$, let the 
associated equilibrium measure, $\mu_{\widetilde{V}}^{\infty}$ (resp., 
$\mu_{\widetilde{V}}^{f})$, and its support, $J_{\infty}$ (resp., $J_{f})$, 
be as described in item~$\pmb{(1)}$ (resp., item~$\pmb{(2)})$ of 
Lemma~\ref{lem3.7}, and, along with the corresponding variational 
constant, $\hat{\ell}$ (resp., $\tilde{\ell})$, satisfy the variational 
conditions~\eqref{eql3.8a} (resp., conditions~\eqref{eql3.8b}$)$$;$ 
moreover, let the associated conditions~{\rm (i)}--{\rm (iv)} of 
item~$\pmb{(1)}$ (resp., item~$\pmb{(2)})$ of Lemma~\ref{lem3.8} be 
valid. Then$:$ $\pmb{(1)}$ for $n \! \in \! \mathbb{N}$ and $k \! \in 
\! \lbrace 1,2,\dotsc,K \rbrace$ such that $\alpha_{p_{\mathfrak{s}}} 
\! := \! \alpha_{k} \! = \! \infty$, $\mathfrak{m} \colon \mathbb{N} 
\times \lbrace 1,2,\dotsc,K \rbrace \times \mathbb{C} \setminus 
\hat{\mathscr{J}} \! \to \! \operatorname{SL}_{2}(\mathbb{C})$, where 
$\hat{\mathscr{J}} \! := \! J_{\infty} \cup \cup_{i=1}^{N}(\hat{a}_{i},
\hat{b}_{i})$, solves the following matrix {\rm RHP}$:$ {\rm (i)} 
$\mathfrak{m}(n,k,z) \! = \! \mathfrak{m}(z)$ is analytic for $z \! \in 
\! \mathbb{C} \setminus \hat{\mathscr{J}}$$;$ {\rm (ii)} the boundary 
values $\mathfrak{m}_{\pm}(z) \! := \! \lim_{\underset{z^{\prime} \in 
\pm \, \mathrm{side \, of} \, \hat{\mathscr{J}}}{z^{\prime} \to z \in 
\hat{\mathscr{J}}}} \mathfrak{m}(z^{\prime})$ satisfy the jump condition
\begin{equation*}
\mathfrak{m}_{+}(z) \! = \! \mathfrak{m}_{-}(z) \daleth (z) \quad 
\mathrm{a.e.} \quad z \! \in \! \hat{\mathscr{J}},
\end{equation*}
where
\begin{align*}
\daleth &\colon \mathbb{N} \times \lbrace 1,2,\dotsc,K \rbrace \times 
\hat{\mathscr{J}} \! \to \! \operatorname{GL}_{2}(\mathbb{C}), \, (n,k,z) 
\! \mapsto \! \daleth (n,k,z) \! =: \! \daleth (z) \\
&= 
\begin{cases}
\mi \sigma_{2}, &\text{$z\! \in \! (\hat{b}_{j-1},\hat{a}_{j}), 
\quad j \! = \! 1,2,\dotsc,N \! + \! 1$,} \\
\me^{-(2 \pi \mi ((n-1)K+k) \int_{\hat{b}_{i}}^{\hat{a}_{N+1}} 
\psi_{\widetilde{V}}^{\infty}(\xi) \, \md \xi) \sigma_{3}}, 
&\text{$z \! \in \! (\hat{a}_{i},\hat{b}_{i}), \quad i \! = \! 
1,2,\dotsc,N;$}
\end{cases}
\end{align*}
{\rm (iii)}
\begin{equation*}
\mathfrak{m}(z) \underset{\overline{\mathbb{C}}_{\pm} \ni z \to \alpha_{k}}{=} 
\mathrm{I} \! + \! \mathcal{O}(z^{-1}), \qquad \qquad \mathfrak{m}(z) 
\underset{\mathbb{C}_{\pm} \ni z \to \alpha_{p_{q}}}{=} \mathcal{O}
(\me^{n(\tilde{\mathscr{P}}_{0}-\tilde{\mathscr{P}}^{\pm}_{1}) \sigma_{3}}), 
\quad q \! = \! 1,2,\dotsc,\mathfrak{s} \! - \! 1;
\end{equation*}
and~$\pmb{(2)}$ for $n \! \in \! \mathbb{N}$ and $k \! \in \! \lbrace 1,2,
\dotsc,K \rbrace$ such that $\alpha_{p_{\mathfrak{s}}} \! := \! \alpha_{k} 
\! \neq \! \infty$, $\mathfrak{m} \colon \mathbb{N} \times \lbrace 1,2,
\dotsc,K \rbrace \times \mathbb{C} \setminus \tilde{\mathscr{J}} \! \to \! 
\operatorname{SL}_{2}(\mathbb{C})$, where $\tilde{\mathscr{J}} \! := \! 
J_{f} \cup \cup_{i=1}^{N}(\tilde{a}_{i},\tilde{b}_{i})$, solves the following 
matrix {\rm RHP}$:$ {\rm (i)} $\mathfrak{m}(n,k,z) \! = \! \mathfrak{m}(z)$ 
is analytic for $z \! \in \! \mathbb{C} \setminus \tilde{\mathscr{J}}$$;$ 
{\rm (ii)} the boundary values $\mathfrak{m}_{\pm}(z) \! := \! 
\lim_{\underset{z^{\prime} \in \pm \, \mathrm{side \, of} \, 
\tilde{\mathscr{J}}}{z^{\prime} \to z \in \tilde{\mathscr{J}}}} 
\mathfrak{m}(z^{\prime})$ satisfy the jump condition
\begin{equation*}
\mathfrak{m}_{+}(z) \! = \! \mathfrak{m}_{-}(z) \daleth (z) \quad 
\mathrm{a.e.} \quad z \! \in \! \tilde{\mathscr{J}},
\end{equation*}
where
\begin{align*}
\daleth &\colon \mathbb{N} \times \lbrace 1,2,\dotsc,K \rbrace \times 
\tilde{\mathscr{J}} \! \to \! \operatorname{GL}_{2}(\mathbb{C}), \, 
(n,k,z) \! \mapsto \! \daleth (n,k,z) \! =: \! \daleth (z) \\
&= 
\begin{cases}
\mi \sigma_{2}, &\text{$z\! \in \! (\tilde{b}_{j-1},\tilde{a}_{j}), 
\quad j \! = \! 1,2,\dotsc,N \! + \! 1$,} \\
\me^{-(2 \pi \mi ((n-1)K+k) \int_{\tilde{b}_{i}}^{\tilde{a}_{N+1}} 
\psi_{\widetilde{V}}^{f}(\xi) \, \md \xi) \sigma_{3}}, 
&\text{$z \! \in \! (\tilde{a}_{i},\tilde{b}_{i}), \quad i \! = \! 
1,2,\dotsc,N;$}
\end{cases}
\end{align*}
{\rm (iii)}
\begin{gather*}
\mathfrak{m}(z) \underset{\mathbb{C}_{\pm} \ni z \to \alpha_{k}}{=} 
(\mathrm{I} \! + \! \mathcal{O}(z \! - \! \alpha_{k})) \mathscr{E}^{\mp 
\sigma_{3}}, \quad \qquad \mathfrak{m}(z) \underset{\overline{
\mathbb{C}}_{\pm} \ni z \to \alpha_{p_{\mathfrak{s}-1}} = \infty}{=} 
\mathcal{O}(\me^{n(\hat{\mathscr{P}}^{\pm}_{0}-\hat{\mathscr{P}}_{1}) 
\sigma_{3}} \mathscr{E}^{\mp \sigma_{3}}), \\
\mathfrak{m}(z) \underset{\mathbb{C}_{\pm} \ni z \to \alpha_{p_{q}}}{=} 
\mathcal{O}(\me^{n(\hat{\mathscr{P}}^{\pm}_{0}-\hat{\mathscr{P}}^{
\pm}_{2}) \sigma_{3}} \mathscr{E}^{\mp \sigma_{3}}), \quad 
q \! = \! 1,2,\dotsc,\mathfrak{s} \! - \! 2.
\end{gather*}
\end{ccccc}
It turns out that, for $n \! \in \! \mathbb{N}$ and $k \! \in \! \lbrace 
1,2,\dotsc,K \rbrace$ such that $\alpha_{p_{\mathfrak{s}}} \! := \! 
\alpha_{k} \! = \! \infty$ (resp., $\alpha_{p_{\mathfrak{s}}} \! := \! 
\alpha_{k} \! \neq \! \infty$), the model matrix RHP $(\mathfrak{m}
(z),\daleth (z),\hat{\mathscr{J}})$ (resp., $(\mathfrak{m}(z),\daleth 
(z),\tilde{\mathscr{J}}))$ formulated in item~$\pmb{(1)}$ (resp., 
item~$\pmb{(2)}$) of Lemma~\ref{lem4.3} is explicitly solvable in 
terms of Riemann theta functions associated with the family of $\hat{K}$ 
(resp., $\tilde{K}$) two-sheeted genus-$N$ hyperelliptic (compact) 
Riemann surfaces $\hat{\mathcal{Y}}$ (resp., $\tilde{\mathcal{Y}})$; 
see, for example, Section~3 of \cite{a67} (see, also, \cite{a60}, and 
Section~4.2 of \cite{a59}): the associated families of solutions and 
their corresponding parametrices are now presented.
\begin{ccccc} \label{lem4.4} 
For $n \! \in \! \mathbb{N}$ and $k \! \in \! \lbrace 1,2,\dotsc,K 
\rbrace$ such that $\alpha_{p_{\mathfrak{s}}} \! := \! \alpha_{k} 
\! = \! \infty$, let $\hat{\gamma} \colon \mathbb{C} \setminus 
((-\infty,\hat{b}_{0}) \cup (\hat{a}_{N+1},+\infty) \cup \cup_{j=1}^{N}
(\hat{a}_{j},\hat{b}_{j})) \! \to \! \mathbb{C}$ be defined by
\begin{equation} \label{eqmaininf10} 
\hat{\gamma}(z) \! := \!
\begin{cases} 
\left(\left(\dfrac{z \! - \! \hat{b}_{0}}{z \! - \! \hat{a}_{N+1}} \right) 
\mathlarger{\prod_{m=1}^{N}} \left(\dfrac{z \! - \! \hat{b}_{m}}{z 
\! - \! \hat{a}_{m}} \right) \right)^{1/4}, &\text{$z \! \in \! 
\mathbb{C}_{+}$ $(\subset \hat{\mathcal{Y}}^{+})$,} \\
-\mi \left(\left(\dfrac{z \! - \! \hat{b}_{0}}{z \! - \! \hat{a}_{N+1}} 
\right) \mathlarger{\prod_{m=1}^{N}} \left(\dfrac{z \! - \! \hat{b}_{
m}}{z \! - \! \hat{a}_{m}} \right) \right)^{1/4}, &\text{$z \! \in \! 
\mathbb{C}_{-}$ $(\subset \hat{\mathcal{Y}}^{-})$,}
\end{cases}
\end{equation}
and set
\begin{equation} \label{eqssabra1} 
\hat{\gamma}(\alpha_{p_{q}}) \! := \! \left(\left(\dfrac{\alpha_{p_{q}} 
\! - \! \hat{b}_{0}}{\alpha_{p_{q}} \! - \! \hat{a}_{N+1}} \right) 
\mathlarger{\prod_{m=1}^{N}} \left(\dfrac{\alpha_{p_{q}} \! - \! 
\hat{b}_{m}}{\alpha_{p_{q}} \! - \! \hat{a}_{m}} \right) \right)^{1/4} 
\, \, (\in \! \mathbb{R}_{+}), \quad q \! = \! 1,2,\dotsc,\mathfrak{s} 
\! - \! 1,
\end{equation}
and, for $n \! \in \! \mathbb{N}$ and $k \! \in \! \lbrace 1,2,\dotsc,
K \rbrace$ such that $\alpha_{p_{\mathfrak{s}}} \! := \! \alpha_{k} 
\! \neq \! \infty$, let $\tilde{\gamma} \colon \mathbb{C} \setminus 
((-\infty,\tilde{b}_{0}) \cup (\tilde{a}_{N+1},+\infty) \cup \cup_{j=1}^{N}
(\tilde{a}_{j},\tilde{b}_{j})) \! \to \! \mathbb{C}$ be defined by
\begin{equation} \label{eqmainfin10} 
\tilde{\gamma}(z) \! := \!
\begin{cases} 
\left(\left(\dfrac{z \! - \! \tilde{b}_{0}}{z \! - \! \tilde{a}_{N+1}} 
\right) \mathlarger{\prod_{m=1}^{N}} \left(\dfrac{z \! - \! 
\tilde{b}_{m}}{z \! - \! \tilde{a}_{m}} \right) \right)^{1/4}, &\text{$z 
\! \in \! \mathbb{C}_{+}$ $(\subset \tilde{\mathcal{Y}}^{+})$,} \\
-\mi \left(\left(\dfrac{z \! - \! \tilde{b}_{0}}{z \! - \! \tilde{a}_{N+1}} 
\right) \mathlarger{\prod_{m=1}^{N}} \left(\dfrac{z \! - \! \tilde{
b}_{m}}{z \! - \! \tilde{a}_{m}} \right) \right)^{1/4}, &\text{$z \! 
\in \! \mathbb{C}_{-}$ $(\subset \tilde{\mathcal{Y}}^{-})$,}
\end{cases}
\end{equation}
and set
\begin{equation} \label{eqssabra2} 
\tilde{\gamma}(\alpha_{p_{q}}) \! := \! \left(\left(\dfrac{\alpha_{
p_{q}} \! - \! \tilde{b}_{0}}{\alpha_{p_{q}} \! - \! \tilde{a}_{N+1}} 
\right) \mathlarger{\prod_{m=1}^{N}} \left(\dfrac{\alpha_{p_{q}} 
\! - \! \tilde{b}_{m}}{\alpha_{p_{q}} \! - \! \tilde{a}_{m}} \right) 
\right)^{1/4} \, \, (\in \! \mathbb{R}_{+}), \quad q \! = \! 1,
\dotsc,\mathfrak{s} \! - \! 2,\mathfrak{s}.
\end{equation}
Then, for $n \! \in \! \mathbb{N}$ and $k \! \in \! \lbrace 1,2,
\dotsc,K \rbrace$ such that $\alpha_{p_{\mathfrak{s}}} \! := \! 
\alpha_{k} \! = \! \infty$,
\begin{equation*}
\lbrace \mathstrut z \! \in \! \mathbb{C}_{\pm} \, (\subset 
\hat{\mathcal{Y}}^{\pm}); \, \hat{\gamma}(z) \! \mp \! (\hat{
\gamma}(z))^{-1} \! = \! 0 \rbrace \! =: \! \lbrace \hat{z}_{
1}^{\pm},\hat{z}_{2}^{\pm},\dotsc,\hat{z}_{N}^{\pm} \rbrace,
\end{equation*} 
where $\hat{z}_{j}^{\pm} \! \in \! (\hat{a}_{j},\hat{b}_{j})^{\pm}$ 
$(\subset \mathbb{C}_{\pm}$ $(\subset \hat{\mathcal{Y}}^{
\pm}))$, $j \! = \! 1,2,\dotsc,N$, and as points on the plane 
$\hat{z}_{j}^{+} \! = \! \hat{z}_{j}^{-} \! = \! \hat{z}_{j} 
\! \in \! (\hat{a}_{j},\hat{b}_{j})$, $j \! = \! 1,2,\dotsc,
N$,\footnote{Strictly speaking, $\boldsymbol{\mathrm{pr}}
(\hat{z}_{j}^{+}) \! = \! \boldsymbol{\mathrm{pr}}(\hat{z}_{j}^{-}) 
\! = \! \hat{z}_{j} \! \in \! (\hat{a}_{j},\hat{b}_{j})$, $j \! = \! 1,2,
\dotsc,N$.} and, for $n \! \in \! \mathbb{N}$ and $k \! \in \! 
\lbrace 1,2,\dotsc,K \rbrace$ such that $\alpha_{p_{\mathfrak{
s}}} \! := \! \alpha_{k} \! \neq \! \infty$,
\begin{equation*}
\lbrace \mathstrut z \! \in \! \mathbb{C}_{\pm} \, (\subset 
\tilde{\mathcal{Y}}^{\pm}); \, \tilde{\gamma}(z) \! \mp \! 
(\tilde{\gamma}(z))^{-1} \! = \! 0 \rbrace \! =: \! \lbrace 
\tilde{z}_{1}^{\pm},\tilde{z}_{2}^{\pm},\dotsc,\tilde{z}_{N}^{\pm} 
\rbrace,
\end{equation*} 
where $\tilde{z}_{j}^{\pm} \! \in \! (\tilde{a}_{j},\tilde{b}_{j})^{
\pm}$ $(\subset \mathbb{C}_{\pm}$ $(\subset \tilde{\mathcal{
Y}}^{\pm}))$, $j \! = \! 1,2,\dotsc,N$, and as points on the 
plane $\tilde{z}_{j}^{+} \! = \! \tilde{z}_{j}^{-} \! = \! \tilde{z}_{j} 
\! \in \! (\tilde{a}_{j},\tilde{b}_{j})$, $j \! = \! 1,2,\dotsc,
N$.\footnote{Strictly speaking, $\boldsymbol{\mathrm{pr}}
(\tilde{z}_{j}^{+}) \! = \! \boldsymbol{\mathrm{pr}}(\tilde{z}_{j}^{
-}) \! = \! \tilde{z}_{j} \! \in \! (\tilde{a}_{j},\tilde{b}_{j})$, $j \! 
= \! 1,2,\dotsc,N$.} Furthermore, for $n \! \in \! \mathbb{N}$ 
and $k \! \in \! \lbrace 1,2,\dotsc,K \rbrace$ such that 
$\alpha_{p_{\mathfrak{s}}} \! := \! \alpha_{k} \! = \! \infty$, 
$\hat{\gamma}(z)$ solves the {\rm RHP}
\begin{enumerate}
\item[\textbullet] $\hat{\gamma}(z)$ is analytic for $z \! \in \! 
\mathbb{C} \setminus ((-\infty,\hat{b}_{0}) \cup (\hat{a}_{N+1},
+\infty) \cup \cup_{j=1}^{N}(\hat{a}_{j},\hat{b}_{j}))$,
\item[\textbullet] $\hat{\gamma}_{+}(z) \! = \! \hat{\gamma}_{-}
(z) \mi$, $z \! \in \! (-\infty,\hat{b}_{0}) \cup (\hat{a}_{N+1},
+\infty) \cup \cup_{j=1}^{N}(\hat{a}_{j},\hat{b}_{j})$,
\item[\textbullet] $\hat{\gamma}(z) \! =_{\overline{\mathbb{
C}}_{\pm} \ni z \to \alpha_{k}} \! (-\mi)^{(1 \mp 1)/2}
(1 \! + \! \mathcal{O}(z^{-1}))$, and $\hat{\gamma}(z) \! 
=_{\mathbb{C}_{\pm} \ni z \to \alpha_{p_{q}}} \! (-\mi)^{(1 \mp 1)/2} 
\hat{\gamma}(\alpha_{p_{q}})(1 \! + \! \mathcal{O}(z \! - \! 
\alpha_{p_{q}}))$, $q \! = \! 1,2,\dotsc,\mathfrak{s} \! - \! 1$,
\end{enumerate}
and, for $n \! \in \! \mathbb{N}$ and $k \! \in \! \lbrace 1,2,
\dotsc,K \rbrace$ such that $\alpha_{p_{\mathfrak{s}}} \! := \! 
\alpha_{k} \! \neq \! \infty$, $\tilde{\gamma}(z)$ solves the {\rm RHP}
\begin{enumerate}
\item[\textbullet] $\tilde{\gamma}(z)$ is analytic for $z \! \in \! 
\mathbb{C} \setminus ((-\infty,\tilde{b}_{0}) \cup (\tilde{a}_{N+1},
+\infty) \cup \cup_{j=1}^{N}(\tilde{a}_{j},\tilde{b}_{j}))$,
\item[\textbullet] $\tilde{\gamma}_{+}(z) \! = \! \tilde{
\gamma}_{-}(z) \mi$, $z \! \in \! (-\infty,\tilde{b}_{0}) \cup 
(\tilde{a}_{N+1},+\infty) \cup \cup_{j=1}^{N}(\tilde{a}_{j},
\tilde{b}_{j})$,
\item[\textbullet] $\tilde{\gamma}(z) \! =_{\overline{\mathbb{
C}}_{\pm} \ni z \to \alpha_{p_{\mathfrak{s}-1}} = \infty} \! 
(-\mi)^{(1 \mp 1)/2}(1 \! + \! \mathcal{O}(z^{-1}))$, and 
$\tilde{\gamma}(z) \! =_{\mathbb{C}_{\pm} \ni z \to \alpha_{p_{q}}} 
\! (-\mi)^{(1 \mp 1)/2} \tilde{\gamma}(\alpha_{p_{q}})(1 \! + \! 
\mathcal{O}(z \! - \! \alpha_{p_{q}}))$, $q \! = \! 1,\dotsc,
\mathfrak{s} \! - \! 2,\mathfrak{s}$.
\end{enumerate}
\end{ccccc}

\emph{Proof}. The proof of this Lemma~\ref{lem4.4} consists 
of two cases: (i) $n \! \in \! \mathbb{N}$ and $k \! \in \! \lbrace 
1,2,\dotsc,K \rbrace$ such that $\alpha_{p_{\mathfrak{s}}} \! 
:= \! \alpha_{k} \! = \! \infty$; and (ii) $n \! \in \! \mathbb{N}$ 
and $k \! \in \! \lbrace 1,2,\dotsc,K \rbrace$ such that 
$\alpha_{p_{\mathfrak{s}}} \! := \! \alpha_{k} \! \neq \! \infty$. 
Notwithstanding the fact that the scheme of the proof is, 
\emph{mutatis mutandis}, similar for both cases, without loss 
of generality, only the proof for case~(ii) is presented in detail, 
whilst case~(i) is proved analogously.

Define $\tilde{\gamma}(z)$ via Equation~\eqref{eqmainfin10}: one 
then notes that $\tilde{\gamma}(z) \! \mp \! (\tilde{\gamma}(z))^{-1} 
\! = \! 0 \! \Leftrightarrow \! (\tilde{\gamma}(z))^{2} \! \mp \! 1 \! 
= \! 0 \! \Rightarrow \! (\tilde{\gamma}(z))^{4} \! - \! 1 \! = \! 0 
\! \Leftrightarrow \! \tilde{\mathscr{Q}}(z) \! := \! (z \! - \! 
\tilde{b}_{0}) \prod_{m=1}^{N}(z \! - \! \tilde{b}_{m}) \! - \! (z \! - 
\! \tilde{a}_{N+1}) \prod_{m=1}^{N}(z \! - \! \tilde{a}_{m}) \! = \! 
0$, whence $\tilde{\mathscr{Q}}(\tilde{a}_{j}) \! = \! (-1)^{N-j+1} 
\tilde{\mathscr{Q}}_{\blacklozenge}(\tilde{a}_{j})$, $j \! = \! 
1,2,\dotsc,N$, where $\tilde{\mathscr{Q}}_{\blacklozenge}
(\tilde{a}_{j}) \! := \! (\tilde{a}_{j} \! - \! \tilde{b}_{0}) \prod_{m_{1}
=1}^{j-1}(\tilde{a}_{j}  \! - \! \tilde{b}_{m_{1}}) \prod_{m_{2}=j}^{N}
(\tilde{b}_{m_{2}} \! - \! \tilde{a}_{j}) \! > \! 0$, and $\tilde{\mathscr{
Q}}(\tilde{b}_{j}) \! = \! (-1)^{N-j} \tilde{\mathscr{Q}}_{\blacklozenge}
(\tilde{b}_{j})$, $j \! = \! 1,2,\dotsc,N$, where $\tilde{\mathscr{Q}}_{
\blacklozenge}(\tilde{b}_{j}) \! := \! (\tilde{a}_{N+1} \! - \! \tilde{b}_{j}) 
\prod_{m_{1}=1}^{j}(\tilde{b}_{j} \! - \! \tilde{a}_{m_{1}}) \prod_{m_{2}
=j+1}^{N}(\tilde{a}_{m_{2}} \! - \! \tilde{b}_{j}) \! > \! 0$; thus, 
$\tilde{\mathscr{Q}}(\tilde{a}_{j}) \tilde{\mathscr{Q}}(\tilde{b}_{j}) 
\! < \! 0$, $j \! = \! 1,2,\dotsc,N$, which shows that: (i) $\tilde{
\mathscr{Q}}(z)$ has a zero, denoted by $\tilde{z}_{j}$, in each 
open interval $(\tilde{a}_{j},\tilde{b}_{j})$, $j \! = \! 1,2,\dotsc,N$; 
and (ii) since $\tilde{\mathscr{Q}}(z)$ is a unital polynomial with 
$\deg (\tilde{\mathscr{Q}}(z)) \! = \! N$, $\lbrace \tilde{z}_{j} 
\rbrace_{j=1}^{N}$ are the only---simple---zeros of $\tilde{
\mathscr{Q}}(z)$. An analysis of the branch cuts shows that, 
for $z \! \in \! \cup_{j=1}^{N}(\tilde{a}_{j},\tilde{b}_{j})^{\pm}$ 
$(\subset \mathbb{C}_{\pm}$ $(\subset \tilde{\mathcal{Y}}^{\pm}))$, 
$\pm (\tilde{\gamma}(z))^{2} \! > \! 0$, whence $\lbrace 
\tilde{z}_{j}^{\pm} \rbrace_{j=1}^{N} \! = \! \lbrace \mathstrut 
z^{\pm} \! \in \! (\tilde{a}_{j},\tilde{b}_{j})^{\pm} \subset 
\mathbb{C}_{\pm} \, (\subset \tilde{\mathcal{Y}}^{\pm}), \, j \! = \! 
1,2,\dotsc,N; \, (\tilde{\gamma}(z) \! \mp \! (\tilde{\gamma}(z))^{-1}) 
\vert_{z=z^{\pm}} \! = \! 0 \rbrace$. Setting $\tilde{J}_{\blacklozenge} 
\! := \! (-\infty,\tilde{b}_{0}) \cup (\tilde{a}_{N+1},+\infty) \cup 
\cup_{j=1}^{N}(\tilde{a}_{j},\tilde{b}_{j})$, one shows, upon performing 
an analysis of the branch cuts, that $\tilde{\gamma}(z)$ solves the 
associated (scalar) RHP $(\tilde{\gamma}(z),\mi,\tilde{J}_{\blacklozenge})$ 
stated in the corresponding item of the lemma. \hfill $\qed$
\begin{ccccc} \label{lem4.5} 
For $n \! \in \! \mathbb{N}$ and $k \! \in \! \lbrace 1,2,\dotsc,K 
\rbrace$ such that $\alpha_{p_{\mathfrak{s}}} \! := \! \alpha_{k} \! 
= \! \infty$ (resp., $\alpha_{p_{\mathfrak{s}}} \! := \! \alpha_{k} \! 
\neq \! \infty)$, let $\mathfrak{m} \colon \mathbb{C} \setminus 
\hat{\mathscr{J}} \! \to \! \mathrm{SL}_{2}(\mathbb{C})$ (resp., 
$\mathfrak{m} \colon \mathbb{C} \setminus \tilde{\mathscr{J}} \! 
\to \! \mathrm{SL}_{2}(\mathbb{C}))$ solve the model {\rm RHP} 
stated in item~{\rm \pmb{(1)}} (resp., item~{\rm \pmb{(2)})} 
of Lemma~\ref{lem4.3}. Then: {\rm \pmb{(1)}} for $n \! \in \! 
\mathbb{N}$ and $k \! \in \! \lbrace 1,2,\dotsc,K \rbrace$ such 
that $\alpha_{p_{\mathfrak{s}}} \! := \! \alpha_{k} \! = \! \infty$,
\begin{equation*}
\mathfrak{m}(z) \! = \! 
\begin{cases}
\mathbb{M}(z), &\text{$z \! \in \! \mathbb{C}_{+}$ $(\subset 
\hat{\mathcal{Y}}^{+})$,} \\
-\mi \mathbb{M}(z) \sigma_{2}, &\text{$z \in \! \mathbb{C}_{-}$ 
$(\subset \hat{\mathcal{Y}}^{-})$,}
\end{cases}
\end{equation*}
with
\begin{equation*}
\mathbb{M}(z) \! = \! 
\tilde{\mathfrak{m}}^{\raise-0.5ex\hbox{$\scriptstyle \infty$}} 
\hat{\mathbb{M}}(z),
\end{equation*}
where $\tilde{\mathfrak{m}}^{\raise-0.5ex\hbox{$\scriptstyle \infty$}}$ 
and $\hat{\mathbb{M}}(z)$ are given by Equations~\eqref{eqmaininf8} 
and~\eqref{eqmaininf9}, respectively, $\hat{\gamma}(z)$ is defined by 
Equation~\eqref{eqmaininf10}, $\hat{\boldsymbol{\Omega}} \! := \! 
(\hat{\Omega}_{1},\hat{\Omega}_{2},\dotsc,\hat{\Omega}_{N})^{
\mathrm{T}}$,\footnote{Note: ${}^{\mathrm{T}}$ denotes transposition.} 
with $\hat{\Omega}_{j} \! = \! 2 \pi \int_{\hat{b}_{j}}^{\hat{a}_{N+1}} 
\psi_{\widetilde{V}}^{\infty}(\xi) \, \md \xi$, $j \! = \! 1,2,\dotsc,N$, 
$\hat{\boldsymbol{d}} \! \equiv \! -\sum_{j=1}^{N} \int_{\hat{a}_{j}}^{
\hat{z}_{j}^{-}} \hat{\boldsymbol{\omega}}$ $(= \! \sum_{j=1}^{N} 
\int_{\hat{a}_{j}}^{\hat{z}_{j}^{+}} \hat{\boldsymbol{\omega}})$, 
where $\lbrace \hat{z}_{j}^{\pm} \rbrace_{j=1}^{N}$ $(\subset 
\hat{\mathcal{Y}}^{\pm})$ are characterised in the corresponding item 
of Lemma~\ref{lem4.4}, $\hat{\boldsymbol{\omega}}$ is the associated 
normalised basis of holomorphic one-forms on $\hat{\mathcal{Y}}$, 
$\hat{\boldsymbol{u}}(z) \! := \! \int_{\hat{a}_{N+1}}^{z} \hat{\boldsymbol{\omega}}$ 
$(\in \! \operatorname{Jac}(\hat{\mathcal{Y}}))$, and $\hat{\boldsymbol{u}}_{+}
(\infty)$ $(= \! \hat{\boldsymbol{u}}(\infty^{+}))$ $:= \! \int_{\hat{a}_{N+1}}^{
\infty^{+}} \hat{\boldsymbol{\omega}}$$;$\footnote{Note: $\alpha_{p_{\mathfrak{s}}}^{\pm} 
\! := \! \alpha_{k}^{\pm} \! = \! \infty^{\pm}$ are the points at infinity in 
$\hat{\mathcal{Y}}^{\pm}$ $(\supset \overline{\mathbb{C}}_{\pm})$.} and {\rm \pmb{(2)}} 
for $n \! \in \! \mathbb{N}$ and $k \! \in \! \lbrace 1,2,\dotsc,K \rbrace$ such that 
$\alpha_{p_{\mathfrak{s}}} \! := \! \alpha_{k} \! \neq \! \infty$,
\begin{equation*}
\mathfrak{m}(z) \! = \! 
\begin{cases}
\mathbb{M}(z), &\text{$z \! \in \! \mathbb{C}_{+}$ $(\subset 
\tilde{\mathcal{Y}}^{+})$,} \\
-\mi \mathbb{M}(z) \sigma_{2}, &\text{$z \in \! \mathbb{C}_{-}$ 
$(\subset \tilde{\mathcal{Y}}^{-})$,}
\end{cases}
\end{equation*}
with
\begin{equation*}
\mathbb{M}(z) \! = \! \widetilde{\mathbb{K}} \, \tilde{\mathbb{M}}(z),
\end{equation*}
where $\widetilde{\mathbb{K}} \! := \! \mathscr{E}^{-\sigma_{3}} 
\tilde{\mathbb{K}}$, with $\tilde{\mathbb{K}}$ and $\mathscr{E}$ 
defined by Equations~\eqref{eqmainfin8} and~\eqref{eqmainfin13}, 
respectively,\footnote{See, also, Equation~\eqref{eqconsfinn1} for 
the definition of $\mathfrak{k}_{0}$.} $\tilde{\mathbb{M}}(z)$ is 
defined by Equation~\eqref{eqmainfin9}, $\tilde{\gamma}(z)$ is 
defined by Equation~\eqref{eqmainfin10}, $\tilde{\boldsymbol{\Omega}} 
\! := \! (\tilde{\Omega}_{1},\tilde{\Omega}_{2},\dotsc,
\tilde{\Omega}_{N})^{\mathrm{T}}$, with $\tilde{\Omega}_{j} \! = \! 
2 \pi \int_{\tilde{b}_{j}}^{\tilde{a}_{N+1}} \psi_{\widetilde{V}}^{f}(\xi) 
\, \md \xi$, $j \! = \! 1,2,\dotsc,N$, $\tilde{\boldsymbol{d}} \! 
\equiv \! -\sum_{j=1}^{N} \int_{\tilde{a}_{j}}^{\tilde{z}_{j}^{-}} 
\tilde{\boldsymbol{\omega}}$ $(= \! \sum_{j=1}^{N} 
\int_{\tilde{a}_{j}}^{\tilde{z}_{j}^{+}} \tilde{\boldsymbol{\omega}})$, 
where $\lbrace \tilde{z}_{j}^{\pm} \rbrace_{j=1}^{N}$ $(\subset \tilde{
\mathcal{Y}}^{\pm})$ are characterised in the corresponding item 
of Lemma~\ref{lem4.4}, $\tilde{\boldsymbol{\omega}}$ is the 
associated normalised basis of holomorphic one-forms on $\tilde{
\mathcal{Y}}$, $\tilde{\boldsymbol{u}}(z) \! := \! \int_{\tilde{a}_{N+1}}^{z} 
\tilde{\boldsymbol{\omega}}$ $(\in \! \operatorname{Jac}(\tilde{
\mathcal{Y}}))$, and $\tilde{\boldsymbol{u}}_{+}(\alpha_{k})$ $(= \! 
\tilde{\boldsymbol{u}}(\alpha_{k}^{+}))$ $:= \! \int_{\tilde{a}_{N+1}}^{
\alpha_{k}^{+}} \tilde{\boldsymbol{\omega}}$.\footnote{Note: 
$\alpha_{p_{\mathfrak{s}}}^{\pm} \! := \! \alpha_{k}^{\pm} \! \in \! 
\tilde{\mathcal{Y}}^{\pm}$ $(\supset \mathbb{C}_{\pm})$.} Furthermore, for 
$n \! \in \! \mathbb{N}$ and $k \! \in \! \lbrace 1,2,\dotsc,K \rbrace$ such 
that $\alpha_{p_{\mathfrak{s}}} \! := \! \alpha_{k} \! = \! \infty$ (resp., 
$\alpha_{p_{\mathfrak{s}}} \! := \! \alpha_{k} \! \neq \! \infty)$, the 
associated solution is unique.
\end{ccccc}

\emph{Proof}. The proof of this Lemma~\ref{lem4.5} consists of two 
cases: (i) $n \! \in \! \mathbb{N}$ and $k \! \in \! \lbrace 1,2,\dotsc,
K \rbrace$ such that $\alpha_{p_{\mathfrak{s}}} \! := \! \alpha_{k} \! 
= \! \infty$; and (ii) $n \! \in \! \mathbb{N}$ and $k \! \in \! \lbrace 
1,2,\dotsc,K \rbrace$ such that $\alpha_{p_{\mathfrak{s}}} \! := \! 
\alpha_{k} \! \neq \! \infty$. Notwithstanding the fact that the scheme 
of the proof is, \emph{mutatis mutandis}, similar for both cases, 
without loss of generality, only the proof for case~(ii) is presented in 
detail, whilst case~(i) is proved analogously.

For $n \! \in \! \mathbb{N}$ and $k \! \in \! \lbrace 1,2,\dotsc,K 
\rbrace$ such that $\alpha_{p_{\mathfrak{s}}} \! := \! \alpha_{k} \! 
\neq \! \infty$, let $\mathfrak{m} \colon \mathbb{C} \setminus 
\tilde{\mathscr{J}} \! \to \! \mathrm{SL}_{2}(\mathbb{C})$ solve 
the associated model (matrix) RHP stated in item~\pmb{(2)} of 
Lemma~\ref{lem4.3}, and define $\mathfrak{m}(z)$ in terms of 
$\mathbb{M}(z)$ as in item~\pmb{(2)} of the lemma. A straightforward 
calculation shows that, for $n \! \in \! \mathbb{N}$ and $k \! \in \! 
\lbrace 1,2,\dotsc,K \rbrace$ such that $\alpha_{p_{\mathfrak{s}}} 
\! := \! \alpha_{k} \! \neq \! \infty$, $\mathbb{M}(z)$ solves the 
following (equivalent) `twisted' matrix RHP: (i) $\mathbb{M}(z)$ is 
analytic for $z \! \in \! \mathbb{C} \setminus \tilde{\mathscr{J}}_{T}$, 
where $\tilde{\mathscr{J}}_{T} \! := \! (-\infty,\tilde{b}_{0}) \cup 
(\tilde{a}_{N+1},+\infty) \cup \cup_{j=1}^{N}(\tilde{a}_{j},\tilde{b}_{j})$; 
(ii) the boundary values $\mathbb{M}_{\pm}(z) \! := \! 
\lim_{\underset{z^{\prime} \, \in \, \pm \, \mathrm{side} \, 
\mathrm{of} \, \tilde{\mathscr{J}}_{T}}{z^{\prime} \to z \, \in \, 
\tilde{\mathscr{J}}_{T}}} \mathbb{M}(z^{\prime})$ satisfy the boundary 
condition $\mathbb{M}_{+}(z) \! = \! \mathbb{M}_{-}(z) \tilde{\daleth}
(z)$ $\text{a.e.}$ $z \! \in \! \tilde{\mathscr{J}}_{T}$, where
\begin{align} \label{eqlmwayv1} 
\tilde{\daleth} &\colon \mathbb{N} \times \lbrace 1,2,\dotsc,K \rbrace 
\times \tilde{\mathscr{J}}_{T} \! \to \! \operatorname{GL}_{2}(\mathbb{C}), 
\, (n,k,z) \! \mapsto \! \tilde{\daleth}(n,k,z) \! =: \! \tilde{\daleth}(z) 
\nonumber \\
&= 
\begin{cases}
\mathrm{I}, &\text{$z\! \in \! J_{f}$,} \\
-\mi \sigma_{2}, &\text{$z \! \in \! (-\infty,\tilde{b}_{0}) \cup 
(\tilde{a}_{N+1},+\infty)$,} \\
-\mi \sigma_{2} \me^{-\mi ((n-1)K+k) \tilde{\Omega}_{j} \sigma_{3}}, 
&\text{$z \! \in \! (\tilde{a}_{j},\tilde{b}_{j}), \quad j \! = \! 1,2,\dotsc,N$,}
\end{cases}
\end{align}
with $\tilde{\Omega}_{j} \! = \! 2 \pi \int_{\tilde{b}_{j}}^{\tilde{a}_{N+1}} 
\psi_{\widetilde{V}}^{f}(\xi) \, \md \xi$, $j \! = \! 1,2,\dotsc,N$; and (iii)
\begin{gather*}
\mathbb{M}(z) \underset{\tilde{\mathcal{Y}}^{+} \supset \mathbb{C}_{+} 
\ni z \to \alpha_{k}}{=} \mathscr{E}^{-\sigma_{3}} \! + \! \mathcal{O}
(z \! - \! \alpha_{k}), \, \quad \qquad \, \mathbb{M}(z) \underset{
\tilde{\mathcal{Y}}^{-} \supset \mathbb{C}_{-} \ni z \to \alpha_{k}}{=} 
\mi \mathscr{E}^{\sigma_{3}} \sigma_{2} \! + \! \mathcal{O}
(z \! - \! \alpha_{k}), \\
\mathbb{M}(z) \underset{\tilde{\mathcal{Y}}^{\pm} \supset \overline{
\mathbb{C}}_{\pm} \ni z \to \alpha_{p_{\mathfrak{s}-1}} = \infty}{=} 
\mathcal{O}(\me^{n(\hat{\mathscr{P}}^{\pm}_{0}-\hat{\mathscr{P}}_{1}) 
\sigma_{3}} \mathscr{E}^{\mp \sigma_{3}}), \, \quad \qquad \mathbb{M}(z) 
\underset{\tilde{\mathcal{Y}}^{\pm} \supset \mathbb{C}_{\pm} \ni z 
\to \alpha_{p_{q}}}{=} \mathcal{O}(\me^{n(\hat{\mathscr{P}}^{\pm}_{0}
-\hat{\mathscr{P}}^{\pm}_{2}) \sigma_{3}} \mathscr{E}^{\mp \sigma_{3}}), 
\quad q \! = \! 1,2,\dotsc,\mathfrak{s} \! - \! 2.
\end{gather*}
The solution of this latter RHP for $\mathbb{M}(z)$ is constructed out 
of the solution of two simpler RHPs: $(\widetilde{\mathscr{N}}(z),-\mi 
\sigma_{2},\tilde{\mathscr{J}}_{T})$ and $(\tilde{\mathfrak{m}}(z),
\tilde{\mathscr{U}}(z),\tilde{\mathscr{J}}_{T})$, where $\tilde{\mathscr{
U}}(z)$ equals $\exp (\mi ((n \! - \! 1)K \! + \! k) \tilde{\Omega}_{j} \, 
\sigma_{3}) \sigma_{1}$ for $z \! \in \! (\tilde{a}_{j},\tilde{b}_{j})$, $j \! 
= \! 1,2,\dotsc,N$, and equals $\sigma_{1}$ for $z \! \in \! (-\infty,
\tilde{b}_{0}) \cup (\tilde{a}_{N+1},+\infty)$.

For $n \! \in \! \mathbb{N}$ and $k \! \in \! \lbrace 1,2,\dotsc,K 
\rbrace$ such that $\alpha_{p_{\mathfrak{s}}} \! := \! \alpha_{k} \! 
\neq \! \infty$, the RHP $(\widetilde{\mathscr{N}}(z),-\mi \sigma_{2},
\tilde{\mathscr{J}}_{T})$ is solved explicitly by diagonalising the jump 
matrix, and thus reduced to two scalar RHPs \cite{a53} (see, also, 
\cite{a51,a61,a67}): the solution is given by {}\footnote{Strictly 
speaking, though, $\widetilde{\mathscr{N}}(z)$ doesn't solve the RHP 
$(\widetilde{\mathscr{N}}(z),-\mi \sigma_{2},\tilde{\mathscr{J}}_{T})$ 
in the heretofore-defined sense, as $\widetilde{\mathscr{N}} \! \! \! 
\upharpoonright_{\mathbb{C}_{\pm}}$ can not be extended continuously 
to $\mathbb{C}_{\pm}$; however, $\widetilde{\mathscr{N}}(\pmb{\cdot} 
\! \pm \! \mi \varepsilon)$ converge in $\mathcal{L}^{2}_{\mathrm{M}_{2}
(\mathbb{C}),\mathrm{loc}}(\mathbb{R})$ as $\varepsilon \! \downarrow 
\! 0$ to $\mathrm{SL}_{2}(\mathbb{C})$-valued functions $\widetilde{
\mathscr{N}}(z)$ in $\mathcal{L}^{2}_{\mathrm{M}_{2}(\mathbb{C})}
(\tilde{\mathscr{J}}_{T})$ that satisfy $\widetilde{\mathscr{N}}_{+}(z) 
\! = \! \widetilde{\mathscr{N}}_{-}(z)(-\mi \sigma_{2})$ $\text{a.e.}$ 
$z \! \in \! \tilde{\mathscr{J}}_{T}$: one then shows that $\widetilde{
\mathscr{N}}(z)$ is the unique solution of the RHP $(\widetilde{
\mathscr{N}}(z),-\mi \sigma_{2},\tilde{\mathscr{J}}_{T})$, where the 
latter boundary/jump condition is interpreted in the $\mathcal{L}^{2}_{
\mathrm{M}_{2}(\mathbb{C}),\mathrm{loc}}(\tilde{\mathscr{J}}_{T})$ sense.}
\begin{equation*}
\widetilde{\mathscr{N}}(z) \! = \! 
\begin{pmatrix}
\frac{1}{2}(\tilde{\gamma}(z) \! + \! (\tilde{\gamma}(z))^{-1}) & 
-\frac{1}{2 \mi}(\tilde{\gamma}(z) \! - \! (\tilde{\gamma}(z))^{-1}) \\
\frac{1}{2 \mi}(\tilde{\gamma}(z) \! - \! (\tilde{\gamma}(z))^{-1}) & 
\frac{1}{2}(\tilde{\gamma}(z) \! + \! (\tilde{\gamma}(z))^{-1})
\end{pmatrix},
\end{equation*}
where $\tilde{\gamma} \colon \mathbb{C} \setminus \tilde{\mathscr{J}}_{T} 
\! \to \! \mathbb{C}$ is described in the corresponding item of 
Lemma~\ref{lem4.4}; furthermore, $\widetilde{\mathscr{N}}(z)$ is 
piecewise (sectionally) analytic for $z \! \in \! \mathbb{C} \setminus 
\tilde{\mathscr{J}}_{T}$, and
\begin{gather*}
\widetilde{\mathscr{N}}(z) \underset{\tilde{\mathcal{Y}}^{+} 
\supset \mathbb{C}_{+} \ni z \to \alpha_{p_{q}}}{=} 
\begin{pmatrix}
\frac{1}{2}(\tilde{\gamma}(\alpha_{p_{q}}) \! + \! (\tilde{\gamma}
(\alpha_{p_{q}}))^{-1}) & -\frac{1}{2 \mi}(\tilde{\gamma}(\alpha_{p_{q}}) 
\! - \! (\tilde{\gamma}(\alpha_{p_{q}}))^{-1}) \\
\frac{1}{2 \mi}(\tilde{\gamma}(\alpha_{p_{q}}) \! - \! (\tilde{\gamma}
(\alpha_{p_{q}}))^{-1}) & \frac{1}{2}(\tilde{\gamma}(\alpha_{p_{q}}) \! 
+ \! (\tilde{\gamma}(\alpha_{p_{q}}))^{-1})
\end{pmatrix} \! + \! \mathcal{O}(z \! - \! \alpha_{p_{q}}), \quad q \! 
= \! 1,\dotsc,\mathfrak{s} \! - \! 2,\mathfrak{s}, \\
\widetilde{\mathscr{N}}(z) \underset{\tilde{\mathcal{Y}}^{-} 
\supset \mathbb{C}_{-} \ni z \to \alpha_{p_{q}}}{=} \mi 
\begin{pmatrix}
\frac{1}{2}(\tilde{\gamma}(\alpha_{p_{q}}) \! + \! (\tilde{\gamma}
(\alpha_{p_{q}}))^{-1}) & -\frac{1}{2 \mi}(\tilde{\gamma}
(\alpha_{p_{q}}) \! - \! (\tilde{\gamma}(\alpha_{p_{q}}))^{-1}) \\
\frac{1}{2 \mi}(\tilde{\gamma}(\alpha_{p_{q}}) \! - \! (\tilde{\gamma}
(\alpha_{p_{q}}))^{-1}) & \frac{1}{2}(\tilde{\gamma}(\alpha_{p_{q}}) \! 
+ \! (\tilde{\gamma}(\alpha_{p_{q}}))^{-1})
\end{pmatrix} \sigma_{2} \! + \! \mathcal{O}(z \! - \! \alpha_{p_{q}}), 
\quad q \! = \! 1,\dotsc,\mathfrak{s} \! - \! 2,\mathfrak{s}, \\
\widetilde{\mathscr{N}}(z) \underset{\tilde{\mathcal{Y}}^{+} \supset 
\overline{\mathbb{C}}_{+} \ni z \to \alpha_{p_{\mathfrak{s}-1}} = 
\infty}{=} \mathrm{I} \! + \! \mathcal{O}(z^{-1}), \qquad \qquad 
\widetilde{\mathscr{N}}(z) \underset{\tilde{\mathcal{Y}}^{-} \supset 
\overline{\mathbb{C}}_{-} \ni z \to \alpha_{p_{\mathfrak{s}-1}} = 
\infty}{=} \mi \sigma_{2} \! + \! \mathcal{O}(z^{-1}),
\end{gather*}
where $\tilde{\gamma}(\alpha_{p_{q}})$, $q \! = \! 1,\dotsc,\mathfrak{s} 
\! - \! 2,\mathfrak{s}$, is defined by Equation~\eqref{eqssabra2}.

For $n \! \in \! \mathbb{N}$ and $k \! \in \! \lbrace 1,2,\dotsc,K 
\rbrace$ such that $\alpha_{p_{\mathfrak{s}}} \! := \! \alpha_{k} \! 
\neq \! \infty$, consider the functions {}\footnote{See, also, the 
construction in \cite{kmo1}.} $\tilde{\boldsymbol{\theta}}
(\tilde{\boldsymbol{u}}(z) \! \pm \! \tilde{\boldsymbol{d}})$, where 
$\tilde{\boldsymbol{u}}(z) \colon \tilde{\mathcal{Y}} \! \to \! 
\operatorname{Jac}(\tilde{\mathcal{Y}})$, $z \! \mapsto \! 
\tilde{\boldsymbol{u}}(z) \! := \! \int_{\tilde{a}_{N+1}}^{z} 
\tilde{\boldsymbol{\omega}}$, with $\tilde{\boldsymbol{\omega}}$ 
the associated normalised basis of holomorphic one-forms on 
$\tilde{\mathcal{Y}}$, $\tilde{\boldsymbol{d}} \! \equiv \! 
-\sum_{j=1}^{N} \int_{\tilde{a}_{j}}^{\tilde{z}_{j}^{-}} \tilde{
\boldsymbol{\omega}} \! = \! \sum_{j=1}^{N} \int_{\tilde{a}_{j}}^{
\tilde{z}_{j}^{+}} \tilde{\boldsymbol{\omega}}$, where $\equiv$ 
denotes equivalence modulo the associated period lattice 
(that is, $\mathbb{C}^{N}/\lbrace \mathrm{N}_{0} \! + \! 
\tilde{\boldsymbol{\tau}} \mathrm{M}_{0} \rbrace, (\mathrm{N}_{0},
\mathrm{M}_{0}) \! \in \! \mathbb{Z}^{N} \times \mathbb{Z}^{N}$, 
where $\tilde{\boldsymbol{\tau}} \! = \! (\tilde{\boldsymbol{\tau}})_{
i,j=1,2,\dotsc,N} \! := \! (\oint_{\tilde{\boldsymbol{\beta}}_{j}} 
\tilde{\omega}_{i})_{i,j=1,2,\dotsc,N}$ is the associated $N \times 
N$ Riemann matrix of $\tilde{\boldsymbol{\beta}}$ periods), 
and $\lbrace \tilde{z}_{j}^{\pm} \rbrace_{j=1}^{N}$ $(\subset 
\tilde{\mathcal{Y}}^{\pm}$ $(\supset \mathbb{C}_{\pm}))$ are 
characterised in the corresponding item of Lemma~\ref{lem4.4}. 
{}From the general theory of theta functions on Riemann surfaces 
(see, for example, \cite{gsp,hmfik}), in the generic case for $z \! \in \! 
\tilde{\mathcal{Y}}$, $\tilde{\boldsymbol{\theta}}(\tilde{\boldsymbol{u}}
(z) \! \pm \! \tilde{\boldsymbol{d}})$ are either identically zero on 
$\tilde{\mathcal{Y}}$ or have precisely $N$ (simple) zeros. In this 
case, since the divisors {}\footnote{Defined slightly differently than 
in Subsection~\ref{subsub1}.} $\sum_{j=1}^{N} \tilde{z}_{j}^{\pm}$ 
are non-special,\footnote{See Section~2 in \cite{jhuim}; see, also, 
\cite{mymo4}.} one uses Lemma~3.27 of \cite{a67} (see, also, 
Lemma~4.2 of \cite{a59}) and the representation \cite{hmfik} 
$\tilde{\boldsymbol{K}} \! = \! \sum_{j=1}^{N} \int_{\tilde{a}_{j}}^{
\tilde{a}_{N+1}} \tilde{\boldsymbol{\omega}}$ for the associated vector 
of Riemann constants, with $\tilde{\boldsymbol{K}}$ a point of order $2$ 
(that is, $2 \tilde{\boldsymbol{K}} \! = \! 0$ and $p \tilde{\boldsymbol{K}} 
\! \not= \! 0$, $0 \! < \! p \! < \! 2)$, to arrive at
\begin{equation*}
\tilde{\boldsymbol{\theta}}(\tilde{\boldsymbol{u}}(z) \! \pm \! 
\tilde{\boldsymbol{d}}) \! = \! \tilde{\boldsymbol{\theta}} \left(
\tilde{\boldsymbol{u}}(z) \! - \! \sum_{j=1}^{N} \int_{\tilde{a}_{
j}}^{\tilde{z}_{j}^{\mp}} \tilde{\boldsymbol{\omega}} \right) \! = 
\! \tilde{\boldsymbol{\theta}} \left(\int_{\tilde{a}_{N+1}}^{z} 
\tilde{\boldsymbol{\omega}} \! - \! \tilde{\boldsymbol{K}} \! - 
\! \sum_{j=1}^{N} \int_{\tilde{a}_{N+1}}^{\tilde{z}_{j}^{\mp}} 
\tilde{\boldsymbol{\omega}} \right) \! = \! 0 \Leftrightarrow z 
\! \in \! \lbrace \tilde{z}_{1}^{\mp},\tilde{z}_{2}^{\mp},\dotsc,
\tilde{z}_{N}^{\mp} \rbrace \, \, (\subset \tilde{\mathcal{Y}}^{\mp}).
\end{equation*}
For $n \! \in \! \mathbb{N}$ and $k \! \in \! \lbrace 1,2,\dotsc,K 
\rbrace$ such that $\alpha_{p_{\mathfrak{s}}} \! := \! \alpha_{k} 
\! \neq \! \infty$, following Lemma~3.21 of \cite{a67}, set
\begin{equation} \label{eqtilfrakm} 
\tilde{\mathfrak{m}}(z) := \! 
\begin{pmatrix}
\frac{\tilde{\boldsymbol{\theta}}(\tilde{\boldsymbol{u}}(z)-
\frac{1}{2 \pi}((n-1)K+k) \tilde{\boldsymbol{\Omega}}+\tilde{
\boldsymbol{d}})}{\tilde{\boldsymbol{\theta}}(\tilde{\boldsymbol{u}}
(z)+\tilde{\boldsymbol{d}})} & \frac{\tilde{\boldsymbol{\theta}}
(-\tilde{\boldsymbol{u}}(z)-\frac{1}{2 \pi}((n-1)K+k) \tilde{
\boldsymbol{\Omega}}+\tilde{\boldsymbol{d}})}{\tilde{
\boldsymbol{\theta}}(-\tilde{\boldsymbol{u}}(z)+\tilde{
\boldsymbol{d}})} \\
\frac{\tilde{\boldsymbol{\theta}}(\tilde{\boldsymbol{u}}(z)-
\frac{1}{2 \pi}((n-1)K+k) \tilde{\boldsymbol{\Omega}}-\tilde{
\boldsymbol{d}})}{\tilde{\boldsymbol{\theta}}(\tilde{\boldsymbol{u}}
(z)-\tilde{\boldsymbol{d}})} & \frac{\tilde{\boldsymbol{\theta}}
(-\tilde{\boldsymbol{u}}(z)-\frac{1}{2 \pi}((n-1)K+k) \tilde{
\boldsymbol{\Omega}}-\tilde{\boldsymbol{d}})}{\tilde{
\boldsymbol{\theta}}(\tilde{\boldsymbol{u}}(z)+\tilde{
\boldsymbol{d}})}
\end{pmatrix},
\end{equation}
where $\tilde{\boldsymbol{\Omega}} \! := \! (\tilde{\Omega}_{1},
\tilde{\Omega}_{2},\dotsc,\tilde{\Omega}_{N})^{\mathrm{T}}$ $(\in 
\! \mathbb{R}^{N})$, with $\tilde{\Omega}_{j}$, $j \! = \! 1,2,\dotsc,N$, 
given above. Using Lemma~3.18 of \cite{a67} (or, equivalently, 
Equations~(4.65) and~(4.66) of \cite{a59}), that is, $\tilde{\boldsymbol{u}}_{+}
(z) \! + \! \tilde{\boldsymbol{u}}_{-}(z) \! \equiv \! -\tilde{\tau}_{j}$ 
$(:= \! -\tilde{\tau}e_{j})$, $z \! \in \! (\tilde{a}_{j},\tilde{b}_{j})$, 
$j \! = \! 1,2,\dotsc,N$, and $\tilde{\boldsymbol{u}}_{+}(z) \! + \! 
\tilde{\boldsymbol{u}}_{-}(z) \! \equiv \! 0$, $z \! \in \! (-\infty,\tilde{b}_{0}) 
\cup (\tilde{a}_{N+1},+\infty)$, where $\tilde{\boldsymbol{u}}_{\pm}(z) \! 
:= \! \int_{\tilde{a}_{N+1}}^{z^{\pm}} \tilde{\boldsymbol{\omega}}$, via the 
evenness and quasi-periodicity properties of the associated Riemann theta 
function $\tilde{\boldsymbol{\theta}}$, one shows that, for $z \! \in \! 
(\tilde{a}_{j},\tilde{b}_{j}) \setminus \lbrace \tilde{z}_{j}^{\pm} \rbrace$, 
$j \! = \! 1,2,\dotsc,N$,
\begin{gather*}
\dfrac{\tilde{\boldsymbol{\theta}}(\tilde{\boldsymbol{u}}_{+}(z) \! - \! 
\frac{1}{2 \pi}((n \! - \! 1)K \! + \! k) \tilde{\boldsymbol{\Omega}} \! + \! 
\tilde{\boldsymbol{d}})}{\tilde{\boldsymbol{\theta}}(\tilde{\boldsymbol{
u}}_{+}(z) \! + \! \tilde{\boldsymbol{d}})} = \me^{-\mi ((n-1)K+k) 
\tilde{\Omega}_{j}} \dfrac{\tilde{\boldsymbol{\theta}}(-\tilde{
\boldsymbol{u}}_{-}(z) \! - \! \frac{1}{2 \pi}((n \! - \! 1)K \! + \! k) 
\tilde{\boldsymbol{\Omega}} \! + \! \tilde{\boldsymbol{d}})}{\tilde{
\boldsymbol{\theta}}(-\tilde{\boldsymbol{u}}_{-}(z) \! + \! \tilde{
\boldsymbol{d}})}, \\
\dfrac{\tilde{\boldsymbol{\theta}}(\tilde{\boldsymbol{u}}_{+}(z) \! - \! 
\frac{1}{2 \pi}((n \! - \! 1)K \! + \! k) \tilde{\boldsymbol{\Omega}} \! - \! 
\tilde{\boldsymbol{d}})}{\tilde{\boldsymbol{\theta}}(\tilde{\boldsymbol{
u}}_{+}(z) \! - \! \tilde{\boldsymbol{d}})} = \me^{-\mi ((n-1)K+k) 
\tilde{\Omega}_{j}} \dfrac{\tilde{\boldsymbol{\theta}}(-\tilde{
\boldsymbol{u}}_{-}(z) \! - \! \frac{1}{2 \pi}((n \! - \! 1)K \! + \! k) 
\tilde{\boldsymbol{\Omega}} \! - \! \tilde{\boldsymbol{d}})}{\tilde{
\boldsymbol{\theta}}(\tilde{\boldsymbol{u}}_{-}(z) \! + \! \tilde{
\boldsymbol{d}})}, \\
\dfrac{\tilde{\boldsymbol{\theta}}(-\tilde{\boldsymbol{u}}_{+}(z) \! - \! 
\frac{1}{2 \pi}((n \! - \! 1)K \! + \! k) \tilde{\boldsymbol{\Omega}} \! + \! 
\tilde{\boldsymbol{d}})}{\tilde{\boldsymbol{\theta}}(-\tilde{\boldsymbol{
u}}_{+}(z) \! + \! \tilde{\boldsymbol{d}})} = \me^{\mi ((n-1)K+k) 
\tilde{\Omega}_{j}} \dfrac{\tilde{\boldsymbol{\theta}}(\tilde{
\boldsymbol{u}}_{-}(z) \! - \! \frac{1}{2 \pi}((n \! - \! 1)K \! + \! k) 
\tilde{\boldsymbol{\Omega}} \! + \! \tilde{\boldsymbol{d}})}{\tilde{
\boldsymbol{\theta}}(\tilde{\boldsymbol{u}}_{-}(z) \! + \! \tilde{
\boldsymbol{d}})}, \\
\dfrac{\tilde{\boldsymbol{\theta}}(-\tilde{\boldsymbol{u}}_{+}(z) \! - \! 
\frac{1}{2 \pi}((n \! - \! 1)K \! + \! k) \tilde{\boldsymbol{\Omega}} \! - \! 
\tilde{\boldsymbol{d}})}{\tilde{\boldsymbol{\theta}}(\tilde{\boldsymbol{
u}}_{+}(z) \! + \! \tilde{\boldsymbol{d}})} = \me^{\mi ((n-1)K+k) 
\tilde{\Omega}_{j}} \dfrac{\tilde{\boldsymbol{\theta}}(\tilde{
\boldsymbol{u}}_{-}(z) \! - \! \frac{1}{2 \pi}((n \! - \! 1)K \! + \! k) 
\tilde{\boldsymbol{\Omega}} \! - \! \tilde{\boldsymbol{d}})}{\tilde{
\boldsymbol{\theta}}(\tilde{\boldsymbol{u}}_{-}(z) \! - \! \tilde{
\boldsymbol{d}})},
\end{gather*}
and, for $z \! \in \! (-\infty,\tilde{b}_{0}) \cup (\tilde{a}_{N+1},+\infty)$, 
one obtains the same relations as above but with $\exp (\mp \mi ((n \! 
- \! 1)K \! + \! k) \tilde{\Omega}_{j}) \! \to \! 1$, $j \! = \! 1,2,\dotsc,N$. 
For $n \! \in \! \mathbb{N}$ and $k \! \in \! \lbrace 1,2,\dotsc,K \rbrace$ 
such that $\alpha_{p_{\mathfrak{s}}} \! := \! \alpha_{k} \! \neq \! \infty$, 
in analogy with Proposition~3.31 of \cite{a67}, let
\begin{equation*}
\tilde{\mathbb{M}}(z) \! := \!
\begin{pmatrix}
\widetilde{\mathscr{N}}_{11}(z) \tilde{\mathfrak{m}}_{11}(z) & 
\widetilde{\mathscr{N}}_{12}(z) \tilde{\mathfrak{m}}_{12}(z) \\
\widetilde{\mathscr{N}}_{21}(z) \tilde{\mathfrak{m}}_{21}(z) & 
\widetilde{\mathscr{N}}_{22}(z) \tilde{\mathfrak{m}}_{22}(z)
\end{pmatrix}.
\end{equation*}
Recalling that, for $n \! \in \! \mathbb{N}$ and $k \! \in \! \lbrace 
1,2,\dotsc,K \rbrace$ such that $\alpha_{p_{\mathfrak{s}}} \! := 
\! \alpha_{k} \! \neq \! \infty$, $\widetilde{\mathscr{N}} \colon 
\mathbb{C} \setminus \tilde{\mathscr{J}}_{T} \! \to \! \mathrm{SL}_{2}
(\mathbb{C})$ solves the RHP $(\widetilde{\mathscr{N}}(z),-\mi 
\sigma_{2},\tilde{\mathscr{J}}_{T})$, via the above theta-functional 
relations and the $\tilde{\mathcal{Y}}^{\pm} \supset 
\mathbb{C}_{\pm} \! \ni \! z \! \to \! \alpha_{p_{q}}$, $q \! = \! 
1,\dotsc,\mathfrak{s} \! - \! 2,\mathfrak{s}$, asymptotics of 
$\tilde{\boldsymbol{u}}(z)$ (see Section~\ref{sek5}, the proof of 
Proposition~\ref{propo5.1}), one shows that, for $n \! \in \! \mathbb{N}$ 
and $k \! \in \! \lbrace 1,2,\dotsc,K \rbrace$ such that $\alpha_{
p_{\mathfrak{s}}} \! := \! \alpha_{k} \! \neq \! \infty$, $\tilde{\mathbb{M}}
(z)$ solves the following RHP: (i) $\tilde{\mathbb{M}}(z)$ is analytic 
for $z \! \in \! \mathbb{C} \setminus \tilde{\mathscr{J}}_{T}$; (ii) 
the boundary values $\tilde{\mathbb{M}}_{\pm}(z) \! := \! 
\lim_{\underset{z^{\prime} \, \in \, \pm \, \mathrm{side} \, \mathrm{of} 
\, \tilde{\mathscr{J}}_{T}}{z^{\prime} \to z \, \in \, \tilde{\mathscr{J}}_{T}}} 
\tilde{\mathbb{M}}(z^{\prime})$ satisfy the jump condition 
$\tilde{\mathbb{M}}_{+}(z) \! = \! \tilde{\mathbb{M}}_{-}(z) \tilde{\daleth}
(z)$ $\text{a.e.}$ $z \! \in \! \tilde{\mathscr{J}}_{T}$, where $\tilde{\daleth}
(z)$ is defined by Equation~\eqref{eqlmwayv1}; and (iii) for $\alpha_{p_{q}} 
\! \in \! (\tilde{a}_{j},\tilde{b}_{j})$, $q \! \in \! \lbrace 1,\dotsc,
\mathfrak{s} \! - \! 2,\mathfrak{s} \rbrace$, $j \! = \! 1,2,\dotsc,N$,
\begin{align*}
\tilde{\mathbb{M}}(z) \underset{\tilde{\mathcal{Y}}^{+} 
\supset \mathbb{C}_{+} \ni z \to \alpha_{k}}{=}& \, 
 \! + \! \mathcal{O}(z \! - \! \alpha_{p_{q}}), \quad 
q \! \in \! \lbrace 1,2,\dotsc,\mathfrak{s} \! - \! 2 \rbrace,
\end{align*}
where $\tilde{\varGamma}_{\triangle (\pm)}(\alpha_{p_{q}}) \! := \! 
\frac{1}{2}(\tilde{\gamma}(\alpha_{p_{q}}) \! \pm \! (\tilde{\gamma}
(\alpha_{p_{q}}))^{-1})$, $q \! \in \! \lbrace 1,\dotsc,\mathfrak{s} \! - \! 2,
\mathfrak{s} \rbrace$, $\tilde{\boldsymbol{u}}_{+}(\alpha_{p_{q}}) \! := \! 
\int_{\tilde{a}_{N+1}}^{\alpha_{p_{q}}^{+}} \tilde{\boldsymbol{\omega}}$, 
and $\tilde{\boldsymbol{u}}_{+}(\infty)$ $(= \! \tilde{\boldsymbol{u}}_{+}
(\alpha_{p_{\mathfrak{s}-1}}))$ $:= \! \int_{\tilde{a}_{N+1}}^{\infty^{+}} 
\tilde{\boldsymbol{\omega}}$, with $\infty^{+} \! = \! \alpha_{p_{
\mathfrak{s}-1}}^{+}$ denoting the point at infinity in $\tilde{
\mathcal{Y}}^{+}$ $(\supset \overline{\mathbb{C}}_{+})$, and, for 
$\alpha_{p_{q}} \! \in \! (-\infty,\tilde{b}_{0}) \cup (\tilde{a}_{N+1},
+\infty)$, $q \! \in \! \lbrace 1,\dotsc,\mathfrak{s} \! - \! 2,
\mathfrak{s} \rbrace$, one obtains the same $\tilde{\mathcal{Y}}^{
\pm} \supset \mathbb{C}_{\pm} \! \ni \! z \! \to \! \alpha_{p_{q}}$ 
asymptotics as above but with $\tilde{\Omega}_{j} \! \to \! 0$, $j \! 
= \! 1,2,\dotsc,N$ (of course, the $\tilde{\mathcal{Y}}^{\pm} \supset 
\overline{\mathbb{C}}_{\pm} \! \ni \! z \! \to \! \alpha_{p_{\mathfrak{s}-1}} 
\! = \! \infty$ asymptotics remain unchanged). Set, now, for $n \! \in \! 
\mathbb{N}$ and $k \! \in \! \lbrace 1,2,\dotsc,K \rbrace$ such that 
$\alpha_{p_{\mathfrak{s}}} \! := \! \alpha_{k} \! \neq \! \infty$, 
$\mathbb{M}(z) \! := \! \widetilde{\mathbb{K}} \, \tilde{\mathbb{M}}(z)$, 
where $\widetilde{\mathbb{K}}$ is defined in item~\pmb{(2)} of the 
lemma: via this definition, and the fact that, if $\alpha_{k} \! \in \! (-\infty,
\tilde{b}_{0}) \cup (\tilde{a}_{N+1},+\infty)$, $\mathscr{E}^{2} \! = \! 1$ 
and $\tilde{\Omega}_{j} \! \to \! 0$, $j \! = \! 1,2,\dotsc,N$, and, if 
$\alpha_{k} \! \in \! (\tilde{a}_{j},\tilde{b}_{j})$, $j \! = \! 1,2,\dotsc,N$, 
$\mathscr{E} \exp (-\mi ((n \! - \! 1)K \! + \! k) \tilde{\Omega}_{j}) 
\! = \! \mathscr{E}^{-1}$, one shows that $\mathbb{M} \colon 
\mathbb{C} \setminus \tilde{\mathscr{J}}_{T} \! \to \! \mathrm{SL}_{2}
(\mathbb{C})$ solves the RHP $(\mathbb{M}(z),\tilde{\daleth}(z),
\tilde{\mathscr{J}}_{T})$, with asymptotics
\begin{equation*}
\mathbb{M}(z) \underset{\tilde{\mathcal{Y}}^{+} \supset \mathbb{C}_{+} 
\ni z \to \alpha_{k}}{=} \mathscr{E}^{-\sigma_{3}} \! + \! \mathcal{O}
(z \! - \! \alpha_{k}), \qquad \qquad \mathbb{M}(z) \underset{\tilde{
\mathcal{Y}}^{-} \supset \mathbb{C}_{-} \ni z \to \alpha_{k}}{=} \mi 
\mathscr{E}^{\sigma_{3}} \sigma_{2} \! + \! \mathcal{O}(z \! - \! 
\alpha_{k}),
\end{equation*}
\begin{align*}
\mathbb{M}(z) \underset{\tilde{\mathcal{Y}}^{+} \supset 
\overline{\mathbb{C}}_{+} \ni z \to \alpha_{p_{\mathfrak{s}-1}} 
= \infty}{=}& \, \widetilde{\mathbb{K}} 
 \! + \! \mathcal{O}(z \! - \! \alpha_{p_{q}}), \quad 
q \! \in \! \lbrace 1,2,\dotsc,\mathfrak{s} \! - \! 2 \rbrace.
\end{align*}
(Note: if $\alpha_{p_{q}} \! \in \! (-\infty,\tilde{b}_{0}) \cup (\tilde{a}_{N+1},
+\infty)$, $q \! \in \! \lbrace 1,2,\dotsc,\mathfrak{s} \! - \! 2 \rbrace$, 
then one obtains the same $\tilde{\mathcal{Y}}^{\pm} \supset 
\mathbb{C}_{\pm} \! \ni \! z \! \to \! \alpha_{p_{q}}$ asymptotics as above 
but with $\tilde{\Omega}_{j} \! \to \! 0$, $j \! = \! 1,2,\dotsc,N$; of course, 
the $\tilde{\mathcal{Y}}^{\pm} \supset \overline{\mathbb{C}}_{\pm} \! \ni \! 
z \! \to \! \alpha_{p_{\mathfrak{s}-1}} \! = \! \infty$ asymptotics remain 
unchanged.) Finally, for $n \! \in \! \mathbb{N}$ and $k \! \in \! \lbrace 
1,2,\dotsc,K \rbrace$ such that $\alpha_{p_{\mathfrak{s}}} \! := \! \alpha_{k} 
\! \neq \! \infty$, defining $\mathfrak{m}(z)$ in terms of $\mathbb{M}(z)$ 
as in item~\pmb{(2)} of the lemma, a straightforward calculation shows 
that $\mathfrak{m} \colon \mathbb{C} \setminus \tilde{\mathscr{J}} 
\! \to \! \mathrm{SL}_{2}(\mathbb{C})$ solves the model RHP 
$(\mathfrak{m}(z),\daleth (z),\tilde{\mathscr{J}})$ stated in 
item~\pmb{(2)} of Lemma~\ref{lem4.3}. For $n \! \in \! \mathbb{N}$ 
and $k \! \in \! \lbrace 1,2,\dotsc,K \rbrace$ such that 
$\alpha_{p_{\mathfrak{s}}} \! := \! \alpha_{k} \! \neq \! \infty$, one 
notes {}from the formula for $\mathbb{M}(z)$ stated in item~\pmb{(2)} 
of the lemma that it is well defined for $\mathbb{C} \setminus 
\mathbb{R} \! \ni \! z$; in particular, it is analytic for $z \! \in \! 
\mathbb{C} \setminus \tilde{\mathscr{J}}_{T}$ (independent of the 
path in $\mathbb{C} \setminus \tilde{\mathscr{J}}_{T}$ chosen to 
evaluate $\tilde{\boldsymbol{u}}(z) \! = \! \int_{\tilde{a}_{N+1}}^{z} 
\tilde{\boldsymbol{\omega}})$. Moreover, since $\lbrace \mathstrut 
z^{\prime} \! \in \! \tilde{\mathcal{Y}}^{\pm} \, (\supset \mathbb{C}_{\pm}); 
\, \tilde{\boldsymbol{\theta}}(\tilde{\boldsymbol{u}}(z^{\prime}) \! \mp 
\! \tilde{\boldsymbol{d}}) \! = \! 0 \rbrace \! = \! \lbrace \tilde{z}_{j}^{\pm} 
\rbrace_{j=1}^{N} \! = \! \lbrace \mathstrut z^{\prime} \! \in \! \mathbb{C}_{
\pm} \, (\subset \tilde{\mathcal{Y}}^{\pm}); \, \tilde{\gamma}(z^{\prime}) 
\! \mp \! (\tilde{\gamma}(z^{\prime}))^{-1} \! = \! 0 \rbrace$, one notes 
that the simple poles of $\mathfrak{m}_{11}(z)$ and $\mathfrak{m}_{22}
(z)$ (resp., $\mathfrak{m}_{12}(z)$ and $\mathfrak{m}_{21}(z))$, that 
is, $\lbrace \mathstrut z^{\prime} \! \in \! \tilde{\mathcal{Y}}^{-} 
\, (\supset \mathbb{C}_{-}); \, \tilde{\boldsymbol{\theta}}
(\tilde{\boldsymbol{u}}(z^{\prime}) \! + \! \tilde{\boldsymbol{d}}) 
\! = \! 0 \rbrace$ (resp., $\lbrace \mathstrut z^{\prime} \! \in 
\! \tilde{\mathcal{Y}}^{+} \, (\supset \mathbb{C}_{+}); \, \tilde{
\boldsymbol{\theta}}(\tilde{\boldsymbol{u}}(z^{\prime}) \! - \! \tilde{
\boldsymbol{d}}) \! = \! 0 \rbrace)$, are exactly cancelled by the simple 
zeros of $\tilde{\gamma}(z) \! + \! (\tilde{\gamma}(z))^{-1} \! = \! 0$ 
(resp., $\tilde{\gamma}(z) \! - \! (\tilde{\gamma}(z))^{-1} \! = \! 0)$; 
hence, $\mathbb{M}(z)$ only has $\tfrac{1}{4}$-root singularities at the 
associated end-points of the intervals, $\tilde{b}_{j-1},\tilde{a}_{j}$, $j \! 
= \! 1,2,\dotsc,N \! + \! 1$, of the support, $J_{f}$, of the equilibrium 
measure, $\mu_{\widetilde{V}}^{f}$. {}From the definition of $\mathfrak{m}
(z)$ in terms of $\mathbb{M}(z)$ (and the explicit formula for $\mathbb{M}
(z))$ stated in item~\pmb{(2)} of the lemma, upon recalling that 
$\mathfrak{m}(z)$ solves the model RHP stated in item~\pmb{(2)} 
of Lemma~\ref{lem4.3}, one learns that, as $\det (\daleth (z)) \! = \! 
1$, $\det (\mathfrak{m}_{+}(z)) \! = \! \det (\mathfrak{m}_{-}(z))$ 
(that is, $\det (\mathfrak{m}(z))$ has no `jumps'), thus $\det 
(\mathfrak{m}(z))$ has, at worst, isolated $\tfrac{1}{2}$-root 
singularities at $\tilde{b}_{j-1},\tilde{a}_{j}$, $j \! = \! 1,2,
\dotsc,N \! + \! 1$, that are removable, which, in turn, implies 
that $\det (\mathfrak{m}(z))$ is entire and bounded; in fact, 
a---tedious---calculation shows that, for $n \! \in \! 
\mathbb{N}$ and $k \! \in \! \lbrace 1,2,\dotsc,K \rbrace$ 
such that $\alpha_{p_{\mathfrak{s}}} \! := \! \alpha_{k} \! 
\neq \! \infty$ (see Section~\ref{sek5}, the proof of 
Proposition~\ref{propo5.1}),\footnote{All square roots are positive.}
\begin{align*}
\det (\mathfrak{m}(z)) \underset{z \to \tilde{a}_{N+1}}{=}& \, 
\dfrac{1}{\mathfrak{k}_{0}} \dfrac{\tilde{\boldsymbol{\theta}}
(-\frac{1}{2 \pi}((n \! - \! 1)K \! + \! k) \tilde{\boldsymbol{
\Omega}} \! + \! \tilde{\boldsymbol{d}})}{\tilde{\boldsymbol{
\theta}}(\tilde{\boldsymbol{d}})} \dfrac{\tilde{\boldsymbol{\theta}}
(-\frac{1}{2 \pi}((n \! - \! 1)K \! + \! k) \tilde{\boldsymbol{\Omega}} 
\! - \! \tilde{\boldsymbol{d}})}{\tilde{\boldsymbol{\theta}}
(-\tilde{\boldsymbol{d}})} \\
\times& \, \left(1 \! + \! \dfrac{\mi}{2}(\tilde{a}_{N+1} \! - \! 
\tilde{b}_{0})^{1/2}(\tilde{\aleph}^{1}_{1}(\tilde{a}_{N+1}) \! - 
\! \tilde{\aleph}^{1}_{-1}(\tilde{a}_{N+1})) \prod_{m=1}^{N} 
\dfrac{(\tilde{a}_{N+1} \! - \! \tilde{b}_{m})^{1/2}}{(\tilde{a}_{N+1} 
\! - \! \tilde{a}_{m})^{1/2}} \right) \! + \! \mathcal{O}
(z \! - \! \tilde{a}_{N+1}), \\
\det (\mathfrak{m}(z)) \underset{z \to \tilde{b}_{0}}{=}& \, 
\dfrac{1}{\mathfrak{k}_{0}} \dfrac{\tilde{\boldsymbol{\theta}}
(\tilde{\boldsymbol{u}}_{+}(\tilde{b}_{0}) \! - \! \frac{1}{2 \pi}
((n \! - \! 1)K \! + \! k) \tilde{\boldsymbol{\Omega}} \! + \! 
\tilde{\boldsymbol{d}})}{\tilde{\boldsymbol{\theta}}(\tilde{
\boldsymbol{u}}_{+}(\tilde{b}_{0}) \! + \! \tilde{\boldsymbol{d}})} 
\dfrac{\tilde{\boldsymbol{\theta}}(\tilde{\boldsymbol{u}}_{+}
(\tilde{b}_{0}) \! - \! \frac{1}{2 \pi}((n \! - \! 1)K \! + \! k) 
\tilde{\boldsymbol{\Omega}} \! - \! \tilde{\boldsymbol{d}})}{
\tilde{\boldsymbol{\theta}}(\tilde{\boldsymbol{u}}_{+}
(\tilde{b}_{0}) \! - \! \tilde{\boldsymbol{d}})} \\
\times& \, \left(1 \! + \! \dfrac{\mi}{2} \me^{\frac{\mi \pi}{2}}
(\tilde{a}_{N+1} \! - \! \tilde{b}_{0})^{1/2}(\tilde{\aleph}^{1}_{1}
(\tilde{b}_{0}) \! - \! \tilde{\aleph}^{1}_{-1}(\tilde{b}_{0})) 
\prod_{m=1}^{N} \dfrac{(\tilde{a}_{m} \! - \! \tilde{b}_{0})^{
1/2}}{(\tilde{b}_{m} \! - \! \tilde{b}_{0})^{1/2}} \right) \! + \! 
\mathcal{O}(z \! - \! \tilde{b}_{0}),
\end{align*}
and, for $j \! = \! 1,2,\dotsc,N$,
\begin{align*}
\det (\mathfrak{m}(z)) \underset{z \to \tilde{a}_{j}}{=}& \, 
\dfrac{\me^{\mi ((n-1)K+k) \tilde{\Omega}_{j}}}{\mathfrak{k}_{0}} 
\dfrac{\tilde{\boldsymbol{\theta}}(\tilde{\boldsymbol{u}}_{+}
(\tilde{a}_{j}) \! - \! \frac{1}{2 \pi}((n \! - \! 1)K \! + \! k) \tilde{
\boldsymbol{\Omega}} \! + \! \tilde{\boldsymbol{d}})}{\tilde{
\boldsymbol{\theta}}(\tilde{\boldsymbol{u}}_{+}(\tilde{a}_{j}) \! 
+ \! \tilde{\boldsymbol{d}})} \dfrac{\tilde{\boldsymbol{\theta}}
(\tilde{\boldsymbol{u}}_{+}(\tilde{a}_{j}) \! - \! \frac{1}{2 \pi}
((n \! - \! 1)K \! + \! k) \tilde{\boldsymbol{\Omega}} \! - \! 
\tilde{\boldsymbol{d}})}{\tilde{\boldsymbol{\theta}}(\tilde{
\boldsymbol{u}}_{+}(\tilde{a}_{j}) \! - \! \tilde{\boldsymbol{d}})} \\
\times& \, \left(1 \! + \! \dfrac{\mi}{2} \dfrac{(\tilde{a}_{j} \! 
- \! \tilde{b}_{0})^{1/2}(\tilde{b}_{j} \! - \! \tilde{a}_{j})^{1/2}}{
(\tilde{a}_{N+1} \! - \! \tilde{a}_{j})^{1/2}}(\tilde{\aleph}^{1}_{1}
(\tilde{a}_{j}) \! - \! \tilde{\aleph}^{1}_{-1}(\tilde{a}_{j})) \prod_{m
=1}^{j-1} \dfrac{(\tilde{a}_{j} \! - \! \tilde{b}_{m})^{1/2}}{(\tilde{a}_{j} 
\! - \! \tilde{a}_{m})^{1/2}} \prod_{m^{\prime}=j+1}^{N} 
\dfrac{(\tilde{b}_{m^{\prime}} \! - \! \tilde{a}_{j})^{1/2}}{
(\tilde{a}_{m^{\prime}} \! - \! \tilde{a}_{j})^{1/2}} \right) 
\! + \! \mathcal{O}(z \! - \! \tilde{a}_{j}), \\
\det (\mathfrak{m}(z)) \underset{z \to \tilde{b}_{j}}{=}& \, 
\dfrac{\me^{\mi ((n-1)K+k) \tilde{\Omega}_{j}}}{\mathfrak{k}_{0}} 
\dfrac{\tilde{\boldsymbol{\theta}}(\tilde{\boldsymbol{u}}_{+}
(\tilde{b}_{j}) \! - \! \frac{1}{2 \pi}((n \! - \! 1)K \! + \! k) \tilde{
\boldsymbol{\Omega}} \! + \! \tilde{\boldsymbol{d}})}{\tilde{
\boldsymbol{\theta}}(\tilde{\boldsymbol{u}}_{+}(\tilde{b}_{j}) \! 
+ \! \tilde{\boldsymbol{d}})} \dfrac{\tilde{\boldsymbol{\theta}}
(\tilde{\boldsymbol{u}}_{+}(\tilde{b}_{j}) \! - \! \frac{1}{2 \pi}
((n \! - \! 1)K \! + \! k) \tilde{\boldsymbol{\Omega}} \! - \! 
\tilde{\boldsymbol{d}})}{\tilde{\boldsymbol{\theta}}(\tilde{
\boldsymbol{u}}_{+}(\tilde{b}_{j}) \! - \! \tilde{\boldsymbol{d}})} \\
\times& \, \left(1 \! + \! \dfrac{\mi}{2} \me^{\frac{\mi \pi}{2}} 
\dfrac{(\tilde{a}_{N+1} \! - \! \tilde{b}_{j})^{1/2}(\tilde{b}_{j} \! 
- \! \tilde{a}_{j})^{1/2}}{(\tilde{b}_{j} \! - \! \tilde{b}_{0})^{1/2}}
(\tilde{\aleph}^{1}_{1}(\tilde{b}_{j}) \! - \! \tilde{\aleph}^{1}_{-1}
(\tilde{b}_{j})) \prod_{m=1}^{j-1} \dfrac{(\tilde{b}_{j} \! - \! 
\tilde{a}_{m})^{1/2}}{(\tilde{b}_{j} \! - \! \tilde{b}_{m})^{1/2}} 
\prod_{m^{\prime}=j+1}^{N} \dfrac{(\tilde{a}_{m^{\prime}} \! - \! 
\tilde{b}_{j})^{1/2}}{(\tilde{b}_{m^{\prime}} \! - \! \tilde{b}_{j})^{1/2}} 
\right) \! + \! \mathcal{O}(z \! - \! \tilde{b}_{j}),
\end{align*}
where $\tilde{\boldsymbol{u}}_{+}(\tilde{a}_{j}) \! := \! \int_{
\tilde{a}_{N+1}}^{\tilde{a}_{j}^{+}} \tilde{\boldsymbol{\omega}}$ 
and $\tilde{\boldsymbol{u}}_{+}(\tilde{b}_{j-1}) \! := \! \int_{
\tilde{a}_{N+1}}^{\tilde{b}_{j-1}^{+}} \tilde{\boldsymbol{\omega}}$, 
$j \! = \! 1,2,\dotsc,N \! + \! 1$ (of course, $\tilde{\boldsymbol{
u}}_{+}(\tilde{a}_{N+1}) \! \equiv \! (0,0,\dotsc,0)$ $(\in \! 
\operatorname{Jac}(\tilde{\mathcal{Y}})))$, and $\tilde{\aleph}^{
\varepsilon_{1}}_{\varepsilon_{2}}(\varsigma)$, $\varepsilon_{1},
\varepsilon_{2} \! = \! \pm 1$, is given by Equations~\eqref{eqmainfin32} 
and~\eqref{eqmainfin34}--\eqref{eqmainfin60}. Thus, via a generalisation 
of Liouville's Theorem and the asymptotics $\det (\mathfrak{m}
(z)) \! =_{\tilde{\mathcal{Y}}^{\pm} \supset \mathbb{C}_{\pm} \ni 
z \to \alpha_{k}} \! 1 \! + \! \mathcal{O}(z \! - \! \alpha_{k})$, one 
arrives at, for $n \! \in \! \mathbb{N}$ and $k \! \in \! \lbrace 1,2,
\dotsc,K \rbrace$ such that $\alpha_{p_{\mathfrak{s}}} \! := \! 
\alpha_{k} \! \neq \! \infty$, $\det (\mathfrak{m}(z)) \! = \! 1$. 
{}From the definition of $\mathfrak{m}(z)$ in terms of $\mathbb{M}
(z)$ (and the explicit formula for $\mathbb{M}(z))$ stated in 
item~\pmb{(2)} of the lemma, it follows that, for $n \! \in \! 
\mathbb{N}$ and $k \! \in \! \lbrace 1,2,\dotsc,K \rbrace$ such 
that $\alpha_{p_{\mathfrak{s}}} \! := \! \alpha_{k} \! \neq \! \infty$, 
both $\mathfrak{m}(z)$ and $(\mathfrak{m}(z))^{-1}$ are uniformly 
bounded functions of $n$ (and $z_{o})$, in the double-scaling limit 
$\mathscr{N},n \! \to \! \infty$ such that $z_{o} \! = \! 1 \! + \! o(1)$, 
for $z$ in compact subsets disjoint {}from $\lbrace \tilde{b}_{j-1},
\tilde{a}_{j} \rbrace_{j=1}^{N+1}$. For $n \! \in \! \mathbb{N}$ 
and $k \! \in \! \lbrace 1,2,\dotsc,K \rbrace$ such that 
$\alpha_{p_{\mathfrak{s}}} \! := \! \alpha_{k} \! \neq \! 
\infty$, let $\mathbb{M}^{\blacklozenge} 
\colon \overline{\mathbb{C}} \setminus \overline{\mathbb{R}} \! 
\to \! \mathrm{SL}_{2}(\mathbb{C})$ be another solution of the RHP 
$(\mathbb{M}(z),\tilde{\daleth}(z),\overline{\mathbb{R}})$, and 
set $\tilde{\Delta}^{o}(z) \! := \! \mathbb{M}^{\blacklozenge}(z)
(\mathbb{M}(z))^{-1}$: then, $\tilde{\Delta}^{o}_{+}(z) \! = \! 
\mathbb{M}^{\blacklozenge}_{+}(z)(\mathbb{M}_{+}(z))^{-1} 
\! = \! \mathbb{M}^{\blacklozenge}_{-}(z) \tilde{\daleth}(z)
(\mathbb{M}_{-}(z) \tilde{\daleth}(z))^{-1} \! = \! \mathbb{M}^{
\blacklozenge}_{-}(z)(\mathbb{M}_{-}(z))^{-1} \! = \! \tilde{
\Delta}^{o}_{-}(z)$ $\text{a.e.}$ $z \! \in \! \overline{\mathbb{R}}$; 
moreover, $\tilde{\Delta}^{o}(z)$ has, at worst, $\mathcal{L}^{1}_{
\mathrm{M}_{2}(\mathbb{C})}(\ast)$-singularities at the associated 
end-points of the intervals, $\tilde{b}_{j-1},\tilde{a}_{j}$, $j \! = \! 
1,2,\dotsc,N \! + \! 1$, of the support, $J_{f}$, of the equilibrium 
measure, $\mu_{\widetilde{V}}^{f}$, which, as per the discussion 
above, are removable; hence, noting that $\tilde{\Delta}^{o}(z) 
\! \to \! \mathrm{I}$ as $(\tilde{\mathcal{Y}}^{\pm} \supset 
\mathbb{C}_{\pm} \! \ni)$ $z \! \to \! \alpha_{p_{q}}$, $q \! = \! 
1,2,\dotsc,\mathfrak{s}$, one concludes that, for $n \! \in \! 
\mathbb{N}$ and $k \! \in \! \lbrace 1,2,\dotsc,K \rbrace$ such 
that $\alpha_{p_{\mathfrak{s}}} \! := \! \alpha_{k} \! \neq \! \infty$, 
$\tilde{\Delta}^{o}(z) \! = \! \mathrm{I}$ $(\Rightarrow \! \mathbb{
M}^{\blacklozenge}(z) \! = \! \mathbb{M}(z))$. \hfill $\qed$

In order to prove that, for $n \! \in \! \mathbb{N}$ and $k \! \in \! 
\lbrace 1,2,\dotsc,K \rbrace$ such that $\alpha_{p_{\mathfrak{s}}} \! 
:= \! \alpha_{k} \! = \! \infty$ (resp., $\alpha_{p_{\mathfrak{s}}} \! 
:= \! \alpha_{k} \! \neq \! \infty)$, there is a solution of the matrix 
RHP $(\mathfrak{M}(z),\mathfrak{v}(z),\hat{\Sigma})$ (resp., 
$(\mathfrak{M}(z),\mathfrak{v}(z),\tilde{\Sigma}))$ stated in 
item~\pmb{(1)} (resp., item~\pmb{(2)}) of Lemma~\ref{lem4.2} 
`close to' the associated parametrix, one needs to establish that the 
parametrix is \emph{uniformly} bounded, that is, by certain general 
theorems (see, for example, \cite{kcligoh}; see, also, Chapter~7 in 
\cite{a51}, and Appendices A and B in \cite{skmclpdm}), one needs to 
show that $\mathfrak{v}(z) \! \to \! \daleth (z)$ in the double-scaling 
limit $\mathscr{N},n \! \to \! \infty$ such that $z_{o} \! = \! 1 \! + \! 
o(1)$ for $z \! \in \! \hat{\Sigma}$ (resp., $z \! \in \! \tilde{\Sigma})$ 
in the $\mathcal{L}^{2}_{\mathrm{M}_{2}(\mathbb{C})}(\hat{\Sigma}) 
\cap \mathcal{L}^{\infty}_{\mathrm{M}_{2}(\mathbb{C})}(\hat{\Sigma})$ 
(resp., $\mathcal{L}^{2}_{\mathrm{M}_{2}(\mathbb{C})}(\tilde{\Sigma}) 
\cap \mathcal{L}^{\infty}_{\mathrm{M}_{2}(\mathbb{C})}(\tilde{\Sigma}))$ 
sense, in other words,
\begin{equation*}
\lim_{\underset{z_{o}=1+o(1)}{\mathscr{N},n \to \infty}} 
\norm{\mathfrak{v}(\pmb{\cdot}) \! - \! \daleth (\pmb{\cdot})}_{
\cap_{p=2,\infty} \mathcal{L}^{p}_{\mathrm{M}_{2}(\mathbb{C})}
(\overset{\ast}{\Sigma})} \! := \! \lim_{\underset{z_{o}=1+o(1)}{
\mathscr{N},n \to \infty}} \sum_{p=2,\infty} \norm{\mathfrak{v}
(\pmb{\cdot}) \! - \! \daleth (\pmb{\cdot})}_{\mathcal{L}^{p}_{
\mathrm{M}_{2}(\mathbb{C})}(\overset{\ast}{\Sigma})} \! = \! 0, 
\quad \ast \! \in \! \lbrace \, \hat \, , \, \tilde \, \rbrace,
\end{equation*}
uniformly; however, for regular $\widetilde{V} \colon \overline{\mathbb{R}} 
\setminus \lbrace \alpha_{1},\alpha_{2},\dotsc,\alpha_{K} \rbrace \! 
\to \! \mathbb{R}$ satisfying conditions~\eqref{eq20}--\eqref{eq22}, 
since, for $n \! \in \! \mathbb{N}$ and $k \! \in \! \lbrace 1,2,\dotsc,
K \rbrace$ such that $\alpha_{p_{\mathfrak{s}}} \! := \! \alpha_{k} \! 
= \! \infty$ (resp., $\alpha_{p_{\mathfrak{s}}} \! := \! \alpha_{k} \! 
\neq \! \infty)$, the strict inequalities $g^{\infty}_{+}(z) \! + \! 
g^{\infty}_{-}(z) \! - \! 2 \tilde{\mathscr{P}}_{0} \! - \! \widetilde{V}(z) 
\! - \! \hat{\ell} \! < \! 0$, $z \! \in \! (-\infty,\hat{b}_{0}) \cup 
(\hat{a}_{N+1},+\infty) \cup \cup_{i=1}^{N}(\hat{a}_{i},\hat{b}_{i})$, 
and $\pm \Re (\mi \int_{z}^{\hat{a}_{N+1}} \psi_{\widetilde{V}}^{\infty}
(\xi) \, \md \xi) \! > \! 0$, $z \! \in \! \mathbb{C}_{\pm} \cap (\cup_{j=
1}^{N+1} \hat{\mathbb{U}}_{j})$ (resp., $g^{f}_{+}(z) \! + \! g^{f}_{-}(z) 
\! - \! \hat{\mathscr{P}}^{+}_{0} \! - \! \hat{\mathscr{P}}^{-}_{0} \! - \! 
\widetilde{V}(z) \! - \! \tilde{\ell} \! < \! 0$, $z \! \in \! (-\infty,\tilde{
b}_{0}) \cup (\tilde{a}_{N+1},+\infty) \cup \cup_{i=1}^{N}(\tilde{a}_{i},
\tilde{b}_{i})$, and $\pm \Re (\mi \int_{z}^{\tilde{a}_{N+1}} \psi_{
\widetilde{V}}^{f}(\xi) \, \md \xi) \! > \! 0$, $z \! \in \! \mathbb{C}_{\pm} 
\cap (\cup_{j=1}^{N+1} \tilde{\mathbb{U}}_{j}))$, fail at the end-points 
of the intervals, $\hat{b}_{j-1},\hat{a}_{j}$ (resp., $\tilde{b}_{j-1},
\tilde{a}_{j})$, $j \! = \! 1,2,\dotsc,N \! + \! 1$, of the support, 
$J_{\infty}$ (resp., $J_{f})$, of the associated equilibrium measure, 
$\mu_{\widetilde{V}}^{\infty}$ (resp., $\mu_{\widetilde{V}}^{f})$, this 
implies that $\mathfrak{v}(z) \! \to \! \daleth (z)$ in the double-scaling 
limit $\mathscr{N},n \! \to \! \infty$ such that $z_{o} \! = \! 1 \! + \! o(1)$ 
\emph{pointwise}, but not uniformly, for $z \! \in \! \hat{\Sigma}$ 
(resp., $z \! \in \! \tilde{\Sigma})$, in which case, one can not conclude 
that $\mathfrak{M}(z) \! \to \! \mathfrak{m}(z)$ in the double-scaling 
limit $\mathscr{N},n \! \to \! \infty$ such that $z_{o} \! = \! 1 \! + \! o(1)$ 
uniformly for $z \! \in \! \hat{\Sigma}$ (resp., $z \! \in \! \tilde{\Sigma})$.
The resolution of this lack of uniformity constitutes, therefore, the essential 
technical obstacle remaining for the analysis of the RHP $(\mathfrak{M}(z),
\mathfrak{v}(z),\hat{\Sigma})$ (resp., $(\mathfrak{M}(z),\mathfrak{v}(z),
\tilde{\Sigma}))$ stated in item~\pmb{(1)} (resp., item~\pmb{(2)}) of 
Lemma~\ref{lem4.2}, and a substantial part of the remaing analysis is 
devoted to overcoming this problem.

For the ensuing discussion, consider, say, and without loss of generality, 
the case $n \! \in \! \mathbb{N}$ and $k \! \in \! \lbrace 1,2,\dotsc,K 
\rbrace$ such that $\alpha_{p_{\mathfrak{s}}} \! := \! \alpha_{k} \! \neq 
\! \infty$ (the case $n \! \in \! \mathbb{N}$ and $k \! \in \! \lbrace 1,2,
\dotsc,K \rbrace$ such that $\alpha_{p_{\mathfrak{s}}} \! := \! \alpha_{k} 
\! = \! \infty$ is analogous). The key ingredient necessary to control the 
above-mentioned technical difficulty is to construct parametrices for the 
solution of the associated RHP $(\mathfrak{M}(z),\mathfrak{v}(z),\tilde{
\Sigma})$ (stated in item~\pmb{(2)} of Lemma~\ref{lem4.2}) in `small' 
(open) neighbourhoods about $\tilde{b}_{j-1},\tilde{a}_{j}$, $j \! = \! 1,
2,\dotsc,N \! + \! 1$, where the convergence of $\mathfrak{v}(z)$ to 
$\daleth (z)$, in the double-scaling limit $\mathscr{N},n \! \to \! \infty$ 
such that $z_{o} \! = \! 1 \! + \! o(1)$, is not uniform, in such a way that 
on the boundaries of these neighbourhoods the parametrices `match' 
with the solution of the associated model RHP, $\mathfrak{m}(z)$, up to 
terms that are $o(1)$, in the double-scaling limit $\mathscr{N},n \! \to \! 
\infty$ such that $z_{o} \! = \! 1 \! + \! o(1)$; moreover, since, for regular 
$\widetilde{V} \colon \overline{\mathbb{R}} \setminus \lbrace \alpha_{1},
\alpha_{2},\dotsc,\alpha_{K} \rbrace \! \to \! \mathbb{R}$ satisfying 
conditions~\eqref{eq20}--\eqref{eq22}, $\psi_{\widetilde{V}}^{f}(x) \! 
=_{x \downarrow \tilde{b}_{j-1}} \! \mathcal{O}((x \! - \! \tilde{b}_{j-1})^{
1/2})$ and $\psi_{\widetilde{V}}^{f}(x) \! =_{x \uparrow \tilde{a}_{j}} \! 
\mathcal{O}((\tilde{a}_{j} \! - \! x)^{1/2})$, $j \! = \! 1,2,\dotsc,N \! + \! 
1$, it is well known \cite{a53} that the parametrices can be expressed in 
terms of Airy functions.\footnote{The general methodology used to construct 
such parametrices is a Vanishing Lemma \cite{rbspdct}.} More precisely, 
for $n \! \in \! \mathbb{N}$ and $k \! \in \! \lbrace 1,2,\dotsc,K \rbrace$ 
such that $\alpha_{p_{\mathfrak{s}}} \! := \! \alpha_{k} \! \neq \! \infty$, 
one surrounds $\tilde{b}_{j-1},\tilde{a}_{j}$, $j \! = \! 1,2,\dotsc,N \! + \! 
1$, by `small', mutually disjoint open discs $\mathscr{O}_{\tilde{\epsilon}^{
b}_{j}}(\tilde{b}_{j-1}) \! := \! \lbrace \mathstrut z \! \in \! \mathbb{C}; \, 
\lvert z \! - \! \tilde{b}_{j-1} \rvert \! < \! \tilde{\epsilon}^{b}_{j} \rbrace$ 
and $\mathscr{O}_{\tilde{\epsilon}^{a}_{j}}(\tilde{a}_{j}) \! := \! \lbrace 
\mathstrut z \! \in \! \mathbb{C}; \, \lvert z \! - \! \tilde{a}_{j} \rvert \! < \! 
\tilde{\epsilon}^{a}_{j} \rbrace$, $j \! = \! 1,2,\dotsc,N \! + \! 1$, respectively, 
where $\tilde{\epsilon}_{j}^{b},\tilde{\epsilon}_{j}^{a}$ are arbitrarily fixed, 
sufficiently small positive real numbers chosen so that $\mathscr{O}_{
\tilde{\epsilon}^{b}_{i}}(\tilde{b}_{i-1}) \cap \mathscr{O}_{\tilde{\epsilon}^{
a}_{j}}(\tilde{a}_{j}) \! = \! \varnothing$ $\forall$ $i,j \! \in \! \lbrace 1,2,
\dotsc,N \! + \! 1 \rbrace$ and $\mathscr{O}_{\tilde{\epsilon}^{b}_{m}}
(\tilde{b}_{m-1}) \cap \lbrace \alpha_{p_{1}},\dotsc,\alpha_{p_{\mathfrak{s}-2}},
\alpha_{p_{\mathfrak{s}}} \rbrace \! = \! \varnothing \! = \! \mathscr{O}_{
\tilde{\epsilon}^{a}_{m}}(\tilde{a}_{m}) \cap \lbrace \alpha_{p_{1}},\dotsc,
\alpha_{p_{\mathfrak{s}-2}},\alpha_{p_{\mathfrak{s}}} \rbrace$ (of course, 
$\mathscr{O}_{\tilde{\epsilon}^{b}_{m}}(\tilde{b}_{m-1}) \cap \lbrace 
\alpha_{p_{\mathfrak{s}-1}} \! = \! \infty \rbrace \! = \! \varnothing \! = \! 
\mathscr{O}_{\tilde{\epsilon}^{a}_{m}}(\tilde{a}_{m}) \cap \lbrace \alpha_{p_{
\mathfrak{s}-1}} \! = \! \infty \rbrace)$, $m \! = \! 1,2,\dotsc,N \! + \! 1$, 
and defines $\tilde{S}_{p}(z)$, the associated parametrix for $\mathfrak{M}
(z)$, to be equal to $\mathfrak{m}(z)$ for $z \! \in \! \mathbb{C} 
\setminus \cup_{j=1}^{N+1}(\mathscr{O}_{\tilde{\epsilon}^{b}_{j}}
(\tilde{b}_{j-1}) \cup \mathscr{O}_{\tilde{\epsilon}^{a}_{j}}(\tilde{a}_{j}))$, and 
by $m_{p}(z)$ (which is to be determined!) for $z \! \in \! \cup_{j=1}^{N+1}
(\mathscr{O}_{\tilde{\epsilon}^{b}_{j}}(\tilde{b}_{j-1}) \cup \mathscr{O}_{
\tilde{\epsilon}^{a}_{j}}(\tilde{a}_{j}))$, and solves the associated local matrix 
RHP for $m_{p}(z)$ on $\cup_{j=1}^{N+1}(\mathscr{O}_{\tilde{\epsilon}^{b}_{j}}
(\tilde{b}_{j-1}) \cup \mathscr{O}_{\tilde{\epsilon}^{a}_{j}}(\tilde{a}_{j}))$ in 
such a way that $m_{p}(z) \! \approx_{\underset{z_{o}=1+o(1)}{\mathscr{N},
n \to \infty}} \! \mathfrak{m}(z)$ for $z \! \in \! \cup_{j=1}^{N+1}
(\partial \mathscr{O}_{\tilde{\epsilon}^{b}_{j}}(\tilde{b}_{j-1}) \cup 
\partial \mathscr{O}_{\tilde{\epsilon}^{a}_{j}}(\tilde{a}_{j}))$, in which case, 
$\tilde{\mathcal{R}}(z) \! := \! \mathfrak{M}(z)(\tilde{S}_{p}(z))^{-1} 
\colon \mathbb{C} \setminus \tilde{\Sigma}_{\tilde{\mathcal{R}}} \! \to 
\! \operatorname{SL}_{2}(\mathbb{C})$, where $\tilde{\Sigma}_{\tilde{
\mathcal{R}}} \! := \! \tilde{\Sigma} \cup \cup_{j=1}^{N+1}(\partial 
\mathscr{O}_{\tilde{\epsilon}^{b}_{j}}(\tilde{b}_{j-1}) \cup \partial 
\mathscr{O}_{\tilde{\epsilon}^{a}_{j}}(\tilde{a}_{j}))$, solves the matrix RHP 
$(\tilde{\mathcal{R}}(z),\tilde{v}_{\tilde{\mathcal{R}}}(z),\tilde{\Sigma}_{
\tilde{\mathcal{R}}})$ with $\norm{\tilde{v}_{\tilde{\mathcal{R}}}(\pmb{\cdot}) 
\! - \! \mathrm{I}}_{\cap_{p=2,\infty} \mathcal{L}^{p}_{\mathrm{M}_{2}
(\mathbb{C})}(\tilde{\Sigma}_{\tilde{\mathcal{R}}})} \! =_{\underset{z_{o}=
1+o(1)}{\mathscr{N},n \to \infty}} \! o(1)$ uniformly; in particular, the error 
term, which is $o(1)$ in the double-scaling limit $\mathscr{N},n \! \to \! 
\infty$ such that $z_{o} \! = \! 1 \! + \! o(1)$, is uniform in $\cap_{p=1,2,
\infty} \mathcal{L}^{p}_{\mathrm{M}_{2}(\mathbb{C})}(\tilde{\Sigma}_{\tilde{
\mathcal{R}}})$. By general Riemann-Hilbert techniques (see, for example, 
\cite{kcligoh}), $\tilde{\mathcal{R}}(z)$, and, thus, $\mathfrak{M}(z)$, can 
be obtained in the double-scaling limit $\mathscr{N},n \! \to \! \infty$ 
such that $z_{o} \! = \! 1 \! + \! o(1)$ via a Neumann series expansion of 
the associated resolvent kernel.

A detailed exposition for the construction of parametrices of the 
above-mentioned type can be found, for example, in \cite{a51,a54,a61,
pz22,a53,a67,a59}; rather than regurgitating, \emph{verbatim}, much of 
the analysis that can be located in these latter references (by no means 
an exhaustive list!), the point of view taken here is that one follows, 
\emph{mutatis mutandis}, the scheme presented therein to arrive at, 
for $n \! \in \! \mathbb{N}$ and $k \! \in \! \lbrace 1,2,\dotsc,K 
\rbrace$ such that $\alpha_{p_{\mathfrak{s}}} \! := \! \alpha_{k} \! = 
\! \infty$ (resp., $\alpha_{p_{\mathfrak{s}}} \! := \! \alpha_{k} \! \neq \! 
\infty)$, the results stated in Lemmata~\ref{lem4.6} and~\ref{lem4.7} (resp., 
Lemmata~\ref{lem4.8} and~\ref{lem4.9}) below for the parametrix of the 
RHP $(\mathfrak{M}(z),\mathfrak{v}(z),\hat{\Sigma})$ (resp., $(\mathfrak{M}(z),
\mathfrak{v}(z),\tilde{\Sigma}))$ stated in item~\pmb{(1)} (resp., item~\pmb{(2)}) 
of Lemma~\ref{lem4.2}. For the reader's convenience, and without loss of 
generality, a terse sketch of a proof, for the case $n \! \in \! \mathbb{N}$ and 
$k \! \in \! \lbrace 1,2,\dotsc,K \rbrace$ such that $\alpha_{p_{\mathfrak{s}}} 
\! := \! \alpha_{k} \! \neq \! \infty$, is presented for the right-most end-points 
(see Lemma~\ref{lem4.9} below) of the intervals, $\tilde{a}_{j}$, $j \! = \! 1,2,
\dotsc,N \! + \! 1$, of the support, $J_{f}$, of the equilibrium measure, 
$\mu_{\widetilde{V}}^{f}$, whilst the left-most end-points (see 
Lemma~\ref{lem4.8} below) of the intervals, $\tilde{b}_{j-1}$, $j \! = \! 1,2,
\dotsc,N \! + \! 1$, of the support, $J_{f}$, of the equilibrium measure, 
$\mu_{\widetilde{V}}^{f}$, are analysed analogously: the proof for the case 
$n \! \in \! \mathbb{N}$ and $k \! \in \! \lbrace 1,2,\dotsc,K \rbrace$ 
such that $\alpha_{p_{\mathfrak{s}}} \! := \! \alpha_{k} \! = \! \infty$ is, 
\emph{mutatis mutandis}, similar (see Lemmata~\ref{lem4.6} 
and~\ref{lem4.7} below).

For $n \! \in \! \mathbb{N}$ and $k \! \in \! \lbrace 1,2,\dotsc,K 
\rbrace$ such that $\alpha_{p_{\mathfrak{s}}} \! := \! \alpha_{k} \! 
= \! \infty$ (resp., $\alpha_{p_{\mathfrak{s}}} \! := \! \alpha_{k} \! 
\neq \! \infty)$, by a parametrix of the RHP $(\mathfrak{M}(z),
\mathfrak{v}(z),\hat{\Sigma})$ (resp., $(\mathfrak{M}(z),\mathfrak{v}
(z),\tilde{\Sigma}))$ stated in item~\pmb{(1)} (resp., item~\pmb{(2)}) 
of Lemma~\ref{lem4.2}, in the neighbourhoods of the end-points 
of the intervals, $\hat{b}_{j-1},\hat{a}_{j}$ (resp., $\tilde{b}_{j-1},
\tilde{a}_{j})$, $j \! = \! 1,2,\dotsc,N \! + \! 1$, of the support, 
$J_{\infty}$ (resp., $J_{f})$, of the equilibrium measure, $\mu_{
\widetilde{V}}^{\infty}$ (resp., $\mu_{\widetilde{V}}^{f})$, is 
meant the solutions of the RHPs stated in Lemmata~\ref{lem4.6} 
and~\ref{lem4.7} (resp., Lemmata~\ref{lem4.8} and~\ref{lem4.9}) 
below. For $n \! \in \! \mathbb{N}$ and $k \! \in \! \lbrace 1,2,\dotsc,
K \rbrace$ such that $\alpha_{p_{\mathfrak{s}}} \! := \! \alpha_{k} \! 
= \! \infty$ (resp., $\alpha_{p_{\mathfrak{s}}} \! := \! \alpha_{k} \! 
\neq \! \infty)$, define the mutually disjoint open discs about 
$\hat{b}_{j-1},\hat{a}_{j}$ (resp., $\tilde{b}_{j-1},\tilde{a}_{j})$, 
$j \! = \! 1,2,\dotsc,N \! + \! 1$, as follows: $\hat{\mathbb{U}}_{
\hat{\delta}_{\hat{b}_{j-1}}} \! := \! \lbrace \mathstrut z \! \in \! 
\mathbb{C}; \, \vert z \! - \! \hat{b}_{j-1} \vert \! < \! \hat{
\delta}_{\hat{b}_{j-1}} \rbrace$ and $\hat{\mathbb{U}}_{\hat{\delta}_{
\hat{a}_{j}}} \! := \! \lbrace \mathstrut z \! \in \! \mathbb{C}; \, \vert 
z \! - \! \hat{a}_{j} \vert \! < \! \hat{\delta}_{\hat{a}_{j}} \rbrace$ 
(resp., $\tilde{\mathbb{U}}_{\tilde{\delta}_{\tilde{b}_{j-1}}} \! := 
\! \lbrace \mathstrut z \! \in \! \mathbb{C}; \, \vert z \! - \! 
\tilde{b}_{j-1} \vert \! < \! \tilde{\delta}_{\tilde{b}_{j-1}} \rbrace$ 
and $\tilde{\mathbb{U}}_{\tilde{\delta}_{\tilde{a}_{j}}} \! := \! \lbrace 
\mathstrut z \! \in \! \mathbb{C}; \, \vert z \! - \! \tilde{a}_{j} \vert 
\! < \! \tilde{\delta}_{\tilde{a}_{j}} \rbrace)$, $j \! = \! 1,2,\dotsc,N 
\! + \! 1$, where, amongst other things (see Lemmata~\ref{lem4.6} 
and~\ref{lem4.7} (resp., Lemmata~\ref{lem4.8} and~\ref{lem4.9}) 
below), $\hat{\delta}_{\hat{b}_{j-1}},\hat{\delta}_{\hat{a}_{j}} \! 
\in \! (0,1)$ (resp., $\tilde{\delta}_{\tilde{b}_{j-1}},\tilde{\delta}_{
\tilde{a}_{j}} \! \in \! (0,1))$ are sufficiently small, positive real numbers 
chosen so that $\hat{\mathbb{U}}_{\hat{\delta}_{\hat{b}_{i-1}}} \cap 
\hat{\mathbb{U}}_{\hat{\delta}_{\hat{a}_{j}}} \! = \! \varnothing$ 
(resp., $\tilde{\mathbb{U}}_{\tilde{\delta}_{\tilde{b}_{i-1}}} \cap 
\tilde{\mathbb{U}}_{\tilde{\delta}_{\tilde{a}_{j}}} \! = \! \varnothing)$ 
$\forall$ $i,j \! \in \! \lbrace 1,2,\dotsc,N \! + \! 1 \rbrace$, and 
$\hat{\mathbb{U}}_{\hat{\delta}_{\hat{b}_{m-1}}} \cap \lbrace 
\alpha_{p_{1}},\alpha_{p_{2}},\dotsc,\alpha_{p_{\mathfrak{s}}} \rbrace 
\! = \! \varnothing \! = \! \hat{\mathbb{U}}_{\hat{\delta}_{\hat{a}_{m}}} 
\cap \lbrace \alpha_{p_{1}},\alpha_{p_{2}},\dotsc,\alpha_{p_{\mathfrak{s}}} 
\rbrace$ (resp., $\tilde{\mathbb{U}}_{\tilde{\delta}_{\tilde{b}_{m-1}}} 
\cap \lbrace \alpha_{p_{1}},\alpha_{p_{2}},\dotsc,\alpha_{p_{\mathfrak{s}}} 
\rbrace \! = \! \varnothing \! = \! \tilde{\mathbb{U}}_{\tilde{\delta}_{
\tilde{a}_{m}}} \cap \lbrace \alpha_{p_{1}},\alpha_{p_{2}},\dotsc,\alpha_{
p_{\mathfrak{s}}} \rbrace)$, $m \! = \! 1,2,\dotsc,N \! + \! 1$: the 
corresponding regions $\hat{\Omega}_{\hat{b}_{j-1}}^{m}$ and 
$\hat{\Omega}_{\hat{a}_{j}}^{m}$ (resp., $\tilde{\Omega}_{\tilde{b}_{j
-1}}^{m}$ and $\tilde{\Omega}_{\tilde{a}_{j}}^{m})$ and arcs $\hat{
\Sigma}_{\hat{b}_{j-1}}^{m}$ and $\hat{\Sigma}_{\hat{a}_{j}}^{m}$ 
(resp., $\tilde{\Sigma}_{\tilde{b}_{j-1}}^{m}$ and $\tilde{\Sigma}_{
\tilde{a}_{j}}^{m})$, $j \! = \! 1,2,\dotsc,N \! + \! 1$, $m \! = \! 
1,2,3,4$, respectively, are defined in Lemmata~\ref{lem4.6} 
and~\ref{lem4.7} (resp., Lemmata~\ref{lem4.8} and~\ref{lem4.9}) 
below; see, also, Figures~\ref{forbhat} and~\ref{forahat} (resp., 
Figures~\ref{forbtil} and~\ref{foratil}), 
respectively.
\begin{eeeee} \label{rem4.4} 
\textsl{In order to simplify the presentation of the results of 
Lemmata~\ref{lem4.6}--\ref{lem4.9} below, it is convenient to introduce 
the following notation: {\rm (i)}
\begin{gather*}
\Psi_{1}(z) \! := \! 
\begin{pmatrix}
\operatorname{Ai}(z) & \operatorname{Ai}(\omega^{2}z) \\
\operatorname{Ai}^{\prime}(z) & \omega^{2} \operatorname{Ai}^{\prime}
(\omega^{2}z) 
\end{pmatrix} \me^{-\frac{\mi \pi}{6} \sigma_{3}}, \qquad \qquad 
\Psi_{2}(z) \! := \! 
\begin{pmatrix}
\operatorname{Ai}(z) & \operatorname{Ai}(\omega^{2}z) \\
\operatorname{Ai}^{\prime}(z) & \omega^{2} \operatorname{Ai}^{\prime}
(\omega^{2}z)
\end{pmatrix} \me^{-\frac{\mi \pi}{6} \sigma_{3}}(\mathrm{I} \! - \! 
\sigma_{-}), \\
\Psi_{3}(z) \! := \! 
\begin{pmatrix}
\operatorname{Ai}(z) & -\omega^{2} \operatorname{Ai}(\omega z) \\
\operatorname{Ai}^{\prime}(z) & -\operatorname{Ai}^{\prime}(\omega z)
\end{pmatrix} \me^{-\frac{\mi \pi}{6} \sigma_{3}}(\mathrm{I} \! + \! 
\sigma_{-}), \qquad \qquad 
\Psi_{4}(z) \! := \! 
\begin{pmatrix}
\operatorname{Ai}(z) & -\omega^{2} \operatorname{Ai}(\omega z) \\
\operatorname{Ai}^{\prime}(z) & -\operatorname{Ai}^{\prime}(\omega z)
\end{pmatrix} \me^{-\frac{\mi \pi}{6} \sigma_{3}},
\end{gather*}
where $\operatorname{Ai}(\pmb{\cdot})$ is the Airy function (cf. 
Subsection~\ref{subsub2}$)$, and $\omega \! = \! \exp (2 \pi \mi/3)$$;$ 
and {\rm (ii)} for $n \! \in \! \mathbb{N}$ and $k \! \in \! \lbrace 1,2,
\dotsc,K \rbrace$ such that $\alpha_{p_{\mathfrak{s}}} \! := \! \alpha_{k} 
\! = \! \infty$,
\begin{equation*}
\hat{\mho}_{j} \! = \!
\begin{cases}
\hat{\Omega}_{j} \! = \! 2 \pi \int_{\hat{b}_{j}}^{\hat{a}_{N+1}} \psi_{
\widetilde{V}}^{\infty}(\xi) \, \md \xi, &\text{$j \! = \! 1,2,\dotsc,N$,} \\
0, &\text{$j \! = \! 0,N \! + \! 1$,}
\end{cases}
\end{equation*}
and, for $n \! \in \! \mathbb{N}$ and $k \! \in \! \lbrace 1,2,\dotsc,K \rbrace$ 
such that $\alpha_{p_{\mathfrak{s}}} \! := \! \alpha_{k} \! \neq \! \infty$,
\begin{equation*}
\tilde{\mho}_{j} \! = \!
\begin{cases}
\tilde{\Omega}_{j} \! = \! 2 \pi \int_{\tilde{b}_{j}}^{\tilde{a}_{N+1}} \psi_{
\widetilde{V}}^{f}(\xi) \, \md \xi, &\text{$j \! = \! 1,2,\dotsc,N$,} \\
0, &\text{$j \! = \! 0,N \! + \! 1$.}
\end{cases}
\end{equation*}}
\end{eeeee}
\begin{ccccc} \label{lem4.6} 
For $n \! \in \! \mathbb{N}$ and $k \! \in \! \lbrace 1,2,\dotsc,K \rbrace$ 
such that $\alpha_{p_{\mathfrak{s}}} \! := \! \alpha_{k} \! = \! \infty$, 
let $\mathfrak{M} \colon \mathbb{C} \setminus \hat{\Sigma} \! 
\to \! \operatorname{SL}_{2}(\mathbb{C})$ solve the {\rm RHP} 
$(\mathfrak{M}(z),\mathfrak{v}(z),\hat{\Sigma})$ stated in 
item~{\rm \pmb{(1)}} of Lemma~\ref{lem4.2}. For $n \! \in \! 
\mathbb{N}$ and $k \! \in \! \lbrace 1,2,\dotsc,K \rbrace$ such that 
$\alpha_{p_{\mathfrak{s}}} \! := \! \alpha_{k} \! = \! \infty$, set 
$\hat{\mathbb{U}}_{\hat{\delta}_{\hat{b}_{i-1}}} \! := \! \lbrace \mathstrut 
z \! \in \! \mathbb{C}; \, \lvert z \! - \! \hat{b}_{i-1} \rvert \! < \! \hat{
\delta}_{\hat{b}_{i-1}} \rbrace$, $i \! = \! 1,2,\dotsc,N \! + \! 1$, and 
let $\hat{\Phi}_{\hat{b}_{j-1}}(z)$ and $\hat{\xi}_{\hat{b}_{j-1}}(z)$, $j \! = 
\! 1,2,\dotsc,N \! + \! 1$, be defined by Equations~\eqref{eqmaininf77} 
and~\eqref{eqmaininf70}, respectively, where, for $z \! \in \! \hat{
\mathbb{U}}_{\hat{\delta}_{\hat{b}_{j-1}}} \setminus (-\infty,\hat{b}_{j-1})$, 
$\hat{\xi}_{\hat{b}_{j-1}}(z) \! = \! \hat{\mathfrak{b}}(z \! - \! \hat{b}_{j-1})^{3/2} 
\hat{G}_{\hat{b}_{j-1}}(z)$, with $\hat{\mathfrak{b}} \! = \! \pm 1$ for 
$z \! \in \! \mathbb{C}_{\pm}$, and $\hat{G}_{\hat{b}_{j-1}}(z)$ is analytic, 
in particular,
\begin{equation*}
\hat{G}_{\hat{b}_{j-1}}(z) \underset{z \to \hat{b}_{j-1}}{=} \dfrac{2}{3}
f(\hat{b}_{j-1}) \! + \! \dfrac{2}{5}f^{\prime}(\hat{b}_{j-1})(z \! - \! 
\hat{b}_{j-1}) \! + \! \dfrac{1}{7}f^{\prime \prime}(\hat{b}_{j-1})(z \! - \! 
\hat{b}_{j-1})^{2} \! + \! \mathcal{O}((z \! - \! \hat{b}_{j-1})^{3}),
\end{equation*}
where {}\footnote{Note: $\mi \! := \! \exp (\mi \pi/2)$.}
\begin{align}
f(\hat{b}_{0}) =& \, \mi (-1)^{N} \hat{h}_{\widetilde{V}}(\hat{b}_{0}) 
\hat{\eta}_{\hat{b}_{0}}, \nonumber \\
f^{\prime}(\hat{b}_{0}) =& \, \mi (-1)^{N} \left(\dfrac{1}{2} \hat{h}_{
\widetilde{V}}(\hat{b}_{0}) \hat{\eta}_{\hat{b}_{0}} \left(\sum_{m=1}^{N} 
\left(\dfrac{1}{\hat{b}_{0} \! - \! \hat{b}_{m}} \! + \! \dfrac{1}{\hat{b}_{0} 
\! - \! \hat{a}_{m}} \right) \! + \! \dfrac{1}{\hat{b}_{0} \! - \! \hat{a}_{N+1}} 
\right) \! + \! (\hat{h}_{\widetilde{V}}(\hat{b}_{0}))^{\prime} \hat{\eta}_{
\hat{b}_{0}} \right), \nonumber \\
f^{\prime \prime}(\hat{b}_{0}) =& \, \mi (-1)^{N} \left(\dfrac{\hat{h}_{
\widetilde{V}}(\hat{b}_{0})(\hat{h}_{\widetilde{V}}(\hat{b}_{0}))^{\prime 
\prime} \! - \! ((\hat{h}_{\widetilde{V}}(\hat{b}_{0}))^{\prime})^{2}}{
\hat{h}_{\widetilde{V}}(\hat{b}_{0})} \hat{\eta}_{\hat{b}_{0}} \! - \! 
\dfrac{1}{2} \hat{h}_{\widetilde{V}}(\hat{b}_{0}) \hat{\eta}_{\hat{b}_{0}} 
\right. \nonumber \\
\times&\left. \, \left(\sum_{m=1}^{N} \left(\dfrac{1}{(\hat{b}_{0} \! - \! 
\hat{b}_{m})^{2}} \! + \! \dfrac{1}{(\hat{b}_{0} \! - \! \hat{a}_{m})^{2}} 
\right) \! + \! \dfrac{1}{(\hat{b}_{0} \! - \! \hat{a}_{N+1})^{2}} \right) 
\right. \nonumber \\
+&\left. \, \left(\dfrac{1}{2} \left(\sum_{m=1}^{N} \left(\dfrac{1}{
\hat{b}_{0} \! - \! \hat{b}_{m}} \! + \! \dfrac{1}{\hat{b}_{0} \! - \! \hat{a}_{m}} 
\right) \! + \! \dfrac{1}{\hat{b}_{0} \! - \! \hat{a}_{N+1}} \right) \! + \! 
\dfrac{(\hat{h}_{\widetilde{V}}(\hat{b}_{0}))^{\prime}}{\hat{h}_{\widetilde{V}}
(\hat{b}_{0})} \right) \right. \nonumber \\
\times&\left. \, \left(\dfrac{1}{2} \hat{h}_{\widetilde{V}}(\hat{b}_{0}) 
\hat{\eta}_{\hat{b}_{0}} \left(\sum_{m=1}^{N} \left(\dfrac{1}{\hat{b}_{0} 
\! - \! \hat{b}_{m}} \! + \! \dfrac{1}{\hat{b}_{0} \! - \! \hat{a}_{m}} \right) 
\! + \! \dfrac{1}{\hat{b}_{0} \! - \! \hat{a}_{N+1}} \right) \! + \! (\hat{h}_{
\widetilde{V}}(\hat{b}_{0}))^{\prime} \hat{\eta}_{\hat{b}_{0}} \right) \right), 
\label{wm24frbohat}
\end{align}
with $\hat{\eta}_{\hat{b}_{0}}$ defined by Equation~\eqref{eqmaininf55}, 
and, for $j \! = \! 1,2,\dotsc,N$,
\begin{align}
f(\hat{b}_{j}) =& \, \mi (-1)^{N-j} \hat{h}_{\widetilde{V}}(\hat{b}_{j}) 
\hat{\eta}_{\hat{b}_{j}}, \nonumber \\
f^{\prime}(\hat{b}_{j}) =& \, \mi (-1)^{N-j} \left(\dfrac{1}{2} \hat{h}_{
\widetilde{V}}(\hat{b}_{j}) \hat{\eta}_{\hat{b}_{j}} \left(\sum_{\substack{m
=1\\m \not= j}}^{N} \left(\dfrac{1}{\hat{b}_{j} \! - \! \hat{b}_{m}} \! + 
\! \dfrac{1}{\hat{b}_{j} \! - \! \hat{a}_{m}} \right) \! + \! \dfrac{1}{\hat{b}_{j} 
\! - \! \hat{a}_{j}} \! + \! \dfrac{1}{\hat{b}_{j} \! - \! \hat{a}_{N+1}} \! + \! 
\dfrac{1}{\hat{b}_{j} \! - \! \hat{b}_{0}} \right) \! + \! (\hat{h}_{\widetilde{V}}
(\hat{b}_{j}))^{\prime} \hat{\eta}_{\hat{b}_{j}} \right), \nonumber \\
f^{\prime \prime}(\hat{b}_{j}) =& \, \mi (-1)^{N-j} \left(\dfrac{\hat{h}_{
\widetilde{V}}(\hat{b}_{j})(\hat{h}_{\widetilde{V}}(\hat{b}_{j}))^{\prime 
\prime} \! - \! ((\hat{h}_{\widetilde{V}}(\hat{b}_{j}))^{\prime})^{2}}{\hat{h}_{
\widetilde{V}}(\hat{b}_{j})} \hat{\eta}_{\hat{b}_{j}} \! - \! \dfrac{1}{2} 
\hat{h}_{\widetilde{V}}(\hat{b}_{j}) \hat{\eta}_{\hat{b}_{j}} \left(
\sum_{\substack{m=1\\m \not= j}}^{N} \left(\dfrac{1}{(\hat{b}_{j} \! - 
\! \hat{b}_{m})^{2}} \! + \! \dfrac{1}{(\hat{b}_{j} \! - \! \hat{a}_{m})^{2}} 
\right) \right. \right. \nonumber \\
+&\left. \left. \, \dfrac{1}{(\hat{b}_{j} \! - \! \hat{a}_{j})^{2}} \! + \! 
\dfrac{1}{(\hat{b}_{j} \! - \! \hat{a}_{N+1})^{2}} \! + \! \dfrac{1}{(\hat{b}_{j} 
\! - \! \hat{b}_{0})^{2}} \right) \! + \! \left(\dfrac{(\hat{h}_{\widetilde{V}}
(\hat{b}_{j}))^{\prime}}{\hat{h}_{\widetilde{V}}(\hat{b}_{j})} \! + \! 
\dfrac{1}{2} \left(\sum_{\substack{m=1\\m \not= j}}^{N} \left(
\dfrac{1}{\hat{b}_{j} \! - \! \hat{b}_{m}} \! + \! \dfrac{1}{\hat{b}_{j} \! 
- \! \hat{a}_{m}} \right) \right. \right. \right. \nonumber \\
+&\left. \left. \left. \, \dfrac{1}{\hat{b}_{j} \! - \! \hat{a}_{j}} \! + \! 
\dfrac{1}{\hat{b}_{j} \! - \! \hat{a}_{N+1}} \! + \! \dfrac{1}{\hat{b}_{j} 
\! - \! \hat{b}_{0}} \right) \right) \left(\dfrac{1}{2} \hat{h}_{\widetilde{V}}
(\hat{b}_{j}) \hat{\eta}_{\hat{b}_{j}} \left(\sum_{\substack{m=1\\m 
\not= j}}^{N} \left(\dfrac{1}{\hat{b}_{j} \! - \! \hat{b}_{m}} \! + \! 
\dfrac{1}{\hat{b}_{j} \! - \! \hat{a}_{m}} \right) \right. \right. \right. 
\nonumber \\
+&\left. \left. \left. \, \dfrac{1}{\hat{b}_{j} \! - \! \hat{a}_{j}} \! + \! 
\dfrac{1}{\hat{b}_{j} \! - \! \hat{a}_{N+1}} \! + \! \dfrac{1}{\hat{b}_{j} 
\! - \! \hat{b}_{0}} \right) \! + \!  (\hat{h}_{\widetilde{V}}(\hat{b}_{j}))^{
\prime} \hat{\eta}_{\hat{b}_{j}} \right) \right), \label{wm24frbjhat}
\end{align}
with $\hat{\eta}_{\hat{b}_{j}}$ defined by Equation~\eqref{eqmaininf56}, 
and $(0,1) \! \ni \! \hat{\delta}_{\hat{b}_{j-1}}$, $j \! = \! 1,2,\dotsc,N \! 
+ \! 1$, are chosen sufficiently small so that $\hat{\Phi}_{\hat{b}_{j-1}}(z)$, 
which are bi-holomorphic, conformal, and non-orientation preserving, 
map (cf. Figure~\ref{forbhat}$)$ $\hat{\mathbb{U}}_{\hat{\delta}_{
\hat{b}_{j-1}}}$, and, thus, the oriented contours $\hat{\Sigma}_{
\hat{b}_{j-1}} \! := \! \cup_{m=1}^{4} \hat{\Sigma}_{\hat{b}_{j-1}}^{m}$, 
$j \! = \! 1,2,\dotsc,N \! + \! 1$, injectively onto open $(n$-, $k$-, and 
$z_{o}$-dependent) neighbourhoods $\hat{\mathbb{U}}_{\hat{\delta}_{
\hat{b}_{j-1}}}^{\ast}$ of the origin, $j \! = \! 1,2,\dotsc,N \! + \! 1$, 
such that $\hat{\Phi}_{\hat{b}_{j-1}}(\hat{b}_{j-1}) \! = \! 0$, 
$\hat{\Phi}_{\hat{b}_{j-1}} \colon \hat{\mathbb{U}}_{\hat{\delta}_{
\hat{b}_{j-1}}} \! \to \! \hat{\mathbb{U}}_{\hat{\delta}_{\hat{b}_{j-1}}}^{
\ast} \! := \! \hat{\Phi}_{\hat{b}_{j-1}}(\hat{\mathbb{U}}_{\hat{\delta}_{
\hat{b}_{j-1}}})$, $\hat{\Phi}_{\hat{b}_{j-1}}(\hat{\mathbb{U}}_{\hat{\delta}_{
\hat{b}_{j-1}}} \cap \hat{\Sigma}_{\hat{b}_{j-1}}^{m}) \! = \! \hat{\Phi}_{
\hat{b}_{j-1}}(\hat{\mathbb{U}}_{\hat{\delta}_{\hat{b}_{j-1}}}) \cap 
\gamma_{\hat{b}_{j-1}}^{\ast,m}$, and $\hat{\Phi}_{\hat{b}_{j-1}}
(\hat{\mathbb{U}}_{\hat{\delta}_{\hat{b}_{j-1}}} \cap \hat{\Omega}_{
\hat{b}_{j-1}}^{m}) \! = \! \hat{\Phi}_{\hat{b}_{j-1}}(\hat{\mathbb{U}}_{
\hat{\delta}_{\hat{b}_{j-1}}}) \cap \hat{\Omega}_{\hat{b}_{j-1}}^{\ast,m}$, 
$m \! = \! 1,2,3,4$, with $\hat{\Omega}_{\hat{b}_{j-1}}^{\ast,1} \! = \! 
\lbrace \mathstrut \zeta \! \in \! \mathbb{C}; \, \arg (\zeta) \! \in \! (0,
2 \pi/3) \rbrace$, $\hat{\Omega}_{\hat{b}_{j-1}}^{\ast,2} \! = \! \lbrace 
\mathstrut \zeta \! \in \! \mathbb{C}; \, \arg (\zeta) \! \in \! (2 \pi/3,\pi) 
\rbrace$, $\hat{\Omega}_{\hat{b}_{j-1}}^{\ast,3} \! = \! \lbrace \mathstrut 
\zeta \! \in \! \mathbb{C}; \, \arg (\zeta) \! \in \! (-\pi,-2 \pi/3) \rbrace$, 
and $\hat{\Omega}_{\hat{b}_{j-1}}^{\ast,4} \! = \! \lbrace \mathstrut 
\zeta \! \in \! \mathbb{C}; \, \arg (\zeta) \! \in \! (-2 \pi/3,0) \rbrace$.

For $n \! \in \! \mathbb{N}$ and $k \! \in \! \lbrace 1,2,\dotsc,K \rbrace$ 
such that $\alpha_{p_{\mathfrak{s}}} \! := \! \alpha_{k} \! = \! \infty$, 
the parametrix for the {\rm RHP} $(\mathfrak{M}(z),\mathfrak{v}(z),
\hat{\Sigma})$, for $z \! \in \! \hat{\mathbb{U}}_{\hat{\delta}_{\hat{b}_{j-1}}}$, 
$j \! = \! 1,2,\dotsc,N \! + \! 1$, is the solution of the following {\rm RHPs} 
for $\hat{\mathcal{X}}^{\hat{b}} \colon \hat{\mathbb{U}}_{\hat{\delta}_{
\hat{b}_{j-1}}} \setminus \hat{\Sigma}_{\hat{b}_{j-1}} \! \to \! \operatorname{SL}_{2}
(\mathbb{C})$, $j \! = \! 1,2,\dotsc,N \! + \! 1$, where $\hat{\Sigma}_{
\hat{b}_{j-1}} \! := \! (\hat{\Phi}_{\hat{b}_{j-1}})^{-1}(\gamma_{\hat{b}_{j-
1}}^{\ast})$, with $(\hat{\Phi}_{\hat{b}_{j-1}})^{-1}$ denoting the inverse 
mapping, and $\gamma_{\hat{b}_{j-1}}^{\ast} \! := \! \cup_{m=1}^{4} 
\gamma_{\hat{b}_{j-1}}^{\ast,m}$$:$ {\rm (i)} $\hat{\mathcal{X}}^{\hat{b}}
(z)$ is analytic for $z \! \in \! \hat{\mathbb{U}}_{\hat{\delta}_{\hat{b}_{j-1}}} 
\setminus \hat{\Sigma}_{\hat{b}_{j-1}}$, $j \! = \! 1,2,\dotsc,N \! + \! 1$$;$ 
{\rm (ii)} the boundary values $\hat{\mathcal{X}}^{\hat{b}}_{\pm}(z) \! := \! 
\lim_{\underset{z^{\prime} \, \in \, \pm \, \mathrm{side} \, \mathrm{of} \, 
\hat{\Sigma}_{\hat{b}_{j-1}}}{z^{\prime} \to z \in \hat{\Sigma}_{\hat{b}_{j-1}}}} 
\hat{\mathcal{X}}^{\hat{b}}(z^{\prime})$ satisfy the jump condition
\begin{equation*}
\hat{\mathcal{X}}^{\hat{b}}_{+}(z) \! = \! \hat{\mathcal{X}}^{\hat{b}}_{-}(z) 
\mathfrak{v}(z), \quad z \! \in \! \hat{\mathbb{U}}_{\hat{\delta}_{\hat{b}_{j-1}}} 
\cap \hat{\Sigma}_{\hat{b}_{j-1}}, \quad j \! = \! 1,2,\dotsc,N \! + \! 1;
\end{equation*}
and {\rm (iii)} uniformly for $z \! \in \! \partial \hat{\mathbb{U}}_{\hat{
\delta}_{\hat{b}_{j-1}}} \! := \! \lbrace \mathstrut z \! \in \! \mathbb{C}; 
\, \vert z \! - \! \hat{b}_{j-1} \vert \! = \! \hat{\delta}_{\hat{b}_{j-1}} 
\rbrace$, $j \! = \! 1,2,\dotsc,N \! + \! 1$,
\begin{equation*}
\mathfrak{m}(z)(\hat{\mathcal{X}}^{\hat{b}}(z))^{-1} 
\underset{\underset{z_{o}=1+o(1)}{\mathscr{N},n \to \infty}}{=} 
\mathrm{I} \! + \! o(1).
\end{equation*}
For $n \! \in \! \mathbb{N}$ and $k \! \in \! \lbrace 1,2,\dotsc,K \rbrace$ 
such that $\alpha_{p_{\mathfrak{s}}} \! := \! \alpha_{k} \! = \! \infty$, 
the solutions of the {\rm RHPs} $(\hat{\mathcal{X}}^{\hat{b}}(z),\mathfrak{v}
(z),\hat{\mathbb{U}}_{\hat{\delta}_{\hat{b}_{j-1}}} \cap \hat{\Sigma}_{
\hat{b}_{j-1}})$, $j \! = \! 1,2,\dotsc,N \! + \! 1$, are: {\rm \pmb{(1)}} for 
$z \! \in \! \hat{\Omega}_{\hat{b}_{j-1}}^{1} \! := \! \hat{\mathbb{U}}_{
\hat{\delta}_{\hat{b}_{j-1}}} \cap (\hat{\Phi}_{\hat{b}_{j-1}})^{-1}
(\hat{\Omega}_{\hat{b}_{j-1}}^{\ast,1})$, $j \! = \! 1,2,\dotsc,N \! + \! 1$,
\begin{equation*}
\hat{\mathcal{X}}^{\hat{b}}(z) \! = \! \sqrt{\smash[b]{\pi}} \, \me^{-
\frac{\mi \pi}{3}} \mathfrak{m}(z) \sigma_{3} \me^{\frac{\mi}{2}((n-1)K+k) 
\hat{\mho}_{j-1} \operatorname{ad}(\sigma_{3})} \begin{pmatrix}
\mi & -\mi \\
1 & 1
\end{pmatrix}(\hat{\Phi}_{\hat{b}_{j-1}}(z))^{\frac{1}{4} \sigma_{3}} 
\Psi_{1}(\hat{\Phi}_{\hat{b}_{j-1}}(z)) \me^{\frac{1}{2}((n-1)K+k) 
\hat{\xi}_{\hat{b}_{j-1}}(z) \sigma_{3}} \sigma_{3},
\end{equation*}
where $\mathfrak{m}(z)$ is given in item~{\rm \pmb{(1)}} of 
Lemma~\ref{lem4.5}$;$ {\rm \pmb{(2)}} for $z \! \in \! \hat{\Omega}_{
\hat{b}_{j-1}}^{2} \! := \! \hat{\mathbb{U}}_{\hat{\delta}_{\hat{b}_{j-1}}} 
\cap (\hat{\Phi}_{\hat{b}_{j-1}})^{-1}(\hat{\Omega}_{\hat{b}_{j-1}}^{\ast,2})$, 
$j \! = \! 1,2,\dotsc,N \! + \! 1$,
\begin{equation*}
\hat{\mathcal{X}}^{\hat{b}}(z) \! = \! \sqrt{\smash[b]{\pi}} \, \me^{-
\frac{\mi \pi}{3}} \mathfrak{m}(z) \sigma_{3} \me^{\frac{\mi}{2}((n-1)K+k) 
\hat{\mho}_{j-1} \operatorname{ad}(\sigma_{3})} 
\begin{pmatrix}
\mi & -\mi \\
1 & 1
\end{pmatrix}(\hat{\Phi}_{\hat{b}_{j-1}}(z))^{\frac{1}{4} \sigma_{3}} \Psi_{2}
(\hat{\Phi}_{\hat{b}_{j-1}}(z)) \me^{\frac{1}{2}((n-1)K+k) 
\hat{\xi}_{\hat{b}_{j-1}}(z) \sigma_{3}} \sigma_{3};
\end{equation*}
{\rm \pmb{(3)}} for $z \! \in \! \hat{\Omega}_{\hat{b}_{j-1}}^{3} \! := \! 
\hat{\mathbb{U}}_{\hat{\delta}_{\hat{b}_{j-1}}} \cap (\hat{\Phi}_{\hat{b}_{j
-1}})^{-1}(\hat{\Omega}_{\hat{b}_{j-1}}^{\ast,3})$, $j \! = \! 1,2,\dotsc,
N \! + \! 1$,
\begin{equation*}
\hat{\mathcal{X}}^{\hat{b}}(z) \! = \! \sqrt{\smash[b]{\pi}} \, \me^{-
\frac{\mi \pi}{3}} \mathfrak{m}(z) \sigma_{3} \me^{-\frac{\mi}{2}((n-1)K+k) 
\hat{\mho}_{j-1} \operatorname{ad}(\sigma_{3})} 
\begin{pmatrix}
\mi & -\mi \\
1 & 1
\end{pmatrix}(\hat{\Phi}_{\hat{b}_{j-1}}(z))^{\frac{1}{4} \sigma_{3}} \Psi_{3}
(\hat{\Phi}_{\hat{b}_{j-1}}(z)) \me^{\frac{1}{2}((n-1)K+k) \hat{\xi}_{\hat{b}_{j
-1}}(z) \sigma_{3}} \sigma_{3};
\end{equation*}
and {\rm \pmb{(4)}} for $z \! \in \! \hat{\Omega}_{\hat{b}_{j-1}}^{4} \! := \! \hat{
\mathbb{U}}_{\hat{\delta}_{\hat{b}_{j-1}}} \cap (\hat{\Phi}_{\hat{b}_{j-1}})^{-1}
(\hat{\Omega}_{\hat{b}_{j-1}}^{\ast,4})$, $j \! = \! 1,2,\dotsc,N \! + \! 1$,
\begin{equation*}
\hat{\mathcal{X}}^{\hat{b}}(z) \! = \! \sqrt{\smash[b]{\pi}} \, \me^{-
\frac{\mi \pi}{3}} \mathfrak{m}(z) \sigma_{3} \me^{-\frac{\mi}{2}((n-1)K+k) 
\hat{\mho}_{j-1} \operatorname{ad}(\sigma_{3})} 
\begin{pmatrix}
\mi & -\mi \\
1 & 1
\end{pmatrix}(\hat{\Phi}_{\hat{b}_{j-1}}(z))^{\frac{1}{4} \sigma_{3}} \Psi_{4}
(\hat{\Phi}_{\hat{b}_{j-1}}(z)) \me^{\frac{1}{2}((n-1)K+k) \hat{\xi}_{\hat{
b}_{j-1}}(z) \sigma_{3}} \sigma_{3}.
\end{equation*}
\end{ccccc}
\begin{ccccc} \label{lem4.7} 
For $n \! \in \! \mathbb{N}$ and $k \! \in \! \lbrace 1,2,\dotsc,K \rbrace$ 
such that $\alpha_{p_{\mathfrak{s}}} \! := \! \alpha_{k} \! = \! \infty$, 
let $\mathfrak{M} \colon \mathbb{C} \setminus \hat{\Sigma} \! 
\to \! \operatorname{SL}_{2}(\mathbb{C})$ solve the {\rm RHP} 
$(\mathfrak{M}(z),\mathfrak{v}(z),\hat{\Sigma})$ stated in 
item~{\rm \pmb{(1)}} of Lemma~\ref{lem4.2}. For $n \! \in \! 
\mathbb{N}$ and $k \! \in \! \lbrace 1,2,\dotsc,K \rbrace$ such that 
$\alpha_{p_{\mathfrak{s}}} \! := \! \alpha_{k} \! = \! \infty$, set $\hat{
\mathbb{U}}_{\hat{\delta}_{\hat{a}_{i}}} \! := \! \lbrace \mathstrut z \! \in 
\! \mathbb{C}; \, \lvert z \! - \! \hat{a}_{i} \rvert \! < \! \hat{\delta}_{
\hat{a}_{i}} \rbrace$, $i \! = \! 1,2,\dotsc,N \! + \! 1$, and let 
$\hat{\Phi}_{\hat{a}_{j}}(z)$ and $\hat{\xi}_{\hat{a}_{j}}(z)$, $j \! = \! 
1,2,\dotsc,N \! + \! 1$, be defined by Equations~\eqref{eqmaininf84} 
and~\eqref{eqmaininf71}, respectively, where, for $z \! \in \! \hat{
\mathbb{U}}_{\hat{\delta}_{\hat{a}_{j}}} \setminus (-\infty,\hat{a}_{j})$, 
$\hat{\xi}_{\hat{a}_{j}}(z) \! = \! (z \! - \! \hat{a}_{j})^{3/2} \hat{G}_{
\hat{a}_{j}}(z)$, and $\hat{G}_{\hat{a}_{j}}(z)$ is analytic, in particular,
\begin{equation*}
\hat{G}_{\hat{a}_{j}}(z) \underset{z \to \hat{a}_{j}}{=} \dfrac{2}{3}f(\hat{a}_{j}) 
\! + \! \dfrac{2}{5}f^{\prime}(\hat{a}_{j})(z \! - \! \hat{a}_{j}) \! + \! \dfrac{
1}{7}f^{\prime \prime}(\hat{a}_{j})(z \! - \! \hat{a}_{j})^{2} \! + \! \mathcal{O}
((z \! - \! \hat{a}_{j})^{3}),
\end{equation*}
where
\begin{align}
f(\hat{a}_{N+1}) =& \, \hat{h}_{\widetilde{V}}(\hat{a}_{N+1}) \hat{\eta}_{
\hat{a}_{N+1}}, \nonumber \\
f^{\prime}(\hat{a}_{N+1}) =& \, \dfrac{1}{2} \hat{h}_{\widetilde{V}}
(\hat{a}_{N+1}) \hat{\eta}_{\hat{a}_{N+1}} \left(\sum_{m=1}^{N} \left(
\dfrac{1}{\hat{a}_{N+1} \! - \! \hat{b}_{m}} \! + \! \dfrac{1}{\hat{a}_{N+1} 
\! - \! \hat{a}_{m}} \right) \! + \! \dfrac{1}{\hat{a}_{N+1} \! - \! \hat{b}_{0}} 
\right) \! + \! (\hat{h}_{\widetilde{V}}(\hat{a}_{N+1}))^{\prime} \hat{\eta}_{
\hat{a}_{N+1}}, \nonumber \\
f^{\prime \prime}(\hat{a}_{N+1}) =& \, \dfrac{\hat{h}_{\widetilde{V}}
(\hat{a}_{N+1})(\hat{h}_{\widetilde{V}}(\hat{a}_{N+1}))^{\prime \prime} 
\! - \! ((\hat{h}_{\widetilde{V}}(\hat{a}_{N+1}))^{\prime})^{2}}{\hat{h}_{
\widetilde{V}}(\hat{a}_{N+1})} \hat{\eta}_{\hat{a}_{N+1}} \! - \! \dfrac{1}{2} 
\hat{h}_{\widetilde{V}}(\hat{a}_{N+1}) \hat{\eta}_{\hat{a}_{N+1}} \nonumber \\
\times& \, \left(\sum_{m=1}^{N} \left(\dfrac{1}{(\hat{a}_{N+1} \! - \! 
\hat{b}_{m})^{2}} \! + \! \dfrac{1}{(\hat{a}_{N+1} \! - \! \hat{a}_{m})^{2}} 
\right) \! + \! \dfrac{1}{(\hat{a}_{N+1} \! - \! \hat{b}_{0})^{2}} \right) 
\nonumber \\
+& \, \left(\dfrac{1}{2} \left(\sum_{m=1}^{N} \left(\dfrac{1}{\hat{a}_{N+1} 
\! - \! \hat{b}_{m}} \! + \! \dfrac{1}{\hat{a}_{N+1} \! - \! \hat{a}_{m}} 
\right) \! + \! \dfrac{1}{\hat{a}_{N+1} \! - \! \hat{b}_{0}} \right) \! + 
\! \dfrac{(\hat{h}_{\widetilde{V}}(\hat{a}_{N+1}))^{\prime}}{\hat{h}_{
\widetilde{V}}(\hat{a}_{N+1})} \right) \nonumber \\
\times& \, \left(\dfrac{1}{2} \hat{h}_{\widetilde{V}}(\hat{a}_{N+1}) \hat{
\eta}_{\hat{a}_{N+1}} \left(\sum_{m=1}^{N} \! \left(\dfrac{1}{\hat{a}_{N+1} 
\! - \! \hat{a}_{m}} \! + \! \dfrac{1}{\hat{a}_{N+1} \! - \! \hat{b}_{m}} 
\right) \! + \! \dfrac{1}{\hat{a}_{N+1} \! - \! \hat{b}_{0}} \right) \right. 
\nonumber \\
+& \left. \, (\hat{h}_{\widetilde{V}}(\hat{a}_{N+1}))^{\prime} 
\hat{\eta}_{\hat{a}_{N+1}} \right), \label{wm24franhat}
\end{align}
with $\hat{\eta}_{\hat{a}_{N+1}}$ defined by Equation~\eqref{eqmaininf57}, 
and, for $j \! = \! 1,2,\dotsc,N$,
\begin{align}
f(\hat{a}_{j}) =& \, (-1)^{N+1-j} \hat{h}_{\widetilde{V}}(\hat{a}_{j}) \hat{
\eta}_{\hat{a}_{j}}, \nonumber \\
f^{\prime}(\hat{a}_{j}) =& \, (-1)^{N+1-j} \left(\dfrac{1}{2} \hat{h}_{
\widetilde{V}}(\hat{a}_{j}) \hat{\eta}_{\hat{a}_{j}} \left(\sum_{
\substack{m=1\\m \not= j}}^{N} \left(\dfrac{1}{\hat{a}_{j} \! - \! 
\hat{b}_{m}} \! + \! \dfrac{1}{\hat{a}_{j} \! - \! \hat{a}_{m}} \right) \! + \! 
\dfrac{1}{\hat{a}_{j} \! - \! \hat{b}_{j}} \! + \! \dfrac{1}{\hat{a}_{j} \! - \! 
\hat{a}_{N+1}} \! + \! \dfrac{1}{\hat{a}_{j} \! - \! \hat{b}_{0}} \right) \! + 
\! (\hat{h}_{\widetilde{V}}(\hat{a}_{j}))^{\prime} \hat{\eta}_{\hat{a}_{j}} 
\right), \nonumber \\
f^{\prime \prime}(\hat{a}_{j}) =& \, (-1)^{N+1-j} \left(\dfrac{\hat{h}_{
\widetilde{V}}(\hat{a}_{j})(\hat{h}_{\widetilde{V}}(\hat{a}_{j}))^{\prime 
\prime} \! - \! ((\hat{h}_{\widetilde{V}}(\hat{a}_{j}))^{\prime})^{2}}{
\hat{h}_{\widetilde{V}}(\hat{a}_{j})} \hat{\eta}_{\hat{a}_{j}} \! - \! 
\dfrac{1}{2} \hat{h}_{\widetilde{V}}(\hat{a}_{j}) \hat{\eta}_{\hat{a}_{j}} 
\left(\sum_{\substack{m=1\\m \not= j}}^{N} \left(\dfrac{1}{(\hat{a}_{j} 
\! - \! \hat{b}_{m})^{2}} \! + \! \dfrac{1}{(\hat{a}_{j} \! - \! \hat{a}_{m})^{2}} 
\right) \right. \right. \nonumber \\
+&\left. \left. \, \dfrac{1}{(\hat{a}_{j} \! - \! \hat{b}_{j})^{2}} \! + \! 
\dfrac{1}{(\hat{a}_{j} \! - \! \hat{a}_{N+1})^{2}} \! + \! \dfrac{1}{(\hat{a}_{j} 
\! - \! \hat{b}_{0})^{2}} \right) \! + \! \left(\dfrac{(\hat{h}_{\widetilde{V}}
(\hat{a}_{j}))^{\prime}}{\hat{h}_{\widetilde{V}}(\hat{a}_{j})} \! + \! 
\dfrac{1}{2} \left(\sum_{\substack{m=1\\m \not= j}}^{N} \left(
\dfrac{1}{\hat{a}_{j} \! - \! \hat{b}_{m}} \! + \! \dfrac{1}{\hat{a}_{j} \! 
- \! \hat{a}_{m}} \right) \right. \right. \right. \nonumber \\
+&\left. \left. \left. \, \dfrac{1}{\hat{a}_{j} \! - \! \hat{b}_{j}} \! + \! 
\dfrac{1}{\hat{a}_{j} \! - \! \hat{a}_{N+1}} \! + \! \dfrac{1}{\hat{a}_{j} 
\! - \! \hat{b}_{0}} \right) \right) \left(\dfrac{1}{2} \hat{h}_{\widetilde{V}}
(\hat{a}_{j}) \hat{\eta}_{\hat{a}_{j}} \left(\sum_{\substack{m=1\\m \not= 
j}}^{N} \left(\dfrac{1}{\hat{a}_{j} \! - \! \hat{b}_{m}} \! + \! \dfrac{1}{
\hat{a}_{j} \! - \! \hat{a}_{m}} \right) \right. \right. \right. \nonumber \\
+&\left. \left. \left. \, \dfrac{1}{\hat{a}_{j} \! - \! \hat{b}_{j}} \! + \! 
\dfrac{1}{\hat{a}_{j} \! - \! \hat{a}_{N+1}} \! + \! \dfrac{1}{\hat{a}_{j} 
\! - \! \hat{b}_{0}} \right) \! + \!  (\hat{h}_{\widetilde{V}}(\hat{a}_{j}))^{
\prime} \hat{\eta}_{\hat{a}_{j}} \right) \right), \label{wm24frajhat}
\end{align}
with $\hat{\eta}_{\hat{a}_{j}}$ defined by Equation~\eqref{eqmaininf58}, 
and $(0,1) \! \ni \! \hat{\delta}_{\hat{a}_{j}}$, $j \! = \! 1,2,\dotsc,N \! 
+ \! 1$, are chosen sufficiently small so that $\hat{\Phi}_{\hat{a}_{j}}(z)$, 
which are bi-holomorphic, conformal, and orientation preserving, map 
(cf. Figure~\ref{forahat}$)$ $\hat{\mathbb{U}}_{\hat{\delta}_{\hat{a}_{j}}}$, 
and, thus, the oriented contours $\hat{\Sigma}_{\hat{a}_{j}} \! := \! 
\cup_{m=1}^{4} \hat{\Sigma}_{\hat{a}_{j}}^{m}$, $j \! = \! 1,2,\dotsc,N \! + \! 1$, 
injectively onto open $(n$-, $k$-, and $z_{o}$-dependent) neighbourhoods 
$\hat{\mathbb{U}}_{\hat{\delta}_{\hat{a}_{j}}}^{\ast}$ of the origin, $j \! = 
\! 1,2,\dotsc,N \! + \! 1$, such that $\hat{\Phi}_{\hat{a}_{j}}(\hat{a}_{j}) \! 
= \! 0$, $\hat{\Phi}_{\hat{a}_{j}} \colon \hat{\mathbb{U}}_{\hat{\delta}_{
\hat{a}_{j}}} \! \to \! \hat{\mathbb{U}}_{\hat{\delta}_{\hat{a}_{j}}}^{\ast} \! 
:= \! \hat{\Phi}_{\hat{a}_{j}}(\hat{\mathbb{U}}_{\hat{\delta}_{\hat{a}_{j}}})$, 
$\hat{\Phi}_{\hat{a}_{j}}(\hat{\mathbb{U}}_{\hat{\delta}_{\hat{a}_{j}}} \cap 
\hat{\Sigma}_{\hat{a}_{j}}^{m}) \! = \! \hat{\Phi}_{\hat{a}_{j}}(\hat{
\mathbb{U}}_{\hat{\delta}_{\hat{a}_{j}}}) \cap \gamma_{\hat{a}_{j}}^{\ast,m}$, 
and $\hat{\Phi}_{\hat{a}_{j}}(\hat{\mathbb{U}}_{\hat{\delta}_{\hat{a}_{j}}} 
\cap \hat{\Omega}_{\hat{a}_{j}}^{m}) \! = \! \hat{\Phi}_{\hat{a}_{j}}
(\hat{\mathbb{U}}_{\hat{\delta}_{\hat{a}_{j}}}) \cap \hat{\Omega}_{
\hat{a}_{j}}^{\ast,m}$, $m \! = \! 1,2,3,4$, with $\hat{\Omega}_{
\hat{a}_{j}}^{\ast,1} \! = \! \lbrace \mathstrut \zeta \! \in \! \mathbb{C}; 
\, \arg (\zeta) \! \in \! (0,2 \pi/3) \rbrace$, $\hat{\Omega}_{\hat{a}_{j}}^{
\ast,2} \! = \! \lbrace \mathstrut \zeta \! \in \! \mathbb{C}; \, \arg (\zeta) 
\! \in \! (2 \pi/3,\pi) \rbrace$, $\hat{\Omega}_{\hat{a}_{j}}^{\ast,3} \! = 
\! \lbrace \mathstrut \zeta \! \in \! \mathbb{C}; \, \arg (\zeta) \! \in \! 
(-\pi,-2 \pi/3) \rbrace$, and $\hat{\Omega}_{\hat{a}_{j}}^{\ast,4} \! = \! 
\lbrace \mathstrut \zeta \! \in \! \mathbb{C}; \, \arg (\zeta) \! \in \! 
(-2 \pi/3,0) \rbrace$.

For $n \! \in \! \mathbb{N}$ and $k \! \in \! \lbrace 1,2,\dotsc,K \rbrace$ 
such that $\alpha_{p_{\mathfrak{s}}} \! := \! \alpha_{k} \! = \! \infty$, 
the parametrix for the {\rm RHP} $(\mathfrak{M}(z),\mathfrak{v}(z),
\hat{\Sigma})$, for $z \! \in \! \hat{\mathbb{U}}_{\hat{\delta}_{\hat{a}_{j}}}$, 
$j \! = \! 1,2,\dotsc,N \! + \! 1$, is the solution of the following {\rm RHPs} 
for $\hat{\mathcal{X}}^{\hat{a}} \colon \hat{\mathbb{U}}_{\hat{\delta}_{
\hat{a}_{j}}} \setminus \hat{\Sigma}_{\hat{a}_{j}} \! \to \! \operatorname{SL}_{2}
(\mathbb{C})$, $j \! = \! 1,2,\dotsc,N \! + \! 1$, where $\hat{\Sigma}_{
\hat{a}_{j}} \! := \! (\hat{\Phi}_{\hat{a}_{j}})^{-1}(\gamma_{\hat{a}_{j}}^{\ast})$, 
with $(\hat{\Phi}_{\hat{a}_{j}})^{-1}$ denoting the inverse mapping, and 
$\gamma_{\hat{a}_{j}}^{\ast} \! := \! \cup_{m=1}^{4} \gamma_{\hat{a}_{j}}^{
\ast,m}$$:$ {\rm (i)} $\hat{\mathcal{X}}^{\hat{a}}(z)$ is analytic for $z \! \in 
\! \hat{\mathbb{U}}_{\hat{\delta}_{\hat{a}_{j}}} \setminus \hat{\Sigma}_{
\hat{a}_{j}}$, $j \! = \! 1,2,\dotsc,N \! + \! 1$$;$ {\rm (ii)} the boundary 
values $\hat{\mathcal{X}}^{\hat{a}}_{\pm}(z) \! := \! \lim_{\underset{
z^{\prime} \, \in \, \pm \, \mathrm{side} \, \mathrm{of} \, \hat{\Sigma}_{
\hat{a}_{j}}}{z^{\prime} \to z \in \hat{\Sigma}_{\hat{a}_{j}}}} \hat{\mathcal{
X}}^{\hat{a}}(z^{\prime})$ satisfy the jump condition
\begin{equation*}
\hat{\mathcal{X}}^{\hat{a}}_{+}(z) \! = \! \hat{\mathcal{X}}^{\hat{a}}_{-}(z) 
\mathfrak{v}(z), \quad z \! \in \! \hat{\mathbb{U}}_{\hat{\delta}_{\hat{a}_{j}}} 
\cap \hat{\Sigma}_{\hat{a}_{j}}, \quad j \! = \! 1,2,\dotsc,N \! + \! 1;
\end{equation*}
and {\rm (iii)} uniformly for $z \! \in \! \partial \hat{\mathbb{U}}_{
\hat{\delta}_{\hat{a}_{j}}} \! := \! \lbrace \mathstrut z \! \in \! 
\mathbb{C}; \, \lvert z \! - \! \hat{a}_{j} \rvert \! = \! \hat{\delta}_{
\hat{a}_{j}} \rbrace$, $j \! = \! 1,2,\dotsc,N \! + \! 1$,
\begin{equation*}
\mathfrak{m}(z)(\hat{\mathcal{X}}^{\hat{a}}(z))^{-1} \underset{\underset{
z_{o}=1+o(1)}{\mathscr{N},n \to \infty}}{=} \mathrm{I} \! + \! o(1).
\end{equation*}
For $n \! \in \! \mathbb{N}$ and $k \! \in \! \lbrace 1,2,\dotsc,K \rbrace$ 
such that $\alpha_{p_{\mathfrak{s}}} \! := \! \alpha_{k} \! = \! \infty$, 
the solutions of the {\rm RHPs} $(\hat{\mathcal{X}}^{\hat{a}}(z),\mathfrak{v}
(z),\hat{\mathbb{U}}_{\hat{\delta}_{\hat{a}_{j}}} \cap \hat{\Sigma}_{\hat{a}_{
j}})$, $j \! = \! 1,2,\dotsc,N \! + \! 1$, are: {\rm \pmb{(1)}} for $z \! \in \! 
\hat{\Omega}_{\hat{a}_{j}}^{1} \! := \! \hat{\mathbb{U}}_{\hat{\delta}_{
\hat{a}_{j}}} \cap (\hat{\Phi}_{\hat{a}_{j}})^{-1}(\hat{\Omega}_{\hat{a}_{
j}}^{\ast,1})$, $j \! = \! 1,2,\dotsc,N \! + \! 1$,
\begin{equation*}
\hat{\mathcal{X}}^{\hat{a}}(z) \! = \! \sqrt{\smash[b]{\pi}} \, \me^{-
\frac{\mi \pi}{3}} \mathfrak{m}(z) \me^{\frac{\mi}{2}((n-1)K+k) 
\hat{\mho}_{j} \operatorname{ad}(\sigma_{3})} 
\begin{pmatrix}
\mi & -\mi \\
1 & 1
\end{pmatrix}(\hat{\Phi}_{\hat{a}_{j}}(z))^{\frac{1}{4} \sigma_{3}} \Psi_{1}
(\hat{\Phi}_{\hat{a}_{j}}(z)) \me^{\frac{1}{2}((n-1)K+k) \hat{\xi}_{\hat{a}_{j}}
(z) \sigma_{3}},
\end{equation*}
where $\mathfrak{m}(z)$ is given in item~{\rm \pmb{(1)}} of 
Lemma~\ref{lem4.5}$;$ {\rm \pmb{(2)}} for $z \! \in \! \hat{\Omega}_{
\hat{a}_{j}}^{2} \! := \! \hat{\mathbb{U}}_{\hat{\delta}_{\hat{a}_{j}}} \cap 
(\hat{\Phi}_{\hat{a}_{j}})^{-1}(\hat{\Omega}_{\hat{a}_{j}}^{\ast,2})$, $j \! 
= \! 1,2,\dotsc,N \! + \! 1$,
\begin{equation*}
\hat{\mathcal{X}}^{\hat{a}}(z) \! = \! \sqrt{\smash[b]{\pi}} \, \me^{-
\frac{\mi \pi}{3}} \mathfrak{m}(z) \me^{\frac{\mi}{2}((n-1)K+k) 
\hat{\mho}_{j} \operatorname{ad}(\sigma_{3})} 
\begin{pmatrix}
\mi & -\mi \\
1 & 1
\end{pmatrix}(\hat{\Phi}_{\hat{a}_{j}}(z))^{\frac{1}{4} \sigma_{3}} \Psi_{2}
(\hat{\Phi}_{\hat{a}_{j}}(z)) \me^{\frac{1}{2}((n-1)K+k) \hat{\xi}_{\hat{a}_{j}}
(z) \sigma_{3}};
\end{equation*}
{\rm \pmb{(3)}} for $z \! \in \! \hat{\Omega}_{\hat{a}_{j}}^{3} \! := \! 
\hat{\mathbb{U}}_{\hat{\delta}_{\hat{a}_{j}}} \cap (\hat{\Phi}_{\hat{a}_{j}})^{-1}
(\hat{\Omega}_{\hat{a}_{j}}^{\ast,3})$, $j \! = \! 1,2,\dotsc,N \! + \! 1$,
\begin{equation*}
\hat{\mathcal{X}}^{\hat{a}}(z) \! = \! \sqrt{\smash[b]{\pi}} \, \me^{-
\frac{\mi \pi}{3}} \mathfrak{m}(z) \me^{-\frac{\mi}{2}((n-1)K+k) 
\hat{\mho}_{j} \operatorname{ad}(\sigma_{3})} 
\begin{pmatrix}
\mi & -\mi \\
1 & 1
\end{pmatrix} 
(\hat{\Phi}_{\hat{a}_{j}}(z))^{\frac{1}{4} \sigma_{3}} \Psi_{3}(\hat{\Phi}_{
\hat{a}_{j}}(z)) \me^{\frac{1}{2}((n-1)K+k) \hat{\xi}_{\hat{a}_{j}}(z) 
\sigma_{3}};
\end{equation*}
and {\rm \pmb{(4)}} for $z \! \in \! \hat{\Omega}_{\hat{a}_{j}}^{4} \! := \! 
\hat{\mathbb{U}}_{\hat{\delta}_{\hat{a}_{j}}} \cap (\hat{\Phi}_{\hat{a}_{j}})^{-1}
(\hat{\Omega}_{\hat{a}_{j}}^{\ast,4})$, $j \! = \! 1,2,\dotsc,N \! + \! 1$,
\begin{equation*}
\hat{\mathcal{X}}^{\hat{a}}(z) \! = \! \sqrt{\smash[b]{\pi}} \, \me^{-
\frac{\mi \pi}{3}}\mathfrak{m}(z) \me^{-\frac{\mi}{2}((n-1)K+k) 
\hat{\mho}_{j} \operatorname{ad}(\sigma_{3})} 
\begin{pmatrix}
\mi & -\mi \\
1 & 1
\end{pmatrix}(\hat{\Phi}_{\hat{a}_{j}}(z))^{\frac{1}{4} \sigma_{3}} \Psi_{4}
(\hat{\Phi}_{\hat{a}_{j}}(z)) \me^{\frac{1}{2}((n-1)K+k) \hat{\xi}_{\hat{a}_{j}}
(z) \sigma_{3}}.
\end{equation*}
\end{ccccc}
\begin{ccccc} \label{lem4.8} 
For $n \! \in \! \mathbb{N}$ and $k \! \in \! \lbrace 1,2,\dotsc,K 
\rbrace$ such that $\alpha_{p_{\mathfrak{s}}} \! := \! \alpha_{k} \! \neq 
\! \infty$, let $\mathfrak{M} \colon \mathbb{C} \setminus \tilde{\Sigma} 
\! \to \! \operatorname{SL}_{2}(\mathbb{C})$ solve the {\rm RHP} 
$(\mathfrak{M}(z),\mathfrak{v}(z),\tilde{\Sigma})$ stated in 
item~{\rm \pmb{(2)}} of Lemma~\ref{lem4.2}. For $n \! \in \! 
\mathbb{N}$ and $k \! \in \! \lbrace 1,2,\dotsc,K \rbrace$ such that 
$\alpha_{p_{\mathfrak{s}}} \! := \! \alpha_{k} \! \neq \! \infty$, set 
$\tilde{\mathbb{U}}_{\tilde{\delta}_{\tilde{b}_{i-1}}} \! := \! \lbrace 
\mathstrut z \! \in \! \mathbb{C}; \, \lvert z \! - \! \tilde{b}_{i-1} 
\rvert \! < \! \tilde{\delta}_{\tilde{b}_{i-1}} \rbrace$, $i \! = \! 1,2,
\dotsc,N \! + \! 1$, and let $\tilde{\Phi}_{\tilde{b}_{j-1}}(z)$ and 
$\tilde{\xi}_{\tilde{b}_{j-1}}(z)$, $j \! = \! 1,2,\dotsc,N \! + \! 1$, be 
defined by Equations~\eqref{eqmainfin78} and~\eqref{eqmainfin72}, 
respectively, where, for $z \! \in \! \tilde{\mathbb{U}}_{
\tilde{\delta}_{\tilde{b}_{j-1}}} \setminus (-\infty,\tilde{b}_{j-1})$, 
$\tilde{\xi}_{\tilde{b}_{j-1}}(z) \! = \! \tilde{\mathfrak{b}}(z \! - \! 
\tilde{b}_{j-1})^{3/2} \tilde{G}_{\tilde{b}_{j-1}}(z)$, with $\tilde{
\mathfrak{b}} \! = \! \pm 1$ for $z \! \in \! \mathbb{C}_{\pm}$, 
and $\tilde{G}_{\tilde{b}_{j-1}}(z)$ is analytic, in particular,
\begin{equation*}
\tilde{G}_{\tilde{b}_{j-1}}(z) \underset{z \to \tilde{b}_{j-1}}{=} \dfrac{2}{3}
f(\tilde{b}_{j-1}) \! + \! \dfrac{2}{5}f^{\prime}(\tilde{b}_{j-1})(z \! - \! 
\tilde{b}_{j-1}) \! + \! \dfrac{1}{7}f^{\prime \prime}(\tilde{b}_{j-1})(z \! 
- \! \tilde{b}_{j-1})^{2} \! + \! \mathcal{O}((z \! - \! \tilde{b}_{j-1})^{3}),
\end{equation*}
where
\begin{align}
f(\tilde{b}_{0}) =& \, \mi (-1)^{N} \tilde{h}_{\widetilde{V}}(\tilde{b}_{0}) 
\tilde{\eta}_{\tilde{b}_{0}}, \nonumber \\
f^{\prime}(\tilde{b}_{0}) =& \, \mi (-1)^{N} \left(\dfrac{1}{2} \tilde{h}_{
\widetilde{V}}(\tilde{b}_{0}) \tilde{\eta}_{\tilde{b}_{0}} \left(\sum_{m=1}^{N} 
\left(\dfrac{1}{\tilde{b}_{0} \! - \! \tilde{b}_{m}} \! + \! \dfrac{1}{\tilde{b}_{0} 
\! - \! \tilde{a}_{m}} \right) \! + \! \dfrac{1}{\tilde{b}_{0} \! - \! \tilde{a}_{N
+1}} \right) \! + \! (\tilde{h}_{\widetilde{V}}(\tilde{b}_{0}))^{\prime} 
\tilde{\eta}_{\tilde{b}_{0}} \right), \nonumber \\
f^{\prime \prime}(\tilde{b}_{0}) =& \, \mi (-1)^{N} \left(\dfrac{\tilde{h}_{
\widetilde{V}}(\tilde{b}_{0})(\tilde{h}_{\widetilde{V}}(\tilde{b}_{0}))^{\prime 
\prime} \! - \! ((\tilde{h}_{\widetilde{V}}(\tilde{b}_{0}))^{\prime})^{2}}{
\tilde{h}_{\widetilde{V}}(\tilde{b}_{0})} \tilde{\eta}_{\tilde{b}_{0}} \! - \! 
\dfrac{1}{2} \tilde{h}_{\widetilde{V}}(\tilde{b}_{0}) \tilde{\eta}_{\tilde{b}_{0}} 
\right. \nonumber \\
\times&\left. \, \left(\sum_{m=1}^{N} \left(\dfrac{1}{(\tilde{b}_{0} \! - \! 
\tilde{b}_{m})^{2}} \! + \! \dfrac{1}{(\tilde{b}_{0} \! - \! \tilde{a}_{m})^{2}} 
\right) \! + \! \dfrac{1}{(\tilde{b}_{0} \! - \! \tilde{a}_{N+1})^{2}} \right) 
\right. \nonumber \\
+&\left. \, \left(\dfrac{1}{2} \left(\sum_{m=1}^{N} \left(\dfrac{1}{
\tilde{b}_{0} \! - \! \tilde{b}_{m}} \! + \! \dfrac{1}{\tilde{b}_{0} \! - \! 
\tilde{a}_{m}} \right) \! + \! \dfrac{1}{\tilde{b}_{0} \! - \! \tilde{a}_{N+1}} 
\right) \! + \! \dfrac{(\tilde{h}_{\widetilde{V}}(\tilde{b}_{0}))^{\prime}}{
\tilde{h}_{\widetilde{V}}(\tilde{b}_{0})} \right) \right. \nonumber \\
\times&\left. \, \left(\dfrac{1}{2} \tilde{h}_{\widetilde{V}}(\tilde{b}_{0}) 
\tilde{\eta}_{\tilde{b}_{0}} \left(\sum_{m=1}^{N} \left(\dfrac{1}{\tilde{b}_{0} 
\! - \! \tilde{b}_{m}} \! + \! \dfrac{1}{\tilde{b}_{0} \! - \! \tilde{a}_{m}} 
\right) \! + \! \dfrac{1}{\tilde{b}_{0} \! - \! \tilde{a}_{N+1}} \right) \! + \! 
(\tilde{h}_{\widetilde{V}}(\tilde{b}_{0}))^{\prime} \tilde{\eta}_{\tilde{b}_{0}} 
\right) \right), \label{wm24frbotil}
\end{align}
with $\tilde{\eta}_{\tilde{b}_{0}}$ defined by Equation~\eqref{eqmainfin57}, 
and, for $j \! = \! 1,2,\dotsc,N$,
\begin{align}
f(\tilde{b}_{j}) =& \, \mi (-1)^{N-j} \tilde{h}_{\widetilde{V}}(\tilde{b}_{j}) 
\tilde{\eta}_{\tilde{b}_{j}}, \nonumber \\
f^{\prime}(\tilde{b}_{j}) =& \, \mi (-1)^{N-j} \left(\dfrac{1}{2} \tilde{h}_{
\widetilde{V}}(\tilde{b}_{j}) \tilde{\eta}_{\tilde{b}_{j}} \left(\sum_{
\substack{m=1\\m \not= j}}^{N} \left(\dfrac{1}{\tilde{b}_{j} \! - \! 
\tilde{b}_{m}} \! + \! \dfrac{1}{\tilde{b}_{j} \! - \! \tilde{a}_{m}} \right) \! 
+ \! \dfrac{1}{\tilde{b}_{j} \! - \! \tilde{a}_{j}} \! + \! \dfrac{1}{\tilde{b}_{j} 
\! - \! \tilde{a}_{N+1}} \! + \! \dfrac{1}{\tilde{b}_{j} \! - \! \tilde{b}_{0}} 
\right) \! + \! (\tilde{h}_{\widetilde{V}}(\tilde{b}_{j}))^{\prime} \tilde{
\eta}_{\tilde{b}_{j}} \right), \nonumber \\
f^{\prime \prime}(\tilde{b}_{j}) =& \, \mi (-1)^{N-j} \left(\dfrac{
\tilde{h}_{\widetilde{V}}(\tilde{b}_{j})(\tilde{h}_{\widetilde{V}}(\tilde{b}_{j}))^{
\prime \prime} \! - \! ((\tilde{h}_{\widetilde{V}}(\tilde{b}_{j}))^{\prime})^{2}}{
\tilde{h}_{\widetilde{V}}(\tilde{b}_{j})} \tilde{\eta}_{\tilde{b}_{j}} \! - \! 
\dfrac{1}{2} \tilde{h}_{\widetilde{V}}(\tilde{b}_{j}) \tilde{\eta}_{\tilde{b}_{j}} 
\left(\sum_{\substack{m=1\\m \not= j}}^{N} \left(\dfrac{1}{(\tilde{b}_{j} \! 
- \! \tilde{b}_{m})^{2}} \! + \! \dfrac{1}{(\tilde{b}_{j} \! - \! \tilde{a}_{m})^{2}} 
\right) \right. \right. \nonumber \\
+&\left. \left. \, \dfrac{1}{(\tilde{b}_{j} \! - \! \tilde{a}_{j})^{2}} \! + \! 
\dfrac{1}{(\tilde{b}_{j} \! - \! \tilde{a}_{N+1})^{2}} \! + \! \dfrac{1}{(\tilde{b}_{j} 
\! - \! \tilde{b}_{0})^{2}} \right) \! + \! \left(\dfrac{(\tilde{h}_{\widetilde{V}}
(\tilde{b}_{j}))^{\prime}}{\tilde{h}_{\widetilde{V}}(\tilde{b}_{j})} \! + \! 
\dfrac{1}{2} \left(\sum_{\substack{m=1\\m \not= j}}^{N} \left(\dfrac{1}{
\tilde{b}_{j} \! - \! \tilde{b}_{m}} \! + \! \dfrac{1}{\tilde{b}_{j} \! - \! 
\tilde{a}_{m}} \right) \right. \right. \right. \nonumber \\
+&\left. \left. \left. \, \dfrac{1}{\tilde{b}_{j} \! - \! \tilde{a}_{j}} \! + \! 
\dfrac{1}{\tilde{b}_{j} \! - \! \tilde{a}_{N+1}} \! + \! \dfrac{1}{\tilde{b}_{j} 
\! - \! \tilde{b}_{0}} \right) \right) \left(\dfrac{1}{2} \tilde{h}_{\widetilde{V}}
(\tilde{b}_{j}) \tilde{\eta}_{\tilde{b}_{j}} \left(\sum_{\substack{m=1\\m 
\not= j}}^{N} \left(\dfrac{1}{\tilde{b}_{j} \! - \! \tilde{b}_{m}} \! + \! 
\dfrac{1}{\tilde{b}_{j} \! - \! \tilde{a}_{m}} \right) \right. \right. \right. 
\nonumber \\
+&\left. \left. \left. \, \dfrac{1}{\tilde{b}_{j} \! - \! \tilde{a}_{j}} \! + \! 
\dfrac{1}{\tilde{b}_{j} \! - \! \tilde{a}_{N+1}} \! + \! \dfrac{1}{\tilde{b}_{j} 
\! - \! \tilde{b}_{0}} \right) \! + \!  (\tilde{h}_{\widetilde{V}}(\tilde{b}_{j}))^{
\prime} \tilde{\eta}_{\tilde{b}_{j}} \right) \right), \label{wm24frbjtil}
\end{align}
with $\tilde{\eta}_{\tilde{b}_{j}}$ defined by Equation~\eqref{eqmainfin58}, 
and $(0,1) \! \ni \! \tilde{\delta}_{\tilde{b}_{j-1}}$, $j \! = \! 1,2,\dotsc,N 
\! + \! 1$, are chosen sufficiently small so that $\tilde{\Phi}_{\tilde{b}_{j-1}}
(z)$, which are bi-holomorphic, conformal, and non-orientation preserving, 
map (cf. Figure~\ref{forbtil}$)$ $\tilde{\mathbb{U}}_{\tilde{\delta}_{
\tilde{b}_{j-1}}}$, and, thus, the oriented contours $\tilde{\Sigma}_{
\tilde{b}_{j-1}} \! := \! \cup_{m=1}^{4} \tilde{\Sigma}_{\tilde{b}_{j-1}}^{m}$, 
$j \! = \! 1,2,\dotsc,N \! + \! 1$, injectively onto open $(n$-, $k$-, and 
$z_{o}$-dependent) neighbourhoods $\tilde{\mathbb{U}}_{\tilde{\delta}_{
\tilde{b}_{j-1}}}^{\ast}$ of the origin, $j \! = \! 1,2,\dotsc,N \! + \! 1$, 
such that $\tilde{\Phi}_{\tilde{b}_{j-1}}(\tilde{b}_{j-1}) \! = \! 0$, 
$\tilde{\Phi}_{\tilde{b}_{j-1}} \colon \tilde{\mathbb{U}}_{\tilde{\delta}_{
\tilde{b}_{j-1}}} \! \to \! \tilde{\mathbb{U}}_{\tilde{\delta}_{\tilde{b}_{j-
1}}}^{\ast} \! := \! \tilde{\Phi}_{\tilde{b}_{j-1}}(\tilde{\mathbb{U}}_{\tilde{
\delta}_{\tilde{b}_{j-1}}})$, $\tilde{\Phi}_{\tilde{b}_{j-1}}(\tilde{\mathbb{U}}_{
\tilde{\delta}_{\tilde{b}_{j-1}}} \cap \tilde{\Sigma}_{\tilde{b}_{j-1}}^{m}) \! 
= \! \tilde{\Phi}_{\tilde{b}_{j-1}}(\tilde{\mathbb{U}}_{\tilde{\delta}_{\tilde{
b}_{j-1}}}) \cap \gamma_{\tilde{b}_{j-1}}^{\ast,m}$, and $\tilde{\Phi}_{
\tilde{b}_{j-1}}(\tilde{\mathbb{U}}_{\tilde{\delta}_{\tilde{b}_{j-1}}} \cap 
\tilde{\Omega}_{\tilde{b}_{j-1}}^{m}) \! = \! \tilde{\Phi}_{\tilde{b}_{j-1}}
(\tilde{\mathbb{U}}_{\tilde{\delta}_{\tilde{b}_{j-1}}}) \cap \tilde{\Omega}_{
\tilde{b}_{j-1}}^{\ast,m}$, $m \! = \! 1,2,3,4$, with $\tilde{\Omega}_{
\tilde{b}_{j-1}}^{\ast,1} \! = \! \lbrace \mathstrut \zeta \! \in \! 
\mathbb{C}; \, \arg (\zeta) \! \in \! (0,2 \pi/3) \rbrace$, $\tilde{\Omega}_{
\tilde{b}_{j-1}}^{\ast,2} \! = \! \lbrace \mathstrut \zeta \! \in \! \mathbb{C}; 
\, \arg (\zeta) \! \in \! (2 \pi/3,\pi) \rbrace$, $\tilde{\Omega}_{\tilde{b}_{j-
1}}^{\ast,3} \! = \! \lbrace \mathstrut \zeta \! \in \! \mathbb{C}; \, \arg 
(\zeta) \! \in \! (-\pi,-2 \pi/3) \rbrace$, and $\tilde{\Omega}_{\tilde{b}_{j-
1}}^{\ast,4} \! = \! \lbrace \mathstrut \zeta \! \in \! \mathbb{C}; \, \arg 
(\zeta) \! \in \! (-2 \pi/3,0) \rbrace$.

For $n \! \in \! \mathbb{N}$ and $k \! \in \! \lbrace 1,2,\dotsc,K \rbrace$ 
such that $\alpha_{p_{\mathfrak{s}}} \! := \! \alpha_{k} \! \neq \! \infty$, 
the parametrix for the {\rm RHP} $(\mathfrak{M}(z),\mathfrak{v}(z),
\tilde{\Sigma})$, for $z \! \in \! \tilde{\mathbb{U}}_{\tilde{\delta}_{
\tilde{b}_{j-1}}}$, $j \! = \! 1,2,\dotsc,N \! + \! 1$, is the solution of the 
following {\rm RHPs} for $\tilde{\mathcal{X}}^{\tilde{b}} \colon \tilde{
\mathbb{U}}_{\tilde{\delta}_{\tilde{b}_{j-1}}} \setminus \tilde{\Sigma}_{
\tilde{b}_{j-1}} \! \to \! \operatorname{SL}_{2}(\mathbb{C})$, $j \! = \! 
1,2,\dotsc,N \! + \! 1$, where $\tilde{\Sigma}_{\tilde{b}_{j-1}} \! := \! 
(\tilde{\Phi}_{\tilde{b}_{j-1}})^{-1}(\gamma_{\tilde{b}_{j-1}}^{\ast})$, with 
$(\tilde{\Phi}_{\tilde{b}_{j-1}})^{-1}$ denoting the inverse mapping, and 
$\gamma_{\tilde{b}_{j-1}}^{\ast} \! := \! \cup_{m=1}^{4} \gamma_{
\tilde{b}_{j-1}}^{\ast,m}$$:$ {\rm (i)} $\tilde{\mathcal{X}}^{\tilde{b}}(z)$ is 
analytic for $z \! \in \! \tilde{\mathbb{U}}_{\tilde{\delta}_{\tilde{b}_{j-1}}} 
\setminus \tilde{\Sigma}_{\tilde{b}_{j-1}}$, $j \! = \! 1,2,\dotsc,N \! + \! 
1$$;$ {\rm (ii)} the boundary values $\tilde{\mathcal{X}}^{\tilde{b}}_{\pm}
(z) \! := \! \lim_{\underset{z^{\prime} \, \in \, \pm \, \mathrm{side} \, 
\mathrm{of} \, \tilde{\Sigma}_{\tilde{b}_{j-1}}}{z^{\prime} \to z \in \tilde{
\Sigma}_{\tilde{b}_{j-1}}}} \tilde{\mathcal{X}}^{\tilde{b}}(z^{\prime})$ satisfy 
the jump condition
\begin{equation*}
\tilde{\mathcal{X}}^{\tilde{b}}_{+}(z) \! = \! \tilde{\mathcal{X}}^{\tilde{b}}_{-}
(z) \mathfrak{v}(z), \quad z \! \in \! \tilde{\mathbb{U}}_{\tilde{\delta}_{
\tilde{b}_{j-1}}} \cap \tilde{\Sigma}_{\tilde{b}_{j-1}}, \quad j \! = \! 1,2,
\dotsc,N \! + \! 1;
\end{equation*}
and {\rm (iii)} uniformly for $z \! \in \! \partial \tilde{\mathbb{U}}_{
\tilde{\delta}_{\tilde{b}_{j-1}}} \! := \! \lbrace \mathstrut z \! \in \! 
\mathbb{C}; \, \vert z \! - \! \tilde{b}_{j-1} \vert \! = \! \tilde{\delta}_{
\tilde{b}_{j-1}} \rbrace$, $j \! = \! 1,2,\dotsc,N \! + \! 1$,
\begin{equation*}
\mathfrak{m}(z)(\tilde{\mathcal{X}}^{\tilde{b}}(z))^{-1} \underset{
\underset{z_{o}=1+o(1)}{\mathscr{N},n \to \infty}}{=} \mathrm{I} \! + \! o(1).
\end{equation*}
For $n \! \in \! \mathbb{N}$ and $k \! \in \! \lbrace 1,2,\dotsc,K \rbrace$ 
such that $\alpha_{p_{\mathfrak{s}}} \! := \! \alpha_{k} \! \neq \! \infty$, 
the solutions of the {\rm RHPs} $(\tilde{\mathcal{X}}^{\tilde{b}}(z),\mathfrak{v}
(z),\tilde{\mathbb{U}}_{\tilde{\delta}_{\tilde{b}_{j-1}}} \cap \tilde{\Sigma}_{
\tilde{b}_{j-1}})$, $j \! = \! 1,2,\dotsc,N \! + \! 1$, are: {\rm \pmb{(1)}} for 
$z \! \in \! \tilde{\Omega}_{\tilde{b}_{j-1}}^{1} \! := \! \tilde{\mathbb{
U}}_{\tilde{\delta}_{\tilde{b}_{j-1}}} \cap (\tilde{\Phi}_{\tilde{b}_{j-1}})^{-1}
(\tilde{\Omega}_{\tilde{b}_{j-1}}^{\ast,1})$, $j \! = \! 1,2,\dotsc,N \! + \! 1$,
\begin{equation*}
\tilde{\mathcal{X}}^{\tilde{b}}(z) \! = \! \sqrt{\smash[b]{\pi}} \, \me^{-
\frac{\mi \pi}{3}} \mathfrak{m}(z) \sigma_{3} \me^{\frac{\mi}{2}((n-1)K+k) 
\tilde{\mho}_{j-1} \operatorname{ad}(\sigma_{3})} \begin{pmatrix}
\mi & -\mi \\
1 & 1
\end{pmatrix}(\tilde{\Phi}_{\tilde{b}_{j-1}}(z))^{\frac{1}{4} \sigma_{3}} 
\Psi_{1}(\tilde{\Phi}_{\tilde{b}_{j-1}}(z)) \me^{\frac{1}{2}((n-1)K+k) 
\tilde{\xi}_{\tilde{b}_{j-1}}(z) \sigma_{3}} \sigma_{3},
\end{equation*}
where $\mathfrak{m}(z)$ is given in item~{\rm \pmb{(2)}} of 
Lemma~\ref{lem4.5}$;$ {\rm \pmb{(2)}} for $z \! \in \! \tilde{
\Omega}_{\tilde{b}_{j-1}}^{2} \! := \! \tilde{\mathbb{U}}_{\tilde{\delta}_{
\tilde{b}_{j-1}}} \cap (\tilde{\Phi}_{\tilde{b}_{j-1}})^{-1}(\tilde{\Omega}_{
\tilde{b}_{j-1}}^{\ast,2})$, $j \! = \! 1,2,\dotsc,N \! + \! 1$,
\begin{equation*}
\tilde{\mathcal{X}}^{\tilde{b}}(z) \! = \! \sqrt{\smash[b]{\pi}} \, \me^{-
\frac{\mi \pi}{3}} \mathfrak{m}(z) \sigma_{3} \me^{\frac{\mi}{2}((n-1)K+k) 
\tilde{\mho}_{j-1} \operatorname{ad}(\sigma_{3})} 
\begin{pmatrix}
\mi & -\mi \\
1 & 1
\end{pmatrix}(\tilde{\Phi}_{\tilde{b}_{j-1}}(z))^{\frac{1}{4} \sigma_{3}} 
\Psi_{2}(\tilde{\Phi}_{\tilde{b}_{j-1}}(z)) \me^{\frac{1}{2}((n-1)K+k) 
\tilde{\xi}_{\tilde{b}_{j-1}}(z) \sigma_{3}} \sigma_{3};
\end{equation*}
{\rm \pmb{(3)}} for $z \! \in \! \tilde{\Omega}_{\tilde{b}_{j-1}}^{3} \! := 
\! \tilde{\mathbb{U}}_{\tilde{\delta}_{\tilde{b}_{j-1}}} \cap (\tilde{\Phi}_{
\tilde{b}_{j-1}})^{-1}(\tilde{\Omega}_{\tilde{b}_{j-1}}^{\ast,3})$, $j \! = \! 
1,2,\dotsc,N \! + \! 1$,
\begin{equation*}
\tilde{\mathcal{X}}^{\tilde{b}}(z) \! = \! \sqrt{\smash[b]{\pi}} \, \me^{-
\frac{\mi \pi}{3}} \mathfrak{m}(z) \sigma_{3} \me^{-\frac{\mi}{2}((n-1)K+k) 
\tilde{\mho}_{j-1} \operatorname{ad}(\sigma_{3})} 
\begin{pmatrix}
\mi & -\mi \\
1 & 1
\end{pmatrix}(\tilde{\Phi}_{\tilde{b}_{j-1}}(z))^{\frac{1}{4} \sigma_{3}} \Psi_{3}
(\tilde{\Phi}_{\tilde{b}_{j-1}}(z)) \me^{\frac{1}{2}((n-1)K+k) \tilde{\xi}_{
\tilde{b}_{j-1}}(z) \sigma_{3}} \sigma_{3};
\end{equation*}
and {\rm \pmb{(4)}} for $z \! \in \! \tilde{\Omega}_{\tilde{b}_{j-1}}^{4} \! 
:= \! \tilde{\mathbb{U}}_{\tilde{\delta}_{\tilde{b}_{j-1}}} \cap (\tilde{\Phi}_{
\tilde{b}_{j-1}})^{-1}(\tilde{\Omega}_{\tilde{b}_{j-1}}^{\ast,4})$, $j \! = \! 
1,2,\dotsc,N \! + \! 1$,
\begin{equation*}
\tilde{\mathcal{X}}^{\tilde{b}}(z) \! = \! \sqrt{\smash[b]{\pi}} \, \me^{-
\frac{\mi \pi}{3}} \mathfrak{m}(z) \sigma_{3} \me^{-\frac{\mi}{2}((n-1)K+k) 
\tilde{\mho}_{j-1} \operatorname{ad}(\sigma_{3})} 
\begin{pmatrix}
\mi & -\mi \\
1 & 1
\end{pmatrix}(\tilde{\Phi}_{\tilde{b}_{j-1}}(z))^{\frac{1}{4} \sigma_{3}} \Psi_{4}
(\tilde{\Phi}_{\tilde{b}_{j-1}}(z)) \me^{\frac{1}{2}((n-1)K+k) \tilde{\xi}_{
\tilde{b}_{j-1}}(z) \sigma_{3}} \sigma_{3}.
\end{equation*}
\end{ccccc}
\begin{ccccc} \label{lem4.9} 
For $n \! \in \! \mathbb{N}$ and $k \! \in \! \lbrace 1,2,\dotsc,K 
\rbrace$ such that $\alpha_{p_{\mathfrak{s}}} \! := \! \alpha_{k} \! \neq 
\! \infty$, let $\mathfrak{M} \colon \mathbb{C} \setminus \tilde{\Sigma} 
\! \to \! \operatorname{SL}_{2}(\mathbb{C})$ solve the {\rm RHP} 
$(\mathfrak{M}(z),\mathfrak{v}(z),\tilde{\Sigma})$ stated in 
item~{\rm \pmb{(2)}} of Lemma~\ref{lem4.2}. For $n \! \in \! 
\mathbb{N}$ and $k \! \in \! \lbrace 1,2,\dotsc,K \rbrace$ such that 
$\alpha_{p_{\mathfrak{s}}} \! := \! \alpha_{k} \! \neq \! \infty$, set 
$\tilde{\mathbb{U}}_{\tilde{\delta}_{\tilde{a}_{i}}} \! := \! \lbrace \mathstrut 
z \! \in \! \mathbb{C}; \, \lvert z \! - \! \tilde{a}_{i} \rvert \! < \! \tilde{
\delta}_{\tilde{a}_{i}} \rbrace$, $i \! = \! 1,2,\dotsc,N \! + \! 1$, and let 
$\tilde{\Phi}_{\tilde{a}_{j}}(z)$ and $\tilde{\xi}_{\tilde{a}_{j}}(z)$, $j \! = \! 
1,2,\dotsc,N \! + \! 1$, be defined by Equations~\eqref{eqmainfin85} 
and~\eqref{eqmainfin73}, respectively, where, for $z \! \in \! \tilde{
\mathbb{U}}_{\tilde{\delta}_{\tilde{a}_{j}}} \setminus (-\infty,\tilde{a}_{j})$, 
$\tilde{\xi}_{\tilde{a}_{j}}(z) \! = \! (z \! - \! \tilde{a}_{j})^{3/2} \tilde{G}_{
\tilde{a}_{j}}(z)$, and $\tilde{G}_{\tilde{a}_{j}}(z)$ is analytic, in particular,
\begin{equation*}
\tilde{G}_{\tilde{a}_{j}}(z) \underset{z \to \tilde{a}_{j}}{=} \dfrac{2}{3}
f(\tilde{a}_{j}) \! + \! \dfrac{2}{5}f^{\prime}(\tilde{a}_{j})(z \! - \! 
\tilde{a}_{j}) \! + \! \dfrac{1}{7}f^{\prime \prime}(\tilde{a}_{j})(z \! 
- \! \tilde{a}_{j})^{2} \! + \! \mathcal{O}((z \! - \! \tilde{a}_{j})^{3}),
\end{equation*}
where
\begin{align}
f(\tilde{a}_{N+1}) =& \, \tilde{h}_{\widetilde{V}}(\tilde{a}_{N+1}) 
\tilde{\eta}_{\tilde{a}_{N+1}}, \nonumber \\
f^{\prime}(\tilde{a}_{N+1}) =& \, \dfrac{1}{2} \tilde{h}_{\widetilde{V}}
(\tilde{a}_{N+1}) \tilde{\eta}_{\tilde{a}_{N+1}} \left(\sum_{m=1}^{N} 
\left(\dfrac{1}{\tilde{a}_{N+1} \! - \! \tilde{b}_{m}} \! + \! \dfrac{1}{
\tilde{a}_{N+1} \! - \! \tilde{a}_{m}} \right) \! + \! \dfrac{1}{\tilde{a}_{N
+1} \! - \! \tilde{b}_{0}} \right) \! + \! (\tilde{h}_{\widetilde{V}}(\tilde{a}_{N
+1}))^{\prime} \tilde{\eta}_{\tilde{a}_{N+1}}, \nonumber \\
f^{\prime \prime}(\tilde{a}_{N+1}) =& \, \dfrac{\tilde{h}_{\widetilde{V}}
(\tilde{a}_{N+1})(\tilde{h}_{\widetilde{V}}(\tilde{a}_{N+1}))^{\prime \prime} 
\! - \! ((\tilde{h}_{\widetilde{V}}(\tilde{a}_{N+1}))^{\prime})^{2}}{\tilde{h}_{
\widetilde{V}}(\tilde{a}_{N+1})} \tilde{\eta}_{\tilde{a}_{N+1}} \! - \! 
\dfrac{1}{2} \tilde{h}_{\widetilde{V}}(\tilde{a}_{N+1}) \tilde{\eta}_{
\tilde{a}_{N+1}} \nonumber \\
\times& \, \left(\sum_{m=1}^{N} \left(\dfrac{1}{(\tilde{a}_{N+1} \! - \! 
\tilde{b}_{m})^{2}} \! + \! \dfrac{1}{(\tilde{a}_{N+1} \! - \! \tilde{a}_{m})^{2}} 
\right) \! + \! \dfrac{1}{(\tilde{a}_{N+1} \! - \! \tilde{b}_{0})^{2}} \right) 
\nonumber \\
+& \, \left(\dfrac{1}{2} \left(\sum_{m=1}^{N} \left(\dfrac{1}{\tilde{a}_{N
+1} \! - \! \tilde{b}_{m}} \! + \! \dfrac{1}{\tilde{a}_{N+1} \! - \! \tilde{a}_{
m}} \right) \! + \! \dfrac{1}{\tilde{a}_{N+1} \! - \! \tilde{b}_{0}} \right) \! 
+ \! \dfrac{(\tilde{h}_{\widetilde{V}}(\tilde{a}_{N+1}))^{\prime}}{
\tilde{h}_{\widetilde{V}}(\tilde{a}_{N+1})} \right) \nonumber \\
\times& \, \left(\dfrac{1}{2} \tilde{h}_{\widetilde{V}}(\tilde{a}_{N+1}) 
\tilde{\eta}_{\tilde{a}_{N+1}} \left(\sum_{m=1}^{N} \! \left(\dfrac{1}{
\tilde{a}_{N+1} \! - \! \tilde{a}_{m}} \! + \! \dfrac{1}{\tilde{a}_{N+1} \! - \! 
\tilde{b}_{m}} \right) \! + \! \dfrac{1}{\tilde{a}_{N+1} \! - \! \tilde{b}_{0}} 
\right) \right. \nonumber \\
+& \left. \, (\tilde{h}_{\widetilde{V}}(\tilde{a}_{N+1}))^{\prime} 
\tilde{\eta}_{\tilde{a}_{N+1}} \right), \label{wm24frantil}
\end{align}
with $\tilde{\eta}_{\tilde{a}_{N+1}}$ defined by 
Equation~\eqref{eqmainfin59}, and, for $j \! = \! 1,2,\dotsc,N$,
\begin{align}
f(\tilde{a}_{j}) =& \, (-1)^{N+1-j} \tilde{h}_{\widetilde{V}}(\tilde{a}_{j}) 
\tilde{\eta}_{\tilde{a}_{j}}, \nonumber \\
f^{\prime}(\tilde{a}_{j}) =& \, (-1)^{N+1-j} \left(\dfrac{1}{2} \tilde{
h}_{\widetilde{V}}(\tilde{a}_{j}) \tilde{\eta}_{\tilde{a}_{j}} \left(\sum_{
\substack{m=1\\m \not= j}}^{N} \left(\dfrac{1}{\tilde{a}_{j} \! - \! 
\tilde{b}_{m}} \! + \! \dfrac{1}{\tilde{a}_{j} \! - \! \tilde{a}_{m}} \right) \! 
+ \! \dfrac{1}{\tilde{a}_{j} \! - \! \tilde{b}_{j}} \! + \! \dfrac{1}{\tilde{a}_{j} 
\! - \! \tilde{a}_{N+1}} \! + \! \dfrac{1}{\tilde{a}_{j} \! - \! \tilde{b}_{0}} 
\right) \! + \! (\tilde{h}_{\widetilde{V}}(\tilde{a}_{j}))^{\prime} \tilde{
\eta}_{\tilde{a}_{j}} \right), \nonumber \\
f^{\prime \prime}(\tilde{a}_{j}) =& \, (-1)^{N+1-j} \left(\dfrac{\tilde{
h}_{\widetilde{V}}(\tilde{a}_{j})(\tilde{h}_{\widetilde{V}}(\tilde{a}_{j}))^{
\prime \prime} \! - \! ((\tilde{h}_{\widetilde{V}}(\tilde{a}_{j}))^{\prime})^{
2}}{\tilde{h}_{\widetilde{V}}(\tilde{a}_{j})} \tilde{\eta}_{\tilde{a}_{j}} \! - \! 
\dfrac{1}{2} \tilde{h}_{\widetilde{V}}(\tilde{a}_{j}) \tilde{\eta}_{\tilde{a}_{j}} 
\left(\sum_{\substack{m=1\\m \not= j}}^{N} \left(\dfrac{1}{(\tilde{a}_{j} 
\! - \! \tilde{b}_{m})^{2}} \! + \! \dfrac{1}{(\tilde{a}_{j} \! - \! \tilde{a}_{
m})^{2}} \right) \right. \right. \nonumber \\
+&\left. \left. \, \dfrac{1}{(\tilde{a}_{j} \! - \! \tilde{b}_{j})^{2}} \! + 
\! \dfrac{1}{(\tilde{a}_{j} \! - \! \tilde{a}_{N+1})^{2}} \! + \! \dfrac{1}{
(\tilde{a}_{j} \! - \! \tilde{b}_{0})^{2}} \right) \! + \! \left(\dfrac{(\tilde{h}_{
\widetilde{V}}(\tilde{a}_{j}))^{\prime}}{\tilde{h}_{\widetilde{V}}(\tilde{a}_{j})} 
\! + \! \dfrac{1}{2} \left(\sum_{\substack{m=1\\m \not= j}}^{N} \left(
\dfrac{1}{\tilde{a}_{j} \! - \! \tilde{b}_{m}} \! + \! \dfrac{1}{\tilde{a}_{j} 
\! - \! \tilde{a}_{m}} \right) \right. \right. \right. \nonumber \\
+&\left. \left. \left. \, \dfrac{1}{\tilde{a}_{j} \! - \! \tilde{b}_{j}} \! + \! 
\dfrac{1}{\tilde{a}_{j} \! - \! \tilde{a}_{N+1}} \! + \! \dfrac{1}{\tilde{a}_{j} 
\! - \! \tilde{b}_{0}} \right) \right) \left(\dfrac{1}{2} \tilde{h}_{\widetilde{V}}
(\tilde{a}_{j}) \tilde{\eta}_{\tilde{a}_{j}} \left(\sum_{\substack{m=1\\m 
\not= j}}^{N} \left(\dfrac{1}{\tilde{a}_{j} \! - \! \tilde{b}_{m}} \! + \! 
\dfrac{1}{\tilde{a}_{j} \! - \! \tilde{a}_{m}} \right) \right. \right. \right. 
\nonumber \\
+&\left. \left. \left. \, \dfrac{1}{\tilde{a}_{j} \! - \! \tilde{b}_{j}} \! + \! 
\dfrac{1}{\tilde{a}_{j} \! - \! \tilde{a}_{N+1}} \! + \! \dfrac{1}{\tilde{a}_{j} 
\! - \! \tilde{b}_{0}} \right) \! + \!  (\tilde{h}_{\widetilde{V}}
(\tilde{a}_{j}))^{\prime} \tilde{\eta}_{\tilde{a}_{j}} \right) \right), 
\label{wm24frajtil}
\end{align}
with $\tilde{\eta}_{\tilde{a}_{j}}$ defined by Equation~\eqref{eqmainfin60}, 
and $(0,1) \! \ni \! \tilde{\delta}_{\tilde{a}_{j}}$, $j \! = \! 1,2,\dotsc,N \! 
+ \! 1$, are chosen sufficiently small so that $\tilde{\Phi}_{\tilde{a}_{j}}
(z)$, which are bi-holomorphic, conformal, and orientation preserving, 
map (cf. Figure~\ref{foratil}$)$ $\tilde{\mathbb{U}}_{\tilde{\delta}_{
\tilde{a}_{j}}}$, and, thus, the oriented contours $\tilde{\Sigma}_{
\tilde{a}_{j}} \! := \! \cup_{m=1}^{4} \tilde{\Sigma}_{\tilde{a}_{j}}^{m}$, 
$j \! = \! 1,2,\dotsc,N \! + \! 1$, injectively onto open $(n$-, $k$-, and 
$z_{o}$-dependent) neighbourhoods $\tilde{\mathbb{U}}_{\tilde{\delta}_{
\tilde{a}_{j}}}^{\ast}$ of the origin, $j \! = \! 1,2,\dotsc,N \! + \! 1$, 
such that $\tilde{\Phi}_{\tilde{a}_{j}}(\tilde{a}_{j}) \! = \! 0$, $\tilde{
\Phi}_{\tilde{a}_{j}} \colon \tilde{\mathbb{U}}_{\tilde{\delta}_{\tilde{a}_{j}}} 
\! \to \! \tilde{\mathbb{U}}_{\tilde{\delta}_{\tilde{a}_{j}}}^{\ast} \! := \! 
\tilde{\Phi}_{\tilde{a}_{j}}(\tilde{\mathbb{U}}_{\tilde{\delta}_{\tilde{a}_{j}}})$, 
$\tilde{\Phi}_{\tilde{a}_{j}}(\tilde{\mathbb{U}}_{\tilde{\delta}_{\tilde{a}_{j}}} 
\cap \tilde{\Sigma}_{\tilde{a}_{j}}^{m}) \! = \! \tilde{\Phi}_{\tilde{a}_{j}}
(\tilde{\mathbb{U}}_{\tilde{\delta}_{\tilde{a}_{j}}}) \cap \gamma_{
\tilde{a}_{j}}^{\ast,m}$, and $\tilde{\Phi}_{\tilde{a}_{j}}(\tilde{\mathbb{
U}}_{\tilde{\delta}_{\tilde{a}_{j}}} \cap \tilde{\Omega}_{\tilde{a}_{j}}^{m}) 
\! = \! \tilde{\Phi}_{\tilde{a}_{j}}(\tilde{\mathbb{U}}_{\tilde{\delta}_{
\tilde{a}_{j}}}) \cap \tilde{\Omega}_{\tilde{a}_{j}}^{\ast,m}$, $m \! = \! 
1,2,3,4$, with $\tilde{\Omega}_{\tilde{a}_{j}}^{\ast,1} \! = \! \lbrace 
\mathstrut \zeta \! \in \! \mathbb{C}; \, \arg (\zeta) \! \in \! (0,
2 \pi/3) \rbrace$, $\tilde{\Omega}_{\tilde{a}_{j}}^{\ast,2} \! = \! 
\lbrace \mathstrut \zeta \! \in \! \mathbb{C}; \, \arg (\zeta) \! \in \! 
(2 \pi/3,\pi) \rbrace$, $\tilde{\Omega}_{\tilde{a}_{j}}^{\ast,3} \! = \! 
\lbrace \mathstrut \zeta \! \in \! \mathbb{C}; \, \arg (\zeta) \! \in \! 
(-\pi,-2 \pi/3) \rbrace$, and $\tilde{\Omega}_{\tilde{a}_{j}}^{\ast,4} \! 
= \! \lbrace \mathstrut \zeta \! \in \! \mathbb{C}; \, \arg (\zeta) \! \in 
\! (-2 \pi/3,0) \rbrace$.

For $n \! \in \! \mathbb{N}$ and $k \! \in \! \lbrace 1,2,\dotsc,K \rbrace$ 
such that $\alpha_{p_{\mathfrak{s}}} \! := \! \alpha_{k} \! \neq \! \infty$, 
the parametrix for the {\rm RHP} $(\mathfrak{M}(z),\mathfrak{v}(z),
\tilde{\Sigma})$, for $z \! \in \! \tilde{\mathbb{U}}_{\tilde{\delta}_{
\tilde{a}_{j}}}$, $j \! = \! 1,2,\dotsc,N \! + \! 1$, is the solution of the 
following {\rm RHPs} for $\tilde{\mathcal{X}}^{\tilde{a}} \colon \tilde{
\mathbb{U}}_{\tilde{\delta}_{\tilde{a}_{j}}} \setminus \tilde{\Sigma}_{\tilde{a}_{j}} 
\! \to \! \operatorname{SL}_{2}(\mathbb{C})$, $j \! = \! 1,2,\dotsc,N \! + \! 1$, 
where $\tilde{\Sigma}_{\tilde{a}_{j}} \! := \! (\tilde{\Phi}_{\tilde{a}_{j}})^{-1}
(\gamma_{\tilde{a}_{j}}^{\ast})$, with $(\tilde{\Phi}_{\tilde{a}_{j}})^{-1}$ 
denoting the inverse mapping, and $\gamma_{\tilde{a}_{j}}^{\ast} \! := 
\! \cup_{m=1}^{4} \gamma_{\tilde{a}_{j}}^{\ast,m}$$:$ {\rm (i)} $\tilde{
\mathcal{X}}^{\tilde{a}}(z)$ is analytic for $z \! \in \! \tilde{\mathbb{U}}_{
\tilde{\delta}_{\tilde{a}_{j}}} \setminus \tilde{\Sigma}_{\tilde{a}_{j}}$, $j 
\! = \! 1,2,\dotsc,N \! + \! 1$$;$ {\rm (ii)} the boundary values $\tilde{
\mathcal{X}}^{\tilde{a}}_{\pm}(z) \! := \! \lim_{\underset{z^{\prime} \, \in 
\, \pm \, \mathrm{side} \, \mathrm{of} \, \tilde{\Sigma}_{\tilde{a}_{j}}}{
z^{\prime} \to z \in \tilde{\Sigma}_{\tilde{a}_{j}}}} \tilde{\mathcal{X}}^{
\tilde{a}}(z^{\prime})$ 
satisfy the jump condition
\begin{equation*}
\tilde{\mathcal{X}}^{\tilde{a}}_{+}(z) \! = \! \tilde{\mathcal{X}}^{\tilde{a}}_{-}
(z) \mathfrak{v}(z), \quad z \! \in \! \tilde{\mathbb{U}}_{\tilde{\delta}_{
\tilde{a}_{j}}} \cap \tilde{\Sigma}_{\tilde{a}_{j}}, \quad j \! = \! 1,2,\dotsc,
N \! + \! 1;
\end{equation*}
and {\rm (iii)} uniformly for $z \! \in \! \partial \tilde{\mathbb{U}}_{
\tilde{\delta}_{\tilde{a}_{j}}} \! := \! \lbrace \mathstrut z \! \in \! 
\mathbb{C}; \, \lvert z \! - \! \tilde{a}_{j} \rvert \! = \! \tilde{\delta}_{
\tilde{a}_{j}} \rbrace$, $j \! = \! 1,2,\dotsc,N \! + \! 1$,
\begin{equation*}
\mathfrak{m}(z)(\tilde{\mathcal{X}}^{\tilde{a}}(z))^{-1} \underset{\underset{
z_{o}=1+o(1)}{\mathscr{N},n \to \infty}}{=} \mathrm{I} \! + \! o(1).
\end{equation*}
For $n \! \in \! \mathbb{N}$ and $k \! \in \! \lbrace 1,2,\dotsc,K \rbrace$ 
such that $\alpha_{p_{\mathfrak{s}}} \! := \! \alpha_{k} \! \neq \! \infty$, 
the solutions of the {\rm RHPs} $(\tilde{\mathcal{X}}^{\tilde{a}}(z),\mathfrak{v}
(z),\tilde{\mathbb{U}}_{\tilde{\delta}_{\tilde{a}_{j}}} \cap \tilde{\Sigma}_{
\tilde{a}_{j}})$, $j \! = \! 1,2,\dotsc,N \! + \! 1$, are: {\rm \pmb{(1)}} for 
$z \! \in \! \tilde{\Omega}_{\tilde{a}_{j}}^{1} \! := \! \tilde{\mathbb{U}}_{
\tilde{\delta}_{\tilde{a}_{j}}} \cap (\tilde{\Phi}_{\tilde{a}_{j}})^{-1}(\tilde{
\Omega}_{\tilde{a}_{j}}^{\ast,1})$, $j \! = \! 1,2,\dotsc,N \! + \! 1$,
\begin{equation*}
\tilde{\mathcal{X}}^{\tilde{a}}(z) \! = \! \sqrt{\smash[b]{\pi}} \, \me^{-
\frac{\mi \pi}{3}} \mathfrak{m}(z) \me^{\frac{\mi}{2}((n-1)K+k) 
\tilde{\mho}_{j} \operatorname{ad}(\sigma_{3})} 
\begin{pmatrix}
\mi & -\mi \\
1 & 1
\end{pmatrix}(\tilde{\Phi}_{\tilde{a}_{j}}(z))^{\frac{1}{4} \sigma_{3}} 
\Psi_{1}(\tilde{\Phi}_{\tilde{a}_{j}}(z)) \me^{\frac{1}{2}((n-1)K+k) 
\tilde{\xi}_{\tilde{a}_{j}}(z) \sigma_{3}},
\end{equation*}
where $\mathfrak{m}(z)$ is given in item~{\rm \pmb{(2)}} of 
Lemma~\ref{lem4.5}$;$ {\rm \pmb{(2)}} for $z \! \in \! \tilde{
\Omega}_{\tilde{a}_{j}}^{2} \! := \! \tilde{\mathbb{U}}_{\tilde{\delta}_{
\tilde{a}_{j}}} \cap (\tilde{\Phi}_{\tilde{a}_{j}})^{-1}(\tilde{\Omega}_{
\tilde{a}_{j}}^{\ast,2})$, $j \! = \! 1,2,\dotsc,N \! + \! 1$,
\begin{equation*}
\tilde{\mathcal{X}}^{\tilde{a}}(z) \! = \! \sqrt{\smash[b]{\pi}} \, \me^{-
\frac{\mi \pi}{3}} \mathfrak{m}(z) \me^{\frac{\mi}{2}((n-1)K+k) 
\tilde{\mho}_{j} \operatorname{ad}(\sigma_{3})} 
\begin{pmatrix}
\mi & -\mi \\
1 & 1
\end{pmatrix}(\tilde{\Phi}_{\tilde{a}_{j}}(z))^{\frac{1}{4} \sigma_{3}} \Psi_{2}
(\tilde{\Phi}_{\tilde{a}_{j}}(z)) \me^{\frac{1}{2}((n-1)K+k) \tilde{\xi}_{
\tilde{a}_{j}}(z) \sigma_{3}};
\end{equation*}
{\rm \pmb{(3)}} for $z \! \in \! \tilde{\Omega}_{\tilde{a}_{j}}^{3} \! := 
\! \tilde{\mathbb{U}}_{\tilde{\delta}_{\tilde{a}_{j}}} \cap (\tilde{\Phi}_{
\tilde{a}_{j}})^{-1}(\tilde{\Omega}_{\tilde{a}_{j}}^{\ast,3})$, $j \! = \! 1,2,
\dotsc,N \! + \! 1$,
\begin{equation*}
\tilde{\mathcal{X}}^{\tilde{a}}(z) \! = \! \sqrt{\smash[b]{\pi}} \, \me^{-
\frac{\mi \pi}{3}} \mathfrak{m}(z) \me^{-\frac{\mi}{2}((n-1)K+k) 
\tilde{\mho}_{j} \operatorname{ad}(\sigma_{3})} 
\begin{pmatrix}
\mi & -\mi \\
1 & 1
\end{pmatrix} 
(\tilde{\Phi}_{\tilde{a}_{j}}(z))^{\frac{1}{4} \sigma_{3}} \Psi_{3}(\tilde{
\Phi}_{\tilde{a}_{j}}(z)) \me^{\frac{1}{2}((n-1)K+k) \tilde{\xi}_{\tilde{a}_{j}}
(z) \sigma_{3}};
\end{equation*}
and {\rm \pmb{(4)}} for $z \! \in \! \tilde{\Omega}_{\tilde{a}_{j}}^{4} \! 
:= \! \tilde{\mathbb{U}}_{\tilde{\delta}_{\tilde{a}_{j}}} \cap (\tilde{\Phi}_{
\tilde{a}_{j}})^{-1}(\tilde{\Omega}_{\tilde{a}_{j}}^{\ast,4})$, $j \! = \! 1,2,
\dotsc,N \! + \! 1$,
\begin{equation*}
\tilde{\mathcal{X}}^{\tilde{a}}(z) \! = \! \sqrt{\smash[b]{\pi}} \, \me^{-
\frac{\mi \pi}{3}} \mathfrak{m}(z) \me^{-\frac{\mi}{2}((n-1)K+k) 
\tilde{\mho}_{j} \operatorname{ad}(\sigma_{3})} 
\begin{pmatrix}
\mi & -\mi \\
1 & 1
\end{pmatrix}(\tilde{\Phi}_{\tilde{a}_{j}}(z))^{\frac{1}{4} \sigma_{3}} \Psi_{4}
(\tilde{\Phi}_{\tilde{a}_{j}}(z)) \me^{\frac{1}{2}((n-1)K+k) \tilde{\xi}_{
\tilde{a}_{j}}(z) \sigma_{3}}.
\end{equation*}
\end{ccccc}

\emph{Proof of Lemma~\ref{lem4.9}.} For $n \! \in \! \mathbb{N}$ 
and $k \! \in \! \lbrace 1,2,\dotsc,K \rbrace$ such that $\alpha_{p_{
\mathfrak{s}}} \! := \! \alpha_{k} \! \neq \! \infty$, let $(\mathfrak{M}
(z),\mathfrak{v}(z),\tilde{\Sigma})$ be the associated RHP stated in 
item~\pmb{(2)} of Lemma~\ref{lem4.2}, and define $\tilde{\mathbb{U}}_{
\tilde{\delta}_{\tilde{a}_{j}}}$, $j \! = \! 1,2,\dotsc,N \! + \! 1$, as in 
Lemma~\ref{lem4.9}. Via the formula for $\mathfrak{v}(z)$ given in 
item~\pmb{(2)} of Lemma~\ref{lem4.2}, one shows, {}from the proof 
of case~\pmb{(A)} of Lemma~\ref{lem4.1}, that, for $n \! \in \! \mathbb{
N}$ and $k \! \in \! \lbrace 1,2,\dotsc,K \rbrace$ such that $\alpha_{
p_{\mathfrak{s}}} \! := \! \alpha_{k} \! \neq \! \infty$: (i) writing 
$2 \pi \mi \int_{z}^{\tilde{a}_{N+1}} \psi_{\widetilde{V}}^{f}(u) \, \md 
u \! = \! 2 \pi \mi (\int_{z}^{\tilde{a}_{j}} \! + \! \int_{\tilde{a}_{j}}^{
\tilde{b}_{j}} \! + \! \int_{\tilde{b}_{j}}^{\tilde{a}_{N+1}}) \psi_{\widetilde{
V}}^{f}(u) \, \md u$ and taking note of the expression for 
$\psi_{\widetilde{V}}^{f}(x)$ given in Equation~\eqref{eql3.7k}, one 
arrives at, upon considering the analytic extension of $2 \pi \mi \int_{z}^{
\tilde{a}_{N+1}} \psi_{\widetilde{V}}^{f}(u) \, \md u$ to $\mathbb{C} 
\setminus \mathbb{R}$ (in particular, to the oriented---open---skeletons 
$\tilde{\mathbb{U}}_{\tilde{\delta}_{\tilde{a}_{j}}} \cap (\tilde{J}_{j}^{\smallfrown} 
\cup \tilde{J}_{j}^{\smallsmile})$, $j \! = \! 1,2,\dotsc,N \! + \! 1)$, $2 \pi 
\mi \int_{z}^{\tilde{a}_{N+1}} \psi_{\widetilde{V}}^{f}(u) \, \md u \! = \! 
-\tilde{\xi}_{\tilde{a}_{j}}(z) \! + \! \mi \tilde{\mho}_{j}$, $j \! = \! 1,2,
\dotsc,N \! + \! 1$, where $\tilde{\xi}_{\tilde{a}_{j}}(z)$ is defined by 
Equation~\eqref{eqmainfin73}, and $\tilde{\mho}_{j}$ is defined in the 
corresponding item~(ii) of Remark~\ref{rem4.4}; and (ii) $g^{f}_{+}
(z) \! + \! g^{f}_{-}(z) \! - \! \hat{\mathscr{P}}^{+}_{0} \! - \! 
\hat{\mathscr{P}}^{-}_{0} \! - \! \widetilde{V}(z) \! - \! \tilde{\ell} \! = \! 
-\left(\tfrac{(n-1)K+k}{n} \right) \int_{\tilde{a}_{j}}^{z}(\tilde{R}(u))^{1/2} 
\tilde{h}_{\widetilde{V}}(u) \, \md u$ $(< \! 0)$, $z \! \in \! (\tilde{a}_{N
+1},+\infty) \cup \cup_{m=1}^{N}(\tilde{a}_{m},\tilde{b}_{m})$. Via the 
latter two formulae, denoting, for $\tilde{\mathbb{U}}_{\tilde{\delta}_{
\tilde{a}_{j}}} \! \ni \! z$, $j \! = \! 1,2,\dotsc,N \! + \! 1$, $\mathfrak{M}
(z)$ by $\tilde{\mathcal{X}}^{\tilde{a}}(z)$, and defining, for $n \! \in \! 
\mathbb{N}$ and $k \! \in \! \lbrace 1,2,\dotsc,K \rbrace$ such that 
$\alpha_{p_{\mathfrak{s}}} \! := \! \alpha_{k} \! \neq \! \infty$, 
\begin{equation*}
\tilde{\mathbb{P}}_{\tilde{a}_{j}}(z) \! := \!
\begin{cases}
\tilde{\mathcal{X}}^{\tilde{a}}(z) \me^{-\frac{1}{2}((n-1)K+k) \tilde{\xi}_{
\tilde{a}_{j}}(z) \sigma_{3}} \me^{\frac{\mi}{2}((n-1)K+k) \tilde{\mho}_{j} 
\sigma_{3}}, &\text{$z \! \in \! \mathbb{C}_{+} \cap \tilde{\mathbb{U}}_{
\tilde{\delta}_{\tilde{a}_{j}}}, \quad j \! = \! 1,2,\dotsc,N \! + \! 1$,} \\
\tilde{\mathcal{X}}^{\tilde{a}}(z) \me^{-\frac{1}{2}((n-1)K+k) \tilde{\xi}_{
\tilde{a}_{j}}(z) \sigma_{3}} \me^{-\frac{\mi}{2}((n-1)K+k) \tilde{\mho}_{j} 
\sigma_{3}}, &\text{$z \! \in \! \mathbb{C}_{-} \cap \tilde{\mathbb{U}}_{
\tilde{\delta}_{\tilde{a}_{j}}}, \quad j \! = \! 1,2,\dotsc,N \! + \! 1$,}
\end{cases}
\end{equation*}
one notes that $\tilde{\mathbb{P}}_{\tilde{a}_{j}} \colon \tilde{\mathbb{U}}_{
\tilde{\delta}_{\tilde{a}_{j}}} \setminus \tilde{J}_{\tilde{a}_{j}} \! \to \! 
\operatorname{GL}_{2}(\mathbb{C})$, where $\tilde{J}_{\tilde{a}_{j}} \! := \!
\tilde{J}_{j}^{\smallfrown} \cup \tilde{J}_{j}^{\smallsmile} \cup (\tilde{a}_{j} 
\! - \! \tilde{\delta}_{\tilde{a}_{j}},\tilde{a}_{j} \! + \! \tilde{\delta}_{\tilde{a}_{
j}})$, $j \! = \! 1,2,\dotsc,N \! + \! 1$, solves the RHP $(\tilde{\mathbb{P}}_{
\tilde{a}_{j}}(z),\tilde{\upsilon}_{\tilde{\mathbb{P}}_{\tilde{a}_{j}}}(z),
\tilde{J}_{\tilde{a}_{j}})$, with $z$-independent jump matrices
\begin{equation*}
\tilde{\upsilon}_{\tilde{\mathbb{P}}_{\tilde{a}_{j}}}(z) \! = \! 
\begin{cases}
\mathrm{I} \! + \! \sigma_{-}, &\text{$z \! \in \! \tilde{\mathbb{U}}_{\tilde{
\delta}_{\tilde{a}_{j}}} \cap (\tilde{J}_{j}^{\smallfrown} \cup \tilde{J}_{j}^{
\smallsmile}) \! = \! \tilde{\Sigma}^{1}_{\tilde{a}_{j}} \cup \tilde{\Sigma}^{
3}_{\tilde{a}_{j}}, \quad  j \! = \! 1,2,\dotsc,N \! + \! 1$,} \\
\mathrm{I} \! + \! \sigma_{+}, &\text{$z \! \in \! \tilde{\mathbb{U}}_{\tilde{
\delta}_{\tilde{a}_{j}}} \cap (\tilde{a}_{j},\tilde{a}_{j} \! + \! \tilde{\delta}_{
\tilde{a}_{j}}) \! = \! \tilde{\Sigma}^{4}_{\tilde{a}_{j}}, \quad  j \! = \! 1,2,
\dotsc,N \! + \! 1$,} \\
\mi \sigma_{2}, &\text{$z \! \in \! \tilde{\mathbb{U}}_{\tilde{\delta}_{
\tilde{a}_{j}}} \cap (\tilde{a}_{j} \! - \! \tilde{\delta}_{\tilde{a}_{j}},\tilde{a}_{j}) 
\! = \! \tilde{\Sigma}^{2}_{\tilde{a}_{j}}, \quad  j \! = \! 1,2,\dotsc,N \! + \! 1$,}
\end{cases}
\end{equation*}
subject, still, to the (uniform for $\partial \tilde{\mathbb{U}}_{\tilde{\delta}_{
\tilde{a}_{j}}} \! \ni \! z$, $j \! = \! 1,2,\dotsc,N \! + \! 1)$ asymptotic 
matching condition $\mathfrak{m}(z)(\tilde{\mathcal{X}}^{\tilde{a}}(z))^{-1} \! 
=_{\underset{z_{o}=1+o(1)}{\mathscr{N},n \to \infty}} \! \mathrm{I} \! + \! o(1)$, 
where $\mathfrak{m}(z)$ is given in item~\pmb{(2)} of Lemma~\ref{lem4.5}. 
For $n \! \in \! \mathbb{N}$ and $k \! \in \! \lbrace 1,2,\dotsc,K \rbrace$ 
such that $\alpha_{p_{\mathfrak{s}}} \! := \! \alpha_{k} \! \neq \! \infty$, let 
$\tilde{\Phi}_{\tilde{a}_{j}}(z)$ be defined by Equation~\eqref{eqmainfin85}: an 
analysis of the branch cuts shows that, for $z \! \in \! \tilde{\mathbb{U}}_{
\tilde{\delta}_{\tilde{a}_{j}}}$, $j \! = \! 1,2,\dotsc,N \! + \! 1$, $\tilde{\Phi}_{
\tilde{a}_{j}}(z)$ and $\tilde{\xi}_{\tilde{a}_{j}}(z)$ possess the properties stated 
in Lemma~\ref{lem4.9}; in particular, for $\tilde{\Phi}_{\tilde{a}_{j}} \colon 
\tilde{\mathbb{U}}_{\tilde{\delta}_{\tilde{a}_{j}}} \! \to \! \tilde{\mathbb{U}}^{
\ast}_{\tilde{\delta}_{\tilde{a}_{j}}}$, $j \! = \! 1,2,\dotsc,N \! + \! 1$, one 
shows that, for $\tilde{\mathbb{U}}_{\tilde{\delta}_{\tilde{a}_{j}}} \! \ni \! z$, 
$j \! = \! 1,2,\dotsc,N \! + \! 1$, $\tilde{\Phi}_{\tilde{a}_{j}}(z) \! = \! (z \! 
- \! \tilde{a}_{j})^{3/2} \tilde{G}_{\tilde{a}_{j}}(z)$, with $\tilde{G}_{\tilde{a}_{j}}
(z)$ as characterised in Lemma~\ref{lem4.9}, and $(\tilde{\Phi}_{\tilde{a}_{j}}
(z))^{\prime} \! \neq \! 0$, with, in particular, $(\tilde{\Phi}_{\tilde{a}_{j}}
(\tilde{a}_{j}))^{\prime} \! > \! 0$. One now chooses, for $n \! \in \! 
\mathbb{N}$ and $k \! \in \! \lbrace 1,2,\dotsc,K \rbrace$ such that 
$\alpha_{p_{\mathfrak{s}}} \! := \! \alpha_{k} \! \neq \! \infty$, $(0,1) \! 
\ni \! \tilde{\delta}_{\tilde{a}_{j}}$ and the oriented---open---skeletons 
$\tilde{J}_{\tilde{a}_{j}}$, $j \! = \! 1,2,\dotsc,N \! + \! 1$, in such a 
way that their images, under the bi-holomorphic, conformal, and 
orientation-preserving mappings $\tilde{\Phi}_{\tilde{a}_{j}}(z)$, are the 
union of the straight-line segments $\gamma_{\tilde{a}_{j}}^{\ast,m}$, 
$m \! = \! 1,2,3,4$. For $n \! \in \! \mathbb{N}$ and $k \! \in \! \lbrace 
1,2,\dotsc,K \rbrace$ such that $\alpha_{p_{\mathfrak{s}}} \! := \! 
\alpha_{k} \! \neq \! \infty$, let $\zeta \! := \! \tilde{\Phi}_{\tilde{a}_{j}}
(z)$, $j \! = \! 1,2,\dotsc,N \! + \! 1$, and consider $\tilde{\mathcal{X}}^{
\tilde{a}}(\tilde{\Phi}_{\tilde{a}_{j}}(z)) \! := \! \tilde{\Psi}(\zeta)$; recalling the 
properties of $\tilde{\Phi}_{\tilde{a}_{j}}(z)$ stated in Lemma~\ref{lem4.9}, 
one shows that, for $j \! = \! 1,2,\dotsc,N \! + \! 1$, $\tilde{\Psi} \colon 
\tilde{\Phi}_{\tilde{a}_{j}}(\tilde{\mathbb{U}}_{\tilde{\delta}_{\tilde{a}_{j}}}) 
\setminus \cup_{m=1}^{4} \gamma_{\tilde{a}_{j}}^{\ast,m} \! \to \! 
\operatorname{GL}_{2}(\mathbb{C})$ solves the RHPs $(\tilde{\Psi}
(\zeta),\tilde{\upsilon}_{\tilde{\Psi}}(\zeta),\cup_{m=1}^{4} \gamma_{
\tilde{a}_{j}}^{\ast,m})$, with $\zeta$-independent jump matrices
\begin{equation*}
\tilde{\upsilon}_{\tilde{\Psi}}(\zeta) \! = \! 
\begin{cases}
\mathrm{I} \! + \! \sigma_{-}, &\text{$\zeta \! \in \! \gamma_{\tilde{a}_{j}}^{
\ast,1} \cup \gamma_{\tilde{a}_{j}}^{\ast,3}, \quad j \! = \! 1,2,\dotsc,N \! + 
\! 1$,} \\
\mathrm{I} \! + \! \sigma_{+}, &\text{$\zeta \! \in \! \gamma_{\tilde{a}_{j}}^{
\ast,4}, \quad j \! = \! 1,2,\dotsc,N \! + \! 1$,} \\
\mi \sigma_{2}, &\text{$\zeta \! \in \! \gamma_{\tilde{a}_{j}}^{\ast,2}, \quad 
j \! = \! 1,2,\dotsc,N \! + \! 1$.}
\end{cases}
\end{equation*}
The solutions of these latter RHPs are well known; in fact, they are expressed in 
terms of Airy functions, and are given by (see, for example, \cite{a51,a54,a61,
pz22,a53,a67,a59})
\begin{equation*}
\tilde{\Psi}(\zeta) \! = \! 
\begin{cases}
\Psi_{1}(\zeta), &\text{$\zeta \! \in \! \tilde{\Omega}_{\tilde{a}_{j}}^{\ast,1}, 
\quad j \! = \! 1,2,\dotsc,N \! + \! 1$,} \\
\Psi_{2}(\zeta), &\text{$\zeta \! \in \! \tilde{\Omega}_{\tilde{a}_{j}}^{\ast,2}, 
\quad j \! = \! 1,2,\dotsc,N \! + \! 1$,} \\
\Psi_{3}(\zeta), &\text{$\zeta \! \in \! \tilde{\Omega}_{\tilde{a}_{j}}^{\ast,3}, 
\quad j \! = \! 1,2,\dotsc,N \! + \! 1$,} \\
\Psi_{4}(\zeta), &\text{$\zeta \! \in \! \tilde{\Omega}_{\tilde{a}_{j}}^{\ast,4}, 
\quad j \! = \! 1,2,\dotsc,N \! + \! 1$,}
\end{cases}
\end{equation*}
where $\Psi_{m}(\pmb{\cdot})$, $m \! = \! 1,2,3,4$, are defined in item~(i) 
of Remark~\ref{rem4.4}. Recalling that $\tilde{\Phi}_{\tilde{a}_{j}}(z)$, $j \! 
= \! 1,2,\dotsc,N \! + \! 1$, are bi-holomorphic, orientation-preserving 
conformal mappings, with $\tilde{\Phi}_{\tilde{a}_{j}}(\tilde{a}_{j}) \! = \! 0$ 
and $\tilde{\Phi}_{\tilde{a}_{j}} \colon \tilde{\mathbb{U}}_{\tilde{\delta}_{
\tilde{a}_{j}}} \! \cap \! \tilde{J}_{\tilde{a}_{j}} \! \to \! \tilde{\mathbb{U}}^{
\ast}_{\tilde{\delta}_{\tilde{a}_{j}}} \cap (\cup_{m=1}^{4} \gamma_{\tilde{a}_{
j}}^{\ast,m})$, one notes that, for any analytic maps $\tilde{E}_{\tilde{a}_{j}} 
\colon \tilde{\mathbb{U}}_{\tilde{\delta}_{\tilde{a}_{j}}} \! \to \! 
\operatorname{GL}_{2}(\mathbb{C})$, $j \! = \! 1,2,\dotsc,N \! + \! 1$, 
$\tilde{\mathbb{U}}_{\tilde{\delta}_{\tilde{a}_{j}}} \setminus \tilde{J}_{
\tilde{a}_{j}} \! \ni \! \zeta \! \mapsto \! \tilde{E}_{\tilde{a}_{j}}(\zeta) 
\tilde{\Psi}(\zeta)$ also solve the latter RHPs $(\tilde{\Psi}(\zeta),
\tilde{\upsilon}_{\tilde{\Psi}}(\zeta),\cup_{m=1}^{4} \gamma_{\tilde{a}_{j}}^{
\ast,m})$: one uses this `degree of freedom' of `multiplying on the left' 
by non-degenerate, analytic, matrix-valued functions in order to satisfy 
the remaining, uniform for $\partial \tilde{\mathbb{U}}_{\tilde{\delta}_{
\tilde{a}_{j}}} \! \ni \! z$, $j \! = \! 1,2,\dotsc,N \! + \! 1$, asymptotic 
matching condition $\mathfrak{m}(z)(\tilde{\mathcal{X}}^{\tilde{a}}(z))^{-1} 
\! =_{\underset{z_{o}=1+o(1)}{\mathscr{N},n \to \infty}} \! \mathrm{I} \! 
+ \! o(1)$ for the parametrix, where $\mathfrak{m}(z)$ is given in 
item~\pmb{(2)} of Lemma~\ref{lem4.5}. For $n \! \in \! \mathbb{N}$ 
and $k \! \in \! \lbrace 1,2,\dotsc,K \rbrace$ such that $\alpha_{
p_{\mathfrak{s}}} \! := \! \alpha_{k} \! \neq \! \infty$, consider, say, the 
regions (cf. Figure~\ref{foratil}) $\tilde{\Omega}_{\tilde{a}_{j}}^{1} \! := \! 
(\tilde{\Phi}_{\tilde{a}_{j}})^{-1}(\tilde{\Omega}_{\tilde{a}_{j}}^{\ast,1})$, 
$j \! = \! 1,2,\dotsc,N \! + \! 1$ (the remaining regions, $\tilde{
\Omega}_{\tilde{a}_{j}}^{m}$, $m \! = \! 2,3,4$, are similar): re-tracing 
the transformations above, one shows that, for $z \! \in \! \tilde{
\Omega}_{\tilde{a}_{j}}^{1}$, $j \! = \! 1,2,\dotsc,N \! + \! 1$, $\tilde{
\mathcal{X}}^{\tilde{a}}(z) \! = \! \tilde{E}_{\tilde{a}_{j}}(z) \tilde{\Psi}
((\tfrac{3}{4}((n \! - \! 1)K \! + \! k) \tilde{\xi}_{\tilde{a}_{j}}(z))^{2/3}) \exp 
(\tfrac{1}{2}((n \! - \! 1)K \! + \! k)(\tilde{\xi}_{\tilde{a}_{j}}(z) \! - \! \mi 
\tilde{\mho}_{j}) \sigma_{3})$, whence, via the expression above for 
$\tilde{\Psi}(\zeta)$, $\zeta \! \in \! \tilde{\Omega}_{\tilde{a}_{j}}^{
\ast,1}$, and the asymptotic expansions for $\operatorname{Ai}
(\pmb{\cdot})$ and $\operatorname{Ai}^{\prime}(\pmb{\cdot})$ given 
in Equation~\eqref{eqairy}, one shows that, for $z \! \in \! \partial 
\tilde{\mathbb{U}}_{\tilde{\delta}_{\tilde{a}_{j}}} \cap \partial \tilde{
\Omega}^{1}_{\tilde{a}_{j}}$, $j \! = \! 1,2,\dotsc,N \! + \! 1$,
\begin{equation*}
\tilde{\mathcal{X}}^{\tilde{a}}(z) \underset{\underset{z_{o}=1+o(1)}{
\mathscr{N},n \to \infty}}{=} \dfrac{1}{\sqrt{2 \pi}} \tilde{E}_{\tilde{a}_{j}}
(z) \left(\left(\dfrac{3}{4}((n \! - \! 1)K \! + \! k) \tilde{\xi}_{\tilde{a}_{j}}
(z) \right)^{2/3} \right)^{-\frac{1}{4} \sigma_{3}} 
\begin{pmatrix}
\me^{-\mi \pi/6} & \me^{\mi \pi/3} \\
-\me^{-\mi \pi/6} & -\me^{4 \pi \mi/3}
\end{pmatrix} \me^{-\frac{\mi}{2}((n-1)K+k) \tilde{\mho}_{j} 
\sigma_{3}}(\mathrm{I} \! + \! o(1));
\end{equation*}
demanding that, for $z \! \in \! \partial \tilde{\mathbb{U}}_{\tilde{\delta}_{
\tilde{a}_{j}}} \cap \partial \tilde{\Omega}_{\tilde{a}_{j}}^{1}$, $j \! = \! 1,2,
\dotsc,N \! + \! 1$, $\mathfrak{m}(z)(\tilde{\mathcal{X}}^{\tilde{a}}(z))^{-1} 
\! =_{\underset{z_{o}=1+o(1)}{\mathscr{N},n \to \infty}} \! \mathrm{I} \! + 
\! o(1)$, one arrives at, for $n \! \in \! \mathbb{N}$ and $k \! \in \! \lbrace 
1,2,\dotsc,K \rbrace$ such that $\alpha_{p_{\mathfrak{s}}} \! := \! \alpha_{k} 
\! \neq \! \infty$,
\begin{equation*}
\tilde{E}_{\tilde{a}_{j}}(z) \! = \! \dfrac{1}{\sqrt{2 \mi}} \mathfrak{m}(z) 
\me^{\frac{\mi}{2}((n-1)K+k) \tilde{\mho}_{j} \sigma_{3}} 
\begin{pmatrix}
\mi & -\mi \\
1 & 1
\end{pmatrix} \left(\left(\dfrac{3}{4}((n \! - \! 1)K \! + \! k) \tilde{\xi}_{
\tilde{a}_{j}}(z) \right)^{2/3} \right)^{\frac{1}{4} \sigma_{3}}, \quad j \! 
= \! 1,2,\dotsc,N \! + \! 1.
\end{equation*}
(Note that $\det (\tilde{E}_{\tilde{a}_{j}}(z)) \! = \! 1$.) One mimics the 
paradigm above for the remaining boundary skeletons $\partial \tilde{
\mathbb{U}}_{\tilde{\delta}_{\tilde{a}_{j}}} \cap \partial \tilde{\Omega}_{
\tilde{a}_{j}}^{m}$, $j \! = \! 1,2,\dotsc,N \! + \! 1$, $m \! = \! 2,3,4$, 
and shows that the same expression for $\tilde{E}_{\tilde{a}_{j}}(z)$ as 
given above is obtained; thus, for $n \! \in \! \mathbb{N}$ and $k \! \in 
\! \lbrace 1,2,\dotsc,K \rbrace$ such that $\alpha_{p_{\mathfrak{s}}} 
\! := \! \alpha_{k} \! \neq \! \infty$, with $\tilde{E}_{\tilde{a}_{j}}(z)$, 
$j \! = \! 1,2,\dotsc,N \! + \! 1$, as given above, one concludes that, 
uniformly for $z \! \in \! \partial \tilde{\mathbb{U}}_{\tilde{\delta}_{
\tilde{a}_{j}}}$, $\mathfrak{m}(z)(\tilde{\mathcal{X}}^{\tilde{a}}(z))^{-1} 
\! =_{\underset{z_{o}=1+o(1)}{\mathscr{N},n \to \infty}} \! \mathrm{I} 
\! + \! o(1)$. There remains, however, the question of unimodularity, 
since (depending on the regions)
\begin{equation*}
\det (\tilde{\mathcal{X}}^{\tilde{a}}(z)) \! = \! 
\left\vert 
\begin{smallmatrix}
\operatorname{Ai}(\tilde{\Phi}_{\tilde{a}_{j}}(z)) & \operatorname{Ai}
(\omega^{2} \tilde{\Phi}_{\tilde{a}_{j}}(z)) \\
\operatorname{Ai}^{\prime}(\tilde{\Phi}_{\tilde{a}_{j}}(z)) & \omega^{2} 
\operatorname{Ai}^{\prime}(\omega^{2} \tilde{\Phi}_{\tilde{a}_{j}}(z))
\end{smallmatrix}
\right\vert \qquad \text{or} \qquad \det (\tilde{\mathcal{X}}^{\tilde{a}}
(z)) \! = \! \left\vert 
\begin{smallmatrix}
\operatorname{Ai}(\tilde{\Phi}_{\tilde{a}_{j}}(z)) & -\omega^{2} 
\operatorname{Ai}(\omega \tilde{\Phi}_{\tilde{a}_{j}}(z)) \\
\operatorname{Ai}^{\prime}(\tilde{\Phi}_{\tilde{a}_{j}}(z)) & 
-\operatorname{Ai}^{\prime}(\omega \tilde{\Phi}_{\tilde{a}_{j}}(z))
\end{smallmatrix}
\right\vert:
\end{equation*}
multiplying $\tilde{\mathcal{X}}^{\tilde{a}}(z)$ `on the left' by a constant, 
$\tilde{c}_{\tilde{\mathcal{X}}}$, say, using the Wronskian relations 
(see, for example, Section~10.4 in \cite{abramsteg}) $\operatorname{W}
(\operatorname{Ai}(x),\operatorname{Ai}(\omega^{2}x)) \! = \! 
(2 \pi)^{-1} \exp (\mi \pi/6)$ or $\operatorname{W}(\operatorname{Ai}
(x),\operatorname{Ai}(\omega x)) \! = \! -(2 \pi)^{-1} \exp (-\mi \pi/6)$, 
and the linear dependence relation for Airy functions, $\operatorname{Ai}
(x) \! + \! \omega \operatorname{Ai}(\omega x) \! + \! \omega^{2} 
\operatorname{Ai}(\omega^{2}x) \! = \! 0$, one shows, for $n \! \in \! 
\mathbb{N}$ and $k \! \in \! \lbrace 1,2,\dotsc,K \rbrace$ such that 
$\alpha_{p_{\mathfrak{s}}} \! := \! \alpha_{k} \! \neq \! \infty$, upon 
imposing the condition $\det (\tilde{\mathcal{X}}^{\tilde{a}}(z)) \! = \! 1$, 
that $\tilde{c}_{\tilde{\mathcal{X}}} \! = \! (2 \pi)^{1/2} \exp (-\mi \pi/12)$. 
\hfill $\qed$

The following Lemma~\ref{lem4.10} states, for $n \! \in \! \mathbb{N}$ 
and $k \! \in \! \lbrace 1,2,\dotsc,K \rbrace$ such that $\alpha_{p_{
\mathfrak{s}}} \! := \! \alpha_{k} \! = \! \infty$ (resp., $\alpha_{p_{
\mathfrak{s}}} \! := \! \alpha_{k} \! \neq \! \infty)$, the extent to which 
the parametrices of Lemmata~\ref{lem4.6} and~\ref{lem4.7} (resp., 
Lemmata~\ref{lem4.8} and~\ref{lem4.9}) and model solution in 
item~\pmb{(1)} (resp., item~\pmb{(2)}) of Lemma~\ref{lem4.5} 
`approximate' the solution of the monic MPC ORF RHP stated in 
item~\pmb{(1)} (resp., item~\pmb{(2)}) of Lemma~\ref{lem4.2}.
\begin{ccccc} \label{lem4.10} 
For $n \! \in \! \mathbb{N}$ and $k \! \in \! \lbrace 1,2,\dotsc,K \rbrace$ 
such that $\alpha_{p_{\mathfrak{s}}} \! := \! \alpha_{k} \! = \! \infty$ 
(resp., $\alpha_{p_{\mathfrak{s}}} \! := \! \alpha_{k} \! \neq \! \infty)$, 
let $\mathfrak{M} \colon \mathbb{C} \setminus \hat{\Sigma} \! \to \! 
\operatorname{SL}_{2}(\mathbb{C})$ (resp., $\mathfrak{M} \colon 
\mathbb{C} \setminus \tilde{\Sigma} \! \to \! \operatorname{SL}_{2}
(\mathbb{C}))$ solve the monic {\rm MPC ORF RHP} $(\mathfrak{M}(z),
\mathfrak{v}(z),\hat{\Sigma})$ (resp., $(\mathfrak{M}(z),\mathfrak{v}(z),
\tilde{\Sigma}))$ stated in item~{\rm \pmb{(1)}} (resp., item~{\rm \pmb{(2)})} 
of Lemma~\ref{lem4.2}. For $n \! \in \! \mathbb{N}$ and $k \! \in \! 
\lbrace 1,2,\dotsc,K \rbrace$ such that $\alpha_{p_{\mathfrak{s}}} \! 
:= \! \alpha_{k} \! = \! \infty$, define
\begin{equation*}
\hat{S}_{p}(z) \! := \! 
\begin{cases}
\mathfrak{m}(z), &\text{$z \! \in \! \mathbb{C} \setminus 
\cup_{j=1}^{N+1}(\hat{\mathbb{U}}_{\hat{\delta}_{\hat{b}_{j-1}}} 
\cup \hat{\mathbb{U}}_{\hat{\delta}_{\hat{a}_{j}}})$,} \\
\hat{\mathcal{X}}(z), &\text{$z \! \in \! \cup_{j=1}^{N+1}
(\hat{\mathbb{U}}_{\hat{\delta}_{\hat{b}_{j-1}}} \cup \hat{
\mathbb{U}}_{\hat{\delta}_{\hat{a}_{j}}})$,}
\end{cases}
\end{equation*}
where
\begin{equation*}
\hat{\mathcal{X}}(z) \! = \! 
\begin{cases}
\hat{\mathcal{X}}^{\hat{b}}(z), &\text{$z \! \in \! \cup_{j=1}^{N+1} 
\hat{\mathbb{U}}_{\hat{\delta}_{\hat{b}_{j-1}}}$,} \\
\hat{\mathcal{X}}^{\hat{a}}(z), &\text{$z \! \in \! \cup_{j=1}^{N+1} 
\hat{\mathbb{U}}_{\hat{\delta}_{\hat{a}_{j}}}$,}
\end{cases}
\end{equation*}
$\mathfrak{m}(z)$ is given in item~{\rm \pmb{(1)}} of Lemma~\ref{lem4.5}, 
and, for $z \! \in \! \hat{\mathbb{U}}_{\hat{\delta}_{\hat{b}_{j-1}}}$ (resp., 
$z \! \in \! \hat{\mathbb{U}}_{\hat{\delta}_{\hat{a}_{j}}})$, $j \! = \! 1,2,
\dotsc,N \! + \! 1$, $\hat{\mathcal{X}}^{\hat{b}} \colon \hat{\mathbb{U}}_{
\hat{\delta}_{\hat{b}_{j-1}}} \setminus \hat{\Sigma}_{\hat{b}_{j-1}} \! \to 
\! \operatorname{SL}_{2}(\mathbb{C})$ (resp., $\hat{\mathcal{X}}^{\hat{a}} 
\colon \hat{\mathbb{U}}_{\hat{\delta}_{\hat{a}_{j}}} \setminus 
\hat{\Sigma}_{\hat{a}_{j}} \! \to \! \operatorname{SL}_{2}(\mathbb{C}))$ solves 
the {\rm RHPs} $(\hat{\mathcal{X}}^{\hat{b}}(z),\mathfrak{v}(z),\hat{\Sigma}_{
\hat{b}_{j-1}})$ (resp., $(\hat{\mathcal{X}}^{\hat{a}}(z),\mathfrak{v}(z),
\hat{\Sigma}_{\hat{a}_{j}}))$ stated in Lemma~\ref{lem4.6} (resp., 
Lemma~\ref{lem4.7}$)$, and, for $n \! \in \! \mathbb{N}$ and $k \! \in \! 
\lbrace 1,2,\dotsc,K \rbrace$ such that $\alpha_{p_{\mathfrak{s}}} \! := \! 
\alpha_{k} \! \neq \! \infty$, define
\begin{equation*}
\tilde{S}_{p}(z) \! := \! 
\begin{cases}
\mathfrak{m}(z), &\text{$z \! \in \! \mathbb{C} \setminus 
\cup_{j=1}^{N+1}(\tilde{\mathbb{U}}_{\tilde{\delta}_{\tilde{b}_{j-1}}} 
\cup \tilde{\mathbb{U}}_{\tilde{\delta}_{\tilde{a}_{j}}})$,} \\
\tilde{\mathcal{X}}(z), &\text{$z \! \in \! \cup_{j=1}^{N+1}
(\tilde{\mathbb{U}}_{\tilde{\delta}_{\tilde{b}_{j-1}}} \cup 
\tilde{\mathbb{U}}_{\tilde{\delta}_{\tilde{a}_{j}}})$,}
\end{cases}
\end{equation*}
where
\begin{equation*}
\tilde{\mathcal{X}}(z) \! = \! 
\begin{cases}
\tilde{\mathcal{X}}^{\tilde{b}}(z), &\text{$z \! \in \! \cup_{j=1}^{N+1} 
\tilde{\mathbb{U}}_{\tilde{\delta}_{\tilde{b}_{j-1}}}$,} \\
\tilde{\mathcal{X}}^{\tilde{a}}(z), &\text{$z \! \in \! \cup_{j=1}^{N+1} 
\tilde{\mathbb{U}}_{\tilde{\delta}_{\tilde{a}_{j}}}$,}
\end{cases}
\end{equation*}
$\mathfrak{m}(z)$ is given in item~{\rm \pmb{(2)}} of Lemma~\ref{lem4.5}, 
and, for $z \! \in \! \tilde{\mathbb{U}}_{\tilde{\delta}_{\tilde{b}_{j-1}}}$ 
(resp., $z \! \in \! \tilde{\mathbb{U}}_{\tilde{\delta}_{\tilde{a}_{j}}})$, 
$j \! = \! 1,2,\dotsc,N \! + \! 1$, $\tilde{\mathcal{X}}^{\tilde{b}} \colon 
\tilde{\mathbb{U}}_{\tilde{\delta}_{\tilde{b}_{j-1}}} \setminus \tilde{\Sigma}_{
\tilde{b}_{j-1}} \! \to \! \operatorname{SL}_{2}(\mathbb{C})$ (resp., 
$\tilde{\mathcal{X}}^{\tilde{a}} \colon \tilde{\mathbb{U}}_{\tilde{\delta}_{
\tilde{a}_{j}}} \setminus \tilde{\Sigma}_{\tilde{a}_{j}} \! \to \! \operatorname{SL}_{2}
(\mathbb{C}))$ solves the {\rm RHPs} $(\tilde{\mathcal{X}}^{\tilde{b}}
(z),\mathfrak{v}(z),\tilde{\Sigma}_{\tilde{b}_{j-1}})$ (resp., $(\tilde{
\mathcal{X}}^{\tilde{a}}(z),\mathfrak{v}(z),\tilde{\Sigma}_{\tilde{a}_{j}}))$ 
stated in Lemma~\ref{lem4.8} (resp., Lemma~\ref{lem4.9}$)$. For $n \! 
\in \! \mathbb{N}$ and $k \! \in \! \lbrace 1,2,\dotsc,K \rbrace$ such 
that $\alpha_{p_{\mathfrak{s}}} \! := \! \alpha_{k} \! = \! \infty$, set
\begin{equation*}
\hat{\mathcal{R}}(z) \! := \! \mathfrak{M}(z)(\hat{S}_{p}(z))^{-1},
\end{equation*}
and define the augmented and oriented contour $\hat{\Sigma}_{\hat{
\mathcal{R}}} \! := \! \hat{\Sigma} \cup \cup_{j=1}^{N+1}(\partial 
\hat{\mathbb{U}}_{\hat{\delta}_{\hat{b}_{j-1}}} \cup \partial \hat{
\mathbb{U}}_{\hat{\delta}_{\hat{a}_{j}}})$ as in Figure~\ref{figforabhat}, 
and, for $n \! \in \! \mathbb{N}$ and $k \! \in \! \lbrace 1,2,\dotsc,K 
\rbrace$ such that $\alpha_{p_{\mathfrak{s}}} \! := \! \alpha_{k} \! 
\neq \! \infty$, set
\begin{equation*}
\tilde{\mathcal{R}}(z) \! := \! \mathfrak{M}(z)(\tilde{S}_{p}(z))^{-1},
\end{equation*}
and define the augmented and oriented contour $\tilde{\Sigma}_{\tilde{
\mathcal{R}}} \! := \! \tilde{\Sigma} \cup \cup_{j=1}^{N+1}(\partial 
\tilde{\mathbb{U}}_{\tilde{\delta}_{\tilde{b}_{j-1}}} \cup \partial \tilde{
\mathbb{U}}_{\tilde{\delta}_{\tilde{a}_{j}}})$ as in Figure~\ref{figforabtil}. 
Then: {\rm \pmb{(1)}} for $n \! \in \! \mathbb{N}$ and $k \! \in \! 
\lbrace 1,2,\dotsc,K \rbrace$ such that $\alpha_{p_{\mathfrak{s}}} \! 
:= \! \alpha_{k} \! = \! \infty$, $\hat{\mathcal{R}} \colon \mathbb{C} 
\setminus \hat{\Sigma}_{\hat{\mathcal{R}}} \! \to \! \operatorname{SL}_{2}
(\mathbb{C})$ solves the following {\rm RHP:} {\rm (i)} $\hat{\mathcal{R}}
(z)$ is analytic for $z \! \in \! \mathbb{C} \setminus \hat{\Sigma}_{
\hat{\mathcal{R}}}$$;$ {\rm (ii)} the boundary values $\hat{\mathcal{R}}_{
\pm}(z) \! := \! \lim_{\underset{z^{\prime} \! \in \, \pm \, \mathrm{side} 
\, \mathrm{of} \, \hat{\Sigma}_{\hat{\mathcal{R}}}}{z^{\prime} \to z \in 
\hat{\Sigma}_{\hat{\mathcal{R}}}}} \hat{\mathcal{R}}(z^{\prime})$ satisfy 
the jump condition $\hat{\mathcal{R}}_{+}(z) \! = \! \hat{\mathcal{R}}_{-}
(z) \hat{v}_{\hat{\mathcal{R}}}(z)$ $\text{a.e.}$ $z \! \in \! \hat{\Sigma}_{
\hat{\mathcal{R}}}$, where
\begin{equation*}
\hat{v}_{\hat{\mathcal{R}}}(z) \! = \! 
\begin{cases}
\mathrm{I} \! + \! \me^{n(g^{\infty}_{+}(z)+g^{\infty}_{-}(z)-2 
\tilde{\mathscr{P}}_{0}-\widetilde{V}(z)-\hat{\ell})} \mathfrak{m}
(z) \sigma_{+}(\mathfrak{m}(z))^{-1}, &\text{$z \! \in \! 
\hat{\Sigma}_{p}^{1}$,} \\
\mathrm{I} \! + \! \me^{-\mi ((n-1)K+k) \hat{\Omega}_{j}} 
\me^{n(g^{\infty}_{+}(z)+g^{\infty}_{-}(z)-2 \tilde{\mathscr{P}}_{0}
-\widetilde{V}(z)-\hat{\ell})} \mathfrak{m}_{-}(z) \sigma_{+}
(\mathfrak{m}_{-}(z))^{-1},  &\text{$z \! \in \! \hat{\Sigma}^{2}_{p,j}, 
\quad j \! = \! 1,2,\dotsc,N$,} \\
\mathrm{I} \! + \! \me^{-2 \pi \mi ((n-1)K+k) \int_{z}^{\hat{a}_{N+1}} 
\psi_{\widetilde{V}}^{\infty}(\xi) \, \md \xi} \mathfrak{m}(z) \sigma_{-}
(\mathfrak{m}(z))^{-1}, &\text{$z \! \in \! \hat{\Sigma}^{3}_{p,j}, 
\quad j \! = \! 1,2,\dotsc,N \! + \! 1$,} \\
\mathrm{I} \! + \! \me^{2 \pi \mi ((n-1)K+k) \int_{z}^{\hat{a}_{N+1}} 
\psi_{\widetilde{V}}^{\infty}(\xi) \, \md \xi} \mathfrak{m}(z) \sigma_{-}
(\mathfrak{m}(z))^{-1}, &\text{$z \! \in \! \hat{\Sigma}^{4}_{p,j}, 
\quad j \! = \! 1,2,\dotsc,N \! + \! 1$,} \\
\hat{\mathcal{X}}(z)(\mathfrak{m}(z))^{-1}, &\text{$z \! \in \! 
\hat{\Sigma}^{5}_{p,j}, \quad j \! = \! 1,2,\dotsc,N \! + \! 1$,} \\
\mathrm{I}, &\text{$z \! \in \! \hat{\Sigma}_{\hat{\mathcal{R}}} 
\setminus \cup_{m=1}^{5} \hat{\Sigma}_{p}^{m}$,}
\end{cases}
\end{equation*}
with
\begin{gather*}
\hat{\Sigma}_{p}^{1} \! := \! (-\infty,\hat{b}_{0} \! - \! \hat{\delta}_{\hat{b}_{0}}) 
\! \cup \! (\hat{a}_{N+1} \! + \! \hat{\delta}_{\hat{a}_{N+1}},+\infty), \qquad 
\hat{\Sigma}^{2}_{p} \! := \! \cup_{j=1}^{N} \hat{\Sigma}^{2}_{p,j} \! := \! 
\cup_{j=1}^{N}(\hat{a}_{j} \! + \! \hat{\delta}_{\hat{a}_{j}},\hat{b}_{j} \! - \! 
\hat{\delta}_{\hat{b}_{j}}), \\
\hat{\Sigma}^{3}_{p} \! := \! \cup_{j=1}^{N+1} \hat{\Sigma}^{3}_{p,j} 
\! := \! \cup_{j=1}^{N+1}(\hat{J}_{j}^{\smallfrown} \setminus 
(\hat{J}_{j}^{\smallfrown} \cap (\hat{\mathbb{U}}_{\hat{\delta}_{\hat{b}_{j-1}}} 
\cup \hat{\mathbb{U}}_{\hat{\delta}_{\hat{a}_{j}}}))), \qquad \hat{\Sigma}^{
4}_{p} \! := \! \cup_{j=1}^{N+1} \hat{\Sigma}^{4}_{p,j} \! := \! \cup_{j=
1}^{N+1}(\hat{J}_{j}^{\smallsmile} \setminus (\hat{J}_{j}^{\smallsmile} \cap 
(\hat{\mathbb{U}}_{\hat{\delta}_{\hat{b}_{j-1}}} \cup \hat{\mathbb{U}}_{
\hat{\delta}_{\hat{a}_{j}}}))), \\
\hat{\Sigma}^{5}_{p} \! := \! \cup_{j=1}^{N+1} \hat{\Sigma}^{5}_{p,j} \! := 
\! \cup_{j=1}^{N+1}(\partial \hat{\mathbb{U}}_{\hat{\delta}_{\hat{b}_{j-1}}} 
\cup \partial \hat{\mathbb{U}}_{\hat{\delta}_{\hat{a}_{j}}});
\end{gather*}
and {\rm (iii)}
\begin{equation*}
\hat{\mathcal{R}}(z) \underset{\overline{\mathbb{C}} \setminus 
\hat{\Sigma}_{\hat{\mathcal{R}}} \ni z \to \alpha_{k}}{=} \mathrm{I} 
\! + \! \mathcal{O}(z^{-1}), \quad \quad \hat{\mathcal{R}}(z) 
\underset{\mathbb{C} \setminus \hat{\Sigma}_{\hat{\mathcal{R}}} 
\ni z \to \alpha_{p_{q}}}{=} \mathcal{O}(1), \quad q \! = \! 
1,2,\dotsc,\mathfrak{s} \! - \! 1;
\end{equation*}
and {\rm \pmb{(2)}} for $n \! \in \! \mathbb{N}$ and $k \! \in \! \lbrace 
1,2,\dotsc,K \rbrace$ such that $\alpha_{p_{\mathfrak{s}}} \! := \! 
\alpha_{k} \! \neq \! \infty$, $\tilde{\mathcal{R}} \colon \mathbb{C} 
\setminus \tilde{\Sigma}_{\tilde{\mathcal{R}}} \! \to \! 
\operatorname{SL}_{2}(\mathbb{C})$ solves the following {\rm RHP:} 
{\rm (i)} $\tilde{\mathcal{R}}(z)$ is analytic for $z \! \in \! \mathbb{C} 
\setminus \tilde{\Sigma}_{\tilde{\mathcal{R}}}$$;$ {\rm (ii)} the boundary 
values $\tilde{\mathcal{R}}_{\pm}(z) \! := \! \lim_{\underset{z^{\prime} 
\! \in \, \pm \, \mathrm{side} \, \mathrm{of} \, \tilde{\Sigma}_{\tilde{
\mathcal{R}}}}{z^{\prime} \to z \in \tilde{\Sigma}_{\tilde{\mathcal{R}}}}} 
\tilde{\mathcal{R}}(z^{\prime})$ satisfy the jump condition $\tilde{
\mathcal{R}}_{+}(z) \! = \! \tilde{\mathcal{R}}_{-}(z) \tilde{v}_{\tilde{
\mathcal{R}}}(z)$ $\text{a.e.}$ $z \! \in \! \tilde{\Sigma}_{\tilde{\mathcal{R}}}$, 
where
\begin{equation*}
\tilde{v}_{\tilde{\mathcal{R}}}(z) \! = \! 
\begin{cases}
\mathrm{I} \! + \! \me^{n(g^{f}_{+}(z)+g^{f}_{-}(z)-\hat{\mathscr{P}}_{0}^{+}
-\hat{\mathscr{P}}_{0}^{-}-\widetilde{V}(z)-\tilde{\ell})} \mathfrak{m}(z) 
\sigma_{+}(\mathfrak{m}(z))^{-1}, &\text{$z \! \in \! 
\tilde{\Sigma}_{p}^{1}$,} \\
\mathrm{I} \! + \! \me^{-\mi ((n-1)K+k) \tilde{\Omega}_{j}} \me^{n(g^{f}_{+}
(z)+g^{f}_{-}(z)-\hat{\mathscr{P}}_{0}^{+}-\hat{\mathscr{P}}_{0}^{-}-
\widetilde{V}(z)-\tilde{\ell})} \mathfrak{m}_{-}(z) \sigma_{+}(\mathfrak{m}_{-}
(z))^{-1},  &\text{$z \! \in \! \tilde{\Sigma}^{2}_{p,j}, \quad j \! = \! 1,2,
\dotsc,N$,} \\
\mathrm{I} \! + \! \me^{-2 \pi \mi ((n-1)K+k) \int_{z}^{\tilde{a}_{N+1}} 
\psi_{\widetilde{V}}^{f}(\xi) \, \md \xi} \mathfrak{m}(z) \sigma_{-}
(\mathfrak{m}(z))^{-1}, &\text{$z \! \in \! \tilde{\Sigma}^{3}_{p,j}, \quad 
j \! = \! 1,2,\dotsc,N \! + \! 1$,} \\
\mathrm{I} \! + \! \me^{2 \pi \mi ((n-1)K+k) \int_{z}^{\tilde{a}_{N+1}} 
\psi_{\widetilde{V}}^{f}(\xi) \, \md \xi} \mathfrak{m}(z) \sigma_{-}
(\mathfrak{m}(z))^{-1}, &\text{$z \! \in \! \tilde{\Sigma}^{4}_{p,j}, \quad 
j \! = \! 1,2,\dotsc,N \! + \! 1$,} \\
\tilde{\mathcal{X}}(z)(\mathfrak{m}(z))^{-1}, &\text{$z \! \in \! 
\tilde{\Sigma}^{5}_{p,j}, \quad j \! = \! 1,2,\dotsc,N \! + \! 1$,} \\
\mathrm{I}, &\text{$z \! \in \! \tilde{\Sigma}_{\tilde{\mathcal{R}}} 
\setminus \cup_{m=1}^{5} \tilde{\Sigma}_{p}^{m}$,}
\end{cases}
\end{equation*}
with
\begin{gather*}
\tilde{\Sigma}_{p}^{1} \! := \! (-\infty,\tilde{b}_{0} \! - \! \tilde{\delta}_{
\tilde{b}_{0}}) \! \cup \! (\tilde{a}_{N+1} \! + \! \tilde{\delta}_{\tilde{a}_{N
+1}},+\infty), \qquad \tilde{\Sigma}^{2}_{p} \! := \! \cup_{j=1}^{N} \tilde{
\Sigma}^{2}_{p,j} \! := \! \cup_{j=1}^{N}(\tilde{a}_{j} \! + \! \tilde{\delta}_{
\tilde{a}_{j}},\tilde{b}_{j} \! - \! \tilde{\delta}_{\tilde{b}_{j}}), \\
\tilde{\Sigma}^{3}_{p} \! := \! \cup_{j=1}^{N+1} \tilde{\Sigma}^{3}_{p,j} 
\! := \! \cup_{j=1}^{N+1}(\tilde{J}_{j}^{\smallfrown} \setminus 
(\tilde{J}_{j}^{\smallfrown} \cap (\tilde{\mathbb{U}}_{\tilde{\delta}_{
\tilde{b}_{j-1}}} \cup \tilde{\mathbb{U}}_{\tilde{\delta}_{\tilde{a}_{j}}}))), 
\qquad \tilde{\Sigma}^{4}_{p} \! := \! \cup_{j=1}^{N+1} \tilde{\Sigma}^{
4}_{p,j} \! := \! \cup_{j=1}^{N+1}(\tilde{J}_{j}^{\smallsmile} \setminus 
(\tilde{J}_{j}^{\smallsmile} \cap (\tilde{\mathbb{U}}_{\tilde{\delta}_{
\tilde{b}_{j-1}}} \cup \tilde{\mathbb{U}}_{\tilde{\delta}_{\tilde{a}_{j}}}))), \\
\tilde{\Sigma}^{5}_{p} \! := \! \cup_{j=1}^{N+1} \tilde{\Sigma}^{5}_{p,j} 
\! := \! \cup_{j=1}^{N+1}(\partial \tilde{\mathbb{U}}_{\tilde{\delta}_{
\tilde{b}_{j-1}}} \cup \partial \tilde{\mathbb{U}}_{\tilde{\delta}_{\tilde{a}_{j}}});
\end{gather*}
and {\rm (iii)}
\begin{gather*}
\tilde{\mathcal{R}}(z) \underset{\mathbb{C} \setminus 
\tilde{\Sigma}_{\tilde{\mathcal{R}}} \ni z \to \alpha_{k}}{=} \mathrm{I} \! + \! 
\mathcal{O}(z \! - \! \alpha_{k}), \qquad \qquad \tilde{\mathcal{R}}(z) 
\underset{\overline{\mathbb{C}} \setminus \tilde{\Sigma}_{\tilde{\mathcal{R}}} 
\ni z \to \alpha_{p_{\mathfrak{s}-1}} = \infty}{=} \mathcal{O}(1), \\
\tilde{\mathcal{R}}(z) \underset{\mathbb{C} \setminus \tilde{\Sigma}_{
\tilde{\mathcal{R}}} \ni z \to \alpha_{p_{q}}}{=} \mathcal{O}(1), 
\quad q \! = \! 1,2,\dotsc,\mathfrak{s} \! - \! 2.
\end{gather*}
\end{ccccc}
\begin{figure}[tbh]
\begin{center}

\end{center}
\vspace{-1.00cm}
\caption{The augmented and oriented contour $\hat{\Sigma}_{\hat{\mathcal{R}}} 
\! := \! \hat{\Sigma} \cup \cup_{j=1}^{N+1}(\partial \hat{\mathbb{U}}_{
\hat{\delta}_{\hat{b}_{j-1}}} \cup \partial \hat{\mathbb{U}}_{\hat{\delta}_{
\hat{a}_{j}}})$.}\label{figforabhat}
\end{figure}
\begin{figure}[tbh]
\begin{center}
%
\end{center}
\vspace{-1.00cm}
\caption{The augmented and oriented contour $\tilde{\Sigma}_{\tilde{
\mathcal{R}}} \! := \! \tilde{\Sigma} \cup \cup_{j=1}^{N+1}(\partial 
\tilde{\mathbb{U}}_{\tilde{\delta}_{\tilde{b}_{j-1}}} \cup \partial \tilde{
\mathbb{U}}_{\tilde{\delta}_{\tilde{a}_{j}}})$.}\label{figforabtil}
\end{figure}
\emph{Proof}. The proof of this Lemma~\ref{lem4.10} consists of two cases: 
(i) $n \! \in \! \mathbb{N}$ and $k \! \in \! \lbrace 1,2,\dotsc,K \rbrace$ 
such that $\alpha_{p_{\mathfrak{s}}} \! := \! \alpha_{k} \! = \! \infty$; and 
(ii) $n \! \in \! \mathbb{N}$ and $k \! \in \! \lbrace 1,2,\dotsc,K \rbrace$ 
such that $\alpha_{p_{\mathfrak{s}}} \! := \! \alpha_{k} \! \neq \! \infty$. 
Notwithstanding the fact that the scheme of the proof is, 
\emph{mutatis mutandis}, similar for both cases, without loss of generality, 
only the proof for case~(ii) is presented, whilst case~(i) is proved analogously.

For $n \! \in \! \mathbb{N}$ and $k \! \in \! \lbrace 1,2,\dotsc,K \rbrace$ 
such that $\alpha_{p_{\mathfrak{s}}} \! := \! \alpha_{k} \! \neq \! \infty$, 
define the augmented and oriented skeleton $\tilde{\Sigma}_{\tilde{
\mathcal{R}}}$ as in item~\pmb{(2)} of the lemma (cf. Figure~\ref{figforabtil}): 
via the definitions of $\tilde{S}_{p}(z)$ and $\tilde{\mathcal{R}}(z)$ given in 
item~\pmb{(2)} of the lemma, the associated RHP $(\tilde{\mathcal{R}}(z),
\tilde{v}_{\tilde{\mathcal{R}}}(z),\tilde{\Sigma}_{\tilde{\mathcal{R}}})$ stated 
in item~\pmb{(2)} of the lemma derives {}from the RHPs $(\mathfrak{M}
(z),\mathfrak{v}(z),\tilde{\Sigma})$ and $(\mathfrak{m}(z),\daleth (z),
\tilde{\mathscr{J}})$ stated in items~\pmb{(2)} of Lemmata~\ref{lem4.2} 
and~\ref{lem4.5}, respectively. \hfill $\qed$
\section{Asymptotics of Solutions of the Monic MPC ORF Families of RHPs 
and MPAs} 
\label{sek5} 
In this section, via the Beals-Coifman (BC) construction \cite{bealscoif}, 
for $n \! \in \! \mathbb{N}$ and $k \! \in \! \lbrace 1,2,\dotsc,K \rbrace$ 
such that $\alpha_{p_{\mathfrak{s}}} \! := \! \alpha_{k} \! = \! \infty$ 
(resp., $\alpha_{p_{\mathfrak{s}}} \! := \! \alpha_{k} \!  \neq \! \infty)$, 
the associated `small-norm' RHP $(\hat{\mathcal{R}}(z),\hat{v}_{\hat{
\mathcal{R}}}(z),\hat{\Sigma}_{\hat{\mathcal{R}}})$ (resp., $(\tilde{
\mathcal{R}}(z),\tilde{v}_{\tilde{\mathcal{R}}}(z),\tilde{\Sigma}_{\tilde{
\mathcal{R}}}))$ stated in item~\pmb{(1)} (resp., item~\pmb{(2)}) of 
Lemma~\ref{lem4.10} is solved asymptotically in the double-scaling limit 
$\mathscr{N},n \! \to \! \infty$ such that $z_{o} \! = \! 1 \! + \! o(1)$; 
and, by re-tracing the sequence of the families of RHP transformations 
$\hat{\mathcal{R}}(z)$ (cf. Lemma~\ref{lem4.10}, item~\pmb{(1)}) $\to 
\! \mathfrak{M}(z)$ (cf. Lemma~\ref{lem4.2}, item~\pmb{(1)}) $\to 
\! \mathcal{M}(z)$ (cf. Lemma~\ref{lem4.1}) $\to \! \mathscr{X}(z)$ 
(cf. Lemma~\ref{lem3.4}) (resp., $\tilde{\mathcal{R}}(z)$ (cf. 
Lemma~\ref{lem4.10}, item~\pmb{(2)}) $\to \! \mathfrak{M}(z)$ (cf. 
Lemma~\ref{lem4.2}, item~\pmb{(2)}) $\to \! \mathcal{M}(z)$ (cf. 
Lemma~\ref{lem4.1}) $\to \! \mathscr{X}(z)$ (cf. Lemma~\ref{lem3.4})) 
$\to \! \mathcal{X}(z)$ (cf. Lemmata~\ref{lem2.1} 
and~$\bm{\mathrm{RHP}_{\mathrm{MPC}}}$), the original monic MPC 
ORF families of RHPs $(\mathcal{X}(z),\upsilon (z),\overline{\mathbb{R}})$ 
are solved uniquely and asymptotically in the double-scaling limit 
$\mathscr{N},n \! \to \! \infty$ such that $z_{o} \! = \! 1 \! + \! o(1)$, 
which then yields the final families of asymptotic results (in the entire 
complex plane) for the associated monic MPC ORFs, $\pmb{\pi}^{n}_{k}
(z)$, the corresponding MPA error term, $\widehat{\pmb{\mathrm{E}}}_{
\tilde{\mu}}(z)$ (resp., $\widetilde{\pmb{\mathrm{E}}}_{\tilde{\mu}}(z))$, 
norming constants, $\mu_{n,\varkappa_{nk}}^{r}(n,k)$, $r \! \in \! \lbrace 
\infty,f \rbrace$, and MPC ORFs, $\phi^{n}_{k}(z)$, stated, respectively, 
in Theorems~\ref{maintheoforinf1}--\ref{maintheoforinf2} (resp., 
Theorems~\ref{maintheoforfin1}--\ref{maintheoforfin2}).
\begin{ccccc} \label{lem5.1} 
For $n \! \in \! \mathbb{N}$ and $k \! \in \! \lbrace 1,2,\dotsc,K 
\rbrace$ such that $\alpha_{p_{\mathfrak{s}}} \! := \! \alpha_{k} \! 
= \! \infty$ (resp., $\alpha_{p_{\mathfrak{s}}} \! := \! \alpha_{k} 
\! \neq \! \infty)$, let $\hat{\mathcal{R}} \colon \mathbb{C} 
\setminus \hat{\Sigma}_{\hat{\mathcal{R}}} \! \to \! \mathrm{SL}_{2}
(\mathbb{C})$ (resp., $\tilde{\mathcal{R}} \colon \mathbb{C} 
\setminus \tilde{\Sigma}_{\tilde{\mathcal{R}}} \! \to \! \mathrm{SL}_{2}
(\mathbb{C}))$ solve the {\rm RHP} stated in item~{\rm \pmb{(1)}} 
(resp., item~{\rm \pmb{(2)}}$)$ of Lemma~\ref{lem4.10}. Then: 
{\rm \pmb{(1)}} for $n \! \in \! \mathbb{N}$ and $k \! \in \! \lbrace 
1,2,\dotsc,K \rbrace$ such that $\alpha_{p_{\mathfrak{s}}} \! := \! 
\alpha_{k} \! = \! \infty$$:$
\begin{enumerate}
\item[$\boldsymbol{\mathrm{(1)}_{i}}$] for $z \! \in \! \hat{\Sigma}_{
p,j}^{2} \! := \! (\hat{a}_{j} \! + \! \hat{\delta}_{\hat{a}_{j}},\hat{b}_{j} 
\! - \! \hat{\delta}_{\hat{b}_{j}})$, $j \! = \! 1,2,\dotsc,N$,
\begin{equation} \label{eqhtvee1} 
\hat{v}_{\hat{\mathcal{R}}}(z) \underset{\underset{z_{o}=1+o(1)}{
\mathscr{N},n \to \infty}}{=} \mathrm{I} \! + \! \mathcal{O} \left(
\hat{\mathfrak{c}}_{\hat{\mathcal{R}},1}(j) \me^{-((n-1)K+k) 
\hat{\lambda}_{\hat{\mathcal{R}},1}(j)(z-\hat{a}_{j}) \chi_{
\hat{A}_{\hat{\mathcal{R}},1}(j)}(z)} \prod_{q \in \hat{Q}_{\hat{
\mathcal{R}},1}(j)} \lvert z \! - \! \alpha_{p_{q}} \rvert^{((n-1)K+k) 
\hat{\mathfrak{c}}_{\hat{\mathcal{R}},2}(j,q) \chi_{\mathscr{O}_{
\hat{\delta}_{\hat{\mathcal{R}},1}(j)}(\alpha_{p_{q}})}(z)} \right),
\end{equation}
where, for $j \! = \! 1,2,\dotsc,N$, $(\mathrm{M}_{2}(\mathbb{C}) 
\! \ni)$ $\hat{\mathfrak{c}}_{\hat{\mathcal{R}},1}(j) \! 
=_{\underset{z_{o}=1+o(1)}{\mathscr{N},n \to \infty}} \! 
\mathcal{O}(1)$, $\hat{\lambda}_{\hat{\mathcal{R}},1}(j) \! = \! 
\hat{\lambda}_{\hat{\mathcal{R}},1}(n,k,z_{o};j) \! =_{\underset{z_{o}
=1+o(1)}{\mathscr{N},n \to \infty}} \! \mathcal{O}(1)$ and $> \! 0$, 
$\hat{\mathfrak{c}}_{\hat{\mathcal{R}},2}(j,q) \! = \! \hat{\mathfrak{
c}}_{\hat{\mathcal{R}},2}(n,k,z_{o};j,q) \! =_{\underset{z_{o}=1+
o(1)}{\mathscr{N},n \to \infty}} \! \mathcal{O}(1)$ and $> \! 0$, 
and $\hat{A}_{\hat{\mathcal{R}},1}(j) \! = \! (\hat{a}_{j} \! + \! 
\hat{\delta}_{\hat{a}_{j}},\hat{b}_{j} \! - \! \hat{\delta}_{\hat{b}_{j}}) 
\setminus \cup_{q \in \hat{Q}_{\hat{\mathcal{R}},1}(j)} \mathscr{O}_{
\hat{\delta}_{\hat{\mathcal{R}},1}(j)}(\alpha_{p_{q}})$, 
with $\hat{Q}_{\hat{\mathcal{R}},1}(j) \! := \! \lbrace 
\mathstrut q^{\prime} \! \in \! \lbrace 1,2,\dotsc,
\mathfrak{s} \! - \! 1 \rbrace; \, \alpha_{p_{q^{\prime}}} \! \in \! 
(\hat{a}_{j} \! + \! \hat{\delta}_{\hat{a}_{j}},\hat{b}_{j} \! - \! \hat{
\delta}_{\hat{b}_{j}}) \rbrace$ and sufficiently small $\hat{\delta}_{
\hat{\mathcal{R}},1}(j) \! > \! 0$ chosen so that $\mathscr{O}_{
\hat{\delta}_{\hat{\mathcal{R}},1}(j)}(\alpha_{p_{q_{1}}}) \cap 
\mathscr{O}_{\hat{\delta}_{\hat{\mathcal{R}},1}(j)}(\alpha_{p_{
q_{2}}}) \! = \! \varnothing$ $\forall$ $q_{1} \! \neq \! q_{2} \! 
\in \! \hat{Q}_{\hat{\mathcal{R}},1}(j)$ and $\mathscr{O}_{
\hat{\delta}_{\hat{\mathcal{R}},1}(j)}(\alpha_{p_{q_{1}}}) \cap 
\lbrace \hat{a}_{j} \! + \! \hat{\delta}_{\hat{a}_{j}} \rbrace \! = \! 
\varnothing \! = \! \mathscr{O}_{\hat{\delta}_{\hat{\mathcal{R}},
1}(j)}(\alpha_{p_{q_{1}}}) \cap \lbrace \hat{b}_{j} \! - \! \hat{
\delta}_{\hat{b}_{j}} \rbrace$$;$\footnote{If, for $j \! = \! 1,2,
\dotsc,N$, $\hat{Q}_{\hat{\mathcal{R}},1}(j) \! = \! \varnothing$, 
that is, $\# \hat{Q}_{\hat{\mathcal{R}},1}(j) \! = \! 0$, then 
$\cup_{q \in \hat{Q}_{\hat{\mathcal{R}},1}(j)} \mathscr{O}_{\hat{
\delta}_{\hat{\mathcal{R}},1}(j)}(\alpha_{p_{q}}) \! := \! \varnothing$, 
in which case $\hat{A}_{\hat{\mathcal{R}},1}(j) \! = \! (\hat{a}_{j} \! + 
\! \hat{\delta}_{\hat{a}_{j}},\hat{b}_{j} \! - \! \hat{\delta}_{\hat{b}_{j}})$, 
and $\prod_{q \in \hat{Q}_{\hat{\mathcal{R}},1}(j)} \lvert z \! - 
\! \alpha_{p_{q}} \rvert^{((n-1)K+k) \hat{\mathfrak{c}}_{\hat{
\mathcal{R}},2}(j,q) \chi_{\mathscr{O}_{\hat{\delta}_{\hat{
\mathcal{R}},1}(j)}(\alpha_{p_{q}})}(z)} \! := \! 1$ (of course, 
$\cup_{q \in \hat{Q}_{\hat{\mathcal{R}},1}(j)} \mathscr{O}_{\hat{
\delta}_{\hat{\mathcal{R}},1}(j)}(\alpha_{p_{q}}) \subset (\hat{a}_{j} 
\! + \! \hat{\delta}_{\hat{a}_{j}},\hat{b}_{j} \! - \! \hat{\delta}_{
\hat{b}_{j}}))$; moreover, $\# \cup_{j=1}^{N} \hat{Q}_{\hat{
\mathcal{R}},1}(j) \! \leqslant \! \mathfrak{s} \! - \! 1$.}
\item[$\boldsymbol{\mathrm{(1)}_{ii}}$] for $z \! \in \! 
\hat{\Sigma}^{1}_{p} \! := \! (-\infty,\hat{b}_{0} \! - \! \hat{
\delta}_{\hat{b}_{0}}) \cup (\hat{a}_{N+1} \! + \! \hat{\delta}_{
\hat{a}_{N+1}},+\infty)$,
\begin{align} 
\hat{v}_{\hat{\mathcal{R}}}(z) \underset{\underset{z_{o}=1+o(1)}{
\mathscr{N},n \to \infty}}{=}& \, \mathrm{I} \! + \! \mathcal{O} 
\left(\hat{\mathfrak{c}}_{\hat{\mathcal{R}},2}(+) \me^{-((n-1)K+
k) \hat{\lambda}_{\hat{\mathcal{R}},2}(+)(z-\hat{a}_{N+1}) \chi_{
\hat{A}_{\hat{\mathcal{R}},2}(+)}(z)} \me^{-((n-1)K+k) \hat{
\lambda}_{\hat{\mathcal{R}},3}(+) \ln (\lvert z \rvert) \chi_{
\mathscr{O}_{\infty}(\alpha_{k})}(z)} \right. \nonumber \\
\times&\left. \, \prod_{q \in \hat{Q}_{\hat{\mathcal{R}},2}(+)} 
\lvert z \! - \! \alpha_{p_{q}} \rvert^{((n-1)K+k) \hat{\mathfrak{c}}_{
\hat{\mathcal{R}},3}(q,+) \chi_{\mathscr{O}_{\hat{\delta}_{\hat{
\mathcal{R}},2}(+)}(\alpha_{p_{q}})}(z)} \right), \quad z \! \in \! 
(\hat{a}_{N+1} \! + \! \hat{\delta}_{\hat{a}_{N+1}},+\infty), 
\label{eqhtvee2} \\
\hat{v}_{\hat{\mathcal{R}}}(z) \underset{\underset{z_{o}=1+o(1)}{
\mathscr{N},n \to \infty}}{=}& \, \mathrm{I} \! + \! \mathcal{O} 
\left(\hat{\mathfrak{c}}_{\hat{\mathcal{R}},2}(-) \me^{-((n-1)K+
k) \hat{\lambda}_{\hat{\mathcal{R}},2}(-) \lvert z-\hat{b}_{0} \rvert 
\chi_{\hat{A}_{\hat{\mathcal{R}},2}(-)}(z)} \me^{-((n-1)K+k) 
\hat{\lambda}_{\hat{\mathcal{R}},3}(-) \ln (\lvert z \rvert) 
\chi_{\mathscr{O}_{\infty}(\alpha_{k})}(z)} \right. \nonumber \\
\times&\left. \, \prod_{q \in \hat{Q}_{\hat{\mathcal{R}},2}(-)} 
\lvert z \! - \! \alpha_{p_{q}} \rvert^{((n-1)K+k) \hat{\mathfrak{c}}_{
\hat{\mathcal{R}},3}(q,-) \chi_{\mathscr{O}_{\hat{\delta}_{\hat{
\mathcal{R}},2}(-)}(\alpha_{p_{q}})}(z)} \right), \quad z \! \in \! 
(-\infty,\hat{b}_{0} \! - \! \hat{\delta}_{\hat{b}_{0}}), \label{eqhtvee3}
\end{align}
where $(\mathrm{M}_{2}(\mathbb{C}) \! \ni)$ $\hat{\mathfrak{c}}_{
\hat{\mathcal{R}},2}(\pm) \! =_{\underset{z_{o}=1+o(1)}{\mathscr{N},
n \to \infty}} \! \mathcal{O}(1)$, $\hat{\lambda}_{\hat{\mathcal{R}},
i_{1}}(\pm) \! = \! \hat{\lambda}_{\hat{\mathcal{R}},i_{1}}(n,k,z_{o};
\pm) \! =_{\underset{z_{o}=1+o(1)}{\mathscr{N},n \to \infty}} \! 
\mathcal{O}(1)$ and $> \! 0$, $i_{1} \! = \! 2,3$, $\hat{\mathfrak{c}}_{
\hat{\mathcal{R}},3}(q,\pm) \! = \! \hat{\mathfrak{c}}_{\hat{\mathcal{R}},
3}(n,k,z_{o};q,\pm) \! =_{\underset{z_{o}=1+o(1)}{\mathscr{N},n \to 
\infty}} \! \mathcal{O}(1)$ and $> \! 0$, $\hat{A}_{\hat{\mathcal{R}},2}
(+) \! = \! (\hat{a}_{N+1} \! + \! \hat{\delta}_{\hat{a}_{N+1}},+\infty) 
\setminus (\mathscr{O}_{\infty}(\alpha_{k}) \cup \cup_{q \in \hat{Q}_{
\hat{\mathcal{R}},2}(+)} \mathscr{O}_{\hat{\delta}_{\hat{\mathcal{R}},2}
(+)}(\alpha_{p_{q}}))$, $\hat{A}_{\hat{\mathcal{R}},2}(-) \! = \! (-\infty,
\hat{b}_{0} \! - \! \hat{\delta}_{\hat{b}_{0}}) \setminus (\mathscr{O}_{
\infty}(\alpha_{k}) \cup \cup_{q \in \hat{Q}_{\hat{\mathcal{R}},2}(-)}  
\mathscr{O}_{\hat{\delta}_{\hat{\mathcal{R}},2}(-)}(\alpha_{p_{q}}))$, 
and $\mathscr{O}_{\infty}(\alpha_{k}) \! = \! \lbrace \mathstrut z \! \in 
\! \overline{\mathbb{C}}; \, \lvert z \rvert \! > \! \hat{\varepsilon}_{
\infty}^{-1} \rbrace$, with $\hat{Q}_{\hat{\mathcal{R}},2}(+) \! := \! 
\lbrace \mathstrut q^{\prime} \! \in \! \lbrace 1,2,\dotsc,\mathfrak{s} 
\! - \! 1 \rbrace; \, \alpha_{p_{q^{\prime}}} \! \in \! (\hat{a}_{N+1} \! 
+ \! \hat{\delta}_{\hat{a}_{N+1}},+\infty) \rbrace$, $\hat{Q}_{\hat{
\mathcal{R}},2}(-) \! := \! \lbrace \mathstrut q^{\prime} \! \in \! \lbrace 
1,2,\dotsc,\mathfrak{s} \! - \! 1 \rbrace; \, \alpha_{p_{q^{\prime}}} \! 
\in \! (-\infty,\hat{b}_{0} \! - \! \hat{\delta}_{\hat{b}_{0}}) \rbrace$, and 
sufficiently small $\hat{\varepsilon}_{\infty} \! > \! 0$, $\hat{\delta}_{
\hat{\mathcal{R}},2}(\pm) \! > \! 0$ chosen so that $\mathscr{O}_{
\hat{\delta}_{\hat{\mathcal{R}},2}(+)}(\alpha_{p_{q^{\prime}_{1}}}) 
\cap \mathscr{O}_{\hat{\delta}_{\hat{\mathcal{R}},2}(+)}(\alpha_{p_{
q^{\prime \prime}_{1}}}) \! = \! \varnothing$ $\forall$ $q^{\prime}_{1} 
\! \neq \! q^{\prime \prime}_{1} \! \in \! \hat{Q}_{\hat{\mathcal{R}},2}
(+)$, $\mathscr{O}_{\hat{\delta}_{\hat{\mathcal{R}},2}(-)}(\alpha_{p_{
q^{\prime}_{2}}}) \cap \mathscr{O}_{\hat{\delta}_{\hat{\mathcal{R}},2}
(-)}(\alpha_{p_{q^{\prime \prime}_{2}}}) \! = \! \varnothing$ $\forall$ 
$q^{\prime}_{2} \! \neq \! q^{\prime \prime}_{2} \! \in \! \hat{Q}_{
\hat{\mathcal{R}},2}(-)$, $\mathscr{O}_{\hat{\delta}_{\hat{\mathcal{R}},
2}(+)}(\alpha_{p_{q^{\prime}_{1}}}) \cap \lbrace \hat{a}_{N+1} \! + \! 
\hat{\delta}_{\hat{a}_{N+1}} \rbrace \! = \! \varnothing \! = \! 
\mathscr{O}_{\hat{\delta}_{\hat{\mathcal{R}},2}(-)}(\alpha_{p_{q^{
\prime}_{2}}}) \cap \lbrace \hat{b}_{0} \! - \! \hat{\delta}_{\hat{b}_{0}} 
\rbrace$, and $\mathscr{O}_{\hat{\delta}_{\hat{\mathcal{R}},2}(+)}
(\alpha_{p_{q^{\prime}_{1}}}) \cap \mathscr{O}_{\infty}(\alpha_{k}) 
\! = \! \varnothing \! = \! \mathscr{O}_{\hat{\delta}_{\hat{\mathcal{R}},
2}(-)}(\alpha_{p_{q^{\prime}_{2}}}) \cap \mathscr{O}_{\infty}
(\alpha_{k})$$;$\footnote{Note that $\hat{\varepsilon}_{\infty}^{-1} 
\! \gg \! \max \lbrace \lvert \hat{b}_{0} \! - \! \hat{\delta}_{
\hat{b}_{0}} \rvert,\lvert \hat{a}_{N+1} \! + \! \hat{\delta}_{
\hat{a}_{N+1}} \rvert \rbrace$. If $\hat{Q}_{\hat{\mathcal{R}},2}(\pm) 
\! = \! \varnothing$, that is, $\# \hat{Q}_{\hat{\mathcal{R}},2}(\pm) 
\! = \! 0$, then $\cup_{q \in \hat{Q}_{\hat{\mathcal{R}},2}(\pm)} 
\mathscr{O}_{\hat{\delta}_{\hat{\mathcal{R}},2}(\pm)}(\alpha_{p_{q}}) 
\! := \! \varnothing$, in which case $\hat{A}_{\hat{\mathcal{R}},2}(+) 
\! = \! (\hat{a}_{N+1} \! + \! \hat{\delta}_{\hat{a}_{N+1}},+\infty) 
\setminus \mathscr{O}_{\infty}(\alpha_{k}) \! = \! (\hat{a}_{N+1} \! 
+ \! \hat{\delta}_{\hat{a}_{N+1}},\hat{\varepsilon}_{\infty}^{-1})$, 
$\hat{A}_{\hat{\mathcal{R}},2}(-) \! = \! (-\infty,\hat{b}_{0} \! - \! 
\hat{\delta}_{\hat{b}_{0}}) \setminus \mathscr{O}_{\infty}(\alpha_{k}) 
\! = \! (-\hat{\varepsilon}_{\infty}^{-1},\hat{b}_{0} \! - \! \hat{
\delta}_{\hat{b}_{0}})$, and $\prod_{q \in \hat{Q}_{\hat{\mathcal{R}},2}
(\pm)} \lvert z \! - \! \alpha_{p_{q}} \rvert^{((n-1)K+k) \hat{
\mathfrak{c}}_{\hat{\mathcal{R}},3}(q,\pm) \chi_{\mathscr{O}_{
\hat{\delta}_{\hat{\mathcal{R}},2}(\pm)}(\alpha_{p_{q}})}(z)} \! := 
\! 1$; moreover, $\cup_{q \in \hat{Q}_{\hat{\mathcal{R}},2}(+)} 
\mathscr{O}_{\hat{\delta}_{\hat{\mathcal{R}},2}(+)}(\alpha_{p_{q}}) 
\subset (\hat{a}_{N+1} \! + \! \hat{\delta}_{\hat{a}_{N+1}},
\hat{\varepsilon}_{\infty}^{-1})$, $\cup_{q \in \hat{Q}_{
\hat{\mathcal{R}},2}(-)} \mathscr{O}_{\hat{\delta}_{\hat{
\mathcal{R}},2}(-)}(\alpha_{p_{q}}) \subset (-\hat{\varepsilon}_{
\infty}^{-1},\hat{b}_{0} \! - \! \hat{\delta}_{\hat{b}_{0}})$, and 
$\# (\hat{Q}_{\hat{\mathcal{R}},2}(+) \cup \hat{Q}_{\hat{
\mathcal{R}},2}(-)) \! \leqslant \! \mathfrak{s} \! - \! 1$.}
\item[$\boldsymbol{\mathrm{(1)}_{iii}}$] for $z \! \in \! \hat{
\Sigma}^{3}_{p,j} \! := \! \hat{J}_{j}^{\smallfrown} \setminus 
(\hat{J}_{j}^{\smallfrown} \cap (\hat{\mathbb{U}}_{\hat{\delta}_{
\hat{b}_{j-1}}} \cup \hat{\mathbb{U}}_{\hat{\delta}_{\hat{a}_{j}}}))$, 
$j \! = \! 1,2,\dotsc,N \! + \! 1$,
\begin{equation} \label{eqhtvee4} 
\hat{v}_{\hat{\mathcal{R}}}(z) \underset{\underset{z_{o}=1+o(1)}{
\mathscr{N},n \to \infty}}{=} \mathrm{I} \! + \! \mathcal{O} \left(
\hat{\mathfrak{c}}_{\hat{\mathcal{R}},3}(j) \me^{-((n-1)K+k) \hat{
\lambda}_{\hat{\mathcal{R}},3}(j) \lvert z-\hat{p}_{+}(j) \rvert} 
\right),
\end{equation}
and, for $z \! \in \! \hat{\Sigma}^{4}_{p,j} \! := \! \hat{J}_{j}^{\smallsmile} 
\setminus (\hat{J}_{j}^{\smallsmile} \cap (\hat{\mathbb{U}}_{\hat{\delta}_{
\hat{b}_{j-1}}} \cup \hat{\mathbb{U}}_{\hat{\delta}_{\hat{a}_{j}}}))$, $j \! 
= \! 1,2,\dotsc,N \! + \! 1$,
\begin{equation} \label{eqhtvee5} 
\hat{v}_{\hat{\mathcal{R}}}(z) \underset{\underset{z_{o}=1+o(1)}{
\mathscr{N},n \to \infty}}{=} \mathrm{I} \! + \! \mathcal{O} \left(
\hat{\mathfrak{c}}_{\hat{\mathcal{R}},4}(j) \me^{-((n-1)K+k) \hat{
\lambda}_{\hat{\mathcal{R}},4}(j) \lvert z-\hat{p}_{-}(j) \rvert} 
\right),
\end{equation}
where, for $j \! = \! 1,2,\dotsc,N \! + \! 1$, $(\mathrm{M}_{2}
(\mathbb{C}) \! \ni)$ $\hat{\mathfrak{c}}_{\hat{\mathcal{R}},i_{i}}(j) \! 
=_{\underset{z_{o}=1+o(1)}{\mathscr{N},n \to \infty}} \! \mathcal{O}
(1)$, $i_{1} \! = \! 3,4$, $\hat{\lambda}_{\hat{\mathcal{R}},i_{1}}(j) \! = \! 
\hat{\lambda}_{\hat{\mathcal{R}},i_{1}}(n,k,z_{o};j) \! =_{\underset{z_{o}
=1+o(1)}{\mathscr{N},n \to \infty}} \! \mathcal{O}(1)$ and $> \! 0$, 
$\hat{p}_{+}(j) \! = \! (\hat{a}_{j} \! + \! \hat{\delta}_{\hat{a}_{j}} \cos 
\sigma_{\hat{a}_{j}}^{+},\hat{\delta}_{\hat{a}_{j}} \sin \sigma_{
\hat{a}_{j}}^{+})$, $\sigma_{\hat{a}_{j}}^{+} \! \in \! (2 \pi/3,\pi)$, and 
$\hat{p}_{-}(j) \! = \! (\hat{a}_{j} \! + \! \hat{\delta}_{\hat{a}_{j}} \cos 
\sigma_{\hat{a}_{j}}^{-},\hat{\delta}_{\hat{a}_{j}} \sin \sigma_{
\hat{a}_{j}}^{-})$, $\sigma_{\hat{a}_{j}}^{-} \! \in \! (-\pi,-2 \pi/3)$$;$ and 
\item[$\boldsymbol{\mathrm{(1)}_{iv}}$] for $z \! \in \hat{\Sigma}_{p,
j}^{5} \! := \! \partial \hat{\mathbb{U}}_{\hat{\delta}_{\hat{b}_{j-1}}} 
\cup \partial \hat{\mathbb{U}}_{\hat{\delta}_{\hat{a}_{j}}}$, $j \! = \! 
1,2,\dotsc,N \! + \! 1$,
\begin{align}
\hat{v}_{\hat{\mathcal{R}}}(z) \underset{\underset{z_{o}=1+o(1)}{
\mathscr{N},n \to \infty}}{=}& \, \mathrm{I} \! + \! \dfrac{1}{((n \! 
- \! 1)K \! + \! k) \hat{\xi}_{\hat{b}_{j-1}}(z)} \mathbb{M}(z) 
\begin{pmatrix}
\mp (s_{1} \! + \! t_{1}) & \mp \mi (s_{1} \! - \! t_{1}) \me^{\mi ((n-1)K+k) 
\hat{\mho}_{j-1}} \\
\mp \mi (s_{1} \! - \! t_{1}) \me^{-\mi ((n-1)K+k) \hat{\mho}_{j-1}} & 
\pm (s_{1} \! + \! t_{1})
\end{pmatrix}(\mathbb{M}(z))^{-1} \nonumber \\
+& \, \mathcal{O} \left(\dfrac{1}{((n \! - \! 1)K \! + \! k)^{2}(\hat{\xi}_{
\hat{b}_{j-1}}(z))^{2}} \mathbb{M}(z) \hat{\mathfrak{c}}^{\triangleright}
(n,k,z_{o};j)(\mathbb{M}(z))^{-1} \right), \quad z \! \in \! \mathbb{C}_{\pm} 
\cap \partial \hat{\mathbb{U}}_{\hat{\delta}_{\hat{b}_{j-1}}}, \label{eqhtvee6} \\
\hat{v}_{\hat{\mathcal{R}}}(z) \underset{\underset{z_{o}=1+o(1)}{
\mathscr{N},n \to \infty}}{=}& \, \mathrm{I} \! + \! \dfrac{1}{((n \! 
- \! 1)K \! + \! k) \hat{\xi}_{\hat{a}_{j}}(z)} \mathbb{M}(z) 
\begin{pmatrix}
\mp (s_{1} \! + \! t_{1}) & \pm \mi (s_{1} \! - \! t_{1}) \me^{\mi ((n-1)K+k) 
\hat{\mho}_{j}} \\
\pm \mi (s_{1} \! - \! t_{1}) \me^{-\mi ((n-1)K+k) \hat{\mho}_{j}} & \pm 
(s_{1} \! + \! t_{1})
\end{pmatrix}(\mathbb{M}(z))^{-1} \nonumber \\
+& \, \mathcal{O} \left(\dfrac{1}{((n \! - \! 1)K \! + \! k)^{2}(\hat{\xi}_{
\hat{a}_{j}}(z))^{2}} \mathbb{M}(z) \hat{\mathfrak{c}}^{\triangleleft}(n,k,
z_{o};j)(\mathbb{M}(z))^{-1} \right), \quad z \! \in \! \mathbb{C}_{\pm} 
\cap \partial \hat{\mathbb{U}}_{\hat{\delta}_{\hat{a}_{j}}}, \label{eqhtvee7}
\end{align}
where, for $j \! = \! 1,2,\dotsc,N \! + \! 1$, $\hat{\xi}_{\hat{b}_{j-1}}(z)$ 
(resp., $\hat{\xi}_{\hat{a}_{j}}(z))$ is defined by Equation~\eqref{eqmaininf70} 
(resp., Equation~\eqref{eqmaininf71}$)$, $\mathbb{M}(z)$ is given in 
item~{\rm \pmb{(1)}} of Lemma~\ref{lem4.5}, $s_{1} \! = \! 5/72$, 
$t_{1} \! = \! -7/72$, $\hat{\mho}_{j}$ is defined in the corresponding 
item~{\rm (ii)} of Remark~\ref{rem4.4}, and $(\mathrm{M}_{2}
(\mathbb{C}) \! \ni)$ $\hat{\mathfrak{c}}^{r}(n,k,z_{o};j) \! 
=_{\underset{z_{o}=1+o(1)}{\mathscr{N},n \to \infty}} \! \mathcal{O}(1)$, 
$r \! \in \! \lbrace \triangleright,\triangleleft \rbrace$$;$
\end{enumerate}
and {\rm \pmb{(2)}} for $n \! \in \! \mathbb{N}$ and $k \! \in \! \lbrace 
1,2,\dotsc,K \rbrace$ such that $\alpha_{p_{\mathfrak{s}}} \! := \! 
\alpha_{k} \! \neq \! \infty$$:$
\begin{enumerate}
\item[$\boldsymbol{\mathrm{(2)}_{i}}$] for $z \! \in \! \tilde{\Sigma}_{
p,j}^{2} \! := \! (\tilde{a}_{j} \! + \! \tilde{\delta}_{\tilde{a}_{j}},\tilde{
b}_{j} \! - \! \tilde{\delta}_{\tilde{b}_{j}})$, $j \! = \! 1,2,\dotsc,N$,
\begin{equation} \label{eqtlvee8} 
\tilde{v}_{\tilde{\mathcal{R}}}(z) \underset{\underset{z_{o}=1+o(1)}{
\mathscr{N},n \to \infty}}{=} \mathrm{I} \! + \! \mathcal{O} \left(
\tilde{\mathfrak{c}}_{\tilde{\mathcal{R}},1}(j) \me^{-((n-1)K+k) 
\tilde{\lambda}_{\tilde{\mathcal{R}},1}(j)(z-\tilde{a}_{j}) \chi_{
\tilde{A}_{\tilde{\mathcal{R}},1}(j)}(z)} \prod_{q \in \tilde{Q}_{\tilde{
\mathcal{R}},1}(j)} \lvert z \! - \! \alpha_{p_{q}} \rvert^{((n-1)K+k) 
\tilde{\mathfrak{c}}_{\tilde{\mathcal{R}},2}(j,q) \chi_{\mathscr{O}_{
\tilde{\delta}_{\tilde{\mathcal{R}},1}(j)}(\alpha_{p_{q}})}(z)} \right),
\end{equation}
where, for $j \! = \! 1,2,\dotsc,N$, $(\mathrm{M}_{2}(\mathbb{C}) 
\! \ni)$ $\tilde{\mathfrak{c}}_{\hat{\mathcal{R}},1}(j) \! 
=_{\underset{z_{o}=1+o(1)}{\mathscr{N},n \to \infty}} \! 
\mathcal{O}(1)$, $\tilde{\lambda}_{\tilde{\mathcal{R}},1}(j) 
\! = \! \tilde{\lambda}_{\tilde{\mathcal{R}},1}(n,k,z_{o};j) \! 
=_{\underset{z_{o}=1+o(1)}{\mathscr{N},n \to \infty}} \! 
\mathcal{O}(1)$ and $> \! 0$, $\tilde{\mathfrak{c}}_{\tilde{
\mathcal{R}},2}(j,q) \! = \! \tilde{\mathfrak{c}}_{\tilde{\mathcal{R}},2}
(n,k,z_{o};j,q) \! =_{\underset{z_{o}=1+o(1)}{\mathscr{N},n \to \infty}} 
\! \mathcal{O}(1)$ and $> \! 0$, and $\tilde{A}_{\tilde{\mathcal{R}},
1}(j) \! = \! (\tilde{a}_{j} \! + \! \tilde{\delta}_{\tilde{a}_{j}},\tilde{
b}_{j} \! - \! \tilde{\delta}_{\tilde{b}_{j}}) \setminus \cup_{q \in 
\tilde{Q}_{\tilde{\mathcal{R}},1}(j)} \mathscr{O}_{\tilde{\delta}_{
\tilde{\mathcal{R}},1}(j)}(\alpha_{p_{q}})$, with $\tilde{Q}_{\tilde{
\mathcal{R}},1}(j) \! := \! \lbrace \mathstrut q^{\prime} \! \in \! 
\lbrace 1,\dotsc,\mathfrak{s} \! - \! 2,\mathfrak{s} \rbrace; \, 
\alpha_{p_{q^{\prime}}} \! \in \! (\tilde{a}_{j} \! + \! \tilde{\delta}_{
\tilde{a}_{j}},\tilde{b}_{j} \! - \! \tilde{\delta}_{\tilde{b}_{j}}) \rbrace$ 
and sufficiently small $\tilde{\delta}_{\tilde{\mathcal{R}},1}(j) \! > 
\! 0$ chosen so that $\mathscr{O}_{\tilde{\delta}_{\tilde{\mathcal{
R}},1}(j)}(\alpha_{p_{q_{1}}}) \cap \mathscr{O}_{\tilde{\delta}_{\tilde{
\mathcal{R}},1}(j)}(\alpha_{p_{q_{2}}}) \! = \! \varnothing$ $\forall$ 
$q_{1} \! \neq \! q_{2} \! \in \! \tilde{Q}_{\tilde{\mathcal{R}},1}(j)$ 
and $\mathscr{O}_{\tilde{\delta}_{\tilde{\mathcal{R}},1}(j)}(\alpha_{
p_{q_{1}}}) \cap \lbrace \tilde{a}_{j} \! + \! \tilde{\delta}_{\tilde{a}_{j}} 
\rbrace \! = \! \varnothing \! = \! \mathscr{O}_{\tilde{\delta}_{\tilde{
\mathcal{R}},1}(j)}(\alpha_{p_{q_{1}}}) \cap \lbrace \tilde{b}_{j} \! - \! 
\tilde{\delta}_{\tilde{b}_{j}} \rbrace$$;$\footnote{If, for $j \! = \! 1,2,
\dotsc,N$, $\tilde{Q}_{\tilde{\mathcal{R}},1}(j) \! = \! \varnothing$, 
that is, $\# \tilde{Q}_{\tilde{\mathcal{R}},1}(j) \! = \! 0$, then 
$\cup_{q \in \tilde{Q}_{\tilde{\mathcal{R}},1}(j)} \mathscr{O}_{\tilde{
\delta}_{\tilde{\mathcal{R}},1}(j)}(\alpha_{p_{q}}) \! := \! \varnothing$, 
in which case $\tilde{A}_{\tilde{\mathcal{R}},1}(j) \! = \! (\tilde{a}_{j} 
\! + \! \tilde{\delta}_{\tilde{a}_{j}},\tilde{b}_{j} \! - \! \tilde{\delta}_{
\tilde{b}_{j}})$, and $\prod_{q \in \tilde{Q}_{\tilde{\mathcal{R}},1}(j)} 
\lvert z \! - \! \alpha_{p_{q}} \rvert^{((n-1)K+k) \tilde{\mathfrak{c}}_{
\tilde{\mathcal{R}},2}(j,q) \chi_{\mathscr{O}_{\tilde{\delta}_{\tilde{
\mathcal{R}},1}(j)}(\alpha_{p_{q}})}(z)} \! := \! 1$ (of course, 
$\cup_{q \in \tilde{Q}_{\tilde{\mathcal{R}},1}(j)} \mathscr{O}_{\tilde{
\delta}_{\tilde{\mathcal{R}},1}(j)}(\alpha_{p_{q}}) \subset (\tilde{a}_{j} 
\! + \! \tilde{\delta}_{\tilde{a}_{j}},\tilde{b}_{j} \! - \! \tilde{\delta}_{
\tilde{b}_{j}}))$; moreover, $\# \cup_{j=1}^{N} \tilde{Q}_{\tilde{
\mathcal{R}},1}(j) \! \leqslant \! \mathfrak{s} \! - \! 1$.}
\item[$\boldsymbol{\mathrm{(2)}_{ii}}$] for $z \! \in \! 
\tilde{\Sigma}^{1}_{p} \! := \! (-\infty,\tilde{b}_{0} \! - \! \tilde{
\delta}_{\tilde{b}_{0}}) \cup (\tilde{a}_{N+1} \! + \! \tilde{\delta}_{
\tilde{a}_{N+1}},+\infty)$,
\begin{align}
\tilde{v}_{\tilde{\mathcal{R}}}(z) \underset{\underset{z_{o}=1+o(1)}{
\mathscr{N},n \to \infty}}{=}& \, \mathrm{I} \! + \! \mathcal{O} 
\left(\tilde{\mathfrak{c}}_{\tilde{\mathcal{R}},2}(+) \me^{-((n-1)K
+k) \tilde{\lambda}_{\tilde{\mathcal{R}},2}(+)(z-\tilde{a}_{N+1}) 
\chi_{\tilde{A}_{\tilde{\mathcal{R}},2}(+)}(z)} \me^{-((n-1)K+k) 
\tilde{\lambda}_{\tilde{\mathcal{R}},3}(+) \ln (\lvert z \rvert) \chi_{
\mathscr{O}_{\infty}(\alpha_{p_{\mathfrak{s}-1}})}(z)} \right. 
\nonumber \\
\times&\left. \, \prod_{q \in \tilde{Q}_{\tilde{\mathcal{R}},2}(+)} 
\lvert z \! - \! \alpha_{p_{q}} \rvert^{((n-1)K+k) \tilde{\mathfrak{c}}_{
\tilde{\mathcal{R}},3}(q,+) \chi_{\mathscr{O}_{\tilde{\delta}_{\tilde{
\mathcal{R}},2}(+)}(\alpha_{p_{q}})}(z)} \right), \quad z \! \in \! 
(\tilde{a}_{N+1} \! + \! \tilde{\delta}_{\tilde{a}_{N+1}},+\infty), 
\label{eqtlvee9} \\
\tilde{v}_{\tilde{\mathcal{R}}}(z) \underset{\underset{z_{o}=1+o(1)}{
\mathscr{N},n \to \infty}}{=}& \, \mathrm{I} \! + \! \mathcal{O} 
\left(\tilde{\mathfrak{c}}_{\tilde{\mathcal{R}},2}(-) \me^{-((n-1)K
+k) \tilde{\lambda}_{\tilde{\mathcal{R}},2}(-) \lvert z-\tilde{b}_{0} 
\rvert \chi_{\tilde{A}_{\tilde{\mathcal{R}},2}(-)}(z)} \me^{-((n-1)K
+k) \tilde{\lambda}_{\tilde{\mathcal{R}},3}(-) \ln (\lvert z \rvert) 
\chi_{\mathscr{O}_{\infty}(\alpha_{p_{\mathfrak{s}-1}})}(z)} \right. 
\nonumber \\
\times&\left. \, \prod_{q \in \tilde{Q}_{\tilde{\mathcal{R}},2}(-)} 
\lvert z \! - \! \alpha_{p_{q}} \rvert^{((n-1)K+k) \tilde{\mathfrak{
c}}_{\tilde{\mathcal{R}},3}(q,-) \chi_{\mathscr{O}_{\tilde{\delta}_{
\tilde{\mathcal{R}},2}(-)}(\alpha_{p_{q}})}(z)} \right), \quad z \! 
\in \! (-\infty,\tilde{b}_{0} \! - \! \tilde{\delta}_{\tilde{b}_{0}}), 
\label{eqtlvee10}
\end{align}
where $(\mathrm{M}_{2}(\mathbb{C}) \! \ni)$ $\tilde{\mathfrak{c}}_{
\tilde{\mathcal{R}},2}(\pm) \! =_{\underset{z_{o}=1+o(1)}{\mathscr{N},
n \to \infty}} \! \mathcal{O}(1)$, $\tilde{\lambda}_{\tilde{\mathcal{R}},
i_{1}}(\pm) \! = \! \tilde{\lambda}_{\tilde{\mathcal{R}},i_{1}}(n,k,z_{o};
\pm) \! =_{\underset{z_{o}=1+o(1)}{\mathscr{N},n \to \infty}} \! 
\mathcal{O}(1)$ and $> \! 0$, $i_{1} \! = \! 2,3$, $\tilde{\mathfrak{c}}_{
\tilde{\mathcal{R}},3}(q,\pm) \! = \! \tilde{\mathfrak{c}}_{\tilde{\mathcal{R}},3}
(n,k,z_{o};q,\pm) \! =_{\underset{z_{o}=1+o(1)}{\mathscr{N},n \to \infty}} 
\! \mathcal{O}(1)$ and $> \! 0$, $\tilde{A}_{\tilde{\mathcal{R}},2}(+) \! 
= \! (\tilde{a}_{N+1} \! + \! \tilde{\delta}_{\tilde{a}_{N+1}},+\infty) 
\setminus (\mathscr{O}_{\infty}(\alpha_{p_{\mathfrak{s}-1}}) \cup 
\cup_{q \in \tilde{Q}_{\tilde{\mathcal{R}},2}(+)} \mathscr{O}_{\tilde{
\delta}_{\tilde{\mathcal{R}},2}(+)}(\alpha_{p_{q}}))$, $\tilde{A}_{
\tilde{\mathcal{R}},2}(-) \! = \! (-\infty,\tilde{b}_{0} \! - \! \tilde{
\delta}_{\tilde{b}_{0}}) \setminus (\mathscr{O}_{\infty}(\alpha_{p_{
\mathfrak{s}-1}}) \cup \cup_{q \in \tilde{Q}_{\tilde{\mathcal{R}},2}
(-)}  \mathscr{O}_{\tilde{\delta}_{\tilde{\mathcal{R}},2}(-)}(\alpha_{
p_{q}}))$, and $\mathscr{O}_{\infty}(\alpha_{p_{\mathfrak{s}-1}}) \! 
= \! \lbrace \mathstrut z \! \in \! \overline{\mathbb{C}}; \, \lvert 
z \rvert \! > \! \tilde{\varepsilon}_{\infty}^{-1} \rbrace$, with 
$\tilde{Q}_{\tilde{\mathcal{R}},2}(+) \! := \! \lbrace \mathstrut 
q^{\prime} \! \in \! \lbrace 1,\dotsc,\mathfrak{s} \! - \! 2,
\mathfrak{s} \rbrace; \, \alpha_{p_{q^{\prime}}} \! \in \! (\tilde{a}_{
N+1} \! + \! \tilde{\delta}_{\tilde{a}_{N+1}},+\infty) \rbrace$, 
$\tilde{Q}_{\tilde{\mathcal{R}},2}(-) \! := \! \lbrace \mathstrut 
q^{\prime} \! \in \! \lbrace 1,\dotsc,\mathfrak{s} \! - \! 2,
\mathfrak{s} \rbrace; \, \alpha_{p_{q^{\prime}}} \! \in \! (-\infty,
\tilde{b}_{0} \! - \! \tilde{\delta}_{\tilde{b}_{0}}) \rbrace$, and 
sufficiently small $\tilde{\varepsilon}_{\infty} \! > \! 0$, $\tilde{
\delta}_{\tilde{\mathcal{R}},2}(\pm) \! > \! 0$ chosen so that 
$\mathscr{O}_{\tilde{\delta}_{\tilde{\mathcal{R}},2}(+)}(\alpha_{
p_{q^{\prime}_{1}}}) \cap \mathscr{O}_{\tilde{\delta}_{\tilde{
\mathcal{R}},2}(+)}(\alpha_{p_{q^{\prime \prime}_{1}}}) \! = \! 
\varnothing$ $\forall$ $q^{\prime}_{1} \! \neq \! q^{\prime 
\prime}_{1} \! \in \! \tilde{Q}_{\tilde{\mathcal{R}},2}(+)$, 
$\mathscr{O}_{\tilde{\delta}_{\tilde{\mathcal{R}},2}(-)}(\alpha_{
p_{q^{\prime}_{2}}}) \cap \mathscr{O}_{\tilde{\delta}_{\tilde{
\mathcal{R}},2}(-)}(\alpha_{p_{q^{\prime \prime}_{2}}}) \! = \! 
\varnothing$ $\forall$ $q^{\prime}_{2} \! \neq \! q^{\prime 
\prime}_{2} \! \in \! \tilde{Q}_{\tilde{\mathcal{R}},2}(-)$, 
$\mathscr{O}_{\tilde{\delta}_{\tilde{\mathcal{R}},2}(+)}(\alpha_{
p_{q^{\prime}_{1}}}) \cap \lbrace \tilde{a}_{N+1} \! + \! \tilde{
\delta}_{\tilde{a}_{N+1}} \rbrace \! = \! \varnothing \! = \! 
\mathscr{O}_{\tilde{\delta}_{\tilde{\mathcal{R}},2}(-)}(\alpha_{p_{
q^{\prime}_{2}}}) \cap \lbrace \tilde{b}_{0} \! - \! \tilde{\delta}_{
\tilde{b}_{0}} \rbrace$, and $\mathscr{O}_{\tilde{\delta}_{\tilde{
\mathcal{R}},2}(+)}(\alpha_{p_{q^{\prime}_{1}}}) \cap \mathscr{O}_{
\infty}(\alpha_{p_{\mathfrak{s}-1}}) \! = \! \varnothing \! = \! 
\mathscr{O}_{\tilde{\delta}_{\tilde{\mathcal{R}},2}(-)}(\alpha_{p_{
q^{\prime}_{2}}}) \cap \mathscr{O}_{\infty}(\alpha_{p_{\mathfrak{s}
-1}})$$;$\footnote{Note that $\tilde{\varepsilon}_{\infty}^{-1} \! \gg 
\! \max \lbrace \lvert \tilde{b}_{0} \! - \! \tilde{\delta}_{\tilde{b}_{0}} 
\rvert,\lvert \tilde{a}_{N+1} \! + \! \tilde{\delta}_{\tilde{a}_{N+1}} 
\rvert \rbrace$. If $\tilde{Q}_{\tilde{\mathcal{R}},2}(\pm) \! = \! 
\varnothing$, that is, $\# \tilde{Q}_{\tilde{\mathcal{R}},2}(\pm) 
\! = \! 0$, then $\cup_{q \in \tilde{Q}_{\tilde{\mathcal{R}},2}(\pm)} 
\mathscr{O}_{\tilde{\delta}_{\tilde{\mathcal{R}},2}(\pm)}(\alpha_{
p_{q}}) \! := \! \varnothing$, in which case $\tilde{A}_{\tilde{
\mathcal{R}},2}(+) \! = \! (\tilde{a}_{N+1} \! + \! \tilde{\delta}_{
\tilde{a}_{N+1}},+\infty) \setminus \mathscr{O}_{\infty}(\alpha_{
p_{\mathfrak{s}-1}}) \! = \! (\tilde{a}_{N+1} \! + \! \tilde{\delta}_{
\tilde{a}_{N+1}},\tilde{\varepsilon}_{\infty}^{-1})$, $\tilde{A}_{
\tilde{\mathcal{R}},2}(-) \! = \! (-\infty,\tilde{b}_{0} \! - \! \tilde{
\delta}_{\tilde{b}_{0}}) \setminus \mathscr{O}_{\infty}(\alpha_{
p_{\mathfrak{s}-1}}) \! = \! (-\tilde{\varepsilon}_{\infty}^{-1},
\tilde{b}_{0} \! - \! \tilde{\delta}_{\tilde{b}_{0}})$, and 
$\prod_{q \in \tilde{Q}_{\tilde{\mathcal{R}},2}(\pm)} \lvert z 
\! - \! \alpha_{p_{q}} \rvert^{((n-1)K+k) \tilde{\mathfrak{c}}_{
\tilde{\mathcal{R}},3}(q,\pm) \chi_{\mathscr{O}_{\tilde{\delta}_{
\tilde{\mathcal{R}},2}(\pm)}(\alpha_{p_{q}})}(z)} \! := \! 1$; 
moreover, $\cup_{q \in \tilde{Q}_{\tilde{\mathcal{R}},2}(+)} 
\mathscr{O}_{\tilde{\delta}_{\tilde{\mathcal{R}},2}(+)}(\alpha_{
p_{q}}) \subset (\tilde{a}_{N+1} \! + \! \tilde{\delta}_{\tilde{a}_{N
+1}},\tilde{\varepsilon}_{\infty}^{-1})$, $\cup_{q \in \tilde{Q}_{
\tilde{\mathcal{R}},2}(-)} \mathscr{O}_{\tilde{\delta}_{\tilde{
\mathcal{R}},2}(-)}(\alpha_{p_{q}}) \subset (-\tilde{\varepsilon}_{
\infty}^{-1},\tilde{b}_{0} \! - \! \tilde{\delta}_{\tilde{b}_{0}})$, 
and $\# (\tilde{Q}_{\tilde{\mathcal{R}},2}(+) \cup \tilde{Q}_{
\tilde{\mathcal{R}},2}(-)) \! \leqslant \! \mathfrak{s} \! - \! 1$.}
\item[$\boldsymbol{\mathrm{(2)}_{iii}}$] for $z \! \in \! \tilde{
\Sigma}^{3}_{p,j} \! := \! \tilde{J}_{j}^{\smallfrown} \setminus 
(\tilde{J}_{j}^{\smallfrown} \cap (\tilde{\mathbb{U}}_{\tilde{\delta}_{
\tilde{b}_{j-1}}} \cup \tilde{\mathbb{U}}_{\tilde{\delta}_{\tilde{a}_{j}}}))$, 
$j \! = \! 1,2,\dotsc,N \! + \! 1$,
\begin{equation} \label{eqtlvee11} 
\tilde{v}_{\tilde{\mathcal{R}}}(z) \underset{\underset{z_{o}=1+o(1)}{
\mathscr{N},n \to \infty}}{=} \mathrm{I} \! + \! \mathcal{O} \left(
\tilde{\mathfrak{c}}_{\tilde{\mathcal{R}},3}(j) \me^{-((n-1)K+k) 
\tilde{\lambda}_{\tilde{\mathcal{R}},3}(j) \lvert z-\tilde{p}_{+}(j) 
\rvert} \right),
\end{equation}
and, for $z \! \in \! \tilde{\Sigma}^{4}_{p,j} \! := \! \tilde{J}_{j}^{
\smallsmile} \setminus (\tilde{J}_{j}^{\smallsmile} \cap (\tilde{
\mathbb{U}}_{\tilde{\delta}_{\tilde{b}_{j-1}}} \cup \tilde{\mathbb{U}}_{
\tilde{\delta}_{\tilde{a}_{j}}}))$, $j \! = \! 1,2,\dotsc,N \! + \! 1$,
\begin{equation} \label{eqtlvee12} 
\tilde{v}_{\tilde{\mathcal{R}}}(z) \underset{\underset{z_{o}=1+o(1)}{
\mathscr{N},n \to \infty}}{=} \mathrm{I} \! + \! \mathcal{O} \left(
\tilde{\mathfrak{c}}_{\tilde{\mathcal{R}},4}(j) \me^{-((n-1)K+k) 
\tilde{\lambda}_{\tilde{\mathcal{R}},4}(j) \lvert z-\tilde{p}_{-}(j) 
\rvert} \right),
\end{equation}
where, for $j \! = \! 1,2,\dotsc,N \! + \! 1$, $(\mathrm{M}_{2}
(\mathbb{C}) \! \ni)$ $\tilde{\mathfrak{c}}_{\tilde{\mathcal{R}},i_{i}}(j) 
\! =_{\underset{z_{o}=1+o(1)}{\mathscr{N},n \to \infty}} \! \mathcal{O}
(1)$, $i_{1} \! = \! 3,4$, $\tilde{\lambda}_{\tilde{\mathcal{R}},i_{1}}(j) 
\! = \! \tilde{\lambda}_{\tilde{\mathcal{R}},i_{1}}(n,k,z_{o};j) \! 
=_{\underset{z_{o}=1+o(1)}{\mathscr{N},n \to \infty}} \! \mathcal{O}(1)$ 
and $> \! 0$, $\tilde{p}_{+}(j) \! = \! (\tilde{a}_{j} \! + \! \tilde{\delta}_{
\tilde{a}_{j}} \cos \sigma_{\tilde{a}_{j}}^{+},\tilde{\delta}_{\tilde{a}_{j}} \sin 
\sigma_{\tilde{a}_{j}}^{+})$, $\sigma_{\tilde{a}_{j}}^{+} \! \in \! (2 \pi/3,
\pi)$, and $\tilde{p}_{-}(j) \! = \! (\tilde{a}_{j} \! + \! \tilde{\delta}_{
\tilde{a}_{j}} \cos \sigma_{\tilde{a}_{j}}^{-},\tilde{\delta}_{\tilde{a}_{j}} 
\sin \sigma_{\tilde{a}_{j}}^{-})$, $\sigma_{\tilde{a}_{j}}^{-} \! \in \! 
(-\pi,-2 \pi/3)$$;$ and 
\item[$\boldsymbol{\mathrm{(2)}_{iv}}$] for $z \! \in \tilde{\Sigma}_{p,
j}^{5} \! := \! \partial \tilde{\mathbb{U}}_{\tilde{\delta}_{\tilde{b}_{j-1}}} 
\cup \partial \tilde{\mathbb{U}}_{\tilde{\delta}_{\tilde{a}_{j}}}$, $j \! 
= \! 1,2,\dotsc,N \! + \! 1$,
\begin{align}
\tilde{v}_{\tilde{\mathcal{R}}}(z) \underset{\underset{z_{o}=1+o(1)}{
\mathscr{N},n \to \infty}}{=}& \, \mathrm{I} \! + \! \dfrac{1}{((n \! 
- \! 1)K \! + \! k) \tilde{\xi}_{\tilde{b}_{j-1}}(z)} \mathbb{M}(z) 
\begin{pmatrix}
\mp (s_{1} \! + \! t_{1}) & \mp \mi (s_{1} \! - \! t_{1}) \me^{\mi ((n-1)K
+k) \tilde{\mho}_{j-1}} \\
\mp \mi (s_{1} \! - \! t_{1}) \me^{-\mi ((n-1)K+k) \tilde{\mho}_{j-1}} & 
\pm (s_{1} \! + \! t_{1})
\end{pmatrix}(\mathbb{M}(z))^{-1} \nonumber \\
+& \, \mathcal{O} \left(\dfrac{1}{((n \! - \! 1)K \! + \! k)^{2}(\tilde{\xi}_{
\tilde{b}_{j-1}}(z))^{2}} \mathbb{M}(z) \tilde{\mathfrak{c}}^{\triangleright}
(n,k,z_{o};j)(\mathbb{M}(z))^{-1} \right), \quad z \! \in \! \mathbb{C}_{\pm} 
\cap \partial \tilde{\mathbb{U}}_{\tilde{\delta}_{\tilde{b}_{j-1}}}, 
\label{eqtlvee13} \\
\tilde{v}_{\tilde{\mathcal{R}}}(z) \underset{\underset{z_{o}=1+o(1)}{
\mathscr{N},n \to \infty}}{=}& \, \mathrm{I} \! + \! \dfrac{1}{((n \! 
- \! 1)K \! + \! k) \tilde{\xi}_{\tilde{a}_{j}}(z)} \mathbb{M}(z) 
\begin{pmatrix}
\mp (s_{1} \! + \! t_{1}) & \pm \mi (s_{1} \! - \! t_{1}) \me^{\mi ((n-1)K
+k) \tilde{\mho}_{j}} \\
\pm \mi (s_{1} \! - \! t_{1}) \me^{-\mi ((n-1)K+k) \tilde{\mho}_{j}} & 
\pm (s_{1} \! + \! t_{1})
\end{pmatrix}(\mathbb{M}(z))^{-1} \nonumber \\
+& \, \mathcal{O} \left(\dfrac{1}{((n \! - \! 1)K \! + \! k)^{2}(\tilde{\xi}_{
\tilde{a}_{j}}(z))^{2}} \mathbb{M}(z) \tilde{\mathfrak{c}}^{\triangleleft}(n,
k,z_{o};j)(\mathbb{M}(z))^{-1} \right), \quad z \! \in \! \mathbb{C}_{\pm} 
\cap \partial \tilde{\mathbb{U}}_{\tilde{\delta}_{\tilde{a}_{j}}}, 
\label{eqtlvee14}
\end{align}
where, for $j \! = \! 1,2,\dotsc,N \! + \! 1$, $\tilde{\xi}_{\tilde{b}_{j-1}}(z)$ 
(resp., $\tilde{\xi}_{\tilde{a}_{j}}(z))$ is defined by Equation~\eqref{eqmainfin72} 
(resp., Equation~\eqref{eqmainfin73}$)$, $\mathbb{M}(z)$ is given in 
item~{\rm \pmb{(2)}} of Lemma~\ref{lem4.5}, $s_{1} \! = \! 5/72$, 
$t_{1} \! = \! -7/72$, $\tilde{\mho}_{j}$ is defined in the corresponding 
item~{\rm (ii)} of Remark~\ref{rem4.4}, and $(\mathrm{M}_{2}(\mathbb{
C}) \! \ni)$ $\tilde{\mathfrak{c}}^{r}(n,k,z_{o};j) \! =_{\underset{z_{o}=1+
o(1)}{\mathscr{N},n \to \infty}} \! \mathcal{O}(1)$, $r \! \in \! \lbrace 
\triangleright,\triangleleft \rbrace$.
\end{enumerate}
\end{ccccc}

\emph{Proof}. The proof of this Lemma~\ref{lem5.1} consists of two cases: 
(i) $n \! \in \! \mathbb{N}$ and $k \! \in \! \lbrace 1,2,\dotsc,K \rbrace$ 
such that $\alpha_{p_{\mathfrak{s}}} \! := \! \alpha_{k} \! = \! \infty$; and 
(ii) $n \! \in \! \mathbb{N}$ and $k \! \in \! \lbrace 1,2,\dotsc,K \rbrace$ 
such that $\alpha_{p_{\mathfrak{s}}} \! := \! \alpha_{k} \! \neq \! \infty$. 
Notwithstanding the fact that the scheme of the proof is, 
\emph{mutatis mutandis}, similar for both cases, without loss of generality, 
only the proof for case~(ii) is presented in detail, whilst case~(i) is proved 
analogously.

For $n \! \in \! \mathbb{N}$ and $k \! \in \! \lbrace 1,2,\dotsc,K \rbrace$ 
such that $\alpha_{p_{\mathfrak{s}}} \! := \! \alpha_{k} \! \neq \! \infty$, 
and for $z \! \in \! \tilde{\Sigma}_{p,j}^{2} \! := \! (\tilde{a}_{j} \! + \! 
\tilde{\delta}_{\tilde{a}_{j}},\tilde{b}_{j} \! - \! \tilde{\delta}_{\tilde{b}_{j}})$, 
$j \! = \! 1,2,\dotsc,N$, define $\tilde{A}_{\tilde{\mathcal{R}},1}(j)$ 
and $\tilde{Q}_{\tilde{\mathcal{R}},1}(j)$, and choose $\tilde{\delta}_{
\tilde{\mathcal{R}},1}(j)$ $(> \! 0)$, as in item~\pmb{(2)}, 
subitem~$\boldsymbol{\mathrm{(2)}_{i}}$ of the lemma: via the definition 
of $\tilde{A}_{\tilde{\mathcal{R}},1}(j)$, it follows that the complement 
of $\tilde{A}_{\tilde{\mathcal{R}},1}(j)$ relative to $(\tilde{a}_{j} \! + \! 
\tilde{\delta}_{\tilde{a}_{j}},\tilde{b}_{j} \! - \! \tilde{\delta}_{\tilde{b}_{j}})$ 
is given by $\tilde{A}_{\tilde{\mathcal{R}},1}^{c}(j) \! = \! \cup_{q \in 
\tilde{Q}_{\tilde{\mathcal{R}},1}(j)} \mathscr{O}_{\tilde{\delta}_{\tilde{
\mathcal{R}},1}(j)}(\alpha_{p_{q}})$, $j \! = \! 1,2,\dotsc,N$ (that is, 
$\tilde{\Sigma}_{p,j}^{2} \! = \! \tilde{A}_{\tilde{\mathcal{R}},1}(j) \cup 
\tilde{A}_{\tilde{\mathcal{R}},1}^{c}(j))$; furthermore, since $\tilde{A}_{
\tilde{\mathcal{R}},1}(j) \cap \tilde{A}_{\tilde{\mathcal{R}},1}^{c}(j) \! = \! 
\varnothing$ and $\mathscr{O}_{\tilde{\delta}_{\tilde{\mathcal{R}},1}(j)}
(\alpha_{p_{q^{\prime}}}) \cap \mathscr{O}_{\tilde{\delta}_{\tilde{
\mathcal{R}},1}(j)}(\alpha_{p_{q^{\prime \prime}}}) \! = \! \varnothing$, 
$j \! = \! 1,2,\dotsc,N$, $q^{\prime} \! \neq \! q^{\prime \prime} 
\! \in \! \tilde{Q}_{\tilde{\mathcal{R}},1}(j)$, it follows that 
$\chi_{\tilde{\Sigma}_{p,j}^{2}}(z) \! = \! \chi_{\tilde{A}_{\tilde{\mathcal{R}},1}
(j) \cup \tilde{A}_{\tilde{\mathcal{R}},1}^{c}(j)}(z) \! = \! \chi_{\tilde{A}_{
\tilde{\mathcal{R}},1}(j)}(z) \! + \! \sum_{q \in \tilde{Q}_{\tilde{\mathcal{R}},1}(j)} 
\chi_{\mathscr{O}_{\tilde{\delta}_{\tilde{\mathcal{R}},1}(j)}(\alpha_{p_{q}})}(z)$, 
whence, via the expression for $\tilde{v}_{\tilde{\mathcal{R}}}(z)$ given in the 
corresponding item~\pmb{(2)} of Lemma~\ref{lem4.10}, one shows that, for 
$z \! \in \! \tilde{\Sigma}_{p,j}^{2}$, $j \! = \! 1,2,\dotsc,N$,
\begin{align} \label{eqlem5.1a} 
\tilde{v}_{\tilde{\mathcal{R}}}(z) =& \, \mathrm{I} \! + \! \me^{-\mi 
((n-1)K+k) \tilde{\Omega}_{j}} \exp \left(n \left(g^{f}_{+}(z) \! + \! 
g^{f}_{-}(z) \! - \! \hat{\mathscr{P}}_{0}^{+} \! - \! \hat{\mathscr{P}}_{
0}^{-} \! - \! \widetilde{V}(z) \! - \! \tilde{\ell} \right) \left(\chi_{
\tilde{A}_{\tilde{\mathcal{R}},1}(j)}(z) \! + \! \sum_{q \in \tilde{Q}_{
\tilde{\mathcal{R}},1}(j)} \chi_{\mathscr{O}_{\tilde{\delta}_{\tilde{
\mathcal{R}},1}(j)}(\alpha_{p_{q}})}(z) \right) \right) \nonumber \\
\times& \, \left(\mathfrak{m}_{-}(z) \left(\chi_{\tilde{A}_{\tilde{
\mathcal{R}},1}(j)}(z) \! + \! \sum_{q \in \tilde{Q}_{\tilde{\mathcal{R}},
1}(j)} \chi_{\mathscr{O}_{\tilde{\delta}_{\tilde{\mathcal{R}},1}(j)}
(\alpha_{p_{q}})}(z) \right) \right) \sigma_{+} \left(\mathfrak{m}_{-}
(z) \left(\chi_{\tilde{A}_{\tilde{\mathcal{R}},1}(j)}(z) \! + \! \sum_{q 
\in \tilde{Q}_{\tilde{\mathcal{R}},1}(j)} \chi_{\mathscr{O}_{\tilde{
\delta}_{\tilde{\mathcal{R}},1}(j)}(\alpha_{p_{q}})}(z) \right) \right)^{-1}.
\end{align}
For $n \! \in \! \mathbb{N}$ and $k \! \in \! \lbrace 1,2,\dotsc,K 
\rbrace$ such that $\alpha_{p_{\mathfrak{s}}} \! := \! \alpha_{k} \! 
\neq \! \infty$, and for $z \! \in \! \tilde{\Sigma}_{p,j}^{2}$, $j \! 
= \! 1,2,\dotsc,N$, if $\tilde{Q}_{\tilde{\mathcal{R}},1}(j) \! = \! 
\varnothing$ (that is, $\# \tilde{Q}_{\tilde{\mathcal{R}},1}(j) \! = \! 
0)$, then one uses the Representation~\eqref{eql3.8AA} for $g^{f}_{+}
(z) \! + \! g^{f}_{-}(z) \! - \! \hat{\mathscr{P}}_{0}^{+} \! - \! \hat{
\mathscr{P}}_{0}^{-} \! - \! \widetilde{V}(z) \! - \! \tilde{\ell}$, or, if 
$\tilde{Q}_{\tilde{\mathcal{R}},1}(j) \! \neq \! \varnothing$, then 
one uses the (equivalent) Representation~\eqref{eql3.8f} for 
$g^{f}_{+}(z) \! + \! g^{f}_{-}(z) \! - \! \hat{\mathscr{P}}_{0}^{+} \! - 
\! \hat{\mathscr{P}}_{0}^{-} \! - \! \widetilde{V}(z) \! - \! \tilde{\ell}$: 
for $j \! = \! 1,2,\dotsc,N$, via the above-mentioned equivalent 
representations for $g^{f}_{+}(z) \! + \! g^{f}_{-}(z) \! - \! \hat{
\mathscr{P}}_{0}^{+} \! - \! \hat{\mathscr{P}}_{0}^{-} \! - \! 
\widetilde{V}(z) \! - \! \tilde{\ell}$, the asymptotics (cf. the proof of 
Lemma~\ref{lem3.8}, case~\pmb{$(\mathrm{A})$}, subitem~\pmb{(2)})
\begin{gather*}
g^{f}_{+}(z) \! + \! g^{f}_{-}(z) \! - \! \hat{\mathscr{P}}_{0}^{+} 
\! - \! \hat{\mathscr{P}}_{0}^{-} \! - \! \widetilde{V}(z) \! - \! 
\tilde{\ell} \underset{\underset{q \in \tilde{Q}_{\tilde{\mathcal{R}},1}
(j)}{\mathscr{O}_{\tilde{\delta}_{\tilde{\mathcal{R}},1}(j)}(\alpha_{p_{q}}) 
\ni z \to \alpha_{p_{q}}=\alpha_{k}}}{=} -\left(\widetilde{V}(z) \! - \! 
\left(\dfrac{\varkappa_{nk} \! - \! 1}{n} \right) \ln (1 \! + \! \lvert z 
\! - \! \alpha_{k} \rvert^{-2}) \right) \! + \! \mathcal{O}(1), \\
g^{f}_{+}(z) \! + \! g^{f}_{-}(z) \! - \! \hat{\mathscr{P}}_{0}^{+} 
\! - \! \hat{\mathscr{P}}_{0}^{-} \! - \! \widetilde{V}(z) \! - \! 
\tilde{\ell} \underset{\underset{q \in \tilde{Q}_{\tilde{\mathcal{R}},1}
(j)}{\mathscr{O}_{\tilde{\delta}_{\tilde{\mathcal{R}},1}(j)}(\alpha_{p_{q}}) 
\ni z \to \alpha_{p_{q}} \neq \alpha_{k}}}{=} -\left(\widetilde{V}(z) 
\! - \! \dfrac{\varkappa_{nk \tilde{k}_{q}}}{n} \ln (1 \! + \! \lvert z \! 
- \! \alpha_{p_{q}} \rvert^{-2}) \right) \! + \! \mathcal{O}(1),
\end{gather*}
the conditions~\eqref{eq20}--\eqref{eq22} for regular $\widetilde{V} 
\colon \overline{\mathbb{R}} \setminus \lbrace \alpha_{1},\alpha_{2},
\dotsc,\alpha_{K} \rbrace \! \to \! \mathbb{R}$,\footnote{Recall that 
(cf. the proof of Lemma~\ref{lem3.1}), for $n \! \in \! \mathbb{N}$ 
and $k \! \in \! \lbrace 1,2,\dotsc,K \rbrace$ such that $\alpha_{p_{
\mathfrak{s}}} \! := \! \alpha_{k} \! \neq \! \infty$, and for $x \! \in 
\! \mathscr{O}_{\tilde{\delta}_{\tilde{\mathcal{R}},1}(j)}(\alpha_{p_{q^{
\prime}}})$, $q^{\prime} \! \in \! \tilde{Q}_{\tilde{\mathcal{R}},1}(j)$, 
$j \! = \! 1,2,\dotsc,N$, there exists $\tilde{\mathfrak{c}}_{q^{\prime}} 
\! = \! \tilde{\mathfrak{c}}_{q^{\prime}}(n,k,z_{o}) \! > \! 0$ and 
$\mathcal{O}(1)$, in the double-scaling limit $\mathscr{N},n \! \to \! 
\infty$ such that $z_{o} \! = \! 1 \! + \! o(1)$, such that $\widetilde{V}
(x) \! \geqslant \! (1 \! + \! \tilde{\mathfrak{c}}_{q^{\prime}}) \ln 
(1 \! + \! (x \! - \! \alpha_{p_{q^{\prime}}})^{-2})$.} the fact that 
$(\tilde{R}(z))^{1/2}$ is continuous and bounded on the compact 
intervals $[\tilde{a}_{j},\tilde{b}_{j}] \supset \tilde{\Sigma}_{p,j}^{2}$, 
$(\tilde{R}(z))^{1/2} \! =_{z \downarrow \tilde{a}_{j}} \! \mathcal{O}
((z \! - \! \tilde{a}_{j})^{1/2})$ and $(\tilde{R}(z))^{1/2} \! 
=_{z \uparrow \tilde{b}_{j}} \! \mathcal{O}((\tilde{b}_{j} \! - \! z)^{1/2})$, 
and is differentiable on the open intervals $\tilde{\Sigma}_{p,j}^{2}$, 
and that $\tilde{h}_{\widetilde{V}}(z)$ is analytic for $\tilde{\Sigma}_{p,
j}^{2} \! \ni \! z$, one arrives at, for $n \! \in \! \mathbb{N}$ and 
$k \! \in \! \lbrace 1,2,\dotsc,K \rbrace$ such that $\alpha_{p_{
\mathfrak{s}}} \! := \! \alpha_{k} \! \neq \! \infty$, and for $z \! 
\in \! \tilde{\Sigma}_{p,j}^{2}$, $j \! = \! 1,2,\dotsc,N$,
\begin{align} \label{eqlem5.1b} 
\me^{n(g^{f}_{+}(z)+g^{f}_{-}(z)-\hat{\mathscr{P}}_{0}^{+}-\hat{
\mathscr{P}}_{0}^{-}-\widetilde{V}(z)-\tilde{\ell})(\chi_{\tilde{A}_{
\tilde{\mathcal{R}},1}(j)}(z)+\sum_{q \in \tilde{Q}_{\tilde{\mathcal{R}},
1}(j)} \chi_{\mathscr{O}_{\tilde{\delta}_{\tilde{\mathcal{R}},1}(j)}
(\alpha_{p_{q}})}(z))} \underset{\underset{z_{o}=1+o(1)}{\mathscr{N},
n \to \infty}}{=}& \, \mathcal{O} \left(\tilde{\mathfrak{c}}_{\tilde{
\mathcal{R}},1}^{\blacklozenge}(j) \me^{-((n-1)K+k) \tilde{
\lambda}_{\tilde{\mathcal{R}},1}(j)(z-\tilde{a}_{j}) \chi_{\tilde{A}_{
\tilde{\mathcal{R}},1}(j)}(z)} \right. \nonumber \\
\times&\left. \, \prod_{q \in \tilde{Q}_{\tilde{\mathcal{R}},1}(j)} \lvert 
z \! - \! \alpha_{p_{q}} \rvert^{((n-1)K+k) \tilde{\mathfrak{c}}_{\tilde{
\mathcal{R}},2}(j,q) \chi_{\mathscr{O}_{\tilde{\delta}_{\tilde{\mathcal{
R}},1}(j)}(\alpha_{p_{q}})}(z)} \right),
\end{align}
where $\tilde{\mathfrak{c}}_{\tilde{\mathcal{R}},1}^{\blacklozenge}(j) 
\! = \! \tilde{\mathfrak{c}}_{\tilde{\mathcal{R}},1}^{\blacklozenge}
(n,k,z_{o};j) \! =_{\underset{z_{o}=1+o(1)}{\mathscr{N},n \to \infty}} 
\! \mathcal{O}(1)$, and $\tilde{\lambda}_{\tilde{\mathcal{R}},1}(j)$ 
$(:= \! \inf_{u \in \tilde{A}_{\tilde{\mathcal{R}},1}(j)} \lvert (\tilde{R}
(u))^{1/2} \tilde{h}_{\widetilde{V}}(u) \rvert)$ and $\tilde{\mathfrak{
c}}_{\tilde{\mathcal{R}},2}(j,q)$ are as described in item~\pmb{(2)}, 
subitem~$\boldsymbol{\mathrm{(2)}_{i}}$ of the 
lemma.\footnote{If, for $j \! = \! 1,2,\dotsc,N$, $\tilde{Q}_{\tilde{
\mathcal{R}},1}(j) \! = \! \varnothing$ (that is, $\# \tilde{Q}_{\tilde{
\mathcal{R}},1}(j) \! = \! 0)$, then $\prod_{q \in \tilde{Q}_{\tilde{
\mathcal{R}},1}(j)} \lvert z \! - \! \alpha_{p_{q}} \rvert^{((n-1)K+k) 
\tilde{\mathfrak{c}}_{\tilde{\mathcal{R}},2}(j,q) \chi_{\mathscr{O}_{
\tilde{\delta}_{\tilde{\mathcal{R}},1}(j)}(\alpha_{p_{q}})}(z)} \! := \! 
1$.} For $n \! \in \! \mathbb{N}$ and $k \! \in \! \lbrace 1,2,\dotsc,
K \rbrace$ such that $\alpha_{p_{\mathfrak{s}}} \! := \! \alpha_{k} 
\! \neq \! \infty$, and for $z \! \in \! \tilde{\Sigma}_{p,j}^{2}$, $j \! 
= \! 1,2,\dotsc,N$, in order to study $\mathfrak{m}_{-}(z) \! := \! 
\mathfrak{m}(z \! - \! \mi 0)$, one notes {}from item~\pmb{(2)} of 
Lemma~\ref{lem4.5} that the $z$-dependent factors $\tilde{
\boldsymbol{\theta}}(\varepsilon_{1} \tilde{\boldsymbol{u}}(z \! - \! 
\mi 0) \! - \! \tfrac{1}{2 \pi}((n \! - \! 1)K \! + \! k) \tilde{\boldsymbol{
\Omega}} \! + \! \varepsilon_{2} \tilde{\boldsymbol{d}})$, $\tilde{
\boldsymbol{\theta}}(\varepsilon_{1} \tilde{\boldsymbol{u}}(z \! - \! 
\mi 0) \! + \! \varepsilon_{2} \tilde{\boldsymbol{d}})$, $\varepsilon_{1},
\varepsilon_{2} \! = \! \pm 1$, $\tilde{\gamma}(z \! - \! \mi 0)$, 
and $(\tilde{\gamma}(z \! - \! \mi 0))^{-1}$ must be analysed (all 
the parameters appearing in the latter formulae are defined in 
item~\pmb{(2)} of Lemma~\ref{lem4.5}). For $n \! \in \! \mathbb{N}$ 
and $k \! \in \! \lbrace 1,2,\dotsc,K \rbrace$ such that $\alpha_{
p_{\mathfrak{s}}} \! := \! \alpha_{k} \! \neq \! \infty$, and for $z \! 
\in \! \tilde{\Sigma}_{p,j}^{2}$ $(= \! \tilde{A}_{\tilde{\mathcal{R}},1}
(j) \cup \tilde{A}_{\tilde{\mathcal{R}},1}^{c}(j))$, $j \! = \! 1,2,\dotsc,
N$, one shows, via the representation for $\tilde{\boldsymbol{\theta}}
(\pmb{\cdot})$ given by Equation~\eqref{eqrmthetafin} and the elementary 
inequality $\lvert \exp (\ast) \rvert \! \leqslant \! \exp (\lvert \ast 
\rvert)$, that $\lvert \tilde{\boldsymbol{\theta}}(\varepsilon_{1} 
\tilde{\boldsymbol{u}}(z \! - \! \mi 0) \! - \! \tfrac{1}{2 \pi}((n \! - \! 
1)K \! + \! k) \tilde{\boldsymbol{\Omega}} \! + \! \varepsilon_{2} 
\tilde{\boldsymbol{d}}) \rvert \! \leqslant \! \sum_{m \in 
\mathbb{Z}^{N}} \me^{-2 \pi \varepsilon_{2}(m,\Im 
(\tilde{\boldsymbol{d}}))} \lvert \me^{\mi \pi (m,\tilde{\boldsymbol{
\tau}}m)} \rvert \lvert \me^{2 \pi \mi \varepsilon_{1}(m,\tilde{
\boldsymbol{u}}(z-\mi 0))} \rvert$: recalling that the associated 
$N \times N$ Riemann matrix of $\tilde{\boldsymbol{\beta}}$-periods, 
$\tilde{\boldsymbol{\tau}} \! = \! (\tilde{\boldsymbol{\tau}})_{j_{1},
j_{2}=1,2,\dotsc,N}$, is non-degenerate, symmetric and pure 
imaginary, and $-\mi \tilde{\boldsymbol{\tau}}$ is positive definite, 
that is, $\mi \pi (m,\tilde{\boldsymbol{\tau}}m) \! = \! \pi 
\sum_{j_{1}=1}^{N} \sum_{j_{2}=1}^{N}m_{j_{1}}(\mi (\tilde{
\boldsymbol{\tau}})_{j_{1},j_{2}})m_{j_{2}} \! < \! 0$, it follows via 
the (coarse, but sufficient!) estimate \footnote{Recall {}from the 
associated Subsection~\ref{subsub1} that, for $n \! \in \! \mathbb{N}$ 
and $k \! \in \! \lbrace 1,2,\dotsc,K \rbrace$ such that $\alpha_{
p_{\mathfrak{s}}} \! := \! \alpha_{k} \! \neq \! \infty$, $(m,\tilde{
\boldsymbol{u}}(z \! - \! \mi 0)) \! = \! \sum_{i_{1}=1}^{N}m_{i_{1}} 
\tilde{\boldsymbol{u}}_{i_{1}}(z \! - \! \mi 0)$, where $\tilde{
\boldsymbol{u}}_{i_{1}}(z \! - \! \mi 0) \! := \! \int_{\tilde{a}_{N+
1}}^{z-\mi 0} \tilde{\omega}_{i_{1}}$, with $\tilde{\omega}_{i_{1}} 
\! = \! \sum_{i_{2}=1}^{N} \tilde{c}_{i_{i}i_{2}}(\tilde{R}(x))^{-1/2}
x^{N-i_{2}} \, \md x$, and $\tilde{c} \! = \! (\tilde{c})_{j_{1},j_{2}=1,2,
\dotsc,N}$ as described in Equations~\eqref{O1} and~\eqref{O2}.}
\begin{equation*}
\lvert \tilde{\boldsymbol{u}}_{i_{1}}(z \! - \! \mi 0) \rvert \! \leqslant 
\! \left(\dfrac{1}{\min \lbrace \tilde{\delta}_{\tilde{a}_{j}},\tilde{\delta}_{
\tilde{b}_{j}} \rbrace} \right)^{N+1}   \left(\dfrac{\lvert \tilde{a}_{N+1} 
\rvert^{N} \! - \! 1}{\lvert \tilde{a}_{N+1} \rvert \! - \! 1} \right)
(\tilde{a}_{N+1} \! - \! \tilde{a}_{j} \! + \! \tilde{\delta}_{\tilde{a}_{j}}) 
\max_{j_{1}=1,2,\dotsc,N} \lbrace \lvert \tilde{c}_{i_{1}j_{j}} \rvert 
\rbrace, \quad z \! \in \! \tilde{\Sigma}_{p,j}^{2}, \quad i_{1},j \! = \! 
1,2,\dotsc,N,
\end{equation*}
and the fact that $\tilde{\boldsymbol{\theta}}(\pmb{\cdot})$ converges 
absolutely and uniformly (on compact subsets of $\mathbb{C}^{N} 
\times \mathfrak{S}_{N}$, where $\mathfrak{S}_{N}$ denotes the Siegel 
upper half space of genus $N$), that, uniformly for $z \! \in \! \tilde{
\Sigma}_{p,j}^{2}$, $j \! = \! 1,2,\dotsc,N$,
\begin{equation*}
\lvert \tilde{\boldsymbol{\theta}}(\varepsilon_{1} \tilde{\boldsymbol{u}}
(z \! - \! \mi 0) \! - \! \tfrac{1}{2 \pi}((n \! - \! 1)K \! + \! k) \tilde{
\boldsymbol{\Omega}} \! + \! \varepsilon_{2} 
\tilde{\boldsymbol{d}}) \rvert \underset{\underset{z_{o}=1+o(1)}{
\mathscr{N},n \to \infty}}{\leqslant} \mathcal{O}(\tilde{\mathfrak{c}}_{
\tilde{\theta},1}(n,k,z_{o};j)), \quad \lvert \tilde{\boldsymbol{\theta}}
(\varepsilon_{1} \tilde{\boldsymbol{u}}(z \! - \! \mi 0) \! + \! 
\varepsilon_{2} \tilde{\boldsymbol{d}}) \rvert \underset{\underset{z_{o}
=1+o(1)}{\mathscr{N},n \to \infty}}{\leqslant} \mathcal{O}(\tilde{
\mathfrak{c}}_{\tilde{\theta},2}(n,k,z_{o};j)),
\end{equation*}
where $\tilde{\mathfrak{c}}_{\tilde{\theta},i_{1}}(n,k,z_{o};j) \! 
=_{\underset{z_{o}=1+o(1)}{\mathscr{N},n \to \infty}} \! \mathcal{O}
(1)$, $i_{1} \! = \! 1,2$. Via the definition of $\tilde{\gamma}(z)$ 
given by Equation~\eqref{eqmainfin10}, one shows that, uniformly 
for $z \! \in \! \tilde{\Sigma}_{p,j}^{2}$, $j \! = \! 1,2,\dotsc,N$,
\begin{equation*}
\lvert \tilde{\gamma}(z \! - \! \mi 0) \rvert \underset{\underset{z_{o}
=1+o(1)}{\mathscr{N},n \to \infty}}{\leqslant} \mathcal{O} \left(
(\tilde{\delta}_{\tilde{a}_{j}})^{-1} \tilde{\mathfrak{c}}_{\tilde{\theta},3}
(n,k,z_{o};j) \right), \quad \quad \lvert (\tilde{\gamma}(z \! - \! \mi 
0))^{-1} \rvert \underset{\underset{z_{o}=1+o(1)}{\mathscr{N},n \to 
\infty}}{\leqslant} \mathcal{O} \left((\tilde{\delta}_{\tilde{b}_{j}})^{-1} 
\tilde{\mathfrak{c}}_{\tilde{\theta},4}(n,k,z_{o};j) \right),
\end{equation*}
where $\tilde{\mathfrak{c}}_{\tilde{\theta},i_{1}}(n,k,z_{o};j) \! 
=_{\underset{z_{o}=1+o(1)}{\mathscr{N},n \to \infty}} \! \mathcal{O}
(1)$, $i_{1} \! = \! 3,4$; hence, for $n \! \in \! \mathbb{N}$ and 
$k \! \in \! \lbrace 1,2,\dotsc,K \rbrace$ such that $\alpha_{
p_{\mathfrak{s}}} \! := \! \alpha_{k} \! \neq \! \infty$, via the 
above-derived bounds, and the (uniform) asymptotics for 
$\mathfrak{m}(z)$ as $\mathbb{C}_{-} \cap \mathscr{O}_{\tilde{
\delta}_{\tilde{\mathcal{R}},1}(j)}(\alpha_{p_{q^{\prime}}}) \! \ni \! z \! 
\to \! \alpha_{p_{q^{\prime}}}$, $q^{\prime} \! \in \! \tilde{Q}_{\tilde{
\mathcal{R}},1}(j)$, $j \! = \! 1,2,\dotsc,N$, given in item~\pmb{(2)} 
of Lemma~\ref{lem4.3}, one arrives at, after a matrix-multiplication 
argument, uniformly for $z \! \in \! \tilde{\Sigma}_{p,j}^{2} \! = \! 
\tilde{A}_{\tilde{\mathcal{R}},1}(j) \cup \tilde{A}_{\tilde{\mathcal{R}},
1}^{c}(j)$, $j \! = \! 1,2,\dotsc,N$,
\begin{align} \label{eqlem5.1c} 
\left(\mathfrak{m}_{-}(z) \left(\chi_{\tilde{A}_{\tilde{\mathcal{R}},1}(j)}
(z) \! + \! \sum_{q \in \tilde{Q}_{\tilde{\mathcal{R}},1}(j)} 
\chi_{\mathscr{O}_{\tilde{\delta}_{\tilde{\mathcal{R}},1}(j)}
(\alpha_{p_{q}})}(z) \right) \right) \sigma_{+} &\left(\mathfrak{m}_{-}
(z) \left(\chi_{\tilde{A}_{\tilde{\mathcal{R}},1}(j)}(z) \! + \! \sum_{q 
\in \tilde{Q}_{\tilde{\mathcal{R}},1}(j)} \chi_{\mathscr{O}_{\tilde{
\delta}_{\tilde{\mathcal{R}},1}(j)}(\alpha_{p_{q}})}(z) \right) \right)^{-1} 
\nonumber \\
\underset{\underset{z_{o}=1+o(1)}{\mathscr{N},n \to \infty}}{\leqslant}& 
\,  \mathcal{O}(\tilde{\mathfrak{c}}_{\tilde{\theta},5}(n,k,z_{o};j)),
\end{align}
where $(\mathrm{M}_{2}(\mathbb{C}) \! \ni)$ $\tilde{\mathfrak{c}}_{
\tilde{\theta},5}(n,k,z_{o};j) \! =_{\underset{z_{o}=1+o(1)}{\mathscr{N},
n \to \infty}} \! \mathcal{O}(1)$; finally, via Equation~\eqref{eqlem5.1a}, 
the Estimates~\eqref{eqlem5.1b} and~\eqref{eqlem5.1c}, and a 
matrix-multiplication argument, one arrives at, for $n \! \in \! 
\mathbb{N}$ and $k \! \in \! \lbrace 1,2,\dotsc,K \rbrace$ such that 
$\alpha_{p_{\mathfrak{s}}} \! := \! \alpha_{k} \! \neq \! \infty$, uniformly 
for $z \! \in \! \tilde{\Sigma}_{p,j}^{2} \! := \! (\tilde{a}_{j} \! + \! 
\tilde{\delta}_{\tilde{a}_{j}},\tilde{b}_{j} \! - \! \tilde{\delta}_{\tilde{b}_{j}})$, 
$j \! = \! 1,2,\dotsc,N$, the asymptotics, in the double-scaling limit 
$\mathscr{N},n \! \to \! \infty$ such that $z_{o} \! = \! 1 \! + \! o(1)$, 
for $\tilde{v}_{\tilde{\mathcal{R}}}(z)$ given in Equation~\eqref{eqtlvee8}.

For $n \! \in \! \mathbb{N}$ and $k \! \in \! \lbrace 1,2,\dotsc,K \rbrace$ 
such that $\alpha_{p_{\mathfrak{s}}} \! := \! \alpha_{k} \! \neq \! \infty$, 
and for $z \! \in \! \tilde{\Sigma}_{p}^{1} \! := \! (-\infty,\tilde{b}_{0} \! - 
\! \tilde{\delta}_{\tilde{b}_{0}}) \cup (\tilde{a}_{N+1} \! + \! \tilde{\delta}_{
\tilde{a}_{N+1}},+\infty)$, consider, say, and without loss of generality, the 
analysis for the case $z \! \in \! (\tilde{a}_{N+1} \! + \! \tilde{\delta}_{
\tilde{a}_{N+1}},+\infty)$ $(\subset \tilde{\Sigma}_{p}^{1})$: the analysis for 
the case $z \! \in \! (-\infty,\tilde{b}_{0} \! - \! \tilde{\delta}_{\tilde{b}_{0}})$ 
$(\subset \tilde{\Sigma}_{p}^{1})$ is analogous. For $z \! \in \! 
(\tilde{a}_{N+1} \! + \! \tilde{\delta}_{\tilde{a}_{N+1}},+\infty)$, define 
$\tilde{A}_{\tilde{\mathcal{R}},2}(+)$ and $\tilde{Q}_{\tilde{\mathcal{R}},
2}(+)$, and choose $\tilde{\varepsilon}_{\infty},\tilde{\delta}_{\tilde{
\mathcal{R}},2}(+) \! > \! 0$ as in the corresponding item~\pmb{(2)}, 
subitem~$\boldsymbol{\mathrm{(2)}_{ii}}$ of the lemma: via the definition 
of $\tilde{A}_{\tilde{\mathcal{R}},2}(+)$, it follows that the complement 
of $\tilde{A}_{\tilde{\mathcal{R}},2}(+)$ relative to $(\tilde{a}_{N+1} \! + 
\! \tilde{\delta}_{\tilde{a}_{N+1}},+\infty)$ is given by $\tilde{A}_{\tilde{
\mathcal{R}},2}^{c}(+) \! = \! \mathscr{O}_{\infty}(\alpha_{p_{\mathfrak{s}
-1}}) \cup \cup_{q \in \tilde{Q}_{\tilde{\mathcal{R}},2}(+)} \mathscr{O}_{
\tilde{\delta}_{\tilde{\mathcal{R}},2}(+)}(\alpha_{p_{q}})$; furthermore, 
since $\tilde{A}_{\tilde{\mathcal{R}},2}(+) \cap \tilde{A}_{\tilde{\mathcal{R}},2}^{c}
(+) \! = \! \varnothing$, $\mathscr{O}_{\tilde{\delta}_{\tilde{\mathcal{R}},2}(+)}
(\alpha_{p_{q^{\prime}}}) \cap \mathscr{O}_{\tilde{\delta}_{\tilde{
\mathcal{R}},2}(+)}(\alpha_{p_{q^{\prime \prime}}}) \! = \! \varnothing$ 
$\forall$ $q^{\prime} \! \neq \! q^{\prime \prime} \! \in \! \tilde{Q}_{
\tilde{\mathcal{R}},2}(+)$, and $\mathscr{O}_{\tilde{\delta}_{\tilde{
\mathcal{R}},2}(+)}(\alpha_{p_{q^{\prime}}}) \cap \mathscr{O}_{\infty}
(\alpha_{p_{\mathfrak{s}-1}}) \! = \! \varnothing$, it follows that 
$\chi_{(\tilde{a}_{N+1}+ \tilde{\delta}_{\tilde{a}_{N+1}},+\infty)}(z) 
\! = \! \chi_{\tilde{A}_{\tilde{\mathcal{R}},2}(+) \cup \tilde{A}_{
\tilde{\mathcal{R}},2}^{c}(+)}(z) \! = \! \chi_{\tilde{A}_{\tilde{
\mathcal{R}},2}(+)}(z) \! + \! \chi_{\mathscr{O}_{\infty}(\alpha_{
p_{\mathfrak{s}-1}})}(z) \! + \! \sum_{q \in \tilde{Q}_{\tilde{
\mathcal{R}},2}(+)} \chi_{\mathscr{O}_{\tilde{\delta}_{\tilde{
\mathcal{R}},2}(+)}(\alpha_{p_{q}})}(z)$, whence, via the expression 
for $\tilde{v}_{\tilde{\mathcal{R}}}(z)$ given in the corresponding 
item~\pmb{(2)} of Lemma~\ref{lem4.10}, one shows that, for $z \! 
\in \! (\tilde{a}_{N+1} \! + \! \tilde{\delta}_{\tilde{a}_{N+1}},+\infty)$,
\begin{align} \label{eqlem5.1d} 
\tilde{v}_{\tilde{\mathcal{R}}}(z) =& \, \mathrm{I} \! + \! \exp \left(n 
\left(g^{f}_{+}(z) \! + \! g^{f}_{-}(z) \! - \! \hat{\mathscr{P}}_{0}^{+} 
\! - \! \hat{\mathscr{P}}_{0}^{-} \! - \! \widetilde{V}(z) \! - \! 
\tilde{\ell} \right) \left(\chi_{\tilde{A}_{\tilde{\mathcal{R}},2}(+)}(z) 
\! + \! \chi_{\mathscr{O}_{\infty}(\alpha_{p_{\mathfrak{s}-1}})}(z) 
\! + \! \sum_{q \in \tilde{Q}_{\tilde{\mathcal{R}},2}(+)} \chi_{
\mathscr{O}_{\tilde{\delta}_{\tilde{\mathcal{R}},2}(+)}(\alpha_{p_{q}})}
(z) \right) \right) \nonumber \\
\times& \, \left(\mathfrak{m}(z) \left(\chi_{\tilde{A}_{\tilde{
\mathcal{R}},2}(+)}(z) \! + \! \chi_{\mathscr{O}_{\infty}(\alpha_{
p_{\mathfrak{s}-1}})}(z) \! + \! \sum_{q \in \tilde{Q}_{\tilde{
\mathcal{R}},2}(+)} \chi_{\mathscr{O}_{\tilde{\delta}_{\tilde{
\mathcal{R}},2}(+)}(\alpha_{p_{q}})}(z) \right) \right) \sigma_{+}
\nonumber \\
\times& \, \left(\mathfrak{m}(z) \left(\chi_{\tilde{A}_{
\tilde{\mathcal{R}},2}(+)}(z) \! + \! \chi_{\mathscr{O}_{\infty}
(\alpha_{p_{\mathfrak{s}-1}})}(z) \! + \! \sum_{q \in \tilde{Q}_{
\tilde{\mathcal{R}},2}(+)} \chi_{\mathscr{O}_{\tilde{\delta}_{
\tilde{\mathcal{R}},2}(+)}(\alpha_{p_{q}})}(z) \right) \right)^{-1}.
\end{align}
For $n \! \in \! \mathbb{N}$ and $k \! \in \! \lbrace 1,2,\dotsc,K 
\rbrace$ such that $\alpha_{p_{\mathfrak{s}}} \! := \! \alpha_{k} 
\! \neq \! \infty$, and for $z \! \in \! (\tilde{a}_{N+1} \! + \! \tilde{
\delta}_{\tilde{a}_{N+1}},+\infty)$ $(\subset \tilde{\Sigma}_{p}^{1})$, 
if $\tilde{Q}_{\tilde{\mathcal{R}},2}(+) \! = \! \varnothing$ (that is, 
$\# \tilde{Q}_{\tilde{\mathcal{R}},2}(+) \! = \! 0)$, then one uses 
the Representation~\eqref{eql3.8BB} for $g^{f}_{+}(z) \! + \! g^{f}_{-}
(z) \! - \! \hat{\mathscr{P}}_{0}^{+} \! - \! \hat{\mathscr{P}}_{0}^{-} 
\! - \! \widetilde{V}(z) \! - \! \tilde{\ell}$, or, if $\tilde{Q}_{\tilde{
\mathcal{R}},2}(+) \! \neq \! \varnothing$, then one uses the 
(equivalent) Representation~\eqref{eql3.8k} for $g^{f}_{+}(z) \! + \! 
g^{f}_{-}(z) \! - \! \hat{\mathscr{P}}_{0}^{+} \! - \! \hat{\mathscr{
P}}_{0}^{-} \! - \! \widetilde{V}(z) \! - \! \tilde{\ell}$: via the 
above-mentioned representations for $g^{f}_{+}(z) \! + \! g^{f}_{-}(z) 
\! - \! \hat{\mathscr{P}}_{0}^{+} \! - \! \hat{\mathscr{P}}_{0}^{-} \! - \! 
\widetilde{V}(z) \! - \! \tilde{\ell}$, the asymptotics (cf. the proof of 
Lemma~\ref{lem3.8}, case~\pmb{$(\mathrm{A})$}, subitem~\pmb{(3)})
\begin{gather*}
g^{f}_{+}(z) \! + \! g^{f}_{-}(z) \! - \! \hat{\mathscr{P}}_{0}^{+} 
\! - \! \hat{\mathscr{P}}_{0}^{-} \! - \! \widetilde{V}(z) \! - \! 
\tilde{\ell} \underset{\underset{q \in \tilde{Q}_{\tilde{\mathcal{R}},2}
(+)}{\mathscr{O}_{\tilde{\delta}_{\tilde{\mathcal{R}},2}(+)}(\alpha_{
p_{q}}) \ni z \to \alpha_{p_{q}}=\alpha_{k}}}{=} -\left(\widetilde{V}(z) 
\! - \! \left(\dfrac{\varkappa_{nk} \! - \! 1}{n} \right) \ln (1 \! + \! 
\lvert z \! - \! \alpha_{k} \rvert^{-2}) \right) \! + \! \mathcal{O}(1), \\
g^{f}_{+}(z) \! + \! g^{f}_{-}(z) \! - \! \hat{\mathscr{P}}_{0}^{+} 
\! - \! \hat{\mathscr{P}}_{0}^{-} \! - \! \widetilde{V}(z) \! - \! 
\tilde{\ell} \underset{\underset{q \in \tilde{Q}_{\tilde{\mathcal{R}},2}
(+)}{\mathscr{O}_{\tilde{\delta}_{\tilde{\mathcal{R}},2}(+)}(\alpha_{
p_{q}}) \ni z \to \alpha_{p_{q}} \neq \alpha_{k}}}{=} -\left(\widetilde{V}
(z) \! - \! \dfrac{\varkappa_{nk \tilde{k}_{q}}}{n} \ln (1 \! + \! \lvert 
z \! - \! \alpha_{p_{q}} \rvert^{-2}) \right) \! + \! \mathcal{O}(1), \\
g^{f}_{+}(z) \! + \! g^{f}_{-}(z) \! - \! \hat{\mathscr{P}}_{0}^{+} 
\! - \! \hat{\mathscr{P}}_{0}^{-} \! - \! \widetilde{V}(z) \! - \! 
\tilde{\ell} \underset{\mathscr{O}_{\infty}(\alpha_{p_{\mathfrak{s}-1}}) 
\ni z \to \alpha_{p_{\mathfrak{s}-1}} = \infty}{=} -\left(\widetilde{V}
(z) \! - \! \left(\dfrac{\varkappa_{nk \tilde{k}_{\mathfrak{s}-1}}^{\infty} 
\! + \! 1}{n} \right) \ln (1 \! + \! z^{2}) \right) \! + \! \mathcal{O}(1),
\end{gather*}
the conditions~\eqref{eq20}--\eqref{eq22} for regular $\widetilde{V} 
\colon \overline{\mathbb{R}} \setminus \lbrace \alpha_{1},\alpha_{2},
\dotsc,\alpha_{K} \rbrace \! \to \! \mathbb{R}$,\footnote{Recall that 
(cf. the proof of Lemma~\ref{lem3.1}), for $n \! \in \! \mathbb{N}$ 
and $k \! \in \! \lbrace 1,2,\dotsc,K \rbrace$ such that $\alpha_{p_{
\mathfrak{s}}} \! := \! \alpha_{k} \! \neq \! \infty$, for $x \! \in \! 
\mathscr{O}_{\infty}(\alpha_{p_{\mathfrak{s}-1}})$, there exists 
$\tilde{\mathfrak{c}}_{\infty} \! = \! \tilde{\mathfrak{c}}_{\infty}(n,k,
z_{o}) \! > \! 0$ and $\mathcal{O}(1)$, in the double-scaling limit 
$\mathscr{N},n \! \to \! \infty$ such that $z_{o} \! = \! 1 \! + \! o(1)$, 
such that $\widetilde{V}(x) \! \geqslant \! (1 \! + \! \tilde{\mathfrak{
c}}_{\infty}) \ln (1 \! + \! x^{2})$, and, for $x \! \in \! \mathscr{O}_{
\tilde{\delta}_{\tilde{\mathcal{R}},2}(+)}(\alpha_{p_{q^{\prime}}})$, 
$q^{\prime} \! \in \! \tilde{Q}_{\tilde{\mathcal{R}},2}(+)$, there exists 
$\tilde{\mathfrak{c}}_{q^{\prime}} \! = \! \tilde{\mathfrak{c}}_{q^{
\prime}}(n,k,z_{o}) \! > \! 0$ and $\mathcal{O}(1)$, in the 
double-scaling limit $\mathscr{N},n \! \to \! \infty$ such that $z_{o} 
\! = \! 1 \! + \! o(1)$, such that $\widetilde{V}(x) \! \geqslant \! (1 
\! + \! \tilde{\mathfrak{c}}_{q^{\prime}}) \ln (1 \! + \! (x \! - \! 
\alpha_{p_{q^{\prime}}})^{-2})$.} the fact that $(\tilde{R}(z))^{1/2}$ 
is differentiable and bounded for $z \! \in \! \tilde{A}_{\tilde{\mathcal{R}},
2}(+)$, and that $\tilde{h}_{\widetilde{V}}(z)$ is analytic for $z \! \in \! 
\tilde{A}_{\tilde{\mathcal{R}},2}(+)$, one arrives at, for $n \! \in \! 
\mathbb{N}$ and $k \! \in \! \lbrace 1,2,\dotsc,K \rbrace$ such that 
$\alpha_{p_{\mathfrak{s}}} \! := \! \alpha_{k} \! \neq \! \infty$, and 
for $z \! \in \! (\tilde{a}_{N+1} \! + \! \tilde{\delta}_{\tilde{a}_{N+1}},
+\infty)$,
\begin{align} \label{eqlem5.1e} 
& \, \me^{n(g^{f}_{+}(z)+g^{f}_{-}(z)-\hat{\mathscr{P}}_{0}^{+}-
\hat{\mathscr{P}}_{0}^{-}-\widetilde{V}(z)-\tilde{\ell})(\chi_{\tilde{
A}_{\tilde{\mathcal{R}},2}(+)}(z)+ \chi_{\mathscr{O}_{\infty}(\alpha_{
p_{\mathfrak{s}-1}})}(z)+\sum_{q \in \tilde{Q}_{\tilde{\mathcal{R}},
2}(+)} \chi_{\mathscr{O}_{\tilde{\delta}_{\tilde{\mathcal{R}},2}(+)}
(\alpha_{p_{q}})}(z))} \underset{\underset{z_{o}=1+o(1)}{\mathscr{N},
n \to \infty}}{=} \mathcal{O} \left(\tilde{\mathfrak{c}}_{\tilde{
\mathcal{R}},2}^{\blacklozenge}(+) \me^{-((n-1)K+k) \tilde{
\lambda}_{\tilde{\mathcal{R}},2}(+)(z-\tilde{a}_{N+1}) \chi_{
\tilde{A}_{\tilde{\mathcal{R}},2}(+)}(z)} \right. \nonumber \\
&\left. \times \, \me^{-((n-1)K+k) \tilde{\lambda}_{\tilde{
\mathcal{R}},3}(+) \ln (\lvert z \rvert) \chi_{\mathscr{O}_{\infty}
(\alpha_{ p_{\mathfrak{s}-1}})}(z)} \prod_{q \in \tilde{Q}_{\tilde{
\mathcal{R}},2}(+)} \lvert z \! - \! \alpha_{p_{q}} \rvert^{((n-1)K+k) 
\tilde{\mathfrak{c}}_{\tilde{\mathcal{R}},3}(q,+) \chi_{\mathscr{O}_{
\tilde{\delta}_{\tilde{\mathcal{R}},2}(+)}(\alpha_{p_{q}})}(z)} \right),
\end{align}
where $\tilde{\mathfrak{c}}_{\tilde{\mathcal{R}},2}^{\blacklozenge}(+) 
\! = \! \tilde{\mathfrak{c}}_{\tilde{\mathcal{R}},2}^{\blacklozenge}
(n,k,z_{o};+) \! =_{\underset{z_{o}=1+o(1)}{\mathscr{N},n \to \infty}} 
\! \mathcal{O}(1)$, and $\tilde{\lambda}_{\tilde{\mathcal{R}},2}(+)$ 
$(:= \! \inf_{u \in \tilde{A}_{\tilde{\mathcal{R}},2}(+)} \lvert (\tilde{R}
(u))^{1/2} \tilde{h}_{\widetilde{V}}(u) \rvert)$, $\tilde{\lambda}_{
\tilde{\mathcal{R}},3}(+)$, and $\tilde{\mathfrak{c}}_{\tilde{
\mathcal{R}},3}(q,+)$ are as described in the corresponding 
item~\pmb{(2)}, subitem~$\boldsymbol{\mathrm{(2)}_{ii}}$ of the 
lemma.\footnote{If $\tilde{Q}_{\tilde{\mathcal{R}},2}(+) \! = \! 
\varnothing$ (that is, $\# \tilde{Q}_{\tilde{\mathcal{R}},2}(+) \! = \! 
0)$, then $\prod_{q \in \tilde{Q}_{\tilde{\mathcal{R}},2}(+)} \lvert z 
\! - \! \alpha_{p_{q}} \rvert^{((n-1)K+k) \tilde{\mathfrak{c}}_{\tilde{
\mathcal{R}},3}(q,+) \chi_{\mathscr{O}_{\tilde{\delta}_{\tilde{
\mathcal{R}},2}(+)}(\alpha_{p_{q}})}(z)} \! := \! 1$.} For $n \! \in \! 
\mathbb{N}$ and $k \! \in \! \lbrace 1,2,\dotsc,K \rbrace$ such 
that $\alpha_{p_{\mathfrak{s}}} \! := \! \alpha_{k} \! \neq \! \infty$, and 
for $z \! \in \! (\tilde{a}_{N+1} \! + \! \tilde{\delta}_{\tilde{a}_{N+1}},
+\infty)$ $(\subset \tilde{\Sigma}_{p}^{1})$, in order to study 
$\mathfrak{m}(z)$, one notes {}from item~\pmb{(2)} of 
Lemma~\ref{lem4.5} that the $z$-dependent factors 
$\tilde{\boldsymbol{\theta}}(\varepsilon_{1} \tilde{\boldsymbol{u}}
(z) \! - \! \tfrac{1}{2 \pi}((n \! - \! 1)K \! + \! k) \tilde{\boldsymbol{
\Omega}} \! + \! \varepsilon_{2} \tilde{\boldsymbol{d}})$, 
$\tilde{\boldsymbol{\theta}}(\varepsilon_{1} \tilde{\boldsymbol{u}}
(z) \! + \! \varepsilon_{2} \tilde{\boldsymbol{d}})$, $\varepsilon_{1},
\varepsilon_{2} \! = \! \pm 1$, $\tilde{\gamma}(z)$, and 
$(\tilde{\gamma}(z))^{-1}$ must be analysed (all the parameters 
appearing in the latter formulae are defined item~\pmb{(2)} of 
Lemma~\ref{lem4.5}). For $n \! \in \! \mathbb{N}$ and $k \! \in \! 
\lbrace 1,2,\dotsc,K \rbrace$ such that $\alpha_{p_{\mathfrak{s}}} \! 
:= \! \alpha_{k} \! \neq \! \infty$, and for $z \! \in \! (\tilde{a}_{N+1} 
\! + \! \tilde{\delta}_{\tilde{a}_{N+1}},+\infty)$, one shows, via the 
representation for $\tilde{\boldsymbol{\theta}}(\pmb{\cdot})$ given by 
Equation~\eqref{eqrmthetafin} and the elementary inequality $\lvert 
\exp (\ast) \rvert \! \leqslant \! \exp (\lvert \ast \rvert)$, that $\lvert 
\tilde{\boldsymbol{\theta}}(\varepsilon_{1} \tilde{\boldsymbol{u}}(z) \! 
- \! \tfrac{1}{2 \pi}((n \! - \! 1)K \! + \! k) \tilde{\boldsymbol{\Omega}} 
\! + \! \varepsilon_{2} \tilde{\boldsymbol{d}}) \rvert \! \leqslant \! 
\sum_{m \in \mathbb{Z}^{N}} \me^{-2 \pi \varepsilon_{2}(m,\Im 
(\tilde{\boldsymbol{d}}))} \lvert \me^{\mi \pi (m,\tilde{\boldsymbol{\tau}}m)} 
\rvert \lvert \me^{2 \pi \mi \varepsilon_{1}(m,\tilde{\boldsymbol{u}}(z))} 
\rvert$: recalling that the associated $N \times N$ Riemann matrix of 
$\tilde{\boldsymbol{\beta}}$-periods, $\tilde{\boldsymbol{\tau}} \! = 
\! (\tilde{\boldsymbol{\tau}})_{i,j=1,2,\dotsc,N}$, is non-degenerate, 
symmetric and pure imaginary, and $-\mi \tilde{\boldsymbol{\tau}}$ is 
positive definite, it follows via the estimate {}\footnote{Recall {}from the 
associated Subsection~\ref{subsub1} that, for $n \! \in \! \mathbb{N}$ 
and $k \! \in \! \lbrace 1,2,\dotsc,K \rbrace$ such that 
$\alpha_{p_{\mathfrak{s}}} \! := \! \alpha_{k} \! \neq \! \infty$, 
$(m,\tilde{\boldsymbol{u}}(z)) \! = \! \sum_{i=1}^{N}m_{i} 
\tilde{\boldsymbol{u}}_{i}(z)$, where $\tilde{\boldsymbol{u}}_{i}(z) \! := \! 
(\int_{\tilde{a}_{N+1}}^{\tilde{a}_{N+1}+ \tilde{\delta}_{\tilde{a}_{N+1}}} 
\! + \! \int_{\tilde{a}_{N+1}+ \tilde{\delta}_{\tilde{a}_{N+1}}}^{z}) 
\tilde{\omega}_{i}$, with $\tilde{\omega}_{i} \! = \! \sum_{j=1}^{N} 
\tilde{c}_{ij}(\tilde{R}(x))^{-1/2}x^{N-j} \, \md x$, and $\tilde{c} \! = \! 
(\tilde{c})_{i,j=1,2,\dotsc,N}$ as described in Equations~\eqref{O1} 
and~\eqref{O2}.}
\begin{equation*}
\lvert \tilde{\boldsymbol{u}}_{i}(z) \rvert \! \leqslant \! \left(
\dfrac{\lvert \tilde{\varepsilon}_{\infty}^{-1} \! - \! (\tilde{a}_{N+1} 
\! + \! \tilde{\delta}_{\tilde{a}_{N+1}}) \rvert \tilde{\varepsilon}_{
\infty}(1 \! - \! \tilde{\varepsilon}_{\infty}^{N})}{(\inf_{x \in \tilde{A}_{
\tilde{\mathcal{R}},2}(+)} \lvert (\tilde{R}(x))^{1/2} \rvert) \tilde{
\varepsilon}_{\infty}^{N}(1 \! - \! \tilde{\varepsilon}_{\infty})} \! + \! 
\dfrac{2(\tilde{\delta}_{\tilde{a}_{N+1}})^{1/2}(\lvert \tilde{a}_{N+1} 
\rvert^{N} \! - \! 1)(1 \! + \! \mathcal{O}(\tilde{\delta}_{\tilde{a}_{N
+1}}))}{\tilde{\eta}_{\tilde{a}_{N+1}}(\lvert \tilde{a}_{N+1} \rvert \! - 
\! 1)} \right) \max_{j=1,2,\dotsc,N} \lbrace \lvert \tilde{c}_{ij} \rvert 
\rbrace, \quad z \! \in \! \tilde{A}_{\tilde{\mathcal{R}},2}(+),
\end{equation*}
$i \! = \! 1,2,\dotsc,N$, with $\tilde{\eta}_{\tilde{a}_{N+1}}$ defined by 
Equation~\eqref{eqmainfin59}, and the fact that $\tilde{\boldsymbol{
\theta}}(\pmb{\cdot})$ converges absolutely and uniformly, that, 
uniformly for $z \! \in \! \tilde{A}_{\tilde{\mathcal{R}},2}(+)$,
\begin{equation*}
\lvert \tilde{\boldsymbol{\theta}}(\varepsilon_{1} \tilde{\boldsymbol{u}}
(z) \! - \! \tfrac{1}{2 \pi}((n \! - \! 1)K \! + \! k) \tilde{\boldsymbol{
\Omega}} \! + \! \varepsilon_{2} \tilde{\boldsymbol{d}}) \rvert 
\underset{\underset{z_{o}=1+o(1)}{\mathscr{N},n \to \infty}}{\leqslant} 
\mathcal{O}(\tilde{\mathfrak{c}}_{\tilde{\theta},1}^{\lozenge}(n,k,z_{o};
+)), \qquad \lvert \tilde{\boldsymbol{\theta}}(\varepsilon_{1} \tilde{
\boldsymbol{u}}(z) \! + \! \varepsilon_{2} \tilde{\boldsymbol{d}}) \rvert 
\underset{\underset{z_{o}=1+o(1)}{\mathscr{N},n \to \infty}}{\leqslant} 
\mathcal{O}(\tilde{\mathfrak{c}}_{\tilde{\theta},2}^{\lozenge}(n,k,z_{o};+)),
\end{equation*}
where $\tilde{\mathfrak{c}}_{\tilde{\theta},i}^{\lozenge}(n,k,z_{o};+) \! 
=_{\underset{z_{o}=1+o(1)}{\mathscr{N},n \to \infty}} \! \mathcal{O}
(1)$, $i \! = \! 1,2$. Via the definition of $\tilde{\gamma}(z)$ given by 
Equation~\eqref{eqmainfin10}, one shows that, uniformly for $z \! \in 
\! \tilde{A}_{\tilde{\mathcal{R}},2}(+)$,
\begin{equation*}
\lvert \tilde{\gamma}(z) \rvert \underset{\underset{z_{o}=1+o(1)}{
\mathscr{N},n \to \infty}}{\leqslant} \mathcal{O} \left((\tilde{\delta}_{
\tilde{a}_{N+1}})^{-1} \tilde{\mathfrak{c}}_{\tilde{\theta},3}^{\lozenge}
(n,k,z_{o};+) \right), \quad \quad \lvert (\tilde{\gamma}(z))^{-1} 
\rvert \underset{\underset{z_{o}=1+o(1)}{\mathscr{N},n \to 
\infty}}{\leqslant} \mathcal{O}(\tilde{\mathfrak{c}}_{\tilde{\theta},
4}^{\lozenge}(n,k,z_{o};+)),
\end{equation*}
where $\tilde{\mathfrak{c}}_{\tilde{\theta},i}^{\lozenge}(n,k,
z_{o};+) \! =_{\underset{z_{o}=1+o(1)}{\mathscr{N},n \to \infty}} \! 
\mathcal{O}(1)$, $i \! = \! 3,4$; hence, for $n \! \in \! \mathbb{N}$ 
and $k \! \in \! \lbrace 1,2,\dotsc,K \rbrace$ such that $\alpha_{
p_{\mathfrak{s}}} \! := \! \alpha_{k} \! \neq \! \infty$, via the 
above-derived bounds, and the (uniform) asymptotics for 
$\mathfrak{m}(z)$ as $\mathbb{C}_{\pm} \cap \mathscr{O}_{
\tilde{\delta}_{\tilde{\mathcal{R}},2}(+)}(\alpha_{p_{q^{\prime}}}) 
\! \ni \! z \! \to \! \alpha_{p_{q^{\prime}}}$, $q^{\prime} \! \in \! 
\tilde{Q}_{\tilde{\mathcal{R}},2}(+)$, and as $\overline{\mathbb{C}}_{
\pm} \cap \mathscr{O}_{\infty}(\alpha_{p_{\mathfrak{s}-1}}) \! \ni 
\! z \! \to \! \alpha_{p_{\mathfrak{s}-1}} \! = \! \infty$ given in 
item~\pmb{(2)} of Lemma~\ref{lem4.3}, one arrives at, after 
a matrix-multiplication argument, uniformly for $z \! \in \! 
(\tilde{a}_{N+1} \! + \! \tilde{\delta}_{\tilde{a}_{N+1}},+\infty) \! 
= \! \tilde{A}_{\tilde{\mathcal{R}},2}(+) \cup \tilde{A}_{\tilde{
\mathcal{R}},2}^{c}(+)$ $(\subset \tilde{\Sigma}_{p}^{1})$,
\begin{align} \label{eqlem5.1f} 
&\, \left(\mathfrak{m}(z) \left(\chi_{\tilde{A}_{\tilde{\mathcal{R}},
2}(+)}(z) \! + \! \chi_{\mathscr{O}_{\infty}(\alpha_{p_{\mathfrak{s}
-1}})}(z) \! + \! \sum_{q \in \tilde{Q}_{\tilde{\mathcal{R}},2}(+)} 
\chi_{\mathscr{O}_{\tilde{\delta}_{\tilde{\mathcal{R}},2}(+)}
(\alpha_{p_{q}})}(z) \right) \right) \sigma_{+} \nonumber \\
\times& \, \left(\mathfrak{m}(z) \left(\chi_{\tilde{A}_{
\tilde{\mathcal{R}},2}(+)}(z) \! + \! \chi_{\mathscr{O}_{\infty}
(\alpha_{p_{\mathfrak{s}-1}})}(z) \! + \! \sum_{q \in \tilde{Q}_{
\tilde{\mathcal{R}},2}(+)} \chi_{\mathscr{O}_{\tilde{\delta}_{
\tilde{\mathcal{R}},2}(+)}(\alpha_{p_{q}})}(z) \right) \right)^{-1} 
\underset{\underset{z_{o}=1+o(1)}{\mathscr{N},n \to \infty}}{\leqslant} 
\mathcal{O}(\tilde{\mathfrak{c}}_{\tilde{\theta},5}^{\lozenge}
(n,k,z_{o};+)),
\end{align}
where $(\mathrm{M}_{2}(\mathbb{C}) \! \ni)$ $\tilde{\mathfrak{c}}_{
\tilde{\theta},5}^{\lozenge}(n,k,z_{o};+) \! =_{\underset{z_{o}=
1+o(1)}{\mathscr{N},n \to \infty}} \! \mathcal{O}(1)$; finally, via 
Equation~\eqref{eqlem5.1d}, the Estimates~\eqref{eqlem5.1e} 
and~\eqref{eqlem5.1f}, and a matrix-multiplication argument, one 
arrives at, for $n \! \in \! \mathbb{N}$ and $k \! \in \! \lbrace 1,2,
\dotsc,K \rbrace$ such that $\alpha_{p_{\mathfrak{s}}} \! := \! 
\alpha_{k} \! \neq \! \infty$, uniformly for $z \! \in \! (\tilde{a}_{N+1} 
\! + \! \tilde{\delta}_{\tilde{a}_{N+1}},+\infty)$ $(\subset 
\tilde{\Sigma}_{p}^{1})$, the asymptotics, in the double-scaling limit 
$\mathscr{N},n \! \to \! \infty$ such that $z_{o} \! = \! 1 \! + \! o(1)$, 
for $\tilde{v}_{\tilde{\mathcal{R}}}(z)$ given in Equation~\eqref{eqtlvee9}. 
For $n \! \in \! \mathbb{N}$ and $k \! \in \! \lbrace 1,2,\dotsc,K \rbrace$ 
such that $\alpha_{p_{\mathfrak{s}}} \! := \! \alpha_{k} \! \neq \! \infty$, 
the case $z \! \in \! (-\infty,\tilde{b}_{0} \! - \! \tilde{\delta}_{\tilde{b}_{0}})$ 
$(\subset \tilde{\Sigma}_{p}^{1})$ is analysed analogously, and leads 
to the asymptotics, in the double-scaling limit $\mathscr{N},n \! \to \! 
\infty$ such that $z_{o} \! = \! 1 \! + \! o(1)$, for $\tilde{v}_{
\tilde{\mathcal{R}}}(z)$ given in Equation~\eqref{eqtlvee10}.

For $n \! \in \! \mathbb{N}$ and $k \! \in \! \lbrace 1,2,\dotsc,K 
\rbrace$ such that $\alpha_{p_{\mathfrak{s}}} \! := \! \alpha_{k} \! 
\neq \! \infty$, and for $z \! \in \! \tilde{\Sigma}^{3}_{p,j} \cup 
\tilde{\Sigma}^{4}_{p,j}$, $j \! = \! 1,2,\dotsc,N \! + \! 1$, where $\tilde{
\Sigma}^{3}_{p,j} \! := \! \tilde{J}_{j}^{\smallfrown} \setminus (\tilde{J}_{j}^{
\smallfrown} \cap (\tilde{\mathbb{U}}_{\tilde{\delta}_{\tilde{b}_{j-1}}} \cup 
\tilde{\mathbb{U}}_{\tilde{\delta}_{\tilde{a}_{j}}}))$ $(\subset \mathbb{C}_{
+})$, and $\tilde{\Sigma}^{4}_{p,j} \! := \! \tilde{J}_{j}^{\smallsmile} 
\setminus (\tilde{J}_{j}^{\smallsmile} \cap (\tilde{\mathbb{U}}_{\tilde{
\delta}_{\tilde{b}_{j-1}}} \cup \tilde{\mathbb{U}}_{\tilde{\delta}_{
\tilde{a}_{j}}}))$ $(\subset \mathbb{C}_{-})$, consider, say, and without loss 
of generality, the case $z \! \in \! \tilde{\Sigma}^{3}_{p,j}$; the analysis for 
the case $z \! \in \! \tilde{\Sigma}^{4}_{p,j}$ is, \emph{mutatis mutandis}, 
analogous. For $n \! \in \! \mathbb{N}$ and $k \! \in \! \lbrace 1,2,
\dotsc,K \rbrace$ such that $\alpha_{p_{\mathfrak{s}}} \! := \! \alpha_{k} 
\! \neq \! \infty$, and for $z \! \in \! \tilde{\Sigma}^{3}_{p,j}$, $j \! 
= \! 1,2,\dotsc,N \! + \! 1$, recall the expression for $\tilde{v}_{\tilde{
\mathcal{R}}}(z)$ given in the corresponding item~\pmb{(2)} of 
Lemma~\ref{lem4.10}:
\begin{equation} \label{eqlem5.1g} 
\tilde{v}_{\tilde{\mathcal{R}}}(z) \! = \! \mathrm{I} \! + \! \me^{-2 \pi 
\mi ((n-1)K+k) \int_{z}^{\tilde{a}_{N+1}} \psi_{\widetilde{V}}^{f}(\xi) \, 
\md \xi} \mathfrak{m}(z) \sigma_{-}(\mathfrak{m}(z))^{-1}, \quad z \! 
\in \! \tilde{\Sigma}^{3}_{p,j}, \quad j \! = \! 1,2,\dotsc,N \! + \! 1,
\end{equation}
where $\psi_{\widetilde{V}}^{f}(\pmb{\cdot})$ is characterised 
completely in item~\pmb{(2)} of Lemma~\ref{lem3.7} (cf. 
Equations~\eqref{eql3.7g}--\eqref{eql3.7l}), $\mathfrak{m}(z)$ is given 
in item~\pmb{(2)} of Lemma~\ref{lem4.5}, and, {}from the corresponding 
item of Lemma~\ref{lem4.1}, $\Re (\mi \int_{z}^{\tilde{a}_{N+1}} 
\psi_{\widetilde{V}}^{f}(\xi) \, \md \xi) \! > \! 0$ for $z \! \in \! \tilde{
\Sigma}^{3}_{p,j}$, $j \! = \! 1,2,\dotsc,N \! + \! 1$. For $n \! \in \! 
\mathbb{N}$ and $k \! \in \! \lbrace 1,2,\dotsc,K \rbrace$ such that 
$\alpha_{p_{\mathfrak{s}}} \! := \! \alpha_{k} \! \neq \! \infty$, and for 
$z \! \in \! \tilde{\Sigma}^{3}_{p,j}$, $j \! = \! 1,2,\dotsc,N \! + \! 1$, 
let $\tilde{\mathbb{U}}_{j}^{\ast} \! := \! \lbrace \mathstrut z \! \in \! 
\mathbb{C}; \, \inf_{\tilde{q} \in (\tilde{b}_{j-1},\tilde{a}_{j})} \lvert z \! 
- \! \tilde{q} \rvert \! < \! \tilde{r}_{j} \rbrace$, $j \! = \! 1,2,\dotsc,
N \! + \! 1$, where $(0,1) \! \ni \! \tilde{r}_{j}$ is chosen so that 
$\tilde{\mathbb{U}}_{j}^{\ast}$ is in the domain of analyticity of 
$\widetilde{V}$, $\tilde{\mathbb{U}}_{i}^{\ast} \cap \tilde{\mathbb{
U}}_{j}^{\ast} \! = \! \varnothing$ $\forall$ $i \! \neq \! j \! \in \! 
\lbrace 1,2,\dotsc,N \! + \! 1 \rbrace$, $\tilde{\mathbb{U}}_{j}^{\ast} 
\cap \lbrace \alpha_{p_{1}},\dotsc,\alpha_{p_{\mathfrak{s}-2}},
\alpha_{k} \rbrace \! = \! \varnothing$,\footnote{Of course, $\tilde{
\mathbb{U}}_{j}^{\ast} \cap \lbrace \alpha_{p_{\mathfrak{s}-1}} \! 
= \! \infty \rbrace \! = \! \varnothing$.} and $\operatorname{int}
(\tilde{\mathbb{U}}_{j}^{\ast}) \cap \operatorname{ext}
(\tilde{\mathbb{U}}_{\tilde{\delta}_{\tilde{b}_{j-1}}} \cup \tilde{
\mathbb{U}}_{\tilde{\delta}_{\tilde{a}_{j}}}) \supset \tilde{\Sigma}^{3}_{p,j} 
\cup \tilde{\Sigma}^{4}_{p,j}$.\footnote{This also implies that, for $j \! 
= \! 1,2,\dotsc,N \! + \! 1$, $\operatorname{int}(\tilde{\mathbb{U}}_{j}^{
\ast}) \cap \operatorname{ext}(\tilde{\mathbb{U}}_{\tilde{\delta}_{\tilde{
b}_{j-1}}} \cup \tilde{\mathbb{U}}_{\tilde{\delta}_{\tilde{a}_{j}}}) \cap 
\mathbb{C}_{+} \supset \tilde{\Sigma}^{3}_{p,j}$, $\tilde{\mathbb{U}}_{j}^{
\ast} \supset \tilde{\mathbb{U}}_{\tilde{\delta}_{\tilde{b}_{j-1}}} \cup 
\tilde{\mathbb{U}}_{\tilde{\delta}_{\tilde{a}_{j}}}$, and $\tilde{\mathbb{
U}}_{j}^{\ast} \supset \tilde{\Sigma}^{3}_{p,j}$.} Let $\tilde{\Sigma}^{3}_{
p,j} \cap \partial \tilde{\mathbb{U}}_{\tilde{\delta}_{\tilde{a}_{j}}} \! := \! 
\lbrace \tilde{p}_{+}(j) \rbrace$, $j \! = \! 1,2,\dotsc,N \! + \! 1$, where 
$\tilde{p}_{+}(j) \! = \! (\tilde{a}_{j} \! + \! \tilde{\delta}_{\tilde{a}_{j}} 
\cos \sigma_{\tilde{a}_{j}}^{+},\tilde{\delta}_{\tilde{a}_{j}} \sin \sigma_{
\tilde{a}_{j}}^{+})$, $\sigma_{\tilde{a}_{j}}^{+} \! \in \! (2 \pi/3,\pi)$. 
Since, for $n \! \in \! \mathbb{N}$ and $k \! \in \! \lbrace 1,2,\dotsc,K 
\rbrace$ such that $\alpha_{p_{\mathfrak{s}}} \! := \! \alpha_{k} \! \neq 
\! \infty$, the oriented skeleton $\tilde{\Sigma}^{3}_{p,j}$, $j \! = \! 1,
2,\dotsc,N \! + \! 1$, is continuous, it can be homotopically deformed 
into, say,\footnote{Other geometries are possible.} an `elliptic lens' 
with parametrisation $x_{j}(\theta) \! = \! \tfrac{1}{2}(\tilde{a}_{j} \! + \! 
\tilde{b}_{j-1}) \! + \! \tfrac{1}{2}(\tilde{a}_{j} \! - \! \tilde{b}_{j-1}) \cos 
\theta$, $y_{j}(\theta) \! = \! \tilde{\eta}_{j} \sin \theta$, $\theta_{0}^{
\tilde{a}}(j) \! \leqslant \! \theta \! \leqslant \! \pi \! - \! \theta_{0}^{
\tilde{b}}(j)$, $j \! = \! 1,2,\dotsc,N \! + \! 1$, where, via an application 
of the Law of Sines and the Law of Cosines,\footnote{Note: $y \! = \! 
\sin^{-1}(x) \Leftrightarrow \sin y \! = \! x$, with $-\tfrac{\pi}{2} \! 
\leqslant \! y \! \leqslant \! \tfrac{\pi}{2}$, where $-1 \! \leqslant \! x 
\! \leqslant \! 1$.}
\begin{gather*}
\theta_{0}^{\tilde{a}}(j) \! = \! \sin^{-1} \left(\dfrac{\tilde{\delta}_{
\tilde{a}_{j}} \sin \sigma^{+}_{\tilde{a}_{j}}}{\tilde{l}_{\tilde{a}_{j}}} \right) 
\! \in \! (0,\pi/2), \quad \quad \theta_{0}^{\tilde{b}}(j) \! = \! \sin^{-1} 
\left(\dfrac{\tilde{\delta}_{\tilde{b}_{j-1}} \sin \sigma^{+}_{\tilde{b}_{j-
1}}}{\tilde{l}_{\tilde{b}_{j-1}}} \right) \! \in \! (0,\pi/2), \\
\tilde{l}_{\tilde{a}_{j}} \! = \! \left((\tilde{\delta}_{\tilde{a}_{j}})^{2} \! + 
\! \left(\tfrac{1}{2}(\tilde{a}_{j} \! - \! \tilde{b}_{j-1}) \right)^{2} \! + 
\! \tilde{\delta}_{\tilde{a}_{j}}(\tilde{a}_{j} \! - \! \tilde{b}_{j-1}) \cos 
\sigma_{\tilde{a}_{j}}^{+} \right)^{1/2} \quad (> \! 0), \\
\tilde{l}_{\tilde{b}_{j-1}} \! = \! \left((\tilde{\delta}_{\tilde{b}_{j-1}})^{2} 
\! + \! \left(\tfrac{1}{2}(\tilde{a}_{j} \! - \! \tilde{b}_{j-1}) \right)^{2} \! 
+ \! \tilde{\delta}_{\tilde{b}_{j-1}}(\tilde{a}_{j} \! - \! \tilde{b}_{j-1}) 
\cos \sigma_{\tilde{b}_{j-1}}^{+} \right)^{1/2} \quad (> \! 0),
\end{gather*}
$\sigma_{\tilde{b}_{j-1}}^{+} \! \in \! (2 \pi/3,\pi)$, and $0 \! < \! \min 
\lbrace \tilde{\delta}_{\tilde{b}_{j-1}} \sin \sigma_{\tilde{b}_{j-1}}^{+},
\tilde{\delta}_{\tilde{a}_{j}} \sin \sigma_{\tilde{a}_{j}}^{+} \rbrace \! < \! 
\tilde{\eta}_{j} \! < \! \tilde{r}_{j} \! < \! 1$.\footnote{Other possible 
choices for $(0,1) \! \ni \! \tilde{\eta}_{j}$, $j \! = \! 1,2,\dotsc,N \! + 
\! 1$, are (by no means an exhaustive list!): $0 \! < \! \tilde{\eta}_{j} \! 
< \! \min \lbrace \tilde{\delta}_{\tilde{b}_{j-1}},\tilde{\delta}_{\tilde{a}_{j}} 
\rbrace \! < \! \tilde{r}_{j} \! < \! 1$; $0 \! < \! \min \lbrace \tilde{
\delta}_{\tilde{b}_{j-1}},\tilde{\delta}_{\tilde{a}_{j}} \rbrace \! < \! \tilde{
\eta}_{j} \! < \! \max \lbrace \tilde{\delta}_{\tilde{b}_{j-1}},\tilde{\delta}_{
\tilde{a}_{j}} \rbrace \! < \! \tilde{r}_{j} \! < \! 1$; or $0 \! < \! \max 
\lbrace \tilde{\delta}_{\tilde{b}_{j-1}},\tilde{\delta}_{\tilde{a}_{j}} \rbrace 
\! < \! \tilde{\eta}_{j} \! < \! \tilde{r}_{j} \! < \! 1$. Note, also, that the 
counter-clockwise sense of traversal of the homotopically deformed 
elliptic lenses is opposite to the clockwise sense of traversal of the 
oriented and continuous skeletons $\tilde{\Sigma}_{p,j}^{3}$, $j \! = \! 
1,2,\dotsc,N \! + \! 1$; but, as all of the ensuing line integrals are 
taken with respect to arc-length parametrisation, and since, as is well 
known, line integrals with respect to arc-length parametrisation are 
invariant under orientation, this doesn't affect any of the results.} For 
$n \! \in \! \mathbb{N}$ and $k \! \in \! \lbrace 1,2,\dotsc,K \rbrace$ 
such that $\alpha_{p_{\mathfrak{s}}} \! := \! \alpha_{k} \! \neq \! \infty$, 
and for $z \! \in \! \tilde{\Sigma}_{p,j}^{3}$, $j \! = \! 1,2,\dotsc,N \! 
+ \! 1$, write
\begin{equation} \label{eqlem5.1h} 
\exp \left(-2 \pi \mi ((n \! - \! 1)K \! + \! k) \int_{z}^{\tilde{a}_{N+1}} 
\psi_{\widetilde{V}}^{f}(\xi) \, \md \xi \right) \! = \! \exp \left(-2 \pi 
\mi ((n \! - \! 1)K \! + \! k) \left(\int_{z}^{\tilde{p}_{+}(j)} \! + \! \int_{
\tilde{C}_{j}^{\ast}} \! + \! \int_{\tilde{a}_{j}+ \tilde{\delta}_{\tilde{a}_{j}}}^{
\tilde{a}_{N+1}} \right) \psi_{\widetilde{V}}^{f}(\xi) \, \md \xi \right),
\end{equation}
where $\tilde{C}_{j}^{\ast} \! := \! \lbrace \mathstrut z(\tau) \! = \! 
x(\tau) \! + \! \mi y(\tau); \, x(\tau) \! = \! \tilde{a}_{j} \! + \! 
\tilde{\delta}_{\tilde{a}_{j}} \cos (\sigma_{\tilde{a}_{j}}^{+} \! - \! \tau), \, 
y(\tau) \! = \! \tilde{\delta}_{\tilde{a}_{j}} \sin (\sigma_{\tilde{a}_{j}}^{+} \! 
- \! \tau), \, 0 \! \leqslant \! \tau \! \leqslant \! \sigma_{\tilde{a}_{j}}^{+} 
\rbrace$, $j \! = \! 1,2,\dotsc,N \! + \! 1$: via 
Equations~\eqref{eql3.7g}--\eqref{eql3.7l} and repeated application 
of the Maximum-Length (ML) Theorem, one shows that, for $z \! \in \! 
\tilde{\Sigma}_{p,j}^{3}$, $j \! = \! 1,2,\dotsc,N \! + \! 1$,
\begin{gather}
\left\vert \me^{-2 \pi \mi ((n-1)K+k) \int_{\tilde{a}_{j}+\tilde{\delta}_{
\tilde{a}_{j}}}^{\tilde{a}_{N+1}} \psi_{\widetilde{V}}^{f}(\xi) \, \md \xi} 
\right\vert \! = \! \left\vert \me^{-2 \pi \mi ((n-1)K+k) \sum_{m=j+
1}^{N+1} \int_{\tilde{b}_{m-1}}^{\tilde{a}_{m}}(2 \pi \mi)^{-1}(\tilde{R}
(\xi))^{1/2}_{+} \tilde{h}_{\widetilde{V}}(\xi) \, \md \xi} \right\vert \! 
= \! 1, \label{eqlem5.1i} \\
\left\vert \me^{-2 \pi \mi ((n-1)K+k) \int_{\tilde{C}_{j}^{\ast}} \psi_{
\widetilde{V}}^{f}(\xi) \, \md \xi} \right\vert \underset{\underset{
z_{o}=1+o(1)}{\mathscr{N},n \to \infty}}{\leqslant} \me^{-((n-1)K
+k)(\tilde{\delta}_{\tilde{a}_{j}})^{N+2} \sigma_{\tilde{a}_{j}}^{+} \inf_{
\tau \in (0,\sigma_{\tilde{a}_{j}}^{+})} \left\vert \tilde{h}_{\widetilde{V}}
(z(\tau)) \right\vert} \quad (\ll \! 1), \label{eqlem5.1j} \\
\left\vert \me^{-2 \pi \mi ((n-1)K+k) \int_{z}^{\tilde{p}_{+}(j)} \psi_{
\widetilde{V}}^{f}(\xi) \, \md \xi} \right\vert \underset{\underset{
z_{o}=1+o(1)}{\mathscr{N},n \to \infty}}{\leqslant} \me^{-((n-1)K+k)
\tilde{\lambda}_{\tilde{\mathcal{R}},3}(j) \lvert z-\tilde{p}_{+}(j) \rvert}, 
\label{eqlem5.1k}
\end{gather}
where
\begin{equation*}
\tilde{\lambda}_{\tilde{\mathcal{R}},3}(j) \! := \! \left(\tilde{\eta}_{j} 
\min \left\lbrace \dfrac{\tilde{\delta}_{\tilde{b}_{j-1}} \sin \sigma_{
\tilde{b}_{j-1}}^{+}}{\tilde{l}_{\tilde{b}_{j-1}}},\dfrac{\tilde{\delta}_{
\tilde{a}_{j}} \sin \sigma_{\tilde{a}_{j}}^{+}}{\tilde{l}_{\tilde{a}_{j}}} 
\right\rbrace \right)^{N+1} \inf_{\theta \in (\theta_{0}^{\tilde{a}}(j),
 \pi -\theta_{0}^{\tilde{b}}(j))} \left\lvert \tilde{h}_{\widetilde{V}}
(x_{j}(\theta) \! + \! \mi y_{j}(\theta)) \right\rvert,
\end{equation*}
with $\tilde{\lambda}_{\tilde{\mathcal{R}},3}(j) \! = \! \tilde{\lambda}_{
\tilde{\mathcal{R}},3}(n,k,z_{o};j) \! =_{\underset{z_{o}=1+o(1)}{
\mathscr{N},n \to \infty}} \! \mathcal{O}(1)$ and $> \! 0$, and 
$x_{j}(\theta),y_{j}(\theta)$ as defined above; hence, for $n \! \in \! 
\mathbb{N}$ and $k \! \in \! \lbrace 1,2,\dotsc,K \rbrace$ such that 
$\alpha_{p_{\mathfrak{s}}} \! := \! \alpha_{k} \! \neq \! \infty$, one 
shows, via Equations~\eqref{eqlem5.1g} and~\eqref{eqlem5.1h} 
and the Estimates~\eqref{eqlem5.1i}--\eqref{eqlem5.1k}, that
\begin{equation} \label{eqlem5.1l} 
\me^{-2 \pi \mi ((n-1)K+k) \int_{z}^{\tilde{a}_{N+1}} \psi_{\widetilde{
V}}^{f}(\xi) \, \md \xi} \underset{\underset{z_{o}=1+o(1)}{\mathscr{N},
n \to \infty}}{=} \mathcal{O} \left(\tilde{\mathfrak{c}}_{\tilde{\mathcal{R}},
3}^{\blacklozenge}(j) \me^{-((n-1)K+k) \tilde{\lambda}_{\tilde{
\mathcal{R}},3}(j) \lvert z-\tilde{p}_{+}(j) \rvert} \right), \quad z \! 
\in \! \tilde{\Sigma}_{p,j}^{3}, \quad j \! = \! 1,2,\dotsc,N \! + \! 1,
\end{equation}
where $\tilde{\mathfrak{c}}_{\tilde{\mathcal{R}},3}^{\blacklozenge}(j) 
\! = \! \tilde{\mathfrak{c}}_{\tilde{\mathcal{R}},3}^{\blacklozenge}
(n,k,z_{o};j) \! =_{\underset{z_{o}=1+o(1)}{\mathscr{N},n \to \infty}} 
\! \mathcal{O}(1)$, $j \! = \! 1,2,\dotsc,N \! + \! 1$. For $n \! \in \! 
\mathbb{N}$ and $k \! \in \! \lbrace 1,2,\dotsc,K \rbrace$ such that 
$\alpha_{p_{\mathfrak{s}}} \! := \! \alpha_{k} \! \neq \! \infty$, and for 
$z \! \in \! \tilde{\Sigma}_{p,j}^{3}$, $j \! = \! 1,2,\dotsc,N \! + \! 1$, 
in order to study $\mathfrak{m}(z)$, one notes {}from item~\pmb{(2)} 
of Lemma~\ref{lem4.5} that the $z$-dependent factors $\tilde{
\boldsymbol{\theta}}(\varepsilon_{1} \tilde{\boldsymbol{u}}(z) \! - \! 
\tfrac{1}{2 \pi}((n \! - \! 1)K \! + \! k) \tilde{\boldsymbol{\Omega}} \! 
+ \! \varepsilon_{2} \tilde{\boldsymbol{d}})$, $\tilde{\boldsymbol{
\theta}}(\varepsilon_{1} \tilde{\boldsymbol{u}}(z) \! + \! 
\varepsilon_{2} \tilde{\boldsymbol{d}})$, $\varepsilon_{1},
\varepsilon_{2} \! = \! \pm 1$, $\tilde{\gamma}(z)$, and 
$(\tilde{\gamma}(z))^{-1}$ must be analysed (all the parameters 
appearing in the latter formulae are defined in item~\pmb{(2)} of 
Lemma~\ref{lem4.5}). For $n \! \in \! \mathbb{N}$ and $k \! \in \! \lbrace 
1,2,\dotsc,K \rbrace$ such that $\alpha_{p_{\mathfrak{s}}} \! := \! 
\alpha_{k} \! \neq \! \infty$, and for $z \! \in \! \tilde{\Sigma}_{p,j}^{3}$, 
$j \! = \! 1,2,\dotsc,N \! + \! 1$, one shows, via the representation 
for $\tilde{\boldsymbol{\theta}}(\pmb{\cdot})$ given by 
Equation~\eqref{eqrmthetafin} and the elementary inequality $\lvert 
\exp (\ast) \rvert \! \leqslant \! \exp (\lvert \ast \rvert)$, that $\lvert 
\tilde{\boldsymbol{\theta}}(\varepsilon_{1} \tilde{\boldsymbol{u}}(z) \! 
- \! \tfrac{1}{2 \pi}((n \! - \! 1)K \! + \! k) \tilde{\boldsymbol{\Omega}} 
\! + \! \varepsilon_{2} \tilde{\boldsymbol{d}}) \rvert \! \leqslant \! 
\sum_{m \in \mathbb{Z}^{N}} \me^{-2 \pi \varepsilon_{2}(m,\Im 
(\tilde{\boldsymbol{d}}))} \lvert \me^{\mi \pi (m,\tilde{\boldsymbol{\tau}}
m)} \rvert \lvert \me^{2 \pi \mi \varepsilon_{1}(m,\tilde{\boldsymbol{u}}
(z))} \rvert$. Recall {}from the associated Subsection~\ref{subsub1} that, 
for $n \! \in \! \mathbb{N}$ and $k \! \in \! \lbrace 1,2,\dotsc,K \rbrace$ 
such that $\alpha_{p_{\mathfrak{s}}} \! := \! \alpha_{k} \! \neq \! \infty$, 
$(m,\tilde{\boldsymbol{u}}(z)) \! = \! \sum_{i=1}^{N}m_{i} \tilde{
\boldsymbol{u}}_{i}(z)$, where $\tilde{\boldsymbol{u}}_{i}(z) \! := \! 
\int_{\tilde{a}_{N+1}}^{z} \tilde{\omega}_{i}$, with $\tilde{\omega}_{i} \! 
= \! \sum_{j=1}^{N} \tilde{c}_{ij}(\tilde{R}(z))^{-1/2}z^{N-j} \, \md z$, 
and $\tilde{c} \! = \! (\tilde{c})_{i,j=1,2,\dotsc,N}$ as described in 
Equations~\eqref{O1} and~\eqref{O2}. For $z \! \in \! \tilde{\Sigma}_{p,
j}^{3}$, $j \! = \! 1,2,\dotsc,N \! + \! 1$, write
\begin{equation*}
\tilde{\boldsymbol{u}}_{i}(z) \! = \! -\left(\int_{z}^{\tilde{p}_{+}(j)} \! + 
\! \int_{\tilde{C}_{j}^{\ast}} \! + \! \int_{\tilde{a}_{j}+\tilde{\delta}_{
\tilde{a}_{j}}}^{\tilde{a}_{N+1}} \right) \tilde{\omega}_{i}, \quad i \! = 
\! 1,2,\dotsc,N;
\end{equation*}
via the identities $x_{1}^{m} \! - \! x_{2}^{m} \! = \! (x_{1} \! - \! x_{2})
(x_{1}^{m-1} \! + \! x_{1}^{m-2}x_{2} \! + \! \dotsb \! + \! x_{1}
x_{2}^{m-2} \! + \! x_{2}^{m-1})$, $m \! \in \! \mathbb{N}$, and 
$\lvert x_{1} \! - \! x_{2} \rvert \! \geqslant \! \lvert x_{1} \rvert \! - 
\! \lvert x_{2} \rvert$, and certain well-known trigonometric identities, 
one arrives at, after repeated application of both the extended triangle 
inequality and the ML Theorem, the following estimates: for $z \! \in 
\! \tilde{\Sigma}_{p,j}^{3}$, $j \! = \! 1,2,\dotsc,N \! + \! 1$, and 
$i \! = \! 1,2,\dotsc,N$,
\begin{align}
\left\lvert \int_{z}^{\tilde{p}_{+}(j)} \tilde{\omega}_{i} \right\rvert 
\underset{\underset{z_{o}=1+o(1)}{\mathscr{N},n \to \infty}}{\leqslant}& 
\, \dfrac{\max_{i_{1}=1,2,\dotsc,N} \lbrace \lvert \tilde{c}_{ii_{1}} \rvert 
\rbrace}{\tilde{\delta}_{\tilde{a}_{j}}^{N+1}} \left(\dfrac{((\tilde{a}_{j} 
\! + \! \tilde{\delta}_{\tilde{a}_{j}} \cos \sigma_{\tilde{a}_{j}}^{+})^{2} \! 
+ \! \tilde{r}_{j}^{2})^{N/2} \! - \! 1}{((\tilde{a}_{j} \! + \! \tilde{\delta}_{
\tilde{a}_{j}} \cos \sigma_{\tilde{a}_{j}}^{+})^{2} \! + \! \tilde{r}_{j}^{2})^{
1/2} \! - \! 1} \right) \nonumber \\
\times& \, \left((\tilde{a}_{j} \! - \! \tilde{b}_{j-1} \! + \! \tilde{\delta}_{
\tilde{a}_{j}} \cos \sigma_{\tilde{a}_{j}}^{+} \! + \! \tilde{\delta}_{\tilde{
b}_{j-1}} \cos \sigma_{\tilde{b}_{j-1}}^{+})^{2} \! + \! (\tilde{\delta}_{
\tilde{a}_{j}} \sin \sigma_{\tilde{a}_{j}}^{+} \! - \! \tilde{\delta}_{\tilde{
b}_{j-1}} \sin \sigma_{\tilde{b}_{j-1}}^{+})^{2} \right)^{1/2}, \label{eqlem5.1m} \\
\left\lvert \int_{\tilde{C}_{j}^{\ast}} \tilde{\omega}_{i} \right\rvert 
\underset{\underset{z_{o}=1+o(1)}{\mathscr{N},n \to \infty}}{\leqslant}& 
\, \dfrac{\sigma_{\tilde{a}_{j}}^{+} \max_{i_{1}=1,2,\dotsc,N} \lbrace 
\lvert \tilde{c}_{ii_{1}} \rvert \rbrace}{\tilde{\delta}_{\tilde{a}_{j}}^{N}} 
\left(\dfrac{(\lvert \tilde{a}_{j} \rvert \! + \! \tilde{\delta}_{\tilde{a}_{j}})^{
N} \! - \! 1}{\lvert \tilde{a}_{j} \rvert \! + \! \tilde{\delta}_{\tilde{a}_{j}} 
\! - \! 1} \right), \label{eqlem5.1n} \\
\left\lvert \int_{\tilde{a}_{N+1}}^{\tilde{a}_{j}+\tilde{\delta}_{\tilde{a}_{j}}} 
\tilde{\omega}_{i} \right\rvert \underset{\underset{z_{o}=1+o(1)}{
\mathscr{N},n \to \infty}}{\leqslant}& \, \dfrac{\max_{i_{1}=1,2,\dotsc,N} 
\lbrace \lvert \tilde{c}_{ii_{1}} \rvert \rbrace}{(\min \lbrace \tilde{\delta}_{
\tilde{a}_{j}},\tilde{\delta}_{\tilde{b}_{j-1}} \rbrace)^{N+1}} \left(\dfrac{
\lvert \tilde{a}_{N+1} \rvert^{N} \! - \! 1}{\lvert \tilde{a}_{N+1} \rvert \! 
- \! 1} \right)(\tilde{a}_{N+1} \! - \! \tilde{a}_{j} \! + \! \tilde{\delta}_{
\tilde{a}_{j}}). \label{eqlem5.1o}
\end{align}
Recalling that the associated $N \times N$ Riemann matrix of $\tilde{
\boldsymbol{\beta}}$-periods, $\tilde{\boldsymbol{\tau}} \! = \! (\tilde{
\boldsymbol{\tau}})_{i,j=1,2,\dotsc,N}$, is non-degenerate, symmetric 
and pure imaginary, and $-\mi \tilde{\boldsymbol{\tau}}$ is positive 
definite, it follows via the Estimates~\eqref{eqlem5.1m}--\eqref{eqlem5.1o}, 
the extended triangle inequality, and the fact that $\tilde{\boldsymbol{
\theta}}(\pmb{\cdot})$ converges absolutely and uniformly, that, for 
$n \! \in \! \mathbb{N}$ and $k \! \in \! \lbrace 1,2,\dotsc,K \rbrace$ 
such that $\alpha_{p_{\mathfrak{s}}} \! := \! \alpha_{k} \! \neq \! \infty$, 
and uniformly for $z \! \in \! \tilde{\Sigma}_{p,j}^{3}$, $j \! = \! 1,2,
\dotsc,N \! + \! 1$,
\begin{equation*}
\lvert \tilde{\boldsymbol{\theta}}(\varepsilon_{1} \tilde{\boldsymbol{u}}
(z) \! - \! \tfrac{1}{2 \pi}((n \! - \! 1)K \! + \! k) \tilde{\boldsymbol{
\Omega}} \! + \! \varepsilon_{2} \tilde{\boldsymbol{d}}) \rvert 
\underset{\underset{z_{o}=1+o(1)}{\mathscr{N},n \to \infty}}{\leqslant} 
\mathcal{O}(\tilde{\mathfrak{c}}_{\tilde{\theta},1}^{\blacklozenge}(n,k,
z_{o};j)), \qquad \lvert \tilde{\boldsymbol{\theta}}(\varepsilon_{1} \tilde{
\boldsymbol{u}}(z) \! + \! \varepsilon_{2} \tilde{\boldsymbol{d}}) \rvert 
\underset{\underset{z_{o}=1+o(1)}{\mathscr{N},n \to \infty}}{\leqslant} 
\mathcal{O}(\tilde{\mathfrak{c}}_{\tilde{\theta},2}^{\blacklozenge}(n,k,
z_{o};j)),
\end{equation*}
where $\tilde{\mathfrak{c}}_{\tilde{\theta},i_{1}}^{\blacklozenge}
(n,k,z_{o};j) \! =_{\underset{z_{o}=1+o(1)}{\mathscr{N},n \to \infty}} \! 
\mathcal{O}(1)$, $i_{1} \! = \! 1,2$. Via the definition of $\tilde{\gamma}
(z)$ given by Equation~\eqref{eqmainfin10}, one shows that, uniformly 
for $z \! \in \! \tilde{\Sigma}_{p,j}^{3}$, $j \! = \! 1,2,\dotsc,N \! + \! 1$,
\begin{align*}
\lvert \tilde{\gamma}(z) \rvert \underset{\underset{z_{o}=1+o(1)}{
\mathscr{N},n \to \infty}}{\leqslant}& \, \dfrac{((\tilde{a}_{j} \! - \! 
\tilde{b}_{j-1} \! + \! \tilde{\delta}_{\tilde{a}_{j}} \cos \sigma_{
\tilde{a}_{j}}^{+})^{2} \! + \! \tilde{r}_{j}^{2})^{1/8}}{((\tilde{b}_{j-1} \! 
- \! \tilde{a}_{j} \! - \! \tilde{\delta}_{\tilde{b}_{j-1}} \cos \sigma_{
\tilde{b}_{j-1}}^{+})^{2} \! + \! (\min \lbrace \tilde{\delta}_{\tilde{a}_{j}} 
\sin \sigma_{\tilde{a}_{j}}^{+},\tilde{\delta}_{\tilde{b}_{j-1}} \sin \sigma_{
\tilde{b}_{j-1}}^{+} \rbrace)^{2})^{1/8}} \\
\times& \, \prod_{\substack{i=1\\i \neq j}}^{N+1} \dfrac{((\tilde{a}_{j} \! - 
\! \tilde{b}_{i-1} \! + \! \tilde{\delta}_{\tilde{a}_{j}} \cos \sigma_{\tilde{
a}_{j}}^{+})^{2} \! + \! \tilde{r}_{j}^{2})^{1/8}}{((\tilde{b}_{j-1} \! - \! 
\tilde{a}_{i} \! - \! \tilde{\delta}_{\tilde{b}_{j-1}} \cos \sigma_{\tilde{b}_{j
-1}}^{+})^{2} \! + \! (\min \lbrace \tilde{\delta}_{\tilde{a}_{j}} \sin 
\sigma_{\tilde{a}_{j}}^{+},\tilde{\delta}_{\tilde{b}_{j-1}} \sin \sigma_{
\tilde{b}_{j-1}}^{+} \rbrace)^{2})^{1/8}} \! =: \! \tilde{\mathfrak{c}}_{
\tilde{\theta},3}^{\blacklozenge}(n,k,z_{o};j), \\
\lvert (\tilde{\gamma}(z))^{-1} \rvert \underset{\underset{z_{o}=1+
o(1)}{\mathscr{N},n \to \infty}}{\leqslant}& \, \dfrac{((\tilde{\delta}_{
\tilde{a}_{j}} \cos \sigma_{\tilde{a}_{j}}^{+})^{2} \! + \! \tilde{r}_{j}^{2})^{
1/8}}{((\tilde{\delta}_{\tilde{b}_{j-1}} \cos \sigma_{\tilde{b}_{j-1}}^{+})^{2} 
\! + \! (\min \lbrace \tilde{\delta}_{\tilde{a}_{j}} \sin \sigma_{\tilde{a}_{j}}^{
+},\tilde{\delta}_{\tilde{b}_{j-1}} \sin \sigma_{\tilde{b}_{j-1}}^{+} \rbrace)^{
2})^{1/8}} \\
\times& \, \prod_{\substack{i=1\\i \neq j}}^{N+1} \dfrac{((\tilde{a}_{j} \! - 
\! \tilde{a}_{i} \! + \! \tilde{\delta}_{\tilde{a}_{j}} \cos \sigma_{\tilde{a}_{
j}}^{+})^{2} \! + \! \tilde{r}_{j}^{2})^{1/8}}{((\tilde{b}_{j-1} \! - \! \tilde{b}_{i
-1} \! - \! \tilde{\delta}_{\tilde{b}_{j-1}} \cos \sigma_{\tilde{b}_{j-1}}^{+})^{2} 
\! + \! (\min \lbrace \tilde{\delta}_{\tilde{a}_{j}} \sin \sigma_{\tilde{a}_{j}}^{+},
\tilde{\delta}_{\tilde{b}_{j-1}} \sin \sigma_{\tilde{b}_{j-1}}^{+} \rbrace)^{2})^{
1/8}} \! =: \! \tilde{\mathfrak{c}}_{\tilde{\theta},4}^{\blacklozenge}(n,k,z_{o};j),
\end{align*}
where $\tilde{\mathfrak{c}}_{\tilde{\theta},i_{1}}^{\blacklozenge}(n,k,z_{o};j) 
\! =_{\underset{z_{o}=1+o(1)}{\mathscr{N},n \to \infty}} \! \mathcal{O}(1)$, 
$i_{1} \! = \! 3,4$; hence, for $n \! \in \! \mathbb{N}$ and $k \! \in \! \lbrace 
1,2,\dotsc,K \rbrace$ such that $\alpha_{p_{\mathfrak{s}}} \! := \! \alpha_{k} 
\! \neq \! \infty$, via the above-derived bounds, one arrives at, after a 
matrix-multiplication argument, uniformly for $z \! \in \! \tilde{\Sigma}_{p,
j}^{3}$, $j \! = \! 1,2,\dotsc,N \! + \! 1$,
\begin{equation} \label{eqlem5.1p} 
\mathfrak{m}(z) \sigma_{-}(\mathfrak{m}(z))^{-1} \underset{\underset{z_{o}
=1+o(1)}{\mathscr{N},n \to \infty}}{=} \mathcal{O}(\tilde{\mathfrak{c}}_{
\tilde{\theta},5}^{\blacklozenge}(n,k,z_{o};j)),
\end{equation}
where $(\mathrm{M}_{2}(\mathbb{C}) \! \ni)$ $\tilde{\mathfrak{c}}_{
\tilde{\theta},5}^{\blacklozenge}(n,k,z_{o};j) \! =_{\underset{z_{o}=
1+o(1)}{\mathscr{N},n \to \infty}} \! \mathcal{O}(1)$. Finally, via 
Equation~\eqref{eqlem5.1g}, the Estimates~\eqref{eqlem5.1l} 
and~\eqref{eqlem5.1p}, and a matrix-multiplication argument, 
one arrives at, for $n \! \in \! \mathbb{N}$ and $k \! \in \! \lbrace 
1,2,\dotsc,K \rbrace$ such that $\alpha_{p_{\mathfrak{s}}} \! := \! 
\alpha_{k} \! \neq \! \infty$, uniformly for $z \! \in \! \tilde{\Sigma}_{p,j}^{3}$, 
$j \! = \! 1,2,\dotsc,N \! + \! 1$, the asymptotics, in the double-scaling limit 
$\mathscr{N},n \! \to \! \infty$ such that $z_{o} \! = \! 1 \! + \! o(1)$, for 
$\tilde{v}_{\tilde{\mathcal{R}}}(z)$ given in Equation~\eqref{eqtlvee11}. For 
$n \! \in \! \mathbb{N}$ and $k \! \in \! \lbrace 1,2,\dotsc,K \rbrace$ such 
that $\alpha_{p_{\mathfrak{s}}} \! := \! \alpha_{k} \! \neq \! \infty$, the case 
$z \! \in \! \tilde{\Sigma}_{p,j}^{4}$, $j \! = \! 1,2,\dotsc,N \! + \! 1$, 
is, \emph{mutatis mutandis}, analysed analogously, and leads to the 
asymptotics, in the double-scaling limit $\mathscr{N},n \! \to \! \infty$ 
such that $z_{o} \! = \! 1 \! + \! o(1)$, for $\tilde{v}_{\tilde{\mathcal{R}}}(z)$ 
given in Equation~\eqref{eqtlvee12}.

For $n \! \in \! \mathbb{N}$ and $k \! \in \! \lbrace 1,2,\dotsc,K \rbrace$ 
such that $\alpha_{p_{\mathfrak{s}}} \! := \! \alpha_{k} \! \neq \! \infty$, 
and for $z \! \in \! \tilde{\Sigma}_{p,j}^{5} \! := \! \partial \tilde{\mathbb{
U}}_{\tilde{\delta}_{\tilde{b}_{j-1}}} \cup \partial \tilde{\mathbb{U}}_{
\tilde{\delta}_{\tilde{a}_{j}}}$, $j \! = \! 1,2,\dotsc,N \! + \! 1$, consider, 
say, and without loss of generality, the case $z \! \in \! \partial \tilde{
\mathbb{U}}_{\tilde{\delta}_{\tilde{a}_{j}}}$; the analysis for the case $z 
\! \in \! \partial \tilde{\mathbb{U}}_{\tilde{\delta}_{\tilde{b}_{j-1}}}$ is, 
\emph{mutatis mutandis}, analogous. For $n \! \in \! \mathbb{N}$ and $k 
\! \in \! \lbrace 1,2,\dotsc,K \rbrace$ such that $\alpha_{p_{\mathfrak{s}}} 
\! := \! \alpha_{k} \! \neq \! \infty$, and for $z \! \in \! \partial \tilde{
\mathbb{U}}_{\tilde{\delta}_{\tilde{a}_{j}}}$ $(\subset \tilde{\Sigma}_{p,
j}^{5})$, $j \! = \! 1,2,\dotsc,N \! + \! 1$, recall the expression for 
$\tilde{v}_{\tilde{\mathcal{R}}}(z)$ given in the corresponding item~\pmb{(2)} 
of Lemma~\ref{lem4.10}: $\tilde{v}_{\tilde{\mathcal{R}}}(z) \! = \! 
\tilde{\mathcal{X}}^{\tilde{a}}(z)(\mathfrak{m}(z))^{-1}$, $z \! \in \! \partial 
\tilde{\mathbb{U}}_{\tilde{\delta}_{\tilde{a}_{j}}}$, $j \! = \! 1,2,\dotsc,N 
\! + \! 1$, where the parametrix, $\tilde{\mathcal{X}}^{\tilde{a}}(z)$, is 
characterised in Lemma~\ref{lem4.9}, and $\mathfrak{m}(z)$ is given in 
item~\pmb{(2)} of Lemma~\ref{lem4.5}. Using the explicit expressions for 
$\tilde{\mathcal{X}}^{\tilde{a}}(z)$ given in Lemma~\ref{lem4.9}, and the 
large-argument asymptotics for the Airy function and its derivative given 
in Equation~\eqref{eqairy}, one shows that, for $n \! \in \! \mathbb{N}$ 
and $k \! \in \! \lbrace 1,2,\dotsc,K \rbrace$ such that $\alpha_{
p_{\mathfrak{s}}} \! := \! \alpha_{k} \! \neq \! \infty$,
\begin{align*}
\tilde{v}_{\tilde{\mathcal{R}}}(z) \underset{\underset{z_{o}=1+o(1)}{\mathscr{N},
n \to \infty}}{=}& \, \mathrm{I} \!  + \! \dfrac{\me^{-\frac{\mi \pi}{3}}}{((n \! 
- \! 1)K \! + \! k) \tilde{\xi}_{\tilde{a}_{j}}(z)} \mathfrak{m}(z) 
\begin{pmatrix}
\mi \me^{\pm \frac{\mi}{2}((n-1)K+k) \tilde{\mho}_{j}} & -\mi \me^{\pm 
\frac{\mi}{2}((n-1)K+k) \tilde{\mho}_{j}} \\
\me^{\mp \frac{\mi}{2}((n-1)K+k) \tilde{\mho}_{j}} & \me^{\mp \frac{\mi}{2}
((n-1)K+k) \tilde{\mho}_{j}}
\end{pmatrix} \\
\times& \, 
\begin{pmatrix}
-s_{1} \me^{-\frac{\mi \pi}{6}} \me^{\mp \frac{\mi}{2}((n-1)K+k) \tilde{
\mho}_{j}} & s_{1} \me^{\frac{\mi \pi}{3}} \me^{\pm \frac{\mi}{2}((n-1)K+k) 
\tilde{\mho}_{j}} \\
t_{1} \me^{-\frac{\mi \pi}{6}} \me^{\mp \frac{\mi}{2}((n-1)K+k) \tilde{
\mho}_{j}} & -t_{1} \me^{-\frac{2 \pi \mi}{3}} \me^{\pm \frac{\mi}{2}
((n-1)K+k) \tilde{\mho}_{j}}
\end{pmatrix}
(\mathfrak{m}(z))^{-1} \\
+& \, \mathcal{O} \left(\dfrac{1}{((n \! - \! 1)K \! + \! k)^{2}
(\tilde{\xi}_{\tilde{a}_{j}}(z))^{2}} \mathfrak{m}(z) \tilde{\mathfrak{c}}^{
\triangleleft}(n,k,z_{o};j)(\mathfrak{m}(z))^{-1} \right), \quad z \! \in 
\! \mathbb{C}_{\pm} \cap \partial \tilde{\mathbb{U}}_{\tilde{\delta}_{
\tilde{a}_{j}}}, \quad j \! = \! 1,2,\dotsc,N \! + \! 1,
\end{align*}
where $\tilde{\xi}_{\tilde{a}_{j}}(z)$, $\tilde{\mho}_{j}$, $s_{1}$, $t_{1}$, 
and $\tilde{\mathfrak{c}}^{\triangleleft}(n,k,z_{o};j)$ are described in the 
corresponding item~\pmb{(2)}, subitem~$\boldsymbol{\mathrm{(2)}_{iv}}$ 
of the lemma, and $\mathfrak{m}(z)$ is given in item~\pmb{(2)} of 
Lemma~\ref{lem4.5}. Using the formula for $\mathfrak{m}(z)$ in terms of 
$\mathbb{M}(z)$ given in item~\pmb{(2)} of Lemma~\ref{lem4.5}, after a 
matrix-multiplication argument, one arrives at, for $n \! \in \! \mathbb{N}$ 
and $k \! \in \! \lbrace 1,2,\dotsc,K \rbrace$ such that $\alpha_{p_{\mathfrak{s}}} 
\! := \! \alpha_{k} \! \neq \! \infty$, and for $z \! \in \! \partial 
\tilde{\mathbb{U}}_{\tilde{\delta}_{\tilde{a}_{j}}}$ $(\subset \tilde{\Sigma}_{p,j}^{5})$, 
$j \! = \! 1,2,\dotsc,N \! + \! 1$, the asymptotics, in the double-scaling limit 
$\mathscr{N},n \! \to \! \infty$ such that $z_{o} \! = \! 1 \! + \! o(1)$, for 
$\tilde{v}_{\tilde{\mathcal{R}}}(z)$ given in Equation~\eqref{eqtlvee13}. For 
$n \! \in \! \mathbb{N}$ and $k \! \in \! \lbrace 1,2,\dotsc,K \rbrace$ such 
that $\alpha_{p_{\mathfrak{s}}} \! := \! \alpha_{k} \! \neq \! \infty$, the case 
$z \! \in \! \partial \tilde{\mathbb{U}}_{\tilde{\delta}_{\tilde{b}_{j-1}}}$ 
$(\subset \tilde{\Sigma}_{p,j}^{5})$, $j \! = \! 1,2,\dotsc,N \! + \! 1$, 
is, \emph{mutatis mutandis}, analysed analogously, and leads to the 
asymptotics, in the double-scaling limit $\mathscr{N},n \! \to \! \infty$ 
such that $z_{o} \! = \! 1 \! + \! o(1)$, for $\tilde{v}_{\tilde{\mathcal{R}}}
(z)$ given in Equation~\eqref{eqtlvee14}. \hfill $\qed$
\begin{aaaaa} \label{def5.1}
\textsl{For an oriented contour $D$ $(\subset \! \overline{\mathbb{C}})$, 
let $\mathfrak{B}_{q}(D)$, $q \! \in \! \lbrace 1,2,\infty \rbrace$, 
denote the set of all bounded linear operators {}from $\mathcal{L}^{q}_{
\mathrm{M}_{2}(\mathbb{C})}(D)$ into $\mathcal{L}^{q}_{\mathrm{M}_{2}
(\mathbb{C})}(D)$.}
\end{aaaaa}

Since the asymptotic analysis, in the double-scaling limit $\mathscr{N},n 
\! \to \! \infty$ such that $z_{o} \! = \! 1 \! + \! o(1)$, that follows relies 
substantially on the BC construction \cite{bealscoif} for the solution of a 
matrix, and suitably normalised, RHP on an oriented and unbounded contour, 
it is prudent to present, at this stage, a succinct synopsis of it. One agrees 
to call a contour $D$ \emph{oriented} if: (i) $\mathbb{C} \setminus D$ has 
countably many open connected components; (ii) $\mathbb{C} \setminus D$ 
is the disjoint union of two, possibly disconnected, open regions, denoted 
by $\boldsymbol{\mathfrak{O}}^{+}$ and $\boldsymbol{\mathfrak{O}}^{-}$; 
and (iii) $D$ may be viewed as either the positively oriented boundary for 
$\boldsymbol{\mathfrak{O}}^{+}$ or the negatively oriented boundary 
for $\boldsymbol{\mathfrak{O}}^{-}$ $(\mathbb{C} \setminus D$ 
is coloured by the two colours $\pm)$. Let $D$, as a closed set, be 
the union of countably many oriented, simple, piecewise-smooth 
arcs. Denote the set of all self-intersections of $D$ by 
$D^{\sharp}$.\footnote{It is assumed that $\mathrm{card}(D^{\sharp}) 
\! < \! \infty$.} Set $D^{\natural} \! := \! D \setminus D^{\sharp}$. The 
BC construction \cite{bealscoif} for the solution of a matrix-valued RHP 
(see, also, \cite{a63,a64,xinzoo1,xinzoo2,xinzoo3}) on an oriented 
contour $D$ consists of finding a matrix-valued function $\mathfrak{W} 
\colon \mathbb{C} \setminus D \! \to \! \mathrm{M}_{2}(\mathbb{C})$ 
such that: (1) $\mathfrak{W}(\zeta)$ is analytic for $\zeta 
\! \in \! \mathbb{C} \setminus D$, $\mathfrak{W}(\zeta) \! \! 
\upharpoonright_{\mathbb{C} \setminus D}$ has a continuous 
extension, {}from `above' and `below', to $D^{\natural}$, and 
$\lim_{\genfrac{}{}{0pt}{2}{\zeta^{\prime} \to \zeta}{\zeta^{\prime} \, 
\in \, \pm \, \mathrm{side} \, \mathrm{of} \, D^{\natural}}} \int_{D^{
\natural}} \lvert \mathfrak{W}(\zeta^{\prime}) \! - \! \mathfrak{W}_{\pm}
(\zeta) \rvert^{2} \, \lvert \md \zeta \rvert \! = \! 0$, where 
$\mathfrak{W}_{\pm}(\zeta) \! := \! \lim_{\underset{\zeta^{\prime} \, \in 
\, \pm \, \mathrm{side} \, \mathrm{of} \, D^{\natural}}{\zeta^{\prime} 
\to \zeta}} \mathfrak{W}(\zeta^{\prime})$; (2) the boundary values 
$\mathfrak{W}_{\pm}(\zeta)$ satisfy the jump condition $\mathfrak{W}_{+}
(\zeta) \! = \! \mathfrak{W}_{-}(\zeta) \upsilon (\zeta)$ $\text{a.e.}$ $\zeta 
\! \in \! D^{\natural}$, for some smooth jump matrix $\upsilon \colon 
D^{\natural} \! \to \! \mathrm{GL}_{2}(\mathbb{C})$; and (3) for an 
arbitrarily fixed $\zeta_{o} \! \in \! \overline{\mathbb{C}}$, $\mathfrak{W}
(\zeta) \! =_{\overline{\mathbb{C}} \setminus D \ni \zeta \to \zeta_{o}} \! 
\mathrm{I} \! + \! o(1)$, where, uniformly with respect to $\zeta$, $o(1) 
\! = \! \mathcal{O}(\zeta \! - \! \zeta_{o})$ if $\zeta_{o}$ is bounded, 
and $o(1) \! = \! \mathcal{O}(\zeta^{-1})$ if $\zeta_{o}$ is the point at 
infinity.\footnote{Condition~(3) is referred to as the \emph{normalisation 
condition}: one says that the RHP is `normalised at $\zeta_{o}$'.} Let 
$\upsilon (\zeta) \! := \! (\mathrm{I} \! - \! w_{-}(\zeta))^{-1}
(\mathrm{I} \! + \! w_{+}(\zeta))$, $\zeta \! \in \! D^{\natural}$, 
be a bounded algebraic factorisation for $\upsilon (\zeta)$, where 
$w_{\pm}(\zeta)$ are upper/lower, or lower/upper, triangular matrices 
(depending on the orientation of $D$), and $w_{\pm}(\zeta) \! \in \! 
\cap_{p=2,\infty} \mathcal{L}^{p}_{\mathrm{M}_{2}(\mathbb{C})}
(D^{\natural})$.\footnote{If $D^{\natural}$ is unbounded, then one 
requires that $w_{\pm}(\zeta) \! =_{\zeta \to \infty} 
\left(
\begin{smallmatrix}
0 & 0 \\
0 & 0
\end{smallmatrix}
\right)$.} Define $w(\zeta) \! := \! w_{+}(\zeta) \! + \! w_{-}
(\zeta)$, and introduce (see Remark~\ref{ouhzphreton} below) 
the---normalised at $\zeta_{o}$---Cauchy (singular integral) 
operators
\begin{equation*}
\mathcal{L}^{2}_{\mathrm{M}_{2}(\mathbb{C})}(D) \! \ni \! f \! \mapsto 
\! (C^{\zeta_{o}}_{\pm}f)(\zeta) := \lim_{\genfrac{}{}{0pt}{2}{\zeta^{\prime} 
\to \zeta}{\zeta^{\prime} \, \in \, \pm \, \mathrm{side} \, \mathrm{of} 
\, D}} \int_{D} \dfrac{(\zeta^{\prime} \! - \! \zeta_{o})f(\xi)}{(\xi \! - \! 
\zeta_{o})(\xi \! - \! \zeta^{\prime})} \, \dfrac{\md \xi}{2 \pi \mi},
\end{equation*}
where $\tfrac{(\zeta-\zeta_{o})}{(\xi -\zeta_{o})(\xi -\zeta)} \, 
\tfrac{\md \xi}{2 \pi \mi}$ is the Cauchy kernel normalised at 
$\zeta_{o}$. Introduce, now, the `modified' BC operator $C^{\zeta_{o}}_{w}$:
\begin{equation*}
\mathcal{L}^{2}_{\mathrm{M}_{2}(\mathbb{C})}(D) \! \ni \! f \! \mapsto \! 
C^{\zeta_{o}}_{w}f \! := \! C^{\zeta_{o}}_{+}(fw_{-}) \! + \! C^{\zeta_{o}}_{-}
(fw_{+}).
\end{equation*}
The solution of the---normalised at $\zeta_{o}$---BC RHP is given by 
the following integral representation \cite{bealscoif}.
\begin{ccccc} \label{lem5.2} 
Let
\begin{equation*}
\mu_{\zeta_{o}}(\zeta) \! = \! \mathfrak{W}_{+}(\zeta)(\mathrm{I} \! + \! 
w_{+}(\zeta))^{-1} \! = \! \mathfrak{W}_{-}(\zeta)(\mathrm{I} \! - \! w_{-}
(\zeta))^{-1}, \quad \zeta \! \in \! D.
\end{equation*}
If $\mu_{\zeta_{o}} \! \in \! \mathrm{I} \! + \! \mathcal{L}^{2}_{\mathrm{
M}_{2}(\mathbb{C})}(D)$ solves the linear singular integral equation
\begin{equation*}
(\id \! - \! C^{\zeta_{o}}_{w})(\mu_{\zeta_{o}}(\zeta) \! - \! \mathrm{I}) \! 
= \! C_{w}^{\zeta_{o}} \mathrm{I} \! = \! C^{\zeta_{o}}_{+}(w_{-}(\zeta)) 
\! + \! C^{\zeta_{o}}_{-}(w_{+}(\zeta)), \quad \zeta \! \in \! D,
\end{equation*}
where $\id$ is the identity operator on $\mathcal{L}^{2}_{\mathrm{M}_{2}
(\mathbb{C})}(D)$, then the solution of the matrix {\rm RHP} $(\mathfrak{W}
(\zeta),\upsilon (\zeta),D)$ is given by
\begin{equation*}
\mathfrak{W}(\zeta) \! = \! \mathrm{I} \! + \! \int\nolimits_{D}\dfrac{(\zeta 
\! - \! \zeta_{o}) \mu_{\zeta_{o}}(\xi)w(\xi)}{(\xi \! -\! \zeta_{o})(\xi \! - 
\! \zeta)} \, \dfrac{\md \xi}{2 \pi \mi}, \quad \zeta \! \in \! \mathbb{C} 
\setminus D,
\end{equation*}
where $\mu_{\zeta_{o}}(\zeta) \! := \! ((\id \! - \! C^{\zeta_{o}}_{w})^{-1} 
\mathrm{I})(\zeta)$.\footnote{The linear singular integral equation for 
$\mu_{\zeta_{o}}$ stated in this Lemma~\ref{lem5.2} is well defined 
provided that $w_{\pm} \! \in \! \cap_{p=2,\infty} \mathcal{L}^{2}_{
\mathrm{M}_{2}(\mathbb{C})}(D)$.}
\end{ccccc}
\begin{eeeee} \label{ouhzphreton} 
\textsl{A scalar- or matrix-valued {\rm RHP} can be defined on the 
Riemann sphere, $\overline{\mathbb{C}}$, in such a way that it is 
invariant in $L^{2}$ under the group of M\"{o}bius transformations 
$z \! \mapsto \! (az \! + \! b)(cz \! + \! d)^{-1}$, $ad \! - \! bc \! 
\neq \! 0$, on $\overline{\mathbb{C}}$, with the spherical metric, 
$2(1 \! + \! \lvert z \rvert^{2})^{-1} \, \lvert \md z \rvert$, induced by 
the stereographic projection {}from the unit sphere $\mathbb{S}^{2}$ 
to $\overline{\mathbb{C}}$. In general, one does not have a similar 
$L^{p}$ invariance for values of $p$ other than $2$$;$ fortunately, 
however, the $L^{2}$ theory is sufficient for the purposes of this lecture 
note. Although it is possible to define a scalar- or matrix-valued 
{\rm RHP} on a fairly arbitrary contour on the sphere, in practice, one 
works only with contours on which the Cauchy operators (see below) are 
bounded in $L^{2}$$:$ such contours have been completely characterised 
(see, for example, {\rm \cite{geedave}}$)$$;$ but, a description of such a 
characterisation is beyond the scope of this monograph. Cauchy operators 
are indispensable in the study of {\rm RHP}s, as scalar- or matrix-valued 
{\rm RHP}s are equivalent to some singular integral equation (or a 
system of singular integral equations) with associated Cauchy operators 
on the contour. When the point at infinity is on the contour, the associated 
Cauchy operators are not bounded in $L^{2}$ with the spherical metric. 
One expedient solution to this problem is to use a M\"{o}bius transformation 
so that the point at infinity is no longer on the contour$;$ another 
approach is to replace the `standard' Cauchy kernel, $\tfrac{1}{\zeta - z} 
\, \tfrac{\md \zeta}{2 \pi \mi}$, with the `generalised' Cauchy kernel, 
$\tfrac{z - \dot{z}_{0}}{(\zeta - \dot{z}_{0})(\zeta - z)} \, \tfrac{\md 
\zeta}{2 \pi \mi}$, where $\dot{z}_{0}$ is any point not on the 
contour (see below$)$$:$ the point $\dot{z}_{0}$ plays the same 
r\^{o}le that the point at infinity does in the case of the standard 
Cauchy operator. One can show that, with this generalised Cauchy 
kernel, the associated Cauchy operators are bounded in $L^{2}$ 
with the spherical metric. When one studies {\rm ORF}s with poles on 
the extended real line, $\overline{\mathbb{R}}$, each typical {\rm RHP} 
is normalised at one of these poles$;$ hence, one needs Cauchy kernels 
normalised at a point (or points) on the contour, namely, the extended 
real line. The generalised Cauchy kernel is clearly bounded {}from the 
weighted $L^{2}$ space, with weight $\lvert z \! - \! \dot{z}_{0} \rvert^{-2}$, 
to itself: when $\dot{z}_{0}$ is not on the contour, this weight gives rise 
to an equivalent norm as given by the spherical metric$;$ but, when 
$\dot{z}_{0}$ is on the contour, the theory relies on the fact that the jump 
matrix, herein denoted as $v$, of the associated {\rm RHP} is normalised 
at $\dot{z}_{0}$ so that $v \! - \! \mathrm{I}$ has sufficient decay to 
balance the growth of the weight $\lvert z \! - \! \dot{z}_{0} \rvert^{-2}$ 
at $\dot{z}_{0}$. The singular integral operators involved are essentially 
compositions of Cauchy operators and the multipliers by $v \! - \! 
\mathrm{I}$, and, therefore, are bounded in $L^{2}$.}
\end{eeeee}
Recall that, for $n \! \in \! \mathbb{N}$ and $k \! \in \! \lbrace 1,2,
\dotsc,K \rbrace$ such that $\alpha_{p_{\mathfrak{s}}} \! := \! \alpha_{k} 
\! = \! \infty$ (resp., $\alpha_{p_{\mathfrak{s}}} \! := \! \alpha_{k} 
\! \neq \! \infty)$, $\hat{\mathcal{R}} \colon \mathbb{C} \setminus 
\hat{\Sigma}_{\hat{\mathcal{R}}} \! \to \! \operatorname{SL}_{2}
(\mathbb{C})$ (resp., $\tilde{\mathcal{R}} \colon \mathbb{C} \setminus 
\tilde{\Sigma}_{\tilde{\mathcal{R}}} \! \to \! \operatorname{SL}_{2}
(\mathbb{C}))$, which solves the associated RHP $(\hat{\mathcal{R}}(z),
\hat{v}_{\hat{\mathcal{R}}}(z),\hat{\Sigma}_{\hat{\mathcal{R}}})$ (resp., 
$(\tilde{\mathcal{R}}(z),\tilde{v}_{\tilde{\mathcal{R}}}(z),\tilde{\Sigma}_{
\tilde{\mathcal{R}}}))$ stated in item~\pmb{(1)} (resp., item~\pmb{(2)}) 
of Lemma~\ref{lem4.10}, is normalised at $\alpha_{p_{\mathfrak{s}}} \! 
:= \! \alpha_{k} \! = \! \infty$ (resp., $\alpha_{p_{\mathfrak{s}}} \! := \! 
\alpha_{k} \! \neq \! \infty)$, that is, $\hat{\mathcal{R}}(\infty) \! = \! 
\mathrm{I}$ (resp., $\tilde{\mathcal{R}}(\alpha_{k}) \! = \! \mathrm{I})$. 
Excising {}from the specification of the associated RHP $(\hat{\mathcal{R}}
(z),\hat{v}_{\hat{\mathcal{R}}}(z),\hat{\Sigma}_{\hat{\mathcal{R}}})$ (resp., 
$(\tilde{\mathcal{R}}(z),\tilde{v}_{\tilde{\mathcal{R}}}(z),\tilde{\Sigma}_{
\tilde{\mathcal{R}}}))$ the oriented skeletons on which the corresponding 
jump matrix, $\hat{v}_{\hat{\mathcal{R}}}(z)$ (resp., $\tilde{v}_{\tilde{
\mathcal{R}}}(z))$, is equal to $\mathrm{I}$, in particular (cf. 
Lemma~\ref{lem4.10}), the oriented skeleton $\hat{\Sigma}_{\hat{\mathcal{R}}} 
\setminus \cup_{j=1}^{5} \hat{\Sigma}_{p}^{j}$ (resp., $\tilde{\Sigma}_{
\tilde{\mathcal{R}}} \setminus \cup_{j=1}^{5} \tilde{\Sigma}_{p}^{j})$, and 
defining $\hat{\Sigma}_{\hat{\mathcal{R}}}^{\sharp} \! := \! \hat{\Sigma}_{
\hat{\mathcal{R}}} \setminus (\hat{\Sigma}_{\hat{\mathcal{R}}} \setminus 
\cup_{j=1}^{5} \hat{\Sigma}_{p}^{j})$ {}\footnote{See Figure~\ref{fighatexcis}.} 
(resp., $\tilde{\Sigma}_{\tilde{\mathcal{R}}}^{\sharp} \! := \! \tilde{\Sigma}_{
\tilde{\mathcal{R}}} \setminus (\tilde{\Sigma}_{\tilde{\mathcal{R}}} \setminus 
\cup_{j=1}^{5} \tilde{\Sigma}_{p}^{j}))$,\footnote{See Figure~\ref{figtilexcis}.} 
one arrives at the associated, equivalent RHP $(\hat{\mathcal{R}}(z),
\hat{v}_{\hat{\mathcal{R}}}(z),\hat{\Sigma}_{\hat{\mathcal{R}}}^{\sharp})$ 
(resp., $(\tilde{\mathcal{R}}(z),\tilde{v}_{\tilde{\mathcal{R}}}(z),\tilde{\Sigma}_{
\tilde{\mathcal{R}}}^{\sharp}))$ for $\hat{\mathcal{R}} \colon \mathbb{C} 
\setminus \hat{\Sigma}_{\hat{\mathcal{R}}}^{\sharp} \! \to \! 
\operatorname{SL}_{2}(\mathbb{C})$ (resp., $\tilde{\mathcal{R}} \colon 
\mathbb{C} \setminus \tilde{\Sigma}_{\tilde{\mathcal{R}}}^{\sharp} \! \to \! 
\operatorname{SL}_{2}(\mathbb{C}))$.\footnote{The respective normalisations 
at $\alpha_{p_{\mathfrak{s}}} \! := \! \alpha_{k}$ remain unchanged.}
\begin{figure}[tbh]
\begin{center}

\end{center}
\vspace{-1.05cm}
\caption{Oriented skeleton $\hat{\Sigma}_{\hat{\mathcal{R}}}^{\sharp} \! := \! 
\hat{\Sigma}_{\hat{\mathcal{R}}} \setminus (\hat{\Sigma}_{\hat{\mathcal{R}}} 
\setminus \cup_{j=1}^{5} \hat{\Sigma}_{p}^{j})$.}\label{fighatexcis}
\end{figure}
\begin{figure}[tbh]
\begin{center}
%
\end{center}
\vspace{-1.05cm}
\caption{Oriented skeleton $\tilde{\Sigma}_{\tilde{\mathcal{R}}}^{\sharp} 
\! := \! \tilde{\Sigma}_{\tilde{\mathcal{R}}} \setminus (\tilde{\Sigma}_{
\tilde{\mathcal{R}}} \setminus \cup_{j=1}^{5} \tilde{\Sigma}_{p}^{j})$.}\label{figtilexcis}
\end{figure}
Via the BC construction \cite{bealscoif}, for $n \! \in \! \mathbb{N}$ and $k 
\! \in \! \lbrace 1,2,\dotsc,K \rbrace$ such that $\alpha_{p_{\mathfrak{s}}} 
\! := \! \alpha_{k} \! = \! \infty$, write, for $\hat{v}_{\hat{\mathcal{R}}} \colon 
\hat{\Sigma}_{\hat{\mathcal{R}}}^{\sharp} \! \to \! \operatorname{SL}_{2}
(\mathbb{C})$, the bounded algebraic factorisation
\begin{equation*}
\hat{v}_{\hat{\mathcal{R}}}(z) \! := \! (\mathrm{I} \! - \! w_{-}^{\Sigma_{\hat{
\mathcal{R}}}}(z))^{-1}(\mathrm{I} \! + \! w_{+}^{\Sigma_{\hat{\mathcal{R}}}}
(z)), \quad z \! \in \! \hat{\Sigma}_{\hat{\mathcal{R}}}^{\sharp},
\end{equation*}
and, for $n \! \in \! \mathbb{N}$ and $k \! \in \! \lbrace 1,2,\dotsc,K \rbrace$ 
such that $\alpha_{p_{\mathfrak{s}}} \! := \! \alpha_{k} \! \neq \! \infty$, write, 
for $\tilde{v}_{\tilde{\mathcal{R}}} \colon \tilde{\Sigma}_{\tilde{\mathcal{R}}}^{
\sharp} \! \to \! \operatorname{SL}_{2}(\mathbb{C})$, the bounded algebraic 
factorisation
\begin{equation*}
\tilde{v}_{\tilde{\mathcal{R}}}(z) \! := \! (\mathrm{I} \! - \! w_{-}^{\Sigma_{\tilde{
\mathcal{R}}}}(z))^{-1}(\mathrm{I} \! + \! w_{+}^{\Sigma_{\tilde{\mathcal{R}}}}
(z)), \quad z \! \in \! \tilde{\Sigma}_{\tilde{\mathcal{R}}}^{\sharp}.
\end{equation*}
Taking, for $n \! \in \! \mathbb{N}$ and $k \! \in \! \lbrace 1,2,\dotsc,K \rbrace$ 
such that $\alpha_{p_{\mathfrak{s}}} \! := \! \alpha_{k} \! = \! \infty$ (resp., 
$\alpha_{p_{\mathfrak{s}}} \! := \! \alpha_{k} \! \neq \! \infty)$, the so-called 
`trivial factorisation' $w^{\Sigma_{\hat{\mathcal{R}}}}_{-}(z) \! \equiv \! 
\left(
\begin{smallmatrix}
0 & 0 \\
0 & 0
\end{smallmatrix}
\right)$ (resp., $w^{\Sigma_{\tilde{\mathcal{R}}}}_{-}(z) \! \equiv \! 
\left(
\begin{smallmatrix}
0 & 0 \\
0 & 0
\end{smallmatrix}
\right))$,\footnote{See, in particular, pp.~293--294 in \cite{a64}; see, also, 
\cite{xinzoo2,xinzoo3}.} whence $\hat{v}_{\hat{\mathcal{R}}}(z) \! = \! \mathrm{I} \! + 
\! w^{\Sigma_{\hat{\mathcal{R}}}}_{+}(z)$, $z \! \in \! \hat{\Sigma}_{\hat{\mathcal{
R}}}^{\sharp}$ (resp., $\tilde{v}_{\tilde{\mathcal{R}}}(z) \! = \! \mathrm{I} \! + 
\! w^{\Sigma_{\tilde{\mathcal{R}}}}_{+}(z)$, $z \! \in \! \tilde{\Sigma}_{\tilde{
\mathcal{R}}}^{\sharp})$, it follows {}from Lemma~\ref{lem5.2} that: (i) for $n 
\! \in \! \mathbb{N}$ and $k \! \in \! \lbrace 1,2,\dotsc,K \rbrace$ such that 
$\alpha_{p_{\mathfrak{s}}} \! := \! \alpha_{k} \! = \! \infty$, upon normalising 
the associated Cauchy operator(s) at $\alpha_{p_{\mathfrak{s}}} \! := \! \alpha_{k} 
\! = \! \infty$, the integral representation for the unique solution of the 
corresponding RHP $(\hat{\mathcal{R}}(z),\hat{v}_{\hat{\mathcal{R}}}(z),\hat{
\Sigma}_{\hat{\mathcal{R}}}^{\sharp})$ is
\begin{equation} \label{eqsek5a} 
\hat{\mathcal{R}}(z) \! = \! \mathrm{I} \! + \! \int_{\hat{\Sigma}_{\hat{
\mathcal{R}}}^{\sharp}} \dfrac{\mu^{\Sigma_{\hat{\mathcal{R}}}}(\xi)w^{
\Sigma_{\hat{\mathcal{R}}}}_{+}(\xi)}{\xi \! - \! z} \, \dfrac{\md \xi}{2 
\pi \mi}, \quad z \! \in \! \mathbb{C} \setminus \hat{\Sigma}_{
\hat{\mathcal{R}}}^{\sharp},
\end{equation}
where $\mu^{\Sigma_{\hat{\mathcal{R}}}}(z) \! = \! \hat{\mathcal{R}}_{-}(z) 
\! = \! \hat{\mathcal{R}}_{+}(z)(\mathrm{I} \! + \! w^{\Sigma_{\hat{
\mathcal{R}}}}_{+}(z))^{-1}$ solves the linear singular integral equation
\begin{equation*}
((\id \! - \! C_{w^{\Sigma_{\hat{\mathcal{R}}}}}) \mu^{\Sigma_{\hat{\mathcal{
R}}}})(z) \! = \! \mathrm{I}, \quad z \! \in \! \hat{\Sigma}_{\hat{\mathcal{R}}}^{
\sharp},
\end{equation*}
with
\begin{equation*}
\mathcal{L}^{2}_{\mathrm{M}_{2}(\mathbb{C})}(\hat{\Sigma}_{\hat{\mathcal{
R}}}^{\sharp}) \! \ni \! f \! \mapsto \! C_{w^{\Sigma_{\hat{\mathcal{R}}}}}f \! 
:= \! \hat{C}_{-}(fw^{\Sigma_{\hat{\mathcal{R}}}}_{+}),
\end{equation*}
where
\begin{equation*}
\mathcal{L}^{2}_{\mathrm{M}_{2}(\mathbb{C})}(\hat{\Sigma}_{\hat{\mathcal{R}}}^{
\sharp}) \! \ni \! g \! \mapsto \! (\hat{C}_{\pm}g)(z) := \lim_{\underset{z^{\prime} 
\, \in \, \pm \, \mathrm{side} \, \mathrm{of} \, \hat{\Sigma}_{\hat{\mathcal{R}}}^{
\sharp}}{z^{\prime} \to z}} \int_{\hat{\Sigma}_{\hat{\mathcal{R}}}^{\sharp}} 
\dfrac{g(\xi)}{\xi \! - \! z^{\prime}} \, \dfrac{\md \xi}{2 \pi \mi};
\end{equation*}
and (ii) for $n \! \in \! \mathbb{N}$ and $k \! \in \! \lbrace 1,2,\dotsc,K \rbrace$ 
such that $\alpha_{p_{\mathfrak{s}}} \! := \! \alpha_{k} \! \neq \! \infty$, upon 
normalising the associated Cauchy operator(s) at $\alpha_{p_{\mathfrak{s}}} \! := 
\! \alpha_{k} \! \neq \! \infty$, the integral representation for the unique solution 
of the corresponding RHP $(\tilde{\mathcal{R}}(z),\tilde{v}_{\tilde{\mathcal{R}}}(z),
\tilde{\Sigma}_{\tilde{\mathcal{R}}}^{\sharp})$ is
\begin{equation} \label{eqsek5b} 
\tilde{\mathcal{R}}(z) \! = \! \mathrm{I} \! + \! \int_{\tilde{\Sigma}_{\tilde{
\mathcal{R}}}^{\sharp}} \dfrac{(z \! - \! \alpha_{k}) \mu^{\Sigma_{\tilde{
\mathcal{R}}}}(\xi)w^{\Sigma_{\tilde{\mathcal{R}}}}_{+}(\xi)}{(\xi \! - \! \alpha_{k})
(\xi \! - \! z)} \, \dfrac{\md \xi}{2 \pi \mi}, \quad z \! \in \! \mathbb{C} 
\setminus \tilde{\Sigma}_{\tilde{\mathcal{R}}}^{\sharp},
\end{equation}
where $\mu^{\Sigma_{\tilde{\mathcal{R}}}}(z) \! = \! \tilde{\mathcal{R}}_{-}(z) 
\! = \! \tilde{\mathcal{R}}_{+}(z)(\mathrm{I} \! + \! w^{\Sigma_{\tilde{\mathcal{
R}}}}_{+}(z))^{-1}$ solves the linear singular integral equation
\begin{equation*}
((\id \! - \! C_{w^{\Sigma_{\tilde{\mathcal{R}}}}}) \mu^{\Sigma_{\tilde{\mathcal{
R}}}})(z) \! = \! \mathrm{I}, \quad z \! \in \! \tilde{\Sigma}_{\tilde{\mathcal{
R}}}^{\sharp},
\end{equation*}
with
\begin{equation*}
\mathcal{L}^{2}_{\mathrm{M}_{2}(\mathbb{C})}(\tilde{\Sigma}_{\tilde{\mathcal{
R}}}^{\sharp}) \! \ni \! f \! \mapsto \! C_{w^{\Sigma_{\tilde{\mathcal{R}}}}}f 
\! := \! \tilde{C}_{-}(fw^{\Sigma_{\tilde{\mathcal{R}}}}_{+}),
\end{equation*}
where
\begin{equation*}
\mathcal{L}^{2}_{\mathrm{M}_{2}(\mathbb{C})}(\tilde{\Sigma}_{\tilde{\mathcal{
R}}}^{\sharp}) \! \ni \! g \! \mapsto \! (\tilde{C}_{\pm}g)(z) := \lim_{\underset{
z^{\prime} \, \in \, \pm \, \mathrm{side} \, \mathrm{of} \, \tilde{\Sigma}_{
\tilde{\mathcal{R}}}^{\sharp}}{z^{\prime} \to z}} \int_{\tilde{\Sigma}_{\tilde{
\mathcal{R}}}^{\sharp}} \dfrac{(z^{\prime} \! - \! \alpha_{k})g(\xi)}{(\xi \! - \! 
\alpha_{k})(\xi \! - \! z^{\prime})} \, \dfrac{\md \xi}{2 \pi \mi}.
\end{equation*}
\begin{bbbbb} \label{propo5.1} 
For $n \! \in \! \mathbb{N}$ and $k \! \in \! \lbrace 1,2,\dotsc,K \rbrace$ 
such that $\alpha_{p_{\mathfrak{s}}} \! := \! \alpha_{k} \! = \! \infty$ (resp., 
$\alpha_{p_{\mathfrak{s}}} \! := \! \alpha_{k} \! \neq \! \infty)$, let 
$\hat{\mathcal{R}} \colon \mathbb{C} \setminus \hat{\Sigma}_{\hat{
\mathcal{R}}}^{\sharp} \! \to \! \operatorname{SL}_{2}(\mathbb{C})$ 
(resp., $\tilde{\mathcal{R}} \colon \mathbb{C} \setminus \tilde{\Sigma}_{
\tilde{\mathcal{R}}}^{\sharp} \! \to \! \operatorname{SL}_{2}(\mathbb{C}))$ solve 
the equivalent {\rm RHP} $(\hat{\mathcal{R}}(z),\hat{v}_{\hat{\mathcal{R}}}(z),
\hat{\Sigma}_{\hat{\mathcal{R}}}^{\sharp})$ (resp., $(\tilde{\mathcal{R}}(z),
\tilde{v}_{\tilde{\mathcal{R}}}(z),\tilde{\Sigma}_{\tilde{\mathcal{R}}}^{\sharp}))$, 
where, in particular, the asymptotic expansion, in the double-scaling limit 
$\mathscr{N},n \! \to \! \infty$ such that $z_{o} \! = \! 1 \! + \! o(1)$, of 
$\hat{v}_{\hat{\mathcal{R}}}(z)$ (resp., $\tilde{v}_{\tilde{\mathcal{R}}}(z))$ for 
$z \! \in \! \hat{\Sigma}_{p,j}^{5} \! := \! \partial \hat{\mathbb{U}}_{\hat{
\delta}_{\hat{b}_{j-1}}} \cup \partial \hat{\mathbb{U}}_{\hat{\delta}_{\hat{a}_{j}}}$ 
(resp., $z \! \in \! \tilde{\Sigma}_{p,j}^{5} \! := \! \partial \tilde{\mathbb{
U}}_{\tilde{\delta}_{\tilde{b}_{j-1}}} \cup \partial \tilde{\mathbb{U}}_{\tilde{
\delta}_{\tilde{a}_{j}}})$, $j \! = \! 1,2,\dotsc,N \! + \! 1$, is given in 
Lemma~\ref{lem5.1}, Equations~\eqref{eqhtvee6} and~\eqref{eqhtvee7} 
(resp., Equations~\eqref{eqtlvee13} and~\eqref{eqtlvee14}$)$. Then$:$ 
{\rm \pmb{(1)}} for $n \! \in \! \mathbb{N}$ and $k \! \in \! \lbrace 1,2,
\dotsc,K \rbrace$ such that $\alpha_{p_{\mathfrak{s}}} \! := \! \alpha_{k} \! 
= \! \infty$, and for $z \! \in \! \partial \hat{\mathbb{U}}_{\hat{\delta}_{
\hat{b}_{j-1}}} \cup \partial \hat{\mathbb{U}}_{\hat{\delta}_{\hat{a}_{j}}}$, 
$j \! = \! 1,2,\dotsc,N \! + \! 1$, with $w^{\Sigma_{\hat{\mathcal{R}}}}_{+}
(z) \! = \! \hat{v}_{\hat{\mathcal{R}}}(z) \! - \! \mathrm{I}$,
\begin{align} \label{eqprophatb} 
w_{+}^{\Sigma_{\hat{\mathcal{R}}}}(z) \underset{\underset{z_{o}=1+o(1)}{
\mathscr{N},n \to \infty}}{=}& \, \dfrac{1}{((n \! - \! 1)K \! + \! k)} \left(
\dfrac{(\hat{\alpha}_{0}(\hat{b}_{j-1}))^{-1}}{(z \! - \! \hat{b}_{j-1})^{2}} 
\hat{\boldsymbol{\mathrm{A}}}(\hat{b}_{j-1}) \! + \! \dfrac{(\hat{\alpha}_{0}
(\hat{b}_{j-1}))^{-2}}{z \! - \! \hat{b}_{j-1}} \left(\hat{\alpha}_{0}(\hat{b}_{j-1}) 
\hat{\boldsymbol{\mathrm{B}}}(\hat{b}_{j-1}) \! - \! \hat{\alpha}_{1}
(\hat{b}_{j-1}) \hat{\boldsymbol{\mathrm{A}}}(\hat{b}_{j-1}) \right) \right. 
\nonumber \\
&\left. \, +(\hat{\alpha}_{0}(\hat{b}_{j-1}))^{-2} \left(\hat{\alpha}_{0}
(\hat{b}_{j-1}) \left(\left(\dfrac{\hat{\alpha}_{1}(\hat{b}_{j-1})}{\hat{\alpha}_{0}
(\hat{b}_{j-1})} \right)^{2} \! - \! \dfrac{\hat{\alpha}_{2}(\hat{b}_{j-1})}{\hat{
\alpha}_{0}(\hat{b}_{j-1})} \right) \hat{\boldsymbol{\mathrm{A}}}(\hat{b}_{j-1}) 
\! - \! \hat{\alpha}_{1}(\hat{b}_{j-1}) \hat{\boldsymbol{\mathrm{B}}}
(\hat{b}_{j-1}) \! + \! \hat{\alpha}_{0}(\hat{b}_{j-1}) \hat{\boldsymbol{
\mathrm{C}}}(\hat{b}_{j-1}) \right) \right) \nonumber \\
& \, +\mathcal{O} \left(\dfrac{1}{((n \! - \! 1)K \! + \! k)} \sum_{m=1}^{
\infty} \hat{c}_{m}^{\triangleright}(n,k,z_{o};\hat{b}_{j-1})(z \! - \! 
\hat{b}_{j-1})^{m} \right) \nonumber \\
& \, + \mathcal{O} \left(\dfrac{1}{((n \! - \! 1)K \! + \! k)^{2}(z \! - \! 
\hat{b}_{j-1})^{3}} \sum_{m=0}^{\infty} \hat{c}_{m}^{\triangleleft}(n,k,z_{o};
\hat{b}_{j-1})(z \! - \! \hat{b}_{j-1})^{m} \right), \quad z \! \in \! \partial 
\hat{\mathbb{U}}_{\hat{\delta}_{\hat{b}_{j-1}}}, \quad j \! = \! 1,2,\dotsc,
N \! + \! 1,
\end{align}
and
\begin{align} \label{eqprophata} 
w_{+}^{\Sigma_{\hat{\mathcal{R}}}}(z) \underset{\underset{z_{o}=1+o(1)}{
\mathscr{N},n \to \infty}}{=}& \, \dfrac{1}{((n \! - \! 1)K \! + \! k)} \left(
\dfrac{(\hat{\alpha}_{0}(\hat{a}_{j}))^{-1}}{(z \! - \! \hat{a}_{j})^{2}} \hat{
\boldsymbol{\mathrm{A}}}(\hat{a}_{j}) \! + \! \dfrac{(\hat{\alpha}_{0}
(\hat{a}_{j}))^{-2}}{z \! - \! \hat{a}_{j}} \left(\hat{\alpha}_{0}(\hat{a}_{j}) 
\hat{\boldsymbol{\mathrm{B}}}(\hat{a}_{j}) \! - \! \hat{\alpha}_{1}(\hat{a}_{j}) 
\hat{\boldsymbol{\mathrm{A}}}(\hat{a}_{j}) \right) \right. \nonumber \\
&\left. \, +(\hat{\alpha}_{0}(\hat{a}_{j}))^{-2} \left(\hat{\alpha}_{0}(\hat{a}_{j}) 
\left(\left(\dfrac{\hat{\alpha}_{1}(\hat{a}_{j})}{\hat{\alpha}_{0}(\hat{a}_{j})} 
\right)^{2} \! - \! \dfrac{\hat{\alpha}_{2}(\hat{a}_{j})}{\hat{\alpha}_{0}
(\hat{a}_{j})} \right) \hat{\boldsymbol{\mathrm{A}}}(\hat{a}_{j}) \! - \! 
\hat{\alpha}_{1}(\hat{a}_{j}) \hat{\boldsymbol{\mathrm{B}}}(\hat{a}_{j}) \! + 
\! \hat{\alpha}_{0}(\hat{a}_{j}) \hat{\boldsymbol{\mathrm{C}}}(\hat{a}_{j}) 
\right) \right) \nonumber \\
& \, +\mathcal{O} \left(\dfrac{1}{((n \! - \! 1)K \! + \! k)} \sum_{m=1}^{
\infty} \hat{c}_{m}^{\triangleright}(n,k,z_{o};\hat{a}_{j})(z \! - \! \hat{a}_{j})^{m} 
\right) \nonumber \\
& \, + \mathcal{O} \left(\dfrac{1}{((n \! - \! 1)K \! + \! k)^{2}(z \! - \! 
\hat{a}_{j})^{3}} \sum_{m=0}^{\infty} \hat{c}_{m}^{\triangleleft}(n,k,z_{o};
\hat{a}_{j})(z \! - \! \hat{a}_{j})^{m} \right), \quad z \! \in \! \partial  
\hat{\mathbb{U}}_{\hat{\delta}_{\hat{a}_{j}}}, \quad j \! = \! 1,2,\dotsc,
N \! + \! 1,
\end{align}
where, for $j \! = \! 1,2,\dotsc,N \! + \! 1$,
\begin{gather}
\hat{\boldsymbol{\mathrm{A}}}(\hat{b}_{j-1}) \! := \! 
\begin{pmatrix}
\tilde{\mathfrak{m}}^{\raise-0.5ex\hbox{$\scriptstyle \infty$}}_{11} & 0 \\
0 & \tilde{\mathfrak{m}}^{\raise-0.5ex\hbox{$\scriptstyle \infty$}}_{22}
\end{pmatrix} 
\begin{pmatrix}
\hat{\mathbb{A}}_{11}(\hat{b}_{j-1}) & \hat{\mathbb{A}}_{12}(\hat{b}_{j-1}) \\
\hat{\mathbb{A}}_{21}(\hat{b}_{j-1}) & \hat{\mathbb{A}}_{22}(\hat{b}_{j-1})
\end{pmatrix} 
\begin{pmatrix}
\tilde{\mathfrak{m}}^{\raise-0.5ex\hbox{$\scriptstyle \infty$}}_{22} & 0 \\
0 & \tilde{\mathfrak{m}}^{\raise-0.5ex\hbox{$\scriptstyle \infty$}}_{11}
\end{pmatrix} \me^{\mi ((n-1)K+k) \hat{\mho}_{j-1}}, \label{eqprop1} \\
\hat{\boldsymbol{\mathrm{B}}}(\hat{b}_{j-1}) \! := \! 
\begin{pmatrix}
\tilde{\mathfrak{m}}^{\raise-0.5ex\hbox{$\scriptstyle \infty$}}_{11} & 0 \\
0 & \tilde{\mathfrak{m}}^{\raise-0.5ex\hbox{$\scriptstyle \infty$}}_{22}
\end{pmatrix} 
\begin{pmatrix}
\hat{\mathbb{B}}_{11}(\hat{b}_{j-1}) & \hat{\mathbb{B}}_{12}(\hat{b}_{j-1}) \\
\hat{\mathbb{B}}_{21}(\hat{b}_{j-1}) & \hat{\mathbb{B}}_{22}(\hat{b}_{j-1})
\end{pmatrix} 
\begin{pmatrix}
\tilde{\mathfrak{m}}^{\raise-0.5ex\hbox{$\scriptstyle \infty$}}_{22} & 0 \\
0 & \tilde{\mathfrak{m}}^{\raise-0.5ex\hbox{$\scriptstyle \infty$}}_{11}
\end{pmatrix} \me^{\mi ((n-1)K+k) \hat{\mho}_{j-1}}, \label{eqprop2} \\
\hat{\boldsymbol{\mathrm{C}}}(\hat{b}_{j-1}) \! := \! 
\begin{pmatrix}
\tilde{\mathfrak{m}}^{\raise-0.5ex\hbox{$\scriptstyle \infty$}}_{11} & 0 \\
0 & \tilde{\mathfrak{m}}^{\raise-0.5ex\hbox{$\scriptstyle \infty$}}_{22}
\end{pmatrix} 
\begin{pmatrix}
\hat{\mathbb{C}}_{11}(\hat{b}_{j-1}) & \hat{\mathbb{C}}_{12}(\hat{b}_{j-1}) \\
\hat{\mathbb{C}}_{21}(\hat{b}_{j-1}) & \hat{\mathbb{C}}_{22}(\hat{b}_{j-1})
\end{pmatrix} 
\begin{pmatrix}
\tilde{\mathfrak{m}}^{\raise-0.5ex\hbox{$\scriptstyle \infty$}}_{22} & 0 \\
0 & \tilde{\mathfrak{m}}^{\raise-0.5ex\hbox{$\scriptstyle \infty$}}_{11}
\end{pmatrix} \me^{\mi ((n-1)K+k) \hat{\mho}_{j-1}}, \label{eqprop3}
\end{gather}
with $\tilde{\mathfrak{m}}^{\raise-0.5ex\hbox{$\scriptstyle \infty$}}_{11}$ 
and $\tilde{\mathfrak{m}}^{\raise-0.5ex\hbox{$\scriptstyle \infty$}}_{22}$ 
given by Equation~\eqref{eqmaininf8}, $\hat{\mho}_{m}$, $m \! = \! 
0,1,\dotsc,N \! + \! 1$, defined in the corresponding item~{\rm (ii)} of 
Remark~\ref{rem4.4},
\begin{gather}
\hat{\mathbb{A}}_{11}(\hat{b}_{j-1}) \! = \! -\hat{\mathbb{A}}_{22}
(\hat{b}_{j-1}) \! = \! -s_{1} \hat{\kappa}_{1}(\hat{b}_{j-1}) \hat{\kappa}_{2}
(\hat{b}_{j-1})(\hat{\mathfrak{Q}}_{0}(\hat{b}_{j-1}))^{-1}, \label{eqprop5} \\
\hat{\mathbb{A}}_{12}(\hat{b}_{j-1}) \! = \! -\mi s_{1}(\hat{\kappa}_{1}
(\hat{b}_{j-1}))^{2}(\hat{\mathfrak{Q}}_{0}(\hat{b}_{j-1}))^{-1}, \qquad 
\hat{\mathbb{A}}_{21}(\hat{b}_{j-1}) \! = \! -\mi s_{1}(\hat{\kappa}_{2}
(\hat{b}_{j-1}))^{2}(\hat{\mathfrak{Q}}_{0}(\hat{b}_{j-1}))^{-1}, \label{eqprop6}
\end{gather}
\begin{align}
\hat{\mathbb{B}}_{11}(\hat{b}_{j-1}) \! = \! -\hat{\mathbb{B}}_{22}
(\hat{b}_{j-1}) =& \, \hat{\kappa}_{1}(\hat{b}_{j-1}) \hat{\kappa}_{2}
(\hat{b}_{j-1}) \left(-s_{1}(\hat{\mathfrak{Q}}_{0}(\hat{b}_{j-1}))^{-1} 
\left(\hat{\daleth}^{1}_{1}(\hat{b}_{j-1}) \! + \! \hat{\daleth}^{1}_{-1}
(\hat{b}_{j-1}) \! - \! \hat{\mathfrak{Q}}_{1}(\hat{b}_{j-1})
(\hat{\mathfrak{Q}}_{0}(\hat{b}_{j-1}))^{-1} \right) \right. \nonumber \\
-&\left. \, t_{1} \left(\hat{\mathfrak{Q}}_{0}(\hat{b}_{j-1}) \! + \! 
(\hat{\mathfrak{Q}}_{0}(\hat{b}_{j-1}))^{-1} \hat{\aleph}^{1}_{1}
(\hat{b}_{j-1}) \hat{\aleph}^{1}_{-1}(\hat{b}_{j-1}) \right) \! + \! \mi 
(s_{1} \! + \! t_{1}) \left(\hat{\aleph}^{1}_{-1}(\hat{b}_{j-1}) \! - \! 
\hat{\aleph}^{1}_{1}(\hat{b}_{j-1}) \right) \right), \label{eqprop7} \\
\hat{\mathbb{B}}_{12}(\hat{b}_{j-1}) =& \, (\hat{\kappa}_{1}
(\hat{b}_{j-1}))^{2} \left(-\mi s_{1}(\hat{\mathfrak{Q}}_{0}(\hat{b}_{j-1}))^{-1} 
\left(2 \hat{\daleth}^{1}_{1}(\hat{b}_{j-1}) \! - \! \hat{\mathfrak{Q}}_{1}
(\hat{b}_{j-1})(\hat{\mathfrak{Q}}_{0}(\hat{b}_{j-1}))^{-1} \right) \right. 
\nonumber \\
+&\left. \, \mi t_{1} \left(\hat{\mathfrak{Q}}_{0}(\hat{b}_{j-1}) \! - \! 
(\hat{\mathfrak{Q}}_{0}(\hat{b}_{j-1}))^{-1}(\hat{\aleph}^{1}_{1}(\hat{b}_{j
-1}))^{2} \right) \! + \! 2(s_{1} \! - \! t_{1}) \hat{\aleph}^{1}_{1}(\hat{b}_{j-1}) 
\right), \label{eqprop8} \\
\hat{\mathbb{B}}_{21}(\hat{b}_{j-1}) =& \, (\hat{\kappa}_{2}
(\hat{b}_{j-1}))^{2} \left(-\mi s_{1}(\hat{\mathfrak{Q}}_{0}(\hat{b}_{j-1}))^{-1} 
\left(2 \hat{\daleth}^{1}_{-1}(\hat{b}_{j-1}) \! - \! \hat{\mathfrak{Q}}_{1}
(\hat{b}_{j-1})(\hat{\mathfrak{Q}}_{0}(\hat{b}_{j-1}))^{-1} \right) \right. 
\nonumber \\
+&\left. \, \mi t_{1} \left(\hat{\mathfrak{Q}}_{0}(\hat{b}_{j-1}) \! - \! 
(\hat{\mathfrak{Q}}_{0}(\hat{b}_{j-1}))^{-1}(\hat{\aleph}^{1}_{-1}(\hat{b}_{j
-1}))^{2} \right) \! - \! 2(s_{1} \! - \! t_{1}) \hat{\aleph}^{1}_{-1}(\hat{b}_{j-1}) 
\right), \label{eqprop9} \\
\hat{\mathbb{C}}_{11}(\hat{b}_{j-1}) \! = \! -\hat{\mathbb{C}}_{22}
(\hat{b}_{j-1}) =& \, \hat{\kappa}_{1}(\hat{b}_{j-1}) \hat{\kappa}_{2}
(\hat{b}_{j-1}) \left(s_{1} \left(-\hat{\mathfrak{Q}}_{0}(\hat{b}_{j-1}) 
\hat{\aleph}^{1}_{1}(\hat{b}_{j-1}) \hat{\aleph}^{1}_{-1}(\hat{b}_{j-1}) \! - \! 
(\hat{\mathfrak{Q}}_{0}(\hat{b}_{j-1}))^{-3}(\hat{\mathfrak{Q}}_{1}(\hat{b}_{j
-1}))^{2} \right. \right. \nonumber \\
+&\left. \left. \, \dfrac{1}{2} \hat{\mathfrak{Q}}_{2}(\hat{b}_{j-1})
(\hat{\mathfrak{Q}}_{0}(\hat{b}_{j-1}))^{-2} \! + \! \hat{\mathfrak{Q}}_{1}
(\hat{b}_{j-1})(\hat{\mathfrak{Q}}_{0}(\hat{b}_{j-1}))^{-2} \left(\hat{
\daleth}^{1}_{1}(\hat{b}_{j-1}) \! + \! \hat{\daleth}^{1}_{-1}(\hat{b}_{j-1}) 
\right) \right. \right. \nonumber \\
-&\left. \left. \, (\hat{\mathfrak{Q}}_{0}(\hat{b}_{j-1}))^{-1} \left(
\hat{\gimel}^{1}_{1}(\hat{b}_{j-1}) \! + \! \hat{\gimel}^{1}_{-1}(\hat{b}_{j-1}) 
\! + \! \hat{\daleth}^{1}_{1}(\hat{b}_{j-1}) \hat{\daleth}^{1}_{-1}(\hat{b}_{j-1}) 
\right) \right) \! + \! t_{1} \left(-\hat{\mathfrak{Q}}_{1}(\hat{b}_{j-1}) \! - \! 
\hat{\mathfrak{Q}}_{0}(\hat{b}_{j-1}) \right. \right. \nonumber \\
\times&\left. \left. \, \left(\hat{\daleth}^{1}_{1}(\hat{b}_{j-1}) \! + \! 
\hat{\daleth}^{1}_{-1}(\hat{b}_{j-1}) \right) \! + \! \hat{\mathfrak{Q}}_{1}
(\hat{b}_{j-1})(\hat{\mathfrak{Q}}_{0}(\hat{b}_{j-1}))^{-2} \hat{\aleph}^{1}_{1}
(\hat{b}_{j-1}) \hat{\aleph}^{1}_{-1}(\hat{b}_{j-1}) \! - \! (\hat{\mathfrak{Q}}_{0}
(\hat{b}_{j-1}))^{-1} \right. \right. \nonumber \\
\times&\left. \left. \, \left(\hat{\aleph}^{1}_{1}(\hat{b}_{j-1}) \hat{\beth}^{1}_{-1}
(\hat{b}_{j-1}) \! + \! \hat{\aleph}^{1}_{-1}(\hat{b}_{j-1}) \hat{\beth}^{1}_{1}
(\hat{b}_{j-1}) \right) \right) \! + \! \mi (s_{1} \! + \! t_{1}) \left(\hat{\beth}^{1}_{-1}
(\hat{b}_{j-1}) \! - \! \hat{\beth}^{1}_{1}(\hat{b}_{j-1}) \right. \right. \nonumber \\
-&\left. \left. \, \hat{\aleph}^{1}_{1}(\hat{b}_{j-1}) \hat{\daleth}^{1}_{-1}
(\hat{b}_{j-1}) \! + \! \hat{\aleph}^{1}_{-1}(\hat{b}_{j-1}) \hat{\daleth}^{1}_{1}
(\hat{b}_{j-1}) \right) \right), \label{eqprop10} \\
\hat{\mathbb{C}}_{12}(\hat{b}_{j-1}) =& \, (\hat{\kappa}_{1}(\hat{b}_{j-1}))^{2} 
\left(\mi s_{1} \left(\hat{\mathfrak{Q}}_{0}(\hat{b}_{j-1})(\hat{\aleph}^{1}_{1}
(\hat{b}_{j-1}))^{2} \! - \! (\hat{\mathfrak{Q}}_{0}(\hat{b}_{j-1}))^{-3}
(\hat{\mathfrak{Q}}_{1}(\hat{b}_{j-1}))^{2} \right. \right. \nonumber \\
+&\left. \left. \, \dfrac{1}{2} \hat{\mathfrak{Q}}_{2}(\hat{b}_{j-1})
(\hat{\mathfrak{Q}}_{0}(\hat{b}_{j-1}))^{-2} \! + \! 2 \hat{\mathfrak{Q}}_{1}
(\hat{b}_{j-1})(\hat{\mathfrak{Q}}_{0}(\hat{b}_{j-1}))^{-2} \hat{\daleth}^{1}_{1}
(\hat{b}_{j-1}) \! - \! (\hat{\mathfrak{Q}}_{0}(\hat{b}_{j-1}))^{-1} \right. \right. 
\nonumber \\
\times&\left. \left. \, \left(2 \hat{\gimel}^{1}_{1}(\hat{b}_{j-1}) \! + \! 
(\hat{\daleth}^{1}_{1}(\hat{b}_{j-1}))^{2} \right) \right) \! + \! \mi t_{1} 
\left(2 \hat{\mathfrak{Q}}_{0}(\hat{b}_{j-1}) \hat{\daleth}^{1}_{1}(\hat{b}_{j-1}) 
\! + \! \hat{\mathfrak{Q}}_{1}(\hat{b}_{j-1}) \right. \right. \nonumber \\
+&\left. \left. \, \hat{\mathfrak{Q}}_{1}(\hat{b}_{j-1})(\hat{\mathfrak{Q}}_{0}
(\hat{b}_{j-1}))^{-2}(\hat{\aleph}^{1}_{1}(\hat{b}_{j-1}))^{2}  \! - \! 
2(\hat{\mathfrak{Q}}_{0}(\hat{b}_{j-1}))^{-1} \hat{\aleph}^{1}_{1}(\hat{b}_{j-1}) 
\hat{\beth}^{1}_{1}(\hat{b}_{j-1}) \right) \right. \nonumber \\
+&\left. \, 2(s_{1} \! - \! t_{1}) \left(\hat{\beth}^{1}_{1}(\hat{b}_{j-1}) \! + \! 
\hat{\aleph}^{1}_{1}(\hat{b}_{j-1}) \hat{\daleth}^{1}_{1}(\hat{b}_{j-1}) \right) 
\right), \label{eqprop11} \\
\hat{\mathbb{C}}_{21}(\hat{b}_{j-1}) =& \, (\hat{\kappa}_{2}(\hat{b}_{j-1}))^{2} 
\left(\mi s_{1} \left(\hat{\mathfrak{Q}}_{0}(\hat{b}_{j-1})(\hat{\aleph}^{1}_{-1}
(\hat{b}_{j-1}))^{2} \! - \! (\hat{\mathfrak{Q}}_{0}(\hat{b}_{j-1}))^{-3}
(\hat{\mathfrak{Q}}_{1}(\hat{b}_{j-1}))^{2} \right. \right. \nonumber \\
+&\left. \left. \, \dfrac{1}{2} \hat{\mathfrak{Q}}_{2}(\hat{b}_{j-1})
(\hat{\mathfrak{Q}}_{0}(\hat{b}_{j-1}))^{-2} \! + \! 2 \hat{\mathfrak{Q}}_{1}
(\hat{b}_{j-1})(\hat{\mathfrak{Q}}_{0}(\hat{b}_{j-1}))^{-2} \hat{\daleth}^{1}_{-1}
(\hat{b}_{j-1}) \! - \! (\hat{\mathfrak{Q}}_{0}(\hat{b}_{j-1}))^{-1} \right. \right. 
\nonumber \\
\times&\left. \left. \, \left(2 \hat{\gimel}^{1}_{-1}(\hat{b}_{j-1}) \! + \! 
(\hat{\daleth}^{1}_{-1}(\hat{b}_{j-1}))^{2} \right) \right) \! + \! \mi t_{1} 
\left(2 \hat{\mathfrak{Q}}_{0}(\hat{b}_{j-1}) \hat{\daleth}^{1}_{-1}(\hat{b}_{j-1}) 
\! + \! \hat{\mathfrak{Q}}_{1}(\hat{b}_{j-1}) \right. \right. \nonumber \\
+&\left. \left. \, \hat{\mathfrak{Q}}_{1}(\hat{b}_{j-1})(\hat{\mathfrak{Q}}_{0}
(\hat{b}_{j-1}))^{-2}(\hat{\aleph}^{1}_{-1}(\hat{b}_{j-1}))^{2}  \! - \! 
2(\hat{\mathfrak{Q}}_{0}(\hat{b}_{j-1}))^{-1} \hat{\aleph}^{1}_{-1}
(\hat{b}_{j-1}) \hat{\beth}^{1}_{-1}(\hat{b}_{j-1}) \right) \right. \nonumber \\
-&\left. \, 2(s_{1} \! - \! t_{1}) \left(\hat{\beth}^{1}_{-1}(\hat{b}_{j-1}) \! + 
\! \hat{\aleph}^{1}_{-1}(\hat{b}_{j-1}) \hat{\daleth}^{1}_{-1}(\hat{b}_{j-1}) 
\right) \right), \label{eqprop12}
\end{align}
$s_{1} \! = \! 5/72$, and $t_{1} \! = \! -7/72$, where, for $\varepsilon_{1},
\varepsilon_{2} \! = \! \pm 1$,\footnote{\, $\pmb{0} \! := \! (0,0,
\dotsc,0)^{\operatorname{T}}$ $(\in \! \mathbb{R}^{N})$.}
\begin{align}
\hat{\kappa}_{1}(\varsigma) =& \, \dfrac{\hat{\boldsymbol{\theta}}
(\hat{\boldsymbol{u}}_{+}(\varsigma) \! - \! \frac{1}{2 \pi}((n \! - \! 1)K \! 
+ \! k) \hat{\boldsymbol{\Omega}} \! + \! \hat{\boldsymbol{d}})}{\hat{
\boldsymbol{\theta}}(\hat{\boldsymbol{u}}_{+}(\varsigma) \! + \! \hat{
\boldsymbol{d}})}, \qquad \quad \hat{\kappa}_{2}(\varsigma) = \dfrac{
\hat{\boldsymbol{\theta}}(\hat{\boldsymbol{u}}_{+}(\varsigma) \! - \! 
\frac{1}{2 \pi}((n \! - \! 1)K \! + \! k) \hat{\boldsymbol{\Omega}} \! - \! 
\hat{\boldsymbol{d}})}{\hat{\boldsymbol{\theta}}(\hat{\boldsymbol{u}}_{+}
(\varsigma) \! - \! \hat{\boldsymbol{d}})}, \label{eqprop13} \\
\hat{\aleph}^{\varepsilon_{1}}_{\varepsilon_{2}}(\varsigma) =& \, 
-\dfrac{\hat{\mathfrak{u}}(\varepsilon_{1},\varepsilon_{2},\bm{0};
\varsigma)}{\hat{\boldsymbol{\theta}}(\varepsilon_{1} 
\hat{\boldsymbol{u}}_{+}(\varsigma) \! + \! \varepsilon_{2} \hat{
\boldsymbol{d}})} \! + \! \dfrac{\hat{\mathfrak{u}}(\varepsilon_{1},
\varepsilon_{2},\hat{\boldsymbol{\Omega}};\varsigma)}{\hat{\boldsymbol{
\theta}}(\varepsilon_{1} \hat{\boldsymbol{u}}_{+}(\varsigma) \! - \! 
\frac{1}{2 \pi}((n \! - \! 1)K \! + \! k) \hat{\boldsymbol{\Omega}} \! + \! 
\varepsilon_{2} \hat{\boldsymbol{d}})}, \label{eqprop14} \\
\hat{\daleth}^{\varepsilon_{1}}_{\varepsilon_{2}}(\varsigma) =& \, 
-\dfrac{\hat{\mathfrak{v}}(\varepsilon_{1},\varepsilon_{2},\bm{0};
\varsigma)}{\hat{\boldsymbol{\theta}}(\varepsilon_{1} 
\hat{\boldsymbol{u}}_{+}(\varsigma) \! + \! \varepsilon_{2} \hat{
\boldsymbol{d}})} \! + \! \dfrac{\hat{\mathfrak{v}}(\varepsilon_{1},
\varepsilon_{2},\hat{\boldsymbol{\Omega}};\varsigma)}{\hat{\boldsymbol{
\theta}}(\varepsilon_{1} \hat{\boldsymbol{u}}_{+}(\varsigma) \! - \! 
\frac{1}{2 \pi}((n \! - \! 1)K \! + \! k) \hat{\boldsymbol{\Omega}} \! + \! 
\varepsilon_{2} \hat{\boldsymbol{d}})} \! - \! \left(\dfrac{\hat{\mathfrak{u}}
(\varepsilon_{1},\varepsilon_{2},\bm{0};\varsigma)}{\hat{\boldsymbol{\theta}}
(\varepsilon_{1} \hat{\boldsymbol{u}}_{+}(\varsigma) \! + \! \varepsilon_{2} 
\hat{\boldsymbol{d}})} \right)^{2} \nonumber \\
+& \, \dfrac{\hat{\mathfrak{u}}(\varepsilon_{1},\varepsilon_{2},\bm{0};\varsigma) 
\hat{\mathfrak{u}}(\varepsilon_{1},\varepsilon_{2},\hat{\boldsymbol{\Omega}};
\varsigma)}{\hat{\boldsymbol{\theta}}(\varepsilon_{1} \hat{\boldsymbol{u}}_{+}
(\varsigma) \! + \! \varepsilon_{2} \hat{\boldsymbol{d}}) \hat{\boldsymbol{
\theta}}(\varepsilon_{1} \hat{\boldsymbol{u}}_{+}(\varsigma) \! - \! \frac{1}{2 
\pi}((n \! - \! 1)K \! + \! k) \hat{\boldsymbol{\Omega}} \! + \! \varepsilon_{2} 
\hat{\boldsymbol{d}})}, \label{eqprop15} \\
\hat{\beth}^{\varepsilon_{1}}_{\varepsilon_{2}}(\xi) =& \, -\dfrac{\hat{
\mathfrak{w}}(\varepsilon_{1},\varepsilon_{2},\bm{0};\varsigma)}{\hat{
\boldsymbol{\theta}}(\varepsilon_{1} \hat{\boldsymbol{u}}_{+}(\varsigma) \! 
+ \! \varepsilon_{2} \hat{\boldsymbol{d}})} \! + \! \dfrac{\hat{\mathfrak{w}}
(\varepsilon_{1},\varepsilon_{2},\hat{\boldsymbol{\Omega}};\varsigma)}{
\hat{\boldsymbol{\theta}}(\varepsilon_{1} \hat{\boldsymbol{u}}_{+}(\varsigma) 
\! - \! \frac{1}{2 \pi}((n \! - \! 1)K \! + \! k) \hat{\boldsymbol{\Omega}} \! + 
\! \varepsilon_{2} \hat{\boldsymbol{d}})} \! + \! \dfrac{2 \hat{\mathfrak{u}}
(\varepsilon_{1},\varepsilon_{2},\bm{0};\varsigma) \hat{\mathfrak{v}}
(\varepsilon_{1},\varepsilon_{2},\bm{0};\varsigma)}{(\hat{\boldsymbol{
\theta}}(\varepsilon_{1} \hat{\boldsymbol{u}}_{+}(\varsigma) \! + \! 
\varepsilon_{2} \hat{\boldsymbol{d}}))^{2}} \nonumber \\
-& \, \dfrac{\hat{\mathfrak{v}}(\varepsilon_{1},\varepsilon_{2},\bm{0};\varsigma) 
\hat{\mathfrak{u}}(\varepsilon_{1},\varepsilon_{2},\hat{\boldsymbol{\Omega}};
\varsigma)}{\hat{\boldsymbol{\theta}}(\varepsilon_{1} \hat{\boldsymbol{u}}_{+}
(\varsigma) \! + \! \varepsilon_{2} \hat{\boldsymbol{d}}) \hat{\boldsymbol{
\theta}}(\varepsilon_{1} \hat{\boldsymbol{u}}_{+}(\varsigma) \! - \! \frac{1}{2 
\pi}((n \! - \! 1)K \! + \! k) \hat{\boldsymbol{\Omega}} \! + \! \varepsilon_{2} 
\hat{\boldsymbol{d}})} \! + \! \left(\dfrac{\hat{\mathfrak{u}}(\varepsilon_{1},
\varepsilon_{2},\bm{0};\varsigma)}{\hat{\boldsymbol{\theta}}(\varepsilon_{1} 
\hat{\boldsymbol{u}}_{+}(\varsigma) \! + \! \varepsilon_{2} \hat{\boldsymbol{
d}})} \right)^{3} \nonumber \\
-& \, \dfrac{\hat{\mathfrak{u}}(\varepsilon_{1},\varepsilon_{2},\bm{0};\varsigma) 
\hat{\mathfrak{v}}(\varepsilon_{1},\varepsilon_{2},\hat{\boldsymbol{\Omega}};
\varsigma)}{\hat{\boldsymbol{\theta}}(\varepsilon_{1} \hat{\boldsymbol{u}}_{+}
(\varsigma) \! + \! \varepsilon_{2} \hat{\boldsymbol{d}}) \hat{\boldsymbol{
\theta}}(\varepsilon_{1} \hat{\boldsymbol{u}}_{+}(\varsigma) \! - \! \frac{1}{
2 \pi}((n \! - \! 1)K \! + \! k) \hat{\boldsymbol{\Omega}} \! + \! \varepsilon_{2} 
\hat{\boldsymbol{d}})} \! - \! \left(\dfrac{\hat{\mathfrak{u}}(\varepsilon_{1},
\varepsilon_{2},\bm{0};\varsigma)}{\hat{\boldsymbol{\theta}}(\varepsilon_{1} 
\hat{\boldsymbol{u}}_{+}(\varsigma) \! + \! \varepsilon_{2} \hat{\boldsymbol{
d}})} \right)^{2} \nonumber \\
\times& \, \dfrac{\hat{\mathfrak{u}}(\varepsilon_{1},\varepsilon_{2},\hat{
\boldsymbol{\Omega}};\varsigma)}{\hat{\boldsymbol{\theta}}(\varepsilon_{1} 
\hat{\boldsymbol{u}}_{+}(\varsigma) \! - \! \frac{1}{2 \pi}((n \! - \! 1)K \! + 
\! k) \hat{\boldsymbol{\Omega}} \! + \! \varepsilon_{2} \hat{\boldsymbol{d}})}, 
\label{eqprop16} \\
\hat{\gimel}^{\varepsilon_{1}}_{\varepsilon_{2}}(\varsigma) =& \, -\dfrac{
\hat{\mathfrak{z}}(\varepsilon_{1},\varepsilon_{2},\bm{0};\varsigma)}{\hat{
\boldsymbol{\theta}}(\varepsilon_{1} \hat{\boldsymbol{u}}_{+}(\varsigma) \! 
+ \! \varepsilon_{2} \hat{\boldsymbol{d}})} \! + \! \dfrac{\hat{\mathfrak{z}}
(\varepsilon_{1},\varepsilon_{2},\hat{\boldsymbol{\Omega}};\varsigma)}{
\hat{\boldsymbol{\theta}}(\varepsilon_{1} \hat{\boldsymbol{u}}_{+}(\varsigma) 
\! - \! \frac{1}{2 \pi}((n \! - \! 1)K \! + \! k) \hat{\boldsymbol{\Omega}} \! + \! 
\varepsilon_{2} \hat{\boldsymbol{d}})} \! + \! \left(\dfrac{\hat{\mathfrak{v}}
(\varepsilon_{1},\varepsilon_{2},\bm{0};\varsigma)}{\hat{\boldsymbol{\theta}}
(\varepsilon_{1} \hat{\boldsymbol{u}}_{+}(\varsigma) \! + \! \varepsilon_{2} 
\hat{\boldsymbol{d}})} \right)^{2} \nonumber \\
-& \, \dfrac{\hat{\mathfrak{v}}(\varepsilon_{1},\varepsilon_{2},\bm{0};\varsigma) 
\hat{\mathfrak{v}}(\varepsilon_{1},\varepsilon_{2},\hat{\boldsymbol{\Omega}};
\varsigma)}{\hat{\boldsymbol{\theta}}(\varepsilon_{1} \hat{\boldsymbol{u}}_{+}
(\varsigma) \! + \! \varepsilon_{2} \hat{\boldsymbol{d}}) \hat{\boldsymbol{
\theta}}(\varepsilon_{1} \hat{\boldsymbol{u}}_{+}(\varsigma) \! - \! \frac{1}{2 
\pi}((n \! - \! 1)K \! + \! k) \hat{\boldsymbol{\Omega}} \! + \! \varepsilon_{2} 
\hat{\boldsymbol{d}})} \! - \! \dfrac{2 \hat{\mathfrak{u}}(\varepsilon_{1},
\varepsilon_{2},\bm{0};\varsigma) \hat{\mathfrak{w}}(\varepsilon_{1},
\varepsilon_{2},\bm{0};\varsigma)}{(\hat{\boldsymbol{\theta}}(\varepsilon_{1} 
\hat{\boldsymbol{u}}_{+}(\varsigma) \! + \! \varepsilon_{2} \hat{\boldsymbol{
d}}))^{2}} \nonumber \\
+& \, \dfrac{\hat{\mathfrak{w}}(\varepsilon_{1},\varepsilon_{2},\bm{0};\varsigma) 
\hat{\mathfrak{u}}(\varepsilon_{1},\varepsilon_{2},\hat{\boldsymbol{\Omega}};
\varsigma)}{\hat{\boldsymbol{\theta}}(\varepsilon_{1} \hat{\boldsymbol{u}}_{+}
(\varsigma) \! + \! \varepsilon_{2} \hat{\boldsymbol{d}}) \hat{\boldsymbol{
\theta}}(\varepsilon_{1} \hat{\boldsymbol{u}}_{+}(\varsigma) \! - \! \frac{1}{2 
\pi}((n \! - \! 1)K \! + \! k) \hat{\boldsymbol{\Omega}} \! + \! \varepsilon_{2} 
\hat{\boldsymbol{d}})} \! + \! \dfrac{3(\hat{\mathfrak{u}}(\varepsilon_{1},
\varepsilon_{2},\bm{0};\varsigma))^{2} \hat{\mathfrak{v}}(\varepsilon_{1},
\varepsilon_{2},\bm{0};\varsigma)}{(\hat{\boldsymbol{\theta}}(\varepsilon_{1} 
\hat{\boldsymbol{u}}_{+}(\varsigma) \! + \! \varepsilon_{2} \hat{\boldsymbol{
d}}))^{3}} \nonumber \\
+& \, \dfrac{\hat{\mathfrak{u}}(\varepsilon_{1},\varepsilon_{2},\bm{0};\varsigma) 
\hat{\mathfrak{w}}(\varepsilon_{1},\varepsilon_{2},\hat{\boldsymbol{\Omega}};
\varsigma)}{\hat{\boldsymbol{\theta}}(\varepsilon_{1} \hat{\boldsymbol{u}}_{+}
(\varsigma) \! + \! \varepsilon_{2} \hat{\boldsymbol{d}}) \hat{\boldsymbol{\theta}}
(\varepsilon_{1} \hat{\boldsymbol{u}}_{+}(\varsigma) \! - \! \frac{1}{2 \pi}
((n \! - \! 1)K \! + \! k) \hat{\boldsymbol{\Omega}} \! + \! \varepsilon_{2} 
\hat{\boldsymbol{d}})} \! + \! \left(\dfrac{\hat{\mathfrak{u}}(\varepsilon_{1},
\varepsilon_{2},\bm{0};\varsigma)}{\hat{\boldsymbol{\theta}}(\varepsilon_{1} 
\hat{\boldsymbol{u}}_{+}(\varsigma) \! + \! \varepsilon_{2} \hat{\boldsymbol{d}})} 
\right)^{4} \nonumber \\
-& \, \dfrac{2 \hat{\mathfrak{u}}(\varepsilon_{1},\varepsilon_{2},\bm{0};\varsigma) 
\hat{\mathfrak{v}}(\varepsilon_{1},\varepsilon_{2},\bm{0};\varsigma) \hat{
\mathfrak{u}}(\varepsilon_{1},\varepsilon_{2},\hat{\boldsymbol{\Omega}};
\varsigma)}{(\hat{\boldsymbol{\theta}}(\varepsilon_{1} \hat{\boldsymbol{u}}_{+}
(\varsigma) \! + \! \varepsilon_{2} \hat{\boldsymbol{d}}))^{2} \hat{\boldsymbol{
\theta}}(\varepsilon_{1} \hat{\boldsymbol{u}}_{+}(\varsigma) \! - \! \frac{1}{2 \pi}
((n \! - \! 1)K \! + \! k) \hat{\boldsymbol{\Omega}} \! + \! \varepsilon_{2} \hat{
\boldsymbol{d}})} \! - \! \left(\dfrac{\hat{\mathfrak{u}}(\varepsilon_{1},
\varepsilon_{2},\bm{0};\varsigma)}{\hat{\boldsymbol{\theta}}(\varepsilon_{1} 
\hat{\boldsymbol{u}}_{+}(\varsigma) \! + \! \varepsilon_{2} \hat{\boldsymbol{d}})} 
\right)^{2} \nonumber \\
\times& \, \dfrac{\hat{\mathfrak{v}}(\varepsilon_{1},\varepsilon_{2},\hat{
\boldsymbol{\Omega}};\varsigma)}{\hat{\boldsymbol{\theta}}(\varepsilon_{1} 
\hat{\boldsymbol{u}}_{+}(\varsigma) \! - \! \frac{1}{2 \pi}((n \! - \! 1)K \! + \! 
k) \hat{\boldsymbol{\Omega}} \! + \! \varepsilon_{2} \hat{\boldsymbol{d}})} 
\! - \! \dfrac{(\hat{\mathfrak{u}}(\varepsilon_{1},\varepsilon_{2},\bm{0};
\varsigma))^{3} \hat{\mathfrak{u}}(\varepsilon_{1},\varepsilon_{2},\hat{
\boldsymbol{\Omega}};\varsigma)}{(\hat{\boldsymbol{\theta}}(\varepsilon_{1} 
\hat{\boldsymbol{u}}_{+}(\varsigma) \! + \! \varepsilon_{2} \hat{\boldsymbol{
d}}))^{3} \hat{\boldsymbol{\theta}}(\varepsilon_{1} \hat{\boldsymbol{u}}_{+}
(\varsigma) \! - \! \frac{1}{2 \pi}((n \! - \! 1)K \! + \! k) \hat{\boldsymbol{
\Omega}} \! + \! \varepsilon_{2} \hat{\boldsymbol{d}})}, \label{eqprop17}
\end{align}
with
\begin{gather}
\hat{\mathfrak{u}}(\varepsilon_{1},\varepsilon_{2},\hat{\boldsymbol{\Omega}};
\varsigma) \! := \! 2 \pi \hat{\Lambda}^{\raise-1.0ex\hbox{$\scriptstyle 1$}}_{0}
(\varepsilon_{1},\varepsilon_{2},\hat{\boldsymbol{\Omega}};\varsigma), \qquad 
\quad \hat{\mathfrak{v}}(\varepsilon_{1},\varepsilon_{2},\hat{\boldsymbol{\Omega}};
\varsigma) \! := \! -2 \pi^{2} \hat{\Lambda}^{\raise-1.0ex\hbox{$\scriptstyle 2$}}_{0}
(\varepsilon_{1},\varepsilon_{2},\hat{\boldsymbol{\Omega}};\varsigma), 
\label{eqprop18} \\
\hat{\mathfrak{w}}(\varepsilon_{1},\varepsilon_{2},\hat{\boldsymbol{\Omega}};
\varsigma) \! := \! 2 \pi \left(\hat{\Lambda}^{\raise-1.0ex\hbox{$\scriptstyle 0$}}_{1}
(\varepsilon_{1},\varepsilon_{2},\hat{\boldsymbol{\Omega}};\varsigma) \! - \! 
\dfrac{2 \pi^{2}}{3} \hat{\Lambda}^{\raise-1.0ex\hbox{$\scriptstyle 3$}}_{0}
(\varepsilon_{1},\varepsilon_{2},\hat{\boldsymbol{\Omega}};\varsigma) \right), 
\label{eqprop19} \\
\hat{\mathfrak{z}}(\varepsilon_{1},\varepsilon_{2},\hat{\boldsymbol{\Omega}};
\varsigma) \! := \! -4 \pi^{2} \left(
\hat{\Lambda}^{\raise-1.0ex\hbox{$\scriptstyle 1$}}_{1}(\varepsilon_{1},
\varepsilon_{2},\hat{\boldsymbol{\Omega}};\varsigma) \! - \! \dfrac{\pi^{2}}{6} 
\hat{\Lambda}^{\raise-1.0ex\hbox{$\scriptstyle 4$}}_{0}(\varepsilon_{1},
\varepsilon_{2},\hat{\boldsymbol{\Omega}};\varsigma) \right), \label{eqprop20} \\
\hat{\Lambda}^{\raise-1.0ex\hbox{$\scriptstyle j_{1}$}}_{j_{2}}(\varepsilon_{1},
\varepsilon_{2},\hat{\boldsymbol{\Omega}};\varsigma) \! = \! \sum_{m \in
\mathbb{Z}^{N}}(\hat{\mathfrak{r}}_{1}(\varsigma))^{j_{1}}(\hat{\mathfrak{r}}_{2}
(\varsigma))^{j_{2}} \me^{2 \pi \mi (m,\varepsilon_{1} \hat{\boldsymbol{u}}_{+}
(\varsigma)-\frac{1}{2 \pi}((n-1)K+k) \hat{\boldsymbol{\Omega}}+ 
\varepsilon_{2} \hat{\boldsymbol{d}})+ \mi \pi (m,\hat{\boldsymbol{\tau}}m)}, 
\quad j_{1},j_{2} \! \in \! \mathbb{N}_{0}, \label{eqprop21} \\
\hat{\mathfrak{r}}_{1}(\varsigma) \! := \! \dfrac{2}{\hat{\leftthreetimes}
(\varsigma)} \sum_{i=1}^{N} \sum_{j=1}^{N}m_{i} \hat{c}_{ij} \varsigma^{N-j}, 
\label{eqprop22} \\
\hat{\mathfrak{r}}_{2}(\varsigma) \! := \! \dfrac{2}{3 \hat{\leftthreetimes}
(\varsigma)} \sum_{i=1}^{N} \sum_{j=1}^{N}m_{i} \hat{c}_{ij} \left(N \! - 
\! j \! - \! \dfrac{\varsigma \hat{\leftthreetimes}^{\prime}(\varsigma)}{
\hat{\leftthreetimes}(\varsigma)} \right) \varsigma^{N-j-1}, \label{eqprop23}
\end{gather}
where {}\footnote{All square roots are positive.}
\begin{align}
\hat{\leftthreetimes}(\hat{b}_{0}) :=& \, \mi (-1)^{N}(\hat{a}_{N+1} \! - \! 
\hat{b}_{0})^{1/2} \prod_{m=1}^{N}(\hat{b}_{m} \! - \! \hat{b}_{0})^{1/2}
(\hat{a}_{m} \! - \! \hat{b}_{0})^{1/2}, \label{eqprop24} \\
\hat{\leftthreetimes}(\hat{b}_{j}) :=& \, \mi (-1)^{N-j}(\hat{b}_{j} \! - \! 
\hat{b}_{0})^{1/2}(\hat{a}_{N+1} \! - \! \hat{b}_{j})^{1/2}(\hat{b}_{j} \! - \! 
\hat{a}_{j})^{1/2} \prod_{m=1}^{j-1}(\hat{b}_{j} \! - \! \hat{b}_{m})^{1/2}
(\hat{b}_{j} \! - \! \hat{a}_{m})^{1/2} \nonumber \\
\times& \, \prod_{m^{\prime}=j+1}^{N}(\hat{b}_{m^{\prime}} \! - \! 
\hat{b}_{j})^{1/2}(\hat{a}_{m^{\prime}} \! - \! \hat{b}_{j})^{1/2}, 
\quad j \! = \! 1,2,\dotsc,N, \label{eqprop25} \\
\dfrac{\hat{\leftthreetimes}^{\prime}(\hat{b}_{0})}{\hat{\leftthreetimes}
(\hat{b}_{0})} =& \, \dfrac{1}{2} \left(\sum_{m=1}^{N} \left(\dfrac{1}{
\hat{b}_{0} \! - \! \hat{b}_{m}} \! + \! \dfrac{1}{\hat{b}_{0} \! - \! \hat{a}_{m}} 
\right) \! + \! \dfrac{1}{\hat{b}_{0} \! - \! \hat{a}_{N+1}} \right), 
\label{eqprop26} \\
\dfrac{\hat{\leftthreetimes}^{\prime}(\hat{b}_{j})}{\hat{\leftthreetimes}
(\hat{b}_{j})} =& \, \dfrac{1}{2} \left(\sum_{\substack{m=1\\m \neq j}}^{N} 
\left(\dfrac{1}{\hat{b}_{j} \! - \! \hat{b}_{m}} \! + \! \dfrac{1}{\hat{b}_{j} \! - 
\! \hat{a}_{m}} \right) \! + \! \dfrac{1}{\hat{b}_{j} \! - \! \hat{a}_{j}} \! + \! 
\dfrac{1}{\hat{b}_{j} \! - \! \hat{b}_{0}} \! + \! \dfrac{1}{\hat{b}_{j} \! - \! 
\hat{a}_{N+1}} \right), \quad j \! = \! 1,2,\dotsc,N, \label{eqprop27}
\end{align}
and $\hat{c}_{ij}$, $i,j \! = \! 1,2,\dotsc,N$, are described in 
Equations~\eqref{E1} and~\eqref{E2}, and
\begin{align}
\hat{\mathfrak{Q}}_{0}(\hat{b}_{0})=& \, -\mi (\hat{a}_{N+1} \! - \! 
\hat{b}_{0})^{-1/2} \prod_{m=1}^{N} \dfrac{(\hat{b}_{m} \! - \! 
\hat{b}_{0})^{1/2}}{(\hat{a}_{m} \! - \! \hat{b}_{0})^{1/2}}, \label{eqprop28} \\
\dfrac{\hat{\mathfrak{Q}}_{1}(\hat{b}_{0})}{\hat{\mathfrak{Q}}_{0}(\hat{b}_{0})}
=& \, \dfrac{1}{2} \left(\sum_{m=1}^{N} \left(\dfrac{1}{\hat{b}_{0} \! - \! 
\hat{b}_{m}} \! - \! \dfrac{1}{\hat{b}_{0} \! - \! \hat{a}_{m}} \right) \! - \! 
\dfrac{1}{\hat{b}_{0} \! - \! \hat{a}_{N+1}} \right), \label{eqprop29} \\
\dfrac{\hat{\mathfrak{Q}}_{2}(\hat{b}_{0})}{\hat{\mathfrak{Q}}_{0}(\hat{b}_{0})}
=& \, \dfrac{1}{2} \left(\sum_{m=1}^{N} \left(\dfrac{1}{(\hat{b}_{0} \! - \! 
\hat{a}_{m})^{2}} \! - \! \dfrac{1}{(\hat{b}_{0} \! - \! \hat{b}_{m})^{2}} \right) 
\! + \! \dfrac{1}{(\hat{b}_{0} \! - \! \hat{a}_{N+1})^{2}} \right) \! + \! 
\dfrac{1}{4} \left(\sum_{m=1}^{N} \left(\dfrac{1}{\hat{b}_{0} \! - \! 
\hat{a}_{m}} \! - \! \dfrac{1}{\hat{b}_{0} \! - \! \hat{b}_{m}} \right) \! + \! 
\dfrac{1}{\hat{b}_{0} \! - \! \hat{a}_{N+1}} \right)^{2}, \label{eqprop30} \\
\hat{\mathfrak{Q}}_{0}(\hat{b}_{j})=& \, -\dfrac{\mi (\hat{b}_{j} \! - \! 
\hat{b}_{0})^{1/2}}{(\hat{a}_{N+1} \! - \! \hat{b}_{j})^{1/2}(\hat{b}_{j} \! 
- \! \hat{a}_{j})^{1/2}} \prod_{m=1}^{j-1} \dfrac{(\hat{b}_{j} \! - \! 
\hat{b}_{m})^{1/2}}{(\hat{b}_{j} \! - \! \hat{a}_{m})^{1/2}} \prod_{m^{\prime}
=j+1}^{N} \dfrac{(\hat{b}_{m^{\prime}} \! - \! \hat{b}_{j})^{1/2}}{(\hat{a}_{
m^{\prime}} \! - \! \hat{b}_{j})^{1/2}}, \quad j \! = \! 1,2,\dotsc,N, 
\label{eqprop31} \\
\dfrac{\hat{\mathfrak{Q}}_{1}(\hat{b}_{j})}{\hat{\mathfrak{Q}}_{0}(\hat{b}_{j})}
=& \, \dfrac{1}{2} \left(\sum_{\substack{m=1\\m \neq j}}^{N} \left(\dfrac{1}{
\hat{b}_{j} \! - \! \hat{b}_{m}} \! - \! \dfrac{1}{\hat{b}_{j} \! - \! \hat{a}_{m}} 
\right) \! + \! \dfrac{1}{\hat{b}_{j} \! - \! \hat{b}_{0}} \! - \! \dfrac{1}{
\hat{b}_{j} \! - \! \hat{a}_{N+1}} \! - \! \dfrac{1}{\hat{b}_{j} \! - \! \hat{a}_{j}} 
\right), \quad j \! = \! 1,2,\dotsc,N, \label{eqprop32} \\
\dfrac{\hat{\mathfrak{Q}}_{2}(\hat{b}_{j})}{\hat{\mathfrak{Q}}_{0}(\hat{b}_{j})}
=& \, -\dfrac{1}{2} \left(\sum_{\substack{m=1\\m \neq j}}^{N} \left(
\dfrac{1}{(\hat{b}_{j} \! - \! \hat{b}_{m})^{2}} \! - \! \dfrac{1}{(\hat{b}_{j} \! - \! 
\hat{a}_{m})^{2}} \right) \! + \! \dfrac{1}{(\hat{b}_{j} \! - \! \hat{b}_{0})^{2}} \! 
- \! \dfrac{1}{(\hat{b}_{j} \! - \! \hat{a}_{N+1})^{2}} \! - \! \dfrac{1}{(\hat{b}_{j} 
\! - \! \hat{a}_{j})^{2}} \right) \nonumber \\
+& \, \dfrac{1}{4} \left(\sum_{\substack{m=1\\m \neq j}}^{N} \left(\dfrac{1}{
\hat{b}_{j} \! - \! \hat{b}_{m}} \! - \! \dfrac{1}{\hat{b}_{j} \! - \! \hat{a}_{m}} 
\right) \! + \! \dfrac{1}{\hat{b}_{j} \! - \! \hat{b}_{0}} \! - \! \dfrac{1}{
\hat{b}_{j} \! - \! \hat{a}_{N+1}} \! - \! \dfrac{1}{\hat{b}_{j} \! - \! \hat{a}_{j}} 
\right)^{2}, \quad j \! = \! 1,2,\dotsc,N, \label{eqprop33}
\end{align}
$\hat{\alpha}_{0}(\hat{b}_{0})$, $\hat{\alpha}_{0}(\hat{b}_{j})$, $\hat{\alpha}_{1}
(\hat{b}_{0})$, and $\hat{\alpha}_{1}(\hat{b}_{j})$, $j \! = \! 1,2,\dotsc,N$, are 
given in Equations~\eqref{eqmaininf47}--\eqref{eqmaininf50}, respectively, 
$\hat{\alpha}_{2}(\hat{b}_{0}) \! = \! f^{\prime \prime}(\hat{b}_{0})/7$ and 
$\hat{\alpha}_{2}(\hat{b}_{j}) \! = \! f^{\prime \prime}(\hat{b}_{j})/7$, $j \! = \! 
1,2,\dotsc,N$, where $f^{\prime \prime}(\hat{b}_{0})$ and $f^{\prime \prime}
(\hat{b}_{j})$ are given in Equations~\eqref{wm24frbohat} and~\eqref{wm24frbjhat}, 
respectively, and $(\mathrm{M}_{2}(\mathbb{C}) \! \ni)$ $\hat{\mathfrak{c}}^{r}_{m}
(n,k,z_{o};\hat{b}_{j-1}) \! =_{\underset{z_{o}=1+o(1)}{\mathscr{N},n \to \infty}} \! 
\mathcal{O}(1)$, $m \! \in \! \mathbb{N}_{0}$, $r \! \in \! \lbrace \triangleright,
\triangleleft \rbrace$, and, for $j \! = \! 1,2,\dotsc,N \! + \! 1$,
\begin{gather}
\hat{\boldsymbol{\mathrm{A}}}(\hat{a}_{j}) \! := \! 
\begin{pmatrix}
\tilde{\mathfrak{m}}^{\raise-0.5ex\hbox{$\scriptstyle \infty$}}_{11} & 0 \\
0 & \tilde{\mathfrak{m}}^{\raise-0.5ex\hbox{$\scriptstyle \infty$}}_{22}
\end{pmatrix} 
\begin{pmatrix}
\hat{\mathbb{A}}_{11}(\hat{a}_{j}) & \hat{\mathbb{A}}_{12}(\hat{a}_{j}) \\
\hat{\mathbb{A}}_{21}(\hat{a}_{j}) & \hat{\mathbb{A}}_{22}(\hat{a}_{j})
\end{pmatrix} 
\begin{pmatrix}
\tilde{\mathfrak{m}}^{\raise-0.5ex\hbox{$\scriptstyle \infty$}}_{22} & 0 \\
0 & \tilde{\mathfrak{m}}^{\raise-0.5ex\hbox{$\scriptstyle \infty$}}_{11}
\end{pmatrix} \me^{\mi ((n-1)K+k) \hat{\mho}_{j}}, \label{eqprop34} \\
\hat{\boldsymbol{\mathrm{B}}}(\hat{a}_{j}) \! := \! 
\begin{pmatrix}
\tilde{\mathfrak{m}}^{\raise-0.5ex\hbox{$\scriptstyle \infty$}}_{11} & 0 \\
0 & \tilde{\mathfrak{m}}^{\raise-0.5ex\hbox{$\scriptstyle \infty$}}_{22}
\end{pmatrix} 
\begin{pmatrix}
\hat{\mathbb{B}}_{11}(\hat{a}_{j}) & \hat{\mathbb{B}}_{12}(\hat{a}_{j}) \\
\hat{\mathbb{B}}_{21}(\hat{a}_{j}) & \hat{\mathbb{B}}_{22}(\hat{a}_{j})
\end{pmatrix} 
\begin{pmatrix}
\tilde{\mathfrak{m}}^{\raise-0.5ex\hbox{$\scriptstyle \infty$}}_{22} & 0 \\
0 & \tilde{\mathfrak{m}}^{\raise-0.5ex\hbox{$\scriptstyle \infty$}}_{11}
\end{pmatrix} \me^{\mi ((n-1)K+k) \hat{\mho}_{j}}, \label{eqprop35} \\
\hat{\boldsymbol{\mathrm{C}}}(\hat{a}_{j}) \! := \! 
\begin{pmatrix}
\tilde{\mathfrak{m}}^{\raise-0.5ex\hbox{$\scriptstyle \infty$}}_{11} & 0 \\
0 & \tilde{\mathfrak{m}}^{\raise-0.5ex\hbox{$\scriptstyle \infty$}}_{22}
\end{pmatrix} 
\begin{pmatrix}
\hat{\mathbb{C}}_{11}(\hat{a}_{j}) & \hat{\mathbb{C}}_{12}(\hat{a}_{j}) \\
\hat{\mathbb{C}}_{21}(\hat{a}_{j}) & \hat{\mathbb{C}}_{22}(\hat{a}_{j})
\end{pmatrix} 
\begin{pmatrix}
\tilde{\mathfrak{m}}^{\raise-0.5ex\hbox{$\scriptstyle \infty$}}_{22} & 0 \\
0 & \tilde{\mathfrak{m}}^{\raise-0.5ex\hbox{$\scriptstyle \infty$}}_{11}
\end{pmatrix} \me^{\mi ((n-1)K+k) \hat{\mho}_{j}}, \label{eqprop36}
\end{gather}
with
\begin{gather}
\hat{\mathbb{A}}_{11}(\hat{a}_{j}) \! = \! -\hat{\mathbb{A}}_{22}(\hat{a}_{j}) 
\! = \! -s_{1} \hat{\kappa}_{1}(\hat{a}_{j}) \hat{\kappa}_{2}(\hat{a}_{j}) 
\hat{\mathfrak{Q}}_{0}(\hat{a}_{j}), \label{eqprop37} \\
\hat{\mathbb{A}}_{12}(\hat{a}_{j}) \! = \! \mi s_{1}(\hat{\kappa}_{1}
(\hat{a}_{j}))^{2} \hat{\mathfrak{Q}}_{0}(\hat{a}_{j}), \qquad 
\hat{\mathbb{A}}_{21}(\hat{a}_{j}) \! = \! \mi s_{1}(\hat{\kappa}_{2}
(\hat{a}_{j}))^{2} \hat{\mathfrak{Q}}_{0}(\hat{a}_{j}), \label{eqprop38}
\end{gather}
\begin{align}
\hat{\mathbb{B}}_{11}(\hat{a}_{j}) \! = \! -\hat{\mathbb{B}}_{22}
(\hat{a}_{j}) =& \, \hat{\kappa}_{1}(\hat{a}_{j}) \hat{\kappa}_{2}
(\hat{a}_{j}) \left(-s_{1} \left(\hat{\mathfrak{Q}}_{1}(\hat{a}_{j}) \! + \! 
\hat{\mathfrak{Q}}_{0}(\hat{a}_{j}) \left(\hat{\daleth}^{1}_{1}(\hat{a}_{j}) \! 
+ \! \hat{\daleth}^{1}_{-1}(\hat{a}_{j}) \right) \right) \right. \nonumber \\
-&\left. \, t_{1} \left((\hat{\mathfrak{Q}}_{0}(\hat{a}_{j}))^{-1} \! + \! 
\hat{\mathfrak{Q}}_{0}(\hat{a}_{j}) \hat{\aleph}^{1}_{1}(\hat{a}_{j}) 
\hat{\aleph}^{1}_{-1}(\hat{a}_{j}) \right) \! + \! \mi (s_{1} \! + \! t_{1}) 
\left(\hat{\aleph}^{1}_{-1}(\hat{a}_{j}) \! - \! \hat{\aleph}^{1}_{1}(\hat{a}_{j}) 
\right) \right), \label{eqprop39} \\
\hat{\mathbb{B}}_{12}(\hat{a}_{j}) =& \, (\hat{\kappa}_{1}(\hat{a}_{j}))^{2} 
\left(\mi s_{1} \left(\hat{\mathfrak{Q}}_{1}(\hat{a}_{j}) \! + \! 2 \hat{
\mathfrak{Q}}_{0}(\hat{a}_{j}) \hat{\daleth}^{1}_{1}(\hat{a}_{j}) \right) 
\right. \nonumber \\
+&\left. \, \mi t_{1} \left(\hat{\mathfrak{Q}}_{0}(\hat{a}_{j})(\hat{\aleph}^{
1}_{1}(\hat{a}_{j}))^{2} \! - \! (\hat{\mathfrak{Q}}_{0}(\hat{a}_{j}))^{-1} 
\right) \! - \! 2(s_{1} \! - \! t_{1}) \hat{\aleph}^{1}_{1}(\hat{a}_{j}) \right), 
\label{eqprop40} \\
\hat{\mathbb{B}}_{21}(\hat{a}_{j}) =& \, (\hat{\kappa}_{2}(\hat{a}_{j}))^{2} 
\left(\mi s_{1} \left(\hat{\mathfrak{Q}}_{1}(\hat{a}_{j}) \! + \! 2 \hat{
\mathfrak{Q}}_{0}(\hat{a}_{j}) \hat{\daleth}^{1}_{-1}(\hat{a}_{j}) \right) 
\right. \nonumber \\
+&\left. \, \mi t_{1} \left(\hat{\mathfrak{Q}}_{0}(\hat{a}_{j})(\hat{\aleph}^{
1}_{-1}(\hat{a}_{j}))^{2} \! - \! (\hat{\mathfrak{Q}}_{0}(\hat{a}_{j}))^{-1} 
\right) \! + \! 2(s_{1} \! - \! t_{1}) \hat{\aleph}^{1}_{-1}(\hat{a}_{j}) 
\right), \label{eqprop41} \\
\hat{\mathbb{C}}_{11}(\hat{a}_{j}) \! = \! -\hat{\mathbb{C}}_{22}(\hat{a}_{j}) 
=& \, \hat{\kappa}_{1}(\hat{a}_{j}) \hat{\kappa}_{2}(\hat{a}_{j}) \left(
-s_{1} \left(\hat{\mathfrak{Q}}_{0}(\hat{a}_{j}) \left(\hat{\gimel}^{1}_{1}
(\hat{a}_{j}) \! + \! \hat{\gimel}^{1}_{-1}(\hat{a}_{j}) \! + \! \hat{\daleth}^{
1}_{1}(\hat{a}_{j}) \hat{\daleth}^{1}_{-1}(\hat{a}_{j}) \right) \right. \right. 
\nonumber \\
+&\left. \left. \, \hat{\mathfrak{Q}}_{1}(\hat{a}_{j}) \left(\hat{\daleth}^{
1}_{1}(\hat{a}_{j}) \! + \! \hat{\daleth}^{1}_{-1}(\hat{a}_{j}) \right) \! + \! 
\dfrac{1}{2} \hat{\mathfrak{Q}}_{2}(\hat{a}_{j}) \! + \! (\hat{\mathfrak{Q}}_{0}
(\hat{a}_{j}))^{-1}\hat{\aleph}^{1}_{1}(\hat{a}_{j}) \hat{\aleph}^{1}_{-1}
(\hat{a}_{j}) \right) \right. \nonumber \\
-&\left. t_{1} \left(\hat{\mathfrak{Q}}_{0}(\hat{a}_{j}) \left(\hat{\aleph}^{
1}_{-1}(\hat{a}_{j}) \hat{\beth}^{1}_{1}(\hat{a}_{j}) \! + \! \hat{\aleph}^{1}_{1}
(\hat{a}_{j}) \hat{\beth}^{1}_{-1}(\hat{a}_{j}) \right) \! + \! (\hat{\mathfrak{
Q}}_{0}(\hat{a}_{j}))^{-1} \right. \right. \nonumber \\
\times&\left. \left. \, \left(\hat{\daleth}^{1}_{1}(\hat{a}_{j}) \! + \! 
\hat{\daleth}^{1}_{-1}(\hat{a}_{j}) \! - \! \hat{\mathfrak{Q}}_{1}
(\hat{a}_{j})(\hat{\mathfrak{Q}}_{0}(\hat{a}_{j}))^{-1} \right) \! + \! 
\hat{\mathfrak{Q}}_{1}(\hat{a}_{j}) \hat{\aleph}^{1}_{1}(\hat{a}_{j}) 
\hat{\aleph}^{1}_{-1}(\hat{a}_{j}) \right) \right. \nonumber \\
+&\left. \, \mi (s_{1} \! + \! t_{1}) \left(\hat{\beth}^{1}_{-1}(\hat{a}_{j}) 
\! - \! \hat{\beth}^{1}_{1}(\hat{a}_{j}) \! + \! \hat{\aleph}^{1}_{-1}
(\hat{a}_{j}) \hat{\daleth}^{1}_{1}(\hat{a}_{j}) \! - \! \hat{\aleph}^{1}_{1}
(\hat{a}_{j}) \hat{\daleth}^{1}_{-1}(\hat{a}_{j}) \right) \right), \label{eqprop42} \\
\hat{\mathbb{C}}_{12}(\hat{a}_{j}) =& \, (\hat{\kappa}_{1}(\hat{a}_{j}))^{2} 
\left(\mi s_{1} \left(\hat{\mathfrak{Q}}_{0}(\hat{a}_{j}) \left(2 \hat{\gimel}^{
1}_{1}(\hat{a}_{j}) \! + \! (\hat{\daleth}^{1}_{1}(\hat{a}_{j}))^{2} \right) 
\! + \! 2 \hat{\mathfrak{Q}}_{1}(\hat{a}_{j}) \hat{\daleth}^{1}_{1}
(\hat{a}_{j}) \right. \right. \nonumber \\
-&\left. \left. \, (\hat{\mathfrak{Q}}_{0}(\hat{a}_{j}))^{-1}(\hat{\aleph}^{
1}_{1}(\hat{a}_{j}))^{2} \! + \! \dfrac{1}{2} \hat{\mathfrak{Q}}_{2}
(\hat{a}_{j}) \right) \! + \! \mi t_{1} \left(2 \hat{\mathfrak{Q}}_{0}
(\hat{a}_{j}) \hat{\aleph}^{1}_{1}(\hat{a}_{j}) \hat{\beth}^{1}_{1}(\hat{a}_{j}) 
\right. \right. \nonumber \\
+&\left. \left. \, \hat{\mathfrak{Q}}_{1}(\hat{a}_{j})(\hat{\aleph}^{1}_{1}
(\hat{a}_{j}))^{2} \! + \! (\hat{\mathfrak{Q}}_{0}(\hat{a}_{j}))^{-1} \left(
\hat{\mathfrak{Q}}_{1}(\hat{a}_{j})(\hat{\mathfrak{Q}}_{0}(\hat{a}_{j}))^{-1} 
\! - \! 2 \hat{\daleth}^{1}_{1}(\hat{a}_{j}) \right) \right) \right. \nonumber \\
-&\left. \, 2(s_{1} \! - \! t_{1}) \left(\hat{\beth}^{1}_{1}(\hat{a}_{j}) \! + \! 
\hat{\aleph}^{1}_{1}(\hat{a}_{j}) \hat{\daleth}^{1}_{1}(\hat{a}_{j}) \right) 
\right), \label{eqprop43} \\
\hat{\mathbb{C}}_{21}(\hat{a}_{j}) =& \, (\hat{\kappa}_{2}(\hat{a}_{j}))^{2} 
\left(\mi s_{1} \left(\hat{\mathfrak{Q}}_{0}(\hat{a}_{j}) \left(2 \hat{\gimel}^{
1}_{-1}(\hat{a}_{j}) \! + \! (\hat{\daleth}^{1}_{-1}(\hat{a}_{j}))^{2} \right) \! 
+ \! 2 \hat{\mathfrak{Q}}_{1}(\hat{a}_{j}) \hat{\daleth}^{1}_{-1}(\hat{a}_{j}) 
\right. \right. \nonumber \\
-&\left. \left. \, (\hat{\mathfrak{Q}}_{0}(\hat{a}_{j}))^{-1}(\hat{\aleph}^{
1}_{-1}(\hat{a}_{j}))^{2} \! + \! \dfrac{1}{2} \hat{\mathfrak{Q}}_{2}
(\hat{a}_{j}) \right) \! + \! \mi t_{1} \left(2 \hat{\mathfrak{Q}}_{0}
(\hat{a}_{j}) \hat{\aleph}^{1}_{-1}(\hat{a}_{j}) \hat{\beth}^{1}_{-1}
(\hat{a}_{j}) \right. \right. \nonumber \\
+&\left. \left. \, \hat{\mathfrak{Q}}_{1}(\hat{a}_{j})(\hat{\aleph}^{1}_{-1}
(\hat{a}_{j}))^{2} \! + \! (\hat{\mathfrak{Q}}_{0}(\hat{a}_{j}))^{-1} \left(
\hat{\mathfrak{Q}}_{1}(\hat{a}_{j})(\hat{\mathfrak{Q}}_{0}(\hat{a}_{j}))^{-1} 
\! - \! 2 \hat{\daleth}^{1}_{-1}(\hat{a}_{j}) \right) \right) \right. \nonumber \\
+&\left. \, 2(s_{1} \! - \! t_{1}) \left(\hat{\beth}^{1}_{-1}(\hat{a}_{j}) \! + \! 
\hat{\aleph}^{1}_{-1}(\hat{a}_{j}) \hat{\daleth}^{1}_{-1}(\hat{a}_{j}) \right) 
\right), \label{eqprop44} 
\end{align}
where, for $\varepsilon_{1},\varepsilon_{2} \! = \! \pm 1$, $\hat{\kappa}_{1}
(\varsigma)$, $\hat{\kappa}_{2}(\varsigma)$, $\hat{\aleph}^{\varepsilon_{
1}}_{\varepsilon_{2}}(\varsigma)$, $\hat{\daleth}^{\varepsilon_{1}}_{
\varepsilon_{2}}(\varsigma)$, $\hat{\beth}^{\varepsilon_{1}}_{\varepsilon_{2}}
(\varsigma)$, and $\hat{\gimel}^{\varepsilon_{1}}_{\varepsilon_{2}}
(\varsigma)$ are defined by Equations~\eqref{eqprop13}--\eqref{eqprop23}, 
and
\begin{align}
\hat{\leftthreetimes}(\hat{a}_{N+1}) :=& \, (\hat{a}_{N+1} \! - \! 
\hat{b}_{0})^{1/2} \prod_{m=1}^{N}(\hat{a}_{N+1} \! - \! \hat{b}_{m})^{1/2}
(\hat{a}_{N+1} \! - \! \hat{a}_{m})^{1/2}, \label{eqprop45} \\
\hat{\leftthreetimes}(\hat{a}_{j}) :=& \, (-1)^{N+1-j}(\hat{a}_{j} \! - \! 
\hat{b}_{0})^{1/2}(\hat{a}_{N+1} \! - \! \hat{a}_{j})^{1/2}(\hat{b}_{j} \! - \! 
\hat{a}_{j})^{1/2} \prod_{m=1}^{j-1}(\hat{a}_{j} \! - \! \hat{b}_{m})^{1/2}
(\hat{a}_{j} \! - \! \hat{a}_{m})^{1/2} \nonumber \\
\times& \, \prod_{m^{\prime}=j+1}^{N}(\hat{b}_{m^{\prime}} \! - \! 
\hat{a}_{j})^{1/2}(\hat{a}_{m^{\prime}} \! - \! \hat{a}_{j})^{1/2}, 
\quad j \! = \! 1,2,\dotsc,N, \label{eqprop46} \\
\dfrac{\hat{\leftthreetimes}^{\prime}(\hat{a}_{N+1})}{\hat{\leftthreetimes}
(\hat{a}_{N+1})} =& \, \dfrac{1}{2} \left(\sum_{m=1}^{N} \left(\dfrac{1}{
\hat{a}_{N+1} \! - \! \hat{b}_{m}} \! + \! \dfrac{1}{\hat{a}_{N+1} \! - \! 
\hat{a}_{m}} \right) \! + \! \dfrac{1}{\hat{a}_{N+1} \! - \! \hat{b}_{0}} \right), 
\label{eqprop47} \\
\dfrac{\hat{\leftthreetimes}^{\prime}(\hat{a}_{j})}{\hat{\leftthreetimes}
(\hat{a}_{j})} =& \, \dfrac{1}{2} \left(\sum_{\substack{m=1\\m \neq j}}^{N} 
\left(\dfrac{1}{\hat{a}_{j} \! - \! \hat{b}_{m}} \! + \! \dfrac{1}{\hat{a}_{j} \! - 
\! \hat{a}_{m}} \right) \! + \! \dfrac{1}{\hat{a}_{j} \! - \! \hat{b}_{j}} \! + \! 
\dfrac{1}{\hat{a}_{j} \! - \! \hat{b}_{0}} \! + \! \dfrac{1}{\hat{a}_{j} \! - \! 
\hat{a}_{N+1}} \right), \quad j \! = \! 1,2,\dotsc,N, \label{eqprop48} \\
\hat{\mathfrak{Q}}_{0}(\hat{a}_{N+1})=& \, (\hat{a}_{N+1} \! - \! 
\hat{b}_{0})^{1/2} \prod_{m=1}^{N} \dfrac{(\hat{a}_{N+1} \! - \! \hat{b}_{
m})^{1/2}}{(\hat{a}_{N+1} \! - \! \hat{a}_{m})^{1/2}}, \label{eqprop49} \\
\dfrac{\hat{\mathfrak{Q}}_{1}(\hat{a}_{N+1})}{\hat{\mathfrak{Q}}_{0}
(\hat{a}_{N+1})}=& \, \dfrac{1}{2} \left(\sum_{m=1}^{N} \left(\dfrac{1}{
\hat{a}_{N+1} \! - \! \hat{b}_{m}} \! - \! \dfrac{1}{\hat{a}_{N+1} \! - \! 
\hat{a}_{m}} \right) \! + \! \dfrac{1}{\hat{a}_{N+1} \! - \! \hat{b}_{0}} \right), 
\label{eqprop50} \\
\dfrac{\hat{\mathfrak{Q}}_{2}(\hat{a}_{N+1})}{\hat{\mathfrak{Q}}_{0}
(\hat{a}_{N+1})}=& \, -\dfrac{1}{2} \left(\sum_{m=1}^{N} \left(\dfrac{1}{
(\hat{a}_{N+1} \! - \! \hat{b}_{m})^{2}} \! - \! \dfrac{1}{(\hat{a}_{N+1} \! - 
\! \hat{a}_{m})^{2}} \right) \! + \! \dfrac{1}{(\hat{a}_{N+1} \! - \! \hat{b}_{0})^{2}} 
\right) \nonumber \\
+& \, \dfrac{1}{4} \left(\sum_{m=1}^{N} \left(\dfrac{1}{\hat{a}_{N+1} 
\! - \! \hat{b}_{m}} \! - \! \dfrac{1}{\hat{a}_{N+1} \! - \! \hat{a}_{m}} \right) 
\! + \! \dfrac{1}{\hat{a}_{N+1} \! - \! \hat{b}_{0}} \right)^{2}, 
\label{eqprop51} \\
\hat{\mathfrak{Q}}_{0}(\hat{a}_{j})=& \, \dfrac{(\hat{a}_{j} \! - \! \hat{b}_{0})^{
1/2}(\hat{b}_{j} \! - \! \hat{a}_{j})^{1/2}}{(\hat{a}_{N+1} \! - \! \hat{a}_{j})^{1/2}} 
\prod_{m=1}^{j-1} \dfrac{(\hat{a}_{j} \! - \! \hat{b}_{m})^{1/2}}{(\hat{a}_{j} 
\! - \! \hat{a}_{m})^{1/2}} \prod_{m^{\prime}=j+1}^{N} \dfrac{(\hat{b}_{
m^{\prime}} \! - \! \hat{a}_{j})^{1/2}}{(\hat{a}_{m^{\prime}} \! - \! 
\hat{a}_{j})^{1/2}}, \quad j \! = \! 1,2,\dotsc,N, \label{eqprop52} \\
\dfrac{\hat{\mathfrak{Q}}_{1}(\hat{a}_{j})}{\hat{\mathfrak{Q}}_{0}(\hat{a}_{j})}
=& \, \dfrac{1}{2} \left(\sum_{\substack{m=1\\m \neq j}}^{N} \left(\dfrac{1}{
\hat{a}_{j} \! - \! \hat{b}_{m}} \! - \! \dfrac{1}{\hat{a}_{j} \! - \! \hat{a}_{m}} 
\right) \! + \! \dfrac{1}{\hat{a}_{j} \! - \! \hat{b}_{0}} \! - \! \dfrac{1}{
\hat{a}_{j} \! - \! \hat{a}_{N+1}} \! + \! \dfrac{1}{\hat{a}_{j} \! - \! \hat{b}_{j}} 
\right), \quad j \! = \! 1,2,\dotsc,N, \label{eqprop53} \\
\dfrac{\hat{\mathfrak{Q}}_{2}(\hat{a}_{j})}{\hat{\mathfrak{Q}}_{0}(\hat{a}_{j})}
=& \, -\dfrac{1}{2} \left(\sum_{\substack{m=1\\m \neq j}}^{N} \left(
\dfrac{1}{(\hat{a}_{j} \! - \! \hat{b}_{m})^{2}} \! - \! \dfrac{1}{(\hat{a}_{j} \! - \! 
\hat{a}_{m})^{2}} \right) \! + \! \dfrac{1}{(\hat{a}_{j} \! - \! \hat{b}_{0})^{2}} \! 
- \! \dfrac{1}{(\hat{a}_{j} \! - \! \hat{a}_{N+1})^{2}} \! + \! \dfrac{1}{(\hat{a}_{j} 
\! - \! \hat{b}_{j})^{2}} \right) \nonumber \\
+& \, \dfrac{1}{4} \left(\sum_{\substack{m=1\\m \neq j}}^{N} \left(\dfrac{1}{
\hat{a}_{j} \! - \! \hat{b}_{m}} \! - \! \dfrac{1}{\hat{a}_{j} \! - \! \hat{a}_{m}} 
\right) \! + \! \dfrac{1}{\hat{a}_{j} \! - \! \hat{b}_{0}} \! - \! \dfrac{1}{
\hat{a}_{j} \! - \! \hat{a}_{N+1}} \! + \! \dfrac{1}{\hat{a}_{j} \! - \! 
\hat{b}_{j}} \right)^{2}, \quad j \! = \! 1,2,\dotsc,N, \label{eqprop54}
\end{align}
$\hat{\alpha}_{0}(\hat{a}_{N+1})$, $\hat{\alpha}_{0}(\hat{a}_{j})$, 
$\hat{\alpha}_{1}(\hat{a}_{N+1})$, and $\hat{\alpha}_{1}
(\hat{a}_{j})$, $j \! = \! 1,2,\dotsc,N$, are given in 
Equations~\eqref{eqmaininf51}--\eqref{eqmaininf54}, respectively, 
$\hat{\alpha}_{2}(\hat{a}_{N+1}) \! = \! f^{\prime \prime}(\hat{a}_{N+1})/7$ 
and $\hat{\alpha}_{2}(\hat{a}_{j}) \! = \! f^{\prime \prime}(\hat{a}_{j})/7$, 
$j \! = \! 1,2,\dotsc,N$, where $f^{\prime \prime}(\hat{a}_{N+1})$ and 
$f^{\prime \prime}(\hat{a}_{j})$, $j \! = \! 1,2,\dotsc,N$, are given in 
Equations~\eqref{wm24franhat} and~\eqref{wm24frajhat}, respectively, 
and $(\mathrm{M}_{2}(\mathbb{C}) \! \ni)$ $\hat{\mathfrak{c}}^{r}_{m}
(n,k,z_{o};\hat{a}_{j}) \! =_{\underset{z_{o}=1+o(1)}{\mathscr{N},n \to \infty}} 
\! \mathcal{O}(1)$, $m \! \in \! \mathbb{N}_{0}$, $r \! \in \! \lbrace \triangleright,
\triangleleft \rbrace$$;$ and {\rm \pmb{(2)}} for $n \! \in \! \mathbb{N}$ and 
$k \! \in \! \lbrace 1,2,\dotsc,K \rbrace$ such that $\alpha_{p_{\mathfrak{s}}} 
\! := \! \alpha_{k} \! \neq \! \infty$, and for $z \! \in \! \partial \tilde{
\mathbb{U}}_{\tilde{\delta}_{\tilde{b}_{j-1}}} \cup \partial \tilde{\mathbb{U}}_{
\tilde{\delta}_{\tilde{a}_{j}}}$, $j \! = \! 1,2,\dotsc,N \! + \! 1$, with 
$w^{\Sigma_{\tilde{\mathcal{R}}}}_{+}(z) \! = \! \tilde{v}_{\tilde{\mathcal{R}}}
(z) \! - \! \mathrm{I}$,
\begin{align} \label{eqproptilb} 
w_{+}^{\Sigma_{\tilde{\mathcal{R}}}}(z) \underset{\underset{z_{o}=1+o(1)}{
\mathscr{N},n \to \infty}}{=}& \, \dfrac{1}{((n \! - \! 1)K \! + \! k)} \left(
\dfrac{(\tilde{\alpha}_{0}(\tilde{b}_{j-1}))^{-1}}{(z \! - \! \tilde{b}_{j-1})^{2}} 
\tilde{\boldsymbol{\mathrm{A}}}(\tilde{b}_{j-1}) \! + \! \dfrac{(\tilde{\alpha}_{0}
(\tilde{b}_{j-1}))^{-2}}{z \! - \! \tilde{b}_{j-1}} \left(\tilde{\alpha}_{0}
(\tilde{b}_{j-1}) \tilde{\boldsymbol{\mathrm{B}}}(\tilde{b}_{j-1}) \! - \! 
\tilde{\alpha}_{1}(\tilde{b}_{j-1}) \tilde{\boldsymbol{\mathrm{A}}}(\tilde{b}_{j-1}) 
\right) \right. \nonumber \\
&\left. \, +(\tilde{\alpha}_{0}(\tilde{b}_{j-1}))^{-2} \left(\tilde{\alpha}_{0}
(\tilde{b}_{j-1}) \left(\left(\dfrac{\tilde{\alpha}_{1}(\tilde{b}_{j-1})}{\tilde{\alpha}_{0}
(\tilde{b}_{j-1})} \right)^{2} \! - \! \dfrac{\tilde{\alpha}_{2}(\tilde{b}_{j-1})}{\tilde{
\alpha}_{0}(\tilde{b}_{j-1})} \right) \tilde{\boldsymbol{\mathrm{A}}}(\tilde{b}_{j
-1}) \! - \! \tilde{\alpha}_{1}(\tilde{b}_{j-1}) \tilde{\boldsymbol{\mathrm{B}}}
(\tilde{b}_{j-1}) \! + \! \tilde{\alpha}_{0}(\tilde{b}_{j-1}) \tilde{\boldsymbol{
\mathrm{C}}}(\tilde{b}_{j-1}) \right) \right) \nonumber \\
& \, +\mathcal{O} \left(\dfrac{1}{((n \! - \! 1)K \! + \! k)} \sum_{m=1}^{
\infty} \tilde{c}_{m}^{\triangleright}(n,k,z_{o};\tilde{b}_{j-1})(z \! - \! 
\tilde{b}_{j-1})^{m} \right) \nonumber \\
& \, + \mathcal{O} \left(\dfrac{1}{((n \! - \! 1)K \! + \! k)^{2}(z \! - \! 
\tilde{b}_{j-1})^{3}} \sum_{m=0}^{\infty} \tilde{c}_{m}^{\triangleleft}(n,k,z_{o};
\tilde{b}_{j-1})(z \! - \! \tilde{b}_{j-1})^{m} \right), \quad z \! \in \! \partial 
\tilde{\mathbb{U}}_{\tilde{\delta}_{\tilde{b}_{j-1}}}, \quad j \! = \! 1,2,\dotsc,
N \! + \! 1,
\end{align}
and
\begin{align} \label{eqproptila} 
w_{+}^{\Sigma_{\hat{\mathcal{R}}}}(z) \underset{\underset{z_{o}=1+o(1)}{
\mathscr{N},n \to \infty}}{=}& \, \dfrac{1}{((n \! - \! 1)K \! + \! k)} \left(
\dfrac{(\tilde{\alpha}_{0}(\tilde{a}_{j}))^{-1}}{(z \! - \! \tilde{a}_{j})^{2}} \tilde{
\boldsymbol{\mathrm{A}}}(\tilde{a}_{j}) \! + \! \dfrac{(\tilde{\alpha}_{0}
(\tilde{a}_{j}))^{-2}}{z \! - \! \tilde{a}_{j}} \left(\tilde{\alpha}_{0}(\tilde{a}_{j}) 
\tilde{\boldsymbol{\mathrm{B}}}(\tilde{a}_{j}) \! - \! \tilde{\alpha}_{1}(\tilde{a}_{j}) 
\tilde{\boldsymbol{\mathrm{A}}}(\tilde{a}_{j}) \right) \right. \nonumber \\
&\left. \, +(\tilde{\alpha}_{0}(\tilde{a}_{j}))^{-2} \left(\tilde{\alpha}_{0}(\tilde{a}_{j}) 
\left(\left(\dfrac{\tilde{\alpha}_{1}(\tilde{a}_{j})}{\tilde{\alpha}_{0}(\tilde{a}_{j})} 
\right)^{2} \! - \! \dfrac{\tilde{\alpha}_{2}(\tilde{a}_{j})}{\tilde{\alpha}_{0}
(\tilde{a}_{j})} \right) \tilde{\boldsymbol{\mathrm{A}}}(\tilde{a}_{j}) \! - \! 
\tilde{\alpha}_{1}(\tilde{a}_{j}) \tilde{\boldsymbol{\mathrm{B}}}(\tilde{a}_{j}) \! 
+ \! \tilde{\alpha}_{0}(\tilde{a}_{j}) \tilde{\boldsymbol{\mathrm{C}}}(\tilde{a}_{j}) 
\right) \right) \nonumber \\
& \, +\mathcal{O} \left(\dfrac{1}{((n \! - \! 1)K \! + \! k)} \sum_{m=1}^{
\infty} \tilde{c}_{m}^{\triangleright}(n,k,z_{o};\tilde{a}_{j})(z \! - \! 
\tilde{a}_{j})^{m} \right) \nonumber \\
& \, + \mathcal{O} \left(\dfrac{1}{((n \! - \! 1)K \! + \! k)^{2}(z \! - \! 
\tilde{a}_{j})^{3}} \sum_{m=0}^{\infty} \tilde{c}_{m}^{\triangleleft}(n,k,z_{o};
\tilde{a}_{j})(z \! - \! \tilde{a}_{j})^{m} \right), \quad z \! \in \! \partial  
\tilde{\mathbb{U}}_{\tilde{\delta}_{\tilde{a}_{j}}}, \quad j \! = \! 1,2,\dotsc,
N \! + \! 1,
\end{align}
where, for $j \! = \! 1,2,\dotsc,N \! + \! 1$,
\begin{gather}
\tilde{\boldsymbol{\mathrm{A}}}(\tilde{b}_{j-1}) \! := \! 
\begin{pmatrix}
\widetilde{\mathbb{K}}_{11} & \widetilde{\mathbb{K}}_{12} \\
\widetilde{\mathbb{K}}_{21} & \widetilde{\mathbb{K}}_{22}
\end{pmatrix} 
\begin{pmatrix}
\tilde{\mathbb{A}}_{11}(\tilde{b}_{j-1}) & \tilde{\mathbb{A}}_{12}
(\tilde{b}_{j-1}) \\
\tilde{\mathbb{A}}_{21}(\tilde{b}_{j-1}) & \tilde{\mathbb{A}}_{22}
(\tilde{b}_{j-1})
\end{pmatrix} 
\begin{pmatrix}
\widetilde{\mathbb{K}}_{22} & -\widetilde{\mathbb{K}}_{12} \\
-\widetilde{\mathbb{K}}_{21} & \widetilde{\mathbb{K}}_{11}
\end{pmatrix} \me^{\mi ((n-1)K+k) \tilde{\mho}_{j-1}}, \label{eqprop55} \\
\tilde{\boldsymbol{\mathrm{B}}}(\tilde{b}_{j-1}) \! := \! 
\begin{pmatrix}
\widetilde{\mathbb{K}}_{11} & \widetilde{\mathbb{K}}_{12} \\
\widetilde{\mathbb{K}}_{21} & \widetilde{\mathbb{K}}_{22}
\end{pmatrix} 
\begin{pmatrix}
\tilde{\mathbb{B}}_{11}(\tilde{b}_{j-1}) & \tilde{\mathbb{B}}_{12}
(\tilde{b}_{j-1}) \\
\tilde{\mathbb{B}}_{21}(\tilde{b}_{j-1}) & \tilde{\mathbb{B}}_{22}
(\tilde{b}_{j-1})
\end{pmatrix} 
\begin{pmatrix}
\widetilde{\mathbb{K}}_{22} & -\widetilde{\mathbb{K}}_{12} \\
-\widetilde{\mathbb{K}}_{21} & \widetilde{\mathbb{K}}_{11}
\end{pmatrix} \me^{\mi ((n-1)K+k) \tilde{\mho}_{j-1}}, \label{eqprop56} \\
\tilde{\boldsymbol{\mathrm{C}}}(\tilde{b}_{j-1}) \! := \! 
\begin{pmatrix}
\widetilde{\mathbb{K}}_{11} & \widetilde{\mathbb{K}}_{12} \\
\widetilde{\mathbb{K}}_{21} & \widetilde{\mathbb{K}}_{22}
\end{pmatrix} 
\begin{pmatrix}
\tilde{\mathbb{C}}_{11}(\tilde{b}_{j-1}) & \tilde{\mathbb{C}}_{12}
(\tilde{b}_{j-1}) \\
\tilde{\mathbb{C}}_{21}(\tilde{b}_{j-1}) & \tilde{\mathbb{C}}_{22}
(\tilde{b}_{j-1})
\end{pmatrix} 
\begin{pmatrix}
\widetilde{\mathbb{K}}_{22} & -\widetilde{\mathbb{K}}_{12} \\
-\widetilde{\mathbb{K}}_{21} & \widetilde{\mathbb{K}}_{11}
\end{pmatrix} \me^{\mi ((n-1)K+k) \tilde{\mho}_{j-1}}, \label{eqprop57}
\end{gather}
with $\widetilde{\mathbb{K}}$ defined in item~{\rm \pmb{(2)}} of 
Lemma~\ref{lem4.5}, $\tilde{\mho}_{m}$, $m \! = \! 0,1,\dotsc,N \! + \! 
1$, defined in the corresponding item~{\rm (ii)} of Remark~\ref{rem4.4},
\begin{gather}
\tilde{\mathbb{A}}_{11}(\tilde{b}_{j-1}) \! = \! -\tilde{\mathbb{A}}_{22}
(\tilde{b}_{j-1}) \! = \! -s_{1} \tilde{\kappa}_{1}(\tilde{b}_{j-1}) \tilde{
\kappa}_{2}(\tilde{b}_{j-1})(\tilde{\mathfrak{Q}}_{0}(\tilde{b}_{j-1}))^{-1}, 
\label{eqprop58} \\
\tilde{\mathbb{A}}_{12}(\tilde{b}_{j-1}) \! = \! -\mi s_{1}(\tilde{\kappa}_{1}
(\tilde{b}_{j-1}))^{2}(\tilde{\mathfrak{Q}}_{0}(\tilde{b}_{j-1}))^{-1}, \qquad 
\tilde{\mathbb{A}}_{21}(\tilde{b}_{j-1}) \! = \! -\mi s_{1}(\tilde{\kappa}_{2}
(\tilde{b}_{j-1}))^{2}(\tilde{\mathfrak{Q}}_{0}(\tilde{b}_{j-1}))^{-1}, 
\label{eqprop59}
\end{gather}
\begin{align}
\tilde{\mathbb{B}}_{11}(\tilde{b}_{j-1}) \! = \! -\tilde{\mathbb{B}}_{22}
(\tilde{b}_{j-1}) =& \, \tilde{\kappa}_{1}(\tilde{b}_{j-1}) \tilde{\kappa}_{2}
(\tilde{b}_{j-1}) \left(-s_{1}(\tilde{\mathfrak{Q}}_{0}(\tilde{b}_{j-1}))^{-1} 
\left(\tilde{\daleth}^{1}_{1}(\tilde{b}_{j-1}) \! + \! \tilde{\daleth}^{1}_{-1}
(\tilde{b}_{j-1}) \! - \! \tilde{\mathfrak{Q}}_{1}(\tilde{b}_{j-1})
(\tilde{\mathfrak{Q}}_{0}(\tilde{b}_{j-1}))^{-1} \right) \right. \nonumber \\
-&\left. \, t_{1} \left(\tilde{\mathfrak{Q}}_{0}(\tilde{b}_{j-1}) \! + \! 
(\tilde{\mathfrak{Q}}_{0}(\tilde{b}_{j-1}))^{-1} \tilde{\aleph}^{1}_{1}
(\tilde{b}_{j-1}) \tilde{\aleph}^{1}_{-1}(\tilde{b}_{j-1}) \right) \! + \! 
\mi (s_{1} \! + \! t_{1}) \left(\tilde{\aleph}^{1}_{-1}(\tilde{b}_{j-1}) \! - \! 
\tilde{\aleph}^{1}_{1}(\tilde{b}_{j-1}) \right) \right), \label{eqprop60} \\
\tilde{\mathbb{B}}_{12}(\tilde{b}_{j-1}) =& \, (\tilde{\kappa}_{1}
(\tilde{b}_{j-1}))^{2} \left(-\mi s_{1}(\tilde{\mathfrak{Q}}_{0}(\tilde{b}_{j
-1}))^{-1} \left(2 \tilde{\daleth}^{1}_{1}(\tilde{b}_{j-1}) \! - \! \tilde{
\mathfrak{Q}}_{1}(\tilde{b}_{j-1})(\tilde{\mathfrak{Q}}_{0}(\tilde{b}_{j-
1}))^{-1} \right) \right. \nonumber \\
+&\left. \, \mi t_{1} \left(\tilde{\mathfrak{Q}}_{0}(\tilde{b}_{j-1}) \! - \! 
(\tilde{\mathfrak{Q}}_{0}(\tilde{b}_{j-1}))^{-1}(\tilde{\aleph}^{1}_{1}
(\tilde{b}_{j-1}))^{2} \right) \! + \! 2(s_{1} \! - \! t_{1}) \tilde{\aleph}^{1}_{1}
(\tilde{b}_{j-1}) \right), \label{eqprop61} \\
\tilde{\mathbb{B}}_{21}(\tilde{b}_{j-1}) =& \, (\tilde{\kappa}_{2}
(\tilde{b}_{j-1}))^{2} \left(-\mi s_{1}(\tilde{\mathfrak{Q}}_{0}(\tilde{b}_{j
-1}))^{-1} \left(2 \tilde{\daleth}^{1}_{-1}(\tilde{b}_{j-1}) \! - \! \tilde{
\mathfrak{Q}}_{1}(\tilde{b}_{j-1})(\tilde{\mathfrak{Q}}_{0}(\tilde{b}_{j-
1}))^{-1} \right) \right. \nonumber \\
+&\left. \, \mi t_{1} \left(\tilde{\mathfrak{Q}}_{0}(\tilde{b}_{j-1}) \! - \! 
(\tilde{\mathfrak{Q}}_{0}(\tilde{b}_{j-1}))^{-1}(\tilde{\aleph}^{1}_{-1}
(\tilde{b}_{j-1}))^{2} \right) \! - \! 2(s_{1} \! - \! t_{1}) \tilde{\aleph}^{
1}_{-1}(\tilde{b}_{j-1}) \right), \label{eqprop62} \\
\tilde{\mathbb{C}}_{11}(\tilde{b}_{j-1}) \! = \! -\tilde{\mathbb{C}}_{22}
(\tilde{b}_{j-1}) =& \, \tilde{\kappa}_{1}(\tilde{b}_{j-1}) \tilde{\kappa}_{2}
(\tilde{b}_{j-1}) \left(s_{1} \left(-\tilde{\mathfrak{Q}}_{0}(\tilde{b}_{j-1}) 
\tilde{\aleph}^{1}_{1}(\tilde{b}_{j-1}) \tilde{\aleph}^{1}_{-1}(\tilde{b}_{j-1}) 
\! - \! (\tilde{\mathfrak{Q}}_{0}(\tilde{b}_{j-1}))^{-3}(\tilde{\mathfrak{Q}}_{1}
(\tilde{b}_{j-1}))^{2} \right. \right. \nonumber \\
+&\left. \left. \, \dfrac{1}{2} \tilde{\mathfrak{Q}}_{2}(\tilde{b}_{j-1})
(\tilde{\mathfrak{Q}}_{0}(\tilde{b}_{j-1}))^{-2} \! + \! \tilde{\mathfrak{Q}}_{1}
(\tilde{b}_{j-1})(\tilde{\mathfrak{Q}}_{0}(\tilde{b}_{j-1}))^{-2} \left(\tilde{
\daleth}^{1}_{1}(\tilde{b}_{j-1}) \! + \! \tilde{\daleth}^{1}_{-1}
(\tilde{b}_{j-1}) \right) \right. \right. \nonumber \\
-&\left. \left. \, (\tilde{\mathfrak{Q}}_{0}(\tilde{b}_{j-1}))^{-1} \left(
\tilde{\gimel}^{1}_{1}(\tilde{b}_{j-1}) \! + \! \tilde{\gimel}^{1}_{-1}
(\tilde{b}_{j-1}) \! + \! \tilde{\daleth}^{1}_{1}(\tilde{b}_{j-1}) \tilde{\daleth}^{
1}_{-1}(\tilde{b}_{j-1}) \right) \right) \! + \! t_{1} \left(-\tilde{\mathfrak{Q}}_{1}
(\tilde{b}_{j-1}) \! - \! \tilde{\mathfrak{Q}}_{0}(\tilde{b}_{j-1}) \right. \right. 
\nonumber \\
\times&\left. \left. \, \left(\tilde{\daleth}^{1}_{1}(\tilde{b}_{j-1}) \! + \! 
\tilde{\daleth}^{1}_{-1}(\tilde{b}_{j-1}) \right) \! + \! \tilde{\mathfrak{Q}}_{1}
(\tilde{b}_{j-1})(\tilde{\mathfrak{Q}}_{0}(\tilde{b}_{j-1}))^{-2} \tilde{\aleph}^{
1}_{1}(\tilde{b}_{j-1}) \tilde{\aleph}^{1}_{-1}(\tilde{b}_{j-1}) \! - \! (\tilde{
\mathfrak{Q}}_{0}(\tilde{b}_{j-1}))^{-1} \right. \right. \nonumber \\
\times&\left. \left. \, \left(\tilde{\aleph}^{1}_{1}(\tilde{b}_{j-1}) \tilde{
\beth}^{1}_{-1}(\tilde{b}_{j-1}) \! + \! \tilde{\aleph}^{1}_{-1}(\tilde{b}_{j-1}) 
\tilde{\beth}^{1}_{1}(\tilde{b}_{j-1}) \right) \right) \! + \! \mi (s_{1} \! + \! 
t_{1}) \left(\tilde{\beth}^{1}_{-1}(\tilde{b}_{j-1}) \! - \! \tilde{\beth}^{1}_{1}
(\tilde{b}_{j-1}) \right. \right. \nonumber \\
-&\left. \left. \, \tilde{\aleph}^{1}_{1}(\tilde{b}_{j-1}) \tilde{\daleth}^{1}_{-1}
(\tilde{b}_{j-1}) \! + \! \tilde{\aleph}^{1}_{-1}(\tilde{b}_{j-1}) \tilde{\daleth}^{
1}_{1}(\tilde{b}_{j-1}) \right) \right), \label{eqprop63} \\
\tilde{\mathbb{C}}_{12}(\tilde{b}_{j-1}) =& \, (\tilde{\kappa}_{1}(\tilde{b}_{j
-1}))^{2} \left(\mi s_{1} \left(\tilde{\mathfrak{Q}}_{0}(\tilde{b}_{j-1})(\tilde{
\aleph}^{1}_{1}(\tilde{b}_{j-1}))^{2} \! - \! (\tilde{\mathfrak{Q}}_{0}(\tilde{
b}_{j-1}))^{-3}(\tilde{\mathfrak{Q}}_{1}(\tilde{b}_{j-1}))^{2} \right. \right. 
\nonumber \\
+&\left. \left. \, \dfrac{1}{2} \tilde{\mathfrak{Q}}_{2}(\tilde{b}_{j-1})
(\tilde{\mathfrak{Q}}_{0}(\tilde{b}_{j-1}))^{-2} \! + \! 2 \tilde{\mathfrak{Q}}_{1}
(\tilde{b}_{j-1})(\tilde{\mathfrak{Q}}_{0}(\tilde{b}_{j-1}))^{-2} \tilde{\daleth}^{
1}_{1}(\tilde{b}_{j-1}) \! - \! (\tilde{\mathfrak{Q}}_{0}(\tilde{b}_{j-1}))^{-1} 
\right. \right. \nonumber \\
\times&\left. \left. \, \left(2 \tilde{\gimel}^{1}_{1}(\tilde{b}_{j-1}) \! + \! 
(\tilde{\daleth}^{1}_{1}(\tilde{b}_{j-1}))^{2} \right) \right) \! + \! \mi t_{1} 
\left(2 \tilde{\mathfrak{Q}}_{0}(\tilde{b}_{j-1}) \tilde{\daleth}^{1}_{1}
(\tilde{b}_{j-1}) \! + \! \tilde{\mathfrak{Q}}_{1}(\tilde{b}_{j-1}) \right. \right. 
\nonumber \\
+&\left. \left. \, \tilde{\mathfrak{Q}}_{1}(\tilde{b}_{j-1})(\tilde{\mathfrak{
Q}}_{0}(\tilde{b}_{j-1}))^{-2}(\tilde{\aleph}^{1}_{1}(\tilde{b}_{j-1}))^{2}  \! - 
\! 2(\tilde{\mathfrak{Q}}_{0}(\tilde{b}_{j-1}))^{-1} \tilde{\aleph}^{1}_{1}
(\tilde{b}_{j-1}) \tilde{\beth}^{1}_{1}(\tilde{b}_{j-1}) \right) \right. \nonumber \\
+&\left. \, 2(s_{1} \! - \! t_{1}) \left(\tilde{\beth}^{1}_{1}(\tilde{b}_{j-1}) 
\! + \! \tilde{\aleph}^{1}_{1}(\tilde{b}_{j-1}) \tilde{\daleth}^{1}_{1}
(\tilde{b}_{j-1}) \right) \right), \label{eqprop64} \\
\tilde{\mathbb{C}}_{21}(\tilde{b}_{j-1}) =& \, (\tilde{\kappa}_{2}(\tilde{b}_{j
-1}))^{2} \left(\mi s_{1} \left(\tilde{\mathfrak{Q}}_{0}(\tilde{b}_{j-1})(\tilde{
\aleph}^{1}_{-1}(\tilde{b}_{j-1}))^{2} \! - \! (\tilde{\mathfrak{Q}}_{0}
(\tilde{b}_{j-1}))^{-3}(\tilde{\mathfrak{Q}}_{1}(\tilde{b}_{j-1}))^{2} \right. 
\right. \nonumber \\
+&\left. \left. \, \dfrac{1}{2} \tilde{\mathfrak{Q}}_{2}(\tilde{b}_{j-1})
(\tilde{\mathfrak{Q}}_{0}(\tilde{b}_{j-1}))^{-2} \! + \! 2 \tilde{\mathfrak{Q}}_{1}
(\tilde{b}_{j-1})(\tilde{\mathfrak{Q}}_{0}(\tilde{b}_{j-1}))^{-2} \tilde{\daleth}^{
1}_{-1}(\tilde{b}_{j-1}) \! - \! (\tilde{\mathfrak{Q}}_{0}(\tilde{b}_{j-1}))^{-1} 
\right. \right. \nonumber \\
\times&\left. \left. \, \left(2 \tilde{\gimel}^{1}_{-1}(\tilde{b}_{j-1}) \! + \! 
(\tilde{\daleth}^{1}_{-1}(\tilde{b}_{j-1}))^{2} \right) \right) \! + \! \mi t_{1} 
\left(2 \tilde{\mathfrak{Q}}_{0}(\tilde{b}_{j-1}) \tilde{\daleth}^{1}_{-1}
(\tilde{b}_{j-1}) \! + \! \tilde{\mathfrak{Q}}_{1}(\tilde{b}_{j-1}) \right. \right. 
\nonumber \\
+&\left. \left. \, \tilde{\mathfrak{Q}}_{1}(\tilde{b}_{j-1})(\tilde{\mathfrak{
Q}}_{0}(\tilde{b}_{j-1}))^{-2}(\tilde{\aleph}^{1}_{-1}(\tilde{b}_{j-1}))^{2}  \! 
- \! 2(\tilde{\mathfrak{Q}}_{0}(\tilde{b}_{j-1}))^{-1} \tilde{\aleph}^{1}_{-1}
(\tilde{b}_{j-1}) \tilde{\beth}^{1}_{-1}(\tilde{b}_{j-1}) \right) \right. 
\nonumber \\
-&\left. \, 2(s_{1} \! - \! t_{1}) \left(\tilde{\beth}^{1}_{-1}(\tilde{b}_{j-1}) 
\! + \! \tilde{\aleph}^{1}_{-1}(\tilde{b}_{j-1}) \tilde{\daleth}^{1}_{-1}
(\tilde{b}_{j-1}) \right) \right), \label{eqprop65}
\end{align}
where, for $\varepsilon_{1},\varepsilon_{2} \! = \! \pm 1$,
\begin{align}
\tilde{\kappa}_{1}(\varsigma) =& \, \dfrac{\tilde{\boldsymbol{\theta}}
(\tilde{\boldsymbol{u}}_{+}(\varsigma) \! - \! \frac{1}{2 \pi}((n \! - \! 1)K \! 
+ \! k) \tilde{\boldsymbol{\Omega}} \! + \! \tilde{\boldsymbol{d}})}{\tilde{
\boldsymbol{\theta}}(\tilde{\boldsymbol{u}}_{+}(\varsigma) \! + \! \tilde{
\boldsymbol{d}})}, \qquad \quad \tilde{\kappa}_{2}(\varsigma) = \dfrac{
\tilde{\boldsymbol{\theta}}(\tilde{\boldsymbol{u}}_{+}(\varsigma) \! - \! 
\frac{1}{2 \pi}((n \! - \! 1)K \! + \! k) \tilde{\boldsymbol{\Omega}} \! - \! 
\tilde{\boldsymbol{d}})}{\tilde{\boldsymbol{\theta}}(\tilde{\boldsymbol{
u}}_{+}(\varsigma) \! - \! \tilde{\boldsymbol{d}})}, \label{eqprop66} \\
\tilde{\aleph}^{\varepsilon_{1}}_{\varepsilon_{2}}(\varsigma) =& \, 
-\dfrac{\tilde{\mathfrak{u}}(\varepsilon_{1},\varepsilon_{2},\bm{0};
\varsigma)}{\tilde{\boldsymbol{\theta}}(\varepsilon_{1} \tilde{
\boldsymbol{u}}_{+}(\varsigma) \! + \! \varepsilon_{2} \tilde{
\boldsymbol{d}})} \! + \! \dfrac{\tilde{\mathfrak{u}}(\varepsilon_{1},
\varepsilon_{2},\tilde{\boldsymbol{\Omega}};\varsigma)}{\tilde{
\boldsymbol{\theta}}(\varepsilon_{1} \tilde{\boldsymbol{u}}_{+}(\varsigma) 
\! - \! \frac{1}{2 \pi}((n \! - \! 1)K \! + \! k) \tilde{\boldsymbol{\Omega}} 
\! + \! \varepsilon_{2} \tilde{\boldsymbol{d}})}, \label{eqprop67} \\
\tilde{\daleth}^{\varepsilon_{1}}_{\varepsilon_{2}}(\varsigma) =& \, 
-\dfrac{\tilde{\mathfrak{v}}(\varepsilon_{1},\varepsilon_{2},\bm{0};
\varsigma)}{\tilde{\boldsymbol{\theta}}(\varepsilon_{1} 
\tilde{\boldsymbol{u}}_{+}(\varsigma) \! + \! \varepsilon_{2} \tilde{
\boldsymbol{d}})} \! + \! \dfrac{\tilde{\mathfrak{v}}(\varepsilon_{1},
\varepsilon_{2},\tilde{\boldsymbol{\Omega}};\varsigma)}{\tilde{
\boldsymbol{\theta}}(\varepsilon_{1} \tilde{\boldsymbol{u}}_{+}(\varsigma) 
\! - \! \frac{1}{2 \pi}((n \! - \! 1)K \! + \! k) \tilde{\boldsymbol{\Omega}} 
\! + \! \varepsilon_{2} \tilde{\boldsymbol{d}})} \! - \! \left(\dfrac{\tilde{
\mathfrak{u}}(\varepsilon_{1},\varepsilon_{2},\bm{0};\varsigma)}{\tilde{
\boldsymbol{\theta}}(\varepsilon_{1} \tilde{\boldsymbol{u}}_{+}(\varsigma) 
\! + \! \varepsilon_{2} \tilde{\boldsymbol{d}})} \right)^{2} \nonumber \\
+& \, \dfrac{\tilde{\mathfrak{u}}(\varepsilon_{1},\varepsilon_{2},\bm{0};
\varsigma) \tilde{\mathfrak{u}}(\varepsilon_{1},\varepsilon_{2},\tilde{
\boldsymbol{\Omega}};\varsigma)}{\tilde{\boldsymbol{\theta}}
(\varepsilon_{1} \tilde{\boldsymbol{u}}_{+}(\varsigma) \! + \! \varepsilon_{2} 
\tilde{\boldsymbol{d}}) \tilde{\boldsymbol{\theta}}(\varepsilon_{1} \tilde{
\boldsymbol{u}}_{+}(\varsigma) \! - \! \frac{1}{2 \pi}((n \! - \! 1)K \! + \! k) 
\tilde{\boldsymbol{\Omega}} \! + \! \varepsilon_{2} \tilde{\boldsymbol{d}})}, 
\label{eqprop68} \\
\tilde{\beth}^{\varepsilon_{1}}_{\varepsilon_{2}}(\xi) =& \, -\dfrac{\tilde{
\mathfrak{w}}(\varepsilon_{1},\varepsilon_{2},\bm{0};\varsigma)}{\tilde{
\boldsymbol{\theta}}(\varepsilon_{1} \tilde{\boldsymbol{u}}_{+}(\varsigma) \! 
+ \! \varepsilon_{2} \tilde{\boldsymbol{d}})} \! + \! \dfrac{\tilde{\mathfrak{w}}
(\varepsilon_{1},\varepsilon_{2},\tilde{\boldsymbol{\Omega}};\varsigma)}{
\tilde{\boldsymbol{\theta}}(\varepsilon_{1} \tilde{\boldsymbol{u}}_{+}
(\varsigma) \! - \! \frac{1}{2 \pi}((n \! - \! 1)K \! + \! k) \tilde{\boldsymbol{
\Omega}} \! + \! \varepsilon_{2} \tilde{\boldsymbol{d}})} \! + \! \dfrac{2 
\tilde{\mathfrak{u}}(\varepsilon_{1},\varepsilon_{2},\bm{0};\varsigma) 
\tilde{\mathfrak{v}}(\varepsilon_{1},\varepsilon_{2},\bm{0};\varsigma)}{
(\tilde{\boldsymbol{\theta}}(\varepsilon_{1} \tilde{\boldsymbol{u}}_{+}
(\varsigma) \! + \! \varepsilon_{2} \tilde{\boldsymbol{d}}))^{2}} \nonumber \\
-& \, \dfrac{\tilde{\mathfrak{v}}(\varepsilon_{1},\varepsilon_{2},\bm{0};\varsigma) 
\tilde{\mathfrak{u}}(\varepsilon_{1},\varepsilon_{2},\tilde{\boldsymbol{\Omega}};
\varsigma)}{\tilde{\boldsymbol{\theta}}(\varepsilon_{1} \tilde{\boldsymbol{u}}_{
+}(\varsigma) \! + \! \varepsilon_{2} \tilde{\boldsymbol{d}}) \tilde{\boldsymbol{
\theta}}(\varepsilon_{1} \tilde{\boldsymbol{u}}_{+}(\varsigma) \! - \! \frac{1}{2 
\pi}((n \! - \! 1)K \! + \! k) \tilde{\boldsymbol{\Omega}} \! + \! \varepsilon_{2} 
\tilde{\boldsymbol{d}})} \! + \! \left(\dfrac{\tilde{\mathfrak{u}}(\varepsilon_{1},
\varepsilon_{2},\bm{0};\varsigma)}{\tilde{\boldsymbol{\theta}}(\varepsilon_{1} 
\tilde{\boldsymbol{u}}_{+}(\varsigma) \! + \! \varepsilon_{2} \tilde{\boldsymbol{
d}})} \right)^{3} \nonumber \\
-& \, \dfrac{\tilde{\mathfrak{u}}(\varepsilon_{1},\varepsilon_{2},\bm{0};\varsigma) 
\tilde{\mathfrak{v}}(\varepsilon_{1},\varepsilon_{2},\tilde{\boldsymbol{\Omega}};
\varsigma)}{\tilde{\boldsymbol{\theta}}(\varepsilon_{1} \tilde{\boldsymbol{u}}_{
+}(\varsigma) \! + \! \varepsilon_{2} \tilde{\boldsymbol{d}}) \tilde{\boldsymbol{
\theta}}(\varepsilon_{1} \tilde{\boldsymbol{u}}_{+}(\varsigma) \! - \! \frac{1}{
2 \pi}((n \! - \! 1)K \! + \! k) \tilde{\boldsymbol{\Omega}} \! + \! \varepsilon_{2} 
\tilde{\boldsymbol{d}})} \! - \! \left(\dfrac{\tilde{\mathfrak{u}}(\varepsilon_{1},
\varepsilon_{2},\bm{0};\varsigma)}{\tilde{\boldsymbol{\theta}}(\varepsilon_{1} 
\tilde{\boldsymbol{u}}_{+}(\varsigma) \! + \! \varepsilon_{2} \tilde{\boldsymbol{
d}})} \right)^{2} \nonumber \\
\times& \, \dfrac{\tilde{\mathfrak{u}}(\varepsilon_{1},\varepsilon_{2},\tilde{
\boldsymbol{\Omega}};\varsigma)}{\tilde{\boldsymbol{\theta}}(\varepsilon_{1} 
\tilde{\boldsymbol{u}}_{+}(\varsigma) \! - \! \frac{1}{2 \pi}((n \! - \! 1)K \! + 
\! k) \tilde{\boldsymbol{\Omega}} \! + \! \varepsilon_{2} \tilde{\boldsymbol{
d}})}, \label{eqprop69} \\
\tilde{\gimel}^{\varepsilon_{1}}_{\varepsilon_{2}}(\varsigma) =& \, -\dfrac{
\tilde{\mathfrak{z}}(\varepsilon_{1},\varepsilon_{2},\bm{0};\varsigma)}{\tilde{
\boldsymbol{\theta}}(\varepsilon_{1} \tilde{\boldsymbol{u}}_{+}(\varsigma) \! 
+ \! \varepsilon_{2} \tilde{\boldsymbol{d}})} \! + \! \dfrac{\tilde{\mathfrak{z}}
(\varepsilon_{1},\varepsilon_{2},\tilde{\boldsymbol{\Omega}};\varsigma)}{
\tilde{\boldsymbol{\theta}}(\varepsilon_{1} \tilde{\boldsymbol{u}}_{+}
(\varsigma) \! - \! \frac{1}{2 \pi}((n \! - \! 1)K \! + \! k) \tilde{\boldsymbol{
\Omega}} \! + \! \varepsilon_{2} \tilde{\boldsymbol{d}})} \! + \! \left(\dfrac{
\tilde{\mathfrak{v}}(\varepsilon_{1},\varepsilon_{2},\bm{0};\varsigma)}{\tilde{
\boldsymbol{\theta}}(\varepsilon_{1} \tilde{\boldsymbol{u}}_{+}(\varsigma) 
\! + \! \varepsilon_{2} \tilde{\boldsymbol{d}})} \right)^{2} \nonumber \\
-& \, \dfrac{\tilde{\mathfrak{v}}(\varepsilon_{1},\varepsilon_{2},\bm{0};\varsigma) 
\tilde{\mathfrak{v}}(\varepsilon_{1},\varepsilon_{2},\tilde{\boldsymbol{\Omega}};
\varsigma)}{\tilde{\boldsymbol{\theta}}(\varepsilon_{1} \tilde{\boldsymbol{u}}_{+}
(\varsigma) \! + \! \varepsilon_{2} \tilde{\boldsymbol{d}}) \tilde{\boldsymbol{
\theta}}(\varepsilon_{1} \tilde{\boldsymbol{u}}_{+}(\varsigma) \! - \! \frac{1}{2 
\pi}((n \! - \! 1)K \! + \! k) \tilde{\boldsymbol{\Omega}} \! + \! \varepsilon_{2} 
\tilde{\boldsymbol{d}})} \! - \! \dfrac{2 \tilde{\mathfrak{u}}(\varepsilon_{1},
\varepsilon_{2},\bm{0};\varsigma) \tilde{\mathfrak{w}}(\varepsilon_{1},
\varepsilon_{2},\bm{0};\varsigma)}{(\tilde{\boldsymbol{\theta}}(\varepsilon_{1} 
\tilde{\boldsymbol{u}}_{+}(\varsigma) \! + \! \varepsilon_{2} \tilde{\boldsymbol{
d}}))^{2}} \nonumber \\
+& \, \dfrac{\tilde{\mathfrak{w}}(\varepsilon_{1},\varepsilon_{2},\bm{0};\varsigma) 
\tilde{\mathfrak{u}}(\varepsilon_{1},\varepsilon_{2},\tilde{\boldsymbol{\Omega}};
\varsigma)}{\tilde{\boldsymbol{\theta}}(\varepsilon_{1} \tilde{\boldsymbol{u}}_{
+}(\varsigma) \! + \! \varepsilon_{2} \tilde{\boldsymbol{d}}) \tilde{\boldsymbol{
\theta}}(\varepsilon_{1} \tilde{\boldsymbol{u}}_{+}(\varsigma) \! - \! \frac{1}{2 
\pi}((n \! - \! 1)K \! + \! k) \tilde{\boldsymbol{\Omega}} \! + \! \varepsilon_{2} 
\tilde{\boldsymbol{d}})} \! + \! \dfrac{3(\tilde{\mathfrak{u}}(\varepsilon_{1},
\varepsilon_{2},\bm{0};\varsigma))^{2} \tilde{\mathfrak{v}}(\varepsilon_{1},
\varepsilon_{2},\bm{0};\varsigma)}{(\tilde{\boldsymbol{\theta}}(\varepsilon_{1} 
\tilde{\boldsymbol{u}}_{+}(\varsigma) \! + \! \varepsilon_{2} \tilde{\boldsymbol{
d}}))^{3}} \nonumber \\
+& \, \dfrac{\tilde{\mathfrak{u}}(\varepsilon_{1},\varepsilon_{2},\bm{0};\varsigma) 
\tilde{\mathfrak{w}}(\varepsilon_{1},\varepsilon_{2},\tilde{\boldsymbol{\Omega}};
\varsigma)}{\tilde{\boldsymbol{\theta}}(\varepsilon_{1} \tilde{\boldsymbol{u}}_{+}
(\varsigma) \! + \! \varepsilon_{2} \tilde{\boldsymbol{d}}) \tilde{\boldsymbol{
\theta}}(\varepsilon_{1} \tilde{\boldsymbol{u}}_{+}(\varsigma) \! - \! \frac{1}{2 
\pi}((n \! - \! 1)K \! + \! k) \tilde{\boldsymbol{\Omega}} \! + \! \varepsilon_{2} 
\tilde{\boldsymbol{d}})} \! + \! \left(\dfrac{\tilde{\mathfrak{u}}(\varepsilon_{1},
\varepsilon_{2},\bm{0};\varsigma)}{\tilde{\boldsymbol{\theta}}(\varepsilon_{1} 
\tilde{\boldsymbol{u}}_{+}(\varsigma) \! + \! \varepsilon_{2} \tilde{\boldsymbol{
d}})} \right)^{4} \nonumber \\
-& \, \dfrac{2 \tilde{\mathfrak{u}}(\varepsilon_{1},\varepsilon_{2},\bm{0};
\varsigma) \tilde{\mathfrak{v}}(\varepsilon_{1},\varepsilon_{2},\bm{0};\varsigma) 
\tilde{\mathfrak{u}}(\varepsilon_{1},\varepsilon_{2},\tilde{\boldsymbol{\Omega}};
\varsigma)}{(\tilde{\boldsymbol{\theta}}(\varepsilon_{1} \tilde{\boldsymbol{u}}_{+}
(\varsigma) \! + \! \varepsilon_{2} \tilde{\boldsymbol{d}}))^{2} \tilde{\boldsymbol{
\theta}}(\varepsilon_{1} \tilde{\boldsymbol{u}}_{+}(\varsigma) \! - \! \frac{1}{2 \pi}
((n \! - \! 1)K \! + \! k) \tilde{\boldsymbol{\Omega}} \! + \! \varepsilon_{2} 
\tilde{\boldsymbol{d}})} \! - \! \left(\dfrac{\tilde{\mathfrak{u}}(\varepsilon_{1},
\varepsilon_{2},\bm{0};\varsigma)}{\tilde{\boldsymbol{\theta}}(\varepsilon_{1} 
\tilde{\boldsymbol{u}}_{+}(\varsigma) \! + \! \varepsilon_{2} \tilde{\boldsymbol{
d}})} \right)^{2} \nonumber \\
\times& \, \dfrac{\tilde{\mathfrak{v}}(\varepsilon_{1},\varepsilon_{2},\tilde{
\boldsymbol{\Omega}};\varsigma)}{\tilde{\boldsymbol{\theta}}(\varepsilon_{1} 
\tilde{\boldsymbol{u}}_{+}(\varsigma) \! - \! \frac{1}{2 \pi}((n \! - \! 1)K \! + \! 
k) \tilde{\boldsymbol{\Omega}} \! + \! \varepsilon_{2} \tilde{\boldsymbol{d}})} 
\! - \! \dfrac{(\tilde{\mathfrak{u}}(\varepsilon_{1},\varepsilon_{2},\bm{0};
\varsigma))^{3} \tilde{\mathfrak{u}}(\varepsilon_{1},\varepsilon_{2},\tilde{
\boldsymbol{\Omega}};\varsigma)}{(\tilde{\boldsymbol{\theta}}(\varepsilon_{1} 
\tilde{\boldsymbol{u}}_{+}(\varsigma) \! + \! \varepsilon_{2} \tilde{\boldsymbol{
d}}))^{3} \tilde{\boldsymbol{\theta}}(\varepsilon_{1} \tilde{\boldsymbol{u}}_{+}
(\varsigma) \! - \! \frac{1}{2 \pi}((n \! - \! 1)K \! + \! k) \tilde{\boldsymbol{
\Omega}} \! + \! \varepsilon_{2} \tilde{\boldsymbol{d}})}, \label{eqprop70}
\end{align}
with
\begin{gather}
\tilde{\mathfrak{u}}(\varepsilon_{1},\varepsilon_{2},\tilde{\boldsymbol{\Omega}};
\varsigma) \! := \! 2 \pi \tilde{\Lambda}^{\raise-1.0ex\hbox{$\scriptstyle 1$}}_{0}
(\varepsilon_{1},\varepsilon_{2},\tilde{\boldsymbol{\Omega}};\varsigma), \qquad 
\quad \tilde{\mathfrak{v}}(\varepsilon_{1},\varepsilon_{2},\tilde{\boldsymbol{
\Omega}};\varsigma) \! := \! -2 \pi^{2} 
\tilde{\Lambda}^{\raise-1.0ex\hbox{$\scriptstyle 2$}}_{0}(\varepsilon_{1},
\varepsilon_{2},\tilde{\boldsymbol{\Omega}};\varsigma), \label{eqprop71} \\
\tilde{\mathfrak{w}}(\varepsilon_{1},\varepsilon_{2},\tilde{\boldsymbol{
\Omega}};\varsigma) \! := \! 2 \pi \left(
\tilde{\Lambda}^{\raise-1.0ex\hbox{$\scriptstyle 0$}}_{1}(\varepsilon_{1},
\varepsilon_{2},\tilde{\boldsymbol{\Omega}};\varsigma) \! - \! \dfrac{2 
\pi^{2}}{3} \tilde{\Lambda}^{\raise-1.0ex\hbox{$\scriptstyle 3$}}_{0}
(\varepsilon_{1},\varepsilon_{2},\tilde{\boldsymbol{\Omega}};\varsigma) 
\right), \label{eqprop72} \\
\tilde{\mathfrak{z}}(\varepsilon_{1},\varepsilon_{2},\tilde{\boldsymbol{
\Omega}};\varsigma) \! := \! -4 \pi^{2} \left(
\tilde{\Lambda}^{\raise-1.0ex\hbox{$\scriptstyle 1$}}_{1}(\varepsilon_{1},
\varepsilon_{2},\tilde{\boldsymbol{\Omega}};\varsigma) \! - \! \dfrac{\pi^{2}}{6} 
\tilde{\Lambda}^{\raise-1.0ex\hbox{$\scriptstyle 4$}}_{0}(\varepsilon_{1},
\varepsilon_{2},\tilde{\boldsymbol{\Omega}};\varsigma) \right), \label{eqprop73} \\
\tilde{\Lambda}^{\raise-1.0ex\hbox{$\scriptstyle j_{1}$}}_{j_{2}}(\varepsilon_{1},
\varepsilon_{2},\tilde{\boldsymbol{\Omega}};\varsigma) \! = \! \sum_{m \in
\mathbb{Z}^{N}}(\tilde{\mathfrak{r}}_{1}(\varsigma))^{j_{1}}(\tilde{\mathfrak{r}}_{2}
(\varsigma))^{j_{2}} \me^{2 \pi \mi (m,\varepsilon_{1} \tilde{\boldsymbol{u}}_{+}
(\varsigma)-\frac{1}{2 \pi}((n-1)K+k) \tilde{\boldsymbol{\Omega}}+ 
\varepsilon_{2} \tilde{\boldsymbol{d}})+ \mi \pi (m,\tilde{\boldsymbol{\tau}}m)}, 
\quad j_{1},j_{2} \! \in \! \mathbb{N}_{0}, \label{eqprop74} \\
\tilde{\mathfrak{r}}_{1}(\varsigma) \! := \! \dfrac{2}{\tilde{\leftthreetimes}
(\varsigma)} \sum_{i=1}^{N} \sum_{j=1}^{N}m_{i} \tilde{c}_{ij} \varsigma^{N-j}, 
\label{eqprop75} \\
\tilde{\mathfrak{r}}_{2}(\varsigma) \! := \! \dfrac{2}{3 \tilde{\leftthreetimes}
(\varsigma)} \sum_{i=1}^{N} \sum_{j=1}^{N}m_{i} \tilde{c}_{ij} \left(N \! - 
\! j \! - \! \dfrac{\varsigma \tilde{\leftthreetimes}^{\prime}(\varsigma)}{
\tilde{\leftthreetimes}(\varsigma)} \right) \varsigma^{N-j-1}, \label{eqprop76}
\end{gather}
where
\begin{align}
\tilde{\leftthreetimes}(\tilde{b}_{0}) :=& \, \mi (-1)^{N}(\tilde{a}_{N+1} \! - \! 
\tilde{b}_{0})^{1/2} \prod_{m=1}^{N}(\tilde{b}_{m} \! - \! \tilde{b}_{0})^{1/2}
(\tilde{a}_{m} \! - \! \tilde{b}_{0})^{1/2}, \label{eqprop77} \\
\tilde{\leftthreetimes}(\tilde{b}_{j}) :=& \, \mi (-1)^{N-j}(\tilde{b}_{j} \! - \! 
\tilde{b}_{0})^{1/2}(\tilde{a}_{N+1} \! - \! \tilde{b}_{j})^{1/2}(\tilde{b}_{j} \! - 
\! \tilde{a}_{j})^{1/2} \prod_{m=1}^{j-1}(\tilde{b}_{j} \! - \! \tilde{b}_{m})^{
1/2}(\tilde{b}_{j} \! - \! \tilde{a}_{m})^{1/2} \nonumber \\
\times& \, \prod_{m^{\prime}=j+1}^{N}(\tilde{b}_{m^{\prime}} \! - \! 
\tilde{b}_{j})^{1/2}(\tilde{a}_{m^{\prime}} \! - \! \tilde{b}_{j})^{1/2}, 
\quad j \! = \! 1,2,\dotsc,N, \label{eqprop78} \\
\dfrac{\tilde{\leftthreetimes}^{\prime}(\tilde{b}_{0})}{\tilde{\leftthreetimes}
(\tilde{b}_{0})} =& \, \dfrac{1}{2} \left(\sum_{m=1}^{N} \left(\dfrac{1}{\tilde{
b}_{0} \! - \! \tilde{b}_{m}} \! + \! \dfrac{1}{\tilde{b}_{0} \! - \! \tilde{a}_{m}} 
\right) \! + \! \dfrac{1}{\tilde{b}_{0} \! - \! \tilde{a}_{N+1}} \right), 
\label{eqprop79} \\
\dfrac{\tilde{\leftthreetimes}^{\prime}(\tilde{b}_{j})}{\tilde{\leftthreetimes}
(\tilde{b}_{j})} =& \, \dfrac{1}{2} \left(\sum_{\substack{m=1\\m \neq j}}^{N} 
\left(\dfrac{1}{\tilde{b}_{j} \! - \! \tilde{b}_{m}} \! + \! \dfrac{1}{\tilde{b}_{j} 
\! - \! \tilde{a}_{m}} \right) \! + \! \dfrac{1}{\tilde{b}_{j} \! - \! \tilde{a}_{j}} 
\! + \! \dfrac{1}{\tilde{b}_{j} \! - \! \tilde{b}_{0}} \! + \! \dfrac{1}{\tilde{b}_{j} 
\! - \! \tilde{a}_{N+1}} \right), \quad j \! = \! 1,2,\dotsc,N, \label{eqprop80}
\end{align}
and $\tilde{c}_{ij}$, $i,j \! = \! 1,2,\dotsc,N$, are described in 
Equations~\eqref{O1} and~\eqref{O2}, and
\begin{align}
\tilde{\mathfrak{Q}}_{0}(\tilde{b}_{0})=& \, -\mi (\tilde{a}_{N+1} \! - \! 
\tilde{b}_{0})^{-1/2} \prod_{m=1}^{N} \dfrac{(\tilde{b}_{m} \! - \! \tilde{b}_{
0})^{1/2}}{(\tilde{a}_{m} \! - \! \tilde{b}_{0})^{1/2}}, \label{eqprop81} \\
\dfrac{\tilde{\mathfrak{Q}}_{1}(\tilde{b}_{0})}{\tilde{\mathfrak{Q}}_{0}
(\tilde{b}_{0})}=& \, \dfrac{1}{2} \left(\sum_{m=1}^{N} \left(\dfrac{1}{\tilde{
b}_{0} \! - \! \tilde{b}_{m}} \! - \! \dfrac{1}{\tilde{b}_{0} \! - \! \tilde{a}_{m}} 
\right) \! - \! \dfrac{1}{\tilde{b}_{0} \! - \! \tilde{a}_{N+1}} \right), 
\label{eqprop82} \\
\dfrac{\tilde{\mathfrak{Q}}_{2}(\tilde{b}_{0})}{\tilde{\mathfrak{Q}}_{0}
(\tilde{b}_{0})}=& \, \dfrac{1}{2} \left(\sum_{m=1}^{N} \left(\dfrac{1}{
(\tilde{b}_{0} \! - \! \tilde{a}_{m})^{2}} \! - \! \dfrac{1}{(\tilde{b}_{0} \! - \! 
\tilde{b}_{m})^{2}} \right) \! + \! \dfrac{1}{(\tilde{b}_{0} \! - \! \tilde{a}_{N
+1})^{2}} \right) \! + \! \dfrac{1}{4} \left(\sum_{m=1}^{N} \left(\dfrac{1}{
\tilde{b}_{0} \! - \! \tilde{a}_{m}} \! - \! \dfrac{1}{\tilde{b}_{0} \! - \! 
\tilde{b}_{m}} \right) \! + \! \dfrac{1}{\tilde{b}_{0} \! - \! \tilde{a}_{N+1}} 
\right)^{2}, \label{eqprop83} \\
\tilde{\mathfrak{Q}}_{0}(\tilde{b}_{j})=& \, -\dfrac{\mi (\tilde{b}_{j} \! - \! 
\tilde{b}_{0})^{1/2}}{(\tilde{a}_{N+1} \! - \! \tilde{b}_{j})^{1/2}(\tilde{b}_{j} 
\! - \! \tilde{a}_{j})^{1/2}} \prod_{m=1}^{j-1} \dfrac{(\tilde{b}_{j} \! - \! 
\tilde{b}_{m})^{1/2}}{(\tilde{b}_{j} \! - \! \tilde{a}_{m})^{1/2}} \prod_{m^{
\prime}=j+1}^{N} \dfrac{(\tilde{b}_{m^{\prime}} \! - \! \tilde{b}_{j})^{1/2}}{
(\tilde{a}_{m^{\prime}} \! - \! \tilde{b}_{j})^{1/2}}, 
\quad j \! = \! 1,2,\dotsc,N, \label{eqprop84} \\
\dfrac{\tilde{\mathfrak{Q}}_{1}(\tilde{b}_{j})}{\tilde{\mathfrak{Q}}_{0}
(\tilde{b}_{j})}=& \, \dfrac{1}{2} \left(\sum_{\substack{m=1\\m \neq j}}^{N} 
\left(\dfrac{1}{\tilde{b}_{j} \! - \! \tilde{b}_{m}} \! - \! \dfrac{1}{\tilde{b}_{j} 
\! - \! \tilde{a}_{m}} \right) \! + \! \dfrac{1}{\tilde{b}_{j} \! - \! \tilde{b}_{0}} 
\! - \! \dfrac{1}{\tilde{b}_{j} \! - \! \tilde{a}_{N+1}} \! - \! \dfrac{1}{\tilde{
b}_{j} \! - \! \tilde{a}_{j}} \right), \quad j \! = \! 1,2,\dotsc,N, \label{eqprop85} \\
\dfrac{\tilde{\mathfrak{Q}}_{2}(\tilde{b}_{j})}{\tilde{\mathfrak{Q}}_{0}
(\tilde{b}_{j})}=& \, -\dfrac{1}{2} \left(\sum_{\substack{m=1\\m \neq j}}^{N} 
\left(\dfrac{1}{(\tilde{b}_{j} \! - \! \tilde{b}_{m})^{2}} \! - \! \dfrac{1}{(\tilde{
b}_{j} \! - \! \tilde{a}_{m})^{2}} \right) \! + \! \dfrac{1}{(\tilde{b}_{j} \! - \! 
\tilde{b}_{0})^{2}} \! - \! \dfrac{1}{(\tilde{b}_{j} \! - \! \tilde{a}_{N+1})^{2}} 
\! - \! \dfrac{1}{(\tilde{b}_{j} \! - \! \tilde{a}_{j})^{2}} \right) \nonumber \\
+& \, \dfrac{1}{4} \left(\sum_{\substack{m=1\\m \neq j}}^{N} \left(\dfrac{1}{
\tilde{b}_{j} \! - \! \tilde{b}_{m}} \! - \! \dfrac{1}{\tilde{b}_{j} \! - \! \tilde{
a}_{m}} \right) \! + \! \dfrac{1}{\tilde{b}_{j} \! - \! \tilde{b}_{0}} \! - \! 
\dfrac{1}{\tilde{b}_{j} \! - \! \tilde{a}_{N+1}} \! - \! \dfrac{1}{\tilde{b}_{j} 
\! - \! \tilde{a}_{j}} \right)^{2}, \quad j \! = \! 1,2,\dotsc,N, \label{eqprop86}
\end{align}
$\tilde{\alpha}_{0}(\tilde{b}_{0})$, $\tilde{\alpha}_{0}(\tilde{b}_{j})$, 
$\tilde{\alpha}_{1}(\tilde{b}_{0})$, and $\tilde{\alpha}_{1}
(\tilde{b}_{j})$, $j \! = \! 1,2,\dotsc,N$, are given in 
Equations~\eqref{eqmainfin49}--\eqref{eqmainfin52}, respectively, and 
$\tilde{\alpha}_{2}(\tilde{b}_{0}) \! = \! f^{\prime \prime}(\tilde{b}_{0})/7$ 
and $\tilde{\alpha}_{2}(\tilde{b}_{j}) \! = \! f^{\prime \prime}(\tilde{b}_{j})/7$, 
$j \! = \! 1,2,\dotsc,N$, where $f^{\prime \prime}(\tilde{b}_{0})$ and 
$f^{\prime \prime}(\tilde{b}_{j})$ are given in Equations~\eqref{wm24frbotil} 
and~\eqref{wm24frbjtil}, respectively, and $(\mathrm{M}_{2}(\mathbb{C}) \! \ni)$ 
$\tilde{\mathfrak{c}}^{r}_{m}(n,k,z_{o};\tilde{b}_{j-1}) \! =_{\underset{z_{o}
=1+o(1)}{\mathscr{N},n \to \infty}} \! \mathcal{O}(1)$, $m \! \in \! \mathbb{N}_{0}$, 
$r \! \in \! \lbrace \triangleright,\triangleleft \rbrace$, and, for $j \! = \! 1,2,
\dotsc,N \! + \! 1$,
\begin{gather}
\tilde{\boldsymbol{\mathrm{A}}}(\tilde{a}_{j}) \! := \! 
\begin{pmatrix}
\widetilde{\mathbb{K}}_{11} & \widetilde{\mathbb{K}}_{12} \\
\widetilde{\mathbb{K}}_{21} & \widetilde{\mathbb{K}}_{22}
\end{pmatrix} 
\begin{pmatrix}
\tilde{\mathbb{A}}_{11}(\tilde{a}_{j}) & \tilde{\mathbb{A}}_{12}
(\tilde{a}_{j}) \\
\tilde{\mathbb{A}}_{21}(\tilde{a}_{j}) & \tilde{\mathbb{A}}_{22}
(\tilde{a}_{j})
\end{pmatrix} 
\begin{pmatrix}
\widetilde{\mathbb{K}}_{22} & -\widetilde{\mathbb{K}}_{12} \\
-\widetilde{\mathbb{K}}_{21} & \widetilde{\mathbb{K}}_{11}
\end{pmatrix} \me^{\mi ((n-1)K+k) \tilde{\mho}_{j}}, \label{eqprop87} \\
\tilde{\boldsymbol{\mathrm{B}}}(\tilde{a}_{j}) \! := \! 
\begin{pmatrix}
\widetilde{\mathbb{K}}_{11} & \widetilde{\mathbb{K}}_{12} \\
\widetilde{\mathbb{K}}_{21} & \widetilde{\mathbb{K}}_{22}
\end{pmatrix} 
\begin{pmatrix}
\tilde{\mathbb{B}}_{11}(\tilde{a}_{j}) & \tilde{\mathbb{B}}_{12}
(\tilde{a}_{j}) \\
\tilde{\mathbb{B}}_{21}(\tilde{a}_{j}) & \tilde{\mathbb{B}}_{22}
(\tilde{a}_{j})
\end{pmatrix} 
\begin{pmatrix}
\widetilde{\mathbb{K}}_{22} & -\widetilde{\mathbb{K}}_{12} \\
-\widetilde{\mathbb{K}}_{21} & \widetilde{\mathbb{K}}_{11}
\end{pmatrix} \me^{\mi ((n-1)K+k) \tilde{\mho}_{j}}, \label{eqprop88} \\
\tilde{\boldsymbol{\mathrm{C}}}(\tilde{a}_{j}) \! := \! 
\begin{pmatrix}
\widetilde{\mathbb{K}}_{11} & \widetilde{\mathbb{K}}_{12} \\
\widetilde{\mathbb{K}}_{21} & \widetilde{\mathbb{K}}_{22}
\end{pmatrix} 
\begin{pmatrix}
\tilde{\mathbb{C}}_{11}(\tilde{a}_{j}) & \tilde{\mathbb{C}}_{12}
(\tilde{a}_{j}) \\
\tilde{\mathbb{C}}_{21}(\tilde{a}_{j}) & \tilde{\mathbb{C}}_{22}
(\tilde{a}_{j})
\end{pmatrix} 
\begin{pmatrix}
\widetilde{\mathbb{K}}_{22} & -\widetilde{\mathbb{K}}_{12} \\
-\widetilde{\mathbb{K}}_{21} & \widetilde{\mathbb{K}}_{11}
\end{pmatrix} \me^{\mi ((n-1)K+k) \tilde{\mho}_{j}}, \label{eqprop89}
\end{gather}
with
\begin{gather}
\tilde{\mathbb{A}}_{11}(\tilde{a}_{j}) \! = \! -\tilde{\mathbb{A}}_{22}
(\tilde{a}_{j}) \! = \! -s_{1} \tilde{\kappa}_{1}(\tilde{a}_{j}) \tilde{\kappa}_{2}
(\tilde{a}_{j}) \tilde{\mathfrak{Q}}_{0}(\tilde{a}_{j}), \label{eqprop90} \\
\tilde{\mathbb{A}}_{12}(\tilde{a}_{j}) \! = \! \mi s_{1}(\tilde{\kappa}_{1}
(\tilde{a}_{j}))^{2} \tilde{\mathfrak{Q}}_{0}(\tilde{a}_{j}), \qquad 
\tilde{\mathbb{A}}_{21}(\tilde{a}_{j}) \! = \! \mi s_{1}(\tilde{\kappa}_{2}
(\tilde{a}_{j}))^{2} \tilde{\mathfrak{Q}}_{0}(\tilde{a}_{j}), \label{eqprop91}
\end{gather}
\begin{align}
\tilde{\mathbb{B}}_{11}(\tilde{a}_{j}) \! = \! -\tilde{\mathbb{B}}_{22}
(\tilde{a}_{j}) =& \, \tilde{\kappa}_{1}(\tilde{a}_{j}) \tilde{\kappa}_{2}
(\tilde{a}_{j}) \left(-s_{1} \left(\tilde{\mathfrak{Q}}_{1}(\tilde{a}_{j}) \! + \! 
\tilde{\mathfrak{Q}}_{0}(\tilde{a}_{j}) \left(\tilde{\daleth}^{1}_{1}(\tilde{a}_{j}) 
\! + \! \tilde{\daleth}^{1}_{-1}(\tilde{a}_{j}) \right) \right) \right. \nonumber \\
-&\left. \, t_{1} \left((\tilde{\mathfrak{Q}}_{0}(\tilde{a}_{j}))^{-1} \! + \! 
\tilde{\mathfrak{Q}}_{0}(\tilde{a}_{j}) \tilde{\aleph}^{1}_{1}(\tilde{a}_{j}) 
\tilde{\aleph}^{1}_{-1}(\tilde{a}_{j}) \right) \! + \! \mi (s_{1} \! + \! t_{1}) 
\left(\tilde{\aleph}^{1}_{-1}(\tilde{a}_{j}) \! - \! \tilde{\aleph}^{1}_{1}
(\tilde{a}_{j}) \right) \right), \label{eqprop92} \\
\tilde{\mathbb{B}}_{12}(\tilde{a}_{j}) =& \, (\tilde{\kappa}_{1}(\tilde{a}_{j}))^{2} 
\left(\mi s_{1} \left(\tilde{\mathfrak{Q}}_{1}(\tilde{a}_{j}) \! + \! 2 \tilde{
\mathfrak{Q}}_{0}(\tilde{a}_{j}) \tilde{\daleth}^{1}_{1}(\tilde{a}_{j}) \right) 
\right. \nonumber \\
+&\left. \, \mi t_{1} \left(\tilde{\mathfrak{Q}}_{0}(\tilde{a}_{j})(\tilde{
\aleph}^{1}_{1}(\tilde{a}_{j}))^{2} \! - \! (\tilde{\mathfrak{Q}}_{0}
(\tilde{a}_{j}))^{-1} \right) \! - \! 2(s_{1} \! - \! t_{1}) \tilde{\aleph}^{1}_{1}
(\tilde{a}_{j}) \right), \label{eqprop93} \\
\tilde{\mathbb{B}}_{21}(\tilde{a}_{j}) =& \, (\tilde{\kappa}_{2}(\tilde{a}_{j}))^{2} 
\left(\mi s_{1} \left(\tilde{\mathfrak{Q}}_{1}(\tilde{a}_{j}) \! + \! 2 \tilde{
\mathfrak{Q}}_{0}(\tilde{a}_{j}) \tilde{\daleth}^{1}_{-1}(\tilde{a}_{j}) \right) 
\right. \nonumber \\
+&\left. \, \mi t_{1} \left(\tilde{\mathfrak{Q}}_{0}(\tilde{a}_{j})(\tilde{
\aleph}^{1}_{-1}(\tilde{a}_{j}))^{2} \! - \! (\tilde{\mathfrak{Q}}_{0}
(\tilde{a}_{j}))^{-1} \right) \! + \! 2(s_{1} \! - \! t_{1}) \tilde{\aleph}^{1}_{-1}
(\tilde{a}_{j}) \right), \label{eqprop94} \\
\tilde{\mathbb{C}}_{11}(\tilde{a}_{j}) \! = \! -\tilde{\mathbb{C}}_{22}
(\tilde{a}_{j}) =& \, \tilde{\kappa}_{1}(\tilde{a}_{j}) \tilde{\kappa}_{2}(\tilde{a}_{j}) 
\left(-s_{1} \left(\tilde{\mathfrak{Q}}_{0}(\tilde{a}_{j}) \left(\tilde{\gimel}^{1}_{1}
(\tilde{a}_{j}) \! + \! \tilde{\gimel}^{1}_{-1}(\tilde{a}_{j}) \! + \! \tilde{\daleth}^{
1}_{1}(\tilde{a}_{j}) \tilde{\daleth}^{1}_{-1}(\tilde{a}_{j}) \right) \right. \right. 
\nonumber \\
+&\left. \left. \, \tilde{\mathfrak{Q}}_{1}(\tilde{a}_{j}) \left(\tilde{\daleth}^{
1}_{1}(\tilde{a}_{j}) \! + \! \tilde{\daleth}^{1}_{-1}(\tilde{a}_{j}) \right) \! + \! 
\dfrac{1}{2} \tilde{\mathfrak{Q}}_{2}(\tilde{a}_{j}) \! + \! (\tilde{\mathfrak{
Q}}_{0}(\tilde{a}_{j}))^{-1} \tilde{\aleph}^{1}_{1}(\tilde{a}_{j}) \tilde{\aleph}^{
1}_{-1}(\tilde{a}_{j}) \right) \right. \nonumber \\
-&\left. t_{1} \left(\tilde{\mathfrak{Q}}_{0}(\tilde{a}_{j}) \left(\tilde{\aleph}^{
1}_{-1}(\tilde{a}_{j}) \tilde{\beth}^{1}_{1}(\tilde{a}_{j}) \! + \! \tilde{\aleph}^{
1}_{1}(\tilde{a}_{j}) \tilde{\beth}^{1}_{-1}(\tilde{a}_{j}) \right) \! + \! 
(\tilde{\mathfrak{Q}}_{0}(\tilde{a}_{j}))^{-1} \right. \right. \nonumber \\
\times&\left. \left. \, \left(\tilde{\daleth}^{1}_{1}(\tilde{a}_{j}) \! + \! 
\tilde{\daleth}^{1}_{-1}(\tilde{a}_{j}) \! - \! \tilde{\mathfrak{Q}}_{1}
(\tilde{a}_{j})(\tilde{\mathfrak{Q}}_{0}(\tilde{a}_{j}))^{-1} \right) \! + \! 
\tilde{\mathfrak{Q}}_{1}(\tilde{a}_{j}) \tilde{\aleph}^{1}_{1}(\tilde{a}_{j}) 
\tilde{\aleph}^{1}_{-1}(\tilde{a}_{j}) \right) \right. \nonumber \\
+&\left. \, \mi (s_{1} \! + \! t_{1}) \left(\tilde{\beth}^{1}_{-1}(\tilde{a}_{j}) 
\! - \! \tilde{\beth}^{1}_{1}(\tilde{a}_{j}) \! + \! \tilde{\aleph}^{1}_{-1}
(\tilde{a}_{j}) \tilde{\daleth}^{1}_{1}(\tilde{a}_{j}) \! - \! \tilde{\aleph}^{1}_{1}
(\tilde{a}_{j}) \tilde{\daleth}^{1}_{-1}(\tilde{a}_{j}) \right) \right), 
\label{eqprop95} \\
\tilde{\mathbb{C}}_{12}(\tilde{a}_{j}) =& \, (\tilde{\kappa}_{1}(\tilde{a}_{j}))^{2} 
\left(\mi s_{1} \left(\tilde{\mathfrak{Q}}_{0}(\tilde{a}_{j}) \left(2 \tilde{
\gimel}^{1}_{1}(\tilde{a}_{j}) \! + \! (\tilde{\daleth}^{1}_{1}(\tilde{a}_{j}))^{2} 
\right) \! + \! 2 \tilde{\mathfrak{Q}}_{1}(\tilde{a}_{j}) \tilde{\daleth}^{1}_{1}
(\tilde{a}_{j}) \right. \right. \nonumber \\
-&\left. \left. \, (\tilde{\mathfrak{Q}}_{0}(\tilde{a}_{j}))^{-1}(\tilde{
\aleph}^{1}_{1}(\tilde{a}_{j}))^{2} \! + \! \dfrac{1}{2} \tilde{\mathfrak{Q}}_{2}
(\tilde{a}_{j}) \right) \! + \! \mi t_{1} \left(2 \tilde{\mathfrak{Q}}_{0}
(\tilde{a}_{j}) \tilde{\aleph}^{1}_{1}(\tilde{a}_{j}) \tilde{\beth}^{1}_{1}
(\tilde{a}_{j}) \right. \right. \nonumber \\
+&\left. \left. \, \tilde{\mathfrak{Q}}_{1}(\tilde{a}_{j})(\tilde{\aleph}^{1}_{1}
(\tilde{a}_{j}))^{2} \! + \! (\tilde{\mathfrak{Q}}_{0}(\tilde{a}_{j}))^{-1} \left(
\tilde{\mathfrak{Q}}_{1}(\tilde{a}_{j})(\tilde{\mathfrak{Q}}_{0}(\tilde{a}_{j}))^{-1} 
\! - \! 2 \tilde{\daleth}^{1}_{1}(\tilde{a}_{j}) \right) \right) \right. \nonumber \\
-&\left. \, 2(s_{1} \! - \! t_{1}) \left(\tilde{\beth}^{1}_{1}(\tilde{a}_{j}) \! + \! 
\tilde{\aleph}^{1}_{1}(\tilde{a}_{j}) \tilde{\daleth}^{1}_{1}(\tilde{a}_{j}) \right) 
\right), \label{eqprop96} \\
\tilde{\mathbb{C}}_{21}(\tilde{a}_{j}) =& \, (\tilde{\kappa}_{2}(\tilde{a}_{j}))^{2} 
\left(\mi s_{1} \left(\tilde{\mathfrak{Q}}_{0}(\tilde{a}_{j}) \left(2 \tilde{
\gimel}^{1}_{-1}(\tilde{a}_{j}) \! + \! (\tilde{\daleth}^{1}_{-1}(\tilde{a}_{j}))^{2} 
\right) \! + \! 2 \tilde{\mathfrak{Q}}_{1}(\tilde{a}_{j}) \tilde{\daleth}^{1}_{-1}
(\tilde{a}_{j}) \right. \right. \nonumber \\
-&\left. \left. \, (\tilde{\mathfrak{Q}}_{0}(\tilde{a}_{j}))^{-1}(\tilde{\aleph}^{
1}_{-1}(\tilde{a}_{j}))^{2} \! + \! \dfrac{1}{2} \tilde{\mathfrak{Q}}_{2}
(\tilde{a}_{j}) \right) \! + \! \mi t_{1} \left(2 \tilde{\mathfrak{Q}}_{0}
(\tilde{a}_{j}) \tilde{\aleph}^{1}_{-1}(\tilde{a}_{j}) \tilde{\beth}^{1}_{-1}
(\tilde{a}_{j}) \right. \right. \nonumber \\
+&\left. \left. \, \tilde{\mathfrak{Q}}_{1}(\tilde{a}_{j})(\tilde{\aleph}^{1}_{-1}
(\tilde{a}_{j}))^{2} \! + \! (\tilde{\mathfrak{Q}}_{0}(\tilde{a}_{j}))^{-1} \left(
\tilde{\mathfrak{Q}}_{1}(\tilde{a}_{j})(\tilde{\mathfrak{Q}}_{0}(\tilde{a}_{j}))^{-1} 
\! - \! 2 \tilde{\daleth}^{1}_{-1}(\tilde{a}_{j}) \right) \right) \right. \nonumber \\
+&\left. \, 2(s_{1} \! - \! t_{1}) \left(\tilde{\beth}^{1}_{-1}(\tilde{a}_{j}) \! + \! 
\tilde{\aleph}^{1}_{-1}(\tilde{a}_{j}) \tilde{\daleth}^{1}_{-1}(\tilde{a}_{j}) 
\right) \right), \label{eqprop97} 
\end{align}
where, for $\varepsilon_{1},\varepsilon_{2} \! = \! \pm 1$, $\tilde{\kappa}_{1}
(\varsigma)$, $\tilde{\kappa}_{2}(\varsigma)$, $\tilde{\aleph}^{\varepsilon_{
1}}_{\varepsilon_{2}}(\varsigma)$, $\tilde{\daleth}^{\varepsilon_{1}}_{
\varepsilon_{2}}(\varsigma)$, $\tilde{\beth}^{\varepsilon_{1}}_{\varepsilon_{2}}
(\varsigma)$, and $\tilde{\gimel}^{\varepsilon_{1}}_{\varepsilon_{2}}
(\varsigma)$ are defined by Equations~\eqref{eqprop66}--\eqref{eqprop76}, 
and
\begin{align}
\tilde{\leftthreetimes}(\tilde{a}_{N+1}) :=& \, (\tilde{a}_{N+1} \! - \! 
\tilde{b}_{0})^{1/2} \prod_{m=1}^{N}(\tilde{a}_{N+1} \! - \! \tilde{b}_{m})^{1/2}
(\tilde{a}_{N+1} \! - \! \tilde{a}_{m})^{1/2}, \label{eqprop98} \\
\tilde{\leftthreetimes}(\tilde{a}_{j}) :=& \, (-1)^{N+1-j}(\tilde{a}_{j}  \! - \! 
\tilde{b}_{0})^{1/2}(\tilde{a}_{N+1} \! - \! \tilde{a}_{j})^{1/2}(\tilde{b}_{j} \! - 
\! \tilde{a}_{j})^{1/2} \prod_{m=1}^{j-1}(\tilde{a}_{j} \! - \! \tilde{b}_{m})^{1/2}
(\tilde{a}_{j} \! - \! \tilde{a}_{m})^{1/2} \nonumber \\
\times& \, \prod_{m^{\prime}=j+1}^{N}(\tilde{b}_{m^{\prime}} \! - \! 
\tilde{a}_{j})^{1/2}(\tilde{a}_{m^{\prime}} \! - \! \tilde{a}_{j})^{1/2}, 
\quad j \! = \! 1,2,\dotsc,N, 
\label{eqprop99} \\
\dfrac{\tilde{\leftthreetimes}^{\prime}(\tilde{a}_{N+1})}{\tilde{\leftthreetimes}
(\tilde{a}_{N+1})} =& \, \dfrac{1}{2} \left(\sum_{m=1}^{N} \left(\dfrac{1}{
\tilde{a}_{N+1} \! - \! \tilde{b}_{m}} \! + \! \dfrac{1}{\tilde{a}_{N+1} \! - \! 
\tilde{a}_{m}} \right) \! + \! \dfrac{1}{\tilde{a}_{N+1} \! - \! \tilde{b}_{0}} 
\right), \label{eqprop100} \\
\dfrac{\tilde{\leftthreetimes}^{\prime}(\tilde{a}_{j})}{\tilde{\leftthreetimes}
(\tilde{a}_{j})} =& \, \dfrac{1}{2} \left(\sum_{\substack{m=1\\m \neq j}}^{N} 
\left(\dfrac{1}{\tilde{a}_{j} \! - \! \tilde{b}_{m}} \! + \! \dfrac{1}{\tilde{a}_{j} 
\! - \! \tilde{a}_{m}} \right) \! + \! \dfrac{1}{\tilde{a}_{j} \! - \! \tilde{b}_{j}} \! 
+ \! \dfrac{1}{\tilde{a}_{j} \! - \! \tilde{b}_{0}} \! + \! \dfrac{1}{\tilde{a}_{j} \! 
- \! \tilde{a}_{N+1}} \right), \quad j \! = \! 1,2,\dotsc,N, \label{eqprop101} \\
\tilde{\mathfrak{Q}}_{0}(\tilde{a}_{N+1})=& \, (\tilde{a}_{N+1} \! - \! 
\tilde{b}_{0})^{1/2} \prod_{m=1}^{N} \dfrac{(\tilde{a}_{N+1} \! - \! \tilde{b}_{
m})^{1/2}}{(\tilde{a}_{N+1} \! - \! \tilde{a}_{m})^{1/2}}, \label{eqprop102} \\
\dfrac{\tilde{\mathfrak{Q}}_{1}(\tilde{a}_{N+1})}{\tilde{\mathfrak{Q}}_{0}
(\tilde{a}_{N+1})}=& \, \dfrac{1}{2} \left(\sum_{m=1}^{N} \left(\dfrac{1}{
\tilde{a}_{N+1} \! - \! \tilde{b}_{m}} \! - \! \dfrac{1}{\tilde{a}_{N+1} \! - \! 
\tilde{a}_{m}} \right) \! + \! \dfrac{1}{\tilde{a}_{N+1} \! - \! \tilde{b}_{0}} 
\right), \label{eqprop103} \\
\dfrac{\tilde{\mathfrak{Q}}_{2}(\tilde{a}_{N+1})}{\tilde{\mathfrak{Q}}_{0}
(\tilde{a}_{N+1})}=& \, -\dfrac{1}{2} \left(\sum_{m=1}^{N} \left(\dfrac{1}{
(\tilde{a}_{N+1} \! - \! \tilde{b}_{m})^{2}} \! - \! \dfrac{1}{(\tilde{a}_{N+1} 
\! - \! \tilde{a}_{m})^{2}} \right) \! + \! \dfrac{1}{(\tilde{a}_{N+1} \! - \! 
\tilde{b}_{0})^{2}} \right) \nonumber \\
+& \, \dfrac{1}{4} \left(\sum_{m=1}^{N} \left(\dfrac{1}{\tilde{a}_{N+1} 
\! - \! \tilde{b}_{m}} \! - \! \dfrac{1}{\tilde{a}_{N+1} \! - \! \tilde{a}_{m}} 
\right) \! + \! \dfrac{1}{\tilde{a}_{N+1} \! - \! \tilde{b}_{0}} \right)^{2}, 
\label{eqprop104} \\
\tilde{\mathfrak{Q}}_{0}(\tilde{a}_{j})=& \, \dfrac{(\tilde{a}_{j} \! - \! 
\tilde{b}_{0})^{1/2}(\tilde{b}_{j} \! - \! \tilde{a}_{j})^{1/2}}{(\tilde{a}_{N+1} 
\! - \! \tilde{a}_{j})^{1/2}} \prod_{m=1}^{j-1} \dfrac{(\tilde{a}_{j} \! - \! 
\tilde{b}_{m})^{1/2}}{(\tilde{a}_{j} \! - \! \tilde{a}_{m})^{1/2}} \prod_{m^{
\prime}=j+1}^{N} \dfrac{(\tilde{b}_{m^{\prime}} \! - \! \tilde{a}_{j})^{
1/2}}{(\tilde{a}_{m^{\prime}} \! - \! \tilde{a}_{j})^{1/2}}, 
\quad j \! = \! 1,2,\dotsc,N, \label{eqprop105} \\
\dfrac{\tilde{\mathfrak{Q}}_{1}(\tilde{a}_{j})}{\tilde{\mathfrak{Q}}_{0}
(\tilde{a}_{j})}=& \, \dfrac{1}{2} \left(\sum_{\substack{m=1\\m \neq j}}^{N} 
\left(\dfrac{1}{\tilde{a}_{j} \! - \! \tilde{b}_{m}} \! - \! \dfrac{1}{\tilde{a}_{j} 
\! - \! \tilde{a}_{m}} \right) \! + \! \dfrac{1}{\tilde{a}_{j} \! - \! \tilde{b}_{0}} 
\! - \! \dfrac{1}{\tilde{a}_{j} \! - \! \tilde{a}_{N+1}} \! + \! \dfrac{1}{\tilde{
a}_{j} \! - \! \tilde{b}_{j}} \right), \quad j \! = \! 1,2,\dotsc,N, 
\label{eqprop106} \\
\dfrac{\tilde{\mathfrak{Q}}_{2}(\tilde{a}_{j})}{\tilde{\mathfrak{Q}}_{0}
(\tilde{a}_{j})}=& \, -\dfrac{1}{2} \left(\sum_{\substack{m=1\\m \neq j}}^{N} 
\left(\dfrac{1}{(\tilde{a}_{j} \! - \! \tilde{b}_{m})^{2}} \! - \! \dfrac{1}{(\tilde{
a}_{j} \! - \! \tilde{a}_{m})^{2}} \right) \! + \! \dfrac{1}{(\tilde{a}_{j} \! - \! 
\tilde{b}_{0})^{2}} \! - \! \dfrac{1}{(\tilde{a}_{j} \! - \! \tilde{a}_{N+1})^{2}} 
\! + \! \dfrac{1}{(\tilde{a}_{j} \! - \! \tilde{b}_{j})^{2}} \right) \nonumber \\
+& \, \dfrac{1}{4} \left(\sum_{\substack{m=1\\m \neq j}}^{N} \left(
\dfrac{1}{\tilde{a}_{j} \! - \! \tilde{b}_{m}} \! - \! \dfrac{1}{\tilde{a}_{j} \! - \! 
\tilde{a}_{m}} \right) \! + \! \dfrac{1}{\tilde{a}_{j} \! - \! \tilde{b}_{0}} \! - \! 
\dfrac{1}{\tilde{a}_{j} \! - \! \tilde{a}_{N+1}} \! + \! \dfrac{1}{\tilde{a}_{j} \! 
- \! \tilde{b}_{j}} \right)^{2}, \quad j \! = \! 1,2,\dotsc,N, \label{eqprop107}
\end{align}
$\tilde{\alpha}_{0}(\tilde{a}_{N+1})$, $\tilde{\alpha}_{0}(\tilde{a}_{j})$, 
$\tilde{\alpha}_{1}(\tilde{a}_{N+1})$, and $\tilde{\alpha}_{1}
(\tilde{a}_{j})$, $j \! = \! 1,2,\dotsc,N$, are given in 
Equations~\eqref{eqmainfin53}--\eqref{eqmainfin56}, respectively, 
$\tilde{\alpha}_{2}(\tilde{a}_{N+1}) \! = \! f^{\prime \prime}(\tilde{a}_{N+1})/7$ 
and $\tilde{\alpha}_{2}(\tilde{a}_{j}) \! = \! f^{\prime \prime}(\tilde{a}_{j})/7$, 
$j \! = \! 1,2,\dotsc,N$, where $f^{\prime \prime}(\tilde{a}_{N+1})$ and 
$f^{\prime \prime}(\tilde{a}_{j})$ are given by Equations~\eqref{wm24frantil} 
and~\eqref{wm24frajtil}, respectively, and $(\mathrm{M}_{2}(\mathbb{C}) \! 
\ni)$ $\tilde{\mathfrak{c}}^{r}_{m}(n,k,z_{o};\tilde{a}_{j}) \! =_{\underset{z_{o}=
1+o(1)}{\mathscr{N},n \to \infty}} \! \mathcal{O}(1)$, $m \! \in \! \mathbb{N}_{0}$, 
$r \! \in \! \lbrace \triangleright,\triangleleft \rbrace$.
\end{bbbbb}

\emph{Proof.} The proof of this Proposition~\ref{propo5.1} consists of two 
cases: (i) $n \! \in \! \mathbb{N}$ and $k \! \in \! \lbrace 1,2,\dotsc,K 
\rbrace$ such that $\alpha_{p_{\mathfrak{s}}} \! := \! \alpha_{k} \! = \! 
\infty$; and (ii) $n \! \in \! \mathbb{N}$ and $k \! \in \! \lbrace 1,2,\dotsc,
K \rbrace$ such that $\alpha_{p_{\mathfrak{s}}} \! := \! \alpha_{k} \! \neq 
\! \infty$. Notwithstanding the fact that the scheme of the proof is, 
\emph{mutatis mutandis}, similar for both cases, case~(ii), nevertheless, 
is the more technically challenging of the two; therefore, without loss of 
generality, only the proof for case~(ii) is presented in detail, whilst case~(i) 
is proved analogously.

For $n \! \in \! \mathbb{N}$ and $k \! \in \! \lbrace 1,2,\dotsc,K \rbrace$ 
such that $\alpha_{p_{\mathfrak{s}}} \! := \! \alpha_{k} \! \neq \! \infty$, 
and for $z \! \in \! \partial \tilde{\mathbb{U}}_{\tilde{\delta}_{\tilde{b}_{j-
1}}} \cup \partial \tilde{\mathbb{U}}_{\tilde{\delta}_{\tilde{a}_{j}}}$, $j \! = 
\! 1,2,\dotsc,N \! + \! 1$, consider, say, and without loss of generality, 
the case $z \! \in \! \partial \tilde{\mathbb{U}}_{\tilde{\delta}_{\tilde{a}_{j}}}$, 
$j \! = \! 1,2,\dotsc,N \! + \! 1$; the analysis for the case $z \! \in \! 
\partial \tilde{\mathbb{U}}_{\tilde{\delta}_{\tilde{b}_{j-1}}}$, $j \! = \! 1,2,
\dotsc,N \! + \! 1$, is analogous. For $n \! \in \! \mathbb{N}$ and $k \! 
\in \! \lbrace 1,2,\dotsc,K \rbrace$ such that $\alpha_{p_{\mathfrak{s}}} 
\! := \! \alpha_{k} \! \neq \! \infty$, via the asymptotic expansion, in the 
double-scaling limit $\mathscr{N},n \! \to \! \infty$ such that $z_{o} 
\! = \! 1 \! + \! o(1)$, for $\tilde{v}_{\tilde{\mathcal{R}}}(z)$ given in 
Equation~\eqref{eqtlvee14}, the expression $\tilde{\xi}_{\tilde{a}_{j}}(z) 
\! = \! (z \! - \! \tilde{a}_{j})^{3/2} \tilde{G}_{\tilde{a}_{j}}(z)$, $z \! \in 
\! \tilde{\mathbb{U}}_{\tilde{\delta}_{\tilde{a}_{j}}} \setminus (-\infty,
\tilde{a}_{j})$, $j \! = \! 1,2,\dotsc,N \! + \! 1$, given in Lemma~\ref{lem4.9}, 
and the formula $w^{\Sigma_{\tilde{\mathcal{R}}}}_{+}(z) \! = \! \tilde{v}_{
\tilde{\mathcal{R}}}(z) \! - \! \mathrm{I}$, it follows that, for $z \! \in \! 
\mathbb{C}_{\pm} \cap \partial \tilde{\mathbb{U}}_{\tilde{\delta}_{\tilde{a}_{j}}}$, 
$j \! = \! 1,2,\dotsc,N \! + \! 1$,
\begin{align} \label{eqlempro1} 
\widetilde{\mathbb{K}}^{-1} w^{\Sigma_{\tilde{\mathcal{R}}}}_{+}(z) 
\widetilde{\mathbb{K}} \underset{\underset{z_{o}=1+o(1)}{\mathscr{N},
n \to \infty}}{=}& \, \dfrac{1}{((n \! - \! 1)K \! + \! k)(z \! - \! 
\tilde{a}_{j})^{3/2} \tilde{G}_{\tilde{a}_{j}}(z)} \tilde{\mathbb{M}}(z) 
\begin{pmatrix}
\mp (s_{1} \! + \! t_{1}) & \pm \mi (s_{1} \! - \! t_{1}) \me^{\mi ((n-1)K
+k) \tilde{\mho}_{j}} \\
\pm \mi (s_{1} \! - \! t_{1}) \me^{-\mi ((n-1)K+k) \tilde{\mho}_{j}} & 
\pm (s_{1} \! + \! t_{1})
\end{pmatrix}(\tilde{\mathbb{M}}(z))^{-1} \nonumber \\
+& \, \mathcal{O} \left(\dfrac{1}{((n \! - \! 1)K \! + \! k)^{2}(z \! - 
\! \tilde{a}_{j})^{3}(\tilde{G}_{\tilde{a}_{j}}(z))^{2}} \tilde{\mathbb{M}}(z) 
\tilde{\mathfrak{c}}^{\triangleleft}(n,k,z_{o};j)(\tilde{\mathbb{M}}(z))^{-1} 
\right),
\end{align}
where $\widetilde{\mathbb{K}}$ and $\tilde{\mathbb{M}}(z)$ are given 
in item~\pmb{(2)} of Lemma~\ref{lem4.5}, and $\tilde{\mho}_{m}$, $m 
\! = \! 0,1,\dotsc,N \! + \! 1$, is defined in the corresponding item~(ii) 
of Remark~\ref{rem4.4}. For $n \! \in \! \mathbb{N}$ and $k \! \in \! 
\lbrace 1,2,\dotsc,K \rbrace$ such that $\alpha_{p_{\mathfrak{s}}} \! 
:= \! \alpha_{k} \! \neq \! \infty$, a matrix-multiplication argument 
shows that, for $z \! \in \! \mathbb{C}_{\pm} \cap \partial \tilde{
\mathbb{U}}_{\tilde{\delta}_{\tilde{a}_{j}}}$, $j \! = \! 1,2,\dotsc,N \! 
+ \! 1$,
\begin{equation} \label{eqlempro2} 
\tilde{\mathbb{M}}(z) 
\begin{pmatrix}
\mp (s_{1} \! + \! t_{1}) & \pm \mi (s_{1} \! - \! t_{1}) \me^{\mi ((n-1)K
+k) \tilde{\mho}_{j}} \\
\pm \mi (s_{1} \! - \! t_{1}) \me^{-\mi ((n-1)K+k) \tilde{\mho}_{j}} & 
\pm (s_{1} \! + \! t_{1})
\end{pmatrix}(\tilde{\mathbb{M}}(z))^{-1} \! = \! 
\begin{pmatrix}
\tilde{\mathfrak{t}}_{11}(z) & \tilde{\mathfrak{t}}_{12}(z) \\
\tilde{\mathfrak{t}}_{21}(z) & \tilde{\mathfrak{t}}_{22}(z)
\end{pmatrix},
\end{equation} 
where
\begin{align}
\tilde{\mathfrak{t}}_{11}(z):=& \, \mp \dfrac{1}{4}(s_{1} \! + \! t_{1}) \left(
\dfrac{(\tilde{\gamma}(z))^{2} \! + \! 1}{\tilde{\gamma}(z)} \right)^{2} 
\tilde{\mathfrak{m}}_{11}(z) \tilde{\mathfrak{m}}_{22}(z) \! \mp \! 
\dfrac{1}{4}(s_{1} \! - \! t_{1}) \left(\dfrac{(\tilde{\gamma}(z))^{4} 
\! - \! 1}{(\tilde{\gamma}(z))^{2}} \right) \tilde{\mathfrak{m}}_{11}(z) 
\tilde{\mathfrak{m}}_{21}(z) \me^{\mi ((n-1)K+k) \tilde{\mho}_{j}} 
\nonumber \\
\mp& \, \dfrac{1}{4}(s_{1} \! + \! t_{1}) \left(\dfrac{(\tilde{\gamma}(z))^{2} 
\! - \! 1}{\tilde{\gamma}(z)} \right)^{2} \tilde{\mathfrak{m}}_{12}(z) 
\tilde{\mathfrak{m}}_{21}(z) \! \mp \! \dfrac{1}{4}(s_{1} \! - \! t_{1}) \left(
\dfrac{(\tilde{\gamma}(z))^{4} \! - \! 1}{(\tilde{\gamma}(z))^{2}} \right) 
\tilde{\mathfrak{m}}_{12}(z) \tilde{\mathfrak{m}}_{22}(z) 
\me^{-\mi ((n-1)K+k) \tilde{\mho}_{j}}, \label{eqlempro3} \\
\tilde{\mathfrak{t}}_{12}(z):=& \, \pm \dfrac{\mi}{2}(s_{1} \! + \! t_{1}) 
\left(\dfrac{(\tilde{\gamma}(z))^{4} \! - \! 1}{(\tilde{\gamma}(z))^{2}} 
\right) \tilde{\mathfrak{m}}_{11}(z) \tilde{\mathfrak{m}}_{12}(z) \! \pm 
\! \dfrac{\mi}{4}(s_{1} \! - \! t_{1}) \left(\dfrac{(\tilde{\gamma}(z))^{2} 
\! + \! 1}{\tilde{\gamma}(z)} \right)^{2}(\tilde{\mathfrak{m}}_{11}(z))^{2} 
\me^{\mi ((n-1)K+k) \tilde{\mho}_{j}} \nonumber \\
\pm& \, \dfrac{\mi}{4}(s_{1} \! - \! t_{1}) \left(\dfrac{(\tilde{\gamma}(z))^{2} 
\! - \! 1}{\tilde{\gamma}(z)} \right)^{2}(\tilde{\mathfrak{m}}_{12}(z))^{2} 
\me^{-\mi ((n-1)K+k) \tilde{\mho}_{j}}, \label{eqlempro4} \\
\tilde{\mathfrak{t}}_{21}(z):=& \, \pm \dfrac{\mi}{2}(s_{1} \! + \! t_{1}) 
\left(\dfrac{(\tilde{\gamma}(z))^{4} \! - \! 1}{(\tilde{\gamma}(z))^{2}} 
\right) \tilde{\mathfrak{m}}_{21}(z) \tilde{\mathfrak{m}}_{22}(z) \! \pm 
\! \dfrac{\mi}{4}(s_{1} \! - \! t_{1}) \left(\dfrac{(\tilde{\gamma}(z))^{2} 
\! - \! 1}{\tilde{\gamma}(z)} \right)^{2}(\tilde{\mathfrak{m}}_{21}(z))^{2} 
\me^{\mi ((n-1)K+k) \tilde{\mho}_{j}} \nonumber \\
\pm& \, \dfrac{\mi}{4}(s_{1} \! - \! t_{1}) \left(\dfrac{(\tilde{\gamma}(z))^{2} 
\! + \! 1}{\tilde{\gamma}(z)} \right)^{2}(\tilde{\mathfrak{m}}_{22}(z))^{2} 
\me^{-\mi ((n-1)K+k) \tilde{\mho}_{j}}, \label{eqlempro5} \\
\tilde{\mathfrak{t}}_{22}(z):=& \, \pm \dfrac{1}{4}(s_{1} \! + \! t_{1}) \left(
\dfrac{(\tilde{\gamma}(z))^{2} \! + \! 1}{\tilde{\gamma}(z)} \right)^{2} 
\tilde{\mathfrak{m}}_{11}(z) \tilde{\mathfrak{m}}_{22}(z) \! \pm \! 
\dfrac{1}{4}(s_{1} \! - \! t_{1}) \left(\dfrac{(\tilde{\gamma}(z))^{4} \! 
- \! 1}{(\tilde{\gamma}(z))^{2}} \right) \tilde{\mathfrak{m}}_{11}(z) 
\tilde{\mathfrak{m}}_{21}(z) \me^{\mi ((n-1)K+k) \tilde{\mho}_{j}} 
\nonumber \\
\pm& \, \dfrac{1}{4}(s_{1} \! + \! t_{1}) \left(\dfrac{(\tilde{\gamma}(z))^{2} 
\! - \! 1}{\tilde{\gamma}(z)} \right)^{2} \tilde{\mathfrak{m}}_{12}(z) 
\tilde{\mathfrak{m}}_{21}(z) \! \pm \! \dfrac{1}{4}(s_{1} \! - \! t_{1}) 
\left(\dfrac{(\tilde{\gamma}(z))^{4} \! - \! 1}{(\tilde{\gamma}(z))^{2}} 
\right) \tilde{\mathfrak{m}}_{12}(z) \tilde{\mathfrak{m}}_{22}(z) 
\me^{-\mi ((n-1)K+k) \tilde{\mho}_{j}}, \label{eqlempro6}
\end{align}
with $\tilde{\gamma}(z)$ and $\tilde{\mathfrak{m}}(z)$ defined by 
Equations~\eqref{eqmainfin10} and~\eqref{eqtilfrakm}, respectively. 
Recall {}from the associated Subsection~\ref{subsub1} that, for $n \! 
\in \! \mathbb{N}$ and $k \! \in \! \lbrace 1,2,\dotsc,K \rbrace$ such 
that $\alpha_{p_{\mathfrak{s}}} \! := \! \alpha_{k} \! \neq \! \infty$, 
$\tilde{\boldsymbol{\omega}} \! := \! (\tilde{\omega}_{1},\tilde{
\omega}_{2},\dotsc,\tilde{\omega}_{N})$, where $\tilde{\omega}_{i} \! 
= \! \sum_{j^{\prime}=1}^{N} \tilde{c}_{ij^{\prime}}(\tilde{R}(z))^{-1/2}
z^{N-j^{\prime}} \linebreak[4]
\pmb{\cdot} \, \md z$, $i \! = \! 1,2,\dotsc,N$, with $c_{ij^{\prime}}$, 
$i,j^{\prime} \! = \! 1,2,\dotsc,N$, described by Equations~\eqref{O1} 
and~\eqref{O2}, and the multi-valued function $(\tilde{R}(z))^{1/2}$ 
is defined by Equation~\eqref{eql3.7j}: a calculation shows that
\begin{equation} \label{eqlempro7} 
\tilde{\omega}_{i} \underset{\underset{j=1,2,\dotsc,N+1}{z \to \tilde{a}_{
j}}}{=} \dfrac{(\tilde{\leftthreetimes}(\tilde{a}_{j}))^{-1}}{(z \! - \! \tilde{a}_{
j})^{1/2}} \left(\tilde{\mathfrak{p}}_{i}(\tilde{a}_{j}) \! + \! \tilde{\mathfrak{
q}}_{i}(\tilde{a}_{j})(z \! - \! \tilde{a}_{j}) \! + \! \tilde{\mathfrak{u}}_{i}
(\tilde{a}_{j})(z \! - \! \tilde{a}_{j})^{2} \! + \! \mathcal{O}((z \! - \! 
\tilde{a}_{j})^{3}) \right) \md z, \quad i \! = \! 1,2,\dotsc,N,
\end{equation}
where
\begin{equation*}
\tilde{\mathfrak{p}}_{i}(\varsigma) \! = \! \sum_{j^{\prime}=1}^{N} 
\tilde{c}_{ij^{\prime}} \varsigma^{N-j^{\prime}}, \qquad \qquad 
\tilde{\mathfrak{q}}_{i}(\varsigma) \! = \! \sum_{j^{\prime}=1}^{N} 
\tilde{c}_{ij^{\prime}} \left(N \! - \! j^{\prime} \! - \! \dfrac{\varsigma 
\tilde{\leftthreetimes}^{\prime}(\varsigma)}{\tilde{\leftthreetimes}
(\varsigma)} \right) \varsigma^{N-j^{\prime}-1},
\end{equation*}
with $\tilde{\leftthreetimes}(\tilde{a}_{j}),\tilde{\leftthreetimes}^{\prime}
(\tilde{a}_{j})$ given by Equations~\eqref{eqprop98}--\eqref{eqprop101}, 
and
\begin{equation*}
\tilde{\mathfrak{u}}_{i}(\varsigma) \! = \! \sum_{m=1}^{N} \tilde{c}_{im} 
\left(\dfrac{1}{2}(N \! - \! m)(N \! - \! m \! - \! 1) \! - \! (N \! - \! m) 
\dfrac{\varsigma \tilde{\leftthreetimes}^{\prime}(\varsigma)}{\tilde{
\leftthreetimes}(\varsigma)} \! + \! \left(\left(\dfrac{\tilde{
\leftthreetimes}^{\prime}(\varsigma)}{\tilde{\leftthreetimes}(\varsigma)} 
\right)^{2} \! - \! \dfrac{\tilde{\leftthreetimes}^{\prime \prime}(\varsigma)}{2 
\tilde{\leftthreetimes}(\varsigma)} \right) \varsigma^{2} \right) 
\varsigma^{N-m-2}, \quad i \! = \! 1,2,\dotsc,N,
\end{equation*}
with
\begin{align*}
\dfrac{\tilde{\leftthreetimes}^{\prime \prime}(\tilde{a}_{N+1})}{
\tilde{\leftthreetimes}(\tilde{a}_{N+1})}=& \, -\dfrac{1}{2} \left(
\sum_{m=1}^{N} \left(\dfrac{1}{(\tilde{a}_{N+1} \! - \! \tilde{b}_{m})^{2}} 
\! + \! \dfrac{1}{(\tilde{a}_{N+1} \! - \! \tilde{a}_{m})^{2}} \right) \! 
+ \! \dfrac{1}{(\tilde{a}_{N+1} \! - \! \tilde{b}_{0})^{2}} \right) \\
+& \, \dfrac{1}{4} \left(\sum_{m=1}^{N} \left(\dfrac{1}{\tilde{a}_{N+1} 
\! - \! \tilde{b}_{m}} \! + \! \dfrac{1}{\tilde{a}_{N+1} \! - \! \tilde{a}_{m}} 
\right) \! + \! \dfrac{1}{\tilde{a}_{N+1} \! - \! \tilde{b}_{0}} \right)^{2}, \\
\dfrac{\tilde{\leftthreetimes}^{\prime \prime}(\tilde{a}_{j^{\prime}})}{
\tilde{\leftthreetimes}(\tilde{a}_{j^{\prime}})}=& \, -\dfrac{1}{2} \left(
\sum_{\substack{m=1\\m \neq j^{\prime}}}^{N} \left(\dfrac{1}{(\tilde{
a}_{j^{\prime}} \! - \! \tilde{b}_{m})^{2}} \! + \! \dfrac{1}{(\tilde{a}_{j^{
\prime}} \! - \! \tilde{a}_{m})^{2}} \right) \! + \! \dfrac{1}{(\tilde{a}_{
j^{\prime}} \! - \! \tilde{b}_{j^{\prime}})^{2}} \! + \! \dfrac{1}{(\tilde{
a}_{j^{\prime}} \! - \! \tilde{a}_{N+1})^{2}} \! + \! \dfrac{1}{(\tilde{a}_{
j^{\prime}} \! - \! \tilde{b}_{0})^{2}} \right) \\
+& \, \dfrac{1}{4} \left(\sum_{\substack{m=1\\m \neq j^{\prime}}}^{N} 
\left(\dfrac{1}{\tilde{a}_{j^{\prime}} \! - \! \tilde{b}_{m}} \! + \! \dfrac{
1}{\tilde{a}_{j^{\prime}} \! - \! \tilde{a}_{m}} \right) \! + \! \dfrac{1}{
\tilde{a}_{j^{\prime}} \! - \! \tilde{b}_{j^{\prime}}} \! + \! \dfrac{1}{
\tilde{a}_{j^{\prime}} \! - \! \tilde{a}_{N+1}} \! + \! \dfrac{1}{\tilde{a}_{
j^{\prime}} \! - \! \tilde{b}_{0}} \right)^{2}, \quad j^{\prime} \! = \! 1,
2,\dotsc,N.
\end{align*}
An integration argument reveals that, for $n \! \in \! \mathbb{N}$ 
and $k \! \in \! \lbrace 1,2,\dotsc,K \rbrace$ such that 
$\alpha_{p_{\mathfrak{s}}} \! := \! \alpha_{k} \! \neq \! \infty$,
\begin{equation} \label{eqlempro8} 
\int_{\tilde{a}_{j}}^{z} \tilde{\omega}_{i} \underset{\underset{j=1,2,\dotsc,
N+1}{z \to \tilde{a}_{j}}}{=} \dfrac{2 \tilde{\mathfrak{p}}_{i}(\tilde{a}_{j})}{
\tilde{\leftthreetimes}(\tilde{a}_{j})}(z \! - \! \tilde{a}_{j})^{1/2} \! + \! 
\dfrac{2 \tilde{\mathfrak{q}}_{i}(\tilde{a}_{j})}{3 \tilde{\leftthreetimes}
(\tilde{a}_{j})}(z \! - \! a_{j})^{3/2} \! + \! \dfrac{2 \tilde{\mathfrak{u}}_{i}
(\tilde{a}_{j})}{5 \tilde{\leftthreetimes}(\tilde{a}_{j})}(z \! - \! \tilde{a}_{j})^{
5/2} \! + \! \mathcal{O}((z \! - \! \tilde{a}_{j})^{7/2}), \quad i \! = \! 1,2,
\dotsc,N.
\end{equation}
Recall {}from the proof of Lemma~\ref{lem4.5} that, for $n \! \in \! 
\mathbb{N}$ and $k \! \in \! \lbrace 1,2,\dotsc,K \rbrace$ such that 
$\alpha_{p_{\mathfrak{s}}} \! := \! \alpha_{k} \! \neq \! \infty$, the 
matrix $\tilde{\mathfrak{m}}(z)$ defined by Equation~\eqref{eqtilfrakm} 
satisfies the jump relation $\tilde{\mathfrak{m}}_{+}(z) \! = \! \tilde{
\mathfrak{m}}_{-}(z)(\exp (-\mi ((n \! - \! 1)K \! + \! k) \tilde{\Omega}_{j^{
\prime}}) \sigma_{-} \! + \! \exp (\mi ((n \! - \! 1)K \! + \! k) \tilde{
\Omega}_{j^{\prime}}) \sigma_{+})$ for $z \! \in \! (\tilde{a}_{j^{\prime}},
\tilde{b}_{j^{\prime}}) \setminus \lbrace \tilde{z}_{j^{\prime}}^{\pm} 
\rbrace$, $j^{\prime} \! = \! 1,2,\dotsc,N$, and $\tilde{\mathfrak{m}}_{+}
(z) \! = \! \tilde{\mathfrak{m}}_{-}(z) 
\left(
\begin{smallmatrix}
0 & 1 \\
1 & 0
\end{smallmatrix}
\right)$ for $z \! \in \! (-\infty,\tilde{b}_{0}) \cup (\tilde{a}_{N+1},+\infty)$: 
via this latter observation, the relation $\tilde{\boldsymbol{u}}(z) \! := 
\! \int_{\tilde{a}_{N+1}}^{z} \tilde{\boldsymbol{\omega}} \! = \! \tilde{
\boldsymbol{u}}(\tilde{a}_{j^{\prime}}) \! + \! (\int_{\tilde{a}_{j^{\prime}}}^{z} 
\tilde{\omega}_{1},\int_{\tilde{a}_{j^{\prime}}}^{z} \tilde{\omega}_{2},\dotsc,
\int_{\tilde{a}_{j^{\prime}}}^{z} \tilde{\omega}_{N})$, $j^{\prime} \! = \! 
1,2,\dotsc,N \! + \! 1$, with $\tilde{\boldsymbol{u}}(\tilde{a}_{N+1}) \! 
\equiv \! (0,0,\dotsc,0)$ $(\in \operatorname{Jac}(\tilde{\mathcal{Y}}))$, 
and the representation for $\tilde{\boldsymbol{\theta}}(z)$ 
given by Equation~\eqref{eqrmthetafin}, it follows {}from the 
Expansion~\eqref{eqlempro8}, and a careful analysis of the branch cuts, 
that, for $n \! \in \! \mathbb{N}$ and $k \! \in \! \lbrace 1,2,\dotsc,K 
\rbrace$ such that $\alpha_{p_{\mathfrak{s}}} \! := \! \alpha_{k} \! \neq 
\! \infty$,
\begin{align}
\tilde{\mathfrak{m}}_{11}(z) \underset{\underset{j=1,2,\dotsc,N+1}{z \to 
\tilde{a}_{j}}}{=}& \, \tilde{\kappa}_{1}(\tilde{a}_{j}) \left(1 \! + \! \mi 
\tilde{\aleph}^{1}_{1}(\tilde{a}_{j})(z \! - \! \tilde{a}_{j})^{1/2} \! + \! 
\tilde{\daleth}^{1}_{1}(\tilde{a}_{j})(z \! - \! \tilde{a}_{j}) \! + \! \mi 
\tilde{\beth}^{1}_{1}(\tilde{a}_{j})(z \! - \! \tilde{a}_{j})^{3/2} \! + \! 
\tilde{\gimel}^{1}_{1}(\tilde{a}_{j})(z \! - \! \tilde{a}_{j})^{2} \! + \! 
\mathcal{O}((z \! - \! \tilde{a}_{j})^{5/2}) \right), \label{eqlempro9} \\
\tilde{\mathfrak{m}}_{12}(z) \underset{\underset{j=1,2,\dotsc,N+1}{z \to 
\tilde{a}_{j}}}{=}& \, \tilde{\kappa}_{1}(\tilde{a}_{j}) \left(1 \! - \! \mi 
\tilde{\aleph}^{-1}_{1}(\tilde{a}_{j})(z \! - \! \tilde{a}_{j})^{1/2} \! + \! 
\tilde{\daleth}^{-1}_{1}(\tilde{a}_{j})(z \! - \! \tilde{a}_{j}) \! - \! \mi 
\tilde{\beth}^{-1}_{1}(\tilde{a}_{j})(z \! - \! \tilde{a}_{j})^{3/2} \! + \! 
\tilde{\gimel}^{-1}_{1}(\tilde{a}_{j})(z \! - \! \tilde{a}_{j})^{2} \! + \! 
\mathcal{O}((z \! - \! \tilde{a}_{j})^{5/2}) \right) \nonumber \\
\times& \, \exp \left(\mi ((n \! - \! 1)K \! + \! k) \tilde{\mho}_{j} \right), 
\label{eqlempro10} \\
\tilde{\mathfrak{m}}_{21}(z) \underset{\underset{j=1,2,\dotsc,N+1}{z \to 
\tilde{a}_{j}}}{=}& \, \tilde{\kappa}_{2}(\tilde{a}_{j}) \left(1 \! + \! \mi 
\tilde{\aleph}^{1}_{-1}(\tilde{a}_{j})(z \! - \! \tilde{a}_{j})^{1/2} \! + \! 
\tilde{\daleth}^{1}_{-1}(\tilde{a}_{j})(z \! - \! \tilde{a}_{j}) \! + \! \mi 
\tilde{\beth}^{1}_{-1}(\tilde{a}_{j})(z \! - \! \tilde{a}_{j})^{3/2} \! + \! 
\tilde{\gimel}^{1}_{-1}(\tilde{a}_{j})(z \! - \! \tilde{a}_{j})^{2} \! + \! 
\mathcal{O}((z \! - \! \tilde{a}_{j})^{5/2}) \right), \label{eqlempro11} \\
\tilde{\mathfrak{m}}_{22}(z) \underset{\underset{j=1,2,\dotsc,N+1}{z \to 
\tilde{a}_{j}}}{=}& \, \tilde{\kappa}_{2}(\tilde{a}_{j}) \left(1 \! - \! \mi 
\tilde{\aleph}^{-1}_{-1}(\tilde{a}_{j})(z \! - \! \tilde{a}_{j})^{1/2} \! + \! 
\tilde{\daleth}^{-1}_{-1}(\tilde{a}_{j})(z \! - \! \tilde{a}_{j}) \! - \! \mi 
\tilde{\beth}^{-1}_{-1}(\tilde{a}_{j})(z \! - \! \tilde{a}_{j})^{3/2} \! + \! 
\tilde{\gimel}^{-1}_{-1}(\tilde{a}_{j})(z \! - \! \tilde{a}_{j})^{2} \! + \! 
\mathcal{O}((z \! - \! \tilde{a}_{j})^{5/2}) \right) \nonumber \\
\times& \, \exp \left(\mi ((n \! - \! 1)K \! + \! K) \tilde{\mho}_{j} \right), 
\label{eqlempro12}
\end{align}
where, for $\varepsilon_{1},\varepsilon_{2} \! = \! \pm 1$, $\tilde{\kappa}_{
1}(\varsigma)$, $\tilde{\kappa}_{2}(\varsigma)$, $\tilde{\aleph}^{\varepsilon_{
1}}_{\varepsilon_{2}}(\varsigma)$, $\tilde{\daleth}^{\varepsilon_{1}}_{
\varepsilon_{2}}(\varsigma)$, $\tilde{\beth}^{\varepsilon_{1}}_{\varepsilon_{2}}
(\varsigma)$, and $\tilde{\gimel}^{\varepsilon_{1}}_{\varepsilon_{2}}
(\varsigma)$ are given by Equations~\eqref{eqprop66}--\eqref{eqprop76} 
and Equations~\eqref{eqprop98}--\eqref{eqprop101}. For $n \! \in \! 
\mathbb{N}$ and $k \! \in \! \lbrace 1,2,\dotsc,K \rbrace$ such that 
$\alpha_{p_{\mathfrak{s}}} \! := \! \alpha_{k} \! \neq \! \infty$, recall the 
expression for $\tilde{\gamma}(z)$ defined by Equation~\eqref{eqmainfin10}: 
a careful analysis of the branch cuts shows that
\begin{equation}
(\tilde{\gamma}(z))^{2} \underset{\tilde{\mathcal{Y}}^{\pm} \supset 
\mathbb{C}_{\pm} \ni z \to \tilde{a}_{j}}{=} \pm \dfrac{\left(\tilde{
\mathfrak{Q}}_{0}(\tilde{a}_{j}) \! + \! \tilde{\mathfrak{Q}}_{1}(\tilde{a}_{j})
(z \! - \! \tilde{a}_{j}) \! + \! \frac{1}{2} \tilde{\mathfrak{Q}}_{2}(\tilde{a}_{j})
(z \! - \! \tilde{a}_{j})^{2} \! + \! \mathcal{O}((z \! - \! \tilde{a}_{j})^{3}) 
\right)}{(z \! - \! \tilde{a}_{j})^{1/2}}, \quad j \! = \! 1,2,\dotsc,N \! + \! 1, 
\label{eqlempro14}
\end{equation}
where $\tilde{\mathfrak{Q}}_{0}(\tilde{a}_{j^{\prime}})$, $\tilde{\mathfrak{
Q}}_{1}(\tilde{a}_{j^{\prime}})$, and $\tilde{\mathfrak{Q}}_{2}(\tilde{a}_{
j^{\prime}})$, $j^{\prime} \! = \! 1,2,\dotsc,N \! + \! 1$, are given in 
Equations~\eqref{eqprop102}--\eqref{eqprop107}. To proceed, and in order 
to eschew a proliferation of overarching notation, consider, for example, 
the $(1 \, 1)$-element of the Expansion~\eqref{eqlempro1} up to terms that are, 
in the double-scaling limit $\mathscr{N},n \! \to \! \infty$ such that $z_{o} \! = \! 
1 \! + \! o(1)$, $\mathcal{O}((((n \! - \! 1)K \! + \! k)^{2}(z \! - \! \tilde{a}_{j})^{3}
(\tilde{G}_{\tilde{a}_{j}}(z))^{2})^{-1} \tilde{\mathbb{M}}(z) \tilde{c}^{\triangleleft}
(n,k,z_{o};\tilde{a}_{j})(\tilde{\mathbb{M}}(z))^{-1})$, $j \! = \! 1,2,\dotsc,N \! 
+ \! 1$, that is (cf. Equations~\eqref{eqlempro2} and~\eqref{eqlempro3}), 
$(((n \! - \! 1)K \! + \! k)(z \! - \! \tilde{a}_{j})^{3/2} \tilde{G}_{\tilde{a}_{j}}(z))^{-1} 
\tilde{\mathfrak{t}}_{11}(z)$.\footnote{Up to an overall minus sign, statements 
made for $(((n \! - \! 1)K \! + \! k)(z \! - \! \tilde{a}_{j})^{3/2} \tilde{G}_{
\tilde{a}_{j}}(z))^{-1} \tilde{\mathfrak{t}}_{11}(z)$ are equally applicable 
to $(((n \! - \! 1)K \! + \! k)(z \! - \! \tilde{a}_{j})^{3/2} \tilde{G}_{
\tilde{a}_{j}}(z))^{-1} \tilde{\mathfrak{t}}_{22}(z)$, since $\operatorname{tr}
(w^{\Sigma_{\tilde{\mathcal{R}}}}_{+}(z)) \! = \! 0$.} Thus, for $n \! \in 
\! \mathbb{N}$ and $k \! \in \! \lbrace 1,2,\dotsc,K \rbrace$ such 
that $\alpha_{p_{\mathfrak{s}}} \! := \! \alpha_{k} \! \neq \! \infty$, 
substituting the Expansions~\eqref{eqlempro9}--\eqref{eqlempro14} 
into Equation~\eqref{eqlempro3} and the prescription $(((n \! - \! 1)K 
\! + \! k)(z \! - \! \tilde{a}_{j})^{3/2} \tilde{G}_{\tilde{a}_{j}}(z))^{-1} 
\tilde{\mathfrak{t}}_{11}(z)$, $j \! = \! 1,2,\dotsc,N \! + \! 1$, collecting 
coefficients of like powers of $(z \! - \! \tilde{a}_{j})^{-m/2}(((n \! - \! 1)
K \! + \! k) \tilde{G}_{\tilde{a}_{j}}(z))^{-1} \exp (\mi ((n \! - \! 1)K \! + \! 
k) \tilde{\mho}_{j})$, $m \! = \! 0,1,2,3,4$, and setting, for economy of 
notation, $\tilde{\mathfrak{Q}}_{i} \! := \! \tilde{\mathfrak{Q}}_{i}(\tilde{
a}_{j})$, $i \! = \! 0,1,2$, $\tilde{\kappa}_{1} \! := \! \tilde{\kappa}_{1}
(\tilde{a}_{j})$, $\tilde{\kappa}_{2} \! := \! \tilde{\kappa}_{2}(\tilde{a}_{j})$, 
$\tilde{\aleph}^{\varepsilon_{1}}_{\varepsilon_{2}} \! := \! \tilde{\aleph}^{
\varepsilon_{1}}_{\varepsilon_{2}}(\tilde{a}_{j})$, $\tilde{\daleth}^{
\varepsilon_{1}}_{\varepsilon_{2}} \! := \! \tilde{\daleth}^{\varepsilon_{1}}_{
\varepsilon_{2}}(\tilde{a}_{j})$, $\tilde{\beth}^{\varepsilon_{1}}_{\varepsilon_{
2}} \! := \! \tilde{\beth}^{\varepsilon_{1}}_{\varepsilon_{2}}(\tilde{a}_{j})$, 
and $\tilde{\gimel}^{\varepsilon_{1}}_{\varepsilon_{2}} \! := \! \tilde{
\gimel}^{\varepsilon_{1}}_{\varepsilon_{2}}(\tilde{a}_{j})$, one arrives at, 
for $j \! = \! 1,2,\dotsc,N \! + \! 1$:
\begin{align}
\mathcal{O} \left(\dfrac{(z \! - \! \tilde{a}_{j})^{-2} \me^{\mi ((n-1)K+k) 
\tilde{\mho}_{j}}}{((n \! - \! 1)K \! + \! k) \tilde{G}_{\tilde{a}_{j}}(z)} \right):& 
\, -\dfrac{(s_{1} \! + \! t_{1}) \tilde{\kappa}_{1} \tilde{\kappa}_{2} \tilde{
\mathfrak{Q}}_{0}}{4} \! - \! \dfrac{(s_{1} \! + \! t_{1}) \tilde{\kappa}_{1} 
\tilde{\kappa}_{2} \tilde{\mathfrak{Q}}_{0}}{4} \! - \! \dfrac{(s_{1} \! - \! t_{1}) 
\tilde{\kappa}_{1} \tilde{\kappa}_{2} \tilde{\mathfrak{Q}}_{0}}{4} \! - \! 
\dfrac{(s_{1} \! - \! t_{1}) \tilde{\kappa}_{1} \tilde{\kappa}_{2} \tilde{
\mathfrak{Q}}_{0}}{4}; \label{eqlempro15} \\
\mathcal{O} \left(\dfrac{(z \! - \! \tilde{a}_{j})^{-3/2} \me^{\mi ((n-1)K
+k) \tilde{\mho}_{j}}}{((n \! - \! 1)K \! + \! k) \tilde{G}_{\tilde{a}_{j}}(z)} 
\right):& \, -\dfrac{\mi (s_{1} \! + \! t_{1}) \tilde{\kappa}_{1} \tilde{
\kappa}_{2} \tilde{\mathfrak{Q}}_{0}(\tilde{\aleph}_{1}^{1} \! - \! \tilde{
\aleph}_{-1}^{-1})}{4} \! - \! \dfrac{\mi (s_{1} \! + \! t_{1}) \tilde{\kappa}_{1} 
\tilde{\kappa}_{2} \tilde{\mathfrak{Q}}_{0}(\tilde{\aleph}^{1}_{-1} \! - \! 
\tilde{\aleph}^{-1}_{1})}{4} \nonumber \\
& \, -\dfrac{\mi (s_{1} \! - \! t_{1}) \tilde{\kappa}_{1} \tilde{\kappa}_{2} 
\tilde{\mathfrak{Q}}_{0}(\tilde{\aleph}^{1}_{-1} \! + \! \tilde{\aleph}^{
1}_{1})}{4} \! + \! \dfrac{\mi (s_{1} \! - \! t_{1}) \tilde{\kappa}_{1} \tilde{
\kappa}_{2} \tilde{\mathfrak{Q}}_{0}(\tilde{\aleph}^{-1}_{-1} \! + \! 
\tilde{\aleph}^{-1}_{1})}{4} \nonumber \\
& \, -\dfrac{(s_{1} \! + \! t_{1}) \tilde{\kappa}_{1} \tilde{\kappa}_{2}}{2} \! 
+ \! \dfrac{(s_{1} \! + \! t_{1}) \tilde{\kappa}_{1} \tilde{\kappa}_{2}}{2}; 
\label{eqlempro16} \\
\mathcal{O} \left(\dfrac{(z \! - \! \tilde{a}_{j})^{-1} \, \me^{\mi ((n-1)K+k)
\tilde{\mho}_{j}}}{((n \! - \! 1)K \! + \! k) \tilde{G}_{\tilde{a}_{j}}(z)} \right):& 
\, -\dfrac{(s_{1} \! + \! t_{1}) \tilde{\kappa}_{1} \tilde{\kappa}_{2}}{4} \left(
\tilde{\mathfrak{Q}}_{1} \! + \! \tilde{\mathfrak{Q}}_{0} \left(\tilde{\daleth}^{
-1}_{-1} \! + \! \tilde{\daleth}^{1}_{1} \! + \! \tilde{\aleph}^{1}_{1} \tilde{
\aleph}^{-1}_{-1} \right) \right) \! - \! \dfrac{(s_{1} \! + \! t_{1}) \tilde{
\kappa}_{1} \tilde{\kappa}_{2}}{4 \tilde{\mathfrak{Q}}_{0}} \nonumber \\
& \, -\dfrac{(s_{1} \! + \! t_{1}) \tilde{\kappa}_{1} \tilde{\kappa}_{2}}{4} \left(
\tilde{\mathfrak{Q}}_{1} \! + \! \tilde{\mathfrak{Q}}_{0} \left(\tilde{\daleth}^{
1}_{-1} \! + \! \tilde{\daleth}^{-1}_{1} \! + \! \tilde{\aleph}^{-1}_{1} \tilde{
\aleph}^{1}_{-1} \right) \right) \! - \! \dfrac{(s_{1} \! + \! t_{1}) \tilde{
\kappa}_{1} \tilde{\kappa}_{2}}{4 \tilde{\mathfrak{Q}}_{0}} \nonumber \\
& \, -\dfrac{(s_{1} \! - \! t_{1}) \tilde{\kappa}_{1} \tilde{\kappa}_{2}}{4} 
\left(\tilde{\mathfrak{Q}}_{1} \! + \! \tilde{\mathfrak{Q}}_{0} \left(
\tilde{\daleth}^{1}_{-1} \! + \! \tilde{\daleth}^{1}_{1} \! - \! \tilde{
\aleph}^{1}_{1} \tilde{\aleph}^{1}_{-1} \right) \right) \! + \! \dfrac{(s_{1} \! 
- \! t_{1}) \tilde{\kappa}_{1} \tilde{\kappa}_{2}}{4 \tilde{\mathfrak{Q}}_{0}} 
\nonumber \\
& \, -\dfrac{(s_{1} \! - \! t_{1}) \tilde{\kappa}_{1} \tilde{\kappa}_{2}}{4} 
\left(\tilde{\mathfrak{Q}}_{1} \! + \! \tilde{\mathfrak{Q}}_{0} \left(\tilde{
\daleth}^{-1}_{-1} \! + \! \tilde{\daleth}^{-1}_{1} \! - \! \tilde{\aleph}^{-
1}_{1} \tilde{\aleph}^{-1}_{-1} \right) \right) \! + \! \dfrac{(s_{1} \! - \! 
t_{1}) \tilde{\kappa}_{1} \tilde{\kappa}_{2}}{4 \tilde{\mathfrak{Q}}_{0}} 
\nonumber \\
& \, -\dfrac{\mi (s_{1} \! + \! t_{1}) \tilde{\kappa}_{1} \tilde{\kappa}_{2}}{2} 
\left(\tilde{\aleph}^{1}_{1} \! - \! \tilde{\aleph}^{-1}_{-1} \right) \! + \! 
\dfrac{\mi (s_{1} \! + \! t_{1}) \tilde{\kappa}_{1} \tilde{\kappa}_{2}}{2} 
\left(\tilde{\aleph}^{1}_{-1} \! - \! \tilde{\aleph}^{-1}_{1} \right); 
\label{eqlempro17} \\
\mathcal{O} \left(\dfrac{(z \! - \! \tilde{a}_{j})^{-1/2} \, \me^{\mi 
((n-1)K+k) \tilde{\mho}_{j}}}{((n \! - \! 1)K \! + \! k) \tilde{G}_{\tilde{a}_{j}}(z)} 
\right):& \, -\dfrac{\mi (s_{1} \! + \! t_{1}) \tilde{\kappa}_{1} \tilde{\kappa}_{
2}}{4} \left(\tilde{\mathfrak{Q}}_{1} \left(\tilde{\aleph}^{1}_{1} \! - \! \tilde{
\aleph}^{-1}_{-1} \right) \! + \! \tilde{\mathfrak{Q}}_{0} \left(\tilde{\beth}^{
1}_{1} \! - \! \tilde{\beth}^{-1}_{-1} \! + \! \tilde{\aleph}^{1}_{1} \tilde{
\daleth}^{-1}_{-1} \! - \! \tilde{\aleph}^{-1}_{-1} \tilde{\daleth}^{1}_{1} 
\right) \right) \nonumber \\
& \, -\dfrac{\mi (s_{1} \! + \! t_{1}) \tilde{\kappa}_{1} \tilde{\kappa}_{2}}{4} 
\left(\tilde{\mathfrak{Q}}_{1} \left(\tilde{\aleph}^{1}_{-1} \! - \! \tilde{
\aleph}^{-1}_{1} \right) \! + \! \tilde{\mathfrak{Q}}_{0} \left(\tilde{\beth}^{
1}_{-1} \! - \! \tilde{\beth}^{-1}_{1} \! + \! \tilde{\aleph}^{1}_{-1} \tilde{
\daleth}^{-1}_{1} \! - \! \tilde{\aleph}^{-1}_{1} \tilde{\daleth}^{1}_{-1} 
\right) \right) \nonumber \\
& \, -\dfrac{\mi (s_{1} \! - \! t_{1}) \tilde{\kappa}_{1} \tilde{\kappa}_{2}}{4} 
\left(\tilde{\mathfrak{Q}}_{1} \left(\tilde{\aleph}^{1}_{-1} \! + \! \tilde{
\aleph}^{1}_{1} \right) \! + \! \tilde{\mathfrak{Q}}_{0} \left(\tilde{\beth}^{
1}_{-1} \! + \! \tilde{\beth}^{1}_{1} \! + \! \tilde{\aleph}^{1}_{1} \tilde{
\daleth}^{1}_{-1} \! + \! \tilde{\aleph}^{1}_{-1} \tilde{\daleth}^{1}_{1} 
\right) \right) \nonumber \\
& \, +\dfrac{\mi (s_{1} \! - \! t_{1}) \tilde{\kappa}_{1} \tilde{\kappa}_{2}}{4} 
\left(\tilde{\mathfrak{Q}}_{1} \left(\tilde{\aleph}^{-1}_{-1} \! + \! \tilde{
\aleph}^{-1}_{1} \right) \! + \! \tilde{\mathfrak{Q}}_{0} \left(\tilde{\beth}^{
-1}_{-1} \! + \! \tilde{\beth}^{-1}_{1} \! + \! \tilde{\aleph}^{-1}_{1} \tilde{
\daleth}^{-1}_{-1} \! + \! \tilde{\aleph}^{-1}_{-1} \tilde{\daleth}^{-1}_{1} 
\right) \right) \nonumber \\
& \, -\dfrac{\mi (s_{1} \! + \! t_{1}) \tilde{\kappa}_{1} \tilde{\kappa}_{2}}{4 
\tilde{\mathfrak{Q}}_{0}} \left(\tilde{\aleph}^{1}_{1} \! - \! \tilde{\aleph}^{
-1}_{-1} \right) \! - \! \dfrac{(s_{1} \! + \! t_{1}) \tilde{\kappa}_{1} \tilde{
\kappa}_{2}}{2} \left(\tilde{\daleth}^{-1}_{-1} \! + \! \tilde{\daleth}^{1}_{1} 
\! + \! \tilde{\aleph}^{1}_{1} \tilde{\aleph}^{-1}_{-1} \right) \nonumber \\
& \, -\dfrac{\mi (s_{1} \! + \! t_{1}) \tilde{\kappa}_{1} \tilde{\kappa}_{2}}{4 
\tilde{\mathfrak{Q}}_{0}} \left(\tilde{\aleph}^{1}_{-1} \! - \! \tilde{\aleph}^{
-1}_{1} \right) \! + \! \dfrac{(s_{1} \! + \! t_{1}) \tilde{\kappa}_{1} \tilde{
\kappa}_{2}}{2} \left(\tilde{\daleth}^{1}_{-1} \! + \! \tilde{\daleth}^{-1}_{1} 
\! + \! \tilde{\aleph}^{-1}_{1} \tilde{\aleph}^{1}_{-1} \right) \nonumber \\
& \, +\dfrac{\mi (s_{1} \! - \! t_{1}) \tilde{\kappa}_{1} \tilde{\kappa}_{2}}{4 
\tilde{\mathfrak{Q}}_{0}} \left(\tilde{\aleph}^{1}_{-1} \! + \! \tilde{\aleph}^{
1}_{1} \right) \! - \! \dfrac{\mi (s_{1} \! - \! t_{1}) \tilde{\kappa}_{1} \tilde{
\kappa}_{2}}{4 \tilde{\mathfrak{Q}}_{0}} \left(\tilde{\aleph}^{-1}_{1} \! + \! 
\tilde{\aleph}^{-1}_{-1} \right); \label{eqlempro18} \\
\mathcal{O} \left(\dfrac{\me^{\mi ((n-1)K+k) \tilde{\mho}_{j}}}{((n \! - \! 
1)K \! + \! k) \tilde{G}_{\tilde{a}_{j}}(z)} \right):& \, -\dfrac{(s_{1} \! + \! t_{1}) 
\tilde{\kappa}_{1} \tilde{\kappa}_{2}}{4} \left(\tilde{\mathfrak{Q}}_{0} \left(
\tilde{\gimel}^{-1}_{-1} \! + \! \tilde{\gimel}^{1}_{1} \! + \! \tilde{\daleth}^{
1}_{1} \tilde{\daleth}^{-1}_{-1} \! + \! \tilde{\aleph}^{-1}_{-1} \tilde{\beth}^{
1}_{1} \! + \! \tilde{\aleph}^{1}_{1} \tilde{\beth}^{-1}_{-1} \right) \! + \! 
\tilde{\mathfrak{Q}}_{1} \left(\tilde{\daleth}^{-1}_{-1} \! + \! \tilde{\daleth}^{
1}_{1} \! + \! \tilde{\aleph}^{1}_{1} \tilde{\aleph}^{-1}_{-1} \right) \right. 
\nonumber \\
&\left. \, +\dfrac{1}{2} \tilde{\mathfrak{Q}}_{2} \right) \! - \! \dfrac{(s_{1} 
\! + \! t_{1}) \tilde{\kappa}_{1} \tilde{\kappa}_{2}}{4} \left(\tilde{\mathfrak{
Q}}_{0} \left(\tilde{\gimel}^{1}_{-1} \! + \! \tilde{\gimel}^{-1}_{1} \! + \! 
\tilde{\daleth}^{-1}_{1} \tilde{\daleth}^{1}_{-1} \! + \! \tilde{\aleph}^{-1}_{1} 
\tilde{\beth}^{1}_{-1} \! + \! \tilde{\aleph}^{1}_{-1} \tilde{\beth}^{-1}_{1} 
\right) \! + \! \dfrac{1}{2} \tilde{\mathfrak{Q}}_{2} \right. \nonumber \\
&\left. \, +\tilde{\mathfrak{Q}}_{1} \left(\tilde{\daleth}^{1}_{-1} \! + \! 
\tilde{\daleth}^{-1}_{1} \! + \! \tilde{\aleph}^{-1}_{1} \tilde{\aleph}^{1}_{-1} 
\right) \right) \! - \! \dfrac{(s_{1} \! - \! t_{1}) \tilde{\kappa}_{1} \tilde{
\kappa}_{2}}{4} \left(\tilde{\mathfrak{Q}}_{0} \left(\tilde{\gimel}^{1}_{-1} 
\! + \! \tilde{\gimel}^{1}_{1} \! + \! \tilde{\daleth}^{1}_{1} \tilde{\daleth}^{
1}_{-1} \! - \! \tilde{\aleph}^{1}_{1} \tilde{\beth}^{1}_{-1} \! - \! \tilde{
\aleph}^{1}_{-1} \tilde{\beth}^{1}_{1} \right) \right. \nonumber \\
&\left. \, +\dfrac{1}{2} \tilde{\mathfrak{Q}}_{2} \! + \! \tilde{\mathfrak{
Q}}_{1} \left(\tilde{\daleth}^{1}_{-1} \! + \! \tilde{\daleth}^{1}_{1} \! - 
\! \tilde{\aleph}^{1}_{1} \tilde{\aleph}^{1}_{-1} \right) \right) \! - \! 
\dfrac{(s_{1} \! - \! t_{1}) \tilde{\kappa}_{1} \tilde{\kappa}_{2}}{4} \left(
\tilde{\mathfrak{Q}}_{0} \left(\tilde{\gimel}^{-1}_{-1} \! + \! \tilde{
\gimel}^{-1}_{1} \! + \! \tilde{\daleth}^{-1}_{1} \tilde{\daleth}^{-1}_{-1} 
\! - \! \tilde{\aleph}^{-1}_{1} \tilde{\beth}^{-1}_{-1} \right. \right. 
\nonumber \\
&\left. \left. \, -\tilde{\aleph}^{-1}_{-1} \tilde{\beth}^{-1}_{1} \right) \! 
+ \! \dfrac{1}{2} \tilde{\mathfrak{Q}}_{2} \! + \! \tilde{\mathfrak{Q}}_{1} 
\left(\tilde{\daleth}^{-1}_{-1} \! + \! \tilde{\daleth}^{-1}_{1} \! - \! 
\tilde{\aleph}^{-1}_{1} \tilde{\aleph}^{-1}_{-1} \right) \right) \! - \! 
\dfrac{(s_{1} \! + \! t_{1}) \tilde{\kappa}_{1} \tilde{\kappa}_{2}}{4 \tilde{
\mathfrak{Q}}_{0}} \left(\tilde{\daleth}^{-1}_{-1} \! + \! \tilde{\daleth}^{
1}_{1} \! + \! \tilde{\aleph}^{1}_{1} \tilde{\aleph}^{-1}_{-1} \right. 
\nonumber \\
&\left. \, -\tilde{\mathfrak{Q}}_{1}(\tilde{\mathfrak{Q}}_{0})^{-1} \right) 
\! - \! \dfrac{(s_{1} \! + \! t_{1}) \tilde{\kappa}_{1} \tilde{\kappa}_{2}}{4 
\tilde{\mathfrak{Q}}_{0}} \left(\tilde{\daleth}^{1}_{-1} \! + \! \tilde{
\daleth}^{-1}_{1} \! + \! \tilde{\aleph}^{-1}_{1} \tilde{\aleph}^{1}_{-1} \! 
- \! \tilde{\mathfrak{Q}}_{1}(\tilde{\mathfrak{Q}}_{0})^{-1} \right) \! - \! 
\dfrac{\mi (s_{1} \! + \! t_{1}) \tilde{\kappa}_{1} \tilde{\kappa}_{2}}{2} 
\nonumber \\
&\times \, \left(\tilde{\beth}^{1}_{1} \! - \! \tilde{\beth}^{-1}_{-1} \! + \! 
\tilde{\aleph}^{1}_{1} \tilde{\daleth}^{-1}_{-1} \! - \! \tilde{\aleph}^{-1}_{
-1} \tilde{\daleth}^{1}_{1} \right) \! + \! \dfrac{\mi (s_{1} \! + \! t_{1}) 
\tilde{\kappa}_{1} \tilde{\kappa}_{2}}{2} \left(\tilde{\beth}^{1}_{-1} \! - 
\! \tilde{\beth}^{-1}_{1} \! + \! \tilde{\aleph}^{1}_{-1} \tilde{\daleth}^{
-1}_{1} \! - \! \tilde{\aleph}^{-1}_{1} \tilde{\daleth}^{1}_{-1} \right) 
\nonumber \\
& \, +\dfrac{(s_{1} \! - \! t_{1}) \tilde{\kappa}_{1} \tilde{\kappa}_{2}}{4 
\tilde{\mathfrak{Q}}_{0}} \left(\tilde{\daleth}^{1}_{-1} \! + \! \tilde{
\daleth}^{1}_{1} \! - \! \tilde{\aleph}^{1}_{1} \tilde{\aleph}^{1}_{-1} \! - 
\! \tilde{\mathfrak{Q}}_{1}(\tilde{\mathfrak{Q}}_{0})^{-1} \right) \! + \! 
\dfrac{(s_{1} \! - \! t_{1}) \tilde{\kappa}_{1} \tilde{\kappa}_{2}}{4 \tilde{
\mathfrak{Q}}_{0}} \left(\tilde{\daleth}^{-1}_{-1} \! + \! \tilde{\daleth}^{
-1}_{1} \right. \nonumber \\
&\left. \, -\tilde{\aleph}^{-1}_{1} \tilde{\aleph}^{-1}_{-1} \! - \! \tilde{
\mathfrak{Q}}_{1}(\tilde{\mathfrak{Q}}_{0})^{-1} \right). \label{eqlempro19}
\end{align}
For $n \! \in \! \mathbb{N}$ and $k \! \in \! \lbrace 1,2,\dotsc,K \rbrace$ 
such that $\alpha_{p_{\mathfrak{s}}} \! := \! \alpha_{k} \! \neq \! \infty$, 
one shows {}from Equations~\eqref{eqprop67}--\eqref{eqprop76} 
and~\eqref{eqprop98}--\eqref{eqprop101} that $\tilde{\aleph}^{-1}_{1}
(\tilde{a}_{j}) \! = \! \tilde{\aleph}^{1}_{1}(\tilde{a}_{j})$, $\tilde{\aleph}^{-
1}_{-1}(\tilde{a}_{j}) \! = \! \tilde{\aleph}^{1}_{-1}(\tilde{a}_{j})$, $\tilde{
\daleth}^{-1}_{1}(\tilde{a}_{j}) \! = \! \tilde{\daleth}^{1}_{1}(\tilde{a}_{j})$, 
$\tilde{\daleth}^{-1}_{-1}(\tilde{a}_{j}) \! = \! \tilde{\daleth}^{1}_{-1}
(\tilde{a}_{j})$, $\tilde{\beth}^{-1}_{1}(\tilde{a}_{j}) \! = \! \tilde{\beth}^{
1}_{1}(\tilde{a}_{j})$, $\tilde{\beth}^{-1}_{-1}(\tilde{a}_{j}) \! = \! \tilde{
\beth}^{1}_{-1}(\tilde{a}_{j})$, $\tilde{\gimel}^{-1}_{1}(\tilde{a}_{j}) \! = \! 
\tilde{\gimel}^{1}_{1}(\tilde{a}_{j})$, and $\tilde{\gimel}^{-1}_{-1}(\tilde{
a}_{j}) \! = \! \tilde{\gimel}^{1}_{-1}(\tilde{a}_{j})$, $j \! = \! 1,2,\dotsc,
N \! + \! 1$; in particular, via these latter relations, one shows 
that Equations~\eqref{eqlempro16} and~\eqref{eqlempro18} are, in 
fact, identically equal to zero, and that Equations~\eqref{eqlempro15}, 
\eqref{eqlempro17}, and~\eqref{eqlempro19} give rise to non-trivial 
contributions (see below). Repeating, \emph{verbatim}, for the $(1 \, 2)$- 
and $(2 \, 1)$-elements of $(((n \! - \! 1)K \! + \! k)(z \! - \! 
\tilde{a}_{j})^{3/2} \tilde{G}_{\tilde{a}_{j}}(z))^{-1} 
\left(
\begin{smallmatrix}
\tilde{\mathfrak{t}}_{11}(z) & \tilde{\mathfrak{t}}_{12}(z) \\
\tilde{\mathfrak{t}}_{21}(z) & -\tilde{\mathfrak{t}}_{11}(z)
\end{smallmatrix}
\right)$, $j \! = \! 1,2,\dotsc,N \! + \! 1$, the analogue of the 
above-described calculation, one shows that, for $n \! \in \! \mathbb{N}$ 
and $k \! \in \! \lbrace 1,2,\dotsc,K \rbrace$ such that $\alpha_{p_{
\mathfrak{s}}} \! := \! \alpha_{k} \! \neq \! \infty$, all the coefficients 
of all the like powers of $(z \! - \! \tilde{a}_{j})^{-m/2} (((n \! - \! 1)K \! 
+ \! k) \tilde{G}_{\tilde{a}_{j}}(z))^{-1} \exp(\mi ((n \! - \! 1)K \! + \! k) 
\tilde{\mho}_{j})$, $m \! = \! 1,3$, $j \! = \! 1,2,\dotsc,N \! + \! 1$, 
are identically equal to zero, and that, via the expansion (cf. 
Lemma~\ref{lem4.9}) $\tilde{G}_{\tilde{a}_{j}}(z) \! =_{z \to \tilde{a}_{j}} \! 
\tilde{\alpha}_{0}(\tilde{a}_{j}) \! + \! \tilde{\alpha}_{1}(\tilde{a}_{j})(z \! - 
\! \tilde{a}_{j}) \! + \! \tilde{\alpha}_{2}(\tilde{a}_{j})(z \! - \! \tilde{a}_{j})^{2} 
\! + \! \mathcal{O}((z \! - \! \tilde{a}_{j})^{3})$, $j \! = \! 1,2,\dotsc,N \! 
+ \! 1$, where $\tilde{\alpha}_{0}(\tilde{a}_{j}) \! := \! 2f(\tilde{a}_{j})/3$, 
$\tilde{\alpha}_{1}(\tilde{a}_{j}) \! := \! 2f^{\prime}(\tilde{a}_{j})/5$, and 
$\tilde{\alpha}_{2}(\tilde{a}_{j}) \! := \! f^{\prime \prime}(\tilde{a}_{j})/7$, 
with $f(\tilde{a}_{j})$, $f^{\prime}(\tilde{a}_{j})$, and $f^{\prime \prime}
(\tilde{a}_{j})$ given in Lemma~\ref{lem4.9}, $w^{\Sigma_{\tilde{\mathcal{R}}}}_{+}(z)$ 
has, after a lengthy algebraic calculation and upon re-inserting explicit $\tilde{a}_{j}$, 
$j \! = \! 1,2,\dotsc,N \! + \! 1$, dependencies, the expansion
\begin{align*}
w_{+}^{\Sigma_{\hat{\mathcal{R}}}}(z) \underset{\underset{z_{o}=1+o(1)}{
\mathscr{N},n \to \infty}}{=}& \, \dfrac{1}{((n \! - \! 1)K \! + \! k)} \left(
\dfrac{(\tilde{\alpha}_{0}(\tilde{a}_{j}))^{-1}}{(z \! - \! \tilde{a}_{j})^{2}} \tilde{
\boldsymbol{\mathrm{A}}}(\tilde{a}_{j}) \! + \! \dfrac{(\tilde{\alpha}_{0}
(\tilde{a}_{j}))^{-2}}{z \! - \! \tilde{a}_{j}} \left(\tilde{\alpha}_{0}(\tilde{a}_{j}) 
\tilde{\boldsymbol{\mathrm{B}}}(\tilde{a}_{j}) \! - \! \tilde{\alpha}_{1}
(\tilde{a}_{j}) \tilde{\boldsymbol{\mathrm{A}}}(\tilde{a}_{j}) \right) \right. \\
&\left. \, +(\tilde{\alpha}_{0}(\tilde{a}_{j}))^{-2} \left(\tilde{\alpha}_{0}
(\tilde{a}_{j}) \left(\left(\dfrac{\tilde{\alpha}_{1}(\tilde{a}_{j})}{\tilde{\alpha}_{0}
(\tilde{a}_{j})} \right)^{2} \! - \! \dfrac{\tilde{\alpha}_{2}(\tilde{a}_{j})}{\tilde{
\alpha}_{0}(\tilde{a}_{j})} \right) \tilde{\boldsymbol{\mathrm{A}}}(\tilde{a}_{j}) 
\! - \! \tilde{\alpha}_{1}(\tilde{a}_{j}) \tilde{\boldsymbol{\mathrm{B}}}
(\tilde{a}_{j}) \! + \! \tilde{\alpha}_{0}(\tilde{a}_{j}) \tilde{\boldsymbol{
\mathrm{C}}}(\tilde{a}_{j}) \right) \right) \\
& \, +\mathcal{O} \left(\dfrac{1}{((n \! - \! 1)K \! + \! k)} \sum_{m \in 
\mathbb{N}} \tilde{c}_{m}^{\triangleright}(n,k,z_{o};\tilde{a}_{j})(z \! - \! 
\tilde{a}_{j})^{m} \right) \! + \! \tilde{\mathbb{E}}_{w^{\Sigma_{\tilde{
\mathcal{R}}}}_{+}}(z), \quad z \! \in \! \partial \tilde{\mathbb{U}}_{
\tilde{\delta}_{\tilde{a}_{j}}}, \quad j \! = \! 1,2,\dotsc,N \! + \! 1,
\end{align*}
where $\tilde{\boldsymbol{\mathrm{A}}}(\tilde{a}_{j})$, $\tilde{
\boldsymbol{\mathrm{B}}}(\tilde{a}_{j})$, and $\tilde{\boldsymbol{
\mathrm{C}}}(\tilde{a}_{j})$ are defined by Equations~\eqref{eqprop87}, 
\eqref{eqprop88}, and~\eqref{eqprop89}, respectively, $(\mathrm{M}_{2}
(\mathbb{C}) \! \ni)$ $\tilde{c}_{m}^{\triangleright}(n,k,z_{o};\linebreak[4] 
\tilde{a}_{j}) \! =_{\underset{z_{o}=1+o(1)}{\mathscr{N},n \to \infty}} \! 
\mathcal{O}(1)$, and
\begin{equation*}
\tilde{\mathbb{E}}_{w^{\Sigma_{\tilde{\mathcal{R}}}}_{+}}(z) 
\underset{\underset{z_{o}=1+o(1)}{\mathscr{N},n \to \infty}}{=} 
\mathcal{O} \left(\dfrac{1}{((n \! - \! 1)K \! + \! k)^{2}(z \! - \! 
\tilde{a}_{j})^{3}(\tilde{G}_{\tilde{a}_{j}}(z))^{2}} \widetilde{\mathbb{K}} 
\, \tilde{\mathbb{M}}(z) \tilde{\mathfrak{c}}^{\triangleleft}(n,k,z_{o};
\tilde{a}_{j})(\widetilde{\mathbb{K}} \, \tilde{\mathbb{M}}(z))^{-1} \right).
\end{equation*}
Finally, proceeding, \emph{verbatim}, as above, one shows, after a lengthy 
calculation, that, for $n \! \in \! \mathbb{N}$ and $k \! \in \! \lbrace 
1,2,\dotsc,K \rbrace$ such that $\alpha_{p_{\mathfrak{s}}} \! := \! 
\alpha_{k} \! \neq \! \infty$, $\tilde{\mathbb{E}}_{w^{\Sigma_{\tilde{
\mathcal{R}}}}_{+}}(z)$, has, for $z \! \in \! \partial \tilde{\mathbb{U}}_{
\tilde{\delta}_{\tilde{a}_{j}}}$, $j \! = \! 1,2,\dotsc,N \! + \! 1$, the 
expansion, in the double-scaling limit $\mathscr{N},n \! \to \! \infty$ such 
that $z_{o} \! = \! 1 \! + \! o(1)$, given in Equation~\eqref{eqproptila}, 
that is, $\tilde{\mathbb{E}}_{w^{\Sigma_{\tilde{\mathcal{R}}}}_{+}}(z) \! 
=_{\underset{z_{o}=1+o(1)}{\mathscr{N},n \to \infty}} \! \mathcal{O}
(((n \! - \! 1)K \! + \! k)^{-2} \sum_{m \in \mathbb{N}_{0}} \tilde{
\mathfrak{c}}^{\triangleleft}_{m}(n,k,z_{o};\tilde{a}_{j}) \linebreak[4] 
\pmb{\cdot}(z \! - \! \tilde{a}_{j})^{m-3})$, $z \! \in \! \partial \tilde{
\mathbb{U}}_{\tilde{\delta}_{\tilde{a}_{j}}}$, $j \! = \! 1,2,\dotsc,N \! + \! 1$, 
where $(\mathrm{M}_{2}(\mathbb{C}) \! \ni)$ $\tilde{c}_{m}^{\triangleleft}
(n,k,z_{o};\tilde{a}_{j}) \! =_{\underset{z_{o}=1+o(1)}{\mathscr{N},n \to 
\infty}} \! \mathcal{O}(1)$. For $n \! \in \! \mathbb{N}$ and $k \! \in 
\! \lbrace 1,2,\dotsc,K \rbrace$ such that $\alpha_{p_{\mathfrak{s}}} 
\! := \! \alpha_{k} \! \neq \! \infty$, the case $z \! \in \! \partial 
\tilde{\mathbb{U}}_{\tilde{\delta}_{\tilde{b}_{j-1}}}$, $j \! = \! 1,2,\dotsc,
N \! + \! 1$, is, \emph{mutatis mutandis}, analysed analogously, and 
leads to the asymptotics, in the double-scaling limit $\mathscr{N},
n \! \to \! \infty$ such that $z_{o} \! = \! 1 \! + \! o(1)$, 
for $w^{\Sigma_{\tilde{\mathcal{R}}}}_{+}(z)$ given in 
Equation~\eqref{eqproptilb}. \hfill $\qed$
\begin{eeeee} \label{rem5.1} 
\textsl{For $n \! \in \! \mathbb{N}$ and $k \! \in \! \lbrace 1,2,\dotsc,K 
\rbrace$ such that $\alpha_{p_{\mathfrak{s}}} \! := \! \alpha_{k} \! = \! 
\infty$ (resp., $\alpha_{p_{\mathfrak{s}}} \! := \! \alpha_{k} \! \neq \! 
\infty)$, note {}from Equations~\eqref{eqprop1}--\eqref{eqprop3} 
and~\eqref{eqprop34}--\eqref{eqprop36} (resp., 
Equations~\eqref{eqprop55}--\eqref{eqprop57} 
and~\eqref{eqprop87}--\eqref{eqprop89}$)$ that $\operatorname{tr}
(\hat{\boldsymbol{\mathrm{A}}}(\hat{b}_{j-1})) \! = \! \operatorname{tr}
(\hat{\boldsymbol{\mathrm{B}}}(\hat{b}_{j-1})) \! = \! \operatorname{tr}
(\hat{\boldsymbol{\mathrm{C}}}(\hat{b}_{j-1})) \! = \! \operatorname{tr}
(\hat{\boldsymbol{\mathrm{A}}}(\hat{a}_{j})) \! = \! \operatorname{tr}
(\hat{\boldsymbol{\mathrm{B}}}(\hat{a}_{j})) \! = \! \operatorname{tr}
(\hat{\boldsymbol{\mathrm{C}}}(\hat{a}_{j})) \! = \! 0$ (resp., 
$\operatorname{tr}(\tilde{\boldsymbol{\mathrm{A}}}(\tilde{b}_{j-1})) \! = \! 
\operatorname{tr}(\tilde{\boldsymbol{\mathrm{B}}}(\tilde{b}_{j-1})) \! = \! 
\operatorname{tr}(\tilde{\boldsymbol{\mathrm{C}}}(\tilde{b}_{j-1})) \! = \! 
\operatorname{tr}(\tilde{\boldsymbol{\mathrm{A}}}(\tilde{a}_{j})) \! = \! 
\operatorname{tr}(\tilde{\boldsymbol{\mathrm{B}}}(\tilde{a}_{j})) \! = \! 
\operatorname{tr}(\tilde{\boldsymbol{\mathrm{C}}}(\tilde{a}_{j})) \! = \! 0)$, 
$j \! = \! 1,2,\dotsc,N \! + \! 1$. Furthermore, it should be noted that, 
for $n \! \in \! \mathbb{N}$ and $k \! \in \! \lbrace 1,2,\dotsc,K \rbrace$ 
such that $\alpha_{p_{\mathfrak{s}}} \! := \! \alpha_{k} \! = \! \infty$ 
(resp., $\alpha_{p_{\mathfrak{s}}} \! := \! \alpha_{k} \! \neq \! \infty)$, 
the expressions for $\hat{\boldsymbol{\mathrm{C}}}(\hat{b}_{j-1})$ and 
$\hat{\boldsymbol{\mathrm{C}}}(\hat{a}_{j})$ (resp., $\tilde{\boldsymbol{
\mathrm{C}}}(\tilde{b}_{j-1})$ and $\tilde{\boldsymbol{\mathrm{C}}}
(\tilde{a}_{j}))$, $j \! = \! 1,2,\dotsc,N \! + \! 1$, defined by 
Equations~\eqref{eqprop3} and~\eqref{eqprop36} (resp., 
Equations~\eqref{eqprop57} and~\eqref{eqprop89}$)$ are necessary in 
order to derive asymptotics, in the double-scaling limit $\mathscr{N},n 
\! \to \! \infty$ such that $z_{o} \! = \! 1 \! + \! o(1)$, in the interior 
of the open discs $\hat{\mathbb{U}}_{\hat{\delta}_{\hat{b}_{j-1}}},
\hat{\mathbb{U}}_{\hat{\delta}_{\hat{a}_{j}}}$ (resp., $\tilde{\mathbb{
U}}_{\tilde{\delta}_{\tilde{b}_{j-1}}},\tilde{\mathbb{U}}_{\tilde{\delta}_{
\tilde{a}_{j}}})$ about the end-points, $\hat{b}_{j-1},\hat{a}_{j}$ (resp., 
$\tilde{b}_{j-1},\tilde{a}_{j})$, of the support, $J_{\infty}$ (resp., $J_{f})$, 
of the associated equilibrium measure, $\mu_{\widetilde{V}}^{\infty}$ 
(resp., $\mu_{\widetilde{V}}^{f})$$:$ see, also, Remark~\ref{rem5.3} below.}
\end{eeeee}
\begin{eeeee}
\textsl{The majority of the parameters appearing in the formulation 
and proof of Proposition~\ref{propo5.2} below have been defined 
heretofore in items~{\rm \pmb{(1)}} and~{\rm \pmb{(2)}} of 
Lemma~\ref{lem5.1}.}
\end{eeeee}
\begin{bbbbb} \label{propo5.2} 
For $n \! \in \! \mathbb{N}$ and $k \! \in \! \lbrace 1,2,\dotsc,K 
\rbrace$ such that $\alpha_{p_{\mathfrak{s}}} \! := \! \alpha_{k} \! 
= \! \infty$ (resp., $\alpha_{p_{\mathfrak{s}}} \! := \! \alpha_{k} 
\! \neq \! \infty)$, let $\hat{\mathcal{R}} \colon \mathbb{C} 
\setminus \hat{\Sigma}_{\hat{\mathcal{R}}}^{\sharp} \! \to \! 
\mathrm{SL}_{2}(\mathbb{C})$ (resp., $\tilde{\mathcal{R}} \colon 
\mathbb{C} \setminus \tilde{\Sigma}_{\tilde{\mathcal{R}}}^{\sharp} 
\! \to \! \mathrm{SL}_{2}(\mathbb{C}))$ solve the equivalent {\rm RHP} 
$(\hat{\mathcal{R}}(z),\hat{v}_{\hat{\mathcal{R}}}(z),\hat{\Sigma}_{
\hat{\mathcal{R}}}^{\sharp})$ (resp., $(\tilde{\mathcal{R}}(z),
\tilde{v}_{\tilde{\mathcal{R}}}(z),\tilde{\Sigma}_{\tilde{\mathcal{R}}}^{\sharp}))$, 
where, for $z \! \in \! \hat{\Sigma}_{\hat{\mathcal{R}}}^{\sharp} \setminus 
\cup_{j=1}^{N+1}(\partial \hat{\mathbb{U}}_{\hat{\delta}_{\hat{b}_{j-1}}} \cup 
\partial \hat{\mathbb{U}}_{\hat{\delta}_{\hat{a}_{j}}})$ (resp., $z \! \in \! 
\tilde{\Sigma}_{\tilde{\mathcal{R}}}^{\sharp} \setminus \cup_{j=1}^{N+1}
(\partial \tilde{\mathbb{U}}_{\tilde{\delta}_{\tilde{b}_{j-1}}} \cup \partial 
\tilde{\mathbb{U}}_{\tilde{\delta}_{\tilde{a}_{j}}}))$, the asymptotics, in the 
double-scaling limit $\mathscr{N},n \! \to \! \infty$ such that $z_{o} \! 
= \! 1 \! + \! o(1)$, of $w^{\Sigma_{\hat{\mathcal{R}}}}_{+}(z) \! = \! 
\hat{v}_{\hat{\mathcal{R}}}(z) \! - \! \mathrm{I}$ (resp., $w^{\Sigma_{
\tilde{\mathcal{R}}}}_{+}(z) \! = \! \tilde{v}_{\tilde{\mathcal{R}}}(z) 
\! - \! \mathrm{I})$ is given by Equations~\eqref{eqhtvee1}--\eqref{eqhtvee5} 
(resp., Equations~\eqref{eqtlvee8}--\eqref{eqtlvee12}$)$, and, for $z \! 
\in \! \partial \hat{\mathbb{U}}_{\hat{\delta}_{\hat{b}_{j-1}}} \cup \partial 
\hat{\mathbb{U}}_{\hat{\delta}_{\hat{a}_{j}}}$ (resp., $z \! \in \! \partial 
\tilde{\mathbb{U}}_{\tilde{\delta}_{\tilde{b}_{j-1}}} \cup \partial \tilde{
\mathbb{U}}_{\tilde{\delta}_{\tilde{a}_{j}}})$, $j \! = \! 1,2,\dotsc,N \! + \! 1$, 
the asymptotics of $w^{\Sigma_{\hat{\mathcal{R}}}}_{+}(z) \! = \! \hat{v}_{
\hat{\mathcal{R}}}(z) \! - \! \mathrm{I}$ (resp., $w^{\Sigma_{\tilde{
\mathcal{R}}}}_{+}(z) \! = \! \tilde{v}_{\tilde{\mathcal{R}}}(z) \! - \! 
\mathrm{I})$ is given by Equations~\eqref{eqprophatb} 
and~\eqref{eqprophata} (resp., Equations~\eqref{eqproptilb} 
and~\eqref{eqproptila}$)$. Then$:$ {\rm \pmb{(1)}} for $n \! \in \! 
\mathbb{N}$ and $k \! \in \! \lbrace 1,2,\dotsc,K \rbrace$ such that 
$\alpha_{p_{\mathfrak{s}}} \! := \! \alpha_{k} \! = \! \infty$$:$
\begin{enumerate}
\item[$\boldsymbol{\mathrm{(1)}_{i}}$] for $z \! \in \! \hat{\Sigma}_{p,j}^{2} 
\! := \! (\hat{a}_{j} \! + \! \hat{\delta}_{\hat{a}_{j}},\hat{b}_{j} \! - \! 
\hat{\delta}_{\hat{b}_{j}})$, $j \! = \! 1,2,\dotsc,N$,
\begin{equation} \label{eqwhtas1} 
\begin{split}
\lvert \lvert w^{\Sigma_{\hat{\mathcal{R}}}}_{+}(\pmb{\cdot}) \rvert \rvert_{
\mathcal{L}^{\infty}_{\mathrm{M}_{2}(\mathbb{C})}(\hat{\Sigma}_{p,j}^{2})} 
&\underset{\underset{z_{o}=1+o(1)}{\mathscr{N},n \to \infty}}{=} \mathcal{O} 
\left(\hat{\mathfrak{c}}_{\hat{\mathcal{R}},1}^{\sharp}(j;\infty) \me^{-((n-1)
K+k) \hat{\lambda}_{\hat{\mathcal{R}},1}^{\sharp}(j)} \right), \\
\lvert \lvert w^{\Sigma_{\hat{\mathcal{R}}}}_{+}(\pmb{\cdot}) \rvert \rvert_{
\mathcal{L}^{q}_{\mathrm{M}_{2}(\mathbb{C})}(\hat{\Sigma}_{p,j}^{2})} 
&\underset{\underset{z_{o}=1+o(1)}{\mathscr{N},n \to \infty}}{=} \mathcal{O} 
\left(\dfrac{\hat{\mathfrak{c}}_{\hat{\mathcal{R}},1}^{\sharp}(j;q) \me^{-((n-
1)K+k) \hat{\lambda}_{\hat{\mathcal{R}},1}^{\sharp}(j)}}{((n \! - \! 1)K \! + \! 
k)^{1/q}} \right), \quad q \! = \! 1,2,
\end{split}
\end{equation}
where $\hat{\lambda}_{\hat{\mathcal{R}},1}^{\sharp}(j) \! = \! \hat{\lambda}_{
\hat{\mathcal{R}},1}^{\sharp}(n,k,z_{o};j) \! =_{\underset{z_{o}=1+o(1)}{
\mathscr{N},n \to \infty}} \! \mathcal{O}(1)$ and $> \! 0$, and $\hat{
\mathfrak{c}}_{\hat{\mathcal{R}},1}^{\sharp}(j;q^{\prime}) \! = \! \hat{
\mathfrak{c}}_{\hat{\mathcal{R}},1}^{\sharp}(n,k,z_{o};j;q^{\prime}) \! 
=_{\underset{z_{o}=1+o(1)}{\mathscr{N},n \to \infty}} \! \mathcal{O}(1)$, 
$q^{\prime} \! \in \! \lbrace 1,2,\infty \rbrace$$;$
\item[$\boldsymbol{\mathrm{(1)}_{ii}}$] for $z \! \in \! \hat{\Sigma}_{p}^{1} \! 
:= \! (-\infty,\hat{b}_{0} \! - \! \hat{\delta}_{\hat{b}_{0}}) \cup (\hat{a}_{N+1} \! 
+ \! \hat{\delta}_{\hat{a}_{N+1}},+\infty)$,
\begin{equation} \label{eqwhtas2} 
\begin{split}
\lvert \lvert w^{\Sigma_{\hat{\mathcal{R}}}}_{+}(\pmb{\cdot}) \rvert \rvert_{
\mathcal{L}^{\infty}_{\mathrm{M}_{2}(\mathbb{C})}((-\infty,\hat{b}_{0}-\hat{
\delta}_{\hat{b}_{0}}))} &\underset{\underset{z_{o}=1+o(1)}{\mathscr{N},n \to 
\infty}}{=} \mathcal{O} \left(\hat{\mathfrak{c}}_{\hat{\mathcal{R}},2}^{\sharp}
(\infty) \me^{-((n-1)K+k) \hat{\lambda}_{\hat{\mathcal{R}},2}^{\sharp}(-)} 
\right), \\
\lvert \lvert w^{\Sigma_{\hat{\mathcal{R}}}}_{+}(\pmb{\cdot}) \rvert \rvert_{
\mathcal{L}^{q}_{\mathrm{M}_{2}(\mathbb{C})}((-\infty,\hat{b}_{0}-\hat{
\delta}_{\hat{b}_{0}}))} &\underset{\underset{z_{o}=1+o(1)}{\mathscr{N},n \to 
\infty}}{=} \mathcal{O} \left(\dfrac{\hat{\mathfrak{c}}_{\hat{\mathcal{R}},2}^{
\sharp}(q) \me^{-((n-1)K+k) \hat{\lambda}_{\hat{\mathcal{R}},2}^{\sharp}
(-)}}{((n \! - \! 1)K \! + \! k)^{1/q}} \right), \quad q \! = \! 1,2,
\end{split}
\end{equation}
\begin{equation} \label{eqwhtas3} 
\begin{split}
\lvert \lvert w^{\Sigma_{\hat{\mathcal{R}}}}_{+}(\pmb{\cdot}) \rvert \rvert_{
\mathcal{L}^{\infty}_{\mathrm{M}_{2}(\mathbb{C})}((\hat{a}_{N+1}+\hat{
\delta}_{\hat{a}_{N+1}},+\infty))} &\underset{\underset{z_{o}=1+o(1)}{
\mathscr{N},n \to \infty}}{=} \mathcal{O} \left(\hat{\mathfrak{c}}_{\hat{
\mathcal{R}},3}^{\sharp}(\infty) \me^{-((n-1)K+k) \hat{\lambda}_{\hat{
\mathcal{R}},2}^{\sharp}(+)} \right), \\
\lvert \lvert w^{\Sigma_{\hat{\mathcal{R}}}}_{+}(\pmb{\cdot}) \rvert \rvert_{
\mathcal{L}^{q}_{\mathrm{M}_{2}(\mathbb{C})}((\hat{a}_{N+1}+\hat{\delta}_{
\hat{a}_{N+1}},+\infty))} &\underset{\underset{z_{o}=1+o(1)}{\mathscr{N},n 
\to \infty}}{=} \mathcal{O} \left(\dfrac{\hat{\mathfrak{c}}_{\hat{\mathcal{R}},
3}^{\sharp}(q) \me^{-((n-1)K+k) \hat{\lambda}_{\hat{\mathcal{R}},2}^{
\sharp}(+)}}{((n \! - \! 1)K \! + \! k)^{1/q}} \right), \quad q \! = \! 1,2,
\end{split}
\end{equation}
where $\hat{\lambda}_{\hat{\mathcal{R}},2}^{\sharp}(\pm) \! = \! \hat{
\lambda}_{\hat{\mathcal{R}},2}^{\sharp}(n,k,z_{o};\pm) \! =_{\underset{z_{o}
=1+o(1)}{\mathscr{N},n \to \infty}} \! \mathcal{O}(1)$ and $> \! 0$, and 
$\hat{\mathfrak{c}}_{\hat{\mathcal{R}},m}^{\sharp}(q^{\prime}) \! = \! 
\hat{\mathfrak{c}}_{\hat{\mathcal{R}},m}^{\sharp}(n,k,z_{o};q^{\prime}) \! 
=_{\underset{z_{o}=1+o(1)}{\mathscr{N},n \to \infty}} \! \mathcal{O}(1)$, 
$m \! = \! 2,3$, $q^{\prime} \! \in \! \lbrace 1,2,\infty \rbrace$$;$
\item[$\boldsymbol{\mathrm{(1)}_{iii}}$] for $z \! \in \! \hat{\Sigma}_{p,j}^{3} 
\! := \! \hat{J}_{j}^{\smallfrown} \setminus (\hat{J}_{j}^{\smallfrown} \cap 
(\hat{\mathbb{U}}_{\hat{\delta}_{\hat{b}_{j-1}}} \cup \hat{\mathbb{U}}_{
\hat{\delta}_{\hat{a}_{j}}}))$, $j \! = \! 1,2,\dotsc,N \! + \! 1$,
\begin{equation} \label{eqwhtas4} 
\begin{split}
\lvert \lvert w^{\Sigma_{\hat{\mathcal{R}}}}_{+}(\pmb{\cdot}) \rvert \rvert_{
\mathcal{L}^{\infty}_{\mathrm{M}_{2}(\mathbb{C})}(\hat{\Sigma}_{p,j}^{3})} 
&\underset{\underset{z_{o}=1+o(1)}{\mathscr{N},n \to \infty}}{=} \mathcal{O} 
\left(\hat{\mathfrak{c}}_{\hat{\mathcal{R}},4}^{\sharp,\smallfrown}(j;\infty) 
\me^{-((n-1)K+k) \hat{\lambda}_{\hat{\mathcal{R}},4}^{\sharp,\smallfrown}(j)} 
\right), \\
\lvert \lvert w^{\Sigma_{\hat{\mathcal{R}}}}_{+}(\pmb{\cdot}) \rvert \rvert_{
\mathcal{L}^{q}_{\mathrm{M}_{2}(\mathbb{C})}(\hat{\Sigma}_{p,j}^{3})} 
&\underset{\underset{z_{o}=1+o(1)}{\mathscr{N},n \to \infty}}{=} \mathcal{O} 
\left(\dfrac{\hat{\mathfrak{c}}_{\hat{\mathcal{R}},4}^{\sharp,\smallfrown}(j;q) 
\me^{-((n-1)K+k) \hat{\lambda}_{\hat{\mathcal{R}},4}^{\sharp,\smallfrown}
(j)}}{((n \! - \! 1)K \! + \! k)^{1/q}} \right), \quad q \! = \! 1,2,
\end{split}
\end{equation}
and, for $z \! \in \! \hat{\Sigma}_{p,j}^{4} \! := \! \hat{J}_{j}^{\smallsmile} 
\setminus (\hat{J}_{j}^{\smallsmile} \cap (\hat{\mathbb{U}}_{\hat{\delta}_{
\hat{b}_{j-1}}} \cup \hat{\mathbb{U}}_{\hat{\delta}_{\hat{a}_{j}}}))$, $j \! = \! 
1,2,\dotsc,N \! + \! 1$,
\begin{equation} \label{eqwhtas5} 
\begin{split}
\lvert \lvert w^{\Sigma_{\hat{\mathcal{R}}}}_{+}(\pmb{\cdot}) \rvert \rvert_{
\mathcal{L}^{\infty}_{\mathrm{M}_{2}(\mathbb{C})}(\hat{\Sigma}_{p,j}^{4})} 
&\underset{\underset{z_{o}=1+o(1)}{\mathscr{N},n \to \infty}}{=} \mathcal{O} 
\left(\hat{\mathfrak{c}}_{\hat{\mathcal{R}},4}^{\sharp,\smallsmile}(j;\infty) 
\me^{-((n-1)K+k) \hat{\lambda}_{\hat{\mathcal{R}},4}^{\sharp,\smallsmile}
(j)} \right), \\
\lvert \lvert w^{\Sigma_{\hat{\mathcal{R}}}}_{+}(\pmb{\cdot}) \rvert \rvert_{
\mathcal{L}^{q}_{\mathrm{M}_{2}(\mathbb{C})}(\hat{\Sigma}_{p,j}^{4})} 
&\underset{\underset{z_{o}=1+o(1)}{\mathscr{N},n \to \infty}}{=} \mathcal{O} 
\left(\dfrac{\hat{\mathfrak{c}}_{\hat{\mathcal{R}},4}^{\sharp,\smallsmile}(j;q) 
\me^{-((n-1)K+k) \hat{\lambda}_{\hat{\mathcal{R}},4}^{\sharp,\smallsmile}
(j)}}{((n \! - \! 1)K \! + \! k)^{1/q}} \right), \quad q \! = \! 1,2,
\end{split}
\end{equation}
where $\hat{\lambda}_{\hat{\mathcal{R}},4}^{\sharp,r_{1}}(j) \! = \! 
\hat{\lambda}_{\hat{\mathcal{R}},4}^{\sharp,r_{1}}(n,k,z_{o};j) \! 
=_{\underset{z_{o}=1+o(1)}{\mathscr{N},n \to \infty}} \! \mathcal{O}(1)$ 
and $> \! 0$, $r_{1} \! \in \! \lbrace \smallfrown,\smallsmile \rbrace$, 
and $\hat{\mathfrak{c}}_{\hat{\mathcal{R}},4}^{\sharp,r_{1}}(j;q^{\prime}) 
\! = \! \hat{\mathfrak{c}}_{\hat{\mathcal{R}},4}^{\sharp,r_{1}}(n,k,z_{o};j;
q^{\prime}) \! =_{\underset{z_{o}=1+o(1)}{\mathscr{N},n \to \infty}} \! 
\mathcal{O}(1)$, $q^{\prime} \! \in \! \lbrace 1,2,\infty \rbrace$$;$ and
\item[$\boldsymbol{\mathrm{(1)}_{iv}}$] for $z \! \in \! \hat{\Sigma}_{p,j}^{5} 
\! := \! \partial \hat{\mathbb{U}}_{\hat{\delta}_{\hat{b}_{j-1}}} \cup \partial 
\hat{\mathbb{U}}_{\hat{\delta}_{\hat{a}_{j}}}$, $j \! = \! 1,2,\dotsc,N \! + \! 1$,
\begin{align}
\lvert \lvert w^{\Sigma_{\hat{\mathcal{R}}}}_{+}(\pmb{\cdot}) \rvert \rvert_{
\mathcal{L}^{q}_{\mathrm{M}_{2}(\mathbb{C})}(\partial \hat{\mathbb{U}}_{
\hat{\delta}_{\hat{b}_{j-1}}})}& \, \underset{\underset{z_{o}=1+o(1)}{\mathscr{
N},n \to \infty}}{=} \mathcal{O} \left(\dfrac{\hat{\mathfrak{c}}_{\hat{\mathcal{
R}}}^{\triangleright}(j;q)}{(n \! - \! 1)K \! + \! k} \right), \quad q \! = \! 1,2,
\infty, \label{eqwhtasbh} \\
\lvert \lvert w^{\Sigma_{\hat{\mathcal{R}}}}_{+}(\pmb{\cdot}) \rvert \rvert_{
\mathcal{L}^{q}_{\mathrm{M}_{2}(\mathbb{C})}(\partial \hat{\mathbb{U}}_{
\hat{\delta}_{\hat{a}_{j}}})}& \, \underset{\underset{z_{o}=1+o(1)}{\mathscr{
N},n \to \infty}}{=} \mathcal{O} \left(\dfrac{\hat{\mathfrak{c}}_{\hat{\mathcal{
R}}}^{\triangleleft}(j;q)}{(n \! - \! 1)K \! + \! k} \right), \quad q \! = \! 1,2,
\infty, \label{eqwhtasah}
\end{align}
where $\hat{\mathfrak{c}}_{\hat{\mathcal{R}}}^{r}(j;q) \! = \! 
\hat{\mathfrak{c}}_{\hat{\mathcal{R}}}^{r}(n,k,z_{o};j;q) \! =_{\underset{z_{o}
=1+o(1)}{\mathscr{N},n \to \infty}} \! \mathcal{O}(1)$, $r \! \in \! \lbrace 
\triangleright,\triangleleft \rbrace$$;$
\end{enumerate}
and {\rm \pmb{(2)}} for $n \! \in \! \mathbb{N}$ and $k \! \in \! \lbrace 1,2,
\dotsc,K \rbrace$ such that $\alpha_{p_{\mathfrak{s}}} \! := \! \alpha_{k} \! 
\neq \! \infty$$:$
\begin{enumerate}
\item[$\boldsymbol{\mathrm{(2)}_{i}}$] for $z \! \in \! \tilde{\Sigma}_{p,j}^{2} 
\! := \! (\tilde{a}_{j} \! + \! \tilde{\delta}_{\tilde{a}_{j}},\tilde{b}_{j} \! - \! \tilde{
\delta}_{\tilde{b}_{j}})$, $j \! = \! 1,2,\dotsc,N$,
\begin{equation} \label{eqwtlas6} 
\begin{split}
\lvert \lvert w^{\Sigma_{\tilde{\mathcal{R}}}}_{+}(\pmb{\cdot}) \rvert \rvert_{
\mathcal{L}^{\infty}_{\mathrm{M}_{2}(\mathbb{C})}(\tilde{\Sigma}_{p,j}^{2})} 
&\underset{\underset{z_{o}=1+o(1)}{\mathscr{N},n \to \infty}}{=} \mathcal{O} 
\left(\tilde{\mathfrak{c}}_{\tilde{\mathcal{R}},1}^{\sharp}(j;\infty) \me^{-((n-1)
K+k) \tilde{\lambda}_{\tilde{\mathcal{R}},1}^{\sharp}(j)} \right), \\
\lvert \lvert w^{\Sigma_{\tilde{\mathcal{R}}}}_{+}(\pmb{\cdot}) \rvert \rvert_{
\mathcal{L}^{q}_{\mathrm{M}_{2}(\mathbb{C})}(\tilde{\Sigma}_{p,j}^{2})} 
&\underset{\underset{z_{o}=1+o(1)}{\mathscr{N},n \to \infty}}{=} \mathcal{O} 
\left(\dfrac{\tilde{\mathfrak{c}}_{\tilde{\mathcal{R}},1}^{\sharp}(j;q) \me^{-((n
-1)K+k) \tilde{\lambda}_{\tilde{\mathcal{R}},1}^{\sharp}(j)}}{((n \! - \! 1)K \! 
+ \! k)^{1/q}} \right), \quad q \! = \! 1,2,
\end{split}
\end{equation}
where $\tilde{\lambda}_{\tilde{\mathcal{R}},1}^{\sharp}(j) \! = \! 
\tilde{\lambda}_{\tilde{\mathcal{R}},1}^{\sharp}(n,k,z_{o};j) \! 
=_{\underset{z_{o}=1+o(1)}{\mathscr{N},n \to \infty}} \! \mathcal{O}(1)$ 
and $> \! 0$, and $\tilde{\mathfrak{c}}_{\tilde{\mathcal{R}},1}^{\sharp}(j;
q^{\prime}) \! = \! \tilde{\mathfrak{c}}_{\tilde{\mathcal{R}},1}^{\sharp}(n,k,
z_{o};j;q^{\prime}) \! =_{\underset{z_{o}=1+o(1)}{\mathscr{N},n \to \infty}} 
\! \mathcal{O}(1)$, $q^{\prime} \! \in \! \lbrace 1,2,\infty \rbrace$$;$
\item[$\boldsymbol{\mathrm{(2)}_{ii}}$] for $z \! \in \! \tilde{\Sigma}_{p}^{1} 
\! := \! (-\infty,\tilde{b}_{0} \! - \! \tilde{\delta}_{\tilde{b}_{0}}) \cup 
(\tilde{a}_{N+1} \! + \! \tilde{\delta}_{\tilde{a}_{N+1}},+\infty)$,
\begin{equation} \label{eqwtlas7} 
\begin{split}
\lvert \lvert w^{\Sigma_{\tilde{\mathcal{R}}}}_{+}(\pmb{\cdot}) \rvert \rvert_{
\mathcal{L}^{\infty}_{\mathrm{M}_{2}(\mathbb{C})}((-\infty,\tilde{b}_{0}-
\tilde{\delta}_{\tilde{b}_{0}}))} &\underset{\underset{z_{o}=1+o(1)}{\mathscr{N},
n \to \infty}}{=} \mathcal{O} \left(\tilde{\mathfrak{c}}_{\tilde{\mathcal{R}},
2}^{\sharp}(\infty) \me^{-((n-1)K+k) \tilde{\lambda}_{\tilde{\mathcal{R}},
2}^{\sharp}(-)} \right), \\
\lvert \lvert w^{\Sigma_{\tilde{\mathcal{R}}}}_{+}(\pmb{\cdot}) \rvert \rvert_{
\mathcal{L}^{q}_{\mathrm{M}_{2}(\mathbb{C})}((-\infty,\tilde{b}_{0}-\tilde{
\delta}_{\tilde{b}_{0}}))} &\underset{\underset{z_{o}=1+o(1)}{\mathscr{N},n \to 
\infty}}{=} \mathcal{O} \left(\dfrac{\tilde{\mathfrak{c}}_{\tilde{\mathcal{R}},
2}^{\sharp}(q) \me^{-((n-1)K+k) \tilde{\lambda}_{\tilde{\mathcal{R}},2}^{
\sharp}(-)}}{((n \! - \! 1)K \! + \! k)^{1/q}} \right), \quad q \! = \! 1,2,
\end{split}
\end{equation}
\begin{equation} \label{eqwtlas8} 
\begin{split}
\lvert \lvert w^{\Sigma_{\tilde{\mathcal{R}}}}_{+}(\pmb{\cdot}) \rvert \rvert_{
\mathcal{L}^{\infty}_{\mathrm{M}_{2}(\mathbb{C})}((\tilde{a}_{N+1}+\tilde{
\delta}_{\tilde{a}_{N+1}},+\infty))} &\underset{\underset{z_{o}=1+o(1)}{
\mathscr{N},n \to \infty}}{=} \mathcal{O} \left(\tilde{\mathfrak{c}}_{\tilde{
\mathcal{R}},3}^{\sharp}(\infty) \me^{-((n-1)K+k) \tilde{\lambda}_{\tilde{
\mathcal{R}},2}^{\sharp}(+)} \right), \\
\lvert \lvert w^{\Sigma_{\tilde{\mathcal{R}}}}_{+}(\pmb{\cdot}) \rvert \rvert_{
\mathcal{L}^{q}_{\mathrm{M}_{2}(\mathbb{C})}((\tilde{a}_{N+1}+\tilde{\delta}_{
\tilde{a}_{N+1}},+\infty))} &\underset{\underset{z_{o}=1+o(1)}{\mathscr{N},n 
\to \infty}}{=} \mathcal{O} \left(\dfrac{\tilde{\mathfrak{c}}_{\tilde{\mathcal{
R}},3}^{\sharp}(q) \me^{-((n-1)K+k) \tilde{\lambda}_{\tilde{\mathcal{R}},
2}^{\sharp}(+)}}{((n \! - \! 1)K \! + \! k)^{1/q}} \right), \quad q \! = \! 1,2,
\end{split}
\end{equation}
where $\tilde{\lambda}_{\tilde{\mathcal{R}},2}^{\sharp}(\pm) \! = \! \tilde{
\lambda}_{\tilde{\mathcal{R}},2}^{\sharp}(n,k,z_{o};\pm) \! =_{\underset{
z_{o}=1+o(1)}{\mathscr{N},n \to \infty}} \! \mathcal{O}(1)$ and $> \! 0$, 
and $\tilde{\mathfrak{c}}_{\tilde{\mathcal{R}},m}^{\sharp}(q^{\prime}) \! = 
\! \tilde{\mathfrak{c}}_{\tilde{\mathcal{R}},m}^{\sharp}(n,k,z_{o};q^{\prime}) 
\! =_{\underset{z_{o}=1+o(1)}{\mathscr{N},n \to \infty}} \! \mathcal{O}(1)$, 
$m \! = \! 2,3$, $q^{\prime} \! \in \! \lbrace 1,2,\infty \rbrace$$;$
\item[$\boldsymbol{\mathrm{(2)}_{iii}}$] for $z \! \in \! \tilde{\Sigma}_{p,
j}^{3} \! := \! \tilde{J}_{j}^{\smallfrown} \setminus (\tilde{J}_{j}^{\smallfrown} 
\cap (\tilde{\mathbb{U}}_{\tilde{\delta}_{\tilde{b}_{j-1}}} \cup \tilde{
\mathbb{U}}_{\tilde{\delta}_{\tilde{a}_{j}}}))$, $j \! = \! 1,2,\dotsc,N \! + \! 1$,
\begin{equation} \label{eqwtlas9} 
\begin{split}
\lvert \lvert w^{\Sigma_{\tilde{\mathcal{R}}}}_{+}(\pmb{\cdot}) \rvert \rvert_{
\mathcal{L}^{\infty}_{\mathrm{M}_{2}(\mathbb{C})}(\tilde{\Sigma}_{p,j}^{3})} 
&\underset{\underset{z_{o}=1+o(1)}{\mathscr{N},n \to \infty}}{=} \mathcal{O} 
\left(\tilde{\mathfrak{c}}_{\tilde{\mathcal{R}},4}^{\sharp,\smallfrown}(j;\infty) 
\me^{-((n-1)K+k) \tilde{\lambda}_{\tilde{\mathcal{R}},4}^{\sharp,\smallfrown}
(j)} \right), \\
\lvert \lvert w^{\Sigma_{\tilde{\mathcal{R}}}}_{+}(\pmb{\cdot}) \rvert \rvert_{
\mathcal{L}^{q}_{\mathrm{M}_{2}(\mathbb{C})}(\tilde{\Sigma}_{p,j}^{3})} 
&\underset{\underset{z_{o}=1+o(1)}{\mathscr{N},n \to \infty}}{=} \mathcal{O} 
\left(\dfrac{\tilde{\mathfrak{c}}_{\tilde{\mathcal{R}},4}^{\sharp,\smallfrown}(j;q) 
\me^{-((n-1)K+k) \tilde{\lambda}_{\tilde{\mathcal{R}},4}^{\sharp,\smallfrown}
(j)}}{((n \! - \! 1)K \! + \! k)^{1/q}} \right), \quad q \! = \! 1,2,
\end{split}
\end{equation}
and, for $z \! \in \! \tilde{\Sigma}_{p,j}^{4} \! := \! \tilde{J}_{j}^{\smallsmile} 
\setminus (\tilde{J}_{j}^{\smallsmile} \cap (\tilde{\mathbb{U}}_{\tilde{\delta}_{
\tilde{b}_{j-1}}} \cup \tilde{\mathbb{U}}_{\tilde{\delta}_{\tilde{a}_{j}}}))$, $j \! 
= \! 1,2,\dotsc,N \! + \! 1$,
\begin{equation} \label{eqwtlas10} 
\begin{split}
\lvert \lvert w^{\Sigma_{\tilde{\mathcal{R}}}}_{+}(\pmb{\cdot}) \rvert \rvert_{
\mathcal{L}^{\infty}_{\mathrm{M}_{2}(\mathbb{C})}(\tilde{\Sigma}_{p,j}^{4})} 
&\underset{\underset{z_{o}=1+o(1)}{\mathscr{N},n \to \infty}}{=} \mathcal{O} 
\left(\tilde{\mathfrak{c}}_{\tilde{\mathcal{R}},4}^{\sharp,\smallsmile}(j;\infty) 
\me^{-((n-1)K+k) \tilde{\lambda}_{\tilde{\mathcal{R}},4}^{\sharp,\smallsmile}
(j)} \right), \\
\lvert \lvert w^{\Sigma_{\tilde{\mathcal{R}}}}_{+}(\pmb{\cdot}) \rvert \rvert_{
\mathcal{L}^{q}_{\mathrm{M}_{2}(\mathbb{C})}(\tilde{\Sigma}_{p,j}^{4})} 
&\underset{\underset{z_{o}=1+o(1)}{\mathscr{N},n \to \infty}}{=} \mathcal{O} 
\left(\dfrac{\tilde{\mathfrak{c}}_{\tilde{\mathcal{R}},4}^{\sharp,\smallsmile}(j;q) 
\me^{-((n-1)K+k) \tilde{\lambda}_{\tilde{\mathcal{R}},4}^{\sharp,\smallsmile}
(j)}}{((n \! - \! 1)K \! + \! k)^{1/q}} \right), \quad q \! = \! 1,2,
\end{split}
\end{equation}
where $\tilde{\lambda}_{\tilde{\mathcal{R}},4}^{\sharp,r_{1}}(j) \! = \! 
\tilde{\lambda}_{\tilde{\mathcal{R}},4}^{\sharp,r_{1}}(n,k,z_{o};j) \! 
=_{\underset{z_{o}=1+o(1)}{\mathscr{N},n \to \infty}} \! \mathcal{O}(1)$ 
and $> \! 0$, $r_{1} \! \in \! \lbrace \smallfrown,\smallsmile \rbrace$, 
and $\tilde{\mathfrak{c}}_{\tilde{\mathcal{R}},4}^{\sharp,r_{1}}(j;q^{\prime}) 
\! = \! \tilde{\mathfrak{c}}_{\tilde{\mathcal{R}},4}^{\sharp,r_{1}}(n,k,z_{o};j;
q^{\prime}) \! =_{\underset{z_{o}=1+o(1)}{\mathscr{N},n \to \infty}} \! 
\mathcal{O}(1)$, $q^{\prime} \! \in \! \lbrace 1,2,\infty \rbrace$$;$ and
\item[$\boldsymbol{\mathrm{(2)}_{iv}}$] for $z \! \in \! \tilde{\Sigma}_{p,j}^{5} 
\! := \! \partial \tilde{\mathbb{U}}_{\tilde{\delta}_{\tilde{b}_{j-1}}} \cup \partial 
\tilde{\mathbb{U}}_{\tilde{\delta}_{\tilde{a}_{j}}}$, $j \! = \! 1,2,\dotsc,
N \! + \! 1$,
\begin{align}
\lvert \lvert w^{\Sigma_{\tilde{\mathcal{R}}}}_{+}(\pmb{\cdot}) \rvert \rvert_{
\mathcal{L}^{q}_{\mathrm{M}_{2}(\mathbb{C})}(\partial \tilde{\mathbb{U}}_{
\tilde{\delta}_{\tilde{b}_{j-1}}})}& \, \underset{\underset{z_{o}=1+o(1)}{\mathscr{
N},n \to \infty}}{=} \mathcal{O} \left(\dfrac{\tilde{\mathfrak{c}}_{\tilde{\mathcal{
R}}}^{\triangleright}(j;q)}{(n \! - \! 1)K \! + \! k} \right), \quad q \! = \! 1,2,
\infty, \label{eqwtlasbt} \\
\lvert \lvert w^{\Sigma_{\tilde{\mathcal{R}}}}_{+}(\pmb{\cdot}) \rvert \rvert_{
\mathcal{L}^{q}_{\mathrm{M}_{2}(\mathbb{C})}(\partial \tilde{\mathbb{U}}_{
\tilde{\delta}_{\tilde{a}_{j}}})}& \, \underset{\underset{z_{o}=1+o(1)}{\mathscr{
N},n \to \infty}}{=} \mathcal{O} \left(\dfrac{\tilde{\mathfrak{c}}_{\tilde{\mathcal{
R}}}^{\triangleleft}(j;q)}{(n \! - \! 1)K \! + \! k} \right), \quad q \! = \! 1,2,
\infty, \label{eqwtlasat}
\end{align}
where $\tilde{\mathfrak{c}}_{\tilde{\mathcal{R}}}^{r}(j;q) \! = \! \tilde{
\mathfrak{c}}_{\tilde{\mathcal{R}}}^{r}(n,k,z_{o};j;q) \! =_{\underset{z_{o}
=1+o(1)}{\mathscr{N},n \to \infty}} \! \mathcal{O}(1)$, $r \! \in \! \lbrace 
\triangleright,\triangleleft \rbrace$.
\end{enumerate}
\end{bbbbb}

\emph{Proof}. The proof of this Proposition~\ref{propo5.2} consists of two 
cases: (i) $n \! \in \! \mathbb{N}$ and $k \! \in \! \lbrace 1,2,\dotsc,K \rbrace$ 
such that $\alpha_{p_{\mathfrak{s}}} \! := \! \alpha_{k} \! = \! \infty$; and 
(ii) $n \! \in \! \mathbb{N}$ and $k \! \in \! \lbrace 1,2,\dotsc,K \rbrace$ 
such that $\alpha_{p_{\mathfrak{s}}} \! := \! \alpha_{k} \! \neq \! \infty$. 
Notwithstanding the fact that the scheme of the proof is, \emph{mutatis mutandis}, 
similar for both cases, without loss of generality, only the proof for case~(ii) 
is presented in detail, whilst case~(i) is proved analogously.

For $n \! \in \! \mathbb{N}$ and $k \! \in \! \lbrace 1,2,\dotsc,K \rbrace$ such 
that $\alpha_{p_{\mathfrak{s}}} \! := \! \alpha_{k} \! \neq \! \infty$, and for $z 
\! \in \! \tilde{\Sigma}_{p,j}^{2} \! := \! (\tilde{a}_{j} \! + \! \tilde{\delta}_{\tilde{
a}_{j}},\tilde{b}_{j} \! - \! \tilde{\delta}_{\tilde{b}_{j}})$, $j \! = \! 1,2,\dotsc,N$, 
recall {}from the proof of the corresponding item of Lemma~\ref{lem5.1} that, 
for $j \! = \! 1,2,\dotsc,N$, $\tilde{\Sigma}_{p,j}^{2} \! = \! \tilde{A}_{\tilde{
\mathcal{R}},1}(j) \cup \tilde{A}_{\tilde{\mathcal{R}},1}^{c}(j)$, where $\tilde{A}_{
\tilde{\mathcal{R}},1}(j) \! = \! (\tilde{a}_{j} \! + \! \tilde{\delta}_{\tilde{a}_{j}},
\tilde{b}_{j} \! - \! \tilde{\delta}_{\tilde{b}_{j}}) \setminus \cup_{q \in \tilde{Q}_{
\tilde{\mathcal{R}},1}(j)} \mathscr{O}_{\tilde{\delta}_{\tilde{\mathcal{R}},1}(j)}
(\alpha_{p_{q}})$, and $\tilde{A}_{\tilde{\mathcal{R}},1}^{c}(j) \! = \! \cup_{q \in 
\tilde{Q}_{\tilde{\mathcal{R}},1}(j)} \mathscr{O}_{\tilde{\delta}_{\tilde{\mathcal{R}},
1}(j)}(\alpha_{p_{q}})$, with $\tilde{Q}_{\tilde{\mathcal{R}},1}(j) \! := \! \lbrace 
\mathstrut q^{\prime} \! \in \! \lbrace 1,\dotsc,\mathfrak{s} \! - \! 2,\mathfrak{s} 
\rbrace; \, \alpha_{p_{q^{\prime}}} \! \in \! (\tilde{a}_{j} \! + \! \tilde{\delta}_{
\tilde{a}_{j}},\tilde{b}_{j} \! - \! \tilde{\delta}_{\tilde{b}_{j}}) \rbrace$, and 
sufficiently small $\tilde{\delta}_{\tilde{\mathcal{R}},1}(j) \! > \! 0$ chosen so 
that $\mathscr{O}_{\tilde{\delta}_{\tilde{\mathcal{R}},1}(j)}(\alpha_{p_{q_{1}}}) 
\cap \mathscr{O}_{\tilde{\delta}_{\tilde{\mathcal{R}},1}(j)}(\alpha_{p_{q_{2}}}) \! 
= \! \varnothing$ $\forall$ $q_{1} \! \neq \! q_{2} \! \in \! \tilde{Q}_{\tilde{
\mathcal{R}},1}(j)$ and $\mathscr{O}_{\tilde{\delta}_{\tilde{\mathcal{R}},1}(j)}
(\alpha_{p_{q_{1}}}) \cap \lbrace \tilde{a}_{j} \! + \! \tilde{\delta}_{\tilde{a}_{j}} 
\rbrace \! = \! \varnothing \! = \! \mathscr{O}_{\tilde{\delta}_{\tilde{\mathcal{R}},
1}(j)}(\alpha_{p_{q_{1}}}) \cap \lbrace \tilde{b}_{j} \! - \! \tilde{\delta}_{\tilde{b}_{j}} 
\rbrace$ (of course, $\tilde{A}_{\tilde{\mathcal{R}},1}(j) \cap \tilde{A}_{\tilde{
\mathcal{R}},1}^{c}(j) \! = \! \varnothing$). Via the formula $w^{\Sigma_{\tilde{
\mathcal{R}}}}_{+}(z) \! = \! \tilde{v}_{\tilde{\mathcal{R}}}(z) \! - \! \mathrm{I}$ 
and the asymptotics, in the double-scaling limit $\mathscr{N},n \! \to \! \infty$ 
such that $z_{o} \! = \! 1 \! + \! o(1)$, of $\tilde{v}_{\tilde{\mathcal{R}}}(z)$ 
given in Equation~\eqref{eqtlvee8}, one shows, via an integration argument, 
that, for $n \! \in \! \mathbb{N}$ and $k \! \in \! \lbrace 1,2,\dotsc,K \rbrace$ 
such that $\alpha_{p_{\mathfrak{s}}} \! := \! \alpha_{k} \! \neq \! \infty$, and 
for $z \! \in \! \tilde{\Sigma}_{p,j}^{2}$, $j \! = \! 1,2,\dotsc,N$,
\begin{equation*}
\lvert \lvert w^{\Sigma_{\tilde{\mathcal{R}}}}_{+}(\pmb{\cdot}) \rvert \rvert_{
\mathcal{L}^{\infty}_{\mathrm{M}_{2}(\mathbb{C})}(\tilde{\Sigma}_{p,j}^{2})} 
\! := \! \max_{i_{1},i_{2}=1,2} \sup_{\xi \in \tilde{A}_{\tilde{\mathcal{R}},1}
(j) \cup \tilde{A}_{\tilde{\mathcal{R}},1}^{c}(j)} \lvert (w^{\Sigma_{\tilde{
\mathcal{R}}}}_{+}(\xi))_{i_{1}i_{2}} \rvert \underset{\underset{z_{o}=1+
o(1)}{\mathscr{N},n \to \infty}}{=} \mathcal{O} \left(\tilde{\mathfrak{c}}_{
\tilde{\mathcal{R}},1}^{\sharp}(j;\infty) \me^{-((n-1)K+k) \tilde{\lambda}_{
\tilde{\mathcal{R}},1}^{\sharp}(j)} \right),
\end{equation*}
and
\begin{align*}
\lvert \lvert w^{\Sigma_{\tilde{\mathcal{R}}}}_{+}(\pmb{\cdot}) \rvert \rvert_{
\mathcal{L}^{q_{1}}_{\mathrm{M}_{2}(\mathbb{C})}(\tilde{\Sigma}_{p,j}^{2})} 
:=& \, \left(\int\limits_{\tilde{A}_{\tilde{\mathcal{R}},1}(j) \cup \tilde{A}_{
\tilde{\mathcal{R}},1}^{c}(j)} \lvert w^{\Sigma_{\tilde{\mathcal{R}}}}_{+}(\xi) 
\rvert^{q_{1}} \, \lvert \md \xi \rvert \right)^{\frac{1}{q_{1}}} \! = \! \left(
\left(\int_{\tilde{A}_{\tilde{\mathcal{R}},1}(j)} \! + \! \int_{\tilde{A}_{\tilde{
\mathcal{R}},1}^{c}(j)} \right) \left(\sum_{i_{1},i_{2}=1,2} \overline{
(w^{\Sigma_{\tilde{\mathcal{R}}}}_{+}(\xi))_{i_{1}i_{2}}} \right. \right. \\
\times&\left. \left. \, (w^{\Sigma_{\tilde{\mathcal{R}}}}_{+}(\xi))_{i_{1}i_{2}} 
\vphantom{M^{M^{M^{M^{M^{M^{M}}}}}}} \right)^{\frac{q_{1}}{2}} \, \lvert 
\md \xi \rvert \right)^{\frac{1}{q_{1}}} \underset{\underset{z_{o}=1+
o(1)}{\mathscr{N},n \to \infty}}{=} \mathcal{O} \left(
\dfrac{\tilde{\mathfrak{c}}_{\tilde{\mathcal{R}},1}^{\sharp}(j;q_{1}) 
\me^{-((n-1)K+k) \tilde{\lambda}_{\tilde{\mathcal{R}},1}^{\sharp}(j)}}{
((n \! - \! 1)K \! + \! k)^{1/q_{1}}} \right), \quad q_{1} \! = \! 1,2,
\end{align*}
where
\begin{equation} \label{eqqlamrj} 
\tilde{\lambda}_{\tilde{\mathcal{R}},1}^{\sharp}(j) \! := \! \min \left\lbrace 
\tilde{\lambda}_{\tilde{\mathcal{R}},1}(j) \min \left\lbrace \tilde{\delta}_{
\tilde{a}_{j}},\min\limits_{q^{\prime} \in \tilde{Q}_{\tilde{\mathcal{R}},1}(j)} 
\lbrace \alpha_{p_{q^{\prime}}} \! - \! \tilde{a}_{j} \! + \! \tilde{\delta}_{
\tilde{\mathcal{R}},1}(j) \rbrace \right\rbrace,K^{-1} \lvert \ln (\tilde{
\delta}_{\tilde{\mathcal{R}},1}(j)) \rvert \min\limits_{q^{\prime} \in \tilde{Q}_{
\tilde{\mathcal{R}},1}(j)} \lbrace \tilde{\mathfrak{c}}_{\tilde{\mathcal{R}},2}
(j,q^{\prime}) \rbrace \right\rbrace,
\end{equation}
with $\tilde{\lambda}_{\tilde{\mathcal{R}},1}(j)$ and $\tilde{
\mathfrak{c}}_{\tilde{\mathcal{R}},2}(j,q^{\prime})$, $q^{\prime} \! 
\in \! \tilde{Q}_{\tilde{\mathcal{R}},1}(j)$, given in item~\pmb{(2)}, 
subitem~$\boldsymbol{\mathrm{(2)}_{i}}$ of Lemma~\ref{lem5.1}, and 
$\tilde{\mathfrak{c}}_{\tilde{\mathcal{R}},1}^{\sharp}(j;q)$, $q \! \in 
\! \lbrace 1,2,\infty \rbrace$, is characterised in item~\pmb{(2)}, 
subitem~$\boldsymbol{\mathrm{(2)}_{i}}$ of the proposition.\footnote{If, 
for $j \! \in \! \lbrace 1,2,\dotsc,N \rbrace$, $\tilde{Q}_{\tilde{\mathcal{
R}},1}(j) \! = \! \varnothing$, that is, $\# \tilde{Q}_{\tilde{\mathcal{R}},1}
(j) \! = \! 0$, then $\cup_{q \in \tilde{Q}_{\tilde{\mathcal{R}},1}(j)} 
\mathscr{O}_{\tilde{\delta}_{\tilde{\mathcal{R}},1}(j)}(\alpha_{p_{q}}) 
\! := \! \varnothing$, in which case $\tilde{A}_{\tilde{\mathcal{R}},1}(j) \! 
= \! \tilde{\Sigma}_{p,j}^{2}$, and the Asymptotics~\eqref{eqwtlas6} are 
applicable provided one defines $\tilde{\lambda}_{\tilde{\mathcal{R}},1}^{
\sharp}(j) \! := \! \tilde{\lambda}_{\tilde{\mathcal{R}},1}(j) \tilde{\delta}_{
\tilde{a}_{j}}$.}

For $n \! \in \! \mathbb{N}$ and $k \! \in \! \lbrace 1,2,\dotsc,K \rbrace$ 
such that $\alpha_{p_{\mathfrak{s}}} \! := \! \alpha_{k} \! \neq \! \infty$, 
and for $z \! \in \! \tilde{\Sigma}_{p}^{1} \! := \! (-\infty,\tilde{b}_{0} \! - 
\! \tilde{\delta}_{\tilde{b}_{0}}) \cup (\tilde{a}_{N+1} \! + \! \tilde{\delta}_{
\tilde{a}_{N+1}},+\infty)$, consider, say, and without loss of generality, 
the analysis for the case $z \! \in \! (\tilde{a}_{N+1} \! + \! \tilde{\delta}_{
\tilde{a}_{N+1}},+\infty)$ $(\subset \tilde{\Sigma}_{p}^{1})$: the analysis 
for the case $z \! \in \! (-\infty,\tilde{b}_{0} \! - \! \tilde{\delta}_{\tilde{
b}_{0}})$ $(\subset \tilde{\Sigma}_{p}^{1})$ is analogous. For $n \! \in \! 
\mathbb{N}$ and $k \! \in \! \lbrace 1,2,\dotsc,K \rbrace$ such that 
$\alpha_{p_{\mathfrak{s}}} \! := \! \alpha_{k} \! \neq \! \infty$, and for $z 
\! \in \! (\tilde{a}_{N+1} \! + \! \tilde{\delta}_{\tilde{a}_{N+1}},+\infty)$, 
recall {}from the proof of the corresponding item of Lemma~\ref{lem5.1} 
that $(\tilde{a}_{N+1} \! + \! \tilde{\delta}_{\tilde{a}_{N+1}},+\infty) \! = 
\! \tilde{A}_{\tilde{\mathcal{R}},2}(+) \cup \tilde{A}_{\tilde{\mathcal{R}},
2}^{c}(+)$, where $\tilde{A}_{\tilde{\mathcal{R}},2}(+) \! = \! (\tilde{a}_{N+1} 
\! + \! \tilde{\delta}_{\tilde{a}_{N+1}},+\infty) \setminus (\mathscr{O}_{\infty}
(\alpha_{p_{\mathfrak{s}-1}}) \cup \cup_{q \in \tilde{Q}_{\tilde{\mathcal{R}},2}(+)} 
\mathscr{O}_{\tilde{\delta}_{\tilde{\mathcal{R}},2}(+)}(\alpha_{p_{q}}))$, 
and $\tilde{A}_{\tilde{\mathcal{R}},2}^{c}(+) \! = \! \mathscr{O}_{\infty}
(\alpha_{p_{\mathfrak{s}-1}}) \cup \cup_{q \in \tilde{Q}_{\tilde{\mathcal{
R}},2}(+)} \mathscr{O}_{\tilde{\delta}_{\tilde{\mathcal{R}},2}(+)}
(\alpha_{p_{q}})$, with $\tilde{Q}_{\tilde{\mathcal{R}},2}(+) \! := \! \lbrace 
\mathstrut q^{\prime} \! \in \! \lbrace 1,\dotsc,\mathfrak{s} \! - \! 2,
\mathfrak{s} \rbrace; \, \alpha_{p_{q^{\prime}}} \! \in \! (\tilde{a}_{
N+1} \! + \! \tilde{\delta}_{\tilde{a}_{N+1}},+\infty) \rbrace$ and 
sufficiently small $\tilde{\varepsilon}_{\infty},\tilde{\delta}_{\tilde{
\mathcal{R}},2}(+) \! > \! 0$ chosen so that $\mathscr{O}_{\tilde{
\delta}_{\tilde{\mathcal{R}},2}(+)}(\alpha_{p_{q^{\prime}_{1}}}) \cap 
\mathscr{O}_{\tilde{\delta}_{\tilde{\mathcal{R}},2}(+)}(\alpha_{p_{q^{\prime 
\prime}_{1}}}) \! = \! \varnothing$ $\forall$ $q^{\prime}_{1} \! \neq \! 
q^{\prime \prime}_{1} \! \in \! \tilde{Q}_{\tilde{\mathcal{R}},2}(+)$, 
$\mathscr{O}_{\tilde{\delta}_{\tilde{\mathcal{R}},2}(+)}(\alpha_{p_{
q^{\prime}_{1}}}) \cap \lbrace \tilde{a}_{N+1} \! + \! \tilde{\delta}_{
\tilde{a}_{N+1}} \rbrace \! = \! \varnothing$, and $\mathscr{O}_{
\tilde{\delta}_{\tilde{\mathcal{R}},2}(+)}(\alpha_{p_{q^{\prime}_{1}}}) \cap 
\mathscr{O}_{\infty}(\alpha_{p_{\mathfrak{s}-1}}) \! = \! \varnothing$ 
(of course, $\tilde{A}_{\tilde{\mathcal{R}},2}(+) \cap \tilde{A}_{\tilde{
\mathcal{R}},2}^{c}(+) \! = \! \varnothing$). Via the formula $w^{\Sigma_{
\tilde{\mathcal{R}}}}_{+}(z) \! = \! \tilde{v}_{\tilde{\mathcal{R}}}(z) \! - \! 
\mathrm{I}$ and the asymptotics, in the double-scaling limit $\mathscr{N},
n \! \to \! \infty$ such that $z_{o} \! = \! 1 \! + \! o(1)$, of $\tilde{v}_{
\tilde{\mathcal{R}}}(z)$ given in Equation~\eqref{eqtlvee9}, one shows, 
via an integration argument, that, for $n \! \in \! \mathbb{N}$ and $k \! 
\in \! \lbrace 1,2,\dotsc,K \rbrace$ such that $\alpha_{p_{\mathfrak{s}}} 
\! := \! \alpha_{k} \! \neq \! \infty$, and for $z \! \in \! (\tilde{a}_{N+1} 
\! + \! \tilde{\delta}_{\tilde{a}_{N+1}},+\infty)$,
\begin{equation*}
\lvert \lvert w^{\Sigma_{\tilde{\mathcal{R}}}}_{+}(\pmb{\cdot}) \rvert \rvert_{
\mathcal{L}^{\infty}_{\mathrm{M}_{2}(\mathbb{C})}((\tilde{a}_{N+1}+
\tilde{\delta}_{\tilde{a}_{N+1}},+\infty))} \! := \! \max_{i_{1},i_{2}=1,2} 
\sup_{\xi \in \tilde{A}_{\tilde{\mathcal{R}},2}(+) \cup \tilde{A}_{\tilde{
\mathcal{R}},2}^{c}(+)} \lvert (w^{\Sigma_{\tilde{\mathcal{R}}}}_{+}(\xi))_{
i_{1}i_{2}} \rvert \underset{\underset{z_{o}=1+o(1)}{\mathscr{N},n \to 
\infty}}{=} \mathcal{O} \left(\tilde{\mathfrak{c}}_{\tilde{\mathcal{R}},3}^{
\sharp}(\infty) \me^{-((n-1)K+k) \tilde{\lambda}_{\tilde{\mathcal{R}},2}^{
\sharp}(+)} \right),
\end{equation*}
and
\begin{align*}
\lvert \lvert w^{\Sigma_{\tilde{\mathcal{R}}}}_{+}(\pmb{\cdot}) \rvert \rvert_{
\mathcal{L}^{q_{1}}_{\mathrm{M}_{2}(\mathbb{C})}((\tilde{a}_{N+1}+\tilde{
\delta}_{\tilde{a}_{N+1}},+\infty))} :=& \, \left(\int\limits_{\tilde{A}_{\tilde{
\mathcal{R}},2}(+) \cup \tilde{A}_{\tilde{\mathcal{R}},2}^{c}(+)} \lvert 
w^{\Sigma_{\tilde{\mathcal{R}}}}_{+}(\xi) \rvert^{q_{1}} \, \lvert \md \xi 
\rvert \right)^{\frac{1}{q_{1}}} \! = \! \left(\left(\int_{\tilde{A}_{\tilde{
\mathcal{R}},2}(+)} \! + \! \int_{\tilde{A}_{\tilde{\mathcal{R}},2}^{c}(+)} \right) 
\left(\sum_{i_{1},i_{2}=1,2} \overline{(w^{\Sigma_{\tilde{\mathcal{R}}}}_{+}
(\xi))_{i_{1}i_{2}}} \right. \right. \\
\times&\left. \left. \, (w^{\Sigma_{\tilde{\mathcal{R}}}}_{+}(\xi))_{i_{1}i_{2}} 
\vphantom{M^{M^{M^{M^{M^{M^{M}}}}}}} \right)^{\frac{q_{1}}{2}} \, \lvert 
\md \xi \rvert \right)^{\frac{1}{q_{1}}} \underset{\underset{z_{o}=1+o(1)}{
\mathscr{N},n \to \infty}}{=} \mathcal{O} \left(\dfrac{\tilde{\mathfrak{c}}_{
\tilde{\mathcal{R}},3}^{\sharp}(q_{1}) \me^{-((n-1)K+k) \tilde{\lambda}_{
\tilde{\mathcal{R}},2}^{\sharp}(+)}}{((n \! - \! 1)K \! + \! k)^{1/q_{1}}} \right), 
\quad q_{1} \! = \! 1,2,
\end{align*}
where
\begin{align} \label{eqqlamrp} 
\tilde{\lambda}_{\tilde{\mathcal{R}},2}^{\sharp}(+) :=& \min \left\lbrace 
\tilde{\lambda}_{\tilde{\mathcal{R}},2}(+) \min \left\lbrace \tilde{\delta}_{
\tilde{a}_{N+1}},\min\limits_{q^{\prime} \in \tilde{Q}_{\tilde{\mathcal{R}},
2}(+)} \lbrace \alpha_{p_{q^{\prime}}} \! - \! \tilde{a}_{N+1} \! + \! \tilde{
\delta}_{\tilde{\mathcal{R}},2}(+) \rbrace \right\rbrace,K^{-1} \lvert \ln 
(\tilde{\varepsilon}_{\infty}) \rvert \tilde{\lambda}_{\tilde{\mathcal{R}},3}
(+), \right. \nonumber \\
&\left. K^{-1} \lvert \ln (\tilde{\delta}_{\tilde{\mathcal{R}},2}(+)) \rvert 
\min\limits_{q^{\prime} \in \tilde{Q}_{\tilde{\mathcal{R}},2}(+)} \lbrace 
\tilde{\mathfrak{c}}_{\tilde{\mathcal{R}},3}(q^{\prime},+) \rbrace 
\right\rbrace,
\end{align}
with $\tilde{\lambda}_{\tilde{\mathcal{R}},2}(+)$, $\tilde{\lambda}_{
\tilde{\mathcal{R}},3}(+)$, and $\tilde{\mathfrak{c}}_{\tilde{\mathcal{R}},3}
(q^{\prime},+)$, $q^{\prime} \! \in \! \tilde{Q}_{\tilde{\mathcal{R}},2}(+)$, 
given in the corresponding item~\pmb{(2)}, 
subitem~$\boldsymbol{\mathrm{(2)}_{ii}}$ of Lemma~\ref{lem5.1}, and 
$\tilde{\mathfrak{c}}_{\tilde{\mathcal{R}},3}^{\sharp}(q)$, $q \! \in \! 
\lbrace 1,2,\infty \rbrace$, is characterised in the corresponding 
item~\pmb{(2)}, subitem~$\boldsymbol{\mathrm{(2)}_{ii}}$ of the 
proposition.\footnote{If $\tilde{Q}_{\tilde{\mathcal{R}},2}(+) \! = \! 
\varnothing$, that is, $\# \tilde{Q}_{\tilde{\mathcal{R}},2}(+) 
\! = \! 0$, then $\cup_{q \in \tilde{Q}_{\tilde{\mathcal{R}},2}(+)} 
\mathscr{O}_{\tilde{\delta}_{\tilde{\mathcal{R}},2}(+)}(\alpha_{p_{q}}) \! 
:= \! \varnothing$, in which case $\tilde{A}_{\tilde{\mathcal{R}},2}(+) 
\! = \! (\tilde{a}_{N+1} \! + \! \tilde{\delta}_{\tilde{a}_{N+1}},+\infty) 
\setminus \mathscr{O}_{\infty}(\alpha_{p_{\mathfrak{s}-1}})$, and the 
Asymptotics~\eqref{eqwtlas8} are applicable provided one defines 
$\tilde{\lambda}_{\tilde{\mathcal{R}},2}^{\sharp}(+) \! := \! \min 
\lbrace \tilde{\lambda}_{\tilde{\mathcal{R}},2}(+) \tilde{\delta}_{
\tilde{a}_{N+1}},K^{-1} \lvert \ln (\tilde{\varepsilon}_{\infty}) \rvert 
\tilde{\lambda}_{\tilde{\mathcal{R}},3}(+) \rbrace$.} For $n \! \in \! 
\mathbb{N}$ and $k \! \in \! \lbrace 1,2,\dotsc,K \rbrace$ such that 
$\alpha_{p_{\mathfrak{s}}} \! := \! \alpha_{k} \! \neq \! \infty$, the 
case $z \! \in \! (-\infty,\tilde{b}_{0} \! - \! \tilde{\delta}_{\tilde{b}_{
0}})$ $(\subset \tilde{\Sigma}_{p}^{1})$ is analysed analogously, and 
leads to the asymptotics, in the double-scaling limit $\mathscr{N},
n \! \to \! \infty$ such that $z_{o} \! = \! 1 \! + \! o(1)$, given in 
Equation~\eqref{eqwtlas7}; in fact, the analogue of $\tilde{\lambda}_{
\tilde{\mathcal{R}},2}^{\sharp}(+)$ for this latter case reads:
\begin{align} \label{eqqlamrm} 
\tilde{\lambda}_{\tilde{\mathcal{R}},2}^{\sharp}(-) :=& \min \left\lbrace 
\tilde{\lambda}_{\tilde{\mathcal{R}},2}(-) \min \left\lbrace \tilde{\delta}_{
\tilde{b}_{0}},\min\limits_{q^{\prime} \in \tilde{Q}_{\tilde{\mathcal{R}},
2}(-)} \lbrace \lvert \alpha_{p_{q^{\prime}}} \! - \! \tilde{b}_{0} \! - \! 
\tilde{\delta}_{\tilde{\mathcal{R}},2}(-) \rvert \rbrace \right\rbrace,
K^{-1} \lvert \ln (\tilde{\varepsilon}_{\infty}) \rvert \tilde{\lambda}_{
\tilde{\mathcal{R}},3}(-), \right. \nonumber \\
&\left. K^{-1} \lvert \ln (\tilde{\delta}_{\tilde{\mathcal{R}},2}(-)) \rvert 
\min\limits_{q^{\prime} \in \tilde{Q}_{\tilde{\mathcal{R}},2}(-)} \lbrace 
\tilde{\mathfrak{c}}_{\tilde{\mathcal{R}},3}(q^{\prime},-) \rbrace 
\right\rbrace,
\end{align}
where $\tilde{\lambda}_{\tilde{\mathcal{R}},2}(-)$, $\tilde{Q}_{\tilde{
\mathcal{R}},2}(-)$, $\tilde{\delta}_{\tilde{\mathcal{R}},2}(-)$, $\tilde{
\lambda}_{\tilde{\mathcal{R}},3}(-)$, and $\tilde{\mathfrak{c}}_{\tilde{
\mathcal{R}},3}(q^{\prime},-)$, $q^{\prime} \! \in \! \tilde{Q}_{\tilde{
\mathcal{R}},2}(-)$, are given in the corresponding item~\pmb{(2)}, 
subitem~$\boldsymbol{\mathrm{(2)}_{ii}}$ of 
Lemma~\ref{lem5.1}.\footnote{If $\tilde{Q}_{\tilde{\mathcal{R}},2}(-) 
\! = \! \varnothing$, that is, $\# \tilde{Q}_{\tilde{\mathcal{R}},2}(-) 
\! = \! 0$, then $\cup_{q \in \tilde{Q}_{\tilde{\mathcal{R}},2}(-)} 
\mathscr{O}_{\tilde{\delta}_{\tilde{\mathcal{R}},2}(-)}(\alpha_{p_{q}}) 
\! := \! \varnothing$, and the Asymptotics~\eqref{eqwtlas7} are 
applicable provided one defines $\tilde{\lambda}_{\tilde{\mathcal{R}},
2}^{\sharp}(-) \! := \! \min \lbrace \tilde{\lambda}_{\tilde{\mathcal{R}},
2}(-) \tilde{\delta}_{\tilde{b}_{0}},K^{-1} \lvert \ln (\tilde{\varepsilon}_{
\infty}) \rvert \tilde{\lambda}_{\tilde{\mathcal{R}},3}(-) \rbrace$.}

For $n \! \in \! \mathbb{N}$ and $k \! \in \! \lbrace 1,2,\dotsc,K \rbrace$ 
such that $\alpha_{p_{\mathfrak{s}}} \! := \! \alpha_{k} \! \neq \! \infty$, 
and for $z \! \in \! \tilde{\Sigma}_{p,j}^{3} \cup \tilde{\Sigma}_{p,j}^{4}$, 
$j \! = \! 1,2,\dotsc,N \! + \! 1$, where $\tilde{\Sigma}_{p,j}^{3} \! := 
\! \tilde{J}_{j}^{\smallfrown} \setminus (\tilde{J}_{j}^{\smallfrown} \cap 
(\tilde{\mathbb{U}}_{\tilde{\delta}_{\tilde{b}_{j-1}}} \cup \tilde{\mathbb{
U}}_{\tilde{\delta}_{\tilde{a}_{j}}}))$ $(\subset \mathbb{C}_{+})$, and 
$\tilde{\Sigma}_{p,j}^{4} \! := \! \tilde{J}_{j}^{\smallsmile} \setminus 
(\tilde{J}_{j}^{\smallsmile} \cap (\tilde{\mathbb{U}}_{\tilde{\delta}_{
\tilde{b}_{j-1}}} \cup \tilde{\mathbb{U}}_{\tilde{\delta}_{\tilde{a}_{j}}}))$ 
$(\subset \mathbb{C}_{-})$, consider, say, and without loss of 
generality, the analysis for the case $z \! \in \! \tilde{\Sigma}_{p,j}^{3}$; 
the analysis for the case $z \! \in \! \tilde{\Sigma}_{p,j}^{4}$ is, 
\emph{mutatis mutandis}, analogous. For $n \! \in \! \mathbb{N}$ 
and $k \! \in \! \lbrace 1,2,\dotsc,K \rbrace$ such that $\alpha_{
p_{\mathfrak{s}}} \! := \! \alpha_{k} \! \neq \! \infty$, via the formula 
$w^{\Sigma_{\tilde{\mathcal{R}}}}_{+}(z) \! = \! \tilde{v}_{\tilde{
\mathcal{R}}}(z) \! - \! \mathrm{I}$, the elementary trigonometric 
inequalities $\sin \theta \! \geqslant \! \tfrac{2 \theta}{\pi}$ and $\cos 
\theta \! \geqslant \! -\tfrac{2 \theta}{\pi} \! + \! 1$ for $0 \! \leqslant 
\! \theta \! \leqslant \! \tfrac{\pi}{2}$, $\sin \theta \! \geqslant \! 
\tfrac{2}{\pi}(\pi \! - \! \theta)$ and $\cos \theta \! \leqslant \! 
-\tfrac{2 \theta}{\pi} \! + \! 1$ for $\tfrac{\pi}{2} \! \leqslant \! \theta 
\! \leqslant \! \pi$, $\sin \theta \! \leqslant \! \tfrac{2 \theta}{\pi}$ 
and $\cos \theta \! \geqslant \! \tfrac{2 \theta}{\pi} \! + \! 1$ for 
$-\tfrac{\pi}{2} \! \leqslant \! \theta \! \leqslant \! 0$, and $\sin 
\theta \! \leqslant \! -\tfrac{2}{\pi}(\pi \! + \! \theta)$ and $\cos 
\theta \! \leqslant \! \tfrac{2 \theta}{\pi} \! + \! 1$ for $-\pi \! \leqslant 
\! \theta \! \leqslant \! -\tfrac{\pi}{2}$, the parametrisation for 
the---elliptic---homotopic deformation of $\tilde{\Sigma}_{p,j}^{3}$ 
given in the proof of the corresponding item of Lemma~\ref{lem5.1}, 
that is, $\tilde{\Sigma}_{p,j}^{3} \! = \! \lbrace (x_{j}(\theta),y_{j}
(\theta)); \, x_{j}(\theta) \! = \! \tfrac{1}{2}(\tilde{a}_{j} \! + \! 
\tilde{b}_{j-1}) \! + \! \tfrac{1}{2}(\tilde{a}_{j} \! - \! \tilde{b}_{j-1}) 
\cos \theta, y_{j}(\theta) \! = \! \tilde{\eta}_{j} \sin \theta, \, 
\theta_{0}^{\tilde{a}}(j) \! \leqslant \! \theta \! \leqslant \! \pi \! 
- \! \theta_{0}^{\tilde{b}}(j) \rbrace$, and the asymptotics, in the 
double-scaling limit $\mathscr{N},n \! \to \! \infty$ such that $z_{o} 
\! = \! 1 \! + \! o(1)$, of $\tilde{v}_{\tilde{\mathcal{R}}}(z)$ given in 
Equation~\eqref{eqtlvee11}, one shows, via an integration-by-parts 
argument, that, for $z \! \in \! \tilde{\Sigma}_{p,j}^{3}$, $j \! = \! 
1,2,\dotsc,N \! + \! 1$,
\begin{equation*}
\lvert \lvert w^{\Sigma_{\tilde{\mathcal{R}}}}_{+}(\pmb{\cdot}) \rvert \rvert_{
\mathcal{L}^{\infty}_{\mathrm{M}_{2}(\mathbb{C})}(\tilde{\Sigma}_{p,j}^{3})} 
\! := \! \max_{i_{1},i_{2}=1,2} \sup_{\theta \in [\theta_{0}^{\tilde{a}}(j),
\pi -\theta_{0}^{\tilde{b}}(j)]} \lvert (w^{\Sigma_{\tilde{\mathcal{R}}}}_{+}
(x_{j}(\theta) \! + \! \mi y_{j}(\theta)))_{i_{1}i_{2}} \rvert \underset{
\underset{z_{o}=1+o(1)}{\mathscr{N},n \to \infty}}{=} \mathcal{O} \left(
\tilde{\mathfrak{c}}_{\tilde{\mathcal{R}},4}^{\sharp,\smallfrown}(j;\infty) 
\me^{-((n-1)K+k) \tilde{\lambda}_{\tilde{\mathcal{R}},4}^{\sharp,
\smallfrown}(j)} \right),
\end{equation*}
and
\begin{align*}
\lvert \lvert w^{\Sigma_{\tilde{\mathcal{R}}}}_{+}(\pmb{\cdot}) \rvert \rvert_{
\mathcal{L}^{q_{1}}_{\mathrm{M}_{2}(\mathbb{C})}(\tilde{\Sigma}_{p,j}^{3})} 
:=& \, \left(\int_{\tilde{\Sigma}_{p,j}^{3}} \lvert w^{\Sigma_{\tilde{\mathcal{
R}}}}_{+}(\xi) \rvert^{q_{1}} \, \lvert \md \xi \rvert \right)^{\frac{1}{q_{1}}} \! 
= \! \left(\left(\int_{\theta_{0}^{\tilde{a}}(j)}^{\pi/2} \! + \! \int_{\pi/2}^{\pi 
- \theta_{0}^{\tilde{b}}(j)} \right) \left(\sum_{i_{1},i_{2}=1,2} \overline{(w^{
\Sigma_{\tilde{\mathcal{R}}}}_{+}(x_{j}(\theta) \! + \! \mi y_{j}(\theta)))_{
i_{1}i_{2}}} \right. \right. \\
\times&\left. \left. \, (w^{\Sigma_{\tilde{\mathcal{R}}}}_{+}(x_{j}(\theta) \! 
+ \! \mi y_{j}(\theta)))_{i_{1}i_{2}} \vphantom{M^{M^{M^{M^{M^{M^{M}}}}}}} 
\right)^{\frac{q_{1}}{2}} \left((x_{j}^{\prime}(\theta))^{2} \! + \! (y_{j}^{\prime}
(\theta))^{2} \right)^{\frac{1}{2}} \, \md \theta \right)^{\frac{1}{q_{1}}} \\
\underset{\underset{z_{o}=1+o(1)}{\mathscr{N},n \to \infty}}{=}& \, 
\mathcal{O} \left(\dfrac{\tilde{\mathfrak{c}}_{\tilde{\mathcal{R}},4}^{
\sharp,\smallfrown}(j;q_{1}) \me^{-((n-1)K+k) \tilde{\lambda}_{\tilde{
\mathcal{R}},4}^{\sharp,\smallfrown}(j)}}{((n \! - \! 1)K \! + \! k)^{1/q_{1}}} 
\right), \quad q_{1} \! = \! 1,2,
\end{align*}
where
\begin{equation} \label{eqqlamrupj} 
\tilde{\lambda}_{\tilde{\mathcal{R}},4}^{\sharp,\smallfrown}(j) \! := 
\! \min \left\lbrace \tilde{\lambda}_{\tilde{\mathcal{R}},3}(j) \tilde{
\mathfrak{x}}_{1}^{\smallfrown}(j) \min \lbrace \tilde{\mathfrak{y}}_{
1}^{\smallfrown}(j),\tilde{\mathfrak{y}}_{2}^{\smallfrown}(j) \rbrace,
\tilde{\lambda}_{\tilde{\mathcal{R}},3}(j) \tilde{\mathfrak{x}}_{2}^{
\smallfrown}(j) \min \lbrace \tilde{\mathfrak{y}}_{3}^{\smallfrown}(j),
\tilde{\mathfrak{y}}_{4}^{\smallfrown}(j) \rbrace \right\rbrace,
\end{equation}
with $\tilde{\lambda}_{\tilde{\mathcal{R}},3}(j)$ given in the corresponding 
item~\pmb{(2)}, subitem~$\boldsymbol{\mathrm{(2)}_{iii}}$ of 
Lemma~\ref{lem5.1}, and {}\footnote{All square roots are positive.}
\begin{gather*}
\tilde{\mathfrak{x}}_{1}^{\smallfrown}(j) \! := \! \left(\dfrac{1}{\pi^{2}}
((\tilde{b}_{j-1} \! - \! \tilde{a}_{j})^{2} \! + \! 4 \tilde{\eta}_{j}^{2})
(\theta_{0}^{\tilde{a}}(j))^{2} \! + \! \tilde{\delta}_{\tilde{a}_{j}}^{2} 
\right)^{1/2}, \\
\tilde{\mathfrak{y}}_{1}^{\smallfrown}(j) \! := \! \left(1 \! - \! 
\dfrac{\frac{2}{\pi} \tilde{\delta}_{\tilde{a}_{j}}((\tilde{b}_{j-1} \! - \! 
\tilde{a}_{j}) \cos \sigma_{\tilde{a}_{j}}^{+} \! + \! 2 \tilde{\eta}_{j} 
\sin \sigma_{\tilde{a}_{j}}^{+}) \theta_{0}^{\tilde{a}}(j)}{\frac{1}{\pi^{2}}
((\tilde{b}_{j-1} \! - \! \tilde{a}_{j})^{2} \! + \! 4 \tilde{\eta}_{j}^{2})
(\theta_{0}^{\tilde{a}}(j))^{2} \! + \! \tilde{\delta}_{\tilde{a}_{j}}^{2}} 
\right)^{1/2}, \\
\tilde{\mathfrak{y}}_{2}^{\smallfrown}(j) \! := \! \left(1 \! - \! 
\dfrac{\tilde{\delta}_{\tilde{a}_{j}}((\tilde{b}_{j-1} \! - \! \tilde{a}_{j}) 
\cos \sigma_{\tilde{a}_{j}}^{+} \! + \! 2 \tilde{\eta}_{j} \sin \sigma_{
\tilde{a}_{j}}^{+})}{\frac{1}{\pi^{2}}((\tilde{b}_{j-1} \! - \! \tilde{a}_{
j})^{2} \! + \! 4 \tilde{\eta}_{j}^{2})(\theta_{0}^{\tilde{a}}(j))^{2} \! + 
\! \tilde{\delta}_{\tilde{a}_{j}}^{2}} \right)^{1/2}, \\
\tilde{\mathfrak{x}}_{2}^{\smallfrown}(j) \! := \! \left(\dfrac{1}{4}
((\tilde{b}_{j-1} \! - \! \tilde{a}_{j})^{2} \! + \! 4 \tilde{\eta}_{j}^{2}) \! 
+ \! \tilde{\delta}_{\tilde{a}_{j}}^{2} \cos^{2} \sigma_{\tilde{a}_{j}}^{+} 
\! + \! (2 \tilde{\eta}_{j} \! - \! \tilde{\delta}_{\tilde{a}_{j}} \sin 
\sigma_{\tilde{a}_{j}}^{+})^{2} \right)^{1/2}, \\
\tilde{\mathfrak{y}}_{3}^{\smallfrown}(j) \! := \! \left(1 \! - \! 
\dfrac{((\tilde{b}_{j-1} \! - \! \tilde{a}_{j}) \tilde{\delta}_{\tilde{a}_{j}} 
\cos \sigma_{\tilde{a}_{j}}^{+} \! + \! 2 \tilde{\eta}_{j}(2 \tilde{\eta}_{j} 
\! - \! \tilde{\delta}_{\tilde{a}_{j}} \sin \sigma_{\tilde{a}_{j}}^{+}))}{
\frac{1}{4}((\tilde{b}_{j-1} \! - \! \tilde{a}_{j})^{2} \! + \! 4 \tilde{
\eta}_{j}^{2}) \! + \! \tilde{\delta}_{\tilde{a}_{j}}^{2} \cos^{2} 
\sigma_{\tilde{a}_{j}}^{+} \! + \! (2 \tilde{\eta}_{j} \! - \! \tilde{
\delta}_{\tilde{a}_{j}} \sin \sigma_{\tilde{a}_{j}}^{+})^{2}} \right)^{1/2}, \\
\tilde{\mathfrak{y}}_{4}^{\smallfrown}(j) \! := \! \left(1 \! - \! 
\dfrac{\frac{2}{\pi}((\tilde{b}_{j-1} \! - \! \tilde{a}_{j}) \tilde{\delta}_{
\tilde{a}_{j}} \cos \sigma_{\tilde{a}_{j}}^{+} \! + \! 2 \tilde{\eta}_{j}(2 
\tilde{\eta}_{j} \! - \! \tilde{\delta}_{\tilde{a}_{j}} \sin \sigma_{
\tilde{a}_{j}}^{+}))(\pi \! - \! \theta_{0}^{\tilde{b}}(j))}{\frac{1}{4}
((\tilde{b}_{j-1} \! - \! \tilde{a}_{j})^{2} \! + \! 4 \tilde{\eta}_{j}^{2}) 
\! + \! \tilde{\delta}_{\tilde{a}_{j}}^{2} \cos^{2} \sigma_{\tilde{a}_{
j}}^{+} \! + \! (2 \tilde{\eta}_{j} \! - \! \tilde{\delta}_{\tilde{a}_{j}} 
\sin \sigma_{\tilde{a}_{j}}^{+})^{2}} \right)^{1/2},
\end{gather*}
where $\sigma_{\tilde{a}_{j}}^{+} \! \in \! (2 \pi/3,\pi)$, and $\tilde{
\mathfrak{c}}_{\tilde{\mathcal{R}},4}^{\sharp,\smallfrown}(j;q)$, $q \! 
\in \! \lbrace 1,2,\infty \rbrace$, is characterised in the corresponding 
item~\pmb{(2)}, subitem~$\boldsymbol{\mathrm{(2)}_{iii}}$ of the 
proposition. For $n \! \in \! \mathbb{N}$ and $k \! \in \! \lbrace 
1,2,\dotsc,K \rbrace$ such that $\alpha_{p_{\mathfrak{s}}} \! := \! 
\alpha_{k} \! \neq \! \infty$, the case $z \! \in \! \tilde{\Sigma}_{p,j}^{4} 
\! = \! \lbrace (x_{j}(\theta),y_{j}(\theta)); \, x_{j}(\theta) \! = \! \tfrac{1}{2}
(\tilde{a}_{j} \! + \! \tilde{b}_{j-1}) \! + \! \tfrac{1}{2}(\tilde{a}_{j} \! - \! 
\tilde{b}_{j-1}) \cos \theta, y_{j}(\theta) \! = \! \tilde{\eta}_{j} \sin \theta, 
\, -\pi \! + \! \theta_{0}^{\tilde{b}}(j) \! \leqslant \! \theta \! \leqslant \! 
-\theta_{0}^{\tilde{a}}(j) \rbrace$, $j \! = \! 1,2,\dotsc,N \! + \! 1$, is, 
\emph{mutatis mutandis}, analysed analogously, and leads to the 
asymptotics, in the double-scaling limit $\mathscr{N},n \! \to \! \infty$ such 
that $z_{o} \! = \! 1 \! + \! o(1)$, given in Equations~\eqref{eqwtlas10}; 
in fact, the analogue of $\tilde{\lambda}_{\tilde{\mathcal{R}},4}^{\sharp,
\smallfrown}(j)$, $j \! = \! 1,2,\dotsc,N \! + \! 1$, for this latter case reads:
\begin{equation} \label{eqqlamrdwj} 
\tilde{\lambda}_{\tilde{\mathcal{R}},4}^{\sharp,\smallsmile}(j) \! := 
\! \min \left\lbrace \tilde{\lambda}_{\tilde{\mathcal{R}},4}(j) \tilde{
\mathfrak{x}}_{1}^{\smallsmile}(j) \min \lbrace \tilde{\mathfrak{y}}_{
1}^{\smallsmile}(j),\tilde{\mathfrak{y}}_{2}^{\smallsmile}(j) \rbrace,
\tilde{\lambda}_{\tilde{\mathcal{R}},4}(j) \tilde{\mathfrak{x}}_{2}^{
\smallsmile}(j) \min \lbrace \tilde{\mathfrak{y}}_{3}^{\smallsmile}(j),
\tilde{\mathfrak{y}}_{4}^{\smallsmile}(j) \rbrace \right\rbrace,
\end{equation}
where $\tilde{\lambda}_{\tilde{\mathcal{R}},4}(j)$ is given in the corresponding 
item~\pmb{(2)}, subitem~$\boldsymbol{\mathrm{(2)}_{iii}}$ of Lemma~\ref{lem5.1}, 
and
\begin{gather*}
\tilde{\mathfrak{x}}_{1}^{\smallsmile}(j) \! := \! \left(\dfrac{1}{4}
((\tilde{b}_{j-1} \! - \! \tilde{a}_{j})^{2} \! + \! 4 \tilde{\eta}_{j}^{2}) 
\! + \! \tilde{\delta}_{\tilde{a}_{j}}^{2} \right)^{1/2}, \\
\tilde{\mathfrak{y}}_{1}^{\smallsmile}(j) \! := \! \left(1 \! - \! 
\dfrac{\tilde{\delta}_{\tilde{a}_{j}}((\tilde{b}_{j-1} \! - \! \tilde{a}_{j}) 
\cos \sigma_{\tilde{a}_{j}}^{-} \! - \! 2 \tilde{\eta}_{j} \sin \sigma_{
\tilde{a}_{j}}^{-})}{\frac{1}{4}((\tilde{b}_{j-1} \! - \! \tilde{a}_{j})^{2} \! 
+ \! 4 \tilde{\eta}_{j}^{2}) \! + \! \tilde{\delta}_{\tilde{a}_{j}}^{2}} 
\right)^{1/2}, \\
\tilde{\mathfrak{y}}_{2}^{\smallsmile}(j) \! := \! \left(1 \! - \! 
\dfrac{\frac{2}{\pi} \tilde{\delta}_{\tilde{a}_{j}}((\tilde{b}_{j-1} \! - \! 
\tilde{a}_{j}) \cos \sigma_{\tilde{a}_{j}}^{-} \! - \! 2 \tilde{\eta}_{j} \sin 
\sigma_{\tilde{a}_{j}}^{-}) \theta_{0}^{\tilde{a}}(j)}{\frac{1}{4}((\tilde{
b}_{j-1} \! - \! \tilde{a}_{j})^{2} \! + \! 4 \tilde{\eta}_{j}^{2}) \! + \! 
\tilde{\delta}_{\tilde{a}_{j}}^{2}} \right)^{1/2}, \\
\tilde{\mathfrak{x}}_{2}^{\smallsmile}(j) \! := \! \left(\dfrac{1}{\pi^{2}} 
((\tilde{b}_{j-1} \! - \! \tilde{a}_{j})^{2} \! + \! 4 \tilde{\eta}_{j}^{2})
(\pi \! - \! \theta_{0}^{\tilde{b}}(j))^{2} \! + \! \tilde{\delta}_{
\tilde{a}_{j}}^{2} \cos^{2} \sigma_{\tilde{a}_{j}}^{-} \! + \! (2 
\tilde{\eta}_{j} \! + \! \tilde{\delta}_{\tilde{a}_{j}} \sin \sigma_{
\tilde{a}_{j}}^{-})^{2} \right)^{1/2}, \\
\tilde{\mathfrak{y}}_{3}^{\smallsmile}(j) \! := \! \left(1 \! - \! 
\dfrac{\frac{2}{\pi}((\tilde{b}_{j-1} \! - \! \tilde{a}_{j}) \tilde{\delta}_{
\tilde{a}_{j}} \cos \sigma_{\tilde{a}_{j}}^{-} \! + \! 2 \tilde{\eta}_{j}
(2 \tilde{\eta}_{j} \! + \! \tilde{\delta}_{\tilde{a}_{j}} \sin \sigma_{
\tilde{a}_{j}}^{-}))(\pi \! - \! \theta_{0}^{\tilde{b}}(j))}{\frac{1}{\pi^{2}}
((\tilde{b}_{j-1} \! - \! \tilde{a}_{j})^{2} \! + \! 4 \tilde{\eta}_{j}^{2})
(\pi \! - \! \theta_{0}^{\tilde{b}}(j))^{2} \! + \! \tilde{\delta}_{
\tilde{a}_{j}}^{2} \cos^{2} \sigma_{\tilde{a}_{j}}^{-} \! + \! (2 
\tilde{\eta}_{j} \! + \! \tilde{\delta}_{\tilde{a}_{j}} \sin \sigma_{
\tilde{a}_{j}}^{-})^{2}} \right)^{1/2}, \\
\tilde{\mathfrak{y}}_{4}^{\smallsmile}(j) \! := \! \left(1 \! - \! 
\dfrac{((\tilde{b}_{j-1} \! - \! \tilde{a}_{j}) \tilde{\delta}_{\tilde{a}_{j}} 
\cos \sigma_{\tilde{a}_{j}}^{-} \! + \! 2 \tilde{\eta}_{j}(2 \tilde{\eta}_{j} 
\! + \! \tilde{\delta}_{\tilde{a}_{j}} \sin \sigma_{\tilde{a}_{j}}^{-}))}{
\frac{1}{\pi^{2}}((\tilde{b}_{j-1} \! - \! \tilde{a}_{j})^{2} \! + \! 4 
\tilde{\eta}_{j}^{2})(\pi \! - \! \theta_{0}^{\tilde{b}}(j))^{2} \! + \! 
\tilde{\delta}_{\tilde{a}_{j}}^{2} \cos^{2} \sigma_{\tilde{a}_{j}}^{-} 
\! + \! (2 \tilde{\eta}_{j} \! + \! \tilde{\delta}_{\tilde{a}_{j}} \sin 
\sigma_{\tilde{a}_{j}}^{-})^{2}} \right)^{1/2},
\end{gather*}
where $\sigma_{\tilde{a}_{j}}^{-} \! \in \! (-\pi,-2 \pi/3)$.

For $n \! \in \! \mathbb{N}$ and $k \! \in \! \lbrace 1,2,\dotsc,K \rbrace$ 
such that $\alpha_{p_{\mathfrak{s}}} \! := \! \alpha_{k} \! \neq \! \infty$, 
and for $z \! \in \! \tilde{\Sigma}_{p,j}^{5} \! := \! \partial \tilde{\mathbb{
U}}_{\tilde{\delta}_{\tilde{b}_{j-1}}} \cup \partial \tilde{\mathbb{U}}_{\tilde{
\delta}_{\tilde{a}_{j}}}$, $j \! = \! 1,2,\dotsc,N \! + \! 1$, consider, say, 
and without loss of generality, the case $z \! \in \! \partial \tilde{
\mathbb{U}}_{\tilde{\delta}_{\tilde{a}_{j}}}$; the analysis for the case $z 
\! \in \! \partial \tilde{\mathbb{U}}_{\tilde{\delta}_{\tilde{b}_{j-1}}}$ is, 
\emph{mutatis mutandis}, analogous. For $n \! \in \! \mathbb{N}$ and 
$k \! \in \! \lbrace 1,2,\dotsc,K \rbrace$ such that $\alpha_{p_{
\mathfrak{s}}} \! := \! \alpha_{k} \! \neq \! \infty$, via the asymptotics, 
in the double-scaling limit $\mathscr{N},n \! \to \! \infty$ such that 
$z_{o} \! = \! 1 \! + \! o(1)$, of $w^{\Sigma_{\tilde{\mathcal{R}}}}_{+}(z)$ 
for $z \! \in \! \partial \tilde{\mathbb{U}}_{\tilde{\delta}_{\tilde{a}_{j}}}$, 
$j \! = \! 1,2,\dotsc,N \! + \! 1$, given in Equation~\eqref{eqproptila}, 
one shows, via a straightforward integration argument, that
\begin{equation*}
\lvert \lvert w^{\Sigma_{\tilde{\mathcal{R}}}}_{+}(\pmb{\cdot}) \rvert \rvert_{
\mathcal{L}^{\infty}_{\mathrm{M}_{2}(\mathbb{C})}(\partial \tilde{\mathbb{
U}}_{\tilde{\delta}_{\tilde{a}_{j}}})} \! := \! \max_{i_{1},i_{2}=1,2} \sup_{\xi 
\in \partial \tilde{\mathbb{U}}_{\tilde{\delta}_{\tilde{a}_{j}}}} \lvert 
(w^{\Sigma_{\tilde{\mathcal{R}}}}_{+}(\xi))_{i_{1}i_{2}} \rvert \underset{
\underset{z_{o}=1+o(1)}{\mathscr{N},n \to \infty}}{=} \mathcal{O} \left(
\dfrac{\tilde{\mathfrak{c}}_{\tilde{\mathcal{R}}}^{\triangleleft}(j;\infty)}{
(n \! - \! 1)K \! + \! k} \right),
\end{equation*}
and
\begin{align*}
\lvert \lvert w^{\Sigma_{\tilde{\mathcal{R}}}}_{+}(\pmb{\cdot}) \rvert \rvert_{
\mathcal{L}^{q_{1}}_{\mathrm{M}_{2}(\mathbb{C})}(\partial \tilde{\mathbb{
U}}_{\tilde{\delta}_{\tilde{a}_{j}}})} :=& \, \left(\int_{\partial \tilde{\mathbb{
U}}_{\tilde{\delta}_{\tilde{a}_{j}}}} \lvert w^{\Sigma_{\tilde{\mathcal{R}}}}_{+}
(\xi) \rvert^{q_{1}} \, \lvert \md \xi \rvert \right)^{\frac{1}{q_{1}}} \! = \! 
\left(\int_{0}^{2 \pi} \left(\sum_{i_{1},i_{2}=1,2} \overline{(w^{\Sigma_{
\tilde{\mathcal{R}}}}_{+}(\tilde{a}_{j} \! + \! \tilde{\delta}_{\tilde{a}_{j}} 
\me^{\mi \theta}))_{i_{1}i_{2}}} \right. \right. \\
\times&\left. \left. \, (w^{\Sigma_{\tilde{\mathcal{R}}}}_{+}(\tilde{a}_{j} 
\! + \! \tilde{\delta}_{\tilde{a}_{j}} \me^{\mi \theta}))_{i_{1}i_{2}} 
\vphantom{M^{M^{M^{M^{M^{M^{M}}}}}}} \right)^{\frac{q_{1}}{2}} \lvert 
\md (\tilde{a}_{j} \! + \! \tilde{\delta}_{\tilde{a}_{j}} \me^{\mi \theta}) 
\rvert \right)^{\frac{1}{q_{1}}} \underset{\underset{z_{o}=1+o(1)}{
\mathscr{N},n \to \infty}}{=} \mathcal{O} \left(\dfrac{\tilde{\mathfrak{c}}_{
\tilde{\mathcal{R}}}^{\triangleleft}(j;q_{1})}{(n \! - \! 1)K \! + \! k} \right), 
\quad q_{1} \! = \! 1,2,
\end{align*}
where $\tilde{\mathfrak{c}}_{\tilde{\mathcal{R}}}^{\triangleleft}(j;q)$, $q \! 
\in \! \lbrace 1,2,\infty \rbrace$, is characterised in the corresponding 
item~\pmb{(2)}, subitem~$\boldsymbol{\mathrm{(2)}_{iv}}$ of the 
proposition. For $n \! \in \! \mathbb{N}$ and $k \! \in \! \lbrace 1,2,
\dotsc,K \rbrace$ such that $\alpha_{p_{\mathfrak{s}}} \! := \! \alpha_{k} 
\! \neq \! \infty$, the case $z \! \in \! \partial \tilde{\mathbb{U}}_{
\tilde{\delta}_{\tilde{b}_{j-1}}}$, $j \! = \! 1,2,\dotsc,N \! + \! 1$, is, 
\emph{mutatis mutandis}, analysed analogously, and leads to the asymptotics, 
in the double-scaling limit $\mathscr{N},n \! \to \! \infty$ such that $z_{o} 
\! = \! 1 \! + \! o(1)$, given in Equations~\eqref{eqwtlasbt}. \hfill $\qed$
\begin{ccccc} \label{lem5.3} 
For $n \! \in \! \mathbb{N}$ and $k \! \in \! \lbrace 1,2,\dotsc,K \rbrace$ 
such that $\alpha_{p_{\mathfrak{s}}} \! := \! \alpha_{k} \! = \! \infty$ 
(resp., $\alpha_{p_{\mathfrak{s}}} \! := \! \alpha_{k} \! \neq \! \infty)$, 
let $\hat{\mathcal{R}} \colon \mathbb{C} \setminus \hat{\Sigma}_{
\hat{\mathcal{R}}}^{\sharp} \! \to \! \mathrm{SL}_{2}(\mathbb{C})$ (resp., 
$\tilde{\mathcal{R}} \colon \mathbb{C} \setminus \tilde{\Sigma}_{
\tilde{\mathcal{R}}}^{\sharp} \! \to \! \mathrm{SL}_{2}(\mathbb{C}))$ solve 
the equivalent {\rm RHP} $(\hat{\mathcal{R}}(z),\hat{v}_{\hat{\mathcal{R}}}
(z),\hat{\Sigma}_{\hat{\mathcal{R}}}^{\sharp})$ (resp., $(\tilde{\mathcal{R}}
(z),\tilde{v}_{\tilde{\mathcal{R}}}(z),\tilde{\Sigma}_{\tilde{\mathcal{R}}}^{
\sharp}))$, with associated {\rm BC} operator $C_{w^{\Sigma_{\hat{
\mathcal{R}}}}}$ (resp., $C_{w^{\Sigma_{\tilde{\mathcal{R}}}}})$. Then$:$ 
{\rm \pmb{(1)}} for $n \! \in \! \mathbb{N}$ and $k \! \in \! \lbrace 1,2,
\dotsc,K \rbrace$ such that $\alpha_{p_{\mathfrak{s}}} \! := \! \alpha_{k} 
\! = \! \infty$,
\begin{equation} \label{eqcwsight} 
\lvert \lvert C_{w^{\Sigma_{\hat{\mathcal{R}}}}} \rvert \rvert_{\mathfrak{B}_{
\infty}(\hat{\Sigma}_{\hat{\mathcal{R}}}^{\sharp})} \underset{\underset{z_{o}
=1+o(1)}{\mathscr{N},n \to \infty}}{=} \mathcal{O} \left(\hat{\mathfrak{c}}_{
\hat{\mathcal{R}},\hat{w}}^{\triangleright} \me^{-((n-1)K+k) \hat{\lambda}_{
\hat{\mathcal{R}},\hat{w}}^{\triangleright}} \right),
\end{equation}
where
\begin{equation} \label{hatlamrwr1} 
\hat{\lambda}_{\hat{\mathcal{R}},\hat{w}}^{\triangleright} \! := \! \min 
\left\lbrace \min_{j=1,2,\dotsc,N} \lbrace \hat{\lambda}_{\hat{\mathcal{R}},
1}^{\sharp}(j) \rbrace,\hat{\lambda}_{\hat{\mathcal{R}},2}^{\sharp}(+),
\hat{\lambda}_{\hat{\mathcal{R}},2}^{\sharp}(-),\min \left\lbrace \min_{j=1,2,
\dotsc,N+1} \lbrace \hat{\lambda}_{\hat{\mathcal{R}},4}^{\sharp,\smallfrown}
(j) \rbrace,\min_{j=1,2,\dotsc,N+1} \lbrace \hat{\lambda}_{\hat{\mathcal{R}},
4}^{\sharp,\smallsmile}(j) \rbrace \right\rbrace \right\rbrace,
\end{equation}
with $\hat{\lambda}_{\hat{\mathcal{R}},\hat{w}}^{\triangleright} \! = \! 
\hat{\lambda}_{\hat{\mathcal{R}},\hat{w}}^{\triangleright}(n,k,z_{o}) \! 
=_{\underset{z_{o}=1+o(1)}{\mathscr{N},n \to \infty}} \! \mathcal{O}(1)$ 
and $> \! 0$,\footnote{Note that $\hat{\lambda}_{\hat{\mathcal{R}},1}^{
\sharp}(j)$, $\hat{\lambda}_{\hat{\mathcal{R}},2}^{\sharp}(\pm)$, 
$\hat{\lambda}_{\hat{\mathcal{R}},4}^{\sharp,\smallfrown}(j)$, and 
$\hat{\lambda}_{\hat{\mathcal{R}},4}^{\sharp,\smallsmile}(j)$ are 
defined in item~\pmb{(1)} of Proposition~\ref{propo5.2}.} and 
$\hat{\mathfrak{c}}_{\hat{\mathcal{R}},\hat{w}}^{\triangleright} \! = \! 
\hat{\mathfrak{c}}_{\hat{\mathcal{R}},\hat{w}}^{\triangleright}(n,k,z_{o}) \! 
=_{\underset{z_{o}=1+o(1)}{\mathscr{N},n \to \infty}} \! \mathcal{O}(1)$, 
and
\begin{equation} \label{eqcwsightt} 
\lvert \lvert C_{w^{\Sigma_{\hat{\mathcal{R}}}}} \rvert \rvert_{\mathfrak{B}_{2}
(\hat{\Sigma}_{\hat{\mathcal{R}}}^{\sharp})} \underset{\underset{z_{o}=1+
o(1)}{\mathscr{N},n \to \infty}}{=} \mathcal{O} \left(\hat{\mathfrak{c}}_{
\hat{\mathcal{R}},\hat{w}}^{\triangleleft} \me^{-\frac{1}{2}((n-1)K+k) 
\hat{\lambda}_{\hat{\mathcal{R}},\hat{w}}^{\triangleright}} \right),
\end{equation}
where $\hat{\mathfrak{c}}_{\hat{\mathcal{R}},\hat{w}}^{\triangleleft} \! = \! 
\hat{\mathfrak{c}}_{\hat{\mathcal{R}},\hat{w}}^{\triangleleft}(n,k,z_{o}) \! 
=_{\underset{z_{o}=1+o(1)}{\mathscr{N},n \to \infty}} \! \mathcal{O}(1)$$;$ 
and {\rm \pmb{(2)}} for $n \! \in \! \mathbb{N}$ and $k \! \in \! \lbrace 1,2,
\dotsc,K \rbrace$ such that $\alpha_{p_{\mathfrak{s}}} \! := \! \alpha_{k} 
\! \neq \! \infty$,
\begin{equation} \label{eqcwsigtl} 
\lvert \lvert C_{w^{\Sigma_{\tilde{\mathcal{R}}}}} \rvert \rvert_{\mathfrak{B}_{
\infty}(\tilde{\Sigma}_{\tilde{\mathcal{R}}}^{\sharp})} \underset{\underset{z_{o}
=1+o(1)}{\mathscr{N},n \to \infty}}{=} \mathcal{O} \left(\tilde{\mathfrak{c}}_{
\tilde{\mathcal{R}},\tilde{w}}^{\triangleright} \me^{-((n-1)K+k) \tilde{
\lambda}_{\tilde{\mathcal{R}},\tilde{w}}^{\triangleright}} \right),
\end{equation}
where
\begin{equation} \label{tillamrwr2} 
\tilde{\lambda}_{\tilde{\mathcal{R}},\tilde{w}}^{\triangleright} \! := \! \min 
\left\lbrace \min_{j=1,2,\dotsc,N} \lbrace \tilde{\lambda}_{\tilde{\mathcal{
R}},1}^{\sharp}(j) \rbrace,\tilde{\lambda}_{\tilde{\mathcal{R}},2}^{\sharp}(+),
\tilde{\lambda}_{\tilde{\mathcal{R}},2}^{\sharp}(-),\min \left\lbrace \min_{j
=1,2,\dotsc,N+1} \lbrace \tilde{\lambda}_{\tilde{\mathcal{R}},4}^{\sharp,
\smallfrown}(j) \rbrace,\min_{j=1,2,\dotsc,N+1} \lbrace \tilde{\lambda}_{
\tilde{\mathcal{R}},4}^{\sharp,\smallsmile}(j) \rbrace \right\rbrace 
\right\rbrace,
\end{equation}
with $\tilde{\lambda}_{\tilde{\mathcal{R}},\tilde{w}}^{\triangleright} \! = \! 
\tilde{\lambda}_{\tilde{\mathcal{R}},\tilde{w}}^{\triangleright}(n,k,z_{o}) \! 
=_{\underset{z_{o}=1+o(1)}{\mathscr{N},n \to \infty}} \! \mathcal{O}(1)$ 
and $> \! 0$,\footnote{Note that $\tilde{\lambda}_{\tilde{\mathcal{R}},1}^{
\sharp}(j)$, $\tilde{\lambda}_{\tilde{\mathcal{R}},2}^{\sharp}(\pm)$, $\tilde{
\lambda}_{\tilde{\mathcal{R}},4}^{\sharp,\smallfrown}(j)$, and $\tilde{
\lambda}_{\tilde{\mathcal{R}},4}^{\sharp,\smallsmile}(j)$ are 
defined in item~\pmb{(2)} of Proposition~\ref{propo5.2}.} and 
$\tilde{\mathfrak{c}}_{\tilde{\mathcal{R}},\tilde{w}}^{\triangleright} \! = \! 
\tilde{\mathfrak{c}}_{\tilde{\mathcal{R}},\tilde{w}}^{\triangleright}(n,k,z_{o}) 
\! =_{\underset{z_{o}=1+o(1)}{\mathscr{N},n \to \infty}} \! \mathcal{O}(1)$, 
and
\begin{equation} \label{eqcwsigtll} 
\lvert \lvert C_{w^{\Sigma_{\tilde{\mathcal{R}}}}} \rvert \rvert_{\mathfrak{B}_{
2}(\tilde{\Sigma}_{\tilde{\mathcal{R}}}^{\sharp})} \underset{\underset{z_{o}=
1+o(1)}{\mathscr{N},n \to \infty}}{=} \mathcal{O} \left(\tilde{\mathfrak{c}}_{
\tilde{\mathcal{R}},\tilde{w}}^{\triangleleft} \me^{-\frac{1}{2}((n-1)K+k) 
\tilde{\lambda}_{\tilde{\mathcal{R}},\tilde{w}}^{\triangleright}} \right),
\end{equation}
where $\tilde{\mathfrak{c}}_{\tilde{\mathcal{R}},\tilde{w}}^{\triangleleft} \! = 
\! \tilde{\mathfrak{c}}_{\tilde{\mathcal{R}},\tilde{w}}^{\triangleleft}(n,k,z_{o}) 
\! =_{\underset{z_{o}=1+o(1)}{\mathscr{N},n \to \infty}} \! \mathcal{O}(1)$.

For $n \! \in \! \mathbb{N}$ and $k \! \in \! \lbrace 1,2,\dotsc,K \rbrace$ 
such that $\alpha_{p_{\mathfrak{s}}} \! := \! \alpha_{k} \! = \! \infty$ 
(resp., $\alpha_{p_{\mathfrak{s}}} \! := \! \alpha_{k} \! \neq \! \infty)$, 
$(\id \! - \! C_{w^{\Sigma_{\hat{\mathcal{R}}}}})^{-1} \! \! \upharpoonright_{
\mathcal{L}^{2}_{\mathrm{M}_{2}(\mathbb{C})}(\hat{\Sigma}_{\hat{\mathcal{
R}}}^{\sharp})}$ (resp., $(\id \! - \! C_{w^{\Sigma_{\tilde{\mathcal{R}}}}})^{-1} 
\! \! \upharpoonright_{\mathcal{L}^{2}_{\mathrm{M}_{2}(\mathbb{C})}
(\tilde{\Sigma}_{\tilde{\mathcal{R}}}^{\sharp})})$ exists; in particular, 
$\lvert \lvert (\id \! - \! C_{w^{\Sigma_{\hat{\mathcal{R}}}}})^{-1} \rvert 
\rvert_{\mathfrak{B}_{2}(\hat{\Sigma}_{\hat{\mathcal{R}}}^{\sharp})} \! 
=_{\underset{z_{o}=1+o(1)}{\mathscr{N},n \to \infty}} \! \mathcal{O}(1)$ 
(resp., $\lvert \lvert (\id \! - \! C_{w^{\Sigma_{\tilde{\mathcal{R}}}}})^{-1} 
\rvert \rvert_{\mathfrak{B}_{2}(\tilde{\Sigma}_{\tilde{\mathcal{R}}}^{\sharp})} 
=_{\underset{z_{o}=1+o(1)}{\mathscr{N},n \to \infty}} \! \mathcal{O}(1))$.
\end{ccccc}

\emph{Proof}. The proof of this Lemma~\ref{lem5.3} consists of two cases: 
(i) $n \! \in \! \mathbb{N}$ and $k \! \in \! \lbrace 1,2,\dotsc,K \rbrace$ 
such that $\alpha_{p_{\mathfrak{s}}} \! := \! \alpha_{k} \! = \! \infty$; and 
(ii) $n \! \in \! \mathbb{N}$ and $k \! \in \! \lbrace 1,2,\dotsc,K \rbrace$ 
such that $\alpha_{p_{\mathfrak{s}}} \! := \! \alpha_{k} \! \neq \! \infty$. 
Notwithstanding the fact that the scheme of the proof is, \emph{mutatis mutandis}, 
similar for both cases, without loss of generality, only the proof for case~(ii) is 
presented in detail, whilst case~(i) is proved analogously.

For $n \! \in \! \mathbb{N}$ and $k \! \in \! \lbrace 1,2,\dotsc,K \rbrace$ 
such that $\alpha_{p_{\mathfrak{s}}} \! := \! \alpha_{k} \! \neq \! \infty$, 
let $\tilde{\mathcal{R}} \colon \mathbb{C} \setminus \tilde{\Sigma}_{
\tilde{\mathcal{R}}}^{\sharp} \! \to \! \mathrm{SL}_{2}(\mathbb{C})$ solve 
the equivalent {\rm RHP} $(\tilde{\mathcal{R}}(z),\tilde{v}_{\tilde{\mathcal{R}}}
(z),\tilde{\Sigma}_{\tilde{\mathcal{R}}}^{\sharp})$, where the action of the 
associated {\rm BC} operator, $C_{w^{\Sigma_{\tilde{\mathcal{R}}}}}$, is 
defined by
\begin{equation*}
\mathcal{L}^{2}_{\mathrm{M}_{2}(\mathbb{C})}(\tilde{\Sigma}_{\tilde{
\mathcal{R}}}^{\sharp}) \! \ni \! f \! \mapsto \! C_{w^{\Sigma_{\tilde{
\mathcal{R}}}}}f \! := \! \lim_{-\tilde{\Sigma}_{\tilde{\mathcal{R}}}^{
\sharp} \ni z^{\prime} \to z} \int_{\tilde{\Sigma}_{\tilde{\mathcal{R}}}^{
\sharp}} \dfrac{(z^{\prime} \! - \! \alpha_{k})(fw_{+}^{\Sigma_{\tilde{
\mathcal{R}}}})(\xi)}{(\xi \! - \! \alpha_{k})(\xi \! - \! z^{\prime})} \, 
\dfrac{\md \xi}{2 \pi \mi},
\end{equation*}
where $-\tilde{\Sigma}_{\tilde{\mathcal{R}}}^{\sharp}$ denotes the 
`minus side' of $\tilde{\Sigma}_{\tilde{\mathcal{R}}}^{\sharp}$, 
$w_{+}^{\Sigma_{\tilde{\mathcal{R}}}}(z) \! = \! \tilde{v}_{\tilde{\mathcal{
R}}}(z) \! - \! \mathrm{I}$, and $\tilde{\Sigma}_{\tilde{\mathcal{R}}}^{
\sharp}$ is the disjoint union $\tilde{\Sigma}_{\tilde{\mathcal{R}}}^{\sharp} 
\! = \! \tilde{\Sigma}_{p}^{1} \cup \cup_{i=1}^{N} \tilde{\Sigma}_{p,i}^{2} 
\cup \cup_{j=1}^{N+1}(\tilde{\Sigma}_{p,j}^{3} \cup \tilde{\Sigma}_{p,j}^{4}) 
\cup \cup_{m=1}^{N+1} \tilde{\Sigma}_{p,m}^{5}$, with $\tilde{\Sigma}_{
p}^{1} \! := \! (-\infty,\tilde{b}_{0} \! - \! \tilde{\delta}_{\tilde{b}_{0}}) 
\cup (\tilde{a}_{N+1} \! + \! \tilde{\delta}_{\tilde{a}_{N+1}},+\infty)$, 
$\tilde{\Sigma}_{p,i}^{2} \! := \! (\tilde{a}_{i} \! + \! \tilde{\delta}_{
\tilde{a}_{i}},\tilde{b}_{i} \! - \! \tilde{\delta}_{\tilde{b}_{i}})$, $\tilde{
\Sigma}_{p,j}^{3} \! := \! \tilde{J}_{j}^{\smallfrown} \setminus 
(\tilde{J}_{j}^{\smallfrown} \cap (\tilde{\mathbb{U}}_{\tilde{\delta}_{
\tilde{b}_{j-1}}} \cup \tilde{\mathbb{U}}_{\tilde{\delta}_{\tilde{a}_{j}}}))$, 
$\tilde{\Sigma}_{p,j}^{4} \! := \! \tilde{J}_{j}^{\smallsmile} \setminus 
(\tilde{J}_{j}^{\smallsmile} \cap (\tilde{\mathbb{U}}_{\tilde{\delta}_{
\tilde{b}_{j-1}}} \cup \tilde{\mathbb{U}}_{\tilde{\delta}_{\tilde{a}_{j}}}))$, 
and $\tilde{\Sigma}_{p,m}^{5} \! := \! \partial \tilde{\mathbb{U}}_{
\tilde{\delta}_{\tilde{b}_{m-1}}} \cup \partial \tilde{\mathbb{U}}_{\tilde{
\delta}_{\tilde{a}_{m}}}$; hence,
\begin{align*}
\lvert \lvert (C_{w^{\Sigma_{\tilde{\mathcal{R}}}}}f)(\pmb{\cdot}) 
\rvert \rvert_{\mathcal{L}^{\infty}_{\mathrm{M}_{2}(\mathbb{C})}
(\tilde{\Sigma}_{\tilde{\mathcal{R}}}^{\sharp})} :=& \, \max_{i_{1},i_{2}
=1,2} \sup_{z \in \tilde{\Sigma}_{\tilde{\mathcal{R}}}^{\sharp}} \lvert 
((C_{w^{\Sigma_{\tilde{\mathcal{R}}}}}f)(z))_{i_{1}i_{2}} \rvert \! = \! 
\max_{i_{1},i_{2}=1,2} \sup_{z \in \tilde{\Sigma}_{\tilde{\mathcal{R}}}^{
\sharp}} \left\lvert \lim_{-\tilde{\Sigma}_{\tilde{\mathcal{R}}}^{\sharp} 
\ni z^{\prime} \to z} \int_{\tilde{\Sigma}_{\tilde{\mathcal{R}}}^{\sharp}} 
\dfrac{(z^{\prime} \! - \! \alpha_{k})((fw_{+}^{\Sigma_{\tilde{\mathcal{
R}}}})(\xi))_{i_{1}i_{2}}}{(\xi \! - \! \alpha_{k})(\xi \! - \! z^{\prime})} \, 
\dfrac{\md \xi}{2 \pi \mi} \right\rvert \\
\leqslant& \, \lvert \lvert f(\pmb{\cdot}) \rvert \rvert_{\mathcal{L}^{
\infty}_{\mathrm{M}_{2}(\mathbb{C})}(\tilde{\Sigma}_{\tilde{\mathcal{
R}}}^{\sharp})} \max_{i_{1},i_{2}=1,2} \sup_{z \in \tilde{\Sigma}_{\tilde{
\mathcal{R}}}^{\sharp}} \left\lvert \lim_{-\tilde{\Sigma}_{\tilde{
\mathcal{R}}}^{\sharp} \ni z^{\prime} \to z} \left(\int_{\tilde{\Sigma}_{
p}^{1}} \! + \! \sum_{i=1}^{N} \int_{\tilde{\Sigma}_{p,i}^{2}} \! + \! 
\sum_{j=1}^{N+1} \left(\int_{\tilde{\Sigma}_{p,j}^{3}} \! + \! \int_{
\tilde{\Sigma}_{p,j}^{4}} \right) \right. \right. \\
+&\left. \left. \, \sum_{m=1}^{N+1} \int_{\tilde{\Sigma}_{p,m}^{5}} 
\right) \dfrac{(z^{\prime} \! - \! \alpha_{k})(w_{+}^{\Sigma_{\tilde{
\mathcal{R}}}}(\xi))_{i_{1}i_{2}}}{(\xi \! - \! \alpha_{k})(\xi \! - \! 
z^{\prime})} \, \dfrac{\md \xi}{2 \pi \mi} \right\rvert \! \leqslant \! 
\lvert \lvert f(\pmb{\cdot}) \rvert \rvert_{\mathcal{L}^{\infty}_{
\mathrm{M}_{2}(\mathbb{C})}(\tilde{\Sigma}_{\tilde{\mathcal{R}}}^{
\sharp})} \max_{i_{1},i_{2}=1,2} \sum_{m=1}^{5} \tilde{\mathbb{I}}_{
m}^{\sharp}(i_{1},i_{2}),
\end{align*}
where
\begin{align}
\tilde{\mathbb{I}}_{1}^{\sharp}(i_{1},i_{2}) \! :=& \sup_{z \in \tilde{
\Sigma}_{\tilde{\mathcal{R}}}^{\sharp}} \left\lvert \sum_{j=1}^{N} 
\lim_{-\tilde{\Sigma}_{\tilde{\mathcal{R}}}^{\sharp} \ni z^{\prime} 
\to z} \int_{\tilde{a}_{j}+\tilde{\delta}_{\tilde{a}_{j}}}^{\tilde{b}_{j}-
\tilde{\delta}_{\tilde{b}_{j}}} \dfrac{(z^{\prime} \! - \! \alpha_{k})
(w_{+}^{\Sigma_{\tilde{\mathcal{R}}}}(\xi))_{i_{1}i_{2}}}{(\xi \! - \! 
\alpha_{k})(\xi \! - \! z^{\prime})} \, \dfrac{\md \xi}{2 \pi \mi} 
\right\rvert, \label{eqiy1} \\
\tilde{\mathbb{I}}_{2}^{\sharp}(i_{1},i_{2}) \! :=& \sup_{z \in \tilde{
\Sigma}_{\tilde{\mathcal{R}}}^{\sharp}} \left\lvert \lim_{-\tilde{
\Sigma}_{\tilde{\mathcal{R}}}^{\sharp} \ni z^{\prime} \to z} \int_{
\tilde{a}_{N+1}+\tilde{\delta}_{\tilde{a}_{N+1}}}^{+\infty} \dfrac{
(z^{\prime} \! - \! \alpha_{k})(w_{+}^{\Sigma_{\tilde{\mathcal{R}}}}
(\xi))_{i_{1}i_{2}}}{(\xi \! - \! \alpha_{k})(\xi \! - \! z^{\prime})} \, 
\dfrac{\md \xi}{2 \pi \mi} \right\rvert, \label{eqiy2} \\
\tilde{\mathbb{I}}_{3}^{\sharp}(i_{1},i_{2}) \! :=& \sup_{z \in \tilde{
\Sigma}_{\tilde{\mathcal{R}}}^{\sharp}} \left\lvert \lim_{-\tilde{
\Sigma}_{\tilde{\mathcal{R}}}^{\sharp} \ni z^{\prime} \to z} 
\int_{-\infty}^{\tilde{b}_{0}-\tilde{\delta}_{\tilde{b}_{0}}} \dfrac{
(z^{\prime} \! - \! \alpha_{k})(w_{+}^{\Sigma_{\tilde{\mathcal{R}}}}
(\xi))_{i_{1}i_{2}}}{(\xi \! - \! \alpha_{k})(\xi \! - \! z^{\prime})} \, 
\dfrac{\md \xi}{2 \pi \mi} \right\rvert, \label{eqiy3} \\
\tilde{\mathbb{I}}_{4}^{\sharp}(i_{1},i_{2}) \! :=& \sup_{z \in \tilde{
\Sigma}_{\tilde{\mathcal{R}}}^{\sharp}} \left\lvert \sum_{j=1}^{N+1} 
\lim_{-\tilde{\Sigma}_{\tilde{\mathcal{R}}}^{\sharp} \ni z^{\prime} 
\to z} \left(\int_{\tilde{\Sigma}_{p,j}^{3}} \! + \! \int_{\tilde{
\Sigma}_{p,j}^{4}} \right) \dfrac{(z^{\prime} \! - \! \alpha_{k})
(w_{+}^{\Sigma_{\tilde{\mathcal{R}}}}(\xi))_{i_{1}i_{2}}}{(\xi \! - \! 
\alpha_{k})(\xi \! - \! z^{\prime})} \, \dfrac{\md \xi}{2 \pi \mi} 
\right\rvert, \label{eqiy4} \\
\tilde{\mathbb{I}}_{5}^{\sharp}(i_{1},i_{2}) \! :=& \sup_{z \in \tilde{
\Sigma}_{\tilde{\mathcal{R}}}^{\sharp}} \left\lvert \sum_{j=1}^{N+1} 
\lim_{-\tilde{\Sigma}_{\tilde{\mathcal{R}}}^{\sharp} \ni z^{\prime} 
\to z} \left(\int_{\partial \tilde{\mathbb{U}}_{\tilde{\delta}_{\tilde{b}_{
j-1}}}} \! + \! \int_{\partial \tilde{\mathbb{U}}_{\tilde{\delta}_{\tilde{
a}_{j}}}} \right) \dfrac{(z^{\prime} \! - \! \alpha_{k})
(w_{+}^{\Sigma_{\tilde{\mathcal{R}}}}(\xi))_{i_{1}i_{2}}}{(\xi \! - \! 
\alpha_{k})(\xi \! - \! z^{\prime})} \, \dfrac{\md \xi}{2 \pi \mi} 
\right\rvert. \label{eqiy5}
\end{align}

For $n \! \in \! \mathbb{N}$ and $k \! \in \! \lbrace 1,2,\dotsc,K \rbrace$ 
such that $\alpha_{p_{\mathfrak{s}}} \! := \! \alpha_{k} \! \neq \! \infty$, 
in order to estimate, asymptotically, in the double-scaling limit 
$\mathscr{N},n \! \to \! \infty$ such that $z_{o} \! = \! 1 \! + \! o(1)$, 
$\tilde{\mathbb{I}}_{1}^{\sharp}(i_{1},i_{2})$, $i_{1},i_{2} \! = \! 1,2$, 
consider, say, for $j \! \in \! \lbrace 1,2,\dotsc,N \rbrace$, the quantity
\begin{equation*}
\tilde{\mathbb{I}}_{1,j}^{\sharp}(i_{1},i_{2}) \! := \! \sup_{z \in \tilde{
\Sigma}_{\tilde{\mathcal{R}}}^{\sharp}} \left\lvert \lim_{-\tilde{\Sigma}_{
\tilde{\mathcal{R}}}^{\sharp} \ni z^{\prime} \to z} \int_{\tilde{a}_{j}+
\tilde{\delta}_{\tilde{a}_{j}}}^{\tilde{b}_{j}-\tilde{\delta}_{\tilde{b}_{j}}} 
\dfrac{(z^{\prime} \! - \! \alpha_{k})(w_{+}^{\Sigma_{\tilde{\mathcal{R}}}}
(\xi))_{i_{1}i_{2}}}{(\xi \! - \! \alpha_{k})(\xi \! - \! z^{\prime})} \, 
\dfrac{\md \xi}{2 \pi \mi} \right\rvert:
\end{equation*}
via the distributional identities $(x \! - \! (x_{0} \! \pm \! \mi 0))^{-1} 
\! = \! (x \! - \! x_{0})^{-1} \! \pm \! \mi \pi \delta (x \! - \! x_{0})$, where 
$\delta (\pmb{\cdot})$ is the Dirac delta function, and $\int_{x_{1}}^{x_{2}}
f(\xi) \delta (\xi \! - \! x) \, \md \xi \! = \! 
\left\{
\begin{smallmatrix}
f(x), \, x \in (x_{1},x_{2}), \\
0, \, x \in \mathbb{R} \setminus (x_{1},x_{2}),
\end{smallmatrix}
\right.$ a partial-fraction decomposition argument, and an elementary 
inequality argument, one shows that, for $j \! \in \! \lbrace 1,2,\dotsc,N 
\rbrace$,
\begin{equation}
\tilde{\mathbb{I}}_{1,j}^{\sharp}(i_{1},i_{2}) \! \leqslant \! \tilde{\mathbb{
I}}_{1,j}^{\sharp,\mathrm{A}}(i_{1},i_{2}) \! + \! \tilde{\mathbb{I}}_{1,j}^{
\sharp,\mathrm{B}}(i_{1},i_{2}) \! + \! \tilde{\mathbb{I}}_{1,j}^{\sharp,
\mathrm{C}}(i_{1},i_{2}), \quad i_{1},i_{2} \! = \! 1,2, \label{eqiy6}
\end{equation} 
where
\begin{gather*}
\tilde{\mathbb{I}}_{1,j}^{\sharp,\mathrm{A}}(i_{1},i_{2}) \! := \! \sup 
\left\lvert \int_{\tilde{a}_{j}+\tilde{\delta}_{\tilde{a}_{j}}}^{\tilde{b}_{j}
-\tilde{\delta}_{\tilde{b}_{j}}} \dfrac{(w_{+}^{\Sigma_{\tilde{\mathcal{R}}}}
(\xi))_{i_{1}i_{2}}}{\xi \! - \! \alpha_{k}} \, \dfrac{\md \xi}{2 \pi \mi} 
\right\rvert, \quad \quad \tilde{\mathbb{I}}_{1,j}^{\sharp,\mathrm{B}}
(i_{1},i_{2}) \! := \! \sup_{z \in \tilde{\Sigma}_{\tilde{\mathcal{R}}}^{\sharp} 
\setminus (\tilde{a}_{j}+\tilde{\delta}_{\tilde{a}_{j}},\tilde{b}_{j}-\tilde{
\delta}_{\tilde{b}_{j}})} \left\lvert \int_{\tilde{a}_{j}+\tilde{\delta}_{
\tilde{a}_{j}}}^{\tilde{b}_{j}-\tilde{\delta}_{\tilde{b}_{j}}} \dfrac{(w_{+}^{
\Sigma_{\tilde{\mathcal{R}}}}(\xi))_{i_{1}i_{2}}}{\xi \! - \! z} \, 
\dfrac{\md \xi}{2 \pi \mi} \right\rvert, \\
\tilde{\mathbb{I}}_{1,j}^{\sharp,\mathrm{C}}(i_{1},i_{2}) \! := \! \sup_{z 
\in (\tilde{a}_{j}+\tilde{\delta}_{\tilde{a}_{j}},\tilde{b}_{j}-\tilde{\delta}_{
\tilde{b}_{j}})} \left\lvert \mathrm{P.V.}
\int_{\tilde{a}_{j}+\tilde{\delta}_{\tilde{a}_{j}}}^{
\tilde{b}_{j}-\tilde{\delta}_{\tilde{b}_{j}}} \dfrac{(w_{+}^{\Sigma_{\tilde{
\mathcal{R}}}}(\xi))_{i_{1}i_{2}}}{\xi \! - \! z} \, \dfrac{\md \xi}{2 \pi \mi} 
\! - \! \dfrac{1}{2}(w_{+}^{\Sigma_{\tilde{\mathcal{R}}}}(z))_{i_{1}i_{2}} 
\right\rvert,
\end{gather*}
with $\mathrm{P.V.}$ denoting the principle-value integral. For $n \! 
\in \! \mathbb{N}$ and $k \! \in \! \lbrace 1,2,\dotsc,K \rbrace$ such 
that $\alpha_{p_{\mathfrak{s}}} \! := \! \alpha_{k} \! \neq \! \infty$, 
recall {}from the proof of the corresponding item of Lemma~\ref{lem5.1} 
that, for $j \! \in \! \lbrace 1,2,\dotsc,N \rbrace$, $\tilde{\Sigma}_{p,j}^{2} 
\! := \! (\tilde{a}_{j} \! + \! \tilde{\delta}_{\tilde{a}_{j}},\tilde{b}_{j} \! - \! 
\tilde{\delta}_{\tilde{b}_{j}}) \! = \! \tilde{A}_{\tilde{\mathcal{R}},1}(j) \cup 
\tilde{A}_{\tilde{\mathcal{R}},1}^{c}(j)$, where $\tilde{A}_{\tilde{\mathcal{R}},
1}(j) \! = \! (\tilde{a}_{j} \! + \! \tilde{\delta}_{\tilde{a}_{j}},\tilde{b}_{j} \! - \! 
\tilde{\delta}_{\tilde{b}_{j}}) \setminus \cup_{q \in \tilde{Q}_{\tilde{\mathcal{
R}},1}(j)} \mathscr{O}_{\tilde{\delta}_{\tilde{\mathcal{R}},1}(j)}(\alpha_{p_{q}})$, 
and $\tilde{A}_{\tilde{\mathcal{R}},1}^{c}(j) \! = \! \cup_{q \in \tilde{Q}_{
\tilde{\mathcal{R}},1}(j)} \mathscr{O}_{\tilde{\delta}_{\tilde{\mathcal{R}},1}(j)}
(\alpha_{p_{q}})$, with $\tilde{Q}_{\tilde{\mathcal{R}},1}(j) \! := \! \lbrace 
\mathstrut q^{\prime} \! \in \! \lbrace 1,\dotsc,\mathfrak{s} \! - \! 2,\mathfrak{s} 
\rbrace; \, \alpha_{p_{q^{\prime}}} \! \in \! (\tilde{a}_{j} \! + \! \tilde{\delta}_{
\tilde{a}_{j}},\tilde{b}_{j} \! - \! \tilde{\delta}_{\tilde{b}_{j}}) \rbrace$, and 
sufficiently small $\tilde{\delta}_{\tilde{\mathcal{R}},1}(j) \! > \! 0$ chosen so 
that $\mathscr{O}_{\tilde{\delta}_{\tilde{\mathcal{R}},1}(j)}(\alpha_{p_{q_{1}}}) 
\cap \mathscr{O}_{\tilde{\delta}_{\tilde{\mathcal{R}},1}(j)}(\alpha_{p_{q_{2}}}) \! 
= \! \varnothing$ $\forall$ $q_{1} \! \neq \! q_{2} \! \in \! \tilde{Q}_{\tilde{
\mathcal{R}},1}(j)$ and $\mathscr{O}_{\tilde{\delta}_{\tilde{\mathcal{R}},1}(j)}
(\alpha_{p_{q_{1}}}) \cap \lbrace \tilde{a}_{j} \! + \! \tilde{\delta}_{\tilde{a}_{j}} 
\rbrace \! = \! \varnothing \! = \! \mathscr{O}_{\tilde{\delta}_{\tilde{\mathcal{R}},
1}(j)}(\alpha_{p_{q_{1}}}) \cap \lbrace \tilde{b}_{j} \! - \! \tilde{\delta}_{\tilde{b}_{j}} 
\rbrace$ (of course, $\tilde{A}_{\tilde{\mathcal{R}},1}(j) \cap \tilde{A}_{\tilde{
\mathcal{R}},1}^{c}(j) \! = \! \varnothing$). Via the formula $w^{\Sigma_{\tilde{
\mathcal{R}}}}_{+}(z) \! = \! \tilde{v}_{\tilde{\mathcal{R}}}(z) \! - \! \mathrm{I}$, 
the asymptotics, in the double-scaling limit $\mathscr{N},n \! \to \! \infty$ 
such that $z_{o} \! = \! 1 \! + \! o(1)$, of $\tilde{v}_{\tilde{\mathcal{R}}}(z)$ 
given in Equation~\eqref{eqtlvee8}, and the inequality $\ln \lvert x \rvert \! 
\leqslant \! \lvert x \rvert \! - \! 1$, one shows, via a tedious integration 
argument, that, for $j \! \in \! \lbrace 1,2,\dotsc,N \rbrace$,
\begin{align*}
\tilde{\mathbb{I}}_{1,j}^{\sharp,\mathrm{A}}(i_{1},i_{2}) =& \, \sup \left\lvert 
\left(\int_{\tilde{A}_{\tilde{\mathcal{R}},1}(j)} \! + \! \int_{\tilde{A}_{\tilde{
\mathcal{R}},1}^{c}(j)} \right) \dfrac{(w_{+}^{\Sigma_{\tilde{\mathcal{R}}}}
(\xi))_{i_{1}i_{2}}}{\xi \! - \! \alpha_{k}} \, \dfrac{\md \xi}{2 \pi \mi} \right\rvert 
\! \underset{\underset{z_{o}=1+o(1)}{\mathscr{N},n \to \infty}}{\leqslant} 
\mathcal{O} \left(\dfrac{\tilde{\mathfrak{c}}_{i_{1}i_{2}}^{\sharp,\mathrm{A}}
(n,k,z_{o};j) \me^{-((n-1)K+k) \tilde{\lambda}_{\tilde{\mathcal{R}},1}^{
\sharp}(j)}}{(n \! - \! 1)K \! + \! k} \right), \\
\tilde{\mathbb{I}}_{1,j}^{\sharp,\mathrm{B}}(i_{1},i_{2}) =& \, \sup_{z \in 
\tilde{\Sigma}_{\tilde{\mathcal{R}}}^{\sharp} \setminus (\tilde{A}_{\tilde{
\mathcal{R}},1}(j) \cup \tilde{A}_{\tilde{\mathcal{R}},1}^{c}(j))} \left\lvert 
\left(\int_{\tilde{A}_{\tilde{\mathcal{R}},1}(j)} \! + \! \int_{\tilde{A}_{\tilde{
\mathcal{R}},1}^{c}(j)} \right) \dfrac{(w_{+}^{\Sigma_{\tilde{\mathcal{R}}}}
(\xi))_{i_{1}i_{2}}}{\xi \! - \! z} \, \dfrac{\md \xi}{2 \pi \mi} \right\rvert \! 
\underset{\underset{z_{o}=1+o(1)}{\mathscr{N},n \to \infty}}{\leqslant} 
\mathcal{O} \left(\dfrac{\tilde{\mathfrak{c}}_{i_{1}i_{2}}^{\sharp,\mathrm{B}}
(n,k,z_{o};j) \me^{-((n-1)K+k) \tilde{\lambda}_{\tilde{\mathcal{R}},1}^{
\sharp}(j)}}{(n \! - \! 1)K \! + \! k} \right), \\
\tilde{\mathbb{I}}_{1,j}^{\sharp,\mathrm{C}}(i_{1},i_{2}) =& \, \sup_{z \in 
\tilde{A}_{\tilde{\mathcal{R}},1}(j) \cup \tilde{A}_{\tilde{\mathcal{R}},1}^{c}(j)} 
\left\lvert \mathrm{P.V.} \left(\int_{\tilde{A}_{\tilde{\mathcal{R}},1}(j)} \! 
+ \! \int_{\tilde{A}_{\tilde{\mathcal{R}},1}^{c}(j)} \right) \dfrac{(w_{+}^{
\Sigma_{\tilde{\mathcal{R}}}}(\xi))_{i_{1}i_{2}}}{\xi \! - \! z} \, \dfrac{\md 
\xi}{2 \pi \mi} \! - \! \dfrac{1}{2}(w_{+}^{\Sigma_{\tilde{\mathcal{R}}}}
(z))_{i_{1}i_{2}} \right\rvert \\
& \, \underset{\underset{z_{o}=1+o(1)}{\mathscr{N},n \to \infty}}{\leqslant} 
\mathcal{O} \left(\tilde{\mathfrak{c}}_{i_{1}i_{2}}^{\sharp,\mathrm{C}}
(n,k,z_{o};j) \me^{-((n-1)K+k) \tilde{\lambda}_{\tilde{\mathcal{R}},1}^{
\sharp}(j)} \right),
\end{align*}
where $\tilde{\lambda}_{\tilde{\mathcal{R}},1}^{\sharp}(j)$ $(> \! 0)$ is 
defined by Equation~\eqref{eqqlamrj}, and $\tilde{\mathfrak{c}}_{i_{1}
i_{2}}^{\sharp,r}(n,k,z_{o};j) \! =_{\underset{z_{o}=1+o(1)}{\mathscr{N},
n \to \infty}} \! \mathcal{O}(1)$, $r \! \in \! \lbrace \mathrm{A},
\mathrm{B},\mathrm{C} \rbrace$, whence, for $i_{1},i_{2} \! = \! 1,2$,
\begin{equation*}
\tilde{\mathbb{I}}_{1,j}^{\sharp}(i_{1},i_{2}) \underset{\underset{z_{o}=
1+o(1)}{\mathscr{N},n \to \infty}}{\leqslant} \mathcal{O} \left(\tilde{
\mathfrak{c}}_{i_{1}i_{2}}^{\sharp,\mathrm{C}}(n,k,z_{o};j) 
\me^{-((n-1)K+k) \tilde{\lambda}_{\tilde{\mathcal{R}},1}^{\sharp}(j)} 
\right), \quad j \! \in \! \lbrace 1,2,\dotsc,N \rbrace:
\end{equation*}
via the latter---asymptotic---estimate, Definition~\eqref{eqiy1}, and 
an elementary inequality argument, one arrives at, for $n \! \in \! 
\mathbb{N}$ and $k \! \in \! \lbrace 1,2,\dotsc,K \rbrace$ such that 
$\alpha_{p_{\mathfrak{s}}} \! := \! \alpha_{k} \! \neq \! \infty$,
\begin{equation} \label{eqiy7} 
\tilde{\mathbb{I}}_{1}^{\sharp}(i_{1},i_{2}) \underset{\underset{z_{o}=1
+o(1)}{\mathscr{N},n \to \infty}}{=} \mathcal{O} \left(\tilde{\mathfrak{
c}}_{i_{1}i_{2}}^{\sharp,1}(n,k,z_{o}) \me^{-((n-1)K+k) \min\limits_{j=1,
2,\dotsc,N} \lbrace \tilde{\lambda}_{\tilde{\mathcal{R}},1}^{\sharp}(j) 
\rbrace} \right), \quad i_{1},i_{2} \! = \! 1,2,
\end{equation}
where $\tilde{\mathfrak{c}}_{i_{1}i_{2}}^{\sharp,1}(n,k,z_{o}) \! =_{
\underset{z_{o}=1+o(1)}{\mathscr{N},n \to \infty}} \! \mathcal{O}(1)$.

For $n \! \in \! \mathbb{N}$ and $k \! \in \! \lbrace 1,2,\dotsc,K \rbrace$ 
such that $\alpha_{p_{\mathfrak{s}}} \! := \! \alpha_{k} \! \neq \! \infty$, 
in order to estimate, asymptotically, in the double-scaling limit 
$\mathscr{N},n \! \to \! \infty$ such that $z_{o} \! = \! 1 \! + \! o(1)$, 
$\tilde{\mathbb{I}}_{2}^{\sharp}(i_{1},i_{2})$, $i_{1},i_{2} \! = \! 1,2$, one uses 
the distributional identities given above, a partial-fraction decomposition 
argument, and an elementary inequality argument, to show that
\begin{equation} \label{eqiy8}
\tilde{\mathbb{I}}_{2}^{\sharp}(i_{1},i_{2}) \! \leqslant \! \tilde{\mathbb{
I}}_{2}^{\sharp,\mathrm{A}}(i_{1},i_{2}) \! + \! \tilde{\mathbb{I}}_{2}^{
\sharp,\mathrm{B}}(i_{1},i_{2}) \! + \! \tilde{\mathbb{I}}_{2}^{\sharp,
\mathrm{C}}(i_{1},i_{2}), \quad i_{1},i_{2} \! = \! 1,2,
\end{equation} 
where
\begin{gather*}
\tilde{\mathbb{I}}_{2}^{\sharp,\mathrm{A}}(i_{1},i_{2}) \! := \! \sup \left\lvert 
\int_{\tilde{a}_{N+1}+\tilde{\delta}_{\tilde{a}_{N+1}}}^{+\infty} \dfrac{(w_{+}^{
\Sigma_{\tilde{\mathcal{R}}}}(\xi))_{i_{1}i_{2}}}{\xi \! - \! \alpha_{k}} \, \dfrac{
\md \xi}{2 \pi \mi} \right\rvert, \quad \quad \tilde{\mathbb{I}}_{2}^{\sharp,
\mathrm{B}}(i_{1},i_{2}) \! := \! \sup_{z \in \tilde{\Sigma}_{\tilde{\mathcal{R}}}^{
\sharp} \setminus (\tilde{a}_{N+1}+\tilde{\delta}_{\tilde{a}_{N+1}},+\infty)} 
\left\lvert \int_{\tilde{a}_{N+1}+\tilde{\delta}_{\tilde{a}_{N+1}}}^{+\infty} 
\dfrac{(w_{+}^{\Sigma_{\tilde{\mathcal{R}}}}(\xi))_{i_{1}i_{2}}}{\xi \! - \! z} \, 
\dfrac{\md \xi}{2 \pi \mi} \right\rvert, \\
\tilde{\mathbb{I}}_{2}^{\sharp,\mathrm{C}}(i_{1},i_{2}) \! := \! \sup_{z 
\in (\tilde{a}_{N+1}+\tilde{\delta}_{\tilde{a}_{N+1}},+\infty)} \left\lvert 
\mathrm{P.V.} \int_{\tilde{a}_{N+1}+\tilde{\delta}_{\tilde{a}_{N+1}}}^{+\infty} 
\dfrac{(w_{+}^{\Sigma_{\tilde{\mathcal{R}}}}(\xi))_{i_{1}i_{2}}}{\xi \! - \! z} 
\, \dfrac{\md \xi}{2 \pi \mi} \! - \! \dfrac{1}{2}(w_{+}^{\Sigma_{\tilde{
\mathcal{R}}}}(z))_{i_{1}i_{2}} \right\rvert.
\end{gather*}
For $n \! \in \! \mathbb{N}$ and $k \! \in \! \lbrace 1,2,\dotsc,K \rbrace$ 
such that $\alpha_{p_{\mathfrak{s}}} \! := \! \alpha_{k} \! \neq \! \infty$, 
recall {}from the proof of the corresponding item of Lemma~\ref{lem5.1} 
that $(\tilde{a}_{N+1} \! + \! \tilde{\delta}_{\tilde{a}_{N+1}},+\infty) \! = 
\! \tilde{A}_{\tilde{\mathcal{R}},2}(+) \cup \tilde{A}_{\tilde{\mathcal{R}},
2}^{c}(+)$, where $\tilde{A}_{\tilde{\mathcal{R}},2}(+) \! = \! (\tilde{a}_{N+1} 
\! + \! \tilde{\delta}_{\tilde{a}_{N+1}},+\infty) \setminus (\mathscr{O}_{\infty}
(\alpha_{p_{\mathfrak{s}-1}}) \cup \cup_{q \in \tilde{Q}_{\tilde{\mathcal{R}},
2}(+)} \mathscr{O}_{\tilde{\delta}_{\tilde{\mathcal{R}},2}(+)}(\alpha_{p_{q}}))$, 
and $\tilde{A}_{\tilde{\mathcal{R}},2}^{c}(+) \! = \! \mathscr{O}_{\infty}
(\alpha_{p_{\mathfrak{s}-1}}) \cup \cup_{q \in \tilde{Q}_{\tilde{\mathcal{
R}},2}(+)} \mathscr{O}_{\tilde{\delta}_{\tilde{\mathcal{R}},2}(+)}
(\alpha_{p_{q}})$, with $\tilde{Q}_{\tilde{\mathcal{R}},2}(+) \! := \! \lbrace 
\mathstrut q^{\prime} \! \in \! \lbrace 1,\dotsc,\mathfrak{s} \! - \! 2,
\mathfrak{s} \rbrace; \, \alpha_{p_{q^{\prime}}} \! \in \! (\tilde{a}_{
N+1} \! + \! \tilde{\delta}_{\tilde{a}_{N+1}},+\infty) \rbrace$ and 
sufficiently small $\tilde{\varepsilon}_{\infty},\tilde{\delta}_{\tilde{
\mathcal{R}},2}(+) \! > \! 0$ chosen so that $\mathscr{O}_{\tilde{
\delta}_{\tilde{\mathcal{R}},2}(+)}(\alpha_{p_{q^{\prime}_{1}}}) \cap 
\mathscr{O}_{\tilde{\delta}_{\tilde{\mathcal{R}},2}(+)}(\alpha_{p_{q^{\prime 
\prime}_{1}}}) \! = \! \varnothing$ $\forall$ $q^{\prime}_{1} \! \neq \! 
q^{\prime \prime}_{1} \! \in \! \tilde{Q}_{\tilde{\mathcal{R}},2}(+)$, 
$\mathscr{O}_{\tilde{\delta}_{\tilde{\mathcal{R}},2}(+)}(\alpha_{p_{
q^{\prime}_{1}}}) \cap \lbrace \tilde{a}_{N+1} \! + \! \tilde{\delta}_{
\tilde{a}_{N+1}} \rbrace \! = \! \varnothing$, and $\mathscr{O}_{
\tilde{\delta}_{\tilde{\mathcal{R}},2}(+)}(\alpha_{p_{q^{\prime}_{1}}}) \cap 
\mathscr{O}_{\infty}(\alpha_{p_{\mathfrak{s}-1}}) \! = \! \varnothing$ 
(of course, $\tilde{A}_{\tilde{\mathcal{R}},2}(+) \cap \tilde{A}_{\tilde{
\mathcal{R}},2}^{c}(+) \! = \! \varnothing$). Via the formula $w^{\Sigma_{
\tilde{\mathcal{R}}}}_{+}(z) \! = \! \tilde{v}_{\tilde{\mathcal{R}}}(z) \! - \! 
\mathrm{I}$, the asymptotics, in the double-scaling limit $\mathscr{N},
n \! \to \! \infty$ such that $z_{o} \! = \! 1 \! + \! o(1)$, of $\tilde{v}_{
\tilde{\mathcal{R}}}(z)$ given in Equation~\eqref{eqtlvee9}, and the 
inequality $\ln \lvert x \rvert \! \leqslant \! \lvert x \rvert \! - \! 1$, one 
shows, via a tedious integration argument, that
\begin{align*}
\tilde{\mathbb{I}}_{2}^{\sharp,\mathrm{A}}(i_{1},i_{2}) =& \, \sup \left\lvert 
\left(\int_{\tilde{A}_{\tilde{\mathcal{R}},2}(+)} \! + \! \int_{\tilde{A}_{\tilde{
\mathcal{R}},2}^{c}(+)} \right) \dfrac{(w_{+}^{\Sigma_{\tilde{\mathcal{R}}}}
(\xi))_{i_{1}i_{2}}}{\xi \! - \! \alpha_{k}} \, \dfrac{\md \xi}{2 \pi \mi} \right\rvert 
\! \underset{\underset{z_{o}=1+o(1)}{\mathscr{N},n \to \infty}}{\leqslant} 
\mathcal{O} \left(\dfrac{\tilde{\mathfrak{c}}_{i_{1}i_{2}}^{\sharp^{\prime},
\mathrm{A}}(n,k,z_{o}) \me^{-((n-1)K+k) \tilde{\lambda}_{\tilde{\mathcal{
R}},2}^{\sharp}(+)}}{(n \! - \! 1)K \! + \! k} \right), \\
\tilde{\mathbb{I}}_{2}^{\sharp,\mathrm{B}}(i_{1},i_{2}) =& \, \sup_{z \in 
\tilde{\Sigma}_{\tilde{\mathcal{R}}}^{\sharp} \setminus (\tilde{A}_{\tilde{
\mathcal{R}},2}(+) \cup \tilde{A}_{\tilde{\mathcal{R}},2}^{c}(+))} \left\lvert 
\left(\int_{\tilde{A}_{\tilde{\mathcal{R}},2}(+)} \! + \! \int_{\tilde{A}_{\tilde{
\mathcal{R}},2}^{c}(+)} \right) \dfrac{(w_{+}^{\Sigma_{\tilde{\mathcal{R}}}}
(\xi))_{i_{1}i_{2}}}{\xi \! - \! z} \, \dfrac{\md \xi}{2 \pi \mi} \right\rvert \! 
\underset{\underset{z_{o}=1+o(1)}{\mathscr{N},n \to \infty}}{\leqslant} 
\mathcal{O} \left(\dfrac{\tilde{\mathfrak{c}}_{i_{1}i_{2}}^{\sharp^{\prime},
\mathrm{B}}(n,k,z_{o}) \me^{-((n-1)K+k) \tilde{\lambda}_{\tilde{\mathcal{
R}},2}^{\sharp}(+)}}{(n \! - \! 1)K \! + \! k} \right), \\
\tilde{\mathbb{I}}_{2}^{\sharp,\mathrm{C}}(i_{1},i_{2}) =& \, \sup_{z \in 
\tilde{A}_{\tilde{\mathcal{R}},2}(+) \cup \tilde{A}_{\tilde{\mathcal{R}},2}^{c}
(+)} \left\lvert \mathrm{P.V.} \left(\int_{\tilde{A}_{\tilde{\mathcal{R}},2}(+)} 
\! + \! \int_{\tilde{A}_{\tilde{\mathcal{R}},2}^{c}(+)} \right) \dfrac{(w_{+}^{
\Sigma_{\tilde{\mathcal{R}}}}(\xi))_{i_{1}i_{2}}}{\xi \! - \! z} \, \dfrac{\md 
\xi}{2 \pi \mi} \! - \! \dfrac{1}{2}(w_{+}^{\Sigma_{\tilde{\mathcal{R}}}}
(z))_{i_{1}i_{2}} \right\rvert \\
& \, \underset{\underset{z_{o}=1+o(1)}{\mathscr{N},n \to \infty}}{\leqslant} 
\mathcal{O} \left(\tilde{\mathfrak{c}}_{i_{1}i_{2}}^{\sharp^{\prime},\mathrm{C}}
(n,k,z_{o}) \me^{-((n-1)K+k) \tilde{\lambda}_{\tilde{\mathcal{R}},2}^{\sharp}
(+)} \right),
\end{align*}
where $\tilde{\lambda}_{\tilde{\mathcal{R}},2}^{\sharp}(+)$ $(> \! 0)$ is 
defined by Equation~\eqref{eqqlamrp}, and $\tilde{\mathfrak{c}}_{i_{1}
i_{2}}^{\sharp^{\prime},r}(n,k,z_{o}) \! =_{\underset{z_{o}=1+o(1)}{
\mathscr{N},n \to \infty}} \! \mathcal{O}(1)$, $r \! \in \! \lbrace \mathrm{A},
\mathrm{B},\mathrm{C} \rbrace$, whence, for $n \! \in \! \mathbb{N}$ and 
$k \! \in \! \lbrace 1,2,\dotsc,K \rbrace$ such that $\alpha_{p_{\mathfrak{s}}} 
\! := \! \alpha_{k} \! \neq \! \infty$,
\begin{equation} \label{eqiy9} 
\tilde{\mathbb{I}}_{2}^{\sharp}(i_{1},i_{2}) \underset{\underset{z_{o}=1+
o(1)}{\mathscr{N},n \to \infty}}{=} \mathcal{O} \left(\tilde{\mathfrak{c}}_{
i_{1}i_{2}}^{\sharp,2}(n,k,z_{o}) \me^{-((n-1)K+k) \tilde{\lambda}_{\tilde{
\mathcal{R}},2}^{\sharp}(+)} \right), \quad i_{1},i_{2} \! = \! 1,2,
\end{equation}
where $\tilde{\mathfrak{c}}_{i_{1}i_{2}}^{\sharp,2}(n,k,z_{o}) \! 
=_{\underset{z_{o}=1+o(1)}{\mathscr{N},n \to \infty}} \! \mathcal{O}(1)$. 
For $n \! \in \! \mathbb{N}$ and $k \! \in \! \lbrace 1,2,\dotsc,K \rbrace$ 
such that $\alpha_{p_{\mathfrak{s}}} \! := \! \alpha_{k} \! \neq \! \infty$, 
the integral $\tilde{\mathbb{I}}_{3}^{\sharp}(i_{1},i_{2})$, $i_{1},i_{2} \! = \! 
1,2$, is analysed analogously (as above for $\tilde{\mathbb{I}}_{2}^{\sharp}
(i_{1},i_{2}))$, and one arrives at
\begin{equation} \label{eqiy10} 
\tilde{\mathbb{I}}_{3}^{\sharp}(i_{1},i_{2}) \underset{\underset{z_{o}=1+
o(1)}{\mathscr{N},n \to \infty}}{=} \mathcal{O} \left(\tilde{\mathfrak{c}}_{
i_{1}i_{2}}^{\sharp,3}(n,k,z_{o}) \me^{-((n-1)K+k) \tilde{\lambda}_{\tilde{
\mathcal{R}},2}^{\sharp}(-)} \right), \quad i_{1},i_{2} \! = \! 1,2,
\end{equation}
where $\tilde{\lambda}_{\tilde{\mathcal{R}},2}^{\sharp}(-)$ $(> \! 0)$ is 
defined by Equation~\eqref{eqqlamrm}, and $\tilde{\mathfrak{c}}_{i_{1}
i_{2}}^{\sharp,3}(n,k,z_{o}) \! =_{\underset{z_{o}=1+o(1)}{\mathscr{N},n 
\to \infty}} \! \mathcal{O}(1)$.

For $n \! \in \! \mathbb{N}$ and $k \! \in \! \lbrace 1,2,\dotsc,K \rbrace$ 
such that $\alpha_{p_{\mathfrak{s}}} \! := \! \alpha_{k} \! \neq \! \infty$, 
in order to estimate, asymptotically, in the double-scaling limit 
$\mathscr{N},n \! \to \! \infty$ such that $z_{o} \! = \! 1 \! + \! o(1)$, 
$\tilde{\mathbb{I}}_{4}^{\sharp}(i_{1},i_{2})$, $i_{1},i_{2} \! = \! 1,2$, consider, 
say, for $j \! \in \! \lbrace 1,2,\dotsc,N \! + \! 1 \rbrace$, the quantity
\begin{equation*}
\tilde{\mathbb{I}}_{4,j}^{\sharp,\smallfrown}(i_{1},i_{2}) \! := \! \sup_{z 
\in \tilde{\Sigma}_{\tilde{\mathcal{R}}}^{\sharp}} \left\lvert \lim_{-\tilde{
\Sigma}_{\tilde{\mathcal{R}}}^{\sharp} \ni z^{\prime} \to z} \int_{\tilde{
\Sigma}_{p,j}^{3}} \dfrac{(z^{\prime} \! - \! \alpha_{k})(w_{+}^{\Sigma_{
\tilde{\mathcal{R}}}}(\xi))_{i_{1}i_{2}}}{(\xi \! - \! \alpha_{k})(\xi \! - \! 
z^{\prime})} \, \dfrac{\md \xi}{2 \pi \mi} \right\rvert:
\end{equation*}
via the distributional identities given above, a partial-fraction decomposition 
argument, and an elementary inequality argument, one shows that, for 
$j \! \in \! \lbrace 1,2,\dotsc,N \! + \! 1 \rbrace$,
\begin{equation}
\tilde{\mathbb{I}}_{4,j}^{\sharp,\smallfrown}(i_{1},i_{2}) \! \leqslant \! 
\tilde{\mathbb{I}}_{4,j}^{\sharp,\smallfrown,\mathrm{A}}(i_{1},i_{2}) \! + \! 
\tilde{\mathbb{I}}_{4,j}^{\sharp,\smallfrown,\mathrm{B}}(i_{1},i_{2}) \! + \! 
\tilde{\mathbb{I}}_{4,j}^{\sharp,\smallfrown,\mathrm{C}}(i_{1},i_{2}), \quad 
i_{1},i_{2} \! = \! 1,2, \label{eqiy11}
\end{equation} 
where
\begin{gather*}
\tilde{\mathbb{I}}_{4,j}^{\sharp,\smallfrown,\mathrm{A}}(i_{1},i_{2}) \! := \! 
\sup \left\lvert \int_{\tilde{\Sigma}_{p,j}^{3}} \dfrac{(w_{+}^{\Sigma_{\tilde{
\mathcal{R}}}}(\xi))_{i_{1}i_{2}}}{\xi \! - \! \alpha_{k}} \, \dfrac{\md \xi}{2 \pi 
\mi} \right\rvert, \quad \quad \tilde{\mathbb{I}}_{4,j}^{\sharp,\smallfrown,
\mathrm{B}}(i_{1},i_{2}) \! := \! \sup_{z \in \tilde{\Sigma}_{\tilde{\mathcal{
R}}}^{\sharp} \setminus \tilde{\Sigma}_{p,j}^{3}} \left\lvert \int_{\tilde{
\Sigma}_{p,j}^{3}} \dfrac{(w_{+}^{\Sigma_{\tilde{\mathcal{R}}}}(\xi))_{i_{1}
i_{2}}}{\xi \! - \! z} \, \dfrac{\md \xi}{2 \pi \mi} \right\rvert, \\
\tilde{\mathbb{I}}_{4,j}^{\sharp,\smallfrown,\mathrm{C}}(i_{1},i_{2}) \! := \! 
\sup_{z \in \tilde{\Sigma}_{p,j}^{3}} \left\lvert \mathrm{P.V.} \int_{\tilde{
\Sigma}_{p,j}^{3}} \dfrac{(w_{+}^{\Sigma_{\tilde{\mathcal{R}}}}(\xi))_{i_{1}
i_{2}}}{\xi \! - \! z} \, \dfrac{\md \xi}{2 \pi \mi} \! - \! \dfrac{1}{2}(w_{+}^{
\Sigma_{\tilde{\mathcal{R}}}}(z))_{i_{1}i_{2}} \right\rvert.
\end{gather*}
Via the formula $w^{\Sigma_{\tilde{\mathcal{R}}}}_{+}(z) \! = \! \tilde{v}_{
\tilde{\mathcal{R}}}(z) \! - \! \mathrm{I}$, the asymptotics, in the 
double-scaling limit $\mathscr{N},n \! \to \! \infty$ such that $z_{o} 
\! = \! 1 \! + \! o(1)$, of $\tilde{v}_{\tilde{\mathcal{R}}}(z)$ given in 
Equation~\eqref{eqtlvee11}, the elementary trigonometric inequalities 
$\sin \theta \! \geqslant \! \tfrac{2 \theta}{\pi}$ and $\cos \theta \! 
\geqslant \! -\tfrac{2 \theta}{\pi} \! + \! 1$ for $0 \! \leqslant \! \theta \! 
\leqslant \! \tfrac{\pi}{2}$, and $\sin \theta \! \geqslant \! \tfrac{2}{\pi}
(\pi \! - \! \theta)$ and $\cos \theta \! \leqslant \! -\tfrac{2 \theta}{\pi} 
\! + \! 1$ for $\tfrac{\pi}{2} \! \leqslant \! \theta \! \leqslant \! \pi$, the 
fact that $\inf \lbrace \mathstrut \lvert \xi \! - \! \alpha_{k} \rvert; \, \xi 
\! \in \! \tilde{\Sigma}_{p,j}^{3} \rbrace \! >\! 0$ (since $\lbrace \alpha_{k} 
\rbrace \cap \tilde{\Sigma}_{p,j}^{3} \! = \! \varnothing)$ and $\inf \lbrace 
\mathstrut \lvert \xi \! - \! z \rvert; \, \xi \! \in \! \tilde{\Sigma}_{p,j}^{3}, 
\, z \! \in \! \tilde{\Sigma}_{\tilde{\mathcal{R}}}^{\sharp} \setminus 
\tilde{\Sigma}_{p,j}^{3} \rbrace \! >\! 0$, and the parametrisation for 
the---elliptic---homotopic deformation of $\tilde{\Sigma}_{p,j}^{3}$ 
given in the proof of the corresponding item of Lemma~\ref{lem5.1}, 
that is, $\tilde{\Sigma}_{p,j}^{3} \! = \! \lbrace \mathstrut (x_{j}(\theta),
y_{j}(\theta)); \, x_{j}(\theta) \! = \! \tfrac{1}{2}(\tilde{a}_{j} \! + \! 
\tilde{b}_{j-1}) \! + \! \tfrac{1}{2}(\tilde{a}_{j} \! - \! \tilde{b}_{j-1}) 
\cos \theta, \, y_{j}(\theta) \! = \! \tilde{\eta}_{j} \sin \theta, \, 
\theta_{0}^{\tilde{a}}(j) \! \leqslant \! \theta \! \leqslant \! \pi \! 
- \! \theta_{0}^{\tilde{b}}(j) \rbrace$, one shows, via a tedious 
integration-by-parts argument, that, for $n \! \in \! \mathbb{N}$ and 
$k \! \in \! \lbrace 1,2,\dotsc,K \rbrace$ such that $\alpha_{p_{
\mathfrak{s}}} \! := \! \alpha_{k} \! \neq \! \infty$, $i_{1},i_{2} \! = \! 
1,2$, and $j \! \in \! \lbrace 1,2,\dotsc,N \! + \! 1 \rbrace$,
\begin{align}
\tilde{\mathbb{I}}_{4,j}^{\sharp,\smallfrown,\mathrm{A}}(i_{1},i_{2}) =& \, 
\sup \left(\left(\int_{\theta_{0}^{\tilde{a}}(j)}^{\pi/2} \! + \! \int_{\pi/2}^{\pi 
-\theta_{0}^{\tilde{b}}(j)} \right) \left\lvert \dfrac{(w_{+}^{\Sigma_{\tilde{
\mathcal{R}}}}(x_{j}(\theta) \! + \! \mi y_{j}(\theta)))_{i_{1}i_{2}}}{x_{j}(\theta) 
\! + \! \mi y_{j}(\theta) \! - \! \alpha_{k}} \right\rvert ((x_{j}^{\prime}
(\theta))^{2} \! + \! (y_{j}^{\prime}(\theta))^{2})^{1/2} \, \dfrac{\md 
\theta}{2 \pi} \right) \nonumber \\
& \, \underset{\underset{z_{o}=1+o(1)}{\mathscr{N},n \to \infty}}{\leqslant} 
\mathcal{O} \left(\dfrac{\tilde{\mathfrak{c}}_{i_{1}i_{2}}^{\sharp,
\smallfrown,\mathrm{A}}(n,k,z_{o};j) \me^{-((n-1)K+k) \tilde{\lambda}_{
\tilde{\mathcal{R}},4}^{\sharp,\smallfrown}(j)}}{(n \! - \! 1)K \! + \! k} 
\right), \label{eqiy12} \\
\tilde{\mathbb{I}}_{4,j}^{\sharp,\smallfrown,\mathrm{B}}(i_{1},i_{2}) =& \, 
\sup_{z \in \tilde{\Sigma}_{\tilde{\mathcal{R}}}^{\sharp} \setminus \tilde{
\Sigma}_{p,j}^{3}} \left(\left(\int_{\theta_{0}^{\tilde{a}}(j)}^{\pi/2} \! + \! 
\int_{\pi/2}^{\pi -\theta_{0}^{\tilde{b}}(j)} \right) \left\lvert \dfrac{(w_{+}^{
\Sigma_{\tilde{\mathcal{R}}}}(x_{j}(\theta) \! + \! \mi y_{j}(\theta)))_{i_{1}
i_{2}}}{x_{j}(\theta) \! + \! \mi y_{j}(\theta) \! - \! z} \right\rvert ((x_{j}^{
\prime}(\theta))^{2} \! + \! (y_{j}^{\prime}(\theta))^{2})^{1/2} \, 
\dfrac{\md \theta}{2 \pi} \right) \nonumber \\
& \, \underset{\underset{z_{o}=1+o(1)}{\mathscr{N},n \to \infty}}{\leqslant} 
\mathcal{O} \left(\dfrac{\tilde{\mathfrak{c}}_{i_{1}i_{2}}^{\sharp,
\smallfrown,\mathrm{B}}(n,k,z_{o};j) \me^{-((n-1)K+k) \tilde{\lambda}_{
\tilde{\mathcal{R}},4}^{\sharp,\smallfrown}(j)}}{(n \! - \! 1)K \! + \! k} 
\right), \label{eqiy13}
\end{align}
where $\tilde{\lambda}_{\tilde{\mathcal{R}},4}^{\sharp,\smallfrown}(j)$ 
$(> \! 0)$ is defined by Equation~\eqref{eqqlamrupj}, and $\tilde{
\mathfrak{c}}_{i_{1}i_{2}}^{\sharp,\smallfrown,r}(n,k,z_{o};j) \! 
=_{\underset{z_{o}=1+o(1)}{\mathscr{N},n \to \infty}} \! \mathcal{O}(1)$, 
$r \! \in \! \lbrace \mathrm{A},\mathrm{B} \rbrace$. An elementary 
inequality argument shows that $\tilde{\mathbb{I}}_{4,j}^{\sharp,
\smallfrown,\mathrm{C}}(i_{1},i_{2})$ can be presented thus:
\begin{align} 
\tilde{\mathbb{I}}_{4,j}^{\sharp,\smallfrown,\mathrm{C}}(i_{1},i_{2}) 
\leqslant& \, \sup_{z \in \tilde{\Sigma}_{p,j}^{3} \cap [\theta_{0}^{
\tilde{a}}(j),\pi/2]} \left(\left\lvert \mathrm{P.V.} \int_{\tilde{\Sigma}_{
p,j}^{3} \cap [\theta_{0}^{\tilde{a}}(j),\pi/2]} \dfrac{(w_{+}^{\Sigma_{
\tilde{\mathcal{R}}}}(\xi))_{i_{1}i_{2}}}{\xi \! - \! z} \, \dfrac{\md \xi}{2 
\pi \mi} \right\rvert \! + \! \left\lvert \int_{\tilde{\Sigma}_{p,j}^{3} \cap 
[\pi/2,\pi -\theta_{0}^{\tilde{b}}(j)]} \dfrac{(w_{+}^{\Sigma_{\tilde{
\mathcal{R}}}}(\xi))_{i_{1}i_{2}}}{\xi \! - \! z} \, \dfrac{\md \xi}{2 \pi \mi} 
\right\rvert \right. \nonumber \\
+&\left. \, \dfrac{1}{2} \left\lvert (w_{+}^{\Sigma_{\tilde{\mathcal{
R}}}}(z))_{i_{1}i_{2}} \right\rvert \right) \! + \! \sup_{z \in \tilde{
\Sigma}_{p,j}^{3} \cap [\pi/2,\pi -\theta_{0}^{\tilde{b}}(j)]} \left(
\left\lvert \int_{\tilde{\Sigma}_{p,j}^{3} \cap [\theta_{0}^{\tilde{a}}(j),\pi/2]} 
\dfrac{(w_{+}^{\Sigma_{\tilde{\mathcal{R}}}}(\xi))_{i_{1}i_{2}}}{\xi \! - \! z} 
\, \dfrac{\md \xi}{2 \pi \mi} \right\rvert \right. \nonumber \\
+&\left. \, \mathrm{P.V.} \left\lvert \int_{\tilde{\Sigma}_{p,j}^{3} \cap 
[\pi/2,\pi -\theta_{0}^{\tilde{b}}(j)]} \dfrac{(w_{+}^{\Sigma_{\tilde{
\mathcal{R}}}}(\xi))_{i_{1}i_{2}}}{\xi \! - \! z} \, \dfrac{\md \xi}{2 \pi 
\mi} \right\rvert \! + \! \dfrac{1}{2} \left\lvert (w_{+}^{\Sigma_{
\tilde{\mathcal{R}}}}(z))_{i_{1}i_{2}} \right\rvert \right). \label{eqiy14}
\end{align}
Via the formula $w^{\Sigma_{\tilde{\mathcal{R}}}}_{+}(z) \! = \! \tilde{v}_{
\tilde{\mathcal{R}}}(z) \! - \! \mathrm{I}$, the asymptotics, in the 
double-scaling limit $\mathscr{N},n \! \to \! \infty$ such that $z_{o} 
\! = \! 1 \! + \! o(1)$, of $\tilde{v}_{\tilde{\mathcal{R}}}(z)$ given in 
Equation~\eqref{eqtlvee11}, the elementary trigonometric inequalities 
$\sin \theta \! \geqslant \! \tfrac{2 \theta}{\pi}$ and $\cos \theta \! 
\geqslant \! -\tfrac{2 \theta}{\pi} \! + \! 1$ for $0 \! \leqslant \! \theta \! 
\leqslant \! \tfrac{\pi}{2}$, and $\sin \theta \! \geqslant \! \tfrac{2}{\pi}
(\pi \! - \! \theta)$ and $\cos \theta \! \leqslant \! -\tfrac{2 \theta}{\pi} 
\! + \! 1$ for $\tfrac{\pi}{2} \! \leqslant \! \theta \! \leqslant \! \pi$, the 
inequality $\ln \lvert x \rvert \! \leqslant \! \lvert x \rvert \! - \! 1$, the 
parametrisation for the homotopic deformation of $\tilde{\Sigma}_{p,j}^{3}$ 
given in the proof of the corresponding item of Lemma~\ref{lem5.1}, that 
is, $\tilde{\Sigma}_{p,j}^{3} \! = \! \lbrace \mathstrut (x_{j}(\theta),y_{j}
(\theta)); \, x_{j}(\theta) \! = \! \tfrac{1}{2}(\tilde{a}_{j} \! + \! \tilde{b}_{j-1}) 
\! + \! \tfrac{1}{2}(\tilde{a}_{j} \! - \! \tilde{b}_{j-1}) \cos \theta, \, y_{j}
(\theta) \! = \! \tilde{\eta}_{j} \sin \theta, \, \theta_{0}^{\tilde{a}}(j) \! 
\leqslant \! \theta \! \leqslant \! \pi \! - \! \theta_{0}^{\tilde{b}}(j) \rbrace$, 
and writing $(\tilde{\Sigma}_{p,j}^{3} \! \ni)$ $z \! = \! x_{j}(\psi) \! + \! 
\mi y_{j}(\psi)$, $\theta_{0}^{\tilde{a}}(j) \! \leqslant \! \psi \! \leqslant 
\! \pi \! - \! \theta_{0}^{\tilde{b}}(j)$, one shows, via a tedious 
integration-by-parts argument, that, for $n \! \in \! \mathbb{N}$ and 
$k \! \in \! \lbrace 1,2,\dotsc,K \rbrace$ such that $\alpha_{p_{
\mathfrak{s}}} \! := \! \alpha_{k} \! \neq \! \infty$, $i_{1},i_{2} \! 
= \! 1,2$, and $j \! \in \! \lbrace 1,2,\dotsc,N \! + \! 1 \rbrace$,
\begin{gather*}
\sup_{z \in \tilde{\Sigma}_{p,j}^{3} \cap [\theta_{0}^{\tilde{a}}(j),\pi/2]} 
\dfrac{1}{2} \left\lvert (w_{+}^{\Sigma_{\tilde{\mathcal{R}}}}(z))_{i_{1}i_{2}} 
\right\rvert \! + \! \sup_{z \in \tilde{\Sigma}_{p,j}^{3} \cap [\pi/2,\pi 
-\theta_{0}^{\tilde{b}}(j)]} \dfrac{1}{2} \left\lvert (w_{+}^{\Sigma_{\tilde{
\mathcal{R}}}}(z))_{i_{1}i_{2}} \right\rvert \! = \! \sup_{\theta_{0}^{\tilde{a}}
(j) \leqslant \psi \leqslant \pi/2} \dfrac{1}{2} \left\lvert (w_{+}^{\Sigma_{
\tilde{\mathcal{R}}}}(x_{j}(\psi) \! + \! \mi y_{j}(\psi)))_{i_{1}i_{2}} \right\rvert \\
+ \sup_{\pi/2 \leqslant \psi \leqslant \pi -\theta_{0}^{\tilde{b}}(j)} 
\dfrac{1}{2} \left\lvert (w_{+}^{\Sigma_{\tilde{\mathcal{R}}}}(x_{j}(\psi) 
\! + \! \mi y_{j}(\psi)))_{i_{1}i_{2}} \right\rvert \underset{\underset{z_{o}
=1+o(1)}{\mathscr{N},n \to \infty}}{\leqslant} \mathcal{O} \left(\tilde{
\mathfrak{c}}_{i_{1}i_{2}}^{\sharp,\smallfrown,\mathrm{C}_{1}}(n,k,
z_{o};j) \me^{-((n-1)K+k) \tilde{\lambda}_{\tilde{\mathcal{R}},4}^{
\sharp,\smallfrown}(j)} \right),
\end{gather*}
\begin{align*}
\sup_{z \in \tilde{\Sigma}_{p,j}^{3} \cap [\theta_{0}^{\tilde{a}}(j),\pi/2]} 
&\left(\left\lvert \int_{\tilde{\Sigma}_{p,j}^{3} \cap [\pi/2,\pi -\theta_{
0}^{\tilde{b}}(j)]} \dfrac{(w_{+}^{\Sigma_{\tilde{\mathcal{R}}}}(\xi))_{i_{1}
i_{2}}}{\xi \! - \! z} \, \dfrac{\md \xi}{2 \pi \mi} \right\rvert \right) \! 
+ \! \sup_{z \in \tilde{\Sigma}_{p,j}^{3} \cap [\pi/2,\pi -\theta_{0}^{
\tilde{b}}(j)]} \left(\left\lvert \int_{\tilde{\Sigma}_{p,j}^{3} \cap [\theta_{
0}^{\tilde{a}}(j),\pi/2]} \dfrac{(w_{+}^{\Sigma_{\tilde{\mathcal{R}}}}(\xi))_{
i_{1}i_{2}}}{\xi \! - \! z} \, \dfrac{\md \xi}{2 \pi \mi} \right\rvert \right) \\
&\leqslant \sup_{\theta_{0}^{\tilde{a}}(j) \leqslant \psi < \pi/2} \left(
\int_{\pi/2}^{\pi -\theta_{0}^{\tilde{b}}(j)} \left\lvert \dfrac{(w_{+}^{
\Sigma_{\tilde{\mathcal{R}}}}(x_{j}(\theta) \! + \! \mi y_{j}(\theta)))_{i_{1}
i_{2}}}{x_{j}(\theta) \! - \! x_{j}(\psi) \! + \! \mi (y_{j}(\theta) \! - \! y_{j}
(\psi))} \right\rvert ((x_{j}^{\prime}(\theta))^{2} \! + \! (y_{j}^{\prime}
(\theta))^{2})^{1/2} \, \dfrac{\md \theta}{2 \pi} \right) \\
&+ \sup_{\pi/2 < \psi \leqslant \pi - \theta_{0}^{\tilde{b}}(j)} \left(
\int_{\theta_{0}^{\tilde{a}}(j)}^{\pi/2} \left\lvert \dfrac{(w_{+}^{
\Sigma_{\tilde{\mathcal{R}}}}(x_{j}(\theta) \! + \! \mi y_{j}(\theta)))_{i_{1}
i_{2}}}{x_{j}(\theta) \! - \! x_{j}(\psi) \! + \! \mi (y_{j}(\theta) \! - \! y_{j}
(\psi))} \right\rvert ((x_{j}^{\prime}(\theta))^{2} \! + \! (y_{j}^{\prime}
(\theta))^{2})^{1/2} \, \dfrac{\md \theta}{2 \pi} \right) \\
&\underset{\underset{z_{o}=1+o(1)}{\mathscr{N},n \to \infty}}{\leqslant} 
\mathcal{O} \left(\dfrac{\tilde{\mathfrak{c}}_{i_{1}i_{2}}^{\sharp,
\smallfrown,\mathrm{C}_{2}}(n,k,z_{o};j) \me^{-((n-1)K+k) \tilde{
\lambda}_{\tilde{\mathcal{R}},4}^{\sharp,\smallfrown}(j)}}{(n \! - \! 1)K 
\! + \! k} \right),
\end{align*}
\begin{align*}
&\sup_{z \in \tilde{\Sigma}_{p,j}^{3} \cap [\theta_{0}^{\tilde{a}}(j),\pi/2]} 
\left(\left\lvert \mathrm{P.V.} \int_{\tilde{\Sigma}_{p,j}^{3} \cap 
[\theta_{0}^{\tilde{a}}(j),\pi/2]} \dfrac{(w_{+}^{\Sigma_{\tilde{\mathcal{R}}}}
(\xi))_{i_{1}i_{2}}}{\xi \! - \! z} \, \dfrac{\md \xi}{2 \pi \mi} \right\rvert 
\right) \! + \! \sup_{z \in \tilde{\Sigma}_{p,j}^{3} \cap [\pi/2,\pi 
-\theta_{0}^{\tilde{b}}(j)]} \left(\left\lvert \mathrm{P.V.} \int_{\tilde{
\Sigma}_{p,j}^{3} \cap [\pi/2,\pi -\theta_{0}^{\tilde{b}}(j)]} \dfrac{(w_{+}^{
\Sigma_{\tilde{\mathcal{R}}}}(\xi))_{i_{1}i_{2}}}{\xi \! - \! z} \, \dfrac{\md 
\xi}{2 \pi \mi} \right\rvert \right) \\
&\underset{\underset{z_{o}=1+o(1)}{\mathscr{N},n \to \infty}}{\leqslant} 
\sup_{\theta_{0}^{\tilde{a}}(j) \leqslant \psi \leqslant \pi/2} \left(\left\lvert 
\left(\int_{\theta_{0}^{\tilde{a}}(j)}^{\psi - n^{-1}} \! + \! \int_{\psi + 
n^{-1}}^{\pi/2} \right) \dfrac{(w_{+}^{\Sigma_{\tilde{\mathcal{R}}}}(x_{j}
(\theta) \! + \! \mi y_{j}(\theta)))_{i_{1}i_{2}}}{x_{j}(\theta) \! - \! x_{j}
(\psi) \! + \! \mi (y_{j}(\theta) \! - \! y_{j}(\psi))} \, \dfrac{\md (x_{j}
(\theta) \! + \! \mi y_{j}(\theta))}{2 \pi \mi} \right\rvert \right) \\
&+ \sup_{\pi/2 \leqslant \psi \leqslant \pi -\theta_{0}^{\tilde{b}}(j)} 
\left(\left\lvert \left(\int_{\pi/2}^{\psi - n^{-1}} \! + \! \int_{\psi + 
n^{-1}}^{\pi - \theta_{0}^{\tilde{b}}(j)} \right) \dfrac{(w_{+}^{\Sigma_{
\tilde{\mathcal{R}}}}(x_{j}(\theta) \! + \! \mi y_{j}(\theta)))_{i_{1}i_{2}}}{
x_{j}(\theta) \! - \! x_{j}(\psi) \! + \! \mi (y_{j}(\theta) \! - \! y_{j}(\psi))} 
\, \dfrac{\md (x_{j}(\theta) \! + \! \mi y_{j}(\theta))}{2 \pi \mi} \right\rvert 
\right) \\
&\underset{\underset{z_{o}=1+o(1)}{\mathscr{N},n \to \infty}}{\leqslant} 
\sup_{\theta_{0}^{\tilde{a}}(j) \leqslant \psi \leqslant \pi/2} \left(\left(
\int_{\theta_{0}^{\tilde{a}}(j)}^{\psi - n^{-1}} \! + \! \int_{\psi + n^{-1}}^{\pi/2} 
\right) \dfrac{\mathcal{O} \left(\lvert (\tilde{\mathfrak{c}}_{\tilde{\mathcal{R}},3}
(j))_{i_{1}i_{2}} \rvert \me^{-((n-1)K+k) \tilde{\lambda}_{\tilde{\mathcal{R}},3}(j) 
\tilde{\mathfrak{x}}_{1}^{\smallfrown}(j) \left(1-\frac{2 \tilde{\delta}_{\tilde{a}_{j}}
((\tilde{b}_{j-1}-\tilde{a}_{j}) \cos \sigma_{\tilde{a}_{j}}^{+}+2 \tilde{\eta}_{j} \sin 
\sigma_{\tilde{a}_{j}}^{+})}{\pi (\tilde{\mathfrak{x}}_{1}^{\smallfrown}(j))^{2}} 
\theta \right)^{1/2}} \right)}{\left\lvert \sin \left(\frac{\theta - \psi}{2} \right) 
\right\rvert \left((\tilde{a}_{j} \! - \! \tilde{b}_{j-1})^{2} \sin^{2} \left(\frac{
\theta_{0}^{\tilde{a}}(j)+ \psi}{2} \right) \! + \! 4 \tilde{\eta}_{j}^{2} \cos^{2} 
\left(\frac{\frac{\pi}{2} + \psi}{2} \right) \right)^{1/2}} \right. \\
&\left. \, \qquad \times  \sqrt{\dfrac{1}{4}(\tilde{a}_{j} \! - \! \tilde{b}_{j-1})^{2} 
\! + \! \tilde{\eta}_{j}^{2}(\cos (\min \lbrace \theta_{0}^{\tilde{a}}(j),
\theta_{0}^{\tilde{b}}(j) \rbrace))^{2}} \, \dfrac{\md \theta}{2 \pi} 
\vphantom{M^{M^{M^{M^{M^{M^{M^{M^{M^{M^{M^{M}}}}}}}}}}}} \right) \\
&+ \sup_{\pi/2 \leqslant \psi \leqslant \pi - \theta_{0}^{\tilde{b}}(j)} 
\left(\left(\int_{\pi/2}^{\psi - n^{-1}} \! + \! \int_{\psi + n^{-1}}^{\pi 
- \theta_{0}^{\tilde{b}}(j)} \right) \dfrac{\mathcal{O} \left(\lvert (\tilde{
\mathfrak{c}}_{\tilde{\mathcal{R}},3}(j))_{i_{1}i_{2}} \rvert \me^{-((n-1)K+k) 
\tilde{\lambda}_{\tilde{\mathcal{R}},3}(j) \tilde{\mathfrak{x}}_{2}^{\smallfrown}
(j) \left(1-\frac{2((\tilde{b}_{j-1}-\tilde{a}_{j}) \tilde{\delta}_{\tilde{a}_{j}} \cos 
\sigma_{\tilde{a}_{j}}^{+}+2 \tilde{\eta}_{j}(2 \tilde{\eta}_{j} - \tilde{\delta}_{
\tilde{a}_{j}} \sin \sigma_{\tilde{a}_{j}}^{+}))}{\pi (\tilde{\mathfrak{x}}_{2}^{
\smallfrown}(j))^{2}} \theta \right)^{1/2}} \right)}{\left\lvert \sin \left(
\frac{\theta - \psi}{2} \right) \right\rvert \left((\tilde{a}_{j} \! - \! 
\tilde{b}_{j-1})^{2} \sin^{2} \left(\frac{\pi -\theta_{0}^{\tilde{b}}(j)+ \psi}{2} 
\right) \! + \! 4 \tilde{\eta}_{j}^{2} \cos^{2} \left(\frac{\frac{\pi}{2} + 
\psi}{2} \right) \right)^{1/2}} \right. \\
&\left. \, \qquad \times \sqrt{\dfrac{1}{4}(\tilde{a}_{j} \! - \! \tilde{b}_{j-1})^{2} 
\! + \! \tilde{\eta}_{j}^{2}(\cos (\min \lbrace \theta_{0}^{\tilde{a}}(j),
\theta_{0}^{\tilde{b}}(j) \rbrace))^{2}} \, \dfrac{\md \theta}{2 \pi} 
\vphantom{M^{M^{M^{M^{M^{M^{M^{M^{M^{M^{M^{M}}}}}}}}}}}} \right) \\
&\underset{\underset{z_{o}=1+o(1)}{\mathscr{N},n \to \infty}}{\leqslant} 
\mathcal{O} \left(\tilde{\mathfrak{c}}_{i_{1}i_{2}}^{\sharp,\smallfrown,
\mathrm{C}_{3}}(n,k,z_{o};j) \me^{-((n-1)K+k) \tilde{\lambda}_{\tilde{
\mathcal{R}},4}^{\sharp,\smallfrown}(j)} \right),
\end{align*}
where $\tilde{\mathfrak{c}}_{i_{1}i_{2}}^{\sharp,\smallfrown,r}(n,k,z_{o};j) 
\! =_{\underset{z_{o}=1+o(1)}{\mathscr{N},n \to \infty}} \! \mathcal{O}(1)$, 
$r \! \in \! \lbrace \mathrm{C}_{1},\mathrm{C}_{2},\mathrm{C}_{3} \rbrace$, 
whence, for $i_{1},i_{2} \! = \! 1,2$, via an elementary inequality argument, 
one shows that
\begin{equation} \label{eqiy15} 
\tilde{\mathbb{I}}_{4,j}^{\sharp,\smallfrown,\mathrm{C}}(i_{1},i_{2}) 
\underset{\underset{z_{o}=1+o(1)}{\mathscr{N},n \to \infty}}{\leqslant} 
\mathcal{O} \left(\tilde{\mathfrak{c}}_{i_{1}i_{2}}^{\sharp,\smallfrown,
\mathrm{C}}(n,k,z_{o};j) \me^{-((n-1)K+k) \tilde{\lambda}_{\tilde{
\mathcal{R}},4}^{\sharp,\smallfrown}(j)} \right),
\end{equation}
where $\tilde{\mathfrak{c}}_{i_{1}i_{2}}^{\sharp,\smallfrown,\mathrm{C}}
(n,k,z_{o};j) \! =_{\underset{z_{o}=1+o(1)}{\mathscr{N},n \to \infty}} \! 
\mathcal{O}(1)$; thus, via the---asymptotic---Estimates~\eqref{eqiy12}, 
\eqref{eqiy13}, and~\eqref{eqiy15}, it follows {}from Equation~\eqref{eqiy11} 
and an elementary inequality argument, that, for $n \! \in \! \mathbb{N}$ 
and $k \! \in \! \lbrace 1,2,\dotsc,K \rbrace$ such that $\alpha_{p_{
\mathfrak{s}}} \! := \! \alpha_{k} \! \neq \! \infty$, $i_{1},i_{2} \! = \! 1,2$, 
and $j \! \in \! \lbrace 1,2,\dotsc,N \! + \! 1 \rbrace$,
\begin{equation} \label{eqiy16} 
\tilde{\mathbb{I}}_{4,j}^{\sharp,\smallfrown}(i_{1},i_{2}) \underset{
\underset{z_{o}=1+o(1)}{\mathscr{N},n \to \infty}}{\leqslant} \mathcal{O} 
\left(\tilde{\mathfrak{c}}_{i_{1}i_{2}}^{\sharp,\smallfrown,\blacklozenge}
(n,k,z_{o};j) \me^{-((n-1)K+k) \tilde{\lambda}_{\tilde{\mathcal{R}},4}^{
\sharp,\smallfrown}(j)} \right),
\end{equation}
where $\tilde{\mathfrak{c}}_{i_{1}i_{2}}^{\sharp,\smallfrown,\blacklozenge}
(n,k,z_{o};j) \! =_{\underset{z_{o}=1+o(1)}{\mathscr{N},n \to \infty}} \! 
\mathcal{O}(1)$. For $n \! \in \! \mathbb{N}$ and $k \! \in \! \lbrace 1,2,
\dotsc,K \rbrace$ such that $\alpha_{p_{\mathfrak{s}}} \! := \! \alpha_{k} 
\! \neq \! \infty$, the integral
\begin{equation*}
\tilde{\mathbb{I}}_{4,j}^{\sharp,\smallsmile}(i_{1},i_{2}) \! := \! \sup_{z 
\in \tilde{\Sigma}_{\tilde{\mathcal{R}}}^{\sharp}} \left\lvert \lim_{-\tilde{
\Sigma}_{\tilde{\mathcal{R}}}^{\sharp} \ni z^{\prime} \to z} \int_{\tilde{
\Sigma}_{p,j}^{4}} \dfrac{(z^{\prime} \! - \! \alpha_{k})(w_{+}^{\Sigma_{
\tilde{\mathcal{R}}}}(\xi))_{i_{1}i_{2}}}{(\xi \! - \! \alpha_{k})(\xi \! - \! 
z^{\prime})} \, \dfrac{\md \xi}{2 \pi \mi} \right\rvert, \quad i_{1},i_{2} 
\! = \! 1,2, \quad j \! \in \! \lbrace 1,2,\dotsc,N \! + \! 1 \rbrace,
\end{equation*}
is analysed analogously (as above for $\tilde{\mathbb{I}}_{4,j}^{\sharp,
\smallsmile}(i_{1},i_{2})$), and one arrives at, for $n \! \in \! \mathbb{N}$ 
and $k \! \in \! \lbrace 1,2,\dotsc,K \rbrace$ such that $\alpha_{p_{
\mathfrak{s}}} \! := \! \alpha_{k} \! \neq \! \infty$, $i_{1},i_{2} \! = \! 
1,2$, and $j \! \in \! \lbrace 1,2,\dotsc,N \! + \! 1 \rbrace$,
\begin{equation} \label{eqiy17} 
\tilde{\mathbb{I}}_{4,j}^{\sharp,\smallsmile}(i_{1},i_{2}) \underset{
\underset{z_{o}=1+o(1)}{\mathscr{N},n \to \infty}}{\leqslant} \mathcal{O} 
\left(\tilde{\mathfrak{c}}_{i_{1}i_{2}}^{\sharp,\smallsmile,\lozenge}(n,k,
z_{o};j) \me^{-((n-1)K+k) \tilde{\lambda}_{\tilde{\mathcal{R}},4}^{\sharp,
\smallsmile}(j)} \right),
\end{equation}
where $\tilde{\lambda}_{\tilde{\mathcal{R}},4}^{\sharp,\smallsmile}(j)$ 
$(> \! 0)$ is defined by Equation~\eqref{eqqlamrdwj}, and $\tilde{
\mathfrak{c}}_{i_{1}i_{2}}^{\sharp,\smallsmile,\lozenge}(n,k,z_{o};j) \! 
=_{\underset{z_{o}=1+o(1)}{\mathscr{N},n \to \infty}} \! \mathcal{O}
(1)$; thus, via the---asymptotic---Estimates~\eqref{eqiy16} 
and~\eqref{eqiy17}, and an elementary inequality argument, one finally 
arrives at, for $n \! \in \! \mathbb{N}$ and $k \! \in \! \lbrace 1,2,\dotsc,
K \rbrace$ such that $\alpha_{p_{\mathfrak{s}}} \! := \! \alpha_{k} \! \neq 
\! \infty$,
\begin{equation} \label{eqiy18} 
\tilde{\mathbb{I}}_{4}^{\sharp,\smallsmile}(i_{1},i_{2}) \underset{
\underset{z_{o}=1+o(1)}{\mathscr{N},n \to \infty}}{=} \mathcal{O} \left(
\tilde{\mathfrak{c}}_{i_{1}i_{2}}^{\sharp,4}(n,k,z_{o}) \me^{-((n-1)K+k) 
\min \left\lbrace \min\limits_{j=1,2,\dotsc,N+1} \lbrace \tilde{\lambda}_{
\tilde{\mathcal{R}},4}^{\sharp,\smallfrown}(j) \rbrace,\min\limits_{j=1,2,
\dotsc,N+1} \lbrace \tilde{\lambda}_{\tilde{\mathcal{R}},4}^{\sharp,
\smallsmile}(j) \rbrace \right\rbrace} \right),
\end{equation}
where $\tilde{\mathfrak{c}}_{i_{1}i_{2}}^{\sharp,4}(n,k,z_{o}) \! 
=_{\underset{z_{o}=1+o(1)}{\mathscr{N},n \to \infty}} \! \mathcal{O}(1)$.

For $n \! \in \! \mathbb{N}$ and $k \! \in \! \lbrace 1,2,\dotsc,K \rbrace$ 
such that $\alpha_{p_{\mathfrak{s}}} \! := \! \alpha_{k} \! \neq \! \infty$, 
in order to estimate, asymptotically, in the double-scaling limit 
$\mathscr{N},n \! \to \! \infty$ such that $z_{o} \! = \! 1 \! + \! o(1)$, 
$\tilde{\mathbb{I}}_{5}^{\sharp}(i_{1},i_{2})$, $i_{1},i_{2} \! = \! 1,2$, 
consider, say, the quantity
\begin{equation*}
\tilde{\mathbb{I}}_{5,\circlearrowright}^{\sharp}(i_{1},i_{2}) \! := \! \sup_{z 
\in \tilde{\Sigma}_{\tilde{\mathcal{R}}}^{\sharp}} \left\lvert \lim_{-\tilde{
\Sigma}_{\tilde{\mathcal{R}}}^{\sharp} \ni z^{\prime} \to z} \int_{\partial 
\tilde{\mathbb{U}}_{\tilde{\delta}_{\tilde{a}_{N+1}}}} \dfrac{(z^{\prime} \! 
- \! \alpha_{k})(w_{+}^{\Sigma_{\tilde{\mathcal{R}}}}(\xi))_{i_{1}i_{2}}}{(\xi 
\! - \! \alpha_{k})(\xi \! - \! z^{\prime})} \, \dfrac{\md \xi}{2 \pi \mi} 
\right\rvert:
\end{equation*}
a straightforward inequality argument shows that
\begin{equation}
\tilde{\mathbb{I}}_{5,\circlearrowright}^{\sharp}(i_{1},i_{2}) \! \leqslant \! 
\tilde{\mathbb{I}}_{5,\circlearrowright}^{\sharp,\mathrm{A}}(i_{1},i_{2}) \! 
+ \! \tilde{\mathbb{I}}_{5,\circlearrowright}^{\sharp,\mathrm{B}}(i_{1},i_{2}) 
\! + \! \tilde{\mathbb{I}}_{5,\circlearrowright}^{\sharp,\mathrm{C}}(i_{1},
i_{2}), \quad i_{1},i_{2} \! = \! 1,2, \label{eqiy19}
\end{equation} 
where
\begin{gather*}
\tilde{\mathbb{I}}_{5,\circlearrowright}^{\sharp,\mathrm{A}}(i_{1},i_{2}) 
\! := \! \sup \left\lvert \int_{\partial \tilde{\mathbb{U}}_{\tilde{\delta}_{
\tilde{a}_{N+1}}}} \dfrac{(w_{+}^{\Sigma_{\tilde{\mathcal{R}}}}(\xi))_{i_{1}
i_{2}}}{\xi \! - \! \alpha_{k}} \, \dfrac{\md \xi}{2 \pi \mi} \right\rvert, \quad 
\quad \tilde{\mathbb{I}}_{5,\circlearrowright}^{\sharp,\mathrm{B}}(i_{1},
i_{2}) \! := \! \sup_{z \in \partial \tilde{\mathbb{U}}_{\tilde{\delta}_{
\tilde{a}_{N+1}}}} \left\lvert \lim_{-\tilde{\Sigma}_{\tilde{\mathcal{R}}}^{
\sharp} \ni z^{\prime} \to z} \int_{\partial \tilde{\mathbb{U}}_{\tilde{\delta}_{
\tilde{a}_{N+1}}}} \dfrac{(w_{+}^{\Sigma_{\tilde{\mathcal{R}}}}(\xi))_{i_{1}
i_{2}}}{\xi \! - \! z^{\prime}} \, \dfrac{\md \xi}{2 \pi \mi} \right\rvert, \\
\tilde{\mathbb{I}}_{5,\circlearrowright}^{\sharp,\mathrm{C}}(i_{1},i_{2}) \! 
:= \! \sup_{z \in \tilde{\Sigma}_{\tilde{\mathcal{R}}}^{\sharp} \setminus 
\partial \tilde{\mathbb{U}}_{\tilde{\delta}_{\tilde{a}_{N+1}}}} \left\lvert 
\lim_{-\tilde{\Sigma}_{\tilde{\mathcal{R}}}^{\sharp} \ni z^{\prime} \to z} 
\int_{\partial \tilde{\mathbb{U}}_{\tilde{\delta}_{\tilde{a}_{N+1}}}} 
\dfrac{(w_{+}^{\Sigma_{\tilde{\mathcal{R}}}}(\xi))_{i_{1}i_{2}}}{\xi \! - \! 
z^{\prime}} \, \dfrac{\md \xi}{2 \pi \mi} \right\rvert.
\end{gather*}
Via the asymptotic expansion, in the double-scaling limit $\mathscr{N},
n \! \to \! \infty$ such that $z_{o} \! = \! 1 \! + \! o(1)$, of $w_{+}^{
\Sigma_{\tilde{\mathcal{R}}}}(z)$ for $z \! \in \! \partial \tilde{\mathbb{U}}_{
\tilde{\delta}_{\tilde{a}_{N+1}}}$ given in Equation~\eqref{eqproptila}, 
H\"{o}lder's Inequality for Integrals, and the indefinite integral \pmb{2.553} 3. 
on p.~180 of \cite{gradryzh}, one shows that, with $p_{0} \! := \! 
(\tilde{a}_{N+1} \! - \! \alpha_{k})^{2} \! + \! \tilde{\delta}_{\tilde{a}_{N+
1}}^{2}$ and $q_{0} \! := \! 2 \tilde{\delta}_{\tilde{a}_{N+1}}(\tilde{a}_{N+1} 
\! - \! \alpha_{k})$ $(\Rightarrow p_{0}^{2} \! - \! q_{0}^{2} \! > \! 0)$, 
for $n \! \in \! \mathbb{N}$ and $k \! \in \! \lbrace 1,2,\dotsc,K 
\rbrace$ such that $\alpha_{p_{\mathfrak{s}}} \! := \! \alpha_{k} 
\! \neq \! \infty$,
\begin{align*}
\tilde{\mathbb{I}}_{5,\circlearrowright}^{\sharp,\mathrm{A}}(i_{1},i_{2}) 
\underset{\underset{z_{o}=1+o(1)}{\mathscr{N},n \to \infty}}{\leqslant}& 
\, \dfrac{\tilde{\delta}_{\tilde{a}_{N+1}}}{\sqrt{2 \pi}} \sup_{\theta_{0} 
\in [0,2 \pi]} \left\lvert (w_{+}^{\Sigma_{\tilde{\mathcal{R}}}}(\tilde{a}_{N+1} 
\! + \! \tilde{\delta}_{\tilde{a}_{N+1}} \me^{\mi \theta_{0}}))_{i_{1}i_{2}} 
\right\rvert \sqrt{\int_{0}^{2 \pi} \dfrac{1}{p_{0} \! + \! q_{0} \cos \theta} 
\, \md \theta} \nonumber \\
\underset{\underset{z_{o}=1+o(1)}{\mathscr{N},n \to \infty}}{\leqslant}& \, 
\dfrac{\tilde{\delta}_{\tilde{a}_{N+1}}}{\sqrt{2 \pi}} \sup_{\theta_{0} 
\in [0,2 \pi]} \left\lvert (w_{+}^{\Sigma_{\tilde{\mathcal{R}}}}(\tilde{a}_{N+1} 
\! + \! \tilde{\delta}_{\tilde{a}_{N+1}} \me^{\mi \theta_{0}}))_{i_{1}i_{2}} 
\right\rvert \sqrt{\dfrac{2}{\sqrt{p_{0}^{2} \! - \! q_{0}^{2}}}} \underbrace{
\sqrt{\left. \tan^{-1} \left(\sqrt{\dfrac{p_{0} \! - \! q_{0}}{p_{0} \! + \! 
q_{0}}} \tan \left(\dfrac{\theta}{2} \right) \right) \right\vert_{0}^{2 
\pi}}}_{= \, 0} \quad \Rightarrow \nonumber
\end{align*}
\begin{equation} \label{eqiy20} 
\tilde{\mathbb{I}}_{5,\circlearrowright}^{\sharp,\mathrm{A}}(i_{1},i_{2}) 
\underset{\underset{z_{o}=1+o(1)}{\mathscr{N},n \to \infty}}{=} 0.
\end{equation}
Via the asymptotic expansion, in the double-scaling limit $\mathscr{N},
n \! \to \! \infty$ such that $z_{o} \! = \! 1 \! + \! o(1)$, of $w_{+}^{
\Sigma_{\tilde{\mathcal{R}}}}(z)$ for $z \! \in \! \partial \tilde{\mathbb{
U}}_{\tilde{\delta}_{\tilde{a}_{N+1}}}$ given in Equation~\eqref{eqproptila}, 
H\"{o}lder's Inequality for Integrals, and the indefinite integral \pmb{2.562} 
1. on p.~185 of \cite{gradryzh}, one shows that, with $A_{\epsilon}(\theta,
\psi) \! := \! 2 \tilde{\delta}_{\tilde{a}_{N+1}} \sin (\tfrac{\theta - \psi}{2}) 
\sin (\tfrac{\theta + \psi}{2}) \! + \! \epsilon \cos (\theta \! - \! \psi)$ and 
$B_{\epsilon}(\theta,\psi) \! := \! -2 \tilde{\delta}_{\tilde{a}_{N+1}} \sin 
(\tfrac{\theta - \psi}{2}) \cos (\tfrac{\theta + \psi}{2}) \! + \! \epsilon \sin 
(\theta \! - \! \psi)$, for $n \! \in \! \mathbb{N}$ and $k \! \in \! \lbrace 
1,2,\dotsc,K \rbrace$ such that $\alpha_{p_{\mathfrak{s}}} \! := \! 
\alpha_{k} \! \neq \! \infty$,
\begin{align*}
\tilde{\mathbb{I}}_{5,\circlearrowright}^{\sharp,\mathrm{B}}(i_{1},i_{2}) 
\underset{\underset{z_{o}=1+o(1)}{\mathscr{N},n \to \infty}}{\leqslant}& 
\, \dfrac{\tilde{\delta}_{\tilde{a}_{N+1}}}{2 \pi} \sup_{\theta_{0} \in [0,
2 \pi]} \left\lvert (w_{+}^{\Sigma_{\tilde{\mathcal{R}}}}(\tilde{a}_{N+1} \! + \! 
\tilde{\delta}_{\tilde{a}_{N+1}} \me^{\mi \theta_{0}}))_{i_{1}i_{2}} \right\rvert 
\sup_{\psi \in [0,2 \pi]} \left\lvert \lim_{\epsilon \downarrow 0} \int_{0}^{2 
\pi} \dfrac{(A_{\epsilon}(\theta,\psi) \! + \! \mi B_{\epsilon}(\theta,\psi))}{
\lvert 2 \mi \tilde{\delta}_{\tilde{a}_{N+1}} \me^{\frac{\mi}{2}(\theta + 
\psi)} \sin (\frac{\theta - \psi}{2}) \! + \! \epsilon \me^{\mi \psi} \rvert^{2}} 
\, \md \theta \right\rvert \nonumber \\
\underset{\underset{z_{o}=1+o(1)}{\mathscr{N},n \to \infty}}{\leqslant}& \, 
\dfrac{\tilde{\delta}_{\tilde{a}_{N+1}}}{\sqrt{2} \pi} \sup_{\theta_{0} \in [0,
2 \pi]} \left\lvert (w_{+}^{\Sigma_{\tilde{\mathcal{R}}}}(\tilde{a}_{N+1} \! + \! 
\tilde{\delta}_{\tilde{a}_{N+1}} \me^{\mi \theta_{0}}))_{i_{1}i_{2}} \right\rvert 
\sup_{\psi \in [0,2 \pi]} \left(\lim_{\epsilon \downarrow 0} \left(
\lim_{\delta \downarrow 0} \left(\int_{0}^{\psi - \delta} \! + \! \int_{\psi 
+ \delta}^{2 \pi} \right) \dfrac{(2 \tilde{\delta}_{\tilde{a}_{N+1}} \! + \! 
\epsilon)}{\epsilon^{2} \! + \! 4 \tilde{\delta}_{\tilde{a}_{N+1}}(\tilde{
\delta}_{\tilde{a}_{N+1}} \! - \! \epsilon) \sin^{2}(\frac{\theta - \psi}{2})} 
\, \md \theta \right) \right) \nonumber \\
\underset{\underset{z_{o}=1+o(1)}{\mathscr{N},n \to \infty}}{\leqslant}& \, 
\dfrac{\sqrt{2} \, \tilde{\delta}_{\tilde{a}_{N+1}}}{\pi} \sup_{\theta_{0} \in 
[0,2 \pi]} \left\lvert (w_{+}^{\Sigma_{\tilde{\mathcal{R}}}}(\tilde{a}_{N+1} 
\! + \! \tilde{\delta}_{\tilde{a}_{N+1}} \me^{\mi \theta_{0}}))_{i_{1}i_{2}} 
\right\rvert \sup_{\psi \in [0,2 \pi]} \left(\lim_{\epsilon \downarrow 0} \left(
\lim_{\delta \downarrow 0} \left(\dfrac{2 \tilde{\delta}_{\tilde{a}_{N+1}} \! 
+ \! \epsilon}{\epsilon (2 \tilde{\delta}_{\tilde{a}_{N+1}} \! - \! \epsilon)} 
\right) \right. \right. \nonumber \\
\times&\left. \left. \, \underbrace{\left(\left. \tan^{-1} \left(\left(
\dfrac{2 \tilde{\delta}_{\tilde{a}_{N+1}} \! - \! \epsilon}{\epsilon} \right) 
\tan u_{1} \right) \right\vert_{u_{1}=-\psi/2}^{u_{1}=-\delta/2} \! + \! 
\left. \tan^{-1} \left(\left(\dfrac{2 \tilde{\delta}_{\tilde{a}_{N+1}} \! - \! 
\epsilon}{\epsilon} \right) \tan u_{2} \right) \right\vert_{u_{2}=\delta/2}^{
u_{2}=\pi -\psi/2} \right)}_{= \, 0} \right) \right) \quad \Rightarrow
\end{align*}
\begin{equation} \label{eqiy21} 
\tilde{\mathbb{I}}_{5,\circlearrowright}^{\sharp,\mathrm{B}}(i_{1},i_{2}) 
\underset{\underset{z_{o}=1+o(1)}{\mathscr{N},n \to \infty}}{=} 0.
\end{equation}
For $n \! \in \! \mathbb{N}$ and $k \! \in \! \lbrace 1,2,\dotsc,K \rbrace$ 
such that $\alpha_{p_{\mathfrak{s}}} \! := \! \alpha_{k} \! \neq \! \infty$, 
one shows that
\begin{equation} \label{eqiy22} 
\tilde{\mathbb{I}}_{5,\circlearrowright}^{\sharp,\mathrm{C}}(i_{1},i_{2}) 
\! \leqslant \! \sum_{m=1}^{5} \tilde{\mathbb{I}}_{5,\circlearrowright}^{
\sharp,\mathrm{C},\flat_{m}}(i_{1},i_{2}), \quad i_{1},i_{2} \! = \! 1,2,
\end{equation}
where
\begin{align*}
\tilde{\mathbb{I}}_{5,\circlearrowright}^{\sharp,\mathrm{C},\flat_{1}}(i_{1},
i_{2}) :=& \, \sup_{z \in \tilde{\Sigma}_{p}^{1}} \left\lvert \lim_{-\tilde{
\Sigma}_{p}^{1} \ni z^{\prime} \to z} \int_{\partial \tilde{\mathbb{U}}_{
\tilde{\delta}_{\tilde{a}_{N+1}}}} \dfrac{(w_{+}^{\Sigma_{\tilde{\mathcal{R}}}}
(\xi))_{i_{1}i_{2}}}{\xi \! - \! z^{\prime}} \, \dfrac{\md \xi}{2 \pi \mi} 
\right\rvert, \\
\tilde{\mathbb{I}}_{5,\circlearrowright}^{\sharp,\mathrm{C},\flat_{2}}(i_{1},
i_{2}) :=& \, \sup_{z \in \cup_{j=1}^{N} \tilde{\Sigma}_{p,j}^{2}} \left\lvert 
\lim_{-\cup_{j=1}^{N} \tilde{\Sigma}_{p,j}^{2} \ni z^{\prime} \to z} \int_{
\partial \tilde{\mathbb{U}}_{\tilde{\delta}_{\tilde{a}_{N+1}}}} \dfrac{(w_{+}^{
\Sigma_{\tilde{\mathcal{R}}}}(\xi))_{i_{1}i_{2}}}{\xi \! - \! z^{\prime}} \, 
\dfrac{\md \xi}{2 \pi \mi} \right\rvert, \\
\tilde{\mathbb{I}}_{5,\circlearrowright}^{\sharp,\mathrm{C},\flat_{3}}(i_{1},
i_{2}) :=& \, \sup_{z \in \cup_{j=1}^{N+1} \tilde{\Sigma}_{p,j}^{3}} \left\lvert 
\lim_{-\cup_{j=1}^{N+1} \tilde{\Sigma}_{p,j}^{3} \ni z^{\prime} \to z} \int_{
\partial \tilde{\mathbb{U}}_{\tilde{\delta}_{\tilde{a}_{N+1}}}} \dfrac{(w_{+}^{
\Sigma_{\tilde{\mathcal{R}}}}(\xi))_{i_{1}i_{2}}}{\xi \! - \! z^{\prime}} \, 
\dfrac{\md \xi}{2 \pi \mi} \right\rvert, \\
\tilde{\mathbb{I}}_{5,\circlearrowright}^{\sharp,\mathrm{C},\flat_{4}}(i_{1},
i_{2}) :=& \, \sup_{z \in \cup_{j=1}^{N+1} \tilde{\Sigma}_{p,j}^{4}} \left\lvert 
\lim_{-\cup_{j=1}^{N+1} \tilde{\Sigma}_{p,j}^{4} \ni z^{\prime} \to z} \int_{
\partial \tilde{\mathbb{U}}_{\tilde{\delta}_{\tilde{a}_{N+1}}}} \dfrac{(w_{+}^{
\Sigma_{\tilde{\mathcal{R}}}}(\xi))_{i_{1}i_{2}}}{\xi \! - \! z^{\prime}} \, 
\dfrac{\md \xi}{2 \pi \mi} \right\rvert, \\
\tilde{\mathbb{I}}_{5,\circlearrowright}^{\sharp,\mathrm{C},\flat_{5}}(i_{1},
i_{2}) :=& \, \sup_{z \in (\cup_{j=1}^{N+1} \tilde{\Sigma}_{p,j}^{5}) \setminus 
\partial \tilde{\mathbb{U}}_{\tilde{\delta}_{\tilde{a}_{N+1}}}} \left\lvert 
\lim_{-((\cup_{j=1}^{N+1} \tilde{\Sigma}_{p,j}^{5}) \setminus \partial \tilde{
\mathbb{U}}_{\tilde{\delta}_{\tilde{a}_{N+1}}}) \ni z^{\prime} \to z} \int_{
\partial \tilde{\mathbb{U}}_{\tilde{\delta}_{\tilde{a}_{N+1}}}} \dfrac{(w_{+}^{
\Sigma_{\tilde{\mathcal{R}}}}(\xi))_{i_{1}i_{2}}}{\xi \! - \! z^{\prime}} \, 
\dfrac{\md \xi}{2 \pi \mi} \right\rvert.
\end{align*}
In order to estimate $\tilde{\mathbb{I}}_{5,\circlearrowright}^{\sharp,
\mathrm{C},\flat_{1}}(i_{1},i_{2})$, consider, say,
\begin{equation*}
\tilde{\mathbb{I}}_{5,\circlearrowright}^{\sharp,\mathrm{C},\flat_{1}}(i_{1},
i_{2};-) \! := \! \sup_{z \in (-\infty,\tilde{b}_{0} - \tilde{\delta}_{\tilde{b}_{0}})} 
\left\lvert \lim_{-(-\infty,\tilde{b}_{0}- \tilde{\delta}_{\tilde{b}_{0}}) \ni 
z^{\prime} \to z} \int_{\partial \tilde{\mathbb{U}}_{\tilde{\delta}_{\tilde{a}_{N
+1}}}} \dfrac{(w_{+}^{\Sigma_{\tilde{\mathcal{R}}}}(\xi))_{i_{1}i_{2}}}{\xi \! - \! 
z^{\prime}} \, \dfrac{\md \xi}{2 \pi \mi} \right\rvert:
\end{equation*}
via the asymptotic expansion, in the double-scaling limit $\mathscr{N},
n \! \to \! \infty$ such that $z_{o} \! = \! 1 \! + \! o(1)$, of $w_{+}^{
\Sigma_{\tilde{\mathcal{R}}}}(z)$ for $z \! \in \! \partial \tilde{\mathbb{
U}}_{\tilde{\delta}_{\tilde{a}_{N+1}}}$ given in Equation~\eqref{eqproptila}, 
H\"{o}lder's Inequality for Integrals, and the indefinite integral \pmb{2.558} 
4. on p.~182 of \cite{gradryzh}, one shows that, with $A_{\epsilon}^{(-)}
(\theta,z) \! := \! (\tilde{a}_{N+1} \! - \! z \! + \! \tilde{\delta}_{\tilde{a}_{N
+1}} \cos \theta) \cos \theta \! + \! (\tilde{\delta}_{\tilde{a}_{N+1}} \sin 
\theta \! + \! \epsilon) \sin \theta$, $B_{\epsilon}^{(-)}(\theta,z) \! := \! 
(\tilde{a}_{N+1} \! - \! z \! + \! \tilde{\delta}_{\tilde{a}_{N+1}} \cos \theta) 
\sin \theta \! - \! (\tilde{\delta}_{\tilde{a}_{N+1}} \sin \theta \! + \! \epsilon) 
\cos \theta$, $r_{\epsilon}^{(-)} \! := \! (\tilde{a}_{N+1} \! - \! z)^{2} \! + \! 
\tilde{\delta}_{\tilde{a}_{N+1}}^{2} \! + \! \epsilon^{2}$, $p_{\epsilon}^{(-)} 
\! := \! 2 \tilde{\delta}_{\tilde{a}_{N+1}} \epsilon$, and $q_{\epsilon}^{(-)} \! 
:= \! 2 \tilde{\delta}_{\tilde{a}_{N+1}}(\tilde{a}_{N+1} \! - \! z)$ $(\Rightarrow 
(r_{\epsilon}^{(-)})^{2} \! - \! (p_{\epsilon}^{(-)})^{2} \! - \! (q_{\epsilon}^{
(-)})^{2} \! = \! ((\tilde{a}_{N+1} \! - \! z)^{2} \! + \! (\epsilon^{2} \! - \! 
\tilde{\delta}_{\tilde{a}_{N+1}}^{2})^{2})^{2} \! > \! 0)$, for $n \! \in \! 
\mathbb{N}$ and $k \! \in \! \lbrace 1,2,\dotsc,K \rbrace$ such that 
$\alpha_{p_{\mathfrak{s}}} \! := \! \alpha_{k} \! \neq \! \infty$,
\begin{align*}
\tilde{\mathbb{I}}_{5,\circlearrowright}^{\sharp,\mathrm{C},\flat_{1}}(i_{1},
i_{2};-) \underset{\underset{z_{o}=1+o(1)}{\mathscr{N},n \to \infty}}{\leqslant}& 
\, \dfrac{\tilde{\delta}_{\tilde{a}_{N+1}}}{2 \pi} \sup_{\theta_{0} \in [0,2 \pi]} 
\left\lvert (w_{+}^{\Sigma_{\tilde{\mathcal{R}}}}(\tilde{a}_{N+1} \! + \! \tilde{
\delta}_{\tilde{a}_{N+1}} \me^{\mi \theta_{0}}))_{i_{1}i_{2}} \right\rvert \sup_{z 
\in (-\infty,\tilde{b}_{0}-\tilde{\delta}_{\tilde{b}_{0}})} \left\lvert \lim_{\epsilon 
\downarrow 0} \int_{0}^{2 \pi} \dfrac{(A_{\epsilon}^{(-)}(\theta,z) \! + \! \mi 
B_{\epsilon}^{(-)}(\theta,z))}{r_{\epsilon}^{(-)} \! + \! p_{\epsilon}^{(-)} \sin 
\theta \! + \! q_{\epsilon}^{(-)} \cos \theta} \, \md \theta \right\rvert 
\nonumber \\
\underset{\underset{z_{o}=1+o(1)}{\mathscr{N},n \to \infty}}{\leqslant}& \, 
\dfrac{\sqrt{2} \, \tilde{\delta}_{\tilde{a}_{N+1}}}{\pi} \sup_{\theta_{0} \in [0,
2 \pi]} \left\lvert (w_{+}^{\Sigma_{\tilde{\mathcal{R}}}}(\tilde{a}_{N+1} \! + \! 
\tilde{\delta}_{\tilde{a}_{N+1}} \me^{\mi \theta_{0}}))_{i_{1}i_{2}} \right\rvert 
\sup_{z \in (-\infty,\tilde{b}_{0}-\tilde{\delta}_{\tilde{b}_{0}})} \left(\lim_{
\epsilon \downarrow 0} \left(\dfrac{\tilde{a}_{N+1} \! - \! z \! + \! 2 \tilde{
\delta}_{\tilde{a}_{N+1}} \! + \! \epsilon}{(\tilde{a}_{N+1} \! - \! z)^{2} \! 
+ \! (\epsilon^{2} \! - \! \tilde{\delta}_{\tilde{a}_{N+1}}^{2})^{2}} \right) 
\right. \nonumber \\
\times&\left. \, \underbrace{\left(\left. \tan^{-1} \left(\dfrac{p_{\epsilon}^{
(-)} \! + \! (r_{\epsilon}^{(-)} \! - \! q_{\epsilon}^{(-)}) \tan(\frac{\theta}{2})}{
(\tilde{a}_{N+1} \! - \! z)^{2} \! + \! (\epsilon^{2} \! - \! \tilde{\delta}_{
\tilde{a}_{N+1}}^{2})^{2}} \right) \right\vert_{0}^{2 \pi} \right)}_{= \, 0} 
\right) \quad \Rightarrow
\end{align*}
\begin{equation*}
\tilde{\mathbb{I}}_{5,\circlearrowright}^{\sharp,\mathrm{C},\flat_{1}}(i_{1},
i_{2};-) \underset{\underset{z_{o}=1+o(1)}{\mathscr{N},n \to \infty}}{=} 0;
\end{equation*}
analogously, one shows that
\begin{align*}
\tilde{\mathbb{I}}_{5,\circlearrowright}^{\sharp,\mathrm{C},\flat_{1}}(i_{1},
i_{2};+) :=& \, \sup_{z \in (\tilde{a}_{N+1}+\tilde{\delta}_{\tilde{a}_{N+1}},
+\infty)} \left\lvert \lim_{-(\tilde{a}_{N+1}+\tilde{\delta}_{\tilde{a}_{N+1}},
+\infty) \ni z^{\prime} \to z} \int_{\partial \tilde{\mathbb{U}}_{\tilde{\delta}_{
\tilde{a}_{N+1}}}} \dfrac{(w_{+}^{\Sigma_{\tilde{\mathcal{R}}}}(\xi))_{i_{1}
i_{2}}}{\xi \! - \! z^{\prime}} \, \dfrac{\md \xi}{2 \pi \mi} \right\rvert \\
\underset{\underset{z_{o}=1+o(1)}{\mathscr{N},n \to \infty}}{\leqslant}& \, 
\dfrac{\sqrt{2} \, \tilde{\delta}_{\tilde{a}_{N+1}}}{\pi} \sup_{\theta_{0} \in [0,
2 \pi]} \left\lvert (w_{+}^{\Sigma_{\tilde{\mathcal{R}}}}(\tilde{a}_{N+1} \! + \! 
\tilde{\delta}_{\tilde{a}_{N+1}} \me^{\mi \theta_{0}}))_{i_{1}i_{2}} \right\rvert 
\sup_{z \in (\tilde{a}_{N+1}+\tilde{\delta}_{\tilde{a}_{N+1}},+\infty)} \left(
\lim_{\epsilon \downarrow 0} \left(\dfrac{\tilde{a}_{N+1} \! - \! z \! + \! 2 
\tilde{\delta}_{\tilde{a}_{N+1}} \! + \! \epsilon}{(\tilde{a}_{N+1} \! - \! z)^{2} 
\! + \! (\epsilon^{2} \! - \! \tilde{\delta}_{\tilde{a}_{N+1}}^{2})^{2}} \right) 
\right. \nonumber \\
\times&\left. \, \underbrace{\left(\left. \tan^{-1} \left(\dfrac{p_{\epsilon}^{
(-)} \! + \! (r_{\epsilon}^{(-)} \! - \! q_{\epsilon}^{(-)}) \tan(\frac{\theta}{2})}{
(\tilde{a}_{N+1} \! - \! z)^{2} \! + \! (\epsilon^{2} \! - \! \tilde{\delta}_{
\tilde{a}_{N+1}}^{2})^{2}} \right) \right\vert_{0}^{2 \pi} \right)}_{= \, 0} 
\right) \quad \Rightarrow
\end{align*}
\begin{equation*}
\tilde{\mathbb{I}}_{5,\circlearrowright}^{\sharp,\mathrm{C},\flat_{1}}(i_{1},
i_{2};+) \underset{\underset{z_{o}=1+o(1)}{\mathscr{N},n \to \infty}}{=} 0,
\end{equation*}
whence, via an elementary inequality argument, one deduces that, for $n 
\! \in \! \mathbb{N}$ and $k \! \in \! \lbrace 1,2,\dotsc,K \rbrace$ such 
that $\alpha_{p_{\mathfrak{s}}} \! := \! \alpha_{k} \! \neq \! \infty$, and 
$i_{1},i_{2} \! = \! 1,2$,
\begin{equation} \label{eqiy23}
\tilde{\mathbb{I}}_{5,\circlearrowright}^{\sharp,\mathrm{C},\flat_{1}}(i_{1},
i_{2}) \underset{\underset{z_{o}=1+o(1)}{\mathscr{N},n \to \infty}}{=} 0.
\end{equation}
Proceeding, \emph{mutatis mutandis}, as above (for $\tilde{\mathbb{I}}_{5,
\circlearrowright}^{\sharp,\mathrm{C},\flat_{1}}(i_{1},i_{2}))$, one also 
shows that, for $n \! \in \! \mathbb{N}$ and $k \! \in \! \lbrace 1,2,\dotsc,
K \rbrace$ such that $\alpha_{p_{\mathfrak{s}}} \! := \! \alpha_{k} \! \neq 
\! \infty$, and $i_{1},i_{2} \! = \! 1,2$,
\begin{equation} \label{eqiy24}
\tilde{\mathbb{I}}_{5,\circlearrowright}^{\sharp,\mathrm{C},\flat_{2}}(i_{1},
i_{2}) \underset{\underset{z_{o}=1+o(1)}{\mathscr{N},n \to \infty}}{=} 0.
\end{equation}
In order to estimate $\tilde{\mathbb{I}}_{5,\circlearrowright}^{\sharp,
\mathrm{C},\flat_{3}}(i_{1},i_{2})$, consider, say, for $j \! \in \! \lbrace 
1,2,\dotsc,N \! + \! 1 \rbrace$,
\begin{equation*}
\tilde{\mathbb{I}}_{5,\circlearrowright}^{\sharp,\mathrm{C},\flat_{3}}(i_{1},
i_{2};\smallfrown \! (j)) \! := \! \sup_{z \in \tilde{\Sigma}_{p,j}^{3}} \left\lvert 
\lim_{-\tilde{\Sigma}_{p,j}^{3} \ni z^{\prime} \to z} \int_{\partial \tilde{
\mathbb{U}}_{\tilde{\delta}_{\tilde{a}_{N+1}}}} \dfrac{(w_{+}^{\Sigma_{
\tilde{\mathcal{R}}}}(\xi))_{i_{1}i_{2}}}{\xi \! - \! z^{\prime}} \, \dfrac{\md 
\xi}{2 \pi \mi} \right\rvert:
\end{equation*}
via the asymptotic expansion, in the double-scaling limit $\mathscr{N},
n \! \to \! \infty$ such that $z_{o} \! = \! 1 \! + \! o(1)$, of $w_{+}^{
\Sigma_{\tilde{\mathcal{R}}}}(z)$ for $z \! \in \! \partial \tilde{\mathbb{
U}}_{\tilde{\delta}_{\tilde{a}_{N+1}}}$ given in Equation~\eqref{eqproptila}, 
H\"{o}lder's Inequality for Integrals, the indefinite integral \pmb{2.558} 
4. on p.~182 of \cite{gradryzh}, and the parametrisation for 
the---elliptic---homotopic deformation of $\tilde{\Sigma}_{p,j}^{3}$ 
given in the proof of the corresponding item of Lemma~\ref{lem5.1}, that 
is, $\tilde{\Sigma}_{p,j}^{3} \! = \! \lbrace \mathstrut (x_{j}(\theta),y_{j}
(\theta)); \, x_{j}(\theta) \! = \! \tfrac{1}{2}(\tilde{a}_{j} \! + \! \tilde{b}_{j
-1}) \! + \! \tfrac{1}{2}(\tilde{a}_{j} \! - \! \tilde{b}_{j-1}) \cos \theta, \, 
y_{j}(\theta) \! = \! \tilde{\eta}_{j} \sin \theta, \, \theta_{0}^{\tilde{a}}(j) \! 
\leqslant \! \theta \! \leqslant \! \pi \! - \! \theta_{0}^{\tilde{b}}(j) \rbrace$, 
one shows that, with $\alpha_{\epsilon}^{\smallfrown}(\theta,\psi) \! := \! 
\tilde{a}_{N+1} \! + \! \tilde{\delta}_{\tilde{a}_{N+1}} \cos \theta \! + \! 
\epsilon \cos \psi \! - \! x_{j}(\psi)$, $\beta_{\epsilon}^{\smallfrown}(\theta,
\psi) \! := \! \tilde{\delta}_{\tilde{a}_{N+1}} \sin \theta \! + \! \epsilon 
\sin \psi \! - \! y_{j}(\psi)$, $A_{\epsilon}^{\smallfrown}(\theta,\psi) \! := \! 
\alpha_{\epsilon}^{\smallfrown}(\theta,\psi) \cos \theta \! + \! \beta_{
\epsilon}^{\smallfrown}(\theta,\psi) \sin \theta$, $B_{\epsilon}^{
\smallfrown}(\theta,\psi) \! := \! \alpha_{\epsilon}^{\smallfrown}(\theta,
\psi) \sin \theta \! - \! \beta_{\epsilon}^{\smallfrown}(\theta,\psi) \cos 
\theta$, $\tilde{A}^{\smallfrown}(\epsilon) \! := \! \tilde{a}_{N+1} \! + \! 
2 \tilde{\delta}_{\tilde{a}_{N+1}} \! + \! \epsilon (1 \! + \! \cos \theta_{
0}^{\tilde{a}}(j)) \! - \! \tfrac{1}{2}(\tilde{a}_{j} \! + \! \tilde{b}_{j-1}) \! + 
\! \tfrac{1}{2}(\tilde{a}_{j} \! - \! \tilde{b}_{j-1}) \cos \theta_{0}^{\tilde{a}}
(j) \! - \! \tilde{\eta}_{j} \min \lbrace \sin \theta_{0}^{\tilde{a}}(j),\sin 
\theta_{0}^{\tilde{b}}(j) \rbrace$, $r_{\epsilon}^{\smallfrown}(\psi) \! := 
\! (\tilde{a}_{N+1} \! - \! x_{j}(\psi))^{2} \! + \! 2 \epsilon (\tilde{a}_{N+1} 
\! - \! x_{j}(\psi)) \cos \psi \! + \! y_{j}^{2}(\psi) \! - \! 2 \epsilon y_{j}(\psi) 
\sin \psi \! + \! \tilde{\delta}_{\tilde{a}_{N+1}}^{2} \! + \! \epsilon^{2}$, 
$p_{\epsilon}^{\smallfrown}(\psi) \! := \! 2 \tilde{\delta}_{\tilde{a}_{N+1}}
(\epsilon \sin \psi \! - \! y_{j}(\psi))$, and $q_{\epsilon}^{\smallfrown}(\psi) 
\! := \! 2 \tilde{\delta}_{\tilde{a}_{N+1}}(\tilde{a}_{N+1} \! + \! \epsilon 
\cos \psi \! - \! x_{j}(\psi))$, for $n \! \in \! \mathbb{N}$ and $k \! \in \! 
\lbrace 1,2,\dotsc,K \rbrace$ such that $\alpha_{p_{\mathfrak{s}}} \! := 
\! \alpha_{k} \! \neq \! \infty$, and $j \! \in \! \lbrace 1,2,\dotsc,N \! 
+ \! 1 \rbrace$,\footnote{If $j \! = \! N \! + \! 1$, then $\psi \! \in \! 
(\theta_{0}^{\tilde{a}}(j),\pi \! - \! \theta_{0}^{\tilde{b}}(j))$.}
\begin{align*}
\tilde{\mathbb{I}}_{5,\circlearrowright}^{\sharp,\mathrm{C},\flat_{3}}(i_{1},
i_{2};\smallfrown \! (j)) \underset{\underset{z_{o}=1+o(1)}{\mathscr{N},n 
\to \infty}}{\leqslant}& \, \dfrac{\tilde{\delta}_{\tilde{a}_{N+1}}}{2 \pi} 
\sup_{\theta_{0} \in [0,2 \pi]} \left\lvert (w_{+}^{\Sigma_{\tilde{\mathcal{R}}}}
(\tilde{a}_{N+1} \! + \! \tilde{\delta}_{\tilde{a}_{N+1}} \me^{\mi \theta_{0}}
))_{i_{1}i_{2}} \right\rvert \sup_{\psi \in [\theta_{0}^{\tilde{a}}(j),\pi -\theta_{0}^{
\tilde{b}}(j)]} \left\lvert \lim_{\epsilon \downarrow 0} \int_{0}^{2 \pi} 
\dfrac{(A_{\epsilon}^{\smallfrown}(\theta,\psi) \! + \! \mi B_{\epsilon}^{
\smallfrown}(\theta,\psi))}{(\alpha_{\epsilon}^{\smallfrown}(\theta,\psi))^{2} 
\! + \! (\beta_{\epsilon}^{\smallfrown}(\theta,\psi))^{2}} \, \md \theta 
\right\rvert \nonumber \\
\underset{\underset{z_{o}=1+o(1)}{\mathscr{N},n \to \infty}}{\leqslant}& \, 
\dfrac{\tilde{\delta}_{\tilde{a}_{N+1}}}{\sqrt{2} \, \pi} \sup_{\theta_{0} \in [0,
2 \pi]} \left\lvert (w_{+}^{\Sigma_{\tilde{\mathcal{R}}}}(\tilde{a}_{N+1} \! + \! 
\tilde{\delta}_{\tilde{a}_{N+1}} \me^{\mi \theta_{0}}))_{i_{1}i_{2}} \right\rvert 
\sup_{\psi \in [\theta_{0}^{\tilde{a}}(j),\pi -\theta_{0}^{\tilde{b}}(j)]} \left(
\lim_{\epsilon \downarrow 0} \left\lvert \tilde{A}^{\smallfrown}(\epsilon) 
\right\rvert \tilde{\mathbb{J}}^{\smallfrown}(\psi,\epsilon) \right),
\end{align*}
where, with $\mathfrak{d}_{\epsilon}^{\smallfrown}(\psi) \! := \! 
(r_{\epsilon}^{\smallfrown}(\psi))^{2} \! - \! (p_{\epsilon}^{\smallfrown}
(\psi))^{2} \! - \! (q_{\epsilon}^{\smallfrown}(\psi))^{2}$,
\begin{align*}
\tilde{\mathbb{J}}^{\smallfrown}(\psi,\epsilon) =& \, 
\begin{cases}
\frac{2}{\sqrt{\mathfrak{d}_{\epsilon}^{\smallfrown}(\psi)}} \left. \tan^{-1} 
\left(\frac{p_{\epsilon}^{\smallfrown}(\psi)+(r_{\epsilon}^{\smallfrown}(\psi)-
q_{\epsilon}^{\smallfrown}(\psi)) \tan (\frac{\theta}{2})}{\sqrt{\mathfrak{d}_{
\epsilon}^{\smallfrown}(\psi)}} \right) \right\vert_{0}^{2 \pi}=0, 
&\text{$\mathfrak{d}_{\epsilon}^{\smallfrown}(\psi) \! > \! 0$,} \\
\frac{1}{\sqrt{-\mathfrak{d}_{\epsilon}^{\smallfrown}(\psi)}} \left. \ln \left(
\frac{p_{\epsilon}^{\smallfrown}(\psi)-\sqrt{-\mathfrak{d}_{\epsilon}^{
\smallfrown}(\psi)}+(r_{\epsilon}^{\smallfrown}(\psi)-q_{\epsilon}^{
\smallfrown}(\psi)) \tan (\frac{\theta}{2})}{p_{\epsilon}^{\smallfrown}(\psi)
+\sqrt{-\mathfrak{d}_{\epsilon}^{\smallfrown}(\psi)}+(r_{\epsilon}^{\smallfrown}
(\psi)-q_{\epsilon}^{\smallfrown}(\psi)) \tan (\frac{\theta}{2})} \right) 
\right\vert_{0}^{2 \pi}=0, &\text{$\mathfrak{d}_{\epsilon}^{\smallfrown}
(\psi) \! < \! 0$,}
\end{cases} \quad \quad \Rightarrow
\end{align*}
\begin{equation*}
\tilde{\mathbb{I}}_{5,\circlearrowright}^{\sharp,\mathrm{C},\flat_{3}}(i_{1},
i_{2};\smallfrown \! (j)) \underset{\underset{z_{o}=1+o(1)}{\mathscr{N},n 
\to \infty}}{=} 0, \quad j \! \in \! \lbrace 1,2,\dotsc,N \! + \! 1 \rbrace,
\end{equation*}
whence, via an elementary inequality argument, one deduces that, for $n 
\! \in \! \mathbb{N}$ and $k \! \in \! \lbrace 1,2,\dotsc,K \rbrace$ such 
that $\alpha_{p_{\mathfrak{s}}} \! := \! \alpha_{k} \! \neq \! \infty$, and 
$i_{1},i_{2} \! = \! 1,2$,
\begin{equation} \label{eqiy25} 
\tilde{\mathbb{I}}_{5,\circlearrowright}^{\sharp,\mathrm{C},\flat_{3}}(i_{1},
i_{2}) \underset{\underset{z_{o}=1+o(1)}{\mathscr{N},n \to \infty}}{=} 0.
\end{equation}
Proceeding, \emph{mutatis mutandis}, as above (for $\tilde{\mathbb{I}}_{5,
\circlearrowright}^{\sharp,\mathrm{C},\flat_{3}}(i_{1},i_{2}))$, one shows 
that, for $n \! \in \! \mathbb{N}$ and $k \! \in \! \lbrace 1,2,\dotsc,K 
\rbrace$ such that $\alpha_{p_{\mathfrak{s}}} \! := \! \alpha_{k} \! \neq 
\! \infty$, and $i_{1},i_{2} \! = \! 1,2$,
\begin{equation} \label{eqiy26} 
\tilde{\mathbb{I}}_{5,\circlearrowright}^{\sharp,\mathrm{C},\flat_{4}}(i_{1},
i_{2}) \underset{\underset{z_{o}=1+o(1)}{\mathscr{N},n \to \infty}}{=} 0.
\end{equation}
In order to estimate $\tilde{\mathbb{I}}_{5,\circlearrowright}^{\sharp,
\mathrm{C},\flat_{5}}(i_{1},i_{2})$, consider, say, for $j \! \in \! \lbrace 
1,2,\dotsc,N \! + \! 1 \rbrace$,
\begin{equation*}
\tilde{\mathbb{I}}_{5,\circlearrowright}^{\sharp,\mathrm{C},\flat_{5}}(i_{1},
i_{2};\circlearrowright \! (j)) \! := \! \sup_{z \in \partial \tilde{\mathbb{U}}_{
\tilde{\delta}_{\tilde{b}_{j-1}}}} \left\lvert \lim_{-\tilde{\mathbb{U}}_{\tilde{
\delta}_{\tilde{b}_{j-1}}} \ni z^{\prime} \to z} \int_{\partial \tilde{\mathbb{
U}}_{\tilde{\delta}_{\tilde{a}_{N+1}}}} \dfrac{(w_{+}^{\Sigma_{\tilde{\mathcal{
R}}}}(\xi))_{i_{1}i_{2}}}{\xi \! - \! z^{\prime}} \, \dfrac{\md \xi}{2 \pi \mi} 
\right\rvert:
\end{equation*}
via the asymptotic expansion, in the double-scaling limit $\mathscr{N},
n \! \to \! \infty$ such that $z_{o} \! = \! 1 \! + \! o(1)$, of $w_{+}^{
\Sigma_{\tilde{\mathcal{R}}}}(z)$ for $z \! \in \! \partial \tilde{\mathbb{
U}}_{\tilde{\delta}_{\tilde{a}_{N+1}}}$ given in Equation~\eqref{eqproptila}, 
H\"{o}lder's Inequality for Integrals, and the indefinite integral \pmb{2.558} 
4. on p.~182 of \cite{gradryzh}, one shows that, with $\alpha_{\epsilon}^{
\circlearrowright}(\theta,\psi) \! := \! \tilde{a}_{N+1} \! - \! \tilde{b}_{j-1} 
\! + \! \tilde{\delta}_{\tilde{a}_{N+1}} \cos \theta \! - \! (\tilde{\delta}_{
\tilde{b}_{j-1}} \! - \! \epsilon) \cos \psi$, $\beta_{\epsilon}^{
\circlearrowright}(\theta,\psi) \! := \! \tilde{\delta}_{\tilde{a}_{N+1}} \sin 
\theta \! - \! (\tilde{\delta}_{\tilde{b}_{j-1}} \! - \! \epsilon) \sin \psi$, 
$A_{\epsilon}^{\circlearrowright}(\theta,\psi) \! := \! \alpha_{\epsilon}^{
\circlearrowright}(\theta,\psi) \cos \theta \! + \! \beta_{\epsilon}^{
\circlearrowright}(\theta,\psi) \sin \theta$, $B_{\epsilon}^{
\circlearrowright}(\theta,\psi) \! := \! \alpha_{\epsilon}^{\circlearrowright}
(\theta,\psi) \sin \theta \! - \! \beta_{\epsilon}^{\circlearrowright}(\theta,
\psi) \cos \theta$, $\tilde{A}^{\circlearrowright}(\epsilon) \! := \! 
\tilde{a}_{N+1} \! - \! \tilde{b}_{j-1} \! + \! 2(\tilde{\delta}_{\tilde{a}_{N+1}} 
\! + \! \tilde{\delta}_{\tilde{b}_{j-1}}) \! - \! 2 \epsilon$, $r_{\epsilon}^{
\circlearrowright}(\psi) \! := \! (\tilde{a}_{N+1} \! - \! \tilde{b}_{j-1})^{2} 
\! - \! 2(\tilde{a}_{N+1} \! - \! \tilde{b}_{j-1})(\tilde{\delta}_{\tilde{b}_{j-1}} 
\! - \! \epsilon) \cos \psi \! + \! \tilde{\delta}_{\tilde{a}_{N+1}}^{2} \! + \! 
\tilde{\delta}_{\tilde{b}_{j-1}}(\tilde{\delta}_{\tilde{b}_{j-1}} \! - \! 2 \epsilon) 
\! + \! \epsilon^{2}$, $p_{\epsilon}^{\circlearrowright}(\psi) \! := \! 2 
\tilde{\delta}_{\tilde{a}_{N+1}}(\epsilon \! - \! \tilde{\delta}_{\tilde{b}_{j-1}}) 
\sin \psi$, and $q_{\epsilon}^{\circlearrowright}(\psi) \! := \! 2 \tilde{
\delta}_{\tilde{a}_{N+1}}(\tilde{a}_{N+1} \! - \! \tilde{b}_{j-1}) \! - \! 2 
\tilde{\delta}_{\tilde{a}_{N+1}}(\tilde{\delta}_{\tilde{b}_{j-1}} \! - \! \epsilon) 
\cos \psi$, for $n \! \in \! \mathbb{N}$ and $k \! \in \! \lbrace 1,2,\dotsc,
K \rbrace$ such that $\alpha_{p_{\mathfrak{s}}} \! := \! \alpha_{k} \! \neq 
\! \infty$, and $j \! \in \! \lbrace 1,2,\dotsc,N \! + \! 1 \rbrace$,
\begin{align*}
\tilde{\mathbb{I}}_{5,\circlearrowright}^{\sharp,\mathrm{C},\flat_{5}}(i_{1},
i_{2};\circlearrowright \! (j)) \underset{\underset{z_{o}=1+o(1)}{\mathscr{N},
n \to \infty}}{\leqslant}& \, \dfrac{\tilde{\delta}_{\tilde{a}_{N+1}}}{2 \pi} 
\sup_{\theta_{0} \in [0,2 \pi]} \left\lvert (w_{+}^{\Sigma_{\tilde{\mathcal{R}}}}
(\tilde{a}_{N+1} \! + \! \tilde{\delta}_{\tilde{a}_{N+1}} \me^{\mi \theta_{0}}
))_{i_{1}i_{2}} \right\rvert \sup_{\psi \in [0,2 \pi]} \left\lvert \lim_{\epsilon 
\downarrow 0} \int_{0}^{2 \pi} \dfrac{(A_{\epsilon}^{\circlearrowright}
(\theta,\psi) \! + \! \mi B_{\epsilon}^{\circlearrowright}(\theta,\psi))}{
(\alpha_{\epsilon}^{\circlearrowright}(\theta,\psi))^{2} \! + \! (\beta_{
\epsilon}^{\circlearrowright}(\theta,\psi))^{2}} \, \md \theta \right\rvert 
\nonumber \\
\underset{\underset{z_{o}=1+o(1)}{\mathscr{N},n \to \infty}}{\leqslant}& \, 
\dfrac{\tilde{\delta}_{\tilde{a}_{N+1}}}{\sqrt{2} \, \pi} \sup_{\theta_{0} \in [0,
2 \pi]} \left\lvert (w_{+}^{\Sigma_{\tilde{\mathcal{R}}}}(\tilde{a}_{N+1} \! + \! 
\tilde{\delta}_{\tilde{a}_{N+1}} \me^{\mi \theta_{0}}))_{i_{1}i_{2}} \right\rvert 
\sup_{\psi \in [0,2 \pi]} \left(\lim_{\epsilon \downarrow 0} \left\lvert 
\tilde{A}^{\circlearrowright}(\epsilon) \right\rvert \tilde{\mathbb{J}}^{
\circlearrowright}(\psi,\epsilon) \right),
\end{align*}
where, with $\mathfrak{d}_{\epsilon}^{\circlearrowright}(\psi) \! := \! 
(r_{\epsilon}^{\circlearrowright}(\psi))^{2} \! - \! (p_{\epsilon}^{\circlearrowright}
(\psi))^{2} \! - \! (q_{\epsilon}^{\circlearrowright}(\psi))^{2}$,
\begin{align*}
\tilde{\mathbb{J}}^{\circlearrowright}(\psi,\epsilon) =& \, 
\begin{cases}
\frac{2}{\sqrt{\mathfrak{d}_{\epsilon}^{\circlearrowright}(\psi)}} \left. \tan^{-1} 
\left(\frac{p_{\epsilon}^{\circlearrowright}(\psi)+(r_{\epsilon}^{\circlearrowright}
(\psi)-q_{\epsilon}^{\circlearrowright}(\psi)) \tan (\frac{\theta}{2})}{\sqrt{
\mathfrak{d}_{\epsilon}^{\circlearrowright}(\psi)}} \right) \right\vert_{0}^{2 \pi}
=0, &\text{$\mathfrak{d}_{\epsilon}^{\circlearrowright}(\psi) \! > \! 0$,} \\
\frac{1}{\sqrt{-\mathfrak{d}_{\epsilon}^{\circlearrowright}(\psi)}} \left. \ln 
\left(\frac{p_{\epsilon}^{\circlearrowright}(\psi)-\sqrt{-\mathfrak{d}_{
\epsilon}^{\circlearrowright}(\psi)}+(r_{\epsilon}^{\circlearrowright}
(\psi)-q_{\epsilon}^{\circlearrowright}(\psi)) \tan (\frac{\theta}{2})}{
p_{\epsilon}^{\circlearrowright}(\psi)+\sqrt{-\mathfrak{d}_{\epsilon}^{
\circlearrowright}(\psi)}+(r_{\epsilon}^{\circlearrowright}(\psi)-q_{\epsilon}^{
\circlearrowright}(\psi)) \tan (\frac{\theta}{2})} \right) \right\vert_{0}^{2 \pi}
=0, &\text{$\mathfrak{d}_{\epsilon}^{\circlearrowright}(\psi) \! < \! 0$,}
\end{cases} \quad \quad \Rightarrow
\end{align*}
\begin{equation*}
\tilde{\mathbb{I}}_{5,\circlearrowright}^{\sharp,\mathrm{C},\flat_{5}}(i_{1},
i_{2};\circlearrowright \! (j)) \underset{\underset{z_{o}=1+o(1)}{\mathscr{N},
n \to \infty}}{=} 0, \quad j \! \in \! \lbrace 1,2,\dotsc,N \! + \! 1 \rbrace;
\end{equation*}
analogously, one shows that
\begin{equation*}
\sup_{z \in \partial \tilde{\mathbb{U}}_{\tilde{\delta}_{\tilde{a}_{j}}}} \left\lvert 
\lim_{-\partial \tilde{\mathbb{U}}_{\tilde{\delta}_{\tilde{a}_{j}}} \ni z^{\prime} 
\to z} \int_{\partial \tilde{\mathbb{U}}_{\tilde{\delta}_{\tilde{a}_{N+1}}}} \dfrac{(
w_{+}^{\Sigma_{\tilde{\mathcal{R}}}}(\xi))_{i_{1}i_{2}}}{\xi \! - \! z^{\prime}} 
\, \dfrac{\md \xi}{2 \pi \mi} \right\rvert \underset{\underset{z_{o}=1+o(1)}{
\mathscr{N},n \to \infty}}{=} 0, \quad j \! \in \! \lbrace 1,2,\dotsc,N \rbrace,
\end{equation*}
whence, via an elementary inequality argument, one arrives at, for $n \! \in \! 
\mathbb{N}$ and $k \! \in \! \lbrace 1,2,\dotsc,K \rbrace$ such that $\alpha_{
p_{\mathfrak{s}}} \! := \! \alpha_{k} \! \neq \! \infty$, and $i_{1},i_{2} \! = \! 
1,2$,
\begin{equation} \label{eqiy27} 
\tilde{\mathbb{I}}_{5,\circlearrowright}^{\sharp,\mathrm{C},\flat_{5}}(i_{1},
i_{2}) \underset{\underset{z_{o}=1+o(1)}{\mathscr{N},n \to \infty}}{=} 0.
\end{equation}
Thus, via Inequality~\eqref{eqiy22} and the 
Estimates~\eqref{eqiy23}--\eqref{eqiy27}, one arrives at, for $n \! \in \! 
\mathbb{N}$ and $k \! \in \! \lbrace 1,2,\dotsc,K \rbrace$ such that 
$\alpha_{p_{\mathfrak{s}}} \! := \! \alpha_{k} \! \neq \! \infty$, and 
$i_{1},i_{2} \! = \! 1,2$,
\begin{equation} \label{eqiy28} 
\tilde{\mathbb{I}}_{5,\circlearrowright}^{\sharp,\mathrm{C}}(i_{1},i_{2}) 
\underset{\underset{z_{o}=1+o(1)}{\mathscr{N},n \to \infty}}{=} 0,
\end{equation}
whence, via Inequality~\eqref{eqiy19} and the Estimates~\eqref{eqiy20}, 
\eqref{eqiy21}, and~\eqref{eqiy28},
\begin{equation} \label{eqiy29} 
\tilde{\mathbb{I}}_{5,\circlearrowright}^{\sharp}(i_{1},i_{2}) \underset{
\underset{z_{o}=1+o(1)}{\mathscr{N},n \to \infty}}{=} 0, \quad i_{1},i_{2} 
\! = \! 1,2.
\end{equation}
Proceeding, \emph{mutatis mutandis}, as above (for $\tilde{\mathbb{I}}_{5,
\circlearrowright}^{\sharp}(i_{1},i_{2})$), one shows that, for $n \! \in \! 
\mathbb{N}$ and $k \! \in \! \lbrace 1,2,\dotsc,K \rbrace$ such that 
$\alpha_{p_{\mathfrak{s}}} \! := \! \alpha_{k} \! \neq \! \infty$, and 
$i_{1},i_{2} \! = \! 1,2$,
\begin{gather*}
\sup_{z \in \tilde{\Sigma}_{\tilde{\mathcal{R}}}^{\sharp}} \left\lvert 
\lim_{-\tilde{\Sigma}_{\tilde{\mathcal{R}}}^{\sharp} \ni z^{\prime} 
\to z} \int_{\partial \tilde{\mathbb{U}}_{\tilde{\delta}_{\tilde{a}_{j}}}} 
\dfrac{(z^{\prime} \! - \! \alpha_{k})(w_{+}^{\Sigma_{\tilde{\mathcal{R}}}}
(\xi))_{i_{1}i_{2}}}{(\xi \! - \! \alpha_{k})(\xi \! - \! z^{\prime})} \, 
\dfrac{\md \xi}{2 \pi \mi} \right\rvert \underset{\underset{z_{o}=1+
o(1)}{\mathscr{N},n \to \infty}}{=} 0, \quad j \! \in \! \lbrace 1,2,\dotsc,
N \rbrace, \\
\sup_{z \in \tilde{\Sigma}_{\tilde{\mathcal{R}}}^{\sharp}} \left\lvert 
\lim_{-\tilde{\Sigma}_{\tilde{\mathcal{R}}}^{\sharp} \ni z^{\prime} 
\to z} \int_{\partial \tilde{\mathbb{U}}_{\tilde{\delta}_{\tilde{b}_{m-1}}}} 
\dfrac{(z^{\prime} \! - \! \alpha_{k})(w_{+}^{\Sigma_{\tilde{\mathcal{R}}}}
(\xi))_{i_{1}i_{2}}}{(\xi \! - \! \alpha_{k})(\xi \! - \! z^{\prime})} \, 
\dfrac{\md \xi}{2 \pi \mi} \right\rvert \underset{\underset{z_{o}=1
+o(1)}{\mathscr{N},n \to \infty}}{=} 0, \quad m \! \in \! \lbrace 1,2,
\dotsc,N \! + \! 1 \rbrace,
\end{gather*}
whence, via an elementary inequality argument, for $n \! \in \! \mathbb{N}$ 
and $k \! \in \! \lbrace 1,2,\dotsc,K \rbrace$ such that $\alpha_{p_{\mathfrak{s}}} 
\! := \! \alpha_{k} \! \neq \! \infty$,
\begin{equation} \label{eqiy30} 
\tilde{\mathbb{I}}_{5}^{\sharp}(i_{1},i_{2}) \underset{\underset{z_{o}=
1+o(1)}{\mathscr{N},n \to \infty}}{=} 0, \quad i_{1},i_{2} \! = \! 1,2.
\end{equation}
Finally, via the Estimates~\eqref{eqiy7}, \eqref{eqiy9}, \eqref{eqiy10}, 
\eqref{eqiy18}, and~\eqref{eqiy30}, one arrives at, for $n \! \in \! 
\mathbb{N}$ and $k \! \in \! \lbrace 1,2,\dotsc,K \rbrace$ such that 
$\alpha_{p_{\mathfrak{s}}} \! := \! \alpha_{k} \! \neq \! \infty$,
\begin{align*}
\lvert \lvert (C_{w^{\Sigma_{\tilde{\mathcal{R}}}}}f)(\pmb{\cdot}) 
\rvert \rvert_{\mathcal{L}^{\infty}_{\mathrm{M}_{2}(\mathbb{C})}
(\tilde{\Sigma}_{\tilde{\mathcal{R}}}^{\sharp})} \underset{\underset{
z_{o}=1+o(1)}{\mathscr{N},n \to \infty}}{\leqslant}& \, \lvert \lvert 
f(\pmb{\cdot}) \rvert \rvert_{\mathcal{L}^{\infty}_{\mathrm{M}_{2}
(\mathbb{C})}(\tilde{\Sigma}_{\tilde{\mathcal{R}}}^{\sharp})} 
\max_{i_{1},i_{2}=1,2} \left\{\mathcal{O} \left(\tilde{\mathfrak{c}}_{i_{1}
i_{2}}^{\sharp,1}(n,k,z_{o}) \me^{-((n-1)K+k) \min\limits_{j=1,2,\dotsc,
N} \lbrace \tilde{\lambda}_{\tilde{\mathcal{R}},1}^{\sharp}(j) \rbrace} 
\right) \right. \\
+&\left. \, \mathcal{O} \left(\tilde{\mathfrak{c}}_{i_{1}i_{2}}^{\sharp,
2}(n,k,z_{o}) \me^{-((n-1)K+k) \tilde{\lambda}_{\tilde{\mathcal{R}},
2}^{\sharp}(+)} \right) \! + \! \mathcal{O} \left(\tilde{\mathfrak{c}}_{
i_{1}i_{2}}^{\sharp,3}(n,k,z_{o}) \me^{-((n-1)K+k) \tilde{\lambda}_{
\tilde{\mathcal{R}},2}^{\sharp}(-)} \right) \right. \\
+&\left. \, \mathcal{O} \left(\tilde{\mathfrak{c}}_{i_{1}i_{2}}^{\sharp,4}
(n,k,z_{o}) \me^{-((n-1)K+k) \min \left\lbrace \min\limits_{j=1,2,
\dotsc,N+1} \lbrace \tilde{\lambda}_{\tilde{\mathcal{R}},4}^{\sharp,
\smallfrown}(j) \rbrace,\min\limits_{j=1,2,\dotsc,N+1} \lbrace 
\tilde{\lambda}_{\tilde{\mathcal{R}},4}^{\sharp,\smallsmile}(j) \rbrace 
\right\rbrace} \right) \right\} \\
\underset{\underset{z_{o}=1+o(1)}{\mathscr{N},n \to \infty}}{\leqslant}& 
\, \lvert \lvert f(\pmb{\cdot}) \rvert \rvert_{\mathcal{L}^{\infty}_{
\mathrm{M}_{2}(\mathbb{C})}(\tilde{\Sigma}_{\tilde{\mathcal{R}}}^{
\sharp})} \mathcal{O} \left(\tilde{\mathfrak{c}}_{\tilde{\mathcal{R}},
\tilde{w}}^{\triangleright} \me^{-((n-1)K+k) \tilde{\lambda}_{\tilde{
\mathcal{R}},\tilde{w}}^{\triangleright}} \right),
\end{align*}
where $\tilde{\lambda}_{\tilde{\mathcal{R}},\tilde{w}}^{\triangleright}$ 
$(> \! 0)$ is defined by Equation~\eqref{tillamrwr2}, and $\tilde{
\mathfrak{c}}_{\tilde{\mathcal{R}},\tilde{w}}^{\triangleright}$ is 
characterised in the corresponding item~\pmb{(2)} of the lemma; 
consequently,
\begin{equation*}
\lvert \lvert C_{w^{\Sigma_{\tilde{\mathcal{R}}}}} \rvert \rvert_{
\mathfrak{B}_{\infty}(\tilde{\Sigma}_{\tilde{\mathcal{R}}}^{\sharp})} 
\underset{\underset{z_{o}=1+o(1)}{\mathscr{N},n \to \infty}}{=} 
\mathcal{O} \left(\tilde{\mathfrak{c}}_{\tilde{\mathcal{R}},\tilde{w}}^{
\triangleright} \me^{-((n-1)K+k) \tilde{\lambda}_{\tilde{\mathcal{R}},
\tilde{w}}^{\triangleright}} \right),
\end{equation*}
which is the Estimate~\eqref{eqcwsigtl}.

For $n \! \in \! \mathbb{N}$ and $k \! \in \! \lbrace 1,2,\dotsc,K \rbrace$ 
such that $\alpha_{p_{\mathfrak{s}}} \! := \! \alpha_{k} \! \neq \! \infty$, 
via the above definitions, H\"{o}lder's Inequality, H\"{o}lder's Inequality for 
Integrals, and a partial-fraction decomposition argument, it follows that
\begin{align}
\lvert \lvert (C_{w^{\Sigma_{\tilde{\mathcal{R}}}}}f)(\pmb{\cdot}) \rvert 
\rvert_{\mathcal{L}^{2}_{\mathrm{M}_{2}(\mathbb{C})}(\tilde{\Sigma}_{
\tilde{\mathcal{R}}}^{\sharp})} =& \, \left(\int_{\tilde{\Sigma}_{\tilde{
\mathcal{R}}}^{\sharp}} \lvert (C_{w^{\Sigma_{\tilde{\mathcal{R}}}}}f)(\xi) 
\rvert^{2} \, \lvert \md \xi \rvert \right)^{1/2} \! = \! \left(\int_{\tilde{
\Sigma}_{\tilde{\mathcal{R}}}^{\sharp}} \sum_{i_{1},i_{2}=1,2} 
\overline{((C_{w^{\Sigma_{\tilde{\mathcal{R}}}}}f)(\xi))_{i_{1}i_{2}}}
((C_{w^{\Sigma_{\tilde{\mathcal{R}}}}}f)(\xi))_{i_{1}i_{2}} \, \lvert \md 
\xi \rvert \right)^{1/2} \nonumber \\
=& \, \left(\int_{\tilde{\Sigma}_{\tilde{\mathcal{R}}}^{\sharp}} \sum_{i_{1},
i_{2}=1,2} \overline{\left(\lim_{-\tilde{\Sigma}_{\tilde{\mathcal{R}}}^{
\sharp} \ni z^{\prime} \to \xi} \int_{\tilde{\Sigma}_{\tilde{\mathcal{R}}}^{
\sharp}} \dfrac{(z^{\prime} \! - \! \alpha_{k})((fw_{+}^{\Sigma_{\tilde{
\mathcal{R}}}})(u))_{i_{1}i_{2}}}{(u \! - \! \alpha_{k})(u \! - \! z^{\prime})} \, 
\dfrac{\md u}{2 \pi \mi} \right)} \right. \nonumber \\
\times&\left. \, \left(\lim_{-\tilde{\Sigma}_{\tilde{\mathcal{R}}}^{\sharp} 
\ni z^{\prime} \to \xi} \int_{\tilde{\Sigma}_{\tilde{\mathcal{R}}}^{\sharp}} 
\dfrac{(z^{\prime} \! - \! \alpha_{k})((fw_{+}^{\Sigma_{\tilde{\mathcal{R}}}})
(v))_{i_{1}i_{2}}}{(v \! - \! \alpha_{k})(v \! - \! z^{\prime})} \, \dfrac{\md 
v}{2 \pi \mi} \right) \, \lvert \md \xi \rvert \right)^{1/2} \nonumber \\
\leqslant& \, \lvert \lvert f(\pmb{\cdot}) \rvert \rvert_{\mathcal{L}^{2}_{
\mathrm{M}_{2}(\mathbb{C})}(\tilde{\Sigma}_{\tilde{\mathcal{R}}}^{
\sharp})} \left(\int_{\tilde{\Sigma}_{\tilde{\mathcal{R}}}^{\sharp}} \left(
\lim_{-\tilde{\Sigma}_{\tilde{\mathcal{R}}}^{\sharp} \ni z^{\prime} \to 
\xi} \left\lvert \left\lvert \dfrac{w_{+}^{\Sigma_{\tilde{\mathcal{R}}}}
(\pmb{\cdot})}{2 \pi \mi (\pmb{\cdot} \! - \! z^{\prime})} \! - \! \dfrac{
w_{+}^{\Sigma_{\tilde{\mathcal{R}}}}(\pmb{\cdot})}{2 \pi \mi (\pmb{\cdot} 
\! - \! \alpha_{k})} \right\rvert \right\rvert_{\mathcal{L}^{2}_{\mathrm{M}_{2}
(\mathbb{C})}(\tilde{\Sigma}_{\tilde{\mathcal{R}}}^{\sharp})}^{2} \right) \, 
\lvert \md \xi \rvert \right)^{1/2}. \label{eqiy31}
\end{align}
For $n \! \in \! \mathbb{N}$ and $k \! \in \! \lbrace 1,2,\dotsc,K \rbrace$ 
such that $\alpha_{p_{\mathfrak{s}}} \! := \! \alpha_{k} \! \neq \! \infty$, 
via the above definitions and the disjointness of the oriented skeleton 
$\tilde{\Sigma}_{\tilde{\mathcal{R}}}^{\sharp}$, a straightforward calculation 
shows that
\begin{equation} \label{eqiy32}
\left\lvert \left\lvert \dfrac{w_{+}^{\Sigma_{\tilde{\mathcal{R}}}}(\pmb{
\cdot})}{2 \pi \mi (\pmb{\cdot} \! - \! \alpha_{k})} \right\rvert \right\rvert_{
\mathcal{L}^{2}_{\mathrm{M}_{2}(\mathbb{C})}(\tilde{\Sigma}_{\tilde{
\mathcal{R}}}^{\sharp})}^{2} \! = \! \tilde{\mathbb{J}}_{A} \! + \! \tilde{
\mathbb{J}}_{B} \! + \! \tilde{\mathbb{J}}_{C} \! + \! \tilde{\mathbb{J}}_{D} 
\! + \! \tilde{\mathbb{J}}_{E},
\end{equation}
where
\begin{gather*}
\tilde{\mathbb{J}}_{A} \! := \! \int_{-\infty}^{\tilde{b}_{0}-\tilde{\delta}_{
\tilde{b}_{0}}} \sum_{i_{1},i_{2}=1,2} \left\lvert \dfrac{(w_{+}^{\Sigma_{
\tilde{\mathcal{R}}}}(\xi))_{i_{1}i_{2}}}{\xi \! - \! \alpha_{k}} \right\rvert^{2} 
\, \dfrac{\lvert \md \xi \rvert}{4 \pi^{2}}, \quad \quad \tilde{\mathbb{J}}_{
B} \! := \! \int_{\tilde{a}_{N+1}+\tilde{\delta}_{\tilde{a}_{N+1}}}^{+\infty} 
\sum_{i_{1},i_{2}=1,2} \left\lvert \dfrac{(w_{+}^{\Sigma_{\tilde{\mathcal{
R}}}}(\xi))_{i_{1}i_{2}}}{\xi \! - \! \alpha_{k}} \right\rvert^{2} \, \dfrac{\lvert 
\md \xi \rvert}{4 \pi^{2}}, \\
\tilde{\mathbb{J}}_{C} \! := \! \sum_{j=1}^{N} \int_{\tilde{a}_{j}+\tilde{
\delta}_{\tilde{a}_{j}}}^{\tilde{b}_{j}-\tilde{\delta}_{\tilde{b}_{j}}} \sum_{i_{1},
i_{2}=1,2} \left\lvert \dfrac{(w_{+}^{\Sigma_{\tilde{\mathcal{R}}}}(\xi))_{
i_{1}i_{2}}}{\xi \! - \! \alpha_{k}} \right\rvert^{2} \, \dfrac{\lvert \md \xi 
\rvert}{4 \pi^{2}}, \quad \quad \tilde{\mathbb{J}}_{D} \! := \! \sum_{j=1}^{
N+1} \left(\int_{\tilde{\Sigma}_{p,j}^{3}} \! + \! \int_{\tilde{\Sigma}_{p,j}^{4}} 
\right) \sum_{i_{1},i_{2}=1,2} \left\lvert \dfrac{(w_{+}^{\Sigma_{\tilde{
\mathcal{R}}}}(\xi))_{i_{1}i_{2}}}{\xi \! - \! \alpha_{k}} \right\rvert^{2} \, 
\dfrac{\lvert \md \xi \rvert}{4 \pi^{2}}, \\
\tilde{\mathbb{J}}_{E} \! := \! \sum_{j=1}^{N+1} \left(\int_{\partial \tilde{
\mathbb{U}}_{\tilde{\delta}_{\tilde{b}_{j-1}}}} \! + \! \int_{\partial \tilde{
\mathbb{U}}_{\tilde{\delta}_{\tilde{a}_{j}}}} \right) \sum_{i_{1},i_{2}=1,2} 
\left\lvert \dfrac{(w_{+}^{\Sigma_{\tilde{\mathcal{R}}}}(\xi))_{i_{1}i_{2}}}{\xi 
\! - \! \alpha_{k}} \right\rvert^{2} \, \dfrac{\lvert \md \xi \rvert}{4 \pi^{2}}.
\end{gather*}
For $n \! \in \! \mathbb{N}$ and $k \! \in \! \lbrace 1,2,\dotsc,K \rbrace$ 
such that $\alpha_{p_{\mathfrak{s}}} \! := \! \alpha_{k} \! \neq \! \infty$, 
recalling {}from the calculations leading to the Estimate~\eqref{eqiy7} that, 
for $j \! \in \! \lbrace 1,2,\dotsc,N \rbrace$, $\tilde{\Sigma}_{p,j}^{2} \! 
:= \! (\tilde{a}_{j} \! + \! \tilde{\delta}_{\tilde{a}_{j}},\tilde{b}_{j} \! - \! 
\tilde{\delta}_{\tilde{b}_{j}}) \! = \! \tilde{A}_{\tilde{\mathcal{R}},1}(j) \cup 
\tilde{A}_{\tilde{\mathcal{R}},1}^{c}(j)$ (with $\tilde{A}_{\tilde{\mathcal{R}},
1}(j) \cap \tilde{A}_{\tilde{\mathcal{R}},1}^{c}(j) \! = \! \varnothing$), where 
$\tilde{A}_{\tilde{\mathcal{R}},1}(j) \! = \! (\tilde{a}_{j} \! + \! \tilde{\delta}_{
\tilde{a}_{j}},\tilde{b}_{j} \! - \! \tilde{\delta}_{\tilde{b}_{j}}) \setminus \cup_{q 
\in \tilde{Q}_{\tilde{\mathcal{R}},1}(j)} \mathscr{O}_{\tilde{\delta}_{\tilde{
\mathcal{R}},1}(j)}(\alpha_{p_{q}})$, and $\tilde{A}_{\tilde{\mathcal{R}},1}^{c}
(j) \! = \! \cup_{q \in \tilde{Q}_{\tilde{\mathcal{R}},1}(j)} \mathscr{O}_{\tilde{
\delta}_{\tilde{\mathcal{R}},1}(j)}(\alpha_{p_{q}})$, with $\tilde{Q}_{\tilde{
\mathcal{R}},1}(j) \! := \! \lbrace \mathstrut q^{\prime} \! \in \! \lbrace 1,
\dotsc,\mathfrak{s} \! - \! 2,\mathfrak{s} \rbrace; \, \alpha_{p_{q^{\prime}}} 
\! \in \! (\tilde{a}_{j} \! + \! \tilde{\delta}_{\tilde{a}_{j}},\tilde{b}_{j} \! - \! 
\tilde{\delta}_{\tilde{b}_{j}}) \rbrace$, and sufficiently small $\tilde{\delta}_{
\tilde{\mathcal{R}},1}(j) \! > \! 0$ chosen so that $\mathscr{O}_{\tilde{
\delta}_{\tilde{\mathcal{R}},1}(j)}(\alpha_{p_{q_{1}}}) \cap \mathscr{O}_{\tilde{
\delta}_{\tilde{\mathcal{R}},1}(j)}(\alpha_{p_{q_{2}}}) \! = \! \varnothing$ 
$\forall$ $q_{1} \! \neq \! q_{2} \! \in \! \tilde{Q}_{\tilde{\mathcal{R}},1}(j)$ 
and $\mathscr{O}_{\tilde{\delta}_{\tilde{\mathcal{R}},1}(j)}(\alpha_{p_{q_{1}}}) 
\cap \lbrace \tilde{a}_{j} \! + \! \tilde{\delta}_{\tilde{a}_{j}} \rbrace \! = \! 
\varnothing \! = \! \mathscr{O}_{\tilde{\delta}_{\tilde{\mathcal{R}},1}(j)}
(\alpha_{p_{q_{1}}}) \cap \lbrace \tilde{b}_{j} \! - \! \tilde{\delta}_{\tilde{b}_{j}} 
\rbrace$, via the formula $w^{\Sigma_{\tilde{\mathcal{R}}}}_{+}(z) \! = \! 
\tilde{v}_{\tilde{\mathcal{R}}}(z) \! - \! \mathrm{I}$, the asymptotics, in 
the double-scaling limit $\mathscr{N},n \! \to \! \infty$ such that $z_{o} 
\! = \! 1 \! + \! o(1)$, of $\tilde{v}_{\tilde{\mathcal{R}}}(z)$ given in 
Equation~\eqref{eqtlvee8}, and the inequality $\ln \lvert x \rvert \! \leqslant 
\! \lvert x \rvert \! - \! 1$, one shows, via an integration argument, that
\begin{equation} \label{eqiy33} 
\tilde{\mathbb{J}}_{C} \! = \! \sum_{j=1}^{N} \left(\int_{\tilde{A}_{\tilde{
\mathcal{R}},1}(j)} \! + \! \int_{\tilde{A}_{\tilde{\mathcal{R}},1}^{c}(j)} 
\right) \sum_{i_{1},i_{2}=1,2} \left\lvert \dfrac{(w_{+}^{\Sigma_{\tilde{
\mathcal{R}}}}(\xi))_{i_{1}i_{2}}}{\xi \! - \! \alpha_{k}} \right\rvert^{2} 
\, \dfrac{\md \xi}{4 \pi^{2}} \! \underset{\underset{z_{o}=1+o(1)}{
\mathscr{N},n \to \infty}}{\leqslant} \! \mathcal{O} \left(\dfrac{\tilde{
\mathfrak{c}}_{\tilde{\mathbb{J}}_{C}}(n,k,z_{o}) \me^{-2((n-1)K+k) 
\min\limits_{j=1,2,\dotsc,N} \lbrace \tilde{\lambda}_{\tilde{\mathcal{R}},
1}^{\sharp}(j) \rbrace}}{(n \! - \! 1)K \! + \! k} \right),
\end{equation}
where $\tilde{\lambda}_{\tilde{\mathcal{R}},1}^{\sharp}(j)$ $(> \! 0)$, 
$j \! \in \! \lbrace 1,2,\dotsc,N \rbrace$, is defined by 
Equation~\eqref{eqqlamrj}, and $\tilde{\mathfrak{c}}_{\tilde{\mathbb{J}}_{C}}
(n,k,z_{o}) \! =_{\underset{z_{o}=1+o(1)}{\mathscr{N},n \to \infty}} \! 
\mathcal{O}(1)$. For $n \! \in \! \mathbb{N}$ and $k \! \in \! \lbrace 1,2,
\dotsc,K \rbrace$ such that $\alpha_{p_{\mathfrak{s}}} \! := \! \alpha_{k} 
\! \neq \! \infty$, recalling {}from the calculations leading to the 
Estimate~\eqref{eqiy9} that $(\tilde{a}_{N+1} \! + \! \tilde{\delta}_{
\tilde{a}_{N+1}},+\infty) \! = \! \tilde{A}_{\tilde{\mathcal{R}},2}(+) \cup 
\tilde{A}_{\tilde{\mathcal{R}},2}^{c}(+)$ (with $\tilde{A}_{\tilde{\mathcal{R}},
2}(+) \cap \tilde{A}_{\tilde{\mathcal{R}},2}^{c}(+) \! = \! \varnothing$), 
where $\tilde{A}_{\tilde{\mathcal{R}},2}(+) \! = \! (\tilde{a}_{N+1} \! + \! 
\tilde{\delta}_{\tilde{a}_{N+1}},+\infty) \setminus (\mathscr{O}_{\infty}
(\alpha_{p_{\mathfrak{s}-1}}) \cup \cup_{q \in \tilde{Q}_{\tilde{\mathcal{R}},2}(+)} 
\mathscr{O}_{\tilde{\delta}_{\tilde{\mathcal{R}},2}(+)}(\alpha_{p_{q}}))$, 
and $\tilde{A}_{\tilde{\mathcal{R}},2}^{c}(+) \! = \! \mathscr{O}_{\infty}
(\alpha_{p_{\mathfrak{s}-1}}) \cup \cup_{q \in \tilde{Q}_{\tilde{
\mathcal{R}},2}(+)} \mathscr{O}_{\tilde{\delta}_{\tilde{\mathcal{R}},2}(+)}
(\alpha_{p_{q}})$, with $\tilde{Q}_{\tilde{\mathcal{R}},2}(+) \! := \! \lbrace 
\mathstrut q^{\prime} \! \in \! \lbrace 1,\dotsc,\mathfrak{s} \! - \! 
2,\mathfrak{s} \rbrace; \, \alpha_{p_{q^{\prime}}} \! \in \! (\tilde{a}_{
N+1} \! + \! \tilde{\delta}_{\tilde{a}_{N+1}},+\infty) \rbrace$ and 
sufficiently small $\tilde{\varepsilon}_{\infty},\tilde{\delta}_{\tilde{
\mathcal{R}},2}(+) \! > \! 0$ chosen so that $\mathscr{O}_{\tilde{
\delta}_{\tilde{\mathcal{R}},2}(+)}(\alpha_{p_{q^{\prime}_{1}}}) \cap 
\mathscr{O}_{\tilde{\delta}_{\tilde{\mathcal{R}},2}(+)}(\alpha_{p_{q^{\prime 
\prime}_{1}}}) \! = \! \varnothing$ $\forall$ $q^{\prime}_{1} \! \neq \! 
q^{\prime \prime}_{1} \! \in \! \tilde{Q}_{\tilde{\mathcal{R}},2}(+)$, 
$\mathscr{O}_{\tilde{\delta}_{\tilde{\mathcal{R}},2}(+)}(\alpha_{p_{
q^{\prime}_{1}}}) \cap \lbrace \tilde{a}_{N+1} \! + \! \tilde{\delta}_{
\tilde{a}_{N+1}} \rbrace \! = \! \varnothing$, and $\mathscr{O}_{
\tilde{\delta}_{\tilde{\mathcal{R}},2}(+)}(\alpha_{p_{q^{\prime}_{1}}}) 
\cap \mathscr{O}_{\infty}(\alpha_{p_{\mathfrak{s}-1}}) \! = \! 
\varnothing$, via the formula $w^{\Sigma_{\tilde{\mathcal{R}}}}_{+}
(z) \! = \! \tilde{v}_{\tilde{\mathcal{R}}}(z) \! - \! \mathrm{I}$, the 
asymptotics, in the double-scaling limit $\mathscr{N},n \! \to \! \infty$ 
such that $z_{o} \! = \! 1 \! + \! o(1)$, of $\tilde{v}_{\tilde{\mathcal{R}}}
(z)$ given in Equation~\eqref{eqtlvee9}, and the inequality $\ln \lvert 
x \rvert \! \leqslant \! \lvert x \rvert \! - \! 1$, one shows, via an 
integration argument, that
\begin{equation} \label{eqiy34} 
\tilde{\mathbb{J}}_{B} \! = \! \left(\int_{\tilde{A}_{\tilde{\mathcal{R}},2}(+)} 
\! + \! \int_{\tilde{A}_{\tilde{\mathcal{R}},2}^{c}(+)} \right) \sum_{i_{1},
i_{2}=1,2} \left\lvert \dfrac{(w_{+}^{\Sigma_{\tilde{\mathcal{R}}}}(\xi))_{
i_{1}i_{2}}}{\xi \! - \! \alpha_{k}} \right\rvert^{2} \, \dfrac{\md \xi}{4 
\pi^{2}} \! \underset{\underset{z_{o}=1+o(1)}{\mathscr{N},n \to \infty}}{
\leqslant} \! \mathcal{O} \left(\dfrac{\tilde{\mathfrak{c}}_{\tilde{\mathbb{
J}}_{B}}(n,k,z_{o}) \me^{-2((n-1)K+k) \tilde{\lambda}_{\tilde{\mathcal{R}},
2}^{\sharp}(+)}}{(n \! - \! 1)K \! + \! k} \right),
\end{equation}
where $\tilde{\lambda}_{\tilde{\mathcal{R}},2}^{\sharp}(+)$ $(> \! 0)$ is 
defined by Equation~\eqref{eqqlamrp}, and $\tilde{\mathfrak{c}}_{\tilde{
\mathbb{J}}_{B}}(n,k,z_{o}) \! =_{\underset{z_{o}=1+o(1)}{\mathscr{N},n 
\to \infty}} \! \mathcal{O}(1)$; proceeding analogously as above (for 
$\tilde{\mathbb{J}}_{B})$, one shows that, for $n \! \in \! \mathbb{N}$ 
and $k \! \in \! \lbrace 1,2,\dotsc,K \rbrace$ such that $\alpha_{p_{
\mathfrak{s}}} \! := \! \alpha_{k} \! \neq \! \infty$,
\begin{equation} \label{eqiy35} 
\tilde{\mathbb{J}}_{A} \! \underset{\underset{z_{o}=1+o(1)}{\mathscr{N},
n \to \infty}}{\leqslant} \! \mathcal{O} \left(\dfrac{\tilde{\mathfrak{c}}_{
\tilde{\mathbb{J}}_{A}}(n,k,z_{o}) \me^{-2((n-1)K+k) \tilde{\lambda}_{
\tilde{\mathcal{R}},2}^{\sharp}(-)}}{(n \! - \! 1)K \! + \! k} \right),
\end{equation}
where $\tilde{\lambda}_{\tilde{\mathcal{R}},2}^{\sharp}(-)$ $(> \! 0)$ is 
defined by Equation~\eqref{eqqlamrm}, and $\tilde{\mathfrak{c}}_{\tilde{
\mathbb{J}}_{A}}(n,k,z_{o}) \! =_{\underset{z_{o}=1+o(1)}{\mathscr{N},n 
\to \infty}} \! \mathcal{O}(1)$. For $n \! \in \! \mathbb{N}$ and $k \! \in 
\! \lbrace 1,2,\dotsc,K \rbrace$ such that $\alpha_{p_{\mathfrak{s}}} \! 
:= \! \alpha_{k} \! \neq \! \infty$, via the formula $w^{\Sigma_{\tilde{
\mathcal{R}}}}_{+}(z) \! = \! \tilde{v}_{\tilde{\mathcal{R}}}(z) \! - \! 
\mathrm{I}$, the asymptotics, in the double-scaling limit $\mathscr{N},
n \! \to \! \infty$ such that $z_{o} \! = \! 1 \! + \! o(1)$, of $\tilde{v}_{
\tilde{\mathcal{R}}}(z)$ given in Equation~\eqref{eqtlvee11}, the 
elementary trigonometric inequalities $\sin \theta \! \geqslant \! 
\tfrac{2 \theta}{\pi}$ and $\cos \theta \! \geqslant \! -\tfrac{2 \theta}{\pi} 
\! + \! 1$ for $0 \! \leqslant \! \theta \! \leqslant \! \tfrac{\pi}{2}$, and 
$\sin \theta \! \geqslant \! \tfrac{2}{\pi}(\pi \! - \! \theta)$ and $\cos 
\theta \! \leqslant \! -\tfrac{2 \theta}{\pi} \! + \! 1$ for $\tfrac{\pi}{2} 
\! \leqslant \! \theta \! \leqslant \! \pi$, the parametrisation $\tilde{
\Sigma}_{p,j}^{3} \! = \! \lbrace \mathstrut (x_{j}(\theta),y_{j}(\theta)); 
\, x_{j}(\theta) \! = \! \tfrac{1}{2}(\tilde{a}_{j} \! + \! \tilde{b}_{j-1}) \! + \! 
\tfrac{1}{2}(\tilde{a}_{j} \! - \! \tilde{b}_{j-1}) \cos \theta, \, y_{j}(\theta) 
\! = \! \tilde{\eta}_{j} \sin \theta, \, \theta_{0}^{\tilde{a}}(j) \! \leqslant 
\! \theta \! \leqslant \! \pi \! - \! \theta_{0}^{\tilde{b}}(j) \rbrace$, $j \! 
\in \! \lbrace 1,2,\dotsc,N \! + \! 1 \rbrace$, for the---elliptic---homotopic 
deformation of $\tilde{\Sigma}_{p,j}^{3}$, and the inequality $\ln \lvert 
x \rvert \! \leqslant \! \lvert x \rvert \! - \! 1$, one shows, via a tedious 
integration-by-parts argument, that
\begin{align*}
& \, \sum_{j=1}^{N+1} \int_{\tilde{\Sigma}_{p,j}^{3}} \sum_{i_{1},i_{2}=1,2} 
\left\lvert \dfrac{(w_{+}^{\Sigma_{\tilde{\mathcal{R}}}}(\xi))_{i_{1}i_{2}}}{
\xi \! - \! \alpha_{k}} \right\rvert^{2} \, \dfrac{\lvert \md \xi \rvert}{4 
\pi^{2}} \! = \! \sum_{j=1}^{N+1} \left(\int_{\theta_{0}^{\tilde{a}}(j)}^{
\pi/2} \! + \! \int_{\pi/2}^{\pi -\theta_{0}^{\tilde{b}}(j)} \right) 
\sum_{i_{1},i_{2}=1,2} \left\lvert \dfrac{(w_{+}^{\Sigma_{\tilde{
\mathcal{R}}}}(x_{j}(\theta) \! + \! \mi y_{j}(\theta)))_{i_{1}i_{2}}}{x_{j}
(\theta) \! + \! \mi y_{j}(\theta) \! - \! \alpha_{k}} \right\rvert^{2} \\
&\times \, ((x_{j}^{\prime}(\theta))^{2} \! + \! (y_{j}^{\prime}
(\theta))^{2})^{1/2} \, \dfrac{\md \theta}{4 \pi^{2}} \! \underset{
\underset{z_{o}=1+o(1)}{\mathscr{N},n \to \infty}}{\leqslant} \! 
\mathcal{O} \left(\dfrac{\tilde{\mathfrak{c}}_{\tilde{\mathbb{J}}_{D}}^{
\smallfrown}(n,k,z_{o}) \me^{-2((n-1)K+k) \min\limits_{j=1,2,\dotsc,
N+1} \lbrace \tilde{\lambda}_{\tilde{\mathcal{R}},4}^{\sharp,\smallfrown}
(j) \rbrace}}{(n \! - \! 1)K \! + \! k} \right),
\end{align*}
where $\tilde{\lambda}_{\tilde{\mathcal{R}},4}^{\sharp,\smallfrown}(j)$ 
$(> \! 0)$ is defined by Equation~\eqref{eqqlamrupj}, and $\tilde{
\mathfrak{c}}_{\tilde{\mathbb{J}}_{D}}^{\smallfrown}(n,k,z_{o}) \! =_{
\underset{z_{o}=1+o(1)}{\mathscr{N},n \to \infty}} \! \mathcal{O}(1)$; 
analogously, one shows that, for $n \! \in \! \mathbb{N}$ and $k \! \in 
\! \lbrace 1,2,\dotsc,K \rbrace$ such that $\alpha_{p_{\mathfrak{s}}} 
\! := \! \alpha_{k} \! \neq \! \infty$,
\begin{equation*}
\sum_{j=1}^{N+1} \int_{\tilde{\Sigma}_{p,j}^{4}} \sum_{i_{1},i_{2}=1,2} 
\left\lvert \dfrac{(w_{+}^{\Sigma_{\tilde{\mathcal{R}}}}(\xi))_{i_{1}i_{2}}}{
\xi \! - \! \alpha_{k}} \right\rvert^{2} \, \dfrac{\lvert \md \xi \rvert}{4 
\pi^{2}} \! \underset{\underset{z_{o}=1+o(1)}{\mathscr{N},n \to 
\infty}}{\leqslant} \! \mathcal{O} \left(\dfrac{\tilde{\mathfrak{c}}_{
\tilde{\mathbb{J}}_{D}}^{\smallsmile}(n,k,z_{o}) \me^{-2((n-1)K+k) 
\min\limits_{j=1,2,\dotsc,N+1} \lbrace \tilde{\lambda}_{\tilde{
\mathcal{R}},4}^{\sharp,\smallsmile}(j) \rbrace}}{(n \! - \! 1)K \! + \! k} 
\right),
\end{equation*}
where $\tilde{\lambda}_{\tilde{\mathcal{R}},4}^{\sharp,\smallsmile}
(j)$ $(> \! 0)$ is defined by Equation~\eqref{eqqlamrdwj}, and $\tilde{
\mathfrak{c}}_{\tilde{\mathbb{J}}_{D}}^{\smallsmile}(n,k,z_{o}) \! =_{
\underset{z_{o}=1+o(1)}{\mathscr{N},n \to \infty}} \! \mathcal{O}(1)$, 
whence
\begin{equation} \label{eqiy36} 
\tilde{\mathbb{J}}_{D} \underset{\underset{z_{o}=1+o(1)}{\mathscr{N},
n \to \infty}}{\leqslant} \mathcal{O} \left(\dfrac{\tilde{\mathfrak{c}}_{
\tilde{\mathbb{J}}_{D}}(n,k,z_{o}) \me^{-2((n-1)K+k) \min \left\lbrace 
\min\limits_{j=1,2,\dotsc,N+1} \lbrace \tilde{\lambda}_{\tilde{\mathcal{
R}},4}^{\sharp,\smallfrown}(j) \rbrace,\min\limits_{j=1,2,\dotsc,N+1} 
\lbrace \tilde{\lambda}_{\tilde{\mathcal{R}},4}^{\sharp,\smallsmile}(j) 
\rbrace \right\rbrace}}{(n \! - \! 1)K \! + \! k} \right),
\end{equation}
where $\tilde{\mathfrak{c}}_{\tilde{\mathbb{J}}_{D}}(n,k,z_{o}) \! =_{
\underset{z_{o}=1+o(1)}{\mathscr{N},n \to \infty}} \! \mathcal{O}(1)$. For 
$n \! \in \! \mathbb{N}$ and $k \! \in \! \lbrace 1,2,\dotsc,K \rbrace$ 
such that $\alpha_{p_{\mathfrak{s}}} \! := \! \alpha_{k} \! \neq \! \infty$, 
proceeding as in the calculations leading to the Estimate~\eqref{eqiy20}, 
one shows, via the asymptotics, in the double-scaling limit $\mathscr{N},
n \! \to \! \infty$ such that $z_{o} \! = \! 1 \! + \! o(1)$, of $w_{+}^{
\Sigma_{\tilde{\mathcal{R}}}}(z)$ for $z \! \in \! \partial \tilde{\mathbb{
U}}_{\tilde{\delta}_{\tilde{b}_{j-1}}}$ and for $z \! \in \! \partial \tilde{
\mathbb{U}}_{\tilde{\delta}_{\tilde{a}_{j}}}$, $j \! \in \! \lbrace 1,2,
\dotsc,N \! + \! 1 \rbrace$, given in Equations~\eqref{eqproptilb} 
and~\eqref{eqproptila}, respectively, and the indefinite integral 
\pmb{2.562} 1. on p.~185 of \cite{gradryzh}, that, for $j \! \in \! 
\lbrace 1,2,\dotsc,N \! + \! 1 \rbrace$,
\begin{align*}
& \, \int_{\partial \tilde{\mathbb{U}}_{\tilde{\delta}_{\tilde{b}_{j-1}}}} 
\sum_{i_{1},i_{2}=1,2} \left\lvert \dfrac{(w_{+}^{\Sigma_{\tilde{
\mathcal{R}}}}(\xi))_{i_{1}i_{2}}}{\xi \! - \! \alpha_{k}} \right\rvert^{2} \, 
\dfrac{\lvert \md \xi \rvert}{4 \pi^{2}} \! \underset{\underset{z_{o}=
1+o(1)}{\mathscr{N},n \to \infty}}{\leqslant} \! \dfrac{\tilde{\delta}_{
\tilde{b}_{j-1}}}{\pi^{2}} \underbrace{\max_{i_{1},i_{2}=1,2} \sup_{
\theta_{0} \in [0,2 \pi]} \left\lvert (w_{+}^{\Sigma_{\tilde{\mathcal{
R}}}}(\tilde{b}_{j-1} \! + \! \tilde{\delta}_{\tilde{b}_{j-1}} \me^{\mi 
\theta_{0}}))_{i_{1}i_{2}} \right\rvert^{2}}_{= \, \mathcal{O}(((n-1)K
+k)^{-2})} \\
&\times \, \int_{0}^{2 \pi} \dfrac{1}{(\tilde{b}_{j-1} \! - \! \alpha_{k})^{
2} \! + \! \tilde{\delta}_{\tilde{b}_{j-1}}^{2} \! + \! 2 \tilde{\delta}_{
\tilde{b}_{j-1}}(\tilde{b}_{j-1} \! - \! \alpha_{k}) \cos \theta} \, \md 
\theta \! \underset{\underset{z_{o}=1+o(1)}{\mathscr{N},n \to \infty}}{
\leqslant} \! \dfrac{\tilde{\delta}_{\tilde{b}_{j-1}}}{\pi^{2}} \mathcal{O} 
\left(\dfrac{1}{((n \! - \! 1)K \! + \! k)^{2}} \right) \dfrac{2}{((\tilde{b}_{
j-1} \! - \! \alpha_{k})^{2} \! - \! \tilde{\delta}_{\tilde{b}_{j-1}}^{2})} \\
&\times \, \underbrace{\left. \tan^{-1} \left(\left(\dfrac{\tilde{b}_{j-1} 
\! - \! \alpha_{k} \! - \! \tilde{\delta}_{\tilde{b}_{j-1}}}{\tilde{b}_{j-1} \! 
- \! \alpha_{k} \! + \! \tilde{\delta}_{\tilde{b}_{j-1}}} \right) \tan \left(
\dfrac{\theta}{2} \right) \right) \right\vert_{0}^{2 \pi}}_{= \, 0} \quad 
\Rightarrow
\end{align*}
\begin{equation*}
\int_{\partial \tilde{\mathbb{U}}_{\tilde{\delta}_{\tilde{b}_{j-1}}}} 
\sum_{i_{1},i_{2}=1,2} \left\lvert \dfrac{(w_{+}^{\Sigma_{\tilde{
\mathcal{R}}}}(\xi))_{i_{1}i_{2}}}{\xi \! - \! \alpha_{k}} \right\rvert^{2} 
\, \dfrac{\lvert \md \xi \rvert}{4 \pi^{2}} \underset{\underset{z_{o}
=1+o(1)}{\mathscr{N},n \to \infty}}{=} 0, 
\end{equation*}
and
\begin{align*}
& \, \int_{\partial \tilde{\mathbb{U}}_{\tilde{\delta}_{\tilde{a}_{j}}}} 
\sum_{i_{1},i_{2}=1,2} \left\lvert \dfrac{(w_{+}^{\Sigma_{\tilde{
\mathcal{R}}}}(\xi))_{i_{1}i_{2}}}{\xi \! - \! \alpha_{k}} \right\rvert^{2} \, 
\dfrac{\lvert \md \xi \rvert}{4 \pi^{2}} \! \underset{\underset{z_{o}=
1+o(1)}{\mathscr{N},n \to \infty}}{\leqslant} \! \dfrac{\tilde{\delta}_{
\tilde{a}_{j}}}{\pi^{2}} \underbrace{\max_{i_{1},i_{2}=1,2} \sup_{
\theta_{0} \in [0,2 \pi]} \left\lvert (w_{+}^{\Sigma_{\tilde{\mathcal{
R}}}}(\tilde{a}_{j} \! + \! \tilde{\delta}_{\tilde{a}_{j}} \me^{\mi 
\theta_{0}}))_{i_{1}i_{2}} \right\rvert^{2}}_{= \, \mathcal{O}(((n-1)K
+k)^{-2})} \\
&\times \, \int_{0}^{2 \pi} \dfrac{1}{(\tilde{a}_{j} \! - \! \alpha_{k})^{
2} \! + \! \tilde{\delta}_{\tilde{a}_{j}}^{2} \! + \! 2 \tilde{\delta}_{
\tilde{a}_{j}}(\tilde{a}_{j} \! - \! \alpha_{k}) \cos \theta} \, \md \theta \! 
\underset{\underset{z_{o}=1+o(1)}{\mathscr{N},n \to \infty}}{\leqslant} 
\! \dfrac{\tilde{\delta}_{\tilde{a}_{j}}}{\pi^{2}} \mathcal{O} \left(
\dfrac{1}{((n \! - \! 1)K \! + \! k)^{2}} \right) \dfrac{2}{((\tilde{a}_{j} \! 
- \! \alpha_{k})^{2} \! - \! \tilde{\delta}_{\tilde{a}_{j}}^{2})} \\
&\times \, \underbrace{\left. \tan^{-1} \left(\left(\dfrac{\tilde{a}_{j} 
\! - \! \alpha_{k} \! - \! \tilde{\delta}_{\tilde{a}_{j}}}{\tilde{a}_{j} \! - 
\! \alpha_{k} \! + \! \tilde{\delta}_{\tilde{a}_{j}}} \right) \tan \left(
\dfrac{\theta}{2} \right) \right) \right\vert_{0}^{2 \pi}}_{= \, 0} \quad 
\Rightarrow
\end{align*}
\begin{equation*}
\int_{\partial \tilde{\mathbb{U}}_{\tilde{\delta}_{\tilde{a}_{j}}}} 
\sum_{i_{1},i_{2}=1,2} \left\lvert \dfrac{(w_{+}^{\Sigma_{\tilde{
\mathcal{R}}}}(\xi))_{i_{1}i_{2}}}{\xi \! - \! \alpha_{k}} \right\rvert^{2} 
\, \dfrac{\lvert \md \xi \rvert}{4 \pi^{2}} \underset{\underset{z_{o}
=1+o(1)}{\mathscr{N},n \to \infty}}{=} 0, 
\end{equation*}
whence
\begin{equation} \label{eqiy37} 
\tilde{\mathbb{J}}_{E} \underset{\underset{z_{o}=1+o(1)}{
\mathscr{N},n \to \infty}}{=} 0.
\end{equation}
Hence, via the Estimates~\eqref{eqiy33}--\eqref{eqiy37} and 
Equation~\eqref{eqiy32}, one arrives at, for $n \! \in \! \mathbb{N}$ 
and $k \! \in \! \lbrace 1,2,\dotsc,K \rbrace$ such that $\alpha_{
p_{\mathfrak{s}}} \! := \! \alpha_{k} \! \neq \! \infty$,
\begin{equation} \label{eqiy38} 
\left\lvert \left\lvert \dfrac{w_{+}^{\Sigma_{\tilde{\mathcal{R}}}}
(\pmb{\cdot})}{2 \pi \mi (\pmb{\cdot} \! - \! \alpha_{k})} \right\rvert 
\right\rvert_{\mathcal{L}^{2}_{\mathrm{M}_{2}(\mathbb{C})}
(\tilde{\Sigma}_{\tilde{\mathcal{R}}}^{\sharp})}^{2} \underset{
\underset{z_{o}=1+o(1)}{\mathscr{N},n \to \infty}}{\leqslant} 
\mathcal{O} \left(\dfrac{\tilde{\mathfrak{c}}_{\tilde{\mathbb{J}}}
(n,k,z_{o}) \me^{-((n-1)K+k) \tilde{\lambda}_{\tilde{\mathcal{R}},
\tilde{w}}^{\triangleright}}}{\left((n \! - \! 1)K \! + \! k \right)^{1/2}} 
\right),
\end{equation}
where $\tilde{\lambda}_{\tilde{\mathcal{R}},\tilde{w}}^{\triangleright}$ 
$(> \! 0)$ is defined by Equation~\eqref{tillamrwr2}, and $\tilde{
\mathfrak{c}}_{\tilde{\mathbb{J}}}(n,k,z_{o}) \! =_{\underset{z_{o}=
1+o(1)}{\mathscr{N},n \to \infty}} \! \mathcal{O}(1)$.

For $n \! \in \! \mathbb{N}$ and $k \! \in \! \lbrace 1,2,\dotsc,K 
\rbrace$ such that $\alpha_{p_{\mathfrak{s}}} \! := \! \alpha_{k} \! 
\neq \! \infty$, via the above definitions, the disjointness of the 
oriented skeleton $\tilde{\Sigma}_{\tilde{\mathcal{R}}}^{\sharp}$, 
H\"{o}lder's Inequality, and H\"{o}lder's Inequality for Integrals, a 
calculation shows that
\begin{equation} \label{eqiy39}
\int_{\tilde{\Sigma}_{\tilde{\mathcal{R}}}^{\sharp}} \left(\lim_{-\tilde{
\Sigma}_{\tilde{\mathcal{R}}}^{\sharp} \ni z^{\prime} \to \xi} \left\lvert 
\left\lvert \dfrac{w_{+}^{\Sigma_{\tilde{\mathcal{R}}}}(\pmb{\cdot})}{2 
\pi \mi (\pmb{\cdot} \! - \! z^{\prime})} \right\rvert \right\rvert_{
\mathcal{L}^{2}_{\mathrm{M}_{2}(\mathbb{C})}(\tilde{\Sigma}_{\tilde{
\mathcal{R}}}^{\sharp})}^{2} \right) \, \lvert \md \xi \rvert \! \leqslant 
\! \tilde{\mathbb{J}}_{\mathbb{A}}^{\sharp} \! + \! \tilde{\mathbb{J}}_{
\mathbb{B}}^{\sharp} \! + \! \tilde{\mathbb{J}}_{\mathbb{C}}^{\sharp} 
\! + \! \tilde{\mathbb{J}}_{\mathbb{D}}^{\sharp} \! + \! \tilde{
\mathbb{J}}_{\mathbb{E}}^{\sharp},
\end{equation}
where
\begin{gather}
\tilde{\mathbb{J}}_{\mathbb{A}}^{\sharp} \! := \! \int_{\tilde{\Sigma}_{
\tilde{\mathcal{R}}}^{\sharp}} \left(\lim_{-\tilde{\Sigma}_{\tilde{
\mathcal{R}}}^{\sharp} \ni z^{\prime} \to \xi} \sum_{j=1}^{N} \int_{
\tilde{a}_{j}+\tilde{\delta}_{\tilde{a}_{j}}}^{\tilde{b}_{j}-\tilde{\delta}_{
\tilde{b}_{j}}} \sum_{i_{1},i_{2}=1,2} \left\lvert \dfrac{(w_{+}^{\Sigma_{
\tilde{\mathcal{R}}}}(u))_{i_{1}i_{2}}}{u \! - \! z^{\prime}} \right\rvert^{2} 
\, \dfrac{\lvert \md u \rvert}{2 \pi^{2}} \right) \lvert \md \xi \rvert, 
\label{eqiy40} \\
\tilde{\mathbb{J}}_{\mathbb{B}}^{\sharp} \! := \! \int_{\tilde{\Sigma}_{
\tilde{\mathcal{R}}}^{\sharp}} \left(\lim_{-\tilde{\Sigma}_{\tilde{
\mathcal{R}}}^{\sharp} \ni z^{\prime} \to \xi} \int_{\tilde{a}_{N+1}+
\tilde{\delta}_{\tilde{a}_{N+1}}}^{+\infty} \sum_{i_{1},i_{2}=1,2} \left\lvert 
\dfrac{(w_{+}^{\Sigma_{\tilde{\mathcal{R}}}}(u))_{i_{1}i_{2}}}{u \! - \! 
z^{\prime}} \right\rvert^{2} \, \dfrac{\lvert \md u \rvert}{2 \pi^{2}} 
\right) \lvert \md \xi \rvert, \label{eqiy41} \\
\tilde{\mathbb{J}}_{\mathbb{C}}^{\sharp} \! := \! \int_{\tilde{\Sigma}_{
\tilde{\mathcal{R}}}^{\sharp}} \left(\lim_{-\tilde{\Sigma}_{\tilde{
\mathcal{R}}}^{\sharp} \ni z^{\prime} \to \xi} \int_{-\infty}^{\tilde{b}_{0}
-\tilde{\delta}_{\tilde{b}_{0}}} \sum_{i_{1},i_{2}=1,2} \left\lvert 
\dfrac{(w_{+}^{\Sigma_{\tilde{\mathcal{R}}}}(u))_{i_{1}i_{2}}}{u \! - \! 
z^{\prime}} \right\rvert^{2} \, \dfrac{\lvert \md u \rvert}{2 \pi^{2}} 
\right) \lvert \md \xi \rvert, \label{eqiy42} \\
\tilde{\mathbb{J}}_{\mathbb{D}}^{\sharp} \! := \! \int_{\tilde{\Sigma}_{
\tilde{\mathcal{R}}}^{\sharp}} \left(\lim_{-\tilde{\Sigma}_{\tilde{
\mathcal{R}}}^{\sharp} \ni z^{\prime} \to \xi} \sum_{j=1}^{N+1} 
\sum_{m=3,4} \int_{\tilde{\Sigma}_{p,j}^{m}} \sum_{i_{1},i_{2}=1,2} 
\left\lvert \dfrac{(w_{+}^{\Sigma_{\tilde{\mathcal{R}}}}(u))_{i_{1}i_{2}}}{u 
\! - \! z^{\prime}} \right\rvert^{2} \, \dfrac{\lvert \md u \rvert}{2 
\pi^{2}} \right) \lvert \md \xi \rvert, \label{eqiy43} \\
\tilde{\mathbb{J}}_{\mathbb{E}}^{\sharp} \! := \! \int_{\tilde{\Sigma}_{
\tilde{\mathcal{R}}}^{\sharp}} \left(\lim_{-\tilde{\Sigma}_{\tilde{
\mathcal{R}}}^{\sharp} \ni z^{\prime} \to \xi} \sum_{j=1}^{N+1} \left(
\int_{\partial \tilde{\mathbb{U}}_{\tilde{\delta}_{\tilde{b}_{j-1}}}} \! + \! 
\int_{\partial \tilde{\mathbb{U}}_{\tilde{\delta}_{\tilde{a}_{j}}}} \right) 
\sum_{i_{1},i_{2}=1,2} \left\lvert \dfrac{(w_{+}^{\Sigma_{\tilde{
\mathcal{R}}}}(u))_{i_{1}i_{2}}}{u \! - \! z^{\prime}} \right\rvert^{2} \, 
\dfrac{\lvert \md u \rvert}{2 \pi^{2}} \right) \lvert \md \xi \rvert. 
\label{eqiy44}
\end{gather}
For $n \! \in \! \mathbb{N}$ and $k \! \in \! \lbrace 1,2,\dotsc,K 
\rbrace$ such that $\alpha_{p_{\mathfrak{s}}} \! := \! \alpha_{k} \! 
\neq \! \infty$, consider, say, and without loss of generality, the 
detailed analysis of the asymptotic estimate, in the double-scaling 
limit $\mathscr{N},n \! \to \! \infty$ such that $z_{o} \! = \! 1 \! + 
\! o(1)$, for the integral $\tilde{\mathbb{J}}_{\mathbb{D}}^{\sharp}$: 
the remaining integrals, that is, $\tilde{\mathbb{J}}_{\mathbb{A}}^{
\sharp}$, $\tilde{\mathbb{J}}_{\mathbb{B}}^{\sharp}$, 
$\tilde{\mathbb{J}}_{\mathbb{C}}^{\sharp}$, and 
$\tilde{\mathbb{J}}_{\mathbb{E}}^{\sharp}$, are estimated analogously. 
For $n \! \in \! \mathbb{N}$ and $k \! \in \! \lbrace 1,2,\dotsc,K \rbrace$ 
such that $\alpha_{p_{\mathfrak{s}}} \! := \! \alpha_{k} \! \neq \! \infty$, 
it follows {}from Definition~\eqref{eqiy43} and the disjointness of the 
oriented skeleton $\tilde{\Sigma}_{\tilde{\mathcal{R}}}^{\sharp}$ that
\begin{equation} \label{eqiy45} 
\tilde{\mathbb{J}}_{\mathbb{D}}^{\sharp} \! = \! \tilde{\mathbb{J}}_{
\mathbb{D},1}^{\sharp} \! + \! \tilde{\mathbb{J}}_{\mathbb{D},2}^{
\sharp} \! + \! \tilde{\mathbb{J}}_{\mathbb{D},3}^{\sharp} \! + \! 
\tilde{\mathbb{J}}_{\mathbb{D},4}^{\sharp} \! + \! \tilde{\mathbb{J}}_{
\mathbb{D},5}^{\sharp},
\end{equation}
where
\begin{gather}
\tilde{\mathbb{J}}_{\mathbb{D},1}^{\sharp} \! := \! \int_{-\infty}^{
\tilde{b}_{0}-\tilde{\delta}_{\tilde{b}_{0}}} \left(\lim_{-\tilde{\Sigma}_{
\tilde{\mathcal{R}}}^{\sharp} \ni z^{\prime} \to \xi} \sum_{j=1}^{N+1} 
\sum_{m=3,4} \int_{\tilde{\Sigma}_{p,j}^{m}} \sum_{i_{1},i_{2}=1,2} 
\left\lvert \dfrac{(w_{+}^{\Sigma_{\tilde{\mathcal{R}}}}(u))_{i_{1}i_{2}}}{u 
\! - \! z^{\prime}} \right\rvert^{2} \, \dfrac{\lvert \md u \rvert}{2 
\pi^{2}} \right) \lvert \md \xi \rvert, \label{eqiy46} \\
\tilde{\mathbb{J}}_{\mathbb{D},2}^{\sharp} \! := \! \int_{\tilde{a}_{N+1}
+\tilde{\delta}_{\tilde{a}_{N+1}}}^{+\infty} \left(\lim_{-\tilde{\Sigma}_{
\tilde{\mathcal{R}}}^{\sharp} \ni z^{\prime} \to \xi} \sum_{j=1}^{N+1} 
\sum_{m=3,4} \int_{\tilde{\Sigma}_{p,j}^{m}} \sum_{i_{1},i_{2}=1,2} 
\left\lvert \dfrac{(w_{+}^{\Sigma_{\tilde{\mathcal{R}}}}(u))_{i_{1}i_{2}}}{u 
\! - \! z^{\prime}} \right\rvert^{2} \, \dfrac{\lvert \md u \rvert}{2 
\pi^{2}} \right) \lvert \md \xi \rvert, \label{eqiy47} \\
\tilde{\mathbb{J}}_{\mathbb{D},3}^{\sharp} \! := \! \sum_{i=1}^{N} 
\int_{\tilde{a}_{i}+\tilde{\delta}_{\tilde{a}_{i}}}^{\tilde{b}_{i}-\tilde{
\delta}_{\tilde{b}_{i}}} \left(\lim_{-\tilde{\Sigma}_{\tilde{\mathcal{R}}}^{
\sharp} \ni z^{\prime} \to \xi} \sum_{j=1}^{N+1} \sum_{m=3,4} 
\int_{\tilde{\Sigma}_{p,j}^{m}} \sum_{i_{1},i_{2}=1,2} \left\lvert 
\dfrac{(w_{+}^{\Sigma_{\tilde{\mathcal{R}}}}(u))_{i_{1}i_{2}}}{u \! - \! 
z^{\prime}} \right\rvert^{2} \, \dfrac{\lvert \md u \rvert}{2 \pi^{2}} 
\right) \lvert \md \xi \rvert, \label{eqiy48} \\
\tilde{\mathbb{J}}_{\mathbb{D},4}^{\sharp} \! := \! \sum_{i=1}^{N+1} 
\sum_{m_{1}=3,4} \int_{\tilde{\Sigma}_{p,i}^{m_{1}}} \left(\lim_{-\tilde{
\Sigma}_{\tilde{\mathcal{R}}}^{\sharp} \ni z^{\prime} \to \xi} \sum_{j
=1}^{N+1} \sum_{m=3,4} \int_{\tilde{\Sigma}_{p,j}^{m}} \sum_{i_{1},
i_{2}=1,2} \left\lvert \dfrac{(w_{+}^{\Sigma_{\tilde{\mathcal{R}}}}
(u))_{i_{1}i_{2}}}{u \! - \! z^{\prime}} \right\rvert^{2} \, \dfrac{\lvert 
\md u \rvert}{2 \pi^{2}} \right) \lvert \md \xi \rvert, \label{eqiy49} \\
\tilde{\mathbb{J}}_{\mathbb{D},5}^{\sharp} \! := \! \sum_{i=1}^{N+1} 
\left(\int_{\partial \tilde{\mathbb{U}}_{\tilde{\delta}_{\tilde{b}_{i-1}}}} 
\! + \! \int_{\partial \tilde{\mathbb{U}}_{\tilde{\delta}_{\tilde{a}_{i}}}} 
\right) \left(\lim_{-\tilde{\Sigma}_{\tilde{\mathcal{R}}}^{\sharp} \ni 
z^{\prime} \to \xi} \sum_{j=1}^{N+1} \sum_{m=3,4} \int_{\tilde{
\Sigma}_{p,j}^{m}} \sum_{i_{1},i_{2}=1,2} \left\lvert \dfrac{(w_{+}^{
\Sigma_{\tilde{\mathcal{R}}}}(u))_{i_{1}i_{2}}}{u \! - \! z^{\prime}} 
\right\rvert^{2} \, \dfrac{\lvert \md u \rvert}{2 \pi^{2}} \right) \lvert 
\md \xi \rvert. \label{eqiy50}
\end{gather}
For $n \! \in \! \mathbb{N}$ and $k \! \in \! \lbrace 1,2,\dotsc,K \rbrace$ 
such that $\alpha_{p_{\mathfrak{s}}} \! := \! \alpha_{k} \! \neq \! \infty$, 
via the formula $w^{\Sigma_{\tilde{\mathcal{R}}}}_{+}(z) \! = \! \tilde{v}_{
\tilde{\mathcal{R}}}(z) \! - \! \mathrm{I}$, the asymptotics, in the 
double-scaling limit $\mathscr{N},n \! \to \! \infty$ such that $z_{o} \! 
= \! 1 \! + \! o(1)$, of $w^{\Sigma_{\tilde{\mathcal{R}}}}_{+}(z)$ for $z 
\! \in \! \tilde{\Sigma}_{p,j}^{3}$ and for $z \! \in \! \tilde{\Sigma}_{p,
j}^{4}$, $j \! \in \! \lbrace 1,2,\dotsc,N \! + \! 1 \rbrace$, given in 
Equations~\eqref{eqtlvee11} and~\eqref{eqtlvee12}, respectively, the 
elementary trigonometric inequalities $\sin \theta \! \geqslant \! 
\tfrac{2 \theta}{\pi}$ and $\cos \theta \! \geqslant \! -\tfrac{2 \theta}{\pi} 
\! + \! 1$ for $0 \! \leqslant \! \theta \! \leqslant \! \tfrac{\pi}{2}$, 
$\sin \theta \! \geqslant \! \tfrac{2}{\pi}(\pi \! - \! \theta)$ and $\cos 
\theta \! \leqslant \! -\tfrac{2 \theta}{\pi} \! + \! 1$ for $\tfrac{\pi}{2} 
\! \leqslant \! \theta \! \leqslant \! \pi$, $\sin \theta \! \leqslant \! 
\tfrac{2 \theta}{\pi}$ and $\cos \theta \! \geqslant \! \tfrac{2 \theta}{\pi} 
\! + \! 1$ for $-\tfrac{\pi}{2} \! \leqslant \! \theta \! \leqslant \! 0$, and 
$\sin \theta \! \leqslant \! -\tfrac{2}{\pi}(\pi \! + \! \theta)$ and $\cos 
\theta \! \leqslant \! \tfrac{2 \theta}{\pi} \! + \! 1$ for $-\pi \! \leqslant 
\! \theta \! \leqslant \! -\tfrac{\pi}{2}$, the parametrisations $\tilde{
\Sigma}_{p,j}^{3} \! = \! \lbrace (x_{j}(\theta),y_{j}(\theta)); \, x_{j}(\theta) 
\! = \! A(j) \! + \! B(j) \cos \theta, y_{j}(\theta) \! = \! \tilde{\eta}_{j} 
\sin \theta, \, \theta_{0}^{\tilde{a}}(j) \! \leqslant \! \theta \! \leqslant 
\! \pi \! - \! \theta_{0}^{\tilde{b}}(j) \rbrace$ and $\tilde{\Sigma}_{p,j}^{4} 
\! = \! \lbrace (x_{j}(\theta),y_{j}(\theta)); \, x_{j}(\theta) \! = \! A(j) \! 
+ \! B(j) \cos \theta, y_{j}(\theta) \! = \! \tilde{\eta}_{j} \sin \theta, \, 
-\pi \! + \! \theta_{0}^{\tilde{b}}(j) \! \leqslant \! \theta \! \leqslant 
\! -\theta_{0}^{\tilde{a}}(j) \rbrace$, $j \! \in \! \lbrace 1,2,\dotsc,N 
\! + \! 1 \rbrace$, with $A(j) \! := \! \tfrac{1}{2}(\tilde{a}_{j} \! + 
\! \tilde{b}_{j-1})$ and $B(j) \! := \! \tfrac{1}{2}(\tilde{a}_{j} \! - \! 
\tilde{b}_{j-1})$, for the---elliptic---homotopic deformations of 
$\tilde{\Sigma}_{p,j}^{3}$ and $\tilde{\Sigma}_{p,j}^{4}$, respectively, 
and the inequality $\ln \lvert x \rvert \! \leqslant \! \lvert x \rvert 
\! - \! 1$, one shows, via a tedious integration argument, that, with 
$\Theta_{M}(j) \! := \! \max \lbrace \cos \theta_{0}^{\tilde{a}}(j),
\cos \theta_{0}^{\tilde{b}}(j) \rbrace$, $j \! \in \! \lbrace 1,2,\dotsc,
N \! + \! 1 \rbrace$,
\begin{align*}
\tilde{\mathbb{J}}_{\mathbb{D},1}^{\sharp} \underset{\underset{z_{o}
=1+o(1)}{\mathscr{N},n \to \infty}}{\leqslant}& \, \mathcal{O} \left(
\me^{-2((n-1)K+k) \min \left\lbrace \min\limits_{j=1,2,\dotsc,N+1} 
\lbrace \tilde{\lambda}_{\tilde{\mathcal{R}},4}^{\sharp,\smallfrown}(j) 
\rbrace,\min\limits_{j=1,2,\dotsc,N+1} \lbrace \tilde{\lambda}_{
\tilde{\mathcal{R}},4}^{\sharp,\smallsmile}(j) \rbrace \right\rbrace} 
\right) \int_{-\infty}^{\tilde{b}_{0}-\tilde{\delta}_{\tilde{b}_{0}}} 
\left(\sum_{j=1}^{N+1} \lim_{\epsilon \downarrow 0} \left(\int_{
\theta_{0}^{\tilde{a}}(j)}^{\pi -\theta_{0}^{\tilde{b}}(j)} \! + \! 
\int_{-\pi +\theta_{0}^{\tilde{b}}(j)}^{-\theta_{0}^{\tilde{a}}(j)} \right) 
\right. \\
&\left. \, \dfrac{((x_{j}^{\prime}(\theta))^{2} \! + \! (y_{j}^{\prime}
(\theta))^{2})^{1/2}}{(x_{j}(\theta) \! - \! \xi)^{2} \! + \! (y_{j}(\theta) \! + 
\! \epsilon)^{2}} \, \md \theta \right) \md \xi \underset{\underset{z_{o}
=1+o(1)}{\mathscr{N},n \to \infty}}{\leqslant} \mathcal{O} \left(
\me^{-2((n-1)K+k) \min \left\lbrace \min\limits_{j=1,2,\dotsc,N+1} 
\lbrace \tilde{\lambda}_{\tilde{\mathcal{R}},4}^{\sharp,\smallfrown}(j) 
\rbrace,\min\limits_{j=1,2,\dotsc,N+1} \lbrace \tilde{\lambda}_{
\tilde{\mathcal{R}},4}^{\sharp,\smallsmile}(j) \rbrace \right\rbrace} 
\right) \\
\times& \, \int_{-\infty}^{\tilde{b}_{0}-\tilde{\delta}_{\tilde{b}_{0}}} 
\left(\sum_{j=1}^{N+1} \lim_{\epsilon \downarrow 0} \left(\int_{
\theta_{0}^{\tilde{a}}(j)}^{\pi -\theta_{0}^{\tilde{b}}(j)} \! + \! 
\int_{-\pi +\theta_{0}^{\tilde{b}}(j)}^{-\theta_{0}^{\tilde{a}}(j)} \right) 
\dfrac{(B^{2}(j) \! + \! \tilde{\eta}_{j}^{2} \Theta_{M}^{2}(j))^{1/2}}{
(A(j) \! - \! B(j) \cos \theta_{0}^{\tilde{b}}(j) \! - \! \xi)^{2} \! + \! 
(-\tilde{\eta}_{j} \! + \! \epsilon)^{2}} \, \md \theta \right) \md 
\xi \\
\underset{\underset{z_{o}=1+o(1)}{\mathscr{N},n \to \infty}}{
\leqslant}& \, \mathcal{O} \left(\me^{-2((n-1)K+k) \min \left\lbrace 
\min\limits_{j=1,2,\dotsc,N+1} \lbrace \tilde{\lambda}_{\tilde{
\mathcal{R}},4}^{\sharp,\smallfrown}(j) \rbrace,\min\limits_{j=1,2,
\dotsc,N+1} \lbrace \tilde{\lambda}_{\tilde{\mathcal{R}},4}^{\sharp,
\smallsmile}(j) \rbrace \right\rbrace} \right) \sum_{j=1}^{N+1}
2(\pi \! - \! \theta_{0}^{\tilde{a}}(j) \! - \! \theta_{0}^{\tilde{b}}(j))
(B^{2}(j) \! + \! \tilde{\eta}_{j}^{2} \Theta_{M}^{2}(j))^{1/2} \\
\times& \, \lim_{\tau \to -\infty} \left(\int_{\tau}^{\tilde{b}_{0}
-\tilde{\delta}_{\tilde{b}_{0}}} \dfrac{1}{(A(j) \! - \! B(j) \cos 
\theta_{0}^{\tilde{b}}(j) \! - \! \xi)^{2} \! + \! \tilde{\eta}_{j}^{2}} \, 
\md \xi \right) \underset{\underset{z_{o}=1+o(1)}{\mathscr{N},
n \to \infty}}{\leqslant} \mathcal{O} \left(\me^{-2((n-1)K+k) 
\min \left\lbrace \min\limits_{j=1,2,\dotsc,N+1} \lbrace \tilde{
\lambda}_{\tilde{\mathcal{R}},4}^{\sharp,\smallfrown}(j) \rbrace,
\min\limits_{j=1,2,\dotsc,N+1} \lbrace \tilde{\lambda}_{\tilde{
\mathcal{R}},4}^{\sharp,\smallsmile}(j) \rbrace \right\rbrace} 
\right) \\
\times& \, \underbrace{\sum_{j=1}^{N+1} \dfrac{2(\pi \! - \! 
\theta_{0}^{\tilde{a}}(j) \! - \! \theta_{0}^{\tilde{b}}(j))(B^{2}(j) \! 
+ \! \tilde{\eta}_{j}^{2} \Theta_{M}^{2}(j))^{1/2}}{\tilde{\eta}_{j}} 
\left(\dfrac{\pi}{2} \! - \! \tan^{-1} \left(\dfrac{A(j) \! - \! B(j) \cos 
\theta_{0}^{\tilde{b}}(j) \! - \! \tilde{b}_{0} \! + \! \tilde{\delta}_{
\tilde{b}_{0}}}{\tilde{\eta}_{j}} \right) \right)}_{= \, \mathcal{O}(1)} 
\quad \Rightarrow
\end{align*}
\begin{equation} \label{eqiy51} 
\tilde{\mathbb{J}}_{\mathbb{D},1}^{\sharp} \underset{\underset{z_{o}
=1+o(1)}{\mathscr{N},n \to \infty}}{\leqslant} \mathcal{O} \left(
\tilde{\mathfrak{c}}_{\mathbb{D},1}^{\sharp}(n,k,z_{o}) \me^{-2((n-1)
K+k) \min \left\lbrace \min\limits_{j=1,2,\dotsc,N+1} \lbrace \tilde{
\lambda}_{\tilde{\mathcal{R}},4}^{\sharp,\smallfrown}(j) \rbrace,
\min\limits_{j=1,2,\dotsc,N+1} \lbrace \tilde{\lambda}_{\tilde{
\mathcal{R}},4}^{\sharp,\smallsmile}(j) \rbrace \right\rbrace} \right),
\end{equation}
where $\tilde{\mathfrak{c}}_{\mathbb{D},1}^{\sharp}(n,k,z_{o}) \! 
=_{\underset{z_{o}=1+o(1)}{\mathscr{N},n \to \infty}} \! \mathcal{O}
(1)$; proceeding analogously as above (for $\tilde{\mathbb{J}}_{
\mathbb{D},1}^{\sharp})$, one shows that, for $n \! \in \! \mathbb{N}$ 
and $k \! \in \! \lbrace 1,2,\dotsc,K \rbrace$ such that $\alpha_{
p_{\mathfrak{s}}} \! := \! \alpha_{k} \! \neq \! \infty$,
\begin{align*}
\tilde{\mathbb{J}}_{\mathbb{D},2}^{\sharp} \underset{\underset{z_{o}
=1+o(1)}{\mathscr{N},n \to \infty}}{\leqslant}& \, \mathcal{O} 
\left(\me^{-2((n-1)K+k) \min \left\lbrace \min\limits_{j=1,2,\dotsc,
N+1} \lbrace \tilde{\lambda}_{\tilde{\mathcal{R}},4}^{\sharp,
\smallfrown}(j) \rbrace,\min\limits_{j=1,2,\dotsc,N+1} \lbrace 
\tilde{\lambda}_{\tilde{\mathcal{R}},4}^{\sharp,\smallsmile}(j) \rbrace 
\right\rbrace} \right) \underbrace{\sum_{j=1}^{N+1} \dfrac{2(\pi \! 
- \! \theta_{0}^{\tilde{a}}(j) \! - \! \theta_{0}^{\tilde{b}}(j))(B^{2}(j) \! 
+ \! \tilde{\eta}_{j}^{2} \Theta_{M}^{2}(j))^{1/2}}{\tilde{\eta}_{j}}}_{= 
\, \mathcal{O}(1)} \\
\times& \, \underbrace{\left(\dfrac{\pi}{2} \! + \! \tan^{-1} \left(
\dfrac{A(j) \! - \! B(j) \cos \theta_{0}^{\tilde{b}}(j) \! - \! \tilde{a}_{
N+1} \! - \! \tilde{\delta}_{\tilde{a}_{N+1}}}{\tilde{\eta}_{j}} \right) 
\right)}_{= \, \mathcal{O}(1)} \quad \Rightarrow
\end{align*}
\begin{equation} \label{eqiy52} 
\tilde{\mathbb{J}}_{\mathbb{D},2}^{\sharp} \underset{\underset{z_{o}
=1+o(1)}{\mathscr{N},n \to \infty}}{\leqslant} \mathcal{O} \left(
\tilde{\mathfrak{c}}_{\mathbb{D},2}^{\sharp}(n,k,z_{o}) \me^{-2((n-1)
K+k) \min \left\lbrace \min\limits_{j=1,2,\dotsc,N+1} \lbrace \tilde{
\lambda}_{\tilde{\mathcal{R}},4}^{\sharp,\smallfrown}(j) \rbrace,
\min\limits_{j=1,2,\dotsc,N+1} \lbrace \tilde{\lambda}_{\tilde{
\mathcal{R}},4}^{\sharp,\smallsmile}(j) \rbrace \right\rbrace} \right),
\end{equation}
where $\tilde{\mathfrak{c}}_{\mathbb{D},2}^{\sharp}(n,k,z_{o}) \! 
=_{\underset{z_{o}=1+o(1)}{\mathscr{N},n \to \infty}} \! \mathcal{O}
(1)$, and
\begin{align*}
\tilde{\mathbb{J}}_{\mathbb{D},3}^{\sharp} \underset{\underset{z_{o}
=1+o(1)}{\mathscr{N},n \to \infty}}{\leqslant}& \, \mathcal{O} \left(
\me^{-2((n-1)K+k) \min \left\lbrace \min\limits_{j=1,2,\dotsc,N+1} 
\lbrace \tilde{\lambda}_{\tilde{\mathcal{R}},4}^{\sharp,\smallfrown}(j) 
\rbrace,\min\limits_{j=1,2,\dotsc,N+1} \lbrace \tilde{\lambda}_{
\tilde{\mathcal{R}},4}^{\sharp,\smallsmile}(j) \rbrace \right\rbrace} 
\right) \underbrace{\sum_{i=1}^{N} \sum_{j=1}^{N+1} \dfrac{2
(\pi \! - \! \theta_{0}^{\tilde{a}}(j) \! - \! \theta_{0}^{\tilde{b}}(j))(B^{2}(j) 
\! + \! \tilde{\eta}_{j}^{2} \Theta_{M}^{2}(j))^{1/2}}{\tilde{\eta}_{j}}}_{= \, 
\mathcal{O}(1)} \\
\times& \, \underbrace{\left(\tan^{-1} \left(\dfrac{A(j) \! - \! B(j) 
\cos \theta_{0}^{\tilde{b}}(j) \! - \! \tilde{a}_{i} \! - \! \tilde{\delta}_{
\tilde{a}_{i}}}{\tilde{\eta}_{j}} \right) \! - \! \tan^{-1} \left(\dfrac{A
(j) \! - \! B(j) \cos \theta_{0}^{\tilde{b}}(j) \! - \! \tilde{b}_{i} \! + \! 
\tilde{\delta}_{\tilde{b}_{i}}}{\tilde{\eta}_{j}} \right) \right)}_{= \, 
\mathcal{O}(1)} \quad \Rightarrow
\end{align*}
\begin{equation} \label{eqiy53} 
\tilde{\mathbb{J}}_{\mathbb{D},3}^{\sharp} \underset{\underset{z_{o}
=1+o(1)}{\mathscr{N},n \to \infty}}{\leqslant} \mathcal{O} \left(
\tilde{\mathfrak{c}}_{\mathbb{D},3}^{\sharp}(n,k,z_{o}) \me^{-2((n-1)
K+k) \min \left\lbrace \min\limits_{j=1,2,\dotsc,N+1} \lbrace \tilde{
\lambda}_{\tilde{\mathcal{R}},4}^{\sharp,\smallfrown}(j) \rbrace,
\min\limits_{j=1,2,\dotsc,N+1} \lbrace \tilde{\lambda}_{\tilde{
\mathcal{R}},4}^{\sharp,\smallsmile}(j) \rbrace \right\rbrace} \right),
\end{equation}
where $\tilde{\mathfrak{c}}_{\mathbb{D},3}^{\sharp}(n,k,z_{o}) \! 
=_{\underset{z_{o}=1+o(1)}{\mathscr{N},n \to \infty}} \! \mathcal{O}
(1)$. A straightforward calculation reveals that $\tilde{\mathbb{J}}_{
\mathbb{D},4}^{\sharp}$ can be presented as
\begin{equation} \label{eqiy54} 
\tilde{\mathbb{J}}_{\mathbb{D},4}^{\sharp} \! = \! \tilde{\mathbb{J}}_{
\mathbb{D},4}^{\sharp,\Ydown} \! + \!  \tilde{\mathbb{J}}_{\mathbb{D},
4}^{\sharp,\Yup} \! + \! \tilde{\mathbb{J}}_{\mathbb{D},4}^{\sharp,
\Yright} \! + \! \tilde{\mathbb{J}}_{\mathbb{D},4}^{\sharp,\Yleft},
\end{equation}
where
\begin{gather*}
\tilde{\mathbb{J}}_{\mathbb{D},4}^{\sharp,\Ydown} \! := \! \sum_{i
=1}^{N+1} \int_{\tilde{\Sigma}_{p,i}^{3}} \left(\lim_{-\tilde{\Sigma}_{
\tilde{\mathcal{R}}}^{\sharp} \ni z^{\prime} \to \xi} \sum_{j=1}^{N
+1} \int_{\tilde{\Sigma}_{p,j}^{4}} \sum_{i_{1},i_{2}=1,2} \left\lvert 
\dfrac{(w_{+}^{\Sigma_{\tilde{\mathcal{R}}}}(u))_{i_{1}i_{2}}}{u \! - \! 
z^{\prime}} \right\rvert^{2} \, \dfrac{\lvert \md u \rvert}{2 \pi^{2}} 
\right) \lvert \md \xi \rvert, \\
\tilde{\mathbb{J}}_{\mathbb{D},4}^{\sharp,\Yup} \! := \! \sum_{i=
1}^{N+1} \int_{\tilde{\Sigma}_{p,i}^{4}} \left(\lim_{-\tilde{\Sigma}_{
\tilde{\mathcal{R}}}^{\sharp} \ni z^{\prime} \to \xi} \sum_{j=1}^{N
+1} \int_{\tilde{\Sigma}_{p,j}^{3}} \sum_{i_{1},i_{2}=1,2} \left\lvert 
\dfrac{(w_{+}^{\Sigma_{\tilde{\mathcal{R}}}}(u))_{i_{1}i_{2}}}{u \! - \! 
z^{\prime}} \right\rvert^{2} \, \dfrac{\lvert \md u \rvert}{2 \pi^{2}} 
\right) \lvert \md \xi \rvert, \\
\tilde{\mathbb{J}}_{\mathbb{D},4}^{\sharp,\Yright} \! := \! \sum_{i
=1}^{N+1} \int_{\tilde{\Sigma}_{p,i}^{3}} \left(\lim_{-\tilde{\Sigma}_{
\tilde{\mathcal{R}}}^{\sharp} \ni z^{\prime} \to \xi} \sum_{j=1}^{N
+1} \int_{\tilde{\Sigma}_{p,j}^{3}} \sum_{i_{1},i_{2}=1,2} \left\lvert 
\dfrac{(w_{+}^{\Sigma_{\tilde{\mathcal{R}}}}(u))_{i_{1}i_{2}}}{u \! - \! 
z^{\prime}} \right\rvert^{2} \, \dfrac{\lvert \md u \rvert}{2 \pi^{2}} 
\right) \lvert \md \xi \rvert, \\
\tilde{\mathbb{J}}_{\mathbb{D},4}^{\sharp,\Yleft} \! := \! \sum_{i
=1}^{N+1} \int_{\tilde{\Sigma}_{p,i}^{4}} \left(\lim_{-\tilde{\Sigma}_{
\tilde{\mathcal{R}}}^{\sharp} \ni z^{\prime} \to \xi} \sum_{j=1}^{N
+1} \int_{\tilde{\Sigma}_{p,j}^{4}} \sum_{i_{1},i_{2}=1,2} \left\lvert 
\dfrac{(w_{+}^{\Sigma_{\tilde{\mathcal{R}}}}(u))_{i_{1}i_{2}}}{u \! - \! 
z^{\prime}} \right\rvert^{2} \, \dfrac{\lvert \md u \rvert}{2 \pi^{2}} 
\right) \lvert \md \xi \rvert.
\end{gather*}
Proceeding as in the calculations leading to the 
Estimates~\eqref{eqiy51}--\eqref{eqiy53}, one shows that, for $n 
\! \in \! \mathbb{N}$ and $k \! \in \! \lbrace 1,2,\dotsc,K \rbrace$ 
such that $\alpha_{p_{\mathfrak{s}}} \! := \! \alpha_{k} \! \neq \! 
\infty$,
\begin{align*}
\tilde{\mathbb{J}}_{\mathbb{D},4}^{\sharp,\Ydown} \underset{
\underset{z_{o}=1+o(1)}{\mathscr{N},n \to \infty}}{\leqslant}& \, 
\mathcal{O} \left(\me^{-2((n-1)K+k) \min \left\lbrace \min\limits_{j
=1,2,\dotsc,N+1} \lbrace \tilde{\lambda}_{\tilde{\mathcal{R}},4}^{
\sharp,\smallfrown}(j) \rbrace,\min\limits_{j=1,2,\dotsc,N+1} \lbrace 
\tilde{\lambda}_{\tilde{\mathcal{R}},4}^{\sharp,\smallsmile}(j) \rbrace 
\right\rbrace} \right) \sum_{i=1}^{N+1} \int_{\theta_{0}^{\tilde{a}}
(i)}^{\pi -\theta_{0}^{\tilde{b}}(i)} \left(\lim_{\epsilon \downarrow 0} 
\sum_{j=1}^{N+1} \int_{-\pi +\theta_{0}^{\tilde{b}}(j)}^{-\theta_{0}^{
\tilde{a}}(j)} \right. \\
&\left. \, \dfrac{((x_{j}^{\prime}(\theta))^{2} \! + \! (y_{j}^{\prime}
(\theta))^{2})^{1/2}}{(x_{j}(\theta) \! - \! x_{i}(\psi) \! + \! \epsilon \cos 
\psi)^{2} \! + \! (y_{j}(\theta) \! - \! y_{i}(\psi) \! + \! \epsilon \sin 
\psi)^{2}} \, \md \theta \right) ((x_{i}^{\prime}(\psi))^{2} \! + \! 
(y_{i}^{\prime}(\psi))^{2})^{1/2} \, \md \psi \\
\underset{\underset{z_{o}=1+o(1)}{\mathscr{N},n \to \infty}}{
\leqslant}& \, \mathcal{O} \left(\me^{-2((n-1)K+k) \min \left\lbrace 
\min\limits_{j=1,2,\dotsc,N+1} \lbrace \tilde{\lambda}_{\tilde{
\mathcal{R}},4}^{\sharp,\smallfrown}(j) \rbrace,\min\limits_{j=1,2,
\dotsc,N+1} \lbrace \tilde{\lambda}_{\tilde{\mathcal{R}},4}^{\sharp,
\smallsmile}(j) \rbrace \right\rbrace} \right) \\
\times& \, \underbrace{\sum_{i=1}^{N+1} \sum_{j=1}^{N+1} 
\dfrac{(B^{2}(j) \! + \! \tilde{\eta}_{j}^{2} \Theta_{M}^{2}(j))^{1/2}
(B^{2}(i) \! + \! \tilde{\eta}_{i}^{2} \Theta_{M}^{2}(i))^{1/2}(\pi \! 
- \! \theta_{0}^{\tilde{a}}(j) \! - \! \theta_{0}^{\tilde{b}}(j))(\pi \! - 
\! \theta_{0}^{\tilde{a}}(i) \! - \! \theta_{0}^{\tilde{b}}(i))}{(A(j) \! 
- \! A(i) \! - \! B(j) \cos \theta_{0}^{\tilde{b}}(j) \! - \! B(i) \cos 
\theta_{0}^{\tilde{a}}(i))^{2} \! + \! (\tilde{\eta}_{j} \! + \! 
\tilde{\eta}_{i})^{2}}}_{= \, \mathcal{O}(1)} \quad \Rightarrow 
\end{align*}
\begin{equation} \label{eqiy55} 
\tilde{\mathbb{J}}_{\mathbb{D},4}^{\sharp,\Ydown} \underset{
\underset{z_{o}=1+o(1)}{\mathscr{N},n \to \infty}}{\leqslant} 
\mathcal{O} \left(\tilde{\mathfrak{c}}_{\mathbb{D},4}^{\sharp,
\Ydown}(n,k,z_{o}) \me^{-2((n-1)K+k) \min \left\lbrace 
\min\limits_{j=1,2,\dotsc,N+1} \lbrace \tilde{\lambda}_{\tilde{
\mathcal{R}},4}^{\sharp,\smallfrown}(j) \rbrace,\min\limits_{j=
1,2,\dotsc,N+1} \lbrace \tilde{\lambda}_{\tilde{\mathcal{R}},
4}^{\sharp,\smallsmile}(j) \rbrace \right\rbrace} \right),
\end{equation}
where $\tilde{\mathfrak{c}}_{\mathbb{D},4}^{\sharp,\Ydown}(n,k,z_{o}) 
\! =_{\underset{z_{o}=1+o(1)}{\mathscr{N},n \to \infty}} \! \mathcal{O}(1)$,
\begin{align*}
\tilde{\mathbb{J}}_{\mathbb{D},4}^{\sharp,\Yup} \underset{\underset{
z_{o}=1+o(1)}{\mathscr{N},n \to \infty}}{\leqslant}& \, \mathcal{O} 
\left(\me^{-2((n-1)K+k) \min \left\lbrace \min\limits_{j=1,2,\dotsc,
N+1} \lbrace \tilde{\lambda}_{\tilde{\mathcal{R}},4}^{\sharp,\smallfrown}
(j) \rbrace,\min\limits_{j=1,2,\dotsc,N+1} \lbrace \tilde{\lambda}_{
\tilde{\mathcal{R}},4}^{\sharp,\smallsmile}(j) \rbrace \right\rbrace} 
\right) \sum_{i=1}^{N+1} \int_{-\pi +\theta_{0}^{\tilde{b}}(i)}^{
-\theta_{0}^{\tilde{a}}(i)} \left(\lim_{\epsilon \downarrow 0} \sum_{j=
1}^{N+1} \int_{\theta_{0}^{\tilde{a}}(j)}^{\pi -\theta_{0}^{\tilde{b}}(j)} 
\right. \\
&\left. \, \dfrac{((x_{j}^{\prime}(\theta))^{2} \! + \! (y_{j}^{\prime}
(\theta))^{2})^{1/2}}{(x_{j}(\theta) \! - \! x_{i}(\psi) \! - \! \epsilon \cos 
\psi)^{2} \! + \! (y_{j}(\theta) \! - \! y_{i}(\psi) \! - \! \epsilon \sin 
\psi)^{2}} \, \md \theta \right) ((x_{i}^{\prime}(\psi))^{2} \! + \! 
(y_{i}^{\prime}(\psi))^{2})^{1/2} \, \md \psi \\
\underset{\underset{z_{o}=1+o(1)}{\mathscr{N},n \to \infty}}{
\leqslant}& \, \mathcal{O} \left(\me^{-2((n-1)K+k) \min \left\lbrace 
\min\limits_{j=1,2,\dotsc,N+1} \lbrace \tilde{\lambda}_{\tilde{
\mathcal{R}},4}^{\sharp,\smallfrown}(j) \rbrace,\min\limits_{j=1,2,
\dotsc,N+1} \lbrace \tilde{\lambda}_{\tilde{\mathcal{R}},4}^{\sharp,
\smallsmile}(j) \rbrace \right\rbrace} \right) \\
\times& \, \underbrace{\sum_{i=1}^{N+1} \sum_{j=1}^{N+1} 
\dfrac{(B^{2}(j) \! + \! \tilde{\eta}_{j}^{2} \Theta_{M}^{2}(j))^{1/2}
(B^{2}(i) \! + \! \tilde{\eta}_{i}^{2} \Theta_{M}^{2}(i))^{1/2}(\pi \! - 
\! \theta_{0}^{\tilde{a}}(j) \! - \! \theta_{0}^{\tilde{b}}(j))(\pi \! - 
\! \theta_{0}^{\tilde{a}}(i) \! - \! \theta_{0}^{\tilde{b}}(i))}{(A(j) \! 
- \! A(i) \! - \! B(j) \cos \theta_{0}^{\tilde{b}}(j) \! - \! B(i) \cos 
\theta_{0}^{\tilde{a}}(i))^{2} \! + \! (\tilde{\eta}_{j} \Theta_{m}(j) 
\! + \! \tilde{\eta}_{i} \Theta_{m}(i))^{2}}}_{= \, \mathcal{O}(1)} 
\quad \Rightarrow 
\end{align*}
\begin{equation} \label{eqiy56} 
\tilde{\mathbb{J}}_{\mathbb{D},4}^{\sharp,\Yup} \underset{
\underset{z_{o}=1+o(1)}{\mathscr{N},n \to \infty}}{\leqslant} 
\mathcal{O} \left(\tilde{\mathfrak{c}}_{\mathbb{D},4}^{\sharp,
\Yup}(n,k,z_{o}) \me^{-2((n-1)K+k) \min \left\lbrace 
\min\limits_{j=1,2,\dotsc,N+1} \lbrace \tilde{\lambda}_{\tilde{
\mathcal{R}},4}^{\sharp,\smallfrown}(j) \rbrace,\min\limits_{j=
1,2,\dotsc,N+1} \lbrace \tilde{\lambda}_{\tilde{\mathcal{R}},
4}^{\sharp,\smallsmile}(j) \rbrace \right\rbrace} \right),
\end{equation}
where $\Theta_{m}(j) \! := \! \min \lbrace \sin \theta_{0}^{\tilde{a}}
(j),\sin \theta_{0}^{\tilde{b}}(j) \rbrace$, $j \! \in \! \lbrace 1,2,
\dotsc,N \! + \! 1 \rbrace$, and $\tilde{\mathfrak{c}}_{\mathbb{D},
4}^{\sharp,\Yup}(n,k,z_{o}) \! =_{\underset{z_{o}=1+o(1)}{\mathscr{N},
n \to \infty}} \! \mathcal{O}(1)$,
\begin{align*}
\tilde{\mathbb{J}}_{\mathbb{D},4}^{\sharp,\Yright} \underset{
\underset{z_{o}=1+o(1)}{\mathscr{N},n \to \infty}}{\leqslant}& \, 
\mathcal{O} \left(\me^{-2((n-1)K+k) \min \left\lbrace \min\limits_{j=
1,2,\dotsc,N+1} \lbrace \tilde{\lambda}_{\tilde{\mathcal{R}},4}^{\sharp,
\smallfrown}(j) \rbrace,\min\limits_{j=1,2,\dotsc,N+1} \lbrace \tilde{
\lambda}_{\tilde{\mathcal{R}},4}^{\sharp,\smallsmile}(j) \rbrace 
\right\rbrace} \right) \left(\sum_{i=1}^{N+1} \int_{\theta_{0}^{\tilde{a}}
(i)}^{\pi -\theta_{0}^{\tilde{b}}(i)} \left(\lim_{\epsilon \downarrow 0} 
\sum_{\substack{j=1\\j \neq i}}^{N+1} \int_{\theta_{0}^{\tilde{a}}
(j)}^{\pi -\theta_{0}^{\tilde{b}}(j)} \right. \right. \\
&\left. \left. \, \dfrac{((x_{j}^{\prime}(\theta))^{2} \! + \! (y_{j}^{
\prime}(\theta))^{2})^{1/2}}{(x_{j}(\theta) \! - \! x_{i}(\psi) \! + \! 
\epsilon \cos \psi)^{2} \! + \! (y_{j}(\theta) \! - \! y_{i}(\psi) \! + \! 
\epsilon \sin \psi)^{2}} \, \md \theta \right) ((x_{i}^{\prime}(\psi))^{2} 
\! + \! (y_{i}^{\prime}(\psi))^{2})^{1/2} \, \md \psi \right. \\
+&\left. \, \sum_{j=1}^{N+1} \int_{\theta_{0}^{\tilde{a}}(j)}^{\pi 
-\theta_{0}^{\tilde{b}}(j)} \left(\lim_{\epsilon \downarrow 0} \int_{
\theta_{0}^{\tilde{a}}(j)}^{\pi -\theta_{0}^{\tilde{b}}(j)} \dfrac{((
x_{j}^{\prime}(\theta))^{2} \! + \! (y_{j}^{\prime}(\theta))^{2})^{1/2}}{
(x_{j}(\theta) \! - \! x_{j}(\psi) \! + \! \epsilon \cos \psi)^{2} \! + \! 
(y_{j}(\theta) \! - \! y_{j}(\psi) \! + \! \epsilon \sin \psi)^{2}} \, \md 
\theta \right) \right. \\
\times&\left. \, ((x_{j}^{\prime}(\psi))^{2} \! + \! (y_{j}^{\prime}
(\psi))^{2})^{1/2} \, \md \psi \right) \underset{\underset{z_{o}=1+
o(1)}{\mathscr{N},n \to \infty}}{\leqslant} \mathcal{O} \left(\me^{-2
((n-1)K+k) \min \left\lbrace \min\limits_{j=1,2,\dotsc,N+1} \lbrace 
\tilde{\lambda}_{\tilde{\mathcal{R}},4}^{\sharp,\smallfrown}(j) 
\rbrace,\min\limits_{j=1,2,\dotsc,N+1} \lbrace \tilde{\lambda}_{
\tilde{\mathcal{R}},4}^{\sharp,\smallsmile}(j) \rbrace \right\rbrace} 
\right) \\
\times& \, \left(\sum_{i=1}^{N+1} \int_{\theta_{0}^{\tilde{a}}(i)}^{
\pi -\theta_{0}^{\tilde{b}}(i)} \left(\lim_{\epsilon \downarrow 0} 
\sum_{\substack{j=1\\j \neq i}}^{N+1} \int_{\theta_{0}^{\tilde{a}}
(j)}^{\pi -\theta_{0}^{\tilde{b}}(j)} \tfrac{(B^{2}(j) + \tilde{\eta}_{j}^{2} 
\Theta_{M}^{2}(j))^{1/2}}{(A(j) - A(i) - B(j) \cos \theta_{0}^{\tilde{b}}(j) 
- B(i) \cos \theta_{0}^{\tilde{a}}(i) - \epsilon \cos \theta_{0}^{\tilde{b}}
(i))^{2} + (\tilde{\eta}_{j} \Theta_{m}(j) - \tilde{\eta}_{i} + \epsilon 
\Theta_{m}(i))^{2}} \, \md \theta \right) \right. \\
\times&\left. \, (B^{2}(i) \! + \! \tilde{\eta}_{i}^{2} \Theta_{M}^{2}(i))^{1/2} 
\, \md \psi \! + \! \sum_{j=1}^{N+1} \int_{\theta_{0}^{\tilde{a}}(j)}^{\pi 
-\theta_{0}^{\tilde{b}}(j)} \left(\lim_{\epsilon \downarrow 0} \int_{
\theta_{0}^{\tilde{a}}(j)}^{\pi -\theta_{0}^{\tilde{b}}(j)} \dfrac{(B^{2}(j) 
\! + \! \tilde{\eta}_{j}^{2} \Theta_{M}^{2}(j))}{r_{j}^{\epsilon} \! + \! 
p_{j}^{\epsilon}(\psi) \sin \theta \! + \! q_{j}^{\epsilon}(\psi) \cos 
\theta} \, \md \theta \right) \md \psi \right),
\end{align*}
where, for $j \! \in \! \lbrace 1,2,\dotsc,N \! + \! 1 \rbrace$, $r_{j}^{
\epsilon} \! := \! 2P_{m}(j) \! - \! 2 \epsilon Q_{M}(j) \! + \! \epsilon^{
2}$, $p_{j}^{\epsilon}(\psi) \! := \! (-2P_{M}(j) \! + \! 2 \epsilon Q_{m}
(j)) \sin \psi$, and $q_{j}^{\epsilon}(\psi) \! := \! (-2P_{M}(j) \! + \! 2 
\epsilon Q_{m}(j)) \cos \psi$, with $P_{m}(j) \! = \! \min \lbrace B^{2}
(j),\tilde{\eta}^{2}_{j} \rbrace$, $Q_{M}(j) \! = \! \max \lbrace B(j),
\tilde{\eta}_{j} \rbrace$, $P_{M}(j) \! = \! \max \lbrace B^{2}(j),
\tilde{\eta}^{2}_{j} \rbrace$, and $Q_{m}(j) \! = \! \min \lbrace B(j),
\tilde{\eta}_{j} \rbrace$, thus, via the indefinite integral \pmb{2.558} 
4. on p.~182 of \cite{gradryzh},
\begin{align*}
\tilde{\mathbb{J}}_{\mathbb{D},4}^{\sharp,\Yright} \underset{
\underset{z_{o}=1+o(1)}{\mathscr{N},n \to \infty}}{\leqslant}& \, 
\mathcal{O} \left(\me^{-2((n-1)K+k) \min \left\lbrace \min\limits_{j=
1,2,\dotsc,N+1} \lbrace \tilde{\lambda}_{\tilde{\mathcal{R}},4}^{\sharp,
\smallfrown}(j) \rbrace,\min\limits_{j=1,2,\dotsc,N+1} \lbrace \tilde{
\lambda}_{\tilde{\mathcal{R}},4}^{\sharp,\smallsmile}(j) \rbrace 
\right\rbrace} \right) \\
\times& \, \left(\underbrace{\sum_{i=1}^{N+1} \sum_{\substack{j=
1\\j \neq i}}^{N+1} \dfrac{(B^{2}(j) \! + \! \tilde{\eta}_{j}^{2} \Theta_{M}^{2}
(j))^{1/2}(B^{2}(i) \! + \! \tilde{\eta}_{i}^{2} \Theta_{M}^{2}(i))^{1/2}(\pi 
\! - \! \theta_{0}^{\tilde{a}}(j) \! - \! \theta_{0}^{\tilde{b}}(j))(\pi \! - \! 
\theta_{0}^{\tilde{a}}(i) \! - \! \theta_{0}^{\tilde{b}}(i))}{(A(j) \! - \! A(i) \! 
- \! B(j) \cos \theta_{0}^{\tilde{b}}(j) \! - \! B(i) \cos \theta_{0}^{\tilde{a}}
(i))^{2} \! + \! (\tilde{\eta}_{j} \Theta_{m}(j) \! - \! \tilde{\eta}_{i})^{2}}}_{= 
\, \mathcal{O}(1)} \right. \\
+&\left. \, \sum_{j=1}^{N+1} \int_{\theta_{0}^{\tilde{a}}(j)}^{\pi 
-\theta_{0}^{\tilde{b}}(j)} \left(\lim_{\epsilon \downarrow 0} \left. 
\tfrac{(B^{2}(j)+ \tilde{\eta}_{j}^{2} \Theta_{M}^{2}(j))}{\sqrt{(p_{j}^{
\epsilon}(\psi))^{2}+(q_{j}^{\epsilon}(\psi))^{2}-(r_{j}^{\epsilon})^{2}}} 
\ln \left(\tfrac{p_{j}^{\epsilon}(\psi)-\sqrt{(p_{j}^{\epsilon}(\psi))^{2}
+(q_{j}^{\epsilon}(\psi))^{2}-(r_{j}^{\epsilon})^{2}}+(r_{j}^{\epsilon}-
q_{j}^{\epsilon}(\psi)) \tan (\theta/2)}{p_{j}^{\epsilon}(\psi)+
\sqrt{(p_{j}^{\epsilon}(\psi))^{2}+(q_{j}^{\epsilon}(\psi))^{2}-
(r_{j}^{\epsilon})^{2}}+(r_{j}^{\epsilon}-q_{j}^{\epsilon}(\psi)) \tan 
(\theta/2)} \right) \right\rvert_{\theta_{0}^{\tilde{a}}(j)}^{\pi 
-\theta_{0}^{\tilde{b}}(j)} \right) \right) \\
\underset{\underset{z_{o}=1+o(1)}{\mathscr{N},n \to \infty}}{\leqslant}& 
\, \mathcal{O} \left(\me^{-2((n-1)K+k) \min \left\lbrace \min\limits_{j=
1,2,\dotsc,N+1} \lbrace \tilde{\lambda}_{\tilde{\mathcal{R}},4}^{\sharp,
\smallfrown}(j) \rbrace,\min\limits_{j=1,2,\dotsc,N+1} \lbrace \tilde{
\lambda}_{\tilde{\mathcal{R}},4}^{\sharp,\smallsmile}(j) \rbrace 
\right\rbrace} \right) \left(\mathcal{O}(1) \! + \! \sum_{j=1}^{N+1} 
\tfrac{(B^{2}(j)+\tilde{\eta}_{j}^{2} \Theta_{M}^{2}(j))}{2(P_{M}^{2}(j)-
P_{m}^{2}(j))^{1/2}} \int_{\theta_{0}^{\tilde{a}}(j)}^{\pi -\theta_{0}^{
\tilde{b}}(j)} \ln \left(\tfrac{\mathrm{F}_{j}(\psi)+\mathrm{G}_{j}
(\psi)}{\mathrm{F}_{j}(\psi)-\mathrm{G}_{j}(\psi)} \right) \md \psi \right),
\end{align*}
where $\mathrm{F}_{j}(\psi) \! := \! P_{M}^{2}(j) \sin^{2} \psi \! - \! 
(P_{M}^{2}(j) \! - \! P_{m}^{2}(j)) \! - \! P_{M}(j) \sin (\psi)(P_{m}(j) \! 
+ \! P_{M}(j) \cos \psi)(\tan (\theta_{0}^{\tilde{a}}(j)/2) \! + \! \cot 
(\theta_{0}^{\tilde{b}}(j)/2)) \! + \! (P_{m}(j) \! + \! P_{M}(j) \cos \psi)^{2} 
\tan (\theta_{0}^{\tilde{a}}(j)/2) \cot (\theta_{0}^{\tilde{b}}(j)/2)$, and 
$\mathrm{G}_{j}(\psi) \! := \! (P_{M}^{2}(j) \! - \! P_{m}^{2}(j))^{1/2}
(P_{m}(j) \! + \! P_{M}(j) \cos \psi)(\cot (\theta_{0}^{\tilde{b}}(j)/2) \! 
- \! \tan (\theta_{0}^{\tilde{a}}(j)/2))$, thus, with $\mathfrak{l}_{0}(j) 
\! := \! P_{m}^{2}(j) \! - \! P_{M}(j) \Theta_{m}(j)(P_{m}(j) \! + \! P_{M}
(j) \cos \theta_{0}^{\tilde{a}}(j))(\tan (\theta_{0}^{\tilde{a}}(j)/2) \! + \! 
\cot (\theta_{0}^{\tilde{b}}(j)/2)) \! + \! (P_{m}(j) \! + \! P_{M}(j) \cos 
\theta_{0}^{\tilde{a}}(j))^{2} \tan (\theta_{0}^{\tilde{a}}(j)/2) \cot 
(\theta_{0}^{\tilde{b}}(j)/2)$, $\mathfrak{l}_{1}(j) \! := \! P_{M}^{2}(j) 
\Theta_{m}^{2}(j) \! - \! (P_{M}^{2}(j) \! - \! P_{m}^{2}(j)) \! - \! P_{M}
(j)(P_{m}(j) \! - \! P_{M}(j) \cos \theta_{0}^{\tilde{b}}(j)) \linebreak[4] 
\pmb{\cdot} (\tan (\theta_{0}^{\tilde{a}}(j)/2) \! + \! \cot (\theta_{
0}^{\tilde{b}}(j)/2)) \! + \! (P_{m}(j) \! - \! P_{M}(j) \cos \theta_{0}^{
\tilde{b}}(j))^{2} \tan (\theta_{0}^{\tilde{a}}(j)/2) \cot (\theta_{0}^{
\tilde{b}}(j)/2)$, and $\mathfrak{l}_{2}(j) \! := \! (P_{M}^{2}(j) 
\! - \! P_{m}^{2}(j))^{1/2}(P_{m}(j) \linebreak[4] 
\! + \! P_{M}(j) \cos \theta_{0}^{\tilde{a}}(j))(\cot (\theta_{0}^{
\tilde{b}}(j)/2) \! - \! \tan (\theta_{0}^{\tilde{a}}(j)/2))$, via the 
monotonicity of $\ln (\pmb{\cdot})$ and an elementary inequality 
argument, one arrives at
\begin{equation*}
\tilde{\mathbb{J}}_{\mathbb{D},4}^{\sharp,\Yright} \underset{
\underset{z_{o}=1+o(1)}{\mathscr{N},n \to \infty}}{\leqslant} \mathcal{O} 
\left(\me^{-2((n-1)K+k) \min \left\lbrace \min\limits_{j=1,2,\dotsc,
N+1} \lbrace \tilde{\lambda}_{\tilde{\mathcal{R}},4}^{\sharp,
\smallfrown}(j) \rbrace,\min\limits_{j=1,2,\dotsc,N+1} \lbrace 
\tilde{\lambda}_{\tilde{\mathcal{R}},4}^{\sharp,\smallsmile}(j) 
\rbrace \right\rbrace} \right) \left(\mathcal{O}(1) \! + \! 
\underbrace{\sum_{j=1}^{N+1} \tfrac{(B^{2}(j)+ \tilde{\eta}_{j}^{2} 
\Theta_{M}^{2}(j))(\pi -\theta_{0}^{\tilde{a}}(j)-\theta_{0}^{\tilde{b}}
(j))}{2(P_{M}^{2}(j)-P_{m}^{2}(j))^{1/2}} \ln \left(\tfrac{\mathfrak{l}_{0}
(j)+ \mathfrak{l}_{2}(j)}{\mathfrak{l}_{1}(j)- \mathfrak{l}_{2}(j)} 
\right)}_{= \, \mathcal{O}(1)} \right) \quad \Rightarrow
\end{equation*}
\begin{equation} \label{eqiy57} 
\tilde{\mathbb{J}}_{\mathbb{D},4}^{\sharp,\Yright} \underset{
\underset{z_{o}=1+o(1)}{\mathscr{N},n \to \infty}}{\leqslant} 
\mathcal{O} \left(\tilde{\mathfrak{c}}_{\mathbb{D},4}^{\sharp,\Yright}
(n,k,z_{o}) \me^{-2((n-1)K+k) \min \left\lbrace \min\limits_{j=1,2,
\dotsc,N+1} \lbrace \tilde{\lambda}_{\tilde{\mathcal{R}},4}^{\sharp,
\smallfrown}(j) \rbrace,\min\limits_{j=1,2,\dotsc,N+1} \lbrace 
\tilde{\lambda}_{\tilde{\mathcal{R}},4}^{\sharp,\smallsmile}(j) 
\rbrace \right\rbrace} \right),
\end{equation}
where $\tilde{\mathfrak{c}}_{\mathbb{D},4}^{\sharp,\Yright}(n,k,
z_{o}) \! =_{\underset{z_{o}=1+o(1)}{\mathscr{N},n \to \infty}} \! 
\mathcal{O}(1)$, and, proceeding as above (for $\tilde{\mathbb{
J}}_{\mathbb{D},4}^{\sharp,\Yright})$,
\begin{align*}
\tilde{\mathbb{J}}_{\mathbb{D},4}^{\sharp,\Yleft} \underset{
\underset{z_{o}=1+o(1)}{\mathscr{N},n \to \infty}}{\leqslant}& \, 
\mathcal{O} \left(\me^{-2((n-1)K+k) \min \left\lbrace \min\limits_{j=
1,2,\dotsc,N+1} \lbrace \tilde{\lambda}_{\tilde{\mathcal{R}},4}^{\sharp,
\smallfrown}(j) \rbrace,\min\limits_{j=1,2,\dotsc,N+1} \lbrace \tilde{
\lambda}_{\tilde{\mathcal{R}},4}^{\sharp,\smallsmile}(j) \rbrace 
\right\rbrace} \right) \left(\sum_{i=1}^{N+1} \int_{-\pi +\theta_{
0}^{\tilde{b}}(i)}^{-\theta_{0}^{\tilde{a}}(i)} \left(\lim_{\epsilon 
\downarrow 0} \sum_{\substack{j=1\\j \neq i}}^{N+1} \int_{-\pi +
\theta_{0}^{\tilde{b}}(j)}^{-\theta_{0}^{\tilde{a}}(j)} \right. \right. \\
&\left. \left. \, \dfrac{((x_{j}^{\prime}(\theta))^{2} \! + \! (y_{j}^{\prime}
(\theta))^{2})^{1/2}}{(x_{j}(\theta) \! - \! x_{i}(\psi) \! - \! \epsilon \cos 
\psi)^{2} \! + \! (y_{j}(\theta) \! - \! y_{i}(\psi) \! - \! \epsilon \sin 
\psi)^{2}} \, \md \theta \right) ((x_{i}^{\prime}(\psi))^{2} \! + \! 
(y_{i}^{\prime}(\psi))^{2})^{1/2} \, \md \psi \right. \\
+&\left. \, \sum_{j=1}^{N+1} \int_{-\pi +\theta_{0}^{\tilde{b}}
(j)}^{-\theta_{0}^{\tilde{a}}(j)} \left(\lim_{\epsilon \downarrow 0} 
\int_{-\pi +\theta_{0}^{\tilde{b}}(j)}^{-\theta_{0}^{\tilde{a}}(j)} \dfrac{
((x_{j}^{\prime}(\theta))^{2} \! + \! (y_{j}^{\prime}(\theta))^{2})^{1/2}}{
(x_{j}(\theta) \! - \! x_{j}(\psi) \! - \! \epsilon \cos \psi)^{2} \! + \! 
(y_{j}(\theta) \! - \! y_{j}(\psi) \! - \! \epsilon \sin \psi)^{2}} \, \md 
\theta \right) \right. \\
\times&\left. \, ((x_{j}^{\prime}(\psi))^{2} \! + \! (y_{j}^{\prime}
(\psi))^{2})^{1/2} \, \md \psi \right) \underset{\underset{z_{o}=1+
o(1)}{\mathscr{N},n \to \infty}}{\leqslant} \mathcal{O} \left(\me^{-2
((n-1)K+k) \min \left\lbrace \min\limits_{j=1,2,\dotsc,N+1} \lbrace 
\tilde{\lambda}_{\tilde{\mathcal{R}},4}^{\sharp,\smallfrown}(j) 
\rbrace,\min\limits_{j=1,2,\dotsc,N+1} \lbrace \tilde{\lambda}_{
\tilde{\mathcal{R}},4}^{\sharp,\smallsmile}(j) \rbrace \right\rbrace} 
\right) \\
\times& \, \left(\underbrace{\sum_{i=1}^{N+1} \sum_{\substack{j
=1\\j \neq i}}^{N+1} \tfrac{(B^{2}(j)+\tilde{\eta}_{j}^{2} \Theta_{M}^{2}
(j))^{1/2}(B^{2}(i)+\tilde{\eta}_{i}^{2} \Theta_{M}^{2}(i))^{1/2}(\pi -
\theta_{0}^{\tilde{a}}(j)-\theta_{0}^{\tilde{b}}(j))(\pi -\theta_{0}^{\tilde{a}}
(i)-\theta_{0}^{\tilde{b}}(i))}{(A(j)-A(i)-B(j) \cos \theta_{0}^{\tilde{b}}(j)-
B(i) \cos \theta_{0}^{\tilde{a}}(i))^{2}+(-\tilde{\eta}_{j}+ \tilde{\eta}_{i} 
\Theta_{m}(i))^{2}}}_{= \, \mathcal{O}(1)} \! + \! \underbrace{\sum_{j=1}^{N+1} 
\tfrac{(B^{2}(j)+ \tilde{\eta}_{j}^{2} \Theta_{M}^{2}(j))(\pi -\theta_{0}^{
\tilde{a}}(j)-\theta_{0}^{\tilde{b}}(j))}{2(P_{M}^{2}(j)-P_{m}^{2}(j))^{1/2}}}_{= 
\, \mathcal{O}(1)} \right. \\
\times&\left. \underbrace{\ln \left(\tfrac{\mathfrak{l}_{0}(j)-(P_{M}^{2}(j)
-P_{m}^{2}(j))^{1/2}(P_{m}(j)-P_{M}(j) \cos \theta_{0}^{\tilde{b}}(j))(\tan 
(\theta_{0}^{\tilde{a}}(j)/2)-\cot (\theta_{0}^{\tilde{b}}(j)/2))}{\mathfrak{l}_{1}
(j)+(P_{M}^{2}(j)-P_{m}^{2}(j))^{1/2}(P_{m}(j)-P_{M}(j) \cos \theta_{0}^{\tilde{b}}
(j))(\tan (\theta_{0}^{\tilde{a}}(j)/2)-\cot (\theta_{0}^{\tilde{b}}(j)/2))} \right)}_{= 
\, \mathcal{O}(1)} \right) \quad \Rightarrow
\end{align*}
\begin{equation} \label{eqiy58} 
\tilde{\mathbb{J}}_{\mathbb{D},4}^{\sharp,\Yleft} \underset{
\underset{z_{o}=1+o(1)}{\mathscr{N},n \to \infty}}{\leqslant} 
\mathcal{O} \left(\tilde{\mathfrak{c}}_{\mathbb{D},4}^{\sharp,
\Yleft}(n,k,z_{o}) \me^{-2((n-1)K+k) \min \left\lbrace \min\limits_{j
=1,2,\dotsc,N+1} \lbrace \tilde{\lambda}_{\tilde{\mathcal{R}},4}^{
\sharp,\smallfrown}(j) \rbrace,\min\limits_{j=1,2,\dotsc,N+1} 
\lbrace \tilde{\lambda}_{\tilde{\mathcal{R}},4}^{\sharp,\smallsmile}
(j) \rbrace \right\rbrace} \right),
\end{equation}
where $\tilde{\mathfrak{c}}_{\mathbb{D},4}^{\sharp,\Yleft}(n,k,z_{o}) 
\! =_{\underset{z_{o}=1+o(1)}{\mathscr{N},n \to \infty}} \! \mathcal{O}
(1)$; thus, via the Estimates~\eqref{eqiy55}--\eqref{eqiy58} and 
Equation~\eqref{eqiy54}, one arrives at, for $n \! \in \! \mathbb{N}$ 
and $k \! \in \! \lbrace 1,2,\dotsc,K \rbrace$ such that $\alpha_{p_{
\mathfrak{s}}} \! := \! \alpha_{k} \! \neq \! \infty$,
\begin{equation} \label{eqiy59} 
\tilde{\mathbb{J}}_{\mathbb{D},4}^{\sharp} \underset{\underset{z_{o}
=1+o(1)}{\mathscr{N},n \to \infty}}{\leqslant} \mathcal{O} \left(\tilde{
\mathfrak{c}}_{\mathbb{D},4}^{\sharp}(n,k,z_{o}) \me^{-2((n-1)K+k) 
\min \left\lbrace \min\limits_{j=1,2,\dotsc,N+1} \lbrace \tilde{
\lambda}_{\tilde{\mathcal{R}},4}^{\sharp,\smallfrown}(j) \rbrace,
\min\limits_{j=1,2,\dotsc,N+1} \lbrace \tilde{\lambda}_{\tilde{
\mathcal{R}},4}^{\sharp,\smallsmile}(j) \rbrace \right\rbrace} \right),
\end{equation}
where $\tilde{\mathfrak{c}}_{\mathbb{D},4}^{\sharp}(n,k,z_{o}) \! 
=_{\underset{z_{o}=1+o(1)}{\mathscr{N},n \to \infty}} \! \mathcal{O}(1)$. 
Proceeding as in the calculations leading to the Estimate~\eqref{eqiy59}, 
one shows, analogously, that, for $n \! \in \! \mathbb{N}$ and $k \! \in 
\! \lbrace 1,2,\dotsc,K \rbrace$ such that $\alpha_{p_{\mathfrak{s}}} \! 
:= \! \alpha_{k} \! \neq \! \infty$:
\begin{equation} \label{eqiy60} 
\tilde{\mathbb{J}}_{\mathbb{D},5}^{\sharp} \! = \! \tilde{\mathbb{J}}_{
\mathbb{D},5}^{\sharp,\Ydown} \! + \! \tilde{\mathbb{J}}_{\mathbb{D},
5}^{\sharp,\Yup} \! + \! \tilde{\mathbb{J}}_{\mathbb{D},5}^{\sharp,
\Yright} \! + \! \tilde{\mathbb{J}}_{\mathbb{D},5}^{\sharp,\Yleft},
\end{equation}
where
\begin{align*}
\tilde{\mathbb{J}}_{\mathbb{D},5}^{\sharp,\Ydown} :=& \, \sum_{i
=1}^{N+1} \int_{\partial \tilde{\mathbb{U}}_{\tilde{\delta}_{\tilde{b}_{i
-1}}}} \left(\lim_{-\tilde{\Sigma}_{\tilde{\mathcal{R}}}^{\sharp} \ni 
z^{\prime} \to \xi} \sum_{j=1}^{N+1} \int_{\tilde{\Sigma}_{p,j}^{3}} 
\sum_{i_{1},i_{2}=1,2} \left\lvert \dfrac{(w_{+}^{\Sigma_{\tilde{
\mathcal{R}}}}(u))_{i_{1}i_{2}}}{u \! - \! z^{\prime}} \right\rvert^{2} 
\, \dfrac{\lvert \md u \rvert}{2 \pi^{2}} \right) \lvert \md \xi \rvert \\
\underset{\underset{z_{o}=1+o(1)}{\mathscr{N},n \to \infty}}{
\leqslant}& \, \mathcal{O} \left(\me^{-2((n-1)K+k) \min \left\lbrace 
\min\limits_{j=1,2,\dotsc,N+1} \lbrace \tilde{\lambda}_{\tilde{\mathcal{
R}},4}^{\sharp,\smallfrown}(j) \rbrace,\min\limits_{j=1,2,\dotsc,N+1} 
\lbrace \tilde{\lambda}_{\tilde{\mathcal{R}},4}^{\sharp,\smallsmile}(j) 
\rbrace \right\rbrace} \right) \underbrace{\sum_{i=1}^{N+1} \sum_{j=
1}^{N+1} \tfrac{2 \pi \tilde{\delta}_{\tilde{b}_{i-1}}(B^{2}(j)+\tilde{
\eta}_{j}^{2} \Theta_{M}^{2}(j))^{1/2}(\pi -\theta_{0}^{\tilde{a}}(j)-\theta_{
0}^{\tilde{b}}(j))}{(A(j)-B(j) \cos \theta_{0}^{\tilde{b}}(j)-\tilde{b}_{i-1}
-\tilde{\delta}_{\tilde{b}_{i-1}})^{2}+(\tilde{\eta}_{j} \Theta_{m}(j)-
\tilde{\delta}_{\tilde{b}_{i-1}})^{2}}}_{= \, \mathcal{O}(1)} \\
\underset{\underset{z_{o}=1+o(1)}{\mathscr{N},n \to \infty}}{
\leqslant}& \, \mathcal{O} \left(\tilde{\mathfrak{c}}_{\mathbb{D},5}^{
\sharp,\Ydown}(n,k,z_{o}) \me^{-2((n-1)K+k) \min \left\lbrace 
\min\limits_{j=1,2,\dotsc,N+1} \lbrace \tilde{\lambda}_{\tilde{
\mathcal{R}},4}^{\sharp,\smallfrown}(j) \rbrace,\min\limits_{j=1,2,
\dotsc,N+1} \lbrace \tilde{\lambda}_{\tilde{\mathcal{R}},4}^{\sharp,
\smallsmile}(j) \rbrace \right\rbrace} \right),
\end{align*}
\begin{align*}
\tilde{\mathbb{J}}_{\mathbb{D},5}^{\sharp,\Yup} :=& \, \sum_{i=
1}^{N+1} \int_{\partial \tilde{\mathbb{U}}_{\tilde{\delta}_{\tilde{b}_{i
-1}}}} \left(\lim_{-\tilde{\Sigma}_{\tilde{\mathcal{R}}}^{\sharp} \ni 
z^{\prime} \to \xi} \sum_{j=1}^{N+1} \int_{\tilde{\Sigma}_{p,j}^{4}} 
\sum_{i_{1},i_{2}=1,2} \left\lvert \dfrac{(w_{+}^{\Sigma_{\tilde{
\mathcal{R}}}}(u))_{i_{1}i_{2}}}{u \! - \! z^{\prime}} \right\rvert^{2} 
\, \dfrac{\lvert \md u \rvert}{2 \pi^{2}} \right) \lvert \md \xi \rvert \\
\underset{\underset{z_{o}=1+o(1)}{\mathscr{N},n \to \infty}}{
\leqslant}& \, \mathcal{O} \left(\me^{-2((n-1)K+k) \min \left\lbrace 
\min\limits_{j=1,2,\dotsc,N+1} \lbrace \tilde{\lambda}_{\tilde{\mathcal{
R}},4}^{\sharp,\smallfrown}(j) \rbrace,\min\limits_{j=1,2,\dotsc,N+1} 
\lbrace \tilde{\lambda}_{\tilde{\mathcal{R}},4}^{\sharp,\smallsmile}(j) 
\rbrace \right\rbrace} \right) \underbrace{\sum_{i=1}^{N+1} \sum_{j=
1}^{N+1} \tfrac{2 \pi \tilde{\delta}_{\tilde{b}_{i-1}}(B^{2}(j)+\tilde{
\eta}_{j}^{2} \Theta_{M}^{2}(j))^{1/2}(\pi -\theta_{0}^{\tilde{a}}(j)-\theta_{
0}^{\tilde{b}}(j))}{(A(j)-B(j) \cos \theta_{0}^{\tilde{b}}(j)-\tilde{b}_{i-1}
-\tilde{\delta}_{\tilde{b}_{i-1}})^{2}+(\tilde{\eta}_{j}+\tilde{\delta}_{
\tilde{b}_{i-1}})^{2}}}_{= \, \mathcal{O}(1)} \\
\underset{\underset{z_{o}=1+o(1)}{\mathscr{N},n \to \infty}}{
\leqslant}& \, \mathcal{O} \left(\tilde{\mathfrak{c}}_{\mathbb{D},5}^{
\sharp,\Yup}(n,k,z_{o}) \me^{-2((n-1)K+k) \min \left\lbrace 
\min\limits_{j=1,2,\dotsc,N+1} \lbrace \tilde{\lambda}_{\tilde{
\mathcal{R}},4}^{\sharp,\smallfrown}(j) \rbrace,\min\limits_{j=1,2,
\dotsc,N+1} \lbrace \tilde{\lambda}_{\tilde{\mathcal{R}},4}^{\sharp,
\smallsmile}(j) \rbrace \right\rbrace} \right),
\end{align*}
\begin{align*}
\tilde{\mathbb{J}}_{\mathbb{D},5}^{\sharp,\Yright} :=& \, \sum_{i=
1}^{N+1} \int_{\partial \tilde{\mathbb{U}}_{\tilde{\delta}_{\tilde{a}_{i}}}} 
\left(\lim_{-\tilde{\Sigma}_{\tilde{\mathcal{R}}}^{\sharp} \ni z^{\prime} 
\to \xi} \sum_{j=1}^{N+1} \int_{\tilde{\Sigma}_{p,j}^{3}} \sum_{i_{1},
i_{2}=1,2} \left\lvert \dfrac{(w_{+}^{\Sigma_{\tilde{\mathcal{R}}}}
(u))_{i_{1}i_{2}}}{u \! - \! z^{\prime}} \right\rvert^{2} \, \dfrac{\lvert 
\md u \rvert}{2 \pi^{2}} \right) \lvert \md \xi \rvert \\
\underset{\underset{z_{o}=1+o(1)}{\mathscr{N},n \to \infty}}{
\leqslant}& \, \mathcal{O} \left(\me^{-2((n-1)K+k) \min \left\lbrace 
\min\limits_{j=1,2,\dotsc,N+1} \lbrace \tilde{\lambda}_{\tilde{\mathcal{
R}},4}^{\sharp,\smallfrown}(j) \rbrace,\min\limits_{j=1,2,\dotsc,N+1} 
\lbrace \tilde{\lambda}_{\tilde{\mathcal{R}},4}^{\sharp,\smallsmile}(j) 
\rbrace \right\rbrace} \right) \underbrace{\sum_{i=1}^{N+1} \sum_{j=
1}^{N+1} \tfrac{2 \pi \tilde{\delta}_{\tilde{a}_{i}}(B^{2}(j)+\tilde{
\eta}_{j}^{2} \Theta_{M}^{2}(j))^{1/2}(\pi -\theta_{0}^{\tilde{a}}(j)-\theta_{
0}^{\tilde{b}}(j))}{(A(j)-B(j) \cos \theta_{0}^{\tilde{b}}(j)-\tilde{a}_{i}-
\tilde{\delta}_{\tilde{a}_{i}})^{2}+(\tilde{\eta}_{j} \Theta_{m}(j)-\tilde{
\delta}_{\tilde{a}_{i}})^{2}}}_{= \, \mathcal{O}(1)} \\
\underset{\underset{z_{o}=1+o(1)}{\mathscr{N},n \to \infty}}{
\leqslant}& \, \mathcal{O} \left(\tilde{\mathfrak{c}}_{\mathbb{D},
5}^{\sharp,\Yright}(n,k,z_{o}) \me^{-2((n-1)K+k) \min \left\lbrace 
\min\limits_{j=1,2,\dotsc,N+1} \lbrace \tilde{\lambda}_{\tilde{
\mathcal{R}},4}^{\sharp,\smallfrown}(j) \rbrace,\min\limits_{j=1,2,
\dotsc,N+1} \lbrace \tilde{\lambda}_{\tilde{\mathcal{R}},4}^{\sharp,
\smallsmile}(j) \rbrace \right\rbrace} \right),
\end{align*}
and
\begin{align*}
\tilde{\mathbb{J}}_{\mathbb{D},5}^{\sharp,\Yleft} :=& \, \sum_{i=
1}^{N+1} \int_{\partial \tilde{\mathbb{U}}_{\tilde{\delta}_{\tilde{a}_{i}}}} 
\left(\lim_{-\tilde{\Sigma}_{\tilde{\mathcal{R}}}^{\sharp} \ni z^{\prime} 
\to \xi} \sum_{j=1}^{N+1} \int_{\tilde{\Sigma}_{p,j}^{4}} \sum_{i_{1},
i_{2}=1,2} \left\lvert \dfrac{(w_{+}^{\Sigma_{\tilde{\mathcal{R}}}}
(u))_{i_{1}i_{2}}}{u \! - \! z^{\prime}} \right\rvert^{2} \, \dfrac{\lvert 
\md u \rvert}{2 \pi^{2}} \right) \lvert \md \xi \rvert \\
\underset{\underset{z_{o}=1+o(1)}{\mathscr{N},n \to \infty}}{
\leqslant}& \, \mathcal{O} \left(\me^{-2((n-1)K+k) \min \left\lbrace 
\min\limits_{j=1,2,\dotsc,N+1} \lbrace \tilde{\lambda}_{\tilde{\mathcal{
R}},4}^{\sharp,\smallfrown}(j) \rbrace,\min\limits_{j=1,2,\dotsc,N+1} 
\lbrace \tilde{\lambda}_{\tilde{\mathcal{R}},4}^{\sharp,\smallsmile}(j) 
\rbrace \right\rbrace} \right) \underbrace{\sum_{i=1}^{N+1} \sum_{j=
1}^{N+1} \tfrac{2 \pi \tilde{\delta}_{\tilde{a}_{i}}(B^{2}(j)+\tilde{
\eta}_{j}^{2} \Theta_{M}^{2}(j))^{1/2}(\pi -\theta_{0}^{\tilde{a}}(j)-\theta_{
0}^{\tilde{b}}(j))}{(A(j)-B(j) \cos \theta_{0}^{\tilde{b}}(j)-\tilde{a}_{i}-
\tilde{\delta}_{\tilde{a}_{i}})^{2}+(\tilde{\eta}_{j}+\tilde{\delta}_{
\tilde{a}_{i}})^{2}}}_{= \, \mathcal{O}(1)} \\
\underset{\underset{z_{o}=1+o(1)}{\mathscr{N},n \to \infty}}{
\leqslant}& \, \mathcal{O} \left(\tilde{\mathfrak{c}}_{\mathbb{D},
5}^{\sharp,\Yleft}(n,k,z_{o}) \me^{-2((n-1)K+k) \min \left\lbrace 
\min\limits_{j=1,2,\dotsc,N+1} \lbrace \tilde{\lambda}_{\tilde{
\mathcal{R}},4}^{\sharp,\smallfrown}(j) \rbrace,\min\limits_{j=1,2,
\dotsc,N+1} \lbrace \tilde{\lambda}_{\tilde{\mathcal{R}},4}^{\sharp,
\smallsmile}(j) \rbrace \right\rbrace} \right),
\end{align*}
where $\tilde{\mathfrak{c}}_{\mathbb{D},5}^{\sharp,r}(n,k,z_{o}) \! 
=_{\underset{z_{o}=1+o(1)}{\mathscr{N},n \to \infty}} \! \mathcal{O}(1)$, 
$r \! \in \! \lbrace \Ydown,\Yup,\Yright,\Yleft \rbrace$, whence
\begin{equation} \label{eqiy61} 
\tilde{\mathbb{J}}_{\mathbb{D},5}^{\sharp} \underset{\underset{z_{o}
=1+o(1)}{\mathscr{N},n \to \infty}}{\leqslant} \mathcal{O} \left(
\tilde{\mathfrak{c}}_{\mathbb{D},5}^{\sharp}(n,k,z_{o}) \me^{-2((n-1)
K+k) \min \left\lbrace \min\limits_{j=1,2,\dotsc,N+1} \lbrace \tilde{
\lambda}_{\tilde{\mathcal{R}},4}^{\sharp,\smallfrown}(j) \rbrace,
\min\limits_{j=1,2,\dotsc,N+1} \lbrace \tilde{\lambda}_{\tilde{
\mathcal{R}},4}^{\sharp,\smallsmile}(j) \rbrace \right\rbrace} \right),
\end{equation}
where $\tilde{\mathfrak{c}}_{\mathbb{D},5}^{\sharp}(n,k,z_{o}) \! 
=_{\underset{z_{o}=1+o(1)}{\mathscr{N},n \to \infty}} \! \mathcal{O}
(1)$. Thus, via Equation~\eqref{eqiy45}, and the Estimates~\eqref{eqiy51}, 
\eqref{eqiy52}, \eqref{eqiy53}, \eqref{eqiy59}, and~\eqref{eqiy61}, one 
arrives at, for $n \! \in \! \mathbb{N}$ and $k \! \in \! \lbrace 1,2,\dotsc,
K \rbrace$ such that $\alpha_{p_{\mathfrak{s}}} \! := \! \alpha_{k} \! 
\neq \! \infty$,
\begin{equation} \label{eqiy62} 
\tilde{\mathbb{J}}_{\mathbb{D}}^{\sharp} \underset{\underset{z_{o}
=1+o(1)}{\mathscr{N},n \to \infty}}{\leqslant} \mathcal{O} \left(
\tilde{\mathfrak{c}}_{\mathbb{D}}^{\sharp}(n,k,z_{o}) \me^{-2((n-1)
K+k) \min \left\lbrace \min\limits_{j=1,2,\dotsc,N+1} \lbrace \tilde{
\lambda}_{\tilde{\mathcal{R}},4}^{\sharp,\smallfrown}(j) \rbrace,
\min\limits_{j=1,2,\dotsc,N+1} \lbrace \tilde{\lambda}_{\tilde{
\mathcal{R}},4}^{\sharp,\smallsmile}(j) \rbrace \right\rbrace} \right),
\end{equation}
where $\tilde{\mathfrak{c}}_{\mathbb{D}}^{\sharp}(n,k,z_{o}) \! =_{
\underset{z_{o}=1+o(1)}{\mathscr{N},n \to \infty}} \! \mathcal{O}
(1)$. Proceeding analogously as in the calculations leading to the 
Estimate~\eqref{eqiy62}, one shows, after several gruelling and 
tedious integration arguments, that, for $n \! \in \! \mathbb{N}$ 
and $k \! \in \! \lbrace 1,2,\dotsc,K \rbrace$ such that $\alpha_{
p_{\mathfrak{s}}} \! := \! \alpha_{k} \! \neq \! \infty$:\footnote{In 
the paradigm of the estimations herewith, definitions, intermediate 
formulae, and final estimates are given provided that the associated 
intermediate formulae are not unwieldy; however, in those cases 
where the intermediate formulae are inordinately unwieldy, 
only associated final estimates are presented.} (i) (cf. 
Definition~\eqref{eqiy40})
\begin{equation} \label{eqiy63} 
\tilde{\mathbb{J}}_{\mathbb{A}}^{\sharp} \! = \! \tilde{\mathbb{J}}_{
\mathbb{A},1}^{\sharp} \! + \! \tilde{\mathbb{J}}_{\mathbb{A},2}^{
\sharp} \! + \! \tilde{\mathbb{J}}_{\mathbb{A},3}^{\sharp} \! + \! 
\tilde{\mathbb{J}}_{\mathbb{A},4}^{\sharp} \! + \! \tilde{\mathbb{J}}_{
\mathbb{A},5}^{\sharp},
\end{equation}
where
\begin{align*}
\tilde{\mathbb{J}}_{\mathbb{A},1}^{\sharp} :=& \, \int_{-\infty}^{
\tilde{b}_{0}-\tilde{\delta}_{\tilde{b}_{0}}} \left(\lim_{-\tilde{\Sigma}_{
\tilde{\mathcal{R}}}^{\sharp} \ni z^{\prime} \to \xi} \sum_{j=1}^{N} 
\int_{\tilde{a}_{j}+\tilde{\delta}_{\tilde{a}_{j}}}^{\tilde{b}_{j}-\tilde{
\delta}_{\tilde{b}_{j}}} \sum_{i_{1},i_{2}=1,2} \left\lvert \dfrac{(w_{+}^{
\Sigma_{\tilde{\mathcal{R}}}}(u))_{i_{1}i_{2}}}{u \! - \! z^{\prime}} 
\right\rvert^{2} \, \dfrac{\lvert \md u \rvert}{2 \pi^{2}} \right) 
\lvert \md \xi \rvert \underset{\underset{z_{o}=1+o(1)}{\mathscr{N},
n \to \infty}}{\leqslant} \mathcal{O} \left(\me^{-2((n-1)K+k) 
\min\limits_{j=1,2,\dotsc,N} \lbrace \tilde{\lambda}_{\tilde{
\mathcal{R}},1}^{\sharp}(j) \rbrace} \right) \\
\times& \, \underbrace{\sum_{j=1}^{N} \ln \left\lvert \dfrac{
\tilde{b}_{j} \! - \! \tilde{b}_{0} \! - \! \tilde{\delta}_{\tilde{b}_{j}} \! 
+ \! \tilde{\delta}_{\tilde{b}_{0}}}{\tilde{a}_{j} \! - \! \tilde{b}_{0} \! 
+ \! \tilde{\delta}_{\tilde{a}_{j}} \! + \! \tilde{\delta}_{\tilde{b}_{0}}} 
\right\rvert}_{= \, \mathcal{O}(1)} \underset{\underset{z_{o}=1+
o(1)}{\mathscr{N},n \to \infty}}{\leqslant} \mathcal{O} \left(\tilde{
\mathfrak{c}}_{\mathbb{A},1}^{\sharp}(n,k,z_{o}) \me^{-2((n-1)K
+k) \min\limits_{j=1,2,\dotsc,N} \lbrace \tilde{\lambda}_{\tilde{
\mathcal{R}},1}^{\sharp}(j) \rbrace} \right), \\
\tilde{\mathbb{J}}_{\mathbb{A},2}^{\sharp} :=& \, \int_{\tilde{a}_{
N+1}+\tilde{\delta}_{\tilde{a}_{N+1}}}^{+\infty} \left(\lim_{-\tilde{
\Sigma}_{\tilde{\mathcal{R}}}^{\sharp} \ni z^{\prime} \to \xi} \sum_{j
=1}^{N} \int_{\tilde{a}_{j}+\tilde{\delta}_{\tilde{a}_{j}}}^{\tilde{b}_{j}-
\tilde{\delta}_{\tilde{b}_{j}}} \sum_{i_{1},i_{2}=1,2} \left\lvert \dfrac{
(w_{+}^{\Sigma_{\tilde{\mathcal{R}}}}(u))_{i_{1}i_{2}}}{u \! - \! 
z^{\prime}} \right\rvert^{2} \, \dfrac{\lvert \md u \rvert}{2 \pi^{2}} 
\right) \lvert \md \xi \rvert \underset{\underset{z_{o}=1+o(1)}{
\mathscr{N},n \to \infty}}{\leqslant} \mathcal{O} \left(\me^{-2((n-1)
K+k) \min\limits_{j=1,2,\dotsc,N} \lbrace \tilde{\lambda}_{\tilde{
\mathcal{R}},1}^{\sharp}(j) \rbrace} \right) \\
\times& \, \underbrace{\sum_{j=1}^{N} \ln \left\lvert \dfrac{
\tilde{a}_{j} \! - \! \tilde{a}_{N+1} \! + \! \tilde{\delta}_{\tilde{a}_{j}} 
\! - \! \tilde{\delta}_{\tilde{a}_{N+1}}}{\tilde{b}_{j} \! - \! \tilde{a}_{N+1} 
\! - \! \tilde{\delta}_{\tilde{b}_{j}} \! - \! \tilde{\delta}_{\tilde{a}_{N+1}}} 
\right\rvert}_{= \, \mathcal{O}(1)} \underset{\underset{z_{o}=1+
o(1)}{\mathscr{N},n \to \infty}}{\leqslant} \mathcal{O} \left(\tilde{
\mathfrak{c}}_{\mathbb{A},2}^{\sharp}(n,k,z_{o}) \me^{-2((n-1)K
+k) \min\limits_{j=1,2,\dotsc,N} \lbrace \tilde{\lambda}_{\tilde{
\mathcal{R}},1}^{\sharp}(j) \rbrace} \right),
\end{align*}
\begin{align*}
\tilde{\mathbb{J}}_{\mathbb{A},3}^{\sharp} :=& \, \sum_{i=1}^{N} 
\int_{\tilde{a}_{i}+\tilde{\delta}_{\tilde{a}_{i}}}^{\tilde{b}_{i}-\tilde{
\delta}_{\tilde{b}_{i}}} \left(\lim_{-\tilde{\Sigma}_{\tilde{\mathcal{
R}}}^{\sharp} \ni z^{\prime} \to \xi} \sum_{j=1}^{N} \int_{\tilde{a}_{j}
+\tilde{\delta}_{\tilde{a}_{j}}}^{\tilde{b}_{j}-\tilde{\delta}_{\tilde{b}_{j}}} 
\sum_{i_{1},i_{2}=1,2} \left\lvert \dfrac{(w_{+}^{\Sigma_{\tilde{
\mathcal{R}}}}(u))_{i_{1}i_{2}}}{u \! - \! z^{\prime}} \right\rvert^{2} \, 
\dfrac{\lvert \md u \rvert}{2 \pi^{2}} \right) \lvert \md \xi \rvert \\
\underset{\underset{z_{o}=1+o(1)}{\mathscr{N},n \to \infty}}{
\leqslant}& \, \mathcal{O} \left(\tilde{\mathfrak{c}}_{\mathbb{A},3}^{
\sharp}(n,k,z_{o}) \me^{-((n-1)K+k) \min\limits_{j=1,2,\dotsc,N} 
\lbrace \tilde{\lambda}_{\tilde{\mathcal{R}},1}^{\sharp}(j) \rbrace} 
\right),
\end{align*}
\begin{align*}
\tilde{\mathbb{J}}_{\mathbb{A},4}^{\sharp} :=& \, \sum_{i=1}^{N+1} 
\sum_{m=3,4} \int_{\tilde{\Sigma}_{p,i}^{m}} \left(\lim_{-\tilde{
\Sigma}_{\tilde{\mathcal{R}}}^{\sharp} \ni z^{\prime} \to \xi} \sum_{j
=1}^{N} \int_{\tilde{a}_{j}+\tilde{\delta}_{\tilde{a}_{j}}}^{\tilde{b}_{j}-
\tilde{\delta}_{\tilde{b}_{j}}} \sum_{i_{1},i_{2}=1,2} \left\lvert \dfrac{
(w_{+}^{\Sigma_{\tilde{\mathcal{R}}}}(u))_{i_{1}i_{2}}}{u \! - \! 
z^{\prime}} \right\rvert^{2} \, \dfrac{\lvert \md u \rvert}{2 \pi^{2}} 
\right) \lvert \md \xi \rvert \\
\underset{\underset{z_{o}=1+o(1)}{\mathscr{N},n \to \infty}}{
\leqslant}& \, \mathcal{O} \left(\me^{-2((n-1)K+k) \min\limits_{j=
1,2,\dotsc,N} \lbrace \tilde{\lambda}_{\tilde{\mathcal{R}},1}^{\sharp}
(j) \rbrace} \right) \underbrace{\sum_{i=1}^{N+1} \sum_{j=1}^{N} 
\dfrac{(B^{2}(i) \! + \! \tilde{\eta}_{i}^{2} \Theta_{M}^{2}(i))^{1/2}}{\tilde{
\eta}_{i}} \ln \left(\left\lvert \dfrac{(1 \! + \! \cos \theta_{0}^{\tilde{b}}
(i)) \sin \theta_{0}^{\tilde{a}}(i)}{(1 \! - \! \cos \theta_{0}^{\tilde{a}}
(i)) \sin \theta_{0}^{\tilde{b}}(i)} \right\rvert \right)}_{= \, \mathcal{O}
(1)} \\
\times& \, \left(\underbrace{\tan^{-1} \left(\dfrac{\tilde{b}_{j} \! - \! 
\tilde{\delta}_{\tilde{b}_{j}} \! - \! A(i) \! + \! B(i) \Theta_{M}(i)}{\tilde{
\eta}_{i} \Theta_{m}(i)} \right) \! - \! \tan^{-1} \left(\dfrac{\tilde{a}_{j} 
\! + \! \tilde{\delta}_{\tilde{a}_{j}} \! - \! A(i) \! - \! B(i) \Theta_{M}(i)}{
\tilde{\eta}_{i}} \right)}_{= \, \mathcal{O}(1)} \right. \\
+&\left. \, \underbrace{\tan^{-1} \left(\dfrac{\tilde{b}_{j} \! - \! 
\tilde{\delta}_{\tilde{b}_{j}} \! - \! A(i) \! + \! B(i) \Theta_{M}(i)}{\tilde{
\eta}_{i}} \right) \! - \! \tan^{-1} \left(\dfrac{\tilde{a}_{j} \! + \! \tilde{
\delta}_{\tilde{a}_{j}} \! - \! A(i) \! - \! B(i) \Theta_{M}(i)}{\tilde{\eta}_{i} 
\Theta_{m}(i)} \right)}_{= \, \mathcal{O}(1)} \right) \\
\underset{\underset{z_{o}=1+o(1)}{\mathscr{N},n \to \infty}}{
\leqslant}& \, \mathcal{O} \left(\tilde{\mathfrak{c}}_{\mathbb{A},4}^{
\sharp}(n,k,z_{o}) \me^{-2((n-1)K+k) \min\limits_{j=1,2,\dotsc,N} 
\lbrace \tilde{\lambda}_{\tilde{\mathcal{R}},1}^{\sharp}(j) \rbrace} 
\right),
\end{align*}
\begin{align*}
\tilde{\mathbb{J}}_{\mathbb{A},5}^{\sharp} :=& \, \sum_{i=1}^{N+1} 
\left(\int_{\partial \tilde{\mathbb{U}}_{\tilde{\delta}_{\tilde{b}_{i-1}}}} 
\! + \! \int_{\partial \tilde{\mathbb{U}}_{\tilde{\delta}_{\tilde{a}_{i}}}} 
\right) \left(\lim_{-\tilde{\Sigma}_{\tilde{\mathcal{R}}}^{\sharp} \ni 
z^{\prime} \to \xi} \sum_{j=1}^{N} \int_{\tilde{a}_{j}+\tilde{\delta}_{
\tilde{a}_{j}}}^{\tilde{b}_{j}-\tilde{\delta}_{\tilde{b}_{j}}} \sum_{i_{1},i_{2}
=1,2} \left\lvert \dfrac{(w_{+}^{\Sigma_{\tilde{\mathcal{R}}}}(u))_{i_{1}
i_{2}}}{u \! - \! z^{\prime}} \right\rvert^{2} \, \dfrac{\lvert \md u 
\rvert}{2 \pi^{2}} \right) \lvert \md \xi \rvert \\
\underset{\underset{z_{o}=1+o(1)}{\mathscr{N},n \to \infty}}{
\leqslant}& \, \mathcal{O} \left(\tilde{\mathfrak{c}}_{\mathbb{A},5}^{
\sharp}(n,k,z_{o}) \me^{-2((n-1)K+k) \min\limits_{j=1,2,\dotsc,N} 
\lbrace \tilde{\lambda}_{\tilde{\mathcal{R}},1}^{\sharp}(j) \rbrace} 
\right),
\end{align*}
where $\tilde{\mathfrak{c}}_{\mathbb{A},r}^{\sharp}(n,k,z_{o}) \! =_{
\underset{z_{o}=1+o(1)}{\mathscr{N},n \to \infty}} \! \mathcal{O}(1)$, 
$r \! \in \! \lbrace 1,2,3,4,5 \rbrace$, whence
\begin{equation} \label{eqiy64}
\tilde{\mathbb{J}}_{\mathbb{A}}^{\sharp} \underset{\underset{z_{o}=
1+o(1)}{\mathscr{N},n \to \infty}}{\leqslant} \mathcal{O} \left(\tilde{
\mathfrak{c}}_{\mathbb{A}}^{\sharp}(n,k,z_{o}) \me^{-((n-1)K+k) 
\min\limits_{j=1,2,\dotsc,N} \lbrace \tilde{\lambda}_{\tilde{
\mathcal{R}},1}^{\sharp}(j) \rbrace} \right),
\end{equation}
where $\tilde{\mathfrak{c}}_{\mathbb{A}}^{\sharp}(n,k,z_{o}) \! =_{
\underset{z_{o}=1+o(1)}{\mathscr{N},n \to \infty}} \! \mathcal{O}(1)$; 
(ii) (cf. Definition~\eqref{eqiy41})
\begin{equation} \label{eqiy65} 
\tilde{\mathbb{J}}_{\mathbb{B}}^{\sharp} \! = \! \tilde{\mathbb{J}}_{
\mathbb{B},1}^{\sharp} \! + \! \tilde{\mathbb{J}}_{\mathbb{B},2}^{
\sharp} \! + \! \tilde{\mathbb{J}}_{\mathbb{B},3}^{\sharp} \! + \! 
\tilde{\mathbb{J}}_{\mathbb{B},4}^{\sharp} \! + \! \tilde{\mathbb{J}}_{
\mathbb{B},5}^{\sharp},
\end{equation}
where
\begin{align*}
\tilde{\mathbb{J}}_{\mathbb{B},1}^{\sharp} :=& \, \int_{-\infty}^{
\tilde{b}_{0}-\tilde{\delta}_{\tilde{b}_{0}}} \left(\lim_{-\tilde{\Sigma}_{
\tilde{\mathcal{R}}}^{\sharp} \ni z^{\prime} \to \xi} \int_{\tilde{a}_{N
+1}+\tilde{\delta}_{\tilde{a}_{N+1}}}^{+\infty} \sum_{i_{1},i_{2}=1,2} 
\left\lvert \dfrac{(w_{+}^{\Sigma_{\tilde{\mathcal{R}}}}(u))_{i_{1}
i_{2}}}{u \! - \! z^{\prime}} \right\rvert^{2} \, \dfrac{\lvert \md u 
\rvert}{2 \pi^{2}} \right) \lvert \md \xi \rvert \\
\underset{\underset{z_{o}=1+o(1)}{\mathscr{N},n \to \infty}}{
\leqslant}& \, \mathcal{O} \left(\tilde{\mathfrak{c}}_{\mathbb{B},1}^{
\sharp}(n,k,z_{o}) \me^{-((n-1)K+k) \tilde{\lambda}_{\tilde{\mathcal{
R}},2}^{\sharp}(+)} \right), \\
\tilde{\mathbb{J}}_{\mathbb{B},2}^{\sharp} :=& \, \int_{\tilde{a}_{N+
1}+\tilde{\delta}_{\tilde{a}_{N+1}}}^{+\infty} \left(\lim_{-\tilde{\Sigma}_{
\tilde{\mathcal{R}}}^{\sharp} \ni z^{\prime} \to \xi} \int_{\tilde{a}_{N
+1}+\tilde{\delta}_{\tilde{a}_{N+1}}}^{+\infty} \sum_{i_{1},i_{2}=1,2} 
\left\lvert \dfrac{(w_{+}^{\Sigma_{\tilde{\mathcal{R}}}}(u))_{i_{1}
i_{2}}}{u \! - \! z^{\prime}} \right\rvert^{2} \, \dfrac{\lvert \md u 
\rvert}{2 \pi^{2}} \right) \lvert \md \xi \rvert \\
\underset{\underset{z_{o}=1+o(1)}{\mathscr{N},n \to \infty}}{
\leqslant}& \, \mathcal{O} \left(\tilde{\mathfrak{c}}_{\mathbb{B},2}^{
\sharp}(n,k,z_{o}) \me^{-((n-1)K+k) \tilde{\lambda}_{\tilde{\mathcal{
R}},2}^{\sharp}(+)} \right),
\end{align*}
\begin{align*}
\tilde{\mathbb{J}}_{\mathbb{B},3}^{\sharp} :=& \, \sum_{j=1}^{N} 
\int_{\tilde{a}_{j}+\tilde{\delta}_{\tilde{a}_{j}}}^{\tilde{b}_{j}-\tilde{
\delta}_{\tilde{b}_{j}}} \left(\lim_{-\tilde{\Sigma}_{\tilde{\mathcal{
R}}}^{\sharp} \ni z^{\prime} \to \xi} \int_{\tilde{a}_{N+1}+\tilde{
\delta}_{\tilde{a}_{N+1}}}^{+\infty} \sum_{i_{1},i_{2}=1,2} \left\lvert 
\dfrac{(w_{+}^{\Sigma_{\tilde{\mathcal{R}}}}(u))_{i_{1}i_{2}}}{u \! - \! 
z^{\prime}} \right\rvert^{2} \, \dfrac{\lvert \md u \rvert}{2 \pi^{2}} 
\right) \lvert \md \xi \rvert \underset{\underset{z_{o}=1+o(1)}{
\mathscr{N},n \to \infty}}{\leqslant} \mathcal{O} \left(\me^{-2((n-1)
K+k) \tilde{\lambda}_{\tilde{\mathcal{R}},2}^{\sharp}(+)} \right) \\
\times& \, \underbrace{\sum_{j=1}^{N} \ln \left\lvert \dfrac{
\tilde{a}_{N+1} \! - \! \tilde{a}_{j} \! + \! \tilde{\delta}_{\tilde{a}_{N+
1}} \! - \! \tilde{\delta}_{\tilde{a}_{j}}}{\tilde{a}_{N+1} \! - \! \tilde{b}_{j} 
\! + \! \tilde{\delta}_{\tilde{a}_{N+1}} \! + \! \tilde{\delta}_{\tilde{b}_{j}}} 
\right\rvert}_{= \, \mathcal{O}(1)} \underset{\underset{z_{o}=1+
o(1)}{\mathscr{N},n \to \infty}}{\leqslant} \mathcal{O} \left(\tilde{
\mathfrak{c}}_{\mathbb{B},3}^{\sharp}(n,k,z_{o}) \me^{-2((n-1)K
+k) \tilde{\lambda}_{\tilde{\mathcal{R}},2}^{\sharp}(+)} \right),
\end{align*}
\begin{align*}
\tilde{\mathbb{J}}_{\mathbb{B},4}^{\sharp} :=& \, \sum_{j=1}^{N+
1} \sum_{m=3,4} \int_{\tilde{\Sigma}_{p,j}^{m}} \left(\lim_{-\tilde{
\Sigma}_{\tilde{\mathcal{R}}}^{\sharp} \ni z^{\prime} \to \xi} 
\int_{\tilde{a}_{N+1}+\tilde{\delta}_{\tilde{a}_{N+1}}}^{+\infty} 
\sum_{i_{1},i_{2}=1,2} \left\lvert \dfrac{(w_{+}^{\Sigma_{\tilde{
\mathcal{R}}}}(u))_{i_{1}i_{2}}}{u \! - \! z^{\prime}} \right\rvert^{2} \, 
\dfrac{\lvert \md u \rvert}{2 \pi^{2}} \right) \lvert \md \xi \rvert 
\underset{\underset{z_{o}=1+o(1)}{\mathscr{N},n \to \infty}}{
\leqslant} \mathcal{O} \left(\me^{-2((n-1)K+k) \tilde{\lambda}_{
\tilde{\mathcal{R}},2}^{\sharp}(+)} \right) \\
\times& \, \underbrace{\sum_{j=1}^{N+1} \dfrac{(B^{2}(j) \! + \! 
\tilde{\eta}_{j}^{2} \Theta_{M}^{2}(j))^{1/2}}{\tilde{\eta}_{j}} \ln \left(
\left\lvert \dfrac{(1 \! - \! \cos \theta_{0}^{\tilde{a}}(j)) \sin \theta_{0}^{
\tilde{b}}(j)}{(1 \! + \! \cos \theta_{0}^{\tilde{b}}(j)) \sin \theta_{
0}^{\tilde{a}}(j)} \right\rvert \right)}_{= \, \mathcal{O}(1)} \left(
\underbrace{\tan^{-1} \left(\dfrac{\tilde{a}_{N+1} \! + \! \tilde{
\delta}_{\tilde{a}_{N+1}} \! - \! A(j) \! - \! B(j) \Theta_{M}(j)}{\tilde{
\eta}_{j}} \right)}_{= \, \mathcal{O}(1)} \right. \\
+&\left. \, \underbrace{\tan^{-1} \left(\dfrac{\tilde{a}_{N+1} \! + \! 
\tilde{\delta}_{\tilde{a}_{N+1}} \! - \! A(j) \! - \! B(j) \Theta_{M}(j)}{
\tilde{\eta}_{j} \Theta_{m}(j)} \right)}_{= \, \mathcal{O}(1)} \right) 
\underset{\underset{z_{o}=1+o(1)}{\mathscr{N},n \to \infty}}{
\leqslant} \mathcal{O} \left(\tilde{\mathfrak{c}}_{\mathbb{B},4}^{
\sharp}(n,k,z_{o}) \me^{-2((n-1)K+k) \tilde{\lambda}_{\tilde{
\mathcal{R}},2}^{\sharp}(+)} \right),
\end{align*}
\begin{align*}
\tilde{\mathbb{J}}_{\mathbb{B},5}^{\sharp} :=& \, \sum_{j=1}^{N+1} 
\left(\int_{\partial \tilde{\mathbb{U}}_{\tilde{\delta}_{\tilde{b}_{j-1}}}} 
\! + \! \int_{\partial \tilde{\mathbb{U}}_{\tilde{\delta}_{\tilde{a}_{j}}}} 
\right) \left(\lim_{-\tilde{\Sigma}_{\tilde{\mathcal{R}}}^{\sharp} \ni 
z^{\prime} \to \xi} \int_{\tilde{a}_{N+1}+\tilde{\delta}_{\tilde{a}_{N+
1}}}^{+\infty} \sum_{i_{1},i_{2}=1,2} \left\lvert \dfrac{(w_{+}^{
\Sigma_{\tilde{\mathcal{R}}}}(u))_{i_{1}i_{2}}}{u \! - \! z^{\prime}} 
\right\rvert^{2} \, \dfrac{\lvert \md u \rvert}{2 \pi^{2}} \right) \lvert 
\md \xi \rvert \\
\underset{\underset{z_{o}=1+o(1)}{\mathscr{N},n \to \infty}}{
\leqslant}& \, \mathcal{O} \left(\tilde{\mathfrak{c}}_{\mathbb{B},
5}^{\sharp}(n,k,z_{o}) \me^{-((n-1)K+k) \tilde{\lambda}_{\tilde{
\mathcal{R}},2}^{\sharp}(+)} \right),
\end{align*}
where $\tilde{\mathfrak{c}}_{\mathbb{B},r}^{\sharp}(n,k,z_{o}) \! =_{
\underset{z_{o}=1+o(1)}{\mathscr{N},n \to \infty}} \! \mathcal{O}(1)$, 
$r \! \in \! \lbrace 1,2,3,4,5 \rbrace$, whence
\begin{equation} \label{eqiy66}
\tilde{\mathbb{J}}_{\mathbb{B}}^{\sharp} \underset{\underset{z_{o}=
1+o(1)}{\mathscr{N},n \to \infty}}{\leqslant} \mathcal{O} \left(\tilde{
\mathfrak{c}}_{\mathbb{B}}^{\sharp}(n,k,z_{o}) \me^{-((n-1)K+k) 
\tilde{\lambda}_{\tilde{\mathcal{R}},2}^{\sharp}(+)} \right),
\end{equation}
where $\tilde{\mathfrak{c}}_{\mathbb{B}}^{\sharp}(n,k,z_{o}) \! =_{
\underset{z_{o}=1+o(1)}{\mathscr{N},n \to \infty}} \! \mathcal{O}(1)$; 
(iii) (cf. Definition~\eqref{eqiy42})
\begin{equation} \label{eqiy67} 
\tilde{\mathbb{J}}_{\mathbb{C}}^{\sharp} \! = \! \tilde{\mathbb{J}}_{
\mathbb{C},1}^{\sharp} \! + \! \tilde{\mathbb{J}}_{\mathbb{C},2}^{
\sharp} \! + \! \tilde{\mathbb{J}}_{\mathbb{C},3}^{\sharp} \! + \! 
\tilde{\mathbb{J}}_{\mathbb{C},4}^{\sharp} \! + \! \tilde{\mathbb{J}}_{
\mathbb{C},5}^{\sharp},
\end{equation}
where
\begin{align*}
\tilde{\mathbb{J}}_{\mathbb{C},1}^{\sharp} :=& \, \int_{-\infty}^{
\tilde{b}_{0}-\tilde{\delta}_{\tilde{b}_{0}}} \left(\lim_{-\tilde{\Sigma}_{
\tilde{\mathcal{R}}}^{\sharp} \ni z^{\prime} \to \xi} \int_{-\infty}^{
\tilde{b}_{0}-\tilde{\delta}_{\tilde{b}_{0}}} \sum_{i_{1},i_{2}=1,2} 
\left\lvert \dfrac{(w_{+}^{\Sigma_{\tilde{\mathcal{R}}}}(u))_{i_{1}
i_{2}}}{u \! - \! z^{\prime}} \right\rvert^{2} \, \dfrac{\lvert \md u 
\rvert}{2 \pi^{2}} \right) \lvert \md \xi \rvert \\
\underset{\underset{z_{o}=1+o(1)}{\mathscr{N},n \to \infty}}{
\leqslant}& \, \mathcal{O} \left(\tilde{\mathfrak{c}}_{\mathbb{C},1}^{
\sharp}(n,k,z_{o}) \me^{-((n-1)K+k) \tilde{\lambda}_{\tilde{\mathcal{
R}},2}^{\sharp}(-)} \right), \\
\tilde{\mathbb{J}}_{\mathbb{C},2}^{\sharp} :=& \, \int_{\tilde{a}_{N+
1}+\tilde{\delta}_{\tilde{a}_{N+1}}}^{+\infty} \left(\lim_{-\tilde{\Sigma}_{
\tilde{\mathcal{R}}}^{\sharp} \ni z^{\prime} \to \xi} \int_{-\infty}^{
\tilde{b}_{0}-\tilde{\delta}_{\tilde{b}_{0}}} \sum_{i_{1},i_{2}=1,2} 
\left\lvert \dfrac{(w_{+}^{\Sigma_{\tilde{\mathcal{R}}}}(u))_{i_{1}
i_{2}}}{u \! - \! z^{\prime}} \right\rvert^{2} \, \dfrac{\lvert \md u 
\rvert}{2 \pi^{2}} \right) \lvert \md \xi \rvert \\
\underset{\underset{z_{o}=1+o(1)}{\mathscr{N},n \to \infty}}{
\leqslant}& \, \mathcal{O} \left(\tilde{\mathfrak{c}}_{\mathbb{C},2}^{
\sharp}(n,k,z_{o}) \me^{-((n-1)K+k) \tilde{\lambda}_{\tilde{\mathcal{
R}},2}^{\sharp}(-)} \right),
\end{align*}
\begin{align*}
\tilde{\mathbb{J}}_{\mathbb{C},3}^{\sharp} :=& \, \sum_{j=1}^{N} 
\int_{\tilde{a}_{j}+\tilde{\delta}_{\tilde{a}_{j}}}^{\tilde{b}_{j}-\tilde{
\delta}_{\tilde{b}_{j}}} \left(\lim_{-\tilde{\Sigma}_{\tilde{\mathcal{
R}}}^{\sharp} \ni z^{\prime} \to \xi} \int_{-\infty}^{\tilde{b}_{0}-
\tilde{\delta}_{\tilde{b}_{0}}} \sum_{i_{1},i_{2}=1,2} \left\lvert 
\dfrac{(w_{+}^{\Sigma_{\tilde{\mathcal{R}}}}(u))_{i_{1}i_{2}}}{u \! - \! 
z^{\prime}} \right\rvert^{2} \, \dfrac{\lvert \md u \rvert}{2 \pi^{2}} 
\right) \lvert \md \xi \rvert \underset{\underset{z_{o}=1+o(1)}{
\mathscr{N},n \to \infty}}{\leqslant} \mathcal{O} \left(\me^{-2((n-1)
K+k) \tilde{\lambda}_{\tilde{\mathcal{R}},2}^{\sharp}(-)} \right) \\
\times& \, \underbrace{\sum_{j=1}^{N} \ln \left\lvert \dfrac{
\tilde{b}_{0} \! - \! \tilde{b}_{j} \! - \! \tilde{\delta}_{\tilde{b}_{0}} 
\! + \! \tilde{\delta}_{\tilde{b}_{j}}}{\tilde{b}_{0} \! - \! \tilde{a}_{j} \! 
- \! \tilde{\delta}_{\tilde{b}_{0}} \! - \! \tilde{\delta}_{\tilde{a}_{j}}} 
\right\rvert}_{= \, \mathcal{O}(1)} \underset{\underset{z_{o}=1+
o(1)}{\mathscr{N},n \to \infty}}{\leqslant} \mathcal{O} \left(\tilde{
\mathfrak{c}}_{\mathbb{C},3}^{\sharp}(n,k,z_{o}) \me^{-2((n-1)K
+k) \tilde{\lambda}_{\tilde{\mathcal{R}},2}^{\sharp}(-)} \right),
\end{align*}
\begin{align*}
\tilde{\mathbb{J}}_{\mathbb{C},4}^{\sharp} :=& \, \sum_{j=1}^{N+
1} \sum_{m=3,4} \int_{\tilde{\Sigma}_{p,j}^{m}} \left(\lim_{-\tilde{
\Sigma}_{\tilde{\mathcal{R}}}^{\sharp} \ni z^{\prime} \to \xi} 
\int_{-\infty}^{\tilde{b}_{0}-\tilde{\delta}_{\tilde{b}_{0}}} \sum_{i_{1},
i_{2}=1,2} \left\lvert \dfrac{(w_{+}^{\Sigma_{\tilde{\mathcal{R}}}}
(u))_{i_{1}i_{2}}}{u \! - \! z^{\prime}} \right\rvert^{2} \, \dfrac{\lvert 
\md u \rvert}{2 \pi^{2}} \right) \lvert \md \xi \rvert \underset{
\underset{z_{o}=1+o(1)}{\mathscr{N},n \to \infty}}{\leqslant} 
\mathcal{O} \left(\me^{-2((n-1)K+k) \tilde{\lambda}_{\tilde{
\mathcal{R}},2}^{\sharp}(-)} \right) \\
\times& \, \underbrace{\sum_{j=1}^{N+1} \dfrac{(B^{2}(j) \! + \! 
\tilde{\eta}_{j}^{2} \Theta_{M}^{2}(j))^{1/2}}{\tilde{\eta}_{j}} \ln \left(
\left\lvert \dfrac{(1 \! + \! \cos \theta_{0}^{\tilde{b}}(j)) \sin \theta_{0}^{
\tilde{a}}(j)}{(1 \! - \! \cos \theta_{0}^{\tilde{a}}(j)) \sin \theta_{
0}^{\tilde{b}}(j)} \right\rvert \right)}_{= \, \mathcal{O}(1)} \left(
\underbrace{\tan^{-1} \left(\dfrac{\tilde{b}_{0} \! - \! \tilde{
\delta}_{\tilde{b}_{0}} \! - \! A(j) \! + \! B(j) \Theta_{M}(j)}{\tilde{
\eta}_{j} \Theta_{m}(j)} \right)}_{= \, \mathcal{O}(1)} \right. \\
+&\left. \, \underbrace{\tan^{-1} \left(\dfrac{\tilde{b}_{0} \! - \! 
\tilde{\delta}_{\tilde{b}_{0}} \! - \! A(j) \! + \! B(j) \Theta_{M}(j)}{
\tilde{\eta}_{j}} \right)}_{= \, \mathcal{O}(1)} \right) \underset{
\underset{z_{o}=1+o(1)}{\mathscr{N},n \to \infty}}{\leqslant} 
\mathcal{O} \left(\tilde{\mathfrak{c}}_{\mathbb{C},4}^{\sharp}
(n,k,z_{o}) \me^{-2((n-1)K+k) \tilde{\lambda}_{\tilde{\mathcal{R}},
2}^{\sharp}(-)} \right),
\end{align*}
\begin{align*}
\tilde{\mathbb{J}}_{\mathbb{C},5}^{\sharp} :=& \, \sum_{j=1}^{N+1} 
\left(\int_{\partial \tilde{\mathbb{U}}_{\tilde{\delta}_{\tilde{b}_{j-1}}}} 
\! + \! \int_{\partial \tilde{\mathbb{U}}_{\tilde{\delta}_{\tilde{a}_{j}}}} 
\right) \left(\lim_{-\tilde{\Sigma}_{\tilde{\mathcal{R}}}^{\sharp} \ni 
z^{\prime} \to \xi} \int_{-\infty}^{\tilde{b}_{0}-\tilde{\delta}_{\tilde{
b}_{0}}} \sum_{i_{1},i_{2}=1,2} \left\lvert \dfrac{(w_{+}^{\Sigma_{
\tilde{\mathcal{R}}}}(u))_{i_{1}i_{2}}}{u \! - \! z^{\prime}} \right\rvert^{2} 
\, \dfrac{\lvert \md u \rvert}{2 \pi^{2}} \right) \lvert \md \xi \rvert \\
\underset{\underset{z_{o}=1+o(1)}{\mathscr{N},n \to \infty}}{
\leqslant}& \, \mathcal{O} \left(\tilde{\mathfrak{c}}_{\mathbb{C},5}^{
\sharp}(n,k,z_{o}) \me^{-((n-1)K+k) \tilde{\lambda}_{\tilde{\mathcal{
R}},2}^{\sharp}(-)} \right),
\end{align*}
where $\tilde{\mathfrak{c}}_{\mathbb{C},r}^{\sharp}(n,k,z_{o}) \! =_{
\underset{z_{o}=1+o(1)}{\mathscr{N},n \to \infty}} \! \mathcal{O}(1)$, 
$r \! \in \! \lbrace 1,2,3,4,5 \rbrace$, whence
\begin{equation} \label{eqiy68}
\tilde{\mathbb{J}}_{\mathbb{C}}^{\sharp} \underset{\underset{z_{o}=
1+o(1)}{\mathscr{N},n \to \infty}}{\leqslant} \mathcal{O} \left(\tilde{
\mathfrak{c}}_{\mathbb{C}}^{\sharp}(n,k,z_{o}) \me^{-((n-1)K+k) 
\tilde{\lambda}_{\tilde{\mathcal{R}},2}^{\sharp}(-)} \right),
\end{equation}
where $\tilde{\mathfrak{c}}_{\mathbb{C}}^{\sharp}(n,k,z_{o}) \! =_{
\underset{z_{o}=1+o(1)}{\mathscr{N},n \to \infty}} \! \mathcal{O}(1)$; 
(iv) (cf. Definition~\eqref{eqiy44})
\begin{equation} \label{eqiy69} 
\tilde{\mathbb{J}}_{\mathbb{E}}^{\sharp} \! = \! \tilde{\mathbb{J}}_{
\mathbb{E},1}^{\sharp} \! + \! \tilde{\mathbb{J}}_{\mathbb{E},2}^{
\sharp} \! + \! \tilde{\mathbb{J}}_{\mathbb{E},3}^{\sharp} \! + \! 
\tilde{\mathbb{J}}_{\mathbb{E},4}^{\sharp} \! + \! \tilde{\mathbb{J}}_{
\mathbb{E},5}^{\sharp},
\end{equation}
where
\begin{align*}
\tilde{\mathbb{J}}_{\mathbb{E},1}^{\sharp} :=& \, \int_{-\infty}^{
\tilde{b}_{0}-\tilde{\delta}_{\tilde{b}_{0}}} \left(\lim_{-\tilde{\Sigma}_{
\tilde{\mathcal{R}}}^{\sharp} \ni z^{\prime} \to \xi} \sum_{j=1}^{N+1} 
\left(\int_{\partial \tilde{\mathbb{U}}_{\tilde{\delta}_{\tilde{b}_{j-1}}}} 
\! + \! \int_{\partial \tilde{\mathbb{U}}_{\tilde{\delta}_{\tilde{a}_{j}}}} 
\right) \sum_{i_{1},i_{2}=1,2} \left\lvert \dfrac{(w_{+}^{\Sigma_{
\tilde{\mathcal{R}}}}(u))_{i_{1}i_{2}}}{u \! - \! z^{\prime}} 
\right\rvert^{2} \, \dfrac{\lvert \md u \rvert}{2 \pi^{2}} \right) \lvert 
\md \xi \rvert \underset{\underset{z_{o}=1+o(1)}{\mathscr{N},n \to 
\infty}}{=} 0, \\
\tilde{\mathbb{J}}_{\mathbb{E},2}^{\sharp} :=& \, \int_{\tilde{a}_{N+
1}+\tilde{\delta}_{\tilde{a}_{N+1}}}^{+\infty} \left(\lim_{-\tilde{\Sigma}_{
\tilde{\mathcal{R}}}^{\sharp} \ni z^{\prime} \to \xi} \sum_{j=1}^{N+1} 
\left(\int_{\partial \tilde{\mathbb{U}}_{\tilde{\delta}_{\tilde{b}_{j-1}}}} 
\! + \! \int_{\partial \tilde{\mathbb{U}}_{\tilde{\delta}_{\tilde{a}_{j}}}} 
\right) \sum_{i_{1},i_{2}=1,2} \left\lvert \dfrac{(w_{+}^{\Sigma_{
\tilde{\mathcal{R}}}}(u))_{i_{1}i_{2}}}{u \! - \! z^{\prime}} \right\rvert^{2} 
\, \dfrac{\lvert \md u \rvert}{2 \pi^{2}} \right) \lvert \md \xi \rvert 
\underset{\underset{z_{o}=1+o(1)}{\mathscr{N},n \to \infty}}{=} 0, \\
\tilde{\mathbb{J}}_{\mathbb{E},3}^{\sharp} :=& \, \sum_{i=1}^{N} 
\int_{\tilde{a}_{i}+\tilde{\delta}_{\tilde{a}_{i}}}^{\tilde{b}_{i}-\tilde{
\delta}_{\tilde{b}_{i}}} \left(\lim_{-\tilde{\Sigma}_{\tilde{\mathcal{
R}}}^{\sharp} \ni z^{\prime} \to \xi} \sum_{j=1}^{N+1} \left(\int_{
\partial \tilde{\mathbb{U}}_{\tilde{\delta}_{\tilde{b}_{j-1}}}} \! + \! 
\int_{\partial \tilde{\mathbb{U}}_{\tilde{\delta}_{\tilde{a}_{j}}}} \right) 
\sum_{i_{1},i_{2}=1,2} \left\lvert \dfrac{(w_{+}^{\Sigma_{\tilde{
\mathcal{R}}}}(u))_{i_{1}i_{2}}}{u \! - \! z^{\prime}} \right\rvert^{2} \, 
\dfrac{\lvert \md u \rvert}{2 \pi^{2}} \right) \lvert \md \xi \rvert 
\underset{\underset{z_{o}=1+o(1)}{\mathscr{N},n \to \infty}}{=} 0, \\
\tilde{\mathbb{J}}_{\mathbb{E},4}^{\sharp} :=& \, \sum_{i=1}^{N+
1} \sum_{m=3,4} \int_{\tilde{\Sigma}_{p,i}^{m}} \left(\lim_{-\tilde{
\Sigma}_{\tilde{\mathcal{R}}}^{\sharp} \ni z^{\prime} \to \xi} \sum_{j
=1}^{N+1} \left(\int_{\partial \tilde{\mathbb{U}}_{\tilde{\delta}_{
\tilde{b}_{j-1}}}} \! + \! \int_{\partial \tilde{\mathbb{U}}_{\tilde{
\delta}_{\tilde{a}_{j}}}} \right) \sum_{i_{1},i_{2}=1,2} \left\lvert 
\dfrac{(w_{+}^{\Sigma_{\tilde{\mathcal{R}}}}(u))_{i_{1}i_{2}}}{u \! 
- \! z^{\prime}} \right\rvert^{2} \, \dfrac{\lvert \md u \rvert}{2 
\pi^{2}} \right) \lvert \md \xi \rvert \underset{\underset{z_{o}=
1+o(1)}{\mathscr{N},n \to \infty}}{=} 0, \\
\tilde{\mathbb{J}}_{\mathbb{E},5}^{\sharp} :=& \, \sum_{i=1}^{N+1} 
\left(\int_{\partial \tilde{\mathbb{U}}_{\tilde{\delta}_{\tilde{b}_{i-1}}}} 
\! + \! \int_{\partial \tilde{\mathbb{U}}_{\tilde{\delta}_{\tilde{a}_{i}}}} 
\right) \left(\lim_{-\tilde{\Sigma}_{\tilde{\mathcal{R}}}^{\sharp} \ni 
z^{\prime} \to \xi} \sum_{j=1}^{N+1} \left(\int_{\partial \tilde{
\mathbb{U}}_{\tilde{\delta}_{\tilde{b}_{j-1}}}} \! + \! \int_{\partial 
\tilde{\mathbb{U}}_{\tilde{\delta}_{\tilde{a}_{j}}}} \right) \sum_{i_{1},
i_{2}=1,2} \left\lvert \dfrac{(w_{+}^{\Sigma_{\tilde{\mathcal{R}}}}
(u))_{i_{1}i_{2}}}{u \! - \! z^{\prime}} \right\rvert^{2} \, \dfrac{\lvert 
\md u \rvert}{2 \pi^{2}} \right) \lvert \md \xi \rvert \underset{
\underset{z_{o}=1+o(1)}{\mathscr{N},n \to \infty}}{=} 0,
\end{align*}
whence
\begin{equation} \label{eqiy70}
\tilde{\mathbb{J}}_{\mathbb{E}}^{\sharp} \underset{\underset{z_{o}=
1+o(1)}{\mathscr{N},n \to \infty}}{=} 0.
\end{equation}
Thus, via Equation~\eqref{eqiy39}, and the Estimates~\eqref{eqiy62}, 
\eqref{eqiy64}, \eqref{eqiy66}, \eqref{eqiy68}, and~\eqref{eqiy70}, 
one shows that, for $n \! \in \! \mathbb{N}$ and $k \! \in \! \lbrace 
1,2,\dotsc,K \rbrace$ such that $\alpha_{p_{\mathfrak{s}}} \! := \! 
\alpha_{k} \! \neq \! \infty$,
\begin{equation} \label{eqiy71} 
\int_{\tilde{\Sigma}_{\tilde{\mathcal{R}}}^{\sharp}} \left(\lim_{-\tilde{
\Sigma}_{\tilde{\mathcal{R}}}^{\sharp} \ni z^{\prime} \to \xi} \left\lvert 
\left\lvert \dfrac{w_{+}^{\Sigma_{\tilde{\mathcal{R}}}}(\pmb{\cdot})}{2 
\pi \mi (\pmb{\cdot} \! - \! z^{\prime})} \right\rvert \right\rvert_{
\mathcal{L}^{2}_{\mathrm{M}_{2}(\mathbb{C})}(\tilde{\Sigma}_{\tilde{
\mathcal{R}}}^{\sharp})}^{2} \right) \, \lvert \md \xi \rvert \underset{
\underset{z_{o}=1+o(1)}{\mathscr{N},n \to \infty}}{\leqslant} 
\mathcal{O} \left(\tilde{\mathfrak{c}}_{\tilde{\mathbb{J}}^{\sharp}}
(n,k,z_{o}) \me^{-((n-1)K+k) \tilde{\lambda}_{\tilde{\mathcal{R}},
\tilde{w}}^{\triangleright}} \right),
\end{equation}
where $\tilde{\lambda}_{\tilde{\mathcal{R}},\tilde{w}}^{\triangleright}$ 
$(> \! 0)$ is defined by Equation~\eqref{tillamrwr2}, and $\tilde{
\mathfrak{c}}_{\tilde{\mathbb{J}}^{\sharp}}(n,k,z_{o}) \! =_{\underset{
z_{o}=1+o(1)}{\mathscr{N},n \to \infty}} \! \mathcal{O}(1)$. Finally, 
via the Estimates~\eqref{eqiy38} and~\eqref{eqiy71}, and an application 
of both H\"{o}lder's and Minkowski's Inequalities for Integrals, one 
arrives at, after an improper integral analysis, for $n \! \in \! 
\mathbb{N}$ and $k \! \in \! \lbrace 1,2,\dotsc,K \rbrace$ such 
that $\alpha_{p_{\mathfrak{s}}} \! := \! \alpha_{k} \! \neq \! \infty$ 
(cf. Inequality~\eqref{eqiy31})
\begin{align*}
\lvert \lvert (C_{w^{\Sigma_{\tilde{\mathcal{R}}}}}f)(\pmb{\cdot}) \rvert 
\rvert_{\mathcal{L}^{2}_{\mathrm{M}_{2}(\mathbb{C})}(\tilde{\Sigma}_{
\tilde{\mathcal{R}}}^{\sharp})} \underset{\underset{z_{o}=1+o(1)}{
\mathscr{N},n \to \infty}}{\leqslant}& \, \lvert \lvert f(\pmb{\cdot}) 
\rvert \rvert_{\mathcal{L}^{2}_{\mathrm{M}_{2}(\mathbb{C})}(\tilde{
\Sigma}_{\tilde{\mathcal{R}}}^{\sharp})} \left(\mathcal{O} \left(
\underbrace{\tilde{\mathfrak{c}}_{\lozenge}^{\sharp}(n,k,z_{o})}_{= 
\, \mathcal{O}(1)} \me^{-2((n-1)K+k) \tilde{\lambda}_{\tilde{
\mathcal{R}},\tilde{w}}^{\triangleright}} \right) \! + \! \mathcal{O} 
\left(\underbrace{\tilde{\mathfrak{c}}_{\blacklozenge}^{\sharp}(n,
k,z_{o})}_{= \, \mathcal{O}(1)} \me^{-\frac{3}{2}((n-1)K+k) \tilde{
\lambda}_{\tilde{\mathcal{R}},\tilde{w}}^{\triangleright}} \right) 
\right. \\
+&\left. \, \mathcal{O} \left(\underbrace{\tilde{\mathfrak{c}}_{
\spadesuit}^{\sharp}(n,k,z_{o})}_{= \, \mathcal{O}(1)} \me^{-((n-1)
K+k) \tilde{\lambda}_{\tilde{\mathcal{R}},\tilde{w}}^{\triangleright}} 
\right) \right)^{1/2} \\
\underset{\underset{z_{o}=1+o(1)}{\mathscr{N},n \to \infty}}{
\leqslant}& \, \lvert \lvert f(\pmb{\cdot}) \rvert \rvert_{\mathcal{
L}^{2}_{\mathrm{M}_{2}(\mathbb{C})}(\tilde{\Sigma}_{\tilde{
\mathcal{R}}}^{\sharp})} \mathcal{O} \left(\tilde{\mathfrak{c}}_{
\tilde{\mathcal{R}},\tilde{w}}^{\triangleleft}(n,k,z_{o}) \me^{-
\frac{1}{2}((n-1)K+k) \tilde{\lambda}_{\tilde{\mathcal{R}},\tilde{w}}^{
\triangleright}} \right), 
\end{align*}
where $\tilde{\mathfrak{c}}_{\tilde{\mathcal{R}},\tilde{w}}^{\triangleleft}
(n,k,z_{o})$ is characterised in the corresponding item~\pmb{(2)} of the 
lemma; consequently,
\begin{equation*}
\lvert \lvert C_{w^{\Sigma_{\tilde{\mathcal{R}}}}} \rvert \rvert_{
\mathfrak{B}_{2}(\tilde{\Sigma}_{\tilde{\mathcal{R}}}^{\sharp})} 
\underset{\underset{z_{o}=1+o(1)}{\mathscr{N},n \to \infty}}{=} 
\mathcal{O} \left(\tilde{\mathfrak{c}}_{\tilde{\mathcal{R}},\tilde{w}}^{
\triangleleft}(n,k,z_{o}) \me^{-\frac{1}{2}((n-1)K+k) \tilde{\lambda}_{
\tilde{\mathcal{R}},\tilde{w}}^{\triangleright}} \right), 
\end{equation*} 
which is the Estimate~\eqref{eqcwsigtll}. Via the Estimate~\eqref{eqcwsigtll}, 
due to a well-known result for bounded linear operators in Hilbert spaces 
(see, for example, Chaper~5 of \cite{joaweid}), it follows that, in the 
double-scaling limit $\mathscr{N},n \! \to \! \infty$ such that $z_{o} \! 
= \! 1 \! + \! o(1)$, $(\id \! - \! C_{w^{\Sigma_{\tilde{\mathcal{R}}}}})^{-1} 
\! \! \upharpoonright_{\mathcal{L}^{2}_{\mathrm{M}_{2}(\mathbb{C})}
(\tilde{\Sigma}_{\tilde{\mathcal{R}}}^{\sharp})}$ exists and can be inverted 
by a Neumann series, with $\lvert \lvert (\id \! - \! C_{w^{\Sigma_{\tilde{
\mathcal{R}}}}})^{-1} \rvert \rvert_{\mathfrak{B}_{2}(\tilde{\Sigma}_{\tilde{
\mathcal{R}}}^{\sharp})} \! \leqslant_{\underset{z_{o}=1+o(1)}{\mathscr{N},
n \to \infty}} \! (1 \! - \! \lvert \lvert C_{w^{\Sigma_{\tilde{\mathcal{R}}}}} 
\rvert \rvert_{\mathfrak{B}_{2}(\tilde{\Sigma}_{\tilde{\mathcal{R}}}^{
\sharp})})^{-1} \! =_{\underset{z_{o}=1+o(1)}{\mathscr{N},n \to \infty}} 
\! \mathcal{O}(1)$.

The analysis for the case $n \! \in \! \mathbb{N}$ and $k \! \in \! \lbrace 
1,2,\dotsc,K \rbrace$ such that $\alpha_{p_{\mathfrak{s}}} \! := \! 
\alpha_{k} \! = \! \infty$ is, \emph{mutatis mutandis}, analogous, and 
leads to the asymptotics, in the double-scaling limit $\mathscr{N},n \! 
\to \! \infty$ such that $z_{o} \! = \! 1 \! + \! o(1)$, for $\lvert \lvert 
C_{w^{\Sigma_{\hat{\mathcal{R}}}}} \rvert \rvert_{\mathfrak{B}_{r}
(\hat{\Sigma}_{\hat{\mathcal{R}}}^{\sharp})}$, $r \! \in \! \lbrace 
2,\infty \rbrace$, stated in Equations~\eqref{eqcwsight} 
and~\eqref{eqcwsightt} (cf. item~\pmb{(1)} of the lemma). \hfill $\qed$
\begin{ccccc} \label{lem5.4} 
For $n \! \in \! \mathbb{N}$ and $k \! \in \! \lbrace 1,2,\dotsc,K \rbrace$ 
such that $\alpha_{p_{\mathfrak{s}}} \! := \! \alpha_{k} \! = \! \infty$ 
(resp., $\alpha_{p_{\mathfrak{s}}} \! := \! \alpha_{k} \! \neq \! \infty)$, 
let $\hat{\mathcal{R}} \colon \mathbb{C} \setminus \hat{\Sigma}_{
\hat{\mathcal{R}}}^{\sharp} \! \to \! \mathrm{SL}_{2}(\mathbb{C})$ (resp., 
$\tilde{\mathcal{R}} \colon \mathbb{C} \setminus \tilde{\Sigma}_{\tilde{
\mathcal{R}}}^{\sharp} \! \to \! \mathrm{SL}_{2}(\mathbb{C}))$ solve the 
equivalent {\rm RHP} $(\hat{\mathcal{R}}(z),\hat{v}_{\hat{\mathcal{R}}}(z),
\hat{\Sigma}_{\hat{\mathcal{R}}}^{\sharp})$ (resp., $(\tilde{\mathcal{R}}
(z),\tilde{v}_{\tilde{\mathcal{R}}}(z),\tilde{\Sigma}_{\tilde{\mathcal{R}}}^{
\sharp}))$ with associated integral representation given by 
Equation~\eqref{eqsek5a} (resp., Equation~\eqref{eqsek5b}$)$, let 
the corresponding operator norms stated in item~{\rm \pmb{(1)}} 
(resp., item~{\rm \pmb{(2)}}$)$ of Lemma~\ref{lem5.3} be valid, 
and set $\hat{\Sigma}_{\circlearrowright} \! := \! \cup_{j=1}^{N+1}
(\partial \hat{\mathbb{U}}_{\hat{\delta}_{\hat{b}_{j-1}}} \cup 
\partial \hat{\mathbb{U}}_{\hat{\delta}_{\hat{a}_{j}}})$ and 
$\hat{\Sigma}_{\scriptscriptstyle \blacksquare} \! := \! \hat{\Sigma}_{
\hat{\mathcal{R}}}^{\sharp} \setminus \hat{\Sigma}_{\circlearrowright}$ 
(resp., $\tilde{\Sigma}_{\circlearrowright} \! := \! \cup_{j=1}^{N+1}
(\partial \tilde{\mathbb{U}}_{\tilde{\delta}_{\tilde{b}_{j-1}}} \cup 
\partial \tilde{\mathbb{U}}_{\tilde{\delta}_{\tilde{a}_{j}}})$ and $\tilde{
\Sigma}_{\scriptscriptstyle \blacksquare} \! := \! \tilde{\Sigma}_{\tilde{
\mathcal{R}}}^{\sharp} \setminus \tilde{\Sigma}_{\circlearrowright})$. 
Then$:$ {\rm \pmb{(1)}} for $n \! \in \! \mathbb{N}$ and $k \! \in \! 
\lbrace 1,2,\dotsc,K \rbrace$ such that $\alpha_{p_{\mathfrak{s}}} 
\! := \! \alpha_{k} \! = \! \infty$, uniformly for compact subsets 
of $\mathbb{C} \setminus \hat{\Sigma}_{\hat{\mathcal{R}}}^{\sharp}$ 
$(\ni \! z)$,
\begin{equation} \label{eqtazh1} 
\hat{\mathcal{R}}(z) \underset{\underset{z_{o}=1+o(1)}{\mathscr{N},
n \to \infty}}{=} \mathrm{I} \! + \! \int_{\hat{\Sigma}_{
\circlearrowright}} \dfrac{w^{\hat{\Sigma}_{\circlearrowright}}_{+}
(\xi)}{\xi \! - \! z} \, \dfrac{\md \xi}{2 \pi \mi} \! + \! \mathcal{O} 
\left(\dfrac{\hat{\mathfrak{c}}_{\hat{\Sigma}_{\circlearrowright}}
(n,k,z_{o}) \me^{-\frac{1}{2}((n-1)K+k) \hat{\lambda}_{\hat{
\mathcal{R}},\hat{w}}^{\triangleright}}}{((n \! - \! 1)K \! + \! k) 
\operatorname{dist}(\hat{\Sigma}_{\hat{\mathcal{R}}}^{\sharp},z)} 
\right), \quad z \! \in \! \mathbb{C} \setminus \hat{\Sigma}_{
\hat{\mathcal{R}}}^{\sharp},
\end{equation}
where $w_{+}^{\hat{\Sigma}_{\circlearrowright}}(z) \! := \! w_{+}^{
\Sigma_{\hat{\mathcal{R}}}}(z) \! \! \upharpoonright_{\hat{\Sigma}_{
\circlearrowright}}$, $\hat{\lambda}_{\hat{\mathcal{R}},\hat{w}}^{
\triangleright}$ $(> \! 0)$ is defined by Equation~\eqref{hatlamrwr1}, 
$\operatorname{dist}(\hat{\Sigma}_{\hat{\mathcal{R}}}^{\sharp},z) 
\! := \! \inf \lbrace \mathstrut \lvert u \! - \! z \rvert; \, u \! \in \! 
\hat{\Sigma}_{\hat{\mathcal{R}}}^{\sharp} \rbrace$, $z \! \in \! 
\mathbb{C} \setminus \hat{\Sigma}_{\hat{\mathcal{R}}}^{\sharp}$, 
and $(\mathrm{M}_{2}(\mathbb{C}) \! \ni)$ $\hat{\mathfrak{c}}_{
\hat{\Sigma}_{\circlearrowright}}(n,k,z_{o}) \! =_{\underset{z_{o}
=1+o(1)}{\mathscr{N},n \to \infty}} \! \mathcal{O}(1)$$;$ and 
{\rm \pmb{(2)}} for $n \! \in \! \mathbb{N}$ and $k \! \in \! \lbrace 
1,2,\dotsc,K \rbrace$ such that $\alpha_{p_{\mathfrak{s}}} \! := \! 
\alpha_{k} \! \neq \! \infty$, uniformly for compact subsets of 
$\mathbb{C} \setminus \tilde{\Sigma}_{\tilde{\mathcal{R}}}^{\sharp}$ 
$(\ni \! z)$,
\begin{equation} \label{eqtazt1} 
\tilde{\mathcal{R}}(z) \underset{\underset{z_{o}=1+o(1)}{\mathscr{N},
n \to \infty}}{=} \mathrm{I} \! + \! \int_{\tilde{\Sigma}_{
\circlearrowright}} \dfrac{(z \! - \! \alpha_{k})w^{\tilde{\Sigma}_{
\circlearrowright}}_{+}(\xi)}{(\xi \! - \! \alpha_{k})(\xi \! - \! z)} \, 
\dfrac{\md \xi}{2 \pi \mi} \! + \! \mathcal{O} \left(\dfrac{\tilde{
\mathfrak{c}}_{\tilde{\Sigma}_{\circlearrowright}}(n,k,z_{o}) 
\me^{-\frac{1}{2}((n-1)K+k) \tilde{\lambda}_{\tilde{\mathcal{R}},
\tilde{w}}^{\triangleright}}}{((n \! - \! 1)K \! + \! k) \min \lbrace 
1,\operatorname{dist}(\tilde{\Sigma}_{\tilde{\mathcal{R}}}^{\sharp},
z) \rbrace} \right), \quad z \! \in \! \mathbb{C} \setminus 
\tilde{\Sigma}_{\tilde{\mathcal{R}}}^{\sharp},
\end{equation}
where $w_{+}^{\tilde{\Sigma}_{\circlearrowright}}(z) \! := \! w_{+}^{
\Sigma_{\tilde{\mathcal{R}}}}(z) \! \! \upharpoonright_{\tilde{\Sigma}_{
\circlearrowright}}$, $\tilde{\lambda}_{\tilde{\mathcal{R}},\tilde{w}}^{
\triangleright}$ $(> \! 0)$ is defined by Equation~\eqref{tillamrwr2}, 
$\operatorname{dist}(\tilde{\Sigma}_{\tilde{\mathcal{R}}}^{\sharp},z) 
\! := \! \inf \lbrace \mathstrut \lvert u \! - \! z \rvert; \, u \! \in \! 
\tilde{\Sigma}_{\tilde{\mathcal{R}}}^{\sharp} \rbrace$, $z \! \in \! 
\mathbb{C} \setminus \tilde{\Sigma}_{\tilde{\mathcal{R}}}^{\sharp}$, 
and $(\mathrm{M}_{2}(\mathbb{C}) \! \ni)$ $\tilde{\mathfrak{c}}_{
\tilde{\Sigma}_{\circlearrowright}}(n,k,z_{o}) \! =_{\underset{z_{o}=
1+o(1)}{\mathscr{N},n \to \infty}} \! \mathcal{O}(1)$.
\end{ccccc}

\emph{Proof}. The proof of this Lemma~\ref{lem5.4} consists of two 
cases: (i) $n \! \in \! \mathbb{N}$ and $k \! \in \! \lbrace 1,2,\dotsc,
K \rbrace$ such that $\alpha_{p_{\mathfrak{s}}} \! := \! \alpha_{k} \! 
= \! \infty$; and (ii) $n \! \in \! \mathbb{N}$ and $k \! \in \! \lbrace 
1,2,\dotsc,K \rbrace$ such that $\alpha_{p_{\mathfrak{s}}} \! := \! 
\alpha_{k} \! \neq \! \infty$. Notwithstanding the fact that the scheme 
of the proof is, \emph{mutatis mutandis}, similar for both cases, 
case~(ii), nevertheless, is the more technically challenging of the two; 
therefore, without loss of generality, only the proof for case~(ii) is 
presented in detail, whilst case~(i) is proved analogously.

Recall that, for $n \! \in \! \mathbb{N}$ and $k \! \in \! \lbrace 1,2,
\dotsc,K \rbrace$ such that $\alpha_{p_{\mathfrak{s}}} \! := \! 
\alpha_{k} \! \neq \! \infty$, $\tilde{\mathcal{R}} \colon \mathbb{C} 
\setminus \tilde{\Sigma}_{\tilde{\mathcal{R}}}^{\sharp} \! \to \! 
\mathrm{SL}_{2}(\mathbb{C})$ solves the equivalent RHP $(\tilde{
\mathcal{R}}(z),\tilde{v}_{\tilde{\mathcal{R}}}(z),\tilde{\Sigma}_{\tilde{
\mathcal{R}}}^{\sharp})$ with associated integral representation (cf. 
Equation~\eqref{eqsek5b})
\begin{equation*}
\tilde{\mathcal{R}}(z) \! = \! \mathrm{I} \! + \! \int_{\tilde{
\Sigma}_{\tilde{\mathcal{R}}}^{\sharp}} \dfrac{(z \! - \! \alpha_{k}) 
\mu^{\Sigma_{\tilde{\mathcal{R}}}}(\xi)w_{+}^{\Sigma_{\tilde{\mathcal{
R}}}}(\xi)}{(\xi \! - \! \alpha_{k})(\xi \! - \! z)} \, \dfrac{\md \xi}{2 \pi 
\mi}, \quad z \! \in \! \mathbb{C} \setminus \tilde{\Sigma}_{\tilde{
\mathcal{R}}}^{\sharp},
\end{equation*}
where $\mu^{\Sigma_{\tilde{\mathcal{R}}}}(z) \! = \! \tilde{\mathcal{R}}_{
-}(z) \! = \! \tilde{\mathcal{R}}_{+}(z)(\mathrm{I} \! + \! w_{+}^{\Sigma_{
\tilde{\mathcal{R}}}}(z))^{-1}$ solves the linear singular integral equation 
$((\id \! - \! C_{w^{\Sigma_{\tilde{\mathcal{R}}}}}) \mu^{\Sigma_{\tilde{
\mathcal{R}}}})(z) \! = \! \mathrm{I}$, $z \! \in \! \tilde{\Sigma}_{\tilde{
\mathcal{R}}}^{\sharp}$, with $\mathcal{L}^{2}_{\mathrm{M}_{2}
(\mathbb{C})}(\tilde{\Sigma}_{\tilde{\mathcal{R}}}^{\sharp}) \! \ni \! 
f \! \mapsto \! (C_{w^{\Sigma_{\tilde{\mathcal{R}}}}}f)(z) \! := \! 
\lim\limits_{-\tilde{\Sigma}_{\tilde{\mathcal{R}}}^{\sharp} \ni z^{\prime} 
\to z} \int_{\tilde{\Sigma}_{\tilde{\mathcal{R}}}^{\sharp}} \tfrac{(z^{\prime}
-\alpha_{k})(fw_{+}^{\Sigma_{\tilde{\mathcal{R}}}})(\xi)}{(\xi -\alpha_{k})
(\xi -z^{\prime})} \, \tfrac{\md \xi}{2 \pi \mi}$, where $-\tilde{\Sigma}_{
\tilde{\mathcal{R}}}^{\sharp}$ denotes the `minus side' of the oriented 
skeleton $\tilde{\Sigma}_{\tilde{\mathcal{R}}}^{\sharp}$. For $n \! \in \! 
\mathbb{N}$ and $k \! \in \! \lbrace 1,2,\dotsc,K \rbrace$ such that 
$\alpha_{p_{\mathfrak{s}}} \! := \! \alpha_{k} \! \neq \! \infty$, let 
$\tilde{\Sigma}_{\circlearrowright} \! := \! \cup_{j=1}^{N+1}(\partial 
\tilde{\mathbb{U}}_{\tilde{\delta}_{\tilde{b}_{j-1}}} \cup \partial \tilde{
\mathbb{U}}_{\tilde{\delta}_{\tilde{a}_{j}}})$ and $\tilde{\Sigma}_{
\scriptscriptstyle \blacksquare} \! := \! \tilde{\Sigma}_{\tilde{
\mathcal{R}}}^{\sharp} \setminus \tilde{\Sigma}_{\circlearrowright}$, 
and write $\tilde{\Sigma}_{\tilde{\mathcal{R}}}^{\sharp} \! = \! 
\tilde{\Sigma}_{\circlearrowright} \cup \tilde{\Sigma}_{\scriptscriptstyle 
\blacksquare}$ (with $\tilde{\Sigma}_{\circlearrowright} \cap \tilde{
\Sigma}_{\scriptscriptstyle \blacksquare} \! = \! \varnothing)$; then, 
via the linearity property of the (normalised at $\alpha_{k}$) Cauchy 
integral operator $C_{w^{\Sigma_{\tilde{\mathcal{R}}}}}$, a calculation 
shows that the action of $C_{w^{\Sigma_{\tilde{\mathcal{R}}}}}$ can be 
presented as $C_{w^{\Sigma_{\tilde{\mathcal{R}}}}} \! = \! C_{w^{\tilde{
\Sigma}_{\circlearrowright}}} \! + \! C_{w^{\tilde{\Sigma}_{\scriptscriptstyle 
\blacksquare}}}$, where $\mathcal{L}^{2}_{\mathrm{M}_{2}(\mathbb{C})}
(\tilde{\Sigma}_{\circlearrowright}) \! \ni \! f_{1} \! \mapsto \! (C_{w^{
\tilde{\Sigma}_{\circlearrowright}}}f_{1})(z) \! := \! \lim\limits_{-\tilde{
\Sigma}_{\circlearrowright} \ni z^{\prime} \to z} \int_{\tilde{\Sigma}_{
\circlearrowright}} \tfrac{(z^{\prime}-\alpha_{k})(f_{1}w_{+}^{\tilde{
\Sigma}_{\circlearrowright}})(\xi)}{(\xi -\alpha_{k})(\xi -z^{\prime})} 
\, \tfrac{\md \xi}{2 \pi \mi}$, with $w_{+}^{\tilde{\Sigma}_{
\circlearrowright}}(z) \! := \! w_{+}^{\Sigma_{\tilde{\mathcal{R}}}}(z) 
\! \! \upharpoonright_{\tilde{\Sigma}_{\circlearrowright}}$, and 
$\mathcal{L}^{2}_{\mathrm{M}_{2}(\mathbb{C})}(\tilde{\Sigma}_{
\scriptscriptstyle \blacksquare}) \! \ni \! f_{2} \! \mapsto \! (C_{w^{
\tilde{\Sigma}_{\scriptscriptstyle \blacksquare}}}f_{2})(z) \! := \! 
\lim\limits_{-\tilde{\Sigma}_{\scriptscriptstyle \blacksquare} \ni 
z^{\prime} \to z} \int_{\tilde{\Sigma}_{\scriptscriptstyle \blacksquare}} 
\tfrac{(z^{\prime}-\alpha_{k})(f_{1}w_{+}^{\tilde{\Sigma}_{\scriptscriptstyle 
\blacksquare}})(\xi)}{(\xi -\alpha_{k})(\xi -z^{\prime})} \, \tfrac{\md \xi}{2 
\pi \mi}$, with $w_{+}^{\tilde{\Sigma}_{\scriptscriptstyle \blacksquare}}(z) 
\! := \! w_{+}^{\Sigma_{\tilde{\mathcal{R}}}}(z) \! \! \upharpoonright_{
\tilde{\Sigma}_{\scriptscriptstyle \blacksquare}}$. Via the latter 
decompositions, iteration of the second resolvent identity 
{}\footnote{For operators $\mathfrak{O}_{1}$ and $\mathfrak{O}_{2}$, 
with $\operatorname{dom}(\mathfrak{O}_{1}) \! = \! \operatorname{dom}
(\mathfrak{O}_{2})$, if $(\id \! - \! \mathfrak{O}_{1})^{-1}$ and $(\id 
\! - \! \mathfrak{O}_{2})^{-1}$ exist, then $(\id \! - \! \mathfrak{O}_{
2})^{-1} \! - \! (\id \! - \! \mathfrak{O}_{1})^{-1} \! = \! (\id \! - \! 
\mathfrak{O}_{2})^{-1}(\mathfrak{O}_{2} \! - \! \mathfrak{O}_{1})
(\id \! - \! \mathfrak{O}_{1})^{-1}$ (see, for example, \cite{joaweid}).} 
shows that
\begin{align*}
\mu^{\Sigma_{\tilde{\mathcal{R}}}}(z) &= \mathrm{I} \! + \! ((\id \! - \! 
C_{w^{\Sigma_{\tilde{\mathcal{R}}}}})^{-1}C_{w^{\Sigma_{\tilde{\mathcal{
R}}}}} \mathrm{I})(z) \! = \! \mathrm{I} \! + \! ((\id \! - \! C_{w^{\tilde{
\Sigma}_{\circlearrowright}}} \! - \! C_{w^{\tilde{\Sigma}_{\scriptscriptstyle
\blacksquare}}})^{-1}(C_{w^{\tilde{\Sigma}_{\circlearrowright}}} \! + \! 
C_{w^{\tilde{\Sigma}_{\scriptscriptstyle \blacksquare}}}) \mathrm{I})(z) \\
&= \mathrm{I} \! + \! ((\id \! - \! C_{w^{\tilde{\Sigma}_{\circlearrowright}}}
\! - \! C_{w^{\tilde{\Sigma}_{\scriptscriptstyle \blacksquare}}})^{-1}
C_{w^{\tilde{\Sigma}_{\circlearrowright}}} \mathrm{I})(z) \! + \! ((\id \! - 
\! C_{w^{\tilde{\Sigma}_{\circlearrowright}}} \! - \! C_{w^{\tilde{\Sigma}_{
\scriptscriptstyle \blacksquare}}})^{-1}C_{w^{\tilde{\Sigma}_{\scriptscriptstyle 
\blacksquare}}} \mathrm{I})(z) \\
&= \mathrm{I} \! + \! (((\id \! - \! C_{w^{\tilde{\Sigma}_{\circlearrowright}}})
(\id \! - \! (\id \! - \! C_{w^{\tilde{\Sigma}_{\circlearrowright}}})^{-1}
C_{w^{\tilde{\Sigma}_{\scriptscriptstyle \blacksquare}}}))^{-1}C_{w^{\tilde{
\Sigma}_{\circlearrowright}}} \mathrm{I})(z) \\
&+ (((\id \! - \! C_{w^{\tilde{\Sigma}_{\scriptscriptstyle \blacksquare}}})
(\id \! - \! (\id \! - \! C_{w^{\tilde{\Sigma}_{\scriptscriptstyle 
\blacksquare}}})^{-1}C_{w^{\tilde{\Sigma}_{\circlearrowright}}}))^{-1}
C_{w^{\tilde{\Sigma}_{\scriptscriptstyle \blacksquare}}} \mathrm{I})(z) \\
&= \mathrm{I} \! + \! ((\id \! - \! (\id \! - \! C_{w^{\tilde{\Sigma}_{
\circlearrowright}}})^{-1}C_{w^{\tilde{\Sigma}_{\scriptscriptstyle
\blacksquare}}})^{-1}(\id \! + \! (\id \! - \! C_{w^{\tilde{\Sigma}_{
\circlearrowright}}})^{-1}C_{w^{\tilde{\Sigma}_{\circlearrowright}}})
C_{w^{\tilde{\Sigma}_{\circlearrowright}}} \mathrm{I})(z) \\
&+ ((\id \! - \! (\id \! - \! C_{w^{\tilde{\Sigma}_{\scriptscriptstyle
\blacksquare}}})^{-1}C_{w^{\tilde{\Sigma}_{\circlearrowright}}})^{-1}
(\id \! + \! (\id \! - \! C_{w^{\tilde{\Sigma}_{\scriptscriptstyle 
\blacksquare}}})^{-1}C_{w^{\tilde{\Sigma}_{\scriptscriptstyle 
\blacksquare}}})C_{w^{\tilde{\Sigma}_{\scriptscriptstyle 
\blacksquare}^{o}}} \mathrm{I})(z) \\
&= \mathrm{I} \! + \! ((\id \! - \! (\id \! - \! C_{w^{\tilde{\Sigma}_{
\circlearrowright}}})^{-1}C_{w^{\tilde{\Sigma}_{\scriptscriptstyle
\blacksquare}}})^{-1}((\id \! - \! C_{w^{\tilde{\Sigma}_{\circlearrowright}}})^{-1}
C_{w^{\tilde{\Sigma}_{\circlearrowright}}})(C_{w^{\tilde{\Sigma}_{
\circlearrowright}}} \mathrm{I}))(z) \\
&+ ((\id \! - \! (\id \! - \! C_{w^{\tilde{\Sigma}_{\scriptscriptstyle
\blacksquare}}})^{-1}C_{w^{\tilde{\Sigma}_{\circlearrowright}}})^{-1}
((\id \! - \! C_{w^{\tilde{\Sigma}_{\scriptscriptstyle \blacksquare}}})^{-1}
C_{w^{\tilde{\Sigma}_{\scriptscriptstyle \blacksquare}}})(C_{w^{\tilde{
\Sigma}_{\scriptscriptstyle \blacksquare}}} \mathrm{I}))(z) \\
&+ ((\id \! - \! (\id \! - \! C_{w^{\tilde{\Sigma}_{\circlearrowright}}})^{-1}
C_{w^{\tilde{\Sigma}_{\scriptscriptstyle \blacksquare}}})^{-1}
(C_{w^{\tilde{\Sigma}_{\circlearrowright}}} \mathrm{I}))(z) \! + \! ((\id 
\! - \! (\id \!- \! C_{w^{\tilde{\Sigma}_{\scriptscriptstyle \blacksquare}}})^{-1}
C_{w^{\tilde{\Sigma}_{\circlearrowright}}})^{-1}(C_{w^{\tilde{\Sigma}_{
\scriptscriptstyle\blacksquare}}} \mathrm{I}))(z) \\
&= \mathrm{I} \! + \! ((\id \! + \! (\id \! - \! (\id \! - \! C_{w^{\tilde{
\Sigma}_{\circlearrowright}}})^{-1}C_{w^{\tilde{\Sigma}_{\scriptscriptstyle
\blacksquare}}})^{-1}(\id \! - \! C_{w^{\tilde{\Sigma}_{\circlearrowright}}})^{-1}
C_{w^{\tilde{\Sigma}_{\scriptscriptstyle \blacksquare}}})(C_{w^{\tilde{
\Sigma}_{\circlearrowright}}} \mathrm{I}))(z) \\
&+ ((\id \! + \! (\id \! - \! (\id \! - \! C_{w^{\tilde{\Sigma}_{\scriptscriptstyle
\blacksquare}}})^{-1}C_{w^{\tilde{\Sigma}_{\circlearrowright}}})^{-1}
(\id\! - \! C_{w^{\tilde{\Sigma}_{\scriptscriptstyle \blacksquare}}})^{-1}
C_{w^{\tilde{\Sigma}_{\circlearrowright}}})(C_{w^{\tilde{\Sigma}_{
\scriptscriptstyle \blacksquare}}} \mathrm{I}))(z) \\
&+ ((\id \! + \! (\id \! - \! (\id \! - \! C_{w^{\tilde{\Sigma}_{
\circlearrowright}}})^{-1}C_{w^{\tilde{\Sigma}_{\scriptscriptstyle
\blacksquare}}})^{-1}(\id \! - \! C_{w^{\tilde{\Sigma}_{\circlearrowright}}})^{-1}
C_{w^{\tilde{\Sigma}_{\scriptscriptstyle \blacksquare}}})(\id \! -\! 
C_{w^{\tilde{\Sigma}_{\circlearrowright}}})^{-1}C_{w^{\tilde{\Sigma}_{
\circlearrowright}}}(C_{w^{\tilde{\Sigma}_{\circlearrowright}}} \mathrm{I}))(z) \\
&+ ((\id \! + \! (\id \! - \! (\id \! - \! C_{w^{\tilde{\Sigma}_{\scriptscriptstyle
\blacksquare}}})^{-1}C_{w^{\tilde{\Sigma}_{\circlearrowright}}})^{-1}
(\id \! - \! C_{w^{\tilde{\Sigma}_{\scriptscriptstyle \blacksquare}}})^{-1}
C_{w^{\tilde{\Sigma}_{\circlearrowright}}})(\id \! - \! C_{w^{\tilde{
\Sigma}_{\scriptscriptstyle \blacksquare}}})^{-1}C_{w^{\tilde{\Sigma}_{
\scriptscriptstyle \blacksquare}}}(C_{w^{\tilde{\Sigma}_{\scriptscriptstyle
\blacksquare}}} \mathrm{I}))(z) \\
&= \mathrm{I} \! + \! (C_{w^{\tilde{\Sigma}_{\circlearrowright}^{o}}} 
\mathrm{I})(z) \! + \! (C_{w^{\tilde{\Sigma}_{\scriptscriptstyle 
\blacksquare}}} \mathrm{I})(z) \! + \! ((\id \! - \! C_{w^{\tilde{\Sigma}_{
\circlearrowright}}})^{-1}C_{w^{\tilde{\Sigma}_{\circlearrowright}}}
(C_{w^{\tilde{\Sigma}_{\circlearrowright}}} \mathrm{I}))(z) \\
&+ ((\id \! - \! C_{w^{\tilde{\Sigma}_{\scriptscriptstyle \blacksquare}}})^{-1}
C_{w^{\tilde{\Sigma}_{\scriptscriptstyle \blacksquare}}}(C_{w^{\tilde{
\Sigma}_{\scriptscriptstyle \blacksquare}}} \mathrm{I}))(z) \! + \! ((\id 
\! - \! (\id \! - \! C_{w^{\tilde{\Sigma}_{\circlearrowright}}})^{-1}
C_{w^{\tilde{\Sigma}_{\scriptscriptstyle \blacksquare}}})^{-1}
(\id \! - \! C_{w^{\tilde{\Sigma}_{\circlearrowright}}})^{-1}
C_{w^{\tilde{\Sigma}_{\scriptscriptstyle \blacksquare}}} \\
&\times (C_{w^{\tilde{\Sigma}_{\circlearrowright}}} \mathrm{I}))(z) \! 
+ \!((\id \! - \! (\id \! - \! C_{w^{\tilde{\Sigma}_{\scriptscriptstyle 
\blacksquare}}})^{-1}C_{w^{\tilde{\Sigma}_{\circlearrowright}}})^{-1}
(\id \! - \! C_{w^{\tilde{\Sigma}_{\scriptscriptstyle \blacksquare}}})^{-1}
C_{w^{\tilde{\Sigma}_{\circlearrowright}}}(C_{w^{\tilde{\Sigma}_{
\scriptscriptstyle \blacksquare}}} \mathrm{I}))(z) \\
&+ ((\id \! - \! (\id \! - \! C_{w^{\tilde{\Sigma}_{\circlearrowright}}})^{-1}
C_{w^{\tilde{\Sigma}_{\scriptscriptstyle \blacksquare}}})^{-1}(\id \! - \! 
C_{w^{\tilde{\Sigma}_{\circlearrowright}}})^{-1}C_{w^{\tilde{\Sigma}_{
\scriptscriptstyle \blacksquare}}}(\id \! - \! C_{w^{\tilde{\Sigma}_{
\circlearrowright}}})^{-1}C_{w^{\tilde{\Sigma}_{\circlearrowright}}}
(C_{w^{\tilde{\Sigma}_{\circlearrowright}}} \mathrm{I}))(z) \\
&+ ((\id \! - \! (\id \! - \! C_{w^{\tilde{\Sigma}_{\scriptscriptstyle
\blacksquare}}})^{-1}C_{w^{\tilde{\Sigma}_{\circlearrowright}}})^{-1}
(\id \! - \! C_{w^{\tilde{\Sigma}_{\scriptscriptstyle \blacksquare}}})^{
-1}C_{w^{\tilde{\Sigma}_{\circlearrowright}}}(\id \! - \! C_{w^{\tilde{
\Sigma}_{\scriptscriptstyle \blacksquare}}})^{-1}C_{w^{\tilde{\Sigma}_{
\scriptscriptstyle \blacksquare}}}(C_{w^{\tilde{\Sigma}_{\scriptscriptstyle
\blacksquare}}} \mathrm{I}))(z) \\
&= \mathrm{I} \! + \! (C_{w^{\tilde{\Sigma}_{\circlearrowright}}} 
\mathrm{I})(z) \! + \! (C_{w^{\tilde{\Sigma}_{\scriptscriptstyle 
\blacksquare}}} \mathrm{I})(z) \! + \! ((\id \! - \! C_{w^{\tilde{\Sigma}_{
\circlearrowright}}})^{-1}C_{w^{\tilde{\Sigma}_{\circlearrowright}}}
(C_{w^{\tilde{\Sigma}_{\circlearrowright}}} \mathrm{I}))(z) \! + \! ((\id 
\! - \! C_{w^{\tilde{\Sigma}_{\scriptscriptstyle \blacksquare}}})^{-1}
C_{w^{\tilde{\Sigma}_{\scriptscriptstyle \blacksquare}}} \\
&\times (C_{w^{\tilde{\Sigma}_{\scriptscriptstyle \blacksquare}}} 
\mathrm{I}))(z) \! + \! ((\id \! - \! (\id \! - \! C_{w^{\tilde{\Sigma}_{
\circlearrowright}}})^{-1}(\id \! - \! C_{w^{\tilde{\Sigma}_{
\scriptscriptstyle \blacksquare}}})^{-1}C_{w^{\tilde{\Sigma}_{
\scriptscriptstyle \blacksquare}}}C_{w^{\tilde{\Sigma}_{
\circlearrowright}}})^{-1}(\id \! - \! C_{w^{\tilde{\Sigma}_{
\circlearrowright}}})^{-1}(\id \! - \! C_{w^{\tilde{\Sigma}_{
\scriptscriptstyle \blacksquare}}})^{-1} \\
&\times C_{w^{\tilde{\Sigma}_{\circlearrowright}}}(C_{w^{\tilde{
\Sigma}_{\scriptscriptstyle \blacksquare}}} \mathrm{I}))(z) \! + \! 
((\id \! - \! (\id \! - \! C_{w^{\tilde{\Sigma}_{\scriptscriptstyle 
\blacksquare}}})^{-1}(\id \! - \! C_{w^{\tilde{\Sigma}_{
\circlearrowright}}})^{-1}C_{w^{\tilde{\Sigma}_{\circlearrowright}}}
C_{w^{\tilde{\Sigma}_{\scriptscriptstyle \blacksquare}}})^{-1}
(\id \! - \! C_{w^{\tilde{\Sigma}_{\scriptscriptstyle \blacksquare}}})^{-1} \\
&\times (\id \! - \! C_{w^{\tilde{\Sigma}_{\circlearrowright}}})^{-1}
C_{w^{\tilde{\Sigma}_{\scriptscriptstyle \blacksquare}}}(C_{w^{\tilde{
\Sigma}_{\circlearrowright}}} \mathrm{I}))(z) \! + \! ((\id \! - \! (\id \! 
- \! C_{w^{\tilde{\Sigma}_{\circlearrowright}}})^{-1}(\id \! - \! C_{w^{
\tilde{\Sigma}_{\scriptscriptstyle \blacksquare}}})^{-1}C_{w^{\tilde{
\Sigma}_{\scriptscriptstyle \blacksquare}}}C_{w^{\tilde{\Sigma}_{
\circlearrowright}}})^{-1} \\
&\times (\id \! - \! C_{w^{\tilde{\Sigma}_{\circlearrowright}}})^{-1}
(\id \! - \! C_{w^{\tilde{\Sigma}_{\scriptscriptstyle \blacksquare}}})^{
-1}C_{w^{\tilde{\Sigma}_{\circlearrowright}}}(\id \! - \! C_{w^{\tilde{
\Sigma}_{\scriptscriptstyle \blacksquare}}})^{-1}C_{w^{\tilde{\Sigma}_{
\scriptscriptstyle \blacksquare}}}(C_{w^{\tilde{\Sigma}_{\scriptscriptstyle
\blacksquare}}} \mathrm{I}))(z) \! + \! ((\id \! - \! (\id \! - \! 
C_{w^{\tilde{\Sigma}_{\scriptscriptstyle \blacksquare}}})^{-1} \\
&\times (\id \! - \! C_{w^{\tilde{\Sigma}_{\circlearrowright}}})^{-1}
C_{w^{\tilde{\Sigma}_{\circlearrowright}}}C_{w^{\tilde{\Sigma}_{
\scriptscriptstyle \blacksquare}}})^{-1}(\id \! - \! C_{w^{\tilde{\Sigma}_{
\scriptscriptstyle \blacksquare}}})^{-1}(\id \! - \! C_{w^{\tilde{\Sigma}_{
\circlearrowright}}})^{-1}C_{w^{\tilde{\Sigma}_{\scriptscriptstyle 
\blacksquare}}}(\id \!- \! C_{w^{\tilde{\Sigma}_{\circlearrowright}}})^{-1}
C_{w^{\tilde{\Sigma}_{\circlearrowright}}}(C_{w^{\tilde{\Sigma}_{
\circlearrowright}}} \mathrm{I}))(z);
\end{align*}
hence, via the above representation for $\mu^{\Sigma_{\tilde{\mathcal{R}}}}
(z)$ and a partial-fraction decomposition argument, one arrives at, for $n \! 
\in \! \mathbb{N}$ and $k \! \in \! \lbrace 1,2,\dotsc,K \rbrace$ such that 
$\alpha_{p_{\mathfrak{s}}} \! := \! \alpha_{k} \! \neq \! \infty$,
\begin{equation} \label{eqtazq1} 
\tilde{\mathcal{R}}(z) \! - \! \mathrm{I} \! - \! \int_{\tilde{\Sigma}_{
\circlearrowright}} \dfrac{(z \! - \! \alpha_{k})w_{+}^{\tilde{\Sigma}_{
\circlearrowright}}(\xi)}{(\xi \! - \! \alpha_{k})(\xi \! - \! z)} \, \dfrac{\md 
\xi}{2 \pi \mi} \! = \! \int_{\tilde{\Sigma}_{\scriptscriptstyle \blacksquare}}
w_{+}^{\tilde{\Sigma}_{\scriptscriptstyle \blacksquare}}(\xi) \left(\dfrac{1}{
\xi \! - \! z} \! - \! \dfrac{1}{\xi \! - \! \alpha_{k}} \right) \dfrac{\md \xi}{2 
\pi \mi} \! + \! \sum_{m=1}^{8} \tilde{\mathbb{J}}_{m}^{\tilde{\mathcal{R}}}
(z), \quad z \! \in \! \overline{\mathbb{C}} \setminus \tilde{\Sigma}_{\tilde{
\mathcal{R}}}^{\sharp},
\end{equation}
where
\begin{align}
\tilde{\mathbb{J}}^{\tilde{\mathcal{R}}}_{1}(z) :=& \, \int_{\tilde{\Sigma}_{
\circlearrowright}}(C_{w^{\tilde{\Sigma}_{\scriptscriptstyle \blacksquare}}} 
\mathrm{I})(\xi)w_{+}^{\tilde{\Sigma}_{\circlearrowright}}(\xi) \left(
\dfrac{1}{\xi \! - \! z} \! - \! \dfrac{1}{\xi \! - \! \alpha_{k}} \right) 
\dfrac{\md \xi}{2 \pi \mi} \! + \! \int_{\tilde{\Sigma}_{\scriptscriptstyle 
\blacksquare}}(C_{w^{\tilde{\Sigma}_{\scriptscriptstyle \blacksquare}}} 
\mathrm{I})(\xi)w_{+}^{\tilde{\Sigma}_{\scriptscriptstyle \blacksquare}}
(\xi) \left(\dfrac{1}{\xi \! - \! z} \! - \! \dfrac{1}{\xi \! - \! \alpha_{k}} 
\right) \dfrac{\md \xi}{2 \pi \mi}, \label{eqtazq2} \\
\tilde{\mathbb{J}}^{\tilde{\mathcal{R}}}_{2}(z) :=& \, \int_{\tilde{\Sigma}_{
\circlearrowright}}(C_{w^{\tilde{\Sigma}_{\circlearrowright}}} \mathrm{I})
(\xi)w_{+}^{\tilde{\Sigma}_{\circlearrowright}}(\xi) \left(\dfrac{1}{\xi \! 
- \! z} \! - \! \dfrac{1}{\xi \! - \! \alpha_{k}} \right) \dfrac{\md \xi}{2 
\pi \mi} \! + \! \int_{\tilde{\Sigma}_{\scriptscriptstyle \blacksquare}}
(C_{w^{\tilde{\Sigma}_{\circlearrowright}}} \mathrm{I})(\xi)w_{+}^{
\tilde{\Sigma}_{\scriptscriptstyle \blacksquare}}(\xi) \left(\dfrac{1}{\xi 
\! - \! z} \! - \! \dfrac{1}{\xi \! - \! \alpha_{k}} \right) \dfrac{\md 
\xi}{2 \pi \mi}, \label{eqtazq3} \\
\tilde{\mathbb{J}}^{\tilde{\mathcal{R}}}_{3}(z) :=& \, \int_{\tilde{
\Sigma}_{\circlearrowright}}(((\id \! - \! C_{w^{\tilde{\Sigma}_{
\scriptscriptstyle \blacksquare}}})^{-1}C_{w^{\tilde{\Sigma}_{
\scriptscriptstyle \blacksquare}}})(C_{w^{\tilde{\Sigma}_{
\scriptscriptstyle\blacksquare}}} \mathrm{I}))(\xi)w_{+}^{\tilde{\Sigma}_{
\circlearrowright}}(\xi) \left(\dfrac{1}{\xi \! - \! z} \! - \! \dfrac{1}{\xi 
\! - \! \alpha_{k}} \right) \dfrac{\md \xi}{2 \pi \mi} \nonumber \\
+& \, \int_{\tilde{\Sigma}_{\scriptscriptstyle \blacksquare}}(((\id \! 
- \! C_{w^{\tilde{\Sigma}_{\scriptscriptstyle \blacksquare}}})^{-1}
C_{w^{\tilde{\Sigma}_{\scriptscriptstyle \blacksquare}}})(C_{
w^{\tilde{\Sigma}_{\scriptscriptstyle\blacksquare}}} \mathrm{I}))(\xi)
w_{+}^{\tilde{\Sigma}_{\scriptscriptstyle \blacksquare}}(\xi) \left(
\dfrac{1}{\xi \! - \! z} \! - \! \dfrac{1}{\xi \! - \! \alpha_{k}} \right) 
\dfrac{\md \xi}{2 \pi \mi}, \label{eqtazq4} \\
\tilde{\mathbb{J}}^{\tilde{\mathcal{R}}}_{4}(z) :=& \, \int_{\tilde{
\Sigma}_{\circlearrowright}}(((\id \! - \! C_{w^{\tilde{\Sigma}_{
\circlearrowright}}})^{-1}C_{w^{\tilde{\Sigma}_{\circlearrowright}}})
(C_{w^{\tilde{\Sigma}_{\circlearrowright}}} \mathrm{I}))(\xi)
w_{+}^{\tilde{\Sigma}_{\circlearrowright}}(\xi) \left(\dfrac{1}{\xi 
\! - \! z} \! - \! \dfrac{1}{\xi \! - \! \alpha_{k}} \right) \dfrac{\md 
\xi}{2 \pi \mi} \nonumber \\
+& \, \int_{\tilde{\Sigma}_{\scriptscriptstyle \blacksquare}}(((\id \! 
- \! C_{w^{\tilde{\Sigma}_{\circlearrowright}}})^{-1}C_{w^{\tilde{
\Sigma}_{\circlearrowright}}})(C_{w^{\tilde{\Sigma}_{\circlearrowright}}} 
\mathrm{I}))(\xi)w_{+}^{\tilde{\Sigma}_{\scriptscriptstyle \blacksquare}}
(\xi) \left(\dfrac{1}{\xi \! - \! z} \! - \! \dfrac{1}{\xi \! - \! \alpha_{k}} 
\right) \dfrac{\md \xi}{2 \pi \mi}, \label{eqtazq5} \\
\tilde{\mathbb{J}}^{\tilde{\mathcal{R}}}_{5}(z) :=& \, \int_{\tilde{
\Sigma}_{\circlearrowright}}(((\id \! - \! (\id \! - \! C_{w^{\tilde{
\Sigma}_{\circlearrowright}}})^{-1}(\id \! - \! C_{w^{\tilde{\Sigma}_{
\scriptscriptstyle \blacksquare}}})^{-1}C_{w^{\tilde{\Sigma}_{
\scriptscriptstyle \blacksquare}}}C_{w^{\tilde{\Sigma}_{
\circlearrowright}}})^{-1}(\id \! - \! C_{w^{\tilde{\Sigma}_{
\circlearrowright}}})^{-1} \nonumber \\
\times& \, (\id \! - \! C_{w^{\tilde{\Sigma}_{\scriptscriptstyle 
\blacksquare}}})^{-1}C_{w^{\tilde{\Sigma}_{\circlearrowright}}})
(C_{w^{\tilde{\Sigma}_{\scriptscriptstyle \blacksquare}}} \mathrm{I}))
(\xi)w_{+}^{\tilde{\Sigma}_{\circlearrowright}}(\xi) \left(\dfrac{1}{\xi 
\! - \! z} \! - \! \dfrac{1}{\xi \! - \! \alpha_{k}} \right) \dfrac{\md 
\xi}{2 \pi \mi} \nonumber \\
+& \, \int_{\tilde{\Sigma}_{\scriptscriptstyle \blacksquare}}(((\id \! 
- \! (\id \! - \! C_{w^{\tilde{\Sigma}_{\circlearrowright}}})^{-1}(\id 
\! - \! C_{w^{\tilde{\Sigma}_{\scriptscriptstyle \blacksquare}}})^{-1}
C_{w^{\tilde{\Sigma}_{\scriptscriptstyle \blacksquare}}}C_{w^{
\tilde{\Sigma}_{\circlearrowright}}})^{-1}(\id \! - \! C_{w^{\tilde{
\Sigma}_{\circlearrowright}}})^{-1} \nonumber \\
\times& \, (\id \! - \! C_{w^{\tilde{\Sigma}_{\scriptscriptstyle 
\blacksquare}}})^{-1}C_{w^{\tilde{\Sigma}_{\circlearrowright}}})
(C_{w^{\tilde{\Sigma}_{\scriptscriptstyle \blacksquare}}} \mathrm{I}))
(\xi)w_{+}^{\tilde{\Sigma}_{\scriptscriptstyle \blacksquare}}(\xi) 
\left(\dfrac{1}{\xi \! - \! z} \! - \! \dfrac{1}{\xi \! - \! \alpha_{k}} 
\right) \dfrac{\md \xi}{2 \pi \mi}, \label{eqtazq6} \\
\tilde{\mathbb{J}}^{\tilde{\mathcal{R}}}_{6}(z) :=& \, \int_{\tilde{
\Sigma}_{\circlearrowright}}(((\id \! - \! (\id \! - \! C_{w^{\tilde{
\Sigma}_{\scriptscriptstyle \blacksquare}}})^{-1}(\id \! - \! 
C_{w^{\tilde{\Sigma}_{\circlearrowright}}})^{-1}C_{w^{\tilde{
\Sigma}_{\circlearrowright}}}C_{w^{\tilde{\Sigma}_{\scriptscriptstyle 
\blacksquare}}})^{-1}(\id \! - \! C_{w^{\tilde{\Sigma}_{\scriptscriptstyle 
\blacksquare}}})^{-1} \nonumber \\
\times& \, (\id \! - \! C_{w^{\tilde{\Sigma}_{\circlearrowright}}})^{-1}
C_{w^{\tilde{\Sigma}_{\scriptscriptstyle \blacksquare}}})(C_{w^{\tilde{
\Sigma}_{\circlearrowright}}} \mathrm{I}))(\xi)w_{+}^{\tilde{\Sigma}_{
\circlearrowright}}(\xi) \left(\dfrac{1}{\xi \! - \! z} \! - \! \dfrac{1}{\xi 
\! - \! \alpha_{k}} \right) \dfrac{\md \xi}{2 \pi \mi} \nonumber \\
+& \, \int_{\tilde{\Sigma}_{\scriptscriptstyle \blacksquare}}(((\id \! 
- \! (\id \! - \! C_{w^{\tilde{\Sigma}_{\scriptscriptstyle \blacksquare}}}
)^{-1}(\id \! - \! C_{w^{\tilde{\Sigma}_{\circlearrowright}}})^{-1}
C_{w^{\tilde{\Sigma}_{\circlearrowright}}}C_{w^{\tilde{\Sigma}_{
\scriptscriptstyle \blacksquare}}})^{-1}(\id \! - \! C_{w^{\tilde{
\Sigma}_{\scriptscriptstyle \blacksquare}}})^{-1} \nonumber \\
\times& \, (\id \! - \! C_{w^{\tilde{\Sigma}_{\circlearrowright}}})^{-1}
C_{w^{\tilde{\Sigma}_{\scriptscriptstyle \blacksquare}}})(C_{w^{
\tilde{\Sigma}_{\circlearrowright}}} \mathrm{I}))(\xi)w_{+}^{\tilde{
\Sigma}_{\scriptscriptstyle \blacksquare}}(\xi) \left(\dfrac{1}{\xi 
\! - \! z} \! - \! \dfrac{1}{\xi \! - \! \alpha_{k}} \right) \dfrac{\md 
\xi}{2 \pi \mi}, \label{eqtazq7} \\
\tilde{\mathbb{J}}^{\tilde{\mathcal{R}}}_{7}(z) :=& \, \int_{\tilde{
\Sigma}_{\circlearrowright}}(((\id \! - \! (\id \! - \! C_{w^{\tilde{
\Sigma}_{\circlearrowright}}})^{-1}(\id \! - \! C_{w^{\tilde{\Sigma}_{
\scriptscriptstyle \blacksquare}}})^{-1}C_{w^{\tilde{\Sigma}_{
\scriptscriptstyle \blacksquare}}}C_{w^{\tilde{\Sigma}_{
\circlearrowright}}})^{-1}(\id \! - \! C_{w^{\tilde{\Sigma}_{
\circlearrowright}}})^{-1}(\id \! - \! C_{w^{\tilde{\Sigma}_{
\scriptscriptstyle \blacksquare}}})^{-1} \nonumber \\
\times& \, C_{w^{\tilde{\Sigma}_{\circlearrowright}}}(\id \! - \! 
C_{w^{\tilde{\Sigma}_{\scriptscriptstyle \blacksquare}}})^{-1}
C_{w^{\tilde{\Sigma}_{\scriptscriptstyle \blacksquare}}})
(C_{w^{\tilde{\Sigma}_{\scriptscriptstyle\blacksquare}}} \mathrm{I}))
(\xi)w_{+}^{\tilde{\Sigma}_{\circlearrowright}}(\xi) \left(\dfrac{1}{
\xi \! - \! z} \! - \! \dfrac{1}{\xi \! - \! \alpha_{k}} \right) 
\dfrac{\md \xi}{2 \pi \mi} \nonumber \\
+& \, \int_{\tilde{\Sigma}_{\scriptscriptstyle \blacksquare}}(((\id 
\! - \! (\id \! - \! C_{w^{\tilde{\Sigma}_{\circlearrowright}}})^{-1}(\id 
\! - \! C_{w^{\tilde{\Sigma}_{\scriptscriptstyle \blacksquare}}})^{-1}
C_{w^{\tilde{\Sigma}_{\scriptscriptstyle \blacksquare}}}C_{w^{
\tilde{\Sigma}_{\circlearrowright}}})^{-1}(\id \! - \! C_{w^{\tilde{
\Sigma}_{\circlearrowright}}})^{-1}(\id \! - \! C_{w^{\tilde{\Sigma}_{
\scriptscriptstyle \blacksquare}}})^{-1} \nonumber \\
\times& \, C_{w^{\tilde{\Sigma}_{\circlearrowright}}}(\id \! - \! 
C_{w^{\tilde{\Sigma}_{\scriptscriptstyle \blacksquare}}})^{-1}
C_{w^{\tilde{\Sigma}_{\scriptscriptstyle \blacksquare}}})
(C_{w^{\tilde{\Sigma}_{\scriptscriptstyle \blacksquare}}} \mathrm{I}))
(\xi)w_{+}^{\tilde{\Sigma}_{\scriptscriptstyle \blacksquare}}(\xi) 
\left(\dfrac{1}{\xi \! - \! z} \! - \! \dfrac{1}{\xi \! - \! \alpha_{k}} 
\right) \dfrac{\md \xi}{2 \pi \mi}, \label{eqtazq8} \\
\tilde{\mathbb{J}}^{\tilde{\mathcal{R}}}_{8}(z) :=& \, \int_{\tilde{
\Sigma}_{\circlearrowright}}(((\id \! - \! (\id \! - \! C_{w^{\tilde{
\Sigma}_{\scriptscriptstyle \blacksquare}}})^{-1}(\id \! - \! 
C_{w^{\tilde{\Sigma}_{\circlearrowright}}})^{-1}C_{w^{\tilde{
\Sigma}_{\circlearrowright}}}C_{w^{\tilde{\Sigma}_{
\scriptscriptstyle \blacksquare}}})^{-1}(\id \! - \! 
C_{w^{\tilde{\Sigma}_{\scriptscriptstyle \blacksquare}}})^{-1}
(\id \! - \! C_{w^{\tilde{\Sigma}_{\circlearrowright}}})^{-1} 
\nonumber \\
\times& \, C_{w^{\tilde{\Sigma}_{\scriptscriptstyle \blacksquare}}}
(\id \! - \! C_{w^{\tilde{\Sigma}_{\circlearrowright}}})^{-1}C_{w^{
\tilde{\Sigma}_{\circlearrowright}}})(C_{w^{\tilde{\Sigma}_{
\circlearrowright}}} \mathrm{I}))(\xi)w_{+}^{\tilde{\Sigma}_{
\circlearrowright}}(\xi) \left(\dfrac{1}{\xi \! - \! z} \! - \! 
\dfrac{1}{\xi \! - \! \alpha_{k}} \right) \dfrac{\md \xi}{2 \pi \mi} 
\nonumber \\
+& \, \int_{\tilde{\Sigma}_{\scriptscriptstyle \blacksquare}}
(((\id \! - \! (\id \! - \! C_{w^{\tilde{\Sigma}_{\scriptscriptstyle 
\blacksquare}}})^{-1}(\id \! - \! C_{w^{\tilde{\Sigma}_{
\circlearrowright}}})^{-1}C_{w^{\tilde{\Sigma}_{\circlearrowright}}}
C_{w^{\tilde{\Sigma}_{\scriptscriptstyle \blacksquare}}})^{-1}
(\id \! - \! C_{w^{\tilde{\Sigma}_{\scriptscriptstyle \blacksquare}}}
)^{-1}(\id \! - \! C_{w^{\tilde{\Sigma}_{\circlearrowright}}})^{-1} 
\nonumber \\
\times& \, C_{w^{\tilde{\Sigma}_{\scriptscriptstyle \blacksquare}}}
(\id \! - \! C_{w^{\tilde{\Sigma}_{\circlearrowright}}})^{-1}C_{w^{
\tilde{\Sigma}_{\circlearrowright}}})(C_{w^{\tilde{\Sigma}_{
\circlearrowright}}} \mathrm{I}))(\xi)w_{+}^{\tilde{\Sigma}_{
\scriptscriptstyle \blacksquare}}(\xi) \left(\dfrac{1}{\xi \! - \! z} 
\! - \! \dfrac{1}{\xi \! - \! \alpha_{k}} \right) \dfrac{\md \xi}{2 
\pi \mi}. \label{eqtazq9}
\end{align}
One now proceeds to estimate asymptotically, in the double-scaling 
limit $\mathscr{N},n \! \to \! \infty$ such that $z_{o} \! = \! 1 \! + \! 
o(1)$, the terms on the right-hand side of Equation~\eqref{eqtazq1}.

An inequality argument shows that, for $z \! \in \! \mathbb{C}
\setminus \tilde{\Sigma}_{\tilde{\mathcal{R}}}^{\sharp}$,
\begin{equation*}
\left\vert \int_{\tilde{\Sigma}_{\scriptscriptstyle \blacksquare}}
w_{+}^{\tilde{\Sigma}_{\scriptscriptstyle \blacksquare}}(\xi) \left(
\dfrac{1}{\xi \! - \! z} \! - \! \dfrac{1}{\xi \! - \! \alpha_{k}} \right) 
\dfrac{\md \xi}{2 \pi \mi} \right\vert \underset{\underset{z_{o}=1
+o(1)}{\mathscr{N},n \to \infty}}{\leqslant} \left\lvert \left\lvert 
\dfrac{w_{+}^{\tilde{\Sigma}_{\scriptscriptstyle \blacksquare}}
(\pmb{\cdot})}{2 \pi \mi (\pmb{\cdot} \! - \! \alpha_{k})} \right\rvert 
\right\rvert_{\mathcal{L}^{1}_{\mathrm{M}_{2}(\mathbb{C})}(\tilde{
\Sigma}_{\scriptscriptstyle \blacksquare})} \! + \! \left\lvert \left\lvert 
\dfrac{w_{+}^{\tilde{\Sigma}_{\scriptscriptstyle \blacksquare}}
(\pmb{\cdot})}{2 \pi \mi (\pmb{\cdot} \! - \! z)} \right\rvert 
\right\rvert_{\mathcal{L}^{1}_{\mathrm{M}_{2}(\mathbb{C})}
(\tilde{\Sigma}_{\scriptscriptstyle \blacksquare})}:
\end{equation*}
proceeding as in the calculations leading to the 
Estimates~\eqref{eqiy38} and~\eqref{eqiy71} (cf. the proof of 
Lemma~\ref{lem5.3}), one shows that
\begin{gather}
\left\lvert \left\lvert \dfrac{w_{+}^{\tilde{\Sigma}_{\scriptscriptstyle 
\blacksquare}}(\pmb{\cdot})}{2 \pi \mi (\pmb{\cdot} \! - \! \alpha_{
k})} \right\rvert \right\rvert_{\mathcal{L}^{1}_{\mathrm{M}_{2}
(\mathbb{C})}(\tilde{\Sigma}_{\scriptscriptstyle \blacksquare})} 
\underset{\underset{z_{o}=1+o(1)}{\mathscr{N},n \to \infty}}{\leqslant} 
\mathcal{O} \left(\dfrac{\tilde{\mathfrak{c}}_{\tilde{\Sigma}}^{1}
(n,k,z_{o}) \me^{-((n-1)K+k) \tilde{\lambda}_{\tilde{\mathcal{R}},
\tilde{w}}^{\triangleright}}}{(n \! - \! 1)K \! + \! k} \right), 
\label{eqtazq10} \\
\left\lvert \left\lvert \dfrac{w_{+}^{\tilde{\Sigma}_{\scriptscriptstyle 
\blacksquare}}(\pmb{\cdot})}{2 \pi \mi (\pmb{\cdot} \! - \! z)} 
\right\rvert \right\rvert_{\mathcal{L}^{1}_{\mathrm{M}_{2}
(\mathbb{C})}(\tilde{\Sigma}_{\scriptscriptstyle \blacksquare})} 
\underset{\underset{z_{o}=1+o(1)}{\mathscr{N},n \to \infty}}{\leqslant} 
\mathcal{O} \left(\dfrac{\tilde{\mathfrak{c}}_{\tilde{\Sigma}}^{2}
(n,k,z_{o}) \me^{-((n-1)K+k) \tilde{\lambda}_{\tilde{\mathcal{R}},
\tilde{w}}^{\triangleright}}}{((n \! - \! 1)K \! + \! k) \operatorname{dist}
(\tilde{\Sigma}_{\tilde{\mathcal{R}}}^{\sharp},z)} \right), \quad z 
\! \in \! \mathbb{C} \setminus \tilde{\Sigma}_{\tilde{\mathcal{
R}}}^{\sharp}, \label{eqtazq11} 
\end{gather}
where $\tilde{\lambda}_{\tilde{\mathcal{R}},\tilde{w}}^{\triangleright}$ 
$(> \! 0)$ is defined by Equation~\eqref{tillamrwr2}, and $\tilde{
\mathfrak{c}}_{\tilde{\Sigma}}^{r}(n,k,z_{o}) \! =_{\underset{z_{o}=
1+o(1)}{\mathscr{N},n \to \infty}} \! \mathcal{O}(1)$, $r \! \in \! 
\lbrace 1,2 \rbrace$, whence, via the Estimates~\eqref{eqtazq10} 
and~\eqref{eqtazq11}, for $n \! \in \! \mathbb{N}$ and $k \! \in \! 
\lbrace 1,2,\dotsc,K \rbrace$ such that $\alpha_{p_{\mathfrak{s}}} \! 
:= \! \alpha_{k} \! \neq \! \infty$,
\begin{equation} \label{eqtazq12} 
\left\vert \int_{\tilde{\Sigma}_{\scriptscriptstyle \blacksquare}}
w_{+}^{\tilde{\Sigma}_{\scriptscriptstyle \blacksquare}}(\xi) \left(
\dfrac{1}{\xi \! - \! z} \! - \! \dfrac{1}{\xi \! - \! \alpha_{k}} \right) 
\dfrac{\md \xi}{2 \pi \mi} \right\vert \underset{\underset{z_{o}=1
+o(1)}{\mathscr{N},n \to \infty}}{\leqslant} \mathcal{O} \left(
\dfrac{\tilde{\mathfrak{c}}_{\tilde{\Sigma}}^{3}(n,k,z_{o}) 
\me^{-((n-1)K+k) \tilde{\lambda}_{\tilde{\mathcal{R}},
\tilde{w}}^{\triangleright}}}{((n \! - \! 1)K \! + \! k) \min \lbrace 1,
\operatorname{dist}(\tilde{\Sigma}_{\tilde{\mathcal{R}}}^{\sharp},
z) \rbrace} \right), \quad z \! \in \! \mathbb{C} \setminus 
\tilde{\Sigma}_{\tilde{\mathcal{R}}}^{\sharp},
\end{equation}
where $\tilde{\mathfrak{c}}_{\tilde{\Sigma}}^{3}(n,k,z_{o}) \! 
=_{\underset{z_{o}=1+o(1)}{\mathscr{N},n \to \infty}} \! 
\mathcal{O}(1)$.

An argument based on H\"{o}lder's Inequality and H\"{o}lder's Inequality 
for Integrals shows that,\footnote{For $\mathrm{F} \! \in \! 
\mathcal{L}^{2}_{\mathrm{M}_{2}(\mathbb{C})}(\mathrm{A})$ and 
$\mathrm{G} \! \in \! \mathcal{L}^{2}_{\mathrm{M}_{2}(\mathbb{C})}
(\mathrm{A})$, where $\mathrm{A} \subseteq \overline{\mathbb{C}}$, 
$\left\lvert \int_{\mathrm{A}} \mathrm{F}(\xi) \mathrm{G}(\xi) \, \md 
\xi \right\rvert \! \leqslant \! \lvert \lvert \mathrm{F}(\pmb{\cdot}) 
\rvert \rvert_{\mathcal{L}^{2}_{\mathrm{M}_{2}(\mathbb{C})}(\mathrm{A})} 
\lvert \lvert \mathrm{G}(\pmb{\cdot}) \rvert \rvert_{\mathcal{L}^{2}_{
\mathrm{M}_{2}(\mathbb{C})}(\mathrm{A})}$.} for $z \! \in \! \mathbb{C} 
\setminus \tilde{\Sigma}_{\tilde{\mathcal{R}}}^{\sharp}$,
\begin{align}
\left\lvert \tilde{\mathbb{J}}_{1}^{\tilde{\mathcal{R}}}(z) \right\rvert 
\underset{\underset{z_{o}=1+o(1)}{\mathscr{N},n \to \infty}}{
\leqslant}& \, \left\lvert \left\lvert (C_{w^{\tilde{\Sigma}_{
\scriptscriptstyle \blacksquare}}} \mathrm{I})(\pmb{\cdot}) 
\right\rvert \right\rvert_{\mathcal{L}^{2}_{\mathrm{M}_{2}
(\mathbb{C})}(\tilde{\Sigma}_{\scriptscriptstyle \blacksquare})} 
\left(\left\lvert \left\lvert \dfrac{w_{+}^{\tilde{\Sigma}_{
\scriptscriptstyle \blacksquare}}(\pmb{\cdot})}{2 \pi \mi 
(\pmb{\cdot} \! - \! \alpha_{k})} \right\rvert \right\rvert_{
\mathcal{L}^{2}_{\mathrm{M}_{2}(\mathbb{C})}(\tilde{\Sigma}_{
\scriptscriptstyle \blacksquare})} \! + \! \left\lvert \left\lvert 
\dfrac{w_{+}^{\tilde{\Sigma}_{\scriptscriptstyle \blacksquare}}
(\pmb{\cdot})}{2 \pi \mi (\pmb{\cdot} \! - \! z)} \right\rvert 
\right\rvert_{\mathcal{L}^{2}_{\mathrm{M}_{2}(\mathbb{C})}
(\tilde{\Sigma}_{\scriptscriptstyle \blacksquare})} \right) 
\nonumber \\
+& \, \left\lvert \left\lvert (C_{w^{\tilde{\Sigma}_{
\scriptscriptstyle \blacksquare}}} \mathrm{I})(\pmb{\cdot}) 
\right\rvert \right\rvert_{\mathcal{L}^{2}_{\mathrm{M}_{2}
(\mathbb{C})}(\tilde{\Sigma}_{\circlearrowright})} \left(
\left\lvert \left\lvert \dfrac{w_{+}^{\tilde{\Sigma}_{
\circlearrowright}}(\pmb{\cdot})}{2 \pi \mi (\pmb{\cdot} \! - 
\! \alpha_{k})} \right\rvert \right\rvert_{\mathcal{L}^{2}_{
\mathrm{M}_{2}(\mathbb{C})}(\tilde{\Sigma}_{\circlearrowright})} 
\! + \! \left\lvert \left\lvert \dfrac{w_{+}^{\tilde{\Sigma}_{
\circlearrowright}}(\pmb{\cdot})}{2 \pi \mi (\pmb{\cdot} \! - \! 
z)} \right\rvert \right\rvert_{\mathcal{L}^{2}_{\mathrm{M}_{2}
(\mathbb{C})}(\tilde{\Sigma}_{\circlearrowright})} \right) \nonumber \\
\underset{\underset{z_{o}=1+o(1)}{\mathscr{N},n \to \infty}}{
\leqslant}& \, \left\lvert \left\lvert C_{w^{\Sigma_{\tilde{
\mathcal{R}}}}} \right\rvert \right\rvert_{\mathfrak{B}_{2}
(\tilde{\Sigma}_{\tilde{\mathcal{R}}}^{\sharp})} \left(\left\lvert 
\left\lvert \dfrac{w_{+}^{\tilde{\Sigma}_{\scriptscriptstyle 
\blacksquare}}(\pmb{\cdot})}{2 \pi \mi (\pmb{\cdot} \! - \! 
\alpha_{k})} \right\rvert \right\rvert_{\mathcal{L}^{2}_{
\mathrm{M}_{2}(\mathbb{C})}(\tilde{\Sigma}_{\scriptscriptstyle 
\blacksquare})} \! + \! \left\lvert \left\lvert \dfrac{w_{+}^{
\tilde{\Sigma}_{\scriptscriptstyle \blacksquare}}(\pmb{\cdot})}{
2 \pi \mi (\pmb{\cdot} \! - \! z)} \right\rvert \right\rvert_{
\mathcal{L}^{2}_{\mathrm{M}_{2}(\mathbb{C})}(\tilde{\Sigma}_{
\scriptscriptstyle \blacksquare})} \right. \nonumber \\
+&\left. \,  \left\lvert \left\lvert \dfrac{w_{+}^{\tilde{\Sigma}_{
\circlearrowright}}(\pmb{\cdot})}{2 \pi \mi (\pmb{\cdot} \! - \! 
\alpha_{k})} \right\rvert \right\rvert_{\mathcal{L}^{2}_{\mathrm{
M}_{2}(\mathbb{C})}(\tilde{\Sigma}_{\circlearrowright})} \! 
+ \! \left\lvert \left\lvert \dfrac{w_{+}^{\tilde{\Sigma}_{
\circlearrowright}}(\pmb{\cdot})}{2 \pi \mi (\pmb{\cdot} \! - \! 
z)} \right\rvert \right\rvert_{\mathcal{L}^{2}_{\mathrm{M}_{2}
(\mathbb{C})}(\tilde{\Sigma}_{\circlearrowright})} \right): 
\label{eqtazq13}
\end{align}
proceeding as in the calculations leading to the 
Estimates~\eqref{eqiy38} and~\eqref{eqiy71} (cf. the 
proof of Lemma~\ref{lem5.3}), one shows that
\begin{gather}
\left\lvert \left\lvert \dfrac{w_{+}^{\tilde{\Sigma}_{\scriptscriptstyle 
\blacksquare}}(\pmb{\cdot})}{2 \pi \mi (\pmb{\cdot} \! - \! \alpha_{
k})} \right\rvert \right\rvert_{\mathcal{L}^{2}_{\mathrm{M}_{2}
(\mathbb{C})}(\tilde{\Sigma}_{\scriptscriptstyle \blacksquare})} 
\underset{\underset{z_{o}=1+o(1)}{\mathscr{N},n \to \infty}}{\leqslant} 
\mathcal{O} \left(\dfrac{\tilde{\mathfrak{c}}_{\tilde{\Sigma}}^{4}
(n,k,z_{o}) \me^{-((n-1)K+k) \tilde{\lambda}_{\tilde{\mathcal{R}},
\tilde{w}}^{\triangleright}}}{((n \! - \! 1)K \! + \! k)^{1/2}} \right), 
\label{eqtazq14} \\
\left\lvert \left\lvert \dfrac{w_{+}^{\tilde{\Sigma}_{\scriptscriptstyle 
\blacksquare}}(\pmb{\cdot})}{2 \pi \mi (\pmb{\cdot} \! - \! z)} 
\right\rvert \right\rvert_{\mathcal{L}^{2}_{\mathrm{M}_{2}
(\mathbb{C})}(\tilde{\Sigma}_{\scriptscriptstyle \blacksquare})} 
\underset{\underset{z_{o}=1+o(1)}{\mathscr{N},n \to \infty}}{\leqslant} 
\mathcal{O} \left(\dfrac{\tilde{\mathfrak{c}}_{\tilde{\Sigma}}^{5}
(n,k,z_{o}) \me^{-((n-1)K+k) \tilde{\lambda}_{\tilde{\mathcal{R}},
\tilde{w}}^{\triangleright}}}{((n \! - \! 1)K \! + \! k)^{1/2} 
\operatorname{dist}(\tilde{\Sigma}_{\tilde{\mathcal{R}}}^{\sharp},z)} 
\right), \quad z \! \in \! \mathbb{C} \setminus \tilde{\Sigma}_{
\tilde{\mathcal{R}}}^{\sharp}, \label{eqtazq15} \\
\left\lvert \left\lvert \dfrac{w_{+}^{\tilde{\Sigma}_{\circlearrowright}}
(\pmb{\cdot})}{2 \pi \mi (\pmb{\cdot} \! - \! \alpha_{k})} \right\rvert 
\right\rvert_{\mathcal{L}^{2}_{\mathrm{M}_{2}(\mathbb{C})}(\tilde{
\Sigma}_{\circlearrowright})} \underset{\underset{z_{o}=1+o(1)}{
\mathscr{N},n \to \infty}}{=} 0, \label{eqtazq16} \\
\left\lvert \left\lvert \dfrac{w_{+}^{\tilde{\Sigma}_{\circlearrowright}}
(\pmb{\cdot})}{2 \pi \mi (\pmb{\cdot} \! - \! z)} \right\rvert 
\right\rvert_{\mathcal{L}^{2}_{\mathrm{M}_{2}(\mathbb{C})}
(\tilde{\Sigma}_{\circlearrowright})} \underset{\underset{z_{o}=1
+o(1)}{\mathscr{N},n \to \infty}}{\leqslant} \mathcal{O} \left(
\dfrac{\tilde{\mathfrak{c}}_{\tilde{\Sigma}}^{6}(n,k,z_{o})}{((n \! - \! 1)
K \! + \! k) \operatorname{dist}(\tilde{\Sigma}_{\tilde{\mathcal{R}}}^{
\sharp},z)} \right), \quad z \! \in \! \mathbb{C} \setminus 
\tilde{\Sigma}_{\tilde{\mathcal{R}}}^{\sharp}, \label{eqtazq17}
\end{gather}
where $\tilde{\mathfrak{c}}_{\tilde{\Sigma}}^{r}(n,k,z_{o}) \! 
=_{\underset{z_{o}=1+o(1)}{\mathscr{N},n \to \infty}} \! 
\mathcal{O}(1)$, $r \! \in \! \lbrace 4,5,6 \rbrace$, 
whence, via the Estimate~\eqref{eqcwsigtll}, the 
Estimates~\eqref{eqtazq14}--\eqref{eqtazq17}, and the 
Inequality~\eqref{eqtazq13}, for $n \! \in \! \mathbb{N}$ and 
$k \! \in \! \lbrace 1,2,\dotsc,K \rbrace$ such that $\alpha_{p_{
\mathfrak{s}}} \! := \! \alpha_{k} \! \neq \! \infty$,
\begin{equation} \label{eqtazq18} 
\left\vert \tilde{\mathbb{J}}^{\tilde{\mathcal{R}}}_{1}(z) \right\rvert 
\underset{\underset{z_{o}=1+o(1)}{\mathscr{N},n \to \infty}}{
\leqslant} \mathcal{O} \left(\dfrac{\tilde{\mathfrak{c}}_{\tilde{
\Sigma}}^{7}(n,k,z_{o}) \me^{-\frac{1}{2}((n-1)K+k) \tilde{
\lambda}_{\tilde{\mathcal{R}},\tilde{w}}^{\triangleright}}}{((n \! 
- \! 1)K \! + \! k) \min \lbrace 1,\operatorname{dist}(\tilde{
\Sigma}_{\tilde{\mathcal{R}}}^{\sharp},z) \rbrace} \right), 
\quad z \! \in \! \mathbb{C} \setminus \tilde{\Sigma}_{
\tilde{\mathcal{R}}}^{\sharp},
\end{equation}
where $\tilde{\mathfrak{c}}_{\tilde{\Sigma}}^{7}(n,k,z_{o}) \! 
=_{\underset{z_{o}=1+o(1)}{\mathscr{N},n \to \infty}} \! 
\mathcal{O}(1)$. Proceeding analogously as in the calculations 
above leading to the Estimate~\eqref{eqtazq18}, one shows that, 
for $n \! \in \! \mathbb{N}$ and $k \! \in \! \lbrace 1,2,\dotsc,K 
\rbrace$ such that $\alpha_{p_{\mathfrak{s}}} \! := \! \alpha_{k} 
\! \neq \! \infty$,
\begin{equation} \label{eqtazq19} 
\left\vert \tilde{\mathbb{J}}^{\tilde{\mathcal{R}}}_{2}(z) \right\rvert 
\underset{\underset{z_{o}=1+o(1)}{\mathscr{N},n \to \infty}}{
\leqslant} \mathcal{O} \left(\dfrac{\tilde{\mathfrak{c}}_{\tilde{
\Sigma}}^{8}(n,k,z_{o}) \me^{-\frac{1}{2}((n-1)K+k) \tilde{
\lambda}_{\tilde{\mathcal{R}},\tilde{w}}^{\triangleright}}}{((n \! 
- \! 1)K \! + \! k) \min \lbrace 1,\operatorname{dist}(\tilde{
\Sigma}_{\tilde{\mathcal{R}}}^{\sharp},z) \rbrace} \right), 
\quad z \! \in \! \mathbb{C} \setminus \tilde{\Sigma}_{
\tilde{\mathcal{R}}}^{\sharp},
\end{equation}
where $\tilde{\mathfrak{c}}_{\tilde{\Sigma}}^{8}(n,k,z_{o}) \! 
=_{\underset{z_{o}=1+o(1)}{\mathscr{N},n \to \infty}} \! 
\mathcal{O}(1)$.

An argument based on H\"{o}lder's Inequality, H\"{o}lder's 
Inequality for Integrals, and a Neumann series inversion (cf. 
Estimate~\eqref{eqcwsigtll}) shows that, for $z \! \in \! \mathbb{C} 
\setminus \tilde{\Sigma}_{\tilde{\mathcal{R}}}^{\sharp}$,
\begin{align}
\left\lvert \tilde{\mathbb{J}}_{3}^{\tilde{\mathcal{R}}}(z) \right\rvert 
\underset{\underset{z_{o}=1+o(1)}{\mathscr{N},n \to \infty}}{
\leqslant}& \, \left\lvert \left\lvert (\id \! - \! C_{w^{\tilde{
\Sigma}_{\scriptscriptstyle \blacksquare}}})^{-1} \right\rvert 
\right\rvert_{\mathfrak{B}_{2}(\tilde{\Sigma}_{\scriptscriptstyle 
\blacksquare})} \left(\left\lvert \left\lvert C_{w^{\tilde{\Sigma}_{
\scriptscriptstyle \blacksquare}}} \right\rvert \right\rvert_{
\mathfrak{B}_{2}(\tilde{\Sigma}_{\scriptscriptstyle 
\blacksquare})} \right)^{2} \left(\left\lvert \left\lvert 
\dfrac{w_{+}^{\tilde{\Sigma}_{\scriptscriptstyle \blacksquare}}
(\pmb{\cdot})}{2 \pi \mi (\pmb{\cdot} \! - \! \alpha_{k})} 
\right\rvert \right\rvert_{\mathcal{L}^{2}_{\mathrm{M}_{2}
(\mathbb{C})}(\tilde{\Sigma}_{\scriptscriptstyle \blacksquare})} 
\! + \! \left\lvert \left\lvert \dfrac{w_{+}^{\tilde{\Sigma}_{
\scriptscriptstyle \blacksquare}}(\pmb{\cdot})}{2 \pi \mi 
(\pmb{\cdot} \! - \! z)} \right\rvert \right\rvert_{\mathcal{
L}^{2}_{\mathrm{M}_{2}(\mathbb{C})}(\tilde{\Sigma}_{
\scriptscriptstyle \blacksquare})} \right) \nonumber \\
+& \, \left\lvert \left\lvert (\id \! - \! C_{w^{\tilde{\Sigma}_{
\scriptscriptstyle \blacksquare}}})^{-1} \right\rvert \right\rvert_{
\mathfrak{B}_{2}(\tilde{\Sigma}_{\circlearrowright})} \left(
\left\lvert \left\lvert C_{w^{\tilde{\Sigma}_{\scriptscriptstyle 
\blacksquare}}} \right\rvert \right\rvert_{\mathfrak{B}_{2}
(\tilde{\Sigma}_{\circlearrowright})} \right)^{2} \left(
\left\lvert \left\lvert \dfrac{w_{+}^{\tilde{\Sigma}_{
\circlearrowright}}(\pmb{\cdot})}{2 \pi \mi (\pmb{\cdot} \! - 
\! \alpha_{k})} \right\rvert \right\rvert_{\mathcal{L}^{2}_{
\mathrm{M}_{2}(\mathbb{C})}(\tilde{\Sigma}_{\circlearrowright})} 
\! + \! \left\lvert \left\lvert \dfrac{w_{+}^{\tilde{\Sigma}_{
\circlearrowright}}(\pmb{\cdot})}{2 \pi \mi (\pmb{\cdot} \! - \! 
z)} \right\rvert \right\rvert_{\mathcal{L}^{2}_{\mathrm{M}_{2}
(\mathbb{C})}(\tilde{\Sigma}_{\circlearrowright})} \right) \nonumber \\
\underset{\underset{z_{o}=1+o(1)}{\mathscr{N},n \to \infty}}{
\leqslant}& \, \dfrac{\left(\left\lvert \left\lvert C_{w^{\Sigma_{
\tilde{\mathcal{R}}}}} \right\rvert \right\rvert_{\mathfrak{B}_{2}
(\tilde{\Sigma}_{\tilde{\mathcal{R}}}^{\sharp})} \right)^{2}}{\left(1 \! - \! 
\left\lvert \left\lvert C_{w^{\Sigma_{\tilde{\mathcal{R}}}}} \right\rvert 
\right\rvert_{\mathfrak{B}_{2}(\tilde{\Sigma}_{\tilde{\mathcal{R}}}^{
\sharp})} \right)} \left(\left\lvert \left\lvert \dfrac{w_{+}^{\tilde{\Sigma}_{
\scriptscriptstyle \blacksquare}}(\pmb{\cdot})}{2 \pi \mi (\pmb{
\cdot} \! - \! \alpha_{k})} \right\rvert \right\rvert_{\mathcal{L}^{2}_{
\mathrm{M}_{2}(\mathbb{C})}(\tilde{\Sigma}_{\scriptscriptstyle 
\blacksquare})} \! + \! \left\lvert \left\lvert \dfrac{w_{+}^{
\tilde{\Sigma}_{\scriptscriptstyle \blacksquare}}(\pmb{\cdot})}{
2 \pi \mi (\pmb{\cdot} \! - \! z)} \right\rvert \right\rvert_{
\mathcal{L}^{2}_{\mathrm{M}_{2}(\mathbb{C})}(\tilde{\Sigma}_{
\scriptscriptstyle \blacksquare})} \right. \nonumber \\
+&\left. \,  \left\lvert \left\lvert \dfrac{w_{+}^{\tilde{\Sigma}_{
\circlearrowright}}(\pmb{\cdot})}{2 \pi \mi (\pmb{\cdot} \! - \! 
\alpha_{k})} \right\rvert \right\rvert_{\mathcal{L}^{2}_{\mathrm{
M}_{2}(\mathbb{C})}(\tilde{\Sigma}_{\circlearrowright})} \! 
+ \! \left\lvert \left\lvert \dfrac{w_{+}^{\tilde{\Sigma}_{
\circlearrowright}}(\pmb{\cdot})}{2 \pi \mi (\pmb{\cdot} \! - \! 
z)} \right\rvert \right\rvert_{\mathcal{L}^{2}_{\mathrm{M}_{2}
(\mathbb{C})}(\tilde{\Sigma}_{\circlearrowright})} \right): 
\label{eqtazq20}
\end{align}
via the Estimate~\eqref{eqcwsigtll}, the 
Estimates~\eqref{eqtazq14}--\eqref{eqtazq17}, and the 
Inequality~\eqref{eqtazq20}, one arrives at, for $n \! \in \! 
\mathbb{N}$ and $k \! \in \! \lbrace 1,2,\dotsc,K \rbrace$ such that 
$\alpha_{p_{\mathfrak{s}}} \! := \! \alpha_{k} \! \neq \! \infty$,
\begin{equation} \label{eqtazq21} 
\left\vert \tilde{\mathbb{J}}^{\tilde{\mathcal{R}}}_{3}(z) \right\rvert 
\underset{\underset{z_{o}=1+o(1)}{\mathscr{N},n \to \infty}}{
\leqslant} \mathcal{O} \left(\dfrac{\tilde{\mathfrak{c}}_{\tilde{
\Sigma}}^{9}(n,k,z_{o}) \me^{-((n-1)K+k) \tilde{\lambda}_{\tilde{
\mathcal{R}},\tilde{w}}^{\triangleright}}}{((n \! - \! 1)K \! + \! k) \min 
\lbrace 1,\operatorname{dist}(\tilde{\Sigma}_{\tilde{\mathcal{R}}}^{
\sharp},z) \rbrace} \right), \quad z \! \in \! \mathbb{C} \setminus 
\tilde{\Sigma}_{\tilde{\mathcal{R}}}^{\sharp},
\end{equation}
where $\tilde{\mathfrak{c}}_{\tilde{\Sigma}}^{9}(n,k,z_{o}) \! 
=_{\underset{z_{o}=1+o(1)}{\mathscr{N},n \to \infty}} \! 
\mathcal{O}(1)$. Proceeding analogously as in the calculations 
above leading to the Estimate~\eqref{eqtazq21}, one shows that, 
for $n \! \in \! \mathbb{N}$ and $k \! \in \! \lbrace 1,2,\dotsc,K 
\rbrace$ such that $\alpha_{p_{\mathfrak{s}}} \! := \! \alpha_{k} 
\! \neq \! \infty$,
\begin{equation} \label{eqtazq22} 
\left\vert \tilde{\mathbb{J}}^{\tilde{\mathcal{R}}}_{4}(z) \right\rvert 
\underset{\underset{z_{o}=1+o(1)}{\mathscr{N},n \to \infty}}{
\leqslant} \mathcal{O} \left(\dfrac{\tilde{\mathfrak{c}}_{\tilde{
\Sigma}}^{10}(n,k,z_{o}) \me^{-((n-1)K+k) \tilde{\lambda}_{
\tilde{\mathcal{R}},\tilde{w}}^{\triangleright}}}{((n \! - \! 1)K \! + \! k) 
\min \lbrace 1,\operatorname{dist}(\tilde{\Sigma}_{\tilde{\mathcal{
R}}}^{\sharp},z) \rbrace} \right), \quad z \! \in \! \mathbb{C} 
\setminus \tilde{\Sigma}_{\tilde{\mathcal{R}}}^{\sharp},
\end{equation}
where $\tilde{\mathfrak{c}}_{\tilde{\Sigma}}^{10}(n,k,z_{o}) \! 
=_{\underset{z_{o}=1+o(1)}{\mathscr{N},n \to \infty}} \! 
\mathcal{O}(1)$.

An argument based on H\"{o}lder's Inequality, H\"{o}lder's 
Inequality for Integrals, and a Neumann series inversion (cf. 
Estimate~\eqref{eqcwsigtll}) shows that, for $z \! \in \! \mathbb{C} 
\setminus \tilde{\Sigma}_{\tilde{\mathcal{R}}}^{\sharp}$,
\begin{align}
\left\lvert \tilde{\mathbb{J}}_{5}^{\tilde{\mathcal{R}}}(z) \right\rvert 
\underset{\underset{z_{o}=1+o(1)}{\mathscr{N},n \to \infty}}{
\leqslant}& \, \left\lvert \left\lvert (\id \! - \! (\id \! - \! 
C_{w^{\tilde{\Sigma}_{\circlearrowright}}})^{-1}(\id \! - \! 
C_{w^{\tilde{\Sigma}_{\scriptscriptstyle \blacksquare}}})^{-1}
C_{w^{\tilde{\Sigma}_{\scriptscriptstyle \blacksquare}}}
C_{w^{\tilde{\Sigma}_{\circlearrowright}}})^{-1} \right\rvert 
\right\rvert_{\mathfrak{B}_{2}(\tilde{\Sigma}_{\scriptscriptstyle 
\blacksquare})} \left\lvert \left\lvert (\id \! - \! C_{w^{\tilde{
\Sigma}_{\circlearrowright}}})^{-1} \right\rvert \right\rvert_{
\mathfrak{B}_{2}(\tilde{\Sigma}_{\scriptscriptstyle \blacksquare})} 
\left\lvert \left\lvert (\id \! - \! C_{w^{\tilde{\Sigma}_{\scriptscriptstyle 
\blacksquare}}})^{-1} \right\rvert \right\rvert_{\mathfrak{B}_{2}
(\tilde{\Sigma}_{\scriptscriptstyle \blacksquare})} \nonumber \\
\times& \, \left\lvert \left\lvert C_{w^{\tilde{\Sigma}_{\circlearrowright}}} 
\right\rvert \right\rvert_{\mathfrak{B}_{2}(\tilde{\Sigma}_{\scriptscriptstyle 
\blacksquare})} \left\lvert \left\lvert C_{w^{\tilde{\Sigma}_{\scriptscriptstyle 
\blacksquare}}} \right\rvert \right\rvert_{\mathfrak{B}_{2}(\tilde{
\Sigma}_{\scriptscriptstyle \blacksquare})} \left(\left\lvert \left\lvert 
\dfrac{w_{+}^{\tilde{\Sigma}_{\scriptscriptstyle \blacksquare}}
(\pmb{\cdot})}{2 \pi \mi (\pmb{\cdot} \! - \! \alpha_{k})} 
\right\rvert \right\rvert_{\mathcal{L}^{2}_{\mathrm{M}_{2}
(\mathbb{C})}(\tilde{\Sigma}_{\scriptscriptstyle \blacksquare})} 
\! + \! \left\lvert \left\lvert \dfrac{w_{+}^{\tilde{\Sigma}_{
\scriptscriptstyle \blacksquare}}(\pmb{\cdot})}{2 \pi \mi 
(\pmb{\cdot} \! - \! z)} \right\rvert \right\rvert_{\mathcal{
L}^{2}_{\mathrm{M}_{2}(\mathbb{C})}(\tilde{\Sigma}_{
\scriptscriptstyle \blacksquare})} \right) \nonumber \\
+& \, \left\lvert \left\lvert (\id \! - \! (\id \! - \! C_{w^{\tilde{\Sigma}_{
\circlearrowright}}})^{-1}(\id \! - \! C_{w^{\tilde{\Sigma}_{\scriptscriptstyle 
\blacksquare}}})^{-1}C_{w^{\tilde{\Sigma}_{\scriptscriptstyle \blacksquare}}}
C_{w^{\tilde{\Sigma}_{\circlearrowright}}})^{-1} \right\rvert \right\rvert_{
\mathfrak{B}_{2}(\tilde{\Sigma}_{\circlearrowright})} \left\lvert \left\lvert 
(\id \! - \! C_{w^{\tilde{\Sigma}_{\circlearrowright}}})^{-1} \right\rvert 
\right\rvert_{\mathfrak{B}_{2}(\tilde{\Sigma}_{\circlearrowright})} 
\left\lvert \left\lvert (\id \! - \! C_{w^{\tilde{\Sigma}_{\scriptscriptstyle 
\blacksquare}}})^{-1} \right\rvert \right\rvert_{\mathfrak{B}_{2}(\tilde{
\Sigma}_{\circlearrowright})} \nonumber \\
\times& \, \left\lvert \left\lvert C_{w^{\tilde{\Sigma}_{\circlearrowright}}} 
\right\rvert \right\rvert_{\mathfrak{B}_{2}(\tilde{\Sigma}_{\circlearrowright})} 
\left\lvert \left\lvert C_{w^{\tilde{\Sigma}_{\scriptscriptstyle \blacksquare}}} 
\right\rvert \right\rvert_{\mathfrak{B}_{2}(\tilde{\Sigma}_{\circlearrowright})} 
\left(\left\lvert \left\lvert \dfrac{w_{+}^{\tilde{\Sigma}_{\circlearrowright}}
(\pmb{\cdot})}{2 \pi \mi (\pmb{\cdot} \! - \! \alpha_{k})} \right\rvert 
\right\rvert_{\mathcal{L}^{2}_{\mathrm{M}_{2}(\mathbb{C})}(\tilde{\Sigma}_{
\circlearrowright})} \! + \! \left\lvert \left\lvert \dfrac{w_{+}^{\tilde{
\Sigma}_{\circlearrowright}}(\pmb{\cdot})}{2 \pi \mi (\pmb{\cdot} \! - \! z)} 
\right\rvert \right\rvert_{\mathcal{L}^{2}_{\mathrm{M}_{2}(\mathbb{C})}
(\tilde{\Sigma}_{\circlearrowright})} \right) \nonumber \\
\underset{\underset{z_{o}=1+o(1)}{\mathscr{N},n \to \infty}}{\leqslant}& 
\, \dfrac{\left(\left\lvert \left\lvert C_{w^{\Sigma_{\tilde{\mathcal{R}}}}} 
\right\rvert \right\rvert_{\mathfrak{B}_{2}(\tilde{\Sigma}_{\tilde{\mathcal{
R}}}^{\sharp})} \right)^{2}}{\left(1 \! - \! \left\lvert \left\lvert C_{w^{\Sigma_{
\tilde{\mathcal{R}}}}} \right\rvert \right\rvert_{\mathfrak{B}_{2}(\tilde{
\Sigma}_{\tilde{\mathcal{R}}}^{\sharp})} \right)^{2} \! - \! \left(\left\lvert 
\left\lvert C_{w^{\Sigma_{\tilde{\mathcal{R}}}}} \right\rvert \right\rvert_{
\mathfrak{B}_{2}(\tilde{\Sigma}_{\tilde{\mathcal{R}}}^{\sharp})} \right)^{2}} 
\left(\left\lvert \left\lvert \dfrac{w_{+}^{\tilde{\Sigma}_{\scriptscriptstyle 
\blacksquare}}(\pmb{\cdot})}{2 \pi \mi (\pmb{\cdot} \! - \! \alpha_{k})} 
\right\rvert \right\rvert_{\mathcal{L}^{2}_{\mathrm{M}_{2}(\mathbb{C})}
(\tilde{\Sigma}_{\scriptscriptstyle \blacksquare})} \! + \! \left\lvert 
\left\lvert \dfrac{w_{+}^{\tilde{\Sigma}_{\scriptscriptstyle \blacksquare}}
(\pmb{\cdot})}{2 \pi \mi (\pmb{\cdot} \! - \! z)} \right\rvert \right\rvert_{
\mathcal{L}^{2}_{\mathrm{M}_{2}(\mathbb{C})}(\tilde{\Sigma}_{
\scriptscriptstyle \blacksquare})} \right. \nonumber \\
+&\left. \,  \left\lvert \left\lvert \dfrac{w_{+}^{\tilde{\Sigma}_{
\circlearrowright}}(\pmb{\cdot})}{2 \pi \mi (\pmb{\cdot} \! - \! 
\alpha_{k})} \right\rvert \right\rvert_{\mathcal{L}^{2}_{\mathrm{
M}_{2}(\mathbb{C})}(\tilde{\Sigma}_{\circlearrowright})} \! 
+ \! \left\lvert \left\lvert \dfrac{w_{+}^{\tilde{\Sigma}_{
\circlearrowright}}(\pmb{\cdot})}{2 \pi \mi (\pmb{\cdot} \! - \! 
z)} \right\rvert \right\rvert_{\mathcal{L}^{2}_{\mathrm{M}_{2}
(\mathbb{C})}(\tilde{\Sigma}_{\circlearrowright})} \right): 
\label{eqtazq23}
\end{align}
via the Estimate~\eqref{eqcwsigtll}, the 
Estimates~\eqref{eqtazq14}--\eqref{eqtazq17}, and the 
Inequality~\eqref{eqtazq23}, one arrives at, for $n \! \in \! 
\mathbb{N}$ and $k \! \in \! \lbrace 1,2,\dotsc,K \rbrace$ such that 
$\alpha_{p_{\mathfrak{s}}} \! := \! \alpha_{k} \! \neq \! \infty$,
\begin{equation} \label{eqtazq24} 
\left\vert \tilde{\mathbb{J}}^{\tilde{\mathcal{R}}}_{5}(z) \right\rvert 
\underset{\underset{z_{o}=1+o(1)}{\mathscr{N},n \to \infty}}{
\leqslant} \mathcal{O} \left(\dfrac{\tilde{\mathfrak{c}}_{\tilde{
\Sigma}}^{11}(n,k,z_{o}) \me^{-((n-1)K+k) \tilde{\lambda}_{\tilde{
\mathcal{R}},\tilde{w}}^{\triangleright}}}{((n \! - \! 1)K \! + \! k) \min 
\lbrace 1,\operatorname{dist}(\tilde{\Sigma}_{\tilde{\mathcal{R}}}^{
\sharp},z) \rbrace} \right), \quad z \! \in \! \mathbb{C} \setminus 
\tilde{\Sigma}_{\tilde{\mathcal{R}}}^{\sharp},
\end{equation}
where $\tilde{\mathfrak{c}}_{\tilde{\Sigma}}^{11}(n,k,z_{o}) 
\! =_{\underset{z_{o}=1+o(1)}{\mathscr{N},n \to \infty}} \! 
\mathcal{O}(1)$. Proceeding analogously as in the calculations 
above leading to the Estimate~\eqref{eqtazq24}, one shows that, 
for $n \! \in \! \mathbb{N}$ and $k \! \in \! \lbrace 1,2,\dotsc,K 
\rbrace$ such that $\alpha_{p_{\mathfrak{s}}} \! := \! \alpha_{k} 
\! \neq \! \infty$,
\begin{equation} \label{eqtazq25} 
\left\vert \tilde{\mathbb{J}}^{\tilde{\mathcal{R}}}_{6}(z) \right\rvert 
\underset{\underset{z_{o}=1+o(1)}{\mathscr{N},n \to \infty}}{
\leqslant} \mathcal{O} \left(\dfrac{\tilde{\mathfrak{c}}_{\tilde{
\Sigma}}^{12}(n,k,z_{o}) \me^{-((n-1)K+k) \tilde{\lambda}_{
\tilde{\mathcal{R}},\tilde{w}}^{\triangleright}}}{((n \! - \! 1)K \! + \! k) 
\min \lbrace 1,\operatorname{dist}(\tilde{\Sigma}_{\tilde{\mathcal{
R}}}^{\sharp},z) \rbrace} \right), \quad z \! \in \! \mathbb{C} 
\setminus \tilde{\Sigma}_{\tilde{\mathcal{R}}}^{\sharp},
\end{equation}
where $\tilde{\mathfrak{c}}_{\tilde{\Sigma}}^{12}(n,k,z_{o}) \! 
=_{\underset{z_{o}=1+o(1)}{\mathscr{N},n \to \infty}} \! 
\mathcal{O}(1)$.

An argument based on H\"{o}lder's Inequality, H\"{o}lder's 
Inequality for Integrals, and a Neumann series inversion (cf. 
Estimate~\eqref{eqcwsigtll}) shows that, for $z \! \in \! \mathbb{C} 
\setminus \tilde{\Sigma}_{\tilde{\mathcal{R}}}^{\sharp}$,
\begin{align}
\left\lvert \tilde{\mathbb{J}}_{7}^{\tilde{\mathcal{R}}}(z) \right\rvert 
\underset{\underset{z_{o}=1+o(1)}{\mathscr{N},n \to \infty}}{
\leqslant}& \, \left\lvert \left\lvert (\id \! - \! (\id \! - \! 
C_{w^{\tilde{\Sigma}_{\circlearrowright}}})^{-1}(\id \! - \! 
C_{w^{\tilde{\Sigma}_{\scriptscriptstyle \blacksquare}}})^{-1}
C_{w^{\tilde{\Sigma}_{\scriptscriptstyle \blacksquare}}}
C_{w^{\tilde{\Sigma}_{\circlearrowright}}})^{-1} \right\rvert 
\right\rvert_{\mathfrak{B}_{2}(\tilde{\Sigma}_{\scriptscriptstyle 
\blacksquare})} \left\lvert \left\lvert (\id \! - \! C_{w^{\tilde{
\Sigma}_{\circlearrowright}}})^{-1} \right\rvert \right\rvert_{
\mathfrak{B}_{2}(\tilde{\Sigma}_{\scriptscriptstyle \blacksquare})} 
\left\lvert \left\lvert (\id \! - \! C_{w^{\tilde{\Sigma}_{\scriptscriptstyle 
\blacksquare}}})^{-1} \right\rvert \right\rvert_{\mathfrak{B}_{2}
(\tilde{\Sigma}_{\scriptscriptstyle \blacksquare})} \nonumber \\
\times& \, \left\lvert \left\lvert C_{w^{\tilde{\Sigma}_{\circlearrowright}}} 
\right\rvert \right\rvert_{\mathfrak{B}_{2}(\tilde{\Sigma}_{\scriptscriptstyle 
\blacksquare})} \left\lvert \left\lvert (\id \! - \! C_{w^{\tilde{\Sigma}_{
\scriptscriptstyle \blacksquare}}})^{-1} \right\rvert \right\rvert_{
\mathfrak{B}_{2}(\tilde{\Sigma}_{\scriptscriptstyle \blacksquare})} 
\left(\left\lvert \left\lvert C_{w^{\tilde{\Sigma}_{\scriptscriptstyle 
\blacksquare}}} \right\rvert \right\rvert_{\mathfrak{B}_{2}(\tilde{
\Sigma}_{\scriptscriptstyle \blacksquare})} \right)^{2} \left(\left\lvert 
\left\lvert \dfrac{w_{+}^{\tilde{\Sigma}_{\scriptscriptstyle \blacksquare}}
(\pmb{\cdot})}{2 \pi \mi (\pmb{\cdot} \! - \! \alpha_{k})} 
\right\rvert \right\rvert_{\mathcal{L}^{2}_{\mathrm{M}_{2}
(\mathbb{C})}(\tilde{\Sigma}_{\scriptscriptstyle \blacksquare})} 
\! + \! \left\lvert \left\lvert \dfrac{w_{+}^{\tilde{\Sigma}_{
\scriptscriptstyle \blacksquare}}(\pmb{\cdot})}{2 \pi \mi 
(\pmb{\cdot} \! - \! z)} \right\rvert \right\rvert_{\mathcal{
L}^{2}_{\mathrm{M}_{2}(\mathbb{C})}(\tilde{\Sigma}_{
\scriptscriptstyle \blacksquare})} \right) \nonumber \\
+& \, \left\lvert \left\lvert (\id \! - \! (\id \! - \! C_{w^{\tilde{\Sigma}_{
\circlearrowright}}})^{-1}(\id \! - \! C_{w^{\tilde{\Sigma}_{\scriptscriptstyle 
\blacksquare}}})^{-1}C_{w^{\tilde{\Sigma}_{\scriptscriptstyle \blacksquare}}}
C_{w^{\tilde{\Sigma}_{\circlearrowright}}})^{-1} \right\rvert \right\rvert_{
\mathfrak{B}_{2}(\tilde{\Sigma}_{\circlearrowright})} \left\lvert \left\lvert 
(\id \! - \! C_{w^{\tilde{\Sigma}_{\circlearrowright}}})^{-1} \right\rvert 
\right\rvert_{\mathfrak{B}_{2}(\tilde{\Sigma}_{\circlearrowright})} 
\left\lvert \left\lvert (\id \! - \! C_{w^{\tilde{\Sigma}_{\scriptscriptstyle 
\blacksquare}}})^{-1} \right\rvert \right\rvert_{\mathfrak{B}_{2}
(\tilde{\Sigma}_{\circlearrowright})} \nonumber \\
\times& \, \left\lvert \left\lvert C_{w^{\tilde{\Sigma}_{\circlearrowright}}} 
\right\rvert \right\rvert_{\mathfrak{B}_{2}(\tilde{\Sigma}_{\circlearrowright})} 
\left\lvert \left\lvert (\id \! - \! C_{w^{\tilde{\Sigma}_{\scriptscriptstyle 
\blacksquare}}})^{-1} \right\rvert \right\rvert_{\mathfrak{B}_{2}(\tilde{
\Sigma}_{\circlearrowright})} \left(\left\lvert \left\lvert C_{w^{\tilde{
\Sigma}_{\scriptscriptstyle \blacksquare}}} \right\rvert \right\rvert_{
\mathfrak{B}_{2}(\tilde{\Sigma}_{\circlearrowright})} \right)^{2} \left(
\left\lvert \left\lvert \dfrac{w_{+}^{\tilde{\Sigma}_{\circlearrowright}}
(\pmb{\cdot})}{2 \pi \mi (\pmb{\cdot} \! - \! \alpha_{k})} \right\rvert 
\right\rvert_{\mathcal{L}^{2}_{\mathrm{M}_{2}(\mathbb{C})}(\tilde{\Sigma}_{
\circlearrowright})} \! + \! \left\lvert \left\lvert \dfrac{w_{+}^{\tilde{
\Sigma}_{\circlearrowright}}(\pmb{\cdot})}{2 \pi \mi (\pmb{\cdot} \! - \! z)} 
\right\rvert \right\rvert_{\mathcal{L}^{2}_{\mathrm{M}_{2}(\mathbb{C})}
(\tilde{\Sigma}_{\circlearrowright})} \right) \nonumber \\
\underset{\underset{z_{o}=1+o(1)}{\mathscr{N},n \to \infty}}{\leqslant}& 
\, \dfrac{\left(\left\lvert \left\lvert C_{w^{\Sigma_{\tilde{\mathcal{R}}}}} 
\right\rvert \right\rvert_{\mathfrak{B}_{2}(\tilde{\Sigma}_{\tilde{\mathcal{
R}}}^{\sharp})} \right)^{3}}{\left(1 \! - \! \left\lvert \left\lvert C_{w^{\Sigma_{
\tilde{\mathcal{R}}}}} \right\rvert \right\rvert_{\mathfrak{B}_{2}(\tilde{
\Sigma}_{\tilde{\mathcal{R}}}^{\sharp})} \right) \left(\left(1 \! - \! \left\lvert 
\left\lvert C_{w^{\Sigma_{\tilde{\mathcal{R}}}}} \right\rvert \right\rvert_{
\mathfrak{B}_{2}(\tilde{\Sigma}_{\tilde{\mathcal{R}}}^{\sharp})} \right)^{2} \! 
- \! \left(\left\lvert \left\lvert C_{w^{\Sigma_{\tilde{\mathcal{R}}}}} \right\rvert 
\right\rvert_{\mathfrak{B}_{2}(\tilde{\Sigma}_{\tilde{\mathcal{R}}}^{\sharp})} 
\right)^{2} \right)} \left(\left\lvert \left\lvert \dfrac{w_{+}^{\tilde{\Sigma}_{
\scriptscriptstyle \blacksquare}}(\pmb{\cdot})}{2 \pi \mi (\pmb{\cdot} \! - \! 
\alpha_{k})} \right\rvert \right\rvert_{\mathcal{L}^{2}_{\mathrm{M}_{2}(\mathbb{C})}
(\tilde{\Sigma}_{\scriptscriptstyle \blacksquare})} \right. \nonumber \\
+&\left. \, \left\lvert \left\lvert \dfrac{w_{+}^{\tilde{\Sigma}_{\scriptscriptstyle 
\blacksquare}}(\pmb{\cdot})}{2 \pi \mi (\pmb{\cdot} \! - \! z)} \right\rvert 
\right\rvert_{\mathcal{L}^{2}_{\mathrm{M}_{2}(\mathbb{C})}(\tilde{\Sigma}_{
\scriptscriptstyle \blacksquare})} \! + \! \left\lvert \left\lvert \dfrac{w_{+}^{
\tilde{\Sigma}_{\circlearrowright}}(\pmb{\cdot})}{2 \pi \mi (\pmb{\cdot} 
\! - \! \alpha_{k})} \right\rvert \right\rvert_{\mathcal{L}^{2}_{
\mathrm{M}_{2}(\mathbb{C})}(\tilde{\Sigma}_{\circlearrowright})} 
\! + \! \left\lvert \left\lvert \dfrac{w_{+}^{\tilde{\Sigma}_{
\circlearrowright}}(\pmb{\cdot})}{2 \pi \mi (\pmb{\cdot} \! - \! 
z)} \right\rvert \right\rvert_{\mathcal{L}^{2}_{\mathrm{M}_{2}
(\mathbb{C})}(\tilde{\Sigma}_{\circlearrowright})} \right): \label{eqtazq26}
\end{align}
via the Estimate~\eqref{eqcwsigtll}, the 
Estimates~\eqref{eqtazq14}--\eqref{eqtazq17}, and the 
Inequality~\eqref{eqtazq26}, one arrives at, for $n \! \in \! 
\mathbb{N}$ and $k \! \in \! \lbrace 1,2,\dotsc,K \rbrace$ such that 
$\alpha_{p_{\mathfrak{s}}} \! := \! \alpha_{k} \! \neq \! \infty$,
\begin{equation} \label{eqtazq27} 
\left\vert \tilde{\mathbb{J}}^{\tilde{\mathcal{R}}}_{7}(z) \right\rvert 
\underset{\underset{z_{o}=1+o(1)}{\mathscr{N},n \to \infty}}{
\leqslant} \mathcal{O} \left(\dfrac{\tilde{\mathfrak{c}}_{\tilde{
\Sigma}}^{13}(n,k,z_{o}) \me^{-\frac{3}{2}((n-1)K+k) \tilde{\lambda}_{
\tilde{\mathcal{R}},\tilde{w}}^{\triangleright}}}{((n \! - \! 1)K \! + \! k) 
\min \lbrace 1,\operatorname{dist}(\tilde{\Sigma}_{\tilde{\mathcal{
R}}}^{\sharp},z) \rbrace} \right), \quad z \! \in \! \mathbb{C} 
\setminus \tilde{\Sigma}_{\tilde{\mathcal{R}}}^{\sharp},
\end{equation}
where $\tilde{\mathfrak{c}}_{\tilde{\Sigma}}^{13}(n,k,z_{o}) \! 
=_{\underset{z_{o}=1+o(1)}{\mathscr{N},n \to \infty}} \! 
\mathcal{O}(1)$. Proceeding analogously as in the calculations 
above leading to the Estimate~\eqref{eqtazq27}, one shows that, 
for $n \! \in \! \mathbb{N}$ and $k \! \in \! \lbrace 1,2,\dotsc,K 
\rbrace$ such that $\alpha_{p_{\mathfrak{s}}} \! := \! \alpha_{k} 
\! \neq \! \infty$,
\begin{equation} \label{eqtazq28} 
\left\vert \tilde{\mathbb{J}}^{\tilde{\mathcal{R}}}_{8}(z) \right\rvert 
\underset{\underset{z_{o}=1+o(1)}{\mathscr{N},n \to \infty}}{
\leqslant} \mathcal{O} \left(\dfrac{\tilde{\mathfrak{c}}_{\tilde{
\Sigma}}^{14}(n,k,z_{o}) \me^{-\frac{3}{2}((n-1)K+k) \tilde{\lambda}_{
\tilde{\mathcal{R}},\tilde{w}}^{\triangleright}}}{((n \! - \! 1)K \! + \! k) 
\min \lbrace 1,\operatorname{dist}(\tilde{\Sigma}_{\tilde{\mathcal{
R}}}^{\sharp},z) \rbrace} \right), \quad z \! \in \! \mathbb{C} 
\setminus \tilde{\Sigma}_{\tilde{\mathcal{R}}}^{\sharp},
\end{equation}
where $\tilde{\mathfrak{c}}_{\tilde{\Sigma}}^{14}(n,k,z_{o}) \! 
=_{\underset{z_{o}=1+o(1)}{\mathscr{N},n \to \infty}} \! 
\mathcal{O}(1)$.

{}From the Estimates~\eqref{eqtazq12}, \eqref{eqtazq18}, 
\eqref{eqtazq19}, \eqref{eqtazq21}, \eqref{eqtazq22}, 
\eqref{eqtazq24}, \eqref{eqtazq25}, \eqref{eqtazq27}, 
and~\eqref{eqtazq28}, one shows that (cf. Equation~\eqref{eqtazq1}), 
for $n \! \in \! \mathbb{N}$ and $k \! \in \! \lbrace 1,2,\dotsc,K 
\rbrace$ such that $\alpha_{p_{\mathfrak{s}}} \! := \! \alpha_{k} 
\! \neq \! \infty$, uniformly for compact subsets of $\mathbb{C} 
\setminus \tilde{\Sigma}_{\tilde{\mathcal{R}}}^{\sharp}$ $(\ni \! z)$,
\begin{equation*}
\left\lvert \int_{\tilde{\Sigma}_{\scriptscriptstyle \blacksquare}}w_{+}^{
\tilde{\Sigma}_{\scriptscriptstyle \blacksquare}}(\xi) \left(\dfrac{1}{\xi 
\! - \! z} \! - \! \dfrac{1}{\xi \! - \! \alpha_{k}} \right) \dfrac{\md \xi}{2 
\pi \mi} \! + \! \sum_{m=1}^{8} \tilde{\mathbb{J}}_{m}^{\tilde{\mathcal{
R}}}(z) \right\rvert \underset{\underset{z_{o}=1+o(1)}{\mathscr{N},n \to 
\infty}}{\leqslant} \mathcal{O} \left(\dfrac{\tilde{\mathfrak{c}}_{\tilde{
\Sigma}}^{\vartriangle}(n,k,z_{o}) \me^{-\frac{1}{2}((n-1)K+k) \tilde{
\lambda}_{\tilde{\mathcal{R}},\tilde{w}}^{\triangleright}}}{((n \! - \! 1)K 
\! + \! k) \min \lbrace 1,\operatorname{dist}(\tilde{\Sigma}_{\tilde{
\mathcal{R}}}^{\sharp},z) \rbrace} \right), \quad z \! \in \! \mathbb{C} 
\setminus \tilde{\Sigma}_{\tilde{\mathcal{R}}}^{\sharp},
\end{equation*}
where $\tilde{\mathfrak{c}}_{\tilde{\Sigma}}^{\vartriangle}(n,k,z_{o}) \! 
=_{\underset{z_{o}=1+o(1)}{\mathscr{N},n \to \infty}} \! \mathcal{O}
(1)$; consequently,
\begin{equation*}
\tilde{\mathcal{R}}(z) \! - \! \mathrm{I} \! - \! \int_{\tilde{\Sigma}_{
\circlearrowright}} \dfrac{(z \! - \! \alpha_{k})w_{+}^{\tilde{\Sigma}_{
\circlearrowright}}(\xi)}{(\xi \! - \! \alpha_{k})(\xi \! - \! z)} \, 
\dfrac{\md \xi}{2 \pi \mi} \underset{\underset{z_{o}=1+o(1)}{
\mathscr{N},n \to \infty}}{=} \mathcal{O} \left(\dfrac{\tilde{
\mathfrak{c}}_{\tilde{\Sigma}_{\circlearrowright}}(n,k,z_{o}) 
\me^{-\frac{1}{2}((n-1)K+k) \tilde{\lambda}_{\tilde{\mathcal{R}},
\tilde{w}}^{\triangleright}}}{((n \! - \! 1)K \! + \! k) \min \lbrace 1,
\operatorname{dist}(\tilde{\Sigma}_{\tilde{\mathcal{R}}}^{\sharp},z) 
\rbrace} \right), \quad z \! \in \! \mathbb{C} \setminus 
\tilde{\Sigma}_{\tilde{\mathcal{R}}}^{\sharp},
\end{equation*}
where $(\mathrm{M}_{2}(\mathbb{C}) \! \ni)$ $\tilde{\mathfrak{c}}_{
\tilde{\Sigma}_{\circlearrowright}}(n,k,z_{o}) \! =_{\underset{z_{o}=
1+o(1)}{\mathscr{N},n \to \infty}} \! \mathcal{O}(1)$, which is 
Equation~\eqref{eqtazt1}.

The analysis for the case $n \! \in \! \mathbb{N}$ and $k \! \in \! 
\lbrace 1,2,\dotsc,K \rbrace$ such that $\alpha_{p_{\mathfrak{s}}} \! 
:= \! \alpha_{k} \! = \! \infty$ is, \emph{mutatis mutandis}, analogous, 
and leads to the asymptotics, in the double-scaling limit $\mathscr{N},
n \! \to \! \infty$ such that $z_{o} \! = \! 1 \! + \! o(1)$, for $\hat{
\mathcal{R}}(z)$ given in Equation~\eqref{eqtazh1}. \hfill $\qed$
\begin{bbbbb} \label{propo5.3} 
For $n \! \in \! \mathbb{N}$ and $k \! \in \! \lbrace 1,2,\dotsc,K 
\rbrace$ such that $\alpha_{p_{\mathfrak{s}}} \! := \! \alpha_{k} \! 
= \! \infty$ (resp., $\alpha_{p_{\mathfrak{s}}} \! := \! \alpha_{k} \! 
\neq \! \infty)$, let $\hat{\mathcal{R}} \colon \mathbb{C} \setminus 
\hat{\Sigma}_{\hat{\mathcal{R}}}^{\sharp} \! \to \! \mathrm{SL}_{2}
(\mathbb{C})$ (resp., $\tilde{\mathcal{R}} \colon \mathbb{C} 
\setminus \tilde{\Sigma}_{\tilde{\mathcal{R}}}^{\sharp} \! \to \! 
\mathrm{SL}_{2}(\mathbb{C}))$ solve the equivalent {\rm RHP} 
$(\hat{\mathcal{R}}(z),\mathrm{I} \! + \! w_{+}^{\Sigma_{\hat{
\mathcal{R}}}}(z),\hat{\Sigma}_{\hat{\mathcal{R}}}^{\sharp})$ (resp., 
$(\tilde{\mathcal{R}}(z),\mathrm{I} \! + \! w_{+}^{\Sigma_{\tilde{
\mathcal{R}}}}(z),\tilde{\Sigma}_{\tilde{\mathcal{R}}}^{\sharp}))$ 
and be given by Equation~\eqref{eqtazh1} (resp., 
Equation~\eqref{eqtazt1}$)$, and let $w_{+}^{\hat{\Sigma}_{
\circlearrowright}}(z) \! := \! w_{+}^{\Sigma_{\hat{\mathcal{R}}}}(z) 
\! \! \upharpoonright_{\hat{\Sigma}_{\circlearrowright}}$ (resp., 
$w_{+}^{\tilde{\Sigma}_{\circlearrowright}}(z) \! := \! w_{+}^{
\Sigma_{\tilde{\mathcal{R}}}}(z) \! \! \upharpoonright_{\tilde{
\Sigma}_{\circlearrowright}})$, where $\hat{\Sigma}_{\circlearrowright} 
\! := \! \cup_{j=1}^{N+1}(\partial \hat{\mathbb{U}}_{\hat{\delta}_{
\hat{b}_{j-1}}} \cup \partial \hat{\mathbb{U}}_{\hat{\delta}_{\hat{a}_{j}}})$ 
(resp., $\tilde{\Sigma}_{\circlearrowright} \! := \! \cup_{j=1}^{N+1}
(\partial \tilde{\mathbb{U}}_{\tilde{\delta}_{\tilde{b}_{j-1}}} \cup 
\partial \tilde{\mathbb{U}}_{\tilde{\delta}_{\tilde{a}_{j}}}))$, have, 
for $z \! \in \! \hat{\Sigma}_{\circlearrowright}$ (resp., $z \! \in 
\! \tilde{\Sigma}_{\circlearrowright})$, the asymptotics given by 
the Expansions~\eqref{eqprophatb} and~\eqref{eqprophata} (resp., 
Expansions~\eqref{eqproptilb} and~\eqref{eqproptila}$)$. Then$:$ 
{\rm \pmb{(1)}} for $n \! \in \! \mathbb{N}$ and $k \! \in \! \lbrace 
1,2,\dotsc,K \rbrace$ such that $\alpha_{p_{\mathfrak{s}}} \! := 
\! \alpha_{k} \! = \! \infty$, uniformly for compact subsets of 
$\mathbb{C} \setminus \hat{\Sigma}_{\hat{\mathcal{R}}}^{\sharp}$ 
$(\ni \! z)$,
\begin{align}
\hat{\mathcal{R}}(z) \underset{\underset{z_{o}=1+o(1)}{\mathscr{N},
n \to \infty}}{=}& \, \mathrm{I} \! + \! \dfrac{1}{((n \! - \! 1)K \! + \! k)} 
\sum_{j=1}^{N+1} \left(\dfrac{(\hat{\alpha}_{0}(\hat{b}_{j-1}))^{-1}}{
(z \! - \! \hat{b}_{j-1})^{2}} \hat{\boldsymbol{\mathrm{A}}}(\hat{b}_{j
-1}) \! + \! \dfrac{(\hat{\alpha}_{0}(\hat{b}_{j-1}))^{-2}}{z \! - \! 
\hat{b}_{j-1}} \left(\hat{\alpha}_{0}(\hat{b}_{j-1}) 
\hat{\boldsymbol{\mathrm{B}}}(\hat{b}_{j-1}) \! - \! \hat{\alpha}_{1}
(\hat{b}_{j-1}) \hat{\boldsymbol{\mathrm{A}}}(\hat{b}_{j-1}) \right) 
\right. \nonumber \\
+&\left. \, \dfrac{(\hat{\alpha}_{0}(\hat{a}_{j}))^{-1}}{(z \! - \! 
\hat{a}_{j})^{2}} \hat{\boldsymbol{\mathrm{A}}}(\hat{a}_{j}) \! + \! 
\dfrac{(\hat{\alpha}_{0}(\hat{a}_{j}))^{-2}}{z \! - \! \hat{a}_{j}} \left(
\hat{\alpha}_{0}(\hat{a}_{j}) \hat{\boldsymbol{\mathrm{B}}}(\hat{a}_{j}) 
\! - \! \hat{\alpha}_{1}(\hat{a}_{j}) \hat{\boldsymbol{\mathrm{A}}}
(\hat{a}_{j}) \right) \! - \! \hat{\mathbb{Y}}_{\hat{b}_{j-1}}(z) 
\chi_{\hat{\mathbb{U}}_{\hat{\delta}_{\hat{b}_{j-1}}}}(z) \! - \! 
\hat{\mathbb{Y}}_{\hat{a}_{j}}(z) \chi_{\hat{\mathbb{U}}_{\hat{
\delta}_{\hat{a}_{j}}}}(z) \vphantom{M^{M^{M^{M^{M^{M^{M}}}}}}} 
\right) \nonumber \\
+& \, \mathcal{O} \left(\dfrac{\hat{\mathfrak{c}}_{\hat{\mathbb{
E}}_{\hat{\mathcal{R}}}^{\ast}}(n,k,z_{o}) \hat{\mathbb{E}}_{\hat{
\mathcal{R}}}^{\ast}(z)}{((n \! - \! 1)K \! + \! k)^{2}} \right), 
\quad z \! \in \! \mathbb{C} \setminus \hat{\Sigma}_{\hat{
\mathcal{R}}}^{\sharp} , \label{eqpropo5.3hA}
\end{align}
where, for $j \! \in \! \lbrace 1,2,\dotsc,N \! + \! 1 \rbrace$, 
$\hat{\boldsymbol{\mathrm{A}}}(\hat{b}_{j-1})$, $\hat{
\boldsymbol{\mathrm{A}}}(\hat{a}_{j})$, $\hat{\boldsymbol{\mathrm{B}}}
(\hat{b}_{j-1})$, $\hat{\boldsymbol{\mathrm{B}}}(\hat{a}_{j})$, 
$\hat{\alpha}_{0}(\hat{b}_{j-1})$, $\hat{\alpha}_{0}(\hat{a}_{j})$, 
$\hat{\alpha}_{1}(\hat{b}_{j-1})$, and $\hat{\alpha}_{1}(\hat{a}_{j})$ 
are defined in item~{\rm \pmb{(1)}} of Proposition~\ref{propo5.1}, 
$\hat{\mathbb{Y}}_{\hat{b}_{j-1}}(z)$ and $\hat{\mathbb{Y}}_{\hat{a}_{j}}
(z)$ are defined by Equations~\eqref{eqmaininf68} and~\eqref{eqmaininf69}, 
respectively, $(\mathrm{M}_{2}(\mathbb{C}) \! \ni)$ $\hat{\mathfrak{c}}_{
\hat{\mathbb{E}}_{\hat{\mathcal{R}}}^{\ast}}(n,k,z_{o}) \! =_{\underset{z_{o}
=1+o(1)}{\mathscr{N},n \to \infty}} \! \mathcal{O}(1)$, and
\begin{equation*}
\hat{\mathbb{E}}_{\hat{\mathcal{R}}}^{\ast}(z) \! = \! 
\begin{cases}
\max\limits_{\underset{m=1,2,3}{j=1,2,\dotsc,N+1}} \left\lbrace 
(z \! - \! \hat{b}_{j-1})^{-m},(z \! - \! \hat{a}_{j})^{-m} \right\rbrace, 
&\text{$z \! \in \! \mathbb{C} \setminus \cup_{j=1}^{N+1}
(\hat{\mathbb{U}}_{\hat{\delta}_{\hat{b}_{j-1}}} \cup \hat{\mathbb{U}}_{
\hat{\delta}_{\hat{a}_{j}}} \cup \partial \hat{\mathbb{U}}_{\hat{\delta}_{
\hat{b}_{j-1}}} \cup \partial \hat{\mathbb{U}}_{\hat{\delta}_{\hat{a}_{j}}})$,} \\
1, &\text{$z \! \in \! \cup_{j=1}^{N+1}(\hat{\mathbb{U}}_{\hat{\delta}_{
\hat{b}_{j-1}}} \cup \hat{\mathbb{U}}_{\hat{\delta}_{\hat{a}_{j}}})$$;$}
\end{cases}
\end{equation*}
and {\rm \pmb{(2)}} for $n \! \in \! \mathbb{N}$ and $k \! \in \! 
\lbrace 1,2,\dotsc,K \rbrace$ such that $\alpha_{p_{\mathfrak{s}}} \! 
:= \! \alpha_{k} \! \neq \! \infty$, uniformly for compact subsets 
of $\mathbb{C} \setminus \tilde{\Sigma}_{\tilde{\mathcal{R}}}^{
\sharp}$ $(\ni \! z)$,
\begin{align}
\tilde{\mathcal{R}}(z) \underset{\underset{z_{o}=1+o(1)}{\mathscr{N},n 
\to \infty}}{=}& \, \mathrm{I} \! + \! \dfrac{1}{((n \! - \! 1)K \! + \! k)} 
\sum_{j=1}^{N+1} \left(\dfrac{(\tilde{\alpha}_{0}(\tilde{b}_{j-1}))^{-1}}{
(z \! - \! \tilde{b}_{j-1})^{2}} \tilde{\boldsymbol{\mathrm{A}}}(\tilde{
b}_{j-1}) \! + \! \dfrac{(\tilde{\alpha}_{0}(\tilde{b}_{j-1}))^{-2}}{z \! - \! 
\tilde{b}_{j-1}} \left(\tilde{\alpha}_{0}(\tilde{b}_{j-1}) \tilde{\boldsymbol{
\mathrm{B}}}(\tilde{b}_{j-1}) \! - \! \tilde{\alpha}_{1}(\tilde{b}_{j-1}) 
\tilde{\boldsymbol{\mathrm{A}}}(\tilde{b}_{j-1}) \right) \right. 
\nonumber \\
+&\left. \, \dfrac{(\tilde{\alpha}_{0}(\tilde{a}_{j}))^{-1}}{(z \! - \! 
\tilde{a}_{j})^{2}} \tilde{\boldsymbol{\mathrm{A}}}(\tilde{a}_{j}) \! + 
\! \dfrac{(\tilde{\alpha}_{0}(\tilde{a}_{j}))^{-2}}{z \! - \! \tilde{a}_{j}} 
\left(\tilde{\alpha}_{0}(\tilde{a}_{j}) \tilde{\boldsymbol{\mathrm{B}}}
(\tilde{a}_{j}) \! - \! \tilde{\alpha}_{1}(\tilde{a}_{j}) \tilde{\boldsymbol{
\mathrm{A}}}(\tilde{a}_{j}) \right) \! + \! \dfrac{((\alpha_{k} \! - \! 
\tilde{b}_{j-1}) \tilde{\alpha}_{1}(\tilde{b}_{j-1}) \! - \! \tilde{\alpha}_{0}
(\tilde{b}_{j-1}))}{(\alpha_{k} \! - \! \tilde{b}_{j-1})^{2}(\tilde{\alpha}_{0}
(\tilde{b}_{j-1}))^{2}} \tilde{\boldsymbol{\mathrm{A}}}(\tilde{b}_{j-1}) 
\right. \nonumber \\
+&\left. \, \dfrac{((\alpha_{k} \! - \! \tilde{a}_{j}) \tilde{\alpha}_{1}
(\tilde{a}_{j}) \! - \! \tilde{\alpha}_{0}(\tilde{a}_{j}))}{(\alpha_{k} 
\! - \! \tilde{a}_{j})^{2}(\tilde{\alpha}_{0}(\tilde{a}_{j}))^{2}} 
\tilde{\boldsymbol{\mathrm{A}}}(\tilde{a}_{j}) \! - \! 
\dfrac{(\tilde{\alpha}_{0}(\tilde{b}_{j-1}))^{-1}}{\alpha_{k} \! - \! 
\tilde{b}_{j-1}} \tilde{\boldsymbol{\mathrm{B}}}(\tilde{b}_{j-1}) \! 
- \! \dfrac{(\tilde{\alpha}_{0}(\tilde{a}_{j}))^{-1}}{\alpha_{k} \! - \! 
\tilde{a}_{j}} \tilde{\boldsymbol{\mathrm{B}}}(\tilde{a}_{j}) \! - \! 
\tilde{\mathbb{Y}}_{\tilde{b}_{j-1}}(z) \chi_{\tilde{\mathbb{U}}_{
\tilde{\delta}_{\tilde{b}_{j-1}}}}(z) \right. \nonumber \\
-&\left. \, \tilde{\mathbb{Y}}_{\tilde{a}_{j}}(z) \chi_{\tilde{
\mathbb{U}}_{\tilde{\delta}_{\tilde{a}_{j}}}}(z) 
\vphantom{M^{M^{M^{M^{M^{M^{M}}}}}}} \right) \! + \! 
\mathcal{O} \left(\dfrac{\tilde{\mathfrak{c}}_{\tilde{\mathbb{E}}_{
\tilde{\mathcal{R}}}^{\ast}}(n,k,z_{o}) \tilde{\mathbb{E}}_{\tilde{
\mathcal{R}}}^{\ast}(z)}{((n \! - \! 1)K \! + \! k)^{2}} \right), 
\quad z \! \in \! \mathbb{C} \setminus \tilde{\Sigma}_{\tilde{
\mathcal{R}}}^{\sharp} , \label{eqpropo5.3tB}
\end{align}
where, for $j \! \in \! \lbrace 1,2,\dotsc,N \! + \! 1 \rbrace$, $\tilde{
\boldsymbol{\mathrm{A}}}(\tilde{b}_{j-1})$, $\tilde{\boldsymbol{
\mathrm{A}}}(\tilde{a}_{j})$, $\tilde{\boldsymbol{\mathrm{B}}}
(\tilde{b}_{j-1})$, $\tilde{\boldsymbol{\mathrm{B}}}(\tilde{a}_{j})$, 
$\tilde{\alpha}_{0}(\tilde{b}_{j-1})$, $\tilde{\alpha}_{0}(\tilde{a}_{j})$, 
$\tilde{\alpha}_{1}(\tilde{b}_{j-1})$, and $\tilde{\alpha}_{1}(\tilde{a}_{j})$ 
are defined in item~{\rm \pmb{(2)}} of Proposition~\ref{propo5.1}, 
$\tilde{\mathbb{Y}}_{\tilde{b}_{j-1}}(z)$ and $\tilde{\mathbb{Y}}_{
\tilde{a}_{j}}(z)$ are defined by Equations~\eqref{eqmainfin70} 
and~\eqref{eqmainfin71}, respectively, $(\mathrm{M}_{2}(\mathbb{C}) 
\! \ni)$ $\tilde{\mathfrak{c}}_{\tilde{\mathbb{E}}_{\tilde{\mathcal{R}}}^{
\ast}}(n,k,z_{o}) \! =_{\underset{z_{o}=1+o(1)}{\mathscr{N},n \to \infty}} 
\! \mathcal{O}(1)$, and
\begin{equation*}
\tilde{\mathbb{E}}_{\tilde{\mathcal{R}}}^{\ast}(z) \! = \! 
\begin{cases}
\max \left\lbrace 1,\max\limits_{\underset{m=1,2,3}{j=1,2,\dotsc,N+1}} 
\left\lbrace (z \! - \! \tilde{b}_{j-1})^{-m},(z \! - \! \tilde{a}_{j})^{-m} 
\right\rbrace \right\rbrace, &\text{$z \! \in \! \mathbb{C} \setminus 
\cup_{j=1}^{N+1}(\tilde{\mathbb{U}}_{\tilde{\delta}_{\tilde{b}_{j-1}}} 
\cup \tilde{\mathbb{U}}_{\tilde{\delta}_{\tilde{a}_{j}}} \cup \partial 
\tilde{\mathbb{U}}_{\tilde{\delta}_{\tilde{b}_{j-1}}} \cup \partial 
\tilde{\mathbb{U}}_{\tilde{\delta}_{\tilde{a}_{j}}})$,} \\
1, &\text{$z \! \in \! \cup_{j=1}^{N+1}(\tilde{\mathbb{U}}_{\tilde{\delta}_{
\tilde{b}_{j-1}}} \cup \tilde{\mathbb{U}}_{\tilde{\delta}_{\tilde{a}_{j}}})$.}
\end{cases}
\end{equation*}
\end{bbbbb}

\emph{Proof.} The proof of this Proposition~\ref{propo5.3} consists of 
two cases: (i) $n \! \in \! \mathbb{N}$ and $k \! \in \! \lbrace 1,2,\dotsc,
K \rbrace$ such that $\alpha_{p_{\mathfrak{s}}} \! := \! \alpha_{k} \! 
= \! \infty$; and (ii) $n \! \in \! \mathbb{N}$ and $k \! \in \! \lbrace 
1,2,\dotsc,K \rbrace$ such that $\alpha_{p_{\mathfrak{s}}} \! := \! 
\alpha_{k} \! \neq \! \infty$. Notwithstanding the fact that the scheme 
of the proof is, \emph{mutatis mutandis}, similar for both cases, 
case~(ii), nevertheless, is the more technically challenging of the two; 
therefore, without loss of generality, only the proof for case~(ii) is 
presented in detail, whilst case~(i) is proved analogously.

For $n \! \in \! \mathbb{N}$ and $k \! \in \! \lbrace 1,2,\dotsc,K 
\rbrace$ such that $\alpha_{p_{\mathfrak{s}}} \! := \! \alpha_{k} 
\! \neq \! \infty$, recall that $\tilde{\mathcal{R}} \colon \mathbb{C} 
\setminus \tilde{\Sigma}_{\tilde{\mathcal{R}}}^{\sharp} \! \to \! 
\mathrm{SL}_{2}(\mathbb{C})$ solves the equivalent RHP 
$(\tilde{\mathcal{R}}(z),\mathrm{I} \! + \! w_{+}^{\Sigma_{\tilde{
\mathcal{R}}}}(z),\tilde{\Sigma}_{\tilde{\mathcal{R}}}^{\sharp})$ with 
the asymptotic, in the double-scaling limit $\mathscr{N},n \! \to \! 
\infty$ such that $z_{o} \! = \! 1 \! + \! o(1)$, representation given 
by Equation~\eqref{eqtazt1}. Recalling {}from Lemma~\ref{lem5.4} 
that $\tilde{\Sigma}_{\circlearrowright} \! := \! \cup_{j=1}^{N+1}
(\partial \tilde{\mathbb{U}}_{\tilde{\delta}_{\tilde{b}_{j-1}}} \cup 
\partial \tilde{\mathbb{U}}_{\tilde{\delta}_{\tilde{a}_{j}}})$, with, 
in particular, $\partial \tilde{\mathbb{U}}_{\tilde{\delta}_{\tilde{b}_{j-1}}} 
\cap \partial \tilde{\mathbb{U}}_{\tilde{\delta}_{\tilde{a}_{j}}} \! = \! 
\varnothing$, and $w_{+}^{\tilde{\Sigma}_{\circlearrowright}}(z) \! 
:= \! w_{+}^{\Sigma_{\tilde{\mathcal{R}}}}(z) \! \! \upharpoonright_{
\tilde{\Sigma}_{\circlearrowright}}$, a straightforward partial-fraction 
decomposition argument shows that, uniformly for compact subsets of 
$\mathbb{C} \setminus \tilde{\Sigma}_{\tilde{\mathcal{R}}}^{\sharp}$ 
$(\ni \! z)$,
\begin{align}
\tilde{\mathcal{R}}(z) \underset{\underset{z_{o}=1+o(1)}{\mathscr{N},
n \to \infty}}{=}& \, \mathrm{I} \! + \! \sum_{j=1}^{N+1} \left(\oint_{
\partial \tilde{\mathbb{U}}_{\tilde{\delta}_{\tilde{b}_{j-1}}}} \! + \! 
\oint_{\partial \tilde{\mathbb{U}}_{\tilde{\delta}_{\tilde{a}_{j}}}} \right) 
\dfrac{w^{\tilde{\Sigma}_{\circlearrowright}}_{+}(\xi)}{\xi \! - \! 
\alpha_{k}} \, \dfrac{\md \xi}{2 \pi \mi} \! - \! \sum_{j=1}^{N+1} \left(
\oint_{\partial \tilde{\mathbb{U}}_{\tilde{\delta}_{\tilde{b}_{j-1}}}} \! + \! 
\oint_{\partial \tilde{\mathbb{U}}_{\tilde{\delta}_{\tilde{a}_{j}}}} \right) 
\dfrac{w^{\tilde{\Sigma}_{\circlearrowright}}_{+}(\xi)}{\xi \! - \! z} \, 
\dfrac{\md \xi}{2 \pi \mi} \nonumber \\
+& \, \mathcal{O} \left(\dfrac{\tilde{\mathfrak{c}}_{\tilde{\Sigma}_{
\circlearrowright}}(n,k,z_{o}) \me^{-\frac{1}{2}((n-1)K+k) \tilde{
\lambda}_{\tilde{\mathcal{R}},\tilde{w}}^{\triangleright}}}{((n \! - \! 
1)K \! + \! k) \min \lbrace 1,\operatorname{dist}(\tilde{\Sigma}_{
\tilde{\mathcal{R}}}^{\sharp},z) \rbrace} \right), \quad z \! \in \! 
\mathbb{C} \setminus \tilde{\Sigma}_{\tilde{\mathcal{R}}}^{\sharp}, 
\label{eqtilrzo1}
\end{align}
where, for $j \! \in \! \lbrace 1,2,\dotsc,N \! + \! 1 \rbrace$, $\oint_{
\partial \tilde{\mathbb{U}}_{\tilde{\delta}_{\tilde{b}_{j-1}}}}$ and 
$\oint_{\partial \tilde{\mathbb{U}}_{\tilde{\delta}_{\tilde{a}_{j}}}}$ are 
counter-clockwise-oriented, closed (contour) integrals about the 
discs of radii $\tilde{\delta}_{\tilde{b}_{j-1}}$ and $\tilde{\delta}_{
\tilde{a}_{j}}$, respectively, surrounding the end-points of the 
intervals, $\tilde{b}_{j-1},\tilde{a}_{j}$, of the support, $J_{f}$, of 
the associated equilibrium measure, $\mu_{\widetilde{V}}^{f}$. 
The $4(N \! + \! 1)$ contour integrals appearing in 
Equation~\eqref{eqtilrzo1} will be evaluated via the Cauchy and 
Residue Theorems.\footnote{In particular, for $j \! \in \! \lbrace 
1,2,\dotsc,N \! + \! 1 \rbrace$ and $r \! \in \! \lbrace \tilde{b}_{j
-1},\tilde{a}_{j} \rbrace$, $\oint_{\partial \tilde{\mathbb{U}}_{
\tilde{\delta}_{r}}} \frac{1}{\xi -z} \, \frac{\md \xi}{2 \pi \mi} \! = \! 
\begin{cases}
0, &\text{$z \! \in \! \operatorname{ext}(\partial \tilde{\mathbb{
U}}_{\tilde{\delta}_{r}}) \! := \! \lbrace \mathstrut z^{\prime} \! \in 
\! \mathbb{C}; \, \oint_{\partial \tilde{\mathbb{U}}_{\tilde{\delta}_{r}}} 
\frac{1}{u-z^{\prime}} \, \frac{\md u}{2 \pi \mi} \! = \! 0 \rbrace$,} \\
1, &\text{$z \! \in \! \operatorname{int}(\partial \tilde{\mathbb{U}}_{
\tilde{\delta}_{r}}) \! := \! \lbrace \mathstrut z^{\prime} \! \in \! 
\mathbb{C}; \, \oint_{\partial \tilde{\mathbb{U}}_{\tilde{\delta}_{r}}} 
\frac{1}{u-z^{\prime}} \, \frac{\md u}{2 \pi \mi} \! \neq \! 0 \rbrace$.}
\end{cases}$} For $n \! \in \! \mathbb{N}$ and $k \! \in \! \lbrace 
1,2,\dotsc,K \rbrace$ such that $\alpha_{p_{\mathfrak{s}}} \! := \! 
\alpha_{k} \! \neq \! \infty$, via the asymptotics~\eqref{eqtlvee13}, 
\eqref{eqtlvee14}, \eqref{eqproptilb}, and~\eqref{eqproptila}, and an 
inordinately lengthy calculation analogous to that of the proof of 
Proposition~\ref{propo5.1}, application of the Cauchy and Residue 
Theorems shows that,\footnote{Recall that, for $j \! \in \! \lbrace 
1,2,\dotsc,N \! + \! 1 \rbrace$, $(\tilde{\mathbb{U}}_{\tilde{\delta}_{
\tilde{b}_{j-1}}} \cup \partial \tilde{\mathbb{U}}_{\tilde{\delta}_{
\tilde{b}_{j-1}}}) \cap \lbrace \alpha_{k} \rbrace \! = \! \varnothing 
\! = \! (\tilde{\mathbb{U}}_{\tilde{\delta}_{\tilde{a}_{j}}} \cup \partial 
\tilde{\mathbb{U}}_{\tilde{\delta}_{\tilde{a}_{j}}}) \cap \lbrace \alpha_{k} 
\rbrace$.} uniformly for compact subsets of $\mathbb{C} \setminus 
\tilde{\Sigma}_{\tilde{\mathcal{R}}}^{\sharp}$ $(\ni \! z)$, for $j \! \in 
\! \lbrace 1,2,\dotsc,N \! + \! 1 \rbrace$,
\begin{align}
\oint_{\partial \tilde{\mathbb{U}}_{\tilde{\delta}_{\tilde{b}_{j-1}}}} 
\dfrac{w^{\tilde{\Sigma}_{\circlearrowright}}_{+}(\xi)}{\xi \! - \! 
\alpha_{k}} \, \dfrac{\md \xi}{2 \pi \mi} =& \, \operatorname{Res} 
\left(\dfrac{w^{\tilde{\Sigma}_{\circlearrowright}}_{+}(z)}{z \! - \! 
\alpha_{k}};\tilde{b}_{j-1} \right) \underset{\underset{z_{o}=1+o(1)}{
\mathscr{N},n \to \infty}}{=} \dfrac{((\alpha_{k} \! - \! \tilde{b}_{j-1}) 
\tilde{\alpha}_{1}(\tilde{b}_{j-1}) \! - \! \tilde{\alpha}_{0}(\tilde{b}_{j-
1}))}{((n \! - \! 1)K \! + \! k)(\alpha_{k} \! - \! \tilde{b}_{j-1})^{2}
(\tilde{\alpha}_{0}(\tilde{b}_{j-1}))^{2}} \tilde{\boldsymbol{\mathrm{A}}}
(\tilde{b}_{j-1}) \nonumber \\
-& \, \dfrac{(\tilde{\alpha}_{0}(\tilde{b}_{j-1}))^{-1}}{((n \! - \! 1)K \! + 
\! k)(\alpha_{k} \! - \! \tilde{b}_{j-1})} \tilde{\boldsymbol{\mathrm{B}}}
(\tilde{b}_{j-1}) \! + \! \mathcal{O} \left(\dfrac{\tilde{\mathfrak{c}}_{
\tilde{\Sigma}_{\circlearrowright}}^{\triangledown}(n,k,z_{o};j)}{((n \! 
- \! 1)K \! + \! k)^{2}} \right), \label{eqtilrzo2} \\
\oint_{\partial \tilde{\mathbb{U}}_{\tilde{\delta}_{\tilde{a}_{j}}}} 
\dfrac{w^{\tilde{\Sigma}_{\circlearrowright}}_{+}(\xi)}{\xi \! - \! 
\alpha_{k}} \, \dfrac{\md \xi}{2 \pi \mi} =& \, \operatorname{Res} 
\left(\dfrac{w^{\tilde{\Sigma}_{\circlearrowright}}_{+}(z)}{z \! - \! 
\alpha_{k}};\tilde{a}_{j} \right) \underset{\underset{z_{o}=1+o(1)}{
\mathscr{N},n \to \infty}}{=} \dfrac{((\alpha_{k} \! - \! \tilde{a}_{j}) 
\tilde{\alpha}_{1}(\tilde{a}_{j}) \! - \! \tilde{\alpha}_{0}(\tilde{a}_{j}))}{
((n \! - \! 1)K \! + \! k)(\alpha_{k} \! - \! \tilde{a}_{j})^{2}(\tilde{
\alpha}_{0}(\tilde{a}_{j}))^{2}} \tilde{\boldsymbol{\mathrm{A}}}
(\tilde{a}_{j}) \nonumber \\
-& \, \dfrac{(\tilde{\alpha}_{0}(\tilde{a}_{j}))^{-1}}{((n \! - \! 1)K \! + 
\! k)(\alpha_{k} \! - \! \tilde{a}_{j})} \tilde{\boldsymbol{\mathrm{B}}}
(\tilde{a}_{j}) \! + \! \mathcal{O} \left(\dfrac{\tilde{\mathfrak{c}}_{
\tilde{\Sigma}_{\circlearrowright}}^{\triangle}(n,k,z_{o};j)}{((n \! - \! 
1)K \! + \! k)^{2}} \right), \label{eqtilrzo3}
\end{align}
\begin{equation} \label{eqtilrzo4}
\oint_{\partial \tilde{\mathbb{U}}_{\tilde{\delta}_{\tilde{b}_{j-1}}}} 
\dfrac{w^{\tilde{\Sigma}_{\circlearrowright}}_{+}(\xi)}{\xi \! - \! z} 
\, \dfrac{\md \xi}{2 \pi \mi} \underset{\underset{z_{o}=1+o(1)}{
\mathscr{N},n \to \infty}}{=}  
\begin{cases}
{\fontsize{9pt}{9pt}\selectfont 
\begin{aligned}[b]
&-\frac{(\tilde{\alpha}_{0}(\tilde{b}_{j-1}))^{-2}}{((n \! - \! 1)K \! + 
\! k)(z \! - \! \tilde{b}_{j-1})} \left(\tilde{\alpha}_{0}(\tilde{b}_{j-1}) 
\tilde{\boldsymbol{\mathrm{B}}}(\tilde{b}_{j-1}) \! - \! \tilde{\alpha}_{1}
(\tilde{b}_{j-1}) \tilde{\boldsymbol{\mathrm{A}}}(\tilde{b}_{j-1}) \right) \\
&-\frac{(\tilde{\alpha}_{0}(\tilde{b}_{j-1}))^{-1}}{((n \! - \! 1)K \! + \! k)
(z \! - \! \tilde{b}_{j-1})^{2}} \tilde{\boldsymbol{\mathrm{A}}}(\tilde{b}_{j
-1}) \! + \! \mathcal{O} \left(\frac{\sum_{m=1}^{3} \tilde{\mathfrak{c}}_{
\tilde{\Sigma}_{\circlearrowright},m}^{\blacktriangle}(n,k,z_{o};j)
(z \! - \! \tilde{b}_{j-1})^{-m}}{((n \! - \! 1)K \! + \! k)^{2}} \right),
\end{aligned} \, \, \, \, \text{$z \! \in \! \mathbb{C} \setminus 
(\tilde{\mathbb{U}}_{\tilde{\delta}_{\tilde{b}_{j-1}}} \cup \partial 
\tilde{\mathbb{U}}_{\tilde{\delta}_{\tilde{b}_{j-1}}})$,}} \\
{\fontsize{9pt}{9pt}\selectfont 
\begin{aligned}[b]
&-\frac{(\tilde{\alpha}_{0}(\tilde{b}_{j-1}))^{-2}}{((n \! - \! 1)K \! + 
\! k)(z \! - \! \tilde{b}_{j-1})} \left(\tilde{\alpha}_{0}(\tilde{b}_{j-1}) 
\tilde{\boldsymbol{\mathrm{B}}}(\tilde{b}_{j-1}) \! - \! \tilde{\alpha}_{1}
(\tilde{b}_{j-1}) \tilde{\boldsymbol{\mathrm{A}}}(\tilde{b}_{j-1}) \right) \\
&-\frac{(\tilde{\alpha}_{0}(\tilde{b}_{j-1}))^{-1}}{((n \! - \! 1)K \! + \! 
k)(z \! - \! \tilde{b}_{j-1})^{2}} \tilde{\boldsymbol{\mathrm{A}}}(\tilde{
b}_{j-1}) \! + \! \frac{1}{((n \! - \! 1)K \! + \! k)} \tilde{\mathbb{Y}}_{
\tilde{b}_{j-1}}(z) \\
&+\mathcal{O} \left(\frac{\sum_{m=0}^{\infty} \tilde{\mathfrak{c}}_{
\tilde{\Sigma}_{\circlearrowright},m}^{\blacktriangledown}(n,k,z_{o};
j)(z \! - \! \tilde{b}_{j-1})^{m}}{((n \! - \! 1)K \! + \! k)^{2}} \right),
\end{aligned} \, \, \, \, \, \text{$z \! \in \! \tilde{\mathbb{U}}_{
\tilde{\delta}_{\tilde{b}_{j-1}}}$,}}
\end{cases}
\end{equation}
\begin{equation} \label{eqtilrzo5}
\oint_{\partial \tilde{\mathbb{U}}_{\tilde{\delta}_{\tilde{a}_{j}}}} 
\dfrac{w^{\tilde{\Sigma}_{\circlearrowright}}_{+}(\xi)}{\xi \! - \! z} 
\, \dfrac{\md \xi}{2 \pi \mi} \underset{\underset{z_{o}=1+o(1)}{
\mathscr{N},n \to \infty}}{=}  
\begin{cases}
{\fontsize{9pt}{9pt}\selectfont 
\begin{aligned}[b]
&-\frac{(\tilde{\alpha}_{0}(\tilde{a}_{j}))^{-2}}{((n \! - \! 1)K \! + \! k)
(z \! - \! \tilde{a}_{j})} \left(\tilde{\alpha}_{0}(\tilde{a}_{j}) \tilde{
\boldsymbol{\mathrm{B}}}(\tilde{a}_{j}) \! - \! \tilde{\alpha}_{1}
(\tilde{a}_{j}) \tilde{\boldsymbol{\mathrm{A}}}(\tilde{a}_{j}) \right) \\
&-\frac{(\tilde{\alpha}_{0}(\tilde{a}_{j}))^{-1}}{((n \! - \! 1)K \! + \! k)
(z \! - \! \tilde{a}_{j})^{2}} \tilde{\boldsymbol{\mathrm{A}}}(\tilde{a}_{j}) 
\! + \! \mathcal{O} \left(\frac{\sum_{m=1}^{3} \tilde{\mathfrak{c}}_{
\tilde{\Sigma}_{\circlearrowright},m}^{\clubsuit}(n,k,z_{o};j)
(z \! - \! \tilde{a}_{j})^{-m}}{((n \! - \! 1)K \! + \! k)^{2}} \right),
\end{aligned} \, \, \, \, \text{$z \! \in \! \mathbb{C} \setminus 
(\tilde{\mathbb{U}}_{\tilde{\delta}_{\tilde{a}_{j}}} \cup \partial 
\tilde{\mathbb{U}}_{\tilde{\delta}_{\tilde{a}_{j}}})$,}} \\
{\fontsize{9pt}{9pt}\selectfont 
\begin{aligned}[b]
&-\frac{(\tilde{\alpha}_{0}(\tilde{a}_{j}))^{-2}}{((n \! - \! 1)K \! + \! k)
(z \! - \! \tilde{a}_{j})} \left(\tilde{\alpha}_{0}(\tilde{a}_{j}) \tilde{
\boldsymbol{\mathrm{B}}}(\tilde{a}_{j}) \! - \! \tilde{\alpha}_{1}
(\tilde{a}_{j}) \tilde{\boldsymbol{\mathrm{A}}}(\tilde{a}_{j}) \right) \\
&-\frac{(\tilde{\alpha}_{0}(\tilde{a}_{j}))^{-1}}{((n \! - \! 1)K \! + \! k)
(z \! - \! \tilde{a}_{j})^{2}} \tilde{\boldsymbol{\mathrm{A}}}(\tilde{a}_{j}) 
\! + \! \frac{1}{((n \! - \! 1)K \! + \! k)} \tilde{\mathbb{Y}}_{\tilde{a}_{j}}
(z) \\
&+\mathcal{O} \left(\frac{\sum_{m=0}^{\infty} \tilde{\mathfrak{c}}_{
\tilde{\Sigma}_{\circlearrowright},m}^{\spadesuit}(n,k,z_{o};j)
(z \! - \! \tilde{a}_{j})^{m}}{((n \! - \! 1)K \! + \! k)^{2}} \right),
\end{aligned} \, \, \, \, \, \text{$z \! \in \! \tilde{\mathbb{U}}_{
\tilde{\delta}_{\tilde{a}_{j}}}$,}}
\end{cases}
\end{equation}
where $\tilde{\boldsymbol{\mathrm{A}}}(\tilde{b}_{j-1})$, $\tilde{
\boldsymbol{\mathrm{A}}}(\tilde{a}_{j})$, $\tilde{\boldsymbol{\mathrm{B}}}
(\tilde{b}_{j-1})$, $\tilde{\boldsymbol{\mathrm{B}}}(\tilde{a}_{j})$, 
$\tilde{\alpha}_{0}(\tilde{b}_{j-1})$, $\tilde{\alpha}_{0}(\tilde{a}_{j})$, 
$\tilde{\alpha}_{1}(\tilde{b}_{j-1})$, $\tilde{\alpha}_{1}(\tilde{a}_{j})$, 
$\tilde{\mathbb{Y}}_{\tilde{b}_{j-1}}(z)$, and $\tilde{\mathbb{Y}}_{
\tilde{a}_{j}}(z)$ are defined in item~\pmb{(2)} of the proposition, 
$(\mathrm{M}_{2}(\mathbb{C}) \! \ni)$ $\tilde{\mathfrak{c}}_{
\tilde{\Sigma}_{\circlearrowright}}^{r_{1}}(n,k,z_{o};j) \! =_{\underset{z_{o}
=1+o(1)}{\mathscr{N},n \to \infty}} \! \mathcal{O}(1)$, $r_{1} \! \in \! 
\lbrace \triangledown,\triangle \rbrace$, and $(\mathrm{M}_{2}
(\mathbb{C}) \! \ni)$ $\tilde{\mathfrak{c}}_{\tilde{\Sigma}_{
\circlearrowright},m}^{r_{2}}(n,k,z_{o};j) \! =_{\underset{z_{o}=1+
o(1)}{\mathscr{N},n \to \infty}} \! \mathcal{O}(1)$, $r_{2} \! \in \! 
\lbrace \blacktriangle,\blacktriangledown,\clubsuit,\spadesuit 
\rbrace$; hence, via Equation~\eqref{eqtilrzo1} and the 
Estimates~\eqref{eqtilrzo2}--\eqref{eqtilrzo5}, one arrives at, for 
$n \! \in \! \mathbb{N}$ and $k \! \in \! \lbrace 1,2,\dotsc,K \rbrace$ 
such that $\alpha_{p_{\mathfrak{s}}} \! := \! \alpha_{k} \! \neq \! 
\infty$, uniformly for compact subsets of $\mathbb{C} \setminus 
\tilde{\Sigma}_{\tilde{\mathcal{R}}}^{\sharp}$ $(\ni \! z)$, the 
Asymptotics~\eqref{eqpropo5.3tB} for $\tilde{\mathcal{R}}(z)$.

The analysis for the case $n \! \in \! \mathbb{N}$ and $k \! \in \! 
\lbrace 1,2,\dotsc,K \rbrace$ such that $\alpha_{p_{\mathfrak{s}}} 
\! := \! \alpha_{k} \! = \! \infty$ is, \emph{mutatis mutandis}, 
analogous, and leads to the (uniform for compact subsets of 
$\mathbb{C} \setminus \hat{\Sigma}_{\hat{\mathcal{R}}}^{\sharp}$ 
$(\ni \! z))$ Asymptotics~\eqref{eqpropo5.3hA} for 
$\hat{\mathcal{R}}(z)$. \hfill $\qed$
\begin{eeeee} \label{rem5.3} 
\textsl{For $n \! \in \! \mathbb{N}$ and $k \! \in \! \lbrace 1,2,
\dotsc,K \rbrace$ such that $\alpha_{p_{\mathfrak{s}}} \! := \! 
\alpha_{k} \! = \! \infty$ (resp., $\alpha_{p_{\mathfrak{s}}} 
\! := \! \alpha_{k} \! \neq \! \infty)$, a perusal of the 
Asymptotics~\eqref{eqpropo5.3hA} (resp., \eqref{eqpropo5.3tB}$)$ 
for $\hat{\mathcal{R}}(z)$ (resp., $\tilde{\mathcal{R}}(z))$ stated 
in item~{\rm \pmb{(1)}} (resp., item~{\rm \pmb{(2)}}$)$ of 
Proposition~\ref{propo5.3} seems to suggest, at first glance, that there 
are second-order poles at $\lbrace \hat{b}_{j-1},\hat{a}_{j} \rbrace_{j
=1}^{N+1}$ (resp., $\lbrace \tilde{b}_{j-1},\tilde{a}_{j} \rbrace_{j=
1}^{N+1})$, which, of course, can not be the case, as a careful 
analysis of the $\mathrm{M}_{2}(\mathbb{C})$-valued factors 
$\hat{\mathbb{Y}}_{\hat{b}_{j-1}}(z)$ and $\hat{\mathbb{Y}}_{
\hat{a}_{j}}(z)$ (resp., $\tilde{\mathbb{Y}}_{\tilde{b}_{j-1}}(z)$ and 
$\tilde{\mathbb{Y}}_{\tilde{a}_{j}}(z))$, $j \! \in \! \lbrace 1,2,\dotsc,
N \! + \! 1 \rbrace$, reveals; in fact, expanding, in the double-scaling 
limit $\mathscr{N},n \! \to \! \infty$ such that $z_{o} \! = \! 1 \! + \! 
o(1)$, $\hat{\mathbb{Y}}_{\hat{b}_{j-1}}(z)$ as $z \! \to \! \hat{b}_{j-
1}$ and $\hat{\mathbb{Y}}_{\hat{a}_{j}}(z)$ as $z \! \to \! \hat{a}_{j}$ 
(resp., $\tilde{\mathbb{Y}}_{\tilde{b}_{j-1}}(z)$ as $z \! \to \! \tilde{
b}_{j-1}$ and $\tilde{\mathbb{Y}}_{\tilde{a}_{j}}(z)$ as $z \! \to \! 
\tilde{a}_{j})$, $j \! \in \! \lbrace 1,2,\dotsc,N \! + \! 1 \rbrace$, as 
in the paradigm of the proof of Proposition~\ref{propo5.1}, one 
shows that: {\rm (i)} for $n \! \in \! \mathbb{N}$ and $k \! \in \! 
\lbrace 1,2,\dotsc,K \rbrace$ such that $\alpha_{p_{\mathfrak{s}}} 
\! := \! \alpha_{k} \! = \! \infty$, uniformly for compact subsets of 
$\mathbb{C} \setminus \hat{\Sigma}_{\hat{\mathcal{R}}}^{\sharp}$ 
$(\ni \! z)$, for $j \! \in \! \lbrace 1,2,\dotsc,N \! + \! 1 \rbrace$,
\begin{align*}
&-\dfrac{1}{((n \! - \! 1)K \! + \! k)} \left(\dfrac{(\hat{\alpha}_{0}
(\hat{b}_{j-1}))^{-1}}{(z \! - \! \hat{b}_{j-1})^{2}} \hat{\boldsymbol{
\mathrm{A}}}(\hat{b}_{j-1}) \! + \! \dfrac{(\hat{\alpha}_{0}(\hat{b}_{j
-1}))^{-2}}{z \! - \! \hat{b}_{j-1}} \left(\hat{\alpha}_{0}(\hat{b}_{j-1}) 
\hat{\boldsymbol{\mathrm{B}}}(\hat{b}_{j-1}) \! - \! \hat{\alpha}_{1}
(\hat{b}_{j-1}) \hat{\boldsymbol{\mathrm{A}}}(\hat{b}_{j-1}) \right) 
\! + \! \hat{\mathbb{Y}}_{\hat{b}_{j-1}}(z) \right) \\
&\underset{\underset{z_{o}=1+o(1)}{\mathscr{N},n \to \infty}}{=} 
\dfrac{(\hat{\alpha}_{0}(\hat{b}_{j-1}))^{-2}}{((n \! - \! 1)K \! + \! k)} 
\left(\hat{\alpha}_{0}(\hat{b}_{j-1}) \left(\left(\dfrac{\hat{\alpha}_{1}
(\hat{b}_{j-1})}{\hat{\alpha}_{0}(\hat{b}_{j-1})} \right)^{2} \! - \! \dfrac{
\hat{\alpha}_{2}(\hat{b}_{j-1})}{\hat{\alpha}_{0}(\hat{b}_{j-1})} \right) 
\hat{\boldsymbol{\mathrm{A}}}(\hat{b}_{j-1}) \! - \! \hat{\alpha}_{1}
(\hat{b}_{j-1}) \hat{\boldsymbol{\mathrm{B}}}(\hat{b}_{j-1}) \! + \! 
\hat{\alpha}_{0}(\hat{b}_{j-1}) \hat{\boldsymbol{\mathrm{C}}}
(\hat{b}_{j-1}) \right) \\
&+\mathcal{O} \left(\dfrac{\sum_{m \in \mathbb{N}} \hat{\mathfrak{
c}}_{m}^{\curlyvee}(n,k,z_{o};j)(z \! - \! \hat{b}_{j-1})^{m}}{(n \! - \! 1)
K \! + \! k} \right), \quad z \! \in \! \hat{\mathbb{U}}_{\hat{\delta}_{
\hat{b}_{j-1}}},
\end{align*}
and
\begin{align*}
&-\dfrac{1}{((n \! - \! 1)K \! + \! k)} \left(\dfrac{(\hat{\alpha}_{0}
(\hat{a}_{j}))^{-1}}{(z \! - \! \hat{a}_{j})^{2}} \hat{\boldsymbol{
\mathrm{A}}}(\hat{a}_{j}) \! + \! \dfrac{(\hat{\alpha}_{0}(\hat{a}_{j})
)^{-2}}{z \! - \! \hat{a}_{j}} \left(\hat{\alpha}_{0}(\hat{a}_{j}) \hat{
\boldsymbol{\mathrm{B}}}(\hat{a}_{j}) \! - \! \hat{\alpha}_{1}
(\hat{a}_{j}) \hat{\boldsymbol{\mathrm{A}}}(\hat{a}_{j}) \right) \! 
+ \! \hat{\mathbb{Y}}_{\hat{a}_{j}}(z) \right) \\
&\underset{\underset{z_{o}=1+o(1)}{\mathscr{N},n \to \infty}}{=} 
\dfrac{(\hat{\alpha}_{0}(\hat{a}_{j}))^{-2}}{((n \! - \! 1)K \! + \! k)} 
\left(\hat{\alpha}_{0}(\hat{a}_{j}) \left(\left(\dfrac{\hat{\alpha}_{1}
(\hat{a}_{j})}{\hat{\alpha}_{0}(\hat{a}_{j})} \right)^{2} \! - \! \dfrac{
\hat{\alpha}_{2}(\hat{a}_{j})}{\hat{\alpha}_{0}(\hat{a}_{j})} \right) 
\hat{\boldsymbol{\mathrm{A}}}(\hat{a}_{j}) \! - \! \hat{\alpha}_{1}
(\hat{a}_{j}) \hat{\boldsymbol{\mathrm{B}}}(\hat{a}_{j}) \! + \! 
\hat{\alpha}_{0}(\hat{a}_{j}) \hat{\boldsymbol{\mathrm{C}}}(\hat{a}_{j}) 
\right) \\
&+\mathcal{O} \left(\dfrac{\sum_{m \in \mathbb{N}} \hat{\mathfrak{
c}}_{m}^{\curlywedge}(n,k,z_{o};j)(z \! - \! \hat{a}_{j})^{m}}{(n \! - \! 1)
K \! + \! k} \right), \quad z \! \in \! \hat{\mathbb{U}}_{\hat{\delta}_{
\hat{a}_{j}}},
\end{align*}
where $\hat{\boldsymbol{\mathrm{A}}}(\hat{b}_{j-1})$, $\hat{
\boldsymbol{\mathrm{A}}}(\hat{a}_{j})$, $\hat{\boldsymbol{
\mathrm{B}}}(\hat{b}_{j-1})$, $\hat{\boldsymbol{\mathrm{B}}}
(\hat{a}_{j})$, $\hat{\boldsymbol{\mathrm{C}}}(\hat{b}_{j-1})$, 
$\hat{\boldsymbol{\mathrm{C}}}(\hat{a}_{j})$, $\hat{\alpha}_{0}
(\hat{b}_{j-1})$, $\hat{\alpha}_{0}(\hat{a}_{j})$, $\hat{\alpha}_{1}
(\hat{b}_{j-1})$, $\hat{\alpha}_{1}(\hat{a}_{j})$, $\hat{\alpha}_{2}
(\hat{b}_{j-1})$, and $\hat{\alpha}_{2}(\hat{a}_{j})$ are defined 
in item~{\rm \pmb{(1)}} of Proposition~\ref{propo5.1}, and 
$(\mathrm{M}_{2}(\mathbb{C}) \! \ni)$ $\hat{\mathfrak{c}}_{m}^{r}
(n,k,z_{o};j) \! =_{\underset{z_{o}=1+o(1)}{\mathscr{N},n \to \infty}} 
\! \mathcal{O}(1)$, $r_{1} \! \in \! \lbrace \curlyvee,\curlywedge 
\rbrace$$;$ and {\rm (ii)} for $n \! \in \! \mathbb{N}$ and $k \! \in 
\! \lbrace 1,2,\dotsc,K \rbrace$ such that $\alpha_{p_{\mathfrak{s}}} 
\! := \! \alpha_{k} \! \neq \! \infty$, uniformly for compact subsets 
of $\mathbb{C} \setminus \tilde{\Sigma}_{\tilde{\mathcal{R}}}^{\sharp}$ 
$(\ni \! z)$, for $j \! \in \! \lbrace 1,2,\dotsc,N \! + \! 1 \rbrace$,
\begin{align*}
&-\dfrac{1}{((n \! - \! 1)K \! + \! k)} \left(\dfrac{(\tilde{\alpha}_{0}
(\tilde{b}_{j-1}))^{-1}}{(z \! - \! \tilde{b}_{j-1})^{2}} \tilde{\boldsymbol{
\mathrm{A}}}(\tilde{b}_{j-1}) \! + \! \dfrac{(\tilde{\alpha}_{0}(\tilde{b}_{j
-1}))^{-2}}{z \! - \! \tilde{b}_{j-1}} \left(\tilde{\alpha}_{0}(\tilde{b}_{j-1}) 
\tilde{\boldsymbol{\mathrm{B}}}(\tilde{b}_{j-1}) \! - \! \tilde{\alpha}_{1}
(\tilde{b}_{j-1}) \tilde{\boldsymbol{\mathrm{A}}}(\tilde{b}_{j-1}) \right) 
\! + \! \tilde{\mathbb{Y}}_{\tilde{b}_{j-1}}(z) \right) \\
&\underset{\underset{z_{o}=1+o(1)}{\mathscr{N},n \to \infty}}{=} 
\dfrac{(\tilde{\alpha}_{0}(\tilde{b}_{j-1}))^{-2}}{((n \! - \! 1)K \! + \! k)} 
\left(\tilde{\alpha}_{0}(\tilde{b}_{j-1}) \left(\left(\dfrac{\tilde{\alpha}_{1}
(\tilde{b}_{j-1})}{\tilde{\alpha}_{0}(\tilde{b}_{j-1})} \right)^{2} \! - \! 
\dfrac{\tilde{\alpha}_{2}(\tilde{b}_{j-1})}{\tilde{\alpha}_{0}(\tilde{b}_{j-1})} 
\right) \tilde{\boldsymbol{\mathrm{A}}}(\tilde{b}_{j-1}) \! - \! \tilde{
\alpha}_{1}(\tilde{b}_{j-1}) \tilde{\boldsymbol{\mathrm{B}}}(\tilde{b}_{j-1}) 
\! + \! \tilde{\alpha}_{0}(\tilde{b}_{j-1}) \tilde{\boldsymbol{\mathrm{C}}}
(\tilde{b}_{j-1}) \right) \\
&+\mathcal{O} \left(\dfrac{\sum_{m \in \mathbb{N}} \tilde{\mathfrak{
c}}_{m}^{\curlyvee}(n,k,z_{o};j)(z \! - \! \tilde{b}_{j-1})^{m}}{(n \! - \! 1)
K \! + \! k} \right), \quad z \! \in \! \tilde{\mathbb{U}}_{\tilde{\delta}_{
\tilde{b}_{j-1}}},
\end{align*}
and
\begin{align*}
&-\dfrac{1}{((n \! - \! 1)K \! + \! k)} \left(\dfrac{(\tilde{\alpha}_{0}
(\tilde{a}_{j}))^{-1}}{(z \! - \! \tilde{a}_{j})^{2}} \tilde{\boldsymbol{
\mathrm{A}}}(\tilde{a}_{j}) \! + \! \dfrac{(\tilde{\alpha}_{0}(\tilde{a}_{j})
)^{-2}}{z \! - \! \tilde{a}_{j}} \left(\tilde{\alpha}_{0}(\tilde{a}_{j}) 
\tilde{\boldsymbol{\mathrm{B}}}(\tilde{a}_{j}) \! - \! \tilde{\alpha}_{1}
(\tilde{a}_{j}) \tilde{\boldsymbol{\mathrm{A}}}(\tilde{a}_{j}) \right) \! 
+ \! \tilde{\mathbb{Y}}_{\tilde{a}_{j}}(z) \right) \\
&\underset{\underset{z_{o}=1+o(1)}{\mathscr{N},n \to \infty}}{=} 
\dfrac{(\tilde{\alpha}_{0}(\tilde{a}_{j}))^{-2}}{((n \! - \! 1)K \! + \! k)} 
\left(\tilde{\alpha}_{0}(\tilde{a}_{j}) \left(\left(\dfrac{\tilde{\alpha}_{1}
(\tilde{a}_{j})}{\tilde{\alpha}_{0}(\tilde{a}_{j})} \right)^{2} \! - \! \dfrac{
\tilde{\alpha}_{2}(\tilde{a}_{j})}{\tilde{\alpha}_{0}(\tilde{a}_{j})} \right) 
\tilde{\boldsymbol{\mathrm{A}}}(\tilde{a}_{j}) \! - \! \tilde{\alpha}_{1}
(\tilde{a}_{j}) \tilde{\boldsymbol{\mathrm{B}}}(\tilde{a}_{j}) \! + \! 
\tilde{\alpha}_{0}(\tilde{a}_{j}) \tilde{\boldsymbol{\mathrm{C}}}
(\tilde{a}_{j}) \right) \\
&+\mathcal{O} \left(\dfrac{\sum_{m \in \mathbb{N}} \tilde{\mathfrak{
c}}_{m}^{\curlywedge}(n,k,z_{o};j)(z \! - \! \tilde{a}_{j})^{m}}{(n \! - \! 1)
K \! + \! k} \right), \quad z \! \in \! \tilde{\mathbb{U}}_{\tilde{\delta}_{
\tilde{a}_{j}}},
\end{align*}
where $\tilde{\boldsymbol{\mathrm{A}}}(\tilde{b}_{j-1})$, $\tilde{
\boldsymbol{\mathrm{A}}}(\tilde{a}_{j})$, $\tilde{\boldsymbol{
\mathrm{B}}}(\tilde{b}_{j-1})$, $\tilde{\boldsymbol{\mathrm{B}}}
(\tilde{a}_{j})$, $\tilde{\boldsymbol{\mathrm{C}}}(\tilde{b}_{j-1})$, 
$\tilde{\boldsymbol{\mathrm{C}}}(\tilde{a}_{j})$, $\tilde{\alpha}_{0}
(\tilde{b}_{j-1})$, $\tilde{\alpha}_{0}(\tilde{a}_{j})$, $\tilde{\alpha}_{1}
(\tilde{b}_{j-1})$, $\tilde{\alpha}_{1}(\tilde{a}_{j})$, $\tilde{\alpha}_{2}
(\tilde{b}_{j-1})$, and $\tilde{\alpha}_{2}(\tilde{a}_{j})$ are defined 
in item~{\rm \pmb{(2)}} of Proposition~\ref{propo5.1}, and 
$(\mathrm{M}_{2}(\mathbb{C}) \! \ni)$ $\tilde{\mathfrak{c}}_{m}^{r}
(n,k,z_{o};j) \! =_{\underset{z_{o}=1+o(1)}{\mathscr{N},n \to \infty}} 
\! \mathcal{O}(1)$, $r_{1} \! \in \! \lbrace \curlyvee,\curlywedge 
\rbrace$.}
\end{eeeee}
\begin{ccccc} \label{lem5.5} 
Let the external field $\widetilde{V} \colon \overline{\mathbb{R}} 
\setminus \lbrace \alpha_{1},\alpha_{2},\dotsc,\alpha_{K} \rbrace \! 
\to \! \mathbb{R}$ satisfy conditions~\eqref{eq20}--\eqref{eq22} and 
be regular. For $n \! \in \! \mathbb{N}$ and $k \! \in \! \lbrace 1,2,
\dotsc,K \rbrace$ such that $\alpha_{p_{\mathfrak{s}}} \! := \! 
\alpha_{k} \! = \! \infty$ (resp., $\alpha_{p_{\mathfrak{s}}} \! := \! 
\alpha_{k} \! \neq \! \infty)$, let the associated equilibrium measure, 
$\mu_{\widetilde{V}}^{\infty}$ (resp., $\mu_{\widetilde{V}}^{f})$, and its 
support, $J_{\infty}$ (resp., $J_{f})$, be as described in item~$\pmb{(1)}$ 
(resp., item~$\pmb{(2)})$ of Lemma~\ref{lem3.7}, and, along with the 
corresponding variational constant, $\hat{\ell}$ (resp., $\tilde{\ell})$, 
satisfy the variational conditions~\eqref{eql3.8a} (resp., 
conditions~\eqref{eql3.8b}$)$$;$ moreover, let the associated 
conditions~{\rm (i)}--{\rm (iv)} of item~$\pmb{(1)}$ (resp., 
item~$\pmb{(2)})$ of Lemma~\ref{lem3.8} be valid. For $n \! \in \! 
\mathbb{N}$ and $k \! \in \! \lbrace 1,2,\dotsc,K \rbrace$, let 
$\mathcal{X} \colon \mathbb{N} \times \lbrace 1,2,\dotsc,K \rbrace 
\times \overline{\mathbb{C}} \setminus \overline{\mathbb{R}} \! \to 
\! \mathrm{SL}_{2}(\mathbb{C})$ be the unique solution of the monic 
{\rm MPC ORF RHP} $(\mathcal{X}(z),\upsilon (z),\overline{\mathbb{R}})$ 
stated in Lemma~$\bm{\mathrm{RHP}_{\mathrm{MPC}}}$. Then$:$ 
{\rm \pmb{(1)}} for $n \! \in \! \mathbb{N}$ and $k \! \in \! \lbrace 
1,2,\dotsc,K \rbrace$ such that $\alpha_{p_{\mathfrak{s}}} \! := \! 
\alpha_{k} \! = \! \infty$,
\begin{align}
\mathcal{X}(z)z^{-\varkappa_{nk} \sigma_{3}} \underset{z \to \alpha_{k}}{=}& 
\, \mathrm{I} \! + \! \me^{\frac{n \hat{\ell}}{2} \mathrm{ad}(\sigma_{3})} 
\left(\dfrac{1}{z} \left(\hat{w}_{0}^{\sharp} \sigma_{3} \! + \! \hat{A}_{0}^{
\sharp}(\alpha_{k}) \! + \! \hat{\mathcal{R}}_{0}^{\hat{A}_{0}^{\sharp}}
(\alpha_{k}) \right) \! + \! \dfrac{1}{z^{2}} \left(\hat{w}_{1}^{\sharp} 
\sigma_{3} \! + \! \dfrac{1}{2}(\hat{w}_{0}^{\sharp})^{2} \mathrm{I} 
\! + \! \hat{w}_{0}^{\sharp} \left(\hat{A}_{0}^{\sharp}(\alpha_{k}) 
\! + \! \hat{\mathcal{R}}_{0}^{\hat{A}_{0}^{\sharp}}(\alpha_{k}) 
\right) \sigma_{3} \right. \right. \nonumber \\
+&\left. \left. \, \hat{B}_{0}^{\sharp}(\alpha_{k}) \! + \! 
\hat{\mathcal{R}}_{0}^{\hat{A}_{0}^{\sharp}}(\alpha_{k}) 
\hat{A}_{0}^{\sharp}(\alpha_{k}) \! + \! \hat{\mathcal{R}}_{0}^{
\hat{B}_{0}^{\sharp}}(\alpha_{k}) \right) \! + \! \mathcal{O} 
\left(\hat{\mathfrak{c}}_{\mathfrak{s}}^{\mathcal{X}}(n,k,z_{o})
z^{-3} \right) \vphantom{M^{M^{M^{M^{M^{M}}}}}} \right), 
\label{eqlem5.5A}
\end{align}
where
\begin{align*}
\hat{w}_{r_{1}}^{\sharp} :=& \, \dfrac{1}{1 \! + \! r_{1}} \left(\sum_{j=
1}^{\mathfrak{s}-1} \varkappa_{nk \tilde{k}_{j}}(\alpha_{p_{j}})^{1+
r_{1}} \! - \! ((n \! - \! 1)K \! + \! k) \int_{J_{\infty}} \xi^{1+r_{1}} \, 
\md \mu_{\widetilde{V}}^{\infty}(\xi) \right), \quad r_{1} \! = \! 0,1, \\
\hat{\mathcal{R}}_{0}^{\hat{A}_{0}^{\sharp}}(\alpha_{k}) \underset{
\underset{z_{o}=1+o(1)}{\mathscr{N},n \to \infty}}{=}& \, \dfrac{1}{((n 
\! - \! 1)K \! + \! k)} \sum_{j=1}^{N+1} \left((\hat{\alpha}_{0}(\hat{b}_{j
-1}))^{-2} \left(\hat{\alpha}_{0}(\hat{b}_{j-1}) \hat{\boldsymbol{
\mathrm{B}}}(\hat{b}_{j-1}) \! - \! \hat{\alpha}_{1}(\hat{b}_{j-1}) 
\hat{\boldsymbol{\mathrm{A}}}(\hat{b}_{j-1}) \right) \right. \\
+&\left. \, (\hat{\alpha}_{0}(\hat{a}_{j}))^{-2} \left(\hat{\alpha}_{0}
(\hat{a}_{j}) \hat{\boldsymbol{\mathrm{B}}}(\hat{a}_{j}) \! - \! \hat{
\alpha}_{1}(\hat{a}_{j}) \hat{\boldsymbol{\mathrm{A}}}(\hat{a}_{j}) \right) 
\right) \! + \! \mathcal{O} \left(\dfrac{\hat{\mathfrak{c}}_{\hat{\mathcal{
R}}_{0}^{\hat{A}_{0}^{\sharp}}}(n,k,z_{o})}{((n \! - \! 1)K \! + \! k)^{2}} 
\right), \\
\hat{\mathcal{R}}_{0}^{\hat{B}_{0}^{\sharp}}(\alpha_{k}) \underset{
\underset{z_{o}=1+o(1)}{\mathscr{N},n \to \infty}}{=}& \, \dfrac{1}{
((n \! - \! 1)K \! + \! k)} \sum_{j=1}^{N+1} \left(\hat{b}_{j-1}
(\hat{\alpha}_{0}(\hat{b}_{j-1}))^{-2} \left(\hat{\alpha}_{0}(\hat{b}_{j-1}) 
\hat{\boldsymbol{\mathrm{B}}}(\hat{b}_{j-1}) \! - \! \hat{\alpha}_{1}
(\hat{b}_{j-1}) \hat{\boldsymbol{\mathrm{A}}}(\hat{b}_{j-1}) \right) 
\right. \\
+&\left. \, \hat{a}_{j}(\hat{\alpha}_{0}(\hat{a}_{j}))^{-2} \left(
\hat{\alpha}_{0}(\hat{a}_{j}) \hat{\boldsymbol{\mathrm{B}}}
(\hat{a}_{j}) \! - \! \hat{\alpha}_{1}(\hat{a}_{j}) \hat{\boldsymbol{
\mathrm{A}}}(\hat{a}_{j}) \right) \! + \! \dfrac{1}{\hat{\alpha}_{0}
(\hat{b}_{j-1})} \hat{\boldsymbol{\mathrm{A}}}(\hat{b}_{j-1}) \! + \! 
\dfrac{1}{\hat{\alpha}_{0}(\hat{a}_{j})} \hat{\boldsymbol{\mathrm{A}}}
(\hat{a}_{j}) \right) \\
+& \, \mathcal{O} \left(\dfrac{\hat{\mathfrak{c}}_{\hat{\mathcal{
R}}_{0}^{\hat{B}_{0}^{\sharp}}}(n,k,z_{o})}{((n \! - \! 1)K \! + \! k)^{2}} 
\right),
\end{align*}
with $\hat{\boldsymbol{\mathrm{A}}}(\hat{b}_{j-1})$, $\hat{
\boldsymbol{\mathrm{A}}}(\hat{a}_{j})$, $\hat{\boldsymbol{
\mathrm{B}}}(\hat{b}_{j-1})$, $\hat{\boldsymbol{\mathrm{B}}}
(\hat{a}_{j})$, $\hat{\alpha}_{0}(\hat{b}_{j-1})$, $\hat{\alpha}_{0}
(\hat{a}_{j})$, $\hat{\alpha}_{1}(\hat{b}_{j-1})$, and $\hat{\alpha}_{1}
(\hat{a}_{j})$, $j \! \in \! \lbrace 1,2,\dotsc,N \! + \! 1 \rbrace$, 
defined in item~{\rm \pmb{(1)}} of Proposition~\ref{propo5.1}, 
and $(\mathrm{M}_{2}(\mathbb{C}) \! \ni)$ $\hat{\mathfrak{c}}_{
\hat{\mathcal{R}}_{0}^{r_{2}}}(n,k,z_{o}) \! =_{\underset{z_{o}=1+
o(1)}{\mathscr{N},n \to \infty}} \! \mathcal{O}(1)$, $r_{2} \! \in 
\! \lbrace \hat{A}_{0}^{\sharp},\hat{B}_{0}^{\sharp} \rbrace$,
\begin{gather*}
\hat{A}_{0}^{\sharp}(\alpha_{k}) \! = \!  
\tilde{\mathfrak{m}}^{\raise-0.5ex\hbox{$\scriptstyle \infty$}} 
\begin{pmatrix}
\hat{\mathbb{M}}_{11}^{1,+}(\alpha_{k}) & \mi \hat{\alpha}_{0}^{
\triangledown} \hat{\mathbb{M}}_{12}^{0,+}(\alpha_{k}) \\
-\mi \hat{\alpha}_{0}^{\triangledown} \hat{\mathbb{M}}_{21}^{0,+}
(\alpha_{k}) & \hat{\mathbb{M}}_{22}^{1,+}(\alpha_{k})
\end{pmatrix}, \\
\hat{B}_{0}^{\sharp}(\alpha_{k}) \! = \! 
\tilde{\mathfrak{m}}^{\raise-0.5ex\hbox{$\scriptstyle \infty$}} 
\begin{pmatrix}
\hat{\mathbb{M}}_{11}^{2,+}(\alpha_{k}) \! + \! \frac{1}{2}
(\hat{\alpha}_{0}^{\triangledown})^{2} \hat{\mathbb{M}}_{11}^{0,+}
(\alpha_{k}) & \mi \hat{\alpha}_{0}^{\triangledown} \hat{\mathbb{
M}}_{12}^{1,+}(\alpha_{k}) \! + \! \mi \hat{\beta}_{0}^{\triangledown} 
\hat{\mathbb{M}}_{12}^{0,+}(\alpha_{k}) \\
-\mi \hat{\alpha}_{0}^{\triangledown} \hat{\mathbb{M}}_{21}^{1,+}
(\alpha_{k}) \! - \! \mi \hat{\beta}_{0}^{\triangledown} \hat{\mathbb{
M}}_{21}^{0,+}(\alpha_{k}) & \hat{\mathbb{M}}_{22}^{2,+}(\alpha_{k}) 
\! + \! \frac{1}{2}(\hat{\alpha}_{0}^{\triangledown})^{2} 
\hat{\mathbb{M}}_{22}^{0,+}(\alpha_{k})
\end{pmatrix},
\end{gather*}
with $\tilde{\mathfrak{m}}^{\raise-0.5ex\hbox{$\scriptstyle \infty$}}$ 
given in Equation~\eqref{eqmaininf8},
\begin{equation*}
\hat{\alpha}_{0}^{\triangledown} \! := \! \dfrac{1}{4} \sum_{j=1}^{N+1}
(\hat{a}_{j} \! - \! \hat{b}_{j-1}), \qquad \qquad \hat{\beta}_{0}^{
\triangledown} \! := \! \dfrac{1}{8} \sum_{j=1}^{N+1}(\hat{a}_{j}^{2} \! 
- \! \hat{b}_{j-1}^{2}),
\end{equation*}
\begin{align*}
\hat{\mathbb{M}}_{11}^{0,+}(\alpha_{k}) =& \, \dfrac{\hat{\mathbb{F}}_{
\hat{\boldsymbol{\Omega}}}^{0,0,+}(1,1)}{\hat{\mathbb{F}}_{
\boldsymbol{0}}^{0,0,+}(1,1)}, \quad \qquad \hat{\mathbb{M}}_{11}^{1,+}
(\alpha_{k}) = \dfrac{\hat{\mathbb{F}}_{\hat{\boldsymbol{\Omega}}}^{0,0,+}
(1,1)}{\hat{\mathbb{F}}_{\boldsymbol{0}}^{0,0,+}(1,1)} \left(\hat{\eta}_{0}^{
\blacklozenge,+}(1,1;\hat{\boldsymbol{\Omega}}) \! - \! \hat{\eta}_{0}^{
\blacklozenge,+}(1,1;\boldsymbol{0}) \right), \\
\hat{\mathbb{M}}_{11}^{2,+}(\alpha_{k}) =& \, \dfrac{\hat{\mathbb{F}}_{
\hat{\boldsymbol{\Omega}}}^{0,0,+}(1,1)}{\hat{\mathbb{F}}_{
\boldsymbol{0}}^{0,0,+}(1,1)} \left(\hat{\eta}_{1}^{\blacklozenge,+}
(1,1;\hat{\boldsymbol{\Omega}}) \! - \! \hat{\eta}_{1}^{\blacklozenge,+}
(1,1;\boldsymbol{0}) \! + \! (\hat{\eta}_{0}^{\blacklozenge,+}(1,1;
\boldsymbol{0}))^{2} \! - \! \hat{\eta}_{0}^{\blacklozenge,+}(1,1;
\hat{\boldsymbol{\Omega}}) \hat{\eta}_{0}^{\blacklozenge,+}(1,1;
\boldsymbol{0}) \right), \\
\hat{\mathbb{M}}_{12}^{0,+}(\alpha_{k}) =& \, \dfrac{\hat{\mathbb{F}}_{
\hat{\boldsymbol{\Omega}}}^{0,0,+}(-1,1)}{\hat{\mathbb{F}}_{
\boldsymbol{0}}^{0,0,+}(-1,1)}, \quad \qquad \hat{\mathbb{M}}_{12}^{1,+}
(\alpha_{k}) = -\dfrac{\hat{\mathbb{F}}_{\hat{\boldsymbol{\Omega}}}^{0,0,+}
(-1,1)}{\hat{\mathbb{F}}_{\boldsymbol{0}}^{0,0,+}(-1,1)} \left(\hat{\eta}_{0}^{
\blacklozenge,+}(-1,1;\hat{\boldsymbol{\Omega}}) \! - \! \hat{\eta}_{0}^{
\blacklozenge,+}(-1,1;\boldsymbol{0}) \right), \\
\hat{\mathbb{M}}_{12}^{2,+}(\alpha_{k}) =& \, -\dfrac{\hat{\mathbb{F}}_{
\hat{\boldsymbol{\Omega}}}^{0,0,+}(-1,1)}{\hat{\mathbb{F}}_{
\boldsymbol{0}}^{0,0,+}(-1,1)} \left(\hat{\eta}_{1}^{\blacklozenge,+}
(-1,1;\hat{\boldsymbol{\Omega}}) \! - \! \hat{\eta}_{1}^{\blacklozenge,+}
(-1,1;\boldsymbol{0}) \! - \! (\hat{\eta}_{0}^{\blacklozenge,+}(-1,1;
\boldsymbol{0}))^{2} \! + \! \hat{\eta}_{0}^{\blacklozenge,+}(-1,1;
\hat{\boldsymbol{\Omega}}) \hat{\eta}_{0}^{\blacklozenge,+}(-1,1;
\boldsymbol{0}) \right), \\
\hat{\mathbb{M}}_{21}^{0,+}(\alpha_{k}) =& \, \dfrac{\hat{\mathbb{F}}_{
\hat{\boldsymbol{\Omega}}}^{0,0,+}(1,-1)}{\hat{\mathbb{F}}_{
\boldsymbol{0}}^{0,0,+}(1,-1)}, \quad \qquad \hat{\mathbb{M}}_{21}^{1,+}
(\alpha_{k}) = \dfrac{\hat{\mathbb{F}}_{\hat{\boldsymbol{\Omega}}}^{0,0,+}
(1,-1)}{\hat{\mathbb{F}}_{\boldsymbol{0}}^{0,0,+}(1,-1)} \left(\hat{\eta}_{0}^{
\blacklozenge,+}(1,-1;\hat{\boldsymbol{\Omega}}) \! - \! \hat{\eta}_{0}^{
\blacklozenge,+}(1,-1;\boldsymbol{0}) \right), \\
\hat{\mathbb{M}}_{21}^{2,+}(\alpha_{k}) =& \, \dfrac{\hat{\mathbb{F}}_{
\hat{\boldsymbol{\Omega}}}^{0,0,+}(1,-1)}{\hat{\mathbb{F}}_{
\boldsymbol{0}}^{0,0,+}(1,-1)} \left(\hat{\eta}_{1}^{\blacklozenge,+}
(1,-1;\hat{\boldsymbol{\Omega}}) \! - \! \hat{\eta}_{1}^{\blacklozenge,+}
(1,-1;\boldsymbol{0}) \! + \! (\hat{\eta}_{0}^{\blacklozenge,+}(1,-1;
\boldsymbol{0}))^{2} \! - \! \hat{\eta}_{0}^{\blacklozenge,+}(1,-1;
\hat{\boldsymbol{\Omega}}) \hat{\eta}_{0}^{\blacklozenge,+}(1,-1;
\boldsymbol{0}) \right), \\
\hat{\mathbb{M}}_{22}^{0,+}(\alpha_{k}) =& \, \dfrac{\hat{\mathbb{F}}_{
\hat{\boldsymbol{\Omega}}}^{0,0,+}(-1,-1)}{\hat{\mathbb{F}}_{
\boldsymbol{0}}^{0,0,+}(-1,-1)}, \quad \qquad \hat{\mathbb{M}}_{22}^{1,+}
(\alpha_{k}) = -\dfrac{\hat{\mathbb{F}}_{\hat{\boldsymbol{\Omega}}}^{0,0,+}
(-1,-1)}{\hat{\mathbb{F}}_{\boldsymbol{0}}^{0,0,+}(-1,-1)} \left(\hat{\eta}_{0}^{
\blacklozenge,+}(-1,-1;\hat{\boldsymbol{\Omega}}) \! - \! \hat{\eta}_{0}^{
\blacklozenge,+}(-1,-1;\boldsymbol{0}) \right), \\
\hat{\mathbb{M}}_{22}^{2,+}(\alpha_{k}) =& \, -\dfrac{\hat{\mathbb{F}}_{
\hat{\boldsymbol{\Omega}}}^{0,0,+}(-1,-1)}{\hat{\mathbb{F}}_{
\boldsymbol{0}}^{0,0,+}(-1,-1)} \left(\hat{\eta}_{1}^{\blacklozenge,+}
(-1,-1;\hat{\boldsymbol{\Omega}}) \! - \! \hat{\eta}_{1}^{\blacklozenge,+}
(-1,-1;\boldsymbol{0}) \! - \! (\hat{\eta}_{0}^{\blacklozenge,+}(-1,-1;
\boldsymbol{0}))^{2} \right. \\
+&\left. \, \hat{\eta}_{0}^{\blacklozenge,+}(-1,-1;\hat{\boldsymbol{
\Omega}}) \hat{\eta}_{0}^{\blacklozenge,+}(-1,-1;\boldsymbol{0}) \right),
\end{align*}
where, for $\varepsilon_{1},\varepsilon_{2} \! = \! \pm 1$,
\begin{gather*}
\hat{\mathbb{F}}_{\hat{\boldsymbol{\Omega}}}^{j_{1},j_{2},\pm}
(\varepsilon_{1},\varepsilon_{2}) \! := \! \sum_{m \in \mathbb{Z}^{N}} 
\me^{2 \pi \mi (m,\varepsilon_{1} \hat{\boldsymbol{u}}_{\pm}(\infty)-
\frac{1}{2 \pi}((n-1)K+k) \hat{\boldsymbol{\Omega}}+ \varepsilon_{2} 
\hat{\boldsymbol{d}})+\mi \pi (m,\hat{\boldsymbol{\tau}}m)}
(\hat{\lambda}_{2}^{\sharp})^{j_{1}}(\hat{\lambda}_{3}^{\sharp})^{j_{2}}, 
\quad j_{1} \! = \! 0,1,2, \quad j_{2} \! = \! 0,1, \\
\hat{\lambda}_{2}^{\sharp} \! = \! -2 \pi \mi \sum_{j=1}^{N}m_{j} 
\hat{c}_{j1}, \qquad \, \, \hat{\lambda}_{3}^{\sharp} \! = \! -\mi \pi 
\sum_{j=1}^{N}m_{j}(\hat{c}_{j2} \! - \! \hat{\lambda}_{0}^{\sharp} 
\hat{c}_{j1}), \qquad \, \, \hat{\lambda}_{0}^{\sharp} \! = \! 
-\dfrac{1}{2} \sum_{j=1}^{N+1}(\hat{b}_{j-1} \! + \! \hat{a}_{j}), 
\end{gather*}
$\hat{c}_{j1}$ and $\hat{c}_{j2}$, $j \! \in \! \lbrace 1,2,\dotsc,N 
\rbrace$, are obtained {}from Equations~\eqref{E1} and~\eqref{E2},
\begin{gather*}
\hat{\eta}_{0}^{\blacklozenge,\pm}(\varepsilon_{1},\varepsilon_{2};
\hat{\boldsymbol{\Omega}}) \! := \! \left(\hat{\mathbb{F}}_{\hat{
\boldsymbol{\Omega}}}^{0,0,\pm}(\varepsilon_{1},\varepsilon_{2}) 
\right)^{-1} \hat{\mathbb{F}}_{\hat{\boldsymbol{\Omega}}}^{1,0,\pm}
(\varepsilon_{1},\varepsilon_{2}), \\
\hat{\eta}_{1}^{\blacklozenge,\pm}(\varepsilon_{1},\varepsilon_{2};
\hat{\boldsymbol{\Omega}}) \! := \! \left(\hat{\mathbb{F}}_{\hat{
\boldsymbol{\Omega}}}^{0,0,\pm}(\varepsilon_{1},\varepsilon_{2}) 
\right)^{-1} \left(\hat{\mathbb{F}}_{\hat{\boldsymbol{\Omega}}}^{0,1,
\pm}(\varepsilon_{1},\varepsilon_{2}) \! \pm \! \dfrac{\varepsilon_{1}}{2} 
\hat{\mathbb{F}}_{\hat{\boldsymbol{\Omega}}}^{2,0,\pm}(\varepsilon_{1},
\varepsilon_{2}) \right),
\end{gather*}
and $(\mathrm{M}_{2}(\mathbb{C}) \! \ni)$ $\hat{\mathfrak{c}}_{
\mathfrak{s}}^{\mathcal{X}}(n,k,z_{o}) \! =_{\underset{z_{o}=1+o(1)}{
\mathscr{N},n \to \infty}} \! \mathcal{O}(1)$, and, for $q \! \in \! \lbrace 
1,2,\dotsc,\mathfrak{s} \! - \! 1 \rbrace$,\footnote{Note: (i) if, for $q \! 
\in \! \lbrace 1,2,\dotsc,\mathfrak{s} \! - \! 1 \rbrace$ and $j \! \in \! 
\lbrace 1,2,\dotsc,N \rbrace$, $\alpha_{p_{q}} \! \in \! (\hat{a}_{j},
\hat{b}_{j})$, then $J_{\infty} \cap \mathbb{R}_{\alpha_{p_{q}}}^{>} \! 
= \! \cup_{i=j}^{N}[\hat{b}_{i},\hat{a}_{i+1}]$, and $\int_{J_{\infty} \cap 
\mathbb{R}_{\alpha_{p_{q}}}^{>}} \md \mu_{\widetilde{V}}^{\infty}(\xi) 
\! = \! \int_{\hat{b}_{j}}^{\hat{a}_{N+1}} \md \mu_{\widetilde{V}}^{\infty}
(\xi) \! = \! \hat{\Omega}_{j}/2 \pi$; (ii) if, for $q \! \in \! \lbrace 1,2,
\dotsc,\mathfrak{s} \! - \! 1 \rbrace$, $\alpha_{p_{q}} \! \in \! (-\infty,
\hat{b}_{0})$, then $J_{\infty} \cap \mathbb{R}_{\alpha_{p_{q}}}^{>} \! = 
\! J_{\infty}$, and $\int_{J_{\infty} \cap \mathbb{R}_{\alpha_{p_{q}}}^{>}} 
\md \mu_{\widetilde{V}}^{\infty}(\xi) \! = \! 1$; and (iii) if, for $q \! \in 
\! \lbrace 1,2,\dotsc,\mathfrak{s} \! - \! 1 \rbrace$, $\alpha_{p_{q}} 
\! \in \! (\hat{a}_{N+1},+\infty)$, then $J_{\infty} \cap \mathbb{R}_{
\alpha_{p_{q}}}^{>} \! = \! \varnothing$, and $\int_{J_{\infty} \cap 
\mathbb{R}_{\alpha_{p_{q}}}^{>}} \md \mu_{\widetilde{V}}^{\infty}
(\xi) \! = \! 0$.}
\begin{align}
\mathcal{X}(z)(z \! - \! \alpha_{p_{q}})^{\varkappa_{nk \tilde{k}_{q}} 
\sigma_{3}} \underset{z \to \alpha_{p_{q}}}{=}& \, \me^{\frac{n \hat{\ell}}{2} 
\mathrm{ad}(\sigma_{3})} \left(\left(\mathrm{I} \! + \! \hat{\mathcal{
R}}_{0}^{\hat{A}_{0}}(\alpha_{p_{q}}) \right) \hat{A}_{0}(\alpha_{p_{q}}) 
\! + \! \left(\hat{w}_{0}^{\triangledown} \left(\mathrm{I} \! + \! 
\hat{\mathcal{R}}_{0}^{\hat{A}_{0}}(\alpha_{p_{q}}) \right) \hat{A}_{0}
(\alpha_{p_{q}}) \sigma_{3} \! + \! \left(\mathrm{I} \! + \! \hat{
\mathcal{R}}_{0}^{\hat{A}_{0}}(\alpha_{p_{q}}) \right) \right. \right. 
\nonumber \\
\times&\left. \left. \, \hat{B}_{0}(\alpha_{p_{q}}) \! + \! 
\hat{\mathcal{R}}_{0}^{\hat{B}_{0}}(\alpha_{p_{q}}) \hat{A}_{0}
(\alpha_{p_{q}}) \right)(z \! - \! \alpha_{p_{q}}) \! + \! \left(\left(
\mathrm{I} \! + \! \hat{\mathcal{R}}_{0}^{\hat{A}_{0}}(\alpha_{p_{q}}) 
\right) \hat{A}_{0}(\alpha_{p_{q}}) \left(\hat{w}_{1}^{\triangledown} 
\sigma_{3} \! + \! \dfrac{1}{2}(\hat{w}_{0}^{\triangledown})^{2} 
\mathrm{I} \right) \right. \right. \nonumber \\
+&\left. \left. \, \hat{w}_{0}^{\triangledown} \left(\left(\mathrm{I} 
\! + \! \hat{\mathcal{R}}_{0}^{\hat{A}_{0}}(\alpha_{p_{q}}) \right) 
\hat{B}_{0}(\alpha_{p_{q}}) \! + \! \hat{\mathcal{R}}_{0}^{\hat{B}_{0}}
(\alpha_{p_{q}}) \hat{A}_{0}(\alpha_{p_{q}}) \right) \sigma_{3} \! + \! 
\left(\mathrm{I} \! + \! \hat{\mathcal{R}}_{0}^{\hat{A}_{0}}(\alpha_{p_{q}}) 
\right) \hat{C}_{0}(\alpha_{p_{q}}) \right. \right. \nonumber \\
+&\left. \left. \, \hat{\mathcal{R}}_{0}^{\hat{B}_{0}}(\alpha_{p_{q}}) 
\hat{B}_{0}(\alpha_{p_{q}}) \! + \! \hat{\mathcal{R}}_{0}^{\hat{C}_{0}}
(\alpha_{p_{q}}) \hat{A}_{0}(\alpha_{p_{q}}) \right)(z \! - \! \alpha_{
p_{q}})^{2} \! + \! \mathcal{O} \left(\hat{\mathfrak{c}}_{q}^{
\mathcal{X}}(n,k,z_{o})(z \! - \! \alpha_{p_{q}})^{3} \right) 
\vphantom{M^{M^{M^{M^{M^{M^{M}}}}}}} \right) \nonumber \\
\times& \, (-1)^{\sum_{j \in \tilde{\Delta}_{\infty}(q)} \varkappa_{nk 
\tilde{k}_{j}}\sigma_{3}} \me^{\hat{\Xi}_{0} \sigma_{3}} \me^{\mi \pi 
((n-1)K+k) \int_{J_{\infty} \cap \mathbb{R}_{\alpha_{p_{q}}}^{>}} \md 
\mu_{\widetilde{V}}^{\infty}(\xi) \sigma_{3}}, \label{eqlem5.5B}
\end{align}
where
\begin{align*}
\hat{w}_{r_{3}}^{\triangledown} :=& \, \dfrac{1}{1 \! + \! r_{3}} \left(
\sum_{\substack{j=1\\j \neq q}}^{\mathfrak{s}-1} \dfrac{\varkappa_{nk 
\tilde{k}_{j}}}{(\alpha_{p_{j}} \! - \! \alpha_{p_{q}})^{1+r_{3}}} \! - \! ((n \! 
- \! 1)K \! + \! k) \int_{J_{\infty}}(\xi \! - \! \alpha_{p_{q}})^{-(1+r_{3})} \, 
\md \mu_{\widetilde{V}}^{\infty}(\xi) \right), \quad r_{3} \! = \! 0,1, \\
\hat{\mathcal{R}}_{0}^{\hat{A}_{0}}(\alpha_{p_{q}}) \underset{\underset{
z_{o}=1+o(1)}{\mathscr{N},n \to \infty}}{=}& \, \dfrac{1}{((n \! - \! 1)K 
\! + \! k)} \sum_{j=1}^{N+1} \left(-\dfrac{(\hat{\alpha}_{0}(\hat{b}_{j-
1}))^{-2}}{\hat{b}_{j-1} \! - \! \alpha_{p_{q}}} \left(\hat{\alpha}_{0}
(\hat{b}_{j-1}) \hat{\boldsymbol{\mathrm{B}}}(\hat{b}_{j-1}) \! - \! 
\hat{\alpha}_{1}(\hat{b}_{j-1}) \hat{\boldsymbol{\mathrm{A}}}
(\hat{b}_{j-1}) \right) \right. \\
-&\left. \, \dfrac{(\hat{\alpha}_{0}(\hat{a}_{j}))^{-2}}{\hat{a}_{j} \! - \! 
\alpha_{p_{q}}} \left(\hat{\alpha}_{0}(\hat{a}_{j}) \hat{\boldsymbol{
\mathrm{B}}}(\hat{a}_{j}) \! - \! \hat{\alpha}_{1}(\hat{a}_{j}) \hat{
\boldsymbol{\mathrm{A}}}(\hat{a}_{j}) \right) \! + \! \dfrac{(\hat{
\alpha}_{0}(\hat{b}_{j-1}))^{-1}}{(\hat{b}_{j-1} \! - \! \alpha_{p_{q}})^{2}} 
\hat{\boldsymbol{\mathrm{A}}}(\hat{b}_{j-1}) \! + \! \dfrac{(\hat{
\alpha}_{0}(\hat{a}_{j}))^{-1}}{(\hat{a}_{j} \! - \! \alpha_{p_{q}})^{2}} 
\hat{\boldsymbol{\mathrm{A}}}(\hat{a}_{j}) \right) \\
+& \, \mathcal{O} \left(\dfrac{\hat{\mathfrak{c}}_{\hat{\mathcal{
R}}_{0}^{\hat{A}_{0}}}(n,k,z_{o})}{((n \! - \! 1)K \! + \! k)^{2}} \right), \\
\hat{\mathcal{R}}_{0}^{\hat{B}_{0}}(\alpha_{p_{q}}) \underset{\underset{
z_{o}=1+o(1)}{\mathscr{N},n \to \infty}}{=}& \, \dfrac{1}{((n \! - \! 1)K 
\! + \! k)} \sum_{j=1}^{N+1} \left(-\dfrac{(\hat{\alpha}_{0}(\hat{b}_{j-
1}))^{-2}}{(\hat{b}_{j-1} \! - \! \alpha_{p_{q}})^{2}} \left(\hat{\alpha}_{0}
(\hat{b}_{j-1}) \hat{\boldsymbol{\mathrm{B}}}(\hat{b}_{j-1}) \! - \! 
\hat{\alpha}_{1}(\hat{b}_{j-1}) \hat{\boldsymbol{\mathrm{A}}}
(\hat{b}_{j-1}) \right) \right. \\
-&\left. \, \dfrac{(\hat{\alpha}_{0}(\hat{a}_{j}))^{-2}}{(\hat{a}_{j} \! - \! 
\alpha_{p_{q}})^{2}} \left(\hat{\alpha}_{0}(\hat{a}_{j}) \hat{\boldsymbol{
\mathrm{B}}}(\hat{a}_{j}) \! - \! \hat{\alpha}_{1}(\hat{a}_{j}) \hat{
\boldsymbol{\mathrm{A}}}(\hat{a}_{j}) \right) \! + \! \dfrac{2(\hat{
\alpha}_{0}(\hat{b}_{j-1}))^{-1}}{(\hat{b}_{j-1} \! - \! \alpha_{p_{q}})^{3}} 
\hat{\boldsymbol{\mathrm{A}}}(\hat{b}_{j-1}) \! + \! \dfrac{2(\hat{
\alpha}_{0}(\hat{a}_{j}))^{-1}}{(\hat{a}_{j} \! - \! \alpha_{p_{q}})^{3}} 
\hat{\boldsymbol{\mathrm{A}}}(\hat{a}_{j}) \right) \\
+& \, \mathcal{O} \left(\dfrac{\hat{\mathfrak{c}}_{\hat{\mathcal{
R}}_{0}^{\hat{B}_{0}}}(n,k,z_{o})}{((n \! - \! 1)K \! + \! k)^{2}} \right), \\
\hat{\mathcal{R}}_{0}^{\hat{C}_{0}}(\alpha_{p_{q}}) \underset{\underset{
z_{o}=1+o(1)}{\mathscr{N},n \to \infty}}{=}& \, \dfrac{1}{((n \! - \! 1)K 
\! + \! k)} \sum_{j=1}^{N+1} \left(-\dfrac{(\hat{\alpha}_{0}(\hat{b}_{j-
1}))^{-2}}{(\hat{b}_{j-1} \! - \! \alpha_{p_{q}})^{3}} \left(\hat{\alpha}_{0}
(\hat{b}_{j-1}) \hat{\boldsymbol{\mathrm{B}}}(\hat{b}_{j-1}) \! - \! 
\hat{\alpha}_{1}(\hat{b}_{j-1}) \hat{\boldsymbol{\mathrm{A}}}
(\hat{b}_{j-1}) \right) \right. \\
-&\left. \, \dfrac{(\hat{\alpha}_{0}(\hat{a}_{j}))^{-2}}{(\hat{a}_{j} \! - \! 
\alpha_{p_{q}})^{3}} \left(\hat{\alpha}_{0}(\hat{a}_{j}) \hat{\boldsymbol{
\mathrm{B}}}(\hat{a}_{j}) \! - \! \hat{\alpha}_{1}(\hat{a}_{j}) \hat{
\boldsymbol{\mathrm{A}}}(\hat{a}_{j}) \right) \! + \! \dfrac{3(\hat{
\alpha}_{0}(\hat{b}_{j-1}))^{-1}}{(\hat{b}_{j-1} \! - \! \alpha_{p_{q}})^{4}} 
\hat{\boldsymbol{\mathrm{A}}}(\hat{b}_{j-1}) \! + \! \dfrac{3(\hat{
\alpha}_{0}(\hat{a}_{j}))^{-1}}{(\hat{a}_{j} \! - \! \alpha_{p_{q}})^{4}} 
\hat{\boldsymbol{\mathrm{A}}}(\hat{a}_{j}) \right) \\
+& \, \mathcal{O} \left(\dfrac{\hat{\mathfrak{c}}_{\hat{\mathcal{
R}}_{0}^{\hat{C}_{0}}}(n,k,z_{o})}{((n \! - \! 1)K \! + \! k)^{2}} \right),
\end{align*}
with $(\mathrm{M}_{2}(\mathbb{C}) \! \ni)$ $\hat{\mathfrak{c}}_{
\hat{\mathcal{R}}_{0}^{r_{4}}}(n,k,z_{o}) \! =_{\underset{z_{o}=1
+o(1)}{\mathscr{N},n \to \infty}} \! \mathcal{O}(1)$, $r_{4} \! \in 
\! \lbrace \hat{A}_{0},\hat{B}_{0},\hat{C}_{0} \rbrace$,
\begin{gather*}
\hat{A}_{0}(\alpha_{p_{q}}) \! = \!  
\tilde{\mathfrak{m}}^{\raise-0.5ex\hbox{$\scriptstyle \infty$}} 
,
\end{gather*}
with
\begin{gather*}
\hat{\mathfrak{h}}_{11}^{0,+} \! := \! \dfrac{1}{2} \left(\hat{\gamma}
(\alpha_{p_{q}}) \! + \! (\hat{\gamma}(\alpha_{p_{q}}))^{-1} \right), 
\qquad \, \, \hat{\mathfrak{h}}_{11}^{1,+} \! := \! \dfrac{\hat{\alpha}^{
\triangledown}}{2} \left(\hat{\gamma}(\alpha_{p_{q}}) \! - \! 
(\hat{\gamma}(\alpha_{p_{q}}))^{-1} \right), \\
\hat{\mathfrak{h}}_{11}^{2,+} \! := \! \dfrac{1}{2} \left(\hat{\beta}^{
\triangledown} \left(\hat{\gamma}(\alpha_{p_{q}}) \! - \! (\hat{\gamma}
(\alpha_{p_{q}}))^{-1} \right) \! + \! \dfrac{(\hat{\alpha}^{\triangledown}
)^{2}}{2} \left(\hat{\gamma}(\alpha_{p_{q}}) \! + \! (\hat{\gamma}
(\alpha_{p_{q}}))^{-1} \right) \right), \\
\hat{\mathfrak{h}}_{12}^{0,+} \! := \! \dfrac{\mi}{2} \left(\hat{\gamma}
(\alpha_{p_{q}}) \! - \! (\hat{\gamma}(\alpha_{p_{q}}))^{-1} \right), \qquad 
\, \, \hat{\mathfrak{h}}_{12}^{1,+} \! := \! \dfrac{\mi \hat{\alpha}^{
\triangledown}}{2} \left(\hat{\gamma}(\alpha_{p_{q}}) \! + \! 
(\hat{\gamma}(\alpha_{p_{q}}))^{-1} \right), \\
\hat{\mathfrak{h}}_{12}^{2,+} \! := \! \dfrac{\mi}{2} \left(\hat{\beta}^{
\triangledown} \left(\hat{\gamma}(\alpha_{p_{q}}) \! + \! (\hat{\gamma}
(\alpha_{p_{q}}))^{-1} \right) \! + \! \dfrac{(\hat{\alpha}^{\triangledown}
)^{2}}{2} \left(\hat{\gamma}(\alpha_{p_{q}}) \! - \! (\hat{\gamma}
(\alpha_{p_{q}}))^{-1} \right) \right),
\end{gather*}
where $\hat{\gamma}(\alpha_{p_{q}})$ is defined by Equation~\eqref{eqssabra1},
\begin{equation*}
\hat{\alpha}^{\triangledown} \! = \! \dfrac{1}{4} \sum_{j=1}^{N+1} 
\left(\dfrac{1}{\hat{a}_{j} \! - \! \alpha_{p_{q}}} \! - \! \dfrac{1}{\hat{b}_{j
-1} \! - \! \alpha_{p_{q}}} \right), \qquad \, \, \hat{\beta}^{\triangledown} 
\! = \! \dfrac{1}{8} \sum_{j=1}^{N+1} \left(\dfrac{1}{(\hat{a}_{j} \! - \! 
\alpha_{p_{q}})^{2}} \! - \! \dfrac{1}{(\hat{b}_{j-1} \! - \! \alpha_{p_{q}}
)^{2}} \right),
\end{equation*}
\begin{align*}
\hat{\mathbb{M}}_{11}^{0,+}(\alpha_{p_{q}}) =& \, \dfrac{\hat{\mathbb{G}}_{
\hat{\boldsymbol{\Omega}}}^{0,0,+}(1,1)}{\hat{\mathbb{G}}_{\boldsymbol{0}}^{
0,0,+}(1,1)}, \quad \qquad \hat{\mathbb{M}}_{11}^{1,+}(\alpha_{p_{q}}) = 
\dfrac{\hat{\mathbb{G}}_{\hat{\boldsymbol{\Omega}}}^{0,0,+}(1,1)}{\hat{
\mathbb{G}}_{\boldsymbol{0}}^{0,0,+}(1,1)} \left(\hat{\eta}_{0}^{\spadesuit,+}
(1,1;\hat{\boldsymbol{\Omega}}) \! - \! \hat{\eta}_{0}^{\spadesuit,+}(1,1;
\boldsymbol{0}) \right), \\
\hat{\mathbb{M}}_{11}^{2,+}(\alpha_{p_{q}}) =& \, \dfrac{\hat{\mathbb{G}}_{
\hat{\boldsymbol{\Omega}}}^{0,0,+}(1,1)}{\hat{\mathbb{G}}_{\boldsymbol{0}}^{
0,0,+}(1,1)} \left(\hat{\eta}_{1}^{\spadesuit,+}(1,1;\hat{\boldsymbol{\Omega}}) 
\! - \! \hat{\eta}_{1}^{\spadesuit,+}(1,1;\boldsymbol{0}) \! + \! (\hat{\eta}_{0}^{
\spadesuit,+}(1,1;\boldsymbol{0}))^{2} \! - \! \hat{\eta}_{0}^{\spadesuit,+}(1,1;
\hat{\boldsymbol{\Omega}}) \hat{\eta}_{0}^{\spadesuit,+}(1,1;\boldsymbol{0}) 
\right), \\
\hat{\mathbb{M}}_{12}^{0,+}(\alpha_{p_{q}}) =& \, \dfrac{\hat{\mathbb{G}}_{
\hat{\boldsymbol{\Omega}}}^{0,0,+}(-1,1)}{\hat{\mathbb{G}}_{\boldsymbol{0}}^{
0,0,+}(-1,1)}, \quad \qquad \hat{\mathbb{M}}_{12}^{1,+}(\alpha_{p_{q}}) = 
-\dfrac{\hat{\mathbb{G}}_{\hat{\boldsymbol{\Omega}}}^{0,0,+}(-1,1)}{\hat{
\mathbb{G}}_{\boldsymbol{0}}^{0,0,+}(-1,1)} \left(\hat{\eta}_{0}^{\spadesuit,+}
(-1,1;\hat{\boldsymbol{\Omega}}) \! - \! \hat{\eta}_{0}^{\spadesuit,+}(-1,1;
\boldsymbol{0}) \right), \\
\hat{\mathbb{M}}_{12}^{2,+}(\alpha_{p_{q}}) =& \, -\dfrac{\hat{\mathbb{G}}_{
\hat{\boldsymbol{\Omega}}}^{0,0,+}(-1,1)}{\hat{\mathbb{G}}_{\boldsymbol{0}}^{
0,0,+}(-1,1)} \left(\hat{\eta}_{1}^{\spadesuit,+}(-1,1;\hat{\boldsymbol{\Omega}}) 
\! - \! \hat{\eta}_{1}^{\spadesuit,+}(-1,1;\boldsymbol{0}) \! - \! (\hat{\eta}_{0}^{
\spadesuit,+}(-1,1;\boldsymbol{0}))^{2} \! + \! \hat{\eta}_{0}^{\spadesuit,+}
(-1,1;\hat{\boldsymbol{\Omega}}) \hat{\eta}_{0}^{\spadesuit,+}(-1,1;
\boldsymbol{0}) \right), \\
\hat{\mathbb{M}}_{21}^{0,+}(\alpha_{p_{q}}) =& \, \dfrac{\hat{\mathbb{G}}_{
\hat{\boldsymbol{\Omega}}}^{0,0,+}(1,-1)}{\hat{\mathbb{G}}_{\boldsymbol{0}}^{
0,0,+}(1,-1)}, \quad \qquad \hat{\mathbb{M}}_{21}^{1,+}(\alpha_{p_{q}}) = 
\dfrac{\hat{\mathbb{G}}_{\hat{\boldsymbol{\Omega}}}^{0,0,+}(1,-1)}{\hat{
\mathbb{G}}_{\boldsymbol{0}}^{0,0,+}(1,-1)} \left(\hat{\eta}_{0}^{\spadesuit,+}
(1,-1;\hat{\boldsymbol{\Omega}}) \! - \! \hat{\eta}_{0}^{\spadesuit,+}(1,-1;
\boldsymbol{0}) \right), \\
\hat{\mathbb{M}}_{21}^{2,+}(\alpha_{p_{q}}) =& \, \dfrac{\hat{\mathbb{G}}_{
\hat{\boldsymbol{\Omega}}}^{0,0,+}(1,-1)}{\hat{\mathbb{G}}_{\boldsymbol{0}}^{
0,0,+}(1,-1)} \left(\hat{\eta}_{1}^{\spadesuit,+}(1,-1;\hat{\boldsymbol{\Omega}}) 
\! - \! \hat{\eta}_{1}^{\spadesuit,+}(1,-1;\boldsymbol{0}) \! + \! (\hat{\eta}_{0}^{
\spadesuit,+}(1,-1;\boldsymbol{0}))^{2} \! - \! \hat{\eta}_{0}^{\spadesuit,+}(1,-1;
\hat{\boldsymbol{\Omega}}) \hat{\eta}_{0}^{\spadesuit,+}(1,-1;\boldsymbol{0}) 
\right), \\
\hat{\mathbb{M}}_{22}^{0,+}(\alpha_{p_{q}}) =& \, \dfrac{\hat{\mathbb{G}}_{
\hat{\boldsymbol{\Omega}}}^{0,0,+}(-1,-1)}{\hat{\mathbb{G}}_{\boldsymbol{0}}^{
0,0,+}(-1,-1)}, \quad \qquad \hat{\mathbb{M}}_{22}^{1,+}(\alpha_{p_{q}}) = 
-\dfrac{\hat{\mathbb{G}}_{\hat{\boldsymbol{\Omega}}}^{0,0,+}(-1,-1)}{\hat{
\mathbb{G}}_{\boldsymbol{0}}^{0,0,+}(-1,-1)} \left(\hat{\eta}_{0}^{\spadesuit,+}
(-1,-1;\hat{\boldsymbol{\Omega}}) \! - \! \hat{\eta}_{0}^{\spadesuit,+}(-1,-1;
\boldsymbol{0}) \right), \\
\hat{\mathbb{M}}_{22}^{2,+}(\alpha_{p_{q}}) =& \, -\dfrac{\hat{\mathbb{G}}_{
\hat{\boldsymbol{\Omega}}}^{0,0,+}(-1,-1)}{\hat{\mathbb{G}}_{\boldsymbol{0}}^{
0,0,+}(-1,-1)} \left(\hat{\eta}_{1}^{\spadesuit,+}(-1,-1;\hat{\boldsymbol{\Omega}}) 
\! - \! \hat{\eta}_{1}^{\spadesuit,+}(-1,-1;\boldsymbol{0}) \! - \! (\hat{\eta}_{0}^{
\spadesuit,+}(-1,-1;\boldsymbol{0}))^{2} \right. \\
+&\left. \, \hat{\eta}_{0}^{\spadesuit,+}(-1,-1;\hat{\boldsymbol{\Omega}}) 
\hat{\eta}_{0}^{\spadesuit,+}(-1,-1;\boldsymbol{0}) \right),
\end{align*}
where, for $\varepsilon_{1},\varepsilon_{2} \! = \! \pm 1$,
\begin{gather*}
\hat{\mathbb{G}}_{\hat{\boldsymbol{\Omega}}}^{j_{1},j_{2},\pm}
(\varepsilon_{1},\varepsilon_{2}) \! := \! \sum_{m \in \mathbb{Z}^{N}} 
\me^{2 \pi \mi (m,\varepsilon_{1} \hat{\boldsymbol{u}}_{\pm}(\alpha_{
p_{q}})-\frac{1}{2 \pi}((n-1)K+k) \hat{\boldsymbol{\Omega}}+ 
\varepsilon_{2} \hat{\boldsymbol{d}})+\mi \pi (m,\hat{\boldsymbol{
\tau}}m)}(\hat{\lambda}_{2})^{j_{1}}(\hat{\lambda}_{3})^{j_{2}}, \quad 
j_{1} \! = \! 0,1,2, \quad j_{2} \! = \! 0,1, \\
\hat{\lambda}_{2} \! = \! 2 \pi \mi \sum_{j=1}^{N}m_{j} \hat{\mathbb{
P}}^{\triangledown}_{j}, \quad \qquad \hat{\lambda}_{3} \! = \! \mi \pi 
\sum_{j=1}^{N}m_{j} \hat{\mathbb{Q}}^{\triangledown}_{j}, \\
\hat{\mathbb{P}}^{\triangledown}_{j} \! := \! \sum_{m=1}^{N} 
\dfrac{\hat{c}_{jm}(\alpha_{p_{q}})^{N-m}}{(-1)^{\hat{\mathfrak{n}}
(\alpha_{p_{q}})} \prod_{i=1}^{N+1}(\lvert \alpha_{p_{q}} \! - \! \hat{b}_{i-1} 
\rvert \lvert \alpha_{p_{q}} \! - \! \hat{a}_{i} \rvert)^{1/2}}, \qquad \, \, 
\hat{\mathbb{Q}}^{\triangledown}_{j} \! := \! \sum_{m=1}^{N} 
\dfrac{\hat{c}_{jm}(\alpha_{p_{q}})^{N-m-1}(N \! - \! m \! - \! \hat{
\lambda}_{0} \alpha_{p_{q}})}{(-1)^{\hat{\mathfrak{n}}(\alpha_{p_{q}})} 
\prod_{i=1}^{N+1}(\lvert \alpha_{p_{q}} \! - \! \hat{b}_{i-1} \rvert 
\lvert \alpha_{p_{q}} \! - \! \hat{a}_{i} \rvert)^{1/2}}, \\
\hat{\lambda}_{0} \! := \! -\dfrac{1}{2} \sum_{j=1}^{N+1} \left(
\dfrac{1}{\hat{b}_{j-1} \! - \! \alpha_{p_{q}}} \! + \! \dfrac{1}{
\hat{a}_{j} \! - \! \alpha_{p_{q}}} \right),
\end{gather*}
with $\hat{c}_{jm}$, $j,m \! \in \! \lbrace 1,2,\dotsc,N \rbrace$, 
obtained {}from Equations~\eqref{E1} and~\eqref{E2}, and $\hat{
\mathfrak{n}}(\alpha_{p_{q}})$ defined in the corresponding item of 
Remark~\ref{remext},
\begin{gather*}
\hat{\eta}_{0}^{\spadesuit,\pm}(\varepsilon_{1},\varepsilon_{2};
\hat{\boldsymbol{\Omega}}) \! := \! \left(\hat{\mathbb{G}}_{\hat{
\boldsymbol{\Omega}}}^{0,0,\pm}(\varepsilon_{1},\varepsilon_{2}) 
\right)^{-1} \hat{\mathbb{G}}_{\hat{\boldsymbol{\Omega}}}^{1,0,\pm}
(\varepsilon_{1},\varepsilon_{2}), \\
\hat{\eta}_{1}^{\spadesuit,\pm}(\varepsilon_{1},\varepsilon_{2};
\hat{\boldsymbol{\Omega}}) \! := \! \left(\hat{\mathbb{G}}_{\hat{
\boldsymbol{\Omega}}}^{0,0,\pm}(\varepsilon_{1},\varepsilon_{2}) 
\right)^{-1} \left(\hat{\mathbb{G}}_{\hat{\boldsymbol{\Omega}}}^{0,1,
\pm}(\varepsilon_{1},\varepsilon_{2}) \! \pm \! \dfrac{\varepsilon_{1}}{2} 
\hat{\mathbb{G}}_{\hat{\boldsymbol{\Omega}}}^{2,0,\pm}(\varepsilon_{1},
\varepsilon_{2}) \right), \\
\hat{\Xi}_{0} \! := \! ((n \! - \! 1)K \! + \! k) \int_{J_{\infty}} \ln (\lvert 
\xi \! - \! \alpha_{p_{q}} \rvert) \, \md \mu_{\widetilde{V}}^{\infty}(\xi) 
\! - \! \sum_{\substack{j=1\\j \neq q}}^{\mathfrak{s}-1} \varkappa_{nk 
\tilde{k}_{j}} \ln (\lvert \alpha_{p_{j}} \! - \! \alpha_{p_{q}} \rvert),
\end{gather*}
$\tilde{\Delta}_{\infty}(q) \! := \! \lbrace \mathstrut j \! \in \! \lbrace 
1,2,\dotsc,\mathfrak{s} \! - \! 1 \rbrace \setminus \lbrace q \rbrace; 
\, \alpha_{p_{j}} \! > \! \alpha_{p_{q}} \rbrace$, and $(\mathrm{M}_{2}
(\mathbb{C}) \! \ni)$ $\hat{\mathfrak{c}}_{q}^{\mathcal{X}}(n,k,z_{o}) \! 
=_{\underset{z_{o}=1+o(1)}{\mathscr{N},n \to \infty}} \! \mathcal{O}(1)$$;$ 
and {\rm \pmb{(2)}} for $n \! \in \! \mathbb{N}$ and $k \! \in \! \lbrace 1,
2,\dotsc,K \rbrace$ such that $\alpha_{p_{\mathfrak{s}}} \! := \! \alpha_{k} 
\! \neq \! \infty$,
\begin{align}
\mathcal{X}(z)(z \! - \! \alpha_{k})^{(\varkappa_{nk}-1) \sigma_{3}} 
\underset{z \to \alpha_{k}}{=}& \, \mathrm{I} \! + \! \me^{\frac{n \tilde{\ell}}{2} 
\mathrm{ad}(\sigma_{3})} \left(\left(w_{0}^{\triangle}(\alpha_{k}) \sigma_{3} \! 
+ \! A_{0}(\alpha_{k}) \! + \! \tilde{\mathcal{R}}_{0}^{A_{0}}(\alpha_{k}) \right)
(z \! - \! \alpha_{k}) \! + \! \left(w_{1}^{\triangle}(\alpha_{k}) \sigma_{3} \! + \! 
\dfrac{1}{2}(w_{0}^{\triangle}(\alpha_{k}))^{2} \mathrm{I} \right. \right. \nonumber \\
+&\left. \left. w_{0}^{\triangle}(\alpha_{k}) \left(A_{0}(\alpha_{k}) \! + 
\! \tilde{\mathcal{R}}_{0}^{A_{0}}(\alpha_{k}) \right) \sigma_{3} \! + \! 
B_{0}(\alpha_{k}) \! + \! \tilde{\mathcal{R}}_{0}^{A_{0}}(\alpha_{k})A_{0}
(\alpha_{k}) \! + \! \tilde{\mathcal{R}}_{0}^{B_{0}}(\alpha_{k}) \right)
(z \! - \! \alpha_{k})^{2} \right. \nonumber \\
+&\left. \, \mathcal{O} \left(\tilde{\mathfrak{c}}_{\mathfrak{s}}^{\mathcal{X}}
(n,k,z_{o})(z \! - \! \alpha_{k})^{3} \right) \vphantom{M^{M^{M^{M^{M^{M}}}}}} 
\right), \label{eqlem5.5C}
\end{align}
and, for $q \! \in \! \lbrace 1,2,\dotsc,\mathfrak{s} \! - \! 2 
\rbrace$,\footnote{Note: (i) if, for $q \! \in \! \lbrace 1,2,\dotsc,
\mathfrak{s} \! - \! 2 \rbrace$ and $j \! \in \! \lbrace 1,2,\dotsc,
N \rbrace$, $\alpha_{p_{q}} \! \in \! (\tilde{a}_{j},\tilde{b}_{j})$, then 
$J_{f} \cap \mathbb{R}_{\alpha_{p_{q}}}^{>} \! = \! \cup_{i=j}^{N}
[\tilde{b}_{i},\tilde{a}_{i+1}]$, and $\int_{J_{f} \cap \mathbb{R}_{
\alpha_{p_{q}}}^{>}} \md \mu_{\widetilde{V}}^{f}(\xi) \! = \! \int_{
\tilde{b}_{j}}^{\tilde{a}_{N+1}} \md \mu_{\widetilde{V}}^{f}(\xi) \! = 
\! \tilde{\Omega}_{j}/2 \pi$; (ii) if, for $q \! \in \! \lbrace 1,2,\dotsc,
\mathfrak{s} \! - \! 2 \rbrace$, $\alpha_{p_{q}} \! \in \! (-\infty,
\tilde{b}_{0})$, then $J_{f} \cap \mathbb{R}_{\alpha_{p_{q}}}^{>} \! 
= \! J_{f}$, and $\int_{J_{f} \cap \mathbb{R}_{\alpha_{p_{q}}}^{>}} \md 
\mu_{\widetilde{V}}^{f}(\xi) \! = \! 1$; and (iii) if, for $q \! \in \! \lbrace 
1,2,\dotsc,\mathfrak{s} \! - \! 2 \rbrace$, $\alpha_{p_{q}} \! \in \! 
(\tilde{a}_{N+1},+\infty)$, then $J_{f} \cap \mathbb{R}_{\alpha_{p_{
q}}}^{>} \! = \! \varnothing$, and $\int_{J_{f} \cap \mathbb{R}_{
\alpha_{p_{q}}}^{>}} \md \mu_{\widetilde{V}}^{f}(\xi) \! = \! 0$.}
\begin{align}
\mathcal{X}(z)(z \! - \! \alpha_{p_{q}})^{\varkappa_{nk \tilde{k}_{q}} 
\sigma_{3}} \underset{z \to \alpha_{p_{q}}}{=}& \, \me^{\frac{n \tilde{\ell}}{2} 
\mathrm{ad}(\sigma_{3})} \left(\vphantom{M^{M^{M^{M^{M^{M}}}}}} 
\left(\mathrm{I} \! + \! \tilde{\mathcal{R}}_{0}^{\tilde{A}_{0}}(\alpha_{p_{q}}) 
\right) \tilde{A}_{0}(\alpha_{p_{q}}) \! + \! \left(\tilde{w}_{0}^{\triangle}
(\alpha_{p_{q}}) \left(\mathrm{I} \! + \! \tilde{\mathcal{R}}_{0}^{\tilde{A}_{0}}
(\alpha_{p_{q}}) \right) \tilde{A}_{0}(\alpha_{p_{q}}) \sigma_{3} \! + \! \left(
\mathrm{I} \! + \! \tilde{\mathcal{R}}_{0}^{\tilde{A}_{0}}(\alpha_{p_{q}}) 
\right) \right. \right. \nonumber \\
\times&\left. \left. \, \tilde{B}_{0}(\alpha_{p_{q}}) \! + \! \tilde{
\mathcal{R}}_{0}^{\tilde{B}_{0}}(\alpha_{p_{q}}) \tilde{A}_{0}(\alpha_{p_{q}}) 
\right)(z \! - \! \alpha_{p_{q}}) \! + \! \left(\left(\mathrm{I} \! + \! 
\tilde{\mathcal{R}}_{0}^{\tilde{A}_{0}}(\alpha_{p_{q}}) \right) \tilde{A}_{0}
(\alpha_{p_{q}}) \left(\tilde{w}_{1}^{\triangle}(\alpha_{p_{q}}) \sigma_{3} 
\! + \! \dfrac{1}{2}(\tilde{w}_{0}^{\triangle}(\alpha_{p_{q}}))^{2} \mathrm{I} 
\right) \right. \right. \nonumber \\
+&\left. \left. \, \tilde{w}_{0}^{\triangle}(\alpha_{p_{q}}) \left(\left(
\mathrm{I} \! + \! \tilde{\mathcal{R}}_{0}^{\tilde{A}_{0}}(\alpha_{p_{q}}) 
\right) \tilde{B}_{0}(\alpha_{p_{q}}) \! + \! \tilde{\mathcal{R}}_{0}^{
\tilde{B}_{0}}(\alpha_{p_{q}}) \tilde{A}_{0}(\alpha_{p_{q}}) \right) \sigma_{3} 
\! + \! \left(\mathrm{I} \! + \! \tilde{\mathcal{R}}_{0}^{\tilde{A}_{0}}
(\alpha_{p_{q}}) \right) \tilde{C}_{0}(\alpha_{p_{q}}) \right. \right. 
\nonumber \\
+&\left. \left. \, \tilde{\mathcal{R}}_{0}^{\tilde{B}_{0}}(\alpha_{p_{q}}) 
\tilde{B}_{0}(\alpha_{p_{q}}) \! + \! \tilde{\mathcal{R}}_{0}^{\tilde{C}_{0}}
(\alpha_{p_{q}}) \tilde{A}_{0}(\alpha_{p_{q}}) \right)(z \! - \! \alpha_{
p_{q}})^{2} \! + \! \mathcal{O} \left(\tilde{\mathfrak{c}}_{q}^{\mathcal{X}}
(n,k,z_{o})(z \! - \! \alpha_{p_{q}})^{3} \right) 
\vphantom{M^{M^{M^{M^{M^{M}}}}}} \right) \nonumber \\
\times& \, (-1)^{(\sum_{j \in \hat{\Delta}_{f}(k)} \varkappa_{nk \tilde{k}_{j}}
+ \sum_{j \in \tilde{\Delta}_{f}(q)} \varkappa_{nk \tilde{k}_{j}}+
(\varkappa_{nk}-1) \epsilon (k,q)) \sigma_{3}} \me^{\tilde{\Xi}_{0}
(\alpha_{p_{q}}) \sigma_{3}} \me^{\mi \pi ((n-1)K+k) \int_{J_{f} \cap 
\mathbb{R}_{\alpha_{p_{q}}}^{>}} \md \mu_{\widetilde{V}}^{f}(\xi) \sigma_{3}} , 
\label{eqlem5.5D}
\end{align}
where
\begin{align*}
\tilde{w}_{r_{5}}^{\triangle}(\alpha_{k}) :=& \, \dfrac{1}{1 \! + \! r_{5}} 
\left(\sum_{j=1}^{\mathfrak{s}-2} \dfrac{\varkappa_{nk \tilde{k}_{j}}}{
(\alpha_{p_{j}} \! - \! \alpha_{k})^{1+r_{5}}} \! - \! ((n \! - \! 1)K \! + \! k) 
\int_{J_{f}}(\xi \! - \! \alpha_{k})^{-(1+r_{5})} \, \md \mu_{\widetilde{V}}^{f}
(\xi) \right), \quad r_{5} \! = \! 0,1, \\
\tilde{w}_{r_{6}}^{\triangle}(\alpha_{p_{q}}) :=& \, \dfrac{1}{1 \! + \! r_{6}} 
\left(\dfrac{\varkappa_{nk} \! - \! 1}{(\alpha_{k} \! - \! \alpha_{p_{q}})^{1
+r_{6}}} \! + \! \sum_{\substack{j=1\\j \neq q}}^{\mathfrak{s}-2} \dfrac{
\varkappa_{nk \tilde{k}_{j}}}{(\alpha_{p_{j}} \! - \! \alpha_{p_{q}})^{1+r_{6}}} 
\! - \! ((n \! - \! 1)K \! + \! k) \int_{J_{f}}(\xi \! - \! \alpha_{p_{q}})^{-(1+
r_{6})} \, \md \mu_{\widetilde{V}}^{f}(\xi) \right), \quad r_{6} \! = \! 0,1, \\
\tilde{\mathcal{R}}_{0}^{A_{0}}(\alpha_{k}) \underset{\underset{z_{o}=1+
o(1)}{\mathscr{N},n \to \infty}}{=}& \, \dfrac{1}{((n \! - \! 1)K \! + \! k)} 
\sum_{j=1}^{N+1} \left(-\dfrac{(\tilde{\alpha}_{0}(\tilde{b}_{j-1}))^{-2}}{
(\tilde{b}_{j-1} \! - \! \alpha_{k})^{2}} \left(\tilde{\alpha}_{0}(\tilde{b}_{j-1}) 
\tilde{\boldsymbol{\mathrm{B}}}(\tilde{b}_{j-1}) \! - \! \tilde{\alpha}_{1}
(\tilde{b}_{j-1}) \tilde{\boldsymbol{\mathrm{A}}}(\tilde{b}_{j-1}) \right) 
\right. \\
-&\left. \, \dfrac{(\tilde{\alpha}_{0}(\tilde{a}_{j}))^{-2}}{(\tilde{a}_{j} \! - 
\! \alpha_{k})^{2}} \left(\tilde{\alpha}_{0}(\tilde{a}_{j}) \tilde{\boldsymbol{
\mathrm{B}}}(\tilde{a}_{j}) \! - \! \tilde{\alpha}_{1}(\tilde{a}_{j}) \tilde{
\boldsymbol{\mathrm{A}}}(\tilde{a}_{j}) \right) \! + \! \dfrac{2(\tilde{
\alpha}_{0}(\tilde{b}_{j-1}))^{-1}}{(\tilde{b}_{j-1} \! - \! \alpha_{k})^{3}} 
\tilde{\boldsymbol{\mathrm{A}}}(\tilde{b}_{j-1}) \! + \! \dfrac{2(\tilde{
\alpha}_{0}(\tilde{a}_{j}))^{-1}}{(\tilde{a}_{j} \! - \! \alpha_{k})^{3}} \tilde{
\boldsymbol{\mathrm{A}}}(\tilde{a}_{j}) \right) \\
+& \, \mathcal{O} \left(\dfrac{\tilde{\mathfrak{c}}_{\tilde{\mathcal{
R}}_{0}^{A_{0}}}(n,k,z_{o})}{((n \! - \! 1)K \! + \! k)^{2}} \right), \\
\tilde{\mathcal{R}}_{0}^{B_{0}}(\alpha_{k}) \underset{\underset{z_{o}=1+
o(1)}{\mathscr{N},n \to \infty}}{=}& \, \dfrac{1}{((n \! - \! 1)K \! + \! k)} 
\sum_{j=1}^{N+1} \left(-\dfrac{(\tilde{\alpha}_{0}(\tilde{b}_{j-1}))^{-2}}{
(\tilde{b}_{j-1} \! - \! \alpha_{k})^{3}} \left(\tilde{\alpha}_{0}(\tilde{b}_{j
-1}) \tilde{\boldsymbol{\mathrm{B}}}(\tilde{b}_{j-1}) \! - \! \tilde{
\alpha}_{1}(\tilde{b}_{j-1}) \tilde{\boldsymbol{\mathrm{A}}}(\tilde{b}_{j
-1}) \right) \right. \\
-&\left. \, \dfrac{(\tilde{\alpha}_{0}(\tilde{a}_{j}))^{-2}}{(\tilde{a}_{j} \! - \! 
\alpha_{k})^{3}} \left(\tilde{\alpha}_{0}(\tilde{a}_{j}) \tilde{\boldsymbol{
\mathrm{B}}}(\tilde{a}_{j}) \! - \! \tilde{\alpha}_{1}(\tilde{a}_{j}) \tilde{
\boldsymbol{\mathrm{A}}}(\tilde{a}_{j}) \right) \! + \! \dfrac{3(\tilde{
\alpha}_{0}(\tilde{b}_{j-1}))^{-1}}{(\tilde{b}_{j-1} \! - \! \alpha_{k})^{4}} 
\tilde{\boldsymbol{\mathrm{A}}}(\tilde{b}_{j-1}) \! + \! \dfrac{3(\tilde{
\alpha}_{0}(\tilde{a}_{j}))^{-1}}{(\tilde{a}_{j} \! - \! \alpha_{k})^{4}} \tilde{
\boldsymbol{\mathrm{A}}}(\tilde{a}_{j}) \right) \\
+& \, \mathcal{O} \left(\dfrac{\tilde{\mathfrak{c}}_{\tilde{\mathcal{
R}}_{0}^{B_{0}}}(n,k,z_{o})}{((n \! - \! 1)K \! + \! k)^{2}} \right), \\
\tilde{\mathcal{R}}_{0}^{\tilde{A}_{0}}(\alpha_{p_{q}}) \underset{\underset{
z_{o}=1+o(1)}{\mathscr{N},n \to \infty}}{=}& \, \dfrac{1}{((n \! - \! 1)K 
\! + \! k)} \sum_{j=1}^{N+1} \left(-\left(\dfrac{1}{\tilde{b}_{j-1} \! - \! 
\alpha_{p_{q}}} \! - \! \dfrac{1}{\tilde{b}_{j-1} \! - \! \alpha_{k}} \right)
(\tilde{\alpha}_{0}(\tilde{b}_{j-1}))^{-2} \left(\tilde{\alpha}_{0}
(\tilde{b}_{j-1}) \tilde{\boldsymbol{\mathrm{B}}}(\tilde{b}_{j-1}) \! - 
\! \tilde{\alpha}_{1}(\tilde{b}_{j-1}) \tilde{\boldsymbol{\mathrm{A}}}
(\tilde{b}_{j-1}) \right) \right. \\
-&\left. \, \left(\dfrac{1}{\tilde{a}_{j} \! - \! \alpha_{p_{q}}} \! - \! 
\dfrac{1}{\tilde{a}_{j} \! - \! \alpha_{k}} \right)(\tilde{\alpha}_{0}
(\tilde{a}_{j}))^{-2} \left(\tilde{\alpha}_{0}(\tilde{a}_{j}) \tilde{\boldsymbol{
\mathrm{B}}}(\tilde{a}_{j}) \! - \! \tilde{\alpha}_{1}(\tilde{a}_{j}) \tilde{
\boldsymbol{\mathrm{A}}}(\tilde{a}_{j}) \right) \! + \! \left(\dfrac{1}{
(\tilde{b}_{j-1} \! - \! \alpha_{p_{q}})^{2}} \! - \! \dfrac{1}{(\tilde{b}_{j-1} 
\! - \! \alpha_{k})^{2}} \right) \right. \\
\times&\left. \, (\tilde{\alpha}_{0}(\tilde{b}_{j-1}))^{-1} \tilde{
\boldsymbol{\mathrm{A}}}(\tilde{b}_{j-1}) \! + \! \left(\dfrac{1}{(\tilde{a}_{j} 
\! - \! \alpha_{p_{q}})^{2}} \! - \! \dfrac{1}{(\tilde{a}_{j} \! - \! \alpha_{
k})^{2}} \right) (\tilde{\alpha}_{0}(\tilde{a}_{j}))^{-1} \tilde{\boldsymbol{
\mathrm{A}}}(\tilde{a}_{j}) \right) \! + \! \mathcal{O} \left(\dfrac{\tilde{
\mathfrak{c}}_{\tilde{\mathcal{R}}_{0}^{\tilde{A}_{0}}}(n,k,z_{o})}{((n \! - \! 
1)K \! + \! k)^{2}} \right), \\
\tilde{\mathcal{R}}_{0}^{\tilde{B}_{0}}(\alpha_{p_{q}}) \underset{\underset{
z_{o}=1+o(1)}{\mathscr{N},n \to \infty}}{=}& \, \dfrac{1}{((n \! - \! 1)K 
\! + \! k)} \sum_{j=1}^{N+1} \left(-\dfrac{(\tilde{\alpha}_{0}(\tilde{b}_{j
-1}))^{-2}}{(\tilde{b}_{j-1} \! - \! \alpha_{p_{q}})^{2}} \left(\tilde{\alpha}_{0}
(\tilde{b}_{j-1}) \tilde{\boldsymbol{\mathrm{B}}}(\tilde{b}_{j-1}) \! - \! 
\tilde{\alpha}_{1}(\tilde{b}_{j-1}) \tilde{\boldsymbol{\mathrm{A}}}
(\tilde{b}_{j-1}) \right) \right. \\
-&\left. \, \dfrac{(\tilde{\alpha}_{0}(\tilde{a}_{j}))^{-2}}{(\tilde{a}_{j} \! - \! 
\alpha_{p_{q}})^{2}} \left(\tilde{\alpha}_{0}(\tilde{a}_{j}) \tilde{\boldsymbol{
\mathrm{B}}}(\tilde{a}_{j}) \! - \! \tilde{\alpha}_{1}(\tilde{a}_{j}) \tilde{
\boldsymbol{\mathrm{A}}}(\tilde{a}_{j}) \right) \! + \! \dfrac{2(\tilde{
\alpha}_{0}(\tilde{b}_{j-1}))^{-1}}{(\tilde{b}_{j-1} \! - \! \alpha_{p_{q}})^{3}} 
\tilde{\boldsymbol{\mathrm{A}}}(\tilde{b}_{j-1}) \! + \! \dfrac{2(\tilde{
\alpha}_{0}(\tilde{a}_{j}))^{-1}}{(\tilde{a}_{j} \! - \! \alpha_{p_{q}})^{3}} 
\tilde{\boldsymbol{\mathrm{A}}}(\tilde{a}_{j}) \right) \\
+& \, \mathcal{O} \left(\dfrac{\tilde{\mathfrak{c}}_{\tilde{\mathcal{
R}}_{0}^{\tilde{B}_{0}}}(n,k,z_{o})}{((n \! - \! 1)K \! + \! k)^{2}} \right), \\
\tilde{\mathcal{R}}_{0}^{\tilde{C}_{0}}(\alpha_{p_{q}}) \underset{\underset{
z_{o}=1+o(1)}{\mathscr{N},n \to \infty}}{=}& \, \dfrac{1}{((n \! - \! 1)K 
\! + \! k)} \sum_{j=1}^{N+1} \left(-\dfrac{(\tilde{\alpha}_{0}(\tilde{b}_{j-
1}))^{-2}}{(\tilde{b}_{j-1} \! - \! \alpha_{p_{q}})^{3}} \left(\tilde{\alpha}_{0}
(\tilde{b}_{j-1}) \tilde{\boldsymbol{\mathrm{B}}}(\tilde{b}_{j-1}) \! - \! 
\tilde{\alpha}_{1}(\tilde{b}_{j-1}) \tilde{\boldsymbol{\mathrm{A}}}
(\tilde{b}_{j-1}) \right) \right. \\
-&\left. \, \dfrac{(\tilde{\alpha}_{0}(\tilde{a}_{j}))^{-2}}{(\tilde{a}_{j} \! - \! 
\alpha_{p_{q}})^{3}} \left(\tilde{\alpha}_{0}(\tilde{a}_{j}) \tilde{\boldsymbol{
\mathrm{B}}}(\tilde{a}_{j}) \! - \! \tilde{\alpha}_{1}(\tilde{a}_{j}) \tilde{
\boldsymbol{\mathrm{A}}}(\tilde{a}_{j}) \right) \! + \! \dfrac{3(\tilde{
\alpha}_{0}(\tilde{b}_{j-1}))^{-1}}{(\tilde{b}_{j-1} \! - \! \alpha_{p_{q}})^{4}} 
\tilde{\boldsymbol{\mathrm{A}}}(\tilde{b}_{j-1}) \! + \! \dfrac{3(\tilde{
\alpha}_{0}(\tilde{a}_{j}))^{-1}}{(\tilde{a}_{j} \! - \! \alpha_{p_{q}})^{4}} 
\tilde{\boldsymbol{\mathrm{A}}}(\tilde{a}_{j}) \right) \\
+& \, \mathcal{O} \left(\dfrac{\tilde{\mathfrak{c}}_{\tilde{\mathcal{
R}}_{0}^{\tilde{C}_{0}}}(n,k,z_{o})}{((n \! - \! 1)K \! + \! k)^{2}} \right),
\end{align*}
with $\tilde{\boldsymbol{\mathrm{A}}}(\tilde{b}_{j-1})$, $\tilde{
\boldsymbol{\mathrm{A}}}(\tilde{a}_{j})$, $\tilde{\boldsymbol{
\mathrm{B}}}(\tilde{b}_{j-1})$, $\tilde{\boldsymbol{\mathrm{B}}}
(\tilde{a}_{j})$, $\tilde{\alpha}_{0}(\tilde{b}_{j-1})$, $\tilde{\alpha}_{0}
(\tilde{a}_{j})$, $\tilde{\alpha}_{1}(\tilde{b}_{j-1})$, and $\tilde{\alpha}_{1}
(\tilde{a}_{j})$, $j \! \in \! \lbrace 1,2,\dotsc,N \! + \! 1 \rbrace$, 
defined in item~{\rm \pmb{(2)}} of Proposition~\ref{propo5.1}, and 
$(\mathrm{M}_{2}(\mathbb{C}) \! \ni)$ $\tilde{\mathfrak{c}}_{
\tilde{\mathcal{R}}_{0}^{r_{7}}}(n,k,z_{o}) \! =_{\underset{z_{o}=1
+o(1)}{\mathscr{N},n \to \infty}} \! \mathcal{O}(1)$, $r_{7} \! \in 
\! \lbrace \tilde{A}_{0},\tilde{B}_{0},\tilde{C}_{0} \rbrace$,
\begin{gather*}
A_{0}(\alpha_{k}) \! = \! \widetilde{\mathbb{K}} 

\right),
\end{gather*}
where $\widetilde{\mathbb{K}} \! := \! \mathscr{E}^{-\sigma_{3}} 
\tilde{\mathbb{K}}$, with $\tilde{\mathbb{K}}$ and $\mathscr{E}$ 
defined by Equations~\eqref{eqmainfin8} and~\eqref{eqmainfin13}, 
respectively,\footnote{See, also, Equation~\eqref{eqconsfinn1} for 
the definition of $\mathfrak{k}_{0}$.} for $q^{\prime} \! \in \! \lbrace 
1,\dotsc,\mathfrak{s} \! - \! 2,\mathfrak{s} \rbrace$,
\begin{gather*}
\mathfrak{h}_{11}^{0,+}(\alpha_{p_{q^{\prime}}}) \! := \! \dfrac{1}{2} 
\left(\tilde{\gamma}(\alpha_{p_{q^{\prime}}}) \! + \! (\tilde{\gamma}
(\alpha_{p_{q^{\prime}}}))^{-1} \right), \qquad \, \, \mathfrak{h}_{11}^{1,+}
(\alpha_{p_{q^{\prime}}}) \! := \! \dfrac{\alpha^{\triangle}}{2} \left(
\tilde{\gamma}(\alpha_{p_{q^{\prime}}}) \! - \! (\tilde{\gamma}
(\alpha_{p_{q^{\prime}}}))^{-1} \right), \\
\mathfrak{h}_{11}^{2,+}(\alpha_{p_{q^{\prime}}}) \! := \! \dfrac{1}{2} 
\left(\beta^{\triangle} \left(\tilde{\gamma}(\alpha_{p_{q^{\prime}}}) 
\! - \! (\tilde{\gamma}(\alpha_{p_{q^{\prime}}}))^{-1} \right) \! + \! 
\dfrac{(\alpha^{\triangle})^{2}}{2} \left(\tilde{\gamma}(\alpha_{p_{
q^{\prime}}}) \! + \! (\tilde{\gamma}(\alpha_{p_{q^{\prime}}}))^{-1} 
\right) \right), \\
\mathfrak{h}_{12}^{0,+}(\alpha_{p_{q^{\prime}}}) \! := \! \dfrac{\mi}{2} 
\left(\tilde{\gamma}(\alpha_{p_{q^{\prime}}}) \! - \! (\tilde{\gamma}
(\alpha_{p_{q^{\prime}}}))^{-1} \right), \qquad \, \, \mathfrak{h}_{12}^{1,+}
(\alpha_{p_{q^{\prime}}}) \! := \! \dfrac{\mi \alpha^{\triangle}}{2} 
\left(\tilde{\gamma}(\alpha_{p_{q^{\prime}}}) \! + \! (\tilde{\gamma}
(\alpha_{p_{q^{\prime}}}))^{-1} \right), \\
\mathfrak{h}_{12}^{2,+}(\alpha_{p_{q^{\prime}}}) \! := \! \dfrac{\mi}{2} 
\left(\beta^{\triangle} \left(\tilde{\gamma}(\alpha_{p_{q^{\prime}}}) 
\! + \! (\tilde{\gamma}(\alpha_{p_{q^{\prime}}}))^{-1} \right) \! + \! 
\dfrac{(\alpha^{\triangle})^{2}}{2} \left(\tilde{\gamma}(\alpha_{p_{
q^{\prime}}}) \! - \! (\tilde{\gamma}(\alpha_{p_{q^{\prime}}}))^{-1} 
\right) \right),
\end{gather*}
where $\tilde{\gamma}(\alpha_{p_{q^{\prime}}})$ is defined by 
Equation~\eqref{eqssabra2},
\begin{equation*}
\alpha^{\triangle} \! = \! \dfrac{1}{4} \sum_{j=1}^{N+1} \left(\dfrac{1}{
\tilde{a}_{j} \! - \! \alpha_{p_{q^{\prime}}}} \! - \! \dfrac{1}{\tilde{b}_{j-1} 
\! - \! \alpha_{p_{q^{\prime}}}} \right), \qquad \, \, \beta^{\triangle} \! 
= \! \dfrac{1}{8} \sum_{j=1}^{N+1} \left(\dfrac{1}{(\tilde{a}_{j} \! - \! 
\alpha_{p_{q^{\prime}}})^{2}} \! - \! \dfrac{1}{(\tilde{b}_{j-1} \! - \! 
\alpha_{p_{q^{\prime}}})^{2}} \right),
\end{equation*}
\begin{align*}
\tilde{\mathbb{M}}_{11}^{0,+}(\alpha_{p_{q^{\prime}}}) =& \, \dfrac{\tilde{
\mathbb{G}}_{\tilde{\boldsymbol{\Omega}}}^{0,0,+}(1,1)}{\tilde{\mathbb{G}}_{
\boldsymbol{0}}^{0,0,+}(1,1)}, \quad \qquad \tilde{\mathbb{M}}_{11}^{1,+}
(\alpha_{p_{q^{\prime}}}) = \dfrac{\tilde{\mathbb{G}}_{\tilde{\boldsymbol{
\Omega}}}^{0,0,+}(1,1)}{\tilde{\mathbb{G}}_{\boldsymbol{0}}^{0,0,+}(1,1)} 
\left(\tilde{\eta}_{0}^{\blacklozenge,+}(1,1;\tilde{\boldsymbol{\Omega}}) 
\! - \! \tilde{\eta}_{0}^{\blacklozenge,+}(1,1;\boldsymbol{0}) \right), \\
\tilde{\mathbb{M}}_{11}^{2,+}(\alpha_{p_{q^{\prime}}}) =& \, \dfrac{
\tilde{\mathbb{G}}_{\tilde{\boldsymbol{\Omega}}}^{0,0,+}(1,1)}{\tilde{
\mathbb{G}}_{\boldsymbol{0}}^{0,0,+}(1,1)} \left(\tilde{\eta}_{1}^{\blacklozenge,
+}(1,1;\tilde{\boldsymbol{\Omega}}) \! - \! \tilde{\eta}_{1}^{\blacklozenge,
+}(1,1;\boldsymbol{0}) \! + \! (\tilde{\eta}_{0}^{\blacklozenge,+}(1,1;
\boldsymbol{0}))^{2} \! - \! \tilde{\eta}_{0}^{\blacklozenge,+}(1,1;
\tilde{\boldsymbol{\Omega}}) \tilde{\eta}_{0}^{\blacklozenge,+}
(1,1;\boldsymbol{0}) \right), \\
\tilde{\mathbb{M}}_{12}^{0,+}(\alpha_{p_{q^{\prime}}}) =& \, 
\dfrac{\tilde{\mathbb{G}}_{\tilde{\boldsymbol{\Omega}}}^{0,0,+}
(-1,1)}{\tilde{\mathbb{G}}_{\boldsymbol{0}}^{0,0,+}(-1,1)}, \quad 
\qquad \tilde{\mathbb{M}}_{12}^{1,+}(\alpha_{p_{q^{\prime}}}) = 
-\dfrac{\tilde{\mathbb{G}}_{\tilde{\boldsymbol{\Omega}}}^{0,0,+}
(-1,1)}{\tilde{\mathbb{G}}_{\boldsymbol{0}}^{0,0,+}(-1,1)} \left(
\tilde{\eta}_{0}^{\blacklozenge,+}(-1,1;\tilde{\boldsymbol{\Omega}}) 
\! - \! \tilde{\eta}_{0}^{\blacklozenge,+}(-1,1;\boldsymbol{0}) \right), \\
\tilde{\mathbb{M}}_{12}^{2,+}(\alpha_{p_{q^{\prime}}}) =& \, 
-\dfrac{\tilde{\mathbb{G}}_{\tilde{\boldsymbol{\Omega}}}^{0,0,+}
(-1,1)}{\tilde{\mathbb{G}}_{\boldsymbol{0}}^{0,0,+}(-1,1)} \left(
\tilde{\eta}_{1}^{\blacklozenge,+}(-1,1;\tilde{\boldsymbol{\Omega}}) 
\! - \! \tilde{\eta}_{1}^{\blacklozenge,+}(-1,1;\boldsymbol{0}) \! - \! 
(\tilde{\eta}_{0}^{\blacklozenge,+}(-1,1;\boldsymbol{0}))^{2} \! + \! 
\tilde{\eta}_{0}^{\blacklozenge,+}(-1,1;\tilde{\boldsymbol{\Omega}}) 
\tilde{\eta}_{0}^{\blacklozenge,+}(-1,1;\boldsymbol{0}) \right), \\
\tilde{\mathbb{M}}_{21}^{0,+}(\alpha_{p_{q^{\prime}}}) =& \, 
\dfrac{\tilde{\mathbb{G}}_{\tilde{\boldsymbol{\Omega}}}^{0,0,+}
(1,-1)}{\tilde{\mathbb{G}}_{\boldsymbol{0}}^{0,0,+}(1,-1)}, \quad 
\qquad \tilde{\mathbb{M}}_{21}^{1,+}(\alpha_{p_{q^{\prime}}}) = 
\dfrac{\tilde{\mathbb{G}}_{\tilde{\boldsymbol{\Omega}}}^{0,0,+}
(1,-1)}{\tilde{\mathbb{G}}_{\boldsymbol{0}}^{0,0,+}(1,-1)} \left(
\tilde{\eta}_{0}^{\blacklozenge,+}(1,-1;\tilde{\boldsymbol{\Omega}}) 
\! - \! \tilde{\eta}_{0}^{\blacklozenge,+}(1,-1;\boldsymbol{0}) \right), \\
\tilde{\mathbb{M}}_{21}^{2,+}(\alpha_{p_{q^{\prime}}}) =& \, \dfrac{
\tilde{\mathbb{G}}_{\tilde{\boldsymbol{\Omega}}}^{0,0,+}(1,-1)}{
\tilde{\mathbb{G}}_{\boldsymbol{0}}^{0,0,+}(1,-1)} \left(\tilde{\eta}_{1}^{
\blacklozenge,+}(1,-1;\tilde{\boldsymbol{\Omega}}) \! - \! \tilde{\eta}_{1}^{
\blacklozenge,+}(1,-1;\boldsymbol{0}) \! + \! (\tilde{\eta}_{0}^{\blacklozenge,+}
(1,-1;\boldsymbol{0}))^{2} \! - \! \tilde{\eta}_{0}^{\blacklozenge,+}(1,-1;
\tilde{\boldsymbol{\Omega}}) \tilde{\eta}_{0}^{\blacklozenge,+}(1,-1;
\boldsymbol{0}) \right), \\
\tilde{\mathbb{M}}_{22}^{0,+}(\alpha_{p_{q^{\prime}}}) =& \, \dfrac{
\tilde{\mathbb{G}}_{\tilde{\boldsymbol{\Omega}}}^{0,0,+}(-1,-1)}{
\tilde{\mathbb{G}}_{\boldsymbol{0}}^{0,0,+}(-1,-1)}, \quad \qquad 
\tilde{\mathbb{M}}_{22}^{1,+}(\alpha_{p_{q^{\prime}}}) = -\dfrac{
\tilde{\mathbb{G}}_{\tilde{\boldsymbol{\Omega}}}^{0,0,+}(-1,-1)}{
\tilde{\mathbb{G}}_{\boldsymbol{0}}^{0,0,+}(-1,-1)} \left(\tilde{
\eta}_{0}^{\blacklozenge,+}(-1,-1;\tilde{\boldsymbol{\Omega}}) \! - 
\! \tilde{\eta}_{0}^{\blacklozenge,+}(-1,-1;\boldsymbol{0}) \right), \\
\tilde{\mathbb{M}}_{22}^{2,+}(\alpha_{p_{q^{\prime}}}) =& \, -\dfrac{
\tilde{\mathbb{G}}_{\tilde{\boldsymbol{\Omega}}}^{0,0,+}(-1,-1)}{
\tilde{\mathbb{G}}_{\boldsymbol{0}}^{0,0,+}(-1,-1)} \left(\tilde{
\eta}_{1}^{\blacklozenge,+}(-1,-1;\tilde{\boldsymbol{\Omega}}) \! - \! 
\tilde{\eta}_{1}^{\blacklozenge,+}(-1,-1;\boldsymbol{0}) \! - \! (\tilde{
\eta}_{0}^{\blacklozenge,+}(-1,-1;\boldsymbol{0}))^{2} \right. \\
+&\left. \, \tilde{\eta}_{0}^{\blacklozenge,+}(-1,-1;\tilde{
\boldsymbol{\Omega}}) \tilde{\eta}_{0}^{\blacklozenge,+}
(-1,-1;\boldsymbol{0}) \right),
\end{align*}
with, for $\varepsilon_{1},\varepsilon_{2} \! = \! \pm 1$,
\begin{gather*}
\tilde{\mathbb{G}}_{\tilde{\boldsymbol{\Omega}}}^{j_{1},j_{2},\pm}
(\varepsilon_{1},\varepsilon_{2}) \! := \! \sum_{m \in \mathbb{Z}^{N}} 
\me^{2 \pi \mi (m,\varepsilon_{1} \tilde{\boldsymbol{u}}_{\pm}
(\alpha_{p_{q^{\prime}}})-\frac{1}{2 \pi}((n-1)K+k) \tilde{\boldsymbol{\Omega}}
+ \varepsilon_{2} \tilde{\boldsymbol{d}})+\mi \pi (m,\tilde{\boldsymbol{
\tau}}m)}(\lambda_{2})^{j_{1}}(\lambda_{3})^{j_{2}}, \quad j_{1} \! = \! 
0,1,2, \quad j_{2} \! = \! 0,1, \\
\lambda_{2} \! = \! 2 \pi \mi \sum_{j=1}^{N}m_{j} \tilde{\mathbb{P}}^{
\triangledown}_{j}, \quad \qquad \lambda_{3} \! = \! \mi \pi \sum_{j=1}^{N}
m_{j} \tilde{\mathbb{Q}}^{\triangledown}_{j}, \\
\tilde{\mathbb{P}}^{\triangledown}_{j} \! := \! \sum_{m=1}^{N} \dfrac{
\tilde{c}_{jm}(\alpha_{p_{q^{\prime}}})^{N-m}}{(-1)^{\tilde{\mathfrak{n}}
(\alpha_{p_{q^{\prime}}})} \prod_{i=1}^{N+1}(\lvert \alpha_{p_{q^{\prime}}} 
\! - \! \tilde{b}_{i-1} \rvert \lvert \alpha_{p_{q^{\prime}}} \! - \! \tilde{a}_{i} 
\rvert)^{1/2}}, \qquad \, \tilde{\mathbb{Q}}^{\triangledown}_{j} \! := \! 
\sum_{m=1}^{N} \dfrac{\tilde{c}_{jm}(\alpha_{p_{q^{\prime}}})^{N-m-1}
(N \! - \! m \! - \! \lambda_{0} \alpha_{p_{q^{\prime}}})}{(-1)^{\tilde{
\mathfrak{n}}(\alpha_{p_{q^{\prime}}})} \prod_{i=1}^{N+1}(\lvert 
\alpha_{p_{q^{\prime}}} \! - \! \tilde{b}_{i-1} \rvert \lvert 
\alpha_{p_{q^{\prime}}} \! - \! \tilde{a}_{i} \rvert)^{1/2}}, \\
\lambda_{0} \! := \! -\dfrac{1}{2} \sum_{j=1}^{N+1} \left(\dfrac{1}{
\tilde{b}_{j-1} \! - \! \alpha_{p_{q^{\prime}}}} \! + \! 
\dfrac{1}{\tilde{a}_{j} \! - \! \alpha_{p_{q^{\prime}}}} \right),
\end{gather*}
where $\tilde{c}_{jm}$, $j,m \! \in \! \lbrace 1,2,\dotsc,N \rbrace$, are 
obtained {}from Equations~\eqref{O1} and~\eqref{O2}, and $\tilde{
\mathfrak{n}}(\alpha_{p_{q^{\prime}}})$ is defined in the corresponding 
item of Remark~\ref{remext},
\begin{gather*}
\tilde{\eta}_{0}^{\blacklozenge,\pm}(\varepsilon_{1},\varepsilon_{2};
\tilde{\boldsymbol{\Omega}}) \! := \! \left(\tilde{\mathbb{G}}_{\tilde{
\boldsymbol{\Omega}}}^{0,0,\pm}(\varepsilon_{1},\varepsilon_{2}) 
\right)^{-1} \tilde{\mathbb{G}}_{\tilde{\boldsymbol{\Omega}}}^{1,0,
\pm}(\varepsilon_{1},\varepsilon_{2}), \\
\tilde{\eta}_{1}^{\blacklozenge,\pm}(\varepsilon_{1},\varepsilon_{2};
\tilde{\boldsymbol{\Omega}}) \! := \! \left(\tilde{\mathbb{G}}_{\tilde{
\boldsymbol{\Omega}}}^{0,0,\pm}(\varepsilon_{1},\varepsilon_{2}) 
\right)^{-1} \left(\tilde{\mathbb{G}}_{\tilde{\boldsymbol{\Omega}}}^{0,
1,\pm}(\varepsilon_{1},\varepsilon_{2}) \! \pm \! \dfrac{\varepsilon_{
1}}{2} \tilde{\mathbb{G}}_{\tilde{\boldsymbol{\Omega}}}^{2,0,\pm}
(\varepsilon_{1},\varepsilon_{2}) \right),
\end{gather*}
for $q \! \in \! \lbrace 1,2,\dotsc,\mathfrak{s} \! - \! 2 \rbrace$, 
$\hat{\Delta}_{f}(k) \! := \! \lbrace \mathstrut j \! \in \! \lbrace 1,2,
\dotsc,\mathfrak{s} \! - \! 2 \rbrace; \, \alpha_{p_{j}} \! > \! \alpha_{k} 
\rbrace$, $\hat{\Delta}_{f}(q) \! := \! \lbrace \mathstrut j \! \in \! \lbrace 
1,2,\dotsc,\mathfrak{s} \! - \! 2 \rbrace \setminus \lbrace q \rbrace; \, 
\alpha_{p_{j}} \! > \! \alpha_{p_{q}} \rbrace$, $\epsilon (k,q) \! = \! 
\left\{
\begin{smallmatrix}
1, \, \, \, \alpha_{p_{q}}< \alpha_{k}, \\
0, \, \, \, \alpha_{p_{q}}> \alpha_{k},
\end{smallmatrix}
\right.$
\begin{equation*}
\tilde{\Xi}_{0}(\alpha_{p_{q}}) \! := \! ((n \! - \! 1)K \! + \! k) \int_{J_{f}} 
\ln \left(\left\lvert \dfrac{\xi \! - \! \alpha_{p_{q}}}{\xi \! - \! \alpha_{k}} 
\right\rvert \right) \, \md \mu_{\widetilde{V}}^{f}(\xi) \! + \! \sum_{
\substack{j=1\\j \neq q}}^{\mathfrak{s}-2} \varkappa_{nk \tilde{k}_{j}} 
\ln \left(\left\lvert \dfrac{\alpha_{p_{j}} \! - \! \alpha_{k}}{\alpha_{p_{j}} 
\! - \! \alpha_{p_{q}}} \right\rvert \right) \! + \! (\varkappa_{nk \tilde{
k}_{q}} \! - \! (\varkappa_{nk} \! - \! 1)) \ln \lvert \alpha_{p_{q}} \! - \! 
\alpha_{k} \rvert,
\end{equation*}
and $(\mathrm{M}_{2}(\mathbb{C}) \! \ni)$ $\tilde{\mathfrak{c}}_{
q^{\prime \prime}}^{\mathcal{X}}(n,k,z_{o}) \! =_{\underset{z_{o}=1+
o(1)}{\mathscr{N},n \to \infty}} \! \mathcal{O}(1)$, $q^{\prime \prime} 
\! \in \! \lbrace 1,\dotsc,\mathfrak{s} \! - \! 2,\mathfrak{s} \rbrace$, 
and
\begin{align}
\mathcal{X}(z)z^{-(\varkappa_{nk \tilde{k}_{\mathfrak{s}-1}}^{\infty}
+1) \sigma_{3}} \underset{z \to \alpha_{p_{\mathfrak{s}-1}} = 
\infty}{=}& \, \me^{\frac{n \tilde{\ell}}{2} \mathrm{ad}(\sigma_{3})} 
\left(\left(\mathrm{I} \! - \! \tilde{\mathcal{R}}_{0}^{\tilde{A}_{0}^{
\sharp}}(\alpha_{p_{\mathfrak{s}-1}}) \right) \tilde{A}_{0}^{\sharp}
(\alpha_{p_{\mathfrak{s}-1}}) \! + \! \dfrac{1}{z} \left(\tilde{w}_{0}^{
\sharp}(\alpha_{p_{\mathfrak{s}-1}}) \left(\mathrm{I} \! - \! \tilde{
\mathcal{R}}_{0}^{\tilde{A}_{0}^{\sharp}}(\alpha_{p_{\mathfrak{s}-1}}) 
\right) \tilde{A}_{0}^{\sharp}(\alpha_{p_{\mathfrak{s}-1}}) \sigma_{3} 
\! + \! \left(\mathrm{I} \! - \! \tilde{\mathcal{R}}_{0}^{\tilde{A}_{0}^{
\sharp}}(\alpha_{p_{\mathfrak{s}-1}}) \right) \right. \right. 
\nonumber \\
\times&\left. \left. \, \tilde{B}_{0}^{\sharp}(\alpha_{p_{\mathfrak{s}-1}}) 
\! + \! \tilde{\mathcal{R}}_{0}^{\tilde{B}_{0}^{\sharp}}(\alpha_{
p_{\mathfrak{s}-1}}) \tilde{A}_{0}^{\sharp}(\alpha_{p_{\mathfrak{s}-1}}) 
\right) \! + \! \dfrac{1}{z^{2}} \left(\left(\mathrm{I} \!  - \! \tilde{
\mathcal{R}}_{0}^{\tilde{A}_{0}^{\sharp}}(\alpha_{p_{\mathfrak{s}-1}}) 
\right) \tilde{A}_{0}^{\sharp}(\alpha_{p_{\mathfrak{s}-1}}) \left(
\tilde{w}_{1}^{\sharp}(\alpha_{p_{\mathfrak{s}-1}}) \sigma_{3} \! + \! 
\dfrac{1}{2}(\tilde{w}_{0}^{\sharp}(\alpha_{p_{\mathfrak{s}-1}}))^{2} 
\mathrm{I} \right) \right. \right. \nonumber \\
+&\left. \left. \, \tilde{w}_{0}^{\sharp}(\alpha_{p_{\mathfrak{s}-1}}) \left(
\left(\mathrm{I} \! - \! \tilde{\mathcal{R}}_{0}^{\tilde{A}_{0}^{\sharp}}
(\alpha_{p_{\mathfrak{s}-1}}) \right) \tilde{B}_{0}^{\sharp}(\alpha_{
p_{\mathfrak{s}-1}}) \! + \! \tilde{\mathcal{R}}_{0}^{\tilde{B}_{0}^{\sharp}}
(\alpha_{p_{\mathfrak{s}-1}}) \tilde{A}_{0}^{\sharp}(\alpha_{p_{\mathfrak{
s}-1}}) \right) \sigma_{3} \! + \! \left(\mathrm{I} \! - \! \tilde{\mathcal{
R}}_{0}^{\tilde{A}_{0}^{\sharp}}(\alpha_{p_{\mathfrak{s}-1}}) \right) 
\tilde{C}_{0}^{\sharp}(\alpha_{p_{\mathfrak{s}-1}}) \right. \right. 
\nonumber \\
+&\left. \left. \, \tilde{\mathcal{R}}_{0}^{\tilde{B}_{0}^{\sharp}}
(\alpha_{p_{\mathfrak{s}-1}}) \tilde{B}_{0}^{\sharp}(\alpha_{
p_{\mathfrak{s}-1}}) \! + \! \tilde{\mathcal{R}}_{0}^{\tilde{C}_{0}^{\sharp}}
(\alpha_{p_{\mathfrak{s}-1}}) \tilde{A}_{0}^{\sharp}(\alpha_{p_{\mathfrak{s}
-1}}) \right) \! + \! \mathcal{O} \left(\tilde{\mathfrak{c}}_{\mathfrak{s}-1}^{
\mathcal{X}}(n,k,z_{o})z^{-3} \right) \vphantom{M^{M^{M^{M^{M^{M}}}}}} 
\right) \nonumber \\
\times& \, (-1)^{\sum_{j \in \hat{\Delta}_{f}(k)} \varkappa_{nk \tilde{k}_{j}} 
\sigma_{3}} \me^{\tilde{\Xi}_{0}^{\sharp}(\alpha_{p_{\mathfrak{s}-1}}) 
\sigma_{3}}, \label{eqlem5.5E}
\end{align}
where
\begin{align*}
\tilde{w}_{r_{8}}^{\sharp}(\alpha_{p_{\mathfrak{s}-1}}) :=& \, \dfrac{1}{1 
\! + \! r_{8}} \left((\varkappa_{nk} \! - \! 1)(\alpha_{k})^{1+r_{8}} \! + \! 
\sum_{j=1}^{\mathfrak{s}-2} \varkappa_{nk \tilde{k}_{j}}(\alpha_{p_{j}})^{1+
r_{8}} \! - \! ((n \! - \! 1)K \! + \! k) \int_{J_{f}} \xi^{1+r_{8}} \, \md \mu_{
\widetilde{V}}^{f}(\xi) \right), \quad r_{8} \! = \! 0,1, \\
\tilde{\mathcal{R}}_{0}^{\tilde{A}_{0}^{\sharp}}(\alpha_{p_{\mathfrak{s}-1}}) 
\underset{\underset{z_{o}=1+o(1)}{\mathscr{N},n \to \infty}}{=}& \, 
\dfrac{1}{((n \! - \! 1)K \! + \! k)} \sum_{j=1}^{N+1} \left(-\dfrac{(\tilde{
\alpha}_{0}(\tilde{b}_{j-1}))^{-2}}{\tilde{b}_{j-1} \! - \! \alpha_{k}} \left(
\tilde{\alpha}_{0}(\tilde{b}_{j-1}) \tilde{\boldsymbol{\mathrm{B}}}
(\tilde{b}_{j-1}) \! - \! \tilde{\alpha}_{1}(\tilde{b}_{j-1}) 
\tilde{\boldsymbol{\mathrm{A}}}(\tilde{b}_{j-1}) \right) \right. \\
-&\left. \, \dfrac{(\tilde{\alpha}_{0}(\tilde{a}_{j}))^{-2}}{\tilde{a}_{j} \! - \! 
\alpha_{k}} \left(\tilde{\alpha}_{0}(\tilde{a}_{j}) \tilde{\boldsymbol{
\mathrm{B}}}(\tilde{a}_{j}) \! - \! \tilde{\alpha}_{1}(\tilde{a}_{j}) \tilde{
\boldsymbol{\mathrm{A}}}(\tilde{a}_{j}) \right) \! + \! \dfrac{(\tilde{
\alpha}_{0}(\tilde{b}_{j-1}))^{-1}}{(\tilde{b}_{j-1} \! - \! \alpha_{k})^{2}} 
\tilde{\boldsymbol{\mathrm{A}}}(\tilde{b}_{j-1}) \! + \! \dfrac{(\tilde{
\alpha}_{0}(\tilde{a}_{j}))^{-1}}{(\tilde{a}_{j} \! - \! \alpha_{k})^{2}} 
\tilde{\boldsymbol{\mathrm{A}}}(\tilde{a}_{j}) \right) \\
+& \, \mathcal{O} \left(\dfrac{\tilde{\mathfrak{c}}_{\tilde{\mathcal{
R}}_{0}^{\tilde{A}_{0}^{\sharp}}}(n,k,z_{o})}{((n \! - \! 1)K \! + \! k)^{2}} 
\right), \\
\tilde{\mathcal{R}}_{0}^{\tilde{B}_{0}^{\sharp}}(\alpha_{p_{\mathfrak{s}-1}}) 
\underset{\underset{z_{o}=1+o(1)}{\mathscr{N},n \to \infty}}{=}& \, 
\dfrac{1}{((n \! - \! 1)K \! + \! k)} \sum_{j=1}^{N+1} \left(\dfrac{1}{
(\tilde{\alpha}_{0}(\tilde{b}_{j-1}))^{2}} \left(\tilde{\alpha}_{0}(\tilde{b}_{j
-1}) \tilde{\boldsymbol{\mathrm{B}}}(\tilde{b}_{j-1}) \! - \! \tilde{
\alpha}_{1}(\tilde{b}_{j-1}) \tilde{\boldsymbol{\mathrm{A}}}(\tilde{b}_{j
-1}) \right) \right. \\
+&\left. \, \dfrac{1}{(\tilde{\alpha}_{0}(\tilde{a}_{j}))^{2}} \left(\tilde{
\alpha}_{0}(\tilde{a}_{j}) \tilde{\boldsymbol{\mathrm{B}}}(\tilde{a}_{j}) 
\! - \! \tilde{\alpha}_{1}(\tilde{a}_{j}) \tilde{\boldsymbol{\mathrm{A}}}
(\tilde{a}_{j}) \right) \right) \! + \! \mathcal{O} \left(\dfrac{\tilde{
\mathfrak{c}}_{\tilde{\mathcal{R}}_{0}^{\tilde{B}_{0}^{\sharp}}}(n,k,
z_{o})}{((n \! - \! 1)K \! + \! k)^{2}} \right), \\
\tilde{\mathcal{R}}_{0}^{\tilde{C}_{0}^{\sharp}}(\alpha_{p_{\mathfrak{s}-
1}}) \underset{\underset{z_{o}=1+o(1)}{\mathscr{N},n \to \infty}}{=}& \, 
\dfrac{1}{((n \! - \! 1)K \! + \! k)} \sum_{j=1}^{N+1} \left(\dfrac{
\tilde{b}_{j-1}}{(\tilde{\alpha}_{0}(\tilde{b}_{j-1}))^{2}} \left(\tilde{
\alpha}_{0}(\tilde{b}_{j-1}) \tilde{\boldsymbol{\mathrm{B}}}(\tilde{b}_{j-
1}) \! - \! \tilde{\alpha}_{1}(\tilde{b}_{j-1}) \tilde{\boldsymbol{\mathrm{
A}}}(\tilde{b}_{j-1}) \right) \right. \\
+&\left. \, \dfrac{\tilde{a}_{j}}{(\tilde{\alpha}_{0}(\tilde{a}_{j}))^{2}} \left(
\tilde{\alpha}_{0}(\tilde{a}_{j}) \tilde{\boldsymbol{\mathrm{B}}}(\tilde{a}_{j}) 
\! - \! \tilde{\alpha}_{1}(\tilde{a}_{j}) \tilde{\boldsymbol{\mathrm{A}}}
(\tilde{a}_{j}) \right) \! + \! \dfrac{1}{\tilde{\alpha}_{0}(\tilde{b}_{j-1})} 
\tilde{\boldsymbol{\mathrm{A}}}(\tilde{b}_{j-1}) \! + \! \dfrac{1}{\tilde{
\alpha}_{0}(\tilde{a}_{j})} \tilde{\boldsymbol{\mathrm{A}}}(\tilde{a}_{j}) 
\right) \\
+& \, \mathcal{O} \left(\dfrac{\tilde{\mathfrak{c}}_{\tilde{\mathcal{
R}}_{0}^{\tilde{C}_{0}^{\sharp}}}(n,k,z_{o})}{((n \! - \! 1)K \! + \! k)^{2}} 
\right), 
\end{align*}
with $(\mathrm{M}_{2}(\mathbb{C}) \! \ni)$ $\tilde{\mathfrak{c}}_{
\tilde{\mathcal{R}}_{0}^{r_{9}}}(n,k,z_{o}) \! =_{\underset{z_{o}=1+
o(1)}{\mathscr{N},n \to \infty}} \! \mathcal{O}(1)$, $r_{9} \! \in \! \lbrace 
\tilde{A}_{0}^{\sharp},\tilde{B}_{0}^{\sharp},\tilde{C}_{0}^{\sharp} \rbrace$,
\begin{gather*}
\tilde{A}_{0}^{\sharp}(\alpha_{p_{\mathfrak{s}-1}}) \! = \! 
\widetilde{\mathbb{K}} 
\begin{pmatrix}
\tilde{\mathbb{M}}_{11}^{0,+}(\alpha_{p_{\mathfrak{s}-1}}) & 0 \\
0 & \tilde{\mathbb{M}}_{22}^{0,+}(\alpha_{p_{\mathfrak{s}-1}}) 
\end{pmatrix}, \qquad \tilde{B}_{0}^{\sharp}(\alpha_{p_{\mathfrak{s}-1}}) 
\! = \! \widetilde{\mathbb{K}} 
\begin{pmatrix}
\tilde{\mathbb{M}}_{11}^{1,+}(\alpha_{p_{\mathfrak{s}-1}}) & \mi 
\tilde{\alpha}_{0}^{\triangle} \tilde{\mathbb{M}}_{12}^{0,+}
(\alpha_{p_{\mathfrak{s}-1}}) \\
-\mi \tilde{\alpha}_{0}^{\triangle} \tilde{\mathbb{M}}_{21}^{0,+}
(\alpha_{p_{\mathfrak{s}-1}}) & \tilde{\mathbb{M}}_{22}^{1,+}
(\alpha_{p_{\mathfrak{s}-1}}) 
\end{pmatrix}, \\
\tilde{C}_{0}^{\sharp}(\alpha_{p_{\mathfrak{s}-1}}) \! = \! 
\widetilde{\mathbb{K}} 
\begin{pmatrix}
\tilde{\mathbb{M}}_{11}^{2,+}(\alpha_{p_{\mathfrak{s}-1}}) \! + \! 
\frac{1}{2}(\tilde{\alpha}_{0}^{\triangle})^{2} \tilde{\mathbb{M}}_{11}^{0,
+}(\alpha_{p_{\mathfrak{s}-1}}) & \mi \tilde{\alpha}_{0}^{\triangle} \tilde{
\mathbb{M}}_{12}^{1,+}(\alpha_{p_{\mathfrak{s}-1}}) \! + \! \mi \tilde{
\beta}_{0}^{\triangle} \tilde{\mathbb{M}}_{12}^{0,+}(\alpha_{p_{\mathfrak{s}
-1}}) \\
-\mi \tilde{\alpha}_{0}^{\triangle} \tilde{\mathbb{M}}_{21}^{1,+}(\alpha_{
p_{\mathfrak{s}-1}}) \! - \! \mi \tilde{\beta}_{0}^{\triangle} \tilde{
\mathbb{M}}_{21}^{0,+}(\alpha_{p_{\mathfrak{s}-1}}) & \tilde{\mathbb{
M}}_{22}^{2,+}(\alpha_{p_{\mathfrak{s}-1}}) \! + \! \frac{1}{2}(\tilde{
\alpha}_{0}^{\triangle})^{2} \tilde{\mathbb{M}}_{22}^{0,+}(\alpha_{
p_{\mathfrak{s}-1}})
\end{pmatrix}, \\
\end{gather*}
where
\begin{equation*}
\tilde{\alpha}_{0}^{\triangle} \! := \! \dfrac{1}{4} \sum_{j=1}^{N+1}
(\tilde{a}_{j} \! - \! \tilde{b}_{j-1}), \quad \qquad \tilde{\beta}_{0}^{
\triangle} \! := \! \dfrac{1}{8} \sum_{j=1}^{N+1}(\tilde{a}_{j}^{2} 
\! - \! \tilde{b}_{j-1}^{2}),
\end{equation*}
\begin{align*}
\tilde{\mathbb{M}}_{11}^{0,+}(\alpha_{p_{\mathfrak{s}-1}}) =& \, 
\dfrac{\tilde{\mathbb{F}}_{\tilde{\boldsymbol{\Omega}}}^{0,0,+}
(1,1)}{\tilde{\mathbb{F}}_{\boldsymbol{0}}^{0,0,+}(1,1)}, \quad 
\qquad \tilde{\mathbb{M}}_{11}^{1,+}(\alpha_{p_{\mathfrak{s}-1}}) = 
\dfrac{\tilde{\mathbb{F}}_{\tilde{\boldsymbol{\Omega}}}^{0,0,+}
(1,1)}{\tilde{\mathbb{F}}_{\boldsymbol{0}}^{0,0,+}(1,1)} \left(
\tilde{\eta}_{0}^{\spadesuit,+}(1,1;\tilde{\boldsymbol{\Omega}}) \! 
- \! \tilde{\eta}_{0}^{\spadesuit,+}(1,1;\boldsymbol{0}) \right), \\
\tilde{\mathbb{M}}_{11}^{2,+}(\alpha_{p_{\mathfrak{s}-1}}) =& \, 
\dfrac{\tilde{\mathbb{F}}_{\tilde{\boldsymbol{\Omega}}}^{0,0,+}
(1,1)}{\tilde{\mathbb{F}}_{\boldsymbol{0}}^{0,0,+}(1,1)} \left(\tilde{
\eta}_{1}^{\spadesuit,+}(1,1;\tilde{\boldsymbol{\Omega}}) \! - \! 
\tilde{\eta}_{1}^{\spadesuit,+}(1,1;\boldsymbol{0}) \! + \! (\tilde{
\eta}_{0}^{\spadesuit,+}(1,1;\boldsymbol{0}))^{2} \! - \! \tilde{\eta}_{
0}^{\spadesuit,+}(1,1;\tilde{\boldsymbol{\Omega}}) \tilde{\eta}_{0}^{
\spadesuit,+}(1,1;\boldsymbol{0}) \right), \\
\tilde{\mathbb{M}}_{12}^{0,+}(\alpha_{p_{\mathfrak{s}-1}}) =& \, 
\dfrac{\tilde{\mathbb{F}}_{\tilde{\boldsymbol{\Omega}}}^{0,0,+}
(-1,1)}{\tilde{\mathbb{F}}_{\boldsymbol{0}}^{0,0,+}(-1,1)}, \quad 
\qquad \tilde{\mathbb{M}}_{12}^{1,+}(\alpha_{p_{\mathfrak{s}-1}}) = 
-\dfrac{\tilde{\mathbb{F}}_{\tilde{\boldsymbol{\Omega}}}^{0,0,+}
(-1,1)}{\tilde{\mathbb{F}}_{\boldsymbol{0}}^{0,0,+}(-1,1)} \left(
\tilde{\eta}_{0}^{\spadesuit,+}(-1,1;\tilde{\boldsymbol{\Omega}}) \! 
- \! \tilde{\eta}_{0}^{\spadesuit,+}(-1,1;\boldsymbol{0}) \right), \\
\tilde{\mathbb{M}}_{12}^{2,+}(\alpha_{p_{\mathfrak{s}-1}}) =& \, 
-\dfrac{\tilde{\mathbb{F}}_{\tilde{\boldsymbol{\Omega}}}^{0,0,+}
(-1,1)}{\tilde{\mathbb{F}}_{\boldsymbol{0}}^{0,0,+}(-1,1)} \left(\tilde{
\eta}_{1}^{\spadesuit,+}(-1,1;\tilde{\boldsymbol{\Omega}}) \! - \! 
\tilde{\eta}_{1}^{\spadesuit,+}(-1,1;\boldsymbol{0}) \! - \! (\tilde{
\eta}_{0}^{\spadesuit,+}(-1,1;\boldsymbol{0}))^{2} \! + \! \tilde{
\eta}_{0}^{\spadesuit,+}(-1,1;\tilde{\boldsymbol{\Omega}}) 
\tilde{\eta}_{0}^{\spadesuit,+}(-1,1;\boldsymbol{0}) \right), \\
\tilde{\mathbb{M}}_{21}^{0,+}(\alpha_{p_{\mathfrak{s}-1}}) =& \, 
\dfrac{\tilde{\mathbb{F}}_{\tilde{\boldsymbol{\Omega}}}^{0,0,+}
(1,-1)}{\tilde{\mathbb{F}}_{\boldsymbol{0}}^{0,0,+}(1,-1)}, \quad 
\qquad \tilde{\mathbb{M}}_{21}^{1,+}(\alpha_{p_{\mathfrak{s}-1}}) = 
\dfrac{\tilde{\mathbb{F}}_{\tilde{\boldsymbol{\Omega}}}^{0,0,+}
(1,-1)}{\tilde{\mathbb{F}}_{\boldsymbol{0}}^{0,0,+}(1,-1)} \left(
\tilde{\eta}_{0}^{\spadesuit,+}(1,-1;\tilde{\boldsymbol{\Omega}}) 
\! - \! \tilde{\eta}_{0}^{\spadesuit,+}(1,-1;\boldsymbol{0}) \right), \\
\tilde{\mathbb{M}}_{21}^{2,+}(\alpha_{p_{\mathfrak{s}-1}}) =& \, 
\dfrac{\tilde{\mathbb{F}}_{\tilde{\boldsymbol{\Omega}}}^{0,0,+}
(1,-1)}{\tilde{\mathbb{F}}_{\boldsymbol{0}}^{0,0,+}(1,-1)} \left(\tilde{
\eta}_{1}^{\spadesuit,+}(1,-1;\tilde{\boldsymbol{\Omega}}) \! - \! 
\tilde{\eta}_{1}^{\spadesuit,+}(1,-1;\boldsymbol{0}) \! + \! (\tilde{
\eta}_{0}^{\spadesuit,+}(1,-1;\boldsymbol{0}))^{2} \! - \! \tilde{\eta}_{0}^{
\spadesuit,+}(1,-1;\tilde{\boldsymbol{\Omega}}) \tilde{\eta}_{0}^{
\spadesuit,+}(1,-1;\boldsymbol{0}) \right), \\
\tilde{\mathbb{M}}_{22}^{0,+}(\alpha_{p_{\mathfrak{s}-1}}) =& \, 
\dfrac{\tilde{\mathbb{F}}_{\tilde{\boldsymbol{\Omega}}}^{0,0,+}
(-1,-1)}{\tilde{\mathbb{F}}_{\boldsymbol{0}}^{0,0,+}(-1,-1)}, \quad 
\qquad \tilde{\mathbb{M}}_{22}^{1,+}(\alpha_{p_{\mathfrak{s}-1}}) = 
-\dfrac{\tilde{\mathbb{F}}_{\tilde{\boldsymbol{\Omega}}}^{0,0,+}
(-1,-1)}{\tilde{\mathbb{F}}_{\boldsymbol{0}}^{0,0,+}(-1,-1)} \left(
\tilde{\eta}_{0}^{\spadesuit,+}(-1,-1;\tilde{\boldsymbol{\Omega}}) \! 
- \! \tilde{\eta}_{0}^{\spadesuit,+}(-1,-1;\boldsymbol{0}) \right), \\
\tilde{\mathbb{M}}_{22}^{2,+}(\alpha_{p_{\mathfrak{s}-1}}) =& \, 
-\dfrac{\tilde{\mathbb{F}}_{\tilde{\boldsymbol{\Omega}}}^{0,0,+}
(-1,-1)}{\tilde{\mathbb{F}}_{\boldsymbol{0}}^{0,0,+}(-1,-1)} \left(
\tilde{\eta}_{1}^{\spadesuit,+}(-1,-1;\tilde{\boldsymbol{\Omega}}) 
\! - \! \tilde{\eta}_{1}^{\spadesuit,+}(-1,-1;\boldsymbol{0}) \! 
- \! (\tilde{\eta}_{0}^{\spadesuit,+}(-1,-1;\boldsymbol{0}))^{2} 
\right. \\
+&\left. \, \tilde{\eta}_{0}^{\spadesuit,+}(-1,-1;
\tilde{\boldsymbol{\Omega}}) \tilde{\eta}_{0}^{\spadesuit,+}
(-1,-1;\boldsymbol{0}) \right),
\end{align*}
with, for $\varepsilon_{1},\varepsilon_{2} \! = \! \pm 1$,
\begin{gather*}
\tilde{\mathbb{F}}_{\tilde{\boldsymbol{\Omega}}}^{j_{1},j_{2},\pm}
(\varepsilon_{1},\varepsilon_{2}) \! := \! \sum_{m \in \mathbb{Z}^{N}} 
\me^{2 \pi \mi (m,\varepsilon_{1} \tilde{\boldsymbol{u}}_{\pm}(\infty)-
\frac{1}{2 \pi}((n-1)K+k) \tilde{\boldsymbol{\Omega}}+\varepsilon_{2} 
\tilde{\boldsymbol{d}})+\mi \pi (m,\tilde{\boldsymbol{\tau}}m)}
(\tilde{\lambda}_{2})^{j_{1}}(\tilde{\lambda}_{3})^{j_{2}}, \quad 
j_{1} \! = \! 0,1,2, \quad j_{2} \! = \! 0,1, \\
\tilde{\lambda}_{2} \! = \! -2 \pi \mi \sum_{j=1}^{N}m_{j} \tilde{c}_{j1}, 
\quad \quad \tilde{\lambda}_{3} \! = \! -\mi \pi \sum_{j=1}^{N}m_{j}
(\tilde{c}_{j2} \! - \! \tilde{\lambda}_{0} \tilde{c}_{j1}), \quad \quad 
\tilde{\lambda}_{0} \! = \! -\dfrac{1}{2} \sum_{j=1}^{N+1}
(\tilde{b}_{j-1} \! + \! \tilde{a}_{j}), 
\end{gather*}
$\tilde{c}_{j1}$ and $\tilde{c}_{j2}$, $j \! \in \! \lbrace 1,2,\dotsc,
N \rbrace$, obtained {}from Equations~\eqref{O1} and~\eqref{O2},
\begin{gather}
\tilde{\eta}_{0}^{\spadesuit,\pm}(\varepsilon_{1},\varepsilon_{2};
\tilde{\boldsymbol{\Omega}}) \! := \! \left(\tilde{\mathbb{F}}_{\tilde{
\boldsymbol{\Omega}}}^{0,0,\pm}(\varepsilon_{1},\varepsilon_{2}) 
\right)^{-1} \tilde{\mathbb{F}}_{\tilde{\boldsymbol{\Omega}}}^{1,0,
\pm}(\varepsilon_{1},\varepsilon_{2}), \nonumber \\
\tilde{\eta}_{1}^{\spadesuit,\pm}(\varepsilon_{1},\varepsilon_{2};
\tilde{\boldsymbol{\Omega}}) \! := \! \left(\tilde{\mathbb{F}}_{\tilde{
\boldsymbol{\Omega}}}^{0,0,\pm}(\varepsilon_{1},\varepsilon_{2}) 
\right)^{-1} \left(\tilde{\mathbb{F}}_{\tilde{\boldsymbol{\Omega}}}^{
0,1,\pm}(\varepsilon_{1},\varepsilon_{2}) \! \pm \! \dfrac{\varepsilon_{
1}}{2} \tilde{\mathbb{F}}_{\tilde{\boldsymbol{\Omega}}}^{2,0,\pm}
(\varepsilon_{1},\varepsilon_{2}) \right), \nonumber \\
\tilde{\Xi}_{0}^{\sharp}(\alpha_{p_{\mathfrak{s}-1}}) \! := \! -((n \! 
- \! 1)K \! + \! k) \int_{J_{f}} \ln (\lvert \xi \! - \! \alpha_{k} \rvert) 
\, \md \mu_{\widetilde{V}}^{f}(\xi) \! + \! \sum_{j=1}^{\mathfrak{s}-2} 
\varkappa_{nk \tilde{k}_{j}} \ln (\lvert \alpha_{p_{j}} \! - \! \alpha_{k} 
\rvert), \label{eqfivecapxiefin1} 
\end{gather} 
and $(\mathrm{M}_{2}(\mathbb{C}) \! \ni)$ $\tilde{\mathfrak{c}}_{
\mathfrak{s}-1}^{\mathcal{X}}(n,k,z_{o}) \! =_{\underset{z_{o}=
1+o(1)}{\mathscr{N},n \to \infty}} \! \mathcal{O}(1)$.
\end{ccccc}

\emph{Proof}. The proof of this Lemma~\ref{lem5.5} consists of two 
cases: (i) $n \! \in \! \mathbb{N}$ and $k \! \in \! \lbrace 1,2,\dotsc,
K \rbrace$ such that $\alpha_{p_{\mathfrak{s}}} \! := \! \alpha_{k} \! 
= \! \infty$; and (ii) $n \! \in \! \mathbb{N}$ and $k \! \in \! \lbrace 
1,2,\dotsc,K \rbrace$ such that $\alpha_{p_{\mathfrak{s}}} \! := \! 
\alpha_{k} \! \neq \! \infty$. Notwithstanding the fact that the scheme 
of the proof is, \emph{mutatis mutandis}, similar for both cases, 
case~(ii), nevertheless, is the more technically challenging of the two; 
therefore, without loss of generality, only the proof for case~(ii) is 
presented in detail, whilst case~(i) is proved analogously.

For $n \! \in \! \mathbb{N}$ and $k \! \in \! \lbrace 1,2,\dotsc,K 
\rbrace$ such that $\alpha_{p_{\mathfrak{s}}} \! := \! \alpha_{k} 
\! \neq \! \infty$, recall that, for $q \! \in \! \lbrace 1,\dotsc,
\mathfrak{s} \! - \! 2,\mathfrak{s} \rbrace$, $\alpha_{p_{q}} \! \in 
\! \mathbb{R} \setminus J_{f} \! = \! (-\infty,\tilde{b}_{0}) \cup 
\cup_{j=1}^{N}(\tilde{a}_{j},\tilde{b}_{j}) \cup (\tilde{a}_{N+1},
+\infty) \subset \tilde{\Upsilon}_{1} \cup \tilde{\Upsilon}_{2}$ (cf. 
Figure~\ref{figsectortil}); hence, via the corresponding invertible 
transformations of Lemmata~\ref{lem3.4}, \ref{lem4.1}, \ref{lem4.2}, 
\ref{lem4.5}, and~\ref{lem4.10}, one shows that, for $z \! \in \! 
\tilde{\Upsilon}_{1} \cup \tilde{\Upsilon}_{2}$, the solution of the 
associated monic MPC ORF RHP $(\mathcal{X}(z),\upsilon (z), \overline{
\mathbb{R}})$ (cf. Lemma~$\bm{\mathrm{RHP}_{\mathrm{MPC}}}$) is 
given by
\begin{equation} \label{eqlem5.5ay} 
\mathcal{X}(z) \! = \! 
\begin{cases}
\me^{\frac{n \tilde{\ell}}{2} \operatorname{ad}(\sigma_{3})} 
\tilde{\mathcal{R}}(z) \widetilde{\mathbb{K}} \, \tilde{\mathbb{M}}(z) 
\mathscr{E}^{\sigma_{3}} \me^{n(g^{f}(z)-\hat{\mathscr{P}}_{0}^{+}) 
\sigma_{3}}, &\text{$z \! \in \! \mathbb{C}_{+} \, (\subset \tilde{
\mathcal{Y}}^{+})$,} \\
-\mi \me^{\frac{n \tilde{\ell}}{2} \operatorname{ad}(\sigma_{3})} 
\tilde{\mathcal{R}}(z) \widetilde{\mathbb{K}} \, \tilde{\mathbb{M}}
(z) \sigma_{2} \mathscr{E}^{-\sigma_{3}} \me^{n(g^{f}(z)-
\hat{\mathscr{P}}_{0}^{-}) \sigma_{3}}, &\text{$z \! \in \! 
\mathbb{C}_{-} \, (\subset \tilde{\mathcal{Y}}^{-})$.}
\end{cases}
\end{equation}
For $n \! \in \! \mathbb{N}$ and $k \! \in \! \lbrace 1,2,\dotsc,K 
\rbrace$ such that $\alpha_{p_{\mathfrak{s}}} \! := \! \alpha_{k} 
\! \neq \! \infty$, recall the expression for $\tilde{\gamma}(z)$ 
defined by Equation~\eqref{eqmainfin10}: a careful analysis 
of the branch cuts shows that, for $q \! \in \! \lbrace 1,\dotsc,
\mathfrak{s} \! - \! 2,\mathfrak{s} \rbrace$,
\begin{equation} \label{eqlem5.5be} 
\tilde{\gamma}(z) \underset{(\tilde{\mathcal{Y}}^{\pm} \supset) \, 
\overline{\mathbb{C}}_{\pm} \ni z \to \alpha_{p_{q}}}{=} (-\mi)^{(1 \mp 1)/2} 
\tilde{\gamma}(\alpha_{p_{q}}) \left(1 \! + \! \alpha^{\triangle}
(z \! - \! \alpha_{p_{q}}) \! + \! \left(\beta^{\triangle}  \! + \! 
\dfrac{(\alpha^{\triangle})^{2}}{2} \right)(z \! - \! \alpha_{p_{q}})^{2} 
\! + \! \mathcal{O} \left(\tilde{\mathfrak{c}}^{1}_{q,\curlywedge}
(n,k,z_{o})(z \! - \! \alpha_{p_{q}})^{3} \right) \right),
\end{equation}
where $\tilde{\gamma}(\alpha_{p_{q}})$ is defined by Equation~\eqref{eqssabra2}, 
$\alpha^{\triangle}$ and $\beta^{\triangle}$ are defined in item~\pmb{(2)} 
of the lemma, and $\tilde{\mathfrak{c}}^{1}_{q,\curlywedge}(n,k,z_{o}) \! 
=_{\underset{z_{o}=1+o(1)}{\mathscr{N},n \to \infty}} \! \mathcal{O}(1)$. For 
$n \! \in \! \mathbb{N}$ and $k \! \in \! \lbrace 1,2,\dotsc,K \rbrace$ such 
that $\alpha_{p_{\mathfrak{s}}} \! := \! \alpha_{k} \! \neq \! \infty$, let, for 
$q \! \in \! \lbrace 1,\dotsc,\mathfrak{s} \! - \! 2,\mathfrak{s} \rbrace$, 
$\tilde{\boldsymbol{u}}(z) \! := \! \int_{\tilde{a}_{N+1}}^{z} 
\tilde{\boldsymbol{\omega}} \! = \! \tilde{\boldsymbol{u}}(\alpha_{p_{q}}) 
\! + \! (\int_{\alpha_{p_{q}}}^{z} \tilde{\omega}_{1},\int_{\alpha_{p_{q}}}^{z} 
\tilde{\omega}_{2},\dotsc,\int_{\alpha_{p_{q}}}^{z} \tilde{\omega}_{N})$, where 
$\tilde{\omega}_{j} \! := \! \sum_{i=1}^{N} \tilde{c}_{ji}(\tilde{R}(z))^{-1/2}
z^{N-i} \, \md z$, $j \! = \! 1,2,\dotsc,N$, with $(\tilde{R}(z))^{1/2}$ defined 
by Equation~\eqref{eql3.7j}, and $(\tilde{c})_{i_{1},i_{2}=1,2,\dotsc,N}$ 
obtained {}from Equations~\eqref{O1} and~\eqref{O2}: a careful analysis 
of the branch cuts shows that, for $q \! \in \! \lbrace 1,\dotsc,\mathfrak{s} 
\! - \! 2,\mathfrak{s} \rbrace$,
\begin{equation} \label{eqlem5.5de} 
\tilde{\boldsymbol{u}}(z) \underset{(\tilde{\mathcal{Y}}^{\pm} 
\supset) \, \mathbb{C}_{\pm} \ni z \to \alpha_{p_{q}}}{=} 
\tilde{\boldsymbol{u}}_{\pm}(\alpha_{p_{q}}) \! + \! \left(
\int_{\alpha_{p_{q}}^{\pm}}^{z} \tilde{\omega}_{1},
\int_{\alpha_{p_{q}}^{\pm}}^{z} \tilde{\omega}_{2},\dotsc,
\int_{\alpha_{p_{q}}^{\pm}}^{z} \tilde{\omega}_{N} \right),
\end{equation} 
where $\tilde{\boldsymbol{u}}_{\pm}(\alpha_{p_{q}}) \! := \! \tilde{
\boldsymbol{u}}(\alpha_{p_{q}}^{\pm}) \! = \! (\int_{\tilde{a}_{N+1}}^{
\alpha_{p_{q}}^{\pm}} \tilde{\omega}_{1},\int_{\tilde{a}_{N+1}}^{\alpha_{
p_{q}}^{\pm}} \tilde{\omega}_{2},\dotsc,\int_{\tilde{a}_{N+1}}^{\alpha_{
p_{q}}^{\pm}} \tilde{\omega}_{N})$, and
\begin{equation} \label{eqlem5.5ee} 
\int_{\alpha_{p_{q}}^{\pm}}^{z} \tilde{\omega}_{j} \underset{
(\tilde{\mathcal{Y}}^{\pm} \supset) \, \mathbb{C}_{\pm} \ni z \to 
\alpha_{p_{q}}}{=} \pm \left(\tilde{\mathbb{P}}^{\triangledown}_{j} \! + 
\! \tilde{\mathbb{Q}}^{\triangledown}_{j}(z \! - \! \alpha_{p_{q}}) \! + \! 
\tilde{\mathbb{R}}^{\triangledown}_{j}(z \! - \! \alpha_{p_{q}})^{2} \! + 
\! \mathcal{O} \left(\tilde{\mathfrak{c}}_{q,\curlywedge}^{2}(n,k,z_{o})
(z \! - \! \alpha_{p_{q}})^{3} \right) \right), \quad j \! = \! 1,2,\dotsc,N,
\end{equation}
where $\tilde{\mathbb{P}}^{\triangledown}_{j}$ and $\tilde{\mathbb{Q}}^{
\triangledown}_{j}$ are defined in item~\pmb{(2)} of the lemma,
\begin{equation*}
\tilde{\mathbb{R}}^{\triangledown}_{j} \! := \! \sum_{m=1}^{N} 
\dfrac{\tilde{c}_{jm}(\alpha_{p_{q}})^{N-m-2} \left(\frac{1}{2}(N \! - \! m)
(N \! - \! m \! - \! 1) \! - \! (N \! - \! m) \lambda_{0} \alpha_{p_{q}} \! - 
\! (\lambda_{1} \! - \! \lambda_{0}^{2}/2) \alpha_{p_{q}}^{2} \right)}{(-1)^{
\tilde{\mathfrak{n}}(\alpha_{p_{q}})} \prod_{i=1}^{N+1}(\lvert \alpha_{p_{q}} 
\! - \! \tilde{b}_{i-1} \rvert \lvert \alpha_{p_{q}} \! - \! \tilde{a}_{i} 
\rvert)^{1/2}},
\end{equation*}
with $\tilde{\mathfrak{n}}(\alpha_{p_{q}})$ defined in the corresponding 
item of Remark~\ref{remext}, $\lambda_{0}$ given in item~\pmb{(2)} 
of the lemma, $\lambda_{1} \! = \! -\tfrac{1}{4} \sum_{i=1}^{N+1}
((\tilde{b}_{i-1} \! - \! \alpha_{p_{q}})^{-2} \! + \! (\tilde{a}_{i} \! - \! 
\alpha_{p_{q}})^{-2})$, and $\tilde{\mathfrak{c}}^{2}_{q,\curlywedge}
(n,k,z_{o}) \! =_{\underset{z_{o}=1+o(1)}{\mathscr{N},n \to \infty}} \! 
\mathcal{O}(1)$. For $n \! \in \! \mathbb{N}$ and $k \! \in \! \lbrace 1,2,
\dotsc,K \rbrace$ such that $\alpha_{p_{\mathfrak{s}}} \! := \! \alpha_{k} 
\! \neq \! \infty$, recall Definition~\eqref{eqtilfrakm} for $\tilde{
\mathfrak{m}}(z)$: via the representation for the Riemann theta 
function $\tilde{\boldsymbol{\theta}}(\pmb{\cdot})$ defined by 
Equation~\eqref{eqrmthetafin}, recalling that the associated $N 
\times N$ Riemann matrix of $\tilde{\boldsymbol{\beta}}$-periods, 
$\tilde{\boldsymbol{\tau}} \! = \! (\tilde{\boldsymbol{\tau}})_{i_{1},i_{2}=
1,2,\dotsc,N}$, is non-degenerate, symmetric and pure imaginary, and 
$-\mi \tilde{\boldsymbol{\tau}}$ is positive definite, that is, $\mi \pi (m,
\tilde{\boldsymbol{\tau}}m) \! = \! \pi \sum_{j_{1}=1}^{N} \sum_{j_{2}=1}^{N}
m_{j_{1}}(\mi (\tilde{\boldsymbol{\tau}})_{j_{1},j_{2}})m_{j_{2}} \! < \! 0$, it 
follows via the Expansions~\eqref{eqlem5.5de} and~\eqref{eqlem5.5ee} 
that, for $\varepsilon_{1},\varepsilon_{2} \! = \! \pm 1$, and $q \! \in 
\! \lbrace 1,\dotsc,\mathfrak{s} \! - \! 2,\mathfrak{s} \rbrace$,
\begin{align} \label{eqlem5.5fe}
\dfrac{\tilde{\boldsymbol{\theta}}(\varepsilon_{1} \tilde{\boldsymbol{u}}
(z) \! - \! \frac{1}{2 \pi}((n \! - \! 1)K \! + \! k) \tilde{\boldsymbol{\Omega}} 
\! + \! \varepsilon_{2} \tilde{\boldsymbol{d}})}{\tilde{\boldsymbol{\theta}}
(\varepsilon_{1} \tilde{\boldsymbol{u}}(z) \! + \! \varepsilon_{2} \tilde{
\boldsymbol{d}})} \underset{(\tilde{\mathcal{Y}}^{\pm} \supset) \, 
\mathbb{C}_{\pm} \ni z \to \alpha_{p_{q}}}{=}& \, \dfrac{\tilde{
\mathbb{G}}_{\tilde{\boldsymbol{\Omega}}}^{0,0,\pm}(\varepsilon_{1},
\varepsilon_{2})}{\tilde{\mathbb{G}}_{\boldsymbol{0}}^{0,0,\pm}
(\varepsilon_{1},\varepsilon_{2})} \left(1 \! \pm \! \varepsilon_{1} \left(
\tilde{\eta}_{0}^{\blacklozenge,\pm}(\varepsilon_{1},\varepsilon_{2};
\tilde{\boldsymbol{\Omega}}) \! - \! \tilde{\eta}_{0}^{\blacklozenge,
\pm}(\varepsilon_{1},\varepsilon_{2};\boldsymbol{0}) \right)
(z \! - \! \alpha_{p_{q}}) \right. \nonumber \\
\pm&\left. \, \varepsilon_{1} \left(\tilde{\eta}_{1}^{\blacklozenge,\pm}
(\varepsilon_{1},\varepsilon_{2};\tilde{\boldsymbol{\Omega}}) \! - \! 
\tilde{\eta}_{1}^{\blacklozenge,\pm}(\varepsilon_{1},\varepsilon_{2};
\boldsymbol{0}) \! \pm \! \varepsilon_{1}(\tilde{\eta}_{0}^{\blacklozenge,
\pm}(\varepsilon_{1},\varepsilon_{2};\boldsymbol{0}))^{2} \right. \right. 
\nonumber \\
\mp&\left. \left. \, \varepsilon_{1} \tilde{\eta}_{0}^{\blacklozenge,\pm}
(\varepsilon_{1},\varepsilon_{2};\tilde{\boldsymbol{\Omega}}) \tilde{
\eta}_{0}^{\blacklozenge,\pm}(\varepsilon_{1},\varepsilon_{2};
\boldsymbol{0}) \right)(z \! - \! \alpha_{p_{q}})^{2} \right. \nonumber \\
+&\left. \, \mathcal{O} \left(\tilde{\mathfrak{c}}_{q,\curlywedge}^{3}
(n,k,z_{o})(z \! - \! \alpha_{p_{q}})^{3} \right) \right),
\end{align}
where $\tilde{\mathbb{G}}_{\tilde{\boldsymbol{\Omega}}}^{j_{1},j_{2},
\pm}(\varepsilon_{1},\varepsilon_{2})$, $j_{1} \! = \! 0,1,2$, $j_{2} \! 
= \! 0,1$, $\tilde{\eta}_{0}^{\blacklozenge,\pm}(\varepsilon_{1},
\varepsilon_{2};\tilde{\boldsymbol{\Omega}})$, and $\tilde{\eta}_{1}^{
\blacklozenge,\pm}(\varepsilon_{1},\varepsilon_{2};\tilde{\boldsymbol{
\Omega}})$ are defined in item~\pmb{(2)} of the lemma, and $\tilde{
\mathfrak{c}}^{3}_{q,\curlywedge}(n,k,z_{o}) \! =_{\underset{z_{o}=
1+o(1)}{\mathscr{N},n \to \infty}} \! \mathcal{O}(1)$. For $n \! \in \! 
\mathbb{N}$ and $k \! \in \! \lbrace 1,2,\dotsc,K \rbrace$ such that 
$\alpha_{p_{\mathfrak{s}}} \! := \! \alpha_{k} \! \neq \! \infty$, recall 
{}from the proof of Lemma~\ref{lem4.5} that the matrix $\tilde{
\mathfrak{m}}(z)$ defined by Equation~\eqref{eqtilfrakm} satisfies the 
jump relation $\tilde{\mathfrak{m}}_{+}(z) \! = \! \tilde{\mathfrak{m}}_{-}
(z)(\exp (-\mi ((n \! - \! 1)K \! + \! k) \tilde{\Omega}_{j^{\prime}}) 
\sigma_{-} \! + \! \exp (\mi ((n \! - \! 1)K \! + \! k) \tilde{\Omega}_{
j^{\prime}}) \sigma_{+})$ for $z \! \in \! (\tilde{a}_{j^{\prime}},
\tilde{b}_{j^{\prime}}) \setminus \lbrace \tilde{z}_{j^{\prime}}^{\pm} 
\rbrace$, $j^{\prime} \! = \! 1,2,\dotsc,N$, and $\tilde{\mathfrak{
m}}_{+}(z) \! = \! \tilde{\mathfrak{m}}_{-}(z) \sigma_{1}$ for $z \! 
\in \! (-\infty,\tilde{b}_{0}) \cup (\tilde{a}_{N+1},+\infty)$; via this 
latter observation, the definitions of $\widetilde{\mathbb{K}}$ and 
$\tilde{\mathbb{M}}(z)$ given in item~\pmb{(2)} of Lemma~\ref{lem4.5}, 
and the relations (cf. the proof of Lemma~\ref{lem4.5})
\begin{gather*}
\widetilde{\mathbb{K}} 
\begin{pmatrix}
\mathfrak{h}_{11}^{0,+}(\alpha_{k}) \tilde{\mathbb{M}}_{11}^{0,+}
(\alpha_{k}) & \mathfrak{h}_{12}^{0,+}(\alpha_{k}) \tilde{\mathbb{
M}}_{12}^{0,+}(\alpha_{k}) \\
-\mathfrak{h}_{12}^{0,+}(\alpha_{k}) \tilde{\mathbb{M}}_{21}^{0,+}
(\alpha_{k}) & \mathfrak{h}_{11}^{0,+}(\alpha_{k}) \tilde{\mathbb{
M}}_{22}^{0,+}(\alpha_{k}) 
\end{pmatrix} \! = \! \mathscr{E}^{-\sigma_{3}}, \\
\widetilde{\mathbb{K}} 
\begin{pmatrix}
\mathfrak{h}_{11}^{0,+}(\alpha_{k}) \tilde{\mathbb{M}}_{11}^{0,+}
(\alpha_{k}) & \mathfrak{h}_{12}^{0,+}(\alpha_{k}) \tilde{\mathbb{
M}}_{12}^{0,+}(\alpha_{k}) \\
-\mathfrak{h}_{12}^{0,+}(\alpha_{k}) \tilde{\mathbb{M}}_{21}^{0,+}
(\alpha_{k}) & \mathfrak{h}_{11}^{0,+}(\alpha_{k}) \tilde{\mathbb{
M}}_{22}^{0,+}(\alpha_{k}) 
\end{pmatrix} 
\sigma_{3} 
\begin{pmatrix}
0 & \me^{\mi ((n-1)K+k) \tilde{\Omega}_{j}} \\
\me^{-\mi ((n-1)K+k) \tilde{\Omega}_{j}} & 0
\end{pmatrix} \! = \! \mi \mathscr{E}^{\sigma_{3}} \sigma_{2},
\end{gather*}
where $\mathfrak{h}_{11}^{0,+}(\alpha_{k})$, $\mathfrak{h}_{12}^{0,+}
(\alpha_{k})$, and $\tilde{\mathbb{M}}_{i_{1}i_{2}}^{0,+}(\alpha_{k})$, 
$i_{1},i_{2} \! = \! 1,2$, are defined in item~\pmb{(2)} of the lemma, 
and $\mathscr{E}$ is defined by Equation~\eqref{eqmainfin13}, one shows, 
via the Expansions~\eqref{eqlem5.5be} and~\eqref{eqlem5.5fe}, and a 
lengthy, but otherwise straightforward, matrix-multiplication argument, 
that
\begin{align}
\widetilde{\mathbb{K}} \, \tilde{\mathbb{M}}(z) \mathscr{E}^{\sigma_{3}} 
\underset{(\tilde{\mathcal{Y}}^{+} \supset) \, \mathbb{C}_{+} \ni z 
\to \alpha_{k}}{=}& \, \mathrm{I} \! + \! A_{0}(\alpha_{k})(z \! - \! 
\alpha_{k}) \! + \! B_{0}(\alpha_{k})(z \! - \! \alpha_{k})^{2} \! + \! 
\mathcal{O} \left(\tilde{\mathfrak{c}}_{\mathfrak{s},\curlywedge}^{4}
(n,k,z_{o})(z \! - \! \alpha_{k})^{3} \right), \label{eqlem5.5ge} \\
-\mi \widetilde{\mathbb{K}} \, \tilde{\mathbb{M}}(z) \sigma_{2} 
\mathscr{E}^{-\sigma_{3}} \underset{(\tilde{\mathcal{Y}}^{-} \supset) 
\, \mathbb{C}_{-} \ni z \to \alpha_{k}}{=}& \, \mathrm{I} \! + \! A_{0}
(\alpha_{k})(z \! - \! \alpha_{k}) \! + \! B_{0}(\alpha_{k})(z \! - \! 
\alpha_{k})^{2} \! + \! \mathcal{O} \left(\tilde{\mathfrak{c}}_{\mathfrak{s},
\curlywedge}^{5}(n,k,z_{o})(z \! - \! \alpha_{k})^{3} \right), \label{eqlem5.5he} 
\end{align}
where $A_{0}(\alpha_{k})$ and $B_{0}(\alpha_{k})$ are defined in 
item~\pmb{(2)} of the lemma, $(\mathrm{M}_{2}(\mathbb{C}) \! \ni)$ 
$\tilde{\mathfrak{c}}^{r}_{\mathfrak{s},\curlywedge}(n,k,z_{o}) \! 
=_{\underset{z_{o}=1+o(1)}{\mathscr{N},n \to \infty}} \! \mathcal{O}(1)$, 
$r \! = \! 4,5$, and, for $q \! \in \! \lbrace 1,2,\dotsc,\mathfrak{s} \! 
- \! 2 \rbrace$,
\begin{align}
\widetilde{\mathbb{K}} \, \tilde{\mathbb{M}}(z) \mathscr{E}^{\sigma_{3}} 
\underset{(\tilde{\mathcal{Y}}^{+} \supset) \, \mathbb{C}_{+} \ni z \to 
\alpha_{p_{q}}}{=}& \, \left(\tilde{A}_{0}(\alpha_{p_{q}}) \! + \! 
\tilde{B}_{0}(\alpha_{p_{q}})(z \! - \! \alpha_{p_{q}}) \! + \! \tilde{C}_{0}
(\alpha_{p_{q}})(z \! - \! \alpha_{p_{q}})^{2} \! + \! \mathcal{O} 
\left(\tilde{\mathfrak{c}}_{q,\curlywedge}^{6}(n,k,z_{o})(z \! - \! 
\alpha_{p_{q}})^{3} \right) \right) \tilde{\Gamma}^{\triangledown}_{+}, 
\label{eqlem5.5ie} \\
-\mi \widetilde{\mathbb{K}} \, \tilde{\mathbb{M}}(z) \sigma_{2} 
\mathscr{E}^{-\sigma_{3}} \underset{(\tilde{\mathcal{Y}}^{-} \supset) 
\, \mathbb{C}_{-} \ni z \to \alpha_{p_{q}}}{=}& \, \left(\tilde{A}_{0}
(\alpha_{p_{q}}) \! + \! \tilde{B}_{0}(\alpha_{p_{q}})(z \! - \! \alpha_{p_{q}}) 
\! + \! \tilde{C}_{0}(\alpha_{p_{q}})(z \! - \! \alpha_{p_{q}})^{2} \! + \! 
\mathcal{O} \left(\tilde{\mathfrak{c}}_{q,\curlywedge}^{7}(n,k,z_{o})(z \! - 
\! \alpha_{p_{q}})^{3} \right) \right) \tilde{\Gamma}^{\triangledown}_{-}, 
\label{eqlem5.5je} 
\end{align}
where $\tilde{A}_{0}(\alpha_{p_{q}})$, $\tilde{B}_{0}(\alpha_{p_{q}})$, and 
$\tilde{C}_{0}(\alpha_{p_{q}})$ are defined in item~\pmb{(2)} of the lemma, 
$\tilde{\Gamma}^{\triangledown}_{+} \! = \! \mathscr{E}^{\sigma_{3}}$, 
$\tilde{\Gamma}^{\triangledown}_{-} \! = \! \mathscr{E}^{-\sigma_{3}} \exp 
(\mi ((n \! - \! 1)K \! + \! k) \tilde{\Omega}_{j} \sigma_{3})$,\footnote{If, for 
$q \! \in \! \lbrace 1,2,\dotsc,\mathfrak{s} \! - \! 2 \rbrace$, $\alpha_{p_{q}} 
\! \in \! (-\infty,\tilde{b}_{0}) \cup (\tilde{a}_{N+1},+\infty)$, then 
$\tilde{\Omega}_{j} \! = \! 0$.} and $(\mathrm{M}_{2}(\mathbb{C}) \! \ni)$ 
$\tilde{\mathfrak{c}}^{r}_{q,\curlywedge}(n,k,z_{o}) \! =_{\underset{z_{o}=1
+o(1)}{\mathscr{N},n \to \infty}} \! \mathcal{O}(1)$, $r \! = \! 6,7$. For $n 
\! \in \! \mathbb{N}$ and $k \! \in \! \lbrace 1,2,\dotsc,K \rbrace$ such 
that $\alpha_{p_{\mathfrak{s}}} \! := \! \alpha_{k} \! \neq \! \infty$, via 
Equation~\eqref{eqpropo5.3tB}, one shows that
\begin{equation} \label{eqlem5.5ke} 
\tilde{\mathcal{R}}(z) \underset{\mathbb{C}_{\pm} \ni z \to \alpha_{k}}{=} 
\mathrm{I} \! + \! \tilde{\mathcal{R}}_{0}^{A_{0}}(\alpha_{k})(z \! - \! 
\alpha_{k}) \! + \! \tilde{\mathcal{R}}_{0}^{B_{0}}(\alpha_{k})(z \! - \! 
\alpha_{k})^{2} \! + \! \mathcal{O} \left(\tilde{\mathfrak{c}}_{\mathfrak{s},
\curlywedge}^{8}(n,k,z_{o})(z \! - \! \alpha_{k})^{3} \right),
\end{equation}
and, for $q \! \in \! \lbrace 1,2,\dotsc,\mathfrak{s} \! - \! 2 \rbrace$,
\begin{equation} \label{eqlem5.5le}
\tilde{\mathcal{R}}(z) \underset{\mathbb{C}_{\pm} \ni z \to \alpha_{p_{q}}}{=} 
\mathrm{I} \! + \! \tilde{\mathcal{R}}_{0}^{\tilde{A}_{0}}(\alpha_{p_{q}}) 
\! + \! \tilde{\mathcal{R}}_{0}^{\tilde{B}_{0}}(\alpha_{p_{q}})(z \! - \! 
\alpha_{p_{q}}) \! + \! \tilde{\mathcal{R}}_{0}^{\tilde{C}_{0}}(\alpha_{p_{q}})
(z \! - \! \alpha_{p_{q}})^{2} \! + \! \mathcal{O} \left(\tilde{\mathfrak{c}}_{q,
\curlywedge}^{9}(n,k,z_{o})(z \! - \! \alpha_{p_{q}})^{3} \right),
\end{equation}
where, in the double-scaling limit $\mathscr{N},n \! \to \! \infty$ 
such that $z_{o} \! = \! 1 \! + \! o(1)$, $\tilde{\mathcal{R}}_{0}^{A_{0}}
(\alpha_{k})$, $\tilde{\mathcal{R}}_{0}^{B_{0}}(\alpha_{k})$, $\tilde{
\mathcal{R}}_{0}^{\tilde{A}_{0}}(\alpha_{p_{q}})$, $\tilde{\mathcal{
R}}_{0}^{\tilde{B}_{0}}(\alpha_{p_{q}})$, and $\tilde{\mathcal{R}}_{
0}^{\tilde{C}_{0}}(\alpha_{p_{q}})$ are given in item~\pmb{(2)} of the lemma, 
and $(\mathrm{M}_{2}(\mathbb{C}) \! \ni)$ $\tilde{\mathfrak{c}}^{r}_{
q^{\prime},\curlywedge}(n,k,z_{o}) \! =_{\underset{z_{o}=1+o(1)}{
\mathscr{N},n \to \infty}} \! \mathcal{O}(1)$, $q^{\prime} \! \in \! \lbrace 1,
\dotsc,\mathfrak{s} \! - \! 2,\mathfrak{s} \rbrace$, $r \! = \! 8,9$. For $n 
\! \in \! \mathbb{N}$ and $k \! \in \! \lbrace 1,2,\dotsc,K \rbrace$ such 
that $\alpha_{p_{\mathfrak{s}}} \! := \! \alpha_{k} \! \neq \! \infty$, via 
Equations~\eqref{eql3.4gee4}, \eqref{eql3.4gee5}, \eqref{eql3.4b}, 
and~\eqref{eql3.4d}--\eqref{eql3.4gee9}, one shows that
\begin{equation} \label{eqlem5.5me} 
\me^{n(g^{f}(z)-\hat{\mathscr{P}}_{0}^{\pm}) \sigma_{3}} 
\underset{(\tilde{\mathcal{Y}}^{\pm} \supset) \, \mathbb{C}_{\pm} 
\ni z \to \alpha_{k}}{=} \me^{-(\varkappa_{nk}-1) \ln (z-\alpha_{k}) 
\sigma_{3}} \me^{\tilde{w}_{0}^{\triangle}(\alpha_{k})(z-\alpha_{k}) 
\sigma_{3}} \me^{\tilde{w}_{1}^{\triangle}(\alpha_{k})(z-\alpha_{k})^{2} 
\sigma_{3}} \me^{\mathcal{O}(\tilde{\mathfrak{c}}_{\mathfrak{s},
\curlywedge}^{10}(n,k,z_{o})(z-\alpha_{k})^{3} \sigma_{3})},
\end{equation}
and, for $q \! \in \! \lbrace 1,2,\dotsc,\mathfrak{s} \! - \! 2 \rbrace$,
\begin{align}
\me^{n(g^{f}(z)-\hat{\mathscr{P}}_{0}^{\pm}) \sigma_{3}} \underset{
(\tilde{\mathcal{Y}}^{\pm} \supset) \, \mathbb{C}_{\pm} \ni z \to 
\alpha_{p_{q}}}{=}& \, \me^{-\varkappa_{nk \tilde{k}_{q}} \ln (z-
\alpha_{p_{q}}) \sigma_{3}} \me^{\tilde{w}_{0}^{\triangle}(\alpha_{p_{q}})
(z-\alpha_{p_{q}}) \sigma_{3}} \me^{\tilde{w}_{1}^{\triangle}(\alpha_{p_{q}})
(z-\alpha_{p_{q}})^{2} \sigma_{3}} \me^{\pm \mi \pi ((n-1)K+k) \int_{J_{f} 
\cap \mathbb{R}_{\alpha_{p_{q}}}^{>}} \md \mu_{\widetilde{V}}^{f}(\xi) 
\sigma_{3}} \nonumber \\
\times& \, \me^{\tilde{\Xi}_{0}(\alpha_{p_{q}}) \sigma_{3}} \mathscr{E}^{
\mp \sigma_{3}}(-1)^{(\sum_{j \in \hat{\Delta}_{f}(k)} \varkappa_{nk 
\tilde{k}_{j}}+\sum_{j \in \tilde{\Delta}_{f}(q)} \varkappa_{nk \tilde{k}_{j}}
+(\varkappa_{nk}-1) \epsilon (k,q)) \sigma_{3}} \me^{\mathcal{O}
(\tilde{\mathfrak{c}}_{q,\curlywedge}^{11}(n,k,z_{o})(z-\alpha_{p_{q}})^{3} 
\sigma_{3})}, \label{eqlem5.5ne}
\end{align}
where $\tilde{w}_{0}^{\triangle}(\alpha_{k})$, $\tilde{w}_{1}^{\triangle}
(\alpha_{k})$, $\hat{\Delta}_{f}(k)$, $\tilde{w}_{0}^{\triangle}(\alpha_{
p_{q}})$, $\tilde{w}_{1}^{\triangle}(\alpha_{p_{q}})$, $\tilde{\Xi}_{0}
(\alpha_{p_{q}})$, $\tilde{\Delta}_{f}(q)$, and $\epsilon (k,q)$ are defined 
in item~\pmb{(2)} of the lemma, and $(\mathrm{M}_{2}(\mathbb{C}) 
\! \ni)$ $\tilde{\mathfrak{c}}^{r}_{q^{\prime},\curlywedge}(n,k,
z_{o}) \! =_{\underset{z_{o}=1+o(1)}{\mathscr{N},n \to \infty}} \! 
\mathcal{O}(1)$, $q^{\prime} \! \in \! \lbrace 1,\dotsc,\mathfrak{s} \! 
- \! 2,\mathfrak{s} \rbrace$, $r \! = \! 10,11$. Hence, for $n \! \in \! 
\mathbb{N}$ and $k \! \in \! \lbrace 1,2,\dotsc,K \rbrace$ such that 
$\alpha_{p_{\mathfrak{s}}} \! := \! \alpha_{k} \! \neq \! \infty$, via 
Equation~\eqref{eqlem5.5ay}, and the Expansions~\eqref{eqlem5.5ge}, 
\eqref{eqlem5.5he}, \eqref{eqlem5.5ke}, and~\eqref{eqlem5.5me}, 
one arrives at, after a matrix-multiplication argument, in the 
double-scaling limit $\mathscr{N},n \! \to \! \infty$ such that $z_{o} 
\! = \! 1 \! + \! o(1)$, the Asymptotics~\eqref{eqlem5.5C}, and, via 
Equation~\eqref{eqlem5.5ay}, the Expansions~\eqref{eqlem5.5ie}, 
\eqref{eqlem5.5je}, \eqref{eqlem5.5le}, and~\eqref{eqlem5.5ne}, and 
the relation
\begin{equation*}
\tilde{\Gamma}^{\triangledown}_{\pm} \mathscr{E}^{\mp \sigma_{3}} 
\me^{\pm \mi \pi ((n-1)K+k) \int_{J_{f} \cap \mathbb{R}_{\alpha_{p_{q}}}^{>}} 
\md \mu_{\widetilde{V}}^{f}(\xi) \sigma_{3}} \! = \! \me^{\mi \pi ((n-1)K+k) 
\int_{J_{f} \cap \mathbb{R}_{\alpha_{p_{q}}}^{>}} \md \mu_{\widetilde{V}}^{f}
(\xi) \sigma_{3}},
\end{equation*}
one arrives at, after a matrix-multiplication argument, in the double-scaling 
limit $\mathscr{N},n \! \to \! \infty$ such that $z_{o} \! = \! 1 \! + \! o(1)$, 
the Asymptotics~\eqref{eqlem5.5D}.

For $n \! \in \! \mathbb{N}$ and $k \! \in \! \lbrace 1,2,\dotsc,K \rbrace$ 
such that $\alpha_{p_{\mathfrak{s}}} \! := \! \alpha_{k} \! \neq \! \infty$, 
the case $z \! \to \! \alpha_{p_{\mathfrak{s}-1}} \! = \! \infty$ is analysed 
similarly; in particular, the analogues of the 
Expansions~\eqref{eqlem5.5be}--\eqref{eqlem5.5ne} read:
\begin{equation} \label{eqlem5.5oe} 
\tilde{\gamma}(z) \underset{(\tilde{\mathcal{Y}}^{\pm} \supset) \, 
\overline{\mathbb{C}}_{\pm} \ni z \to \alpha_{p_{\mathfrak{s}-1}} 
= \infty}{=} (-\mi)^{(1 \mp 1)/2} \left(1 \! + \! \dfrac{1}{z} 
\tilde{\alpha}_{0}^{\triangle} \! + \! \dfrac{1}{z^{2}} \left(\tilde{\beta}_{
0}^{\triangle}  \! + \! \dfrac{(\tilde{\alpha}_{0}^{\triangle})^{2}}{2} \right) 
\! + \! \mathcal{O} \left(\tilde{\mathfrak{c}}^{1}_{\mathfrak{s}-1,
\curlyvee}(n,k,z_{o})z^{-3} \right) \right),
\end{equation}
where $\tilde{\alpha}_{0}^{\triangle}$ and $\tilde{\beta}_{0}^{\triangle}$ 
are defined in item~\pmb{(2)} of the lemma, and $\tilde{\mathfrak{c}}^{
1}_{\mathfrak{s}-1,\curlyvee}(n,k,z_{o}) \! =_{\underset{z_{o}=1
+o(1)}{\mathscr{N},n \to \infty}} \! \mathcal{O}(1)$,
\begin{equation} \label{eqlem5.5qe} 
\tilde{\boldsymbol{u}}(z) \underset{(\tilde{\mathcal{Y}}^{\pm} \supset) 
\, \overline{\mathbb{C}}_{\pm} \ni z \to \alpha_{p_{\mathfrak{s}-1}} = 
\infty}{=} \tilde{\boldsymbol{u}}_{\pm}(\infty) \! + \! \left(\int_{\infty^{
\pm}}^{z} \tilde{\omega}_{1},\int_{\infty^{\pm}}^{z} \tilde{\omega}_{2},
\dotsc,\int_{\infty^{\pm}}^{z} \tilde{\omega}_{N} \right),
\end{equation} 
where $\tilde{\boldsymbol{u}}_{\pm}(\infty)$ $(= \! \tilde{\boldsymbol{
u}}_{\pm}(\alpha_{p_{\mathfrak{s}-1}}))$ $:= \! \tilde{\boldsymbol{u}}
(\infty^{\pm}) \! = \! (\int_{\tilde{a}_{N+1}}^{\infty^{\pm}} \tilde{
\omega}_{1},\int_{\tilde{a}_{N+1}}^{\infty^{\pm}} \tilde{\omega}_{2},
\dotsc,\int_{\tilde{a}_{N+1}}^{\infty^{\pm}} \tilde{\omega}_{N})$, and
\begin{equation} \label{eqlem5.5re} 
\int_{\infty^{\pm}}^{z} \tilde{\omega}_{j} \underset{(\tilde{\mathcal{Y}}^{\pm} 
\supset) \, \overline{\mathbb{C}}_{\pm} \ni z \to \alpha_{p_{\mathfrak{s}-1}} 
= \infty}{=} \pm \left(-\dfrac{1}{z} \tilde{c}_{j1} \! - \! \dfrac{1}{z^{2}} 
\dfrac{(\tilde{c}_{j2} \! - \! \tilde{\lambda}_{0} \tilde{c}_{j1})}{2} \! + \! 
\mathcal{O} \left(\tilde{\mathfrak{c}}_{\mathfrak{s}-1,\curlyvee}^{2}(n,k,z_{o})
z^{-3} \right) \right), \quad j \! = \! 1,2,\dotsc,N,
\end{equation}
with $\tilde{c}_{j1}$, $\tilde{c}_{j2}$, and $\tilde{\lambda}_{0}$ defined in 
item~\pmb{(2)} of the lemma, and $\tilde{\mathfrak{c}}^{2}_{\mathfrak{s}
-1,\curlyvee}(n,k,z_{o}) \! =_{\underset{z_{o}=1+o(1)}{\mathscr{N},n \to 
\infty}} \! \mathcal{O}(1)$,
\begin{align} \label{eqlem5.5se}
\dfrac{\tilde{\boldsymbol{\theta}}(\varepsilon_{1} \tilde{\boldsymbol{u}}
(z) \! - \! \frac{1}{2 \pi}((n \! - \! 1)K \! + \! k) \tilde{\boldsymbol{\Omega}} 
\! + \! \varepsilon_{2} \tilde{\boldsymbol{d}})}{\tilde{\boldsymbol{\theta}}
(\varepsilon_{1} \tilde{\boldsymbol{u}}(z) \! + \! \varepsilon_{2} \tilde{
\boldsymbol{d}})} \underset{(\tilde{\mathcal{Y}}^{\pm} \supset) \, 
\overline{\mathbb{C}}_{\pm} \ni z \to \alpha_{p_{\mathfrak{s}-1}} = 
\infty}{=}& \, \dfrac{\tilde{\mathbb{F}}_{\tilde{\boldsymbol{
\Omega}}}^{0,0,\pm}(\varepsilon_{1},\varepsilon_{2})}{\tilde{\mathbb{
F}}_{\boldsymbol{0}}^{0,0,\pm}(\varepsilon_{1},\varepsilon_{2})} \left(
1 \! \pm \! \dfrac{1}{z} \varepsilon_{1} \left(\tilde{\eta}_{0}^{\spadesuit,
\pm}(\varepsilon_{1},\varepsilon_{2};\tilde{\boldsymbol{\Omega}}) \! - 
\! \tilde{\eta}_{0}^{\spadesuit,\pm}(\varepsilon_{1},\varepsilon_{2};
\boldsymbol{0}) \right) \right. \nonumber \\
\pm&\left. \, \dfrac{1}{z^{2}} \varepsilon_{1} \left(\tilde{\eta}_{1}^{
\spadesuit,\pm}(\varepsilon_{1},\varepsilon_{2};\tilde{\boldsymbol{\Omega}}) 
\! - \! \tilde{\eta}_{1}^{\spadesuit,\pm}(\varepsilon_{1},\varepsilon_{2};
\boldsymbol{0}) \! \pm \! \varepsilon_{1}(\tilde{\eta}_{0}^{\spadesuit,\pm}
(\varepsilon_{1},\varepsilon_{2};\boldsymbol{0}))^{2} \right. \right. 
\nonumber \\
\mp&\left. \left. \, \varepsilon_{1} \tilde{\eta}_{0}^{\spadesuit,\pm}
(\varepsilon_{1},\varepsilon_{2};\tilde{\boldsymbol{\Omega}}) \tilde{
\eta}_{0}^{\spadesuit,\pm}(\varepsilon_{1},\varepsilon_{2};\boldsymbol{0}) 
\right) \! + \! \mathcal{O} \left(\tilde{\mathfrak{c}}_{\mathfrak{s}-1,
\curlyvee}^{3}(n,k,z_{o})z^{-3} \right) \right),
\end{align}
where, for $\varepsilon_{1},\varepsilon_{2} \! = \! \pm 1$, 
$\tilde{\mathbb{F}}_{\tilde{\boldsymbol{\Omega}}}^{j_{1},j_{2},\pm}
(\varepsilon_{1},\varepsilon_{2})$, $j_{1} \! = \! 0,1,2$, $j_{2} \! = \! 0,1$, 
$\tilde{\eta}_{0}^{\spadesuit,\pm}(\varepsilon_{1},\varepsilon_{2};
\tilde{\boldsymbol{\Omega}})$, and $\tilde{\eta}_{1}^{\spadesuit,\pm}
(\varepsilon_{1},\varepsilon_{2};\tilde{\boldsymbol{\Omega}})$ are defined 
in item~\pmb{(2)} of the lemma, and $\tilde{\mathfrak{c}}^{3}_{\mathfrak{s}
-1,\curlyvee}(n,k,z_{o}) \! =_{\underset{z_{o}=1+o(1)}{\mathscr{N},n \to 
\infty}} \! \mathcal{O}(1)$,
\begin{align}
\widetilde{\mathbb{K}} \, \tilde{\mathbb{M}}(z) \mathscr{E}^{\sigma_{3}} 
\underset{(\tilde{\mathcal{Y}}^{+} \supset) \, \overline{\mathbb{C}}_{+} 
\ni z \to \alpha_{p_{\mathfrak{s}-1}} = \infty}{=}& \, \left(\tilde{A}_{0}^{
\sharp}(\alpha_{p_{\mathfrak{s}-1}}) \! + \! \dfrac{1}{z} \tilde{B}_{0}^{
\sharp}(\alpha_{p_{\mathfrak{s}-1}}) \! + \! \dfrac{1}{z^{2}} \tilde{C}_{
0}^{\sharp}(\alpha_{p_{\mathfrak{s}-1}}) \! + \! \mathcal{O} \left(
\tilde{\mathfrak{c}}_{\mathfrak{s}-1,\curlyvee}^{4}(n,k,z_{o})z^{-3} 
\right) \right) \mathscr{E}^{\sigma_{3}}, \label{eqlem5.5te} \\
-\mi \widetilde{\mathbb{K}} \, \tilde{\mathbb{M}}(z) \sigma_{2} 
\mathscr{E}^{-\sigma_{3}} \underset{(\tilde{\mathcal{Y}}^{-} \supset) \, 
\overline{\mathbb{C}}_{-} \ni z \to \alpha_{p_{\mathfrak{s}-1}} = \infty}{
=}& \, \left(\tilde{A}_{0}^{\sharp}(\alpha_{p_{\mathfrak{s}-1}}) \! + \! 
\dfrac{1}{z} \tilde{B}_{0}^{\sharp}(\alpha_{p_{\mathfrak{s}-1}}) \! + \! 
\dfrac{1}{z^{2}} \tilde{C}_{0}^{\sharp}(\alpha_{p_{\mathfrak{s}-1}}) \! + 
\! \mathcal{O} \left(\tilde{\mathfrak{c}}_{\mathfrak{s}-1,\curlyvee}^{5}
(n,k,z_{o})z^{-3} \right) \right) \mathscr{E}^{-\sigma_{3}}, \label{eqlem5.5ue}
\end{align}
where $\tilde{A}_{0}^{\sharp}(\alpha_{p_{\mathfrak{s}-1}})$, 
$\tilde{B}_{0}^{\sharp}(\alpha_{p_{\mathfrak{s}-1}})$, and 
$\tilde{C}_{0}^{\sharp}(\alpha_{p_{\mathfrak{s}-1}})$ are defined in 
item~\pmb{(2)} of the lemma, and $(\mathrm{M}_{2}(\mathbb{C}) \! 
\ni)$ $\tilde{\mathfrak{c}}^{r}_{\mathfrak{s}-1,\curlyvee}(n,k,z_{o}) \! 
=_{\underset{z_{o}=1+o(1)}{\mathscr{N},n \to \infty}} \! \mathcal{O}(1)$, 
$r \! = \! 4,5$,
\begin{align}
\tilde{\mathcal{R}}(z) \underset{\overline{\mathbb{C}}_{\pm} \ni z 
\to \alpha_{p_{\mathfrak{s}-1}} = \infty}{=}& \, \mathrm{I} \! - \! 
\tilde{\mathcal{R}}_{0}^{\tilde{A}_{0}^{\sharp}}(\alpha_{p_{\mathfrak{s}
-1}}) \! + \! \dfrac{1}{z} \tilde{\mathcal{R}}_{0}^{\tilde{B}_{0}^{\sharp}}
(\alpha_{p_{\mathfrak{s}-1}}) \! + \! \dfrac{1}{z^{2}} \tilde{\mathcal{
R}}_{0}^{\tilde{C}_{0}^{\sharp}}(\alpha_{p_{\mathfrak{s}-1}}) \! + \! 
\mathcal{O} \left(\tilde{\mathfrak{c}}_{\mathfrak{s}-1,\curlyvee}^{6}
(n,k,z_{o})z^{-3} \right), \label{eqlem5.5ve}
\end{align}
where, in the double-scaling limit $\mathscr{N},n \! \to \! \infty$ such 
that $z_{o} \! = \! 1 \! + \! o(1)$, $\tilde{\mathcal{R}}_{0}^{\tilde{A}_{
0}^{\sharp}}(\alpha_{p_{\mathfrak{s}-1}})$, $\tilde{\mathcal{R}}_{0}^{
\tilde{B}_{0}^{\sharp}}(\alpha_{p_{\mathfrak{s}-1}})$, and $\tilde{
\mathcal{R}}_{0}^{\tilde{C}_{0}^{\sharp}}(\alpha_{p_{\mathfrak{s}-1}})$ 
are given in item~\pmb{(2)} of the lemma, and $(\mathrm{M}_{2}
(\mathbb{C}) \! \ni)$ $\tilde{\mathfrak{c}}^{6}_{\mathfrak{s}-1,
\curlyvee}(n,k,z_{o}) \! =_{\underset{z_{o}=1+o(1)}{\mathscr{N},n 
\to \infty}} \! \mathcal{O}(1)$,
\begin{align} \label{eqlem5.5we} 
\me^{n(g^{f}(z)-\hat{\mathscr{P}}_{0}^{\pm}) \sigma_{3}} \underset{
(\tilde{\mathcal{Y}}^{\pm} \supset) \, \overline{\mathbb{C}}_{\pm} \ni 
z \to \alpha_{p_{\mathfrak{s}-1}} = \infty}{=}& \, \me^{(\varkappa_{nk 
\tilde{k}_{\mathfrak{s}-1}}^{\infty}+1) \ln (z) \sigma_{3}} \me^{\tilde{
w}_{0}^{\sharp}(\alpha_{p_{\mathfrak{s}-1}})z^{-1} \sigma_{3}} 
\me^{\tilde{w}_{1}^{\sharp}(\alpha_{p_{\mathfrak{s}-1}})z^{-2} 
\sigma_{3}} \me^{\tilde{\Xi}_{0}^{\sharp}(\alpha_{p_{\mathfrak{s}
-1}}) \sigma_{3}} \mathscr{E}^{\mp \sigma_{3}} \nonumber \\
\times& \, (-1)^{\sum_{j \in \hat{\Delta}_{f}(k)} \varkappa_{nk 
\tilde{k}_{j}} \sigma_{3}} \me^{\mathcal{O}(\tilde{\mathfrak{c}}_{
\mathfrak{s}-1,\curlyvee}^{7}(n,k,z_{o})z^{-3} \sigma_{3})},
\end{align}
where $\tilde{w}_{0}^{\sharp}(\alpha_{p_{\mathfrak{s}-1}})$, 
$\tilde{w}_{1}^{\sharp}(\alpha_{p_{\mathfrak{s}-1}})$, and 
$\tilde{\Xi}_{0}^{\sharp}(\alpha_{p_{\mathfrak{s}-1}})$ are defined in 
item~\pmb{(2)} of the lemma, and $(\mathrm{M}_{2}(\mathbb{C}) 
\! \ni)$ $\tilde{\mathfrak{c}}^{7}_{\mathfrak{s}-1,\curlyvee}
(n,k,z_{o}) \! =_{\underset{z_{o}=1+o(1)}{\mathscr{N},n \to \infty}} \! 
\mathcal{O}(1)$. Hence, for $n \! \in \! \mathbb{N}$ and $k \! \in \! 
\lbrace 1,2,\dotsc,K \rbrace$ such that $\alpha_{p_{\mathfrak{s}}} \! 
:= \! \alpha_{k} \! \neq \! \infty$, via Equation~\eqref{eqlem5.5ay}, 
and the Expansions~\eqref{eqlem5.5te}--\eqref{eqlem5.5we}, 
one arrives at, after a matrix-multiplication argument, in the 
double-scaling limit $\mathscr{N},n \! \to \! \infty$ such that 
$z_{o} \! = \! 1 \! + \! o(1)$, the Asymptotics~\eqref{eqlem5.5E}.

The analysis for the case $n \! \in \! \mathbb{N}$ and $k \! \in \! 
\lbrace 1,2,\dotsc,K \rbrace$ such that $\alpha_{p_{\mathfrak{s}}} 
\! := \! \alpha_{k} \! = \! \infty$ is, \emph{mutatis mutandis}, 
analogous, and leads to, in the double-scaling limit $\mathscr{N},
n \! \to \! \infty$ such that $z_{o} \! = \! 1 \! + \! o(1)$, the 
Asymptotics~\eqref{eqlem5.5A} and~\eqref{eqlem5.5B}. \hfill $\qed$
\begin{ccccc} \label{lem5.6} 
Let the external field $\widetilde{V} \colon \overline{\mathbb{R}} 
\setminus \lbrace \alpha_{1},\alpha_{2},\dotsc,\alpha_{K} \rbrace \! 
\to \! \mathbb{R}$ satisfy conditions~\eqref{eq20}--\eqref{eq22} and 
be regular. For $n \! \in \! \mathbb{N}$ and $k \! \in \! \lbrace 1,2,
\dotsc,K \rbrace$ such that $\alpha_{p_{\mathfrak{s}}} \! := \! 
\alpha_{k} \! = \! \infty$ (resp., $\alpha_{p_{\mathfrak{s}}} \! := \! 
\alpha_{k} \! \neq \! \infty)$, let the associated equilibrium measure, 
$\mu_{\widetilde{V}}^{\infty}$ (resp., $\mu_{\widetilde{V}}^{f})$, and its 
support, $J_{\infty}$ (resp., $J_{f})$, be as described in item~$\pmb{(1)}$ 
(resp., item~$\pmb{(2)})$ of Lemma~\ref{lem3.7}, and, along with the 
corresponding variational constant, $\hat{\ell}$ (resp., $\tilde{\ell})$, 
satisfy the variational conditions~\eqref{eql3.8a} (resp., 
conditions~\eqref{eql3.8b}$)$$;$ moreover, let the associated 
conditions~{\rm (i)}--{\rm (iv)} of item~$\pmb{(1)}$ (resp., 
item~$\pmb{(2)})$ of Lemma~\ref{lem3.8} be valid. For $n \! \in \! 
\mathbb{N}$ and $k \! \in \! \lbrace 1,2,\dotsc,K \rbrace$, let 
$\mathcal{X} \colon \mathbb{N} \times \lbrace 1,2,\dotsc,K \rbrace 
\times \overline{\mathbb{C}} \setminus \overline{\mathbb{R}} \! \to 
\! \mathrm{SL}_{2}(\mathbb{C})$ be the unique solution of the monic 
{\rm MPC ORF RHP} $(\mathcal{X}(z),\upsilon (z),\overline{\mathbb{R}})$ 
stated in Lemma~$\bm{\mathrm{RHP}_{\mathrm{MPC}}}$, where, in 
particular, for $n \! \in \! \mathbb{N}$ and $k \! \in \! \lbrace 1,2,
\dotsc,K \rbrace$ such that $\alpha_{p_{\mathfrak{s}}} \! := \! \alpha_{k} 
\! = \! \infty$ (resp., $\alpha_{p_{\mathfrak{s}}} \! := \! \alpha_{k} \! 
\neq \! \infty)$, it has the $z \! \to \! \alpha_{p_{\mathfrak{s}}} \! := \! 
\alpha_{k} \! = \! \infty$ (resp., $z \! \to \! \alpha_{p_{\mathfrak{s}}} 
\! := \! \alpha_{k} \! \neq \! \infty)$ Asymptotics~\eqref{eqlem5.5A} 
(resp.,  Asymptotics~\eqref{eqlem5.5C}$)$ stated in item~{\rm \pmb{(1)}} 
(resp., item~{\rm \pmb{(2)}}$)$ of Lemma~\ref{lem5.5}, with 
associated norming constant, $\mu_{n,\varkappa_{nk}}^{\infty}
(n,k)$ (resp., $\mu_{n,\varkappa_{nk}}^{f}(n,k))$, 
defined as in Subsection~\ref{subsubsec1.2.1} 
(resp., Subsection~\ref{subsubsec1.2.2}$)$. Then$:$ {\rm \pmb{(1)}} 
for $n \! \in \! \mathbb{N}$ and $k \! \in \! \lbrace 1,2,\dotsc,K \rbrace$ 
such that $\alpha_{p_{\mathfrak{s}}} \! := \! \alpha_{k} \! = \! \infty$,
\begin{equation} \label{eqlem5.6A} 
\mu_{n,\varkappa_{nk}}^{\infty}(n,k) \! = \! \lvert \lvert \pmb{\pi}_{k}^{n}
(\pmb{\cdot}) \rvert \rvert_{\mathscr{L}}^{-1} \underset{\underset{z_{o}
=1+o(1)}{\mathscr{N},n \to \infty}}{=} \sqrt{\dfrac{1}{-2 \pi \mi 
\me^{n \hat{\ell}} \left(\hat{A}_{0}^{\sharp}(\alpha_{k}) \! + \! 
\hat{\mathcal{R}}_{0}^{\hat{A}_{0}^{\sharp}}(\alpha_{k}) \right)_{12}}},
\end{equation}
where, in the double-scaling limit $\mathscr{N},n \! \to \! \infty$ 
such that $z_{o} \! = \! 1 \! + \! o(1)$, $\hat{A}_{0}^{\sharp}(\alpha_{k})$ 
and $\hat{\mathcal{R}}_{0}^{\hat{A}_{0}^{\sharp}}(\alpha_{k})$ are 
given in item~{\rm \pmb{(1)}} of Lemma~\ref{lem5.5}$;$ and 
{\rm \pmb{(2)}} for $n \! \in \! \mathbb{N}$ and $k \! \in \! \lbrace 
1,2,\dotsc,K \rbrace$ such that $\alpha_{p_{\mathfrak{s}}} \! := \! 
\alpha_{k} \! \neq \! \infty$,
\begin{equation} \label{eqlem5.6B} 
\mu_{n,\varkappa_{nk}}^{f}(n,k) \! = \! \lvert \lvert \pmb{\pi}_{k}^{n}
(\pmb{\cdot}) \rvert \rvert_{\mathscr{L}}^{-1} \underset{\underset{z_{o}
=1+o(1)}{\mathscr{N},n \to \infty}}{=} \sqrt{\dfrac{1}{2 \pi \mi 
\me^{n \tilde{\ell}} \left(A_{0}(\alpha_{k}) \! + \! \tilde{\mathcal{R}}_{0}^{A_{0}}
(\alpha_{k}) \right)_{12}}},
\end{equation}
where, in the double-scaling limit $\mathscr{N},n \! \to \! \infty$ 
such that $z_{o} \! = \! 1 \! + \! o(1)$, $A_{0}(\alpha_{k})$ and 
$\tilde{\mathcal{R}}_{0}^{A_{0}}(\alpha_{k})$ are given in 
item~{\rm \pmb{(2)}} of Lemma~\ref{lem5.5}.
\end{ccccc}

\emph{Proof.} The proof of this Lemma~\ref{lem5.6} consists of two 
cases: (i) $n \! \in \! \mathbb{N}$ and $k \! \in \! \lbrace 1,2,\dotsc,
K \rbrace$ such that $\alpha_{p_{\mathfrak{s}}} \! := \! \alpha_{k} \! 
= \! \infty$; and (ii) $n \! \in \! \mathbb{N}$ and $k \! \in \! \lbrace 
1,2,\dotsc,K \rbrace$ such that $\alpha_{p_{\mathfrak{s}}} \! := \! 
\alpha_{k} \! \neq \! \infty$. Notwithstanding the fact that the scheme 
of the proof is, \emph{mutatis mutandis}, similar for both cases, 
without loss of generality, only the proof for case~(ii) is presented in 
detail, whilst case~(i) is proved analogously.

For $n \! \in \! \mathbb{N}$ and $k \! \in \! \lbrace 1,2,\dotsc,K \rbrace$ 
such that $\alpha_{p_{\mathfrak{s}}} \! := \! \alpha_{k} \! \neq \! \infty$, 
recall {}from the corresponding representation~\eqref{intrepfin} for 
$\mathcal{X}(z)$ that, for $z \! \in \! \mathbb{C} \setminus \mathbb{R}$, 
$(\mathcal{X}(z))_{11} \! = \! (z \! - \! \alpha_{k}) \pmb{\pi}_{k}^{n}(z)$ 
and $(\mathcal{X}(z))_{12} \! = \! (z \! - \! \alpha_{k}) \int_{\mathbb{R}} 
\tfrac{(\mathcal{X}(\xi))_{11} \me^{-n \widetilde{V}(\xi)}}{(\xi - \alpha_{k})
(\xi -z)} \, \tfrac{\md \xi}{2 \pi \mi}$; via the latter formulae, the orthogonality 
conditions~\eqref{eq16}--\eqref{eq19}, the relation $\phi_{k}^{n}(z) 
\! = \! \mu_{n,\varkappa_{nk}}^{f}(n,k) \pmb{\pi}_{k}^{n}(z)$, where 
$\mu_{n,\varkappa_{nk}}^{f}(n,k)$ $(> \! 0)$ is the associated MPC 
ORF norming constant, and the representation of the monic MPC ORF given 
in the corresponding item of Remark~\ref{remcoeffs}, one proceeds as 
follows (e.g.,  $0 \! < \! \lvert z \! - \! \alpha_{k} \rvert \! < \! \min \lbrace 
\min_{i \neq j \in \lbrace 1,\dotsc,\mathfrak{s}-2,\mathfrak{s} \rbrace} 
\lbrace \lvert \alpha_{p_{i}} \! - \! \alpha_{p_{j}} \rvert \rbrace,
\min_{\underset{q \in \lbrace 1,\dotsc,\mathfrak{s}-2,\mathfrak{s} 
\rbrace}{m=1,2,\dotsc,N+1}} \lbrace \lvert \tilde{b}_{m-1} \! - \! \alpha_{p_{q}} 
\rvert,\lvert \tilde{a}_{m} \! - \! \alpha_{p_{q}} \rvert \rbrace \rbrace)$:
\begin{align*}
(\mathcal{X}(z))_{12} \underset{\mathbb{C} \setminus \mathbb{R} 
\ni z \to \alpha_{k}}{=}& \, (z \! - \! \alpha_{k}) \int_{\mathbb{R}} 
\dfrac{\pmb{\pi}_{k}^{n}(\xi)}{\xi \! - \! \alpha_{k}} \left(
1 \! + \! \dfrac{z \! - \! \alpha_{k}}{\xi \! - \! \alpha_{k}} 
\! + \! \dotsb \! + \! \dfrac{(z \! - \! \alpha_{k})^{\varkappa_{nk}-1}}{(\xi 
\! - \! \alpha_{k})^{\varkappa_{nk}-1}} \! + \! \dfrac{(z \! - \! \alpha_{k})^{
\varkappa_{nk}}}{(\xi \! - \! \alpha_{k})^{\varkappa_{nk}}} \! + \! 
\dotsb \right) \me^{-n \widetilde{V}(\xi)} \, \dfrac{\md \xi}{2 \pi \mi} \\
\underset{\mathbb{C} \setminus \mathbb{R} \ni z \to \alpha_{k}}{=}& 
\, (z \! - \! \alpha_{k})^{\varkappa_{nk}} \int_{\mathbb{R}} 
\dfrac{\pmb{\pi}_{k}^{n}(\xi) \me^{-n \widetilde{V}(\xi)}}{(\xi \! 
- \! \alpha_{k})^{\varkappa_{nk}}} \, \dfrac{\md \xi}{2 \pi \mi} 
\! + \! \mathcal{O} \left(\tilde{\mathfrak{c}}_{\mathcal{X}}^{\ast}
(n,k,z_{o})(z \! - \! \alpha_{k})^{\varkappa_{nk}+1} \right) \\
\underset{\mathbb{C} \setminus \mathbb{R} \ni z \to \alpha_{k}}{=}& 
\, (z \! - \! \alpha_{k})^{\varkappa_{nk}} \int_{\mathbb{R}} 
\pmb{\pi}_{k}^{n}(\xi) \left(\dfrac{1}{(\xi \! - \! \alpha_{k})^{
\varkappa_{nk}}} \! + \! \dfrac{1}{\mu^{f}_{n,\varkappa_{nk}}(n,k)} 
\sum_{m=1}^{\varkappa_{nk}-1} \dfrac{\mu^{f}_{n,m}(n,k)}{(\xi \! 
- \! \alpha_{k})^{m}} \! + \! \dfrac{1}{\mu^{f}_{n,\varkappa_{nk}}(n,k)} 
\sum_{l=1}^{\varkappa_{nk \tilde{k}_{\mathfrak{s}-1}}^{\infty}} 
\hat{\nu}^{f}_{n,l}(n,k) \xi^{l} \right. \\
+&\left. \, \dfrac{1}{\mu^{f}_{n,\varkappa_{nk}}(n,k)} \sum_{q=1}^{
\mathfrak{s}-2} \sum_{r=1}^{\varkappa_{nk \tilde{k}_{q}}} \dfrac{
\tilde{\nu}^{f}_{r,q}(n,k)}{(\xi \! - \! \alpha_{p_{q}})^{r}} \! + \! \dfrac{
\phi^{f}_{0}(n,k)}{\mu^{f}_{n,\varkappa_{nk}}(n,k)} \right) \me^{-n 
\widetilde{V}(\xi)} \, \dfrac{\md \xi}{2 \pi \mi} \! + \! \mathcal{O} 
\left(\tilde{\mathfrak{c}}_{\mathcal{X}}^{\ast}(n,k,z_{o})
(z \! - \! \alpha_{k})^{\varkappa_{nk}+1} \right) \\
\underset{\mathbb{C} \setminus \mathbb{R} \ni z \to \alpha_{k}}{=}& 
\, \dfrac{(z \! - \! \alpha_{k})^{\varkappa_{nk}}}{2 \pi \mi} \int_{\mathbb{R}} 
\pmb{\pi}_{k}^{n}(\xi) \pmb{\pi}_{k}^{n}(\xi) \me^{-n \widetilde{V}(\xi)} 
\, \md \xi \! + \! \mathcal{O} \left(\tilde{\mathfrak{c}}_{\mathcal{X}}^{
\ast}(n,k,z_{o})(z \! - \! \alpha_{k})^{\varkappa_{nk}+1} \right) \\
\underset{\mathbb{C} \setminus \mathbb{R} \ni z \to \alpha_{k}}{=}& 
\, \dfrac{(z \! - \! \alpha_{k})^{\varkappa_{nk}}}{2 \pi \mi (\mu^{f}_{n,
\varkappa_{nk}}(n,k))^{2}} \underbrace{\int_{\mathbb{R}} \phi_{k}^{n}
(\xi) \phi_{k}^{n}(\xi) \me^{-n \widetilde{V}(\xi)} \, \md \xi}_{= \, 1} 
\! + \mathcal{O} \left(\tilde{\mathfrak{c}}_{\mathcal{X}}^{\ast}(n,k,z_{o})
(z \! - \! \alpha_{k})^{\varkappa_{nk}+1} \right) \quad \Rightarrow
\end{align*}
\begin{equation} \label{eqlem5.6C} 
\left(\mathcal{X}(z)(z \! - \! \alpha_{k})^{(\varkappa_{nk}-1) \sigma_{3}} 
\right)_{12} \underset{\mathbb{C} \setminus \mathbb{R} \ni z \to 
\alpha_{k}}{=} \left(\dfrac{1}{2 \pi \mi (\mu^{f}_{n,\varkappa_{nk}}
(n,k))^{2}} \right)(z \! - \! \alpha_{k}) \! + \! \mathcal{O} \left(
\tilde{\mathfrak{c}}_{\mathcal{X}}^{\ast}(n,k,z_{o})(z \! - \! 
\alpha_{k})^{2} \right),
\end{equation}
where $\tilde{\mathfrak{c}}_{\mathcal{X}}^{\ast}(n,k,z_{o}) \! =_{
\underset{z_{o}=1+o(1)}{\mathscr{N},n \to \infty}} \! \mathcal{O}(1)$. 
For $n \! \in \! \mathbb{N}$ and $k \! \in \! \lbrace 1,2,\dotsc,K 
\rbrace$ such that $\alpha_{p_{\mathfrak{s}}} \! := \! \alpha_{k} \! 
\neq \! \infty$, it follows {}from the Asymptotics~\eqref{eqlem5.5C} 
that
\begin{equation} \label{eqlem5.6D} 
\left(\mathcal{X}(z)(z \! - \! \alpha_{k})^{(\varkappa_{nk}-1) \sigma_{3}} 
\right)_{12} \underset{\mathbb{C} \setminus \mathbb{R} \ni z \to 
\alpha_{k}}{=} \me^{n \tilde{\ell}} \left(A_{0}(\alpha_{k}) \! + \! 
\tilde{\mathcal{R}}_{0}^{A_{0}}(\alpha_{k}) \right)_{12}(z \! - \! \alpha_{k}) 
\! + \! \mathcal{O} \left(\tilde{\mathfrak{c}}_{\mathcal{X}}^{\ast^{\prime}}
(n,k,z_{o})(z \! - \! \alpha_{k})^{2} \right),
\end{equation}
where, in the double-scaling limit $\mathscr{N},n \! \to \! \infty$ such 
that $z_{o} \! = \! 1 \! + \! o(1)$, $A_{0}(\alpha_{k})$ and $\tilde{
\mathcal{R}}_{0}^{A_{0}}(\alpha_{k})$ are given in item~{\rm \pmb{(2)}} 
of Lemma~\ref{lem5.5}, and $\tilde{\mathfrak{c}}_{\mathcal{X}}^{
\ast^{\prime}}(n,k,z_{o}) \! =_{\underset{z_{o}=1+o(1)}{\mathscr{N},
n \to \infty}} \! \mathcal{O}(1)$; hence, via the asymptotic 
expansions~\eqref{eqlem5.6C} and~\eqref{eqlem5.6D}, one arrives 
at, for $n \! \in \! \mathbb{N}$ and $k \! \in \! \lbrace 1,2,\dotsc,
K \rbrace$ such that $\alpha_{p_{\mathfrak{s}}} \! := \! \alpha_{k} 
\! \neq \! \infty$, Equation~\eqref{eqlem5.6B}.

The analysis for the case $n \! \in \! \mathbb{N}$ and $k \! \in \! 
\lbrace 1,2,\dotsc,K \rbrace$ such that $\alpha_{p_{\mathfrak{s}}} 
\! := \! \alpha_{k} \! = \! \infty$ is, \emph{mutatis mutandis}, 
analogous, and leads to Equation~\eqref{eqlem5.6A}. \hfill $\qed$
\begin{ccccc} \label{lemmpainffin} 
Let the external field $\widetilde{V} \colon \overline{\mathbb{R}} 
\setminus \lbrace \alpha_{1},\alpha_{2},\dotsc,\alpha_{K} \rbrace \! 
\to \! \mathbb{R}$ satisfy conditions~\eqref{eq20}--\eqref{eq22} and 
be regular. For $n \! \in \! \mathbb{N}$ and $k \! \in \! \lbrace 1,2,
\dotsc,K \rbrace$ such that $\alpha_{p_{\mathfrak{s}}} \! := \! 
\alpha_{k} \! = \! \infty$ (resp., $\alpha_{p_{\mathfrak{s}}} \! := \! 
\alpha_{k} \! \neq \! \infty)$, let the associated equilibrium measure, 
$\mu_{\widetilde{V}}^{\infty}$ (resp., $\mu_{\widetilde{V}}^{f})$, and its 
support, $J_{\infty}$ (resp., $J_{f})$, be as described in item~$\pmb{(1)}$ 
(resp., item~$\pmb{(2)})$ of Lemma~\ref{lem3.7}, and, along with the 
corresponding variational constant, $\hat{\ell}$ (resp., $\tilde{\ell})$, 
satisfy the variational conditions~\eqref{eql3.8a} (resp., 
conditions~\eqref{eql3.8b}$)$$;$ moreover, let the associated 
conditions~{\rm (i)}--{\rm (iv)} of item~$\pmb{(1)}$ (resp., 
item~$\pmb{(2)})$ of Lemma~\ref{lem3.8} be valid. For $n \! \in \! 
\mathbb{N}$ and $k \! \in \! \lbrace 1,2,\dotsc,K \rbrace$, let 
$\mathcal{X} \colon \mathbb{N} \times \lbrace 1,2,\dotsc,K \rbrace 
\times \overline{\mathbb{C}} \setminus \overline{\mathbb{R}} \! \to 
\! \mathrm{SL}_{2}(\mathbb{C})$ be the unique solution of the monic 
{\rm MPC ORF RHP} $(\mathcal{X}(z),\upsilon (z),\overline{\mathbb{R}})$ 
stated in Lemma~$\bm{\mathrm{RHP}_{\mathrm{MPC}}}$. Define the 
Markov-Stieltjes transform of the probability measure $\widetilde{\mu}$ 
as follows:\footnote{This is Equation~\eqref{mvssinf1} transformed 
according to (cf. Remark~\ref{rem1.3.2}) $\md \mu \! \to \! \md 
\widetilde{\mu}$.}
\begin{equation} \label{eqmvsstildemu} 
\mathrm{F}_{\tilde{\mu}}(z) \! := \! \int_{\mathbb{R}}(z \! - \! \xi)^{-1} 
\me^{-n \widetilde{V}(\xi)} \, \md \xi.
\end{equation}
For $n \! \in \! \mathbb{N}$ and $k \! \in \! \lbrace 1,2,\dotsc,K 
\rbrace$ such that $\alpha_{p_{\mathfrak{s}}} \! := \! \alpha_{k} \! = \! 
\infty$, define the corresponding associated $\mathrm{R}$-function 
as {}\footnote{This is Equation~\eqref{mvssinf4} transformed 
according to (cf. Remark~\ref{rem1.3.2}) $\md \mu \! \to \! \md 
\widetilde{\mu}$.}
\begin{equation} \label{eqlemmvssinfmpa1} 
\widehat{\pmb{\mathrm{R}}}_{\tilde{\mu}}(z) \! := \! \int_{\mathbb{R}} 
\left(\dfrac{\pmb{\pi}^{n}_{k}(\xi) \! - \! \pmb{\pi}^{n}_{k}(z)}{\xi 
\! - \! z} \right) \me^{-n \widetilde{V}(\xi)} \, \md \xi,
\end{equation}
and, for $n \! \in \! \mathbb{N}$ and $k \! \in \! \lbrace 1,2,\dotsc,K 
\rbrace$ such that $\alpha_{p_{\mathfrak{s}}} \! := \! \alpha_{k} \! \neq 
\! \infty$, define the corresponding associated $\mathrm{R}$-function 
as {}\footnote{This is Equation~\eqref{mvssfin3} transformed 
according to (cf. Remark~\ref{rem1.3.2}) $\md \mu \! \to \! \md 
\widetilde{\mu}$.}
\begin{equation} \label{eqlemmvssfinmpa1} 
\widetilde{\pmb{\mathrm{R}}}_{\tilde{\mu}}(z) \! := \! \int_{\mathbb{R}} 
\left(\dfrac{\pmb{\pi}^{n}_{k}(\xi) \! - \! \pmb{\pi}^{n}_{k}(z)}{\xi \! 
- \! z} \right) \me^{-n \widetilde{V}(\xi)} \, \md \xi.
\end{equation}
Then$:$ {\rm \pmb{(1)}} for $n \! \in \! \mathbb{N}$ and $k \! \in \! 
\lbrace 1,2,\dotsc,K \rbrace$ such that $\alpha_{p_{\mathfrak{s}}} 
\! := \! \alpha_{k} \! = \! \infty$, in the double-scaling limit 
$\mathscr{N},n \! \to \! \infty$ such that $z_{o} \! = \! 1 \! + \! o(1)$,
\begin{equation} \label{eqlemmvssinfmpa2} 
\dfrac{\widehat{\pmb{\mathrm{R}}}_{\tilde{\mu}}(z)}{\pmb{\pi}^{n}_{k}
(z)} \! = \! \dfrac{\widehat{\mathrm{U}}_{\tilde{\mu}}(z)}{\widehat{
\mathrm{V}}_{\tilde{\mu}}(z)} \! := \! \dfrac{\sum_{j=0}^{(n-1)K+k-1} 
\hat{r}_{j}z^{j}}{\sum_{j=0}^{(n-1)K+k} \hat{t}_{j}z^{j}}
\end{equation}
is the {\rm MPA} of type $((n \! - \! 1)K \! + \! k \! - \! 1,(n \! - \! 1)K \! + 
\! k)$ for the Markov-Stieltjes transform, where, with $\hat{t}_{(n-1)K+k} 
\! = \! \hat{t}_{(n-1)K+k}(n,k,z_{o}) \! = \! 1$, $\hat{r}_{j} \! = \! \hat{r}_{j}
(n,k,z_{o})$ and $\hat{t}_{j} \! = \! \hat{t}_{j}(n,k,z_{o})$, $j \! = \! 0,1,
\dotsc,(n \! - \! 1)K \! + \! k \! - \! 1$, solve the linear inhomogeneous 
algebraic system of equations~\eqref{wm24frevhmpa20} below, and 
the corresponding {\rm MPA} error term is given by {}\footnote{This 
is Equations~\eqref{mvssinf9} and~\eqref{mvssinf10} transformed 
according to (cf. Remark~\ref{rem1.3.2}) $\md \mu \! \to \! \md 
\widetilde{\mu}$.}
\begin{equation} \label{eqlemmvssinfmpa3} 
\widehat{\pmb{\mathrm{E}}}_{\tilde{\mu}}(z) \! := \! 
\dfrac{\widehat{\pmb{\mathrm{R}}}_{\tilde{\mu}}(z)}{\pmb{\pi}^{n}_{k}
(z)} \! - \! \mathrm{F}_{\tilde{\mu}}(z) \! = \! 2 \pi \mi 
\dfrac{(\mathcal{X}(z))_{12}}{(\mathcal{X}(z))_{11}},
\end{equation}
where $(\mathcal{X}(z))_{11}$ and $(\mathcal{X}(z))_{12}$ are the 
$(1 \, 1)$- and $(1 \, 2)$-elements, respectively, of the matrix 
representation~\eqref{intrepinf}$;$ and {\rm \pmb{(2)}} for $n \! \in 
\! \mathbb{N}$ and $k \! \in \! \lbrace 1,2,\dotsc,K \rbrace$ such 
that $\alpha_{p_{\mathfrak{s}}} \! := \! \alpha_{k} \! \neq \! \infty$, 
in the double-scaling limit $\mathscr{N},n \! \to \! \infty$ such 
that $z_{o} \! = \! 1 \! + \! o(1)$,
\begin{equation} \label{eqlemmvssfinmpa2} 
\dfrac{\widetilde{\pmb{\mathrm{R}}}_{\tilde{\mu}}(z)}{\pmb{\pi}^{n}_{k}
(z)} \! = \! \dfrac{\widetilde{\mathrm{U}}_{\tilde{\mu}}(z)}{\widetilde{
\mathrm{V}}_{\tilde{\mu}}(z)} \! := \! \dfrac{\sum_{j=0}^{(n-1)K+k-1} 
\tilde{r}_{j}z^{j}}{\sum_{j=0}^{(n-1)K+k} \tilde{t}_{j}z^{j}}
\end{equation}
is the {\rm MPA} of type $((n \! - \! 1)K \! + \! k \! - \! 1,(n \! - \! 1)K \! + \! 
k)$ for the Markov-Stieltjes transform, where $\tilde{r}_{j} \! = \! \tilde{r}_{j}
(n,k,z_{o})$ and $\tilde{t}_{j} \! = \! \tilde{t}_{j}(n,k,z_{o})$, $j \! = \! 0,1,
\dotsc,(n \! - \! 1)K \! + \! k \! - \! 1$, solve the linear inhomogeneous 
algebraic system of equations~\eqref{wm24frevhmpa12} below,
\begin{align} \label{eqlemmvssfinmpa3}
\tilde{t}_{(n-1)K+k} \! = \! \tilde{t}_{(n-1)K+k}(n,k,z_{o}) 
\! = \! \widetilde{\nu}_{\mathfrak{s}-1,\varkappa_{nk \tilde{k}_{
\mathfrak{s}-1}}^{\infty}}^{\raise-1.0ex\hbox{$\scriptstyle f$}}(n,k) 
\underset{\underset{z_{o}=1+o(1)}{\mathscr{N},n \to \infty}}{=}& \, 
\left(\tilde{t}_{(n-1)K+k} \right)_{0} \! + \! \dfrac{1}{(n \! - \! 1)K 
\! + \! k} \left(\tilde{t}_{(n-1)K+k} \right)_{1} \nonumber \\
+& \, \mathcal{O} \left(\dfrac{\tilde{\mathfrak{c}}_{\tilde{t}}
(n,k,z_{o})}{((n \! - \! 1)K \! + \! k)^{2}} \right),
\end{align}
where
\begin{align}
(\tilde{t}_{(n-1)K+k})_{0} =& \, \mathscr{E}^{-1} \tilde{\mathbb{K}}_{11} 
\dfrac{\tilde{\boldsymbol{\theta}}(\tilde{\boldsymbol{u}}_{+}(\infty) \! 
- \! \frac{1}{2 \pi}((n \! - \! 1)K \! + \! k) \tilde{\boldsymbol{\Omega}} 
\! + \! \tilde{\boldsymbol{d}})}{\tilde{\boldsymbol{\theta}}(\tilde{\boldsymbol{
u}}_{+}(\infty) \! + \! \tilde{\boldsymbol{d}})} \me^{\tilde{\Xi}_{0}^{\sharp}
(\alpha_{p_{\mathfrak{s}-1}})}(-1)^{\sum_{j \in \hat{\Delta}_{f}(k)} 
\varkappa_{nk \tilde{k}_{j}}}, \label{eqlemmvssfinmpa4} \\
(\tilde{t}_{(n-1)K+k})_{1} =& \, \dfrac{\tilde{\boldsymbol{\theta}}
(\tilde{\boldsymbol{u}}_{+}(\infty) \! - \! \frac{1}{2 \pi}((n \! - \! 1)K \! 
+ \! k) \tilde{\boldsymbol{\Omega}} \! + \! \tilde{\boldsymbol{d}})}{
\tilde{\boldsymbol{\theta}}(\tilde{\boldsymbol{u}}_{+}(\infty) \! + \! 
\tilde{\boldsymbol{d}})} \left(\mathscr{E}^{-1} \tilde{\mathbb{K}}_{11} 
\left(\tilde{\mathcal{R}}_{0,\Delta}^{\tilde{A}_{0}^{\sharp}}
(\alpha_{p_{\mathfrak{s}-1}}) \right)_{11} \! + \! \mathscr{E} 
\tilde{\mathbb{K}}_{21} \left(\tilde{\mathcal{R}}_{0,\Delta}^{
\tilde{A}_{0}^{\sharp}}(\alpha_{p_{\mathfrak{s}-1}}) \right)_{12} 
\right) \nonumber \\
\times& \, \me^{\tilde{\Xi}_{0}^{\sharp}(\alpha_{p_{\mathfrak{s}-1}})}
(-1)^{\sum_{j \in \hat{\Delta}_{f}(k)} \varkappa_{nk \tilde{k}_{j}}}, 
\label{eqlemmvssfinmpa5}
\end{align}
with $\mathscr{E}$ defined by Equation~\eqref{eqmainfin13}, 
$\tilde{\mathbb{K}}$ defined by Equations~\eqref{eqmainfin8} 
and~\eqref{eqconsfinn1}, $\tilde{\Xi}_{0}^{\sharp}(\alpha_{
p_{\mathfrak{s}-1}})$ defined by Equation~\eqref{eqfivecapxiefin1}, 
$\hat{\Delta}_{f}(k) \! := \! \lbrace \mathstrut j \! \in \! \lbrace 
1,2,\dotsc,\mathfrak{s} \! - \! 2 \rbrace; \, \alpha_{p_{j}} \! > \! 
\alpha_{k} \rbrace$, and
\begin{align}
\tilde{\mathcal{R}}_{0,\Delta}^{\tilde{A}_{0}^{\sharp}}
(\alpha_{p_{\mathfrak{s}-1}}) :=& \, \sum_{j=1}^{N+1} \left(
\dfrac{(\tilde{\alpha}_{0}(\tilde{b}_{j-1}))^{-2}}{\tilde{b}_{j-1} 
\! - \! \alpha_{k}} \left(\tilde{\alpha}_{0}(\tilde{b}_{j-1}) 
\tilde{\boldsymbol{\mathrm{B}}}(\tilde{b}_{j-1}) \! - \! \tilde{\alpha}_{1}
(\tilde{b}_{j-1}) \tilde{\boldsymbol{\mathrm{A}}}(\tilde{b}_{j-1}) \right) 
\! - \! \dfrac{(\tilde{\alpha}_{0}(\tilde{b}_{j-1}))^{-1}}{(\tilde{b}_{j-1} 
\! - \! \alpha_{k})^{2}} \tilde{\boldsymbol{\mathrm{A}}}(\tilde{b}_{j-1}) 
\right. \nonumber \\
+&\left. \, \dfrac{(\tilde{\alpha}_{0}(\tilde{a}_{j}))^{-2}}{\tilde{a}_{j} \! - 
\! \alpha_{k}} \left(\tilde{\alpha}_{0}(\tilde{a}_{j}) \tilde{\boldsymbol{
\mathrm{B}}}(\tilde{a}_{j}) \! - \! \tilde{\alpha}_{1}(\tilde{a}_{j}) \tilde{
\boldsymbol{\mathrm{A}}}(\tilde{a}_{j}) \right) \! - \! \dfrac{(\tilde{
\alpha}_{0}(\tilde{a}_{j}))^{-1}}{(\tilde{a}_{j} \! - \! \alpha_{k})^{2}} 
\tilde{\boldsymbol{\mathrm{A}}}(\tilde{a}_{j}) \right), 
\label{eqlemmvssfinmpa6}
\end{align}
where $\tilde{\boldsymbol{\mathrm{A}}}(\tilde{b}_{j-1})$, $\tilde{
\boldsymbol{\mathrm{A}}}(\tilde{a}_{j})$, $\tilde{\boldsymbol{
\mathrm{B}}}(\tilde{b}_{j-1})$, $\tilde{\boldsymbol{\mathrm{B}}}
(\tilde{a}_{j})$, $\tilde{\alpha}_{0}(\tilde{b}_{j-1})$, $\tilde{\alpha}_{0}
(\tilde{a}_{j})$, $\tilde{\alpha}_{1}(\tilde{b}_{j-1})$, and $\tilde{\alpha}_{1}
(\tilde{a}_{j})$, $j \! \in \! \lbrace 1,2,\dotsc,N \! + \! 1 \rbrace$, are 
defined in item~{\rm \pmb{(2)}} of Proposition~\ref{propo5.1}, and 
$\tilde{\mathfrak{c}}_{\tilde{t}}(n,k,z_{o}) \! =_{\underset{z_{o}
=1+o(1)}{\mathscr{N},n \to \infty}} \! \mathcal{O}(1)$, and the 
corresponding {\rm MPA} error term is given by {}\footnote{This 
is Equations~\eqref{mvssfin9} and~\eqref{mvssfin10} transformed 
according to (cf. Remark~\ref{rem1.3.2}) $\md \mu \! \to \! \md 
\widetilde{\mu}$.}
\begin{equation} \label{eqlemmvssfinmpa7} 
\widetilde{\pmb{\mathrm{E}}}_{\tilde{\mu}}(z) \! := \! 
\dfrac{\widetilde{\pmb{\mathrm{R}}}_{\tilde{\mu}}(z)}{\pmb{\pi}^{n}_{k}
(z)} \! - \! \mathrm{F}_{\tilde{\mu}}(z) \! = \! 2 \pi \mi 
\dfrac{(\mathcal{X}(z))_{12}}{(\mathcal{X}(z))_{11}},
\end{equation}
where $(\mathcal{X}(z))_{11}$ and $(\mathcal{X}(z))_{12}$ are the 
$(1 \, 1)$- and $(1 \, 2)$-elements, respectively, of the matrix 
representation~\eqref{intrepfin}.
\end{ccccc}

\emph{Proof.} The proof of this Lemma~\ref{lemmpainffin} consists 
of two cases: (i) $n \! \in \! \mathbb{N}$ and $k \! \in \! \lbrace 
1,2,\dotsc,K \rbrace$ such that $\alpha_{p_{\mathfrak{s}}} \! := \! 
\alpha_{k} \! = \! \infty$; and (ii) $n \! \in \! \mathbb{N}$ and $k \! 
\in \! \lbrace 1,2,\dotsc,K \rbrace$ such that $\alpha_{p_{\mathfrak{s}}} 
\! := \! \alpha_{k} \! \neq \! \infty$. Notwithstanding the fact that the 
scheme of the proof is, \emph{mutatis mutandis}, similar for both 
cases, without loss of generality, only the proof for case~(ii) is 
presented in detail, whilst case~(i) is proved analogously.

\pmb{(1)} For $n \! \in \! \mathbb{N}$ and $k \! \in \! \lbrace 1,2,
\dotsc,K \rbrace$ such that $\alpha_{p_{\mathfrak{s}}} \! := \! \alpha_{k} 
\! \neq \! \infty$, define the associated $\mathrm{R}$-function as 
in Equation~\eqref{eqlemmvssfinmpa1}: via the representation of 
the monic MPC ORF given in the corresponding item of 
Remark~\ref{remcoeffs},\footnote{Recall that $\widetilde{\nu}_{
\mathfrak{s},\varkappa_{nk}}^{\raise-1.0ex\hbox{$\scriptstyle f$}}
(n,k) \! = \! 1$.} the fact that $\widetilde{\mu} \! \in \! \mathscr{M}_{1}
(\mathbb{R})$, the identity $y_{1}^{m} \! - \! y_{2}^{m} \! = \! (y_{1} \! 
- \! y_{2})(y_{1}^{m-1} \! + \! y_{1}^{m-2}y_{2} \! + \! \dotsb \! + \! 
y_{1}y_{2}^{m-2} \! + \! y_{2}^{m-1})$, and the moment integrals 
$c^{(q)}_{j}(\alpha_{p_{q}}) \! := \! -\int_{\mathbb{R}}(\xi \! - \! 
\alpha_{p_{q}})^{-(1+j)} \, \md \widetilde{\mu}(\xi)$, $(j,q) \! \in 
\! \mathbb{N}_{0} \times \lbrace 1,\dotsc,\mathfrak{s} \! - \! 2,
\mathfrak{s} \rbrace$, and $c^{(\infty)}_{i} \! := \! \int_{\mathbb{R}} 
\xi^{i-1} \, \md \widetilde{\mu}(\xi)$, $i \! \in \! \mathbb{N}$, 
with $c^{(\infty)}_{1} \! = \! 1$, one shows that 
\begin{align*}
\int_{\mathbb{R}} \left(\dfrac{\pmb{\pi}^{n}_{k}(\xi) \! - \! 
\pmb{\pi}^{n}_{k}(z)}{\xi \! - \! z} \right) \md \widetilde{\mu}
(\xi) =& \, \sum_{q=1}^{\mathfrak{s}-2} \sum_{m=1}^{
\varkappa_{nk \tilde{k}_{q}}} \sum_{j=1}^{m} 
\dfrac{\widetilde{\nu}_{q,m}^{\raise-1.0ex\hbox{$\scriptstyle f$}}
(n,k)c^{(q)}_{m-j}(\alpha_{p_{q}})}{(z \! - \! \alpha_{p_{q}})^{j}} \! 
+ \! \sum_{m=1}^{\varkappa_{nk \tilde{k}_{\mathfrak{s}-1}}^{
\infty}} \sum_{j=0}^{m-1} 
\widetilde{\nu}_{\mathfrak{s}-1,m}^{\raise-1.0ex\hbox{$\scriptstyle f$}}
(n,k)c^{(\infty)}_{m-j}z^{j} \\
+& \, \sum_{m=1}^{\varkappa_{nk}} \sum_{j=1}^{m} 
\dfrac{\widetilde{\nu}_{\mathfrak{s},m}^{\raise-1.0ex\hbox{$\scriptstyle f$}}
(n,k)c^{(\mathfrak{s})}_{m-j}(\alpha_{k})}{(z \! - \! \alpha_{k})^{j}},
\end{align*}
whence, via Equation~\eqref{fincount}, one arrives at (the improper 
fraction)
\begin{equation} \label{wm24frevhmpa1} 
\widetilde{\pmb{\mathrm{R}}}_{\tilde{\mu}}(z) \! = \! 
\dfrac{\widetilde{\mathrm{U}}_{\tilde{\mu}}(z)}{\prod_{q=1}^{
\mathfrak{s}-2}(z \! - \! \alpha_{p_{q}})^{\varkappa_{nk \tilde{k}_{q}}}
(z \! - \! \alpha_{k})^{\varkappa_{nk}}} \! := \! \dfrac{\sum_{j=0}^{
(n-1)K+k-1} \tilde{r}_{j}z^{j}}{\prod_{q=1}^{\mathfrak{s}-2}(z \! - \! 
\alpha_{p_{q}})^{\varkappa_{nk \tilde{k}_{q}}}(z \! - \! \alpha_{k})^{
\varkappa_{nk}}},
\end{equation}
where $\tilde{r}_{j} \! = \! \tilde{r}_{j}(n,k,z_{o})$, $j \! = \! 0,1,\dotsc,
(n \! - \! 1)K \! + \! k \! - \! 1$.\footnote{The $z \! \to \! 
\alpha_{p_{\mathfrak{s}-1}} \! = \! \infty$ asymptotic analysis of 
Equation~\eqref{wm24frevhmpa1} shows that $\tilde{r}_{(n-1)K+k-1} 
\! = \! \widetilde{\nu}_{\mathfrak{s}-1,\varkappa_{nk \tilde{k}_{
\mathfrak{s}-1}}^{\infty}}^{\raise-1.0ex\hbox{$\scriptstyle f$}}(n,k)$.} For 
$n \! \in \! \mathbb{N}$ and $k \! \in \! \lbrace 1,2,\dotsc,K \rbrace$ such 
that $\alpha_{p_{\mathfrak{s}}} \! := \! \alpha_{k} \! \neq \! \infty$, a 
calculation based on Equation~\eqref{fincount}, and the representation of the 
monic MPC ORF given in the corresponding item of Remark~\ref{remcoeffs}, 
shows that $\pmb{\pi}^{n}_{k}(z)$ can be presented as (the improper fraction)
\begin{equation} \label{wm24frevhmpa2} 
\pmb{\pi}^{n}_{k}(z) \! = \! \dfrac{\widetilde{\mathrm{V}}_{\tilde{\mu}}
(z)}{\prod_{q=1}^{\mathfrak{s}-2}(z \! - \! \alpha_{p_{q}})^{\varkappa_{nk 
\tilde{k}_{q}}}(z \! - \! \alpha_{k})^{\varkappa_{nk}}} \! := \! \dfrac{
\sum_{j=0}^{(n-1)K+k} \tilde{t}_{j}z^{j}}{\prod_{q=1}^{\mathfrak{s}-2}
(z \! - \! \alpha_{p_{q}})^{\varkappa_{nk \tilde{k}_{q}}}(z \! - \! 
\alpha_{k})^{\varkappa_{nk}}},
\end{equation}
where $\tilde{t}_{j} \! = \! \tilde{t}_{j}(n,k,z_{o})$, $j \! = \! 0,1,
\dotsc,(n \! - \! 1)K \! + \! k$, with $\tilde{t}_{(n-1)K+k} 
\! = \! \widetilde{\nu}_{\mathfrak{s}-1,\varkappa_{nk 
\tilde{k}_{\mathfrak{s}-1}}^{\infty}}^{\raise-1.0ex\hbox{$\scriptstyle f$}}
(n,k)$. For $n \! \in \! \mathbb{N}$ and $k \! \in \! \lbrace 1,2,
\dotsc,K \rbrace$ such that $\alpha_{p_{\mathfrak{s}}} \! := \! \alpha_{k} 
\! \neq \! \infty$, the $z \! \to \! \alpha_{p_{\mathfrak{s}-1}} \! = \! \infty$ 
asymptotic analysis of Equation~\eqref{wm24frevhmpa2} shows that 
$\pmb{\pi}^{n}_{k}(z) \! =_{z \to \alpha_{p_{\mathfrak{s}-1}} = \infty} \! 
\tilde{t}_{(n-1)K+k}z^{\varkappa_{nk \tilde{k}_{\mathfrak{s}-1}}^{\infty}}
(1 \! + \! \mathcal{O}(z^{-1}))$, whence, via the $(1 \, 1)$-element of 
the matrix representation~\eqref{intrepfin}, $(\mathcal{X}(z))_{11}
z^{-(\varkappa_{nk \tilde{k}_{\mathfrak{s}-1}}^{\infty}+1)} \linebreak[4] 
\! =_{z \to \alpha_{p_{\mathfrak{s}-1}} = \infty} \! \tilde{t}_{(n-1)K+k}
(1 \! + \! \mathcal{O}(z^{-1}))$; hence, {}from the relevant $(1 \, 1)$-element 
of the asymptotics~\eqref{eqlem5.5E}, one shows that, in the double-scaling 
limit $\mathscr{N},n \! \to \! \infty$ such that $z_{o} \! = \! 1 \! + \! o(1)$, 
$\tilde{t}_{(n-1)K+k}$ has the asymptotics~\eqref{eqlemmvssfinmpa3}. For 
$n \! \in \! \mathbb{N}$ and $k \! \in \! \lbrace 1,2,\dotsc,K \rbrace$ such 
that $\alpha_{p_{\mathfrak{s}}} \! := \! \alpha_{k} \! \neq \! \infty$, a 
straightforward calculation based on Equations~\eqref{wm24frevhmpa1} 
and~\eqref{wm24frevhmpa2} (see, also, \cite{n2,n1}) shows that 
the proper fraction given by Equation~\eqref{eqlemmvssfinmpa2} 
is the corresponding MPA of type $((n \! - \! 1)K \! + \! k \! - 
\! 1,(n \! - \! 1)K \! + \! k)$ for the Markov-Stieltjes transform 
defined by Equation~\eqref{eqmvsstildemu}, in the sense that the 
interpolation conditions (cf. Subsection~\ref{subsubsec1.2.2}, 
Equations~\eqref{mvssfin6}--\eqref{mvssfin8})
\begin{gather}
\dfrac{\widetilde{\mathrm{U}}_{\tilde{\mu}}(z)}{\widetilde{\mathrm{V}}_{
\tilde{\mu}}(z)} \! - \! \sum_{j=0}^{2 \varkappa_{nk \tilde{k}_{q}}-1}
c^{(q)}_{j}(\alpha_{p_{q}})(z \! - \! \alpha_{p_{q}})^{j} \underset{z \to 
\alpha_{p_{q}}}{=} \mathcal{O} \left(\tilde{\mathfrak{c}}^{(q)}_{\tilde{\mu}}
(n,k,z_{o})(z \! - \! \alpha_{p_{q}})^{2 \varkappa_{nk \tilde{k}_{q}}} \right), 
\quad q \! = \! 1,2,\dotsc,\mathfrak{s} \! - \! 2, \label{wm24frevhmpa3} \\
\dfrac{\widetilde{\mathrm{U}}_{\tilde{\mu}}(z)}{\widetilde{\mathrm{V}}_{
\tilde{\mu}}(z)} \! - \! \sum_{j=0}^{2 \varkappa_{nk}-1}c^{(\mathfrak{s})}_{j}
(\alpha_{k})(z \! - \! \alpha_{k})^{j} \underset{z \to \alpha_{k}}{=} 
\mathcal{O} \left(\tilde{\mathfrak{c}}^{(\mathfrak{s})}_{\tilde{\mu}}
(n,k,z_{o})(z \! - \! \alpha_{k})^{2 \varkappa_{nk}} \right), 
\label{wm24frevhmpa4} \\
\dfrac{\widetilde{\mathrm{U}}_{\tilde{\mu}}(z)}{\widetilde{\mathrm{V}}_{
\tilde{\mu}}(z)} \! - \! \sum_{j=1}^{2 \varkappa^{\infty}_{nk \tilde{k}_{
\mathfrak{s}-1}}}c^{(\infty)}_{j}z^{-j} \underset{z \to \alpha_{
p_{\mathfrak{s}-1}} = \infty}{=} \mathcal{O} \left(
\tilde{\mathfrak{c}}^{(\mathfrak{s}-1)}_{\tilde{\mu}}(n,k,z_{o})
z^{-(2 \varkappa^{\infty}_{nk \tilde{k}_{\mathfrak{s}-1}}+1)} \right), 
\label{wm24frevhmpa5}
\end{gather}
where $\tilde{\mathfrak{c}}^{(q)}_{\tilde{\mu}}(n,k,z_{o}) \! =_{
\underset{z_{o}=1+o(1)}{\mathscr{N},n \to \infty}} \! \mathcal{O}
(1)$, $q \! \in \! \lbrace 1,2,\dotsc,\mathfrak{s} \rbrace$, are satisfied. 
Note that Equations~\eqref{wm24frevhmpa3}, \eqref{wm24frevhmpa4}, 
and~\eqref{wm24frevhmpa5}, respectively, give rise to 
$2 \sum_{q=1}^{\mathfrak{s}-2} \varkappa_{nk \tilde{k}_{q}}$, 
$2 \varkappa_{nk}$, and $2 \varkappa_{nk \tilde{k}_{
\mathfrak{s}-1}}^{\infty}$ interpolation conditions, for a combined 
total of (cf. Equation~\eqref{fincount}) $2((n \! - \! 1)K \! + \! k)$ 
conditions, which is precisely the number necessary in order to 
determine---uniquely---the $2((n \! - \! 1)K \! + \! k)$ coefficients 
$\tilde{r}_{0},\tilde{r}_{1},\dotsc,\tilde{r}_{(n-1)K+k-1},\tilde{t}_{0},
\tilde{t}_{1},\dotsc,\tilde{t}_{(n-1)K+k-1}$ (recall that, in the 
double-scaling limit $\mathscr{N},n \! \to \! \infty$ such that 
$z_{o} \! = \! 1 \! + \! o(1)$, $\tilde{t}_{(n-1)K+k}$ is given by 
Equation~\eqref{eqlemmvssfinmpa3}). Substituting the expressions 
for the polynomials $\widetilde{\mathrm{U}}_{\tilde{\mu}}(z)$ 
and $\widetilde{\mathrm{V}}_{\tilde{\mu}}(z)$ defined by 
Equations~\eqref{wm24frevhmpa1} and~\eqref{wm24frevhmpa2}, 
respectively, into the interpolation condition~\eqref{wm24frevhmpa5}, 
one shows that, in the double-scaling limit $\mathscr{N},n \! \to 
\! \infty$ such that $z_{o} \! = \! 1 \! + \! o(1)$, upon equating 
coefficients of like powers of $z^{-j}$, $j \! = \! 1,2,\dotsc,
2 \varkappa_{nk \tilde{k}_{\mathfrak{s}-1}}^{\infty}$:
\begin{gather}
\dfrac{\tilde{r}_{(n-1)K+k-1}}{\tilde{t}_{(n-1)K+k}} \! = \! 
c^{(\infty)}_{1} \! = \! 1, \label{wm24frevhmpa6} \\
\dfrac{\tilde{r}_{(n-1)K+k-1}}{\tilde{t}_{(n-1)K+k}} \left(
\dfrac{\tilde{r}_{(n-1)K+k-2}}{\tilde{r}_{(n-1)K+k-1}} \! - \! 
\dfrac{\tilde{t}_{(n-1)K+k-1}}{\tilde{t}_{(n-1)K+k}} \right) 
\! = \! c^{(\infty)}_{2}, \label{wm24frevhmpa7} \\
\dfrac{\tilde{r}_{(n-1)K+k-1}}{\tilde{t}_{(n-1)K+k}} 
\widetilde{\varpi}_{\tilde{\mu}}^{\ast}(m \! - \! 1) \! = \! 
c^{(\infty)}_{m}, \quad m \! = \! 3,4,\dotsc,2 \varkappa_{nk 
\tilde{k}_{\mathfrak{s}-1}}^{\infty}, \label{wm24frevhmpa8}
\end{gather}
where, for $m \! = \! 2,3,\dotsc,2 \varkappa_{nk 
\tilde{k}_{\mathfrak{s}-1}}^{\infty} \! - \! 1$,
\begin{equation*}
\widetilde{\varpi}_{\tilde{\mu}}^{\ast}(m) \! := \! 
\mathlarger{\sum_{\substack{j_{1},j_{2}=1,2,\dotsc,2 
\varkappa_{nk \tilde{k}_{\mathfrak{s}-1}}^{\infty}-1\\j_{1}+j_{2}=m}}} 
\, \dfrac{\tilde{r}_{(n-1)K+k-1-j_{1}}}{\tilde{r}_{(n-1)K+k-1}} 
\sum_{m^{\prime}=1}^{j_{2}}(-1)^{m^{\prime}} \widetilde{C}_{
\tilde{\mu}}^{\ast}(m^{\prime},j_{2}) \! + \! \sum_{j=1}^{m}
(-1)^{j} \widetilde{C}_{\tilde{\mu}}^{\ast}(j,m) \! + \! 
\dfrac{\tilde{r}_{(n-1)K+k-1-m}}{\tilde{r}_{(n-1)K+k-1}},
\end{equation*}
with
\begin{equation*}
\widetilde{C}_{\tilde{\mu}}^{\ast}(j,m) \! = \! 
\mathlarger{\sum_{\substack{i_{1},i_{2},\dotsc,i_{j}=1,2,\dotsc,2 
\varkappa_{nk \tilde{k}_{\mathfrak{s}-1}}^{\infty}-(j-1)\\i_{1}+i_{2}
+ \dotsb + i_{j}=m}}} \, \prod_{l=1}^{j} \dfrac{\tilde{t}_{(n-1)
K+k-i_{l}}}{\tilde{t}_{(n-1)K+k}};
\end{equation*}
e.g., for $m \! = \! 2$,
\begin{gather*}
\widetilde{C}_{\tilde{\mu}}^{\ast}(1,1) \! = \! \mathlarger{
\sum_{\substack{i_{1}=1,2,\dotsc,2 \varkappa_{nk 
\tilde{k}_{\mathfrak{s}-1}}^{\infty}\\i_{1}=1}}} \, \dfrac{\tilde{t}_{
(n-1)K+k-i_{1}}}{\tilde{t}_{(n-1)K+k}} \! = \! \dfrac{\tilde{t}_{(n-
1)K+k-1}}{\tilde{t}_{(n-1)K+k}}, \qquad \quad \widetilde{C}_{
\tilde{\mu}}^{\ast}(1,2) \! = \! \mathlarger{\sum_{\substack{i_{1}
=1,2,\dotsc,2 \varkappa_{nk \tilde{k}_{\mathfrak{s}-1}}^{\infty}\\
i_{1}=2}}} \, \dfrac{\tilde{t}_{(n-1)K+k-i_{1}}}{\tilde{t}_{(n-1)K+k}} 
\! = \! \dfrac{\tilde{t}_{(n-1)K+k-2}}{\tilde{t}_{(n-1)K+k}}, \\
\widetilde{C}_{\tilde{\mu}}^{\ast}(2,2) \! = \! \mathlarger{
\sum_{\substack{i_{1},i_{2}=1,2,\dotsc,2 \varkappa_{nk 
\tilde{k}_{\mathfrak{s}-1}}^{\infty}-1\\i_{1}+i_{2}=2}}} \, 
\dfrac{\tilde{t}_{(n-1)K+k-i_{1}}}{\tilde{t}_{(n-1)K+k}} 
\dfrac{\tilde{t}_{(n-1)K+k-i_{2}}}{\tilde{t}_{(n-1)K+k}} \! = 
\! \left(\dfrac{\tilde{t}_{(n-1)K+k-1}}{\tilde{t}_{(n-1)K+k}} 
\right)^{2},
\end{gather*}
in which case
\begin{align*}
\widetilde{\varpi}_{\tilde{\mu}}^{\ast}(2) =& \, -\dfrac{
\tilde{r}_{(n-1)K+k-2}}{\tilde{r}_{(n-1)K+k-1}} \widetilde{C}_{
\tilde{\mu}}^{\ast}(1,1) \! - \! \widetilde{C}_{\tilde{\mu}}^{\ast}
(1,2) \! + \! \widetilde{C}_{\tilde{\mu}}^{\ast}(2,2) \! + \! 
\dfrac{\tilde{r}_{(n-1)K+k-3}}{\tilde{r}_{(n-1)K+k-1}} \\
=& \, -\dfrac{\tilde{r}_{(n-1)K+k-2}}{\tilde{r}_{(n-1)K+k-1}} 
\dfrac{\tilde{t}_{(n-1)K+k-1}}{\tilde{t}_{(n-1)K+k}} \! - \! 
\dfrac{\tilde{t}_{(n-1)K+k-2}}{\tilde{t}_{(n-1)K+k}} \! + \! 
\left(\dfrac{\tilde{t}_{(n-1)K+k-1}}{\tilde{t}_{(n-1)K+k}} 
\right)^{2} \! + \! \dfrac{\tilde{r}_{(n-1)K+k-3}}{\tilde{r}_{
(n-1)K+k-1}},
\end{align*}
for $m \! = \! 3$, with $\widetilde{C}_{\tilde{\mu}}^{\ast}
(1,1)$, $\widetilde{C}_{\tilde{\mu}}^{\ast}(1,2)$, and 
$\widetilde{C}_{\tilde{\mu}}^{\ast}(2,2)$ given above,
\begin{gather*}
\widetilde{C}_{\tilde{\mu}}^{\ast}(1,3) \! = \! \mathlarger{
\sum_{\substack{i_{1}=1,2,\dotsc,2 \varkappa_{nk 
\tilde{k}_{\mathfrak{s}-1}}^{\infty}\\i_{1}=3}}} \, \dfrac{\tilde{t}_{
(n-1)K+k-i_{1}}}{\tilde{t}_{(n-1)K+k}} \! = \! \dfrac{\tilde{t}_{(n
-1)K+k-3}}{\tilde{t}_{(n-1)K+k}}, \\
\widetilde{C}_{\tilde{\mu}}^{\ast}(2,3) \! = \! \mathlarger{
\sum_{\substack{i_{1},i_{2}=1,2,\dotsc,2 \varkappa_{nk 
\tilde{k}_{\mathfrak{s}-1}}^{\infty}-1\\i_{1}+i_{2}=3}}} \, 
\dfrac{\tilde{t}_{(n-1)K+k-i_{1}}}{\tilde{t}_{(n-1)K+k}} 
\dfrac{\tilde{t}_{(n-1)K+k-i_{2}}}{\tilde{t}_{(n-1)K+k}} \! = 
\! 2 \dfrac{\tilde{t}_{(n-1)K+k-1}}{\tilde{t}_{(n-1)K+k}} 
\dfrac{\tilde{t}_{(n-1)K+k-2}}{\tilde{t}_{(n-1)K+k}}, \\
\widetilde{C}_{\tilde{\mu}}^{\ast}(3,3) \! = \! \mathlarger{
\sum_{\substack{i_{1},i_{2},i_{3}=1,2,\dotsc,2 \varkappa_{nk 
\tilde{k}_{\mathfrak{s}-1}}^{\infty}-2\\i_{1}+i_{2}+i_{3}=3}}} 
\, \dfrac{\tilde{t}_{(n-1)K+k-i_{1}}}{\tilde{t}_{(n-1)K+k}} 
\dfrac{\tilde{t}_{(n-1)K+k-i_{2}}}{\tilde{t}_{(n-1)K+k}} 
\dfrac{\tilde{t}_{(n-1)K+k-i_{3}}}{\tilde{t}_{(n-1)K+k}} \! = \! 
\left(\dfrac{\tilde{t}_{(n-1)K+k-1}}{\tilde{t}_{(n-1)K+k}} \right)^{3},
\end{gather*}
in which case,
\begin{align*}
\widetilde{\varpi}_{\tilde{\mu}}^{\ast}(3) =& \, \dfrac{
\tilde{r}_{(n-1)K+k-2}}{\tilde{r}_{(n-1)K+k-1}} \left(
\widetilde{C}_{\tilde{\mu}}^{\ast}(2,2) \! - \! \widetilde{C}_{
\tilde{\mu}}^{\ast}(1,2) \right) \! - \! \dfrac{\tilde{r}_{(n-1)K
+k-3}}{\tilde{r}_{(n-1)K+k-1}} \widetilde{C}_{\tilde{\mu}}^{
\ast}(1,1) \! - \! \widetilde{C}_{\tilde{\mu}}^{\ast}(1,3) \! + \! 
\widetilde{C}_{\tilde{\mu}}^{\ast}(2,3) \! - \! \widetilde{C}_{
\tilde{\mu}}^{\ast}(3,3) \! + \! \dfrac{\tilde{r}_{(n-1)K+k-
4}}{\tilde{r}_{(n-1)K+k-1}} \\
=& \, \dfrac{\tilde{r}_{(n-1)K+k-2}}{\tilde{r}_{(n-1)K+k-1}} 
\left(\left(\dfrac{\tilde{t}_{(n-1)K+k-1}}{\tilde{t}_{(n-1)K+k}} 
\right)^{2} \! - \! \dfrac{\tilde{t}_{(n-1)K+k-2}}{\tilde{t}_{(n-
1)K+k}} \right) \! - \! \dfrac{\tilde{r}_{(n-1)K+k-3}}{\tilde{r}_{
(n-1)K+k-1}} \dfrac{\tilde{t}_{(n-1)K+k-1}}{\tilde{t}_{(n-1)K+
k}} \! - \! \dfrac{\tilde{t}_{(n-1)K+k-3}}{\tilde{t}_{(n-1)K+k}} \\
+& \, 2 \dfrac{\tilde{t}_{(n-1)K+k-1}}{\tilde{t}_{(n-1)K+k}} 
\dfrac{\tilde{t}_{(n-1)K+k-2}}{\tilde{t}_{(n-1)K+k}} \! - \! \left(
\dfrac{\tilde{t}_{(n-1)K+k-1}}{\tilde{t}_{(n-1)K+k}} \right)^{3} 
\! + \! \dfrac{\tilde{r}_{(n-1)K+k-4}}{\tilde{r}_{(n-1)K+k-1}},
\end{align*}
etc. Despite the non-linear appearance of 
Equations~\eqref{wm24frevhmpa6}--\eqref{wm24frevhmpa8}, 
the incredulous fact is that, for $n \! \in \! \mathbb{N}$ and 
$k \! \in \! \lbrace 1,2,\dotsc,K \rbrace$ such that 
$\alpha_{p_{\mathfrak{s}}} \! := \! \alpha_{k} \! \neq \! \infty$, via 
tedious algebraic and concomitant factorisation, re-arrangement, 
recursion, and induction arguments, one shows that they are 
linearised to {}\footnote{Note the convention $\sum_{m=1}^{0} 
\pmb{\ast} \! := \! 0$.}
\begin{equation} \label{wm24frevhmpa9} 
\tilde{r}_{(n-1)K+k-j} \! - \! \sum_{m=1}^{j-1}c^{(\infty)}_{j-m} 
\tilde{t}_{(n-1)K+k-m} \! = \! c^{(\infty)}_{j} \tilde{t}_{(n-1)K+k}, 
\quad j \! = \! 1,2,\dotsc,2 \varkappa_{nk \tilde{k}_{\mathfrak{s}
-1}}^{\infty}.
\end{equation}
Substituting the expressions for the polynomials $\widetilde{
\mathrm{U}}_{\tilde{\mu}}(z)$ and $\widetilde{\mathrm{V}}_{
\tilde{\mu}}(z)$ defined by Equations~\eqref{wm24frevhmpa1} 
and~\eqref{wm24frevhmpa2}, respectively, into the interpolation 
conditions~\eqref{wm24frevhmpa3} and~\eqref{wm24frevhmpa4}, 
and taking note of the differential identity $\tfrac{\partial}{\partial 
\alpha_{p_{q}}}c^{(q)}_{j}(\alpha_{p_{q}}) \! = \! (1 \! + \! j)
c^{(q)}_{j+1}(\alpha_{p_{q}})$, $(j,q) \! \in \! \mathbb{N}_{0} 
\times \lbrace 1,\dotsc,\mathfrak{s} \! - \! 2,\mathfrak{s} 
\rbrace$, one shows that, for $n \! \in \! \mathbb{N}$ and 
$k \! \in \! \lbrace 1,2,\dotsc,K \rbrace$ such that 
$\alpha_{p_{\mathfrak{s}}} \! := \! \alpha_{k} \! \neq \! \infty$, 
in the double-scaling limit $\mathscr{N},n \! \to \! \infty$ such 
that $z_{o} \! = \! 1 \! + \! o(1)$, upon equating coefficients of 
like powers of $(z \! - \! \alpha_{p_{q}})^{i_{1}}$ and $(z \! - \! 
\alpha_{k})^{i_{2}}$, $q \! = \! 1,2,\dotsc,\mathfrak{s} \! - \! 2$, 
$i_{1} \! = \! 0,1,\dotsc,2 \varkappa_{nk \tilde{k}_{q}} \! - \! 1$, 
and $i_{2} \! = \! 0,1,\dotsc,2 \varkappa_{nk} \! - \! 1$,\footnote{In 
Equations~\eqref{wm24frevhmpa10} and~\eqref{wm24frevhmpa11}, 
$\tilde{r}_{j},\tilde{t}_{j}$, $j \! = \! 0,1,\dotsc,(n \! - \! 1)K \! + \! k 
\! - \! 1$, and $\tilde{t}_{(n-1)K+k}$ are not differentiated; note, also, 
the conventions $(\tfrac{\partial}{\partial \alpha_{p_{q}}})^{0} \! := 
\! 0$ and, for $m \! \in \! \mathbb{N}$, $(\tfrac{\partial}{\partial 
\alpha_{p_{q}}})^{m} \! := \! \underbrace{\tfrac{\partial}{\partial 
\alpha_{p_{q}}} \tfrac{\partial}{\partial \alpha_{p_{q}}} \dotsb 
\tfrac{\partial}{\partial \alpha_{p_{q}}}}_{m}$.}
\begin{gather}
\left(\dfrac{\partial}{\partial \alpha_{p_{q}}} \right)^{m_{1}} 
\left(\sum_{j=0}^{(n-1)K+k-1} \tilde{r}_{j}(\alpha_{p_{q}})^{j} 
\! - \! c^{(q)}_{0}(\alpha_{p_{q}}) \sum_{j=0}^{(n-1)K+k} 
\tilde{t}_{j}(\alpha_{p_{q}})^{j} \right) \! = \! 0, \quad m_{1} 
\! = \! 0,1,\dotsc,2 \varkappa_{nk \tilde{k}_{q}} \! - \! 1, 
\quad q \! = \! 1,2,\dotsc,\mathfrak{s} \! - \! 2, 
\label{wm24frevhmpa10} \\
\left(\dfrac{\partial}{\partial \alpha_{k}} \right)^{m_{2}} \left(
\sum_{j=0}^{(n-1)K+k-1} \tilde{r}_{j}(\alpha_{k})^{j} \! - \! 
c^{(\mathfrak{s})}_{0}(\alpha_{k}) \sum_{j=0}^{(n-1)K+k} 
\tilde{t}_{j}(\alpha_{k})^{j} \right) \! = \! 0, \quad m_{2} 
\! = \! 0,1,\dotsc,2 \varkappa_{nk} \! - \! 1.
\label{wm24frevhmpa11}
\end{gather}
For $n \! \in \! \mathbb{N}$ and $k \! \in \! \lbrace 1,2,\dotsc,
K \rbrace$ such that $\alpha_{p_{\mathfrak{s}}} \! := \! \alpha_{k} 
\! \neq \! \infty$, in the double-scaling limit $\mathscr{N},
n \! \to \! \infty$ such that $z_{o} \! = \! 1 \! + \! o(1)$, 
Equations~\eqref{wm24frevhmpa9}--\eqref{wm24frevhmpa11} 
constitute the system of $2((n \! - \! 1)K \! + \! k)$ linear 
inhomogeneous algebraic equations for the $2((n \! - \! 1)K \! 
+ \! k)$ coefficients $\tilde{r}_{0},\tilde{r}_{1},\dotsc,\tilde{r}_{
(n-1)K+k-1},\tilde{t}_{0},\tilde{t}_{1},\dotsc,\tilde{t}_{(n-1)K+k-1}$, 
which can be written in the following ordered block-matrix 
form:
\begin{equation} \label{wm24frevhmpa12} 
\left(
, 
\end{equation}
where, for $q \! \in \! \lbrace 1,2,\dotsc,\mathfrak{s} \! - \! 2 \rbrace$, 
$q \! = \! \mathfrak{s} \! - \! 1$, and $q \! = \! \mathfrak{s}$, the 
corresponding $2 \varkappa_{nk \tilde{k}_{q}} \times ((n \! - \! 1)K \! 
+ \! k)$, $2 \varkappa_{nk \tilde{k}_{\mathfrak{s}-1}}^{\infty} \times 
((n \! - \! 1)K \! + \! k)$, and $2 \varkappa_{nk} \times ((n \! - \! 1)K 
\! + \! k)$ non-degenerate sub-block matrices $\widetilde{\mathrm{
W}}_{\mathrm{L}}^{\ast}(q)$ and $\widetilde{\mathrm{W}}_{\mathrm{
R}}^{\ast}(q)$, $\widetilde{\mathrm{W}}_{\mathrm{L}}^{\ast}
(\mathfrak{s} \! - \! 1)$ and $\widetilde{\mathrm{W}}_{\mathrm{R}}^{
\ast}(\mathfrak{s} \! - \! 1)$, and $\widetilde{\mathrm{W}}_{\mathrm{
L}}^{\ast}(\mathfrak{s})$ and $\widetilde{\mathrm{W}}_{\mathrm{R}}^{
\ast}(\mathfrak{s})$, respectively, are given by
\begin{gather*}
\left(\widetilde{\mathrm{W}}_{\mathrm{L}}^{\ast}(q) \right)_{ij} 
\underset{q=1,2,\dotsc,\mathfrak{s}-2}{:=} \left(\dfrac{\partial}{
\partial \alpha_{p_{q}}} \right)^{i-1} \alpha_{p_{q}}^{j-1}, \quad 
i \! = \! 1,\dotsc,2\varkappa_{nk \tilde{k}_{q}}, \, j \! = \! 1,\dotsc,
(n \! - \! 1)K \! + \! k, \\
\left(\widetilde{\mathrm{W}}_{\mathrm{R}}^{\ast}(q) \right)_{ij} 
\underset{q=1,2,\dotsc,\mathfrak{s}-2}{:=} -\left(\dfrac{\partial}{
\partial \alpha_{p_{q}}} \right)^{i-1} \left(\dfrac{c^{(q)}_{0}
(\alpha_{p_{q}})}{\alpha_{p_{q}}^{(n-1)K+k-j+1}} \right), \quad i \! 
= \! 1,\dotsc,2 \varkappa_{nk \tilde{k}_{q}}, \, j \! = \! (n \! - \! 
1)K \! + \! k \! + \! 1,\dotsc,2((n \! - \! 1)K \! + \! k), \\
\left(\widetilde{\mathrm{W}}_{\mathrm{L}}^{\ast}(\mathfrak{s} \! 
- \! 1) \right)_{ij} := \delta_{i \, \, (n-1)K+k-j+1}, \quad i \! = \! 
1,\dotsc,2\varkappa_{nk \tilde{k}_{\mathfrak{s}-1}}^{\infty}, \, 
j \! = \! 1,\dotsc,(n \! - \! 1)K \! + \! k, \\
\left(\widetilde{\mathrm{W}}_{\mathrm{R}}^{\ast}(\mathfrak{s} 
\! - \! 1) \right)_{ij} := 
\begin{cases}
0, &\text{$i \! = \! 1, \, j \! = \! (n \! - \! 1)K \! + \! k \! + \! 1,
\dotsc,2((n \! - \! 1)K \! + \! k)$,} \\
-c^{(\infty)}_{i+j-2((n-1)K+k)-1}, &\text{$i \! = \! 2,\dotsc,2 
\varkappa_{nk \tilde{k}_{\mathfrak{s}-1}}^{\infty}, \, j \! = \! 2((n \! - 
\! 1)K \! + \! k) \! - \! i \! + \! 2,\dotsc,2((n \! - \! 1)K \! + \! k)$,}
\end{cases} \\
\left(\widetilde{\mathrm{W}}_{\mathrm{L}}^{\ast}(\mathfrak{s}) 
\right)_{ij} := \left(\dfrac{\partial}{\partial \alpha_{k}} \right)^{i-1} 
\alpha_{k}^{j-1}, \quad i \! = \! 1,\dotsc,2\varkappa_{nk}, \, 
j \! = \! 1,\dotsc,(n \! - \! 1)K \! + \! k, \\
\left(\widetilde{\mathrm{W}}_{\mathrm{R}}^{\ast}(\mathfrak{s}) 
\right)_{ij} := -\left(\dfrac{\partial}{\partial \alpha_{k}} \right)^{i-1} 
\left(\dfrac{c^{(\mathfrak{s})}_{0}(\alpha_{k})}{\alpha_{k}^{(n-1)K
+k-j+1}} \right), 
\quad i \! = \! 1,\dotsc,2 \varkappa_{nk}, \, j \! = \! (n \! - \! 1)K 
\! + \! k \! + \! 1,\dotsc,2((n \! - \! 1)K \! + \! k),
\end{gather*}
and, for $q \! \in \! \lbrace 1,2,\dotsc,\mathfrak{s} \! - \! 2 \rbrace$, 
$q \! = \! \mathfrak{s} \! - \! 1$, and $q \! = \! \mathfrak{s}$, 
the corresponding $2 \varkappa_{nk \tilde{k}_{q}} \times 1$, $2 
\varkappa_{nk \tilde{k}_{\mathfrak{s}-1}}^{\infty} \times 1$, and $2 
\varkappa_{nk} \times 1$ column tuples $\widetilde{\P}^{\ast}(q)$, 
$\widetilde{\P}^{\ast}(\mathfrak{s} \! - \! 1)$, and $\widetilde{\P}^{
\ast}(\mathfrak{s})$, respectively, are given by
\begin{gather*}
\left(\widetilde{\P}^{\ast}(q) \right)_{i} \underset{q=1,2,\dotsc,
\mathfrak{s}-2}{:=} (i \! - \! 1)! \alpha_{p_{q}}^{(n-1)K+k}
c^{(q)}_{i-1}(\alpha_{p_{q}}) \tilde{t}_{(n-1)K+k}, \quad i \! 
= \! 1,2,\dotsc,2 \varkappa_{nk \tilde{k}_{q}}, \\
\left(\widetilde{\P}^{\ast}(\mathfrak{s} \! - \! 1) \right)_{i} := 
c^{(\infty)}_{i} \tilde{t}_{(n-1)K+k}, \quad i \! = \! 1,2,\dotsc,
2 \varkappa_{nk \tilde{k}_{\mathfrak{s}-1}}^{\infty}, \\
\left(\widetilde{\P}^{\ast}(\mathfrak{s}) \right)_{i} := (i \! - \! 1)! 
\alpha_{k}^{(n-1)K+k}c^{(\mathfrak{s})}_{i-1}(\alpha_{k}) 
\tilde{t}_{(n-1)K+k}, \quad i \! = \! 1,2,\dotsc,2 \varkappa_{nk}.
\end{gather*}
For $n \! \in \! \mathbb{N}$ and $k \! \in \! \lbrace 1,2,\dotsc,K 
\rbrace$ such that $\alpha_{p_{\mathfrak{s}}} \! := \! \alpha_{k} 
\! \neq \! \infty$, via the definition of the Markov-Stieltjes 
transform and the associated $\mathrm{R}$-function given by 
Equations~\eqref{eqmvsstildemu} and~\eqref{eqlemmvssfinmpa1}, 
respectively, one proceeds as follows:
\begin{gather}
\mathrm{F}_{\tilde{\mu}}(z) \! = \! \int_{\mathbb{R}}(z \! - \! \xi)^{-1} 
\, \md \widetilde{\mu}(\xi) \quad \Rightarrow \quad \mathrm{F}_{
\tilde{\mu}}(z) \pmb{\pi}^{n}_{k}(z) \! = \! \underbrace{\int_{
\mathbb{R}} \left(\dfrac{\pmb{\pi}^{n}_{k}(\xi) \! - \! \pmb{\pi}^{n}_{k}
(z)}{\xi \! - \! z} \right) \md \widetilde{\mu}(\xi)}_{= \, \widetilde{
\pmb{\mathrm{R}}}_{\tilde{\mu}}(z)} + \! \int_{\mathbb{R}} 
\dfrac{\pmb{\pi}^{n}_{k}(\xi)}{z \! - \! \xi} \, \md \widetilde{\mu}
(\xi) \quad \nonumber \Rightarrow \\
\mathrm{F}_{\tilde{\mu}}(z) \! = \! \dfrac{\widetilde{\pmb{\mathrm{R}}}_{
\tilde{\mu}}(z)}{\pmb{\pi}^{n}_{k}(z)} \! + \! \dfrac{1}{(z \! - \! 
\alpha_{k}) \pmb{\pi}^{n}_{k}(z)} \left((z \! - \! \alpha_{k}) \int_{
\mathbb{R}} \dfrac{((\xi \! - \! \alpha_{k}) \pmb{\pi}^{n}_{k}
(\xi))}{(\xi \! - \! \alpha_{k})(z \! - \! \xi)} \, \md \widetilde{\mu}
(\xi) \right); \label{wm24frevhmpa13}
\end{gather} 
hence, via the $(1 \, 1)$- and $(1 \, 2)$-elements, respectively, of the matrix 
representation~\eqref{intrepfin}, and the definition of the corresponding 
MPA error term $\widetilde{\pmb{\mathrm{E}}}_{\tilde{\mu}}(z) \! 
:= \! (\pmb{\pi}^{n}_{k}(z))^{-1} \widetilde{\pmb{\mathrm{R}}}_{
\tilde{\mu}}(z) \! - \! \mathrm{F}_{\tilde{\mu}}(z)$,\footnote{In this 
context, `error term' refers to the fact that, for $n \! \in \! \mathbb{N}$ 
and $k \! \in \! \lbrace 1,2,\dotsc,K \rbrace$ such that 
$\alpha_{p_{\mathfrak{s}}} \! := \! \alpha_{k} \! \neq \! \infty$, 
in the double-scaling limit $\mathscr{N},n \! \to \! \infty$ such 
that $z_{o} \! = \! 1 \! + \! o(1)$, $\widetilde{\pmb{\mathrm{E}}}_{
\tilde{\mu}}(z) \! =_{z \to \alpha_{p_{q}}} \! \mathcal{O}(\tilde{
\mathfrak{c}}^{(q)}_{\tilde{\mu}}(n,k,z_{o})(z \! - \! \alpha_{p_{q}})^{
2 \varkappa_{nk \tilde{k}_{q}}})$, $q \! = \! 1,2,\dotsc,\mathfrak{s} 
\! - \! 2$, $\widetilde{\pmb{\mathrm{E}}}_{\tilde{\mu}}(z) \! =_{z 
\to \alpha_{k}} \! \mathcal{O}(\tilde{\mathfrak{c}}^{(\mathfrak{s})}_{
\tilde{\mu}}(n,k,z_{o})(z \! - \! \alpha_{k})^{2 \varkappa_{nk}})$, and 
$\widetilde{\pmb{\mathrm{E}}}_{\tilde{\mu}}(z) \! =_{z \to \alpha_{
p_{\mathfrak{s}-1} = \infty}} \! \mathcal{O} \left(\tilde{\mathfrak{c}}^{
(\mathfrak{s}-1)}_{\tilde{\mu}}(n,k,z_{o})z^{-(2 \varkappa_{nk 
\tilde{k}_{\mathfrak{s}-1}}^{\infty}+1)} \right)$, where $\tilde{
\mathfrak{c}}^{(\mathfrak{q})}_{\tilde{\mu}}(n,k,z_{o}) \! =_{
\underset{z_{o}=1+o(1)}{\mathscr{N},n \to \infty}} \! \mathcal{O}
(1)$, $\mathfrak{q} \! \in \! \lbrace 1,2,\dotsc,\mathfrak{s} 
\rbrace$.} one arrives at Equation~\eqref{eqlemmvssfinmpa7}.

\pmb{(2)} The proof of this case, that is, $n \! \in \! \mathbb{N}$ 
and $k \! \in \! \lbrace 1,2,\dotsc,K \rbrace$ such that 
$\alpha_{p_{\mathfrak{s}}} \! := \! \alpha_{k} \! = \! \infty$, is 
virtually identical to the proof of \pmb{(1)} above; one mimics, 
\emph{verbatim}, the scheme of the calculations presented in 
\pmb{(1)} above in order to arrive at the corresponding claims 
stated in item~\pmb{(1)} of the lemma; in order to do so, however, 
the analogues of 
Equations~\eqref{wm24frevhmpa1}--\eqref{wm24frevhmpa5} and 
\eqref{wm24frevhmpa9}--\eqref{wm24frevhmpa13} are necessary, 
which, in the present case, read:
\begin{gather}
\widehat{\pmb{\mathrm{R}}}_{\tilde{\mu}}(z) \! = \! \dfrac{
\widehat{\mathrm{U}}_{\tilde{\mu}}(z)}{\prod_{q=1}^{\mathfrak{s}-1}
(z \! - \! \alpha_{p_{q}})^{\varkappa_{nk \tilde{k}_{q}}}} \! := \! \dfrac{
\sum_{j=0}^{(n-1)K+k-1} \hat{r}_{j}z^{j}}{\prod_{q=1}^{\mathfrak{s}
-1}(z \! - \! \alpha_{p_{q}})^{\varkappa_{nk \tilde{k}_{q}}}}, 
\label{wm24frevhmpa14} \\
\pmb{\pi}^{n}_{k}(z) \! = \! \dfrac{\widehat{\mathrm{V}}_{\tilde{\mu}}
(z)}{\prod_{q=1}^{\mathfrak{s}-1}(z \! - \! \alpha_{p_{q}})^{\varkappa_{nk 
\tilde{k}_{q}}}} \! := \! \dfrac{\sum_{j=0}^{(n-1)K+k} \hat{t}_{j}z^{j}}{
\prod_{q=1}^{\mathfrak{s}-1}(z \! - \! \alpha_{p_{q}})^{\varkappa_{nk 
\tilde{k}_{q}}}}, \label{wm24frevhmpa15} \\
\dfrac{\widehat{\mathrm{U}}_{\tilde{\mu}}(z)}{\widehat{\mathrm{V}}_{
\tilde{\mu}}(z)} \! - \! \sum_{j=0}^{2 \varkappa_{nk \tilde{k}_{q}}-1}
c^{(q)}_{j}(\alpha_{p_{q}})(z \! - \! \alpha_{p_{q}})^{j} \underset{z \to 
\alpha_{p_{q}}}{=} \mathcal{O} \left(\hat{\mathfrak{c}}^{(q)}_{\tilde{\mu}}
(n,k,z_{o})(z \! - \! \alpha_{p_{q}})^{2 \varkappa_{nk \tilde{k}_{q}}} \right), 
\quad q \! = \! 1,2,\dotsc,\mathfrak{s} \! - \! 1, \label{wm24frevhmpa16} \\
\dfrac{\widehat{\mathrm{U}}_{\tilde{\mu}}(z)}{\widehat{\mathrm{V}}_{
\tilde{\mu}}(z)} \! - \! \sum_{j=1}^{2 \varkappa_{nk}}c^{(\infty)}_{j}z^{-j} 
\underset{z \to \alpha_{k}}{=} \mathcal{O} \left(\hat{\mathfrak{c}}^{
(\mathfrak{s})}_{\tilde{\mu}}(n,k,z_{o})z^{-(2 \varkappa_{nk}+1)} \right), 
\label{wm24frevhmpa17}
\end{gather}
where $c^{(q)}_{j}(\alpha_{p_{q}}) \! := \! -\int_{\mathbb{R}}(\xi \! - \! 
\alpha_{p_{q}})^{-(1+j)} \, \md \widetilde{\mu}(\xi)$, $(j,q) \! \in \! 
\mathbb{N}_{0} \times \lbrace 1,2,\dotsc,\mathfrak{s} \! - \! 1 \rbrace$, 
and $c^{(\infty)}_{i} \! := \! \int_{\mathbb{R}} \xi^{i-1} \, \md 
\widetilde{\mu}(\xi)$, $i \! \in \! \mathbb{N}$, with $c^{(\infty)}_{1} 
\! = \! 1$, and $\hat{\mathfrak{c}}^{(m)}_{\tilde{\mu}}(n,k,z_{o}) \! 
=_{\underset{z_{o}=1+o(1)}{\mathscr{N},n \to \infty}} \! \mathcal{O}
(1)$, $m \! \in \! \lbrace 1,2,\dotsc,\mathfrak{s} \rbrace$, $\hat{r}_{j} 
\! = \! \hat{r}_{j}(n,k,z_{o})$ and $\hat{t}_{j} \! = \! \hat{t}_{j}(n,k,z_{o})$, 
$j \! = \! 0,1,\dotsc,(n \! - \! 1)K \! + \! k \! - \! 1$, with $\hat{t}_{(n
-1)K+k} \! = \! \hat{t}_{(n-1)K+k}(n,k,z_{o}) \! = \! 1$, solve the 
system of $2((n \! - \! 1)K \! + \! k)$ linear inhomogeneous algebraic 
equations {}\footnote{In Equation~\eqref{wm24frevhmpa19}, $\hat{r}_{j},
\hat{t}_{j}$, $j \! = \! 0,1,\dotsc,(n \! - \! 1)K \! + \! k \! - \! 1$, and 
$\hat{t}_{(n-1)K+k} \! = \! 1$ are not differentiated.}
\begin{gather}
\hat{r}_{(n-1)K+k-j} \! - \! \sum_{m=1}^{j-1}c^{(\infty)}_{j-m} 
\hat{t}_{(n-1)K+k-m} \! = \! c^{(\infty)}_{j}, \quad j \! = \! 1,2,
\dotsc,2 \varkappa_{nk}, \label{wm24frevhmpa18} \\
\left(\dfrac{\partial}{\partial \alpha_{p_{q}}} \right)^{m} \left(
\sum_{j=0}^{(n-1)K+k-1} \hat{r}_{j}(\alpha_{p_{q}})^{j} \! - \! 
c^{(q)}_{0}(\alpha_{p_{q}}) \sum_{j=0}^{(n-1)K+k} \hat{t}_{j}
(\alpha_{p_{q}})^{j} \right) \! = \! 0, \quad m \! = \! 0,1,\dotsc,
2 \varkappa_{nk \tilde{k}_{q}} \! - \! 1, \quad q \! = \! 1,2,\dotsc,
\mathfrak{s} \! - \! 1, \label{wm24frevhmpa19}
\end{gather}
which can be presented, equivalently, as the ordered block-matrix form
\begin{equation} \label{wm24frevhmpa20} 
\left(
\begin{smallmatrix}
{} & {} \\
\boxed{
\begin{matrix}
{} & {} & {} \\
{} & \widehat{\mathrm{W}}_{\mathrm{L}}^{\ast}(1) & {} \\
{} & {} & {}
\end{matrix}} & 
\boxed{
\begin{matrix}
{} & {} & {} \\
{} & \widehat{\mathrm{W}}_{\mathrm{R}}^{\ast}(1) & {} \\
{} & {} & {}
\end{matrix}} \\
{} & {} \\
\boxed{
\begin{matrix}
{} & {} & {} \\
{} & \widehat{\mathrm{W}}_{\mathrm{L}}^{\ast}(2) & {} \\
{} & {} & {}
\end{matrix}} & 
\boxed{
\begin{matrix}
{} & {} & {} \\
{} & \widehat{\mathrm{W}}_{\mathrm{R}}^{\ast}(2) & {} \\
{} & {} & {}
\end{matrix}} \\
{} & {} \\
\vdots & \vdots \\
{} & {} \\
\boxed{
\begin{matrix}
{} & {} & {} \\
{} & \widehat{\mathrm{W}}_{\mathrm{L}}^{\ast}(\mathfrak{s} \! - \! 1) & {} \\
{} & {} & {}
\end{matrix}} & 
\boxed{
\begin{matrix}
{} & {} & {} \\
{} & \widehat{\mathrm{W}}_{\mathrm{R}}^{\ast}(\mathfrak{s} \! - \! 1) & {} \\
{} & {} & {}
\end{matrix}} \\
{} & {} \\
\boxed{
\begin{matrix}
{} & {} & {} \\
{} & \widehat{\mathrm{W}}_{\mathrm{L}}^{\ast}(\mathfrak{s}) & {} \\
{} & {} & {}
\end{matrix}} & 
\boxed{
\begin{matrix}
{} & {} & {} \\
{} & \widehat{\mathrm{W}}_{\mathrm{R}}^{\ast}(\mathfrak{s}) & {} \\
{} & {} & {}
\end{matrix}} \\
{} & {}
\end{smallmatrix}
\right)
\begin{pmatrix}
\hat{r}_{0} \\
\hat{r}_{1} \\
\vdots \\
\vdots \\
\hat{r}_{(n-1)K+k-2} \\
\hat{r}_{(n-1)K+k-1} \\
\hat{t}_{0} \\
\hat{t}_{1} \\
\vdots \\
\vdots \\
\hat{t}_{(n-1)K+k-2} \\
\hat{t}_{(n-1)K+k-1}
\end{pmatrix} = 
\begin{pmatrix}
\widehat{\P}^{\ast}(1) \\
\widehat{\P}^{\ast}(2) \\
\vdots \\
\vdots \\
\vdots \\
\vdots \\
\vdots \\
\vdots \\
\widehat{\P}^{\ast}(\mathfrak{s} \! - \! 1) \\
\widehat{\P}^{\ast}(\mathfrak{s})
\end{pmatrix}, 
\end{equation}
where, for $q \! \in \! \lbrace 1,2,\dotsc,\mathfrak{s} \! - \! 1 \rbrace$ 
and $q \! = \! \mathfrak{s}$, the corresponding $2 \varkappa_{nk 
\tilde{k}_{q}} \times ((n \! - \! 1)K \! + \! k)$ and $2 \varkappa_{nk} 
\times ((n \! - \! 1)K \! + \! k)$ non-degenerate sub-block matrices 
$\widehat{\mathrm{W}}_{\mathrm{L}}^{\ast}(q)$ and $\widehat{
\mathrm{W}}_{\mathrm{R}}^{\ast}(q)$, and $\widehat{\mathrm{W}}_{
\mathrm{L}}^{\ast}(\mathfrak{s})$ and $\widehat{\mathrm{W}}_{
\mathrm{R}}^{\ast}(\mathfrak{s})$, respectively, are given by
\begin{gather*}
\left(\widehat{\mathrm{W}}_{\mathrm{L}}^{\ast}(q) \right)_{ij} 
\underset{q=1,2,\dotsc,\mathfrak{s}-1}{:=} \left(\dfrac{\partial}{
\partial \alpha_{p_{q}}} \right)^{i-1} \alpha_{p_{q}}^{j-1}, \quad 
i \! = \! 1,\dotsc,2\varkappa_{nk \tilde{k}_{q}}, \, j \! = \! 1,\dotsc,
(n \! - \! 1)K \! + \! k, \\
\left(\widehat{\mathrm{W}}_{\mathrm{R}}^{\ast}(q) \right)_{ij} 
\underset{q=1,2,\dotsc,\mathfrak{s}-1}{:=} -\left(\dfrac{\partial}{
\partial \alpha_{p_{q}}} \right)^{i-1} \left(\dfrac{c^{(q)}_{0}
(\alpha_{p_{q}})}{\alpha_{p_{q}}^{(n-1)K+k-j+1}} \right), \quad i \! 
= \! 1,\dotsc,2 \varkappa_{nk \tilde{k}_{q}}, \, j \! = \! (n \! - \! 1)
K \! + \! k \! + \! 1,\dotsc,2((n \! - \! 1)K \! + \! k), \\
\left(\widehat{\mathrm{W}}_{\mathrm{L}}^{\ast}(\mathfrak{s}) 
\right)_{ij} := \delta_{i \, \, (n-1)K+k-j+1}, \quad i \! = \! 1,\dotsc,
2 \varkappa_{nk}, \, j \! = \! 1,\dotsc,(n \! - \! 1)K \! + \! k, \\
\left(\widehat{\mathrm{W}}_{\mathrm{R}}^{\ast}(\mathfrak{s}) 
\right)_{ij} := 
\begin{cases}
0, &\text{$i \! = \! 1, \, j \! = \! (n \! - \! 1)K \! + \! k \! + \! 
1,\dotsc,2((n \! - \! 1)K \! + \! k)$,} \\
-c^{(\infty)}_{i+j-2((n-1)K+k)-1}, &\text{$i \! = \! 2,\dotsc,2 
\varkappa_{nk}, \, j \! = \! 2((n \! - \! 1)K \! + \! k) \! - \! i \! 
+ \! 2,\dotsc,2((n \! - \! 1)K \! + \! k)$,}
\end{cases}
\end{gather*}
and, for $q \! \in \! \lbrace 1,2,\dotsc,\mathfrak{s} \! - \! 1 \rbrace$ 
and $q \! = \! \mathfrak{s}$, the corresponding $2 \varkappa_{nk 
\tilde{k}_{q}} \times 1$ and $2 \varkappa_{nk} \times 1$ column tuples 
$\widehat{\P}^{\ast}(q)$ and $\widehat{\P}^{\ast}(\mathfrak{s})$, 
respectively, are given by
\begin{gather*}
\left(\widehat{\P}^{\ast}(q) \right)_{i} \underset{q=1,2,\dotsc,
\mathfrak{s}-1}{:=} (i \! - \! 1)! \alpha_{p_{q}}^{(n-1)K+k}
c^{(q)}_{i-1}(\alpha_{p_{q}}), \quad i \! = \! 1,2,\dotsc,2 \varkappa_{nk 
\tilde{k}_{q}}, \\
\left(\widehat{\P}^{\ast}(\mathfrak{s}) \right)_{i} := c^{(\infty)}_{i}, 
\quad i \! = \! 1,2,\dotsc,2 \varkappa_{nk},
\end{gather*}
and
\begin{equation} \label{wm24frevhmpa21} 
\mathrm{F}_{\tilde{\mu}}(z) \! = \! \dfrac{\widehat{\pmb{\mathrm{R}}}_{
\tilde{\mu}}(z)}{\pmb{\pi}^{n}_{k}(z)} \! + \! \dfrac{1}{\pmb{\pi}^{n}_{k}
(z)} \int_{\mathbb{R}} \dfrac{\pmb{\pi}^{n}_{k}(\xi)}{z \! - \! \xi} 
\, \md \widetilde{\mu}(\xi),
\end{equation} 
hence, via the $(1 \, 1)$- and $(1 \, 2)$-elements, respectively, of the matrix 
representation~\eqref{intrepinf}, and the definition of the corresponding 
MPA error term $\widehat{\pmb{\mathrm{E}}}_{\tilde{\mu}}(z) \! := \! 
(\pmb{\pi}^{n}_{k}(z))^{-1} \widehat{\pmb{\mathrm{R}}}_{\tilde{\mu}}(z) 
\! - \! \mathrm{F}_{\tilde{\mu}}(z)$,\footnote{In this context, `error term' 
refers to the fact that, for $n \! \in \! \mathbb{N}$ and $k \! \in \! \lbrace 
1,2,\dotsc,K \rbrace$ such that $\alpha_{p_{\mathfrak{s}}} \! := \! \alpha_{k} 
\! = \! \infty$, in the double-scaling limit $\mathscr{N},n \! \to \! \infty$ 
such that $z_{o} \! = \! 1 \! + \! o(1)$, $\widehat{\pmb{\mathrm{E}}}_{
\tilde{\mu}}(z) \! =_{z \to \alpha_{p_{q}}} \! \mathcal{O}(\hat{
\mathfrak{c}}^{(q)}_{\tilde{\mu}}(n,k,z_{o})(z \! - \! \alpha_{p_{q}})^{
2 \varkappa_{nk \tilde{k}_{q}}})$, $q \! = \! 1,2,\dotsc,\mathfrak{s} 
\! - \! 1$, and $\widehat{\pmb{\mathrm{E}}}_{\tilde{\mu}}(z) \! 
=_{z \to \alpha_{k}} \! \mathcal{O}(\hat{\mathfrak{c}}^{(\mathfrak{s})}_{
\tilde{\mu}}(n,k,z_{o})z^{-(2 \varkappa_{nk}+1)})$, where 
$\hat{\mathfrak{c}}^{(\mathfrak{q})}_{\tilde{\mu}}(n,k,z_{o}) 
\! =_{\underset{z_{o}=1+o(1)}{\mathscr{N},n \to \infty}} \! 
\mathcal{O}(1)$, $\mathfrak{q} \! \in \! \lbrace 1,2,\dotsc,
\mathfrak{s} \rbrace$.} one arrives at Equation~\eqref{eqlemmvssinfmpa3}. 
\hfill $\qed$

All the necessary ingredients are now available in order to state the 
final results of this monograph, namely, for $n \! \in \! \mathbb{N}$ 
and $k \! \in \! \lbrace 1,2,\dotsc,K \rbrace$ such that $\alpha_{
p_{\mathfrak{s}}} \! := \! \alpha_{k} \! = \! \infty$ or $\alpha_{p_{
\mathfrak{s}}} \! := \! \alpha_{k} \! \neq \! \infty$, asymptotics, in the 
double-scaling limit $\mathscr{N},n \! \to \! \infty$ such that $z_{o} 
\! = \! 1 \! + \! o(1)$, for the monic MPC ORF, $\pmb{\pi}_{k}^{n}
(z)$, $z \! \in \! \mathbb{C}$, the corresponding MPA error terms, 
$\widehat{\pmb{\mathrm{E}}}_{\tilde{\mu}}(z)$ (for $\alpha_{
p_{\mathfrak{s}}} \! := \! \alpha_{k} \! = \! \infty)$ and $\widetilde{
\pmb{\mathrm{E}}}_{\tilde{\mu}}(z)$ (for $\alpha_{p_{\mathfrak{s}}} 
\! := \! \alpha_{k} \! \neq \! \infty)$, $z \! \in \! \mathbb{C}$, the 
associated norming constants, $\mu_{n,\varkappa_{nk}}^{r}(n,k)$, 
$r \! \in \! \lbrace \infty,f \rbrace$, and the MPC ORF, $\phi_{k}^{n}
(z)$,  $z \! \in \! \mathbb{C}$. In order to derive the above-mentioned 
asymptotics for the monic MPC ORF, $\pmb{\pi}_{k}^{n}(z)$, however, 
one must re-trace the sequence of RHP transformations (all of which 
are invertible) stated in Lemmata~\ref{lem3.4}, \ref{lem4.1}, \ref{lem4.2}, 
and~\ref{lem4.5}--\ref{lem4.10}, and Proposition~\ref{propo5.3}, for 
the solution of the monic MPC ORF RHP $(\mathcal{X}(z),\upsilon (z),
\overline{\mathbb{R}})$ stated in Lemma~$\bm{\mathrm{RHP}_{\mathrm{MPC}}}$: 
(1) for $n \! \in \! \mathbb{N}$ and $k \! \in \! \lbrace 1,2,\dotsc,K 
\rbrace$ such that $\alpha_{p_{\mathfrak{s}}} \! := \! \alpha_{k} \! 
= \! \infty$, and $j \! = \! 1,2,\dotsc,N \! + \! 1$,
\begin{equation} \label{finmpcorfmonichat} 
\mathcal{X}(z) \! = \! 
\begin{cases}
\me^{\frac{n \hat{\ell}}{2} \operatorname{ad}(\sigma_{3})} \hat{\mathcal{R}}
(z) \tilde{\mathfrak{m}}^{\raise-0.5ex\hbox{$\scriptstyle \infty$}} \hat{
\mathbb{M}}(z) \me^{n(g^{\infty}(z)-\tilde{\mathscr{P}}_{0}) \sigma_{3}}, 
&\text{$z \! \in \! \hat{\Upsilon}_{1}$,} \\
-\mi \me^{\frac{n \hat{\ell}}{2} \operatorname{ad}(\sigma_{3})} \hat{\mathcal{R}}
(z) \tilde{\mathfrak{m}}^{\raise-0.5ex\hbox{$\scriptstyle \infty$}} \hat{
\mathbb{M}}(z) \sigma_{2} \me^{n(g^{\infty}(z)-\tilde{\mathscr{P}}_{0}) 
\sigma_{3}}, &\text{$z \! \in \! \hat{\Upsilon}_{2}$,} \\
\me^{\frac{n \hat{\ell}}{2} \operatorname{ad}(\sigma_{3})} \hat{\mathcal{R}}
(z) \tilde{\mathfrak{m}}^{\raise-0.5ex\hbox{$\scriptstyle \infty$}} \hat{
\mathbb{M}}(z) \left(\mathrm{I} \! + \! \me^{-2 \pi \mi ((n-1)K+k) 
\int_{z}^{\hat{a}_{N+1}} \psi_{\widetilde{V}}^{\infty}(\xi) \, \md \xi} 
\sigma_{-} \right) \me^{n(g^{\infty}(z)-\tilde{\mathscr{P}}_{0}) \sigma_{3}}, 
&\text{$z \! \in \! \hat{\Upsilon}_{3}$,} \\
-\mi \me^{\frac{n \hat{\ell}}{2} \operatorname{ad}(\sigma_{3})} \hat{\mathcal{R}}
(z) \tilde{\mathfrak{m}}^{\raise-0.5ex\hbox{$\scriptstyle \infty$}} \hat{
\mathbb{M}}(z) \sigma_{2} \left(\mathrm{I} \! - \! \me^{2 \pi \mi ((n-1)K+k) 
\int_{z}^{\hat{a}_{N+1}} \psi_{\widetilde{V}}^{\infty}(\xi) \, \md \xi} 
\sigma_{-} \right) \me^{n(g^{\infty}(z)-\tilde{\mathscr{P}}_{0}) \sigma_{3}}, 
&\text{$z \! \in \! \hat{\Upsilon}_{4}$,} \\
\me^{\frac{n \hat{\ell}}{2} \operatorname{ad}(\sigma_{3})} \hat{\mathcal{R}}(z) 
\hat{\mathcal{X}}^{\hat{b}}(z) \me^{n(g^{\infty}(z)-\tilde{\mathscr{P}}_{0}) 
\sigma_{3}}, &\text{$z \! \in \! \hat{\Omega}_{\hat{b}_{j-1}}^{1}$,} \\
\me^{\frac{n \hat{\ell}}{2} \operatorname{ad}(\sigma_{3})} \hat{\mathcal{R}}(z) 
\hat{\mathcal{X}}^{\hat{a}}(z) \me^{n(g^{\infty}(z)-\tilde{\mathscr{P}}_{0}) 
\sigma_{3}}, &\text{$z \! \in \! \hat{\Omega}_{\hat{a}_{j}}^{1}$,} \\
\me^{\frac{n \hat{\ell}}{2} \operatorname{ad}(\sigma_{3})} \hat{\mathcal{R}}(z) 
\hat{\mathcal{X}}^{\hat{b}}(z) \left(\mathrm{I} \! + \! \me^{-2 \pi \mi ((n-1)K
+k) \int_{z}^{\hat{a}_{N+1}} \psi_{\widetilde{V}}^{\infty}(\xi) \, \md \xi} 
\sigma_{-} \right) \me^{n(g^{\infty}(z)-\tilde{\mathscr{P}}_{0}) \sigma_{3}}, 
&\text{$z \! \in \! \hat{\Omega}_{\hat{b}_{j-1}}^{2}$,} \\
\me^{\frac{n \hat{\ell}}{2} \operatorname{ad}(\sigma_{3})} \hat{\mathcal{R}}(z) 
\hat{\mathcal{X}}^{\hat{a}}(z) \left(\mathrm{I} \! + \! \me^{-2 \pi \mi ((n-1)K
+k) \int_{z}^{\hat{a}_{N+1}} \psi_{\widetilde{V}}^{\infty}(\xi) \, \md \xi} 
\sigma_{-} \right) \me^{n(g^{\infty}(z)-\tilde{\mathscr{P}}_{0}) \sigma_{3}}, 
&\text{$z \! \in \! \hat{\Omega}_{\hat{a}_{j}}^{2}$,} \\
\me^{\frac{n \hat{\ell}}{2} \operatorname{ad}(\sigma_{3})} \hat{\mathcal{R}}(z) 
\hat{\mathcal{X}}^{\hat{b}}(z) \left(\mathrm{I} \! - \! \me^{2 \pi \mi ((n-1)K
+k) \int_{z}^{\hat{a}_{N+1}} \psi_{\widetilde{V}}^{\infty}(\xi) \, \md \xi} 
\sigma_{-} \right) \me^{n(g^{\infty}(z)-\tilde{\mathscr{P}}_{0}) \sigma_{3}}, 
&\text{$z \! \in \! \hat{\Omega}_{\hat{b}_{j-1}}^{3}$,} \\
\me^{\frac{n \hat{\ell}}{2} \operatorname{ad}(\sigma_{3})} \hat{\mathcal{R}}(z) 
\hat{\mathcal{X}}^{\hat{a}}(z) \left(\mathrm{I} \! - \! \me^{2 \pi \mi ((n-1)K
+k) \int_{z}^{\hat{a}_{N+1}} \psi_{\widetilde{V}}^{\infty}(\xi) \, \md \xi} 
\sigma_{-} \right) \me^{n(g^{\infty}(z)-\tilde{\mathscr{P}}_{0}) \sigma_{3}}, 
&\text{$z \! \in \! \hat{\Omega}_{\hat{a}_{j}}^{3}$,} \\
\me^{\frac{n \hat{\ell}}{2} \operatorname{ad}(\sigma_{3})} \hat{\mathcal{R}}(z) 
\hat{\mathcal{X}}^{\hat{b}}(z) \me^{n(g^{\infty}(z)-\tilde{\mathscr{P}}_{0}) 
\sigma_{3}}, &\text{$z \! \in \! \hat{\Omega}_{\hat{b}_{j-1}}^{4}$,} \\
\me^{\frac{n \hat{\ell}}{2} \operatorname{ad}(\sigma_{3})} \hat{\mathcal{R}}(z) 
\hat{\mathcal{X}}^{\hat{a}}(z) \me^{n(g^{\infty}(z)-\tilde{\mathscr{P}}_{0}) 
\sigma_{3}}, &\text{$z \! \in \! \hat{\Omega}_{\hat{a}_{j}}^{4}$;}
\end{cases}
\end{equation}
and (2) for $n \! \in \! \mathbb{N}$ and $k \! \in \! \lbrace 1,2,\dotsc,K \rbrace$ 
such that $\alpha_{p_{\mathfrak{s}}} \! := \! \alpha_{k} \! \neq \! \infty$, and 
$j \! = \! 1,2,\dotsc,N \! + \! 1$,
\begin{equation} \label{finmpcorfmonictil} 
\mathcal{X}(z) \! = \! 
\begin{cases}
\me^{\frac{n \tilde{\ell}}{2} \operatorname{ad}(\sigma_{3})} \tilde{\mathcal{R}}
(z) \widetilde{\mathbb{K}} \, \tilde{\mathbb{M}}(z) \mathscr{E}^{\sigma_{3}} 
\me^{n(g^{f}(z)-\hat{\mathscr{P}}_{0}^{+}) \sigma_{3}}, &\text{$z \! \in \! 
\tilde{\Upsilon}_{1}$,} \\
-\mi \me^{\frac{n \tilde{\ell}}{2} \operatorname{ad}(\sigma_{3})} 
\tilde{\mathcal{R}}(z) \widetilde{\mathbb{K}} \, \tilde{\mathbb{M}}(z) 
\sigma_{2}\mathscr{E}^{-\sigma_{3}} \me^{n(g^{f}(z)-\hat{\mathscr{P}}_{0}^{-}) 
\sigma_{3}}, &\text{$z \! \in \! \tilde{\Upsilon}_{2}$,} \\
\me^{\frac{n \tilde{\ell}}{2} \operatorname{ad}(\sigma_{3})} \tilde{\mathcal{R}}
(z) \widetilde{\mathbb{K}} \, \tilde{\mathbb{M}}(z) \left(\mathrm{I} \! + \! 
\me^{-2 \pi \mi ((n-1)K+k) \int_{z}^{\tilde{a}_{N+1}} \psi_{\widetilde{V}}^{f}
(\xi) \, \md \xi} \sigma_{-} \right) \mathscr{E}^{\sigma_{3}} \me^{n(g^{f}(z)-
\hat{\mathscr{P}}_{0}^{+}) \sigma_{3}}, &\text{$z \! \in \! \tilde{\Upsilon}_{3}$,} \\
-\mi \me^{\frac{n \tilde{\ell}}{2} \operatorname{ad}(\sigma_{3})} 
\tilde{\mathcal{R}}(z) \widetilde{\mathbb{K}} \, \tilde{\mathbb{M}}(z) 
\sigma_{2} \left(\mathrm{I} \! - \! \me^{2 \pi \mi ((n-1)K+k) \int_{z}^{
\tilde{a}_{N+1}} \psi_{\widetilde{V}}^{f}(\xi) \, \md \xi} \sigma_{-} \right) 
\mathscr{E}^{-\sigma_{3}} \me^{n(g^{f}(z)-\hat{\mathscr{P}}_{0}^{-}) 
\sigma_{3}}, &\text{$z \! \in \! \tilde{\Upsilon}_{4}$,} \\
\me^{\frac{n \tilde{\ell}}{2} \operatorname{ad}(\sigma_{3})} \tilde{\mathcal{R}}
(z) \tilde{\mathcal{X}}^{\tilde{b}}(z) \mathscr{E}^{\sigma_{3}} \me^{n(g^{f}(z)
-\hat{\mathscr{P}}_{0}^{+}) \sigma_{3}}, &\text{$z \! \in \! \tilde{\Omega}_{
\tilde{b}_{j-1}}^{1}$,} \\
\me^{\frac{n \tilde{\ell}}{2} \operatorname{ad}(\sigma_{3})} \tilde{\mathcal{R}}
(z) \tilde{\mathcal{X}}^{\tilde{a}}(z) \mathscr{E}^{\sigma_{3}} \me^{n(g^{f}(z)
-\hat{\mathscr{P}}_{0}^{+}) \sigma_{3}}, &\text{$z \! \in \! \tilde{\Omega}_{
\tilde{a}_{j}}^{1}$,} \\
\me^{\frac{n \tilde{\ell}}{2} \operatorname{ad}(\sigma_{3})} \tilde{\mathcal{R}}
(z) \tilde{\mathcal{X}}^{\tilde{b}}(z) \left(\mathrm{I} \! + \! \me^{-2 \pi \mi 
((n-1)K+k) \int_{z}^{\tilde{a}_{N+1}} \psi_{\widetilde{V}}^{f}(\xi) \, \md \xi} 
\sigma_{-} \right) \mathscr{E}^{\sigma_{3}} \me^{n(g^{f}(z)-\hat{\mathscr{
P}}_{0}^{+}) \sigma_{3}}, &\text{$z \! \in \! \tilde{\Omega}_{\tilde{b}_{j-1}}^{2}$,} \\
\me^{\frac{n \tilde{\ell}}{2} \operatorname{ad}(\sigma_{3})} \tilde{\mathcal{R}}
(z) \tilde{\mathcal{X}}^{\tilde{a}}(z) \left(\mathrm{I} \! + \! \me^{-2 \pi \mi 
((n-1)K+k) \int_{z}^{\tilde{a}_{N+1}} \psi_{\widetilde{V}}^{f}(\xi) \, \md \xi} 
\sigma_{-} \right) \mathscr{E}^{\sigma_{3}} \me^{n(g^{f}(z)-\hat{\mathscr{
P}}_{0}^{+}) \sigma_{3}}, &\text{$z \! \in \! \tilde{\Omega}_{\tilde{a}_{j}}^{2}$,} \\
\me^{\frac{n \tilde{\ell}}{2} \operatorname{ad}(\sigma_{3})} \tilde{\mathcal{R}}
(z) \tilde{\mathcal{X}}^{\tilde{b}}(z) \left(\mathrm{I} \! - \! \me^{2 \pi \mi 
((n-1)K+k) \int_{z}^{\tilde{a}_{N+1}} \psi_{\widetilde{V}}^{f}(\xi) \, \md \xi} 
\sigma_{-} \right) \mathscr{E}^{-\sigma_{3}} \me^{n(g^{f}(z)-\hat{\mathscr{
P}}_{0}^{-}) \sigma_{3}}, &\text{$z \! \in \! \tilde{\Omega}_{\tilde{b}_{j-1}}^{3}$,} \\
\me^{\frac{n \tilde{\ell}}{2} \operatorname{ad}(\sigma_{3})} \tilde{\mathcal{R}}
(z) \tilde{\mathcal{X}}^{\tilde{a}}(z) \left(\mathrm{I} \! - \! \me^{2 \pi \mi 
((n-1)K+k) \int_{z}^{\tilde{a}_{N+1}} \psi_{\widetilde{V}}^{f}(\xi) \, \md \xi} 
\sigma_{-} \right) \mathscr{E}^{-\sigma_{3}} \me^{n(g^{f}(z)-\hat{\mathscr{
P}}_{0}^{-}) \sigma_{3}}, &\text{$z \! \in \! \tilde{\Omega}_{\tilde{a}_{j}}^{3}$,} \\
\me^{\frac{n \tilde{\ell}}{2} \operatorname{ad}(\sigma_{3})} \tilde{\mathcal{R}}
(z) \tilde{\mathcal{X}}^{\tilde{b}}(z) \mathscr{E}^{-\sigma_{3}} \me^{n(g^{f}(z)
-\hat{\mathscr{P}}_{0}^{-}) \sigma_{3}}, &\text{$z \! \in \! \tilde{\Omega}_{
\tilde{b}_{j-1}}^{4}$,} \\
\me^{\frac{n \tilde{\ell}}{2} \operatorname{ad}(\sigma_{3})} \tilde{\mathcal{R}}
(z) \tilde{\mathcal{X}}^{\tilde{a}}(z) \mathscr{E}^{-\sigma_{3}} \me^{n(g^{f}(z)
-\hat{\mathscr{P}}_{0}^{-}) \sigma_{3}}, &\text{$z \! \in \! \tilde{\Omega}_{
\tilde{a}_{j}}^{4}$.}
\end{cases}
\end{equation}
For $n \! \in \! \mathbb{N}$ and $k \! \in \! \lbrace 1,2,\dotsc,K \rbrace$ 
such that $\alpha_{p_{\mathfrak{s}}} \! := \! \alpha_{k} \! = \! \infty$ 
(resp., $\alpha_{p_{\mathfrak{s}}} \! := \! \alpha_{k} \! \neq \! \infty)$, 
via a matrix-multiplication argument and collecting $(1 \, 1)$- and 
$(1 \, 2)$-elements, respectively, in Equation~\eqref{finmpcorfmonichat} 
(resp., Equation~\eqref{finmpcorfmonictil}), one arrives at the asymptotic, 
in the double-scaling limit $\mathscr{N},n \! \to \! \infty$ such 
that $z_{o} \! = \! 1 \! + \! o(1)$, results for the monic MPC ORF, 
$\pmb{\pi}_{k}^{n}(z) \! = \! (\mathcal{X}(z))_{11}$, and $\int_{\mathbb{R}} 
\tfrac{\pmb{\pi}_{k}^{n}(\xi) \me^{-n \widetilde{V}(\xi)}}{\xi -z} \, 
\tfrac{\md \xi}{2 \pi \mi} \! = \! (\mathcal{X}(z))_{12}$ (resp., 
$\pmb{\pi}_{k}^{n}(z) \! = \! (z \! - \! \alpha_{k})^{-1}(\mathcal{X}(z))_{11}$, 
and $\int_{\mathbb{R}} \tfrac{((\xi -\alpha_{k}) \pmb{\pi}_{k}^{n}(\xi)) 
\me^{-n \widetilde{V}(\xi)}}{(\xi -\alpha_{k})(\xi -z)} \, \tfrac{\md \xi}{2 
\pi \mi} \! = \! (z \! - \! \alpha_{k})^{-1}(\mathcal{X}(z))_{12})$, in the 
entire complex plane, stated in Theorem~\ref{maintheoforinf1} (resp., 
\ref{maintheoforfin1}); furthermore, for $n \! \in \! \mathbb{N}$ and 
$k \! \in \! \lbrace 1,2,\dotsc,K \rbrace$ such that $\alpha_{
p_{\mathfrak{s}}} \! := \! \alpha_{k} \! = \! \infty$ (resp., 
$\alpha_{p_{\mathfrak{s}}} \! := \! \alpha_{k} \! \neq \! \infty)$, via 
Equation~\eqref{eqlemmvssinfmpa3} (resp., \eqref{eqlemmvssfinmpa7}), 
one arrives at the asymptotic, in the double-scaling limit $\mathscr{N},n 
\! \to \! \infty$ such that $z_{o} \! = \! 1 \! + \! o(1)$, result for the 
corresponding MPA error term $\widehat{\pmb{\mathrm{E}}}_{\tilde{\mu}}(z)$ 
(resp., $\widetilde{\pmb{\mathrm{E}}}_{\tilde{\mu}}(z))$, in the entire 
complex plane, stated in Theorem~\ref{maintheompainf} (resp., 
\ref{maintheompafin}).

For $n \! \in \! \mathbb{N}$ and $k \! \in \! \lbrace 1,2,\dotsc,K \rbrace$ 
such that $\alpha_{p_{\mathfrak{s}}} \! := \! \alpha_{k} \! = \! \infty$ 
(resp., $\alpha_{p_{\mathfrak{s}}} \! := \! \alpha_{k} \! \neq \! \infty)$, 
via Equation~\eqref{eqlem5.6A} (resp., \eqref{eqlem5.6B}), one arrives at 
the asymptotic, in the double-scaling limit $\mathscr{N},n \! \to \! \infty$ 
such that $z_{o} \! = \! 1 \! + \! o(1)$, result for the associated norming 
constant, $\mu_{n,\varkappa_{nk}}^{\infty}(n,k)$ (resp., 
$\mu_{n,\varkappa_{nk}}^{f}(n,k))$ stated in Theorem~\ref{maintheoforinf2} 
(resp., \ref{maintheoforfin2}); moreover, via the asymptotics for 
$\pmb{\pi}_{k}^{n}(z)$ and $\mu_{n,\varkappa_{nk}}^{\infty}(n,k)$ (resp., 
$\pmb{\pi}_{k}^{n}(z)$ and $\mu_{n,\varkappa_{nk}}^{f}(n,k))$ stated in 
Theorems~\ref{maintheoforinf1} and~\ref{maintheoforinf2} (resp., 
\ref{maintheoforfin1} and~\ref{maintheoforfin2}) and the relation $\phi_{k}^{n}
(z) \! = \! \mu_{n,\varkappa_{nk}}^{\infty}(n,k) \pmb{\pi}_{k}^{n}(z)$ (resp., 
$\phi_{k}^{n}(z) \! = \! \mu_{n,\varkappa_{nk}}^{f}(n,k) \pmb{\pi}_{k}^{n}(z))$, 
one arrives at the asymptotic, in the double-scaling limit $\mathscr{N},n 
\! \to \! \infty$ such that $z_{o} \! = \! 1 \! + \! o(1)$, results for the 
MPC ORF, $\phi_{k}^{n}(z)$, in the entire complex plane, stated in 
Theorem~\ref{maintheoforinf2} (resp., \ref{maintheoforfin2}).
\clearpage
\section*{Appendix A: MPC ORF Moments} \label{appamom}
\addcontentsline{toc}{section}{Appendix A: MPC ORF Moments}
\setcounter{section}{1}
\setcounter{equation}{0}
\renewcommand{\thesection}{\Alph{section}}
\renewcommand{\theequation}{\Alph{section}.\arabic{equation}}
For $n \! \in \! \mathbb{N}$ and $k \! \in \! \lbrace 1,2,\dotsc,K \rbrace$ 
such that $\alpha_{p_{\mathfrak{s}}} \! := \! \alpha_{k} \! = \! \infty$ 
(resp., $\alpha_{p_{\mathfrak{s}}} \! := \! \alpha_{k} \! \neq \! \infty)$, 
explicit formulae for the various moments of the associated equilibrium 
measure $\mu_{\widetilde{V}}^{\infty}$ (resp., $\mu_{\widetilde{V}}^{f})$, 
that is, $\int_{\mathbb{R}} \xi^{m} \, \md \mu_{\widetilde{V}}^{\infty}
(\xi)$ and $\int_{\mathbb{R}}(\xi \! - \! \alpha_{p_{q}})^{-m} \, \md 
\mu_{\widetilde{V}}^{\infty}(\xi)$, $m \! \in \! \mathbb{N}$, $q \! = \! 
1,2,\dotsc,\mathfrak{s} \! - \! 1$ (resp., $\int_{\mathbb{R}} \xi^{m} \, 
\md \mu_{\widetilde{V}}^{f}(\xi)$ and $\int_{\mathbb{R}}(\xi \! - \! 
\alpha_{p_{q}})^{-m} \, \md \mu_{\widetilde{V}}^{f}(\xi)$, $m \! \in \! 
\mathbb{N}$, $q \! = \! 1,\dotsc,\mathfrak{s} \! - \! 2,\mathfrak{s})$, in 
terms of the regular external field $\widetilde{V} \colon \overline{\mathbb{R}} 
\setminus \lbrace \alpha_{1},\alpha_{2},\dotsc,\alpha_{K} \rbrace \! \to \! 
\mathbb{R}$ satisfying conditions~\eqref{eq20}--\eqref{eq22}, and the 
multi-valued function $(\hat{R}(z))^{1/2}$ (resp., $(\tilde{R}(z))^{1/2})$ 
defined by Equation~\eqref{eql3.7d} (resp., \eqref{eql3.7j}), where 
$\lbrace \hat{b}_{j-1},\hat{a}_{j} \rbrace_{j=1}^{N+1}$ (resp., $\lbrace 
\tilde{b}_{j-1},\tilde{a}_{j} \rbrace_{j=1}^{N+1})$ satisfy, in the 
double-scaling limit $\mathscr{N},n \! \to \! \infty$ such that $z_{o} \! = 
\! 1 \! + \! o(1)$, the locally solvable system of $2(N \! + \! 1)$ moment 
equations~\eqref{eql3.7a}--\eqref{eql3.7c} (resp., \eqref{eql3.7g}--\eqref{eql3.7i}), 
can be derived. Without loss of generality, consider the case $n \! \in \! 
\mathbb{N}$ and $k \! \in \! \lbrace 1,2,\dotsc,K \rbrace$ such that 
$\alpha_{p_{\mathfrak{s}}} \! := \! \alpha_{k} \! \neq \! \infty$, and calculate, 
say, the corresponding `first three' moments of the associated equilibrium 
measure, $\mu_{\widetilde{V}}^{f}$, that is, $\int_{J_{f}} \xi^{r_{1}} 
\psi_{\widetilde{V}}^{f}(\xi) \, \md \xi$ and $\int_{J_{f}}(\xi \! - \! 
\alpha_{p_{q}})^{-r_{1}} \psi_{\widetilde{V}}^{f}(\xi) \, \md \xi$, $r_{1} \! 
= \! 1,2,3$, $q \! = \! 1,\dotsc,\mathfrak{s} \! - \! 2,\mathfrak{s}$. (Of 
course, the scheme of the calculations below generalises to the higher-order 
moments $\int_{J_{f}} \xi^{r_{2}+3} \psi_{\widetilde{V}}^{f}(\xi) \, \md \xi$ 
and $\int_{J_{f}}(\xi \! - \! \alpha_{p_{q}})^{-(r_{2}+3)} \psi_{\widetilde{V}}^{
f}(\xi) \, \md \xi$, $r_{2} \! \in \! \mathbb{N}$, $q \! = \! 1,\dotsc,
\mathfrak{s} \! - \! 2,\mathfrak{s}$.) For $n \! \in \! \mathbb{N}$ and $k 
\! \in \! \lbrace 1,2,\dotsc,K \rbrace$ such that $\alpha_{p_{\mathfrak{s}}} 
\! := \! \alpha_{k} \! \neq \! \infty$, recall the representations for 
$\mathfrak{F}_{f}(z)$ given in Equations~\eqref{eql3.4a} and~\eqref{eql3.7a12}. 
{}From these representations for $\mathfrak{F}_{f}(z)$, one derives the following 
asymptotic expansions: (i) for $\xi \! \in \! J_{f}$ and $z \! \notin \! J_{f}$ 
such that $\lvert (z \! - \! \alpha_{k})/(\xi \! - \! \alpha_{k}) \rvert \! \ll \! 1$ 
(e.g., $0 \! < \! \lvert z \! - \! \alpha_{k} \rvert \! \ll \! \min \lbrace 
\min_{i \neq j \in \lbrace 1,\dotsc,\mathfrak{s}-2,\mathfrak{s} \rbrace} 
\lbrace \lvert \alpha_{p_{i}} \! - \! \alpha_{p_{j}} \rvert \rbrace,\inf_{\xi 
\in J_{f}} \lbrace \lvert \xi \! - \! \alpha_{k} \rvert \rbrace,\min_{j=1,2,
\dotsc,N+1} \lbrace \lvert \lvert \tilde{b}_{j-1} \! - \! \tilde{a}_{j} \rvert \! - \! 
\alpha_{k} \rvert \rbrace \rbrace)$, via the expansions $\tfrac{1}{(z-\alpha_{k})
-(\xi -\alpha_{k})} \! = \! -\sum_{j=0}^{l} \tfrac{(z-\alpha_{k})^{j}}{(\xi 
-\alpha_{k})^{j+1}} \! + \! \tfrac{(z-\alpha_{k})^{l+1}}{(\xi -\alpha_{k})^{l+1}
(z- \xi)}$, $l \! \in \! \mathbb{N}_{0}$, and $\ln (1 \! - \! \blacklozenge) \! 
= \! -\sum_{m=1}^{\infty} \blacklozenge^{m}/m$, $\lvert \blacklozenge 
\rvert \! \ll \! 1$,
\begin{align} \label{eqappa1} 
\mathfrak{F}_{f}(z) \underset{z \to \alpha_{k}}{=}& \, -\dfrac{1}{\mi \pi} 
\left(\dfrac{(\varkappa_{nk} \! - \! 1)}{n(z \! - \! \alpha_{k})} \! + \! 
\sum_{j=0}^{3} \left(-\sum_{q=1}^{\mathfrak{s}-2} \dfrac{\varkappa_{nk 
\tilde{k}_{q}}}{n(\alpha_{p_{q}} \! - \! \alpha_{k})^{j+1}} \! + \! 
\left(\dfrac{(n \! - \! 1)K \! + \! k}{n} \right) \int_{J_{f}} 
\dfrac{\psi_{\widetilde{V}}^{f}(\xi)}{(\xi \! - \! \alpha_{k})^{j+1}} \, 
\md \xi \right)(z \! - \! \alpha_{k})^{j} \right) \nonumber \\
+& \, \mathcal{O} \left((z \! - \! \alpha_{k})^{4} \right),
\end{align}
and
\begin{align} \label{eqappa2} 
\mathfrak{F}_{f}(z) \underset{z \to \alpha_{k}}{=}& \, -\dfrac{1}{\mi \pi} 
\dfrac{(\varkappa_{nk} \! - \! 1)}{n(z \! - \! \alpha_{k})} \! + \! 
\dfrac{1}{\mi \pi} \sum_{q=1}^{\mathfrak{s}-2} \dfrac{\varkappa_{nk 
\tilde{k}_{q}}}{n(\alpha_{p_{q}} \! - \! \alpha_{k})} \! + \! 
(-1)^{\tilde{\mathfrak{n}}(\alpha_{k})} \left(\prod_{j=1}^{N+1} \lvert 
\tilde{b}_{j-1} \! - \! \alpha_{k} \rvert \lvert \tilde{a}_{j} \! - \! 
\alpha_{k} \rvert \right)^{1/2} \tilde{\mathscr{Q}}_{0}^{\sharp} 
\nonumber \\
+& \, \left(\dfrac{1}{\mi \pi} \sum_{q=1}^{\mathfrak{s}-2} 
\dfrac{\varkappa_{nk \tilde{k}_{q}}}{n(\alpha_{p_{q}} \! - \! \alpha_{k})^{2}} 
\! + \! (-1)^{\tilde{\mathfrak{n}}(\alpha_{k})} \left(\prod_{j=1}^{N+1} 
\lvert \tilde{b}_{j-1} \! - \! \alpha_{k} \rvert \lvert \tilde{a}_{j} \! - \! 
\alpha_{k} \rvert \right)^{1/2} \left(\tilde{\mathscr{Q}}_{1}^{\sharp} \! + 
\! \tilde{\alpha}_{\triangle}^{\sharp} \tilde{\mathscr{Q}}_{0}^{\sharp} 
\right) \right)(z \! - \! \alpha_{k}) \nonumber \\
+& \, \left(\dfrac{1}{\mi \pi} \sum_{q=1}^{\mathfrak{s}-2} 
\dfrac{\varkappa_{nk \tilde{k}_{q}}}{n(\alpha_{p_{q}} \! - \! \alpha_{k})^{3}} 
\! + \! (-1)^{\tilde{\mathfrak{n}}(\alpha_{k})} \left(\prod_{j=1}^{N+1} 
\lvert \tilde{b}_{j-1} \! - \! \alpha_{k} \rvert \lvert \tilde{a}_{j} \! - \! 
\alpha_{k} \rvert \right)^{1/2} \left(\tilde{\mathscr{Q}}_{2}^{\sharp} \! + \! 
\tilde{\alpha}_{\triangle}^{\sharp} \tilde{\mathscr{Q}}_{1}^{\sharp} \! + \! 
\left(\tilde{\beta}_{\triangle}^{\sharp} \! + \! \dfrac{1}{2}(\tilde{\alpha}_{
\triangle}^{\sharp})^{2} \right) \tilde{\mathscr{Q}}_{0}^{\sharp} \right) 
\right)(z \! - \! \alpha_{k})^{2} \nonumber \\
+& \, \mathcal{O} \left((z \! - \! \alpha_{k})^{3} \right),
\end{align}
where $\tilde{\mathfrak{n}}(\alpha_{k})$ is given in the corresponding 
item of Remark~\ref{remext},
\begin{gather*}
\tilde{\mathscr{Q}}_{j}^{\sharp} \! := \! \int_{J_{f}} \dfrac{(\tilde{R}
(\xi))^{-1/2}_{+}}{(\xi \! - \! \alpha_{k})^{j+1}} \left(\dfrac{2}{\mi \pi} 
\left(\dfrac{(\varkappa_{nk} \! - \! 1)}{n(\xi \! - \! \alpha_{k})} \! + \! 
\sum_{q=1}^{\mathfrak{s}-2} \dfrac{\varkappa_{nk \tilde{k}_{q}}}{n(\xi \! - 
\! \alpha_{p_{q}})} \right) \! + \! \dfrac{\widetilde{V}^{\prime}(\xi)}{\mi 
\pi} \right) \dfrac{\md \xi}{2 \pi \mi}, \quad j \! = \! 0,1,2, \\
\tilde{\alpha}_{\triangle}^{\sharp} \! := \! -\dfrac{1}{2} \sum_{j=1}^{N+1} 
\left(\dfrac{1}{\tilde{b}_{j-1} \! - \! \alpha_{k}} \! + \! \dfrac{1}{\tilde{a}_{j} 
\! - \! \alpha_{k}} \right), \qquad \qquad \tilde{\beta}_{\triangle}^{\sharp} 
\! := \! -\dfrac{1}{4} \sum_{j=1}^{N+1} \left(\dfrac{1}{(\tilde{b}_{j-1} \! 
- \! \alpha_{k})^{2}} \! + \! \dfrac{1}{(\tilde{a}_{j} \! - \! \alpha_{k})^{2}} 
\right);
\end{gather*}
and (ii) for $\xi \! \in \! J_{f}$ and $z \! \notin \! J_{f}$ such that 
$\lvert \xi/z \rvert \! \ll \! 1$ (e.g., $\lvert z \rvert \! \gg \! \max 
\lbrace \max_{i \neq j \in \lbrace 1,\dotsc,\mathfrak{s}-2,\mathfrak{s} 
\rbrace} \lbrace \lvert \alpha_{p_{i}} \! - \! \alpha_{p_{j}} \rvert \rbrace,
\max_{q=1,\dotsc,\mathfrak{s}-2,\mathfrak{s}} \lbrace \lvert 
\alpha_{p_{q}} \rvert \rbrace,\linebreak[4]
\max_{j=1,2,\dotsc,N+1} \lbrace \lvert \tilde{b}_{j-1} \! - \! 
\tilde{a}_{j} \rvert \rbrace \rbrace)$, via the expansions $\tfrac{1}{\xi -z} \! 
= \! -\sum_{j=0}^{l} \tfrac{\xi^{j}}{z^{j+1}} \! + \! \tfrac{\xi^{l+1}}{z^{l+1}
(\xi -z)}$, $l \! \in \! \mathbb{N}_{0}$, and $\ln (1 \! - \! \blacklozenge) \! 
= \! -\sum_{m=1}^{\infty} \blacklozenge^{m}/m$, $\lvert \blacklozenge 
\rvert \! \ll \! 1$, and the associated $2(N \! + \! 1)$ moment 
equations~\eqref{eql3.7g}--\eqref{eql3.7i},
\begin{align} \label{eqappa3} 
\mathfrak{F}_{f}(z) \underset{z \to \alpha_{p_{\mathfrak{s}-1}} = \infty}{=}& \, 
\dfrac{(\varkappa_{nk \tilde{k}_{\mathfrak{s}-1}}^{\infty} \! + \! 1)}{\mi \pi 
nz} \! - \! \dfrac{1}{\mi \pi z} \sum_{m=1}^{3} \left(\left(\dfrac{\varkappa_{nk} 
\! - \! 1}{n} \right)(\alpha_{k})^{m} \! + \! \sum_{q=1}^{\mathfrak{s}-2} 
\dfrac{\varkappa_{nk \tilde{k}_{q}}}{n}(\alpha_{p_{q}})^{m} \! - \! 
\left(\dfrac{(n \! - \! 1)K \! + \! k}{n} \right) \int_{J_{f}} \xi^{m} 
\psi_{\widetilde{V}}^{f}(\xi) \, \md \xi \right)z^{-m} \nonumber \\
+& \, \mathcal{O}(z^{-5}),
\end{align}
and
\begin{align} \label{eqappa4} 
\mathfrak{F}_{f}(z) \underset{z \to \alpha_{p_{\mathfrak{s}-1}} = \infty}{=}& \, 
\dfrac{(\varkappa_{nk \tilde{k}_{\mathfrak{s}-1}}^{\infty} \! + \! 1)}{\mi \pi nz} 
\! + \! \left(-\tilde{\mathscr{Q}}_{0}^{\flat} \! + \! \dfrac{\tilde{\alpha}_{
\triangle}^{\flat}}{\mi \pi} \left(\dfrac{(n \! - \! 1)K \! + \! k}{n} \right) \! - \! 
\dfrac{1}{\mi \pi} \left(\left(\dfrac{\varkappa_{nk} \! - \! 1}{n} \right) \alpha_{k} 
\! + \! \sum_{q=1}^{\mathfrak{s}-2} \dfrac{\varkappa_{nk \tilde{k}_{q}}}{n} 
\alpha_{p_{q}} \right) \right)z^{-2} \nonumber \\
+& \, \left(-\tilde{\mathscr{Q}}_{1}^{\flat} \! - \! \tilde{\alpha}_{\triangle}^{
\flat} \tilde{\mathscr{Q}}_{0}^{\flat} \! + \! \dfrac{(\tilde{\beta}_{\triangle}^{
\flat} \! + \! (\tilde{\alpha}_{\triangle}^{\flat})^{2}/2)}{\mi \pi} \left(\dfrac{(n 
\! - \! 1)K \! + \! k}{n} \right) \! - \! \dfrac{1}{\mi \pi} \left(\left(\dfrac{
\varkappa_{nk} \! - \! 1}{n} \right)(\alpha_{k})^{2} \! + \! \sum_{q=1}^{
\mathfrak{s}-2} \dfrac{\varkappa_{nk \tilde{k}_{q}}}{n}(\alpha_{p_{q}})^{2} 
\right) \right)z^{-3} \nonumber \\
+& \, \left(-\tilde{\mathscr{Q}}_{2}^{\flat} \! - \! \tilde{\alpha}_{\triangle}^{
\flat} \tilde{\mathscr{Q}}_{1}^{\flat} \! - \! \left(\tilde{\beta}_{\triangle}^{
\flat} \! + \! \dfrac{1}{2}(\tilde{\alpha}_{\triangle}^{\flat})^{2} \right) 
\tilde{\mathscr{Q}}_{0}^{\flat} \! + \! \dfrac{(\tilde{\gamma}_{\triangle}^{
\flat} \! + \! \tilde{\alpha}_{\triangle}^{\flat} \tilde{\beta}_{\triangle}^{\flat} 
\! + \! (\tilde{\alpha}_{\triangle}^{\flat})^{3}/3!)}{\mi \pi} \left(\dfrac{(n 
\! - \! 1)K \! + \! k}{n} \right) \! - \! \dfrac{1}{\mi \pi} \left(\left(
\dfrac{\varkappa_{nk} \! - \! 1}{n} \right)(\alpha_{k})^{3} \right. \right. 
\nonumber \\
+&\left. \left. \, \sum_{q=1}^{\mathfrak{s}-2} \dfrac{\varkappa_{nk 
\tilde{k}_{q}}}{n}(\alpha_{p_{q}})^{3} \right) \right)z^{-4} \! + \! 
\mathcal{O}(z^{-5}),
\end{align}
where
\begin{gather*}
\tilde{\mathscr{Q}}_{j}^{\flat} \! := \! \int_{J_{f}} \dfrac{\xi^{N+2+j}}{
(\tilde{R}(\xi))^{1/2}_{+}} \left(\dfrac{2}{\mi \pi} \left(\dfrac{
(\varkappa_{nk} \! - \! 1)}{n(\xi \! - \! \alpha_{k})} \! + \! \sum_{q=1}^{
\mathfrak{s}-2} \dfrac{\varkappa_{nk \tilde{k}_{q}}}{n(\xi \! - \! 
\alpha_{p_{q}})} \right) \! + \! \dfrac{\widetilde{V}^{\prime}(\xi)}{\mi 
\pi} \right) \dfrac{\md \xi}{2 \pi \mi}, \quad j \! = \! 0,1,2, \\
\tilde{\alpha}_{\triangle}^{\flat} \! := \! -\dfrac{1}{2} \sum_{j=1}^{N+1}
(\tilde{b}_{j-1} \! + \! \tilde{a}_{j}), \quad \quad \tilde{\beta}_{\triangle}^{
\flat} \! := \! -\dfrac{1}{4} \sum_{j=1}^{N+1}((\tilde{b}_{j-1})^{2} \! + \! 
(\tilde{a}_{j})^{2}), \quad \quad \tilde{\gamma}_{\triangle}^{\flat} \! := \! 
-\dfrac{1}{3!} \sum_{j=1}^{N+1}((\tilde{b}_{j-1})^{3} \! + \! (\tilde{a}_{j})^{3}).
\end{gather*}
For $n \! \in \! \mathbb{N}$ and $k \! \in \! \lbrace 1,2,\dotsc,K \rbrace$ 
such that $\alpha_{p_{\mathfrak{s}}} \! := \! \alpha_{k} \! \neq \! \infty$, 
equating coefficients of like powers of $(z \! - \! \alpha_{k})^{i}$, $i \! = \! -1,
0,1,2$, in the asymptotic expansions~\eqref{eqappa1} and~\eqref{eqappa2} 
for $\mathfrak{F}_{f}(z)$, one arrives at
\begin{align*}
\int_{J_{f}} \dfrac{\psi_{\widetilde{V}}^{f}(\xi)}{\xi \! - \! \alpha_{k}} \, 
\md \xi =& \, \dfrac{(-1)^{\tilde{\mathfrak{n}}(\alpha_{k})}(\prod_{j=1}^{N+1} 
\lvert \tilde{b}_{j-1} \! - \! \alpha_{k} \rvert \lvert \tilde{a}_{j} \! - \! 
\alpha_{k} \rvert)^{1/2}}{2((n \! - \! 1)K \! + \! k)/n} \int_{J_{f}} 
\dfrac{(\tilde{R}(\xi))^{-1/2}_{+}}{(\xi \! - \! \alpha_{k})} \left(\dfrac{2 
\mi}{\pi} \left(\dfrac{(\varkappa_{nk} \! - \! 1)}{n(\xi \! - \! \alpha_{k})} 
\! + \! \sum_{q=1}^{\mathfrak{s}-2} \dfrac{\varkappa_{nk \tilde{k}_{q}}}{n(\xi 
\! - \! \alpha_{p_{q}})} \right) \right. \\
+&\left. \, \dfrac{\mi \widetilde{V}^{\prime}(\xi)}{\pi} \right) \md \xi, \\
\int_{J_{f}} \dfrac{\psi_{\widetilde{V}}^{f}(\xi)}{(\xi \! - \! \alpha_{k})^{
2}} \, \md \xi =& \, \dfrac{(-1)^{\tilde{\mathfrak{n}}(\alpha_{k})}(\prod_{j=
1}^{N+1} \lvert \tilde{b}_{j-1} \! - \! \alpha_{k} \rvert \lvert \tilde{a}_{j} 
\! - \! \alpha_{k} \rvert)^{1/2}}{2((n \! - \! 1)K \! + \! k)/n} \left(
\int_{J_{f}} \dfrac{(\tilde{R}(\xi))^{-1/2}_{+}}{(\xi \! - \! \alpha_{k})^{2}} 
\left(\dfrac{2 \mi}{\pi} \left(\dfrac{(\varkappa_{nk} \! - \! 1)}{n(\xi \! - 
\! \alpha_{k})} \! + \! \sum_{q=1}^{\mathfrak{s}-2} \dfrac{\varkappa_{nk 
\tilde{k}_{q}}}{n(\xi \! - \! \alpha_{p_{q}})} \right) \right. \right. \\
+&\left. \left. \, \dfrac{\mi \widetilde{V}^{\prime}(\xi)}{\pi} \right) \md 
\xi \! - \! \int_{J_{f}} \dfrac{(\tilde{R}(\xi))^{-1/2}_{+}}{(\xi \! - \! 
\alpha_{k})} \left(\dfrac{2 \mi}{\pi} \left(\dfrac{(\varkappa_{nk} \! - \! 
1)}{n(\xi \! - \! \alpha_{k})} \! + \! \sum_{q=1}^{\mathfrak{s}-2} 
\dfrac{\varkappa_{nk \tilde{k}_{q}}}{n(\xi \! - \! \alpha_{p_{q}})} \right) \! 
+ \! \dfrac{\mi \widetilde{V}^{\prime}(\xi)}{\pi} \right) \md \xi \right. \\
\times&\left. \, \dfrac{1}{2} \sum_{m=1}^{N+1} \left(\dfrac{1}{\tilde{b}_{m-1} 
\! - \! \alpha_{k}} \! + \! \dfrac{1}{\tilde{a}_{m} \! - \! \alpha_{k}} 
\right) \right), \\
\int_{J_{f}} \dfrac{\psi_{\widetilde{V}}^{f}(\xi)}{(\xi \! - \! \alpha_{k})^{
3}} \, \md \xi =& \, \dfrac{(-1)^{\tilde{\mathfrak{n}}(\alpha_{k})}(\prod_{j=
1}^{N+1} \lvert \tilde{b}_{j-1} \! - \! \alpha_{k} \rvert \lvert \tilde{a}_{j} 
\! - \! \alpha_{k} \rvert)^{1/2}}{2((n \! - \! 1)K \! + \! k)/n} \left(
\int_{J_{f}} \dfrac{(\tilde{R}(\xi))^{-1/2}_{+}}{(\xi \! - \! \alpha_{k})^{3}} 
\left(\dfrac{2 \mi}{\pi} \left(\dfrac{(\varkappa_{nk} \! - \! 1)}{n(\xi \! - 
\! \alpha_{k})} \! + \! \sum_{q=1}^{\mathfrak{s}-2} \dfrac{\varkappa_{nk 
\tilde{k}_{q}}}{n(\xi \! - \! \alpha_{p_{q}})} \right) \right. \right. \\
+&\left. \left. \, \dfrac{\mi \widetilde{V}^{\prime}(\xi)}{\pi} \right) \md 
\xi \! - \! \int_{J_{f}} \dfrac{(\tilde{R}(\xi))^{-1/2}_{+}}{(\xi \! - \! 
\alpha_{k})^{2}} \left(\dfrac{2 \mi}{\pi} \left(\dfrac{(\varkappa_{nk} \! 
- \! 1)}{n(\xi \! - \! \alpha_{k})} \! + \! \sum_{q=1}^{\mathfrak{s}-2} 
\dfrac{\varkappa_{nk \tilde{k}_{q}}}{n(\xi \! - \! \alpha_{p_{q}})} \right) \! 
+ \! \dfrac{\mi \widetilde{V}^{\prime}(\xi)}{\pi} \right) \md \xi \right. \\
\times&\left. \, \dfrac{1}{2} \sum_{m=1}^{N+1} \left(\dfrac{1}{\tilde{b}_{m-1} 
\! - \! \alpha_{k}} \! + \! \dfrac{1}{\tilde{a}_{m} \! - \! \alpha_{k}} 
\right) \! + \! \int_{J_{f}} \dfrac{(\tilde{R}(\xi))^{-1/2}_{+}}{(\xi \! 
- \! \alpha_{k})} \left(\dfrac{2 \mi}{\pi} \left(\dfrac{(\varkappa_{nk} 
\! - \! 1)}{n(\xi \! - \! \alpha_{k})} \! + \! \sum_{q=1}^{\mathfrak{s}-2} 
\dfrac{\varkappa_{nk \tilde{k}_{q}}}{n(\xi \! - \! \alpha_{p_{q}})} \right) \! 
+ \! \dfrac{\mi \widetilde{V}^{\prime}(\xi)}{\pi} \right) \md \xi \right. \\
\times&\left. \, \left(\dfrac{1}{8} \left(\sum_{m=1}^{N+1} \left(\dfrac{1}{
\tilde{b}_{m-1} \! - \! \alpha_{k}} \! + \! \dfrac{1}{\tilde{a}_{m} \! - \! 
\alpha_{k}} \right) \right)^{2} \! - \! \dfrac{1}{4} \sum_{m=1}^{N+1} \left(
\dfrac{1}{(\tilde{b}_{m-1} \! - \! \alpha_{k})^{2}} \! + \! \dfrac{1}{
(\tilde{a}_{m} \! - \! \alpha_{k})^{2}} \right) \right) \right),
\end{align*}
and, equating coefficients of like powers of $z^{-j}$, $j \! = \! 1,2,3,4$, 
in the asymptotic expansions~\eqref{eqappa3} and~\eqref{eqappa4} for 
$\mathfrak{F}_{f}(z)$, one arrives at
\begin{align*}
\int_{J_{f}} \xi \psi_{\widetilde{V}}^{f}(\xi) \, \md \xi =& \, \dfrac{1}{2
((n \! - \! 1)K \! + \! k)/n} \int_{J_{f}} \dfrac{\xi^{N+2}}{(\tilde{R}
(\xi))^{1/2}_{+}} \left(\dfrac{2 \mi}{\pi} \left(\dfrac{(\varkappa_{nk} \! 
- \! 1)}{n(\xi \! - \! \alpha_{k})} \! + \! \sum_{q=1}^{\mathfrak{s}-2} 
\dfrac{\varkappa_{nk \tilde{k}_{q}}}{n(\xi \! - \! \alpha_{p_{q}})} \right) 
\! + \! \dfrac{\mi \widetilde{V}^{\prime}(\xi)}{\pi} \right) \md \xi \\
-& \, \dfrac{1}{2} \sum_{j=1}^{N+1}(\tilde{b}_{j-1} \! + \! \tilde{a}_{j}), \\
\int_{J_{f}} \xi^{2} \psi_{\widetilde{V}}^{f}(\xi) \, \md \xi =& \, 
\dfrac{1}{2((n \! - \! 1)K \! + \! k)/n} \left(\int_{J_{f}} 
\dfrac{\xi^{N+3}}{(\tilde{R}(\xi))^{1/2}_{+}} \left(\dfrac{2 \mi}{\pi} \left(
\dfrac{(\varkappa_{nk} \! - \! 1)}{n(\xi \! - \! \alpha_{k})} \! + \! 
\sum_{q=1}^{\mathfrak{s}-2} \dfrac{\varkappa_{nk \tilde{k}_{q}}}{n(\xi \! - \! 
\alpha_{p_{q}})} \right) \! + \! \dfrac{\mi \widetilde{V}^{\prime}(\xi)}{\pi} 
\right) \md \xi \right. \\
+&\left. \, \dfrac{1}{2} \left(\sum_{j=1}^{N+1}(\tilde{b}_{j-1} \! + \! 
\tilde{a}_{j}) \right) \int_{J_{f}} \dfrac{\xi^{N+2}}{(\tilde{R}(\xi))^{
1/2}_{+}} \left(\dfrac{2 \mi}{\pi} \left(\dfrac{(\varkappa_{nk} \! - \! 
1)}{n(\xi \! - \! \alpha_{k})} \! + \! \sum_{q=1}^{\mathfrak{s}-2} \dfrac{
\varkappa_{nk \tilde{k}_{q}}}{n(\xi \! - \! \alpha_{p_{q}})} \right) \! + \! 
\dfrac{\mi \widetilde{V}^{\prime}(\xi)}{\pi} \right) \md \xi \right) \\
-& \, \dfrac{1}{4} \sum_{j=1}^{N+1}((\tilde{b}_{j-1})^{2} \! + \! 
(\tilde{a}_{j})^{2}) \! + \! \dfrac{1}{8} \left(\sum_{j=1}^{N+1}
(\tilde{b}_{j-1} \! + \! \tilde{a}_{j}) \right)^{2}, \\
\int_{J_{f}} \xi^{3} \psi_{\widetilde{V}}^{f}(\xi) \, \md \xi =& \, 
\dfrac{1}{2((n \! - \! 1)K \! + \! k)/n} \left(\int_{J_{f}} 
\dfrac{\xi^{N+4}}{(\tilde{R}(\xi))^{1/2}_{+}} \left(\dfrac{2 \mi}{\pi} 
\left(\dfrac{(\varkappa_{nk} \! - \! 1)}{n(\xi \! - \! \alpha_{k})} \! + \! 
\sum_{q=1}^{\mathfrak{s}-2} \dfrac{\varkappa_{nk \tilde{k}_{q}}}{n(\xi \! - \! 
\alpha_{p_{q}})} \right) \! + \! \dfrac{\mi \widetilde{V}^{\prime}(\xi)}{\pi} 
\right) \md \xi \right. \\
-&\left. \, \dfrac{1}{2} \left(\sum_{j=1}^{N+1}(\tilde{b}_{j-1} \! + \! 
\tilde{a}_{j}) \right) \int_{J_{f}} \dfrac{\xi^{N+3}}{(\tilde{R}(\xi))^{
1/2}_{+}} \left(\dfrac{2 \mi}{\pi} \left(\dfrac{(\varkappa_{nk} \! - \! 
1)}{n(\xi \! - \! \alpha_{k})} \! + \! \sum_{q=1}^{\mathfrak{s}-2} 
\dfrac{\varkappa_{nk \tilde{k}_{q}}}{n(\xi \! - \! \alpha_{p_{q}})} \right) \! 
+ \! \dfrac{\mi \widetilde{V}^{\prime}(\xi)}{\pi} \right) \md \xi \right. \\
-&\left. \, \left(\dfrac{1}{4} \sum_{j=1}^{N+1}((\tilde{b}_{j-1})^{2} \! 
+ \! (\tilde{a}_{j})^{2}) \! - \! \dfrac{1}{8} \left(\sum_{j=1}^{N+1}
(\tilde{b}_{j-1} \! + \! \tilde{a}_{j}) \right)^{2} \right) \int_{J_{f}} 
\dfrac{\xi^{N+2}}{(\tilde{R}(\xi))^{1/2}_{+}} \left(\dfrac{2 \mi}{\pi} \left(
\dfrac{(\varkappa_{nk} \! - \! 1)}{n(\xi \! - \! \alpha_{k})} \right. \right. 
\right. \\
+&\left. \left. \left. \, \sum_{q=1}^{\mathfrak{s}-2} \dfrac{\varkappa_{nk 
\tilde{k}_{q}}}{n(\xi \! - \! \alpha_{p_{q}})} \right) \! + \! \dfrac{\mi 
\widetilde{V}^{\prime}(\xi)}{\pi} \right) \md \xi \right) \! - \! \dfrac{1}{8 
\cdot 3!} \left(\sum_{j=1}^{N+1}(\tilde{b}_{j-1} \! + \! \tilde{a}_{j}) 
\right)^{3} \\
-& \, \dfrac{1}{3!} \sum_{j=1}^{N+1}((\tilde{b}_{j-1})^{3} \! + \! 
(\tilde{a}_{j})^{3}) \! + \! \dfrac{1}{8} \left(\sum_{j=1}^{N+1}
(\tilde{b}_{j-1} \! + \! \tilde{a}_{j}) \right) \left(\sum_{m=1}^{N+1}
((\tilde{b}_{m-1})^{2} \! + \! (\tilde{a}_{m})^{2}) \right).
\end{align*}
Similarly, for the case $n \! \in \! \mathbb{N}$ and $k \! \in \! \lbrace 1,2,
\dotsc,K \rbrace$ such that $\alpha_{p_{\mathfrak{s}}} \! := \! \alpha_{k} 
\! = \! \infty$, via the representations for $\mathfrak{F}_{\infty}(z)$ given 
in Equations~\eqref{eql3.4i} and~\eqref{eql3.7b19}, one proceeds as per the 
paradigm above to show that, for the first three moments of the associated 
equilibrium measure, $\mu_{\widetilde{V}}^{\infty}$, for $q \! \in \! \lbrace 
1,2,\dotsc,\mathfrak{s} \! - \! 1 \rbrace$,
\begin{align*}
\int_{J_{\infty}} \dfrac{\psi_{\widetilde{V}}^{\infty}(\xi)}{\xi \! - \! 
\alpha_{p_{q}}} \, \md \xi =& \, \dfrac{(-1)^{\hat{\mathfrak{n}}
(\alpha_{p_{q}})}(\prod_{j=1}^{N+1} \lvert \hat{b}_{j-1} \! - \! 
\alpha_{p_{q}} \rvert \lvert \hat{a}_{j} \! - \! \alpha_{p_{q}} \rvert)^{
1/2}}{2((n \! - \! 1)K \! + \! k)/n} \int_{J_{\infty}} \dfrac{(\hat{R}
(\xi))^{-1/2}_{+}}{(\xi \! - \! \alpha_{p_{q}})} \left(\dfrac{2 \mi}{\pi} 
\sum_{q^{\prime}=1}^{\mathfrak{s}-1} \dfrac{\varkappa_{nk \tilde{k}_{
q^{\prime}}}}{n(\xi \! - \! \alpha_{p_{q^{\prime}}})} \right. \\
+&\left. \, \dfrac{\mi \widetilde{V}^{\prime}(\xi)}{\pi} \right) 
\md \xi, \\
\int_{J_{\infty}} \dfrac{\psi_{\widetilde{V}}^{\infty}(\xi)}{(\xi \! - \! 
\alpha_{p_{q}})^{2}} \, \md \xi =& \, \dfrac{(-1)^{\hat{\mathfrak{n}}
(\alpha_{p_{q}})}(\prod_{j=1}^{N+1} \lvert \hat{b}_{j-1} \! - \! 
\alpha_{p_{q}} \rvert \lvert \hat{a}_{j} \! - \! \alpha_{p_{q}} \rvert)^{
1/2}}{2((n \! - \! 1)K \! + \! k)/n} \left(\int_{J_{\infty}} \dfrac{(\hat{R}
(\xi))^{-1/2}_{+}}{(\xi \! - \! \alpha_{p_{q}})^{2}} \left(\dfrac{2 \mi}{\pi} 
\sum_{q^{\prime}=1}^{\mathfrak{s}-1} \dfrac{\varkappa_{nk \tilde{k}_{
q^{\prime}}}}{n(\xi \! - \! \alpha_{p_{q^{\prime}}})} \right. \right. \\
+&\left. \left. \, \dfrac{\mi \widetilde{V}^{\prime}(\xi)}{\pi} \right) \md 
\xi \! - \! \int_{J_{\infty}} \dfrac{(\hat{R}(\xi))^{-1/2}_{+}}{(\xi \! - \! 
\alpha_{p_{q}})} \left(\dfrac{2 \mi}{\pi} \sum_{q^{\prime}=1}^{\mathfrak{s}
-1} \dfrac{\varkappa_{nk \tilde{k}_{q^{\prime}}}}{n(\xi \! - \! 
\alpha_{p_{q^{\prime}}})} \! + \! \dfrac{\mi \widetilde{V}^{\prime}(\xi)}{\pi} 
\right) \md \xi \right. \\
\times&\left. \, \dfrac{1}{2} \sum_{m=1}^{N+1} \left(\dfrac{1}{\hat{b}_{m-1} 
\! - \! \alpha_{p_{q}}} \! + \! \dfrac{1}{\hat{a}_{m} \! - \! \alpha_{p_{q}}} 
\right) \right), \\
\int_{J_{\infty}} \dfrac{\psi_{\widetilde{V}}^{\infty}(\xi)}{(\xi \! - \! 
\alpha_{p_{q}})^{3}} \, \md \xi =& \, \dfrac{(-1)^{\hat{\mathfrak{n}}
(\alpha_{p_{q}})}(\prod_{j=1}^{N+1} \lvert \hat{b}_{j-1} \! - \! 
\alpha_{p_{q}} \rvert \lvert \hat{a}_{j} \! - \! \alpha_{p_{q}} \rvert)^{
1/2}}{2((n \! - \! 1)K \! + \! k)/n} \left(\int_{J_{\infty}} \dfrac{(\hat{R}
(\xi))^{-1/2}_{+}}{(\xi \! - \! \alpha_{p_{q}})^{3}} \left(\dfrac{2 \mi}{\pi} 
\sum_{q^{\prime}=1}^{\mathfrak{s}-1} \dfrac{\varkappa_{nk \tilde{k}_{
q^{\prime}}}}{n(\xi \! - \! \alpha_{p_{q^{\prime}}})} \right. \right. \\
+&\left. \left. \, \dfrac{\mi \widetilde{V}^{\prime}(\xi)}{\pi} \right) \md 
\xi \! - \! \int_{J_{\infty}} \dfrac{(\hat{R}(\xi))^{-1/2}_{+}}{(\xi \! - \! 
\alpha_{p_{q}})^{2}} \left(\dfrac{2 \mi}{\pi} \sum_{q^{\prime}=1}^{
\mathfrak{s}-1} \dfrac{\varkappa_{nk \tilde{k}_{q^{\prime}}}}{n(\xi \! - \! 
\alpha_{p_{q^{\prime}}})} \! + \! \dfrac{\mi \widetilde{V}^{\prime}(\xi)}{\pi} 
\right) \md \xi \right. \\
\times&\left. \, \dfrac{1}{2} \sum_{m=1}^{N+1} \left(\dfrac{1}{\hat{b}_{m-1} 
\! - \! \alpha_{p_{q}}} \! + \! \dfrac{1}{\hat{a}_{m} \! - \! \alpha_{p_{q}}} 
\right) \! + \! \int_{J_{\infty}} \dfrac{(\hat{R}(\xi))^{-1/2}_{+}}{(\xi \! - 
\! \alpha_{p_{q}})} \left(\dfrac{2 \mi}{\pi} \sum_{q^{\prime}=1}^{\mathfrak{s}
-1} \dfrac{\varkappa_{nk \tilde{k}_{q^{\prime}}}}{n(\xi \! - \! 
\alpha_{p_{q^{\prime}}})} \! + \! \dfrac{\mi \widetilde{V}^{\prime}(\xi)}{\pi} 
\right) \md \xi \right. \\
\times&\left. \, \left(\dfrac{1}{8} \left(\sum_{m=1}^{N+1} \left(\dfrac{1}{
\hat{b}_{m-1} \! - \! \alpha_{p_{q}}} \! + \! \dfrac{1}{\hat{a}_{m} \! - \! 
\alpha_{p_{q}}} \right) \right)^{2} \! - \! \dfrac{1}{4} \sum_{m=1}^{N+1} 
\left(\dfrac{1}{(\hat{b}_{m-1} \! - \! \alpha_{p_{q}})^{2}} \! + \! 
\dfrac{1}{(\hat{a}_{m} \! - \! \alpha_{p_{q}})^{2}} \right) \right) \right),
\end{align*}
where $\hat{\mathfrak{n}}(\alpha_{p_{q}})$ is given in the corresponding 
item of Remark~\ref{remext}, and
\begin{align*}
\int_{J_{\infty}} \xi \psi_{\widetilde{V}}^{\infty}(\xi) \, \md \xi =& \, 
\dfrac{1}{2((n \! - \! 1)K \! + \! k)/n} \int_{J_{\infty}} \dfrac{\xi^{N+2}}{
(\hat{R}(\xi))^{1/2}_{+}} \left(\dfrac{2 \mi}{\pi} \sum_{q=1}^{\mathfrak{s}-1} 
\dfrac{\varkappa_{nk \tilde{k}_{q}}}{n(\xi \! - \! \alpha_{p_{q}})} \! + \! 
\dfrac{\mi \widetilde{V}^{\prime}(\xi)}{\pi} \right) \md \xi \! - \! 
\dfrac{1}{2} \sum_{j=1}^{N+1}(\hat{b}_{j-1} \! + \! \hat{a}_{j}), \\
\int_{J_{\infty}} \xi^{2} \psi_{\widetilde{V}}^{\infty}(\xi) \, \md \xi =& \, 
\dfrac{1}{2((n \! - \! 1)K \! + \! k)/n} \left(\int_{J_{\infty}} \dfrac{
\xi^{N+3}}{(\hat{R}(\xi))^{1/2}_{+}} \left(\dfrac{2 \mi}{\pi} \sum_{q=1}^{
\mathfrak{s}-1} \dfrac{\varkappa_{nk \tilde{k}_{q}}}{n(\xi \! - \! 
\alpha_{p_{q}})} \! + \! \dfrac{\mi \widetilde{V}^{\prime}(\xi)}{\pi} \right) 
\md \xi \right. \\
+&\left. \, \dfrac{1}{2} \left(\sum_{j=1}^{N+1}(\hat{b}_{j-1} \! + \! 
\hat{a}_{j}) \right) \int_{J_{\infty}} \dfrac{\xi^{N+2}}{(\hat{R}(\xi))^{
1/2}_{+}} \left(\dfrac{2 \mi}{\pi} \sum_{q=1}^{\mathfrak{s}-1} \dfrac{
\varkappa_{nk \tilde{k}_{q}}}{n(\xi \! - \! \alpha_{p_{q}})} \! + \! 
\dfrac{\mi \widetilde{V}^{\prime}(\xi)}{\pi} \right) \md \xi \right) \\
-&\, \dfrac{1}{4} \sum_{j=1}^{N+1}((\hat{b}_{j-1})^{2} \! + \! 
(\hat{a}_{j})^{2}) \! + \! \dfrac{1}{8} \left(\sum_{j=1}^{N+1}
(\hat{b}_{j-1} \! + \! \hat{a}_{j}) \right)^{2}, \\
\int_{J_{\infty}} \xi^{3} \psi_{\widetilde{V}}^{\infty}(\xi) \, \md \xi =& \, 
\dfrac{1}{2((n \! - \! 1)K \! + \! k)/n} \left(\int_{J_{\infty}} \dfrac{
\xi^{N+4}}{(\hat{R}(\xi))^{1/2}_{+}} \left(\dfrac{2 \mi}{\pi} \sum_{q=1}^{
\mathfrak{s}-1} \dfrac{\varkappa_{nk \tilde{k}_{q}}}{n(\xi \! - \! 
\alpha_{p_{q}})} \! + \! \dfrac{\mi \widetilde{V}^{\prime}(\xi)}{\pi} \right) 
\md \xi \right. \\
-&\left. \, \dfrac{1}{2} \left(\sum_{j=1}^{N+1}(\hat{b}_{j-1} \! + \! 
\hat{a}_{j}) \right) \int_{J_{\infty}} \dfrac{\xi^{N+3}}{(\hat{R}(\xi))^{
1/2}_{+}} \left(\dfrac{2 \mi}{\pi} \sum_{q=1}^{\mathfrak{s}-1} \dfrac{
\varkappa_{nk \tilde{k}_{q}}}{n(\xi \! - \! \alpha_{p_{q}})} \! + \! 
\dfrac{\mi \widetilde{V}^{\prime}(\xi)}{\pi} \right) \md \xi \right. \\
-&\left. \, \left(\dfrac{1}{4} \sum_{j=1}^{N+1}((\hat{b}_{j-1})^{2} \! + \! 
(\hat{a}_{j})^{2}) \! - \! \dfrac{1}{8} \left(\sum_{j=1}^{N+1}(\hat{b}_{j-1} 
\! + \! \hat{a}_{j}) \right)^{2} \right) \int_{J_{\infty}} \dfrac{\xi^{N+
2}}{(\hat{R}(\xi))^{1/2}_{+}} \right. \\
\times&\left. \, \left(\dfrac{2 \mi}{\pi} \sum_{q=1}^{\mathfrak{s}-1} \dfrac{
\varkappa_{nk \tilde{k}_{q}}}{n(\xi \! - \! \alpha_{p_{q}})} \! + \! 
\dfrac{\mi \widetilde{V}^{\prime}(\xi)}{\pi} \right) \md \xi \right) \! 
- \! \dfrac{1}{8 \cdot 3!} \left(\sum_{j=1}^{N+1}(\hat{b}_{j-1} \! + \! 
\hat{a}_{j}) \right)^{3} \\
-& \, \dfrac{1}{3!} \sum_{j=1}^{N+1}((\hat{b}_{j-1})^{3} \! + \! 
(\hat{a}_{j})^{3}) \! + \! \dfrac{1}{8} \left(\sum_{j=1}^{N+1}(\hat{b}_{j-1} 
\! + \! \hat{a}_{j}) \right) \left(\sum_{m=1}^{N+1}((\hat{b}_{m-1})^{2} 
\! + \! (\hat{a}_{m})^{2}) \right).
\end{align*}
Since, for $\xi \! \in \! J_{\infty}$ (resp., $\xi \! \in \! J_{f})$, 
$(\hat{R}(\xi))^{1/2}_{+} \! = \! \mi (\lvert \hat{R}(\xi) \rvert)^{1/2} \! 
\in \! \mi \mathbb{R}$ (resp., $(\tilde{R}(\xi))^{1/2}_{+} \! = \! \mi 
(\lvert \tilde{R}(\xi) \rvert)^{1/2} \! \in \! \mi \mathbb{R})$, and, for $j 
\! = \! 1,2,\dotsc,N \! + \! 1$, $(\hat{R}(\xi))^{1/2} \! =_{\xi \downarrow 
\hat{b}_{j-1}} \! \mathcal{O}((\xi \! - \! \hat{b}_{j-1})^{1/2})$ and 
$(\hat{R}(\xi))^{1/2} \! =_{\xi \uparrow \hat{a}_{j}} \! \mathcal{O}
((\hat{a}_{j} \! - \! \xi)^{1/2})$ (resp., $(\tilde{R}(\xi))^{1/2} \! 
=_{\xi \downarrow \tilde{b}_{j-1}} \! \mathcal{O}((\xi \! - \! 
\tilde{b}_{j-1})^{1/2})$ and $(\tilde{R}(\xi))^{1/2} \! =_{\xi \uparrow 
\tilde{a}_{j}} \! \mathcal{O}((\tilde{a}_{j} \! - \! \xi)^{1/2}))$, that 
is, there are integrable singularities at the end-points of the intervals, 
$\hat{b}_{j-1},\hat{a}_{j}$ (resp., $\tilde{b}_{j-1},\tilde{a}_{j})$, of the 
support, $J_{\infty}$ (resp., $J_{f})$, of the associated equilibrium measure, 
$\mu_{\widetilde{V}}^{\infty}$ (resp., $\mu_{\widetilde{V}}^{f})$, it follows 
that, for $n \! \in \! \mathbb{N}$ and $k \! \in \! \lbrace 1,2,\dotsc,K 
\rbrace$ such that $\alpha_{p_{\mathfrak{s}}} \! := \! \alpha_{k} \! = \! 
\infty$ (resp., $\alpha_{p_{\mathfrak{s}}} \! := \! \alpha_{k} \! \neq \! 
\infty)$, all of the above moment integrals are real valued and bounded.

\end{document}